\documentclass{svmono}
\usepackage{amsmath}
\usepackage{amssymb}

\usepackage{amsthm}

\usepackage{color} 
 
\newcommand{\red}[1]{{\color{red} #1}}

\newcommand{\bred}[1]{{\color{red}\textbf{#1}}}

\input{infmacs.sty}

\newcommand{\Ch}{\mathop{\rm Ch}\nolimits}
\newcommand{\graph}{\mathop{\rm graph}\nolimits}
\newcommand{\heis}{\mbox{$\mathfrak{heis}$}}
\newcommand{\gf}{\mbox{$\mathfrak{gf}$}}
\newcommand{\gs}{\mbox{$\mathfrak{gs}$}}
\newcommand{\gh}{\mbox{$\mathfrak{gh}$}}
\newcommand{\cJ}{{\mathcal J}}
\renewcommand{\be}{e}

\newcommand{\bt}{{\bf t}}
\newcommand{\bq}{{\bf q}}

\renewcommand{\dd}{d} 
\renewcommand{\mlabel}{\label} 
\DeclareMathOperator{\Emb}{Emb}

\DeclareMathOperator{\lcx}{lcx}
\DeclareMathOperator{\DL}{dl}
\DeclareMathOperator{\fbx}{fb}

\usepackage{multicol}

\begin{document}

\frontmatter

\begin{titlepage}
  \begin{center}
    {\huge\bf  Infinite-Dimensional Lie Groups} \\[2cm]
    {\Large\bf  Karl-Hermann Neeb, Helge Gl\"ockner} \\[1cm]
{\Large\bf Version of January 2, 2026} \\[1.5cm]
  \end{center}

\author{Karl-Hermann Neeb, Helge Gl\"ockner}

\newpage

\end{titlepage}

\pagenumbering{roman}

\tableofcontents

\newpage

\setcounter{chapter}{-1}
\mainmatter

\chapter{Introduction} \mlabel{ch:intro} 

Symmetries play a central role in the natural sciences and 
throughout mathematics. In mathematics, abstract 
symmetries are represented by groups. If symmetries 
are viewed dynamically, by the development 
of a system with respect to a time-parameter, 
then they correspond to one-parameter groups. 
Dealing with symmetries depending on many parameters leads to 
the concept of a Lie group. 
The core topic of this monograph 
is the theory of infinite-dimensional Lie groups, i.e., 
of symmetries depending on infinitely many
parameters. Such symmetries may be studied on an
infinitesimal, local or global level, which amounts 
to studying Lie algebras, local Lie groups
or global Lie groups, respectively.

Finite-dimensional Lie theory was created in the late
19th century by 
\textsc{Marius Sophus Lie} and \textsc{Friedrich Engel}, who showed that in
finite dimensions the local and the infinitesimal 
theory are equivalent (\cite{Lie80, Lie95, LE93}). 
The global theory had to wait until the early 20th century for its language 
of topological spaces and smooth manifolds to be properly developed. 
A crucial feature of the finite-dimensional theory
is that finiteness conditions permit a full-fledged
structure theory of finite-dimensional Lie groups to 
be developed in terms
of the Levi splitting and the fine structure of semisimple 
Lie groups. We refer to \cite{HiNe12} for the theory
  of finite-dimensional Lie groups.

In infinite dimensions,
the passage from the infinitesimal to the local
and from the local to the global level is not always possible, 
so that the theory splits
into three properly distinct levels.
A substantial part of the literature
on infinite-dimensional Lie
theory deals exclusively with Lie algebras (the infinitesimal level), 
their structure, and their
representations. 
However, only special classes of groups, such as Kac--Moody groups  
or certain direct limit groups, can
be approached by purely algebraic methods, combined with 
finite-dimensional Lie theory 
\cite{Kum02, Ma18}. 
The algebraic side of Lie theory is highly 
relevant for many applications in Mathematical Physics, 
where the infinitesimal approach 
is convenient for calculations, but a global perspective is required 
to understand geometric and topological phenomena.

In this monograph, we focus on
the local and global level of infinite-dimensional Lie theory,
as well as the mechanisms  required to pass from 
the infinitesimal to the local and from there to the global level. 
Our approach is  based on a notion
of a Lie group which is both simple and general:
\begin{center}
A {\bf Lie group}  is a smooth manifold, endowed with 
a group structure
such that multiplication and inversion are smooth maps. 
\end{center}
The main difference compared
to the finite-dimensional theory concerns
the notion of a manifold: The manifolds we consider
are modeled on (not necessarily complete) locally convex spaces.
It is quite useful to approach Lie groups
from such a general perspective,
because this enables a unified discussion of
all basic aspects of the theory.
Although we simply call them 
{\it Lie groups}, a more specific terminology is {\it locally convex Lie groups}. Depending on the type 
of the model space, we obtain in particular the classes of 
finite-dimensional, Banach--, Fr\'echet--, LF-- and Silva--Lie groups. 

There are weaker concepts of Lie groups and  infinite-dimensional manifolds. 
One is based on the ``convenient setting'' for global analysis developed 
by \textsc{Fr\"olicher, Kriegl and Michor} (\cite{FB66, Mi84, FK88, KM97}). 
In the context of Fr\'echet 
manifolds, this setting does not differ from ours, but for more general 
model spaces it provides a concept of a smooth map which does not imply 
continuity, hence leads to ``Lie groups'' which are not topological groups. 
Another approach is based on \textsc{J.-M.~Souriau}'s 
concept of a diffeological space 
which can be used to study spaces like quotients of $\R$ by 
non-discrete subgroups in a differential geometric context 
\cite{So84, So85, DI85, Lk92} (see \cite{HeMa02}
for applications to diffeomorphism groups). It has the 
advantage that the category of diffeological spaces is cartesian closed 
and that any quotient of a diffeological space carries a natural diffeology, 
but, as a consequence, diffeological structures are rather weak.
This incredible freedom makes it
harder to distinguish regular objects from non-regular ones.

In this context, we also mention the more recent concept
of a half--Lie group (\cite{KMR16}, \cite{MNe18}, \cite{BHM23}), which  consists
of a topological group with a smooth Banach manifold structure for which only
left multiplications are smooth.
This subsumes  the groups  of $C^k$-diffeomorphisms playing 
an important role in geometric analysis (cf.\ \cite{AK98}, \cite{EMi99}).
Its main feature is that it permits to keep track of quantitative
information and to provide a context for perturbation theory.
\\ 

To demonstrate  the diversity of infinite-dimensional 
Lie groups, let us briefly discuss 
several types of infinite-dimensional 
Lie groups that motivate the development of the corresponding calculus for 
related classes of locally convex spaces: 
\begin{itemize}
\item[$\bullet$] Linear Lie groups 
\item[$\bullet$] Mapping groups
\item[$\bullet$] Direct limit groups 
\item[$\bullet$] Diffeomorphism groups 
\end{itemize}

\nin {\bf Linear Lie groups.} 
In finite-dimensional Lie theory, a natural approach to Lie groups proceeds via 
matrix groups, i.e., subgroups of the group 
$\GL_n(\R)$ of invertible real $n \times n$-matrices. 
Since every finite-dimensional (associative) algebra can be embedded into a matrix algebra, 
this is equivalent to considering subgroups of the unit groups 
\[ \cA^\times := \{ a \in \cA \: (\exists b \in \cA)\ ab = ba = \1\}\] 
of 
finite-dimensional unital associative algebras~$\cA$. The advantage of this approach is 
that 
one can define the exponential function quite directly via the 
exponential series and thus take a shortcut 
to several deeper results on Lie groups \cite{HiNe12}. 
This approach also works quite well in the context of Banach--Lie groups. 
Here the linear Lie groups are subgroups of unit groups of Banach algebras, 
but this setting is too restrictive for many applications 
of infinite-dimensional Lie theory. 

Dealing with locally convex spaces, it is natural to consider 
for a locally convex space $V$ the unital algebra 
$\cA := {\cal L}(V)$ of continuous linear endomorphisms of~$V$. 
Its  unit group $\GL(V)$ carries no 
natural Lie group structure if $V$ is not a Banach space.  
In particular, it is not open in any vector  topology 
(cf.\ \cite{Ms63}): 
If $A \in \cL(V)$ has unbounded spectrum,  then there exists a sequence 
$t_n \to 0$ for which no operator $\1+ t_n A$ is contained in $\GL(V)$. 
One therefore needs a class of algebras which behaves better than~$\cL(V)$. The most natural class of algebras for 
infinite-dimensional Lie theory are 
{\it continuous inverse algebras} (cias). 
These are unital locally convex algebras $\cA$ with continuous multiplication 
such that the unit group $\cA^\times$ 
is open and the inversion is a continuous map $\cA^\times \to \cA$. 
Typical examples of cias are: 
\begin{itemize}
\item[$\bullet$] Unital Banach algebras:
The convergence of the {\it Neumann series} \break 
$(\1 - x)^{-1} = \sum_{k = 0}^\infty x^k$ for $\|x\| < 1$ 
implies that $\cA^\times$ is open and that inversion is continuous.
\item[$\bullet$]   For a compact  
finite-dimensional smooth manifold $M$, 
the  algebra $C^\infty(M,\C)$ of smooth functions 
is a unital Fr\'echet continuous inverse algebra.
If $M$ is non-compact and $\sigma$-compact, then
$C^\infty_c(M,\C)$ is a non-metrizable LF space and it has no algebra unit.
The corresponding unital extension is a continuous inverse algebra. 
\item[$\bullet$] For a compact subset $K \subeq \C^n$, the algebra 
$\cO(K)$ of germs of holomorphic functions defined on a neighborhood of $K$ 
is a cia which is a Silva space, i.e., a direct limit of Banach algebras 
with compact connecting maps. 
\end{itemize}
The unit group $\cA^\times$ of a cia $\cA$ is a Lie group  
when endowed with its natural manifold structures as an open subset; 
this applies in particular to the unit group 
$\GL_n(\cA)$ of $M_n(\cA)$. 
We think of ``Lie subgroups'' of these 
groups as {\it linear Lie groups}. 
Most classical Lie groups are defined as centralizers of certain 
matrices or as the set of fixed points for a group of automorphisms. All 
these constructions have natural generalizations to matrices with 
entries in cias. \\

\nin {\bf Mapping groups.} 
In the context of Banach space calculus, one constructs Lie groups of mappings as follows. 
For a compact space $X$ and a Banach--Lie group $K$, the group 
$C(X,K)$ of continuous maps is a Banach--Lie group with Lie algebra 
$C(X,\fk)$, where $\fk := \Lie(K)$ is the Lie algebra of $K$. 

In the larger context of locally convex Lie groups, one also obtains 
for each Lie group $K$ with Lie algebra $\fk$ 
and a $\sigma$-compact smooth finite-dimensional 
manifold $M$, a Lie group structure on the group 
$C_c^\infty(M,K)$ of smooth compactly supported maps from $M$ to $K$ 
such that $C^\infty_c(M,\fk)$ is its Lie algebra. 
The latter is Fr\'echet if $M$ is compact and $K$ is finite-dimensional, 
but if $M$ is non-compact, then the model space is only LF. 

The passage from continuous maps to smooth maps 
is motivated in particular by the existence of derivations on algebras 
of smooth functions. In particular, we can form semidirect products 
of Lie groups such as 
\[ C^\infty(M,K) \rtimes \Diff(M),\] 
which acts by bundle automorphisms on the trivial bundle $M \times K \to M$. 
There is no Lie group analog of this in the $C^0$-context of Banach--Lie groups. \\

\nin {\bf Direct limit groups.} 
One way to obtain infinite-dimensional 
groups from finite-dimensional ones is to consider 
a sequence $(G_n)_{n \in\N}$ of 
finite-dimensional Lie groups and morphisms 
$\phi_n \: G_n \to G_{n+1}$, so that we can define a direct limit group 
$G := \indlim G_n$ whose representations correspond to compatible 
sequences of representations of the groups $G_n$. 
We shall see that the direct limit group $G$ can always be 
endowed with a natural Lie group 
structure (Chapter~\ref{ch:dirlim}). Its Lie algebra $\L(G)$ is the 
countably-dimensional 
direct limit space $\indlim \L(G_n)$, endowed with the finest locally convex 
topology. This provides an interesting class of infinite-dimensional Lie groups
which is still quite close to finite-dimensional groups and has a very rich 
representation theory (\cite{DiPe99, NRW01, Wol05}). \\

\nin {\bf Groups of diffeomorphisms.} 
Many infinite-dimensional 
groups arise naturally from geometric or other structures on 
manifolds as their automorphism groups. 
In \textsc{Felix Klein}'s Erlangen Program \cite{Kl1872}, 
geometric structures are even defined in terms of their automorphism groups. 
Starting with the differentiable structure  on a 
$\sigma$-compact finite-dimensional smooth manifold~$M$, the automorphism group 
is the full diffeomorphism group $\Diff(M)$. To endow this group 
with a manifold structure, let us first assume that $M$ is compact. 
Then $\Diff(M)$ has a natural Lie group structure modeled on the space 
${\cal V}(M)$ of (smooth) vector fields on $M$, which is 
the Lie algebra of this group 
(\cite{Les67, Omo70, EM69, EM70, Gu77, Mr80, Ham82}). 
Then smooth left actions 
$G \times M \to M$ correspond to Lie group homomorphisms 
$G \to \Diff(M)$. For $G = \R$, we obtain in particular the 
well-known correspondence between 
smooth flows on $M$, smooth vector fields on $M$, 
and one-parameter groups of $\Diff(M)$. 
Other important groups of diffeomorphisms arise as subgroups of $\Diff(M)$ 
stabilizing a volume form $\mu$, a symplectic structure $\omega$, or a contact 
form $\alpha$ (cf.\ \cite{KM97}). 

If $M$ is not compact, then 
it is still possible to turn $\Diff(M)$ into a Lie group, but then 
it has to be modeled on the Lie algebra 
${\cal V}_c(M)$ of smooth vector fields with compact support, 
endowed with its natural LF structure 
(cf.\ \cite{Mr80, Mil82, Gl03c}). 
This has the disadvantage that a smooth flow generated by a vector 
field with non-compact support defines a discontinuous homomorphism 
$\R \to \Diff(M)$.  
For this Lie group structure, the normal subgroup 
$\Diff_c(M)$ of all diffeomorphisms which coincide with $\id_M$ outside 
a compact set is an open subgroup of $\Diff(M)$. 

Of a different nature, but also locally convex Lie groups, are 
groups of formal diffeomorphism as studied by 
\textsc{Lewis} \cite{Lew39}, \textsc{Sternberg} \cite{St61} 
and \textsc{Kuranishi} \cite{Kur59}, 
groups of germs of smooth and analytic diffeomorphisms of $\R^n$ 
fixing $0$ (\cite{RK97, Rob02}), 
and also germs of biholomorphic maps of $\C^n$ fixing $0$ 
(\cite{Pis76, Pis77, Pis79}). \\ 

\nin {\bf Similarity with the infinite-dimensional unitary group.} The 
situation for non-compact manifolds is similar to the situation we encounter 
in the theory of unitary group representations. Let $\cH$ be an 
infinite-dimensional complex Hilbert space and 
$\U(\cH)$ be its unitary group. This group has two natural topologies. 
The norm topology on $\U(\cH)$ inherited from the Banach algebra ${\cal L}(\cH)$ turns it 
into a Banach--Lie group $\U(\cH)_n$, but this topology is 
rather fine. The strong operator topology (the topology of pointwise convergence) 
turns $\U(\cH)$ into a topological group $\U(\cH)_s$, whose one-parameter 
groups are generated by (unbounded) selfadjoint operators. 
There is a striking analogy with diffeomorphism groups if $M$ is non-compact, 
where the ``non-Lie group'' $\Diff(M)$ corresponds to 
the topological group $\U(\cH)_s$, 
smooth vector fields on $M$ correspond to selfadjoint operators, 
and compactness of the support corresponds to boundedness. 
Accordingly, the Lie group structure on $\Diff_c(M)$ compares to the Lie 
group structure on $\U(\cH)_n$. \\

\nin{\bf Diffeomorphisms in infinite dimensions.} 
We have already seen above that, for a non-compact smooth manifold $M$, 
its diffeomorphism group $\Diff(M)$ carries no natural Lie group  structure. 
However, even if $M$ is an infinite-dimensional  
manifold, we can still make 
sense of ``smooth'' maps $f \: N \to \Diff(M)$, where $N$ is a smooth manifold, 
by requiring  the corresponding map 
\[  N \times M \to M^2, \quad (n,m) \mapsto (f(n)(m), f(n)^{-1}(m)) \] 
to be smooth. Then a smooth action of a Lie group $G$ on $M$ 
is a {\it smooth homomorphism} $G \to \Diff(M)$. 
This trick of ``moving arguments to the left'' underlies many constructions 
in the convenient calculus \cite{KM97}. 
We  find it useful in 
numerous contexts, 
as a tool  to apply rudiments of Lie theoretic arguments to 
``non--Lie groups'' such as $\Diff(M)$, or to $\GL(V)$ for a locally convex space. \\

\nin {\bf Content of the book.} 
The discussion of these classes of Lie groups shows that 
locally convex Lie groups occur in quite different types: Banach--Lie groups, 
groups of diffeomorphisms (modeled on Fr\'echet and LF spaces), 
groups of germs (modeled on Silva spaces) and formal groups (modeled 
on Fr\'echet spaces such as $\R^\N$). This diversity has to  be taken 
into account in the development of the corresponding analytical 
tools. Accordingly, this monograph is divided into two parts, 
Part A deals with analysis on locally convex spaces 
and Part B with infinite-dimensional Lie groups. 

{\bf Part A} consists of four chapters. The first two introduce 
basic and more specific aspects of differential calculus, 
the third one introduces manifolds modeled on locally convex spaces and 
the fourth chapter  describes constructions such as manifolds of mappings, 
box products and direct limits. 

In Chapter~\ref{chapcalcul}, we explain our setting of
differential calculus in locally convex spaces~$E$. 
The approach is a very natural one: Partial derivatives,
as familiar from finite-dimensional calculus, are replaced by directional
derivatives, which are required to be continuous both in the point and the direction.
Let $E$ and $F$ be locally convex spaces, $U
\subeq E$ open and $f \: U \to F$ a map. Then the {\it derivative
  of $f$ at $x$ in the direction $h$} is defined as 
\[  df(x,h) := (D_h f)(x) := \derat0 f(x + t h) 
= \lim_{t \to 0} \frac{1}{t}(f(x+th) -f(x)) \] 
whenever it exists. The function $f$ is called {\it differentiable at
  $x$} if $df(x,h)$ exists for all $h \in E$. It is called {\it
  continuously differentiable} (or $C^1$), if it is differentiable at all
points of $U$ and $df  \: U \times E \to F$ 
is a continuous map. 
We call $f$ a {\it $C^k$-map}, $k \in \N \cup \{\infty\}$, 
if it is continuous, the iterated directional derivatives 
\[  d^{\,(j)}f(x,h_1,\ldots, h_j)
:= (D_{h_j} \cdots D_{h_1}f)(x) \] 
exist for all integers $1\leq j \leq k$, $x \in U$ and $h_1,\ldots, h_j \in E$, 
and all maps $d^{\,(j)} f \: U \times E^j \to F$ are continuous. 
As usual, $C^\infty$-maps are called {\it smooth}. 
This approach to differential calculus in locally convex spaces
goes back to \textsc{Andr\'{e}e Bastiani}~\cite{Ba64}.

We neither assume completeness of $E$ of $F$, nor do we restrict 
to mappings on open subsets, but develop differential calculus
for mappings on ``locally convex'' subsets with dense interior, 
where differentiability means that the maps $d^{\,(j)}f$ extend continuously 
to the boundary.
This enables us to speak about smooth maps on compact intervals,
squares and simplices. 
We encounter analogs of many classical results of differential calculus
in the infinite-dimensional setting, like the 
Fundamental Theorem of Calculus,
the Chain Rule, and Taylor's Theorem.
As a prerequisite for many important constructions in infinite-dimensional Lie theory, 
we also discuss the continuity and differentiability
properties of typical mappings between
spaces of $C^k$-functions.

In Chapter~\ref{chapfurcalc}, 
we discuss further concepts and tools of infinite-dimensional calculus
which are relevant for specific examples and applications, such as 
real and complex analytic mappings, ordinary differential equations 
and specifics of calculus on metrizable spaces, Silva spaces 
and locally convex direct sums. 

Chapters~\ref{chapcalcul} and \ref{chapfurcalc} may also be useful 
in other branches of mathematics, like dynamical systems or geometry, 
irrespective of Lie theory. They lay a broad foundation for applications 
of infinite-dimensional differential calculus. 

In Chapter~\ref{chapmanif}  we turn to manifolds modeled
on locally convex spaces and the basic concepts going along with them:
smooth (and $C^r$-) maps between manifolds, tangent maps,
the Lie algebra of smooth vector fields, and differential forms. 
We also prepare  basic constructions of infinite-dimensional
manifolds by discussing vector bundles, principal bundles, sprays and 
local additions. 

In Chapter~\ref{chap-manifold-constructions}
  we construct major classes of infinite-dimensional manifolds. 
Notably, for each compact smooth manifold $M$
and paracompact finite-dimensional smooth manifold~$N$,
we turn the set $C^k(M,N)$
of $N$-valued $C^k$-maps on~$M$
into a smooth manifold.
More generally, $N$ can be any smooth manifold
admitting a local addition.
The construction is essential for
infinite-dimensional Lie theory:
First, for each Lie group~$K$,
we obtain a Lie group structure on $C^k(M,K)$. 
Second, the group $\Diff(M)$
of smooth diffeomorphisms of~$M$
turns out to be an open subset of $C^\infty(M,M)$,
and the smooth manifold structure as an open
subset makes it a Lie group. 
As a starting point,
we define a topology on $C^k(M,N)$
for arbitrary $C^k$-manifolds $M$ and~$N$
modeled on locally convex spaces. 
If $N:=F$ is a locally convex space, the latter topology makes $C^k(M,F)$
a locally convex space. Similarly,
taking $N:=E$ for a given vector bundle $E\to M$,
the induced topology on the subset $\Gamma_{C^k}(E)\sub C^k(M,E)$
of $C^k$-sections is  a locally convex vector topology.
We also study differentiability properties of non-linear mappings between such spaces of mappings or spaces of sections. 
Notably, we obtain exponential laws
for function spaces on products of manifolds. 
On this foundation, we carry out the following constructions:\medskip

\nin \emph{Fine box products of manifolds.}
For an arbitrary sequence $(M_n)_{n\in \N}$
of smooth manifolds, we construct a smooth manifold structure on the cartesian
product
\[
\prod_{n\in \N}M_n,
\]
endowed with a certain topology which is finer than the product
topology.  Such manifolds
are called \emph{fine box products}. They are denoted 
${\prod^{\fbx}_{n\in \N}} M_n$. 
Note that the product topology on $\T^\N$ is not compatible with any manifold 
structure because this space is not locally contractible. 
\\ 

\noindent
\emph{Manifolds of mappings on non-compact manifolds.}
For every $\sigma$-compact finite-dimensional
smooth manifold $M$ and smooth manifold $N$ admitting a local addition,
we show that the image of the map
\[\rho\colon C^k(M,N)\to{\prod_{n\in \N}}^{\fbx} C^k(M_n,N),\quad f\mto (f|_{M_n})_{n\in \N}
\]
is a submanifold of the fine box product, for each locally finite cover $(M_n)_{n\in\N}$
of $M$ by compact full submanifolds~$M_n$. We give $C^k(M,N)$
the smooth manifold structure making $\rho$ a $C^\infty$-diffeomorphism
onto the image. \\ 

\noindent
\emph{Direct limits of finite-dimensional manifolds.}
We construct a natural smooth manifold structure on
the union $\bigcup_{n\in \N}M_n$ for each
ascending sequence
\[
M_1\sub M_2\sub\cdots
\]
of finite-dimensional $C^\infty$-manifolds such that all inclusion maps
$M_n\to M_{n+1}$ are smooth immersions.

{\bf Part B} on Lie groups consists of eleven  chapters. 
Chapter~\ref{ch:3}  introduces 
locally convex Lie groups and their basic structural features. 
The other chapters deal with important properties 
of Lie groups, such as regularity (Chapter~\ref{ch:4}) and local exponentiality 
(Chapter~\ref{ch:5}). This is followed by Chapter~\ref{ch:7} dealing with 
local Lie theory and Chapter~\ref{ch:6} on subgroups and quotients. 
Chapters~\ref{ch:lingrp}--\ref{ch:diffeo}
 deal with the four important classes of Lie groups 
introduced above: 
linear Lie groups, mapping groups, direct limit Lie groups and 
diffeomorphism groups.
Homotopy groups of various classes of infinite-dimensional
Lie groups are studied in Chapter~\ref{ch:top}. 
The final Chapter~\ref{ch:sectop} is a collection of 
smaller sections on more specialized topics. 

As in finite dimensions, the Lie algebra $\g = \L(G)$ of a Lie group 
$G$ is identified with the tangent space $T_\be(G)$ in the identity $\be$, 
where the Lie bracket is obtained by identification with the space of left invariant vector fields. 
This turns $\L(G)$ into a locally convex topological Lie algebra. 
Associating to a morphism $\phi$ of Lie groups its tangent map 
$\L(\phi) := T_\be(\phi)$, 
we obtain the {\it Lie functor $\L$} from the category of locally convex Lie 
groups to the category of locally convex topological Lie algebras. 
The core of Lie theory now consists in determining how much information the Lie 
functor forgets and how much information on the group level
  can be reconstructed from it.

As a general rule, all constructions from finite-dimensional Lie theory 
that require only ``differentiation'' carry over to the infinite-dimensional 
context, but as soon as integration processes and solutions of differential 
equations are involved, one has to refine the axiomatic setup. 
For example, an  important tool in the finite-dimensional and Banach 
context is the exponential map, 
but as vector fields on locally convex manifolds need not possess integral 
curves, there is no general theorem that guarantees the existence of a smooth 
exponential map, i.e., a smooth function 
\[ \exp_G \: \L(G) \to G\] 
for which the curves $\gamma_x(t) := \exp_G(tx)$ are homomorphisms 
$(\R,+) \to G$ with $\gamma_x'(0)= x$. Therefore the existence of an 
exponential function has to be treated as an additional requirement. 
Even stronger 
is the requirement of {\it regularity}, meaning that, for each 
smooth curve $\xi \: [0,1] \to \L(G)$, the initial value problem 
\[  \gamma'(t) =\gamma(t).\xi(t) := T_\be(\lambda_{\gamma(t)})\xi(t), 
\quad \gamma(0) = \be \]
has a solution $\gamma_\xi \:[0,1] \to G$ and that 
$\gamma_\xi(1)$ depends smoothly on $\xi$. 
Regularity, discussed in Chapter~\ref{ch:4}, 
 is a natural assumption that provides a good deal of methods to 
pass from the infinitesimal to the global level. This regularity 
concept is due to \textsc{Milnor} \cite{Mil84}. 
It weakens the $\mu$-regularity 
introduced by \textsc{Omori} et al. (see \cite{KYM85} for a survey), 
but it is still strong enough for the essential 
Lie-theoretic applications.  
Presently, we do not know of any Lie group modeled on a complete space 
which is not regular. For all major concrete classes discussed below,  
one can prove regularity, but there is 
no general theorem asserting that each locally convex Lie group with a complete 
model space is regular or merely  that it has an exponential function. 
To prove or disprove such a theorem is a fundamental open problem of the theory. 

An assumption of a different nature than regularity, 
and which can be used to develop a profound 
Lie theory, is that $G$ is {\it locally exponential} in the sense that 
it has an exponential function which is a local diffeomorphism in $0$. 
Groups with this property are studied in Chapter~\ref{ch:5}. 
Even stronger is the assumption that $G$ is analytic and that 
the exponential function is an analytic local diffeomorphism in $0$. 
Groups with this property are called {\it BCH--Lie groups}, because the local 
multiplication in canonical local coordinates is given by the 
Baker--Campbell--Hausdorff (BCH) series 
\[ x * y = x + y + \frac{1}{2}[x,y]+\frac{1}{12}[x,[x,y]]+\frac{1}{12}[y,[y,x]]
+\ldots \] 
This class contains in particular all Banach--Lie groups. 

For any locally exponential Lie group, much information on the group 
$G$ is already contained in the local multiplication defined by 
the BCH series on a suitable open $0$-neighborhood $U \subeq \g$. 
The so-obtained structure is a {\it locally exponential local Lie group} 
because the straight lines $(-\eps,\eps) \to U, t \mapsto tx$ are 
local one-parameter groups. This motivates the discussion of local 
Lie groups and the particular class associated to 
{\it locally exponential Lie algebras} in Chapter~\ref{ch:7}. 

In Chapter~\ref{ch:6} we turn to subgroups of Lie groups. This is an area
where  even finite-dimensional Lie theory becomes quite subtle. 
In a finite-dimensional Lie group $G$, one basically considers two types 
of subgroups: The best-behaved class are 
{\it Lie subgroups}. These are the closed subgroups $H \subeq G$; they 
are automatically embedded submanifolds and the coset space $G/H$ carries a 
natural manifold structure. The {\it integral subgroups} 
$H = \la \exp \fh \ra$ are generated by  the exponential image of a Lie 
subalgebra $\fh \subeq \g$. These subgroups also carry a canonical 
Lie group structure for which the inclusion $\iota_H \: H \to G$ 
is an immersion, but in general not an embedding. 
A~unifying picture is obtained by the theorem asserting that 
{\bf every} subgroup~$H$ of a finite-dimensional Lie group 
carries a unique Lie group structure for which the inclusion $\iota_H$ 
defines an {\it initial submanifold}. By  Yamabe's Theorem,  
ar\-cwise connectedness of $H$ as a topological subspace of $G$ 
is equivalent to connectedness of the corresponding Lie group~$H$ 
(\cite{HiNe12}).

Most of these extremely strong results break down in the infinite 
dimensional context, first of all Yamabe's Theorem, which heavily 
rests on Brouwer's Fixed Point Theorem. 
There are subgroups $H$ of infinite-di\-men\-sio\-nal Lie groups $G$ 
which are arcwise connected, but in which all smooth arcs are constant. 
They also carry a canonical initial manifold structure, but it corresponds 
to the discrete topology. Likewise, the closedness of a subgroup 
does in general not imply the existence of a submanifold structure. 
These defects are typical in the sense that, in the infinite-dimensional 
context, topological requirements are not enough and one has to 
consider properties taking the smooth structure into account. 
This leads to the categorical concept of an {\it initial Lie subgroup} 
and the notion of a {\it split Lie subgroup}, which specifies those 
subgroups $H$ for which $G/H$ carries a natural manifold structure 
and right multiplication turns $G$ into an $H$-principal bundle. 

Besides this ``bad news'', there 
is also some good news, namely that many results on subgroups of 
Banach--Lie groups carry over to locally exponential Lie groups. 
First of all, locally exponential subalgebras $\fh$ of the 
Lie algebra~$\g$ of a locally exponential 
Lie group generate integral subgroups (Integral Subgroup 
Theorem), and good criteria are available for integral 
subgroups to be initial. The finite-dimensional 
result that closed subgroups are Lie groups survives in the form 
that locally compact subgroups of locally exponential Lie groups 
are Lie subgroups; but local compactness implies that these
subgroups are finite-dimensional. 

We do not go deeper into 
Chapters~\ref{ch:lingrp}-\ref{ch:diffeo}
on linear Lie groups, mapping groups, direct limit Lie groups and 
diffeomorphism groups.

Chapter~\ref{ch:top} deals with the topology of infinite-dimensional Lie 
groups. For a Lie group $G$, the most important topological information 
is contained 
in the first three homotopy groups  $\pi_0(G)$ (the group of connected components), 
$\pi_1(G)$ (the fundamental group), and the second homotopy group $\pi_2(G)$. 
The importance of $\pi_0(G)$ is clear 
because one often needs to know whether a concretely 
given group is connected and the influence of the Lie algebra 
does not reach beyond the identity component. 
Information on the fundamental group 
is important for the integration of Lie 
algebra homomorphisms to group homomorphisms and hence in particular for 
representation theory (Theorem~\ref{thm3.2.11}). 
It also shows up in the integration of 
$1$-forms satisfying the Maurer--Cartan equation (Theorem~\ref{thm-fundamental}). 
The interest in $\pi_2(G)$ stems 
from the crucial role this group plays for enlargeability of Lie algebras and for extensions of $G$ (Subsection~\ref{subsec:per-grp}). 
After introducing higher homotopy groups and some techniques 
for their calculation, such as a suitable long exact homotopy sequence, 
we briefly recall some of the key 
results on homotopy groups of finite-dimensional Lie groups 
and groups of diffeomorphisms. 
The core part of this chapter deals with 
various types of groups of operators on Hilbert space, 
direct limit groups and mapping groups. \\

The long final Chapter~\ref{ch:sectop} is a collection of 
smaller sections dealing with more specialized topics. 
Without going into detail, let us highlight one topic of particular 
importance that is treated in this chapter, namely the existence 
of Lie groups  for a given Lie algebra. 
We call a locally exponential Lie algebra 
$\g$ {\it enlargeable} if there exists a locally exponential 
Lie group $G$ with $\L(G) = \g$. 
Every finite-dimensional Lie algebra is enlargeable 
by Lie's Third Theorem 
and examples by \textsc{van Est} and \textsc{Korthagen} showed 
in the 1960s that not every Banach--Lie algebra $\g$ is enlargeable. 
The obstruction lies in the non-discreteness of a subgroup 
$\Gamma(\g)$  of the center $\fz(\g)$ (the {\it period group}). 
This criterion is generalized to locally exponential Lie groups 
in Section~\ref{sec:11.4}. 
As local exponentiality of a Lie algebra $\g$ already implies the 
existence of a local group $(U,*) \subeq \g$ with Lie algebra~$\g$, 
this result deals with the local to global passage 
 which can be addressed with well-known constructions from topology 
using suitable path and loop groups. 
It is much harder to develop criteria for more general locally 
convex Lie algebras to correspond to a local or global Lie group. \\

To make the book more self-contained, we include long appendices on 
point set topology (Appendix A) and locally convex spaces (Appendix B),
a short one on topological groups (Appendix C) and Appendix D 
on ``smooth maps with values in non-Lie groups'', dealing with maps 
with values in diffeomorphism groups of infinite-dimensional 
manifolds. \\

\nin {\bf Some history.} 
To put the theory of locally convex Lie groups in perspective, we take a 
brief look at the historical development of infinite-dimensional Lie theory. 
Infinite-dimensional Lie algebras, such as 
Lie algebras of vector fields, where present in Lie theory right 
from the beginning, when \textsc{Sophus Lie} started to study (local) Lie groups 
as groups ``generated'' by finite-dimensional Lie algebras 
of vector fields \cite{Lie80}. 
The global theory 
of finite-dimensional Lie groups started to develop in the 
late 19th century, driven substantially by \textsc{\'E.~Cartan}'s work 
on symmetric spaces \cite{CaE98}. 
The first exposition of a global theory, 
including the description of all connected groups with a given Lie algebra 
and integral subgroups, was given by \textsc{Mayer} and \textsc{Thomas} 
\cite{MaTh35}. After the combination 
with the structure theory of Lie algebras, the theory reached 
its mature form in the middle of the 20th century, which is exposed in 
the fundamental books of \textsc{Chevalley} \cite{Ch46})
 and \textsc{Hochschild} \cite{Ho65} (see also \textsc{Pontrjagin} \cite{Po39} 
for an early textbook situated on the borderline between topological groups 
and Lie groups). 

Already in the work of \textsc{Sophus Lie},  
infinite-dimensional groups occur as groups of 
(local) diffeomorphisms of open domains in $\R^n$ \cite{Lie95}.  
Later, \textsc{\'E.~Cartan} undertook a more systematic study of certain types of infinite-dimensional 
Lie algebras, resp., groups of diffeomorphisms preserving geometric structures 
on a manifold, such as symplectic, contact or volume forms \cite{CaE04}. 
The advent of Quantum Mechanics in the 1920s created 
a need to understand the structure of 
groups of operators on Hilbert spaces, which is a quite different class of 
infinite-dimensional groups (cf.\ \cite{De32}). 

The first steps toward a theory of 
 infinite-dimensional groups as smooth manifolds    
were undertaken by \textsc{Birkhoff} in \cite{Bir36, Bir38}, 
where he developed the local 
Lie theory of Banach--Lie groups, resp., Banach--Lie algebras (see also 
\cite{MiE37} 
for first steps in extending Lie's theory of local transformation groups to the 
Banach setting). 
In particular, he proved that local $C^1$-Banach--Lie groups 
admit exponential coordinates, which lead to analytic local 
Lie group structures, 
and showed that continuous homomorphisms are analytic and  
that, for every Banach--Lie 
algebra, the BCH series defines an analytic local group structure. 
He also defined the Lie algebra of a local group, 
derived its functoriality properties and established the correspondence 
between closed subalgebras/ideals and the corresponding local subgroups. 
Even product integrals, which play a central role in the modern theory, 
appear in his work as solutions of left invariant differential equations. 
The local theory of Banach--Lie groups was continued 
by \textsc{Dynkin} \cite{Dy47, Dy53} who developed the algebraic theory of the 
BCH series further and by \textsc{Laugwitz} 
\cite{Lau55, Lau56} who developed 
a differential-geometric perspective, which is quite close in spirit to the 
theory of locally exponential Lie groups. 
Put in modern terms, he used the Maurer--Cartan form and integrability conditions 
on (partial) differential equations on Banach spaces, developed by \textsc{Michal}
and \textsc{Elconin} \cite{MiE37, MicA48} 
to derive the existence of the local group structure 
from the Maurer--Cartan form, which in turn is obtained from the Lie bracket. 
In the finite-dimensional case, this strategy is due to \textsc{F.~Schur} 
\cite{SchF90a} and quite close to \textsc{Lie}'s original approach. 
In \cite{Lau55}, \textsc{Laugwitz} showed in particular that 
the center and any locally compact subgroup of a Banach--Lie group 
are Banach--Lie subgroups. 
Formal Lie groups in  infinitely many parameters were introduced by \textsc{Ritt} a few years 
earlier \cite{Ri50}. 

The global theory of Banach--Lie groups started in the early 1960s with 
\textsc{Maissen}'s paper 
\cite{Ms62} which contains the first basic results on the Lie functor 
on the global level, such as the existence of integral subgroups for closed Lie 
subalgebras and the integrability of Lie algebra homomorphisms for simply connected 
groups. Later \textsc{van Est} and \textsc{Korthagen} 
studied the integrability problem for Banach--Lie algebras 
and found the first example 
of a non-integrable Banach--Lie algebra \cite{EK64}. 
Based on Kuiper's Theorem 
that the unitary group of an infinite-dimensional Hilbert space is contractible 
\cite{Ku65}, simpler examples were constructed 
later by \textsc{Douady} and \textsc{Lazard} \cite{DL66}. 
Chapters 2 and 3 in \textsc{Bourbaki}'s ``Lie Groups and Lie Algebras'' 
contain in particular the basic local theory of Banach--Lie groups and Lie algebras 
and also some global aspects \cite{Bou89}. 
The material in \textsc{K.H. Hofmann}'s Tulane Lecture Notes 
\cite{Hof68}, approaching the subject from a topological group perspective, 
was an important source for people working on Banach--Lie theory 
(see also \cite{Hof72, Hof75}); most of it was published in 
\cite{HoM98}. 

In the early 1970s, \textsc{de la Harpe} extended \textsc{\'E.~Cartan}'s classification of 
Riemannian symmetric spaces to Hilbert manifolds associated to a certain class of 
Hilbert--Lie algebras, called $L^*$-algebras, 
and studied different classes of operator groups related to Schatten ideals. 
Another context, where a structure-theoretic approach leads quite far is the
theory of bounded symmetric domains in Banach spaces and the related theory 
of (normed) symmetric spaces, developed by \textsc{Kaup} and \textsc{Upmeier} 
(cf.\ \cite{Ka81, Ka83, Up85}). 
For a more general approach to Banach symmetric spaces in the sense of \textsc{Loos}
\cite{Lo69}, extending the class of all finite-dimensional symmetric spaces, not 
only Riemannian ones,  we refer to \cite{Ne02c} and \cite{Kl11,Kl12} (cf.\ also 
\cite{La99} for the corresponding 
basic differential geometry). In the context of symplectic geometry, resp., 
Hamiltonian flows, Banach manifolds were introduced by \textsc{Marsden} 
\cite{Mar67}, and 
\textsc{Weinstein} proved a Darboux Theorem in this context \cite{Wei69}.
\textsc{Schmid}'s monograph \cite{Sch87} is a nice introduction to 
infinite-dimensional Hamiltonian systems. 
For more recent results on Banach--K\"ahler manifolds and their connections 
to representation theory, we refer to 
\cite{Ne04b, Bel06} and for Banach--Poisson manifolds to the work of 
\textsc{Ratiu}, \textsc{Odzijewicz} and \textsc{Beltita} 
(\cite{RO03, RO04}, \cite{BR05, BR07}). 
Some ideas on Poisson manifolds beyond the Banach context can be found in 
\cite{NTS15}.

Although \textsc{Birkhoff} was already aware of the fact that his theory 
covered groups of operators on Banach spaces, but not groups of diffeomorphisms, 
it took 30 years until infinite-dimensional Lie groups modeled on 
(complete) locally convex spaces occurred for the first time,  
in the context of Lie group structures  
on the group $\Diff(M)$ of diffeomorphisms 
of a compact manifold~$M$, in the work of \textsc{Leslie} [Les67] and 
\textsc{Omori} \cite{Omo70}. 
This theory was developed further by \textsc{Omori} 
in the context of strong ILB--Lie groups \cite{Omo74}. 
A large part of \cite{Omo74} 
is devoted to the construction of a strong ILB--Lie group 
structure on various types of groups of diffeomorphisms. In the 1980s, this theory 
was refined substantially by imposing and proving additional regularity conditions 
for such groups \cite{OMY82, OMY83a, KYM85}. 

A different type of Lie group was studied by \textsc{Pisanelli} in 
\cite{Pis76, Pis77, Pis79}, 
namely the group $\Gh_n(\C)$ of germs of biholomorphic maps 
of $\C^n$ fixing~$0$. This group carries the structure of a Silva--Lie 
group and is one of the first non-Fr\'echet--Lie groups studied 
systematically in a Lie theoretic context (see also \cite{Da11}). 
In \cite{BCR81}, \textsc{Boseck}, \textsc{Czichowski} and \textsc{Rudolph} approached infinite-dimensional 
Lie groups from a topological group perspective. They use the same concept of an 
infinite-dimensional manifold as we do here, but a stronger Lie group concept. This 
leads to a natural setting for mapping groups on non-compact manifolds 
modeled on spaces of rapidly decreasing functions 
(see also \cite{Nik15, Wa12}). 

In his lecture notes \cite{Mil84}, \textsc{Milnor} 
outlined a general theory of Lie groups modeled on
sequentially complete locally convex spaces, which already contained important cornerstones. 
This paper and the earlier preprint \cite{Mil82} 
had a strong influence on the development of infinite-dimensional Lie theory. 
Both contain precise formulations of several problems, 
some of which have been solved in the meantime and some of which are 
still open, as we shall see in more detail below (see \cite{Gl06b} 
for a survey on some of these problems). 

In the middle of the 1980s, groups of smooth maps, and in particular 
groups of smooth loops became popular because of their intimate connection 
with Kac--Moody theory, topology and string theory (cf.\ \cite{PS86}, 
\cite{Mick87, Mick89}). 
The interest in direct limits of finite-dimensional Lie groups grew in 
the 1990s (cf.\ \cite{NRW91, NRW93, NRW94, NRW01}). They show 
up naturally in the structure and representation theory of Lie algebras 
(cf.\ \cite{Ne98, Ne01b}, \cite{DiPe99}, \cite{NRW01}, 
\cite{NS01}, \cite{Wol05}). The Lie theory of 
these groups was put into its definitive form in 
\cite{Gl05c}. 

Our discussion of smoothness of maps with values in diffeomorphism groups 
of (possibly infinite-dimensional) manifolds 
is inspired by the diffeological approach. We shall see in particular 
in Appendix~\ref{app:nonlie} that, to some extent, one can use differential methods to deal with 
groups with no Lie group structure, such as groups of diffeomorphisms 
of non-compact manifolds or groups of linear automorphisms of 
locally convex spaces, and that this provides a natural framework 
for a Lie theory of smooth actions on manifolds and smooth linear 
representations. \\

In this monograph, we present our 
personal view of the current state of several aspects of the Lie theory of locally 
convex Lie groups. 
We had to make choices, and as a result, 
we could not take up many interesting 
directions such as the modern theory of 
symmetries of differential equations as exposed 
in \textsc{Olver}'s beautiful book \cite{Olv93} and  
the fine structure and the geometry of specific groups of 
diffeomorphisms, such as the group $\Diff(M,\omega)$ of symplectomorphisms 
of a symplectic manifold $(M,\omega)$ 
(\cite{Ban97, MDS98, Pol01} are textbooks on this topic). We do not go into 
(unitary) representation theory (cf.\ \cite{AHM93, Is96, DP03, Pic00a, 
Pic00b, Ki05, Ne14, Ne17, NSZ17})  
and connections to physics, which are nicely described in 
 surveys of \textsc{Goldin} 
\cite{Go04} and \textsc{Schmid} \cite{Sch04}. 
Other topics are only mentioned very briefly, such as 
the ILB and ILH-theory of Lie groups of diffeomorphisms 
which plays an important role in geometric analysis 
(cf.\ \cite{AK98, EMi99}) 
and the group of invertible Fourier integral operators of order zero, 
whose Lie group property 
was the main goal of the series of papers 
by \textsc{Omori} and coauhtors (cf.\ \cite{OMY80, MOK85}). 
An alternative approach to these groups is described in 
\cite{ARS85,ARS86a, ARS86b}. 
Very interesting results concerning diffeomorphism 
groups and Fourier integral operators on non-compact manifolds 
(with bounded geometry) have been obtained by \textsc{Eichhorn} and \textsc{Schmid} 
\cite{ES96, ES01}). We also do not go into Riemannian metrics 
on Fr\'echet homogeneous spaces and shape theory (\cite{MM13, BBM14, BBM16}) 
or the rather recent topic of coarse geometry (\cite{Ro13, Ro18}, \cite{ADM20}). 
Last, but not least, we mention the nice results on the existence of slices in the infinite-dimensional context by 
\textsc{Diez} and \textsc{Rudolph} \cite{DR19}.

\nin {\bf Notation.} Throughout we shall use the notation $\K$ for the fields 
$\R$ and $\C$ of real and complex numbers, 
 and all vector spaces are real or complex. 
For two topological vector spaces $V,W$, we write 
${\cal L}(V,W)$ for the space of continuous linear operators 
$V \to W$ and put ${\cal L}(V) := {\cal L}(V,V)$.

\part{Analysis on locally convex spaces}

\chapter{Differential calculus in locally convex spaces}\label{chapcalcul}
\chaptermark{Differential calculus\hfill
\copyright{} H. Gl\"{o}ckner and K.-H. Neeb}
In this chapter,
we explain our setting of
differential calculus in locally convex topological
vector spaces.
The approach is a very natural one: Partial derivatives,
as familiar from finite-dimensional calculus, are replaced with directional
derivatives, which are required to be continuous both in the point and the direction.
We shall not restrict attention
to locally convex spaces satisfying
certain completeness conditions,
but work with
arbitrary
locally convex spaces.
Moreover,
we shall not only consider
mappings on open subsets
of locally convex spaces,
but develop differential calculus
for mappings on certain
``locally convex'' subsets
with dense interior.
This enables us to
speak about smooth maps on compact intervals,
squares and simplices,
which will play a role
in our
development of 
infinite-dimensional Lie theory.

As the chapter progresses,
we shall encounter analogs of many classical results of differential calculus
in the infinite-dimensional setting, like the\linebreak
Fundamental Theorem of Calculus,
the Chain Rule, and Taylor's Theorem.
As a prerequiste for many important constructions in infinite-dimensional Lie theory
(like the construction of a Lie group structure on the group
$C([0,1],G)$ of continuous paths in an infinite-dimensional Lie group~$G$),
we also discuss the continuity and differentiability
properties of typical mappings between
spaces of continuous functions (and spaces of $C^k$-functions).

The chapter compiles core results of differential calculus which are a prerequisite for the subsequent
introduction to infinite-dimensional Lie groups.
Sections~\ref{seccuint}--\ref{secTay} contain essential foundations,
and may well be read (or taught) in one piece.
The results in Section~\ref{secCspaces} are also essential for Lie theory,
but shall be used only at a later stage; they might be looked up when they are needed.
More specific results of differential calculus,
which are only important for the discussion of certain examples
or Lie groups with particular properties (like analytic Lie groups),
have been relegated to Chapter~\ref{chapfurcalc}.\\[2.2mm]
The letter $\K$ always stands
for~$\R$ or~$\C$.
All vector spaces will be $\K$-vector
spaces and all linear maps will
be $\K$-linear, unless the contrary
is stated.\\[2.2mm]
{\bf Prerequisites for Chapter~\ref{chapcalcul}.}
The reader should be familiar with
basic facts
concerning
locally convex topological vector spaces,
as compiled (with proofs) in Appendix~\ref{chaplcx}.
Also some basic facts from point set topology will be used,
like the Wallace Lemma~\ref{Wallla}.
Section~\ref{secCspaces} presumes basic facts on the compact-open topology
(as recalled in Appendix~\ref{appcotop}).
\section{Curves and integrals in locally convex spaces}\label{seccuint}
Before we turn to differentiable
mappings between subsets of locally convex spaces,
let us discuss the simpler special case
of curves (which also is an important
preliminary for the general case).
\begin{defn}\label{defnC1curve}
Let $E$ be a locally convex space.
A continuous mapping\linebreak
$\gamma\colon  I\to E$
on an
interval $I\sub \R$ is called a \emph{$C^0$-curve}. \index{curve!$C^0$}
If $I$ is non-degenerate,\footnote{That is, $I$ has more than
one element. Henceforth, we shall always assume this.}
then a $C^0$-curve $\gamma\colon I\to E$
is called a
\emph{$C^1$-curve}   \index{curve!$C^1$}
if the limit
\begin{equation}\label{dfdercu}
\gamma'(t)\,:=\,\lim_{s\to 0}\;
{\textstyle \frac{1}{s}(\gamma(t+s)-\gamma(t))}
\end{equation}
exists for all $t\in I$, and $\gamma'\colon
I\to E$, $t\mto \gamma'(t)$
is continuous.\footnote{Thus, we require
that $\lim_{n\to\infty}
\frac{1}{s_n}(\gamma(t+s_n)-\gamma(t))$
exists, for each sequence $s_n\in \R \setminus\{0\}$
such that $t+s_n\in I$
for each $n$, and $\lim_{n\to\infty}s_n=0$.}
We also write $\gamma^{(1)}:=\gamma'$.
Let $\gamma^{(0)}:=\gamma$. Recursively,
given $k\in \N$, we call $\gamma$ a \emph{$C^k$-curve} 
if \index{curve!$C^k$}
$\gamma$ is a $C^{k-1}$-curve and $\gamma^{(k-1)}$ is a $C^1$-curve;
we define $\gamma^{(k)}:=(\gamma^{(k-1)})'$.
If $\gamma$ is $C^k$ for each $k\in \N$, then
the curve $\gamma$ is called $C^\infty$ or
\emph{smooth}. \index{curve!smooth}
Occasionally, we also
write $\frac{d^k\gamma}{dt^k}(t):=\gamma^{(k)}(t)$ or
$\frac{d^k}{dt^k}\big|_{t=t_0}\gamma(t):=\gamma^{(k)}(t_0)$.
\end{defn}
\begin{rem}
It is clear (by a trivial induction)
that $\gamma$ is a $C^k$-curve if and only
if $\gamma$ is a $C^1$-curve
and $\gamma'$ is a $C^{k-1}$-curve.
Furthermore, $\gamma^{(k)}=(\gamma')^{(k-1)}$.
\end{rem}
\begin{rem}\label{remderthcts}
For a map $\gamma\colon I\to E$,
the existence of the derivative (\ref{dfdercu})
at~$t$ implies continuity of $\gamma$ at $t$,
since
\[
\gamma(t+s)-\gamma(t)=s\,\frac{\gamma(t+s)-\gamma(t)}{s}\to 0\cdot\gamma'(t)=0\;\;
\mbox{as $\,s\to 0$.}
\]
\end{rem}
Many results concerning curves in
locally convex spaces
can be reduced to well-known facts
concerning real-valued functions,
by applying continuous linear functionals.
The following simple lemma facilitates this
procedure.
\begin{lem}\label{babychain}
If $E$ and $F$ are locally convex spaces,
$\lambda\colon E\to F$ a continuous linear
map and $\gamma\colon  I\to E$ a $C^1$-curve,
then also $\lambda\circ \gamma$ is a $C^1$-curve,
and
\begin{equation}\label{eqbabych}
(\lambda\circ \gamma)'\; =\; \lambda\circ \gamma'\,.
\end{equation}
\end{lem}
\begin{prf}
For $t\in I$ and $0\not= s\in\R$
such that $t+s\in I$,
using the linearity and continuity of~$\lambda$
we find that
\begin{eqnarray*}
\lim_{s\to 0}\; {\ts \frac{1}{s}(\lambda(\gamma(t+s))-\lambda(\gamma(t)))}
& = & 
\lim_{s\to 0}\; {\ts \lambda\! \left(
\frac{1}{s}(\gamma(t+s)-\gamma(t))\right)}\\
& = &
\lambda\! \left(
\lim_{s\to 0}\;
{\ts \frac{1}{s}(\gamma(t+s)-\gamma(t))}\right)
=\lambda(\gamma'(t))\, .
\end{eqnarray*}
Since $\lambda\circ \gamma'$ is continuous,
$\lambda\circ\gamma$ is~$C^1$,
with the asserted derivative.
\end{prf}
We want to establish (one half of) the Fundamental Theorem of Calculus:
\begin{prop}[Fundamental Theorem of Calculus: First Part]\label{fundamental}
Let $\gamma\colon I\to E$ be a $C^1$-curve
in a locally convex space~$E$,
and $a,b\in I$. Then
\[
\gamma(b)-\gamma(a)\;=\; \int_a^b \gamma'(t)\, d t\,.
\]
\end{prop}
The following terminology and notation is used here:
\begin{defn}\label{defnweakint}
Let $\gamma\colon I\to E$ be a $C^0$-curve in a locally convex $\K$-vector
space~$E$,
and $a,b\in I$. If there exists an element
$z\in E$ such that
\[
\lambda(z)\;=\; \int_a^b \lambda(\gamma(t))\, dt\quad \mbox{for
all $\lambda\in E'$,}
\]
then $z\in E$ is called the \emph{weak integral of
$\gamma$} from $a$ to $b$, and denoted \index{weak integral} 
\[
\int_a^b\gamma(t)\, dt\; :=\; z\,.
\]
\end{defn}
\begin{rem}
Note that the element $z\in E$ in Definition~\ref{defnweakint}
is uniquely determined if it exists, as the dual space $E'$
of all continuous linear functionals
separates points on~$E$ by the Hahn--Banach Theorem
(Theorem~\ref{dualsep}).
\end{rem}
\noindent{\em Proof of  Proposition}~\ref{fundamental}.
Let $\lambda\in E'$.
By Lemma~\ref{babychain}, $\lambda\circ \gamma\colon I\to \K$
is a $C^1$-curve, and $(\lambda\circ \gamma)'=\lambda\circ
\gamma'$.
The standard Fundamental Theorem
of Calculus yields
\[
\lambda(\gamma(b)-\gamma(a))=
\lambda(\gamma(b))-\lambda(\gamma(a))=\int_a^b(\lambda\circ \gamma)'(t)\, dt
=\int_a^b \lambda(\gamma'(t))\, dt\,.
\]
Hence $z=\gamma(b)-\gamma(a)$
satisfies the defining property of the
weak integral
$\int_a^b\gamma'(t)\, dt$, as required.\qed
\begin{rem}\label{remlk}\label{remkonst}
We make some simple observations concerning the
existence of weak integrals. Further basic facts are compiled in Exercises~\ref{excbasicwint}
and~\ref{excprodpartint}(d).\medskip

\noindent(a)
If the weak integrals of $\gamma\colon I\to E$
and $\eta\colon I\to E$ from $a$ to~$b$ exist,
and $s \in \K$, then also the weak integral of
$\gamma+s\eta$ exists, and
\begin{equation}\label{lincob}
\int_a^b(\gamma(t)+s\eta(t))\, dt=\int_a^b\gamma(t)\, dt
+s\int_a^b\eta(t)\, dt\,.
\end{equation}
In fact, it is readily verified that
the right hand side of (\ref{lincob})
satisfies the defining property
of the weak integral $\int_a^b(\gamma(t)+s\eta(t))\, dt$.\medskip

\noindent(b)
If $\gamma\colon I\to E$ is constant, say $\gamma(t)=K$,
then $\int_a^b\gamma(t)\, dt$ exists and is equal to
$K\cdot (b-a)$. In fact, it is easily verified
that $K\cdot (b-a)$ satisfies the defining property
of $\int_a^b\gamma(t)\, dt$.\medskip

\noindent(c)
If $\gamma\colon I\to E$ is a $C^0$-curve in a complex locally
convex space~$E$ and $z:=\int_a^b\gamma(t)\, dt$
exists in $E$ considered as a real
vector space, then~$z$ also is the weak integral
of~$\gamma$ from $a$ to~$b$ in $E$,
considered as a complex vector space.
To verify this, let
$\lambda\colon E\to \C$ be a continuous complex linear functional.
Then
\[
\lambda(x)=\Repart(\lambda(x))+i\,\Impart(\lambda(x))=
\Repart(\lambda(x))-i\,\Repart(\lambda(ix))=\rho(x)-i\sigma(x)\,,
\]
where $\rho\colon E\to \R$,
$\rho(x):=\Repart(\lambda(x))$
and $\sigma\colon E\to \R$, $\sigma(x):=\Repart(\lambda(ix))$
are continuous real linear functionals.
Hence (as required)
\begin{eqnarray*}
\lambda(z) & = & \rho(z)-i\sigma(z)=\int_a^b\rho(\gamma(t))\, dt
-i\int_a^b \sigma(\gamma(t))\, dt\\
& = & 
\int_a^b(\rho(\gamma(t))
-i\sigma(\gamma(t)))\, dt
=
\int_a^b\lambda(\gamma(t))\, dt\,.
\end{eqnarray*}
\end{rem}
We now generalize the
estimate $\bigl|\int_a^bf(t)\, dt\bigr|
\leq |b-a|\cdot \sup\{|f(t)|\colon t\in [a,b]\}$
for Riemann integrals of continuous real-valued functions.
\begin{lem}\label{weakinballb}
Let $a\leq b$ be real numbers and
$\gamma\colon [a,b]\to E$ be a continuous curve
in a locally convex space~$E$. If the weak integral
$\int_a^b\gamma(t)\, dt$
exists,
then
\[
\left\| \int_a^b\gamma(t)\, dt\right\|_p
\,\leq\,   \int_a^b\|\gamma(t)\|_p\,dt
\, \leq\,  (b-a)\cdot
\max\{\|\gamma(t)\|_p\colon t\in [a,b]\}
\]
for each continuous seminorm $\|\cdot\|_p$
on~$E$.
\end{lem}
\begin{prf}
We may assume that $\K=\R$.
Let $z:=\int_a^b\gamma(t)\,dt$.
By the Hahn--Banach Extension Theorem (Corollary~\ref{corHBext}),
there is $\lambda\in E'$ such that
\[
\mbox{$|\lambda(x)|\leq \|x\|_p$ for all $x\in E$ and $\lambda(z)=\|z\|_p$.}
\]
Hence
\[
\left\| \int_a^b\gamma(t)\, dt\right\|_p=\lambda(z)=\int_a^b\lambda(\gamma(t))\,dt
\leq \int_a^b|\lambda(\gamma(t))|\,dt
\leq\int_a^b \|\gamma(t)\|_p\,dt,
\]
from which the assertions follow.
\end{prf}
The following lemma
makes it easy to see that
weak integrals depend continuously on parameters.
\begin{lem}\label{comparg1}
Let $P$ and $X$ be topological spaces,
$K\sub X$ be compact
and $F$ be a topological vector space.
Let $f\colon P\times X \to F$
be a continuous map,
$U\sub F$ be a
$0$-neighborhood,
and $p\in P$.
Then there exists an open neighborhood $Q\sub P$
of~$p$ and an open subset $Y\sub X$
containing~$K$,
such that
\[
f(q,y)-f(p,y)\;\in U\quad
\mbox{for all $\, q\in Q$ and $y\in Y$.}
\]
\end{lem}
\begin{prf}
Let $p\in P$.
The map $g\colon P\times X\to F$, $g(q,y):=f(q,y)-f(p,y)$
is continuous and $g(p,y)=0$ for each $y\in X$.
Hence
$g^{-1}(U)$
is an open subset of $P\times X$ which contains $\{p\}\times K$.
Since $\{p\}$ and $K$ are compact, the Wallace Lemma~\ref{Wallla}
provides open subsets $Q\sub P$ and $Y\sub X$ such that
\[
\{p\}\times K\, \sub \, Q\times Y\, \sub \, g^{-1}(U).
\]
Then $f(q,y)-f(p,y)=g(q,y)\in U$ for all $q\in Q$ and $y\in Y$.
\end{prf}
\begin{lem}[Continuity of parameter-dependent integrals]\label{intpar}
Let $P$ be a topological space
and $a< b$ be real numbers. Let $f\colon P\times [a,b]\to E$
be a continuous map to a locally convex
space~$E$.
If the weak integral
\[
g(p)\; :=\; \int_a^b f(p,t)\; dt
\]
exists in~$E$ for each $p\in P$,
then $g\colon  P\to E$ is continuous.
\end{lem}
\begin{prf}
To prove continuity at $p\in P$, let $\|\cdot\|_s$ be a continuous seminorm
on~$E$ and $\ve>0$.
By Lemma~\ref{comparg1},
there is a neighborhood $Q\sub P$ of~$p$ such that
\[
\|f(q,t)-f(p,t)\|_s\leq \frac{\ve}{b-a}
\quad \mbox{for all $q\in Q$ and $t\in[a,b]$.}
\]
Using Lemma~\ref{weakinballb}, we deduce that
\begin{eqnarray*}
\|g(q)-g(p)\|_s&=& \left\|\int_a^b (f(q,t)-f(p,t))\, dt\right\|_s \leq \int_a^b\|f(q,t)-f(p,t)\|_s\, dt\\
&\leq & (b-a)\frac{\ve}{b-a}\leq\ve.
\end{eqnarray*}
Thus $g$ is continuous at $p$.
\end{prf}
We now formulate the second half of the Fundamental
Theorem of Calculus.
See Proposition~\ref{rieman}
for a criterion ensuring the
existence of
weak integrals.
\begin{prop}[Second Part of the Fundamental Theorem]\label{secpart}
Let $\gamma\colon I\to E$ be a $C^0$-curve, $a\in I$, and
assume that the weak integral
\[
\eta(t):=\int_a^t \gamma(s)\, ds
\]
exists for all $t\in I$. Then $\eta\colon I\to E$
is a $C^1$-curve in~$E$, and $\eta'=\gamma$.
\end{prop}
\begin{prf}
Let $t_0\in I$ and $t\in I-t_0$ such that $t\not=0$.
Then
\begin{eqnarray}
\frac{\eta(t_0+t)-\eta(t_0)}{t} &=& \frac{1}{t}\left(
\int_a^{t_0+t}\gamma(s)\, ds-\int_a^{t_0}\gamma(s)\, ds\right) \notag \\
&=& \frac{1}{t}\int_{t_0}^{t_0+t}\gamma(s)\,ds=\int_0^1 \gamma(t_0+rt)\, dr, \label{secpafd}
\end{eqnarray}
using Exercise~\ref{excbasicwint}(e) to combine the integrals
and then the Substitution Rule
with $s=t_0+rt$, $ds=t\, dr$ (see Exercise~\ref{excbasicwint}(d)).
Note that the final integrand also makes sense for $t=0$, and defines
a continuous map
\[
f \colon (I-t_0)\times [0,1]\to E,\quad f(t,r) := \gamma(t_0+rt).
\]
Since $f(0,r)=\gamma(t_0)$ is independent of~$r$,
the weak integral $\int_0^1 f(t,r)\,dr$ exists also for $t=0$
and is given by $\gamma(t_0)$ then.
Hence
\[
\frac{\eta(t_0+t)-\eta(t_0)}{t}  =\int_0^1 f(t,r)\,dr\to\int_0^1 f(0,r)\, dr=\gamma(t_0)\;\,\mbox{as $t\to 0$,}
\]
by Proposition~\ref{intpar}.
Thus $\eta'=\gamma$, which is a continuous map.
Since $\eta$ is continuous (see Remark~\ref{remderthcts}),
we find that $\eta$ is a $C^1$-curve.
\end{prf}
We can also interchange weak integrals and
uniform limits.
\begin{defn}\label{defnunifcon}
Let $X$ be a set and $E$ be a locally convex space.
We say that a sequence $(\gamma_n)_{n\in \N}$
of maps $\gamma_n\colon X \to E$
\emph{converges
uniformly} to a map \index{uniform convergence} \index{convergence!uniform} 
$\gamma\colon X\to E$ if, for each continuous seminorm $\|\cdot\|_p$ on~$E$,
\[
\lim_{n\to\infty}\, \sup \big\{\|\gamma_n(x)-\gamma(x)\|_p\colon
x\in X \big\}\, =\, 0\,.
\]
\end{defn}
\begin{lem}\label{lemcurvesconv}
Let $E$ be a locally convex space,
$a\leq b$ be real numbers
and
$(\gamma_n)_{n\in \N}$
be a sequence of $C^0$-curves
$\gamma_n\colon  [a,b] \to E$
that converges uniformly
to a $C^0$-curve $\gamma\colon [a,b] \to E$.
If the weak integrals
$\int_a^b \gamma_n(t)\, dt$
and
$\int_a^b \gamma(t)\, dt$ exist in~$E$, then
\begin{equation}\label{versionpar}
\lim_{n\to\infty}
\int_a^b \gamma_n(t)\, dt\;=\; \int_a^b \gamma(t)\, dt\,.
\end{equation}
\end{lem}
\begin{prf}
Let $\|\cdot\|_p$ be a continuous seminorm on~$E$.
Then
\begin{eqnarray*}
\Big\|\int_a^b\gamma_n(t)\, dt-\int_a^b\gamma(t)\,dt\Big\|_p
&=&
\Big\|\int_a^b(\gamma_n(t)-\gamma(t))\,dt\Big\|_p\\
&\leq & (b-a)\,\sup\{\|\gamma_n(t)-\gamma(t)\|_p\colon t\in [a,b]\}
\end{eqnarray*}
tends to~$0$
as $n\to\infty$,
using Lemma~\ref{weakinballb}
and the uniform convergence.
\end{prf}
It is natural to ask
for conditions ensuring that
weak integrals exist.
It turns out that
weak integrals in
\emph{sequentially complete} \index{sequentially complete space} 
locally convex spaces (as in Definition~\ref{defnCauchy}(b))
always exist.
Notably, weak integrals
in complete locally convex spaces
(see Definition~\ref{defnCauchy}(d))
always exist.
\begin{prop}\label{rieman}
Let $E$ be a sequentially complete locally convex
space and $\gamma\colon [a,b] \to E$
be a $C^0$-curve.
Then the weak integral
$\int_a^b \gamma(t)\, dt$
exists~in~$E$. 
\end{prop}
\begin{prf}
Using that~$E$ is sequentially
complete, the weak integral can be constructed
as a limit of a sequence of Riemann sums
(as in the case of real-valued functions).
Details can be looked
up in
Section~\ref{secapp-ch2},
the appendix to Chapter~\ref{chapcalcul}
(after Lemma~\ref{unif}).
\end{prf}
\begin{rem}\label{remMackeyc}
Recall from
Proposition~\ref{cmplexists}
that every locally convex space~$E$
can be completed.
This is very useful, because
the weak integral $z:=\int_a^b \gamma(t)\, dt$
of a $C^0$-curve
$\gamma\colon [a,b]\to E$
always exists in the completion $\wt{E}$,
by Proposition~\ref{rieman}.
Consider $E$ as a vector subspace of~$\wt{E}$.
As we shall see in Exercise~\ref{excmorewint}(b),
the weak integral $w:=\int_a^b\gamma(t)\,dt$ in~$E$
exists if and only if $z\in E$,
in which case $z=w$.
Frequently, this makes
it possible to discuss a problem
concerning an $E$-valued mapping~$f$
in two steps: First, we consider
$f$ as a map to~$\wt{E}$,
where weak integrals are available, and solve
the problem there.
In the next step, we verify that
all weak integrals of interest
lie in~$E$, and verify that the
original problem
is also solved, for~$f$ as a map to~$E$.
\end{rem}
For many purposes,
sequential completeness can be replaced by
``Mackey completeness,''
a weaker completeness property
which we recall~now.
\begin{defn}\label{defncomptness}
Let $E$ be a locally convex space.
\begin{description}[(D)]
\item[(a)]
Let $a<b$ be real numbers.
A map $\gamma\colon [a,b]\to E$
is called
a \emph{Lipschitz curve} \index{curve!Lipschitz} \index{Lipschitz curve} 
if the set
\[
\Big\{\frac{\gamma(s)-\gamma(t)}{s-t}\colon
\mbox{$s,t\in [a,b]$, $s\not=t$}\Big\}
\]
is bounded in~$E$.
Thus, $\gamma$ is a Lipschitz curve if and only if,
for each continuous seminorm~$\|\cdot\|_p$ on~$E$,
there exists $L \in [0,\infty[$ such that
\[
\|\gamma(s)-\gamma(t)\|_p\;\leq\; L \, |s-t|\quad
\mbox{for all $s,t\in [a,b]$.}
\]
In particular, every Lipschitz curve is continuous.
\item[(b)]
The space~$E$
is called
\emph{Mackey complete} \index{Mackey complete l.c. space}
\footnote{Mackey complete locally convex spaces are also called
``convenient,'' ``$c^\infty$-complete,''
or ``locally complete'' in the literature.}
if the weak integral $\int_a^b \gamma(t)\, dt$
exists in~$E$
for all real numbers $a<b$
and Lipschitz curves $\gamma\colon [a,b]\to E$.
\end{description}
\end{defn}
\begin{rem}
(a)
Every sequentially complete locally convex
space is Mackey complete.
In fact, if
$E$ is sequentially
complete,
then the weak integral $\int_a^b\gamma(t)\, dt$
exists in~$E$ for each $C^0$-curve
$\gamma\colon [a,b]\to E$,
by Proposition~\ref{rieman}.\vspace{-.4mm}
Hence $\int_a^b\gamma(t)\, dt$
exists a fortiori
for each Lipschitz curve
$\gamma\colon [a,b]\to E$
and thus~$E$ is Mackey complete.\medskip

\noindent(b)
Using Proposition~\ref{fundamental}
and Lemma~\ref{weakinballb},
it is easy to see
that every $C^1$-curve on a compact interval
is Lipschitz (Exercise~\ref{excC1lip}).\medskip

\noindent(c)
Replacing
Lipschitz curves by $C^k$-curves
for a fixed
$k\in \N\cup \{\infty\}$
in the definition of Mackey
completeness, one obtains an equivalent
property (cf.\ \cite[Thm.\,2.14]{KM97}).
Lipschitz curves have the advantage
that a curve $\gamma\colon [a,b]\to E$
in a locally convex space~$E$
whose image is contained in
a vector subspace $E_0\sub E$
is a Lipschitz curve in~$E$ if and only
if it is a Lipschitz curve in~$E_0$
(cf.\ Lemma~\ref{lembdsub}).
While all difference quotients automatically
lie in~$E_0$ if so does the image of~$\gamma$,
the corresponding differential quotients
involve limits which may exist
in~$E$ but not in~$E_0$.
Examples show that
a $C^1$-curve~$\gamma$
in~$E$ with image in a vector subspace
$E_0\sub E$ need not be a $C^1$-curve
in~$E_0$ (see Exercise~\ref{excpathcurvsub}).
We shall use Lipschitz curves
only once, to prove the implication
``(e)$\impl$(a)'' in Proposition~\ref{analytsingle}.\medskip

\noindent
(d)
Mackey completeness
ensures the existence
of weak integrals
of particular importance.
Notably, as we shall
see in Section~\ref{seccxan},
the existence of weak integrals
substantially simplifies the theory of complex
analytic functions to Mackey complete locally
convex spaces.\medskip

\noindent
(e)
It can be shown that a locally convex space~$E$
is Mackey complete if and only if every
\emph{Mackey--Cauchy sequence} 
\index{Mackey--Cauchy sequence} 
$(x_n)_{n\in \N}$
in~$E$ is convergent,
i.e., every sequence
for which there exists a bounded subset $B\sub E$
and a family $(r_{k,\ell})_{k,\ell\in \N}$
of real numbers $r_{k,\ell} \geq 0$
such that $x_k-x_\ell \in r_{k,\ell}B$
for all $k,\ell \in \N$,
and such that
$r_{k,\ell} \to 0$ as both $k,\ell \to \infty$,
in the sense that for each $\ve >0$,
we find $N\in \N$ such that $r_{k,\ell}\leq \ve$
for all $k,\ell\geq N$
(cf.\ \cite[Lem.\,2.2 and Thm.\,2.14]{KM97},
where a long list of alternative characterizations
of Mackey completeness can be found. For further
characterizations, see~\cite[Thm.\,A.2]{GE92}).\medskip

\noindent
(f)
In the appendix to this chapter,
we show that Mackey completeness
also ensures the existence of weak
integrals of $C^1$-maps
on higher-dimensional simplices
(Proposition~\ref{intsimplex}).
The argument should only be read after
Proposition~\ref{diffpar}.
See also \cite{BG14} and \cite{NSZ14} for parameter dependence of weak integrals
over higher-dimensional (or more general) sets.\medskip

\noindent
(g)
Let us say that a locally convex space $E$ is 
\index{integral complete l.c. space}
\emph{integral complete}\vspace*{-.4mm}
if the weak integral $\int_0^1\gamma(t)\,dt$ exists in $E$ for each $C^0$-curve $\gamma\colon [0,1]\to E$.
Every sequentially complete locally convex space is integral complete;
every integral complete locally convex space is Mackey complete.
It is known that a locally convex space $E$ is integral complete if and only if it has the
\index{metric convex compactness property}  
\emph{metric convex compactness property}, 
i.e., the closed convex hull
$\wb{\conv(K)}$
is compact for each metrizable compact subset $K\sub E$
(\cite{Wz12}, cf.\ \cite{Voi92}). See \cite[Ex.~A.3]{GE92}
for a Mackey complete locally convex space which fails to be integral complete.
\end{rem}
\begin{small}
\subsection*{Exercises for Section~\ref{seccuint}}
\begin{exer}\label{excintRn}
Show: The weak integral
of a $C^0$-curve
\[ \gamma=(\gamma_1,\ldots, \gamma_n)\colon [a,b]\to\R^n \] 
is given by
$\int_a^b\gamma(t)\, dt=\bigl(\int_a^b\gamma_1(t)\,dt,\ldots,
\int_a^b\gamma_n(t)\,dt\bigr)$,
where $\int_a^b\gamma_j(t)\, dt\in \R$
is the usual Riemann integral.
\end{exer}

\begin{exer}\label{excbasicwint}
Let $E$ be a locally convex space,
$\gamma,\eta\colon I\to E$ be $C^0$-curves,\break 
{$a,b,c\in I$,} $t_0\in \R$
and $r\in \K$.
Assume that the weak integral
on the right hand side exists
and deduce that the weak
integral on the left hand side
exists and equality holds:
\begin{description}[(D)]
\item[(a)]
$\int_a^b (\gamma(t)+r\eta(t))\, dt =
\int_a^b \gamma(t)\, dt\,+\, r\int_a^b \eta(t)\, dt$;
\item[(b)]
$\int_a^b \gamma(t)\, dt =
\int_{a-t_0}^{b-t_0} \gamma(s+t_0)\, ds$.
\item[(c)] $\int_b^a\gamma(t)\, dt=-\int_a^b\gamma(t)\, dt$;
\item[(d)]
$\int_\alpha^\beta \gamma(\phi(t))\, \phi'(t)\, dt
=\int_a^b \gamma(s)\, ds$,
where
$\phi\colon [\alpha,\beta]\to I$
is a $C^1$-function such that
$\phi(\alpha)=a$ and $\phi(\beta)=b$ (Substitution Rule).
\end{description}
In (b), (c), and (d), also show that existence of the integral on the right hand side implies
existence of the integral on the left.
Finally, show that if two of the following weak integrals exist in~$E$, then also the third,
and equality holds:
\begin{description}[(D)]
\item[(e)] 
$\int_a^c\gamma(t)\, dt =
\int_a^b\gamma(t)\, dt \,+\, \int_b^c\gamma(t)\, dt$.
\end{description}
\end{exer}

\begin{exer}\label{excmorewint}
Let $E$ be a locally convex space,
$\gamma\colon I\to E$ be a $C^0$-curve
and $a,b\in I$.
\begin{description}[(D)]
\item[(a)]
Show that
if the weak integral
$z:=\int_a^b \gamma(t)\, dt$ exists in~$E$
and $\alpha\colon E\to F$
is a continuous linear map to a
locally convex space~$F$,
then $\alpha(z)$ is the weak integral
of $\alpha\circ \gamma$
from~$a$ to~$b$.
\item[(b)]
If $Y\sub E$ is a vector subspace such that $\gamma(I)\sub Y$
and~$y$ an element of~$Y$, show that $y=\int_a^b\gamma|^Y(t)\,dt$ in~$Y$ (endowed
with the topology induced by~$E$) if and only if $y=\int_a^b\gamma(t)\,dt$ in~$E$.\\[.3mm]
[Recall that each continuous linear functional $\lambda\in Y'$
extends to a continuous linear functional on~$E$,
by the Hahn-Banach extension theorem.]
\item[(c)]
Let $\wt{E}$
be a completion of~$E$
such that $E\sub \wt{E}$.
Let
$w:=\int_a^b \gamma(t)\, dt$
be the weak integral of
$\gamma$ in~$\wt{E}$
(i.e.,
the weak integral of
$\kappa \circ \gamma$,
where
$\kappa \colon E\to \wt{E}$
is the inclusion map).
Show that the weak integral
$z:=\int_a^b\gamma(t)\, dt$
exists in~$E$ if and only
if $w\in E$,
in which case
$z=w$.
\item[(d)]
Let $a<b$ now and
assume that
$z:=\int_a^b\gamma(t)\, dt$ exists in~$E$.
Let $A\sub E$ be a closed convex set
such that $\gamma([a,b])\sub A$.
Since $z$ coincides
with the weak integral in~$\wt{E}$,
the proof of Proposition~\ref{rieman}
provides a sequence
$(S_n)_{n\in \N}$
of Riemann sums
converging
to $\int_a^b\gamma(t)\, dt$
in $\wt{E}$ and hence in~$E$.
Show that $S_n\in (b-a)A$
for each $n\in \N$ and
deduce that
$\int_a^b \gamma(t)\, dt\in (b-a)A$.
\end{description}
\end{exer}

\begin{exer}\label{exccurveinsub}
Show that if
$\gamma\colon I\to E$
is a $C^1$-curve in a locally
convex space~$E$
and $E_0\sub E$ is a closed vector subspace
such that
$\gamma(I)\sub E_0$,
then $\gamma$ also is $C^1$
when considered as a curve in~$E_0$.
The same conclusion holds if~$E_0$ is not necessarily closed but
both $\gamma(I)\sub E_0$ and $\gamma'(I)\sub E_0$.
\end{exer}

\begin{exer}\label{excpsepwint}
Let $E$ be a locally convex space, $a<b$ and $\gamma\colon [a,b]\to E$
be a continuous curve such that the weak integral $\int_a^b\gamma(t)\, dt$ exists in~$E$.
Let $\Lambda\sub E'$ be a set of continuous linear functionals which separate
points on~$E$, i.e., for all $x\not=y$ in~$E$ we find $\lambda\in\Lambda$ such that
$\lambda(x)\not=\lambda(y)$. Let $z\in E$ be an element such that
\[
\lambda(z)=\int_a^b\lambda(\gamma(t))\, dt\quad\mbox{for all $\lambda\in\Lambda$.}
\]
Show that $z=\int_a^b\gamma(t)\, dt$.
\end{exer}

\begin{exer}\label{excpathcurvsub}
Let $E:=\R^\N$,
equipped with the product topology and
let $E_0:=\ell^\infty\sub E$
be the space of bounded real-valued
sequences, equipped with the
topology induced by~$E$.
Consider $\gamma\colon \R\to E$,
$\gamma(t) := (\cos(kt))_{k\in \N}$
and $\eta\colon \R\to E$, $\eta(t):=(\sin(kt))_{k\in\N}$.
Since~$E$ is endowed with the product topology,
$\gamma$ and~$\eta$ are continuous.
\begin{description}[(D)]
\item[(a)]
Since~$E$ is complete, the weak integral $\int_0^t\gamma(s)\, ds$ exists in~$E$ for each $t\in\R$.
Using the continuous linear point evaluations $\ve_k\colon E\to \R$, $(x_n)_{n\in\N}\mapsto x_k$,
deduce with Exercise~\ref{excpsepwint} that $\eta(t)=\int_0^t\gamma(s)\, ds$. 
\item[(b)]
Deduce that $\eta$ is a $C^1$-curve in~$E$, with $\eta'=\gamma$.
\item[(c)]
Show that although $\eta(\R)\sub E_0$,
the co-restriction $\eta|^{E_0}\colon \R\to E_0$
is not a $C^1$-curve in~$E_0$.
\end{description}
\end{exer}

\begin{exer}\label{excpathnonlcx}
Given $0<p<1$,
let $L^p[0,1]$
be the topological vector space
of equivalence classes $[\gamma]$
of measurable functions
$\gamma\colon [0,1]\to \R$
such that
\[
\int_0^1 |\gamma(s)|^p\, ds<\infty,
\]
modulo functions vanishing almost
everywhere.
The topology on
$L^p[0,1]$
comes from the metric
\[
d\colon L^p[0,1] \times L^p[0,1]\to [0,\infty[,
\quad d([\gamma],[\eta]):=\int_0^1 |\gamma(s)-\eta(s)|^p\, ds.
\]
Using characteristic functions,
we define a curve
\[
\beta \colon [0,1] \to L^p[0,1]\,,\qquad
\beta (t):=[\one_{[0,t[}]\,.
\]
Show that $\beta$
is injective and differentiable at each point,
with $\beta'(t)=0$ for each
$t\in [0,1]$.
Deduce that $L^p[0,1]$ is not locally
convex.
\end{exer}

\begin{exer}\label{excprecia}
Let~$\cA$ be a locally convex,
unital, associative
topological algebra (see
Definition~\ref{deftopalg})
and
$\gamma\colon I\to \cA$ be a $C^1$-curve.
Assume that $\gamma(I)$ is contained in the set $\cA^\times\sub\cA$
of invertible elements
and that the inversion map $\iota\colon \cA^\times\to\cA$, $x\mto x^{-1}$ is continuous.
Show that
\[
\eta:=\iota\circ \gamma\colon I\to\cA, \quad t\mto \gamma(t)^{-1}
\]
is $C^1$ with $\eta'(t)=-\gamma(t)^{-1}\gamma'(t)\gamma(t)$ for all $t\in I$.
\,[Hint: Write
\[
\frac{\gamma(t+s)^{-1}-\gamma(t)^{-1}}{s}=
-\gamma(t+s)^{-1}\left(\frac{\gamma(t+s)-\gamma(t)}{s}\right)\gamma(t)^{-1}
\]
for $t\in I$ and $0\not=s\in\R$ such that $t+s\in I$].
\end{exer}

\begin{exer}\label{excunifk}
Let $F$ be a locally convex space,
$X$ be a topological space,
$(\gamma_n)_{n\in \N}$
be a sequence of maps
$\gamma_n\colon X\to F$
and $\gamma\colon X\to F$
be a map. Show the following:
\begin{description}[(D)]
\item[(a)]
If each $\gamma_n$ is continuous
and $\gamma_n\to\gamma$ uniformly,
then also $\gamma$ is continuous.
\item[(b)]
If $X$ is metrizable, then a map
$f \colon X\to Y$ to a topological space~$Y$
is continuous if and only if $f|_K$
is continuous for each compact subset $K\sub X$
(i.e., $X$ is a so-called ``$k$-space".)\footnote{See \ref{defnkkR}, also Exercise~\ref{exc-kviamap}.}
Hint: Together with its limit,
the elements of a convergent sequence
form a compact set.
\item[(c)]
If $X$ is metrizable,
each $\gamma_n$ is continuous
and $\gamma_n|_K$ converges uniformly
to $\gamma|_K$ for each compact subset
$K\sub X$,
then $\gamma$ is continuous.
\end{description}
\end{exer}

\begin{exer}\label{excbilincts}
Let $E_1$, $E_2$ and $F$ be locally convex spaces,
$\beta\colon E_1\times E_2\to F$ be a continuous
bilinear map and
$W\sub F$ be a $0$-neighborhood.
\begin{description}[(D)]
\item[(a)]
Observe that there exist
$0$-neighborhoods $U_1\sub E_1$ and $U_2\sub E_2$
such that $\beta(U_1\times U_2)\sub W$.
For each~$r>0$, find a $0$-neighborhood
$V\sub E_1$ such that $\beta(V\times rU_2)\sub W$.
\item[(b)]
Show that, for each bounded subset
$B\sub E_2$, there exists a $0$-neighborhood
$V\sub E_1$ such that $\beta(V\times B)\sub W$.
\end{description}
\end{exer}

\begin{exer}\label{excrefinedwint}
Let $E$ be a locally convex space
and $I=[a,b]$.
\begin{description}[(D)]
\item[(a)]
Verify that the definition
of $\int_a^b \gamma(t)\, dt$
can be extended to the case
that\linebreak
$\gamma\colon I\to E$ is merely
``weakly continuous,''
i.e.
$\lambda\circ \gamma\colon I\to \K$
is continuous for each $\lambda\in E'$.
Also the estimate
\[
\left\|  \int_a^b\gamma(t)\, dt\right\|_p\leq (b-a)\sup\{\|\gamma(t)\|_p\colon t\in [a,b]\}
\]
(cf.\ Lemma~\ref{weakinballb})
remains valid, for each continuous seminorm $\|\cdot\|_p$~on~$E$.
\item[(b)]
Let $P$ be a topological space,
$\alpha\colon P\times I\to \K$ be
a continuous function and\linebreak
$\beta\colon I\to E$
be a weakly continuous map.
Write $E_w$ for~$E$,
equipped with the weak topology
(see Definition~\ref{defweaktop}).
Then $\beta(I)$ is compact in $E_w$,
hence bounded in~$E_w$ and hence bounded
in~$E$, by Mackey's Theorem (Theorem~\ref{Mackey}).
We assume that the weak integral
$g(p) := \int_a^b \alpha(p,t)\, \beta(t)\,dt$
exists in~$E$, for each $p \in P$.
Using Exercise~\ref{excbilincts}(b),
show that $g\colon P\to E$ is continuous.
\end{description}
\end{exer}

\begin{exer}\label{excC1lip}
Show that each $C^1$-curve
$\gamma\colon [a,b]\to E$
on a compact interval is Lipschitz.
\end{exer}

\begin{exer}\label{excunifcnets}
\begin{description}[(D)]
\item[(a)]
Replacing $\N$ with a directed set $(A,\leq)$ in Definition~\ref{defnunifcon},
we obtain a notion of uniform convergence for nets of functions.
Show that the conclusions of Lemma~\ref{lemcurvesconv}
and Exercise~\ref{excunifk} remain valid if sequences are replaced with
nets.
\item[(b)]
Let $(\gamma_\alpha)_{\alpha\in A}$
be a net of $C^1$-curves $\gamma_\alpha\colon [a,b]\to E$ in a locally convex space~$E$
such that $\gamma_\alpha\to\gamma$ and $\gamma_\alpha'\to\eta$ for certain
continuous curves $\gamma,\eta\colon [a,b]\to E$.
Show that $\gamma$ is a $C^1$-curve and $\gamma'=\eta$.\\[.3mm]
[For $x\in[a,b]$, by Proposition~\ref{fundamental} we have
$\gamma_\alpha(x)=\gamma_\alpha(a)+\int_a^x\gamma_\alpha'(t)\,dt$
for all $\alpha\in A$, by Proposition~\ref{fundamental}. Using (a), passing to the limit
we get $\gamma(x)=\gamma(a)+\int_a^x\eta(t)\,dt$. The assertion now follows from Proposition~\ref{secpart}].
\end{description}
\end{exer}
\end{small}
\section[The differential calculus of $C^1$-maps]{Differential
calculus: {\boldmath $C^1$\/}-maps}\label{secC1}
We now discuss $C^1$-maps
on open subsets of locally convex spaces.
\begin{defn}\label{defnC1map}
Let $E$ and $F$ be locally convex $\K$-vector spaces
and\linebreak
$f\colon U\to F$ be a mapping on an open
subset $U\sub E$.
The \emph{derivative of $f$ at $x\in U$
in the direction $y\in E$} 
\index{directional derivative} 
\index{derivative!directional} 
is defined as
\[
df(x,y)\;:=\; (D_y\, f)(x)\;:=\;
\lim_{t\to 0}
\, \frac{1}{t}\, \Big(f(x+ty)-f(x)\Big)\,,
\]
whenever the limit exists.\footnote{Here $0\not= t\in \K$
such that $x+ty\in U$.}
We say that $f$ is
$C^1_\K$ (or simply $C^1$, if $\K$ is understood)
if $f$ is continuous,
the directional derivative
$df(x,y)$
exists for all $(x,y)\in U\times E$, and the map
$df\colon U\times E\to F$ so obtained is continuous.
In later chapters,
we shall also write $\dd f(x)(y):=df(x,y)$.
\end{defn}
\begin{ex}\label{exlinC1}
If $\lambda\colon E\to F$ is a continuous linear map
between
locally convex spaces~$E$ and~$F$,
then~$\lambda$ is $C^1$.
In fact, for all $x,y\in E$
and $t\in \K^\times$,
using the linearity of~$\lambda$, we obtain
\[
\frac{1}{t}\big( \lambda(x+ty)-\lambda(x)\big)\;=\;
\lambda(y)\,,
\]
which converges to $\lambda(y)$ as $t\to 0$.
Hence the directional derivative
$d\lambda(x,y)$ exists and is given
by
\begin{equation}\label{difflin}
d\lambda(x,y)\; =\; \lambda(y)\qquad
\mbox{for all $\, x,y\in E$.}
\end{equation}
The map $d\lambda$ being continuous,
$\lambda$ is~$C^1$.
For later use,
note that $d\lambda\colon E\times E\to F$ is also
a continuous linear map.
\end{ex}
\begin{ex}\label{exmultilinC1}
Let $E_1,E_2,F$
be locally convex spaces and $\beta\colon E_1\times
E_2\to F$ be a continuous bilinear map.
Then $\beta$ is $C^1$, and
\begin{equation}\label{diffbilin}
d\beta((x_1,x_2), (y_1,y_2))\;=\;
\beta(x_1,y_2)+\beta(y_1,x_2)
\end{equation}
for all $x_1,y_1\in E_1$ and $x_2,y_2\in E_2$,
because
\[
\frac{1}{t}\big(\beta(x_1+ty_1,x_2+ty_2)-\beta(x_1,x_2)\big)
\;=\; \beta(x_1,y_2)+\beta(y_1,x_2)+t\beta(y_1,y_2)
\]
converges to (\ref{diffbilin})
as $t\to 0$, and $d\beta$ is continuous.
A similar argument shows that
each continuous $n$-linear
map $\beta\colon E_1\times \cdots\times E_n\to F$
is $C^1$, with
\[
d\beta((x_1,\ldots, x_n), (y_1,\ldots, y_n))\;=\;
\beta(y_1, x_2,\ldots, x_n)+\cdots+
\beta(x_1,\ldots, x_{n-1}, y_n)\,.
\]
\end{ex}
Let us verify that $df(x,\cdot)\colon E\to F$
is homogeneous.
\begin{lem}\label{homog}
Let $f\colon E\supseteq U\to F$ be a $C^1$-map.
Then
\begin{equation}\label{eqhom}
df(x,sy)\; =s\; df(x,y)
\end{equation}
for all $x\in U$, $y\in E$ and
$s\in \K$.
\end{lem}
\begin{prf}
Because $df(x,0 y)=df(x,0)=0=0 df(x,y)$,
we may assume that
$s\not=0$.
Then
\begin{eqnarray*}
df(x,sy)&=&\vspace{-1.3mm}
\lim_{t\to 0}
\frac{1}{t}\big(f(x+tsy)-f(x)\big)
= s \, \lim_{t\to 0}
\frac{1}{ts}\big(f(x+tsy)-f(x)\big)\\
&=& sdf(x,y).\hspace*{7.4cm}\qedhere
\end{eqnarray*}
\end{prf}
Similar arguments show:
\begin{lem}\label{samecurve}
A $C^0$-curve $\gamma\colon I\to E$
on an open interval $I\sub\R$
is a $C^1$-curve
if and only if $\gamma$
is a $C^1$-map.
In this case,
\begin{equation}\label{dfcurve}
d\gamma(s,r)\; = \; r \, \gamma'(s)\quad
\mbox{for all $\,s\in I$ and $r\in \R$.}
\end{equation}
\end{lem}
\begin{prf}
Let $\gamma$ be a $C^1$-curve. If we can prove~(\ref{dfcurve}),
then $\gamma$ will be~$C^1$.
But (\ref{dfcurve}) can be shown as in the proof of Lemma~\ref{homog}:
Let $s\in I$, $r\in \R$.
Since trivially $d\gamma(s,0)=0$,
we may assume that $r\not=0$.
Then
\[
\frac{1}{t}\big(\gamma(s+tr)-\gamma(s)\big)\,=\,
r\Big( \frac{1}{tr}\big( \gamma(s+tr)-\gamma(s)\big)\Big)
\to r \gamma'(s)
\]
as $t\to 0$ (since this entails $tr\to 0$) and thus $d\gamma(s,r)=r\gamma'(s)$,
as required.\\[3mm]
If, conversely, $\gamma$ is~$C^1$
as a map, then the $C^0$-curve $\gamma$ is differentiable at each $s\in I$ with
\[
\gamma'(s)=\lim_{t\to 0}\frac{\gamma(s+t)-\gamma(s)}{t}=
d\gamma(s,1).
\]
Since $d\gamma(s,1)$ is continuous in~$s$, we deduce that $\gamma$ is a $C^1$-curve.
\end{prf}
The
Mean Value Theorem
carries over to maps
between locally convex spaces:
\begin{prop}[\textbf{Mean Value Theorem}]\label{lemMVTop}
Let $E$ and $F$ be locally\linebreak
convex spaces
and $f\colon U\to F$
be a $C^1$-map on an open subset $U\sub E$.
Then
\begin{equation}\label{formMVTop}
f(y)-f(x)\;=\; \int_0^1 df(x+t(y-x),y-x)\; dt
\end{equation}
for all $x,y\in U$ such that
$U$ contains
$[x,y]:=\{x+t(y-x)\colon t\in [0,1]\}$,
the line segment joining $x$ and~$y$.
\end{prop}
\begin{prf}
Since
$f$ admits all directional
derivatives, we see that
$\gamma\colon [0,1]\to F$, $\gamma(t):=f(x+t(y-x))$
is differentiable at each $t\in [0,1]$,
with derivative
\begin{eqnarray*}
\gamma'(t) \! &\! = \! &\! \lim_{s\to 0} \frac{\gamma(t+s)-\gamma(t)}{s}
= \lim_{s\to 0} \frac{f(x+t(y-x)+s(y-x))-f(x+t(y-x))}{s}\\
& = & df(x+t(y-x),y-x)\,.
\end{eqnarray*}
The preceding formula shows that $\gamma'\colon I\to F$
is continuous, and thus $\gamma'$
is a $C^1$-curve.
Since $\gamma'(t)=df(x+t(y-x),\,y-x)$,
using the the Fundamental Theorem
(Proposition~\ref{fundamental})
to calculate $\gamma(1)-\gamma(0)$,
we obtain~(\ref{formMVTop}).
\end{prf}
Exercise~\ref{excpathnonlcx}
shows that
local convexity of the range~$F$ is essential
for the following property of~$C^1$-maps.
\begin{lem}\label{locconst}
Let $E$ and $F$ be locally convex spaces
and
$f\colon U\to F$ be a $C^1$-map
on an open subset $U\sub E$.
If $df=0$, then $f$ is locally constant.
\end{lem}
\begin{prf}
Given $x\in U$, let $V$ be a convex
neighborhood of~$x$ in~$U$.
For each $y\in V$, we obtain
\[
f(y)-f(x)\;=\; \int_0^1 df(x+t(y-x),\, y-x)\; dt\;=\; 0\,,
\]
using Proposition~\ref{lemMVTop} and the fact that $df=0$.
Thus $f|_V$ is constant.
\end{prf}
\begin{prop}[\textbf{Rule on partial differentials}]\label{rulepartial} 
Let $E_1$, $E_2$ and $F$ be locally convex spaces,
$U \sub E_1\times E_2$ be an open subset
and $f \: U \to F$ be continuous.
Then $f$ is $C^1$ if and only if 
the limits
\[
d_1 f(x_1, x_2, h_1) \, := \, \lim_{t \to 0} \frac{1}{t} 
\big( f(x_1 + t h_1,x_2) - f(x_1, x_2)\big)
\]
and
\[ 
d_2 f(x_1, x_2, h_2) \, := \, \lim_{t \to 0} \frac{1}{t} 
\big( f(x_1,x_2+ t h_2) - f(x_1, x_2)\big)
\]
exist for all $(x_1,x_2)\in U$,
$h_1\in E_1$
and $h_2\in E_2$,
and define continuous mappings
$d_jf\colon U\times E_j\to F$
$($for $j\in \{1,2\})$.
In this case, we have
\begin{equation}\label{eqnpartial}
df(x_1, x_2, h_1,h_2) \; = \; d_1 f(x_1, x_2, h_1) + d_2 f(x_1,
x_2, h_2)
\end{equation}
for all $(x_1,x_2)\in U$ and $h_1\in E_1$, $h_2\in E_2$.
\end{prop}
\begin{prf}
It is obvious that $d_1f$ and $d_2f$ exist
if $f$ is~$C^1$; they are given by
\begin{equation}\label{d1foru}
d_1f(x_1,x_2,h_1)=df(x_1,x_2,h_1,0)\quad\mbox{for $(x_1,x_2)\in U$, $h_1\in E_1$}
\end{equation}
and
\begin{equation}\label{d2foru}
d_2f(x_1,x_2,h_2)=df(x_1,x_2,0,h_2)\quad\mbox{for $(x_1,x_2)\in U$, $h_2\in E_2$.}
\end{equation}
Conversely, assume that $d_1f$ and $d_2f$ exist.
Given $x=(x_1,x_2)\in U$ and
$(h_1,h_2)\in E_1\times E_2$,
there exists $\ve>0$ such that
$x+ \bD_\ve h_1\times \bD_\ve h_2 \sub U$, where
$\bD_\ve:=\{z\in \K\colon |z|\leq \ve\}$.
Then
$x+[0,1]th_1\times [0,1]th_2\sub U$
for each $0\not=t\in \bD_\ve$.
We can write
\begin{eqnarray}
\!\!
\lefteqn{\frac{f((x_1,x_2)+t(h_1,h_2))-f(x_1,x_2)}{t}}\qquad\nonumber\\
&=&
\!
 \frac{f(x_1+th_1, x_2)-f(x_1, x_2)}{t}
\, + \,
\frac{f(x_1+t h_1, x_2+th_2)-f(x_1+th_1, x_2)}{t}.\nonumber \\
& & \label{applpar1}
\end{eqnarray}
The first summand in (\ref{applpar1})
converges to $d_1f(x_1,x_2,h_1)$ as $t\to 0$,\linebreak
by definition of the partial differential.
The second summand in (\ref{applpar1})
can be written as a weak integral,
\begin{equation}\label{rwrpara}
\frac{f(x_1+t h_1, x_2+th_2)-f(x_1+th_1, x_2)}{t}
=
\int_0^1 d_2f(x_1+th_1, x_2+sth_2, h_2)ds,
\end{equation}
as we can
apply Proposition~\ref{lemMVTop} (and Lemma~\ref{homog})
to~$f$ with fixed first argument,
which is~$C^1$ as a consequence of the hypotheses.
Note that the integral in (\ref{rwrpara})
makes sense also for $t=0$
(the integrand is then constant),
and hence defines a function $I_2 \colon \bD_\ve\to F$ of $t$.
The map
\[
\bD_\ve\times [0,1]\to F,\;\,
(t,s)\mto d_2f(x_1+th_1, x_2+sth_2, h_2)
\]
being
continuous, we deduce from Lemma~\ref{intpar}
that $I_2$ is continuous.
Hence the right hand side of (\ref{applpar1})
converges as $t\to 0$,
with limit
\[
d_1f(x_1,x_2,h_1)+I_2(0)=d_1f(x_1,x_2,h_1)+d_2f(x_1,x_2,h_2).
\]
Thus $df$ exists and is given by the right hand side
of (\ref{eqnpartial}) and hence continuous,
whence $f$ is~$C^1$.
\end{prf}
\begin{rem}
The mappings $d_1f\colon U\times E_1\to F$
and $d_2f\colon U\times E_2\to F$ in
Proposition~\ref{rulepartial}
are called the \emph{first} and \emph{second
partial differentials} of~$f$, \index{partial differential} 
respectively.
As in finite-dimensional
analysis, where it is often easier
to calculate partial derivatives
than to find $f'(x)$ directly,
the ``Rule on Partial Differentials''~(\ref{eqnpartial})
is an invaluable tool to calculate
$df$ for complicated mappings on products
(because it is often not too hard to see
what the partial differentials are).
We shall encounter various examples.
\end{rem}
The following characterization of $C^1$-maps
will turn the proof of the Chain Rule
into a triviality.
\begin{lem}\label{linkBGN}
Let $E$ and $F$ be locally convex $\K$-vector spaces,
$U\sub E$ be an open subset and $f\colon U\to F$ be a continuous
map.
Then $f$ is $C^1$ if and only if
the directional difference quotients
\[
\frac{1}{t}\,\big( f(x+ty)-f(x)\big)\,,
\]
which make sense
for all $(x,y,t)\in U\times E\times \K^\times$
such that $x+ty\in U$,
admit a continuous extension to $t=0$.
More precisely,
$f$ is $C^1$ if and only if
there exists a continuous
map $f^{[1]}\colon U^{[1]}\to F$ on
\begin{equation}\label{defnU1}
U^{[1]}\; :=\; \{(x,y,t)\in U\times E\times \K\colon x+ty\in U\}
\end{equation}
such that
\[
f^{[1]}(x,y,t)\; =\; \frac{1}{t}\,\big( f(x+ty)-f(x)\big)
\]
for all $(x,y,t)\in U^{[1]}$ such that $t\not=0$.
\end{lem}
\begin{prf}
If $f^{[1]}$ exists, then
\[
df(x,y)
\,=\, \lim_{t\to 0} \frac{1}{t}
\big(f(x+ty)-f(x)\big)
\,=\, \lim_{t\to 0}
f^{[1]}(x,y,t)\,=\,f^{[1]}(x,y,0)
\]
exists for each $(x,y)\in U\times E$
and defines a continuous map $U\times E\to F$,
as $f^{[1]}$ is continuous.
Hence $f$ is~$C^1$,
with
\begin{equation}\label{reusabl}
df(x,y)\;=\; f^{[1]}(x,y,0)\qquad\mbox{for all $\,(x,y)\in U\times E$.}
\end{equation}
If, conversely, $f$ is~$C^1$, we define
\[
f^{[1]}\colon U^{[1]}\to F\,,\qquad
f^{[1]}(x,y,t)\;:=\;
\left\{
\begin{array}{cl}
\frac{1}{t}\,\big(f(x+ty)-f(x)\big) &\;\mbox{if $\,t\not=0$;}\\
df(x,y) &\; \mbox{if $\, t=0$.}
\end{array}
\right.
\]
Then $f^{[1]}$ is continuous.
In fact, since $f$ is continuous,
the map $f^{[1]}$ is continuous
at each $(x_0,y_0,t_0)\in U^{[1]}$
such that $t_0\not=0$.
Given $(x_0,y_0)\in U\times E$,
there exist $\ve>0$ and neighborhoods
$X$ of~$x_0$ and $Y$ of~$y_0$
in~$E$ such that $X+ \bD_\ve Y\sub U$,
where $\bD_\ve:=\{t\in \K\colon |t|\leq \ve \}$.
Then $Q:=X\times Y\times \bD_\ve\sub U^{[1]}$.
Moreover,
$[x,x+ty]\sub U$ for all $(x,y,t)\in Q$;
indeed, $x+sty\in U$ for all
$s\in [0,1]$ as $(x,y,st)\in X\times Y\times \bD_\ve$.
Using Proposition~\ref{lemMVTop} and Lemma~\ref{homog}, we obtain
\begin{equation}\label{enterscene}
f^{[1]}(x,y,t)\;=\;\int_0^1 df(x+sty,y)\;ds
\end{equation}
for all $(x,y,t)\in Q$ such that
$t\not=0$; if $t=0$ on the other hand,
then (\ref{enterscene})
also holds because the integrand is the constant function
$s\mto df(x,y)$.
The map $Q\times [0,1]\to F$, $(x,y,t,s)\mto df(x+sty,y)$ being
continuous, Lemma~\ref{intpar}
shows that the right hand side
of (\ref{enterscene}) depends continuously
on $(x,y,t)\in Q$, whence
$f^{[1]}$ is continuous at $(x_0,y_0,0)$ in particular.
Thus $f^{[1]}$ is continuous at each point and
hence continuous.
\end{prf}
\begin{lem}\label{difctslin}
If $f\colon E\supseteq U\to F$ is $C^1$, then
$f'(x):=df(x,\cdot)\colon E\to F$ is
a continuous linear map, for each $x\in U$.
\end{lem}
\begin{prf}
The map
$f'(x)\colon E\to F$
is continuous, being a partial map of $df$.
We already verified in Lemma~\ref{homog}
that $f'(x)$ is homogeneous.\\[3mm]
\emph{Additivity}:
Given $x\in U$ and
$y_1,y_2\in E$,
there exists $\ve>0$ such that\linebreak
$x+\bD_\ve y_1+ \bD_\ve y_2\sub U$.
For each $0\not=t\in \bD_\ve$,
we obtain
\begin{eqnarray}
\!\!\!\!\!\!\!\!
\lefteqn{\frac{f(x+t(y_1+y_2))-f(x)}{t}}\qquad\nonumber\\
&=&
\frac{f(x+ty_1)-f(x)}{t}
\; + \;
\frac{f(x+t y_1+t y_2)-f(x+ty_1)}{t} \nonumber \\
&=& f^{[1]}(x,y_1,t) \, +\, f^{[1]}(x+ty_1,y_2,t)\,.
\label{applpar}
\end{eqnarray}
The term in the final line
of (\ref{applpar}) also makes sense
for $t=0$ and is continuous
in $t\in \bD_\ve$ by continuity of~$f^{[1]}$.
The limit $df(x,y_1+y_2)$ of the difference quotients in the first line
of (\ref{applpar}) coincides with the limit of the last line
as $t\to 0$, i.e., with
$f^{[1]}(x,y_1,0)+f^{[1]}(x,y_2,0)
=df(x,y_1)+df(x,y_2)$.
\end{prf}
\begin{rem}\label{doublestnd}
The preceding definition of $f'$
leads to a double meaning of $f'(t)$ if $f\colon I\to F$
is a $C^1$-curve,
corresponding to the linear isomorphism
\[
\cL(\R,F)\to F\,, \quad \alpha\mto \alpha(1)
\]
with inverse $x\mto (\R \ni s\mto s x)$
(where $\cL(\R,F)$ is the space
of linear maps from $\R$ to~$F$).
In fact, $f'(t)$
either means the derivative
$f^{(1)}(t)$ (as in Definition~\ref{defnC1curve}),
or it means $df(t,\cdot)$
(as in Lemma~\ref{difctslin}).
Here the linear map $df(t,\cdot)\colon \R\to F$
is determined by $df(t,s)=sf^{(1)}(t)$
(see (\ref{dfcurve})),
and conversely
$f^{(1)}(t)$ can be recovered from
$df(t,\cdot)$ as
$f^{(1)}(t)=df(t,1)$.
It will always
be clear from the context
which meaning of $f'(t)$ is intended.
\end{rem}
\begin{prop}[\textbf{Chain Rule}]\label{chainC1}
Let $E$, $F$ and $G$ be locally convex spaces,
$U\sub E$, $V\sub F$ be open subsets
and $f\colon U\to F$, $g\colon V\to G$
be $C^1$-maps with $f(U)\sub V$.
Then the composition $g\circ f\colon U\to G$ is $C^1$,
and
\begin{equation}\label{chn1}
d(g\circ f)(x,y)\;=\; dg(f(x), df(x,y))\quad
\mbox{for all $\,(x,y)\in U\times E$,}
\end{equation}
that is,
\begin{equation}\label{chn2}
(g\circ f)'(x)\;=\; g'(f(x))\circ f'(x)\quad\mbox{for all $\, x\in U$.}
\end{equation}
\end{prop}
\begin{prf}
For each $(x,y,t)\in U^{[1]}$ such that $t\not=0$,
we calculate
\begin{eqnarray}
\frac{g(f(x+ty))-g(f(x))}{t} &=&
\frac{g\big(f(x)+t\, \frac{f(x+ty)-f(x)}{t}\big)-g(f(x))}{t}\nonumber\\
&=& g^{[1]}(f(x), f^{[1]}(x,y,t), t)\,,\label{easychain}
\end{eqnarray}
where $f^{[1]}$ and $g^{[1]}$ are as in Lemma~\ref{linkBGN}.
The function $h\colon U^{[1]}\to G$,
$h(x,y,t):=g^{[1]}(f(x), f^{[1]}(x,y,t), t)$
is continuous and extends
the right hand side of (\ref{easychain}).
Hence Lemma~\ref{linkBGN}
shows that $g\circ f$ is $C^1$,
with
\begin{equation}\label{gf1}
(g\circ f)^{[1]}(x,y,t)\;=\;
g^{[1]}(f(x), f^{[1]}(x,y,t), t)\quad\mbox{for all
$(x,y,t)\in U^{[1]}$.}
\end{equation}
In particular, $d(g\circ f)(x,y)=(g\circ f)^{[1]}(x,y,0)=
dg(f(x), df(x,y))$.\vspace{2mm}
\end{prf}
\begin{small}
\subsection*{Exercises for Section~\ref{secC1}}

\begin{exer}\label{excmultilindiff}
Show that every continuous $n$-linear map
$\beta$ is $C^1$ and that $d\beta$
is of the form asserted in Example~\ref{exmultilinC1}.
\end{exer}

\begin{exer}\label{excC1wocts}
Let $E$ and $F$ be locally convex spaces,
$U\sub E$ be open
and\linebreak
$f\colon U\to F$ be a map
such that the directional derivative
$df(x,y)$ exists for all $(x,y)\in U\times E$
and $df\colon U\times E\to F$ is continuous.
Observe that the conclusion of the Mean Value Theorem remains valid and deduce that~$f$ is continuous
(and hence~$C^1$).
\end{exer}

\begin{exer}\label{excprodpartint}
\begin{description}[(D)]
\item[(a)]
(Product Rule).
Let $\beta\colon E_1\times E_2\to F$
be a continuous bilinear map
between locally convex spaces
and $\gamma_1\colon I\to E_1$
as well as $\gamma_2\colon I\to E_2$ be $C^1$-curves, defined on an open interval $I\sub\R$.
Show that
\[
\gamma\colon I\to F\,,\quad
\gamma(t)\, :=\, \beta(\gamma_1(t),\gamma_2(t))
\]
is a $C^1$-curve and $\gamma'(t)=\beta(\gamma_1'(t),\gamma_2(t))
+\beta(\gamma_1(t),\gamma_2'(t))$.
\item[(b)]
Proceed in the same way
if $\gamma:=\beta\circ (\gamma_1,\ldots, \gamma_n)$
for $C^1$-curves
$\gamma_k\colon I\to E_k$
and a continuous $n$-linear map
$\beta\colon E_1\times \cdots\times E_n \to F$.
\item[(c)]
Let~$\cA$ be a locally convex,
unital, associative
topological algebra (see
Definition~\ref{deftopalg}),
$\gamma\colon I\to \cA$ be a $C^1$-curve
and $n\in \N$.
Find $\eta'(t)$ for\linebreak
$\eta\colon I\! \to\!  \cA$, $\eta(t)\! :=\! (\gamma(t))^n$.
When can the formula be simplified as expected\,?
\item[(d)] (Partial Integration).
In the situation of~(a), show that the weak integral
$\int_a^b\beta(\gamma_1(t),\gamma'_2(t))\, dt$ exists in $F$
for given $a,b\in I$ if and only if the weak integral $\int_a^b\beta(\gamma_1'(t),\gamma_2(t))\, dt$
exists in~$F$. In this case,\footnote{As usual, $[f(t)]_a^b:=f(b)-f(a)$.}
\[
\int_a^b\beta(\gamma_1(t),\gamma_2'(t))\,dt =\big[\beta(\gamma_1(t),\gamma_2(t))\big]_a^b
-\int_a^b\beta(\gamma_1'(t),\gamma_2(t))\, dt.
\]
\end{description}
[Openness of $I$ is irrelevant for (a)--(d) (cf.\ Proposition~\ref{chainno} and Lemma~\ref{samecurveno})].
\end{exer}

\begin{exer}\label{chain-individual}
Let $E$, $F$, and~$G$ be locally convex spaces, and $U\sub E$ as well as $V\sub F$ be open subsets.
Let $f\colon U\to F$ be a mapping such that $f(U)\sub V$ and $g\colon V\to G$
be a $C^1$-map. Let $x\in U$ and $y\in E$ be given.
Show that if the directional derivative $df(x,y)=(D_yf)(x)$ of~$f$ at~$x$ exists,
then $d(g\circ f)(x,y)$ exists and is given by $d(g\circ f)(x,y)=dg(f(x),df(x,y))$.\\[.7mm]
[For $0\not=t\in\K$ such that $x+ty\in U$, we have $(g(f(x+ty))-g(f(x))/t=g^{[1]}(f(x),(f(x+ty)-f(x))/t,t)$,
which converges to $g^{[1]}(f(x),df(x,y),0)$ as $t\to 0$.]
\end{exer}

\begin{exer}\label{chain-affine}
Let $E$, $F$, and $G$ be locally convex spaces, and $U\sub E$ as well as $V\sub F$ be open subsets.
Let $f\colon U\to F$ be a mapping such that $f(U)\sub V$,
which is the restriction of a continuous affine-linear map.
Let $x\in U$, $y\in E$ and $g\colon V\to G$
be a mapping whose directional derivative $dg(f(x),df(x,y))$
exists. Show that the directional derivative $d(g\circ f)(x,y)$ exists
and coincides with $dg(f(x),df(x,y))$.\\[.7mm]
[By affine linearity, $f(x+ty)=f(x)+tdf(x,y)$. Hence
$(g(f(x+ty))-g(f(x)))/t=(g(f(x)+tdf(x,y))-g(f(x))/t$.]
\end{exer}

\begin{exer}\label{excfirstpushf}
Given a norm on~$\R^k$,
the supremum norm makes
$C([0,1],\R^k)$
a Banach space.
\begin{description}[(D)]
\item[(a)]
Let
$f\colon \R\to\R$ be continuous
and $(\gamma_n)_{n\in \N}$ a
uniformly convergent
sequence of continuous functions
$\gamma_n\colon [0,1] \to \R$,
with limit $\gamma$.
Exploiting that~$f$ is uniformly
continuous on each compact interval,
deduce that
$f\circ  \gamma_n\to f\circ \gamma$
uniformly. Hence
$\phi := C([0,1], f)\colon
C([0,1],\R)\to C([0,1], \R)$,
$\phi(\gamma):=f \circ \gamma$
is continuous.
\end{description}
We want to see that $\phi$ is $C^1$
if $f$ is $C^1$.
To this end, let
$\gamma,\eta\in C([0,1],\R)$.
\begin{description}[(D)]
\item[\rm (b)]
Assume that
\begin{equation}\label{hint}
d\phi(\gamma,\eta)\;:=\;
\lim_{t\to 0}\frac{1}{t}\big(\phi(\gamma+t\eta)-\phi(\gamma)\big)
\end{equation}
exists.
The
point evaluation
$\ev_x\colon C([0,1],\R)\to\R$,
$\zeta\mto \zeta(x)$
is continuous linear,
for each $x\in [0,1]$.
Apply $\ev_x$ to both sides of
(\ref{hint}) and find the only possible
candidate $\kappa(\gamma,\eta)$
for $d\phi(\gamma,\eta)$
in this way.
\item[\rm(c)]
Show that $\kappa$ is continuous.
\item[\rm (d)]
Applying point evaluations,
verify that $\frac{f(\gamma+t\eta)-f(\gamma)}{t}=
\int_0^1\kappa(\gamma+st\eta,\eta)\, ds$.
\item[\rm (e)]
Show that
$\kappa(\gamma,\eta)$
from~(b)
is indeed the directional derivative
$d\phi(\gamma,\eta)$
and deduce that~$\phi$ is~$C^1$.
\end{description}
\end{exer}
\end{small}
\section{Differentiability of higher order}\label{secCk}
Differentiability of higher order
is defined recursively.
\begin{defn}\label{defnhigher}
Let $E$ and $F$ be locally convex spaces,
$U\sub E$ be an open subset and
$f\colon U\to F$ be a map.
If $f$ is continuous, then we also say that
$f$ is $C^0$. \index{map!$C^k$} 
Let $k\in \N$. We say that $f$ is a $C^k$-map
if $f$ is $C^1$ and
\[
df\colon U\times E\to F
\]
is~$C^{k-1}$ on the open subset $U\times E\sub E\times E$.
We say that $f$ is $C^\infty$ or
\emph{smooth} \index{smooth map} \index{map!smooth} 
if $f$ is $C^k$ for each $k\in \N$.
If we wish to emphasize the ground field $\K$,
we also speak of $C^k_\K$-maps.
\end{defn}
\begin{ex}\label{exctslin}
Each continuous linear map
$\lambda\colon E\to F$
between locally convex spaces
is smooth.
In fact, we know from
Example~\ref{exmultilinC1}
that $\lambda$ is $C^1$,
with $d\lambda$ a continuous linear map.
Given $k\in \N$, the map
$d\lambda$
is $C^{k-1}$ by induction,
and hence $\lambda$ is~$C^k$.
\end{ex}
\begin{lem}[\textbf{Mappings to products}]\label{lemprod}\label{lemprodconv}
Let $E$ be a locally convex space,
$(F_j)_{j\in J}$
be a family of locally convex spaces,
$f\colon U\to \prod_{j\in J}F_j=:F$
be a map on an open subset
$U\sub E$, and
$r\in \N \cup \{\infty\}$.
Let $f_j:=\pr_j\circ f\colon U\to F_j$
for $j\in J$,
where $\pr_j\colon F\to F_j$
is the projection onto the $j$th component.
Then $f$ is $C^r$ if and only if $f_j$ is $C^r$ for each $j\in J$.
In this case,
\begin{equation}\label{diffprod}
df(x,y)\;=\; \big(df_j(x,y)\big)_{j\in J}\quad
\mbox{for all $\,x\in U$ and $y\in E$.}
\end{equation}
\end{lem}
\begin{prf}
Assume that each $f_j$ is $C^r$;
to see that $f$ is $C^r$,
we may assume that $r$ is finite
and proceed by induction.
If $r=1$, $x\in U$ and $y\in E$,
then
\begin{equation}\label{dncequot}
{\textstyle\frac{1}{t}}
\big( f(x+ty)-f(x)\big)
\;=\; \big({\textstyle\frac{1}{t}}(f_j(x+ty)-f_j(x))\big)_{j\in J}
\end{equation}
for each $t\in \K^\times$ such that $x+ty\in U$.
For each $j\in J$,
the $j$th component
$\frac{1}{t}\big(f_j(x+ty)-f_j(x)\big)$
converges to $df_j(x,y)$ as $t\to 0$.
Hence the difference quotient
in (\ref{dncequot})
converges, with limit
$df(x,y)$ as described in~(\ref{diffprod}).
Since $\pr_j\circ df=df_j$
is continuous for each $j\in J$, we conclude that
$df\colon U\times E\to F$ is
continuous.
Hence
$f$ is~$C^1$.
Since $\pr_j\circ df=df_j$
is $C^{r-1}$ for each $j\in J$, the map $df$ is $C^{r-1}$ by the inductive
hypothesis.
Thus $f$ is $C^r$.

If, conversely, $f$ is $C^r$, then $f_j=\pr_j\circ f$ is $C^1$ by the Chain Rule (Proposition~\ref{chainC1})
with $df_j=\pr_j\circ df$ (using that $\pr_j$ is continuous linear).
If $r\geq 2$, then $df_j=\pr_j\circ df$ is $C^{r-1}$ by the inductive hypothesis
(as $df$ is $C^{r-1}$), and thus $f_j$ is~$C^r$.
\end{prf}
\begin{prop}\label{chainCk}
Let $E$, $F$ and $G$ be locally convex spaces,
$U\sub E$ and $V\sub F$ be open subsets, $k\in \N\cup\{\infty\}$
and $f\colon U\to F$ as well as $g\colon V\to G$
be $C^k$-maps.
If $f(U)\sub V$,
then also $g\circ f\colon U\to G$ is $C^k$.
\end{prop}
\begin{prf}
It suffices to prove the assertion
for $k \in \N$.
The proof is by induction.
The case $k=1$ having been settled
in Proposition~\ref{chainC1} (Chain Rule),
we may assume that $k\geq 2$
and that the assertion holds
when~$k$ is replaced with $k-1$.
By Proposition~\ref{chainC1},
$g\circ f$ is $C^1$, and $d(g\circ f)(x,y)=dg(f(x), df(x,y))$
for all $(x,y)\in U\times E$.
Thus
\begin{equation}\label{niceform}
d(g\circ f)\; =\; dg\circ (f\circ \pr_1, df)\, ,
\end{equation}
where $\pr_1\colon U\times E\to U$, $\pr_1(x,y):=x$
is $C^\infty$ as the restriction of a
continuous linear map (Example~\ref{exctslin}).
By induction, $f\circ \pr_1$ is $C^{k-1}$
and hence so is
$(f\circ \pr_1, df)$, by Lemma~\ref{lemprod}.
Applying the inductive hypothesis
to the composition in (\ref{niceform}),
we see that $d(g\circ f)$ is $C^{k-1}$.
Hence $g\circ f$ is~$C^k$.
\end{prf}
\begin{ex}\label{exmultilin}
Every continuous $n$-linear map
$\beta\colon E_1\times\cdots\times E_n\to F$
between locally convex spaces is smooth.\\[3mm]
Indeed, by Example~\ref{exmultilinC1}
the map $\beta$ is~$C^1$, and
\begin{equation}\label{enbind}
d\beta\; = \; \sum_{j=1}^n \beta\circ \lambda_j\, ,
\end{equation}
where
$\lambda_j\colon (E_1\times \cdots\times E_n)^2\to E_1\times\cdots\times E_n$,
\[
\lambda_j(x_1,\ldots, x_n, y_1,\ldots, y_n)\;:=\;
(x_1,\ldots,x_{j-1}, y_j,x_{j+1},\ldots,x_n)
\]
is continuous linear and hence smooth
(Example~\ref{exctslin}).
If~$\beta$ is $C^k$,
using Proposition~\ref{chainCk},
we deduce from (\ref{enbind})
that $d\beta$ is $C^k$,
and thus $\beta$ is $C^{k+1}$.
\end{ex}
We now give an application
to topological algebras (as in Definition~\ref{deftopalg}).
\begin{prop}\label{invalongCk}
Let $\cA$ be a locally convex,
unital, associative
topological algebra
such that the inversion map
$\iota\colon \cA^\times \to \cA$, $\iota(x):=x^{-1}$
is continuous.
Let $E$ be a locally convex space and
$f\colon U\to \cA$ be $C^k$-map
on an open subset $U\sub E$
such that $f(U)\sub \cA^\times$.
Then
$\iota\circ f\colon U\to \cA$, $x\mto f(x)^{-1}$
is~$C^k$.
\end{prop}
\begin{prf}
Assume that $f$ is $C^k$, where $k\in \N$.
Let $x\in U$, $y\in E$.
Using that
\begin{equation}\label{handoninv}
b^{-1}-a^{-1}=b^{-1}(a-b)a^{-1}\;\;
\mbox{for all $a,b\in \cA^\times$,}
\end{equation}
we obtain
for all
$0\not=t\in \K$ such that $x+ty\in U$:
\[
\frac{1}{t}\big( f(x+ty)^{-1}-f(x)^{-1}\big) \; =\; 
-\, f(x+ty)^{-1}\,\frac{f(x+ty)-f(x)}{t} \, f(x)^{-1}\,.
\]
The algebra multiplication being continuous,
the right hand side converges
to $-f(x)^{-1}df(x,y)f(x)^{-1}$
as $t\to 0$. Thus $d(\iota\circ f)(x,y)$ exists
and is given~by
\begin{equation}\label{simplindu}
\hspace*{-2mm}d(\iota \circ\! f)(x,y) \! = 
-f(x)^{-1}df(x,y)f(x)^{-1}\!
=-\tau\big(\iota(f(x)),df(x,y), \iota(f(x))\big),
\end{equation}
where $\tau\colon \cA\times \cA\times \cA\to \cA$, $\tau(a,b,c):=abc$
is continuous trilinear and hence smooth
(Example~\ref{exmultilin}).
By (\ref{simplindu}),
$d(\iota\circ f)$ is continuous,
and hence $\iota\circ f$ is~$C^1$.
We assume now that
$\iota\circ f$ is~$C^{k-1}$,
by induction.
Then
(\ref{simplindu})
shows that $d(\iota\circ f)$ is $C^{k-1}$
as a composition of $C^{k-1}$-maps
(see Proposition~\ref{chainCk} and Lemma~\ref{lemprod}).
As $\iota\circ f$ is~$C^1$ and $d(\iota\circ f)$ is~$C^{k-1}$,
the map $\iota \circ f$ is~$C^k$.
\end{prf}
\begin{defn}\label{defncia}
A locally convex, unital associative topological
$\K$-algebra $\cA$ is called a 
\index{continuous inverse algebra} 
\emph{continuous inverse algebra} (or cia) if its unit group
$\cA^\times$ is open in~$\cA$ and the inversion map
$\iota\colon \cA^\times\to \cA$, $\iota(x):=x^{-1}$
is continuous.
\end{defn}
For example, every unital Banach algebra
is a continuous inverse algebra (Exercise~\ref{excBancia}). Further
examples are discussed in Chapter~\ref{ch:lingrp}.
\begin{cor}\label{invsmoocia}
If $\cA$ is a continuous inverse algebra,
then the inversion map $\iota\colon \cA^\times\to \cA$
is smooth, and
\begin{equation}\label{diffinve}
d\iota(x,y)\; =\; -x^{-1}yx^{-1}
\qquad
\mbox{for all $\, (x,y)\in \cA^\times\times \cA\,$.}
\end{equation}
\end{cor}
\begin{prf}
Applying Proposition~\ref{invalongCk}
to $f\colon \cA^\times \to \cA$, $f(x):=x$,
we see that $\iota=\iota\circ f$
is smooth. The formula for $d\iota$
follows from~(\ref{simplindu}).
\end{prf}
\begin{rem}
The algebra multiplication $\cA\times \cA\to \cA$
in a continuous inverse algebra~$\cA$ being continuous bilinear
and thus smooth, also
the group multiplication $\cA^\times\times \cA^\times\to \cA^\times$
is smooth. Hence, multiplication and
inversion being smooth,
$\cA^\times$ is a $\K$-Lie group
(in the terminology of Chapter~\ref{ch:3}).
\end{rem}
\begin{prop}\label{higherdiff}
Let $E$ and $F$ be locally convex topological $\K$-vector spaces,
$f\colon U\to F$ be a continuous map
on an open subset $U\sub E$
and $r\in \N\cup \{\infty\}$.
Then $f$ is $C^r$ if and only
if the iterated directional derivatives
\[
d^{\,(k)}f(x,y_1,\ldots, y_k)\;:=\; (D_{y_k}\cdots D_{y_1}f)(x)
\]
exist for all $k\in \N$ such that $k\leq r$,
$x\in U$ and $y_1,\ldots, y_k\in E$,
and define
continuous mappings $d^{\,(k)}f\colon U\times E^k\to F$.
\end{prop}
Thus $d^{\, (1)}f=df$ in particular.
We set $d^{\,(0)}f:=f$.
\begin{defn}
$\, d^{\,(k)}f$ is called the 
\index{$k$th differential $d^{(k)}$} 
\index{differential $d^{(k)}$} 
\emph{$k$th differential}
of~$f$.
\end{defn}
The following lemma will help us to prove Proposition~\ref{higherdiff}.
\begin{lem}\label{prehighdf}
Let $E$, $F$ be locally convex spaces, $\ell\in\N$ and $(W_j)_{1\leq j\leq \ell}$ be a family
of locally convex spaces $W_j$ for
$j\in\{1,\ldots,\ell\}$.
Abbreviate $W:=W_1\times\cdots\times W_\ell$.
Let $U\sub E$ be an open subset,
$r\in\N\cup\{\infty\}$
and
\[
f\colon U\times W\to F
\]
be a continuous mapping such that
\[
f(x,\cdot)\colon W_1\times\cdots\times W_\ell\to F
\]
is $\ell$-linear for each $x\in U$,
the iterated directional derivatives
\[
d_1^{\,(k)}f(x,w,y_1,\ldots, y_k):=(D_{(y_k,0)}\cdots D_{(y_1,0)}f)(x,w)
\]
exist for all $k\in\N$ with $k\leq r$,
$x\in U$, $w\in W$ and $y_1,\ldots, y_k\in E$,
and
\[
d_1^{\,(k)}f\colon U\times W\times E^k\to F
\]
is continuous.
Then $f$ is a $C^r$-map.
\end{lem}
\begin{prf}
We may assume that $r<\infty$;
the proof is by induction on $r\in \N$.
If $r=1$, then $d_1f$ exists and is continuous, by hypothesis.
Since $f(x,\cdot)$ is $\ell$-linear and continuous, also
$d_2f(x,w,h)$
exists for all $x\in U$, $w=(w_1,\ldots,w_\ell)\in W$
and $h=(h_1,\ldots, h_\ell)\in W$, and is given by
\begin{equation}\label{d2fformul}
d_2f(x,w,h)=\sum_{j=1}^\ell f(x,w_1,\ldots, w_{j-1},h_j,w_{j+1},\ldots, w_\ell)
\end{equation}
(see Example~\ref{exmultilin}).
Thus
\begin{equation}\label{eninducto}
d_2f=\sum_{j=1}^\ell f\circ (\id_U\times \lambda_j)
\end{equation}
with the continuous linear functions
\[
\lambda_j\colon W\times W\to W,\quad (w,h)\mto(w_1,\ldots, w_{j-1},h_j,w_{j+1},\ldots, w_\ell).
\]
Since each of the maps $\id_U\times \lambda_j$
is continuous, we deduce from (\ref{eninducto})
that $d_2f$ is continuous. Hence $f$ is $C^1$, by
the Rule on Partial Differentials
(Proposition~\ref{rulepartial}), with
\begin{eqnarray}
df(x,w,y,h)&=&d_1f(x,w,y)+d_2f(x,w,h)\nonumber\\
&=&d_1f(x,w,y)+\sum_{j=1}^\ell f(x,\lambda_j(w,h)).\label{thatsdf}
\end{eqnarray}
If $r\geq 2$, then $f$ is $C^{r-1}$ be the inductive hypothesis.
Since
\[
\id_U\times\lambda_j
=(\pr_1,\lambda_j\circ \pr_2)
\]
is smooth (where $\pr_1$ and $\pr_2$ is the projection from $U\times W$
onto $U$ and $W$, respectively),
we deduce with Proposition~\ref{chainCk}
that each of the summands $f\circ (\id_U\times \lambda_j)$ in (\ref{thatsdf})
is $C^{r-1}$.
Now
\[
d_1f(x,w,y)=
\lim_{t\to 0}\frac{1}{t}(f(x+ty,w)-f(x,w))
\]
is $\ell$-linear in $w=(w_1,\ldots, w_\ell)$ as a pointwise limit
of functions with this property. Moreover, $d_1f(x,w,y)=d(f(\cdot,w))(x,y)$ is linear
in~$y$. Hence
\[
d_1f\colon U\times (W_1\times \cdots\times W_\ell\times E)\to F
\]
is a continuous function such that $d_1f(x,\cdot)$ is $(\ell+1)$-linear for each
$x\in U$,
\begin{equation}\label{thscts}
d_1^{\,(k)}(d_1f)(x,(w,y), y_1,\ldots, y_k)=d_1^{\,(k+1)}f(x,w,y,y_1,\ldots, y_k)
\end{equation}
exists for all $k\in \N$ with $k\leq r-1$, $x\in U$, $(w,y)\in W\times E$ and $y_1,\ldots, y_k\in E$,
and such that $d^{(k)}_1(d_1f)$ is continuous (as is clear from (\ref{thscts})).
By the inductive hypothesis, $d_1f$ is $C^{r-1}$.
As continuous linear maps are $C^\infty$, we deduce with Proposition~\ref{chainCk} that
\[
U\times W\times E\times W\to F,\quad (x,w,y,h)\mto d_1f(x,w,y)
\]
is $C^{r-1}$. Since all summands are $C^{r-1}$, we now conclude from (\ref{thatsdf})
that $df$ is $C^{r-1}$. Hence $f$ is $C^r$, which completes the inductive proof.
\end{prf}
\noindent
{\em Proof of Proposition}~\ref{higherdiff}.
We may assume that $r<\infty$.
Let us show first that the higher differentials up to order $r$ exist and are continuous
if $f$ is $C^r$,
by induction on $r\in \N$.
If $r=1$, then $d^{\,(1)}f=df$ exists and is continuous.
If $r\geq 2$, then $df$ is $C^{r-1}$
and hence $df$ has continuous higher differentials $d^{\,(k)}(df)$
for all $k\in\N$ such that $k\leq r-1$.
If $k\in \{2,\ldots, r\}$, $x\in U$ and $y_1,\ldots, y_k\in E$,
then
\[
(D_{y_k}\cdots D_{y_1}f)(x)=D_{y_k}\cdots D_{y_2}(df(\cdot,y_1))(x)
=(D_{(y_k,0)}\cdots D_{(y_2,0)}(df))(x,y_1)
\]
and thus
\begin{equation}\label{dodgehgr}
d^{\,(k)}f(x,y_1,\ldots, y_k)=d^{\,(k-1)}(df)((x,y_1),(y_2,0),\ldots, (y_k,0)),
\end{equation}
which is a continuous $F$-valued function of $(x,y_1,\ldots,y_k)\in U\times E^k$.\\[2.3mm]
Conversely, assume that $f$ is continuous and that continuous higher differentials
exist up to order $r\in \N$. Let us show that $f$ is $C^r$, by induction on $r\in\N$.
If $r=1$, then $f$ is continuous and $df=d^{\,(1)}f$ exists and is continuous,
whence $f$ is $C^1$.
If $r\geq 2$, then $f$ is $C^1$ by the base of the induction
and
\[
df\colon U\times E\to F
\]
is a continuous function which is linear in its second argument
and such that
\[
d_1^{\,(k)}(df)=d^{\,(k+1)}f
\]
exists and is continuous, for all $k\in \N$ such that $k\leq r-1$.
Hence $df$ is $C^{r-1}$ (by Lemma~\ref{prehighdf})
and thus $f$ is $C^r$.\qed
\begin{rem}
We mention that
\[
d^{\,(k)}f\colon U\times E^k\to F
\]
is a $C^{r-k}$-map if $f\colon U\to F$ is a $C^r$-map on an open subset $U\sub E$
with $r\in \N$, and $k\in \{1,\ldots, r\}$.\\[2mm]
[By the recursive definition of a $C^r$-map, $df$ is $C^{r-1}$
and so the claim holds if $k=1$. Thus $d^{\,(k-1)}(df)$ is $C^{r-1-(k-1)}=C^{r-k}$
by induction, if $k\geq 2$. As
\[
d^{\,(k)}f(x,y_1,\ldots,y_k)=d^{\,(k-1)}(df)((x,y_1),(y_2,0),\ldots, (y_k,0))
\]
by (\ref{dodgehgr}), we conclude with Proposition~\ref{chainCk} that $d^{\,(k)}f$ is $C^{r-k}$.]
\end{rem}
\begin{rem}\label{lemhigherdprod}
In the situation of Lemma~\ref{lemprod},
we have
\[
d^{\,(k)}f=(d^{\,(k)}f_j)_{j\in J}
\]
for all $k\in\N$ such that $k\leq r$.
In fact, the lemma subsumes the case $k=1$.
By induction, we have
\[
d^{\,(k)}(df)=(d^{\,(k)}(df_j))_{j\in J}
\]
for each $k\in\N$ such that $k\leq r-1$ and thus
\begin{eqnarray*}
d^{\,(k+1)}f(x,y_1,\ldots, y_{k+1})&=& d^{\,(k)}(df)((x,y_1),(y_2,0),\ldots, (y_{k+1},0))\\
&=&(d^{\,(k)}(df_j)((x,y_1),(y_2,0),\ldots,(y_{k+1},0)))_{j\in J}\\
&=& (d^{\,(k+1)}f_j(x,y_1,\ldots, y_{k+1}))_{j\in J},
\end{eqnarray*}
which completes the inductive argument.
\end{rem}
Our next goal is a version of the Theorem of Hermann Amandus Schwarz
for $C^k$-maps. The following lemma on differentiability of parameter-dependent integrals
will be useful in the proof.
\begin{lem}\label{babydiffpar}
Let $E$ be a locally convex space, $U\sub \K$ be a convex open subset, $t_0\in U$ and
$f\colon U\times [a,b]\to E$ be a continuous function, where $a<b$.
Assume that the partial derivative
\[
\frac{\partial f}{\partial t}(t,s):=(D_{(1,0)}f)(t,s)
\]
exists and is a continuous $E$-valued function of $(t,s)\in U\times [a,b]$.
Moreover, assume that
the weak integral
\begin{equation}\label{intder}
\int_a^b\frac{\partial f}{\partial t}(t_0,s)\, ds
\end{equation}
exists in~$E$, as well as the weak integrals
$\gamma(t):=\int_a^b f(t,s)\, ds$
for all $t\in U$. Then $\gamma\colon U\to E$
is differentiable at~$t_0$, with
$\gamma'(t_0)=\int_a^b\frac{\partial f}{\partial t}(t_0,s)\, ds$.
\end{lem}
\begin{prf}
For $0\not=r\in U-t_0$, we have
\begin{eqnarray*}
\Delta_r &:=& \frac{\gamma(t_0+r)-\gamma(t_0)}{r}
=\int_a^b\frac{f(t_0+r,s)-f(t_0,s)}{r}\,ds\\
&=&\int_a^b\int_0^1 \frac{\partial f}{\partial t}(t_0+\theta r,s)\,d\theta\, ds
\end{eqnarray*}
by the Mean Value Theorem (applied to the functions $f(\cdot, s)$).\footnote{If $\K=\C$, compare
Lemma~\ref{singlevarcx} with $k=1$ for further details.}
For $r=0$, the inner integral exists as well and coincides with the value of its integrand
$\frac{\partial f}{\partial t}(t_0,s)$, which is independent of~$\theta$.
Hence also the outer integral $\Delta_0$ exists\linebreak
and coincides with
(\ref{intder}). By Lemma~\ref{intpar},
the parameter-dependent integral $U-t_0\to E$, $r\mto \Delta_r$ is continuous.
We deduce that $\gamma'(t_0)=\lim_{r\to 0}\Delta_r=\Delta_0$ exists and
coincides with (\ref{intder}).
\end{prf}
Note that $f$ need not be $C^2$ in the following lemma
(since neither existence nor continuity of $\frac{\partial^2 f}{\partial s^2}$
and $\frac{\partial^2 f}{\partial t^2}$
is required).
\begin{lem}\label{stro2nd}
Let $E$ be a locally convex space, $U\sub \K^2$ be an open subset
and $f\colon U\to E$, $(s,t)\mto f(s,t)$
be a continuous function such that the partial derivatives
\[
\frac{\partial f}{\partial s},\quad
\frac{\partial f}{\partial t}\quad\mbox{and}\quad\;
\frac{\partial^2 f}{\partial s\partial t}
\]
exist and are continuous functions from $U$ to~$E$.
Then also $\frac{\partial^2 f}{\partial t\partial s}$ exists\vspace{-.4mm}
and coincides with $\frac{\partial^2 f}{\partial s\partial t}$.
\end{lem}
\begin{prf}
If $(s,t)\in U$, there is $\ve>0$ such that $(s,t+r)\in U$
for all $r\in\K$ with $|r|<\ve$. For $r\not=0$,
we have
\begin{equation}\label{toshorten}
\frac{f(s,t+r)-f(s,t)}{r}=\int_0^1\frac{\partial f}{\partial t}(s,t+\theta r)\,d\theta
\end{equation}
by the Mean Value Theorem.
Considering both sides of (\ref{toshorten}) as elements in a completion $\wt{E}$ of~$E$ with $E\sub\wt{E}$
to ensure the existence of weak integrals, Lemma~\ref{babydiffpar} enables us to differentiate under the integral sign; we obtain
\begin{equation}\label{bothsgoo}
\frac{\frac{\partial f}{\partial s}(s,t+r)-\frac{\partial f}{\partial s}(s,t)}{r}
=\int_0^1\frac{\partial^2 f}{\partial s \partial t}(s,t+\theta r)\,d\theta.
\end{equation}
As the left hand side is in~$E$, so is the right hand side and thus the weak integral also exists in~$E$.
For $r=0$, the integrand on the right hand side is independent of~$\theta$,
whence the weak integral exists in~$E$ also in this case.
Now the right hand side of (\ref{bothsgoo})
converges to $\int_0^1 \frac{\partial^2 f}{\partial s \partial t}(s,t)\,d\theta
=\frac{\partial^2 f}{\partial s \partial t}(s,t)$ as $r\to 0$.
Hence also the left hand side converges, showing that
$\frac{\partial^2 f}{\partial t \partial s}(s,t)$ exists and coincides with the limit
$\frac{\partial^2 f}{\partial s\partial t}(s,t)$ of the right hand side.
\end{prf}
Interchanging adjacent directional derivatives in turn, we readily deduce:
\begin{prop}[\textbf{Schwarz' Theorem}]\label{schwarz}
Let $E$ and $F$ be locally convex topological $\K$-vector spaces,
$U\sub E$ be an open subset, $r\in \N_0\cup\{\infty\}$
and $f\colon U\to F$ be a $C^r_\K$-map.
Let $k\in \N$ with $k\leq r$.
Then
\[
d^{\,(k)}f(x,\cdot)\colon E^k\to F
\]
is a continuous, symmetric $k$-linear map, for each $x\in U$.
\end{prop}
\begin{prf}
We observe first that $d^{\,(k)}f(x,\cdot)$
is continuous as a partial map of $d^{\,(k)}f$.
Now
$d^{\,(k)}f(x,y_1,\ldots, y_k)
=d (d^{\,(k-1)}f(\cdot,y_1,\ldots,y_{k-1}))(x,y_k)$
is linear in~$y_k$,
by Lemma~\ref{difctslin}.
Hence $d^{\,(k)}f(x,\cdot)$
will be $k$-linear
if we can show that $d^{\,(k)}f(x,\cdot)$
is symmetric. We show this by induction on~$k$.
If $k=2$, then
\begin{eqnarray*}
d^{\,(2)}f(x,y_1,y_2)&=&
\frac{\partial^2}{\partial s\partial t}\Big|_{s=t=0}f(x+ty_1+sy_2)
=\frac{\partial^2}{\partial t\partial s}\Big|_{s=t=0}f(x+ty_1+sy_2)\\
&=& d^{\,(2)}f(x,y_2,y_1)
\end{eqnarray*}
by Lemma~\ref{stro2nd}.
If $k>2$, then
\begin{eqnarray*}
d^{\,(k)}f(x,y_1,\ldots, y_{k-2},y_k,y_{k-1})&=&
(D_{y_{k-1}}D_{y_k})(d^{(k-2)}f(\cdot,y_1,\ldots,y_{k-2}))(x)\\
&=&
(D_{y_k}D_{y_{k-1}})(d^{(k-2)}f(\cdot,y_1,\ldots,y_{k-2}))(x)\\
&=& d^{\,(k)}f(x,y_1,\ldots, y_k)
\end{eqnarray*}
by the case just treated, using that $d^{\,(k-2)}f(\cdot,y_1,\ldots,y_{k-2})$ is $C^2$.
By the inductive hypothesis,
$d^{(k-1)}f(x,y_1,\ldots,y_{k-1})$ is unchanged  if $y_j$ and $y_{j+1}$
are swapped for some $j\in \{1,\ldots, k-2\}$.
Applying $D_{y_k}$, we find that also $d^{\,(k)}f(x,y_1,\ldots, y_k)$
remains unchanged. Since every permutation $\pi$ of $\{1,\ldots, k\}$
can be written as a product of transpositions of adjacent numbers,
we deduce that $d^{\,(k)}f(x,y_{\pi(1)},\ldots, y_{\pi(k)})=d^{\,(k)}f(x,y_1,\ldots,y_k)$.
\end{prf}
We mention that there is a simple explicit formula for the higher differentials
$d^{\,(k)}(g\circ f)$ of a composition, in term of those of~$f$ and~$g$.
To formulate it, for $k\in\N$ and $j\in\{1,\ldots, k\}$, we write $P_{k,j}$ for the set
of all partitions $P=\{I_1,\ldots, I_j\}$ of the set $\{1,\ldots, k\}$
into $j$ non-empty disjoint subsets $I_1,\ldots, I_j\sub\{1,\ldots, k\}$.
Thus $I_a\cap I_b=\emptyset$ if $a\not=b$ and $\{1,\ldots, k\}=\bigcup_{a=1}^j I_a$.
We write $|I_a|$ for the number of elements of the set~$I_a$.
If $E$ is a vector space, $y=(y_1,\ldots,y_k)\in E^k$ and $I\sub \{1,\ldots, k\}$ a non-empty subset,
say $I=\{i_1,\ldots, i_\ell\}$ with $i_1<\cdots<i_\ell$,
we abbreviate
\[
y_I:=(y_{i_1},\ldots, y_{i_\ell})\in E^\ell.
\]
\begin{thm}[Fa\`{a} di Bruno's Formula]\label{faadk}
Let $E$, $F$ and $H$ be locally convex spaces,
$U\sub E$ and $V\sub F$ be open,
$k\in \N$ and
$f\colon U\to V\sub F$,
$g\colon V\to H$ be $C^k$-maps.
Then
\[
d^{\,(k)}(g\circ f)(x,y)=\sum_{j=1}^k\sum_{P\in P_{k,j}}\! d^{\,(j)}g\big(f(x),
d^{\,(|I_1|)}f(x,y_{I_1}),\ldots,d^{\,(|I_j|)}f(x,y_{I_j})\big)
\]
for all $x\in U$ and $y=(y_1,\ldots, y_k)\in E^k$, where $P=\{I_1,\ldots, I_j\}$
and the right hand side is well defined, independent of the order of $I_1,\ldots, I_j$.
\end{thm}
\begin{prf}
The right hand side of Fa\`{a} di Bruno's formula is well defined because, for each $j\in\{1,\ldots, k\}$,
the $j$-linear map $d^{\,(j)}g(f(x),\cdot)$ is
symmetric. 
For $k=1$, we have $d(g\circ f)(x,y)=dg(f(x),df(x,y))$ by the Chain Rule,
which is of the desired form (as $P_{1,1}$ is a singleton and only contains $P=\{\{1\}\}$).
Now assume that $f$ and $g$ are $C^{k+1}$ and assume that $d^{\,(k)}(g\circ f)$
is of the asserted form. If $y=(y_1,\ldots, y_k)\in E^k$ and $y_{k+1}\in E$,
we differentiate for $x\in U$ the summand
\[
d^{\,(j)}g(f(x),
d^{\,(|I_1|)}f(x,y_{I_1}),\ldots,d^{\,(|I_j|)}f(x,y_{I_j}))
\]
indexed by $P=\{I_1,\ldots, I_j\}\in P_{k,j}$ at $x$ in the direction~$y_{k+1}$.
By the Rule on Partial Differentials (Proposition~\ref{rulepartial}),
the Chain Rule (Proposition~\ref{chainC1}) and equations (\ref{diffprod}) as well as (\ref{enbind}), we obtain a sum of $j+1$ terms,
namely
\begin{eqnarray*}
\hspace*{-3cm}\lefteqn{d^{\,(j)}g\big(f(x),
d^{\,(|I_1|)}f(x,y_{I_1}),\ldots,d^{(|I_{a-1}|)}f(x,y_{I_{a-1}}),
d^{\,(|I_a|+1)}f(x,y_{I_a},y_{k+1}),}\qquad\qquad\qquad \quad\\
& & d^{\,(|I_{a+1}|)}f(x,y_{I_{a+1}}),\ldots
d^{\,(|I_j|)}f(x,y_{I_j})\big)
\end{eqnarray*}
for $a\in\{1,\ldots, j\}$ and
\[
d^{\,(j+1)}g\big(f(x),
d^{\,(|I_1|)}f(x,y_{I_1}),\ldots,d^{\,(|I_j|)}f(x,y_{I_j}), d^{\,(1)}f(x,y_{k+1})\big).
\]
Note that $P=\{I_1,\ldots, I_j\}\in P_{k,j}$ gives rise to a set $P'$ of $j+1$ pairwise distinct
partitions of $\{1,\ldots, k+1\}$, namely
\[
\{I_1,\ldots, I_{a-1}, I_a\cup\{k+1\},I_{a+1},\ldots, I_j\}\in P_{k+1,j}
\]
for $a\in \{1,\ldots, j\}$ and $P\cup\{\{k+1\}\}\in P_{k+1,j+1}$.
It is clear that $P'\cap Q'=\emptyset$ if $Q$ is a partition of $\{1,\ldots, k\}$ such that $Q\not=P$.
Moreover, each partition of $\{1,\ldots,k+1\}$ is an element of~$P'$ for some $j\in\{1,\ldots,k\}$
and $P\in P_{k,j}$. The preceding summands are therefore exactly the summands in
\[
\sum_{j=1}^{k+1}\sum_{P\in P_{k+1,j}}d^{\,(j)}g\big(f(x),
d^{\,(|I_1|)}f(x,y_{I_1}),\ldots,d^{\,(|I_j|)}f(x,y_{I_j})\big),
\]
which completes the proof.
\end{prf}

We now record various simple,
but very useful
observations.\\[3mm]
Recall that a subset $A\sub X$ of a topological space~$X$
is called 
\emph{sequentially closed}
if $\lim_{n\to \infty}x_n\in A$
for each sequence $(x_n)_{n\in \N}$
in~$A$ which converges in~$X$.
We shall use the following
lemma mainly in the case of closed vector
subspaces, but mere sequential
closedness is needed for the proof.
\begin{lem}\label{corestr}
Let $E$ and $F$ be locally convex spaces, $F_0\sub F$ be a sequentially closed
vector subspace, $U\sub E$ an open subset
and $f\colon U\to F$
a map
such that $f(U)\sub F_0$.
Let $k\in \N\cup \{\infty\}$.
Then $f$ is $C^k$ if and only
if the corestriction $f|^{F_0}\colon U \to F_0$
is $C^k$.
\end{lem}
\begin{prf}
If $f|^{F_0}$ is~$C^k$,
then also $f=\lambda\circ f|^{F_0}$ is~$C^k$,
the inclusion map $\lambda\colon F_0\to F$ being continuous
linear and hence~$C^k$.
For the proof of the converse,
we may assume that~$k$ is finite,
and proceed by induction.

If $f$ is $C^1$,
given $x\in U$ and $y\in E$,
we pick a sequence $(t_n)_{n\in \N}$ in~$\K^\times$
such that $t_n\to 0$ as $n\to\infty$
and $x+t_ny\in U$
for each~$n$. Then
\[
df(x,y)\;=\;\lim_{n\to \infty}\, \frac{1}{t_n}\big(f(x+t_ny)-f(x)\big)\;
\in \, F_0\,,
\]
because
each of the difference quotients
is contained in~$F_0$ and~$F_0$ is sequentially closed in~$F$.
It readily follows that $\frac{1}{t}(f(x+ty)-f(x))\to df(x,y)$
in~$F_0$.
The map $(df)|^{F_0}$ being continuous,
we deduce that $f|^{F_0}$ is $C^1$
with
$d(f|^{F_0})=(df)|^{F_0}$.\\[3mm]
If $f$ is~$C^k$,
$d(f|^{F_0})=(df)|^{F_0}$ is $C^{k-1}$ by induction
and hence $f$ is~$C^k$.
\end{prf}
Our next lemma
deals with mappings to projective limits
of locally convex spaces (see Definition~\ref{defprosys}).
We shall use it frequently as a tool.
\begin{lem}\label{lemPL}
Let $E$ and $F$ be locally convex spaces,
$U\sub E$ be an open subset,
$f\colon U\to F$ be a map,
and $k\in \N\cup\{\infty\}$.
Assume that $F=\pl F_j$\vspace{-.7mm} for a projective
system $((F_i)_{i\in I},(q_{ij})_{i\leq j})$
of locally convex spaces
and continuous linear maps $q_{ij}\colon F_j\to F_i$,
with limit maps $q_i\colon F\to F_i$.
Then $f$ is $C^k$ if and only if $q_i\circ f\colon U\to F_i$
is $C^k$ for each $i\in I$.
\end{lem}
\begin{prf}
If $f$ is $C^k$, then so is $q_i\circ f$,
the limit map $q_i$ being continuous linear and hence smooth.
For the proof of the converse,
we may assume that $F$ is realized
as a closed vector subspace of $\prod_{i\in I}F_i=:P$
and $q_i=\pr_i|_F$.
If $q_i\circ f=\pr_i\circ f$ is $C^k$
for each $i\in I$,
then~$f$ is $C^k$ as a map into~$P$
(by Lemma~\ref{lemprod})
and hence $C^k$ as a map into~$F$,
by Lemma~\ref{corestr}.
\end{prf}
The reader may wish to consult
Appendix~\ref{secfinaltop}
for the definition
and basic properties of quotient maps,
as well as Proposition~\ref{baconlcx}(c).
\begin{lem}\label{lemquot}
Let $E$ and $F$ be locally convex spaces,
$N$ be a closed vector subspace of~$E$
and $q\colon E\to E/N=:E_1$ be the quotient map.
Let $k\in \N_0\cup\{\infty\}$
and $f_1\colon U_1\to F$
be a map, defined on an open subset $U_1\sub E_1$.
Let $U\sub E$ be an open subset such that $q(U)=U_1$.
Then $f_1$ is $C^k$ if and only if $f:=f_1\circ q|_U^{U_1}\colon
U\to F$ is~$C^k$.
\end{lem}
\begin{prf}
The continuous linear map~$q$ being smooth,
$f$ will be $C^k$ if so is~$f_1$.
For the converse, we may assume that $k\in \N_0$;
the proof is by induction.\\[3mm]
\emph{The case $k=0$.}
The map $q|_U^{U_1}\colon U\to U_1$ is a continuous
open surjection and hence a quotient map.
Thus, if $f=f_1\circ q|_U^{U_1}$ is continuous,
then so is~$f_1$.\\[3mm]
\emph{Induction step.}
Assume that $f$ is $C^k$ for some $k\geq 1$.
Given $x_1\in U_1$ and $y_1\in E_1$,
we find $x\in U$ and $y\in E$ such that $q(x)=x_1$ and $q(y)=y_1$.
For $t\in \K^\times$ such that $x+ty\in U$, we have
$x_1+ty_1=q(x+ty)\in U_1$ and
\[
\frac{1}{t}\big(f_1(x_1+ty_1)-f_1(x_1)\big)
\,=\,
\frac{1}{t}\big(f(x+ty)-f(x)\big)\,.
\]
As the right hand side converges to $df(x,y)$ as $t\to0$,
we see that $df_1(x_1,y_1)$ exists,
and $df_1(q(x),q(y))=df(x,y)$.
Since $q\times q\colon E\times E\to E_1\times E_1$
is a quotient map and
$(q\times q)(U\times E)=U_1\times E_1$,
the preceding formula enables us to apply
the inductive hypotheses to~$df_1$:
Thus $df_1$ is $C^{k-1}$ and hence~$f_1$
is $C^k$.
\end{prf}
\begin{rem}
(a)
Note that Lemma~\ref{lemquot}
does not require that $N$
be complemented in~$E$ (in the sense of Definition~\ref{defcplsub});
if~$N$ is complemented, the assertion is trivial.
Although every closed vector subspace~$N$
of a Hilbert space~$E$
is complemented by its orthogonal complement $N^\perp:=\{x\in E
\colon (\forall y\in N)\,\langle
x,y\rangle=0\}$, 
already for more general Banach spaces
this is no longer true.
For example, $c_0(\N,\R)$ is uncomplemented
in $\ell^\infty(\N,\R)$
(see \cite[Satz~IV.6.5]{Wer95}).\medskip

\noindent
(b)
If $N\sub E$ is uncomplemented,
then the quotient map $q\colon E\to E/N$
does not admit a local $C^1$-section
around~$0$
(a $C^1$-map $\sigma\colon U\to E$
on an open $0$-neighborhood $U\sub E/N$
such that $q\circ \sigma=\id_U$).
In fact, otherwise the image
of the linear map $\sigma'(0)\colon E/N\to E$
would be a complement for~$N$ in~$E$.
Nonetheless, if $E$ is a Fr\'{e}chet space,
then $q\colon E\to E/N$ always admits a
\emph{continuous} global section
$\sigma\colon E/N\to E$
and hence is a \emph{topological} $N$-principal bundle,
by Michael's Selection Theorem
(see \cite{MicE59} for the Banach case,
\cite[Ch.\,II, \S4.7, Prop.\,12]{Bou87}
for the general~result).
\end{rem}
Parameter-dependent integrals can be differentiated as expected.
The next result (which will mainly be used in Section~\ref{seccxan} for the study of complex analytic maps)
follows from Lemma~\ref{babydiffpar} by a straightforward induction.
\begin{prop}[\textbf{Differentiation under the integral sign}]\label{diffpar}
Let $E$ and $F$ be locally convex
spaces, $U\sub E$ be an open subset,
$I\sub \R$ be an interval, $a,b\in I$,
$k\in \N$
and $f\colon U\times I \to F$ be a map
such that $f(\cdot,t)\colon U\to F$
is $C^k$, for each $t\in I$.
We assume that the mappings
\begin{equation}\label{ctsreq}
d_1^{\,(j)}f \colon
U\times I \times E^j  \to F\,,\quad
d_1^{\,(j)} f(x,t,y):=d^{\,(j)} f(\cdot,t)(x,y)
\end{equation}
are continuous for each $j\leq k$,
and we assume
that the weak integral
\[
g(x)\;:=\; \int_a^b f(x,t)\, dt
\quad\mbox{exists in~$F$}
\]
for each $x\in U$,
as well as the weak integral
\[
\int_a^b d_1^{\, (j)} f(x,t, y)\, dt,\quad\mbox{for all $j\leq k$, $x\in U$ and $y\in E^j$.}
\]
Then $g\colon U\to F$ is a $C^k$-map and
\begin{equation}\label{intononc}
d^{\, (j)}g (x,y)\, =\, \int_a^b d_1^{\,(j)} f(x,t,y)\, dt\;\; \,
\mbox{for all $\, j\leq k$, $x\in U$ and $y\in E^j$.}
\end{equation}
\end{prop}
\begin{rem}\label{impspecsit}
Note that
$f(\cdot,t)\colon U\to F$
is $C^k$ and the mappings
in (\ref{ctsreq}) are continuous
in the situation
of Proposition~\ref{diffpar},
in the following cases:\medskip

\noindent
(a)
$I\sub\R$ is open, $\K=\R$ and the
map $f\colon U\times I\to F$
is $C^k$; then
\[
d_1^{\,(j)}f(x,t,y_1,\ldots, y_j)
\, =\, d^{\,(j)}f((x,t),\, (y_1,0),\, \ldots, \,(y_j,0))
\]
for all $j\leq k$, $x\in U$ and $y_1,\ldots, y_j \in E$.\medskip

\noindent(b)
$\K=\C$
and $f=h|_{U\times I}$
for a $C^k_\C$-map $h\colon U\times W\to F$,
where $W\sub \C$
is an open subset
such that $I\sub W$.
Then
\[
d_1^{\,(j)}f(x,t,y_1,\ldots, y_j)
\, =\, d^{\,(j)} h((x,t),\, (y_1,0),\, \ldots,\, (y_j,0))\,.
\]
\noindent
(c)
$f(x,t)=\beta(c(t),h(x,t))$,
where $X$ and $Y$ are
locally
convex spaces,\linebreak
$c\colon  I \to X$ is a continuous map,
$\beta\colon  X\times Y\to F$
is continuous bilinear,
and $h\colon U\times W\to Y$
is a $C^k$-map, where $W$ is an
open subset of~$\C$ with
$I\sub W$
if $\K=\C$,
resp., $W=I$ open in~$\R$ if $\K=\R$.
Then
\[
d_1^{\,(j)}f(x,t,y_1,\ldots, y_j)
\, =\, \beta\Big( c(t),\,
d^{\,(j)} h\big((x,t),\, (y_1,0),\, \ldots,\, (y_j,0)\big)\Big)\,.
\]
Proposition~\ref{diffpar} is particularly easy to use
if~$F$ is sequentially complete,
as the
weak integrals
automatically exist in this case
(by Proposition~\ref{rieman}).
\end{rem}
\begin{small}
\subsection*{Exercises for Section~\ref{secCk}}

\begin{exer}\label{exchigherdlin}
Let $\lambda\colon E\to F$
be a continuous linear map
and $\beta\colon E_1\times E_2\to F$ be a continuous
bilinear map,
where $E$, $E_1$, $E_2$ and $F$
are locally convex spaces.
Using the formulas
for $d\lambda$ and $d\beta$
provided above,
calculate the following:
\begin{description}[(D)]
\item[\rm (a)]
$d^{\,(2)}\lambda$
and finally $d^{(k)}\lambda$
for all $k\geq 2$.
\item[\rm (b)]
$d^{\,(2)}\beta$,
$d^{\,(3)}\beta$ and finally $d^{\,(k)}\beta$
for all $k\geq 2$.
\end{description}
\end{exer}

\begin{exer}\label{excfindimCk}
Let $f\colon U\to\R$
be a function on an open subset
$U\sub \R^n$, and $k\in \N$.
Show by induction that $f$ is $C^k$
in the sense of Definition~\ref{defnhigher}
if and only if
all partial derivatives of~$f$
up to order~$k$
exist and are continuous.
In this case,
\[
d^{\,(k)}f(x,y_1,\ldots, y_k)\;=\;
\sum_{j_1,\ldots,j_k=1}^n
y_{1,j_1}\cdots y_{k,j_k} \, \frac{\partial^k f}{\partial x_{j_k}
\cdots \partial x_{j_1}} (x)
\]
for all
$x=(x_1,\ldots, x_n)\in U$
and $y_i=(y_{i,1},\ldots, y_{i,n})\in \R^n$
for $i=1,\ldots, k$.
Can the range $\R$ be replaced by an arbitrary
locally convex space\,?
\end{exer}

\begin{exer}\label{excCkmaptoFo}
Let $E$ and $F$ be locally convex spaces,
$F_0\sub F$ be a vector subspace
and $f\colon U\to F$
be a $C^k$-map on an open
subset $U\sub E$
such that
$\im(d^{\,(j)}f)\sub F_0$
for all $j\in \N_0$ such that $j\leq k$.
Using the characterization from Proposition~\ref{higherdiff},
show that the co-restriction
$f|^{F_0}\colon U \to F_0$
is $C^k$, with
$d^{(j)}(f|^{F_0})=(d^{(j)}f)|^{F_0}$
for all $j\in \N$ with $j\leq k$.
\end{exer}

\begin{exer}\label{excpassfindim}
Let $E$ and $F$ be locally convex spaces,
$f\colon U\to F$ be a $C^k$-map
on an open subset $U\sub E$,
where $k\in \N$,
and $x\in U$, $y_1,\ldots, y_k\in E$.
There is an open $0$-neighborhood
$V\sub \K^k$ such that $s(z):=x+z_1y_1+\cdots+z_ky_k\in U$
for all $z=(z_1,\ldots, z_k)\in V$.
Verify that $h\colon V\to F$, $h(z):=
f(s(z))$ is $C^k$
and
$\frac{\partial^j h}{\partial z_{i_j}\cdots \partial z_{i_1}}
(z)\;=\; d^{\,(j)}f(s(z),y_{i_1},\ldots, y_{i_j})$
for all $j\in\{1,\ldots, k\}$ and $i_1,\ldots,i_j\in \{1,\ldots,k\}$.
In particular,
\begin{equation}\label{forlater}
\frac{\partial^k h}{\partial z_k\cdots \partial z_1}
(0)\;=\; d^{\,(k)}f(x,y_1,\ldots, y_k)\,.
\end{equation}
\end{exer}

\begin{exer}\label{excCkloc}
Let $E$ and $F$ be locally convex spaces,
$U\sub E$ be an open subset,
$f\colon U\to F$ be a map
and $(U_i)_{i\in I}$ be an
open cover for~$U$
(a family
of open subsets $U_i$ of~$U$
with union~$U$).
Show that $f$ is $C^k$ if and only
$f|_{U_i}$ is $C^k$ for each $i\in I$.
Thus being $C^k$ is a local property.
\end{exer}

\begin{exer}\label{exccartpromap}
Let $E_1$, $E_2$, $F_1$ and $F_2$ be locally
convex space, $U_1\sub E_1$ and $U_2\sub E_2$
be open subsets
and $f_1\colon U_1\to F_1$,
$f_2\colon U_2\to F_2$ be $C^k$-maps.
Show
that $f_1\times f_2 \colon U_1\times U_2\to F_1\times F_2$,
$(x_1,x_2)\mto (f_1(x_1), f_2(x_2))$ is $C^k$.
\end{exer}

\begin{exer}\label{exc-explicithdmult}
Given $m\in\N$, let $W_1,\ldots, W_m$ and $F$ be locally convex spaces
and $\beta\colon W_1\times\cdots\times W_m\to F$
be a continuous $m$-linear map.
For $j\in\N$, let $\Theta_{m,j}$ be the set of all functions
\[
\theta\colon \{1,\ldots,m\}\to \{0,\ldots,j\}
\]
such that $\theta^{-1}(i)$
is a singleton for all $i\in\{1,\ldots,j\}$.
Show that
\[
d^{\,(j)}\beta(y_0,\ldots, y_j)=\sum_{\theta\in\Theta_{m,j}}\beta(y_{\theta(1),1},\ldots, y_{\theta(m),m})\vspace{-.4mm}
\]
for all $j\in \N$ and $y_0,\ldots, y_j\in W_1\times\cdots\times W_m$,
where $y_i=(y_{i,1},\ldots, y_{i,m})$ with $y_{i,a}\in W_a$
for $i\in\{0,\ldots,j\}$ and $a\in\{1,\ldots,m\}$.
In particular, $\Theta_{m,j}=\emptyset$ and hence $d^{\,(j)}\beta=0$ for all $j>m$.
\end{exer}

\begin{exer}\label{excbilinctso}
Let $E_1, E_2,\ldots, E_k$, $E$ and $F$
be topological vector spaces.
\begin{description}[(D)]
\item[(a)]
Using that
$\beta(y_1,y_2)-\beta(x_1,x_2)=\beta(y_1-x_1,y_2)
+\beta(x_1,y_2-x_2)$
for all $x_1,y_1\in E_1$,
$x_2,y_2\in E_2$,
show that a bilinear map $\beta\colon E_1\times E_2\to F$
is continuous if and only
if it is continuous at~$(0,0)$.
(Further hints: Use that each $0$-neighborhood
is absorbing and
$\beta(tz_1 ,z_2)=\beta(z_1,tz_2)$ for all $(z_1,z_2)\in E_1\times E_2$).
\item[(b)]
Show that a $k$-linear map
$\beta\colon E_1\times\cdots \times E_k\to F$
is continuous if and only if it is continuous
at~$0$.
\item[(c)]
Now assume that $(E_1,\|\cdot\|_{E_1}),\ldots , (E_k,\|\cdot\|_{E_k})$ and $(F,\|\cdot\|_F)$ are normed spaces.
Show that a $k$-linear map
$\beta\colon E_1\times \cdots\times E_k\to F$ is continuous
if and only if
\[
\|\beta\|_{\op} \, :=\, \sup\{ \|\beta(x_1,\ldots, x_k)\|_F \colon
\mbox{$x_j\in E_j$, $\|x_j\|_{E_j}\leq 1$}\}\; <\; \infty\,.
\]
Then $\|\cdot\|_{\op}$ is a norm on the space $\cL^k(E_1,\ldots, E_k;F)$
of all continuous $k$-linear maps from $E_1\times \cdots\times
E_k$ to~$F$,
and
\[
\|\beta(x_1,\ldots , x_k)\|_F \leq \|\beta\|_{\op} \|x_1\|_{E_1}
\cdots\|x_k\|_{E_k}
\]
holds for all $\beta\in \cL^k(E_1,\ldots, E_k;F)$
and $x_j\in E_j$ (cf.\ Lemma~\ref{bf-multi-cts}).
\end{description}
\end{exer}

\begin{exer}\label{excBancia}
Every unital
Banach algebra~$\cA$
(as in Remark~\ref{whatisban})
is a cia:
\begin{description}[(D)]
\item[(a)]
Show that multiplication is continuous.
Hence $\cA$ is a topological algebra.
\item[(b)]
Show that the ``Neumann series''
$\sum_{k=0}^\infty x^k$ (going back to Carl Gottfried Neumann)
converges uniformly on $B_1^\cA(0)$,
whence $f\colon B_1^\cA(0)\to \cA$,
$f(x):=\sum_{k=0}^\infty x^k$ is continuous.
Calculate $f(x)(\one -x)$ and $(\one -x)f(x)$.
Infer $\one -x\in \cA^\times$ and
\begin{equation}\label{neum1}
(\one -x)^{-1}\, =\; f(x)\;=\; \sum_{k=0}^\infty \,x^k
\;\;\,\mbox{for all $x\in B_1^\cA(0)$.}
\end{equation}
Hence $B_1^\cA(\one )\sub \cA^\times$
and $\iota\colon \cA^\times\to \cA$,
$\iota(x):=x^{-1}$ is continuous on $B_1^\cA(\one )$.
It now easily follows that
$\cA^\times$ is open and $\iota$ is continuous
(see Lemma~\ref{spotcia}).
\item[(c)]
Using (\ref{neum1}), show that $\|(\one-x)^{-1}\|\leq \frac{1}{1-\|x\|}$
for all $x\in B^\cA_1(0)$.
\end{description}
\end{exer}

\begin{exer}\label{exclininvctssmoo}
Let $E$ and $F$ be locally convex spaces, $U\sub E$ open, $k\in \N$
and $A_x\colon F\to F$ be an isomorphism of topological vector spaces
for $x\in U$ such that
\[
f\colon U\times F \to F,\quad f(x,y):=A_x(y)
\]
is $C^k$. Show that if $g\colon U\times F\to F$, $g(x,y):=(A_x)^{-1}(y)$ is continuous,
then $g$ is~$C^k$ [reuse (\ref{handoninv}) and other
ideas from the proof of Proposition~\ref{invalongCk}].
\end{exer}
\end{small}
\section{Differential calculus on non-open domains}\label{secnoop}
Except for the case of curves, we only considered $C^k$-maps
on open domains so far.
In this section, we extend the theory
and define and study $C^k$-maps
on suitable
not necessarily open subsets of a locally convex space~$E$.
In particular, our approach subsumes
$C^k$-maps
on cubes
\[
[0,1]^n\sub \R^n,
\]
half-spaces
\[
\R^{n-1}\times [0,\infty[\; \sub \R^n,
\]
and similar sets. This is indispensable
for infinite-dimensional
Lie theory. For example,
our discussions of
``regular'' Lie groups
(a well-behaved class of infinite-dimensional
Lie groups)
will involve mappings on $[0,1]$ and
$[0,1]^2$.
Also mappings on open subsets
of closed half-spaces
are of importance,
since ``manifolds with boundary''
are modeled on such sets.\\[2.3mm]
More precisely, we consider $C^k$-maps $f\colon U\to F$
to a locally convex space~$F$, defined
on a locally convex subset $U\sub E$ with dense interior.
The main point is that a version of the Mean Value Theorem
is still available for such maps. As a consequence,
the Chain Rule and other basic facts can be established
essentially as in the case of open domains.
Let us first define and discuss the new type of domain,
and then turn to the functions thereon.
\begin{defn}\label{lcxset}
Let $E$ be a locally convex space.
A subset $U\sub E$ is called 
\index{locally convex subset}\index{subset!locally convex}
\emph{locally convex}
if every point $x\in U$ has a convex neighborhood~$V$
in~$U$.
\end{defn}
\begin{rem}\label{basiclcxs}
For example, every convex subset
of~$E$ is locally convex,
and so is any open subset.
Furthermore, intersections
of finitely many locally convex
subsets are locally convex.
For example, $C \cap U$ is locally convex,
for every convex subset $C\sub E$
and open subset $U\sub E$.
In particular, every (relatively)
open subset
of a closed half-space in~$\R^n$
is locally convex.
We remark that every connected, closed,
locally convex
subset $U\sub E$ is convex (see Exercise~\ref{excclosedlcx}).
\end{rem}
\begin{rem}\label{useflrm}
In Definition~\ref{lcxset},
we can always achieve that~$V$ is open in~$U$.
In fact, $E$ being locally convex,
there is an open, convex neighborhood~$W$
of~$x$ in~$E$ such that $U\cap W\sub V$.
Then $U\cap W=V\cap W$ is a convex, open neighborhood
of~$x$ in~$U$. Taking $W$ sufficiently
small, we actually see that
each neighborhood $P\sub U$ of~$x$
contains an open, convex neighborhood.
\end{rem}
If $E$ is a locally convex space, we write $U^0$
for the interior (relative to~$E$) of a subset $U\sub E$.
We say that $U$ has dense interior if $U^0$ is dense in~$U$.
\begin{defn}\label{defnCkno}
Let $E$ and $F$ be locally convex topological $\K$-vector spaces
and $U\sub E$ be a locally convex subset
with dense interior.
A map $f\colon U\to F$ is called
$C^1$ if it is continuous,
$f|_{U^0}$ is $C^1$ and $d(f|_{U^0})\colon U^0\times E\to F$
extends to a (necessarily unique)
continuous map
\[
df\colon U\times E\to F.
\]
We say that $f$ is $C^0$ if $f$ is continuous.
If $k\in\N$ with $k\geq 2$, we say that $f$ is $C^k$ if $f$ is $C^1$
and $df$ is $C^{k-1}$ (noting that $U\times E$ is a locally convex subset with dense interior
in $E\times E$). If $f$ is $C^k$ for all $k\in\N_0$,
then we say that $f$ is smooth or $C^\infty$.
\end{defn}
\begin{lem}\label{lemhidffsno}
Let $E$ and $F$ be locally convex spaces, $U\sub E$ be a locally convex subset
with dense interior and $k\in \N\cup\{\infty\}$.
Then $f$ is $C^k$ if and only if
$f|_{U^0}$ is $C^k$ and $d^{\,(j)}(f|_{U^0})$
admits a $($necessarily unique$)$ continuous extension
\[
d^{\,(j)}f\colon U\times E^j\to F
\]
for all $j\in\N$ such that $j\leq k$.
\end{lem}
\begin{prf}
We may assume that $k\in\N$. If $k=1$,
then the statement holds by definition.
For $k\geq 2$, we now prove the necessity of the second condition
(the sufficiency will be proved after Proposition~\ref{chainno}).
To this end, assume that $f$ is $C^k$.
Then $f|_{U^0}$ is $C^1$ with $d(f|_{U_0})=df|_{U^0\times E}$
a $C^{k-1}$-map and thus $f|_{U^0}$ is $C^k$.
Since
$f$ is $C^{k-1}$ in particular, there exist continuous extensions $d^{\,(j)}f$
for $j\in\{1,\ldots, k-1\}$, by induction.
Also $df$ is $C^{k-1}$, whence (by induction) there is a continuous
extension $d^{\,(k-1)}(df)$ of $d^{\,(k-1)}(df|_{U^0\times E})$.
Recalling that
\[
d^{\,(k)}(f|_{U^0})(x,y_1,\ldots, y_k)=d^{\,(k-1)}(d(f|_{U^0}))((x,y_1),(y_2,0),\ldots,(y_k,0))
\]
for all $x\in U^0$ and $y_1,\ldots, y_k\in E$, we find that
$d^{\,(k)}f\colon U\times E^k\to F$,
\[
(x,y_1,\ldots, y_k)\mto d^{\,(k-1)}(df)((x,y_1),(y_2,0),\ldots,(y_k,0))
\]
is a continuous extension for $d^{\,(k)}(f|_{U^0})$.
\end{prf}
A mapping to a product is $C^r$ if and only if so is each of its components.
\begin{lem}[\textbf{Mappings to products}]\label{lemprodno}
Let $E$ be a locally convex space,
$(F_j)_{j\in J}$
be a family of locally convex spaces,
and $f\colon U\to \prod_{j\in J}F_j=:F$
be a map on a locally convex subset
$U\sub E$ with dense interior.
For $j\in J$, let $\pr_j\colon F\to F_j$
be the projection onto the $j$-th component.
Let $r\in \N \cup \{\infty\}$.
Then $f$ is $C^r$ if and only if $f_j:=\pr_j\circ f\colon U\to F_j$
is $C^r$ for each $j\in J$.
In this case, we have
\begin{equation}\label{dfincomp}
d^{\,(k)}f=(d^{\,(k)}f_j)_{j\in J}\;\;\mbox{for all $k\in \N$ such that $k\leq r$.}
\end{equation}
\end{lem}
\begin{prf}
We may assume that~$r$ is finite.
If each $f_j$ is $C^r$, then $f=(f_j)_{j\in J}$ is continuous.
Moreover, $f_j|_{U^0}$ is $C^r$
and hence $f|_{U^0}=(f_j|_{U^0})_{j\in J}$ is $C^r$,
by Lemma~\ref{lemprod}, with
\[
d(f|_{U^0})=(d(f_j|_{U^0}))_{j\in J}=((df_j)|_{U^0\times E})_{j\in J}.
\]
As the continuous map $df:=(df_j)_{j\in J}$ extends $d(f|_{U^0})$,
we see that $f$ is $C^1$. If $r\geq 2$, then $df=(d(f_j))_{j\in J}$ is $C^{r-1}$
by induction and hence $f$ is $C^r$.\\[2.3mm]
If, conversely, $f$ is $C^r$, then each $f_j$ is continuous.
Moreover, $f|_{U^0}$ is $C^1$, so $f_j|_{U^0}$ is $C^1$. Since
$(df)|_{U^0\times E}=d(f|_{U^0})=(d(f_j|_{U^0}))_{j\in J}$
by Lemma~\ref{lemprod}, we see that the continuous maps $\pr_j\circ df$
extend $d(f_j|_{U^0})$. Hence $f_j$ is $C^1$ and
\begin{equation}\label{indhypapp}
df=(df_j)_{j\in J}.
\end{equation}
Since $df$ is $C^{r-1}$,
we deduce from (\ref{indhypapp})
and the inductive hypothesis that $df_j$ is
$C^{r-1}$. Hence $f_j$ is $C^r$.
\end{prf}
\begin{lem}\label{homogno}
If $E$ and $F$ are locally convex spaces, $U\sub E$ a locally convex subset
with dense interior, $r\in \N\cup\{\infty \}$
and $f\colon U\to F$ is a $C^r$-map,
then
\[
d^{\,(k)}f(x,\cdot)\colon E^k\to F
\]
is a continuous, symmetric $k$-linear map for each $k\in\N$ with $k\leq r$ and~$x\in U$.
\end{lem}
\begin{prf}
Since $d^{\,(k)}f$ is continuous, also $d^{\,(k)}f(x,\cdot)$ is continuous.
As $U^0$ is dense in~$U$ and $d^{\,(k)}f(\cdot,y)$ is continuous
for each $y\in E^k$, the multilinearity and symmetry of $d^{\,(k)}f(x,\cdot)$
for $x\in U$ follow from the corresponding properties for $x\in U^0$.
We give more details for
the linearity of $d^{\,(k)}(x,y_1,\ldots, y_k)$ in~$y_1$;
all other properties are established analogously.
Let $y_j\in E_j$ for $j\in\{1,\ldots,k\}$,
$z_1\in E_1$ and $t\in\K$. Then
\begin{eqnarray*}
x&\mto& d^{\,(k)}f(x,y_1+tz_1,y_2,\ldots, y_k)\qquad\mbox{and}
\\
x&\mto & d^{\,(k)}f(x,y_1,\ldots, y_k)+td^{\,(k)}f(x,z_1,y_2,\ldots,y_k)
\end{eqnarray*}
are continuous functions from $U$ to the Hausdorff space $F$ which coincide on the dense
subset $U^0$ of~$U$. Therefore both functions coincide and thus
$d^{\,(k)}f(x,y_1+tz_1,y_2,\ldots, y_k)
=d^{\,(k)}f(x,y_1,\ldots, y_k)+td^{\,(k)}f(x,z_1,y_2,\ldots,y_k)$
for all $x\in U$.
\end{prf}
It is essential that the
Mean Value Theorem
remains valid for mappings on locally convex subsets
with dense interior.
\begin{prop}\label{lemMVT}
Let $E$ and $F$ be locally convex spaces
and $f\colon U\to F$
be a $C^1$-map on a locally convex subset $U\sub E$
with dense interior.
Then
\begin{equation}\label{formMVT}
f(y)-f(x)\;=\; \int_0^1 df(x+t(y-x),y-x)\; dt
\end{equation}
for all $x,y\in U$ such that
$U$ contains the line segment
$[x,y]$.
\end{prop}
\begin{prf}
Let us first consider the special case
that $U$ is a convex subset
of~$E$ whose interior is non-empty (and hence dense in~$U$,
by Lemma~\ref{baseconvex}(c)).
We choose a
completion $\tilde{F}$ of~$F$ with $F\sub\wt{F}$.
By Lemma~\ref{intpar}
on parameter-dependent integrals,
the function
\[
h\colon U\times U\to \tilde{F}\,,\quad
h(u,v)\, :=\, \int_0^1 df(u+t(v-u), v-u)\, dt
\]
is continuous. By Proposition~\ref{lemMVTop}, we have $h(u,v)=f(v)-f(u)$
for all $u,v\in U^0$.
Since $U^0\times U^0$ is dense
in $U\times U$ and both $h$ and the mapping\linebreak
$(u,v)\mto f(v)-f(u)$ are continuous,
we deduce that $h(u,v)=f(v)-f(u)$
for all $u,v\in U$. In particular, $h(u,v)\in F$
and (\ref{formMVT}) holds for all $x,y\in U$.\\[2.3mm]
\emph{General case}: Let $U\sub E$ be as described in the proposition,
and $x,y\in U$ such that $[x,y]\sub U$.
Each $z\in [x,y]$ has a convex open neighborhood
in~$U$.
A compactness argument
yields a partition
$0=t_0<t_1< \cdots< t_{n-1}<t_n=1$
of $[0,1]$ such that, setting $z_j:=x+t_j(y-x)$,
we have $[z_j,z_{j+1}]\sub V_j$
for some open convex subset~$V_j$ of~$U$.
Since $U$ has dense interior, also each $V_j$ has 
dense interior.
By the special case already treated,
\begin{equation}\label{ifshow}
f(z_{j+1})-f(z_j)\; =\; \int_0^1
df(z_j+s(z_{j+1}-z_j),z_{j+1}-z_j)\,ds
\end{equation}
for each~$j$.
Re-writing $z_{j+1}-z_j=(t_{j+1}-t_j)(y-x)$
in (\ref{ifshow}) yields
\begin{eqnarray*}
f(z_{j+1})-f(z_j) &= &
\,(t_{j+1}-t_j)\!
\int_0^1 \!\! df\bigl(x+(t_j+s(t_{j+1}-t_j))\,(y-x),\, y-x
\bigr) ds\\
&= &
\,\int_{t_j}^{t_{j+1}}df(x+t(y-x), y-x)\,dt\,,
\end{eqnarray*}
where we used the substitution $t=t_j+s(t_{j+1}-t_j)$
to pass to the last line
(see Exercise~\ref{excbasicwint}(d)
for the Substitution Rule)
and used the homogeneity from Lemma~\ref{homogno}.
Since
\[
f(y)-f(x)\; =\; \sum_{j=0}^{n-1}(f(z_{j+1})-f(z_j)),
\]
we now obtain (\ref{formMVT}).
\end{prf}
\begin{lem}\label{linkBGNno}
Let $E$, $F$ be locally convex $\K$-vector spaces,
$U\sub E$ be a locally convex set with
dense interior, and $f\colon U\to F$ be a continuous
map.
Then $f$ is $C^1$ if and only if
there exists a continuous
map $f^{[1]}\colon U^{[1]}\to F$ on
\begin{equation}\label{defnU1no}
U^{[1]}\; :=\; \{(x,y,t)\in U\times E\times \K\colon x+ty\in U\}
\end{equation}
such that
\[
f^{[1]}(x,y,t)\; =\; \frac{1}{t}\,\big( f(x+ty)-f(x)\big)
\]
for all $(x,y,t)\in U^{[1]}$ such that $t\not=0$.
\end{lem}
\begin{prf}
If $f^{[1]}$ exists, then its restriction
to $(U^0)^{[1]}$ is a continuous map $(f|_{U^0})^{[1]}$
extending the directed difference quotients of $f|_{U^0}$.
Hence $f|_{U^0}$ is $C^1$, by Lemma~\ref{linkBGN},
with
\begin{equation}\label{nwsam}
df(x,y)=f^{[1]}(x,y,0)
\end{equation}
for all $(x,y)\in U^0\times E$.
Since (\ref{nwsam}) can be used to define
a continuous function $df\colon U\times E\to F$ which extends
$df|_{U^0}$, we see that $f$ is $C^1$. \\[2.3mm]
If, conversely, $f$ is~$C^1$, we define
\[
f^{[1]}\colon U^{[1]}\to F\,,\qquad
f^{[1]}(x,y,t)\;:=\;
\left\{
\begin{array}{cl}
\frac{1}{t}\,\big(f(x+ty)-f(x)\big) &\;\mbox{if $\,t\not=0$;}\\
df(x,y) &\; \mbox{if $\, t=0$.}
\end{array}
\right.
\]
Then $f^{[1]}$ is continuous.
In fact, since $f$ is continuous,
the map $f^{[1]}$ is continuous
at each $(x_0,y_0,t_0)\in U^{[1]}$
such that $t_0\not=0$.
Given $(x_0,y_0)\in U\times E$,
there exists a convex neighborhood~$V$ of~$x_0$
in~$U$. Let $W$ be a neighborhood
of~$x_0$ in~$E$ such that $W\cap U\sub V$.
There exist neighborhoods
$X$ of~$x_0$ and $Y$ of~$y_0$
in~$E$, and $\ve>0$ such that $X+ \bD_\ve Y\sub W$,
where $\bD_\ve:=\{t\in \K\colon |t|\leq \ve \}$.
Then $X\sub X+\bD_\ve Y\sub W$
and $(X+ \bD_\ve Y)\cap U\sub W\cap U\sub V$.
Hence, for each
$(x,y,t)$
in the neighborhood
$Q:=U^{[1]}\cap (X\times Y\times \bD_\ve)$
of $(x_0,y_0,0)$ in $U^{[1]}$, we have $x,
x+ty\in (X+ \bD_\ve Y)\cap U\sub V$
and hence $[x, x+t y]\sub V$,
since $V$ is convex.
Now (\ref{enterscene})
and the remainder of the proof of Lemma~\ref{linkBGN}
can be copied verbatim.
\end{prf}
As a consequence, also the Chain Rule
remains valid for mappings on locally convex
subsets with dense interior:
\begin{prop}\label{chainno}
Let $E,F,G$ be locally convex spaces,
$U\sub E$, $V\sub F$ be locally convex sets
with dense interior,
and $f\colon U\to F$, $g\colon V\to G$
be $C^k$-maps such that $f(U)\sub V$,
where $k\in \N\cup\{\infty\}$.
Then also $g\circ f\colon U\to G$ is $C^k$, and
\[
d(g\circ f)(x,y)=dg(f(x),df(x,y))\;\;\mbox{for all $(x,y)\in U\times E$.}
\]
\end{prop}
\begin{prf}
We may assume that $k\in \N$ and proceed
by induction.
The case $k=1$ can be proved like
Proposition~\ref{chainC1},
using Lemma~\ref{linkBGNno} instead of Lemma~\ref{linkBGN}.
The induction step can be performed as in the proof of Proposition~\ref{chainCk},
using Lemma~\ref{lemprodno} instead of Lemma~\ref{lemprod}.
\end{prf}
\begin{prop}[\textbf{Rule on partial differentials}]\label{rulepartialno} 
Let $E_1$, $E_2$ and $F$ be locally convex spaces,
$U \sub E_1\times E_2$ be a locally convex subset
with dense interior
and $f \colon  U \to F$ be continuous.
Then $f$ is $C^1$ if and only if 
the limits
\[
d_1 f(x_1, x_2, h_1) \, := \, \lim_{t \to 0} \frac{1}{t} 
\big( f(x_1 + t h_1,x_2) - f(x_1, x_2)\big)
\]
and
\[ 
d_2 f(x_1, x_2, h_2) \, := \, \lim_{t \to 0} \frac{1}{t} 
\big( f(x_1,x_2+ t h_2) - f(x_1, x_2)\big)
\]
exist for all $(x_1,x_2)\in U^0$,
$h_1\in E_1$
and $h_2\in E_2$,
and extend to continuous mappings
$d_jf\colon U\times E_j\to F$
$($for $j\in \{1,2\})$.
In this case, we have
\begin{equation}\label{eqnpartialno}
df(x_1, x_2, h_1,h_2) \; = \; d_1 f(x_1, x_2, h_1) + d_2 f(x_1,
x_2, h_2)
\end{equation}
for all $(x_1,x_2)\in U$ and $h_1\in E_1$, $h_2\in E_2$.
\end{prop}
\begin{prf}
If $f$ is $C^1$, then the limits described in the proposition exist;
moreover, (\ref{d1foru}) and (\ref{d2foru})
define continuous extensions $d_1f$ and $d_2f$,
respectively, and (\ref{eqnpartialno}) holds.
If, conversely, the limits and continuous extensions $d_1f$ and $d_2f$ exist,
then $f|_{U^0}$ is $C^1$
and
\[
d(f|_{U^0})(x_1, x_2, h_1,h_2) \; = \; d_1 f(x_1, x_2, h_1) + d_2 f(x_1,
x_2, h_2)
\]
for all $(x_1,x_2)\in U^0$ and $h_1\in E_1$, $h_2\in E_2$,
by Proposition~\ref{rulepartial}.
Hence, if we define $df\colon U\times (E_1\times E_2)\to F$ via
(\ref{eqnpartialno}), then $df$ is continuous and
extends $d(f|_{U^0})$. Thus $f$ is $C^1$ and (\ref{eqnpartialno}) holds.
\end{prf}
\begin{rem}\label{prehighdfno}
Lemma~\ref{prehighdf}
on mappings $f\colon U\times (W_1\times\cdots \times W_\ell)\to F$
that are $\ell$-linear in the second argument
extends to the case when $U\sub E$ is a locally convex set with dense interior.\\[2.3mm]
Indeed, $d_2f(x,w,h)$ exists and is given by (\ref{d2fformul})
for all $(x,w,h)\in U^0\times W\times W$\linebreak
in this case,
and (\ref{eninducto}) provides the continuous extension $d_2f$.
Using now\linebreak
Lemma~\ref{rulepartialno}
instead of Lemma~\ref{rulepartial},
we see that $f$ is $C^1$ and (\ref{thatsdf}) holds.
As in the proof of Lemma~\ref{prehighdf},
we see that the maps $d_1f(x,\cdot)$ are $\ell$-linear for all $x$
in the dense subset $U^0\sub U$
and hence for all $x\in U$ by continuity of $d_1f$ (cf.\ proof of Lemma~\ref{homogno}).
Replacing Proposition~\ref{chainCk} with Proposition~\ref{chainno},
the remainder of the proof of Lemma~\ref{prehighdf}
carries over.
\end{rem}
\noindent
\emph{Proof of Lemma} \ref{lemhidffsno}, \emph{completed.}
Let $k\geq 2$ and assume that the continuous extensions $d^{\,(j)}f$ exist for
$j\in \{1,\ldots, k\}$.
Then $f$ is $C^1$ by the base of the induction.
The continuous map $df\colon U\times E\to F$ is linear in the second argument
and admits continuous differentials $d_1^{\,(j)}(df)$
for $j\in \{1,\ldots, k-1\}$. Hence $df$ is $C^{k-1}$ (see Remark~\ref{prehighdfno})
and thus~$f$ is~$C^k$.\vspace{2.3mm}\qed

\noindent
Our concept of $C^k$-map is compatible with the terminology for curves:
\begin{lem}\label{samecurveno}
Let $E$ be a locally convex space, $k\in \N$ and $I\sub\R$ be a non-degenerate interval.
Then a map $\gamma\colon I\to E$
is a $C^k$-curve $($as in Definition~{\rm\ref{defnC1curve}}$)$
if and only if $\gamma$
is a $C^k$-map in the sense of
Definition~{\rm \ref{defnCkno}}.
In this case,
\begin{equation}\label{dfcurveno}
d\gamma(s,r)\; = \; r \, \gamma'(s)\quad
\mbox{for all $\,s\in I$ and $r\in \R$.}
\end{equation}
\end{lem}
\begin{prf}
Let $\gamma$ be a $C^1$-curve. By Lemma~\ref{samecurve},
$\gamma|_{I_0}$ is a $C^1$-map with $d(\gamma|_{I^0})(s,r)=r\gamma'(s)$
for all $(s,r)\in I^0\times\R$. As the continuous map $I\times\R\to F$, $(s,r)\mto
r\gamma'(s)$ extends $d(\gamma|_{I^0})$, we see that $\gamma$ is a $C^1$-map and
(\ref{dfcurveno}) holds.\\[3mm]
Conversely, assume that the map~$\gamma$ is~$C^1$.
For each $s\in I$ and
$0\not= t\in I-s$, we then have
$s+[0,1]t\sub I$ and hence
\[
\frac{\gamma(s+t)-\gamma(s)}{t}\;=\;
\int_0^1 d\gamma(s+rt,1)\; dr
\]
by Proposition~\ref{lemMVT}
and Lemma~\ref{homogno}.
The map $h\colon (I-s)\times [0,1]\to E$,
$h(t,r):=d\gamma(s+rt,1)$
being continuous,
the map $g\colon I-s\to E$, $g(t):=\int_0^1 h(t,r)\, dr$
is continuous by Lemma~\ref{intpar}
(noting that the weak integral also
exists for $t=0$, as $h(0,\cdot)$ is constant).
Therefore $\frac{1}{t}(\gamma(s+t)-\gamma(s))=g(t)\to g(0)=h(0,0)=d\gamma(s,1)$
as $t\to 0$, entailing that
$\gamma'(s)$ exists and is given by
$\gamma'(s)=d\gamma(s,1)$ and hence continuous.
Thus $\gamma$ is a $C^1$-curve.
\end{prf}
The following lemma will be used to juxtapose solutions
to differential equations (see, e.g., Remark~\ref{glue-sol-ode}(a)).
In Section~\ref{secmetrcalc}, it will also help
us to find smooth parametrizations
for certain infinite polygons.
\begin{lem}\label{C1glueing}
Let $I\sub \R$ be a non-degenerate interval, $t_0\in I$, $k\in\N\cup\{\infty\}$
and $\gamma\colon I \to E$
be a continuous map to a locally convex space~$E$ such that $\gamma|_{I \setminus\{t_0\}}$
is $C^k$ and
$\lambda_j:=
\lim_{t\to t_0}\gamma^{(j)}(t)$ for $t\not=t_0$ exists
for all $j\in\N$ such that $j\leq k$.
Then $\gamma$ is $C^k$, and $\gamma^{(j)}(t_0)=\lambda_j$
for all $j\in\N$ such that $j\leq k$.
\end{lem}
\begin{prf}
Assume that $t_0$ is an interior point of~$I$ (the excluded cases are similar, and easier).
We may assume that $k<\infty$, and proceed by induction.
The case $k=1$:
By the hypotheses, $\eta:=\gamma|_{I\cap [t_0,\infty[}$
is a $C^1$-map in the sense of
Definition~\ref{defnCkno}
and thus a $C^1$-curve (Lemma~\ref{samecurveno}).
Then $\eta'(t_0)=\lim_{t\to t_0}\gamma'(t)=\lambda_1$.
Likewise,
$\zeta:=\gamma|_{I\cap\, ]{-\infty},t_0]}$ is a $C^1$-curve
such that $\zeta'(t_0)=\lambda_1$.
Since both the left and right one-sided derivatives of~$\gamma$ at~$t_0$
exist and have the same value~$\lambda_1$, the derivative $\gamma'(t_0)$
exists and is given by~$\lambda_1$.
Hence $\gamma'$ is continuous
and thus~$\gamma$ is~$C^1$.
Induction step: Let $k>1$ and assume that the assertion holds for $k-1$ in place of~$k$.
By the case $k=1$, $\gamma$ is $C^1$ with $\gamma'(t_0)=\lambda_1$.
The inductive hypothesis applies to $\gamma'$ in place of $\gamma$ and $\lambda_{j+1}$
in place of $\lambda_j$, for $j\in\{1,\ldots,k-1\}$. Hence $\gamma'$ is $C^{k-1}$ (whence $\gamma$ is
$C^k$) and $\gamma^{(j+1)}(t_0)=(\gamma')^{(j)}(t_0)=\lambda_{j+1}$ for all $j\in\{1,\ldots, k-1\}$,
which completes the inductive proof.
\end{prf}
\begin{rem}\label{nopecarryov}
Many results from calculus on open domains carry over directly
to locally convex domains with dense interior,
or only require spelling out an obvious candidate for a continuous extension.
For example:\medskip

\noindent
(a)
The characterization of locally constant $C^1$-functions
from Lemma~\ref{locconst}
remains valid if $U\sub E$
is a locally convex subset with dense
interior (the proof can be repeated verbatim).\medskip

\noindent
(b)
The conclusions of Proposition~\ref{diffpar} (concerning differentiation under the integral sign)
remain valid if $U\sub E$ is merely a locally convex subset with dense
interior. Indeed, (\ref{intononc}) describes the higher differentials of $g|_{U^0}$;
as the integral on the right hand side of (\ref{intononc})
provides a continuous extension $d^{\,(j)}g\colon U\times E^j\to F$,
Lemma~\ref{lemhidffsno} shows that~$g$~is~$C^k$.\medskip

\noindent
(c) Fa\`{a} di Bruno's Formula
(as in Theorem~\ref{faadk}) remains valid if $U$ and $V$
are merely locally convex subsets with dense interior
(as the formula holds for $(x,y)\in U^0\times E^k$ and its right hand side provides
a continuous extension to all of $U\times E^k$, which has to coincide with $d^{\,(k)}(g\circ f)$).
\end{rem}
The following analog of Lemma~\ref{corestr} is available.
\begin{lem}\label{nonocorestr}
Let $E$ and $F$ be locally convex spaces,
$f\colon U\to F$ be a map on a
locally convex subset
$U\sub E$ with dense
interior, $F_0\sub F$
be a sequentially closed vector subspace
such that $f(U)\sub F_0$,
and $k\in \N\cup \{\infty\}$.
Then $f$ is $C^k$ if and only
if the corestriction $f|^{F_0}\colon U \to F_0$
is $C^k$.
\end{lem}
\begin{prf}
If $f|^{F_0}$ is~$C^k$,
then also $f=\lambda\circ f|^{F_0}$ is~$C^k$,
the inclusion map $\lambda\colon F_0\to F$ being continuous
linear and hence~$C^k$.
If, conversely, $f$ is~$C^k$, then $f|_{U^0}^{F_0}$
is~$C^k$ by Lemma~\ref{corestr}, with
$d^{\,(j)}(f|^{F_0}_{U^0})=(d^{\,(j)}(f|_{U^0}))|^{F_0}$
for all $j\in\N$ such that $j\leq k$.
Given $(x,y)\in U\times E^j$,
we pick a convex neighborhood $V\sub U$ of~$x$.
Since $U^0$ is dense in~$U$, we have
$V^0\not=\emptyset$. Let
$z\in V^0$.
Then \[ x_n:=z+(1-2^{-n})(x-z)\in V^0\sub U^0 \] 
for each $n\in \N$ (cf.\ proof of Lemma~\ref{baseconvex}(b)),
and thus
\[
d^{\,(j)}f(x,y)\; =\; \lim_{n\to\infty} d^{\,(j)}f(x_n,y)\; \in \, F_0
\]
as $d^{\,(j)}f(x_n,y)\in F_0$ for each~$n$ and~$F_0$ is sequentially
closed. Hence $d^{\,(j)}f$ takes its values in~$F_0$
and thus $(d^{\,(j)}f)|^{F_0}$ is a continuous extension
for $d^{\,(j)}(f|^{F_0}_{U^0})$. Hence $f|^{F_0}$ is $C^k$.
\end{prf}
The remainder of this section compiles auxiliary results which are more
spezialized; we recommend to read these results only later,
when they are actually used in the text.

The next lemma links real and complex
differentiability. We shall use it as a tool in the proofs of
Lemmas~\ref{Ckinlimit} and~\ref{twotoplem}, and also in Section~\ref{secCspaces}.
The topic will be taken up in Proposition~\ref{analytsingle}, Corollary~\ref{corcxcurve2}
and Theorem~\ref{charcxcompl}.
See Section~\ref{seccxan} for a more comprehensive discussion of complex differentiable
maps.
\begin{lem}\label{realvscxearly}
Let $E$ and $F$ be complex locally convex spaces, $U\sub E$ be a locally convex subset with dense interior
and $k\in \N\cup\{\infty\}$. Then a map
$f\colon U\to F$ is $C^k_\C$ if and only if~$f$ is~$C^k_\R$
and $df(x,\cdot)\colon E\to F$ is complex linear for each $x\in U$.
\end{lem}
\begin{prf}
It is clear from the definition that every $C^k_\C$-map~$f$ is $C^k_\R$,
with the same differentials $d^{\,(j)}f$ for $j\in\N_0$ such that $j\leq k$.
In particular, $df(x,\cdot)$ is complex linear by Lemma~\ref{homogno}.

Conversely, assume that $f$ is $C^k_\R$ and $df(x,\cdot)$ is $\C$-linear for each $x\in U$.
Then an easy induction shows that $d^{\,(j)}f(x,y_1,\ldots, y_j)$
is $\C$-linear in~$y_1$ for all $j$ as before, $x\in U$ and $y_1,\ldots, y_j\in E$
(being a limit of difference quotients which are $\C$-linear in $y_1$).
Since $d^{\,(j)}f(x,\cdot)\colon E^j\to F$ is symmetric, we deduce that $d^{\,(j)}f(x,\cdot)$
is complex $j$-linear. It remains to show that 
\[
d^{\,(j)}f(x,y_1,\ldots, y_j)=D_{y_j}(d^{\,(j-1)}f(\cdot,y_1,\ldots, y_{j-1}))(x)
\]
as a complex directional derivative, for all
$x\in U^0$ and $y_1,\ldots,y_j\in E$. We have $x+\bD_r y_j\sub U$ for some $r>0$.
By the Mean Value Theorem,
\begin{eqnarray*}
\Delta_z & :=& \frac{d^{\,(j-1)}f(x+zy_j,y_1,\ldots y_{j-1})-d^{\,(j-1)}f(x,y_1,\ldots, y_{j-1})}{z}\\
&=&\frac{1}{z}\int_0^1d^{\,(j)}f(x+tzy_j,y_1,\ldots, y_{j-1},zy_j)\,dt
\end{eqnarray*}
for $z\in\bD_r\setminus\{0\}$.
By $\C$-linearity of $d^{\,(j)}f$ in its final argument, $z$ cancels.
Thus
\[
\Delta_z=\int_0^1d^{\,(j)}f(x+tzy_j,y_1,\ldots, y_j)\,dt,
\]
which converges to $\int_0^1d^{\,(j)}f(x,y_1,\ldots, y_j)\, dt=d^{\,(j)}f(x,y_1,\ldots, y_j)$
as $z\to 0$, by continuous parameter-dependence of integrals (Lemma~\ref{intpar}).
\end{prf}
Let us check that the $C^k$-property passes to limits,
in certain situations. This will be used in the proof of Lemma~\ref{Ckellpartial}.
\begin{lem}\label{Ckinlimit}
Let $X$ be a topological space, $x_0\in X$ be an accumulation point,\footnote{Thus $U\setminus\{x_0\}\not=\emptyset$
for each neighborhood $U\sub X$ of~$x_0$.}
$E$ and $F$ be locally convex spaces, $V\sub E$ be a locally convex subset with dense
interior, $k\in\N\cup\{\infty\}$ and $f\colon X \times V\to F$
be a continuous map. For $x\in X$, abbreviate $f_x:=f(x,\cdot)\colon V\to F$.
Assume that
\begin{description}[(D)]
\item[{\rm(a)}]
For each $x\in X \setminus\{x_0\}$, the map $f_x$ is $C^k$; and
\item[{\rm(b)}]
For each $j\in\N$ such that $j\leq k$, there exists a continuous mapping
$g_j\colon X \times V\times E^j\to F$ such that $g_j(x,y,w)=d^{\,(j)}(f_x)(y,w)$
holds for all $x\in X \setminus\{x_0\}$, $y\in V$ and $w\in E^j$.
\end{description}
Then also $f_{x_0}$ is $C^k$, and $d^{\,(j)}(f_{x_0})=g_j(x_0,\cdot)$ for all $j\in\N$ with $j\leq k$.
\end{lem}
\begin{prf}
If $\K=\C$ and the assertion holds in the real case, then $f_{x_0}$ is $C^k_\R$
and $d(f_{x_0})(y,w)=g_1(x_0,y,w)$ is complex linear in $w\in E$
(being the limit of $g_1(x,y,w)=d(f_x)(y,w)$ as $x\to x_0$ with $x\in X\setminus \{x_0\}$).
Hence $f_{x_0}$ is $C^k_\C$, by Lemma~\ref{realvscxearly}.

We may therefore assume $\K=\R$ now. We may also assume that $k\in\N$, and proceed by induction.
If $k=1$, $y\in V^0$ and $w\in E$,
there is $r>0$ such that $y+\,]{-r},r[\, w\sub V^0$.
Let $\wt{F}$ be a completion of~$F$ with $F\sub\wt{F}$.
The map
\[
h_1\colon X \times \,]{-r},r[\, \to\wt{F},\quad
(x,t)\mto f(x,y)+\int_0^1 g_1(x,y+stw,tw)\, ds
\]
is continuous, by Lemma~\ref{intpar}.
Also
\[
h_2\colon X\times \,]{-r},r[\,\to\wt{F},\quad (x,t)\mto f(x,y+tw)
\]
is continuous.
Now $h_1(x,t)=h_2(x,t)$ for all $x\in X\setminus\{x_0\}$ and $t\in\,]{-r},r[$
by the Mean Value Theorem. Since $(X\setminus\{x_0\})\times \,]{-r},r[$ is dense
in $X\times \,]{-r},r[$, this implies that $h_1=h_2$.
Hence
\begin{eqnarray*}
f_{x_0}(y+tw)&=&f(x_0,y)+\int_0^1 g_1(x_0,y+stw,tw)\,ds\\
&=& f(x_0,y)+\int_0^t g_1(x_0,y+\tau w,w)\,d\tau,
\end{eqnarray*}
where we substituted $\tau=st$.
The second part of the Fundamental Theorem (Proposition~\ref{secpart})
now shows that $]{-r},r[\,\to F$, $t\mto f_{x_0}(y+tw)$ is $C^1$,
with $\frac{d}{dt}\big|_{t=0}f_{x_0}(y+tw)=g_1(x_0,y,w)$.
In particular, $d(f_{x_0})(y,w)$ exists and coincides with $g_1(x_0,y,w)$.
Since $g_1$ is continuous,
the assertion follows.

If $k\geq 2$ and the assertion holds for $k-1$ in place of~$k$, we apply the same argument to
the function
\[
d^{\,(k-1)}(f_x)(y,w_1,\ldots, w_{k-1})=g_{k-1}(x,y,w_1,\ldots, w_{k-1})
\]
of $x\in X$ and $y\in V^0$ (with fixed $w_1,\ldots, w_{k-1}\in E$)
to see that the directional derivative $d^{\,(k)}(f_{x_0})(y,w_1,\ldots, w_k)$ exists for each $w_k\in E$,
and coincides with $g_k(x_0,y,w_1,\ldots, w_k)$. Since $g_k$ is continuous, the assertion follows.
\end{prf}
The next lemma will only be used once
(in the proof of Proposition~\ref{companotions}).
\begin{lem}\label{twotoplem}
Let $\cO_1$ and $\cO_2$ be two topologies
on a vector space~$F$ making it a locally convex space,
such that $\cO_1\sub\cO_2$ and $(F,\cO_2)$ is sequentially complete
$($or integral complete$)$.
\begin{description}[(D)]
\item[\rm(a)]
Let $\gamma\colon [a,b]\to (F,\cO_1)$ be a $C^1$-curve
such that $\gamma'$ is continuous as a map $[a,b]\to (F,\cO_2)$.
Then $\gamma$ is also $C^1$ as a map to $(F,\cO_2)$,
with the same derivative.
\item[\rm(b)]
Let $E$ be a locally convex space, $U\sub  E$ be a locally convex subset with dense interior,
$k\in\N$ and $f\colon U\to (F,\cO_1)$
be a $C^k$-map which is $C^{k-1}$ as a map to $(F,\cO_2)$.
If $d^{\,(k)}f$ is continuous as a map to $(F,\cO_2)$,
then $f$ is $C^k$ as a map to $(F,\cO_2)$.
\end{description}
\end{lem}
\begin{prf}
(a) Let $x\in [a,b]$. Since $\gamma'$ is continuous as a map to $(F,\cO_2)$ and $(F,\cO_2)$
is integral complete, the weak integral
\[
\eta(x):=\int_a^x \gamma'(t)\, dt
\]
exists in $(F,\cO_2)$. The element $z:=\gamma(x)-\gamma(a)$ satisfies
$\lambda(z)=\int_a^x\lambda(\gamma'(t))\,dt$ for all $\lambda\in (F,\cO_1)'$,
by Proposition~\ref{fundamental}.
Since $(F,\cO_1)'$ is a subset of $(F,\cO_2)'$,
we have $\lambda(\eta(t))=\int_a^x\lambda(\gamma'(t))\,dt$
for all $\lambda\in (F,\cO_1)'$. As $(F,\cO_1)'$ separates points on~$F$,
we infer $z=\eta(t)$ (see 
Exercise~\ref{excpsepwint}).
Thus $\gamma(x)=\gamma(a)+\eta(x)$.
Now $\eta$ is $C^1$ as a map to $(F,\cO_2)$, by Proposition~\ref{secpart},
with $\eta'=\gamma'$.
Hence $\gamma=\gamma(a)+\eta$ is $C^1$ as a map to $(F,\cO_2)$,
with $\frac{d}{dt}(\gamma(a)+\eta)=\eta'=\gamma'$.

(b) Write $g$ for~$f$, considered as a mapping to $(F,\cO_2)$.
Then $f=\lambda\circ g$ with the continuous linear map $\lambda\colon (F,\cO_2)\to(F,\cO_1)$,
$y\mto y$ and thus
\begin{equation}\label{usepresen}
d^{\,(k-1)}f=\lambda\circ d^{\,(k-1)}g =d^{\,(k-1)}g,
\end{equation}
by the Chain Rule.
If $\K=\C$ and we can show that $g$ is $C^k_\R$, then $dg(x,\cdot)=df(x,\cdot)$
will be complex linear for all $x\in U$ and thus $g$ will be $C^k_\C$,
by Lemma~\ref{realvscxearly}.
We may therefore assume that $\K=\R$ now.
It suffices to show that $d^{(k)}g(x,y_1,\ldots, y_k)$
exists for $x\in U^0$ and $y_1,\ldots, y_k\in E$ and coincides with $d^{\,(k)}f(x,y_1,\ldots, y_k)$;
then $d^{\,(k)}g:=d^{\,(k)}f\colon U\times E^k\to (F,\cO_2)$ will be a continuous extension
and so~$g$ will be~$C^k$.
There is $r>0$ such that $x+[{-r},r]y_k\sub U$.
Then
$\gamma\colon [{-r},r]\to (F,\cO_1)$,
\[
t\mto d^{\,(k-1)}g(x+ty_k,y_1,\ldots, y_{k-1})=d^{\,(k-1)}f(x+ty_k,y_1,\ldots, y_{k-1})
\]
is a $C^1$-curve, noting that $\gamma'(t)=d^{\,(k)}f(x+ty_k,y_1,\ldots, y_k)$.
Since $d^{\,(k)}f$ is continuous as a map to $(F,\cO_2)$, so is~$\gamma'$.
Thus (a) shows that $\gamma$ is $C^1$ as a map to $(F,\cO_2)$.
In particular,
\[
d^{\,(k)}g(x,y_1,\ldots, y_k)=\gamma'(0)
\]
exists in $(F,\cO_2)$ and coincides with $d^{\,(k)}f(x,y_1,\ldots, y_k)$.
\end{prf}

\begin{small}
\subsection*{Exercises for Section~\ref{secnoop}}

\begin{exer}\label{excclosedlcx}
Let $E$ be a locally convex space and
$A\sub E$ be a closed, connected,
locally convex set. Show that $A$ is convex,
as follows:
\begin{description}[(D)]
\item[(a)]
Let $E=\R^2$ first
and assume that $[0,e_1]\cup [0,e_2]\sub A$,
where $e_1=(1,0)$,
$e_2=(0,1)$,
and where $[u,v]$ is the line segment
joining $u,v \in E$.
Show that, for each $\theta\in [0,1[$
such that $[\theta e_1, e_2]\sub A$,
there is $\ve\in \; ]0,1-\theta]$
such that $[(\theta+r)e_1,e_2]\sub A$
for $r\in [0,\ve]$.
Show that $\{\theta\in [0,1]\colon
[\theta e_1, e_2]\sub A\}=:T$
is closed in~$\R$. Deduce that $1\in T$
and thus $[e_1,e_2]\sub A$.
\item[(b)]
Given $x\in E$,
let $S_x:=\{y\in A\colon [x,y]\sub A\}$.
Show that $S_x$ is closed.
Using (a), show that $S_x$ is convex.
Let $x,z\in A$ such that
$S_x\cap S_z\not=\emptyset$.
Using~(a),
show that $z\in S_x$ and in fact $S_x=S_z$.
Deduce that $A=S_x$ for each $x\in A$.
\end{description}
\end{exer}

\begin{exer}\label{excLeibru}
Let $I$ be a non-degenerate interval, $\beta\colon E_1\times E_2\to F$ be a\linebreak
continuous
bilinear map between locally convex spaces, $n\in\N_0$ and $\gamma_1\colon I\to E_1$
as well as $\gamma_2\colon I\to E_2$ be $C^n$-curves. Show the Leibniz Rule:
\[
\frac{d^n}{dt^n} \,\beta(\gamma_1(t),\gamma_2(t))=\sum_{k=0}^n\left(\begin{array}{c}
n\\
k
\end{array}\right) \beta(\gamma^{(k)}_1(t),\gamma_2^{(n-k)}(t)).
\]
\end{exer}

\begin{exer}\label{seq-cts-Leib}
Show that the conclusion of Exercise~\ref{excLeibru}
remains valid if~$\beta$ is merely sequentially continuous
(compare Corollary~\ref{hypo-seq}, Corollary~\ref{bar-hypo},
and Example~\ref{exa-eval} for such maps).
\end{exer}

\begin{exer}\label{fin-diff-cts}
Let $E$ be a finite-dimensional $\K$-vector space with basis $v_1,\ldots, v_n$
and $F$ be a locally convex space. Show that the map
$f'\colon U\to\cL(E,F)_b$, $x\mto f'(x)$ is continuous for each $C^1$-map $f\colon U\to F$
on a locally convex subset $U\sub E$ with dense interior.\\[.7mm]
[Using the isomorphism $\Phi\colon \cL(E,F)_b\to F^n$ of locally convex spaces
from Exercise~\ref{exc-on-fin}, we have $(\Phi\circ f')(x)=(df(x,v_j))_{j\in\{1,\ldots, n\}}$.]
\end{exer}

\begin{exer}\label{exc-chainpt}
Let $E$ and $F$ be locally convex spaces, $U\sub E$ be a locally convex subset with dense interior
and $f\colon U\to F$ be a $C^1$-map. Let $I\sub\R$ be a non-degenerate interval, $t\in I$
and $\gamma\colon I\to U$ be a map which is differentiable at~$t$ in the sense that
\[
\gamma'(t)=\lim_{s\to t}\frac{\gamma(s)-\gamma(t)}{s-t}
\]
(with $s\in I\setminus\{t\}$) exists in~$F$. Using that
\[
\frac{f(\gamma(s))-f(\gamma(t))}{s-t}=f^{[1]}\left(
\gamma(t),\frac{\gamma(s)-\gamma(t)}{s-t},s-t\right)
\]
for $s\in I\setminus \{t\}$, deduce that $f\circ\gamma\colon I\to F$ is differentiable at~$t$ and
$(f\circ\gamma)'(t)=df(\gamma(t),\gamma'(t))$.
\end{exer}

\begin{exer}\label{exc-patho-nonop}
\begin{description}[(D)]
%
%
\item[(a)]
We consider the map $f\colon\R^2\to\R$, 
$(x,y)\mapsto x^2+y^2$. Then $[1,4]\sub\R$
is a locally convex subset with dense interior;
show that $f^{-1}([1,4])$ is not locally convex.\smallskip
\item[(b)]
Let $\alpha\colon E\to F$ be a continuous linear map between
locally convex spaces and $Q\sub F$ be a locally convex subset.
Show that $\alpha^{-1}(Q)$ is locally convex.\\[1.2mm]
[Given $x\in\alpha^{-1}(Q)$, let $P\sub F$ be a convex $x$-neighborhood such that
$Q\cap P$ is convex. Then $\alpha^{-1}(P\cap Q)=\alpha^{-1}(Q)\cap \alpha^{-1}(P)$
is convex and an $x$-neighborhood in $\alpha^{-1}(Q)$.]\smallskip
\item[(c)]
Consider the continuous linear map $\lambda\colon\R\to\R^2$, $x\mto (x,0)$.
The disk $K:=\{(x,y)\in\R^2\colon (x-1)^2+y^2\leq 1\}$ is a locally convex
subset of $E:=\R^2$ and has
dense interior. Show that the locally convex set
$\lambda^{-1}(K)$ has empty interior.
Moreover, the intersection $K\cap F$ with the closed vector subspace
$F:=\{0\}\times\R$ has empty interior in~$F$.
\end{description}
\end{exer}
\end{small}
\section{Lipschitz continuity and related topics}\label{sec-lip}
Lipschitz continuity and Lipschitz conditions play an important role
in differential calculus, for example in the study of fixed points and their parameter dependence;
in connection with the inverse function theorem; and in the\linebreak
theory of ordinary
differential equations.
In this section, we first record some simple basic facts
concerning the relation between Lipschitz constants for a function
and the operator norm of its derivatives (Lemma~\ref{lipviaprime}).
We then introduce Lipschitz maps between subsets of locally convex spaces
and prove several more specialized results,
which may be skipped
on a first reading.
\begin{defn}
A map $f\colon X\to Y$ between metric spaces $(X,d_X)$ and $(Y,d_Y)$
is called (globally) 
\emph{Lipschitz} \index{Lipschitz map} 
if there exists $L\in [0,\infty[$ (a so-called \emph{Lipschitz constant})
such that
\begin{equation}\label{the-lip}
d_Y(f(x),f(y))\leq L\, d_X(x,y)\quad\mbox{for all $\, x,y\in X$.}
\end{equation}
The latter holds if and only if
\[
\Lip(f):=\sup\left\{
\frac{d_Y(f(x),f(y))}{d_X(x,y)} \colon \mbox{$x,y\in X$ with $x\not=y$} \right\}\in[0,\infty]
\]
is finite, and $\Lip(f)$ is the minimum Lipschitz constant in this case.
The map~$f$ is called 
\index{locally Lipschitz map}
\emph{locally Lipschitz}
if each $x_0\in X$ has a neighborhood $U\sub X$ such that $f|_U$ is Lipschitz.
\end{defn}
Every Lipschitz map is continuous, and so is every locally Lipschitz map.
\begin{rem}
We are most interested in the special case that
$X$ and $Y$ are subsets of normed spaces $(E,\|\cdot\|_E)$ and $(F,\|\cdot\|_F)$,
respectively, $d_X(x,y):=\|y-x\|_E$ for $x,y\in X$, and~$d_Y$ is defined analogously using~$\|\cdot\|_F$.
Then (\ref{the-lip}) turns into
\[
\|f(y)-f(x)\|_F\leq L \|y-x\|_E\quad\mbox{for all $x,y\in X$,}
\]
and $\Lip(f)$ is the supremum over $\frac{\|f(y)-f(x)\|_F}{\|y-x\|_E}$ with $x,y\in X$ such that $x\not=y$.
\end{rem}
The following lemma is essential for differential calculus.  
\begin{lem}\label{lipviaprime}
Let $(E,\|\cdot\|_E)$ and $(F,\|\cdot\|_F)$ be normed spaces
and $f\colon U\to F$ be a mapping on a locally convex subset $U\sub E$
with dense interior.
\begin{description}
\item[\rm(a)]
If the directional derivative $df(x,y)=(D_yf)(x)$
exists for some $x\in U^0$ and $y\in E$, then
$\|df(x,y)\|_F\leq\Lip(f) \|y\|_E$.
\item[\rm(b)]
If $f$ is $C^1$, then
$\|f'(x)\|_{\op} \leq \Lip(f)$ for all $x\in U$.
\item[\rm(c)]
If $f$ is $C^1$ and~$U$ is convex, then
\begin{equation}\label{linkderlip}
\Lip(f)=\sup\{\|f'(x)\|_{\op}\colon x\in U\}.
\end{equation}
\item[\rm(d)]
If $f$ is $C^1$, then $f$ is locally Lipschitz.
\end{description}
\end{lem}
\begin{prf}
(a) We have $x+ty\in U$ for small $t\in\K\setminus\{0\}$ and
\[\|f(x+ty)-f(x)\|_F\leq \Lip(f)\|ty\|_E=|t|\Lip(f)\|y\|_E \] entails
\[
\|df(x,y)\|_F=\lim_{t\to 0}\frac{1}{|t|}\|f(x+ty)-f(x)\|_F\leq \Lip(f)\|y\|_E .
\]
To establish (b),
we have to show that $\|df(x,y)\|_F\leq \Lip(f)$ for all $x\in U$
and $y\in\wb{B}^E_1(0)$. Since $df$ is continuous and $U^0$ is dense in~$U$,
it suffices to have the estimate for $x\in U^0$, and this was established in~(a).

(c)  The inequality ``$\geq$'' in (\ref{linkderlip}) holds by~(b).
For the converse inequality, excluding a trival case we may assume that
$L:=\sup\{\|f'(x)\|_{\op}\colon x\in U\}<\infty$. Then
$\Lip(f)\leq L$ since
\begin{eqnarray*}
\|f(y)-f(x)\|_F &=& \left\|\int_0^1 f'(x+t(y-x))(y-x)\, dt\right\|_F\\
&\leq &\int_0^1\|f'(x+t(y-x))\|_{\op}\|y-x\|_E\,dt\leq L\|y-x\|_E
\end{eqnarray*}
for all $x,y\in U$, using
the Mean Value Theorem (Proposition~\ref{lemMVT}).

(d) Given $x_0\in U$, we have $df(x_0,0)=0$. Since $df$ is continuous,
we find a convex $x_0$-neighborhood $U_0\sub U$ and $\ve>0$ such that
$df(U_0\times B^E_\ve(0))\sub B^F_1(0)$. Thus $\|f'(x)\|_{\op}\leq \frac{1}{\ve}$
for all $x\in U_0$, by Lemma~\ref{like-op-no-mult},
and hence $\Lip(f|_{U_0})\leq\frac{1}{\ve}$,
by~(c).
\end{prf}
Analogs for mappings between
subsets of locally convex spaces can also be of interest.
\begin{defn}\label{lip-lcx}
Let $E$ and $F$ be locally convex spaces and $f\colon X\to F$ be a function,
defined on a subset $X\sub E$.
\begin{description}[(D)]
\item[(a)] We say that~$f$ is (globally) 
\emph{Lipschitz} \index{Lipschitz map!on l.c. space}
if, for each continuous seminorm~$p$ on~$F$, there exist a continuous seminorm~$q$ on~$E$
and $L\in [0,\infty[$ such that
\begin{equation}\label{lcxlip}
p(f(y)-f(x))\leq L\, q(y-x)\;\;\mbox{for all $x,y\in X$.}
\end{equation}
\item[(b)] We say that $f$ is \emph{locally Lipschitz}
\index{locally Lipschitz map!on l.c. space}
if each $x_0\in X$ has a neighborhood $U$ in~$X$ such that~$f|_U$ is Lipschitz.
\item[(c)]
Let $p$ be a continuous seminorm on~$F$ and $q$ a continuous seminorm on~$E$.
If there exist $x,y\in X$ such that $q(y-x)=0$ and $p(f(y)-f(x))\not=0$,
we define $\Lip_{p,q}(f):=\infty$.
If $q(y-x)=0$ implies $p(f(y)-f(x))=0$ for all $x,y\in X$, we define
\[
\Lip_{p,q}(f):=\sup\left\{\frac{p(f(y)-f(x))}{q(y-x)}\colon \mbox{$x,y\in X$ with $q(y-x)\not=0$}\right\}
\in [0,\infty].
\]
\end{description}
\end{defn}
\begin{rem}\label{rem-lip-lcx}
(a) A map $f$ as in Definition~\ref{lip-lcx}
is Lipschitz if and only if, for each continuous seminorm~$p$ on~$F$,
there exists a continuous seminorm~$q$ on~$E$ such that $\Lip_{p,q}(f)<\infty$.
For such $p$ and~$q$,
\[
p(f(y)-f(x))\,\leq \, \Lip_{p,q}(f)\, q(y-x)\quad \mbox{for all $\, x,y\in X$,}
\]
and $\Lip_{p,q}(f)$ is the minimum choice for $L$ in (\ref{lcxlip}).\medskip

\noindent
(b) Every Lipschitz map is continuous, and so is every locally Lipschitz map.\medskip

\noindent
(c)
If $(E,\|\cdot\|)$ is a normed space in Definition~\ref{lip-lcx}, then Lipschitz continuity of~$f$ means
that for each continuous seminorm~$p$ on~$F$, there exists $L\geq 0$ such that
\[
p(f(y)-f(x))\leq L\|y-x\|\quad\mbox{for all $x,y\in X$}
\]
(cf.\ Exercise~\ref{exer-semin-norm}).\medskip

\noindent
(d) If $p$ is a continuous seminorm on a locally convex space~$F$, let $\alpha_p\colon F\to F_p$, $x\mto x+p^{-1}(\{0\})$
be the canonical map to the associated normed space $(F_p,\|\cdot\|_p)$, as in~\ref{modout}.
By definition, $\|\alpha_p(x)\|_p=p(x)$ for all $x\in F$.
As a consequence, a map
$f$ as in Definition~\ref{lip-lcx} is Lipschitz if and only if
$\alpha_p\circ f\colon X\to (F_p,\|\cdot\|_p)$ is Lipschitz for each
continuous seminorm~$p$ on~$F$.\medskip

\noindent
(e) Definition~\ref{lip-lcx} and the preceding remarks (a)--(d)
apply just as well if $F$ (and~$E$) may not be Hausdorff.
This will be useful for the refined information provided by
Remark~\ref{back-to-normed}(a).
\end{rem}
The next two lemmas
will be used in the Appendix for Chapter~\ref{chapcalcul}
to prove the existence of weak integrals for functions on
higher-dimensional sets.
\begin{lem}\label{cplipaut}
Let $K$ be a compact convex subset with non-empty interior
in a finite-dimensional normed space $(E,\|\cdot\|_E)$
and $f\colon K\to F$ be a $C^1$-map to a locally convex space~$F$. Then $f$ is Lipschitz.
\end{lem}
\begin{prf}
By Remark~\ref{rem-lip-lcx}(c), we may assume that $(F,\|\cdot\|_F)$ is a normed
space. Since $f'\colon K\to (\cL(E,F),\|\cdot\|_{\op})$ is continuous by Exercise~\ref{fin-diff-cts},
we have
\[
\|f'\|_\infty:=\sup\{\|f'(x)\|_{\op}\colon x\in K\}<\infty
\]
and deduce with Lemma~\ref{lipviaprime}(c)
that $\Lip(f)=\|f'\|_\infty<\infty$.
\end{prf}
Generalizing the case of functions on subsets of normed spaces
(as in Remark~\ref{rem-lip-lcx}(b)), given a metric space $(X,d)$,
we call a function $f\colon X\to F$ to a locally convex space~$F$ Lipschitz if,
for each continuous seminorm~$p$ on~$F$, there exists $L\geq 0$ such that
$p(f(y)-f(x))\leq L\, d(x,y)$ for all $x,y\in X$.

If $(X,d_X)$ and $(Y,d_Y)$ are metric spaces, we shall use the maximum metric
on $X\times Y$, given by $d((x_1,y_1),(x_2,y_2)):=\max\{d_X(x_1,x_2),d_Y(y_1,y_2)\}$
for $(x_1,y_1),(x_2,y_2)\in X\times Y$.
\begin{lem}\label{lippardepi}
Let $(X,d)$ be a metric space, $a<b$ be real numbers and
$f\colon X\times [a,b]\to F$ be a Lipschitz map to a locally convex space~$F$ such that the weak
integral
\[
g(x):=\int_a^bf(x,t)\, dt
\]
exists in $F$ for each $x\in X$. Then also $g\colon X\to F$, $x\mto g(x)$
is Lipschitz.
\end{lem}
\begin{prf}
For each continuous seminorm $p$ on~$F$, there is $L\in[0,\infty[$
such that
\[
p(f(x_2,t_2)-f(x_1,t_1))\leq L \max\{d(x_1,x_2),|t_2-t_1|\}
\]
for all $(x_1,t_1),(x_2,t_2)\in X\times [a,b]$. Hence
\[
p(g(x_2)-g(x_1))\leq\int_a^b p(f(x_2,t)-f(x_1,t))\,dt\leq L(b-a)d(x_1,x_2)
\]
for all $x_1,x_2\in X$ (by Lemma~\ref{weakinballb}), showing that $g$ is Lipschitz.
\end{prf}
The following lemma
will be used in the proof of Lemma~\ref{dl-beyond-ana}
to establish continuity for certain non-linear
mappings on subsets of locally convex direct limits.
Lemma~\ref{sim-old} (with $k=1$) will be used
in our development of differential calculus on normed spaces
(in the proof of Lemma~\ref{prepacompa}),
and also to establish
a refined result concerning
regularity properties of infinite-dimensional Lie groups
(Lemma~\ref{Evol-cts-dir}). 
A special case of Lemma~\ref{lem-der-to-ops} (singled out in Remark~\ref{back-to-normed}(c))
will be used in Exercise~\ref{exc-loc-quadr}
to discuss the asymptotics of the Newton method,
and in the proof of Theorem~\ref{thm-inv-para} (a version of the Inverse Function Theorem).
\begin{lem}\label{lipviadergen}
Let $E$ and $F$ be locally convex spaces,
$U\sub E$ be a convex subset with non-empty interior
and $f\colon U\to F$ be a~$C^1$-map.
Let $p$ be a continuous seminorm on~$F$
and $q$ be a continuous seminorm on~$E$
such that
\[
L:=\sup\{p(df(x,y))\colon x\in U, y\in \wb{B}^q_1(0)\}\,<\,\infty.
\]
Then $p(f(y)-f(x))\leq Lq(y-x)$ for all $x,y\in U$ and thus $\Lip_{p,q}(f)\leq L$.
\end{lem}
\begin{prf}
By Lemma~\ref{like-op-no-mult}, we have $p(df(x,y))\leq Lq(y)$ for all $x\in U$ and $y\in E$.
Using the Mean Value Theorem,
we deduce that
\begin{eqnarray*}
p(f(y)-f(x))&=&p\left(\int_0^1df(x+t(y-x),y-x)\,dt\right)\\
&\leq & \int_0^1p(df(x+ty,y-x))\,dt
\leq \int_0^1 Lq(y-x)\,dt = Lq(y-x)
\end{eqnarray*}
for all $x,y\in U$.
\end{prf}
\begin{lem}\label{sim-old}
Let $E$, $E_1,\ldots, E_k$ and~$F$ be locally convex spaces,
$U\sub E$ be a locally convex subset with dense interior
and
$g\colon U\times (E_1\times\cdots\times E_k)\to F$ be a continuous map such that
$g(x,\cdot)\colon E_1\times \cdots\times E_k\to F$ is $k$-linear for each $x\in X$.
Assume that $g$ admits a continuous partial differential
\[
d_1g\colon U\times (E_1\times\cdots\times E_k)\times E\to F
\]
with respect to the first variable $($as in Proposition~\emph{\ref{rulepartialno})}.
Then
\[
g^\vee\colon U\to \cL^k(E_1,\ldots, E_k;F)_b,\;\, x\mto g(x,\cdot)
\]
is continuous.
For each $x_0\in U$ and each continuous seminorm $p$ on~$F$,
there exist an $x_0$-neighborhood $U_0\sub U$,
continuous seminorms $q_j$ on~$E_j$ for $j\in\{1,\ldots, k\}$
and a continuous seminorm~$q$ on~$E$ such that
\begin{equation}\label{morethnno}
p(g(x_2,y)-g(x_1,y))\leq
q_1(y_1)\ldots q_k(y_k)q(x_2-x_1)
\end{equation}
for all $x_1,x_2\in U_0$ and $y=(y_1,\ldots, y_k)\in E_1\times\cdots\times E_k$.
\end{lem}
\begin{prf}
Note that $d_1g(x,\cdot)\colon E_1\times \cdots\times E_k\times E\to F$
is $(k+1)$-linear for each $x\in U$: By continuity, it suffices to check this if $x\in U^0$.
For $y\in E_1\times\cdots\times E_k$ and $z\in E$,
$d_1g(x,y,z)=d(g(\cdot,y))(x,z)$ is linear in~$z$ as $g(\cdot,y)$ is a $C^1$-map.
For $z\in E$, the map $d_1g(x,\cdot,z)$ is a pointwise limit of the $k$-linear maps
$(g(x+tz,\cdot)-g(x,\cdot))/t$, whence $d_1g(x,\cdot,z)$ is $k$-linear.

Let~$p$ be a continuous seminorm on~$F$ and $x_0\in U$.
Since $d_1g(x,0)=0$ and $d_1g$ is continuous, we find a convex $x_0$-neighborhood
$U_0\sub U$, continuous seminorms~$q_j$ on~$E_j$ for $j\in\{1,\ldots, k\}$
and a continuous seminorm~$q$ on~$E$ such that
\[
d_1g(U_0\times B^{q_1}_1(0)\times\cdots\times B^{q_k}_1(0)
\times B^q_1(0))\sub B^p_1(0).
\]
Then
\begin{equation}\label{prepreppp}
p(d_1g(x,y_1,\ldots,y_k,z))\leq q_1(y_1)\cdots q_k(y_k)q(z)
\end{equation}
for all $x\in U_0$ and $(y_1,\ldots, y_k,z)\in E_1\times \cdots\times E_k\times E$,
by Lemma~\ref{like-op-no-mult}.
Using the Mean Value Theorem and (\ref{prepreppp}), we deduce that
\begin{eqnarray*}
p(g(x_2,y)-g(x_1,y))&=&p\left(\int_0^1d_1g(x_1+t(x_2-x_1),y,x_2-x_1)\, dt\right)\\
&\leq & \int_0 ^1 p(d_1g(x_1+t(x_1-x_1),y,x_2-x_1))\,dt\\
&\leq & q_1(y_1)\ldots q_k(y_k)q(x_2-x_1)
\end{eqnarray*}
for all $x_1,x_2\in U_0$ and $y=(y_1,\ldots, y_k)\in E_1\times\cdots\times E_k$.
If $B_j\sub E_j$ is a bounded subset for $j\in\{1,\ldots, k\}$
and $B:=B_1\times\cdots\times B_k$,
then
\[
C:=\sup (q_1(B_1)\cdots q_k(B_k))<\infty.
\]
By the preceding, we have
\begin{equation}\label{even-lip}
\|g^\vee(x_2)-g^\vee(x_1)\|_{B,p}\leq C\, q(x_2-x_1)\;\, \mbox{for all $x_1,x_2\in U_0$.}
\end{equation}
Hence $g^\vee\colon U\to (\cL^k(E_1,\ldots, E_k;F),\|\cdot\|_{B,p})$
is locally Lipschitz and thus continuous. As the seminorms of the form $\|\cdot\|_{B,p}$ define the
locally convex vector topology on
$\cL^k(E_1,\ldots, E_k;F)_b$, we deduce that~$g^\vee$ is also continuous as a map to this space.
\end{prf}
\begin{rem}
In view of the $k$-linearity of~$g$ in the second argument,
the differentiability condition in Lemma~\ref{sim-old}
is equivalent to~$g$ being~$C^1$ (as follows from Example~\ref{exmultilinC1}
and Proposition~\ref{rulepartialno}; cf.\ also Exercise~\ref{exc-multCkell}).
Given $\ell\in\N$, one can show that $g^\vee\colon U\to \cL(E_1,\ldots, E_k;F)_b$ is $C^\ell$
whenever $g$ is~$C^{\ell+1}$ (see \cite{Gl07e}).
\end{rem}
\begin{lem}\label{lem-der-to-ops}
Let $E$ and $F$ be locally convex spaces, $U\sub E$ be a locally convex subset with dense interior,
$k\in \N$ and $f\colon U\to F$ be a $C^{k+1}$-map.
Then the following map is continuous:
\[
(d^{\,(k)}f)^\vee\colon U\to\cL^k(E,\ldots, E;F)_b,\;\, x\mto d^{\,(k)}f(x,\cdot).
\]
\end{lem}
\begin{prf} Apply Lemma~\ref{sim-old} to $g:= d^{\,(k)}f\colon U\times E^k\to F$.
\end{prf}
\begin{rem}\label{back-to-normed}
(a) The proof of Lemma~\ref{sim-old} shows that the mapping
$g^\vee\colon U\to (\cL^k(E_1,\ldots, E_k,;F),\|\cdot\|_{B,p})$
is locally Lipschitz in the situation of the lemma, for each
continuous seminorm~$p$ on~$F$ and each bounded subset
$B\sub E_1\times\cdots\times E_k$ (see (\ref{even-lip})).\medskip

\noindent
(b) As a consequence,
$(d^{\,(k)}f)^\vee\colon U\to(\cL^k(E,\ldots, E;F),\|\cdot\|_{B,p})$
is locally Lipschitz in the situation of Lemma~\ref{lem-der-to-ops},
for each continuous seminorm~$p$ on~$F$ and each bounded subset
$B\sub E^k$.\medskip

\noindent
(c) In particular, if $(E,\|\cdot\|_E)$ and $(F,\|\cdot\|_F)$
are normed spaces and $f\colon U\to F$ is a $C^2$-map on a locally convex subset $U\sub E$ with dense interior, then $f'\colon U\to \cL(E,F)$ is locally Lipschitz
as a map to $(\cL(E,F),\|\cdot\|_{\op})$ (and hence continuous).
\end{rem}
\begin{small}
\subsection*{Exercises for Section~\ref{sec-lip}}

\begin{exer}\label{exc-lip-for-op}
Let $(E,\|\cdot\|_E)$ and $(F,\|\cdot\|_F)$ be normed spaces and
$\alpha\colon E\to F$ be a continuous linear map.
Show that~$\alpha$ is Lipschitz and $\Lip(\alpha)=\|\alpha\|_{\op}$.
\end{exer}

\begin{exer}\label{exc-lip-com-pro}
Let $(X,d_X)$, $(Y,d_Y)$ and $(Z,d_Z)$ be
metric spaces.
\begin{description}[(D)]
\item[(a)]
If $f\colon X\to Y$ and $g\colon Y\to Z$ are Lipschitz maps,
show that $g\circ f\colon X\to Z$ is Lipschitz with
$\Lip(g\circ f)\leq\Lip(g)\Lip(f)$.
\item[(b)]
Endow $Y\times Z$ with the metric $d((y_1, z_1), (y_2,z_2)) = \max(d_Y(y_1, y_2), 
d_Z(z_1, z_2))$. 
Show that a function $f=(f_1,f_2)\colon X\to Y\times Z$ is Lipschitz if and only if both of
its components~$f_1$ and~$f_2$ are Lipschitz.
Moreover, $\Lip(f)=\max\{\Lip(f_1),\Lip(f_2)\}$.
\end{description}
\end{exer}

\begin{exer}
Let $(K,d_K)$ be a compact metric space and $(Y,d_Y)$ be a metric space.
Show that a function $f\colon K\to Y$ is Lipschitz if and only if it is locally Lipschitz.
\end{exer}

\begin{exer}\label{exc-lcx-loc-lip}
Let $E$ and~$F$ be locally convex süpaces, $U\sub E$ be a locally convex subset with dense interior and $f\colon U\to F$
be a function.
Show: If $f$ is~$C^1$, then
$f\colon U\to (F,p)$ is locally Lipschitz for each continuous seminorm~$p$ on~$F$.\\[1mm]
[Use Lemma~\ref{lipviadergen}
and a variant of the proof of Lemma~\ref{lipviaprime}(d).]
\end{exer}

\end{small}
\section{Taylor expansions and homogeneous polynomials}\label{secTay}
We now complete
our exposition of $C^k$-maps
with a discussion of Taylor expansions.
As in the finite-dimensional
case, Taylor's Theorem will allow us to
write a $C^k$-function $f$ close to $x$
in the form $f(x+y)=p(y)+R(y)$, where $p$ is
a polynomial of degree $\leq k$
and $\lim_{t\to 0}\frac{1}{t^k}R(ty)=0$.
\subsection*{Taylor's Formula for curves}
\begin{prop}\label{taycurve}
Let $F$ be a real locally convex space, $k\in \N$,
$\gamma\colon I\to F$ be a $C^k$-curve
on a non-degenerate interval $I\sub \R$,
and $x\in I$.
\begin{description}[(D)]
\item[\rm (a)]
Then
\begin{equation}\label{tay1}
\gamma(x+t)\, =\, \gamma(x)+t\gamma'(x)+\cdots
+\frac{t^{k-1}}{(k-1)!}\gamma^{(k-1)}(x)+\rho_k(t)
\end{equation}
for all $t\in I-x$, with remainder term $\rho_k(t)$ given by
\begin{equation}\label{tay2}
\rho_k(t)=\frac{t^k}{(k-1)!}\int_0^1(1-r)^{k-1}\gamma^{(k)}(x+rt)\, dr\,.
\end{equation}
\item[\rm(b)]
Furthermore,
\begin{equation}\label{taynicer}
\gamma(x+t)=\gamma(x)+t\gamma'(x)+\cdots
+\frac{t^k}{k!}\gamma^{(k)}(x)+t^k R_k(t)
\end{equation}
for $t\in I-x$, for a uniquely determined
continuous map
$R_k\colon I-x\to F$. Here $R_k(0)=0$,
and $R_k(t)$ can be written as the weak integral
\begin{equation}\label{betterrem}
R_k(t)=
\frac{1}{(k-1)!} \int_0^1(1-r)^{k-1} (\gamma^{(k)}(x+rt)
-\gamma^{(k)}(x))\, dr\,.
\end{equation}
\end{description}
\end{prop}
\begin{prf}
(a) If $k=1$, then the Mean Value Theorem shows that
\[
\rho_1(t)=\gamma(x+t)-\gamma(x)=\int_0^1d\gamma(x+rt,t)\,dr=t\int_0^1\gamma'(x+rt)\,dr,
\]
which
is of the asserted form. Now assume that the result holds for $k$ and $\gamma$ is $C^{k+1}$.
Noting that $\int (1-r)^{k-1}\, dr=-\frac{1}{k}(1-r)^k$
and $\frac{d}{dr}(\gamma^{(k)}(x+rt))=t\gamma^{(k+1)}(x+rt)$,
Partial Integration
(as in Exercise~\ref{excprodpartint}) yields
\begin{eqnarray*}
\rho_k(t)&=& \frac{t^k}{(k-1)!}\int_0^1(1-r)^{k-1}\gamma^{(k)}(x+rt)\, dr\\
&=&\frac{t^k}{k!}\gamma^{(k)}(x)+
\frac{t^{k+1}}{k!}\int_0^1(1-r)^k\gamma^{(k+1)}(x+rt)\, dr,
\end{eqnarray*}
whence $\gamma(x+t)=\sum_{j=0}^k\frac{t^j}{j!}\gamma^{(j)}(x)+\rho_{k+1}(t)$
with $\rho_{k+1}(t)$ of the required form.

(b) Since $k\int_0^1 (1-r)^{k-1} dr= [-(1-r)^k]_0^1=1$,
we can rewrite the remainder term $\rho_k$ in~(a)
as
\begin{eqnarray*}
\rho_k(t) &= &\frac{t^k}{k!}\gamma^{(k)}(x)
+\rho_k(t)-\frac{t^k}{k!}k \int_0^1 (1-r)^{k-1} \gamma^{(k)}(x)
\, dr\\
& = & \frac{t^k}{k!}\gamma^{(k)}(x)
+t^k \underbrace{\frac{1}{(k-1)!}\int_0^1(1-r)^{k-1} (\gamma^{(k)}(x+rt)
-\gamma^{(k)}(x))\, dr}_{=:R_k(t)}\,.
\end{eqnarray*}
Note that $R_k$ is continuous as a parameter-dependent
integral with continuous integrand, and $R_k(0)=0$.
For $t\not=0$, apparently $R_k(t)$ is uniquely
determined by $\gamma$, because we can then solve
(\ref{taynicer})
for $R_k(t)$.
Then also $R_k(0)=\lim_{0\not=t\to 0}R_k(t)$
is uniquely determined, by continuity.
\end{prf}
\subsection*{Taylor's Theorem for $C^k$-maps}
In multi-variable calculus,
Taylor's Theorem
explains how to approximate a function
close to a given point as good as possible
by a polynomial (polynomial function) of several variables.
Our goal is to do the same
for functions between locally convex
spaces, using polynomials in the following sense.
\begin{defn}
Let $E$ and $F$ be $\K$-vector spaces.
A \emph{homogeneous polynomial} \index{homogeneous polynomial}
of degree~$k\in\N_0$ from $E$ to $F$
is a map
$p  \colon E\to F$
of the form
\begin{equation}\label{betaviamu}
p (x)\, =\, \beta\big(\underbrace{x,x,\ldots,x}_{k}\big)\vspace{-2.5mm}
\end{equation}
for some $k$-linear map
$\beta \colon E^k\to F$.
\end{defn}
Note that homogeneous polynomials
of degree~$0$ are constant maps.\footnote{A $0$-linear
map $E^0\to F$ is a (necessarily
constant) map
$\beta\colon E^0=\{0\}\to F$.
Reading ``$\underbrace{x,\ldots, x}_{0}$''
as $0$,
we see that each
homogeneous\vspace{-3mm}
polynomial of degree $0$ is constant.}
By definition,
a homogeneous polynomial
$p\colon E\to F$ of degree~$1$ is
a linear map.
\begin{rem}\label{ishomog}
If $p\colon E\to F$
is a homogeneous polynomial of degree~$k$,
then $p(tx)=\beta (tx,tx,\ldots,tx)=t^k\beta (x,\ldots,x)=t^kp(x)$
for all $x\in E$ and $t\in \K$.
\end{rem}
\begin{defn}\label{defgenpol}
Let $k\in\N_0$.
A function $p\colon E\to F$ between $\K$-vector spaces is called a 
\index{polynomial of degree $\leq k$}
\emph{polynomial} of degree $\leq k$
if it can be written as a sum $p=\sum_{j=0}^kp_j$ of homogeneous polynomials
$p_j\colon E\to F$
of degree~$j$, for $j\in\{0,1,\ldots, k\}$.
\end{defn}
\begin{rem}
The homogeneous components $p_j$ of a polynomial $p=\sum_{j=0}^kp_j$ (as in Definition~\ref{defgenpol})
are uniquely determined. To see this,
it suffices to show that $\lambda\circ p_j$ is determined by $\lambda\circ p$
for each linear functional\linebreak
$\lambda\colon F\to\K$
(as the latter separate points
on~$F$). However,
for each $y\in E$, the function
$\K\to \K$, $t\mto \lambda(p(ty))=\sum_{j=0}^k\lambda(p_j(y))t^j$
is an ordinary polynomial and hence determines its coefficients $\lambda(p_j(y))$.
\end{rem}
\begin{rem}
The $k$-linear map $\beta\colon E^k\to F$
in the definition of a homogeneous
polynomial can always be chosen
as a \emph{symmetric} $k$-linear map.
In fact, given $\beta$,
the map
\[
\wb{\beta} \colon E^k\to F\,,\quad
\wb{\beta}(x_1,\ldots, x_k)\; :=\;
\frac{1}{k!}\sum_{\sigma\in S_k} \beta(x_{\sigma(1)},\ldots, x_{\sigma(k)})\vspace{-2mm}
\]
is symmetric $k$-linear, and
$p(x)=\beta(x,\ldots, x )=\wb{\beta}(x,\ldots, x)$
for each $x\in E$. If $E$ and $F$ are topological
vector spaces and $\beta$ is continuous, then also~$\wb{\beta}$ is continuous.
We shall see later that $\wb{\beta}$ is uniquely
determined by~$p$.
\end{rem}
The following lemma provides examples of homogeneous polynomials
and shows that compositions of polynomials are polynomials.
\begin{lem}\label{compopol}
Let $E$, $E_1,\ldots, E_n$, $F$, and $H$ be vector spaces.
Then the\linebreak
following holds:
\begin{description}[(D)]
\item[\rm(a)]
If $\beta\colon E_1\times \cdots\times E_n\to F$
is an $n$-linear map, then $\beta$ is a homogeneous polynomial of degree~$n$.
\item[\rm (b)]
If $\beta\colon E_1\times \cdots \times E_n\to H$ is an $n$-linear map
and $(p_j)_{j\in\{1,\ldots, n\}}$ is an $n$-tuple of homogeneous polynomials
$p_j\colon E\to E_j$
of degree~$k_j$,
then
\[
p\colon E\to F, \quad p(x):=\beta(p_1(x),\ldots, p_n(x))
\]
is a homogeneous polynomial
of degree $k_1+\cdots+ k_n$.
\item[\rm (c)]
If
$q \colon F\to H$ is
a homogeneous polynomial of degree~$k$
and
$p\colon E\to F$ a homogeneous
polynomial
of degree $\ell$,
then $q\circ p\colon E\to H$
is a homogeneous polynomial
of degree~$k\ell$.
\item[\rm (d)]
If
$q\colon F\to H$ is a
polynomial
of degree $\leq k$ and $p\colon E\to F$
a polynomial of degree~$\leq \ell$,
then $q\circ p \colon E\to H$
is a polynomial
of degree~$\leq k\ell$.
\end{description}
\end{lem}
\begin{prf}
(a) Define
$\alpha\colon (E_1\times\cdots\times E_n)^n\to F$,
\[
\alpha(x_1,\ldots, x_n):= \beta(x_{1,1},\ldots,x_{n,n})\,,
\]
where
$x_j=(x_{j,1},\ldots, x_{j,n})$.
Then $\alpha$ is $n$-linear and
$\beta(x)=\alpha(x,\ldots, x)$
for all $x=(x_1,\ldots, x_n)\in E_1\times \cdots\times
E_n$.

(b) For $j\in \{1,\ldots, n\}$,
let $\beta_j\colon E^{k_j}\to F$ be a $k_j$-linear
map such that $p_j(x)=\beta_j(x,\ldots, x)$.
Then
$\alpha:=
\beta\circ (\beta_1\times \cdots\times \beta_n)
\colon E^{k_1+\cdots+ k_n}\to H$\linebreak
is $(k_1+\cdots+k_n)$-linear
and $p(x)=\alpha(x,\ldots, x)$,
whence indeed $p$ is a homogeneous polynomial
of degree $k_1+\cdots+k_n$.

(c) We have $q(x)=\beta(x,\ldots, x)$
for a $k$-linear map $\beta\colon F^k\to H$.
Since $q(p(x))=\beta(p(x),\ldots, p(x))$,
Part\,(b) applies with
$p_1,\ldots, p_k:=p$.

(d) We have $q=\sum_{j=0}^k q_j$
and $p=\sum_{i=0}^\ell p_i$,
where $p_i\colon E\to F$
is a homogeneous polynomial
of degree~$i$ and
$q_j(x)=\beta_j(x,\ldots, x)$
for a $j$-linear map $\beta_j\colon F^j\to H$.
Then
\[
q\circ p\;=\; \sum_{j=0}^k \; \sum_{i_1, \ldots, i_j=1}^\ell
\beta_j\circ (p_{i_1},\ldots, p_{i_j})\,,
\]
where each summand is a homogeneous
polynomial of degree
$i_1+\cdots+i_j\leq j\ell\leq k\ell$,
by~(b). The assertion follows.
\end{prf}
G\^{a}teaux differentials are an important source of homogeneous polynomials,
and account for the importance of the latter in differential calculus.
\begin{defn}
Let $E$ and $F$ be locally convex spaces,
$f\colon U\to F$ be a map on an open
subset $U\sub E$, $x\in U$, $y\in E$  and
$k\in \N$.
If there exists an open $0$-neighborhood $W\sub \K$ such that
$x+Wy\sub U$
and
$W\to F$,
$t\mto f(x+ty)$
is $k$ times differentiable,\footnote{Thus, we require that $k$ real derivatives
(as in \ref{dfdercu}) can be formed if $\K=\R$, and $k$ complex derivatives (as in~\ref{defcxdersing}) if $\K=\C$.}
then we define
\[
\delta_x^kf(y)\;:=\; \frac{d^k}{dt^k}\Big|_{t=0}f(x+ty)\,.
\]
If $\delta_x^kf(y)$ exists for each $y\in E$
and
$\delta_x^kf\colon E\to F$,
$y\mto \delta_x^kf(y)$
is a homogeneous polynomial
of degree~$k$,
then $\delta_x^kf\colon E\to F$
is called the
\emph{$k$th G\^{a}teaux differential of $f$ at~$x$.}
\index{G\^{a}teaux differential} 
We let $\delta_x^0f$ be the constant map $E\to F$,
$y\mto f(x)$.
\end{defn}
\begin{lem}\label{examppoly}
Let $E$ and $F$ be locally convex spaces,
$k\in \N\cup\{\infty\}$
and $f\colon U\to F$ a $C^k$-map
on an open subset $U\sub E$.
Then $\delta_x^jf$ exists for each $j\in \N_0$ such that $j\leq k$
and each $x\in E$, and is given by
\[
\delta_x^jf(y)\;=\; d^{\,(j)}f(x,y,\ldots, y)\,.
\]
Thus
$\delta_x^jf\colon E\to F$ is a continuous
homogeneous polynomial of degree~$j$.
\end{lem}
\begin{prf}
There is $\ve>0$ such that $x+\bD_\ve^0  \, y\sub U$.
Then $\gamma\colon \bD_\ve^0 \to F$, $\gamma(t):=
f(x+ty)$ is~$C^k$.
We have
\[
\gamma'(t)\,=\,\lim_{s\to 0}\,\frac{1}{s}\left(f(x+ty+sy)-f(x+ty)\right)
\,=\, df(x+ty,y)
\]
and inductively $\gamma^{(j)}(t)=d^{\,(j)}f(x+ty,y,\ldots,y)$
for all $j=1,\ldots, k$.
Thus $\delta^j_xf(y)\stackrel{\text{def}}{=}
\frac{d^j}{dt^j}\big|_{t=0}f(x+ty)
=\gamma^{(j)}(0)
=d^{\,(j)}f(x,y,\ldots,y)$,
as required.
\end{prf}
The G\^{a}teaux differentials
of homogeneous
polynomials (and hence of arbitrary
polynomials) can be calculated easily.
\begin{lem}\label{Gatpoly}
If $E$ and $F$ are locally convex spaces,
$\!n\in \N_0$ and
$\beta\colon \!E^n \to F$ is a symmetric
$n$-linear map, then
$p\colon E\to F$, $p(x):=\beta(x,\ldots, x)$
admits a $k$th G\^{a}teaux
differential $\delta_x^kp$,
for all $k\in \N_0$ and $x\in E$.
It is given by
\begin{equation}\label{gatp1}
\delta_x^kp(y)\; =\;
\left\{
\begin{array}{cl}
\frac{n!}{(n-k)!}\,
\beta(\underbrace{x,\ldots, x}_{n-k},\underbrace{y,\ldots, y}_{k}) &
\; \mbox{if $\,k\leq n$;}\\[2mm]
0 & \; \mbox{if $\, k>n$.}
\end{array}
\right.
\end{equation}
In particular,
\begin{equation}\label{gatp2}
\delta_0^kp(y)\; = \; \left\{
\begin{array}{cl}
n!\, p(y) &\;\mbox{if $\, k=n$;}\\
0 & \; \mbox{if $\, k\not=n$.}
\end{array}
\right.
\end{equation}
If $\beta$ is continuous, then $p$ is smooth.
\end{lem}
\begin{prf}
Given $x,y\in E$,
define $\gamma\colon \K\to F$,
$\gamma(t):=p(x+ty)$.
Then
\[
\gamma(t)\,=\,\beta(x+ty,\ldots, x+ty)\,=\,
\sum_{j=1}^n t^j \Big({n\atop j}\Big)
\beta(\underbrace{x,\ldots,x}_{n-j},\,
\underbrace{y,\ldots,y}_{j})\,,\vspace{-2mm}
\]
as $\beta$ is $n$-linear and symmetric.
It is clear from this formula
that $\gamma$ is smooth.
Differentiating this formula
as usual, we obtain (\ref{gatp1})
for $\delta_x^kp(y)=\gamma^{(k)}(0)$
(which apparently is a homogeneous
polynomial of degree~$k$ in~$y$).
Taking $x:=0$, (\ref{gatp1}) implies
(\ref{gatp2}).
If $\beta$ is continuous,
then~$p$ is smooth,
as it is
the composition of the smooth map
$\beta$ (see Example~\ref{exmultilin})
and the continuous linear and hence smooth
map $E\to E^n$, $x\mto (x,\ldots, x)$.
\end{prf}
\begin{rem}
If $p=\sum_{k=0}^\ell p_k\colon E\to F$
is a polynomial of degree $\leq \ell$
between locally convex spaces,
with homogeneous polynomials
$p_k\colon E\to F$ of degree~$k$,
then (\ref{gatp2}) implies that
\begin{equation}\label{enabtayl}
\frac{\delta^k_0 p}{k!}\;=\;
\left\{
\begin{array}{cl}
p_k & \;\mbox{if $\, k\leq \ell$;}\\
0 & \;\mbox{if $\, k>\ell$.}
\end{array}
\right.
\end{equation}
\end{rem}
\begin{defn}
Let $E$ and $F$ be locally convex spaces, $U\sub E$ be an open subset, $k\in\N_0$
and $f\colon U\to F$ a $C^k$-map. For $x\in U$, the map
\[
P^k_xf\colon E\to F,\quad P^k_xf(y):=\sum_{j=0}^k\frac{1}{j!}\delta^j_xf(y)
\]
is called the \emph{$k$th order Taylor polynomial} 
\index{Taylor polynomial!of order $k$} 
\index{$k$-jet} 
of $f$ at~$x$, or also the \emph{$k$-jet} of~$f$ at~$x$.
\end{defn}
\begin{rem}
By Lemma~\ref{examppoly},
$P_x^kf=\sum_{j=0}^k\frac{1}{j!}\delta^j_xf$\vspace{-.5mm} is a polynomial of degree $\leq k$ with continuous
homogeneous components. Hence $P^k_xf\colon E\to F$ is a smooth map,
by Lemma~\ref{Gatpoly}.
\end{rem}
Taylor's Theorem provides information on the remainder term which\linebreak
occurs
if we approximate $f(x+y)$ by $P^k_xf(y)$.
\begin{thm}[Taylor's Theorem]\label{thmtay}
Let $E$ and
$F$ be locally convex spaces,
$k\in \N_0 $ and
$f\colon U\to F$ be a $C^k$-map
on an open subset $U\sub E$.
Given $x\in U$, write
\begin{equation}\label{detRkvar}
f(x+y) = P^k_xf(y)+
R(y)
\end{equation}
for $y\in U-x$.
Then
\begin{equation}\label{explicitvar}
R(y)\, =\, \frac{1}{(k-1)!}\int_0^1(1-s)^{k-1} (\delta_{x+sy}^kf(y)
-\delta_x^kf(y))\, ds
\end{equation}
for all $y\in U-x$ such that $x+[0,1]y\sub U$.
\end{thm}
\begin{prf}
For $y\in U-x$ such that $x+[0,1]y\sub U$, the map
$\gamma\colon [0,1]\to F$, $\gamma(t):=f(x+ty)$
is a $C^k$-curve,
by Proposition~\ref{chainCk}.
As in the proof of Lemma~\ref{examppoly},
we see that
\begin{equation}\label{useme}
\gamma^{(j)}(t)  = d^{\,(j)}f(x+ty,y,\ldots,y)
=\delta_{x+ty}^jf(y)
\end{equation}
for $j=1,\ldots, k$.
Using
Proposition~\ref{taycurve}(b), we deduce that
\begin{eqnarray*}
f(x+y)\; = \; \gamma(1) & = &\gamma(0)+\gamma'(0)+\cdots +
{\ts \frac{1}{k!}}\gamma^{(k)}(0)\\
& & \quad +{\ts \frac{1}{(k-1)!}}\int_0^1
(1-s)^{k-1}(\gamma^{(k)}(s)-\gamma^{(k)}(0))\, ds\\
& = & f(x)+\delta^1_xf(y)+\cdots + {\ts \frac{1}{k!}} \delta^k_xf(y)\\
& & \quad
+ {\ts \frac{1}{(k-1)!}}
\int_0^1 (1-s)^{k-1}(\delta^k_{x+sy}f(y)-\delta^k_xf(y))\,
ds,
\end{eqnarray*}
from which (\ref{explicitvar}) follows.
\end{prf}
The $k$th order Taylor polynomial can be characterized in several ways.
\begin{prop}\label{charTaypo}
Let $E$ and
$F$ be locally convex spaces,
$k\in \N_0 $ and
$f\colon U\to F$ be a $C^k$-map
on an open subset $U\sub E$.
Then the following holds for each $x\in U$:
\begin{description}[(D)]
\item[\rm (a)]
$p:=P_x^kf\colon E\to F$ is the unique polynomial of degree $\leq k$
such that
$\delta_0^j(p)=\delta_x^jf$
for each $j\leq k$
$($and thus $P_0^k(p)=P_x^kf)$.
\item[\rm (b)]
$p:=P^k_x f$
satisfies
\begin{equation}\label{char2tp}
\lim_{t\to 0} \frac{f(x+ty)-p(ty)}{t^k}\; =\; 0\quad
\mbox{for each $y\in E$}
\end{equation}
and is the unique polynomial
of degree $\leq k$ with this property.
\end{description}
\end{prop}
\begin{prf}
(a) is a trivial consequence of (\ref{enabtayl}).

(b) The map $P^k_x f:=\sum_{j=0}^k \frac{1}{j!}\delta^j_xf$
is a continuous polynomial of degree $\leq k$.
Given $x\in U$ and $y\in E$,
there is $\ve>0$ such that $x+[0,\ve]y\sub U$.
Applying (\ref{explicitvar})
to $f(x+ty)-P^k_xf(ty)=R(ty)$ with $t\in [0,\ve]$, we see that
\[
R(ty)
= {\ts \frac{1}{(k-1)!}}
\int_0^1 (1-s)^{k-1}(\delta^k_{x+sty}f(ty)-\delta^k_xf(ty))\,
ds= t^k \int_0^1 h(t,s)\, ds
\]
with $h\colon [0,\ve]\times [0,1]\to F$,
\[
h(t,s):=
{\ts \frac{1}{(k-1)!}} (1-s)^{k-1}(\delta^k_{x+sy}f(y)-\delta^k_xf(y)).
\]
As $h$ is continuous, $t^{-k}(f(x+ty)-P^k_xf(ty))=\int_0^1h(t,s)\,ds\to \int_0^1h(0,s)=0$
for $t\to 0$ by the continuity of parameter-dependent integrals
(Lemma~\ref{intpar}).

Now let $p\colon E\to F$
be any polynomial of degree $\leq k$ such that (\ref{char2tp}) holds.
We have $p=\sum_{j=0}^k p_j$ with homogeneous polynomials $p_j\colon E\to F$
of degree~$j$.
Then, for each~$y\in E$ and
each continuous linear functional~$\lambda\in F'$,
\[
q\colon \K\to \K,\quad q(t)\, :=\, \lambda(P_x^k f(ty)-p(ty))\, =\,
\sum_{j=0}^kt^j\lambda(\delta^j_xf(y)- p_j(y))
\]
is an ordinary polynomial satisfying
$\lim_{t\to 0} \,\frac{1}{t^k}q(t)=0$.
As~$q$ has degree $\leq k$, this is only possible
if~$q$ vanishes identically.
Since~$F'$ separates points on $F$,
we deduce that $P_x^kf(ty)=p(ty)$ for all $t\in \K$,
whence $P_x^k f(y)=p(y)$ in particular.
As $y$ was arbitrary, we obtain $P_x^k f=p$.
\end{prf}
\begin{ex}\label{extayofpol}
If $E$ and $F$ are locally convex spaces and $p=\sum_{j=0}^\ell p_j$
is a polynomial with continuous homogeneous components $p_j\colon E\to F$,
then (\ref{enabtayl}) entails that
\[
P^k_0(p)=\left\{
\begin{array}{cl}
\sum_{j=0}^kp_j &\mbox{if $\, k\leq \ell$;}\\
p &\mbox{if $\,k\geq \ell$.}
\end{array}
\right.
\]
\end{ex}
\begin{ex}\label{truncatedco}
Let $E$, $F$, and $H$ be locally convex spaces,
$p\colon E\to F$ and $q\colon F\to H$ be polynomials
of degree $\leq k$,
with continuous homogeneous components.
Then
\[
q\circ p=\sum_{j=0}^{k^2} r_j
\]
with continuous homogeneous polynomials $r_j\colon E\to F$ of degree~$j$,
by Lemma~\ref{compopol}(d) and its proof.
Example~\ref{extayofpol} now shows that
\[
P^k_0(q\circ p)=\sum_{j=0}^k r_j
\]
is the ``truncated composition" of $p$ and $q$, obtained by chopping off
all\linebreak
homogeneous components of degree $>k$.
\end{ex}
\begin{lem}\label{polyshift}
Let $f\colon U\to F$
be a $C^k$-map on an open $0$-neighborhood
$U\sub E$,
and $a\in E$.
Then $P_a^k(f(\cdot -a))=P^k_0 f$.
\end{lem}
\begin{prf}
Applying Proposition~\ref{charTaypo}(b)
to the map $f$ around~$0$,
we find that $t^{-k}(f(ty)-P^k_0f(ty))\to 0$
as $t\to 0$.
Hence $t^{-k}(f((a+ty)-a)-P^k_0f(ty))\to 0$,
whence $P_a^k(f(\cdot -a))=P_0^kf$
by uniqueness
in Proposition~\ref{charTaypo}(b).
\end{prf}
\subsection*{The Polarization Formula}
Generalizing Example~\ref{truncatedco}, we shall discuss Taylor polynomials of
compositions of $C^k$-maps. Before, let us analyze how a symmetric
multilinear map can be recovered from the associated homogeneous polynomial.
For a symmetric bilinear map $\beta\colon E\times E\to F$,
this is easy: If $p(x):=\beta(x,x)$, then
\[
p(x\pm y)=\beta(x,x)\pm 2\beta(x,y)+\beta(y,y)
\]
for all $x,y\in E$, whence
\[
\beta(x,y)=\frac{1}{4}(p(x+y)-p(x-y))=\frac{1}{8}(p(x+y)-p(x-y)-p(-x+y)+p(-x-y)).
\]
We now state an appropriate generalization for symmetric $k$-linear maps.
\begin{prop}[\textbf{Polarization Formula}]\label{proppolarvar}
Let $E$ and $F$ be vector spaces,
$\beta\colon E^k\to F$ be
symmetric $k$-linear, and $p \colon E\to F$,
$p(x):=\beta(x,\ldots, x)$ be the associated
homogeneous polynomial.
Then
\begin{equation}\label{polarformulvar}
\beta(x_1,\ldots,x_k)\,=\,
\frac{1}{k!\, 2^k}\sum_{\ve_1,\ldots,\ve_k\in \{1,-1\}}
\ve_1\cdots\ve_k\,
p(\ve_1x_1+\cdots+\ve_kx_k)
\end{equation}
for all $x_1,\ldots, x_k\in E$.
Thus $\beta$ is uniquely
determined by~$p$, and if $E$ and $F$ are topological
vector spaces, then $\beta$ is continuous if and only
if so is~$p$.
\end{prop}
\begin{prf}
Let $M$ be the set of all functions $\sigma\colon \{1,\ldots, k\}\to\{1,\ldots, k\}$
and $S_k\sub M$ be the symmetric group of all permutations of $\{1,\ldots, k\}$.
If $x_1,\ldots, x_k\in E$, let $\theta(x_1,\ldots, x_k)$
be the right-hand side of (\ref{polarformulvar}). Since
\begin{eqnarray*}
p(\ve_1x_1+\cdots+\ve_kx_k)
&=& \beta(\ve_1x_1+\cdots+\ve_kx_k,\ldots, \ve_1x_1+\cdots+\ve_kx_k)\\
&=& \sum_{\sigma\in M}
\ve_{\sigma(1)}\cdots\ve_{\sigma(k)}\, \beta(x_{\sigma(1)},\ldots, x_{\sigma(k)}),
\end{eqnarray*}
we have
\[
\theta(x_1,\ldots, x_k)
= \frac{1}{k!\,2^k}\sum_{\sigma\in M} A_\sigma\,
\beta(x_{\sigma(1)},\ldots, x_{\sigma(k)})
\]
with
\[
A_\sigma :=
\sum_{\ve_1,\ldots,\ve_k=\pm 1}
\ve_{\sigma(1)}\cdots\ve_{\sigma(k)}\, \ve_1\cdots\ve_k\,.
\]
If $\sigma\in S_k$, then
$A_\sigma
=\sum_{\ve_1,\ldots,\ve_k=\pm 1}1=2^k$.
If $\sigma\in M\setminus S_k$,
then there exists $j\in \{1,\ldots, k\}$
outside the image of~$\sigma$.
Hence
\[
A_\sigma = (1-1)
\sum_{\ve_i=\pm 1, i\not=j}
\ve_{\sigma(1)}\cdots\ve_{\sigma(k)}\prod_{i\not=j}\ve_i=0\,.
\]
Using the symmetry of $\beta$, we deduce that
\[
\theta(x_1,\ldots, x_k)
= \frac{1}{k!\,2^k}\sum_{\sigma\in S_k} 2^k
\beta(x_{\sigma(1)},\ldots, x_{\sigma(k)})
=\beta(x_1,\ldots,x_k).\qedhere
\]
\end{prf}
\begin{cor}\label{applpolze}
Let $E$ and $F$ be locally convex spaces
and $f\colon U\to F$
be a $C^k$-map on an open subset
$U\sub E$, where $k\in \N$.
Let $F_0\sub F$ be a vector subspace,
and $x\in U$. Then
$\delta_x^k f(E)\sub F_0$
if and only if $d^{\,(k)}f(\{x\}\times E^k)\sub F_0$.
\end{cor}
\begin{prf}
Since $\delta_x^kf(y)=d^{\,(k)}f(x,y, \ldots, y)$,
it is obvious that $F_0$ will contain
the image of $\delta_x^kf$ if it contains
the image of $d^{\,(k)}f(x,\sbull)$.
Since the symmetric $k$-linear map
$d^{\,(k)}f(x,\sbull)$ can be recovered
from $\delta_x^kf$ by the Polarization
Formula (\ref{polarformulvar}),
it is clear that also
the converse holds.
\end{prf}
\subsection*{Chain Rule for Taylor Polynomials}
\begin{prop}\label{chainrtay}
Let $E,F$ and $H$ be locally convex spaces,
$U\sub E$ and $V\sub F$ be open,
$k\in \N_0$ and
$f\colon U\to V\sub F$,
$g\colon V\to H$ be $C^k$-maps.
Given $x\in U$, set $z:=f(x)$. Then
\begin{equation}\label{chainjet}
P_x^k(g\circ f)\;=\; P_0^k\big( (P_z^kg)\circ (P_x^k f \, -z )\big),
\end{equation}
where the right hand side can be calculated as in Example~{\em\ref{truncatedco}}. Thus
\begin{eqnarray}
P^k_x(g\circ f)(y) &=&\qquad\qquad\hspace*{7.5cm}\label{prenmfaa}\\
\hspace*{-3mm}\lefteqn{g(f(x))+\sum_{j=1}^k\frac{1}{j!}\sum_{i_1+\cdots+ i_j\leq k}\frac{1}{i_1!\cdots i_j!}
d^{\,(j)}g\big(f(x),\delta^{i_1}_xf(y),\ldots, \delta^{i_j}_xf(y)\big).}\qquad\qquad \notag
\end{eqnarray}
\end{prop}
\begin{prf}
If $f_1\colon U_1\to V_1$
and $g_1\colon V_1\to H$ are $C^k$-maps
on open neighborhoods
$U_1\sub E$ of~$x$ and $V_1\sub F$ of~$z$
such that
$P_x^k(f_1)=P_x^k(f)$
and $P_z^k(g_1)=P_z^k(g)$,
then $\delta^j_xf=\delta^j_xf_1$ and thus $d^{\,(j)}f(x,\cdot)=d^{\,(j)}f_1(x,\cdot)$
(by polarization), for each $j\in\{0,1,\ldots, k\}$.
Likewise, $d^{\,(j)}g(z,\cdot)=d^{\,(j)}g_1(z,\cdot)$ for each~$j$.
Hence $d^{\,(j)}(g\circ f)(x,\cdot)=d^{\,(j)}(g_1\circ f_1)(x,\cdot)$
for each~$j$, by Fa\`{a} di Bruno's Formula (Theorem~\ref{faadk})
and thus
\begin{equation}\label{givesformjet}
P_x^k(g\circ f)\;=\; P_x^k ( g_1\circ f_1 ).
\end{equation}
We can
apply (\ref{givesformjet})
to $f_1\colon E\to F$, $v\mto (P_x^kf)(v-x)$
and $g_1\colon F\to H$,\linebreak
$w\mto (P_z^kg)(w-z)$ since
$P^k_x f_1 =P^k_0(P_x^kf)=P_x^kf$
and
$P_z^k g_1=P_0^k(P_z^k g)= P^k_zg$
by Lemma~\ref{polyshift} and Proposition~\ref{charTaypo}(a).
Hence
\begin{eqnarray*}
P_x^k(g\circ f)
& = & P_x^k(g_1\circ f_1)
\; =\;
P_x^k\big((P_z^kg)(\cdot -z)
\circ (P_x^kf)(\cdot  -x)\big)\\
&=&
P_0^k\big((P_z^kg)(\cdot -z)
\circ P_x^kf \big)\,,
\end{eqnarray*}
using (\ref{givesformjet}) to obtain
the first equality
and Lemma~\ref{polyshift}
for the last.
As the final term is
$P_0^k\big( (P_z^kg)\circ (P_x^k f \, -z )\big)$,
we have established (\ref{chainjet}).

Given $y\in E$, we substitute
$v:=P^k_xf(y)-z=\sum_{i=1}^k\frac{1}{i!}\delta^i_xf(y)$ in
\[
P^k_zg(v)=g(z)+\sum_{j=1}^k\frac{1}{j!}d^{\,(j)}g(f(x),v,\ldots, v).
\]
Since $d^{\,(j)}g(f(x),\cdot )$ is $j$-linear, we find that $(P^k_zg)(P^k_xf(y)-z)$
equals
\[
g(f(x))+\sum_{j=1}^{k^2}\frac{1}{j!}\sum_{i_1+\cdots+ i_j\leq k}\frac{1}{i_1!\cdots i_j!}
d^{\,(j)}g\big(f(x),\delta^{i_1}_xf(y),\ldots, \delta^{i_j}_xf(y)\big).
\]
Discarding the summands with $j>k$ (which are homogeneous polynomials in~$y$
of degree $j>k$), we obtain (\ref{prenmfaa}).
\end{prf}
\begin{rem}
Equation (\ref{prenmfaa}) for $P^k_x(g\circ f)$ entails a version of
Fa\`{a} di Bruno's formula for G\^{a}teaux differentials of a composition $g\circ f$:
The homogeneous polynomial of degree $k\in\N$ in $P^k_x(g\circ f)$
is $\frac{1}{k!}\delta^k_x(g\circ f)$; comparing with (\ref{prenmfaa}),
we deduce that
\[
\delta^k_x(g\circ f)(y)=\sum_{j=0}^k\frac{k!}{j!} \sum_{i_1+\cdots+ i_j=k}\frac{1}{i_1!\cdots i_j!}
\,d^{\,(j)}g(f(x),\delta^{i_1}_xf(y),\ldots, \delta^{i_j}_x(y))
\]
for all $y\in E$. (This also follows from Theorem~\ref{faadk}
with $y_1,\ldots, y_k:=y$.)\footnote{However,
the combinatorics to explain the constants would be more involved.}
\end{rem}
\begin{convent}
Let $E$ and $F$ be locally convex spaces,
$U\sub E$ be an open $0$-neighborhood
and $f,g\colon U\to
F$ be two $C^k$-maps.
If $k$ is understood,
occasionally we shall simply write $f(x)=g(x)+\cdots$
if $P_0^k(f-g)=0$.
Hence $f(x+y)=P_x^kf(y)+\cdots$ in particular,
for each $C^k$-map~$f$ and~$x$ in its domain.
\end{convent}
\subsection*{Continuity of homogeneous components}
A polynomial is continuous if and only if all of its homogeneous\linebreak
components are continuous.
\begin{prop}\label{contcompo}
Let $p\!=\!\sum_{j=0}^kp_j \colon \!E\to F$ be a polynomial of degree \mbox{$\leq \! k$} between
locally convex spaces, with homogeneous components $p_j$ of degree~$j$.
If $p$ is continuous at some~$x_0\in E$, then $p$ and all the~$p_j$ are continuous.
\end{prop}
For the proof, let us recall how the coefficients of a polynomial $p\colon \K\to F$
of degree $\leq k$
can be recovered from the values of~$p$
on a suitable finite set.

Given a $\K$-vector space~$F$, let us write $P(\K,F)$
for the space of all polynomials $p\colon \K\to F$
and $P_{\leq k}(\K,F)\sub P(\K,F)$ for the vector subspace of all polynomials
of degree $\leq k$.
\begin{lem}\label{coeffviamatr}
For each $k\in\N_0$, there exists a $(k+1)\times (k+1)$-matrix $(b_{i,j})_{i,j\in\{0,1,\ldots, k\}}$
of rational numbers~$b_{i,j}$ such that, for each vector space~$F$ and each polynomial
$p\colon\K\to F$, $t\mto\sum_{i=0}^ka_it^i$ of degree $\leq k$ with
$a_i\in F$, we have\vspace{-1mm}
\[
a_i=\sum_{j=0}^k b_{i,j}p(j)\quad\mbox{for all $\, i\in \{0,1,\ldots, k\}$.}
\]
\end{lem}
\begin{prf}
For fixed $k\in\N_0$,
define the Lagrange interpolation polynomial
$p^{1,k}_j:=p_j  \in P_{\leq k}(\K,\K)$ via
\[
p_j(x):=\prod_{\ell\not=j}^k\frac{x-\ell}{j-\ell}
\]
for $j\in\{0,1\ldots, k\}$ (and $\ell$ in the same set).
Then
\[
p_j(i)=\delta_{j,i}\quad\mbox{for all $\, i\in \{0,\ldots, k\}$}
\]
and $p_j$ has rational coefficients.
Consider the linear map
\[
\phi_F\colon P_{\leq k}(\K,F)\to F^{k+1}\quad \phi_F(p):=(p(i))_{i\in\{0,\ldots, k\}}.
\]
Then $\phi_\K(p_j)=\delta_{j,\cdot}=e_j$ is the $j$th standard basis vector,
whence $\phi_\K$ is surjective and hence an isomorphism,
as both $P_{\leq k}(\K,\K)$ and $\K^{k+1}$ have dimension~$k+1$.
For $i\in\{0,\ldots, k\}$, consider the map
\[
\psi_i\colon P_{\leq k}(\K,\K)\to\K, \quad p\mto a_i
\]
for $p$ of the form $p(x)=\sum_{j=0}^ka_j\,x^j$.
Since
\[
b_{i,j}:=(\psi_i\circ(\phi_\K)^{-1})(e_j)=\psi_i(p_j)
\]
is the $i$th coefficient of~$p_j$,
we have $b_{i,j}\in\Q$.
As $\psi_i\circ (\phi_\K)^{-1}$ is linear, we get
\[
(\psi_i\circ (\phi_\K)^{-1})(y)=\sum_{j=0}^k b_{i,j}\, y_j
\]
for all $y=(y_0,\ldots, y_k)\in \K^{k+1}$.
Hence
\[
a_i=\psi_i(p)=(\psi_i\circ(\phi_\K)^{-1})(\phi_\K(p))=\sum_{j=0}^k b_{i,j}p(j)
\]
for $p\in P_{\leq k}(\K,\K)$ as before.
If $p\in P_{\leq k}(\K,F)$ with $p(x)=\sum_{j=0}^k a_j\, x^j$,
then $(\lambda\circ p)(x)=\sum_{j=0}^k \lambda(a_j)\, x^j$
for each linear functional $\lambda\colon F\to\K$.
Now $\lambda\circ p$ has $i$th coefficient
\[
\lambda(a_i)=\sum_{j=0}^k b_{i,j}\, (\lambda\circ p)(j)
=\lambda\Big(\sum_{j=0}^k b_{i,j}\, p(j)\Big).
\]
As the linear functionals~$\lambda$ separate points, $a_i=\sum_{j=0}^k b_{i,j}\, p(j)$ follows.
\end{prf}
\noindent{\em Proof of Proposition}~\ref{contcompo}.
Assume first that $p=\sum_{j=0}^kp_j$ is continuous at~$x_0=0$,
where the $p_j\colon E\to F$ are homogeneous polynomials of degree~$j$.
Let $(b_{i,j})_{i,j\in\{0,1,\ldots, k\}}$ be as in Lemma~\ref{coeffviamatr}.
For each $y\in E$, the map $\K\to F$, $t\mto p(ty)=\sum_{j=0}^kt^j p_j(y)$
is a polynomial of degree $\leq k$ and hence
\begin{equation}\label{henceconts}
p_i(y)=\sum_{j=0}^k b_{i,j}\, p(jy)
\end{equation}
for all $i\in \{0,1,\ldots, k\}$. Since $p$ is continuous at~$0$,
we deduce from~(\ref{henceconts}) that $p_i$ is continuous
at~$0$ for each~$i$.
Hence~$p_i$ is continuous, by Exercise~\ref{excctspolpt0}.

Now assume that $x_0$ is arbitrary. Then $q\colon E\to F$, $q(x):=p(x-x_0)$
is a polynomial which is continuous at~$0$, whence~$q$ is continuous by the preceding.
Hence also $p$ is continuous and thus each~$p_j$ is continuous.\qed
\subsection*{Estimates for the remainder term}
We prove estimates concerning the remainder terms in Taylor expansions.
The results can be skipped on a first reading.
\begin{lem}\label{est-if-k+1}
Let $E$ and~$F$ be locally convex spaces,
$U\sub E$ be an open subset, $k\in\N_0$
and $f\colon U\to F$ be a $C^{k+1}$-map.
If $x_0\in U$ and $p$ is a continuous seminorm on~$F$,
then there exists a continuous seminorm~$q$ on~$E$
and $r>0$ such that
$B^q_{2r}(x_0)\sub U$ and such that the Taylor remainder
\[
R_x(y):=\frac{1}{(k-1)!}\int_0^1(1-s)^{k-1}(\delta_{x+sy}^kf(y)-\delta^k_xf(y))\,ds
\]
satisfies $p(R_x(y))\leq (q(y))^{k+1}$ for all $x\in B^q_r(x_0)$ and $y\in B^q_r(0)$.
\end{lem}
\begin{prf}
The map $U\!\times \! E\to F$, $(x,y)\mto\delta^{k+1}_x(y)$ is continuous and $\delta_{x_0}^{k+1}f(0)=0$.
Thus, there is $\ve>0$ and a continuous seminorm~$Q$ on~$E$ such that $B^Q_{2\ve}(x_0)\sub U$
and $p(\delta^{k+1}_zf(y))\leq 1$ for all $z\in B^Q_{2\ve}(x_0)$ and $y\in B^Q_\ve(0)$. 
Then
\[
p(\delta^{k+1}_zf(y))\leq (Q(y)/\ve)^{k+1}\;\,
\mbox{for all $z\in B^Q_{2\ve}(x_0)$ and $y\in E$,}
\]
by Lemma~\ref{like-op-no}. For $x\in B^Q_\ve(x_0)$ and $y\in B^Q_\ve(0)$, we deduce that
\begin{eqnarray*}
p(R_x(y))&=& \frac{1}{(k-1)!}p\left(\int_0^1(1-s)^{k-1}\int_0^s\delta_{x+ty}^{k+1}(y)\,dt\,ds\right)\\
& \leq & \frac{1}{(k-1)!}\int_0^1\int_0^sp(\delta^{k+1}_{x+ty}f(y))\, dt\,ds
\leq \frac{(Q(y)/\ve)^{k+1}}{(k-1)!\,}.
\end{eqnarray*}
The assertion follows with $q:=\frac{1}{\sqrt[k+1]{(k-1)!}\,}\,Q/\ve$ and $r:=\frac{1}{\sqrt[k-1]{(k-1)!}}$.
\end{prf}
\begin{rem}\label{est-if-rem}
If $(E,\|\cdot\|_E)$ and $(F,\|\cdot\|_F)$ are normed spaces in the situation of Lemma~\ref{est-if-k+1}
and $p:=\|\cdot\|_F$, we may assume that $q=\!\sqrt[k+1]{C}\,\|\cdot\|_E$ for some $C>0$
and deduce that
\[
\|R_x(y)\|_F\leq C(\|y\|_E)^{k+1}\;\,\mbox{for all $\,x\in B^E_r(x_0)$ and $\,y\in B^E_r(0)$.}
\]
\end{rem}
\begin{lem}\label{est-if-on-lines}
Let $E$, $F$, and~$Z$ be locally convex spaces,
$U\sub E$ and $V\sub F$ be open subsets,
$(x_0,y_0)\in U\times V$ and $f\colon U\times V\to Z$ be a $C^2$-function.
Let~$p$ be a continuous seminorm on~$Z$.
If there exist continuous linear functions $\alpha\colon E\to Z$ and $\beta\colon F\to Z$
such that
\[
f(x,y_0)=f(x_0,y_0)+\alpha(x-x_0)\;\,\mbox{and}\;\,
f(x_0,y)=f(x_0,y_0)+\beta(y-y_0)
\]
for all $x\in U$ and $y\in V$,
then there exist $r>0$ and continuous seminorms~$q_1$ and~$q_2$ on~$E$ and~$F$, respectively,
such that
$B^{q_1}_r(x_0)\sub U$, $B^{q_2}_r(y_0)\sub V$, and
\[
p(f(x_0+x,y_0+y)-f(x_0,y_0)-\alpha(x)-\beta(y))\, \leq \, q_1(x)\,q_2(y)
\]
for all $x\in B^{q_1}_r(0)$ and $y\in B^{q_2}_r(0)$.
\end{lem}
\begin{prf}
After replacing $f$ with 
\[ g\colon U\times V\to Z,\quad g(x,y):=f(x,y)-f(x_0,y_0)-\alpha(x-x_0)-\beta(y-y_0),\]
we may assume that $f(x_0,y_0)=0$, $\alpha=0$, and $\beta=0$.
Since 
\[ h \colon U\times V\times E\times F\to Z,\quad h(x,y,x_1,y_1):=d^{\,(2)}f((x,y),(0,y_1),(x_1,0))\]
is continuous,
there exist continuous seminorms~$q_1$ and~$q_2$
on~$E$ and~$F$, respectively, such that $B^{q_1}_1(x_0)\sub U$, $B^{q_2}_1(y_0)\sub V$
and
$p(h(x,y,x_1,y_1))\leq 1$ for all $(x,y,x_1,y_1)\in B^{q_1}_1(x_0)\times B^{q_2}_1(y_0)\times B^{q_1}_1(0)\times B^{q_2}_1(0)$.
Since $h(x,y,\cdot)\colon E\times F\to Z$ is bilinear, Lemma~\ref{like-op-no-mult}
shows that
\begin{equation}\label{size-intg}
p(h(x,y,x_1,y_1))\leq q_1(x_1)q_2(y_1)
\end{equation}
for all $(x,y,x_1,y_1)\in B^{q_1}_1(x_0)\times B^{q_2}_1(y_0)\times E\times F$.
As $f(x_0,y)=0$ for~$y\in V$,
\begin{equation}\label{for-sec-id}
d_2f(x_0,y_0+ty;y)=\frac{d}{dt}f(x_0,y_0+ty)=0
\end{equation}
holds for all $(y,t)\in B^{q_2}_1(0)\times [0,1]$.
Using $f(x_0+x,y_0)=0$ and (\ref{for-sec-id}), we~get
\begin{eqnarray*}
f(x_0+x,y_0+y)& =& \int_0^1d_2f(x_0+x,y_0+ty;y)\,dt\\
&=& \int_0^1\int_0^1d^{\,(2)}f((x_0+sx,y_0+ty),(0,y),(x,0))\,
ds\,dt
\end{eqnarray*}
for all $x\in B^{q_1}_1(0)$ and $y\in B^{q_2}_1(0)$. Estimating the integrand with (\ref{size-intg}),
\[
p(f(x_0+x,y_0+y))\leq q_1(x)q_2(y)\;\,\mbox{for all $x\in B^{q_1}_1(0)$ and $y\in B^{q_2}_1(0)$}
\]
follows. The assertion is established with $r:=1$.
\end{prf}

\begin{small}
\subsection*{Exercises for Section~\ref{secTay}.}

\begin{exer}\label{exc-homog-then-pol}
Let $k\in\N$ and $f\colon E\to F$ be a $C^k$-map between locally convex spaces
which is homogeneous of degree~$k$ in the sense that $f(tx)=t^kf(x)$ for all
$t\in\R$ and $x\in E$. Show that $f(x)=\delta^k_0f(x)$ for all $x\in E$. In particular, $f$ is a continuous homogeneous polynomial.\\[1mm]
[Use that $\delta_0^kf(x)=\frac{d^k}{dt^k}\big|_{t=0}f(tx)=\frac{d^k}{dt^k}\big|_{t=0}t^kf(x)$.\,]
\end{exer}
\begin{exer}\label{excTaycia}
Let $\cA$ be a continuous inverse algebra.
We determine the G\^{a}teaux differentials
and Taylor polynomials
of the inversion map
$\iota\colon \cA^\times \to \cA$, $\iota(x):=x^{-1}$.
\begin{description}[(D)]
\item[(a)]
Given
$x\in \cA^\times$ and $y\in \cA$,
consider
$\gamma\colon W\to \cA^\times$,
$\gamma(z) := \iota(x+zy)$
on the open subset $W:=\{z\in \C\colon x+zy\in \cA^\times\}\sub \C$.
Then $\gamma'(z)=d\iota(x+zy ,y)= -\gamma(z)\, y\,
\gamma(z)$
by (\ref{diffinve}).
Show that
\[
\gamma^{(k)}(z) \; = \; ({-1})^k \, k!\, \gamma(z)
\, \big(y\, \gamma(z) \big)^k
\]
for each $z\in W$ and infer that
$\delta^k_x\iota(y)=
({-1})^k k! \, x^{-1}\big(yx^{-1}\big)^k$.
\item[(b)]
Determine $P^k_x\iota(y)$
for each $k$. Simplify your formulas for $x=\one$.
\end{description}
\end{exer}

\begin{exer}\label{excPolfindim}
Let $F$ be a vector space and
$n,k\in \N$.
Show that a mapping
$p\colon \K^n\to F$
is a polynomial of degree $\leq k$
in the sense of
Definition~\ref{defgenpol}
if and only if
$p$ is an ordinary polynomial
of degree $\leq k$ with coefficients
in~$F$, viz.\
there are elements $a_\alpha\in F$
indexed by multi-indices
$\alpha\in \N_0^n$
of order $|\alpha|:=\sum_{j=1}^n\alpha_j\leq k$
such that
$p(x)= \sum_{|\alpha| \leq k} x^\alpha\, a_\alpha$
for all $x=(x_1,\ldots,x_n)\in \K^n$,
where $x^\alpha:=x_1^{\alpha_1}\cdots x_n^{\alpha_n}$. 
\end{exer}

\begin{exer}\label{excnormpolvsmulti}
Let $(E,\|\cdot\|_E)$ and $(F,\|\cdot\|_F)$ be normed spaces.
The norm of a
continuous homogeneous
polynomial $p\colon E\to F$ of degree~$k$
is defined as
\[
\|p\|_{\op}:=\sup\{\|p(x)\|_F\colon x\in E, \|x\|_E\leq 1\}.
\]
Thus $\|p(x)\|_F\leq \|p\|_{\op}(\|x\|_E)^k$
for all $x\in E$.
Using polarization,
show that if\linebreak
$\beta\colon E^k\to F$ is a continuous
symmetric $k$-linear map and $p\colon E\to F$,
$x\mto\beta(x,\ldots, x)$
the associated homogeneous polynomial,
then
\begin{equation}\label{linknormpm}
\|p\|_{\op} \leq \; \|\beta\|_{\op} \leq \; \frac{k^k}{k!}\, \|p\|_{\op}\,.
\end{equation}
\end{exer}

\begin{exer}\label{excctspolpt0}
Show that
a homogeneous polynomial
$p\colon E\to F$
between
topological vector spaces
is continuous
if and only if it is continuous at~$0$
(use the Polarization Formula and Exercice~\ref{excbilinctso}(b)).
\end{exer}

\begin{exer}\label{excCkGatcts}
Let $k\in \N$ and
$f\colon E\subseteq U\to F$ be a $C^k$-map.
Show that the map $U\times E\to F$, $(x,y)\mto\delta^k_x(y)$
is continuous.
\end{exer}

\begin{exer}\label{excbabyBCH}
Given a unital Banach algebra~$\cA$,
define
$\exp\colon \cA\to \cA$, $\exp(x):=\sum_{k=0}^\infty
\frac{1}{k!}x^k$\vspace{-.4mm}
and $\log\colon B_1^\cA(\one )\to \cA$,
$\log(x):=\sum_{k=1}^\infty \frac{({-1})^{k+1}}{k}\, (x-\one)^k$.
\begin{description}[(D)]
\item[(a)]
Show that these series
converge
uniformly
on each ball $B_r^\cA(0)$,
resp., uniformly on
$B_r^\cA(\one)$ for each $r<1$.
Deduce that $\exp$ and $\log$
are continuous.
\item[(b)]
It will become clear later
that
$\exp$ and $\log$ are $C^\infty$, with
$n$th order Taylor polynomials
$P_0^n(\exp)(x)=\sum_{k=0}^n \frac{1}{k!}x^k$
and $P_{{\tiny\one} }^n(\log)(x)=
\sum_{k=1}^n \frac{({-1})^{k+1}}{k}\, x^k$
(cf.\ Proposition~\ref{singlpoint}).
Taking this for granted,
consider the $C^\infty$-map
\[
f\colon \cA\times \cA\supseteq U\to \cA\, ,\quad
f(x,y)\, :=\, \log (\exp(x)\exp(y) )\,,
\]
defined on some open
$(0,0)$-neighborhood $U\sub \cA\times \cA$.
Justified by Proposition~\ref{chainrtay},
determine the Taylor polynomial
$P^2_{(0,0)}f$
as follows:
Replace all functions
by their
second order Taylor
polynomials, multiply out
and discard all terms of order exceeding~$2$.
Calculate the 2nd order Taylor polynomial
of $\log(\exp(x)\exp(y)\exp({-x})\exp({-y}))$
at $(0,0)$ in the same way.
\end{description}
\end{exer}

\begin{exer}
Let $E$ and $F$ be $\K$-vector spaces,
$\beta\colon E^k\to F$ be a symmetric $k$-linear map and $p\colon E\to F$, $x\mto\beta(x,\ldots,x)$
be the associated polynomial.
\begin{description}[(D)]
\item[(a)]
Show that $p(E_0)$ has finite-dimensional linear span~$F_0$ in~$F$,
for each finite-dimensional vector subspace $E_0\sub E$.
Since $F_0$ can be endowed with the unique Hausdorff vector topology,
this enables us to speak of limits in~$F_0$.
\item[(b)]
Show that the directional derivatives $d^{\,(j)}p(x,y_1,\ldots,y_j):=
(D_{y_j}\cdots D_{y_1}p)(x)$ exist for all $j\in\N$, $x\in E$ and $y_1,\ldots, y_j\in E$,
and are given by
\[
d^{\,(j)}p(x,y_1,\ldots, y_j)\; =\;
\left\{
\begin{array}{cl}
\frac{k!}{(k-j)!}\,
\beta(\underbrace{x,\ldots, x}_{k-j},y_1,\ldots, y_j) &
\; \mbox{if $\,j\leq k$;}\\[2mm]
0 & \; \mbox{if $\, j>k$.}
\end{array}
\right.
\]
In particular, this yields an alternative proof that $\beta=d^{\,(k)}p(x,\cdot )$ is uniquely determined by~$p$ 
(without recourse to the Polarization Formula).
\end{description}
\end{exer}

\begin{exer}
Let $E$ and $F$ be topological vector spaces, $X$ be a topological space, $k\in\N_0$
and $f\colon X\times E\to F$ be a continuous map such that the partial map $f_x:=f(x,\cdot )\colon E\to F$
is a polynomial of degree $\leq k$ for each $x\in X$. For $j\in \{0,\ldots, k\}$,
let $(f_x)_j\colon E\to F$ be the homogeneous component of degree~$j$ of~$f_x$.
Show that $f_j\colon X\times E\to F$, $(x,y)\mto (f_x)_j(y)$
is continuous for all $j\in\{0,1,\ldots, k\}$.
\end{exer}

\begin{exer}\label{exc-equipol}
Let $E$ and $F$ be locally convex spaces,
$n\in \N$, and $P_{\leq n}(E,F)$ be the space of all
polynomals $p\colon E\to F$ of degree $\leq n$.
\begin{description}[(D)]
\item[(a)]
Let $\Gamma\sub F^E$ be a set of homogeneous
polynomials of degree~$n$. Show that $\Gamma$ is equicontinuous
if and only if~$\Gamma$ is equicontinuous at~$0$.\\[1mm]
[Use the Polarization Formula, Exercise~\ref{exc-equimult}
and Lemma~\ref{basics-equi}(a).]\vspace{1mm}
\item[(b)]
Varying the proof of Proposition~\ref{contcompo},
show that a set $\Gamma\sub P_{\leq n}(E,F)$
is equicontinuous if and only if
$\Gamma$ is equicontinuous at some~$x_0\in E$.
\end{description}
\end{exer}

\begin{exer}
Given a $\K$-vector space~$F$ and $n\in\N$, let us write $P(\K^n,F)$
for the space of all polynomials $p\colon \K^n\to F$.
\begin{description}[(D)]
\item[(a)]
For $p\in P(\K^2,F)$, show that $p^\vee(x_1):=p(x_1,\cdot)\in P(\K,F)$
for all $x_1\in \K$, and $p^\vee\in P(\K,P(\K,F))$. [Hint: If
$p(x_1,x_2)=\sum_{i_1,i_2=0}^k a_{i_1,i_2}\, x_1^{i_1}x_2^{i_2}$, then
\[
p(x_1,x_2)=\sum_{i_2=0}^k\Big(\sum_{i_1=0}^k a_{i_1,i_2}\,x_1^{i_1}\Big)x_2^{i_2}
=\sum_{i_1=0}^kq_{i_1}(x_2)\, x_1^{i_1}
\]
with $q_{i_1}(x_2):=\sum_{i_2=0}^ka_{i_1,i_2}\, x_2^{i_2}$.]
\item[(b)]
For $k\in\N_0$, let $(b_{i,j})_{i,j=0}^k$ be the $(k+1)\times (k+1)$-matrix with rational
entries from Lemma~\ref{coeffviamatr}.
Show that
\[
a_{i_1,i_2}=\sum_{j_1,j_2=0}^k b_{i_1,j_1}b_{i_2,j_2}\,p(j_1,j_2)
\]
for all $i_1,i_2\in\{0,1,\ldots, k\}$ and $p\in P(\K^2,F)$ of multi-degree $\leq (k,k)$ as in~(a).
\item[(c)]
Find analogs for $n>2$.
\item[(d)]
Show that $\Phi\colon P(\K^2,F)\to P(\K,P(\K,E))$, $p\mto p^\vee$
is an isomorphism of $\K$-vector spaces.
\end{description}
\end{exer}

\begin{exer}
For $n\in\N$ and $k\in\N_0$, abbreviate $\Delta_{n,k}:=\{\alpha\in \N_0^n\colon |\alpha|\leq k\}$,
where $|\alpha|:=\alpha_1+\cdots+\alpha_n$ for $\alpha=(\alpha_1,\ldots,\alpha_n)\in\N_0^n$.
Given a $\K$-vector space $F$, let $P_{\leq k}(\K^n,F)$ be the space of all $F$-valued polynomials
on $\K^n$ of degree $\leq k$, i.e, $p(x)=\sum_{\alpha\in \Delta_{n,k}} a_\alpha \, x^\alpha$
with suitable $a_\alpha\in F$.
\begin{description}[(D)]
\item[(a)]
Let us show by induction on~$n\in\N$ that, for each $k\in\N_0$,
there are polynomials
$p^{n,k}_\alpha\in P_{\leq k}(\K^n,\K)$ for $\alpha\in\Delta_{n,k}$
such that
\begin{equation}\label{willdetp}
p^{n,k}_\alpha(\beta)=\delta_{\alpha,\beta}\quad
\mbox{for all $\,\beta\in\Delta_{n,k}$},
\end{equation}
using Kronecker's delta.
If $n=1$, for $p^{1,k}_j$
we can take the Lagrange interpolation polynomial encountered in the proof of Lemma~\ref{coeffviamatr},
for all $k\geq 1$ and $j\in\{0,\ldots, k\}$.
If $n\geq 2$ and polynomials $p^{n-1,k}_\alpha$ have been found for $k\in\N_0$
and $\alpha\in \Delta_{n-1,k}$, we find $p^{n,k}_\alpha$ for $\alpha\in\Delta_{n,k}$
by induction on $k\in\N_0$, as follows:
If $k=0$, we take $p^{n,0}_0:=1$.
If $k\geq 1$ and the $p^{n,k-1}_\beta$ have been found for $\beta\in\Delta_{n,k-1}$,
let $\alpha\in\Delta_{n,k}$ and verify that $p^{n,k}_\alpha$ with the required properties can
be obtained in the following fashion: If $|\alpha|<k$, set
\[
p^{n,k}_\alpha(x):=\frac{k-x_1-\cdots - x_n}{k-\alpha_1-\cdots - \alpha_n}p^{n,k-1}_\alpha(x)
\quad\mbox{for $x=(x_1,\ldots, x_n)\in\K^n$.}
\]
If $|\alpha|=k$, then $\alpha$ lies in the hyperplane $H:=\{x\in\K^n\colon x_1+\cdots+x_n=k\}$
and $\alpha_n=k-\alpha_1-\cdots-\alpha_{n-1}$ is determined by $(\alpha_1,\ldots,\alpha_{n-1})\in\Delta_{n-1,k}$.
The map
\[
\psi\colon \K^{n-1}\to H,\quad (x_1,\ldots,x_{n-1})\mto \Big(x_1,\ldots, x_{n-1},k-\sum_{j=1}^{n-1}x_j\Big)
\]
is a bijection which takes $\Delta_{n-1,k}$ onto $\Delta_{n,k}\cap H$.
Moreover, $\psi^{-1}=\pr|_H$, where $\pr\colon \K^n\to\K^{n-1}$, $(x_1,\ldots, x_n)\mto (x_1,\ldots, x_{n-1})$.
Then $q:=p^{n-1,k}_{\psi^{-1}(\alpha)}\circ \pr\in P_{\leq k}(\K^n,\K)$
and we can take
\[
p^{n,k}_\alpha\; :=\; q\; -\!\!\! \sum_{\beta\in \Delta_{n,k}\setminus H}q(\beta)\, p^{n,k}_\beta,
\]
where the right hand side is already defined as $\Delta_{n,k}\setminus H=\{\beta\in\Delta_{n,k}\colon
|\beta|<k\}$.
Show that all polynomials constructed have coefficients in the field~$\Q$.
\item[(b)]
For $n\in\N$ and $k\in\N_0$, consider the linear map
\[
\phi_F\colon P_{\leq k}(\K^n,F)\to F^{\Delta_{n,k}},\quad \phi_F(p):=(p(\alpha))_{\alpha\in \Delta_{n,k}}.
\]
Show as in the proof of Lemma~\ref{coeffviamatr}
that $\phi_F$ is an isomorphism and that there are families
$(b^{n,k}_{\alpha,\beta})_{\alpha,\beta\in \Delta_{n,k}}$
of rationals such that
\begin{equation}\label{gtcoef}
a_\alpha=\sum_{\beta\in \Delta_{n,k}}b^{n,k}_{\alpha,\beta}\, p(\beta)
\end{equation}
for all $p\in P_{\leq k}(\K^n,F)$ as above.
Using that $\phi_\K$ is injective,
also show that each $p^{n,k}_\alpha\in P_{\leq k}(\K^n,\K)$ is uniquely
determined by~(\ref{willdetp}).
\end{description}
\end{exer}

\end{small}

\section[Spaces of $C^k$-functions and mappings between them]{Spaces of {\boldmath$C^k$-}functions and mappings between
them}\label{secCspaces}
Function spaces and differentiability properties of typical mappings
between such spaces are essential both for the theory of infinite-dimensional Lie groups,
and for many important examples.
We now study such mappings, in three stages.

First, we discuss differentiability properties of mappings of the form
\[
C(K,f)\colon C(K,U)\to C(K,F),\quad \gamma\mto f\circ \gamma
\]
(and related maps), if $K$ is a compact topological space, $E$ and $F$ are locally convex spaces
and $f\colon U\to F$ a $C^k$-map on an open subset $U\sub E$.
The results obtained (Propositions~\ref{superctisCk} and \ref{pfctisCk}) can be used to turn $C(K,G)$
into a Lie group for each Lie group~$G$,
once the local description of Lie groups (Theorem~\ref{thm:locglob}) is available.
Besides unit groups of continuous inverse algebras,
these are among the easiest examples of infinite-dimensional Lie groups.
A variant of the cited propositions involving parameters (Proposition~\ref{pfctisCkpar})
will be used (in Section~\ref{provis-Ban}) to study the parameter-dependence of solutions to
differential equations in Banach spaces.
Moreover, it will be used to show that both Banach--Lie groups
and diffeomorphism groups are well-behaved Lie groups, 
so-called $C^0$-regular 
Lie groups (\red{insert references}),
in which continuous Lie algebra-valued curves admit well-behaved evolutions
in the group (see Definition~\ref{def:regular0} for details).

We then introduce a natural topology on spaces $C^k(U,F)$ of $C^k$-functions
for $U$ an open subset of a locally convex space and show that 
\[
C^k(U,f)\colon C^k(U,E)\to C^k(U,F),\quad\gamma\mto f\circ\gamma
\]
is continuous for each $C^k$-map $f\colon E\to F$ (Proposition~\ref{Crpushf}). This result (and Proposition~\ref{Crpuba})
will be used when we introduce topologies on spaces of sections
in vector bundles (see Section~\ref{secspacemfd}). 
We also study properties like completeness and metrizability
for the function spaces $C^k(U,E)$.

The high point of this section is the exponential law, which shows that smooth functions $f\colon U\to C^\infty(V,F)$
with values in space of smooth functions simply correspond to smooth functions $U\times V \to F$
in two variables (Corollary~\ref{cinfty-explaw}), in good cases.
More generally, we study when $C^k$-functions $U\to C^\ell(V,E)$
correspond to suitably-defined $C^{k,\ell}$-functions $U\times V\to F$
with different degrees of differentiability in its two arguments (Theorem~\ref{explawCkell}).
Exponential laws are powerful tools of infinite-dimensional calculus.
We shall use them frequently, both in the development of infinite-dimensional Lie theory
(e.g., for the discussion of $C^k$-regularity), and for the discussion of examples
(like the construction of the Lie group structure on diffeomorphism groups
and the proof of their regularity).

Some of the results are technical and it may be advisable to read them
only
once they are needed, to have motivating applications directly at hand.
\subsection*{Mappings between spaces of continuous functions}
If $K$ is a compact topological space
and $E$ a (real or complex) locally convex space,
then the compact-open topology turns the vector space $C(K,E)$
of all continuous $E$-valued functions on~$K$
into a locally convex topological vector space (see Lemma~\ref{sammelsu}(g)).
We now discuss differentiability properties of important mappings between
such function spaces and open subsets thereof.
In Appendix~\ref{appcotop}, we already established \emph{continuity} properties
for mappings between function spaces;
these facts are the basis for the following discussion, and we urge the reader
to consult the appendix for the necessary background.

We begin with a result concerning differentiability properties
of superposition operators (the proof of which is covered by the following more general theory,
as explained in Remark~\ref{remvarspec}(a)):
\begin{prop}\label{superctisCk}
Let $K$ be a compact topological space, $E$ and $F$ be locally convex spaces,
$U\sub E$ be an open subset, $k\in\N_0\cup\{\infty\}$ and
$g\colon U\to F$ be a $C^k$-map. Then also the following map is~$C^k$:
\[
C(K,g)\colon C(K,U)\to C(K,F),\quad \gamma\mto g\circ \gamma.
\]
\end{prop}
Recall that $C(K,U)=\lfloor K,U\rfloor$ is open in $C(K,E)$ if $U\sub E$ is open;
moreover, the topology induced by $C(K,E)$ on $C(K,U)$ coincides with the compact-open topology
(see Remark~\ref{reminduco}).

If $f\colon K\times U\to F$ is a continuous map in the preceding situation,
then
\[
f_*(\gamma):=f\circ (\id_K,\gamma)\in C(K,F)
\]
for all $\gamma\in C(K,U)$
and Lemma~\ref{ctspfapp}
shows that the map
\[
f_*\colon C(K,U)\to C(K,F),\quad \gamma\mto f_*(\gamma)
\]
so obtained is continuous
(the so-called ``pushforward").
More explicitly,
\[
f_*(\gamma)(x):=f(x,\gamma(x))\quad \mbox{for $\,x\in K$.}
\]
We want to show that $f_*$ is $C^k$ if $f$ is $C^k$.
In fact, this conclusion will hold if $f$ is
merely a so-called $C^{0,k}$-map;
this is important for some of our applications.
\begin{defn}\label{defnC0k}
Let $X$ be a Hausdorff topological space, $F$ be a locally convex space,
$U$ be an open subset of a locally convex space~$E$, and $k\in \N_0\cup\{\infty\}$.
A mapping $f\colon X\times U\to F$ is called $C^{0,k}$
if $f$ is continuous, the iterated directional derivatives
\[
d^{\,(0,j)}f(x,y,w_1,\ldots, w_j):=(D_{w_j}\cdots D_{w_1}f_x)(y)
\]
of $f_x:=f(x,\cdot)\colon U\to F$ exist for all $x\in X$, $y\in U$, $j\in\N$ such that $j\leq k$ and $w_1,\ldots, w_j\in E$,
and the mappings $d^{\,(0,j)}f\colon X\times U\times E^j\to F$
are continuous.\footnote{Thus
$d^{\,(0,j)}f(x,y,w_1,\ldots, w_j):=(D_{(0,w_j)}\cdots D_{(0,w_1)}f)(x,y)$
if~$X$ is a subset of a locally convex space.}
If $U\sub E$ is merely a locally convex subset with dense interior,
then a map $f$ (as before) is called $C^{0,k}$ if~$f$ is continuous,
$f|_{X\times U^0}$ is $C^{0,k}$, and $d^{\,(0,j)}(f|_{X\times U^0})\colon X\times U^0\times E^j\to F$
admits a (necessarily unique) continuous extension $d^{\,(0,j)}f\colon X\times U\times E^j\to F$,
for all $j\in\N$ such that $j\leq k$.
We mention that $C^{k,0}$-maps $f\colon U\times X\to F$
and
their partial differentials
$d^{\,(j,0)}f\colon U\times X\times E^j\to F$
can be defined following the same pattern
(reversing the roles of $x$ and $y$).
\end{defn}
\begin{rem}
(a) We already encountered $C^{k,0}$-maps in Proposition~\ref{diffpar}
on differentiable dependence of integrals on parameters (without attaching a name
to such maps), and also in Lemma~\ref{prehighdf}.
In the earlier notation,
\[
d^{\,(i,0)}f=d^{\,(i)}_1f\quad\mbox{and}\quad
d^{\,(0,j)}f=d^{\,(j)}_2f.
\]
Our new notation prepares the consideration of $C^{k,\ell}$-maps
on products with different orders $k$ and $\ell$ of differentiability
in the two variables, and their differentials $d^{\,(i,j)}f$ for $i,j\in \N_0$
such that $i\leq k$ and $j\leq \ell$ (see Definition~\ref{defnCkell}).\medskip

\noindent
(b)
If $E_1$, $E_2$, and $F$ are locally convex spaces and $U_1\sub E_1$, $U_2\sub E_2$
open subsets,
then every $C^k$-map $f\colon U_1\times U_2\to F$ is $C^{0,k}$,
with
\begin{equation}\label{rhsisnice}
d^{\,(0,j)}f(x,y,w_1,\ldots, w_j)=d^{\,(j)}f(x,y,(0,w_1),\ldots, (0,w_j))
\end{equation}
for all $x\in U_1$, $y\in U_2$, $j\in\N$ with $j\leq k$ and $w_1,\ldots, w_j\in E_2$.
This is clear from the definitions. If $U_1$ and $U_2$ are merely locally convex subsets
with dense interior, then the right hand side of (\ref{rhsisnice}) provides a continuous extension
for $d^{\,(0,j)}(f|_{U_1\times U_2^0})$, whence again $f$ is $C^{0,k}$.
Likewise, every $C^k$-map on a product is $C^{k,0}$.
\end{rem}
Now our result on pushforwards reads as follows.
\begin{prop}\label{pfctisCk}
Let $K$ be a compact topological space, $E$ and $F$ be locally convex spaces,
$U\sub E$ be an open subset, $k\in\N_0\cup\{\infty\}$, and
$f\colon K\times U\to F$ be a $C^{0,k}$-map. Then the following map is~$C^k$:
\[
f_*\colon C(K,U)\to C(K,F),\quad \gamma\mto f\circ (\id_K,\gamma).
\]
\end{prop}
We shall deduce Proposition~\ref{pfctisCk} from a
variant with parameters.
In the following result, we abbreviate
\begin{equation}\label{fstapreu}
f^p(x,y):=f(x,y,p)\quad \mbox{for $(x,y,p)\in K\times U\times P$}
\end{equation}
and $(f^p)_*(\gamma)(x)=f(x,\gamma(x),p)$ for $x\in K$, $\gamma\in C(K,U)$ and
$p\in P$.
\begin{prop}\label{pfctisCkpar}
Let $E$, $F$, and $Z$ be locally convex spaces, $U\sub E$ be an open subset,
$P\sub Z$ be a locally convex subset with dense interior,
$K$ be a compact topological
space, $k\in\N_0\cup\{\infty\}$ and
$f\colon K\times (U\times P) \to F$ be 
a $C^{0,k}$-map.
Then the following map is $C^k$:
\[
\Phi\colon C(K,U)\times P\to C(K,F), \quad (\gamma,p)\mto (f^p)_*(\gamma).
\]
\end{prop}
Three lemmas will help us to prove Proposition~\ref{pfctisCkpar}.
In the first, we use $U^{[1]}\sub U\times E\times\K$ as in (\ref{defnU1}) and (\ref{defnU1no}).
\begin{lem}\label{f1withpar}
Let $P$ be a Hausdorff topological space, $E$ and $F$ be locally convex spaces,
$U\sub E$ be a locally convex subset with dense interior and $f\colon P\times U\to F$ be a $C^{0,1}$-map.
Then the following map is continuous:
\[
f^{[0,1]}\colon P\times U^{[1]}\to F,\;\,
f^{[0,1]}(p,x,y,t):=\left\{
\begin{array}{cl}
{\displaystyle\frac{1}{t}}\big(f(p,x+ty)-f(p,x)\big) & \mbox{if $\,t\not=0$;}\\[1mm]
{\displaystyle d^{\,(0,1)}f(p,x,y)}&\mbox{if $\,t=0$.}
\end{array}
\right.
\]
\end{lem}
We mention that
Lemma~\ref{f1withpar} can be proved like Lemmas~\ref{linkBGN} and~\ref{linkBGNno};
likewise, the following lemma
can be proved like Lemma~\ref{prehighdf}
(see also Remark~\ref{prehighdfno}):
we only need to insert the parameter $p\in P$ in all
formulas, and check continuous dependence also on~$p$.
\begin{lem}\label{multidiffauto}
Let $E$ and $F$ be locally convex spaces, $\ell\in\N$, and $(W_j)_{1\leq j\leq \ell}$ be a family
of locally convex spaces $W_j$ for
$j\in\{1,\ldots,\ell\}$.
Abbreviate $W:=W_1\times\cdots\times W_\ell$.
Let $U\sub E$ be a locally convex subset with dense interior, $P$ be a Hausdorff
topological space,
$r\in\N\cup\{\infty\}$
and
\[
f\colon U\times (W\times P) \to F
\]
be a $C^{r,0}$-map such that
\[
f(x,\cdot ,p)\colon W_1\times\cdots\times W_\ell\to F
\]
is $\ell$-linear for all $x\in U$ and $p\in P$.
Then $f$ is $C^{r,0}$ also as a map from $(U\times W)\times P$ to~$F$.
\end{lem}
We record a typical application, for later use.
\begin{lem}\label{diffphigherd}
Let $X$ be a Hausdorff topological space, $E$ and $F$ be locally convex spaces, $U\sub E$ be
a locally convex subset with dense interior, $k\in\N_0\cup\{\infty\}$
and
$f\colon X \times U \to F$ be a $C^{0,k}$-map. Then $d^{\,(0,j)}f\colon X \times (U\times E^j)\to F$
is $C^{0,k-j}$, for all $j\in \N$ such that $j\leq k$.
\end{lem}
\begin{prf}
By Exercise~\ref{exc-babychainC0k}, we need only show that
$g\colon X\times (E^j\times U)\to F$, $g(x,w,y):=d^{\,(0,j)}f(x,y,w)$
is $C^{0,k-j}$ (where $x\in X$, $y\in U$, $w\in E^j$).
It is clear that the restriction of~$g$ to $(X\times E^j)\times U^0$
is $C^{0,k-j}$, with
\begin{equation}\label{rhswiext}
d^{\,(0,i)}g((x,w),y,z_1,\ldots, z_i)=d^{\,(0,j+i)}f(x,y,w,z_1,\ldots , z_i)
\end{equation}
for all $x\in X$, $y\in U^0$, $w\in E^j$ and $z_1,\ldots, z_i\in E$.
As the right had side of~(\ref{rhswiext}) defines a continuous
$F$-valued map~$h_i$ on $X\times E^j\times U\times E^i$,
we see that $g$ is $C^{0,k-j}$ on $(X\times E^j)\times U$, with $d^{\,(0,i)}g=h_i$ for all~$i$.
Since $g(x,w,y)$ is $j$-linear in $w\in E^j$,
Lemma~\ref{multidiffauto}
implies that $g$ is $C^{0,k-j}$ also as a map from
$X\times (E^j\times U)$ to~$F$.
\end{prf}
\noindent
\emph{Proof of Proposition}~\ref{pfctisCkpar}.
We may assume that $k\in\N_0$ and proceed by induction.
By Lemma~\ref{pushpar}, the map~$\Phi$ is
continuous. Now let $k\geq 1$ and assume the assertion holds for $k-1$ in place of~$k$.
Since the mapping $d^{\,(0,1)}f\colon K\times (U\times P\times E\times Z)\to F$ is $C^{0,k-1}$ (see Lemma~\ref{diffphigherd}),
also
\[
g\colon K\times (U\times E\times P\times Z)\to F,\quad g(x,y,z,p,q):=d^{\,(0,1)}f(x,y,p,z,q)
\]
is $C^{0,k-1}$ (see Exercise~\ref{exc-babychainC0k}). Hence
\[
\Psi\colon C(K,U\times E) \times (P\times Z)\to C(K,F),\quad (\gamma,\eta,p,q)\mto (g^{p,q})_*(\gamma,\eta)
\]
is $C^{k-1}$, by the inductive hypothesis (here and in the following,
$C(K,U\times E)$ is identified with $C(K,U)\times C(K,E)$ as in Lemma~\ref{cotopprod}).
We claim that the directional derivative
of $\Phi$ exists at each
$(\gamma,p)\in C(K,U)\times P^0$
in each direction $(\eta,q)\in C(K,E)\times Z$,
and is given by
\begin{equation}\label{thedirrdhere}
d\Phi((\gamma,p),(\eta,q))=\Psi(\gamma,\eta,p,q).
\end{equation}
As the right hand side of (\ref{thedirrdhere}) defines a continuous function of $(\gamma,p,\eta,q)\in C(K,U)\times P\times C(K,E)\times Z$, this implies that $\Phi$ is $C^1$ and that $d\Phi$
is given by (\ref{thedirrdhere}) on all of its domain; hence $d\Phi$ will be a $C^{k-1}$-map and
so $\Phi$ will be $C^k$.
To prove the claim, pick $\ve>0$
such that $\gamma(K)+\bD_\ve \eta(K)\sub U$
and $p+\bD_\ve q\sub P$.
Then $(x,y,p,z,q,t)\in K\times (U \times P)^{[1]}$
for all $(x,y,z,t)\in K\times \gamma(K)\times\eta(K)\times \bD_\ve$
and the map
\[
h\colon K\times \gamma(K)\times \eta(K)\times \bD_\ve\to F,\quad
h(x,y,z,t):=f^{[0,1]}(x,y,p,z,q,t)
\]
is continuous. For $0\not=t\in\bD_\ve$, consider the difference quotient
\[
\Delta_t:= \frac{\Phi(\gamma+t\eta,p+tq)-\Phi(\gamma,p)}{t}\in C(K,F).
\]
For $x\in K$, we have
\begin{eqnarray*}
\Delta_t(x)&=&\frac{f(x,\gamma(x)+t\eta(x),p+tq)-f(x,\gamma(x),p)}{t}\\
&=& f^{[0,1]}(x,\gamma(x),p,\eta(x),q,t) \,=\, (h^t)_*(\gamma,\eta)(x).
\end{eqnarray*}
Hence $\Delta_t=(h^t)_*(\gamma,\eta)$ for $t\in\bD_\ve\setminus\{0\}$,
where the right hand side is a continuous $C(K,F)$-valued function of $t\in\bD_\ve$
by Lemma~\ref{pushpar}. The difference quotient $\Delta_t=(h^t)_*(\gamma,\eta)$
therefore converges to $(h^0)_*(\gamma,\eta)$ as $t\to 0$.
Since $(h^0)_*(\gamma,\eta)(x)=d^{\,(0,1)}f(x,\gamma(x),p,\eta(x),q)=\Psi(\gamma,\eta,p,q)(x)$,
we have established~(\ref{thedirrdhere}). This completes the proof.\qed

\begin{rem}\label{remvarspec} (a)
Taking $P$ as a singleton, Proposition~\ref{pfctisCk} follows from Proposition~\ref{pfctisCkpar}.
Taking $f(x,y):=g(y)$, Proposition~\ref{superctisCk} becomes a special case of Proposition~\ref{pfctisCk}.\medskip

\noindent
(b) We mention that Proposition~\ref{pfctisCkpar}
(and hence also Proposition~\ref{pfctisCk} and Proposition~\ref{superctisCk})
remains valid if~$U$ is, instead, a convex subset of~$E$ with non-empty interior.
In fact, $C(K,U)^0=C(K,U^0)$ in this case (Exercise~\ref{exc-intmpset}).
In particular, the convex set $C(K,U)$ has non-empty interior,
whence $C(K,U)\times P$ is a locally convex subset of $C(K,E)\times Z$ with dense
interior. Now
the case $k=0$ in the inductive proof of Proposition~\ref{pfctisCkpar}
does not require openness of~$U$.
In the induction step, the directional derivatives still have the form described
in (\ref{thedirrdhere}) for $(\gamma,p,\eta,q)\in C(K,U^0)\times P^0\times C(K,E)\times Z$
(by the proof as it stands). As~$\Psi$ is~$C^{k-1}$ by induction and hence
continuous, we deduce that~$\Phi$ is~$C^1$ with~$d\Phi$ given by (\ref{thedirrdhere})
on all of its domain.
Hence~$d\Phi$ is~$C^{k-1}$ and~$\Phi$ is~$C^k$.
\end{rem}
\subsection*{Continuity of pullbacks and global pushforwards between spaces
of {\boldmath$C^k$}-functions}
We now turn the space $C^k(U,F)$ of $C^k$-maps $\gamma\colon U\to F$ into a locally convex space,
when $U$ is a locally convex subset with dense interior in a locally convex space~$E$ and $k\in\N_0\cup\{\infty\}$.
Continuity is established for two types of mappings between such spaces; these
are needed later to discuss spaces of~$C^k$-maps on manifolds, and spaces of $C^k$-sections in vector bundles.
\begin{defn}\label{firstdefCktop}
Let $E$ and $F$ be locally convex spaces, $U\sub E$ be a locally convex subset with dense
interior and $k\in\N_0\cup\{\infty\}$.
For $j\in\N_0$ such that $j\leq k$, endow $C(U\times E^j,F)$ with the compact-open topology.
The initial topology on $C^k(U,F)$ with respect to the mappings
\[
d^{\,(j)}\colon C^k(U,F)\to C(U\times E^j,F),\quad \gamma\mto d^{\,(j)}\gamma,
\]
(for $j$ as before) is called the 
\emph{compact-open $C^k$-topology}.
\index{compact-open $C^k$-topology} 
As the $d^{\,(j)}$ form a point-separating family of linear maps
from $C^k(U,E)$ to locally convex spaces,
also $C^k(U,E)$ is a locally convex topological vector space.
\end{defn}
\begin{prop}\label{Crpuba}
Let $E_1$, $E_2$, and $F$ be locally convex spaces, $U\sub E_1$ and $V\sub E_2$
be locally convex subsets with dense interior, $r\in\N_0\cup\{\infty\}$
and $f\colon U\to E_2$ be a $C^r$-map such that $f(U)\sub V$.
Then the following map is continuous and linear:
\[
f^*\colon C^r(V,F)\to C^r(U,F),\quad\gamma\mto\gamma\circ f.
\]
\end{prop}
\begin{prf}
The map $f^*$ is continuous to $C(U,F)$ by Lemma~\ref{pubas}. It remains to show that
$d^{\,(k)}(f^*(\gamma))=d^{\,(k)}(\gamma\circ f)\in C(U\times E_1^k,F)$
is continuous in $\gamma\in C^r(V,F)$, for all $k\in\N$ such that $k\leq r$.
By Fa\`{a} di Bruno's formula (Theorem~\ref{faadk}),
$d^{\,(k)}(\gamma\circ f)$ is the sum of summands $h_P(\gamma)$ given by
\[
h_P(\gamma)(x,y)=d^{\,(j)}\gamma(f(x),d^{(|I_1|)}f(x,y_{I_1}),\ldots,d^{\,(|I_j|)}f(x,y_{I_j}))
\]
for $j\in\{1,\ldots, k\}$, $P=\{I_1,\ldots, I_j\}\in P_{k,j}$
and $(x,y)\in U\times E_1^k$.
To see that each $h_P$ (and hence $d^{\,(k)}\circ f^*$) is continuous, note that
\[
g_P\colon U\times E_1^k\to V\times E_2^j,\quad (x,y)\mto (f(x),d^{(|I_1|)}f(x,y_{I_1}),\ldots,d^{\,(|I_j|)}f(x,y_{I_j}))
\]
is a continuous map. Hence also $g_P^*\colon C(V\times E_2^j,F)\to C(U\times E_1^k,F)$, $\eta\mto \eta\circ g_P$
is continuous (see Lemma~\ref{pubas}). Since $d^{\,(j)}\colon C^r(V,F)\to C(V\times E_2^j,F)$, $\gamma\mto d^{\,(j)}\gamma$
is continuous, indeed $h_P=g_P^*\circ d^{\,(j)}$
is continuous.
\end{prf}
\begin{prop}\label{Crpushf}
Let $Z$, $E$, and $F$ be locally convex spaces, $U\sub Z$ be a locally convex subset
with dense interior, $r\in\N_0\cup\{\infty\}$ and
$f\colon U\times E\to F$ be a $C^r$-map. Then the following map is continuous:
\[
f_*\colon C^r(U,E)\to C^r(U,F),\quad \gamma\mto f\circ (\id_U,\gamma).
\]
\end{prop}
\begin{prf}
Let $k\in\N_0$ such that $k\leq r$; we have to show that
$d^{\,(k)} \circ f_*\colon C^r(U,E)\to C(U\times Z^k,F)$ is continuous.
For $k=0$, this holds by Lemma~\ref{ctspfapp}.
If $k>0$, 
we have $d^{\,(k)}(f_*(\gamma))=d^{\,(k)}(f\circ\wt{\gamma})$
with $\wt{\gamma}:=(\id_U,\gamma)\colon U\to U\times E$.
Hence Fa\`{a} di Bruno's Formula (Theorem~\ref{faadk}) shows that
$d^{\,(k)}(f_*(\gamma))(x,y_1,\ldots, y_k)$ (for $x\in U$ and $y=(y_1,\ldots, y_k)\in Z^k$)
is the sum of the summands $h_P(\gamma)\in C(U\times Z^k,F)$ given by
\[
h_P(\gamma)(x,y):=d^{\,(j)}f(\wt{\gamma}(x),d^{\,(|I_1|)}\wt{\gamma}(x,y_{I_1}),\ldots,d^{\,(|I_j|)}\wt{\gamma}(x,y_{I_j}))
\]
for $j\in\{1,\ldots, k\}$ and $P=\{I_1\ldots, I_j\}\in P_{k,j}$ (with notation as in
Theorem~\ref{faadk}). To see that $h_P$ is continuous,
we consider for $I\sub \{1,\ldots, k\}$ the continuous map
\[
g_I\colon U\times Z^k\to U\times Z^{|I|},\quad (x,y)\mto (x,y_I);
\]
thus $g_\emptyset$ is the projection $U\times Z^k\to U$
$(x,y)\mto x$.
Then
\[
g_I^*\colon C(U\times Z^{|I|},Z\times E)\to C(U\times Z^k,Z\times E),\quad\eta\mto\eta\circ g_I
\]
is continuous (by Lemma~\ref{pubas}). 
We can now write
\[
h_P=C(U\times Z^k,d^{\,(j)}f)\circ (g_\emptyset^*,g_{I_1}^*\circ d^{\,(|I_1|)},\cdots, g_{I_j}^*\circ d^{\,(|I_j|)})\circ \widetilde{\cdot},
\]
where the map
$C^r(U,E)\to C^r(U,U\times E)$, $\gamma \mto \wt{\gamma}=(\id_U,\gamma)$
is continuous,\footnote{As $d^{\,(k)}(\id_U,\gamma)=(d^{\,(k)}\id_U,d^{\,(k)}\gamma)\in C(U\times Z^k,Z\times E)$
depends continuously on $\gamma\in C^r(U,E)$.} the map
$d^{\,(i)}\colon C^r(U,U\times E)\to C(U\times Z^i,Z\times E)$
is continuous for $i\in\N$ such that $i\leq r$,
and, finally,
\[
C(U\times Z^k,d^{\,(j)}f)\colon C(U\times Z^k,U\times E\times(Z\times E)^j)\to C(U\times Z^k,F)
\]
is continuous by Lemma~\ref{covsuppo}. Hence $h_P$ is continuous, as required.
\end{prf}
With $r:=k$ and $f(x,y):=g(y)$, we deduce:
\begin{cor}\label{Cksuppo}
Let $Z$, $E$, and $F$ be locally convex spaces, $U\sub Z$ be a locally convex subset
with dense interior, $k\in\N_0\cup\{\infty\}$ and
$g\colon E\to F$ be a $C^k$-map. Then the following map is continuous:
\[
C^k(U,g)\colon C^k(U,E)\to C^k(U,F),\quad \gamma\mto g\circ \gamma.
\]
\end{cor}
\begin{cor}\label{Ckrespprod}
Let $Z$, $E_1$, and $E_2$ be locally convex spaces, $U\sub Z$ be a locally convex subset
with dense interior, $k\in\N_0\cup\{\infty\}$ and $\pr_j\colon E_1\times E_2\to E_j$ be the projection
onto the $j$th component, for $j\in\{1,2\}$.
Then the map 
\[
\Phi:=(C^k(U,\pr_1),C^k(U,\pr_2))\colon C^k(U,E_1\times E_2)\to C^k(U,E_1)\times
C^k(U,E_2)
\]
taking $\gamma\in C^k(U,E_1\times E_2)$ to the pair $(\pr_1\circ \gamma,\pr_2\circ\gamma)$
of its components is an isomorphism of topological vector spaces.
\end{cor}
\begin{prf}
Let $\lambda_1\colon E_1\to E_1\times E_2$, $x\mto(x,0)$ and $\lambda_2\colon E_2\to E_1\times E_2$,
$y\mto(0,y)$ be the inclusion maps. For $j\in \{1,2\}$, let $\pi_j\colon C^k(U,E_1)\times C^k(U,E_2)\to C^k(U,E_j)$
be the projection onto the $j$th component.
Since $\lambda_1$ and $\lambda_2$ are
continuous linear and hence smooth, we deduce with Corollary~\ref{Cksuppo} that
\[
\Phi^{-1}=C^k(U,\lambda_1)\circ \pi_1+C^k(U,\lambda_2)\circ \pi_2
\]
is continuous.
\end{prf}
\subsection*{Elementary properties of the spaces {\boldmath$C^k(U,F)$}}
We now compile some properties of the function spaces $C^k(U,F)$,
notably if $U$ is a subset of a finite-dimensional
$\K$-vector space.
\begin{lem}\label{ima-clo}
Let $E$ and $F$ be locally convex spaces, $k\in\N_0\cup\{\infty\}$
and $U\sub E$ be a locally convex subset with dense interior.
Then the linear map
\[
\Phi\colon C^k(U,F)\to\prod_{j\leq k}C(U\times E^j,F),\quad f\mto (d^jf)_{j\leq k}
\]
$($where $j\in\N_0$ with $j\leq k)$
is a topological embedding with closed image.
\end{lem}
\begin{prf}
It is clear that~$\Phi$ is linear and it is a topological
embedding by definition of the compact-open $C^k$-topology.
To see that $\Phi$ has closed image,
let $(f_\alpha)_{\alpha\in A}$ a net in $C^k(U,F)$
such that
\[
\Phi(f_\alpha)\to h
\]
for some $h=(h_j)_{j\leq k}\in\prod_{j\leq k}C(U\times E^j,F)$.
Abbreviate $f:=h_0$.
We claim that the iterated directional derivative $d^{\,(j)}f(x,y_1,\ldots,y_j)$
exists for all $j\in\N_0$ with $j\leq k$, $x\in U^0$ and $(y_1,\ldots,y_j)\in E^j$.
If this is true, then $h_j$ will be a continuous extension for $d^{\,(j)}(f|_{U^0})$
and thus $f$ will be $C^k$ with $d^{\,(j)}f=h_j$ for all $j$, whence
$\Phi(f)=h$.

To prove the claim, we may assume that~$U$ is open.\footnote{An alternative idea of proof would
be a reduction to Lemma~\ref{Ckinlimit}. However, we prefer the following
arguments, which are not longer and more elementary.}
It suffices to work over the ground field~$\R$
(if $k\geq 1$ and each $f_\alpha$ is $C^k_\C$ and $d^{(j)}f=h_j$
in the real sense, then the complex linearity of $df(x,\cdot)=h_1(x,\cdot)$
implies that $f$ is $C^k_\C$, see Lemma~\ref{realvscxearly}).
We now show by induction on $j\in\N_0$ with $j\leq k$
that $f$ is $C^j$ with $d^{(j)}f=h_j$.
The case $j=0$ holds by definition.
If the assertion holds for $0\leq j<k$,
let $x\in U$ and $y_1,\ldots, y_{j+1}\in E$.
Pick $r>0$ such that $x+[{-r},r]y_{j+1}\sub U$.
By the inductive hypothesis, $f$ is $C^j$ with $d^{\,(j)}f=h_j$ and thus
\[
\gamma(t):=d^{(j)}f(x+ty_{j+1},y_1,\ldots,y_j)=h_j(x+ty_{j+1},y_1,\ldots,y_j)
\]
for $t\in[{-r},r]$.
Define $\gamma_\alpha\colon [{-r},r]\to E$ and $\eta\colon [{-r},r]\to E$
via $\gamma_\alpha(t):=d^{\,(j)}f_\alpha(x+ty_{j+1},y_1,\ldots,y_j)$
and $\eta(t):=h_{j+1}(x+ty_{j+1},y_1,\ldots,y_j,y_{j+1})$.
Since
\[
\gamma_\alpha(t)\to \gamma(t)\quad\mbox{and}\quad \gamma_\alpha'(t)
=d^{\,(j+1)}f_\alpha(x+ty_{j+1},y_1,\ldots, y_{j+1})\to\eta(t)
\]
uniformly in $t\in[{-r},r]$, Exercise~\ref{excunifcnets}(b)
shows that $\gamma$ is a $C^1$-curve with $\gamma'=\eta$.
In particular,
$d^{\,(j+1)}f(x,y_1,\ldots,y_{j+1})=\gamma'(0)$ exists and is given by
\[
d^{(j+1)}f(x,y_1,\ldots, y_{j+1})=\gamma'(0)=\eta(0)=h_{j+1}(x,y_1,\ldots,y_{j+1}),
\]
which completes the proof.
\end{prf}
\begin{prop}\label{compl-and-metr}
Let $E$ be a metrizable locally convex space,
$F$ be a locally convex space, $k\in\N_0\cup\{\infty\}$
and $U\sub E$ be a locally convex subset with dense interior.
Then the following holds:
\begin{description}[(D)]
\item[\rm(a)]
If $F$ is complete, quasi-complete, sequentially complete and Mackey complete,
respectively, then also $C^k(U,F)$ is complete,
quasi-complete, sequentially complete and Mackey complete,
respectively.
\item[\rm(b)]
If $F$ is metrizable and $U$ is locally compact, then $C^k(U,F)$ is metrizable.
\item[\rm(c)]
If $F$ is normable, $U$ is compact and $k\in\N_0$, then
$C^k(U,F)$ is normable.
\end{description}
\end{prop}
\begin{prf}
(a) The topological space $U\times E^j$ is metrizable for each $j\in \N_0$
and hence a $k_\R$-space. The
the locally convex space $\prod_{j\leq k}C(U\times E^j,F)$ (with $j\in\N_0$)
is therefore complete,
quasi-complete, sequentially complete and Mackey complete, respectively,
by Lemma~\ref{sammelsu} (d), (e), (h) and (i), respectively.
Using Lemma~\ref{ima-clo}, the assertions follow.

(b) If $U$ is locally compact, then also $E$ is locally compact and hence
a finite-dimensional vector space. Let $p_1\leq p_2\leq\cdots$ be an ascending
sequence of seminorms on~$F$
defining its locally convex vector topology,
and $(K_n)_{n\in\N}$ be an exhaustion of~$U$ by compact sets (see Exercise~\ref{lcp-hemi}).
Then the seminorms $\|\cdot\|_{C^j,K_n,p_m}$ described
in Exercise~\ref{exc-the-topf} define the locally convex vector topology on $C^k(U,F)$,
for $n,m\in\N$ and $j\in\N_0$ such that $j\leq k$.
As the set of seminorms is countable, $C^k(U,F)$ is metrizable by Corollary~\ref{charmetriz}.

(c) If $p$ is a norm on~$F$ defining its topology, then
$\|\cdot\|_{C^k,U,p}$ defines the topology of $C^k(U,F)$, by the proof of~(b).
\end{prf}
\subsection*{{\boldmath$C^{k,\ell}$}-functions and the exponential law}
\begin{defn}\label{defnCkell}
Let $E_1$, $E_2$, and $F$ be locally convex spaces, $U\sub E_1$ and $V\sub E_2$ be locally convex subsets
with dense interior and $k,\ell\in\N_0\cup\{\infty\}$.
A map $f \colon U\times V\to F$ is called a 
\emph{$C^{k,\ell}$-map}
\index{$C^{k,\ell}$-map}\index{map!$C^{k,\ell}$}
if $d^{\,(0,0)}f:=f$ is continuous, the iterated directional derivatives
\begin{eqnarray*}
\lefteqn{d^{\,(i,j)}f(x,y,v_1,\ldots,v_i,w_1,\ldots, w_j)}\qquad\qquad\qquad \\
&:=& (D_{(v_i,0)}\cdots D_{(v_1,0)}D_{(0,w_j)}\cdots D_{(0,w_1)}f)(x,y)
\end{eqnarray*}
exist for all $i,j\in\N_0$ such that $i\leq k$ and $j\leq\ell$,
$(x,y)\in U^0\times V^0$ and $v_1,\ldots, v_i\in E_1$,
$w_1, \ldots, w_j\in E_2$, and admit continuous extensions
\[
d^{\,(i,j)}f\colon U\times V\times E_1^i\times E_2^j\to F.
\]
We endow the space $C^{k,\ell}(U\times V,F)$ of all $C^{k,\ell}$-maps $f\colon U\times V\to F$
with the so-called 
\index{compact-open $C^{k,\ell}$-topology} 
\emph{compact-open $C^{k,\ell}$-topology}, i.e.,
the initial topology with respect to the linear mappings
\[
d^{\,(i,j)}\colon C^{k,\ell}(U\times V,F)\to C(U\times V\times E_1^i\times E_2^j,F)
\]
for all $i,j$ as before, using the compact-open topology on the right hand side.
\end{defn}
In the preceding situation, we have:
\begin{lem}\label{Ckellpartial}
If $f\colon U\times V\to F$ is a $C^{k,\ell}$-map,
then $f_x:=f(x,\cdot)\colon V\to F$ is a $C^\ell$-map for
each $x\in U$, and
\begin{equation}\label{whatneenow}
d^{\,(j)}(f_x)(y,w_1,\ldots, w_j)=d^{\,(0,j)}f(x,y,w_1,\ldots,w_j)
\end{equation}
for all $j\in\N$ such that $j\leq\ell$ and $(y,w_1,\ldots, w_j)\in V\times E_2^j$.
\end{lem}
\begin{prf}
If $x\in U^0$, then $f_x|_{V^0}$ is $C^\ell$
with
\begin{equation}\label{babyargu}
d^{\,(j)}(f_x|_{V^0})(y,w)=d^{\,(0,j)}f(x,y,w)
\end{equation}
for all
$j\in\N$ such that $j\leq \ell$, $y\in V^0$ and $w=(w_1,\ldots, w_j)\in E_2^j$.
As the right hand side of (\ref{babyargu})
defines a continuous $F$-valued function of $(y,w)\in V\times E_2^j$,
we see that $f_x$ is $C^\ell$ and (\ref{whatneenow}) holds
for $j$ as before, for each $x\in U^0$.

If $x\in U\setminus U^0$, we can apply Lemma~\ref{Ckinlimit}
to $f|_{X\times V}$ with $X:=U^0\cup\{x\}$ and $g_j:=d^{\,(0,j)}f|_{X\times V\times E_2^j}$
to see that again $f_x$ is $C^k$ and (\ref{whatneenow}) holds.\footnote{Note that $x$ is an accumulation point of $X$ since $U^0$ is dense in~$X$ and $x\not\in U^0$.}
\end{prf}
\begin{rem}\label{remCkellobs}
(a)
In particular, the $C^{0,\ell}$-maps on $U\times V\sub E_1\times E_2$
in the sense of Definition~\ref{defnC0k} coincide with those of Definition~\ref{defnCkell}
(and likewise for $C^{k,0}$-maps).\medskip

\noindent(b) For $f\in C^{k,\ell}(U\times V,F)$,
the same reasoning as in Lemma~\ref{Ckellpartial} shows that $U\to F$, $x\mto d^{\,(0,j)}f(x,y,w_1,\ldots, w_j)$
is a $C^k$-map for all $j\in\N_0$ with $j\leq\ell$, $y\in V$ and $w_1,\ldots, w_j\in E_2$.
\end{rem}
Mappings to products behave as expected.
\begin{prop}\label{Cklprod}
Let $E_1$ and $E_2$ be locally convex spaces, $U\sub E_1$ and $V\sub E_2$
be locally convex subsets with dense interior, and $k,\ell\in\N_0\cup\{\infty\}$.
Let $(F_a)_{a\in A}$ be a family of locally convex spaces.
Then a map
\[
f=(f_a)_{a\in A}\colon U\times V\to \prod_{a\in A} F_a\vspace{-.3mm}
\]
is $C^{k,\ell}$
if and only if each component $f_a\colon U\times V\to F_a$ is $C^{k,\ell}$.
In this case,
\begin{equation}\label{thederprod}
d^{\,(i,j)}f=(d^{(i,j)}f_a)_{a\in A}\;\,
\mbox{for all $i,j\in \N_0$ such that $i\leq k$ and $j\leq \ell$.}
\end{equation}
\end{prop}
\begin{prf}
Let $F:=\prod_{a\in A}F_a$.
If $f$ is $C^{k,\ell}$, then each component $f_a$ is continuous.
Moreover, the iterated directional derivative
\[
D_{(v_1,0)}\cdots D_{(v_1,0)}D_{(0,w_j)}\cdots D_{(0,w_1)}f_a(x,y)
\]
exists for all $i,j\in \N_0$ such that $i\leq k$ and $j\leq \ell$,
all $(x,y)\in U^0\times V^0$\linebreak
and all $v=(v_1,\ldots, v_i)\in E_1^i$, $w=(w_1,\ldots, w_j)\in E_2^j$,
and it equals
$\pr_a(d^{\,(i,j)}f(x,y,v,w))$,
as we can apply the continuous linear projection\linebreak
$\pr_a\colon F\to F_a$
to the corresponding derivatives of $f$. Since the mapping $\pr_a\circ \,d^{\,(i,j)}f\colon U\times V\times E_1^i\times E_2^j\to F$
is continuous, it serves as the continuous extension $d^{\,(i,j)}f_a$
and we deduce that $f_a$ is $C^{k,\ell}$.
If, conversely,  each $f_a$ is $C^{k,\ell}$, then $f$ is continuous.
Moreover,
for all $i,j\in \N_0$ such that $i\leq k$ and $j\leq \ell$,
all $(x,y)\in U^0\times V^0$ and all $v=(v_1,\ldots, v_i)\in E_1^i$, $w=(w_1,\ldots, w_j)\in E_2^j$,
\[
D_{(v_1,0)}\cdots D_{(v_1,0)}D_{(0,w_j)}\cdots D_{(0,w_1)}f(x,y)
\]
exists and equals $(d^{\,(i,j)}f_a(x,y,v,w))_{a\in A}$,
as limits can be formed component\-wise.
Since $(d^{\,(i,j)}f_a)_{a\in A}\colon U\times V\times E_1^i\times E_2^j\to F$ is continuous,
we see that $f$ is $C^{k,\ell}$ and (\ref{thederprod}) holds.
\end{prf}
We record a version of
Schwarz' Theorem for $C^{k,\ell}$-maps.
\begin{prop}\label{schwarzCkell}
Let $E_1$, $E_2$, and $F$ be locally convex spaces, $U\sub E_1$ and $V\sub E_2$
be open subsets, $k,\ell\in\N_0\cup\{\infty\}$ and $f\colon U\times V\to F$ be a $C^{k,\ell}$-map.
Let $i,j\in\N_0$ with $i\leq k$ and $j\leq\ell$.
For $v_1,\ldots, v_i\in E_1$, $w_1,\ldots, w_j\in E_2$,
define
\[
u_a:=\left\{
\begin{array}{cl}
(v_a,0) & \mbox{if $a\in\{1,\ldots,i\}$;}\\
(0,w_{a-i}) & \mbox{if $a\in \{i+1,\ldots, i+j\}$.}
\end{array}
\right.
\]
Let $\pi\in S_{i+j}$ be a permutation of $\{1,\ldots, i+j\}$.
Then the iterated directional derivative
\[
(D_{u_{\pi(1)}}\cdots D_{u_{\pi(i+j)}}f)(x,y)
\]
exists for all $(x,y)\in U\times V$, and it equals $d^{\,(i,j)}f(x,y,v_1,\ldots,v_i,
w_1,\ldots, w_j)$.
\end{prop}
\begin{prf}
Since $f(x,\cdot )$ is $C^\ell$ and $f(\cdot ,y)$ is $C^k$
(whence Proposition~\ref{schwarz} applies to these maps), we may assume that both $i,j\geq 1$
(and hence $k,\ell\geq 1$). Let $j=1$ first;
we show by induction on $i\in \N$ with $i\leq k$ that
\[
(D_{(0,w_1)}D_{(v_i,0)}\cdots D_{(v_1,0)}f)(x,y)
\]
exists for all $(x,y)\in U\times V$ and equals $d^{\,(i,1)}f(x,y,v_1,\ldots, v_i,w_1)$.
For $i=1$, this follows from Lemma~\ref{stro2nd}, applied to $(s,t)\mto f(x+sv_1,y+tw_1)$
with $(s,t)$ in a $0$-neighborhood in~$\K^2$.
Now let $i\geq 2$ and assume that the assertion holds for $i-1$ in place of~$i$;
thus the following derivatives exist and coincide:
\begin{equation}\label{thusnw}
\hspace*{-.5mm}D_{(0,w_1)}D_{(v_{i-1},0)}\cdots D_{(v_1,0)}f(x,y)=
D_{(v_{i-1},0)}\cdots D_{(v_1,0)}D_{(0,w_1)}f(x,y).\!\!
\end{equation}
Then
\[
g\colon U\times V\to F,\quad (x,y)\mto D_{(v_{i-1},0)}\cdots D_{(v_1,0)}f(x,y)
\]
is $C^{(1,0)}$ with $d^{(1,0)}g(x,y,v_i)=d^{(i,0)}f(x,y,v_1,\ldots, v_i)$.
By (\ref{thusnw}),
\[
D_{(0,w_1)}g(x,y)
\]
exists and equals $d^{\,(i-1,1)}f(x,y,v_1,\ldots, v_{i-1},w_1)$;
it is therefore continuous in $(x,y,w_1)\in U\times V\times E_2$.
The last identity entails that
\[
D_{(v_i,0)}D_{(0,w_1)}g(x,y)
\]
exists and equals $d^{\,(i,1)}f(x,y,v_1,\ldots, v_i,w_1)$,
which is a continuous $F$-valued function of $(x,y,v_i,w_1)\in U\times V\times E_1\times E_2$.
Hence $(s,t)\mto g(x+sv_i,y+tw_1)$
satisfies the hypotheses of Lemma~\ref{stro2nd} and we deduce that also
\[
D_{(0,w_1)}D_{(v_i,0)}g(x,y)
= D_{(0,w_1)}D_{(v_i,0)}\cdots D_{(v_1,0)}f(x,y)
\]
exists and equals
\begin{eqnarray*}
D_{(v_i,0)}D_{(0,w_1)}g(x,y)&=& D_{(v_i,0)}D_{(0,w_1)}D_{(v_{i-1},0)}\cdots
D_{(v_1,0)}f(x,y)\\
&=&D_{(v_i,0)}\cdots D_{(v_1,0)}D_{(0,w_1)}f(x,y),
\end{eqnarray*}
using (\ref{thusnw}) to pass to the second line. This completes the induction step.

Now let $\pi\in S_{i+j}$ and $u_1,\ldots, u_{i+j}$ be as in the proposition.
The proof is by induction on $i+j$ with $i,j\geq 1$.
The case $i+j=2$ has already been settled. Now assume $i+j\geq 3$.

If $\pi(1)\in\{1,\ldots,i\}$,
then $(D_{u_{\pi(2)}}\cdots D_{u_{\pi(i+j)}}f)(x,y)$ exists by induction and equals
\[
d^{\,(i-1,j)}f(x,y,v_1,\ldots,v_{\pi(1)-1},v_{\pi(1)+1},\ldots,
v_i,w_1,\ldots,w_j).
\]
This function of~$(x,y)$
can be differentiated in the direction $u_{\pi_1}=(v_{\pi(1)},0)$. Therefore the following derivatives exist and coincide:
\begin{eqnarray*}
\lefteqn{(D_{u_{\pi(1)}}\cdots D_{u_{\pi(i+j)}}f)(x,y)}\qquad\\
&=& d^{\,(i,j)}f(x,y,v_1,\ldots,v_{\pi(1)-1},v_{\pi(1)+1},\ldots,
v_i,v_{\pi(1)},w_1,\ldots,w_j)\\
&=&
d^{\,(i,j)}f(x,y,v_1,\ldots,v_i,w_1,\ldots,w_j);
\end{eqnarray*}
the final equality holds by Proposition~\ref{schwarz},
as $x\mto d^{\,(0,j)}f(x,y,w_1,\ldots, w_j)$ is~$C^k$.

If $\pi(1)\in\{i+1,\ldots, i+j\}$,
then the following derivatives exist by induction and coincide:
\begin{eqnarray*}
\lefteqn{(D_{u_{\pi(2)}}\cdots D_{u_{\pi(i+j)}}f)(x,y)}\qquad\\
&=&d^{\,(i,j-1)}f(x,y,v_1,\ldots,v_i,w_1,\ldots,w_{\pi(1)-1},
w_{\pi(1)+1},\ldots, w_j).
\end{eqnarray*}
Now $h\colon U\times V\to F$, $(x,y)\mto d^{\,(0,j-1)}f(x,y,w_1,\ldots,w_{\pi(1)-1},w_{\pi(1)+1},\ldots, w_j)$
is a $C^{k,1}$-function with
\begin{eqnarray*}
\lefteqn{d^{(i,1)}h(x,y,v_1,\ldots, v_i,w_{\pi(1)})}\qquad\\
&=& d^{\,(i,j)}f(x,y,v_1,\ldots,v_i,w_1,\ldots,w_{\pi(1)-1},w_{\pi(1)+1},\ldots, w_j,w_{\pi(1)})\\
&=&
d^{\,(i,j)}f(x,y,v_1,\ldots, v_i,w_1,\ldots, w_j),
\end{eqnarray*}
using that $f(x,\cdot )$ is $C^\ell$ to reorder the differentiations in the second
variable (with Proposition~\ref{schwarz}).
By the case $j=1$ treated above,
\[
D_{u_{\pi(1)}}\cdots D_{u_{\pi(i+j)}}f(x,y)=
D_{w_{\pi(1)}}D_{(v_i,0)}\cdots D_{(v_1,0)}h(x,y)
\]
exists and equals 
\[
d^{(i,1)}h(x,y,v_1,\ldots, v_i,w_{\pi(1)})
=d^{\,(i,j)}f(x,y,v_1,\ldots, v_i,w_1,\ldots, w_j).\qedhere
\]
\end{prf}
\begin{cor}\label{CkellCellk}
Let $E_1$, $E_2$, and $F$ be locally convex spaces,
$U\sub E_1$ and $V\sub E_2$ be locally convex subsets with dense interior,
$k,\ell\in \N_0\cup\{\infty\}$ and
$f\colon U\times V\to F$ be a $C^{k,\ell}$-map.
Then
\[
g\colon V\times U\to F,\quad g(y,x):=f(x,y)
\]
is a $C^{\ell,k}$-map, and
\begin{equation}\label{mixedinter}
d^{(j,i)}g(y,x,w_1,\ldots, w_j,v_1,\ldots, v_i)\!=\!
d^{(i,j)}f(x,y,v_1,\ldots, v_i,w_1,\ldots, w_j)
\end{equation}
for all
$i,j\in\N_0$ with $i\leq k$ and $j\leq \ell$, $x\in U$, $y\in V$,
$v=(v_1,\ldots, v_i)\in E_1^i$ and $w=(w_1,\ldots, w_j)\in E_2^j$.
\end{cor}
\begin{prf}
By Proposition~\ref{schwarzCkell},
the directional derivative $d^{\,(j,i)}g(y,x,w,v)$
exists if $(x,y)\in V^0\times U^0$ and is given by~(\ref{mixedinter}).
Since the right hand side of (\ref{mixedinter})
provides a continuous extension to all $(y,x,w,v)\in V\times U\times E_2^j\times E_1^i$,
we see that $g$ is $C^{\ell,k}$ and (\ref{mixedinter}) holds.
\end{prf}
\begin{lem}\label{linCkell}
Let $E_1$, $E_2$, and $F$ be locally convex spaces, $U\sub E_1$
be a locally convex subset with dense interior, $k\in\N_0\cup\{\infty\}$
and $f\colon U\times E_2\to F$ be a $C^{k,0}$-map
such that $f_x:=f(x,\cdot )\colon E_2\to F$ is linear for each
$x\in U$. Then $f$ is $C^{k,\infty}$.
\end{lem}
\begin{prf}
By definition, $i$-fold directional derivatives in the first variable
exist at $(x,y)\in U^0\times E_2$
and a continuous extension $d^{\,(i,0)}f$,
for all $i\in\N$ such that $i\leq k$.
Using that $f_x$ is continuous linear, we find that
\begin{equation}\label{reftoths}
d^{\,(0,1)}f(x,y,w_1)=f(x,w_1)
\end{equation}
for all $(x,y,w_1)\in U^0\times E_2\times E_2$ and (as (\ref{reftoths}) is independent of~$y$)
\[
d^{\,(0,j)}f(x,y,w_1,\ldots, w_j)=0\;\,\mbox{for all $(x,y,w_1,\ldots, w_j)\in U^0\times E_2^{j+1}$,}
\]
for all $j\in \N$ with $j\geq 2$.
In both cases, we can form directional derivatives in the first variable;
for $i$ as before and $v_1,\ldots, v_i\in E_1$,
we obtain
\[
d^{\,(i,1)}f(x,y,v_1,\ldots, v_i,w_1)=d^{\,(i,0)}f(x,w_1,v_1,\ldots, v_i)\quad\mbox{and}\vspace{-1mm}
\]
\[
d^{\,(i,j)}f(x,y,v_1,\ldots,v_i,w_1,\ldots, w_j)= 0.
\]
As the right hand sides of the last and penultimate equation define continuous
$F$-valued functions on $U\times E_2\times E_1^i\times E_2$ and
$U\times E_2\times E_1^i\times E_2^j$, respectively,
we see that $f$ is $C^{k,\infty}$.
\end{prf}
\begin{rem}\label{linCkellvar}
Let $E_1$, $E_2$, and $F$ be locally convex spaces, $V\sub E_2$ be a locally
convex subset with dense interior, $\ell\in \N_0\cup\{\infty\}$ and $g\colon E_1\times V\to F$
be a $C^{0,\ell}$-map such that $g^y:=g(\cdot ,y)\colon E_1\to F$ is linear
for each $y\in V$. Then $g$ is $C^{\infty,\ell}$
(as we can combine Corollary~\ref{CkellCellk} and Lemma~\ref{linCkell}).
\end{rem}
\begin{lem}\label{evaldiffprop}
Let $E$ be a finite-dimensional vector space,
$F$ be a locally convex space, $\ell\in \N_0\cup\{\infty\}$
and $U\sub E$ be a locally convex subset with dense interior.
If $U$ is locally compact, then the evaluation map
\[
\ve\colon C^\ell(U,F)\times U\to F,\quad (f,x)\mto f(x)
\]
is a $C^{\infty,\ell}$-map.
\end{lem}
\begin{prf}
As the inclusion map $C^\ell(U,F)\to C(U,F)$ is continuous,
we deduce from Lemma~\ref{laeval} that $\ve$ is continuous and hence~$C^{0,0}$.
Now assume $\ell\geq 1$. For fixed $f\in C^\ell(U,F)$, we have
\[
\ve(f,\cdot )=f.
\]
It is clear from this that the derivatives
$d^{\,(0,j)}\ve(f,x,y_1,\ldots,y_j)$ exist
for all $j\in\N$ with $j\leq \ell$, $x\in U^0$ and $y_1,\ldots, y_j\in E$,
and are given by
\[
d^{\,(0,j)}\ve(f,x,y_1,\ldots, y_j)=d^{\,(j)}f(x,y_1,\ldots,y_j).
\]
An extension to $(f,x,y_1,\ldots,y_j)\in C^\ell(U,F)\times U\times E^j$
is given by the
function
\[
h_j\colon C^\ell(U,E)\times U\times E^j\to F,\quad
(f,x,y_1,\ldots, y_j)\mto \ve_j(d^{\,(j)}f,(x,y_1,\ldots,y_j)),
\]
where $\ve_j\colon C(U\times E^j,F)\times (U\times E^j)\to F$ is the evaluation map.
Now $\ve_j$
is continuous by Lemma~\ref{laeval} and the map
\[
d^{\,(j)}\colon C^\ell(U,F)\to C(U\times E^j,F),\quad f\mto d^{\,(j)}f
\]
is continuous by definition of the compact-open $C^\ell$-topology.
Hence $h_j$ is continuous and we deduce that $\ve$ is $C^{0,\ell}$
with
\[
d^{\,(0,j)}\ve(f,x,y_1,\ldots, y_j)=d^{\,(j)}f(x,y_1,\ldots, y_j)
\]
for all $j\in \N$ such that $j\leq\ell$, $x\in U$ and $y_1,\ldots, y_j\in E$.
Since $\ve(f,x)$ is linear in $f$,
we deduce that~$\ve$ is $C^{\infty,\ell}$ (see Remark~\ref{linCkellvar}).
\end{prf}
\begin{prop}[Chain Rule 1 for {\boldmath$C^{r,s}$}-maps]\label{chainR1}
Let $E_1$, $E_2$, $F_1$, $F_2$, and $F$ be locally convex spaces,
$U_1\sub E_1$, $U_2\sub E_2$, $V_1\sub F_1$, and $V_2\sub F_2$
be locally convex subsets with dense interior, $r,s\in\N_0\cup\{\infty\}$,
$g\colon V_1\times V_2\to F$ be a $C^{r,s}$-map,
$f_1\colon U_1\to F_1$ be a $C^r$-map with $f_1(U_1)\sub V_1$
and $f_2\colon U_2\to F_2$ be a $C^s$-map with $f_2(U_2)\sub V_2$.
Then also the following map is $C^{r,s}$:
\[
g\circ (f_1\times f_2)\colon U_1\times U_2\to F,\quad (x,y)\mto g(f_1(x),f_2(y)).
\]
\end{prop}
\begin{prf}
Let $h:=g\circ(f_1\times f_2)$;
assume first that~$U_1$ and~$U_2$ are open.
Let $k,\ell\in\N_0$ with $k\leq r$ and $\ell\leq s$.
We show that the iterated directional derivatives
\[
d^{\,(k,\ell)}h(x,y,v,w)
\]
exist for all $x\in U_1$, $y\in U_2$, $v=(v_1,\ldots, v_k)\in E_1^k$
and $w=(w_1,\ldots, w_\ell)\in E_2^\ell$.
Moreover, we shall obtain explicit formulas for $d^{\,(k,\ell)}h(x,y,v,w)$,
which show that it is a continuous $F$-valued function of $(x,y,v,w)$.

If $\ell\geq 1$, then $g(f_1(x),\cdot)\colon V_2\to F$, $z\mto g(f_1(x),z)$
is a $C^s$-map by Lemma~\ref{Ckellpartial}, with $j$th differential $d^{\,(0,j)}g(f_1(x),\cdot)$
for $j\in\N$ with $j\leq s$.
Hence $U_2\to F$, $y\mto g(f_1(x),f_2(y))$ is a $C^s$-map, by the Chain Rule (Proposition~\ref{chainno}),
with $\ell$th differential taking $(y,w)\in U_2\times E_2^\ell$ to
\begin{equation}\label{FDB1}
\sum_{j=1}^\ell\sum_{P\in P_{\ell,j}}d^{\,(0,j)}g(f_1(x),f_2(y),d^{\,(|I_1|)}f_2(y,w_{I_1}),\ldots,
d^{\,(|I_j|)}f_2(y,w_{I_j})),
\end{equation}
by Fa\`{a} di Bruno's Formula (Theorem~\ref{faadk}
and Remark~\ref{nopecarryov}(c)),
using the notation introduced there, with $P=\{I_1,\ldots, I_j\}$.
Thus $d^{\,(0,\ell)}h(x,y,w)$ is given by (\ref{FDB1}),
and so $d^{\,(0,\ell)}h$ is continuous.

Likewise, if $k\geq 1$ and $\ell=0$, we find that $d^{\,(k,0)}h(x,y,v)$ exists and equals
\begin{equation}\label{FDB2}
\sum_{i=1}^k\sum_{Q\in P_{k,i}}d^{\,(i,0)}g(f_1(x),f_2(y),d^{\,(|J_1|)}f_1(x,v_{J_1}),\ldots,
d^{\,(|J_i|)}f_1(x,v_{J_i}))
\end{equation}
with $Q=\{J_1,\ldots, J_i\}$. Thus $d^{\,(k,0)}h$ is continuous.

If $k,\ell\geq 1$, we use that $V_1\to F$, $z\mto
d^{\,(0,j)}g(z,u)$ is a $C^r$-map for $u\in V_2\times F_2^j$ given by
\[
u=(u_0,u_1,\ldots, u_j):=(f_2(y),d^{\,(|I_1|)}f_2(y,w_{I_1}),\ldots,
d^{\,(|I_j|)}f_2(y,w_{I_j}))
\]
(with notation as in \hspace*{-.4mm}(\ref{FDB1})),
whose $i$th differential takes $(z_0,z_1,\ldots, z_i)\!\in\! V_1\!\times \! F_1^i$\linebreak
to
$d^{\,(i, j)}g(z_0,u_0,z_1,\ldots,z_i,u_1,\ldots,u_j)$ (see Remark~\ref{remCkellobs}(b)).\vspace{.2mm}
Hence $U_1\to F$, $x\mto d^{\,(0,\ell)}h(x,y,w)$ is a $C^r$-map, by the Chain Rule (Proposition~\ref{chainno}),
with $k$th differential taking $(x,v)\in U_1\times E_1^k$ to
\[
\sum_{i=1}^k\sum_{j=1}^\ell \sum_{Q\in P_{k,i}}\sum_{P\in P_{\ell,j}}
d^{\,(i,j)}g\big(f_1(x),f_2(y),
d^{\,(|J_1|}f_1(x,v_{J_1}),\ldots,d^{\,(|J_i|)}f_1(x,v_{J_i}),\vspace{-1mm}
\]
\begin{equation}\label{FDB3}
d^{\,(|I_1|)}f_2(y,w_{I_1}),\ldots,d^{\,(|I_j|)}f_2(y,w_{I_j})\big).
\end{equation}
So $d^{\,(k,\ell)}h(x,y,v,w)$ exists and is given by (\ref{FDB3});
thus $d^{\,(k,\ell)}h$ is continuous.

In the case of general $U_1$ and $U_2$, for all $k,\ell\in \N_0$ with $k\leq r$ and $\ell\leq s$,
the iterated directional derivative
$d^{\,(k,\ell)}h(x,y,v,w)$ exists and is given by (\ref{FDB1}), (\ref{FDB2}) and (\ref{FDB3}),
respectively,
for all $(x,y,v,w)\in U_1^0\times U_2^0\times E_1^k\times E_2^\ell$.
Since (\ref{FDB1}), (\ref{FDB2}) and (\ref{FDB3}) define continuous maps
\[
d^{\,(k,\ell)}h\colon U_1\times U_2\times E_1^k\times E_2^\ell\to F
\]
that extend $d^{\,(k,\ell)}(h|_{U_1^0\times U_2^0})$, we see that~$h$ is~$C^{r,s}$.
\end{prf}
\begin{rem}\label{faa-chain-1}
The expressions for $d^{\,(k,\ell)}h(x,y,v,w)$ provided by (\ref{FDB1}), (\ref{FDB2}), and (\ref{FDB3})
are the analogs of Fa\`{a} di Bruno's Formula
for the differentials $d^{\,(k,\ell)}h$ of $h:=g\circ (f_1\times f_2)$
in the cases $k=0$, resp., $\ell=0$, resp., $k,\ell\not=0$.
\end{rem}
\begin{prop}[Chain Rule 2 for {\boldmath$C^{r,s}$}-maps]\label{chainR2}
Let $E_1$, $E_2$, $F$, and $H$ be locally convex spaces, $U\sub E_1$, $V\sub E_2$,
and $W\sub F$ be locally convex subsets with dense interior, $r,s\in\N_0\cup\{\infty\}$,
$f\colon U\times V\to F$ be a $C^{r,s}$-map such that $f(U\times V)\sub W$,
and $g\colon W\to H$ be a $C^{r+s}$-map.
Then also $g\circ f\colon U\times V\to H$ is a $C^{r,s}$-map.
\end{prop}
\begin{prf}
Let $h:=g\circ f$.
For $k,\ell\in\N_0$ with $k\leq r$ and $\ell\leq s$,
we show that
\[
d^{\,(k,\ell)}h(x,y,v,w)
\]
exists for all $x\in U^0\!$, $y\in V^0\!$, $v=(v_1,\ldots, v_k)\in E_1^k$
and $w\!=\!(w_1,\ldots, w_\ell)\!\in\! E_2^\ell$.

If $\ell\geq 1$, then $f(x,\cdot )\colon V^0 \to F$, $y\mto f(x,y)$
is a $C^s$-map by Lemma~\ref{Ckellpartial}, with $j$th differential $d^{\,(0,j)}f(x,\cdot )$
for $j\in\N$ with $j\leq s$.
Hence $V^0\to F$, $y\mto g(f(x,y))$ is a $C^s$-map, by the Chain Rule (Proposition~\ref{chainno}),
with $\ell$th differential taking $(y,w)\in V^0\times E_2^\ell$ to
\begin{equation}\label{FDB4}
\sum_{j=1}^\ell\sum_{P\in P_{\ell,j}}d^{\,(j)}g(f(x,y),d^{\,(0,|I_1|)}f(x,y,w_{I_1}),\ldots,
d^{\,(0,|I_j|)}f(x,y,w_{I_j})),
\end{equation}
by Fa\`{a} di Bruno's Formula (Theorem~\ref{faadk}
and Remark~\ref{nopecarryov}(c))
using the notation introduced there, with $P=\{I_1,\ldots, I_j\}$.
Thus $d^{\,(0,\ell)}h(x,y,w)$ is given by (\ref{FDB4}).

Likewise, if $k\geq 1$ and $\ell=0$, we find that $d^{\,(k,0)}h(x,y,v)$ exists and equals
\begin{equation}\label{FDB5}
\sum_{i=1}^k\sum_{Q\in P_{k,i}}d^{\,(i)}g(f(x,y),d^{\,(|J_1|,0)}f(x,y,v_{J_1}),\ldots,
d^{\,(|J_i|,0)}f(x,y,v_{J_i}))
\end{equation}
with $Q=\{J_1,\ldots, J_i\}$.

If $k,\ell\geq 1$, we use that $\theta_P^{y,w} \colon U\to W\times F^j$,
\[
x \mto \big(f(x,y),d^{\,(0,|I_1|)}f(x,y,w_{I_1}),\ldots, d^{\,(0,|I_j|)}f(x,y,w_{I_j})\big)
\]
is a $C^r$-map (see Remark~\ref{remCkellobs}(b)), with
\begin{eqnarray}
d^{\,(m)}\theta_P^{y,w}(x,u)
&=&
\big(d^{\,(m,0)}f(x,y,u),d^{\,(m,|I_1|)}f(x,y,u,w_{I_1}),\qquad \notag\\
& & \qquad\qquad\qquad\qquad  \ldots, d^{\,(m,|I_j|)}f(x,y,u,w_{I_j})\big)\label{preFDB6}
\end{eqnarray}
for $x\in U$, $m\in\N$ with $m\leq r$, and $u\in E_1^m$.
Hence $U^0\to F$,
\[
x\mto d^{\,(0,\ell)}h(x,y,w)=\sum_{j=1}^\ell \sum_{P\in P_{\ell,j}}
d^{\,(j)}g(\theta_P^{y,w}(x))
\]
is a $C^r$-map, by the Chain Rule (Proposition~\ref{chainno}),
with $k$th differential taking $(x,v)\in U^0\times E_1^k$ to
\begin{equation}\label{FDB6}
\sum_{i=1}^k\sum_{j=1}^\ell \sum_{Q\in P_{k,i}}\!\sum_{P\in P_{\ell,j}}
\!\!\!d^{\,(i+j)}g\big(\theta_P^{y,w}(x),
d^{\,(|J_1|)}\theta_P^{y,w}(x,v_{J_1}),\ldots,d^{\,(|J_i|)}\theta_P^{y,w}(x,v_{J_i})\big).
\end{equation}
Thus $d^{\,(k,\ell)}h(x,y,v,w)$ exists for all $(x,y,v,w)\in U^0\times V^0\times E_1^k\times E_2^\ell$
and is given by (\ref{FDB6}).
Now (\ref{preFDB6}) shows that $d^{(|J_a|)}\theta_P^{y,w}(x,v_{J_a})$
is a continuous function
of $(x,y,v,w)\in U\times V\times E_1^k\times E_2^\ell$, for all $a\in\{1,\ldots, i\}$.
Therefore (\ref{FDB6}) defines a continuous function
$d^{\,(k,\ell)}h$ of $(x,y,v,w)\in U_1\times U_2\times E_1^k\times E_2^\ell$
which extends $d^{\,(k,\ell)}(h|_{U_1^0\times U_2^0})$.
Likewise, (\ref{FDB4}) and (\ref{FDB5}) provide\vspace{-.3mm} continuous extensions $d^{\,(k,\ell)}h$ if $k=0$ and $\ell=0$, respectively.
Thus $h$ is $C^{r,s}$.
\end{prf}
\begin{rem}
(a) The expressions for $d^{\,(k,\ell)}(g\circ f)(x,y,v,w)$ given by (\ref{FDB4}), (\ref{FDB5}), and (\ref{FDB6})
(together with (\ref{preFDB6}))
are the analogs of Fa\`{a} di Bruno's Formula
for the differentials $d^{\,(k,\ell)}(g\circ f)$
in the cases $k=0$, resp., $\ell=0$, resp., $k,\ell\not=0$.\medskip

(b) If $\ell=0$, then $U_2$ and $V_2$ (resp., $V$) can be replaced with Hausdorff topological spaces
in Propositions \ref{chainR1} and \ref{chainR2}, respectively.
Likewise, if $k=0$, then~$U_1$ and~$V_1$ (resp., $U$) can be replaced with Hausdorff
topological spaces in the cited propositions.
\end{rem}
On a first reading, one may concentrate on case~(a)
of the next theorem. See \ref{defnkkR}(c) in the appendix for the notion of a $k_\R$-space.
\begin{thm}[Exponential law for spaces of {\boldmath$C^{k,\ell}$}-functions]\label{explawCkell}
Let $E_1$, $E_2$, and $F$ be locally convex spaces,
$U\sub E_1$ and $V\sub E_2$ be locally convex subsets with dense interior
and $k,\ell\in\N_0\cup\{\infty\}$.
If $f \in C^{k,\ell}(U\times V,F)$, set $f^\vee(x):=f_x:=f(x,\cdot )$
for $x\in U$. Then
\[
f^\vee\colon U\to C^\ell(V,F)
\]
is a $C^k$-map and the mapping
\[
\Phi\colon C^{k,\ell}(U\times V,F)\to C^k(U,C^\ell(V,F)),\quad f \mto f^\vee
\]
is linear and a topological embedding. If
\begin{description}
\item[\rm(a)]
$V$ is locally compact; or
\item[\rm(b)]
$U\times V\times E_1^i \times E_2^j$ is a $k_\R$-space,
for all $i,j\in \N_0$ with $i\leq k$ and $j\leq \ell$,
\end{description}
then $\Phi$ is an isomorphism of topological vector spaces.
\end{thm}
\begin{rem}\label{inwashexp}
(a) If $k=0$, then $U$ may be replaced with a Hausdorff topological space
in the setting of Theorem~\ref{explawCkell} (ignoring the symbol $E_1^0$).
Likewise, $V$ may be replaced with a Hausdorff
topological space if $\ell=0$.\medskip

\noindent
(b) Let $g\colon U\to C^\ell(V,F)$ be a map.
If $\Phi$ is an isomorphism, then $g$ is $C^k$ if and only if
\[
g^\wedge \colon U\times V\to F,\quad g^\wedge(x,y):=g(x)(y)
\]
is $C^{k,\ell}$. [If $g$ is $C^k$, then $g^\wedge=\Phi^{-1}(g)\in C^{k,\ell}(U\times V,F)$.
If $g^\wedge\in C^{k,\ell}(U\times V,F)$, then $g=\Phi(g^\wedge)\in C^k(U,C^\ell(V,F))$.]\medskip

\noindent(c)
Let $f\colon U\times V\to F$ be a $C^{k,\ell}$-map.
We shall see in the proof of Theorem~\ref{explawCkell} that
\begin{equation}\label{mightbeuse}
(d^{\,(j)}\circ (d^{\,(i)}f^\vee))(x,v)(y,w)=d^{\,(i,j)}f(x,y,v,w)
\end{equation}
for all $i,j\in\N_0$ with $i\leq k$ and $j\leq \ell$, and all
$(x,y,v,w)\in U\times V\times E_1^i\times E_2^j$,
\end{rem}
\noindent
\emph{Proof of Theorem}~\ref{explawCkell}.
We know from Lemma~\ref{Ckellpartial} that $f^\vee(x)=f_x\in C^\ell(V,F)$
for each $x\in U$.

\emph{The map $f^\vee\colon U\to C^\ell(V,F)$ is continuous.}
This will hold if
\[
d^{\,(j)}\circ f^\vee\colon U\to C(V\times E_2^j,F)
\]
is continuous for each
$j\in \N_0$ with $j \leq \ell$, as the topology on $C^\ell(V,F)$ is initial with respect
to the maps $d^{\,(j)}\colon C^\ell(V,F)\to C(V\times E_2^j,F)$
for~$j$ as before.
If $x\in U$, $y\in V$ and $w\in E_2^j$, we have
\[
(d^{\,(j)}\circ f^\vee)(x)(y,w) = d^{\,(j)}(f_x)(y,w)=d^{\,(0,j)}f(x,y,w)
=(d^{\,(0,j)}f)^\vee(x)(y,w).
\]
Hence
\[
d^{\,(j)}\circ f^\vee = (d^{\,(0,j)}f)^\vee,
\]
if we consider $d^{\,(0,j)}f$ as a continuous function $U\times (V\times E_2^j)\to F$.
By Proposition~\ref{ctsexp}(a), the map $(d^{\,(0,j)}f)^\vee\colon U\to C(V\times E_2^j,F)$
(and hence also $d^{\,(j)}\circ f^\vee$) is continuous.

\emph{The map $f^\vee\colon U\to C^\ell(V,F)$ is $C^k$.}
To see this, we show by induction that $f^\vee\colon U\to C^\ell(V,F)$ is $C^i$
for each $i\in \N$ with $i\leq k$, and
\begin{equation}\label{theformul}
d^{\,(i)}(f^\vee)(x,v)= (d^{\,(i,0)}f)(x,\cdot,v)
\end{equation}
for all $x\in U$ and $v\in E_1^i$.
Using
Corollary~\ref{CkellCellk} twice, we see that the map
\[
h_i\colon (U\times E_1^i)\times V\to F,\quad ((x,v),y)\mto d^{\,(i,0)}f(x,y,v)
\]
is $C^{0,\ell}$, with
\[
d^{\,(0,j)}h_i((x,v),y,w)=d^{\,(i,j)}f(x,y,v,w)
\]
for all $(x,v,y,w)\in U\times E_1^i \times V\times E_2^j$.
Hence
\[
h_i^\vee\colon U\times E_1^i\to C^\ell(V,F),\quad (x,v)\mto (d^{\,(i,0)}f)(x,\cdot ,v)
\]
is continuous, by the case $i=0$ treated above, with
\begin{equation}\label{notherll}
(d^{\,(j)}\circ (h_i^\vee))(x,v)(y,w)=d^{\,(i,j)}f(x,y,v,w).
\end{equation}
It therefore suffices to show that the iterated directional
derivative on the left hand side of (\ref{theformul}) exists and
coincides with $h_i^\vee(x,v)$, for all $x\in U^0$ and $v\in E_1^i$;
then $f^\vee$ will be~$C^i$ and\footnote{Note that (\ref{notherll})
and (\ref{latertopo}) imply (\ref{mightbeuse}).}
\begin{equation}\label{latertopo}
d^{\,(i)}(f^\vee)=h_i^\vee.
\end{equation}
If $i=1$, we have $x+\bD_\ve  v\sub U^0$ for some $\ve>0$.
For $t\in \bD_\ve\setminus\{0\}$,
consider
\[
\Delta_t:=
\frac{f^\vee(x+tv)-f^\vee(x)}{t} \in C^\ell(V,F).
\]
We show that
\[
d^{\,(j)}(\Delta_t)\to d^{\,(j)}( h^\vee_i(x,v))
\]
as $t\to 0$.
Since $d^{\,(0,j)}f\colon U\times (V\times E_2^j)\to F$ is $C^{1,0}$,
there is a continuous map
\[
(d^{\,(0,j)}f)^{[1,0]}\colon U^{[1]}\times (V\times E_2^j)\to F
\]
extending the difference quotient map, as in Lemma~\ref{f1withpar}.
By Proposition~\ref{ctsexp}, the corresponding map
\[
((d^{\,(0,j)}f)^{[1,0]})^\vee\colon U^{[1]}\to C(V\times E_2^j,F)
\]
is continuous.
Now
\begin{eqnarray*}
(d^{\,(j)}(\Delta_t))(y,w)&=& \frac{(d^{\,(j)}\circ  f^\vee)(x+tv)-(d^{\,(j)}\circ f^\vee)(x)}{t}(y,w)\\
&=&\frac{(d^{\,(0,j)}f)^\vee(x+tv)-(d^{\,(0,j)}f)^\vee(x,v)}{t}(y,w)\\
&=&
\frac{d^{\,(0,j)}f(x+tv,y,w)-d^{\,(0,j)}f(x,y,w)}{t}\\
&=& (d^{\,(0,j)}f)^{[1,0]}((x,v,t),(y,w)).
\end{eqnarray*}
Hence
\[
d^{\,(j)}(\Delta_t)=((d^{\,(0,j)}f)^{[1,0]})^\vee(x,v,t),
\]
which converges to
\[
((d^{\,(0,j)}f)^{[1,0]})^\vee(x,v,0)\colon (y,w)\mto d^{\,(1,j)}f(x,y,v,w)
\]
(which equals $(d^{\,(j)}\circ h^\vee_1)(x,v)$) as $t\to 0$.

Now assume that $2\leq i\leq k$ and that (\ref{theformul}) has been established for $i-1$
in place of~$i$.
Let $x\in U^0$ and $v=(v_1,\ldots, v_i)\in E_1^i$.
Set $v':=(v_1,\ldots,v_{i-1})$.
We have $x+\bD_\ve  v_i\sub U^0$ for some $\ve>0$.
For $t\in \bD_\ve\setminus\{0\}$,
consider
\[
\Delta_t:=
\frac{d^{\,(i-1)}(f^\vee)(x+tv_i,v')-d^{\,(i-1)}(f^\vee)(x,v')}{t} \in C^\ell(V,F).
\]
We show that
\[
d^{\,(j)}(\Delta_t)\to d^{\,(j)}(h^\vee_i(x,v))
\]
as $t\to 0$.
Since $d^{\,(i-1,j)}f\colon U\times (V\times E_1^{i-1}\times E_2^j)\to F$ is $C^{1,0}$,
there is a continuous map
\[
(d^{\,(i-1,j)}f)^{[1,0]}\colon U^{[1]}\times (V\times E_1^{i-1}\times E_2^j)\to F
\]
extending the difference quotient map, as in Lemma~\ref{f1withpar}.
Hence also
\[
g\colon U^{[1]}\times (V\times E_2^j)\to F,\quad (x,u,t,y,w)\mto (d^{\,(i-1,j)}f)^{[1,0]}((x,u,t),(y,v',w))
\]
is continuous.
By Proposition~\ref{ctsexp}, the corresponding map
\[
g^\vee\colon U^{[1]}\to C(V\times E_2^j,F)
\]
is continuous.
Now
\begin{eqnarray*}
d^{\,(j)}(\Delta_t)(y,w)&=& \frac{(d^{\,(j)}\circ  h_{i-1}^\vee)(x+tv_i,v')-(d^{\,(j)}\circ h_{i-1}^\vee)(x,v')}{t}(y,w)\\
&=&
\frac{d^{\,(i-1,j)}f(x+tv_i,y,v',w)-d^{\,(i-1,j)}f(x,y,v',w)}{t}\\
&=& g^\vee(x,v_i,t)(y,w)
\end{eqnarray*}
and thus
\[
d^{\,(j)}(\Delta_t)=g^\vee(x,v_i,t),
\]
which converges to
\[
g^\vee(x,v_i,0)\colon (y,w)\mto d^{\,(i,j)}f(x,y,v,w)
\]
(i.e, to $(d^{\,(j)}\circ h^\vee_i)(x,v)$) as $t\to 0$.
By the preceding, $f^\vee\in C^k(U,C^\ell(V,F))$ and thus $\Phi$ can be defined.
It is clear that $\Phi$ is linear and injective.

\emph{The map $\Phi$ is a topological embedding.}
To se this, let $\cO$ be the compact-open $C^{k,\ell}$-topology on $C^{k,\ell}(U\times V,F)$
and $\cT$ be the initial topology on\linebreak
$C^{k,\ell}(U\times V, F)$ with respect to~$\Phi$,
which turns~$\Phi$ into a topological embedding.
By transitivity of initial topologies, $\cT$ is initial with repect to the maps
\[
d^{\,(i)}\circ\Phi\colon C^{k,\ell}(U\times V,F)\to C(U\times E_1^i,C^\ell(V,F)).
\]
As the topology on $C^\ell(V,F)$ is initial with respect to the mappings\linebreak
$d^{\,(j)}\colon
C^\ell(V,F)\to C(V\times E_2^j,F)$, we deduce from Lemma~\ref{inipush}
that the topology on $C(U\times E_1^i,C^\ell(V,F))$ is initial with respect to the maps
\[
C(U\times E_1^i,d^{\,(j)})\colon C(U\times E_1^i,C^\ell(V,F))\to C(U\times E_1^i,C(V\times E_2^j,F)),\,
\gamma\mto d^{\,(j)}\circ\gamma.
\]
Hence $\cT$ is initial with respect to the maps
\[
\phi_{i,j}:=C(U\times E_1^i,d^{\,(j)})\circ d^{\,(i)}.
\]
Now $\cO$ is initial with respect to the maps
\[
d^{\,(i,j)}\colon C^{k,\ell}(U\times V,F)\to C(U\times V\times E_1^i\times E_2^j,F).
\]
Lemma~\ref{pubas} implies that the map
\[
\Theta_{i,j}\colon
C(U\times V\times E_1^i\times E_2^j,F)\to C(U\times E_1^i\times V\times E_2^j,F)
\]
given by $\Theta_{i,j}(\gamma)(x,v,y,w):=\gamma(x,y,v,w)$ is a homeomorphism,
whence its domain carries the initial topology with respect to $\Theta_{i,j}$.
By Proposition~\ref{ctsexp}, the map
\[
\Psi_{i,j}\colon C(U\times E_1^i\times V\times E_2^j,F)\to C(U\times E_1^i,C(V\times E_2^j,F)),\;\;
\gamma\mto\gamma^\vee
\]
is a topological embedding, whence its domain carries the initial topology with respect to $\Psi_{i,j}$.
Hence, by transitivity of initial topologies, $\cO$ is the initial topology
with respect to the mappings $\Psi_{i,j}\circ\Theta_{i,j}\circ d^{\,(i,j)}$.
Since
\[
\Psi_{i,j}\circ\Theta_{i,j}\circ d^{\,(i,j)}=\phi_{i,j},
\]
we deduce that $\cO=\cT$. Hence $\Phi$ is a topological embedding.

\emph{If condition} (a) \emph{is satisfied, then $\Phi$ is surjective}.
To see this, let $g\in C^k(U,C^\ell(V,F))$.
As we assume that $V$ is locally compact (condition~(a)), the evaluation map
$\ve\colon C^\ell(V,F)\times V\to F$ is $C^{\infty,\ell}$ (Lemma~\ref{evaldiffprop}), whence
\[
g^\wedge :=\ve\circ (g\times \id_V)\colon U\times V\to F,\quad (x,y)\mto g(x)(y)
\]
is $C^{k,\ell}$ by Chain Rule~1 (Proposition~\ref{chainR1}).
Now $\Phi(g^\wedge)=g$ since $(g^\wedge)^\vee(x)(y)$ $=(g^\wedge)_x(y)=g^\wedge(x,y)=g(x)(y)$
for all $x\in U$ and $y\in V$.

\emph{If condition}~(b) \emph{is satisfied, then $\Phi$ is surjective.}
To see this, let $g\in C^k(U,C^\ell(V,F))$.
Each of the maps $\Psi_{i,j}$ from above is a homeomorphism (see Proposition~\ref{ctsexp}),
whence also
\[
\Psi_{i,j}\circ \Theta_{i,j}\colon C(U\times V\times E_1^i\times E_2^j,F)\to C(U\times E_1^i,C(V\times E_2^j,F))
\]
is a homeomorphism. Hence, there is a unique $f_{i,j}\in C(U\times V\times E_1^i\times E_2^j,F)$ such that
\begin{equation}\label{hereverbatim}
(\Psi_{i,j}\circ\Theta_{i,j})(f_{i,j})=d^{\,(j)}\circ (d^{\,(i)}g).
\end{equation}
Then $f:=f_{0,0}\in C(U\times V,F)$ satisfies $f^\vee=g$.
Thus $\Phi(f)=g$ if we can show that $f$ is $C^{k,\ell}$.
It suffices to show that the iterated directional derivative $d^{\,(i,j)}f(x,y,v,w)$
exists and equals $f_{i,j}(x,y,v,w)$ for all $i,j\in\N_0$ with $i\leq k$ and $j\leq \ell$,
and all $v\in E_1^i$ as well as $w\in E_2^j$, if $(x,y)\in U^0\times V^0$.
This will hold if we can show that $f|_{U^0\times (V^0\cap H)}$ is $C^{k,\ell}$
for each finite-dimensional vector subspace $H\sub E_2$, and
\begin{equation}\label{fijwant}
d^{\,(i,j)}(f|_{U^0\times (V^0\cap H)})=f_{i,j}|_{U^0\times (V^0\cap H)\times E_1^i\times H^j}
\end{equation}
for all $i,j$ as before. Let $\lambda\colon H\to E_2$ be the inclusion map.
Since $\lambda|_{V^0\cap H}$ is $C^\ell$, the restriction map
\[
\rho:=C^\ell(\lambda|_{V^0\cap H},F)\colon C^\ell(V,F)\to C^\ell(V^0\cap H,F),\quad \eta\mto\eta|_{V^0\cap H}
\]
is continuous linear (see Proposition~\ref{Crpuba}) and hence smooth, entailing that
\[
g_H:=\rho\circ g|_{U^0}\colon U^0\to C^\ell(V^0\cap H,F)
\]
is~$C^k$. Since $V^0\cap H$ is locally compact, we know that
\[
(g_H)^\wedge \colon U^0\times (V^0\cap H)\to F
\]
is a $C^{k,\ell}$-map, with
\begin{equation}\label{ingredi}
(d^{\,(j)}\circ (d^{\,(i)}g_H))(x,v)(y,w)=d^{\,(i,j)}g_H^\wedge(x,y,v,w)
\end{equation}
(cf.\ (\ref{mightbeuse})).
Then $g_H^\wedge(x,y)=g_H(x)(y)=g(x)(\lambda(y))=g(x)(y)=f(x,y)$ for all $(x,y)\in U^0\times (V^0\times H)$,
whence $f|_{U^0\times (V^0\cap H)}=g_H^\wedge$ is $C^{k,\ell}$ and
\begin{equation}\label{ingredii}
(d^{\,(i,j)}f)|_{U^0\times (V^0\cap H)\times E_1^i\times H^j}=d^{\,(i,j)}(g_H^\wedge).
\end{equation}
We have $d^{\,(i)}(g|_{U^0})=(d^{\,(i)}g)|_{U^0\times E_1^i}$
and hence
\begin{equation}\label{ingrediii}
d^{\,(i)}g_H=\rho\circ d^{\,(i)}(g|_{U^0})=\rho\circ (d^{\,(i)}g)|_{U^0\times E_1^i},
\end{equation}
using that $\rho$ is continuous and linear. Moreover,
\[
d^{\,(j)}(\rho(\gamma))=d^{\,(j)}(\gamma|_{V^0\cap H})=(d^{\,(j)}\gamma)|_{(V^0\cap H)\times H^j}
\]
for each $\gamma\in C^ \ell(V,F)$. Applying this to $\gamma:=(d^{\,(i)}g_H)(x,v)$ (as in (\ref{ingrediii}))
with $x\in U^0$ and $v \in E_1^i$,
we obtain
\[
(d^{\,(j)}\circ d^{\,(i)}g_H)(x,v)(y,w)=d^{\,(j)}(d^{\,(i)}g(x,v))(y,w)=f_{i,j}(x,y,v,w)
\]
for all $(y,w)\in (V^0\cap H)\times H^j$.
Together with (\ref{ingredi}) and (\ref{ingredii}), this implies (\ref{fijwant}).
The proof is complete.\qed
\begin{rem}\label{defversexp}
It is possible to replace condition (b)
in Theorem~\ref{explawCkell} by the following condition:
\begin{description}[(D)]
\item[\rm(i)]
$U\times V\times E_1\times E_2$ is a $k_\R$-space; or
\item[\rm(ii)] $\ell=0$ and $U\times V\times E_1$ is a $k_\R$-space; or
\item[\rm(iii)]
$k=0$ and $U\times V\times E_2$ is a $k_\R$-space.
\end{description}
[Assume that (i) holds. Let $g\in C^k(U,C^\ell(V,F))$.
For $i,j\in\N_0$ such that $i\leq k$ and $j\leq \ell$, define $f_{i,j}\colon U\times V\times E_1^i\times E_2^j\to F$ via
\[
f_{i,j}(x,y,v,w):=(d^{\,(j)}\circ (d^{\,(i)}g))(x,v)(y,w).
\]
Consider
\[
\theta_i\colon V\times E_2\to V\times E_2^i,\quad (y,w)\mto (y,w,\ldots,w).
\]
Then $C(\theta_i,F)\colon C(V\times E_2^i,F)\to C(V\times E_2,F)$, $\gamma\mto\gamma\circ \theta_i$
is a continuous map (see Proposition~\ref{pubas}), whence
\[
g_{i,j}\colon U\times E_1\to C(V\times E_2,F),
\]
\[
(x,v)\mto C(\theta_j,F)(d^{\,(j)}(\delta^i_xg(v)))=d^{\,(j)} (\delta^i_xg(v))\, \circ\,\theta_j
\]
is continuous. Using hypothesis (i), Proposition~\ref{ctsexp} shows that
\[
g_{i,j}^\wedge\colon U\times E_1\times V\times E_2\to F
\]
is continuous. Now $f_{i,j}(x,y,v,w)$ is symmetric $i$-linear in $v\!=\!(v_1,\ldots, v_i)\!\in\! E_1^i$\linebreak
for fixed other arguments. Hence, if we define $h_{i,j}\colon U\times V\times E_1^i\times E_2\to F$
using the corresponding homogeneous polynomials,
\[
h_{i,j}(x,y,u,w):=f_{i,j}(x,y,(u,\ldots, u),w)
\]
(with $i$ occurences of $u\in E_1$),
then
\[
f_{i,j}(x,y,v,w)=
\frac{1}{i!\, 2^i}\sum_{\ve_1,\ldots,\ve_i\in \{1,-1\}}
\ve_1\cdots\ve_i\,
h_{i,j}(x,y,\ve_1v_1+\cdots+\ve_iv_i,w)
\]
by the Polarization Formula (Proposition~\ref{proppolarvar}). But $f_{i,j}(x,y,v,w)$
(and hence also $h_{i,j}(x,y,u,w)$) is symmetric $j$-linear in $w\in E_2^j$,
and taking the corresponding homogeneous polynomials yields the map $g_{i,j}^\wedge$.
Using the Polarization Formula again, we find that
\[
f_{i,j}(x,y,v,w)=
\frac{1}{i!j!\, 2^{i+j}}\sum_{\ve_1,\ldots,\ve_i\in \{1,-1\}}\sum_{\sigma_1,\ldots, \sigma_j\in\{1,-1\}}
\ve_1\cdots\ve_i\,\sigma_1\cdots\sigma_j\quad\vspace{-1mm}
\]
\[
\qquad\qquad\qquad g_{i,j}^\wedge(x,y,\ve_1v_1+\cdots+\ve_iv_i,\sigma_1w_1+\cdots+\sigma_j w_j).
\]
Hence $f_{i,j}$ is continuous and we can now show that $f:=f_{0,0}\in C^{k,\ell}(U\times V,F)$
and $\Phi(f)=g$ as in the proof of Theorem~\ref{explawCkell}(b), starting at (\ref{hereverbatim}).]

The cases (ii) and (iii) are similar (see Exercise~\ref{exc-smallerp});
only one polarization argument is needed then.
\end{rem}
\begin{prop}\label{CkvsCkk}
Let $E_1$, $E_2$, and $F$ be locally convex spaces,
$U\sub E_1$ and $V\sub E_2$ be locally convex subsets with dense interior
and $k,\ell\in\N_0\cup\{\infty\}$.
Then the following holds:
\begin{description}
\item[\rm(a)]
Every $C^{k+\ell}$-map $f\colon U\times V\to F$ is $C^{k,\ell}$, and the inclusion map
\[
\lambda\colon C^{k+\ell}(U\times V,F)\to C^{k,\ell}(U\times V,F)
\]
is continuous and linear.
\item[\rm(b)]
Every $C^{k,k}$-map $f\colon U\times V\to F$ is $C^k$, and the inclusion map
\[
\mu\colon C^{k,k}(U\times V,F)\to C^k(U\times V,F)
\]
is continuous and linear.
\end{description}
\end{prop}
\begin{prf}
(a) For all $i,j \in\N_0$ such that $i\leq k$ and $j\leq \ell$,
the following iterated
directional derivatives exist
and coincide:
\begin{equation}\label{oarhs}
d^{\,(i,j)}f(x,y,v,w)=d^{\,(i+j)}f(x,y,(0,w_1),\ldots, (0,w_j),(v_1,0),\ldots, (v_i,0))
\end{equation}
for $v=(v_1,\ldots, v_i)\in E_1^i$, $w=(w_1,\ldots, w_j)\in E_2^j$ and $(x,y)\in U^0\times V^0$.
As the right hand side of (\ref{oarhs}) defines a continuous $F$-valued function of
$(x,y,v,w)\in U\times V\times E_1^i\times E_2^j$, we see that $f$ is $C^{k,\ell}$
and that (\ref{oarhs}) holds for all $i,j$ as before and all $(x,y,v,w)\in U\times V\times E_1^i\times E_2^j$.
Note that\linebreak
$h_{i,j}\colon U\times V\times E_1^i\times E_2^j\to U\times V\times (E_1\times E_2)^{i+j}$,
\[
(x,y,v,w)\mto (x,y,(0,w_1),\ldots, (0,w_j), (v_1,0),\ldots, (v_i,0))
\]
is a continuous map. As a consequence, also the map
\[
C(h_{i,j},F)\colon C(U\times V\times (E_1\times E_2)^{i+j},F)\to C(U\times V\times E_1^i\times E_2^j,F),\;\gamma\mto
\gamma\circ h_{i,j}
\]
is continuous.
Now
\[
d^{\,(i,j)}\circ \lambda=C(h_{i,j},F)\circ d^{\,(i+j)}
\]
by (\ref{oarhs}), which is a continuous function since $C(h_{i,j},F)$ and
\[
d^{\,(i+j)}\colon C^{k+\ell}(U\times V,F)\to C(U\times V\times (E_1\times E_2)^{i+j},F)
\]
are continuous. As the topology on $C^{k,\ell}(U\times V,F)$ is initial with respect
to the maps $d^{\,(i,j)}\colon C^{k,\ell}(U\times V,F)\to C(U\times V\times E_1^i\times E_2^j,F)$,
we deduce that $\lambda$ is continuous.

(b) We show by induction on $m \in\N$ that every $f\in C^{m,m}(U\times V,F)$ is $C^m$, with
\begin{equation}\label{prepbtop}
d^{\,(m)}f(x,y,w)=\!\!\sum_{I\sub\{1,\ldots, m\}} d^{\,(a,m-a)}
f(x,y,w_{i_1,1},\ldots, w_{i_a,1},w_{j_1,2},\ldots,w_{j_{m-a},2})
\end{equation}
for $(x,y,w)\in U\times V\times (E_1\times E_2)^m$, where $w=(w_1,\ldots, w_m)$
with $w_b=(w_{b,1},w_{b,2})\in E_1\times E_2$
for $b\in\{1,\ldots, m\}$,
\[
I=\{i_1,\ldots, i_a\}\quad\mbox{with}\quad i_1<i_2<\cdots < i_a
\]
and $\{1,\ldots,m\}\setminus I=\{j_1,\ldots, j_{m-a}\}$ with $j_1<\ldots < j_{m-a}$.
If $m=1$, then $f$ is $C^1$ by the Rule on Partial Differentials (Proposition~\ref{rulepartialno})
and
\[
df(x,y,w_1)=d^{\,(1,0)}f(x,y,w_{1,1})+d^{\,(0,1)}f(x,y,w_{1,2}),
\]
whence (\ref{prepbtop}) holds.
Let $f$ be a $C^{m+1,m+1}$-map now with $m\in\N$ and
assume that (\ref{prepbtop}) holds for~$m$. Fix $I\sub\{1,\ldots, m\}$ and $w\in (E_1\times E_2)^m$,
and let $z:=(w,w_{m+1})$ with $w_{m+1}=(w_{m+1,1},w_{m+1,2})\in E_1\times E_2$.
Abbreviate
\[
w^I:=(w_{i_1,1},\ldots, w_{i_a,1},w_{j_1,2},\ldots,w_{j_{m-a},2}).
\]
As a consequence of the Schwarz Theorem (Proposition~\ref{schwarzCkell}), the map
\[
h_{I,w}\colon U\times V\to F,\quad (x,y)\mto d^{\,(a,m-a)}f(x,y,w^I)
\]
is $C^{1,1}$ with
\[
d^{\,(1,0)}h_{I,w}(x,y,w_{m+1,1})=d^{\,(a+1,m-a)}f(x,y,(w,w_{m+1})^{I\cup\{m+1\}})
\]
and
\[
d^{\,(0,1)}h_{I,w}(x,y,w_{m+1,2})=d^{\,(a,m-a+1)}f(x,y,(w,w_{m+1})^{I}).
\]
By the Rule on Partial Differentials, $h_{I,w}$ is $C^1$ with $d^{\,(1)}h_{I,w}(x,y,w_{m+1})$ equal to
\[
d^{\,(a+1,m-a)}f(x,y,(w,w_{m+1})^{I\cup\{m+1\}})+
d^{\,(a,m-a+1)}f(x,y,(w,w_{m+1})^{I}).
\]
Hence $d^{\,(m+1)}f(x,y,z)$ exists and is given by
\[
\sum_{I\sub\{1,\ldots m\}}\big(d^{\,(a+1,m-a)}f(x,y,z^{I\cup\{m+1\}})
+ d^{\,(a,m-a+1)}f(x,y,z^I)\big).
\]
This is of the form (\ref{prepbtop}), with $m+1$ in place of~$m$ and $z$
in place of~$w$.

To see that $\mu$ is continuous, note that the map
\[
g_{m,I}\colon
U\times V\times (E_1\times E_2)^m\to U\times V\times E_1^a\times E_2^{m-a},\quad (x,y,w)\mto (x,y,w^I)
\]
is continuous for $I\sub \{1,\ldots, m\}$, with notation as in (\ref{prepbtop}).
Hence also
\[
C(g_{m,I},F)\colon C(U\times V\times E_1^a\times E_2^{m-a},F)\to C(U\times V\times (E_1\times E_2)^m),\;
\gamma\mto\gamma\circ g_{m,I}
\]
is continuous (see Lemma~\ref{pubas}). Note that
\[
d^{\,(m)}f=\bigcup_{I\sub \{1,\ldots,m\}}C(g_{m,I},F)(d^{\,(a,m-a)}f)\quad\mbox{for $\,f\in C^{k,k}(U\times V,F)$}
\]
for $m\in \N_0$ such that $m\leq k$, by (\ref{prepbtop}).
Since
\[
d^{\,(a,m-a)}\colon C^{k,k}(U\times V,F)\to C(U\times V\times E_1^a\times E_2^{m-a},F)
\]
is continuous, we deduce that $d^{\,(m)}\circ\mu$ is continuous.
Hence $\mu$ is continuous, using that the topology on $C^k(U\times V,F)$ 
is initial with respect to the mappings $d^{\,(m)}\colon C^k(U\times V,F)\to C(U\times V\times (E_1\times E_2)^m,F)$.
\end{prf}
Taking $k=\ell=\infty$, Proposition~\ref{CkvsCkk} subsumes the following:
\begin{cor}
Let $E_1$, $E_2$, and $F$ be locally convex spaces,
$U\sub E_1$ and $V\sub E_2$ be locally convex subsets with dense interior.
Then a mapping\linebreak
$f\colon U\times V\to F$ is $C^\infty$ if and only if it is $C^{\infty,\infty}$.
Moreover,
\[
C^\infty(U\times V,F)=C^{\infty,\infty}(U\times V,F)
\]
as a locally convex space, i.e.\ the compact-open $C^\infty$-topology and the compact-open
$C^{\infty,\infty}$-topology coincide.\qed
\end{cor}
Taking $k=\ell=\infty$ in Theorem~\ref{explawCkell} (and Remark~\ref{defversexp}), and using
that $C^{\infty,\infty}(U\times V,F)=C^\infty(U\times V,F)$,
we obtain an important special case:
\begin{cor}[Exponential law for spaces of smooth functions]\label{cinfty-explaw}
Let $E_1$, $E_2$, and $F$ be locally convex spaces,
and $U\sub E_1$ as well as $V\sub E_2$ be locally convex subsets with dense interior.
If $f \in C^\infty(U\times V,F)$, set $f^\vee(x):=f_x:=f(x,\cdot )$
for $x\in U$. Then
\[
f^\vee\colon U\to C^\infty(V,F)
\]
is a $C^\infty$-map and the mapping
\[
\Phi\colon C^\infty(U\times V,F)\to C^\infty(U,C^\ell(V,F)),\quad f \mto f^\vee
\]
is linear and a topological embedding. If
$V$ is locally compact or
$U\times V\times E_1\times E_2$ is a $k_\R$-space,
then $\Phi$ is an isomorphism of topological vector spaces.\qed
\end{cor}
We close this section with
five lemmas which can be skipped on a first reading.
Lemmas~\ref{C0kviaC01},
\ref{C0chainy}, and
\ref{explaw-vari}
will be used for
technical results which eventually feed into
the proof of $C^0$-regularity for Banach-Lie groups and diffeomorphism groups.
Lemmas~\ref{fstamult} and~\ref{fstaFCk} (or Exercise~\ref{exc-ifFCk})
are needed for a result concerning differential
equations in Banach spaces (Corollary~\ref{paraFCk}), which will not be used elsewhere.
\begin{lem}\label{C0kviaC01}
Let $X$ be a topological space, $E$ and $F$ be locally convex spaces,
$U\sub E$ a locally convex subset with dense interior and $f\colon X\times U\to F$ be a
$C^{0,1}$-function. If $k\in\N$
and $d^{(0,1)}f\colon X\times (U\times E)\to F$ is $C^{0,k}$,
then~$f$ is~$C^{0,k+1}$.
\end{lem}
\begin{prf}
For $x\in X$, $y\in U^0$, $j\in\{2,\ldots, k+1\}$, and $w_1,\ldots, w_j\in E$,
the iterated directional derivative $d^{(0,j)}f(x,y,w_1,\ldots, w_j)$
exists as it is given by
\begin{eqnarray*}
(D_{w_j}\cdots D_{w_1}(f_x))(y)&=&(D_{(w_j,0)}\cdots D_{(w_2,0)}(d^{(0,1)}f)_x)(y,w_1)\\
&=& d^{(0,j-1)}(d^{(0,1)}f)(x,(y,w_1),(w_2,0),\ldots,(w_j,0)).
\end{eqnarray*}
As $f$ is $C^{0,1}$ and the right hand side of the preceding formula defines a continuous function of
$(x,y,w_1,\ldots,w_j)\in X\times U\times E^j$ for all $j\in\{2,\ldots,k+1\}$,
we see that $f$ is~$C^{0,k+1}$.
\end{prf}
Here is  another version of the Chain Rule.
\begin{lem}\label{C0chainy}
Let $g_1\colon X\to Y$ be a continuous map between topological spaces.
Let $E$, $F$, and $H$ be locally convex spaces,
$U\sub E$ and $V\sub F$ be locally convex subsets with dense interior,
$r\in\N\cup\{\infty\}$ and $f\colon Y\times V\to H$
as well as
$g_2\colon X\times U\to F$ be $C^{0,r}$-maps.
If $g_2(X\times U)\sub V$, then also the following map is~$C^{0,r}$:
\[
h\colon X\times U\to H,\quad (x,y)\mto f(g_1(x),g_2(x,y)).
\]
\end{lem}
\begin{prf}
We may assume that $r\in\N$.
Holding $x$ fixed, both $f(g_1(x),\cdot)$ and $g_2(x,\cdot)$ are~$C^1$-maps
with differentials $d^{(0,1)}f(g_1(x),\cdot)$ and $d^{(0,1)}g_2(x,\cdot)$,
respectively (see Lemma~\ref{Ckellpartial}).
Thus $h_x:=h(x,\cdot)$ is $C^1$ with $d^{(0,1)}h(x,y,w)=dh_x(y,w)=d^{(0,1)}f(g_1(x),g_2(x,y),
d^{(0,1)}g_2(x,y,w))$, which depends continuously on $(x,y,w)\in X\times U\times E$.
Hence $h$ is $C^{0,1}$. Moreover,
\[
d^{(0,1)}h(x,y,w)=d^{(0,1)}f(g_1(x),g_2(x,y),d^{(0,1)}g_2(x,y,w))
\]
is a $C^{0,r-1}$-function of $(x,y,w)\in X\times (U\times E)$
by induction, exploiting that $d^{(0,1)}f\colon Y\times (V\times F)\to H$ is $C^{0,r-1}$
by Lemma~\ref{diffphigherd} and $X\times (U\times E)\to F\times F$, $(x,y,w)\mto (g_2(x,y),d^{(0,1)}g_2(x,y,w))$
is $C^{0,r-1}$ since its components are so, by Exercise~\ref{exc-babychainC0k}
and Lemma~\ref{diffphigherd}.
Thus~$h$ is $C^{0,r}$, by Lemma~\ref{C0kviaC01}.
\end{prf}
Finally, we need another variant of the exponential law.
\begin{lem}\label{explaw-vari}
Let $Y$ and $Z$ be locally compact topological spaces,
$E$ and $F$ be locally convex spaces, $U\sub E$ be a locally convex subset with dense interior,
and $k\in\N_0\cup\{\infty\}$.
Endow $C(Z,F)$ with the compact-open topology.
If $f\colon U\times Y\to C(Z,F)$ is a $C^{k,0}$-map,
then also the following map is $C^{k,0}$:
\[
f^\wedge\colon U\times (Y\times Z)\to F,\quad (x,y,z)\mto f(x,y)(z).
\]
\end{lem}
\begin{prf}
By Theorem~\ref{explawCkell} and Remark~\ref{inwashexp}(a),
the mapping
\[
f^\vee\colon U\to C(Y,C(Z,F)), \quad x \mto f(x,\cdot)
\]
is~$C^k$.
The linear map
$\Phi\colon C(Y\times Z,F)\to C(Y,C(Z,F))$, $\eta\mto\eta^\vee$
(with $\eta^\vee(y)(z):=\eta(y,z)$) is a homeomorphism,
by Proposition~\ref{ctsexp}.
Hence
\[
\Psi:=\Phi^{-1}\colon C(Y,C(Z,F))\to C(Y\times Z,F),\quad \gamma\mto \gamma^\wedge
\]
is a continuous linear mapping (where $\gamma^\wedge(y,z):=\gamma(y)(z)$).
Thus \break $g:=\Psi\circ f^\vee\colon$ $U\to C(Y\times Z,F)$ is~$C^k$, whence
\[
g^\wedge\colon U\times (Y\times Z)\to F,\quad g^\wedge(x,y,z):=g(x)(y,z)
\]
is $C^{k,0}$, by Theorem~\ref{explawCkell} and Remark~\ref{inwashexp}(a).
But
$g^\wedge=f^\wedge$ as $g^\wedge(x,y,z)$ $=g(x)(y,z)\!=\!\Psi(f^\vee(x))(y,z)\!=
\!\Psi(f(x,\cdot))(y)(z)\!=\!
f(x,y)(z)\!=\!f^\wedge(x,y,z)$.
\end{prf}
The final lemmas use
notation from Section~\ref{provis-Ban},
and Proposition~\ref{companotions}.
\begin{lem}\label{fstamult}
Let $K$ be a compact topological space, $P$ be a topological space,
$k\in \N$, and $(E_j,\|\cdot\|_j)$ be normed spaces for $j\in\{1,\ldots, k\}$.
Let $(E,\|\cdot\|_E)$ and $(F,\|\cdot\|_F)$ be normed spaces,
$U\sub E$ be an open subset and
\[
f\colon K\times U\times (E_1\times\cdots\times E_k)\times P\to F
\]
be a continuous mapping such that $g(x,y,p):=f(x,y,\cdot,p)\colon E_1\times \cdots\times E_k\to F$
is $k$-linear for all $(x,y,p)\in K\times U\times P$ and the map
\[
g\colon K\times U\times P\to \cL^k(E_1,\ldots, E_k;F)_b
\]
so obtained is continuous. Then also the map
\[
\Theta \colon C(K,U)\times P\to \cL^k(C(K,E_1),\ldots,C(K,E_k);C(K,F))_b
\]
determined by $\Theta(\gamma,p)(\gamma_1,\ldots,\gamma_k)(x):=
f(x,\gamma(x),\gamma_1(x),\ldots,\gamma_k(x),p)$ is continuous.
\end{lem}
\begin{prf}
In the usual way, we identify $C(K,U\times E_1\times \cdots \times E_k)$ and the product\linebreak $C(K,U)\times C(K,E_1)\times\cdots \times C(K,E_k)$.
By Lemma~\ref{pushpar}, the map
\[
\Phi\colon C(K,U)\times C(K,E_1)\times\cdots\times C(K,E_k)\times P\to C(K,F)
\]
determined by $\Phi(\gamma,\gamma_1,\ldots,\gamma_k,p)(x)
:=f(x,\gamma(x),\gamma_1(x),\dots,\gamma_k(x),p)$
is continuous.
Therefore, by Proposition~\ref{ctsexp}, $\Theta$ is continuous as a mapping to\linebreak
$\cL^k(C(K,E_1),\ldots, C(K,E_k);C(K,F))_c$.
To see that~$\Theta$ is also continuous with respect to the operator norm on its range, let $\gamma\in C(K,U)$, $p\in P$, and $\ve>0$. Since $\gamma(K)\sub U$ is compact,
there exists $r>0$ with $\gamma(K)+B^E_r(0)\sub U$.
As
\[
h(x,y,q):=g(x,\gamma(x)+y,q)-g(x,\gamma(x),p)
\]
is continuous in $(x,y,q)\in K\times B^E_r(0)\times P$
and $h(x,0,p)=0$, using the Wallace Lemma we find an open $p$-neighborhood $P_0\sub P$
and $s\in \,]0,r]$ with
\[
\|h(x,y,q)\|_{\op}\leq \ve \quad\mbox{for all $\,(x,y,q)\in K\times B^E_s(0)\times P_0$.}
\]
For $\eta\in C(K,B^E_s(0))$ and $q\in P_0$, this entails that
\begin{eqnarray*}
\hspace*{-3mm}\lefteqn{\|(\Theta(\gamma+\eta,q)-\Theta(\gamma,p))(\gamma_1,\ldots,\gamma_k)(x)\|_F}\qquad
\qquad\qquad\\
&=&\|h(x,\eta(x),q)(\gamma_1(x),\ldots,\gamma_k(x))\|_F
\leq \|h(x,\eta(x),q)\|_{\op}\leq \ve
\end{eqnarray*}
for all $\gamma_j\in C(K,E_j)$ of supremum norm
$\|\gamma_j\|_\infty\leq 1$ for $j\in\{1,\ldots, k\}$
and $x\in K$. Taking the supremum over~$x$, we get
\[
\|\big(\Theta(\gamma+\eta,q)-\Theta(\gamma,p)\big)(\gamma_1,\ldots,\gamma_k)\|_{\infty}
\leq \ve.
\]
Taking the supremum over $\gamma_1,\ldots,\gamma_k$ as before, we get
\[
\|\Theta(\gamma+\eta,q)-\Theta(\gamma,p)\|_{\op}\leq\ve
\]
for all $q\in P_0$ and $\eta\in C(K,B^E_s(0))$. Thus~$\Theta$ is continuous at~$(\gamma,p)$.
\end{prf}
With notation as in Proposition~\ref{pfctisCk}, 
we have:
\begin{lem}\label{fstaFCk}
Let $(E,\|\cdot\|_E)$, $(F,\|\cdot\|_F)$, and $(Z,\|\cdot\|_Z)$
be normed spaces, $U\sub E$ be an open subset,
$P\sub Z$ be a locally convex subset with dense interior,
$K$ be a compact topological
space and $k\in\N_0\cup\{\infty\}$. Moreover,
let
$f\colon K\times (U\times P) \to F$ be 
a $C^{0,k}$-map such that
\begin{equation}\label{partbcont}
(d^{\,(0,j)}f)^\vee\colon K\times U\times P\to\cL^j(E\times Z,\ldots,E\times Z;F)_b
\end{equation}
is continuous for all $j\in\N$ with $j\leq k$.
Then the following map is $FC^k$:
\[
\Phi\colon C(K,U)\times P\to C(K,F), \quad 
(\gamma,p)\mto (f^p)_*(\gamma) 
= f^p \circ (\id_K,\gamma). 
\]
\end{lem}
\begin{prf}
Let $\lambda\colon Z\to C(K,Z)$ be the map taking $p\in Z$ to the constant function
$\lambda(p)\colon x\mto p$.
Thus~$\lambda$ is linear and an isometry. By Proposition~\ref{pfctisCkpar},
$\Phi$ is~$C^k$ and by Exercise~\ref{exc-higherdfsta},
we have
\begin{eqnarray}
\lefteqn{d^{\,(j)}\Phi(\gamma,p,\gamma_1,p_1,\ldots,\gamma_j,p_j)(x)}\qquad\qquad\notag\\
&=&d^{\,(0,j)}f(x,\gamma(x),p,\gamma_1(x),p_1,\ldots,\gamma_j(x),p_j)\label{thedofphi}
\end{eqnarray}
for all $j\in \N$ such that $j\leq k$, $\gamma\in C(K,U)$,
$\gamma_1,\ldots,\gamma_j\in C(K,E)$, and $p\in U$, $p_1,\ldots, p_j\in Z$.
By Lemma~\ref{fstamult}, the map
\[
\Theta_j\colon C(K,U)\times P\to\cL^j(C(K,E\times Z),\ldots, C(K,E\times Z); C(K,F))_b
\]
determined by $\Theta_j(\gamma,p)(\eta_1,\ldots,\eta_j)(x):=
d^{\,(0,j)}f(x,\gamma(x),p,\eta_1(x),\ldots,\eta_j(x))$ is continuous (thanks to~(\ref{partbcont})).
By~(\ref{thedofphi}), the map
\[
(d^{(j)}\Phi)^\vee\colon C(K,U)\times P\to \cL^k((C(K,E)\times Z)^j;C(K,F))_b
\]
is given by
$(d^{\,(j)}\Phi)^\vee(\gamma,p)(\gamma_1,p_1,\ldots\gamma_j,p_j)=
\Theta_j(\gamma,p)(\gamma_1,\lambda(p_1),\ldots,\gamma_j,\lambda(p_j))$,
whence $(d^{\,(j)}\Phi)^\vee$ is continuous.
Thus~$\Phi$ is $FC^k$, by Proposition~\ref{companotions}(c).
\end{prf}
\begin{rem}\label{fstaFCk2}
If $U$, $P$ and~$F$ are as in Lemma~\ref{fstaFCk},
$K$ is a compact convex subset with dense interior in a finite-dimensional vector space
and $f\colon K\times U\times P\to F$ is~$FC^k$,
then all of the maps in~(\ref{partbcont}) are continuous
(as a consequence of Proposition~\ref{companotions}(c))
and the lemma applies. An alternative, easier proof for this special case
is sketched in Exercise~\ref{exc-ifFCk}.
\end{rem}
The following lemma will be useful later, when
we discuss a
topology on  $C^{r,s}(M\times N,F)$
for manifolds $M$ and $N$.
\begin{lem}\label{pb-prod-Ckl}
Let $r,s\in\N_0\cup\{\infty\}$,
$E_1$, $E_2$, $F_1$, $F_2$, and $F$
be locally convex spaces, $U_j\sub E_j$ and $V_j\sub F_j$
be locally convex subsets with dense interior
for $j\in \{1,2\}$,
$g_1\colon U_j\to F_1$ be a $C^r$-map with $g_1(U_1)\sub V_1$
and $g_2\colon U_2\to F_2$ be a $C^s$-map
with $g_2(U_2)\sub V_2$.
Then the following linear map is continuous:
\[
C^{r,s}(g_1\times g_2,F)\colon C^{r,s}(V_1\times V_2,F)\to C^{r,s}(U_1,\times U_2,F),\;\;
f\mto f\circ (g_1\times g_2).
\]
\end{lem}
\begin{prf}
From Remark~\ref{faa-chain-1},
we known that for all $(k,\ell)\in \N_0\times \N_0$ with $k\leq r$ and $\ell\leq s$,
there exist a finite set $A_{k,\ell}$, a pair
$(i(a),j(a))\in \N_0\times\N_0$ with $i(a)\leq k$ and $j(a)\leq \ell$
and a continuous function
\[
\psi_a\colon U_1\times U_2\times E_1^k\times E_2^\ell\to
V_1\times V_2\times F_1^{i(a)}\times F_2^{j(a)}
\]
 for $a\in A_{k,\ell}$ such that, for all $f\in C^{r,s}(V_1\times V_2,F)$, we have
\[
d^{\,(k,\ell)}(f\circ (g_1\times g_2))=\sum_{a\in A_{k,\ell}}
(C(\psi_a,F)\circ d^{\,(i(a),j(a))})(f).
\]
Thus $d^{\,(k,\ell)}\circ C^{r,s}(g_1\times g_2,F)$
is continuous (using Lemma~\ref{pubas}),
whence $C^{r,s}(g_1\times g_2,F)$ is continuous.
\end{prf}
\begin{lem}\label{Ckl-cover}
Let $r,s\in \N_0\cup\{\infty\}$,
$E_1$, $E_2$, and $F$ be locally convex spaces,
$U\sub E_1$ and $V\sub E_2$ be
locally convex subsets with dense interior,
$I$ be a set,
$(U_i)_{i\in I}$ an open cover of~$U$
and $(V_i)_{i\in I}$ be an open cover of $V$.
Then the compact-open $C^{r,s}$-topology $\cO$ on
$C^{r,s}(U\times V,F)$
is initial with respect to the family of restriction maps
$\rho_i\colon C^{r,s}(U\times V,F)\to C^{r,s}(U_i\times V_i,F)$, $f\mto f|_{U_i\times V_i}$
for $i\in I$.
\end{lem}
\begin{prf}
For all $k,\ell\in \N_0$ such that $k\leq r$ and $\ell\leq s$,
the compact-open topology on $C(U\times V\times E_1^k\times E_2^\ell,F)$
is initial with respect to the restriction maps
$\rho_{i,k,\ell}\colon C(U\times V\times E_1^k\times E_2^\ell)\to
C(U_i\times V_i\times E_1^k\times E_2^\ell,F)$, by Lemma~\ref{coveremb}.
As a consequence, $\cO$ is initial with respect to the mappings
$\rho_{i,k,\ell}\circ d^{\,(k,\ell)}$ (see Lemma~\ref{transinit}).
For each $i\in I$, the compact-open $C^{r,s}$-topology
on the set $C^{r,s}(U_i\times V_i,F)$ is initial with respect
to the mappings
\[
\tau_{i,k,\ell}\colon
C^{r,s}(U_i\times V_i,F)\to C(U_i\times V_i\times E_1^k\times E_2 ^\ell,F),\;\;
f\mto d^{\,(k,\ell)}f
\]
for $k,\ell\in \N_0$ with $k\leq r$ and $\ell\leq s$.
Since $\rho_{i,k,\ell}\circ d^{\,(k,\ell)}=\tau_{i,k\ell}\circ \rho_i$,
we deduce with Lemma~\ref{transinit} that $\cO$ is initial with respect to the mappings
$\rho_i$ for $i\in I$.
\end{prf}
\subsection*{Exercises for Section~\ref{secCspaces}}

\begin{small}
\begin{exer}\label{exc-babychainC0k}
Let $X$ be a Hausdorff topological space, $E_1$, $E_2$ and $F$ be locally
convex spaces, $\alpha\colon E_1\to E_2$ be a continuous linear map,
and $k\in\N$.
Let $U\sub E_1$ and $V\sub E_2$ be locally convex subsets with dense
interior such that $\alpha(U)\sub V$, and $f\colon X\times V\to F$ be a $C^{0,k}$-map.
Show that also the map
\[
g\colon X\times U\to F, \quad (x,y)\mto f(x,\alpha(y))
\]
is $C^{0,k}$, and $d^{\,(0,k)}g(x,y,z_1,\ldots,z_k)=d^{\,(0,k)}f(x,\alpha(y),\alpha(z_1),\ldots, \alpha(z_k))$
for all\linebreak
$x\in X$, $y\in U$, and $z_1,\ldots, z_k\in E_1$.
\end{exer}

\begin{exer}\label{exc-cxC0k}
Let $X$ be a topological space, $E$ and $F$ be complex locally convex spaces,
$V\sub E$ be a locally convex subset with dense interior and $k\in\N_0\cup\{\infty\}$.
Show: If a map $f\colon X\times U\to F$ is $C^{0,k}_\R$ (i.e., $C^{0,k}$
when $E$ and~$F$ are considered as real vector spaces) and $f(x,\cdot)\colon U\to F$
is~$C^k_\C$ for each $x\in X$, then~$f$ is $C^{0,k}_\C$.
\end{exer}
\begin{exer}\label{exc-intmpset}
Let $K$ be a compact topological space, $E$ be a locally convex space and $U\sub E$ be a convex subset
with non-empty interior. Show that the subset $C(K,U)$ of $C(K,E)$
has interior $C(K,U)^0=C(K,U^0)$.\\[2mm]
[Hint: If there was $\gamma\in C(K,U)^0\setminus C(K,U^0)$,
then we could find $x\in K$ such that $\gamma(x)\in U\setminus U^0$.
Derive a contradiction, using that the point evaluation
$C(K,E)\to E$, $\eta\mto \eta(x)$ is continuous
and $E\to C(K,E)$, $v\mto c_v$ (with $c_v(x):=v$ for all $v\in E$)
a continuous map (see Lemma~\ref{tocon})].
\end{exer}

\begin{exer}\label{exc-multCkell}
Let $E_1$, $E_2$, and $F$ be locally convex spaces, $U\sub E_1$ be a locally convex subset
with dense interior and $f\colon U\times E_2\to F$ be a map such that
$E_2=W_1\times\cdots\times W_m$ with locally convex spaces
$W_1,\ldots, W_m$.
Show that if $f$ is a $C^{k,0}$ and $f_x:=f(x,\cdot)$ is $m$-linear for each $x\in U$,
then $f$ is $C^{k,\infty}$.\\[1mm]
[Use Exercise~\ref{exc-explicithdmult} to find an explicit formula for
$d^{\,(0,j)}f(x,y_0,y_1,\ldots,y_j)$ for $x\in U$ and $y_0,y_1,\ldots, y_j\in E_2$.
Then proceed as in the proof of Lemma~\ref{linCkell}.]
\end{exer}

\begin{exer}\label{exc-higherdfsta}
In the situation of Proposition~\ref{pfctisCkpar}, show that
\[
d^{(j)}\Phi(\gamma,p,\gamma_1,p_1,\ldots,\gamma_j,p_j)(x)=
d^{(0,j)}f(x,\gamma(x),p,\gamma_1(x),p_1,\ldots,\gamma_j(x),p_j)
\]
for all $j\in \N$ such that $j\leq k$ and all $\gamma\in C(K,U)$,
$\gamma_1,\ldots, \gamma_j\in C(K,E)$,
$p\in P$, $p_1,\ldots, p_j\in Z$, and $x\in K$.\\[1mm]
[The point evaluation $\ve_x\colon C(K,F)\to F$, $\eta\mto\eta(x)$ is continuous linear, whence
$\ve_x\circ (d^{(j)}\Phi)=d^{(j)}(\ve_x\circ \Phi)$, where
$(\ve_x\circ\Phi)(\gamma,p)=f(x,\gamma(x),p)$.]
\end{exer}
 
\begin{exer}\label{exc-smallerp}
Prove the
cases (ii) and (iii) described in Remark~\ref{defversexp}.
\end{exer}

\begin{exer}\label{exc-convol}
Let $E$ be a sequentially complete locally convex space,
$f\colon \R^n\to\R$ be continuous and $g\colon \R^n\to E$
be a continuous function with compact support.
\begin{description}[(D)]
\item[(a)]
Show that $(f*g)(x):=\int_{\R^n} f(y)g(x-y)\,dy=\int_{\R^n} f(x-y)g(y)$
defines a continuous function $f*g\colon\R^n\to E$.\\[.3mm]
[The integral is as in Definition~\ref{def-wint-gen}
and exists by Proposition~\ref{rieman},
as it can be understood as an iterated weak integral over compact intervals.]
\item[(b)]
If the directional derivatives $D_vf(x)$ (resp., $D_vg(x)$) exist for some $v\in\R^n$
and all $x\in\R^n$,
and $D_vf\colon\R^n\to\R$ (resp., $D_vg\colon\R^n\to E$) is continuous, then
$D_v(f*g)(x)$ exists for all $x\in\R^n$ and
$D_v(f*g)=(D_vf)*g$ (resp., $f*D_vg$).
\item[(c)]
If $f$ is smooth, then also $f*g$ is smooth.
\item[(d)]
We let $h\colon \R^n\to \R$ be a smooth function
with $h(\R^n)\sub [0,\infty[$, $\Supp(h)\sub[{-1},1]^n$ and $\int_{\R^n}h(x)\,dx=1$.
For $\ve>0$, we define
\[
h_\ve\colon\R^n\to\R,\quad x\mto \frac{1}{\ve^n} h\Big(\frac{1}{\ve}x\Big).
\]
Show that $h_\ve*f\to f$ in $C(\R^n,E)$ as $\ve\to 0$.
\end{description}
\end{exer}

\begin{exer}\label{exc-the-topf}
Let $E$ be a locally convex space, $U\sub\R^n$ be a locally convex subset with dense
interior and $k\in\N_0\cup\{\infty\}$.
Show that the compact-open $C^k$-topology on $C^k(U,E)$
can be described by the seminorms $\|f\|_{C^j,K,p}$
for $j\in\N_0$ such that $j\leq k$, compact
subsets $K\sub U$ and continuous seminorms $p\colon E\to[0,\infty[$,
where
\[
\|f\|_{C^j,K,p}:=\max_{|\alpha|\leq j}\sup_{x\in K}p\Big(\frac{\partial^{|\alpha|}f}{\partial x^\alpha}(x)\Big)\quad
\mbox{for $f\in C^k(U,E)$.}
\]
\end{exer}

\begin{exer}\label{exc-Ck-complete}
Let $E$ and $F$ be locally convex spaces, $k\in \N_0\cup\{\infty\}$
and $U\sub E$ be a locally convex
subset with dense interior. Given $j\in\N_0$
with $j\leq k$ and $\gamma\in C^k(U,F)$, we define $\Delta_j\gamma\colon U\times E\to F$
via $\Delta_j\gamma(x,y):=d^{\,(j)}\gamma(x,y,\ldots,y)$ (the $j$th G\^{a}teaux-differential
$\delta^j_x\gamma(y)$).
\begin{description}[(D)]
\item[(a)]
Deduce from Lemma~\ref{ima-clo} that the map
\[
\Delta\colon C^k(U,F)\to C(U,F)\times \prod_{\N\ni j\leq k}C(U\times E,F),\;
\gamma\mto \left(\gamma,(\Delta_j(\gamma))_{\N\ni j\leq k}\right)
\]
is a linear topological embedding with closed image.\\[2mm]
[For $j\in\N$ with $j\leq k$, we have $\Delta_j\gamma=(d^{\,(j)}\gamma)\circ g_j$
using the continuous map $g_i\colon U\times E\to U\times E^j$, $(x,y)\mto (x,y,\ldots,y)$.
By the Polarization Formula, conversely $d^{\,(j)}\gamma=\sum_{i\in I_j}c_{ij}
(\Delta_j\gamma)\circ g_{ij}$ with suitable $c_{ij}\in\Q$ and continuous functions
$g_{ij}\colon U\times E^j\to U\times E$.
Lemma~\ref{pubas} enables us to pass from limits of nets
$\Delta_j\gamma_\alpha$ to limits of the $d^{\,(j)}\gamma_\alpha$,
and conversely.\,]\vspace{1mm}
\item[(b)]
For $k\in \N\cup\{\infty\}$, show: If $U\times E$ is a $k_\R$-space
and~$F$ is complete, quasi-complete, sequentially complete and
Mackey complete, respectvely, then also $C^k(U,F)$ is complete, quasi-complete, sequentially complete
and Mackey complete, respectively.\\[1.3mm]
[Argue as in the proof of Proposition~\ref{compl-and-metr}(a). For $k=0$, see
Lemma~\ref{sammelsu}.\,] 
\end{description}
\end{exer}

\begin{exer}\label{exc-gat-Ck}
Let $E$ be a locally convex space, $U\sub\R^n$ be an open subset,
$k\in\N\cup\{\infty\}$ and $f\colon U\to E$ be a continuous function whose G\^{a}teaux differentials
\[
\delta^{(j)}f(x,y):=\delta^{(j)}_xf(y):=((D_y)^jf)(x)
\]
exist for all $j\in \N$ such that $j\leq k$ and $(x,y)\in U\times\R^n$, and define
continuous functions $\delta^{(j)}f\colon U\times\R^n\to E$.
Our goal is to show that $f$ is~$C^k$.
As it suffices to see that $f$ is $C^k$ on an open neighborhood of each point $x\in U$, after multiplication with a function $\chi \in C^\infty_c(U,\R)$ which is~$1$ on a neighborhood of~$x$, and extending $\chi f$ by $0$
on $\R^n\setminus\Supp(h)$, we may assume that $U=\R^n$ and $f$ has compact support.
If we can show that $f$ is $C^k$ as a map to a completion~$\wt{E}$,
then $\im(d^{(j)}f)\sub E$ because $\im(\delta^{(j)}f)\sub E$
and $d^{(j)}f(x,\cdot)$ can be recovered from $\delta^{(j)}_xf$ using the Polarization
Formula. As a consequence, $f$ will also be $C^k$ as a map to~$E$. We may therefore assume now
that~$E$ is complete. We let $h\colon \R^n\to \R$ be a smooth function
with $h(\R^n)\sub [0,\infty[$, $\Supp(h)\sub[{-1},1]^n$ and $\int_{\R^n}h(x)\,dx=1$.
Let $h_\ve\colon\R^n\to\R$ be as in Exercise~\ref{exc-convol}(b).
Then $f_\ve:=h_\ve*f\colon\R^n\to E$ is a smooth function such that $(D_y)^jf_\ve=
h_\ve*(D_y^j)f$ for all $j\in\N$ such that $j\leq k$ and $y\in\R^n$,
by Exercise~\ref{exc-convol}(b).
\begin{description}[(D)]
\item[(a)]
Show that
$\frac{\partial^{|\alpha|}f_\ve}{\partial x^\alpha}$
converges in $C(\R^n,E)$ to a continuous function $g_\alpha$ as $\ve\to0$,
for each multi-index $\alpha\in\N_0^n$ such that $|\alpha|\leq k$.\\[.3mm]
[By the Polarization Formula and Exercise \ref{exc-convol}(b),
$\frac{\partial^{|\alpha|}f_\ve}{\partial x^\alpha}$
is a linear combination
of $(D_y)^{|\alpha|}f_\ve=h_\ve*(D_y)^{|\alpha|}f$ for suitable $y\in\R^n$,
with coefficients independent of~$\ve$.]
\item[(b)]
Using Exercise~\ref{exc-the-topf},
deduce from~(a) that $\big(\frac{\partial^{|\alpha|}f_\ve}{\partial x^\alpha}\big)_{\ve>0}$
is a Cauchy net in the complete locally convex space
$C^k(\R^n,E)$.
Deduce that $f=g_0\in C^k(\R^n,E)$.
\end{description}
\end{exer}

\begin{exer}\label{exc-prepa-pfwd}
Let $K$ be a compact topological space; let $(Y,\|\cdot\|_Y)$
and $(F,\|\cdot\|_F)$ be normed spaces. We endow spaces of continuous maps
from $K$ to normed spaces with the supremum norm; we endow
spaces of continuous linear operators between normed spaces with the operator norm.
Show that $\Xi(\gamma)(\eta)(x):=\gamma(x)(\eta(x))$
for $\gamma\in C(K,\cL(Y,F))$, $\eta\in C(K,Y)$ and $x\in K$
defines a continuous map $\Xi(\gamma)(\eta)\in C(K,F)$
which is continuous linear in $\eta\in C(K,Y)$.
Show that the map
\[
\Xi\colon C(K,\cL(Y,F))\to \cL(C(K,Y),C(K,F)),\;\;
\gamma\mto\Xi(\gamma)
\]
so obtained is continuous and linear, with $\|\Xi\|_{\op}\leq 1$.
\end{exer}

\begin{exer}\label{exc-ifFCk}
Let $(E,\|\cdot\|_E)$, $(F,\|\cdot\|_F)$, and $(Z,\|\cdot\|_Z)$
be normed spaces, $U\sub E$ be an open subset,
$P\sub Z$ be a locally convex subset with dense interior,
$k\in\N_0\cup\{\infty\}$ and
$K$ be a compact convex subset with non-empty interior
in a finite-dimensional vector space~$X$.
If $f\colon K\times U\times P\to F$
is an $FC^k$-map, show that
\[
\Phi\colon C(K,U)\times P\to C(K,F), \quad (\gamma,p)\mto (f^p)_*(\gamma)=f(\cdot,p)\circ (\id_K,\gamma)
\]
is $FC^k$. You may assume that $k$ is finite and proceed by induction;
the case $k=0$ is subsumed by Lemma~\ref{pushpar}. Now let $k\in\N$
and assume the assertion holds for $k-1$ in place of~$k$.
The map $\Phi$ is $C^k$ and hence $C^1$, by Proposition~\ref{pfctisCkpar}.
We write $f'(x,y,p):=df(x,y,p,\cdot)\in\cL(X\times E\times Z,F)$
for $(x,y,p)\in K\times U\times P$ and
$\Phi'(\gamma,p):=d\Phi(\gamma,p,\cdot)\in\cL(C(K,E)\times Z,C(K,F))$
for $\gamma\in C(K,U)$ and $p\in P$.
\begin{description}[(D)]
\item[(a)]
Show that the calculation leading to (\ref{thedirrdhere}) also yields the formula
\[
\Phi'(\gamma,p)(\eta,q)(x)=f'(x,\gamma(x),p)(0,\eta(x),q)
\]
for all $\gamma\in C(K,U)$, $p\in P$, $x\in K$, $\eta\in C(K,E)$, and $q\in Z$.
\item[(b)]
For $Y:=X\times E\times Z$, define $\Xi$ as in Exercise~\ref{exc-prepa-pfwd}.
By induction, the map
\[
\Psi\colon C(K,U)\times P\to C(K,\cL(Y,F)),\quad (\gamma,p)\mto
f'(\cdot,p)\circ (\id_K,\gamma)
\]
is $FC^{k-1}$. The map $\lambda\colon Z\to C(K,Z)$ taking $q\in Z$ to the constant function
$\lambda(q)\colon K\to Z$, $x\mto q$ is continuous and linear.
Show that the formula from (a) can be rewritten as
\[
\Phi'= \cL(\Lambda, C(K,F))\circ \Xi\circ \Psi,
\]
where $\Lambda\colon C(K,E)\times Z\to C(K,Y)\cong C(K,X)\times C(K,E)
\times C(K,Z)$, $(\eta,q)\mto (0,\eta,\lambda(q))$ is continuous
linear. Deduce that $\Phi'$ is $FC^{k-1}$
and infer that~$\Phi$~is~$FC^k$.
\end{description}
\end{exer}

\end{small}

\section{Appendix to Chapter~\ref{chapcalcul}}\label{secapp-ch2}
In this appendix to Chapter~\ref{chapcalcul},
we discuss the existence of weak integrals for curves
and suitable mappings on higher-dimensional sets.
\subsection*{Proof of Proposition~\ref{rieman}}
We prove Proposition~\ref{rieman}
on the existence of weak integrals for curves in sequentially complete,
locally convex spaces.
To construct the weak integral of a continuous curve,
we approximate the curve by step functions.
\begin{defn}\label{defnEstep}
Let $E$ be a locally convex space
and $a<b$ in~$\R$.
\begin{description}[(D)]
\item[(a)]
A function $\gamma \colon
[a,b]\to E$
is called an ($E$-valued) \emph{step function} \index{step function} 
if there exist $n\in \N$ and
$a=t_0<t_1<\cdots< t_n=b$
such that $\gamma |_{]t_{k-1},t_k[}$
is constant for all $k\in \{1,\ldots, n\}$.
\item[(b)]
If $\gamma \colon
[a,b]\to E$ is a step function
and $a=t_0<t_1<\cdots< t_n=b$
a partition such that
$\gamma|_{]t_{k-1},t_k[}=\alpha_k$
for suitable $\alpha_k\in E$,
we define
\[
\int_a^b \gamma (t)\, dt\; :=\;
\sum_{k=1}^n\, (t_k-t_{k-1})\, \alpha_k\, \in \, E\,.
\]
\end{description}
\end{defn}
\begin{rem}
(a)
As in the case of real-valued functions,
using that any two partitions
of $[a,b]$ admit a common refinement,
we see that $\int_a^b \gamma (t)\, dt$ is independent
of the choice of $t_0,\ldots, t_k$.\medskip

\noindent
(b)
Using the existence of common refinements,
we see that the set $\cS([a,b],E)$
of all $E$-valued step functions
is a vector subspace of $E^{[a,b]}$,
and that $\cS([a,b],E)\to E$,
$\gamma \mto \int_a^b \gamma (t)\, dt$ is
a linear map.\medskip

\noindent
(c)
For each step function $\gamma \colon [a,b]\to E$
and seminorm $\|\cdot\|_p$ on~$E$, we have
\begin{equation}\label{prepaX}
\Big\| \int_a^b \gamma (t)\, dt\Big\|_p\; \leq\;
(b-a)\, \max\{\|\gamma (t)\|_p \colon t\in [a,b]\}\,.
\end{equation}
To see this, let $t_k$ and $\alpha_k$
be as in Definition~\ref{defnEstep}(b).
We define
$M:=$\linebreak
$\max\{\|\alpha_k\|_p\colon k=1,\ldots, n\}$.
Then
$\|\int_a^b \gamma (t)\, dt\|_p=\|\sum_{k=1}^n(t_k-t_{k-1})\alpha_k\|_p$
$\leq \sum_{k=1}^n(t_k-t_{k-1})\|\alpha_k\|_p
\leq M\,\sum_{k=1}^n(t_k-t_{k-1})
=(b-a)\, M$.
\end{rem}
To see that weak integrals of continuous
curves exist, it is useful to show first that
any such curve is uniformly continuous,
in the following sense:
\begin{defn}
A mapping $f\colon X\to E$ from a metric
space $(X,d)$ to a topological vector space
$E$
is called \emph{uniformly continuous} \index{uniformly continuous map}
if for each $0$-neighborhood $U\sub E$,
there exists $\delta>0$ such that
\[
f(x)-f(y)\, \in \, U\quad \mbox{for all $\, x,y\in X$ of distance $d(x,y)<\delta$.}
\]
\end{defn}
Uniform continuity of curves
now follows from the next lemma.
\begin{lem}\label{unif}
If $f\colon X\to E$ is a continuous
map from a compact metric
space $(X,d)$ to a topological
vector space~$E$, then $f$ is uniformly continuous.
\end{lem}
\begin{prf}
Given a $0$-neighborhood
$U$ in~$E$, let $V\sub E$ be a $0$-neighborhood
such that $V-V\sub U$.
For $x\in X$, by continuity of~$f$
we find $\delta_x>0$
such that
\[
f(y)-f(x)\in V\quad \mbox{for all $y\in B_{\delta_x}(x)$,}
\]
where $B_{\delta_x}(x)$ denotes the open ball of radius
$\delta_x$ around $x$ in the metric space~$X$.
Then $(B_{\frac{\delta_x}{2}}(x))_{x\in X}$
is an open cover of~$X$. The space~$X$ being compact,
we find a finite subset $F\sub X$ such that
$X=\bigcup_{x\in F}B_{\frac{\delta_x}{2}}(x)$.
Define
$\delta:=\min\left\{\frac{\delta_x}{2}\colon x\in F\right\}$.
For $y,z\in X$ such that $d(y,z)<\delta$,
we find $x\in F$ such that $y\in B_{\frac{\delta_x}{2}}(x)$.
Then $d(x,z)\leq d(x,y)+d(y,z)<\frac{\delta_x}{2}+\delta\leq
\delta_x$ shows that also $z\in B_{\delta_x}(x)$
and thus
$f(z)-f(y)=f(z)-f(x)-(f(y)-f(x))\in V-V\sub U$.
We have shown that $f$ is uniformly continuous.
\end{prf}
\noindent
{\em Proof of Proposition}~\ref{rieman}.
We may assume that
$a<b$.
Let $\gamma \colon [a,b]\to E$
be a $C^0$-curve.
For each $n\in \N$,
using the partition
$(a+k2^{-n}(b-a))_{k=0}^{2^n}$,
we define a
step function $\gamma_n\colon [a,b]\to E$
via $\gamma_n(a):=\gamma(a)$ and
\[
\gamma_n(t)\;:=\; \gamma( a+k\, 2^{-n}(b-a))
\]
for $k\in \{0,\ldots, 2^n-1\}$
and $t\in \; ]a+k\frac{b-a}{2^n}, a+(k+1)\frac{b-a}{2^n}]$.
We show that
\[
S_n\; :=\; \int_a^b \gamma_n(t)\, dt\;=\;
\sum_{k=0}^{2^n-1}\frac{b-a}{2^n}\,
\gamma(a+k\, 2^{-n}(b-a))\, \in \, E
\]
is a Cauchy sequence.
To this end, let $U\sub E$ be a $0$-neighborhood.
After shrinking~$U$, we may
assume that $U=\{x\in E\colon \|x\|_p\leq 1\}$
is the closed unit ball
of some continuous seminorm $\|\cdot\|_p$ on~$E$
(see Proposition~\ref{Minkowsk}).
The interval $[a,b]$ being compact and $\gamma$
continuous, Lemma~\ref{unif}
provides $\delta>0$ such that
$\gamma(r)-\gamma(s)\in\frac{1}{b-a}\, U$
for all $r,s\in [a,b]$ such that $|r-s|<\delta$.
Then
\begin{equation}\label{forall}
\|\gamma(r)-\gamma(s)\|_p\; \leq \; (b-a)^{-1}
\quad
\mbox{whenever $|r-s|<\delta$.}
\end{equation}
Let $N\in \N$ such that $2^{-N}\cdot (b-a)<\delta$.
Then
\begin{equation}\label{forall2}
\max \{\|\gamma_n(t)-\gamma_m(t)\|_p\colon t\in [a,b]\}\;\leq\;
(b-a)^{-1}
\end{equation}
for all $n,m> N$. In fact,
by construction of $\gamma_n$,
for each $t\in [a,b]$,
we have $\gamma_n(t)=\gamma(v)$
for some $v\in [a,b]$ such that
$|t-v|\leq 2^{-n}(b-a)$.
Likewise, $\gamma_m(t)=\gamma(w)$
for some $w\in [a,b]$ such that $|t-w|\leq 2^{-m}(b-a)$.
Then $|v-w|\leq |v-t|+|t-w|\leq (2^{-n}+2^{-m})(b-a)
\leq 2^{-N}(b-a)<\delta$,
whence $\|\gamma_n(t)-\gamma_m(t)\|_p=
\|\gamma(v)-\gamma(w)\|_p\leq (b-a)^{-1}$, by (\ref{forall}).
Thus (\ref{forall2}) holds.

Combining (\ref{forall2}) and (\ref{prepaX}),
we see that $\|S_n-S_m\|_p\leq 1$
for all ${n,m\geq N}$ and hence
$S_m-S_n\in U$.
Thus $(S_n)_{n\in \N}$ is a Cauchy sequence
in~$E$ and hence convergent,
by sequential completeness.
Let $S:=\lim_{n\to\infty} S_n$.
If $\lambda\in E'$, then
$\lambda(S_n)=
\sum_{k=0}^{2^n-1}\frac{b-a}{2^n}\,
(\lambda\circ\gamma) (a+k\, 2^{-m}(b-a))\in \K$
is a Riemann sum for
the continuous function $\lambda\circ \gamma\colon
[a,b]\to \K$. Thus
$\lambda(S)= \lim_{n\to\infty} \lambda(S_n)=
\int_a^b\lambda(\gamma(t))\,dt$,
as is well known from single variable calculus.
Hence $S$ satisfies the defining property
of the weak integral $\int_a^b\gamma(t)\, dt$.\qed
\subsection*{Vector-valued
integrals over higher-dimensional sets}
We now discuss weak integrals
over $[0,1]^n$, simplices
and more general sets, for later
use.
We find it convenient to
work in the framework
of Lebesgue integrals.
Readers lacking the necessary background
may replace $X$ by a Jordan measurable
subset of~$\R^n$, assume that $f$ is continuous,
and replace Lebesgue integrals by Riemann
integrals. This suffices for our actual
applications.
\begin{defn}\label{def-wint-gen}
Let $(X,\Sigma,\mu)$ be a measure
space, $E$ be a locally convex space
and $f \colon X\to F$
be a map such that $\lambda\circ f\colon X\to\K$
is Lebesgue integrable over~$X$
with respect to the measure~$\mu$, for each $\lambda\in E'$.
If there exists $z\in E$ such that $\lambda(z)=\int_X \lambda\circ f\,
d\mu$ for each $\lambda\in E'$,
then we call $z$ the \emph{weak integral} \index{weak integral} 
of~$f$ over~$X$
and write $\int_X f\, d\mu:=z$.
If $X$ is a Borel subset of~$\R^n$
and $\mu$ the restriction
of Lebesgue measure to~$X$,
we simply write $\int_X f(x)\,dx:=\int_X f\, d\mu$.
\end{defn}
Since $E'$ separates points on~$E$, the weak integral
is unique (if it exists).
In the remainder of this section, we endow $\R^n$ and its subsets with the metric $(x,y)\mto\|y-x\|_\infty$
arising from the maximum norm.
We endow the product $X\times Y$ of metric spaces $(X,d_X)$ and $(Y,d_Y)$
with the maximum metric, $d((x_1,y_1),(x_2,y_2)):=\max\{d_X(x_1,x_2),d_Y(y_1,y_2)\}$
for $(x_1,y_1),(x_2,y_2)\in X\times Y$.
\begin{prop}\label{intsquare}
Let $E$ be a Mackey complete locally convex space
and $0<a<b<\infty$.
Then $\int_{[a,b]^n}f(x)\, dx$
exists in~$E$, for each $n\in \N$ and each
Lipschitz map $f\colon [a,b]^n\to E$.
\end{prop}
For example, $f\colon [0,1]^n\to E$ might be any $C^1$-map (see Lemma~\ref{cplipaut}).
\begin{prf}
The proof is by induction.
The case $n=1$ holds by definition of Mackey completeness.
Now let $n\geq 2$.
For each $y\in [a,b]^{n-1}$,
the curve
\[ [a,b]\to E, \quad t\mto f(y,t)\]  is Lipschitz, whence
the weak integral
$g(y):=\int_{[a,b]}f(y,t)\, dt$ exists.
By Lemma~\ref{lippardepi}, the map
$g\colon [a,b]^{n-1}\to E$
is Lipschitz.
Hence
\[ z:=\int_{[a,b]^{n-1}}g(y)\, dy \] 
exists in~$E$, by induction.
For each $\lambda\in E'$, we have
$\lambda(z)=\int_{[a,b]^{n-1}}\lambda(g(y))\, dy=
\int_{[a,b]^{n-1}}\int_{[a,b]}\lambda(f(y,t))\, dt\,dy=\int_{[a,b]^n}
\lambda(f(x))\, dx$, using Fubini's Theorem
for the last equality.
Thus $z=\int_{[a,b]^n}f(x)\, dx$.
\end{prf}
\begin{cor}\label{intcompsup}
Let $U\sub \R^n$ be open,
$E$ be a Mackey complete
locally convex space
and $f\colon U\to E$ be
a $C^1$-map whose support
\[ \Supp(f):=\wb{\{x\in U\colon
f(x)\not=0\}} \] 
is a compact subset of $U$.
Then the weak integral $\int_Uf(x)\,dx$ exists in~$E$.
\end{cor}
\begin{prf}
Since $f$ can be extended by~$0$
to a $C^1$-map on all of~$\R^n$,
we may assume that $U=\R^n$.
There is $r>0$ such that $\Supp(f)\sub [{-r},r]^n$.
Then $z:=\int_{[{-r},r]^n}f(x)\,dx$
exists in~$E$ by Proposition~\ref{intsquare}.
Applying continuous linear functionals,
we see that~$z$
satisfies the defining property of the weak integral
$\int_Uf(x)\,dx$.
\end{prf}
\begin{prop}\label{intsimplex}
Let $E$ be a Mackey complete locally convex space.
Then the weak integral $\int_{\Delta_n}f(x)\, dx$
exists in~$E$, for each $n\in \N$ and
Lipschitz map $f\colon \!\Delta_n \!\to E$
on the simplex $\Delta_n\!:=\!\{(x_1,\ldots, x_n)
\in [0,\infty[^n\colon \!\sum_{k=1}^n x_k\leq 1\}$.
\end{prop}
For example, $f\colon \Delta_n\to E$ might be any $C^1$-map (see Lemma~\ref{cplipaut}).
\begin{prf}
To reduce to Proposition~\ref{intsquare},
we use the surjective smooth map
\[
\phi\colon [0,1]^n\to\Delta_n ,\;\;
\phi(t_1,\ldots, t_n):=(t_1,(1-t_1)t_2,\ldots,
(1-t_1)\cdots (1-t_{n-1})t_n)
\]
with Jacobian
\[
\delta(t)\; :=\;
\dt \phi'(t)=(1-t_1)^{n-1}(1-t_2)^{n-2}\cdots (1-t_{n-1})
\]
for $t=(t_1,\ldots, t_n)\in [0,1]^n$.
By Lemma~\ref{cplipaut}, $\phi$ is Lipschitz.
Then $\delta(t)>0$
for $t\in\; ]0,1[^n$,
and $\phi$ takes
$]0,1[^n$ diffeomorphically
onto $\Delta_n^0$.
Given $f$, we now consider the
map
\[
g\colon [0,1]^n\to E\,,\quad
g(t)\; :=\; \delta(t)\,
f(\phi(t)),
\]
which is Lipschitz as it is the product of two bounded Lipschitz maps (see Exercise~\ref{excprodlip}).
By Proposition~\ref{intsquare},
the weak integral $z:=\int_{[0,1]^n}g(t)\,dt$
exists in~$E$.
Then
\begin{eqnarray*}
\lambda(z) &= &\int_{[0,1]^n}\delta(t)
(\lambda\circ f)(\phi(t))\,dt\; =\;
\int_{]0,1[^n}\delta(t) (\lambda\circ f)(\phi(t))\,dt\\
&=& \int_{\Delta_n^0}(\lambda\circ f)(x)\, dx
\; =\; \int_{\Delta_n}(\lambda\circ f)(x)\, dx
\end{eqnarray*}
for each $\lambda\in E'$, by Transformation
of Integrals.
Thus $z=\int_{\Delta_n}f(x)\, dx$.
\end{prf}
\begin{small}
\section*{Exercises for Section~\ref{secapp-ch2}}

\begin{exer}\label{excprodlip}
Let $(X,d)$ be a metric space, $E_1$, $E_2$, and $F$ be locally convex spaces,
$\beta\colon E_1\times E_2\to F$
be a continuous bilinear map and $f_1\colon X\to E_1$ as well as $f_2\colon X\to E_2$
be Lipschitz maps such that $f_1(X)\sub E_1$ and $f_2(X)\sub E_2$ are bounded.
Prove that $\beta\circ (f_1,f_2)\colon X\to F$ is Lipschitz. Show by example that the conclusion can go wrong
if $f_1$ or $f_2$ is unbounded.
\end{exer}

\end{small}
\section{Notes and comments on Chapter~\ref{chapcalcul}}\label{notes-ch2}

\textbf{Differential calculus of smooth and {\boldmath$C^k$}-maps.}
The approach to differential calculus
in locally convex spaces
presented here
goes back to Andr\'{e}e Bastiani~\cite{Ba64}.
Together with the
so-called
``Convenient
Differential Calculus''
developed by A. Fr\"{o}licher, A. Kriegl and P.\,W.
Michor (see \cite{FK88} and \cite{KM97}),
it is one of the main
approaches to infinite-dimensional
differential calculus.
It was widely popularized in the 1980s
by Milnor's survey on infinite-dimensional
Lie groups~\cite{Mil84}.
The $C^k$-maps
we consider are also
known in the literature
as ``Keller's $C^k_c$-maps,''
because they were called $C^k_c$-maps
in Keller's comparative study of differential
calculi in locally convex spaces~\cite{Kr74}.
Expositions of the differential calculus of $C^k_c$-maps
(on open subsets)
can be found in
\cite{Mr80},
\cite{BCR81}, \cite{Ham82},
\cite{Mil82}, \cite{Mil84},
\cite{Ne01a} and \cite{Gl02c},
in varying generality.
In contrast to these texts,
Keller \cite{Kr74}
based his investigations on
convergence structures,
rather than topological spaces.
Our formulation of Fa\`{a} di Bruno's Formula (Theorem~\ref{faadk})
was inspired by~\cite{ClH12}.
See \cite[Lemma~2.3]{GN17} for a proof of Exercise~\ref{exclininvctssmoo},
which generalizes \cite[Thm.~5.3.1]{Ham82}.

\textbf{Convenient differential calculus.}
Smoothness
in the sense of convenient differential
calculus is a weakened
smoothness property,
which is also easy to work with
and easy to verify.
In this context, one calls
a map $f\colon E\to F$ between
locally convex spaces ``smooth'' if
$f\circ \gamma \colon \R\to F$
is smooth for each smooth
curve $\gamma\colon \R\to E$.
While every smooth map in our sense
is also smooth in the convenient sense,
the converse is not true:
A conveniently smooth map
need not even be continuous.
For mappings on Fr\'{e}chet spaces,
smoothness and convenient smoothness
are equivalent
(cf.\ \cite[Thms.\ 12.8 and 4.11(1)]{KM97}),
and so the techniques
of convenient calculus
can also be used in our setting
in this special case.

\textbf{Completeness properties.}
The observation that completeness
properties are unnecessary for the basic
development of differential calculus
has been made in \cite[Appendix]{Kr74}
and~\cite{Gl02c}.
This is important for quotient constructions,
because the quotient space~$E/F$
of a complete locally convex space~$E$
by a closed vector subspace $F\sub E$
need not be complete (nor Mackey complete),
unless~$E$ is metrizable.
An example for this pathology (based on~\cite{SN81}) is given in Example~\ref{bad-quot};
see also \cite[\S\,31.6]{Koe69}
and \cite[Ex., p.\,22]{KM97}.
As already mentioned,
Mackey completeness is essential
for convenient differential calculus,
and refined results and a wealth of information
concerning this concept 
can be found in~\cite{KM97}.
Mackey completeness of the range space also simplifies the theory
of complex analytic functions, as we recall in Section~\ref{seccxan}.

\textbf{Classical differential calculus in Banach spaces.}
For mappings between Banach spaces,
there also is the classical
concept of continuously
Fr\'{e}chet differentiable maps
(see, e.g., \cite{Di60} and~\cite{CaH67}).
Additional results
from finite-dimensional
calculus are available for such maps
which become false beyond Banach spaces
(like the Inverse and Implicit
Function Theorems, or uniqueness
and existence of solutions
to ordinary differential equations).
Recall that a continuous map $f\colon E\supseteq U\to F$
from an open subset of a Banach space~$E$
to a Banach space~$F$
is called \emph{continuously Fr\'{e}chet differentiable} 
\index{continuously Fr\'{e}chet differentiable, $FC^1$} 
(or $FC^1$)
if it is totally differentiable
at each $x\in U$
(i.e., there exists a continuous
linear map $f'(x)\in \cL(E,F)$
such that
\[
f(x+y)-f(x)-f'(x)(y)=o(\|y\|)
\]
holds),
and the map $f'\colon U\to \cL(E,F)$
is continuous, where $\cL(E,F)$ is
equipped with the operator norm.
Inductively, $f$ is called $FC^k$ if it is $FC^1$
and $f'\colon U\to \cL(E,F)$
is $FC^{k-1}$
(see Section~\ref{provis-Ban} for more details).
It can be shown that every $FC^k$-map is $C^k$,
and that every $C^{k+1}$-map
between open subsets of Banach spaces is $FC^k$
(see Proposition~\ref{companotions},
or \cite{Kr74},
\cite{Mil82}, \cite[Appendix]{Gl07e},
\cite[Appendix A.3]{Wa12},
and a sketch in~\cite{Ne01a}).
This gap of one order of differentiability
does not matter much for practical
purposes.
In particular, a map between
open subsets of Banach spaces is $C^\infty$
if and only if it is $FC^\infty$.

\textbf{Differential calculus on non-open subsets.}
We found it essential to acquaint
the reader early on with
analysis on non-open subsets.
We believe that
the class
of locally convex sets with dense interior
is an ideal choice of domains for $C^k$-maps,
which is both general enough (for Lie-theoretic applications)
and facilitates
a simple and convincing theory.
In many cases, such $C^k$-maps simply are restrictions
of $C^k$-maps on open subsets of vector spaces
(as traditionally considered in differential geometry),
notably in the case of finite-dimensional domains.
In fact,
for any locally convex space~$F$, open subset $U\sub\R^n$
and $k\in\N_0\cup\{\infty\}$, any $C^k$-map
$f\colon U\cap [0,\infty[^n\to F$
extends to an $F$-valued $C^k$-map on~$U$,
and actually there is a continuous linear extension operator
$C^k(U\cap [0,\infty[^n,F)\to C^k(U,F)$
(see \cite{Hn23} and \cite{Gl20b},
based on ideas of \cite{See64}; cf.\ also \cite{Mit61}).
For a compact convex subset $K\sub\R^n$ with non-empty interior,
$C^k(K,\R)$ can be identified with the space $\cE^k(K)$
of Whitney jets of class $C^k$, which admit $C^k$-extensions to $\R^n$ by Whitney's
Extension Theorem~\cite{Wh34}.
For $k=\infty$, the existence of a continuous linear extension operator $C^\infty(K,\R)\cong
\cE^\infty(K)\to C^\infty(\R^n,\R)$
is classical (cf.\ \cite[p.\,181]{Stn70},
cited from \cite[p.\,123]{Frk07}). 
For finite~$k$, the continuous linear restriction map $\rho\colon C^k(\R^n,\R)\to C^k(K,\R)$
is surjective as just recalled
and hence a quotient map, by the Open Mapping Theorem for Fr\'{e}chet spaces.
Thus \cite{Fef07} provides a continuous linear extension operator
$C^k(K,\R)\to C^k(\R^n,\R)$. For any $k\in \N_0\cup\{\infty\}$
and any sequentially complete locally convex space~$F$,
a continuous linear extension operator $C^k(K,F)\to C^k(\R^n,F)$
can now be obtained from the preceeding ones using $\ve$-tensor products
(see \cite{Gl20b}).
We mention related results:
For general compact subsets $K\sub\R^n$,
the existence of continuous linear operators $\cE^\infty(K)\to C^\infty(\R^n,\R)$
can be characterized (see \cite{Tn79} and \cite{Frk07}).
Easier characterizations are available
for the existence of tame linear extension operators \cite{FJW11},
and linear extension operators such that $C^k$-norms of extensions
can be estimated in terms of $C^k$-norms of jets
(see \cite{FJW16}).
Extension operators on spaces of Whitney jets are championed
in \cite{Mi20}, with a view towards manifolds of mappings between finite-dimensional manifolds.
Extension results for mappings on subsets of normed spaces
can be found in~\cite{MO92};
parameters in locally convex spaces
are covered by~\cite{Hn23}.
Compare also
\cite{Wo06} and
\cite{RS18} for related results.

\textbf{{\boldmath$C^{k,\ell}$}-functions and exponential laws.}
The differential calculus of $C^{k,\ell}$-maps was developed in~\cite{Alz13},
\cite{Alz19},
and~\cite{AlS15}; the first two works even discuss
$C^\alpha$-mappings $f\colon U_1\times\cdots\times U_n\to F$ on products of subsets $U_1,\ldots, U_n$
of locally convex spaces $E_1,\ldots, E_n$,
where $\alpha=(\alpha_1,\ldots,\alpha_n)$
is an $n$-tuple whose entry $\alpha_j\in\N_0\cup\{\infty\}$
indicates the order of differentiability in the $j$th variable.
The exponential law for $C^{k,\ell}$-maps (Theorem~\ref{explawCkell} and Remark~\ref{defversexp})
and its special case of $C^\infty$-maps (Corollary~\ref{cinfty-explaw})
can be found in~\cite{Alz13}, \cite{Alz19}, and~\cite{AlS15},
with $k$-spaces in place of $k_\R$-spaces.
More limited versions of the exponential law
were available earlier (see, e.g., \cite[\S12]{Gl04b}).
Inspiration also came from
\cite{Th96}, where the product topology $\cO$ on products
(like those in the definition of $C^k$-maps)
is replaced by the coarsest $k$-space topology on the product which is finer than~$\cO$
(as in Exercise~\ref{exc-Kelleyfic});
cf.\ also~\cite{Sei72}.
Further aspects of $C^\alpha$-maps are treated in~\cite{GS22}.

In convenient differential calculus,
a different topology is used on spaces $C^\infty(U,F)$
of smooth maps (in the convenient sense) on an open (or more general) subset $U$
of a Mackey complete locally convex space~$E$,
namely the initial topology with respect to the linear
maps
\[
C^\infty(\gamma,F)\colon C^\infty(U,F)\to C^\infty(\R,F),\quad f\mto f\circ\gamma
\]
to the space of smooth curves in~$F$ (endowed with the compact-open $C^\infty$-topology),
for $\gamma\colon \R\to U$ ranging through the set of all smooth curves.
Then an exponential law for smooth functions always holds (without
extra conditions on the Mackey complete locally convex spaces involved)~\cite{KM97},
but only in the sense that the map
\[
\Phi\colon C^\infty(U\times V,F)\to C^\infty(U,C^\infty(V,F)),\quad f\mto f^\vee
\]
is an isomorphism of bornological vector spaces (i.e., both $\Phi$ and $\Phi^{-1}$
take bounded sets to bounded sets).
Already for $E_1=E_2=U=V=\R$, the map $\Phi$
is discontinuous for the topologies used in the convenient
setting~\cite{Bil02}.

Beyond the cases of smooth or $C^{k,\ell}$-functions,
exponential laws are available for some other relevant spaces of differentiable functions
(e.g., for Schwartz spaces of vector-valued rapidly decreasing
smooth functions), see \cite{KMR16} and \cite{Nik15}.

\textbf{Further generalizations of differential calculus.}
Let~$E$ and $F$ be
topological $\K$-vector spaces
over a non-discrete, Hausdorff topological field~$\K$
and $U\sub E$ be an open subset.
Following~\cite{BGN04},
a map $f\colon U\to F$ is called $C^1$
if it is continuous and there exists
a continuous extension
$f^{[1]}\colon U^{[1]}\to F$
of the directional difference quotient map,
as described in Lemma~\ref{linkBGN}.
Inductively,
$f$ is called $C^k$ if $f$ is~$C^1$
and $f^{[1]}$ is $C^{k-1}$.
Lemma~\ref{linkBGN} ensures
that a map between open subsets
of locally convex spaces over $\K\in\{\R,\C\}$
is $C^1$ in the sense of \cite{BGN04}
if and only if it is $C^1$
in the sense of Definition~\ref{defnC1map},
and also the $C^k$-properties are equivalent (see \cite{BGN04}).
Surprisingly large parts of differential calculus
remain intact for $C^k$-maps over topological fields,
and hence also for $C^k$-maps between
real or complex topological vector spaces which are not
necessarily locally convex.
However, difference calculus can be much more cumbersome than differential
calculus. We therefore use difference calculus only
occasionally, when it is advantageous;
an example is the short proof of the Chain Rule
(Proposition~\ref{chainC1}) with the help of~$f^{[1]}$.
See~\cite{Ber08}, \cite{Gl04b}, \cite{Gl06a}, and
\cite{Gl07d}
for more information on this
general approach to differential calculus.
Compare also
earlier
literature on non-archimedean analysis
(like~\cite{Sf84} and the references therein),
where mostly functions of a single scalar variable were considered.

It is also possible to consider vector-valued $C^k$-functions $f\colon U\to F$
on an open subset $U$ of topological group~$G$ (instead of a locally convex space~$E$),
and corresponding $C^{k,\ell}$-maps.
Directional derivatives at $x\in U$ in a direction $y\in E$ are replaced with derivatives
\[
df(x,\gamma):=(D_\gamma f)(x):=\frac{d}{dt}\Big|_{t=0}f(x\gamma(t))
\]
along continuous one-parameter groups $\gamma\colon\R\to G$ here,
and again exponential laws are available under suitable hypotheses (see
\cite{BCR81}, \cite{BNi15}, \cite{Nik15} and the seminal paper
\cite{Rs53}).

\chapter[Further techniques  of differential calculus]{Further techniques of differential calculus}\label{chapfurcalc}

In this chapter, we discuss further concepts and tools of infinite-dimensional calculus
which are relevant for specific examples and applications,
but may be skipped in a first course on infinite-dimensional Lie groups.
\section{Complex analytic mappings}\label{seccxan}
In this section, we shall
refer to $C^k$-maps
over $\K\in \{\R,\C\}$
as $C^k_\K$-maps,
for added clarity.
Given complex locally convex spaces $E$ and $F$,
we show that a map
$f\colon U\to F$ on an open subset $U\sub E$
is $C^\infty_\C$
if and only if it is complex analytic,
i.e., a continuous map
which can be written
close to each point
as a series of continuous
homogeneous polynomials.\\[3mm]
It is useful to
discuss mappings on open subsets
of~$\C$ first and pass to
domains $U\sub E$ only in a second step.
As in our earlier discussions of curves
on real intervals,
a simpler
description of the
$C^k_\C$-property is available in this case.
\begin{lem}\label{singlevarcx}
Let $f\colon U\to E$ be a map from
an open subset $U\sub \C$ to a complex locally
convex space~$E$, and $k\in\N\cup\{\infty\}$.
Then $f$ is a $C^k_\C$-map if and only if $f$
admits the complex derivatives
\begin{equation}\label{defcxdersing}
f^{(1)}(z)\, :=\, f'(z)\, :=\, \lim_{w\to 0}
\frac{1}{w}\big( f(z+w)-f(z)\big)
\end{equation}
and $f^{(j)}(z):=(f^{(j-1)})'(z)$
for all $z\in U$ and $j\in\N$ such that $j\leq k$ $($where $f^{(0)}:=f)$,
and $f^{(j)}\colon U\to E$ is continuous for all $j\in\N_0$ with $j\leq k$.
In this~case,
\begin{equation}\label{dervsdiff}
f^{(j)}(z)\;=\; d^{\,(j)}f(z,1,\ldots,1)\;=\;
\delta^j_z f(1)
\end{equation}
and
\begin{equation}\label{dervsdiff2}
d^{\,(j)}f(z, w_1,\ldots, w_j)
\;=\; w_1\cdot\ldots\cdot w_j\cdot f^{(j)}(z)\,,
\end{equation}
for all $j\in \N$ such that $j\leq k$,
$z\in U$, and $w_1,\ldots, w_j\in\C$.
\end{lem}
\begin{prf}
We may assume that $k<\infty$.
If $f^{(j)}$ exists and is continuous for all $j\in \N$ such that
$j\leq k$,
then we see as in the
proof of Lemma~\ref{homog}
that $df(z,w_1)$ exists for
all $(z,w_1)\in U\times \C $
and is given by
\begin{equation}\label{simplversi}
df(z,w_1)=w_1 df(z,1)=w_1 f'(z),
\end{equation}
which is a continuous $E$-valued function of $(z,w_1)\in U\times \C$.
Hence $f$ is $C^1_\C$.
If $k\geq 2$, then
$f'\colon U\to E$ admits continuous complex derivatives
$(f')^{(j)}=f^{(j+1)}$ for all $j\in\N$ such that $j\leq k-1$,
whence $f'\colon U\to E$ is a $C^{k-1}_\C$-map by the inductive hypothesis,
with
\begin{equation}\label{simplverso}
d(f')(z,w_2,\ldots,w_k)=w_2\cdots w_k(f')^{(k-1)}(z)
\end{equation}
for $z\in U$ and $w_2,\ldots, w_k\in\C$.
Now (\ref{simplversi}) shows that $df$ is $C^{k-1}_\C$
and thus $f$ is $C^k_\C$.
Holding $w_1\in\C$ fixed, using (\ref{simplverso}) we can form $k-1$ directional derivatives in the
$z$-variable
in (\ref{simplversi}) and obtain (\ref{dervsdiff2}).

If, conversely, $f$ is a $C^k_\C$-map,
then apparently
$f^{(j)}(z)$ exists and is given by $d^{\,(j)}f(z,1,\ldots, 1)$,
for all $j\in \N$ such that $j\leq k$. It only remains to observe that $d^{\,(j)}f(z,1,\ldots,1)$
is continuous in~$z$.
\end{prf}
Remarks analogous to
Remark~\ref{doublestnd}
apply to $C^1_\C$-maps $f\colon \C\supseteq U\to E$.
\begin{defn}
Let $E$ and $F$ be complex locally convex spaces
and $U\sub E$ be open.
A map $f\colon U\to F$
is called  \emph{complex analytic}\index{complex analytic} 
\index{$\C$-analytic} 
(\emph{$\C$-analytic}, or $C^\omega_\C$)
if it is continuous
and, for each $x\in U$,
there exists a sequence $(p_k)_{k\in \N_0}$
of continuous homogeneous polynomials
$p_k\colon E\to F$ of degree~$k$
such that
\begin{equation}\label{cxexpa}
f(x+y)\;=\; \sum_{k=0}^\infty p_k(y)\quad
\mbox{as a pointwise limit,}
\end{equation}
for all $y$ in a $0$-neighborhood
$Y\sub E$ such that~$x+ Y\sub U$.
\end{defn}
\begin{rem}
Taking $E=\C$ in the preceding definition,
we see that
a continuous map
$f\colon \C\supseteq U\to F$
is complex analytic if and only if,
for each $x\in U$, there exists
a sequence $(a_k)_{k\in \N_0}$
of elements $a_k\in F$
such that\footnote{The continuity of $f$ is actually
automatic here, because
for each continuous seminorm
$\|\cdot\|_p$ on $F$ and $r>1$ such that
$\bD_r\sub Y$, we have
\[ \|f(x+z)-f(x)\|_p\leq |z|\,C\sum_{k=1}^\infty (1/2)^{k-1} \] 
for each $z\in \bD_{r/2}$,
where
$C:=\sup\,\{\|r^ka_k\|_p\colon k\in \N\}<\infty$
because the convergence
of $\sum_{k=1}^\infty r^ka_k$
entails that $\lim_{k\to \infty} r^ka_k=0$.}
\begin{equation}\label{cxexpa2}
f(x+z)\;=\; \sum_{k=0}^\infty z^k\, a_k
\end{equation}
for all $z$ in some $0$-neighborhood $Y\sub \C$
such that $x+Y\sub U$.
\end{rem}
Thus,
a function of a single complex
variable is complex analytic if and only
if it is given locally by a power series.
Notably, a
function $\C\supseteq U\to\C$
is complex analytic if and only if it is holomorphic.
\begin{defn}
Given a $C^1$-curve $\gamma \colon [a,b]\to \C$
and a continuous map $f\colon \gamma([a,b])\to E$
to a complex locally convex space~$E$,
we define
\[
\int_\gamma f(\zeta)\, d\zeta
\; :=\; \int_a^b f(\gamma(t))\, \gamma'(t)\; dt\,,
\]
if the weak integral on
the right hand side exists in~$E$.
If $\gamma\colon [0,2\pi]\to\C$,
$\gamma(t)=a+re^{it}$
where
$a\in \C$ and $r>0$, 
let
\[
\int_{|\zeta -a |=r} f(\zeta)\, d\zeta
:=\int_\gamma f(\zeta)\, d\zeta.
\]
\end{defn}
\begin{rem}
(a)
It is clear that $\int_\gamma f(\zeta)\, d\zeta$
exists if and only if there is a
(necessarily unique)
element $z\in E$ such that $\lambda(z)=\int_\gamma\lambda(f(\zeta))\, d\zeta$
for each $\lambda\in E'$.
In this case, $z=\int_\gamma f(\zeta)\, d\zeta$.\medskip

\noindent
(b)
Suppose that
$\gamma$ is merely
piecewise $C^1$, i.e., a $C^0$-curve
$\gamma\colon [a,b]\to \C$
such that, for suitable $a=t_0<t_1<\cdots<t_n=b$,
the restriction $\gamma|_{[t_{j-1},t_j]}$
is $C^1$ for each $j\in \{1,\ldots, n\}$.
In this case, we define
$\int_\gamma f(\zeta)\, d\zeta$
as the unique element $z\in E$
such that $\lambda(z)=\sum_{j=1}^n \int_{t_{j-1}}^{t_j}
\lambda(f(\gamma(t)))\,\gamma'(t)\,dt$
for each $\lambda\in E'$,
whenever~$z$ exists.
It is easy to see that
this definition is independent
of the choice of~$n$ and $t_0,\ldots, t_n$.
\end{rem}
We shall use the notation
$\bD_r=\{z\in \C\colon |z|\leq r\}$
and $\bD=\bD_1$.
\begin{prop}\label{analytsingle}
Let $U\sub \C$ be an open subset
and $f\colon U\to E$ be a map to
a Mackey complete complex locally convex space~$E$.
Then the following conditions are equivalent:
\begin{description}[(D)]
\item[\rm (a)]
$f$ is $C^\infty_\C$;
\item[\rm (b)]
$f$ is $C^1_\C$;
\item[\rm (c)]
$f$ is $C^1_\R$ and $df(x,\cdot ) \colon \C \to E$
is complex linear, for each $x\in U$;
\item[\rm (d)]
$f$ is continuous and
``weakly analytic,''
i.e.,
for each $\lambda\in E'$,
the map
$\lambda\circ f\colon \C\supseteq
U\to \C$ is complex analytic.
\item[\rm (e)]
$f$ is continuous and,
for all $x\in U$
and $r>0$ with
$x+\bD_r \sub U$,
we have
\begin{equation}\label{Cauchyf}
f(z)\;=\;
\frac{1}{2\pi i}
\int_{|\zeta-x|=r} \frac{f(\zeta)}{\zeta - z }\; d\zeta
\end{equation}
for all $\,z\in \C$ such that $|z-x|<r$.
\item[\rm (f)]
$f$ is continuous and for each $x\in U$ there exists
$r$ as in {\rm (e)} such that {\rm (\ref{Cauchyf})}
holds.
\item[\rm (g)]
$f$ is complex analytic.
\end{description}
\end{prop}
\begin{prf}
The implications (a)$\impl$(b),
(b)$\impl$(c) and (e)$\impl$(f)
are trivial.

(c)$\impl$(d): Being $C^1_\R$, the map $f$ is continuous.
For each $\lambda\in E'$, the composition
$\lambda\circ f\colon U\to \C$
is $C^1_\R$;
its real differential is
$(\lambda\circ f)'(x)=\lambda\circ f'(x)
\colon \C\to \C$ by the Chain Rule
and hence complex linear,
the map $f'(x)$ being complex linear
by hypothesis. Hence $\lambda\circ f$ is
holomorphic and thus complex analytic,
being a $C^1_\R$-map
which satisfies the Cauchy--Riemann equations.

(d)$\impl$(e):
$\lambda\circ f\colon \C\supseteq
U\to \C$ is holomorphic
for each $\lambda\in E'$,
whence
\[
\lambda(f(z)) = (\lambda\circ f)(z)=
\frac{1}{2\pi i}
\int_{|\zeta-x|=r} \!\!\frac{(\lambda\circ f)(\zeta)}{\zeta - z }\,d\zeta
=
\int_{|\zeta-x|=r} \!\lambda\Big(\frac{1}{2\pi i}\,\frac{f(\zeta)}{\zeta - z }
\Big)
d\zeta
\]
for $x,r$ and $z$ as described in~(e),
by Cauchy's Integral Formula.
Hence $f(z)$ is the weak integral
$\int_{|\zeta-x|=r} \frac{1}{2\pi i}\,\frac{f(\zeta)}{\zeta - z }\,
d\zeta$.

(f)$\impl$(g):
Let $r>0$ be as described in (f),
and $s\in \; ]0,r[$.
For each $z\in \bD_s$,
we have $x+z\in x+\bD_r^0$ and hence
\[
f(x+z)\,=\,
\frac{1}{2\pi i}
\int_{|\zeta-x|=r}\frac{f(\zeta )}{\zeta-(x+z)}\, d\zeta
\,=\,
\frac{1}{2\pi}\int_0^{2\pi}
\frac{f(x+re^{it})}{re^{it}- z}\, r e^{it}\, dt\,,
\]
by hypothesis.
Here
\[
\frac{r e^{it}}{r e^{it}-z}\;=\;
\frac{1}{1-\frac{z}{r e^{it}}}
\;=\; \sum_{n=0}^\infty
\Big(\frac{z}{re^{it}}\Big)^n ,
\]
as $|\frac{z}{re^{it}}|\leq \frac{s}{r}<1$.
We claim that
\[
\sum_{k=0}^n f(x+re^{it})
\Big(\frac{z}{ re^{it}}\Big)^k
\;\to\;
\frac{f(x+re^{it})}{re^{it}-z}\, re^{it}
\]
uniformly in $t\in [0,2\pi]$ as $n\to \infty$.
If this is true,
Lemma~\ref{lemcurvesconv}
shows~that
\begin{eqnarray*}
f(x+z) &= &
\frac{1}{2\pi}\int_0^{2\pi}
\frac{f(x+re^{it})}{re^{it}-z}\, re^{it}\, dt
\; =\;
\sum_{n=0}^\infty
z^n\,
\frac{1}{2\pi}\int_0^{2\pi}
\frac{f(x+r e^{it})}{(re^{it})^n }\, dt\\
&=&
\sum_{n=0}^\infty
z^n\, a_n\vspace{-1mm}
\end{eqnarray*}
with $a_n:=\frac{1}{2\pi i}\int_{|\zeta-x|=r}
\frac{f(\zeta)}{(\zeta -x)^{n+1} }\, d\zeta$,
and hence~$f$ is complex analytic, as\linebreak
required. To prove the claim, let~$p$
be a continuous seminorm on~$E$.
Then
\[
C\; :=\; \sup\big\{ \|f(x+re^{it})\|_p\colon t\in [0,2\pi]\big\}
\; <\; \infty\,,
\]
since $[0,2\pi]\to \R$, $t\mto
\|f(x+re^{it})\|_p$ is a continuous function on
a compact set.
Fix~$z$,
and abbreviate $K:=\frac{C}{1-\frac{s}{r}}$.
For each $n\in \N$, we obtain
\begin{eqnarray*}
\!\!\!\!\! \lefteqn{\Big\|\frac{f(x+re^{it})}{re^{it}-z}\, re^{it}-
\sum_{k=0}^n f(x+r e^{it})
{\textstyle \big(\frac{z}{ re^{it}}\big)^k}\Big\|_p}\qquad\\
&=&
\, \left\|\sum_{k=n+1}^\infty f(x+ re^{it})
{\textstyle \big(\frac{z}{ re^{it}}\big)^k}\right\|_p
\, \leq\,  C \sum_{k=n+1}^\infty (s/r)^k \,=\,  K \,(s/r)^{n+1}
\end{eqnarray*}
for each $t\in [0,2\pi]$,
whence also the supremum over all~$t$
is $\leq K \, (s/r)^{n+1}$.
The desired uniform convergence is established.

(g)$\impl$(e):
If $f$ is complex analytic,
given $x\in U$
there exist elements $a_n\in E$
such that
$f(z)=\sum_{n=0}^\infty (z-x)^n \, a_n$
for all $z$ sufficiently close
to~$x$
and hence $\lambda(f(z))=\sum_{n=0}^\infty
(z-x)^n\, \lambda(a_n)$
for each $\lambda\in E'$,
showing that the function $\lambda\circ f\colon \C\supseteq
U\to \C$ is complex analytic
and hence holomorphic.
By Cauchy's Integral Formula,
for each $x\in U$ and
$r>0$ such that $x+\bD_r\sub U$,
we have
\[
\lambda(f(z))\, =\,(\lambda\circ f)(z)\,=\,
\frac{1}{2\pi i}\int_{|\zeta-x|=r}\frac{\lambda(f(\zeta))}{\zeta-z}\,
d\zeta
\]
for each $z\in x+\bD_r$.
As $\lambda\in E'$ was arbitrary,
we see as in ``(d)$\impl$(e)''
that $f(z)$ is the weak integral
$\frac{1}{2\pi i} \int_{|\zeta-x|=r}
\frac{f(\zeta)}{\zeta-z}\, d\zeta$,
for each $z\in x+\bD_r$.
Hence~(e)~holds.

(e)$\impl$(a):
Given $x\in U$, let $r>0$ such that
$x+\bD_r\sub U$.
The map $g\colon \C\to\C$,
$g(w):=x+re^{iw}$
being continuous, the pre-image $W:=g^{-1}(U\setminus (x+\bD_{r/2}))$
is open, and it contains $[0,2\pi]$
by hypothesis.
The map
\[
\phi \colon
(x+\bD_{r/2}^0) \times [0,2\pi]
\to E\,,\quad
\phi(z,t):=\frac{1}{2\pi} \frac{f(x+re^{it})}{re^{it}-z}\, r e^{it}
\]
is of the form
$\phi(t,z)=\beta(h(z,t), c(t))$,
where $\beta\colon \C\times E\to E$,
$\beta(s,v):=sv$ is the scalar multiplication
which is continuous bilinear,
$c\colon [0,2\pi]\to E$,
$c(t):=\frac{f(x+re^{it})re^{it}}{2\pi}$
is continuous, and
$h\colon (x+\bD_{r/2}^0) \times W\to \C$,
$h(z,w):=\frac{1}{re^{iw}-z}$
is $C^\infty_\C$.
Hence, we are in the situation of
Remark~\ref{impspecsit}(c).
If we assume that $E$ is sequentially
complete,
then
all relevant weak
integrals exist automatically,
enabling us to apply Proposition~\ref{diffpar}
to the parameter-dependent integral
$f(z)= \int_0^{2\pi} \phi(z,t)\, dt$.
Hence $f$ is~$C^\infty_\C$ in this case, with
\begin{eqnarray}
\!\!f^{(j)}(z) & = & d^{\, (j)}f(z,1,\ldots,1) \; =\,
\int_0^{2\pi}\beta(d_1^{\,(j)}h(z,t,1,\ldots,1),
c(t))\, dt\nonumber \\
&= &
\frac{j!}{2\pi}
\int_0^{2\pi}\!\!
\frac{f(x+re^{it})re^{it}}{(re^{it}-z)^{j+1}}\, dt
=
\frac{j!}{2\pi i}
\int_{|\zeta-x|=r}\!
\frac{f(\zeta)}{(\zeta -z)^{j+1}}\, d\zeta.\; \label{Cauchys}
\end{eqnarray}
In the general case, if $E$ is merely Mackey complete,
the preceding argument ensures that~$f$
is $C^\infty_\C$ as a mapping to a completion $\wt{E}$,
entailing that
\begin{equation}\label{singleLip}
[0,2\pi]\to \wt{E}\,, \quad
t\mto
\frac{j! f(x+re^{it})re^{it}}{2\pi (re^{it} -z)^{j+1}}
\end{equation}
is $C^\infty_\R$
and hence a Lipschitz curve,
for each $z\in x+\bD_{r/2}$
and $j\in \N$.
But then (\ref{singleLip})
also is a Lipschitz curve in~$E$,
and hence all of the weak integrals
in (\ref{Cauchys}) exist in~$E$,
the latter being Mackey complete.
Now $f$ being a $C^\infty_\C$-map
$\C\supseteq U\to \wt{E}$
with $f^{(j)}(U)\sub E$
for each $j\in \N_0$,
we deduce from (\ref{dervsdiff2}) that $d^{\,(j)}f(U\times \C^j)\sub E$
for each $j\in\N_0$.
Hence~$f$ is also $C^\infty_\C$
as a map to~$E$, by Exercise~\ref{excCkmaptoFo}.
\end{prf}
\begin{rem}
The continuity assumption in
Proposition~\ref{analytsingle}(d)
can be omitted.
Indeed, if $f\colon \C\supseteq U\to E$ is
weakly analytic, $x\in U$ and $K\sub U$
a compact neighborhood of~$x$,
then $\lambda(f(K))\sub \C$
is compact and thus bounded
for each $\lambda\in E'$,
the function $\lambda\circ f\colon U\to \C$
being complex analytic and thus continuous.
Hence $f(K)$ is weakly bounded
and thus bounded in~$E$,
by Mackey's Theorem (Theorem~\ref{Mackey}).
For each $r>0$ such that $x+\bD_r\sub U$
and $z\in \bD_r^0$, by Cauchy's Integral Formula we have
\[
\lambda(f(z))\; =\; \frac{1}{2\pi}
\int_0^{2\pi} \frac{\lambda(f(x+re^{it}))re^{it}}{re^{it}-z}\, dt
\]
for each $\lambda\in E'$,
whence $f(z)=\int_0^{2\pi}
\frac{re^{it}}{2\pi (re^{it}-z)}\,f(x+re^{it})\, dt$
as a weak integral (in the extended sense
of Exercise~\ref{excrefinedwint}(a)).
Here, the first factor of the integrand
is a continuous complex-valued function,
while the second factor
is bounded, since $f(x+\bD_r)$ is bounded
by compactness of $x+\bD_r$
(as explained above).
Exercise~\ref{excrefinedwint}(b)
now shows
that the parameter-dependent
integral $f$ is continuous.
\end{rem}
\begin{cor}\label{corcxcurve1}
Let $U\sub \C$ be open,
$E$ be a complex locally convex
space and $f\colon U\to E$ be a complex
analytic map.
Then $f$ is $C^\infty_\C$, and
the following holds
for each $x\in U$ and $r>0$ such that
$x+\bD_r \sub U$:
\begin{description}[(D)]
\item[\rm (a)]
The coefficients $a_n$ in {\rm (\ref{cxexpa2})}
are uniquely determined;
they are given by $a_n=\frac{1}{n!}f^{(n)}(x)$.
\item[\rm (b)]
We have $f(x+z)=\sum_{n=0}^\infty z^n \frac{f^{(n)}(x)}{n!}$
for each $z\in \bD_r^0 $.
\item[\rm (c)]
Cauchy's Integral Formulas hold:
For each $z\in x+\bD_r^0$ and $n\in \N_0$,
\begin{equation}\label{Cauchys2}
f^{(n)}(z)
\;=\; \frac{n!}{2\pi i}\int_{|\zeta-x|=r}\frac{f(\zeta)}{(\zeta-z)^{n+1}}
\; d\zeta\,.
\end{equation}
As a consequence, for each continuous
seminorm $\|\cdot\|_p$ on~$E$, we have
the ``Cauchy Estimates''
\begin{equation}\label{Cestimate}
\|f^{(n)}(x)\|_p\; \leq \; \frac{n!}{r^n}\,
\max\, \big\{\|f(\zeta)\|_p\colon \mbox{$|\zeta-x|=r$}
\big\}\,.
\end{equation}
\end{description}
\end{cor}
\begin{prf}
(a) Let $Y\sub \C$ be a $0$-neighborhood
such that $x+Y\sub U$ and
$f(x+z)=\sum_{n=0}^\infty z^na_n$
for all $z\in Y$.
Then $\lambda(f(x+z))=\sum_{n=0}^\infty
z^n\lambda(a_n)$, whence
\[
h \colon Y\to \C,
\quad h(z):=\lambda(f(x+z))
\]
is a holomorphic
function and thus $\lambda(a_n)=\frac{h^{(n)}(0)}{n!}$
for each $n\in \N_0$,
as is well known from complex analysis.
If $E$ is Mackey complete,
then $f$ is $C^\infty_\C$ by Proposition~\ref{analytsingle}
and thus $h'(z)=\lambda(f'(x+z))$
and more generally $h^{(n)}(z)=\lambda(f^{(n)}(x+z))$
for each $n\in \N_0$, by induction.
Thus $\lambda(a_n)=\lambda\big(\frac{f^{(n)}(x)}{n!} \big)$\vspace{-.4mm}
for each $\lambda\in E'$,
whence indeed
$a_n=\frac{f^{(n)}(x)}{n!}$.
If $E$ is not Mackey complete,
we consider $f$ as a complex analytic function
to a completion $\wt{E}$.
Then $f$ is $C^\infty_\C$
as a function to $\wt{E}$, with
$f^{(n)}(x) =n!\, a_n\in E$ for each $n\in \N_0$.
Since~$f$ is a $C^\infty_\C$-map $\C\supseteq U\to \wt{E}$
such that $f^{(n)}(U)\sub E$
for each $n\in \N_0$, we find (as at the end of the proof of Proposition~\ref{analytsingle})
that $f$ is also~$C^\infty_\C$
as a map to~$E$.

(b) and (c): As $f(x+z)\in E$ in~(b)
and all summands are in~$E$,
the sequence in question converges
to $f(x+z)$ in~$E$ if and only if it
does so in~$\wt{E}$.
Likewise,
since $f^{(n)}(z)\in E$ in\,(c)
and the integrand takes values in~$E$,
the formulas in\,(c) will hold in~$E$
if we can prove them for $f$ as a map to~$\wt{E}$.
We may therefore assume that~$E$ is complete.
The proof of ``(f)$\impl$(g)''
of Proposition~\ref{analytsingle}
shows that $f(x+z)=\sum_{n=0}^\infty z^na_n$
for each $z\in \bD_r^0$,
where $a_n=
\frac{1}{2\pi i}\int_{|\zeta-x|=r}\frac{f(\zeta)}{(\zeta-x)^{n+1}}\,d\zeta$.
By (a) already established,
we have $a_n=\frac{f^{(n)}(x)}{n!}$
here, whence Cauchy's Integral Formula
holds for $f^{(n)}(x)$.
Applying Lemma~\ref{weakinballb}
to (\ref{Cauchys2}),
we obtain the Cauchy Estimates.
\end{prf}
\begin{cor}\label{corcxcurve2}
Let $U\sub \C$ be open,
$E$ be a complex locally convex
space, and $f\colon U\to E$ be a map.
Then the following conditions are equivalent:
\begin{description}[(D)]
\item[\rm (a)]
$f$ is $C^\infty_\C$;
\item[\rm (b)]
$f$ is $C^\infty_\R$ and $df(x,\cdot )\colon \C \to E$
is complex linear, for each $x\in U$;
\item[\rm (c)]
$f$ is complex analytic.
\end{description}
\end{cor}
\begin{prf}
The implication (a)$\impl$(b)
is trivial, and (c)$\impl$(a)
has been established in Lemma~\ref{corcxcurve1}.

(b)$\impl$(c): By Proposition~\ref{analytsingle},
$f$ is complex analytic as a map
to a\linebreak
completion $\wt{E}$ of~$E$.
Given $x\in U$, choose $r>0$ such that $x+\bD_r\sub U$.
By Corollary~\ref{corcxcurve1},
we have $f(x+z)=\sum_{n=0}^\infty z^na_n$
for all $z\in \bD_r^0$,
where $a_n:=\frac{f^{(n)}(x)}{n!}
=\frac{1}{n!}d^{\,(n)}f(x,1,\ldots,1)\in E$
for each $n\in \N_0$,
since $f$ is $C^\infty_\R$ as a map to~$E$.
Thus $f(x+z)=\sum_{n=0}^\infty z^na_n$ with
$a_n\in E$ for each $n$,
and thus $f$ is complex analytic as a map to~$E$.
\end{prf}
\begin{lem}\label{diffan}
Let $r>0$ and $f\colon \bD^0_r\to E$
be a complex
analytic map to a
complex locally convex
space~$E$. Then also $f'\colon U\to E$ is complex analytic.
If $f(z)=\sum_{n=0}^\infty z^na_n$,
then $f'(z)=\sum_{n=1}^\infty n z^{n-1} a_n$.
\end{lem}
\begin{prf}
Since $f$ is complex
analytic, $f$ is $C^\infty_\C$,
whence also $f'$ is $C^\infty_\C$ and hence complex
analytic. By Corollary~\ref{corcxcurve1}(a),
we have $a_n=\frac{f^{(n)}(0)}{n!}$.
Again by
Corollary~\ref{corcxcurve1}(a),
we have $f'(z)=\sum_{n=0}^\infty
z^n b_n$ with $b_n=\frac{(f')^{(n)}(0)}{n!}=(n+1)\frac{f^{(n+1)}(0)}{(n+1)!}
=(n+1)a_{n+1}$, from which the desired formula follows.
\end{prf}
\begin{defn}\label{defnGato}
Let $E$, $F$ be complex locally convex spaces,
$U\sub E$ be open
and $f\colon U\to F$ be a map.
We call $f$ \emph{G\^{a}teaux analytic}\index{G\^{a}teaux analytic}
if it is
complex analytic along each affine line,
i.e., for each $x\in U$ and $y\in E$,
the function
\[
Z\to F\,, \quad z\mto f(x+zy)
\]
on the
open subset $Z:=\{z\in \C\colon x+zy\in U\}\sub \C$
is complex analytic.
\end{defn}
\begin{thm}\label{charcxcompl}
Let $E$ and $F$ be complex locally convex spaces,
$U\sub E$ be open, and
$f\colon U\to F$ be a map.
Then the following conditions are
equivalent:
\begin{description}[(D)]
\item[\rm (a)]
$f$ is $C^\infty_\C$;
\item[\rm (b)]
$f$ is $C^\infty_\R$ and $f'(x)\colon E\to F$
is complex linear, for each $x\in U$;
\item[\rm (c)]
$f$ is continuous and
G\^{a}teaux analytic;
\item[\rm (d)]
$f$ is complex analytic.
\end{description}
If $F$ is Mackey complete,
then also each of {\rm (e)} and {\rm (f)}
is equivalent to {\rm (a)--(d)}:
\begin{description}[(D)]
\item[\rm (e)]
$f$ is $C^1_\C$;
\item[\rm (f)]
$f$ is $C^1_\R$ and $f'(x)\colon E\to F$
is complex linear, for each $x\in U$.
\end{description}
\end{thm}
\begin{prf}
The implications (a)$\impl$(b),
(a)$\impl$(e) and (e)$\impl$(f)
are trivial.

(b)$\impl$(c):
Given $x\in U$ and $y\in E$, let $Z:=\{z\in \C\colon
x+zy\in U\}$.
The map $Z\to U$, $z\mto x+zy$
being the restriction of a continuous affine-linear
map and hence~$C^\infty_\R$, the composition
$h\colon Z\to F$, $h(z):=f(x+zy)$
is~$C^\infty_\R$,
with differential $h'(z)$ at $z\in Z$ given by
$h'(z)(w)=f'(x+zy)(wy)$ for $w\in \C$,
by the Chain Rule. Since $f'(x+zy)$
is complex linear, the preceding formula
shows that $h'(z)\colon \C\to F$ is complex linear,
and thus $h$ is complex analytic
by (b)$\impl$(c) in Corollary~\ref{corcxcurve2}.
Thus $f$ is continuous and G\^{a}teaux analytic.

(c)$\impl$(a): We show that $f$ is $C^k_\C$
for each $k\in \N$,
by induction.
Given $x_0\in U$,
there is an open neighborhood
$X\sub U$ of~$x_0$ and a balanced $0$-neighborhood
$Y\sub E$ such that $X+2 Y\sub U$.
Given $x\in X$ and $y\in Y$,
the map
$h\colon \bD_2^0\to F$, $h(z):=f(x+zy)$ is complex
analytic by G\^{a}teaux analyticity of~$f$,
entailing that
$df(x,y)=h'(0)$ exists.
Since $h'(0)=\frac{1}{2\pi i}\int_{|\zeta|=1}\frac{h(\zeta)}{\zeta^2}\,d\zeta
=\frac{1}{2\pi i}\int_{|\zeta|=1}\frac{f(x+\zeta y)}{\zeta^2}\, d\zeta$
by Corollary~\ref{corcxcurve1}(c),
we find that
\begin{equation}\label{goodindu}
df(x,y)\;=\; \frac{1}{2\pi i}\int_{|\zeta|=1}\frac{f(x+\zeta y)}{\zeta^2}
\; d\zeta \quad \mbox{for all $\,x\in X$ and $y\in Y$.}
\end{equation}
As in Lemma~\ref{homog}, we see that $df(x,y)$ actually
exists for all $x\in X$ and $y\in E=\bigcup_{s>0}sY$,
and is homogeneous in~$y$.
Thus $df(x,y)$
exists for all $(x,y)\in U\times E$,
as~$x_0$ was arbitrary.
Since~$f$ is continuous, we infer with Lemma~\ref{intpar}
from (\ref{goodindu}) that $df|_{X\times Y}$
is continuous. As $E=\bigcup_{s>0}sY$
and
\begin{equation}\label{expanddomain}
df(x,y)\; =\; s\, df|_{X\times Y}(x,s^{-1}y)\qquad
\mbox{for each $x\in X$ and $y\in sY$,}
\end{equation}
we deduce that $df|_{X\times E}$ is continuous
(and thus $df$, as $x_0$ was arbitrary).
Hence $f$ is $C^1_\C$.
Let $Q\sub \C$ be an open subset
such that $[0,2\pi]\sub Q$
and $|e^{is}|\leq 2$ for all $s\in Q$.
If now $f$ is $C^k_\C$ by induction
for some $k\in \N$,
then (\ref{goodindu})
expresses $df(x,y)$
as a parameter-dependent integral,
$df(x,y)=
\frac{1}{2\pi}\int_0^{2\pi} g(x,y,t)\,dt$,
where $g\colon X\times Y\times Q \to F$,
$g(x,y,s):=e^{-is}f(x+e^{is}y)$
is $C^k_\C$.
In view of Remark~\ref{impspecsit}(b),
Proposition~\ref{diffpar} shows that
$df|_{X\times Y}$
is $C^k_\C$ as a map
to $\wt{F}$,
whence so is $df|_{X\times E}$
by (\ref{expanddomain})
and hence $df$ (as $x_0$ was arbitrary).
Thus $f$ is $C^{k+1}_\C$
as a map to $\wt{F}$.
Since $\delta_x^jf(y)=\frac{d^j}{dz^j}\big|_{z=0}f(x+zy)\in
F$ by complex analyticity of $z\mto f(x+zy)$
as a map into~$F$,
we deduce that $\delta^j_xf(E)\sub F$
for each $x\in U$, $y\in E$
and $j\leq k+1$,
whence also $d^{(j)}f(U\times E^j)\sub F$
by polarization
(see Corollary~\ref{applpolze}).
Hence $f$ is $C^{k+1}_\C$ also as
a map to~$F$, which completes the inductive proof.

(c)$\impl$(d):
By ``(c)$\impl$(a)''
already established,
$f$ is~$C^\infty_\C$.
Given $x\in U$, let
$W\sub E$ be a balanced
open $0$-neighborhood such that $x+W\sub U$.
For each $y\in W$,
the set $Z:=\{z\in \C\colon zy  \in W\}$
is a balanced open $0$-neighborhood
in~$\C$.
By hypothesis, the map $h\colon Z\to F$,
$h(z):=f(x+zy)$ is complex analytic.
As $1\in Z$ and $Z$ is open,
we have $r\in Z$ for some $r>1$
and thus
$\bD_r \sub Z$, since $Z$ is balanced.
Now Corollary~\ref{corcxcurve1}(b)
shows that
\[
h(z)\;=\; \sum_{n=0}^\infty z^n \frac{h^{(n)}(0)}{n!}
\;=\;
\sum_{n=0}^\infty z^n \, \frac{\delta^n_xf(y)}{n!}
\quad \mbox{for each $z\in \bD_r^0\,$.}
\]
Thus $f(x+y)=h(1)=\sum_{n=0}^\infty \frac{\delta^n_xf(y)}{n!}$,
where $\frac{\delta^n_x f}{n!}\colon E\to F$
is a continuous homogeneous polynomial
of degree~$n$ by Lemma~\ref{examppoly}.
Hence~$f$ is complex analytic.

(d)$\impl$(c):
Given $x\in U$ and $y\in E$,
the map
$h\colon Z\to F$, $z\mto f(x+zy)$
on $Z:=\{z\in \C\colon x+zy\in U\}\sub \C$ is continuous.
Given $z\in Z$,
there exists an open $0$-neighborhood
$W\sub E$ such that $(x+zy)+W\sub U$
and continuous
homogeneous polynomials
$p_k\colon E\to F$ of degree~$k$ such that
$f((x+zy)+w)=\sum_{k=0}^\infty p_k(w)$
for each $w\in W$.
There exists an open $0$-neighborhood $Y\sub \C$
such that $Yy\sub W$. 
Then $h(z+s)=f((x+zy)+sy)
=\sum_{k=0}^\infty s^k\,p_k(y)$ for each $s\in Y$.
Hence $h$ is complex analytic
and thus $f$ is G\^{a}teaux analytic.

(f)$\impl$(c): Assume that $f$ is $C^1_\R$ and
each $f'(x)$ is complex linear.
Then $f$ is continuous.
Given $x\in U$ and $y\in E$,
let $h\colon Z\to F$, $h(z):=f(x+zy)$ be as in
Definition~\ref{defnGato}.
Then $h$ is $C^1_\R$ and its differential
$dh(z,\cdot )\colon
\C\to F$, $w\mto df(x+zy ,w y)$
is $\C$-linear for each $z\in Z$,
whence $h$ is complex analytic by
Proposition~\ref{analytsingle}.
Hence $f$ is continuous
and G\^{a}teaux analytic.
\end{prf}
\begin{rem}
Having shown that complex analytic maps
are the same thing as $C^\infty_\C$-maps,
we deduce from Proposition~\ref{chainCk}
that compositions of composable complex
analytic maps are complex analytic.
\end{rem}
Every complex analytic map is given
locally by its Taylor series.
\begin{cor}\label{corcxmap}
Let $f\colon E\supseteq U\to F$ be complex analytic,
and $x\in U$.
Then the following holds:
\begin{description}[(D)]
\item[\rm (a)]
The homogeneous
polynomials $p_k\colon E\to F$ in {\rm (\ref{cxexpa})}
are uniquely determined;
they are given by
$p_k = \delta^{k}_x f/k!\,$,
where $\delta^k_xf$ is the $k$th G\^{a}teaux
differential of~$f$ at~$x$.
\item[\rm (b)]
If $V\sub E$ is
a balanced, open
$0$-neighborhood
such that $x+V\sub U$, then
\[
f(x+z)\;=\; \sum_{k=0}^\infty \frac{\delta_x^kf(z)}{k!}\qquad
\mbox{for all $\, z\in V$.}
\]
\end{description}
\end{cor}
\begin{prf}
Assume that $f(x+z)=\sum_{k=0}^\infty p_k(z)$
for each $z\in Y$,
as in (\ref{cxexpa}).
Given $y\in E$, there exists $r>0$
such that $\bD_r^0\, y\sub Y$.
Then $h\colon \bD_r^0\to F$, $h(s)=f(x+sy)=\sum_{k=0}^\infty
s^k\, p_k(y)$ is complex analytic.
By Corollary~\ref{corcxcurve1}(a), we have $p_k(y)=\frac{h^{(k)}(0)}{k!}=
\frac{\delta^k_xf(y)}{k!}$ for each $k\in \N_0$.
This shows~(a);
assertion~(b) was already established
when we proved ``(c)$\impl$(d)"
in Theorem~\ref{charcxcompl}.
\end{prf}
\begin{ex}\label{invcxcia}
Let $\cA$ be a complex continuous
inverse algebra.
Being $C^\infty_\C$ (by Corollary~\ref{invsmoocia}),
the inversion map $\iota\colon \cA^\times \to \cA^\times$,
$\iota(x):=x^{-1}$
is complex analytic.
Let us expand $\iota$
into its Taylor series around
a point $x\in \cA^\times$.
From Exercise~\ref{excTaycia},
we know that $\delta^k_x\iota(y)=
({-1})^k k! \, x^{-1}\big(yx^{-1}\big)^k$
for each $y\in \cA$.
Now Corollary~\ref{corcxmap} shows that,
for each absolutely
convex, open $0$-neighborhood $\cY\sub A$ such that $x+\cY\sub \cA^\times$,
we have
\[
\iota(x+y)\;=\; \sum_{k=0}^\infty ({-1})^k x^{-1}\big(yx^{-1}\big)^k
\quad\mbox{for all $\,y\in Y$.}
\]
If $x=\one $ in particular, then
$\iota(\one +y)=\sum_{k=0}^\infty ({-1})^k y^k$
and hence
\begin{equation}\label{neumcia}
(\one -y)^{-1}\; =\; \sum_{k=0}^\infty y^k\,.
\end{equation}
Thus $(\one -y)^{-1}$ is given by Neumann's series,
as in the case of Banach algebras already encountered in
Exercise~\ref{excBancia}.
\end{ex}
\begin{thm}[\textbf{Identity Theorem}]\label{identitycx}
Let $E$ and $F$ be complex locally convex spaces,
$U\sub E$ be an open,
connected subset and $f\colon U\to F$
be complex analytic.
Assume that one of the
following conditions
is satisfied:
\begin{description}[(D)]
\item[\rm (a)]
There exists an open, non-empty subset $V\sub U$
such that $f|_V=0$;
\item[\rm (b)]
There exists $x\in U$
such that $\delta_x^kf =0$ for all
$k\in \N_0$; or:
\item[\rm (c)]
There exists a real vector subspace
$H\sub E$ such that $H+iH=E$,
and $x\in U$ such that $f|_{(x+H)\cap U}=0$.
\end{description}
Then $f=0$.
\end{thm}
\begin{prf}
Note first that (a)$\impl$(b)
as $f|_V=0$ entails that
$d^{\,(k)}f(V\times E^k)=0$
for each $k\in \N$ (by an apparent induction),
whence $\delta^k_xf=0$ for all $k\in \N$, $x\in V$.

Also (c)$\impl$(b). In fact, given $x\in U$
such that $f|_{(x+H)\cap U}=0$
and $y\in H$,
there is $r>0$ such that
$x+\bD_r^0\, y\sub U$.
Then $h\colon \bD_r^0\to F$,
$h(z):=f(x+zy)$ is $C^\infty_\C$.
Here $h(t)=0$ for each $t\in \;]{-r},r[$,
by hypothesis,
whence also $h'(t)=\lim_{s\to 0}\frac{h(t+s)-h(t)}{s}=0$
for each
$t\in \;]{-r},r[$,
as we can use a sequence of real
numbers~$s$ to calculate the limit.
Repeating the argument, we find that
$h^{(k)}(t)=0$ for each $k\in \N_0$
and each $t\in \,]{-r},r[$,
whence $\delta_x^kf(y)=h^{(k)}(0)=0$
for each $k\in \N_0$ and each $y\in H$
in particular. Now the Polarization
Formula (\ref{polarformulvar}) shows that $d^{\,(k)}f(x,y_1,\ldots, y_k)=0$
for all $y_1,\ldots, y_k\in H$,
whence $d^{\,(k)}f(x,u_1+iv_1,\ldots,u_k+iv_k)=0$
for all $u_1,\ldots, u_k,v_1,\ldots,v_k\in H$,
using the complex $k$-linearity of $d^{\,(k)}f(x,\cdot)$.
Since $E=H+iH$, we deduce that $d^{\,(k)}f(x,\cdot)=0$
and thus $\delta^k_x f=0$.

Now assume that (b) holds.
Using Exercise~\ref{excCkGatcts}, we see that
\[
W\, :=\, \{x\in U\colon (\forall k\in \N_0)\; \delta^k_xf=0\}
\,=\,
\bigcap_{k\in \N_0}\bigcap_{y\in E} \{x\in U\colon \delta^k_xf(y)=0\}\vspace{-1mm}
\]
is a closed subset of~$U$.
By our hypothesis, $W$ is non-empty.
Since $U$ is connected,
we will have $W=U$ (which finishes the proof)
if we can show that $W$ is also open.
To this end, let $x\in W$;
thus $\delta^k_xf=0$
for each $k\in \N_0$.
There exists a balanced, open $0$-neighborhood
$V\sub E$ such that $x+V\sub U$.
Then $f(x+z)=\sum_{k=0}^\infty \frac{\delta^k_xf(z)}{k!}=0$
for each $z\in V$ (by Corollary~\ref{corcxmap}),
whence $f|_{x+V}=0$ and thus $x+V\sub W$,
by the proof of (a)$\impl$(b).
Hence $W$ is a neighborhood of each of its
points and thus open.
\end{prf}
In the definition of complex
analyticity, we required
the existence of an expansion
into a series of homogeneous
polynomials \emph{around each point}.
In good cases,
it suffices to have such an expansion
around a single point.
\begin{prop}\label{singlpoint}
Let $E$ and $F$
be complex locally convex spaces,
$U\sub E$ be a balanced, open $0$-neighborhood
and $(p_k)_{k\in \N_0}$ be a sequence
of continuous homogeneous polynomials
$p_k\colon E\to F$ of degree~$k$.
If the limit
\[
f(z)\;:=\; \sum_{k=0}^\infty \, p_k(z)
\]
exists for each $z\in U$,
the map $f\colon U\to F$ is continuous
and $F$ is Mackey complete,
then $f\colon U\to F$ is complex analytic.
\end{prop}
\begin{prf}
Since $f$ is continuous,
in view of ``(c)$\impl$(d)'' in Theorem~\ref{charcxcompl}
we only need to show that $f$ is
G\^{a}teaux analytic.
As a consequence of Proposition~\ref{analytsingle}(d),
this is the case if and only if
$\lambda\circ f$ is G\^{a}teaux analytic
for each $\lambda\in F'$,
where $(\lambda\circ f)(z)=\sum_{k=0}^\infty (\lambda\circ p_k)(z)$.
We may therefore assume that~$F=\C$.
We claim that,
for each compact subset $K\sub U$,
\begin{equation}\label{goodconv}
\sum_{k=0}^\infty \, \sup\, \{|p_k(x)|\colon x\in K\}\; < \; \infty\,.
\end{equation}
To see this, choose $r>1$ such that $\bD_r K\sub U$
(which is possible since $K$ is compact and $U$
is open and balanced).
Given $x\in K$,
the map $h\colon \bD_r^0\to \C$,
$h(z):=f(zx)=\sum_{k=0}^\infty z^k p_k(x)$
is holomorphic, with
$p_k(x)=\frac{h^{(k)}(0)}{k!}$.
Then $C:=\sup\,\{|f(y)|\colon y\in \bD_rK \}<\infty$
by compactness of $\bD_r K$ and continuity of~$f$.
Using the Cauchy Estimates~(\ref{Cestimate}),
we obtain
\[
|k!\,p_k(x)|\;=\;
|h^{(k)}(0)| \; \leq \; \frac{k!}{r^k}\,
\sup\,\Big\{|h(z)|\colon z\in \bD_r \Big\}
\; \leq\; \frac{k!}{r^k}\, C\,. 
\]
Hence $|p_k(x)|\leq C/r^k$ for each $x\in K$
and thus
\[ \sum_{k=0}^\infty \, \sup\, \{|p_k(x)|\colon x\in K\}
\leq \sum_{k=0}^\infty C/r^k=\frac{C}{1-\frac{1}{r}}<\infty,\]
whence (\ref{goodconv}) holds.

Given $x\in U$ and $y\in E$,
define $Z:=\{z \in \C\colon x+z y\in U\}$
and
$g \colon Z\to \C$,
$g(z):=f(x+zy)$.
For each $z_0\in Z$, we find $r>0$
such that $z_0+\bD_r\sub Z$.
The set $K:=x+ (z_0+\bD_r)y \sub U$ being
compact, (\ref{goodconv})
entails that the series
$g(z)=\sum_{k=0}^\infty p_k(x+zy)$
converges
uniformly on $z_0+\bD_r$.
Hence $g|_{z_0+\bD_r^0}$
is holomorphic,
being a uniform limit
of the holomorphic (polynomial)
functions $z_0 + \bD_r^0\to\C$,
$z\mto \sum_{k=0}^n p_k(x+zy)$.
As $z_0\in Z$ was arbitrary,
we see that $g$ is holomorphic,
and thus $f$ is G\^{a}teaux holomorphic,
as required.
\end{prf}
\begin{ex}\label{explogban}
Let $\cA$ be a unital complex Banach algebra.
We have seen in Exercise~\ref{excbabyBCH}
that the mappings
$\exp\colon \cA\to \cA$,\vspace{-.6mm}
$\exp(x):=\sum_{k=0}^\infty\frac{1}{k!}x^k$
and $\log \colon  B_1^\cA(\one )\to \cA$,
$\log(x):=\sum_{k=1}^\infty \frac{({-1})^{k+1}}{k}(x-\one)^k$
are continuous.
Hence Proposition~\ref{singlpoint}
shows that
$\exp$ and $\log$
are complex analytic.
Corollary~\ref{corcxmap}
entails that $P_0^n(\exp)(x)=\sum_{k=0}^n\frac{1}{k!}\, x^k$\vspace{-.6mm}
and $P_{\tiny\one}^n(\log)(x)=\sum_{k=1}^n\frac{({-1})^{k+1}}{k}\, x^k$,
for each $n\in \N$.
\end{ex}
\subsection*{Complex analytic functions on direct limits}
We close this section with a more specialized result concerning complex analytic functions
on open subsets of locally convex direct limits.
The following theorem is a slight generalization of the main result in~\cite{Da10}.
It is useful for some constructions of complex Lie groups, and for the proof of their regularity
(see, e.g., \cite{Da11} and \cite{Da14}). 
\begin{thm}[Dahmen's Theorem]\label{thm-dahmen}
Let $E_1\sub E_2\sub\cdots$ be a direct sequence of complex locally convex spaces
whose locally convex direct limit $E:=\bigcup_{n\in\N} E_n$ is Hausdorff.
Let $U_n\sub E_n$ be a convex open subset for $n\in\N$ such that $U_1\sub U_2\sub\cdots$.
Then $U:=\bigcup_{n\in\N}U_n$ is open in~$E$.
If $F$ is a complex locally convex space and $f\colon U\to F$ a function such that
$f_n:=f|_{U_n}\colon U_n\to F$ is complex analytic for each $n\in\N$ and $f(U_n)\sub F$
is bounded, then~$f$ is complex analytic.
\end{thm}
In the next lemma and its proof, we use terminology and notation from Definition~\ref{lip-lcx}.
\begin{lem}\label{dl-beyond-ana}
Let $E_1\sub E_2\sub\cdots$ be a direct sequence of locally convex spaces,
with locally convex direct limit $E:=\bigcup_{n\in\N} E_n$.
Let $U_n\sub E_n$ be a convex open subset for $n\in\N$ such that $U_1\sub U_2\sub\cdots$.
If $F$ is a locally convex space and $f\colon U\to F$ a function such that
$f_n:=f|_{U_n}\colon U_n\to F$ is Lipschitz on $U_n\sub E_n$
for each $n\in\N$,
then~$f$ is continuous.
\end{lem}
\begin{prf}
Given $x\in U$, let us show that $f$ is continuous at~$x$.
We have $x\in U_m$ for some $m\in\N$; after passing to the direct sequence
$E_m\sub E_{m+1}\sub\cdots$, we may assume that $x\in E_1$.
After replacing $U_n$ with $U_n-x$ for each $n\in \N$ and $f$ with
$g\colon U-x\to F$, $g(y):=f(x+y)$, we may assume that $x=0$.
Let~$p$ be a continuous seminorm on~$F$ and $\ve>0$.
For each $n\in\N$,
there exists a continuous seminorm~$q_n$ on~$E_n$ such that $\Lip_{p,q_n}(f_n)<\infty$.
Pick $r_n>0$ with
\[
\sum_{n=1}^\infty r_n\Lip_{p,q_n}(f_n)\,<\,\ve;
\]
for example, take $r_n:=\ve 2^{-n}/(\Lip_{p,q_n}(f_n)+1)$ for $n\in\N$.
Then
\[
V:=\sum_{k\in\N} B^{q_k}_{r_k}(0)\cap 2^{-k}U_k
\]
is an open $0$-neighborhood in~$E$ (see Lemma~\ref{lemnbhdlcxDL}(c)).
Moreover,
\[
\sum_{k=1}^n 2^{-k}U_k\sub
\sum_{k=1}^n 2^{-k}U_n
\sub U_n\quad\mbox{for each $n\in\N$,}
\]
as $U_n$ is convex and $0\in U_n$.
If $y\in V$, then $y=\sum_{k=1}^n z_k$
for some $n\in\N$ and elements $z_k\in B^{q_k}_{r_k}(0)\cap 2^{-k}U_k$
for $k\in\{1,\ldots, n\}$. Abbreviate
\[
y_k:=\sum_{j=1}^k z_j\;\,\mbox{for $k\in\{0,1,\ldots, n\}$;}
\]
thus $y_0=0$ and $y_n=y$. Moreover, $\{y_{k-1},y_k\}\sub U_k$ for all
$k\in\{1,\ldots, n\}$ and
$q_k(y_k-y_{k-1})=q_k(z_k)<r_k$.
Hence
\begin{eqnarray*}
p(f(y)-f(0)) &\leq& \sum_{k=1}^n p(f(y_k)-f(y_{k-1}))
=\sum_{k=1}^n p(f_k(y_k)-f_k(y_{k-1}))\\
&\leq& \sum_{k=1}^n q_k(y_k-y_{k-1})\Lip_{p,q_k}(f_k)
\leq \sum_{k=1}^n r_k\Lip_{p,q_k}(f_k)<\ve.
\end{eqnarray*}
Thus $p(f(y)-f(0))<\ve$ for all $y\in V$, and so
$f$ is continuous at~$0$.
\end{prf}
\begin{lem}\label{excgenesti}
Let $E, F$ be complex locally convex spaces,
$U\sub E$ be open,
$f\colon U\to F$ be a complex
analytic map and $x\in U$.
Let $p$ be a continuous seminorm on~$F$
and $q$ be a continuous seminorm on~$E$
such that $V:=B_r^q(x)\sub U$
and $p(f(V))$ is bounded.
Then
\begin{equation}\label{est-via-CE}
p(\delta_x^kf(y))
\leq \frac{k!}{r^k}\sup \{p(f(v))\colon v\in V\}\;\, \mbox{for all $y\in \wb{B}^q_1(0)$ and $k\in \N_0$.}
\end{equation}
Moreover, $\Lip_{p,q}(f|_{B^q_\rho(x)})<\infty$ for each $\rho\in\,]0,r[$.
\end{lem}
\begin{prf}
Abbreviate $C:=\sup p(f(V))<\infty$.
Given $y\in \wb{B}^q_1(0)$,
we define $Z:=\{z\in\C\colon x+zy\in V\}$. Then $h\colon Z\to F$,
$h(z):=f(x+zy)$ is complex analytic and bounded.
For each $s\in\,]0,r[$, we have $\bD_s\sub Z$.
Using the Cauchy Estimates (\ref{Cestimate})
recorded in Corollary~\ref{corcxcurve1},
for $k\in\N_0$ we find that
\[
p(\delta_x^kf(y))=p(h^{(k)}(0))\leq
\frac{k!}{s^k}C.
\]
Letting $s\to r$, (\ref{est-via-CE}) follows.
Given $\rho\in \,]0,r[$ and $v\in B^q_\rho(x)$,
we have $B^q_{r-\rho}(v)\sub B^q_r(x)=V$. Applying (\ref{est-via-CE})
with $v$ in place of~$x$ and $r-\rho$ in place of $r$, we deduce that
\[
p(df(v,y))=p(\delta^1_vf(y))\leq \frac{C}{r-\rho}.
\]
As a consequence, $\Lip_{p,q}(f|_{B^q_\rho(x)})\leq C/(r-\rho)<\infty$, by
Lemma~\ref{lipviadergen}.
\end{prf}
\noindent
{\em Proof of Theorem}~\ref{thm-dahmen}.
To see that~$f$ is G\^{a}teaux analytic, let $x\in U$ and $y\in E$.
There exists $m\in \N$ such that $x\in U_m$ and $y\in E_m$.
Then $Z_n:=\{z\in\C\colon x+zy\in U_n\}$ is an open subset of~$\C$
for each $n\geq m$, and the map $Z_n\to F$, $z\mto f(x+zy)=f_n(x+zy)$
is complex analytic.
Moreover, $Z:=\{z\in \C\colon x+ty\in U\}$ is an open subset of~$\C$.
The map $h\colon Z\to F$, $z\mto f(x+zy)$ is complex analytic, as
$Z=\bigcup_{n\geq m}Z_n$ and $h|_{Z_n}$ is complex analytic for each
$n\geq m$, as just observed. Thus $f$ is G\^{a}teaux analytic.

By ``(c)$\impl$(d)'' in Theorem~\ref{charcxcompl},
the map~$f$ will be complex analytic if we can show that~$f$ is continuous.
To see that $f$ is continuous at $x\in U$, after passing to $E_m\sub E_{m+1}\sub\cdots$
for some $m$
we may assume that $x\in U_1$. After replacing each $U_n$ with $U_n-x$ and $f$
with $f(x+\cdot)$, we may assume that $x=0$.
Thus each $U_n$ is a convex, open $0$-neighborhood in~$E_n$.
By Exercise~\ref{ac-kernel},
\[
U_n':=\{y\in U_n\colon \bD y\sub U_n\}
\]
is an absolutely convex open $0$-neighborhood in~$E_n$ for each $n\in\N$,
and $U_1'\sub U_2'\sub\cdots$.
If $q_n\colon E_n\to[0,\infty[$ is the Minkowski functional of~$U_n'$,
then $U_n'=B^{q_n}_1(0)$.
We define $W_n:=\frac{1}{2}U_n'=B^{q_n}_{1/2}(0)$.
Then $W_1\sub W_2\sub\cdots$ and $W:=\bigcup_{n\in\N}W_n$
is an open $0$-neighborhood in~$E$.
By Lemma~\ref{excgenesti},
we have
$\Lip_{p,q_n}(f_n|_{W_n})<\infty$
for each continuous seminorm~$p$ on~$F$.
Notably, $f_n|_{W_n}\colon W_n\to F$ is Lipschitz on $W_n\sub E_n$
for each $n\in\N$.
Hence $f|_W$ is continuous, by Lemma~\ref{dl-beyond-ana}.\qed
\begin{small}
\subsection*{Exercises for Section \ref{seccxan}.}

\begin{exer}\label{excfthvec}
Let $E$ be a Mackey complete, complex
locally convex space
and $f\colon U\to E$
be a map on an open subset $U\sub \C$.
Show the following:
\begin{description}[(D)]
\item[(a)]
(Cauchy's Integral Theorem).
If $f$ is complex analytic
and $U$ is star-shaped
around some $z_0\in U$
(or at least simply connected),
then
$\int_\gamma f(\zeta)\, d\zeta = 0$\vspace{-.7mm}
for each closed, piecewise $C^1$ curve $\gamma$ in~$U$.
\item[(b)]
(Morera's Theorem).
If $f$ is continuous
and
$\int_{\partial \Delta} f(\zeta)\, d\zeta=0$
for each (solid) triangle
$\Delta\sub U$,
then $f$ is complex analytic.
\item[(c)]
(Liouville's Theorem).
Every bounded $C^\infty_\C$-map
$f\colon \C\to E$ is constant.
\end{description}
\end{exer}

\begin{exer}\label{excposep}
Let $E$ and $F$ be complex locally convex spaces,
$U\sub E$ be open and $f\colon U\to F$
be a continuous mapping such that $\lambda\circ f\colon U\to\C$
is complex analytic for all $\lambda$ in a subset $\Lambda\sub F'$
which separates points on~$F$ (i.e.,
for all $x\not=y$ in~$F$, there exists $\lambda\in\Lambda$ such that
$\lambda(x)\not=\lambda(y)$). Assume also that~$F$
is sequentially complete.
Using ``(e)$\impl$(g)" from Proposition~\ref{analytsingle} and
``(c)$\impl$(d)" from Theorem~\ref{charcxcompl},
show that $f$ is complex analytic.
\end{exer}

\begin{exer}\label{excmackincompl}
Let $E$ and $F$ be complex locally convex spaces,
$U\sub E$ be open and $f\colon U\to F$
be a map.
Using ``(e)$\impl$(g)" from Proposition~\ref{analytsingle} and
``(c)$\impl$(d)" from Theorem~\ref{charcxcompl},
show that if~$F$ is Mackey complete and~$f$ is complex analytic as a map to a completion $\wt{F}$ of~$F$ with $F\sub \wt{F}$,
then $f$ is also complex analytic as a map to~$F$.
\end{exer}

\begin{exer}\label{excspeccia}
Let $\cA$ be a complex unital continuous inverse algebra;
we assume that
$\cA\not=\{0\}$
(entailing that $\one\not= 0$).
Given $x\in \cA$ we define its
\emph{resolvent set} as \index{resolvent} 
$\rho(x):= \{z\in \C\colon z\one -x\in \cA^\times\}$
and its \emph{spectrum} as \index{spectrum} 
$\sigma(x):=\C\setminus \rho(x)$.
\begin{description}[(D)]
\item[(a)]
Note that
$\rho(x)$ is open
in~$\C$ and $\sigma(x)$ is closed.
\item[(b)]
For each $z\in \C^\times$, we have
\begin{equation}\label{preinvert}
z\one-x\; = \; z\, \Big(\one-\frac{1}{z}\, x\Big).
\end{equation}
Deduce that there exists $r>0$
such that $z\in \rho(x)$
for all $z\in \C$ such that
$|z|>r$. Infer that $\sigma(x)$ is compact.
\item[(c)]
Calculate $(z\one-x)^{-1}$
using (\ref{preinvert})
and deduce that
$\lim_{|z|\to\infty}(z\one-x)^{-1}=0$.
\item[(d)]
Infer from facts already established
that
$R \colon \rho(x)\to \cA$,
$R (z):=(z\one-x)^{-1}$
is complex analytic.
Given $r\geq 0$, abbreviate
$M_r:=\{z\in \C\colon |z|\geq r\}$.
Deduce from (c) that
$\|R(M_r)\|_p\sub \R$
is bounded for each continuous seminorm
$\|\cdot\|_p$ on~$\cA$. Hence
$R(M_r)\sub \cA$ is bounded.
\item[(e)]
Show that $\sigma(x)\not=\emptyset$.
(Hint: Otherwise $R\colon \C\to \cA$
would be a bounded, complex analytic map,
to which Liouville's Theorem from Exercise~\ref{excfthvec}(c) applies).
\end{description}
\end{exer}

\begin{exer}\label{exclimcxan}
This exercise is a continuation
of Exercise~\ref{excunifk}.
Let $F$ be a locally convex space,
$X$ be a topological space,
$(\gamma_n)_{n\in \N}$
be a sequence of maps
$\gamma_n\colon X\to F$
and $\gamma\colon X\to F$
be a map.
Show the following:
\begin{description}[(D)]
\item[(a)]
If $F$ is Mackey complete,
$X$ is an open subset of~$\C$,
each $\gamma_n$ is complex analytic
and $\gamma_n\to \gamma$ uniformly,
then also $\gamma$ is complex analytic.
\item[(b)]
Given a complex Banach algebra~$\cA$,
set
$\eta_n(x):=\sum_{k=0}^n \frac{1}{k!}x^k$
for $x\in \cA$ and $n\in \N$.
Show that~$\eta_n$ converges uniformly
on each open ball
$B_r(0)\sub \cA$. Deduce that
$\exp\colon \cA\to \cA$, $\exp(x):=\sum_{n=0}^\infty \frac{1}{n!}x^n$
is complex analytic.
\item[(c)]
If $F$ is sequentially complete,
$(a_k)_{k\in \N}$
a sequence in~$F$ and
$0\not=x \in \C$ a complex number such that
$\sum_{k=0}^\infty x^k\,a_k$
converges, then
$\{x^k\, a_k\colon k\in \N_0\}$
is bounded in~$F$
and $\sum_{k=0}^\infty  z^k\, a_k$
converges uniformly on $\bD^0_r$
to a complex analytic map,
for each $r\in \; ]0,|x|[$.
(This remains valid if~$F$ is merely
Mackey complete, since the relevant
partial sums form Mackey--Cauchy sequences).
\item[(d)]
If $F$ is Mackey complete,
$X$ an open subset
of a metrizable locally convex space~$E$,
each $\gamma_n$ is complex analytic
and $\gamma_n|_K$ converges uniformly
to $\gamma|_K$ for each compact subset
$K\sub X$,
then $\gamma$ is complex analytic.
\end{description}
\end{exer}

\begin{exer}\label{exccxpolcts}
Let $p=\sum_{j=0}^k p_j \colon E\to F$
be a polynomial between complex locally
convex spaces, with homogeneous
components~$p_j$.
Show that~$p$ is G\^{a}teaux analytic.
Deduce that if~$p$ is continuous, then~$p$ is complex analytic and $p_j\colon E\to F$
is continuous for each $j\in \{0,\ldots, k\}$.
Compare with Proposition~\ref{contcompo}, which was established by entirely different arguments.
\end{exer}

\begin{exer}\label{excholtoBan}
Let $E, F$ be complex Banach spaces
and $f\colon U\to F$ be a map on an open subset
$U\sub E$. Show that
the following conditions are equivalent:
\begin{description}[(D)]
\item[(a)]
$f$ is complex analytic;
\item[(b)]
$f$ is G\^{a}teaux analytic
and locally bounded, i.e.
each $x\in U$ has an open neighborhood
$V\sub U$ such that $f(V)$ is bounded
in~$F$;
\item[(c)]
$f$ is locally bounded and weakly
analytic, i.e., $\lambda\circ f\colon U\to\C$
is complex analytic for each $\lambda\in F'$.
\end{description}
\end{exer}

\begin{exer}\label{excholonfindm}
Show that a map $f\colon U\to F$
from an open subset $U\sub \C^n$
to a Mackey complete,
complex locally convex space~$F$
is complex analytic if and only if it is weakly
analytic. (Hint: Each $x\in U$ has a compact neighborhood~$K$.
Show  that, if $f$ is weakly analytic,
then $f(K)$ is bounded).
\end{exer}

\begin{exer}\label{exccompBou1}
Let $E, F$ be complex Banach spaces,
$U\sub E$ be open,
$f\colon U\to F$ be a complex
analytic map and $x\in U$.
There is $r>0$ such that $V:=B_r^E(x)\sub U$
and $f$ is bounded on $V$.
\begin{description}[(D)]
\item[(a)]
Using the Cauchy Estimates,
prove that
$\|\delta_x^kf\|_{\op}\leq \frac{k!}{r^k}\|f|_V\|_\infty$,
where $\|\delta_x^kf\|_{\op}$
is the norm of the
homogeneous polynomial $\delta_x^kf$
(as in Exercise~\ref{excnormpolvsmulti})
and $\|\cdot\|_\infty$ the supremum norm.
Hence
$\sum_{k=0}^\infty \frac{s^k}{k!}
\|\delta^k_xf\|_{\op}< \infty$
for each $s<r$.
\item[(b)]
Using that $\lim_{k\to\infty}k/\sqrt[k]{k!}=e$
(Euler's constant)
as a consequence of Stirling's
formula,
deduce that
$\sum_{k=0}^\infty \frac{s^k}{k!}\,
\|d^{\,(k)}f(x,\cdot)\|_{\op}
< \infty$ for all $s< \frac{r}{e}$.
\end{description}
Thus, for maps between open subsets of Banach spaces,
complex analyticity in our sense coincides
with complex analyticity in the sense
of Bourbaki~\cite{Bou67}.
\end{exer}

\begin{exer}\label{exccompBou2}
Let $E:=\ell^1(\N,\C)$ be the Banach space
of absolutely summable complex sequences,
with the norm given by $\|x\|_1:=\sum_{k=1}^\infty
|x_k|$ for $x=(x_k)_{k\in\N}\in E$.
We consider $f\colon E \to \C$,
$f(x):=\sum_{k=1}^\infty 2^k \, (x_k)^{2k}$.
\begin{description}[(D)]
\item[(a)]
Exploiting that $\sum_{k=n}^\infty 2^k|x_k|^{2k}\leq
\sum_{k=n}^\infty 2^{-k}=2^{1-n}$
for each sequence $(x_k)_{k\geq n}$
of complex numbers $x_k\in \bD_{1/2}$,
show that the series defining~$f$
converges uniformly
on some neighborhood of each $x\in E$.
Deduce that $f$ is complex analytic.
\item[(b)]
Show that $f$ is unbounded on $\wb{B}_1^E(x)$, for each $x\in E$.
\item[(c)]
Define $g\colon E\to \C^\N$, $g(x):=(f(kx))_{k\in \N}$.
Show that $g$ is complex analytic
but unbounded on each $0$-neighborhood
in~$E$.
\end{description}
Being not locally bounded at~$0$,
the function $g$ from (c) is not complex analytic
in the sense of Bourbaki~\cite{Bou67}.
In particular, \cite[3.3.1]{Bou67}
is false.
\end{exer}

\begin{exer}\label{excmulticauch}
Let $f\colon U\to F$ be complex analytic
map from an open subset $U\sub \C^n$
to a complex locally convex space~$E$.
Given $r=(r_1,\ldots, r_n)\in \;]0,\infty[^n$
and $x=(x_1,\ldots,x_n)\in \C^n$,
we let
$P_r(x):=B_{r_1}^\C(x_1)\times \cdots\times B_{r_n}^\C(x_n)\sub \C^n$
be the open polycylinder of multi-radius~$r$
in~$\C^n$.
\begin{description}
\item[(a)]
(Cauchy's Integral Formula in several variables).
Show by induction on~$n$: If $U$ contains
the closure $\wb{P}_r(x)$ of $P_r(x)$,
then
\[
\!\!\! \partial^\alpha \! f(z)= \frac{\alpha!}{(2\pi i)^n}
\int_{|\zeta_1-x_1|=r_1}\!\!\!\!\!\!\! \cdots\int_{|\zeta_n-x_n|=r_n}
\!\frac{f(\zeta_1,\ldots, \zeta_n )\, d\zeta_n
\cdots d\zeta_1}{(\zeta_1-z_1)^{\alpha_1+1}
\cdots(\zeta_n-z_n)^{\alpha_n+1}}
\]
for all $z=(z_1,\ldots, z_n)\in P_r(x)$
and all $\alpha\in \N_0^n$.
\item[(b)]
(Cauchy estimates).
Deduce that
\begin{equation}\label{cauchymult}
\|\partial^\alpha \! f(x)\|_p \; \leq \; \frac{\alpha!}{r^\alpha}\,
\sup\big\{\|f(\zeta_1,\ldots,\zeta_n)\|_p\colon |\zeta_j-x_j|=r_j\big\}
\end{equation}
for each continuous seminorm $\|\cdot\|_p$ on~$E$
and each $\alpha\in \N_0^n$.
\end{description}
\end{exer}

\begin{exer}\label{excmultpowser}
This exercise is directed to readers familiar
with absolutely summable
families in Banach spaces
(see Dieudonn\'{e}~\cite[Ch.\,V, \S\,2]{Di60}).
Using the Cauchy estimates (\ref{cauchymult}),
show that a map $f\colon U\to E$
from an open subset $U\sub \C^n$ to
a Banach space~$E$ is complex analytic if and only
if it is given close to each $x\in U$
by an absolutely convergent power series,
i.e.\
$f(x+z) = \sum_{\alpha\in \N_0^n}z^\alpha a_\alpha$
for $z$ in some $0$-neighborhood in~$\C^n$
and suitable coefficients $a_\alpha\in E$.
Dieudonn\'{e} uses
this property to define complex
analytic maps.
\end{exer}

\begin{exer}\label{exc-cx-notCk}
We give $\C^\N$ the product topology
and consider
$f\colon \C \to \C^\N$, $f(z):=(e^{kz})_{k\in \N}$.
\begin{description}[(b)]
\item[(a)]
Show that $f$ is $C^\infty_\C$
and calculate $f'(z) \in \C^\N$
as well as $f^{(2)}(z)$
for $z\in \C$.
\item[(b)]
Equip $E:=\Spann_\C\big(f(\C)\cup f'(\C)\big)\sub \C^\N$
with the induced topology.
Show that the corestriction $g:=f|^E\colon \C \to E$
is $C^1_\C$ but not $C^2_\C$
(whence $g$ is not complex analytic).
Hint: Use without proof that the functions
$h_{n,z}\colon \N\to\C$, $h_{n,z}(k):=k^n z^k$
are linearly
independent for $n\in \N_0$, $z\in \C^\times$.
\item[(c)]
Given $k\in \N$,
find a vector subspace $F\sub \C^\N$
containing~$f(\C)$
such that $f$ is $C^k_\C$ as a map
to~$F$, but not $C_\C^{k+1}$.
\end{description}
\end{exer}

\begin{exer}\label{exchlindepdop}
Show that the functions $h_{k,z}$ $($for $k\in \N_0$, $z\in \C^\times)$
are linearly independent in~$\C^\N$,
as follows:
\begin{description}[(b)]
\item[(a)]
Check that,
for each fixed
$z\in \C^\times$, the functions
$h_{k,z}$ ($k\in \N_0$)
are linearly independent.
\item[(b)]
Now consider the shift operator
$S\colon \C^\N\to\C^\N$, $S(f)(n)=f(n+1)$.
Then $S$ is a linear endomorphism
of the space $\C^\N$ of all complex sequences.
Given $z\in \C$, we let
$V^z$ be the generalized
eigenspace of~$S$,
consisting of all $f\in \C^\N$ such that
$(S-z)^kf=0$ for some $k\in \N$.
Show that $h_{k,z}\in V^z$
for each~$z\in \C^\times$ 
and $k\in \N_0$,
by induction on~$k$.
\end{description}
The assertion follows, as the sum $\sum_{z\in \C^\times}V^z$
of generalized eigenspaces is direct.
\end{exer}
\end{small}
\section{Real analytic mappings}\label{secrealan}
We now define real analytic maps.
These are certain well-behaved functions
which are, in particular,
smooth and
given locally by their Taylor series
around each point
(the actual definition will require
a slightly stronger property).
To state the precise definition,
we need
complexifications
of topological vector spaces.
\begin{defn}\label{defcxvsp}
If $E$ is a real topological vector space,
then $E_\C:=E\times E$, equipped
with the product topology,
becomes a complex topological vector space
if we define
\[
(x+iy).(u,v)\;:=\; (xu-yv, xv+yu)\quad\mbox{for $x,y\in \R$, $u,v\in E$.}
\]
The complex vector space $E_\C$ so obtained is
called the \emph{complexification} \index{complexification!of l.c. space} 
of~$E$.
Note that $E_\C$ is locally convex if so is~$E$.
Identifying $E$ with $E\times \{0\}\sub E_\C$,
we have $E_\C=E\oplus iE$ as a real vector space.
\end{defn}
More generally, a complex topological
vector space $F$, together with a continuous
$\R$-linear map
$\gamma\colon E\to F$, is called
a complexification of~$E$ if there exists
an isomorphism $\phi \colon F\to E_\C=E\times E$
of complex topological vector spaces such
that $\phi(\gamma(x))=(x,0)$ for each
$x\in E$. Again, we write $E_\C:=F$.\\[3mm]
Now real analytic mappings are defined as follows:
\begin{defn}
Let $E$ and $F$ be real locally convex topological
vector spaces,
$U\sub E$ be open, and $f\colon U\to F$ be a map.
Then $f$ is called \emph{real analytic}\index{real analytic}\index{$\R$ analytic}
(\emph{$\R$-analytic},
$C^\omega_\R$, or simply
$C^\omega$)
if it extends to a complex analytic mapping
$V\to F_\C$,
defined on some open neighborhood $V$ of~$U$ in $E_\C$.
\end{defn}
Real analyticity of a map $f\colon U\to F$ is a \emph{local}
property:
\begin{lem}\label{localprop}
If there exists
an open cover $(U_j)_{j\in J}$
of $U$ such that $f|_{U_j}$
is real analytic for each $j\in J$,
then $f$ is real analytic.
\end{lem}
\begin{prf}
If $f|_{U_j}$ is real analytic for each~$j$, then for each
$x\in U$ we find a convex, open $0$-neighborhood
$W_x\sub E$ such that $x+W_x\sub U$
and $f|_{x+W_x}=g_x|_{x+W_x}$
for a complex analytic map $g_x\colon V_x:=(x+W_x)+iW_x\to F_\C$.
Given $x,y\in U$,
the set
$V_x\cap V_y$ is convex as
an intersection of convex
sets,
and hence connected.
If $V_x\cap V_y\not=\emptyset$,
then
the formula
$V_x\cap V_y=((x+W_x)\cap (y+W_y))+i(W_x\cap W_y)$
shows that
$\emptyset\not= E\cap (V_x\cap V_y)=
(x+W_x)\cap (y+W_y)= :C$.
Since $g_x$ and $g_y$ coincide
on~$C$, they coincide on all of $V_x\cap V_y $
by Theorem~\ref{identitycx}(c).
Thus $g_x|_{V_x\cap V_y}=g_y|_{V_x\cap V_y}$ for all $x,y\in U$, whence
$g:=\bigcup_{x\in U}g_x\colon \bigcup_{x\in U} V_x\to F_\C$
is a well-defined complex analytic mapping extending~$f$.
\end{prf}
\begin{prop}\label{chainrealana}
If $E,F$ and~$H$ are real locally convex spaces,
$U\sub E$ and $V\sub F$ are open
and $f\colon U\to V\sub F$
and $g\colon V\to H$ are real analytic maps,
then also $g\circ f\colon U\to H$ is real analytic.
\end{prop}
\begin{prf}
Since $g$ is real analytic,
there is an open subset $V_1\sub F_\C$
and a complex analytic map
$g_1\colon V_1\to H_\C$ such that
$V\sub V_1$ and $g_1|_V=g$.
Since $f$ is real analytic,
there is an open subset $U_1\sub E_\C$
and a complex analytic map
$f_1\colon U_1\to F_\C$ such that
$U\sub U_1$ and $f_1|_U=f$.
Then $U_2:=(f_1)^{-1}(V_1)$ is an open subset
of~$U_1$ which contains~$U$.
After replacing $U_1$ with $U_2$ if necessary,
we may assume that $f_1(U_1)\sub V_1$.
Then $g_1\circ f_1 \colon U_1\to H_\C$
is complex analytic and $(g_1\circ f_1)|_U=g\circ f$,
whence $g\circ f$ is real analytic.
\end{prf}
\begin{lem}\label{realansmooth}
Every real analytic map
$f\colon E\supseteq U\to F$
is smooth.
For each complex analytic extension
$f_1\colon E_\C\supseteq
U_1\to F_\C$ of~$f$ with $U\sub U_1$,
we have $d^{\,(k)}f(x,y_1,\ldots,y_k)=
d^{\,(k)}(f_1)(x,y_1,\ldots,y_k)$
for all $k\in \N$, $x\in U$ and $y_1,\ldots, y_k\in E$.
\end{lem}
\begin{prf}
The map $f_1\colon U_1\to F_\C$
being complex analytic, it is $C^\infty_\C$ and
hence also $C^\infty_\R$.
The inclusion map $\gamma\colon E\to E_\C$,
$\gamma(x)=x$ being continuous linear and thus
$C^\infty_\R$, we deduce that
$f=f_1\circ \gamma|_U^{U_1}$
is $C^\infty_\R$ as a map to $F_\C$.
Since $f(U)\sub F$ where
$F=F\times \{0\}$ is closed in $F\times F=F_\C$,
Lemma~\ref{corestr} shows that $f$ is also
$C^\infty_\R$ as a map to~$F$.
If $x\in U$ and $y\in E$, choose a sequence $(s_n)_{n\in\N}$
of non-zero real numbers such that $x+s_ny\in U$ for all~$n$.
Then
\[
\frac{1}{s_n}(f_1(x+s_ny)-f_1(x))=\frac{1}{s_n}(f(x+s_ny)-f(x))
\]
converges to both $df_1(x,y)$ and $df(x,y)$, whence $df_1(x,y)=df(x,y)$.
This verifies the final assertion for $k=1$, and the general case follows by induction, repeating the preceding argument for the iterated directional derivatives.
\end{prf}
Every real analytic map is given locally by its Taylor
series.
\begin{lem}\label{lemsimphalf}
If $f\colon E\supseteq U\to F$ is real analytic
and $x\in U$,
then there exists an open $0$-neighborhood
$Y\sub E$ such that $x+Y \sub U$
and $f(x+y)=\sum_{k=0}^\infty \frac{\delta_x^kf(y)}{k!}$
for all $y\in Y$.
\end{lem}
\begin{prf}
Let $f_1 \colon U_1\to F_\C$
be a complex analytic extension of
$f$, defined on an open neighborhood
$U_1$ of~$U$ in $E_\C$.
Let $Y_1\sub E_\C$ be an open, balanced
$0$-neighborhood
such that $x+Y_1\sub U_1$.
Then $f_1(x+y)=\sum_{k=0}^\infty \frac{\delta_x^kf_1(y)}{k!}$
for all $y\in Y_1$ by Corollary~\ref{corcxmap}
and hence
$f(x+y)=\sum_{k=0}^\infty \frac{\delta_x^kf_1(y)}{k!}
=\sum_{k=0}^\infty \frac{\delta_x^kf(y)}{k!}$
for all $y\in Y:=Y_1\cap (U - x)\sub E$,
using Lemma~\ref{realansmooth}.
\end{prf}
\begin{rem}
If $E$ is a Fr\'{e}chet space,
then also the converse
of Lemma~\ref{lemsimphalf} holds:
\emph{$f\colon E\supseteq U\to F$
is real analytic if and only
if $f$ is $C^\infty_\R$ and given by its
Taylor series close to each point}.
This follows from
the fact
that, if $E$ is a Fr\'{e}chet
space, then
a map $f\colon E\supseteq U\to F$
to a real locally convex space~$F$
is real analytic if and only if
it is continuous, and
locally of the form
$f(x+y)=\sum_{k=0}^\infty p_k(y)$
for a pointwise convergent series
of continuous homogeneous
polynomials $p_k\colon E\to F$
of degree~$k$
(which necessarily coincide
with $\frac{1}{k!}\delta_x^kf$):
cf.\ \cite[Thm.\,5.2]{BS71b}.
We shall not use this fact,
nor prove~it.
\end{rem}
\begin{thm}[\textbf{Identity Theorem}]\label{identityreal}
Let $E$ and $F$ be real locally convex spaces,
$U\sub E$ be an open,
connected subset and $f\colon U\to F$
be real analytic.
If
$f|_V=0$
for some open, non-empty subset $V\sub U$,
then $f=0$.
More generally, if there is $x\in U$
such that $\delta_x^kf=0$ for all
$k\in \N_0$, then $f=0$.
\end{thm}
\begin{prf}
It is clear that the first condition
implies the second.
We therefore assume now that there exists
$x\in U$ such that $\delta_x^kf=0$
for each $k\in \N_0$.
There exists an open subset
$U_1\sub E_\C$ such that $U\sub U_1$
and a complex analytic map
$f_1\colon U_1\to F_\C$ such that $f_1|_U=f$.
After replacing $U_1$ with its connected
component containing~$U$,
we may assume that $U_1$ is connected.
Then $f_1(x)=f(x)=0$
and $d^{\,(k)}(f_1)(x,y_1,\ldots,y_k)=
d^{\,(k)}f(x,y_1,\ldots,y_k)=0$
for each $k\in \N$ and $y_1,\ldots, y_k\in E$,
by Lemma~\ref{realansmooth},
whence $d^{\,(k)}f(x,y_1,\ldots,y_k)=0$
for all $y_1,\ldots,y_k\in E_\C$
by complex $k$-linearity of $d^{\,(k)}(f_1)(x,\cdot)$.
Hence $f_1=0$ by case~(b)
of Theorem~\ref{identitycx}
and hence also $f=f_1|_U=0$.
\end{prf}
We want to show that every
complex analytic map is also real analytic.
The following two lemmas facilitate this.
\begin{lem}\label{cxofcx}
Given a complex vector space
$E$, put $J\colon E\to E$, $J(x):=ix$.
Let $E_\R$ be the real vector space
underlying $E$ and $(E_\R)_\C=E_\R\oplus iE_\R$
be its complexification.
Then the map
\[
\phi\colon (E_\R)_\C\to E\oplus \wb{E}\,,
\quad \phi(x+iy):=(x+Jy,x-Jy)
\]
is an isomorphism of complex
vector spaces,
where $\wb{E}$ denotes~$E$,
equipped with the ``opposite
complex structure''
with scalar
multiplication $\C\times \wb{E}\to \wb{E}$,
$(z,x)\mto z*x:=\wb{z}x$.
\end{lem}
\begin{prf}
We realize $(E_\R)_\C$ as $E_\R\times E_\R$,
as in Definition~\ref{defcxvsp}.
Given $u,v\in E$,
the equations
$u=x+Jy$ and $v=x-Jy$
have a unique solution $(x,y)\in E^2$,
given by
$x=(u+v)/2$ and $y=J(v-u)/2$.
Hence $\phi$ is a bijection, with
\begin{equation}\label{invphi}
\phi^{-1}\colon E\oplus \wb{E}\to (E_\R)_\C\,,\qquad
\phi^{-1}(u,v)\;=\;
\Big(\frac{u+v}{2},\, \frac{J(v-u)}{2}\Big)\,.
\end{equation}
Clearly $\phi$ is real linear.
It is complex linear, since also
$\phi(i(x,y))=\phi(-y,x)=({-y}+Jx, {-y}-Jx)=(J(x+Jy),
(-J)(x-Jy))=i\phi(x,y)$.
\end{prf}
It is clear from the definition of~$\phi$
and (\ref{invphi}) that $\phi$
is an isomorphism of topological
vector spaces if~$E$ is a complex topological
vector space.
\begin{lem}
If $f\colon E\supseteq U\to F$
is complex analytic, then so
is $f$, considered as a mapping from
$U\sub \wb{E}$ to $\wb{F}$
$($using the opposite complex
structures$)$.
\end{lem}
\begin{prf}
The continuity of~$f$ carries over.
Furthermore, given $x\in U$,
we have
\begin{equation}\label{rettoexp}
f(x+z)\;=\;\sum_{k=0}^\infty
p_k(z)
\end{equation}
for all $z$ in a $0$-neighborhood
in~$E$ and suitable continuous homogeneous
polynomials $p_k\colon E\to F$
of degree~$k$.
Here $p_k(z)=\beta_k(z,\ldots, z)$
for a suitable
continuous,
symmetric $k$-linear map $\beta_k\colon E^k\to F$.
Then also $\beta_k\colon \wb{E}^k\to \wb{F}$
is $k$-linear, as $\beta_k(z*y_1,y_2,\ldots, y_k)=
\beta_k(\wb{z}y_1,y_2,\ldots, y_k)=\wb{z}
\beta_k(y_1,y_2,\ldots, y_k)=z*
\beta_k(y_1,y_2,\ldots, y_k)$.
Therefore each $p_k$ also is a continuous
homogeneous polynomial of degree~$k$
from $\wb{E}$ to $\wb{F}$,
and thus (\ref{rettoexp})
shows that $f$ is complex analytic
from $U\sub \wb{E}$ to $\wb{F}$.
\end{prf}
\begin{prop}
Every complex analytic map
is also real analytic.
\end{prop}
\begin{prf}
Let $E$ and $F$ be complex locally
convex spaces, $U\sub E$ be open and $f\colon
U\to F$ be a complex analytic map.
We write $\wb{U}$ for~$U$,
considered as an open subset of~$\wb{E}$.
Since $f$ is complex analytic
as a map from $\wb{U}$
to $\wb{F}$ by the preceding lemma,
also the map $f\times f\colon U\times \wb{U}\to F\times
\wb{F}$, $(x,y)\mto (f(x), f(y))$
is complex analytic.
For each $x_0\in U$,
there exists an open $x_0$-neighborhood
$V\sub E$ and an open $0$-neighborhood
$W\sub E$ such that $V+J(W)\sub U$
and $V-J(W)\sub U$.
If $\phi\colon (E_\R)_\C\to E\times \wb{E}$
is the isomorphism of complex
locally convex spaces from Lemma~\ref{cxofcx}
and $\psi\colon (F_\R)_\C\to F\times \wb{F}$
is defined analogously,
then
\[
g\, :=\, \psi^{-1}\circ (f\times f)\circ \phi|_{V+iW}^{U\times \wb{U}}
\colon V+iW\to (F_\R)_\C
\]
is a complex analytic map on $V+iW\sub (E_\R)_\C
=E_\R\oplus i E_\R $
which extends $f|_V$
as $g(x)=\psi^{-1}(f(x),f(x))=f(x)$
for each $x\in V$.
Hence $f$ is locally real analytic and
hence real analytic, by Lemma~\ref{localprop}.
\end{prf}
\begin{prop}\label{extend-real-to-complex}
Let $E$ be a real locally convex space, $F$ be a complex locally convex space,
$U\sub E$ be an open subset and $f\colon U\to F$ be a real analytic map.
Then there exists an open subset $V\sub E_\C$ and a complex analytic mapping
$g\colon V\to F$ such that $U\sub V$ and $g|_U=f$.
\end{prop}
\begin{prf}
Let $F_\R$ be the real locally convex space underlying~$F$
and $(F_\R)_\C$ be a complexification of $F_\R$.
Since $f\colon U\to F_\R$ is real analytic,
there exists an open subset $V\sub E_\C$ and a complex analytic mapping $h\colon V\to (F_\R)_\C$
such that $U\sub V$ and $h|_U=f$.
The real linear map $\id\colon F_\R\to F$ induces a complex linear map
$\alpha\colon (F_\R)_\C\to F$, $v+iw\mapsto v+iw$ for $v,w\in F_\R$,
using multiplication with~$i$ in $(F_\R)_\C$ and $F$, respectively.
Then $\alpha\circ f=f$ and $g:=\alpha\circ h\colon V\to F$
is a complex analytic extension of~$f$.
\end{prf}
\subsection*{Example: Inversion in a real continuous
inverse algebra}
To illustrate the concept of a real
analytic map,
we show that inversion
in a (unital) real continuous inverse
algebra is not only smooth
(as shown in Corollary~\ref{invsmoocia}),
but real analytic.
\begin{prop}\label{invrealana}
Let $\cA$ be a real continuous inverse
algebra.
Then the inversion map
$\iota\colon \cA^\times\to \cA^\times$, $\iota(x):=x^{-1}$
is real analytic.
\end{prop}
Two lemmas are useful for the proof.
\begin{lem}\label{spotcia}
Let $\cA$ be a unital, associative
topological algebra. Then the following holds:
\begin{description}[(D)]
\item[\rm (a)]
If $\cA^\times$ is an identity neighborhood,
then $\cA^\times$ is open in~$\cA$.
\item[\rm (b)]
If the inversion map $\iota\colon \cA^\times\hspace*{-.5mm}\to \cA^\times $
is continuous at~$\one$, then $\iota$ is continuous.
\end{description}
\end{lem}
\begin{prf}
(a) Given $x\in \cA$, consider
the left multiplication map
$\lambda_x\colon \cA\to \cA$, $\lambda_x(a):=xa$. The algebra multiplication being
continuous, $\lambda_x\colon \cA\to \cA$
is continuous. Given $x,y\in \cA$,
we have $(\lambda_x\circ \lambda_y)(a)=xya=\lambda_{xy}(a)$
for all $a\in \cA$ and thus $\lambda_x\circ \lambda_y=\lambda_{xy}$.
Furthermore, clearly $\lambda_{\tiny\one} =\id_\cA$.
If $x\in \cA^\times$, then $\lambda_x\circ \lambda_{x^{-1}}=\lambda_{\tiny\one}
=\id_\cA$
and similarly $\lambda_{x^{-1}}\circ \lambda_x=\id_\cA$,
showing that $\lambda_x\colon \cA\to \cA$ is an invertible continuous map,
with continuous inverse $(\lambda_x)^{-1}=\lambda_{x^{-1}}$.
Thus $\lambda_x$ is a homeomorphism.
Consequently, we get for the interior:
\[
(\lambda_x(\cA^\times))^0\;=\; \lambda_x((\cA^\times)^0)\; \ni\;
\lambda_x(\one )\; =\; x\,.
\]
Now, every $x\in \cA^\times$
being an interior point, we deduce that $\cA^\times$ is open.

(b) Given $x\in \cA^\times$, let
$\ell_x:=\lambda_x|_{\cA^\times}^{\cA^\times}$
and define $\rho_x\colon \cA^\times\to \cA^\times$,
$\rho_x(y)=yx$.
The formula
\[
\iota(\ell_x(y))\, =\, (\ell_x(y))^{-1}\, =\,
(xy)^{-1}\, =\, y^{-1}x^{-1}
\, =\, \rho_{x^{-1}}(y^{-1})\, =\, \rho_{x^{-1}}(\iota(y))
\]
shows that $\iota\circ \ell_x$
is continuous at $\one $,
entailing that $\iota=(\iota\circ \ell_x)\circ (\ell_x)^{-1}$
is continuous at $\ell_x(\one )=x$.
\end{prf}
\begin{defn}
If $\cA$ is a real algebra,
we make $\cA_\C$ a complex algebra
via
\begin{equation}\label{cxmult}
(x_1+iy_1)\cdot (x_2+iy_2)\;:=\;
(x_1x_2-y_1y_2)+i (x_1y_2+y_1x_2)
\end{equation}
for all $x_1,x_2,y_1,y_2\in \cA$.
\end{defn}
\begin{rem}\label{remeasycheck}
It is easy to check that the map
$\cA_\C\times \cA_\C\to \cA_\C$
defined via (\ref{cxmult})
is complex bilinear.
Furthermore, if $\cA$ is associative,
then so is $\cA_\C$,
and if $\cA$ is unital,
then so is $\cA_\C$, with the same unit (exercise). 
\end{rem}
\begin{lem}\label{compliscia}
If $\cA$ is a real continuous inverse
algebra, then also $\cA_\C$ is a continuous
inverse algebra.
\end{lem}
\begin{prf}
Since $\cA$ is a continuous inverse algebra, we find
an open identity neighborhood $U\sub \cA^\times$
and an open $0$-neighborhood $V\sub \cA$
such that
\[
\one +(a^{-1}b)^2\in \cA^\times\;\;\,
\mbox{for all $\, a\in U$, $b\in V$.}
\]
For $a,b$ as before, we have
$a+ib=a(\one +ia^{-1}b)$ in $\cA_\C$,
where $a$ is invertible
and so is $\one +ia^{-1}b$.
In fact, abbreviating $c:=a^{-1}b$,
we observe that $\one +ic$ and $\one -ic$ commute,
whence $\one +ic$ also commutes with $\one +c^2=(\one +ic)(\one -ic)$
and its inverse.
Now
$(\one +ic)(\one +c^2)^{-1}(\one -ic)=
(\one +c^2)^{-1}(\one -ic)(\one +ic)=\one$
implies that $(\one +c^2)^{-1}(\one -ic)$
is the inverse of $\one +ic$.
We have shown that the open subset $U+iV$
is contained in $\cA_\C^\times$
and
\[
(a+ib)^{-1}=(\one +(a^{-1}b)^2)^{-1}(\one -ia^{-1}b)a^{-1}
\;\;\,
\mbox{for all $\,(a,b)\in U\times V$,}
\]
which depends continuously on~$(a,b)$.
Using Lemma~\ref{spotcia}, we deduce that
$\cA_\C$ is a continuous
inverse algebra.
\end{prf}
\noindent{\em Proof
of Proposition} \ref{invrealana}.
By Lemma~\ref{compliscia},
$\cA_\C$ is a continuous
inverse algebra.
The inversion map
$j\colon (\cA_\C)^\times\to \cA_\C$, $j(x):=x^{-1}$
is complex analytic (see Example~\ref{invcxcia})
and extends~$\iota$.
Hence $\iota$ is real analytic.\qed
\subsection*{Complex analytic extension of maps on real analytic manifolds}
Let us briefly discuss complexifications and complex analytic extensions in a global setting,
for mappings on real analytic and complex analytic manifolds
modeled on locally convex spaces
(as defined in Sections~\ref{secmanlcx} and~\ref{sec:2.2}).
\begin{defn}
Let $M$ be a real analytic manifold modeled on a real locally convex space~$E$.
A complex analytic manifold $M^*$ modeled on $E_\C$ is called a \emph{complexification}
of $M$ if $M\sub M^*$ holds and each $x\in M$ is contained in the domain $U^*$ of some chart
$\phi\colon U^*\to U'\sub E_\C$ of $M^*$ such that $\phi(M\cap U^*)=E\cap U'$
and $\phi|_{M\cap U^*}\colon M\cap U^*\to E\cap U'$ is a chart for~$M$.
\end{defn}
\begin{prop}\label{global-analyt-extension}
Let $M$ be a real analytic manifold modeled on a real locally convex space~$E$ and $M^*$ be a
complexification of~$M$. If $M$ is closed in $M^*$ and $M^*$ is paracompact,
then the following holds:
For each complex analytic manifold~$L$ modeled on a complex locally convex space and each
real analytic map $f\colon M\to L$, there exists an open subset $U\sub M^*$ such that $M\sub U$
and a complex analytic mapping $g\colon U\to L$ such that $g|_M=f$.
\end{prop}
\begin{prf}
For each $x\in M$,
there exists a chart $\phi_x\colon U_x^*\to U_x'\sub E_\C$
of $M^*$ around~$x$ such that $\phi_x(M\cap U_x^*)=E\cap U'_x$ holds and 
\[ \psi_x:=\phi_x|_{M\cap U_x^*}\colon M\cap U_x^*\to E\cap U'_x \] 
 is a chart for~$M$.
Given $f\colon M\to L$, let $F$ be the modeling space of~$L$.
For $x\in M$, let $\kappa_x\colon P_x\to Q_x\sub F$
be a chart of~$L$ around $f(x)$.
After shrinking $U_x^*$ and $U_x'$,
we may assume that
$\kappa_x\circ f\circ \psi_x^{-1}$ is defined on all of $E\cap U_x'$.
Being real analytic, $\kappa_x\circ f\circ \psi_x^{-1}$
has a complex analytic extension $f_x\colon W_x\to F$
to some open subset $W_x\sub E_\C$ such that $E\cap U_x'\sub W_x$
(see Proposition~\ref{extend-real-to-complex}).
After shrinking $W_x$, we may assume that $f_x(W_x)\sub Q_x$.
After replacing $W_x$ with $W_x\cap U_x'$, we may assume that
$W_x\sub U_x'$ and
\begin{equation}\label{hence-wrx}
E\cap W_x
=E\cap U_x'.
\end{equation}
Now $V_x:=\phi_x^{-1}(W_x)$ is an open subset of $U_x^*$
and
\begin{equation}\label{wrx-2}
M\cap V_x=M\cap U_x^*,
\end{equation}
by(\ref{hence-wrx}).
The function $g_x:=\kappa_x^{-1}\circ f_x\circ \phi_x|_{V_x}\colon V_x\to L$
is complex analytic and $g_x|_{M\cap V_x}=f|_{M\cap V_x}$,
using~(\ref{wrx-2}).
If $x,y\in M$ and $a\in M\cap V_x\cap V_y$,
we have $g_y(a)=f(a)=g_x(a)$ and find a connected, open
$a$-neighborhood $O\sub V_x\cap V_y$
such that $g_y(O)\sub P_x$.
Thus $f_x$ and $\kappa_x\circ g_y\circ \phi_x^{-1}|_{\phi_x(O)}$
are complex analytic functions which coincide
for all $z\in E\cap \phi_x(O)$
with $(\kappa_x\circ f\circ \phi_x^{-1})(z)$.
By the Identity Theorem (Theorem~\ref{identitycx}(c)),
$f_x|_{\phi_x(O)}=\kappa_x\circ g_y\circ \phi_x^{-1}|_{\phi_x(O)}$.
Thus $O(a,x,y):=O$
is an open $a$-neighborhood in $M^*$
such that $g_x|_{O(a,x,y)}=g_y|_{O(a,x,y)}$.
By Lemma~\ref{pre-godement},
there exists an open subset $U\sub M^*$ such that $M\sub U$
and a function
$g\colon U\to L$
such that each
$z\in U$
has an open neighborhood $W\sub U$
such that $W\sub V_x$ for some $x\in M$ and
$g|_W=g_x|_W$.
As a consequence, $g$ is complex analytic and $g|_M=f$.
\end{prf}
\begin{small}
\subsection*{Exercises for Section~\ref{secrealan}}

\begin{exer}
Show that the complexification $\gamma\colon E\to E_\C$ of a real locally convex space~$E$
has the following universal property:
For each continuous $\R$-linear map $\alpha\colon E\to F$ to a complex locally convex space~$F$,
there is a unique continuous $\C$-linear map $\wt{\alpha}\colon E_\C\to F$
such that
\begin{equation}\label{sooneasia}
\wt{\alpha}\circ\gamma=\alpha.
\end{equation}
In particular, for each continuous $\R$-linear map $\beta\colon E_1\to E_2$ between real locally convex spaces
with complexifications $\gamma_1\colon E_1\to (E_1)_\C$ and $\gamma_2\colon E_2\to (E_2)_\C$,
there is a unique continuous $\C$-linear map $\beta_\C\colon (E_1)_\C\to (E_2)_\C$ such that
$\beta_\C\circ \gamma_1=\gamma_2\circ \beta$.
If $E\sub E_\C$ and $\gamma$ is the inclusion map, then (\ref{sooneasia}) simply reads $\wt{\alpha}|_E=\alpha$.
\end{exer}

\begin{exer}\label{excravsTay}
Let $E$ be a sequentially complete
real locally convex space.
Show that a function
$f\colon U\to E$
on an open subset $U\sub \R$
is real analytic if and only if
$f$ is a $C^\infty_\R$-map
which is given locally
around each $x\in U$
by its Taylor series,
viz.\
$f(x+z)=\sum_{k=0}^\infty z^k\, \frac{f^{(k)}(x)}{k!}$
for all $z\in \R$ close to~$0$.
(Hint: Use the Taylor series
to extend $f$ locally
to a complex analytic map).
Remark: As the partial sums
form a Mackey--Cauchy sequence,
the result remains valid
for Mackey complete~$E$.
\end{exer}

\begin{exer}\label{excpolrea}
Show that each continuous linear and each
continuous $n$-linear map between
real locally convex spaces is real analytic.
Deduce that every continuous homogeneous
polynomial between real locally convex spaces
is real analytic.
\end{exer}

\begin{exer}\label{excreatosub}
Let $E$ and $F$ be real locally convex spaces,
$U\sub E$ be open,
$F_0\sub F$ be a sequentially closed
vector subspace and $f\colon U\to F$ be
a map such that $f(U)\sub F_0$.
Show that $f$ is real analytic if and only if
its corestriction $f|^{F_0}\colon U\to F_0$
is real analytic.
\end{exer}

\begin{exer}\label{excreatoprod}
Let $E$ and $F_1,\ldots, F_n$
be real locally convex spaces (where $n\in \N$)
and $U\sub E$ be open. Show that a map
$f=(f_j)_{j=1}^n \colon U\to F_1\times\cdots \times F_n$
is real analytic if and only if
its components $f_j\colon U\to F_j$
are real analytic for all $j\in \{1,\ldots, n\}$.
\end{exer}

\begin{exer}\label{excreatoprodfail}
\begin{description}[(D)]
\item[(a)]
Let $f_1\colon \C\supseteq U\to \C$
be a complex analytic extension
of the real analytic map $f\colon \R\to \R$,
$f(x):=\frac{1}{1+x^2}$.
Show that $\bD\not\sub U$.
\item[(b)]
Find a sequence $(f_k)_{k\in \N}$
of real analytic functions $f_k\colon \R\to\R$
such that the map $f:=(f_k)_{k\in \N}\colon \R\to \R^{\N}$
is not real analytic.
\end{description}
\end{exer}

\begin{exer}
Show that if a map
$f\colon E\supseteq U\to F$ is real
analytic, then also its differentials $d^{\,(j)}f\colon U\times E^j\to F$
are real analytic, for all $j\in \N$.
\end{exer}

\begin{exer}\label{exccheckrm}
Verify the details of Remark~\ref{remeasycheck}.
\end{exer}

\begin{exer}\label{excreaquot}
Let $E$ and $F$ be real locally convex spaces,
$N$ be a closed vector subspace of~$E$
and $q\colon E\to E/N=:E_1$ be the quotient map.
Let $f_1\colon U_1\to F$
be a map, defined on an open subset $U_1\sub E_1$.
Let $U\sub E$ be an open subset such that $q(U)=U_1$.
Show that $f_1$ is real analytic
if and only if $f:=f_1\circ q|_U^{U_1}\colon
U\to F$ is real analytic
(reduce to Lemma~\ref{lemquot} via complex analytic extension).
\end{exer}

\begin{exer}\label{excreaTaynotrea}
Let $E\sub \R^\N$ be the space
of all sequences $(x_n)_{n\in \N}$
such that $\lim_{n\to\infty}\frac{x_n}{n^k}=0$
for some $k\in \N$. Give $E$
the topology induced by~$\R^\N$. Show that
\[
f\colon \R\to E\, , \quad f(t)\, :=\, (\sin(nt))_{n\in \N}
\]
is $C^\infty_\R$
and given by its Taylor series
around each point (cf.\ Exercise~\ref{exccurveinsub}).
Yet, $f$~is not real analytic.
\end{exer}
\end{small}
\section{Differential calculus on normed spaces}\label{provis-Ban}
In this section, we develop the differential calculus of $k$ times continuously Fr\'{e}chet
differentiable mappings between open (and more general) subsets of normed spaces,
including versions of the Inverse Function Theorem and Implicit Function Theorem.
We call such mappings \emph{$FC^k$-maps}, for short.
Each $FC^k$-map is $C^k$, but not conversely (see
Exercise~\ref{exc-C1notFC1} for
a $C^1$-map which fails to be~$FC^1$).
At least, every $C^{k+1}$-map between subsets of normed spaces is $FC^k$,
whence $C^\infty$-maps and $FC^\infty$-maps coincide
(see Proposition~\ref{companotions}).
The topic shall be continued in Sections~\ref{sec-ode} and~\ref{sec-ode-mfd},
which subsume the basic theory of ordinary differential equations (ODEs) in Banach spaces
and Banach manifolds.
Results concerning the parameter dependence of fixed points, inverse functions, and solutions to initial value problems
are important for our ends. Notably, some more specialized results concerning
the parameter dependence of solutions to ODEs will be needed later to see that Lie groups 
modeled on Banach spaces
are well-behaved ($C^0$-regular) Lie groups
(and, likewise,
diffeomorphism groups of $\sigma$-compact finite-dimensional smooth manifolds).
We therefore emphasize parameter dependence throughout our discussions.
\subsection*{Continuously Fr\'{e}chet differentiable maps and {\boldmath$FC^k$}-maps}
Throughout this section,
$(E,\|\cdot\|_E)$, $(F,\|\cdot\|_F)$ and $(E_1,\|\cdot\|_1)$, $(E_2,\|\cdot\|_2)$, etc. are normed spaces over $\K\in\{\R,\C\}$.
\begin{defn}\label{defn-todiff}
A mapping $f\colon U\to F$ on an open subset $U\sub E$ is called
\emph{totally differentiable} \index{totally differentiable map}
at~$x\in U$
if there exist a continuous linear map $f'(x)\colon E\to F$
such that the remainder term $R(y):=f(y)-f(x)-f'(x)(y-x)$
in the (affine) linear approximation
\begin{equation}\label{linapprox}
f(y)=f(x)+f'(x)(y-x)+R(y)\quad\mbox{for $y\in U$}
\end{equation}
satisfies
\begin{equation}\label{little-o}
\lim_{y\to x}\frac{R(y)}{\|y-x\|_E}=0.
\end{equation}
\end{defn}
\begin{rem}
(a) If $f\colon U\to E$ is totally differentiable at~$x$, then $f$ is continuous at~$x$ since
\[
f(y)-f(x)=f'(x)(y-x)+ \|y-x\|_E \frac{R(y)}{\|y-x\|_E} \to 0
\]
as $y\not=x$ tends to~$x$.\medskip

\noindent
(b) The linear map $f'(x)$ is uniquely determined by (\ref{linapprox}) and (\ref{little-o})
as total differentiability at~$x$ implies the existence of
the directional derivative $df(x,y)$ for all $y\in E$, and
\begin{equation}\label{primevsd}
f'(x)(y)=df(x,y).
\end{equation}
It suffices to show this if $y\not=0$.
For $0\not= t\in\K$ with $x+ty\in U$, we have
\[
\frac{1}{t}(f(x+ty)-f(x))=f'(x)(y)+\frac{1}{t}R(x+ty),
\]
where $\frac{\|R(x+ty)\|_F}{|t|}=\|y\|_E \frac{\|R(x+ty)\|_F}{\|ty\|_E}\to 0$ as $t\to 0$, by (\ref{little-o}).
\end{rem}
\begin{numba}
In this section, we endow the space $\cL(E,F)$
of continuous linear maps $\alpha\colon E\to F$ with
the operator norm, unless the contrary is stated.
It defines the topology of 
bounded convergence, whence we also write $\cL(E,F)_b$
for emphasis. If, instead, the compact-open topology is used
(which coincides with the topology of compact convergence), we write $\cL(E,F)_c$.
\end{numba}
\begin{defn}\label{defnFCk}
A mapping $f\colon U\to F$ on an open subset $U\sub E$ is called
\emph{continuously Fr\'{e}chet differentiable}
\index{continuously Fr\'{e}chet differentiable, $FC^1$}
(or $FC^1$) if it is totally differentiable
at each $x\in U$ and the map
\[
f'\colon U\to \cL(E,F)
\]
is continuous.
If $f$ is $FC^1$ and $f'\colon U\to\cL(E,F)$ is $FC^k$ for some $k\in\N$,
we say that $f$ is $FC^{k+1}$.
If $f$ is continuous, we say that $f$ is $FC^0$;
if $f$ is $FC^k$ for each $k\in\N_0$, we say that $f$ is $FC^\infty$.
\end{defn}
A straightforward induction shows that every $FC^{k+1}$-map is $FC^k$.
\begin{exs}\label{extodiff}
(a) Every constant map $f\colon E\to F$ between normed spaces is $FC^1$
with $f'=0$, whence $f$ is $FC^k$ for all $k\in\N$ by induction.
Hence $f$ is~$FC^\infty$.

\nin (b) Every continuous linear map $\lambda\colon E\to F$ is $FC^1$ with $\lambda'(x)=\lambda$
at each $x\in E$ and remainder $R=0$. Since $\lambda'$ is constant and hence $FC^\infty$, we deduce that
$\lambda$ is $FC^\infty$.

\nin (c) Linear combinations of $FC^k$-maps are $FC^k$.

\nin (d) Every continuous bilinear map $\beta\colon E_1\times E_2\to F$
is $FC^\infty$, with
\[
\beta'(x_1,x_2)(y_1,y_2)\!=\!\beta(x_1,y_2)+\beta(y_1,x_2)\;\mbox{for all $(x_1,x_2), (y_1,y_2)\in E_1\!
\times\! E_2$.}
\]
To see this, endow $E_1\times E_2$ with the maximum norm,
\[ \|(x_1,x_2)\|:=\max\{\|x_1\|_1,\|x_2\|_2\} \quad \mbox{ for } \quad 
(x_1,x_2)\in E_1\times E_2.\]
Using the preceding continuous linear map $\beta'(x_1,x_2)$,
the remainder $R(y)$ in the linear approximation at $x=(x_1,x_2)$ is given by
\[
\beta(y_1,y_2)-\beta(x_1,x_2)-\beta(x_1,y_2-x_2)-\beta(y_1-x_1,x_2)=
\beta(y_1-x_1,y_2-x_2)
\]
for $y=(y_1,y_2)\in E_1\times E_2$, with
\[
\frac{\,\|R(y)\|_F }{\|y-x\|}\leq\frac{\|\beta\|_{\op}\|y_1-x_1\|_1\|y_2-x_2\|_2}{\|y-x\|}
\leq\|\beta\|_{\op}\|y-x\|\to 0
\]
as $y\not=x$ tends to~$x$. Hence $\beta$ is $FC^1$. Note that $\beta'(x_1,x_2)$ is linear in $(x_1,x_2)$.
Since 
\begin{align*}
 \|\beta'(x_1,x_2)(y_1,y_2)\|_F
&\leq\|\beta\|_{\op}(\|x_1\|_1\|y_2\|_2+\|y_1\|_1\|x_2\|_2) \\
&\leq
2\|\beta\|_{\op}\|(x_1,x_2)\|\,\|(y_1,y_2)\|,
\end{align*}
we have $\|\beta'(x_1,x_2)\|_{\op}\leq2\|\beta\|_{\op}\|(x_1,x_2)\|<\infty$,
whence $\beta'$ is continuous linear and hence $FC^\infty$. Thus $\beta$ is $FC^\infty$.
\end{exs}
\begin{rem}\label{impbil}
It is essential for calculus in normed spaces that the following bilinear maps are continuous
(and hence $FC^\infty$):
\begin{description}[(D)]
\item[(a)]
The evaluation map $\ve\colon \cL(E,F)\times E\to F$, $(\alpha,x)\mto\alpha(x)$;
\item[(b)]
The composition map 
\[ \Gamma\colon \cL(E_2,E_3)\times\cL(E_1,E_2)\to\cL(E_1,E_3), \quad 
(\alpha,\beta)\mto\alpha\circ\beta. \] 
\end{description}
In fact, $\|\ve(\alpha,x)\|_F \leq\|\alpha\|_{\op}\|x\|_E$
and $\|\Gamma(\alpha,\beta)\|_{\op}\leq\|\alpha\|_{\op}\|\beta\|_{\op}$,
whence $\|\ve\|_{\op}\leq 1$ and $\|\Gamma\|_{\op}\leq 1$.
\end{rem}
\begin{rem}\label{L-bifunct}
Given a continuous linear map $\alpha\colon E_1\to E_2$ between normed spaces
and a normed space~$E$,
we shall also use the linear maps
\[
\cL(E,\alpha)\colon\cL(E,E_1)\to\cL(E,E_2),\quad \beta\mto\alpha\circ \beta
\]
and
\[
\cL(\alpha,E)\colon \cL(E_2,E)\to\cL(E_1,E),\quad \beta\mto\beta\circ\alpha
\]
which are continuous with
\begin{equation}\label{est-opno-bif}
\|\cL(E,\alpha)\|_{\op}\leq\|\alpha\|_{\op}\quad\mbox{and}\quad
\|\cL(\alpha, E)\|_{\op}\leq\|\alpha\|_{\op}
\end{equation}
(by the estimates for composition maps as just discussed).
\end{rem}
\begin{numba}
In this section, we endow the space $\cL^k(E_1,\ldots, E_k;F)$ of continuous $k$-linear maps
$E_1\times\cdots\times E_k\to F$ with the norm $\|\cdot\|_{\op}$
(as in Exercise~\ref{excbilinctso}), and we write $\cL^k(E_1,\ldots,E_k;F)_b$ for emphasis.
Thus $\cL^1(E;F)=\cL(E,F)$.
When the compact open topology induced by $C(E_1\times\cdots\times E_k;F)$ is used
instead, we write $\cL^k(E_1,\ldots, E_k;F)_c$
(cf.\ also Definition~\ref{defn-top-multi} and Exercise~\ref{exer-semin-norm}).
\end{numba}
\begin{rem}\label{multileval}
The evaluation map
\[
\ve\colon \cL^k(E_1,\ldots ,E_k;F)\times E_1\times\cdots\times E_k\to F,\;
(\beta,x_1,\ldots, x_k)\mto\beta(x_1,\ldots, x_k)
\]
is $(k+1)$-linear, and continuous with $\|\ve\|_{\op}\leq 1$ as
\[
\|\ve(\beta,x_1,\ldots, x_k)\|_F=\|\beta(x_1,\ldots, x_k)\|_F\leq\|\beta\|_{\op}\|x_1\|_1\cdots\|x_k\|_k
\]
for all $\beta\in\cL^k(E_1,\ldots, E_k;F)$ and $x_j\in E_j$ for $j\in\{1,\ldots, k\}$.
\end{rem}
\begin{rem}
Exponential laws like
\[
\cL^2(E_1,E_2;F)_b\cong \cL(E_1,\cL(E_2,F)_b)_b
\]
(as discussed in Exercise~\ref{exc-multil-exp})
and corresponding identifications
play a certain role in the differential calculus in normed spaces.
As the recursive Definition~\ref{defnFCk}
avoids recourse to higher Fr\'{e}chet derivatives,
we can minimize the use of such identifications;
exponential laws for spaces of multilinear maps will only be used
for the proofs of Proposition~\ref{companotions}(b) (via Lemma~\ref{prepacompa})
and~(c).
\end{rem}
With a view towards right-hand-sides of differential equations, which are often defined
on non-open subsets, we also consider non-open domains.
\begin{defn}
A mapping $f\colon U\to F$ on a locally convex subset $U\sub E$ with dense interior
is called
\emph{continuously Fr\'{e}chet differentiable} (or $FC^1$) if it is continuous,
its restriction $f|_{U^0}$ to the interior is $FC^1$ and $(f|_{U^0})'$ 
admits a (necessarily unique) continuous extension
\[
f'\colon U\to \cL(E,F).
\]
If $f$ is $FC^1$ and $f'\colon U\to\cL(E,F)$ is $FC^k$ for some $k\in\N$,
we say that $f$ is $FC^{k+1}$.
If $f$ is continuous, we say that $f$ is $FC^0$;
if $f$ is $FC^k$ for each $k\in\N_0$, we say that $f$ is $FC^\infty$.
\end{defn}
Besides $FC^k$-maps, we can consider $C^k$-maps between
normed spaces (as in Chapter~\ref{chapcalcul}).
We now compare the two concepts, for~$k=1$.
General $k\in\N$ will be considered in Propositions \ref{FCk-findim}
and~\ref{companotions}.
\begin{lem}\label{FC1vsC1}
Let $U\sub E$ be a locally convex subset with dense interior.
For a mapping $f\colon U\to F$, the following conditions are equivalent:
\begin{description}[(D)]
\item[\rm(a)]
$f$ is continuously Fr\'{e}chet differentiable;
\item[\rm(b)]
$f$ is $C^1$ and the map $U\to\cL(E,F)$, $x\mto df(x,\cdot)$
is continuous.
\end{description}
If the conditions are satisfied, then
\begin{equation}\label{primevsd2}
df(x,y)=f'(x)(y)\quad\mbox{for all $x\in U$ and $y\in E$.}
\end{equation}
\end{lem}
\begin{prf}
(a)$\impl$(b):
If $f$ is $FC^1$, then $f$ is continuous
with directional derivatives $df(x,y)=f'(x)(y)$ for all $x\in U^0$ and $y\in E$
(see (\ref{primevsd})), whence
\[
df(x,y):=f'(x)(y)=\ve(f'(x),y)\;\;\mbox{for $(x,y)\in U\times E$}
\]
(with $\ve$ as in Remark~\ref{impbil}(a)) defines a continuous extension $df$ to $U\times E$. Thus~$f$ is $C^1$
and (\ref{primevsd2}) holds.

(b)$\impl$(a): Writing $f'(x):=df(x,\cdot)$ for $x\in U$, we obtain a continuous map $f'\colon
U\to\cL(U,F)$.
Given $x\in U$ and $\ve>0$,
there exists a convex $x$-neighborhood $Y\sub U$
such that
\begin{equation}\label{subit2}
\|f'(y)-f'(x)\|_{op}\leq\ve\quad\mbox{for all $y\in Y$.}
\end{equation}
Using the Mean Value Theorem, we deduce that
\begin{eqnarray*}
\lefteqn{\|f(z)-f(x)-f'(x)(z-x)\|_F}\qquad\\
&=&\left\|\int_0^1 (f'(x+t(z-x))-f'(x))(z-x)\,dt\right\|_F\\
&\leq&\max\{\|f'(x+t(z-x))-f'(x)\|_{op}\colon t\in[0,1]\}\|z-x\|_E
\leq \ve\|z-x\|_E
\end{eqnarray*}
for all $z\in Y$, using (\ref{subit2}) with $y:=x+t(z-x)$.
If $x\in U^0$, this means that $f$ is totally differentiable at~$x$
with derivative $f'(x)$. Hence $f$ is $FC^1$.
\end{prf}
See Exercise~\ref{exc-C1notFC1}
for an example of a $C^1$-map between Banach spaces which fails to be $FC^1$.
\begin{prop}[Chain Rule for {\boldmath$FC^k$}-maps]\label{chainFCk}
Let $U\sub E_1, V\sub E_2$ be locally convex subsets with dense interior,
$k\in\N\cup\{\infty\}$ and $f\colon U\to E_2$ as well as
$g\colon V\to E_3$ be $FC^k$-maps such that $f(U)\sub V$.
Then also the mapping
$g\circ f\colon U\to E_3$, $x\mto g(f(x))$ is $FC^k$, and
\begin{equation}\label{CRtotal}
(g\circ f)'(x)=g'(f(x))\circ f'(x)\quad\mbox{for all $x\in U$.}
\end{equation}
\end{prop}
\begin{prf}
Being $FC^1$, both $f$ and $g$ are $C^1$ and $df(x,\cdot)=f'(x)$
as well as $dg(z,\cdot)=g'(z)$ depend continuously on $x\in U$ and $z\in V$,
respectively (see Lemma~\ref{FC1vsC1}).
Hence $g\circ f$ is $C^1$ with
\[
d(g\circ f)(x,y)=dg(f(x),df(x,y))=(g'(f(x))\circ f'(x))(y)
\]
for all $x\in U$ and $y\in E$ (see Proposition~\ref{chainno}). Using the continuous bilinear
composition map
\[
\Gamma\colon \cL(E_2,E_3)\times \cL(E_1,E_2)\to\cL(E_1,E_3),\;
(\alpha,\beta)\mto\alpha\circ\beta,
\]
we have
\[
d(g\circ f)(x,\cdot)=g'(f(x))\circ f'(x)=\Gamma(g'(f(x)),f'(x)),
\]
which is continuous in $x\in U$. 
By Lemma~\ref{FC1vsC1}, $g\circ f$ is $FC^1$ and (\ref{CRtotal})
holds. If $k\geq 2$ and the assertions hold for $k-1$ in place of~$k$,
consider the continuous linear (and hence $FC^{k-1}$-) maps
\[
\lambda_1\colon \cL(E_2,E_2)\to\cL(E_2,E_3)\times \cL(E_1,E_2),\;\; \alpha\mto (\alpha, 0)
\]
and $\lambda_2\colon \cL(E_1,E_2)\to\cL(E_2,E_3)\times\cL(E_1,E_2)$, $\beta\mto (0,\beta)$.
Then
\[
(g'\circ f,f')=\lambda_1\circ g'\circ f+\lambda_2\circ f'
\]
is $FC^{k-1}$ as $\lambda_1\circ g'\circ f$ and $\lambda_2\circ f'$ are
$FC^{k-1}$ by the inductive hypothesis.
Since $\Gamma$ is $FC^{k-1}$,
using the inductive hypothesis again, we see that
\[
(g\circ f)'=\Gamma\circ (g'\circ f, f')
\]
is $FC^{k-1}$. As $g\circ f$ is $FC^1$ and $(g\circ f)'$ is $FC^{k-1}$,
the map $g\circ f$ is $FC^k$.
\end{prf}
\begin{rem}
Let $U\sub E$ be a locally convex subset with dense interior
and
$k\in\N_0\cup\{\infty\}$.
Then a map $f=(f_1,f_2)\colon U\to E_1\times E_2$ is $FC^k$ if and only if both
$f_1$ and $f_2$ are $FC^k$.\\[2mm]
[To see this, let $\pi_j\colon E_1\times E_2\to E_j$ be the projection onto the $j$th factor
for $j\in\{1,2\}$
and $\lambda_j\colon E_j \to E_1\times E_2$ be the map taking $x$ to $(x,0)$ (if $j=1$)
and $(0,x)$ (if $j=2$). Then $\pi_1$, $\pi_2$, $\lambda_1$, and $\lambda_2$
are continuous linear and thus $FC^\infty$.
Using Proposition~\ref{chainFCk}, the assertion follows from the identities
\begin{equation}\label{maybereus}
f_1=\pi_1\circ f,\quad f_2=\pi_2\circ f\quad\mbox{and}\quad
f=\lambda_1\circ f_1+\lambda_2\circ f_2.
\end{equation}
\end{rem}
\begin{lem}\label{strdiff}
Let $U\sub E$ be a locally convex subset with dense interior
and $f\colon U\to F$ be continuously Fr\'{e}chet
differentiable. Let $x\in U$ and $R\colon U\to F$ 
be the remainder term in the affine
linear approximation
\[
f(y)=f(x)+f'(x)(y-x)+R(y)
\]
of $f$ around~$x$. Then $R|_{U\cap B^E_r(x)}$ is Lipschitz for small $r>0$, and
\begin{equation}\label{liptozero}
\lim_{r\to 0}\Lip(R|_{U\cap B^E_r(0)})=0.
\end{equation}
\end{lem}
\begin{prf}
Given $\ve>0$, there exists a convex $x$-neighborhood $Y\sub U$
such that
\[
\|f'(y)-f'(x)\|_{\op}\leq\ve\quad\mbox{for all $y\in Y$.}
\]
For all $y,z\in Y$, using the Mean Value Theorem,
we deduce that
\begin{eqnarray*}
\|R(z)-R(y)\|_F&=&\|f(z)-f(y)-f'(x)(z-y)\|_F\\
&=&\left\|\int_0^1(f'(y+t(z-y))-f'(x))(z-y)\,dt\right\|_F\\
&\leq & \int_0^1\|f'(y+t(z-y))-f'(x)\|_{\op}\|z-y\|_E\, dt
\leq\ve\|z-y\|_E,
\end{eqnarray*}
from which the assertion follows.
\end{prf}
\begin{rem}\label{str-as-lip-pert}
In the situation of Lemma~\ref{strdiff},
define $g\colon U\to F$ via $g(y):=f(y)-f'(x)(y)$ for $y\in U$.
Assume that~$U$ is open in~$E$.
As $g(y)=R(y)+f(x)-f'(x)(x)$ differs from $R(y)$ only by a constant,
the map
\[
f|_{B^E_r(x)}=f'(x)|_{B^E_r(x)}+g|_{B^E_r(x)}
\]
is a Lipschitz perturbation of the
continuous linear map $f'(x)$ for small $r>0$, with
$\Lip(g|_{B^E_r(x)})=\Lip(R|_{B^E_r(x)})\to 0$ as $r\to 0$.
This will be useful when we prove the Inverse Function Theorem
(Theorem~\ref{inv-fct-class}).
\end{rem}
\begin{prop}\label{companotions}
Let $U\sub E$ be a locally convex subset with dense
interior, $f\colon U\to F$ be a function and $k\in\N$.
Then the following holds:
\begin{description}[(D)]
\item[\rm(a)]
If $f$ is $FC^k$, then $f$ is~$C^k$.
\item[\rm(b)]
If $f$ is $C^{k+1}$, then $f$ is $FC^k$.
\item[\rm(c)]
$f$ is $FC^k$ if and only if $f$ is $C^k$ and the map
$U\to\cL^k(E,\ldots, E;F)_b$, $x\mto d^{\,(k)}f(x,\cdot)$
is continuous.
\end{description}
\end{prop}
\begin{prf}
(a) If $f$ is $FC^k$, then $f$ is $C^1$ and
\begin{equation}\label{enabchai}
df=\ve \circ (f'\times \id_E),
\end{equation}
where $\ve\colon \cL(E,F)\times E\to F$,
$(\alpha,x)\mto\alpha(x)$ is the evaluation map (see Lemma~\ref{FC1vsC1}).
Being continuous bilinear, $\ve$ is $C^{k-1}$.
The $FC^{k-1}$-map $f'$ is $C^{k-1}$ by induction (resp., continuous if $k=1$)
and
the continuous linear map $\id_E$ is $C^{k-1}$,
whence $f'\times \id_E$ is $C^{k-1}$ (cf.\ Exercise~\ref{exccartpromap}).
Using the Chain Rule (Proposition~\ref{chainno}), we deduce from (\ref{enabchai}) that
$df$ is $C^{k-1}$. Hence $f$ is~$C^k$.

(b) if $k=1$: We show that $f$ is $FC^1$ if it is~$C^2$.
Let $x\in U$. Since $d^{\,(2)}f\colon U\times E\times E\to F$
is continuous and $d^{\,(2)}f(x,0,0)=0$,
given $\ve>0$,
there is a convex neighborhood $Y \sub U$ of~$x$
and $r>0$ such that
\[
\|d^{\,(2)}f(y,v,w)\|_F\leq\ve\quad\mbox{for all $y\in Y$ and $v,w\in B^E_r(0)$.}
\]
After shrinking~$Y$, we may assume that $Y-Y\sub B^E_{r^2}(0)$.
Applying the Mean Value Theorem to the $C^1$-function $df$,
we obtain for all $y\in Y$ and $v\in E$ such that $\|v\|_E\leq 1$:
\begin{eqnarray*}
\|(f'(y)-f'(x))(v)\|_F &=& \|df(y,v)-df(x,v)\|_F\\
&=&\left\|\int_0^1 d^{\,(2)}f(x+t(y-x),v,y-x)\, dt\right\|_F\\
&\leq &\int_0^1 \Big\|d^{\,(2)}f\big(x+t(y-x),rv,\frac{1}{r}(y-x)\big)\Big\|_F\,dt
\leq\ve.
\end{eqnarray*}
Thus
\begin{equation}\label{subit}
\|f'(y)-f'(x)\|_{\op}\leq\ve \quad\mbox{for all $y\in Y$,}
\end{equation}
showing that $f'$ is continuous at~$x$ and hence continuous.
Therefore $f$ is $FC^1$, by Lemma~\ref{FC1vsC1}.

(c) for $k=1$ is the content of Lemma~\ref{FC1vsC1}.

For $k\geq 2$, the assertions (b) and (c) are less important for this book,
and we recommend to skip their proofs for a first reading. The proof of (b) uses another lemma.
\begin{lem}\label{prepacompa}
Let $k\in\N$, $\ell\in \N_0$ and $U\sub E$ be a locally convex subset
with dense interior.
If a map
$f\colon U\to\cL^k(E_1,\ldots, E_k;F)$ is $C^{\ell+1}$ as a map to $\cL^k(E_1,\ldots, E_k;F)_c$,
then it is $FC^\ell$ as a map to $\cL^k(E_1,\ldots, E_k;F)_b$. 
\end{lem}
\begin{prf}
The proof is by induction on $\ell\in\N_0$.
If $\ell=0$, we
consider $f$ as a $C^1$-map
to $C(E_1\times\cdots \times E_k,F)$
and deduce from Theorem~\ref{explawCkell} and Remark~\ref{inwashexp}(b)
that the map
\[
f^\wedge\colon U\times (E_1\times\cdots \times E_k)\to F,\;
(x,y_1,\ldots, y_k)\mto f(x)(y_1,\ldots, y_k)
\]
is $C^{1,0}$.
Hence $f=(f^\wedge)^\vee$ is continuous (and thus $FC^0$) as a
map to $\cL^k(E_1,\ldots, E_k;F)_b$, by Lemma~\ref{sim-old}.

If $f$ is $C^{\ell+1}$ as a map to $\cL^k(E_1,\ldots, E_k;F)_c$
for some $\ell\in\N$, then
\[
df\colon U\times E\to \cL^k(E_1,\ldots,E_k;F)_c\sub C(E_1\times\cdots\times E_k;F)
\]
is $C^\ell$, whence
\[
\wh{df}\colon (U\times E)\times (E_1\times\cdots\times E_k)\to F,\quad
(x,y,x_1,\ldots, x_k)\mto df(x,y)(x_1,\ldots, x_k)
\]
is a $C^{\ell,0}$-map.
Interpreting $\wh{df}$ as a $C^{\ell,0}$-map
$U\times (E\times E_1\times\cdots\times E_k)\to F$ (which is a weaker property),
Theorem~\ref{explawCkell}
provides a $C^\ell$-function
\[
g:=(\wh{df})^\vee\colon U\to C(E\times E_1,\ldots\times E_k,F).
\]
As $g(U)$ is a subset of the closed vector subspace $\cL^{k+1}(E,E_1,\ldots, E_k;F)$
of $C(E\times E_1\times\cdots\times E_k,F)$, Lemma~\ref{nonocorestr}
allows us to interpret $g$ as a $C^\ell$-function
to
$\cL^{k+1}(E,E_1,\ldots,E_k;F)_c$. By the inductive hypothesis, $g$ is $FC^{\ell-1}$
as a map to $\cL^{k+1}(E,E_1,\ldots, E_k;F)_b$.
Consider the map
\[
\Phi\colon \cL^{k+1}(E,E_1,\ldots,E_k;F)_b\to\cL(E,\cL^k(E_1,\ldots, E_k;F)_b)_b
\]
from Exercise~\ref{exc-multil-exp}, which is an isomorphism of topological vector spaces.
If we can show that
\[
f'=\Phi\circ g,
\]
then $f'$ will be $FC^{\ell-1}$ and thus $f$ will be~$FC^\ell$ (as required).
For $x\in U$, $y\in E$ and $x_j\in E_j$ for $j\in\{1,\ldots, k\}$,
we calculate
\begin{eqnarray*}
\Phi(g(x))(y)(x_1,\ldots, x_k)&=&g(x)(y,x_1,\ldots,x_k)=\wh{df}(x,y,x_1,\ldots,x_k)\\
&=&df(x,y)(x_1,\ldots,x_k),
\end{eqnarray*}
whence $\Phi(g(x))(y)=df(x,y)=f'(x)(y)$ and hence $\Phi(g(x))=f'(x)$.
\end{prf}
\noindent
\emph{Proof of Proposition}~\ref{companotions}, \emph{completed}.
(b) If $f\colon U\to F$ is $C^{k+1}$ with $k\geq 1$, then $f$ is $FC^1$ (by the case $k=1$
already settled).
Since $df\colon U\times E\to F$ is $C^k$, the map
\[
(df)^\vee\colon U\to C(E,F),\quad x\mto df(x,\cdot)
\]
is $C^k$, by Theorem~\ref{explawCkell}.
As the image of $(df)^\vee$ is contained in the closed vector subspace
$\cL(E,F)$ of $C(E,F)$, we deduce that also the corestriction $f'\colon U\to\cL(E,F)_c$
of $(df)^\vee$ is~$C^k$. Then $f'\colon U\to\cL(E,F)_b$ is $FC^{k-1}$,
by Lemma~\ref{prepacompa}, and thus $f$ is $FC^k$.

(c) Let $k\geq 2$ and suppose the assertion holds for $k-1$ in place of~$k$.
Let $f\colon U\to F$ be a $C^k$-map such that $U\to\cL^k(E\times\cdots\times E;F)_b$,
$x\mto d^{\,(k)}f(x,\cdot)$ is continuous.
By (b), $f$ is $FC^{k-1}$, whence $f'\colon U\to\cL(E,F)_b$ is $FC^{k-2}$
and hence $C^{k-2}$, by~(a).
As $df\colon U\times E\to F$ is $C^{k-1}$, the map $(df)^\vee\colon U\to C(E,F)$,
$x\mto df(x,\cdot)$
and hence also
\[
g:=f'\colon U\to\cL(E,F)_c
\]
is $C^{k-1}$.
We claim that
\begin{equation}
  \label{eq:122}
d^{\,(k-1)}g(x,y_1,\ldots, y_{k-1})(y)=d^{\,(k)}f(x,y_1,\ldots, y_{k-1},y)
\end{equation}
for all
$x\in U$ and $y,y_1,\ldots, y_{k-1}\in E$.
If this is true, then
\[
d^{\,(k-1)}g(x,\cdot)=\Phi(d^{\,(k)}f(x,\cdot))\quad\mbox{for all $x\in U$,}
\]
using the isomorphism
\begin{equation}\label{2usespsi}
\Phi\colon \cL^k(E,\ldots,E;F)_b\to\cL^{k-1}(E,\ldots, E;\cL(E,F)_b)_b
\end{equation}
of topological vector spaces discussed in Exercise~\ref{exc-multil-exp}.
Thus
\begin{equation}\label{uugly}
U\to\cL^{k-1}(E,\ldots, E;\cL(E,F)_b)_b,\;\,
x\mto d^{\,(k-1)}g(x,\cdot)=\Phi(d^{\,(k)}f(x,\cdot)) 
\end{equation}
is continuous. As the evaluation map
\[
\ve\colon \cL^{k-1}(E,\ldots, E;\cL(E,F)_b)_b\times E^{k-1}\to\cL(E,F)_b
\]
is continuous (see \ref{multileval}), we deduce that
\[
d^{\,(k-1)}g(x,y_1,\ldots, y_{k-1})=\ve(d^{\,(k-1)}g(x,\cdot),y_1,\ldots, y_{k-1})\in\cL(E,F)_b
\]
depends continuously on $(x,y_1,\ldots, y_{k-1})\in U\times E^{k-1}$.
Hence $f'=g$ is $C^{k-1}$ as a map to $\cL(E,F)_b$, by Lemma~\ref{twotoplem}(b).
Since the map in (\ref{uugly})
is continuous, the map
$f'$ is $FC^{k-1}$, by induction. Hence $f$ is~$FC^k$.

To verify \eqref{eq:122}, 
let $y\in E$. As the evaluation map $\ve_y\colon\cL(E,F)_c\to F$, $\alpha\mto\alpha(x)$
is continuous linear, we have
\begin{eqnarray}
d^{\,(k-1)}(\ve_y\circ g)(x,y_1,\ldots,y_{k-1})& =& \ve_y(d^{\,(k-1)}g(x,y_1,\ldots,y_{k-1}))\label{prelimviacha}\\
&=& d^{\,(k-1)}g(x,y_1,\ldots, y_{k-1})(y)\notag
\end{eqnarray}
for all $(x,y_1,\ldots,y_{k-1})\in U\times E^{k-1}$.
Since $\ve_y\circ g=df(\cdot,y)$, the first term in (\ref{prelimviacha}) is given by
$d^{\,(k-1)}(\ve_y\circ g)(x,y_1,\ldots,y_{k-1})=d^{\,(k)}f(x,y,y_1,\ldots, y_{k-1})
=d^{\,(k)}f(x,y_1,\ldots, y_k,y)$,
whence the claim ist true.

If, conversely, $f$ is $FC^k$, then $f$ is $C^k$, by~(a).
Moreover, $f'\colon U\to\cL(E,F)_b$ is $FC^{k-1}$,
whence $f'$ is $C^{k-1}$ by induction and
\[
U\to\cL^{k-1}(E,\ldots, E;\cL(E,F)_b)_b,\quad x\mto d^{\,(k-1)}(f')(x,\cdot)
\]
is a continuous map. Given $y\in E$, the evaluation map $\ve_y\colon\cL(E,F)_b\to F$, $\alpha\mto\alpha(y)$
is continuous linear and $\ve_y\circ f'=df(\cdot,y)$, entailing that
\begin{eqnarray*}
d^{\,(k-1)}(f')(x,y_1,\ldots,y_{k-1})(y) &=& d^{\,(k-1)}(\ve_y\circ f')(x,y_1,\ldots, y_{k-1})\\
&=&d^{\,(k)}f(x,y,y_1,\ldots,y_{k-1})
\end{eqnarray*}
for all $(x,y_1,\ldots,y_{k-1})\in U\times E^{k-1}$.
Thus, re-using the map $\Phi$ from (\ref{2usespsi}),
we have $d^{\,(k)}f(x,\cdot)=\Phi^{-1}(d^{\,(k-1)}(f')(x,\cdot))\in\cL^k(E,\ldots, E;F)_b$,
which depends continuously on $x\in U$.
\end{prf}
\begin{prop}\label{FCk-findim}
Let $E$ be a finite-dimensional vector space,
$U\sub E$ be a locally convex subset with dense interior,
$F$ be a normed space, $f\colon U\to F$ be a mapping and
$k\in\N\cup\{\infty\}$.
Then~$f$ is~$C^k$ if and only if~$f$ is~$FC^k$.
\end{prop}
\begin{prf}
If $f$ is $FC^k$, then $f$ is~$C^k$ by
Proposition~\ref{companotions}(a).
To prove the converse, we may assume that $k\in\N$.
The proof is by induction.
The case $k=1$:
Let $v_1,\ldots, v_n$ be a basis for~$E$.
Then
\[
\phi \colon \cL(E,F)_b\to F^n,\quad \alpha\mto (\alpha(v_1),\ldots,\alpha(v_n))
\]
is an isomorphism of topological vector spaces (see Exercise~\ref{exc-on-fin}).
The map
\begin{equation}\label{hence-indc1}
\phi\circ f'=(df(\cdot,v_1),\ldots,df(\cdot,v_n))\colon U\to F^n
\end{equation}
is $C^{k-1}$ and hence continuous, entailing that $f'=\phi^{-1}\circ(\phi\circ f')$
is continuous. Thus $f$ is $FC^1$, by Lemma~\ref{FC1vsC1}.
If $f$ is $C^k$ with $k\geq 2$, then we deduce from~(\ref{hence-indc1}) that
$\phi\circ f'$ is $C^{k-1}$.
Then $f'=\phi^{-1}\circ (\phi\circ f')$ is $C^{k-1}$
and thus $FC^{k-1}$, by the inductive hypothesis.
Hence~$f$ is~$FC^k$.
\end{prf}
\subsection*{Fixed points, inverse and implicit functions}
Among other things, we shall prove
the classical inverse function theorem
for $FC^k$-maps now:
\begin{thm}[Inverse Function Theorem for {\boldmath$FC^k$}-maps]\label{inv-fct-class}
Let $E$ be a Banach space
and $f\colon U\to E$ be an $FC^k$-map on an open subset $U\sub E$, where $k\in\N\cup\{\infty\}$.
Let $x_0\in U$. If $f'(x_0)\colon E\to E$ is invertible,
then there exists an open $x_0$-neighborhood $V\sub U$ such that $f(V)$ is open in~$E$ and $f|_V\colon V\to f(V)$ is an $FC^k$-diffeomorphism.
\end{thm}
We shall also prove the classical implicit function theorem
(and slightly more, as we need not assume that $E$ is complete and
$U\sub E$ an open subset):
\begin{thm}[Implicit Function Theorem for {\boldmath$FC^k$}-maps]\label{impl-fct-class}
Let $E$ be a normed space and $F$ be a Banach space,
$U\sub E$ be a locally convex subset with dense interior and $V\sub F$ be an open subset.
Let $f\colon U\times V\to F$ be an $FC^k$-map with $k\in\N\cup\{\infty\}$
and $(x_0,y_0)\in U\times V$ such that $f(x_0,y_0)=0$ and $(f_{x_0})'(y_0)\colon F\to F$ is invertible,
where $f_x:=f(x,\cdot)\colon V\to F$ for $x\in U$.
Then there exists an open $x_0$-neighborhood $U_0\sub U$, an open $y_0$-neighborhood
$V_0\sub V$ and an $FC^k$-map $\phi\colon U_0\to V_0$ such that
\begin{equation}\label{sol-is-graph}
\{(x,y)\in U_0\times V_0\colon f(x,y)=0\}=\graph(\phi).
\end{equation}
\end{thm}
As before, vector spaces (and differentiability properties) are considered over a base
field $\K\in\{\R,\C\}$ here and in the remainder of the section.
The proofs for Theorems~\ref{inv-fct-class} and \ref{impl-fct-class}
require some preparations; they shall be given after
Remarks~\ref{L-not-lip} and \ref{rem:2.3.44},
respectively. 

A version of the implicit function theorem is even available if~$E$
is an arbitrary locally convex space (Corollary~\ref{impl-Ck}) and we shall 
derive Theorem~\ref{impl-fct-class} from this result.
As a tool for the proofs, we study the dependence
of fixed points of contractions on parameters (which is
also of independent
interest). Let us recall the concept.
\begin{defn}
Let $(X,d)$ be a metric space. A mapping $f\colon X\to X$ is said to be
a \emph{contraction} \index{contraction}
if there exists $L\in[0,1[$ such that
\[
d(f(x),f(y))\leq L \, d(x,y)\quad\mbox{for all $x,y\in X$,}
\]
viz.\ $f$ is Lipschitz continuous with $\Lip(f)<1$.
Lipschitz constants $L\in[0,1[$ are also called \emph{contraction constants}. \index{contraction constant}
\end{defn}
\begin{lem}\label{fp-uniqu}
If a contraction $f\colon X\to X$ of a metric space $(X,d)$
has a fixed point~$x_\infty$,
then the latter is necessarily unique.
\end{lem}
\begin{prf}
Let $\wb{x}\in X$ with $f(\wb{x})=\wb{x}$.
If $\wb{x}\not= x_\infty$, then $d(\wb{x},x_\infty)>0$
and thus
\[
d(\wb{x},x_\infty)=d(f(\wb{x}),f(x_\infty))\leq L\, d(\wb{x},x_0)<d(\wb{x},x_\infty),
\]
which is absurd.
\end{prf}
We recall Banach's Fixed Point Theorem (also known as the Contraction Mapping
Principle).
\begin{lem}[Banach's Fixed Point Theorem]\label{ban-fix}
Let $(X,d)$ be a complete metric space with $X\not=\emptyset$
and $f\colon X\to X$ be a contraction with contraction constant $L\in[0,1[$.
Then the following holds:
\begin{description}[(D)]
\item[\rm(a)]
$f$ has a unique fixed point $x_\infty$.
\item[\rm(b)]
For every $x_0\in X$,
\begin{equation}\label{conv-to}
\lim_{n\to\infty}f^n(x_0)=x_\infty.
\end{equation}
\item[\rm(c)]
The following a priori estimate holds: For all $n\in\N_0$,
\begin{equation}\label{a-prio-est}
d(f^n(x_0),x_\infty)\leq \frac{L^n}{1-L}\, d(f(x_0),x_0).
\end{equation}
\end{description}
\end{lem}
\begin{prf}
Given $x_0\in X$,
set $x_n:=f^n(x_0)$ for $n\in\N_0$. Then
\begin{equation}\label{fist-banf}
d(x_{n+1},x_n)\leq L^nd(x_1,x_0)
\end{equation}
for all $n\in\N_0$, as the estimate is trivial for $n=0$
and for $n\in\N$
\[
d(x_{n+1},x_n)=d(f(x_n),f(x_{n-1}))\leq L\,
d(x_n,x_{n-1})\leq L^nd(x_1,x_0)
\]
as $d(x_n,x_{n-1})\leq L^{n-1}d(x_1,x_0)$
by induction.
For $n\in\N_0$ and $m\in \N$, this implies that
\begin{equation}\label{sec-banf}
d(x_{n+m},x_n)\leq \frac{L^n}{1-L}\,d(x_1,x_0);
\end{equation}
in fact, $d(x_{n+m},x_n)\leq \sum_{k=n}^{n+m-1}d(x_{k+1},x_k)
\leq\sum_{k=n}^{n+m-1}L^kd(x_1,x_0)=L^nd(x_1,x_0)\sum_{j=0}^{m-1}L^j\leq \frac{L^n}{1-L}d(x_1,x_0)$
using the triangle inequality, (\ref{fist-banf}), and the summation formula for the
geometric series. By (\ref{sec-banf}), $(x_n)_{n\in\N_0}$ is a Cauchy sequence
and hence convergent. Let $x_\infty\in X$ be its limit.
Then $x_\infty$ is fixed by~$f$, as
\[
f(x_\infty)=f\left(\lim_{n\to\infty}x_n\right)=\lim_{n\to\infty}f(x_n)=\lim_{n\to\infty}x_{n+1}=x_\infty.
\]
Letting $m\to \infty$ in (\ref{sec-banf}), we obtain (\ref{a-prio-est}).
See Lemma~\ref{fp-uniqu} for uniqueness of fixed points.
\end{prf}
The following hypothesis enables meaningful results concerning
fixed points for families of contractions.
\begin{defn}
Let $(X,d)$ be a metric space
and $(f_j)_{j\in J}$ be a family of contractions
$f_j\colon X\to X$. A number $\theta\in [0,1[$
is called a \emph{uniform contraction constant}
\index{uniform contraction constant}
for $(f_j)_{j\in J}$
if it is a contraction constant for each~$f_j$.
We say that $(f_j)_{j\in J}$ is a \emph{uniform family of contractions}
\index{uniform family of contractions}
if it admits a uniform contraction constant, which holds if and only if
\begin{equation}\label{unif-contra-sup}
\sup_{j\in J} \,\Lip(f_j)<1.
\end{equation}
\end{defn}
\begin{prop}\label{fp-cts-dep}
Let $(X,d)$ be a complete metric space with $X\not=\emptyset$.
Let $P$ be a topological space and $f\colon P\times X\to X$
be a mapping such that the maps
\[
f_p:=f(p,\cdot)\colon X\to X
\]
form a uniform family $(f_p)_{p\in P}$ of contractions.
Let~$\theta$ be a uniform contraction constant.
For $p\in P$, let $x_p\in X$ be the unique fixed point of~$f_p$.
Then the following holds:
\begin{description}[(D)]
\item[\rm(a)]
If $f$ is continuous, then $\phi\colon P\to X$, $p\mto x_p$
is continuous.
\item[\rm(b)]
If there exists a metric $d_P$ on~$P$ defining its topology and such that
$f^x:=f(\cdot,x)\colon P\to X$, $p\mto f(p,x)$ is Lipschitz for each
$x\in X$ with
\[
L:=\sup_{x\in X}\Lip(f^x)<\infty,
\]
then $\phi$ is Lipschitz with $\Lip(\phi)\leq \frac{L}{1-\theta}$.
\end{description}
\end{prop}
\begin{prf}
For all $p,q\in P$, (\ref{a-prio-est}) shows that
\begin{equation}\label{basisfound}
d(x_p,x_q)\leq\frac{1}{1-\theta}d(f_q(x_p),x_p)=\frac{1}{1-\theta}
d(f(q,x_p),f(p,x_p)).
\end{equation}
\indent
(a) Let $p\in P$. Given $\ve>0$, there exists a
$p$-neighborhood
$Q\sub P$ such that $d(f(q,x_p),f(p,x_p))< (1-\theta)\ve$ for all $q\in Q$
and thus $d(x_p,x_q)<\ve$, by (\ref{basisfound}). Hence $\phi$ is continuous at~$p$.

(b) Since $d(f(q,x_p),f(p,x_p))\leq Ld_P(q,p)$, (\ref{basisfound}) shows
that $\phi$ is Lipschitz with $\Lip(\phi)\leq L/(1-\theta)$.
\end{prf}
\begin{prop}\label{fp-Ck-dep}
Let $E$ be a locally convex space, $P\sub E$ be a locally convex subset with dense interior,
$(F,\|\cdot\|)$ be a Banach space and $U\sub F$ be a closed convex subset with dense interior.
Let $f\colon P\times U\to U$ be a map
such that the maps $f_p:=f(p,\cdot)\colon U\to U$ form a uniform family
$(f_p)_{p\in P}$ of contractions. Let $k\in\N\cup\{\infty\}$ and
let $x_p$ be the fixed point of~$f_p$ for $p\in P$.
Then the following holds:
\begin{description}[(D)]
\item[\rm(a)]
If $f$ is $C^k$, then also the map $\phi\colon P\to U$, $p\mto x_p$ is $C^k$.
\item[\rm(b)]
If $E$ is a normed space and~$f$ is $FC^k$, then also $\phi\colon P\to U$ is~$FC^k$.
\end{description}
\end{prop}
\begin{prf} Let $\theta\in[0,1[$ be a uniform contraction constant for $(f_p)_{p\in P}$.
We may assume that $k\in\N$.

(a) The case $k=1$:
We claim that, for each $n\in\N_0$, there is a continuous map
$h_n\colon P^{[1]}\to F$ such that
\begin{equation}\label{defn-hk}
h_n(p,q,t)=\frac{(f_{p+tq})^{n+1}(x_p)-(f_{p+tq})^n(x_p)}{t}
\;\mbox{for all $(p,q,t)\in P^{[1]}$ with $t\not=0$,}
\end{equation}
where $P^{[1]}$ is as in (\ref{defnU1no}).
For all $(p,q,t)\in P^{[1]}$ with $t\not=0$, we have
\[
\phi(p+tq)=x_{p+tq}=\lim_{N\to\infty} (f_{p+tq})^N(x_p),
\]
whence the left hand side of
\begin{eqnarray}
\frac{(f_{p+tq})^{N+1}(x_p)-x_p}{t}&=&
\sum_{n=0}^N\frac{(f_{p+tq})^{n+1}(x_p)-(f_{p+tq})^n(x_p)}{t}\notag\\
&=&\sum_{n=0}^N h_n(p,q,t)\label{expl-strat}
\end{eqnarray}
converges to
\[
\frac{\phi(p+tq)-\phi(p)}{t}.
\]
If we can show that $\sum_{n=0}^\infty h_n$ converges to a continuous function
$h\colon P^{[1]}\to F$, then $\phi$ will be $C^1$ with $\phi^{[1]}=h$ (see
Lemma~\ref{linkBGNno}).
We prove the claim by induction.
If $n=0$, using $x_p=f_p(x_p)=f(p,x_p)$ we find that 
\[
\frac{f_{p+tq}(x_p)-x_p}{t}=\frac{f(p+tq,x_p)-f(p,x_p)}{t}=f^{[1]}(p,x_p,q,0,t)
\]
for all $(p,q,t)\in P^{[1]}$ such that $t\not=0$.
Hence $h_0\colon P^{[1]}\to F$,
\[
(p,q,t)\mto f^{[1]}(p,x_p,q,0,t)
\]
is as required.
If $n\in\N$ and $h_{n-1}$ has already been found,
writing
\begin{eqnarray*}
(f_{p+tq})^n(x_p)&=&(f_{p+tq})^{n-1}(x_p)+t\frac{(f_{p+tq})^n(x_p)-(f_{p+tq})^{n-1}(x_p)}{t}\\
&=&(f_{p+tq})^{n-1}(x_p)+th_{n-1}(p,q,t)
\end{eqnarray*}
we obtain that
\begin{eqnarray}
\lefteqn{\frac{(f_{p+tq})^{n+1}(x_p)-(f_{p+tq})^n(x_p)}{t}}\qquad\notag\\
&=&
\frac{f(p+tq,(f_{p+tq})^n(x_p))-f(p+tq,(f_{p+tq})^{n-1}(x_p))}{t}\notag\\
&=&f^{[1]}(p+tq,(f_{p+tq})^{n-1}(x_p),0,h_{n-1}(p,q,t),t)\label{sem-defn}
\end{eqnarray}
for all $(p,q,t)\in P^{[1]}$ such that $t\not=0$.
Note that the final term in (\ref{sem-defn})
makes sense for all $(p,q,t)\in P^{[1]}$
and defines a continuous function $h_n\colon P^{[1]}\to F$.
This completes the recursive construction.

Since $h_0$ is continuous, every $(p_0,q_0,t_0)\in P^{[1]}$
has an open neighborhood $W\sub P^{[1]}$ such that $h_0(W)$ is bounded;
thus
\[
M:= \|h_0|_W\|_\infty <\infty.
\]
Using (\ref{fist-banf}), we find that
\begin{eqnarray*}
\|h_n(p,q,t)\|&=&\frac{1}{|t|} \|(f_{p+tq})^{n+1}(x_p)-(f_{p+tq})^n(x_p)\|
\leq \frac{\theta^n}{|t|}\|f_{p+tq}(x_p)-x_p\|\\
&=& \theta^n\|h_0(p,q,t)\|\leq \theta^nM
\end{eqnarray*}
for all $(p,q,t)\in W$ such that $t\not=0$.
Hence
\[
\|h_n(p,q,t)\|\leq\theta^n M\;\; \mbox{for all $\, (p,q,t)\in W$,}
\]
since $h_n$ is continuous and $\{(p,q,t)\in W\colon t\not=0\}$ is dense in~$W$.
Thus
\[
\sum_{n=0}^\infty \|h_n|_W\|_\infty\leq\sum_{n=0}^\infty \theta^n M=\frac{M}{1-\theta}<\infty,
\]
showing that $\sum_{n=0}^\infty h_n|_W$ converges
uniformly to a continuous function\linebreak
$W\to F$.
As $(p_0,q_0,t_0)$
was arbitrary, we deduce that $\sum_{n=0}^\infty h_n$ converges
locally uniformly
to a continuous function $h\colon P^{[1]}\to F$.

Induction step. If $f$ is $C^k$ for some $k\in\N$ with $k\geq 2$,
then~$\phi$ is $C^{k-1}$ (and hence $C^1$), by the inductive hypothesis.
Since
\[
f(p,\phi(p))=\phi(p)
\]
for $p\in P$, using the Chain Rule we get that
\begin{equation}\label{indu-fixed}
df(p,\phi(p),q,d\phi(p,q))=d\phi(p,q)\;\; \mbox{for all $\, (p,q)\in P\times E$.}
\end{equation}
Now $g\colon P\times E\times F\to F$, $(p,q,y)\mto df(p,\phi(p),q,y)$
is a $C^{k-1}$-map. Setting $g_{p,q}:=g(p,q,\cdot)\colon F\to F$,
we can rewrite (\ref{indu-fixed})
as
\begin{equation}\label{fixed-lev2}
g_{p,q}(d\phi(p,q))=d\phi(p,q)\;\,\mbox{for all $(p,q)\in P\times E$.}
\end{equation}
Since $g_{p,q}=d_1f(p,\phi(p),q)+d_2f(p,\phi(p),\cdot)$,
using Lemma~\ref{lipviaprime}(b) we find that
\begin{eqnarray*}
\Lip(g_{p,q})&=&\Lip(d_2f(p,\phi(p),\cdot))
=\|d_2f(p,\phi(p),\cdot)\|_{\op}\\
&=&\|(f_p)'(\phi(p))\|_{\op}\leq\Lip(f_p)\leq\theta.
\end{eqnarray*}
Thus $(g_{p,q})_{(p,q)\in P\times E}$ is a uniform family of contractions.
As $d\phi(p,q)$ is the fixed point of $g_{p,q}$ by (\ref{fixed-lev2}),
$d\phi$ is $C^{k-1}$ by the inductive hypothesis. So~$\phi$ is~$C^k$.

(b) By \eqref{indu-fixed} in (a), the map $\phi$ is~$C^1$ and
\[
\phi'(p)=f'(p,\phi(p))\circ \lambda_1+ f'(p,\phi(p)) \circ \lambda_2\circ \phi'(p)
\quad \mbox{ for } \quad p\in P,
\] 
using the continuous linear mappings $\lambda_1\colon E\to E\times F$, $q\mto (q,0)$
and\linebreak
$\lambda_2\colon F\to E\times F$, $y\mto (0,y)$.
Thus
\begin{equation}\label{fixed-lev3}
\phi'(p)=h_p(\phi'(p))
\end{equation}
in terms of the map
\[
h\colon P\times \cL(E,F)\to\cL(E,F),\;\, (p,\alpha)\mto
f'(p,\phi(p))\circ\lambda_1 +f'(p,\phi(p))\circ \lambda_2\circ\alpha
\]
which is $FC^{k-1}$ by the Chain Rule and Remark~\ref{impbil}(b),
writing as before $h_p:=h(p,\cdot)$. Here we use 
that $\phi$ is $FC^{k-1}$;
in fact, $\phi$ is continuous (hence $FC^0$)
and if $k\geq 2$, then the $FC^{k-1}$-property
can be assumed by induction.
Using (\ref{est-opno-bif}), we see that
\[
\Lip(h_p)=\|\cL(E,f_p'(\phi(p))\circ\lambda_2)\|_{\op}
\leq\|f_p'(\phi(p))\circ\lambda_2\|_{\op}
\leq\|f_p'(\phi(p))\|_{\op}\leq\theta,
\]
whence $(h_p)_{p\in P}$ is a uniform family of contractions.
Since $\phi'(p)$ is the fixed point of $h_p$ by (\ref{fixed-lev3}),
the inductive hypothesis shows that $\phi'$ is $FC^{k-1}$
and hence continuous, whence~$\phi$ is $FC^1$ by Lemma~\ref{FC1vsC1}.
Thus~$\phi$ is $FC^k$.
\end{prf}
Results concerning real analytic dependence
of fixed points on parameters are also available
(see Exercices~\ref{exerc-fp-an1} and \ref{exerc-fp-an2}).
The following lemma concerning $C^{0,k}$-dependence
can be skipped on a first reading; it is only needed for further technical
results which eventually enable the proof of $C^0$-regularity
for Banach-Lie groups and diffeomorphism groups.
\begin{lem}\label{fp-C0k-dep}
Let $X$ be a topological space,
$E$ be a locally convex space, $P\sub E$ be a locally convex subset with dense interior,
$(F,\|\cdot\|)$ be a Banach space,
$U\sub F$ a closed convex subset with dense interior, and $k\in \N_0\cup\{\infty\}$.\linebreak
Let $f\colon X\times (P\times U)\to U$ be a $C^{0,k}$-map
such that the mappings $f_{x,p}:=f(x,p,\cdot)\colon U\to U$ form a uniform family
$(f_{x,p})_{(x,p)\in X\times P}$ of contractions.
Let $\phi(x,p)\in U$ be the fixed point of~$f_{x,p}$ for $(x,p)\in X\times P$.
Then $\phi\colon X\times P\to F$, $(x,p)\mto\phi(x,p)$ is a $C^{0,k}$-map.
\end{lem}
\begin{prf}
We may assume that $k<\infty$. The proof is by induction on $k\in\N_0$.
For $k=0$, the assertion holds by Proposition~\ref{fp-cts-dep}(a).
Now let $k\in\N$.
For fixed $x\in X$, the map $f_x:=f(x,\cdot)\colon P\times U\to U$ is
$C^k$ (see Lemma~\ref{Ckellpartial}) and defines a uniform family of contractions, whence
$\phi_x:=\phi(x,\cdot)\colon P\to F$ is a $C^k$-map, by Proposition~\ref{fp-Ck-dep}(a).
Holding $x$ fixed, applying the Chain Rule to
\[
\phi_x(p)=f_x(p,\phi_x(p)), 
\]
we obtain that $d\phi_x(p,q)=d(f_x)(p,\phi_x(p),q,d\phi_x(p,q))$ for all $(p,q)\in P\times E$
and thus
\[
d^{(0,1)}\phi(x,p,q)=d^{(0,1)}f(x,p,\phi(x,p),q,d^{(0,1)}\phi(x,p,q)).
\]
By the preceding, $d^{(0,1)}\phi(x,p,q)$
is the fixed point of the map $g_{x,p,q}:=g(x,p,q,\cdot)\colon F\to F$ for $(x,p,q)\in X\times P\times E$, 
if we define
\[
g\colon X\times (P\times E\times F)\to F,\quad
g(x,p,q,y):=d^{(0,1)}f(x,p,\phi(x,p),q,y).
\]
Now $d^{(0,1)}f\colon X\times (P\times U \times E\times F)\to F$ is $C^{0,k-1}$ by
Lemma~\ref{diffphigherd}
and the map
\[
X\times (P\times E\times F)\to P\times U\times E\times F,\quad
(x,p,q,y)\mto (p,\phi(x,p),q,y)
\]
is $C^{0,k-1}$, as~$\phi$ is~$C^{0,k-1}$ by induction (cf.\ also Exercise~\ref{exc-babychainC0k}).
Hence~$g$ is $C^{0,k-1}$, by Lemma~\ref{C0chainy}.
Now
\[
g(x,p,q,y)=
d_2f(x,p,\phi(x,p);q)+
d_3f(x,p,\phi(x,p);y)
\]
in terms of partial differentials.
If $\theta$ is a uniform contraction constant for $(f_{x,p})_{(x,p)\in X\times P}$,
then $\|d_3f(x,p,\phi(x,p);\cdot)\|_{\op}\leq\theta$ for all $(x,p)\in X\times P$, by
Lemma~\ref{lipviaprime}.
As a consequence, $\Lip(g_{x,p,q})\leq\theta$ (again by Lemma~\ref{lipviaprime}).
Hence $(g_{x,p,q})_{(x,p,q)\in X\times P\times E}$
is a uniform family of contractions.
Thus $d^{(0,1)}\phi$ is $C^{0,k-1}$, by the inductive hypothesis.
Hence $\phi$ is $C^{0,k}$, by Lemma~\ref{C0kviaC01}.
\end{prf}
\begin{defn}
Let $k\in\N\cup\{\infty\}$.
Let $E$ and~$F$ be locally convex spaces over $\K\in\{\R,\C\}$
and $f\colon U\to V$ be a mapping between open subsets
$U\sub E$ and $V\sub F$;
when speaking about $FC^k$-maps, we assume that both~$E$ and $F$
are normable.\medskip

\noindent
(a) $f\colon U\to V$ is called a \emph{$C^k$-diffeomorphism},
\emph{$FC^k$-diffeomorphism}, and \emph{$\K$-analytic diffeomorphism},
respectively,
if $f$ is a bijection and both $f$ and $f^{-1}$ are $C^k$-maps (resp., $FC^k$-maps,
resp., $\K$-analytic maps).
In the same fashion, we define $C^k$-diffeomorphisms and $FC^k$-diffeomorphisms
if, instead, $U\sub E$ and $V\sub F$ are locally convex subsets with dense interior.\medskip

\index{$C^k$-diffeomorphism} \index{diffeomorphism!$C^k$} 
\index{local $C^k$-diffeomorphism} 
\index{$FC^k$-diffeomorphism} \index{diffeomorphism!$FC^k$} 
\index{analytic diffeomorphism} \index{diffeomorphism!analytic} 

\noindent
(b) $f$ is called a \emph{local $C^k$-diffeomorphism at $x\in U$}
if there exists an open $x$-neighborhood $W\sub U$ such that $f(W)$
is open in~$F$ and $f|_W\colon W\to f(W)$ is a $C^k$-diffeomorphism.
If $f$ is a local $C^k$-diffeomorphism at each $x\in U$,
then $f$ is called a \emph{local $C^k$-diffeomorphism}.
Local $FC^k$-diffeomorphisms and local $\K$-analytic diffeomorphisms
(as well as such at~$x$) are defined analogously,
replacing $C^k$-maps with $FC^k$-maps and $\K$-analytic mappings, respectively.
\end{defn}
\begin{rem}
If $f\colon U\to V$ (as before) is a $C^1$-diffeomorphism,
applying the Chain Rule to $\id_U=f^{-1}\circ f$ and $\id_V=f\circ f^{-1}$
we find that
$\id_E=(f^{-1})'(f(x))\circ f'(x)$ for all $x\in U$ and
\begin{equation}\label{then-repl}
\id_F=f'(f^{-1}(y))\circ (f^{-1})'(y)
\end{equation}
for all $y\in V$. Taking $y:=f(x)$ in (\ref{then-repl}), we obtain $\id_F=f'(x)\circ (f^{-1})'(f(x))$.
Hence $f'(x)$ is invertible for all $x\in U$, and
\begin{equation}
(f^{-1})'(f(x))=f'(x)^{-1}.
\end{equation}
\end{rem}
\begin{rem}
Let $(E,\|\cdot\|)$ be a normed space with $E\not=\{0\}$ and \break $\alpha\colon E\to E$
be an automorphism of~$E$ as a topological vector space.
In the following, the number
\[
\frac{1}{\|\alpha^{-1}\|_{\op}}
\]
will be encountered repeatedly.
It can be interpreted as the least expansion factor of~$\alpha$,
in the sense that
\begin{equation}\label{least-facto}
\frac{1}{\|\alpha^{-1}\|_{\op}}=\inf\left\{ \frac{\|\alpha(x)\|}{\|x\|}\colon 0\not=x\in E\right\}
\in\; ]0,\infty[.
\end{equation}
To see this, rewrite the right hand side of (\ref{least-facto}) as
\[
\frac{1}{\sup\{\|x\|/\|\alpha(x)\|\colon 0\not=x\in E\}}
=
\frac{1}{\sup\{\|\alpha^{-1}(y)\|/\|y\|\colon 0\not=y\in E\}}.
\]
\end{rem}
We now prove a quantitive version of the inverse function theorem
for Lipschitz perturbations $f=\alpha+g$ of a linear automorphism~$\alpha$.
Its proof is based on a simplified Newton iteration
which does not involve derivatives of~$f$
but uses the given~$\alpha$ instead.
\begin{thm}[Quantitative Inverse Function Theorem]\label{lip-inv-fct}
\,Let\linebreak
$(E,\|\cdot\|)$ be a Banach space with $E\not=\{0\}$
and $\alpha\in\GL(E)$. Let $x\in E$, $r>0$ and $g\colon B^E_r(x)\to E$
be a Lipschitz map such that
\[
\Lip(g)<\frac{1}{\|\alpha^{-1}\|_{\op}}.
\]
Then $f\colon B^E_r(x)\to E$, $y\mto \alpha(y)+g(y)$ has the following
properties:
\begin{description}[(D)]
\item[\rm(a)]
$f$ has open image and is a homeomorphism onto its image.
\item[\rm(b)]
$f^{-1}\colon f(B^E_r(x))\to B^E_r(x)$ is Lipschitz with
\[
\Lip(f^{-1})\leq a^{-1} \quad \mbox{ for } \quad a := 
\frac{1}{\|\alpha^{-1}\|_{\op}}-\Lip(g) > 0.
\]
\item[\rm(c)]
If we set $h:=f^{-1}-\alpha^{-1}$, then $f^{-1}=\alpha^{-1}+h$
and $h\colon f(B^E_r(x))\to E$ is Lipschitz with
\begin{equation}\label{the-lip-h}
\Lip(h)\leq\frac{\|\alpha^{-1}\|_{\op}\Lip(g)}{a}.
\end{equation}
\item[\rm(d)]
Abbreviate
$b:=\|\alpha\|_{\op}+\Lip(g)$.
Then
\[
a\|z-y\|\leq \|f(z)-f(y)\|\leq b\|z-y\|\quad\mbox{for all $y,z\in B^E_r(x)$.}
\]
\item[\rm(e)]
The following estimates for the images of balls are available:
\[
B^E_{ar}(f(x))\sub f(B^E_r(x))\sub B^E_{br}(f(x))
\]
and, more generally,
\begin{equation}\label{ball-ba}
B^E_{as}(f(y))\sub f(B^E_s(y))\sub B^E_{bs}(f(y))
\end{equation}
for all $y\in B^E_r(x)$ and $0 < s \leq r-\|y-x\|$.
\item[\rm(f)]
If $f$ is $C^k$ $($resp., $FC^k)$ with $k\in\N\cup\{\infty\}$,
then also $f^{-1}\colon f(B^E_r(x))\to E$ is $C^k$ $($resp., $FC^k)$.
\end{description}
\end{thm}
\begin{prf}
Let $L:=\Lip(g)$; then
\begin{equation}\label{neenee}
a+L=\frac{1}{\|\alpha^{-1}\|_{\op}}.
\end{equation}
\indent
(d) For all $y,z\in B^E_r(x)$, we have
\[\|f(z)-f(y)\|=\|\alpha(z-y)+g(z)-g(y)\|, \] 
whence $\|f(z)-f(y)\|\leq (\|\alpha\|_{\op}+L)\|z-y\|$
by the triangle inequality
and $\|f(z)-f(y)\|\geq\|\alpha(z-y)\|-\|g(z)-g(y)\|\geq a\|z-y\|$,
using~(\ref{least-facto}).

(b) The lower estimate in (d) entails that~$f$ is injective. If $v,w\in f(B^E_r(x))$,
then $v=f(y)$ and $w=f(z)$ with $y:=f^{-1}(v)$ and $z:=f^{-1}(w)$.
By (d), we have
$\|f^{-1}(w)-f^{-1}(v)\|=\|z-y\|\leq \frac{1}{a}\|f(z)-f(y)\|=\frac{1}{a}\|w-v\|$.
Thus $\Lip(f^{-1})\|\leq\frac{1}{a}$.

(e) Let $y$ and $s$ be as in (\ref{ball-ba}).
The second inclusion in~(\ref{ball-ba}) holds by~(d).
To prove the first, we need only show that
\[
\wb{B}^E_{at}(f(y))\sub f(\wb{B}^E_t(y))
\]
for all $0 < t < r-\|y-x\|$. To this end, let $c\in \wb{B}^E_{at}(f(y))$
and consider
\[
h_c\colon \wb{B}^E_t(y)\to E,\quad z\mto z-\alpha^{-1}(f(z)-c).
\]
Then $h_c(z)\in \wb{B}_t^E(y)$ for all $z\in\wb{B}_t^E(y)$, as
\begin{eqnarray*}
\|h_c(z)-y\|&=& \|\alpha^{-1}(\alpha(z)-\alpha(y)-f(z)+c)\|\\
&=&\|\alpha^{-1}(-g(z)+g(y)-f(y)+c)\|\\
&\leq & \|\alpha^{-1}\|_{\op}(L\|z-y\|+\|c-f(y)\|)
\leq\|\alpha^{-1}\|_{\op}(Lt+at)\leq t,
\end{eqnarray*}
using~(\ref{neenee}).
Moreover, $h_c\colon \wb{B}_t^E(y)\to\wb{B}_t^E(y)$
is a contraction, as
\begin{eqnarray*}
\|h_c(z)-h_c(v)\|&=&\|\alpha^{-1}(\alpha(z)-\alpha(v)-f(z)+f(v))\|\\
&=&\|\alpha^{-1}(g(v)-g(z))\|\leq\|\alpha^{-1}\|_{\op}L\|z-v\|
\end{eqnarray*}
shows that
\begin{equation}\label{thus-uni-co}
\Lip(h_c) \leq \|\alpha^{-1}\|_{\op}L<1.
\end{equation}
By Banach's Fixed Point Theorem, there exists $z\in
\wb{B}^E_t(y)$ such that
\[
z=h_c(z)=z-\alpha^{-1}(f(z)-c)
\]
and thus $c=f(z)$.

(a) By~(b), $f$ is a homeomorphism onto its image.
The first inclusion in (\ref{ball-ba}) shows that $f(B^E_r(x))$
is open.

(c) We have $\id=f^{-1}\circ f=(\alpha^{-1}+h)\circ (\alpha+g)=
\id+\alpha^{-1}\circ g+h\circ f$ and hence
$h\circ f=-\, \alpha^{-1}\circ g$, i.e.,
\[
h=-\,\alpha^{-1}\circ g\circ f^{-1}.
\]
Hence $\Lip(h)\leq\|\alpha^{-1}\|_{\op}\Lip(g)\Lip(f^{-1})$.
Estimating $\Lip(f^{-1})$ as in~(b),
we obtain~(\ref{the-lip-h}).

(f) If $f$ is $C^k$ (resp., $FC^k$),
we show that $f^{-1}$ is $C^k$ (resp., $FC^k$)
on an open neighborhood of $f(y)$ for each $y\in B^E_r(x)$.
To this end, let $0 < t < r-\|y-x\|$. Then
\[
\tilde h \colon \wb{B}^E_{at}(f(y))\times \wb{B}_t^E(y)\to \wb{B}^E_t(y),\quad z\mto z-\alpha^{-1}(f(z)-c)
\]
is $C^k$ (resp., $FC^k$).
By (\ref{thus-uni-co}), the maps $\tilde h_c:=\tilde h(c,\cdot)$
for a uniform family of contractions for $c\in \wb{B}^E_{at}(f(y))$,
and we also saw that the fixed point~$z$ of~$\tilde h_c$ coincides with
$f^{-1}(c)$. Hence $f^{-1}$ is $C^k$ (resp., $FC^k$) on $\wb{B}^E_{at}(f(y))$,
by Proposition~\ref{fp-Ck-dep}.
\end{prf}
\begin{rem}\label{L-not-lip}
Repeating the preceding proof with $L\in [\Lip(g),\frac{1}{\|\alpha^{-1}\|_{\op}}[$
in place of $\Lip(g)$, we obtain statements analogous to (b)--(e)
in Theorem~\ref{lip-inv-fct}, except that $L$ replaces~$\Lip(g)$.
\end{rem}
\noindent\emph{Proof of Theorem}~\ref{inv-fct-class}.
We have $B^E_r(x_0)\sub U$ for some $r>0$.
By Lemma~\ref{strdiff}, after shrinking~$r$ we can achieve that
the function
\[
g\colon B^E_r(x_0)\to E,\quad x\mto f(x)-f'(x_0)(x)
\]
is Lipschitz with
\[
\Lip(g)<\frac{1}{\|f'(x_0)^{-1}\|_{\op}}
\]
(see Remark~\ref{str-as-lip-pert}).
Then $f|_{B^E_r(x_0)}=f'(x_0)+g$ has open image and is an $FC^k$-diffeomorphism
onto the latter, by Theorem~\ref{lip-inv-fct}.\qed\medskip

As a tool, we now study families of local inverses.
\begin{thm}[Inverse Function Theorem with Parameters]\label{thm-inv-para}
Let $E$ be a locally convex space, $(F,\|\cdot\|)$ be a Banach space,
$U\sub E$ be a locally convex subset with dense interior,
$V\sub F$ be an open subset and $f\colon U\times V\to F$
be a $C^k$-function with $k\in\{2,3,\ldots\}\cup\{\infty\}$.
Abbreviate $f_x:=f(x,\cdot)\colon V\to F$ for $x\in U$.
Assume that $(x_0,y_0)\in U\times V$ with $f_{x_0}'(y_0)\in\GL(F)$.
Then there exist an open $x_0$-neighborhood $U_0\sub U$
and an open $y_0$-neighborhood $V_0\sub V$ with the following properties:
\begin{description}[(D)]
\item[\rm(a)]
$f_x(V_0)$ is open in~$F$ and $f_x|_{V_0}\colon V_0\to f_x(V_0)$ is
a $C^k$-diffeomorphism, for each $x\in U_0$;
\item[\rm(b)]
The set $W:=\bigcup_{x\in U_0}\{x\}\times f_x(V_0)$ is open in $U_0\times F$
and the mapping
\[
\psi\colon U_0\times V_0\to W,\quad (x,y)\mto (x,f_x(y))
\]
is a $C^k$-diffeomorphism.
\end{description}
\end{thm}
The proof uses part~(a) of the following lemma.
\begin{lem}\label{est-if-C2}
Let $E$ be a locally convex space, $(F,\|\cdot\|)$ be a normed space,
$U\sub E$ be a locally convex subset with dense interior,
$V\sub F$ be an open subset and $f\colon U\times V\to F$
be a function.
Abbreviate $f_x:=f(x,\cdot)\colon V\to F$ for $x\in U$.
Assume that
\begin{description}[(D)]
\item[\rm(a)]
$f$ is $C^2$; or
\item[\rm(b)]
$F$ is finite-dimensional and $f$ is~$C^1$ $($or at least $C^{0,1})$; or
\item[\rm(c)]
$E$ is normable and $f$ is $FC^1$.
\end{description}
Then $g\colon U\times V\to \cL(F)_b$, $(x,y)\mto (f_x)'(y)$
is a continuous map.
\end{lem}
\begin{prf}
(a) See Remark~\ref{back-to-normed}.

(b)
Let $v_1,\ldots, v_n$ be a basis for~$F$.
By Exercise~\ref{exc-on-fin}(b), the map
\[
\phi \colon \cL(F)_b\to F^n,\quad \alpha\mto (\alpha(v_1),\ldots,\alpha(v_n))
\]
is an isomorphism of topological vector spaces.
As the map
\begin{equation}\label{hence-indc}
\phi\circ g=(d_2f(\cdot,v_1),\ldots,d_2f(\cdot,v_n))\colon U\times V \to F^n
\end{equation}
is continuous, also $g=\phi^{-1}\circ(\phi\circ g)$
is continuous.

(c) Since $\lambda\colon F\to E\times F$, $y\mto (0,y)$ is continuous linear, the maps
\[ \cL(\lambda,F)\colon \cL(E\times F,F)_b\to\cL(F)_b,\quad 
\alpha\mto\alpha\circ\lambda \] 
and $g=\cL(\lambda,F)\circ f'$ are continuous.
\end{prf}
\noindent\emph{Proof of Theorem}~\ref{thm-inv-para}.
We may assume that $F\not=\{0\}$,
excluding a trivial case.
Let $\alpha:=(f_{x_0})'(y_0)$ and $0 < L < \frac{1}{\|\alpha^{-1}\|_{\op}}$;
thus
\[
a:=\frac{1}{\|\alpha^{-1}\|_{\op}}-L>0.
\]
By Lemma~\ref{est-if-C2}(a),
there exists an open $x_0$-neighborhood $U_0\sub U$
and $r>0$ with $V_0:=B^F_r(y_0)\sub V$ such that
\begin{equation}\label{need-ths}
\big( \forall \, (x,y)\in U_0\times V_0\big)\;\;
\|(f_x)'(y)-\alpha\|_{\op} = 
\|(f_x)'(y)-(f_{x_0})'(y_0)\|_{\op}\leq L .
\end{equation}
Then $g_x\colon V_0\to F$, $y\mto f_x(y)-\alpha(y)$ satisfies
$\|(g_x)'(y)\|_{\op}\leq L$ for all $y$ in the convex set~$V_0$, whence
$\Lip(g_x)\leq L$ by Lemma~\ref{lipviaprime}.

(a) Since $f_x|_{V_0}=\alpha+g_x$,
Theorem~\ref{lip-inv-fct} shows that $f_x(V_0)$ is open in~$F$ and $f_x|_{V_0}\colon V_0\to f_x(V_0)$
is a $C^k$-diffeomorphism, for each $x\in U_0$. As a consequence, the map $\psi$ in (b) is injective.

(b) Let $\wb{x}\in U_0$ and $\wb{y}\in V_0$.
We show that $W$ contains a neighborhood~$\Omega$ of $(\wb{x},f_x(\wb{y}))$
in $U_0\times F$ and that $\psi^{-1}$ is $C^k$ on this neighborhood.
Let \break $0 < t < r-\|\wb{y}-y_0\|$
and $Q$ be an open neighborhood of~$\wb{x}$ in~$U_0$
such that
\begin{equation}\label{gives-inside}
f_x(\wb{y})\in \wb{B}^F_{at/2}(f_{\wb{x}}(\wb{y}))\quad\mbox{for all $\,x\in Q$.}
\end{equation}
Then
\[
\wb{B}^F_{at}(f_x(\wb{y}))\supseteq\wb{B}^F_{at/2}(f_{\wb{x}}(\wb{y}))
\quad\mbox{for all $\,x\in Q$,}
\]
by (\ref{gives-inside}) and the triangle inequality.
Hence
\[
f_x(V_0)\supseteq f_x(\wb{B}^F_t(\wb{y}))\supseteq
\wb{B}^F_{at}(f_x(\wb{y}))\supseteq \wb{B}^F_{at/2}(f_{\wb{x}}(\wb{y}))
\]
for $x\in Q$,
using Theorem~\ref{lip-inv-fct}(e)
and Remark~\ref{L-not-lip} for the second inclusion.
Thus
\[
\Omega:=Q\times \wb{B}^F_{at/2}(f_{\wb{x}}(\wb{y}))\sub W.
\]
To see that $\psi^{-1}|_\Omega$ is~$C^k$,
recall from the proof of Theorem~\ref{lip-inv-fct}(e)
that $(f_x|_{V_0})^{-1}(c)$ is the unique fixed point of the contraction
\[
h_{x,c}\colon \wb{B}_t^F(\wb{y})\to \wb{B}^F_t(\wb{y}),
\quad z\mto z-\alpha^{-1}(f_x(z)-c)
\]
for $x\in Q$ and $c\in \wb{B}^F_{at}(f_x(\wb{y}))$,
which satisfies $\Lip(h_{x,c})\leq L\|\alpha^{-1}\|_{\op}$.
By the preceding,
$(h_{x,c})_{(x,c)\in \Omega}$
is a uniform family of contractions with uniform contraction constant $L \|\alpha^{-1}\|_{\op}<1$.
As the map
\[
h\colon \Omega\times \wb{B}^F_t(\wb{y})\to
\wb{B}^F_t(\wb{y}),\quad (x,c,z)\mto
z-\alpha^{-1}(f(x,z)-c)
\] 
is~$C^k$, Proposition~\ref{fp-Ck-dep}(a)
shows that the fixed point $(f_x|_{V_0})^{-1}(c)\in F$
is a $C^k$-function of $(x,c)\in \Omega$.\qed

\begin{rem}\label{prepa-impl}
Since $W$ is open in Theorem~\ref{thm-inv-para}(b),
it contains $Q\times Z$ for an open $x_0$-neighborhood $Q\sub U_0$
and an open neighborhood $Z\sub F$
of $f_{x_0}(y_0)$. After replacing~$U_0$ with~$Q$,
we may therefore assume that
\[U_0\times Z\sub W.\]
Thus $Z\sub f_x(V_0)$ for all $x\in U_0$, and the map
\[
U_0\times Z\to V_0,\quad (x,z)\mto (f_x|_{V_0})^{-1}(z)
\]
is $C^k$ (as it is the 2nd component of the $C^k$-map $\psi^{-1}|_{U_0\times Z}\colon
U_0\times Z\to U_0\times V_0$).
\end{rem}
\begin{rem}\label{variants-i-para}
(a) Theorem~\ref{thm-inv-para}
also holds for $k=1$ if
there exist an open $x_0$-neighborhood $U_0\sub U$
and $r>0$ with $V_0:=B^F_r(y_0)\sub V$ such that
\begin{equation}\label{need-ths2}
L:=\sup_{(x,y)\in U_0\times V_0}\|(f_x)'(y)-(f_{x_0})'(y_0)\|_{\op}
<\frac{1}{\|(f_{x_0})'(y_0)^{-1}\|_{\op}}.\vspace{-.7mm}
\end{equation}
This condition is automatic if~$F$ has finite dimension
or~$E$ is normable and~$f$ is~$FC^1$
(see Lemma~\ref{est-if-C2}(b) and (c), respectively).
If we assume (\ref{need-ths2}),
then (\ref{need-ths}) holds, 
and the remainder of the proof of Theorem~\ref{thm-inv-para}
can be repeated verbatim.\medskip

\noindent
(b) If~$E$ is normable and $k\in\N\cup\{\infty\}$,
then we can replace~$C^k$-maps with $FC^k$-maps in all hypotheses and conclusions
of Theorem~\ref{thm-inv-para}. In fact, we just observed
that (\ref{need-ths}) holds if~$f$ is $FC^k$ (and hence~$FC^1$).
The remainder of the proof of Theorem~\ref{thm-inv-para}
can now be repeated, writing $FC^k$ instead of~$C^k$,
and using Proposition~\ref{fp-Ck-dep}(b) instead of Proposition~\ref{fp-Ck-dep}(a).\medskip

\noindent
(c) If $U\sub E$ is open in Theorem~\ref{thm-inv-para}, then we can replace $C^k$-maps
with $\K$-analytic maps in all hypotheses and conclusions of the theorem.
This is trivial if $\K=\C$, as a mapping on an open subset
is $\K$-analytic if and only if it is~$C^\infty_\C$ (see Theorem~\ref{charcxcompl}).
If $\K=\R$ and $f$ is real analytic, after shrinking~$U$ and $V$ we may assume
that $f$ has a complex analytic extension
\[
\wt{f}\colon \wt{U}\times\wt{V}\to F_\C,
\]
where $\wt{U}\sub E_\C$ is an open $x_0$-neighborhood
and $\wt{V}\sub F_\C$ is an open $y_0$-neighborhood.
Since $(\wt{f}_{x_0})'(y_0)=(f_{x_0})'(y_0)_\C\in\GL(F_\C)$,
Theorem~\ref{thm-inv-para} provides an open $x_0$-neighborhood
$\wt{U}_0\sub\wt{U}$ and an open $y_0$-neighborhood $\wt{V}_0\sub\wt{V}$
such that
$\wt{W}:=\bigcup_{x\in\wt{U}_0}\{x\}\times \wt{f}_x(\wt{V}_0)$ is open in
$E_\C\times F_\C$ and
\[
\wt{\psi}\colon \wt{U}_0\times\wt{V}_0\to\wt{W},\quad (x,y)\mto\wt{f}_x(y)
\]
is a complex analytic diffeomorphism.
Applying Theorem~\ref{thm-inv-para} to the $C^\infty$-map
$f|_{(U\cap \wt{U}_0)\times (V\cap\wt{V}_0)}$,
we find an open $x_0$-neighborhood $U_0\sub U\cap \wt{U}_0$
and an open $y_0$-neighborhood $V_0\sub V\cap\wt{V}_0$
as described in the theorem.
Since $\wt{\psi}$ is a complex analytic extension of~$\psi$
and $\wt{\psi}^{-1}$ is a complex analytic extension of~$\psi^{-1}$,
we see that~$\psi$ and~$\psi^{-1}$ are real analytic.
Hence~$\psi$ is a real analytic diffeomorphism.
As $(\wt{f}_x|_{\wt{V}_0})^{-1}$ is a complex analytic extension
of $(f_x|_{V_0})^{-1}$, we see that $(f_x|_{V_0})^{-1}$
is real analytic for all $x\in U_0$, whence $f_x|_{V_0}\colon V_0\to f_x(V_0)$
is a real analytic diffeomorphism.
\end{rem}
Choosing $E=\{0\}$ and $x_0=0$ in Theorem~\ref{thm-inv-para}
and its analog for analytic functions described in Remark~\ref{variants-i-para}(c),
we obtain inverse function theorems for $C^k$-maps and
$\K$-analytic maps:
\begin{cor}\label{inv-fct-Ck}
Let $F$ be a Banach space over $\K\in\{\R,\C\}$
and $f\colon V\to F$ be a mapping defined on an open subset $V\sub F$.
Assume that $f$ is $C^k$ for some $k\in\{2,3,\ldots\}\cup\{\infty\}$
$($resp., $\K$-analytic$)$.
Let $y_0\in V$ such that $f'(y_0)\colon F\to F$ is invertible.
Then there exists an open $y_0$-neighborhood $V_0\sub V$ such that $f(V_0)$ is open in~$F$
and $f|_{V_0}\colon V_0\to f(V_0)$ is a $C^k$-diffeomorphism $($resp., a $\K$-analytic
diffeomorphism$)$.\qed
\end{cor}
\begin{rem}
Corollary~\ref{inv-fct-Ck} remains valid for $k=1$
whenever there exists a $y_0$-neighborhood $V_0\sub V$
with 
\[ \sup\{\|f'(y)-f'(y_0)\|_{\op}\colon y\in V_0\}<\frac{1}{\|f'(y_0)^{-1}\|_{\op}}\]
(cf.\ Remark~\ref{variants-i-para}(a)).
Without such an extra condition, the inverse function
theorem may fail for $C^1$-maps, see Exercise~\ref{exc-C1notFC1}.
\end{rem}
Implicit function theorems also follow immediately from Theorem~\ref{thm-inv-para}.
\begin{thm}\label{impl-Ck}
Let $E$ be a locally convex space and $F$ be a Banach space over $\K\in\{\R,\C\}$,
$U\sub E$ be a locally convex subset with dense interior and $V\sub F$ be an open subset.
Let $f\colon U\times V\to F$ be a $C^k$-map with $k\in\{2,3,\ldots\} \cup\{\infty\}$
$($resp., a $\K$-analytic map$)$
and $(x_0,y_0)\in U\times V$ such that $f(x_0,y_0)=0$ and $(f_{x_0})'(y_0)\colon F\to F$ is invertible,
where $f_x:=f(x,\cdot)\colon V\to F$ for $x\in U$.
Then there exists an open $x_0$-neighborhood $U_0\sub U$, an open $y_0$-neighborhood
$V_0\sub V$ and a $C^k$-function $($resp., $\K$-analytic function$)$ $\phi\colon U_0\to V_0$ such that
\begin{equation}\label{sol-is-graph-2}
\{(x,y)\in U_0\times V_0\colon f(x,y)=0\}=\graph(\phi).
\end{equation}
\end{thm}
\begin{prf}
Let $U_0$, $V_0$, $W$, and $\psi$ be as in Theorem~\ref{thm-inv-para};
as observed in Remark~\ref{prepa-impl}, we may assume that
$U_0\times Z\sub W$ for some open $0$-neighborhood $Z\sub F$.
Hence $f_x|_{V_0}$ has a zero $\phi(x)\in V_0$ for each $x\in U_0$,
and the latter is unique as $f_x|_{V_0}$ is injective.
Since
\[
\psi(x,\phi(x))=(x,0),
\]
we see that $\phi(x)$ is the second component of $\psi^{-1}(x,0)\in U_0\times V_0$
and hence (like $\psi^{-1}(x,0)$) a $C^k$-function
of $x\in U_0$ (resp., a $\K$-analytic function, cf.\ Remark~\ref{variants-i-para}(c)).
\end{prf}
\begin{rem} \label{rem:2.3.44}
Theorem~\ref{impl-Ck}
also holds for $k=1$ if we assume that
there exists an open $x_0$-neighborhood $U_0\sub U$
and $r>0$ with $V_0:=B^F_r(y_0)\sub V$ such that
(\ref{need-ths2}) holds.
In fact, the preceding proof can be repeated in this case,
using the $C^1$-version of Theorem~\ref{thm-inv-para}
described in Remark~\ref{variants-i-para}(a).
\end{rem}
\noindent\emph{Proof of Theorem}~\ref{impl-fct-class}.
Replacing Theorem~\ref{thm-inv-para} with its $FC^k$-variant
described in Remark~\ref{variants-i-para}(b),
we can repeat the proof of Theorem~\ref{impl-Ck} verbatim.\qed
\subsection*{Global inverse function theorems}
Theorem~\ref{lip-inv-fct}
immediately entails a global inverse function theorem.
\begin{thm}\label{glob-lip-inv}
Let $(E,\|\cdot\|)$ be a Banach space,
$\alpha\in\GL(E)$ and $g\colon E\to E$ be Lipschitz
with
\[
\Lip(g)<\frac{1}{\|\alpha^{-1}\|_{\op}}.
\]
Then $f:=\alpha+g\colon E\to E$ is a homeomorphism
and the inverse mapping\linebreak
$f^{-1}\colon E\to E$ is Lipschitz with
\begin{equation}\label{norber}
\Lip(f^{-1})\leq a^{-1} \quad \mbox{ for } \quad a:= {\frac{1}{\|\alpha^{-1}\|_{\op}}-\Lip(g)}.
\end{equation}
Moreover, $h:=f^{-1}-\alpha^{-1}$ is Lipschitz, $f^{-1}=\alpha^{-1}+h$ and
\begin{equation}\label{nobbi2}
\Lip(h)\leq\frac{\|\alpha^{-1}\|_{\op}\Lip(g)}{a}.
\end{equation}
If $g$ is $C^k$ for some $k\in\N$ $($resp., $FC^k$, resp., $\K$-analytic$)$,
then also $f^{-1}$ is $C^k$ $($resp., $FC^k$, resp., $\K$-analytic$)$.
\end{thm}
\begin{prf}
Let $L:=\Lip(g)$, $g_r:=g|_{B^E_r(0)}$ and $f_r:=f|_{B^E_r(0)}$ for $r>0$.
Applying Theorem~\ref{lip-inv-fct} to $g_r$ and $f_r$,
we see that $f_r$ is injective with
\begin{equation}\label{lippi}
\Lip((f_r)^{-1})\leq a^{-1}
\end{equation}
and $\im(f_r)\supseteq B^E_{ar}(f(0))$ (see also Remark~\ref{L-not-lip}).
Given $x\not=y$ in~$E$, we have $x,y\in B^E_r(0)$ for some $r>0$
and obtain $f(x)=f_r(x)\not=f_r(y)=f(y)$; thus $f$ is injective.
As the sets $B^E_{ar}(f(0))$ cover~$E$ for $r > 0$, 
we see that~$f$ is surjective.
Given $x,y\in E$, we have $x,y\in \im(f_r)$ for some $r>0$.
Since
\[
\|f^{-1}(x)-f^{-1}(y)\|=\|(f_r)^{-1}(x)-(f_r)^{-1}(y)\|\leq\Lip((f_r)^{-1})\|x-y\|,
\]
we deduce with (\ref{lippi}) that (\ref{norber}) holds.
Likewise, (\ref{nobbi2}) follows from the estimate
\[
\Lip((f_r)^{-1}-\alpha^{-1})
\leq\frac{\|\alpha^{-1}\|_{\op}L}{a} 
\]
(the analog of (\ref{the-lip-h}) described in
Remark~\ref{L-not-lip}).
If $g$ is $C^k$ and $FC^k$, respectively, then $(f_r)^{-1}$
is $C^k$ and $FC^k$, respectively (see Theorem~\ref{lip-inv-fct})
for each $r>0$, whence $f^{-1}$ is $C^k$ (resp., $FC^k$).
If $g$ is real analytic, then the invertibility of $f'(x)$
implies that~$f$ is a local analytic diffeomorphism at~$x$ (see Corollary~\ref{inv-fct-Ck}),
whence $f^{-1}$ is real analytic on an open neighborhood of~$f(x)$.
As a consequence, $f^{-1}$ is real analytic.
\end{prf}
If $g$ is $C^1$ in the situation of Theorem~\ref{glob-lip-inv},
then
\[
f'(x)=\alpha+g'(x)=\alpha\circ (\id_E+\alpha^{-1}\circ g'(x))
\]
is invertible for all $x\in E$ as
\[
\|\alpha^{-1}\circ g'(x)\|_{\op}\leq\|\alpha^{-1}\|_{\op}\|g'(x)\|_{\op}
\leq\|\alpha^{-1}\|_{\op}\Lip(g)<1.
\]
Moreover,
\[
\|f'(x)^{-1}\|_{\op}\leq\|(\id_E+\alpha^{-1}\circ g'(x))^{-1}\|_{\op}\|\alpha^{-1}\|_{\op}
\leq \frac{\|\alpha^{-1}\|_{\op}}{1-\|\alpha^{-1}\|_{\op}\Lip(g)}
\]
is a bounded function of~$x$
(using Exercise~\ref{excBancia}(c) for the second estimate).
The main conclusion of Theorem~\ref{glob-lip-inv}
(that $f$ is a bijection)
therefore is a special case of the following
version of Hadamard's Global Inverse Function Theorem,
which we mention for completeness.
We recommend to skip the proof on a first reading,
as Theorem~\ref{glob-lip-inv} usually suffices
in Lie-theoretic studies of diffeomorphism groups.
\begin{thm}\label{glob-hadamard}
Let $(E,\|\cdot\|)$ be a Banach space
and $q\colon E\to E$ be a local $C^1$-diffeomorphism
such that
\begin{equation}\label{prereq-hadam}
L:=\sup_{x\in E}\|q'(x)^{-1}\|_{\op}<\infty.
\end{equation}
Then $q(E)=E$ and $q$ is a $C^1$-diffeomorphism.
\end{thm}
We shall use a fact concerning uniqueness of lifts.
Recall that a mapping $q\colon X\to Y$ between topological spaces
is said to be a \emph{local homeomorphism} \index{local homeomorphism}
if each $x\in X$ has an open neighborhood $V\sub X$
such that $q(V)$ is open in~$Y$ and $q|_V\colon V\to q(V)$
is a homeomorphism. 
\begin{lem}\label{uni-lft}
Let $q\colon Z\to Y$ be a local homeomorphism between topological spaces.
Let $X$ be a topological space and
$\gamma\colon X\to Z$ as well as $\eta\colon X\to Z$ be continuous mappings
such that
\[
q\circ \gamma=q\circ\eta.
\]
If $X$ is connected, $Z$ is Hausdorff and $\gamma(x_0)=\eta(x_0)$ for some $x_0\in X$,
then $\gamma=\eta$.
\end{lem}
\begin{prf}
The set $E:=\{x\in X\colon \gamma(x)=\eta(x)\}$ is closed (as $\gamma$ and $\eta$
are continuous and $Z$ is Hausdorff) and not empty, since $x_0\in E$.
If we can show that~$E$ is open, then $E=X$ as~$X$ is connected, and thus $\gamma=\eta$.
Given $x\in E$, let $V\sub Z$ be a neighborhood of $\gamma(x)=\eta(x)$ such that
$q|_V$ is injective. Since~$\gamma$ and~$\eta$ are continuous,
there exists an $x$-neighborhood $U\sub X$ such that $\gamma(U)\sub V$
and $\eta(U)\sub V$. Now $q|_V \circ \gamma|_U=q|_V\circ\eta|_U$ implies
$\gamma|_U=\eta|_U$, whence $U\sub E$. Thus~$E$ is a neighborhood
of~$x$ in~$X$ and thus $E$ is open.
\end{prf}
\begin{lem}\label{smooth-lft}
Let $X$, $Y$, and $Z$ be locally convex spaces,
$U\sub X$ be a locally convex subset with dense interior,
$V\sub Y$ and $W\sub Z$ be open subsets and $q\colon W\to V$ be a $C^k$-diffeomorphism,
where $k\in\N\cup\{\infty,\omega\}$
and $U$ is assumed open in~$X$ if $k=\omega$.
If $\gamma\colon U\to W$ is a continuous map such that $q\circ \gamma$ is~$C^k$,
then also $\gamma$ is $C^k$, and
\begin{equation}\label{lift-will-lip}
\gamma'(x)=q'(\gamma(x))^{-1} \circ (q\circ \gamma)'(x)\quad
\mbox{for all $x\in U$.}
\end{equation}
\end{lem}
\begin{prf}
Given $x\in U$, let $W_0\sub W$ be an open neighborhood of $\gamma(x)$
such that $q(W_0)$ is open in~$Y$ and $q|_{W_0}\colon W_0\to q(W_0)$ is a $C^k$-diffeomorphism.
Let $U_0\sub U$ be an open neighborhood of~$x$ such that $\gamma(U_0)\sub W_0$.
Then
\begin{equation}\label{pre-equ}
\gamma|_{U_0}=(q|_{W_0})^{-1}\circ (q \circ \gamma)|_{U_0} 
\end{equation}
is $C^k$ (and hence also~$\gamma$). Applying the Chain Rule to (\ref{pre-equ})
and using $y:=\gamma(x)$ in
\begin{equation}
((q|_{W_0})^{-1})'(q(y))=q'(y)^{-1}\quad\mbox{for all $y\in W_0$,}
\end{equation}
we obtain~(\ref{lift-will-lip}).
\end{prf}
The key to Proposition~\ref{glob-hadamard}
is the existence of lifts for $C^1$-curves and
homotopies between such.
\begin{lem}\label{lifts-for-hadam}
In the situation of Theorem~\emph{\ref{glob-hadamard}},
the following holds:
\begin{description}[(DD)]
\item[\rm(a)]
If $\gamma\colon [0,1]\to E$ is a $C^1$-curve
and $z_0\in E$ such that $q(z_0)=\gamma(0)$,
then there exists a $C^1$-curve $\eta\colon [0,1]\to E$
such that $q\circ \eta=\gamma$ and $\eta(0)=z_0$.
\item[\rm(b)]
If $F \colon [0,1]\times [0,1]\to E$ is a $C^1$-map
and $z_0\in E$ such that $q(z_0)=F(0,0)$,
then there exists a $C^1$-map $G \colon [0,1]\times [0,1] \to E$
such that $q\circ G= F$ and $G(0,0)=z_0$.
\end{description}
\end{lem}
\begin{prf}
(a) Let $r\in [0,1[$ and $\eta \colon [0,r] \to E$ be a continuous map
such that $q\circ\eta=\gamma|_{[0,r]}$ and $\eta(0)=z_0$
(this situation can always be achieved if we take $r:=0$).
There exists an open neighborhood $W_0\sub E$ of $\eta(r)$
such that $q(W_0)$ is open in~$E$ and $q|_{W_0}\colon W_0\to q(W_0)$
is a homeomorphism. Since $q(W_0)$ is a neighborhood of $q(\eta(r))=\gamma(r)$,
there exists $\rho \in\;]r,1]$ such that
$\gamma([r,\rho])\sub q(W_0)$.
Setting
\[
\eta(t):=(q|_{W_0})^{-1}(\gamma(t))
\]
for $t\in [r,\rho]$, we obtain a continuous extension of $\eta$ to a
mapping $[0,\rho]\to E$ such that $\eta(0)=z_0$ and $q\circ \eta=\gamma|_{[0,\rho]}$.

By the preceeding, for some $r\in\;]0,1]$
there exists a continuous map $\eta_r\colon [0,r[\;\to E$
such that $q\circ\eta_r=\gamma|_{[0,r[}$ and $\eta_r(0)=z_0$.
Lemma~\ref{uni-lft} entails that~$r$ can be chosen maximal.
Let
\[
M:=\max\{\|\gamma'(t)\|\colon t\in[0,r]\}.
\]
Lemma~\ref{smooth-lft} entails that $\eta_r$ is~$C^1$, with
\[
\|\eta_r'(t)\|=\|(q(\eta_r(t))')^{-1}(\gamma'(t))\|\leq LM
\]
for all $t\in[0,r[$.
As a consequence, $\eta_r$ is Lipschitz (cf.\ Lemma~\ref{lipviaprime})
and hence uniformly continuous. As $(E,\|\cdot\|)$ is complete,
$\eta_r$ has a uniformly continuous extension $\eta\colon [0,r]\to E$.
Then $q\circ \eta=\gamma|_{[0,r]}$, since equality holds on the dense subset $[0,r[$.
If we had $r<1$, then~$\eta$ could be extended to $[0,\rho]$ for some
$\rho \in\;]r,1]$ (as shown at the beginning of the proof),
contradicting the maximality of~$r$. Hence $r=1$,
and it only remains to observe that~$\eta$ is~$C^1$,
by Lemma~\ref{smooth-lft}.

(b) Let $B_\ve(s,t):=\{(a,b)\in[0,1]^2\colon |a-s|,|b-t|<\ve\}$
for $(s,t)\in [0,1]^2$ and $\ve>0$.
Let $r\in [0,1[$ and $G \colon [0,r]\times [0,1]\; \to E$ be a continuous map
such that $q\circ G =F |_{[0,r]\times [0,1]}$ and $G(0,0)=z_0$
(by~(a), this situation can always be achieved if we take $r:=0$).
For each $t\in[0,1]$,
there exists an open neighborhood $W_t\sub E$ of $G(r,t)$
such that $q(W_t)$ is open in~$E$ and $q|_{W_t}\colon W_t\to q(W_t)$
is a homeomorphism. Since $q(W_t)$ is a neighborhood of $q(G(r,t))=F(r,t)$,
there exists $\ve_t>0$ such that $F(V_t)\sub q(W_t)$
holds for $V_t:=B_{\ve_t}(r,t)$.
Then
\[
G^t:=(q|_{W_t})^{-1}\circ F|_{V_t}
\]
is a continuous map which coincides with~$G$ on the convex (and hence connected) set
$V_t\cap ([0,r]\times[0,1])$
as $G^t(r,t)=G(r,t)$
(see Lemma~\ref{uni-lft}).
If $t,\tau\in[0,1]$ and $V_t\cap V_\tau\not=\emptyset$,
then $V_t\cap V_\tau$ coincides with
the convex set
\[
\{a\in [0,1]\colon |a-r|<\min\{\ve_t,\ve_\tau\}\}\times
\{b\in[0,1]\colon |b-t|<\ve_t\;\mbox{and}\;|b-\tau|<\ve_\tau\},
\]
which contains $(r,b)$ for some $b\in[0,1]$.
Since $G^t(r,b)=G(r,b)=G^\tau(r,b)$, Lemma~\ref{uni-lft}
shows that $G^t|_{V_t\cap V_\tau}=G^\tau|_{V_t\cap V_\tau}$.
Using the Wallace Lemma (Lemma~\ref{Wallla}), we find $\rho\in \;]r,1]$ such that
\[
[0,\rho]\times [0,1]\sub ([0,r]\times [0,1])\cup \bigcup_{t\in[0,1]}V_t=([0,r[\;\times [0,1])
\cup\bigcup_{t\in[0,1]}V_t.
\]
Then $G(a,b):=G^t(a,b)$ for $t\in[0,1]$ and $(a,b)\in V_t\cap ([0,\rho]\times[0,1])$
defines a continuous extension of~$G$ to a map $[0,\rho]\times[0,1]\to E$
such that $q\circ G=F|_{[0,\rho]\times[0,1]}$.

By the preceeding, for some $r\in\;]0,1]$
there exists a continuous map $G_r\colon [0,r[\,\times [0,1]\to E$
such that $q\circ G_r=F |_{[0,r[\,\times[0,1]}$ and $G_r(0,0)=z_0$.
Lemma~\ref{uni-lft} entails that~$r$ can be chosen maximal.
Let
\[
M:=\max\{\|F'(s,t)\|_{\op}\colon (s,t)\in[0,r]\times[0,1]\}.
\]
Lemma~\ref{smooth-lft} entails that $G_r$ is~$C^1$, with
\[
\|G_r'(s,t)\|_{\op}=\|q'(G_r(s,t))^{-1}\circ F'(s,t)\|_{\op}\leq LM
\]
for all $(s,t)\in[0,r[\, \times [0,1]$.
As a result, $G_r$ is Lipschitz (see Lemma~\ref{lipviaprime})
and hence uniformly continuous. As $(E,\|\cdot\|)$ is complete,
$G_r$ has a uniformly continuous extension $G\colon [0,r]\times[0,1]\to E$.
Then $q\circ G=F|_{[0,r]\times[0,1]}$, since equality holds on the dense subset $[0,r[\,\times[0,1]$.
If we had $r<1$, then~$G$ could be extended to $[0,\rho]\times[0,1]$ for some
$\rho\in\;]r,1]$ (as already shown),
contradicting the maximality of~$r$. Hence $r=1$,
and it only remains to observe that~$G$ is~$C^1$,
by Lemma~\ref{smooth-lft}.
\end{prf}
\noindent\emph{Proof of Theorem}~\ref{glob-hadamard}.
To see that $q$ is surjective, let $y\in E$.
Let $z_0\in E$ and consider the $C^1$-curve
\[
\gamma\colon [0,1]\to E,\quad t\mto q(z_0)+t(y-q(z_0)).
\]
Let $\eta\colon [0,1]\to E$ be as in Lemma~\ref{lifts-for-hadam}(a).
Then $q\circ\eta=\gamma$, whence $y=\gamma(1)=q(\eta(1))$ is in the image of~$q$.

To see that $q$ is injective, let $z_0,z_1\in E$ such that $q(z_0)=q(z_1)$.
Then
\[
\eta\colon [0,1] \to E,\quad t\mto z_0+t(z_1-z_0)
\]
is a $C^1$-curve from $z_0$ to~$z_1$,
entailing that $\gamma:=q\circ\eta$ is a loop.
Using $x_0:=\gamma(0)=\gamma(1)$, consider the $C^1$-map
\[
F\colon [0,1]\times[0,1]\to E,\quad (s,t)\mto x_0+s(\gamma(t)-x_0).
\]
By Lemma~\ref{lifts-for-hadam}(b), there is a $C^1$-map
$G\colon [0,1]\times[0,1]\to E$ such that
\[
q\circ G=F\quad\mbox{and}\quad G(0,0)=z_0.
\]
Since $F(0,t)=x_0$ for all $t\in[0,1]$, we must have $G(0,t)=z_0$ for all~$t$,
by Lemma~\ref{uni-lft}.
Likewise, $F(s,0)=F(s,1)=x_0$ and $G(0,0)=G(0,1)=z_0$ entail that $G(s,0)=G(s,1)=z_0$ for all $s\in[0,1]$.
But $F(1,t)=\gamma(t)$ for all $t\in[0,1]$ and $G(1,0)=z_0$
entail $G(1,t)=\eta(t)$ for all $t\in[0,1]$.
Thus $z_1=\eta(1)=G(1,1)=z_0$.\qed
\begin{small}
\subsection*{Exercises for Section~\ref{provis-Ban}}

\begin{exer}\label{exc-multil-exp}
Given $n,m\in\N$, let $(E_j,\|\cdot\|_j)$ for $j\in\{1,\ldots, n+m\}$
and $(F,\|\cdot\|_F)$ be normed spaces.
\begin{description}[(D)]
\item[(a)]
Show that $\beta^\vee(x):=\beta(x,\cdot)\in\cL^m(E_{n+1},\ldots, E_{n+m};F)$
for all
$x=(x_1,\ldots, x_n)$ in $E_1\times\cdots\times E_n$ and
$\beta\in\cL^{n+m}(E_1,\ldots, E_{n+m};F)$, with
\[
\|\beta^\vee(x)\|_{\op}\leq \|\beta\|_{\op}\|x_1\|_1\ldots\|x_n\|_n.
\]
Show that $\beta^\vee\in\cL^n(E_1,\ldots, E_n;\cL^m(E_{n+1},\ldots, E_{n+m};F))$
and $\|\beta^\vee\|_{\op}=\|\beta\|_{\op}$.
\item[(b)]
Let
$\ve_n$ be the evaluation map with domain
\[
\cL^n(E_1,\ldots, E_n;\cL^m(E_{n+1},\ldots, E_{n+m};F))\times E_1\times\cdots
\times E_n
\]
and range
$\cL^m(E_{n+1},\ldots, E_{n+m};F)$, which
is a continuous $(n+1)$-linear map by Remark~\ref{multileval}.
Let
$\delta_m\colon \cL^m(E_{n+1},\ldots,E_{n+m};F)\times E_{n+1}\times\cdots\times E_{n+m}\to F$
be the evaluation map, which is continuous $(m+1)$-linear.
Define
\begin{eqnarray*}
\wh{\beta}(x_1,\ldots,x_{n+m}) &:=& \beta(x_1,\ldots, x_n)(x_{n+1},\ldots x_{n+m})\\
&=&
\delta_m(\beta(x_1,\ldots, x_n),x_{n+1},\ldots, x_{n+m})\\
&=&\delta_m(\ve_n(\beta, x_1,\ldots, x_n),x_{n+1},\ldots, x_{n+m})
\end{eqnarray*}
for $\beta$ in $\cL^n(E_1,\ldots, E_n;\cL^m(E_{n+1},\ldots,E_{n+m};F))$ and $(x_1,\ldots,x_{n+m})$ in
the\linebreak
set $E_1\times\cdots\times E_{n+m}$.
Show that $\wh{\beta}\in\cL^{n+m}(E_1,\ldots, E_{n+m};F)$
and $(\wh{\beta})^\vee=\beta$.
\end{description}
By (a) and (b), the map
\[
\Phi\colon \cL^{n+m}(E_1,\ldots, E_{n+m};F)\to\cL^n(E_1,\ldots, E_n;\cL^m(E_{n+1},\ldots, E_{n+m});F),
\; \beta\mto\beta^\vee
\]
is an isometric isomorphism whose inverse mapping is given by $\Phi^{-1}(\beta)=\wh{\beta}$
for\linebreak
$\beta\in\cL^n(E_1,\ldots, E_n;\cL^m(E_{n+1},\ldots, E_{n+m};F))$.
Now let $(F, \|\cdot\|_F)$ be a Banach space.
\begin{description}[(D)]
\item[(c)]
Using the well-known fact that $(\cL(E,F),\|\cdot\|_{\op})$
is a Banach space for each normed space $(E,\|\cdot\|_E)$,
show by induction that $(\cL^n(E_1,\ldots, E_n;F),\|\cdot\|_{\op})$
is a Banach space for all $n\in\N$ and all normed spaces $(E_1,\|\cdot\|_1),\ldots,(E_n,\|\cdot\|_n)$.
\end{description}
\end{exer}

\begin{exer}
Let $(E,\|\cdot\|_E)$ and $(F,\|\cdot\|_F)$ be normed spaces and $U\sub E$
be a locally convex subset with dense interior.
Let us define total differentiability of~$f$ at~$x\in U$
by the requirements of Definition~\ref{defn-todiff}
also in this more general situation.
\begin{description}[(D)]
\item[(a)]
Show: If $f$ is totally differentiable at $x\in U$,
then $f$ is continuous at~$x$
and $f'(x)\in\cL(E,F)$ is uniquely determined.
\item[(b)]
Show that $f$ is $FC^1$ if and only if~$f$ is totally
differentiable at each $x\in U$ and the map $f'\colon U\to\cL(E,F)$,
$x\mto f'(x)$ is continuous.\\[1mm]
[Compare the proof of Lemma~\ref{FC1vsC1} for the necessity of this condition.]
\end{description}
\end{exer}

\begin{exer}
Let $(E,\|\cdot\|_E)$
and $(F,\|\cdot\|_F)$
be normed spaces and $U\sub E$
be an arbitrary subset with dense interior.
Let us say that a map $f\colon U\to F$ is $FC^1$
if there exists a continuous map $f'\colon U\to \cL(E,F)$
such that, for each $x\in U$, the map~$f$ admits a linear approximation~(\ref{linapprox})
satisfying~(\ref{little-o}), as in
Definition~\ref{defn-todiff}.
\begin{description}[(D)]
\item[(a)]
Show that, if $f\colon U\to F$ is $FC^1$, then
$f'$ is unique. Moreover, $f$ is continuous at each $x\in U$.
\item[(b)]
Assume that also $(Y,\|.\|_Y)$ is a normed space, $V\sub F$ a subset with dense
interior such that $f(U)\sub V$ and both $f$ and a map $g\colon V\to Y$ are $FC^1$.
Show that $g\circ f\colon U\to Y$ is $FC^1$ with $(g\circ f)'(x)=g'(f(x))\circ f'(x)$
for all $x\in U$.
\end{description}
\end{exer}

\begin{exer}\label{uniq-fixed-gener} (Uniqueness of fixed points).
Let $(X,d)$ be a metric space, $U\sub X$ be a subset,
$f\colon U\to X$ be a mapping with $\theta:=\Lip(f)<1$
and $x,y\in U$ such that $f(x)=x$ and $f(y)=y$. Show that $x=y$.
\end{exer}

\begin{exer}\label{contr-onlyself}
Let $(X,d)$ be a metric space, $U\sub X$ be an open subset,
$P$ be a topological space and $f\colon P\times U\to X$
be a continuous mapping such that the maps $f_p:=f(p,\cdot)\colon U\to X$
satisfy
\[
\theta:=\sup\{\Lip(f_p)\colon p\in P\}<1.
\]
Let $P_0$ be the set of all $p\in P$ such that
$f_p$ has a fixed point~$x_p$
(which is necessarily unique, by the preceding exercise).
\begin{description}[(D)]
\item[(a)] (Reduction to Self-Maps).
Given $p\in P_0$,
let $\ve>0$ such that $V:=\wb{B}^d_\ve(x_p)\sub U$.
Let $Q\sub P$ be a neighborhood of~$p$ such that $d(f(q,x_p),x_p)\leq (1-\theta)\ve$
for all $q\in Q$.
Show that $f_q(V)\sub V$ for all $q\in Q$.
\item[(b)]
If $(X,d)$ complete, deduce from~(a) that~$P_0$ is open in~$P$
and use Proposition~\ref{fp-cts-dep}(a)
to see that the map $\phi\colon P_0\to U$, $p\mto x_p$ 
is continuous.
\end{description}
\end{exer}

\begin{exer}\label{exerc-fp-an1}
Let $P$ be an open subset of a locally convex space~$E$.
Let~$V$ be an open subset of a Banach space $(F,\|\cdot\|)$
and $f\colon P\times V\to V$ be a real analytic map.
Write $f_p:=f(p,\cdot)\colon V\to V$ for $p\in P$ and assume that
$(f_p)_{p\in P}$ is a uniform family of contractions.
Let~$P_0$ be the set of all $p\in P$ for which $f_p$ has a fixed point~$x_p$.
By the preceding exercise, $P_0$ is open in~$P$.
Show that the map $\phi\colon P_0\to V$, $p\mto x_p$ is real
analytic.\\[1.3mm]
[We have $\|f'_p(x)\|_{\op}<1$ for all $p\in P$ and $x\in V$.
For $p\in P_0$,
the fixed point $x_p$ is the unique solution $x$ to the equation
$x-f(p,x)=0$, which can be discussed using the real analytic version of
the Implicit Function Theorem.\,]
\end{exer}

\begin{exer}\label{exerc-cx-norm}
If $E$ is a real vector space and~$\|\cdot\|$
a norm thereon, let $B:=B^E_1(0)$ be the open unit ball.
Recall that $\bD:=\{z\in\C\colon |z|\leq 1\}$.
Within the complex vector space~$E_\C$, define
\[
U:=\abscnv(B)=\cnv (\bD B).
\]
\begin{description}[(D)]
\item[(a)]
Since $\bD B\sub B+iB$ and the latter set is convex,
we have
\begin{equation}\label{half-prod}
U\sub B+ iB.
\end{equation}
Deduce that $U\cap E\sub B$ and thus $U\cap E=B$.
\item[(b)]
To see that $U$ is open, let $x\in U$.
There are $n\in\N$, $z_1,\ldots, z_n\in \bD$,
$x_1,\ldots, x_n\in B$
and $t_1,\ldots, t_n\in[0,1]$ with $t_1+\cdots+t_n=1$
such that
$x=t_1z_1x_1+\cdots+ t_nz_nx_n$. Since~$B$ is open, there is $r>1$ such that
$rx_1,\ldots, r x_n\in B$. Deduce that $rx\in U$ and infer with Lemma~\ref{baseconvex}(d)
that $x\in U^0$.
\item[(c)]
Show that $\frac{1}{2}B+\frac{1}{2}iB\sub U$.
\item[(d)]
Deduce from (\ref{half-prod})
and (c) that the Minkowski functional $\|\cdot\|_\C:=\mu_U$
of~$U$ is a norm on the complex vector space~$E_\C$
which is equivalent to the maximum norm on the real vector space
$E_\C=E\oplus iE$ (and hence defines the product topology).
\item[(e)] Show that $\|x\|_\C=\|x\|$
for each $x\in E$.
\item[(f)]
Let $\|\cdot\|'$ be any norm
on the complex vector space $E_\C$
such that $\|x\|'=\|x\|$ for each $x\in E$.
Let $C$ be the open unit ball in $(E,\|\cdot\|')$.
Show that $B\sub C$, $U\sub C$, and $\|\cdot\|'\leq\|\cdot\|_\C$.\\[.8mm]
We call $\|\cdot\|_\C$ the \emph{maximal complexification}
of~$\|\cdot\|$.
\item[(g)]
Let also $(F,\|\cdot\|_F)$ be a normed vector space over~$\R$
and $\alpha\colon E\to F$ be a continuous linear map.
Show that $\|\alpha_\C\|_{\op}=\|\alpha\|_{\op}$
if we consider $\alpha_\C\colon E_\C\to F_\C$, $x+iy\mto \alpha(x)+i\alpha(y)$
as a complex linear map from $(E_\C,\|\cdot\|_\C)$
to $(F_\C,(\|\cdot\|_F)_\C)$.
\item[(h)]
Let $\beta\colon E_1\times E_2\to F$ be a continuous bilinear map between normed vector spaces
over~$\R$ and $\beta_\C\colon (E_1)_\C\times (E_2)_\C\to F_\C$
be its unique complex bilinear extension. Show that $\|\beta_\C\|_{\op}=\|\beta\|_{\op}$
holds when the maximal complexified norms are used on the complexified vector spaces.
Formulate and prove an analogous result for continuous $n$-linear maps with $n\in\N$. 
\end{description}
\end{exer}

\begin{exer}\label{exerc-fp-an2}
Let $P$ be an open subset of a locally convex space~$E$.
Let $V$ be an open subset of a Banach space $(F,\|\cdot\|)$,
and $f\colon P\times V\to F$ be a real analytic map.
Let $p\in P$ and $x_p\in V$ be a fixed point of $f_p:=f(p,\cdot)$
such $\|f_p'(x_p)\|_{\op}<1$.
After shrinking the $p$-neighborhood $P\sub E$ and the $x_p$-neighborhood
$V\sub F$, we may assume that $f$ has a complex analytic extension
$\wt{f}\colon\wt{P}\times\wt{V}\to F_\C$
for an open neighborhood $\wt{P}$ of~$P$ in~$E_\C$
and an open neighborhood $\wt{V}$ of~$V$ in~$F_\C$.
\begin{description}[(D)]
\item[(a)]
Show that, after shrinking the neighborhoods if necessary,
we can assume that there exists $\theta\in [0,1[$
such that the map $f_q:=f(q,\cdot)\colon V\to F$
satisfies $\|f_q'(x)\|_{\op}\leq\theta$ for all $(q,x)\in P\times V$
and the map $\wt{f}_q:=\wt{f}(q,\cdot)\colon \wt{V}\to F_\C$
satisfies $\|\wt{f}_q'(x)\|_{\op}\leq\theta$
for all $(q,x)\in \wt{P}\times\wt{V}$, using the maximal complexification
$\|\cdot\|_\C$ of the norm~$\|\cdot\|$ on~$F_\C$.
\item[(b)]
Show that, moreover, we can achieve that there is $r>0$
with $W:=\wb{B}^F_r(x_p)\sub V$ and $\wt{W}:=\wb{B}^{F_\C}_r(x_p)\sub\wt{V}$
such that $f_q(W)\sub W$ for all $q\in P$ and $\wt{f}_q(\wt{W})\sub\wt{W}$
for all $q\in\wt{P}$.
\item[(c)]
As a consequence of~(b), $f_q|_W$ has a unique fixed point $\phi(q):=x_q$
for each $q\in P$ and $\wt{f}_q|_{\wt{W}}$ has a unique fixed point $\wt{\phi}(q)$
for each $q\in\wt{P}$.
Use Proposition~\ref{fp-Ck-dep} to see that $\wt{\phi}\colon\wt{P}\to F_\C$ is $C^\infty_\C$
and thus complex analytic,
and deduce that $\phi\colon P\to F$ is real analytic.
\end{description}
Compare the current alternative approach with Exercise~\ref{exerc-fp-an1}.
\end{exer}

\begin{exer}\label{exer-stri1}
Let $(E,\|\cdot\|)$ be a normed space
and $U\sub E$ be an open subset. A function $f\colon U\to F$ to a normed space $(F,\|\cdot\|_F)$
is called \emph{strictly differentiable}
\index{strictly differentiable map} 
\index{map!strictly differentiable} 
at a point $x\in U$
if there exists a continuous linear map $f'(x)\colon E\to F$
such that the unique function $R\colon U\to F$ with
\begin{equation}\label{strict-condi}
(\forall y\in U)\;\, f(y)=f(x)+f'(x)(y-x)+R(y)
\end{equation}
satisfies $\Lip(R|_{B^E_r(x)})\to 0$ as $r\to 0$. Then $R(x)=0$
and $f$ is totally differentiable at~$x$ with derivative~$f'(x)$, whence $f'(x)$ is unique.
Moreover,
\begin{equation}\label{using-g}
f(y)=\alpha(y)+g(y)\quad\mbox{for $\,y\in U$}
\end{equation}
with $\alpha:=f'(x)$ and $g(y):=f(x)-f'(x)(x)+R(y)$ satisfying $\Lip(g|_{B^E_r(x)})=
\Lip(R|_{B^E_r(x)})\to 0$ as $r\to 0$.
If $f$ is strictly differentiable at each $x\in U$, then $f$ is called \emph{strictly differentiable}.\\[2.3mm]
Now let $(E,\|\cdot\|)$ be a Banach space, $x\in E$, $r>0$
and $f\colon B^E_r(x)\to E$ be a function which is strictly
differentiable at~$x$, with $f'(x)\in \GL(E)$.
Show that, after shrinking~$r$ if necessary,
$f|_{B^E_r(x)}$ is a homeomorphism onto its image,
$f(B^E_r(x))$ is open in~$E$, and the map $f^{-1}\colon f(B^E_r(x))\to E$
is strictly differentiable at $f(x)$ with derivative $f'(x)^{-1}$.\\[1mm]
[After shrinking~$r$, we are in the situation of Theorem~\ref{lip-inv-fct}.
Using the estimate (\ref{the-lip-h}) for~$\Lip(h)$
in Theorem~\ref{lip-inv-fct}, which tends to~$0$ as $r\to 0$,
we see that $f^{-1}$ is strictly differentiable at~$f(x)$.]
\end{exer}

\begin{exer}\label{exer-stri2}
Let $(E,\|\cdot\|)$ and $(F,\|\cdot\|_F)$ be normed spaces
and $U\sub E$ be an open subset.
\begin{description}[(D)]
\item[(a)]
Show that a function $f\colon U\to F$
is strictly differentiable if and only if $f$ is~$FC^1$.\\[1mm]
[If $f$ is strictly differentiable, then $f$ is totally differentiable at each $y\in U$;
differentiating (\ref{strict-condi}) with respect to $y$, we obtain
\[
f'(y)=f'(x)+R'(y)\quad\mbox{for all $y\in U$,}
\]
where $\|f'(y)-f'(x)\|_{\op}=\|R'(y)\|_{\op}\leq \Lip(R|_{B^E_r(x)})$ for
$y\in B^E_r(x)$. If, conversely, $f$ is $FC^1$, then $f$
is strictly differentiable by Lemma~\ref{strdiff}.]
\item[(b)]
Combining (a) with Exercise~\ref{exer-stri1},
obtain an alternative proof for the $FC^1$-case of
Theorem~\ref{inv-fct-class}.
\end{description}
\end{exer}

\begin{exer}
Let $E$ be a locally convex space, $U\sub E$ be a locally convex subset with
dense interior,
$F$ be a Banach space and $V\sub U\times F$ be relatively open.
Let $k\in\N\cup\{\infty\}$ with $k\geq 2$
and $f\colon V\to F$ be a $C^k$-map such that
$d_2f(x,y,\cdot)\in \GL(F)$ for all $(x,y)\in V$.
Show:
\begin{description}[(D)]
\item[(a)]
If $\gamma\colon U\to F$ is a continuous function with graph
in~$V$ such that $f(x,\gamma(x))=0$ for all $x\in U$, then
$\gamma$ is~$C^k$.
\item[(b)]
If $U$ is connected and $\gamma,\eta\colon U\to F$
are continuous functions with graph in~$V$
such that
\[
f(x,\gamma(x))=f(x,\eta(x))=0
\]
for all $x\in U$, and $\gamma(x_0)=\eta(x_0)$ for some
$x_0\in U$, then $\gamma=\eta$.
\end{description}
\end{exer}

\begin{exer}\label{exc-C1notFC1}
In this exercise, $L^1[0,1]$
denotes the Banach space of equivalence classes $[\gamma]$
(modulo functions vanishing almost everywhere)
of measurable functions $\gamma\colon[0,1]\to\R$
which are integrable with respect to Lebesgue--Borel measure.
Also, $L^\infty[0,1]$ is the space of equivalence classes of essentially bounded
measurable functions $\gamma\colon [0,1]\to\R$. As usual, we shall not distinguish functions and equivalence classes in
the notation.
It is well known that $L^1[0,1]'\cong L^\infty[0,1]$;
the isomorphism takes $\theta\in L^\infty[0,1]$ to the continuous linear functional
\[
\lambda_\theta\colon L^1[0,1]\to\R,\quad \gamma\mto \int_{[0,1]}\theta(x)\gamma(x)\, dx.
\]
\begin{description}[(D)]
\item[(a)]
Show that
$f(\gamma):=\sin \circ \, \gamma \in L^1[0,1]$
for all $\gamma\in L^1[0,1]$.
\item[(b)]
Show that $f\colon L^1[0,1]\to L^1[0,1]$
is Lipschitz continuous with $\Lip(f)\leq 1$.
Does the same reasoning apply to $g(\gamma):=\cos\circ \,\gamma$?
\item[(c)]
For $\gamma\in L^\infty[0,1]$, show that the linear map
\[
m_\gamma\colon L^1[0,1]\to L^1[0,1],\quad
\eta\mapsto \gamma\eta
\]
with $(\gamma\eta)(x):=\gamma(x)\eta(x)$
is continuous, and
$\|m_\gamma\|_{\op}\leq\|\gamma\|_{L^\infty}$.
\item[(d)]
Show that the mapping
\[
h\colon L^1[0,1]\times L^1[0,1]\to L^1[0,1],\quad h(\gamma,\eta):=
m_{\cos\circ \gamma}(\eta)=(\cos\circ \,\gamma)\eta
\]
is continuous. [Hint: In view of (c) and Exercise~\ref{exer-wopbd}(c), it suffices that
\[
L^1[0,1]\to  L^1[0,1],\;\; \gamma\mto h(\gamma,\eta)=m_\eta(g(\gamma))
\]
is continuous for each $\eta$ in the dense vector subspace $L^\infty[0,1]$ of $L^1[0,1]$.\,]
\item[(e)]
Let $\gamma,\eta\in L^1[0,1]$ und $t\in \R\setminus\{0\}$.
Verify that
\[
\frac{f(\gamma+t\eta)-f(\gamma)}{t}=\int_0^1 h(\gamma+st\eta,\eta)\,ds;
\]
to this end, apply the linear functionals $\lambda_\theta$
to both sides, for $\theta\in L^\infty[0,1]$.
\item[(f)]
Letting $t\to 0$ in (e), show that $df(\gamma,\eta)=h(\gamma,\eta)$.
Deduce that $f$ is $C^1$.
\item[(g)]
We have $f(0)=0$. Show that $f(L^1[0,1])$
is not a $0$-neighborhood in $L^1[0,1]$ although $f'(0)=\mbox{id}_{L^1[0,1]}$.
Deduce that $f$ is not $FC^1$.
\end{description}
\end{exer}

\begin{exer}\label{exc-C0k-impl}
Assume that $U\sub E$ is replaced with $U=P\times Q$
in the situation of Theorem~\ref{thm-inv-para}, where~$P$ is a topological space and~$Q$ a locally convex
subset with dense interior in a locally convex space~$Y$. Also,
assume that $f\colon P\times (Q \times V)\to F$ is not~$C^k$, but merely $C^{0,k}$
with $k\in\N\cup\{\infty\}$ (requiring~(\ref{need-ths2}) from Remark~\ref{variants-i-para}(a)
if $F$ has infinite dimension).
\begin{description}[(D)]
\item[(a)]
Using Lemma~\ref{fp-C0k-dep}, show that the conclusions of Theorem~\ref{thm-inv-para}
remain valid, except that
$\phi\colon P_0\times (Q_0\times V_0)\to W$ is only $C^{0,k}$
with~$C^{0,k}$ inverse (if we choose $U_0=P_0\times Q_0$
with an open $p_0$-neighborhood~$P_0$ in~$P$
and a relatively open $q_0$-neighborhood~$Q_0$ in~$Q$, where $x_0=(p_0,q_0)$).
\item[(b)]
Show that also the conclusion of the Implicit Function Theorem (Theorem~\ref{impl-Ck})
remains valid, except that $\phi\colon U_0=P_0\times Q_0\to V_0$
is only~$C^{0,k}$.
\end{description}
\end{exer}

\begin{exer}
Let $Z$ be a topological space, $Y$ be a set and
$q\colon Z\to Y$ be a mapping which is \emph{locally injective}
in the sense that each $z\in Z$ has a neighborhood $V\sub Z$ such that $q|_V$
is injective.
Let $X$ be a topological space and
$\gamma\colon X\to Z$ as well as $\eta\colon X\to Z$ be continuous mappings
such that
$q\circ \gamma=q\circ\eta$. Show:
If $X$ is connected, $Z$ is Hausdorff and $\gamma(x_0)=\eta(x_0)$ for some $x_0\in X$,
then $\gamma=\eta$.
\end{exer}

\begin{exer}[Simplified Newton Iteration]\label{exc-simpl-nwt}
Let $(E,\|\cdot\|)$ be a Banach space with $E\not=\{0\}$. Let $x_0\in E$, $r>0$, and
$f\colon B^E_r(x_0)\to E$ be an $FC^1$-map with $f'(x_0)\in\GL(E)$,
\begin{equation}\label{enab-sNewton}
a:=\frac{1}{\|f'(x_0)^{-1}\|_{\op}}-\Lip(f-f'(x_0))>0,\;\,\mbox{and}\;\,
\|f(x_0)\|< ar.
\end{equation}
\begin{description}[(D)]
\item[(a)]
Show that $f$ has a unique zero $x_\infty$ in $B^E_r(x_0)$.
\item[(b)]
Show that $x_n:=x_{n-1}-f'(x_0)^{-1}(f(x_{n-1}))\in B^E_r(x_0)$ holds for all
$n\in\N$ and $x_n\to x_\infty$ as $n\to\infty$.
\item[(c)]
Let $L:=\|f'(x_0)^{-1}\|_{\op}\Lip(f-f'(x_0))$ (which is an element of $[0,1[$)
and $C:=\|x_1-x_0\|/(1-L)$.
Show that $\|x_n-x_\infty\|\leq CL^n$ for all $n\in\N_0$.
\end{description}
[See Theorem~\ref{lip-inv-fct} and its proof.]
\end{exer}

\noindent
Every zero $x_\infty$ with $f'(x_\infty)\in\GL(E)$ can be
calculated by a simplified Newton iteration,
if we start with an initial approximation~$x_0$ sufficiently close to~$x_\infty$:

\begin{exer}
Let $(E,\|\cdot\|)$ be a Banach space such that $E\not=\{0\}$. Let $U\sub E$ be an open subset,
$f\colon U\to E$ be an $FC^1$-function and $x_\infty\in U$ be a
zero of~$f$ such that $f'(x_\infty)\in\GL(E)$.
Show that there exist $r\geq s >0$ with
$B^E_{2r}(x_\infty)\sub U$ such that $f|_{B^E_r(x_0)}$ satisfies the hypotheses (\ref{enab-sNewton})
of Exercise~\ref{exc-simpl-nwt} for each $x_0\in B^E_s(x_\infty)$.\\[2mm]
[Let $\ve>0$ such that $4\ve\leq 1/\|f'(x_\infty)^{-1}\|_{\op}$.
Choosing $r>0$ with $B^E_{2r}(x_\infty)\sub U$ small enough, we
can achieve that $f'(x_0)\in\GL(E)$ and
\[
\left|\frac{1}{\|f'(x_0)^{-1}\|_{\op}}-\frac{1}{\|f'(x_\infty)^{-1}\|_{\op}}\right|\leq\ve\;\,
\mbox{for all $x_0\in B^E_r(x_\infty)$.}
\]
We can also achieve that $\|f'(x)-f'(x_\infty)\|_{\op}\leq\ve$
for all $x\in B^E_{2r}(x_\infty)$, whence $\|f'(x)-f'(x_0)\|_{\op}\leq 2\ve$
for all $x_0\in B^E_r(x_\infty)$ and $x\in B^E_{2r}(x_\infty)$.
As a consequence,
$\Lip\big(f|_{B^E_r(x_0)}-f'(x_0)\big)\leq 2\ve$ for all $x_0\in B^E_r(x_\infty)$ and thus
\[
\frac{1}{\|f'(x_0)^{-1}\|_{\op}}-\Lip\big(f|_{B^E_r(x_0)}-f'(x_0)\big)
\geq
\frac{1}{\|f'(x_\infty)^{-1}\|_{\op}}-3\ve\geq \ve>0.
\]
It remains to choose $s\in\;]0,r]$ so small that $\|f(x_0)\|<\ve r$ for all
$x_0\in B^E_s(x_\infty)$.]
\end{exer}

\noindent
Also the Newton iterates tend to a zero $x_\infty$ with $f'(x_\infty)\in\GL(E)$, for~$x_0$ near~$x_\infty$.
\begin{exer}[Newton Method]\label{exc-loc-quadr}
Let $(E,\|\cdot\|)$ be a Banach space, $f\colon U\to E$ be an $FC^1$-function on an open subset
$U\sub E$ and $x_\infty\in U$ such that $f(x_\infty)=0$ and $f'(x_\infty)\in\GL(E)$.
Assume $E\not=\{0\}$.
\begin{description}[(D)]
\item[(a)]
Let $r>0$ with $B^E_r(x_\infty)\sub U$. After shrinking~$r$, we may assume that
$f'(x)\in\GL(E)$ for all $x\in B^E_r(x_\infty)$, and
$M:=\sup\{\|f'(x)^{-1}\|_{\op}\colon x\in B^E_r(x_\infty)\}<\infty$.
Shrinking~$r$ further, we may also assume that
\[
L:=M\sup\{\|f'(x)-f'(y)\|_{\op}\colon x,y\in B^E_r(x_\infty)\}<1.
\]
For each $x\in B^E_r(x_\infty)$, the map
\[
g_x\colon B^E_r(x_\infty)\to E,\quad y\mto y-f'(x)^{-1}(f(y))
\]
satisfies $g_x(x_\infty)=x_\infty$ and
\begin{equation}\label{coarse-contra}
\Lip(g_x)\leq \|f'(x)^{-1}\|_{\op}\Lip\big(f|_{B^E_r(x_\infty)}-f'(x)\big)\leq L,
\end{equation}
whence it is a self-map of $B^E_r(x_\infty)$ and a contraction.
Given $x_0\in B^E_r(x_\infty)$, this enables us to define 
\[
x_n:=g_{x_{n-1}}(x_{n-1})=x_{n-1}-f'(x_{n-1})^{-1}(f(x_{n-1}))\in B^E_r(x_\infty)
\]
recursively for all $n\in\N$.
Deduce from~(\ref{coarse-contra}) that $\|x_n-x_\infty\|\leq L^n\|x_0-x_\infty\|$
for all $n\in\N$. Hence $x_n\to x_\infty$ for $n\to\infty$.
\item[(b)]
To get better estimates for the asymptotic behavior as $n\to\infty$, define
\[
\theta_s:=M\sup\{\|f'(x)-f'(y)\|_{\op}\colon x,y\in \wb{B}^E_s(x_\infty)\}\;\,\mbox{for $s\in [0,r[$.}
\]
Then $\theta_s\leq L<1$ for all $s\in [0,r[$, and $\theta_s\to 0$ as $s\to 0$.
Check that
\[
\Lip\big(g_x|_{\wb{B}^E_{\|x-x_\infty\|}(x_\infty)}\big)\leq\theta_{\|x-x_\infty\|}
\]
for all $x\in B^E_r(x_\infty)$. Deduce that, in the situation of~(a),
we have the estimate $\|x_n-x_\infty\|\leq \theta_{\|x_{n-1}-x_\infty\|}\cdots
\theta_{\|x_0-x_\infty\|}\|x_0-x_\infty\|$, for all $n\in\N$.
\item[(c)] Deduce from~(b): For each $\theta>0$, given $x_0\in B^E_r(x_\infty)$ there exists $C>0$ such that
$\|x_n-x_\infty\|\leq C\hspace*{.2mm}\theta^n$ for all $n\in\N$.
\item[(d)]
Now assume that~$f$ is~$C^2$.
By Remark~\ref{back-to-normed}(c),
after shrinking~$r$ if necessary, we can achieve that
there exists $K>0$ such that $\|f'(y)-f'(x)\|_{\op}\leq \frac{1}{2}K\|y-x\|$ for all $x,y\in B^E_r(x_\infty)$
and thus
\[
\theta_s\leq MKs\;\,\mbox{for all $s\in [0,r[$.}
\]
Hence $\|x_n-x_\infty\|\leq \theta_{\|x_{n-1}-x_\infty\|}\|x_{n-1}-x_\infty\|\leq
MK \|x_{n-1}-x_\infty\|^2$ and thus
\[
MK\|x_n-x_\infty\|\leq (MK\|x_{n-1}-x_\infty\|)^2
\]
for all $n\in\N$. Abbreviating $\delta_n:= MK\|x_n-x_\infty\|$
for $n\in \N_0$, this becomes
\[
\delta_n\leq(\delta_{n-1})^2\;\,\mbox{for all $n\in\N$}
\]
(an estimate referred to as ``locally quadratic convergence'' in the literature).
Given $\theta\in\,]0,1[$, there is $n_\theta\in\N$ such that
$\delta_{n_\theta}\leq\theta$. Deduce that
\[
\delta_{n_\theta+k}\leq \theta^{(2^k)}\;\,\mbox{for all $k\in \N_0$.}
\]
\end{description}
\end{exer}
\noindent
We mention a criterion for convergence of the Newton iteration.
\begin{exer}\label{exc-class-newt}
Let $(E,\|\cdot\|)$ be a Banach space with $E\not=\{0\}$
and $f\colon U\to E$ be an $FC^1$-map on an open subset $U\sub E$.
Verify the details of the following sketch:
\begin{description}[(D)]
\item[(a)]
Let $x_0\in U$ and $r>0$ such that $B^E_{2r}(x_0)\sub U$,
$f'(x)\in \GL(E)$ for all $x\in B^E_{2r}(x_0)$,
$m:=\sup\{\|f'(x)^{-1}\|_{\op}\colon x\in B^E_{2r}(x_0)\}<\infty$ and
\[
\ell\, :=\, m\,\sup\{\|f'(y)-f'(x)\|_{\op}\colon x,y\in B^E_{2r}(x_0)\}<1.
\]
Then $a:=\frac{1}{\|f'(x_0)^{-1}\|_{\op}}-\Lip\big(f|_{B^E_{2r}(x_0)}-f'(x_0)\big)\geq
\frac{1}{m}(1-\ell)\in \;]0,1]$. If $\|f(x_0)\|<ar$,
then $f|_{B^E_r(x_0)}$ has a zero~$x_\infty$, by Theorem~\ref{lip-inv-fct}.
Also, $f|_{B^E_r(x_\infty)}$ satisfies the hypotheses of
Exercise~\ref{exc-loc-quadr}(a),
and $x_0\in B^E_r(x_\infty)$. Hence
$x_n:=x_{n-1}-f'(x_{n-1})^{-1}(f(x_{n-1}))$ is defined for all~$n$ and
$x_n\to x_\infty$ as $n\to\infty$.
\item[(b)]
If $x_\infty$ is any zero of~$f$ such that $f'(x_\infty)\in\GL(E)$, then
there exists $r>0$ such that $B^E_{3r}(x_\infty)\sub U$, $f'(x)\in\GL(E)$ for all
$x\in B^E_{3r}(x_\infty)$,
\[
M:=\sup\{\|f'(x)^{-1}\|_{\op}\colon x\in B^E_{3r}(x_\infty)\}<\infty,
\]
and
$L:=M\sup\{\|f'(y)-f'(x)\|_{\op}\colon x,y\in B^E_{3r}(x_\infty)\}<1$.
Choose $s\in \;]0,r]$ so small that $\|f(x_0)\|\leq \frac{1}{M}(1-L)r$
for all $x_0\in B^E_s(x_\infty)$.
Then $f|_{B^E_{2r}(x_0)}$\vspace{-.5mm} satisfies the hypotheses of~(a),
for each $x_0\in B^E_s(x_\infty)$.
\end{description}
\end{exer}

\begin{exer}
Define $g\colon \R\to\R$ via $g(x)=0$ if $x\leq 1/2$, $g(x)=x/2-1/4$ if $x\in [1/2,3/2]$,
$g(x)=1/2$
if $x\geq 3/2$. Let $x_0:=2$ and consider $f\colon \R\to\R$, $x\mto x+\int_0^x g(t)\,dt$.
\begin{description}[(D)]
\item[(a)]
Let $x_n:=x_{n-1}-f'(x_0)^{-1}f(x_{n-1})$ for $n\in\N$. Show that $x_n=(1/3)^n$
for all $n\in\N$. Deduce that the conclusion of Exercise~\ref{exc-loc-quadr}(c) does not hold for the \emph{simplified}
Newton iteration in general.
\item[(b)]
Now consider the Newton iterates $y_n:=y_{n-1}-f'(y_{n-1})^{-1}f(y_{n-1})$
for $n\in\N$, with $y_0:=x_0$. Show that $y_n=0$, the unique zero of $f$, for all $n\geq 2$.
\end{description}
\end{exer}
To speed up the convergence to the fixed point,
one can frequently use a Newton iteration instead of the iteration
described in Banach's Fixed Point Theorem.
\begin{exer}
Let $(E,\|\cdot\|)$ be a Banach space, $U\sub E$ be an open subset and $f\colon U\to E$ be an $FC^1$-map
such that $\sup\{\|f'(x)\|_{\op}\colon x\in U\}<1$.
Assume that~$f$ has a fixed point~$x_\infty$.
Then $x_\infty$ is a zero of the $FC^1$-map $h\colon U\to E$, $h(x):=x-f(x)$.
Show that $h'(x)\in\GL(E)$ for all $x\in U$; notably, $h'(x_\infty)\in\GL(E)$.
Deduce that~$x_\infty$ can be obtained as the limit of a Newton iteration for the function~$h$,
if we start with an initial value~$x_0$ sufficiently close to $x_\infty$.  
\end{exer}

\end{small}
\section{Differential equations in locally convex spaces}\label{sec-ode}
We now study differential equations on subsets of locally convex spaces,
including the classical case of differential equations in Banach spaces.

All vector spaces considered in this section are vector spaces
over the field of real numbers, unless the contrary is stated.
\begin{defn}
Let $E$ be a locally convex space, $W\sub \R\times E$ be a subset and
$f\colon W\to E$ be a function.\medskip

\noindent
(a) We say that a function $\gamma \colon I\to E$
on a non-degenerate interval $I\sub \R$
is a \emph{solution to the differential equation}
\index{solution to a differential equation}
\begin{equation}\label{theode}
y'(t)=f(t,y(t))
\end{equation}
if $\gamma$ is $C^1$,
\begin{equation}\label{theode2}
(t,\gamma(t))\in W\quad\mbox{and}\quad \gamma'(t)=f(t,\gamma(t))\;\;
\mbox{for all $\,t\in I$.}
\end{equation}
(b) If $(t_0,y_0)\in W$ and $\gamma\colon I\to E$ is a solution to~(\ref{theode})
such that $t_0\in I$ and $\gamma(t_0)=y_0$,
then $\gamma$ is called a 
\emph{solution to the initial value problem}
\index{solution to an initial value problem}
\begin{equation}\label{theiniprob}
\left\{
\begin{array}{rcl}
y'(t) &=& f(t,y(t))\\
y(t_0)&=&y_0.
\end{array}
\right.
\end{equation}
(c) We say that the differential equation (\ref{theode})
satisfies 
\emph{local uniqueness of solutions} 
\index{local uniqueness of solutions} 
if the following condition holds:
For all solutions $\gamma_1\colon I_1\to E$ and $\gamma_2\colon I_2\to E$ to (\ref{theode})
such that $\gamma_1(t_0)=\gamma_2(t_0)$ at some $t_0\in I_1\cap I_2$,
there exists a neighborhood~$I$ of~$t_0$ in $I_1\cap I_2$ such that
$\gamma_1|_I=\gamma_2|_I$.\medskip

\noindent
(d) If $W=J\times U$ for a non-degenerate interval $J\sub\R$ and some subset $U\sub E$
(or if $W\sub J\times U$ with $J$ and $U$ as before),\footnote{Later, $W$ will be relatively open in $J\times U$ and hence locally of this form.}
we say that the differential equation (\ref{theode}) satisfies 
\emph{local existence of solutions} 
if the following condition holds:\footnote{More precisely, we should speak of local existence
of solutions \emph{with respect to~$J$}; however, the choice of~$J$
will always be clear from the context.}
For all $(t_0,y_0)\in W$, there exists a solution $\gamma\colon I\to E$ to the initial value
problem (\ref{theiniprob}) whose domain~$I$ is relatively open in~$J$.
\end{defn}
If~$E$ is a Banach space, then there are well-known easy criteria for local existence and uniqueness
of solutions (which we shall recall soon).
For locally convex spaces~$E$ which are not Banach spaces
and suitable continuous linear mappings $A\colon E\to E$,
not even linear differential equations of the form
\begin{equation}\label{lin-ode}
y'(t)=Ay(t)
\end{equation}
need to satisfy local uniqueness;
it is also possible to choose $E$ and~$A$ in such a way that
local existence of solutions fails.\footnote{The right-hand-side
$f\colon \R\times E\to E$, $f(t,y):=A(y)$ is used here.}
We recall well-known examples:
\begin{ex}\label{ode-no}
Consider the locally convex direct sum
$E:=\R^{(\N)}=\bigoplus_{n\in\N}\R$ (see Examples~\ref{theKinfty}
and \ref{Rinftysum})
and the right shift
\[
A \colon E\to E,\quad (x_1,x_2,\ldots)\mto(0,x_1,x_2,\ldots )
\]
which is continuous (as a consequence of Lemma~\ref{ctsonsteps}(ii))
and linear. If the initial value problem
\[
y'(t)=Ay(t), \quad y(0)=(1,0,0,\ldots)
\]
had a solution $\gamma=(\gamma_n)_{n\in\N}\colon I\to E$,
then $\gamma_1$ would satisfy $\gamma_1(0)=1$ and
\[
\gamma_1'(t)=0\;\;\mbox{for all $t\in I$,}
\]
whence $\gamma_1(t)=1$ for all $t\in I$.
For $n>1$, we would have $\gamma_n(0)=0$ and
\[
\gamma_n'(t)=\gamma_{n-1}(t)\quad\mbox{for all $t\in I$,}
\]
whence $\gamma_n(t)=t^{n-1}/(n-1)!$ for all $t\in I$ by induction on~$n$.
Since
\[
(t^{n-1}/(n-1)!)_{n\in\N}\not\in \R^{(\N)}\quad\mbox{for all $t\in\R\setminus\{0\}$,}
\]
we deduce that $I=\{0\}$, contradicting the hypothesis that~$I$ is non-degenerate.
Thus~$\gamma$ cannot exist.
\end{ex}
\begin{ex}\label{ode-many}
Let $E:=\R^{\N_0}=\prod_{n\in\N_0}\R$
with the product topology, which makes it a non-normable Fr\'{e}chet space
(cf.\ Exercise~\ref{exc-allseq}). We consider the left shift
\[
A\colon E\to E,\quad (x_n)_{n\in\N_0}\mto (x_{n+1})_{n\in\N_0}.
\]
Then the $0$-function $\gamma\colon \R\to E$, $t\mto 0$
solves the initial value problem
\begin{equation}\label{ini-here}
y'(t)=Ay(t), \quad y(0)=0.
\end{equation}
Let $h\colon\R\to\R$ be a smooth function
such that $h^{(n)}(0)=0$ for all $n\in\N_0$
but $h(t)\not=0$ for all $t\not=0$
(for example, $h(t):=e^{-1/t^2}$ for $t\not=0$).
Then also the function $\eta=(\eta_n)_{n\in\N_0}\colon\R\to\R^{\N_0}$
with $n$th component $\eta_n:=h^{(n)}$
is a solution to~(\ref{ini-here}).
Looking at $\eta_0$, we see that $\gamma$ and $\eta$ do not coincide
on any $0$-neighborhood.
\end{ex}
For later use, we mention:
\begin{lem}\label{via-integral-eq}
Let $E$ be a locally convex space, $W\sub \R\times E$ be a subset and
$f\colon W\to E$ be a continuous function.
Let $\gamma\colon I\to E$ be a function
on a non-degenerate interval $I\sub \R$ such that
$(t,\gamma(t))\in W$ for all $t\in I$.
Let $t_0\in I$ and $y_0:=\gamma(t_0)$.
Then the following conditions are equivalent:
\begin{description}[(D)]
\item[\rm(a)]
$\gamma$ is~$C^1$ and solves the inital value problem \emph{(\ref{theiniprob})};
\item[\rm(b)]
$\gamma$ is continuous, the weak integral $\int_{t_0}^tf(s,\gamma(s))\,ds$
exists in~$E$
for all $t\in I$, and $\gamma(t)=y_0+\int_{t_0}^tf(s,\gamma(s))\,ds$.
\end{description}
\end{lem}
\begin{prf}
If (a) holds, then also~(b), by Proposition~\ref{fundamental} and (\ref{theiniprob}).
If (b) holds, then also~(a), by Proposition~\ref{secpart}.
\end{prf}
\begin{rem}\label{glue-sol-ode}
Let~$E$ be a locally convex space, $W\sub \R\times E$ be a subset and
$f\colon W\to E$ be a continuous function.\medskip

(a) If $I\sub \R$ is a non-degenerate interval, $t_0\in I$
and $\gamma\colon I\to E$
a continuous map with $\{(t,\gamma(t))\colon t\in I\}\sub W$
whose restrictions to the connected components of $I\setminus\{t_0\}$
solve~(\ref{theode}), then also~$\gamma$ satisfies~(\ref{theode})
as the function $I\to E$, $t\mto f(t,\gamma(t))$ is continuous
and Lemma~\ref{C1glueing} applies.\medskip

(b) If $I_1$ and~$I_2$ are non-degenerate intervals in~$\R$ such that
$t_0:=\max(I_1)=\min(I_2)$ and $\gamma_j\colon I_j\to E$ are solutions to
(\ref{theode}) for $j\in\{1,2\}$ such that $\gamma_1(t_0)=\gamma_2(t_0)$,
then also
\[
\gamma\colon I_1\cup I_2\to E,\quad t\mto\left\{
\begin{array}{cc}
\gamma_1(t) &\mbox{if $t\in I_1$;}\\
\gamma_2(t)&\mbox{if $t\in I_2$}
\end{array}
\right.
\]
is a solution to~(\ref{theode}), by~(a).\medskip

(c) If $W\sub\R\times E$ is a locally convex subset with dense interior, $k\in\N_0$
and $f\colon W\to E$ is~$C^k$, then every solution
$\gamma\colon I\to E$ to~(\ref{theode}) is~$C^{k+1}$. For $k=0$,
this holds by definition. If $k\geq 1$ and the assertion holds for~$k-1$,
then $\gamma$ is $C^k$. Since
\[
\gamma'(t)=f(t,\gamma(t)),
\]
$\gamma'$ is $C^k$ by the Chain Rule and thus~$\gamma$ is $C^{k+1}$.
\end{rem}
\begin{lem}\label{lem-loc-glob}
Consider a differential equation \emph{(\ref{theode})}
that satisfies local uniqueness of solutions.
If $\gamma_1\colon I_1\to E$ and $\gamma_2\colon I_2\to E$
are solutions to \emph{(\ref{theode})}
such that $\gamma_1(t_0)=\gamma_2(t_0)$
for some $t_0\in I_1\cap I_2$, then $\gamma_1|_{I_1\cap I_2}=\gamma_2|_{I_1\cap I_2}$.
\end{lem}
\begin{prf}
The subset $K:=\{t\in I_1\cap I_2\colon\gamma_1(t)=\gamma_2(t)\}$
is closed in $I_1\cap I_2$ since $E$ is Hausdorff and the functions $\gamma_1$
and~$\gamma_2$ are continuous. Since (\ref{theode}) satisfies local uniqueness
of solutions, $K$ is relatively open in~$I_1\cap I_2$. By hypothesis, $K\not=\emptyset$.
As $I_1\cap I_2$ is an interval and hence connected, $K=I_1\cap I_2$ follows.
\end{prf}
\begin{lem}\label{ex-maxi-sol}
Consider a differential equation
\begin{equation}\label{goodode}
y'(t)=f(t,y(t))
\end{equation}
satisfying both local uniqueness and local existence
of solutions, where $E$ is a locally convex space, $J\sub \R$ a non-degenerate interval,
$U\sub E$
and $f\colon W\to E$ a function on a subset $W\sub J\times U$.
For each $(t_0,y_0)\in W$, there exists
a solution $\gamma\colon I\to E$ to the initial value problem \emph{(\ref{theiniprob})}
with the following property:
For each solution $\eta\colon I_\eta\to E$ of \emph{(\ref{theiniprob})},
\begin{equation}\label{proper-max}
I_\eta\sub I\quad\mbox{and}\quad \eta=\gamma|_{I_\eta}
\end{equation}
holds.
The function $\gamma$ is uniquely determined by the latter property,
and~$I$ is relatively open in~$J$.
\end{lem}
\begin{prf}
For fixed $(t_0,y_0)\in W$,
let~$S$ be the set of all solutions $\eta\colon I_\eta\to E$ of (\ref{theiniprob}).
Let $I:=\bigcup_{\eta\in S}I_\eta$.
Then
\[
\gamma\colon I\to E,\quad \gamma(t):=\eta(t)\;\;\mbox{if $\eta\in S$ and $t\in I_\eta$}
\]
is a well-defined function, as
$\eta|_{I_\eta\cap I_\zeta}=\zeta|_{I_\eta\cap I_\zeta}$ for all $\eta,\zeta\in S$
by Lemma~\ref{lem-loc-glob}. 
By construction, $I_\eta\sub I$ and $\gamma|_{I_\eta}=\eta$ for all
$\eta\in S$, whence $\gamma$ is a solution to (\ref{theiniprob})
and~(\ref{proper-max}) holds. The uniqueness assertion made in the lemma is clear.
If $I$ was not open in~$J$, then we would have $t_1:=\sup(I)\in I$ but $t_1<b$ for some $b\in J$
(which we assume now), or $t_1:=\inf(I)\in I$ but $a<t_1$ for some $a\in J$
(a case which can be discussed by analogous arguments).
By local existence, there is a solution $\theta\colon K\to E$ to the differential equation (\ref{theode})
with $\theta(t_1)=\gamma(t_1)$, whose domain~$K$ is relatively open in~$J$.
Then $[t_1,t_2]\sub K$ for some $t_2>t_1$. Now $\eta\colon I\cup [t_1,t_2]\to E$,
$\eta(t):=\gamma(t)$ if $t\in I$, $\eta(t):=\theta(t)$ if $t\in[t_1,t_2]$
is a solution to the initial value problem~(\ref{theiniprob}).
Hence $\eta\in S$ and $I\cup[t_1,t_2]=I_\eta\sub I$,
contrary to $t_2>\sup(I)$.
\end{prf}
\begin{defn}\label{flow-nonauto}
The solution $\gamma\colon I\to E$ constructed in Lemma~\ref{ex-maxi-sol}
for $(t_0,y_0)\in W$ is called the \emph{maximal solution}
to the initial value problem (\ref{theiniprob});
to emphasize its dependence on $t_0$ and $y_0$,
we also write $\gamma_{t_0,y_0}\colon I_{t_0,y_0}\to E$ 
in place of~$\gamma$.
The subset $\Omega\sub J\times J\times U$ given by
\[
\Omega:=\bigcup_{(t_0,y_0)\in W}I_{t_0,y_0}\times \{(t_0,y_0)\}
\]
is the domain of the so-called (maximal) \emph{flow}
\[
\Fl\colon \Omega\to E,\quad \Fl(t,t_0,y_0):=\gamma_{t_0,y_0}(t).
\]
Given $t,t_0\in J$, we define $\Omega_{t,t_0}:=\{y_0\in E\colon
(t_0,y_0)\in W\;\mbox{and}\;
t\in I_{t_0,y_0}\}$
and consider the partial map
\[
\Fl_{t,t_0} \colon\Omega_{t,t_0}\to E,\quad y_0\mto\Fl(t,t_0,y_0).
\]
\end{defn}
In the situation of Lemma~\ref{ex-maxi-sol}
and Definition~\ref{flow-nonauto}, we have:
\begin{lem}\label{basics-flow}
\begin{description}[(D)]
\item[\rm(a)]
If $(t_0,y_0)\in W$ and $t_1\in I_{t_0,y_0}$,
then $\gamma_{t_0,y_0}=\gamma_{t_1,y_1}$ with $y_1:=\gamma_{t_0,y_0}(t_1)
=\Fl_{t_1,t_0}(y_0)$.\vspace{.8mm}
\item[\rm(b)]
If $t_2\in I_{t_1,y_1}$ in \emph{(a)},
then $t_2\in I_{t_0,y_0}$ and
$\Fl_{t_2,t_0}(y_0)=\Fl_{t_2,t_1}(\Fl_{t_1,t_0}(y_0))$.\vspace{.8mm}
\item[\rm(c)]
For all $t_1,t_0\in J$, the map $\Fl_{t_1,t_0}$ is injective,
$\Fl_{t_1,t_0}(\Omega_{t_1,t_0})=\Omega_{t_0,t_1}$, and $\Fl_{t_0,t_1}=(\Fl_{t_1,t_0})^{-1}$.
\end{description}
\end{lem}
\begin{prf}
(a) Abbreviate $y_1:=\Fl_{t_1,t_0}(y_0)$.
The functions $\gamma_{t_0,y_0}$ and $\gamma_{t_1,y_1}$
are solutions to the differential equation~(\ref{theode})
whose values at~$t_1$ coincide.
Thus $I_{t_0,y_0}\sub I_{t_1,y_1}$
by maximality of $\gamma_{t_1,y_1}$, 
and $\gamma_{t_0,y_0}$ is the restriction of $\gamma_{t_1,y_1}$
to~$I_{t_0,y_0}$. Hence $t_0\in I_{t_1,y_1}$ and $\gamma_{t_1,y_1}(t_0)=\gamma_{t_0,y_0}(t_0)=y_0$.
Thus $I_{t_1,y_1}\sub I_{t_0,y_0}$ and $\gamma_{t_1,y_1}=\gamma_{t_0,y_0}|_{I_{t_1,y_1}}$,
by maximailty of $\gamma_{t_0,y_0}$.
Hence $\gamma_{t_0,y_0}=\gamma_{t_1,y_1}$.\smallskip

(b) By~(a), we have $t_2\in I_{t_0,y_0}$ and $\Fl_{t_2,t_0}(y_0)=\gamma_{t_0,y_0}(t_2)
=\gamma_{t_1,y_1}(t_2)=
\Fl_{t_2,t_1}(\Fl_{t_1,t_0}(y_0))$.\smallskip

(c) If $y_0\in\Omega_{t_1,t_0}$, abbreviate $y_1:=\gamma_{t_0,y_0}(t_1)=\Fl_{t_1,t_0}(y_0)$.
Then $\gamma_{t_0,y_0}=\gamma_{t_1,y_1}$, by~(a).
Hence $t_0\in I _{t_1,y_1}$,
whence $y_1\in\Omega_{t_0,t_1}$ and thus
\begin{equation}\label{one-incls}
\Fl_{t_1,t_0}(\Omega_{t_1,t_0})\sub \Omega_{t_0,t_1}.
\end{equation}
Moreover,
$\Fl_{t_0,t_1}(\Fl_{t_1,t_0}(y_0))=\Fl_{t_0,t_0}(y_0)=y_0$, by~(b).
As a result, $\Fl_{t_1,t_0}$ is injective
and $\Fl_{t_0,t_1}\circ\Fl_{t_1,t_0}=\id_{\Omega_{t_1,t_0}}$,
whence $\Fl_{t_0,t_1}(\Omega_{t_0,t_1})\supseteq \Omega_{t_1,t_0}$.
Reversing the roles of $t_0$ and $t_1$, we also have
$\Fl_{t_1,t_0}(\Omega_{t_1,t_0})\supseteq \Omega_{t_0,t_1}$.
Hence $\Fl_{t_1,t_0}(\Omega_{t_1,t_0})=\Omega_{t_0,t_1}$, using~(\ref{one-incls}).
\end{prf}
We now discuss three conditions ensuring local uniqueness
of solutions.
\begin{defn}
Let $X=\R$ (or any topological space),
$(E,\|\cdot\|)$ be a normed space
and $f\colon W\to E$ be a function on a subset
$W\sub X\times E$.
Let $L\in[0,\infty[$. We say that~$f$ satisfies a \emph{global Lipschitz condition
with constant~$L$}~if
\[
\|f(t,y_2)-f(t,y_1)\|\leq L\|y_2-y_1\|
\]
for all $t\in X$ and $y_1,y_2\in E$ such that $(t,y_1),(t,y_2)\in W$.
If $f$ satisfies a global Lipschitz condition with constant~$L$ for some~$L$,
we say that $f$ satisfies a \emph{global Lipschitz condition}.
If each $(t_0,y_0)\in W$ has a neighborhood $Y\sub W$ such that
$f|_Y$ satisfies a global Lipschitz condition, then $f$ is said to satisfy a \emph{local
Lipschitz condition}. We also speak of a global or local Lipschitz condition
\emph{in the second argument}, or \emph{in the $y$-variable}, for added clarity.
\end{defn}
If also $P$ is a topological space and $W$ a subset of $X\times E\times P$,
we say that a function $f\colon W\to E$ satisfies a global (or local)
Lipschitz condition in the second argument if the map $((t,p),y)\mto f(t,y,p)$
on the corresponding subset of $(X\times P)\times E$ does so.
\begin{prop}\label{uniq-normd}
Let $(E,\|\cdot\|)$ be a normed space and $f\colon W\to E$
be a function on a subset $W\sub \R\times E$.
If $f$ satisfies a local Lipschitz condition,
then the differential equation
$y'(t)=f(t,y(t))$ satisfies local uniqueness
of solutions.
\end{prop}
\begin{prf}
Let $\gamma_1\colon I_1\to E$ and $\gamma_2\colon I_2\to E$
be solutions of (\ref{theode}) and $\gamma_1(t_0)=\gamma_2(t_0)$
for some $t_0\in I_1\cap I_2$. Abbreviate $y_0:=\gamma_1(t_0)=\gamma_2(t_0)$.
Since $f$ satisfies a local Lipschitz condition,
we find $\ve>0$, $\delta>0$ and $L\in[0,\infty[$
such that $f|_Y$ satisfies a global Lipschitz condition with constant~$L$
on
\[
Y:=\{(t,y)\in W\colon |t-t_0|\leq \delta\;\mbox{and}\;\|y-y_0\|\leq \ve\}.
\]
After shrinking~$\delta$, we may assume that $\delta L<1$
and $\|\gamma_j(t)-y_0\|\leq \ve$ for all $j\in \{1,2\}$ and $t$ in
\[ U:=\{t\in I_1\cap I_2\colon |t-t_0|\leq\delta \},\]
which is a neighborhood of~$t_0$ in $I_1\cap I_2$. For all $s\in U$, we have
\[
\|f(s,\gamma_2(s))-f(s,\gamma_1(s))\|\leq L\|\gamma_2(s)-\gamma_1(s)\|
\leq L\|\gamma_1|_U-\gamma_2|_U\|_\infty,
\]
where $\|\cdot\|_\infty$ denotes the supremum norm on $C(U,E)$. Using Lemma~\ref{via-integral-eq}(b),
this implies that
\begin{eqnarray*}
\|\gamma_2(t)-\gamma_1(t)\| & =& \left\|
\int_{t_0}^t\big(f(s,\gamma_2(s))-f(s,\gamma_1(s))\big)\,ds\right\| \\
&\leq & \left|\int_{t_0}^t\|f(s,\gamma_2(s))-f(s,\gamma_1(s))\|\,ds\right|
\leq \delta L\|\gamma_2|_U-\gamma_1|_U\|_\infty
\end{eqnarray*}
for all $t\in U$, entailing that
\[
\|\gamma_2|_U-\gamma_2|_U\|_\infty\leq \delta L \|\gamma_2|_U-\gamma_1|_U\|_\infty.
\]
Since $\delta L < 1$, this implies that $\|\gamma_2|_U-\gamma_1|_U\|_\infty=0$
and thus $\gamma_1|_U=\gamma_2|_U$.
\end{prf}
\begin{ex}\label{gives-lip-con}
If $X$ is a topological space, $(E,\|\cdot\|)$
a Banach space, $U\sub E$ a locally convex subset with dense interior
and $f\colon X\times U\to E$ a $C^{0,1}$-map, then
$f$ satisfies a local Lipschitz condition in the second argument.\\[2mm]
[Given $(x_0,y_0)\in X\times U$, there exists a neighborhood $X_0$ of $x_0$ in~$X$,
a convex neighborhood $U_0$ of~$y_0$ in~$U$
and $r>0$ such that
\[
d_2f(X_0\times U_0 \times \wb{B}^E_r(0))\sub \wb{B}^E_1(0).
\]
Then $\|d_2f(x,y,\cdot)\|_{\op}\leq\frac{1}{r}$ for all $(x,y)\in X_0\times U_0$
and thus $\Lip(f(x,\cdot)|_{U_0})\leq\frac{1}{r}$ for all $x\in X_0$, by Lemma~\ref{lipviaprime}.]
\end{ex}
\begin{defn}
Let $(E,\|\cdot\|)$ be a locally convex space over $\K\in\{\R,\C\}$.
A \emph{linear differential equation} in~$E$
is a differential equation of the form
\[
y'(t)=A(t)y(t)+b(t),
\]
where $J\sub\R$ is a non-degenerate interval,
$b\colon J\to E$ is a continuous function
and $A\colon J\to\cL(E)$ is a function such that
\[
\wh{A}\colon J\times E\to E,\quad (t,y)\mto A(t)(y)
\]
is continuous.
If $b=0$, then the linear differential equation is called \emph{homogeneous};
if $b\not=0$, it is called \emph{inhomogeneous}.
\end{defn}
\begin{rem}
Continuity
of $\wh{A}\colon J\times E\to E$ implies continuity of $A$ as a map
to $\cL(E)$, endowed with the compact-open topology;
if $E$ is metrizable (or, more generally, if $J\times E$ is a $k_\R$-space),
then continuity of $\wh{A}$ and continuity of $A\colon J\to\cL(E)_c$ are equivalent
(cf.\ Proposition~\ref{ctsexp}).
The same conclusion holds if~$E$ is barreled.
To see this, assume that
$A\colon J\to\cL(E)_c$ is continuous.
The bilinear evaluation map
$\beta\colon E\times\cL(E)_c\to E$, $(y,\alpha)\mto \alpha(y)$
is separately continuous and hence hypocontinuous if~$E$ is barreled
(see Corollary~\ref{bar-hypo}).
Now $A(K)\sub\cL(E)_c$ is compact (and hence bounded) for each compact subset $K\sub J$.
Hence
$\beta|_{E\times A(K)}$
is continuous (see Proposition~\ref{three-hypo}(a)),
whence
$K\times E\to E$, $(t,y)\mto \wh{A}(t,y)=\beta(y,A(t))$ is continuous
and hence also $\wh{A}$.
\end{rem}
\begin{ex}\label{lin-loclip}
For every linear differential equation
$y'(t)=A(t)y(t)+b(t)$ in a normed space $(E,\|\cdot\|)$,
the continuous function
\[
f\colon J\times E\to E,\quad f(t,y):=A(t)(y)+b(t)
\]
satisfies a local Lipschitz condition in the second argument;
if~$J$ is compact, then~$f$ satisfies a global Lipschitz condition.\\[2mm]
[It suffices to prove the final assertion.
As $\wh{A}\colon J\times E\to E$ is continuous
and $\wh{A}|_{J\times\{0\}}=0$, there exists $r>0$ such that $\wh{A}(J\times \wb{B}^E_r(0))\sub\wb{B}^E_1(0)$,
by the Wallace Lemma. Thus $\Lip f(t,\cdot)=\Lip A(t)=\|A(t)\|_{\op}\leq\frac{1}{r}$
for all $t\in J$.]
\end{ex}
\begin{prop}\label{uniq-pl}
Let $E$ be a locally convex space and $f\colon W\to E$
be a map on a subset $W\sub \R\times E$.
Let $\Lambda$ be a set of continuous linear maps
$\lambda\colon E\to E_\lambda$ to locally convex spaces~$E_\lambda$
satisfying the following conditions:
\begin{description}[(D)]
\item[\rm(a)]
$\Lambda$ separates points on~$E$;
\item[\rm(b)]
For all $\lambda\in\Lambda$, there exists a function $f_\lambda\colon W_\lambda\to E_\lambda$
on a subset \break $W_\lambda\sub \R \times E_\lambda$ such that
$(\id_\R\times \lambda)(W)\sub W_\lambda$,
\[
\lambda\circ f=f_\lambda\circ (\id_\R\times \lambda)|_W
\]
and such that the differential equation $y'(t)=f_\lambda(t,y(t))$
satisfies local uniqueness of solutions.
\end{description}
Then~$f$ satisfies local uniqueness of solutions.
\end{prop}
\begin{prf}
Let $\gamma_1\colon I_1\to E$ and $\gamma_2\colon I_2\to E$
be solutions to (\ref{theode}) and $t_0\in I_1\cap I_2$
such that $\gamma_1(t_0)=\gamma_2(t_0)$.
For $\lambda\in\Lambda$ and $j\in\{1,2\}$, we then have
\[
(\lambda\circ \gamma_j)'(t)=\lambda(\gamma_j'(t))=\lambda(f(t,\gamma_j(t)))
=f_\lambda(t,(\lambda\circ \gamma_j)(t))
\]
for all $t\in I_j$, whence $\lambda\circ \gamma_j$
is a solution to $y'(t)=f_\lambda(t,y(t))$.
As the values of
$\lambda\circ \gamma_1$ and $\lambda\circ\gamma_2$ agree at $t_0$,
and $y'(t)=f_\lambda(t,y(t))$ satisfies local uniqueness of solutions
by~(b), we have $\lambda\circ \gamma_1|_{I_1\cap I_2}=
\lambda\circ\gamma_2|{I_1\cap I_2}$ by Lemma~\ref{lem-loc-glob}.
As the $\lambda\in\Lambda$ separate points on~$E$,
$\gamma_1|_{I_1\cap I_2}=\gamma_2|_{I_1\cap I_2}$ follows.
\end{prf}
\begin{rem}\label{typic-unq}
In typical applications of Proposition~\ref{uniq-pl},
each~$E_\lambda$ is a normed space and $f_\lambda$
satisfies a local Lipschitz condition; the final condition
of Proposition~\ref{uniq-pl}(b) is then automatically satisfied,
by Proposition~\ref{uniq-normd}.
\end{rem}
\begin{ex}\label{reallynee}
Let $F$ be a locally convex space, $E:=F^n$ with $n\in\N$ and $J\sub\R$ be a non-degenerate
interval. Let $a_{jk}\colon J\to\R$ be continuous functions
for $j,k\in\{1,\ldots,n\}$. Moreover, let $b=(b_1,\ldots, b_n)\colon J\to F^n$
be a continuous function.
For $t\in J$, let $A(t)\colon F^n\to F^n$ be the continuous linear mapping
taking $y=(y_1,\ldots, y_n)\in F^n$ to
\[
A(t)(y):=\left(\sum_{k=1}^n a_{1k}(t)y_k,\ldots, \sum_{k=1}^n a_{nk}(t)y_k\right).
\]
Then the linear differential equation
\begin{equation}\label{speclin}
y'(t)=A(t)y(t)+b(t)
\end{equation}
satisfies local uniqueness of solutions.\\[2mm]
To see this, let $\Lambda$ be the set of all $\lambda^n:=\lambda\times\cdots\times \lambda\colon
F^n\to\R^n$, $(y_1,\ldots, y_n)\mto(\lambda(y_1),\ldots,\lambda(y_n))$,
for $\lambda$ ranging through the dual space~$F'$ of all continuous linear functionals
of~$F$.
Then~$\Lambda$ separates points on~$E$.
For $t\in J$, define a matrix $B(t)\in \R^{n\times n}$ via
\[
B(t):=(a_{jk}(t))_{j,k=1}^n.
\]
Then $B\colon J\to \R^{n\times n}$ and $b_\lambda:=\lambda^n\circ b\colon J\to \R^n$
are continuous mappings.
For the right hand side $f\colon J\times E\to E$, $f(t,y):=A(t)(y)+b(t)$
of~(\ref{speclin}), we have
\[
\lambda^n\circ f=f_{\lambda^n}\circ (\id_J\times \lambda^n)
\]
with $f_{\lambda^n}\colon J\times\R^n\to \R^n$, $f_{\lambda^n}(t,y)=B(t)y+b_\lambda(t)$.
Since
\[
y'(t)=B(t)y+b_\lambda(t)
\] 
satisfies local uniqueness
of solutions by Proposition~\ref{uniq-normd} and Example~\ref{lin-loclip},
also $y'(t)=f(t,y(t))$ satisfies local uniqueness of solutions,
by Proposition~\ref{uniq-pl}.
\end{ex}
The third
criterion for local uniqueness -- existence of local flows --
is of particular importance.
\begin{defn}\label{defn-loc-flo}
Let $J\sub \R$ be a non-degenerate interval, $E$ be a locally convex space,
$U\sub E$ be a locally convex subset with dense interior
and $f\colon W\to E$ be a function on an open subset $W\sub J\times U$.
Let $k\in\N_0\cup\{\infty\}$.
We say that the differential equation $y'(t)=f(t,y(t))$
\emph{admits local $C^k$-flows} if, for all $(\wb{t},\wb{y}) \in W$,
there exist a relatively open interval $I\sub J$ with $\wb{t}\in I$,
an open neighborhood~$V$ of $\wb{y}$ in~$U$
with $I\times V\sub W$
and a $C^k$-function
\[
\Phi\colon I\times I\times V\to E
\]
with the following properties:
\begin{description}[(D)]
\item[(a)]
For all $(t_0,y_0)\in I\times V$, the function
$I\to E$, $t\mto \Phi_{t,t_0}(y_0):=\Phi(t,t_0,y_0)$ is a solution to
the initial value problem (\ref{theiniprob});
\item[(b)]
There is an open $\wb{y}$-neighborhood $Y\sub V$
such that $\Phi_{t_1,t_0}(Y)\sub V$ for all $t_0,t_1\in I$ and
\[
\Phi_{t_2,t_1}(\Phi_{t_1,t_0}(y_0))=\Phi_{t_2,t_0}(y_0)
\;\;\mbox{for all $\,t_0,t_1\in I$ and $y_0\in Y$.}
\]
\end{description}
If $J$ is an open interval, $U\sub E$ an open subset, $f$ is real analytic
and $\Phi$ as before can be chosen as a real analytic map
for all $(\wb{t},\wb{y})\in W$,
then we say that $y'(t)=f(t,y(t))$ \emph{admits real analytic local flows}.
\end{defn}
\begin{prop}\label{locflow-uniq}
Let $J\sub \R$ be a non-degenerate interval, $E$ be a locally convex space,
$U\sub E$ be a locally convex subset with dense interior
and\linebreak $f\colon W\to E$ be a function on an open subset $W\sub J\times U$.
If the differential equation $y'(t)=f(t,y(t))$
admits local $C^1$-flows, then it satisfies local uniqueness
of solutions.
\end{prop}
\begin{prf}
Let $\gamma_j\colon I_j\to E$
be solutions to $y'(t)=f(t,y(t))$
for $j\in \{1,2\}$ and $\wb{t}\in I_1\cap I_2$ such that
$\wb{y}:=\gamma_1(\wb{t})=\gamma_2(\wb{t})$.
To see that $\gamma_1$ and $\gamma_2$ coincide on a neighborhood of~$\wb{t}$
in $I_1\cap I_2$, we may assume that $I_1\cap I_2\not=\{t_0\}$
(excluding only a trivial case). Thus $I_1\cap I_2$ is a non-degenerate interval.
Let $I$, $V$, $\Phi$, and~$Y$
be as in Definition~\ref{defn-loc-flo}.
For $(t,t_0,y_0)\in I\times I\times V$, we write
\[
(\partial_2\Phi)(t,t_0,y_0):=\frac{\partial\Phi}{\partial t_0}(t,t_0,y_0)
\]
for the partial derivative of~$\Phi$ with respect to the second variable.
There exists a relatively open interval $K\sub I_1\cap I_2\cap I$
with $\wb{t}\in K$ such that $\gamma_1(K)\sub Y$,
$\gamma_2(K)\sub Y$ and $\Phi_{t,\wb{t}}(\wb{y})\in Y$ for all $t\in K$.
After shrinking~$K$ if necessary, we can also assume that
\[
\theta_j(t):=\Phi_{\wb{t},t}(\gamma_j(t))\in Y \;\;\mbox{for all
$t\in K$ and $j\in\{1,2\}$.}
\]
It suffices to show that
\[
\gamma_j(t)=\Phi_{t,\wb{t}}(\wb{y})\;\;\mbox{for
$t\in K$ and $j\in\{1,2\}$.}
\] 
Since $\Phi_{t,\wb{t}}\circ \Phi_{\wb{t},t}|_Y=\id_Y$ for all
$t\in I$, the map $\Phi_{\wb{t},t}|_Y$ is injective.
Hence $\gamma_1|_K=\gamma_2|_K$
will hold if we can show that both $\theta_1$ and $\theta_2$ coincide
with
\[
\theta\colon K\to E,\quad t\mto \Phi_{\wb{t},t}(\Phi_{t,\wb{t}}(\wb{y}))=\wb{y}.
\]
Since $\theta_j(\wb{t})=\wb{y}=\theta(\wb{t})$ for $j\in\{1,2\}$,
the latter will hold if we can show that
\[
\theta_j'(t)=\theta'(t)=0\;\;\mbox{for  $t\in K$.}
\]
For all $t\in I$ and $z\in Y$, we have
\[
z=\Phi_{\wb{t},t}(\Phi_{t,\wb{t}}(z))
\]
and hence, differentiating with respect to~$t$,
\begin{equation}\label{tool-flw}
0=\partial_2\Phi(\wb{t},t,\Phi_{t,\wb{t}}(z))+d\Phi_{\wb{t},t}(\Phi_{t,\wb{t}}(z),f(t,\Phi_{t,\wb{t}}(z)))\;\;
\mbox{for } \quad t\in I, z\in Y,
\end{equation}
exploiting that $\frac{d}{dt}\Phi_{t,\wb{t}}(z)=f(t,\Phi_{t,\wb{t}}(z))$.
For $t\in K$ and $j\in \{1,2\}$,
we have $\gamma_j(t)=\Phi_{t,\wb{t}}(z)$ with $z:=\theta_j(t)\in Y$.
Since $\gamma_j'(t)=f(t,\gamma_j(t))=f(t,\Phi_{t,\wb{t}}(z))$, we get
\begin{align*}
\theta_j'(t) 
&=\partial_2\Phi(\wb{t},t,\gamma_j(t))+d\Phi_{\wb{t},t}(\gamma_j(t),\gamma_j'(t))\\
&=\partial_2\Phi(\wb{t},t,\Phi_{t,\wb{t}}(z))+
d\Phi_{\wb{t},t}(\Phi_{t,\wb{t}}(z),f(t,\Phi_{t,\wb{t}}(z)))=0
\end{align*}
as a special case of~(\ref{tool-flw}).
\end{prf}
\begin{rem}
Let $J\sub \R$ be a non-degenerate interval, $E$ be a locally convex space,
$U\sub E$ be a locally convex subset with dense interior
and\linebreak $f\colon W\to E$ be a function on an open subset $W\sub J\times U$.
Let $k\in\N_0\cup\{\infty,\omega\}$; if $k=\omega$, assume that~$J$ and~$U$ are open.
If the differential equation $y'(t)=f(t,y(t))$
satisfies local existence and local uniqueness
of solutions, then the following conditions are equivalent:
\begin{description}[(D)]
\item[(a)]
$y'(t)=f(t,y(t))$ admits local $C^k$-flows;
\item[(b)]
The domain $\Omega$ of the maximal flow $\Fl\colon \Omega\to E$
associated with $y'(t)=f(t,y(t))$ is a neighborhood of
$(\wb{t},\wb{t}, \wb{y})$ in $J\times J\times U$ for all $(\wb{t},\wb{y})\in W$,
and is $C^k$ on a (possibly smaller) open neighborhood.
\end{description}
In fact, if $\Phi\colon I\times I\times V\to E$ is as in Definition~\ref{defn-loc-flo},
then $I\times I\times V\sub\Omega$ and $\Fl|_{I\times I\times V}=\Phi$ is~$C^k$;
thus (a)$\Rightarrow$(b). Conversely, assume that $\Omega$ contains an open neighborhood
$I\times I\times V$
of $(\wb{t},\wb{t},\wb{y})$ for all $(\wb{t},\wb{y})\in W$ such that $\Phi:=\Fl|_{I\times I\times V}$
is~$C^k$. Since $\Phi(\wb{t},\wb{t},\wb{y})=\wb{y}$ and~$\Phi$ is continuous, after shrinking~$I$
if necessary we may assume that there exists an open $\wb{y}$-neighborhood $Y\sub V$
such that $\Phi(I\times I\times Y)\sub V$. Condition~(b) of Definition~\ref{defn-loc-flo} is satisfied by
Lemma~\ref{basics-flow}(b).
Thus (b)$\Rightarrow$(a).
\end{rem}
\begin{lem}\label{power-fix}
Consider a self-map
$f\colon X\to X$ of a set~$X$.
If there exists $n\in\N$ such that the iterate $g:=f^n\colon X\to X$ has a unique fixed point~$x_0$,
then $x_0$ also is a fixed point of~$f$.
If, moreover, $X$ is a topological space, $f$ is continuous at~$x_0$
and $y\in X$ is an element such that $g^k(y)\to x_0$ as $k\to\infty$,
then also $f^k(y)\to x_0$.
\end{lem}
\begin{prf}
The identity $f^n(f(x_0))=f(f^n(x_0))=f(x_0)$
shows that $f(x_0)$ is
a fixed point of~$f^n$. Hence $f(x_0)=x_0$, by uniqueness of the fixed point of~$f^n$.

To establish the second assertion, it suffices to observe that $f^{j+nk}(y)=f^j(g^k(y))\to f^j(x_0)=x_0$
as $k\to\infty$, for each $j\in\{0,1,\ldots,n-1\}$.
\end{prf}
\begin{thm}[Quantitative Existence Theorem]\label{quant-PL}
Let $(E,\|\cdot\|)$ be a Banach space,
$a<b$ be real numbers and $t_0\in [a,b]$.
Let $y_0\in E$, $R>0$ and $f\colon [a,b]\times \wb{B}^E_R(y_0)\to E$
be a bounded continuous function satisfying a global Lipschitz condition in its
second argument. Abbreviate
\[
M:=\|f\|_\infty:=\sup\{\|f(t,y)\|\colon (t,y)\in [a,b]\times \wb{B}^E_R(y_0)\}\in [0,\infty[
\]
and define\footnote{Read $R/0:=\infty$ if $M=0$.}
\[
\ve_1:=\min\{t_0-a,R/M\},\quad \ve_2:=\min\{b-t_0,R/M\}.
\]
Then the following holds:
\begin{description}[(D)]
\item[{\rm(a)}]
There is a unique $C^1$-function
$\gamma\colon I \to \wb{B}^E_R(y_0)$
on the non-degenerate interval $I:=[t_0-\ve_1,t_0+\ve_2]$
which solves the initial value problem
\begin{equation}\label{theinivp}
y'(t)=f(t,y(t)),\quad y(t_0)=y_0.
\end{equation}
\item[{\rm(b)}]
If we define $\gamma_0\in C(I,E)$ via $\gamma_0(t):=y_0$
and $\gamma_n\in C(I,E)$ for $n\in\N$ via
\[
\gamma_n(t):=y_0+\int_{t_0}^t f(s,\gamma_{n-1}(s))\, ds\;\;\mbox{for $\,t\in I$,}
\]
then $\gamma_n\to \gamma$ uniformly.
\item[{\rm(c)}]
If $R/M\geq\max\{t_0-a,b-t_0\}$, then $I=[a,b]$ and thus $\gamma$
is defined on all of $[a,b]$.
\end{description}
\end{thm}
\begin{prf}
(c) is clear from the definition of~$\ve_1$ and~$\ve_2$.

(a) and~(b). The solution~$\gamma$ is unique (if it exists), by Proposition~\ref{uniq-normd}.
By hypothesis, we have
\[
L:=\sup\{\Lip(f(t,\cdot))\colon t\in [a,b]\}<\infty.
\]
Let $I:=[t_0-\ve_1,t_0+\ve_2]$ and $F:=C(I,E)$, endowed with the supremum norm.
We consider the map
\[
g\colon \{\eta\in C(I,E)\colon \|\eta\|_\infty\leq R\}\to C(I,E),\quad
g(\eta)(t):=\int_{t_0}^t f(s,y_0+\eta(s))\, ds.
\]
We claim that
\begin{description}[(D)]
\item[(i)] $g$ is a self-map of $\wb{B}^F_R(0)\sub C(I,E)$; and
\item[(ii)]
$g^n\colon \wb{B}^F_R(0)\to \wb{B}^F_R(0)$ is Lipschitz
for each $n\in\N$ with
\begin{equation}\label{thelipn}
\Lip(g^n)\leq\frac{L^n}{n!}\ve^n,
\end{equation}
where $\ve:=\max\{\ve_1,\ve_2\}\leq R/M$.
Hence $g^n$ is a contraction
for all $n\in\N$ such that\footnote{Note that this condition is satisfied for large~$n$
since $\sum_{n=0}^\infty\frac{L^n}{n!}\ve^n=e^{L\ve}<\infty$.}
$\frac{L^n}{n!}\ve^n<1$.
\end{description}
If this is true, fix~$n$ as in~(ii); then $g^n$ has a unique fixed point~$\eta$
(by Banach's Fixed Point Theorem, Lemma~\ref{ban-fix}), and $\eta=\lim_{k\to\infty}g^{nk}(0)$.
Then also~$g$ has~$\eta$ as a fixed point and
\begin{equation}\label{iterg}
\eta=\lim_{k\to\infty}g^k(0)\;\;\mbox{in $C(I,E)$,}
\end{equation}
by Lemma~\ref{power-fix}. 
By Lemma~\ref{via-integral-eq}, $\eta$ solves
\begin{equation}\label{theinivp2}
z'(t)=f(t,z(t)+y_0),\quad z(t_0)=0.
\end{equation}
Identify $y_0$ with the constant function $I\to E$, $t\mto y_0$.
By the preceding, $\gamma:=\eta+y_0$ solves~(\ref{theinivp})
and since $\gamma_k=g^k(y_0)+y_0$,
we see that $\gamma_k\to\gamma$ uniformly as $k\to\infty$.
To complete the proof of (a) and~(b), it therefore only remains to prove the claim.

(i) Let $\zeta\in \wb{B}^F_R(0)$. For each $t\in I$, we have
\[
\|g(\zeta)(t)\|=\left\|\int_{t_0}^tf(s,\zeta(s)+y_0)\, ds\right\|
\leq |t-t_0|M\leq R
\]
as $\ve_1,\ve_2\leq R/M$. Hence $g(\zeta)\in\wb{B}^F_R(0)$.

(ii) will follow if, for each $n\in\N$, we can show that
\begin{equation}\label{powers-esti}
\|(g^n(\eta_2)-g^n(\eta_1))(t)\|\leq \frac{|t-t_0|^n}{n!}L^n\|\eta_2-\eta_2\|_\infty
\end{equation}
for all $\eta_1,\eta_2\in \wb{B}^F_R(0)$ and $t\in I$
(since $|t-t_0|\leq\max\{\ve_1,\ve_2\}=\ve$). We proceed by induction;
for $n=1$, we have
\begin{eqnarray*}
\|g(\eta_2)(t)-g(\eta_1)(t)\|&=&
\left\|\int_{t_0}^t \big(f(s,\eta_2(s)+y_0)-f(s,\eta_2(s)+y_0)\big)\,ds\right\|\\
&\leq & \left|\int_{t_0}^t\|f(s,\eta_2(s)+y_0)-f(s,\eta_1(s)+y_0)\|\, ds\right|\\
&\leq& |t-t_0|L\|\eta_2-\eta_1\|_\infty,
\end{eqnarray*}
as required.
If (\ref{powers-esti}) holds for~$n$, we get
\begin{eqnarray*}
\lefteqn{\|g^{n+1}(\eta_2)(t)-g^{n+1}(\eta_1)(t)\|}\qquad\\
&=&
\left\|\int_{t_0}^t \big(f(s,g^n(\eta_2)(s)+y_0)-f(s,g^n(\eta_1)(s)+y_0)\big)\,ds\right\|\\
&\leq & \left|\int_{t_0}^t\|f(s,g^n(\eta_2)(s)+y_0)-f(s,g^n(\eta_1)(s)+y_0)\|\, ds\right|\\
&\leq &\left|\int_{t_0}^t \frac{|s-t_0|^n}{n!}L^{n+1}\|\eta_2-\eta_1\|_\infty\,ds\right|
= \frac{|t-t_0|^{n+1}}{(n+1)!}L^{n+1}\|\eta_2-\eta_1\|_\infty,
\end{eqnarray*}
using the inductive hypothesis.
This completes the induction step.
\end{prf}
\begin{rem}
The proof of Theorem~\ref{quant-PL} can be shortend
if one uses the Fixed Point Theorem of Banach--Weissinger~\cite{Wsg52}
instead of Banach's Fixed Point Theorem;
the former applies directly to the self-map~$g$,
not only to a power~$g^n$ (see Exercise~\ref{exc-ban-weis}).
We chose to work nonetheless with the standard Banach Fixed Point Theorem,
as we wish to re-use the proof later in this section when we discuss parameters. We can then employ
the results concerning parameter-dependence of fixed points
from Section~\ref{provis-Ban}.
\end{rem}
\begin{cor}\label{qual-exi}
Let $J\sub\R$ be a non-degenerate interval, $(E,\|\cdot\|)$
be a Banach space, $W\sub J\times E$ be an open subset
and $f\colon W\to E$ be a continuous function which satisfies a local Lipschitz condition
in its second argument.
Then $y'(t)=f(t,y(t))$ satisfies local existence of solutions.
\end{cor}
\begin{prf}
Given $(t_0,y_0)\in W$, there exist real numbers $a<b$ such that $[a,b]$
is a $t_0$-neighborhood in~$J$, and $R>0$
such that $[a,b]\times \wb{B}^E_R(y_0)\sub W$.
After shrinking $R$ and the interval $[a,b]$, we may assume
that
\[
L:=\sup_{t\in[a,b]}\Lip f(t,\cdot)|_{\wb{B}^E_R(y_0)}<\infty
\]
and $M:=\sup\{\|f(t,y)\|\colon t\in[a,b],\, y\in\wb{B}^E_R(y_0)\}<\infty$.
By Theorem~\ref{quant-PL}, there exist $\ve_1,\ve_2\geq 0$
such that $I:=[t_0-\ve_1,t_0+\ve_2]$ is a neighborhood of~$t_0$ in $[a,b]$,
and a solution $\gamma\colon I\to E$ to the initial value problem
$y'(t)=f(t,y(t))$, $y(t_0)=y_0$. It only remains to note that $I$ is also a $t_0$-neighborhood
in~$J$.
\end{prf}
\begin{prop}\label{ex-lin-ode}
Let $(E,\|\cdot\|)$ be a Banach space, $J\sub\R$ be a non-degenerate interval,
$b\colon J\to E$ a continuous function and $A\colon J\to\cL(E)$ be a function such that
$\wh{A}\colon J\times E\to E$ is continuous.
Then
\begin{equation}\label{our-lin}
y'(t)=A(t)y(t)+b(t)
\end{equation}
satisfies local existence and local uniqueness of solutions.
For each $(t_0,y_0)\in J\times E$, the maximal solution to
$y'(t)=A(t)y(t)+b(t)$, $y(t_0)=y_0$ is defined on all of~$J$.
\end{prop}
\begin{prf}
By Example~\ref{lin-loclip}, the map
$f \colon J\times E\to E$, $(t,y)\mto A(t)(y)+b(t)$
satisfies a local Lipschitz condition.
Hence (\ref{our-lin}) satisfies local uniqueness
of solutions, by Proposition~\ref{uniq-normd}.
The remaining assertions follow if we can show that,
for $(t_0,y_0)\in J\times E$,
the initial value problem
$y'(t)=A(t)y(t)+b(t)$, $y(t_0)=y_0$
has a solution on each non-degenerate compact subinterval $I\sub J$ with $t_0\in I$.
Fix~$I=[a,b]$. By Example~\ref{lin-loclip}, we have
\[
L:=\sup_{t\in I}\Lip f(t,\cdot)<\infty.
\]
We define
$g\colon C(I,E)\to C(I,E)$ via $g(\gamma)(t):=y_0+\int_{t_0}^tf(s,\gamma(s))\,ds$.
As in the proof of Theorem~\ref{quant-PL}, an induction on $n\in\N_0$ shows that
\[
\|g^n(\gamma)(t)- g^n(\eta)(t)\|\leq \frac{L^n|t-t_0|^n}{n!}\|\gamma-\eta\|_\infty
\]
for all $\gamma,\eta\in C(I,E)$ and $t\in I$, whence
$\Lip(g^n)\leq \frac{(b-a)^nL^n}{n!}$.
We deduce that $\Lip(g^n)<1$ for large~$n$,
whence~$g$ has a fixed point~$\gamma$, by Lemma~\ref{power-fix}
and Banach's Fixed Point Theorem.
Then~$\gamma$ is defined on all of~$I$ and solves the initial value problem
we consider, by Lemma~\ref{via-integral-eq}.
\end{prf}
\begin{rem}
Note that $A$ need not be continuous as a map $J\to\cL(E)_b$
in Proposition~\ref{ex-lin-ode} (see Exercise~\ref{exc-generato}).
\end{rem}
We now study the dependence of solutions to differential equations on parameters
(and continue to study the dependence on initial conditions).
We begin with terminology.
\begin{defn}\label{defn-ode-par-m}
Let $E$ be a locally convex space, $P$ be a set,
$J\sub \R$ be a non-degenerate interval,
$U\sub E$ be a subset
and $f\colon W\to E$ be a function on a subset $W\sub J\times U\times P$.
For $p\in P$, consider the differential equation
\begin{equation}\label{ode-with-pr}
y'(t)=f(t,y(t),p)
\end{equation}
whose right-hand side is given by the function $W^p\to E$, $(t,y)\mto f(t,y,p)$
on the subset $W^p:=\{(t,y)\in J\times U\colon (t,y,p)\in W\}$ of $J\times U$.
If the differential equations (\ref{ode-with-pr})
satisfy both local existence and local uniqueness of solutions for all $p\in P$,
we let $\gamma_{t_0,y_0,p}\colon I_{t_0,y_0,p}\to E$ be the maximal solution to the
initial value problem
\begin{equation}\label{ivp-manif}
y'(t)=f(t,y(t),p),\quad y(t_0)=y_0
\end{equation}
for $(t_0,y_0,p)\in W$.
We define
\[
\Omega:=\bigcup_{(t_0,y_0,p)\in W}I_{t_0,y_0,p}\times\{(t_0,y_0,p)\}
\sub J\times J\times E\times P
\]
in this case and call the map
\[
\Fl\colon \Omega\to E,\quad \Fl(t,t_0,y_0,p):=\gamma_{t_0,y_0,p}(t)
\]
the associated (maximal) flow.
\end{defn}
We now formulate a prototypical result concerning
the dependence of solutions to differential equations on parameters and inital conditions.
\begin{thm}\label{nice-spec-para}
Let $J\sub \R$ be a non-degenerate interval,
$(E,\|\cdot\|)$ be a Banach space, $F$ be a locally convex space,
$U\sub E$ be an open subset, and $P\sub F$ be a locally convex subset with dense interior.
Let $W\sub J\times U\times P$ be an open subset, $k\in\N_0\cup\{\infty\}$
and
\[
f\colon W\to E,\quad (t,y,p)\mto f(t,y,p)
\]
be a $C^k$-function.
If~$k=0$, assume that~$f$ satisfies a local Lipschitz condition in the $y$-variable.
Then the following holds:
\begin{description}[(D)]
\item[\rm(a)]
The differential equation
\begin{equation}\label{odewithpa}
y'(t)=f(t,y(t),p)
\end{equation}
satisfies local existence and local uniqueness of solutions, for each $p\in P$;
\item[\rm(b)]
The domain $\Omega$ of the flow of {\rm(\ref{odewithpa})}
is open in $J\times J\times U\times P$;
\item[\rm(c)]
The flow $\Fl\colon \Omega\to E$ of~{\rm(\ref{odewithpa})}
is~$C^k$.
\end{description}
If~$J$ is an open interval, $P\sub F$ is an open subset and~$f$ is real analytic,
then also $\Fl$ is real analytic.
\end{thm}
To establish Theorem~\ref{nice-spec-para},
we first develop certain results concerning the local behavior of solutions.
Moreover, tools are introduced to tackle differential equations with analytic right-hand sides.
We then globalize the investigations in the following Section~\ref{sec-ode-mfd},
in the more general setting of differential equations on manifolds.
As we shall see, Theorem~\ref{nice-spec-para} does
not require a separate proof; it will be a special case of
Corollary~\ref{pardep-ode-ban-m}.

The following lemma will help us to gain real analytic solutions
to real differential equations from
complex analytic solutions to complex differential equations
(in which case $E=F_\C$).
\begin{lem}\label{ode-in-subspace}
Let $J\sub\R$ be a non-degenerate interval,
$E$ be a locally convex space, $U\sub E$ be a subset,
$W\sub J\times U$ be a subset and $f\colon W\to E$ be a function such that
the differential equation $y'(t)=f(t,y(t))$ satisfies local uniqueness of solutions.
Let $F\sub E$ be a closed vector subspace
such that
\[
g(t,y):=f(t,y)\in F\quad\mbox{for all $\, (t,y)\in W\cap (J\times F)$.}
\]
Then the differential equation $y'(t)=g(t,y(t))$ satisfies local uniqueness of solutions.
If it also satisfies local existence of solutions, then 
the following holds: If $\gamma\colon I\to E$ is a solution
to $y'(t)=f(t,y(t))$ and $y(t_0)\in F$ for some $t_0\in I$,
then $\gamma(I)\sub F$ and $\gamma$ solves $y'(t)=g(t,y(t))$.
\end{lem}
\begin{prf}
Since solutions to $y'(t)=g(t,y(t))$ also solve $y'(t)=f(t,y(t))$,
the first differential equation inherits local uniqueness of solutions from the second
one.
Now let $\gamma\colon I\to E$ be a solution to $y'(t)=f(t,y(t))$ such that $\gamma(t_0)\in F$ for some
$t_0\in I$. We may assume that~$\gamma$ is the maximal solution
for the initial data $(t_0,\gamma(t_0))$.
As~$\gamma$ is continuous and~$F$ is closed in~$E$,
the preimage~$\gamma^{-1}(F)$
is closed in~$I$. By hypothesis, $\gamma^{-1}(F)$ is non-empty.
Since~$I$ is connected, we shall have $\gamma^{-1}(F)=I$
(which completes the proof) if we can show that $\gamma^{-1}(F)$ is open in~$I$.
To this end, let $t_1\in \gamma^{-1}(F)$.
By local existence of solutions, the initial value problem
\[
y'(t)=g(t,y(t)),\quad y(t_1)=\gamma(t_1)
\]
has a solution $\eta\colon K\to F$ on a relatively open subinterval $K\sub J$.
Since~$\gamma$ is maximal, $K\sub I$ and $\eta=\gamma|_K$
follows, whence $K\sub\gamma^{-1}(F)$.
Thus $\gamma^{-1}(F)$ is a neighborhood of~$t_1$ and hence open in~$I$.
\end{prf}
As a tool for the discussion of real analytic solutions to
differential equations, let us briefly discuss complex differential equations.
\begin{numba}\label{def-cx-ode}
Let $E$ be a complex locally convex space, $W\sub \C\times E$ be a subset and
$f\colon W\to E$ be a function.\medskip

\noindent
(a) We say that a function $\gamma \colon I\to E$
on a convex subset $I\sub\C$ with non-empty interior
is a \emph{solution to the complex differential equation}
\begin{equation}\label{theode-cx}
y'(z)=f(z,y(z))
\end{equation}
if $\gamma$ is a $C^1_\C$-function,
$(z,\gamma(z))\in W$ and $\gamma'(z)=f(z,\gamma(z))$
for all $z\in I$,
where $\gamma'(z)$ is understood as the complex derivative $\frac{d\gamma}{dz}(z)$
($:=d\gamma(z,1)$, cf.\ Definition~\ref{singlevarcx}).\medskip

\nin (b) If $(z_0,y_0)\in W$ and $\gamma\colon I\to E$ is a solution to the complex differential equation~(\ref{theode-cx})
such that $z_0\in I$ and $\gamma(z_0)=y_0$,
then $\gamma$ is called a \emph{solution to the initial value problem}
\begin{equation}\label{theiniprob-cx}
\left\{
\begin{array}{rcl}
y'(z) &=& f(z,y(z))\\
y(z_0)&=&y_0.
\end{array}
\right.
\end{equation}
\end{numba}
\begin{rem}
(a) We shall only discuss complex differential equations in a complex Banach space~$E$.
Of course, we are mostly interested in the case that $I\sub \C$ and $W\sub \C\times E$
are open subsets, in which case the $C^1_\C$-property of a solution~$\gamma$
is equivalent to~$\gamma$ being complex analytic
(see Proposition~\ref{analytsingle}).
However, our discussion of complex differential equations in Banach spaces entirely avoids the use of
power series and is based throughout on the $C^1_\C$-property (or $C^k_\C$-property)
of the functions at hand, rather than complex analyticity.
In our main result in this context,
Proposition~\ref{local-ivp-dep-Ck}, openness of~$P$ and~$J$
would be an unnecessary extra condition.
We therefore stated~\ref{def-cx-ode} in full generality, including the case
of non-open sets.\medskip

(b) In the following, we shall frequently use the latter ``$t$'' (or others) in place of~$z$
when dealing with complex differential equations,
to enable a uniform notation for both real and complex differential equations. Yet, all derivatives
with respect to~$t$ are intended as complex derivatives in the case of a complex differential equation.
\end{rem}
It can be useful to pass to normalized initial value problems.
\begin{numba}\label{normal-ivp}
Let
$J\sub\K$ be a convex subset with non-empty interior,
$t_0\in J$ and $\gamma\colon J\to E$ be a solution to the 
initial value problem
\[
y'(t)=f(t,y(t)),\quad y(t_0)=y_0
\]
over~$\K$. Let $t\in J$. Then
\begin{equation}\label{gamm-via-et}
\gamma(t)=\eta(1)+y_0
\end{equation}
for the $C^1_\R$-curve
\[
\eta\colon [0,1]\to E,\quad \eta(\tau):=\gamma(t_0+\tau(t-t_0))-y_0
\]
which solves the initial value problem
\begin{equation}\label{normalized-prob}
w'(\tau)=(t-t_0)f(t_0+\tau(t-t_0),w(t)+y_0),\quad w(0)=0
\end{equation}
which is normalized in the sense that the initial conditions
are $\tau_0=0$ and $w(\tau_0)=0$.
If we can find the solution~$\eta$ to~(\ref{normalized-prob}),
then we also know $\gamma(t)$ (because of~(\ref{gamm-via-et})).
This idea is the backbone of the following proof.
\end{numba}
\begin{thm}\label{local-ivp-dep-Ck}
Let $(E,\|\cdot\|)$ be a Banach space and~$F$ be a locally convex space
over $\K\in\{\R,\C\}$.
Given $(\wb{t},\wb{y},\wb{p})\in \K\times E\times F$,
let $J\sub\K$ be a convex subset with non-empty interior
such that $\wb{t}\in J$;
let $U\sub E$ be an open $\wb{y}$-neighborhood and $P\sub F$ be a locally convex subset with dense
interior such that $\wb{p}\in P$.
Let $k\in \N_0\cup\{\infty\}$ and
$f\colon J\times U\times P\to E$
be a $C^k_\K$-function;
if $k=0$, assume that~$f$ satisfies a local Lipschitz condition in its second argument;
if $\K=\C$, assume that $k\geq 1$.
Then there exist a convex, relatively open $\wb{t}$-neighborhood $J_0\sub J$,
an open $\wb{y}$-neighborhood $U_0\sub U$, and
a relatively open $\wb{p}$-neighborhood $P_0\sub P$,
with the following properties:
\begin{description}[(D)]
\item[\rm(a)]
For all $(t_0,y_0,p)\in J_0\times U_0\times P_0$, the initial value problem
\begin{equation}\label{ini-over-K}
y'(t)=f(t,y(t),p),\quad y(t_0)=y_0
\end{equation}
over~$\K$ has a unique $C^k_\K$-solution $\gamma_{t_0,y_0,p}\colon J_0\to E$;\vspace{.8mm}
\item[\rm(b)]
The map
$\Phi\colon J_0\times J_0\times U_0\times P_0\to E$, $(t,t_0,y_0,p)\mto\gamma_{t_0,y_0,p}(t)$
is~$C^k_\K$.
\end{description}
If $\K=\R$, all of $J\sub\R$, $U\sub E$, and $P\sub F$ are open subsets
and $f$ is real analytic, then we can also achieve that~$\Phi$ is real analytic.
\end{thm}
\begin{prf}
If $\K=\R$, the uniqueness in (a) follows from Proposition~\ref{uniq-normd}.
If $\K=\C$, then $\gamma_{t_0,y_0,p}(t)$ in~(a) is unique
since $\gamma_{t_0,y_0,p}(t)=\eta(1)$ if we define
$\eta\colon [0,1]\to E$, $\eta(\tau):=\gamma_{t_0,y_0,p}(t_0+\tau(t-t_0))-y_0$,
which is a solution to the initial value problem $w'(\tau)=(t-t_0)f(t_0+\tau(t-t_0),w(t)+y_0,p)$, $w(0)=0$
and hence unique as $(\tau,w)\mto (t-t_0)f(t_0+\tau(t-t_0),w+y_0,p)$
satisfies a local Lipschitz condition in the $w$-variable (cf.\ \ref{normal-ivp}).

Let us construct solutions now if $k\in \N_0\cup\{\infty\}$.
After shrinking~$U$, we may assume that $U=B^E_{2R}(\wb{y})$
for some $R>0$. After shrinking $J$, $R$, and $P$, we may assume that
there exist $M,L\in[0,\infty[$ such that
\begin{equation}\label{expli-1-ode}
(\forall t\in J)\,(\forall p\in P)\quad \Lip(f(t,\cdot,p))\leq L,\quad\mbox{and $\,\|f\|_\infty\leq M$.}
\end{equation}
After shrinking~$J$, we may assume that its diameter $r:=\diam(J):=$\linebreak
$\sup\{|s-t|\colon
s,t\in J\}$ is so small that
\begin{equation}\label{expli-2-ode}
rL<1\quad\mbox{and}\quad rM\leq R.
\end{equation}
We now establish (a) and~(b) with $U_0:=B^E_R(\wb{y})$
and $J_0:=J$, $P_0:=P$.
Consider the $C^k_\R$-function $h\colon [0,1]\times (B^E_{2R}(\wb{y})\times
J_0\times J_0\times P_0)\to E$,
\[
h(\tau, y,t,t_0,p):=(t-t_0)f(t_0+\tau(t-t_0),y,p).
\]
Then~$h$ is $C^{0,k}_\R$, in particular.
Moreover, $h$ is $C^{0,k}_\C$ if $\K=\C$, as the mapping
$h(\tau,\cdot)\colon
B^E_{2R}(\wb{y})\times
J_0\times J_0\times P_0\to E$ is~$C^k_\C$ for each $\tau\in[0,1]$
(whence Exercise~\ref{exc-cxC0k} applies).
Now $\|h\|_\infty\leq rM$ and
\[
\Lip h(\tau ,\cdot,t,t_0,p)\leq |t-t_0|\Lip(f(t_0+\tau(t-t_0),\cdot,p))
\leq rL
\]
for all $(\tau,t,t_0,p)\in [0,1]\times J_0\times J_0\times P_0$.
Consider $Y:=C([0,1],E)$ with the supremum norm.
By the Quantitative Existence Theorem (Theorem~\ref{quant-PL}) and its proof
(with $\ve_1=0,\ve=\ve_2=1$ and $n=1$),
the initial value problem
\begin{equation}\label{the-norma}
w'(\tau)=h(\tau,w(\tau)+y_0,t,t_0,p),\quad w(0)=0
\end{equation}
has a unique solution $\eta_{t,t_0,y_0,p}\colon [0,1]\to E$
for all $(t,t_0,y_0,p)$ in the set\linebreak
$J_0\!\times\! J_0\!\times\! B^E_R(\wb{y})\! \times \! P_0$,
and $\eta_{t,t_0,y_0,p}$ is the unique fixed point of the contraction
\[
g_{t,t_0,y_0,p}\colon \wb{B}^Y_R(0)\to \wb{B}^Y_R(0)
\]
with contraction constant~$rL$,
where $g\colon J_0\times J_0\times B^E_R(\wb{y})\times P_0\times \wb{B}^Y_R(0)
\to \wb{B}^Y_R(0)$ is given by
\[
g(t,t_0,y_0,p,\psi)(\tau):=\int_0^\tau h(s,\psi(s)+y_0,t,t_0,p)\,ds
\]
for $\tau\in [0,1]$ (and $g_{t,t_0,y_0,p}:=g(t,t_0,y_0,p,\cdot)$, as usual).
Thus
\begin{equation}\label{thusbFck}
g(t,t_0,y_0,p,\psi)=I\big((h^{t,t_0,p})_*(\psi+y_0)),
\end{equation}
where the map
$I\colon C([0,1],E)\to C([0,1],E)$ determined by $I(\psi)(\tau):=\int_0^\tau\psi(s)\,ds$
is continuous $\K$-linear and the map
\[
C([0,1],B^E_{2R}(\wb{y}))\times J_0\times J_0\times P\to C([0,1],E),\;\,
(\psi,t,t_0,p)\mto (h^{t,t_0,p})_*(\psi)
\]
is~$C^k_\K$, by Proposition~\ref{pfctisCkpar}
(where $C([0,1],B^E_{2R}(\wb{y}))=B^Y_{2R}(\wb{y})$).
As a consequence, $g$ is~$C^k_\K$, whence also
\[
\Psi\colon J_0\times J_0\times B^E_R(\wb{y})\times P_0\to C([0,1],E),\;\,
(t,t_0,y_0,p)\mto \eta_{t,t_0,y_0,p}
\]
is~$C^k_\K$, by Proposition~\ref{fp-Ck-dep}(a). 
The point evaluation
\[
\ev_1\colon C([0,1],E)\to E,\quad \psi\mto\psi(1)
\]
at~$1$ being continuous $\K$-linear, we deduce that
\[
\ev_1\circ\, \Psi \colon J_0\times J_0\times B^E_R(\wb{y})\times P_0\to E,\;\,
(t,t_0,y_0,p)\mto \eta_{t,t_0,y_0,p}(1)
\]
is~$C^k_\K$, and hence also the map
\[
\Phi\colon J_0\times J_0\times B^E_R(\wb{y})\times P_0\to E,\;\,
(t,t_0,y_0,p)\mto \eta_{t,t_0,y_0,p}(1)+y_0.
\]
Given $(t,t_0,y_0,p,\tau)\in J_0\times J_0\times B^E_R(\wb{y})\times P_0\times [0,1]$,
both $\eta_{t_0+s(t-t_0),t_0,y_0,p}$
and the function $[0,1]\to E$, $\tau\mto\eta_{t,t_0,y_0,p}(s\tau)$
solve the initial value problem
$w'(\tau)=h(\tau,w(\tau)+y_0,t_0+s(t-t_0),t_0,p)$, $w(0)=0$,
whence both functions coincide. Taking $\tau=1$, we deduce that
\begin{eqnarray}
\eta_{t,t_0,y_0,p}(s)+y_0&=&\eta_{t_0+s(t-t_0),t_0,y_0,p}(1)+y_0\nonumber\\
&=&\Phi(t_0+s(t-t_0),t_0,y_0,p)\label{enaderiv}
\end{eqnarray}
for all
$(t,t_0,y_0,p,s)\in J_0\times J_0\times B^E_R(\wb{y})\times P_0\times [0,1]$.

The case $\K=\R$:
Let $(t_0,y_0,p)\in J_0\times B^E_R(\wb{y})\times P_0$.
The function $\gamma\colon J_0\to E$, $\gamma(t):=\Phi(t,t_0,y_0,p)$
satisfies $\gamma(t_0)=y_0$.
In order that~$\gamma$ is a solution to (\ref{ini-over-K}),
it suffices to show that the restriction $\gamma|_K$
is so for each $t\in J_0\setminus \{t_0\}$,
with $K:=\{t_0+s(t-t_0)\colon s\in [0,1]\}$.
Now
\[
\gamma(t_0+s(t-t_0))=\Phi(t_0+s(t-t_0),t_0,y_0,p)=
\eta_{t,t_0,y_0,p}(s)+y_0
\]
is a $C^1$-function of $s\in [0,1]$.
Since $[0,1]\to K$, $s\mto t_0+s(t-t_0)$
is a $C^1$-diffeomorphism,
we deduce that $\gamma|_K$ is $C^1$.
Moreover,
\begin{eqnarray*}
(t-t_0)\gamma'(t_0+s(t-t_0))&=&\frac{d}{ds}\gamma(t_0+s(t-t_0))
=\eta_{t,t_0,y_0,p}'(s)\\
&=&
h(s,\eta_{t,t_0,y_0,p}(s)+y_0,t,t_0,p)\\
&=&h(s,\gamma(t_0+s(t-t_0)),t,t_0,p)\\
&=& (t-t_0)f(t_0+s(t-t_0),\gamma(t_0+s(t-t_0)),p),
\end{eqnarray*}
whence $\gamma'(x)=f(x,\gamma(x),p)$ for all $x\in K$.
Thus $\gamma_{t_0,y_0,p}:=\gamma$
solves~(\ref{ini-over-K}), and we note that~(a) and~(b) are established.\medskip

The case $\K=\C$:
Then $\Phi$ is $C^k_\C$ with $k\geq 1$.
Hence, for $(t_0,y_0,p)\in J_0\times B^E_R(\wb{y})\times P_0$, the function
\[
\gamma\colon J_0\to E,\quad t\mto \Phi(t,t_0,y_0,p)
\]
is $C^k_\C$. We have
$\gamma(t_0)=y_0+\eta_{t_0,t_0,y_0,p}(1)=y_0$.
The complex derivative $\gamma'\colon J_0\to E$
being continuous, $\gamma$ will satisfy the initial value problem~(\ref{ini-over-K})
if we can show that
\begin{equation}\label{slightly-less}
\gamma'(t)=f(t,\gamma(t),p)\;\,\mbox{for all $t\in J_0\setminus\{t_0\}$.}
\end{equation}
But~(\ref{enaderiv}) and the Chain Rule yield
\begin{eqnarray*}
(t-t_0)\,\gamma'(t)& =& \frac{d}{ds}\Big|_{s=1}\gamma(t_0+s(t-t_0))=
\eta'_{t,t_0,y_0,p}(1)\\
&=& (t-t_0)f(t,\eta_{t,t_0,y_0,p}(1)+y_0,p) =(t-t_0)f(t,\gamma(t),p)
\end{eqnarray*}
(with $s\in [0,1]$), from which~(\ref{slightly-less}) follows.\medskip

The real analytic case: Assume that $J\sub\R$, $U\sub E$ and $P\sub F$
are open and $f$ is real analytic. Then $f$ has a complex analytic extension
$\wt{f}\colon Y\to E_\C$,
defined on an open neighborhood~$Y$ of $J\times U\times P$ in
$\C\times E_\C\times F_\C$.
Let $\wt{J}\sub\C$ be a convex open neighborhood of~$\wb{t}$,
and $\wt{U}\sub E_\C$ as well as  $\wt{P}\sub F_\C$ be open neighborhoods of~$\wb{y}$
and~$\wb{p}$, respectively, such that $\wt{J}\times \wt{U}\times \wt{P}\sub Y$.
After shrinking the latter neighborhoods, we may assume that
$\wt{J}\cap \R\sub J$, $\wt{U}\cap E\sub U$ and $\wt{P}\cap F\sub P$.
After replacing $J$, $U$ and $P$ with the preceding intersections, we may assume that
\[
J=\wt{J}\cap \R,\quad U=\wt{U}\cap E,\;\,\mbox{and}\;\, \wt{P}\cap F=P.
\]
By the complex case of the proposition already established,
there exist a convex open $\wb{t}$-neighborhood
$\wt{J}_0\sub \wt{J}$ and open neighborhoods
$\wt{U}_0\sub\wt{U}$ as well as $\wt{P}_0\sub\wt{P}$
of~$\wb{y}$ and $\wb{p}$, respectively, such that the initial value problem
\[
y'(t)=\wt{f}(t,y(t),p),\quad y(t_0)=y_0
\]
has a complex analytic solution $\wt{\gamma}_{t_0,y_0,p}\colon \wt{J}_0\to E_\C$
for all $(t_0,y_0,p)$ in\linebreak
$\wt{J}_0\times\wt{U}_0\times\wt{P}_0$
and the map
\[
\wt{\Phi}\colon \wt{J}_0\times\wt{J}_0\times \wt{U_0}\times\wt{P}\to E_\C,\;\,
(t,t_0,y_0,p)\mto\wt{\gamma}_{t_0,y_0,p}(t)
\]
is complex analytic.
Define $J_0:=\wt{J}_0\cap\R$, $U_0:=\wt{U}_0\cap E$ and $P_0:=\wt{P}_0\cap F$.
For all $(t_0,y_0,p)\in J_0\times U_0\times P_0$, Lemma~\ref{ode-in-subspace}
shows that
$\wt{\gamma}_{t_0,y_0,p}(J_0)\sub E$
and that $\gamma_{t_0,y_0,p}:=\wt{\gamma}_{t_0,y_0,p}|_{J_0}$
is a solution to~(\ref{ini-over-K}).
Moreover, the map~$\Phi$ defined in~(b) is real analytic,
as $\wt{\Phi}$ is a complex analytic extension.
\end{prf}
\begin{rem}
The proof of Theorem~\ref{local-ivp-dep-Ck} provides quantitative information:
If $U=B^E_{2R}(\wb{y})$
and $J$, $P$ are chosen such that
(\ref{expli-1-ode}) and
(\ref{expli-2-ode}) are satisfied,
then $U_0:=B^E_R(\wb{y})$, $J_0:=J$
and $P_0:=P$ can be chosen
if $\K=\C$ or $\K=\R$ and $k\in\N_0\cup\{\infty\}$.\medskip

If, moreover, $f$ is real analytic, then $\Phi$ (as in Theorem~\ref{local-ivp-dep-Ck}(b))
is real analytic on all of
$J\times J\times B^E_R(\wb{y})\times P$. However, the latter
will only become clear later (once Theorem~\ref{nice-spec-para}
has been established);
in the real analytic case, our proof of Theorem~\ref{local-ivp-dep-Ck}
is \emph{not} quantitative and may involve smaller neighborhoods.
\end{rem}
\begin{rem}
If $k=0$, then $P$ can be replaced with an arbitrary topological space
in Theorem~\ref{local-ivp-dep-Ck}.
Likewise in Theorem~\ref{nice-spec-para},
Proposition~\ref{special-ini-par},
Definition~\ref{defn-ode-par-m},
Definition~\ref{defn-loc-flo-m},
Theorem~\ref{ode-loc-glob-par},
and Corollary~\ref{pardep-ode-ban-m}.
\end{rem}
We also have a version of
Theorem~\ref{local-ivp-dep-Ck}
in the setting of~$FC^k$-maps.
\begin{cor}\label{paraFCk}
Let $(E,\|\cdot\|_E)$ be a real Banach space, $(F,\|\cdot\|_F)$
be a normed space over~$\R$,
and $(\wb{t},\wb{y},\wb{p})\in \R\times E\times F$.
Let
$U\sub E$ be an open $\wb{y}$-neighborhood,
$J\sub\R$ be a non-degenerate interval containing~$\wb{t}$
and $P\sub F$ be a locally convex subset with dense interior
such that $\wb{p}\in P$.
Let $k\in \N \cup\{\infty\}$ and
$f\colon J\times U\times P\to E$
be an $FC^k$-map.
Then there exist a relatively open subinterval $J_0\sub J$ containing~$\wb{t}$,
an open $\wb{y}$-neighborhood $U_0\sub U$, and
a relatively open $\wb{p}$-neighborhood $P_0\sub P$,
with the following properties:
\begin{description}[(D)]
\item[\rm(a)]
For all $(t_0,y_0,p)\in J_0\times U_0\times P_0$, the initial value problem
\begin{equation}\label{ini-over-K2}
y'(t)=f(t,y(t),p),\quad y(t_0)=y_0
\end{equation}
has a unique solution $\gamma_{t_0,y_0,p}\colon J_0\to E$;\vspace{.8mm}
\item[\rm(b)]
The map
$\Phi\colon J_0\times J_0\times U_0\times P_0\to E$, $(t,t_0,y_0,p)\mto\gamma_{t_0,y_0,p}(t)$
is~$FC^k$.
\end{description}
\end{cor}
\begin{prf}
Since $f$ is $FC^k$, also the map $h$ in the proof of Theorem~\ref{local-ivp-dep-Ck} is~$FC^k$.
In view of (\ref{thusbFck}), we deduce from
Lemma~\ref{fstaFCk} and Remark~\ref{fstaFCk2}
(or the alternative discussion in Exercise~\ref{exc-ifFCk})
that the map~$g$ in the proof of Theorem~\ref{local-ivp-dep-Ck} is~$FC^k$.
Using Proposition~\ref{fp-Ck-dep}(b), we see that the fixed point
$\eta_{t,t_0,y_0,p}$ of the contraction $g_{t,t_0,y_0,p}$
is an $FC^k$-function of $(t,t_0,y_0,p)$,
whence also $\Phi(t,t_0,y_0,p)=\eta_{t,t_0,y_0,p}(1)+y_0$ is $FC^k$ in $(t,t_0,y_0,p)$.
\end{prf}
We now turn to an analog of Theorem~\ref{local-ivp-dep-Ck}
for the case when the right-hand side of a differential equation with parameters
is a $C^{0,k}$-function. We recommend to skip this more specialized result
on a first reading, which will only be used twice
(for the proofs of $C^0$-regularity for Banach-Lie groups
and diffeomorphism groups).
\begin{prop}\label{special-ini-par}
Let $(E,\|\cdot\|)$ be a Banach space, $F$ be a locally convex space, $P\sub F$ be a locally convex subset with dense interior,
$a<b$ be real numbers, $\wb{y}\in E$,
$R>0$, $k\in\N_0\cup\{\infty\}$
and
\[
f\colon [a,b]\times (B^E_{2R}(\wb{y})\times P)\to E,\;\, (t,y,p)\mto f(t,y,p)
\]
be a $C^{0,k}$-function satisfying a global Lipschitz condition in the~$y$-variable. Assume that
\[
M:=\|f\|_\infty:=\sup\{\|f(t,y,p)\|\colon (t,y,p)\in [a,b]\times B^E_{2R}(\wb{y})\times P\}<\infty
\]
and $R/M\geq b-a$.
Then, for each $(t_0,y_0,p)\in [a,b]\times B^E_R(\wb{y})\times P$,
the initial value problem
\[
y'(t)=f(t,y(t),p),\quad y(t_0)=y_0
\]
has a unique solution $\gamma_{t_0,y_0,p}\colon [a,b]\to E$, and the
following map is $C^{0,k}$:
\[
\Phi\colon ([a,b]\times [a,b])\times (B^E_R(\wb{y})\times P)\to E,\quad
\Phi(t,t_0,y_0,p):=\gamma_{t_0,y_0,p}(t).
\]
\end{prop}
\begin{prf}
Abbreviate $F:=C([a,b],E)$.
By hypothesis,
\[
L:=\sup\{\Lip(f(t,\cdot,p)\colon (t,p)\in [a,b]\times P\}\in[0,\infty[.
\]
Fix $n\in\N$ such that
\[
\theta:=\frac{L^n}{n!}(b-a)^n<1.
\]
The existence of~$\gamma_{t_0,y_0,p}$ is guaranteed by
Quantitative Existence Theorem (Theorem~\ref{quant-PL}),
and the proof of the latter shows that
\[
\eta_{t_0,y_0,p}:=\gamma_{t_0,y_0,p}-y_0
\]
is the unique fixed point of the contraction $(g_{t_0,y_0,p})^n$ of $\wb{B}^F_R(0)$
with contraction constant~$\theta$,
where
\begin{equation}\label{monsta}
g\colon [a,b]\times (B^E_R(\wb{y})\times P\times \wb{B}^F_R(0))\to \wb{B}^F_R(0)
\end{equation}
is given by $g(t_0,y_0,p,\psi)(t):=\int_{t_0}^tf(s,\psi(s)+y_0,p)\,ds$ for
$(t_0,y_0,p,\psi)$ in the left-hand side of~(\ref{monsta}) and $t\in[a,b]$
(and $g_{t_0,y_0,p}(\psi):=g(t_0,y_0,p,\psi)$).
Thus
\begin{equation}\label{aufgedroeselt}
g(t_0,y_0,p,\psi)
=I\big(t_0,(f^p)_*(\psi+y_0)\big)
\end{equation}
using the map
\[
I\colon  [a,b]\times C([a,b],E)\to C([a,b],E)
\]
given by $I(t_0,\psi)(t):=\int_{t_0}^t\psi(s)\,ds$
for $(t_0,\psi,t)\in [a,b]\times C([a,b],E)\times [a,b]$
and the map
\[
(f^p)_*\colon C([a,b],B^E_{2R}(\wb{y}))\to C([a,b],E)
\]
given by $(f^p)_*(\psi)(t):=f(t,\psi(t),p)$,
where $C([a,b],B^E_{2R}(\wb{y}))=B^F_{2R}(\wb{y})$
if we identify~$\wb{y}$ with the constant function $[a,b]\to E$,
$t\mto \wb{y}$.
As the map $f\colon [a,b]\times (B^E_{2R}(\wb{y})\times P)\to E$
is $C^{0,k}$, Proposition~\ref{pfctisCkpar} shows that
\[
P\times B^F_{2R}(\wb{y})\to C([a,b],E),\quad (p,\psi)\mto (f^p)_*(\psi)
\]
is a $C^k$-map.
Moreover, $I$ is continuous as
\begin{eqnarray*}
\lefteqn{\|I(t_1,\psi_1)-I(t_0,\psi_0)\|_\infty}\qquad\qquad\\
&=&
\|I(t_1,\psi_1-\psi_0)+I(t_1,\psi_0)-I(t_0,\psi_0)\|_\infty\\
&=& \left\|t\mto \int_{t_1}^t (\psi_1(s)-\psi_0(s))\,ds +\int_{t_1}^{t_0}\psi_0(s)\,ds\right\|_\infty\\
&\leq&  (b-a)\|\psi_1-\psi_0\|_\infty+|t_1-t_0|\, \|\psi_0\|_\infty\to 0
\end{eqnarray*}
as $(t_1,\psi_1)\to (t_0,\psi_0)$ in $[a,b]\times C([a,b],E)$.
Since~$I$ is linear in the second argument, we deduce from Lemma~\ref{linCkell} that~$I$ is $C^{0,\infty}$
(and hence $C^{0,k}$). Applying Proposition~\ref{chainR1} to~(\ref{aufgedroeselt}),
we infer that~$g$ is~$C^{0,k}$.
Hence
\[
h\colon
[a,b]\times (B^E_R(\wb{y})\times P\times \wb{B}^F_R(0))\to \wb{B}^F_R(0),\;\,
h(t_0,y_0,p,\psi):=(g_{t_0,y_0,p})^n(\psi)
\]
is~$C^{0,k}$, by repeated use of Lemma~\ref{C0chainy}.
As a consequence, the mapping
\[
[a,b]\times(\wb{B}^E_R(\wb{y})\times P)\to C([a,b],E),\quad (t_0,y_0,p)\mto\eta_{t_0,y_0,p}
\]
assigning to $(t_0,y_0,p)$ the fixed point of the contraction $h_{t_0,y_0,p}=(g_{t_0,y_0,p})^n$
is~$C^{0,k}$ (see Lemma~\ref{fp-C0k-dep}).
Then also
\[
[a,b]\times(\wb{B}^E_R(\wb{y})\times P)\to C([a,b],E),\quad (t_0,y_0,p)\mto\gamma_{t_0,y_0,p}
=\eta_{t_0,y_0,p}+y_0
\]
is $C^{0,k}$ and thus
$\Phi$ is~$C^{0,k}$, by Lemma~\ref{explaw-vari}. 
\end{prf}
\begin{small}
\subsection*{Exercises for Section~\ref{sec-ode}}

\begin{exer}\label{exc-ban-weis}
Let $(X,d)$ be a complete metric space with $X\not=\emptyset$
and $f\colon X\to X$ be a Lipschitz map whose iterates satisfy
\[
\sum_{n=1}^\infty\Lip(f^n)<\infty.
\]
Prove the Fixed Point Theorem of Banach--Weissinger:
\begin{description}[(D)]
\item[(a)]
$d(f^{n+m}(x),f^n(x))\leq d(f(x),x)\sum_{k=n}^{n+m-1}\Lip(f^k)$
for all $n,m\in\N$ and $x\in X$;
\item[(b)]
For each $x\in X$, the sequence $(f^n(x))_{n\in\N}$
converges and its limit is a fixed point of~$f$;
\item[(c)]
The fixed point $x_\infty$ of~$f$ is unique;
\item[(d)]
$d(x_\infty,f^n(x))\leq d(f(x),x)\sum_{k=n}^\infty \Lip(f^k)$, for each $x\in X$ and $n\in\N$.
\item[(e)]
Show that the preceding fixed point theorem applies to the self-map~$g$
in the proof of Theorem~\ref{quant-PL} (not only to a power~$g^n$).
\end{description}
\end{exer}

\begin{exer}\label{exc-ode-rea}
Let $E$ be a locally convex space and $f\colon W\to E$ be a real analytic
mapping on an open subset $W\sub \R\times E$.
\begin{description}[(D)]
\item[(a)]
Let $\gamma$ be a real analytic solution of
\begin{equation}\label{ana-equ}
y'(t)=f(t,y(t)), \quad y(t_0)=y_0.
\end{equation}
Then
\[
\gamma''(t)=\frac{\partial f}{\partial t}(t,\gamma(t))+d_2f(t,\gamma(t),\gamma'(t))
=f_2(t,\gamma(t))
\]
with $f_2(t,y):=\frac{\partial f}{\partial t}(t,y)+d_2f(t,y,f(t,y))$.
Proceeding in this way, show that there are real analytic functions $f_n\colon W\to E$
for all integers $n\geq 2$ such that each real analytic solution $\gamma\colon I\to E$ of~(\ref{ana-equ})
satisfies
\[
\gamma^{(n)}(t)=f_n(t,\gamma(t))\quad\mbox{for all $n\in\N$ and $t\in I$.}
\]
\item[(b)] Deduce that if $\gamma_j\colon I_j\to E$ are real analytic solutions of~(\ref{ana-equ})
for $j\in\{1,2\}$,
then $\gamma_1$ and~$\gamma_2$ agree on some neighborhood of~$t_0$.
Infer that $\gamma_1|_{I_1\cap I_2}=\gamma_2|_{I_1\cap I_2}$.
\item[(c)]
Show that if (\ref{ana-equ}) has a real analytic solution,
then it also has a maximal real analytic solution $\gamma^{\omega}_{t_0,y_0}\colon I_{t_0,y_0}^\omega
\to E$
(such that every real analytic soluton is a restriction of the former).
\item[(d)]
Assuming that~(\ref{ana-equ}) has a real analytic solution for all $(t_0,y_0)\in W$,
define $\Omega^\omega:=\bigcup_{(t_0,y_0)\in W}I^\omega_{t_0,y_0}\times\{(t_0,y_0)\}$
and $\Fl^\omega\colon \Omega^\omega\to E$, $\Fl^\omega(t,t_0,y_0):=\Fl^\omega_{t,t_0}(y_0):=
\gamma^\omega_{t_0,y_0}(t)$.
Show that if $(t_0,y_0)\in W$ is given, $t_1\in I_{t_0,y_0}$ and $t_2\in I_{t_1,\Fl^\omega(t_1,t_0,y_0)}$,
then $t_2\in I_{t_0,y_0}$ and
\begin{equation}
\Fl_{t_2,t_0}^\omega(y_0)=\Fl_{t_2,t_1}^\omega(\Fl_{t_1,t_0}^\omega(y_0)).
\end{equation}
\end{description}
\end{exer}

\begin{exer}\label{exc-intlip}
Let $J\sub \R$ be a non-degenerate interval,
$(E,\|\cdot\|)$ be a normed space, $U\sub E$ be a subset
and $f\colon J\times U \to E$ be a function.
Assume that there exists a locally Lebesgue integrable function
$g\colon J\to [0,\infty]$ such that
\[
\Lip f(t,\cdot)\leq g(t)\quad\mbox{for all $\,t\in J$.}
\]
Show that $y'(t)=f(t,y(t))$ satisfies local uniqueness of solutions.
\end{exer}

\begin{exer}\label{exc-generato}
Let $E:=\ell^2$ be the Hilbert space of square summable complex sequences.
For $t\in\R$, define $A(t)\in \cL(E)$ with $\|A(t)\|_{\op}=1$
via
\[
A(t)(x):=(e^{int}x_n)_{n\in\N}\quad\mbox{for $x=(x_n)_{n\in\N}\in E$.}
\]
Show that the map
$\wh{A}\colon \R\times E\to E$, $(t,x)\mto A(t)(x)$ is continuous.
Show that the map $A\colon \R\to \cL(E)_b$, $t\mto A(t)$ is not continuous.
\end{exer}

\begin{exer}\label{exc-C1-nec}
Let $U\sub\C$ be an open subset, $E$ be a complex Banach space
and $f\colon U\to E$ be a continuous function.
Show that if the differential equation $y'(z)=f(z)$ admits a $C^1_\C$-solution
$\gamma\colon J\to E$ on a non-empty open connected subset $J\sub U$,
then $f|_J$ is complex analytic and hence~$C^1_\C$.\\[2mm]
The condition $k\geq 1$ therefore cannot be omitted if $\K=\C$
in Theorem~\ref{local-ivp-dep-Ck}.
\end{exer}

\begin{exer}\label{exc-lie-der-loc}
Let $J$ be a non-degenerate interval, $E$ be a locally convex space,
$U\sub E$ be a locally convex subset with dense interior,
$W\sub J\times U$ be an open subset
and $f\colon W\to E$ be a function such that
$y'(t)=f(t,y(t))$ admits local $C^1$-flows.
Let $\Fl\colon \Omega\to U$, $(t,t_0,y_0)\mto\Fl(t,t_0,y_0)=\Fl_{t,t_0}(y_0)$
be the maximal flow; write
\[
\partial_1\!\Fl(t,t_0,y_0):=\frac{\partial}{\partial t}\Fl(t,t_0,y_0):=
d\Fl(t,t_0,y_0,1,0,0)
\]
and $\partial_2\!\Fl(t,t_0,y_0):=\frac{\partial}{\partial t_0}\Fl(t,t_0,y_0):=
d\Fl(t,t_0,y_0,0,1,0)$ for $(t,t_0,y_0)\in \Omega$.
Since $I_{t_0,y_0}\to E$, $t\mto\Fl_{t,t_0}(y_0)$ solves
$y'(t)=f(t,y(t))$, $y(t_0)=y_0$, we have
\begin{equation}\label{flo-basic-2}
\Fl(t_0,t_0,y_0)=y_0\;\;\mbox{for all $(t_0,y_0)\in W$}
\end{equation}
and
\begin{equation}\label{start-flo-basic}
\partial_1\!\Fl(t,t_0,y_0)=f(t,\Fl(t,t_0,y_0))\;\;\mbox{for all $(t,t_0,y_0)\in\Omega$.}
\end{equation}
Notably,
$\partial_1\!\Fl(t_0,t_0,y_0)=f(t_0,y_0)$ for all $(t_0,y_0)\in W$.
\begin{description}[(a)]
\item[(a)]
Differentiating (\ref{flo-basic-2}) with respect to $t_0$, show that
\[
\partial_2\!\Fl(t_0,t_0,y_0)=-\partial_1\!\Fl(t_0,t_0,y_0)=-f(t_0,y_0)
\]
for all $(t_0,y_0)\in W$. Also deduce that
\begin{equation}\label{basic-flo-3}
d\Fl_{t_0,t_0}(y_0,z)=z\;\;\mbox{for all $\,(t_0,y_0)\in W$ and $z\in E$.}
\end{equation}
\item[(b)]
If $g\colon U\to F$ is a $C^1$-function to a locally convex space~$F$,
show that
\[
\frac{d}{d t}\Big|_{t=t_0}g\big(\Fl_{t,t_0}(y_0)\big)=dg\big(y_0,f_{t_0}(y_0)\big)
\]
for all $(t_0,y_0)\in W$, where $f_{t_0}(y_0):=f(t_0,y_0)$.
\item[(c)]
If $f$ admits local $C^2$-flows, deduce from (\ref{basic-flo-3})
that $d^{\,(2)}\Fl_{t_0,t_0}(y_0,z,w)=0$
for all $(t_0,y_0)\in W$ and $z,w\in E$.
\item[(d)]
If $f$ admits local $C^2$-flows and $g\colon U\to E$ is a $C^1$-map, show that
\[
\frac{d}{dt}\Big|_{t=t_0}d\Fl_{t_0,t}\big(\Fl_{t,t_0}(y_0),g(\Fl_{t,t_0}(y_0))\big)
=dg\big(y_0,f_{t_0}(y_0)\big)-df_{t_0}\big(y_0,g(y_0)\big)
\]
for all $(t_0,y_0)\in W$.\\[1mm]
[Using the Schwarz Theorem, we get
$d^{\,(2)}\Fl((t_0,t_0,y_0),(0,0,g(y_0)),(0,1,0))
=d^{\,(2)}\Fl((t_0,t_0,y_0),(0,1,0),(0,0,g(y_0)))
=d\,\partial_2\!\Fl((t_0,t_0,y_0),(0,0,g(y_0)))$,
which equals $-df_{t_0}(y_0,g(y_0))$.]
\end{description}
\end{exer}

\begin{exer}
Show that the flow $\Fl\colon \Omega\to M$ in Theorem~\ref{nice-spec-para}
is not only $C^k$, but that also
$\frac{\partial}{\partial t} \Fl\colon \Omega\to E$ is~$C^k$.\\[1mm]
[Note that $\frac{\partial}{\partial t} \Fl
(t,t_0,y_0,p)=f(t,\Fl(t,t_0,y_0,p),p)$.]
\end{exer}

\begin{exer}\label{exc-par-also}
Let $J\sub\R$ be a non-degenerate interval, $E$ and $F$ be locally convex spaces,
$U\sub E$ and $P\sub F$ be subsets,
$f\colon J\times U\times P\to E$ be a function and $(t_0,y_0,p_0)\in J\times U\times P$.
Show that a $C^1$-function $\gamma\colon J\to E$ is a solution to
\[
y'(t)=f(t,y(t),p_0),\quad y(t_0)=y_0
\]
if and only if $\eta\colon J\to E\times F$,
$\eta(t):=(\gamma(t),p_0)$ is a solution to
\[
w'(t)=g(t,w(t)), \quad w(t_0)=(y_0,p_0)
\]
with $g\colon J\times U\times P\to E\times F$,
$g(t,y,p):=(f(t,y,p),0)$.
\end{exer}
\begin{exer}\label{exc-ders-flow}
Let $E$ be a locally convex space, $U\sub E$ be a locally convex subset with
dense interior, $J\sub\R$ be a non-degenerate interval,
$W\sub J\times U$ be an open subset and $f\colon W\to E$ be
a $C^1$-map such that $y'(t)=f(t,y)$ admits local $C^1$-flows.
Let $\Fl\colon \Omega\to E$ be the flow;
thus $\frac{d}{dt}\Fl_{t,t_0}(x)=f(t,\Fl_{t,t_0}(x))$
for all $(t,t_0,x)\in\Omega$
and $\Fl_{t_0,t_0}(x)=x$ for all $(t_0,x)\in W$.
\begin{description}[(D)]
\item[(a)]
Verify that $\frac{d}{dt}\Fl_{t,t_0}(x)=f(t,\Fl_{t,t_0}(x))$
is $C^1$ and hence $C^{0,1}$ in $(t,t_0,x)\in \Omega \sub J\times
(J\times U)$. Conclude that $\Fl$ is $C^{1,1}$
on $\Omega\sub J\times (J\times U)$.
\item[(b)]
By~(a), the derivative $\frac{d}{dt}(d\Fl_{t,t_0}(x,y))$
exists for all $(t,t_0,x)\in\Omega$ and $y\in E$. Using Proposition~\ref{schwarzCkell},
show that
\[
\frac{d}{dt}d\Fl_{t,t_0}(x,y)=d_2f(t,\Fl_{t,t_0}(x),d\Fl_{t,t_0}(x,y)).
\]
Notably, $\frac{d}{dt}\big|_{t=t_0}\Fl_{t,t_0}(x,y)=d_2f(t_0,x,y)$
for all $(t_0,x)\in W$ and $y\in E$.
\end{description}
\end{exer}

Let us discuss Euler's method for approximate solutions to
initial value problems, and Peano's existence theorem.
In both Exercise~\ref{exc-euler} and Exercise~\ref{exc-peano},
we let $a<b$ be real numbers, $(E,\|\cdot\|_E)$ be a Banach space, $y_0\in E$
and
\[
f\colon [a,b]\times \wb{B}^E_R(y_0)\to E
\]
be a continuous function such that
\[
M:=\sup\{\|f(t,y)\|_\infty\colon (t,y)\in [a,b]\times \wb{B}^E_R(y_0)\}<\infty\;\;
\mbox{and}\;\;
(b-a)M\leq R.
\]
Given a subdivision $Z=(t_0,\ldots, t_\ell)$ of $[a,b]$ with $a=t_0<t_1<\cdots<t_\ell=b$,
we write $m(Z)$ for the maximum of $t_j-t_{j-1}$ for $j\in\{1,\ldots,\ell\}$. Recursively,
we define $y_j\in E$ with $\|y_j-y_0\|_E\leq (t_j-a)M\leq (b-a)M\leq R$ via
\[
y_j:=y_{j-1}+(t_j-t_{j-1})f(t_{j-1},y_{j-1})\;\,\mbox{for $j\in\{1,\ldots, \ell\}$.}
\]
The Euler polygon corresponding to~$Z$ is the piecewise affine-linear function
$\gamma_Z\colon [a,b]\to \wb{B}^E_R(y_0)$ which maps
$t\in [t_{j-1},t_j]$ to $\gamma_Z(t):=y_{j-1}+\frac{t-t_{j-1}}{t_j-t_{j-1}}
(y_j-y_{j-1})$, for $j\in\{1,\ldots, \ell\}$.
Thus, if we define $\tau_Z\colon [a,b]\to Z$ via
$\tau_Z(t):=t_{j-1}$ if $t\in [t_{j-1},t_j[$ and $j\in\{1,\ldots,\ell-1\}$
or $j=\ell$ and $t\in [t_{\ell-1}, b]$, then
\begin{equation}\label{intfo-eul}
\gamma_Z(t)=y_0+\int_a^t f(\tau_Z(s),\gamma_Z(\tau_Z(s)))\, ds
\end{equation}
for all $t\in [a,b]$, where $\gamma_Z(\tau_Z(s))=y_{j-1}$ if $\tau_Z(s)=t_{j-1}$.
Note that
\begin{eqnarray}
\Lip(\gamma_Z)&\leq & M\quad\mbox{and}\label{bas-eul1}\\
|\tau_Z(s)-s|&\leq & m(Z)\quad\mbox{for all $s\in [a,b]$.}\label{bas-eul2}
\end{eqnarray}
We let $(Z_k)_{k\in\N}$ be a sequence of subdivisions of~$[a,b]$
such that $m(Z_k)\to 0$ as $k\to\infty$, and abbreviate $\gamma_k:=\gamma_{Z_k}$
and $\tau_k:=\tau_{Z_k}$.
Let $F:=C([a,b],E)$ with the supremum norm $\|\cdot\|_\infty$ and consider the map
$g\colon \wb{B}^F_R(y_0)\to F$ determined by
\begin{equation}\label{the-g-io}
g(\eta)(t):=y_0+\int_a^tf(s,\eta(s))\,ds\;\;\mbox{for $\,\eta\in\wb{B}^F_R(y_0)\,$ and $\,t\in[a,b]$,}
\end{equation}
writing $y_0$ for $[a,b]\to E$, $t\mto y_0$.
Since $(b-a)M\leq R$, we have $g(\eta)\in \wb{B}^F_R(y_0)$ for each $\eta\in \wb{B}^F_R(y_0)$.
Give $\R\times E$ a norm via $\|(t,y)\|:=\max\{|t|,\|y\|_E)$.
\begin{exer}[Euler Method]\label{exc-euler}
\begin{description}[(D)]
\item[(a)]
If $E$ is finite-dimensional, then $[a,b]\times \wb{B}^E_R(y_0)$ is compact, whence
$f$ is uniformly continuous. Given $\ve>0$, pick $\delta>0$ such that
$\|f(s_2,v_2)-f(s_1,v_2)\|_E\leq\ve$ for all $(s_1,y_1)$, $(s_2,y_2)\in [a,b]\times\wb{B}^E_R(y_0)$
with $|s_2-s_1|\leq\delta$ and $\|v_2-v_1\|_E\leq\delta$.~Using
\[
g(\gamma_k)(t)-\gamma_k(t)=\int_a^t\Big(f\big(s,\gamma_k(s)\big)-f\big(\tau_k(s),\gamma_k(\tau_k(s))\big)\Big)\, ds,
\]
show that $\|g(\gamma_k)-\gamma_k\|_\infty\leq \ve(b-a)$
if~$k$ is so large that $\max\{1,M\}m(Z_k)\leq\delta$.
\end{description}
Now assume that $f$ as before satisfies a global Lipschitz condition in its second argument, i.e.,
$L:=\sup\{\Lip(f(t,\cdot))\colon t\in[a,b]\}<\infty$. Choose $n\in \N$ so large that
$\theta:=\frac{L^n}{n!}(b-a)^n<1$; then $g^n\colon \wb{B}^F_R(y_0)\to \wb{B}^F_R(y_0)$
is a contraction with $\Lip(g^n)\leq\theta$; moreover, $\Lip(g)\leq (b-a)L$ (cf.\ proof of
Theorem~\ref{quant-PL}).
Let $\gamma\colon [a,b]\to \wb{B}^E_R(y_0)$
be the solution to $y'(t)=f(t,y(t))$, $y(t_0)=y_0$, which exists by Theorem~\ref{quant-PL}.
\begin{description}[(D)]
\item[(b)]
Using the telescopic sum $g^n(\gamma_k)-\gamma_k=\sum_{j=1}^n(g^j(\gamma_k)-g^{j-1}(\gamma_k))$,
deduce that $\|g^n(\gamma_k)-\gamma_k\|_\infty\leq C\ve$ for~$k$
as in~(a),
with $C:=(b-a)\sum_{j=1}^n(\Lip(g))^j$. Employing the a priori estimate~(\ref{a-prio-est})
from Lemma~\ref{ban-fix},
deduce that\linebreak
$\|\gamma-\gamma_k\|_\infty\leq \frac{C}{1-\theta}\,\ve$.
Thus $\gamma_k\to\gamma$ uniformly as $k\to\infty$.
\item[(c)]
Show that all conclusions of~(a) and~(b) remain valid for an arbitrary Banach space~$E$,
if we assume as an extra hypothesis that~$f$ is uniformly continuous.
The latter holds, e.g., if $f$ is Lipschitz continuous; show that
$\|g(\gamma_k)-\gamma_k\|_\infty\leq \Lip(f)\max\{1,M\}m(Z_k)$ and
\[
\|\gamma-\gamma_k\|_\infty\leq \frac{C\Lip(f)\max\{1,M\}}{1-\theta} m(Z_k)\;\;
\mbox{for all $k\in\N$}
\]
in this case, with~$C$ as in~(b).
\end{description}
\end{exer}

\begin{exer}[Peano Existence Theorem]\label{exc-peano}
If $E$ is finite-dimensional, then $\wb{B}^E_R(y_0)$ is compact.
Deduce from $\Lip(\gamma_k)\leq M$ for all $k\in \N$ that $\{\gamma_k\colon k\in \N\}$
is equicontinuous. Since $\gamma_k(t)\in\wb{B}^E_R(y_0)$ for each $t\in [a,b]$
and $k\in\N$, Ascoli's Theorem shows that $\{\gamma_k\colon k\in\N\}$
is relatively compact in $(C([a,b],E),\|\cdot\|_\infty)$.
Hence, after passing to a subsequence, we may assume that $\gamma_k\to \gamma$ 
uniformly for some $\gamma\in C([a,b],E)$ (and then $\gamma\in \wb{B}^F_R(y_0)$).
Write~$g$ from (\ref{the-g-io}) as
\[
g(\eta)=y_0+I(f_*(\eta)),
\]
where $f_*\colon C([a,b],\wb{B}^E_R(y_0))\to C([a,b],E)$, $\eta\mto f\circ (\id_{[a,b]},\eta)$
is continuous by Lemma~\ref{ctspfapp} and $I\colon C([a,b],E)\to C([a,b],E)$
is the continuous linear map given by $I(\eta)(t):=\int_a^t\eta(s)\, ds$.
We know from Exercise~\ref{exc-euler}(a) that
$\Delta_k:=g(\gamma_k)-\gamma_k\to 0$. Letting $k\to\infty$ in
\[
\gamma_k=g(\gamma_k)-\Delta_k=y_0 + I(f_*(\gamma_k))-\Delta_k,
\]
deduce that $\gamma=g(\gamma)$. Thus $\gamma$ solves
$y'(t)=f(t,y(t))$, $y(a)=y_0$.
\end{exer}

\end{small}
\section{Differential equations on manifolds}\label{sec-ode-mfd}
In Section~\ref{sec-ode},
we already considered
differential equations on subsets of locally convex spaces,
with an emphasis on local results.
We now turn to global aspects, which are more adequately
studied in the framework of differential equations on manifolds.
On the reader's part, this section only presupposes
knowledge of the most basic definitions concerning manifolds:
The notion of a 
$C^k$-manifold modeled on a locally convex space
(possibly with boundary);\footnote{All results in this section apply if
``boundary'' is interpreted in the most general form of a ``rough boundary''.}
the notion of a $C^k$-map
$f\colon M\to N$ between such manifolds; and the tangent map
$Tf\colon TM\to TN$. All manifolds in this section are manifolds over the
ground field $\K=\R$.
\begin{defn}
If $M$ is a $C^1$-manifold (possibly with boundary), $I\sub \R$ a non-degenerate interval and
$\gamma\colon I\to M$ a $C^1$-map, we define
$\dot{\gamma}\colon I\to TM$ via
\[
\dot{\gamma}(t):=T\gamma(t,1)\in T_{\gamma(t)}M\quad\mbox{for $t\in I$,}
\]
identifying $TI$ with $I\times\R$ as usual. Thus, if $M=U$ is an open subset of a locally convex space~$E$,
then
\[
\dot{\gamma}=(\gamma,\gamma')\colon I\to TU=U\times E.
\]
It is useful to use different notation for $\dot{\gamma}$ and $\gamma'$
in the current section; in later parts of the book, we shall often write
$\gamma'$ in place of $\dot{\gamma}$, and the intended meaning
will be clear from the context.
\end{defn}
\begin{defn}\label{ode-on-m}
Let $M$ be a $C^1$-manifold modeled on a locally convex space~$E$.
Let $W\sub \R\times M$ be a subset and
$f\colon W\to TM$ be a function
such that $f(t,y)\in T_yM$ for all $(t,y)\in W$.\medskip

\noindent
(a) We say that a function $\gamma \colon I\to M$
on a non-degenerate interval $I\sub \R$
is a \emph{solution to the differential equation}
\begin{equation}\label{theode-m}
\dot{y}(t)=f(t,y(t))
\end{equation}
if $\gamma$ is a $C^1$-map,
$(t,\gamma(t))\in W$ for all $t\in I$ and $\dot{\gamma}(t)=f(t,\gamma(t))$.\medskip

\noindent
(b) If $(t_0,y_0)\in W$ and $\gamma\colon I\to M$ is a solution to~(\ref{theode-m})
such that $t_0\in I$ and $\gamma(t_0)=y_0$,
then $\gamma$ is called a \emph{solution to the initial value problem}
\begin{equation}\label{theiniprob-m}
\left\{
\begin{array}{ccl}
\dot{y}(t) &=& f(t,y(t))\\
y(t_0)&=&y_0.
\end{array}
\right.
\end{equation}
(c) We say that the differential equation (\ref{theode-m})
satisfies \emph{local uniqueness of solutions} if the following condition holds:
For all solutions $\gamma_1\colon I_1\to M$ and $\gamma_2\colon I_2\to M$ to (\ref{theode-m})
such that $\gamma_1(t_0)=\gamma_2(t_0)$ for some $t_0\in I_1\cap I_2$,
there exists a neighborhood~$I$ of~$t_0$ in $I_1\cap I_2$ such that
$\gamma_1|_I=\gamma_2|_I$.\medskip

\noindent
(d) If $W$ is a subset of $J\times M$ for a non-degenerate interval
$J\sub\R$, we say that the differential equation (\ref{theode-m}) satisfies \emph{local existence
of solutions} (with respect to~$J$) if the following condition holds:
For all $(t_0,y_0)\in W$, there exists a solution $\gamma\colon I\to M$ to the initial value
problem (\ref{theiniprob-m}) whose domain~$I$ is relatively open in~$J$.
\end{defn}
\begin{rem}\label{rem-ode-m}
(a) In most applications, we have $W=J\times M$ for a non-degenerate interval
$J\sub\R$. If~$M$ is a $C^{k+1}$-manifold
and $f\colon J\times M\to TM$ as above is $C^k$ for some $k\in\N_0$,
then we can consider $f$ as a \emph{time-dependent vector field of class~$C^k$}
on~$M$, in the sense that~$f$ is~$C^k$ and $f_t:=f(t,\cdot)\colon M\to TM$
is a $C^k$-vector field on~$M$ for each $t\in J$.\medskip

(b) If $f\colon W\to TM$ is as in Definition~\ref{ode-on-m}
and $\phi\colon U_\phi\to V_\phi\sub E$ is a chart for~$M$,
define $W_\phi:=\{(t,\phi(y))\colon (t,y)\in W\cap (\R\times U_\phi)\}\sub\R\times V_\phi$
and
\[
f_\phi\colon W_\phi\to E,\quad f_\phi(t,y):=d\phi f(t,\phi^{-1}(y)).
\]
If $\gamma\colon I\to M$ is a $C^1$-map on a non-degenerate interval $I\sub\R$
such that $\gamma(I)\sub U_\phi$, then $\gamma$ solves~(\ref{theode-m})
if and only if $\phi\circ\gamma$ solves $y'(t)=f_\phi(t,y(t))$.
\end{rem}
Replacing~$E$ with~$M$ in the proofs of Lemmas~\ref{lem-loc-glob}
and~Remark~\ref{glue-sol-ode}, 
we obtain:
\begin{lem}\label{lem-loc-glob-m}
Consider a differential equation \emph{(\ref{theode-m})}
with local uniqueness of solutions.
If $\gamma_1\colon I_1\to M$ and $\gamma_2\colon I_2\to M$
are solutions to \emph{(\ref{theode-m})}
such that $\gamma_1(t_0)=\gamma_2(t_0)$
for some $t_0\in I_1\cap I_2$, then $\gamma_1|_{I_1\cap I_2}
=\gamma_2|_{I_1\cap I_2}$.\qed
\end{lem}
\begin{lem}\label{ex-maxi-sol-m}
Consider a differential equation
$\dot{y}(t)=f(t,y(t))$
satisfying both local uniqueness and local existence
of solutions, where $M$ is a manifold $($possibly with boundary$)$ modeled
on a locally convex space, $J\sub \R$ a non-degenerate interval
and $f\colon W\to TM$ a function on a subset~$W$ of $J\times U$
such that $f(t,y)\in T_y(M)$ for all $(t,y)\in W$.
For each $(t_0,y_0)\in W$, there exists
a solution $\gamma\colon I\to E$ to the initial value problem \emph{(\ref{theiniprob-m})}
such that
\begin{equation}\label{proper-max-m}
I_\eta\sub I\quad\mbox{and}\quad \eta=\gamma|_{I_\eta}
\end{equation}
for each solution $\eta\colon I_\eta\to M$ of~\emph{(\ref{theiniprob-m})}.
The function $\gamma$ is uniquely determined by the latter property,
and~$I$ is relatively open in~$J$.\qed
\end{lem}
\begin{defn}\label{flow-nonauto-m}
The solution $\gamma\colon I\to M$ described in Lemma~\ref{ex-maxi-sol-m}
for $(t_0,y_0)\in W$ is called the \emph{maximal solution}
to the initial value problem (\ref{theiniprob-m});
to emphasize its dependence on $t_0$ and $y_0$,
we also write $\gamma_{t_0,y_0}\colon I_{t_0,y_0}\to M$ 
in place of~$\gamma$.
The subset $\Omega\sub J\times J\times M$ given by
\[
\Omega:=\bigcup_{(t_0,y_0)\in W}I_{t_0,y_0}\times \{(t_0,y_0)\}
\]
is the domain of the so-called (maximal) \emph{flow}
$\Fl\colon \Omega\to M$, $\Fl(t,t_0,y_0):=\gamma_{t_0,y_0}(t)$.
Given $t,t_0\in J$, we let $\Omega_{t,t_0}:=\{y_0\in M\colon
\mbox{$(t_0,y_0)\in W$ and $t\in I_{t_0,y_0}$}\}$
and consider the partial map
\[
\Fl_{t,t_0} \colon\Omega_{t,t_0}\to M,\quad y_0\mto\Fl(t,t_0,y_0).
\]
\end{defn}
In the situation of Lemma~\ref{ex-maxi-sol-m}
and Definition~\ref{flow-nonauto-m}, we can repeat the proof
of Lemma~\ref{basics-flow} and obtain:
\begin{lem}\label{basics-flow-m}
\begin{description}[(D)]
\item[\rm(a)]
If $(t_0,y_0)\in W$ and $t_1\in I_{t_0,y_0}$,
then $\gamma_{t_0,y_0}=\gamma_{t_1,y_1}$ with $y_1:=\gamma_{t_0,y_0}(t_1)
=\Fl_{t_1,t_0}(y_0)$.\vspace{.8mm}
\item[\rm(b)]
If $t_2\in I_{t_1,y_1}$ in \emph{(a)},
then $t_2\in I_{t_0,y_0}$ and
$\Fl_{t_2,t_0}(y_0)=\Fl_{t_2,t_1}(\Fl_{t_1,t_0}(y_0))$.\vspace{.8mm}
\item[\rm(c)]
For all $t_1,t_0\in J$, the map $\Fl_{t_1,t_0}$ is injective,
$\Phi_{t_1,t_0}(\Omega_{t_1,t_0})=\Omega_{t_0,t_1}$, and $\Fl_{t_0,t_1}
=(\Fl_{t_1,t_0})^{-1}$.\qed
\end{description}
\end{lem}
The following proposition can be established like Lemma~\ref{ode-in-subspace};
we only need to replace~$y'$, $E$, and~$F$ with $\dot{y}$, $M$, and~$N$,
respectively,
in its proof.
\begin{prop}\label{ode-in-submfd}
Let~$M$ be a $C^1$-manifold modeled on a locally convex space~$E$.
Let $N\sub M$ be a closed subset which is a submanifold modeled on a closed vector subspace $F\sub E$,
and $J\sub\R$ be a non-degenerate interval.
Let $W\sub J\times M$ be a subset and $f\colon W\to TM$ be a
function such that $f(t,y)\in T_yM$ for all $(t,y)\in W$
and $f(t,y)\in T_yN$ for all $(t,y)\in W\cap(\R\times N)$.
Define $g\colon W\cap (\R\times N)\to TN$, $g(t,y):=f(t,y)$.
Assume that $\dot{y}(t)=f(t,y(t))$ satisfies local uniqueness of solutions; 
then also $\dot{y}(t)=g(t,y(t))$ satisfies local uniqueness of solutions.
If, moreover, $\dot{y}(t)=g(t,y(t))$ satisfies local existence of solutions,
then the following holds: If $\gamma\colon I\to E$ is a solution
to $\dot{y}(t)=f(t,y(t))$ and $y(t_0)\in N$ for some $t_0\in I$,
then $\gamma(I)\sub N$ and $\gamma$ solves $\dot{y}(t)=g(t,y(t))$.\qed
\end{prop}
\begin{defn}\label{defn-loc-flo-m}
Let $J\sub \R$ be a non-degenerate interval and $k\in \N_0\cup\{\infty\}$.
Let $M$ be a~$C^k$-manifold (if $k\geq 1$), resp., a $C^1$-manifold (if $k=0$),
possibly with boundary,
modeled on a locally convex space~$E$.
Let $f\colon W\to TM$ be a function on an open subset $W\sub J\times M$
such that $f(t,y)\in T_yM$ for all $(t,y)\in W$.
We say that the differential equation $\dot{y}(t)=f(t,y(t))$
\emph{admits local $C^k$-flows} if, for all $(\wb{t},\wb{y}) \in W$,
there exist a relatively open interval $I\sub J$ with $\wb{t}\in I$,
an open neighborhood~$V$ of $\wb{y}$ in~$M$
and a $C^k$-function
\[
\Phi\colon I\times I\times V\to M
\]
with the following properties:
\begin{description}[(D)]
\item[(a)]
For all $(t_0,y_0)\in I\times V$, the function
$I\to M$, $t\mto \Phi_{t,t_0}(y_0):=\Phi(t,t_0,y_0)$ is a solution to
the initial value problem (\ref{theiniprob-m});
\item[(b)]
There is an open $\wb{y}$-neighborhood $Y\sub V$
such that $\Phi_{t_1,t_0}(Y)\sub V$ for all $t_0,t_1\in I$ and
\[
\Phi_{t_2,t_1}(\Phi_{t_1,t_0}(y_0))=\Phi_{t_2,t_0}(y_0)
\;\;\mbox{for all $\,t_0,t_1,t_2\in I$ and $y_0\in Y$.}
\]
\end{description}
If~$J$ is an open interval, $M$ a real analytic manifold modeled on a locally convex space,
$f$ $($as above$)$ is real analytic and~$\Phi$ with {\rm(a)} and {\rm(b)}
can be chosen real analytic, then we say that $\dot{y}(t)=f(t,y(t))$
\emph{admits real analytic local
flows}.
\end{defn}
\begin{lem}\label{flo-cons}
If $\dot{y}(t)=f(t,y(t))$ $($as in
Definition~\emph{\ref{defn-loc-flo-m}}$)$
admits local $C^1$-flows,
then it satisfies local uniqueness of solutions.
\end{lem}
\begin{prf}
Note that
$y'(t)=f_\phi(t,y(t))$ (as in Remark~\ref{rem-ode-m}(b))
admits local $C^1$-flows
for each chart $\phi$ of~$M$
and hence satisfies local uniqueness of solutions,
by Proposition~\ref{locflow-uniq}.
The assertion follows.
\end{prf}
\begin{rem}\label{rems-ode-md}
If $M$, $k$, and $f\colon W\to TM$ are as in Definition~\ref{defn-loc-flo-m}
and
$\dot{y}(t)=f(t,y(t))$ satisfies local existence and local uniqueness of solutions,
then it admits local $C^k$-flows if and only if the domain
$\Omega$ of the corresponding maximal flow $\Fl\colon \Omega\to M$
is a neighborhood of $(t,t,y)$ in $J\times J\times M$
for all $(t,y)\in W$ and $\Fl$ is~$C^k$ on some open neighborhood
of $(t,t,y)$ in~$\Omega$.
In fact, for $\Phi\colon I\times I\times V\to M$ as in Definition~\ref{defn-loc-flo-m},
we have $I\times I\times V\sub \Omega$ and $\Fl|_{I\times I\times V}=\Phi$.
Conversely, we can define $\Phi$ as a restriction of~$\Fl$.
\end{rem}
To illustrate the general theory,
let us consider two examples, namely
1.) time-dependent
left invariant vector fields on a Lie group~$G$; and
2.) time-dependent fundamental vector fields
on a right $G$-space~$M$.
The
examples are directed to readers
who
are familiar with
the concept of a Lie group
modeled on a locally convex space, and basic Lie-theoretic facts
(as in Section~\ref{sec:3.1}).\medskip

In the following two lemmas, let $G$ be a Lie group
modeled on a locally convex space and
$\g:=\L(G):=T_\be(G)$ be the tangent space of~$G$ at the neutral element~$\be$.
We shall use the smooth left action
\begin{equation}\label{canleft}
G\times T(G)\to T(G),\;\; (g,v)\mto gv:=T(\lambda_g)v
\end{equation}
of~$G$ on its tangent bundle,
where $\lambda_g\colon G\to G$, $x\mto gx$ is left translation by~$g$.
Each $v\in \g$ defines a left invariant vector field
\[
v_l\colon G\to T(G),\quad v_l(g):=gv.
\]
Let $\cV(G)^l$ be the vector space of left invarant vector fields on~$G$.
Given a $C^k$-curve $\gamma\colon [0,1]\to\g$ with $k\in \N_0\cup\{\infty\}$,
we can consider the map
\[
[0,1]\to \cV(G)^l,\quad t\mto \gamma(t)_l
\]
and the corresponding time-dependent vector field
\[ [0,1]\times G\to T(G), \quad (t,g)\mto \gamma(t)_l(g)=g\gamma(t) \] 
which is a $C^k$-map as $\gamma$ is $C^k$ and the action~(\ref{canleft}) is smooth.
A $C^1$-curve $\eta\colon [0,1]\to G$ solves the differential equation
\[
\dot{y}(t)=\gamma(t)_l(y(t))=y(t)\gamma(t)
\]
if and only if its left logarithmic derivative $\delta(\eta)\colon [0,1]\to\g$,
$t\mto \eta(t)^{-1}\dot{\eta}(t)$ satisfies
$\delta(\eta)=\gamma$.
In this case, $\eta$ is $C^{k+1}$ (cf.\ Remark~\ref{glue-sol-ode}(c)).
\begin{lem}\label{liegpflow}
If $\gamma\colon [0,1]\to\g$ is a $C^k$-curve and there exists
a $C^{k+1}$-curve $\eta\colon [0,1]\to G$ with left logarithmic derivative
$\delta(\eta)=\gamma$,
then
\[
\dot{y}(t)=\gamma(t)_l(y(t))=y(t)\gamma(t)
\]
satisfies local existence and local uniqueness of solutions. The corresponding maximal
flow is defined on all of
$[0,1]\times [0,1]\times G$; it is the $C^{k+1}$-map
\[
\Fl\colon [0,1]\times [0,1]\times G\to G,\;\;
\Fl(t,t_0,y_0):=y_0\eta(t_0)^{-1}\eta(t).
\]
\end{lem}
\begin{prf}
As the multiplication map $\sigma\colon G\times G\to G$ is a smooth right
action of $G$ on itself and $v_l(g)=T_\be(\lambda_g)v=T\sigma(g,\cdot)v$
for all $g\in G$ and $v\in\g$, the lemma is a special case of the subsequent
Lemma~\ref{Gspaceflow}.
\end{prf}
For $G$ as before, let $M$ be a smooth manifold modeled
on a locally convex space (possibly with rough boundary)
and $\sigma\colon M\times G\to M$ be a smooth right action of $G$ on~$M$.
Given $v\in\g$, we define the corresponding \emph{fundamental vector field} $v_\sharp$ on $M$
as
\[
v_\sharp\colon M\to TM,\quad v_\sharp(p):=T(\sigma_p)(v)=T(\sigma)(0_p,v),
\]
where $M\to TM$, $p\mto 0_p\in T_pM$ is the zero-section
(which is smooth), $\sigma_p:=\sigma(p,\cdot)\colon G\to M$ is the orbit map at~$p\in M$
and $T(M\times G)$ is identified with $T(M)\times T(G)$, as usual.
If $\gamma\colon [0,1]\to \g$ is a $C^k$-curve, then the time-dependent vector field
\[
[0,1]\times M\to TM,\;\; (t,p)\mto\gamma(t)_\sharp(p)=T\sigma(0_p,\gamma(t))
\]
is $C^k$.
\begin{lem}\label{Gspaceflow}
If $\gamma\colon [0,1]\to\g$ is a $C^k$-curve and there exists a $C^{k+1}$-curve $\eta\colon [0,1]\to G$
with left logarithmc derivative $\delta(\eta)=\gamma$,
then
\[
\dot{y}(t)=\gamma(t)_\sharp(y(t))
\]
satisfies local existence and local uniqueness of solutions. The corresponding maximal flow is defined on all of
$[0,1]\times [0,1]\times M$; it is the $C^{k+1}$-map
\begin{equation}\label{Mfloweq}
\Fl\colon [0,1]\times [0,1]\times M\to M,\;\;
\Fl(t,t_0,y_0):=\sigma(y_0,\eta(t_0)^{-1}\eta(t)).
\end{equation}
\end{lem}
\begin{prf}
Abbreviate $\sigma(p,g)=:pg$ for $p\in M$ and $g\in G$.
We define
\[
\Phi\colon [0,1]\times[0,1]\times M\to M,\quad (t,t_0,y_0)\mto y_0\eta(t_0)^{-1}\eta(t)
\]
and write $\Phi_{t,t_0}(y_0):=\Phi(t,t_0,y_0)$.
Then~$\Phi$ is a $C^{k+1}$-map.
Given $y_0\in M$ and $t_0\in [0,1]$, we define $\zeta\colon [0,1]\to G$ via
$\zeta(t):=\eta(t_0)^{-1}\eta(t)$ and
consider $\phi\colon [0,1]\to M$, $\phi(t):=y_0\zeta(t)=\Phi_{t,t_0}(y_0)$.
Then $\delta(\zeta)=\delta(\eta)=\gamma$ and
\begin{eqnarray*}
\phi'(t)&=&\frac{d}{ds}\Big|_{s=0}\phi(t+s)=\frac{d}{ds}\Big|_{s=0}y_0\zeta(t+s)\\
&=& \frac{d}{ds}\Big|_{s=0}
\sigma(y_0\zeta(t),\zeta(t)^{-1}\zeta(t+s))=T(\sigma_{y_0\zeta(t)})\delta(\zeta)(t)
=\gamma(t)_\sharp(\phi(t)),
\end{eqnarray*}
whence $\phi$ solves (\ref{Mfloweq}).
For all
$t,t_0,t_1\in [0,1]$ and $y_0\in M$, setting $y_1:=\Phi_{t_1,t_0}(y_0)=y_0\eta(t_0)^{-1}\eta(t_1)$
we have
\begin{eqnarray*}
\Phi_{t,t_1}(\Phi_{t_1,t_0}(y_0))&=& \Phi_{t,t_1}(y_1)=y_1\eta(t_1)^{-1}\eta(t)\\
&=&y_0\eta(t_0)^{-1}\eta(t_1)\eta(t_1)^{-1}\eta(t)=y_0\eta(t_0)^{-1}\eta(t)=\Phi_{t,t_0}(y_0).
\end{eqnarray*}
Thus $\Phi$ is a local flow for (\ref{Mfloweq}).
Hence (\ref{Mfloweq}) satisfies local uniqueness of solutions (see
Lemma~\ref{flo-cons})
and local existence. Moreover, $\Fl=\Phi$.
\end{prf}
\begin{rem}
(a) If $k=\infty$, then the uniqueness result in Lemma~\ref{Gspaceflow} (and hence in Lemma~\ref{liegpflow})
could also be deduced
from Lemma~\ref{unique-intcurve}.\medskip

\noindent
(b) For another proof of the uniqueness statement in Lemma~\ref{liegpflow}
based on Lie-theoretic arguments, cf.\ Lemma~\ref{lem:c.12b}.
\end{rem}
\begin{thm}\label{ode-loc-glob}
Let $J\sub \R$ be a non-degenerate interval and $k\in \N_0\cup\{\infty\}$.
Let $M$ be a~$C^k$-manifold $($if $k\geq 1)$, resp., a $C^1$-manifold $($if $k=0)$,
possibly with boundary,
modeled on a locally convex space~$E$.
Let $f\colon W\to TM$ be a function on an open subset $W\sub J\times M$
such that $f(t,y)\in T_yM$ for all $(t,y)\in W$.
Assume that $\dot{y}(t)=f(t,y(t))$
admits local $C^k$-flows;
if $k=0$, assume moreover that $\dot{y}(t)=f(t,y(t))$
satisfies local uniqueness of solutions.
Let $\Fl\colon \Omega\to M$ be the maximal flow of $\dot{y}(t)=f(t,y(t))$.
Then we have:
\begin{description}[(D)]
\item[\rm(a)]
The domain $\Omega$ of the maximal flow~$\Fl$ is open in $J\times J\times M$
and $\Fl$ is~$C^k$;\vspace{.7mm}
\item[\rm(b)]
$\Omega_{t,t_0}$ is open in~$M$
for all $t_0,t\in J$ and the map
$\Phi_{t,t_0}\colon\Omega_{t,t_0}\to \Omega_{t_0,t}$
is a $C^k$-diffeomorphism.\footnote{If $k=0$, then $\Phi_{t,t_0}$ is a homeomorphism.}
\end{description}
If~$J$ is an open interval,
$M$ a real analytic manifold modeled on a locally convex space,
$f$ is real analytic and $\dot{y}(t)=f(t,y(t))$ admits
real analytic local flows, then also $\Fl\colon\Omega\to M$ is real analytic and the maps
$\Fl_{t,t_0}\colon\Omega_{t,t_0}\to\Omega_{t_0,t}$ are real analytic diffeomorphisms. 
\end{thm}
\begin{prf}
In the real analytic case, set $k:=\omega$.
Taking a singleton set of parameters,
we obtain
openness of~$\Omega$ in $J\times J\times M$ and the~$C^k$-property of~$\Fl$
as a special case of the next result (Theorem~\ref{ode-loc-glob-par}).
Openness of~$\Omega$ implies that the sets~$\Omega_{t,t_0}$ are open in~$M$.
We know from Lemma~\ref{basics-flow-m}(c)
that $\Fl_{t,t_0}\colon\Omega_{t,t_0}\to\Omega_{t_0,t}$
is a bijection, with inverse $\Fl_{t_0,t}$. Since~$\Fl$ is~$C^k$,
also its partial maps $\Fl_{t,t_0}$ and $\Fl_{t_0,t}$ are~$C^k$.
Hence $\Fl_{t,t_0}$ is a $C^k$-diffeomorphism.
\end{prf}
In the following two definitions, let
$J\sub \R$ be a non-degenerate interval, $P$ be a $C^k$-manifold
modeled on a locally convex space~$F$, with $k\in \N_0\cup\{\infty,\omega\}$,
and
$M$ be a~$C^k$-manifold (if $k\geq 1$) or a $C^1$-manifold (if $k=0$),
modeled on a locally convex space~$E$. If $k\not=\omega$,
both~$P$ and~$M$ may have a boundary.
If $k=\omega$, we assume that~$J$ is an open interval.
Let $W\sub J\times M\times P$ be an open subset and $f\colon W\to TM$ be a function
such that $f(t,y,p)\in T_yM$ for all $(t,y,p)\in W$.
\begin{defn}\label{defn-ode-par-mfd}
If $\dot{y}(t)=f(t,y(t),p)$ satisfies local existence and local uniqueness of solutions for all $p\in P$,
we let $\gamma_{t_0,y_0,p}\colon I_{t_0,y_0,p}\to M$ be the maximal solution to the
initial value problem
\begin{equation}\label{ivp-manifp}
\dot{y}(t)=f(t,y(t),p),\quad y(t_0)=y_0
\end{equation}
for $(t_0,y_0,p)\in W$.
We define
\[
\Omega:=\bigcup_{(t_0,y_0,p)\in W}I_{t_0,y_0,p}\times\{(t_0,y_0,p)\}
\sub J\times J\times M\times P
\]
in this case and call the map
$\Fl\colon \Omega\to M$, $\Fl(t,t_0,y_0,p):=\gamma_{t_0,y_0,p}(t)$
the associated (maximal) flow.
\end{defn}
\begin{defn}\label{defn-loc-flo-m-para}
We say that the differential equation $\dot{y}(t)=f(t,y(t),p)$
\emph{admits local $C^k$-flows with parameters} if, for all $(\wb{t},\wb{y},p) \in W$,
there exist a relatively open interval $I\sub J$ with $\wb{t}\in I$,
an open $\wb{y}$-neighborhood $V\sub M$,
an open $p$-neighborhood $Q\sub P$
and a $C^k$-function
\[
\Phi\colon I\times I\times V\times Q\to M
\]
with the following properties:
\begin{description}[(D)]
\item[(a)]
For all $(t_0,y_0,q)\in I\times V\times Q$, the function
$I\to M$, $t\mto \Phi_{t,t_0,q}(y_0):=\Phi(t,t_0,y_0,q)$ is a solution to
the initial value problem (\ref{theiniprob-m});
\item[(b)]
There exists an open $\wb{y}$-neighborhood $Y\sub V$
such that $\Phi_{t_1,t_0,q}(Y)\sub V$ for all $t_0,t_1\in I$ and $q\in Q$,
and moreover
\[
\Phi_{t_2,t_1,q}(\Phi_{t_1,t_0,q}(y_0))=\Phi_{t_2,t_0,q}(y_0)
\;\mbox{for all $\,t_0,t_1,t_2\in I$, $y_0\in Y$, and $q\in Q$.}
\]
\end{description}
\end{defn}
\begin{thm}\label{ode-loc-glob-par}
Let $k\in \N_0\cup\{\infty\}$ and $J\sub \R$ be a non-degenerate interval.
Let $M$ be a~$C^k$-manifold $($if $k\geq 1)$, resp., a $C^1$-manifold $($if $k=0)$,
possibly with boundary,
modeled on a locally convex space. Let~$P$ be a $C^k$-manifold, possibly with boundary,
modeled on a locally convex space.
Let $f\colon W\to TM$ be a function on an open subset $W\sub J\times M\times P$
such that $f(t,y,p)\in T_yM$ for all $(t,y,p)\in W$.
Assume that $\dot{y}(t)=f(t,y(t),p)$
admits local $C^k$-flows with parameters;
if $k=0$, assume moreover that $\dot{y}(t)=f(t,y(t),p)$
satisfies local uniqueness of solutions, for all $p\in P$.
Then we have:
\begin{description}[(D)]
\item[\rm(a)]
The domain $\Omega$ of the corresponding flow
$\Fl\colon \Omega\to M$ with parameters
is open in $J\times J\times M\times P$
and $\Fl\colon \Omega\to M$ is a $C^k$-map.
\item[\rm(b)]
If~$J$ is an open interval, $M$ and~$P$ are real analytic manifolds modeled on locally convex spaces,
$f$ is real analytic and $\dot{y}(t)=f(t,y(t),p)$ admits real analytic
local flows with parameters,
then~$\Fl$ is real analytic.
\end{description}
\end{thm}
\begin{prf}
In the real analytic case, let us set $k:=\omega$.
Given $(\tau,t_0,y_0,p)\in \Omega$, let~$K$ be the set of
all $t_1\in I_{t_0,y_0,p}$ for which~$\Omega$ contains a neighborhood~$Q(t_1)$
of $(t_1,t_0,y_0,p)$ in $J\times J\times M\times P$
such that~$\Fl|_{Q(t_1)}$ is~$C^k$. By definition, $K$ is relatively open in~$I_{t_0,y_0,p}$.
As the differential equation admits local $C^k$-flows with parameters,
we have $t_0\in K$.
If we can show that~$K$ is closed in~$I_{t_0,y_0,p}$,
then $K=I_{t_0,y_0,p}$ by connectedness of the latter interval. Notably, $\tau\in K$,
whence $\Omega$ (which contains~$Q(\tau)$) is a neighborhood of
$(\tau,t_0,y_0,p)$ in
$J\times J\times M\times P$ and $\Fl|_{Q(\tau)}$ is~$C^k$ (from which all assertions follow).
To see that~$K$ is closed, let~$\theta$ be in the closure
of~$K$ in~$I_{t_0,y_0,p}$. As the differential equation admits local flows with
parameters, there exists an open neighborhood
$S\sub J\times J\times M\times P$ of $(\theta,\theta,\gamma_{t_0,y_0,p}(\theta),p)$
such that $S\sub\Omega$ and $\Fl|_S$ is~$C^k$.
After shrinking~$S$, we may assume that
$S=I_1\times Y$, where $I_1\sub J$ is a relatively open
interval with $\theta\in I_1$ and $Y$ is an open subset of
$J\times M\times P$.
There exists $t_2\in K$ such that $(t_2,t_2,\gamma_{t_0,y_0,p}(t_2),p)\in S$.
After shrinking $Q:=Q(t_2)$, we may assume that
$Q=I_0\times Z$ for a relatively open interval
$I_0\sub J$ with $t_2\in I_0$ and an open $(t_0,y_0,p)$-neighborhood
$Z\sub J\times M\times P$.
Moreover, we may assume that
\[
(t_2,\gamma_{t_1,y_1,q}(t_2),q)\in Y \;\,\mbox{for all $(t_1,y_1,q)\in Z$.}
\]
Now Lemma~\ref{basics-flow-m}
shows that $(t,t_1,y_1,q)\in\Omega$ for all $(t,t_1,y_1,q)\in I_1\times Z$ and
\[
\Fl(t,t_1,y_1,q)=\Fl|_S(t,t_2,\Fl|_Q(t_2,t_1,y_1,q),q),
\]
which is a $C^k$-function of $(t,t_1,y_1,q)\in I_1\times Z$. Hence $\theta\in K$,
and we deduce that $K$ is closed.
\end{prf}
\begin{defn}\label{def-lipcond-mfd}
Let $M$ and $N$ be $C^1$-manifolds (possibly with boundary) modeled on normed spaces
$(E,\|\cdot\|_E)$ and $(F,\|\cdot\|_F)$, respectively,
and~$X$ be a topological space. A continuous mapping $f\colon W\to N$
on a subset $W\sub X\times M$
is said to satisfy a \emph{local Lipschitz condition in the second
argument} if it has the following property:
For each $(t,y)\in W$, there exists
a $t$-neighborhood $X_0\sub X$, a chart $\phi\colon U_\phi\to V_\phi\sub E$
of~$M$ around~$y$ and a chart \break $\psi\colon U_\psi\to V_\psi\sub F$
of~$N$ around $f(t,y)$ with $f(W\cap (X_0\times U_\phi))\sub U_\psi$
such that the map
\[
(\id_X\times \phi)(W\cap (X_0\times U_\phi))\to F,\;\, (a,b)\mto \psi(f(a,\phi^{-1}(b))),
\]
which is defined on a subset of $X\times E$, satisfies a local Lipschitz condition in the
second argument.\footnote{We assume that $M$ and $N$ are $C^1$-manifolds
to ensure that the property just defined is independent
of the choice of charts (see Exercise~\ref{exc-loc-lip-indep}).}
\end{defn}
If also $P$ is a topological space and $W$ a subset of $X\times M\times P$,
we say that a mapping $f\colon W\to N$
satisfies a local Lipschitz condition in the second argument if
the mapping
$((t,p),y)\mto f(t,y,p)$ on the corresponding subset of
$(X\times P)\times M$ does so.
\begin{cor}\label{pardep-ode-ban-m}
Let $J\sub \R$ be a non-degenerate interval, $k\in \N_0\cup\{\infty\}$,
and $P$ be a $C^k$-manifold $($possibly with boundary$)$ modeled
on a locally convex space.
Let $M$ be a $C^{k+1}$-manifold $($if $k\geq 1)$,
resp., a $C^2$-manifold $($if $k=0)$, modeled
on a Banach space. 
Let $W\sub J\times M\times P$ be an open subset and
\[
f\colon W\to TM,\quad (t,y,p)\mto f(t,y,p)
\]
be a $C^k$-function such that $f(t,y,p)\in T_yM$ for all $(t,y,p)\in W$.
If~$k=0$, assume that~$f$ satisfies a local Lipschitz condition in the second variable.
Then the following holds:
\begin{description}[(D)]
\item[\rm(a)]
The differential equation
\begin{equation}\label{odewithpa-m}
\dot{y}(t)=f(t,y(t),p)
\end{equation}
satisfies local existence and local uniqueness of solutions, for each $p\in P$;
\item[\rm(b)]
The domain $\Omega$ of the maximal flow of {\rm(\ref{odewithpa-m})}
is open in $J\times J\times M\times P$;
\item[\rm(c)]
The maximal flow $\Fl\colon \Omega\to M$ of~{\rm(\ref{odewithpa-m})}
is~$C^k$.
\end{description}
If $J$ is an open interval, $P$ is a real analytic manifold,
$M$ a real analytic Banach manifold and $f\colon W\to TM$
is real analytic, then~$\Fl$ is real analytic.
\end{cor}
\begin{prf}
Let~$E$ and~$F$ be the modeling spaces of~$M$
and~$P$, respectively.
In the real analytic case, let $k:=\omega$.
Given $a:=(\wb{t},\wb{y},p)\in W$,
let $\phi\colon U_\phi\to V_\phi\sub E$
and $\psi\colon U_\psi\to V_\psi\sub F$
be charts of~$M$ and~$P$ around~$\wb{y}$ and~$p$,
respectively. Define $W_a:=\{(t,y,q)\in J\times V_\phi\times V_\psi\colon
(t,\phi^{-1}(y),\psi^{-1}(q))\in W\}$ and
\[
f_a\colon W_a\to E,\quad (t,y,q)\mto d\phi (f(t,\phi^{-1}(y),\psi^{-1}(q))).
\]
Now $y'(t)=f_a(t,y(t),q)$ satisfies local existence and local uniqueness
of\linebreak
solutions
and admits local $C^k$-flows with parameters for each $a\in W$,
by Proposition~\ref{uniq-normd} and Theorem~\ref{local-ivp-dep-Ck}.
In view of Remark~\ref{rem-ode-m}(b),
assertion~(a) follows from the preceding
and $\dot{y}(t)=f(t,y(t),p)$ admits local $C^k$-flows with parameters.
Thus~(b) and~(c) hold, by Theorem~\ref{ode-loc-glob-par}.
\end{prf}
Corollary~\ref{pardep-ode-ban-m}
(with $P$ a singleton)
and Theorem~\ref{ode-loc-glob} imply:
\begin{cor}\label{niceflow-ban-m}
Let $J\sub \R$ be a non-degenerate interval and $k\in \N_0\cup\{\infty\}$.
Let $M$ be a $C^{k+1}$-manifold $($if $k\geq 1)$,
resp., a $C^2$-manifold $($if $k=0)$ which is modeled
on a Banach space. 
Let $W\sub J\times M$ be an open subset and
$f\colon W\to TM$
be a $C^k$-function such that $f(t,y)\in T_yM$ for all $(t,y)\in W$.
If~$k=0$, assume that~$f$ satisfies a local Lipschitz condition in the second variable.
Then the following holds:
\begin{description}[(D)]
\item[\rm(a)]
The differential equation
\begin{equation}\label{odewithpa-4}
\dot{y}(t)=f(t,y(t))
\end{equation}
satisfies local existence and and local uniqueness of solutions.
\item[\rm(b)]
The domain $\Omega$ of the maximal flow of {\rm(\ref{odewithpa-4})}
is open in $J\times J\times M$.
\item[\rm(c)]
The maximal flow $\Fl\colon \Omega\to M$ of~{\rm(\ref{odewithpa-4})}
is~$C^k$.
\item[\rm(d)]
$\Omega_{t,t_0}$ is open in~$M$
for all $t_0,t\in J$ and the map
$\Phi_{t,t_0}\colon\Omega_{t,t_0}\to \Omega_{t_0,t}$
is a $C^k$-diffeomorphism.
\end{description}
If $J$ is an open interval, $P$ is a real analytic manifold,
$M$ a real analytic Banach manifold and $f\colon W\to TM$
real analytic, then $\Fl\colon \Omega\to M$ is real analytic
and the maps $\Fl_{t,t_0}$ are real analytic diffeomorphisms.\qed
\end{cor}
\begin{rem}
For $k\geq 1$, an analog of Corollary~\ref{niceflow-ban-m}
also holds if $C^k$-maps and $C^{k+1}$-manifolds
are replaced with $FC^k$-maps and $FC^{k+1}$-manifolds,
respectively,
in all hypotheses and conclusions.
Replacing, moreover, $P$ with an $FC^k$-manifold
modeled on a normed space, $\Fl$ will be $FC^k$
also in Corollary~\ref{pardep-ode-ban-m}.
We only need to replace Theorem~\ref{local-ivp-dep-Ck}
with Corollary~\ref{paraFCk} in the proof.
\end{rem}
\begin{defn}\label{flow-auton}
Let $M$ be a $C^1$-manifold (possibly with boundary) modeled
on a locally convex space and $X\colon M\to TM$ a continuous
vector field such that $\dot{y}(t)=X(y(t))$ satisfies local existence
and local uniqueness of solutions.
Write $\gamma_{y_0}\colon I_{y_0}\to M$ for the maximal solution to $\dot{y}(t)=X(y(t))$, $y(0)=y_0$,
for $y_0\in M$, and define
\[
\Omega^X:=\bigcup_{y_0\in M}I_{y_0}\times\{y_0\}
\]
and $\Fl^X\colon \Omega^X\to M$, $\Fl^X(t,y_0):=\gamma_{y_0}(t)$.
In this context, the maximal solution $\gamma_{y_0}$
is also called the \emph{integral curve} passing through~$y_0$.
Write $\Fl^X_t:=\Fl(t,\cdot)\colon \Omega_t^X \to M$ for $t\in\R$,
with $\Omega_t^X:=\{y_0\in M\colon (t,y_0)\in\Omega^X\}=\{y_0\in M\colon t\in I_{y_0}\}$.
\end{defn}
\begin{rem}\label{bridgetoautono}
(a) In the setting of Definition~\ref{flow-auton},
define a function $f\colon  \R \times  M\to TM$ via $f(t,y):=X(y)$.
Using the notation $\gamma_{t_0,y_0}$, $I_{t_0,y_0}$, $\Omega$, $\Fl\colon \Omega\to M$,
$\Fl_{t,t_0}$ and $\Omega_{t_1,t_0}$ introduced for~$f$
in Definition~\ref{flow-nonauto-m},
we have
$\gamma_{y_0}=\gamma_{0,y_0}$, $I_{y_0}=I_{0,y_0}$,
\[
\Omega^X=\{(t,y)\in \R\times M\colon (t,0,y)\in \Omega\},\;\,
\Fl^X(t,y)=\Fl(t,0,y),
\]
$\Omega_t^X=\Omega_{t,0}$, and $\Fl^X_t=\Fl_{t,0}$.\medskip

\noindent
(b) Note that
\[
I_{t_0+t_1,y_0}=I_{t_0,y_0}+t_1\;\,\mbox{and $\,\gamma_{t_0+t_1,y_0}(t)=\gamma_{t_0,y_0}(t-t_1)$
for all $t\in I_{t_0+t_1,y_0}$.}
\]
In fact,
$I_{t_0,y_0}+t_1\to M$, $t\mto \gamma_{t_0,y_0}(t-t_1)$
solves $\dot{y}(t)=X(y(t))$, $y(t_0+t_1)=y_0$ for all $t_0,t_1\in\R$ and $y_0\in M$,
whence $I_{t_0,y_0}+t_1\sub I_{t_0+t_1,y_0}$ and $\gamma_{t_0+t_1,y_0}(t)=\gamma_{t_0,y_0}(t-t_1)$
for all $t\in I_{t_0,y_0}+t_1$. Replacing $(t_0,t_1)$ with $(t_0+t_1,-t_1)$ here, we find that also
$I_{t_0+t_1,y_0}-t_1\sub I_{t_0,y_0}$, whence $I_{t_0+t_1,y_0}\sub I_{t_0,y_0}+t_1$
and thus $I_{t_0+t_1,y_0}=I_{t_0,y_0}+t_1$.\medskip

\noindent
(c) If $y_0\in M$ and $s,t\in\R$ such that $(s,y_0)\in \Omega^X$, then $(s+t,y_0)\in\Omega^X$
if and only if $(t,\Fl^X_s(y_0))\in \Omega^X$; in this case,
and $\Fl^X_t(\Fl^X_s(y_0))=\Fl^X_{s+t}(y_0)$.\\[1mm]
In fact, $I_{0,y_0}=I_{s,y_1}$ and $\gamma_{0,y_0}=\gamma_{s,y_1}$ by
Lemma~\ref{basics-flow-m}(a), with $y_1:=\gamma_{0,y_0}(s)=\Fl^X_s(y_0)$.
Now $(s+t,y_0)\in \Omega^X$ implies $s+t\in I_{0,y_0}=I_{s,y_1}=I_{0,y_1}+s$.
Thus $t\in I_{0,y_1}$ (whence $(t,y_1)\in \Omega^X$)
and $\Fl^X_t(\Fl^X_s(y_0))=\gamma_{0,y_1}(t)=\gamma_{s,y_1}(s+t)=\gamma_{0,y_0}(s+t)
=\Fl^X_{s+t}(y_0)$, using part~(b) of this remark.
Conversely, $(t,y_1)\in\Omega^X$ implies $t\in I_{0,y_1}$.
Using~(b), we get $s+t\in I_{s,y_1}=I_{0,y_0}$
and thus $(s+t,y_0)\in \Omega^X$.
\end{rem}
\begin{lem}
Let $k\in\N_0\cup\{\infty,\omega\}$.
If $k=0$, let~$M$ be a $C^1$-manifold modeled on a locally convex space;
otherwise, let $M$ be a $C^k$-manifold modeled on a locally convex space;
if $k\not=\omega$, we may allow~$M$ to have a boundary.
Let $X\colon M\to TM$ be a function such that $X(y)\in T_yM$ for each $y\in M$.
Consider the following conditions:
\begin{description}[(D)]
\item[\rm(a)]
$\dot{y}(t)=X(y(t))$ admits local $C^k$-flows.
\item[\rm(b)]
For each $\overline{y}\in M$, there exist
$\ve>0$, an open $\overline{y}$-neighborhood $V\sub M$,
and a $C^k$-function
\[
\Psi\colon \,]{-\ve},\ve[\,\times V\to M,\quad (t,y)\mto \Psi(t,y)=:\Psi_t(y)
\]
with the following properties:
\begin{description}[(D)]
\item[\rm(i)]
For each $y_0\in V$, the function $\Psi(\cdot,y_0)$
solves the initial value problem $\dot{y}(t)=X(y(t))$, $y(0)=y_0$.
\item[\rm(ii)]
There exists a $\oline{y}$-neighborhood $Y\sub V$ such that
$\Psi_t(y)\in V$ for all $(t,y)\in\,]{-\ve},\ve[\,\times Y$ and
$\Psi_s(\Psi_t(y))=\Psi_{s+t}(y)$ for all $y\in Y$ and $s,t\in\,]{-\ve},\ve[$
with $s+t\in\,]{-\ve},\ve[$.
\end{description}
\end{description}
Then \emph{(b)} implies \emph{(a)}.
If $\dot{y}(t)=X(y(t))$ satisfies local uniqueness
of solutions or $k\in\N\cup\{\omega\}$,
then \emph{(a)} and \emph{(b)} are equivalent.
\end{lem}
\begin{prf}
(b)$\Rightarrow$(a):
If $\Psi$ and~$Y$
are as described, given $\overline{t}\in\R$ we define $I:=\,]\overline{t}-\ve/2,\overline{t}+\ve/2[$ and note that
$\Phi\colon I\times I\times V\to M$, $\Phi(t,t_0,y_0):=\Psi(t-t_0,y_0)$
and the given~$Y$ satisfy the conditions in Definition~\ref{defn-loc-flo-m}.

(a)$\impl$(b):
Assume that $\dot{y}(t)=X(y(t))$ has local $C^k$-flows.
If $k\in\N\cup\{\infty\}$, this implies that $\dot{y}(t)=X(y(t))$
satisfies local uniqueness of solutions; if $k=0$,
we assume the latter as a hypothesis.
Hence, a maximal flow $\Fl$ can be associated with $\dot{y}=X(y(t))$ as in
Definition~\ref{flow-nonauto-m}.
Let $\Phi\colon I\times I\times V\to M$ and~$Y$ be as in
Definition~\ref{defn-loc-flo-m}, with $\wb{t}:=0$ and a given element $\wb{y}\in M$.
There exists $\ve>0$ with $]{-\ve}, \ve[\,\sub I$.
Then $\Psi\colon \,]{-\ve},\ve[\,\times V\to M$, $\Psi(t,y):=\Phi_{t,0}(y)$ is~$C^k$ and
satisfies~(i) and~(ii),
using that $\Phi=\Fl|_{I\times I\times V}$ and thus $\Psi(t,y)=\Fl(t,0,y)=\Fl^X_t(y)$
with properties as in Remark~\ref{bridgetoautono}(c).
\end{prf}
\begin{prop}\label{leave-cp}
Let $J\sub\R$ be a non-degenerate
interval,
$M$ be a $C^1$-manifold $($possibly with boundary$)$
modeled on a locally convex space,  and $f\colon J\times M\to TM$
be a function with $f(t,y)\in T_yM$ for all $(t,y)\in J\times M$.
Assume that $\dot{y}(t)=f(t,y(t))$ satisfies local existence
and local uniqueness of solutions, and admits local $C^0$-flows.
Let $(t_0,y_0)\in J\times M$ and assume that
the maximal solution $\gamma\colon I_{t_0,y_0}\to M$
to $\dot{y}(t)=f(t,y(t))$, $y(t_0)=y_0$
is not defined for all possible times in the future, i.e.,
$I_{t_0,y_0}\cap [t_0,\infty[$
is a proper subset of $J\cap [t_0,\infty[$.
Then the following holds:
For each compact subset $K\sub M$,
there exists $t_K\in I_{t_0,y_0}\cap [t_0,\infty[$
such that $\gamma(t)\not\in K$ for all $t\in I_{t_0,y_0}\cap\,]t_K,\infty[$.
\end{prop}
\begin{prf}
Since $I_{t_0,y_0}\cap [t_0,\infty[$ is a proper subinterval of $J\cap [t_0,\infty[$,
the interval $J\cap [t_0,\infty[$ has to be non-degenerate.
In view of local existence of solutions, $I_{t_0,y_0}\cap [t_0,\infty[$ cannot be a closed
interval; thus $I_{t_0,y_0}\cap[t_0,\infty[\,=[t_0,\overline{t}[$
for some~$\overline{t}>t_0$. Then $\overline{t}\in J$. If the assertion was wrong,
we could find a compact subset $K\sub M$ and a sequence $(t_n)_{n\in\N}$ in
$[t_0,\overline{t}[$ such that $\gamma(t_n)\in K$ for all $n\in\N$ and $t_n\to \overline{t}$ as $n\to\infty$.
Let $\overline{y}$ be a cluster point of $(\gamma(t_n))_{n\in\N}$ in the compact topological space~$K$.
Let $I\sub J$ be a relatively open interval containing~$\overline{t}$
and $V\sub M$ be an open $\overline{y}$-neighborhood such that $I\times I\times V$
is contained in the domain $\Omega\sub J\times J\times M$ of the maximal
flow~$\Fl$ of $\dot{y}(t)=f(t,y(t))$
(cf.\ Definition~\ref{defn-loc-flo-m}). There exists $n_0\in\N$ such that $t_n\in I$ for all
$n\geq n_0$. By definition of a cluster point, we find $m\geq n_0$ such that $\gamma(t_m)\in V$.
Then $I\times \{(t_m,\gamma(t_m))\}\sub I\times I\times V\sub \Omega$ shows that the
domain of definition $I_{t_m,\gamma(t_m)}$
of the maximal solution $\gamma_{t_m,\gamma(t_m)}$
to the initial value problem $\dot{y}(t)=f(t,y(t))$, $y(t_m)=\gamma(t_m)$
contains~$I$. But $\gamma_{t_m,\gamma(t_m)}=\gamma_{t_0,y_0}=\gamma$,
whence $\overline{t}\in I\sub I_{t_0,y_0}=[t_0,\overline{t}[$, contradiction.
\end{prf}
\begin{rem}\label{backwd}
Of course, an analogous result applies if $I_{t_0,y_0}\cap \,]{-\infty},t_0]$ is a proper subinterval of $J\cap \,]{-\infty},t_0]$.
\end{rem}
\begin{cor}\label{ode-cp-mfd-1}
Let $M$ be a compact $C^1$-manifold $($possibly with boundary$)$
$J\sub\R$ be a non-degenerate
interval and $f\colon J\times M\to TM$
be a function with $f(t,y)\in T_yM$ for all $(t,y)\in J\times M$.
Assume that $\dot{y}(t)=f(t,y(t))$ satisfies local existence
and uniqueness of solutions, and admits local $C^0$-flows.
Then the initial value problem $\dot{y}(t)=f(t,y(t))$, $y(t_0)=y_0$
has a solution defined on all of~$J$, for all $(t_0,y_0)\in J\times M$.
\end{cor}
\begin{prf}
For all $(t_0,y_0)\in J\times M$,
we have $I_{t_0,y_0}\cap [t_0,\infty[=J\cap [t_0,\infty[$ by Proposition~\ref{leave-cp}
and also $I_{t_0,y_0}\cap \,]{-\infty},t_0]=J\cap\,]{-\infty},t_0]$ (see Remark~\ref{backwd}),
whence $I_{t_0,y_0}=J$.
\end{prf}
\noindent
Combining Corollaries~\ref{niceflow-ban-m} and~\ref{ode-cp-mfd-1}
(see also Remark~\ref{rems-ode-md}(b)), we~get:
\begin{cor}\label{ode-cp-mfd-2}
Let $J\sub\R$ be a non-degenerate
interval,
$M$ be a compact $C^2$-manifold without boundary,
and $f\colon J\times M\to TM$
be a continuous function with $f(t,y)\in T_yM$ for all $(t,y)\in J\times M$
which is locally Lipschitz in its second argument.
Then the initial value problem $\dot{y}(t)=f(t,y(t))$, $y(t_0)=y_0$
has a solution defined on all of~$J$, for all $(t_0,y_0)\in J\times M$.\qed
\end{cor}
Compactly supported vector fields can be treated in the same way.
\begin{cor}
Let $J\sub\R$ be a non-degenerate
interval,
$M$ be a finite-dimensional $C^2$-manifold without boundary,
and $f\colon J\times M\to TM$
be a continuous function with $f(t,y)\in T_yM$ for all $(t,y)\in J\times M$
which is locally Lipschitz in its second argument.
If there exists a compact subset $K\sub M$ such that $f(t,y)=0$ for all $t\in J$
and $y\in M\setminus K$,
then the initial value problem $\dot{y}(t)=f(t,y(t))$, $y(t_0)=y_0$
has a solution defined on all of~$J$, for all $(t_0,y_0)\in J\times M$.\qed
\end{cor}
\begin{small}
\subsection*{Exercises for Section~\ref{sec-ode-mfd}}

\begin{exer}\label{exc-flow-in-chart}
Let $M$ be a $C^1$-manifold (possibly with boundary) modeled on a locally convex space~$E$.
Let $J\sub \R$ be non-degenerate interval and
$f\colon W\to TM$ be a function on a subset $W\sub J\times TM$
such that $f(t,x)\in T_xM$ for all $(t,x)\in W$.
Let $\phi\colon U_\phi\to V_\phi$ be a chart for~$M$
and $f_\phi\colon W_\phi\to E$ be as in Remark~\ref{rem-ode-m}(b).
Assume that $\dot{y}(t)=f(t,y(t))$ satisfies local existence and local uniqueness of solutions.
\begin{description}[(D)]
\item[(a)]
Show that $y'(t)=f_\phi(t,y(t))$ satisfies local existence and local uniqueness of solutons.
\item[(b)]
Let $\Fl\colon \Omega\to M$ be the maximal flow of $\dot{y}(t)=f(t,y(t))$ and
$\Phi\colon Q\to E$ be the maximal flow of $y'(t)=f_\phi(t,y(t))$.
Show that $(\id_J\times\id_J\times\phi^{-1})(Q)\sub\Omega$ and
\[
\Fl_{t,t_0}(\phi^{-1}(x))=\phi^{-1}\big(\Phi_{t,t_0}(x)\big)\quad\mbox{for all
$\,(t,t_0,x)\in Q$.}
\]
\end{description}
\end{exer}

\begin{exer}\label{exc-flo-opensub}
Let $M$ be a $C^1$-manifold, $J\sub\R$ be a non-degenerate interval,
$W\sub J\times M$ be a subset
and $f\colon W\to TM$ be a function such that $f(t,y)\in T_yM$ for all $(t,y)\in W$.
Assume that $\dot{y}(t)=f(t,y(t))$ satisfies
local existence and local uniqueness of solutions,
and let $\Fl\colon \Omega\to M$ be its maximal flow.
Given an open subset $U\sub M$, consider
$f_U:=f|_{W\cap( J\times U)}\colon W\cap(J\times U)\to TU$. Prove the following:
\begin{description}[(D)]
\item[(a)]
The differential equation $\dot{y}(t)=f_U(t,y(t))$ satisfies local existence and uniqueness of
solutions; its maximal flow $\Fl^U$ has domain
\[
\Omega^U:=\{(t,t_0,y_0)\in \Omega\colon (\forall \tau\in[0,1])\;\,
\Phi(t_0+\tau(t-t_0),t_0,y_0)\in U\}
\]
and is given by $\Fl^U=\Fl|_{\Omega^U}\colon \Omega^U\to U$.
\item[(b)]
If $\Omega$ is open in $J\times J\times M$ and $\Fl$ is continuous,
then $\Omega^U$ is open in $J\times J\times U$.
\end{description}
\end{exer}

\begin{exer}\label{exc-flow-related}
Given a $C^1$-map $g\colon M\to N$ between $C^1$-manifolds $M$ and $N$,
let\linebreak
$X\colon M\to TM$ and $Y\colon N\to TN$
be continuous vector fields which are $g$-related in the sense that $Y\circ g=Tg\circ X$.
Let $\gamma\colon I\to M$ be a $C^1$-map on a non-degenerate interval $I\sub \R$.
Show:
\begin{description}[(a)]
\item[(a)]
If $\gamma$ is a solution to the differential equation $\dot{y}(t)=X(y(t))$,
then $g\circ \gamma$ solves $\dot{y}(t)=Y(y(t))$.
\item[(b)]
If $T_pg\colon T_pM\to T_{g(p)}N$ is injective for all $p\in M$
and $g\circ\gamma$ is a solution to the differential equation
$\dot{y}(t)=Y(y(t))$, then $\gamma$ solves $\dot{y}(t)=X(y(t))$.
\item[(c)]
If $g$ is a $C^1$-diffeomorphism, then $\dot{y}(t)=X(y(t))$
admits local $C^1$-flows if and only if $\dot{y}(t)=Y(y(t))$
admits local $C^1$-flows.
In this case, we have $\Omega^Y=(\id_\R\times g)(\Omega^X)$
and $\Fl^Y_t(g(y_0))=g(\Fl^X_t(y_0))$ for all $(t,y_0)\in\Omega^X$.
\item[(d)]
Formulate analogs for time-dependent vector fields.
\end{description}
\end{exer}

\begin{exer}\label{exc-auton-flow-bij}
Let $M$ be a $C^1$-manifold (possibly with boundary) modeled on a locally convex space
and $X\colon M \to TM$ be a continuous vector field.
Consider $\Fl^X\colon\Omega^X\to M$ and $\Omega^X_t$ as in
Definition~\ref{flow-auton}.
Show that $\Fl^X_t(\Omega_t^X)=\Omega_{-t}^X$ for all $t\in\R$.
Also show that $\Fl^X_t\colon \Omega_t^X\to\Omega_{-t}^X$
is bijective with inverse $(\Fl^X_t)^{-1}=\Fl^X_{-t}$.
\end{exer}

\begin{exer}\label{exc-details-ode-mfd}
Fill in details in Remark~\ref{rems-ode-md} and the proofs of Lemma~\ref{flo-cons}
and Corollary~\ref{pardep-ode-ban-m}.
\end{exer}

\begin{exer}\label{exc-loc-lip-indep}
Let $M$ and $N$ be $C^1$-manifolds (possibly with boundary)
modeled on normed space $(E,\|\cdot\|_E)$ and $(F,\|\cdot\|_F)$,
respectively.
Let $X$ be a topological space, $W\sub X\times M$ be a subset
and $f\colon W\to N$ be a continuous function
which satisfies a local Lipschitz condition in its second argument.
Let $Y\sub X$ be an open subset,
$\kappa\colon P_\kappa\to Q_\kappa\sub E$ be a chart of~$M$
and $\theta\colon P_\theta\to Q_\theta\sub F$
be a chart of~$N$ such that
$f(W\cap (Y\times P_\kappa))\sub P_\theta$.
Show that the map
\[
h\colon (\id_Y\times \kappa)(W\cap (Y\times P_\kappa))\to F,\;\, (a,b)\mto \theta\big(f\big(a,\kappa^{-1}(b)
\big)\big)
\]
satisfies a local Lipschitz condition in its second argument.\\[1mm]
[Insert charts $\phi$ and $\psi$ as in Definition~\ref{def-lipcond-mfd} and use Lemma~\ref{lipviaprime}(a)
to see that the transition maps between the charts are locally Lipschitz.]
\end{exer}

\begin{exer}\label{exc-lie-der-mfd}
Let $X\colon M\to TM$
be a continuous vector field on a $C^1$-manifold~$M$
such that $\dot{y}(t)=X(y(t))$ admits local $C^1$-flows.
In (a), $M$ may have a boundary.
\begin{description}[(a)]
\item[(a)]
Show that
\[
\frac{d}{dt}\Big|_{t=0}f(\Fl^X_t(p))=df(X(p))
\]
for each $C^1$-function $f\colon M\to F$ to a locally convex space~$F$
and $p\in M$.
\item[(b)]
If $M$ and $X$ are smooth and $X$ admits local $C^2$-flows, then
\[
\frac{d}{dt}\Big|_{t=0}d\Fl_{-t}^X\big(Y\big(\Fl^X_t(p)\big)\big)=[X,Y](p)
\]
for each smooth vector field
$Y\colon M\to TM$ and $p\in M$
(where $[X,Y]$ is the Lie bracket of the smooth vector fields $X$ and $Y$,
as in Theorem~\ref{Liebrvf}).\\[1mm]
[Use Exercises~\ref{exc-flo-opensub} and \ref{exc-flow-related}(c)
to reduce to Exercise~\ref{exc-lie-der-loc}.]
\end{description}
\end{exer}

\end{small}
\section{Differential calculus on metrizable spaces}\label{secmetrcalc}
In this section, we provide certain tools
which simplify
differential calculus
on metrizable locally convex spaces.
We prove the following facts:
\begin{thm}\label{mayn1}
Let $E$ be a metrizable real locally convex space,
$U\sub E$ be open
and $f\colon U\to X$ be a map
to a topological space~$X$.
Then $f$ is continuous 
if and only if
$f\circ \gamma\colon \R\to X$ is continuous for each smooth curve
$\gamma\colon \R\to U$.
\end{thm}
\begin{thm}\label{mayn2}
Let $E$ be a metrizable real locally convex space,
$U\sub E$ be open
and $f\colon U\to F$ be a map
to a locally
convex space~$F$. Let~$k\in \N_0$.
Then~$f$
is~$C^k$ 
if and only if $f\circ \gamma\colon \R^{k+1}\to F$
is~$C^k$ for each smooth map $\gamma\colon \R^{k+1}\to U$.
\end{thm}
\begin{rem}\label{linkconven}
Theorem~\ref{mayn1}
can be interpreted as follows:
If $E$ is a metrizable locally convex space
and $U\sub E$ an open subset, then
the topology induced by~$E$ on~$U$
coincides with the final topology
on~$U$ with respect to the set $C^\infty(\R,U)$
of smooth curves $\gamma\colon \R\to U$.
\end{rem}
Before we establish these theorems, let us prove
an easier result first:
\begin{prop}\label{contarcs}
Let $E$ be a metrizable
locally convex space, $U\sub E$ be open
and $f\colon U\to X$ be a map to a topological space~$X$.
Then $f$ is continuous if and only if $f\circ \gamma\colon \R\to
X$ is continuous for each continuous curve $\gamma\colon \R\to U$.
\end{prop}
\begin{prf}
The necessity of the condition is obvious.
To see that the condition is also sufficient,
we argue by contraposition. Thus, suppose
that $f$ is discontinuous. The space~$E$ being metrizable,
we find a point $x\in U$ and a sequence
$(x_n)_{n\in \N}$ in $U$ converging to~$x$
such that $(f(x_n))_{n\in \N}$ does not converge to~$f(x)$.
Let $V\sub U$ be a convex neighborhood
of~$x$; after passage to a subsequence,
we may assume that $x_n\in V$ for all $n\in \N$.
We now define $\gamma\colon \R\to U$
via
\[
\gamma(t)\; :=\; 
\left\{
\begin{array}{cl}
x \; & \;\mbox{if $\,t\leq 0$}\\
x_1 &\;\mbox{if $\,t\geq 1$}\\
x_{n+1}+ {\textstyle \frac{t-\frac{1}{n+1}}{\frac{1}{n}-\frac{1}{n+1}}
(x_n-x_{n+1})} & \;\mbox{if $\,t\in
[\frac{1}{n+1},\frac{1}{n}]$}
\end{array}
\right.
\]
and observe that $\gamma(\frac{1}{n}) = x_n$.
Note that the image of this map is the `infinite polygon'
obtained by drawing a line segment from $x_n$ to $x_{n+1}$
for each~$n$, together with $\{x\}$.
Because 
$\{x\}\cup\{x_n\colon n\in \N\}\sub
V$ and $V$ is convex, we deduce that $\,\im\,\gamma\sub V\sub U$.
Then $\gamma$ is apparently
continuous on $\R\setminus\{0\}$. It is also
continuous at~$0$:
In fact, given a convex neighborhood
$W\sub U$ of~$x$,
there exists $N\in \N$
such that $x_n \in W$ for all $n\geq N$.
For each $t\in \,]0,\frac{1}{N}]$, there exists $n\geq N$
such that $t\in \;]\frac{1}{n+1},\frac{1}{n}]$.
Then $\gamma(t)=\tau x_{n+1}+(1-\tau)x_n$
for some $\tau\in [0,1]$ and hence
$\gamma(t)\in W$,
using that $W$ is convex and $x_n,x_{n+1}\in W$.
Thus $\gamma$ is continuous.
It remains to note that $\gamma(\frac{1}{n})=x_n$,
whence $(f\circ \gamma)(\frac{1}{n})$ does not converge to $(f\circ\gamma)(0)$
although $\frac{1}{n}\to 0$. Hence $f\circ
\gamma$ is not continuous.
\end{prf}
To be able to join the points $x_n$
of a converging sequence $(x_n)_{n\in \N}$
not only by a continuous curve, but a smooth
curve, we need to assume that the sequence $(x_n)_{n\in \N}$
converges sufficiently fast.
\begin{defn}
A sequence $(x_n)_{n\in \N}$
in a locally convex space~$E$ is said
to \emph{converge fast}
to $x\in E$ if
${\displaystyle \lim_{n\to \infty}n^k\cdot (x_n-x)=0}$,
for each $k\in \N_0$.
\end{defn}
\begin{lem}[Special Curve Lemma]\label{specialcurve}
Let $E$ be a real
locally convex space, $x\in E$,
and $(x_n)_{n\in \N}$ be a sequence in~$E$ which converges
fast to~$x$. Then there exists a smooth curve
$\gamma\colon \R\to E$
with image $\gamma(\R)$ contained in
$\{x\}\cup\{tx_n+(1-t)x_{n+1}\colon n\in \N,\, t\in [0,1]\}$
such that $\gamma(\frac{1}{n})=x_n$ for each $n\in\N$ and
$\gamma(0)=x$.
\end{lem}
\begin{prf}
We choose
a smooth map $\tau\colon \R\to \R$
such that $\im\,\tau\sub [0,1]$, $\tau(t)=0$ if $t\leq \frac{1}{3}$
and
$\tau(t)=1$ if $t\geq \frac{2}{3}$.
We now define a map $\gamma\colon \R\to U$
via
\[
\gamma(t) \;:=\;
\left\{
\begin{array}{cl}
x & \; \mbox{if $\;t\leq 0$}\\
x_1 &\; \mbox{if $\;t\geq 1$}\\
x_{n+1}+ {\textstyle \tau\left(
\frac{t-\frac{1}{n+1}}{\frac{1}{n}-\frac{1}{n+1}}\right)
(x_n-x_{n+1})} & \; \mbox{if $\;
t\in [\frac{1}{n+1},\frac{1}{n}]$.}
\end{array}
\right.
\]
Then $\gamma$ is either constant
or given by a closed formula defining a smooth function
on an open neighborhood of each given $t\in \R\setminus\{0\}$,
whence $\gamma|_{\R\setminus\{0\}}$ is smooth.
To see that $\gamma$ is continuous at~$0$
and hence continuous, we can argue as
in the proof of Proposition~\ref{contarcs}.
Let $k\in \N$ now;
we want to show that $\gamma^{(k)}(t)\to 0$
as $t\to 0^+$ (and hence also if $0\not=t\to 0$).
To this end, note that for $t\in\;]0,1]$,
say $t\in [\frac{1}{n+1}, \frac{1}{n}]$,
we have
\begin{eqnarray*}
\gamma^{(k)} (t)
&=&
{\textstyle \tau^{(k)}
\left(
\frac{t-\frac{1}{n+1}}{\frac{1}{n}-\frac{1}{n+1}}
\right)
\left(\frac{1}{\frac{1}{n}-\frac{1}{n+1}}\right)^k
(x_n-x_{n+1})}\\
&=&
{\textstyle \tau^{(k)}
\left(
\frac{t-\frac{1}{n+1}}{\frac{1}{n}-\frac{1}{n+1}}
\right)
\left( n(n+1) \right)^k
(x_n-x_{n+1})}\,,
\end{eqnarray*}
by the Chain Rule.
To see that this converges to $0$, let $\|\cdot\|_p\colon E\to \R$
be a continuous seminorm on~$E$.
The set $\tau^{(k)}(\R)=\tau^{(k)}([\frac{1}{3},\frac{2}{3}])$
is compact and 
hence contained in $[-M,M]$
for some $M\geq 0$.
We now estimate for $n \geq 2$: 
\begin{eqnarray*}
\|\gamma^{(k)}(t)\|_p &=&
{\textstyle \left|\tau^{(k)}\left(
\frac{t-\frac{1}{n+1}}{\frac{1}{n}-\frac{1}{n+1}}\right)\right|
\left(n(n+1)\right)^k
\|x_n-x_{n+1}\|_p}\\
&\leq& M
\left(n(n+1)\right)^k
\|x_n-x_{n+1}\|_p\\
&=& M
\left(n(n+1)\right)^k
\|(x_n-x)-(x_{n+1}-x)\|_p\\
&\leq &
M\left(n(n+1)\right)^k
\|x_n-x\|_p +
M\left(n(n+1)\right)^k
\| x_{n+1}-x\|_p\\
&\leq &
Mn^{3k}
\|x_n-x\|_p +
M(n+1)^{2k}
\|x_{n+1}-x\|_p \\
&=& M\| n^{3k}(x_n-x)\|_p +
M\|(n+1)^{2k}(x_{n+1}-x)\|_p \,.
\end{eqnarray*}
Hence
$\lim_{t\to 0}\|\gamma^{(k)}(t)\|_p=0$,
using that $n=n(t)\to\infty$;
we exploit here that $x_n\to x$ fast.

By the preceding,
$\gamma$ is continuous, $\gamma|_{\R\setminus\{0\}}$
is smooth, and $\gamma^{(k)}$ extends from
$\R\setminus\{0\}$ to a continuous map $\R\to\R$
by $0\mto 0$. Using Lemma~\ref{C1glueing},
we see that
$\gamma$ is $C^1$ with $\gamma'(0)=0$.
Repeating the argument, we see that $\gamma$
is $C^k$ with $\gamma^{(k)}(0)=0$,
for each $k\in \N$.
\end{prf}
\begin{lem}\label{picknose}
Let $(x_n)_{n\in \N}$ be a convergent
sequence in a metrizable locally convex space~$E$,
with limit $x\in E$. Then $(x_n)_{n\in \N}$ has a subsequence
$(x_{n_j})_{j\in\N}$ which converges fast to~$x$.
\end{lem}
\begin{prf}
Since $E$ is locally convex and metrizable,
there exists a sequence of continuous seminorms 
$\|\cdot\|_1\leq \|\cdot\|_2\leq \cdots$
on $E$ defining its locally convex vector
topology.
As $x_n\to x$,
for each $j\in \N$ we find $n_j\in \N$ such that
\begin{equation}\label{eqnjs}
\|x_n-x\|_j\leq \frac{1}{j^j}\quad \mbox{for all $n\geq n_j$.}
\end{equation}
We may assume that $n_1<n_2<\cdots$.
We claim that $x_{n_j}$ converges fast to~$x$ as $j\to\infty$.
To see this, fix $k\in \N_0$. In order that $\lim_{j\to\infty}
j^k(x_{n_j}-x)=0$, we have to show that $\lim_{j\to\infty}
\|j^k(x_{n_j}-x)\|_N=0$ for each $N\in \N$.
But indeed, for each $j\geq N$, using (\ref{eqnjs})
we find that
$\|j^k(x_{n_j}-x)\|_N\leq \|j^k(x_{n_j}-x)\|_j
=j^k\|x_{n_j}-x\|_j \leq j^k\frac{1}{j^j}=\frac{1}{j^{j-k}}\to 0$
as $j\to \infty$.
\end{prf}
\noindent
{\em Proof of Theorem}~\ref{mayn1}.
If $f$ is continuous, then apparently also $f\circ \gamma$
is continuous for each smooth curve $\gamma\colon \R\to U$.

To see the converse, we argue by contraposition.
Thus, assume that $f\colon U\to X$ is discontinuous
at some point~$x\in U$.
Let $V\sub U$ be a convex neighborhood of~$x$.
Then there exist
a neighborhood
$W\sub F$ of~$f(x)$ and
a sequence $(x_n)_{n\in \N}$ in~$V$
such that $x_n\to x$
but $f(x_n)\not\in W$ for each
$n\in \N$.
By Lemma~\ref{picknose},
after passing to a subsequence
we may assume that $x_n\to x$
fast.
Hence the Special Curve Lemma (Lemma~\ref{specialcurve})
provides a smooth curve $\gamma\colon \R\to E$
such that $\gamma(\frac{1}{n})=x_n$
for all $n\in \N$, $\gamma(0)=x$ and such that $\,\im\,\gamma$
is contained in the convex hull
$\,\conv(\{x\}\cup\{x_n\colon n\in \N\})$,
which is a subset of~$V$ and hence of~$U$.
Thus $\,\im\,\gamma\sub U$, enabling us to consider
the composition $f\circ \gamma$.
Since $(f\circ \gamma)(0)=f(x)\in W$
but $(f\circ \gamma)(\frac{1}{n})=f(x_n)\not\in W$
for all $n\in \N$, where $\frac{1}{n}\to 0$,
we see that $f\circ \gamma$ is discontinuous.\vspace{2mm}\qed

The following two lemmas will help us
to prove Theorem~\ref{mayn2}:
\begin{lem}\label{alsoopen}
Let~$k\in \N$, $r\in \N_0\cup\{\infty\}$,
$E$ be a metrizable real locally convex space,
$U\sub E$ be open
and $f\colon U\to F$ be a map to a locally
convex space~$F$
such that $f\circ \gamma\colon \R^k\to F$
is~$C^r$ for each smooth map $\gamma\colon \R^k\to U$.
Then
$f\circ \gamma\colon W \to F$
is~$C^r$ for each open subset $W\sub \R^k$
and smooth map $\gamma\colon W\to U$.
\end{lem}
\begin{prf}
Let $\gamma\colon W\to U$ be a smooth map, where
$W\sub \R^k$ is open. In order that $f\circ \gamma$
be $C^r$, it suffices to show that
it is $C^r$ on some open neighborhood
of each given point $x\in W$.
Thus, let $x\in W$.
Then
the open ball~$B_\ve(x)$
with respect to $\|\cdot\|_\infty$
is contained in~$W$,
for some $\ve>0$.
There is a smooth function
$h\colon \R^k\to W$ such that $\,\im\, h\sub B_\ve(x)$
and $h(y)=y$ for all $y\in B_{\frac{\ve}{2}}(x)$.
Now $\gamma\circ h\colon\R^k\to U$ being smooth,
by hypothesis $f\circ (\gamma\circ h)$ is~$C^r$.
Thus $(f\circ \gamma\circ h)|_{B_{\frac{\ve}{2}}(x)}=
f\circ \gamma|_{B_{\frac{\ve}{2}}(x)}$
is~$C^r$,
which completes the proof.
\end{prf}
\begin{lem}\label{lowerdm}
Let~$k\in \N$, $r\in \N_0\cup\{\infty\}$,
$E$ be a metrizable real locally convex space,
$U\sub E$ be open
and $f\colon U\to F$ be a map to a locally
convex space~$F$
such that $f\circ \gamma\colon \R^k\to F$
is~$C^r$ for each smooth map $\gamma\colon \R^k\to U$.
Then
$f\circ \gamma\colon \R^j\to F$
is~$C^r$ for each $j<k$ and
each smooth map $\gamma\colon \R^j\to U$.
\end{lem}
\begin{prf}
Given $j<k$ and a smooth map $\gamma\colon \R^j\to U$,
we define a smooth map $\eta\colon \R^k\to U$
via
$\eta(t_1,t_2,\ldots, t_k):=\gamma(t_1,t_2,\ldots, t_j)$
for $(t_1,\ldots, t_k)\in \R^k$.
Then $f\circ \eta$ is of class $C^r$ by hypothesis.
Since $f\circ \gamma=(f\circ \eta)\circ \lambda$,
where
\[
\lambda\colon \R^j\to \R^k,\quad
\lambda(t_1,\ldots, t_j)\, :=\, (t_1,\ldots, t_j,0,\ldots, 0)
\]
is continuous linear and thus~$C^r$,
we see that also $f\circ\gamma$ is~$C^r$.
\end{prf}
We are now well-prepared to prove Theorem~\ref{mayn2}.\\[3mm]
{\em Proof of Theorem}~\ref{mayn2}.
The proof is by induction on
$k\in \N_0$. The case $k=0$ being covered by Theorem~\ref{mayn1},
let us assume that $k\in \N$ now,
and assume that the assertion of Theorem~\ref{mayn2} is correct if
$k$ is replaced with~$k-1$. Let $f\colon U\to F$
be a map such that $f\circ \gamma$ is~$C^k$
for all smooth maps $\gamma\colon \R^{k+1}\to U$.
Then $f\circ \gamma$ is~$C^k$ and hence~$C^{k-1}$
for all smooth maps $\gamma\colon \R^k\to U$
(Lemma~\ref{lowerdm}),
whence $f$ is~$C^{k-1}$, by induction.

Let $x\in U$ and $y_1,\ldots, y_k\in E$;
we claim that $d^{\,(k)}f(x,y_1,\ldots,y_k)$ exists.
To see this, choose $r>0$ such that
$\gamma(t_1,\ldots, t_k):=x+t_1y_1+\cdots+t_ky_k\in U$
for all $(t_1,\ldots,t_k)\in \; ]{-r},r[^k=:W$.
By Lemmas~\ref{alsoopen} and \ref{lowerdm}, 
the map $f\circ \gamma\colon W\to F$ is~$C^k$,
where $\gamma\colon W\to U$ is defined as just described.
Let $e_i$ be the $i$-th standard basis vector
of~$\R^k$, for $i\in\{1,\ldots, k\}$.
Then, for each $t=(t_1,\ldots,t_k)\in W$, we calculate
\begin{eqnarray*}
d(f\circ \gamma)(t,e_1)&=&{ \frac{d}{ds}}\Big|_{s=0}f(x+t_1y_1+\cdots
+t_ky_k+sy_1)\\
&=&
df(x+t_1y_1+\cdots
+t_ky_k, y_1)=df(\gamma(t),y_1)\,,
\end{eqnarray*}
and inductively, by the same argument,
\[
d^{\,(j)}(f\circ \gamma)(t,e_1,\ldots, e_j)
\,=\, d^{\,(j)}f(\gamma(t),y_1,\ldots, y_j)
\quad\mbox{for all $j\in\{1,\ldots, k\}$.}
\]
In particular, we have shown that the limit
defining $d^{\,(k)}f(x,y_1,\ldots,y_k)$
exists; it is given by
$d^{\,(k)}f(x,y_1,\ldots, y_k)=d^{\,(k)}(f\circ \gamma)
(0,e_1,\ldots, e_k)$.

It remains to show that $d^{\,(k)}f\colon U\times E^k\to F$
is continuous. The set $U\times E^k$
being open in~$E^{k+1}$ and~$E^{k+1}$ being metrizable,
we only need to show that $d^{\,(k)}f\circ \gamma\colon\R\to F$
is continuous for each smooth curve $\gamma\colon\R\to
U\times E^k$. This will be the case
if we can show that $d^{\,(k)}f\circ \gamma$ is
continuous on some neighborhood
of each element $t_0\in \R$.
To verify this property, let us write
$\gamma=(\gamma_0,\gamma_1,\ldots,\gamma_k)\colon \R\to
U\times \R^k$.
By continuity, there exists
an open interval $I\sub \R$
containing~$t_0$
and
$\ve>0$ such that
\[
\eta(t,t_1,\ldots,t_k)\, :=\, \gamma_0(t)+t_1\gamma_1(t)+\cdots
+t_k\gamma_k(t)\, \in \, U
\]
for all $(t,t_1,\ldots,t_k)\in I\times \, ]{-\ve},
\ve[^k=:Y$.
Note that $\eta\colon Y\to U$ is a smooth map
on the open subset $Y\sub \R^{k+1}$.
Hence $f\circ \eta$ is~$C^k$
(see Lemma~\ref{alsoopen}). Labelling
the standard basis vectors
of $\R^{k+1}$ by $0,1,\ldots,k$ now,
we obtain
\begin{eqnarray*}
d(f\circ \eta)((t,t_1,\ldots,t_k), e_1)
&=&\!{\frac{d}{ds}}\Big|_{s=0}f(
\gamma_0(t)+t_1\gamma_1(t)+\cdots
+t_k\gamma_k(t) + s\gamma_1(t))\\
&=&df(\gamma_0(t)+t_1\gamma_1(t)+\cdots
+t_k\gamma_k(t),\gamma_1(t))\\
&=&
df(\eta(t,t_1,\ldots, t_k),\gamma_1(t))
\end{eqnarray*}
for all $(t,t_1,\ldots, t_k)\in Y$.
Repeating the argument, we find that
\[
d^{\,(j)}(f\circ \eta)((t,t_1,\ldots,t_k),\, e_1,\ldots, e_j)
=d^{\,(j)}f(\eta(t,t_1,\ldots, t_k),\gamma_1(t),\ldots, \gamma_j(t))
\]
for all $j=1,\ldots, k$. In particular,
\begin{eqnarray*}
d^{\, (k)}(f\circ \eta)((t,0,\ldots,0),\, e_1,\ldots, e_k)
&=&d^{\,(k)}f(\eta(t,0,\ldots, 0),\, \gamma_1(t),\ldots,\gamma_k(t))\\
&=&d^{\,(k)}f(\gamma_0(t),\gamma_1(t),\ldots, \gamma_k(t))\\
&=& (d^{\,(k)}f\circ \gamma)(t)\,,
\end{eqnarray*}
for all $t\in I$. The map
$d^{\, (k)}
(f\circ\eta)$ being continuous,
the preceding formula shows that
also $(d^{\, (k)}f)\circ \gamma|_I$ is continuous,
which completes the proof.\qed
\begin{rem}
Lemma~\ref{specialcurve}
and the proof of
Theorem~\ref{mayn1}
have been adapted from \cite[\S2]{KM97}.
Theorem~\ref{mayn2} is a special
case of \cite[Thm.~12.4]{BGN04}.
\end{rem}
Theorem~\ref{mayn2}
implies that $f$ is smooth if and only
if $f\circ \gamma$ is smooth
for each $k\in \N$ and smooth map
$\gamma\colon \R^k\to U$.
We mention
a closely related fact:
\begin{rem}
Let $f\colon E\supseteq
U\to F$ be a map
from an open subset
of a metrizable locally convex space~$E$ to a locally
convex space~$F$.
It can be shown that
$f$ is $C^\infty$
if and only if $f\circ \gamma\colon \R\to F$
is smooth for each smooth curve $\gamma\colon \R\to U$
(cf.\ Theorems 4.11 and 12.8 in \cite{KM97}; see also \cite[Prop.~E.3]{Gl04b}).
For mappings $f\colon \R^n\to \R$,
the preceding fact is known
as \emph{Boman's Theorem}~\cite{Bm67}.
\end{rem}
\begin{rem}\label{two-var-enough}
Theorem~\ref{mayn2} can be strengthened as follows:\\[2.3mm]
\emph{Let $E$ be a metrizable real locally convex space,
$k\in\N_0\cup\{\infty\}$
and $f\colon U\to F$
be a map from an open subset $U\sub E$ to a locally convex space~$F$.
Then $f$ is $C^k$ if and only if $f\circ \gamma\colon\R^2\to F$ is $C^k$
for each smooth map $\gamma\colon \R^2\to U$.}\\[2mm]
To see the sufficiency for $k\geq 1$,
let $\gamma_1,\gamma_2\colon \,]{-\ve},\ve[\,\to E$
be $C^\infty$-curves and
\[
\gamma(s,t):=\gamma_1(t)+s\gamma_2(t)\in U
\]
for all $s,t\in\,]{-\ve },\ve[$.
Then $f\circ\gamma$ is $C^k$ and, for all $j\in\N$ such that $j\leq k$,
\[
\frac{\partial^j(f\circ \gamma)}{\partial s^j}(s,t)=\delta^{(j)}_{\gamma(s,t)}f(\gamma_2(t)).
\]
Setting $s=0$, we see that
the G\^{a}teaux differential $\delta^{(j)}_{\gamma_1(t)}f(\gamma_2(t))$
exists for all $t\in \,]{-\ve},\ve[$, entailing that $\delta^{(j)}f\colon U\times E\to F$
exists and is continuous (by Theorem~\ref{mayn1}), since
\[
(\delta^{(j)}f)\circ(\gamma_1,\gamma_2)=\frac{\partial^j(f\circ \gamma)}{\partial s^j})(0,\cdot)
\]
is continuous. Hence $f$ is $C^k$, by Exercise~\ref{exc-gat-Ck}.
\end{rem}

\begin{small}
\subsection*{Exercises for Section~\ref{secmetrcalc}}

\begin{exer}
In the proof of Lemma~\ref{alsoopen},
a certain smooth mapping~$h$ was used.
Show that such a map $h$ exists.
\end{exer}

\begin{exer}
Verify that Theorem~\ref{mayn1} and its proof remain valid if $U\sub E$
is not necessarily open, but a locally convex subset with dense interior.
\end{exer}

\end{small}
\section{Differential calculus on Silva spaces}\label{sec-silva-calc}
Some infinite-dimensional Lie grous of interest
are modeled on Silva spaces,
like direct limits of finite-dimensional Lie groups (see Chapter~\ref{ch:dirlim}) 
or the Lie group $\Diff^\omega(M)$ of real analytic diffeomorphisms
of a compact real analytic manifold~$M$ (see \cite{KM97} or \cite{DS15},
cf.\ also \cite{Les82, Les83}).
In this section, we provide tools of differential calculus which
are useful for the construction of the Lie group structure on such groups.
The required background on Silva spaces, direct limits of topological
spaces and direct limits of locally convex spaces can be found
in Appendices \ref{app-basic-DL} and \ref{sec-appDLvec}.
\begin{prop}\label{diff-silva}
Let $E_1\sub E_2\sub\cdots$ be a direct sequence of Banach spaces $(E_n,\|\cdot\|_n)$
such that all inclusion maps $E_n\to E_{n+1}$ are compact operators.
Let $E:=\bigcup_{n\in\N}E_n$ be the corresponding Silva space,
endowed with the locally convex direct limit topology~$\cO$.
Let $U_1\sub U_2\sub\cdots$ be a sequence of open subsets~$U_n$ in $(E_n,\|\cdot\|_n)$.
Then $U:=\bigcup_{n\in\N}U_n$ is open in $(E,\cO)$. If $k\in\N_0\cup\{\infty\}$
and $f\colon U\to F$ is a function to a locally convex space~$F$
such that
\[
f|_{U_n}\colon U_n\to F\quad\mbox{is $C^k$ on $U_n\sub E_n$
for each $n\in\N$,}
\]
then $f$ is $C^k$.
\end{prop}
The proof uses a lemma.
\begin{lem}\label{diff-steps}
Let $E_1\sub E_2\sub \cdots$ be a direct sequence of locally convex
spaces and $\cO$ be a locally convex vector topology on $E:=\bigcup_{n\in\N}E_n$
making each inclusion map $E_n\to E$ continuous.
Let $U_1\sub U_2\sub\cdots$ be an ascending sequence of open subsets
$U_n\sub E_n$ such that $U:=\bigcup_{n\in\N}U_n$ is open in $(E,\cO)$.
If $k\in\N_0\cup\{\infty\}$ and $f\colon U\to F$ is a function to a locally convex space~$F$ such that
$f|_{U_n}\colon U_n\to F$ is a $C^k$-function on the open subset $U_n\sub E_n$
for each $n\in\N$,
then the iterated directional derivative
\[
d^{\,(j)}f(x,y_1,\ldots, y_j)
\]
exists in~$F$ for all $j\in\N$ with $j\leq k$, $x\in U$ and $y_1,\ldots, y_j\in E$;
it satisfies
\begin{equation}\label{der-pieces}
d^{\,(j)}f|_{U_n\times (E_n)^j}=d^{\,(j)}(f|_{U_n})\quad\mbox{for all $n\in\N$.}
\end{equation}
\end{lem}
\begin{prf}
If $x\in U_n$, $y_1,\ldots, y_j\in E_n$, we have
\begin{eqnarray*}
df(x,y_1) &=&
\lim_{t\to 0} \frac{1}{t} (f(x+ty_1)-f(x))=\lim_{t\to0}
\frac{1}{t} (f|_{U_n}(x+ty_1)-f|_{U_n}(x))\\
&=& d(f|_{U_n})(x,y_1)
\end{eqnarray*}
and recursively
\begin{eqnarray*}
d^{\,(j)} f(x,y_1,\ldots,y_j)
&=& \frac{d}{dt} \Big|_{t=0} d^{\,(j-1)}f(x+ty_j ,y_1,\ldots, y_{j-1})\\
&=& \frac{d}{dt} \Big|_{t=0} d^{\,(j-1)}(f|_{U_n})(x+ty_j,y_1,\ldots, y_{j-1})\\
&=&
d^{\,(j)}(f|_{U_n})(x,y_1,\ldots, y_j)
\end{eqnarray*}
for all integers $0<j\leq k$.
\end{prf}
{\em Proof of Proposition}~\ref{diff-silva}.
By Proposition~\ref{silvahaveDL}(a),
the locally convex direct limit topology $\cO$ coincides with the topology making~$E$ the direct limit of the topological spaces
$E_1\sub E_2\sub\cdots$. Hence~$U$ is open in $(E,\cO)$, by Lemma~\ref{basicDL}(a),
and~$f$ is continuous as $f|_{U_n}$ is continuous for each $n$ (see Lemmas~\ref{basicDL}(b)
and \ref{ctsonsteps}(ii)).
By Lemma~\ref{diff-steps},
the $j$-fold iterated directional derivatives of~$f$ exist,
and are given by (\ref{der-pieces}) for all $j\in\N$ such that $j\leq k$.
Now $U\times E^j$ is an open subset of $E^{j+1}$
which is the locally convex (and topological) direct limit of
the directed sequence
\[
(E_1)^{j+1}\sub (E_2)^{j+1}\sub\cdots
\]
of Banach spaces and compact inclusion operators.
Since $d^{\,(j)}f|_{U_n\times (E_n)^j}$ is continuous for each $n\in\N$ by~(\ref{der-pieces}),
Lemmas~\ref{basicDL}(b) and \ref{ctsonsteps}(ii) show
that $d^{\,(j)}f$ is continuous. Hence $f$ is~$C^k$.\qed
\section{Differential calculus on locally convex direct sums}\label{sec-sums}
We compile some results concerning continuity, differentiability properties and analyticity
for a well-behaved class of mappings between locally convex direct sums.
The results will be used in connection with spaces of compactly supported
smooth sections in vector bundles over a $\sigma$-compact finite-dimensional base manifold~$M$,
and also in the construction of a Lie group structure on the group $\Diff(M)$
of all $C^\infty$-diffeomorphisms of~$M$.
\begin{prop}\label{diff-sum}
Let $(E_n)_{n\in\N}$ and $(F_n)_{n\in\N}$
be sequences of locally convex topological $\K$-vector spaces.
Let $k\in\N_0\cup\{\infty\}$
and $(f_n)_{n\in\N}$ be a sequence of $C^k_\K$-maps $f_n\colon U_n\to F_n$
on open subsets $U_n\sub E_n$ such that $0\in U_n$ and $f_n(0)=0$
for all but finitely many $n\in\N$. Then $\bigoplus_{n\in\N}U_n$
is open in the locally convex direct sum $\bigoplus_{n\in \N}E_n$
and
\[
f:=\bigoplus_{n\in\N}f_n\colon \bigoplus_{n\in\N}U_n\to \bigoplus_{n\in\N}F_n,
\;\, (x_n)_{n\in\N}\mto (f_n(x_n))_{n\in\N}
\]
is a $C^k_\K$-map. If $\K=\R$ and each $f_n$ is real analytic,
then $f$ is real analytic.
\end{prop}
\begin{prf}
We know from Remark~\ref{firstremsums}(a) that $U:=\bigoplus_{n\in\N}U_n$ is open in
$E:=\bigoplus_{n\in\N}E_n$.
To see that~$f$ is $C^k_\K$, we may assume that $k<\infty$ and proceed by induction.

If $k=0$, we have to show that $f$ is continuous, i.e., continuous at each $x=(x_n)_{n\in\N}
\in U$. Let $V$ be an open neighborhood of $f(x)$ in $F:=\bigoplus_{n\in\N}F_n$;
after shrinking~$V$, we may assume that $V=\bigoplus_{n\in\N}V_n$
with $V_n$ an open neighborhood of $f_n(x_n)$ in~$F_n$
(see Remark~\ref{firstremsums}(a)).
Since $f_n^{-1}(V_n)$ is an open neighborhood of~$x_n$ in~$E_n$,
we deduce that
\[
f^{-1}(V)=\bigoplus_{n\in\N}f_n^{-1}(V_n)
\]
is an open neighborhood of~$x$ in~$E$. Hence $f$ is continuous at~$x$.

If $k>0$, then for $x=(x_n)_{n\in\N}\in U$ and $y=(y_n)_{n\in\N}\in E$,
we find $N\in\N$ such that $x,y\in E_1\times\cdots\times E_N$
and
\[
(\forall n>N) \quad 0\in U_n\quad\mbox{and}\quad f_n(0)=0.
\]
For $t\in\K\setminus\{0\}$ close to~$0$,
we have $x+ty\in U_1\times\cdots\times U_N$
and see that
\begin{eqnarray*}
\frac{1}{t}(f(x+ty)-f(x))&=& 
\Big(\frac{1}{t}(f_n(x_n+ty_n)-f_n(x_n))\Big)_{n\in\{1,\ldots, N\}}\\
&\to & (df_n(x_n,y_n))_{n\in\{1,\ldots, N\}}=(df_n(x_n,y_n))_{n\in\N}
\end{eqnarray*}
as $t\to0$ in $F_1\times\cdots \times F_N$ and hence in~$F$
(identifying the former with a vector subspace of~$F$ as in Remark~\ref{firstremsums}(c)).
Thus $df(x,y)$ exists and
\begin{equation}\label{sumindu}
df =\Big(\bigoplus_{n\in\N}df_n\Big)\circ \Phi|_{U\times E},
\end{equation}
where $\Phi\colon E\times E\to \bigoplus_{n\in\N} (E_n\times E_n)$
is the isomorphism of topological vector spaces sending $(x,y)$
with $x=(x_n)_{n\in\N}$, $y=(y_n)_{n\in\N}$ to
$\Phi(x,y)=(x_n,y_n)_{n\in\N}$ 
(see Lemma~\ref{firstlasum}(c)).
Now $\oplus_{n\in\N}df_n$ (and hence $df$, by (\ref{sumindu})) is $C^{k-1}_\K$ by induction and hence
continuous, whence $f$ is $C^1_\K$ with $df$ a $C^{k-1}_\K$-map.
Hence $f$ is $C^k_\K$.

If $\K=\R$ and each $f_n$ is real analytic,
we find complex analytic maps $\wt{f}_n\colon \wt{U}_n\to (F_n)_\C$
on open subsets $\wt{U}_n\sub (E_n)_\C$ such that $U_n\sub\wt{U}_n$ and $\wt{f}|_{U_n}=f_n$.
Then $\wt{f}_n$ is $C^\infty_\C$, by Theorem~\ref{charcxcompl}, whence
\[
\wt{f}:=\bigoplus_{n\in\N} \wt{f}_n\colon\bigoplus_{n\in\N}\wt{U}_n\to \bigoplus_{n\in\N}(F_n)_\C
\]
is $C^\infty_\C$ and hence complex analytic, using Theorem~\ref{charcxcompl} again.
Identifying $E_\C$ with $\bigoplus_{n\in\N}(E_n)_\C$
and $F_\C$ with $\bigoplus_{n\in\N}(F_n)_\C$ (cf.\ Lemma~\ref{firstlasum}(c)),
we can interpret $\wt{f}$ as a complex analytic extension of~$f$.
Hence $f$ is real analytic.
\end{prf}

\section{Notes and comments on Chapter~\ref{chapfurcalc}}
The results presented profited from many sources.\medskip

\textbf{Analytic mappings.}
In the complex analytic case,
we have taken the expository papers
\cite{BS71a} and \cite{BS71b}
by J. Bochnak and J. Siciak
as a guideline,
removing however the completeness conditions
imposed there.
Analytic maps between non-complete spaces
have been discussed in~\cite{GE92} (cf.\ also~\cite{GE04}),
\cite{Gl02c}, and~\cite{Gl07c};
the examples of $C^k_\C$-maps which are not~$C^{k+1}_\C$
discussed in Exercises~\ref{exchlindepdop} and~\ref{exc-cx-notCk}
were taken from \cite{Gl07c} (compare \cite[Ex.~2.3]{GE92}
for functions which are $k$ times complex differentiable but not $k+1$ times,
disregarding continuity of the derivatives).
Further sources concerning complex analytic mappings include~\cite{He71}
and~\cite{He89}. Analytic mappings between Banach spaces
are also subsumed by the outline in~\cite{Bou67},
and interesting tools can be found in~\cite[Appendix]{La75}
and~\cite[Appendix~III]{Ne00}.

As far as real analytic mappings are concerned,
we favor
the approach of John Milnor~\cite{Mil84} (who defined real analytic maps via the existence of a complex
analytic extension),
as it is not 
clear whether
the class of real analytic maps
in the weaker sense of~\cite{BS71b}
is closed under composition. The exposition expands~\cite{Gl02c}.
For further material concerning complexifications of infinite-dimensional
manifolds and complex analytic extensions
of real analytic mappings in this context, we refer to~\cite{DGS14}
(cf.\ also~\cite{PS09}); for finite-dimensional manifolds,
classical references are~\cite{BrW59} and~\cite{Gr58b}.\medskip

{\bf Calculus in Banach spaces, inverse and implicit functions.}
Relations between
$C^k$-maps between normed spaces and classical $FC^k$-maps
were already described in the notes for Chapter~\ref{chapcalcul}.
The implication ``$\Leftarrow$'' from Proposition~\ref{companotions}(c)
was also recorded in \cite[Prop.~A.3.2]{Wa12}.
As we already have the theory of~$C^k$-maps at our disposal,
our discussion of $FC^k$-maps can deviate somewhat from classical
discussions as in \cite{Bou67}, \cite{CaH67}, \cite{Di60}, or~\cite{La99}.
For example, we can do without a Chain Rule for
compositions of mappings which are totally differentiable merely at a point.
Also, we are interested in the dependence on parameters not only in
normed spaces, but in general locally convex spaces (and related topics),
which do not have a counterpart in the classical works.
Exercise~\ref{exc-C1notFC1} was contributed by the authors.
Further examples of $C^1$-maps which are not~$FC^1$ can be found in
\cite[Ex.~6.8]{Mil82} and~\cite{WHO18}.

Continuous, Lipschitz, and Fr\'{e}chet differentiable
parameter-dependence of fixed points is a standard topic in the theory
of dynamical systems (see, e.g., \cite{Ir72}, \cite{We76}, and \cite[Thm.~C.7]{Ir80};
cf.\ also \cite[10.1.1]{Di60}).
Our discussion of $C^k$-dependence draws on~\cite{Gl06d} and~\cite{Gl07d}.
The Quantitative Inverse Function Theorem we present
varies \cite[pp.\,285--286]{We76}
(see also \cite{HrP70}), \cite{Gl06d}, and~\cite{Gl07d}).
As a simplified Newton iteration is used in the proof, we
discuss further aspects of the Newton method and variants in the exercises.
A criterion for convergence of the Newton method was first given in~\cite{KA64}.

The strategy to prove an Inverse Function Theorem with Parameters
using quantitative information on the images of balls combined with $C^k$-dependence of
fixed points is taken from~\cite{Gl06d} and~\cite{Gl07d};
theorems on implicit and inverse functions are then obtained as corollaries
(compare \cite{Hi99}, \cite{Te01b}, and \cite{Gl06a}
for related earlier results).
In contrast to the classical approach which seeks to reduce the implicit function theorem
to the inverse function theorem (and hence only applies to functions on open subsets of Banach spaces),
our strategy does not require completeness of~$E$ in Theorem~\ref{impl-fct-class},
nor openness of~$U$.
It also allows Theorem~\ref{impl-Ck} to be established with~$E$ an arbitrary locally convex space,
and variants for $C^{0,k}$-maps, as in Exercise~\ref{exc-C0k-impl}(b)
(which shall be valuable for the proof Proposition~\ref{special-ini-par},
a result concerning parameter-dependence of solutions to ordinary
differential equations which shall feed into the proof of $C^0$-regularity
of Banach--Lie groups and diffeomorphism groups).
Strict differentiability of mappings between Banach spaces (as broached in Exercises
\ref{exer-stri1} and \ref{exer-stri2})
was introduced by Leach (see \cite{Lea61})
and used to get an inverse function theorem
assuming only differentiability at one point.
See also~\cite{Bou67}, \cite{CaH67}, and ~\cite{Nij74}
for further information. The $FC^1$-analog of the
Global Inverse Function Theorem, Theorem~\ref{glob-hadamard},
is also discussed in~\cite{CH82} and~\cite{CDD85};
for generalizations, see~\cite{JLM17} and the references therein.

We refrain from proving advanced inverse function theorems,
like the Inverse Function Theorem of Nash and Moser
(as presented in~\cite{Ham82} and \cite[\S51]{KM97})
or Poppenberg's inverse function theorem~\cite{Pop99}
which imposes properties from the modern theory of Fr\'{e}chet space
(cf.\ \cite{MV97}). Unlike the preceding results, the inverse function theorems
in~\cite{Mu08} and~\cite{Gl06d} apply also to Fr\'{e}chet spaces which do not admit
a continuous norm (but their hypotheses are quite restrictive).
Another framework for inverse functions was presented in~\cite{Ma01},
which exploits generalizations of mapping degrees.\medskip

{\bf Differential equations in locally convex spaces.}
It is well-known that differential equations on non-normable locally convex spaces
need not be well-behaved; Examples \ref{ode-no} and \ref{ode-many}
of initial value problems
without (or with infinitely many) local solutions
are taken from \cite[Exs.~6.1/2]{Mil82}.\footnote{Compare
\cite{Dei77} for countable systems of ordinary differential equations in Banach spaces.}
To deal with this situation, we decided to turn local uniqueness
and existence (and the existence of local $C^k$-flows)
into axioms which may or may not be satisfied by a given differential equation.
When the axioms are satisfied, general results follow
(like existence of maximal solutions,
existence of a maximal flow, or differentiability properties of the latter) -- irrespective
of the classical theory of differential equations in Banach spaces.
Of course, many of the situations in which the axioms can be verified rely
on calculus in Banach spaces in some form.
For example, local uniqueness is guaranteed
if a differential equation in a projective limit of Banach spaces arises from a compatible
system of differential equations on the Banach spaces (whose right hand sides satisfiy a
local Lipschitz condition).
This idea was generalized in Proposition~\ref{uniq-pl}.
A typical case, Example~\ref{reallynee}, was used implicitly in the discussion
of a Frobenius theorem for finite-dimensional
vector distributions  on infinite-dimensional manifolds
in~\cite{Te01b} and made explicit in~\cite{Eyn12};
similar reasonings are also needed for a treatment of Banach distributions,
as in~\cite{Eyn14}.\footnote{Frobenius theorems for co-Banach distributions are available as well,
see~\cite{Hi00} and~\cite{Eyn16}.}
Compare also~\cite{DGV16} for a general theory of projective limits of Banach manifolds.

In Section~\ref{sec-ode}, we complemented the formal, axiomatic approach by a discussion of 
the classical case of initial value problems in Banach spaces,
which exploits the results concerning non-linear mappings
between function spaces developed in Section~\ref{secCspaces}.
Our starting point
is quantitative version of the Picard--Lindel\"{o}f Theorem,
which constructs the solution as the fixed point of a contraction.
We then used results concerning parameter dependence
of fixed points from Section~\ref{provis-Ban} to get
Theorem~\ref{local-ivp-dep-Ck} on $C^k$-dependence
of solutions on initial conditions and parameters
(as well as real analytic dependence), and
Proposition~\ref{special-ini-par} on $C^{0,k}$-dependence.
Our treatment on $C^k$-dependence is a variant of
the method in~\cite{Gl06d} (where an implicit function theorem is applied instead),
while~\cite{AlS15} (where~$t_0$ is fixed) and~\cite{GN17}
are precursors of our treatment of $C^{0,k}$-dependence.
For classical discussions, see \cite[Thm.\,3.7.1]{CaH67}
(showing $FC^k$-dependence on $(t,y_0,p)$ for fixed $t_0$ and parameter~$p$
in a Banach space), or \cite[10.8.2]{Di60} (for $FC^k$-dependence on $(t,t_0,y_0)$)
and \cite[10.7.4]{Di60} (for $FC^k$-dependence on $(t,p)$ for fixed $(t_0,y_0)$
and~$p$ in a Banach space). See also \cite[Thm.\,IV.1.11]{La99}
for $FC^k$-dependence on $(t,y_0)$ (for fixed~$t_0$)
and Thm.\,IV.1.16 in loc.\ cit. for the $FC^k$-property of the global flow
(as a function of $(t,y_0)$) of a time-independent vector field.
The general idea to study the dependence of solutions to initial value problems
by means of the implicit function theorem is due to~\cite{Ro68} (cf.\ also~\cite{CH82}).
To prove real analytic dependence on initial conditions and parameters,
we constructed a complex analytic extension and therefore had to consider
differential equations in a complex variable as well.
This general strategy was also used in~\cite{Di60}.
Likewise, we reduced the real analytic implicit function theorem
to the complex analytic case, which simply amounts to the $C^\infty$-case over~$\C$
(cf.\ Remark~\ref{variants-i-para}(c)).
In both cases, this strategy (which is known from~\cite{Di60}) makes it superfluous
to prove the convergence of power series expansions of solutions
by hand, or any estimates to insure the convergence
(in contrast to more computational approaches as 
in~\cite{KP02a, KP02b}).
We mention that -- valuable as it can be -- we did not need to use
Groenwall's Lemma in our discussion of differential equations.
We also omit a discussion of differential equations with 
discontinuous right hand sides; see \cite{Son98} or~\cite{Scr97}
(and the references therein) for differential equations in finite-dimensional or Banach spaces,
\cite{Gl15b} and \cite{NR17} for generalizations and Lie theoretic applications
(cf.\ also~\cite{MNe18}).

In Exercises~\ref{exerc-cx-norm} and \ref{exerc-fp-an2}, we used certain complexification norms
(the maximal ones). For general discussions of
complexifications of norms, see~\cite{MST99}
and~\cite{Kir01}. Compare~\cite{Da14} for Lie theoretic applications.

{\bf Calculus on metrizable spaces.}
It is well known that a map~$f$ on a metrizable locally convex space is smooth
(in the sense of Definition~\ref{defnhigher}) if and only if~$f$ is smooth along each smooth curve,
viz., a smooth map in the sense of Convenient Differential Calculus
(compare \cite[Thm.\,12.8]{KM97};
for mappings between finite-dimensional spaces, see~\cite{Bm67}).
The idea to check the $C^k$-property along smooth maps of $k+1$ variables
(as in Theorem~\ref{mayn2} is taken from~\cite{BGN04};
our proof is simpler as we assume local convexity.
Already for mappings between finite-dimensional spaces, the $C^k$-property
for finite~$k$ cannot be checked along smooth curves , see~\cite{Bm67}.
But two variables suffice, as noted in~\cite{Bm67}.
That two variables suffice also for maps on metrizable spaces (as explained in
Remark~\ref{two-var-enough})
is enabled by an argument from \cite[p.\,27]{KM97}
(as worked out in Exercise~\ref{exc-gat-Ck}).
Theorem~\ref{mayn1} is a variant of \cite[Thm.\,4.11(1)]{KM97}.
The idea to use functions of finitely many variables (so-called ``plots'') to single out
functions of interest is also at the heart of the theory of diffeological spaces
(cf.\ \cite{So84}, \cite{DI85}, and the references therein).
We mention that the Special Curve Lemma (Lemma~\ref{specialcurve})
is a standard tool of Convenient Differential calculus
(see \cite[p.\,18]{KM97} or~\cite{FK88}).

{\bf Calculus on Silva spaces.}
Our results concerning $C^k$-maps on Silva spaces are based on~\cite{Gl07a}.
We mention that smoothness in the sense of Definition~\ref{defnhigher}
and smoothness along smooth curves (as in Convenient Differential calculus)
are equivalent for mappings on open subsets of Silva spaces (see, e.g., \cite{CaP14}).
For manifolds modeled on metrizable spaces, and for manifolds modeled on Silva spaces,
smooth maps can be characterized as those which pull back smooth functions (on open subsets)
to such (see~\cite{EW17}).
For a certain concept of Silva differentiable mappings between Silva spaces,
compare \cite[Defn.\,2.2]{Les82} and the references therein.

{\bf Mappings between direct sums.}
The results concerning mappings between locally convex direct sums are based on~\cite{Gl03b}
(where also uncountable direct sums are considered)
and~\cite{Gl04b} (which includes results with parameters).

As a counterpart to the explicit description of a basis of $0$-neighborhoods
in the space $C^\infty_c(\R^n)$ of test functions given in \cite[Ch.\,III, \S1]{Schw57},
we shall see later that spaces of vector-valued compactly supported smooth functions (or compactly supported smooth sections in vector bundles) can be identified with closed vector subspaces
of suitable locally convex direct sums
(as, e.g., in \cite{Gl04b}, \cite{Gl03c}, \cite{Gl05d}, and~\cite{BG14}).
Proposition~\ref{diff-sum} can then serve as a tool to establish continuity and differentiability properties
for non-linear mappings between such spaces of compactly supported functions (or sections).

\chapter[Basic theory of infinite-dimensional manifolds]{Basic\hspace*{-.3mm} theory\hspace*{-.4mm} of \hspace*{-.3mm}infinite-dimensional\hspace*{-.3mm} manifolds}\label{chapmanif}

In this chapter,
we discuss manifolds modeled
on locally convex spaces
and the basic concepts going along with them:
smooth (and $C^r$-) maps between
manifolds, tangent maps,
vector fields, and differential forms.

To assist readers without background
knowledge concerning differentiable manifolds,
we explain all of the basic ideas
and constructions in detail.
More experienced readers are invited
to skip these explanations,
and direct their attention to
the differences compared to
the finite-dimensional case.

At the end of the chapter, we compile
preparatory material needed in Chapter~\ref{chap-manifold-constructions} for basic constructions of infinite-dimensional
manifolds. Thus,
Section~\ref{tools-mfd} covers aspects of vector bundles and principal bundles;
Section~\ref{sec:local-add} is devoted to the construction of local additions and tubular
neighborhoods using sprays and their associated exponential maps.
\section{Manifolds modeled on locally convex spaces}\label{secmanlcx}
So far, we have mainly been dealing with
$C^r$-maps between open
subsets of locally convex spaces.
We now introduce
the corresponding ``global''
objects, namely $C^r$-manifolds
modeled on locally convex spaces.
While classical manifolds
locally look like
open subsets of~$\R^n$,
the manifolds considered here
locally look like open subsets of a
given locally convex space.
As we shall see later, $C^k$-maps
between $C^r$-manifolds (for $k\leq r$)
can be defined in a natural way,
and also many other concepts of analysis can be adapted
to the case of manifolds.

Smooth manifolds ($r=\infty$) are of main interest,
but we find it useful to present
the foundations more generally for $C^r$-manifolds
with $r\in \N_0\cup\{\infty,\omega\}$
over a ground field $\K\in\{\R,\C\}$,
using ``$C^\omega$'' as a shorthand for ``analytic''.\footnote{Concerning the symbols
$\infty$ and $\omega$, our conventions are that $n<\infty<\omega$ for each $n\in\N_0$.
Moreover, $\infty-n:=\infty$ and $\omega-n:=\omega$.}
All locally convex spaces are locally convex spaces over~$\K$.
If we wish to emphasize the chosen ground field, we write $C^r_\K$
in place of~$C^r$.

The following simplification will be used:
\begin{convent}
Let $f\colon X\to Y$ and $g\colon A\to B$
be any maps.
We shall write $g\circ f$
as a shorthand for
the composition
$g\circ f|_{f^{-1}(A)} \colon f^{-1}(A)\to B$,
whenever this helps to avoid
clumsy notation.
\end{convent}
\begin{rem}
We shall mostly use the
preceding convention
if both $Y$ and~$A$
are subsets of a given topological space~$Z$,
$A$ is open in~$Z$, and~$f$~continuous.
Then the domain
$f^{-1}(A)=f^{-1}(Y\cap A)$
of $g\circ f$ is open in~$X$.
\end{rem}
\begin{defn}\label{defnatlas}
Let $M$ be a Hausdorff topological space
and $E$ be a locally convex topological $\K$-vector
space.
An \emph{$E$-chart} for~$M$
is a homeomorphism
$\phi\colon U_\phi\to V_\phi$
from an
open subset
of $U_\phi \sub M$ onto an open subset
$V_\phi\sub E$.
Let $r\in \N_0\cup\{\infty,\omega \}$.
A \emph{$C^r$-atlas for $M$} (modeled on~$E$)
is a set $\cA$ of $E$-charts for~$M$
satisfying the following
conditions:
\begin{description}[(DA)]
\item[(A1)]
$M=\bigcup_{\phi\in \cA} U_\phi$,
that is,
$M$ is covered by the domains $U_\phi$;
\item[(A2)]
All
$\phi,\psi\in \cA$
are \emph{$C^r$-compatible}
in the sense
that the \emph{transition maps}
\begin{equation}\label{transi}
\phi\circ \psi^{-1} \colon
\psi(U_\phi\cap U_\psi)\to V_\phi\sub E
\end{equation}
(which are mappings between open subsets of $E$)
are~$C^r$.
\end{description}
\end{defn}
Note that the set of $C^r$-atlases for $M$ can be ordered
by inclusion.
\begin{defn}\label{defnmfd}
A \emph{$C^r$-manifold modeled on~$E$}
is a pair $(M,\cA)$,
where~$M$ is
a Hausdorff topological space
and $\cA$ a maximal $C^r$-atlas for~$M$.
We shall mainly be concerned with
$C^\infty$-manifolds,
which are also called \emph{smooth mani\-folds}.
The $C^\omega_\K$-manifolds are also called
\emph{$\K$-analytic manifolds}.
If $(M,\cA)$ is a manifold,
we reserve the word ``chart''
for the $E$-charts $\phi\in \cA$.
Given $x\in M$,
a \emph{chart around~$x$}
is a chart $\phi\colon U_\phi\to V_\phi$
in~$\cA$ such that $x\in U_\phi$.
We write $\cA_x$ for the set of all charts
of~$M$ around~$x$.
A manifold
modeled on a Banach space
(Hilbert space, or Fr\'{e}chet space)
is called a \emph{Banach manifold},
\emph{Hilbert manifold} and
\emph{Fr\'{e}chet manifold},
respectively.
\end{defn}
\begin{rem}\label{precrem}
Every $C^r$-atlas $\cA$ for~$M$
is contained in a (uniquely determined)
maximal $C^r$-atlas, namely the
set $\cM$
of all $E$-charts $\phi\colon U_\phi\to V_\phi$
which are compatible with $\cA$
in the sense that
the transition maps
\[
\phi\circ \kappa^{-1}|_{\kappa(U_\phi\cap U_\kappa)}
\qquad\mbox{ and }\qquad
\kappa\circ \phi^{-1}|_{\phi(U_\kappa\cap U_\phi)}
\]
are~$C^r$ for all $\kappa\in \cA$.
To see that $\cM$ is a $C^r$-atlas,
one has to check that all
$\phi,\psi\in \cM$ are $C^r$-compatible
(in the sense of Definition~\ref{defnatlas},
(A2)).
We verify this via
``insertion of charts,'' a simple
standard argument based
on the Chain Rule which will
usually be left to the reader
in the sequel.
To see that $\phi\circ \psi^{-1}$
is $C^r$ on a neighborhood
of a
given element $x\in \psi(U_\phi\cap U_\psi)$,
the idea is to pick a chart
$\kappa\in \cA$ such that
$\psi^{-1}(x) \in U_\kappa$.
Then $W:=\psi(U_\phi\cap U_\psi\cap U_\kappa)$
is an open neighborhood of~$x$
in $\psi(U_\phi\cap U_\psi)\sub E$,
and the definition of~$\cM$ implies that
\[
\phi\circ \psi^{-1}|_W
\;=\;
(\phi\circ \kappa^{-1})
\circ
(\kappa \circ \psi^{-1})|_W
\]
is a $C^r$-map on the open neighborhood~$W$
of~$x$, being a composition of
$C^r$-maps. As $x$ was arbitrary,
we deduce that $\phi\circ \psi^{-1}$
is $C^r$ on all of its domain $\psi(U_\phi\cap U_\psi)$.
Thus
$\phi$ and $\psi$ are indeed $C^r$-compatible.
\end{rem}
\begin{rem}
Let
$(M,\cA)$ be a
$C^r$-manifold
with modeling space~$E$.\medskip

\noindent
(a) Occasionally, we shall refer to the charts $\phi\in\cA$
also as \emph{local charts}.
A chart $\phi\colon U_\phi\to V_\phi$
is called a \emph{global chart} for~$M$ if $U_\phi=M$.
Of course, global charts need not exist.\medskip

\noindent
(b) Let $F$ be a locally convex space
isomorphic to~$E$.
Occasionally,
we shall find it useful
to loosen our terminology
and call also an $F$-chart $\phi$ of~$M$
a ``chart of~$M$''
if~$\phi$ is $C^r$-compatible with~$\cA$
(in the sense of
Remark~\ref{precrem}).
This situation
will occur only rarely,
and the abuse will be clear from the context.
\end{rem}
\begin{rem}
Usually, for a given topological space~$M$
which we want to turn into a $C^r$-manifold,
we do not have a maximal $C^r$-atlas
right away; as a rule,
we only know some (non-maximal)
atlas.
We can then extend this atlas to a maximal atlas
in order to consider~$M$ as a $C^r$-manifold.
\end{rem}
Here are some
elementary
examples of manifolds.
\begin{ex}\label{exmfd1}
Every locally convex $\K$-vector space $E$
is a $C^\infty$-manifold modeled on~$E$ in a natural way,
equipped with the maximal $C^\infty$-atlas containing
the $C^\infty$-atlas $\{\id_E\colon E\to E\}$.
More generally, each open subset $U\sub E$ becomes
a $C^\infty$-manifold when equipped with the maximal
$C^\infty$-atlas containing $\{\id_U\}$.
Likewise, $E$ (and each open subset thereof)
can be considered as a $\K$-analytic manifold.
\end{ex}
\begin{ex}\label{exmfd2}
If $(M,\cA)$ is a $C^r$-manifold
and
$U\sub M$ an open subset,
then $U$ can be made a $C^r$-manifold:
We equip $U$ with the induced
topology and note that
$\{\phi\in \cA\colon U_\phi\sub U\}$
is a maximal $C^r$-atlas for $U$
(exercise\,!).
\end{ex}
More interesting examples of manifolds
are (not necessarily open)
\emph{submanifolds} of given manifolds.
%
%
%
\begin{defn}\label{defsubm}
Let $(M,\cA)$ be a $C^r$-manifold
modeled on a locally convex space~$E$,
and $F\sub E$ be a
closed vector subspace.
A \emph{sub\-mani\-fold of~$M$ modeled on~$F$}
is a subset $N\sub M$
such that, for each $x\in N$,
there exists a chart $\phi\colon U_\phi\to V_\phi$
of~$M$ around~$x$ such that
$\phi(U_\phi\cap N)=V_\phi\cap F$.
In this case, $\phi_N:=\phi|_{U_\phi\cap N}^{V_\phi\cap F}$
is an $F$-chart for~$N$,
and the set $\cB$ of all $F$-charts
$\phi_N$ so obtained is a $C^r$-atlas
for~$N$,
exploiting that $F$ is closed.\begin{footnote}
{Here it would suffice that $F$ is sequentially closed, but we are 
not aware of natural situations, where the model space of a submanifold 
is sequentially closed and not closed.} 
\end{footnote}
In fact, given
$\phi_N, \psi_N\in \cB$,
their transition map
\[
\phi_N\circ \psi_N^{-1}\; =\;
\phi\circ \psi^{-1}|_{\psi(U_\phi\cap U_\psi\cap N)}
\]
is $C^r$ as a map to~$E$
(because $\phi\circ \psi^{-1}$ is~$C^r$).
Since $\phi_N\circ\psi_N^{-1}$
takes its values in the
closed vector subspace~$F$,
it is also $C^r$ as a map to~$F$,
by Lemma~\ref{corestr}
(resp., Exercise~\ref{excreatosub}).\\[2.5mm]
We now equip $N$ with the maximal atlas
containing~$\cB$.
The charts $\phi_N\in \cB$ are called
\emph{submanifold charts} and the corresponding charts~$\phi$ of~$M$
are called \emph{adapted to~$N$}.
If $N$ is
a submanifold of~$M$ modeled
on~$F$ and $F$ is a
complemented
vector subspace of~$E$
(as in Definition~\ref{defcplsub}),
then~$N$ is called a \emph{split} submanifold.
\end{defn}
\begin{rem}
Knowing that a submanifold
$N\sub M$ is split in the preceding sense
can be useful
if~$M$ is a Banach manifold
(whence Inverse-
and Implicit
Function Theorems are available).
For non-Banach manifolds,
the concept is less useful.
We shall
therefore define ``split Lie subgroups''
later in a different way
(see Definition~\ref{defc.4}).
\end{rem}
\begin{rem}
We shall see
in Lemma~\ref{basinitmfd}
that $N$
determines~$F$ up to isomorphism
and that
the $C^r$-manifold structure induced
on~$N$ is essentially unique,
if $N\not=\emptyset$ and $r\geq 1$.
\end{rem}
\begin{rem}
For the purposes
of Lie theory,
more general subsets
than submanifolds
(so-called initial submanifolds)
will be needed.
We shall
introduce these
in Definition~\ref{definitmfd}.
\end{rem}
\begin{rem}\label{underlyingmf}
If $(M,\cA)$ is a $C^r$-manifold
and $k\in \N_0\cup\{\infty \}$ is such that $k\leq r$,
then the maximal $C^r$-atlas $\cA$
also is a $C^k$-atlas, and hence extends to
a maximal $C^k$-atlas $\cB$.
Thus $(M,\cA)$ has an underlying $C^k$-manifold
$(M,\cB)$. By abuse of notation, we shall mostly
write $M$ as a shorthand for $(M,\cA)$,
and use the same symbol for $(M,\cB)$
whenever we want to consider $M$ as a $C^k$-manifold.
In the same way, every $C^r_\C$-manifold $(M,\cA)$ has an underlying $C^r_\R$-manifold,
obtained by passing to the maxmal $C^r_\R$-atlas containing~$\cA$.
\end{rem}
\begin{ex}
Submanifolds of locally convex spaces
already provide a rich supply of interesting
manifolds.
For instance, we shall see later that each closed
subgroup of $\GL_n(\R)$ is a smooth submanifold
of~$M_n(\R)$.
Some elementary examples of submanifolds
will be discussed presently in
the exercises, like the circle
$\bS^1$ in~$\R^2$, the $(n-1)$-sphere
$\bS^{n-1}$ in~$\R^n$
and the special linear group $\SL_2(\R)$
inside $M_2(\R)\isom \R^4$.
\end{ex}
Direct products of manifolds
will often be used.
\begin{defn}\label{exmprod}
Let $(M_1,\cA_1)$ and $(M_2,\cA_2)$
be $C^r$-manifolds
modeled on topological $\K$-vector spaces $E_1$
and $E_2$, respectively.
Then
\[ \cA:=\{\phi_1\times \phi_2\colon \phi_1\in \cA_1,\phi_2\in \cA_2\} \] 
is a $C^r$-atlas for
$M_1\times M_2$, endowed with the product topology (as we presently
explain).
Hence the maximal atlas containing~$\cA$
makes $M_1\times M_2$ a $C^r$-manifold modeled on $E_1\times E_2$,
called the \emph{direct product} of~$M_1$ and~$M_2$.
\end{defn}
To see that $\cA$ is an atlas,
the main point is to show that
all transition maps are $C^r$.
Given
$\phi_1\times \phi_2$, $\psi_1\times \psi_2\in \cA$,
the map
\[ (\phi_1\times \phi_2)\circ (\psi_1\times\psi_2)^{-1}
=
(\phi_1\circ \psi_1^{-1}) \times (\phi_2\circ \psi_2^{-1})\] 
is $C^r$, being a direct product
of two $C^r$-maps (cf.\ Lemma~\ref{lemprod}
and Exercise~\ref{excreatoprod}).
Thus (A2) holds, and also
(A1) is easily verified.
\begin{rem}
In classical textbooks on differential geometry
dealing with finite-dimen\-sional manifolds,
all manifolds are modeled on~$\R^n$,
and it is usually assumed
that manifolds satisfy the second axiom of countability
(or equivalently, that they are
$\sigma$-compact).
This makes manifolds more manageable:
Every $\sigma$-compact
finite-dimensional smooth manifold is
paracompact (see Lemma~\ref{pre-lcp-para})
and actually it admits
smooth partitions
of unity (see Proposition~\ref{findim-smoothly-para}).
The latter facilitate localization
arguments, and are important tools for
differential geometry and global analysis.
We mention that
every $\sigma$-compact
finite-dimensional smooth manifold can be realized as
a smooth submanifold of some $\R^n$,
by Whitney's Embedding Theorem \cite{Wh36}. 
Paracompact finite-dimensional
manifolds are merely disjoint unions
of open $\sigma$-compact submanifolds
(see Proposition~\ref{lcp-parac})
and therefore equally well tractable.\\[2.3mm]
It would not make sense to include
conditions like $\sigma$-compactness,
paracompactness or second countability
in the definition of general manifolds modeled
on locally convex spaces
(because these are not satisfied
by many important examples).
Rather, whenever such an
additional property is required,
we shall say so explicitly.
\end{rem}
\begin{rem}\label{remregul}
We remark in this connection
that manifolds modeled on
locally convex spaces
need not be regular topological spaces,
in contrast to finite-dimensional manifolds,
which are regular because they are locally
compact
(recall that a topological space~$X$ is
\emph{regular} if it is Hausdorff
and every neighborhood~$U$ of
any given point $x\in X$ contains
a closed neighborhood of~$x$).
An example of such
a non-regular manifold
can be found in Exercise~\ref{exer:2.1.6} 
(adapted from \cite{Mil82}).
For Lie groups modeled on locally convex spaces
(as defined later), the pathology cannot occur:
Each Lie group being a topological group,
its underlying topological space
is regular,
and even completely regular
(see \cite[Ch.\,II, \S\,8]{HR79}).
\end{rem}
\begin{small}
\subsection*{Exercises for Section~\ref{secmanlcx}} 
\label{Ex3.1}
\begin{exer}
  \label{exer:2.1.1}
Define
$\phi \colon \;]0 ,2\pi[\;\to \bS^1\setminus \{1\}$,
$\phi(t):=e^{it}$ and $\psi\colon
\;]{-\pi},\pi[\;\to \bS^1\setminus \{{-1}\}$,
$\psi(t):=e^{it}$,
where $\bS^1:=\{z\in \C\colon |z|=1\}$ is the unit
circle.
Show that $\{\phi^{-1}, \psi^{-1}\}$
is a $C^\infty$-atlas for~$\bS^1$.
\end{exer}

\begin{exer}
  \label{exer:2.1.2} Let $M_1$ and $M_2$ be $C^r$-manifolds
and $N_1\sub M_1$
as well as $N_2\sub M_2$
be submanifolds.
Show that
$N_1\times N_2$ is
a submanifold of~$M_1\times M_2$
whose manifold structure
(as a submanifold)
coincides with the direct product
of~$N_1$ and~$N_2$.
\end{exer}

\begin{exer}
  \label{exer:2.1.6} Let $n,d\in \N$
with $d<n$ and $M\sub \R^n$ be a subset
with the following property:
For each $p\in M$, there exists
a $C^r$-map $f\colon W\to\R^{n-d}$
on an open neighborhood $W\sub \R^n$ of~$p$
such that $M\cap W=f^{-1}(\{0\})$
and such that $f'(x)$ has full rank
for each $x\in W$.
Show that $M$ is a $C^r$-submanifold
of~$\R^n$ modeled on~$\R^d\times \{0\}$.
(Hint: For $p$ and $f$ as before,
after a permutation of the
coordinates we may assume that
$h\colon \R^n\supseteq
W\to\R^n$, $h(x_1,\ldots, x_n):=(x_1,\ldots, x_d,f(x))$
has an invertible differential at~$p$
and hence restricts to a local
diffeomorphism $\phi$ on some open neighborhood
of~$p$, by the Inverse Function Theorem).
\end{exer}

\begin{exer}
  \label{exer:2.1.7} Using Exercise~\ref{exer:2.1.6}, show that the sphere
\[ \bS^{n-1}:=\Big\{(x_1,\ldots, x_n)\in \R^n\colon
\sum_{j=1}^nx_j^2=1\Big\}\] 
is an $(n-1)$-dimensional
smooth submanifold of~$\R^n$.
Also show that
$\SL_2(\R):=\{x\in M_2(\R)\colon \det(x)=1\}$
is a $3$-dimensional smooth submanifold
of $M_2(\R)\isom \R^4$.
\end{exer}

\begin{exer}\label{exc-mfd-local}
Let $r\in\N_0\cup\{\infty,\omega\}$,
$E$ be a locally convex space, $M$ be a Hausdorff topological space
and $(U_j)_{j\in J}$
be an open cover of~$M$. Assume that a $C^r$-manifold structure
modeled on~$E$ is given on $U_j$ for each $j\in J$,
such that $U_i$ and $U_j$ induce the same $C^r$-manfold structure
on $U_i\cap U_j$ for all $i,j\in J$.
Show that $M$ can be made a $C^r$-manifold modeled
on~$E$ in such a way that $M$ induces the given $C^r$-manifold
structure on each of the open substs $U_j$.
The latter $C^r$-manifold structure is unique.
\end{exer}

\begin{exer}
  \label{exer:2.1.8} (A non-regular manifold) 
Let $M$ be the subset of the Fr\'{e}chet space
$\R^\N$ consisting of all sequences $(x_k)_{k\in \N}$
such that all but finitely many $x_k$ belong to the open interval
$J :=\;]{-1},1[$.
Thus $M=\bigcup_{n\in \N_0}M_n$,
where $M_n:=\R^n\times J^{\{n+1,n+2,\ldots\}}$.
We equip $M_n$ with the product topology
and give $M$ the final topology with
respect to the inclusion maps $M_n\to M$
(see Definition~\ref{defnfinal}
in Appendix~\ref{appA}).
Thus $U\sub M$ is open if and only
if, for each $n\in \N_0$,
the intersection $U\cap M_n$ is open in $M_n$.
\begin{description}[(D)]
\item[(a)]
Show that $M_n$ is open in $M$,
and that $M$ induces the given topology on~$M_n$,
for each $n\in \N_0$.
Deduce that~$M$ is Hausdorff.
\item[(b)]
Find the closure of a set of the form
$U_1\times \cdots\times U_n\times J^{\{n+1,n+2,\ldots\}}$
in~$M$, where $n\in \N_0$
and $U_1,\ldots, U_n$ are any non-empty subsets of~$\R$.
\item[(c)]
Show that every neighborhood~$U$ of $(0,0,\ldots)$
in~$M$ contains a set of the form described
in~(b) and deduce that $\wb{U}\not\sub J^\N$.
Noting that $J^\N$ is a neighborhood of
$(0,0,\ldots)$, conclude that the topological
space $M$ is not regular.
\item[(d)]
The map
$f\colon J \to \R$, $f(x):=x/(1-x^2)$
is a $C^\infty$-diffeomorphism.
Show that the map 
\[ \phi_n\colon M_n\to \R^\N, \quad 
(x_k)_{k\in \N} \mapsto (x_1,\ldots, x_n,f(x_{n+1}),f(x_{n+2}),\ldots)\] 
is a homeomorphism,
for each $n\in \N_0$.
Show that $\{\phi_n\colon n\in \N_0\}$
is a $C^\infty$-atlas making $M$ a
smooth manifold modeled on $\R^\N$
(and actually a real analytic manifold).
\end{description}
\end{exer}
\begin{exer} \label{exer:2.1.9} (Banach--Grassmann manifolds) 
Given a Banach space~$E$ and
complemented vector subspace
$F\sub E$
(as in Definition~\ref{defcplsub}),
let $\Gr_F(E)$ be the set of all
complemented vector subspaces $H\sub E$
such that $F\isom H$ and $E/F\isom E/H$
as topological vector spaces.
The goal of this exercise
is to make $\Gr_F(E)$ an
analytic manifold modeled
on $\cL(F,E/F)$,
a so-called \emph{Grassmann manifold}.
\begin{description}[(D)]
\item[(a)]
Given $H\in \Gr_F(E)$
and a complementary subspace
$Y \sub E$ for~$H$,
we let $\pr_1^{H,Y}\colon E=H\oplus Y\to H$
and $\pr_2^{H,Y}\colon E=H\oplus Y\to Y$
be the projections onto the first
and second component.
Show that, for each continuous
linear map $\alpha\colon H \to Y$,
the graph 
\[ \Gamma^{H,Y}(\alpha):=\{x+\alpha(x)\colon x\in H\} \] 
is a closed
vector subspace of~$E$
which has~$Y$ as a complement.
Also show that
$\pr_1^{H,Y}$
restricts to
an isomorphism of topological
vector spaces
$\Gamma^{H,Y}(\alpha)\to H$,
with inverse
$\theta^{H,Y}_\alpha
\colon H\to
\Gamma^{H,Y}(\alpha)$,
$\theta^{H,Y}_\alpha(x):=x+\alpha(x)$.
As a consequence,
$\Gamma^{H,Y}(\alpha)\in \Gr_F(E)$.
We set 
\[ U_Y:=\{\Gamma^{H,Y}(\alpha) \colon \alpha\in \cL(H,Y)\}.\] 
Show that the map $\Gamma^{H,Y}\colon \cL(H,Y)\to U_Y$,
$\alpha\mto \Gamma^{H,Y}(\alpha)$
is a bijection.
Also show that $U_Y$ is
the set of all
vector complements of~$Y$ in~$E$.
\item[(b)]
Let~$H$ and~$Y$
be as before,
and $G\in U_Y$.
Then $H=\Gamma^{G,Y}(\alpha)$
for some $\alpha\in \cL(G,Y)$.
Show that $(\Gamma^{G,Y})^{-1}(\Gamma^{H,Y}(\beta))
=\alpha+ \beta \circ \theta^{G,Y}_\alpha$
for each $\beta\in \cL(H,Y)$.
Hence $(\Gamma^{G,Y})^{-1}\circ\Gamma^{H,Y}
=\alpha+ \cL(\theta^{G,Y}_\alpha,Y)$
is a continuous
affine-linear (and hence analytic) map.
\item[(c)]
Let $H\in \Gr_F(E)$ and both~$Y$ and~$Z$
be vector complements for~$H$ in~$E$.
Given $\alpha\in \cL(H,Y)$,
show that $\Gamma^{H,Y}(\alpha)\in U_Z$
if and only if the linear map
$\kappa_\alpha\colon
\pr_1^{H,Z}\! \circ \,\theta^{H,Y}_\alpha\colon
H\to H$ is invertible.
Using that $\GL(H)$ is open in $\cL(H)$,
deduce that $(\Gamma^{H,Y})^{-1}(U_Z)$
is open in $\cL(H,Y)$.
For
$\alpha\in (\Gamma^{H,Y})^{-1}(U_Z)$,
show that $\Gamma^{H,Y}(\alpha)=\Gamma^{H,Z}(\beta)$
with $\beta=
\kappa_\alpha^{-1}-\alpha\circ\kappa_\alpha^{-1}-\id_H$.
Thus
$(\Gamma^{H,Z})^{-1}(\Gamma^{H,Y}(\alpha))=
\kappa_\alpha^{-1}-\alpha\circ\kappa_\alpha^{-1}-\id_H$,
showing that
$(\Gamma^{H,Z})^{-1}\circ \Gamma^{H,Y}$
is an analytic map.
\item[(d)]
Give $\Gr_F(E)$
the final topology with respect to
the maps $\theta^{H,Y}\colon \cL(H,Y)\to\Gr_F(E)$,
for $H\in \Gr_F(E)$ and $Y$ a vector complement
to~$H$ in~$E$.
Show that each $U_Y$ is open in $\Gr_F(E)$
and that $\theta^{H,Y}$ is a homeomorphism
onto~$U_Y$.
\item[(e)]
$\GL(E)$ acts on $\Gr_F(E)$ via
$\sigma(\Lambda,H):=\Lambda(H)$ for $\Lambda\in \GL(E)$,
$H\in \Gr_F(E)$.
The stabilizer $S:=\{\Lambda\in \GL(E)\colon
\Lambda(F)=F\}$ of~$F$ is a closed
subgroup of $\GL(E)$ and hence $\GL(E)/S$
is Hausdorff (Exercise~\ref{exer:6.1.1}). 
Since the action is transitive,
the map
$h\colon \GL(E)/S\to \Gr_F(E)$, $h(\Lambda S):=\Lambda(F)$
is a bijection.
Hence $\Gr_F(E)$ will be Hausdorff
if we can show that~$h$ is an open map.
We proceed in steps:
Show that $\sigma_\Lambda \colon
\Gr_F(E)\to\Gr_F(E)$, $H\mto \Lambda(H)$
is a homeomorphism,
for each $\Lambda\in \GL(E)$.
Fix a vector complement $Y$ for
$F$ in~$E$. Show that 
\[ \tau(\Gamma^{F,Y}(\alpha)):=\id_E +\alpha\circ \pr^{F,Y}_1
\in \GL(E) \] 
for $\Gamma^{F,Y}(\alpha)$
in some neighborhood $W\sub U_Y$ of~$F$,
and that $\tau\colon W\to \GL(E)$
is a continuous local section
to $\sigma_\Lambda$
(i.e., $\sigma_\Lambda\circ \tau=\id_W$)
such that $\tau(F)=\one$.
Conclude that $\sigma_\Lambda$
is open at~$\one$ (i.e.,
it takes $\one$-neighborhoods to
neighborhoods of~$F$).
Infer that $\sigma_\Lambda$
is open and hence also~$h$.
\item[(f)]
For each $H$ and~$Y$ as before,
there are isomorphisms
of topological vector spaces
$\ve\colon F\to H$
and $\delta\colon Y\to E/F$.
Define $\phi_{H,Y,\ve,\delta}\colon
U_Y\to \cL(F,E/F)$
via $\phi_{H,Y,\ve,\delta}:=
\cL(\ve,\delta)\circ (\theta^{H,Y})^{-1}$,
where $\cL(\ve,\delta)\colon \cL(H,Y)\to\cL(F,E/F)$,
$\alpha\mto \delta\circ \alpha\circ \ve$.
Show that the set of these maps
$\phi_{H,Y,\ve,\delta}$
is an analytic atlas of
$\cL(F,E/F)$-charts for $\Gr_F(E)$.
Passing to a maximal atlas,
we obtain the desired
analytic manifold structure
on $\Gr_F(E)$.
\end{description}
\end{exer} 
\begin{exer}
  \label{exer:2.1.10} Let us discuss phenomena
related to Grassmann manifolds.
\begin{description}[(D)]
\item[(a)]
Let $E:=\ell^2$,
$F:=\{x \in E \colon x_1=0\}$
und $H:=\{x\in E \colon x_1=x_2=0\}$.
Show that $F\isom H$ but
$E/F\not\isom E/H$
and hence $H\not\in \Gr_F(E)$.
\item[(b)]
Show that every $1$-dimensional
vector subspace $F$ of a
locally convex space~$E$
is complemented (if $F\not=E$,
then $\ker(\lambda)$ is a vector complement for~$F$,
for a suitable linear functional $\lambda\in E'$).
Using induction on~$n$,
deduce that every $n$-dimensional
vector subspace $F\sub E$ is complemented.
\item[(c)]
Let $E$ be a Banach space, $n\in \N$
and both $F\sub E$ and $H\sub E$
be $n$-dimensional
vector subspaces.
Show that $E/F\isom E/H$ and hence
$H\in \Gr_F(E)$.
\,[Hint: Let $Y$ and $Z$ be vector complements
for $F$ and~$H$, respectively.
Then $Y\cap Z$
has finite codimension in~$E$,
hence the same codimension in both
$Y$ and~$Z$].\\[1.2mm]
By the preceding,
$\Gr_n(E):=\Gr_F(E)$ is independent
of the choice of~$F$
(up to isomorphism of the modeling space
of $\Gr_n(E)$).
The manifold
$\bP(E):=\Gr_1(E)$ of all $1$-dimensional
vector subspaces of~$E$ is of particular interest;
it is called
the \emph{projective space} of~$E$.
\end{description}
\end{exer}
\end{small} 
\section{$C^r$-maps between manifolds}
\label{sec:2.2}
We now define $C^r$-maps between $C^r$-manifolds
and compile some basic properties and examples.
\begin{defn}\label{defnCrmf}
Let $M$ and $N$
be $C^r$-manifolds
modeled on locally convex spaces $E$ and $F$,
respectively.
A map $f\colon M\to N$ is called $C^r$
if $f$ is continuous and, for every chart $\phi$ of $M$ and every chart $\psi$ of
$N$, the map
\begin{equation}\label{uggly}
\psi\circ f \circ \phi^{-1}\colon E\, \supseteq \,
\phi(f^{-1}(U_\psi)\cap U_\phi)\to F
\end{equation}
is~$C^r$.
Note that the domain
$\phi(f^{-1}(U_\psi)\cap U_\phi)$ of the composition is an open subset
of~$V_\phi$ and hence also
open in~$E$.
We write $C^r(M,N)$ for the set
of all $C^r$-maps from~$M$ to~$N$.
\end{defn}
\begin{rem}\label{oneatlasenough}
Let $f\colon M\to N$ be a continuous map between $C^r$-manifolds
$(M,\cA_1)$ and $(N,\cA_2)$.
Let $\cB_1\sub \cA_1$ and $\cB_2\sub \cA_2$
be any atlases.
Using the technique
of ``insertion of charts,''
we see that $f$ is
$C^r$
if and only if
the map (\ref{uggly})
is $C^r$ for all $\phi\in \cB_1$ and $\psi\in \cB_2$.
Thus,
the $C^r$-property can be tested
using non-maximal atlases,
which is quite useful in practice.
It even suffices to check that,
for each $x\in M$,
there exists a chart $\phi\in \cA_1$
around~$x$ and a chart $\psi\in\cA_2$
around $f(x)$ such that
$\psi\circ f\circ \phi^{-1}$ is~$C^r$.
\end{rem}
Here are simple examples:
\begin{ex}\label{exmpCr1}
Let $E$ and $F$ be locally convex spaces
and $f\colon U\to F$ be a map
on an open subset $U\sub E$.
Then $f$ is $C^r$
as a map between open subsets of locally
convex spaces
(as in Chapter~\ref{chapcalcul})
if and only if $f$ is $C^r$ as a map
between $U$ and $F$,
considered as $C^r$-manifolds
modeled on $E$, resp., $F$,
equipped with the maximal $C^r$-atlas
containing the atlas $\{\id_U\}$ and $\{\id_F\}$,
respectively.
\end{ex}
Thus, our previous terminology is compatible
with the current one.
\begin{ex}\label{coprojcr}
Let $M_1$ and $M_2$
be $C^r$-manifolds
modeled on~$E_1$ and
$E_2$, respectively.
Then the coordinate projections
\[
\pr_i\colon M_1\times M_2\to M_i
\]
are $C^r$. In fact, for each chart $\phi_1\times \phi_2$
in the (non-maximal) atlas $\cA$ of $M_1\times M_2$
from Definition~\ref{exmprod}
and each chart $\psi$ of $M_1$,
we have
\[
\psi\circ \pr_1\circ (\phi_1\times \phi_2)^{-1}
=\psi\circ \phi_1^{-1}\circ \pi_1\,,
\]
where $\pi_1\colon E_1\times E_2\to E_1$
is the coordinate projection.
Being composed of $C^r$-maps,
this map is~$C^r$.
Thus $\pr_1$ is~$C^r$.
Similarly, $\pr_2$ is~$C^r$.
\end{ex}
Using the technique of ``insertion of charts,''
we see:
\begin{prop}\label{compmfds}
Let $M$, $N$ and $X$ be $C^r$-manifolds.
If $f\colon M\to N$ and $g\colon N\to X$
are $C^r$-maps,
then also
their composition
$g\circ f \colon M\to X$
is $C^r$.
\end{prop}
\begin{prf}
The map $g\circ f$ is continuous.
Given charts $\phi\colon U_\phi\to V_\phi$ of $M$ and
$\psi\colon U_\psi\to V_\psi$
of $X$,
let $x\in \phi((g\circ f)^{-1}(U_\psi)\cap U_\phi)$.
Let $\kappa$ be a chart of~$N$
around $f(\phi^{-1}(x))$.
Then
\[
\psi\circ (g\circ f)\circ\phi^{-1}=(\psi\circ g\circ \kappa^{-1})
\circ (\kappa\circ f\circ \phi^{-1})
\]
holds on some open neighborhood of~$x$
in~$V_\phi$,
and this is a $C^r$-map
as it is a composition of $C^r$-maps between
open subsets of topological vector spaces.
Thus
$\psi\circ (g\circ f)\circ \phi^{-1}$
is locally $C^r$ and hence~$C^r$.
Hence $g\circ f$ is~$C^r$.
\end{prf}
\begin{lem}\label{mapinprodmfd}
Let $M_1$, $M_2$, and $N$ be
$C^r$-manifolds.
Then a map $f=(f_1,f_2)\colon N\to M_1\times M_2$ is~$C^r$ if and only if
its components $f_1\colon N\to M_1$
and $f_2\colon N\to M_2$ are~$C^r$.
\end{lem}
\begin{prf}
If $f$ is $C^r$, then so are
$f_1=\pr_1\circ f$
and $f_2=\pr_2\circ f$,
being compositions of $C^r$-maps.
Conversely, assume that $f_1$ and $f_2$ are~$C^r$.
Consider
a chart~$\phi$ of~$N$ and
a chart of $M_1\times M_2$
of the form $\phi_1\times \phi_2$
(which it suffices to consider
by Remark~\ref{oneatlasenough}).
Then
$(\phi_1\times \phi_2)\circ f\circ \phi^{-1}
=(\phi_1\circ f_1\circ \phi^{-1},\phi_2\circ f_2\circ \phi^{-1})$
is $C^r$
(see Lemma~\ref{lemprod} and
Exercise~\ref{excreatoprod}),
as required.
\end{prf}
\begin{lem}\label{mapsubmfd}
Let $M$ be a $C^r$-manifold
and $N\sub M$ a submanifold. 
Then the inclusion map $\iota\colon N\to M$,
$x\mto x$ is $C^r$.
Furthermore, a map
$f\colon P\to N$ from a $C^r$-manifold
to~$N$
is $C^r$ if and only if $\iota\circ f\colon P\to M$
is $C^r$.
\end{lem}
\begin{prf}
Let $E$ and $F\sub E$
be the modeling spaces of~$M$ and $N$,
respectively.
For each chart $\phi\colon U_\phi\to V_\phi$
of~$M$ adapted to~$N$ (as in Definition~\ref{defsubm})
and the corresponding submanifold chart $\phi_N$ of~$N$,
the composition
$\phi\circ \iota\circ \phi_N^{-1}$ is the inclusion map
$V_\phi\cap F\to V_\phi$, which is $C^r$ as the restriction
of the continuous linear inclusion map $F\to E$.
Hence $\iota$ is $C^r$
(see Remark~\ref{oneatlasenough}).

If $f$ is $C^r$, then so is $\iota\circ f$,
by Proposition~\ref{compmfds}.
Conversely, assume that $\iota\circ f$ is~$C^r$.
Let $\phi$ and $\phi_N$ be as before and $\psi\colon U_\psi\to V_\psi$
be a chart for~$P$. Abbreviate $W:=\psi(U_\psi\cap f^{-1}(U_\phi))$.
Then $\phi\circ \iota\circ f\circ \psi^{-1} \colon W\to E$
is a $C^r$-map taking its values in the
closed vector subspace $F\sub E$,
and thus $(\phi\circ \iota\circ
f\circ \psi^{-1})|^F=\phi_N\circ f\circ \psi^{-1}\colon
W\to F$ is $C^r$, by Lemma~\ref{corestr}
(resp., Exercise~\ref{excreatosub}).
Hence $f$ is $C^r$.
\end{prf}
\begin{defn}\label{secdefdiff}
A \emph{$C^r$-diffeomorphism}
is a $C^r$-map $f\colon M\to N$
between $C^r$-manifolds such that there exists
a $C^r$-map $g\colon N\to M$ with
$g\circ f=\id_M$ and $f\circ g=\id_N$.
\end{defn}
In other words, $f$ is a bijective
$C^r$-map with
$f^{-1}\in C^r(N,M)$.
\begin{defn}\label{secdeflocdiff}
A mapping $f\colon M\to N$ between $C^r$-manifolds is called a
\emph{local $C^r$-diffeomorphism} if each point $x\in M$ has an open
neighborhood~$U$ in~$M$ such that
$f(U)$ is open in~$N$ and $f|_U\colon U\to f(U)$ is a $C^r$-diffeomorphism.
Then $f$ is $C^r$, in particular.
\end{defn}
If the choice of~$r$ is clear from the context, $C^r$-diffeomorphisms
and local $C^r$-diffeomorphisms
may simply be called \emph{diffeomorphisms} and \emph{local diffeomorphisms},
respectively.
\subsection*{Covering manifolds} 
Manifold structures can be pushed forward or pulled back along local homeo\-morphisms
under appropriate hypotheses.

Recall that a map $q\colon X\to Y$ between topological spaces is called a \emph{local homeomorphism}
if each $x\in X$ has  an open neighborhood $U\sub X$ such that
$q(U)$ is open in~$Y$ and $q|_U\colon U\to q(U)$ a homeomorphism.
\begin{thm}\label{coverdiff}
Let $M$ and $N$ be Hausdorff spaces and $q\colon M\to N$ be a local homeomorphism.
Let $E$ be a locally convex space over $\K\in\{\R,\C\}$ and $r\in \N_0\cup\{\infty,\omega\}$.
Then the following holds.
\begin{description}[(D)]
\item[\rm(a)] 
  If $N$ is a $C^r_\K$-manifold modeled on~$E$, then there is a unique
  structure of a $C^r_\K$-manifold on~$M$
modeled on~$E$ for which the map $q$ is a local $C^r_\K$-diffeomorphism.
\item[\rm(b)]
If $q$ is surjective and $M$ is a $C^r_\K$-manifold modeled on~$E$, then the following conditions are equivalent:
\begin{description}[(D)]
\item[\rm(i)]
$N$ admits a $C^r_\K$-manifold structure modeled on~$E$ which makes~$q$ a local $C^r_\K$-diffeo\-mor\-phism;
\item[\rm(ii)]
For all $x,y\in M$ with $q(x)=q(y)$, there exists a $C^r_\K$-diffeomorphism $h\colon U\to V$ between an open $x$-neighborhood $U\sub M$ and an open $y$-neighborhood $V\sub M$ such that $q\circ h=q|_U$.
\end{description}
The $C^r_\K$-manifold structure on~$N$ described in {\rm(i)} is unique, if it exists.
\end{description}
\end{thm}
\begin{prf}
(a) Uniqueness: Let $\cA_1$ and $\cA_2$ be maximal $C^r_\K$-atlases on~$M$
modeled on~$E$ which turn~$q$ into a local $C^r_\K$-diffeomorphism.
Given $x\in M$, let $U_j\sub M$ be an open $x$-neighborhood for $j\in \{1,2\}$
such that $q(U_j)\sub N$
is open and $q|_{U_j}\colon U_j\to q(U_j)$ is a $C^r_\K$-diffeomorphism for the manifold structure
induced by $(M,\cA_j)$ on~$U_j$.
After shrinking~$U_1$ and $U_2$, we may assume that there is a chart $\phi\colon P\to Q\sub E$ of~$N$
such that $q(U_1)\cup q(U_2)\sub P$. Then $U:=U_1\cap U_2$ is an open $x$-neighborhood in~$M$
and $\phi_x:=\phi\circ q|_U\colon U\to \phi(q(U))$ is an element of both $\cA_1$ and $\cA_2$.
Thus both $\cA_1$ and $\cA_2$ coincide with the maximal $C^r_\K$-atlas containing the $C^r_\K$-atlas $\{\phi_x\colon x\in M\}$.

Existence: Each $x\in M$ has an open neighborhood $U_x\sub M$ such that $q(U_x)$ is open in~$N$
and $q|_{U_x}\colon U_x\to q(U_x)$ is
a homeomorphism. After shrinking~$U_x$, we may assume that $q(U_x)\sub P_x$ for some chart $\psi_x\colon P_x\to Q_x\sub E$ of~$N$. After shrinking~$P_x$ and~$Q_x$, we may assume that $P_x=q(U_x)$.
Then $\phi_x:=\psi_x\circ q|_{U_x}\colon U_x \to \psi_x(q(U_x))$ is a homeomorphism
between open sets. Now $\{\phi_x\colon x\in M\}$
is a $C^r_\K$-atlas, noting that
\[
\phi_x\circ\phi_y^{-1}=\psi_x\circ q|_{U_x}\circ (q|_{U_y})^{-1}\circ \psi_y^{-1}
\]
coincides with $\psi_x\circ \psi_y^{-1}$ on its open domain $\phi_y(U_x\cap U_y)$.
We endow~$M$ with the corresponding maximal $C^r_\K$-atlas~$\cA$.
Then~$q$ is a local $C^r_\K$-diffeomorphism on $(M,\cA)$, since $\psi_x\circ q|_{U_x}=\phi_x\in\cA$
implies that $\phi_x$, and hence also $q|_{U_x}^{q(U_x)}=\psi_x^{-1}\circ \phi_x$,
is a $C^r_\K$-diffeomorphism.

(b) Assuming~(i), let us deduce~(ii). Let $x,y\in M$ with $q(x)=q(y)$.
Since~$q$ is a local diffeomorphism, we find an open $x$-neighborhood $U\sub M$ and an open $y$-neighborhood $V\sub M$ such that $q(U)$ and $q(V)$ are open in~$N$ and the mappings
$q|_U\colon U\to q(U)$ and $q|_V\colon V\to q(V)$ are $C^r_\K$-diffeomorphisms.
After replacing $U$ with $(q|_U)^{-1}(q(U)\cap q(V))$ and $V$
with $(q|_V)^{-1}(q(U)\cap q(V))$, we may assume that $q(U)=q(V)$.
Now $h:=(q|_V)^{-1}\circ q|_U\colon U\to V$ is a $C^r_\K$-diffeomorphism such that
$q\circ h=q|_U$.

Conversely, assume that (ii) holds. Given $y\in M$,
let $\psi_y\colon P_y\to V_y$ be
a chart of~$M$ with $y\in P_y$.
By hypothesis, there is an open $y$-neighborhood $W_y\sub M$ such that $U_y:=q(W_y)$ is open
in~$N$ and $q|_{W_y}\colon W_y\to q(W_y)$ is a homeomorphism. After shrinking~$W_y$,
we may assume that $W_y\sub P_y$. After shrinking~$P_y$ and~$V_y$,
we may assume that $W_y=P_y$. Now
\[
\phi_y:=\psi_y\circ (q|_{W_y})^{-1}\colon U_y\to V_y
\]
is a homeomorphism between open subsets of~$N$ and~$E$, respectively.
To see that $\{\phi_y\colon y\in M\}$ is a $C^r_\K$-atlas for~$N$,
note that $x\in N$ satisfies $y\in U_y$ if we pick $y\in M$ with $q(y)=x$.
Let $y,y'\in M$ and set $x:=q(y)$, $x':=q(y')$;
we show that $\phi_{y'}\circ\phi_y^{-1}$ is $C^r_\K$
on an open neighborhood of a given point~$v$ in its domain $\phi_y(U_y\cap U_{y'})$
(whence $\phi_{y'}\circ \phi_y^{-1}$ is $C^r_\K$).
Let $z:=\phi_y^{-1}(v)\in U_y\cap U_{y'}$.
There exist $a\in W_y$ and $b\in W_{y'}$ such that $q(a)=q(b)=z$.
By hypothesis~(ii), there exists an open $a$-neighborhood $U\sub M$, an open $b$-neighborhood $V\sub M$
and a $C^r_\K$-diffeomorphism $h\colon U\to V$ such that $q\circ h=q|_U$.
Now
\[
(q|_{W_y})^{-1}(U_y\cap U_{y'})\cap U
\]
is an open neighborhood of~$a$ in~$W_y$ and hence in~$M$,
whence also its connected component~$C$ containing~$a$ is open in~$M$ (as $M$ is a manifold and hence locally connected). Now
\[
g\colon C\to M,\quad c\mto (q|_{W_{y'}})^{-1}(q(c))
\]
and $h|_C\colon C\to M$ are continuous maps such that $q\circ g=q|_C$
and $q\circ h|_C=q|_C$. As $g(a)=h(a)=b$, we have
$g=h|_C$ by Lemma~\ref{uni-lft}.
As a consequence,
\[
\phi_{y'}\circ \phi_y^{-1}|_{\psi_y(C)}=\psi_{y'}\circ (q|_{W_{y'}})^{-1}\circ q|_{W_y}\circ \psi_y^{-1}|_{\psi_y(C)}
=\psi_{y'}\circ h|_C\circ \psi_y^{-1}|_{\psi_y(C)}
\]
is $C^r_\K$ on the open neighborhood $\psi_y(C)$ of~$v$ in~$V_y$.

Let $\cA$ be the maximal $C^r_\K$-atlas on~$N$ which contains $\{\phi_y\colon y\in M\}$.
For each $y\in M$, the image $q(W_y)=U_y$ is open in~$N$
and $q|_{W_y}^{U_y}=\phi_y^{-1}\circ \psi_y$ is a $C^r_\K$-diffeomorphism.
Thus $q\colon M\to (N,\cA)$ is a local $C^r_\K$-diffeomorphism.\vspace{.3mm}

Uniqueness: For $j\in\{1,2\}$, let $\cA_j$ be a maximal $C^r_\K$-atlas
on~$N$ modeled on~$\!E$
such that $q\colon \!M\!\to \!(N,\cA_j)$ is a local $C^r_\K$-diffeomorphism.
Given $x\!\in\! N$, pick $y\in M$ such that $q(y)=x$.
For $j\in \{1,2\}$, there exists an open $y$-neighborhood $U_j\sub M$ such that $q(U_j)$ is open in~$N$
and $q|_{U_j}\colon U_j\to q(U_j)$ is a $C^r_\K$-diffeomorphism,
considering $q(U_j)$ as
an open submanifold of $(N,\cA_j)$. We may assume that $W_x:=U_1=U_2$,
after replacing both $U_1$ and $U_2$ with $U_1\cap U_2$.
After shrinking~$W_x$, we may assume that~$W_x$ is the domain
of a chart $\psi_x\colon W_x\to V_x\sub E$ of~$M$.
Then $\phi_x:=\psi_x\circ (q|_{W_x})^{-1}\colon q(W_x)\to V_x$
is a $C^r_\K$-diffeomorphism defined on the open submanifold $q(W_x)$ of $(M,\cA_j)$ for $j\in\{1,2\}$,
and thus $\phi_x\in \cA_1\cap\cA_2$.
As a consequence, both $\cA_1$ and $\cA_2$ coincide with the maximal $C^r_\K$-atlas on~$M$ which contains $\{\phi_x\colon x\in N\}$. Thus $\cA_1=\cA_2$.
\end{prf}
Let us recall some basic terminology concerning covering maps.
\begin{defn}
A mapping $q\colon M\to N$ between topological spaces
is called a \emph{covering map} if $q$ is surjective and each $x\in N$ has an open neighborhood
$V\sub N$ which is \emph{evenly covered} in the sense that
\[
q^{-1}(V)=\bigcup_{j\in J} W_j
\]
for a family $(W_j)_{j\in J}$ of disjoint open subsets $W_j\sub M$ such that $q(W_j)=V$ and
$q|_{W_j}\colon W_j\to V$
is a homeomorphism, for all $j\in J$.
Given a covering map $q\colon M\to N$, a homeomorphism $h\colon M\to M$ is called
a \emph{deck transformation} if $q\circ h=q$. The set $\Deck(q)$
of all deck transformations is a group under composition.
If $\Deck(q)$ acts transitively on $q^{-1}(\{x\})$
for all $x\in N$, then $q$ is called a \emph{normal} covering.
\end{defn}
We record an important special case of Theorem~\ref{coverdiff}.
\begin{cor}\label{cov-coverdiff}
Let $q\colon M\to N$ be a covering map between Hausdorff spaces,
$r\in\N_0\cup\{\infty,\omega\}$ and~$E$ be a locally convex space over $\K\in\{\R,\C\}$.
\begin{description}[(D)]
\item[\rm(a)]
If $N$ is a $C^r_\K$-manifold modeled on~$E$, then~$M$ admits a unique
$C^r_\K$-manifold structure modeled on~$E$
that makes~$q$ a local $C^r_\K$-diffeomorphism.
\item[\rm(b)]
If $M$ is a $C^r_\K$-manifold modeled on~$E$, the covering~$q$ is normal
and each $h\in \Deck(q)$
is a $C^r_\K$-diffeomorphism,
then~$N$ admits a unique $C^r_\K$-manifold structure modeled on~$E$
which turns~$q$ into a local $C^r_\K$-diffeomorphism. \qed
\end{description}
\end{cor}
\subsection*{Pure manifolds and sets of modeling spaces}
So far, we considered $C^r$-manifolds modeled on a single
locally convex space~$F$, so-called \emph{pure}
manifolds.
For some purposes (like manifold structures on manifolds
of mappings treated in Section~\ref{sec-mfdmps}),
this is not sufficient.
\begin{defn}\label{non-pure}
Given a Hausdorff space~$M$ and a set $\cE$ of locally
convex spaces, we now consider charts
$(E,\phi)$ where $E\in \cE$ and $\phi\colon U_\phi\to V_\phi$
is a homeomorphism from an open subset $U_\phi$ of~$M$
onto an open subset $V_\phi$ of~$E$.
We usually write $E_\phi:=E$;
if $E_\phi$ is understood,
we simply call $\phi\colon U_\phi\to V_\phi\sub E_\phi$
a chart.
Given $r\in \N_0\cup\{\infty,\omega\}$, we say that
charts $\phi\colon U_\phi\to V_\phi\sub
E_\phi$ and $\psi\colon U_\psi\to V_\psi\sub E_\psi$
are $C^r$-compatible if $\psi\circ \phi^{-1}\colon \phi(U_\phi\cap U_\psi)\to \psi(U_\phi\cap U_\psi)$
is a $C^r$-diffeomorphsm. A set of $C^r$-compatible charts
whose domains cover~$M$
is a called a $C^r$-atlas on~$M$ modeled on~$\cE$;
a \emph{$C^r$-manifold modeld on~$\cE$}
is a Hausdorff space, together with a maximal $C^r$-atlas modeled
on~$\cE$. If $\cE$ is a singleton $\{E\}$,
then $M$ is a $C^r$-manifold modeled on~$E$ as before.
\end{defn}
Using charts $\phi\colon U_\phi\to V_\phi\sub E_\phi$
and $\psi\colon U_\psi\to V_\psi\sub F_\psi$
in Definition~\ref{defnCrmf}
instead of $E$-charts and $F$-charts with fixed locally convex spaces $E$ and~$F$,
we can define $C^r$-maps between $C^r$-manifolds $M$ and $N$
modeled on sets $\cE$ and $\cF$ of locally convex spaces.
Also further concepts can be adapted to non-pure manifolds
in a straightforward fashion, which we leave to the reader.
Submanifolds are defined as follows:
\begin{defn}
Let $r\in\N_0\cup\{\infty,\omega\}$
and $M$ be a $C^r$-manifold
modeled on a set $\cE$ of locally convex spaces.
A subset $N\sub M$ is called
a (not necessarily pure$)$ submanifold if, for each $x\in N$,
there exists a chart $\phi\colon U_\phi\to V_\phi\sub E_\phi$
and a closed vector subspace $F_\phi\sub E_\phi$
such that $\phi(N\cap U_\phi)=F_\phi \cap V_\phi$.
\end{defn}
We endow $N$ with the maximal $C^r$-atlas modeled
on the set $\cF$ of the $F_\phi$ which occur,
that contains all of charts $\phi|_{N\cap U_\phi}\colon U_\phi\cap N\to V_\phi\cap F_\phi\sub F_\phi$.
\begin{convent}
The manifolds considered in this book
are pure manifolds, unless the contrary is stated.
\end{convent}
Non-pure manifolds are needed only in Sections~\ref{sec-mfdmps},
\ref{sec-box-mfd}, and \ref{sec-mfd-map-sigma}.
%
We also consider not necessarily pure manifolds
when we discuss submersions, immersions
and embeddings in Section~\ref{sec:2.3},
as constructions with preimages
may not give rise to pure submanifolds
in infinite dimensions.
\begin{small}
\subsection*{Exercises for Section~\ref{sec:2.2}} 
\label{Exer3.2}
\begin{exer}\label{exc-chartdiffeo}
Let $M$ be a $C^r$-manifold modeled on a locally convex space~$E$,
where $r\in\N_0\cup\{\infty,\omega\}$.
Let $\phi\colon U\to V$ be a map between open subsets $U\sub M$ and $V\sub E$.
Show that $\phi$ is a chart for~$M$ if and only if~$\phi$ is a $C^r$-diffeomorphism.
\end{exer}

\begin{exer}
Given $\K\in\{\R,\C\}$ and $r\in\N\cup\{\infty,\omega\}$,
let $M$ and~$N$ be $C^r_\K$-manifolds and
$q\colon M\to N$ be a local $C^r_\K$-diffeomorphim
which is a normal covering map.
Show that each $h\in \Deck(q)$ is a $C^r_\K$-diffeomorphism.
\end{exer}

\begin{exer}
  \label{exer:2.2.1} Let $G$ be a Lie group modeled
on a locally convex space~$E$;
i.e., $G$ is a group and a smooth manifold,
and both the group multiplication
$\mu\colon G\times G\to G$, $\mu(x,y)=xy$ and
the inversion map $\eta\colon G\to G$, $x\mto x^{-1}$
are smooth.
\begin{description}[(D)]
\item[(a)]
Show that, for each $g\in G$,
the left translation map
$\lambda_g\colon G\to G$, $\lambda_g(x):=gx$
and the right translation map
$\rho_g\colon G\to G$, $\rho_g(x):=xg$ are
$C^\infty$-diffeomorphisms.
\item[(b)]
Let $\phi\colon U_\phi\to V_\phi\sub E$
be a chart for~$G$,
and $g\in G$ be any group element.
Show that also $\phi_g \colon g U_\phi\to V_\phi$,
$x\mto \phi(g^{-1}x)$ is a chart.
If $\phi$ is a chart around~$\one$,
then $\phi_g$ is a chart around~$g$.
\item[(c)]
Let $F\sub E$ be a
closed
vector subspace.
Show that a subgroup $H\sub G$ is a submanifold of~$G$
modeled on $F$ if and only if there
is a chart $\phi\colon U_\phi\to V_\phi$
of~$G$ around~$\be$ such that
$\phi(U_\phi\cap H)=V_\phi\cap F$.
\end{description}
\end{exer}

\begin{exer}
  \label{exer:2.2.2} 
Specify a suitable chart of~$\GL_2(\R)$
around~${\bf 1}$ as in Exercise~\ref{exer:2.2.1}(c) 
and deduce that $\SL_2(\R)$ (as in Exercise~\ref{exer:2.1.7}) 
is a submanifold of~$\GL_2(\R)$
modeled on~$\fsl_2(\R)=\{x\in M_2(\R)\colon \tr(x)=0\}$.
Hint: The matrix exponential function
$\exp\colon M_2(\R)\to \GL_2(\R)$
satisfies $\det (\exp x)=e^{\tr x}$
(cf.\ also Example~\ref{explogban}).
\end{exer}

\begin{exer}
  \label{exer:2.2.3} Show that if $G$ is a Lie group
and $H\sub G$ a subgroup
and a submanifold,
then the submanifold structure
makes~$H$ a Lie group.
Deduce that $\SL_2(\R)$
is a Lie group.
\end{exer}

\begin{exer}
Give $\{0,1\}$ the discrete topology
and endow $M:=\{0,1\}\times\R$ with the product topology.
Let $\theta\colon \R\to\R$ be a homeomorphism which is not a $C^\infty$-diffeomorphism.
Then
\[
\phi\colon \{0\}\times\R\to\R,\;\; (0,x)\mto x
\]
and $\psi\colon \{1\}\times\R$, $(1,x)\mto \theta(x)$ are homeomorphisms between
open subsets of~$M$ and~$\R$, and $\{\phi,\psi\}$ is a smooth atlas for~$M$.
Let~$\cA$ be the corresponding maximal smooth atlas.
Show that $\pr_2\colon M\to\R$, $(j,x)\mto x$ is a covering map.
Show that there is no smooth manifold structure~$\cB$ on~$\R$ turning
$\pr_2\colon (M,\cA)\to (\R,\cB)$ into a local $C^\infty$-diffeomorphism.
\end{exer}
\end{small}
\section{Tangent bundles and tangent maps}
\label{sec:2.3}%
In this section, we define the tangent
space $T_p(M)$ of a $C^r$-manifold~$M$
at a point~$p\in M$,
if $r\geq 1$.
Every $C^r$-map $f\colon M\to N$ between
manifolds will give rise
to a  ``tangent map'' $T_p(f)\colon T_p(M)\to T_{f(p)}(N)$,
which is a certain analog
of the derivative $f'(x)$ of
a $C^1$-map between locally convex spaces.
Before we give a more abstract
definition of tangent spaces
(which, later, will also work for manifolds
with boundary), 
we start with a discussion of geometric tangent spaces,
which only make sense for ordinary
manifolds (without boundary), as currently considered.
The latter motivate the formalism
and provide the necessary geometric intuition.
We fix $r\in \N\cup\{\infty,\omega\}$.
\subsection*{Geometric tangent spaces}
If $M\sub E$ is a submanifold
of a locally convex space,
then a vector subspace~$V$ of~$E$
can be associated to each
$p\in M$,
namely the set
of all velocity vectors $v=\gamma'(0)$
of
$C^1$-curves $\gamma\colon \, ]{-\ve},\ve[\;
\to M$ which ``pass through $p$''
in the sense that $\gamma(0)=p$
(see Exercise~\ref{exer:2.3.5}(a)).
Then $p+V$ is an affine vector subspace
of~$E$ passing through~$p$,
which is tangent to~$M$ in~$p$.\\[2.5mm]
For example, if $E=\R^n$,
$M=\bS^{n-1}$ and $p\in \bS^{n-1}$,
then $p+V=p+p^\perp$,
where $p^\perp=\{x\in \R^n\colon \langle p,x\rangle=0\}$
(see Exercise~\ref{exer:2.3.5}(b)).\\[2.5mm]
Unfortunately,
since the ambient locally convex space~$E$
is needed to calculate the derivatives $\gamma'(0)$,
the preceding idea cannot be used to
define tangent spaces for general
manifolds (beyond submanifolds
of vector spaces).
To find a replacement for velocity vectors
in the general case,
we note that
$v:=\gamma'(x)$
determines (and is uniquely determined by)
the set $[\gamma]$ of all $C^1$-curves
$\eta$ passing through~$p$
such that $\eta'(0)=v$, in the case $M\sub E$.
This motivates the following definition.
\begin{defn}\label{def-geom-tang}
Let $M$ be a $C^r$-manifold
modeled on a locally convex space~$E$,
and $p\in M$.
For $C^1$-curves $\gamma,\eta$
passing through~$p$,
we write $\gamma\sim \eta$~if
\begin{equation}\label{specsim}
(\phi\circ \gamma)'(0)\; =\; (\phi\circ\eta)'(0)
\end{equation}
for some chart~$\phi$ of~$M$
around~$p$.
Then (\ref{specsim})
holds for each chart $\psi$
around~$p$ in place of~$\phi$, using that
$(\psi\circ\gamma)'(0)
=d(\psi\circ\phi^{-1})(\phi(p),(\phi\circ\gamma)'(0))$
by the Chain Rule.
Hence equality in (\ref{specsim})
can be checked with any given chart,
entailing that $\sim$ is an equivalence relation
on the set of all $C^1$-curves $\gamma$
passing through~$p$.
The equivalence class $[\gamma]$
of a $C^1$-curve $\gamma$ passing through~$p$
is called a
\emph{geometric tangent vector} (of~$M$ at~$p$).
The \emph{geometric tangent space of~$M$ at~$p$}
is defined as the set
$\cT_p(M)$ of all geometric tangent vectors
of~$M$ at~$p$.
\end{defn}
The following lemma facilitates $\cT_p(M)$
to be turned into a locally convex space
isomorphic to~$E$, in a canonical way.
\begin{lem}\label{changehphis}
\begin{description}[(D)]
\item[\rm (a)]
Given a chart $\phi$ of~$M$ around~$p$,
set $x_\phi:=\phi(p)$. Then
\[
h_\phi\colon E\to \cT_p(M)\,,\quad
h_\phi(y)\, :=\, [t\mto \phi^{-1}(x_\phi+ty)]
\]
is a bijection,
with inverse $\ell_\phi\colon \cT_p(M)\to E$, $[\gamma]\mto
(\phi\circ \gamma)'(0)$.
\item[\rm (b)]
For all charts $\phi,\psi$ of~$M$ around~$p$, we have
$h_\psi^{-1}\circ h_\phi = (\psi\circ\phi^{-1})'(x_\phi)$,
which is an automorphism of the topological
vector space~$E$.
\end{description}
\end{lem}
\begin{prf}
(a) It is clear from the definition
of tangent vectors that $\ell_\phi$ is injective.
Given $y\in E$, we have
$\ell_\phi (h_\phi (y))=
\frac{d}{dt}\big|_{t=0}\phi(\phi^{-1}(x_\phi+ty))=y$.
Hence $\ell_\phi$ is also surjective
and $h_\phi=\ell_\phi^{-1}$.

(b) Given $y\in E$,
we have 
\[ h_\psi^{-1}(h_\phi(y))=
\frac{d}{dt}\Big|_{t=0}\psi(\phi^{-1}(x_\phi+ty)) = (\psi\circ\phi^{-1})'(x_\phi)y.
\qedhere\]
\end{prf}
We equip $\cT_p(M)$ with the locally convex topological
vector space structure making the bijection
$h_\phi$ an isomorphism of topological vector spaces
for some $($and hence each$)$ chart $\phi$ of~$M$ around~$p$.
\subsection*{Abstract definition of tangent spaces and tangent bundles}
We now give an alternative definition of tangent spaces
$T_p(M)$,
in a purely formal way.
To motivate our
definition,
let us return to a geometric tangent vector $v\in \cT_p(M)$ at $p\in M$.
By Lemma~\ref{changehphis}\,(a),
for each chart $\phi$ of~$M$ around~$p$,
we have $v=h_\phi(v_\phi)$
for a uniquely determined vector $v_\phi\in E$.
We can thus identify $v$ with the family
$(v_\phi)_{\phi\in \cA_p}$.
Given $\phi,\psi\in \cA_p$,
Lemma~\ref{changehphis}(b) also tells us that
\begin{equation}\label{pretraf2}
v_\psi\; =\; d(\psi\circ\phi^{-1})(x_\phi,v_\phi)\, ,
\end{equation}
where $x_\phi:=\phi(p)$.
Conversely, every family $(v_\phi)_{\phi\in \cA_p}$
of vectors $v_\phi\in E$ showing the
transformation behavior described in (\ref{pretraf2})
arises from a unique geometric tangent vector~$v\in \cT_p(M)$
which is given by $v=h_\phi(v_\phi)$
for one (and hence any) $\phi\in \cA_p$.
We can therefore identify geometric tangent vectors
with families of vectors satisfying
a specified transformation behavior.
This is the basic idea underlying the abstract
definition of tangent spaces, which we describe
in this section.
\begin{defn}\label{defTefX}
Let $E$ and $F$ be locally convex spaces and
$U\sub E$ and $V\sub F$ be open
subsets.
If $f\colon U\to V$ is $C^1$,
we define
\[
Tf\colon U\times E\to V\times F\,,\quad
Tf(x,y)\,:=\, (f(x),df(x,y))\,.
\]
\end{defn}
If $f$ is $C^r$, then $Tf$ is $C^{r-1}$.
We also write $T(f)$ instead of $Tf$.
\begin{rem}\label{remsimpachn}
Using the notation just introduced, the Chain Rule
can be reformulated in a convenient way:
Given $C^1$-maps $f\colon U\to V$ and $g\colon V\to W$,
where $U\sub E$, $V\sub F$ and $W\sub H$ are open
subsets of locally convex spaces, we have
\begin{equation}\label{simplchain1}
d(g\circ f)\;=\; dg\circ Tf
\end{equation}
and hence
\begin{equation}\label{simplchain2}
T(g\circ f)\;=\; Tg\circ Tf\,.
\end{equation}
\end{rem}
Using this formalism,
(\ref{pretraf2})
can be rewritten as
\[ (x_\psi,v_\psi)=T(\psi\circ\phi^{-1})(x_\phi,v_\phi).\] 
\begin{defn}\label{taspa}
Let $M$ be a $C^r$-manifold
and $p\in M$.
We define a relation $\sim$ on the set
of triples $(\phi, x, v)$,
where
$\phi\colon U_\phi\to V_\phi$ is a chart of~$M$ around~$p$,
$x:=\phi(p)\in V_\phi$ and $v\in E$,
as follows: Given two such triples
$(\phi,x,v)$ and $(\psi,y,w)$,
we declare that $(\phi,x,v)\sim (\psi,y,w)$
if and only if $(y,w)=T(\psi\circ \phi^{-1})(x,v)$.
It readily follows from (\ref{simplchain2})
that $\sim$ is an equivalence
relation. The $\sim$-equivalence class
$[\phi, x,v]$ of $(\phi, x,v)$
is called a
{\em tangent vector of~$M$ at~$p$}.
The set $T_p(M)$ (or $T_pM$)
of all tangent vectors of~$M$ at~$p$ is called
the {\em tangent space of~$M$ at~$p$.}
Given $p\in M$, we pick a chart $\phi$ of~$M$ around $p$
and give $T_pM$
the unique locally convex topological vector
space structure which makes the bijection
\begin{equation}\label{secstage}
k_\phi\colon E\to T_pM\,,\quad
k_\phi(v)\, :=\, T\phi^{-1}(x,v)\,:=\, [\phi,x,v]
\end{equation}
an isomorphism of topological vector spaces,
where $x:=\phi(p)$.
\end{defn}
Since
\begin{equation}\label{cmbto}
k_\psi^{-1}\circ k_\phi\, =\, (\psi\circ \phi^{-1})'(x)
\end{equation}
is an isomorphism of topological vector spaces,
the vector topology on $T_pM$ is
well defined, independent of the choice of~$\phi$.
\begin{rem}\label{tangonopsub}
If $U\sub M$ is an open submanifold and
$p\in U$, then every chart $\phi$
for~$U$ around~$p$ also is a chart
for~$M$, entailing that $T_pU\to T_pM$,
$[\phi,x,y]\mto [\phi,x,y]$
is an isomorphism of topological
vector spaces. We shall frequently
identify $T_pU$ with $T_pM$, without further
mention.
\end{rem}
\begin{defn}\label{tabu}
Let $(M,\cA)$ be a
$C^r$-manifold
modeled on a locally convex space~$E$.
We call $TM:=\bigcup_{p\in M}T_pM$
the \emph{tangent bundle} of~$M$.
The tangent spaces $T_pM$ being mutually
disjoint, the map $\pi_{TM}\colon TM\to M$
sending a tangent vector
$v\in T_pM$ to $\pi_{TM}(v):=p$
is well defined;
it is called the \emph{bundle projection}.
We equip $TM$ with the final topology
with respect to the family $(T\phi^{-1})_{\phi\in \cA}$
of the mappings
\[
T\phi^{-1}\colon V_\phi\times E\to TM\,,\quad
(T\phi^{-1})(x,y)\,:=\, [\phi,x,y]\,.
\]
We\hspace*{-.07mm}
shall\hspace*{-.07mm}
see\hspace*{-.07mm}
presently\hspace*{-.07mm}
that\hspace*{-.07mm}
$TU_\phi=\pi_{TM}^{-1}(U_\phi)$\hspace*{-.07mm}
is\hspace*{-.07mm}
open\hspace*{-.07mm}
in\hspace*{-.07mm}
$TM$\hspace*{-.07mm}
for\hspace*{-.07mm}
all\hspace*{-.07mm}
$\phi\in \cA$, that
\begin{equation}\label{defTphi}
T\phi\, :=\, (T\phi^{-1})^{-1}\colon TU_\phi\to V_\phi\times E
\end{equation}
is a homeomorphism for each $\phi\in \cA$, and that
$\cB:=\{T\phi\colon \phi\in \cA\}$ is a $C^{r-1}$-atlas
for~$TM$ modeled on~$E\times E$.
We equip $TM$ with the maximal $C^{r-1}$-atlas
containing~$\cB$.
In this way, $TM$ becomes a $C^{r-1}$-manifold.
\end{defn}
We need to check various details.
\begin{lem}\label{chckdt}
In the situation of {\rm Definition~\ref{tabu}},
the following holds:
\begin{description}[(D)]
\item[\rm (a)]
For all $\phi,\psi\in \cA$, the diagram
\[
\begin{array}{ccc}
T(U_\phi\cap U_\psi)
&\mapright{\id} & T(U_\phi\cap U_\psi)\\[.6mm]
\mapdown{T\phi} &&\mapdown{T\psi}\\
\phi(U_\phi\cap U_\psi)\times E &
\mapright{T(\psi\circ \phi^{-1})} &
\psi(U_\phi\cap U_\psi)\times E
\end{array}
\]
is commutative
$($where $T\phi$ and $T\psi$
have to be restricted to $T(U_\phi\cap U_\psi))$.
\item[\rm (b)]
$TU_\phi$ is open in $TM$, for each
chart $\phi\colon U_\phi\to V_\phi$ of~$M$.
\item[\rm (c)]
$T\phi\colon TU_\phi\to V_\phi\times E$
is a homeomorphism.
\item[\rm (d)]
$\cB:=\{T\phi\colon \phi\in \cA\}$ is a $C^{r-1}$-atlas.
\item[\rm (e)]
$\pi_{TM}\colon TM \to M$ is a $C^{r-1}$-map.
\item[\rm (f)]
$TM$ induces on each tangent space $T_pM$
its natural vector topology.
\end{description}
\end{lem}
\begin{prf}
(a) is essentially a reformulation of (\ref{cmbto}).

(b) Given $\phi\in \cA$, for each $\psi\in \cA$ we have
the preimage $(T\psi^{-1})^{-1}(TU_\phi)=
T\psi (TU_\phi\cap TU_\psi)=
\psi(U_\phi\cap U_\psi)\times E$,
which is open in $V_\psi\times E$.
Hence $TU_\phi$ is open in $TM$,
by definition of the final topology on~$TM$.

(c) $T\phi^{-1}$ being continuous by definition
of the topology on~$TM$, we deduce that
$T\phi$ is an open map.
To see that $T\phi$ is continuous,
let $U\sub V_\phi\times E$ be open.
Let $\psi\in \cA$.
Then $W:=U\cap \, (\phi(U_\phi\cap U_\psi)\! \times \! E)$
is open in $V_\phi\times E$, whence
\begin{eqnarray*}
(T\psi^{-1})^{-1}((T\phi)^{-1}(U))
& = & T\psi(TU_\psi\cap (T\phi)^{-1}(U))\\
& = &T\psi((T\phi)^{-1}(W))
\; = \; T(\psi\circ\phi^{-1})(W)
\end{eqnarray*}
is open in $V_\psi\times E$,
the map $T(\psi\circ \phi^{-1})\colon
\phi(U_\phi\cap U_\psi)\times E \to \psi(U_\phi\cap U_\psi)\times E$
being a homeomorphism from an open subset
of $V_\phi\times E$ onto an open subset
of $V_\psi\times E$.
Hence $(T\phi)^{-1}(U)$ is open in~$TM$.

(d) By (b) and (c), each $T\phi$ is a homeomorphism
from an open subset of $TM$ onto an open subset
of $E\times E$
and thus $T\phi$ is an $E\times E$-chart.
It is clear that the sets $TU_\phi$ cover~$TM$,
and (a) shows that the transition
maps are of the form
$T\psi\circ (T\phi)^{-1}=T(\psi\circ\phi^{-1})$
and hence~$C^{r-1}$.

(e) $\phi\circ\pi_{TM}\circ T\phi^{-1}\colon V_\phi\times E\to V_\phi$
is the projection onto the first component and hence $C^{r-1}$, for each $\phi\in\cA$.

(f) is immediate from (c) and the definition of the
vector topology in Definition~\ref{taspa}.
\end{prf}
\begin{ex}\label{anothtriv}
Let $E$ be a locally convex space
and $U\sub E$ be an open subset.
Then $\id_U$ is a global chart for~$U$,
whence
$T(\id_U)\colon TU\to U\times E$
(defined as in (\ref{defTphi}))
is a global chart for $TU$ and hence
a $C^r$-diffeomorphism
for each $r\in \N\cup\{\infty,\omega\}$.
We shall frequently identify $TU$
with $U\times E$ by means of $T(\id_U)$,
and suppress $T(\id_U)$
in the formulas. The abuse will always
be clear from the context.
Under this identification, $\pi_{TU}$
corresponds to the projection
$\pr_1\colon U\times E\to U$
onto the first component, i.e.,
the following diagram commutes:
\[
\begin{array}{rcl}
TU \; &\mapright{T(\id_U)} & U\times E\\
\mapdown{\pi_{TU}}\; & &\;\;\;  \mapdown{\pr_1} \\
U\  &\mapright{\id_U} & \ \ U.
\end{array}
\]
\end{ex}
\begin{rem}\label{agnsotrv}
Similarly, for each chart $\phi\colon U_\phi\to V_\phi$
of a $C^r$-manifold~$M$ modeled on~$E$
we have the commutative diagram
\[
\begin{array}{rcl}
TU_\phi  &\mapright{T\phi} & V_\phi\times E\\
\mapdown{\pi_{TM}|_{TU_\phi}}\;\, & &\;\;\;\;  \mapdown{\pr_1}\\
U_\phi &\mapright{\phi} & \ \  \; V_\phi\,,
\end{array}
\]
where $\pr_1\colon V_\phi\times E\to V_\phi$
is the projection onto the first component.
\end{rem}
\subsection*{Tangent maps}
\begin{defn}\label{deftangmap}
Let $f\colon M\to N$ be a $C^r$-map
between $C^r$-manifolds.
Then there exists a uniquely determined map
$Tf\colon TM\to TN$ (also denoted $T(f)$)
with the following properties:
\begin{description}[(D)]
\item[(a)]
$Tf(T_pM)\sub T_{f(p)}N$ for each $p\in M$,
i.e., $\pi_{TN}\circ Tf=f\circ \pi_{TM}$.
\item[(b)]
For each chart $\psi\colon U_\psi\to V_\psi$ of~$N$ and each chart
$\phi\colon U_\phi\to V_\phi$ of $M$
such that $f(U_\phi)\sub U_\psi$, the following diagram
is commutative:
\[
\begin{array}{ccc}
TU_\phi & \mapright{Tf|_{TU_\phi}^{TU_\psi}}&
TU_\psi\\[1.1mm]
\mapdown{T\phi}&& \mapdown{T\psi}\\
V_\phi\times E &
\mapright{T(\psi\circ f\circ \phi^{-1})} &
V_\psi\times F\,.
\end{array}
\]
\end{description}
The map $Tf$
is called the \emph{tangent map} of~$f$.
Given $p\in M$, we define
$T_pf:=T_p(f):=Tf\big|_{T_pM}^{T_{f(p)}N}\colon T_pM\to T_{f(p)}N$.
\end{defn}
To see that $Tf$ exists,
simply set
\[
Tf(v):=\bigl((T\psi)^{-1}\circ T(\psi\circ f\circ \phi^{-1})
\circ T\phi\bigr)(v)\quad \mbox{  for } \quad v\in TM,
\] 
where $\phi$ and $\psi$ are any charts
as in~(b) such that $\pi_{TM}(v)\in U_\phi$.
Using insertion of charts,
Lemma~\ref{chckdt}(a),
and (\ref{simplchain2}),
we see that $Tf$ is well defined.
By definition, the desired
diagrams are commutative.
Finally, $Tf$ is $C^{r-1}$,
as it can be written locally as
the composition
\[ Tf|_{TU_\phi}=(T\psi)^{-1}\circ T(\psi\circ f\circ \phi^{-1})
\circ T\phi \] 
of $C^{r-1}$-maps.
\begin{prop}[Chain Rule]\label{halffunct}
Let $M,N$ and $K$ be $C^r$-manifolds
and
$f\colon M\to N$ as well as $g\colon N\to K$
be $C^r$-maps. Then
$T(g\circ f)=Tg\circ Tf$.
\end{prop}
\begin{prf}
Given $p\in M$, we find charts $\theta$ of~$K$,
$\psi$ of~$N$ and $\phi$ of~$M$
around~$g(f(p))$,
$f(p)$ and $p$, respectively,
such that $f(U_\phi)\sub U_\psi$ and $g(U_\psi)\sub U_\theta$.
Then
\[
T(\theta \circ (g\circ f)\circ \phi^{-1})
=
T((\theta \circ g\circ \psi^{-1})\circ (\psi\circ f \circ \phi^{-1}))
=
T(\theta \circ g\circ \psi^{-1})\circ
T(\psi\circ f \circ \phi^{-1})
\]
by (\ref{simplchain2}).
Thus
\begin{eqnarray*}
T(g\circ f)|_{TU_\phi} & = &
(T\theta)^{-1}\circ T(\theta\circ g\circ f\circ \phi^{-1})\circ T\phi\\
&=&
(T\theta)^{-1}\circ
T(\theta \circ g\circ \psi^{-1})\circ T\psi\circ (T\psi)^{-1}
\circ
T(\psi\circ f \circ \phi^{-1})
\circ T\phi\\
&= & Tg\circ Tf|_{TU_\phi}\, .
\end{eqnarray*}
The assertion follows since $p$ was arbitrary.
\end{prf}
\begin{rem}\label{tanfunct}
Proposition~\ref{compmfds}
enables us to consider the class
of all \break {$C^r$-manifolds}
as a category, with
$C^r$-maps as the morphisms.
Since also $T(\id_M)=\id_{TM}$,
Proposition~\ref{halffunct}
shows that $T$ is a functor
from the category of $C^r$-manifolds to
the category of $C^{r-1}$-manifolds.
\end{rem}
\begin{lem}\label{tanprod}
Let $M_1$ and $M_2$ be $C^r$-manifolds
and $\pr_j\colon M_1\times M_2\to M_j$ be the canonical
projections.
Then we have:
\begin{description}[(D)]
\item[\rm (a)]
$\theta:=
(T(\pr_1), T(\pr_2))\colon T(M_1\times M_2)\to T(M_1)\times T(M_2)$
is a $C^{r-1}$-diffeomorphism.
\item[\rm (b)]
$\theta(T_x(M_1\times M_2))\sub T_{x_1}M_1\times T_{x_2}M_2$
holds
for each $x=(x_1,x_2)\in M_1\times M_2$,
and the restriction
$\theta\colon T_x(M_1\times M_2)\to
T_{x_1}M_1\times T_{x_2}M_2$
is an isomorphism of topological vector spaces.
\end{description}
\end{lem}
\begin{prf}
The first half of (b) is clear.
For any charts $\phi_1$ of $M_1$ and $\phi_2$ of $M_2$,
using that
$\big(T(\phi_1\times \phi_2)\big)^{-1}=T(
\phi_1^{-1}\times \phi_2^{-1})$,
we see that
\begin{eqnarray*}
\!\!\!\!
\lefteqn{(T(\phi_1)\times T(\phi_2))\circ \theta\circ
\big(T(\phi_1\times \phi_2)\big)^{-1}}\qquad\\
&=&
\big(
T(\phi_1) \circ T\pr_1 \circ
\big(T(\phi_1\times \phi_2)\big)^{-1},\;
T(\phi_2)\circ T\pr_2 \circ
\big(T(\phi_1\times \phi_2)\big)^{-1}\big)\\
&=&
(T(\phi_1\circ \pr_1\circ \, (\phi_1^{-1}\times \phi_2^{-1})),\,
T(\phi_2\circ \pr_2\circ (\phi_1^{-1}\times \phi_2^{-1})))\\
&=& (T\pi_1,T\pi_2)\,,
\end{eqnarray*}
where $\pi_j\colon V_{\phi_1}\times V_{\phi_2}\to V_{\phi_j}$
is the projection.
Here $(T\pi_1,T\pi_2)$ is the map
\[
V_{\phi_1}\times V_{\phi_2}\times
E_1\times E_2\to 
V_{\phi_1}\times E_1\times V_{\phi_2}\times E_2\,,\quad
(x_1,x_2,y_1,y_2)\mto (x_1,y_1,x_2,y_2)
\]
which is a $C^{r-1}$-diffeomorphism. This entails
that $\theta$ is a $C^{r-1}$-diffeomorphism.
Notably, $\theta$ is a homeomorphism and hence restricts
to a homeomorphism
from $T_x(M_1\times M_2)$ onto $T_{x_1}M_1\times T_{x_2}M_2$
for each $x=(x_1,x_2)\in M_1\times M_2$.
Being also linear (as it equals $(T_x\pr_1, T_x\pr_2)$),
the restriction is an isomorphism of topological vector spaces.
\end{prf}
\begin{rem} \label{rem:tangbun-prod} 
(a) We shall frequently identify the tangent bundle\linebreak
$T(M_1\times M_2)$ with
$TM_1\times TM_2$ by means of the diffeomorphism~$\theta$
described in Lemma~\ref{tanprod}.
Iterating Lemma~\ref{tanprod},
we see that
\[ T(M_1\times\cdots \times M_n)\isom TM_1\times\cdots\times
TM_n \quad \mbox{  via } \quad (T\pr_1,\ldots,T\pr_n),\]
for any finite product of $C^r$-manifolds.\medskip

\noindent
(b) Contrary to (a), in the case of open subsets
$U\sub E$ and $V\sub F$ of locally convex spaces $E$ and $F$,
following the convention in Example~\ref{anothtriv} we can also identify $T(U\times V)$ with
$(U\times V)\times (E\times F)$, considering $U\times V$ as an open subset of the locally
convex space $E\times F$. The choice of identification will always be clear from the context.
\end{rem}
The following simple observation is useful:
\begin{lem}\label{hestd}
If $f_1\colon M_1\to N_1$ and $f_2\colon M_2\to N_2$
are $C^1$-maps between $C^1$-manifolds,
then $T(f_1\times f_2)\colon T(M_1\times M_2)\to T(N_1\times N_2)$
corresponds to the map $Tf_1\times Tf_2\colon TM_1\times TM_2\to
TN_1\times TN_2$
under the preceding identifications.
\end{lem}
\begin{prf}
We shall use the projections
$\pr_j\colon M_1\times M_2\to M_j$,
$\pi_j\colon N_1\times N_2\to N_j$
and $\Pi_j\colon TN_1\times TN_2\to TN_j$.
For $j\in \{1,2\}$, we then have
\begin{eqnarray*}
\Pi_j\circ (T\pi_1, T\pi_2)\circ T(f_1\times f_2)
&=& T\pi_j\circ T(f_1\times f_2)
\,=\, T(\pi_j\circ (f_1\times f_2))\\
&=& T(f_j\circ\pr_j)=T(f_j)\circ
T(\pr_j)\\
&=& \Pi_j\circ (Tf_1\times Tf_2)\circ (T\pr_1, T\pr_2)
\end{eqnarray*}
and hence
$(T\pi_1,T\pi_2)\circ T(f_1\times f_2)
=(Tf_1\times Tf_2)\circ (T\pr_1, T\pr_2)$.
\end{prf}
\begin{rem}\label{simrk}
Similarly, $T(f_1,f_2)$ corresponds to $(Tf_1,Tf_2)$
for a map $(f_1,f_2)\colon M\to N_1\times N_2$
(exercise).
\end{rem}
The Rule on Partial
Differentials generalizes
to maps between manifolds in the following form.
\begin{lem}[Rule on Partial Differentials]\label{pardifmfd}
Let $M_1,M_2$ and $N$ be $C^1$-manifolds
and $f\colon M_1\times M_2\to N$ be
a $C^1$-map. Let $\pr_j\colon M_1\times M_2\to M_j$
be the projections for $j\in \{1,2\}$.
Then
\[
T_pf(v)\;=\; T_xf(\cdot,y)(T\pr_1(v))+T_yf(x,\cdot)(T\pr_2(v))
\]
for all $p=(x,y)\in M_1\times M_2$
and $v\in T_p(M_1\times M_2)$.
Thus, identifying $T(M_1\times M_2)$ with
$TM_1\times TM_2$ as above, we have
\[
Tf_p(v_x,v_y)\;=\;T_x(f(\cdot,y))(v_x)+T_y(f(x,\cdot))(v_y)
\]
for all $(x,y)\in M_1\times M_2$, $v_x\in T_xM_1$,
and $v_y\in T_yM_2$.
\end{lem}
\begin{prf}
The desired formula follows from the Rule
on Partial Derivatives (Proposition~\ref{rulepartial}),
applied in local charts.
The details are left to the reader as an exercise.
\end{prf}
\subsection*{Tangent maps of mappings to vector spaces}
\begin{defn}\label{defndfinE}
Let $E$ be a locally convex space,
$U\sub E$ be open,
$M$ be a $C^r$-manifold
and $f\colon M\to U$ be a $C^r$-map.
Identifying $TU$ with $U\times E$
as in Example~\ref{anothtriv},
we can consider the first and second
component of $Tf$.
By Definition~\ref{deftangmap}(a),
we have
$\pi_{TU}\circ Tf=f\circ \pi_{TM}$,
where $\pi_{TU}\colon U\times E\to U$
is the projection onto the first coordinate
(in view of the identification made).
Thus, writing $df\colon TM\to E$
for the second component of~$Tf$,
we have
\begin{equation}\label{compsTf}
Tf \; =\; (f\circ \pi_{TM}, \, df)\,.
\end{equation}
Given $p\in M$, we shall also
write 
\[ df(p):=(df)_p:=df|_{T_p(M)};\] 
this is a continuous linear map $T_p(M)\to E$.
\end{defn}
\begin{rem}
During our exposition of the basic theory,
we shall strictly distinguish between
$Tf$ and $df$.
In the literature,
a clear distinction is not always made.
\end{rem}
\begin{numba}
Let $U$ be an open subset of a locally convex space,
$f\colon U\to M$ a $C^r$-map to a $C^r$-manifold with $r\in \N\cup\{\infty,\omega\}$.
Given $x\in U$, we identify $T_xU$ with $\{x\}\times E$,
whence $T_xf$ is a map $\{x\}\times E\to T_{f(x)}M$.
We often relax notation
and write $T_xf$ also for the map $T_xf(x,\cdot)\colon E\to T_{f(x)}M$.
Using this convention, $(T_xf)^{-1}=d(f^{-1})|_{T_{f(x)}M}$
if $f$ is a $C^r$-diffeomorphism.
\end{numba}
\subsection*{Initial submanifolds}
Using tangent maps as a tool,
we are able to see now that
the modeling space of a submanifold
is unique up to isomorphism.
More generally, this conclusion will hold
for so-called \emph{initial}
submanifolds,
which are very useful generalizations
of submanifolds.
\begin{defn}\label{definitmfd}
Let $M$ be a $C^r$-manifold
and $N$ be a $C^r$-manifold
such that $N\sub M$ as a set.
We call $N$ an \emph{initial submanifold of~$M$}
if the following conditions are satisfied:
\begin{description}[(D)]
\item[(a)]
The inclusion map
$\iota\colon N\to M$
is $C^r$ and $T_p\iota\colon T_pN\to T_pM$
is injective 
for each $p\in N$.
\item[(b)]
A mapping $f\colon P\to N$
from a $C^r$-manifold~$P$
to~$N$ is $C^r$ if and only if
$\iota\circ f\colon
P\to M$ is $C^r$.
\end{description}
\end{defn}
\begin{rem}\label{newrminiti}
The topology on~$N$ need not be the induced topology here;
it may be properly finer.
Using the injective continuous linear map $T_p\iota$, we identify $T_pN$
with the vector subspace $T_p\iota(T_pN)$ of~$T_pM$.
Given a $C^r$-map $f\colon P\to M$
with $f(P)\sub N$,
we then have
$T_xf(T_xP)\sub T_{f(x)}N$
for each $x\in P$
(because $T_xf=T_{f(x)}(\iota)\circ T_x(f|^N)$).
\end{rem}
\begin{lem}\label{basinitmfd}
Let $M$ be a $C^r$-manifold
and $N\sub M$ be a subset.
\begin{description}[(D)]
\item[\rm (a)]
If $\cA_1$ and $\cA_2$ are atlases
making $N$ an initial submanifold of~$M$,
then the identity map
$\id\colon (N,\cA_1)\to (N,\cA_2)$
is a $C^r$-diffeomorphism.
In this sense, the manifold structure
on an initial submanifold
is unique.
Furthermore, its modeling space
is unique up to an isomorphism of
topological vector spaces, if $N\not=\emptyset$.
\item[\rm (b)]
Every submanifold $N\sub M$ is an initial
submanifold.
\item[\rm (c)]
Let $\cA_1$ and $\cA_2$ be atlases making
$N\not=\emptyset$ a submanifold of~$M$, modeled on
vector subspaces $F_1$ and $F_2$
of the modeling space~$E$ of~$M$, respectively.
Then there is an isomorphism
$\alpha\colon E\to E$ of topological
vector spaces such that
$\alpha(F_1)=F_2$.
\end{description}
\end{lem}
\begin{prf}
(a) Let $\iota_j\colon (N,\cA_j)\to M$
be the inclusion map for $j\in\{1,2\}$, and $f\colon (N,\cA_1)\to (N,\cA_2)$
be the identity map.
Then $\iota_2\circ f=\iota_1$ is $C^r$ and hence $f$
is $C^r$, since $(N,\cA_2)$ is an initial
submanifold. A similar argument shows that also
$f^{-1}$ is $C^r$.
Pick $p\in N$.
Then $T_pf\colon T_p(N,\cA_1)\to T_p(N,\cA_2)$
is an isomorphism of topological vector spaces,
where $T_p(N,\cA_j)\isom F_j$.

(b) Pick $p\in N$.
Let
$\phi\colon U_\phi\to V_\phi$ be a chart of~$M$ adapted to~$N$
such that $p\in U_\phi$,
and $\phi_N$ be the corresponding submanifold chart for~$N$,
as in Definition~\ref{defsubm}.
Set $x_0:=\phi(p)$.
Then $\phi\circ \iota\circ \phi_N^{-1}$
is the restriction $\lambda|_{V_\phi\cap F} \colon V_\phi\cap F\to V_\phi$
of
the inclusion map $\lambda\colon F\to E$ of
the modeling space of~$N$ into the modeling space of~$M$,
which is linear and a topological embedding.
Therefore
$T_p\phi\circ T_p\iota\circ (T_p\phi_N)^{-1}
=T_{x_0}(\phi\circ \iota\circ \phi_N^{-1})$
sends $(x_0,y)\in T_{x_0}(V_\phi\cap F)=\{x_0\}\times F$
to $(x_0,y)\in \{x_0\}\times E=T_{x_0}V_\phi$
and hence corresponds to
the embedding $\lambda\colon F\to E$.
Thus $T_p\iota$ is a topological embedding.
All other requirements have been verified in Lemma~\ref{mapsubmfd}.

(c) Given $p\in N$, let $\phi$ be a
chart for~$M$ around~$p$ adapted to the submanifold~$N$ modeled on~$F_1$,
and $\phi_N\in\cA_1$ be the corresponding submanifold chart,
as in Definition~\ref{defsubm}.
Let~$\psi$ be a chart for~$M$ around~$p$ adapted to the submanifold~$N$ modeled on~$F_2$,
and $\psi_N\in\cA_2$ be the corresponding submanifold chart.
After replacing the domains $U_\phi$ and $U_\psi$ with $U_\phi\cap U_\psi$,
we may assume that $U_\phi=U_\psi$.
Then $\psi\circ \phi^{-1}$ is a $C^r$-diffeomorphism,
so $\alpha:=$ $(\psi\circ\phi^{-1})'(\phi(p))$
is an automorphism
of the topological vector space~$E$.
As
\[
\psi\circ \phi^{-1}\big\vert_{V_\phi\cap F_1}^{V_\psi\cap F_2}=
\psi_N\circ \phi_N^{-1},
\]
calculating directional derivatives we see that $\alpha(F_1)\sub F_2$ holds
and its inverse $\alpha^{-1}=(\phi\circ\psi^{-1})'(\psi(p))$
maps $F_2$ inside~$F_1$.
As a consequence, $\alpha(F_1)=F_2$.
%
\end{prf}
\subsection*{Immersions, submersions and embeddings}
It is useful to consider mappings which locally look
like a linear projection $E\times F \to F$ onto a direct summand,
or a linear inclusion map $E\to E\oplus F$.
\begin{defn}\label{defn-imm-subm}
Let $r\in \N_0\cup\{\infty,\omega\}$ and $f\colon M\to N$
be a $C^r$-map between
$C^r$-manifolds modeled on sets $\cE$ and $\cF$
of locally convex spaces,
respectively.
\begin{description}[(D)]
\item[(a)]
$f$ is called a \emph{submersion}
if, for each $x\in M$, there exist a chart \break $\phi\colon U_\phi\to V_\phi\sub E_\phi$
of~$M$ around $x$, a chart $\psi\colon U_\psi\to V_\psi\sub F_\psi$
of~$N$ around $f(x)$
and a continuous linear map $q\colon E_\phi\to F_\psi$
which has a continuous linear right inverse,
such that $f(U_\phi)\sub U_\psi$, $q(V_\phi)\sub V_\psi$,
and $\psi\circ f\circ \phi^{-1}=q|_{V_\phi}$.
\item[(b)]
$f$ is called an \emph{immersion}
if, for each $x\in M$, there exist a chart \break $\phi\colon U_\phi\to V_\phi\sub E_\phi$
of~$M$ around $x$, a chart $\psi\colon U_\psi\to V_\psi\sub F_\psi$
of~$N$ around $f(x)$
and a continuous linear map $i\colon E_\phi\to F_\psi$
which has a continuous linear left inverse,
such that $f(U_\phi)\sub U_\psi$, $i(V_\phi)\sub V_\psi$,
and $\psi\circ f\circ \phi^{-1}=i|_{V_\phi}$.
\item[(c)]
$f$ is called an \emph{embedding of $C^r$-manifolds}
if $f$ is an immersion and a topological embedding.
\end{description}
\end{defn}
\begin{rem}
We here use tangent spaces and tangent maps
also for {manifolds} modeled on sets of locally convex spaces, as in Definition~\ref{non-pure}.
The definition of geometric tangent spaces carries over directly:
simply use charts $\phi\colon U_\phi\to V_\phi\sub E_\phi$ with $E_\phi$ in the set $\cE$
of modeling spaces.
In the definition of abstract tangent spaces, use quadrupels $(E,\phi,x,v)$ with $E\in\cE$
and $\phi\colon U_\phi\to V_\phi\sub E$; give $T_pM$ the locally convex topological vector space structure making $k_\phi\colon E\to T_pM$, $v\mto [E,\phi,x , v]$ an isomorphism of topological vector spaces,
for each chart $\phi\colon U_\phi\to V_\phi\sub E$ of $M$ around~$p$, with $x:=\phi(p)$.
\end{rem}
Note that every submersion is an open map.
\begin{prop}\label{subm-sec}
Let $r\in\N_0\cup\{\infty,\omega\}$
and $f\colon M\to N$ be a surjective $C^r$-submersion between
$C^r$-manifolds modeled on sets $\cE$ and $\cF$
of locally convex spaces,
respectively. Then the following holds:
\begin{description}[(D)]
\item[\rm(a)]
For each $x\in M$, there exists an open $f(x)$-neighborhood
$W\sub N$ and a $C^r$-map $\sigma\colon W\to M$
such that $f\circ \sigma=\id_W$.
\item[\rm(b)]
For each $C^r$-manifold $L$, a map $g\colon N\to L$
is $C^r$ if and only if $g\circ f$ is~$C^r$.
\item[\rm(c)]
For every $($not necessarily pure$)$ submanifold $S\sub N$,
the pre-image $f^{-1}(S)$ is a
$($not necessarily pure$)$
submanifold of~$M$.
The restriction $f|_{f^{-1}(S)}\colon f^{-1}(S)\to S$
is a submersion.
\end{description}
\end{prop}
\begin{prf}
(a) Let $\phi$, $\psi$, and $q$ be as in Definition~\ref{defn-imm-subm}(a).
After replacing $\phi$ with $\phi-\phi(x)$ and
$\psi$ with $\psi-\psi(f(x))=\psi-q(\phi(x))$,
we may assume that $\phi(x)=0$ and $\psi(f(x))=0$.
Let $\alpha\colon F\to E$ be a continuous linear
map with $q\circ \alpha=\id_F$.
Then $W:=\psi^{-1}(\alpha^{-1}(V_\phi)\cap V_\psi)$
is an open $f(x)$-neighborhood in~$N$
and $\sigma:=\phi^{-1}\circ \alpha\circ \psi|_W$
is a $C^r$-map such that $\sigma(f(x))=x$
and $\psi \circ f\circ \sigma=q\circ\phi\circ \sigma=\psi|_W$,
entailing that $f\circ \sigma=\id_W$.

(b) If $g$ is $C^r$, then $g\circ f$ is $C^r$.
If, conversely,
$g\circ f$ is $C^r$, let $y\in N$.
Since $f$ is surjective, there is $x\in M$
such that $f(x)=y$. Let $\sigma$ be as in~(a).
Then $g|_W=g|_W\circ \id_W=g \circ f\circ \sigma$
is $C^r$. Thus $g$ is $C^r$ on an open neighborhood~$W$
of each point~$y$, and thus $g$ is~$C^r$.

(c)  Given $x\in f^{-1}(S)$,
let $\phi$, $\psi$, and $q$ be as in Definition~\ref{defn-imm-subm}(a);
we may assume that $\phi(x)=0$ and $\psi(f(x))=0$.
Let $\alpha\colon F\to E$ be a continuous linear
map such that $q\circ\alpha=\id_F$.
Then
\[
\beta\colon E\to F\times \ker(q),\;\, y\mto (q(y),y-\alpha(q(y)))
\]
is an isomorphism of locally convex space and
$\pr_1\circ \beta=q$ using the projection $\pr_1\colon F\times\ker(q)\to F$,
$(v,w)\mto v$.
Thus $\beta(V_\phi)$
is an open $0$-neighborhood in $F\times \ker(q)$;
after shrinking $V_\phi$,
we may assume that $\beta(V_\phi)=P\times Q$
for open $0$-neighborhoods $P\sub F$ and $Q\sub \ker(q)$.
Then $P\sub V_\phi$; after replacing $V_\psi$ with $P$,
we have the $C^r$-diffeomoprhisms
$\theta:=\beta\circ \phi\colon U_\phi\to P\times Q$ and $\psi\colon U_\psi\to P$
such that $\pr_1\circ \beta=\psi$.
Then $\psi(S\cap U_\psi)$ is a submanifold
of $P$ and $A:=\pr_1^{-1}(\psi(S\cap U_\phi))=\psi(S\cap U_\phi)\times Q$
is a submanifold of $P\times Q$ such that $\pr_1|_A$ is a submersion.
Thus $\theta^{-1}(A)$ is a submanifold of $U_\phi$ and hence of~$M$
such that $f|_{\theta^{-1}(A)}=\psi^{-1}\circ \pr_1|_A\circ \theta|_{\theta^{-1}(A)}$
is a submersion.
For $y\in U_\phi$,
we have $f(y)\in S$ if and only if $\pr_1(\theta(y))=q(\phi(y))=\psi(f(y))\in \psi(S\cap U_\psi)$.
Thus $U_\phi\cap f^{-1}(S)=\theta^{-1}(A)$ is a submanifold.
Being a submanifold locally, $f^{-1}$ is a submanifold.
\end{prf}
\begin{cor}\label{unique-quot-mfd}
Let $r\in \N_0\cup\{\infty,\omega\}$,
$M$ be a $C^r$-manifold modeled on a set~$\cE$
of locally convex spaces and
$q\colon M\to N$ be a surjective map.
For $j\in\{1,2\}$ let us write $N_j$ for~$N$,
endowed with a $C^r$-manifold structure modeled on a set $\cE_j$
of locally convex spaces such that $q\colon M\to N_j$
is a $C^r$-submersion.
Then $N_1\to N_2$, $x\mto x$
is a $C^r$-diffeomorphism.
If $\cE_1=\cE_2$, then $N_1=N_2$.
\end{cor}
\begin{prf}
Write $q_j$ for $q$ as a map to $N_j$.
Let $f\colon N_1\to N_2$ be the map $x\mto x$.
Since $f\circ q_1=q_2$ is $C^r$, the map $f$ is $C^r$ by Proposition~\ref{subm-sec}(b).
Likewise, $f^{-1}\circ q_2=q_1$ implies that $f^{-1}$ is $C^r$.
\end{prf}
\begin{prop}\label{naive-imm-subm}
Let $r\in \N\cup\{\infty,\omega\}$ and $f\colon M\to N$
be a $C^r$-map between
$C^r$-manifolds $M$ and $N$ modeled on sets $\cE$ and $\cF$
of locally convex spaces, respectively.
\begin{description}[(D)]
\item[\rm(a)]
If each $F\in \cF$ has finite dimension, then $f$ is a submersion
if and only if $T_xf\colon T_xM\to T_{f(x)}N$
is surjective for all $x\in M$.
\item[\rm(b)]
If each $E\in \cE$ has finite dimension,
then $f$ is an immersion
if and only if $T_xf\colon T_xM\to T_{f(x)}N$
is injective for all $x\in M$.
\end{description}
\end{prop}
\begin{prf}
(a) If $f$ is a submersion
and $x\in M$, let $\phi$, $\psi$, and $q$
be as in Definition~\ref{defn-imm-subm}(a).
Then $f=\psi^{-1}\circ q|_{V_\phi}\circ \phi$,
whence $T_xf=T_{\psi(f(x))}\psi^{-1}\circ q\circ d\phi|_{T_xM}$
is surjective.
Conversely,
assume that each tangent map is surjective.
Given $x\in M$, there exists a chart
$\phi\colon U_\phi\to V_\phi\sub E_\phi$ for $M$ around~$x$
and a chart $\psi\colon U_\psi\to V_\psi\sub F_\psi$
for~$N$ around $f(x)$, such that $f(U_\phi)\sub U_\psi$.
By hypothesis,
\[
q:=(\psi\circ f\circ \phi^{-1})'(\phi(x))\colon E_\phi \to F_\psi
\]
is surjective. Let $b_1,\ldots, b_n$ be a basis for $F_\psi$;
we find $a_1,\ldots, a_n\in E_\phi$ such that $q(a_j)=b_j$
for $j\in\{1,\ldots, n\}$. Let $\alpha\colon F_\psi\to E_\phi$
be the continuous linear map such that $\alpha(b_j)=a_j$
for all $j\in\{1,\ldots, n\}$.
Then
\[
\beta\colon E_\phi\to \ker(q)\times F_\psi,\;\,
v\mto (v-\alpha(q(v)),q(v))
\]
is an isomorphism of topological vector spaces
with $\beta^{-1}(a,b)=a+\alpha(b)$.
After shrinking $V_\phi$, we may assume that
$\beta(V_\phi)=A\times B$ with open
$0$-neighborhoods $A\sub \ker(q)$ and $B\sub F_\psi$.
Then the map $h\colon A\times B\to \ker(q)\times F_\psi$,
\[
(a,b)\mto
(a,(\psi\circ f\circ \phi^{-1}\circ \beta^{-1})(a,b))
=(a,(\psi\circ f\circ \phi^{-1})(a+\alpha(b)))
\]
is $C^r$ and $d_2(\psi\circ f\circ \phi^{-1}\circ\beta^{-1})(0,0,c)
=q(\alpha(c))=c$, whence
\[
d_2(\psi\circ f\circ\phi^{-1}\circ\beta^{-1})(0,0,\cdot)=\id_{F_\psi}.
\]
By the Inverse Function Theorem with Parameters,
after shrinking $A$ and $B$, the image $h(A\times B)$
is open and $h$ is a $C^r$-diffeomorphism
onto $h(A\times B)$
(see Theorem~\ref{thm-inv-para} as well as (a) and (c) in Remark~\ref{variants-i-para} 
for the finite-dimensional and the analytic case).
Then $W:=\beta^{-1}(h(A\times B))$ is open in~$E_\phi$
and $\theta:=\beta^{-1}\circ h\circ\beta\circ \phi\colon U_\phi\to W$
is a chart for~$M$.
Using the projection $\pr_2\colon \ker(q)\times F_\psi\to F_\psi$, $(v,w)\mto w$,
we have $\psi\circ f\circ \phi^{-1}=\pr_2\circ h\circ \beta|_{V_\phi}$,
whence
$\psi\circ f\circ \theta^{-1}
=\psi\circ f\circ \phi^{-1}\circ \beta^{-1}\circ h^{-1}\circ \beta|_W
=\pr_2\circ h \circ h^{-1}\circ \beta=q|_W$.

(b) If $f$ is an immersion
and $x\in M$, let $\phi$, $\psi$, and $i$
be as in Definition~\ref{defn-imm-subm}(b).
Then $f=\psi^{-1}\circ i|_{V_\phi}\circ \phi$,
whence $T_xf=T_{\psi(f(x))}\psi^{-1}\circ i\circ d\phi|_{T_xM}$
is injective.
Assume that, conversely, each tangent map is injective.
Let $x_0\in M$.
We find a chart
$\phi\colon U_\phi\to V_\phi\sub E_\phi$ of~$M$
around~$x_0$
and a chart $\psi\colon U_\psi\to V_\psi\sub F_\psi$
of~$N$ around $f(x_0)$ such that $f(U_\phi)\sub U_\psi$,
$\phi(x_0)=0$, and $\psi(f(x_0))=0$.
Then
\[
\alpha:=(\psi\circ f\circ \phi^{-1})'(\phi(x))\colon E_\phi\to F_\psi
\]
is injective. As $\alpha(E_\phi)$ has finite dimension,
it is complemented in the topological vector space $F_\psi$;
let $Y\sub F_\psi$ be a topological complement.
Let $\pr_1\colon F_\psi\to Y$ and $\pr_2\colon F_\psi\to\alpha(E_\phi)$
be the continuous linear mappings such that $y=\pr_1(y)+\pr_2(y)$
Then $g:=\pr_2\circ \psi\circ f\circ \phi^{-1}\colon U_\phi\to \alpha(F_\phi)$
has derivative $g'(0)=\pr_2\circ \alpha=\alpha|^{\alpha(E_\phi)}$
at $0$, which is an isomorphism between finite-dimensional
vector spaces. By the Inverse Function Theorem,
after shrinking $V_\phi$ we may assume that $g(V_\phi)$ is open in $\alpha(E_\phi)$
and $g\colon V_\phi\to g(V_\phi)$ is a $C^r$-diffeomorphism.
Then $Q:=g'(0)^{-1}(g(V_\phi))$ is an open $0$-neighborhood in~$E_\phi$
and $\kappa:=g'(0)^{-1}\circ g\circ \phi\colon U_\phi\to Q$
is a chart for $M$ around~$x_0$.
Moreover, $h:=\psi\circ f\circ \phi^{-1}\circ g^{-1}\colon g(V_\phi)\to F_\psi$
is a $C^r$-map with $\pr_2\circ \, h=g\circ g^{-1}=\id_{V_\phi}$.
Now $Y+V_\phi$
and $Y+g(V_\phi)$ are open $0$-neighborhoods in $F_\psi$,
and
\[
\Xi\colon Y+g(V_\phi)\to Y+g(V_\phi),\;\,
v+w\mto v+ (\pr_1\circ \, h)(w)+w=v+h(w)
\]
for $v\in Y$, $w\in g(V_\phi)$
is a $C^r$-diffeomorphism with 
inverse $a+b\mto a + 2b- h(b)$ for $a\in Y$, $b\in g(V_\phi)$.
Then $P:=\Xi^{-1}(V_\psi)$ is an open subset of $E_\phi$
and $\theta:=\Xi^{-1}\circ \psi\colon U_\psi\to P$ is a chart of~$N$
such that $(\theta\circ f\circ \phi^{-1})(x)=g(x)$ for all
$x\in V_\phi$ and hence $(\theta\circ f\circ \kappa^{-1})(x)=g'(0)(x)
=\alpha(x)$ for all $x\in Q$.
Since $\alpha$ has a continuous linear left inverse,
$f$ is an immersion.
\end{prf}
\begin{prop}\label{embedding-submfd}
Let $r\in \N_0\cup\{\infty,\omega\}$ and $f\colon M\to N$
be a $C^r$-map between
$C^r$-manifolds $M$ and $N$, which need not be pure.
The following conditions are equivalent.
\begin{description}[(D)]
\item[\rm(a)]
$f$ is an embedding of $C^r$-manifolds.
\item[\rm(b)]
$f(M)$ is a $($not necessarily pure$)$
submanifold of~$N$
and $f\colon M\to f(M)$ is a $C^r$-diffeomorphism.
\end{description}
\end{prop}
\begin{prf}
If (b) holds, we find for $x_0\in N$
a chart $\psi\colon U_\psi\to V_\psi\sub E_\psi$
which is adapted to~$f(M)$ in the sense
that $\psi(f(M)\cap U_\psi)=F\cap V_\psi$
for some closed vector subspace $F\sub E_\psi$.
Let $\psi_{f(M)}:=\psi|_{f(M)\cap U_\psi}\colon f(M)\cap U_\psi\to F\cap V_\psi\sub F$
be the corresponding submanifold chart.
Then $W:=f^{-1}(f(M)\cap U_\psi)$ is an open subset of~$M$
and $\phi:=\psi\circ f|_W\colon W\to F\cap V_\psi$
is a $C^r$-diffeomorphism and hence an
$F$-chart for~$M$.
If $i\colon F\to E_\psi$
is the inclusion map, then $\psi\circ f\circ \psi^{-1}=i|_{F\cap V_\psi}$,
entailing that $f$ is an immersion. As the inclusion map $j\colon f(M)\to M$
is a topological embedding and $f|^{f(M)}$ a homeomorphism,
we deduce that $f=j\circ f|^{f(M)}$ is a topological embedding.
Hence $f$ is an embedding of $C^r$-manifolds.

If (a) holds and $y\in f(M)$,
we let $x\in M$ with $y=f(x)$.
Since $f$ is an immersion, there exist a chart
$\phi\colon U_\phi\to V_\phi\sub E_\phi$
of~$M$ around~$x$ and a chart $\psi\colon U_\psi\to V_\psi\sub F_\phi$
of~$N$ around $f(x)$ such that $f(U_\phi)\sub U_\psi$
and $\psi\circ f\circ \phi^{-1}\colon V_\phi\to  V_\psi$ equals
$i|_{V_\phi}$ for some continuous linear map $i\colon E_\phi\to F_\psi$
admitting a continuous linear left inverse $q\colon F_\psi\to E_\phi$.
Then \break $F_\psi=i(E_\phi)\oplus \ker(q)$ as a locally convex space.
Now
\[
i(V_\phi)=(\psi\circ f\circ \phi^{-1})(V_\phi)\sub V_\psi.
\]
After replacing $V_\psi$ with $V_\psi\cap(i(V_\phi)\oplus \ker(q))$
and adapting $U_\psi$,
we may assume that
\[
i(E_\phi)\cap V_\psi=i(V_\phi).
\]
Since $f$ is a topological embedding,
$f(U_\phi)$ is relatively open in $f(M)$.
Thus $f(U_\phi)=f(M)\cap Q$ for some open subset $Q\sub N$.
After replacing $Q$ with $Q\cap U_\psi$,
we may assume that $Q\sub U_\psi$.
After replacing $U_\psi$ with~$Q$
and $V_\psi$ with $\psi(Q)$,
we may assume that $U_\psi=Q$, and thus
\[
f(M)\cap U_\psi=f(U_\phi).
\]
Then
\begin{eqnarray*}
\psi(f(M)\cap U_\psi)& =& \psi(f(U_\phi))=(\psi\circ f\circ \phi^{-1})(V_\phi)
=i(V_\phi)= i(E_\phi)\cap V_\psi,
\end{eqnarray*}
showing that $\psi$ is a chart adapted to $f(M)$.
Thus $f(M)$ a submanifold of~$N$.
Since $f$ is an embedding of $C^r$-manifolds,
$f|^{f(M)}\colon M\to f(M)$ is a bijection.
For $\phi$ and $\psi$ as before, $\psi$ is adapted to $f(M)$.
Let $\psi_{f(M)}:= \psi|_{f(M)\cap U_\psi}\colon f(M)\cap U_\psi\to i(E_\phi)\cap
V_\psi$ be the corresponding submanifold chart.
Then
\[
\psi_{f(M)}\circ f|^{f(M)}\circ \phi^{-1}=\psi\circ f\circ \phi^{-1})|^{i(E_\phi)\cap V_\phi}
=i_{V_\phi}^{i(V_\phi)}
\]
is a $C^r$-diffeomorphism. Being a bijection and a local
$C^r$-diffeomorphism, $f|^{f(M)}$ is a $C^r$-diffeomorphism.
\end{prf}
\begin{small}
\subsection*{Exercises for Section~\ref{sec:2.3}} 
\label{Exer3.3}
\begin{exer}
  \label{exer:2.3.1} 
Let $r\geq 1$
and $M$ be a $C^r$-manifold
modeled on a locally convex space~$E$.
Find an isomorphism
$\theta_p\colon \cT_p(M)\to T_pM$
for $p\in M$.
If $f\colon M\to N$ is a $C^r$-map,
show that the linear map $\cT_p(f)\colon
\cT_p(M)\to \cT_{f(p)}(N)$
corresponding to $T_p(f)$
(i.e., $\cT_p(f):=\theta_{f(p)}^{-1}\circ T_p(f)
\circ \theta_p$)
is the map sending a geometric tangent
vector (equivalence class of curves)
$[\gamma]$ to $[f\circ \gamma]$.
\end{exer}

\begin{exer}
  \label{exer:2.3.3}
Let $r\in \N\cup\{\infty,\omega\}$
and
$M$ be a $C^r$-manifold modeled on a locally convex space~$E$.
\begin{description}[(D)]
\item[(a)]
Show that $T_pM$ is a split submanifold
of~$TM$ modeled on $\{0\}\times E$ ($\isom E$),
for each $p\in M$.
\item[(b)]
Show that the ``zero-section''
$\sigma\colon M\to TM$, $\sigma(p):=0_p\in T_p(M)$
is an injective $C^{r-1}$-map.
\item[(c)]
Show that $\sigma(M)$ is a split
submanifold of~$TM$ modeled on
$E\times \{0\}$ ($\isom E$).
\item[(d)]
Show that $\sigma|^{\sigma(M)}\colon M\to \sigma(M)$
is a $C^{r-1}$-diffeomorphism.
\end{description}
\end{exer}

\begin{exer}
  \label{exer:2.3.4} Verify the Rule on Partial Differentials
for mappings on products of manifolds
(Lemma~\ref{pardifmfd}), using local charts.
Also fill in the details in Remark~\ref{simrk}.
\end{exer}

\begin{exer}
  \label{exer:2.3.5} Let $E$ be a locally convex space,
$M\sub E$ be a $C^r$-submanifold
modeled on $F\sub E$
(where $r\geq 1$) and $p\in M$.
Let $\iota\colon M\to E$
be the inclusion map
and $V\sub E$ be the set of all
velocity vectors $\gamma'(0)$
of $C^1$-curves~$\gamma$ in~$M$ which are passing through~$p$.
\begin{description}[(D)]
\item[(a)]
Show that $V$ is the image
of $T_p(\iota)$,
and that $T_p(\iota)\colon T_pM\to V$
is an isomorphism of topological
vector spaces.
In particular,
$V$ is a vector subspace
of~$E$ isomorphic to~$F$.
\item[(b)]
If $E=\R^n$ and $M=\bS^{n-1}$,
show that $V=p^\perp=\{x\in \R^n\colon\langle x,p\rangle=0\}$.
For $n\in \{2,3\}$,
make a sketch showing $\bS^{n-1}$, $p$ and $p+V$.
\end{description}
\end{exer}

\begin{exer}
  \label{exer:2.3.6a} (Product Rule) Let
$\beta\colon E_1\times E_2\to F$
be a continuous bilinear map
between locally convex spaces,
$M$ be a $C^1$-manifold
and $f_1\colon M\to E_1$
as well as $f_2\colon M\to E_2$
be $C^1$-maps.
Let $\phi:=\beta\circ (f_1,f_2)\colon M\to F$.
Show that 
\[ d\phi=\beta\circ (df_1,f_2\circ \pi_{TM})
+\beta\circ (f_1\circ\pi_{TM},df_2),\] 
where $\pi_{TM}\colon TM\to M$
is the bundle projection.
\end{exer}

\begin{exer}\label{exc-unique-model-component}
Let $r\in \N\cup\{\infty,\omega\}$ and $M$
be a $C^r$-manifold modeled on a set $\cE$
of locally convex spaces. For $x,y\in M$, write $x\sim y$
if $T_xM\cong T_yM$ as a topological vector space.
Verify that $\sim$ is an equivalence relation on~$M$
and show that the equivalence classes~$[x]$ are open.
Deduce that, for each connected component $C$ of~$M$, we have
$T_xM\cong T_yM$ for all $x,y\in C$.
\end{exer}
  
\end{small}
\section{The Lie algebra of vector fields on a manifold}
\label{sec:2.4} 
In this section,
we define vector fields
and the Lie bracket of vector fields
on a smooth manifold~$M$.
Vector fields act as differential operators
on smooth functions,
and they can be used to define
differential equations
on a manifold.
The Lie bracket of vector fields
turns the space $\cV(M)$
of all vector fields into a Lie algebra.
Later on, the Lie bracket of vector fields will
enable us to turn
the tangent space $T_{\be}(G)$ of
a Lie group~$G$ into a Lie algebra.\\[2.5mm]
To increase the readability,
we shall focus
on smooth vector fields
on smooth manifolds over the ground field~$\R$.
Variants and generalizations
are described in Remark~\ref{othervectorf}.
\begin{defn}\label{defVF}
Let $M$ be a smooth manifold.
A (smooth) \emph{vector field on~$M$}
is a smooth map $X \colon M\to TM$
assigning to each $p\in M$ a tangent vector
$X(p)\in T_pM$.
Thus a vector field is a smooth map
$X\colon M\to TM$ such that
$\pi_{TM}\circ X=\id_M$
(that is, a smooth section of the bundle projection
$\pi_{TM}\colon TM\to M$).
We let $\cV(M)$ denote the set of all
vector fields on~$M$.\\[2.5mm]
Given $X,Y\in \cV(M)$ and
$r,s\in \R$, we define
a section $rX+sY\colon M\to TM$ of~$\pi_{TM}$
via
$(rX+ sY)(p):=rX(p)+sY(p)\in T_p(M)$ for $p\in M$.
Then $rX+sY\in \cV(M)$, since it is
easy to verify using charts
that $rX+sY$ is smooth.
Clearly $\cV(M)$
becomes an $\R$-vector space in this way.
We can also multiply vector
fields by functions:
If $f\in C^\infty(M,\R)$ and
$X\in \cV(M)$,
we define $fX\in \cV(M)$
via pointwise scalar multiplication,
$(fX)(p):=f(p)X(p)$.
In this way, $\cV(M)$ becomes
a $C^\infty(M,\R)$-module,
if the space $C^\infty(M,\R)$ of real-valued
smooth functions on~$M$
is considered as an associative
$\R$-algebra with pointwise
multiplication of functions
as the algebra multiplication.
\end{defn}
\begin{ex}\label{exXE}
Let us consider the simplest case where
our manifold is an open subset $U\sub E$ of a topological
vector space~$E$.
In this case, we have $TU=U\times E$,
and $\pi_{TU}\colon  U\times E\to U$
is the projection onto the first
component. A vector field is a smooth map of the form
\[
X \, =\, (\id_U,X_E)\colon U\to U\times E, \quad  x\mto (x,X_E(x)), \]
where $X_E\colon U\to E$ is a smooth map.
\end{ex}
\begin{numba}\label{numnewappr}
If $M$ is a smooth manifold
modeled on~$E$
and $X$ a vector field on~$M$,
then we still have an analog of $X_E$
locally, for each given chart $\phi\colon U_\phi\to V_\phi$
of~$M$. In fact, $T(\phi)\circ X\circ \phi^{-1}\colon V_\phi\to
V_\phi\times E$
is a vector field on~$V_\phi$ and hence of the form
\begin{equation}\label{newappr}
T(\phi)\circ X\circ \phi^{-1}\,=\,(\id_{V_\phi},X_\phi)
\end{equation}
with
the smooth map
$X_\phi:=d\phi\circ X\circ \phi^{-1}\colon V_\phi\to E$.
We call $X_\phi \in C^\infty(V_\phi, E)$ the \emph{local representative} of~$X$
with respect to the chart~$\phi$.
\end{numba}
\begin{rem}\label{passloclin}
The linearity of the maps $d\phi|_{T_pM}$ implies
that the map
\[
\cV(M)\to C^\infty(V_\phi,E)\,, \quad
X \mto X_\phi
\]
is $\R$-linear, for each chart $\phi\colon  U_\phi\to V_\phi\sub E$.
\end{rem}
By definition of $X_\phi$, we have the commutative
diagram
\begin{equation}\label{prereltd}
\begin{array}{rcl}
TU_\phi &\mapright{T\phi}\;\;\, & V_\phi\times E\\
\mapup{X|_{U_\phi}} \, \ & &\;\;\;  \mapup{(\id,X_\phi)}\\
U_\phi &\mapright{\phi} & \; V_\phi\,.
\end{array}
\end{equation}
It is useful to introduce terminology
describing such situations.
\begin{defn}\label{defrelvefi}
Let $f\colon  M\to N$ be a smooth map
between smooth manifolds.
Two vector fields
$X\in \cV(M)$ and $Y\in \cV(N)$ are called \emph{$f$-related}
if
\[
Y\circ f\; =\; T(f)\circ X\, .
\]
\end{defn}
Thus, we require commutativity
of the diagram
\[
\begin{array}{rcl}
TM & \mapright{Tf} & TN\\
\mapup{X} \, & & \mapup{Y}\\
M & \mapright{f} & N\,.
\end{array}
\]
\begin{rems}\label{rmsrel}
\begin{description}[(D)]
\item[(a)]
If $f\colon M\to N$ is a smooth
map and $X\in \cV(M)$, then there
need not exist a vector field
$Y\in \cV(N)$ which is $f$-related
to~$X$ (nor conversely).
\item[(b)]
If $f\colon M\to N$ and $g\colon N\to K$
are smooth mappings, $X\in \cV(M)$ and
$Y\in \cV(N)$ are $f$-related
and $Y$ and $Z\in \cV(K)$ are
$g$-related, then $X$ and~$Z$ are
$g\circ f$-related,
since $Z \circ g\circ f
= Tg\circ Y\circ f=Tg\circ Tf\circ X=
T(g\circ f) \circ X$.
\item[(c)]
By~(\ref{prereltd}),
the vector fields $X|_{U_\phi}$ and $(\id,X_\phi)$
are $\phi$-related.
\item[(d)]
If $f\colon M\to N$ is a smooth map, $U\sub N$ is open
and $X\in \cV(M)$, $Y\in \cV(N)$
are $f$-related, then $X|_{f^{-1}(U)}$ and $Y|_U$
are $f|_{f^{-1}(U)}^U$-related.
\item[(e)]
Combining (b)--(d),
we see: If $\phi\colon U_\phi\to V_\phi$
and $\psi\colon U_\psi\to V_\psi$
are charts for a smooth manifold~$M$
and $X\in \cV(M)$,
then $(\id,X_\phi)$
and $(\id, X_\psi)$ (restricted to
$\phi(U_\phi\cap U_\psi)$, resp.,
$\psi(U_\phi\cap U_\psi)$)
are related via the transition map
$\tau:=\psi \circ \phi^{-1}$.
The latter means that
\begin{equation}\label{relinloc}
X_\psi(\tau(x)) \; =\; d\tau (x,X_\phi(x))\quad
\mbox{for all $x\in \phi(U_\phi\cap U_\psi)$.}
\end{equation}
\end{description}
\end{rems}
Conversely, compatible
families of smooth maps combine
to a vector field:
\begin{lem}\label{ladealwvf1}
Let $(M,\cA)$ be a smooth manifold modeled on a
locally convex space~$E$,
and let $(\xi_\phi)_{\phi\in \cA}$
be a family of smooth maps $\xi_\phi\colon 
V_\phi\to E$ which are compatible
in the sense that
\[
\xi_\psi\circ (\psi\circ \phi^{-1})\; =\;
d(\psi\circ \phi^{-1})\circ (\id,\,\xi_\phi)
\quad \mbox{ on }\quad \phi(U_\phi\cap U_\psi)
\]
{\rm(}i.e., $\xi_\phi$ and $\xi_\psi$ represent $\psi \circ \phi^{-1}$-related vector 
fields{\rm)} for all charts $\phi,\psi\in \cA$.
Then there is a uniquely determined
vector field $X\in \cV(M)$
with the given family of local representatives,
i.e., $X_\phi=\xi_\phi$
for each $\phi \in \cA$.
\end{lem}
\begin{prf}
\!\!Uniqueness:
$d\phi|_{T_pM}$ being injective,
the condition $d\phi(X(p))\!=\!\xi_\phi(\phi(p))$ determines $X(p)$
for $p\in U_\phi$. As to existence,
we define $X \colon  M\to TM$ via
\[ X(p) := T(\phi^{-1})\big(\phi(p),\xi_\phi(\phi(p))\big) 
\quad \mbox{ for } \quad p\in U_\phi.\] 
By compatibility of the maps $\xi_\phi$,
this yields a well-defined map $X\colon M\to TM$.
By construction, $X$ is a section for~$\pi_{TM}$.
The defining formula for
$X|_{TU_\phi}$ shows that $X$ is smooth
on $TU_\phi$ for each chart $\phi\colon U_\phi\to V_\phi$
and hence smooth on all of~$M$. By construction, $X_\phi=\xi_\phi$.
\end{prf}
\begin{rem}\label{vecfifin}
We mention that vector fields on an open subset
$U\sub \R^n=:E$
can be described more explicitly.
In fact, let $X\in \cV(U)$ be given
and write $X_E=(f_1,\ldots, f_n)$
with smooth functions $f_1,\ldots, f_n\colon
U\to \R^n$, using the notation from
Example~\ref{exXE}.
Let $e_j\in \R^n$ be the $j$th standard
unit vector and
write $\frac{\partial}{\partial x_j}$
for the smooth vector field
on~$U$ such that $\big(\frac{\partial}{\partial x_j}\big)_E(x)=e_j$
for all $x\in U$ (again using the
notation from Example~\ref{exXE}).
Then $X$ is given by
\[
X\;=\; f_1\, \frac{\partial}{\partial x_1}+ \cdots +
f_n\, \frac{\partial}{\partial x_n}\,,
\]
and the smooth functions
$f_1,\ldots, f_n$ are uniquely
determined by this identity.
\end{rem}
The preceding notation already suggests that
vector fields act on functions
as differential operators.
This is indeed the case
(also in infinite dimensions).
\begin{defn}\label{defvfact}
Let $f\colon  M\to F$ be a smooth map
from a smooth manifold~$M$ to
a locally convex space~$F$,
and $X\in \cV(M)$ be a vector field.
We define a smooth map $X .f\colon  M\to F$ as the composition
\[ X .f\, :=\, df\circ X \in C^\infty(M,F). \]
If $f\colon U\to F$ is a smooth map on an open subset
$U\sub M$, we write $X. f:=(X|_U).f=df\circ (X|_U)$.
\end{defn}
\begin{rem}\label{rembil}
Note that
$\cV(M)\times C^\infty(M,F)\to C^\infty(M,F)$,
$(X,f)\mto X.f$ is an $\R$-bilinear map,
since $df$ depends linearly on~$f$ and $df|_{T_pM}\colon  T_pM\to F$
is linear for each $p\in M$.
\end{rem}
\begin{rem}\label{vfonope}
Let us consider the special case where
$M=U$ is an open subset of~$E$. 
In this case, we have $TU=U\times E$,
the vector field $X$ is of the form
$X=(\id_U,X_E)$ for some smooth map
$X_E\colon  U\to E$, and $df\colon TU=U\times E\to F$ is the usual
differential. Thus
\begin{equation}\label{actopen}
(X.f)(x)\, =\, df(X(x))=df(x,X_E(x))
\end{equation}
for $x\in U$, i.e.,
$(X.f)(x)$ is the directional derivative
of $f$ at~$x$ in the direction $X_E(x)$
provided by the vector field~$X$ at~$x$.
In the special case described here,
we abbreviate
\begin{equation}\label{abbrnot}
X_E.f\, :=\, (\id_U,X_E).f\, =\, X .f\,.
\end{equation}
\end{rem}
In the following, we use
Lie algebras and
various related notions (ideals,
Lie algebra homomorphisms,
derivations, and modules).
Readers unfamiliar with these
concepts are referred to any introductory text on Lie algebras, 
such as \cite{HiNe12}. 
\begin{rem}\label{introdDX}
Consider $C^\infty(M,\R)$ as a unital, associative
$\R$-algebra.
Then the map
\begin{equation}\label{firstl}
\cL_X\colon C^\infty(M,\R)\to C^\infty(M,\R)\,,\qquad
\cL_X(f)\,:=\,X.f
\end{equation}
is a derivation of~$C^\infty(M,\R)$,
for each vector field $X\in \cV(M)$.
Indeed, $\cL_X$ is $\R$-linear by Remark~\ref{rembil}.
To verify the Product Rule,
let $f,g\in C^\infty(M,\R)$.
Since
$d(f\cdot g)= (df)\cdot (g\circ\pi_{TM})
+(f\circ \pi_{TM})\cdot dg$
by Exercise~\ref{exer:2.3.6a}, 
we have 
\[ \cL_X(f\cdot g)=d(f\cdot g)\circ X
=(df\circ X)\cdot g+f\cdot (dg\circ X)=
\cL_X(f)\cdot g+f\cdot \cL_X(g).\]
\end{rem}
\begin{thm}\label{Liebrvf}
For each smooth manifold $(M,\cA)$,
the following holds:
\begin{description}[(D)]
\item[\rm(a)]
For all $X, Y\in \cV(M)$,
there exists a uniquely determined
smooth vector field
$[X,Y]\in \cV(M)$
with local representations
\begin{equation}\label{formzeta}
[X,Y]_\phi\, =\, X_\phi.Y_\phi-Y_\phi.X_\phi\quad\mbox{for
all $\phi\in \cA$,}
\end{equation}
using the notation from {\rm (\ref{abbrnot})}.
\item[\rm(b)]
The map $[\cdot,\cdot]\colon  \cV(M)\times \cV(M)\to \cV(M)$
obtained from {\rm (a)}
is $\R$-bilinear and makes $\cV(M)$ a Lie
algebra.
\item[\rm(c)]
The $\R$-linear map
\[
\cL\colon \cV(M)\to\der(C^\infty(M,\R))\,,\quad
X\mto \cL_X
\]
$($with $\cL_X$ as in {\rm (\ref{firstl}))} is a homomorphism
of Lie algebras, i.e., 
\[ \cL_{[X,Y]}=\cL_X\circ \cL_Y-\cL_Y\circ \cL_X.\]
\end{description}
\end{thm}
To prepare the proof
of Theorem~\ref{Liebrvf},
we first study vector fields
on an open subset~$U$ of a locally
convex space~$E$.
Given $\xi,\eta\in C^\infty(U,E)$,
we abbreviate
\[
[\xi,\eta]\;:=\; \xi.\eta-\eta.\xi\,,
\]
using notation as in~(\ref{abbrnot}).\\[2.5mm]
Given open subsets $U\sub E$ and $V\sub F$
of locally convex spaces and a smooth map
$\phi\colon U\to V$, it is convenient to
say that two smooth functions
$\xi\colon U\to E$ and $\eta\colon V\to F$
are \emph{$\phi$-related}
if the corresponding vector fields
are so, or equivalently, if
\begin{equation}\label{conce}
\eta\circ \phi \;=\; d\phi \circ (\id_U,\xi)\,.
\end{equation}
It is useful to record
a formula for differentials
of maps like the ones in~(\ref{conce}).
\begin{lem}\label{lemhelpstud}
Let $E$, $F$ and $H$ be locally convex spaces,
$U\sub E$ and $V\sub H$ be open subsets,
$f\colon U\to F$ be a $C^2$-map
and $g\colon V\to U$ and $h\colon V\to E$ be~$C^1$.
Then the differential of the $C^1$-map
\[ \phi:=df\circ (g,h)\colon V\to F, \quad
  \phi(x)=df(g(x),h(x)) \]  is given by
\begin{equation}\label{helpstudents}
d\phi(x,y)\;=\; d^{\, (2)}f\big(g(x),h(x),dg(x,y)\big)+
df\big(g(x),dh(x,y)\big)
\end{equation}
for all $x\in V$ and $y\in H$.
\end{lem}
\begin{prf}
We have
\begin{eqnarray*}
d\phi(x,y) & = & d(d f)(g(x),h(x),dg(x,y),dh(x,y))\\
&=& d_1(d f)(g(x),h(x),dg(x,y))+d_2(df)(g(x),h(x), dh(x,y))\\
&=& d^{\, (2)}f(g(x),h(x),dg(x,y))
+df(g(x),dh(x,y))\,,
\end{eqnarray*}
using the Chain Rule (Proposition~\ref{chainC1}),
the formula $d(g,h)=(dg,dh)$
from Lemma~\ref{lemprod}
and the fact that
\[
d(df)(x,y_1;y_2,y_3)=d^{(2)}f(x,y_1,y_2)+df(x,y_3)
\] 
by~(\ref{dodgehgr}) in the proof of Proposition~\ref{higherdiff}.
\end{prf}
The next lemma captures
much of
Theorem~\ref{Liebrvf},
in the case of open subsets.
The following notation will be useful for the proof:
If $\fg$ is a vector space and
$\fg\times\fg\to\fg$, $(u,v)\mto [u,v]$
a bilinear map, we write
\begin{equation}\label{jacexpr}
\Jac(u,v,w):=
[u,[v,w]]+[v,[w,u]]+[w,[u,v]]\,.
\end{equation}
The Jacobi identity of a Lie algebra
now simply reads $\Jac(u,v,w)=0$.
\begin{lem}\label{preall}
Let
$E$ be a locally convex space and
$U\sub E$ be open.
Define $\cL_\xi(f):=df\circ (\id_U,\xi)$
for $\xi \in C^\infty(U,E)$
and $f\in C^\infty(U,\R)$.
Then we have:
\begin{description}[(D)]
\item[\rm(a)]
The map $\cL_\xi\colon C^\infty(U,\R)\to C^\infty(U,\R)$
is a derivation of $C^\infty(U,\R)$,
for each $\xi\in C^\infty(U,E)$.
\item[\rm(b)]
$\cL_{[\xi,\eta]}=
\cL_\xi\circ \cL_\eta-\cL_\eta\circ \cL_\xi$,
for all $\xi,\eta\in C^\infty(U,E)$.
\item[\rm(c)]
The map $\cL\colon C^\infty(U,E)\to\der(C^\infty(U,\R))$,
$\xi\mto \cL_\xi$
is linear and injective.
\item[\rm(d)]
The map $[\cdot,\cdot]\colon C^\infty(U,E)\times C^\infty(U,E)
\to C^\infty(U,E)$,
$(\xi,\eta)\mto [\xi,\eta]$ makes
$C^\infty(U,E)$ a Lie algebra.
\end{description}
\end{lem}
\begin{prf}
(a) has been checked in Remark~\ref{introdDX}.

(b) Let $f\in C^\infty(U,\R)$.
Since
\[
d(df\circ (\id_U,\eta))(x,y)=
d^{(2)}f(x,\eta(x),y)+
df(x,d\eta(x,y)) \;\mbox{ for all } \;(x,y)\in U\times E
\] 
by Lemma~\ref{lemhelpstud},
we deduce that
\[
\cL_\xi(\cL_\eta(f))
\, =\,  \cL_\xi(df\circ (\id_U,\eta))(x)
\, =\, d(df\circ (\id_U,\eta))
\circ (\id_U,\xi)
\]
is given by
\[
\cL_\xi(\cL_\eta(f))(x)\,=\,
d^{(2)}f(x,\eta(x),\xi(x))+
df(x,d\eta(x,\xi(x)))\,.
\]
Using this formula
and the corresponding formula
with $\xi$ and $\eta$ interchanged,
we obtain
\begin{eqnarray*}
(\cL_\xi\circ \cL_\eta-\cL_\eta\circ\cL_\xi)(f)(x)
&=& d^{(2)}f(x,\eta(x),\xi(x))+
df(x,d\eta(x,\xi(x)))\\
& & \;\;
- d^{(2)}f(x,\xi(x),\eta(x))-
df(x,d\xi(x,\eta(x)))\\
&=&
df(x,d\eta(x,\xi(x)))- df(x,d\xi(x,\eta(x)))\\
&=&
df(x,[\xi,\eta](x))
\;=\; \cL_{[\xi,\eta]}(f)(x)\,,
\end{eqnarray*}
as required (the second order terms
cancel by Schwarz' Theorem).

(c) $\cL_\xi$ is linear in $\xi$
as so is $df(x,\xi(x))$,
for all $f\in C^\infty(U,E)$ and $x\in U$.\linebreak
Since $\cL$ is linear,
it will be injective if
its kernel is trivial.
To verify the latter, let $0\not=\xi\in C^\infty(U,E)$.
There exists $x\in U$ such that $\xi(x)\not=0$.
By the Hahn--Banach theorem,
there exists
$\lambda\in E'$ such that
$\lambda(\xi(x))\not=0$ (see Theorem~\ref{dualsep}).
Thus $\cL_\xi(\lambda)(x)=d\lambda(x,\xi(x))=\lambda(\xi(x))\not=0$,
whence $\cL_\xi\not=0$.

(d) It is clear that $[\xi,\eta]$
is bilinear in $(\xi,\eta)$
(cf.\ Remark~\ref{rembil}).
Hence $(C^\infty(U,E),[\cdot,\cdot])$
is an algebra. Since $[\xi,\xi]=\xi.\xi-\xi.\xi=0$,
it only remains to check the Jacobi identity,
$\Jac(\xi,\eta,\zeta)=0$
for all $\xi,\eta,\zeta\in C^\infty(U,E)$
(using the notation from~(\ref{jacexpr})).
However,
\[
\cL_{\Jac(\xi,\eta,\zeta)}\;=\;
\Jac(\cL_\xi,\cL_\eta,\cL_\zeta)\;=\;0\,,
\]
using~(b) to obtain the first equality and
the fact that $\der(C^\infty(U,\R))$
is a Lie algebra to obtain the second.
Since $\cL\colon C^\infty(U,E)\to\der(C^\infty(U,\R))$
is injective, it follows that
$\Jac(\xi,\eta,\zeta)=0$.
\end{prf}
The following lemma is an efficient
tool, which will be applied
frequently. It accounts
for the usefulness of the concept
of $f$-relatedness.
\begin{lem}[Lemma on related vector fields on open subsets]\label{larelloc}
Let $U\sub E$ and $V\sub F$ be open subsets
of locally convex spaces. Let
$\phi\colon U\to V$ as well as
$\xi_1,\xi_2\colon U\to E$ and
$\eta_1,\eta_2\colon V\to F$ be smooth
functions.
If $\xi_1$ is $\phi$-related to~$\eta_1$
and $\xi_2$ is $\phi$-related to~$\eta_2$,
then $[\xi_1,\xi_2]$ is $\phi$-related
to $[\eta_1,\eta_2]$.
\end{lem}
\begin{prf} We claim that $\xi \in C^\infty(U,E)$ and $\eta \in C^\infty(V,F)$ 
are $\phi$-related if and only if 
\begin{equation}
  \label{eq:lie-rel}
\cL_\xi \circ \phi^* = \phi^* \circ \cL_\eta \: 
C^\infty(V,\R) \to C^\infty(U,\R), \quad \mbox{ where }  \quad 
\phi^*f = f \circ \phi.
\end{equation}
As 
\[ \cL_\xi(f \circ \phi)(x) 
= d(f \circ \phi)(x, \xi(x)) 
= df(\phi(x), d\phi(x,\xi(x))) \] 
and 
\[ (\cL_\eta f)(\phi(x)) = d f(\phi(x), \eta(\phi(x))),\] 
we see that if $\xi$ and $\eta$ are $\phi$-related, then 
\eqref{eq:lie-rel} holds. Conversely, we may apply 
\eqref{eq:lie-rel} to restrictions of continuous linear functionals to see that 
it implies \eqref{conce}, i.e., that $\xi$ is $\phi$-related to~$\eta$. 

Now our assumptions imply
\[ \cL_{\xi_1}\circ \cL_{\xi_2} \circ \phi^* 
= \cL_{\xi_1} \circ \phi^* \circ \cL_{\eta_2} 
= \phi^* \circ \cL_{\eta_1} \circ \cL_{\eta_2}, \] 
so that we obtain with Lemma~\ref{preall}(b) the relation 
\[ \cL_{[\xi_1, \xi_2]} \circ \phi^* 
= [\cL_{\xi_1}, \cL_{\xi_2}] \circ \phi^* 
= \phi^* \circ [\cL_{\eta_1}, \cL_{\eta_2}] 
= \phi^* \circ \cL_{[\eta_1, \eta_2]},\] 
and our preceding discussion shows that this means that 
$[\xi_1,\xi_2]$ is $\phi$-related to $[\eta_1, \eta_2]$. 
\end{prf}
\emph{Proof of Theorem}~\ref{Liebrvf}.
(a)
We only need to show that
the family $([X_\phi,Y_\phi])_{\phi\in \cA}$
is compatible
in the sense described in
Lemma~\ref{ladealwvf1}.
To this end, let $\phi,\psi\in \cA$ be charts.
Since $X_\phi$ and $X_\psi$
(restricted to $\phi(U_\phi\cap U_\psi)$,
resp., $\psi(U_\phi\cap U_\psi)$)
are $\tau:=\psi\circ \phi^{-1}$-related
and also the corresponding restrictions
of $Y_\phi$ and~$Y_\psi$
are $\tau$-related,
it follows from Lemma~\ref{larelloc}
that also the restrictions of
$[X_\phi,Y_\phi]$ and $[X_\psi, Y_\psi]$
are $\tau$-related.
This is the required compatibility condition.

(b) Since $[X,X]_\phi=[X_\phi,X_\phi]=0$
for each $\phi\in\cA$ by
Lemma~\ref{preall},
we have $[X,X]=0$.
Likewise,
the Jacobi identity can be tested
in charts.
There, it is valid by Lemma~\ref{preall}.

(c)
We have
\begin{eqnarray*}
\cL_X(f)\circ \phi^{-1}
& = & df\circ X \circ\phi^{-1}
\; =\; df\circ T\phi^{-1}\circ (\id_{V_\phi},X_\phi)\\
& = & d(f\circ \phi^{-1})\circ (\id_{V_\phi},X_\phi)
\; =\; \cL_{X_\phi}(f\circ \phi^{-1})
\end{eqnarray*}
for all $f\in C^\infty(M,\R)$,
$X\in \cV(M)$ and $\phi\in \cA$.
Hence
\[ \cL_{[X,Y]}(f)\circ \phi^{-1}
=\cL_{[X_\phi,Y_\phi]}(f\circ \phi^{-1})
=[\cL_{X_\phi},\cL_{Y_\phi}](f\circ \phi^{-1})
=[\cL_X,\cL_Y](f)\circ \phi^{-1},\]
for all
$X,Y\in \cV(M)$,
$f\in C^\infty(M,\R)$
and $\phi\in \cA$.
Since~$\phi$ was arbitrary,
it follows that
$\cL_{[X,Y]}(f)=[\cL_X,\cL_Y](f)$,
as required.\qed

There also is an (even more important) global version
of Lemma~\ref{larelloc}.
\begin{lem}[Lemma on related vector fields]\label{larelglob}
Let
$f\colon M\to N$ be a smooth map
between smooth manifolds
and $X^1,X^2\in \cV(M)$
as well as $Y^1,Y^2\in \cV(N)$ be
vector fields.
If $X^1$ is $f$-related to~$Y^1$
and $X^2$ is $f$-related to~$Y^2$,
then $[X^1,X^2]$ is $f$-related
to $[Y^1,Y^2]$.
\end{lem}
\begin{prf} We have to show that 
\[ Tf \circ [X^1, X^2] = [Y^1, Y^2] \circ f \]
holds pointwise, at every $m \in M$. To this end, we choose charts 
$\phi$ of $M$ around $m$ and $\psi$ of $N$ around $f(m)$. 
We may assume that $f(U_\phi) \subeq U_\psi$. 
Then $X^j|_{U_\phi}$ is $\phi$-related to $X^j_\phi$, the map
$X^j_\phi$ is $(\psi \circ f \circ \phi^{-1})$-related to $Y^j_\psi$, and 
$Y^j|_{U_\psi}$ is $\psi$-related to $Y^j_\psi$. 
By construction of the Lie bracket on vector fields, 
$[X^1, X^2]$ is $\phi$-related to 
$[X^1_\phi, X^2_\phi]$ and 
$[Y^1, Y^2]$ is $\psi$-related to 
$[Y^1_\psi, Y^2_\psi]$. 
Further, Lemma~\ref{larelloc} shows that 
$[X^1_\phi, X^2_\phi]$ is $(\psi \circ f \circ \phi^{-1})$-related 
to $[Y^1_\psi, Y^2_\psi]$. Combining these relations, we see that 
$[Y^1, Y^2]$ is $f$-related to $[X^1, X^2]$. 
\end{prf}

\begin{rem}\label{fdvfvsder}
If $M$ is a finite-dimensional smooth manifold
over~$\R$,
then the linear map
\[
\cL\colon \cV(M)\to \der(C^\infty(M,\R))
\,,\quad
X\mto \cL_X
\]
is a bijection, whence the Lie algebra
structure of $\der(C^\infty(M,\R))$ can be transported
to $\cV(M)$.
Unfortunately,
$\cL$ is not a bijection in general
if $M$ is a smooth manifold modeled
on a locally convex space,
whence we cannot turn $\cV(M)$
as easily into a Lie algebra in this case
as in the case where $M=U$ is an open
subset of a locally convex space
(discussed in Lemma~\ref{preall}).
See \cite[Ex.~8.3]{BGN04}
for an example where~$\cL$ is not surjective
(based on \cite[Lemma~28.4]{KM97}).
\end{rem}
We do not know the answer to the following problem.
\begin{probl}\label{mfd-prob-1}
Is the map
\begin{equation}\label{givesinje}
\cV(M)\to \der(C^\infty(M,\R)),\;\;
X\mto \cL_X
\end{equation}
injective, for every
smooth manifold~$M$ modeled on a locally
convex space\,?
\end{probl}
\begin{rem}
Let $M$ be a smooth manifold such that,
for each $p\in M$, the differentials $df|_{T_pM}$ with $f\in C^\infty(M,\R)$
separate points on $T_pM$.
Then the map (\ref{givesinje}) is injective.\footnote{If $0\not=X\in \cV(M)$,
then $X(p)\not=0$ for some $p\in M$. The hypothesis provides $f\in C^\infty(M,\R)$ with
$0\not=df(X(p))=\cL_X(f)(p)$.}
\end{rem}
\begin{probl}\label{mfd-prob-2}
Do the differentials $df|_{T_pM}$ for $f\in C^\infty(M,\R)$ separate points on~$T_pM$,
for each smooth manifold~$M$ and each $p\in M$\,?
\end{probl}
A related problem is the following.
\begin{probl}\label{mfd-prob-3}
Does $C^\infty(M,\R)$ separate points on~$M$,
for each smooth manifold~$M$\,?
\end{probl}
For smoothly regular smooth manifolds
discussed in Section~\ref{sec-mfd-boundary} (like finite-dimensional ones),
all of Problems~\ref{mfd-prob-1}--\ref{mfd-prob-3}
have positive solutions.
\begin{rem}\label{othervectorf}
%
%
(a) If $M$ is a $C^r$-manifold with $r\in\N\cup\{\infty\}$
and $k\in \N_0$ with $k\leq r-1$,
we say that a $C^k$-map $X\colon M\to TM$ is a \emph{$C^k$-vector field}
if $\pi_{TM}\circ X=\id_{TM}$.
The set of all $C^k$-vector fields is a real vector space $\cV_{C^k}(M)$
under pointwise operations,
as in Definition~\ref{defVF}.
Their local representatives are $C^k$-functions $X_\phi\colon V_\phi\to E$
to the modeling space~$E$ of~$M$.
If $f\colon M\to N$ is a $C^r$-map to a $C^r$-manifold~$N$,
we say that $C^k$-vector fields $X\in \cV_{C^k}(M)$ and $Y\in \cV_{C^k}(N)$ are \emph{$f$-related}
if $Y\circ f=Tf\circ X$. An analog of Remrak~\ref{rmsrel}(b) holds.
Given $X\in \cV_{C^k}(M)$ and a $C^r$-function $f\colon M\to F$ to a locally convex space,
we can define a $C^k$-map $X.f\colon M\to F$ via $X.f:=df\circ X$.

(b) If $\K\in\{\R,\C\}$ and $M$ is a $\K$-analytic manifold,
we say that a \break $\K$-analytic map $X\colon M\to TM$ is a
\emph{$\K$-analytic vector field} if $\pi_{TM}\circ X=\id_{TM}$.
The set $\cV_{C^\omega}(M)$ of all $\K$-analytic vector fields is a $\K$-vector
space under pointwise operations and a Lie subalgebra of the Lie algebra
$\cV(M)$ of all smooth vector fields, as $X_\phi.Y_\phi-Y_\phi.X_\phi$
is $\K$-analytic for all $X,Y\in \cV_{C^\omega}(M)$ and
each chart $\phi$ of~$M$. Thus $\cV_{C^\omega}(M)$
is a Lie algebra over~$\K$ (noting that the Lie bracket is complex bilinear
if $\K=\C$ as $X_\phi.Y_\phi-Y_\phi.X_\phi$
is complex bilinear in $(X,Y)$ for each chart~$\phi$).
If $f\colon M\to F$ is a $\K$-analytic map to a locally convex topological
$\K$-vector space~$F$, then also the map $X.f:=df\circ X$ is $\K$-analytic
(and $\K$-linear in~$f$).
If $F:=\K$, then $X.f=\cL_X(f)$ defines a derivation of $C^\omega(M,\K)$
as a function of $f\in C^\omega(M,\K)$,
and we obtain a Lie algebra homomorphism $\cV_{C^\omega}(M)\to \der(C^\omega(M,\K))$.
\end{rem}
\begin{small}
\subsection*{Exercises for Section~\ref{sec:2.4}} 
\label{Exer3.4}

\begin{exer}
  \label{exer:2.4.1a} Let $f\colon M\to N$
be a $C^r$-map between $C^r$-manifolds
and let $X\colon M\to TM$ and $Y\colon N\to TN$
be $C^r$-vector fields. Show:
\begin{description}[(D)]
\item[(a)]
If $X$ and $Y$ are $f$-related,
then
$(Y.\gamma)\circ f= X.(\gamma\circ f)$
for each $C^r$-map $\gamma\colon N\to\K$.
\item[(b)]
$X$ and $Y$ are $f$-related
if and only if
$(Y.\gamma)\circ f|_{f^{-1}(U)}= X.(\gamma\circ f|_{f^{-1}(U)})$
for each open subset $U\sub N$
and each
$C^r$-map $\gamma\colon U\to\K$.
\end{description}
\end{exer}

\begin{exer}
  \label{exer:2.4.2} 
Given a vector field
$X\in \cV(M)$, define $\Supp(X)$
as the closure of the set 
$\{p\in M\colon X(p)\not= 0\}$.
Show that, for a finite-dimensional smooth manifold $M$, the set $\cV(M)_c$
of all vector fields with compact support
(``compactly supported vector fields'')
is an ideal of the Lie algebra
$\cV(M)$.
\end{exer}

\begin{exer}
  \label{exer:2.4.3}
Let $E$ be a locally
convex space, $U\sub E$ be open
and let {$\xi,\eta,\zeta\colon U\to E$}
be smooth maps
(which we interpret as
local representatives of smooth
vector fields on~$U$).
Then
\[ [\xi,\eta](x):=(\xi.\eta-\eta.\xi)(x)
=d\eta(x,\xi(x))-d\xi(x,\eta(x)) \] 
defines $[\xi,\eta]\colon U\to E$.
Verify by direct calculation
that $[.,.]$ makes $C^\infty(U,E)$
a Lie algebra
(in particular, verify the Jacobi identity).
\end{exer}

\begin{exer}
  \label{exer:2.4.4} 
Let $U\sub \R^2$ be an open subset,
$X=f_1\frac{\partial}{\partial x_1}
+f_2\frac{\partial}{\partial x_2}$ and
$Y=g_1\frac{\partial}{\partial x_1}
+g_2\frac{\partial}{\partial x_2}$
be smooth vector fields
on~$U$ (see Remark~\ref{vecfifin}).
Then $[X,Y]=h_1\frac{\partial}{\partial x_1}
+h_2\frac{\partial}{\partial x_2}$
with suitable smooth functions
$h_1,h_2\colon U\to \R$.
Calculate $h_1$ and $h_2$ explicitly
in terms of $f_1,f_2,g_1$ and $g_2$.
\end{exer}
\end{small}
\section{Manifolds with boundary and generalizations}\label{sec-mfd-boundary}
Various objects of interest
(like the closed unit ball in~$\R^n$,
or $[0,1]^2$ in~$\R^2$)
are not manifolds.
In order to capture these
examples, one has to loosen
the definition of a manifold.
As we have a differential calculus
available for mappings on non-open subsets,
this generalization does not pose any problems.
In this section,
we introduce
``manifolds with rough boundary''
modeled on locally convex subsets
with dense interior
of a given locally convex space.
These subsume manifolds modeled on locally convex spaces
(as in the preceding sections),
convex subsets of locally convex spaces with non-empty interior,
and traditional manifolds with boundary or corners (as recalled in Remark~\ref{remtradi}).
\begin{defn}\label{defroughbd}
Let $E$ be a locally convex space
and $r\in \N_0\cup\{\infty\}$.
\begin{description}
\item[(a)]
A \emph{rough $E$-chart} for
a Hausdorff topological space~$M$
is a homeomorphism
$\phi\colon U_\phi\to V_\phi$
from an open subset $U_\phi\sub M$
onto a locally convex subset $V_\phi\sub E$
with dense interior.
If $r\geq 1$, then
a \emph{$C^r$-atlas of rough $E$-charts}
for~$M$ is a set $\cA$ of rough
$E$-charts $\phi\colon U_\phi\to V_\phi$
for~$M$ such that
(A1) and (A2) from Definition~\ref{defnatlas}
hold.
If $r=0$, we assume in addition:\vspace{1mm}
\begin{description}[(DD)]
\item[(A3)]
For all $\phi,\psi\in \cA$ and
$x\in U_\phi\cap U_\psi$,
we have $\phi(x)\in \partial V_\phi$
if and only if $\psi(x)\in \partial V_\psi$
(where $\partial V_\phi=\wb{V_\phi}\setminus V_\phi^0$
is the topological boundary of $V_\phi$
in~$E$).\vspace{.5mm}
\end{description}
\item[(b)]
A \emph{$C^r$-manifold with rough boundary}
is a Hausdorff topological space~$M$,
together with a maximal
$C^r$-atlas~$\cA$ of rough
$E$-charts for~$M$.
Once~$\cA$ has been fixed,
we refer to the
rough $E$-charts $\phi\in \cA$ simply as
\emph{charts}.
\item[(c)]
If $(M,\cA)$ is a $C^r$-manifold
with rough boundary modeled
on~$E$, we call $x\in M$
a (formal) \emph{boundary
point}
if $\phi(x)\in \partial V_\phi$
for some $\phi\colon U_\phi\to V_\phi$
in~$\cA$ such that $x\in U_\phi$.
The set
$\partial M$
of all boundary points
is called the (formal) \emph{boundary of~$M$.}
\end{description}
\end{defn}
\begin{rem}
Note that $U_\phi\cap U_\psi$ is open in~$U_\phi$
in the situation of (A2),
whence $\psi(U_\phi\cap U_\psi)$ is open
in $V_\psi$.
Consequently, $\psi(U_\phi\cap U_\psi)$ is
a locally convex subset of~$E$ with dense interior.
Therefore
the $C^r$-property
of the transition map in (\ref{transi})
makes sense.
\end{rem}
\begin{ex}\label{exarough}
Manifolds with rough boundary arise
naturally:
\begin{description}[(D)]
\item[(a)]
It is clear that
every $C^r$-manifold $(M,\cA)$
(where $r\not=\omega$) has an underlying
$C^r$-manifold with rough boundary,
obtained by replacing~$\cA$
by the maximal $C^r$-atlas $\wh{\cA}$ of rough $E$-charts
containing~$\cA$.
%
\item[(b)]
If~$E$ is a locally convex space,
then every locally convex subset $U\sub E$
with non-empty interior
is a $C^\infty$-manifold with
rough boundary when equipped
with the maximal $C^\infty$-atlas
of rough $E$-charts
containing $\{\id_U\}$.
Thus $[0,1]^n\sub \R^n$ is a smooth
manifold with rough
boundary modeled on~$\R^n$, and so is
the standard simplex 
\[ \Delta_n:=\Big\{(x_1,\ldots, x_n)\in [0,1]^n \colon \sum_{k=1}^n x_k\leq 1\Big\}.\] 
\end{description}
\end{ex}
\begin{ex}\label{exmfd3}
If $(M,\cA)$ is a $C^r$-manifold
with rough boundary, then every open subset
$U\sub M$
can be made a $C^r$-manifold
with rough boundary,
as in Example~\ref{exmfd2}.
\end{ex}
\begin{ex}\label{exmprodvar}
Let $(M_1,\cA_1)$ and $(M_2,\cA_2)$
be $C^r$-manifolds with rough boundary,
modeled on locally convex spaces $\K$-vector spaces $E_1$
and $E_2$, respectively.
As in the case of $C^r$-manifolds without boundary,
we see that
\[ \cA:=\{\phi_1\times \phi_2\colon \phi_1\in \cA_1,\phi_2\in \cA_2\}\] 
is a $C^r$-atlas of rough $(E_1\times E_2)$-charts for
$M_1\times M_2$,
endowed with the product topology
(if $r=0$, then the additional condition (A3)
is easily verified).
Hence the maximal $C^r$-atlas of rough $(E_1\times E_2)$-charts
containing~$\cA$
makes $M_1\times M_2$ a $C^r$-manifold
with rough boundary,
the \emph{direct product} of~$M_1$ and~$M_2$.
\end{ex}
The next lemma shows that condition
(A3) is automatic if $r\geq 1$.
As a consequence, boundaries are what
they should be.
\begin{lem}\label{goodbdy}
Let $(M,\cA)$
be a $C^r$-manifold
with rough boundary $($where
$r\in \N\cup\{\infty\})$,
and $x\in M$.
Let $\phi\colon U_\phi\to V_\phi$
and $\psi\colon U_\psi\to V_\psi$
be charts for~$M$
around~$x$.
Then $\phi(x)\in \partial V_\phi$
if and only if $\psi(x)\in \partial V_\psi$.
\end{lem}
\begin{prf}
We may assume that $\K=\R$
and $E\not=\{0\}$.
Since the roles of~$\phi$
and~$\psi$ can be interchanged,
it suffices to show
that if $z:=\psi(x)\in \partial V_\psi$,
then also $y:=\phi(x)\in \partial V_\phi$.
We assume to the contrary
that $z \in \partial V_\psi$
but $y$ is in the interior
$V_\phi^0$ of~$V_\phi$,
and derive a contradiction.
Since~$V_\psi$ is locally
convex, $z$ has
a convex neighborhood $W\sub V_\psi$.
Then $z\in \partial W$ and hence
$z\not\in W^0$, where $W^0$
is the interior of~$W$ in~$E$.
Applying the Hahn--Banach Theorem
(Theorem~\ref{geomHBreal}) to the convex
open subset $W^0-z$, we find
$\lambda\in E'$ such that $\lambda(W^0-z)\sub \;]0,\infty[$,
whence $\lambda(w)>\lambda(z)$
for each $w\in W^0$
and $\lambda(w)\geq \lambda(z)$
for each $w\in W$.
There exists $v\in E$ such that $\lambda(v)<0$.
Since
$\kappa := \psi\circ \phi^{-1}\colon
\phi(U_\phi\cap U_\psi)\to V_\psi$
and $\kappa^{-1}=\phi\circ \psi^{-1}$
are~$C^r$, the Chain
Rule entails that $\kappa'(y)\colon E\to E$
is an invertible linear map.
Hence there exists
$u\in E$ such that $\kappa'(y)u=v$.
Since $\kappa^{-1}(W)$ is a neighborhood
of~$y$ in~$V_\phi$
and hence also in~$E$
(as $y\in V_\phi^0$),
we find $\ve>0$ such that
$y+tu\in \kappa^{-1}(W)$
for each $t\in \;]{-\ve},\ve[$.
Then $\lambda(\kappa(y+tu))
\geq \lambda(z)$
for  $|t| < \ve$,
whence $\lambda(\kappa'(y)u)
=
\frac{d}{dt}\big|_{t=0}
\lambda(\kappa(y+tu))= 0$.
But
$\lambda(\kappa'(y)u)
=\lambda(v)<0$, by choice of~$u$.
We have reached a contradiction.
\end{prf}
\begin{rem}
The preceding lemma entails
that $\partial M=\emptyset$
(and $\wh{\cA}=\cA$)
in the situation of Example~\ref{exarough}(a).
It also ensures that the formal boundary
of $U$
arises from the topological
boundary as $\partial U\cap U$
in the situation of
Example~\ref{exarough}(b).
In particular, the formal
boundary of~$U$ coincides with
its topological boundary $\partial U$ in~$E$
if~$U$ is closed.
\end{rem}
\begin{rem}
We recall \emph{Invariance of Domain} for $E=\R^n$:
\emph{If $U\sub E$ is an open subset and $f\colon U\to E$ an injective continuous map,
then $f(U)$ is open in~$E$ and $f\colon U\to f(U)$ is a homeomorphism}
(see, e.g., \cite[Thm.~2B.3]{Ha02}).
As a consequence, condition (A3) in Definition~\ref{defroughbd}
is automatically satisfied (and could be omitted)
if $r=0$ and~$E$ is finite-dimensional.
As many examples of infinite-dimensional normed spaces are homeomorphic to a proper vector subspace,
invariance of domain becomes false for such (although it may still hold for special classes of continuous functions $f$, see \cite{Sr29} and \cite{Tr72}).
In~\cite{Kr76}, one finds an easy example of
a continuous bijection from the open unit ball  in
$\ell^2$ onto itself which fails to be an open map.
\end{rem}
\begin{rem}
If $r=0$ and $E$ is infinite-dimensional, then condition (A3) cannot be omitted, because of
well-known pathologies.\medskip

\noindent
(a) Every closed convex subset $U$ with non-empty interior in an infinite-dimensional Fr\'{e}chet space~$E$
is homeomorphic to~$E$ (see \cite[Ch.~VI, Thm.~62]{BP75}).\footnote{That $\ell^2$ is homeomorphic to a half space in $\ell^2$ was already shown in~\cite{Kl53}.}
As every $x\in E$ can be mapped to any $y\in E$ using a translation, we
deduce that also $U$ is a homogeneous
topological space (i.e., for all $x,y\in U$, there exists a homeomorphism
$\phi\colon U\to U$ such that $\phi(x)=y$). In particular, for each $x\in\partial U$ and $y$ in the interior of
$U$ we have $\phi(x)=y$ for some homeomorphism $\phi\colon U\to U$.\medskip

\noindent
(b)
In particular, (a) applies to the closed unit ball $U$ in an infinite-dimensional Banach space $(E,\|\cdot\|)$.
We mention that the latter is homeomorphic to the open ball (as the latter is also homeomorphic to~$E$,
by \cite[Ch.~III, Thm.~6.2]{BP75}).
\end{rem}
\begin{rem}\label{remtradi}
The traditional manifolds with boundary
and manifolds with corners
from finite-dimensional
differential geometry
can be considered as special
cases of manifolds with rough boundary:
\begin{description}[(D)]
\item[(a)]
A \emph{manifold with boundary}
is a $C^r$-manifold~$M$
with rough boundary
modeled on some $\R^n$,
whose atlas contains
an atlas of charts
whose ranges are open
subsets of $\R^{n-1}\times [0,\infty[$.
\item[(b)]
A \emph{manifold with corners}
is a $C^r$-manifold~$M$
with rough boundary
modeled on some $\R^n$,
whose atlas contains
an atlas of charts
$\phi\colon U_\phi\to V_\phi$
whose ranges $V_\phi$ are open
subsets of $[0,\infty[^n$.
\end{description}
\end{rem}
We define $C^r$-maps between $C^r$-manifolds with rough boundary
as in the case of ordinary manifolds (without boundary).
\begin{defn}\label{defnCrmfbd}
Let $M$ and $N$
be $C^r$-manifolds
with rough boundary
modeled on locally convex spaces $E$ and $F$,
respectively.
A map $f\colon M\to N$ is called $C^r$
if $f$ is continuous and, for all charts $\phi$ of $M$ and $\psi$ of
$N$, the map
\begin{equation}\label{ugglybd}
\psi\circ f \circ \phi^{-1}\colon E\, \supseteq \,
\phi(f^{-1}(U_\psi)\cap U_\phi)\to F
\end{equation}
is~$C^r$.
Note that the domain
$\phi(f^{-1}(U_\psi)\cap U_\phi)$ of the composition is an open subset
of~$V_\phi$ and hence a locally convex subset of~$E$
with dense interior.\footnote{Thus, pathologies as in Exercise~\ref{exc-patho-nonop}(a) cannot occur here.}
We write $C^r(M,N)$ for the set
of all $C^r$-maps from~$M$ to~$N$.
\end{defn}
As in the case of ordinary manifolds, compositions of $C^r$-maps are~$C^r$.
\begin{prop}\label{compmfdsbd}
Let $M$, $N$ and $X$ be $C^r$-manifolds
with rough boundary.
If $f\colon M\to N$ and $g\colon N\to X$
are $C^r$-maps,
then also
their composition
$g\circ f \colon M\to X$
is $C^r$.
\end{prop}
\begin{prf}
The map $g\circ f$ is continuous.
Given charts $\phi\colon U_\phi\to V_\phi$ of $M$ and
$\psi\colon U_\psi\to V_\psi$
of $X$,
let $x$ be an element of the domain $\phi((g\circ f)^{-1}(U_\psi)\cap U_\phi)$
of the composition $\psi\circ g\circ f \circ \phi^{-1}$,
which is an open subset of~$V_\phi$ and hence a locally convex subset of the modeling
space with dense interior.
Let $\kappa$ be a chart of~$N$
around $f(\phi^{-1}(x))$.
Then
\[
\psi\circ (g\circ f)\circ\phi^{-1}=(\psi\circ g\circ \kappa^{-1})
\circ (\kappa\circ f\circ \phi^{-1})
\]
holds on some open neighborhood of~$x$
in~$V_\phi$,
and this is a $C^r$-map
as it is a composition of $C^r$-maps between
locally convex subsets of locally convex spaces
with dense interior.
Therefore the mapping
$\psi\circ (g\circ f)\circ \phi^{-1}$
is locally~$C^r$ and hence~$C^r$.
Thus $g\circ f$ is~$C^r$.
\end{prf}
\begin{rem}\label{alsowithbd}
Many results concerning $C^r$-maps formulated in
Section~\ref{sec:2.2} for ordinary manifolds (without boundary)
carry over directly to $C^r$-maps between manifolds with rough boundary,
for $r\not=\omega$.
Thus, Remark~\ref{oneatlasenough}, Example~\ref{coprojcr},
Lemma~\ref{mapinprodmfd} remain valid if
$M$, $N$, $M_1$, and $M_2$, respectively, are $C^r$-manifolds with rough boundary.
Example~\ref{exmpCr1} remains valid if we consider a locally convex subset $U\sub E$
with dense interior as a $C^r$-manifold with rough boundary.
We mention that $C^r$-diffeomorphisms and local \break {$C^r$-diffeomorphisms} between $C^r$-manifolds with rough boundary can be defined as in Definitions~\ref{secdefdiff} and \ref{secdeflocdiff}, respectively,
replacing $M$ and $N$ with \break $C^r$-manifolds with rough boundary there.
Exercise~\ref{exc-chartdiffeo} remains valid if~$M$ is a $C^r$-manifold with
rough boundary and $V\sub E$ a locally convex subset with dense interior.
\end{rem}
For manifolds with rough boundary, we define submanifolds as follows.
\begin{defn}\label{submfdbd}
Let $r\in\N_0\cup\{\infty\}$ and~$M$ be a $C^r$-manifold with rough boundary modeled
on a locally convex space~$E$. Let $F\sub E$ be a closed vector subspace. A subset $N\sub M$ is called a \emph{$C^r$-submanifold of~$M$ modeled on~$F$}
if it has the following properties:
\begin{description}[(S2)]
\item[\rm(S1)]
For each $x\in N$, there exists a chart $\phi\colon U_\phi\to V_\phi\sub E$ of~$M$ around~$x$
which is \emph{adapted to~$N$} in the sense that
\[
\phi(U_\phi\cap N)=V_\phi\cap F
\]
and the interior of $V_\phi\cap F$ relative~$F$ is dense in $V_\phi\cap F$.
\item[\rm(S2)]
For all $x\in N$ and charts $\phi\colon U_\phi\to V_\phi$ and $\psi\colon U_\psi\to V_\psi$
as in {\rm(S1)}, we have $\phi(x)\in \partial (V_\phi\cap F)$ relative~$F$
if and only if $\psi(x)\in \partial (V_\psi\cap F)$.
\end{description}
\end{defn}
Then $N$ is a $C^r$-manifold with rough boundary when endowed with the topology induced by~$M$
and the maximal $C^r$-atlas of rough $F$-charts that contains
$\phi_N:=\phi|_{U_\phi\cap N}\colon U_\phi\cap N\to V_\phi\cap F$
for all charts~$\phi$ of~$M$ as in (S1).\\[2.3mm]
%
For example, every open subset $N\sub M$ is a submanifold modeled on~$E$.
\begin{numba} ($C^0$-submanifolds with prescribed boundary).
In the case $r=0$, assume that (S1) is satisfied.
Let $B$ be a subset of $N$ and assume the set
$\cA_B$ of all rough $F$-charts of the form $\phi_N$ satisfying
\[
(\forall x\in U_\phi\cap N)\;\, x\in B\;\,\Leftrightarrow \;\, \mbox{$\phi_N(x)\in\partial (V_\phi\cap F)$
relative~$F$}
\]
is a rough $C^0$-atlas for~$N$. Then $N$, endowed with the topology induced by~$M$ and the maximal rough $C^0$-atlas containing $\cA_B$, is a $C^0$-manifold with rough boundary.
We refer to it as a \emph{submanifold with boundary $B$}.
\end{numba}
For example, consider a $C^0$-manifold with rough boundary modeled on~$E$,
a submanifold $N$ of $M$ and a locally convex space~$Y$.
Then $N\times Y$ can be considered as a submanifold of $M\times Y$
modeled on $F\times Y$ with boundary $(\partial N)\times \id_Y$,
using the atlas of rough $F\times Y$-charts $\phi_N\times\id_Y$.
\begin{rem}
(a) Condition~(S2) is automatic if $r\geq 1$ and can be omitted in this case,
by Lemma~\ref{goodbdy}. If $r=0$, we have to spell it out.
By Invariance of Domain, (S2) can also be omitted if $r=0$ and $F$
has finite dimension.
\medskip

\noindent
(b) Note that $V_\phi\cap F$ is locally convex in the situation of (S1),
being an intersection of locally convex sets (see Remark~\ref{basiclcxs}).\medskip

\noindent
(c) If $V_\phi$ is open in~$E$ (e.g., if $M$ does not have a boundary), then
$V_\phi\cap F$ is open in~$F$,
whence it has dense interior (automatically).\medskip

\noindent
(d) If $V_\phi$ is not open, then the density of the interior of $V_\phi\cap F$ relative~$F$
is a proper requirement which is not automatic
(cf.\ Exercise~\ref{exc-patho-nonop}(c)).\medskip

\noindent
(e) Lemma~\ref{mapsubmfd} remains valid if $r\in \N_0\cup\{\infty\}$ and $M$, $N$ and $P$
are $C^r$-manifolds with rough boundary,
replacing the references to Definition~\ref{defsubm} and Proposition~\ref{compmfds}
in its proof by Definition~\ref{submfdbd} and Proposition~\ref{compmfdsbd}.
\end{rem}
Graphs of mappings to manifolds without boundary
provide important examples of
submanifolds.
\begin{prop}\label{graphsubm}
Given $r\in\N_0\cup\{\infty\}$, let $M$ be a $C^r$-manifold $($possibly with rough boundary$)$,
modeled on a locally convex space~$E$, and $N$ be a $C^r$-manifold $($without boundary$)$,
modeled on a locally convex space~$F$. Let $f\colon M\to N$ be a $C^r$-map.
Then the graph
\[
\graph(f):=\{(x,f(x))\colon x\in M\}
\]
is a submanifold of $M\times N$ modeled on the complemented vector subspace $E\times \{0\}$
of $E\times F$, with boundary $\partial\graph(f)=\{(x,f(x))\colon x\in\partial M\}$.
%
\end{prop}
\begin{proof}
Given $x\in M$, let $\psi\colon U_\psi\to V_\psi$ be a chart for~$N$ with
$f(x)\in U_\psi$, such that $\psi(f(x))=0$; here $V_\psi$ is an open subset of~$F$.
There exists a chart
$\phi\colon U_\phi\to V_\phi\sub E$ for~$M$ such that $x\in U_\phi$,
$\phi(x)=0$ and $f(U_\phi)\sub U_\psi$; here $V_\phi$ is a locally convex subset of~$E$
with dense interior. Then $V_\phi\times F$ is a locally convex subset of $E\times F$
with dense interior. The map
\[
h\colon V_\phi\to F,\;\; z\mto \psi(f(\phi^{-1}(z)))
\]
is $C^r$.  As a consequence, the map
\[
g\colon V_\phi\times F\to V_\phi\times F,\;\, (z,y)\mto (z,y-h(z))
\]
is $C^r$ and in fact a $C^r$-diffeomorphism, as the $C^r$-map
\[
V_\phi\times F\to V_\phi\times F,\;\, (z,y)\mto (z,y+h(z))
\]
is its inverse. Thus
\[
W:=g(V_\phi\times V_\psi)
\]
is an open subset of $V_\phi\times F$ and hence a locally convex subset of $E\times F$ with
dense interior. Moreover,
\[
\theta:=g\circ (\phi\times \psi)\colon U_\phi\times U_\psi\to W
\]
is a $C^r$-diffeomorphism and hence a chart for $M\times N$ (see Exercise~\ref{exc-chartdiffeo}
and Remark~\ref{alsowithbd}).
By construction, we have
\[
\theta((U_\phi\times U_\psi)\cap\graph(f))=W\cap (E\times \{0\})
=V_\phi\times\{0\},
\]
which is a locally convex subset of $E\times \{0\}$
with dense interior $V_\phi^0\times\{0\}$ relative $E\times\{0\}$.
Thus $\theta$ is adapted to $\graph(f)$,
whence condition (S1) is satisfied.
Relative $E\times\{0\}$, we have $\partial (V_\phi\times\{0\})=(\partial V_\phi)\times\{0\}$.
For $(x,y)\in (U_\phi\times U_\psi)\cap \graph(f)$,
we therefore have $\theta(x,y)\in \partial(V_\phi\times \{0\})$
if and only if $\phi(x)\in \partial V_\phi$,
which holds if and only if $x\in \partial M$.
If
$r\geq 1$, this shows that  $\graph(f)$
is a submanifold of $M\times N$ modeled on $E\times\{0\}$,
whose formal boundary is as asserted. If $r=0$, the preceding shows that
$\graph(f)$
is a submanifold of $M\times N$ modeled on $E\times\{0\}$,
with prescribed boundary $\{(x,f(x))\colon x\in \partial M\}$.
\end{proof}
Occasionally, results concerning ordinary manifolds do not carry over to
manifolds with rough boundary, and caution is needed; e.g., graphs need not be submanifolds
for $C^r$-maps to
$C^r$-manifolds with rough boundary.
\begin{ex}\label{bad-grph}
The convex subset $N:=[0,\infty[$ of $E:=\R$
with dense interior can be considered as a smooth manifold
with boundary, with global chart $\id_N$.
Moreover, $M:=\R$ is a smooth manifold.
The map $f\colon M\to N$, $x\mto x^2$
is smooth; its graph $\graph(f)$
is not a submanifold of the smooth manifold with boundary $M\times N$
(see Exercise~\ref{surpeasier}).\\[2mm]
The proof of Proposition~\ref{graphsubm} does not carry over because the map
\[
M\times N\to \R\times\R,\quad (x,y)\mto (x,y-f(x))
\]
(which corresponds to $g\circ(\phi\times\psi)$ in the proof of Proposition~\ref{graphsubm},
applied with $\phi:=\id_M$ and $\psi:=\id_N$) has image
$\{(x,z) \in \R \times \R \colon \, z\geq {-x^2}\}$,
which fails to be locally convex.
\end{ex}
A second type of submanifold will be used repeatedly.
\begin{defn}\label{full-subm}
Let $r\in \N_0\cup\{\infty\}$
and $M$ be a $C^r$-manifold with rough boundary
modeled on a locally convex space~$E$.
A subset $N\sub M$ is called a \emph{full submanifold}
if, for each $x\in N$, there exists a chart $\phi\colon U_\phi\to V_\phi\sub E$
of~$M$ around~$x$ such that $W_\phi:=\phi(U_\phi\cap N)$
is a locally convex subset of~$E$ with dense interior.
\end{defn}
Note that for boundaries relative $E$ we have
\[
W_\phi\cap \partial W_\phi = (W_\phi\cap \partial V_\phi)\cup
(W_\phi\setminus W_\phi^0),
\]
using the interior $W_\phi^0$ of $W_\phi$ relative $V_\phi$.
As a consequence, for $x\in W_\phi$
\begin{equation}\label{ensures-A3}
x\in \partial W_\phi\quad\Leftrightarrow\quad
\phi^{-1}(x)\in \partial M\cup (N\setminus N^0),
\end{equation}
using the formal boundary of~$M$.
Thus $N$ is a $C^r$-manifold with rough boundary modeled
on~$E$, with the ma\-xi\-mal $C^r$-atlas
containing the rough charts $\phi|_{U_\phi\cap N}\colon U_\phi\cap N\to W_\phi$ for $\phi$ as before;
(\ref{ensures-A3}) guarantees condition~(A3).
\begin{ex}
Each closed ball $\wb{B}_\ve(x)$ with respect to a given norm on $\R^n$
is a full submanifold of $\R^n$.
\end{ex}
\begin{ex}
Let $M$ be a locally compact $C^r$-manifold with rough boundary,
modeled on a locally convex space~$E$
(for example, any finite-dimensional manifold without boundary or with
$C^r$-boundary, or any finite-dimensional $C^r$-manifold with corners).
Then $E$ is locally compact and hence of finite
dimension. If $\phi\colon U_\phi\to V_\phi\sub E$
is a chart of~$M$ around some $x\in M$, then $V_\phi$ contains a convex $\phi(x)$-neighborhood $W$.
Then $\phi^{-1}(W)$ contains a compact $x$-neighborhood $K$.
We choose a norm $\|\cdot\|$ on $E$ and find $\ve>0$
such that $V_\phi\cap \wb{B}_\ve(\phi(x))\sub \phi(K)$.
Then $\phi(K)\cap\wb{B}_\ve(\phi(x))$ is compact.
Also, $C=W\cap \wb{B}_\ve(\phi(x))$ is convex and has non-empty (and hence dense)
interior relative~$E$. Thus $\phi^{-1}(C)$ is a compact, full submanifold
of~$M$ such that $x\in (\phi^{-1}(C))^0$
relative~$M$.
\end{ex}
\begin{numba}
For $r\geq 1$, the definition of a tangent bundle $TM$
works just as well if $M$ is a $C^r$-manifold
with rough boundary modeled on a locally convex space~$E$,
with tangent spaces $T_pM\cong E$ for $p\in M$.
Likewise, tangent maps can be defined with the usual properties.
In fact, Definition~\ref{defTefX} and Remark~\ref{remsimpachn}
make sense if $U\sub E$, $V\sub F$ (and $W\sub H$)
are locally convex subsets with dense interior.
Definition~\ref{taspa}
and Remark~\ref{tangonopsub}
can be copied if $M$ is a $C^r$-manifold with
rough boundary
and also Definition~\ref{tabu},
replacing charts with rough charts
and
creating a $C^{r-1}$-manifold $TM$
with rough boundary.
Note that, as $(V_\phi\times E)^0=V_\phi^0\times E$,
the maps $T\phi$ in (\ref{defTphi})
satisfy
property (A3) of a rough $C^{r-1}$-atlas
(as stated in Definition~\ref{defroughbd}) .
They are rough $E\times E$ charts
with image the locally convex subset $V_\phi \times E\sub E\times E$
with dense interior
(replacing ordinary charts onto open subsets in the proof of Lemma~\ref{chckdt}(d),
which now produces a rough $C^{r-1}$-atlas).
Likewise, in all of Example~\ref{anothtriv} through Definition~\ref{defndfinE},
manifolds can be replaced by manifolds with rough boundary, open subsets of locally convex spaces by
locally convex subsets with dense interior, and charts by rough charts.
We can also define smooth vector fields
on a $C^\infty$-manifold with smooth boundary and their Lie bracket.
In all of Definition~\ref{defVF} through Lemma~\ref{larelglob},
simply replace smooth manifolds by smooth manifolds manifolds with rough boundary,
open subsets of locally convex spaces by
locally convex subsets with dense interior, and charts by rough charts.
\end{numba}
\begin{lem}\label{sub-TM}
Let $r\in \N\cup\{\infty,\omega\}$
and $M$ be a $C^r$-manifold
which may have a rough boundary if $r\not=\omega$.
If $N\sub M$ is a submanifold $($or full submanifold,
if $r\not=\omega)$ and $\pi_{TN}\colon TN\to N$
the bundle projection,
then $TN$ is a submanifold of the $C^{r-1}$-manifold
$TM$ with boundary $\pi_{TN}^{-1}(\partial N)$
$($resp., a full submanifold of $TM)$.
Moreover, the $C^{r-1}$-manifold structure as a submanifold
coincides with the $C^{r-1}$-manifold structure on $TN$
as the tangent bundle of~$N$.
\end{lem}
\begin{prf}
See Exercise \ref{tangent-vs-sub}.
\end{prf}
\begin{defn}\label{defCklmfd}
Let
$M_1$ be a $C^k$-manifold with rough boundary,
$M_2$ be a $C^\ell$-manifold with rough boundary,
and $N$ be a $C^{k+\ell}$-manifold with rough boundary,
for $k,\ell\in \N_0\cup\{\infty\}$.
We say that a map $f\colon M_1\times M_2\to N$
is $C^{k,\ell}$ if, for each $x=(x_1,x_2)\in M_1\times M_2$
there exist charts $\phi_j\colon U_j\to V_j$ of $M_j$ around $x_j$
for $j\in \{1,2\}$ and a chart $\psi\colon U_\psi\to V_\psi$ of~$N$
such that $f(U_1\times U_2)\sub U_\psi$
and $\psi\circ f\circ (\phi_1^{-1}\times \phi_2^{-1})\colon V_1\times V_2\to V_\psi$
is $C^{k,\ell}$.
\end{defn}
\noindent
Every $C^{k,\ell}$-map $f$ is continuous, the restrictions $f|_{U_1\times U_2}$
being continuous.
\begin{rem}
It is good enough for Definition~\ref{defCklmfd}
that $M_1\times M_2$
is a topological space
and we have the atlases for $M_1$, $M_2$, and $N$ available.
If both $M_1$ and $M_2$ are $C^r$-manifolds for some $r\geq \max\{k,\ell\}$,
we can consider $M_1\times M_2$ as a $C^r$-manifold;
but this is not required for the purpose.
Of course, in most applications both $M_1$ and $M_2$ will be smooth manifolds.
\end{rem}
\begin{rem}\label{rem-simpl-Ckl}
Let $k,\ell\in \N_0\cup\{\infty\}$.
For manifolds which may have a rough boundary, we have:
\medskip

\noindent
(a) If $M_1$ is a $C^k$-manifold, $M_2$ a $C^\ell$-manifold,
$N$ a $C^{k+\ell}$-manifold,\break 
${f\colon M_1\times M_2\to N}$
a $C^{k,\ell}$-map,
$\phi_j\colon U_j\to V_j$ a chart for $M_j$
for $j\in\{1,2\}$ and $\psi\colon U_\psi\to V_\psi$ a chart for $N$
such that $f(U_1\times U_2)\sub U_\psi$,
then $\psi\circ f\circ (\phi_1^{-1}\times \phi_2^{-1})$
is $C^{k,\ell}$.
Inserting charts as in Definition~\ref{defCklmfd},
this follows from Propositions~\ref{chainR1} and \ref{chainR2}.\medskip

\noindent
Proposition~\ref{CkvsCkk} implies:\medskip

\noindent
(b) If $M_1$, $M_2$, and $N$ are $C^{k+\ell}$-manifolds
and $f\colon M_1\times M_2\to N$ is $C^{k+\ell}$,
then $f$ is $C^{k,\ell}$.\medskip

\noindent
(c) If $M_1$ and $M_2$ are $C^k$-manifolds, $N$ a $C^{k+k}$-manifold
and a mapping $f\colon M_1\times M_2\to N$ is $C^{k,k}$,
then $f$ is $C^k$.\medskip

\noindent
(d) If $M_1$, $M_2$, and $N$ are smooth manifolds, then
$f\colon M_1\times M_2\to N$ is $C^\infty$ if and only if $f$
is $C^{\infty,\infty}$.\medskip

\noindent
Proposition~\ref{Cklprod} implies, using~(a):\medskip

\noindent
(e) If $M_1$ is a $C^k$-manifold, $M_2$ a $C^\ell$-manifold and
$N=N_1\times \cdots\times N_n$ a product of $C^{k+\ell}$-manifolds,
then a map $f=(f_1,\ldots, f_n)\colon M_1\times M_2\to N$
is $C^{k,\ell}$ if and only if all of its components
$f_1,\ldots, f_n$ are $C^{k,\ell}$.
\end{rem}
\noindent
Applying Propositions~\ref{chainR1} and \ref{chainR2}
in local charts,
we get the following versions of the Chain Rule:
\begin{prop}\label{Cklchainmfd}
Let $M_1$ be a $C^k$-manifold, $M_2$ a $C^\ell$-manifold,
and $N$ a $C^{k+\ell}$-manifold with $k,\ell\in\N_0\cup\{\infty\}$,
possibly with rough boundary. Let  $f\colon M_1\times M_2\to N$ be
a $C^{k+\ell}$-map. Then the following holds:
\begin{description}[(D)]
\item[\rm(a)]
If $L_1$ is a $C^k$-manifold with rough boundary, $L_2$ a $C^\ell$-manifold with rough boundary,
$g_1\colon L_1\to M_1$ a $C^k$-map, and $g_2\colon L_2\to M_2$ a $C^\ell$-map;
then
$f\circ (g_1\times g_2)\colon L_1\times L_2\to N$ is a $C^{k,\ell}$-map.
\item[\rm(b)]
If $L$ is a $C^{k+\ell}$-manifold with rough boundary and $g\colon N\to L$
a $C^{k+\ell}$-map, then
$g\circ f \colon M_1\times M_2\to L$
is a $C^{k,\ell}$-map.
\end{description}
\end{prop}
\begin{rem}\label{subm-sec-bd}
We mention that $C^k$-manifolds
with rough boundary modeled on a set $\cE$ of locally convex spaces
can be defined in an obvious fashion for $k\in\N_0\cup\{\infty\}$,
using rough $E$-charts with $E\in\cE$
(cf.\ Definition~\ref{non-pure} for the case without boundary).
Also $C^k$-maps between such can be defined in the straightforward way;
the details are left to the reader.
Likewise, $C^{k,\ell}$-maps $M\times N\to L$
can be defined for smooth
manifolds $M$, $N$, and $L$ with rough boundary which need not be pure.
Proposition~\ref{subm-sec}(b) remains valid if~$L$
has a rough boundary.
\end{rem}
\subsection*{Smooth regularity and partitions of unity}
It can be useful to know that a given smooth manifold
admits a certain supply of smooth functions.
We briefly discuss two concepts: smooth regularity and smooth paracompactness.
The ground field is $\K=\R$ in these considerations.
\begin{prop}\label{equiv-ck-reg}
For $r\in \N_0\cup\{\infty\}$,
let $M$ be a $C^r$-manifold $($possibly with rough boundary$)$
modeled on a real locally convex space. Then the following
statements are equivalent:
\begin{description}[(D)]
\item[\rm(a)]
The topology on~$M$ is initial with respect to
$C^r(M,\R)$.
\item[\rm(b)]
For every $x\in M$ and neighborhood
$U\sub M$ of~$x$,
there exists
a $C^r$-function $f \colon M\to\R$
such that $f(x)\not=0$ and $\Supp(f)\sub U$.
\item[\rm(c)]
For every $x\in M$ and $x$-neighborhood
$U\sub M$,
there exists
a $C^r$-function $f\colon M\to \R$
such that $f(M)\sub [0,1]$,
$\Supp(f)\sub U$ and $f|_V=1$ for some $x$-neighborhood $V\sub M$.
\end{description}
\end{prop}
\begin{prf}
(a)$\Rightarrow$(c):
Every $x$-neighborhood $U\sub M$ in the initial topology contains
a neighborhood of the form $f_1^{-1}(V_1)\cap\cdots\cap f^{-1}_n(V_n)$
with $n\in\N$, $f_1,\ldots, f_n\in C^r(M,\R)$ and $f_j(x)$-neighborhoods $V_j\sub\R$ for
$j\in\{1,\ldots, n\}$. After shrinking~$V_j$, we may assume that, for some $\ve>0$, we have
$V_j=\,]f_j(x)-\ve,f_j(x)+\ve[$ for all $j\in\{1,\ldots, n\}$.
Then $g:=1+\sum_{j=1}^n(f_j-f_j(p))^2\in C^r(M,\R)$ satisfies $g(y)\geq 1$ for all $y\in M$,
$g(x)=1$,
and $g^{-1}([1,1+\ve^2[)\sub U$. If $h\colon \R\to \R$ is a smooth function such that
$h(t)=1$ for $t$ in some \break  {$1$-neighborhood,} $h(\R)\sub [0,1]$ and $\Supp(h)\sub \,]{-\infty},1+\ve^2[$,
then $f:=h\circ g$ is~$1$ on some $x$-neighborhood, $\Supp(f)\sub U$ and $f(M)\sub[0,1]$.

The implications ``(c)$\Rightarrow$(b)'' and ``(b)$\Rightarrow$(a)''
are trivial.
\end{prf}
\begin{defn}\label{defn-smoothly-reg}
If $M$ satisfies the equivalent conditions (a)--(c)
of Proposition~\ref{equiv-ck-reg},
then $M$ is called a \emph{$C^r$-regular}
manifold (or also \emph{smoothly regular}, if $r=\infty$).
Functions as in (b) are called \emph{bump functions};
functions as in (c) are called \emph{cut-off functions}.
\end{defn}
\begin{prop}\label{from-space-to-mfd}
Let $r\in \N_0\cup\{\infty\}$
and $M$ be a $C^r$-manifold $($possibly with rough boundary$)$
modeled on a real locally convex space~$E$.
If $M$ is a regular topological space and $E$ is $C^r$-regular,
then $M$ is $C^r$-regular.
Notably, every locally compact $C^r$-manifold with rough boundary
is $C^r$-regular.
\end{prop}
\begin{prf}
Let $x\in M$ and $U\sub M$ be an open $x$-neighborhood.
After shrinking~$U$, we may assume that~$U$ is the domain of a chart
$\phi\colon U\to V\sub E$ of~$M$ around~$x$.
Since $M$ is regular, there exists an $x$-neighborhood $A\sub U$
which is closed in~$M$.
Then $\phi(A)$ is an $x$-neighborhood and contains
a relatively open $x$-neighborhood $Q\sub V$.
Thus $Q=V\cap P$ for some open
$x$-neighborhood~$P$ in~$E$.
By Proposition~\ref{equiv-ck-reg}(b),
there exists $f\in C^r(E,\R)$ such that $f(x)\not=0$ and $\Supp(f)\sub P$.
Thus $\Supp(f)\cap V\sub Q\sub \phi(A)$.
We define a function $g\colon M\to\R$ piecewise via $g(y):=0$
for $y$ in the open set $M\setminus A$ and $g(y):=f(\phi(y))$
for $y\in U$. Then $g$ is $C^r$, $g(x)\not=0$ and $\Supp(g)\sub A\sub U$.
Thus~$M$ satisfies condition~(b) of Proposition~\ref{equiv-ck-reg}.
If $M$ is locally compact, then $M$ is regular (see Proposition~\ref{paracomp-reg})
and its modeling space $E$ is locally compact, whence $E$ is finite-dimensional.
Thus $E$ is $C^r$-regular (since $E\cong \R^n$ for some $n$
and it is a classical fact that $\R^n$ admits smooth cut-off functions).
\end{prf}
$C^r$-regularity is inherited by finite direct products
and by submanifolds.
\begin{defn}\label{defn-smoo-para}
Let $r\in \N_0\cup\{\infty\}$
and $M$ be a $C^r$-manifold
modeled on a real locally convex space
$($which may have a rough boundary$)$.
A \emph{$C^r$-partition of unity}
on~$M$ is a family $(h_j)_{j\in J}$
of $C^r$-functions
$h_j\colon M\to\R$
such that
\begin{description}[(D)]
\item[(a)]
$h_j(M)\sub [0,1]$ for all $j\in J$;
\item[(b)]
The family
$(\Supp(h_j))_{j\in J}$ of supports is locally finite;
and
\item[(c)]
For each $x\in M$, we have $\sum_{j\in J}h_j(x)=1$.
\end{description}
We say that $(h_j)_{j\in J}$ is \emph{subordinate}
to
an open cover $(U_i)_{i\in I}$ of~$M$
if, for each $j\in J$, there exists $i\in I$ such that $\Supp(h_j)\sub U_i$.
We say that $M$ is \emph{$C^r$-paracompact}
if, for each open cover $(U_i)_{i\in I}$ of~$M$,
there exists a $C^r$-partition of unity subordinate to it.
If $M$ is $C^r$-paracompact with $r=\infty$, we also say that~$M$
is \emph{smoothly paracompact}.
\end{defn}
By (b), each $x\in M$ has a neighborhood $W$ in~$M$
such that
\[
J_0:=\{j\in J\colon \Supp(h_j)\cap W\not=\emptyset\}
\]
is finite. Thus $\sum_{j\in J_0}h_j(y)=1$
for all $y\in W$ and all other summands in
$\sum_{j\in J}h_j(y)$ vanish.
\begin{rem}
Note that every $C^k$-paracompact $C^r$-manifold is paracompact.
In fact, for an open cover $(U_i)_{i\in I}$
and a subordinate $C^r$-partition of unity $(h_j)_{j\in J}$,
the preimages $h_j^{-1}(]0,\infty[)$
form a locally finite open cover subordinate to $(U_i)_{i\in I}$.
\end{rem}
\begin{prop}\label{findim-smoothly-para}
Every locally compact,
paracompact smooth manifold~$M$ with rough boundary
is smoothly paracompact.
\end{prop}
\begin{prf}
Let $(U_i)_{i\in I}$ be an open cover of~$M$.\\[1mm]
Step~1. If $M$ is $\sigma$-compact,
let $L_n$ and $O_n$ be as in Lemma~\ref{lem-rings}, applied with $A:=X:=M$.
For each $x\in X$, there exists $i(x)\in I$ such that $x\in U_{i(x)}$.
If $n\in\N$
and $x\in L_n$,
then Proposition~\ref{from-space-to-mfd} provides a smooth function
$g_{n,x}\in C^\infty_c(M,\R)$
such that $g_{n,x}(M)\sub [0,1]$,
$\Supp(g_{n,x})\sub U_{i(x)}\cap O_n$
and $g_{n,x}(x)=1$.
For $x\in L_n$, the set
$g_{n,x}^{-1}(]0,\infty[)$
form an open cover of~$L_n$ in~$O_n$.
By compactness of $L_n$, we find a finite subset $\Phi_n\sub L_n$
such that
\begin{equation}\label{hence-pos}
L_n\sub \bigcup_{x\in \Phi_n}
g_{n,x}^{-1} (]0,\infty[).
\end{equation}
We define $J:=\bigcup_{n\in\N}\{n\}\times \Phi_n$.
Then $(\Supp(g_{n,x}))_{(n,x)\in J}$ is a locally finite
family (cf.\ proof of Lemma~\ref{pre-lcp-para}),
entailing that
\[
g\colon M\to \R,\quad y\mto\sum_{(n,x)\in J}g_{n,x}(y)
\]
is a finite sum on some neighborhood
of each point and hence smooth.
As a consequence of~(\ref{hence-pos}),
we have $g(x)>0$ for all $x\in M$.
Then $h_{n,x}(y):=g_{n,x}(y)/g(y)\in [0,1]$
for all $y\in M$ and $h_{n,x}\colon M\to\R$
is a smooth function for all $(n,x)\in J$.
For all $y\in M$, we have $\sum_{(n,x)\in J}
h_{n,x}(y)=\frac{1}{g(y)}\sum_{(n,x)\in J}g_{n,x}(y)=
\frac{1}{g(y)}g(y)=1$.\\[1mm]
Step~2.
In general,
$M$ is the union of
of a family $(M_a)_{a\in A}$
of $\sigma$-compact open subsets $M_a\sub M$
(see Proposition~\ref{lcp-parac}).
By Step~1, $M_a$ is smoothly paracompact,
whence there exists a smooth partition of unity $(h_{a,b})_{b\in B_a}$
on~$M_a$ which is subordinate to $(U_i\cap M_a)_{i\in I}$.
Let $g_{a,b}\colon M\to\R$ the extension of $h_{a,n}$
by~$0$.
Define $J:=\bigcup_{a\in A}\{a\}\times B_a$.
Then $(g_j)_{j\in J}$
is a smooth partition of unity on~$M$ which is subordinate to $(U_i)_{i\in I}$.
\end{prf}
\begin{prop}\label{lcp-mfd-parac}
Let $r\in \N_0\cup\{\infty\}$
and $M$ be a $C^r$-manifold $($possibly with rough boundary$)$.
Then the following holds:
\begin{description}[(D)]
\item[\rm(a)]
Every connected component of~$M$
is open.
\item[\rm(b)]
If $M$ is locally compact, then $M$ is paracompact
if and only if all of its connected components are $\sigma$-compact.
\end{description}
\end{prop}
\begin{prf}
(a) For every $x\in M$, there exists a chart $\phi\colon U\to V\sub E$.
The set $V$ is locally convex; after shrinking~$V$, we may assume that~$V$
is convex and hence path-connected. Thus $U=\phi^{-1}(V)$ is path-connected
as well and hence connected. As a consequence, each connected component
$C$ of~$M$ is a neighborhood in~$M$ for each of its points.
Thus $C$ is open in~$M$.

(b) If each connected component is $\sigma$-compact, then $M$ is paracompact, by
Proposition~\ref{lcp-parac}.
If $M$ is paracompact, then $M$ the union of a family $(M_j)_{j\in J}$
of pairwise disjoint, open subsets $M_j\sub M$
which are $\sigma$-compact (see Proposition~\ref{lcp-parac}).
Each connected component $C$ of~$M$ is contained in some~$M_j$.
Since $C$ is closed in~$M$ and hence in $M_j$, we deduce that~$C$
is $\sigma$-compact.
\end{prf}
\begin{small}
\subsection*{Exercises for Section~\ref{sec-mfd-boundary}}

\begin{exer}
  \label{exer:2.1.3} Let $M$ and~$N$ be $C^r$-manifolds
with rough boundary.
Show that $\partial (M\times N)$
is the union of $(\partial M) \times N$
and $M \times (\partial N)$.
\end{exer}

\begin{exer}
  \label{exer:2.1.4} Show that each finite-dimensional
$C^r$-manifold with boundary
also is a $C^r$-manifold
with corners.
\end{exer}

\begin{exer}
  \label{exer:2.1.5} Let $M$ and $N$ be finite-dimensional
$C^r$-manifolds with boun\-dary.
Show that $M\times N$ is a $C^r$-manifold with corners.
If $\partial N=\emptyset$, show that $M\times N$ is a $C^r$-manifold with boundary.
\end{exer}

\begin{exer}
  \label{exer:2.2.4} Let $f\colon M\to N$
be a $C^r$-diffeomorphism between
$C^r$-manifolds with
rough boundary, where $r\geq 1$.
Show that $f(\partial M)=\partial N$.
\end{exer}

\begin{exer}
  \label{exer:2.3.2} Show that
$\partial (TM)=(\pi_{TM})^{-1}(\partial M)$.
\end{exer}

\begin{exer}\label{localbsub}
Let $M$ be a $C^r$-manifold with rough boundary modeled on
a locally convex space~$E$. Let $F\sub E$ be a closed vector subspace,
$N\sub M$ be a subset and $\phi\colon U_\phi\to V_\phi$ of~$M$
be a chart for~$M$ such that
\[
\phi(U_\phi\cap N)=V_\phi\cap F.
\]
Show that $(V_\phi\cap F)\cap \partial(V_\phi\cap F)\sub V_\phi\cap \partial V_\phi$,
where the boundary on the left is formed relative~$F$.\\[2mm]
[If $x\in (V_\phi\cap F)\cap \partial(V_\phi\cap F)$, then $x\in V_\phi$.
For each $x$-neighborhood $U\sub E$, the $x$-neighborhood $U\cap F$ relative~$F$
contains an element $y\in F\setminus (V_\phi\cap F)$. Then $y\not\in V_\phi$
and $y\in U$, whence $x\in \partial V_\phi$.]
\end{exer}

\begin{exer}
Let $M$ be a $C^r$-manifold with rough boundary
and $N\sub M$ be a submanifold. Using Exercise~\ref{localbsub},
show that $\partial N\sub \partial M$.
Show that $\partial N$ can be a proper subset of
$N\cap\partial M$.
\end{exer}


\begin{exer}
Let $M$ and $N$ be analytic manifolds over $\K\in\{\R,\C\}$,
modeled on locally convex spaces~$E$ and~$F$, respectively. Let $f\colon M\to N$
be an analytic map. Show that the graph of~$f$ is a submanifold of the analytic manifold $M\times N$,
modeled on $E\times\{0\}$.
\end{exer}

\begin{exer}\label{surpeasier}
Given $k\in\N_0\cup\{\infty\}$, we consider  the $C^k$-map
$f\colon \R \to [0,\infty[$, $x\mto x^2$.
Show that $G:=\mbox{graph}(f)$
is not a submanifold of the $C^k$-manifold with boundary
$H:=\R \times [0,\infty[$, as follows:
Suppose, to the contrary, that we could find
a rough $\R^2$-chart $\phi\colon U\to V\sub\R^2$ of~$H$
such that $(0,0)\in U$ and
\[
\phi(U\cap G)=V\cap F
\]
for a vector subspace $F\sub \R^2$; we shall derive a contradiction.
After a translation in the range, we may assume that $\phi(0,0)=(0,0)$.
\begin{description}[(D)]
\item[(a)]
$G$ is not a $(0,0)$-neighborhood in $\R\times [0,\infty[$;
deduce that $F\not=\R^2$. Show that $F\not=\{0\}$;
thus $F$ has dimension~$1$. After a rotation in the range,
we may assume that $F=\R\times\{0\}$.
\item[(b)]
We may assume that $U=\,]{-r},r[\,\times[0,r[$ for some $r\in\,]0,1]$,
since $U$ contains such a subset as a relatively open $(0,0)$-neighborhood.
As~$V$ is locally convex, it contains a convex $(0,0)$-neighborhood $C$
which is relatively open in~$V$.
Thus $G\cap U=\{(x,x^2)\colon x\in \,]{-r},r[\}$ is connected.
Thus $V\cap F$ is a connected subset of $F=\R\times \{0\}\sim\R$
and hence of the form $I\times \{0\}$ for some interval $I\sub\R$
which contains~$0$.
Since $(G\cap U)\setminus \{(0,0)\}$ has two connected components,
also $I\setminus\{0\}$ has two connected components.
We therefore find a real number $t>0$ such that $[{-t},t]\sub I$;
after shrinking~$t$, we may assume that $[{-t},t]\times \{0\}\sub C$.
\item[(c)]
We have $t=\phi(x,x^2)$ and $-t=\phi(y,y^2)$ for certain $x,y\in \,]{-r},r[\,\setminus\{0\}$.
Thus $(x,x^2)$ and $(y,y^2)$ are in the open subset $]{-r},r[\,\times\, ]0,r[$.
If $k\geq 1$, then $\phi$ is $C^k$ (as also $\id_H$
is a chart for $H$);
using the inverse function theorem (or invariance of domain if $k=0$),
show that~$V$ (and hence also~$C$) is a neighborhood of $({-t},0)$ and $(t,0)$ in~$\R^2$.
We therefore find $\ve>0$ such that
\[
({-t},\ve), ({-t},-\ve), (t,\ve), (t,{-\ve})\,\in\, C.
\]
Using the convexity of~$C$, deduce that $[{-t},t]\times[{-\ve},\ve]\sub C\sub V$.
\item[(d)]
Using the inverse function theorem (or invariance of domain), deduce that
$U=\phi^{-1}(V)$ is a neighborhood of $(0,0)$ in~$\R^2$.
This is a contradiction as $U=]{-r},r[\times[0,r[$;
hence~$\phi$ cannot exist.
\item[(e)]
Using more complicated arguments, for $k\in\N\cup\{\infty\}$
one can show that $\{(x,x^2)\colon x\in [0,\infty[\}$
is not a submanifold of $[0,\infty[^2$, considered as a \break
$C^k$-manifold with rough boundary
(see \cite{Gl21b}).
Show that
\[ \{(x,x^2)\colon x\in [0,\infty[\}\] is a submanifold of
$[0,\infty[^2$, considered as a $C^0$-manifold with rough boundary.
\end{description}
\end{exer}

\begin{exer}\label{tangent-vs-sub}
Prove Lemma~\ref{sub-TM}.
\end{exer}

\begin{exer}\label{exc-smoothreg-3}
Let $r=\infty$ in the situation of
Exercise~\ref{exer:2.4.1a} and assume that~$N$
is smoothly regular.
Show that $X$ and $Y$ are $f$-related if and only
if $(Y.\gamma)\circ f=X.(\gamma\circ f)$
for each $\gamma\in C^\infty(N,\R)$.
\end{exer}

\begin{exer}\label{exc-smoothreg-2}
\begin{description}[(D)]
\item[(a)]
Using facts from
analysis, show that
$\R^n$ is smoothly regular.
\item[(b)]
Let $(\cH, \langle \cdot,\cdot \rangle)$
be a real Hilbert space.
Show that $\cH$ is smoothly regular.\\[.5mm]
[Hint: The map
$q\colon \cH\to \R$, $q(x):=\|x\|^2:=\langle x,x\rangle$
is smooth.\,]
\item[(c)]
Part\,(b) can be pushed
further:
If $E$ is a real locally convex space
whose vector topology is defined
by a set $\cP$
of seminorms $p\colon E\to\R$
such that $p(x)=\sqrt{\langle x,x\rangle_p\,}$
for a positive semidefinite
symmetric bilinear form $\langle \cdot, \cdot \rangle_p\colon
E\times E\to\R$ on~$E$
($p$ is a Hilbert seminorm),
then $E$ admits smooth bump functions
because the argument from (b)
can be applied to the pre-Hilbert spaces
$E_p:=E/N_p$, and $f\circ \alpha_p\colon E\to\R$
is smooth for each smooth map $f\colon E_p\to \R$
(where $N_p$ and
$\alpha_p\colon E\to E_p$, $\alpha_p(x):=x+N_p$
are as in Appendix~\ref{chaplcx},
Section~\ref{modout}).\vspace{-.5mm}
\end{description}
\end{exer}%
\nin All ``nuclear'' locally convex spaces
have the property described
in~(c),
and thus many non-Banach spaces of relevance
for analysis (see, e.g., \cite{Pi72} and \cite{Tr67}).

\begin{exer}\label{exc-cutoff}
Let $k\in \N_0\cup\{\infty\}$ and
$M$ be a $C^k$-manifold (possibly with rough boundary)
which is $C^k$-regular. Show that, if $K\sub M$
is a compact subset and $U\sub M$ an open subset with $K\sub U$,
there exists a $C^k$-function $h\colon M\to [0,1]$
such that $\Supp(h)\sub U$ and $h|_W=1$ for some open subset $W\sub M$
with $K\sub W$.\\[1mm]
[Let $g\colon \R\to\R$ be monotonically increasing and $C^k$ with
$g(t)=0$ for all $t\leq 0$ and $g(t)=1$ for all $t\geq 1/2$.
For $z\in K$, let $h_z\in C^k(M,\R)$
wth $h_z(z)=1$ and $\Supp(h_z)\sub U$. Then $K\sub \bigcup_{z\in \Phi}h_z^{-1}(]1/2,\infty[)$
for some finite subset $\Phi\sub K$. Set $h(x):=g(\sum_{z\in \Phi}h_z(x))$.\,]
\end{exer}

\end{small}
\section{Differential forms}
\label{sec:diff-form}

Differential forms play a significant role throughout this book; 
either as differential forms on Lie groups, or as differential
forms on manifolds on which certain Lie groups act. 
In the present section, we describe a natural approach to
differential forms on manifolds modeled on locally convex spaces.
The main differences compared to the finite-dimensional case
are that there is no natural coordinate description
for differential forms in local charts; and also that
smooth partitions of unity are unavailable
for general locally convex manifolds.
Therefore, one has to be careful with localization arguments
(even in the case of Banach manifolds).

As shown above, the space 
${\cal V}(M)$ of smooth vector fields on $M$ carries a natural Lie 
algebra structure,
for each smooth manifold~$M$.
We shall see below that
each smooth $p$-form $\omega \in \Omega^p(M,E)$ with values in a locally 
convex space~$E$ gives rise to an alternating $p$-linear map 
\[
{\cal V}(M)^p \to C^\infty(M,E), \quad 
(X_1,\ldots, X_p) \mapsto \omega(X_1, \ldots, X_p).
\]
If $M$ has the property that each tangent vector 
extends to a smooth vector field, which is always the case
locally, then this leads to an inclusion of $\Omega^p(M,E)$
into the space of Lie algebra cochains for ${\cal V}(M)$
with values in the ${\cal V}(M)$-module $C^\infty(M,E)$.
We shall define the exterior derivative on differential  
forms in such a way that
it corresponds to the Lie algebra differential
(as in Definition~\ref{def:C.1}) 
with respect to the preceding identification.
The connection between differential
forms and Lie algebra cohomology
will be exploited throughout the book.
In the current section,
we use it to reduce assertions concerning
geometric objects (differential forms,
Lie derivatives and exterior derivatives)
to simple assertions concerning Lie algebra cochains,
the proofs of which are purely formal.
They can be looked up quickly in
Appendix~\ref{app:liealg-cohom}  once they are needed.\\[3mm]
{\bf Convention.}
Although we are mainly interested
in the real case, we can work over $\K\in \{\R,\C\}$
in this section.
All manifolds are smooth manifolds
modeled on locally convex $\K$-vector
spaces, with or without (rough) boundary.
%
%
\begin{defn} \label{def2.3.4a}
If $M$ is a smooth manifold and $E$ a locally convex space,
then an \emph{$E$-valued $p$-form} $\omega$ on $M$ is a function 
$\omega$ which associates to each $x \in M$ a $k$-linear 
alternating map $\omega_x \colon T_x(M)^p \to E$
such that in local charts the map 
$(x,v_1, \ldots, v_p) \mapsto \omega_x(v_1, \ldots, v_p)$
is smooth.
More explicitly, we require that the map
$\omega_\phi\colon V \times Z^p\to E$,
%
\begin{equation}\label{ominloc}
\omega_\phi(x,v_1,\ldots, v_p):=
\omega_{\phi^{-1}(x)}(T_x\phi^{-1}(v_1),\ldots,
T_x\phi^{-1}(v_p))
\end{equation}
is smooth for each chart $\phi\colon U \to V$
of~$M$, where $Z$ is the modeling space
of~$M$.
We write $\Omega^p(M,E)$ for the space of
$E$-valued $p$-forms on $M$ and identify
$\Omega^0(M,E)$ with the space $C^\infty(M,E)$ 
of smooth $E$-valued functions on $M$.
\end{defn}
%
%
\begin{rem}\label{firstencount}
We shall mostly deal with $p$-forms for $p\in \{0,1,2\}$.
As already mentioned, $E$-valued
$0$-forms simply are smooth $E$-valued
functions on~$M$.
Typical examples
of 1-forms arise
from the differentials
of smooth $E$-valued functions $f\colon M\to E$.
For any such function, $\omega_x:=df|_{T_xM}$
defines a smooth $E$-valued $1$-form~$\omega$,
where
$df\colon TM\to E$
is as in Definition~\ref{defndfinE}
(we shall come back to this
example in Remark~\ref{combck} below).
Typical examples of differential forms
of second and higher order arise as wedge
products of lower order differential
forms, which we now define.
\end{rem}
%
%
\begin{defn}\label{def2.3.4b}
Let $E_1, E_2, E_3$ be locally convex spaces and 
\[ \beta \: E_1 \times E_2 \to E_3 \] 
be a continuous bilinear map. 
Given $p,q\in \N_0$,
let $S_{p+q}$ be the symmetric group
of all permutations of $\{1,\ldots, p+q\}$.
We define the wedge product
$\omega \wedge \eta \in \Omega^{p+q}(M,E_3)$
of $\omega \in \Omega^p(M,E_1)$ and $\eta\in \Omega^q(M,E_2)$
via
$(\omega \wedge \eta)_x := \omega_x \wedge \eta_x$,
where  
\begin{eqnarray*}
&&(\omega_x \wedge \eta_x)(v_1,\ldots, v_{p+q}) \\
&&:= \frac{1}{p!q!} \sum_{\sigma \in S_{p+q}} 
\sgn(\sigma) 
\beta\big(\omega_x(v_{\sigma(1)}, \ldots, v_{\sigma(p)}), 
\eta_x(v_{\sigma(p+1)}, \ldots, v_{\sigma(p+q)})\big)
\end{eqnarray*}
is alternating in $v_1,\ldots, v_{p+q}\in T_xM$
(see Exercise~\ref{exc-omeg-1}).
\end{defn}
Then
\[
\Omega^p(M,E_1) \times \Omega^q(M,E_2) 
\to \Omega^{p+q}(M,E_3),
\quad (\omega, \eta) \mapsto \omega \wedge \eta
\]
is a bilinear map.
For $p = q = 1$, we simply have
\[
(\omega \wedge \eta)_x(v_1, v_2) 
=
\beta(\omega_x(v_1), \eta_x(v_2)) - \beta(\omega_x(v_2), \eta_x(v_1)).
\]
\begin{rem}
Wedge products are used
in particular in
the following
important cases:
\begin{description}[(D)]
\item[(a)]
$\beta \: \K \times E \to E$ is the scalar multiplication of $E$. 
\item[(b)]
$\beta \: \cA \times \cA \to \cA$ is the multiplication map of
an associative locally convex topological algebra. 
\item[(c)]
$\beta \: \g \times \g \to \g$ is the Lie bracket of a
locally convex topological Lie algebra. 
In this case we also write $[\omega,\eta] := \omega \wedge \eta$. 
\end{description}
\end{rem}
\begin{rem}
We can interpret a $\K$-valued
$0$-form as a smooth function
$f\colon M\to \K$.
If $\beta\colon \K\times E\to E$
is the scalar multiplication,
then the wedge product of~$f$ and a
$p$-form $\omega\in \Omega^p(M,E)$
is given by $(f\wedge \omega)_x=
f(x)\omega_x$, using the scalar multiplication
in $\Alt^p(T_xM,E)$.
We simply write $f\omega:=f\wedge \omega$
in this case and note that $\Omega^p(M,E)$
becomes an $C^\infty(M,\K)$-module
in this way.
If furthermore $E=\K$,
then
the wedge product of
two $\K$-valued $0$-forms
(interpreted as $\K$-valued functions)
is given by pointwise
scalar multiplication,
$(f_1\wedge f_2)(x)=f_1(x)f_2(x)$,
and thus coincides with the product
in the associative algebra $C^\infty(M,\K)$.
\end{rem}
We mention
that various
consistent definitions of the wedge product
of differential forms are possible
(as already in the finite-dimensional
case). We here follow Lang's conventions,
which are particularly natural
(see \cite[Remark preceding Defn.\,V.3.3]{La99}
for further discussions of this point).
%
%
\begin{defn}\label{def2.3.4c}
The pullback $\phi^*\omega$ of $\omega \in \Omega^p(M,E)$
with respect to a smooth map $\phi \: N \to M$ is the smooth
$p$-form in $\Omega^p(N,E)$ defined 
by $(\phi^*\omega)_x := (T_x\phi)^*\omega_{\phi(x)}$, i.e., 
\[
(\phi^*\omega)_x(v_1, \ldots, v_p) 
:= \omega_{\phi(x)}(T_x\phi(v_1), \ldots, T_x\phi(v_p)).
\]
\end{defn}
Note that the Chain Rule implies that 
%
\begin{eqnarray}\label{eq:pb1}
\id_M^*\omega = \omega \quad \hbox{ and } \quad 
\phi_1^*(\phi_2^*\omega) = (\phi_2 \circ \phi_1)^*\omega 
\end{eqnarray}
holds for compositions of smooth maps. Moreover, 
%
\begin{eqnarray}\label{eq:pb2}
\phi^*(\omega \wedge \eta) = \phi^*\omega \wedge \phi^*\eta 
\end{eqnarray}
follows directly from the definitions.
For $f = \omega \in \Omega^0(M,E)$, we simply have
$\phi^*f = f \circ \phi$.\\[3mm]
The definition of the exterior differential 
\[
d \: \Omega^p(M,E) \to \Omega^{p+1}(M,E)
\]
is slightly more subtle than in finite dimensions.
We shall use the following
notation:
Given $\omega\in \Omega^p(M,E)$
and vector fields $X_1,\ldots, X_p\in \cV(U)$
on an open subset $U\sub M$,
we define a smooth map $\omega(X_1,\ldots, X_p)\colon U\to E$
via
\[
\omega(X_1,\ldots, X_p)(x)\,:=\,
\omega_x(X_1(x),\ldots, X_p(x)).
\]
%
%
\begin{prop}\label{propextd}
For $\omega \in \Omega^p(M,E)$, $x \in M$ and
$v_0,\ldots, v_p \in T_x(M)$, we choose
smooth vector fields $X_i$
defined on a neighborhood of~$x$ such that
$X_i(x) = v_i$.
Then
\begin{eqnarray} \label{eq:3.3.0} 
\hspace*{-14mm}\lefteqn{(d\omega)_x(v_0, \ldots, v_p) 
\,:= \,\sum_{i = 0}^p (-1)^i
\big(X_i.\omega(X_0, \ldots, \hat X_i, \ldots, X_p)\big)(x)}\qquad \notag \\ 
&&\ \ \ \ \;\; + \sum_{i < j} (-1)^{i+j}
\omega([X_i, X_j], X_0, \ldots, \hat X_i,\ldots, \hat X_j, \ldots, X_p)(x)
\label{eq:3.3.1} 
\end{eqnarray}
does not depend on the choice of the vector fields~$X_i$
and defines a $(p+1)$-form $d\omega \in \Omega^{p+1}(M,E)$.   
\end{prop}
%
%
\begin{rem}\label{remextd}
As usual, the hat $\wh{\;}$ indicates
that the corresponding
entry has to be omitted.
The definition of the differential is designed in such a way that 
for $X_0, \ldots, X_p \in {\cal V}(M)$, we have 
\begin{eqnarray*}
(d\omega)(X_0, \ldots, X_p)
&= & 
\sum_{i = 0}^p (-1)^i X_i.\omega(X_0, \ldots, \hat X_i, \ldots,X_p) \cr
& &\ \  + \sum_{i < j} (-1)^{i+j} \omega([X_i, X_j], X_0, \ldots, \hat X_i,
\ldots, \hat X_j, \ldots, X_p). 
\end{eqnarray*}
\end{rem}
\nin {\em Proof of Proposition}~\ref{propextd}.
We have to verify that the right hand side of
(\ref{eq:3.3.1}) does not depend on
the choice of the vector fields $X_k$ and that it is alternating in
the $v_k$. First we show that $d\omega$ does not depend on the choice
of the vector fields $X_k$, which amounts to showing that if one
vector field $X_k$ vanishes in $x$, then the right hand side of
(\ref{eq:3.3.1}) vanishes.\footnote{In fact,
if also $Y_1,\ldots, Y_p$
are vector fields with $Y_k(x)=v_k$,
we only need to show
that the difference
of the right hand side of (\ref{eq:3.3.1})
using $X_1,\ldots, X_k,Y_{k+1},\ldots, Y_p$
and the same expression for
$X_1,\ldots, X_{k-1}, Y_k,\ldots, Y_p$
vanishes, for $k=1,\ldots,p$.}

Suppose that $X_k(x) =0$. Then the only terms not obviously
vanishing~in $x$ are
\begin{eqnarray}
&& \sum_{i \not=k}^p (-1)^i \big(X_i.\omega(X_0, \ldots, \hat X_i, \ldots,
X_p)\big)(x), \label{eq:3.3.2} \\
&&  \sum_{i < k} (-1)^{i+k} \omega([X_i, X_k], X_0, \ldots, \hat X_i,
\ldots, \hat X_k, \ldots, X_p)(x), \;\;\mbox{and} \quad \label{eq:3.3.3} \\
&& \sum_{k < i} (-1)^{i+k} \omega([X_k, X_i], X_0, \ldots, \hat X_k,
\ldots, \hat X_i, \ldots, X_p)(x). \label{eq:3.3.4} 
\end{eqnarray} 
In a local chart, we have 
\begin{eqnarray*}
\hspace*{-6mm}\lefteqn{\big(X_i.\omega(X_0, \ldots, \hat X_i, \ldots,X_p)\big)(x)}\qquad \\
&=& (d_1\omega)(x,X_i(x))
(X_1(x), \ldots, \hat X_i(x), \ldots, X_p(x))\\[.5mm]
&&
+ \sum_{j < i} \omega_x(X_0(x), \ldots, d X_j(x)X_i(x), \ldots, \hat
X_i(x),\ldots, X_p(x))\cr
&&
+ \sum_{j > i} \omega_x(X_0(x), \ldots, \hat X_i(x),\ldots, dX_j(x)X_i(x),
\ldots, X_p(x)).
\end{eqnarray*}
Here, we interpret $\omega$
as a smooth function of $p+1$ variables
(as in~\ref{ominloc});
the direction of differentiation
in the partial differential $d_1\omega$
has been written as the second argument,
for better readability (cf.\ also Exercise~\ref{exc-omeg-3}).
For a fixed $i > k$, the assumption $X_k(x) = 0$ implies that
only the term 
\[
\omega_x(X_0(x), \ldots, d X_k(x)X_i(x), \ldots, \hat X_i(x),\ldots, X_p(x))
\]
contributes. In view of $X_k(x) = 0$, we have 
\[
dX_k(x) X_i(x) = d X_k(x) X_i(x) - dX_i(x)X_k(x) = [X_i,X_k](x).
\]
This leads to 
\begin{eqnarray*}
\hspace*{-12mm}\lefteqn{\omega_x(X_0(x), \ldots, dX_k(x)X_i(x), \ldots,
\hat X_i(x),\ldots, X_p(x))}\qquad \\
&=& 
-(-1)^k \omega([X_k, X_i], X_0, \ldots, \hat X_k,
\ldots, \hat X_i, \ldots, X_p)(x), 
\end{eqnarray*}
so that corresponding terms in (\ref{eq:3.3.2}) and (\ref{eq:3.3.4}) 
cancel, and the same happens for $i < k$ with the terms in 
(\ref{eq:3.3.2}) and (\ref{eq:3.3.3}).  
This proves that $d\omega$ is well defined by (\ref{eq:3.3.0}). 

To see that we obtain a smooth $(p+1)$-form, we use a local chart
and choose the vector fields $X_i$ as constant vector fields. Then 
%
\begin{eqnarray}
\label{eq:ext-diff}
(d\omega)_x(v_0,\ldots, v_p) 
= \sum_{i=0}^p (-1)^{i} (d_1\omega)(x,v_i)(v_0, \ldots, \hat v_i, \ldots,
v_p) 
\end{eqnarray}
is a smooth function of $(x,v_0, \ldots, v_p)$. 

It remains to show that $(d\omega)_x$ is alternating. 
To this end, assume
that $v_i = v_j$ for some $i < j$.
Since each map $(d_1\omega)(x,v)$ is alternating,
(\ref{eq:ext-diff}) implies that
\begin{eqnarray*}
&& \hspace*{-4mm}(d\omega)_x(v_0, v_1, \ldots, v_p) \\[.5mm] 
&=&(-1)^i (d_1\omega)(x,v_i)(v_0, \ldots, \hat v_i, \ldots,v_p)\\
&&\qquad\qquad + (-1)^j (d_1\omega)(x,v_j)(v_0, \ldots, \hat v_j, \ldots,v_p)\cr
&=& (-1)^i (d_1\omega)(x,v_i)(v_0, \ldots, \hat v_i, \ldots,v_p)\\
&&\qquad\qquad- (-1)^{i} (d_1\omega)(x,v_i)(v_0, \ldots, \hat v_i, \ldots,v_p)\\
&=&0.
\end{eqnarray*}
Thus $(d\omega)_x$ is alternating,
and the proof is complete.\qed
%
%
\begin{prop}\label{propDD}
For each $\omega \in \Omega^p(M,E)$, we have $d^2 \omega :=d(d\omega)= 0$. 
\end{prop}
\begin{prf} It clearly suffices to verify this for the case where 
$M$ has a global chart, as $M$ is covered by open chart domains.
Each $p$-form $\omega \in \Omega^p(M,E)$ defines an
alternating $p$-linear map
%
\begin{equation}\label{omegg}
\omega_\g \colon  {\cal V}(M)^p \to C^\infty(M,E)\,\quad
(X_1,\ldots, X_p)\mto \omega(X_1,\ldots, X_p)\,.
\end{equation}
Thus $\omega_\g$ is a $p$-cochain for the 
Lie algebra $\g := {\cal V}(M)$ with values in the
${\cal V}(M)$-module 
$C^\infty(M,E)$,
where the module structure is the natural one given by 
$X.f := df \circ X$. The
map $\omega \mapsto \omega_\g$ is injective, as we see by 
evaluating $p$-forms on vector fields which are constant in the global chart.
By definition of $d$,
we have $d_\g \omega_\g = (d\omega)_\g$ (see Remark~\ref{remextd}).
Now $(d^2\omega)_\g = d_\g^2\omega_\g = 0$  
by Lemma~\ref{lem:C.2}(4) in Appendix~\ref{app:liealg-cohom},
entailing that $d^2\omega = 0$.
\end{prf}
\begin{rem}
Another way to verify that $d^2\omega = 0$ is to calculate directly  
in local charts using formula (\ref{eq:ext-diff}). 
Then $d^2\omega = 0$ easily follows from the symmetry of second derivatives 
of $\omega$ (Theorem~\ref{schwarz}); see Exercise~\ref{exer:2.6.4}.
\end{rem}
\begin{defn}
Extending $d$ to a linear self-map of the space
\[
\Omega(M,E) := \bigoplus_{p \in \N_0} \Omega^p(M,E)
\]
of all $E$-valued differential forms on $M$, 
the relation $d^2 = 0$ implies that the space
\[
Z^p_{\rm dR}(M,E) := \ker(d\res_{\Omega^p(M,E)})
\]
of \emph{closed forms} contains the space 
$B^p_{\rm dR}(M,E) := d(\Omega^{p-1}(M,E))$
of \emph{exact forms} (where $B^0_{\rm dR}(M,E):=\{0\}$),
so that the \emph{$E$-valued de Rham cohomology space} 
\[
H^p_{\rm dR}(M,E) := Z^p_{\rm dR}(M,E) / B^p_{\rm dR}(M,E)
\]
can be defined.
\end{defn}
%
%
%
\begin{rem}\label{combck}
Considering a smooth function $f \: M \to E$ as a differential
form of degree $0$,
we can form its exterior differential
$df\in \Omega^1(M,E)$.
\begin{description}[(D)]
\item[(a)]
Then the exterior differential
$df$ is the $1$-form given by $(df)_x=df|_{TxM}$
(as already encountered in Remark~\ref{firstencount}),
where $df\colon TM\to E$ on the right hand side
is the differential of $f$,
as introduced in
Definition~\ref{defndfinE}.\footnote{The multiple
use of the symbol ``$df$''
will not create confusion.}
\item[(b)]
Since $M$ is locally convex, vanishing of $df$ 
means that the function $f$ is locally constant. 
Thus $Z^0_{\rm dR}(M,E) \cong H^0_{\rm dR}(M,E)$ is the space of locally 
constant functions on $M$. If $M$ has $d$ connected components, then 
$H^0_{\rm dR}(M,E) \cong E^d$. 
\end{description}
\end{rem}
%
%
\begin{lem} \label{lem-pullback}
If $\phi \: N \to M$ is a smooth map and $\omega \in \Omega^p(M,E)$, then 
$d(\phi^*\omega) = \phi^*d\omega.$
\end{lem}
\begin{prf} First we assume that $\phi$ is a diffeomorphism. 
Let $X_0, \ldots, X_p \in {\cal V}(N)$ and 
define $Y_0, \ldots, Y_p \in {\cal V}(M)$ by 
$Y_i(\phi(x)) := T_x\phi(X_i(x))$, so that 
$Y_i \circ \phi = T\phi \circ X_i$.
Then
$[Y_i, Y_j] \circ \phi = T\phi \circ [X_i, X_j]$ for 
all $i,j \in \{0,\ldots, p\}$,
by Lemma~\ref{larelglob}.
Moreover, we have 
\[
\phi^*(\omega(Y_0, \ldots, \hat Y_i, \ldots, Y_p))
\,= \, (\phi^*\omega)(X_0, \ldots, \hat X_i, \ldots, X_p). 
\]
We further have for each smooth function $f$ on $M$ the relation 
\[
\phi^*(Y_i.f)(x) 
= df_{\phi(x)}(Y_i(\phi(x))) 
= df_{\phi(x)}(T_x(\phi)X_i(x)) = (X_i.(\phi^*f))(x)
\]
(cf.\ also Exercise~\ref{exer:2.4.1a}(a)), 
so that we obtain with (\ref{eq:3.3.0}) 
\[
\phi^*(d\omega)(X_0, \ldots, X_p) 
\, = \, d(\phi^*\omega)(X_0, \ldots, X_p).
\]
Since this relation also holds on each open subset of $M$, resp., $N$, 
we conclude that $d(\phi^*\omega) = \phi^*(d\omega)$. 
The preceding argument applies in particular to local diffeomorphisms 
defined by charts. 

To complete the proof of the general case, we may now 
assume w.l.o.g.\ that $M$ and $N$ are open subsets of 
locally convex spaces (or locally
convex subsets with dense interior,
in the case of manifolds with rough boundary).
By (\ref{eq:ext-diff}),
we then have 
\[
(d\omega)_y(v_0,\ldots, v_p) 
= \sum_{i = 0}^p (-1)^i
(d_1\omega(y,v_i))(v_0, \ldots, \hat v_i, \ldots, v_p)
\]
for all $y\in M$, and therefore, for $x\in N$,
\begin{eqnarray*}
\hspace*{-3mm}\lefteqn{(\phi^*(d\omega))_x(v_0,\ldots, v_p)} \quad \\
&=& \sum_{i = 0}^p (-1)^i (d_1\omega(\phi(x),d\phi(x)v_i))
(d\phi(x)v_0, \ldots, d\phi(x)\hat v_i, \ldots, d\phi(x)v_p). 
\end{eqnarray*}
On the other hand, the Chain Rule leads to 
\begin{eqnarray*}
&& \hspace*{-3mm}(d(\phi^*\omega))_x
(v_0, \ldots, v_p) \\[.5mm]
&=& \sum_{i = 0}^p (-1)^i
(d_1\omega(\phi(x),d\phi(x)v_i))(d\phi(x)v_0, \ldots, 
d\phi(x)\hat v_i, \ldots, d\phi(x)v_p)\cr
&& \!+ \sum_{i = 0}^p (-1)^i \sum_{j < i} \omega_{\phi(x)}
(d\phi(x)v_0, \ldots, d^{(2)}\phi(x,v_j,v_i),\ldots, 
d\phi(x)\hat v_i, \ldots, d\phi(x)v_p)\cr
&& \! + \sum_{i = 0}^p (-1)^i \sum_{j > i} \omega_{\phi(x)} 
(d\phi(x)v_0,\ldots, d\phi(x)\hat v_i,
\ldots, d^{(2)}\phi(x,v_j,v_i), \ldots, d\phi(x)v_p)\cr
&=& \sum_{i = 0}^p (-1)^i
\big(d_1\omega(\phi(x),d\phi(x)v_i)\big)(d\phi(x)v_0, \ldots, 
d\phi(x)\hat v_i, \ldots, d\phi(x)v_p),
\end{eqnarray*}
where the terms in the penultimate line and the one preceding it cancel because of the
symmetry of the bilinear maps $d^{(2)}\phi(x,\cdot)$
(Proposition~\ref{schwarz}).
This proves the assertion. 
\end{prf}
For any smooth manifold $M$ and locally convex space~$E$, 
we now define a natural representation of 
the Lie algebra ${\cal V}(M)$ on the space $\Omega^p(M,E)$ of 
$E$-valued $p$-forms on $M$, given by the 
Lie derivative.
\begin{defn}
For $Y \in {\cal V}(M)$ and
$\omega\in \Omega^p(M,E)$,
we define
the \emph{Lie derivative} 
${\cal L}_Y \omega\in \Omega^p(M,E)$
via
\begin{eqnarray}
&&({\cal L}_Y\omega)_x(v_1,\ldots, v_p) \notag \\
&&= (Y.\omega(X_1, \ldots, X_p))(x) - \sum_{j = 1}^p \omega(X_1, \ldots,
[Y, X_j], \ldots, X_p)(x) \label{good-summands} \label{goodsummands} \\
&&= Y.\omega(X_1, \ldots, X_p)(x) 
+ \sum_{j = 1}^p (-1)^j \omega([Y, X_j], X_1, \ldots, \hat X_j, 
\ldots, X_p)(x), \notag
\end{eqnarray}
where $X_1, \ldots, X_p$ are smooth vector fields on an open
neighborhood of $x$ such that $X_i(x) = v_i$.
\end{defn}
To see that the right hand side does not depend on the choice of the 
vector fields $X_i$, suppose that 
$X_i(x) = 0$ for some $i$.
In a local chart,
evaluation of the right hand side in $x$ yields
\begin{eqnarray*}
\hspace*{-6mm}\lefteqn{(Y.\omega(X_1, \ldots, X_p))(x)
- \omega(X_1, \ldots,[Y, X_i], \ldots, X_p)(x)}\quad \\
&=& \omega(x)(X_1(x),\ldots, dX_i(x)Y(x),\ldots, X_p(x)) \\
&& \; - \omega(X_1(x), \ldots,dX_i(x)Y(x) - dY(x)X_i(x), \ldots, X_p(x))\\
&=&  0. 
\end{eqnarray*}
Therefore $({\cal L}_Y\omega)_x(v_1,\ldots, v_p)$ is well
defined, and it is now easy to see that $\cL_Y\omega$
is an $E$-valued $p$-form on~$M$.
In fact, we may assume that~$M$ has a global chart and $X_1,\ldots, X_p$
are constant in the chart; then $(\cL_Y\omega)_x(v_1,\ldots,v_p)$
is given by the first term in~(\ref{good-summands}),
hence alternating in $v_1,\ldots, v_p$.
It is clear that
the map $\cL_Y\colon \Omega^p(M,E)\to\Omega^p(M,E)$,
$\omega\mto \cL_Y\omega$
is linear, for each $p\in P$.
We also write $\cL_Y$
for the unique linear map
$\Omega(M,E)\to\Omega(M,E)$
which restricts to the self-map
$\cL_Y$ of $\Omega^p(M,E)$
for each~$p$.
It is clear that the map
$Y\mto \cL_Y$ from $\cV(M)$
to the space of linear self-maps
of~$\Omega(M,E)$ is linear.\\[3mm]
It is useful to introduce a further standard notation.
\begin{defn}
Given
$X \in {\cal V}(M)$ and $p \geq 1$,
we define $i_X\omega\in \Omega^{p-1}(M,E)$
for $\omega\in \Omega^p(M,E)$ via
\[
(i_X\omega)_x(v_1,\ldots, v_{p-1}) :=
\omega_x(X(x), v_1,\ldots, v_{p-1}). 
\]
Then
$i_X \colon \Omega^p(M,E) \to \Omega^{p-1}(M,E)$,
$\omega\mto  i_X\omega$
is a linear map.
For $\omega \in \Omega^0(M,E) = C^\infty(M,E)$,
we put $i_X \omega := 0$.
We also consider $i_X$ as a linear self-map
of $\Omega (M,E)$.
\end{defn}
%
%
%
\begin{prop} \label{prop-cart-for}
For $X, Y \in {\cal V}(M)$, we have on $\Omega(M,E)$: 
\begin{description}[(D)]
\item[\rm(a)] $[{\cal L}_X, {\cal L}_Y] = {\cal L}_{[X,Y]}$, i.e., 
the Lie derivative defines a representation of the Lie algebra
${\cal V}(M)$ 
on $\Omega(M,E)$, and also on each $\Omega^p(M,E)$.
\item[\rm(b)] $[{\cal L}_X, i_Y] = i_{[X,Y]}$. 
\item[\rm(c)] ${\cal L}_X = d \circ i_X + i_X \circ d$ $($Cartan's formula$)$.
\item[\rm(d)] ${\cal L}_X \circ d = d \circ {\cal L}_X$. 
\item[\rm(e)] ${\cal L}_X(Z^p_{\rm dR}(M,E)) \subeq B^p_{\rm dR}(M,E)$. 
\end{description}
\end{prop}
\begin{prf}
(a)--(d):
It suffices to verify these formulas on an open cover of~$M$, so that we may
assume that $M$ has a global chart.
Then everything follows from the corresponding formulas 
in Lemma~\ref{lem:C.2} from Appendix~\ref{app:liealg-cohom}, 
applied to the Lie algebra\footnote{We give $\K$, $\g$ and $V$ the discrete topology
here, as their locally convex topologies need not make $C^\infty(M,E)$
a topological $\cV(M)$-module, see
Exercise~\ref{exc-not-top-module}.}
$\g := {\cal V}(M)$ and the $\g$-module $V:=C^\infty(M,E)$.
To see this, recall from the proof of 
Proposition~\ref{propDD}
that the map
$\phi\colon \Omega(M,E)\to C_c^\bullet(\g,V)$, $\omega\mto \omega_\g$
is injective and intertwines $d$ and the Lie algebra differential
$d_\g$, i.e., $(d\omega)_\g=d_\g(\omega_\g)$.
Likewise, $\phi$ intertwines
$\cL$ and the corresponding Lie derivative
on $C_c^\bullet(\g,V)$,
and it also intertwines $i_X$ and the corresponding insertion
operator on $C_c^\bullet(\g,V)$.
Hence, after applying $\phi$ to all identities,
we can verify them in $C_c^\bullet(\g,V)$.

(e) follows from~(c).
\end{prf}
%
%
\begin{rem}\label{liedban}
We defined the Lie derivative
$\cL_Y\omega$ of a vector-valued
$p$-form $\omega\in \Omega^p(M,E)$
along a vector field~$Y$
in a way which works for
manifolds~$M$
modeled
on arbitrary locally convex spaces.
In the case of Banach manifolds
(without boundary), another
description of $\cL_Y\omega$
is possible, which provides a better
intuition for this object
(but will not be used in the present book).
In this case, a flow $\Fl^Y$ is associated to~$Y$ (see Definition~\ref{flow-auton}
and Corollary~\ref{niceflow-ban-m}).
Locally around each given point,
$\Fl_t^Y$ defines a local diffeomorphism
for small $t$, which can be used
to move the differential
form around. Using a local chart,
one easily verifies that
%
\begin{equation}\label{interprtLd}
\frac{d}{dt}\Big|_{t=0}((\Fl_t^Y)^*\omega)_x(v_1,\ldots, v_p)\;=\;(\cL_Y\omega)_x(v_1,\ldots, v_p)
\end{equation}
(cf.\ \cite[Ch.\,V,
Prop.\ 5.1 and 5.2]{La99}
for the case of scalar-valued forms).
The formula remains valid if~$M$ is modeled
on a locally convex space and $Y$ is assumed
to admit local $C^1$-flows (see Exercise~\ref{exc-flow-ext-der}).
\end{rem}
\begin{rem}
Of course, we can just as well define
real analytic differential forms
on real analytic manifolds
modeled on real locally
convex spaces, by requiring that the maps
(\ref{ominloc}) be real analytic.
All of the preceding constructions,
results and proofs carry over directly
to this case, if we replace
the word ``smooth'' everywhere
by ``real analytic'' and the symbol $C^\infty$
by $C^\omega$.
\end{rem}
\noindent
In the remainder of this section,
we have a closer look at differential
forms on finite-dimensional
smooth manifolds.
Some classical facts will be
used without proof, but precise
references will be given.\\[2.5mm]
To get a more explicit description
of differential forms on subsets
of~$\R^n$,
it is convenient to recall classical
notation.
\begin{defn}
Let $U$ be an open subset of $\R^n$
(or a locally convex subset $U\sub \R^n$
with dense interior).
To introduce standard notation, let
$\lambda_j\colon U\to \R$, $(x_1,\ldots, x_n)\mto x_j$
be the $j$-th coordinate projection
for $j\in \{1,\ldots, n\}$.
One writes
$dx_j:=d\lambda_j$
for the exterior differential of~$\lambda_j$;
thus $dx_j$ is a $1$-form in $\Omega^1(U,\R)$,
and
$(dx_j)_x(v_1,\ldots, v_n)=v_j$.
For scalar-valued differential forms
$\omega_1,\ldots, \omega_k$ on $U$,
let us define
$\omega_1\wedge\cdots\wedge \omega_k$
recursively via
\[
\omega_1\wedge\cdots\wedge \omega_k:=
(\omega_1\wedge \cdots\wedge \omega_{k-1})\wedge \omega_k
\]
(the choice of brackets actually is irrelevant
because the wedge product of scalar-valued
forms is associative, see Exercise~\ref{exer:2.6.6}). 
Given $\omega\in \Omega^p(U,\R)$
and a smooth vector-valued function
$f\in C^\infty(M,E)$, we define
$f\omega\in \Omega^p(M,E)$ via $(f\omega)_x(v_1,\ldots, v_p)
=\omega_x(v_1,\ldots, v_p)f(x)$
(which is $f\wedge \omega$,
formed using the scalar multiplication
$\R\times E\to E$).
\end{defn}
It is clear that $\Omega^p(U,\R)=\{0\}$
for $p>n$.
For $p\leq n$,
the forms \break $dx_{i_1}\wedge\cdots\wedge dx_{i_k}$ 
can be used to describe arbitrary $p$-forms.
%
%
\begin{lem}\label{relclass}
Let $U\sub \R^n$ be an open subset
$($or a locally convex subset with dense interior$)$,
$E$ be a locally convex space, and $p\in \{1,\ldots, n\}$.
\begin{description}[(D)]
\item[\rm (a)]
Then each $E$-valued $p$-form $\omega\in \Omega^p(U,E)$
can be written as
%
\begin{equation}\label{shpform}
\omega\;=\; \sum_{i_1<\cdots<i_p}f_{i_1,\ldots, i_p}\,
dx_{i_1}\wedge\cdots\wedge dx_{i_p}\,,
\end{equation}
with unique smooth functions
$f_{i_1,\ldots, i_p}\colon U\to E$
for $1\leq i_1<\cdots<i_p\leq n$.
\item[\rm (b)]
In particular, every $E$-valued
$n$-form $\omega\in \Omega^p(U,E)$
can be written as
\[
\omega\,=\; f\, dx_1\wedge\cdots\wedge dx_n\,,
\]
where $f\colon U\to E$ is a uniquely determined
smooth function.
\end{description}
\end{lem}
\begin{prf}
The proof uses the following fact:
If $1\leq i_1<\cdots<i_p\leq n$
and $1\leq j_1<\cdots< j_p\leq n$, then
\[
(dx_{i_1}\wedge\cdots\wedge dx_{i_p})\Big(
\frac{\partial}{\partial x_{j_1}},
\ldots, \frac{\partial}{\partial x_{j_p}}\Big)
\;=
\left\{
\begin{array}{cl}
1 &\;\mbox{if $\,(i_1,\ldots, i_p)=(j_1,\ldots, j_p)$}\\
0 & \;\mbox{otherwise.}
\end{array}
\right.
\]
This is well known and can be quickly checked
by induction (e.g., using Exercise~\ref{exer:2.6.5} with $q=1$).
We prove (b);
the proof of (a) (which we shall not use)
is similar and left as an exercise.
We show uniqueness of $f$ first.
If $\omega=f\, dx_1\wedge\cdots\wedge dx_n$,
then
\[
\omega_x
({\textstyle
\frac{\partial}{\partial x_1}(x), \ldots, \frac{\partial}{\partial x_p}(x)})
=f(x)\,
(dx_1\wedge\cdots\wedge dx_n)\big(
{\textstyle \frac{\partial}{\partial x_1}(x), \ldots,
\frac{\partial}{\partial x_p}(x)}\big)
=f(x)\,,
\]
whence $f(x)$ is determined by~$\omega$.
To see that $f$ exists, we simply
take the preceding formula as the definition
of $f(x)$.
Then $f=\omega(\frac{\partial}{\partial x_1},
\ldots, \frac{\partial}{\partial x_p})$ is smooth,
and satisfies our needs.
\end{prf}
\begin{defn}
Recall that a $C^\infty$-atlas~$\cA$ of a
finite-dimensional smooth manifold
(with or without smooth boundary)
is called an \emph{oriented} atlas
if $\det (\phi\circ\psi^{-1})'(x)>0$
for all $\phi,\psi\in\cA$ and all $x$
in the domain of the transition map $\phi\circ\psi^{-1}$.
An oriented atlas~$\cA$ which is maximal
among oriented atlases is called
an \emph{orientation}, and $(M,\cA)$
is then called an \emph{oriented manifold}.
Connected finite-dimensional manifolds
admitting an orientation compatible
with the given atlas are called
\emph{orientable}.
It is easy to
see that an orientable manifold
admits exactly two orientations.
(See \cite[Ch.\,XVI, \S\,3]{La99}
and \cite[Ch.\,0, \S\,4]{dCa92}
for further information).
\end{defn}
Integration of differential forms $\omega \in \Omega^p(M,E)$ 
only makes sense if $M$ is an oriented
manifold
(possibly with smooth boundary)
of dimension~$p$ 
and $E$ is Mackey complete.
We need the Mackey completeness to ensure 
that each smooth function $f \: Q \to E$ on a cube $Q := \prod_{i=1}^p 
[a_i, b_i] \subeq \R^p$ has an iterated integral
\[
\int_Q f(x)\; dx \, =\,  \int_{a_1}^{b_1} \cdots 
\int_{a_p}^{b_p} f(x_1, \ldots, x_p)\; dx_1\cdots dx_p
\]
(see Proposition~\ref{intsquare} and its proof).
In the remainder of this section,
we only consider
$\sigma$-compact finite-dimensional
smooth manifolds, for simplicity.
As is well known,
$\sigma$-compactness ensures
the existence of smooth partitions of unity
(cf.\ \cite[Ch.\,II, Cor.\,3.4]{La99} for manifolds
without boundary),
which is useful for us.
\begin{defn}
Let $M$ be an oriented $p$-dimensional smooth manifold
with or without smooth boundary
and $\omega\in \Omega^p(M,E)$
be a compactly supported $E$-valued $p$-form,
i.e., $\supp(\omega):=\wb{\{x\in M\colon \omega_x\not=0\}}$
is compact.
If there exists a chart
$\phi \colon  U \to \phi(U)\sub \R^p$ of $M$ compatible with the
orientation
such that $\supp(\omega)\sub U$,
then we define 
\[
\int_M \omega := \int_{\phi(U)} (\phi^{-1})^*\omega 
= \int_{\phi(U)} f\ dx,
\]
where $f \in C^\infty(\phi(U),E)$ is the compactly supported
function determined by 
\[
((\phi^{-1})^*\omega)(x) = f(x)\ dx_1 \wedge \ldots \wedge dx_p\,.
\]
The existence of the integral is guaranteed
by Corollary~\ref{intcompsup}.
If $\supp(\omega)$ is not necessarily
contained in a chart domain,
we choose
a smooth partition of unity
$(\chi_i)_{i \in I}$
such that each
$\supp(\chi_i)$ is contained in a chart domain,
and define 
%
\begin{equation}\label{defeqifo}
\int_M \omega  \; := \,  \sum_{i \in I} \int_M \chi_i \; \omega\,.
\end{equation}
\end{defn}
Note that the right hand side of (\ref{defeqifo})
is a finite sum and that
each integral in the sum has already been defined,
because $\supp(\chi_i \, \omega)$
is contained in a chart domain.
Using the transformation formula for $p$-dimensional integrals, 
it is easy to see that the integral $\int_M \omega$
is well defined -- it does
not depend on the choice of
charts, nor the partition of unity.\\[3mm]
We shall use the following version
of Stokes' Theorem:
\begin{thm}
Let $M$ be an oriented,
$\sigma$-compact
$p$-dimensional smooth manifold
$($possibly with smooth boundary$)$,
$E$ be a Mackey complete locally convex space,
and
$\omega$ be a compactly supported
$E$-valued $(p-1)$-form on~$M$. Then
%
\begin{equation}\label{eqstoks}
\int_M d\omega  \; = \; \int_{\partial M} \omega\,,
\end{equation}
where the boundary $\partial M$ is given
the induced orientation.
\end{thm}
\noindent
More precisely,
the integrand on the right
is $j^*\omega$,
where $j\colon \partial M\to M$ is the inclusion map.
\begin{prf}
We can verify the validity of (\ref{eqstoks})
by applying continuous linear functionals.
Thus (\ref{eqstoks}) holds if and only if
$\int_M \lambda_*(d\omega)=\int_{\partial M}\lambda_*\omega$
for each $\lambda\in E'$,
where $(\lambda_*\omega)_x:=\lambda\circ \omega_x$.
But $\lambda_*(d\omega)=d(\lambda_*\omega)$
(see Exercise~\ref{exer:2.6.1}),
and
$\int_M d(\lambda_*\omega)=\int_{\partial M}\lambda_*\omega$
holds by
Stokes' Theorem for scalar-valued forms (see
\cite[Ch.\,XVII, Thm.\,2.1]{La99}).
\end{prf}
\begin{rem}
The preceding considerations
do not really depend on $\sigma$-compactness.
In fact,
each compact subset $K$ of
a finite-dimensional manifold~$M$
is contained in a $\sigma$-compact
open subset $N\sub M$,
as is easily seen.
Then $N$ admits smooth partitions of unity,
and this suffices for our purposes.\vspace{2mm}
\end{rem}

\subsection*{The Poincar\'e Lemma} 

As a consequence of Lemma~\ref{lem-pullback}, pullbacks of
closed, resp., exact forms are closed, resp., exact. 
We conclude that each smooth
map $f \: M \to N$ leads to well-defined linear maps (functoriality
of de Rham cohomology)
\[
f^* \: H^k_{\mathrm{dR}}(N,E) \to H^k_{\mathrm{dR}}(M,E), 
\quad [\alpha] \mapsto
[f^*\alpha].\]

\begin{lem}[Homotopy Lemma]
\index{theorem!homotopy!lemma} \label{lem:homot}
Let $M$ and $N$ be smooth manifolds, $F \: [0,1]\times M \to N$ be a
smooth map and $E$ be a Mackey complete locally convex 
space. Then the induced smooth maps $F_0, F_1 \: M \to N$ satisfy
\[
F_0^* = F_1^* \: H^k_{\mathrm{dR}}(N,E) \to H^k_{\mathrm{dR}}(M,E) 
\quad \mbox{ for } \quad k \in \N_0.\]
\end{lem}

\begin{prf} Put $J:= [0,1]$.
We write $i_t\colon M\to \{t\}\times M\subeq J\times M$ for the
canonical embeddings.
For $\omega\in\Omega^{k}(J\times M,E)$ with $k \in\N$, we define
the fiber integral $I(\omega) \in\Omega^{k-1}(M,E)$ by
\[
I(\omega)_x(v_1,\ldots, v_{k-1}) := \int_0^1 \omega_{(t,x)}
\bigl({\textstyle\frac{\partial}{\partial t}}, v_1,\ldots, v_{k-1} \bigr)\, dt,
\]
i.e.,
\[
I(\omega)_x = \int_0^1 i_{\frac{\partial}{\partial t}}\omega_{(t,x)}\, dt
\]
in the space $\cL^{k-1}(T_xM,E)$ of continuous $E$-valued
$(k-1)$-linear maps,
endowed with the topology of pointwise convergence.
Let $X_0, \ldots, X_{k-1} \in\cV(M)$ and extend these vector fields
in the canonical fashion to vector fields $\tilde{X}_i$ on $J\times M$, constant zero
in the first component. Then we have $[\tilde{X}_i,
\tilde{X}_j] = [X_i, X_j]\;\tilde{}$, and from Cartan's formula 
(Proposition~\ref{prop-cart-for}) we further get
\begin{align*}
& \bigl(d_M I(\omega)\bigr) (X_0,\ldots, X_{k-1})(x)
\\
& = \int_0^1 \! \bigl(d_{J\times M} i_{\frac{\partial}{\partial
t}}\omega \bigr) (\tilde{X}_0, \ldots,
\tilde{X}_{k-1})(t,x)\, dt
\\
& = \int_0^1 \! \bigl(\cL_{\frac{\partial}{\partial t}}\omega\bigr)
(\tilde{X}_0, \ldots, \tilde{X}_{k-1})(t,x)\, dt - \!
\int_0^1 \! i_{\frac{\partial}{\partial
t}}(d_{J\times M}\omega) (\tilde{X}_0, \ldots, \tilde{X}_{k-1})(t,x)\, dt
\\
& = (i_1^*\omega) (X_0, \ldots, X_{k-1})(x) - \bigl(i_0^*\omega
\bigr) (X_0, \ldots, X_{k-1})(x)\\
& \hspace*{40.6mm}
- I(d_{J\times M}\omega) (X_0,\ldots,
X_{k-1})(x).
\end{align*}
This means that we have the homotopy formula
\begin{equation}
\label{eq:homtop-form}
d_M I(\omega) + I(d_{J\times M}\omega) = i_1^*
\omega- i_0^*\omega.
\end{equation}
We apply this formula to $\omega= F^*\alpha$ for a closed form
$\alpha$ of degree $k\geq1$ on~$N$ and obtain
\[
\bigl[F_1^*\alpha- F_0^*\alpha\bigr] = \bigl[i_1^*\omega- i_0^*\omega\bigr] =
\bigl[I(d_{J\times M}\omega)\bigr] = \bigl[I\bigl(F^*d\alpha\bigr)\bigr] =
0.
\]
For degree $k = 0$, the space $H^0_{\mathrm{dR}}(N,E)$ consists of locally
constant functions~$f$, and since $F_t^*f = f \circ F_t$ does not
depend on $t$, we also get $F_0^* = F_1^*$ for $k =0$.
\end{prf}

\begin{cor}\label{cor:prod-poin}
If $M$ is a smooth manifold and $E,F$ locally convex spaces, 
where $E$ is assumed Mackey complete, then the projection
$p_M \: M \times F \to M$ defines an isomorphism
\[
p_M^* \: H^k_{\mathrm{dR}}(M,E) \to H^k_{\mathrm{dR}}\bigl(M \times F,E\bigr),
\quad [\alpha] \mapsto\bigl[p_M^*\alpha\bigr]
\]
for each $k \in\N_0$.
\end{cor}

\begin{prf} We consider the smooth map
\[
F \: M \times F \times [0,1] \to M \times F, \quad (m,x,t) \mapsto(m,tx).
\]
Then $F_1 = \id_{M \times F}$, and
$F_0$ is the projection onto $M \times\{0\}$.

The inclusion $i_0 \: M \to M \times F, m \mapsto(m,0)$
satisfies $p_M \circ i_0 = \id_M$ and $i_0 \circ p_M = F_0$. In view
of Lemma~\ref{lem:homot}, we have $F_0^*=F_1^*=\id$, so that the
pull-back maps $p_M^*$ and $i_0^*$ induce mutually inverse
isomorphisms between the spaces $H^k_{\mathrm{dR}}(M,E)$ 
and $H^k_{\mathrm{dR}}(M \times F,E)$.
\end{prf}

\begin{defn}\label{smoothcontractibledef}
A smooth manifold $M$ is called {\it smoothly
contractible}\index{contractible!smoothly}\index{smooth!contraction}
to a point $p\in M$ if there exists a smooth map $F\colon M\times
[0,1] \to M$ with
\[
F(x,0)=x \quad\mbox{and} \quad F(x,1)=p\;\;\,\mbox{for all $\, x\in M$.}
\]
\end{defn}

\begin{thm}[Poincar\'{e} Lemma]\label{poincarelemma}\index
{theorem!Poincar\'{e}
Lemma}\index{Poincar\'{e}!lemma}\index{Poincar\'{e}, Henri
(1854--1912)} Suppose that $M$ is a smoothly contractible manifold.
Then any closed form of degree at least one with values in a Mackey 
complete space is exact.
\end{thm}

\begin{prf} Let $F\colon M\times [0,1]\to M$ be a smooth contraction of $M$
to $p\in M$ and $\alpha\in\Omega^k(M,E)$, $k > 0$. Since $F_0 =
\id_M$ and $F_1 = p$ is the constant map, the Homotopy Lemma implies
that
$[\alpha] = \bigl[F_0^*\alpha\bigr] = \bigl[F_1^*\alpha\bigr]$. 
But since $F_1$ is constant and $k > 0$, we have $F_1^*\alpha= 0$.
\end{prf}

\begin{cor} \label{cor:poincare-lemma} 
Suppose that $M$ is open and star-shaped 
in the locally convex space~$F$. Then any closed
form of degree at least one on $M$ with values in a Mackey complete 
space is exact.
\end{cor}

\begin{prf} \hspace*{-.8mm}Contract\hspace*{-.2mm} linearly\hspace*{-.2mm} to\hspace*{-.2mm} any
\hspace*{-.2mm}point\hspace*{-.2mm} with\hspace*{-.2mm} respect\hspace*{-.2mm} to\hspace*{-.2mm} which\hspace*{-.2mm} $M$\hspace*{-.2mm}
is\hspace*{-.2mm} star-shaped.\hspace*{-1.9mm}
\end{prf}
\begin{small}
\subsection*{Exercises for Section~\ref{sec:diff-form}} 
\label{exer:3.3.1a}

\begin{exer}\label{exc-omeg-1}
Show that $(\omega\wedge\eta)_x$ is alternating
in the situation of Definition~\ref{def2.3.4b}.
\end{exer}
\begin{exer}
  \label{exer:2.6.2} 
Let $(\fg,[.,.])$ be a locally convex topological
Lie algebra and $\omega\in \Omega^1(M,\fg)$.
Show that
$[\omega,\omega]:=\omega\wedge\omega$ is given by
$[\omega,\omega]_x(v_1,v_2)=2[\omega_x(v_1),\omega_x(v_2)]$.
\end{exer}
\begin{exer}\label{exc-omeg-3}
Let $\omega\in \Omega^p(M,E)$,
$\phi\colon U\to V$ be a chart for~$M$
and $X_0^\phi,\ldots, X_p^\phi\in C^\infty(V,Z)$
be the local representatives of
$X_0,\ldots, X_p\in\cV(M)$, where~$M$
is modeled on~$Z$. Let $x\in M$.
Using Exercise~\ref{exer:2.4.1a},
show that
$(d\omega)_x(X_0(x),\ldots,X_p(x))$
equals
\begin{eqnarray*}
\lefteqn{\sum_{i=0}^pX_i^\phi . (\omega_\phi\circ (\id_V,X_0^\phi,\ldots, \wh{X_i^\phi},\ldots,
X_p^\phi))}\qquad\qquad\\
& & \;+
\sum_{i<j}(-1)^{i+j}\omega_\phi\circ
(\id_V,[X_i^\phi,X_j^\phi],X_0^\phi,\ldots\wh{X_i^\phi},\ldots,
\wh{X_j^\phi},\ldots, X_p^\phi)
\end{eqnarray*}
evaluated at $\phi(x)$,
with notation as in (\ref{ominloc}).
\end{exer}
\begin{exer}
  \label{exer:2.6.1} 
Let $\alpha\colon E\to F$ be a continuous
linear map between locally convex spaces,
$M$ be a smooth manifold, and $p\in \N_0$.
Check that $(\alpha_*\omega)_x:=\alpha\circ \omega_x$
defines an $F$-valued $p$-form
$\alpha_*\omega$ for each $\omega\in \Omega^p(M,E)$.
Verify that $\alpha_*(d\omega)=d(\alpha_*\omega)$.
\end{exer}
\begin{exer}
  \label{exer:2.6.4} 
Verify that $d^2\omega = 0$ for the exterior differential 
on $\Omega^p(M,E)$ ($M$ a smooth manifold modeled on $Z$, 
$E$ a locally convex space) 
directly in local charts, using formula (\ref{eq:ext-diff}).
Hint: For each $x \in M$, the bilinear map 
\[
Z^2 \to \Alt^p(Z, E), \quad (v,w)
\mapsto d_1^{\,(2)}\omega_x(v,w)
\]
(arising from the
second derivative
with respect to the first argument~$x$
of $\omega$, as in Proposition~\ref{diffpar}) is symmetric by
Proposition~\ref{schwarz}.
\end{exer}
\begin{exer}\label{exc-omeg-6}
Let $\phi\colon M\to N$ be a smooth map between
smooth manifolds, $X\in \cV(M)$ and $Y\in \cV(N)$ be $\phi$-related
vector fields, $E$ be a locally convex space, $p\in\N_0$ and $\omega\in \Omega^p(N,E)$.
Show that
\[
\phi^*(i_Y\omega)=i_X(\phi^*\omega)\quad\mbox{and}\quad
\phi^*(\cL_Y\omega)=\cL_X(\phi^*\omega).
\]
\end{exer}
\begin{exer}
  \label{exer:2.6.5}
In the situation of
Definition~\ref{def2.3.4b},
let us write $S_{p,q}$
for the set of all permutations $\sigma\in S_{p+q}$
such that $\sigma(1)<\cdots <\sigma(p)$
and $\sigma(p+1)<\cdots<\sigma(p+q)$.
Show that $(\omega\wedge \eta)_x(v_1,\ldots, v_{p+q})$ can be written as
%
\begin{equation}\label{ridfactorial}
\sum_{\sigma\in S_{p,q}} \sgn(\sigma)\,
\beta\big(\omega_x(v_{\sigma(1)},\ldots,v_{\sigma(p)}),
\eta_x(v_{\sigma(p+1)},\ldots, v_{\sigma(p+q)})\big)
\end{equation}
and note that the elements $\sigma\in S_{p,q}$
are in one-to-one correspondence
with pairs $(A,B)$ of disjoint subsets of $\{1,\ldots, p+q\}$
such that $A$ has $p$ elements,
$B$ has $q$ elements, and $A\cup B=\{1,\ldots, p+q\}$
(the mapping taking
$\sigma\in S_{p.q}$ to
$\big(\{\sigma(1),\ldots,\sigma(p)\},
\{\sigma(p+1),\ldots,\sigma(p+q)\}\big)$
is a suitable bijection).
We can therefore think of
the sum in (\ref{ridfactorial})
as a sum over all such partitions $(A,B)$.
\end{exer}
\begin{exer}
  \label{exer:2.6.6} 
Let $\cA$ be a locally convex,
associative topological algebra,
and $M$ be a smooth manifold.
Show that $(\omega \wedge \eta)\wedge\zeta=
\omega\wedge (\eta\wedge \zeta)$
holds for  all $p,q,r\in \N_0$ and $\omega\in \Omega^p(M,\cA)$,
$\eta\in \Omega^q(M,\cA)$, $\zeta\in \Omega^r(M,\cA)$,
if the algebra
multiplication $\beta\colon \cA\times \cA\to \cA$,
$\beta(a,b)= a b$
is used to define the wedge products.
Deduce that the
wedge product makes $\Omega(M,\cA)$
an associative algebra.
Hint: It is convenient
to use the interpretation of wedge products
provided by the previous exercise.
Show that both
$(\omega_x\wedge\eta_x)\wedge \zeta_x$
and $\omega_x\wedge (\eta_x\wedge \zeta_x)$
take\linebreak
$(v_1,\ldots, v_{p+q+r})$
to
\[ \sum_{\sigma\in S_{p,q,r}}
\!\sgn(\sigma)\,
\omega_x(v_{\sigma(1)},\ldots,v_{\sigma(p)})
\eta_x(v_{\sigma(p+1)},\ldots,v_{\sigma(p+q)})
\zeta_x(v_{\sigma(p+q+1)},\ldots,v_{p+q+r}),\]
using notation analogous to the one in
the preceding exercise for partitions into three sets.
\end{exer}
\begin{exer}
  \label{exer:2.6.3} 
Let $U\sub \R^n$ be an open subset
(or a locally convex subset with dense
interior), and $E$ be a locally
convex space.
\begin{description}[(D)]
\item[(a)]
Show that if $f\colon U\to E$ is a smooth map, then
$df=\frac{\partial f}{\partial x_1}dx_1+\cdots+
\frac{\partial f}{\partial x_n}dx_n$.
\item[(b)]
Verify the fact mentioned in the
proof of Lemma~\ref{relclass}
and prove part\,(a) of the lemma.
\end{description}
\end{exer}
\begin{exer}
  \label{exer:2.6.7} 
Let $M$ be a smooth manifold (possibly with rough
boundary), $E$ be a locally
convex space, and $U\sub \R^n$
be a locally convex subset with dense interior.
\begin{description}[(D)]
\item[(a)]
Let $f\colon M\to \R$ be a smooth function.
Using formula
(\ref{eq:ext-diff}),
show that $d(f\omega)=df\wedge \omega+ fd\omega$
for each $E$-valued
$1$-form $\omega$ on~$M$,
and more generally for each
$\omega\in \Omega^p(M,E)$.
\item[(b)]
Show that
$d(\sum_{j=1}^n f_j\, dx_j)=
\sum_{j\not=i} \frac{\partial f_j}{\partial x_i}\, dx_i\wedge
dx_j$
for each smooth $1$-form
$\omega=\sum_{j=1}^n f_j\, dx_j$ on~$U$,
using part (a) and Exercise~\ref{exer:2.6.3}(b). 
\item[(c)]
We now use a continuous bilinear map
$\beta\colon E_1\times E_2\to E_3$
between locally convex spaces
to define wedge products.
Prove that
\[
d(\omega\wedge \eta)\; =\; d\omega\wedge \eta\, +\, 
(-1)^p \omega\wedge d\eta\,,
\]
for all $p,q\in \N_0$ and
$\omega\in \Omega^p(M,E_1)$,
$\eta\in \Omega^q(M,E_2)$.
\end{description}
\end{exer}
\begin{exer}
  \label{exer:2.6.8} 
Check that the integral in (\ref{defeqifo})
is well defined.
\end{exer}
\begin{exer}
  \label{exer:2.6.9} 
Verify formula (\ref{interprtLd})
from Remark~\ref{liedban}.
\end{exer}
%
%
%
%
\begin{exer}\label{exer:2.x.x} 
Let $M$ be a smooth manifold, $E$ a locally convex space 
and $f \in C^\infty(M,E)$. Suppose that $df(m) = 0$ for some 
$m \in M$. Show that there exists a symmetric continuous bilinear map
\[
d^2 f(m) \: T_m(M)\times T_m(M)  \to E
\] 
with the property that for $X,Y \in \cV(M)$, we have 
\[ 
(X.(Y. f)(m) = d^2 f(m)(X(m), Y(m)).
\] 
\end{exer}
\begin{exer}\label{exc-not-top-module}
Let $M$ be an open subset of a locally convex space~$Z$.
Show that the map
\[
\cV(M)\times C^\infty(M,\R)\to C^\infty(M,\R),\quad (X,f)\mto X.f
\]
is continuous if and only if $Z$ has finite dimension.\\[2.3mm]
[If $Z$ has infinite dimension, use the continuous linear maps
$Z'\to C^\infty(M,\R)$, \break  $\lambda\mto\lambda|_M$
and $Z\to C^\infty(M,Z)\cong \cV(M)$, $v\mto (x\mto v)$
to create a link to Exercise~\ref{exc-eval-co}].
\end{exer}
\begin{exer}\label{exc-flow-ext-der}
Let $M$ be a smooth manifold modeled on a locally convex space.
Let $E$ be a locally convex space, $p\in\N_0$, $\omega\in\Omega^p(M,E)$
and $Y\in\cV(M)$ such that~$Y$ admits local $C^1$-flows.
For $x\in M$ and $v_1,\ldots, v_p\in T_xM$,
show that $(\cL_Y\omega)_x(v_1,\ldots, v_p)$
is given by (\ref{interprtLd}). [Use Exercises~\ref{exc-flow-in-chart}
and \ref{exc-ders-flow}.]
\end{exer}
\end{small}
%
%

%
\section{Vector bundles and principal bundles}\label{tools-mfd}
This section and the next
prepare Chapter~\ref{chap-manifold-constructions}, which is devoted to
the construction of major classes of examples of infinite-dimensional manifolds.
We define vector bundles and principal bundles
modeled on locally convex spaces.
Some basic facts are described, focussing on what is needed
for basic constructions of infinite-dimensional manifolds (like manifolds of mappings and
direct limits of ascending sequences of finite-dimensional smooth manifolds).
\subsection*{Vector bundles}
We consider vector spaces and manifolds
over~$\K\in \{\R,\C\}$.
\begin{defn}\label{defn-vbdl}
Let $F$ be a locally convex space, $r\in \N_0\cup\{\infty,\omega\}$
and~$M$ be a $C^r$-manifold.
A \emph{$C^r$-vector bundle over~$M$}
with typical fibre~$F$
is a $C^r$-manifold~$E$
together with a surjective $C^r$-map $\pi\colon E\to M$
and a vector space structure on $E_x:=\pi^{-1}(\{x\})$
for each $x\in M$ with the following local triviality property:
For each $x\in M$, there exists an open neighborhood $U\sub M$
and a $C^r$-diffeomorphism $\theta=(\theta_1,\theta_2)\colon E|_U\to U\times F$
with domain $E|_U:=\pi^{-1}(U)$ such that
\[
\theta_1=\pi|_{E|_U}\;\,\mbox{and}\;\, \mbox{$\theta_2|_{E_x}\colon E_x\to F$ is an isomorphism
for each $x\in U$.}
\]
We call $\theta$ a \emph{local trivialization} of~$E$ (and a \emph{trivialization} of $E|_U$).
If all of~$E$ admits a trivialization, $E$ is called \emph{trivializable}.
\end{defn}
\begin{rem}\label{firstremvbdl}
(a) Since $\theta$ is a homeomorphism, we deduce that $E_x$ is a locally convex topological
vector space in the topology induced by~$E$ and $\theta_2|_{E_x}$
an isomorphism of topological vector spaces.

(b) For~$M$ and $F$ as before, $E:=M\times F$ becomes a $C^r$-vector bundle
over~$M$ if we let $\pi:=\pr_1$ be the projection onto the first component
and endow $E_x=\{x\}\times F$ with the vector space structure making
the bijection $\{x\}\times F\to F$, $(x,y)\mto y$ an isomorphism of vector spaces.
Such vector bundles are called \emph{trivial} vector bundles. Notably, $E=M\times F$
is trivializable with $\id_E$ as a trivialization.

(c) A map $f \colon E\to H$ between $C^r$-vector bundles $\pi_E\colon
E\to M$ and $\pi_H\colon H\to N$
on bases $M$ and~$N$, respectively, is called a \emph{vector bundle map}
if $f(E_x)\sub E_{g(x)}$ for some $g(x)\in N$, for each $x\in M$.
Then $g\colon M\to N$ is a $C^r$-map (as $g=\pi_H\circ f\circ 0_E$ for the 
$0$-section $0_E \: M \to E, x \mapsto 0_x$),
and we also say that $f$ is a \emph{vector bundle map over~$g$}.
If $M=N$ and $f$ is a bijective vector bundle map over~$\id_M$
such that also $f^{-1}$ is $C^r$, then $f$ is called an \emph{isomorphism
of vector bundles} over~$\id_M$. For example,
a map $\theta\colon E|_U\to U\times F$
is a local trivialization if and only if it is a vector bundle isomorphism
over $\id_U$,
and thus vector bundles are locally isomorphic to trivial vector bundles.
\end{rem}
\begin{ex}\label{extanbun}
If $M$ is a $C^r$-manifold modeled on~$F$ with $r\in\N\cup\{\infty,\omega\}$,
then its tangent bundle~$TM$ is a $C^{r-1}$-vector bundle over~$M$
together with the bundle projection $\pi_{TM}\colon TM\to M$ taking all
$v\in T_xM$ to $x\in M$. In fact,
\[
\theta_\phi:=(\pi_{TM}|_{TU},d\phi)    \colon T(U)\to U\times F
\]
is a local trivialization for each chart $\phi\colon U\to V\sub F$ of~$M$,
with inverse $U\times F\to TU$, $(x,y)\mto T\phi^{-1}(\phi(x),y)$.
\end{ex}
\begin{rem}\label{rem-vbdl-rough}
If $r\not=\omega$, we may allow that $M$ and $E$ are $C^r$-manifolds
with rough boundary in Definition~\ref{defn-vbdl}.
Remark~\ref{firstremvbdl}
and Example~\ref{extanbun} carry over to this situation.
\end{rem}
We briefly discuss restrictions and pullbacks of vector bundles.
\begin{lem}\label{restrict-vbdl}
Let $F$ be a locally convex space, $r\in\N_0\cup\{\infty,\omega\}$
and $E$ be a $C^r$-vector bundle with typical fiber~$F$
over a $C^r$-manifold $N$
modeled on a locally convex space~$Z$,
with bundle projection $\pi\colon E\to N$.
Let $S\sub N$ be a submanifold modeled on a closed vector subspace
$Y\sub Z$. Then $E|_S:=\pi^{-1}(S)$
is a submanifold of~$E$
and a $C^r$-vector bundle over~$S$ with fibre~$F$
if we retain the vector space structure on $E_x$ for $x\in S$,
and restrict $\pi$ to a bundle projection $E|_S\to S$.
If $r\not=\omega$, we may assume that $N$
has a rough boundary and find that $E|_S$
is a submanfold with boundary $\pi^{-1}(\partial S)$;
likewise if $S$ is a full submanifold of~$N$.
\end{lem}
\begin{prf}
For $x\in S$, there exist an open $x$-neighborhood $U\sub N$
and a local trivialization $\theta\colon E|_U\to U\times F$.
Let $\theta_2\colon E|_U\to F$ be its second component.
After shrinking~$U$, we may assume that there exists
a chart $\phi\colon U\to V\sub Z$ of~$M$ adapted to~$S$.
Then $\psi:=(\phi\times \id_F)\circ \theta\colon E|_U\to V\times F$
is a chart for~$E$ and $\psi(E|_S\cap E|_U)=(V\cap Y)\times F=:W$
a locally convex subset with dense interior of $Y\times F$
(resp., an open subset, if $r=\omega$).
The boundary of~$W$ relative $Y\times F$
is $(\partial(V\cap Y))\times F$. For $V\in E|_S\cap E|_U$,
this entails that $\psi(v)\in \partial W$ if and only if $\pi(v)\in \partial S$,
in the case $r\not=\omega$. Thus $E|_S$ is a submanifold of~$S$
with boundary as asserted.
By the preceding, $\theta$ restricts to a $C^r$-diffeomorphism
$E|_S\cap E|_U\to (S\cap U)\times F$. The latter is a local trivialization, whence
$E|_S$ is a vector bundle.
\end{prf}
\begin{lem}\label{pullbackbdl}
Let $F$ be a locally convex space, $r\in\N_0\cup\{\infty,\omega\}$
and $E$ be a $C^r$-vector bundle with typical fiber~$F$
over a $C^r$-manifold $N$
modeled on a locally convex space~$Z$,
with bundle projection $\pi\colon E\to N$.
Let $M$ be a $C^r$-manifold
modeled on a locally convex space~$Y$, and $f\colon M\to N$ be a $C^r$-map.
If $r\not=\omega$, both $M$ and $N$ may have a rough boundary.
For $x\in M$, endow $\{x\}\times E_{f(x)}$ with
the vector space structure making $\{x\}\times E_{f(x)}\to E_{f(x)}$,
$(x,y)\mto y$
an isomorphism. Then
\[
f^*(E):=\bigcup_{x\in M}\{x\}\times E_{f(x)}
\]
admits a unique $C^r$-vector bundle structure over~$M$ with typical fibre~$F$
modeled on $Y\times F$
such that $\pi_{f^*(E)}\colon f^*(E)\to M$, $(x,v)\mto x$
for $x\in M$, $v\in E_{f(x)}$ is the bundle projection,
$\partial f^*(E)=(\pi_{f^*(E)})^{-1}(\partial M)$
and
\[
\theta_f\colon f^*(E)|_{f^{-1}(U)}\to f^{-1}(U)\times F,\;\,
(x,v)\mto (x,\theta_2(v))
\]
a local trivialization of $f^*(E)$\hspace*{-.3mm} for any local trivialization
$\theta \colon \! E|_U\! \to U\!\times \! F$ of~$E$.
\end{lem}
\begin{prf}
For convenience, let us introduce an index and write
$\{\theta_i\colon i\in I\}$ for the set of all local trivializations
$\theta_i\colon E|_{U_i}\to U_i\times F$ of~$E$. Let
$\theta_{i,2}\colon E|_{U_i}\to F$ be the second component of $\theta_i$.
For $i,j\in I$, we have
\[
(\theta_i\circ \theta_j^{-1})(x,y)=(x,g_{ij}(x,y))\quad\mbox{for all $\,(x,y)\in (U_i\cap U_j)\times F$}
\]
for a smooth $C^r$-function $g_{ij}\colon U_i\cap U_j\times F\to F$
which is linear in~$y$.
We give $f^*(E)$ the
final topology with respect to the mappings
\[
h_i\colon f^{-1}(U_i)\times F\to f^*(E),\;\;
(x,y)\mto (x,\theta_i^{-1}(f(x),y)).
\]
Note that $h_i$ is injective and has image $\pi_{f^*(E)}^{-1}(f^{-1}(U_i))=:f^*(E)|_{f^{-1}(U_i)}$;
the inverse function sends $(x,v)$ to $(x,\theta_{i,2}(v))$.
Thus
$h_i^{-1}(\im(f_i)\cap\im(f_j))=f^{-1}(U_i\cap U_j)\times F$
is open in the domain of~$h_i$ and the map
$h_i^{-1}\circ h_j$,
\[
(x,y)\mto (x, \theta_{i,2}(\theta_j^{-1}(f(x),y)))=
(x,g_{ij}(f(x),y))
\]
is a $C^r$-diffeomorphism of $f^{-1}(U_i\cap U_j)\times F$
onto itself with inverse $h_j^{-1}\circ h_i$.
By Exercise~\ref{topforbun}, each $h_i$ has open image and is a homeomorphism
onto its image. In particular,
distinct points in the image of $h_i$ can be separated by open neighborhoods.
Since $ \pi_{f^*(E)}\circ h_i\colon f^{-1}(U_i)\times F\to f^{-1}(U_i)$,
$(x,y)\mto x$ is continuous for each $i\in I$, the map $\pi_{f^*(E)}$
is continuous, whence all $v,w \in f^*(E)$
in different fibres of $\pi_{f^*(E)}$ can be separated. Thus $f^*(E)$
is Hausdorff. Since $h_i^{-1}\circ h_j$ is a $C^r$-diffeomorphism,
$\im(h_i)$ and $\im(h_j)$ induce the same $C^r$-manifold
structure on their intersection which is open in both
(with the same boundary, if $r\not=\omega$).
Using Exercise~\ref{exc-mfd-local}
and an analogous argument for manifolds
with rough boundary,
we see that $f^*(E)$ is a $C^r$-manifold
with the asserted boundary. Now
$\Theta_i:=h_i^{-1}$ is a $C^r$-diffeomorphism
and we readily see it is a local trivialization,
whence $f^*(E)$ is a $C^r$-vector bundle over~$M$
with\linebreak
typical fibre~$F$. The uniqueness assertion
is clear, as local trivializations determine
the $C^r$-manifold structure on $f^*(E)$
(cf.\ Exercise~\ref{exc-mfd-local}).
\end{prf}
\begin{rem}\label{restr-as-pb}
Of course, $E|_S\cong j^*(E)$ for the
inclusion map $j\colon S\to N$ in the situation
of Lemma~\ref{restrict-vbdl}.
\end{rem}
\noindent
For use in Section~\ref{dl-construction},
we show that all
$C^\infty$-vector bundles over a smoothly
contractible, $\sigma$-compact finite-dimensional
base manifold are trivial.\\[1.3mm]
%
%
Recall that a smooth manifold $M$ is called \emph{smoothly contractible}
if there exists a homotopy
$H\colon [0,1]\times M\to M$
from $\id_M$ to the constant function
$c_{x_0}\colon M\to M$, $x\mto x_0$
for some $x_0\in M$, such that~$H$ is a smooth map.
\begin{rem}
If a $\sigma$-compact, finite-dimensional smooth manifold~$M$
is contractible in the sense
that there exists a homotopy
$H\colon [0,1]\times M\to M$
from $\id_M$ to the constant function
$c_{x_0}$
for some $x_0\in M$,
then $M$ is smoothly contractible
(see, e.g., \cite[Ex.~15.4(a)]{Gl20a}).
We shall not use this fact.
\end{rem}
Two lemmas are helpful.
\begin{lem}\label{hat1}
Let $M$ be a $C^\infty$-manifold,
$a<b$ in~$\R$
and $\pi\colon E\to [a,b]\times M$
be a smooth vector bundle over $[a,b]\times M$
whose typical fibre is a locally convex space~$F$.
If there exist $\alpha<\beta$ in $\,]a,b[$
such that $E|_{[a,\beta[\,\times M}$
and $E|_{]\alpha,b]\times M}$ are trivializable
smooth vector bundles, then~$E$
is smoothly trivializable.
\end{lem}
\begin{prf}
Let $\theta_1\colon E|_{[a,\beta[\,\times M}\to
[\alpha,\beta[\,\times M\times F$
and $\theta_2\colon E|_{]\alpha,b]\times M}\to
\;]\alpha,b]\times M\times F$ be $C^\infty$-trivializations,
with 2nd components $\theta_{1,2}$ and $\theta_{2,2}$,
respectively.
Then
\begin{equation}\label{trivcha}
]\alpha,\beta[\,\times M\times F\to \; ]\alpha,\beta[\,\times M\times F,
\;(t,x,y)\mto \theta_i(\theta_j^{-1}(t,x,y))
\end{equation}
is a $C^\infty$-diffeomorphism for $(i,j)\in\{(1,2),(2,1)\}$
and linear in the final argument.
We let $g_{i,j}\colon \,]\alpha,\beta[\,\times M\times F\to F$
be the third component of the diffeomorphism in~(\ref{trivcha});
thus
\[
\theta_1(\theta_2^{-1}(t,x,y))=(t,x,g_{1,2}(t,x,y)).
\]
Pick $r<s$ in $]\alpha,\beta[$.
There is a $C^\infty$-map $\tau\colon \,]\alpha,b]\to\R$
such that $\tau$ is monotonically increasing,
$\tau(t)=s$ for all $t\in [s,b]$
and $\tau(t)=t$ for all $t\in \,]\alpha,r]$.
Define $\theta\colon E\to [a,b]\times M\times F$
via
\[
\theta(v):=\left\{
\begin{array}{cl}
\theta_1(v) &\mbox{ if $t\in [a,r[$;}\\
(t,x,g_{1,2}(\tau(t),x,\theta_{2,2}(v))) &\mbox{ if $t\in \,]\alpha,b]$}
\end{array}\right.
\]
for $v\in E$,
with
$(t,x):=\pi(v)$.
Then $\theta$ is a $C^\infty$-diffeomorphism,
as we readily check that the inverse is the map taking
$(t,x,y)\in [a,b]\times M\times F$ to
$\theta_1^{-1}(t,x,y)$
and
$\theta_2^{-1}(t,x,g_{2,1}(\tau(t),x,y))$
if $t\in[a,r[$ and $t\in \,]\alpha,b]$,
respectively. Then $\theta$ is a global
$C^\infty$-trivialization for~$E$.
\end{prf}
\begin{lem}\label{hat2}
Let
$M$ be a $C^\infty$-manifold,
$a<b$ in $\R$,
and $\pi\colon E\to [a,b]\times M$
be a smooth vector bundle over $[a,b]\times M$
whose typical fibre is a locally convex space~$F$.
Then each $x\in M$ has an open neighborhood
$U\sub M$ such that $E|_{[a,b]\times U}$
is smoothly trivializable.
\end{lem}
\begin{prf}
For each $s\in [a,b]$,
there exists an open neighborhood $V_s$ of $(s,x)$ in
$[a,b]\times M$ such that $E|_{V_s}$
is smoothly trivializable. We may assume
that $V_s=J_s\times U_s$
for open subsets $J_s\sub [a,b]$
and $U_s\sub M$.
Let $\delta>0$
be a Lebesgue number for the open cover $(J_s)_{s\in[a,b]}$
of the compact metric space $[a,b]$.
Pick $a=t_0<t_1<\cdots<t_n=b$ such that
$t_j-t_{j-1}<\delta/3$ for all $j\in\{1,\ldots,n\}$.
Set $t_{-1}:=t_0$ and $t_{n+1}:=t_n$.
For each $j\in\{1,\ldots,n\}$,
we find $s_j\in [a,b]$ such that
\[
[t_{j-2},t_{j+1}]\sub J_{s_j},
\]
as the interval on the left has length $<\delta$.
Then
\[
U:=U_{s_1}\cap\cdots\cap U_{s_n}
\]
is an open neighborhood of~$x$ in~$M$.
A straightforward induction based
on Lemma~\ref{hat1} shows that
$E|_{[a,t_{j+1}]\times U}$ is smoothly trivializable
for all $j\in\{1,\ldots, n\}$.
Notably, $E|_{[a,b]\times U}$
is smoothly trivializable.
\end{prf}
\begin{prop}\label{eq-boundary}
Let $F$ be a locally convex space, $M$ be a $\sigma$-compact,
finite-dimensional $C^\infty$-manifold,
$a<b$ in $\R$ and $\pi\colon E\to [a,b]\times M$
be a smooth vector bundle with typical fibre~$F$.
For $t\in [a,b]$, let $\lambda_t\colon M\to [a,b]\times M$
be the map $x\mto (t,x)$.
Then there exists an isomorphism
$\lambda_a^*(E)\to\lambda_b^*(E)$
of smooth vector bundles over $\id_M$. 
\end{prop}
\begin{prf}
Abbreviate $I:=[a,b]$.
By Lemma~\ref{hat2},
$M$ admits a cover $\cU$ by open sets $U\sub M$
such that $E|_{I\times U}$
is smoothly trivializable.
By smooth paracompactness (Proposition~\ref{findim-smoothly-para}),
we find a smooth partition of
unity $(h_j)_{j\in J}$ on~$M$
such that $S(j):=\Supp(h_j)\sub U(j)$ for some $U(j)\in\cU$,
for each $j\in J$. Since~$M$ is $\sigma$-compact,
$h_j\not=0$ for only countably many $j\in J$.
We may therefore assume that $J$ is countable
and actually that $J=\N$.
For $n\in \N_0$, define the smooth map
\[
g_n:=a+(b-a)\sum_{j=1}^nh_j\colon M\to I.
\]
Consider the pullback bundle
$f_n^*(E)$
over~$M$ determined by the smooth map $f_n\colon
M\to I \times M$, $x\mto (g_n(x),x)$.
Thus $f_0=\lambda_a$.
For $n\in\N_0$, the interior
\[
W_n:=\{x\in M\colon g_n(x)=b\}^0
\]
is an open subset of~$M$ such that the open subset
$f_n^*(E)|_{W_n}$ of $f_n^*(E)$
and the open subset $\lambda_b^*(E)|_{W_n}$ of $\lambda_b^*(E)$
coincide as smooth vector bundles; thus
\[
f_n^*(E)|_{W_n}=\lambda_b^*(E)|_{W_n}.
\]
Likewise,
$f_n^*(E)|_{M\setminus S_n}=\lambda_a^*(E)|_{M\setminus S_n}$
with $S_n:=\bigcup_{j=1}^n\Supp(h_j)$.
We now construct a sequence $(\alpha_n)_{n\in\N_0}$
of $C^\infty$-vector bundle isomorphisms
\[
\alpha_n\colon \lambda_a^*(E)\to f_n^*(E)
\]
over $\id_M$
such that, for all
$m\geq n$ in $\N_0$,
we have
\begin{equation}\label{loc-statio}
\alpha_m(v)=\alpha_n(v)\mbox{ for all $v\in \lambda_a^*(E)|_{W_n}$.}
\end{equation}
Once this is accomplished, we get a well-defined map
\[
\alpha\colon \lambda_a^*(E)\to \lambda_b^*(E)
\]
if we send $v\in \lambda_a^*(E)$ to
$\alpha(v):=\alpha_n(v)$,
independent of the choice of $n\in\N_0$
such that $v\in \lambda_a^*(E)|_{W_n}$.
By construction, $\alpha$ is an isomorphism
of $C^\infty$-vector bundles over~$\id_M$.\\[2.3mm]
For $j\in \N$, we let
$\theta_j\colon E|_{I\times U(j)}\to I \times U(j)\times F$
be a $C^\infty$-trivialization of $E|_{I\times U(j)}$,
with second component $\theta_{j,2}\colon E|_{I\times U(j)}\to F$.
We let $\alpha_0$ be the identity map $\lambda_a^*(E)\to \lambda_a^*(E)$.
To define a $C^\infty$-vector bundle isomorphism
\begin{equation}\label{stepbyst}
\alpha_{n,n-1}\colon f_{n-1}^*(E)\to f_n^*(E)
\end{equation}
over $\id_M$ for $n\in\N$, recall that
the local trivialization $\theta_n$
of~$E$ yields local trivializations
\[
\Theta_n \colon f_n^*(E)|_{U(n)}\to U(n)\times F,\;\,
(x,y)\mto (x,\theta_{n,2}(y))
\]
and
\[
\Xi_n \colon f_{n-1}^*(E)|_{U(n)}\to U(n)\times F,\;\,
(x,y)\mto (x,\theta_{n,2}(y))
\]
of the pullback bundles $f_n^*(E)$
and $f_{n-1}^*(E)$, respectively
(where $x\in U(n)$ and $y\in E_{f_n(x)}$,
resp., $y\in E_{f_{n-1}(x)}$).
Note that
\[
\Theta_n^{-1}(x,z)=(x,\theta_n^{-1}(g_n(x),x,z))
\]
for $(x,z)\in U(n)\times F$,
and
$\Xi_n^{-1}(x,z)=(x,\theta_n^{-1}(g_{n-1}(x),x,z))$.
Since $g_n(x)=g_{n-1}(x)$
if $x\in U(n)\setminus S(n)$, we see that
\begin{equation}\label{hencewelldf}
\Theta_n^{-1}|_{(U(n)\setminus S(n))\times F}
=
\Xi_n^{-1}|_{(U(n)\setminus S(n))\times F}.
\end{equation}
We define $\alpha_{n,n-1}$ as in (\ref{stepbyst})
via $\alpha_{n,n-1}(v):=v$ if $v\in f_{n-1}^*(E)|_{M\setminus S(n)}$
and
\[
\alpha_{n,n-1}(v):=(\Theta_n^{-1}\circ\Xi_n)(v)
\]
for $v\in f_{n-1}^*(E)|_{U(n)}$;
the map is well defined by~(\ref{hencewelldf}).
By construction, $\alpha_{n,n-1}$ is a $C^\infty$-vector
bundle isomorphism over~$\id_M$.
Moreover,
\begin{equation}\label{thusgood}
Y:=f_n^*(E)|_{M\setminus S(n)}=f_{n-1}^*(E)|_{M\setminus S(n)}
\mbox{ and } \alpha_{n,n-1}|_Y=\id_Y.
\end{equation}
We get a $C^\infty$-vector bundle
isomorphism $\alpha_n\colon \lambda_a^*(E)\to f_n^*(E)$ over $\id_M$ via
\[
\alpha_n:=\alpha_{n,n-1}\circ \cdots\circ \alpha_{1,0}.
\]
As $S_k$ and $W_n$ are disjoint for $k>n$,
the maps $\alpha_m=\alpha_{m,m-1}\circ \cdots\circ \alpha_{n+1,n}\circ\alpha_n$
and $\alpha_n$ coincide on $\lambda_a^*(E)|_{W_n}$ for all $m>n$.
\end{prf}
\begin{cor}\label{contractible-trivial}
Let $M$ be a smoothly contractible, $\sigma$-compact
smooth manifold of finite dimension.
Then every smooth vector bundle
$E$ over~$M$ $($with typical fibre
a locally convex space~$F)$
is trivializable.
\end{cor}
\begin{prf}
Let $H\colon [0,1]\times M\to M$
be a $C^\infty$-homotopy
from $\id_M$ to $c_{x_0}\colon M\to M$
with $x_0\in M$.
By Proposition~\ref{eq-boundary},
the $C^\infty$-vector bundles
$E=H(0,\cdot)^*(E)$
and $H(1,\cdot)^*(E)\cong M\times F$
are $C^\infty$-isomorphic over $\id_M$.
\end{prf}
\subsection*{Sections, frames, and normal bundles}
We introduce sections of vector bundles
(see Section~\ref{secspacemfd} for further studies).
Normal bundles are a tool
in the construction of tubular neighborhoods.
\begin{defn}
Consider a $C^r$-vector bundle $\pi\colon E\to M$, with \break 
 ${r\in\N_0\cup\{\infty,\omega\}}$.
We say that a mapping $\sigma\colon M\to E$ is a \emph{$C^r$-section}
of~$E$ if~$\sigma$ is $C^r$ and $\sigma(x)\in E_x$ for all $x\in M$, i.e.,
$\pi\circ\sigma=\id_M$.
If $W$ is an open subset of~$M$
and $\sigma\colon W\to E$
a $C^r$-map with $\pi\circ\sigma=\id_W$,
then $\sigma$ is called a \emph{local $C^r$-section}.
\end{defn}
\begin{rem}
(a) The map $0_E\colon M\to E$ taking $x\in M$ to the zero-vector $0_x$ in the vector space~$E_x$
is called the \emph{zero-section} of~$E$; it is $C^r$ as $0_E|_U=\theta^{-1}(x,0)$
for each local trivialization $\theta\colon E|_U\to U\times F$.
Thus $0_E$ is a $C^r$-section of~$E$
(sometimes, also $0_E(M)\sub E$ is called the zero-section).\medskip

\noindent
(b) Consider a map $\sigma\colon M\to E$ with $\pi\circ\sigma=\id_M$.
If $\sigma$ is a $C^r$-section, then
\begin{equation}\label{locrepisCk}
\sigma_\theta:=\theta_2\circ \sigma|_U\in C^r(U,F)
\end{equation}
for each local trivialization $\theta\colon E|_U\to U\times F$ of~$E$.
If, conversely, each point $x\in M$ has an open neighborhood $U\sub M$
such that there exists a trivialization $\theta\colon E|_U\to U\times F$ for which
(\ref{locrepisCk}) holds, then $\sigma$ is $C^r$ (and hence a $C^r$-section)
as $\sigma|_U=\theta^{-1}\circ (\id_U,
\theta_2\circ\sigma|_U)$ for all~$U$ and $\theta$ as before.
\end{rem}
\begin{defn}
Let $F$ be a locally convex space, $r\in \N_0\cup\{\infty,\omega\}$
and $\pi\colon E\to M$ be a $C^r$-vector bundle with typical fibre~$F$
over a $C^r$-manifold~$M$ modeled on a locally convex space~$Z$,
which may have a rough boundary if $r\not=\omega$.
Let $Y\sub F$ be a closed vector subspace.
A subset $H\sub E$ is called a \emph{vector subbundle}
of~$E$ with typical fibre~$Y$
if, for each $x\in M$,
there exists a local trivalization $\theta\colon E|_U\to U\times F$
of~$E$ with $x\in U$ such that $\theta(H)=U\times Y$.
\end{defn}

We endow $H_x:=H\cap E_x$ with the vector space structure
making the second component $\theta_2|_{H_x}\colon H_x\to Y$
an isomorphism of vector spaces.
Then $H|_U:=\theta^{-1}(U\times Y)$
has a unique $C^r$-manifold structure modeled on $Z\times Y$
(with boundary $H_U\cap \pi^{-1}(\partial M)$)
making $\theta_H\colon H|_U\to U\times Y$, $v\mto \theta(y)$
a $C^r$-diffeomorphism.
Then $H$ admits a unique $C^r$-manfold structure
(with boundary $H\cap \pi^{-1}(\partial M)$)
making each $H|_U$ an open submanifold
(cf.\ Exercise~\ref{exc-mfd-local}).
Then $\theta_H$ s a local trivialization and hence $H$ a $C^r$-vector bundle
over~$M$ with typical fibre~$Y$ and bundle projection $\pi|_H$.
\begin{lem}\label{subbun-via-frame}
Let $F$ be a locally convex space, $r\in \N_0\cup\{\infty,\omega\}$
and ${\pi\colon E\to M}$ be a $C^r$-vector bundle with typical fibre~$F$
over a $C^r$-manfold~$M$ modeled on a locally convex space~$Z$,
which may have a rough boundary if $r\not=\omega$.
Let $n\in\N$ and $H\sub E$ be a subset such that
$H_x:=H\cap E_x$ is an $n$-dimensional vector subspace of
$E_x$ for each $x\in M$.
Then {\rm(a)} and {\rm(b)} are equivalent:
\begin{description}[(D)]
\item[\rm(a)]
$H$ is a vector subbundle of~$E$;
\item[\rm(b)]
For each $x_0\in M$, there exist an open $x_0$-neighborhood
$U\sub M$ and local $C^r$-sections $\sigma_1,\ldots,\sigma_n\colon U\to E$
of $M$ such that $\sigma_1(x),\ldots,\sigma_n(x)$
is a basis of $H_x$ for all $x\in U$.
\end{description}
\end{lem}
\noindent
In the situation of~(b), $(\sigma_1,\ldots,\sigma_n)$
is called a \emph{local frame} for~$H$.\smallskip

\noindent
\begin{prf}
Assume that $H$ is a vector subbundle of~$E$
with typical fibre an $n$-dimensional vector subspace $V$ of~$F$.
Let $b_1,\ldots, b_n$ be a basis for~$V$.
If $\theta\colon E|_U\to U\times F$
is a local trivialization with $x_0\in U$ and $\theta(H|_U)=U\times V$,
then the functions $\sigma_j\colon U\to H$, $x\mto \theta^{-1}(x,b_j)$
define a local frame $(\sigma_1,\ldots, \sigma_n)$
for $H$. Thus (a)$\Rightarrow$(b).

Conversely, assume (b). For $x_0\in M$,
let $U$ and $\sigma_1,\ldots,\sigma_n\colon U\to E$
be as in~(b).
After shrinking $U$,
there exists a local trivialization
$\theta\colon E|_U\to U\times F$ of~$E$;
let $\theta_2\colon U\to F$ be its second component.
Let $V$ be the $\K$-linear span
of the $b_j:=\theta_2(\sigma_j(x_0))$ for $j\in\{1,\ldots, n\}$.
Let $b_1^*,\ldots, b_n^*$ be the dual basis of $V^*$
determined by $b_i^*(b_j)=\delta_{ij}$
for $i,j\in\{1,\ldots, n\}$,
using Kronecker's delta.
The finite-dimensional vector subspace
$V$ of~$F$ has a topological complement $W$ in~$F$.
Let $\pr_V(v+w):=v$ and $\pr_W(v+w):=w$ for $v\in V$, $w\in W$.
For $x\in U$, we obtain $\K$-linear maps
$\alpha_x\colon V\to V$ and
 $\beta_x\colon V\to W$ via
via
\[
\alpha_x(v):=\sum_{j=1}^nb_j^*(v)\pr_V(\theta_2(\sigma_j(x)))
\;\;\mbox{and}\;\;
\beta_x(v):=\sum_{j=1}^n b_j^*(v)\pr_W(\theta_2(\sigma_j(x)))
\]
for $v\in V$. Then $\alpha_{x_0}=\id_V$. Since $U \ni x \mto\alpha_x\in \cL(V)$
is continuous, after shrinking~$U$ we may assume that $\alpha_x\in\GL(V)$
for all $x\in U$. Note that
the maps $U\times V\to V$, $(x,v)\mto\alpha_x(v)$
and $U\times V\to V$, $(x,v)\mto\alpha_x^{-1}(v)$ are $C^r$,
as well as the map $U\times V\to W$, $(x,v)\mto \beta_x(v)$.
Now
\[
h\colon U\times F\to E,\quad (x,v+w)\mto
\sum_{j=1}^nb_j^*(v)\sigma_j(x)+\theta^{-1}(x,w)
\]
for $y\in U$, $v\in V$, $w\in W$
is a $C^r$-map. To see that $h(U\times F)=E|_U$ and $h\colon U\times F\to E|_U$
is a $C^r$-diffeomorphism, it suffices to show that ${\theta(h(U\times F))=U\times F}$
and $\theta\circ h$ is a $C^r$-diffeomorphism. For $a\in V$ and $b\in W$
and $x\in U$, we have
\[
\theta(h(x,v+w))=\left( x,
\alpha_x(v)+\beta_x(v)+w\right)=(x,a+b)
\]
for $(v,w)\in V\times W$ if and only if
\begin{equation}\label{thusinvCr}
v=\alpha_x^{-1}(a)\quad
\mbox{and}\quad w=b-\beta_x(v).
\end{equation}
Thus $\theta\circ h$ is a bijection with
$(\theta\circ h)^{-1}(x,a+b)=(v,w)$ as in (\ref{thusinvCr})
a $C^r$-function of $(x,a+b)$.
We deduce that $h^{-1}\colon E|_U\to U\times F$
is a local trivialization for~$E$. For each $x\in U$, we have
$h(\{x\}\times V)=H_x$ by construction, whence $h^{-1}(H|_U)=U\times V$.
Thus $H$ is a vector subbundle.
\end{prf}
\begin{numba}\label{riemann-exists}
Let $M$ be a paracompact finite-dimensional smooth manifold, of dimension~$n$.
Then $M$ admits a \emph{Riemannian metric}, i.e., a
family $(g_x)_{x\in M}$
of positive definite bilinear form $g_x\colon T_xM\times T_xM\to \R$
such that
\begin{equation}\label{diff-prop-riemann}
V\times \R^n\times \R^n\to \R, \;\;
(x,v,w)\mto g_{\phi^{-1}(x)}(v,w)
\end{equation}
is smooth for each chart $\phi\colon U\to V\sub\R^n$ of~$M$.
In fact, there exists a smooth partition of unity $(h_j)_{j\in J}$
on~$M$ subordinate to the open cover of chart domains (see Proposition~\ref{findim-smoothly-para}).
Thus $\Supp(h_j)\sub U_j$ for some chart $\phi_j\colon U_j\to V_j\sub\R^n$ of~$M$.
Let $\langle\cdot,\cdot\rangle$ be a scalar product on $\R^n$.
Then
\[
g_x(v,w):=\sum_{j\in J} h_j(x)\,\langle d\phi_j(v),d\phi_j(w)\rangle
\]
defines a Riemannian metric on~$M$, reading summands
with $h_j(x)=0$ as $0$.
\end{numba}
\begin{prop}\label{normal-bundle-subbundle}
Let $M$ be a paracompact $C^\infty$-manifold of finite dimension,
$(g_x)_{x\in M}$ a Riemannan metric on~$M$
and $N\sub M$ a submanifold.
For $x\in N$, let $(T_xN)^\perp:=\{v\in T_xM\colon (\forall w\in T_xN)\;g_x(v,w)=0\}$.
Then
\[
TN^\perp := \bigcup_{x\in N} (T_xN)^\perp
\]
is a vector subbundle of $(TM)|_N$.
\end{prop}
\noindent
We call $TN^\perp$ the \emph{normal bundle} of~$N$.\smallskip

\noindent
\begin{prf}
Let $m$ and $n$ be the dimension of~$M$ and $N$, respectively.
For $x_0\in N$,
there exists a chart $\phi\colon U\to V\sub\R^m$ of $M$
around~$x_0$ such that $\phi(U\cap N)=V\cap(\R^n\times\{0\})$.
Consider the local trivialization
\[
\theta\colon TU\to U\times\R^n,\quad v\mto(\pi_{TM}(v),d\phi(v))
\]
of $TM$ associated with $\phi$.
Then $\theta$ restricts to a local trivialization
\[
\theta_N\colon TU\cap TM|_N\to (U\cap N)\times \R^m
\]
of $TM|_N$ with $\theta_N(T(U\cap N))=(U\cap N)\times (\R^n\times\{0\})$.
Using the standard basis vectors
$e_1,\ldots, e_m$ in $\R^m$, get a local frame $(\tau_1,\ldots,\tau_m)$
for $TM|_N$ via
\[
\tau_j(x):=\theta_N^{-1}(x,e_j)\;\;\mbox{for $j\in\{1,\ldots, m\}$ and $x\in U$.}
\]
We obtain local $C^\infty$-sections $\xi_j\colon U\cap N\to (TM)|_N$
and $\sigma_j\colon U\cap N \to (TM)|_N$ for $j\in\{1,\ldots, m\}$
using Gram--Schmidt orthonormalization:
\begin{eqnarray*}
\xi_j(x) &:= &\tau_j(x)-\sum_{i=1}^{j-1}g_x(\sigma_i(x),\tau_j(x))\, \sigma_i(x),\\
\sigma_j(x) & :=&
\frac{1}{\sqrt{g_x(\xi_j(x),\xi_j(x))}}\, \xi_j(x)
\end{eqnarray*}
for $x\in U\cap N$.
Then $\sigma_1(x),\ldots,\sigma_m(x)$ is an orthonormal basis
of $T_xM$ for all $x\in U\cap N$
and $\sigma_1(x),\ldots,\sigma_n(x)$ is a basis
of $T_xN$, entaling that $\sigma_{n+1}(x),\ldots,\sigma_m(x)$
is a basis for $(T_xN)^\perp$.
Thus $(\sigma_{n+1},\ldots,\sigma_m)$
is a local frame for $(TN)^\perp$.
Hence $(TN)^\perp$ is a vector subbundle of $(TM)|_N$,
by Lemma~\ref{subbun-via-frame}.
\end{prf}
\begin{numba}\label{sett-whit-sum}
Let $r\in \N_0\cup\{\infty,\omega\}$,
$M$ be a $C^r$-manifold which may have a rough boundary
if $r\not=\omega$.
Let $F_k$ be a locally convex topological $\K$-vector space
and $\pi_k\colon E_j\to M$
be a $C^r$-vector bundle over~$M$
for $k\in \{1,2\}$,
with typical fibre~$F_k$.
Let $E_{k,x}:=\pi_k^{-1}(\{x\}$ for $x\in M$.
Define
\[
E_1\oplus E_2:=\bigcup_{x\in M}(E_{1,x}\times E_{2,x})
\]
as a set and note that $(E_{1,x}\times E_{2,x})_{x\in M}$
is a family of pairwise disjoint sets.
We define $\pi(v_1,v_2):=x$ for $x\in M$ and $(v_1,v_2)\in E_{1,x}\times E_{2,x}$.
Moreover, we give $(E_1\oplus E_2)_x=E_{1,x}\times E_{2,x}$
the vector space structure as a direct product.
If $\theta_k\colon E_k|_U\to U\times F_k$
are trivializations of $E_k|_U$ for $k\in \{1,2\}$
for the same open subset $U\sub M$,
setting $(E_1\oplus E_2)|_U:=\pi^{-1}(U)$ we have a bijection
\[
\Theta\colon (E_1\oplus E_2)\to U\times F_1\times F_2,\;\;
(v_1,v_2)\mto (\theta_{1,2}(v_1),\theta_{2,2}(v_2)),
\]
where $\theta_{2,2}\colon E_2|_U\to F_2$ is the second component
of~$\theta_2$. Then
\[
\Theta^{-1}(x,y_1,y_2)=(\theta_1^{-1}(x,y_1),\theta_2^{-1}(x,y_2))\;\;
\mbox{for all $(x,y_1,y_2)\in U\times F_1\times F_2$.}
\]
We let $\cA$ be the set of all the functions $\Theta$
as before. We let $(\Theta_i)_{i\in I}$
be a family of bijections $\Theta_i\colon (E_1\oplus E_2)|_{U_i}\to U_i\times F_1\times F_2$
such that $\cA=\{\Theta_i\colon i \in I\}$.
Let $\Theta_i$ be obtained from the trivializations
$\theta_{i,k}$ of $E_k|_{U_i}$ for $k\in \{1,2\}$.
Let $Z$ be the modeling space of~$M$.
\end{numba}
\begin{prop}\label{whitney-sum}
In the situation of {\rm\ref{sett-whit-sum}},
there is a unique $C^r$-manifold structure
on $E_1\oplus E_2$ with boundary $\pi^{-1}(\partial M)$
which is modeled on $Z\times F_1\times F_2$,
makes $(E_1\oplus E_2)|_{U_i}$
an open subset
and $\Theta_i$ a $C^r$-diffeomorphism,
for each $i\in I$.
The structure makes $E_1\oplus E_2$
a $C^r$-vector bundle over~$M$
with typical fibre~$F_1\times F_2$.
\end{prop}
$E_1\oplus E_2$ is called the \emph{Whitney sum} of $E_1$ and $E_2$.
\begin{prf}
We give $E_1\oplus E_2$ the final topology $\cO$ with respect to
the mappings $\Theta_i^{-1}\colon U_i\times F_1\times F_2\to E_1\oplus E_2$
for $i\in I$. Then 
\[ \Theta_i^{-1}(\im(\Theta_j^{-1}))=(U_i\times U_j)\times 
F_1\times F_2\] is open for all $i,j\in I$.
We have $\theta_{i,1}(\theta_{j,1}^{-1}(x,y_1))=(x,g_{ij}(x,y_1))$
and $\theta_{i,2}(\theta_{j,2}^{-1}(x,y_2))=(x,h_{ij}(x,y_2))$
for all $x\in U_i\cap U_j$ and $y_1\in F_1$, $y_2\in F_2$,
with $C^r$-functions $g_{ij}\colon U_i\cap U_j\times F_1\to F_1$
and $h_{ij}\colon U_i\cap U_j\times F_2\to F_2$.
Then
\[
(\Theta_i\circ\Theta_j^{-1})(x,y_1,y_2)=(x,g_{ij}(x,y_1),h_{ij}(x,y_2))
\]
is $C^r$ in $(x,y_1,y_2)\in U_i\cap U_2\times F_1\times F_2$
and hence a $C^r$-diffeomorphism, as the inverse map is
$\Theta_j\circ\Theta_i^{-1}$. By Exercise~\ref{topforbun}, each of the maps $\Theta_i^{-1}$
has open image in $(E_1\oplus E_2,\cO)$ and is a homeomorphism
onto its image. Since $\pi\circ \Theta_i^{-1}$
is the projection $(U_1\times U_2)\times F_1\times F_2\to M$,
$(x,y_1,y_2)\mto x$, hence $C^r$ and thus continuous,
we can separate $v\not=w$ in $ E_1\oplus E_2$
by open sets if $\pi(v)\not=\pi(w)$.
If $\pi(v)=\pi(w)$,
then $v,w$ are in the open subset
$(E_1\oplus E_2)|_{U_i}$ for some $i$, which is Hausdorff
as it is homeomorphic to $U_i\times F_1\times F_2$.
Thus $\cO$ is Hausdorff and $E_1\oplus E_2$
admits a unique $C^r$-manifold structure
making $(E_1\oplus E_2)|_{U_i}$
an open subset and $\Theta_i$
a $C^r$-diffeomorphism for each $i\in I$
(cf.\ Exercise~\ref{exc-mfd-local}).
Then each $\Theta_i$ is a local trivialization and
$E_1\oplus E_2$ a vector bundle.
\end{prf}
\subsection*{Principal bundles}
We define principal bundles which may be infinite-dimensional and
record basic facts. For brevity, we restrict attention to
smooth principal bundles.
\begin{numba}
Let $G$ be a group.
If $X$ is a set, then a map 
\[ \tau\colon X\times G\to X, \quad (x,g)\mto x.g\] 
 is called a \emph{right $G$-action}
(and $(X,\tau)$ a right $G$-set)
if $x.e=x$ and $x.(g_1g_2)=(x.g_1).g_2$ for all $x\in X$ and $g_1,g_2\in G$.
If $(X_1,\tau_1)$ and $(X_2,\tau_2)$ are right $G$-sets, then a map
$f\colon X_1\to X_2$ is called \emph{equivariant} if $f
(x.g)=f(x).g$ for all $(x,g)\in X_1\times G$.
If $X$ is a smooth manifold (possibly with rough boundary)
and $\tau$ is smooth, then $(X,\tau)$ is called a smooth $G$-manifold.
\end{numba}
\begin{numba} Let $G$ be a Lie group
and $U$ be a smooth manifold.
The $U\times G$ is a smooth $G$-manifold with the smooth right $G$-action
$\tau\colon (U\times G)\times G\to U\times G$ defined via
$(x,g).h:=(x,gh)$ for $x\in X$, $g,h\in G$.
Smooth $G$-manifolds of this form are called
\emph{trivial smooth $G$-principal bundles}
(with rough boundary, if $U$ is a smooth manifold with rough boundary).
\end{numba}
General principal bundles locally look like trivial ones.
\begin{defn}
Let $G$ be a Lie group and $M$
be a smooth manfold. A~smooth $G$-manifold~$P$,
together
with a surjective smooth mapping \break ${\pi\colon P\to M}$,
is called a \emph{$G$-principal bundle over~$M$}
if, for each $x\in M$,
there exist an open $x$-neighborhood $U\sub M$
and a $C^\infty$-diffeomorphism $\theta\colon P|_U\to U\times G$
on $P|_U:=\pi^{-1}(U)$ which
satisfies $\pr_1(\theta(x))=\pi(x)$ for all $x\in P|_U$
and is $G$-equivariant in the sense that
$\theta(x.g)=\theta(x).g$ for all $x\in P|_U$ and $g\in G$.
We call $\theta$ a local trivialization of~$P$.
If $P$ admits a trivialization $\theta\colon P\to M\times G$,
then $P$ is called \emph{trivializable}.
\end{defn}
\begin{rem}\label{principal-subm}
(a) $\pi$ is a submersion as
it locally looks like projections $U\times G\to U$.\smallskip

(b) In the same way, we can define smooth $G$-principal bundles
with rough boundary, replacing $M$ and $P$
by smooth manifolds with rough boundary.\smallskip

(c) A $C^\infty$-diffeomorphism $\theta\colon P|_U\to U\times G$
with second component $\theta_2\colon P|_U\to G$
and first component $\pi|_{P|_U}$ is a local trivialization
if and only if $\theta_2\colon P|_U\to G$
is equivariant if we use the smooth group multiplication
$\mu\colon G\times G\to G$
as the right $G$-action of $G$ on itself.
\end{rem}
Homogeneous spaces of Lie groups are a basic source of examples 
(cf.\ Proposition~\ref{prop:split-Lie}).
\begin{lem}\label{GH-principal-sub}
Let $G$ be a Lie group and $H$ be a closed subgroup of~$G$.
Assume that $G/H$ is endowed with a smooth manifold structure
which turns the canonical map $q\colon G\to G/H$, $g\mto gH$
into a smooth map.
Then the following conditions are equivalent:
\begin{description}[(D)]
\item[\rm(a)]
$G$, with the right action $G\times H\to G$, $(g,h)\mto gh$
and $q\colon G\to G/H$
is a smooth $H$-principal bundle.
\item[\rm(b)]
$q$ is a submersion.
\end{description}
In this case, $\sigma(W)H$ is open in~$G$
for each smooth local section $\sigma\colon W\to G$ of~$q$
and the map $\psi\colon W\times H\to\sigma(W)H$, $(x,h)\mto\sigma(x)h$
is a $C^\infty$-diffeomorphism.
\end{lem}
\begin{prf}
By Remark~\ref{principal-subm}, (a) implies~(b).
If~(b) holds,
then $q$ admits smooth local sections (see Proposition~\ref{subm-sec}(a)).
For $\sigma$ and $\psi$ as in the lemma, $\sigma(W)H=q^{-1}(W)$
is open in~$G$. By basic group theory, $\psi$ is a bijection.
Since $q(\sigma(x)h)=x$,
we readily see that the smooth map
$\theta\colon \sigma(W)H\to W\times H$, $g\mto (q(g),\sigma(q(g))^{-1}g)$
is the inverse of~$\psi$, whence $\psi$ is a $C^\infty$-diffeomorphism
and $\theta$ a local trivialization.
\end{prf}
\begin{rem}\label{like-vbdl}
Let $G$ be a Lie group.
Taking $r=\infty$,
we can define the restrictions $P|_S$
and pullback $f^*(P)$
of a $G$-principal bundle $\pi\colon P\to N$
as in Lemmas~\ref{restrict-vbdl}
and \ref{pullbackbdl},
replacing $E$ with $P$, $E_x$ with $P_x:=\pi^{-1}(\{x\})$,
$F$~with~$G$, vector spaces with smooth $G$-manfolds
and linearity with equivariance
in the proofs. Remark~\ref{restr-as-pb},
Lemmas~\ref{hat1} and \ref{hat2},
Proposition~\ref{eq-boundary}
and Corollary~\ref{contractible-trivial}
carry over as well, with the same
modifications. Notably:
\end{rem}
%
%
\begin{prop}\label{contr-principal-trivial}
If $M$ is a smoothly contractible, $\sigma$-compact
smooth manfold of finite dimension
and $G$ a Lie group,
then every smooth $G$-principal bundle
$\pi\colon P\to M$
is trivializable.\qed
\end{prop}
\begin{small}

\subsection*{Exercises for Section~\ref{tools-mfd}}

\begin{exer}
Show that $f^*(E)\to E$, $(x,v)\mto v$ is a $C^r$-vector bundle map
over~$f$ in the situation of Lemma~\ref{pullbackbdl}.
\end{exer}

\begin{exer}
Prove the assertion of Remark~\ref{restr-as-pb}.
\end{exer}
\begin{exer}\label{pullback-is-submfd}
If $N$ does not have a boundary in the situation of
Lemma~\ref{pullbackbdl}, show that $f^*(E)$
is a split submanifold of the $C^r$-manifold $M\times E$
modeled on the closed vector subspace $Y\times \{0\}\times F\cong Y\times F$
of $Y\times Z\times F$. The $C^r$-manifold structure on $f^*(E)$
as a submanifold coincides with the $C^r$-manifold structure constructed in
Proposition~\ref{pullbackbdl}.
\end{exer}

\begin{exer}
Check the details of Remark~\ref{like-vbdl}.
\end{exer}

\end{small}
\section{Sprays, local additions, and tubular neighborhoods}\label{sec:local-add}
\subsection*{Differential equations of second order and sprays}
Sprays on a smooth Banach manifold~$M$
are vector fields on $TM$
which can be used to formulate second order
differential equations on~$M$
with the additional property
that also $t\mto \gamma(ts)$
is a solution if $\gamma$ is so, for all $s\in\R$.
In this book, we are only interested in sprays
as a tool to obtain so-called local additions
(as in Definition~\ref{defn-loc-add}),
and we limit our discussion
to those aspects which are necessary for this goal.
\begin{numba}
Let $E$ be a Banach space and $V\sub E$ be an open subset;
as usual, we identify $TV$ with $V\times E$,
which is an open subset of $E\times E$.
Accordingly, we identify $T^2(V):=T(TV)$
with the open subset $T(V\times E)=V\times E\times E\times E$
of $E^4$.
If $I\sub\R$ is an open interval containing~$0$ 
and $\gamma\colon I\to V$ a $C^1$-curve, then
\[ \dot{\gamma}(t)=(\gamma(t),\gamma'(t))\in V\times E 
\quad \mbox{ for } \quad t\in I.\]
If $\gamma\colon I\to V$ is a $C^2$-curve, then
\[
\ddot{\gamma}(t)=(\dot{\gamma})\dot{}\hspace*{.1mm}(t)=(\gamma(t),\gamma'(t),\gamma'(t),\gamma''(t))
\in TV\times(E\times E)=T^2V.
\]
If $X\colon TV=V\times E\to T(TV)=TV\times (E\times E)$ is a smooth vector field
and for all $(x_0,v_0)\in V\times E$
the initial value problem
\begin{equation}\label{awpsev}
\dot{y}(t)=X(y(t)),\quad y(0)=(x_0,v_0)
\end{equation}
has a solution $(\gamma,\eta)$ with
$\eta=\gamma'$
(whence $\ddot{\gamma}(t)=X(\gamma(t),\gamma'(t))$,
$\gamma(0)=x_0$ und $\gamma'(0)=v_0$),
then we must have
\begin{equation}\label{spr1}
X(x,v)=(x,v,v,f(x,v))\;\,\mbox{for all $v\in V$}
\end{equation}
for some smooth function $f\colon V\times  E\to E$.
\end{numba}
\begin{numba}\label{doppel-frag}
For $(x_0,v_0)\in V\times E$, let
\[
(\gamma,\eta)\colon I\to V\times E
\]
be a solution to (\ref{awpsev}).
If (\ref{spr1}) holds, then
\[
(\gamma(t),\eta(t),\gamma'(t),\eta'(t))=
(\gamma(t),\eta(t),\eta(t),f(\gamma(t),\eta(t)))
\]
for all $t\in I$, whence
\[
\eta(t)=\gamma'(t)
\]
(i.e., $(\gamma,\eta)=\dot{\gamma}$)
and
\begin{equation}\label{sosec}
\gamma''(t)=f(\gamma(t),\gamma'(t)).
\end{equation}
If, conversely, $(\gamma(0),\gamma'(0))=(x_0,v_0)$
and (\ref{sosec}) holds,
then (\ref{awpsev})
holds for the function $\dot{\gamma}$.
\end{numba}
\begin{numba}
Assuming (\ref{spr1}), let
$(\gamma,\gamma')\colon I\to V\times\R^m$
be a solution to~(\ref{awpsev}) on an open interval~$I$.
For $s\in\R$,
\[
I_s:=\{t\in\R\colon ts\in I\}
\]
is an open interval in $\R$ such that $0\in I_s$.
Consider the smooth function
\[
\gamma_s\colon I_s\to V,\quad t\mto \gamma(ts).
\]
Then $(\gamma_s)'(t)=s\gamma'(st)$
and $(\gamma_s)''(t)=s^2\gamma''(st)$;
thus $(\gamma_s,(\gamma_s)')$ satisfies the differential equation
$\dot{y}(t)=X(y(t))$ if and only if
\[
s^2\gamma''(st)=(\gamma_s)''(t)=f(\gamma(st),s\gamma'(st)),
\]
i.e., $s^2 f(\gamma(st),\gamma'(st))=f(\gamma(st),s\gamma'(st))$.
Hence (a) and (b) are equivalent:
\begin{description}[(D)]
\item[(a)]
For all $(x_0,v_0)\in TV$,
using the maximal solution
\[ (\gamma,\gamma')=\gamma_{0,x_0,v_0}\colon I\to TV \] 
to~(\ref{awpsev}),
we have that $(\gamma_s,(\gamma_s)') \colon I_s\to TV$
solves the differential equation $\dot{y}(t)=X(y(t))$
(and hence the initial value problem $\dot{y}(t)=X(y(t))$,
$y(0)=(x_0,sv_0)$),
for all $s\in \R$.
\item[(b)] The following identity is satisfied:
\begin{equation}\label{spr2}
f(x,sv)=s^2f(x,v)\;\,\mbox{for all $(x,v)\in V\times\R^m$
and $s\in\R$.}
\end{equation}
\end{description}
\end{numba}
\begin{defn}\label{defn-spray-local}
If $E$ is a Banach space and $V\sub E$ an open subset,
we call a smooth vector field
$X\colon TV\to T^2V$
a \emph{spray on~$V$}
if it is of the form~(\ref{spr1})
and satisfies (\ref{spr2}).
\end{defn}
The following lemma will enable us to pass to manifolds.
\begin{lem}\label{change-chart-spray}
Let $\tau\colon V\to W$ be a $C^\infty$-diffeomorphism
between open subsets $V$ and $W$
of a Banach space~$E$.
If $X\colon TV\to T(TV)$
is a spray on~$V$, then
$Y:=T(T\tau)\circ X\circ T\tau^{-1}$
is a spray on~$W$.
\end{lem}
\begin{prf}
We have
\[
T\tau\colon V\times E\to W\times E,
\;\,
(x,v)\mto (\tau(x),d\tau(x,v))\vspace{-.3mm}
\]
and $T(T\tau)\colon (V\times E)\times (E\times E)\to
(W\times E)\times(E\times E)$,
\begin{eqnarray*}
(x,v,a,b) &\mto &
(T\tau(x,v),dT\tau(x,v,a,b))\\
&=& (\tau(x),d\tau(x,v),d\tau(x,a),
d^{(2)}\tau(x,v,a)+d\tau(x,b)).
\end{eqnarray*}
For $(y,w)\in W\times E$
let $(x,v):=T\tau^{-1}(y,w)$;
thus $(y,w)=T\tau(x,v)=(\tau(x),d\tau(x,v))$.
Then
\[
X(T\tau^{-1}(y,w))=X(x,v)=(x,v,v,f(x,v))
\]
using a smooth function $f\colon V\times E\to E$
as in~(\ref{spr1}).
Accordingly,
\begin{eqnarray*}
Y(y,w) &=&  T^2\tau(X(x,v))=
(\tau(x),d\tau(x,v),d\tau(x,v),g(y,w))\\
&=& (y,w,w,g(y,w))
\end{eqnarray*}
with $g(y,w):=d^{(2)}\tau(x,v,v)+d\tau(x,f(x,v))$,
whence $Y$ is of the form~(\ref{spr1}).
If we replace $w$ with $sw$ for some $s\in\R$,
then $v$ has to be replaced with $sv$.
Thus
\begin{eqnarray*}
g(y,sw) &=& d^{(2)}\tau(x,sv,sv)+d\tau(x,f(x,sv))\\
&=& s^2(d^{(2)}\tau(x,v,v)+d\tau(x,f(x,v)))
=s^2g(y,w);
\end{eqnarray*}
as a consequence, $Y$ satisfies (\ref{spr2}),
and thus $Y$ is a spray.
\end{prf}
\begin{lem}\label{conv-spray}
If $X_1,\ldots, X_n$ are sprays on an open subset $V$
of a Banach space~$E$
and $h_1,\ldots,h_n\colon V\to\R$
smooth functions such that
$\sum_{j=1}^n h_j=1$,
then also $\sum_{j=1}^n (h_j\circ \pi_{TV})X_j$
is a spray on~$V$.
\end{lem}
\begin{prf}
For each $j$, we
have $X_j(x,v)=(x,v,v,f_j(x,v))$
for a smooth function $f_j$ satisfying $f_j(x,sv)=s^2f_j(x,v)$. Thus
\begin{eqnarray*}
\Big(\sum_{j=1}^n (h_j\circ \pi_{TV})X_j(x,v)\Big) &=&
\Big(x,v,\sum_{j=1}^n h_j(x)(v,f_j(x,v))\Big)\\[-.6mm]
&=& \Big(x,v,\sum_{j=1}^nh_j(x)v,\sum_{j=1}^nh_j(x)f_j(x,v)\Big)\\[-.3mm]
&=& (x,v,v,f(x,v))\vspace{-1.2mm}
\end{eqnarray*}
with $f(x,v):=\sum_{j=1}^n h_j(x)f_j(x,v)$.
Since $f_j(x,sv)=s^2f_j(x,v)$ for each $j$,
we have $f(x,sv)=\sum_{j=1}^nh_j(x)f_j(x,sv)
=s^2f(x,v)$.
\end{prf}
\begin{rem}\label{trivial-rem-spray}
If $X\colon TV\to T^2V$, $(x,v)\mto (x,v,v,f(x,v))$
is a spray on an open subset $V\sub E$,
then, for each $x\in V$, the maximal solution
$\eta$ of the initial value problem
\[
\dot{y}(t)=X(y(t)),\quad y(0)=(x,0)
\]
is defined on all of~$\R$,
and
$\eta(t)=(x,0)$ for all $t\in\R$.
\end{rem}
In fact, if we set
$\eta(t):=(x,0)$
for $t\in\R$,
then $\eta(0)=(x,0)$
and $\dot{\eta}(t)=(x,0,0,0)$.
Since $f(x,0)=f(x,00)=0^2f(x,0)=0$, we have
$X(\eta(t))=
X(x,0)=(x,0,0,f(x,0))=(x,0,0,0)=\dot{\eta}(t)$.
\begin{defn}\label{defn-global-spray}
Let $M$ be a smooth manifold modeled on a Banach space~$E$.
A smooth vector field
$X\colon TM\to T(TM)$
is called a {\it  spray} on $M$ if
$T^2\phi\circ X\circ T\phi^{-1}$
is a spray on $V_\phi$
for each chart
$\phi\colon U_\phi\to V_\phi\sub E$
of~$M$.
\end{defn}
By Lemma~\ref{change-chart-spray},
Definitions \ref{defn-spray-local} and \ref{defn-global-spray}
are equivalent if~$M=V$ is an open subset of~$E$.
\begin{prop}\label{sprays-exist}
Every smoothly paracompact
smooth manifold modeled on a Banach space
admits a spray.
\end{prop}
\begin{prf}
Let $\cA$
be the maximal $C^\infty$-atlas of~$M$,
comprised of charts of the form $\phi\colon U_\phi\to V_\phi\sub E$.
We choose a smooth partition of unity $(h_j)_{j\in J}$
on~$M$ which is subordinate to $(U_\phi)_{\phi\in\cA}$.
For each $j\in J$; there exists a chart
$\phi_j\colon U_j\to V_j$
with $\Supp(h_j)\sub U_j$.
Then $Y_j \colon TV_j\to T(TV_j)$,
$(x,v)\mto (x,v,v,0)$
is a spray on~$V_j$, whence
$X_j:=T^2(\phi^{-1}_j)\circ Y_j\circ T\phi_j$
is a spray on~$U_j$. We show that
\begin{equation}\label{try}
X\colon TM\to T(TM),\quad
v\mto \sum_{j\in J}h_j(\pi_{TM}(v))X_j(v)
\end{equation}
is a finite sum for each $v\in TM$
(if $h_j(\pi_{TM}(v))=0$, the summand
should be read as~$0$),
and that~$X$ is a spray.
Given $v_0\in TM$, the element
$x_0:=\pi_{TM}(v_0)$
has an open neighborhood~$W$ in~$M$
such that
\[
J_0:=\{j\in J\colon W\cap\Supp(h_j)\not=\emptyset\}
\]
is finite. If $j\in J_0$ and
$x_0\not\in\Supp(h_j)$, we can replace~$W$ by its intersection with
$X\setminus\Supp(h_j)$; we may therefore assume that
$x_0\in\Supp(h_j)\sub U_j$
for all $j\in J_0$, whence
$U:=W\cap \bigcap_{j\in J_0}U_j$
is an open $x_0$-neighborhood in~$M$.
Then all summands in (\ref{try}) with $j\in J\setminus J_0$
vanish for all
$v$ in the open
$v_0$-neighborhood $TU$ in $TM$;
thus
\[
X(v)=\sum_{j\in J_0}h_j(\pi_{TM}(v))X_j(v).
\]
Hence $X|_{TU}$ is smooth and Lemma~\ref{conv-spray}
implies that
$X|_{TU}$ is a spray on~$U$
(for $k\in J_0$, at each point 
$T^2(\phi_k|_U)\circ X|_{TU}\circ T(\phi_k|_U)^{-1}$
is an affine 
combination of the sprays
$T^2(\phi_k|_U)\circ X_j|_{TU}\circ T(\phi_k|_U)^{-1}$
with $j\in J_0$).
\end{prf}
\begin{prop}\label{expfu-spray}
Let $M$ be a $C^\infty$-manifold modeled on a Banach space~$E$
and
$X\colon TM\to T^2M$ be a spray.
Then the following holds:
\begin{description}
\item[\rm(a)]
If $\eta\colon I\to TM$ is a solution to the differential equation
$\dot{y}(t)=X(y(t))$,
then
$\eta=\dot{\gamma}$
holds for $\gamma:=\pi_{TM}\circ\,\eta$.
\item[\rm (b)]
For the flow
$\Fl\colon \Omega\to TM$
of the differential equation $\dot{y}(t)=X(y(t))$
on $TM$,
the set $\Omega_{1,0}=\{v\in TM\colon (1,0,v)\in\Omega\}$
is an open subset of~$TM$
and $0_M(M)=\{0_x\in T_xM\colon x\in M\}\sub \Omega_{1,0}$.
\item[\rm(c)]
The map
$\exp_X\colon \Omega_{1,0}\to M$, $\,v\mto \pi_{TM}(\Fl_{1,0}(v))$
is smooth and satisfies
$T_0\exp_X|_{\Omega\cap T_xM}=\id_{T_xM}$
for all $x\in M$.
\end{description}
\end{prop}
Here
$\Omega\sub\R\times\R\times TM$
is the domain of the flow;
for the notations $\Omega$, $\Fl$, $\Fl_{1,0}$,
and $\Omega_{1,0}$, see Definition~\ref{flow-nonauto-m}.
\begin{prf}
(a) We may assume that
$\eta(I)\sub TU_\phi$
for a chart $\phi\colon U_\phi\to V_\phi\sub E$
of~$M$, as $I$ can be covered by such intervals.
Now
\[ Y:=T^2\phi\circ X|_{TU_\phi}\circ T\phi^{-1} \] 
is a spray on~$V_\phi$.
Since $X|_{TU_\phi}$ and $Y$ are $T\phi$-related,
we deduce that
${\zeta:=T\phi\circ \eta}$
solves the differential equation
$\dot{y}(t)=Y(y(t))$.
Thus
$\zeta=\dot{\theta}$ with ${\theta:=\pi_{TV_\phi}\circ \zeta=\pr_1\circ\, \zeta}$
(see \ref{doppel-frag}).
For
\[
\gamma:=\pi_{TM}\circ \eta=\pi_{TM}\circ T\phi^{-1}\circ\zeta
=\phi^{-1}\circ\pi_{TV_\phi}\circ\zeta=\phi^{-1}\circ\theta,
\]
we have 
$\dot{\gamma}=T\phi^{-1}\circ \dot{\theta}=T\phi^{-1}\circ\zeta=\eta$.\\[1mm]
(b) and (c):
For each $x\in M$, we have $0_x\in \Omega_{1,0}$
and $\exp_X(x)=x$ (see
Remark~\ref{trivial-rem-spray}; hence
$0_M(M)\sub \Omega_{1,0}$.
Now $\Omega_{1,0}$
is open in $TM$
and $\Fl$ is smooth (see Corollary~\ref{niceflow-ban-m}),
whence $\exp_X$ is smooth.
For $x\in M$ and $v\in T_xM$,
we have $rv\in \Omega_{1,0}$
for small $r>0$.
Let $\eta\colon I\to M$ with $[0,1]\sub I$
be a solution to the initial value problem
$\dot{y}(t)=X(y(t))$,
$y(0)=rv$. Let
$\gamma:=\pr_{\!TM}\circ\,\eta$,
so that $\eta=\dot{\gamma}$, by~(a).
For all $s\in [0,1]$,
$t\mto \eta(st)$ then is the solution to the initial value problem
with $srv$ in place of~$rv$.
In particular, taking $t=1$ we obtain that
\[
\gamma(s)=\pi_{TM}(\eta(st))=\exp_X(srv).
\]
Hence
\begin{eqnarray*}
T_{0_x}(\exp_X|_{T_xM\cap \Omega_{1,0}})(rv)
&=& [s\mto\exp_X(srv)]=[s\mto \gamma(s)]\\[.3mm]
&=& \dot{\gamma}(0)=\eta(0)=rv.\vspace{-3mm}
\end{eqnarray*}
Thus
$T_{0_x}(\exp_X|_{T_xM\cap \Omega_{1,0}})(v)=v$,
entailing that
$T_{0_x}(\exp_X|_{T_xM\cap \Omega_{1,0}})=\id_{T_xM}$.
\end{prf}
\noindent
The map $\exp_X$ is called the \emph{exponential function}
associated with the spray.
\subsection*{Local additions and tubular neighborhoods}
Local additions are an important tool for the creation
of infinite-dimensional manifold structures
(see Section~\ref{sec-mfdmps}). We shall also use them
in the construction of tubular neighborhoods (Theorem~\ref{thm-tubular}).
\begin{defn}\label{defn-loc-add}\label{local-add}
Let $M$ be a smooth manifold.
A smooth map
$\Sigma\colon U\to M$
on an open subset $U\sub TM$
is called a \emph{local addition}
if
\begin{description}[(D)]
\item[(a)]
$0_x\in U$
for all $x\in M$ and $\Sigma(0_x)=x$
(for $0_x\in T_xM$); and
\item[(b)]
The set $U':=\{(\pi_{TM}(v),\Sigma(v))\colon v\in U\}$
is open in $M\times M$ and the map
$\theta\colon U\to U'$, $v\mto (\pi_{TM}(v),\Sigma(v))$
is a $C^\infty$-diffeomorphism.
\end{description}
We say that $\Sigma$ is \emph{normalized} if, moreover,
\begin{equation}\label{bettersigma}
T_{0_x}(\Sigma|_{T_xM\cap Q})=\id_{T_xM}\;\,\mbox{for all $x\in M$.}
\end{equation}
If $M$ is a $\K$-analytic manifold with $\K\in\{\R,\C\}$
and $\Sigma$ a local addition for which
$\theta$ is a diffeomorphism of $\K$-analytic
manifolds, then $\Sigma$ is called a $\K$-analytic
local addition.
\end{defn}
Lie groups are a first source of examples.
\begin{prop}\label{lie-local-add}\label{exliegploa}
Let $G$ be a Lie group, $\g:=T_eG$
and $G\times TG\to TG$, $(g,v)\mto g.v:=T\lambda_g(v)$
be the smooth left action of $G$ on $TG$ defined by the group multiplication. 
If $\phi\colon V\to W$ a $C^\infty$-diffeomorphism
from an open $0$-neighborhood $V\sub \g$
onto an open identity neigborhood $W\sub G$
such that $\phi(0)=e$,
then $U:=\{v\in TG\colon \pi_{TG}(v)^{-1}.v\in V\}$
is an open subset of $TG$ with $0_G(G)\sub U$ and the smooth map
\[
\Sigma\colon U\to G,\;\;
v\mto \pi_{TG}(v)\phi(\pi_{TG}(v)^{-1}.v)
\]
is a local addition. If $T_0\phi=\id_\g$,
then $\Sigma$ is normalized.
If $G$ is a $\K$-analytic Lie group and $\phi$
is a diffeomorphism of $\K$-analytic manifolds,
then $\Sigma$ is a $\K$-analytic local addition.
\end{prop}
\begin{prf}
The set $U':=\{(\pi_{TG}(v),\Sigma(v))\colon v\in U\}$
equals 
\[ \bigcup_{g\in G}(\{g\}\times gW)=\{(g,h)\in G\times G\colon g^{-1}h\in W\},\]
whence it is open in $G\times G$.
The map $\theta$ is a $C^\infty$-diffeomorphism. 
We readily check that the map 
\[ U'\to TG, \quad (g,h)\mto g.\phi^{-1}(g^{-1}h) \] 
is inverse to $\theta$, whence $\theta\colon U\to U'$
is a $C^\infty$-diffeomorphism.
If $T_e\phi=\id_{\g}$, applying the chain rule for tangent maps to
\[
\Sigma|_{U\cap T_gG}=\lambda_g\circ \phi\circ T\lambda_{g^{-1}},
\]
we obtain
$T_{0_g}(\Sigma|_{U\cap T_gG})=T\lambda_g\circ T_0\phi \circ T_{0_g}\lambda_{g^{-1}}
=\id_{T_gG}$.
For the $\K$-analytic case, replace ``$C^\infty$'' with ``$\K$-analytic''
in the proof.
\end{prf}
By Proposition~\ref{findim-smoothly-para}, the following
proposition applies to every finite-dimensional
smooth manifold which is paracompact.
\begin{prop}\label{local-add-parac}\label{exfindimloa}
Let $M$ be a smoothly paracompact smooth manifold
modeled on a Banach space.
Then $M$ admits a normalized local addition.
\end{prop}
The following lemma will be used in the proof.
\begin{lem}\label{tool-loca-exists}
Let
$M$ be a paracompact smooth manifold modeled on a Banach
space
and $\Sigma\colon V\to M$ be a smooth map
on an open subset $V\sub TM$ with $0_M(M)\sub V$
such that
\begin{equation}\label{norma-gives-loca}
T_{0_x}\Sigma|_{T_xM\cap V}=\id_{T_xM}\;\,
\mbox{for all $x\in M$.}
\end{equation}
Then there exists an open subset $U\sub V$
with $0_M(M)\sub U$ such that $\Sigma|_U$
is a normalized local addition.
\end{lem}
\begin{prf}
If $M$ is finite dimensional, then $M$
is a topological sum of $\sigma$-compact
open subsets, whence so is $M\times M$,
entailing that $M\times M$ is paracompact
(see Proposition~\ref{lcp-parac}).
In the general case, every $x\in M$ has an open neighborhood $U_x$
which is homeomorphic to an open subset
of $E$ and hence metrizable.
Since $M$, being paracompact, is regular (see Proposition~\ref{paracomp-reg}),
$U_x$ contains a neighborhood $A_x$
of~$x$ which is closed in~$M$.
By paracompactness, there is a locally finite
open cover $(Q_j)_{j\in J}$ of~$M$
which is subordinate to $(A_x^0)_{x\in M}$.
Then also $(\wb{Q_j})_{j\in J}$
is locally finite (see Lemma~\ref{ops-locfin}(a))
and subordinate to $(A_x)_{x\in M}$
and hence subordinate to $(U_x)_{x\in M}$,
entailing that each of the sets $\wb{Q_j}$
is metrizable. Hence $M$ is metrizable, by
\cite[Thm.~4.4.19]{Ek89}.
As a consequence, $M\times M$ is metrizable,
whence $M \times M$ is paracompact, by Stone's
Theorem (see \cite[Thm.~4.4.1]{Ek89}).
The smooth map
$h:=(\pi_{TM}|_V,\Sigma)$
takes $0_M(M)=\{0_x\colon x\in M\}$
bijectively and continuously onto the diagonal
$\Delta_M:=\{(x,x)\colon x\in M\}$
in $M\times M$
and also the inverse function
$(x,x)\mto 0_x=0_M(x)$ is continuous, as
$0_M\colon M\to TM$, $x\mto 0_x$
is continuous.
Given $x_0\in M$, let $\phi\colon U_\phi\to V_\phi\sub E$
be a chart of~$M$ around $x_0$ such that $\phi(x_0)=0$.
The map
\[
\alpha\colon E\times E\to E\times E,\quad (y,z)\mto (y,y+z)
\]
is an isomorphism of topological vector spaces with
inverse $(y,w)\mto (y,w-y)$.
The open set $T\phi(TU_\phi\cap h^{-1}(U_\phi\times U_\phi))$
contains $W\times W$ for some open $0$-neighborhood $W\sub E$.
Thus
\[
f:=\alpha^{-1}\circ (\phi\times\phi)\circ h\circ T(\phi^{-1})\colon
W\times W\to E\times E 
\]
is a smooth map. Since $f(y,0)=(y,0)$ for all $y\in W$, we have 
\[ df((y,0),(v,0))=d_1f(y,0,v)=(v,0)\quad \mbox{  for all } \quad v\in E.\] 
As a consequence of (\ref{norma-gives-loca}), $df((y,0),(0,w))=d_2f((y,0),w)=(0,w)$
for all $w\in E$. Hence $f'(0,0)=df((0,0),\cdot)=\id_{E\times E}$.
Using the Inverse Function Theorem,
we deduce that $f$ is a local $C^\infty$-diffeomorphism 
at $0$, whence $h$ is a local $C^\infty$-diffeomorphism at~$0_x$.
Thus, there exists an open $0_x$-neighborhood $Q_x\sub V$
such that $h(V)$ is open in $M\times M$
and $h|_{Q_x}\colon Q_x\to h(Q_x)$ is a $C^\infty$-diffeomorphism.
After replacing $V$ with $\bigcup_{x\in M} Q_x$,
we may assume that $h$ is a local $C^\infty$-diffeomorphism.
Theorem~\ref{godement-globalize}
now yields an open subset $U\sub V$
with $0_M(M)\sub U$
such that $U':=h(U)$ is open in $M\times M$
and $\theta :=h|_U\colon U\to U'$
is a homeomorphism.
Since $\theta$ also is a local $C^\infty$-diffeomorphism,
the homeomorphism $\theta$ is a $C^\infty$-diffeomorphism.
\end{prf}
\noindent
\emph{Proof of Proposition}~\ref{local-add-parac}.
By Proposition~\ref{sprays-exist},
there exists a spray $X$ on~$M$.
By Proposition~\ref{expfu-spray},
its exponential function $\exp_X$
satisfies the hypothesis of Lemma ~\ref{tool-loca-exists}.
Thus $\exp_X$ restricts to a local addition on~$M$.\vspace{1.2mm}\qed

\noindent
Similar arguments can be used to prove the following theorem.
\begin{thm}\label{thm-tubular}
Let $M$ be a paracompact, finite-dimensional smooth manifold
and $N\sub M$ be a closed submanifold of~$M$.
Then there exists a smooth vector bundle
$\pi\colon E\to N$ over~$N$,
an open subset $Q\sub E$ with $0_N(N)\sub Q$
and a $C^\infty$-diffeomorphism
$\psi\colon Q\to P$
for some open subset $P\sub M$ with $N\sub P$,
such that $\psi(0_x)=x$ for all $x\in N$.
\end{thm}
\noindent
In the previous situation,
$P$ is called a \emph{tubular neighborhood} of~$N$
in~$M$.
The proof shows that $E$ can be chosen as the normal bundle
$(TN)^\perp\sub TM|_N$ for a Riemannan metric on~$M$.
\begin{prf}
Let $\Sigma\colon U\to M$
be a local addition for~$M$ (see Proposition~\ref{local-add-parac})
and $g$ be a Riemannan metric on~$M$ (see \ref{riemann-exists}).
Then the normal bundle
$(TN)^\perp$ 
is a submanifold of $TM|_N$
(see Proposition~\ref{normal-bundle-subbundle}
which is a submanifold of $TM$ (see Lemma~\ref{restrict-vbdl}),
whence
$W:=(TN)^\perp\cap U$ is open in $(TN)^\perp$.
Moreover,
$0_N(N)=(TN)^\perp\cap 0_M(M)\sub W$.
The restriction 
\[ h:=\Sigma|_W\colon W\to M \] is smooth and maps
$0_N(N)=\{0_x\colon x\in N\}$
bijectively and continuously onto~$N$.
The inverse map
$0_N\colon N\to 0_N(N)$ is continuous.
Let $i\colon N\to M$
be the inclusion map.
Since $h\circ O_N=i$, we see that for $x\in N$ the image
$\im(T_{0_x}h)$ contains
$\im(T_x(h\circ O_N))=\im T_xi=T_xN$.
Using the smooth inclusion map $j\colon (T_xN)^\perp\cap W\to W\sub TM$
we have $(h\circ j)(y)=\Sigma|_{(T_xN)^\perp\cap Q}(y)$,
whence $\im(T_xh)$ also contains $\im(T_0 (h\circ j))=
\id_{T_xM}((T_xN)^\perp)
=(T_xN)^\perp$.
As $T_xM=T_xN \oplus (T_xN)^\perp$, we see that the linear map
\[ T_{0_x}h\colon T_{0_x}(W)=T_{0_x}((TN)^\perp)\to T_xM \] 
is surjective and hence bijective,
the domain and range having the same finite dimension.
For each $x\in N$, the Inverse Function Theorem 
yields an open {$0_x$-neighborhood} $W_x\sub W$
such that $h(W_x)$ is open in~$M$
and \break {$h|_{W_x}\colon W_x\to h(W_x)$}
is a $C^\infty$-diffeomorphism.
After replacing $W$ with $\bigcup_{x\in M}W_x$,
we may assume that $h$ is a local $C^\infty$-diffeomorphism.
Exploiting the paracompactness of~$M$,
Theorem~\ref{godement-globalize}
provides an open subset $Q\sub W$
with $0_N(N)\sub Q$ such that $P:=h(Q)$ is open in~$M$
and $\psi:=h|_Q\colon Q\to P$
is a homeomorphism.
Since $\psi$ also is a local $C^\infty$-diffeomorphism,
$\psi$
is a \break {$C^\infty$-diffeomorphism.}
\end{prf}
\begin{rem}\label{numbetsigm}
If a smooth manifold $M$ admits a local addition, then it also admits
a normalized local addition.\medskip

In fact, if $\Sigma\colon U\to M$ is arbitrary,
then the diffeomorphism
$\theta:=(\pi_{TM},\Sigma)\colon U\to U'\sub M\times M$
takes the open set $T_xM\cap U\sub T_xM$ to the submanifold
$(\{x\}\times M)\cap U'$ of~$M\times M$ for each $x\in M$
and restricts to a \break {$C^\infty$-diffeomorphism} between these sets.
Hence $\Sigma|_{T_xM\cap U}$ is a {$C^\infty$-diffeo\-morphism}
onto an open subset of~$M$, whence
$\alpha_x:=T_{0_x}(\Sigma|_{T_xM})\in\GL(T_xM)$ for each $x\in M$ and
\[
f\colon TM\to TM,\;\, f(v):= \alpha_{\pi_{TM}(v)}^{-1}(v)
\]
is a smooth map. Now $\Sigma\circ f\colon f^{-1}(U)\to M$
is a normalized local addition.
\end{rem}

\begin{lem}\label{prodlocadd}
For $j\in\{1,2\}$, let
$M_j$ be a smooth manifold which
admits a local addition $\Sigma_j\colon U_j\to M_j$.
Identifying $T(M_1\times M_2)$ with $TM_1\times TM_2$,
we can consider $\Sigma_1\times \Sigma_2\colon U_1\times U_2\to M_1\times M_2$
as a local addition for $M_1\times M_2$.
\end{lem}
\begin{prf}
In fact, the map $\psi\colon M_1\times M_2\times M_1\times M_2\to M_1\times M_1\times M_2\times M_2$,
$(a,b,c,d)\mto (a,c,b,d)$ is a $C^\infty$-diffeomorphism
and $\psi\circ (\pi_{TM},\Sigma_1\times\Sigma_2)
=(\pi_{TM_1},\Sigma_1)\times (\pi_{TM_2},\Sigma_2)$
corestricts to a $C^\infty$-diffeomorphism onto its open image $U_1'\times U_2'$,
where $M:=M_1\times M_2$.
\end{prf}
\subsection*{Appendix: Alternative constructions}
If $f\colon M\to N$ is a local diffeomorphism
between manifolds and $f|_A$ is injective
for a subset $A\sub M$,
one would like to find an open neighborhood $U$ of~$A$ in~$M$
such that $f|_U$ is injective and hence a diffeomorphism onto an open subset of~$N$.
We already discussed such situations, exploiting paracompactness.
In this appendix, we describe
complementary situations in which~$U$ exists.
As a preliminary, we study locally injective maps between
topological spaces.
\begin{prop}\label{around-compact-homeo}
Let $f\colon X\to Y$ be a continuous map between
Hausdorff topological spaces and
$K$ be a compact subset of~$X$.
Assume that $f|_K$ is injective
and each $x\in K$ has a neighborhood
$W(x)$ in~$X$ on which $f$ is injective.
Then $f|_U$ is injective for
some open subset $U\sub X$ such that $K\sub U$.
\end{prop}
\begin{prf}
We may assume that $K\not=\emptyset$.
Let $J$ be the set of all open subsets $U\sub X$
such that $K\sub U$.
Writing $U\leq V$ for $U,V\in J$ if and only if $V\sub U$,
we obtain a directed set $(J,\leq)$.
If the conclusion of the proposition was wrong,
then we could find $x_U\not=y_U\in U$
for each $U\in J$ such that $f(x_U)=f(y_U)$.
We claim that the net $(x_U,y_U)_{U\in J}$
in $X\times X$ has a subnet, say $(x_{U(i)},y_{U(i)})_{i\in I}$,
which converges to some
$(x,y)\in K\times K$.
If this is true, then
\[
f(x)=\lim f(x_{U(i)})=\lim f(y_{U(i)})=f(y)
\]
by continuity, whence $x=y$.
For $i\in I$ sufficently large, we have
$x_{U(i)},y_{U(i)}\in W(x)$ and would deduce from $f(x_{U(i)})=f(y_{U(i)})$
that $x_{U(i)}=y_{U(i)}$, which is absurd. Thus, the conclusion
of the proposition must hold.

To establish the claim, suppose it was wrong.
Then, for each $(x,y)\in$ $K\times K$ there exist an open $(x,y)$-neighborhood
$V(x,y)$ in $X\times X$ and $U(x,y)\in J$
such that $(x_U,y_U)\not\in V(x,y)$ for all $U\in J$ with
$U\geq U(x,y)$. In fact, otherwise
\[
I:=\{(V, U)\in \cV(x,y)\times J\colon (x_U,y_U) \in V\}
\]
would be a directed set via $(V,U)\leq (V',U')$
if $V'\sub V$ and $U'\sub U$,
writing $\cV(x,y)$ for the set of open neighborhoods of $(x,y)$ in $X\times X$.
Moreover, $(x_U,y_U)_{(V,U)\in I}$
would be a subnet of $(x_U,y_U)_{U\in J}$
which converges to $(x,y)$; but we already ruled out the existence
of such subnets.
Thus $U(x,y)$ and $V(x,y)$ do exist.
By compactness of $K\times K$,
there exists a non-empty finite subset $\Phi\sub K\times K$ such that
\[
K\times K\sub \bigcup_{(x,y)\in\Phi}V(x,y)=:V.\vspace{-2mm}
\]
By the Wallace Lemma,
there are open subsets $P,Q\sub X$ with $K\sub P$,
$K\sub Q$ and $P\times Q\sub V$.
Then $W:=\bigcap_{(x,y)\in\Phi}U(x,y)\in J$.
Also $U:=W\cap P\cap Q\in J$ and $U\geq W\geq U(x,y)$
for all $(x,y)\in \Phi$
whence $(x_U,y_U)\not\in V(x,y)$ and thus $(x_U,y_U)\not\in V$. But $(x_U,y_U)\in U\times U\sub P\times Q\sub V$,
contradiction.
\end{prf}
If $f$ is, moreover, a local homeomorphism
in Proposition~\ref{around-compact-homeo},
then $f(U)$ is open in~$Y$ and $f|_U\colon U\to f(U)$
is a homeomorphism.
Likewise, we conclude:
\begin{cor}\label{around-compact-diffeo}
Let $r\in\N\cup\{\infty,\omega\}$
and $f\colon M\to N$ be a local \break $C^r$-diffeomorphism between
$C^r$-manifolds
modeled on locally convex spaces
$($or between $C^r$-manifolds with rough boundary, if $r\not=\omega)$.
If $K\sub M$ is a compact subset such that $f|_K$ is injective,
then there exists an open subset $U\subset M$
with $K\sub U$ such that $f(U)$ is open in~$N$ and $f|_U\colon U\to f(U)$
is a \break $C^r$-diffeomorphism.\qed
\end{cor}
\begin{prop}\label{gives-inj-around}
Let $X$ be a Hausdorff topological space, $Y$ be a locally compact space,
$f\colon X\to Y$ be a continuous map and $A\sub X$ be a
$\sigma$-compact
subset such that $f(A)$ is closed in~$Y$ and $f|_A$ a topological
embedding.
Assume that each $x\in A$ has  a neighborhood $W(x)$
in~$X$ such that $f|_{W(x)}$ is injective.
Then $f|_U$ is injective for some open subset $U\sub X$ with $A\sub U$.
\end{prop}
\begin{prf}
By Lemma~\ref{lem-rings},
there exist compact subsets $L_j\sub f(A)$
and open subsets $O_j\sub Y$ with $L_j\sub O_j$ such that
$f(A)=\bigcup_{j\in\N}L_j$ and
$|i-j|\leq 1$ for all $i,j\in\N$ such that $O_i\cap O_j\not=\emptyset$.
Set $K_j:=(f|_A)^{-1}(L_j)$
for $j\in\N$ and $K_0:=\emptyset$.
Let $j\in\N$.
Applying Proposition~\ref{around-compact-homeo}
to $f$
and the compact set $K_{j-1}\cup K_j\cup K_{j+1}$,
we find an open subset $V_j\sub X$
such that
\[
K_{j-1}\cup K_j\cup K_{j+1}\sub V_j
\]
and $f|_{V_j}$ is injective. Let $V_0:=X$.
For each $j\in\N$,
\[
U_j:=V_{j-1}\cap V_j\cap V_{j-1}\cap f^{-1}(O_j)
\]
then is an open subset of~$X$ such that $K_j\sub U_j$
and thus $A\sub \bigcup_{j\in \N} U_j=:U$.
If $x,y\in U$ and $f(x)=f(y)$, then there exist $i,j\in\N$
such that $x\in U_i$ and $y\in U_j$.
Then $f(x)=f(y)\in O_i\cap O_j$
and thus $|i-j|\leq 1$, whence $\{x,y\}\sub U_i\cup U_j\sub V_j$
and thus $x=y$ by injectivity of $f|_{V_j}$.
\end{prf}
Proposition~\ref{gives-inj-around}
readily entails the following.
\begin{cor}\label{gives-diffeo-around}
Let $r\in \N\cup\{\infty,\omega\}$ and $f\colon M\to N$
be a $C^r$-map between finite-dimensional $C^r$-manifolds
$($or between locally compact $C^r$-manfolds
with rough boundary, if $r\not=\omega)$.
Let $A\sub M$ be a $\sigma$-compact
subset such that $f(A)$ is closed in~$N$ and $f|_A$ a topological
embedding. Moreover, assume
that each $x\in A$ has  an open neighborhood~$W$
in~$M$ such that $f(W)$ is open in~$N$ and $f|_W\colon W\to f(W)$ a $C^r$-diffeomorphism.
Then there exists an open subset $U\sub M$ with $K\sub U$ such that $f(U)$ is open in~$N$
and $f|_U\colon U\to f(U)$ is a $C^r$-diffeomorphism.\qed
\end{cor}
\begin{small}
\subsection*{Exercises for Section~\ref{sec:local-add}}

\begin{exer}
Let $M$ be a real analytic finite-dimensional manifold
and $(g_x)_{x\in M}$ be a Riemannian metric on~$M$
which is real analytic in the sense that the map in (\ref{diff-prop-riemann})
is real analytic for each chart $\phi\colon U\to V$ of~$M$.
Show that $TN^\perp$ is a vector subbundle of the real
analytic vector bundle $(TM)|_N$,
for each submanifold $N$ of~$M$ which is closed in~$M$.
\end{exer}

\noindent
Any $\sigma$-compact, finite-dimensional real analytic manifold
can be embedded in~$\R^n$ for some $n\in\N$ (see \cite{Gr58b});
pulling back the euclidean metric,
a real analytic Riemannian metric on~$M$
is obtained. To the latter, a spray $X\colon TM\to T^2(M)$
can be associated in a canonical way (see \cite[Ch.\,V, Thm.\,7.1]{La99}),
and the latter is real analytic.

\begin{exer}
Let $M$ be a finite-dimensional real
analytic manifold and \break {$X\colon TM\to T^2(M)$} be a real analytic spray.
Show that the associated exponential map (as in Proposition~\ref{expfu-spray})
is real analytic. If $M$ is paracompact,
show that the associated local addition
$\Sigma\colon U\to M$ (as in the proof of Proposition~\ref{local-add-parac})
is a real analytic local addition.
Using such a local addition, show that the tubular neighborhood
around a closed submanifold $N\sub M$
constructed in Theorem~\ref{thm-tubular}
is real analytic in the sense that the map $\psi$ is a $C^\omega$-diffeomorphism.
\end{exer}

\begin{exer}\label{submfd-spray}
Let $M$ be a $C^\infty$-paracompact smooth manifold
modeled on a Banach space.
Let $N\sub M$ be a submanifold.
Let $X\colon TM\to T^2(M)$ be a spray on~$M$
as constructed in the proof of Proposition~\ref{sprays-exist},
using now charts of~$M$ adapted to
the submanifold~$N$ instead of general charts of~$M$.
Show that $X(TN)\sub T^2(N)$ holds and $X|_{TN}\colon TN\to T^2(N)$
is a spray for~$N$. Also show that the exponential map of $X|_{TN}$
is a restriction of the exponential map $\exp_X$.
Let $\Sigma\colon U\to M$ be an associated local addition, as in
the proof of Proposition~\ref{local-add-parac},
which is a restriction of $\exp_X$.
Consider $\theta:=(\pi_{TM},\Sigma)$ as a $C^\infty$-diffeomorphism onto
its open image~$U'$.
If $N$ is closed in~$M$, show that, after shrinking~$U$,
one can achieve that $\Sigma|_{U\cap TN}$ is a local addition
for~$N$. After replacing~$U'$
with an open subset intersecting $N\times N$ in the relatively open set $\theta(U\cap TN)$,
one can achieve that, moreover, $U'\cap (N\times N)=\theta(U\cap TN)$.
\end{exer}

\end{small}
\mlabel{app:liealg-cohom} 
\section{Appendix: Lie algebra 
cohomology} 

The cohomology of Lie algebras is the natural tool to understand how 
we can build new Lie algebras $\hat \g$ from given Lie algebras $\g$ and $\fa$ in 
such a way that $\fa \trile \hat\g$ and $\hat\g/\fa \cong \g$. An important special 
case of this situation arises if $\fa$ is assumed to be abelian. We will see in particular 
how the abelian extensions of Lie algebras can be parametrized by a certain cohomology space.

Let $\K$ be a topological field of characteristic zero (all field operations are assumed to be 
continuous). A {\it topological Lie algebra $\g$} is a $\K$-Lie algebra which is a topological 
vector space for which the Lie bracket is a continuous bilinear map. 
A {\it topological $\g$-module} is a $\g$-module $V$ which is a topological vector 
space for which the module structure, viewed as a map $\g \times  V\to V$, 
$(x,v) \mapsto x.v$ is continuous. 
Note that every module $V$ of a Lie algebra $\g$ over a field $\K$ 
becomes a topological module if we endow $\K$, $\g$ and $V$ with the 
discrete topology. In this sense, all the following applies in particular to general 
modules of Lie algebras over fields of characteristic zero. 

\begin{defn} \mlabel{def:C.1} 
Let $V$ be a {\it topological module of the
topological Lie algebra $\g$}. 
For $p \in \N_0$, let $C^p_c(\g,V)$ denote the space of continuous 
alternating maps $\g^p \to V$, i.e., 
the {\it Lie algebra $p$-cochains with values in the module $V$}. 
We write
\[ C^\bullet_c(\g,V) := \bigoplus_{p \in \N_0} C^p_c(\g,V).\] 
Note that $C^1_c(\g,V) = {\cal L}(\g,V)$ is the space of continuous linear
maps $\g \to V$. We use the convention $C^0_c(\g,V) = V$. 
We then obtain a chain complex with the differential 
$$ d_\g \: C^p_c(\g,V) \to C^{p+1}_c(\g,V) $$
given on $f \in C^p_c(\g,V)$ by 
\begin{align*}
(d_\g f)(x_0, \ldots, x_p) 
&:= \sum_{j = 0}^p (-1)^j x_j.f(x_0, \ldots, \hat x_j, \ldots, 
x_p) \cr
& + \sum_{i < j} (-1)^{i + j} f([x_i, x_j], x_0, \ldots, \hat
x_i, \ldots, \hat x_j, \ldots, x_p), 
\end{align*}
where $\hat x_j$ indicates omission of $x_j$. Note that the continuity of the bracket 
on $\g$ and the action on $V$ imply that $d_\g f$ is continuous 
and hence an element of $C^{p+1}_c(\g,V)$,
as $d_\g f$ is alternating (see Exercise~\ref{exc-e-1}).

For elements of low degree, we have in particular: 
\begin{align*}
p=0: & \quad  d_\g f(x) = x.f \cr
p=1: & \quad d_\g f(x,y) = x.f(y) - y.f(x) - f([x,y])
\cr 
p=2: & \quad d_\g f(x,y,z) \\ 
& = x.f(y,z) - y.f(x,z) +
                                 z.f(x,y) \\
  &\qquad - f([x,y],z) + f([x,z],y) - f([y,z],x)\cr 
& = \sum_\cyc x.f(y,z) - f([x,y],z),  
\end{align*}
where we have used the notation 
$$ \sum_\cyc \gamma(x,y,z) := \gamma(x,y,z) +\gamma(y,z,x)
+\gamma(z,x,y). $$
In this sense, the Jacobi identity reads $\sum_\cyc [[x,y],z] = 0.$

Below we shall show that $d_\g^2 = 0$, so that the space 
$Z^p_c(\g,V) := \ker(d_\g\res_{C^p_c(\g,V)})$
of {\it $p$-cocycles} contains the space 
$B^p_c(\g,V) :=  d_\g(C^{p-1}_c(\g,V))$
of {\it $p$-coboundaries}. The quotient 
$$ H^p_c(\g,V) := Z^p_c(\g,V)/B^p_c(\g,V) $$
is the {\it $p$-th continuous cohomology space of $\g$ with values in
the $\g$-module $V$}. We write $[f] := f + B^p_c(\g,V)$ for the 
cohomology class of the cocycle $f$. 
\end{defn} 

On $C^\bullet_c(\g,V)$, 
we have a natural representation of $\g$,
given for $x \in \g$ and $f \in C^p_c(\g,V)$ by the {\it Lie derivative} 
\begin{align*}
({\cal L}_x f)(x_1, \ldots, x_p) 
&= x.f(x_1, \ldots, x_p) - \sum_{j = 1}^p f(x_1, \ldots,
[x, x_j], \ldots, x_p) \cr
&= x.f(x_1, \ldots, x_p) 
+ \sum_{j = 1}^p (-1)^j f([x, x_j], x_1, \ldots, \hat x_j, 
\ldots, x_p)
\end{align*}
(see Exercise~\ref{exc-e-2}).
We further have for each $x \in \g$ an insertion map 
$$ i_x \: C^p_c(\g,V) \to C^{p-1}_c(\g,V), 
\quad 
\big(i_xf\big)(x_1, \ldots, x_{p-1}) = f(x, x_1, \ldots,
x_{p-1}), $$
where we define $i_x$ to be $0$ on $C^0_c(\g,V) \cong V$. 

\begin{lem}
  \mlabel{lem:C.2} For $x,y \in \g$, we have the following identities: 
\begin{description}
\item[\rm(1)]  ${\cal L}_x = d_\g \circ i_x + i_x \circ d_\g$ {\rm(Cartan Formula)}. 
\item[\rm(2)] $[{\cal L}_x, i_y] = i_{[x,y]}$. 
\item[\rm(3)] $[{\cal L}_x, d_\g] = 0$. 
\item[\rm(4)] $d_\g^2 = 0$. 
\item[\rm(5)] ${\cal L}_x(Z^p_c(\g,V)) \subeq B^p_c(\g,V)$. In particular, 
the natural $\g$-action on $H^p_c(\g,V)$ is trivial. 
\end{description}
\end{lem}

\begin{prf} (1) Using the insertion map $i_{x_0}$, we 
can rewrite the formula for the coboundary operator as 
\begin{align*}
\big(i_{x_0}.d_\g f\big)(x_1, \ldots, x_p)  
= &x_0.f(x_1, \ldots, x_p) - 
\sum_{j = 1}^p (-1)^{j-1} x_j.f(x_0, \ldots, \hat x_j, \ldots, 
x_p) \cr
& + \sum_{j = 1}^p (-1)^j f([x_0, x_j], x_1, 
\ldots, \hat x_j, \ldots, x_p)  \cr
&+ \sum_{1 \leq i < j} (-1)^{i + j} f([x_i, x_j], x_0, \ldots, \hat
x_i, \ldots, \hat x_j, \ldots, x_p) \cr
=& x_0.f(x_1, \ldots, x_p) - 
\sum_{j = 1}^p (-1)^{j-1} x_j.f(x_0, \ldots, \hat x_j, \ldots, 
x_p) \cr
& - \sum_{j = 1}^p f(x_1, \ldots, x_{j-1}, [x_0, x_j], x_{j+1}, 
\ldots, x_p)  \cr
&- \sum_{1 \leq i < j} (-1)^{i + j} f(x_0, [x_i, x_j], \ldots, \hat
x_i, \ldots, \hat x_j, \ldots, x_p) \cr
=& ({\cal L}_{x_0}f)(x_1, \ldots, x_p) 
- d_\g\big(i_{x_0} f\big)(x_1, \ldots, x_p). 
\end{align*}
This proves our assertion. 

(2) The explicit formula for ${\cal L}_x$  
implies that for $y = x_1$, we have 
$i_y {\cal L}_x = {\cal L}_x i_y - i_{[x,y]}.$

(3),(4) Let $\phi \: C^\bullet_c(\g,V) \to C^\bullet_c(\g,V)$ be a linear map 
for which there exists an $\eps \in \{\pm 1\}$ with 
$\phi \circ i_x = \eps i_x \circ \phi$ for all $x \in \g$
and a $k \in \N$ with $\phi(C^p_c(\g,V)) \subeq C^{p+k}_c(\g,V)$ for each 
$p \in \N_0$. We claim that $\phi = 0$. 
Since the operators $i_x \: C^p_c(\g,V) \to C^{p-1}_c(\g,V)$, $x \in \g$, separate the points, 
it suffices to show that 
$i_x \circ \phi = \eps \phi \circ i_x$ vanishes for each $x \in \g$. 
On $C^0_c(\g,V)$, this follows from the definition of $i_x$, and 
on $C^p_c(\g,V)$, $p \in \N$, we obtain it by induction. 

Now we prove (3). From (1) and (2) we get 
\begin{align*}
{\cal L}_{[x,y]}
&= [{\cal L}_x, {\cal L}_y] 
= [d_\g \circ i_x, {\cal L}_y] + [i_x \circ d_\g, {\cal L}_y] \cr
&= [d_\g, {\cal L}_y] \circ i_x + d_\g \circ i_{[x,y]} 
+ i_{[x,y]} \circ d_\g + i_x \circ [d_\g, {\cal L}_y] \cr
&= [d_\g, {\cal L}_y] \circ i_x + {\cal L}_{[x,y]} + i_x \circ [d_\g, {\cal L}_y],
\end{align*}
so that $\phi := [d_\g, {\cal L}_y]$ anticommutes with the operators 
$i_x$ ($\eps = -1$ and $k = 1$). Therefore the argument in the preceding paragraph 
shows that $\phi$ vanishes, which is (3). 

To obtain (4), we consider the operator $\phi = d_\g^2$. 
Combining (3) with the Cartan Formula, we get 
$$ 0 = [d_\g, {\cal L}_x] = d_\g^2 \circ i_x - i_x \circ d_\g^2, \leqno(C.2) $$
so that the argument above applies with $\eps = 1$ and $k =2$. 
This proves that $d_\g^2 = 0$. 

(5) follows immediately from the Cartan Formula~(1). 
\end{prf} 

\begin{defn} \mlabel{def:C.3} 
A linear subspace $W$ of a topological vector space $V$ is called 
{\it (topologically) split} if it is closed and there is a continuous linear map 
$\sigma \: V/W \to V$ for which the map 
$$ W \times V/W \to V, \quad (w,x) \mapsto w + \sigma(x) $$
is an isomorphism of topological vector spaces. Note that the closedness of 
$W$ guarantees that the quotient topology turns $V/W$ into a Hausdorff space which 
is a topological vector space with respect to 
the induced vector space structure. 
A continuous linear map $f \: V \to W$ between topological vector spaces  
is said to be {\it (topologically) split} if the subspaces 
$\ker(f) \subeq V$ and $\im(f) \subeq W$ are topologically split
and the induced map $V/\ker(f)\to\im(f)$ is an isomorphism
of topological vector spaces. 
\end{defn} 

\begin{rem}
  \mlabel{rem:C.4} Let $\g$ be a topological Lie algebra and 
$$  0 \to V_1 \sssmapright{\alpha} V_2 \sssmapright{\beta} V_3 \to 0
$$
be a topologically split short exact sequence of topological $\g$-modules.
Identifying $V_1$ with $\alpha(V_1) \subeq V_2$, we obtain injective maps 
$\alpha_p \: C^p_c(\g,V_1) \to C^p_c(\g,V_2)$ and surjective maps
$\beta_p \: C^p_c(\g,V_2) \to C^p_c(\g,V_3)$ which lead to a short exact
sequence 
$$  0 \to C^\bullet_c(\g,V_1) \sssmapright{\alpha_*} C^\bullet_c(\g,V_2)
\sssmapright{\beta_*} C^\bullet_c(\g,V_3)  \to 0 $$
of cochain complexes. These maps can be combined to a long exact sequence 
\begin{align*}
 0 
&\to H^0_c(\g,V_1) \to  H^0_c(\g,V_2) \to H^0_c(\g,V_3) \\
&\to H^1_c(\g,V_1) \to  H^1_c(\g,V_2) \to H^1_c(\g,V_3) \to \ldots, 
\end{align*}
where, for $p \in \N_0$, the connecting map 
$$ \delta \: H^p_c(\g,V_3) \to H^{p+1}_c(\g,V_1) $$
is defined by $\delta([f]) = [d_\g \tilde f]$, where
$\tilde f \in C^p_c(\g,V_2)$ satisfies $\beta \circ \tilde f = f$, 
which implies that $\im(d_\g \tilde f) \subeq V_1$ if $f$ is a cocycle. 
\end{rem}

\begin{small}
\subsection*{Exercises for Appendix~\ref{app:liealg-cohom}}
\begin{exer}\label{exc-e-1}
Let $f \in C^p_c(\g,V)$ and $x_0,\ldots,x_p\in \g$.
Show that
if there exist $a<b$ in $\{0,\ldots, p\}$ with $x_a=x_b$,
then $d_\g f(x_0,\ldots,x_p)=0$. 
Deduce that $d_\g f\in C^{p+1}_c(\g,V)$.
\end{exer}
\begin{exer}\label{exc-e-2}
  Show that
  \[ ({\cal L}_x({\cal L}_yf)-{\cal L}_y({\cal L}_xf))(x_1,\ldots,x_p)
=({\cal L}_{[x,y]}f)(x_1,\ldots,x_p) \] 
for all $f\in C^p_c(\g,V)$
and $x,y,x_1,\ldots,x_p\in\g$.
\end{exer}
\end{small}

\section{Notes and comments on Chapter~\ref{chapmanif}}
In this chapter,
we adapted the relevant concepts
from the theory of finite-dimensional manifolds
(tangent spaces and tangent maps;
vector fields; differential forms)
to the case of manifolds modeled
on locally convex spaces.
%
Frequently, we only needed to pick the appropriate
definition among the classical ones.
For example, both the geometric
(kinematic)
and formal definitions of tangent vectors
turned out to be suitable for
infinite-dimensional manifolds
(while the popular interpretation
as point derivations would be unsuitable
because locally convex spaces usually
admit too many of these,  and  smooth partitions of unity need not exist).
However, we also encountered
various differences compared to the finite-dimensional
case. We now recall the most relevant ones
and provide some more background information.
\begin{itemize}
\item
Although
we still have local charts $\phi$ available
for manifolds modeled on locally convex spaces,
we do not have
``local coordinates'' $\phi(p)=(x_1,\ldots, x_n)\in \R^n$
in the strict sense of the word.
Therefore all concepts have to be
introduced and discussed
in a coordinate-free way.
For Banach manifolds, such a coordinate-free
formulation has been spelled out
in Lang's classical text books
(like \cite{La99}).
However, further deviations
from classical definitions
are necessary here
which are not required in the
Banach case.
\item
While the definition of a tangent bundle $TM$
does not cause any problems,
it is not possible in general to define
a \emph{canonical} smooth vector bundle structure
on the dual bundle $T^*M:=\bigcup_{x\in M} (T_xM)'$
of a smooth manifold~$M$ modeled on a non-normable
locally convex space
(where the topological dual space $(T_xM)'$
is equipped, say, with the topology of uniform convergence
on bounded sets).
Behind this is the well-known
difficulty that the evaluation map
$\cL(E)\times E\to E$, $(A,x)\mto A(x)$
is discontinuous
for each non-normable locally convex space~$E$
and each vector topology on the space
$\cL(E)$ of continuous endomorphisms of~$E$
(cf.\ \cite{Ms63} and \cite[p.\,2]{KM97}).
An explicit counterexample
is given in \cite{Gl07e};
in this case, $TM$ is a trivial bundle,
but different trivializations of~$TM$
give rise to incompatible vector bundle
structures on $T^*M$.
Because of such pathologies,
it is not possible in general to define
$1$-forms (and higher differential forms)
as smooth sections in suitable smooth
vector bundles, as customary in the
finite-dimensional theory.
However, as first observed by Beggs \cite{Bg87},
this problem can be circumvented.
We have chosen an alternative definition
of differential forms
which is as elementary and
simple as possible.
Equivalently,
the reader may think of $E$-valued smooth
$p$-forms on~$M$ as smooth $E$-valued
maps $\omega\colon \bigoplus_{j=1}^p TM\to E$
on the $p$-fold Whitney sum of the vector
bundle $TM$, such that $\omega|_{(T_xM)^p}$
is an alternating $p$-linear map
for each $x\in M$.
We mention that the problems
just described
are absent in the Convenient Setting of
Analysis. In this framework,
the formation of cotangent
and tensor bundles
is always possible (see \cite{KM97}).
\item
Also our definition of the
Lie derivative $\cL_Y\omega$ of a differential form
had to deviate from the conventional
definition for finite-dimensional
and Banach manifolds
(as alluded to in Remark~\ref{liedban}),
because not every vector field~$Y$
admits a flow due to well-known pathologies
(as encountered in Exercises \ref{ode-no} and \ref{ode-many}).
\end{itemize}
Let us remark that
the submanifolds of a Banach manifold~$M$
considered in the literature are usually
what we call split submanifolds,
modeled
on a split vector subspace of the modeling space of~$M$
(see, e.g., the works of Bourbaki,
where non-split submanifolds are called
quasi-submanifolds). As the splitting condition
is very much less useful in the case
of general locally convex manifolds,
we prefer not to make it part of the definition of
a submanifold, but a mere extra condition
which can be satisfied or not.

For a more extensive discussion of immersions and submersions between
manifolds modeled on locally convex spaces, the reader
is referred to \cite{Gl15a} (or \cite{Ham82}
in the case of Fr\'{e}chet manifolds). 
In particular  \cite{Jo82} contains the observation that the map 
$f \: \R \to \R^2, f(x) = (x^2,x^3)$ has the universal property 
of immersions but its differential is not everywhere non-zero. 
Notably, variants of Proposition~\ref{naive-imm-subm}
are available when finite-dimensional vector spaces are replaced with Banach spaces.
The discussion of Proposition~\ref{eq-boundary} and its corollary
concerning vector bundles over smoothly contractible bases
vary \cite[pp.\ 20--21]{Ha02},
where finite-dimensional topological vector bundles
over paracompact topological spaces are considered.
They are a special case of \cite[\S15]{Gl20a}.

See \cite{KM97} for a thorough discussion of
smooth regularity and smooth paracompactness,
in the setting of convenient differential calculus.
Notably,
every Hilbert space is smoothly paracompact (see \cite[Cor.\,16.16]{KM97}
and the references given there), and so is every nuclear
Fr\'{e}chet space, and every nuclear Silva space (see \cite[Thm.\,16.10]{KM97}).
Also direct limits $\bigcup_{n\in\N}M_n$ of finite-dimensional smooth manifolds
$M_1\sub M_2\sub\cdots$ (as in Theorem~\ref{dl-mfd})
are known to be smoothly paracompact (see \cite[Prop.\,3.6]{Gl05c}).

For smoothing techniques in algebraic topology, see \cite{MW09, Wo09, Gl20a};
cf.\ also~\cite{Hr76} and references therein.

Manifolds with corners go back to \cite{Ce61} and \cite{Dou61},
and were also used in \cite{Mr80}; Banach manifolds with corners were
considered in \cite{MO92}.
Manifolds with rough boundary are introduced here as an overarching
framework. Every strongly convex, compact
subset of a Riemannian manifold
is a full submanifold and hence a manifold with rough boundary~\cite{RS18}.

Both our discussion of sprays, and the discussion of tubular neighborhoods
took inspiration from~\cite{La99}. For local additions
in convenient differential calculus, see~\cite[42.4]{KM97}.
\chapter[Basic constructions of infinite-dimensional manifolds]{\hspace*{.2mm}Basic constructions of infinite-dimensional\\ \hspace*{4.9mm}manifolds}\label{chap-manifold-constructions}
In this chapter, we construct major classes of infinite-dimensional manifolds.
Notably,
for each compact smooth manifold $M$
and paracompact finite-dimensional smooth manifold~$N$,
we turn the set $C^k(M,N)$
of $N$-valued $C^k$-maps on~$M$
into a smooth manifold.
More generally, $N$ can be any smooth manifold
admitting a local addition.
The construction is essential for
infinite-dimensional Lie theory:
First, for each Lie group~$K$,
we obtain a Lie group structure on $C^k(M,K)$,
notably on the loop group
$C^k(\bS,K)$ (see Proposition~\ref{gpCkmps}).
Second, the group $\Diff(M)$
of smooth diffeomorphisms of~$M$
turns out to be an open subset of $C^\infty(M,M)$,
and the smooth manifold structure as an open
subset makes it a Lie group
(see Theorem~\ref{diffKsmooth}).
As a starting point,
we define a topology on $C^k(M,N)$
for arbitrary $C^k$-manifolds $M$ and~$N$
modeled on locally convex spaces,
If $N:=F$ is a locally convex space, the latter topology makes $C^k(M,F)$
a locally convex space. Similarly,
taking $N:=E$ for a given vector bundle $E\to M$,
the induced topology on the subset $\Gamma_{C^k}(E)\sub C^k(M,E)$
of $C^k$-sections is  a locally convex vector topology.
We study differentiability properties of non-linear mappings between such spaces of mappings or spaces of sections, continuing the discussions started in Section~\ref{secCspaces}. Notably, we obtain exponential laws
now also for function spaces on products of manifolds (see Section~\ref{secspacemfd}).
On this foundation, we carry out the following constructions:\medskip

\noindent
\emph{Manifold of mappings on compact manifolds.}
The smooth manfold structure on $C^k(M,N)$
announced above is constructed in Section~\ref{sec-mfdmps},
for compact~$M$.\medskip

\noindent
\emph{Fine box products of manifolds.}
For an arbitrary sequence $(M_n)_{n\in \N}$
of smooth manifolds, we construct a smooth manifold structure on the cartesian
product
\[
\prod_{n\in \N}M_n,
\]
endowed with a certain topology which is finer than the product
topology (and finer than the box topology). Such manifolds, 
denoted $\prod_{n\in \N}^{\rm fb} M_n$,
are called \emph{fine box products} (see Section~\ref{sec-box-mfd}).\medskip

\noindent
\emph{Manifolds of mappings on non-compact manifolds.}
For every $\sigma$-compact finite-dimensional
smooth manifold $M$ and smooth manifold $N$ admitting a local addition,
we show that the image of the map
\[
\rho\colon C^k(M,N)\to{\prod_{n\in \N}}^{\fbx}C^k(M_n,N),\quad f\mto (f|_{M_n})_{n\in \N}
\]
is a submanifold, for each locally finite cover $(M_n)_{n\in\N}$
of $M$ by compact full submanifolds~$M_n$. We give $C^k(M,N)$
the smooth manifold structure making $\rho$ a $C^\infty$-diffeomorphism
onto the image; it is independent of the choice of
$(M_n)_{n\in\N}$ (see Section~\ref{sec-mfd-map-sigma}).\medskip

\noindent
\emph{Direct limits of finite-dimensional manifolds.}
We construct a natural smooth manifold structure on
the union $\bigcup_{n\in \N}M_n$ for each
ascending sequence
\[
M_1\sub M_2\sub\cdots
\]
of finite-dimensional $C^\infty$-manifolds such that all inclusion maps
$M_n\to M_{n+1}$ are smooth immersions (see Section~\ref{dl-construction}).\medskip

Readers interested in the corresponding
Lie groups $C^k(M,K)$ of mappings
and
Lie groups $\Diff(M)$ of diffeomorphisms for compact~$M$
can pass directly from Section~\ref{sec-mfdmps}
to Chapters \ref{ch:mapgrp} and \ref{ch:diffeo};
for non-compact~$M$,
Sections~\ref{sec-box-mfd} and \ref{sec-mfd-map-sigma}
are prerequisites.
The construction of the Lie group structure does
not require additional theory from other chapters.\footnote{For the additional information
that well-behaved (regular) Lie groups
are obtained, basic results from Sections \ref{sec:4.3}
and \ref{tools-regularity}
are used.}
Likewise, the construction of the Lie group structure on
an ascending union $\bigcup_{n\in \N}G_n$
of finite-dimensional Lie groups $G_1\sub G_2\sub\cdots$
in Chapter~\ref{ch:dirlim}
only presupposes Section~\ref{dl-construction}.
We already encountered unit groups of continuous inverse
algebras as another basic source of examples of infinite-dimensonal Lie groups.
Also parts of Chapter~\ref{ch:lingrp} can be read
without additional Lie theory. Results concerning the compact-open
$C^k$-topology from the current chapter are useful.
\section{Function spaces on manifolds and spaces of sections}
\label{secspacemfd}
We now endow the set $C^k(M,N)$ of $C^k$-maps between $C^k$-manifolds~$M$ and~$N$
with a topology, the so-called \emph{compact-open $C^k$-topology}.
As a special case, we obtain topologies on spaces of vector-valued $C^k$-functions
on manifolds, and also topologies on spaces of sections in vector bundles.
\begin{numba} {\bf General convention.}
We let $\K\in\{\R,\C\}$
and consider manifolds over the ground field $\K$.
All manifolds considered in Section~\ref{secspacemfd}
are modeled on arbitrary locally convex spaces and may have
a rough boundary, unless the contrary is stated.
\end{numba}
\subsection*{The compact-open {\boldmath$C^k$}-topology on {\boldmath$C^k(M,N)$}}
If $M$ is a $C^k$-manifold, we define its iterated tangent bundles
via $T^0M:=M$ and $T^jM:=T(T^{j-1}M)$ for $j\in\N$ such that $j\leq k$.
\begin{defn}\label{defncoCk}\label{def:smooth-co-top}
Let $M$ and $N$ be $C^k$-manifolds, where $k \in \N_0\cup\{\infty\}$.
We let $C^k(M,N)$ be the set of all $C^k$-maps $f \colon M\to N$,
endowed with the initial topology with respect to the mappings
\[
T^j\colon C^k(M,N)\to C(T^kM,T^kN),\;\, f\mto T^k f ,
\]
where $C(T^kM,T^kN)$ is endowed with the compact-open topology,
$T^0f:=f$, and $T^j f:=T(T^{j-1}f )$ if $j\geq 1$.
\end{defn}
\begin{rem}\label{fiobtopoCk}
(a) As $C(M,N)$ is Hausdorff and has continuous point evaluations
$\ve_x\colon C(M,N)\to N$, $f \mto f(x)$ for all $x\in M$,
also $C^k(M,N)$ is Hausdorff and the point evaluations $\ve_x\colon C^k(M,N)\to N$
are continuous.

(b) If $M$ and $N$ are smooth manifolds, then $C^\infty(M,N)=\bigcap_{k\in\N_0}C^k(M,N)$
and the topology on $C^\infty(M,N)$ is initial with respect to the inclusion maps
$C^\infty(M,N)\to C^k(M,N)$ for $k\in\N_0$ (in view of Lemma~\ref{transinit}).
Hence
\[
C^\infty(M,N)=\pl \, C^k(M,N)
\]
is the projective limit of the projective system $C^0(M,N)\leftarrow C^1(M,N)\leftarrow\cdots$
(all limit maps and all bonding maps are the respective inclusion maps).
\end{rem}
\begin{lem}\label{Cktoppu}
Let $N$, $M$, and $L$ be $C^k$-manifolds with $k\in\N_0\cup\{\infty\}$
and $g \colon M\to L$ be a $C^k$-map. Then the following holds:
\begin{description}[(D)]
\item[\rm(a)]
The map $C^k(g,N)\colon C^k(L,N)\to C^k(M,N)$, $f\mto f\circ g$
is continuous.
\item[\rm(b)]
The map $C^k(N,g)\colon C^k(N,M)\to C^k(N,L)$, $f\mto g\circ f$ is continuous.
\end{description}
\end{lem}
\begin{prf}
(a) holds since
$T^j(f\circ g)=T^jf\circ T^jg=C(T^jg,T^jN)(T^j f)\in C(T^jM,T^jN)$
depends continuously on $f\in C^k(L,N)$ for each $j\in\N_0$ such that $j\leq k$,
by Lemma~\ref{pubas}.

(b) holds since $T^j(g\circ f)=T^jg\circ T^jf=C(T^jN,T^jg)(T^jf)
\in C(T^jN,T^jL)$ depends continuously on $f\in C^k(N,M)$ for each $j\in\N_0$
such that $j\leq k$, by Lemma~\ref{covsuppo}.
\end{prf}
\begin{rem}\label{earlyremCk}
Let $M$ and $N$ be $C^k$-manifolds.

(a) Assume that $S\sub M$ is endowed with a $C^k$-manifold structure
turning the inclusion map $j\colon S\to M$ into a $C^k$-map
(e.g., $S$ could be any open subset of~$M$, or $M$ a manifold without boundary
and $S$ a submanifold of~$M$).
Then the restriction map
\[
C^k(M,N)\to C^k(S,N),\;\, f \mto f|_S
\]
is continuous (as it coincides with $C^k(j,N)$).

(b) If $S\sub N$ is a submanifold,
then the compact-open $C^k$-topology on
$C^k(M,S)$ coincides with the topology induced by $C^k(M,N)$.
This follows from the fact that $T^jS$ carries the topology induced by $T^jN$
(see Lemma~\ref{sub-TM}),
whence $C(T^jM,T^jN)$ induces the compact-open topology on $C(T^jM,T^jS)$,
by Lemma~\ref{ctsemb}, for each $j\in\N_0$ such that $j\leq k$.
Likewise for full submanifolds $S\sub N$.
\end{rem}
\begin{lem}\label{onlyT}
If $k\in\N_0$, then $T^j\colon C^k(M,N)\to C^{k-j}(T^jM,T^jN)$
is a topological embedding for all $j\in \{0,1,\ldots, k\}$.
In particular, the map\linebreak
$T^k\colon C^k(M,N)\to C(T^kM,T^kN)$, $f\mto T^kf$
is a topological embedding and if $k\geq 1$, then $T\colon C^k(M,N)\to C^{k-1}(TM,TN)$,
$f\mto Tf$ is a topological embedding.
\end{lem}
\begin{prf}
It suffices to prove the final assertion.
As a consequence of Definition~\ref{defncoCk},
the map
\[
C^k(M,N)\to C(M,N)\times C^{k-1}(TM,TN),\;\, f\mto (f,Tf)
\]
is a topological embedding (cf.\ Lemma~\ref{transinit}).
Let $\cO$ be the compact-open $C^k$-topology on $C^k(M,N)$ and $\cT$
be the topology on $C^k(M,N)$ turning $T\colon C^k(M,N)\to C^{k-1}(TM,TN)$
into a topological embedding. By the preceding, $\cT\sub \cO$.
If we can show that $\iota\colon (C^k(M,N),\cT)\to C(M,N)$ is continuous,
$\cO\sub \cT$ (and thus $\cO=\cT$) will follow.
Let $\pi_N\colon TN\to N$ be the bundle projection
and $0_M\colon M\to TM$, $x\mto 0\in T_xM$ be the zero-section
(both of which are $C^{k-1}$). Then $\pi_N(Tf(0_M(x)))=\pi_N(0_N(f(x)))=f(x)$
for all $x\in M$ and thus
\[
\iota(f)=f=(C(M,\pi_N)\circ C(0_M,TN))(Tf),
\]
which depends continuously on $f\in (C^k(M,N),\cT)$
(see Lemma~\ref{pubas} and Lemma~\ref{covsuppo}).
\end{prf}
\begin{lem}\label{prodmapsp}
Let $M$, $N_1$, and $N_2$ be $C^k$-manifolds
with $k\in\N_0\cup\{\infty\}$ and $\pr_i\colon N_1\times N_2\to N_i$
be the projection onto the $i$th component for $i\in \{1,2\}$.
Then the map
\[
\Phi:=(C^k(M,\pr_1),C^k(M,\pr_2))\colon C^k(M,N_1\times N_2)\to C^k(M,N_1)\times C^k(M,N_2)
\]
taking $f$ to $(\pr_1\circ f,\pr_2\circ f)$
is a homeomorphism.
\end{lem}
\begin{prf}
In view of Remark~\ref{fiobtopoCk}(b), it suffices to consider
$k\in\N_0$. If $k=0$, then the assertion holds by Lemma~\ref{cotopprod}.
Let $k\in\N$ now and assume that the assertion holds for $k-1$ in place of~$k$.
By Lemma~\ref{mapinprodmfd}, $\Phi$ is a bijection;
it thus suffices to show that $\Phi$ is a topological embedding,
which holds if and only if $(T\times T)\circ \Phi\colon C^k(M,N_1\times N_2)\to C^{k-1}(TM,TN_1)\times
C^{k-1}(TM,TN_2)$ is a topological embedding, by Lemma~\ref{onlyT}.
Now
\[ \psi:=(T\pr_1,T\pr_2)\colon T(N_1\times N_2)\to TN_1\times TN_2 \] 
is a $C^{k-1}$-diffeomorphism. Hence
\begin{eqnarray*}
C^{k-1}(TM,T(N_1\times N_2)) & \sim &  C^{k-1}(TM,TN_1\times TN_2)\\
& \sim & C^{k-1}(TM,TN_1)\times
C^{k-1}(TM,TN_2),
\end{eqnarray*}
where the first homeomorphism is $C^{k-1}(TM,\psi)$ (cf.\ Lemma~\ref{Cktoppu}(b))
and the second homeomorphism is provided by the inductive hypothesis. Let $\Theta\colon
C^{k-1}(TM,T(N_1\times N_2))\to C^{k-1}(TM,TN_1)\times
C^{k-1}(TM,TN_2)$ be the resulting homeomorphism.
Using Lemma~\ref{onlyT}, we deduce that $\Theta\circ T=(T\times T)\circ \Phi$
is a topological embedding. This completes the induction step.
\end{prf}
\begin{lem}
Let $M$, $N_1$, and $N_2$ be $C^k$-manifolds with $k\in\N_0\cup\{\infty\}$
and $g\colon M\times N_1\to N_2$ be a $C^k$-map.
Then 
\[
g_*\colon C^k(M,N_1)\to C^k(M,N_2),\;\, f\mto g\circ (\id_M, f)
\]
is continuous.
\end{lem}
\begin{prf}
Identifying $C^k(M,M\times N_1)$ with $C^k(M,M)\times C^k(M,N_1)$
as a topological space as in Lemma~\ref{prodmapsp},
we have $g_*(f)=C^k(M,g)(\id_M, f)\in C^k(M,N_2)$, which is continuous
in $f\in C^k(M,N_1)$ by Lemma~\ref{Cktoppu}(b).
\end{prf}
\begin{lem}\label{family-restrictions}
Given $k\in\N_0\cup\{\infty\}$, let $M$ and $N$ be $C^k$-manifolds.
Let $(V_i)_{i\in I}$ be a family of full submanifolds of~$M$
whose interiors $V_i^0$ relative~$M$ cover~$M$.
Then
\[
\rho\colon C^k(M,N)\to\prod_{i\in I} C^k(V_i,N),\;\, f\mto (f|_{V_i})_{i\in I}
\]
is a topological embedding with closed image.
\end{lem}
\begin{prf}
The image of $\rho$ is the set of all $(\gamma_i)_{i\in I}\in\prod_{i\in I}C^k(V_i,N)$
such that $\gamma_i(x)=\gamma_j(x)$ for all $i,k\in I$ and $x\in V_i\cap V_k$.
The latter set is closed in $\prod_{i\in I}C^k(V_i,N)$ (cf.\ Remark~\ref{fiobtopoCk}(a)).

It remains to show that the compact-open $C^k$-topology~$\cO$ coincides with the initial topology~$\cT$
with respect to the restriction maps $\rho_i:=C^k(\eta_i,N)\colon C^k(M,N)\to C^k(V_i,N)$
for $i\in I$,
where $\eta_i\colon V_i\to M$ is the inclusion map (which is $C^k$).
As the maps $C^k(\eta_i,N)$ are continuous by Lemma~\ref{Cktoppu}(a),
we have $\cT\sub\cO$. The restriction maps
\[ r_i\colon C^k(V_i,N)\to C^k(V_i^0,N) \] 
being continuous, we have $\cS\sub\cT$ for the initial topology~$\cS$ on $C^k(M,N)$ with respect to
the mappings $r_i\circ \rho_i\colon$ $C^k(M,N)\to C^k(V_i^0,N)$, $\gamma\mto\gamma|_{V_i^0}$.
Hence $\cT=\cO$ will hold if we can show that $\cS=\cO$.
After replacing each $V_i$ with $V_i^0$, we may therefore assume that each $V_i$ is an open submanifold
of~$M$.

For each $j\in \N_0$ such that $j\leq k$,
the sets $T^jV_i$ form an open cover of $T^jM$ for $i\in I$,
whence the compact-open topology on $C(T^jM,T^jN)$ is initial with respect to
the restriction maps $\rho_{i,j}\colon C(T^jM,T^jN)\to C(T^jV_i,T^jN)$ for $i\in I$ (see
Lemma~\ref{coveremb}).
By transitivity of initial topologies (Lemma~\ref{transinit}),
the topology~$\cO$ on $C^k(M,N)$ is initial with respect
to the mappings $\rho_{i,j}\circ T^j$ for $i\in I$ and $j\in\N_0$ with $j\leq k$.
Again by transitivity of initial topologies, the initial topology on $C^k(M,N)$ with respect
to the maps $\rho_{i,j}\circ T^j=T^j\circ\rho_i$ (which is $\cO$ by the preceding)
coincides with the initial topology~$\cT$ with respect to the mappings~$\rho_i$. 
\end{prf}
We record a further simple fact.
\begin{lem}\label{toconstcts}
Let $N$ and~$M$ be $C^k$-manifolds, with $k\in\N_0\cup\{\infty\}$.
Consider the mapping $c\colon N\to C^k(M,N)$ taking $x\in N$
to the constant function $c_x\colon M\to N$, $y\mto x$.
Then~$c$ is continuous.
\end{lem}
\begin{prf}
In view of Remark~\ref{fiobtopoCk}(b), we may assume that $k$ is finite.
The assertion holds for $k=0$, as $c^{-1}(\lfloor K, V\rfloor)$
is open (namely $V$ or $N$) for each compact subset $K\sub M$ and open subset $V\sub N$.
Now assume that the map $C\colon TN\to C^{k-1}(TM,TN)$
taking $v\in TN$ to the constant mapping\linebreak
$C_v\colon TM\to TN$, $w\mto v$ is continuous.
Let $0_N\colon N\to TN$ be the zero section taking $x\in N$ to $0_N(x)=0\in T_xN$.
We know that $0_N$ is~$C^{k-1}$.
As $T(c_x)\in C^{k-1}(TM,TN)$ is the constant map
$C_{0_N(x)}\colon TM\to TN$ for $x\in N$, the map
$T\circ c\colon N\to C^{k-1}(TM,TN)$ concides with $C\circ 0_N$,
whence it is continuous. The map $T\colon C^k(M,N)\to C^{k-1}(TM,TN)$ being a topological
embedding, we deduce that~$c$ is continuous. 
\end{prf}
\subsection*{Spaces of vector-valued $C^k$-functions on manifolds}
We now consider $C^k(M,N)$ in the special case that $N$ is a locally convex topological
vector space.
\begin{prop}
Given $k\in\N_0\cup\{\infty\}$, let $M$ be a $C^k$-manifold
and $F$ be a locally convex space.
Then the compact-open $C^k$-topology $($as in Definition~\emph{\ref{defncoCk}}$)$
makes $C^k(M,F)$ a locally convex topological vector space.
\end{prop}
\begin{prf}
Let $J$ be the set of all $j\in\N_0$ with $j\leq k$.
Now $TF=F\times F$ is a locally convex space and 
$T\colon C^k(M,F)\to C(TM,TF)$, $f \mto (f \circ \pi_M,df)$
is linear (where $\pi_M\colon TM\to M$ is the bundle projection). More generally,
$T^jF$ is a locally convex space (isomorphic to $F^{2^j}$)
for each $j\in J$ and the mapping $T^j\colon C^k(M,F)\to C(T^jM,T^jF)$ is linear,
whence the topological embedding
\[
\Phi\colon C^k(M,F)\to\prod_{j\in J}C(T^jM,T^jF)
\]
is linear. Like the image of $\Phi$, also $C^k(M,F)$ is a locally convex space.
\end{prf}
We hasten to check that the topologies described in
Definitions~\ref{firstdefCktop} and \ref{defncoCk} coincide whenever both are defined.
\begin{lem}\label{df-vs-Tf}
Let $E$ and $F$ be locally convex spaces, $U\sub E$ be a locally convex subset
with dense interior and $k\in\N_0\cup\{\infty\}$.
Then the compact-open $C^k$-topology~$\cO$ on
$C^k(U,F)$ $($as in Definition~\emph{\ref{defncoCk}}$)$ is
initial with respect to the maps $d^{\,(j)}\colon C^k(U,F)\to C(U\times E^j,F)$,
$f \mto d^{\,(j)}f $ for $j\in \N_0$ such that $j\leq k$.
\end{lem}
\begin{prf}
Let $\cT$ be the initial topology on $C^k(U,F)$ with respect to the maps $d^{\,(j)}$ for $j\in\N_0$
such that $j\leq k$.

To verify that $\cT\sub\cO$,
we recursively define continuous linear mappings
\[
\alpha_j\colon F^{2^j}\to F\quad\mbox{and}\quad
\beta_j\colon E^{j+1}\to E^{2^j}
\]
with $\beta_j(U\times E^j)\sub U\times E^{2^j-1}=T^jU$ for $j\in\N$ with $j\leq k$
as follows: We let
$\alpha_1\colon F\times F\to F$, $(v,w)\mto w$ be the projection onto the second component
and $\beta_1\colon E\times E\to E\times E$ be the identity map.
Recursively, we define
\[
\alpha_j:=\alpha_{j-1}\circ\pr_2,
\]
where $\pr_2\colon F^{2^{j-1}}\times F^{2^{j-1}}\to F^{2^{j-1}}$
is the projection onto the second component, and
\[
\beta_j(x,y_1,\ldots,y_j):=(T\beta_{j-1})((x,y_1,\ldots,y_{j-1}),(y_j,0,\cdots,0))
\]
for $x,y_1,\ldots,y_j\in E$. Since $d^{\,(0)}f=f=T^0f$, the topology~$\cO$
makes~$d^{\,(0)}$ continuous.
We claim that
\begin{equation}\label{explformu}
d^{\,(j)}\gamma=\alpha_j\circ (T^j\gamma) \circ\beta_j|_{U\times E^j}
\end{equation}
for all $j\in\N$ with $j\leq k$. If this is true, then
\[
d^{\,(j)}f=
(C(U\times E^j,\alpha_j) \circ C(\beta_j|_{U\times E^j},T^jF))
(T^jf)\in C(U\times E^j,F)
\]
depends continuously on $f\in (C^k(U,F),\cO)$,
and we deduce that $\cT\sub\cO$.
Now
\[
d^{\,(1)}f=df=\pr_2\circ Tf=\alpha_1\circ Tf\circ\beta_1|_U.
\]
Recall that if $\theta\colon X\to Y$ is a continuous linear map between
locally convex spaces, then $d\theta=\theta\circ\pr_2$ using
the projection $\pr_2\colon X\times X\to X$
onto the second factor, and $T\theta=\theta\times\theta$;
notably, $d\theta$ and $T\theta$ are continuous linear.
Hence,
if
$d^{\,(j-1)}f=\alpha_{j-1}\circ (T^{j-1}f)\circ\beta_{j-1}|_{U\times E^{j-1}}$
holds, then
\[
d^{\,(j)}f(x,y_1,\ldots,y_j)=d(d^{\,(j-1)}f)((x,y_1,\ldots,y_{j-1})(y_j,0,\ldots,0))
\]
with
\begin{eqnarray*}
d(d^{\,(j-1)}f) &=&d(\alpha_{j-1}\circ T^{j-1}f\circ\beta_{j-1}|_{U\times E^{j-1}})\\
&=&(d\alpha_{j-1})\circ T^jf\circ (T\beta_{j-1})|_{U\times E^{2j-1}}\\
&=& (\alpha_{j-1}\circ\pr_2)\circ T^jf\circ (T\beta_{j-1})|_{U\times E^{2j-1}},
\end{eqnarray*}
using the Chain Rule.
Thus (\ref{explformu}) also holds for~$j$.

To see that $\cO\sub \cT$ (and thus $\cO=\cT$), note first that
$T^0f=f=d^{\,(0)}f\in C(U,F)$ depends continuously on $f\in
(C^k(U,F),\cT)$. We claim that, for each $j\in\N$ such that $j\leq k$, there exists
$m_j\in\N$, $i_a\in \{1,\ldots, j\}$ for $a\in \{1,\ldots, m_j\}$ and continuous linear mappings
\[
\alpha_a\colon F\to F^{2^j}\;\,\mbox{as well as}\;\,
\beta_a\colon E^{2^j}\to E^{i_a+1}
\]
such that $\beta_a(U\times E^{2^j-1})\sub U\times E^{i_a}$ and
\begin{equation}\label{themonsta}
T^jf=\sum_{a=1}^{m_j}\alpha_a\circ d^{\,(i_a)}f\circ \beta_a|_{T^jU}\;\,
\mbox{for all $f\in C^k(U,F)$.}
\end{equation}
If this is true, then $T^jf$ is a finite sum of contributions of
the form $(C(T^jU,\alpha_a)\circ C(\beta_a,F))  (d^{(i_a)}f)\in C(T^jU,T^jF)$,
which depend continuously on $f\in (C^k(U,F),\cT)$. Thus $\cO\sub\cT$.
The claim holds for $j=1$ since
\[
T^1f=Tf=(f\circ\pr_1,df)
=\alpha_1\circ (d^{\,(0)}f) \circ \beta_1+\alpha_2\circ (d^{\,(1)}f)\circ \beta_2
\]
using the linear maps $\alpha_1\colon F\to F\times F$, $v\mto (v,0)$ and $\alpha_2\colon F\to F\times F$,
$v\mto (0,v)$ as well as
$\beta_1\colon E\times E\to E$, $(x,y)\mto x$ and $\beta_2:=\id_{E\times E}$.
Now assume that $j\in\N$ with $j<k$ and assume that the claim holds for~$j$.
Write $T^jf$ as in~(\ref{themonsta}).
Applying~$T$ to (\ref{themonsta}), the Chain Rule yields
\begin{eqnarray}
T^{j+1}f &=& \sum_{a=1}^{m_j}
(T\alpha_a)\circ T(d^{\,(i_a)}f)\circ T\beta_a\label{stpx}\\
&=&\sum_{a=1}^{m_j}
\big( (\alpha_a,0)\circ (d^{\,(i_a)}f)\circ \pr_1 \circ T\beta_a
+(0,\alpha_a)\circ d(d^{\,(i_a)}f)\circ T\beta_a\big)\notag
\end{eqnarray}
as $T\alpha_a=\alpha_a\times\alpha_a$,
where $\pr_1\colon E^{i_a+1}\times E^{i_a+1}\to E^{i_a+1}$
is the projection onto the first factor.
Now $d(d^{\,(i_a)}f)=d^{\,(1)}f$ if $i_a=0$. If $i_a\geq 1$, we have
\begin{equation}\label{stpxx}
d(d^{\,(i_a)}f)=d^{\,(i_a+1)}f\circ \theta_{i_a}|_{U\times E^{2i_a+1}}
+\sum_{\ell=1}^{i_a}
(d^{\,(i_a)}f)\circ \theta_{i_a,\ell}|_{U\times E^{2i_a+1}},
\end{equation}
where $\theta_{i_a}\colon E^{2i_a+1}\to E^{i_a+2}$,
$(x,y_1,\ldots,y_{i_a},z,w_1,\ldots,w_{i_a})\mto
(x,y_1,\ldots,y_{i_a},z)$
and $\theta_{i_a,\ell}\colon E^{2i_a+2}\to E^{i_a+1}$ is the map given by
\[
(x,y_1,\ldots,y_{i_a},z,w_1,\ldots,w_{i_a})\mto
(x,y_1,\ldots,y_{\ell-1},w_\ell,y_{\ell+1},\ldots,y_{i_a}).
\]
Substituting (\ref{stpxx}) (and $d(d^{\,(0)}f)=d^{\,(1)}f$)
into (\ref{stpx}) and expanding, we obtain a finite sum of terms
as in (\ref{themonsta}) (with $j+1$ in place of~$j$).
If $m_{j+1}$ is the number of summands, we can enumerate them
by $1,\ldots,m_{j+1}$ to complete the induction step.
\end{prf}
\begin{lem}\label{fun-modules}
For each $k\in \N_0\cup\{\infty\}$
and $C^k$-manifold~$M$ modeled on a topological
$\K$-vector space,
the following hold:
\begin{description}[(D)]
\item[\rm(a)]
Using pointwise multiplication,
$C^k(M,\K)$ is an associative unital
topological $\K$-algebra.
\item[\rm(b)]
For each locally convex topological $\K$-vector space~$E$,
pointwise multiplication makes
$C^k(M,E)$ a topological
$C^k(M,\K)$-module.
In particular, for each $f\in C^k(M,\K)$ we get a continuous linear
multiplication operator
\[
m_f\colon C^k(M,E)\to C^k(M,E),\;\;
\gamma\mto f\gamma.
\]
\end{description}
\end{lem}
\begin{prf}
(a) Let $\beta \colon \K\times\K\to \K$, $(s,t)\mto st$ be the multiplication map.
Then $C^k(M,\beta)\colon C^k(M,\K\times \K)\to C^k(M,\K)$
is continuous, by Lemma~\ref{Cktoppu}(b). Identifying the domain with $C^k(M,\K)\times C^k(M,\K)$
as in Lemma~\ref{prodmapsp},
we obtain the multiplication map $C^k(M,\K)\times C^k(M,\K)\to
C^k(M,\K)$.
The proof of (b) follows the same lines,
replacing $\beta$ with the scalar multiplication map \break 
$\K\times E \to E$.
\end{prf}
\begin{lem}\label{space-lin-init}
Let $M$ be a $C^k$-manifold with $k\in\N_0\cup\{\infty\}$
and $F$ be a locally convex space whose topology is initial with respect
to a family $(\lambda_i)_{i\in I}$
of linear mappings $\lambda_i\colon F\to F_i$
to locally convex spaces~$F_i$.
Then the compact-open $C^k$-topology on $C^k(M,F)$ is initial with respect to
the mappings $C^k(M,\lambda_i)\colon C^k(M,F)\to C^k(M,F_i)$
for $i\in I$.
\end{lem}
\begin{prf}
Let $\cO$ be the compact-open $C^k$-topology on $C^k(M,F)$.
For $j\in \N_0$ such that $j\leq k$, the topology on $T^jF=F^{2^{j-1}}$
is initial with respect to the linear maps $T^j\lambda_i=\lambda_i^{2^{j-1}}$,
whence the compact-open topology on $C(T^jM,T^jF)$ is initial with respect to the maps
$C(T^jM,T^j\lambda_i)$ for $i\in I$
(see Lemma~\ref{inipush}). Thus $\cO$ is inital with respect to the maps
$C(T^jM,T^j\lambda_i)\circ T^j$ with $T^j\colon C^k(M,F)\to C(T^jM,T^jF)$.
Since $T^j(\lambda_i\circ f)=T^j\lambda_i\circ T^j f$ for $f\in C^k(M,F)$,
the latter maps coincide with $\tau_{i,j}\circ C^k(M,\lambda_i)$,
writing $\tau_{i,j}\colon C^k(M,F_i)\to C(T^jM,T^jF_i)$, $g\mto T^jg$.
The topology on $C^k(M,F_i)$ being initial with respect to the maps $\tau_{i,j}$
with $k\geq j\in\N_0$, transitivity of initial topologies
shows that $\cO$ is initial with respect to the maps $C^k(M,\lambda_i)$.
\end{prf}
\begin{prop}\label{lem:smooth-eval-vec}\label{prop:smooth-eval}
Let $k\in\N_0\cup\{\infty\}$, $M$ be a locally compact
$C^k$-manifold
and
$F$ be a locally convex space.
Then the evaluation map
\[
\ve\colon C^k(M,F)\times M\to F,\;\; (f,x)\mto f(x)
\]
is $C^{\infty,k}$ and thus $C^k$.
Notably, $\ve$ is smooth if $k=\infty$.
\end{prop}
\begin{prf}
Let $E$ be the modeling space of~$M$.
Given $x_0\in M$, consider a chart $\phi\colon U_\phi\to V_\phi\sub E$ 
of~$M$ around $x_0$.
Considering $\phi^{-1}$ as a $C^k$-map to~$M$,
Lemma~\ref{Cktoppu}(a)
provides a continuous linear map 
\[
(\phi^{-1})^*:=C^k(\phi^{-1},F)\colon C^k(M,F)\to C^k(V_\phi, F).
\]
The evaluation map $\ev\colon C^k(V_\phi,F)\times V_\phi\to F$ is $C^{\infty,k}$
by Lemma~\ref{evaldiffprop}.
Thus
\[
\ve(f,x)=f(x)=(f\circ \phi^{-1})(\phi(x))=\ev((\phi^{-1})^*(f),\phi(x))
\]
is $C^{\infty,k}$ in $(f,x)\in C^k(M,F)\times U_\phi$,
by Proposition~\ref{Cklchainmfd}(a). Being $C^{\infty,k}$ on the sets $C^k(M,F)\times U_\phi$
which form an open cover of $C^k(M,F)\times M$, the function~$\ve$ is~$C^{\infty,k}$.
\end{prf}
\subsection*{{\boldmath$C^{k,\ell}(M\times N,F)$} as a locally convex space; exponential law}
We define a topology on spaces of vector-valued
$C^{k,\ell}$-functions and establish an exponential law.
\begin{defn}\label{defn-ckl-top-mfd}
Let $F$ be  a locally convex space, $M$ be a $C^k$-manifold
and $N$ a $C^\ell$-manifold, with $k,\ell\in \N_0\cup\{\infty\}$.
We endow the vector subspace $C^{k,\ell}(M\times N,F)\sub F^{M\times N}$
of all $C^{k,\ell}$-maps $f\colon M\times N\to F$
with the inital topology with respect to the mappings
\[
(\phi^{-1}\times\psi^{-1})^*\colon C^{k,\ell}(M\times N,F)\to C^{k,\ell}(V_\phi\times V_\psi,F),
\]
for $\phi\colon U_\phi\to V_\phi\sub E_1$
ranging through the maximal $C^k$-atlas of~$M$ and $\psi\colon U_\psi\to V_\psi\sub E_2$
through the maximal $C^\ell$-atlas of~$N$.
\end{defn}
Here $C^{k,\ell}(V_\phi\times, V_\psi, F)$
is endowed with the compact-open $C^{k,\ell}$-topology
as in Definition~\ref{defnCkell}.\\[1mm]
In the situation of Definition~\ref{defn-ckl-top-mfd},
we have:
\begin{lem}\label{embed-prod-ckl}
If $\cA$ is a $C^k$-atlas for~$M$ and $\cB$ a $C^\ell$-atlas
for~$N$,
then
\[
\Phi\colon C^{k,\ell}(M\times N,F)\to\! \prod_{(\phi,\psi)\in \cA\times\cB}
C^{k,\ell}(V_\phi\times V_\psi, F),\;\;\,
f\mto f\circ (\phi^{-1}\times \psi^{-1})
\]
is a linear topological embedding with closed image.
\end{lem}
\begin{prf}
The linearity and injectivity of $\Phi$ are obvious.
The image of $\Phi$ consists of all families
of $C^{k,\ell}$-functions $f_{\phi,\psi}\colon V_\phi\times V_\psi\to F$
such that, for all $\phi_1,\phi_2\in~\cA$, $\psi_1,\psi_2\in~\cB$ and
$(x,y)\in (U_{\phi_1}\cap U_{\phi_2})\times (U_{\psi_1}\cap U_{\psi_2})$,
we have
\[
f_{\phi_1,\psi_1}(\phi_1(x),\psi_1(x))=f_{\phi_2,\psi_2}(\phi_2(x),\psi_2(x)).
\]
As $F$ is Hausdorff and point evaluations are continuous, we deduce that $\im(\Phi)$ is closed.
By definition of the compact-open $C^{k,\ell}$-topology, the map $(\phi^{-1}\times \psi^{-1})^*\colon
C^{k,\ell}(M\times N,F)\to C^{k,\ell}(V_\phi\times V_\psi,F)$
is continuous for all $(\phi,\psi)\in \cA\times \cB$. Hence $\Phi$
is continuous. Let $\tilde{\cA}$ be the maximal $C^k$-atlas containing $\cA$
and $\tilde{\cB}$ be the maximal $C^\ell$-atlas containing~$\cB$.
To conclude that $\Phi$ is a topological embedding, let $\cO$
be the initial topology on $C^{k,\ell}(M\times N,F)$ with respect to $\{\Phi\}$.
It remains to show that $\cO$ makes the mappings
\begin{equation}\label{muma}
(\sigma^{-1}\times\tau^{-1})^*\colon C^{k,\ell}(M\times N,F)\to C^{k,\ell}(V_\sigma\times V_\tau, F)
\end{equation}
continuous for all $(\sigma,\tau)\in \tilde{\cA}\times\tilde{\cB}$.
By Lemma~\ref{Ckl-cover}, the topology on the right-hand side of $(\ref{muma})$
is initial with respect to the restriction maps $\rho_{\sigma,\tau,\phi,\psi}$ to
$C^{k,\ell}(\sigma(U_\sigma\cap U_\phi)\times\tau(U_\tau\cap U_\psi),F)$
for $(\phi,\psi)\in\cA\times\cB$.
By Lemma~\ref{pb-prod-Ckl},
$C^{k,\ell}((\tau\circ \phi^{-1})\times (\sigma\circ\psi^{-1}),F)$
is continuous as a mapping
\[
C^{k,\ell}(V_\phi\times  V_\psi,F)\to
C^{k,\ell}(\sigma(U_\sigma\cap U_\phi)\times\tau(U_\tau\cap U_\psi),F).
\]
Since $\rho_{\sigma,\tau,\phi,\psi}\circ (\sigma^{-1}\times\tau^{-1})^*
=
C^{k,\ell}((\tau\circ \phi^{-1})\times (\sigma\circ\psi^{-1}),F)\circ
(\phi^{-1}\times\psi^{-1})^*$
is continuous, so is
$(\sigma^{-1}\times\tau^{-1})^*$.
\end{prf}
In the next two lemmas, $F$ is a locally convex space
and $k,\ell\in \N_0\cup\{\infty\}$.
\begin{lem}\label{CplusCkl}
For all $C^{k+\ell}$-manifolds $M$ and $N$, the inclusion mapping\linebreak
$j\colon C^{k+\ell}(M\times N,F)\to C^{k,\ell}(M\times N,F)$
is continuous and linear.
\end{lem}
\begin{prf}
If $\phi\colon U_\phi\to V_\phi$ is a chart for $M$ and $\psi\colon U_\psi\to V_\psi$
a chart for $N$, then\linebreak
$f\circ(\phi^{-1}\times\psi^{-1})\colon V_\phi\times V_\psi\to F$
is $C^{k+\ell}$ and thus $C^{k,\ell}$,
by Proposition~\ref{CkvsCkk}(a); hence $f$ is $C^{k,\ell}$.
The inclusion $j_{\phi,\psi}\colon C^{k+\ell}(V_\phi\times V_\psi,F)\to
C^{k,\ell}(V_\phi\times V_\psi,F)$ is continuous and linear
by Proposition~\ref{CkvsCkk}(a).
We now use that the topology on $C^{k,\ell}(M\times N,F)$ is initial
with respect to the mappings 
\[ (\phi^{-1}\times\psi^{-1})^*\colon 
C^{k,\ell}(M\times N,F)\to C^{k,\ell}(V_\phi\times V_\psi,F).\]
Since
\[
(\phi^{-1}\times\psi^{-1})^*\circ j=j_{\phi,\psi}\circ C^{k+\ell}(\phi^{-1}\times \psi^{-1},F)
\]
is continuous by Lemma~\ref{Cktoppu}(a), the linear map $j$ is continuous.
\end{prf}
\begin{lem}\label{fromCkk2Ck}
For all $C^k$-manifolds $M$ and $N$, the inclusion mapping\linebreak
$j\colon C^{k,k}(M\times N,F)\to C^k(M\times N,F)$
is continuous and linear.
\end{lem}
\begin{prf}
If $\phi\colon U_\phi\to V_\phi$ is a chart for $M$ and $\psi\colon U_\psi\to V_\psi$
a chart for $N$, then $f\circ(\phi^{-1}\times\psi^{-1})\colon V_\phi\times V_\psi\to F$
is $C^{k,k}$ and thus $C^k$,
by Proposition~\ref{CkvsCkk}(b); hence $f$ is~$C^k$.
The inclusion map $j_{\phi,\psi}\colon C^{k,k}(V_\phi\times V_\psi,F)\to
C^k(V_\phi\times V_\psi,F)$ is continuous and linear
by Proposition~\ref{CkvsCkk}(b).
By definition of the compact-open $C^{k,k}$-topology,
the map
\[
(\phi^{-1}\times\psi^{-1})^*\colon
C^{k,k}(M\times N,F)\to C^{k,k}(V_\phi\times V_\psi,F)
\]
is continuous.
We now use that the topology on $C^k(M\times N,F)$ is initial
with respect to the maps $C^k(\phi^{-1}\times\psi^{-1},F)\colon
C^k(M\times N,F)\to C^{k,\ell}(V_\phi\times V_\psi,F)$,
as a consequence of Lemmas~\ref{Cktoppu}(a) and \ref{family-restrictions}.
Since
\[
C^k(\phi^{-1}\times\psi^{-1},F)\circ j=j_{\phi,\psi}\circ
(\phi^{-1}\times \psi^{-1})^*
\]
is continuous, the linear map $j$ is continuous.
\end{prf}
\begin{numba}\label{infinfinf}
Notably, $C^{\infty,\infty}(M\times N,F)=C^\infty(M\times N,F)$ as a topological vector space
for each locally convex space $F$ and all smooth manifolds $M$ and~$N$,
by Lemmas~\ref{CplusCkl} and~\ref{fromCkk2Ck}.
\end{numba}
\begin{prop}\label{explawfuncmfd}\label{prop:cartes-closed}
Let $M$ be a $C^k$-manifold
and $N$ be a $C^\ell$-manifold with $k,\ell\in \N_0\cup\{\infty\}$.
Let $E_1$ and $E_2$ be the modeling spaces of $M$ and $N$, respectively,
and $F$ be a locally convex space.
Then the following holds:
\begin{description}[(D)]
\item[\rm(a)]
If $f\colon M\times N\to F$ is $C^{k,\ell}$, then $f^\vee(x):=f_x:=f(x,\cdot)\in C^\ell(N,F)$
for each $x\in M$ and the map $f^\vee\colon M\to C^\ell(N,F)$
is~$C^k$.
\item[\rm(b)]
The map $\Phi\colon C^{k,\ell}(M\times N,F)\to C^k(M,C^\ell(N,F))$, $f\mto f^\vee$
is linear and a topological embedding.
\item[\rm(c)]
If $N$ is locally compact or both $E_1$ and $E_2$ are metrizable,
then the map $\Phi$ in {\rm(b)} is an isomorphism of topological vector spaces.
Moreover, a map $g\colon M\to C^\ell(N,F)$ then is $C^k$ if and only if
the map $g^\wedge\colon M\times N\to F$ defined via
$g^\wedge(x,y):=g(x)(y)$ is~$C^{k,\ell}$.
\item[\rm(d)]
For $k=\ell=\infty$, the map $\Phi\colon C^\infty(M\times N,F)\to C^\infty(M,C^\infty(N,F))$
is linear and a topological embedding.
If $N$ is locally compact or $E_1$ and $E_2$ are metrizable,
then $\Phi$ is an isomorphism of topological vector spaces.
\end{description}
\end{prop}
\begin{prf}
(a)
In view of Definition~\ref{defCklmfd} in terms of local charts,
Lemma~\ref{Ckellpartial} implies that
$f^\vee(x)\colon N\to F$ is a $C^\ell$-map.
Let $\cA$ be an atlas of charts of~$M$ and $\cB$ be an atlas for~$N$.
For charts
\[ \phi\colon U_\phi\to V_\phi\sub E_1\quad \mbox{ in } \quad \cA
  \quad \mbox{ and } \quad \psi\colon U_\psi\to V_\psi\sub E_2
  \quad \mbox{ in } \quad \cB,\]
the map $f_{\phi,\psi}:=f\circ
(\phi^{-1}\times\psi^{-1})\colon V_\phi\times V_\psi\to F$
is $C^{k,\ell}$, whence the map $f_{\phi,\psi}^\vee\colon V_\phi\to C^\ell(V_\psi,F)$
is $C^k$. Then also $\psi^*\circ f_{\phi,\psi}^\vee\circ \phi$
is $C^k$, using the isomorphism $\psi^*\colon C^\ell(V_\phi,F)\to C^\ell(U_\phi,F)$
of topological vector spaces (cf.\ Lemma~\ref{Cktoppu}(b)).
By Lemma~\ref{family-restrictions}, the map
\[
\rho\colon C^\ell(N,F)\to\prod_{\psi\in \cB}C^\ell(U_\psi,F), \;\; g\mto (g|_{U_\psi})_{\psi\in \cB}
\]
is a linear topological embedding with closed image.
As $\psi^*\circ f_{\phi,\psi}^\vee\circ\phi$ is equal to the $\psi$-component
of $\rho\circ f^\vee$, we deduce that $f^\vee$ is $C^k$.

(b) We readily check that $\Phi$ is linear and injective.
The compact-open $C^{k,\ell}$-topology $\cO$
on $C^{k,\ell}(M\times N,F)$
is initial with respect to the mappings
$(\phi^{-1}\times\psi^{-1})^*\colon C^{k,\ell}(M\times N,F)\to C^{k,\ell}(V_\phi\times V_\psi,F)$,
$f\mto f\circ (\phi^{-1}\circ \psi^{-1})$
for $(\phi,\psi)\in\cA\times\cB$, by Lemma~\ref{embed-prod-ckl}
and hence with respect to the mappings
\[
\Theta_{\phi,\psi}\circ (\phi^{-1}\times\psi^{-1})^*,
\]
using the topological embeddings
\[ \Theta_{\phi,\psi}\colon C^{k,\ell}(V_\phi\times V_\psi, F)
\to C^k(V_\phi, C^\ell(V_\psi, F)), \quad g\mto g^\vee \] 
as in Theorem~\ref{explawCkell}. 
Since
\[
C^k(\phi^{-1},C^\ell(N,F))\colon C^k(U_\phi,C^\ell(N,F))\to C^k(V_\phi,C^\ell(N,F))
\]
is an isomorphism of topological vector spaces (cf.\ Lemma~\ref{Cktoppu}(a))
and the linear map $C^k(M,C^\ell(N,F))\to\prod_{\phi\in \cA}C^k(U_\phi,C^\ell(N,F))$,
$g\mto (g|_{U_\phi})_{\phi\in\cA}$ is a topological embedding
(by Lemma~\ref{family-restrictions}),
the topology on $C^k(M,C^k(N,F))$ is initial with respect
to the linear mappings
\[
(\phi^{-1})^*\colon C^k(M,C^\ell(N,F))\to C^k(V_\phi, C^k(N,F)),\;\;
g\mto g\circ \phi^{-1}
\]
for $\phi\in \cA$.
Likewise,
the topology on $C^\ell(N,F)$ is initial with respect to the linear
mappings $(\psi^{-1})^*\colon C^\ell(N,F)\to C^\ell(V_\psi,F)$ for $\psi\in\cB$,
entailing that the topology on $C^k(V_\phi,C^\ell(N,F))$
is initial with respect to the mappings
$C^k(V_\phi,(\psi^{-1})^*)\colon C^k(V_\phi,C^\ell(N,F))\to C^k(V_\phi, C^\ell(V_\psi,F))$
(see Lemma~\ref{space-lin-init}).
As a consequence, the topology $\cT$ on $C^{k,\ell}(M\times N,F)$
making $\Phi$ a topological embedding is initial with respect to the maps
\begin{equation}\label{thusO=T}
C^k(V_\phi,(\psi^{-1})^*)\circ (\phi^{-1})^*\circ \Phi
=\Theta_{\phi,\psi}\circ (\phi^{-1}\times\psi^{-1})^*
\end{equation}
for $(\phi,\psi)\in \cA\times\cB$. Thus $\cT=\cO$.

(c) If $g$ is $C^k$, then also $g_{\phi,\psi}:=C^\ell(\psi^{-1},F)\circ g\circ \phi^{-1}$
for all $(\phi,\psi)\in\cA\times \cB$.
Then $g^\wedge\circ (\phi^{-1}\times\psi^{-1})=(g_{\phi,\psi})^\wedge$ is $C^{k,\ell}$
by Theorem~\ref{explawCkell}, whence $g^\wedge$ is $C^{k,\ell}$.
Notably, $g=(g^\wedge)^\vee=\Phi(g^\wedge)\in \im(\Phi)$,
whence $\Phi$ is surjective and hence an isomorphism of topological vector spaces.
If, conversely, $g^\wedge$ is $C^{k,\ell}$,
then $g=(g^\wedge)^\vee$ is $C^k$ by~(a).

(d) Is a special case of (b) and (c),
as $C^\infty(M\times N,F)=C^{\infty,\infty}(M\times N,F)$ by \ref{infinfinf}.
\end{prf}
\begin{rem}
The conclusions of~(c) and the final conclusion of~(d)
remain valid if $E_1$ and $E_2$ are $k_\omega$-spaces
and $M$ and $N$ admit atlases of charts whose ranges $V_\phi$ and $V_\psi$
are open subsets of closed subsets $A_\phi$ and $A_\psi$ of $E_1$ and $E_2$,
respectively.\\[1mm]
In fact, by part (a) of Exercise~\ref{exc-k-omeg},
the closed subsets $A_\phi$ and $A_\psi$ are \break $k_\omega$-spaces.
By part~(b) of the exercise, each point $x$ in the open subset $V_\phi\sub A_\phi$
has an open neighborhood $W_x$ which is a $k_\omega$-space.
Then $\phi$ restricts to a chart $\phi^{-1}(W_x)\to W_x$
whose range is a $k_\omega$-space.
Likewise for $\psi$. Thus, we may assume that all ranges of the charts
are $k_\omega$-spaces.
The proof of Proposition~\ref{exc-k-omeg}
now carries over if we work with atlases $\cA$ and $\cB$ of such charts,
exploiting that finite direct products of $k_\omega$-spaces
are $k_\omega$-spaces (by Remark~\ref{prod-of-k-omeg})
and hence $k$-spaces.
\end{rem}
We record typical applications.
\begin{prop}\label{pushf-cp}
Let
$k,\ell\in \N_0\cup\{\infty\}$
and $M$ be a
compact $C^{k+\ell}$-manifold
over $\K\in \{\R,\C\}$.
Let $E$ and $F$ be locally convex
topological $\K$-vector spaces, $U\sub E$
be an open subset, and
$g\colon M\times U\to F$
a $C^{k+\ell}$-map.
Then $g_*\colon C^\ell(M,U)\to C^\ell(M,F)$, $f\mto g\circ (\id_M,f)$
is~$C^k$.
\end{prop}
\begin{prf}
The map $(g_*)^\wedge\colon C^\ell(M,U)\times M\to F$,
$(f,x)\mto g(x,f(x))=g(x,\ve(f,x))$
is $C^{k,\ell}$.
In fact, $g$ is $C^{k+\ell}$,
the evaluation map $\ve \colon C^\ell(M,U)\times M\to U$
is $C^{k,\ell}$ (by Proposition~\ref{lem:smooth-eval-vec}),
and the map $C^\ell(M,U)\times M\to M$, $(f,x)\mto x$ is $C^{k,\ell}$
(being~$C^{k+\ell}$).
Thus
Proposition~\ref{Cklchainmfd}(b) applies.
Now $g_*$ is~$C^k$, by Proposition~\ref{explawfuncmfd}(a).
\end{prf}
\begin{prop}\label{superpo-cp}
Let $k\in \N_0\cup\{\infty,\omega\}$,
$\ell\in \N_0\cup\{\infty\}$
and $M$ be a
compact $C^\ell$-manifold
over~$\K$. If $\K=\C$, let $\bL:=\C$;
if $\K=\R$, let
$\bL\in\{\R,\C\}$. Let $E$ and $F$ be locally convex
topological $\bL$-vector spaces
and $U\sub E$
be an open subset.
Then the following holds:
\begin{description}[(D)]
\item[\rm(a)]
If $g\colon U\to F$ is a $C^{k+\ell}_{\bL}$-map,
then the map
\[ C^\ell(M,g)\colon C^\ell(M,U)\to C^\ell(M,F),\quad f\mto g\circ f \] 
is $C^k_{\bL}$.
\item[\rm(b)]
If $g\colon M\times U\to F$
is of the form $g=h\circ (q\times \id_U)$
for some $C^\ell_\K$-map $q\colon M\to N$ to a $C^{k+\ell}_{\bL}$-manifold~$N$
without boundary and a $C^{k+\ell}_{\bL}$-map $h\colon N\times U\to F$,
then $g_*\colon C^\ell(M,U)\to C^\ell(M,F)$, $f\mto g\circ (\id_M,f)$
is~$C^k_{\bL}$.
\end{description}
\end{prop}
\begin{prf}
(a) The map 
\[ C^\ell(M,g)^\wedge\colon C^\ell(M,U)\times M\to F, \quad 
(f,x)\mto g(f(x))=g(\ve(f,x))\] 
is $C^{k,\ell}_\K$, as $g$ is $C^{k+\ell}_\K$
and the evaluation map $\ve \colon C^\ell(M,U)\times M\to U$
is $C^{k,\ell}_\K$ (by Proposition~\ref{lem:smooth-eval-vec}), whence
Proposition~\ref{Cklchainmfd}(b) applies.
Hence $C^\ell(M,g)$ is $C^k_\K$, by Proposition~\ref{explawfuncmfd}(a).
It remains to show that $C^\ell(M,g)$ is $C^k_\C$
if $\K=\R$ and $\bL=\C$. If $k\geq 1$, let $x\in M$
and consider the continuous linear point evaluation
$\ve_x\colon C^\ell(M,F)\to F$, $f\mto f(x)$.
For $f\in C^\ell(M,U)$
and $h\in C^\ell(M,E)$,
\[
\ve_x\circ (dC^\ell(M,g)(f,h))=d(\ve_x\circ C^\ell(M,g))(f,h)=
dg(f(x),h(x))
\]
is $\C$-linear in $h$, whence so is $d C^\ell(M,g)(f,h)$.
The $C^k_\R$-map $C^\ell(M,g)$
therefore is $C^k_\C$, by Lemma~\ref{realvscxearly}
(and this trivially holds also if $k=0$).

(b) Let $H$ be the modeling space of~$N$.
For $x\in M$, let $\phi_x\colon U_x\to V_x\sub H$ be a chart of~$N$
around~$q(x)$. Then there exists a compact full submanifold {$M_x\sub M$}
such that $x\in M_x^0$ and $q(M_x)\sub U_x$.
Then $M=\bigcup_{x\in \Phi}M_x^0=M$ for some finite subset $\Phi\sub M$.
Let $\rho_x\colon C^\ell(M,F)\to C^\ell(M_x,F)$
be the restriction map for $x\in \Phi$. The map
\[
(\rho_x)_{x\in \Phi}\colon C^\ell(M,F)\to\prod_{x\in \Phi}C^\ell(M_x,F)
\]
being a $\bL$-linear topological embedding with closed image,
it suffices to show that $\rho_x\circ g_*\colon C^\ell(M,U)\to C^\ell(M_x,F)$
is $C^k_\bL$ for each $x\in \Phi$.
Identify $C^\ell(M_x,H\times E)$ with $C^\ell(M_x,H)\times C^\ell(M_x,E)$
as in Lemma~\ref{prodmapsp}.
Now $a_x\colon \colon V_x\times U \to F$, $(z,y)\mto h(\phi_x^{-1}(z),y)$
is $C^{k+\ell}_\bL$ and
\[
(\rho_x\circ g_*)(f)=C^\ell(M_x,a_x)(\phi\circ q|_{M_x},f|_{M_x})
\]
is $C^k_\bL$ in $f$, as $C^\ell(M_x,a_x)\colon C^\ell(M_x,V_x\times U)\to
C^\ell(M_x,F)$ is $C^k_\bL$ by (a)
and the restriction map $C^\ell(M,U)\to C^\ell(M_x,U)$
is $C^\infty_\bL$.
\end{prf}
Properties of the function spaces $C^k(M,F)$
are subsumed by Proposition~\ref{properties-spaces-sections} below devoted to spaces of sections
in vector bundles (which include function spaces in
the special case of trivial bundles, by Corollary~\ref{sectrivializable}).
\subsection*{Spaces of sections in vector bundles}

\begin{defn}
Let $\pi\colon E\to M$ be a $C^k$-vector bundle with typical fibre~$F$.
\begin{description}[(D)]
\item[(a)]
Write $\Gamma_{C^k}(E)$ for the set of all
$C^k$-sections $\sigma\colon M\to E$.
We give $C^k(M,E)$ the compact-open $C^k$-topology and
endow $\Gamma_{C^k}(E)$ with the topology induced
by $C^k(M,E)$.
\item[(b)]
If $\sigma,\tau\in\Gamma^{C^k}(E)$, we define a map $\sigma+\tau\colon M\to E$ via
\[
(\sigma+\tau)(x):=\sigma(x)+\tau(x),
\]
using addition in the vector space~$E_x$. Since $(\sigma+\tau)_\theta=\sigma_\theta+\tau_\theta
\in C^k(U,F)$ for each local trivialization $\theta\colon E|_U\to U\times F$,
the map $\sigma+\tau$ is $C^k$ and hence a $C^k$-section of~$E$
(see Remark~\ref{firstremvbdl}(b)). Likewise, $(s\sigma)(x):=s\sigma(x)\in E_x$
defines a $C^k$-section $s\sigma$ of~$E$, for each $\sigma\in\Gamma_{C^k}(E)$
and $s\in\K$.
\end{description}
\end{defn}
\begin{prop}\label{secsviafctns}
Let $\pi\colon E\to M$ be a $C^k$-vector bundle with typical fibre~$F$
and $(\theta_i)_{i\in I}$ be a family of local trivializations $\theta_i\colon E|_{U_i}\to U_i\times F$
of~$E$ such that $\bigcup_{i\in I}U_i=M$. Then the mapping
\[
\Phi\colon \Gamma_{C^k}(E)\to\prod_{i\in I} C^k(U_i,F),\;\, \sigma\mto (\sigma_{\theta_i})_{i\in I}
\]
is linear and a topological embedding with closed image \break $($where
$\sigma_{\theta_i}={\pr_2\circ \theta_i\circ\sigma|_{U_i}}$ is as in Remark~\emph{\ref{firstremvbdl}(b)).}
\end{prop}
\begin{prf}
By construction of the vector space structure on~$\Gamma_{C^k}(E)$,
each of the maps $\sigma\mto \sigma_{\theta_i}$ is linear, and hence also~$\Phi$.
Lemma~\ref{Cktoppu}(a) and Remark~\ref{earlyremCk}(b)
entail that the map $\rho\colon\Gamma_{C^k}(E)\to\prod_{i\in I}\Gamma(E|_{U_i})$,
$\sigma\mto\sigma|_{U_i}$ is a topological embedding.
Now
\[
\psi_i\colon \Gamma_{C^k}(E|_{U_i})\to C^k(U_i,F),\;\,\sigma\mto (\pr_2\circ \, \theta_i)\circ\sigma
\]
is a continuous bijection (by Lemma~\ref{Cktoppu}(b))
with inverse $f \mto \theta_i^{-1}\circ (\id_{U_i},f)$,
which is continuous by Lemmas~\ref{Cktoppu}(b) and \ref{prodmapsp}.
Thus $\Phi=(\prod_{i\in I}\psi_i)\circ\rho$ is a topological embedding, being
the composition of a homeomorphism and a topological embedding.
Given $i,j\in I$, the map
\[
g_{ij}\colon U\times F\to F,\;\, (x,y)\mto
(\pr_2\circ \, \theta_i\circ\theta_j^{-1})(x,y)
\]
is~$C^k$ and $g_{ij}(x,\cdot)\colon F\to F$ is linear.
The image of~$\Phi$ coincides with the closed vector subspace
\[
\Big\{(f_i)_{i\in I}\in\prod_{i\in I}C^k(U_i,F)\colon (\forall i,j\in I)\,(\forall x\in U_i\cap U_j)\,
f_i(x)=g_{ij}(x,f_j(x))\Big\}
\]
of $\prod_{i\in I}C^k(U_i,F)$, as it is contained in it by direct verification.
If, conversely, $(f_i)_{i\in I}$ is an element of the latter vector subspace,
then $\sigma(x):=\theta_j^{-1}(x,f_j(x))$ for $x\in U_j$
is well-defined (as $\theta_i(\theta_j^{-1}(x,\sigma_j(x)))
=(x,g_{ij}(x,\sigma_j(x)))=(x,\sigma_i(x))$
if also $x\in U_i$), and one easily checks that $\sigma\in\Gamma_{C^k}(E)$ and $\Phi(\sigma)=
(f_i)_{i\in I}$.
\end{prf}
Taking $J$ as a singleton, we deduce:
\begin{cor}\label{sectrivializable}
If $\pi\colon E\to M$ is a trivializable $C^k$-vector bundle with typical fibre $f$
and $\theta=(\pi,\theta_2)\colon E\to M\times F$ a trivialization, then
\[
\Gamma_{C^k}(E)\to C^k(M,F),\;\,\sigma\mto \sigma_\theta:=\theta_2\circ\sigma
\]
is an isomorphism of topological vector spaces.\,\qed
\end{cor}
\begin{lem}\label{family-sections}
Let $\pi\colon E\to M$ be a $C^k$-vector bundles
with typical fibre~$F$ and $k\in\N_0\cup\{\infty\}$.
Let $(M_i)_{i\in I}$ be a family of full submanifolds
$M_i$ of~$M$ whose interiors relative~$M$
cover~$M$.
Then
\[
\rho\colon \Gamma_{C^k}(E)\to\prod_{i\in I}\Gamma_{C^k}(E|_{M_i}),\;\;
\sigma\mto (\sigma|_{M_i})_{i\in I}
\]
is linear and a topological embedding with closed image.
If $E|_{M_i}$ is trivializable and $\theta_i\colon E|_{M_i}\to M_i\times F$
a trivialization, then also the following linear map is a topological
embedding with closed image:
\[
\Psi\colon \Gamma_{C^k}(E)\to\prod_{i\in I}C^k(M_i,F),\;\, \sigma\mto
\pr_2\circ \, \theta_i\circ\sigma|_{U_i}.
\]
\end{lem}
\begin{prf}
The map $C^k(M,E)\to \prod_{i\in I}C^k(M_i,E)$,
$f\mto (f|_{M_i})_{i\in I}$ is a topological embedding
with closed image,
by Lemma~\ref{family-restrictions}.
It therefore restricts to
a topological embedding of~$\Gamma_{C^k}(E)$.
Since $\Gamma_{C^k}(E)$ is closed in $C^k(M,E)$,
the image is closed in $\prod_{i\in I}C^k(M_i,E)$,
hence also closed in the subset $\prod_{i\in I}C^k(M,E|_{M_i})$
(on which it induces the product of the compact-open $C^k$-topologies
by Lemma~\ref{restrict-vbdl}
and Remark~\ref{fiobtopoCk}(b))
and its subset $\prod_{i\in I}\Gamma_{C^k}(E|_{M_i})$.

To establish the final assertion, note that $h_i\colon \Gamma_{C^k}(E|_{M_i})\to C^k(M_i,F)$,
$\tau\mto \pr_2\circ\, \theta_i\circ \tau$
is an isomorphism of topological vector spaces, by Corollary~\ref{sectrivializable}.
Since $\rho$ is a linear topological embedding with closed image, so is
$\Psi=\big(\prod_{i\in I} h_i\big)\circ \rho$.
\end{prf}
\begin{prop}\label{properties-spaces-sections}
Let $k\in \N_0\cup\{\infty\}$ and $\pi\colon E\to M$
be a $C^k$-vector bundle over a locally compact $C^k$-manifold~$M$,
with typical fibre~$F$.
\begin{description}[(D)]
\item[\rm(a)]
If $F$ is complete,
quasi-complete, sequentially complete, and Mackey-complete,
respectively, then also $\Gamma_{C^k}(E)$
is complete, quasi-complete, sequentially complete,
resp., Mackey complete.
\item[\rm(b)]
If $M$ is $\sigma$-compact and $F$ is metrizable,
then $\Gamma_{C^k}(E)$ is metrizable.
\item[\rm(c)]
If $M$ is compact, $F$ is normable and $k\in\N_0$,
then $\Gamma_{C^k}(E)$ is normable.
\end{description}
\end{prop}
\begin{prf}
Let $Z$ be the modeling space of~$M$,
which has finite dimension.
Let $(\theta_i)_{i\in I}$
be a family of local trivializations $\theta_i\colon E|_{U_i}\to U_i\times F$
of~$M$ such that $(U_i)_{i\in I}$
is an open cover of~$M$
and there exists a chart $\phi_i\colon U_i\to V_i\sub Z$;
let $\Phi$ be as in Lemma~\ref{family-sections}.
The map $C^k(\phi_i,F)\colon C^k(V_i,F)\to C^k(U_i,F)$
is an isomorphism of topological vector spaces
(cf.\ Proposition~\ref{Cktoppu}(a)).

(a)
By Proposition~\ref{compl-and-metr}(a),
$C^k(V_i,F)$ inherits the respective completeness
property from~$F$,
hence also $C^k(U_i,F)$, the direct product $\prod_{i\in I}C^k(U_i,F)$,
its closed vector subspace $\im(\Phi)$,
and $\Gamma_{C^k}(E)$.

(b) We write~$M$ as a countable union of compact sets.
As each of these is covered by finitely many $U_i$,
we find a countable subset $J\sub I$ such that
$M=\bigcup_{i\in J}U_i$. After replacing $I$ with~$J$,
we may assume that~$I$ is countable.
Then $\prod_{i\in I}C^k(U_i,F)$ is metrizable
as a countable product of metrizable spaces (by Proposition~\ref{compl-and-metr}(b)),
whence so are $\im(\Phi)$ and $\Gamma_{C^k}(E)$.

(c) If $M$ is compact, then each element $x\in M$
is contained in $U_{i(x)}$ for some $i(x) \in I$;
we can find a compact convex $\phi_{i(x)}(x)$-neighborhood
$W_x\sub V_{i(x)}$ and obtain a compact, full submanifold $M_x:=\phi_{i(x)}^{-1}(W_x)$
which contains $x$ in its interior $M_x^0$ relative~$M$.
Then $M=\bigcup_{x\in  A}M_x^0$ for some finite subset
$A\sub M$.
By Lemma~\ref{family-sections},
we have a linear topological embedding
\[
\Psi\colon \Gamma_{C^k}(E)\to\prod_{x\in A}C^k(M_x,F)
\]
taking $\sigma$ to $(\pr_2\circ \, \theta_{i(x)}\circ \sigma|_{M_x})_{x\in A}$.
Now $C^k(W_x,F)$
is normable by Proposition~\ref{compl-and-metr}(b),
whence so are $C^k(M_x,F)$, $\prod_{x\in A}C^k(M_x,F)$, and $\Gamma_{C^k}(E)$.
\end{prf}
\subsection*{Spaces of sections in pullback-bundles}
\noindent
Spaces of sections of pullback bundles
are of particular interest with a view towards manifolds
of mappings. In Lemma~\ref{sections-pb}
through
Proposition~\ref{explaw-pb},
we consider the following setting:
\begin{numba}\label{def-gamma-f}
Let $\pi\colon E\to N$ be a smooth vector bundle
over a smooth manifold~$N$ without boundary, with typical fibre~$F$.
Let $M$ be a $C^\ell$-manifold for some $\ell\in \N_0\cup\{\infty\}$ (possibly with rough boundary)
and $f\colon M\to N$ be a $C^\ell$-map.
We consider the pullback bundle $f^*(E)$
as in Lemma~\ref{pullbackbdl},
with bundle projection $\pi_{f^*(E)}\colon f^*(E)\to M$.
Recall from Exercise~\ref{pullback-is-submfd}
that $f^*(E)$ is a submanifold
of the $C^\ell$-manifold $M\times E$.
We let $\pr_2\colon M\times E\to E$ be the projection onto
the second factor.
Then
\[
\Gamma_f:=\{\tau\in C^\ell(M,E)\colon \pi\circ \tau=f\}
\]
is a vector subspace of $\prod_{x\in M}E_{f(x)}$.
We give $\Gamma_f$ the topology induced by $C^\ell(M,E)$,
endowed with the compact-open $C^\ell$-topology.
\end{numba}
\begin{lem}\label{sections-pb}
$\Gamma_f$ is a locally convex space.
The map $\Phi\colon \Gamma_{C^\ell}(f^*(E))\to \Gamma_f$, $\sigma\mto  \pr_2\circ \, \sigma$
is an isomorphism of topological vector spaces.
\end{lem}
\begin{prf}
The map $C^\ell(M,\pr_2)\colon C^\ell(M,M\times E)\to C^\ell(M,E)$
is continuous by Lemma~\ref{Cktoppu}(b). The also its restriction to $\Gamma_{C^\ell}(f^*(E))$
is continuous (using Remark~\ref{earlyremCk}(b)) and also the co-restriction $\Phi$ of the latter
map to $\Gamma_f$. We readily check that
$\Gamma_f\to\Gamma_{C^\ell}(f^*(E))$, $\tau\mto (\id_M,\tau)$
is the inverse map of~$\Phi$.
This mapping is continuous as it is continuous to $C^\ell(M,M\times E)\sim C^\ell(M,M)\times
C^\ell(M,E)$. Thus $\Phi$ is a homeomorphism.
As $\Phi$ is linear and $\Gamma_{C^\ell}(f^*(E))$
is a locally convex space, $\Gamma_f$
is a locally convex space.
\end{prf}
%
%
\begin{prop}\label{eval-pb}
If $M$ is locally compact, then the evaluation mapping\linebreak
$\ve\colon \Gamma_f\times M\to E$,
$(\tau,x)\mto \tau(x)$ is $C^{\infty,\ell}$.
\end{prop}
\begin{prf}
Given $x_0\in N$,
there exists a local trivialization $\theta\colon E|_U\to U\times F$
of~$E$ for some open $f(x_0)$-neighborhood $U\sub N$.
Let $\pr_F\colon U\times F\to F$
be the projection onto the second factor and $\theta_2:=\pr_F\circ \theta$.
Then $f^{-1}(U)$ is an open $x$-neighborhood in~$M$;
we find a chart $\phi\colon U_\phi\to V_\phi$
of $M$ around $x_0$ such that $U_\phi\sub f^{-1}(U)$.
We show that $\ve$ is $C^{\infty,\ell}$
on the open subset $\Gamma_f\times U_\phi$ of $\Gamma_f\times M$.
It suffices to show that $(\theta\circ \ve)(\tau,\phi^{-1}(x))$
is $C^{k,\ell}$ in $(\tau,x)\in \Gamma_f\times V_\phi$.
By Lemma~\ref{evaldiffprop}, the evaluation map $\ve\colon C^\ell(V_\phi,F)\times V_\phi\to F$
is $C^{\infty,\ell}$.
Now
\begin{eqnarray*}
\theta(\ve(\tau,\phi^{-1}(x))) &=& \theta(\tau(\phi^{-1}(x)))=
(f(\phi^{-1}(x)),\theta_2(\tau(\phi^{-1}(x))))\\
& =& (f(\phi^{-1}(x)),\ev(\theta_2\circ \tau\circ \phi^{-1},x)).
\end{eqnarray*}
The map $h\colon \Gamma_f\to C^\ell(V_\phi,F)$, $\tau\mto \theta_2\circ \tau\circ\phi^{-1}$
is continuous by (a) and (b) in~Lemma~\ref{Cktoppu},
and $h$ is linear. By Proposition~\ref{Cklchainmfd}(a), $\ev\circ (h\times \id_{V_\phi})$
is $C^{\infty,\ell}$. The map $\Gamma_f\times V_\phi\to M$,
$(\tau,x)\mto f(\phi^{-1}(x))$ is $C^\ell$
and thus $C^{\infty,\ell}$, being independent of~$\tau$.
Using Remark~\ref{rem-simpl-Ckl}(e), $\theta\circ \ve\circ (\id_{\Gamma_f}\circ \phi^{-1})$
is $C^{k,\ell}$.
\end{prf}
\begin{rem}
Note that, since $f^*(E)$ merely is a $C^\ell$-manifold,
for $\ell<\infty$ an analog of Proposition~\ref{eval-pb}
does not make sense for the evaluation map
$\Gamma_{C^\ell}(f^*(E))\times M\to f^*(E)$.
We could only ask whether this map is $C^{\infty,k}$ if $f^*(E)$
was a smooth manifold.
\end{rem}
\begin{prop}\label{explaw-pb}
If $M$ is locally compact, then the following holds
for each $k\in \N_0\cup\{\infty\}$
and each $C^k$-manifold~$L$:
A map $g \colon L\to\Gamma_f$
is $C^k$ if and only if $g^\wedge\colon L\times M\to E$,
$(x,y)\mto g(x)(y)$
is $C^{k,\ell}$.
\end{prop}
\begin{prf}
By Proposition~\ref{eval-pb},
the evaluation map $\ve\colon \Gamma_f\times M\to E$ is $C^{\infty,\ell}$
and thus $C^{k,\ell}$.
If $g$ is $C^k$, then $g^\wedge=\ve\circ (g\times \id_M)$
is $C^{k,\ell}$, by Proposition~\ref{Cklchainmfd}(a).
Conversely, assume $g^\wedge$ is~$C^{k,\ell}$.
Let $\cA$ be the set of all charts $\phi\colon U_\phi\to V_\phi$
of~$M$ for which $U_\phi\sub f^{-1}(W(\phi))$
for some trivialization $\theta_\phi\colon E|_{W(\phi)}\to W(\phi)\times F$
of~$E$. Let $\theta_{\phi,2}\colon E|_{W(\phi)}\to F$
be the second component of $\theta_\phi$.
Then
\[
\Theta_\phi\colon \colon f^*(E)|_{U_\phi}\colon f^*(E)|_{U_\phi}\to U_\phi\times F,\;\;
v\mto(\pr_{f^*(E)},(\theta_{\phi,2}\circ \pr_2)(v))
\]
is a local trivialization of $f^*(E)$
and
$\cA$ an atlas for~$M$. Let $\Theta_{\phi,2}\colon \! f^*(E)|_{U_\phi}\!\! \to F$\linebreak
be the second
component of $\Theta_\phi$.
By Proposition~\ref{secsviafctns}, the linear map
\[
\Phi\colon \Gamma_{C^\ell}(f^*(E))\to\prod_{\phi\in \cA} C^\ell(U_\phi,F),\;\;
\sigma\mto (\Theta_{\phi,2}\circ\sigma|_{U_\phi})_{\phi\in\cA}
\]
is a homeomorphism onto its closed image.
As
\[ (\phi^{-1})^*\colon  C^\ell(U_\phi,F) \to C^\ell(V_\phi,F), \quad 
g\mto g\circ\phi^{-1} \] 
is an isomorphism of topological vector spaces, also
\[
\Big(\prod_{\phi\in\cA} (\phi^{-1})^*\Big)\circ \Phi\colon \Gamma_{C^\ell}(f^*(E))\to\prod_{\phi\in \cA}C^\ell(V_\phi,F)
\]
is a linear topological embedding with closed image.
Let $\Psi\colon \Gamma_f\to\Gamma_{C^\ell}(f^*(E))$
be the isomorphism of topological
vector spaces taking $\tau$ to $(\id_M,\tau)$
(cf.\ Lemma~\ref{sections-pb}). The map $g$ will be $C^k$ if
$\prod_{\phi\in \cA}(\phi^{-1})^*\circ \Phi\circ \Psi\circ g$
is $C^k$, which holds if we can show that $(\phi^{-1})^*\circ\, \Phi\circ\Psi\circ g$
is $C^k$. This will hold if $(\phi^{-1})^*\circ \, \Phi\circ\Psi\circ g\circ \psi^{-1}$
is $C^k$ for each chart $\psi\colon U_\psi\to V_\psi$ for~$L$.
Since $g(x)(y)\in E|_{W(\phi)}$ for all $x\in L$ and $y\in U_\phi$,
we can form the composition
\[
h:=\theta_2\circ g^\wedge\circ (\psi^{-1}\times\phi^{-1})\colon V_\psi\times V_\phi\to F,
\]
which is $C^{k,\ell}$.
Then $h^\wedge\colon V_\psi\to C^\ell(V_\phi, F)$
is $C^k$, by Theorem~\ref{explawCkell}.
It only remains to observe
that $(\phi^{-1})^*\circ \Phi\circ\Psi\circ g\circ \psi^{-1}=h^\wedge$.
\end{prf}
\subsection*{Spaces of compactly supported sections}
%
\begin{numba}
From Definition~\ref{defn-cp-supp-sec}
until Lemma~\ref{embed-sections-sum},
let $F$ be a locally convex space,
$k\in \N_0\cup\{\infty\}$ and $\pi\colon E\to M$ be a $C^k$-vector bundle
with typical fibre $F$ over a $C^k$-manifold~$M$.
\end{numba}
\begin{defn}\label{defn-cp-supp-sec}
The support $\Supp(\sigma)$ of a $C^k$-section $\sigma\colon M\to E$
of~$E$ is defined as the closure in~$M$
of $\{x\in M\colon \sigma(x)\not=0_x\}$.
\begin{description}[(D)]
\item[(a)]
For a closed subset $L\sub M$, we endow the closed vector subspace
\begin{eqnarray*}
\Gamma_{C^k_L}(E) &:=& \{\sigma\in \Gamma_{C^k}(E)\colon \Supp(\sigma)\sub L\}\\
&=& \{\sigma\in \Gamma_{C^k}(E)\colon (\forall x\in M\setminus L)\,\sigma(x)=0_x\}
\end{eqnarray*}
of $\Gamma_{C^k}(E)$ with the induced topology.
\item[(b)]
If $M$ is locally compact and $\sigma$-compact,
we let
\[
\Gamma_{C^k_c}(E):=\{\sigma\in \Gamma_{C^k}(E)\colon \mbox{$\Supp(\sigma)$ is compact}\}.
\]
The set $\cK(M)$ of compact subsets of~$M$ is directed under inclusion
and $\Gamma_{C^k_K}(E)\sub \Gamma_{C^k_L}(E)$ for $K,L\in\cK(M)$ with $K\sub L$.
We give
\[
\Gamma_{C^k_c}(E)=\bigcup_{L\in\cK(M)}\Gamma_{C^k_L}(E)=\dl \Gamma_{C^k_L}(E)\vspace{-.3mm}
\]
the locally convex direct limit topology.
\end{description}
\end{defn}
\begin{rem}\label{rem-supported-sections}
In the situation of Definition~\ref{defn-cp-supp-sec}(b), we have:\\
(a) If $K_1\sub K_2\sub\cdots$ is a compact exhaustion of $M$,
then $(K_n)_{n\in \N}$ is cofinal in $\cK(M)$,
whence
\[
\Gamma_{C^k_c}(E)=\dl \Gamma_{C^k_{K_n}}(E)\vspace{-.3mm}
\]
as a locally convex space.\\ 
(b) The inclusion map $\Gamma_{C^k_L}(E)\to \Gamma_{C^k}(E)$ is continuous
and linear for each $L\in \cK(M)$, entailing that the inclusion map
$\Gamma_{C^k_c}(E)\to\Gamma_{C^k}(E)$ is continuous (see Remark~\ref{firstremlcxDL}(d)).
As a consequence, $\Gamma_{C^k_c}(E)$ is Hausdorff and
induces the given topology on $\Gamma_{C^k_L}(E)$ for each $L\in \cK(M)$.\\ 
(c) Of course, $\Gamma_{C^k_c}(E)$ is only of interest if $\K=\R$.
\end{rem}
\begin{lem}\label{compactinopen-restr}
If $L\sub M$ is a closed subset and $U\sub M$
an open subset such that $L\sub U$
$($or a full submanifold such that $L\sub U^0)$,
then the
restriction map
$\rho\colon \Gamma_{C^k_L}(E)\to \Gamma_{C^k_L}(E|_U)$
is an isomorphism of topological vector spaces.
\end{lem}
\begin{prf}
The map $\rho$ is linear. We readily check that it is injective.
To see that $\rho$ is surjective, let $\sigma\in \Gamma_{C^k_L}(U)$.
Then $\tau(x):=\sigma(x)$ if $x\in U$,
$\tau(x):=0_x\in E_x$ for $x\in M\setminus L$
defines an element $\tau\in \Gamma_{C^k_L}(E)$ such that $\rho(\tau)=\sigma$.
By Lemmas~\ref{family-sections} and \ref{transinit},
the topology on $\Gamma_{C^k_L}(E)$
is initial with respect to $\rho$ and the restriction map
$r\colon \Gamma_{C^k_L}(E)\to\Gamma_{C^k}(E|_{M\setminus L})$.
Since $r$ is constant, we can omit it without changing the
initial topology. Thus, the topology on $\Gamma_{C^k_L}(E)$
is initial with respect to the bijection~$\rho$,
whence $\rho$ is a homeomorphism.
\end{prf}
\begin{lem}\label{vbdl-map-induce}
Let also $\pi_H\colon H\to M$ be a $C^k$-vector bundle
over $M$ with typical fibre a locally convex space~$Y$.
If $f\colon E\to H$ is a vector bundle map
of class $C^k$ over~$\id_M$, then the following holds:
\begin{description}[(D)]
\item[\rm(a)]
The mapping
\[
\Gamma_{C^k}(f)\colon \Gamma_{C^k}(E)\to\Gamma_{C^k}(H),\;\;
\sigma\mto f\circ \sigma
\]
is continuous linear, as well as its restriction
$\Gamma_{C^k_L}(f)\colon \Gamma_{C^k_L}(E)\to \Gamma_{C^k_L}(H)$
for each closed subset $L\sub M$.
\item[\rm(b)]
If $M$ is locally compact and $\sigma$-compact,
then $\Gamma_{C^k_c}(f)\colon \Gamma_{C^k_c}(E)\to\Gamma_{C^k_c}(H)$,
$\sigma\mto f\circ\sigma$ is a continuous linear map.
\end{description}
\end{lem}
\begin{prf}
(a) is immediate from Lemma~\ref{Cktoppu}(b).

(b) For $L\in \cK(M)$,
let $j_L\colon \Gamma_{C^k_L}(H)\to \Gamma_{C^k}(H)$
be the continuous linear inclusion map.
The map $\Gamma_{C^k_c}(f)$ is linear and its restriction to $\Gamma_{C^k_L}(E)$
is the continuous map $j_L\circ \Gamma_{C^k_L}(f)$, for each $L\in \cK(M)$.
Hence $\Gamma_{C^k_c}(f)$ is continuous, by Remark~\ref{firstremlcxDL}(d). 
\end{prf}
\begin{lem}\label{secs-top-module}
Pointwise multiplication turns
$\Gamma_{C^k}(E)$ into a topological
$C^k(M,\K)$-module.
In particular, for each $f\in C^k(M,\K)$ we get a continuous linear
multiplication operator
\[
m_f\colon \Gamma_{C^k}(E)\to\Gamma_{C^k}(E),\;\;
\sigma\mto f\sigma.
\]
\end{lem}
\begin{prf}
For $t\in \K$, $x\in M$ and $v\in E_x$, we can form the product $tv$ in the vector space
$E_x$. We obtain a multiplication map
$\mu\colon \K\times E\to E$ which is $C^k$ (as we readily check
using local trivializations).
Let 
\[ \Phi\colon C^k(M,\K\times E)\to C^k(M,\K)\times C^k(M,E) \] 
be the homeomorphism from Lemma~\ref{prodmapsp}.
By Lemma~\ref{Cktoppu}(b), the map $C^k(M,\mu)\colon C^k(M,\K\times E)\to C^k(M,E)$
is continuous.
Hence also the map $C^k(M,\mu)\circ \Phi^{-1}\colon C^k(M,\K)\times C^k(M,E)
\to C^k(M,E)$ is continuous, which restricts to the module
multiplication $m\colon C^k(M,\K)\times \Gamma_{C^k}(E)\to\Gamma_{C^k}(E)$.
Then also $m_f=m(f,\cdot)$ is continuous.
\end{prf}
In the following lemma, we assume that the ground field is $\K=\R$.
\begin{lem}\label{embed-sections-sum}
If $M$ is $\sigma$-compact and locally compact,
let $(V_i)_{i\in I}$ be a countable
locally finite cover of~$M$
by relatively compact, open subsets which cover~$M$
$($or full submanifolds
whose interiors $V_i^0$ cover~$M)$
Then
\[
\rho\colon \Gamma_{C^k_c}(E)\to\bigoplus_{i\in I}\Gamma_{C^k}(E|_{V_i}),\;\,
\sigma\mto (\sigma|_{V_i})_{i\in I}
\]
is a continuous linear map and has a continuous linear left inverse.
Notably, $\rho$ is
a linear topological embedding with closed image.
\end{lem}
\begin{prf}
Note first that $\rho$ is linear.
Being continuous on $\Gamma_{C^k_L}(E)$
for each $L\in \cK(M)$ by Remark~\ref{earlyremCk}(a), $\rho$ is continuous
(see Remark~\ref{firstremlcxDL}(d)).
By Proposition~\ref{findim-smoothly-para}, we find a smooth partition $(h_j)_{j\in J}$ of unity
subordinate to $(V_i)_{i\in I}$.
Thus $\Supp(h_j)\sub V_{i(j)}$ for some $i(j)\in I$.
Let $J_i:=\{j\in J\colon i(j)=i\}$ for $i\in I$
and define $g_i:=\sum_{j\in J_i}h_j$.
Then $(g_i)_{i\in I}$ is a smooth partition of unity
such that $L_i:=\Supp(g_i)\sub V_i$
(cf.\ Lemma~\ref{ops-locfin}(c)).
By Lemma~\ref{compactinopen-restr}, the restriction operator
\[
\rho_i\colon \Gamma_{C^k_{L_i}}(E)\to \Gamma_{C^k_{L_i}}(E|_{V_i})
\]
is an isomorphism of topological vector spaces.
Thus $\lambda_i:=\rho_i^{-1}\colon \Gamma_{C^k_{L_i}}(E|_{V_i})\to
\Gamma_{C^k_{L_i}}(E)\sub \Gamma_{C^k_c}(E)$ is a continuous
linear map. By Lemma~\ref{secs-top-module}, the multiplication operator
\[
m_i\colon \Gamma_{C^k}(E|_{V_i})\to\Gamma_{C^k_{L_i}}(E|_{V_i}),\;\;
\sigma\mto g_i|_{V_i}\sigma
\]
is continuous and linear. Thus $\lambda_i\circ m_i\colon \Gamma_{C^k}(E|_{V_i})\to\Gamma_{C^k_c}(E)$
is a continuous linear map for each $i\in I$. Then also the map
\[
\alpha\colon \bigoplus_{i\in I} \Gamma_{C^k}(E|_{V_i})\to\Gamma_{C^k_c}(E),\;\;
(\sigma_i)_{i\in I}\mto\sum_{i\in I}(\lambda_i\circ m_i)(\sigma_i)
\]
is continuous and linear, by Lemma~\ref{firstlasum}(a).
It remains to observe that $\alpha(\rho(\sigma))(x)=\sum_{i\in I}g_i(x)\sigma(x)=\sigma(x)$
for all $\sigma\in\Gamma_{C^k_c}(E)$ and $x\in M$,
whence $(\alpha\circ \rho)(\sigma)=\sigma$ and $\alpha\circ \rho$
is the identity map on $\Gamma_{C^k_c}(E)$.
\end{prf}
\begin{numba}\label{the-set-sum-sections}
Let
$k\in\N_0\cup\{\infty\}$.
Let $F_j$ be a locally convex space
and $\pi_j\colon E_j\to M$ be $C^k$-vector bundles with typical fibre~$F_j$
over a $C^k$-manifold~$M$,
for $j\in \{1,2\}$.
Using the local trivializations of the Whitney sum
described in \ref{sett-whit-sum}, we readily verify: The map
\[
p_j\colon E_1\oplus E_2\to E_j,\quad (v_1,v_2)\mto v_j
\]
is a vector bundle map over~$\id_M$ of class $C^k$,
for $j\in \{1,2\}$, and likewise
\[
\lambda_1\colon E_1\to E_1\oplus E_2,\;\; v\mto (v,0_x)
\]
with $x:=\pi_1(v)$ and $\lambda_2\colon E_2\to E_1\oplus E_2$,
$v\mto (0_x,v)$ with $x:=\pi_2(v)$.
\end{numba}
\begin{lem}\label{sections-whitney-sum}
In the situation of {\rm\ref{the-set-sum-sections}},
we have:
\begin{description}[(D)]
\item[\rm(a)]
The mapping
\[
\Phi:=(\Gamma_{C^k}(p_1),\Gamma_{C^k}(p_2))\colon
\Gamma_{C^k}(E_1\oplus E_2)\to
\Gamma_{C^k}(E_1)\times\Gamma_{C^k}(E_2)
\]
is an isomorphism of topological vector spaces,
and so is its restriction
$\Phi_L:=(\Gamma_{C^k_L}(p_1),\Gamma_{C^k_L}(p_2))\colon
\Gamma_{C^k_L}(E_1\oplus E_2)\to
\Gamma_{C^k_L}(E_1)\times\Gamma_{C^k_L}(E_2)$,
for each closed subset $L\sub M$.
\item[\rm(b)]
The map $\Phi_c:=(\Gamma_{C^k_c}(p_1),\Gamma_{C^k_c}(p_2))\colon
\Gamma_{C^k_c}(E_1\oplus E_2)\to
\Gamma_{C^k_c}(E_1)\times\Gamma_{C^k_c}(E_2)$
is an isomorphism of topological vector spaces
if $M$ is $\sigma$-compact and locally compact.
\end{description}
\end{lem}
\begin{prf}
The projection $\pr_j\colon \Gamma_{C^k}(E_1)\times \Gamma_{C^k}(E_2)\to
\Gamma_{C^k}(E_j)$ is continuous linear for $j\in \{1,2\}$.
To prove the first assertion, we consider the continuous linear map
$\Psi:=\Gamma_{C^k}(\lambda_1)\circ \pr_1+\Gamma_{C^k}(\lambda_2)\circ \pr_2$
and check that $\Phi\circ \Psi$ and $\Psi\circ \Phi$
are the identity maps on $\Gamma_{C^k}(E_1)\times \Gamma_{C^k}(E_2)$
and $\Gamma_{C^k}(E_1\oplus E_2)$, respectively.
The other asertions are proved along the same lines, replacing
the symbol $C^k$ with $C^k_L$ and $C^k_c$, respectively.
\end{prf}
\begin{ex}\label{sections-trivia}
If $F$ is a locally convex space, $k\in \N_0\cup\{\infty\}$
and $M$ be a $C^k$-manifold,
we can consider the trivial vector bundle $E:=M\times F$
with bundle projection $\pi:=\pr_1\colon M\times F\to M$.
Let $\pr_2\colon M\times F\to F$ be the projection onto the second factor.
Then $\id_E$ is a global trivialization and
\begin{equation}\label{sec-vs-fun}
\Gamma_{C^k}(E)\to C^k(M,F),\;\; \sigma\mto \pr_2\circ\, \sigma
\end{equation}
is an isomorphism of topological vector spaces, by Corollary~\ref{sectrivializable}.
If we endow $C^k_L(M,F):=\{f\in C^k(M,F)\colon f|_{M\setminus L}=0\}$
with the topology induced by $C^k(M,F)$
for a closed subset $L\sub M$, then the map in (\ref{sec-vs-fun})
restricts to an isomorphism
\[
\Gamma_{C^k_L}(E)\to C^k_L(M,F)
\]
of topological vector spaces. If $M$ is $\sigma$-compact and
locally compact, we can make $C^k_c(M,F)=\bigcup_{L\in \cK(M)} C^k_L(M,F)$
the locally convex direct limit $\dl C^k_L(M,F)$
and deduce that the map
\[
\Gamma_{C^k_c}(E)\to C^k_c(M,F),\;\; \sigma\mto \pr_2\circ \, \sigma
\]
is an isomorphism of topological vector spaces.
In each case, the inverse function takes $f$ to $(\id_M,f)$.
As a consequence of Lemma~\ref{embed-sections-sum}, in the case $\K=\R$
the map
\[
C^k_c(M,F)\to\bigoplus_{i\in I}C^k(M_i,F),\;\; f\mto (f|_{M_i})_{i\in I}
\]
is a linear topological embedding with closed image,
for each countable locally finite family $(M_i)_{i\in I}$
of relatively compact, open subsets of $M$ which cover~$M$
(or relatively compact, full submanifolds whose interiors $M_i^0$
cover~$M$).
\end{ex}
\begin{defn}\label{mixed-bundles}
Let $F$ is a locally convex space over the ground field~$\C$.
Let $M$ be a $C^k$-manifold
over the ground field~$\R$, with $k\in \N_0\cup\{\infty,\omega\}$.
Let $\pi\colon E\to M$ is a $C^k$-vector bundle with typical fibre~$F$,
endowed with a complex vector space structure on each fibre $E_x$.
We call $E$ a $C^k$-vector bundle with
typical fibre the complex locally convex space~$F$
if, for each $x_0\in M$, we find an open $x_0$-neighborhood $U\sub M$
and a local trivialization $\theta\colon E|_U\to U\times F$
whose second component is $\C$-linear on $E_x$ for each $x\in U$.
\end{defn}
\begin{rem}
The results and proofs of this section (and the
preceding one)
remain valid if $\K=\R$ and $F$ is a complex vector space,
if we read ``vector space'' as complex vector space
and ``linear map'' as ``complex linear map.''
Analogous to Lemma~\ref{fun-modules}(a), we can take $\K=\R$ as the ground field
but conclude that $C^k(M,\C)$ is a topological $\C$-algebra.
Analogous to Lemma~\ref{secs-top-module}, we get a topological $C^k(M,\C)$-module structure on $\Gamma_{C^k}(E)$
if $E$ is a $C^k$-vector bundle over the ground field~$\R$
which has as its typical fibre a complex locally convex space.
\end{rem}
\subsection*{{\boldmath$C^k$}-maps between spaces of compactly supported sections}
Let $k,\ell\in \N_0\cup\{\infty\}$,
$M$ be a $\sigma$-compact
finite-dimensional $C^\ell$-manifold
and $\pi_j\colon E_j\to M$ be $C^\ell$-vector bundles over~$M$.
Frequently, we have to deal with a non-linear map
$f\colon \Gamma_{C^\ell_c}(E_1)\to \Gamma_{C^\ell_c}(E_2)$
and know that the restriction to $\Gamma_{C^\ell_K}(E_1)$
is $C^\ell$ for each compact subset $K\sub M$.
We would like to deduce that $f$ is $C^k$,
but unfortunately this is not possible without additional hypotheses.
\begin{ex}
The bilinear map $C^\infty_c(\R,\R)\times C^\infty_c(\R,\R)\to C^\infty_c(\R\times\R,\R)$,
$(f,g)\mto f\otimes g$ with $(f\otimes g)(x,y):=f(x)g(y)$
is discontinuous on its domain $C^\infty_c(\R,\R\hspace*{-.3mm}\times\hspace*{-.3mm}\R)$
although its restriction to $C^\infty_K(\R,\R\hspace*{-.3mm}\times \hspace*{-.3mm}\R)\cong
C^\infty_K(\R,\R)\hspace*{-.3mm}\times\hspace{-.3mm}
C^\infty_K(\R,\R)$\linebreak
is continuous bilinear and thus $C^\infty$ for all compact sets $K\sub\R$,
see \cite{HST01}.
\end{ex}
\begin{ex}
The map $C^\infty_c(\R)\to C^\infty_c(\R)$, $f\mto f\circ f - f(0)$
is discontinuous, although its restriction to $C^\infty_K(\R)$
is smooth for each compact subset $K\sub\R$ (see \cite{Gl06b}).
\end{ex}
We now describe additional hypotheses which prevent
such pathologies.
Consider the following setting, using the ground field $\K=\R$.
\begin{numba}\label{sett-for-almloc}
Let $\ell\in \N_0\cup\{\infty\}$ and
$M_j$
be a locally compact $C^\ell$-manifold
for $j\in\{1,2\}$.
Let $F_j$ be a locally convex space over $\bL\in \{\R,\C\}$
and $\pi_j\colon E_j\to M_j$ be a $C^\ell$-vector bundle over~$M_j$
with typical fibre the locally convex topological $\bL$-vector space $F_j$, for $j\in \{1,2\}$.
Let $f\colon \Omega\to\Gamma_{C^\ell_c}(E_2)$ be a mapping
on an open subset $\Omega\sub \Gamma_{C^\ell_c}(E_1)$.
\end{numba}
\begin{defn}\label{def-almloc}
Consider a map $f$ as in~\ref{sett-for-almloc}.
\begin{description}[(D)]
\item[(a)]
If $M:=M_1=M_2$,
the map $f$ is called \emph{local}
if $f(\sigma)(x)=f(\tau)(x)$ for all $\sigma,\tau\in \Omega$
and all $x\in M$ such that $\sigma|_V=\tau|_V$ for some $x$-neighborhood
$V\sub M$.
\item[(b)]
$f$ is called \emph{almost local}
if there exist locally finite countable covers $(V_i)_{i\in I}$ of $M_1$
and $(W_i)_{i\in I}$ of~$M_2$ by relatively compact, open subsets
such that
\[
(\forall \sigma,\tau\in \Omega)(\forall i\in I)\;\;\sigma|_{V_i}=\tau|_{V_i}\;\Rightarrow\;
f(\sigma)|_{W_i}=f(\tau)|_{W_i}.
\]
\item[(c)]
If each $\sigma\in \Omega$ has an open neighborhood
$Q\sub\Omega$ such that $f|_Q$ is almost local, then
$f$ is called \emph{locally almost local}.
\end{description}
\end{defn}
\begin{rem}
(a) Every local map $f\colon \Omega\to \Gamma_{C^\ell_c}(E_2)$
is almost local. In fact, let $(V_i)_{i\in I}$
be any locally finite cover of~$M$ by relatively compact, open
subsets of~$M$ and set $W_i:=V_i$.\medskip

(b) Every almost local map is locally almost local.
\end{rem}
\begin{thm}\label{almloc-thm}
Let $k\in \N_0\cup\{\infty,\omega\}$
and $f\colon \Omega\to \Gamma_{C^\ell_c}(E_2)$
be a mapping as in {\rm\ref{sett-for-almloc}}.
If $f$ is locally almost local and its restriction to a mapping\linebreak
$\Omega\cap \Gamma_{C^\ell_K}(E_1)\to\Gamma_{C^\ell_c}(E_2)$
is $C^k_{\bL}$ for each compact set $K\sub M_1$,
then $f$ is~$C^k_{\bL}$.
\end{thm}
In the proof, we use that the relevant covers can be inflated.
\begin{lem}\label{inflate-cover}
Let $X$ be a $\sigma$-compact locally compact topological space and
$(V_i)_{i\in I}$ be a locally finite family of relatively compact,
open subsets $V_i\sub X$ such that $\bigcup_{i\in I}V_i=X$.
Then there exist relatively compact, open subsets $U_i\sub M$ containing the
closure $\wb{V}_i$
such that $(U_i)_{i\in I}$
is locally finite.
\end{lem}
\begin{prf}
Let $K_1\sub K_2\sub\cdots$ be a compact exhaustion of~$X$
(see Lemma~\ref{exhaustions-exist})
and $K_0:=\emptyset$.
For $i\in I$,
let
\[
n(i):=\max\{n\in \N_0\colon V_i\cap K_n=\emptyset\}.
\]
By Remark~\ref{trivi-lcp},
there exists a compact subset $L_i\sub X\setminus K_{n(i)}$
such that $\wb{V}_i\sub L_i^0=:U_i$.
For each $x\in X$, we have $x\in K_n^0$ for some $n\in \N$.
By Lemma~\ref{ops-locfin}(b), the set $I_0:=\{i\in \colon V_i\cap K_n\not=\emptyset\}$
is finite. For $i\in I\setminus I_0$, we have $n(i)\geq n$,
whence $U_i\sub X\setminus K_{n(i)}\sub X\setminus K_n$
and hence $U_i\cap K_n=\emptyset$.
\end{prf}
\emph{Proof of Theorem}~\ref{almloc-thm}.
We show that $f$ is $C^k_{\bL}$ on an open neighborhood
of a given section $\sigma_0\in \Omega$. After replacing $\Omega$
with a smaller open $\sigma_0$-neighborhood, we may assume that $f$
is almost local. After replacing $\Omega$ with $\Omega-\sigma_0$
and $f$ with $\sigma\mto f(\sigma+\sigma_0)$,
we may assume that $\sigma_0=0$.
We let $(V_i)_{i\in I}$ and $(W_i)_{i\in I}$
be as in Definition~\ref{def-almloc}(b).
By Lemma~\ref{inflate-cover},
there exists a locally finite cover $(U_i)_{i\in I}$
of $M_1$ by relatively compact, open subsets $U_i$
such that $\wb{V}_i\sub U_i$. By Exercise~\ref{exc-cutoff},
there exist $C^\ell$-functions $h_i\colon M\to [0,1]$
such that $L_i:=\Supp(h_i)\sub U_i$
and $h_i|_{\wb{V}_i}=1$.
By Lemma~\ref{embed-sections-sum}, the maps
\[
\rho_1\colon \Gamma_{C^\ell_c}(E_1)\to \bigoplus_{i\in I}\Gamma_{C^\ell}(E_1|_{U_i}),\;\;
\sigma\mto (\sigma|_{U_i})_{i\in I}
\]
and 
\[ \rho_2\colon \Gamma_{C^\ell_c}(E_2)
\to \bigoplus_{i\in I}\Gamma_{C^\ell}(E_2|_{W_i}),\quad 
\sigma\mto (\sigma|_{W_i})_{i\in I} \] 
are $\bL$-linear topological embeddings with closed image.
For each $i\in I$, the map
\[
r_i\colon \Gamma_{C^\ell_{L_i}}(E_1)\to \Gamma_{C^\ell_{L_i}}(E_1|_{U_i}),\;\;
\sigma\mto \sigma|_{U_i}
\]
is an isomorphism of topological $\bL$-vector spaces, by Lemma~\ref{compactinopen-restr}.
Hence $\lambda:=\oplus_{i\in I}r_i^{-1}\colon \bigoplus_{i\in I}\Gamma_{C^\ell_{L_i}}(E_1|_{U_i})
\to\bigoplus_{i\in I}\Gamma_{C^\ell_{L_i}}(E_1)$
is an isomorphism of topological $\bL$-vector spaces.
We let
\[
m_i\colon \Gamma_{C^\ell}(E_1|_{U_i})\to\Gamma_{C^\ell_{L_i}}(E_1|_{U_i})
\]
be the multiplication operator $\sigma\mto h_i|_{U_i}\sigma$,
which is continuous and \break {$\bL$-linear} by Lemma~\ref{secs-top-module}.
Then also the linear mapping 
\[ m:=\oplus_{i\in I} m_i\colon \bigoplus_{i\in I}\Gamma_{C^\ell}(E_1|_{U_i})
\to\bigoplus_{i\in I}\Gamma_{C^\ell_{L_i}}(E_1|_{U_i}) \] 
is continuous.
For each $i\in I$, the set
\[
\Omega_i:=\Omega\cap \Gamma_{C^\ell_{L_i}}(E_1)
\]
is an open $0$-neighborhood in $\Gamma_{C^\ell_{L_i}}(E_1)$.
Hence $Q:=(\lambda\circ m\circ \rho_1)^{-1}(\bigoplus_{i\in I}\Omega_i)$
is an open $0$-neighborhood in
$\Gamma_{C^\ell_c}(E_1)$. For each $i\in I$, the map
\[
f_i\colon \Omega_i\to \Gamma_{C^\ell}(E_2|_{W_i}),\;\;
\sigma\mto f(\sigma)|_{W_i}
\]
is $C^\ell_{\bL}$, whence
$\oplus_{i\in I}f_i \colon \bigoplus_{i\in I}\Omega_i\to \bigoplus_{i\in I}
\Gamma_{C^\ell}(E_2|_{W_i})$ is $C^k_{\bL}$ by Proposition~\ref{diff-sum}.
For each $\sigma\in \Omega\cap Q$ and $i\in I$, we have
\[
f(\lambda_i(h_i|_{U_i}\sigma|_{U_i}))|_{W_i}
=f(h_i\sigma)|_{W_i}=f(\sigma)|_{W_i}
\]
since $(h_i\sigma)|_{U_i}=\sigma|_{U_i}$. Hence
\[
\rho_2\circ f|_{\Omega\cap Q}=\oplus_{i\in I}f_i\circ \lambda\circ m\circ \rho_1|_{\Omega\cap Q},
\]
showing that $\rho_2 \circ f|_{\Omega\cap Q}$ is $C^k_{\bL}$.
Then also $f|_{\Omega\cap Q}$ is $C^k_{\bL}$, by Lemma~\ref{corestr}
and Exercise~\ref{excreatosub}, respectively.\vspace{.7mm}
\qed

Results like the following one
are frequently referred to as an ``$\omega$-Lemma.''
\begin{prop}\label{omega-lemma}
Let $k,\ell\in \N_0\cup\{\infty\}$,
$M$ be a $C^{k+\ell}$-manifold over the ground field~$\R$
and $\pi_j\colon E_j\to M$ be a $C^{k+\ell}$-vector bundle
with typical fibre a locally convex real vector space~$F_j$,
for $j\in\{1,2\}$. Let $\Omega\sub E_1$ be
an open subset and $f\colon \Omega\to E_2$
be a $C^{k+\ell}$-map such that
$f(\Omega\cap (E_1)_x)\sub (E_2)_x$ for all $x\in M$.
Assume that there exists a compact subset $B\sub M$ such that
$f(0_x)=0_x$ for all $x\in M\setminus B$.
Then $\Gamma_{C^\ell_c}(E_1)\cap C^\ell(M,\Omega)$
is an open subset of $\Gamma_{C^\ell_c}(E_1)$ and the following
map is $C^k$:
\[
f_*\colon \Gamma_{C^\ell_c}(E_1)\cap C^\ell(M,\Omega)\to \Gamma_{C^\ell_c}(E_2),\;\;
\sigma\mto f\circ \sigma.
\]
\end{prop}
Let us single out an argument of the proof.
\begin{lem}\label{domain-open}
Let $k\in \N_0\cup\{\infty\}$
and $\pi\colon E\to M$ be a $C^k$-vector bundle
over a $\sigma$-compact, locally compact $C^k$-manifold~$M$.
If $\Omega\sub E$ is open,
then\linebreak
$\{\sigma\in \Gamma_{C^k_c}(E)\colon \sigma(M)\sub \Omega\}=$
$\Gamma_{C^k_c}(E)\cap C^k(M,\Omega)$
is an open~subset~of~$\Gamma_{C^k_c}(E)$.
\end{lem}
\begin{prf}
Let $(M_i)_{i\in I}$ be a countable
locally finite family of compact full submanifolds $M_i\sub M$
whose interiors cover~$M$.
By Lemma~\ref{embed-sections-sum}, the map
\[
\rho\colon \Gamma_{C^k_c}(E)\to\bigoplus_{i\in I}\Gamma_{C^k}(E|_{M_i}),\;\;
\sigma\mto (\sigma|_{M_i})_{i\in I}
\]
is continuous.
Now $C^k(M_i,\Omega)=C^k(M_i,E)\cap \lfloor M_i,\Omega\rfloor$
is open in $C^k(M_i,E)$, whence $\Gamma_{C^k}(E|_{M_i})\cap C^k(M_i,\Omega)$
is open in $\Gamma_{C^k}(E|_{M_i})$.
So $\Gamma_{C^k_c}(E)\cap C^k(M,\Omega)$\linebreak
$=\rho^{-1}\big(\bigoplus_{i\in I}(\Gamma_{C^k}(E|_{M_i})\cap C^k(M_i,\Omega))\big)$ is open in $\Gamma_{C^k_c}(E)$.
\end{prf}
\noindent
\emph{Proof of Proposition}~\ref{omega-lemma}.
The domain is open by Lemma~\ref{domain-open}.
The map $f$ is local,
whence we only need to show
that its restriction to a map $g\colon \Gamma_{C^\ell_K}(E_1)\cap C^\ell(M,\Omega)\to
\Gamma_{C^\ell_c}(E_2)$
is~$C^k$. As the image of $g$ is contained in the closed
vector subspace $\Gamma_{C^\ell_{K\cup B}}(E_2)$
which is also a closed vector subspace of $\Gamma_{C^\ell}(E_2)$,
it suffices to show that $g$ is $C^k$ as a map to
$\Gamma_{C^\ell}(E_2)$ (see Lemma~\ref{corestr}).
We find a countable locally finite family $(M_i)_{i\in I}$
of compact full submanifolds of~$M$ whose interiors cover~$M$,
such that local trivializations
$\theta_i\colon E_1|_{M_i}\to M_i\times F_1$
and $\Theta_i\colon E_2|_{M_i}\to M_i\times F_2$
exist for each $i\in I$. Let $\theta_{i,2}\colon E_1|_{M_i}\to F_1$
and $\Theta_{i,2}\colon E_2|_{M_i}\to F_2$ be the second component of $\theta_i$
and $\Theta_i$, respectively.
The map
\[
\Gamma_{C^\ell}(E_2)\to\prod_{i\in I}C^\ell(M_i,F_2),\;\;
\sigma\mto (\Theta_{i,2}\circ\sigma|_{M_i})_{i\in I}
\]
is a linear topological embedding with closed
image by Lemma~\ref{family-sections}.
Hence, by Lemma~\ref{corestr},
it suffices to show that
\[
g_i\colon \Gamma_{C^\ell_L}(E_1)\cap C^\ell(M,\Omega)
\to C^\ell(M_i,F_2),\;\; \sigma\mto \Theta_{i,2}\circ f\circ \sigma|_{M_i}
\]
is $C^k$ for each $i\in I$.
Note that $\Omega_i:=\{(x,y)\in M_i \times F_1)\colon \theta_i^{-1}(x,y)\in \Omega\}$
is an open subset of $M_i\times F_1$. We have
\[
\Theta_i(f\circ \theta_i^{-1}(x,y))=(x,h_i(x,y))\;\;\mbox{for $(x,y)\in\Omega_i$}
\]
for a $C^{k+\ell}$-map $h_i\colon \!\Omega_i\hspace*{-.2mm}\to\hspace*{-.2mm} F_2$.
Now $Q_i\hspace*{-.3mm}:=\hspace*{-.3mm}\{\phi\hspace*{-.2mm}\in\hspace*{-.2mm} C^\ell(M_i,F_1)\colon \! (\id_{M_i},\phi)(M_i)
\!\sub \! \Omega_i\}$\linebreak
is an open subset of $C^\ell(M_i,F_1)$ and we show that the map
\[
(h_i)_*\colon Q_i\to C^\ell(M_i, F_2);\;\; \phi\mto h_i\circ (\id_{M_i}, \phi)
\]
is $C^k$. In fact, the evaluation map $\ve_i\colon C^\ell(M_i,F_1)\times M_i\to F_1$
is $C^{\infty, \ell}$ by Proposition~\ref{lem:smooth-eval-vec} and thus $C^{k,\ell}$, whence
\[
((h_i)_*)^\wedge(\phi,x)=h_i(x,\phi(x))=h_i(x,\ve_i(\phi,x))
\]
is $C^{k,\ell}$ in $(\phi,x)\in Q_i\times M_i$,
by Remark~\ref{rem-simpl-Ckl}(e) and Proposition~\ref{Cklchainmfd}(b).
Hence $h_i$ is $C^k$, by Proposition~\ref{explawfuncmfd}(c).\vspace{1mm}\qed

For trivial vector bundles, we deduce:
\begin{cor}\label{open-in-tefu}
Let $\ell\in \N_0\cup\{\infty\}$,
$M$ be a $\sigma$-compact, locally compact $C^\ell$-manifold
over~$\R$
and $E$ be a locally convex
space. If $U\sub E$
is open then $C^\ell_c(M,U):=\{f\in C^\ell_c(M,E)\colon f(M)\sub U\}$
is an open subset of $C^\ell_c(M,E)$.
\end{cor}
\begin{prf}
The subset $M\times U$ of the trivial vector bundle $M\times E$
is open. We identify $C^\ell_c(M,E)$ with $\Gamma_{C^\ell_c}(M\times E)$
and $C^\ell_c(M,F)$ with $\Gamma_{C^\ell_c}(M\times F)$,
as in Example~\ref{sections-trivia}.
Then $C^\ell_c(M,U)$ corresponds to the set
\[ \Gamma_{C^\ell_c}(M\times E)\cap C^\ell(M,M\times U),\]
which is open in $\Gamma_{C^\ell_c}(M\times E)$ by
Lemma~\ref{domain-open}.
\end{prf}
\begin{cor}\label{pushforwards-cp-supp}
Let $k,\ell\in \N_0\cup\{\infty\}$,
$M$ be a $\sigma$-compact, locally compact $C^{k+\ell}$-manifold
over the ground field~$\R$
and $E$ as well as $F$ be locally convex
spaces. Let $U\sub E$
be an open subset
and $g\colon M\times U\to F$ be a $C^{k+\ell}$-map.
If $M$ is not compact, we assume that $0\in U$
and that there exists a compact subset $B\sub M$
such that $f(x,0)=0$ for all $x\in M\setminus B$.
Then the following map is $C^k$:
\[
g_*\colon C^\ell_c(M,U)\to C^\ell_c(M,F),\quad f\mto g\circ (\id_M,f).
\]
\end{cor}
\begin{prf}
Identify $C^\ell_c(M,E)$ with $\Gamma_{C^\ell_c}(M\times E)$
and $C^\ell_c(M,F)$ with the space $\Gamma_{C^\ell_c}(M\times F)$,
as in Example~\ref{sections-trivia}.
Using the $C^{k+\ell}$-map
$h\colon M\times U\to M\times F$, \break {$(x,y)\mto (x,g(x,y))$}, the mapping
$g_*$ corresponds to the map
\[
h_*\colon \Gamma_{C^\ell_c}(M\times E)\cap C^\ell(M,M\times U)\to \Gamma_{C^\ell_c}(M\times F),
\;\; \sigma\mto h\circ \sigma,
\]
which is $C^k$ by Proposition~\ref{omega-lemma}.
\end{prf}
Also the following variants are useful.
\begin{prop}\label{variant-pushforwards}
Let $k\in \N_0\cup\{\infty,\omega\}$,
$\ell\in \N_0\cup\{\infty\}$,
$M$ be a $\sigma$-compact, locally compact $C^\ell$-manifold over~$\R$.
For
$\bL\in\{\R,\C\}$, let $E$ and $F$ be locally convex
topological $\bL$-vector spaces, and $U\sub E$
be open.
Then we have:
\begin{description}[(D)]
\item[\rm(a)]
If $g\colon U\to F$ is a $C^{k+\ell}_{\bL}$-map, then
$C^\ell_c(M,g)\colon C^\ell_c(M,U)\to C^\ell_c(M,F)$, $f\mto g\circ f$
is a $C^k_{\bL}$-map;
if $M$ is not compact, we assume that $0\in U$
and $g(0)=0$ here.
\item[\rm(b)]
Consider a mapping $g\colon M\times U\to F$
of the form $g=h\circ (q\times \id_U)$
for some $C^\ell_\K$-map $q\colon M\to N$ to a $C^{k+\ell}_{\bL}$-manifold~$N$
without boundary and a
$C^{k+\ell}_{\bL}$-map $h\colon N\times U\to F$.
If $M$ is not compact, assume that $0\in U$
and assume that there exists a compact subset $B\sub M$ with $h(q(x),0)=0$
for all $x\in M\setminus B$.
Then $g_*\colon C^\ell_c(M,U)\to C^\ell_c(M,F)$, $f\mto g\circ (\id_M,f)$
is~$C^k_{\bL}$.
\end{description}
\end{prop}
\begin{prf}
(a) We first assume that $k\not=\omega$.
Let $(M_i)_{i\in I}$ be a locally finite, countable family of compact
full submanifolds of~$M_i$.
For each $i\in I$,
the map
$C^\ell(M_i,g)\colon C^\ell(M_i,U)\to C^\ell(M_i,F)$
is $C^k_\bL$, by Proposition~\ref{superpo-cp}(a).
Since $\rho\colon C^\ell_c(M,F)\to\bigoplus_{i\in I}C^\ell(M_i,F)$,
$f\mto (f|_{M_i})_{i\in I}$
is a linear topological embedding with closed image,
it only remains to show that $\rho\circ C^\ell(M,g)$ is~$C^k_\bL$.
But $r\colon C^\ell_c(M,U)\to\bigoplus_{i\in I} C^\ell(M_i,U)$,
$f\mto (f|_{M_i})_{i\in I}$ is $C^\infty_{\bL}$ and 
\[
\rho\circ C^\ell_c(M,g)=\big(\oplus_{i\in I} C^\ell(M_i,g)\big)\circ r
\]
is $C^k_{\bL}$ as the first map on the right-hand side is $C^k_\bL$ by
Proposition~\ref{diff-sum}.

The case $k=\omega$, $\bL=\C$ coincides with the case $k=\infty$, $\bL=\C$.

If $k=\omega$ and $\bL=\R$,
there is  a komplex analytic extension $h\colon W\to F_\C$
of~$g$ to an open set $W\sub E_\C$ with $U\sub W$.
Then $C^\ell_c(M,W)$ is open in
$C^\ell_c(M,E_\C)=C^\ell_c(M,E)_\C$
and $C^\ell_c(M,h)\colon C^\ell_c(M,W)\to C^\ell_c(M,F_\C)=C^\ell_c(M,F)_\C$
is a complex analytic extension for $C^\ell_c(M,g)$,
whence $C^\ell_c(M,g)$ is real analytic.

(b) Let $(M_i)_{i\in I}$, $r$ and $\rho$ be as in the proof of~(a).
For $i\in I$, let $g_i$ be the restriction of $g$ to a map
$M_i\times U\to F$. Then $(g_i)_*\colon C^\ell(M_i,U)\to C^\ell(M_i,F)$ is
$C^k_\bL$ for each $i \in I$, by Proposition~\ref{superpo-cp}(b).
Since $\rho\circ g_*=\big(\oplus_{i\in I}(g_i)_*\big)\circ r$,
we deduce that $g_*$ is~$C^k_\bL$.
\end{prf}
\begin{cor}\label{bilin-pushf}
Let $M$ be a $\sigma$-compact, locally compact $C^\ell$-manifold over~$\R$
with $\ell\in \N_0\cup\{\infty\}$.
Let $E_1$, $E_2$, and $F$ be locally convex spaces and $\beta\colon E_1\times E_2\to F$
be a continuous bilinear map, Then also the bilinear map
$C^\ell_c(M,E_1)\times C^\ell_c(M,E_2)\to C^\ell_c(M,F)$,
$(f,g)\mto \beta\circ (f,g)$ is continuous.
\end{cor}
\begin{prf}
Identifying the domain with $C^\ell_c(M,E_1\times E_2)$,
the mapping in contention corresponds to $C^\ell_c(M,\beta)$,
which is continuous as $\beta$ is smooth and hence $C^\ell$
(see Proposition~\ref{variant-pushforwards}(a)).
\end{prf}
\begin{ex}\label{tefu-algs}
Let $\ell\in \N_0\cup\{\infty\}$
and $M$ be a $\sigma$-compact, locally compact $C^\ell$-manifold over~$\R$.
Let $\bL\in\{\R,\C\}$.
Corollary~\ref{bilin-pushf} implies
the following, for the usual pointwise operations:
\begin{description}[(D)]
\item[(a)]
$C^\ell_c(M,\bL)$ is an
associative topological $\bL$-algebra
(without unit element, if $M$ is not compact).
\item[(b)]
If $\cA$ is an associative, locally convex topological $\bL$-algebra, then
$C^\ell_c(M,\cA)$ is an associative topological
$\bL$-algebra (and topological $C^\ell_c(M,\bL)$-algebra);
\item[(c)]
If $\g$ is a locally convex topological Lie algebra over~$\bL$,
then also $C^\ell_c(M,\g)$ is topological Lie algebra over~$\bL$.
\end{description}
\end{ex}
We also have the following.
\begin{prop}\label{bundles-cp-supp-modules}
Let $\ell\in \N_0\cup\{\infty\}$
and $M$ be a $\sigma$-compact, locally compact $C^\ell$-manifold over~$\R$.
Let $\bL\in\{\R,\C\}$ and $\pi\colon E\to M$ be a $C^\ell$-vector bundle over~$M$,
with typical fibre a locally convex topological $\bL$-vector space~$F$.
Then pointwise multplication turns the topological $\bL$-vector space
$\Gamma_{C^\ell_c}(E)$ into a topological $C^\ell_c(M,\bL)$-module.
\end{prop}
\begin{prf}
For $f\in C^\ell(M,\bL)$ and $\sigma\in \Gamma_{C^\ell_c}(E)$,
setting $\mu(f,\sigma)(x):=(f\sigma)(x):=f(x)\sigma(x)$,
we get a $C^\ell$-section $f\sigma$ of~$E$ (cf.\ Lemma~\ref{secs-top-module}),
and we readily check it has compact support.
Using Example~\ref{sections-trivia} and Lemma~\ref{sections-whitney-sum}(b),
we identify $C^\ell_c(M,\bL)\times \Gamma_{C^\ell_c}(E)$
with $\Gamma_{C^\ell_c}((M\times\bL)\oplus E)$.
Then $\mu$ is a local map and we only need to check
its restriction to $C^\ell_K(M,\bL)\times \Gamma_{C^\ell_K}(E)$
is continuous for each compact subset $K\sub M$.
This follows from the fact that the module multiplication
$C^\ell(M,\bL)\times \Gamma_{C^\ell}(E)\to\Gamma_{C^\ell}(E)$
is continuous, as a consequence of
Lemma~\ref{secs-top-module}.
\end{prf}
\begin{small}
\subsection*{Exercises for Section~\ref{secspacemfd}}

\begin{exer} \mlabel{exer:4.1.5} Show that, for $k \in \N_0$, the map 
$$ C^k(M,N) \into C(T^k(M), T^k(N))_{c.o.}, \quad 
f \mapsto T^k(f), $$ 
is a topological embedding if $C^k(M,N)$ is endowed with the compact open \break 
$C^k$-topology. 
\end{exer}

\begin{exer} \mlabel{exer:4.1.6} Let $M_1, M_2, N_1$ and $N_2$ be $C^k$-manifolds
with $k\in \N_0\cup\{\infty\}$
and $\phi \: M_1 \to M_2$ and $\psi \:  N_2 \to N_1$ be $C^k$-diffeomorphisms. 
Show that the map 
\[  C^k(M_2,N_2) \to C^k(M_1,N_1), \quad f \mapsto \psi \circ f \circ \phi \] 
is a homeomorphism with respect to the compact-open $C^k$-topology.
\end{exer}

\begin{exer}
Along the lines of Lemma~\ref{prodmapsp}, show that $C^k(M,\prod_{i\in I}N_i)$ is homeomorphic to
$\prod_{i\in I}C^k(M,N_i)$
for each $k\in\N_0\cup\{\infty\}$, each $C^k$-manifold~$M$ and each 
finite family $(N_i)_{i\in I}$
of $C^k$-manifolds.
\end{exer}
\begin{exer}
Given $k\in\N_0\cup\{\infty\}$, let $M$, $N$, $L$, and $P$ be $C^k$-manifolds
and $g\colon M\times N\times P\to L$ be a $C^k$-map.
\begin{description}[(D)]
\item[(a)]
Show that the map
$\Phi\colon C^k(M,N)\times P\to C^k(M,L)$, $(f,p)\mto g^p\circ (\id_M,f)$
is continuous, where $g^p:=g(\cdot,p)$ [Note that $\Phi(f)=C^k(M,g)(\id_M,g,c_p)$ with notation as in Lemma~\ref{toconstcts}).
\item[(b)]
Show that $\Phi$ is continuous also if
$g\colon (M\times N)\times P\to N$
is $C^{k,0}$, by induction on $k\in\N_0$ [use Lemma~\ref{pushpar} and the identity
\[ T(g^p\circ (\id_M,f))=T(g^p)\circ (\id_{TM},Tf). \] 
\end{description}
\end{exer}

\begin{exer}
Let $M$ be a finite-dimensional $C^k$-manifold (possibly with rough boundary),
$E$ and $F$ be locally convex spaces, $U\sub E$ be an open subset,
$K\sub M$ be a compact subset with non-empty interior,
$k,\ell\in\N_0\cup\{\infty\}$ and $g\colon K^0\times U\to F$
be a $C^{k+\ell}$-function. Then $C^k(M,E)\cap\lfloor K,U\rfloor=
\{f\in C^k(M,E)\colon f(K)\sub U\}$ is an open subset of $C^k(M,E)$.
Using the Exponential Law \ref{explawCkell}, show that
\[
\{f\in C^k(M,E)\colon f(K)\sub U\}
\to C^k(K^0,F),\;\, f\mto g\circ (\id_{K^0},f|_{K^0})
\]
is a $C^\ell$-map.
\end{exer}

\end{small}
\section{Manifolds of mappings on compact manifolds}\label{sec-mfdmps}
%
%
%
We now construct a smooth manifold structure
on $C^\ell(M,N)$ if $N$ a smooth manifold admitting a local addition,
$\ell\in \N_0\cup\{\infty\}$,
and $M$ is a compact $C^\ell$-manifold,
which may have a rough boundary.
The manifold of mappings is not a pure manifold,
but modeled on a set $\cE$ of locally convex spaces (as in Definition~\ref{non-pure}),
which we introduce first.
\begin{numba}\label{special-gamma}
For $f\in C^\ell(M,N)$,
we define $\Gamma_f:=\{\tau\in C^\ell(M,TN)\colon \pi_{TN}\circ \tau=f\}$
as in \ref{def-gamma-f}, with $E:=TN$.
Thus $\Gamma_f\cong \Gamma_{C^\ell}(f^*(TN))$
(see Lemma~\ref{sections-pb}). We let $\cE:=\{\Gamma_f\colon f\in C^\ell(M,N)\}$.
\end{numba}
\begin{thm}\label{thmmfdmps}
Let $N$ be a smooth manifold admitting a local addition,
$\ell\in\N_0\cup\{\infty\}$,
and $M$ be a compact $C^\ell$-manifold
which may have a rough boundary.
Then the following holds:
\begin{description}[(D)]
\item[\rm(a)]
There is a unique smooth manifold structure on the set $C^\ell(M,N)$
which is modeled on $\cE$ and for which the following Exponential Law holds for all $k\in\N_0\cup\{\infty\}$:
For each $C^k$-manifold~$L$, possibly with rough boundary,
a map $g\colon L\to C^\ell(M,N)$ is~$C^k$ if and only if
\[ g^\wedge\colon L\times M\to N,\quad  g^\wedge(x,y):=f(x)(y) \] 
 is a $C^{k,\ell}$-map.
\item[\rm(b)]
The topology underlying the smooth manifold $C^\ell(M,N)$ from \emph{(a)}
is the compact-open $C^\ell$-topology.
\item[\rm(c)]
The evaluation map $\ve\colon C^\ell(M,N)\times M\to N$, $(f,x)\mto f(x)$ is $C^{\infty,\ell}$
and thus~$C^\ell$. Notably, the point evaluation $\ve_x\colon C^\ell(M,N)\to N$,
$f\mto f(x)$ is smooth for each $x\in M$.
\item[\rm(d)]
For each $f\in C^\ell(M,N)$
and $v\in T_f (C^\ell(M,N))$, we have $(T\ve_x(v))_{x\in M}\in \Gamma_f$
and the map
$\Theta\colon T_f(C^\ell(M,N))\to \Gamma_f$, $v\mto (T\ve_x(v))_{x\in M}$
is an isomorphism of topological vector spaces.
\end{description}
\end{thm}
\begin{prf}
Consider a local addition $\Sigma\colon U\to N$ for~$N$,
the open subset $U':=(\pi_{TN},\Sigma)(U)\sub N\times N$
and the associated $C^\infty$-diffeomorphism $\theta\colon U\to U'$,\linebreak
$v\mto (\pi_{TN}(v),\Sigma(v))$. 
For $f\in C^\ell(M,N)$,
\[
O_f := \Gamma_f \cap C^\ell(M,U)
\]
is an open subset of~$\Gamma_f$,
\[
O_f':=\{g\in C^\ell(M,N)\colon (f,g)(M)\sub U'\}
\]
is an open subset of $C^\ell(M,N)$ and the map
\[
\psi_f \colon O_f \to O'_f,\quad\tau\mto \Sigma\circ \tau
\]
is a homeomorphism with inverse $\phi_f\colon O_f'\to O_f$,
$h \mto \theta^{-1}\circ (f,h)$.
If also $g\in C^\ell(M,N)$, then $\psi_g^{-1}\circ \psi_f$
has an open domain, $\psi_f^{-1}(O_f'\cap O_g')$.
It is smooth there by Proposition~\ref{explaw-pb}, as the map
\begin{equation}\label{eq1}
(\tau,x)  \mto  (\psi_g^{-1}\circ\psi_f)(\tau)(x)=\theta^{-1}(g(x),\Sigma(\tau(x)))
\end{equation}
is $C^{\infty,\ell}$. Hence $C^\ell(M,N)$ has a smooth manifold structure
such that each of the maps $\phi_f$ is a local chart.

Note that the smooth manifold structure on $C^\ell(M,N)$ just constructed satifies~(b),
by construction. It will satisfy~(c) if we can show that
$\ve(\psi_f(\tau),x)$ is $C^{\infty,\ell}$
in $(\tau,x)\in O_f\times M$ for all $f\in C^\ell(M,N)$.
But
\[
\ve(\psi_f(\tau),x)=\Sigma(\tau(x))=\Sigma(\ev(\tau,x)),
\]
where $\ev \colon \Gamma_f \times M\to TN$, $(\tau,x)\mto\tau(x)$
is $C^{\infty,\ell}$ by Proposition~\ref{eval-pb}.
To see that the smooth manifold structure we constructed satisfies the exponential law
described in~(a),
let $k\in\N_0\cup\{\infty\}$ and $g\colon L\to C^\ell(M,N)$ be a map, where
$L$ is a $C^k$-manifold which may have
a rough boundary.
If~$g$ is~$C^k$, then $g^\wedge=\ve\circ (g\times\id_M)$
is~$C^{k,\ell}$.
If, conversely, $g^\wedge$ is $C^{k,\ell}$,
let $x_0\in L$ and $f:=g(x_0)$.
Then $g$ is continuous as a map to $C(M,N)$,
by Proposition~\ref{ctsexp},
whence $L\to C(M,N\times N)$, $x\mto (f,g(x))$
is continuous.
We therefore find an open $x_0$-neighborhood
$W\sub L$ such that $(f,g(x)) \in C(M,U')$
for all $x\in W$, whence $g(W)\sub O_f'$.
Abbreviate $h:=\phi_f\circ g|_W\colon W\to O
_f\sub\Gamma_f$.
Then
\[
h^\wedge(x,y)=\phi_f(g(x))(y)=(\theta^{-1}\circ (f,g(x)))(y)=\theta^{-1}(f(y),g^\wedge(x,y))
\]
is $C^{k,\ell}$, whence $h$ is~$C^k$ by Proposition~\ref{explaw-pb}.
Thus $g|_W=\psi_f\circ h$ is $C^k$ and hence so is~$g$.

Uniqueness in~(a): We write $C^\ell(M,N)'$
for $C^\ell(M,N)$, endowed with another smooth manifold structure
for which the exponential law from~(a) holds.
We prove uniqueness of the smooth manifold structures in the sense that
$g:=\id\colon C^\ell(M,N)'\to C^\ell(M,N)$ is a $C^\infty$-diffeomorphism.
Hence, if we assume that also $C^\ell(M,N)'$
is modeled on $\cE$, then the smooth manifold structures will coincide.
Since $\id\colon C^\ell(M,N)'\to C^\ell(M,N)'$ is a smooth map,
the exponential law entails that the evaluation map
$\ve=\id^\wedge\colon C^\ell(M,N)'\times M\to N$ is~$C^{\infty,\ell}$.
As~$g$ satisfies
$g^\wedge=\ve$, the exponential law entails that~$g$ is smooth.
Replacing the roles of $C^\ell(M,N)$ and $C^\ell(M,N)'$,
we see that also $g^{-1}$ is smooth.
Thus $g$ is a $C^\infty$-diffeomorphism.

(d) By the preceding, the smooth manifold structure on $C^\ell(M,N)$
is independent of the choice of local addition.
We now assume that the local addition~$\Sigma$
is normalized. Fix $f\in C^\ell(M,N)$. It suffices to show that
\begin{equation}\label{chncha}
(\Theta\circ T\psi_f)(0,\tau)
=\tau \;\,\mbox{for each $\tau\in \Gamma_f$,}
\end{equation}
i.e., $(T\ve_x\circ T\psi_f)(0,\tau)=\tau(x)$
for all $\tau$ as before and all $x\in M$.
We here use the natural identification of $T\Gamma_f$ with $\Gamma_f\times \Gamma_f$
which identifies $T_0\Gamma_f$ with $\{0\}\times \Gamma_f$.
To verify (\ref{chncha}), we find it convenient to work
with geometric tangent vectors
(as in Definition~\ref{def-geom-tang} and Exercise~\ref{exer:3.2.1}).
Thus $(0,\tau)$ corresponds to the geometric
tangent vector $[t\mto t\tau]$.
We calculate
\begin{eqnarray*}
T\ve_x T\psi_f([t\mto t\tau])
&=&
T\ve_x([t\mto \psi_f(t\tau)])=T\ve_x([t\mto\Sigma\circ t\tau])\\
&=& [t\mto \Sigma(t\tau(x))]
=
[t\mto \Sigma|_{U\cap T_{f(x)}N}(t\tau(x))]\\
&=&
T_0\Sigma|_{U \cap T_{f(x)}N}([t\mto t \tau(x)])=\tau(x),
\end{eqnarray*}
using (\ref{bettersigma}) for the final equality.
\end{prf}
\begin{defn}
For $M$, $N$, and $\ell$ as in Theorem~\ref{thmmfdmps},
the unique smooth manifold structure provided by Theorem~\ref{thmmfdmps}(a)
will be called the \emph{canonical} smooth manifold structure on $C^\ell(M,N)$.
\end{defn}
\begin{rem}
If $M$ and~$N$ are as in Theorem~\ref{thmmfdmps}, $L$ is a $C^k$-manifold (possibly with rough
boundary) and $h\colon L\times M\to N$ is a map such that $h_x:=h(x,\cdot)\in C^\ell(M,N)$
for all $x\in L$, then~$h$ is $C^{k,\ell}$ if and only if the map
\[
h^\vee\colon L\to C^\ell(M,N), \;\, x\mto h_x
\]
is~$C^k$ (as
Theorem~\ref{thmmfdmps}(a) applies to the map $g:= h^\vee$ with $g^\wedge=h$).
\end{rem}
\begin{lem}\label{inCkopensub}
Let $N$ be a smooth manifold admitting a local addition,
$\ell\in\N_0\cup\{\infty\}$, and
$M$ be a compact smooth manifold $($possibly with rough boundary$)$.
If $W\sub N$ is an open subset, then $W$ admits a local addition
and the smooth manifold structure induced by $C^\ell(M,N)$
on its open subset $C^\ell(M,W)=C^\ell(M,N)\cap\lfloor M,W\rfloor$
coincides with the canonical
manifold structure on $C^\ell(M,W)$.
\end{lem}
\begin{prf}
If $\Sigma \colon U \to N$ is a local addition for~$N$
on an open neighborhood $U \sub TN$ of the $0$-section, then
$U' :=\{ (\pi_N(v),\Sigma(v))\colon v\in U\}$
is open in $N\times N$ and $\theta:=(\pi_N,\Sigma)\colon U \to U'$
is a $C^\infty$-diffeomorphism. Then $U_W:=U\cap TW$ is an open neighborhood of the $0$-section
in~$TW$, $U_W'= \{ (\pi_W(v), \Sigma_W(v))\colon v\in U_W\}=\theta(U\cap TW)$
is open in $N\times N$ (hence also in $W\times W$),
and $\Sigma_W:=(\pi_W,\Sigma_W)\colon U_W\to U_W'$ is a $C^\infty$-diffeomorphism
(as it coincides with $\theta|_{U_W}$). Thus~$\Sigma_W$ is a local addition for~$W$.
As the charts for $C^\ell(M,W)$ constructed using $\Sigma_W'$
coincide with restrictions of the charts $\phi_f$ for $C^\ell(M,N)$
(as in the proof of Theorem~\ref{thmmfdmps}) to open subsets,
the canonical manifold structure on $C^\ell(M,W)$ coincides with
the smooth manifold structure as an open submanifold of $C^\ell(M,N)$.
\end{prf}
\begin{lem}\label{ClmapCk}
Let $N_1$ and~$N_2$ be smooth manifolds admitting local additions,
$k,\ell\in\N_0\cup\{\infty\}$
and $M$ be a compact $C^\ell$-manifold which may have a rough boundary.
If $h\colon N_1\to N_2$ is a $C^{k+\ell}$-map,
then
\[
C^\ell(M,h)\colon C^\ell(M,N_1)\to C^\ell(M,N_2),\;\, f \mto h\circ f
\]
is a $C^k$-map. Notably, $C^\ell(M,h)$ is smooth whenever $h$ is smooth.
\end{lem}
\begin{prf}
By the Exponential Law (as in Theorem~\ref{thmmfdmps}(a)), we need only show that
$(C^\ell(M,h))^\wedge(f,x)=h(f(x))\in N_2$ is $C^{k,\ell}$
in $(f,x)\in C^\ell(M,N_1)\times M$. Now
$f(x)$ is $C^{k,\ell}$ in $(f,x)$, by Theorem~\ref{thmmfdmps}(c).
Hence $C^\ell(M,h)^\wedge$ is $C^{k,\ell}$, using the Chain Rule
(Proposition~\ref{Cklchainmfd}).
\end{prf}
%
%
\begin{lem}\label{fstarCkmfd}
Let
$N_1$ and~$N_2$ be smooth manifolds admitting local additions,
$k,\ell\in \N_0\cup\{\infty\}$
and $M$ be a compact $C^{k+\ell}$-manifold $($possibly with rough boundary$)$.
If $g\colon M\times N_1\to N_2$ is a $C^{k+\ell}$-map,
then
\[
g_*\colon C^\ell(M,N_1)\to C^\ell(M,N_1),\;\, f\mto g\circ (\id_M,f) 
\]
is a $C^k$-map. Notably, $g_*$ is smooth whenever $g$ is smooth.
\end{lem}
\begin{prf}
By the Exponential Law (as in Theorem~\ref{thmmfdmps}(a)), we need only show that
$(g_*)^\wedge(f,x)=g(x,f(x))\in N_2$ is $C^{k,\ell}$
in $(f,x)\in C^\ell(M,N_1)\times M$. Now $x$ is $C^{k,\ell}$
in $(f,x)$ (being $C^{k+\ell}$)
and $f(x)$ is $C^{k,\ell}$ in $(f,x)$, by Theorem~\ref{thmmfdmps}(c).
Hence $(g_*)^\wedge$ is $C^{k,\ell}$, using the Chain Rule
(Proposition~\ref{Cklchainmfd})
and Remark~\ref{rem-simpl-Ckl}(e).
\end{prf}
%
%
%
\begin{lem}\label{prodasmfd}
For $j\hspace{-.2mm}\in\hspace*{-.2mm}\{1,2\}$, let
$N_j$ be a $C^\infty$\hspace*{-.2mm}-manifold which
admits a local addition. Let~$M$ be a compact $C^\ell$-manifold $($possibly with rough boundary$)$
with $\ell\in \N_0\cup\{\infty\}$.
Then $\Phi\colon C^\ell(M,N_1\times N_2)\to C^\ell(M,N_1)\times C^\ell(M,N_2)$,
$f\mto (\pr_1\circ f,\pr_2\circ f)$ is a
$C^\infty$-diffeomorphism.
\end{lem}
\begin{prf}
We know that $\Phi$ is a homeomorphism (Lemma~\ref{prodmapsp}).
By Lemma~\ref{ClmapCk}, the maps $C^\ell(M,\pr_j)$ are smooth for $j\in\{1,2\}$,
whence also the bijection~$\Phi$ is smooth. The inverse
$\Phi^{-1}\colon C^\ell(M,N_1)\times C^\ell(M,N_2)\to C^\ell(M,N_1\times N_2)$
will be smooth if we can show that
\[
(\Phi^{-1})^\wedge\colon C^\ell(M,N_1)\times C^\ell(M,N_2)\times M\to N_1\times N_2
\]
is $C^{\infty,\ell}$,
by the Exponential Law in Theorem~\ref{thmmfdmps}(a). But
\[
(\Phi^{-1})^\wedge(f,g,x)=(f(x),g(x))
\]
is $C^{\infty,\ell}$ as its components are evaluation maps
as in Theorem~\ref{thmmfdmps}(c).
\end{prf}
\begin{small}
\subsection*{Exercises for Section~\ref{sec-mfdmps}}

\begin{exer}
Let
$N$ be a smoothy paracompact Banach manifold, $S$
be a submanifold of~$N$ which is a closed subset,
$\ell\in \N_0\cup\{\infty\}$, and
$M$ be a compact $C^\ell$-manifold
which may have a rough boundary.
Using a local addition for~$N$ which is adapted to~$S$
as in Exercise~\ref{submfd-spray}, show that
$C^\ell(M,S)$ is a submanifold
of $C^\ell(M,N)$.
\end{exer}

\end{small}
\section{Fine box products of manifolds}\label{sec-box-mfd}
As a tool, we create
manifold structures
on countable direct products
of manifolds.
Such ``fine box products''
are useful, for instance,
for the study of manifolds of mappings on non-compact manifolds
(see Section~\ref{sec-mfd-map-sigma}).
In the case of finite products,
we get the usual product manifold.
In the case of infinite products,
the topology we use may 
be finer than the product topology.
Our setting is as follows:
\begin{numba}
Let $\K\in\{\R,\C\}$,
$r\in\N_0\cup\{\infty,\omega\}$,
$I$ be a countable non-empty set
and $M_i$, for $i\in I$, be
a $C^r_\K$-manifold
modeled on locally convex spaces
which is not necessarily pure.
Consider the cartesian product
$M:=\prod_{i\in I}M_i$ as a set.
\end{numba}
\begin{numba}\label{the-fine-top}
For $\phi:=(\phi_i)_{i\in I}$ ranging through the
families of charts $\phi_i\colon U_i\to V_i\sub E_i$
of~$M_i$ such that $0\in V_i$,
we endow $E_\phi:=\bigoplus_{i\in I}E_i$
with the locally convex direct sum topology.
Then $\bigoplus_{i\in I}V_i:=
E_\phi\cap\prod_{i\in I}V_i$
is an open $0$-neighborhood in~$E_\phi$
(see Remark~\ref{firstremsums}(a)),
and we give it the topology induced by~$E_\phi$.
Let
$\cO_{\fbx}$ be the final topology on~$M$
with respect to the mappings
\begin{equation}\label{para-box}
\Theta_\phi \colon \bigoplus_{i\in I} V_i \to M,
\;\, (x_i)_{i\in I}\mto(\phi_i^{-1}(x_i))_{i\in I},\vspace{-.7mm}
\end{equation}
for all $\phi$ as before, and
$U_\phi:=\Theta_\phi(V_\phi)$.
Thus
\begin{equation}\label{U-phi-concrete}
U_\phi=\Big\{(y_i)_{i\in I}\in \prod_{i\in I}U_i\colon
\mbox{$y_i\not=\phi_i^{-1}(0)$ for only finitely many $i\in I$}\Big\}.\vspace{-.7mm}
\end{equation}
\end{numba}
\noindent
We call $\cO_{\fbx}$
the \emph{fine box topology } on~$M$.
The fine box topology is Hausdorff
since it is finer than the product topology.
In fact, consider the projection $\pr_i\colon M\to M_i$
onto the $i$th component for $i\in I$.
Using the continuous linear projection
$\pi_i\colon E_\phi\to E_i$ onto the $i$th component,
we deduce from the continuity of $\pr_i\circ\, \Theta_\phi=\phi_i^{-1}\circ \pi_i|_{V_\phi}$
for each $\phi$
that $\pr_i$ is continuous.
\begin{numba}\label{basics-fine-box}
Let $\phi$ be as before and consider an analogous family~$\psi$
of charts
$\psi_i\colon R_i\to S_i\sub F_i$.
If $\phi_i^{-1}(0)=\psi_i^{-1}(0)$ for all
but finitely many $i\in I$, then
\[
(\Theta_\phi)^{-1}(U_\phi\cap U_\psi)=\bigoplus_{i\in I}\phi_i(U_i\cap R_i),\vspace{-.7mm}
\]
which is an
open $0$-neighbourhood in $\bigoplus_{i\in I}E_i$.
By Proposition~\ref{diff-sum}, the transition map
\[
(\Theta_\phi)^{-1}\circ \, \Theta_\psi\colon \bigoplus_{i\in I}\psi_i(U_i\cap R_i)
\to \bigoplus_{i\in I}\phi_i(U_i\cap R_i),\;\,
(x_i)_{i\in I}\mto ((\phi_i\circ \psi_i^{-1})(x_i))_{i\in I}\vspace{-.7mm}
\]
is $C^r_\K$
and in fact a $C^r_\K$-diffeomorphism,
and hence a homeomorphism,
since $\Theta_\psi^{-1}\circ \, \Theta_\phi$
is the inverse map.
If $\phi_i^{-1}(0)\not=\psi_i^{-1}(0)$ for
infinitely many $i\in I$,
then $(\Theta_\phi)^{-1}(U_\phi\cap U_\psi)=\emptyset$
and the transition map trivially is a homeomorphism.
Hence
$U_\phi=\Theta_\phi(V_\phi)$ is open in $(M,\cO_{\fbx})$
for all~$\phi$
and $\Theta_\phi$ is a homeomorphism onto its image
(see Exercise~\ref{topforbun}).
By the preceding,
the maps $\Phi_\phi:=(\Theta_\phi|^{U_\phi})^{-1}\colon U_\phi\to V_\phi\sub E_\phi$
are $C^r_\K$-compatible 
and hence form an atlas for a $C^r_\K$-manifold
structure on~$M$.
\end{numba}
\begin{defn}
We write $M^{\fbx}$ for~$M$, endowed with the topology~$\cO_{\fbx}$
and the $C^r_\K$-manifold structure just described,
and call $M^{\fbx}$ the \emph{fine box product}
of the family $(M_i)_{i\in I}$
of $C^r_\K$-manifolds.
\end{defn}
\begin{lem}\label{maps-between-box}
Let $I$ be a countable set, $r\in \N_0\cup\{\infty\}$
and $f_i\colon M_i\to N_i$
be a $C^r_\K$-map between
$C^r_\K$-manifolds
for $i\in I$.
Then the following map is $C^r_\K$:
\[
f:=\prod_{i\in I}f_i\colon {\prod_{i\in I}}^{\fbx}M_i\to {\prod_{i\in I}}^{\fbx}N_i,\;\,
(x_i)_{i\in I}\mto (f_i(x_i))_{i\in I}.
\]
\end{lem}
\begin{prf}
Let $x=(x_i)_{i\in I}$ be an element of $M:=\prod_{i\in I}^{\fbx}M_i$,
$\psi_i\colon P_i\to Q_i\sub F_i$
be a chart for $N_i$ around $y_i :=f(x_i)$ for $i\in I$
such that $\psi_i(y_i)=0$,
and \break $\phi_i\colon U_i\to V_i\sub E_i$
be a chart of $M_i$ around $x_i$ such that $\phi_i(x_i)=0$.
After shrinking $U_i$,
we may assume that $f_i(U_i)\sub P_i$.
Let $\phi:=(\phi_i)_{i \in I}$
and $\Phi_\phi\colon U_\phi\to \bigoplus_{i\in I}V_i\sub \bigoplus_{i\in I}E_i$
be the corresponding chart of $M$ around~$x$,
as in \ref{basics-fine-box}
Likewise, consider $\psi:=(\psi_i)_{i\in I}$
and the corresponding chart $\Phi_\psi\colon U_\psi\to \bigoplus_{i\in I}Q_i\sub\bigoplus_{i\in I}F_i$
of $\prod_{i\in I}^{\fbx} N_i$ around $f(x)$.
Then $\Phi_\psi\circ f\circ (\Phi_\phi)^{-1}$ is the map
\[
\oplus_{i\in I}(\psi_i\circ f_i\circ \phi_i^{-1})\colon
\bigoplus_{i\in I}V_i\to \bigoplus_{i\in }Q_i,\;\,
(z_i)_{i\in I}\mto ((\psi_i\circ f_i\circ \phi_i^{-1})(z_i))_{i\in I},
\]
which is $C^r_\K$ by Proposition~\ref{diff-sum}.
\end{prf}
\begin{lem}\label{finebox-product}
Let $I$ be a countable set, $r\in \N_0\cup\{\infty\}$,
and $M_{i,1}$, $M_{i,2}$
be $C^r_\K$-manifolds
for $i\in I$. For $j\in\{1,2\}$,
let $\pr_{i,j}\colon M_{i,1}\times M_{i,2}\to M_{i,j}$
be the projection onto the $j$th component.
Then the following map is a $C^r_\K$-diffeomorphism:
\[
\Psi:=\left(\prod_{i \in I}\pr_{i,1},\prod_{i \in I}\pr_{i,2}\right)\colon 
{\prod_{i\in I}}^{\fbx}(M_{i,1}\times M_{i,2})\to
\Big({\prod_{i\in I}}^{\fbx}M_{i,1}\Big)\times\Big({\prod_{i\in I}}^{\fbx}M_{i,2}\Big).
\]
\end{lem}
\begin{prf}
Let $x=(x_i)_{i\in I}\in M_1:=\prod_{i\in I}^{\fbx}M_{i,1}$
and $y=(y_i)_{i\in I}\in M_2:=\prod_{i\in I}^{\fbx}M_{i,2}$.
For each $i\in I$, let $\phi_i\colon U_i\to V_i\sub E_i$
be a chart of $M_{i,1}$ around $x_i$
such that $\phi_i(x_i)=0$;
set $\phi:=(\phi_i)_{i\in I}$
and let $\Phi_\phi\colon U_\phi\to\bigoplus_{i\in I}V_i\sub\bigoplus_{i\in I}E_i$
be the corresponding chart of $M_1$.
Likewise, let $\psi_i\colon P_i\to Q_i\sub F_i$
be a chart of $M_{i,2}$ around $y_i$
such that $\psi_i(y_i)=0$;
set $\psi:=(\psi_i)_{i\in I}$
and let $\Phi_\psi\colon U_\psi\to\bigoplus_{i\in I}Q_i\sub\bigoplus_{i\in I}F_i$
be the corresponding chart of $M_2$.
Let $\theta:=(\phi_i\times \psi_i)_{i\in I}$
and $\Phi_\theta\colon U_\theta\to\bigoplus_{i\in I}(V_i\times Q_i)\sub \bigoplus_{i\in I}
(E_i\times F_i)$
be the corresponding chart of $M:=\prod_{i \in I}^{\fb}(M_{i,1}\times M_{i,2})$.
By Lemma~\ref{firstlasum}(c), the following map is an isomorphism of topological
vector spaces:
\[
\alpha\colon \bigoplus_{i\in I}(E_i\times F_i)\to\bigoplus_{i\in I}E_i\times\bigoplus_{i\in I}F_i,
\;\,
(v_i,w_i)_{i\in I}\mto((v_i)_{i\in I},(w_i)_{i\in I}).
\]
It restricts to a $C^r_\K$-diffeomorphism
from $W:=\bigoplus_{i\in I}(V_i\times Q_i)$
onto $V\times Q$ with $V:=\bigoplus_{i\in I}V_i$,
$Q:=\bigoplus_{\in I}Q_i$.
Now $\Psi(U_\theta)=U_\phi\times U_\psi$
and $(\Phi_\phi\times \Phi_\psi)\circ \Psi\circ \Phi_\theta^{-1}=\alpha|_W^{V\times Q}$
is a $C^r_\K$-diffeomorphism.
Thus $\Psi$ is a bijective map which is a
local $C^r_\K$-diffeomorphism at each point, and thus $\Psi$
is a $C^r_\K$-diffeomorphism.
\end{prf}
\section{Manifolds of mappings on non-compact manifolds}\label{sec-mfd-map-sigma}
We now construct smooth manifold structures
on sets of mappings on \break $\sigma$-compact, locally compact manifolds
(which may have a rough boundary).
Using countable locally finite covers by
compact full submanifolds, much of the construction can be reduced
to the case of compact domains treated in Section~\ref{sec-mfdmps},
with the help of embeddings in suitable fine box products.
The manifolds will
be modeled on certain locally convex spaces $\Gamma_{f,c}$
isomorphic to the space $\Gamma_{C^\ell_c}(f^*(TN))$
of compactly supported $C^\ell$-sections
in the pullback-bundle $f^*(TN)$,
see Lemma~\ref{details-gamma-nonc}. We prove the following theorem.
\begin{thm}\label{map-mfd-noncomp}
Let $N$ be a smooth manifold admitting
a local addition,
$\ell\in \N_0\cup\{\infty\}$
and $M$
be a $\sigma$-compact, locally compact
$C^\ell$-manifold which may have a rough boundary.
Let $\pi_{TN}\colon TN\to N$ be the bundle projection.
Then $C^\ell(M,N)$
admits a unique smooth manifold structure
modeled on the set $\cE:=\{\Gamma_{f,c}\colon f\in C^\ell(M,N)\}$
of locally convex spaces such that,
for each $f\in C^\ell(M,N)$
and local addition $\Sigma\colon TN\supseteq U \to N$
of~$N$, the map
\[
\Gamma_{f,c}\cap C^\ell(M,U)\to C^\ell(M,N),\;\,
\tau\mto \Sigma\circ\tau
\]
is a $C^\infty$-diffeomorphism onto an open subset of $C^\ell(M,N)$.
For each locally finite countable family $(K_i)_{i\in I}$ of compact full submanifolds
$K_i$ of~$M$ whose interiors cover~$M$, we have an embedding of smooth manifolds
\begin{equation}\label{rho-in-the box}
\rho\colon C^\ell(M,N)\to{\prod_{i\in I}}^{\fbx}C^\ell(K_i,N),\;\;
f\mto (f|_{K_i})_{i\in I}.
\end{equation}
\end{thm}
We also show:
\begin{prop}\label{the-tangents}
In the situation of Theorem~{\rm\ref{map-mfd-noncomp}},
the point evaluation $\ve_x\colon C^\ell(M,N)\to N$
is smooth for each $x\in M$. For each $f\in C^\ell(M,N)$,
we have $(T\ve_x (v))_{x\in M}\in\Gamma_{f,c}$
for each $v\in T_f(C^\ell(M,N))$ and the map
\[
\Theta\colon T_f(C^\ell(M,N))\to\Gamma_{f,c},\;\;
v\mto (T\ve_x(v))_{x\in M}
\]
is an isomorphism of topological vector spaces.
\end{prop}
Some auxiliary results are helpful. We use notation as in~\ref{def-gamma-f}.
\begin{lem}\label{details-gamma-nonc}
In the situation of Theorem~{\rm\ref{map-mfd-noncomp}},
let $f\in C^\ell(M,N)$.
Then, for each compact subset $K$ of~$M$,
\[
\Gamma_{f,K}:=
\{\tau\in \Gamma_f\colon (\forall x\in M\setminus K)\;
\tau(x)=0
\in T_{f(x)}N\}
\]
is a closed vector subspace of $\Gamma_f$
and hence a locally convex space in the induced topology.
The locally convex direct limit topology on
$\Gamma_{f,c}=\bigcup_K\Gamma_{f,K}$ is Hausdorff.
\end{lem}
\begin{prf}
By Proposition~\ref{eval-pb},
the point evaluation $\ev_x\colon \!\Gamma_f\!\to\! TN$, $\tau\!\mto\! \tau(x)$
is~$C^\infty$ for $x\in M$,
whence it is continuous linear as a map to $T_{f(x)}N$.
Thus $\Gamma_{f,K}$ is a closed vector subspace of $\Gamma_f$.
The locally convex direct limit topology on $\Gamma_{f,c}$ makes the linear
inclusion map $\Gamma_{f,c}\to \Gamma_f$
continuous
(see Remark~\ref{firstremlcxDL}(d)). Since $\Gamma_f$ is Hausdorff,
so is $\Gamma_{f,c}$.
\end{prf}
Henceforth, we endow $\Gamma_{f,c}$ with the locally convex direct
limit topology.
%
%
\begin{lem}\label{ops-with-Gamma}
Let
$N$ be a smooth manifold, $\ell \in \N_0\cup\{\infty\}$ and
$M$ be a $\sigma$-compact, locally compact $C^\ell$-manifold
which may have a rough boundary.
For $f\in C^\ell(M,N)$, the following holds:
\begin{description}[(D)]
\item[\rm(a)]
The bilinear map $C^\ell(M,\R)\times \Gamma_f\to\Gamma_f$,
$(h,\tau)\mto h\tau$ is continuous, where $(h\tau)(x)=h(x)\tau(x)$.
\item[\rm(b)]
If $K\sub M$
is a compact full submanifold,
then the linear map
$\Gamma_f\to \Gamma_{f|_K}$, $\tau\mto \tau|_K$ is continuous.
\item[\rm(c)]
If $K$
a compact full submanifold of $M$
and $L\sub K$ compact,
then
$r\colon \Gamma_{f,L}\to \Gamma_{f|_K,L}$,
$\tau\mto \tau|_K$ is an isomorphism of topological vector spaces.
\end{description}
\end{lem}
\begin{prf}
(a) Using the isomorphism $\Gamma_{C^\ell}(f^*(TN))\to \Gamma_f$
from Lemma~\ref{sections-pb},
this follows from the fact that $\Gamma_{C^\ell}(f^*(TN))$
is a topological $C^\ell(M,\R)$-module (see Lemma~\ref{secs-top-module}).

(b) The mapping in contention is a restriction
of the restriction map $C^\ell(M,TN)\to C^\ell(K,TN)$,
which is continuous (see Remark~\ref{earlyremCk}(a)).

(c) In view of Lemma~\ref{sections-pb},
this follows from
Lemma~\ref{compactinopen-restr}.
\end{prf}
\noindent
\emph{Proof of Theorem}~\ref{map-mfd-noncomp}.
Let $(K_i)_{i\in I}$ be a
locally finite family of 
compact, full submanifolds $K_i$
of~$M$ whose interiors cover~$M$.
Since $M$ is $\sigma$-compact, after omitting empty sets
we may assume that~$I$ is countable.
The map
\[
\rho\colon C^\ell(M,N)\to\prod_{i\in I}C^\ell(K_i,N),\;\,
f\mto (f|_{K_i})_{i\in I}
\]
is injective and its image $\im(\rho)$ equals
\begin{equation}\label{get-image}
\Big\{(f_i)_{i\in I}\in
\prod_{i\in I}C^\ell(K_i,N)
\colon (\forall i,j\in I)\,(\forall x\in K_i\cap K_j)\,
f_i(x)=f_j(x)\Big\}.
\end{equation}
In fact, the inclusion ``$\sub$''
is obvious. If $(f_i)_{i\in I}$
is in the set on the right-hand side,
then a piecewise definition,
$f(x):=f_i(x)$ if $x\in K_i$,
gives a well-defined function $f\colon M\to N$
which is $C^\ell$ since $f|_{(K_i)^0}=f_i|_{(K_i)^0}$
is $C^\ell$ for each $i\in I$.
Then $\rho(f)=(f_i)_{i\in I}$.\\[2mm]
For each $i\in I$,
endow $C^\ell(K_i,N)$
with the canonical smooth manifold structure,
as in Theorem~\ref{thmmfdmps}, modeled
on the set $\{\Gamma_f\colon f\in C^\ell(K_i,N)\}$
of the locally convex spaces
$\Gamma_f:=\{\tau\in C^\ell(K_i,TN)\colon \pi_{TN}\circ
\tau=f\}$
for $f\in C^\ell(K_i,N)$.
Let
$\Sigma\colon TN\supseteq U\to N$
be a local addition for~$N$;
as in Definition~\ref{defn-loc-add}, write
$U':=\{(\pi_{TN}(v),\Sigma(v))\colon v\in U\}$
and $\theta:=(\pi_{TN}|_U,\Sigma)\colon U\to U'$.
For $f\in C^\ell(K_i,N)$,
consider
$O_f:=\Gamma_f\cap C^\ell(K_i,U)$,
$O_f':=\{g\in C^\ell(K_i,N)\colon (f,g)\in C^\ell(K_i,U')\}$,
and $\psi_f\colon O_f\to O_f'$, $\tau\mto \Sigma\circ \tau$
as in Theorem~\ref{thmmfdmps} and its proof.
For $f\in C^\ell(M,N)$,
let
$\Gamma_{f,c}$
be the set of all $\tau\in C^\ell(M,TN)$
such that $\pi_{TN}\circ
\tau=f$ and
\[
\{x\in M\colon \tau(x)\not= 0_{f(x)}\in T_{f(x)}N\}
\]
is relatively compact in~$M$.
Define $O_f:=\Gamma_{f,c}\cap C^\ell(M,U)$
and let~$O_f'$ be the set of all
$g\in C^\ell(M,N)$ such that
\[
(f,g)\in C^\ell(M,U')\;\;\mbox{and}\;\; \mbox{$g|_{M\setminus K}=f|_{M\setminus K}$
for some compact subset $K\sub M$.}
\]
Then $\psi_f\colon O_f\to O_f'$, $\tau\mto \Sigma\circ \tau$
is a bijection with $(\psi_f)^{-1}(g)=\theta^{-1}\circ (f,g)$.
The linear map
\[
s\colon \Gamma_{f,c}\to \bigoplus_{i\in I}\Gamma_{f|_{K_i}},\quad
\tau\mto (\tau|_{K_i})_{i\in I}
\]
is continuous on $\Gamma_{f,L}$ for each
compact subset $L\sub M$ (see Lemma~\ref{ops-with-Gamma}(b))
and hence continuous on the locally convex direct
limit $\Gamma_{f,c}$ by its universal property (see Remark~\ref{firstremlcxDL}(d)).
As above,
\begin{equation}\label{id-image-2}
\im(s)\hspace*{-.2mm}=\hspace*{-.2mm}
\Big\{(\tau_i)_{i\in I}\!\in\bigoplus_{i\in I}\Gamma_{f|_{K_i}}\!\!\colon\!
(\forall i,j\in I)(\forall x\in K_i\cap K_j)\, \tau_i(x)=\tau_j(x)\Big\}\! ,\!\!\!
\end{equation}
which is a closed vector subspace of $\bigoplus_{i\in I}
\Gamma_{f|_{K_i}}$.
We now show that $s$ is a homeomorphism onto its
image. In fact, the corestriction to the image
admits a continuous
linear right inverse.
To see this,
pick a $C^\ell$-partition of unity
$(h_i)_{i\in I}$ on~$M$
subordinate to $(K_i^0)_{i\in I}$;
then $L_i:=\Supp(h_i)$ is a closed subset of $K_i$
and thus compact.
The multiplication operator
$\beta_i\colon \Gamma_{f|_{K_i}}\to \Gamma_{f|_{K_i},L_i}$,
$\tau\mto h_i\tau$
is continuous linear (by
Lemma~\ref{ops-with-Gamma}\,(a)).
Moreover, the restriction operator
$s_i\colon \Gamma_{f,L_i}\to\Gamma_{f|_{K_i},L_i}$
is an isomorphism of topological vector spaces
(Lemma~\ref{ops-with-Gamma}\,(c)).
Thus $s_i^{-1}\circ \beta_i\colon \Gamma_{f|_{K_i}}\!\to\hspace{.2mm}\Gamma_{f,L_i}\sub \Gamma_{f,c}$
is a continuous linear map. By Lemma~\ref{firstlasum}(a),
also the linear map
\[
\sigma\colon \bigoplus_{i\in I}\Gamma_{f|_{K_i}}\to \Gamma_{f,c},\;\,
(\tau_i)_{i\in I}\mto \sum_{i\in I}(s_i^{-1}\circ \beta_i)(\tau_i)
\]
is continuous. Hence $\sigma|_{\im(s)}$
is continuous and linear. We easily verify that
$s\circ \sigma|_{\im(s)}=\id_{\im(s)}$.\\[2mm]
Abbreviate $\phi_i:=(\psi_{f|_{K_i}})^{-1}$
and $\phi:=(\phi_i)_{i\in I}$.
We now use the $C^\infty$-diffeomorphism
\[
\Theta_\phi\colon
\bigoplus_{i\in I} O_{f|_{K_i}}\! \to \,U_\phi,\;\;
(\tau_i)_{i\in I}\mto (\phi_i^{-1}(\tau_i))_{i\in I}=
(\Sigma \circ \tau_i)_{i\in I}
\]
as in (\ref{para-box}) and~\ref{basics-fine-box},
the inverse of which is the chart
\[
\Phi_\phi\colon U_\phi\to
\bigoplus_{i\in I}O_{f|_{K_i}}, \;\,
(g_i)_{i\in I}\mto
(\phi_i(g_i))_{i\in I}
\]
of $\prod_{i\in I}^{\fbx}C^\ell(K_i,N)$
around $(f|_{K_i})_{i\in I}$.
For $(\tau_i)_{i\in I}\in\bigoplus_{i\in I}O_{f|_{K_i}}$,
we have
\[
\Theta_\phi((\tau_i)_{i\in I})
\in \im(\rho)\;\;\Leftrightarrow\;\;
(\tau_i)_{i\in I}\in  \im(s).
\]
In fact, for $i,j\in I$ and $x\in K_i\cap K_j$
we have
$\Sigma(\tau_i(x))=\Sigma(\tau_j(x))$
if and only if $\tau_i(x)=\tau_j(x)$, from
which the assertion follows in view of~(\ref{get-image})
and~(\ref{id-image-2}).
Thus
\[
\Phi_\phi(\im(\rho)\cap U_\phi)
=\im(s)\cap \bigoplus_{i\in I}O_{f|_{K_i}},
\]
showing that $\im(\rho)$
is a submanifold of $\prod_{i\in I}^{\fbx}C^\ell(K_i,N)$.
Let
\[
\Psi_\phi\colon \im(\rho)\cap U_\phi
\to \im(s)\cap \bigoplus_{i\in I}O_{f|_{K_i}},\;\;
(g_i)_{i\in I}\mto \Phi_\phi((g_i)_{i\in I})
\]
be the corresponding submanifold
chart for $\im(\rho)$.
Then
\[
\rho(O_f')=\im(\rho)\cap U_\phi\quad\mbox{and}\quad
s(O_f)=\im(s)\cap \bigoplus_{i\in I}O_{f|_{K_i}}.
\]
Hence $\phi_f=s^{-1}\circ \Psi_\phi\circ \rho|_{O_f'}\colon O_f'\to O_f$
is a chart for the smooth manifold structure on $C^\ell(M,N)$
modeled on $\cE$ (the set of all $\Gamma_f$)
which makes
$\rho\colon C^\ell(M,N)\to \im(\rho)$
a $C^\infty$-diffeomorphism.
Note that the smooth manifold structure
on $C^\ell(M,N)$ modeled on~$\cE$ 
that makes $\rho$ a $C^\infty$-diffeomorphism
is uniquely determined by these properties.
Thus, it is independent of the choice of~$\Sigma$.
On the other hand, the $\phi_f$
form a $C^\infty$-atlas for a given
local addition~$\Sigma$.
As the definition of the $\phi_f$
does not involve the cover $(K_i)_{i\in I}$,
the smooth manifold structure just constructed
is independent of the choice of $(K_i)_{i\in I}$. \medskip\qed

\emph{Proof of Proposition}~\ref{the-tangents}.
With notation as in the preceding proof, given $x\in M$ we have $x\in K_i$ for some $i\in I$.
The restriction map $\rho_i\colon C^\ell(M,N)\to C^\ell(K_i,N)$
is smooth (being a component of~$\rho$).
Moreover, the point evaluation $\ev_x\colon C^\ell(K_i,N)\to N$, $f\mto f(x)$
is smooth, as shown in Theorem~\ref{thmmfdmps}(c).
Hence $\ve_x=\ev_x\circ \rho_i$ is smooth.
We can now repeat the proof of Theorem~\ref{thmmfdmps}(d) verbatim,
replacing $\Gamma_f$ with $\Gamma_{f,c}$. \qed
\begin{prop}\label{pushf-mfd-maps}
Let
$N_1$ and $N_2$ be smooth manifolds admitting local additions
and $k,\ell\in \N_0\cup\{\infty\}$.
\begin{description}[(D)]
\item[\rm(a)]
If $M$ is a $\sigma$-compact, locally compact $C^{k+\ell}$-manifold which may have a rough boundary
and $g\colon M\times N_1\to N_2$ a $C^{k+\ell}$-map, then the mapping
$g_*\colon C^\ell(M,N_1)\to C^\ell(M,N_2)$, $f\mto g\circ (\id_N,f)$ is~$C^k$.
\item[\rm(b)]
If $M$ is a $\sigma$-compact, locally compact $C^\ell$-manifold which may have a rough boundary
and $h\colon N_1\to N_2$ a $C^{k+\ell}$-map,
then the mapping $C^\ell(M,h)\colon C^\ell(M,N_1)\to C^\ell(M,N_2)$,
$f\mto h\circ f$ is~$C^k$.
\end{description}
\end{prop}
\begin{prf}
(a) Let $(K_i)_{i\in I}$ be a countable locally finite family of compact full submanifolds
of~$M$ whose interiors cover~$M$. Then
$H_i\colon
C^\ell(K_i,N_1)\to C^\ell(K_i,N_2)$, $f\mto g\circ (\id_{K_i},f)$
is $C^k$ for each $i\in I$, by Lemma~\ref{fstarCkmfd}.
Hence also $H:=\prod_{i\in I}H_i\colon \prod_{i\in I}^{\fbx}
C^\ell(K_i,N_1)\to\prod_{i\in I}^{\fbx} C^\ell(K_i,N_2)$
is $C^k$, by Lemma~\ref{maps-between-box}. Now
\[
\rho_j\colon C^\ell(M,N_j)\to{\prod_{i\in I}}^{\fbx}C^\ell(K_i,N_j),\;\;
f\mto (f|_{K_i})_{i\in I}
\]
is an embedding of smooth manifolds for $j\in\{1,2\}$.
Hence $g_*$ will
be $C^k$ if $\rho_2\circ g_*$ is $C^k$.
But this is the case as $\rho_2\circ g_*=H\circ \rho_1$.

(b) can be proved like~(a), replacing $g_*$ with $C^\ell(M,h)$
and taking $H_i:=C^\ell(K_i,h)$ as in Lemma~\ref{ClmapCk}.
\end{prf}
\begin{prop}\label{mapmfd-compa-prod}
Let $M$ be a $\sigma$-compact, locally compact $C^\ell$-manifold which may have a rough boundary,
with $\ell\in \N_0\cup\{\infty\}$. Let
$N_1$ and $N_2$ be smooth manifolds admitting local additions.
Then also $N_1\times N_2$
admits a local addition. Let $\pr_j\colon N_1\times N_2\to N_j$
be the projection onto the $j$th component for $j\in \{1,2\}$.
Then the following map is a $C^\infty$-diffeomorphism:
\[
\Phi:=(C^\ell(M,\pr_1),C^\ell(M,\pr_2))\colon C^\ell(M,N_1\times N_2)\to
C^\ell(M,N_1)\times C^\ell(M,N_2).
\]
\end{prop}
\begin{prf}
By Proposition~\ref{pushf-mfd-maps}(b), $\Phi$ is smooth;
moreover, $\Phi$ is a bijection. To see that $\Phi^{-1}$ is continuous,
let $(K_i)_{i\in I}$ be a countable locally finite family of compact full submanifolds
of~$M$ whose interiors cover~$M$.
For $i\in I$, the map
\[
\Phi_i:=
(C^\ell(K_i,\pr_1), C^\ell(K_i, \pr_2))\colon C^\ell(K_i, N_1\times N_2)\to
C^\ell(K_i,N_1)\times C^\ell(K_i, N_2)
\]
is a $C^\infty$-diffeomorphism (see Lemma~\ref{prodasmfd}). Then
\[
h:=\prod_{i\in I}\Phi_i^{-1}\colon
{\prod_{i\in I}}^{\fbx} (C^\ell(K_i,N_1)\times C^\ell(K_i, N_2))\to
{\prod_{i\in I}}^{\fbx} C^\ell(K_i, N_1\times N_2)
\]
is smooth, by Lemma~\ref{maps-between-box}.
By Theorem~\ref{map-mfd-noncomp},
\[
\rho_j\colon C^\ell(M,N_j)\to {\prod_{i\in I}}^{\fbx}C^\ell(K_i,N_j),\;\;
f\mto (f|_{K_i})_{i\in I}
\]
is an embedding of smooth manifolds for $j\in\{1,2\}$
and also the mapping $\rho\colon C^\ell(M,N_1\times N_2)\to {\prod_{i\in I}}^{\fbx} C^\ell(K_i,N_1\times N_2)$,
$f\mto (f|_{K_i})_{i\in I}$ is an embedding.
Let
$\Psi$ be the $C^\infty$-diffeomorphism
\[
{\prod_{i\in I}}^{\fbx} (C^\ell(K_i,N_1)\times C^\ell(K_i, N_2))
\to \left( {\prod_{i\in I}}^{\fbx} C^\ell(K_i, N_1)\right)
\times
\left(
{\prod_{i\in I}}^{\fbx} C^\ell(K_i, N_2)
\right)
\]
from
Lemma~\ref{finebox-product}.
Then $\rho\circ \Phi^{-1}=\Psi\circ h\circ \Psi^{-1}\circ (\rho_1\times \rho_2)$
is smooth, entailing smoothness of~$\Phi^{-1}$.
\end{prf}
\section{Direct limits of finite-dimensional manifolds}
\label{dl-construction}
Consider an ascending sequence
$M_1\sub M_2\sub\cdots$ of finite-dimensional smooth manifolds,
such that the inclusion map $\lambda_{n,m}\colon M_m\to M_n$
is a smooth immersion for all $n,m\in\N$ such that $n\geq m$.
Let 
\[ d:=\sup\{\dim(M_n)\colon n\in\N\}\in \N_0\cup\{\infty\} \quad \mbox{ and } \quad 
\R^\infty:=\R^{(\N)}:=\bigoplus_{n\in \N}\R,\]
endowed with the locally convex direct sum topology.
In this section, we prove the following theorem:
\begin{thm}\label{dl-mfd}
There is a unique smooth manifold structure on 
\[ M:=\bigcup_{n\in\N}M_n\quad \mbox{ modeled on } \quad \R^d,\] with the following property:
For each $k\in \N_0\cup\{\infty\}$ and $C^k$-manifold
$N$ which may have a rough boundary, a map
$f\colon M\to N$ is $C^k$ if and only if $f|_{M_n}\colon M_n\to N$
is $C^k$ for each $n\in \N$.
For each $n\in\N$, the inclusion map $\lambda_n\colon M_n\to M$ is
smooth for the latter smooth manifold structure on~$M$. 
The topology underlying $M$ is the final topology with respect to the inclusion maps
$\lambda_n\colon M_n\to M$ for $n\in\N$.
\end{thm}
We shall use the following terminology:
\begin{defn}
Let $M$ be a smooth manifold
and $N\sub M$ be a submanifold.
A map $r\colon M\to N$
is called a \emph{smooth retraction}
from $M$ to~$N$
if $r$ is smooth and $r|_N=\id_N$.
\end{defn}
In the following lemma, for $n<m$ we identify
$\R^n$ with the vector subspace $\R^n\times\{0\}\sub\R^n\times\R^{m-n}=\R^m$.
\begin{lem}[Extension of charts]\label{extend-charts}
Let $M$ be an $m$-dimensional smooth manifold,
$N\sub M$ be an $n$-dimensional immersed submanifold and 
\break ${\phi\colon U_\phi\to V_\phi\sub\R^n}$
be a chart for~$N$
such that $U_\phi$ is relatively compact in~$M$
and smoothly contractible.
Then there exists a chart $\psi\colon U_\psi\to V_\psi\sub \R^m$
of $M$ with the following properties:
\begin{description}[(D)]
\item[\rm(a)]
$U_\psi$ is relatively compact in~$M$
and smoothly contractible;
\item[\rm(b)]
$N\cap U_\psi=U_\phi$,
$\R^n\cap V_\psi=V_\phi$,
and $\phi=\psi|_{U_\phi}$;
\item[\rm(c)]
There exists a smooth retraction from $U_\psi$ to~$U_\phi$.
\end{description}
\end{lem}
\begin{prf}
Let $\wb{U}_\phi$ be the closure of $U_\phi$ in~$N$.
Since $\wb{U}_\phi$ is compact and~$M$ Hausdorff,
the inclusion map $\wb{U}_\phi\to M$
is a topological embedding,
entailing that $M$ and $N$ induce the same topology
on $\wb{U}_\phi$ and also
$j|_{U_\phi}\colon U_\phi\to M$
is a topological embedding.
Since~$M$ is locally compact, we find a compact subset $K\sub M$
such that $\wb{U}_\phi\sub K^0$.
Since $K$ is covered by finitely many chart domains which are separable,
$K^0$
is separable. Since $U_\phi$ is open in $\wb{U}_\phi$,
there exists an open subset $W\sub M$ such that
\begin{equation}\label{thus-rel-closed}
U_\phi=\wb{U}_\phi \cap W.
\end{equation}
After replacing $W$ with $W\cap K^0$,
we may assume that $W\sub K^0$.
By~\eqref{thus-rel-closed},
$U_\phi$ is closed in~$W$.
The above reasoning shows that
$j|_{U_\phi}\colon U_\phi\to W$
is a topological embedding.
Also being an immersion, $j|_{U_\phi}\colon U_\phi\to W$
is an embedding of smooth manifolds.
Thus $U_\phi$ is a submanifold of~$W$ (see Proposition~\ref{embedding-submfd})
and a closed subset of~$W$.
By Theorem~\ref{thm-tubular},
there exists a $C^\infty$-diffeomorphism
$f\colon P\to Q$ from
an open subset $P\sub W$
with $U_\phi\sub P$
onto an open subset $Q$ of a smooth vector bundle
$E$ over~$U_\phi$, such that $Q$ contains
the zero-section and $f(x)=0_x\in E_x$
for each $x\in U_\phi$. Let $\pi\colon E\to U_\phi$
be the bundle projection.
By Corollary~\ref{contractible-trivial}, $E$ is trivializable;
we may therefore assume that $E=U_\phi\times\R^{m-n}$.
Since $(x,0)\in Q$, for each $x\in U_\phi$
there exist open subsets $U_x\sub U_\phi$ and $Y_x\sub \R^{m-n}$
such that $x\in U_\phi$, $0\in Y_x$, and $U_x\times Y_x\sub Q$.
After replacing $Y_x$ with an open convex $0$-neighborhood,
we may assume that $ty\in Y_x$ for all $y\in Y_x$ and $t\in [0,1]$.
After replacing $Q$ with $\bigoplus_{x\in U_\phi}(U_x\times Y_x)$
and $P$ with the preimage of this open set under~$f$,
we may assume that
\[
F(t,x,y):=(x,(1-t)y)\in Q\;\,\mbox{for all $(x,y)\in Q$ and $t\in [0,1]$.}
\]
The map $F\colon [0,1]\times Q\to Q$ is smooth.
Let $h\colon [0,1]\to[0,1]$ be a smooth function
such that, for some $\ve\in]0,1/2[$, $h(t)=0$ for all $t\in [0,\ve[$
and $h(t)=1$ for all $t\in \;]1-\ve,1]$.
Since $U_\phi$ is smoothly contractible,
there exists $x_0\in U_\phi$ and a smooth map $G\colon [0,1]\times U_\phi$
such that $G(0,x)=x$ and $G(1,x)=x_0$ for all $x\in U_\phi$.
Consider the map
\[
H\colon [0,1]\times Q\to Q,\;\,
(t,x,y)\mto \left\{
\begin{array}{rl}
F(h(2t),x,y) &\mbox{if $t\in [0,1/2]$;}\\
(G(h(2t-1),x),0) & \mbox{if $t\in [1/2,1]$.}
\end{array}
\right.
\]
Since $H(t,x,y)=(x,0)$ for all $t\in \;]1/2-\ve/2,1/2+\ve/2[$,
we easily see that $H$ is smooth. Moreover, $H(0,x,y)=(x,y)$
and $H(1,x,y)=(x_0,0)$ for all $(x,y)\in Q$.
Hence~$Q$ is smoothly contractible and thus so is $U_\psi:=P$.
As a consequence of (\ref{thus-rel-closed}), we have $N\cap U_\psi=U_\phi$.
Moreover, $U_\psi$ is relatively compact as $U_\psi\sub W\sub K$.
Finally, $V_\psi:=(\phi\times \id_{\R^{m-n}})(Q)\sub \R^m$
is an open set and $(\phi\times \id_{\R^{m-n}})\circ f\colon U_\phi\to V_\psi$
is a $C^\infty$-diffeomorphism;
moreover,
\[ V_\psi\cap (\R^n\times \{0\})=(\phi\times\id_{\R^{m-n}})(Q\cap (U_\phi\times\{0\}))=
\phi(U_\phi)\times\{0\}=V_\phi\times\{0\},\]
Since $f(x)=(x,0)$ for all $x\in U_\phi$, we have $\psi(x)=(\phi(x),0)$
and thus $\psi|_{U_\phi}=\phi$, if we identify
$\R^n$ with $\R^n\times\{0\}\sub\R^m$.

(c) By~(b), $U_\phi$ is a submanifold of~$U_\psi$.
The map $\pi\circ f\colon U_\psi\to U_\phi$ is a smooth retract.
\end{prf}
\emph{Proof of Theorem}~\ref{dl-mfd}.
We endow $M:=\bigcup_{n\in \N}M_n$
with the final topology with respect to the inclusion maps
$\lambda_n\colon M_n\to M$,
making it the direct limit $\dl\, M_n$\vspace{-.7mm} as a topological space (see
Definition~\ref{defnDLtop}).
Then $M$ is Hausdorff (see Lemma~\ref{basicDL}(h)).
Let $d_n:=\dim(M_n)$ for $n\in\N$.
Given $x_0\in M$, there exists $n_0\in \N$
such that $x_0\in M_{n_0}$.
Let $\phi_{n_0}\colon U_{n_0}\to V_{n_0}\sub\R^{d_{n_0}}$
be a chart for $M_{n_0}$
around $x_0$.
After shrinking $V_{n_0}$,
we may assume that $V_{n_0}$
is a ball $B_r(0)\sub\R^{d_{n_0}}$
for some $r>0$, whence $V_{n_0}$
(and hence also $U_{n_0}$) is smoothly contractible.
After replacing~$r$ with $r/2$,
we may assume that $U_{n_0}$
is relatively compact in $M_{n_0}$.
Recursively, using Lemma~\ref{extend-charts},
we find charts $\phi_n\colon U_n\to V_n\sub\R^{d_n}$
of $M_n$ for all $n\in \N$ with $n>n_0$
such that $U_n$ is relatively compact in~$M_n$ and smoothly
contractible for all $n\in\N$,
and $M_n\cap U_{n+1}=U_n$,
$\phi_n=\phi_{n+1}|_{U_n}$,
and $\R^{d_n}\cap V_{n+1}=V_n$.
Then
\[
U:=\bigcup_{n\geq n_0}U_n
\]
is open in ${\dl}_{n\geq n_0}M_n=\dl\,M_n$,\vspace{-.7mm}
and the topology induced by $\dl\, M_n$
on $U$ makes the latter the topological space
${\dl}_{n\geq n_0}U_n$\vspace{-.7mm}
(see (a) and (b) in Lemma~\ref{basicDL},
and Remark~\ref{firstremDL}(d)).
\phantom{\Big(} Likewise, $V:=\bigcup_{n\geq n_0}V_n$
is open in $\bigcup_{n\geq n_0}\R^{d_n}=\R^d$
and the topology induced by $\R^d$ on $V$
makes the latter the direct limit topological space
$\dl_{n\geq n_0}V_n$.
As $\dl_{n\geq n_0}\phi_n$
and its inverse map $\dl_{n\geq n_0}(\phi_n^{-1})$
are continuous (see Definition~\ref{DLmaps}),
the map $\phi:=\dl_{n\geq n_0}\phi_n\colon U\to V$
is a homeomorphism.
Let $\cA$ be the set of all homeomorphisms $\phi$
so obtained; let us verify that $\cA$ is a $C^\infty$-atlas
for~$M$.
To this end, let also $\psi=\dl_{n\geq m_0}\psi_n$
be in $\cA$, with charts $\psi_n\colon P_n\to Q_n\sub\R^{d_n}$.
such that $M_n\cap P_{n+1}=P_n$,
$\R^{d_n}\cap Q_{n+1}=Q_n$
and $\psi_{n+1}|_{P_n}=\psi_n$
for all integers $n\geq m_0$.
Let $P:=\bigcup_{n\geq m_0}P_n$.
Then $U\cap P$ is open in~$M$
and $M_n\cap U\cap P=U_n\cap P_n$
for all $n\geq \max\{n_0,m_0\}=:k_0$,
entailing that $\phi(U\cap P)=\bigcup_{n\geq k_0}\phi_n(U_n\cap P_n)$.
Since
\[
(\psi\circ\phi^{-1})|_{\phi_n(U_n\cap P_n)}=(\psi_n\circ \phi_n^{-1})|_{\phi_n(U_n\cap P_n)}
\]
is smooth for each $n\geq k_0$,
Proposition~\ref{diff-silva} shows that the bijection
\[ \psi\circ \phi^{-1}\colon \phi(U\cap P)\to \psi(U\cap P) \] is smooth.
Likewise, the inverse map
$\phi\circ \psi^{-1}\colon \psi(U\cap P)\to \phi(U\cap P)$ 
is smooth, whence $\psi\circ \phi^{-1}$
is a $C^\infty$-diffeomorphism.
Hence $\cA$ is a $C^\infty$-atlas for~$M$
modeled on $\R^d$
and we endow~$M$
with the corresponding maximal $C^\infty$-atlas.
By definition of the final topology,
$\lambda_{n_0}\colon M_{n_0}\to M$
is continuous for each $n_0\in\N$.
For each $x_0\in M_{n_0}$
and chart $\phi=\dl_{n\geq n_0}\phi_n$
as above, we have $\lambda_n(U_n)\sub U$
and $\phi\circ \lambda_n\circ \phi_{n_0}^{-1}=j|_{V_{n_0}}$
for the inclusion map $j\colon \R^{d_{n_0}}\to\R^d$,
whence $\lambda_n$ is $C^\infty$.
Now consider a map $f\colon M\to N$ to a $C^k$-manifold~$N$
with rough boundary, for $k\in \N_0\cup\{\infty\}$.
If $f$ is $C^k$, then $f|_{M_n}=f\circ\lambda_n$
is $C^k$ for each $n\in\N$.
Conversely, assume that $f|_{M_n}$
is $C^k$ for each $n\in \N$.
Then $f$ is continuous, by Lemma~\ref{ctsonsteps}(a).
Let $x_0\in M$ and $\psi\colon U_\psi\to V_\psi\sub F$
be a chart for~$N$ around $f(x_0)$.
Let $\phi=\dl_{n\geq n_0}\colon U\to V$
be a chart for $M$ around $x_0$
as above. Then $Q:=\phi(U\cap f^{-1}(U_\psi))$ is open
in $\R^d$
and
\[
\psi\circ f\circ \phi^{-1}|_{Q\circ V_n}=\psi\circ f|_{M_n}\circ \phi_n^{-1}|_{Q\circ V_n}
\]
is $C^k$ for each $n\geq n_0$, whence $\psi\circ f\circ \phi^{-1}|_Q\colon Q\to F$
is $C^k$, by Proposition~\ref{diff-silva}.
Hence $f$ is $C^k$. \qed
\begin{lem}\label{dl-of-ck-maps}
Let $M_1\sub M_2\sub \cdots$ and $N_1\sub N_2\sub\cdots$
be ascending sequences of finite-dimensional smooth manifolds,
such that all inclusion maps $M_n\to M_{n+1}$ and $N_n\to N_{n+1}$
are immersions. Endow $M:=\bigcup_{n\in \N}M_n$
and $N:=\bigcup_{n\in \N}N_n$
with the smooth manifold structures described in Theorem~{\rm\ref{dl-mfd}}.
Let $k\in \N_0\cup\{\infty\}$
and $(f_n)_{n\in\N}$
be a sequence of $C^k$-maps
$f_n\colon M_n\to N_n$
such that $f_{n+1}|_{M_n}=f_n$
for all $n\in\N$. Then $f\colon M\to N$,
\[
f(x):=f_n(x)\;\,\mbox{if $x\in M_n$}
\]
is a well-defined map and $C^k$.
\end{lem}
\begin{numba}
We write $\dl f_n:=f$
for the map in Lemma~\ref{dl-of-ck-maps}.
\end{numba}
\begin{prf}
For each $n\in\N$,
the inclusion maps $\lambda_n\colon M_n\to M$ and $j_n\colon N_n\to N$
are smooth,
by Theorem~\ref{dl-mfd}.
Thus $f\circ \lambda_n=j_n\circ f_n$ is $C^k$
for each $n\in \N$. Hence $f$ is $C^k$, by Theorem~\ref{dl-mfd}.
\end{prf}
\begin{lem}\label{prod-of-dl-mfd}
Let $M_1\sub M_2\sub \cdots$ and $N_1\sub N_2\sub\cdots$
be ascending sequences of finite-dimensional smooth manifolds,
such that all inclusion maps $M_n\to M_{n+1}$ and $N_n\to N_{n+1}$
are immersions. Endow $M:=\bigcup_{n\in \N}M_n$,
$N:=\bigcup_{n\in \N}N_n$,
and
\[
P=\bigcup_{n\in\N} (M_n\times N_n)
\]
with the smooth manifold structures described in Theorem~{\rm\ref{dl-mfd}}.
Then the map $\Psi\colon P\to M\times N$, $(x,y)\mto (x,y)$
to the product manifold is a $C^\infty$-diffeomorphism.
\end{lem}
\begin{prf}
By Lemma~\ref{DLcompaprod}, $\Psi$ is a homeomorphism. 
For $(x_0,y_0)\in P$, there is $n_0\in \N$
such that $x_0\in M_{n_0}$ and $y_0\in N_{n_0}$.
Let\vspace{-.3mm} $\phi=\dl_{\!n\geq n_0\,}\phi_n\colon U\to V\sub \R^d$
be a chart as in the proof of Theorem~\ref{dl-mfd},
with $\phi_n\colon U_n\to V_n\sub \R^{d_n}$.
Let $\psi=\dl_{\!\!n\geq n_0\,}\psi_n \: R \to Q\sub \R^{e}$\vspace{-.2mm}
be a corresponding chart of $N$, with ${\psi_n\colon R_n\to Q_n\sub\R^{e_n}}$.
For integers $n>n_0$, let $a_n:=d_n-d_{n-1}$
and $b_n:=e_n-e_{n-1}$.
Let $\alpha_n$ be the isomorphism of vector spaces from
\[
\R^{d_n}\times\R^{e_n}=\left(\R^{d_{n_0}}\oplus \bigoplus_{j=n_0+1}^n\R^{a_n}\right)
\times
\left(\R^{e_{n_0}}\oplus \bigoplus_{j=n_0+1}^n\R^{b_n}\right
)
\]
onto $(\R^{d_{n_0}}\times \R^{e_{n_0}})\oplus\bigoplus_{j=n_0+1}^n
(\R^{a_n}\times \R^{b_m})=\R^{d_n+e_n}$
which maps an element 
$(x_{n_0},\ldots, x_n,y_{n_0},\ldots, y_n)$
to $(x_{n_0}
,y_{n_0},\ldots, x_n,y_n)$.
Then $\theta_n:=\alpha_n\circ (\phi_n\circ \psi_n)$
are compatible charts for the ascending sequence of product manifolds
${M_n\times N_n}$, and we obtain a chart
$\theta:=\dl_{n\geq n_0}\theta_n\colon T\to S\sub \R^{d+e}$
of~$P$ from $T:=\bigcup_{n\geq n_0}(U_n\times R_n)$
onto $S:=\bigcup_{n\geq n_0}\theta_n(V_n\times Q_n)\sub\R^{d+e}$.
The map $\beta$ from
$
\R^{d+e}=(\R^{d_{n_0}}\times\R^{e_{n_0}})\oplus\bigoplus_{j>n_0}(\R^{a_j}\times \R^{b_j})$
to
\[
\left(\R^{d_{n_0}}\oplus \bigoplus_{j>n_0}\R^{a_j}\right)
\times \left(\R^{e_{n_0}}\oplus \bigoplus_{j>n_0}\R^{b_j}\right)
= \R^d\times \R^e
\]
taking $(x_j,y_j)_{j\geq n_0}$ to $((x_j)_{j\geq n_0},(y_j)_{j\geq n_0})$
is an isomorphism of topological vector spaces.
By construction, we have $\beta(S)=V\times Q$,
whence the map ${\beta|_S\colon S\to V\times Q}$ is a $C^\infty$-diffeomorphism.
Since
$(\phi\times \psi)\circ \Psi\circ \theta^{-1}=\beta|_S$,
we deduce that $\Psi|_R\colon R\to U\times P$
is a $C^\infty$-diffeomorphism. The assertion follows.
\end{prf}
\section{Notes and comments on Chapter~\ref{chap-manifold-constructions}}
Spaces of $C^k$-functions,
$C^k$-sections and manifolds of $C^k$-maps
are central objects in infinite-dimensional analysis and geometry.
For a compact smooth manifold~$M$,
a finite-dimensional smooth manifold~$N$
and $\ell\in \N_0$,
a Banach manifold structure on $C^\ell(M,N)$ was first constructed
by Eells \cite{Ee58, Ee66}.
Various classes of sections in fibre bundles and mappings between spaces of sections
were studied by Palais \cite{Pa68}.
Omori and collaborators considered
Lie groups of smooth mappings as intersections of descending
sequences of Banach manifolds (see \cite{OMY82, Omo97}
and further works).
For a paracompact smooth manifold $M$ with corners and
a finite-dimensional smooth manifold~$N$,
the smooth manifold structure on $C^\infty(M,N)$
was first constructed by Michor
(see \cite{Mr80}). For a discussion of manifolds of mappings in
the convenient setting of analysis, see~\cite{KM97}.
As in \cite{Ne06a}, use
a compact-open $C^k$-topology on $C^k(M,N)$
as the starting point, irrespective of a manifold structure on
$C^k(M,N)$. Proposition~\ref{explawfuncmfd}
varies the exponential laws for spaces of $C^{k,\ell}$-maps
in \cite{Alz13, AlS15, Alz19}.
Our construction of the smooth manifold structure on $C^\ell(M,N)$
for compact~$M$ varies the one in \cite{AGS20}.
For $\sigma$-compact $M$, we proceed as in
\cite{Gl21b} and \cite{GS22}, where also smooth manifolds
$C^{k,\ell}(M_1\times M_2,N)$
of $C^{k,\ell}$-maps are discussed and analogs
for functions of $n$ variables. On compact domains,
an exponential law of the form $C^{k,\ell}(M_1\times M_2,N)\cong
C^k(M_1,C^\ell(M_2,N))$ is then available.
Fine box products were developed in~\cite{Gl21b}
and used in \cite{Gl20a, GS22}.

Let us briefly recall the concepts of
large box products and small box products in general topology.
For a family $(X_j)_{j\in J}$
of topological spaces
one can endow $\prod_{j\in J}X_j$ with the box topology,
which has products $\prod_{j\in J}U_j$
of open subsets $U_j\sub X_j$
as a basis. If $(X_j,x_j)_{j\in J}$
is a family of pointed topological spaces,
one can consider the subset $\square_{j\in J}(X_j,x_j)$
of all families $(y_j)_{j\in J}$ such that $y_j=x_j$
for all but finitely many $j\in J$,
and endow it with the induced topology.
Compare~\cite{BR10}
for further information.
The fine box topology on a product $\prod_{j\in J} M_j$ of manifolds
is finer than the box topology and can be properly finer,
even for countable~$J$.
For example, the path component of~$0$ in
$\prod_{n\in \N}^{\fbx}\R$ equals $\bigoplus_{n\in \N}\R$,
while $\prod_{n\in \N}\R$ is path connected in the box topology.
These phenomena are analogous to the passage from the $\cD$-topology to the
fine $\cD$-topology
on $C^\infty(M,N)$ in \cite[\S4]{Mr80}.
Certain `` weak direct products''
of Lie groups $G_j$ were introduced in~\cite{Gl03b}.
We shall use them for countable $J$
and use the notation $\bigoplus_{j\in J}G_j$ (see Proposition~\ref{box-prod-regular});
as topological groups, they then coincide with $\square_{j\in J}(G_j,e)$.
Compare
\cite{HS17} and
\cite{Mr80} for
other descriptions of
the topology on $C^\ell(M,N)$ used in Theorem~\ref{map-mfd-noncomp},
notably for $\ell=\infty$
(cf.\ also~\cite{Il03} for a related coarser topology).

It is helpful to know simple additional hypotheses
which ensure that a nonlinear mapping
$f\colon C^\infty_c(\R)\to C^\infty_c(\R)$
is smooth if its restrictions to $C^\infty_K(\R)$
are smooth for each compact subset $K\sub\R$;
likewise for mappings between open subsets in spaces of compactly supported smooth
sections in vector bundles.
The concepts of local and almost local maps were developed in \cite{Gl05d, Gl03c, Gl04b};
see \cite{Gl21c} for a variant with parameters
(applied in \cite{PzS21}).
As we have seen in Theorem~\ref{almloc-thm}, the main point is the linear embedding of a space of
compactly supported sections onto a closed vector subspace of a suitable locally convex direct sum
(as in Lemma~\ref{embed-sections-sum}).
Such embeddings are also used in the alternative approach of ``patched mappings''
developed in \cite{Gl03c, Gl04b} and applied in \cite{Sme15}.
Without such tools,
mappings between spaces of sections were
discussed in~\cite{Gl13b}, under
less restrictive conditions concerning $f$ as in the above Proposition~\ref{omega-lemma}.
An interesting example of a smooth mapping $C^\infty_c(\R)\times C^\infty_c(\R)\to C^\infty_c(\R)$
which is neither local nor almost local on its domain $C^\infty_c(\R, \R\times \R)$
is the convolution of test functions, $(f,g)\mto f*g$,
which is continuous bilinear (see \cite{HST01, Gl12a}).
For analogous questions concerning convolution on finite-dimensional Lie groups, see~\cite{BG14}.
We mention that for a local map $f\colon \Gamma_{C^k}(E_1)\to\Gamma_{C^k}(E_2)$
which is $C^1$, for $\sigma$ in the domain and $x$ in the base manifold~$M$, the element
$f(\sigma)(x)\in E_2$
only depends on the G\^{a}teaux derivatives
$\delta^j_{\phi(x)}(\theta_2\circ \sigma\circ\phi^{-1})$ for $j\in \N_0$ with $j\leq k$,
where $\theta\colon E_1|_U\to U\times F_1$
is a local trivialization of~$E_1$ with $x\in U$
and $\phi\colon U_\phi\to V_\phi$ a chart of~$M$
around~$x$. Using the Fundamental Theorem of Calculus,
this follows from the fact that
$df(\tau,\cdot)$ is a linear local map, to which an appropriate version of
Peetre's Theorem as in~\cite{Da15} applies (see \cite{Ptr59}
for the original result).

In rare cases, smooth manifold structures compatible with the compact-open $C^\ell$-topology
are available on $C^\ell(M,N)$ also for non-compact~$M$, e.g.\ on $C^\infty(\R,K)$
if $K$ is a regular Lie group~\cite{NeW08b}
(see \cite{Alz13, Alz21, GS22} for generalizations).

Topological spaces which are ascending unions of finite-dimensional topological manifolds $M_1\sub M_2\sub\cdots$ (with the direct limit topology) were first
discussed by \cite{Han71}. Ascending unions of finite-dimensional Lie groups $G_1\sub G_2\sub\cdots$
were considered in \cite{NRW91} and \cite{KM97}
in important cases (see Chapter~\ref{ch:dirlim} for more details).
Ascending unions of finite-dimensional manifolds were turned into
manifolds in \cite{Gl03a, Gl05c}.

\part{General Lie theory} 

\chapter{Locally convex Lie groups} \mlabel{ch:3}
In this monograph we follow the philosophy that the definition of an 
infinite-dimensional Lie group should be as simple as possible. Therefore 
it is natural to define a Lie group as a locally convex manifold $G$, 
endowed with a group structure,  
for which multiplication and inversion are smooth maps. This is the familiar 
definition from finite-dimensional Lie theory, the only difference is that 
we admit a larger class of model spaces. In finite dimensions, and even 
in the Banach context, it would lead to the same objects if we require the 
manifold structure and the group operations to be analytic. This is no 
longer true in general and it actually turns out that many interesting 
infinite-dimensional groups carry no analytic Lie group structure. 

For any Lie group 
$G$, the tangent space $T_\be(G)$ in the identity element carries a natural 
Lie algebra structure, obtained from the Lie bracket of left invariant 
vector fields on~$G$. We thus obtain the Lie functor 
$\L$ from the category of Lie groups to the category of 
locally convex topological Lie algebras. This process is explained in 
Section~\ref{sec:3.1}. 

In Section~\ref{sec:3.2}, we then turn to the adjoint action $\Ad$ of the Lie group $G$ 
on its Lie algebra, which is best understood as obtained from the restriction of 
the conjugation 
action of the tangent group $T(G)$ on itself. Here we use that the tangent 
bundle $T(G)$ carries a natural Lie group structure defined by the tangent map 
of the multiplication on $G$. We also discuss semidirect products and show that 
$T(G) \cong \L(G) \rtimes_{\Ad} G$.  

Section~\ref{sec:3.3} provides a collection of tools to obtain Lie group, 
resp., manifold structures on groups. Here a natural strategy is 
to construct a chart around the identity in which the group 
operations are smooth. This suffices in many situations 
to specify a global Lie group structure. 

In Section~\ref{sec:3.4}, we introduce the logarithmic derivative of a 
smooth function $f \: M \to G$ from a manifold $M$ with values in a Lie group 
$G$. If $M$ is connected, then such a function is uniquely determined 
by its logarithmic derivative and the value in one point of $M$. 
This has several immediate implications for homomorphisms of 
Lie groups, in particular, that a morphism $\phi$ of connected Lie groups 
is determined by $\L(\phi)$. Compared to the classical finite-dimensional 
theory, we put a stronger emphasize on logarithmic derivatives because 
an exponential function is not available in general. 

We devote Section~\ref{sec:3.5} to the exponential function 
$\exp_G \: \L(G) \to G$ in the general context of locally convex Lie groups. 
We discuss a variety of exponential functions of different 
kinds of groups. As these examples show, the exponential function need not be 
a local diffeomorphism, so that there are many Lie groups without 
natural charts. This is a serious difficulty one has to face 
in infinite-dimensional Lie theory. In Chapter~\ref{ch:5}, 
we shall investigate the class of 
locally exponential Lie groups, for which this pathology does not occur. 


{\bf Prerequisites:} The basics from infinite-dimensional manifolds, 
such as the 
tangent bundle and 
tangent maps. In Section~\ref{sec:3.3} we further need some facts from 
the theory of covering spaces, such as the existence of a universal covering 
space and the existence of liftings. We shall also need the basic facts on 
integral curves of time-dependent vector fields on manifolds.

\section{Lie groups and their Lie algebras} \mlabel{sec:3.1} 

\nin In this section, we define locally convex Lie groups.  
We explain how the Lie algebra and the corresponding Lie functor are defined 
and describe some basic properties, such as the relation of the Lie bracket 
with the Taylor expansion of the group multiplication in local charts. 

\subsection{The definition of a Lie group} 

There are two types of additional structures 
on groups. The first level consists of a topological structure 
compatible with the group structure, which leads to the concept 
of a topological group,  and the second level is a 
differentiable structure, which leads to the concept 
of a Lie group. 

\begin{defn} \mlabel{def:3.1.1} 
(a) A {\it topological group} \index{topological group} 
is a Hausdorff space $G$ 
 endowed with a group structure such that the multiplication 
map and the inversion map are continuous. 

(b) A {\it Lie group} \index{Lie group} 
is a 
locally convex manifold $G$, endowed with a group structure such that the multiplication map and the inversion map are smooth. 
Since smooth maps are continuous, every Lie group is 
in particular a topological group.

If, in addition, $G$ is a $\K$-analytic manifold for $\K = \R$ or $\C$,
  and the group operations are $\K$-analytic, we call $G$ a
  {\it $\K$-analytic Lie group}.  \index{Lie group!$\K$-analytic}

(c) If $G$ and $H$ are topological groups, then a group homomorphism 
$\phi \: G \to H$ is called a {\it morphism of topological groups} 
if $\phi$ is continuous. If, in addition, $G$ and $H$ are Lie groups,
then a group homomorphism 
$\phi \: G \to H$ is called a {\it morphism of Lie groups} 
if $\phi$ is smooth. 
\index{morphism of topological/Lie groups} 
\end{defn} 

Throughout this book, we use the following notation: 
We write $\be$ for the identity element of $G$. 
For $g \in G$ we write $\lambda_g \: G \to G, x \mapsto gx$ for the
{\it left multiplication} 
\index{left multiplication maps $\lambda_g$} 
\index{right multiplication maps $\rho_g$} 
\index{conjugation maps!$c_g(x) = gxg^{-1}$} 
\index{inversion!$\eta_G(g) = g^{-1}$} 
by $g$ and 
$\rho_g \: G \to G, x \mapsto xg$ for the {\it right multiplication} by 
$g$. Both are diffeomorphisms of $G$ (Exercise~\ref{exer:2.2.1}). We write 
$$m_G \: G \times G \to G,\quad (x,y) \mapsto xy$$ 
for the {\it multiplication map} and 
$$\eta_G \: G \to G, \quad x \mapsto x^{-1}$$ for the {\it inversion}. 

\begin{ex} (Vector groups) Each locally convex space 
$E$ is an abelian Lie group with respect to addition and 
the obvious manifold structure. 
\end{ex}

Vector groups $(E,+)$ form the most elementary Lie groups. 
The next natural class are unit groups of algebras. This leads us to the 
concept of a continuous inverse algebra, which came up in the 1950s 
(cf.\ \cite{Wa54a}, \cite{Wa54b}): 

\begin{defn} (a) A {\it locally convex algebra} \index{locally convex algebra} 
is a locally convex space 
$\cA$, endowed with a continuous bilinear multiplication 
$\cA \times \cA \to \cA, (a,b) \mapsto ab$ 
which is associative. 
A unital locally convex algebra $\cA$ is called a 
{\it continuous inverse algebra} \index{continuous inverse algebra, cia} 
(cia for short) 
if its unit group $\cA^\times$ 
is open and the inversion is a continuous map $\cA^\times \to \cA, 
a \mapsto a^{-1}$. 

(b) If $\cA$ is a locally convex algebra which is 
not unital, then we obtain a monoid structure on $\cA$ 
by $x * y := x + y + xy$, for which $0$ is the identity element. 
In this case we write $\cA^\times$ for the unit group of $(\cA,*)$ 
and say that 
$\cA$ is a {\it non-unital continuous inverse algebra} if $\cA^\times$ is open 
and the inversion map 
$\eta_\cA \: \cA^\times \to \cA$ is continuous. 

If $\cA_+ := \cA \times \K$ is the unital locally convex algebra with the 
multiplication $(x,t)(x',t') := (xx'+tx'+t'x,tt')$, then the map 
$(\cA,*) \to \cA \times \{1\}, a \mapsto (a,1)$ is an isomorphism of monoids, 
and it is easy to see that $\cA_+$ is a continuous inverse algebra if and only 
if $\cA$ is a (not necessarily unital) continuous inverse algebra (Exercise~\ref{exer:3.1.4}). 
\end{defn}

\begin{ex} \mlabel{ex:cia-grp} (Unit groups as Lie groups) Let 
$\cA$ be a continuous inverse algebra over $\K$ and $\cA^\times$ 
be its unit group. 
As an open subset of $\cA$, the group $\cA^\times$ carries a natural manifold 
structure. The multiplication on $\cA$ is bilinear and continuous, 
hence a smooth map. 
Therefore the multiplication of $\cA^\times$ is smooth and it remains to 
see that the inversion $\eta \: \cA^\times \to \cA^\times$ 
is smooth. Since $\cA$ is a continuous inverse algebra, $\eta$ is continuous. 

For $a,b \in \cA^\times$, we have 
$b^{-1} - a^{-1} = a^{-1}(a - b)b^{-1},$
so that for $t \in \R$ we get 
\[ \eta(a + th) - \eta(a) = (a + th)^{-1} - a^{-1} = a^{-1}(-th)(a+th)^{-1} = -ta^{-1}h(a+th)^{-1}.\]
Therefore $\eta$ is everywhere differentiable with 
\begin{eqnarray}
  \label{eq:diff-inv}d\eta(a)(h) = -a^{-1}ha^{-1}. 
\end{eqnarray}
Now the continuity of $\eta$ implies that 
$\dd\eta \: \cA^\times \times \cA \to \cA$ is 
continuous, hence that $\eta$ is a $C^1$-map. 
From (\ref{eq:diff-inv}), we further derive that if $\eta$ is $C^k$, then 
$d\eta$ is also $C^k$, so that $\eta$ is $C^{k+1}$. Inductively, it 
follows that $\eta$ is smooth.

In some cases it is also possible to obtain a Lie group structure on the unit group 
$\cA^\times$ of a unital locally convex algebra if $\cA^\times$ is not open 
(cf.\ Example~\ref{ex:exotic-units}). 
\end{ex}

\begin{rem}
  \label{rem:hausdorff} 
Since it also makes sense to consider manifolds without 
assuming that they are Hausdorff (cf.\ \cite{Pa57}, 
\cite{La99}), it is worthwhile to observe that 
this does not lead to a larger class of Lie groups. 

In fact, let $G$ be a ``Lie group'', which is not necessarily Hausdorff. Then $G$  
is in particular a topological group which possesses an identity 
neighborhood $U$ homeomorphic to an open subset of a locally convex space. 
As $U$ is Hausdorff, and since the subgroup $\oline{\{\be\}}$ of $G$ coincides 
with the intersection of all $\be$-neighborhoods 
(Exercise~\ref{exer:3.1.5}), the closedness of 
$\{\be\}$ in $U$ implies that $\oline{\{\be\}} \cap U = \{\be\}$ 
is closed and hence that $G$ is a Hausdorff topological 
group (Exercise~\ref{exer:3.1.6}). 
\end{rem}

\subsection{The Lie algebra of a Lie group} 
We now turn to the definition of the Lie algebra of a Lie group. 

\begin{defn} {\rm(The Lie algebra of a Lie group)} 
\mlabel{def:liealg}
A vector field $X$ on the Lie group $G$ is called 
 {\it left invariant} if 
\index{vector field!left invariant} 
 \begin{eqnarray}
   \label{eq:left-invar}
X \circ \lambda_g = T(\lambda_g) \circ X \: G \to T(G) 
 \end{eqnarray}
holds for each $g \in G$, i.e., $X$ is $\lambda_g$-related to itself for each $g \in G$. 
We write ${\cal V}(G)^l$ for the set of left invariant vector fields in 
${\cal V}(G)$. 
The left invariance of a vector field $X$ implies in particular that, for each 
$g \in G$, we have $X(g) = T_\be(\lambda_g)X(\be)$. 

For each $x \in \g$, we define a vector field $x_l \in {\cal V}(G)$ 
by $x_l(g) := T_\be(\lambda_g)x$. For $h \in G$ we then have 
$$ x_l(gh) 
= T_\be(\lambda_{gh})x
= T_h(\lambda_{g})T_\be(\lambda_{h})x 
= T_h(\lambda_{g}) x_l(h), $$
i.e., $x_l \in {\cal V}(G)^l$ is the unique left invariant vector field with 
$x_l(\be) = x$. This means that the map 
$$ \ev_\be \: {\cal V}(G)^l \to T_\be(G), \quad X \mapsto X(\be) $$ 
is a linear bijection. If $X,Y$ are left invariant, then (\ref{eq:left-invar}) means that 
they are 
$\lambda_g$-related to themselves, and 
their Lie bracket $[X,Y]$ inherits this property. 
We conclude that the set 
${\cal V}(G)^l$ of left invariant vector fields on $G$ is a Lie subalgebra of the 
Lie algebra ${\cal V}(G)$ of all smooth vector fields on $G$. Since the map 
$$ T_\be(G) \to {\cal V}(G)^l, \quad x \mapsto x_l $$
is a linear isomorphism with inverse $\ev_\be \: X \mapsto X(\be)$, we obtain a Lie bracket on 
$T_\be(G)$ by 
$$ [x,y] := [x_l, y_l](\be). $$
By definition, this Lie bracket satisfies 
\begin{eqnarray}
  \label{eq:left-brack}
[x_l, y_l] = [x,y]_l.
\end{eqnarray}

To show that the Lie bracket on $T_\be(G)$ is continuous, 
let $E := T_\be(G)$ and choose a 
local $E$-chart $(\phi,U)$ of $G$ with $\phi(\be) = 0$ and $T_\be(\phi) = \id_E$. 
For $x \in T_\be(G)$ we then obtain a smooth vector field 
\begin{eqnarray}
  \label{eq:left-vec} 
\tilde x_l := T(\phi) \circ x_l \circ \phi^{-1}
\end{eqnarray}
on $V := \phi(U)$. We identify 
$T(V) \cong V \times E$ and $T^2(V) \cong T(V) \times T(E) \cong V \times E^3$. 
Then the Related Vector Field Lemma~\ref{larelglob} implies 
\begin{eqnarray}
  \label{eq:loc-brack}
[x,y] 
&=& [x_l, y_l](\be) = [\tilde x_l, \tilde y_l](0) 
= d\tilde y_l(0)\tilde x_l(0) - d\tilde x_l(0)\tilde y_l(0) \notag \\
&=& d\tilde y_l(0)x - d\tilde x_l(0)y. 
\end{eqnarray}
Clearly, the function 
$$ \psi \: E \times V \to E, \quad 
(x,z) \mapsto \tilde x_l(z) 
= T(\phi) x_l(\phi^{-1}(z)) 
= T(\phi) T(m_G)(0_{\phi^{-1}(z)}, x) $$
is smooth, so that 
$(x,y) \mapsto d\tilde x_l(0)y$
is continuous bilinear, and hence the bracket on $E \cong T_\be(G)$ is continuous. 

The Lie algebra 
$$ \L(G) := (T_\be(G),[\cdot,\cdot]) $$
is called the {\it Lie algebra of the Lie group $G$}. 
\index{Lie algebra of $G$, $\L(G)$,$\g$} 
To simplify the notation, we shall 
also use $\g$ instead of $\L(G)$. 

We call a locally convex Lie algebra $\g$ {\it integrable} 
if \index{Lie algebra!integrable} 
$\g \cong \L(G)$ holds for some  Lie group $G$. 
\end{defn}

\begin{lem} \mlabel{lem:second-order} Let $E,F$ and $X$ be locally convex spaces, 
$U \subeq E \times F$ an open $0$-neighborhood and 
$f \: U \to X$ a $C^3$-function, such that 
$$ f(x,0) = \alpha(x) \quad \mbox{ and } \quad  
f(0,y) = \beta(y) $$
holds for $(x,0),(0,y) \in U$ and linear functions 
$\alpha \: E \to X$, resp., $\beta \: F \to Y$. Then 
the second order Taylor polynomial of $f$ in $(0,0)$ is of the form 
$$ (x,y) \mapsto \alpha(x)  + \beta(y) + \gamma(x,y), $$
where $\gamma \: E \times F \to X$ is a continuous bilinear map, and for 
$(x,y) \in E \times F$ we have 
$$ \gamma(x,y) = 
\hbox{$\frac{\partial^2}{\partial s \partial t}$\vrule}_{s,t=0}  \ 
f(sx,ty). $$ 
\end{lem}

\begin{prf} For the first order term of the Taylor polynomial 
of $f$ in $(0,0)$, our assumptions lead to 
$$ d f(0,0)(x,y) =  d f(0,0)(x,0) + d f(0,0)(0,y) = \alpha(x) + \beta(y). $$

The second order term $\gamma \: E \times F \to X$ is quadratic, 
and our assumptions imply that 
$$ \gamma(x,0) = \gamma(0,y) = 0 \quad\mbox{ for } \quad 
x \in E, y \in F. $$
Let $B \: (E \times  F) \times (E \times F) \to X$ be the unique 
continuous bilinear map with 
$$ \gamma(x,y) = B((x,y),(x,y)). $$
Then $0 = B((x,0),(x,0)) = B((0,y),(0,y))$ shows that 
$$ \gamma(x,y) = B((x,0), (0,y)) + B((0,y),(x,0)) = 2 B((x,0), (0,y)), $$
is bilinear. 
We conclude in particular that the second order Taylor polynomial of 
$f(sv,tw)$ is $s\alpha(v) + t \beta(w) + ts \gamma(v,w)$, which implies the assertion. 
\end{prf}

\begin{rem} {\rm(Lie bracket and Taylor expansions)} \mlabel{rem:brack-taylor} 

Let $\g := \L(G)$ and pick a $\g$-chart $(\phi,U)$ of $G$ with 
$\phi(\be) = 0$ and $T_\be(\phi) = \id_\g$. Write $V := \phi(U)$. 
Then we obtain a smooth map
\begin{eqnarray*}
D&:=& \{ (x,y) \in V \times V \: \phi^{-1}(x)\phi^{-1}(y) \in U\} \to \g, \\ 
 && (x,y) \mapsto x * y := m(x,y) := \phi\big(\phi^{-1}(x)\phi^{-1}(y)\big). 
\end{eqnarray*}

(a) We consider the Taylor expansion of $m$ in $(0,0) 
\in \g \times \g$. In view of $0 * 0 =0$, the constant term vanishes, so that 
\[  x * y = b_1(x,y) + b_2(x,y) + b_3(x,y) + \ldots 
\quad \mbox{ where  } \quad b_k = \frac{1}{k!} \delta^k_{(0,0)} m \]
are continuous homogeneous polynomials of degree $k$. 

In view of $x = x * 0 =0 * x$, Lemma~\ref{lem:second-order} implies that 
$$ b_1(x,y) = x+y $$
and that $b_2$ is bilinear with 
$$ b_2(x,y) = \hbox{$\frac{\partial^2}{\partial s \partial t}$\vrule}_{s,t=0}  \ 
sx * ty. $$
To evaluate the right hand side, we first note that 
\begin{eqnarray*}
\hbox{$\frac{\partial}{\partial t}$\vrule}_{t=0} \ x * ty 
&=&  \hbox{$\frac{\partial}{\partial t}$\vrule}_{t=0} \  \phi(\phi^{-1}(x)\phi^{-1}(ty))\\ 
&=& T(\phi) T_\be(\lambda_{\phi^{-1}(x)})y 
= T(\phi) y_l(\phi^{-1}(x)) = \tilde y_l(x), 
\end{eqnarray*}
(cf.\ \eqref{eq:left-vec})
and hence 
$$ b_2(x,y) = \hbox{$\frac{\partial}{\partial t}$\vrule}_{t=0} \ \tilde y_l(tx) 
= d\tilde y_l(0)x. $$ 
We thus arrive with (\ref{eq:loc-brack}) at 
$$ [x,y] = d \tilde y_l(0) x - d\tilde x_l(0) y = b_2(x,y) - b_2(y,x), $$
i.e., the Lie bracket of $\g$ is the skew-symmetric part of the bilinear 
map~$b_2$.

(b) Let $\eta(x) = x^{-1} = s_1(x) + s_2(x) +\cdots$ denote the Taylor expansion of 
the inversion map in $0$. From $x * x^{-1} = 0$ we obtain 
$$ 0 = x * (s_1(x) + s_2(x) + \cdots) = \big(x + s_1(x)\big) + 
\big(s_2(x) + b_2(x, s_1(x))\big) + \cdots. $$
This leads to 
\begin{eqnarray}
  \label{eq:4.1.4}
s_1(x) = -x \quad \hbox{ and } \quad s_2(x) = -b_2(x,s_1(x)) =
b_2(x,x) 
\end{eqnarray}
because $b_2$ is bilinear. 

(c) If $x$ and $y$ are sufficiently close to $0$, we obtain with (b) 
the Taylor expansion of the conjugation map from the Chain Rule for 
Taylor polynomials: 
\begin{eqnarray} \label{eq:conjug}
(x * y) * x^{-1} 
&&= (x + y + b_2(x,y) + \cdots) * (-x + b_2(x,x) + \cdots) \cr
&&= (x + y -x) + b_2(x,y) + b_2(x,x) + b_2(x + y,-x) + \cdots \cr
&&= y + b_2(x,y) - b_2(y,x) + \cdots = y + [x,y] + \cdots 
\end{eqnarray}
For the commutator $(x,y) := x * y * x^{-1} * y^{-1}$ we therefore get 
\begin{eqnarray} \label{eq:commut}
&& x * y * x^{-1} * y^{-1} 
= (y + [x,y] + \cdots) * (-y + b_2(y,y) + \cdots) \cr
&&= [x,y] + b_2(y,y) + b_2(y,-y)  + \cdots =  [x,y]  + \cdots   
\end{eqnarray}
\end{rem}

We now take a look at the Lie algebras of the Lie groups from the examples discussed above. 
\begin{ex} \mlabel{ex:abgrp} (a) If $G$ is an abelian Lie group, then the bilinear 
map \break $b_2\: \g \times \g \to \g$ in Remark~\ref{rem:brack-taylor} 
is symmetric, which implies that $\L(G) = \g$ is abelian. 
This applies in particular to the additive Lie group $(E,+)$ of a locally 
convex space~$E$. 

(b) Let $\cA$ be a continuous inverse algebra. Then 
\[ \phi \: \cA^\times \to \cA, x \mapsto x - \1 \] 
is a global chart of $\cA^\times$, satisfying $\phi(\1) = 0$. 
In this chart, the group multiplication is given by 
\[  x * y := \phi(\phi^{-1}(x)\phi^{-1}(y)) = (x +\1)(y+\1) - \1 
= x + y + xy. \]
In the terminology of Remark~\ref{rem:brack-taylor}, we then have 
$b_2(x,y) = xy$ and therefore 
$\L(\cA^\times) = (\cA,[\cdot,\cdot])$, where 
$[x,y] =xy - yx$ is the commutator bracket on the associative algebra~$\cA$.

As we shall see in Example~\ref{ex:expo-cia}, any continuous inverse algebra has an exponential 
function 
$$ \exp \: \cA \to \cA^\times, \quad x \mapsto 
e^x := \sum_{n = 0}^\infty \frac{1}{n!} x^n, $$
defined by a converging power series defining a smooth map. 
Moreover, the power series 
$$ \log(g)  = \sum_{n=1}^\infty \frac{(-1)^{n+1}}{n} (g-\1)^n $$
defines a smooth function in a neighborhood of $\1$, 
and $\log \circ \exp\res_U = \id_U$ holds on some open 
$0$-neighborhood, so that $\log$ defines a chart of the unit group 
$\cA^\times$. One can show that the corresponding product 
$$ x * y := \log(\exp x \exp y) $$
is given by the {\it Hausdorff series} \index{Hausdorff series} 
\begin{eqnarray*}
 x * y 
&&=x + \sum_{k,m\ge 0\atop p_i+q_i>0}(-1)^k \frac{
(\ad x)^{p_1}(\ad y)^{q_1}\ldots
(\ad x)^{p_k}(\ad y)^{q_k}(\ad x)^m}{ 
(k+1)(q_1+\ldots+q_k+1)p_1!q_1!\ldots p_k!q_k!m!}y \\
&&=x+y+\frac{1}{2}[x,y]+\frac{1}{12}[x,[x,y]]+\frac{1}{12}[y,[y,x]]+\ldots 
\end{eqnarray*}
In particular, $b_2(x,y) = \frac{1}{2}[x,y]$ is skew-symmetric in 
this case. 
\end{ex}

We now turn to homomorphisms of local Lie groups, which leads to 
the interpretation of the assignment of the Lie algebra as a functor. 

\begin{prop} If $\phi \: G \to H$ is a morphism 
of Lie groups, then 
$$ \L(\phi) := T_{\be}(\phi) \: \L(G) \to \L(H) $$
is a continuous homomorphism of Lie algebras. 
\end{prop}

\begin{prf} Let $x \in \g$ and $x_l(g) = T_\be(\lambda_g)(x)$ be 
the corresponding left invariant 
vector field on $G$. For each $g \in G$, we then have 
$\phi \circ \lambda_g = \lambda_{\phi(g)} \circ \phi.$
This implies that 
$$ T_g(\phi)x_l(g) = T_g(\phi) T_\be(\lambda_g)x  
= T_\be(\lambda_{\phi(g)}) T_\be(\phi)(x)
= \big(\L(\phi)(x)\big)_l(\phi(g)), $$
so that the vector fields $x_l$ on $G$ 
and $\big(\L(\phi)x\big)_l$ on $H$ are $\phi$-related. 
We conclude with Lemma~\ref{larelglob} %
that for $x,y \in \L(G)$ the vector fields 
$[x_l,y_l]$ and 
$[(\L(\phi)x)_l,(\L(\phi)y)_l]$ 
are also $\phi$-related. Evaluating in $\be$, we obtain 
$$ [\L(\phi)x,\L(\phi)y] = [\L(\phi)x,\L(\phi)y]_l(\be) 
= T_\be(\phi)([x_l,y_l](\be)) = \L(\phi)[x,y]. $$
Hence $\L(\phi)$ is a homomorphism of Lie algebras. Its continuity follows from 
the smoothness of $\phi$. 
\end{prf}

\begin{cor} \mlabel{cor:liefunct} If 
we assign to a Lie group $G$ its Lie algebra $\L(G)$ and to a 
morphism of Lie groups $\phi \: G \to H$ the linear map 
$\L(\phi) := T_\be(\phi)$, then $\L$ defines a functor 
from the category of Lie groups to 
the category of topological Lie algebras. 
\end{cor}

\begin{prf} We have already seen that $\L(\phi)$ is a morphism of topological 
Lie algebras. As a direct consequence of the Chain Rule, we get 
$\L(\id_G)  = \id_{\L(G)}$ and $\L(\phi_1 \circ \phi_2) = \L(\phi_1) \circ \L(\phi_2)$, 
and the assertion follows. 
\end{prf}

\begin{rem} The functoriality of $\L$ implies in particular that for 
each isomorphism of Lie groups $\phi \: G \to H$, the map 
$\L(\phi) \: \L(G) \to \L(H)$ is an isomorphism of Lie algebras. 
In particular, we obtain for each Lie group $G$ a group homomorphism 
$$ \L \: \Aut(G) \to \Aut(\L(G)), \quad \phi \mapsto \L(\phi). $$
Composing with the conjugation action 
$G \to \Aut(G),g \mapsto c_g$, of $G$ on itself, we thus obtain a 
homomorphism  $\Ad \: G \to \Aut(\L(G))$ defining the 
{\it adjoint representation of $G$}, \index{adjoint representation} 
which is studied in more detail in the following section.   
\end{rem}

  \begin{small}
    \subsection*{Exercises for Section~\ref{sec:3.1}} 

    \begin{exer} \mlabel{exer:3.1.1} Let $G$ be a Lie group. 
      A vector field $X$ on the Lie group $G$ is called 
      {\it right invariant} if  \index{vector field!right invariant}
      $X \circ \rho_g = T(\rho_g) \circ X \: G \to T(G) $ 
      holds for each $g \in G$.  We write ${\cal V}(G)^r$ for the set of right 
      invariant vector fields on $G$. Show that: 
      \begin{enumerate}
      \item[\rm(1)] The space ${\cal V}(G)^r$ of right invariant vector fields on $G$ 
        is a Lie subalgebra of ${\cal V}(G)$. 
      \item[\rm(2)] $(\eta_G)_* \: {\cal V}(G)^l \to {\cal V}(G)^r, X \mapsto 
        T(\eta_G) \circ X \circ \eta_G^{-1}$ is an isomorphism 
        of Lie algebras. 
      \item[\rm(3)] The evaluation map $\ev_\be \: {\cal V}(G)^r \to T_\be(G)$ 
        is a linear isomorphism which induces an anti-isomorphism of Lie algebras, i.e., 
        $$ \ev_\be([X,Y]) = - [\ev_\be(X),\ev_\be(Y)]. $$
        Hint: Use (2) and $T_\be(\eta_G)x = -x$. 
      \item[\rm(4)] Let $x_r$ denote the unique right invariant vector field with 
        $x_r(\be) = x$. Then $[x_r, y_r] = - [x,y]_r$ for $x,y \in T_\be(G)$. 
      \end{enumerate}
    \end{exer}

    \begin{exer} \mlabel{exer:3.1.2} Show that each unital Banach algebra $\cA$ 
      is a continuous inverse algebra. 
      Hint: For $\|x\| < 1$ the Neumann series $\sum_{n = 0}^\infty x^n$ provides an inverse 
      of the element $\be - x$.   
    \end{exer}

    \begin{exer}
      Let $G$ be a Lie group. Show that 
      all connected components of $G$ are diffeomorphic. 
      Hint: The maps $\lambda_g \: G \to G, x \mapsto gx$ are 
      diffeomorphisms. 
    \end{exer}
    
    \begin{exer}\mlabel{exer:3.1.4b} 
      Let $G$ and $H$ be topological (Lie) groups 
      and $\phi \: G \to H$ be a group homomorphism. 
      Show that $\phi$ is continuous (smooth) if there exists an open 
      identity neighborhood $U \subeq G$ on which $\phi$ is 
      continuous (smooth). 
    \end{exer}
    
    \begin{exer} \mlabel{exer:3.1.4} 
      Let $\cA$ be a locally convex $\K$-algebra and define a product on 
      $\cA_+ := \cA \times \K$ by $(x,t)(x',t') := (xx'+tx'+t'x,tt')$. 
      Show that for the multiplication 
      $$ a * b := a + b + ab $$
      on $\cA$, the map 
      $(\cA,*) \to \cA \times \{1\}, a \mapsto (a,1)$, is an isomorphism of monoids 
      and that $\cA_+$ is a continuous inverse algebra if and only 
      if $\cA$ is a continuous inverse algebra. 
    \end{exer}

    \begin{exer} \mlabel{exer:3.5.2}
      Let $\phi \: G \to H$ be a surjective morphism of Banach--Lie groups. 
      Show that $\phi$ is a covering morphism if and only if 
      $\L(\phi)$ is an isomorphism of Banach--Lie groups. 
    \end{exer}

    \begin{exer}
      \mlabel{exer:3.1.7} 
      (a) Let $m \: G \times G \to G$ be a smooth
      associative multiplication on the manifold $G$ with identity element
      $\be$. Show that the differential in $(\be,\be)$ is given by 
      $$ T_{(\be,\be)}(m) \: T_\be(G) \times T_\be(G) \to T_\be(G), \quad 
      (v,w) \mapsto v + w. $$

      (b) Show that the smoothness of
      the inversion in the definition of a Banach--Lie group 
      is redundant because the Inverse Function Theorem 
      can be applied to the map 
      $$G \times G \to G \times G, \quad 
      (x,y) \mapsto (x, xy) $$
      whose differential in $(\be,\be)$ is given by the map $(v,w) 
      \mapsto (v, v + w)$.
    \end{exer}

\begin{exer}\mlabel{exer:III.7} Let $G$ be a Lie group with Lie algebra $\g$ and 
$\phi \: U_G \to \g$ a local chart with $\phi(\be) =0$. Show that:
\begin{description}[(D)]
\item[(1)] For the local multiplication $x * y := \phi(\phi^{-1}(x)\phi^{-1}(y))$,  
the second order Taylor polynomial of $x*y*x^{-1}*y^{-1}$ in 
$(0,0)$ is the Lie bracket $[x,y]$. 
\item[(2)] Use (1) to show that for each morphism of Lie groups 
$\phi \: G \to H$, the map $d\phi(\be)$ is a homomorphism of Lie algebras. 
Hint: Compare the second order Taylor polynomials of 
$\phi(x)*\phi(y)*\phi(x)^{-1}*\phi(y)^{-1}$ 
and $\phi(x*y*x^{-1}*y^{-1})$ by using the Chain Rule for 
Taylor polynomials. 
\end{description}
\end{exer}

  \end{small} 

  \section{Adjoint action and semidirect products} \mlabel{sec:3.2}

  The adjoint representation $\Ad \: G \to \Aut(\L(G)), g \mapsto \L(c_g)$ 
  is an important structural feature of a Lie group. 
  It is trivial for abelian Lie groups, and in general it can be considered 
  as providing a linear picture of the non-commutativity of $G$. To understand 
  the adjoint action properly, it is convenient to introduce the natural group 
  structure on the tangent bundle $T(G)$, obtained from the tangent map of the 
  multiplication of $G$. The tangent bundle is a 
  semidirect product of the additive group $(\L(G),+)$ with $G$ acting by the 
  adjoint representation. In this section, we present the basic facts concerning 
  the adjoint representation, the tangent bundle and semidirect products of 
  Lie groups. 

  \subsection{The adjoint representation} 

  Before we discuss the adjoint representation in particular, it 
  is useful to start with some basic observations concerning actions 
  of Lie groups in general. The following proposition describes how to 
  associate to a smooth action of a Lie group $G$ on a manifold $M$ 
  a Lie algebra homomorphism $\L(G) \to {\cal V}(M)$, the corresponding 
  {\it infinitesimal action}. \index{infinitesimal action} 
We further associate to each 
  smooth representation $\pi \: G \to \GL(E)$ on a locally convex space 
  its {\it derived representation}. \index{derived representation} 

  \begin{prop} [The derived action] 
    \mlabel{prop:3e.1.5} Let $G$ be a Lie group, $M$ a smooth manifold and 
    $E$ a locally convex space. 
    \begin{enumerate}
    \item[\rm(i)] If $\sigma \: M \times G \to M$ is a smooth right action of $G$
      on $M$, then 
      $$ \dot\sigma\:  \L(G) \to {\cal V}(M), \qquad 
      \dot\sigma(x)_p :=  T_{(p,\be)}(\sigma)(0,x) $$ 
      defines a homomorphism of Lie algebras. 
    \item[\rm(ii)] If $\sigma \: G \times M \to M$ is a smooth left action of $G$
      on $M$, then 
      $$ \dot\sigma\:  \L(G) \to {\cal V}(M), \qquad 
      \dot\sigma(x)_p :=  -T_{(\be,p)}(\sigma)(x,0) $$ 
      defines a homomorphism of Lie algebras. 
    \item[\rm(iii)] If $\pi \: G \to \GL(E)$ is a smooth representation of $G$ on $E$, 
      i.e., $\sigma(g,v) := \pi(g)v$ defines a smooth action of $G$ on $E$, then 
      $$ \L(\pi)(x) v := T_{(\be,v)}(\sigma)(x,0) =  - \dot\sigma(x)(v)$$ 
      defines a homomorphism of Lie algebras $\L(\pi) \: \L(G) \to \gl(E).$
    \end{enumerate}
  \end{prop} 

  \begin{prf} (i) We pick $p \in M$ and write 
    $\sigma^p \: G \to M, g \mapsto p.g := \sigma(p,g)$ for the smooth orbit map of $p$. 
    Then $\sigma^p \circ \lambda_g = \sigma^{p.g}$ leads to 
    \[ T(\sigma^p)(x_l(g)) = T_\be(\sigma^{p.g})x = \dot\sigma(x)_{p.g} 
    = \dot\sigma(x)_{\sigma^p(g)}, \] 
    which means that the vector fields $x_l$ and $\dot\sigma(x)$ are 
    $\sigma^p$-related. We conclude that, for $x,y \in \L(G)$, the vector fields 
    $[x_l, y_l] = [x,y]_l$ and $[\dot\sigma(x), \dot\sigma(y)]$ are also $\phi^p$-related 
    (Lemma~\ref{larelglob}). This  leads to 
    $$ [\dot\sigma(x), \dot\sigma(y)](p) 
    = T_\be(\sigma^p)[x,y]_l(\be) 
    = T_{(p,\be)}(\sigma)(0,[x,y]) 
    = \dot\sigma([x,y])(p). $$

    (ii) If $\sigma$ is a left action, then $\sigma^\vee(m,g) := \sigma(g^{-1},m)$ 
    defines a right action, and $T_\be(\eta_G)x = -x$ implies that 
    $$ \dot\sigma(x) 
    =  -T_{(\be,p)}(\sigma)(x,0) 
    =  T_{(p,\be)}(\sigma^\vee)(0,x), $$
    so that the assertion follows from (i). 

    (iii) For linear vector fields 
    $X_A(x) = Ax$, $X_B(x) = Bx$, $A,B \in {\cal L}(E)$, we have 
    $$ [X_A,X_B](x) = d X_B(x)X_A(x)-d X_A(x)X_B(x) 
    = BAx-ABx =-[A,B]x, $$
    so that the corresponding map ${\cal L}(E) \to {\cal V}(E), A \mapsto -X_A$ is a 
    homomorphism of Lie algebras. Therefore (iii) follows from (ii). 
  \end{prf}

  \begin{prop} \mlabel{prop:tangentgrp} 
    Let $G$ be a Lie group and $TG$ its tangent bundle. 
    Then the tangent map 
    $$T(m_G) \: TG \times TG \to TG, \quad 
    (v_g, w_h) \mapsto T(\rho_h)v_g + T(\lambda_g)w_h, $$
    for $v_g \in T_g(G), w_h \in T_h(G),$ defines a Lie group 
    structure on $TG$ with identity element $0 \in T_\be(G)$ and inversion 
    $T(\eta_G)$. 

    Furthermore, the 
    zero section $\sigma \: G \to TG, g \mapsto 0_g,$ and the bundle projection 
    $\pi \: TG \to G$ are morphisms of Lie groups with $\pi \circ \sigma = \id_G$.
  \end{prop}

  \begin{prf} We identify $T(G \times G)$ canonically with 
    $T(G) \times T(G)$ (cf.\ Remark~\ref{rem:tangbun-prod}). 
    Since the multiplication map $m_G \: G \times G \to G$ is smooth, the
    same holds for its tangent map 
    $$ T(m_G) \: T(G \times G) \cong TG \times TG \to TG. $$

    Let $\eps_G \: G  \to G, g \mapsto \be$, be the constant homomorphism. 
    Then the group axioms are encoded in the relations 
    \begin{enumerate}
    \item[\rm(a)] $m_G \circ (m_G \times \id_G) = m_G \circ (\id_G \times m_G)$ (associativity), 
    \item[\rm(b)] $m_G \circ (\eta_G,\id_G) = m_G \circ (\id_G,\eta_G) = \eps_G$ 
      (inversion), and 
    \item[\rm(c)] $m_G \circ (\eps_G, \id_G) = m_G \circ (\id_G, \eps_G) = \id_G$ (unit element).
    \end{enumerate}

    \nin Applying the functor $T$ to these relations, 
    we see that $T(m_G)$ defines a Lie group structure on 
    $T(G)$ for which $T(\eta_G)$ is the inversion and $0_\be \in T_\be(G)$ is the identity. 

    The product of two elements $v_g \in T_g(G)$ and $w_h \in T_h(G)$ can be written as 
    \begin{eqnarray*}
      v_g\cdot w_h 
      &=& T(m_G)(v_g, w_h) 
      = T(m_G)(v_g, 0_h) + T(m_G)(0_g, w_h) \\
      &=& T(\rho_h)v_g + T(\lambda_g)w_h \in T_{gh}(G).
    \end{eqnarray*}
    From that we derive in particular that $\pi$ is a group homomorphism, and 
    from  $0_g \cdot 0_h = T(\rho_h)0_g + T(\lambda_g)0_h = 0_{gh}$ it also follows 
    that the zero section $\sigma$ is a group homomorphism. Both are smooth 
    and satisfy $\pi \circ \sigma = \id_G$ by definition. 
  \end{prf} 

  \begin{defn} \mlabel{def:3.1.2a} 
    Let $G$ be a Lie group. From the Lie group structure on $T(G)$, we obtain 
    a smooth left action 
    $$ G \times T(G) \to T(G), \quad (g,v) \mapsto gv := T(\lambda_g)v $$ 
    and a smooth right action 
    $$ T(G) \times G \to T(G), \quad (v,g) \mapsto vg := T(\rho_g)v. $$ 
    Restricting to tangent vectors in $\be$, we 
    obtain diffeomorphisms 
    $$ G \times \g \to T(G), \quad (g,x) \mapsto gx 
    \quad \mbox{ and } \quad 
    \g \times G \to T(G), \quad (x,g) \mapsto xg. $$
  \end{defn}

  \begin{defn} \mlabel{def:3.1.2} 
    Let $G$ be a Lie group with Lie algebra $\g$. Then, for each $g \in G$, the map 
    $$ c_g \: G \to G, \quad x \mapsto gxg^{-1},$$ 
    is a smooth automorphism, hence induces a continuous linear
    automorphism 
    $$ \Ad(g) := \L(c_g)  \: \g \to \g $$
    and we thus obtain the {\it adjoint representation} \index{adjoint representation} 
    $\Ad \: G \to \Aut(\g).$
    Note that 
    $$ \Ad(g_1 g_2) = \L(c_{g_1 g_2}) = \L(c_{g_1} c_{g_2}) =\L(c_{g_1}) \L(c_{g_2}) 
    = \Ad(g_1) \Ad(g_2) $$
    follows immediately from the functoriality of $\L$. 

    If we identify $G$ via the zero section with a subset of $TG$, then 
    we have in terms of the multiplication on $TG$ (Proposition~\ref{prop:tangentgrp}) 
    the relation 
    $$ \Ad(g)x = gxg^{-1} \quad \mbox{ for } \quad g \in G, x \in \g = T_\be(G), $$
    which implies in particular that the {\it adjoint action} \index{adjoint action} 
    $ \sigma_{\Ad} \: G \times \g \to \g, \break (g, x) \mapsto \Ad(g)x $ 
    defined by the adjoint representation is smooth. 
  \end{defn}

  \begin{prop} \mlabel{prop:der-Ad} 
    The derived representation $\L(\Ad)$ of the adjoint representation is given by 
    $$ \ad x(y) = [x,y]. $$
    In particular, we have for $g \in G$ and $x,y  \in \g$ the relation 
    $$ T_{(g,y)}(\sigma_{\Ad})(gx,z) = \Ad(g)([x,y]+z). $$
  \end{prop}

  \begin{prf} To calculate the linear maps $\L(\Ad)(x)$, 
    we consider a local chart $\phi \: U \to G$ with 
    $\phi(0) = \be$ and $d\phi(\be) = \id_\g$ and the local multiplication 
    $x * y = \phi^{-1}(\phi(x)\phi(y))$. From Remark~\ref{rem:brack-taylor} 
    and Lemma~\ref{lem:second-order} we derive 
    $$ \L(\Ad)(x)y 
    = \hbox{$\frac{\partial}{\partial s}$\vrule}_{s=0} 
    \hbox{$\frac{\partial}{\partial t}$\vrule}_{t=0} 
    sx * ty * (sx)^{-1} 
    = \hbox{$\frac{\partial^2}{\partial s \partial t}$\vrule}_{s,t=0} 
    sx * ty * (sx)^{-1} 
    = [x,y]. $$
    Since $\Ad$ is a group homomorphism, it follows that 
    \[ T_{(g,y)}(\sigma_{\Ad})(gx,z) = \Ad(g) \circ T_{(\be,y)}(\sigma_{\Ad})(x,z) 
    = \Ad(g)([x,y]+z).\qedhere \]
  \end{prf} 

  \begin{ex} (a) If $G = \cA^\times$ is the unit group of a continuous inverse algebra 
    (cf.\ Example~\ref{ex:abgrp}(b)), then the conjugation 
    maps $c_g(a) = gag^{-1}$ are linear in the canonical chart given 
    by the inclusion $\cA^\times \into \cA$, which implies that 
    $$\Ad(c_g)x = gxg^{-1}\quad \mbox{ for } \quad g \in \cA^\times, x \in \cA 
    \cong \L(G). $$

    (b) If $E$ is a Banach space, then the algebra $\cA = \cL(E)$ of 
    bounded linear operators on $E$ is a Banach algebra, hence in 
    particular a continuous inverse algebra. Its unit group is the group 
    $\GL(E)$ of invertible continuous linear operators on $E$. 
    Its Lie algebra is $\gl(E) = (\cL(E), [\cdot, \cdot])$, 
    and by specializing from (a), we see that the adjoint action 
    of this group is given by $\Ad(g)x = gxg^{-1}$. 

    (c) Specializing further to $E= \K^n$ for $\K \in \{\R,\C\}$, 
    the preceding argument also applies to the matrix groups 
    $G = \GL_n(\R)$ and $G = \GL_n(\C)$. 
  \end{ex}

  \subsection{Semidirect products} 

  The easiest way to construct a new Lie group 
  from two given Lie groups $G$ and $H$,  
  is to endow the product manifold $G \times H$ with the multiplication 
  $$ (g_1, h_1) (g_2, h_2) := (g_1 g_2, h_1 h_2). $$
  The resulting group is called the {\it direct product} \index{direct product} 
of the Lie groups 
  $G$ and~$H$. Here $G$ and $H$ can be identified with normal subgroups 
  of $G \times H$ for which the multiplication map is bijective. 
  Relaxing this condition in the sense that only one factor is assumed to 
  be normal, leads to the concept of a semidirect product of Lie groups, 
  that we discuss in this subsection. 

  \begin{defn} \mlabel{def:semdir} Let $N$ and $G$ be Lie groups and 
    $\alpha \: G \to \Aut(N)$
    be a group homomorphism defining a smooth action 
    $(g,n) \mapsto \alpha(g)(n)$ of $G$ on $N$. Then the product manifold 
    $N \times G$ is a group with respect to the product 
    $$ (n,g) (n',g') := (n \alpha(g)(n'), gg') $$
    and the inversion 
    $$ (n,g)^{-1} = (\alpha(g^{-1})(n^{-1}), g^{-1}). $$
    Since multiplication and inversion are smooth, this group is a Lie group, 
    called the {\it semidirect product of $N$ and $G$ with respect to $\alpha$.} 
    We denote this group by \index{semidirect product of Lie groups} 
    $N \rtimes_\alpha G.$
  \end{defn}

  \begin{prop} \mlabel{prop:liealg-semdir} The Lie 
    algebra of the semidirect product 
    group $N \rtimes_\alpha G$ is given by 
    $$ \L(N\rtimes_\alpha G) \cong \L(N) \rtimes_{\L(\alpha^L)} \L(G), $$
    where $\L(\alpha^L) \: \L(G) \to \der(\L(N))$ is the derived representation 
    of $\L(G)$ on $\L(N)$ corresponding to the smooth representation 
    $\alpha^L \: G \to \Aut(N)$, defined by 
    $\alpha_g^L := \L(\alpha_g) \in \Aut(\L(N))$ for $g \in G$. 
  \end{prop}

  \begin{prf} The smooth action of $G$ on $N$ induces a smooth action 
    of $G$ on the tangent bundle $T(N)$ (Exercise~\ref{exer:3.2.1}), 
    and on $T_\be(N) \cong \L(N)$ we thus obtain 
    the smooth representation $g \mapsto \alpha_g^L$. 
    First we show that $\im(\L(\alpha^L)) \subeq \der(\L(N))$. 
    For $x \in \L(G)$, let $\gamma \: [0,1] \to G$ be a smooth curve with 
    $\gamma(0) = \be$ and $\gamma'(0) = x$. For $y,z \in \L(N)$, we then have 
    $$ \alpha_{\gamma(t)}^L[y,z] = [\alpha_{\gamma(t)}^L y,\alpha_{\gamma(t)}^L z] $$ 
    for each $t \in [0,1]$. For the derivative 
    $D_x := \L(\alpha^L)x = \derat0 \alpha_{\gamma(t)}^L$ in $t = 0$, 
    we thus get with the Product Rule
    $$  D_x[y,z] = [D_x y,z] + + [y, D_x z], $$ 
    so that $D_x \in \der(\L(N))$. 

    Next we choose local charts $(\phi_N, U_N)$ of $N$ and $(\phi_G, U_G)$ of $G$ with 
    $$ \phi_N(0) = \be, \quad T_0(\phi_N) = \id_{\L(N)}, 
    \quad 
    \phi_G(0) = \be\quad \mbox{ and } \quad T_0(\phi_G) = \id_{\L(G)}, $$ 
    so that the multiplication 
    $$ (n,g) (n',g') := (n \alpha_g(n'), gg') $$
    can be expressed in local coordinates as 
    $$ (x,y) * (x',y') = (x * \tilde\alpha_y (x'), y * y'), $$
    where $\tilde\alpha_y(x) = \phi_N(\alpha_{\phi_G(y)}(\phi_N^{-1}(x))).$
    For the function $f(x,y) = \tilde\alpha_y(x)$ we then have 
    $$ f(x,0) = \tilde\alpha_0(x) = x 
    \quad \mbox{ and } \quad 
    f(0,y) = \tilde\alpha_y(0) = 0. $$
    Therefore Lemma~\ref{lem:second-order} and 
    $$ \hbox{$\frac{\partial^2}{\partial s \partial t}$\vrule}_{s,t=0}  \ 
    \tilde\alpha_{ty}(sx)  
    = \hbox{$\frac{\partial}{\partial t}$\vrule}_{t=0}  \ 
    \tilde\alpha(ty)^L x 
    = D_yx $$ 
    imply that the second order 
    Taylor polynomial of $f$ in $(0,0)$ is given by 
    $$ x + D_y x. $$
    We conclude that the second order Taylor polynomial of 
    $(x,y) * (x',y')$ in $((0,0), (0,0))$ is of the form 
    $$ (x + x' + D_y x' + b_2^N(x,x'), y + y' + b_2^G(y,y')), $$
    so that we derive from Remark~\ref{rem:brack-taylor} 
    \begin{eqnarray*}
      &&  [(x,y), (x',y')] \\
      &=& (D_yx' + b_2^N(x,x'), b_2^G(y,y'))
      - (D_{y'}x + b_2^N(x',x), b_2^G(y',y))\\
      &=& (D_y x' - D_{y'}x + [x,x'], [y,y']). 
    \end{eqnarray*}
    This proves the assertion. 
  \end{prf}

  \begin{ex}
    \mlabel{ex:affgrp}
    Let $E$ be a Banach space and $\Aff(E)$ be its affine group. 
    The subgroup of translations is isomorphic to $(E,+)$, which 
    leads to a group isomorphism 
    $$ \Aff(E) \cong E \rtimes_\alpha \GL(E), $$
    where $\alpha(g)v = gv$ denotes the canonical action of $\GL(E)$ on $E$. 
    Since this action is smooth, $\Aff(E)$ carries a natural Lie group
    structure. Its Lie algebra is 
    $$ \aff(E) := \L(\Aff(E)) \cong E \rtimes_{\L(\alpha)} \gl(E), 
    \quad \mbox{ where } \quad 
    \L(\alpha)(X)v = Xv. $$
  \end{ex}

  We conclude this section with some observations concerning the 
  recognition of semidirect products. 

  \begin{rem} If $\hat G := N \rtimes_\alpha G$ is a semidirect product, then 
    $$ \pi \: \hat G \to G, \quad (n,g) \mapsto g, \qquad 
    \sigma \: G \to \hat G, \quad g \mapsto (\be,g) $$
    and  $\iota \: N \to \hat G, n \mapsto (n,\be)$
    are morphisms of Lie groups with $\pi \circ \sigma = \id_G$ and 
    $\iota$ is an isomorphism of~$N$ onto the initial submanifold $\ker \pi$ of $\hat G$. 
  \end{rem}

  \begin{prop} \mlabel{prop:semdir-crit} Let $N, G$ and $\hat G$ be Lie groups and 
    $$ \pi \: \hat G \to G, \quad \sigma \: G \to \hat G \quad 
    \mbox{ and } \quad  \iota \: N \to \hat G $$
    morphisms of Lie groups with $\pi \circ \sigma = \id_G$ such that 
    $\iota$ is an isomorphism of $N$ onto $\ker \pi$ which is an initial submanifold of 
    $G$. Then 
    $$ \alpha \: G \to \Aut(N), \quad \alpha_g(n) := \iota^{-1}(\sigma(g)\iota(n)\sigma(g)^{-1}) $$
    defines a smooth action of $G$ on $N$ and the map 
    $$ \Phi \: N \rtimes_\alpha G \to \hat G, \quad (n, g) \mapsto \iota(n) \sigma(g) $$
    is an isomorphism of Lie groups. 
  \end{prop} 

  \begin{prf} Since $\sigma$ is smooth, the assignment 
    $$ \sigma \: G \times \hat G \to \hat G, \quad (g,x) \mapsto \sigma(g)x \sigma(g)^{-1} $$ 
    defines a smooth action of $G$ on $\hat G$, and the normal subgroup 
    $\im(\iota) = \ker \pi$ is invariant under this action. As it is an initial submanifold, 
    the restriction 
    $$\sigma\res_{G \times \ker \pi} \: G \times \ker \pi \to \ker \pi $$ 
    is smooth. Hence $\alpha$ defines a smooth action of $G$ on $N$. 
    We can thus construct the semidirect product group $N \rtimes_\alpha G$. 
    A direct calculation shows that $\Phi$ is a group homomorphism. 
    Further $(\pi \circ \Phi)(n,g) = g$ implies that $\Phi$ is bijective with inverse 
    given by 
    $\Phi^{-1}(x) = \big(\iota^{-1}(x \sigma(\pi(x))^{-1}) , \pi(x)\big). $
    It follows in particular that $\Phi$ is a diffeomorphism, hence an isomorphism of 
    Lie groups.   
  \end{prf}

  \begin{ex} \mlabel{ex:tangentgrp} 
    Let $G$ be a Lie group and $T(G)$ its tangent Lie group 
    (Proposition~\ref{prop:tangentgrp}). Further, let 
    $\sigma \: G \to T(G)$ be the zero section and $\pi \: T(G) \to G$ the bundle projection. 
    Then $\sigma$ and $\pi$ are morphisms of Lie groups and 
    $T_\be(G)= \g = \ker \pi$ is an initial submanifold of $T(G)$, and the group structure on 
    this subgroup is given by 
    $$T_{(\be,\be)}(m_G)(x,y) = T(\rho_\be)x + T(\lambda_\be)y = x + y$$
    (Proposition~\ref{prop:tangentgrp}). 
    From Proposition~\ref{prop:semdir-crit} we now derive that 
    $T(G)$ is a semidirect product $\g \rtimes_\alpha G$, where 
    $$ \alpha_g x = \sigma(g)x\sigma(g)^{-1} = \L(c_g)x = \Ad(g)x $$
    (Definition~\ref{def:3.1.2}).

    For its Lie algebra, we obtain with Propositions~\ref{prop:liealg-semdir} 
    and \ref{prop:der-Ad}  
    $$ \L(T(G)) 
    \cong |\g| \rtimes_{\L(\Ad)} \g \cong |\g| \rtimes_{\ad} \g, $$
    where we can also write $|\g|$ for the locally convex space underlying 
    $\g$, considered as an abelian Lie algebra. 

    Using the algebra $\R[\eps] = \R \oplus \R \eps$ of 
{\it dual numbers}, \index{dual numbers}  
    where $\eps^2 = 0$, we obtain 
    $$ \L(T(G)) \cong T(\g) := \g \otimes_\R \R[\eps], $$
    where the bracket on the right is given by the $\R[\eps]$-bilinear 
    extension of the bracket on $\g$: 
    $$ [x \otimes a, y \otimes b] = [x,y] \otimes ab. $$
  \end{ex}

  \begin{small}
\subsection*{Exercises for Section~\ref{sec:3.2}} 

We have seen above that tangent bundles of Lie groups carry a natural 
Lie group structure. The following exercise extends this construction to 
smooth actions. 

\begin{exer} \mlabel{exer:3.2.1} Let $\sigma \: G \times M \to M$ be a smooth action of the Lie group 
$G$ on the smooth manifold $M$. Then the tangent map $T\sigma \: TG \times TM \to TM$ is a
smooth action of the Lie group $TG$ on $TM$. Hint: The defining properties of a group 
action can be written as 
$$\sigma \circ (m_G \times \id_M) = \sigma \circ (\id_G \times \sigma) 
\: G \times G \times M \to M $$
and 
$$\sigma \circ (\eps_G \times \id_M) = \pr_M \: G \times M \to M, \quad 
(g,m) \mapsto m. $$
Now apply the functor $T$. 
\end{exer} 

\begin{exer} \mlabel{exer:3.2.2} (Multiplication on $T^2(G)$) 
Show that the adjoint representation of $T(G)$ on its Lie algebra 
$\L(T(G)) \cong \L(G) \times \L(G)$ is given by 
$$ \Ad(x,g)(v,w) = (\Ad(g)v + [x, \Ad(g)w], \Ad(g)w). $$
Conclude that the multiplication of 
$T^2(G) := T(T(G))$, where we write 
$T(G) \cong \L(G) \rtimes_{\Ad} G$, satisfies the commutator formula 
\begin{eqnarray*}
&& \big((0,0),(v,\be)\big) \big((0,w),(0,\be)\big) \big((0,0),(v,\be)\big)^{-1} 
\big((0,w),(0,\be)\big)^{-1} \\
&=& \big(\Ad(v,\be)(0,w),(0,\be)\big) \big((0,-w),(0,\be)\big)
= \big(([v,w], 0),(0,\be)\big), 
\end{eqnarray*}
relating the group commutator on $T^2(G)$ to the Lie bracket on $\L(G)$. 
\end{exer}

\begin{exer} \mlabel{exer:3.2.3} Let $\cA$ be a continuous inverse algebra and $\cA^\times$ its unit group, which is a Lie group. 
Identify $T(\cA^\times)$ in the canonical way with the product set $\cA^\times \times \cA$ 
and show that the Lie group structure on $T(\cA^\times)$ is given by 
$$ (a,x)(b,y) = (ab, ay + xb). $$
\end{exer}

\begin{exer} \mlabel{exer:3.2.3b} Let $\cA$ be a continuous inverse algebra 
and $\R[\eps]$ be the algebra of dual numbers. Show that the tensor 
product algebra 
$$ T(\cA) := \cA \otimes \R[\eps] \cong \cA \oplus \cA \eps $$
is a continuous inverse algebra with $T(\cA)^\times \cong T(\cA^\times)$. 
\end{exer}

\begin{exer} \mlabel{exer:3.2.4} Let $G$ and $N$ be Lie groups and 
$\alpha \: G \to \Aut(N)$ be a homomorphism defining a smooth action of 
$G$ on $N$. Show that we obtain a Lie group structure on the product manifold 
$G \times N$ by 
$$ (g,n) (g',n') := (gg', \alpha_{g'}^{-1}(n)n'). $$
Call this group $G \ltimes_\alpha N$. 
Show furthermore that 
$$ \Phi \: N \rtimes_\alpha G\to G \ltimes_\alpha N, \quad 
(n,g) \mapsto (g, \alpha_g^{-1}n) $$
is an isomorphism of Lie groups. 
\end{exer} 

\begin{exer} \mlabel{exer:smooth-crit}  
Let $E$ be a locally convex space and $G$ be a Lie group. 
We consider a homomorphism $\sigma \: G \to \GL(E), g \mapsto \sigma_g$. 
Show that the corresponding action 
of $G$ on $E$ is smooth if it is smooth on an open neighborhood 
$U_G \times U_E$ of $(\be,0)$ in~$G \times E$. \\ 
Hint: Show first the smoothness on $U_G \times E$ and use that all maps 
$\sigma_g$ are smooth. 
\end{exer}

\end{small}

\section{Constructing Lie group structures on groups} \mlabel{sec:3.3} 

In this subsection we describe some methods to construct 
Lie group structures on groups, starting from a manifold structure 
on some ``identity neighborhood'' for which the group operations 
are smooth close to $\be$. 
The results presented here are classical for finite-dimensional 
groups and the corresponding proofs carry over quite easily. 

\subsection{Lie groups from local data} 

The following theorem is the smooth version of Lemma~\ref{lem:toplocglob}  
from the topological context. Throughout this book, 
it is our main tool to construct Lie group structures on groups. 

\begin{thm} \mlabel{thm:locglob} 
Let $G$ be a group and 
$U = U^{-1}$ be a symmetric subset containing $\be$. 
We further assume that $U$ is a smooth manifold and that 
\begin{description}
\item[\rm(L1)] $D := \{ (x,y) \in U \times U \: xy\in U\}$ is an open subset 
and the multiplication $m_U \: D \to U, (x,y) \mapsto xy$ is smooth, 
\item[\rm(L2)] the inversion map $\eta_U \: U \to U, u \mapsto u^{-1}$ is
smooth, and, 
\item[\rm(L3)] for each $g \in G$, there exists an open $\be$-neighborhood $U_g \subeq
U$ with $c_g(U_g) \subeq U$ and such that the conjugation map 
$c_g \: U_g \to U, x \mapsto gxg^{-1}$
is smooth. 
\end{description}
Then there exists a unique structure of a
Lie group on $G$ such that the inclusion map $U \into G$ is a diffeomorphism 
onto an open subset of $G$. 

If, in addition, $U$ generates $G$, then {\rm(L1/2)} imply {\rm(L3)}.
\end{thm}

\begin{prf}  From Lemma~\ref{lem:toplocglob},  
we obtain a unique group topology on $G$ for which 
the inclusion map $U \into G$ is an open embedding.  

Now we turn to the manifold structure. 
Let $V = V^{-1}\subeq U$ be an open $\be$-neighborhood with 
$VV \times VV \subeq D$, for which there exists a locally convex space $E$ and 
an $E$-chart $(\phi,V)$ of $U$.  
For $g \in G$ we consider the map 
$$ \phi_g \: gV \to E, \quad \phi_g(x) = \phi(g^{-1}x) $$
which is a homeomorphisms of $gV$ onto the open subset 
$\phi(V)\subeq E$. We claim that
$(\phi_g, gV)_{g \in G}$ is a smooth atlas of $G$. 

Let $g_1, g_2 \in G$ and put $W := g_1 V \cap g_2 V$. If $W \not=
\eset$, then $g_2^{-1} g_1 \in V V^{-1} = VV$. The smoothness of the 
map 
$$ \psi := \phi_{g_2} \circ \phi_{g_1}^{-1}\res_{\phi_{g_1}(W)}
\: \phi_{g_1}(W) \to \phi_{g_2}(W) $$
given by 
$$ \psi(x) 
= \phi_{g_2}(\phi_{g_1}^{-1}(x))
= \phi_{g_2}(g_1 \phi^{-1}(x))
= \phi(g_2^{-1} g_1 \phi^{-1}(x)) $$
follows from the smoothness of the
multiplication $VV \times VV \to U$. This proves that the charts 
$(\phi_g, gU)_{g \in G}$ form a smooth atlas of $G$. Moreover, the
construction implies that all left translations of $G$ are smooth
maps because $\phi_g \circ \lambda_h = \phi_{h^{-1}g}$ for 
$g, h \in G$. 
 
The construction also shows that, for $g \in G$, the conjugation $c_g \: G
\to G$ is smooth in a neighborhood of $\be$. Since all left
translations are smooth, and 
$$ c_g \circ \lambda_x = \lambda_{c_g(x)} \circ c_g, $$
the smoothness of $c_g$ in a neighborhood of  $x \in G$ follows. 
Therefore all conjugations and hence also all
right multiplications are smooth. The smoothness of the inversion
follows from its smoothness on $V$ and the fact that left and right
multiplications are smooth. Finally, the smoothness of the
multiplication follows from the smoothness in $\be \times \be$ because 
$$ g_1 x g_2 y 
= g_1 g_2 c_{g_2^{-1}}(x) y. $$ 
We conclude that $G$ is a Lie group. 

Next we show that the inclusion $U \into G$ of $U$ is a diffeomorphism. 
So let $x \in U$ and recall the open subset $U_{x} = U \cap x^{-1}U$. 
 Then $\lambda_x$ restricts to a smooth map 
$U_x \to U$ with image $U_{x^{-1}}$. Its 
inverse is also smooth. We thus obtain a diffeomorphism 
$\lambda_x^U \: U_x \to U_{x^{-1}}$. 
Since $\lambda_x \: G \to G$ also is a diffeomorphism, 
the inclusion 
$\lambda_x \circ \lambda_{x^{-1}}^U \: U_{x^{-1}} \to G$ 
is a diffeomorphism. As $x$ was arbitrary, the inclusion of $U$ in $G$ 
is a diffeomorphic embedding. 

The uniqueness of the Lie group structure is clear because each locally
diffeomorphic bijective homomorphism between Lie groups is a
diffeomorphism (cf.~Exercise~\ref{exer:3.1.4b}). 

Finally, we assume that $G$ is generated by $U$. 
We show that in this case (L3) is a consequence of (L1) and (L2). 
For each $g \in U$, there exists an open $\be$-neighborhood $U_g$ 
with $gU_g \times \{g^{-1}\} \subeq D$. 
Then $c_g(U_g) \subeq U$, and the smoothness of $m_U$ 
implies that $c_g\res_{U_g} \: U_g\to U$ is smooth. 
Hence, for each $g \in U$, the conjugation $c_g$ is smooth in a neighborhood of $\be$. 
Since the set of all these
$g$ is a submonoid of $G$ containing $U$, it contains $U^n$ for each
$n \in \N$, hence all of $G$ because $G$ is generated by $U= U^{-1}$. 
Therefore (L3) is satisfied. 
\end{prf}

\begin{cor} \mlabel{cor:open-lie} 
Let $G$ be a group and $N \trile G$ be a normal subgroup of $G$ that 
carries a Lie group structure. Then there exists a unique Lie group structure on $G$ 
for which $N$ is an open subgroup if and only if, for each $g \in G$, 
the restriction $c_g\res_N$ is a smooth automorphism of $N$. 
\end{cor}

\begin{prf} If $N$ is an open normal subgroup of the Lie group $G$, then clearly 
all inner automorphisms of $G$ restrict to smooth automorphisms of $N$. 

Suppose, conversely, that $N$ is a normal subgroup of the group $G$ which is a 
Lie group, and that all inner automorphisms of $G$ restrict to smooth automorphisms of 
$N$. Then we can apply Theorem~\ref{thm:locglob} with $U = N$ and obtain a unique Lie group 
structure on $G$ for which the inclusion $N \to G$ is a 
diffeomorphism onto an open subgroup of $G$.   
\end{prf}
\begin{rem}
If $\K\in \{\R,\C\}$ and $U$ is endowed with a $\K$-analytic
manifold structure in Theorem~\ref{thm:locglob},
replacing smoothness with $\K$-analyticity in all
hypotheses we obtain a $\K$-analytic Lie group
structure on~$G$ with $U$
as an open submanifold, by the same reasoning.
\end{rem}

\subsection{Coverings of Lie groups}

In this subsection we show that, for a Lie group $G$,  
coverings and quotients by discrete subgroups carry natural Lie group 
structures. This provides the means to classify
Lie groups with a given Lie algebra as quotients of a universal covering group. 

\begin{cor} \mlabel{cor:cover-lie} 
Let $\phi \: G \to H$ be a covering of topological groups. 
If $G$ or $H$ is a Lie group, then the other group has a unique Lie group 
structure for which $\phi$ is a morphism of Lie groups which is a local diffeomorphism. 
\end{cor}

\begin{prf} Let $U_G \subeq G$ be an open symmetric $\be$-neighborhood 
for which $\phi\res_{U_G}$ is a homeomorphism onto an open 
subset $U_H$ of $H$. In addition, we assume thaty $U_G^3 \cap \ker \phi = \{1\}$. 

Suppose first that $G$ is a Lie group. Then we apply 
Theorem~\ref{thm:locglob} to $U_H$, endowed with the manifold structure for 
which $\phi\res_{U_G}$ is a diffeomorphism. Then (L2) follows from 
$\phi(x)^{-1} = \phi(x^{-1})$. To verify the smoothness of the multiplication map 
$$ m_{U_H} \: D_H := \{ (a,b) \in U_H \times U_H \: ab \in U_H \} \to U_H, $$
we first observe that, if $x,y \in {U_G}$ satisfy 
$(\phi(x),\phi(y)) \in D_H$, i.e., $\phi(xy) \in U_H$, then there exists a $z \in {U_G}$ with 
$\phi(xy) = \phi(z)$, and
\[ xyz^{-1} \in {U_G}^3 \cap \ker(\phi)=\{\be\} \]  yields 
$z = xy$. We thus have 
$D_H = (\phi \times \phi)(D_G)$  
for 
$$ D_G := \{ (x,y) \in {U_G} \times {U_G} \: xy \in {U_G}\} $$
and the smoothness of $m_H$ follows from the smoothness of the multiplication 
$m_{U_G} \: D_G \to {U_G}$ and 
$m_{U_H} \circ (\phi \times \phi) = \phi \circ m_{U_G}.$ 
To verify 
(L3), we note that the surjectivity of $\phi$ 
implies that, for each $h \in H$, there is an element $g \in G$ with 
$\phi(g) = h$. Now we choose an open $\be$-neighborhood $U_g \subeq U_G$ with 
$c_{g}(U_g) \subeq U_G$ and put $U_h := \phi(U_g)$. 

If, conversely, $H$ is a Lie group, then we apply 
Theorem~\ref{thm:locglob} to $U_G$, endowed with the manifold structure for 
which $\phi\res_{U_G}$ is a diffeomorphism onto $U_H$. 
Again, (L2) follows right away, and (L1) follows from 
$(\phi \circ \phi)(D_G) \subeq D_H$ and the smoothness of 
$m_{U_H} \circ (\phi \circ \phi) = \phi \circ m_{U_G}.$
For (L3), we choose $U_g$ as any open $\be$-neighborhood in $U_G$ with 
$c_g(U_g) \subeq U$. Then the smoothness of $c_g\res_{U_g}$ follows 
from the smoothness the maps 
of $\phi \circ c_g = c_{\phi(g)} \circ \phi.$
\end{prf}

\begin{cor} \mlabel{cor:unicov} If 
$G$ is a connected Lie group and $q_G \: \tilde G \to G$ its universal 
covering space, then $\tilde G$ carries a unique Lie group structure for which 
$q_G$ is a smooth covering homomorphism.   
\end{cor}

\begin{prf} We first have to construct a (topological) group structure on the universal 
covering space~$\tilde G$. Pick an element $\tilde\be \in q_G^{-1}(\be)$. 
Then the multiplication map $m_G \: G \times G \to G$ 
lifts uniquely to a continuous map $\tilde m_G \: \tilde G \times 
\tilde G \to \tilde G$
with $\tilde m_G(\tilde\be,\tilde\be) = \tilde\be$ (\cite[Thm.~III.4.1]{Bre93}). 
To see that the multiplication map 
$\tilde m_G$ is associative, we observe that 
\begin{eqnarray*}
 && q_G \circ \tilde m_G \circ (\id_{\tilde G} \times \tilde m_G) 
=  m_G \circ (q_G \times q_G) \circ (\id_{\tilde G} \times \tilde m_G) \\
&=&  m_G \circ (\id_{G} \times m_G) \circ (q_G \times q_G \times q_G)
=  m_G \circ (m_G \times \id_{G}) \circ (q_G \times q_G \times q_G)\\
&=&  q_G \circ \tilde m_G \circ (\tilde m_G \times \id_{\tilde G}),
\end{eqnarray*}
so that the two continuous maps 
$$ \tilde m_G \circ (\id_{\tilde G} \times \tilde m_G), \ \ 
\tilde m_G \circ (\tilde m_G \times \id_{\tilde G})  \: \tilde G^3 \to G, $$
are lifts of the same map $G^3 \to G$ 
and both map $(\tilde\be,\tilde\be,\tilde\be)$ to $\tilde\be$. 
Hence the uniqueness of lifts implies that $\tilde m_G$ is associative. 
We likewise obtain that the unique lift $\tilde \eta_G \: \tilde G \to \tilde G$
of the inversion map $\eta_G \: G \to G$ with $\tilde\eta_G(\tilde\be) = \tilde\be$ 
satisfies 
$$ \tilde m_G \circ (\eta_G \times \id_{\tilde G}) = \tilde\be 
=  \tilde m_G \circ (\id_{\tilde G} \times \eta_G). $$
Therefore $\tilde m_G$ defines on $\tilde G$ a topological group structure 
such that \break $q_G \: \tilde G \to G$ is a covering morphism of topological groups. 
Now Corollary~\ref{cor:cover-lie} applies. 
\end{prf}

\begin{prop} \mlabel{prop:3.3.8} 
If $\phi \: G \to H$ is a morphism of connected 
Lie groups and $q_G \: \tilde G \to G$, $q_H \: \tilde H \to H$ 
simply connected coverings, then there exists a unique 
morphism of Lie groups $\tilde \phi \: \tilde G \to \tilde H$ 
with $q_H \circ \tilde \phi = \phi \circ q_G$. 
\end{prop}

\begin{prf} According to \cite[Thm.~III.4.1]{Bre93}, the 
map $\phi \circ q_G \: \tilde G \to H$ has a unique continuous 
lift $\tilde \phi \: \tilde G \to \tilde H$ 
with $\tilde\phi(\be) = \be$. Then 
\begin{equation}
  \label{eq:3.3.8}
m_{\tilde H} \circ (\tilde \phi  \times \tilde \phi) \: 
\tilde G \times \tilde G \to \tilde H 
\end{equation}
is continuous and maps $(\be,\be)$ to $\be$. 
The same holds for the continuous map 
\begin{equation}
  \label{eq:3.3.9}
\tilde\phi \circ m_{\tilde G} \: \tilde G \times \tilde G \to \tilde H. 
\end{equation}
In view of 
\begin{align*}
q_{H} \circ \tilde \phi \circ m_{\tilde G} 
&=  \phi \circ q_{G} \circ m_{\tilde G} 
=  \phi \circ m_G \circ (q_{G} \times q_{G}) \\
&=  m_H \circ (\phi \times \phi) \circ (q_{G} \times q_{G}) 
=  m_H \circ (q_H \circ \tilde\phi \times q_H \circ \tilde\phi)\\ 
&=  q_{H} \circ m_{\tilde H} \circ (\tilde\phi \times \tilde\phi),    
\end{align*}
the maps \eqref{eq:3.3.8} and \eqref{eq:3.3.9} 
are lifts of the same map to $H$, 
so that they are equal. This means that $\tilde\phi$ 
is a group homomorphism. 

By definition, $\tilde \phi$ is continuous, and since 
$q_{H}$ is a local diffeomorphism and $q_H \circ \tilde\phi$  
is smooth, $\tilde \phi$ is smooth 
in an identity neighborhood, hence smooth 
(Exercise~\ref{exer:3.1.4b}; see also Exercise~\ref{exer:3.3.7}). 
\end{prf}

\begin{thm}[Lie groups with the same covering group] \mlabel{thm:3.3.6} 
\index{Lie groups!with same covering group}
Let $q_G \: \tilde G \to G$ be the universal covering morphism 
of the connected Lie group~$G$. Then the following assertions hold: 
\begin{enumerate}
\item[\rm(i)] $\ker q_G \cong \pi_1(G)$ is a discrete central subgroup and in 
particular abelian. 
\item[\rm(ii)] $G \cong \tilde G/\Gamma$ for a discrete central subgroup 
$\Gamma \subeq Z(\tilde G)$. 
\item[\rm(iii)] For two discrete subgroups $\Gamma_1, \Gamma_2 \subeq Z(\tilde G)$, 
the groups $\tilde G/\Gamma_1$ and $\tilde G/\Gamma_2$ are isomorphic if and only 
if there exists an automorphism $\phi \in \Aut(\tilde G)$ with 
$\phi(\Gamma_1) = \Gamma_2$. 
\item[\rm(iv)] The action of $\Aut(\tilde  G)$ on $Z(\tilde G)$ factors through 
an action of the group $\Out(\tilde G) := \Aut(\tilde G)/\Inn(\tilde G)$ 
of out automorphosms, where \[ \Inn(\tilde G) := \{c_g \: g \in \tilde G\} \] 
is the normal subgroup of inner automorphisms. 
The orbits of this group in the set of discrete subgroups of 
$Z(\tilde G)$ parametrize the isomorphy classes of connected Lie groups 
with the universal covering group $\tilde G$. 
\item[\rm(v)] Lifting automorphisms of $G$ to $\tilde G$ yields an embedding 
$\Aut(G) \into \Aut(\tilde G)$ whose image coincides with 
the normalizer of $\ker q_G$: 
\[ N_{\Aut(\tilde G)}(\ker q_G) = \{ \phi \in \Aut(\tilde G) \: \phi(\ker q_G) = \ker q_G\}.\] 
\end{enumerate}
\end{thm}

\begin{prf} (i) Since $\ker q_G$ is a discrete normal subgroup of the connected group $\tilde G$, 
it is central by Exercise~\ref{exer:3.3.4}. Left multiplications 
by elements of this group lead to deck transformations of the covering $\tilde G \to G$, 
and this yields a group isomorphism $\pi_1(G) \cong \ker q_G$. 
In particular, we see that $\pi_1(G)$ is abelian. 

(ii) follows from $G \cong \tilde G/\ker q_G$ and (i). 

(iii) follows from Proposition~\ref{prop:3.3.8}. 

(iv) Since the normal subgroup 
$\Inn(\tilde G) := \{ c_g \: g \in \tilde G\}$ of 
inner automorphisms acts trivially on $Z(\tilde G)$, the action of 
$\Aut(\tilde G)$ on $Z(\tilde G)$ factors through an action of 
$\Out(\tilde G)$. 

(v) That each automorphism $\phi \in \Aut(G)$ lifts to a unique automorphism 
$\tilde\phi \in \Aut(\tilde G)$ follows from Exercise~\ref{exer:3.3.5}. 
Conversely, an automorphism of $\tilde G$ is such a lift if and only 
if it normalizes $\ker q_G$. 
\end{prf}

\begin{exs} \mlabel{exs:ab-quot} 
Let $E$ be a locally convex space, so that 
$(E,+)$ is a Lie group. Then, for each discrete subgroup $\Gamma \subeq E$, 
the quotient group $E/\Gamma$ carries a natural Lie group structure 
(Corollary~\ref{cor:cover-lie}). 

We shall see in Proposition~\ref{prop:reg-ab} below that all connected abelian Lie groups 
with an exponential function are of this form. 

(a) As special cases we obtain in particular the finite-dimensional tori 
$$ \T^d \cong \R^d/\Z^d. $$

(b) We consider the countably-dimensional vector space $E := \R^{(\N)}$ 
with basis $(e_n)_{n \in \N}$, endowed with the finest locally convex topology, 
i.e., each seminorm on $E$ is continuous. 
Then the subgroup 
$\Gamma := \Z^{(\N)}$ (the free abelian group generated by the basis)  
is discrete because $]-1,1[^\N \cap E$ is a convex $0$-neighborhood in $E$,  
intersecting $\Gamma$ only in $\{0\}$. Therefore 
the quotient 
$$ \T^{(\N)} := \R^{(\N)}/\Z^{(\N)} $$
carries a natural Lie group structure. This group is a directed union of 
finite-dimensional tori. 

(c) If $E = \R^\N$ is endowed with the product topology, then 
the subgroup $\Gamma := \Z^\N$ is not discrete because the sequence 
$(e_n)_{n \in \N}$ in $\Gamma$ converges to $0$. The quotient 
group $\T^\N \cong \R^\N/\Z^\N$ is a compact group which obviously is 
not a Lie group because each compact Lie group is modeled on a locally 
compact space, hence finite-dimensional (\cite[Thm.~1.22]{Ru91}). 

However, the topology on the group $\T^\N$ can be refined 
to a Lie group topology. In the Banach space 
$E := \ell^\infty(\N,\R)$ of all bounded sequences, endowed with the sup-norm 
$\|\cdot\|_\infty$, the subgroup 
$\Gamma := \ell^\infty(\N,\Z)$ of all bounded $\Z$-valued sequences  
is discrete, so that the quotient group 
$$ \T^\N \cong  \ell^\infty(\N,\R)/\ell^\infty(\N,\Z)$$
inherits a Lie group structure. The corresponding topology on $\T^\N$ 
is much finer than the product topology. 
\end{exs}

\begin{ex} \mlabel{ex:exotic-units} (a) Let $\cA$ be a commutative unital locally convex algebra 
with an exponential function 
$$ \exp_\cA \: \cA \to \cA^\times, $$
i.e., which means that $\exp_\cA \: (\cA,+) \to (\cA^\times, \cdot)$ 
is a group homomorphism with 
$T_0(\exp_\cA)= \id_\cA$. 

Then $\Gamma_\cA := \ker(\exp_\cA)$ is a closed subgroup of $\cA$, not containing any 
line. Suppose that $\Gamma_\cA$ is discrete. Then 
$Q := \cA/\Gamma_\cA$ carries a natural Lie group structure (Example~\ref{exs:ab-quot})  
and the exponential function factors through an injection $Q \into \cA^\times$.
We may therefore use Corollary~\ref{cor:open-lie} to define a Lie group structure on the group 
$\cA^\times$ for which the identity component is $\exp_\cA(\cA) \cong Q$ and 
$\exp_\cA$ is a covering morphism.

(b) Let $M$ be a finite-dimensional $\sigma$-compact manifold. Then 
the space $\cA := C(M,\C)_{c.o.}$ of all continuous functions $M \to \C$, 
endowed with the topology of uniform convergence on compact subsets, defined by the seminorms 
$$ p_K(f) :=  \sup\{ |f(k)| \: k \in K\}, 
\quad K \subeq M \ \mbox{ compact},  $$
is a locally convex algebra. 

If $M$ is non-compact, then there are unbounded continuous functions on $M$, 
and this implies that the unit group $\cA^\times$ is not open because for each unbounded function $f \: M \to \C$, the element 
$\be + \lambda f$ is not invertible for some 
$\lambda \in \C$ arbitrarily close to $0$. 
We conclude that $\cA$ is a continuous inverse algebra if and only if $M$ is compact, in which case it 
even is a unital Banach algebra. 

We write 
$$ \exp_\cA \: \cA \to \cA^\times = C(M,\C^\times), \quad f \mapsto e^f $$
for the exponential function of $\cA$, which is easily seen to be smooth. 
Its kernel 
$$\Gamma_\cA = \ker(\exp_\cA) = C(M,2\pi i \Z) \cong C(M,\Z)$$ 
is discrete if and only if $M$ has only finitely many connected components 
(Exercise~\ref{exer:3.3.6}). 
If this is the case, (a) implies that $\cA^\times$ carries a Lie group structure for which 
$\exp_\cA$ is a covering map. 

A typical example is $M = \R$ and $\cA = C^\infty(\R,\C)$ with 
$$ \cA^\times = C^\infty(\R,\C^\times) 
\cong C^\infty_*(\R,\C^\times) \times \C^\times,  $$
where $C^\infty_*$ denotes functions mapping $0$ to $1$, resp., to $0$. 
Since all elements of $C^\infty_*(\R,\C^\times)$ have a unique lift to an 
element of $C^\infty_*(\R,\C)$, and this induces a homeomorphism for the 
compact open topology, we get 
$$ \cA^\times \cong C^\infty_*(\R,\C) \times \C^\times. $$
For $M = \N$, we have $\cA \cong \C^\N$, and $\Gamma_\cA \cong \Z^\N$ is not discrete. 

If $M$ is connected and $\sigma$-compact, it is not hard to see that the logarithmic derivative 
\[ \delta \: \cA^\times \cong C^\infty(M,\C^\times) \to Z^1_{\rm dR}(M,\C), \quad 
f \mapsto \frac{d f}{f}\] 
induces a topological isomorphism of $\cA^\times/\C^\times$ onto the 
group $Z^1_{\rm dR}(M,\Z)$ of closed $1$-forms whose periods are contained in 
$2\pi i \Z$, and the arc-component $\cA^\times_a$ of the identity is mapped onto 
the set of exact $1$-forms (cf.\ Corollary~\ref{cor:fundamental} below). 
We conclude that, as topological groups, 
$$ \pi_0(\cA^\times) \cong \cA^\times/\cA^\times_a \cong H^1_{\rm dR}(M,\Z) 
\cong \Hom(\pi_1(M),\Z) \cong \Hom(H_1(M),\Z), $$
and this group is discrete with respect to the topology of pointwise 
convergence if and only if the singular cohomology group 
$$H^1_{\rm sing}(M,\Z) \cong \Hom(H_1(M),\Z)$$ 
is finitely generated 
(cf.\ \cite{NeW08b}). This shows that the arc-component of the identity in 
$\cA^\times$ is open if and only if the abelian group 
$H^1_{\rm sing}(M,\Z)$ is finitely generated. 

For $M := \C\setminus \N$, the group  
$H_1(M) \cong \Z^{(\N)}$ is of infinite rank, 
$H^1_{\rm sing}(M,\Z) \cong \Z^\N$ is not discrete, 
but $M$ is connected, so that $\cA^\times$ carries a Lie group structure whose underlying 
topology is finer than the subspace topology of $\cA^\times$ induced from $\cA$. 
\end{ex}

\begin{small}
\subsection*{Exercises for Section~\ref{sec:3.3}}

\begin{exer} \mlabel{exer:3.3.3}
Let $G$ be an abelian group and $N \leq G$ a subgroup carrying a Lie 
group structure. Then there exists a unique Lie group structure on $G$ for which 
$N$ is an open subgroup. 
\end{exer}

\begin{exer} \mlabel{exer:3.3.6} Let $M$ be a topological finite-dimensional 
manifold. Show that the group $C(M,\Z)_{c.o.}$ (carrying the compact open 
topology) is discrete if and only if 
$M$ has only finitely many connected components. 
\end{exer} 

\begin{exer} \mlabel{exer:3.3.7} 
Let $M$, $N$ and $X$ be smooth manifolds, 
$q \: M \to N$ be a smooth covering map and 
$F \: X \to M$ a continuous map. Show that 
$F$ is smooth if and only if $q \circ F$ is smooth. 
\end{exer}

%
\end{small}

\section{Groups of maps and diffeomorphisms} 
\mlabel{sec:grpsmapsanddiffeo}

In this section, we introduce, for a smooth manifold and a Lie group 
$K$, the {\it compact open $C^r$-topology} \index{ compact open $C^r$-topology} 
on the group $C^r(M,K)$, 
which is a suitable refinement of the compact open topology. 
If $M$ is compact (possibly with boundary) 
and $K$ is a Lie group, this topology is compatible 
with a natural Lie group structure on this group.
The current outline is independent of the
more exhaustive discussion in Chapter~\ref{chap-manifold-constructions},
which is needed for more technical results like
the discussion of regularity properties in Section~\ref{tools-regularity}, Chapter~\ref{ch:mapgrp}, and
Chapter~\ref{ch:diffeo}. Moreover, the local-to-global approach
is an important alternative, which avoids putting a manifold structure
on general sets of maps.
\begin{defn} \mlabel{def:top-mapgrp} (Groups of differentiable maps as topological groups)  
Let $M$ be a smooth manifold,  
$K$ a Lie group with Lie algebra $\fk$ and $r \in \N_0 \cup \{\infty\}$. 
We endow the group $G := C^r(M,K)$ with the compact open $C^r$-topology 
(Definition~\ref{def:smooth-co-top}). 

We know already that the tangent bundle $TK$ of $K$ is a Lie group (Proposition~\ref{prop:tangentgrp}). 
Iterating this procedure, we obtain a Lie group structure on all higher tangent bundles $T^n K$. 
For each $n \in \N_0$, we thus obtain topological groups $C(T^n M, T^n K)_{c.o.}$ 
by using the compact open topology (Definition~\ref{def:smooth-co-top}). 
We also observe that, for two smooth maps $f_1, f_2 \: M \to K$,  
the functoriality of $T$ yields 
$$ T(f_1\cdot f_2) = T(m_G \circ (f_1 \times f_2)) 
=T(m_G) \circ (Tf_1 \times Tf_2) = Tf_1 \cdot Tf_2. $$
Therefore the inclusion map 
$$ C^r(M,K) \into \prod_{n = 0}^r  C(T^n M, T^n K)_{c.o.}, \quad 
f \mapsto (T^n f)_{0 \leq n \leq r} $$
is a group homomorphism, so that the inverse image of the product 
topology from the right hand side is a 
group topology on $C^r(M,K)$. Hence 
the compact open $C^r$-topology turns $C^r(M,K)$ into a topological group, 
even if $M$ and $K$ are infinite dimensional. 
\end{defn}

\begin{thm} \mlabel{thm:mapgro-Lie} 
  Let $K$ be a Lie group with Lie  algebra $\fk$,  
$M$ a compact manifold (possibly with boundary),  
and $r \in \N_0\cup \{\infty\}$. 
Then 
$C^r(M,K)$ carries a Lie group structure compatible with the compact open 
$C^r$-topology, 
and its Lie algebra is $C^r(M,\fk)$, endowed with the pointwise bracket. 
\end{thm}

\begin{prf} (Sketch) Let $G := C^r(M,K)$ and $\g:= C^r(M,\fk)$. 
The Lie group structure on $G$ can be constructed with 
Theorem~\ref{thm:locglob} 
as follows. 
Let $\phi_K \: U_K \to\fk$ be a chart of $K$. Then the set 
$U_G := \{ f \in G \: f(M) \subeq U_K\}$ is an open subset of 
$G$. Assume, in addition, that $\be \in U_K$ and $\phi_K(\be)=0$. 
Then the map 
$$ \phi_G \: U_G \to \g, \quad f \mapsto \phi_K \circ f $$
defines a chart $(\phi_G, U_G)$ of $G$ (cf.\ Exercise~\ref{exer:4.1.6}). 
To apply Theorem~\ref{thm:locglob}, one 
has to verify that in this chart 
the inversion is a smooth map, that the multiplication map 
$$ D_G := \{ (f,g) \in U_G \times U_G \: fg \in U_G \} \to U_G $$
is smooth and that conjugation maps are smooth 
in some $\be$-neighborhood of $U_G$ 
(Propositions~\ref{pushf-cp}). 
%
%
To calculate the Lie algebra of this group, we observe that, for 
$m \in M$, we have for the multiplication in local coordinates 
\begin{eqnarray*}
(f *_G g)(m) 
&:=& \phi_G\Big(\phi_G^{-1}(f)\phi_G^{-1}(g)\Big)(m) 
= \phi_K\big(\phi_K^{-1}(f(m))\phi_K^{-1}(g(m))\big) \cr
&=& f(m) *_K g(m)= f(m) + g(m) + b_\fk(f(m), g(m)) + \cdots.
\end{eqnarray*}
In view of Remark~\ref{rem:brack-taylor}, this implies that 
$b_\g(f,g)(m) = b_\fk(f(m),g(m)),$
and hence that 
\begin{eqnarray*}
[f,g](m) &=& b_\g(f,g)(m) - b_\g(g,f)(m) \\ 
&=& b_\fk(f(m),g(m))-b_\fk(g(m),f(m)) = [f(m),g(m)].
\end{eqnarray*}
Therefore $\L(C^r(M,K)) = C^r(M,\fk)$, endowed 
with the pointwise defined Lie bracket. 
\end{prf}

\begin{ex}
  \mlabel{ex:III.1.19} Let $M$ be a compact manifold and 
$\g = {\cal V}(M)$, the Lie algebra of smooth vector fields on $M$. 
We now explain how the group 
$\Diff(M)$ can be turned into a Lie group, modeled on~$\g$. 

We shall see in Chapter~\ref{ch:diffeo} 
 that, although $\Diff(M)$ has a smooth exponential function, 
it is not a local diffeomorphism 
of a $0$-neighborhood in $\g$ onto an identity neighborhood in $G$. Therefore we cannot 
use it to define charts for $G$. But there is an easy way around this problem. 

Let $g$ be a Riemannian metric on $M$ and 
$\Exp \: TM \to M$
be its exponential function, which assigns to $v \in T_m(M)$ the point 
$\gamma(1)$, where $\gamma \: [0,1] \to M$ is the geodesic segment with 
$\gamma(0) = m$ and $\gamma'(0) = v$. We then obtain a smooth map 
$$ \Phi \: TM \to M \times M, \quad v \mapsto (m, \Exp v), \quad v \in T_m(M). $$
There exists an open neighborhood $U \subeq TM$ of the zero section such that 
$\Phi$ maps $U$ diffeomorphically onto an open neighborhood of the diagonal in $M \times M$.
Now 
$$ U_\g := \{ X \in {\cal V}(M) \: X(M) \subeq U\} $$
is an open subset of the Fr\'echet space ${\cal V}(M)$, and we define a map 
$$ \phi \: U_\g \to C^\infty(M,M), \quad \phi(X)(m) := \Exp(X(m)). $$
It is clear that $\phi(0) = \id_M$. 
One can show that after shrinking $U_\g$ to a sufficiently small 
$0$-neighborhood in the $C^1$-topology on ${\cal V}(M)$, 
we may achieve that $\phi(U_\g) \subeq \Diff(M)$. 
To see that $\Diff(M)$ carries a Lie group 
structure for which $\phi$ is a chart, one has to verify that the 
group operations are smooth in a $0$-neighborhood 
when transfered to $U_\g$ via $\phi$, so that  
Theorem~\ref{thm:locglob} below applies. 
We thus obtain a Lie group structure on $\Diff(M)$. 
    
From the smoothness of 
$U_\g \times M \to M, (X,m) \mapsto \phi(X)(m) = \Exp(X(m))$ it follows 
that the canonical left action $\sigma \: \Diff(M) \times M \to M, 
(\phi,m) \mapsto \phi(m)$ is smooth in an identity neighborhood of 
$\Diff(M)$, and hence smooth, because it is an action by smooth maps. 
The corresponding homomorphism of Lie algebras 
\[ \dot\sigma \: \Lie(\Diff(M)) \to {\cal V}(M) \] 
is given by 
$$ \dot\sigma(X)(m) = -T\sigma(X,0_m) = -(d\Exp)_{0_m}(X(m)) = -X(m), $$
i.e., $\dot\sigma = -\id_{{\cal V}(M)}$. 
This leads to 
$$ \Lie(\Diff(M)) = ({\cal V}(M),[\cdot,\cdot])^{\rm op}. $$

This ``wrong'' sign is caused by the fact that we consider 
$\Diff(M)$ as a group acting on $M$ from the left. If we consider 
it as a group acting on the right, we obtain the opposite multiplication 
$$ \phi * \psi := \psi \circ \phi, $$
and 
$$ \Lie(\Diff(M)^{\rm op}) \cong  ({\cal V}(M),[\cdot,\cdot]) $$
follows from Proposition~\ref{prop:3e.1.5}. 
The tangent bundle of $\Diff(M)$ can be identified with the set 
$$ T(\Diff(M)) := \{ X \in C^\infty(M,TM) \: \pi_{TM} \circ X \in
\Diff(M)\},  $$
where the map 
$$ \pi \: T(\Diff(M)) \to \Diff(M), \quad X \mapsto \pi_{TM}
\circ X $$
is the bundle projection. Then 
$$T_\phi(\Diff(M)) := \pi^{-1}(\phi)
= \{ X \in C^\infty(M,TM) \: \pi_{TM} \circ X = \phi\} $$
is the tangent space in the diffeomorphism $\phi$. 
The multiplication in the group $T(\Diff(M))$ is given by the formula 
$$ X \cdot Y := \pi_{T^2M} \circ TX \circ Y, $$
where $\pi_{T^2M} \: T^2M \to TM$ is the natural projection. 
Note that 
$$ \pi_{TM} \circ (X \cdot Y) 
= \pi_{TM} \circ \pi_{T^2M} \circ TX \circ Y
= \pi_{TM} \circ X \circ \pi_{TM} \circ Y $$
shows that $\pi$ is a group homomorphism. 
Identifying $\phi \in \Diff(M)$ with the origin in $T_\phi(\Diff(M))$, we get 
$$ X \cdot \phi 
= \pi_{T^2M} \circ TX \circ \phi
= X \circ \phi \quad \hbox{ and } \quad \phi \cdot X 
= \pi_{T^2M} \circ T\phi \circ X 
= T\phi \circ X. $$
In particular, this leads to the formula 
\begin{equation}
  \label{eq:adjoint-vect}
  \Ad(\phi)X = T\phi \circ X \circ \phi^{-1}
\end{equation}
for the adjoint action of 
$\Diff(M)$ on $T_0(\Diff(M))= {\cal V}(M)$. 
\end{ex}

\begin{rem} The manifold structure on diffeomorphism groups 
can be made particularly explicit for the group 
$G :=\Diff_c(M)$ of compactly supported diffeomorphisms of an 
open domain $M \subeq \R^n$. Here we call a diffeomorphism 
$\phi \in \Diff(M)$ compactly supported if the set 
$\{ x \in M \: \phi(x) \not=x\}$ is relatively compact in $M$. 
The manifold structure of this group is modeled on the Lie algra 
$\g := \cV_c(M)$ of compactly supported smooth vector fields. 
Now 
\[ U  := \{ X \in \cV_c(M) \: \id_M + X \in \Diff(M) \} \] 
turns out to be an open subset of $\cV_c(M)$, 
so that we obtain a chart 
\[ \phi \:  U \to \Diff_c(M), \quad X \mapsto \id_M + X.\] 
In this chart the multiplication has the particularly simple form 
\begin{align*}
X * Y &= \phi^{-1}(\phi(X)\circ \phi(Y)) 
= (\id_M + X) \circ (\id_M + Y) - \id_M \\
&= X \circ (\id_M + Y) + Y = X + Y + dX\cdot Y + \cdots,
\end{align*}
so that the Lie bracket on $\g$ is given by 
\[ [X,Y] = d X \cdot Y - dY \cdot X, \quad 
[X,Y](m) =dX(m)Y(m) - dY(m)X(m), m \in M \] 
(cf.\ Remark~\ref{rem:brack-taylor}).
The special case $M=\R^n$ was treated in \cite{Gl05d}.
The discussion of the general case is similar;
compare \cite{GN17} for an analogous discussion
of the Lie group $\Diff_{\partial M}(M)$
of all $C^\infty$-diffeomorphisms of a compact convex subset
$M\sub \R^n$ fixing the boundary pointwise,
or~\cite{Gl23} for the diffeomorphism group of a convex polytope. 
\end{rem}

\section{The Maurer--Cartan form} \mlabel{sec:3.4} 

In this section, we introduce some tools that play a key role throughout this 
book: the Maurer--Cartan form and logarithmic derivatives. We have already 
seen that the tangent bundle $T(G)$ of a Lie group $G$ carries the structure 
of a Lie group and that the zero section $\sigma \: G\to T(G)$ and the 
bundle projection $\pi \: T(G) \to G$ are homomorphisms 
of Lie groups (Example~\ref{ex:tangentgrp}). 
This implies that the multiplication map 
$$ G \times \L(G) \to T(G), \quad (g,x) \mapsto gx $$
is a diffeomorphism whose inverse can be written as a pair of maps 
$(\pi,\kappa_G)$, where $\kappa_G \: T(G) \to \L(G)$ is an $\L(G)$-valued 
$1$-form, called the {\it (left) Maurer--Cartan form}. \index{Maurer--Cartan form} 
Composing with the 
Maurer--Cartan form, we associate to each smooth function $f \: M \to G$ 
its logarithmic derivative 
$$\delta(f) := \kappa_G \circ T(f) = f^* \kappa_G \: TM \to \L(G) $$  
in $\Omega^1(M,\L(G))$. 

The importance of these concepts is due to the fact that they permit us 
to develop some basic tools for a calculus of group-valued smooth functions 
on manifolds. In particular, there is a product and a quotient rule from 
which we 
derive powerful uniqueness statements, such as: if $M$ is connected, 
then 
$f$ is uniquely determined by $\delta(f)$ and its value in a single point of $M$. 
This in turn has interesting applications to Lie groups, namely that a homomorphism 
$\phi \: G \to H$ of connected Lie groups is uniquely determined by $\L(\phi)$. 

Here we only deal with the ``differentiation'' aspects of these concepts. 
The more 
refined ``integration'' aspects are dealt with in Chapter~\ref{ch:4} on 
regular Lie groups. We postpone these issues 
 because the integration of Lie algebra-valued $1$-forms requires 
some extra properties of the target group which are naturally encoded in the 
regularity concept.

\subsection{Logarithmic derivatives} \mlabel{subsec3.2.1}

\begin{defn} \mlabel{def:c.11} Restricting the multiplication of the tangent bundle 
$T(G)$ 
to the subset $G \times T(G)$, where we identify $G$ with the zero section,   
we obtain in particular a smooth action of $G$ on $T(G)$ which we simply write 
$(g,v) \mapsto gv$. 

The {\it (left) Maurer--Cartan form} \index{Maurer--Cartan form, $\kappa_G$}
$\kappa_G \in \Omega^1(G,\L(G))$ 
of the Lie group $G$ is defined by 
$$ \kappa_G(v) := \pi(v)^{-1}v, $$
where $\pi \: T(G) \to G$ is the bundle projection. This means that 
$\kappa_G(gx) = x$ for $g \in G$ and $x \in \L(G)\cong T_\be(G)$. The smoothness of this 
$1$-form follows from the smoothness of the action of $G$ on $T(G)$. 
For left invariant vector fields $x_l(g) = gx$, we 
obtain constant functions $\kappa_G(x_l(g)) = x$. 

Let $M$ be a smooth manifold and $G$ a Lie group with
Lie algebra $\L(G)$. 
For an element $f \in C^\infty(M,G)$ we define a smooth $\L(G)$-valued 
$1$-form $\delta(f) \in \Omega^1(M,\g)$ by 
$$ \delta(f) = \kappa_G \circ T(f) = f^*\kappa_G. $$
This means that for each $m \in M$ and $v \in T_m(M)$ we have 
$$ \delta(f)_m(v) := f(m)^{-1} T_m(f)(v), $$
which is written symbolically as 
$$ \delta(f) = f^{-1}.d f. $$
This $1$-form is called the 
{\it (left) logarithmic derivative of $f$} 
\index{logarithmic derivative (left,right)} 
Note that $\kappa_G = \delta(\id_G)$; sometimes this is written as 
``$\kappa_G = g^{-1}.d g$,'' but we won't use this highly ambiguous 
notation. 

We define the {\it right logarithmic derivative} 
by $\delta^r(f) = d f.f^{-1} = f^*\kappa_G^r$ for 
the {\it right Maurer--Cartan form} 
\index{Maurer--Cartan form!right} 
$\kappa_G^r(v) = v\pi(v)^{-1}$, where the 
product refers to the group structure of $T(G)$. 
\end{defn}

\begin{exs} (a) If $G = (E,+)$ is the additive group of a locally convex 
space $E$, then 
we identify $T(E)$ in the usual way with $E \times E$, and in these coordinates we 
have 
$$ \kappa_E(x,v) = v. $$
For any smooth function $f \: M \to E$, the logarithmic derivative 
$\delta(f)$ therefore coincides with the differential $d f$ of $f$. 

(b) If $G = \cA^\times$ is the unit group of a continuous inverse algebra and we identify 
$T(\cA^\times) \cong \cA^\times \times \cA$, then 
$$ \kappa_{\cA^\times}(g,v) = g^{-1}v. $$
For any smooth function $f \: M \to \cA^\times$ we therefore have 
$\delta(f) = f^{-1}.d f$ in the pointwise sense. 
\end{exs}

\begin{lem} {\rm(Quotient and Product Rule)} \mlabel{lem:c.12} 
For $f,g \in C^\infty(M,G)$,  
the following assertions hold: 
\begin{enumerate}
\item[\rm(i)] The map $f^{-1} \: M \to G, m \mapsto f(m)^{-1}$ is smooth with 
  \begin{eqnarray*}
    \delta(f^{-1}) &=& - \Ad(f)\delta(f) = - \delta^r(f).
   \end{eqnarray*}
\item[\rm(ii)] We have the following product and quotient rules: 
$$ \delta(fg) 
= \delta(g) + \Ad(g^{-1})\delta(f), \qquad 
\delta^r(fg) = \delta^r(f) + \Ad(f)\delta^r(g), $$
and 
$$ \delta(fg^{-1}) 
= \Ad(g)(\delta(f) - \delta(g)). $$
\end{enumerate}
\end{lem}

\begin{prf} Although the lemma follows from 
Lemma~\ref{lem:e.1.2} in Appendix~E, 
we give a direct proof using less conceptual background. 

Writing $fg = m_G \circ (f,g)$, we obtain from 
$$ T_{(g,h)}(m_G)(v,w) = vh + gw $$ 
(Proposition~\ref{prop:tangentgrp}) the relation 
$$ T(fg) = T(m_G) \circ (T(f), T(g)) = T(f)\cdot g + f \cdot T(g) 
\: M \to T(G), $$
which immediately leads to 
$$ \delta(fg) 
= (fg)^{-1}\cdot (T(f)\cdot g + f \cdot T(g)) 
= g^{-1}\cdot \delta(f)\cdot g + \delta(g) 
= \Ad(g)^{-1}.\delta(f) + \delta(g). $$
For the right logarithmic derivative, we similarly obtain 
$$ \delta^r(fg) 
= (T(f)\cdot g + f\cdot T(g))\cdot (fg)^{-1}  
= \delta^r(f) + \Ad(f).\delta^r(g). $$

Next we note that 
$$ \delta^r(f) = T(f)\cdot f^{-1} = f\cdot (f^{-1}\cdot T(f))\cdot f^{-1} 
= \Ad(f).\delta(f). $$
From the Product Rule and $\delta(ff^{-1}) = 0$, the formula for 
$\delta(f^{-1})$ follows, and by combining this with the Product Rule, 
we get the Quotient Rule. 
\end{prf}

\begin{lem} {\rm(Uniqueness Lemma)} \mlabel{lem:c.12b} 
Suppose that $M$ is connected. 
For $f_1, f_2 \in C^\infty(M,G)$, the relation $\delta(f_1)= \delta(f_2)$ 
is equivalent to the
existence of a $g \in G$ with $f_2 = \lambda_g \circ f_1$. 
In particular 
$$ \delta(f_1) = \delta(f_2) \quad \hbox{ and } \quad 
f_1(m_0) = f_2(m_0) \quad \hbox{ for some } \quad m_0 \in M $$
imply $f_1 = f_2$. 
\end{lem}

\begin{prf} First we note that 
$\delta(\lambda_g \circ f) = f^*\lambda_g^*\kappa_G 
= f^*\kappa_G = \delta(f)$ for any smooth function 
$f : M \to G$. 

If $\delta(f_1) = \delta(f_2)$, then the Quotient Rule implies that 
$\delta(f_1 f_2^{-1}) = 0$, so that $f_1 f_2^{-1}$ is locally constant 
(Lemma~\ref{locconst}). 
\end{prf}

In the following, we define for a manifold $M$ and a locally convex Lie
algebra $\g$ the bracket 
$$ \Omega^1(M,\g) \times \Omega^1(M,\g) \to \Omega^2(M,\g), \quad 
(\alpha, \beta) \mapsto [\alpha,\beta] $$
by 
$$[\alpha,\beta]_p(v,w) 
:= [\alpha_p(v), \beta_p(w)] - [\alpha_p(w), \beta_p(v)] 
\quad \mbox{ for } \quad  v,w \in T_p(M). $$
Note that $[\alpha,\beta]= [\beta,\alpha]$. 

\begin{lem} \mlabel{lem:MC} 
For each $f \in C^\infty(M,G)$, the $1$-form 
$\alpha :=  \delta(f)$ satisfies the Maurer--Cartan equation 
\[ d \alpha + \shalf[\alpha,\alpha] =0. \tag{MC} \]
\end{lem}

\begin{prf} First we show that $\kappa_G = \delta(\id_G)$ satisfies the 
Maurer--Cartan equation. It suffices to evaluate $d \alpha$ on 
left invariant vector fields $x_l, y_l$, where $x,y \in \g$. 
Since $\kappa_G(x_l) = x$ is constant, we have 
\begin{eqnarray*}
d\kappa_G(x_l, y_l) 
&=& x_l\kappa_G(y_l) - y_l\kappa_G(x_l) - \kappa_G([x_l, y_l])
= -\kappa_G([x,y]_l) = -[x,y] \\
&=& -\shalf[\kappa_G, \kappa_G](x_l, y_l). 
\end{eqnarray*}
Therefore $\alpha = f^*\kappa_G$ satisfies 
$$ d \alpha = f^* d\kappa_G 
= - \shalf f^*[\kappa_G, \kappa_G] 
= - \shalf [f^*\kappa_G, f^*\kappa_G] 
= - \shalf [\alpha, \alpha], $$
which is the Maurer--Cartan equation. 
\end{prf} 

\subsection{Lie group homomorphisms} 

\begin{lem} \mlabel{lem:c.13}  Let $\phi \: G \to H$ be a morphism of Lie groups. 
Then 
$$ \delta(\phi) = \phi^*\kappa_H = \L(\phi) \circ \kappa_G,  $$
and for each smooth map $f \: M \to G$ we have 
$$ \delta(\phi \circ f) = \L(\phi) \circ \delta(f). $$
\end{lem}

\begin{prf} For $g \in G$ we have $\phi \circ \lambda_g =
\lambda_{\phi(g)} \circ \phi$, which implies 
for $g \in G$ and $x \in \L(G)$ that 
$$ \delta(\phi)(gx) 
= \phi(g)^{-1}T_g(\phi)(gx) = T_\be(\phi)x = \L(\phi)x. $$
This means that $\delta(\phi) = \L(\phi) \circ \kappa_G$. 

The second assertion follows from 
\[(\phi \circ f)^*\kappa_G 
= f^*\phi^*\kappa_G 
= f^*(\L(\phi) \circ \kappa_G) \\
= \L(\phi) \circ (f^*\kappa_G) 
= \L(\phi) \circ \delta(f). \qedhere\]
\end{prf}

\begin{prop} \mlabel{propc.14} 
If $G$ is connected and $\phi_1, \phi_2 : G \to H$ are
morphisms of Lie groups with $\L(\phi_1) = \L(\phi_2)$, then 
$\phi_1 = \phi_2$. 
\end{prop}

\begin{prf} In view of Lemma~\ref{lem:c.13}, 
$\delta(\phi_1) 
=  \L(\phi_1) \circ \kappa_G
=  \L(\phi_2) \circ \kappa_G = \delta(\phi_2),$
so that the uniqueness assertion follows from the Uniqueness 
Lemma~\ref{lem:c.12b}. 
\end{prf}

The preceding observation has some interesting consequences: 

\begin{cor} \mlabel{cor:3.2.8} 
If $G$ is a connected Lie group, then 
$\ker \Ad = Z(G)$. 
\end{cor}

\begin{prf} Let $c_g(x) = gxg^{-1}$. 
In view of Proposition~\ref{propc.14}, for $g \in G$ the conditions $c_g =
\id_G$ and $\L(c_g) = \Ad(g) = \id_{\L(G)}$ are
equivalent. This implies the assertion. 
\end{prf}

\begin{prop} \mlabel{prop:ab} 
A connected Lie group $G$ is abelian if and only 
if its Lie algebra is abelian. 
\end{prop}

\begin{prf} That the Lie algebra of an abelian Lie group is abelian has already 
been observed in Example~\ref{ex:abgrp}(a). 
Suppose, conversely, that $\L(G)$ is abelian. 

In view of the preceding corollary, we have to show that $\Ad(g) = \be$ 
holds for each $g \in G$. 
Let $x \in \L(G)$ and consider a smooth curve 
\break $\gamma \: [0,1] \to G$
 with $\gamma(0) = \be$ and $\gamma(1) = g$. For 
$\eta(t) := \Ad(\gamma(t))x$ we then have by Proposition~\ref{prop:der-Ad} 
$$ \eta'(t) 
= T(\sigma_{\Ad})(\gamma(t)\delta(\gamma)_t, 0) 
= \Ad(\gamma(t))[\delta(\gamma)_t, x] = 0, $$
so that $\eta$ is constant. This shows that 
$\Ad(g)x = \eta(1) = \eta(0) = x$. 
\end{prf}

In Section~\ref{sec:3.6} below, 
we shall generalize the preceding result 
to nilpotent and solvable Lie groups. 
The following lemma links the multiplicativity of a map between 
Lie groups to the Maurer--Cartan form. 

\begin{lem} \mlabel{lem:c.15} 
Let $U$ be an open $\be$-neighborhood of the Lie group $G$,
$H$ a Lie group, $\psi \: \L(G) \to \L(H)$ a continuous morphism of
topological Lie algebras, 
and $\phi \: U \to H$ a smooth map with 
$$ \phi(\be_G) = \be_H \quad \hbox { and } 
\quad \phi^*\kappa_H = \psi \circ \kappa_G. $$
Let $V \subeq U$ be an open connected symmetric
$\be$-neighborhood with $V\cdot V \subeq U$. Then 
$$ \phi(xy) = \phi(x) \phi(y) \quad \hbox{ for all} \quad x,y \in V. $$
\end{lem}

\begin{prf} Let $x \in V$. In view of $V \cdot V \subeq U$, we obtain a smooth map 
$$ f = \lambda_{\phi(x)^{-1}} \circ \phi \circ \lambda_x 
\: V \to H, \quad y \mapsto \phi(x)^{-1} \phi(xy). $$
Then $f(\be_G) = \be_H$ and 
\begin{eqnarray*}
\delta(f)
&=& \delta(\lambda_{\phi(x)^{-1}} \circ \phi \circ \lambda_x) 
= \delta(\phi \circ \lambda_x) 
= \lambda_x^* \phi^* \kappa_H \\ 
&=& \lambda_x^* (\psi \circ \kappa_G) 
= \psi \circ (\lambda_x^*\kappa_G) 
= \psi \circ \kappa_G = \delta(\phi).  
\end{eqnarray*}
Since $V$ is connected and $f(\be_G) = \phi(\be_G)$, we obtain 
$f = \phi\res_V$ from the Uniqueness Lemma~\ref{lem:c.12b}. This means that 
$\phi(xy) = \phi(x)\phi(y)$ for $x,y~\in~V$. 
\end{prf}

\begin{prop} \mlabel{prop:homocrit} 
Let $G$ and $H$ be Lie groups. Assume that $G$ is connected and that 
$\psi \: \L(G) \to \L(H)$ is a morphism of topological Lie algebras. 
If $\phi \: G \to H$ is a smooth map with 
\begin{equation}
  \label{eq:3.2.2}
\phi(\be_G) = \be_H \quad \hbox { and } 
\quad \phi^*\kappa_H = \psi \circ \kappa_G, 
\end{equation}
then $\phi$ is a homomorphism of Lie groups. 
\end{prop}

\begin{prf} We apply the preceding lemma with $V = U = G$. 
\end{prf}

\begin{thm} \mlabel{thm:c.16} 
Let $G$ be a $1$-connected Lie group 
and $H$ a Lie group. Further,  let $\psi \: \L(G) \to \L(H)$ be a
morphism of locally convex Lie algebras. Suppose that 
there exists an open $\be$-neighborhood in $G$ 
and a smooth map $\phi \: U \to H$ with 
\begin{equation}
  \label{eq:3.2.1}
\phi(\be_G) = \be_H \quad \hbox { and } 
\quad \phi^*\kappa_H = \psi \circ \kappa_G\res_U. 
\end{equation}
Then there exists a morphism of Lie groups $\eta \: G \to H$ 
which coincides with $\phi$ on each connected open symmetric
$\be$-neighborhood $V$ with $V\cdot V \subeq U$. In particular 
we have 
$$ \L(\eta) = \psi. $$
\end{thm}

\begin{prf} Let $V \subeq G$ be a connected open symmetric
$\be$-neighborhood with $V\cdot V \subeq U$. Then Lemma~\ref{lem:c.15} implies
$\phi(xy) = \phi(x)\phi(y)$ for $x,y \in V$. Now 
Proposition~\ref{prop:lochomo-ext} yields an extension 
$\eta \: G \to H$ of $\phi\res_V$ which is a group homomorphism. Since
$\phi$ is smooth, the homomorphism $\eta$ is smooth in an identity
neighborhood, hence smooth, i.e., a morphism of Lie groups. 
Its differential satisfies 
$$ \L(\eta) = T_\be(\phi) = (\phi^*\kappa_H)(\be) = \psi \circ \kappa_G(\be) = \psi. $$

If $V'$ is another connected open symmetric $\be$-neighborhood in $G$
with $V'V' \subeq U$, then the extension of $\phi\res_{V'}$ leads to a
Lie group morphism $\eta' \: G \to H$ with $\L(\eta') = \psi$. 
Now the connectedness of $G$ implies $\eta = \eta'$ (Proposition~\ref{propc.14}), so that 
$\eta\res_{V'} = \eta'\res_{V'} = \phi\res_{V'}.$ 
\end{prf}

\begin{small}
\subsection*{Exercises for Section~\ref{sec:3.4}} 

\begin{exer} \mlabel{exer:3.4.1} Let $G$, $H$ be Lie groups, 
asume that $G$ is connected, and let 
$\alpha \: G \to \Aut(H), g \mapsto \alpha_g$ 
be a homomorphism defining a smooth action
of $G$ on $H$. For a smooth map $f \: G \to H$, the following are equivalent:
\begin{enumerate}
\item[\rm(1)] $f$ is a {\it (left) crossed homomorphism}, \index{crossed homomorphism} 
$$ f(xy) = f(x)\cdot \alpha_x(f(y)) \quad \mbox{ for } \quad x,y \in G. $$
\item[\rm(2)] $f(\be) = \be$ and $\delta(f) \in \Omega^1(G,\L(H))$ is
{\it equivariant}, i.e., \index{differential form!equivariant}
$$\lambda_g^*\delta(f) = \L(\alpha_g) \circ \delta(f)$$
for each $g \in G$.
\item[\rm(3)] $(f,\id_G) \: G \to H \rtimes_\alpha G$ is a homomorphism.
\end{enumerate}
Hint: Apply the Uniqueness Lemma to functions of the
form $f \circ \lambda_x$, $x \in G$.
\end{exer}

\begin{exer} \mlabel{exer:3.4.2} Let $G$, $H$ be Lie groups, 
asume that $G$ is connected, and let 
$\alpha \: G \to \Aut(H), g \mapsto \alpha_g$ be a homomorphism defining a smooth action
of $G$ on $H$. For a smooth map $f \: G \to H$, the following are equivalent:
\begin{enumerate}
\item[\rm(1)] $f$ is a {\it right crossed homomorphism}, \index{crossed homomorphism} 
$$ f(xy) = \alpha_y^{-1}(f(x)) \cdot f(y)\quad \mbox{ for } 
\quad x,y \in G. $$
\item[\rm(2)] $f(\be) = \be$ and $\delta(f) \in \Omega^1(G,\L(H))$ is
{\it right equivariant}, i.e.,
$$\rho_g^*\delta(f) = \L(\alpha_g^{-1}) \circ \delta(f)$$
for each $g \in G$.
\item[\rm(3)] $(\id_G, f) \: G \to G \ltimes_\alpha H$ is a homomorphism.
\end{enumerate}
Hint: Apply the Uniqueness Lemma to functions of the
form $f \circ \rho_y$, $y \in G$.
\end{exer}

\begin{exer} \mlabel{exer:local-one-para} 
Let $G$ be a Lie group, $\eps > 0$ and 
$\gamma \: [0,\eps] \to G$ be a $C^1$-curve with 
$\delta(\gamma)_t = x$ for $0 \leq t \leq \eps$. Show that: 
\begin{description}
\item[\rm(a)] $\gamma(t+s) = \gamma(t) \gamma(s)$ for $0 \leq t,s, t + s \leq \eps$. 
\item[\rm(b)] $\gamma(t) := \lim_{n \to \infty} \gamma(t/n)^n$ for $t > 0$ 
defines a one-parameter semigroup of $G$, i.e., 
$\gamma(t+s) = \gamma(t)\gamma(s)$ for $t,s \geq 0$. 
\item[\rm(c)] $\gamma(t) := \gamma(-t)$ for $t \leq 0$ defines a smooth 
one-parameter group with $\gamma'(0) = x$. 
\end{description}
\end{exer}

\end{small}

\section{The exponential function of a Lie group} \mlabel{sec:3.5} 

The exponential function of a Lie group is the key tool in the 
finite-dimensional Lie theory and also for Banach--Lie groups. 
For locally convex Lie groups, exponential functions are much less 
well-behaved. In general, they are not 
local diffeomorphisms, so that there are many Lie groups without 
natural charts. This is a serious difficulty one has to face 
in infinite-dimensional Lie theory. In Chapter~\ref{ch:5}, we shall study 
the locally exponential Lie groups 
for which this pathology does not occur. 
In the present section, we take a first look at the exponential function. 
In particular, we derive an explicit formula for its 
logarithmic derivative. 

\subsection{Some examples of exponential functions} 

\begin{defn} \mlabel{def:efunc}
For a Lie group $G$ with Lie algebra $\g = \L(G)$, we call a 
smooth function $\exp_G \: \g \to G$ an 
{\it exponential function for $G$} \index{exponential function, of a Lie group}
if, for each $x \in \g$, the curve 
$\gamma_x(t) := \exp_G(tx)$ is a smooth one-parameter group 
with $\gamma_x'(0) = x$. 
\end{defn} 

From the definition we immediately see that any exponential function 
satisfies 
\begin{eqnarray}
  \label{eq:der-ext0} 
T_0(\exp_G) = \id_{\L(G)}. 
\end{eqnarray}

\begin{lem} \mlabel{lem:5.1.2} If $\exp_G \: \L(G) \to G$ 
is a function which is smooth in an identity neighborhood 
and, for each $x \in \g$, the curve 
$\gamma_x(t) := \exp_G(tx)$ is a one-parameter group 
with $\gamma_x'(0) = x$, then $\exp_G$ is an exponential function.
\end{lem}

\begin{prf} Let $x \in \g$ and $U$ be as in Definition \ref{def:5.1.1}. Then there
exists an $n \in \N$ with $\frac{1}{n} x \in U$. Now the function 
$$ n U \to G, \quad y \mapsto \Big(\exp_G\big(\frac{y}{n}\big)\Big)^n = \exp_G y $$
is smooth in a neighborhood of $x$, and therefore $\exp_G$ is smooth. 
\end{prf}

Since any morphism of Lie groups $\gamma \: \R \to G$ is uniquely determined 
by $\L(\gamma) = \gamma'(0)$ (Proposition~\ref{propc.14}), we have: 

\begin{prop} \mlabel{prop:exp-unique} A Lie group $G$ has at most one exponential function. 
\end{prop}

So much about the uniqueness of an exponential function; the existence is 
a more difficult issue (cf.\ Example~\ref{ex:no-expfct} below). 

The relation $T_0(\exp_G) = \id_{\L(G)}$ 
is not as useful in the locally convex context 
as it is in the finite-dimensional or Banach context. For Banach--Lie groups, it 
follows from the Inverse Function Theorem that 
$\exp_G$ restricts to a diffeomorphism of some open $0$-neighborhood in $\g$ to an 
open $\be$-neighborhood in $G$, so that we can use the exponential function to obtain 
charts around~$\be$. We shall see below (Example~\ref{ex:bad-expfct}) 
that this conclusion does not work for 
general Fr\'echet--Lie groups because there is no general Inverse Function 
Theorem and the image of $\exp_G$ need not be an identity neighborhood in $G$.

\begin{rem} \mlabel{rem:one-par} If $\gamma \: \R \to G$ is a smooth homomorphism and 
$x = \gamma'(0) = \L(\gamma)1$, then $\gamma'(t) = \gamma(t)\gamma'(0)$ 
implies that 
$\delta(\gamma)$ is constant (cf.\ Lemma~\ref{lem:c.13}), 
so that $\gamma$ is the unique solution of the initial value problem 
\begin{eqnarray}
  \label{eq:ivp}
\gamma(0) = \be\quad \mbox{ and }\quad \gamma'(t) = \gamma(t)x. 
\end{eqnarray}

Assume, conversely, that a smooth curve $\gamma \: \R \to G$ satisfies 
$\delta(\gamma) = x$ for some $x \in \L(G)$. In terms of differential forms, 
this means that $\gamma^*\kappa_G = \delta(\gamma) = \kappa_\R \cdot x$, where 
$\kappa_\R$ is the Maurer--Cartan form of $\R$. 
This in turn can be interpreted as 
$\gamma^*\kappa_G = \psi \circ \kappa_\R$ for the linear map 
$\psi \:\R \to \L(G), t \mapsto tx$, and Proposition~\ref{prop:homocrit} implies 
that $\gamma$ is a homomorphism of Lie groups. 

We conclude that the existence of an exponential function is equivalent to 
the existence of a solution $\gamma_x$ to (\ref{eq:ivp}) for $x \in \g$ (which is unique 
by Lemma~\ref{lem:c.12b}), such that $\gamma_x(1)$ depends smoothly on $x$. 
\end{rem}

\begin{ex} \mlabel{ex:3.5.5} 
If $G = (E,+)$ is a locally convex space, then 
$\exp_E := \id_E$ is an exponential function because, for each 
$x \in E$, $\gamma_x(t) := tx$ is a smooth one-parameter group. 

If $G = E/\Gamma$ holds for a discrete subgroup of $\Gamma$, 
then the preceding observation implies that the quotient map 
$E \to E/\Gamma$ is an exponential function of~$G$. 
\end{ex} 

\begin{ex} \mlabel{ex:exp-lingrp} If 
$\cA$ is a unital Banach algebra and $\cA^\times$ its unit group, 
then the exponential series 
$$ \exp_\cA(x) := \sum_{n = 0}^\infty \frac{x^n}{n!} $$
converges uniformly on each ball, hence defines an analytic 
(hence smooth) function $\cA \to \cA$. 

As a consequence of the Cauchy Product Formula for absolutely convergent 
series (Exercise~\ref{exer:2.1.3c}), we obtain the relation 
$$ \exp_\cA(x+y) = \exp_\cA(x) \exp_\cA(y) $$
for commuting elements $x,y \in \cA$. This implies in particular  
$$\exp_\cA(x)\exp_\cA(-x) = \be = \exp_\cA(-x)\exp_\cA(x), $$
so that $\im(\exp_\cA) \subeq \cA^\times$, and it also follows that the 
curves $\gamma_x(t) := \exp_\cA(tx)$ are homomorphisms. Hence 
$$ \exp_{\cA^\times} \: \cA \to \cA^\times, \quad x \mapsto 
\sum_{n = 0}^\infty \frac{x^n}{n!} $$
is the exponential function of the Banach--Lie group $\cA^\times$ 
(see Example~\ref{ex:expo-cia} for a generalization to cias). 
\end{ex} 

\begin{ex}\mlabel{ex:mapgrp} 
Let $K$ be a Lie group with Lie  algebra $\fk$,  
$M$ a compact manifold (possibly with boundary),  
and $r \in \N_0\cup \{\infty\}$. 
In Theorem~\ref{thm:mapgro-Lie} we have seen that 
$C^r(M,K)$ carries a Lie group structure compatible with the compact open 
$C^r$-topology and that its Lie algebra is the space 
$C^r(M,\fk)$, endowed with the pointwise bracket. 

Suppose that $\exp_K \: \fk \to K$ is a smooth 
exponential function of $K$. We claim that 
\[ \exp \:   C^r(M,\fk) \to C^r(M,K), \quad \xi \mapsto \exp_K \circ \, \xi \] 
is an exponential function of $C^r(M,K)$. Clearly, for every 
$\xi \in C^r(M,\fk)$, the curve 
$\gamma_\xi(t) := \exp(t\xi) = \exp_K \circ \, t\xi$ is a 
one-parameter group of $C^r(M,K)$. It remains to show that 
$\exp$ is smooth in a $0$-neighborhood of $C^r(M,\fk)$ 
(Lemma~\ref{lem:5.1.2}). 

In the proof of Theorem~\ref{thm:mapgro-Lie} we obtain 
charts of $G := C^r(M,K)$ as follows. 
Let $\phi_K \: U_K \to\fk$ be a chart of $K$ with 
$\phi_K(\be) = 0$. Then the set 
$U_G := \{ f \in G \: f(M) \subeq U_K\}$ is an open subset of~$G$ 
and $\phi_G \: U_G \to \g, f \mapsto \phi_K \circ f$ 
defines a chart $(\phi_G, U_G)$ of $G$. 

Clearly, $V := \{\xi \in C^r(M,\fk) \:  \xi(M) \subeq \exp_K^{-1}(U_K)\}$ 
is an open $0$-neighborhood in $C^r(M,\fk)$ mapped by 
$\exp$ into $U_G$. Further, 
$\phi_G(\exp\xi) = \phi_K \circ \exp_K \circ \, \xi$, 
so that the smoothness of $\phi_G \circ \exp\res_V$ follows from
Proposition~\ref{pushf-cp}.
\end{ex}

\nin {\bf Open Problem:} 
Presently, it is an important open question if each Lie group modeled 
on a Mackey complete space has an exponential function. The following 
example shows that, without the completeness hypothesis, it is easy to 
find examples of Lie groups without an exponential function. This is 
not surprising because any method for solving ordinary differential 
equations is based on completeness assumptions. 

\begin{ex} \mlabel{ex:no-expfct} (A Lie group without exponential function)  
Let $\cA$ denote the subalgebra of all rational functions in the 
Banach algebra  $C([0,1],\C)$,  
endowed with the sup norm $\|f\| := \sup_{0 \leq t \leq 1} |f(t)|$.  
Since an element $f \in \cA$ is invertible if and only if 
$\im(f) \subeq \C^\times$, the unit group $\cA^\times$ is open, and since 
$C([0,1],\C)$ is a Banach algebra, $\cA$ is a continuous inverse algebra, which in turn implies that 
$\cA^\times$ is a Lie group. 

If $\gamma \: \R \to \cA^\times$ is a smooth one-parameter group, then 
it also is a one-parameter group of $C([0,1],\C^\times)$, hence of the form 
$\gamma(t) = e^{tf}$ for some continuous function $f = \gamma'(0) \in \cA$ 
(Example~\ref{ex:exp-lingrp}). 
For $f \in \cA$, the function $e^f$ is a rational function if and only 
if $f$ is constant. This implies that the Lie group $\cA^\times$ has no 
exponential function, resp., only the constant functions 
generate smooth one-parameter groups. 
\end{ex}
  
\begin{ex} \mlabel{ex:expo-cia} Let $G := \cA^\times$ be the unit group of a 
Mackey complete continuous inverse algebra. Then the fact that $\cA^\times$ is open implies that, for each 
$a \in \cA$, the spectrum $\Spec(a)$ is a compact subset (which also is non-empty), 
and it is shown in \cite{Gl02b} (see also Chapter~\ref{ch:lingrp}) that the 
holomorphic functional calculus works as for Banach algebras. We thus obtain 
an analytic exponential function 
$$ \exp_\cA \: \cA \to \cA^\times, 
\quad x \mapsto \frac{1}{2\pi i} \oint_\Gamma e^\zeta \cdot (\zeta\1 - x)^{-1}\, 
d\zeta, $$
where $\Gamma$ is any piecewise smooth contour around $\Spec(x)$. 
Here we consider $\cA$ as a closed subalgebra of its complexification 
$\cA_\C = \cA \otimes_\R \C$, which also is a continuous inverse algebra 
(Corollary~\ref{cor:8.1.5}). 
From 
$$\exp_\cA(tx) = \frac{1}{2\pi i} \oint_\Gamma e^{t\zeta} 
\cdot (\zeta\1 - x)^{-1}\, d\zeta, $$
the fact that 
$e_t(z) := e^{tz}$ defines a continuous homomorphism $\R\to {\cal O}(\C)$ 
(the Fr\'echet algebra of holomorphic functions on $\C$), 
and the fact that the holomorphic functional calculus yields a continuous algebra homomorphism 
$$ {\cal O}(\C) \to \cA_\C, \quad f \mapsto 
 \frac{1}{2\pi i} \oint_\Gamma f(\zeta) \cdot(\zeta\1 - x)^{-1}\, d\zeta, $$
it follows that $\gamma_x(t) := \exp_\cA(tx)$ is a one-parameter group, so 
that $\exp_\cA$ is indeed an exponential function of the Lie group $\cA^\times$. 
In addition, the convergence of the exponential series in the Fr\'echet algebra $\cO(\C)$ 
implies that 
\[ \exp_\cA(x) = \sum_{n = 0}^\infty \frac{x^n}{n!} \] 
also holds in $\cA$. 
\end{ex}

The following proposition provides the exponential function for an interesting 
class of examples that is still quite simple, but it also displays already many 
of the problems arising for the exponential function of non-Banach--Lie groups. 

\begin{prop} \mlabel{prop:semdir-exp} 
Let $E$ be a Mackey complete space, $K$ be a Lie group with 
exponential function $\exp_K$ and let 
$\alpha \: K \to \GL(E)$ define a smooth linear action of $K$ on~$E$. 
Then the semi-direct product Lie group $G := E \rtimes_\alpha K$ has an 
exponential function given by 
\[  \exp_G(v,x) := (\beta(x)v, \exp_K(x)) \quad \mbox{ with } \quad 
\beta(x)v = \int_0^1 \alpha_{\exp_K(sx)}v\, ds. \] 
\end{prop} 

\begin{prf} Considering the one-parameter group 
$t \mapsto \exp_K(tx)$ of $K$, we see that it suffices to verify 
the assertion in the special case where $K = \R$. 

First we observe that a curve 
$\gamma(t) = (\xi(t),\eta(t))$ in $G$ is a smooth homomorphism 
if and only if there exists an $a \in \R$ with $\eta(t) = at$ and 
$$ \xi(t_1 + t_2) = \xi(t_1) + \alpha_{at_1}\xi(t_2) 
\quad \mbox{ for } \quad t_1, t_2 \in \R,  $$
i.e., $\xi \: \R \to E$ is a crossed homomorphism 
for the action of $\R$ on $E$ by $t.v := \alpha_{at}v$.  
In view of Exercise~\ref{exer:3.4.1}, this is equivalent to 
the equivariance of the derivative $\xi'\: \R \to E$, which is 
equivalent to 
$$ \xi'(t) = \alpha_{at}\xi'(0). $$
This is in turn is equivalent to 
\[ \xi(t) = \int_0^t \alpha_{as}v\, ds \quad \mbox{ for } \quad v = \xi'(0). \] 
We thus obtain a description of the smooth one-parameter groups 
of $G$ by 
$$ \gamma(t) = \Big(\int_0^t \alpha_{as}v\, ds, at\Big), \quad 
\mbox{ where } \quad 
\gamma'(0) = (v, a). $$
This in turn implies that the map 
$$\exp_G(v,a) := \Big(\int_0^1 \alpha_{as}v\, ds, a\Big)$$
is an exponential map of $G$ because its smoothness follows from 
Proposition~\ref{prop:cartes-closed}. 
\end{prf}

Note that, for $t \not=0$, we can also write 
\[\beta(t) = \frac{1}{t} \int_0^t \alpha_s \, ds,\] 
so that $\beta(t)$ is the mean value of $\alpha$ 
on the interval between $0$ and $t$.

\begin{ex} \mlabel{ex:bad-expfct} Let $G =E \rtimes_\alpha \R$ be as 
in the preceding proposition. 
From the explicit formula for $\exp_G$ 
it is clear that, for $t \not=0$,  $(w,t) \in \im(\exp_G)$ is equivalent to 
$w \in \im(\beta(t))$. We conclude that 
$\exp_G$ is injective on some $0$-neighborhood 
if and only if $\beta(t)$ is injective for $t$ close to $0$, 
and it is surjective onto some $\be$-neighborhood in $G$ if 
and only if $\beta(t)$ is surjective for $t$ close to $0$ 
(cf.\ Problem~\ref{prob:5.1}). 

Note that the eigenvector equation 
$Dv = \lambda v$ for $t\lambda \not=0$ implies that 
\[  \beta(t)v = \int_0^1 e^{st\lambda}v\, ds  = \frac{e^{t\lambda} - 1}{t\lambda}v, \]
so that $\beta(t)v = 0$ is equivalent to $t\lambda \in 2 \pi i \Z \setminus \{0\}$. 

(a) For the Fr\'echet space $E = \C^\N$ 
and the diagonal operator $D$ given by 
$D(z_n) = (2\pi in z_n)$, we see that 
$\beta(\frac{1}{n}) e_n = 0$ holds for $e_n = (\delta_{mn})_{m \in \N}$, 
and $e_n \not\in \im\big(\beta(\frac{1}{n})\big)$. 
We conclude that $(e_n, \frac{1}{n})$ is not contained in the image of $\exp_G$, and since 
$(e_n, \frac{1}{n})\to (0,0)$, the identity of $G$, 
$\im(\exp_G)$ does not contain any identity neighborhood of $G$. 
Hence the exponential function of the Fr\'echet--Lie group 
$G = E \rtimes_\alpha \R$ is neither locally injective nor locally surjective in $0$. 

(b) For the Fr\'echet space $E = \R^\N$ 
and the diagonal operator $D$ given by 
$D(z_n) = (n z_n)$, all operators $\beta(t)$ are invertible, and 
the map 
$$\tilde\beta \: \R \times E \times E, \quad 
(t,v) \mapsto (\beta(t)v,\beta(t)^{-1}v)$$ 
is smooth. This implies that $\exp_G \: \g \to G$ is a diffeomorphism. 
\end{ex}

\subsection{Naturality of the exponential function} 

We now return to the general discussion of the exponential function.

\begin{prop} \mlabel{prop:exp-diag} If $\phi \: G \to H$ is a morphism of Lie groups with an exponential 
function, then the following diagram commutes: 
$$ \begin{matrix}
 G & \smapright{\phi} & H \\ 
\mapup{\exp_G} & & \mapup{\exp_H} \\ 
\L(G)& \smapright{\L(\phi)}& \L(G). 
\end{matrix} $$
\end{prop}

\begin{prf} For each $x\in \L(G)$ and $\gamma^G_x(t) := \exp_G(tx)$, the map 
$\phi \circ \gamma_x^G \: \R \to H$ is a smooth group homomorphism 
with $(\phi \circ \gamma_x^G)'(0) = \L(\phi)x$. Since this determines 
$\phi \circ \gamma_x^G$ uniquely, it follows that 
$\phi \circ \gamma_x^G = \gamma_{\L(\phi)x}^H$, so that evaluating in 
$t = 1$ implies $\phi(\exp_G(x)) = \exp_H(\L(\phi)x)$. 
\end{prf}

\begin{prop} \mlabel{prop:exp-cover} 
Let $\phi \: G \to H$ be a covering morphism of Lie groups. 
  \begin{enumerate}
  \item[\rm(i)] If $\exp_G$ is an exponential function of $G$, then 
$\exp_H := \phi \circ \exp_G \circ \L(\phi)^{-1}$ 
is an exponential function of $H$. 
  \item[\rm(ii)] If $H$ has an exponential function $\exp_H$, then 
$G$ has an exponential $\exp_G$ with $\phi \circ \exp_G = \exp_H \circ 
\L(\phi)$. 
  \end{enumerate}
\end{prop}

\begin{prf} (i) is an immediate consequence of the definitions and the fact that 
$\L(\phi)$ is a topological isomorphism. 

(ii) Let $\exp_G \: \L(G) \to G$ be the unique lift 
of the map \break $\exp_H \circ \L(\phi) \: \L(G) \to H$ with 
$\exp_G(0) = \be$. Then $\exp_G$ is smooth (Exercise~\ref{exer:3.3.7}) 
and satisfies 
$$ \phi  \circ \exp_G = \exp_H \circ \L(\phi). $$
Moreover, for each $x \in \L(G)$, the curve 
$\tilde\gamma_x(t) := \exp_G(t x)$ is the unique continuous lift 
of the one-parameter group $\phi(\tilde\gamma_x(t)) = \exp_H(t \L(\phi)x)$ 
mapping $0$ to $\be$. Therefore Proposition~\ref{prop:3.3.8} 
implies that $\tilde\gamma_x$ is a one-parameter group. 
From $\L(\phi)\tilde\gamma_x'(0) = \L(\phi)x$ we derive 
that $\tilde\gamma_x'(0) = x$, and this implies that 
$\exp_G$ is an exponential function of~$G$. 
\end{prf}

\begin{thm} \mlabel{thm:ban-expfct} Each 
Banach--Lie group has a (smooth) exponential function. 
\end{thm}

\begin{prf} {\bf Step 1:} Each left invariant vector field $X$ on $G$ is complete. 

Let $g \in G$ and $\gamma\: I \to G$ be the unique maximal integral curve 
of $X \in {\cal V}(G)^l$ with $\gamma(0) = g$ (cf.\ \cite{La99}). 

For each $h \in G$ we have the relation $T(\lambda_h) \circ X \circ \lambda_h^{-1}= X$, 
which implies that $\eta := \lambda_h \circ \gamma$ also is an integral curve of $X$. Put 
$h = \gamma(s)g^{-1}$ for some $s > 0$. Then 
$$\eta(0)= (\lambda_h \circ \gamma)(0) = h \gamma(0) = hg = \gamma(s),$$ 
and the uniqueness of integral curves implies that 
$\gamma(t+s) = \eta(t)$ for all $t$ in the interval $I \cap (I-s)$. 
In view of the maximality of $I$, it now follows that $I - s \subeq I$, and hence that 
$I - ns \subeq I$ for each $n \in \N$, so that 
the interval $I$ is unbounded from below. 
Applying the same argument to some $s < 0$, we see that 
$I$ is also unbounded from above. Hence $I = \R$, which means that $X$ is complete. 

{\bf Step 2:} We now define the exponential function  
$$ \exp_G \: \L(G) \to G, \quad  \exp_G(x) := \gamma_x(1), $$
where $\gamma_x \: \R \to G$ is the unique maximal integral curve of the left invariant 
vector field $x_l$, satisfying $\gamma_x(0) = \be$. 
This means that $\gamma_x$ is the unique solution 
of the initial value problem 
$$ \gamma(0) = \be, \quad \gamma'(t) = x_l(\gamma(t)) = \gamma(t)x \quad 
\mbox{ for all } \quad t \in \R, $$
hence a one-parameter group (Remark~\ref{rem:one-par}). 

{\bf Step 3:} $\exp_G$ is smooth. 

The map $\Psi \: \L(G) \times G \to T(G), (x,g) \mapsto x_l(g) = gx$ 
is smooth because the multiplication on $T(G)$ is smooth. 
It can be interpreted as a smooth vector field on the product 
manifold $\L(G) \times G$, if we identify the tangent bundle 
$T(\L(G) \times G)$ with $T(\L(G)) \times T(G)$. 
Then the integral curves of $\Psi$ are of the form 
$\gamma(t) = (x, \gamma_G(t))$, where $\gamma_G$ is an integral curve 
of the vector field $x_l$. Therefore Step $1$ implies that 
$\Psi$ is a complete vector field, hence that its global flow 
$$ \R \times \L(G) \times G \to \L(G) \times G, \quad 
(t,x,g) \mapsto (x,g \exp_G(tx)) $$
is smooth
(see Corollary~\ref{pardep-ode-ban-m}). 
Restricting to 
{$\{1\} \times \L(G) \times \{\be\}$}, the smoothness of $\exp_G$ follows. 
\end{prf}

\begin{prop} If $G$ is a Lie group with an exponential function 
whose image is an identity neighborhood, then 
the subgroup $\la \exp_G(\L(G))\ra$ of $G$ generated by 
$\exp_G(\L(G))$ coincides with the identity component $G_0$ of~$G$.
\end{prop}

\begin{prf} Since $\im(\exp_G)$ is an identity neighborhood, 
the subgroup $H := \la \exp_G(\L(G))\ra$ generated 
by the exponential image is a $\be$-neighborhood. That it actually is a subgroup 
follows from $\exp_G(-x) = \exp_G(x)^{-1}$ for each $x \in \L(G)$. 
From Exercise~\ref{exer:3.1.3}(3) we derive that 
$H$ is an open closed subgroup of $G$, hence contains $G_0$. 

On the other hand, $\exp_G$ is continuous, so that it maps the 
connected space $\L(G)$ into the identity component $G_0$ of $G$, 
which leads to $H \subeq G_0$, and hence to equality. 
\end{prf}

\begin{lem} \mlabel{lem:4.1.4} Let $G$ be Lie groups with a smooth 
exponential function. 
  \begin{enumerate} 
  \item[\rm(i)] If $x,y \in \L(G)$ commute, then $\exp_G(x+y)= \exp_G(x)\exp_G(y)$. 
  \item[\rm(ii)] $\exp_G \circ \Ad(g)  = c_g \circ \exp_G$ for each $g \in G$.
  \end{enumerate}
\end{lem}

\begin{prf} (i) If $x, y \in \L(G)$ commute, then we first consider the smooth curve 
$$ \eta \: \R \to \g, \quad t \mapsto \Ad(\gamma_y(t))x. $$
According to Proposition~\ref{prop:der-Ad}, its derivative is given by 
$$ \eta'(t) 
= T_{(\gamma_y(t),x)}(\sigma_{\Ad})(\gamma_y(t)y, 0) 
= \Ad(\gamma_y(t))[y,x] = 0. $$
Hence $\Ad(\gamma_y(t))x = x$ for each $t$, so that 
Lemma~\ref{lem:c.12}(ii) implies  
$$ \delta(\gamma_x \cdot \gamma_y) 
= \delta(\gamma_y) + \Ad(\gamma_y)^{-1}\delta(\gamma_x) 
= y + \Ad(\gamma_y)^{-1}x = y + x, $$
which leads to $\gamma_{x + y} = \gamma_x \cdot \gamma_y$ 
(Remark~\ref{rem:one-par}), and finally to 
$$ \exp_G(x + y) = \gamma_{x+y}(1) = \gamma_x(1)\gamma_y(1) = \exp_G x \exp_G y. $$

(ii) follows from Proposition~\ref{prop:exp-diag} with $H = G$ and $\phi = c_g$. 
\end{prf}

\subsection{Abelian Lie groups with exponential function}

The following proposition characterizes connected 
abelian Lie groups with an exponential function. 
Since all known Lie groups with a 
Mackey complete Lie algebra are regular and 
therefore have an exponential function, it applies 
to all these groups. 

\begin{prop} \mlabel{prop:reg-ab}
Let $A$ be a connected abelian Lie group modeled on a Mackey complete 
space $\fa$. Then $A$ has an exponential function if and only if 
$A \cong \fa/\Gamma_A$  holds for a discrete subgroup $\Gamma_A$ of $\fa$. 
\end{prop}

\begin{prf}  If $\Gamma_A$ is a discrete subgroup of the 
Mackey complete locally convex space $\fa$ and 
$A := \fa/\Gamma_A$ the topological quotient group, then 
$A$ carries a natural Lie group structure (Corollary~\ref{cor:cover-lie}) 
and the quotient map 
$q_A \: \fa \to A$ can be viewed as the exponential function of $A$ 
and also as the universal covering map. In particular, 
$A$ has an  exponential function. 

Suppose, conversely, that $A$ is a connected abelian Lie group modeled on a Mackey complete 
space with an exponential 
function $\exp_A$. First we claim that $\exp_A$ is surjective. In view of 
Lemma \ref{lem:4.1.4}(i),  
$\exp_A \: (\fa,+) \to A$ is a morphism of Lie groups.
Let $a \in A$ and consider a smooth path $\gamma \: [0,1] \to A$ with 
$\gamma(0) = \be$ and $\gamma(1) = a$. Then the logarithmic 
derivative $\xi := \delta(\gamma)$ is a smooth map $[0,1] \to \fa$ and we consider 
the smooth path 
$$ \eta(t) := \exp_A\Big(\int_0^t \xi(s)\, ds\Big) $$
that also satisfies $\delta(\eta) = \xi$ because 
$\delta(\exp_A) = \id_\fa$ (Lemma~\ref{lem:c.13}). Here we have used the Mackey 
completeness of $\fa$ to ensure the existence of the Riemann integral 
of the smooth curve $\xi$. Now $\eta(0) = \gamma(0) = \be$ implies that 
$$ a = \gamma(1) = \eta(1) 
= \exp_A\Big(\int_0^1 \xi(s)\, ds\Big) \in \im(\exp_A) $$
(Lemma~\ref{lem:c.12b}). 

Let $q_A \: \tilde A \to A$ denote a universal covering homomorphism with 
$\L(q_A) = \id_\fa$.  Then Proposition~\ref{prop:exp-cover} implies that 
$\tilde A$ has an  exponential function 
$\exp_{\tilde A} \: \fa \to \tilde A$ with $\exp_A = q_A \circ \exp_{\tilde A}$. 
Since $\tilde A$ is simply connected, the Lie algebra homomorphism 
$\id_\fa \: \fa\to \fa$ integrates to a Lie group homomorphism 
$L \: \tilde A \to \fa$ with $\L(L) = \id_\fa$. In fact, 
the closedness of the Maurer--Cartan form 
of $\tilde A$ implies the existence of a unique smooth function 
$L \: \tilde A \to \fa$ with $L(\be) = 0$ and 
$d L = \kappa_{\tilde A}$, so that Proposition~\ref{prop:homocrit} implies
that $L$ is a homomorphism. 
We now have 
$$ L \circ \exp_{\tilde A} = \exp_\fa \circ \L(f) = \id_\fa \circ \id_\fa = \id_\fa, $$
and hence $\exp_{\tilde A} \circ L$ restricts to the identity on 
$\im(\exp_{\tilde A}) = \fa$, which also leads to 
$\exp_{\tilde A} \circ L = \id_{\tilde A}.$
Hence $\tilde A \cong \fa$ as Lie groups, which implies that 
$\exp_A$ is a covering morphism and therefore that 
$\Gamma_A := \ker(\exp_A) \subeq \fa$ is discrete with 
$A \cong \fa/\Gamma_A$. 
\end{prf}

\begin{rem} If we do not know that the exponential function \break 
$\exp_A \: \fa \to A$ of a connected simply connected abelian 
Lie group $A$ is smooth, the proof of the preceding proposition does not work. 
Since $\fa$ is abelian, it follows from 
Lemma \ref{lem:4.1.4}(i) that 
$$ \exp_A \: (\fa,+) \to A $$
is a group homomorphism. 
If $\fa$ is Mackey complete, then the Lie algebra homomorphism 
$\id_\fa \: \fa\to \fa$ integrates to a smooth group homomorphism 
$L \: A \to \fa$ with $\L(L) = \id_\fa$ and $L \circ \exp_{A} = \id_\fa$. 
Hence $\exp_{A} \circ L$ restricts to the identity on 
$\im(\exp_{A})$, but there is no reason for $\exp_{A}$ to be smooth 
or a diffeomorphism. It is interesting that we always get the smooth function 
$L$, which can be viewed as a logarithm function, but that the smoothness of 
$\exp_{A}$ seems to lie deeper. 
\end{rem}

\begin{rem} Proposition~\ref{prop:reg-ab} applies in particular 
to finite-dimensional connected abelian Lie groups $A$, so that 
$A \cong \fa/\Gamma$, where $\Gamma \subeq \fa$ is a discrete 
subgroup. In view of Exercise~\ref{exer:3.5.1}, there exists a 
basis $v_1, \ldots, v_n$ of $\fa$ with 
$\Gamma = \Z v_1 + \cdots + \Z v_k$ for some $k \leq n$. 
We thus obtain the product structure 
$$ A \cong \fa/\Gamma \cong \R^n/\Z^k \cong \T^k \times \R^{n-k}. $$
\end{rem}

\subsection{The derivative of the exponential function} 

To obtain a formula for the (logarithmic) derivative 
of the exponential function, we first observe that 
$T(\exp_G)$ is an exponential function of the tangent 
group $TG$. 

\begin{lem} \mlabel{lem:Texp} If $\exp_G \: \g \to G$ is an exponential 
function, then 
\[ T(\exp_G) \: T\g \cong \g \oplus \g \to TG \] 
is an exponential function of the tangent group $TG$. 
\end{lem}

\begin{prf} For $t \in \R$ we consider the smooth maps 
$E_t \: \g \to G, E_t(x)  = \exp_G(tx).$
They satisfy the relation 
\[ E_t \cdot E_s = E_{t+s} \quad \mbox{ for } \quad t,s \in\R.\]Applying the tangent functor, we obtain the relation 
\begin{equation}
  \label{eq:tet}
T(E_t) \cdot T(E_s) = T(E_{t+s}) 
\quad \mbox{ for } \quad t,s \in\R.
\end{equation}
We conclude that, for $(x,y) \in T\g \cong \g \oplus \g$, 
the smooth curve 
\[ \gamma \: \R \to TG, \quad \gamma(t) 
= T(E_t)(x,y) = T(\exp_G)(tx, ty)\] 
is a one-parameter group. Now 
\begin{align*}
\gamma'(0) 
&= \derat0 T(\exp_G)(tx,ty) 
=  \derat0 \exp_G(tx) + \derat0 T_0(\exp_G)ty \\
&= (x,0) + (0,y) = (x,y)
\end{align*}
implies that $T(\exp_G)$ is an exponential 
function of $T(G)$ (cf.\ Definition~\ref{def:efunc}). 
\end{prf}

\begin{defn} \mlabel{def:kappa-liealg}
Suppose that $G$ is a Lie group with an exponential 
function and $\g = \L(G)$. 
Then the one-parameter group  
$\Ad(\exp_G(tx))$ of $\GL(\g)$ has the infinitesimal generator 
$\ad x$ (Proposition~\ref{prop:der-Ad}), so that we may also write 
\begin{eqnarray}
  \label{eq:2.5.4}
\Ad(\exp_G(tx)) = e^{t\ad x}
\end{eqnarray}
(cf.\ Proposition~\ref{prop:e.2.6}). 
If, in addition, $\g$ is Mackey complete, then the operator-valued integral 
\begin{eqnarray}
  \label{eq:2.5.5}
\kappa_\g(x) := \int_0^1 e^{-t\ad x}\, dt 
\end{eqnarray}
exists in the pointwise sense 
because the curves $t \mapsto e^{-t\ad x}y$ are smooth 
and the integral denotes the linear map 
$y \mapsto  \int_0^1 e^{-t\ad x}y\, dt.$ 
Note that $\kappa_\g(x)$ 
commutes with $\ad x$, so that 
\begin{equation}
  \label{eq:kappa-comm}
\kappa_\g(x) \circ \ad x 
= \ad x \circ \kappa_\g(x) 
= \int_0^1 e^{-t\ad x}\ad x\, dt 
= \1- e^{-\ad x}.   
\end{equation}
We also observe that 
\begin{equation}\label{eq:erel2} 
e^{\ad x} \kappa_\g(x) 
= \int_0^1 e^{(1-t)\ad x}\, dt 
= \int_0^1 e^{t\ad x}\, dt = \kappa_\g(-x). 
\end{equation}
\end{defn}

\begin{thm} \mlabel{thm:exp-logder} 
If $\g$ is Mackey complete and 
$\exp_G \: \g \to G$ is an exponential function, then 
\begin{eqnarray}
  \label{eq:2.5.2}
\delta(\exp_G)(x) = \kappa_\g(x) = \int_0^1 e^{-t\ad x}\, dt. 
\end{eqnarray}
\end{thm}

\begin{prf} {\bf 1st proof:} 
First we recall that  
$T(G) \cong (\g,+) \rtimes_{\Ad} G$ is a semidirect product 
of the abelian Lie group $(\g,+)$ and $G$, where $G$ acts 
on $\g$ by the adjoint representation 
(Example~\ref{ex:tangentgrp}). Provided 
$\g$ is Mackey complete, 
Proposition~\ref{prop:semdir-exp}, thus 
provides the explicit formula 
\[  T(\exp_G)(y,x) = \exp_{T(G)}(y,x) = (\beta(x)y, \exp_G(x))\] with 
\[ \beta(x) 
= \int_0^1 \Ad(\exp_G(sx))\, ds
= \int_0^1 e^{s\ad x}\, ds = \kappa_\g(-x). \] 
We conclude with \eqref{eq:erel2}  that 
\begin{align*}
\delta(\exp_G)_x y 
&= \exp_G(x)^{-1} \cdot T(\exp_G)(x,y) 
= \exp_G(-x) \cdot \kappa_\g(-x)y\cdot \exp_G(x) \\
&= e^{-\ad x}\kappa_\g(-x)y 
= \kappa_\g(x)y.
\end{align*}

\nin{\bf 2nd proof:} 
For $t \in \R$ the functions 
$E_t \: \g \to G, x \mapsto \exp_G(tx)$, 
satisfy $E_{t+s} = E_t E_s$ pointwise on $\g$. The 
Product Rule (Lemma~\ref{lem:c.12}) implies that 
$$ \delta(E_{t+s}) = \delta(E_s) + \Ad(E_s)^{-1} 
\delta(E_t). $$
For the smooth curve 
\[ \psi \: \R \to \g, \quad 
\psi(t) := \delta(\exp_G)_{tx}(ty) = \delta(E_t)_x y \] 
we therefore obtain 
\begin{eqnarray*}
\psi(t+s) 
&=& \delta(E_{t+s})_x(y) 
= \delta(E_s)_x(y) + \Ad(E_s)^{-1}\delta(E_t)_x(y) \\
&=& \psi(s) + \Ad(\exp_G(-sx))\psi(t) 
= \psi(s) + e^{-s\ad x}\psi(t). 
\end{eqnarray*}
We have $\psi(0) = 0$ and 
$$\psi'(0) = \lim_{t \to 0} \delta(\exp_G)_{tx}(y) 
= \delta(\exp_G)_0(y) = y, $$ 
so that taking derivatives with respect to $t$ in $0$, 
leads to 
$\psi'(s) = e^{-s\ad x}y.$
Now the assertion follows from 
$\delta(\exp_G)_x(y) = \psi(1) = \int_0^1 \psi'(s)\, ds$. 
\end{prf}

\begin{rem} \mlabel{rem:diffexp-alg} If 
$G$ is the unit group $\cA^\times$ of a Mackey complete continuous inverse algebra 
(cf.\ Example~\ref{ex:expo-cia}), then 
we identify $T(G)$ with $\cA^\times \times \cA$ and note that 
$\ad x = \lambda_x - \rho_x$ and $e^{\ad x}y = e^x ye^{-x}.$
Therefore (\ref{eq:2.5.2}) can be written as 
$$ T_x(\exp_G)y 
= e^x \int_0^1 e^{-tx} y e^{tx}\, dt 
= \int_0^1 e^{(1-t)x} y e^{tx}\, dt. $$
\end{rem}

\begin{rem} \mlabel{rem:ban-dexp} 
If $\g$ is a Banach--Lie algebra, then $\kappa_\g(x)$ can be represented by 
a convergent power series 
$$ \kappa_\g(x) 
= \int_0^1 e^{-t\ad x}\, dt 
= \sum_{k = 0}^\infty \frac{(-1)^k}{(k+1)!} (\ad x)^k. $$
This means that $\kappa_\g(x) = f(\ad x)$ holds for the entire function 
\[  f(z) :=  \sum_{k = 0}^\infty \frac{(-1)^k}{(k+1)!} z^k = \frac{1 - e^{-z}}{z}. \]
As $ f^{-1}(0) = 2\pi i \Z \setminus \{0\}$, 
and $\Spec(\kappa_\g(x)) = f(\Spec(\ad x))$
by the Spectral Mapping Theorem (\cite[Thm.~10.28]{Ru91}), 
$$\kappa_\g(x) \in \GL(\g) \quad \Longleftrightarrow \quad 
\Spec(\ad x) \cap 2 \pi i \Z \subeq \{0\}.$$

Part of this observation can be secured in the general case. 
If $\g$ is Mackey complete, 
$\kappa_\g(x)$ is not injective if and only if there exists 
some $n \in \N$ with 
\begin{equation}
  \label{eq:singexp-ev}
 \ker((\ad x)^2 + 4 \pi^2 n^2\1) \not= \{0\}
\end{equation}
(Proposition~\ref{prop:kerexp-locconv}). 
If $\g$ is a complex Lie algebra, this means that some  
$2\pi i n \in 2 \pi i \Z \setminus \{0\}$ is an eigenvalue of $\ad x$.
\end{rem}

A closer inspection of \eqref{eq:2.5.5} leads to the 
following results on the behavior of the exponential function 
(cf.\ \cite{LaTi66} for the finite-dimensional case): 

\begin{prop} \mlabel{prop:exp-noninj} Let $G$ be a Lie group 
with an exponential function and Lie algebra~$\g$. 
Then the following assertions hold for 
$x,y \in \g$: 
\begin{enumerate}
\item[\rm(i)] If $\kappa_\g(x)y = 0$, 
then 
$$ \exp_G(e^{t\ad y}x) = \exp_G(x) \quad \mbox{ for all } \quad t \in \R.$$
\item[\rm(ii)] If $\kappa_\g(x)$ is not injective and $\g$ is Mackey complete, 
then $\exp_G$ is not injective in any neighborhood of $x$. 
\item[\rm(iii)] If $\kappa_\g(x)$ is injective, then 
  \begin{enumerate}
  \item[\rm(a)] $\exp_G(y) = \exp_G(x)$ implies $[x,y]=0$ and $\exp_G(x-y) = \be$. 
  \item[\rm(b)] $\Ad(\exp_G(x)) = \id \Longleftrightarrow x \in \z(\g)$. 
  \item[\rm(c)] $\exp_G(x) = \be$ implies $x \in \z(\g)$. 
  \end{enumerate}
\item[\rm(iv)] Suppose that $0$ is isolated in $\exp_G^{-1}(\be)$. 
Then $x$ is isolated in $\exp_G^{-1}(\exp_G(x))$ 
if and only if $\kappa_\g(x)$ is injective. 
\item[\rm(v)] If $\fa \subeq \g$ is an abelian subalgebra, then 
$\exp_G \res_\fa \: \fa \to G$ is a morphism of Lie groups. 
Its kernel $\Gamma_\fa := \ker(\exp_G\res_\fa)$ is a closed subgroup of $\fa$ 
in which all $C^1$-curves are constant. It intersects each finite-dimensional 
subspace of $\fa$ in a discrete subgroup. 
\end{enumerate}
\end{prop}

\begin{prf} (i) We consider the smooth curve $\gamma \: \R \to G$, given by 
\begin{eqnarray*}
\gamma(t) &:=& \exp_G(e^{t\ad y}x)
= \exp_G(\Ad(\exp_G(ty))x) \\
&=& \exp_G(ty) \exp_G(x) \exp_G(-ty). 
\end{eqnarray*}
We show that $\gamma'(t) = 0$ for each $t \in \R$, which implies 
that $\gamma$ is constant. 
Let $\beta(t) := \Ad(\exp_G(ty))x = e^{t \ad y}x$. Then 
$$ \kappa_\g(\beta(t)) 
= \Ad(\exp_G(ty))\circ \kappa_\g(x)\circ \Ad(\exp_G(ty))^{-1} $$
(Exercise~\ref{exer:3.5.7}) implies that 
\begin{eqnarray*}
\delta(\gamma)_t 
&=& \delta(\exp_G)_{\beta(t)} \beta'(t) 
= \delta(\exp_G)_{\beta(t)}\Ad(\exp_G(ty))[y,x] \\ 
&=& \kappa_\g(\beta(t))\Ad(\exp_G(ty))[y,x]  
= \Ad(\exp_G(ty))\kappa_\g(x)[y,x] \\
&=& -\Ad(\exp_G(ty)) \circ \ad x \circ \kappa_\g(x)(y)= 0 
\end{eqnarray*}
(Exercise~\ref{exer:3.5.7}). 

(ii) If $\kappa_\g(x)$ is not injective, then the Mackey completeness 
of $\g$ implies the existence of 
$n \in \N$ and $0\not= y \in \g$ with 
$$ (\ad x)^2y = - 4 \pi^2 n^2 y. $$
(Proposition~\ref{prop:kerexp-locconv}). 
Now $[x,y] \not=0$ and $\kappa_\g(x)(y)= 0$, so that 
(i) implies that the non-constant curve $e^{\R \ad y}x$ is contained 
in $\exp_G^{-1}(\exp_G(x))$. Hence $\exp_G$ is not injective 
in any neighborhood of $x$. 

(iii) (a) We consider the curve $\gamma(t) := \exp_G(e^{t\ad y}x)$ as in (i), 
and recall that
$$ \delta(\gamma)_t  = \Ad(\exp_G(ty))\kappa_\g(x)[y,x]. $$
Further, 
$$ \gamma(t) 
= \exp_G(ty) \exp_G(x) \exp_G(-ty)  
= \exp_G(ty) \exp_G(y) \exp_G(-ty) = \be $$
implies that $\gamma'(t) =0$, so that the injectivity of 
$\kappa_\g(x)$ leads to $[y,x]=0$. 
Then Lemma~\ref{lem:4.1.4} yields 
$$ \exp_G(x - y) = \exp_G(x) \exp_G(y)^{-1} = \be. $$

(iii) (b) From \eqref{eq:kappa-comm} we derive 
$\kappa_\g(x) \circ \ad x = \be - e^{-\ad x},$
so that $e^{\ad x} = \be$ is equivalent to  $\ad x = 0$ 
because $\kappa_\g(x)$ is invertible. 

(iii) (c) follows from (b). 

(iv) If $\kappa_\g(x)$ is not injective, then the argument under 
(ii) shows that $x$ is not isolated in $\exp_G^{-1}(\exp_G(x))$.
Suppose, conversely, that $\kappa_\g(x)$ is injective. In view of (iii), 
$\exp_G(y) = \exp_G(x)$ 
implies $\exp_G(x - y) = \be$. If $0$ is isolated in $\exp_G^{-1}(\be)$, 
it follows that $x$ is isolated in $\exp_G^{-1}(\exp_G(x))$. 

(v) Since $\fa$ is abelian, the restriction of $\exp_G$ to $\fa$ is a smooth 
group homomorphism $\fa \to G$ whose kernel is $\Gamma_\fa$ 
(Lemma~\ref{lem:4.1.4}). 
Hence $\Gamma_\fa$ is a closed subgroup of $\fa$. Since 
$\L(\exp_G\res_\fa)$ is the inclusion map $\fa \into \g$, it is injective, 
which in turn implies that the derivative of each curve 
$\gamma \: [0,1] \to \Gamma_\fa$ lies in $\ker \L(\exp_G\res_\fa)= \{0\}$. 
Hence each $C^1$-curve in $\Gamma_\fa$ is constant. If $E \subeq \fa$ is 
finite-dimensional, then $E \cap \Gamma_\fa$ is a closed subgroup of $E$ 
containing no non-trivial vector subspace, so that it 
must be discrete 
(Exercise~\ref{exer:3.5.6}). 
\end{prf}

\begin{cor} \mlabel{cor:expinj} 
Let $G$ be a Lie group with an exponential function, Lie algebra $\g$, 
and $U \subeq \g$ be a symmetric subset such that $\kappa_\g(x)$ is injective for 
each $x \in U + U$. Then, for $x,y \in U$, 
\[  \exp_G x = \exp_G y \quad \Longleftrightarrow \quad 
x - y \in \z(\g) \quad \mbox{ and } \quad \exp_G(x - y) = \be. \] 
\end{cor}

\begin{prf} If $x - y \in \z(\g)$ satisfies $\exp_G(x - y) = \be$, then we have 
$$ \exp_G y = \exp_G(x + (y-x)) = \exp_G x \exp_G(x - y) = \exp_G x$$
(Lemma~\ref{lem:4.1.4}(i)). 

Suppose, conversely, that $x,y \in U$ satisfy $\exp_G x = \exp_G y$. 
Then Proposition~\ref{prop:exp-noninj} implies $\exp_G(x-y) = \be$. 
Since $\kappa_\g(x-y)$ is invertible, we further obtain $x-y \in \z(\g)$. 
\end{prf}

\begin{rem} (a) 
Note that, in the situation of the preceding corollary, we may always assume that 
$U + \z(\g) = U$ because the subset $V := U + \z(\g)$ also is symmetric and 
the relation 
$\kappa_\g(x+z) = \kappa_\g(x)$ implies that  
$\kappa_\g(x)$ is injective for each $x \in V + V = U + U + \z(\g)$. 

If we assume, in addition, that the closed subgroup 
$\Gamma_{\z(\g)} := \ker(\exp_G\res_{\z(\g)})$ intersects 
$U + U$ only in~$\{0\}$, it follows that $\exp_G$ is injective on $U$. 

(b) If $G$ is a Banach--Lie group and 
$\g = \L(G)$ carries a norm with $\|[x,y]\| \leq \|x\|\cdot\|y\|$, 
then $\| \ad x\| \leq \|x\|$. Therefore $\|x\| < 2\pi$ implies 
that $\kappa_\g(x)$ is invertible (Remark~\ref{rem:ban-dexp}). 
 If $\exp_G\res_{\z(\g)}$ is injective, 
i.e., $Z(G)$ is simply connected, Corollary~\ref{cor:expinj} 
implies 
that $\exp_G$ is injective on the open ball $B_\pi := \{ x \in \g \: \|x\| < \pi\}$ 
(cf.\ \cite{LaTi66}). In general, we may put 
$$ \delta_G := \inf \{ \|x\|\: 0 \not= x \in \Gamma_{\z(\g)}\} $$
to obtain that $\exp_G$ is injective on the ball of radius 
$r := \min \{ \pi, \frac{\delta_G}{2}\}$
(cf.\ \cite{GN03}; \cite[Rem.~2.3]{Bel04}). 
\end{rem} 

\begin{rem} It is an important problem to decide for a given 
Lie group $G$ with exponential function 
if $0$ is isolated in $\exp_G^{-1}(\be)$. This is certainly the case 
if $\exp_G$ is injective on some $0$-neighborhood in $\g$. 
Presently, we are note aware of any Lie group $G$ for which 
$0$ is not isolated in $\exp_G^{-1}(\be)$. 

This property can also be interpreted as some rigidity 
for morphisms $\gamma \: \T\cong \R/\Z \to G$. For any such 
morphism $\gamma$, we have $\exp_G(\L(\gamma)1) = \be$, 
and, conversely, any $x \in \exp_G^{-1}(\be)$, defines by 
$\gamma_x(t+\Z) := \exp_G(tx)$ a morphism $\gamma_x \: \T \to G$. 
Therefore $0$ is isolated in $\exp_G^{-1}(\be)$ if any morphism 
$\T \to G$ which is $C^1$-close to the trivial one is actually 
trivial. 

If $0$ is isolated in $\exp_G^{-1}(\be)$, then 
it follows in particular, that the closed subgroup 
$$ \Gamma_Z := \{ x \in \z(\g) \: \exp_G x = \be\} $$
of $\z(\g)$ is closed, a property that we shall encounter again 
and again in other contexts. In view of 
Proposition~\ref{prop:exp-noninj}, for $x \in \g\setminus 
\z(\g)$ satisfying $\exp_G x =\be$, the operator 
$\kappa_\g(x)$ is not injective. 
\end{rem}

\subsection*{Sternberg's example} 
\mlabel{subsec:sternberg}

Important examples of {\it pro-nilpotent Lie groups}, i.e., 
\index{Lie group!pro-nilpotent} 
projective limits of nilpotent Lie groups, arise 
as groups of formal diffeomorphisms. 
For $\K \in \{\R,\C\}$, we write $\Gf_n(\K)$ for the group 
of {\it formal diffeomorphisms of 
$\K^n$} fixing $0$. \index{formal diffeomorphism} 
The elements of this group are represented by formal power series of the form 
\[  \phi(x) = gx + \sum_{|{\bf m}| > 1} c_{\bf m} x^{\bf m}, \quad \mbox{ where } \quad 
g \in \GL_n(\K),\] 
${\bf m} = (m_1,\ldots, m_n) \in \N_0^n,$ and 
$$ |{\bf m}| := m_1 + \ldots + m_n, \quad 
x^{\bf m} := x_1^{m_1} \cdots x_n^{m_n}, \quad c_{{\bf m}} \in \K^n. $$
The group operation is given by composition of power series. 
We call $\phi$ {\it pro-unipotent} \index{formal diffeomorphism!pro-unipotent} 
if $g = E_n$ is the unit matrix. 
It is easy to see that the pro-unipotent formal diffeomorphisms form a pro-nilpotent 
Lie group $\Gf_n(\K)_1 = \prolim G_k$, where 
$G_k$ is the finite-dimensional nilpotent Lie group obtained by 
composing polynomials of the form 
$$ \phi(x) = x + \sum_{1 < |{\bf m}| \leq k} c_{\bf m} x^{\bf m} $$
modulo terms of order $> k$. The group $\Gf_n(\K)$ of all formal diffeomorphisms 
of $\K^n$ fixing $0$ is a semidirect product 
$$ \Gf_n(\K) \cong \Gf_n(\K)_1 \rtimes \GL_n(\K),  $$
where the group $\GL_n(\K)$ of linear automorphisms acts by conjugation. 
As this action is smooth, $\Gf_n(\K)$ is a Fr\'echet--Lie group (Definition~\ref{def:semdir}). 

\begin{ex} The  group $\Gf_n(\C)$ has been studied by Sternberg 
in \cite{St61}, 
where he shows in particular 
that, for $n = 1$, the elements 
$$ \phi_m(x) = e^{\frac{2 \pi i}{m}} x + p x^{m+1}, 
\quad m \in \N\setminus \{1\}, p \in \C^\times, $$
are not contained in the image of the exponential function. This is of particular 
interest because $\phi_m \to \be$ in the Lie group $\Gf_n(\C)$, 
so that the image of the exponential function in this group 
is not an identity neighborhood. 
A detailed analysis of the exponential function 
of this group can also be found in Lewis' paper \cite{Lew39}; 
see also \cite{Rob02}.  

To see that $\phi_m$ is not in the image of the exponential function 
of $\Gf_n(\C)$, 
it suffices to verify this in the finite-dimensional solvable quotient group 
\break $G_{m+1} \rtimes \C^\times$, i.e., modulo terms of order $m+2$. 
The subgroup $H := \C x^{m+1} \rtimes \C^\times$ 
is isomorphic to 
$\C \rtimes \C^\times$ with the multiplication 
$$ (z,w)(z',w') = (z + w^{-m}z', ww') $$
because 
$w(w^{-1}x)^{m+1} = w^{-m} x^{m+1}. $
Its  exponential function is  given by 
$$ \exp_H(z,w) = \Big( \frac{e^{-wm} - 1}{(-wm)}z, e^w\Big)= \Big(\frac{1 - (e^w)^{-m}}{wm}z, 
e^w\Big) $$ 
(cf.\ Proposition~\ref{prop:semdir-exp}). 
We conclude that 
$\phi_m \in H$ is not contained in the exponential image. 

Any element $\xi \in \g_{m+1} \rtimes \C$ 
with $\exp\xi = \phi_m$ must be contained in 
the plane $\Spann_\C \{x, x^{m+1}\}$. In fact, 
$Z := \C x^{m+1}$ is central in the nilpotent group $G_{m+1}$ and 
the eigenvalues for the action of 
$w_0 := e^{\frac{2\pi i}{m}}$ on $\g_m \cong \g_{m+1}/\C x^{m+1}$ are 
$$ w_0^{m-1}, w_0^{m-2}, \ldots, w_0, $$
all different from $1$. Therefore  the exponential 
function of the solvable Lie algebra $\g_m \rtimes \C$ 
is regular in $(0, \frac{2\pi i}{ m})$ 
(cf.\ Remark~\ref{rem:ban-dexp}). Hence any $z$ with 
\[ \exp_{G_m}(z) = \exp_{G_m}\Big(0, \frac{2\pi i}{m}\Big) \] 
is of the form $(a,b) \in \g_m \rtimes \C$ with 
$(a,b)$ commuting with $w_0$ (Proposition~\ref{prop:exp-noninj}), 
but this implies $a = 0$. Hence 
any $\xi \in \g_{m+1} \rtimes \C$ with $\exp_{G_{m+1}} \xi = \phi_m$ 
is contained in $\Spann_\C \{ x, x^{m+1}\}$. 
This completes the proof. 
\end{ex}

\begin{small}
\subsection*{Exercises for Section~\ref{sec:3.5}} 

\begin{exer}[Cauchy Product Formula] \mlabel{exer:2.1.3c}
\index{Cauchy!product formula} \index{Cauchy, Augustin Louis
(1789--1857)} Let $X, Y, Z$ be Banach spaces and \break 
$\beta \: X \times Y\to Z$ a continuous bilinear map. Suppose that $x :=
\sum_{n=0}^\infty x_n$ is absolutely convergent in~$X$ and that $y
:= \sum_{n=0}^\infty y_n$ is absolutely convergent in~$Y$. Then
\[ \beta(x,y) = \sum_{n = 0}^\infty \sum_{k=0}^n \beta(x_k, y_{n-k}). \]
\end{exer}

\begin{exer}[Cauchy Product Formula] \mlabel{exer:2.1.3b}
Let $X, Y, Z$ be locally convex spaces where $Z$ is sequentially complete and 
$\beta \: X \times Y\to Z$ be a continuous bilinear map. Suppose that $x :=
\sum_{n=0}^\infty x_n$ is absolutely convergent in~$X$ and that $y
:= \sum_{n=0}^\infty y_n$ is absolutely convergent in~$Y$. Then
\[ \beta(x,y) = \sum_{n = 0}^\infty \sum_{k=0}^n \beta(x_k, y_{n-k}). \]
Hint: Reduce this verification to the Banach case (Exercise~\ref{exer:2.1.3c}). 
One can w.l.o.g.\ assume that $Z$ is a Banach space. Then show that $\beta$ factors through 
a bilinear maps $\tilde X \times \tilde Y \to Z$, where $\tilde X$ and $\tilde Y$ 
are Banach spaces. 
\end{exer}

\begin{exer} \mlabel{exer:5.6.2b} 
Let $E$ be a real topological vector space.  Show that every continuous
group homomorphism $\gamma \: (\R,+) \to (E,+)$ can be written as 
$\gamma(t) = t v$ for some  $v \in E$. 
\end{exer}

\begin{exer} \mlabel{exer:3.5.1} Show that, for each discrete subgroup 
$D \subeq \R^n$, there exist linearly independent elements 
$v_1, \ldots, v_k \in \R^n$ with $D = \sum_{i=1}^k
\Z v_i$. Use induction on $\dim(\Spann D)$ and proceed along the 
following steps: 
\begin{enumerate}
\item[\rm(1)] Show that $D$ is closed. 
\item[\rm(2)] Reduce to $\Spann D = \R^n$. 
\item[\rm(3)] Every compact subset $C \subeq \R^n$ intersects $D$ in a finite subset. 
\item[\rm(4)] Assume that $\Spann D = \R^n$ and 
that there exists a basis $f_1,\ldots, f_n$ of $\R^n$, contained in $D$, such that 
the hyperplane $F := \Spann \{f_1,\ldots, f_{n-1}\}$ satisfies 
$F \cap D = \Z f_1 + \ldots + \Z f_{n-1}$. 
Show that 
$$ \delta := \inf \Big\{ \lambda_n > 0 \: (\exists \lambda_1, \ldots, \lambda_{n-1}\in \R)\ 
\sum_{i = 1}^n \lambda_i f_i \in D \Big\} > 0. $$
Hint: It suffices to assume $0 \leq \lambda_i \leq 1$ for $i =1,\ldots, n$. 
\item[\rm(5)] Apply induction on $n$ to find $f_1,\ldots, f_n$ as in (4) 
and pick $f_n' := \sum_{i = 1}^n \lambda_i f_i \in D$ with $\lambda_n = \delta$. 
Show that $D = \Z f_1 + \ldots + \Z f_{n-1} + \Z f_n'$. 
\end{enumerate}
\end{exer} 

\begin{exer} \mlabel{exer:3.5.3} (The structure 
of connected abelian Banach--Lie groups) 
(cf.\ Proposition~\ref{prop:reg-ab}) 
Let $A$ be a connected abelian Banach--Lie group. Show that 
\begin{enumerate}
\item[\rm(1)] $\exp_A \: (\L(A),+) \to A$ is a covering morphism of Lie groups. 
\item[\rm(2)] $\Gamma_A := \ker \exp_A$ is a discrete subgroup of $(\L(A),+)$ and 
$A \cong \L(A)/\Gamma_A$. 
\item[\rm(3)] If $\dim A < \infty$, then $A \cong \R^k \times \T^m$ 
for some $k,m \geq 0$. Hint: Exercise~\ref{exer:3.5.1}. 
\end{enumerate}
\end{exer}

\begin{exer} \mlabel{exer:3.5.4} (Divisible groups) An 
abelian group $D$ is called {\it divisible} if, for 
each $d \in D$ and \index{group!divisible abelian} 
$n \in \N$, 
there exists an $a \in D$ with $a^n = d$. Show that: 
\begin{enumerate}
\item[\rm(1)] If $G$ is an abelian group, $H$ a subgroup and 
$f \: H \to D$ a homomorphism into a divisible group $D$, then there exists 
an extension of $f$ to a homomorphism 
$\tilde f \: G \to D$. Hint: Use Zorn's Lemma to reduce the situation to the case where 
$G$ is generated by $H$ and one additional element. 
\item[\rm(2)] If $G$ is an abelian group and $D$ a divisible subgroup, then 
$G \cong D \times H$ for some subgroup $H$ of $G$. Hint: Extend $\id_D \: D \to D$ to a 
homomorphism $f \: G \to D$ and put $H := \ker f$. 
\end{enumerate}
\end{exer}

\begin{exer} \mlabel{exer:3.5.5} (General abelian Banach--Lie groups) 
Let $A$ be an abelian Banach--Lie group. Show that: 
\begin{enumerate}
\item[\rm(1)] The identity component of $A_0$ of $A$ 
is divisible (cf.\ Proposition~\ref{prop:reg-ab}).  
\item[\rm(2)] $A \cong A_0 \times \pi_0(A)$, where $\pi_0(A) := A/A_0$. 
\item[\rm(3)] If $\dim A< \infty$, then there exists a discrete abelian group $D$ and 
$k,m \in \N_0$ with $A \cong \R^k \times \T^m \times D$. 
\end{enumerate}
\end{exer}

\begin{exer} \mlabel{exer:3.5.6} Let $E$ be a finite-dimensional space and 
$\Gamma \subeq E$ be a closed subgroup. \
Show that, if $\Gamma$ is not discrete, then $\Gamma$ contains a non-trivial vector 
subspace. Hint: If $\Gamma$ is not discrete, then there exists a sequence 
$(\gamma_n)$ of non-zero elements in $\Gamma$ converging to $0$. 
Pick a norm $\|\cdot\|$ on $E$ and show that we may assume that 
$v_n := \frac{\gamma_n}{\|\gamma_n\|}$ converges in the unit sphere of $(E,\|\cdot\|)$ 
to some $v$. Then show that $\R v \subeq \Gamma$. 
\end{exer}

\begin{exer}
  \mlabel{exer:3.5.7} 
Let $\g$ be a Mackey complete Lie algebra, $x \in \g$ such that 
$\ad x$ generated a smooth one-parameter group $e^{t\ad x}$ 
of automorphisms of $\g$
 and $\phi \in \Aut(\g)$. Show that: 
\begin{enumerate}
\item[\rm(1)] $\ad(\phi(x))$ also generates a smooth one-parameter group on $\g$, 
and we have 
$$ \phi \circ e^{\ad x} \circ \phi^{-1} = e^{\ad \phi(x)}. $$
\item[\rm(2)] $\kappa_\g(\phi(x)) = \phi \circ \kappa_\g(x) \circ \phi^{-1}$ holds for 
$\kappa_\g(x) := \int_0^1 e^{-t\ad x}\, dt$. 
\end{enumerate}
\end{exer}

\begin{exer}\mlabel{exer:I.6} Let $X$ be a Banach space and let 
$G$ be a finite-dimensional connected Lie group,  
and $\pi \: G \to \GL(X)$ be a faithful representation which is continuous when 
$\GL(X)$ carries the uniform topology inherited from the Banach algebra 
${\cal L}(X)$ and for which $\pi(G)$ is bounded. 
Show that $\g := \Lie(G)$ is a compact Lie algebra by using the following steps: 
\begin{enumerate}
\item[\rm(1)] $\pi$ is a smooth homomorphism of Lie groups. In particular,  
we have a representation of the Lie algebra $\Lie(\pi) \: \g \to {\cal L}(X)$. 
\item[\rm(2)] $\|x\| := \|\Lie(\pi)(x)\|$ defines a norm on $\g$, and 
$\Ad(G)$ is bounded with respect to this norm. 
\item[\rm(3)] $\Ad(G)$ has compact closure, so that $\g$ is a compact Lie algebra. 
\end{enumerate}

If, in addition, $X$ is a Hilbert space, then one can even show that there exists a 
scalar product compatible with the topology which is invariant under $G$, 
so that $\pi$ becomes a unitary representation with respect to this scalar 
product. This can be achieved by showing 
that the set of all compatible scalar products is a Bruhat--Tits space and 
then applying the Bruhat--Tits Fixed Point Theorem (\cite{La99}). 
\end{exer} 

\begin{exer} Let $\fn$ be a $2$-step nilpotent Lie algebra, endowed 
with  the BCH multiplication 
$$ x * y := x + y + \frac{1}{2} [x,y]. $$
Further let $D \in \der(\fn)$ and $\g := \fn \rtimes_D \R$. 
Calculate the exponential function of the group 
$G := (\fn,*) \rtimes_\alpha \R$ for $\alpha(t) := e^{tD}$ explicitly. 
\end{exer}

\begin{exer} (Torsion elements in Banach--Lie groups) 
Let $G$ be a Banach--Lie group and 
\[ E_n(G) := \{ g \in G \:  g^n = \be \} \] 
denote the closed subset of those elements whose order divides~$n$. 
Show that each $E_n(G)$ carries the structure of a Banach manifold 
and that the identity component $G_0$ of $G$ acts transitively on each connected 
component of $E_n(G)$\begin{footnote}{See \cite{AW08} for the 
case of finite-dimensional Lie groups.}\end{footnote}\\
Hint: Verify the following assertion: 
\begin{enumerate}
\item[\rm(1)] For $g \in E_n(G)$, consider the closed subspaces 
\[ \fa := \ker(\Ad(g)-\1) \quad \mbox{ and } \quad 
\fb := \ker(\1 + \Ad(g) + \cdots + \Ad(g)^{n-1}) = \im(\Ad(g)-\1). \]
Show that $\g = \fa \oplus \fb$. Hint: $(z^n-1) = (z-1)(z^{n-1} + \cdots + z + 1)$ for 
$z \in\C$ and the complex linear extension of $\Ad(g)$ to $\g_\C$ 
is diagonalizable.
\item[\rm(2)] For the map 
\[ F \: \fa \times \fb \to G, \quad F(x,y) := \exp(y) g \exp(x) \exp(-y) \] 
there exist open $0$-neighborhoods $U_\fa \subeq \fa$ and $U_\fb \subeq \fb$ 
such that $F\res_{U_\fa \times U_\fb}$ is a diffeomorphism onto an open subset 
of $G$. Hint: Show that 
\[ T_{(x,y)}(F)(v,w) = g.(v + \Ad(g)^{-1}w - w)\] 
and use the Inverse Function Theorem.  
\item[\rm(3)] If $U_\fa$ is sufficiently small and $(x,y) \in U_\fa \times U_\fb$, 
then 
$F(x,y) \in E_n(G)$ is equivalent to $x = 0$. 
\item[\rm(4)]Conclude that $E_n(G)$ is a 
submanifold of $G$ and that the $G_0$-conjugation orbits are the connected 
components of $E_n(G)$. 
\end{enumerate}
\end{exer}

\end{small}

\section{Notes and comments on Chapter~\ref{ch:3}} 

In our discussion of Lie groups, we essentially follow \cite{Mil82, Mil84}, 
but, as for manifolds, we do not assume that the model space 
of a Lie group is complete (\cite{Gl02a}). 

\nin{\bf Section~\ref{sec:3.3}:} Theorem~\ref{thm:locglob} is an infinite-dimensional 
version of a classical result for finite-dimensional groups 
(cf.\ \cite{Ch46}, \S 14, Prop.~2; or \cite{Ti83}, p.14). 

\nin{\bf Section~\ref{sec:3.5}:} 
Proposition~\ref{prop:reg-ab} is due to Michor/Teichmann (\cite{MT99}). 

Our presentation of the derivative of the exponential 
function is inspired by Grabowski's paper \cite{Gr93}. 
The formula in Remark~\ref{rem:diffexp-alg} for the differential of 
the exponential function of a Banach algebra can already be found in 
Michael's paper \cite{MicA45} and earlier versions already appeared
in \cite{Kn30}.

The global behavior of the exponential function and in 
particular the question of its surjectivity is a quite complicated issue, depending 
very much on specific properties of the groups under consideration  
(cf.\ \cite{Wue03, Wue05}). 

For finite-dimensional Lie groups, the most basic general result is that 
if $G$ is a connected Lie group with compact Lie algebra $\g$, then 
$\exp_G$ is surjective (\cite{Ho65}). Since the compactness of $\g$ is equivalent to the 
existence of an $\Ad(G)$-invariant scalar product on $\g$, which in turn leads to a 
biinvariant Riemannian metric on $G$, the surjectivity of $\exp_G$ can be 
derived from the Hopf--Rinow Theorem in Riemannian geometry or the 
{\it Maximal Torus Theorem}, asserting that each compact Lie group is the union 
of its maximal tori. 

A natural generalization of the notion of a compact Lie algebra to 
the Banach context is obtained by calling 
a real Banach--Lie algebra $(\g,\|\cdot\|)$  
{\it elliptic} if the norm on $\g$ is \index{Banach--Lie algebra!elliptic} 
invariant under the group $\Inn(\g) := \la e^{\ad \g}\ra \subeq \Aut(\g)$ 
of {\it inner automorphisms} (cf.\ \cite[Def.~IV.3]{Ne02c}). 
A finite-dimensional Lie algebra $\g$ is elliptic with respect to some
norm if and only if it is compact. 
In this case the requirement of an invariant scalar product
leads to the same class of Lie algebras, but in the
infinite-dimensional context this is different. Here the existence of 
an invariant scalar product turning $\g$ into a real Hilbert space
leads to the structure of a 
\index{Hilbert--Lie algebra} 
{\it Hilbert--Lie algebra}, 
resp., an $L^*$-algebrea of compact type. 
Simple Hilbert--Lie algebras can be
classified, and each Hilbert--Lie algebra 
is a Hilbert space direct sum of
simple ideals and its center (cf.\ \cite{Sc60}, \cite{Sc61}, 
\cite{dlH72}, \cite{CGM90}, \cite{Neh93}, \cite{St99}). In particular, the
classification shows that every Hilbert--Lie algebra can be
realized as a closed subalgebra of the 
Hilbert--Lie algebra $\fu_2({\cal H})$ of 
skew-hermitian Hilbert--Schmidt operators 
on a complex Hilbert space ${\cal H}$. 

The class of elliptic Lie algebras is much bigger. It contains
the algebra $\fu(\cA)$ of skew-hermitian elements of any $C^*$-algebra
$\cA$ and in particular the Lie algebra $\fu({\cal H})$ of the full unitary
group $\U({\cal H})$ of a Hilbert space ${\cal H}$. 

Although finite-dimensional connected 
Lie groups with compact Lie algebra have a surjective exponential
function, this is no longer true for 
connected Banach--Lie groups with elliptic Lie
algebra. This is a quite remarkable phenomenon discovered  by Putnam and
Winter in \cite{PW52}: the orthogonal group $\OO({\cal H})$ of a real 
infinite-dimensional Hilbert 
space is a connected Banach--Lie group with elliptic Lie algebra, but
its exponential function is {\sl not} surjective. 
In fact, a reflection $g \in \OO(\cH)$ in a 
hyperplane $E \subeq \cH$ is not contained in the exponential 
image because any $X \in \fo(\cH)$ with $e^X = g$ commutes 
with $g$, so that $X^\top = - X$ and $\dim E^\bot = 1$ 
lead to $X\res_{E^\bot} = 0$, contradicting 
$g\res_{E^\bot} = - \id$. This contrasts the fact that 
the exponential function of the unitary group $\U({\cal H})$ of a complex Hilbert space 
is always surjective. 

\nin{\bf Discrete subgroups of Banach spaces:} 
Discsrete central subgroups of Lie groups and, in particular, 
discrete subgroups of Banach spaces play an important role
in Lie theory. A systematic exploration of closed and in particular 
discrete subgroups of Banach spaces has been undertaken in 
\cite{ADG94}. Of fundamental importance in this context is the following: 
\begin{thm} {\rm(\cite[Thm.~1.1]{ADG94})} 
For a closed subgroup $\Gamma$ of a separable complete metric topological 
vector space $E$, the following are equivalent: 
\begin{enumerate}
\item[\rm(a)] $\Gamma$ is discrete. 
\item[\rm(b)] $\Gamma$ is locally compact and contains no lines. 
\item[\rm(c)] $\Gamma$ is countable. 
\item[\rm(d)] $\Gamma \cong \Z^{(J)}$ for some set $J$, i.e., $\Gamma$ is a 
free abelian group. 
\end{enumerate}  
\end{thm}
The paper \cite{ADG94} also contains results on closed subgroups 
of the sequence space $\R^\N$ which have later been embedded into a 
systematic duality theory for weakly complete spaces (spaces of the form 
$\R^J$, $J$ a set) by K.H.\,Hofmann and S.\,Morris in the context of projective 
limits of finite-dimensional Lie groups (\cite{HoM07}). 
We also mention \cite[Thm.~4]{ADG94}, asserting that a Banach space 
$E$ contains a closed, non-discrete, totally disconnected subgroup 
if and only if it contains a copy of $c_0 = c_0(\N)$, 
the space of $0$-sequences. One implication is easy to see, because
for every sequence $a = (a_n)_{n \in \N} \in c_0$, 
the subgroup 
\[ \Gamma := \{ b = (b_n)_{n \in \N} \in c_0 \: 
(\forall n \in \N)\, b_n \in \Z a_n \} \subeq c_0\] 
is a non-discrete, weakly closed closed, totally disconnected 
subgroup of~$c_0$.

\nin{\bf Epimorphisms:} A morphism of Banach--Lie groups 
$f \: G \to H$ is called an {\it epimorphism} if, 
for any two morphisms $h_1, h_2 \: H \to L$ of Banach--Lie groups, 
the equality $h_1 \circ f = h_2 \circ f$ implies $h_1 = h_2$. 
Every morphism with dense range is an epimorphism, but there are many 
interesting examples of epimorphisms with non-dense range, such as 
the inclusion of the Borel subgroup in $\SL_2(\R)$ (\cite[Thm.~6.3]{PU21}). 
For a recent discussion of epimorphisms in several categories of 
infinite-dimensional groups, we refer to \cite{PU21}. 

\chapter{Regular Lie groups} \mlabel{ch:4}
%
%

\nin Since the basic theory of ordinary differential equations works perfectly 
for Banach manifolds, for any Banach--Lie group $G$ and any continuous 
curve $\xi \: [0,1] \to \L(G)$, the initial value problem 
\begin{eqnarray}
  \label{eq:ivp0}
\gamma(0) = \be, \quad  \gamma'(t) = \gamma(t)\xi(t) 
\end{eqnarray}
has a unique solution (see Definition~\ref{def:3.1.2a} for the notation). 
Beyond the class of Banach--Lie groups, 
no general result exists that guarantees the existence of a solution to 
this kind of initial value problem, not even for constant curves~$\xi$. 
On the one hand, there are many 
specific classes of Lie groups for which 
one can prove the existence of solutions to (\ref{eq:ivp0}) by direct 
means, depending on the group under consideration. 
On the other hand, we do not know of any Lie group $G$, modeled on a 
Mackey complete space, for which (\ref{eq:ivp0}) does not have a solution  
for some smooth curve $\xi$. The notion of a regular Lie group, 
introduced in this chapter, provides a natural axiomatic background 
to deal with this problem. 

In particular, we shall see how regularity of a Lie group can be used to 
obtain a generalization of the Fundamental Theorem of Calculus to Lie group-valued 
smooth functions, which provides a characterization of all those 
$\L(G)$-valued $1$-forms on a connected manifold $M$ which arise 
as logarithmic derivatives of $G$-valued smooth functions. 

This leads to solutions to many integrability questions. 
For example, for each
continuous Lie algebra 
homomorphism $\psi \:  \L(G) \to \L(H)$ 
from the Lie algebra of 
a $1$-connected Lie group $G$ to the Lie algebra of a regular 
Lie group $H$, there exists a unique morphism of Lie groups $\phi$ with 
$\L(\phi) = \psi$. As an important consequence, all $1$-connected 
regular Lie groups with the same Lie algebra are isomorphic as 
Lie groups. 

All Banach--Lie groups and all Lie groups 
which are quotients of Mackey complete vector spaces by discrete subgroups are regular 
and, conversely, the Lie algebra of a regular Lie group is Mackey complete. 
However, presently no non-regular Lie group modeled on a Mackey complete space is known. 

{\bf Notation:} Throughout this  chapter, we write $I := [0,1]$ and 
$C^k_*(I,G)$ denotes the set of all $C^k$-paths in $G$ 
starting at~$\be$. 

{\bf Prerequisites:} Sections~\ref{seccuint}-\ref{secCk} and 
\ref{secTay}, and Chapters~\ref{chapmanif} and~\ref{ch:3}.

\section{Regularity and the Fundamental Theorem} \mlabel{sec:4.3} 

In this section, we shall see how regularity of a Lie group can be used to 
obtain a Fundamental Theorem of Calculus for Lie group-valued smooth functions. 
This implies solutions to many integrability questions. 
For example, for each homomorphism $\psi \:  \L(G) \to \L(H)$ from the Lie algebra of 
a connected simply connected Lie group $G$ into the Lie algebra of a regular 
Lie group $H$, there exists a unique morphism of Lie groups $\phi$ with 
$\L(\phi) = \psi$.

\subsection{Regular Lie groups} 

We shall introduce a hierarchy of regularity conditions, called 
$C^k$-regularity, where, for $k \leq r$, the $C^r$-regularity follows 
from $C^k$-regularity. The most important of these notions is the 
$C^\infty$-regularity, introduced by J.~Milnor in \cite{Mil82}. 
The main point in discussing the stronger concepts of $C^k$-regularity is 
that they provide additional information
which is particularly useful in applications.
The Uniqueness Lemma~\ref{lem:c.12b} 
as well as the Product and Quotient Rules
for logarithmic derivatives (Lemma~\ref{lem:c.12})
will be essential. They apply just as well to $C^{k+1}$-functions
(with $k\in \N_0\cup\{\infty\}$), with identical proof.
In this case, the logarithmic derivative is a Lie algebra-valued
$1$-form of class $C^k$.
When applied to a $C^{k+1}$-path $\gamma\colon [0,1]\to G$, the logarithmic derivative
yields a  $C^k$-form which can be identified with
the $C^k$-path $\delta(\gamma)\colon [0,1]\to \g$, $t\mto \gamma(t)^{-1}\dot{\gamma}(t)$,
using the smooth left action of $G$ on $TG$ inherited from left translations.
\begin{defn} \mlabel{def:regular0} 
Let $k \in \N_0 \cup \{\infty\}$.
A Lie group $G$ with Lie algebra $\g$ is called 
\index{Lie group!$C^k$-regular} 
\index{Lie group!regular} 
\index{regular Lie group} 
{\it $C^k$-semiregular}, if for each
$\xi \in C^k(I,\g)$, the initial value problem 
\begin{eqnarray}
  \label{eq:lie-ivp}
\gamma(0) = \be, \quad \delta(\gamma) = \gamma^{-1}\dot\gamma= \xi 
\end{eqnarray}
has a solution $\gamma_\xi$, which is then contained in 
$C^{k+1}(I,G)$.
If, moreover,
the corresponding evolution map 
$$ \evol_G \: C^k(I,\g) \to G, \quad \xi \mapsto \gamma_\xi(1) $$
is smooth, then $G$ is called \emph{$C^k$-regular}.
We recall from the Uniqueness Lemma~\ref{lem:c.12b} that the solutions of 
(\ref{eq:lie-ivp})  are unique whenever they exist. 
If $G$ is $C^k$-semiregular, we write 
$$ \Evol_G \: C^k(I,\g) \to C^{k+1}(I,G), \quad \xi \mapsto \gamma_\xi $$
for the corresponding map on the level of Lie group-valued curves.
If $G$ is understood from the context,
we relax notation and also write $\evol$ and $\Evol$
in place of $\evol_G$ and $\Evol_G$.
The group $G$ is called {\it regular} if it is $C^\infty$-regular. 
\end{defn}

\begin{rem} \mlabel{rem:reg1a} (a) Any regular Lie group $G$ 
has a smooth  exponential function 
$$ \exp_G \: \g\to G \quad \hbox{ by } \quad \exp_G(x) := \gamma_x(1), $$
where $x \in \g$ is considered as a constant function $I \to \g$. 
As a restriction of the smooth function $\evol_G$ to the topological subspace 
$\g \subeq C^\infty(I,\g)$ of constant functions, 
the exponential function is smooth. 

(b) Assume that $G$ is $C^k$-semiregular. 
For $\xi \in C^k(I,\g)$, $0 \leq s \leq 1$, and $\eta(t) 
:= \gamma_\xi(st)$, 
we have 
$\delta(\eta)(t) = s\xi(st)$. Therefore we obtain with 
$S(\xi,s)(t) := s\xi(st)$ the relation 
$$ \Evol_G(\xi)(s) = \gamma_\xi(s) = \evol_G(S(\xi,s)). $$
If $k=\infty$ and $G$ is regular, let us show we obtain a smooth mapping
\[ C^\infty(I,\g) \times I \to G, \quad 
(\xi,s) \mapsto \Evol(\xi)(s).\] 
To verify the assertion, we have to verify the smoothness of the map 
\[  S\: C^\infty(I,\g) \times I \to C^\infty(I,\g), \quad 
S(\xi,s)(t) := s\xi(st) \] 
Since the evaluation map $\ev \: I \times C^\infty(I,\g) \to \g$ is smooth 
(Proposition~\ref{prop:smooth-eval}), the map 
$$ I \times C^\infty(I,\g) \times I \to \g, \quad 
(s,\xi,t) \mapsto s \xi(st) $$ 
is also smooth. This in turn implies that
$S$ is smooth (Proposition~\ref{prop:cartes-closed}).

(c) For $k \leq r$, the $C^k$-regularity of a Lie group $G$ implies its 
$C^r$-regularity because the inclusion map 
$C^r(I,\g) \to C^k(I,\g)$ is continuous linear, hence smooth. 
\end{rem} 

We shall see in Theorem~\ref{thm:banach-regular} 
below that Banach--Lie 
groups are regular and in Proposition~\ref{prop:reg-Mack} that the Lie algebra
of a regular Lie group is Mackey complete. 
Since the proofs of these facts requires some more information 
on mapping groups developed in Section~\ref{sec:4.1}, 
we first take a closer look at the simpler 
case of abelian groups. 

\begin{prop} \mlabel{prop:mc-reg} 
Let $E$ be a locally convex space and $k \in \N \cup \{\infty\}$. 
Then the Lie group $(E,+)$ is $C^k$-regular if and only if $E$ is Mackey complete. 
\end{prop}

\begin{prf} If $(E,+)$ is $C^k$-regular, then, 
for each smooth curve $\xi \: [0,1] \to E$, there exists a smooth curve 
$\eta \: [0,1] \to E$ with $\eta' = \xi$, so that $E$ is Mackey complete. 

If, conversely,  $E$ is Mackey complete, 
then each $C^1$-curve $\xi \: [0,1] \to E$ 
has an integral, which leads to a continuous linear map 
$$ {\cal I} \: C^k(I,E) \to C^{k+1}_*(I,E), 
\quad \xi \mapsto \Big(t \mapsto \int_0^t\xi(s)\, ds\Big). $$
Now $\evol_E(\xi) = {\cal I}(\xi)(t)$ 
implies the $C^k$-regularity of~$(E,+)$. 
\end{prf}

The following lemma is an effective tool to verify the 
regularity of a Lie group in terms of local data. 

\begin{lem} \mlabel{lem:reg-crit} 
{\rm(Local regularity criterion)} 
Let $G$ be a Lie group with 
Lie algebra $\g$ and $k \in \N_0 \cup \{\infty\}$. Suppose that 
{\rm(\ref{eq:lie-ivp})} has a solution for each $\xi$ in an 
open $0$-neighborhood $U \subeq C^k([0,1],\g)$. 
Then it has a solution for each $\xi \in C^k([0,1],\g)$. 
If the evolution map 
$$ \evol_G \: C^k([0,1],\g) \to G, \quad \xi \mapsto \gamma_\xi(1) $$
is smooth in $U$, then it is smooth on all of $C^k([0,1],\g)$. 
\end{lem}

\begin{prf} Let $\xi \in C^k([0,1],\g)$. 
For $n \in \N$ and $i \in \{0,\ldots, n-1\}$ we define 
$$\xi_i^n \in C^k([0,1],\g) \quad \mbox{ by } \quad 
\xi_i^n(t) := \frac{1}{n}\xi\Big(\frac{i+t}{n}\Big). $$
We then have for each $m \leq k$: 
$$ \frac{d^m}{dt^m}\xi_i^n(t) = \frac{1}{n^{m+1}}\xi^{(m)}\Big(\frac{i+t}{n}\Big). $$
We conclude that, for $n \to \infty$, the sequence 
$\frac{d^m}{dt^m}\xi_i^n$ tends to $0$ in $C^0([0,1],\g)$, uniformly in $i$, and hence that 
$\xi_i^n$ tends to $0$ in $C^k([0,1],\g)$, uniformly in $i$. 
In particular, there exists some $N \in \N$ for which 
$\xi_i^N \in U$ for $i = 0,1,\ldots, N-1$. 

We define a path $\gamma_\xi \: [0,1] \to G$ by 
$$ \gamma_\xi(t) := 
\gamma_{\xi_0^N}(1)\cdots \gamma_{\xi_{i-1}^N}(1)
\gamma_{\xi_i^N}(nt-i) 
\quad \mbox{ for } \quad \frac{i}{n} \leq t \leq \frac{i+1}{n} $$
and observe that $\delta(\gamma_\xi)$ exists on all of $[0,1]$ and 
equals $\xi$. We now put 
$$\Evol_G(\xi) := \gamma_\xi \quad \mbox{ and } \quad 
\evol_G(\xi) := \gamma_\xi(1). $$

Since, for each $i$, 
the assignment $\xi \mapsto \xi_i^N$ is linear and continuous, 
there exists an open neighborhood $V$ of $\xi$ such that 
$\eta_i^N \in U$ holds for each $\eta \in V$ and $i =0,1,\ldots, N-1$. 
It suffices to show that $\evol_G$ is smooth on $V$, 
but this follows from the fact that 
\[ \evol_G(\eta) = \evol_G(\eta_0^N) \cdots \evol_G(\eta_{N-1}^N)   \] 
is a product of $N$ smooth functions. 
\end{prf}

The preceding lemma shows that regularity of a Lie group is a local 
property. Therefore it immediately implies: 

\begin{prop} \mlabel{prop:4.1.5} 
Let $\phi \: G \to H$ be a covering morphism of Lie groups. 
If one of these groups is $C^k$-regular for some $k \in \N_0$, then 
the other group is also $C^k$-regular. 
In particular, the universal covering group of a $C^k$-regular 
Lie group is $C^k$-regular. 
\end{prop}

\subsection[The Fundamental Theorem] 
{The Fundamental Theorem for Lie group-valued functions} 

\begin{defn} Let $G$ be a Lie group with Lie algebra $\g = \L(G)$. 
We call a $\g$-valued $1$-form 
$\alpha \in \Omega^1(M,\g)$ {\it integrable} 
\index{integrable $1$-form, on Lie group} 
\index{locally integrable $1$-form, on Lie group} 
if there exists a smooth 
function $f \: M \to G$ with $\delta(f) = \alpha$. 
The $1$-form $\alpha$ is said to be 
{\it locally integrable} if each point $m \in M$ has an open 
neighborhood $U$ such that $\alpha\res_U$ is integrable. 
\end{defn}

We now turn to the preparations of the proof 
of the Fundamental Theorem for functions with values in a 
Lie group $G$ with Lie algebra $\g$. 

\begin{lem} \mlabel{lem:laxuni}
If $G$ is $C^k$-regular, 
$x\in \g$ and $\xi \in C^k([0,1],\g)$, then the initial value problem 
\begin{equation}
  \label{eq:e1}
\eta'(t) = [\eta(t),\xi(t)], \quad \eta(0) = x \tag{E1}
\end{equation}
has a unique solution given by 
\begin{equation}
  \label{eq:e2}
\eta(t) =\Ad(\gamma_\xi(t))^{-1}x. \tag{E2}
\end{equation}
\end{lem}

\begin{prf} For 
$\gamma(t) := \gamma_\xi(t)$ we get with Lemma~\ref{lem:c.12}(i)  
$$ \delta(\gamma^{-1})_t = - \Ad(\gamma(t))\delta(\gamma)_t  
= -\Ad(\gamma(t))\xi(t). $$
We define $\eta$ by (E2). Then $\eta$ is a $C^{k+1}$-curve 
with 
\begin{eqnarray*}
\eta'(t) 
&&= \Ad(\gamma(t))^{-1}[- \Ad(\gamma(t))\xi(t), x] 
= [\Ad(\gamma(t))^{-1}x, \xi(t)] = [\eta(t), \xi(t)] 
\end{eqnarray*} 
(Proposition~\ref{prop:der-Ad}). 

Now let $\beta$ be any solution of (E1) and consider the curve 
$$ \zeta(t) :=  \Ad(\gamma(t))\beta(t). $$
Then $\zeta(0) = \beta(0) = x$ and 
Proposition~\ref{prop:der-Ad} lead to 
\begin{eqnarray*}
\zeta'(t) 
&&= \Ad(\gamma(t))[\delta(\gamma)_t, \beta(t)] 
+  \Ad(\gamma(t))\beta'(t) \\
&&= \Ad(\gamma(t))\big([\xi(t), \beta(t)] + \beta'(t)\big) = 0.
\end{eqnarray*}
Therefore $\zeta$ is constant equal to $x$, and we obtain 
\[\beta(t) = \Ad(\gamma(t))^{-1}\zeta(t) = \Ad(\gamma(t))^{-1}x = \eta(t).
\qedhere\] 
\end{prf}

Before we turn to more general situations, we first take a closer  
look at smooth functions $I^2 \to G$, where $I = [0,1]$ is the unit 
interval. 

\begin{lem} \mlabel{lem:mc-square}
  \begin{enumerate}
  \item[\rm(i)] A smooth $\g$-valued $1$-form 
\[  \alpha = v d x + w d y \in \Omega^1(I^2,\g) \quad \mbox{ with } \quad 
v,w \in C^\infty(I^2, \g) \]
satisfies the Maurer--Cartan equation if and only if 
\begin{eqnarray}
  \label{eq:3.3.6}
\frac{\partial v}{\partial y} - \frac{\partial w}{\partial x}  = [v,w]. 
\end{eqnarray}
  \item[\rm(ii)] Suppose that $\alpha$ satisfies the 
Maurer--Cartan equation. 
\begin{enumerate}
\item[\rm(a)] If $f \: I^2 \to G$ is smooth
with $\delta(f)(\partial_y) = w$, and 
$\delta(f)(\partial_x)(x,0) = v(x,0)$ holds for $x \in I$, then 
$\delta(f) = \alpha$. 
\item[\rm(b)] The smooth function $f \: I^2 \to G$, defined by 
$$ f(x,0) := \gamma_{v(\cdot,0)}(x) \quad \mbox{and} \quad  
f(x,y) := f(x,0) \cdot \gamma_{w(x,\cdot)}(y) $$ 
satisfies $\delta(f) = \alpha$. 
\end{enumerate}
\end{enumerate}
\end{lem}

\begin{prf} (i) To evaluate the Maurer--Cartan equation 
for $\alpha$, we first observe that 
$$ \shalf[\alpha, \alpha] 
\big(\frac{\partial}{\partial x}, \frac{\partial}{\partial y}\big)
 = \Big[ \alpha\Big(\frac{\partial}{\partial x}\Big), 
          \alpha\Big(\frac{\partial}{\partial y}\Big)\Big] = [v,w], $$
and obtain 
\begin{align*}
d \alpha + \shalf [\alpha, \alpha] 
&= \frac{\partial v}{\partial y} d y  \wedge d x 
+ \frac{\partial w}{\partial x} d x  \wedge d y + [v,w] d x  \wedge d y  \\
&= \big(\frac{\partial w}{\partial x} 
-\frac{\partial v}{\partial y} + [v,w]\Big)d x \wedge d y. 
\end{align*}

(ii) (a) We write 
$$\delta f = \hat v \dd x + w \dd y \quad \mbox{ with } \quad 
\hat v(x,0) = v(x,0) \quad \mbox{ for} \quad x \in I. $$
The Maurer--Cartan equation for $\delta f$ reads 
$$ \frac{\partial \hat v}{\partial y} - \frac{\partial w}{\partial x}  = [\hat v,w], $$
so that subtraction of this equation from (\ref{eq:3.3.6}) leads to 
$$ \frac{\partial (v- \hat v)}{\partial y} = [v-\hat v,w]. $$
As $(v - \hat v)(x,0) = 0$, the uniqueness assertion of Lemma~\ref{lem:laxuni}, 
applied to $\eta_x(t) := (v - \hat v)(x,t)$, 
implies that $(v - \hat v)(x,y) = 0$ for all $x,y \in I$, hence that 
$v = \hat v$, which means that $\delta(f) = v \dd x + w \dd y = \alpha$.  
  
(b) Since the map 
$I \to C^\infty(I,\g), x \mapsto w(x,\cdot)$
is smooth (Proposition~\ref{prop:cartes-closed}), 
Remark~\ref{rem:reg1a}(b) implies that $f$ is a smooth function. 
By construction, it satisfies the assumptions of (a), so that 
$\delta(f) =\alpha$. 
\end{prf}

\begin{lem} \mlabel{lem:locexist} 
Let $U$ be an open convex subset of the locally convex space~$E$, 
$G$ a regular Lie group with Lie algebra $\g$ 
and $\alpha \in \Omega^1(U,\g)$ satisfying the Maurer--Cartan 
equation. Then $\alpha$ is integrable. 
\end{lem}

\begin{prf} We may w.l.o.g.\ assume that $x_0 = 0 \in U$. 
For $x \in U$ we then consider the smooth curve 
$$ \xi_x \: I \to \g, \quad t \mapsto \alpha_{tx}(x). $$
Then the map $U \to C^\infty(I,\g)$, $x \mapsto \xi_x$ is smooth because the map 
$$ U \times I \to \g, \quad (x,t) \mapsto \xi_x(t) $$
is smooth (Proposition~\ref{prop:cartes-closed}). Hence  the function 
$$ f \: U \to G, \quad x \mapsto \evol(\xi_x) $$
is smooth. 

First we show that $f(sx) = \gamma_{\xi_x}(s)$ holds for each $s \in I$. 
From Remark~\ref{rem:reg1a} we derive that 
$$ S(\xi_x,s)(t) = s \xi_x(st) = \alpha_{stx}(sx) = \xi_{sx}(t), $$
which leads to 
$f(sx) 
= \gamma_{\xi_{sx}}(1) 
= \gamma_{\xi_{x}(s\cdot)}(1) 
= \gamma_{\xi_x}(s)$. 

For $x, x + h \in U$, we consider the smooth map 
$$\beta \: I \times I \to U, \quad (s,t) \mapsto t (x + sh) $$ 
and the smooth function $F := f \circ \beta$. Then the preceding considerations imply 
$F(s,0) = f(0) = \be$, 
\begin{eqnarray*}
\frac{\partial F}{\partial t}(s,t) 
&&= \frac{d}{dt} f(t(x + sh)) 
= \frac{d}{dt} \gamma_{\xi_{x + sh}}(t) = F(s,t)\xi_{x + sh}(t) \\
&&= F(s,t)\alpha_{t(x + sh)}(x + sh) 
= F(s,t)(\beta^*\alpha)_{(s,t)}\big(\frac{\partial}{\partial t}\big). 
\end{eqnarray*}
As we have seen in Lemma~\ref{lem:mc-square}(ii), these relations lead to 
$$ \delta(F) = \beta^* \alpha \quad \mbox{ on } \quad I \times I. $$
We therefore obtain 
$$ \frac{\partial}{\partial s} f(x + sh) 
= \frac{\partial}{\partial s} F(s,1) 
= F(s,1)\alpha_{x + sh}(h) 
= f(x + sh)\alpha_{x + sh}(h), $$
and for $s = 0$ this leads to 
$T_x(f)(h) = f(x)\alpha_x(h),$
so that $\delta(f)=~\alpha$. 
\end{prf}

\begin{prop} \mlabel{prop:cov} 
  Let $M$ be a connected manifold, $G$ a Lie group with Lie algebra $\g$, and
  $\alpha \in \Omega^1(M,\g)$ a $1$-form. If $\alpha$ is locally integrable, 
then there exists a
connected covering $q \: \hat M \to M$ such that $q^*\alpha$ is integrable. 
If, in addition, $M$ is simply connected, then $\alpha$ is integrable. 
\end{prop}

\begin{prf} We consider the product set $P := M \times G$ with the two projection 
maps $F \: P \to G$ and $q\: P \to M$. We define a topology on $P$ as follows. 
For each pair $(U,f)$, consisting of an open subset $U \subeq M$ and a smooth 
function $f \: U \to G$ with $\delta(f) = \alpha\res_U$,  
the graph $\Gamma(f,U) := \{(x,f(x)) \: x \in U\}$ 
is a subset of $P$. 

{\bf Step 1:} We claim that the sets $\Gamma(f,U)$ form a basis for a topology $\tau$ on $P$. 
In fact, let $p = (m,g)\in \Gamma(f_1, U_1) \cap \Gamma(f_2, U_2)$.  
We choose a connected open neighborhood $U$ of $m$, contained in 
$U_1 \cap U_2$. Then 
$f_1(m) = f_2(m) = g$ and $\delta(f_1\res_U) = \delta(f_2\res_U) = \alpha\res_{U}$ 
imply that $f_1 = f_2$ on $U$ (Lemma~\ref{lem:c.12b}), so that 
$p \in \Gamma(f_1\res_U,U) \subeq \Gamma(f_1, U_1) \cap \Gamma(f_2, U_2)$. 
This proves our claim. We endow $P$ with the topology $\tau$. 

{\bf Step 2:} $\tau$ is Hausdorff: 
It is clear that the projection $q \: P \to M$ is continuous with respect 
to $\tau$. It therefore suffices to show that 
two different points $(m,g_1), (m,g_2) \in P$ can be separated by open subsets. 
If $U$ is a connected open neighborhood of $m$ on which 
$\alpha$ is integrable, we find a smooth function $f_1 \: U \to G$ with 
$f_1(m) = g_1$ and $\delta(f_1) = \alpha$. Now $f_2 := g_2g_1^{-1} 
f_1$ satisfies 
$f_2(m) = g_2$ and $\delta(f_2) = \alpha$. Hence 
$\Gamma(f_1, U)$ is an open neighborhood of $(m,g_1)$ and $\Gamma(f_2, U)$ an 
open neighborhood of $(m,g_2)$. If both intersect, there exists a 
point $x \in U$ with $f_1(x) = f_2(x)$, so that the connectedness of 
$U$ implies  $f_1 = f_2$, contradicting $g_1 \not= 
g_2$ (Lemma~\ref{lem:c.12b}). 
Hence $P$ is Hausdorff. 

{\bf Step 3:} The restrictions $q\res_{\Gamma(f,U)}$ are homeomorphisms. 
We know already that $q$ is continuous. 
Let $\sigma := (\id_U, f) \: U \to \Gamma(f,U) \subeq P$. 
Then, for each basic neighborhood 
$\Gamma(g,V)$, we have 
$$ \sigma^{-1}(\Gamma(g,V)) = \{ x \in U \cap V \: f(x) = g(x)\}. $$
In view of the uniqueness assertion of 
Lemma~\ref{lem:c.12b}, the latter set is open, so that $\sigma$ 
is continuous. Now $q \circ \sigma = \id_U$ implies that 
$q\res_{\Gamma(f,U)}$ is a homeomorphism onto $U$. 

{\bf Step 4:} The mapping $q \: P \to M$ is a covering:  
Let $x \in M$. Since $\alpha$ is 
integrable, there exists a connected open neighborhood 
$U_x$ of $x$ such that $\alpha\res_{U_x}$ is integrable. 
Then there exists for each $g \in G$ a smooth function 
$f_g \: U_x \to G$ with $f_g(x) = g$ and 
$\delta(f_g) = \alpha\res_U$. 
Now $q^{-1}(U) = U \times G = \bigcup_{g \in G} \Gamma(f_{g},U)$ 
is a disjoint union of open subsets of $P$, as we have seen in Step~$1$, 
and in Step~$3$ we have seen that it restricts to homeomorphisms 
on each of these sets. Hence $q$ is a covering. 

{\bf Step 5:} We conclude with Theorem~\ref{coverdiff}(a)
that $P$ carries a natural manifold structure for which $q$ is a local diffeomorphism. 
For this manifold 
structure, the function $F \: P \to G$ is smooth, because 
for each basic open subset $\Gamma(f,U)$, $U$ connected, the inverse of the corresponding 
restriction of $U$ is given by 
$$ \sigma_{(f,U)} := (q\res_{\Gamma(f,U)})^{-1} \: U \to \Gamma(f,U), \quad 
x \mapsto (x,f(x)), $$
and $F \circ \sigma_{(f,U)} = f$ is smooth. Moreover, 
$f \circ q = F\res_{\Gamma(f,U)}$ leads to 
$$ \delta(F) = q^*\delta(f) = q^*\alpha $$
on each set $\Gamma(f,U)$, hence on all of $P$.

Fix a point $m_0 \in M$. Then the 
connected component $\hat M$ of $(m_0, \be)$ in~$P$ is a connected covering manifold 
of $M$ with the required properties. 
If, in addition, $M$ is simply connected, then $q$ is a trivial covering, hence 
a diffeomorphism, and therefore $\alpha$ is integrable. 
\end{prf}

Let $\alpha \in \Omega^1(M,\g)$. If $\gamma \: I = [0,1] \to M$ 
is a piecewise smooth loop, then $\gamma^*\alpha \in \Omega^1(I,\g) 
\cong C^\infty(I,\g)$, so that 
$\evol_G(\gamma^*\alpha) \in G$ is defined because $G$ was assumed to 
be regular. 

\begin{lem} \mlabel{per-invar} Suppose that $M$ is connected and 
$\alpha \in \Omega^1(M,\g)$ is locally integrable. Pick $m_0 \in M$. 
Then the element $\evol_G(\gamma^*\alpha)$ does not change under homotopies with 
fixed endpoints and 
$$ \per^{m_0}_\alpha \: \pi_1(M,m_0) \to G,\quad 
[\gamma] \mapsto \evol_G(\gamma^*\alpha) $$
is a group homomorphism.   
\end{lem}

\begin{prf} Let $q_M \: \tilde M \to M$ denote a $1$-connected 
covering manifold of $M$ 
and choose a base point $\tilde m_0 \in \tilde M$ with 
$q_M(\tilde m_0) = m_0$. Then the $\g$-valued \break $1$-form 
$q_M^*\alpha$ on $\tilde M$ also satisfies the Maurer--Cartan 
equation, so that Proposition~\ref{prop:cov} 
implies the existence of a unique 
smooth function $\tilde f \: \tilde M \to G$ with 
$\delta(\tilde f) = q_M^*\alpha$ and $\tilde f(\tilde m_0) = \be$. 

We write 
$$ \sigma \: \pi_1(M, m_0) \times 
\tilde M \to \tilde M, \quad (d,m) \mapsto d.m  = \sigma_d(m) $$
for the left 
action of the fundamental group $\pi_1(M,m_0)$ on $\tilde M$. 
Then 
$$\delta(\tilde f \circ \sigma_d) = \sigma_d^*q_M^*\alpha = q_M^*\alpha
= \delta(\tilde f)$$ 
for each $d \in \pi_1(M,m_0)$ implies the existence of a function 
$$ \chi\: \pi_1(M,m_0) \to G 
\quad \hbox{ with } \quad 
\tilde f \circ \sigma_d = \chi(d) \cdot \tilde f, 
\quad d \in \pi_1(M,m_0). $$
For $d_1, d_2 \in \pi_1(M,m_0)$,  we then have 
$$ \tilde f \circ \sigma_{d_1 d_2} 
= \tilde f \circ \sigma_{d_1} \circ \sigma_{d_2} 
= (\chi(d_1) \cdot \tilde f) \circ \sigma_{d_2} 
= \chi(d_1) \cdot (\tilde f \circ \sigma_{d_2}) 
= \chi(d_1)\chi(d_2) \cdot \tilde f,$$
hence $\chi$ is a group homomorphism. 

We now pick a piecewise smooth $\tilde \gamma \: I \to \tilde M$ 
with $q_M \circ \tilde\gamma = \gamma$ and observe that 
$$ \delta(\tilde f \circ \tilde\gamma) 
= \tilde\gamma^*\delta(\tilde f) 
= \tilde\gamma^*q_M^*\alpha = \gamma^*\alpha, $$
so that 
$$ \chi([\gamma]) = \tilde f([\gamma].\tilde m_0) 
=  \tilde f(\tilde\gamma(1)) = \evol_G(\gamma^*\alpha). $$
This completes the proof. 
\end{prf} 

\begin{defn} \mlabel{def:period-homo}
For each locally integrable $\alpha \in \Omega^1(M,\g)$, the homomorphism 
$$ \per_\alpha := \per_\alpha^{m_0} \: \pi_1(M,m_0) \to G 
\quad \mbox{ with } \quad  \per_\alpha^{m_0}([\gamma]) = \evol_G(\gamma^*\alpha), $$
for each piecewise smooth loop $\gamma \: I \to M$ in $m_0$,  
is called the {\it period homomorphism of $\alpha$ with respect to $m_0$}. 
\index{period homomorphism of $1$-form} 
\end{defn}

The following theorem is one version of the Fundamental Theorem of calculus 
for functions with values in regular Lie groups. 

\begin{thm} 
[Fundamental Theorem for Lie group-valued functions]  
\index{Fundamental Theorem for Lie group-valued functions} 
\mlabel{thm-fundamental} 
Let $M$ be a smooth manifold, $G$ a Lie group with Lie algebra $\g$, and 
$\alpha \in \Omega^1(M,\g)$. Then the following assertions hold: 
\begin{enumerate}
\item[\rm(i)] If $G$ is regular, then $\alpha$ is locally integrable 
if and only if it satisfies the Maurer--Cartan equation. 
\item[\rm(ii)] If $\alpha$ is locally integrable, then
  \begin{enumerate}
\item[\rm(a)] if $M$ is $1$-connected, then $\alpha$ is integrable. 
\item[\rm(b)] if $M$ is connected and $m_0 \in M$, then $\per_\alpha^{m_0}$ 
vanishes if and only if $\alpha$ is integrable. 
  \end{enumerate}
\end{enumerate}
\end{thm}

\begin{prf} (i) If $\alpha$ is locally integrable, then Lemma~\ref{lem:MC} 
implies that it satisfies the Maurer--Cartan equation. If, converse, it 
satisfies the Maurer--Cartan equation, then Lemma~\ref{lem:locexist} implies its local 
integrability, provided $G$ is regular. 

(ii) (a) is Proposition~\ref{prop:cov}. 

(ii) (b) The function $\tilde f$ constructed in the proof of 
Lemma~\ref{per-invar} above factors through a smooth function on $M$ if and only if 
the period homomorphism is trivial. This implies (b). 
\end{prf}

\begin{cor} \mlabel{cor:fundamental} 
Let $M$ be a smooth connected manifold, $E$ a Mackey complete locally convex space, 
$\Gamma \subeq E$ a discrete subgroup, $G := E/\Gamma$, and $\alpha \in \Omega^1(M,E)$. 
Then the following assertions hold: 
\begin{enumerate}
\item[\rm(i)] $\alpha$ is locally integrable if and only if $\alpha$ is closed. 
\item[\rm(ii)] $\alpha$ is integrable to a smooth $G$-valued function 
if and only if all integrals of $\alpha$ over piecewise smooth loops are contained in $\Gamma$. 
\end{enumerate}
\end{cor}

\begin{prf} For the proof we only have to note that 
Proposition~\ref{prop:mc-reg} implies that $G$ is regular with 
$\evol_G(\xi) = q_G\big(\int_0^1 \xi(s)\, ds\big)$. 
\end{prf}

\begin{rem} (a) If $M$ is one-dimensional, 
then each $\g$-valued $2$-form on $M$ vanishes, so that 
$$ [\alpha, \beta] = 0 = \dd \alpha $$ 
for $\alpha, \beta \in \Omega^1(M,\g)$. Therefore 
all $1$-forms trivially solve the Maurer--Cartan equation. 

(b) If $M$ is a complex manifold, 
$G$ is a complex regular Lie group and 
$\alpha \in \Omega^1(M,\g)$ is a holomorphic $1$-form, then 
for any smooth function $f \: M \to G$ with $\delta(f) = \alpha$ 
the differential of $f$ is complex linear in each point, so that 
$f$ is holomorphic. Conversely, the logarithmic derivative of any holomorphic 
function $f$ is a holomorphic $1$-form. 

If, in addition, $M$ is a complex curve, i.e., a one-dimensional complex manifold, 
$G$ is a complex regular Lie group, and 
$\alpha \in \Omega^1(M,\g)$ is a holomorphic $1$-form, then 
the same argument as in (a) shows that the holomorphic $2$-forms 
$[\alpha, \alpha]$ and $\dd \alpha$ vanish, so that the Maurer--Cartan equation 
is automatically satisfied. 
\end{rem}


The following theorem was one of the central motivations to introduce 
the concept of a regular Lie group. 

\begin{thm}[Integration Theorem--regular case]\mlabel{thm3.2.11} 
\index{Integration Theorem!regular case}
 If $H$ is a regular Lie group, $G$ is a $1$-connected Lie group, 
and $\psi \: \L(G) \to \L(H)$ is a continuous homomorphism of Lie algebras,
then there exists a unique Lie group homomorphism $\phi \: G \to H$
with $\L(\phi) = \psi$. 
\end{thm}

\begin{prf} The uniqueness assertion follows from Proposition~\ref{propc.14} 
and does not require the regularity of~$H$. 

On $G$ we consider the smooth $\L(H)$-valued $1$-form given by 
$\alpha := \psi \circ \kappa_G$. That it satisfies the Maurer--Cartan equation 
follows from 
$$ \dd \alpha = \psi \circ \dd\kappa_G 
= -\shalf \psi \circ [\kappa_G, \kappa_G]  
= -\shalf [\psi \circ \kappa_G, \psi \circ \kappa_G] 
= -\shalf [\alpha, \alpha]. $$
Therefore the Fundamental Theorem implies the existence of a unique 
smooth function $\phi \: G \to H$ with $\delta(\phi) = \alpha$ and $\phi(\be_G) = \be_H$. 
In view of Proposition~\ref{prop:homocrit}, the function $\phi$ is a morphism of 
Lie groups, and we clearly have 
$\L(\phi) = \alpha_\be = \psi$. 
\end{prf}

\begin{cor} \mlabel{cor:to-existhom} If $G_1$ and $G_2$ are 
regular $1$-connected 
Lie groups with isomorphic Lie algebras, then $G_1$ and $G_2$ are isomorphic. 
\end{cor}

\begin{prf} Let $\psi\: \L(G_1) \to \L(G_2)$ be an isomorphism of Lie algebras. 
Then there exists a unique morphism of Lie groups 
$\phi \: G_1 \to G_2$ with $\L(\phi) = \psi$, and likewise 
there exists a morphism 
$\hat\phi \: G_2 \to G_1$ with $\L(\hat\phi) = \psi^{-1}$. 
Then $\L(\phi\circ \hat\phi) = \id_{\L(G_2)}$ and 
$\L(\hat\phi\circ \phi) = \id_{\L(G_1)}$ 
lead to $\hat \phi \circ \phi = \id_{G_1}$ and 
$\phi \circ \hat\phi = \id_{G_2}$ 
(Proposition~\ref{propc.14}). 
Thus $\phi$ is an isomorphism of Lie groups. 
\end{prf}

\begin{cor} \mlabel{cor:4.2.16} If $G$ is a $1$-connected regular Lie group 
with Lie algebra~$\g$, then 
\[ \L \: \Aut(G) \to \Aut(\g), \quad \phi \mapsto \L(\phi) \] 
is an isomorphism of groups.
\end{cor}

\begin{thm}
[Regular Lie groups with the same Lie algebra]\mlabel{thm:4.2.17} 
\index{Lie groups!regular, with same Lie algebra}Let $\g$ be the Lie algebra of  a regular Lie group. 
Then there exists a $1$-connected Lie group $G$ with Lie algebra~$\g$ 
which is unique up to isomorphism. 
Every connected regular Lie group with Lie algebra $\g$ is isomorphic 
to a quotient $G/\Gamma$, where $\Gamma\subeq Z(G)$ is a discrete subgroup. 
For two discrete subgroups $\Gamma_1, \Gamma_2 \subeq Z(G)$, 
the groups $G/\Gamma_1$ and $G/\Gamma_2$ are isomorphic if and only 
if there exists an automorphism $\phi \in \Aut(G)$ with 
$\phi(\Gamma_1) = \Gamma_2$. 
\end{thm}

\begin{prf} Proposition~\ref{prop:4.1.5} implies the existence of a 
$1$-connected regular Lie group $G$ with Lie algebra~$\g$ 
and the uniqueness up to isomorphism follows  from Corollary~\ref{cor:4.2.16}. 
For any connected regular Lie group $H$ whose Lie algebra 
is isomorphic to $\g$, Corollary~\ref{cor:4.2.16} further implies the 
the universal covering group $\tilde H$ is isomorphic to $G$. 
Now the remaining assertions follow from Theorem~\ref{thm:3.3.6}. 
\end{prf}

\section{More on Lie algebra-valued $1$-forms} 
\mlabel{sec:4.1}

In the following, we write $\Omega^1_{C^k}(M,E)$ for the
space of $E$-valued $p$-forms on $M$ of class $C^k$ ($k \in \N_0 \cup \{\infty\}$),
for a locally convex space~$E$,
i.e., the corresponding function $T(M) \to E$ is of class $C^k$ 
(Definition~\ref{def2.3.4a}). 
We topologize this space using the linear embedding
of $\Omega^1_{C^k}(M,E)$ into the locally convex space $C^k(TM,E)$,
which has closed image.

\begin{lem} \mlabel{lem:smoothprop} 
Let $M$ be a compact smooth manifold (possibly with boundary) 
and $K$ a Lie group with Lie algebra~$\fk$. 
Then, for each $r \in \N_0 \cup \{\infty\}$, the action of the Lie group 
$C^r(M,K)$ on $\Omega^1_{C^r}(M,\fk)$ by 
$(g,\alpha) \mapsto \Ad(g).\alpha$ is smooth. 
\end{lem}

\begin{prf} Assume that $d = \dim M$. 
Since each covering $(U_i)_{i \in I}$ of $M$ by compact submanifolds with 
boundary, diffeomorphic to $d$-dimensional balls, whose interiors define an 
atlas, yields a linear topological embedding
\[ \Omega^1_{C^r}(M,\fk) \into \prod_{i \in I} \Omega^1_{C^r}(U_i,\fk) 
\cong  \prod_{i \in I} C^r(U_i,\fk)^d\]
with closed image (cf.\ Lemma~\ref{family-restrictions}),
it suffices to show that the action of $C^r(M,K)$ on each space $C^r(U_i,\fk)$, given by 
\[ (g,f) \mapsto \Ad(g).f = \sigma_{\Ad} \circ (g,f), \quad 
(\Ad(g).f)(m) := \Ad(g(m))f(m)\]  
is smooth. This action factors through the Lie group morphisms 
$$ C^r(M,K) \onto C^r(U_i,K) $$ 
(Exercise~\ref{exer:4.1.7}), 
and for the Lie groups $G := C^r(U_i,K)$, it coincides with the adjoint action 
on $\L(G) \cong C^r(U_i,\fk)$, which is smooth. 
This proves the lemma. 
\end{prf} 
\begin{lem} \mlabel{lem:smoothprop2} 
Let $U$ be an open subset of a locally convex space $E$, 
$F$ a locally convex space, 
$M$ a compact manifold (with boundary), $r \in \N \cup \{ \infty\}$, and 
$\alpha \in \Omega^1(U,F)$. Then the map 
$$ C^r(M,U) \to \Omega^1_{C^{r-1}}(M,F), \quad 
f \mapsto f^*\alpha $$
is smooth if $C^r(M,U)$ is considered as an open subset of $C^r(M,E)$. 
\end{lem}

\begin{prf}
Both components of the map 
\begin{eqnarray*}
C^r(M,U) &\to& C^{r-1}(TM,TU) \cong C^{r-1}(TM,U) \times C^{r-1}(TM,E), \\
f &\mapsto & Tf = (f \circ \pi, \dd f) 
\end{eqnarray*}
are restrictions of continuous linear maps, hence smooth. 
Since $f^*\alpha = \alpha \circ Tf$, the assertion follows from the 
smoothness of $\alpha$, which implies that the map 
$$ \alpha_* \: C^{r-1}(TM,TU) \to C^{r-1}(TM, F), 
\quad h \mapsto \alpha \circ h 
$$
is smooth (Lemma~\ref{ClmapCk}) 
if we topologize 
$\Omega^1_{C^{r-1}}(M,F)$ as a closed subspace of $C^{r-1}(TM,F)$.  
\end{prf}

\begin{prop} \mlabel{prop:smooth-logder0} For any Lie group $K$ with Lie algebra 
$\fk$, any compact manifold $M$ 
(with boundary) and any $r \in \N \cup \{\infty\}$, the left logarithmic derivative 
$$ \delta \: C^r(M,K) \to  \Omega^1_{C^{r-1}}(M,\fk) $$ 
is a smooth map with respect to the Lie group structure on $C^r(M,K)$, 
and 
\begin{eqnarray}
  \label{eq:diff-delta}
T_\be(\delta) = \dd \: C^r(M,\fk) \to \Omega^1_{C^{r-1}}(M,\fk), \quad 
\xi \mapsto \dd \xi. 
\end{eqnarray}
\end{prop}

\begin{prf} From Lemma~\ref{lem:smoothprop}, we already know that the action 
of the Lie group $G := C^r(M,K)$ on $\Omega^1_{C^r}(M,\fk)$ by $f.\alpha := \Ad(f).\alpha$ 
is smooth. Since the inclusion map 
$C^r(M,K) \to C^{r-1}(M,K)$ is a smooth morphism of Lie groups  
(Exercise~\ref{exer:4.1.3b}), the action of $G$ on 
on $\Omega^1_{C^{r-1}}(M,\fk)$ is also smooth. 

The product rule 
$\delta(\eta\gamma) 
= \delta(\gamma) + \Ad(\gamma)^{-1}.\delta(\eta)$
means that $\delta$ is a right crossed homomorphism for the smooth action of 
$C^r(M,K)$ on $\Omega^1_{C^{r-1}}(M,\fk)$ (cf.\ Exercises~\ref{exer:3.4.2} 
and \ref{exer:4.1.4}). 
It therefore suffices to verify its smoothness in a neighborhood 
of the identity. 

Let $(\phi_K, U_K)$ be a chart of an identity neighborhood of $K$ with 
$T_\be(\phi_K) = \id_\fk$, so 
that $(\phi_G,U_G)$ with 
$$ U_G := \lfloor M,U_K\rfloor = \{ \gamma \in G \: \gamma(M) \subeq U_K \}, 
\quad \phi_G(\gamma) := \phi_K \circ \gamma $$
is a chart of an identity neighborhood of the Lie group $G$.  Then we have a map 
\begin{align*}
C^r(M, \phi_K(U_K)) &\to \Omega^1_{C^{r-1}}(M,\fk), \\ 
\xi &\mapsto \delta(\phi_K^{-1} \circ \xi) 
= (\phi_K^{-1} \circ \xi)^*\kappa_K 
= \xi^*(\phi_K^{-1})^*\kappa_K 
\end{align*}
whose smoothness follows from Lemma~\ref{lem:smoothprop2}. 

Now (\ref{eq:diff-delta}) follows from the fact that 
for each $m \in M$ and $v \in T_m(M)$, we have for 
$\beta := (\phi_K^{-1})^*\kappa_K$ 
the relations $\beta_0 = \id_\fk$ and 
\begin{align*}
\derat0 (t\xi)^*\beta v  
&= \derat0 \beta_{t\xi(m)}(t \dd\xi(m)v)
= \lim_{t \to 0} \beta_{t\xi(m)}\dd\xi(m)v \\
&= \beta_0 \dd\xi(m)v 
= \dd\xi(m)v. 
\qedhere\end{align*}
\end{prf}

For every Lie algebra $\g$, 
the closed commutator algebra 
$\oline D^1(\g) := \oline{[\g,\g]}$ is an ideal for which 
$\g_{\rm ab} := \g/\oline D^1(\g)$ is abelian and the quotient map 
$p \: \g \to \g_{\rm ab}$ is a continuous morphism of Lie algebras. 
If $G$ is a $1$-connected Lie group with Lie algebra $\g$, we can now ask 
for a homomorphism $p_G \: G \to (\g_{\rm ab},+)$ integrating 
$p$. If $G$ is locally exponential, then the existence of $p_G$ 
follows from Theorem~\ref{thm:int-thm} below, but in general this theorem 
does not apply. Even if $\g$ is Mackey complete, 
there is no guarantee that $\g_{\rm ab}$ is Mackey complete and 
Theorem~\ref{thm3.2.11} does also not apply. Nevertheless, we have: 

\begin{prop} If $G$ is a $1$-connected Lie group with a Mackey complete 
Lie algebra $\g$, then there exists a 
unique homomorphism $p_G \: G \to (\g_{\rm ab},+)$ for which 
$\L(p_G) = p \:  \g \to \g_{\rm ab}$ is the quotient map. 
If $\gamma \: [0,1] \to G$ is a smooth path starting in $\be$, then 
\begin{equation}
  \label{eq:int-abquot}
p_G(\gamma(1)) = p\Big(\int_0^1 \delta(\gamma)_t\, dt\Big).
\end{equation}
\end{prop}

\begin{prf} Let $E$ be the completion of $\g_{\rm ab}$. Then 
$E$ is in particular Mackey complete, so that 
Theorem~\ref{thm3.2.11} implies the existence of a unique 
Lie group morphism $p_G\: G \to E$ with $\L(p_G) = p$. 
If $\gamma \: [0,1] \to G$ is a smooth path starting in $\be$, then 
\[ \frac{d}{dt} p_G(\gamma(t)) = p(\delta(\gamma)_t),\] 
and since $\int_0^1 \delta(\gamma)_t\, dt$ exists in $\g$, the 
Fundamental Theorem implies \eqref{eq:int-abquot}. This shows that 
$p_G(G) \subeq \g_{\rm ab}$. To see that $p_G$ is actually smooth as a 
$\g_{\rm ab}$-valued map, let $\phi \: U \to G$ be a chart of $G$, where 
$U \subeq \g$ is an open $0$-neighborhood and $\phi(0) = \be$. 
Then 
\[ \sigma \: U \to C^\infty_*([0,1],G), \quad 
\sigma(x)(t) := \phi(tx) \] 
is a smooth map, and now the explicit formula 
\[ p_G(\phi(x)) = p\Big(\int_0^1 \delta(\sigma(x))_t\, dt\Big)\] 
and the smoothness of $\delta$ (Proposition~\ref{prop:smooth-logder0}) 
shows that $p_G$ is smooth on $\phi(U)$, hence smooth. 
\end{prf}

\begin{cor} If $G$ is a $1$-connected Lie group with a Mackey complete 
Lie algebra $\g$ and exponential fuction, then 
\[ \{ x \in \g \: \exp_G(x) = \be\} \subeq 
\oline D^1(\g).\] 
\end{cor} 

\begin{prf} We apply the preceding proposition with $\gamma(t) = \exp(tx)$ and 
obtain $0 = p_G(\exp x) = p(x)$, hence $x \in \ker p = \oline D^1(\g)$. 
\end{prf}

\begin{cor} \mlabel{cor:5.2.7} 
If $G$ is a $1$-connected Banach--Lie group, then 
\[ \{ x \in \g \: \exp_G(x) = \be\} \subeq 
\bigcap_{n \in \N} \oline D^n(\g).\] 
\end{cor} 

\begin{prf} In this case $p_G \: G \to \g_{\rm ab}$ is a quotient map and 
$\g_{\rm ab}$ is contractible. Therefore the kernel $N := \ker(q_G)$ 
satisfies 
\[ \pi_{k-1}(N) \cong \pi_{k}(\g_{\rm ab}) = \{0\} \quad \mbox{ for } \quad 
k = 1,2\] 
(Corollary~\ref{cor10.1.2b} and the subsequent remark). 
We conclude that the Lie group $N$ with Lie algebra 
$\oline D^1(\g)$ is also $1$-connected. We may therefore iterate the 
use of Corollary~\ref{cor:5.2.7} to obtain the assertion. 
\end{prf}

\section{Paths in groups and Lie algebras} 
 \mlabel{sec:4.2b} 

After the technical preparations in the preceding section, 
we now turn to some of the finer points of regular Lie groups. 
First we compare the smoothness requirements on $\evol_G$ and $\Evol_G$ 
and show that, for a regular Lie group~$G$ with Lie algebra $\g$, the map 
$\Evol_G \:  C^\infty(I,\g) \to C^\infty_*(I,G)$ is a diffeomorphism. 
We also determine the structure of connected abelian regular Lie groups 
and show that Banach--Lie groups are $C^0$-regular. 

\begin{rem} \mlabel{rem:reg1} 
Suppose that $G$ is a $C^k$-regular Lie group. 
Then the evolution map 
$$ \Evol_G \: C^k(I,\g) \to C^{k+1}(I,G), \quad \xi \mapsto \gamma_\xi $$
has values in the Lie group of $G$-valued $C^{k+1}$-curves. 

If $\Evol_G$ is smooth, then $\evol_G$ also is smooth 
because the evaluation map 
$$\ev \: C^{k+1}(I,G) \to G, \quad f \mapsto f(1)$$ 
is smooth (Exercise~\ref{exer:4.1.7} or Theorem~\ref{thmmfdmps}(c)). 
Suppose, conversely, that $G$ is regular. 
From Remark~\ref{rem:reg1a} we recall that 
$$ \Evol_G(\xi)(s) = \gamma_\xi(s) = \evol_G(S(\xi,s))
\quad \mbox{ for } \quad S(\xi,s)(t) = s \xi(st) $$ 
and that 
$$ S \: C^\infty(I,\g) \times I \to C^\infty(I,\g) $$
is smooth. 
Hence the map $(\xi,s) \mapsto \Evol_G(\xi)(s)$ is smooth, 
and this in turn implies the smoothness of 
$$ \Evol_G \:  C^\infty(I,\g) \to C^\infty(I,G) $$
(see Theorem~\ref{thmmfdmps}(a)).
\end{rem} 

The following proposition gives a natural necessary condition for regularity. 
Presently, the authors are not aware of any  non-regular Lie group $G$ 
modeled on a Mackey complete space. To find such an example is one of 
the central open problems of the theory. 

For the Lie group structure on the group $C^r_*(I,G)$ we refer to 
Exercise~\ref{exer:4.2.2}. 

\begin{prop} \mlabel{prop:reg-Mack} If 
$G$ is a regular Lie group, 
then its  Lie algebra $\g$ is Mackey complete and 
$$ \Evol_G \: C^\infty(I,\g) \to C^\infty_*(I,G) $$
is a diffeomorphism with $\Evol_G^{-1} = \delta$. 
\end{prop} 

\begin{prf} In view of Proposition~\ref{prop:smooth-logder0}, the map 
$$ \delta \: C^\infty(I,G) \to \Omega^1(I,\g) \cong C^\infty(I,\g) $$
is smooth with $T_\be(\delta)\xi = \xi'$. 
Since $G$ is regular, the evolution map 
$$ \Evol_G \: C^\infty(I,\g) \to C^\infty_*(I,G) $$
is also smooth (Remark~\ref{rem:reg1}) 
with 
\[ \Evol_G \circ \delta = \id_{C^\infty_*(I,G)} \quad \mbox{ and } \quad 
\delta \circ \Evol_G = \id_{C^\infty(I,\g)}.\] 
Hence, for each $\xi \in C^\infty(I,\g)$, the smooth curve 
$\eta := T_0(\Evol_G)\xi \in C^\infty(I,\g)$ satisfies 
$\eta' = T_\be(\delta)\eta = \xi$, and this implies that $\g$ is Mackey complete. 
\end{prf}

\begin{thm}[Connected regular abelian Lie groups] 
\mlabel{thm:reg-subcover}
\index{Lie groups!regular, connected abelian} 
A connected abelian Lie group $A$ is regular if and only if it 
is isomorphic to a quotient $\fa/\Gamma$ of a Mackey complete space 
$\fa$ by a discrete subgroup~$\Gamma$. 
\end{thm}

\begin{prf} Regularity is trivially inherited by all quotient groups 
$\fa/\Gamma$, where $\Gamma \subeq \fa$ is a discrete subgroup of a 
the Mackey complete space~$\fa$. 

If, conversely, $A$ is a connected regular abelian Lie group 
with Lie algebra~$\fa$, then 
Proposition~\ref{prop:reg-Mack} implies 
that $\fa$ is Mackey complete. Since $A$ has an exponential function 
$\exp_A \: \fa \to A$ (Remark~\ref{rem:reg1a}), 
Proposition~\ref{prop:reg-ab} now implies that 
it is of the form $\fa/\Gamma_A$, 
where $\Gamma_A = \ker(\exp_A) \subeq \fa$ is a discrete subgroup. 
\end{prf}

Now we are ready to verify that all Banach--Lie groups are 
$C^0$-regular.

\begin{thm} \mlabel{thm:banach-regular} 
\index{Banach--Lie groups!regularity} 
Each Banach--Lie group $G$ is 
$C^0$-regular. Moreover, for each $k \in \N_0$, 
$$ \Evol_G \: C^k(I,\g) \to C^{k+1}_*(I,G), \quad \xi \mapsto \gamma_\xi,  $$
is a  $C^\infty$-diffeomorphism. 
\end{thm}

\begin{prf}

  \nin{\bf Step 1:} Fix $k \in \N_0$. Then the logarithmic derivative
$$ \delta \: C^{k+1}(I,G) \to C^k(I,\g), \quad \gamma \mapsto 
\gamma^{-1}.\dot{\gamma} = \delta(\gamma) 
= \gamma^*\kappa_G $$
is smooth, as we saw in Proposition~\ref{prop:smooth-logder0}. 
Here we give an alternative, more direct argument. We have seen 
in Proposition~\ref{prop:tangentgrp} that for any  Lie group 
$G$, the tangent bundle $T(G)$ inherits a natural Lie group structure. 
We further have a map 
$$ C^{k+1}(I,G) \to C^{k}(I,T(G)), \quad \gamma \mapsto \dot{\gamma}, $$ 
which is easily seen to be a morphism of Lie groups (Exercise~\ref{exer:4.2.1}). 
We write $\pi \: T(G) \to G$ for the bundle projection. 
Then the Maurer--Cartan form $\kappa_G \: T(G) \to \g$ is the projection 
onto the second factor in the semidirect product 
$T(G) \cong G \ltimes \g$, satisfying the cocycle identity 
$$ \kappa_G(v \cdot w) = \kappa_G(w) + \Ad(\pi(w)^{-1})\kappa_G(v), $$
hence it is a right crossed homomorphism of Lie groups. It therefore induces a 
smooth crossed homomorphism of Lie groups 
$(\kappa_G)_* \: C^{k}(I,T(G)) \to C^k(I,\g)$ (Exercise~\ref{exer:4.1.4}). 
Now the assertion follows from $\delta(\gamma) = \kappa_G \circ \dot{\gamma}$. 

\nin{\bf Step 2:} We obtain a smooth bijective map 
$\delta \: C^{k+1}_*(I,G) \to C^k(I,\g).$ 
From Proposition~\ref{prop:smooth-logder0}, we further know that 
$T_\be(\delta)\xi = \xi',$
and this defines a topological isomorphism 
$C^{k+1}_*(I,\g) \to C^k(I,\g)$ (cf.\ Proposition~\ref{prop:reg-Mack}),  
whose inverse is the 
integration map 
$$ {\cal I} \: C^k(I,\g) \to C^{k+1}_*(I,\g), \quad {\cal I}(\xi)(t) := \int_0^t 
\xi(s)\, ds. $$
From the product rule 
$\delta(\gamma_1 \gamma_2) 
= \delta(\gamma_2) + \Ad(\gamma_2)^{-1}.\delta(\gamma_1),$
we further derive that 
$$ T_{\gamma}(\delta)(\xi.\gamma) = \Ad(\gamma)^{-1}.\xi', $$
and hence that all differentials $T_\gamma(\delta)$ are linear topological 
isomorphisms, so that $\delta$ is a $C^\infty$-diffeomorphism by 
the Inverse Function Theorem. Hence its inverse 
also is a $C^\infty$-diffeomorphism.
By construction, this inverse is the evolution map
$$ \Evol_G \: C^k(I,\g) \to C^{k+1}_*(I,G). $$
Hence $G$ is $C^k$-regular.
\end{prf}

\section{Tools to establish regularity of Lie groups}\label{tools-regularity}
We now provide results which frequently make it easy
to establish regularity for a Lie group at hand.
We shall apply the tools to mapping groups,
diffeomorphism groups,
box products of Lie groups, and direct limits
of finite-dimensional Lie groups.
The following characterization is useful.
We shall use the implication (b)$\Rightarrow$(a)
to deduce $C^0$-regularity in the examples.
\begin{thm}\label{thm-cts-0-reg}
Let $G$ be a Lie group with Lie algebra~$\g$.
Then {\rm(a)} and {\rm(b)} are equivalent:
\begin{description}[(D)]
\item[\rm(a)]
$G$ is $C^0$-regular.
\item[\rm(b)]
$G$ is $C^0$-semiregular, $\g$ is integral complete,
and $\evol\colon C(I,\g)\to G$
is continuous at~$0$.
\end{description}
%
%
%
\end{thm}
Subgroups inherit regularity
properties from the surrounding Lie group
in good cases. We shall only apply the following
result to Lie subgroups. But it holds, more generally,
for so-called strictly initial Lie subgroups (as in Definition~\ref{def-init-Lie}),\footnote{We require
to include manifolds $M$ \emph{with boundary} in
the definition of a $C^k$-initial
morphism of Lie groups.}
in the following form.
\begin{prop}\label{reg-equalizer}
Let $G$ be a Lie group, $\g$ be its Lie algebra,
and $k\in \N_0\cup\{\infty\}$.
Let $H$ be a subgroup of~$G$
admitting a Lie group structure turning
the inclusion map $\iota\colon H\to G$
both into a $C^\infty$-initial morphism
of Lie groups and a $C^\ell$-initial morphism
for some $\ell\in \N\cup\{\infty\}$ such that $\ell\leq k+1$.
Assume that there exist a set $J$, Lie groups $G_j$ for $j\in J$
and smooth group homomorphisms
$\alpha_j\colon G\to G_j$ and $\beta_j\colon G\to G_j$
such that
\begin{equation}
H=\{g\in G\colon (\forall j\in J)\; \alpha_j(g)=\beta_j(g)\}.
\end{equation}
Then the following holds:
\begin{description}[(D)]
\item[\rm(a)]
If $G$ is $C^k$-semiregular, then also $H$
is $C^k$-semiregular.
\item[\rm(b)]
If $G$ is $C^k$-regular, then also $H$
is $C^k$-regular.
\end{description}
\end{prop}
To prepare the proofs,
we now formulate and establish various auxiliary results.
As a byproduct, these show:
\begin{prop}
For $k\in \N_0\cup\{\infty\}$,
a $C^k$-semiregular Lie group $G$ is $C^k$-regular
if and only if the map $\Evol\colon C^k(I,\g)\to C^{k+1}(I,G)$
is smooth.
\end{prop}
\begin{prf}
Smoothness of $\Evol$ implies smoothness
of $\evol=\ev_1\circ \Evol$, as the evaluation mapping  $\ev_1\colon C^{k+1}(I,G)\to G$,
$\theta\mto\theta(1)$ is smooth (see Theorem~\ref{thmmfdmps}(c)); thus $G$ is $C^k$-regular.
If, conversely, $G$ is $C^k$-regular, then $\evol\colon C^k(I,\g)\to G$
is smooth and hence also $\Evol\colon C^k(I,\g)\to C^{k+1}(I,G)$
will be smooth, by part (f) of the following lemma.
\end{prf}
\begin{rem}
Note that, if a Lie group $G$ is $C^k$-semiregular
for some $k\in \N_0\cup\{\infty\}$, then
the map
\[
\delta\colon C^{k+1}_*(I,G)\to C^k(I,\g)
\]
is a bijection. We let $\odot$ be the group multiplication on $C^k(I,\g)$
turning $\delta$ into an isomorphism of groups.
For all $\gamma,\eta\in C^k(I,\g)$, we have
\[
(\gamma\odot \eta)(t)=\eta(t)+\Ad_{\Evol(\eta)(t)^{-1}}(\gamma(t))
\]
for all $t\in I$ (cf.\ Lemma~\ref{lem:c.12}) and thus
\[
\gamma\odot\eta =\eta +\Ad_{\Evol(\eta)^{-1}}(\gamma)
\]
if we write
\[
\Ad_{\theta}(\gamma):= C^k(I, \Ad^\wedge)(\theta,\gamma)
\]
for $\gamma$ as before and $\theta\in C^k(I,G)$,
using the adjoint action $\Ad^\wedge\colon G\times\g\to\g$, $(g,y)\mto \Ad_g(y)$
and the map
\[
C^k(I,\Ad^\wedge)\colon C^k(I,G)\times C^k(I,\g)\to C^k(I,\g),\;
(f,g)\mto \Ad^\wedge\circ \, (f,g)
\]
which is smooth by Proposition~\ref{variant-pushforwards}(b).
In particular, the right translation
\[
\rho_\eta\colon C^k(I,\g)\to C^k(I,\g),\;\; \gamma\mto \gamma\odot \eta
\]
is a continuous affine-linear map and hence smooth.
In fact, $\rho_\eta$ is a $C^\infty$-diffeomorphism,
as the smooth map $\rho_{\eta^{-1}}$ is its inverse.
Moreover,
\begin{equation}\label{line-part}
d\rho_\eta(\gamma,\gamma_1)=\Ad_{\Evol(\eta)^{-1}}(\gamma_1)
\end{equation}
is given by the linear part of the preceding affine-linear map.
These observations will be essential in the following discussions.
\end{rem}
\begin{rem}\label{pre-S}
We shall also use the map
\[
S\colon C^k(I,\g)\times I\to C^k(I,G)
\]
given by $S(\gamma,s)(t)=s\gamma(st)$.
Using the evaluation map $\ev\colon C^k(I,\g)\times I\to \g$,
the map
\[
S^\wedge\colon (C^k(I,\g)\times I)\times I\to \g,\quad ((\gamma,s),t)\mto s\gamma(st)
\]
can be written as
\[
S^\wedge((\gamma,s),t)=s\ev(\gamma,st).
\]
Since $\ev$ is $C^k$ (see Lemma~\ref{evaldiffprop}
and Proposition~\ref{CkvsCkk}),
this implies that $S^\wedge$ is $C^k$
and thus
$C^{0,k}$
(see Proposition~\ref{CkvsCkk}(a)),
whence $S$ is continuous as a map to $C^k(I,\g)$
using the Exponential Law (Theorem~\ref{explawCkell}(a)).
Being linear in the first argument,
the continuous map $S\colon C^k(I,\g)\times I\to C^k(I,\g)$
is $C^{\infty,0}$
(see Lemma~\ref{linCkell} and Corollary~\ref{CkellCellk}).
\end{rem}
\begin{rem}\label{henc-diffeom}
Let $G$ be a Lie group and $\phi\colon U\to V$
is $C^\infty$-diffeomorphism
from an open $\be$-neighborhood $U\sub G$ onto
an open $0$-neighbourhood $V\sub \g$
such that $\phi(\be)=0$ and $d\phi|_\g=\id_g$.
Let $\Sigma$
be the local addition on~$G$ provided by Proposition~\ref{exliegploa},
applied with
$\phi^{-1}$ in place of $\phi$.
Let $k\in \N_0\cup\{\infty\}$.
Write $f\colon I\to G$ for the constant map $t\mto\be$.
Then $O_f=C^k(I,V)$, $O_f'=C^k(I,U)$
and the chart $\phi_f$ for $C^k(I,G)$ constructed in the proof of
Theorem~\ref{thmmfdmps}(a) is the mapping
$C^k(I,\phi)\colon C^k(I,U)\to C^k(I,V)$, $\gamma\mto \phi\circ\gamma$.
\end{rem}
\begin{lem}\label{reg-in-steps}
Let $k\in \N_0\cup\{\infty\}$, $G$ be a $C^k$-semiregular Lie group 
and $\g$ be its Lie algebra.
Then the following holds.
\begin{description}[(D)]
\item[\rm(a)]
If $\evol\colon C^k(I,\g)\to G$
is continuous at~$0$,
then $\Evol\colon C^k(I,\g)\to C(I,G)$
is a continuous at~$0$.
\item[\rm(b)]
If the mapping $\Evol\colon C^k(I,\g)\to C(I,G)$
is continuous at~$0$,
then $\Evol\colon C^k(I,\g)\to C(I,G)$
is continuous.
\item[\rm(c)]
If $j\in \N_0\cup\{\infty\}$ with $j\leq k+1$
and $U\sub C^k(I,\g)$ is an open $0$-neighborhood
such that the restriction $\Evol|_U\colon U\to C^j(I,G)$
is $C^\ell$ for some $\ell\in \N_0\cup\{\infty\}$, then
$\Evol\colon C^k(I,\g)\to C^j(I,G)$
is $C^\ell$.
\item[\rm(d)]
If $j\in \N_0$ with $j<k$
and $\Evol \colon C^k(I,\g) \to C^j(I,G)$
is smooth, then $\Evol\colon C^k(I,\g)\to C^{j+1}(I,G)$
is smooth.
\item[\rm(e)]
If $\Evol|_U\colon U\to C(I,G)$ is smooth for
some open $0$-neighborhood $U\sub C^k(I,\g)$,
then $\Evol\colon C^k(I,\g)\to C^{k+1}(I,G)$
is smooth.
\item[\rm(f)]
If $\evol$ is smooth, then $\Evol\colon C^k(I,\g)\to C^{k+1}(I,G)$
is smooth.
\end{description}
\end{lem}
\begin{proof}
(a) Let $P$ be an identity neighborhood in $C(I,G)$.
After shrinking $P$, we may assume that $P=C(I,Q)$
for an identity neighborhood $Q\sub G$.
Since $\evol\colon C^k(I,\g)\to G$ is continuous at~$0$,
there exists an open $0$-neighborhood $V\sub C^k(I,\g)$
such that $\evol(V)\sub Q$.
As the map $S\colon C^k(I,\g)\times I\to C^k(I,\g)$
is continuous and $S(0,s)=0\in V$ for each $s\in I$, we see that
$\{0\}\times I$ is a subset of the open set $S^{-1}(V)$.
The Wallace Lemma (Lemma~\ref{Wallla})
provides an open $0$-neighborhood $W\sub C^k(I,\g)$
such that
\[
W\times I\sub S^{-1}(V).
\]
For each $\gamma\in W$, we then have
\[
\Evol(\gamma)(s)=\evol(S(\gamma,s))\in \evol(V)\sub Q
\]
for each $s\in I$ and thus $\Evol(\gamma)\in C(I,Q)=P$.

(b) For $\zeta\in C(I,G)$, let $\rho_\zeta\colon C(I,G)\to C(I,G)$
be the right translation map. Since
\[
\Evol\colon (C^k(I,\g),\odot)\to C^{k+1}_*(I,G)\sub C(I,G)
\]
is a group homomorphism, we have
$\Evol\circ \, \rho_\eta=\rho_\zeta\circ \Evol$ 
with $\zeta:=\Evol(\eta)$. Thus
$\Evol=\rho_\zeta\circ\Evol\circ \, \rho_\eta^{-1}$
is continuous at~$\eta$ if it is continuous at~$0$.

(c) For each $\eta\in C^k(I,\g)$, the set $U\eta:=\rho_\eta(U)$ is an open $\eta$-neighborhood,
since $\rho_\eta$ is a $C^\infty$-diffeomorphism and hence a homeomorphism.
Writing $\zeta:=\Evol(\eta)$, the formula
\[
\Evol|_{U\eta}=\rho_\zeta\circ \Evol|_U\circ \rho_\eta^{-1}|_{U\eta}
\]
shows that $\Evol|_{U \eta}$ is~$C^\ell$.

(d) Let $\phi\colon U\to V$ be a $C^\infty$-diffeomorphism from an
open $e$-neighborhood $U\sub G$ onto
an open $0$-neighborhood $V\sub \g$
such that $\phi(e)=0$ and $d\phi|_\g=\id_\g$.
Let $W_1\sub U$ be an open $e$-neighborhood
such that $W_1W_1\sub  U$ and $W:=\phi(W_1)$.
We then get a smooth mapping
\[
\mu\colon W\times W\to V,\;\,
(x,y) \mto \phi(\phi^{-1}(x)\phi^{-1}(y))
\]
which expresses the group multiplication
in the local chart. Using the partial differential in the second variable,
we obtain a smooth map
\[
\tau\colon W\times \g\to \g,\;\, (x,y)\mto d_2\mu(x,0,y).
\]
As a consequence, the map
\[
h:=C^j(I,\tau)\colon C^j(I,W)\times C^j(I,\g)\to C^j(I,\g),\;\,
(f,g)\mto \tau\circ (f,g)
\]
is smooth (see Proposition~\ref{superpo-cp}(a)).
Since $\Evol\colon C^k(I,\g)\to C^j(I,G)$
is continuous, the set $O:=\{\gamma\in C^k(I,\g)\colon \Evol(\gamma)\in C^j(I,W_1)\}$
is an open identity neighborhood in $C^k(I,\g)$.
As observed in Remark~\ref{henc-diffeom}, the map
\[
C^{j+1}(I,\phi)\colon C^{j+1}(I,W_1)\to C^j(I,W),\;\, f\mto \phi\circ f
\]
is a $C^\infty$-diffeomorphism. Hence $\Evol|_O$
will be smooth as a map to $C^{j+1}(I,G)$
(and thus also $\Evol$, by (c))
if we can show that $C^{j+1}(I,\phi)\circ \Evol|_O$
is smooth.
Consider the inclusion map $\incl\colon C^{j+1}(I,\g)\to C(I,\g)$
and the map $D\colon C^{j+1}(I,\g)\to C^j(I,\g)$, $f\mto f'$.
We know that the linear mapping
\[
(\incl,D)\colon C^{j+1}(I,\g)\to C(I,\g)\times C^j(I,\g)
\]
is a topological embedding with closed image.
Hence $C^{j+1}(I,\phi)\circ \Evol|_O$ will be smooth as a map to $C^{j+1}(I,\g)$
if it is smooth as a map to $C(I,\g)$
(which holds as $\Evol$ is smooth to $C(I,G)$)
and $D\circ C^{j+1}(I,\phi)\circ \Evol|_O$
is smooth.
We claim that the latter map coincides with
\[
h\circ ((C^j(I,\phi)\circ \Evol|_O)\times \incl),
\]
using the smooth map $\Evol\colon C^k(I,\g)\to C^j(I,G)$
on the right-hand side. Being a composition of smooth maps,
this map is smooth.

It remains to prove the claim.
For $\gamma\in O$, abbreviate $\eta:=\Evol(\gamma)$
and $\zeta:=\phi\circ \eta$.
Let $x\in W$ and $z:=\phi^{-1}(x)$.
Writing $\lambda_x:=\mu(x,\cdot)\colon W\to V$,
we have $\lambda_x=\phi\circ \lambda_z \circ\phi^{-1}|_W$
and thus
\[
\tau(x,y)=    d\lambda_x(0,y)=d\phi T\lambda_z T\phi(0,y)=
d\phi (zy).
\]
Given $t\in I$, we take $x:=\zeta(t)$ and $y:=\gamma(t)$.
Then $z=\eta(t)$.
We obtain that
\begin{eqnarray*}
h(\zeta,\gamma)(t)&=& (\tau\circ (\zeta,\gamma))(t)=
\tau(\zeta(t),\gamma(t))=d\phi(\eta(t)\gamma(t)) =d\phi(\dot{\eta}(t))\\
&=& \frac{d}{dt}(\phi\circ \eta)(t)=\zeta'(t).
\end{eqnarray*}
Thus $D(\phi\circ \eta)=D(\zeta)=h(\zeta,\gamma)$
with $\eta=\Evol(\gamma)$ and $\zeta=\phi\circ \eta$, 
as required.

(e) By (c), the map $\Evol\colon C^k(I,\g)\to C(I,G)$ is smooth.
As a consequence, $\Evol\colon C^k(I,\cG)\to C^j(I,G)$
is smooth for each $j\in \N_0$ such that $j\leq k+1$,
by (d) and induction. If $k$ is finite, we can choose $j=k+1$.
If $k=\infty$, we deduce using the Exponential Law (Theorem~\ref{explawCkell}(a))
that
\[
(\Evol)^\wedge\colon C^\infty(I,\g)\times I\to G,\;\,
(\gamma,t)\mto\Evol(\gamma)(t)
\]
is a $C^{\infty,j}$-map for each $j\in\N_0$
and hence $C^{\infty,\infty}$.
As a consequence, the map $\Evol\colon C^\infty(I,\g)\to C^\infty(I,G)$ is $C^\infty$, using
Theorem~\ref{explawCkell}(a) again.

(f) By Remark~\ref{rem:reg1a}(a), $\Evol(\gamma)(s)=\evol(S(\gamma,s))$ holds
for all $\gamma\in C^k(I,\g)$ and $s\in I$, whence
\[
\Evol^\wedge\colon C^k(I,\g)\times I\to G,\;\; (\gamma,s)\mto \Evol(\gamma)(s)
\]
is of the form $\Evol^\wedge=\evol\circ  \,S$.
Since $S$ is a $C^{\infty,0}$-map (see Remark~\ref{pre-S}) and $\evol$ is smooth, the Chain Rule shows that
$\Evol^\wedge=\evol\circ \, S$ is $C^{\infty,0}$ (see Proposition~\ref{Cklchainmfd}(a)).
Hence $\Evol\colon  C^k(I,\g)\to C(I,G)$
is smooth (by the Exponential Law, Theorem~\ref{thmmfdmps}(a)).
Then also $\Evol\colon C^k(I,\g)\to C^{k+1}(I,g)$
is smooth, by~(d).
\end{proof}
\begin{lem}\label{Evol-cts-dir}
Let~$G$ be a Lie group which is $C^0$-semiregular
and
such that the evaluation map $\Evol\colon C(I,\g)\to C(I,G)$
is continuous. If $\g$ is integral
complete, then $\Evol\colon C(I,\g)\to C(I,G)$
is a smooth map.
\end{lem}
\begin{prf}
Step 1. Let
$\phi\colon U\to V$ as well as $W_1$, $W$ and $\mu\colon W\times W\to V$
be as in the proof of Lemma~\ref{reg-in-steps}(d).
Since $\Evol\colon C(I,\g)\to C(I,G)$
is continuous, we find an open $0$-neighborhood
$O\sub C(I,\g)$ such that $\Evol(\gamma)\in C(I,U)$
for all $\gamma\in O$.
The map $C(I,\phi)\colon C(I,U)\to C(I,V)$, $\theta\mto\phi\circ\theta$
is a $C^\infty$-diffeomorphism (see Remark~\ref{henc-diffeom}).
Hence, if we can show that
\[
h:=C(I,\phi)\circ \Evol|_O\colon O\to C(I,F)
\]
is $C^1$, then also $\Evol|_O\colon O\to C(I,G)$
will be $C^1$ and thus $\Evol\colon C(I,\g)\to C(I,G)$
will be $C^1$, by Lemma~\ref{reg-in-steps}(c).
We claim that $h$ has all directional derivatives at~$0$
and $\alpha:=dh(0,\cdot)\colon C(I,\g)\to C(I,\g)$
is continuous.
If this is true,
then $h$ will be~$C^1$.
To see this, let $\gamma\in O$;
Since $\Evol$ is a group homomorphism from
$(C(I,\g),\odot)$ to $C(I,G)$,
we have
\[ \rho_{\Evol(\gamma)}\circ \Evol\circ \,\rho_{\gamma^{-1}}=\Evol \] 
and thus, on some open $\gamma$-neighborhood $O'\sub O$,
\begin{eqnarray*}
h|_{O'} &= &
C(I,\phi)\circ \rho_{\Evol(\gamma)} \circ C(I,\phi^{-1})\circ
C(I,\phi)\circ \Evol\circ \rho_{\gamma^{-1}}|_{O'}\\
&=&
(C(I,\phi)\circ \rho_{\Evol(\gamma)} \circ C(I,\phi^{-1}))\circ
h \circ \rho_{\gamma^{-1}}|_{O'}.
\end{eqnarray*}
Since $\rho_{\gamma^{-1}}$
is a continuous affine-linear map,
we deduce with Exercises~\ref{chain-affine}
and \ref{exc-chainpt} that the directional derivatives
$dh(\gamma,\gamma_1)$ at $\gamma$ exist for all $\gamma_1\in C(I,\g)$,
and are given by
\[
d\big(
C(I,\phi)\circ \rho_{\Evol(\gamma)} \circ C(I,\phi^{-1})
\big)
\big( 0,
\alpha(T\rho_{\gamma^{-1}}(\gamma,\gamma_1))
\big) ,
\]
which is a continuous function of $(\gamma,\gamma_1)$.\smallskip
 
To prove the claim, we identify $U$ with $V$ by means of the chart~$\phi$,
to ease notation. Thus $0$ now is the neutral element of~$G$.
Let $\alpha\colon C(I,\g)\to C(I,\g)$ be the continuous linear map given by
\[
\alpha(\gamma)(t)=\int_0^t\gamma(s)\, ds.
\]
Consider the smooth map $W\times \g\to\g$, $(x,y)\mto x.y:=d\lambda_x(0)(y)$,
where $\lambda_x\colon W\to V$, $z\mto xz$ is left translation by~$x$.
For $\gamma\in C(I,\g)$ and $0\not=\tau\in\R$ close to~$0$, we have 
\begin{eqnarray*}
\Evol(\tau\gamma)(t)&=& \int_0^t(\Evol(\tau\gamma))'(s)\,ds
=\int_0^t\Evol(\tau\gamma)(s).(\tau\gamma(s))\,ds\\
&=& \tau\int_0^t\Evol(\tau\gamma)(s).\gamma(s)\,ds,
\end{eqnarray*}
using the definition of $\Evol$. Hence
\[
\Delta_\tau := \frac{\Evol(\tau \gamma)-\Evol(0)}{\tau}-\alpha(\gamma)
=\frac{\Evol(\tau\gamma)}{\tau}-\alpha(\gamma)
\]
satisfies
\begin{eqnarray*}
\Delta_\tau(t) &=& 
\int_0^t\big(\Evol(\tau\gamma)(s).\gamma(s)-\gamma(s)\big)\,ds
= \int_0^t f(\Evol(\tau\gamma)(s),\gamma(s))\,ds, 
\end{eqnarray*}
with $f\colon W\times \g\to\g$, $f(x,y):=x.y-y$.
To see that $\Delta_\tau\to 0$ in $C(I,\g)$
(which completes the proof), let us show that, for each continuous seminorm
$p$ on $\g$, there exists $\rho>0$ such that 
 $[{-\rho},\rho] \gamma\subeq \Omega$ 
and
\begin{equation}\label{gives-diff-ev}
\|\Delta_\tau\|_{p,\infty}:=\sup_{t\in I}p(\Delta_\tau(t))\leq 1
\end{equation}
for all $\tau\in\R\setminus\{0\}$ such that $|\tau|\leq \rho$.
Now $f(x,0)=0$ and $f(0,y)=0$ for all $x\in W$ and $y\in\g$, whence there exist
continuous seminorms $q_1$ and $q_2$ on~$\g$ such that $B^{q_1}_1(0)\sub W$ and
and $p(f(x,y))\leq q_1(x)q_2(y)$ for all $x\in B^{q_1}_1(0)$
and $y\in\g$ (cf.\ Lemma~\ref{est-if-on-lines};
in view of the linearity of~$f$ in the second argument,
the latter need not be restricted to a $q_2$-ball around~$0$).
%
Since $q_2(\gamma([0,1]))$ is compact, we find $C\geq 1$ such that
$q_2(\gamma(s))\leq C$ for all $s\in I$.
Since~$\Evol$ is continuous at~$0$, we find a $0$-neighborhood $Q\sub\Omega$
such that
$\Evol(\eta)\in C(I,B^{q_1}_{1/C}(0))$ for all $\eta\in Q$.
Shrinking $\rho$, we may assume that $[{-\rho},\rho]\gamma\sub Q$.
Then
\[
p(\Delta_\tau(t))\leq\int_0^tp(f(\Evol(\tau\gamma)(s),\gamma(s))\,ds\leq 1\;\,
\mbox{for all $\tau
\in[{-\rho},\rho]\setminus\{0\}$,}
\]
as $p\big(f(\Evol(\tau\gamma)(s),\gamma(s))\big)\leq
q_1(\Evol(\tau\gamma)(s))q_2(\gamma(s))\leq\frac{1}{C}C\leq 1$.\smallskip

Step 2. We show by induction that $\Evol\colon C(I,\g)\to C(I,G)$
is $C^j$ for each $j\in \N$ and hence smooth.
The case $j=1$ was established in Step~2.
Since $\Evol$ as a group homomorphism from
$(C(I,\g),\odot)$ to $C(I,G)$,
we have $\rho_{\Evol(\gamma)}\circ \Evol\circ \,\rho_{\gamma^{-1}}=\Evol$
for each $\gamma\in C(I,\g)$.
Since $\Evol$ is $C^1$, taking tangent maps
we find that
$T\Evol
=T\rho_{\Evol(\gamma)}\circ T\Evol\circ 
\, T\rho_{\gamma^{-1}}$
and thus
\[
T\Evol(\gamma,\gamma_1)=
(T\rho_{\Evol(\gamma)}\circ T\Evol)(0,C(I,(\Ad)^\wedge)(\Evol(\gamma),\gamma_1))
\]
for all $\gamma,\gamma_1\in C(I,\g)$,
using~(\ref{line-part}).
If $\Evol$ is $C^j$ for some $j\in \N$, the formula shows that $T\Evol$ is $C^j$.
Now $\Evol$ being $C^1$ with $T\Evol$ a $C^j$-map,
it follows that $\Evol$ is a $C^{j+1}$-map.
\end{prf}
\noindent
\emph{Proof of Theorem}~\ref{thm-cts-0-reg}.
(a)$\Rightarrow$(b):
If $G$ is $C^0$-regular, then $G$ is $C^0$-semiregular
and $\evol=\ev_1\circ \Evol$ is continuous,
as the evaluation map $C^1(I,G)\to G$, $f\mto f(1)$
is continuous (Theorem~\ref{thmmfdmps}(c)).
To see that $\g$ is integral complete,
consider $\Evol$ as a map to $C^1(I,g)_*$.
Applying the Chain Rule to
\[
\delta\circ \Evol=\id_{C(I,\g)},
\]
we deduce from Proposition~\ref{prop:smooth-logder0}
that $\gamma=T_e\delta(T_0\Evol(\gamma)))=(T_0\Evol(\gamma))'$
for each $\gamma\in C^1(I,\g)_*$, whence $\eta:=T_0\Evol(\gamma)$
is an element of $C^1(I,\g)$ such that $\eta'=\gamma$.
Therefore $\eta(t)$ is the weak integral $\int_0^t\gamma(s)\,ds$
for all $t\in I$. Notably, $\int_0^1\gamma(s)\,ds$
exists in~$\g$.\smallskip

(b)$\Rightarrow$(a):
If $G$ is $C^0$-semiregular and $\evol\colon C(I,\g)\to G$ is continuous at~$0$,
then $\Evol\colon C(I,\g)\to C(I,G)$
is continuous by part (a) and (b) of Lemma~\ref{reg-in-steps}.
As we assume that $\g$ is integral complete,
Lemma~\ref{Evol-cts-dir} shows that $\Evol\colon C(I,\g)\to C(I,G)$
is smooth. Hence $\evol=\ev_1\circ \Evol$
is smooth, using the smoothness of the evaluation map $\ev_1\colon C(I,G)\to G$,
$f\mto f(1)$ (see Theorem~\ref{thmmfdmps}(c)).
Thus $G$ is $C^0$-regular.\qed\smallskip

To prepare the proof of Proposition~\ref{reg-equalizer},
we record an elementary fact.
\begin{lem}\label{morph-comp}
Let $\alpha\colon G\to H$ be a morphism of Lie groups.
Then $\delta(\alpha\circ\eta)=\Lie(\alpha)\circ \delta(\eta)$
for each $\eta\in C^1(I,G)$.
\end{lem}
\begin{prf}
Abbreviate $\zeta:=\alpha\circ \eta$.
For each $g\in G$, we have $\alpha\circ \lambda_g=\lambda_{\alpha(g)}\circ \alpha$
and hence $T\alpha\circ T\lambda_g=T\lambda_{\alpha(g)}\circ T\alpha$.
Applying this linear map to the tangent vector $\dot{\eta}(t)$ for $t\in I$
with $g:=\eta(t)^{-1}$, we get $T\alpha \delta(\eta)(t)= \alpha(\eta(t))^{-1}(T\alpha \dot{\eta}(t))=\zeta(t)^{-1}
\dot{\zeta}(t)=\delta(\zeta)(t)$.
\end{prf}
\noindent
\emph{Proof of Proposition}~\ref{reg-equalizer}.
(a) Let $\iota\colon H\to G$
be the inclusion map, which is a morphism of Lie groups
with injective tangent map $\L(\iota)\colon \h\to\g$
between $\h:=\L(G)$ and $\g:=\L(G)$.
Given $\gamma\in C^k(I,\h)$,
let $\eta:=\Evol_G(\L(\iota)\circ \gamma)\in C^{k+1}(I,G)$.
For all $j\in J$, we have
\[
\alpha_j\circ \iota=\beta_j\circ \iota
\]
and hence $\L(\alpha_j)\circ \L(\iota)=\L(\alpha_j\circ \iota)=\L(\beta_j\circ\iota)=
\L(\beta_j)\circ\L(\iota)$. As a consequence,
\[
\delta(\alpha_j\circ \eta)=\L(\alpha_j)\circ\delta(\eta)=\L(\alpha_j)\circ \L(\iota)\circ \gamma
=\b(\beta_j)\circ\L(\iota)\circ\gamma,
\]
whence $\alpha_j\circ\eta=\beta_j\circ\eta$.
Hence $\eta(t)\in H$ for all $t\in I$.
Then $\eta|^H\colon I\to H$ is a $C^\ell$-map
by $C^\ell$-initiality. Since
\[
\L(\iota)\circ\gamma=\delta(\eta)=\delta(\iota\circ \eta|^H)=\L(\iota)\circ \delta(\eta|^H),
\]
the injectivity of $L(\iota)$ implies that $\delta(\eta|^H)=\gamma$
and thus $\eta|^H=\Evol(\gamma)$.\smallskip

(b) We know from (a) that $\evol_H=(\evol_G\circ \, C^k(I,\L(\iota)))|^H$.
Using that $\evol_G\circ \, C^k(I,\L(\iota))\colon C^k(I,\h)\to G$
is a smooth mapping with image in~$H$,
using $C^\infty$-initiality we see that $\evol_H$ is smooth.\medskip\qed

Since $\L(\alpha_j\circ \,\eta)=\L(\alpha_j)\circ \, \delta(\eta)$
by Lemma~\ref{morph-comp}, we also infer
the following.
%
\begin{lem}\label{sep-points-evol}
Let $G$ be a Lie group, $\gamma\colon I\to \L(G)$
be a continuous path and $\eta\colon I\to G$ be $C^1$.
Assume that
\[
\delta(\alpha_j\circ \eta)=\L(\alpha_j)\circ \gamma\quad
\mbox{for all $\,j\in J$}
\]
for a family $(\alpha_j)_{j\in J}$
of morphisms $\alpha_j\colon G\to G_j$
of Lie groups such that the maps $\L(\alpha_j)\colon \L(G)\to\L(G_j)$
separate points on~$\L(G)$.
Then $\delta(\eta)=\gamma$.\qed
\end{lem}
We record a simple consequence.
The topic will be taken up again in Theorem~\ref{thm:reg-ext-prop}.
\begin{lem}\label{prod-reg}
Let $G_1$ and $G_2$ be Lie groups and $k\in \N_0\cup\{\infty\}$.
\begin{description}[(D)]
\item[\rm(a)]
If $G_1$ and $G_2$ are $C^k$-semiregular,
then the Lie group $G:=G_1\times G_2$ is also 
$C^k$-semiregular.
Abbreviate $\g:=\L(G)$.
If $\pr_j\colon G\to G_j$, $(g_1,g_2)\mto g_j$
denotes the projection for $j\in \{1,2\}$, then
\[
\Evol_G(\gamma)=(\Evol_{G_1}(\gamma_1),\Evol_{G_2}(\gamma_2))
\]
for all $\gamma\in C^k(I,\g)$,
with $\gamma_j:=\L(\pr_j)\circ\gamma$.
\item[\rm(b)]
If $G_1$ and $G_2$ are $C^k$-regular, then also
$G$ is $C^k$-regular.
\end{description}
\end{lem}
\begin{prf}
Abbreviate $\g_j:=\L(G_j)$ for $j\in\{1,2\}$.

(a) Given $\gamma\in C^k(I,\g)$, we have
$\gamma_j:=\L(\pr_j)\circ \gamma\in C^k(I,\g_j)$.
Then $\eta:=(\Evol(\gamma_1),\Evol(\gamma_2))\colon I\to G_1\times G_2=G$,
\[
t\mto(\Evol(\gamma_1)(t),\Evol(\gamma_2)(t))
\]
is a $C^k$-map such that $\delta(\pr_j\circ\eta)=\delta(\Evol(\gamma_j))=\gamma_j=\L(\pr_j)\circ\gamma$
and thus $\eta=\Evol(\gamma)$, by Lemma~\ref{sep-points-evol}.

(b) By (a), we have $\evol_G(\gamma)=(\evol_{G_1}(\L(\pr_1)\circ \gamma,\L(\pr_2)\circ\gamma)$
for all $\gamma\in C^k(I,\g)$. Since $\evol_{G_j}\colon C^k(I,\g_j)\to G_j$
is smooth and the linear mappings
\[
C^k(I,\L(\pr_j))\colon C^k(I,\g)\to C^k(I,\g_j),\quad \gamma\mto \L(\pr_j)\circ \gamma
\]
are continuous by Corollary~\ref{Cksuppo} and hence smooth,
we deduce that $\evol_G\colon C^k(I,\g)\to G$ is smooth.
\end{prf}
\begin{rem}\label{lie-produ}
By Lemma~\ref{tanprod}(b), the Lie algebra homomorphism
$\theta:=(\L(\pr_1),\L(\pr_2))\colon \L(G)\to \L(G_1)\times \L(G_2)$
is an isomorphism of topological vector spaces. If we identify $\L(G)$ with $\L(G_1)\times \L(G_2)$
via~$\theta$, then $\L(\pr_j)$ corresponds to the projection $\pi_j\colon \L(G_1)\times \L(G_2)\to\L(G_j)$
onto the $j$th component. Accordingly, we can identify $\gamma\in C^k(I,\g)$
with the pair $(\gamma_1,\gamma_2)\in C^k(I,\g_1)\times C^k(I,\g_2)$
and write $\Evol_G(\gamma)=(\Evol_{G_1}(\gamma_1),\Evol_{G_2}(\gamma_2))$
for $\gamma=(\gamma_1,\gamma_2)\in C^k(I,\g)$
in the situation of Lemma~\ref{prod-reg}(a).
\end{rem}
\begin{small}
\subsection*{Exercises for Chapter~\ref{ch:4}} 

\begin{exer} \mlabel{exer:reg-fixed-points} 
Let $G$ be a regular Lie group and 
$\Gamma \subeq \Aut(G)$ be a subgroup. Suppose that the subgroup 
\[ G^\Gamma := \{ g \in G \: (\forall \gamma \in \Gamma) \ \gamma(g) = g \} \] 
is a Lie subgroup. Show that $G^\Gamma$ is regular. 
\end{exer}

\begin{exer} \mlabel{exer:4.1.2b} Let $M$ be a finite-dimensional manifold, 
$E$ a locally convex space, and $\Omega^p_{C^r}(M,E)$ the 
space of $E$-valued $p$-forms of class $C^r$. Show that the 
map 
$$ \Phi \: \Omega^p_{C^r}(M,E) \to C^r(M,E)^{{\cal V}(M)^p}, 
\quad \Phi(\omega)(X_1,\ldots, X_p) := \omega(X_1,\ldots, X_p) $$
is a topological embedding of locally convex spaces. 
\end{exer}

\begin{exer} \mlabel{exer:4.1.3b} Let 
$r \leq s \in \N_0 \cup \{\infty\}$, 
$M$ be a compact manifold and $K$ a Lie group. 
Show that the inclusion $C^s(M,K) \into C^r(M,K)$ is a morphism of Lie groups. 
\end{exer}

\begin{exer} \mlabel{exer:4.1.4} 
Let $G$ and $N$ be a Lie groups and suppose that the homomorphism 
$\alpha \: G \to \Aut(N)$ defines a smooth action of $G$ on~$N$. 
Show that for each compact manifold $M$ and any 
$r \in \N_0 \cup \{\infty\}$, the 
prescription $(f.\xi)(m) := f(m).\xi(m)$ defines a smooth action of 
the Lie group $C^r(M,G)$ on the Lie group $C^r(M,N)$, and that 
any smooth crossed homomorphism 
$\chi \: G \to N$ induces a smooth crossed homomorphism 
$$ \chi_* \: C^k(M,G) \to C^k(M,N),  \quad f \mapsto \chi \circ f $$ 
(cf.\ Exercise~\ref{exer:3.4.2}). 
\end{exer}

\begin{exer} \mlabel{exer:4.1.7} Let $G$ be a Lie group and 
$\phi \: M_1 \to M_2$ a smooth map between compact manifolds with boundary. 
Show that, for each $r \in \N_0 \cup \{\infty\}$, the map 
$$ \phi^* \: C^r(M_2,K) \to C^r(M_1, K),\quad f \mapsto f \circ \phi $$
is a morphism of Lie groups. 
\end{exer}
  
\begin{exer} \mlabel{exer:4.2.1} Show that, 
for each Lie group $G$ and 
$k \in \N$, the map 
$$ C^k(I,G) \to C^{k-1}(I,TG), \quad \gamma \mapsto \gamma' $$
is a morphism of Lie groups if $TG$ is endowed with its canonical Lie group 
structure (Proposition~\ref{prop:tangentgrp}). 
\end{exer}

\begin{exer} \mlabel{exer:4.2.2} Let $K$ be a Lie group with Lie algebra $\fk$, 
$M$ be a finite-dimensional manifold 
and $X \subeq M$ a subset with compact closure. 
Show that for each $r \in \N_0 \cup \{\infty\}$, the group 
$$ C^r_X(M,K) := \{ f \in C^r(M,K)\: (\forall m \in M \setminus X)\ f(m) = \be\} $$
carries a natural Lie group structure with Lie algebra 
$C^r_X(M,\fk)$. Hint: Construct charts as for the case $X = M$, $M$ compact. 
\end{exer}
  
\begin{exer}\mlabel{exer:IV.2} Let $M$ be a smooth manifold, 
$H$ a regular Lie group and $\alpha  \in \Omega^1(M,\fh)$ a solution 
of the MC equation. Show that: 
\begin{enumerate} 
\item[\rm(1)] For any diffeomorphism $\phi \in \Diff(M)$, we have 
$$ \per_\alpha^{m_0}(\phi^*\alpha) = \per_\alpha^{\phi(m_0)}(\alpha) \circ 
\pi_1(\phi,m_0) \: \pi_1(M,m_0) \to H. $$
\item[\rm(2)] Let $G$ be a Lie group, acting smoothly on $M$ from the left 
by $g.m = \sigma_g(m)$ and 
also on $H$, resp., $\h$, by automorphisms $\rho_H(g)$, resp., $\rho_\h(g)$. 
We call $\alpha$ an equivariant form if 
$$ \sigma_g^* \alpha = \rho_\h(g) \circ \alpha $$
holds for each $g \in G$. Show that if $\alpha$ is equivariant, then 
$$ \rho_H(g) \circ \per_\alpha^{m_0}(\alpha) = \per_\alpha^{g.m_0}(\alpha) \circ 
\pi_1(\sigma_g,m_0) \: \pi_1(M,m_0) \to G. $$
If, in addition, $m_0$ is fixed by $G$ and $G$ is connected, then 
$$ \im(\per_\alpha^{m_0}) \subeq H^G. $$
\end{enumerate}
\end{exer} 

\begin{exer}
Show that Lemma~\ref{lem:laxuni}
also follows directly from
Lemma~\ref{Gspaceflow},
applied with the smooth right action
$\g\times G\to\g$, $(y,g)\mto\Ad_{g^{-1}}(y)$.
\end{exer}

\end{small}

\section{Notes and comments on Chapter~\ref{ch:4}} 

Theorem~\ref{thm3.2.11} is Theorem 8.1 in \cite{Mil84}. 
Actually this theorem was one of the guiding motivations 
to define the concept of regularity.
A version of the Fundamental Theorem (Theorem~\ref{thm-fundamental})
for $G$-valued $C^k$-maps can be found in~\cite{Alz13, Alz21}.
For a discussion of 
regular Lie groups in the convenient setting of 
Kriegl--Michor, we refer to \cite{KM97}.
Regularity properties for many classes of examples were established in
\cite{KM97b} and \cite{Gl12b}.
Stronger regularity properties, such as 
$L^1$-and $L^\infty$-regularity, have recently been investigated in~\cite{Gl15b}
(see also \cite{Nik21}, \cite{GH23}, and \cite{GSu23});
for $L^1$-regularity,
a smooth evolution $\Evol\colon L^1(I,\g)\to C(I,G)$
is required with values in the group of absolutely
continuous $G$-valued paths.
In the case of a Banach--Lie group,
it is even possible to replace $L^1(I,\g)$
with a Banach space of non-atomic $\g$-valued vector measures
and consider evolutions which are continuous $G$-valued functions
of bounded variation (see \cite{GSS23}).
Hanusch~\cite{Hn22} introduced many new ideas into regularity theory.
Notably, characterizations of $C^k$-regularity were obtained.
These involve certain completeness properties of the Lie group.
Completeness of infinite-dimensional Lie groups in their left uniformity was
studied in \cite{Gl17, Gl19a}.
Theorem~\ref{thm-cts-0-reg}
was obtained in \cite{Hn22},
partially based on previous work in~\cite{Gl12b}.
The Trotter product formula and commutator formula (as familiar in finite-dimensional Lie theory)
carry over to $L^\infty$-regular Lie groups in a strong form~\cite{Gl15b},
and also to $C^0$-regular Lie groups~\cite{Hn20}.
Every $C^0$-semiregular Fr\'{e}chet--Lie group is $C^0$-regular~\cite{Hn19}.

Regularity is a natural assumption that provides a good deal of methods to 
pass from the infinitesimal to the global level. This regularity 
concept is due to Milnor \cite{Mil84}. 
It weakens the $\mu$-regularity (see also \cite{Ne06a} introduced 
by Omori et al.\ in \cite{OMY82, OMY83a} (see \cite{KYM85} for a survey 
and \cite{Te02} for a discussion in the context of 
Fr\"olicher groups), 
but it is still strong enough for the essential 
Lie theoretic applications. 

\chapter{Locally exponential Lie groups} \mlabel{ch:5}

In this chapter, we turn to Lie groups $G$ with an exponential function 
which is 
well-behaved in the sense that it maps some $0$-neighborhood in $\L(G)$ 
diffeomorphically onto a $\be$-neighborhood in~$G$. We call 
such Lie groups {\it locally exponential}. 
As a consequence of the Inverse Function Theorem, 
all Banach--Lie groups and, 
in particular, all finite-dimensionale Lie groups are locally
exponential. The class of locally exponential Lie groups contains a huge variety of 
other types of interesting groups, such as 
the matrix groups $\GL_n(\cA)$ for a Mackey complete 
cia $\cA$ and the groups 
$C^\infty(M,K)$ of smooth maps of a compact manifold with values in any 
locally exponential Lie group $K$. Even the gauge groups $\Gau(P)$ of 
principal bundles over compact manifolds with locally exponential structure 
group $K$ are locally exponential. These types of groups, 
their integral subgroups 
(cf.\ Section~\ref{sec:6.2}), and their central extensions 
have been our principal motivation for a systematic 
study of locally exponential Lie groups. 
As we shall see in Chapter~\ref{ch:6}, 
the assumption of local exponentiality permits us to develop a 
powerful Lie theory for subgroups and there even is a characterization 
of those normal subgroups for which we may form Lie group quotients. 

There exist regular Fr\'echet--Lie groups 
which are not locally exponential. 
A prominent example is the diffeomorphisms group 
$\Diff(\bS^1)$ of the circle, but there are simpler examples 
such as the semidirect product $G = \C^\N \rtimes_\alpha \R$ 
with $\alpha(t)((z_n)) = (e^{itn}z_n)$ (Example~\ref{ex:bad-expfct}). 

A Lie group is called {\it exponential} if the exponential map 
$\exp_G \: \L(G) \to G$ 
is a global diffeomorphism, a requirement which is much stronger 
than local exponentiality. 
Important classes of such groups arise 
from simply connected finite-dimensional nilpotent Lie groups 
(\cite{HiNe12}), but the class of exponential 
Lie groups turns out to be closed under projective limits, so that it also contains 
the interesting class of simply connected groups associated to pro-nilpotent 
Lie algebras, i.e., projective limits of nilpotent Lie algebras.

{\bf Preliminaries:} Chapter~\ref{ch:3}

\section{Examples of locally exponential Lie groups} \mlabel{sec:5.1} 

In this first section, we present some of the most important 
classes of locally exponential Lie groups: groups of differentiable 
maps on compact manifolds and unit groups of continuous inverse algebras. We also provide some 
concrete examples that will be used repeatedly later on to 
display various kinds of bad behavior. In particular, we shall 
see an exponential Lie group which is analytic, but for which 
the corresponding global multiplication on the 
Lie algebra is not analytic. We also encounter 
an exponential Lie group 
whose Lie algebra contains closed ideals not invariant 
under the adjoint action. 

\subsection{Locally exponential Lie groups} 

\index{Lie group!locally exponential}
\begin{defn} \mlabel{def:5.1.1} We call a Lie group $G$ {\it locally exponential} if 
it has an exponential function $\exp_G \:  \L(G) \to G$ 
and there exists an open $0$-neighborhood $U \subeq \L(G)$ such that 
$\exp_G\res_U$ is a diffeomorphism onto an open
$\be$-neighborhood of~$G$. 

\index{Lie group!exponential}
A Lie group $G$ is called {\it exponential} if it has an exponential function 
which is a global diffeomorphism $\exp_G \: \L(G) \to G$. 
\end{defn}

\begin{prop} \mlabel{prop:5.1.3} Each Banach--Lie group is locally exponential. 
In particular, each finite-dimensional Lie group is locally exponential. 
\end{prop}

\begin{prf} In view of $T_0(\exp_G) = \id_{\L(G)}$, 
the assertion follows from the Inverse Function Theorem for Banach manifolds. 
\end{prf}

\begin{prop} \mlabel{prop:5.1.3b} If $E$ is a locally convex space 
and $\Gamma \subeq E$ a discrete subgroup,  then 
$E/\Gamma$ is a locally exponential Lie group. 
Conversely, for every connected abelian locally exponential Lie group 
$G$, the exponential function $\exp_G \: \g \to G$ is a smooth 
homomorphism and $\ker(\exp_G)$ is a discrete subgroup with 
$G \cong \g/\ker \exp_G$. 
\end{prop}

\begin{prf} That the quotients $E/\Gamma$ are locally exponential 
follows from the fact that, for these groups, $\L(E/\Gamma) \cong E$ 
and the quotient map $q \: E \to E/\Gamma$ 
is an exponential function (Example~\ref{ex:3.5.5}). 

If, conversely, $G$ is a connected abelian locally exponential 
Lie group, then its Lie algebra $\g$ is also abelian 
(Example~\ref{ex:abgrp}), so that $\exp_G \: (\g,+) \to G$ 
is a group homomorphism (Lemma~\ref{lem:4.1.4}). The local exponentiality 
of $G$ implies that $\exp_G$ is a covering morphism with discrete kernel, 
hence that $\exp_G$ factors through an isomorphism 
$\g/\ker(\exp_G) \to G$ of Lie groups. 
\end{prf}

\begin{exs} \mlabel{exs:5.1.4} 

(a) Unit groups of Mackey complete continuous inverse algebras are locally exponential 
(cf.~Theorem~\ref{thm:IV.1.11} and its proof): 
We recall from Example~\ref{ex:expo-cia} that the unit group $\cA^\times$ of 
each complex Mackey complete continuous inverse algebra has an analytic exponential function 
\[ \exp_\cA \: \cA \to \cA^\times,\]
defined by holomorphic functional calculus. 
The subset 
$$ \Omega :=  \{ a \in \cA \: \Spec(a) \cap ]-\infty,0] = \eset\} $$
is an open $\1$-neighborhood in $\cA^\times$, and with the complex 
logarithm function $\log \: \C \setminus ]-\infty,0] \to \C$ satisfying 
$\log(1) =0$, we get the holomorphic function 
\[  \log_\cA \: \Omega \to \cA, \quad x \mapsto \frac{1}{2\pi i} \oint_\Gamma 
\log(\zeta)(\zeta\1 - x)^{-1}\, d\zeta, \] 
where $\Gamma$ is a contour around $\Spec(x)$, lying in $\C \setminus 
]-\infty,0]$. Now the Spectral Mapping Theorem implies that 
$(\log_\cA \circ \exp_\cA)(x) = x$ holds if 
the spectral radius $\rho(x)$ is sufficiently small, and 
$(\exp_\cA \circ \log_\cA)(x) = x$ for each $x \in \Omega$. We conclude that 
the unit group $\cA^\times$ is locally exponential.

If $\cA$ is a real continuous inverse algebra, then one uses the fact that its complexification 
$\cA_\C$ is a continuous inverse algebra (Corollary~\ref{cor:8.1.5}) 
to see that $\log_{\cA_\C}(\Omega \cap \cA^\times) \subeq \cA$, and that 
$\log_\cA := \log_{\cA_\C} \res_\Omega$ is a smooth local inverse of 
$\exp_\cA = \exp_{\cA_\C} \res_\cA$.

(b) If $M$ is a compact manifold (possibly with boundary) 
and $K$ a Lie group with 
an exponential function, 
then Example~\ref{ex:mapgrp} implies that 
the Lie group $G := C^r(M,K)$, $r \in \N_0 \cup \{\infty\}$, has an 
exponential function given by 
$$ \exp_G \: C^r(M,\fk) \to C^r(M,K), \quad \xi \mapsto \exp_K \circ \xi. $$
If, in addition, 
$K$ is locally exponential, then the proof of 
Theorem~\ref{thm:mapgro-Lie} shows that 
we can use $\exp_G$ to define a chart of $G$, and hence that $G$ is locally exponential.

\index{gauge group}
(c) A larger class of groups of smooth maps 
is obtained as {\it gauge groups} of principal bundles. 
If $q \: P \to M$ is a smooth principal bundle with structure group $K$ 
and $\sigma \: P \times K \to P, (p,k) \mapsto \sigma_k(p) = p.k$ 
denotes the right action of $K$ on $P$, then 
\[  \Gau(P) := \{ \phi \in \Diff(P) \: q \circ \phi = q,\ 
(\forall k \in k)\ \phi \circ \sigma_k = \sigma_k \circ \phi\} \] 
is called the {\it gauge group of the bundle} and its elements 
are called {\it gauge transformations}.  
\index{gauge transformation}
In view of $q \circ \phi = q$, each gauge 
transformation $\phi$ can be written as 
$\phi(p) = p.f(p)$
for some smooth function $f \: P \to K$, and from 
$\phi \circ \sigma_k = \sigma_k \circ \phi$ we derive that 
$k f(p.k) = f(p)k$, i.e., 
\begin{eqnarray}
  \label{eq:5.1.1}
f(p.k) = k^{-1} f(p)k \quad \mbox{ for } \quad p \in P, k \in K.
\end{eqnarray}
Conversely, every smooth function $f \: P \to K$ satisfying (\ref{eq:5.1.1}) 
defines a gauge transformation by $\phi_f(p) := p.f(p)$. 
Moreover, 
\begin{eqnarray*}
\phi_{f_1}(\phi_{f_2}(p)) 
&= &\phi_{f_2}(p).f_1\big(\phi_{f_2}(p)\big)
= p.(f_2(p)f_1(p.f_2(p))) \cr
&=& p.(f_1(p)f_2(p)) = \phi_{f_1 f_2}(p) 
\end{eqnarray*}
implies that we have a group isomorphism 
$C^\infty(P,K)^K \to \Gau(P), f \mapsto \phi_f$, where 
$$ C^\infty(P,K)^K := \{ f \in C^\infty(P,K)\: 
(\forall p \in P)(\forall k \in K)\ f(p.k) = k^{-1} f(p)k\}. $$

We may therefore view $\Gau(P)$ as a subgroup of the group 
$C^\infty(P,K)$, endowed with the pointwise product. 
If the bundle $P$ is trivial, then there exists a smooth global section 
$\sigma \: M \to P$, and the map 
$$ C^\infty(P,K)^K \to C^\infty(M,K), \quad f \mapsto f \circ \sigma $$
is an isomorphism of groups. 

Results of Ch.~Wockel (\cite{Wo07a})
imply that (b) generalizes even to gauge 
groups: If $K$ is locally exponential and $M$ is compact, 
then $\Gau(P)$ carries a natural 
Lie group structure, turning it into a locally exponential Lie group
(cf.\ also \cite{SJ13}).
This is verified by showing that for 
\[  C^\infty(P,\fk)^K := \{ \xi \in C^\infty(P,\fk)\: 
(\forall k \in K)\ \xi \circ \sigma_k = \Ad(k)^{-1} \circ \xi\} \]
the map 
\[  \gau(P) := C^\infty(P,\fk)^K \to C^\infty(P,K)^K \cong \Gau(P), \quad 
\xi \mapsto \exp_K \circ \xi \]
is a local homeomorphism, 
and that it can be used to define a locally exponential 
Lie group structure on $\Gau(P)$ (see also \cite[Prop.~6.6]{OMY83a}). 

If $M$ is not compact, then the subgroup 
\[  \Gau_c(P) \trile \Gau(P) \] 
consisting of all elements $\phi \in \Gau(P)$ for which 
the image of the subset $\{ p \in P \: \phi(p) \not=p\}$ under~$q$ 
is relatively compact, carries a locally exponential Lie group 
structure if $K$ is locally exponential. 
\end{exs}

\subsection{Some exponential Lie groups} 

Next we take a closer look at exponential Lie groups. 
For any exponential Lie group $G$ with Lie algebra $\g$, we may use 
the diffeomorphism $\exp_G \: \g \to G$ 
to obtain a Lie group structure on $\g$ by 
$$ x * y := \exp_G^{-1}(\exp_G x\exp_G y). $$
Then $(\g,*)$ is a Lie group with $\L(\g,*) = \g$, and the multiplication 
on $\g$ satisfies 
\begin{eqnarray}
  \label{eq:5.1.2}
tx * sx = (t+s) x\quad \mbox{ for} \quad t,s \in \R, x \in \g.    
\end{eqnarray}
This implies that $\id_\g$ is an exponential 
function for this group. If, conversely, $*$ defines a Lie group structure on a Lie 
algebra $\g$ with neutral element $0$, then $\id_\g$ is an exponential function 
if and only if (\ref{eq:5.1.2}) is satisfied. 
In this sense exponential Lie groups can be identified with 
locally convex Lie algebras carrying a Lie group structure 
satisfying~(\ref{eq:5.1.2}), and Theorem~\ref{thm:5.2.13} below implies that the Lie algebra structure determines the exponential group 
structure on $\g$ completely. 
We call these Lie algebras {\it exponential}. 
\index{Lie algebra!exponential} 
In Theorem~\ref{thm:BCH-nil} below, we shall see in particular that 
each locally convex nilpotent Lie algebra is exponential. 

With the tools from Section~\ref{sec:3.4}, we can immediately 
show that 
each Lie algebra $\g$ carries at most one corresponding exponential Lie group structure. 
To this end we need: 

\begin{prop} \mlabel{prop:5.2.13b}  
Let $(\g_1, *)$ and $(\g_2, *)$ be two exponential Lie algebras and 
$\psi \: \g_1 \to \g_2$ be a continuous homomorphism of Lie algberas. 
Then \break $\psi \: (\g_1, *) \to (\g_2, *)$ is a homomorphism of Lie groups. 
\end{prop}

\begin{prf} In view of
\[ \delta(\psi)(x) = \kappa_{\g_2}(\psi(x)) \circ \psi 
= \int_0^1 e^{-t\ad \psi(x)}\, dt \circ \psi 
= \psi \circ \int_0^1 e^{-t\ad x}\, dt 
= \psi \circ \kappa_{\g_1}(x), \]
Proposition~\ref{prop:homocrit} implies that 
$\psi$ is a homomorphism of Lie groups. 
\end{prf}

As a consequence, we obtain the following uniqueness result: 
\begin{thm} 
[Uniqueness Theorem for exponential Lie groups]   \mlabel{thm:5.2.13}  
\index{exponential Lie groups!Uniqueness Theorem} 
Two exponential Lie group structures $(\g,*_1)$ and $(\g,*_2)$ 
on the same locally convex Lie algebra $\g$ coincide. 
\end{thm}

\begin{ex} \mlabel{ex:5.1.9} A typical example of a real exponential Lie group which is not nilpotent 
is the $2$-dimensional group $G = \R \rtimes \R$, with the multiplication
$$ (x,y) (x',y') = (x + e^y x', y + y'). $$
It is isomorphic to the identity component 
$\Aff(\R)_0$ of the affine group of the real line 
(cf.\ Example~\ref{ex:affgrp}). 
Its exponential map is given by 
$$ \exp_G \: \R^2 \to \R^2, \quad \exp_G(x,y) 
= \Big( \frac{e^y - 1}{y}x, y\Big)$$ 
(cf.\ Proposition~\ref{prop:semdir-exp}). 
It is a global diffeomorphism whose inverse is given by 
$$ \log_G \: \R^2 \to \R^2, \quad \log_G(x,y) 
= \Big( \frac{y}{e^y - 1}x,y\Big) $$
(Exercise~\ref{exer:5.6.1}). 
On the Lie algebra level, we have 
$$ [(x,y),(x',y')] = (yx' - y'x, 0), $$
so that the generators $e_1 = (1,0)$ and $e_2 = (0,1)$ 
satisfy $[e_2,e_1] = e_1$, which immediately 
implies that $\g = \L(G)$ is solvable but not nilpotent. 
\end{ex}

\begin{ex} \mlabel{ex:5.1.10} We describe an exponential 
Fr\'echet--Lie group $G$ which is
analytic, whose exponential function is analytic, but 
its inverse is not analytic (there is no analytic Inverse 
Function Theorem for Fr\'echet spaces), and the corresponding 
multiplication on $\g$ is not analytic. 

Let $H = \R \rtimes \R\cong \Aff(\R)_0$ 
be as in Example~\ref{ex:5.1.9} and 
put 
$$ G := H^\N = (\R\rtimes \R)^\N, $$
endowed with the direct product group structure. 
We endow $G$ with the manifold structure obtained by identifying it
with the Fr\'echet space $(\R^2)^\N$. 
This turns $G$ into an analytic manifold. 
As the power series defining the multiplication 
converges globally, the multiplication of $G$ is analytic, and the
same holds for the inversion map, which in each component is given by 
$$ (x,y)^{-1} = (-e^{-y}x ,-y). $$
In this sense, $G$ is an analytic Lie group. 

The exponential map of $G$ is given by 
$$ \exp_G(x,y)_n= \Big( \frac{e^{y_n} - 1}{y_n}x_n, y_n\Big), $$
and this map is analytic because the corresponding
power series converges globally. For the inverse function, we obtain 
$$ \exp_G^{-1}(x,y)_n = \Big( \frac{y_n}{e^{y_n} - 1}x_n, y_n\Big), $$
but this map is not analytic because the power series of the real
analytic function 
$$ y \mapsto \frac{y}{e^y - 1} $$
converges only on the interval from $-2\pi$ to $2\pi$, and the product
of infinitely many such intervals is not an open subset in $\R^\N$. 

The multiplication $*$ on the Lie algebra $\fh$ obtained from the exponential 
chart is given by 
\begin{eqnarray*}
(x,y) * (x',y') 
&=& \log(\exp_G(x,y)\exp_G(x',y')) \\
&=& \log\Big( \frac{e^y -1}{y}x + e^y\frac{e^{y'} -1}{ y'}x', y + y'\Big) \\ 
&=& \Big( \frac{y + y'}{e^{y + y'} - 1}
\Big(\frac{e^y -1}{y}x + e^y \frac{e^{y'}-1}{y'}x'\Big), y + y'\Big)
\end{eqnarray*}
and in particular 
$$ (0,y) * (1,0) = \log(e^y, y) = \Big(\frac{y e^y}{e^y-1}, y\Big)
= \Big(\frac{y}{1 - e^{-y}}, y\Big).$$ 
Therefore the argument from above shows that the multiplication 
on the product Lie algebra $\g \cong \fh^\N$ is not analytic. 
\end{ex}

\begin{rem} Let $E$ be a Mackey complete locally convex space 
and $\alpha \: \R \to \GL(E)$ be a smooth $\R$-action on $E$, so that 
we can form the semidirect product Lie group 
$G := E \rtimes_\alpha \R$. 
The group $G$ has an exponential 
function given by 
$$ \exp_G(v,t) := (\beta(t)v, t) \quad \mbox{ with } \quad 
\beta(t) = \int_0^1 \alpha_{st}\, ds $$
(Proposition~\ref{prop:semdir-exp}). 
From this formula, it follows immediately that $\exp_G$ is bijective if and only if, for each $t \in \R$, 
the operator $\beta(t) \in \cL(E)$ is bijective, and that 
in this case the inverse of $\exp_G$ is given by 
$$ \log_G(v,t) = (\beta(t)^{-1}v, t). $$
Therefore $\exp_G$ is a diffeomorphism if and only if $\beta(\R) \subeq \GL(E)$ 
and the curve $\beta \: \R \to \GL(E)$ is smooth in the sense that 
$$ \tilde\beta \: \R \times E \to E \times E, \quad 
(t,v) \mapsto (\beta(t)v, \beta(t)^{-1}v) $$
is smooth. In view of Proposition~\ref{prop:e.2.10}(2), 
the map $\tilde\beta_1 \: (t,v) \mapsto \beta(t)v$ is smooth, and 
Theorem~\ref{thm:autsmooth-lin} implies that 
the map $\tilde\beta_2 \: (t,v) \mapsto \beta(t)^{-1}v$ is smooth 
if it is continuous. 
\end{rem}

We now turn to a concrete example of an exponential Lie algebra of this form. 
 
\begin{ex} \mlabel{ex:5.1.11} Let 
$E \subeq C^\infty(\R,\C)$ be the closed subspace of $1$-periodic functions, 
$\mu \in \R^\times$, and 
consider the homomorphism $\alpha \: \R \to \GL(E)$ given by 
$$ (\alpha_t f)(x) := e^{\mu t} f(x + t). $$
That the corresponding $\R$-action on $E$ is smooth follows from the smoothness of 
the map 
$(t,f,x) \mapsto e^{\mu t} f(x + t)$ 
(Proposition~\ref{prop:cartes-closed}), which in turn follows from 
the smoothness of the evaluation map of $C^\infty(\R,\C)$ 
(Lemma~\ref{lem:smooth-eval-vec}). 
Understanding the exponential in the sense of 
Remark~\ref{rem:e.2.2}, we have $\alpha_t = e^{tD}$ for $Df = \mu f + f'$. 
The functions $e_n(x) = e^{2\pi i nx}$, $n \in \Z$, are eigenfunctions 
of $D$ and $\alpha$ with 
$$ D e_n = (\mu + 2\pi i n) e_n 
\quad \mbox{ and } \quad \alpha_t e_n = e^{t(\mu + 2 \pi i n)} e_n. $$
For the operators $\beta(t)$, this leads to 
$\beta(t)e_n = \beta_n(t) e_n$ with 
$$\beta_n(t) 
= \int_0^1 e^{st(\mu + 2 \pi i n)}\, ds 
= \left\{ 
  \begin{array}{cl} 
    1 & \mbox{for $t = 0$} \\           
 \frac{e^{t(\mu + 2 \pi i n)}-1}{t(\mu + 2 \pi i n)} & \mbox{ otherwise.} 
\end{array} \right. $$

The functions $\beta_n \: \R \to \C$ are analytic, hence in particular smooth. 
To see that the corresponding operators $\beta(t)$ are invertible, 
we consider the Fourier transformed picture of $E$, showing that 
$$ E = \Big\{ \sum_{n \in \Z} a_n e_n \: (\forall k \in\N) \lim_{n \to \infty} n^k a_n =0\Big\} $$
whose topology is defined by the seminorms 
$$ p_k\Big(\sum_n a_n e_n\Big) := \sup_{n \in \Z} (|n|+ 1)^k |a_n|. $$ 
The corresponding Banach spaces 
$$E_k := \{ a =(a_n)_{n \in \N} \in \C^\N \: p_k(a) < \infty\}$$ 
satisfy $E \cong \prolim E_k.$ 
We write $q_k \:E \to E_k$ for the canonical map, which is continuous and injective. 

For a sequence $c = (c_n)_{n \in \Z}$, the prescription 
$T_c(e_n) := c_n e_n$ extends to a continuous linear operator on $E$ 
if and only if, for each $k \in \N$, there exists some $m \in \N$ and $C_{k,m} > 0$ with 
$$ p_k\Big(\sum_n a_n c_n e_n\Big) \leq C_{k,m} p_m\Big(\sum_n a_n e_n\Big). $$
Applying it to the basis elements $e_n$, this inequality leads to 
$$  |c_n| (1 + |n|)^k \leq C_{k,m} (1 + |n|)^m, $$
which implies that $(c_n)$ is of polynomial growth. 
If, conversely, $(c_n)$ is of polynomial growth: 
$|c_n| \leq M(1 + |n|^m)$ for some $m \in \N$, $M > 0$, then we have for each $k \in \N$
$$  |c_n| |n|^k \leq M(1 + |n|^m) |n|^k \leq M' |n|^{m+k} $$
for some $M' > 0$, and this implies that $T_c$ extends to a bounded operator on $E$.
The operator $T_c$ is invertible if each $c_n$ is non-zero and the sequence 
$c^{-1}  = (c_n^{-1})_{n \in \Z}$ also is of polynomial growth. 

To see that $\beta(t) \in \GL(E)$, we therefore have 
to show that each sequence $(\beta_n(t))_{n \in \Z}$ and its inverse is of polynomial 
growth. For $t \not=0$ with $r \leq |t| \leq R$ we have 
$$ |\beta_n(t)| 
=  \frac{|e^{t(\mu + 2 \pi i n)}-1|}{|t(\mu + 2 \pi i n)|}
\leq  \frac{e^{t\mu} + 1}{|t|\mu} 
\leq  \frac{e^{R\mu} + 1}{r\mu} $$
and 
$$ |\beta_n(t)^{-1}| 
=  \frac{|t(\mu + 2 \pi i n)|}{|e^{t(\mu + 2 \pi i n)}-1|}
\leq  \frac{|t|(|\mu| + 2 \pi |n|)}{|e^{t\mu} - 1|} 
\leq  \frac{R(|\mu| + 2 \pi |n|)}{\min(|e^{R\mu} - 1|,|e^{r\mu}-1|)}, $$
which implies that $\beta(t) \in \GL(E)$ and that we have uniform 
estimates for $r \leq |t| \leq R$, showing that $\beta(t)$ extends to a continuous 
linear operator in $\GL(E_k)$, $k \in \N$, and that 
$\beta(t)^{-1}$ extends 
to a continuous linear operator in ${\cal L}(E_{k+1},E_k)$. 
To verify the continuity of the map $\hat\beta_2$, it suffices to verify, 
 for each $k$, the continuity of the corresponding map 
$q_k \circ \hat\beta_2 \: \R \times E \to E_k,$
which, as we have just observed, factors through a map 
$\gamma \: \R \times E_{k+1} \to E_k.$
For  
$\delta(t) := \beta(t)^{-1}\res_{E_{k+1}}^{E_k}$
we have $\gamma(t,v) = \delta(t)v$. 
The estimates above show that $\delta$ is locally bounded. On the other hand, 
the maps $t \mapsto \delta(t)e_n$ are continuous for each $n \in \Z$, and 
the $e_n$ span a dense subspace of $E_{k+1}$. From that it follows easily that 
$\gamma$ is continuous (Exercise~\ref{exer:5.6.11}). 
This in turn entails that $\tilde\beta_2$ is continuous, 
and hence that $\tilde\beta_2$ is smooth (Theorem~\ref{thm:autsmooth-lin}). 
We conclude that $\beta \: \R \to \GL(E)$ is smooth, so that 
$G$ is an exponential Lie group. 
\end{ex}

\begin{small}
\subsection*{Exercises for Section~\ref{sec:5.1}} 

\begin{exer} \mlabel{exer:5.6.1} 
Let $G$ be the $2$-dimensional group obtained by endowing $\R^2$ with the 
multiplication 
$$ (x,y) (x',y') = (x + e^y x', y + y'). $$
Show that
\begin{enumerate}
\item[\rm(a)] $G \cong \Aff(\R)_+$, the group of orientation preserving 
affine maps of $\R$. 
\item[\rm(b)] $G$ is exponential. Its exponential function is given by 
$$ \exp_G \: \R^2 \to \R^2, \quad \exp_G(x,y) = \Big( \frac{e^y - 1}{y}x, y\Big),
$$
and that its inverse is given by 
$$ \log_G \: \R^2 \to \R^2, \quad \log(x,y) = \Big( \frac{y}{e^y - 1}x,y\Big). $$
\end{enumerate}
\end{exer} 

\begin{exer} \mlabel{exer:5.6.1b} 
Let $G$ be the $3$-dimensional group obtained by endowing $\C \times \R$ with the 
multiplication 
$$ (x,y) (x',y') = (x + e^{iy} x', y + y'). $$
Show that this is a $3$-dimensional solvable Lie group which is not 
exponential. Its exponential function is given by 
$$ \exp_G \: \C \times \R \to \C \times \R, \quad \exp_G(x,y) 
= \Big( \frac{e^{iy} - 1}{iy}x, y\Big). $$
Determine the image of the exponential function, its singular points, 
and all pairs $((x,y), (x',y'))$ with 
$\exp_G(x,y) = \exp_G(x',y')$. 

The group $G$ occurs quite naturally in geometry as the universal covering 
of the identity component $\Mot_2(\R)_0$ of the 
euclidean motion group $\Mot_2(\R)$, consisting of all (affine) isometries 
$\R^2 \to \R^2$ with respect to the euclidean metric. The covering map 
$$ q \: G \to \Mot_2(\R)_0 \quad \mbox{ is given by }\quad 
q(x,y)z := x + e^{iy} z, $$
if we identify the euclidean plane $\R^2$ in the canonical way with $\C$.  
\end{exer} 

\begin{exer} \mlabel{exer:5.6.11} Let $X$ be a topological space, 
$E,F$ Banach spaces, and consider a locally bounded 
map $\phi \: X \to {\cal L}(E,F)$. 
Suppose that the set of all $v \in E$ for which the map 
$\phi^v \: X \to F, x \mapsto \phi(x)v$ is continuous is dense in $E$. 
Then the map $\hat\phi \: X \times E \to F, (x,v) \mapsto \phi(x)(v)$ is continuous. \\
Hint: $\|\phi(x)(v) - \phi(x')(v')\| \leq \|\phi^v(x) - \phi^v(x')\| + 
\|\phi(x')\|\|v - v'\|$.   
\end{exer}
  
\end{small}

\section{Homomorphisms of locally exponential Lie groups} \mlabel{sec:5.5} 

The main result of the present section is 
Theorem~\ref{thm:5.4.4} which asserts  that 
continuous homomorphisms of locally exponential Lie groups are smooth 
and that any homomorphism $\L(G) \to \L(H)$ of the Lie algebra 
$\L(G)$ of a \break {$1$-connected} locally exponential Lie group to the Lie algebra 
of a group $H$ with exponential function can be integrated 
to a Lie group homomorphism. Since locally exponential Lie groups need not be 
regular, these results complement the corresponding results on 
regular Lie groups such as Theorem~\ref{thm3.2.11}. 

\subsection{Automatic Smoothness} 

We start with the case of one-parameter groups. 

\begin{lem} \mlabel{lem:5.4.1} 
If $G$ is locally exponential and 
$\gamma \: \R \to G$ a continuous one-parameter group, then 
$\gamma$ is smooth and there exists a unique $x \in \L(G)$ with 
$\gamma(t) = \exp_G(tx)$ for all $t \in \R$. 
\end{lem}

\begin{prf} Let $U = - U$ be a convex
$0$-neighborhood in $\g$ for which $\exp_G\res_U$ is a diffeomorphism
onto an open subset of $G$ and put
$U_1 := \frac{1}{2} U$.
Since $\gamma$ is continuous in $0$, there exists an $\eps > 0$ such that
$\gamma([-\eps,\eps]) \subeq  \exp_G(U_1)$.
Then $\alpha(t) := (\exp_G\res_U)^{-1}(\gamma(t))$ defines a continuous curve
$\alpha \: [-\eps,\eps] \to U_1$
with $\exp(\alpha(t)) = \gamma(t)$ for $|t| \leq \eps$.
For any such $t$, we then have
$$ \exp_G\big(2\alpha({\textstyle \frac t2})\big) = \exp_G(\alpha({\textstyle \frac t2})\big)^2
= \gamma({\textstyle \frac t2})^2 = \gamma(t)
= \exp_G(\alpha(t)), $$
so that the injectivity of $\exp_G$ on $U$ yields
$\alpha(\textstyle \frac{t}{2}) 
= \textstyle \frac{1}{2}\alpha(t)$ for $|t| \leq \eps.$
Inductively, we thus obtain
\begin{equation}
  \label{eq:2.2.3}
\alpha(\textstyle \frac{t}{2^k}) = \textstyle \frac{1}{2^k}\alpha(t) \quad \hbox{ for } \quad |t| \leq \eps,
k \in \N.
\end{equation}
In particular, we obtain
$ \alpha(t) \in \frac{1}{2^k} U_1$ for $|t| \leq \frac{\eps}{2^k}.$
For $n \in \Z$ with $|n| \leq 2^k$ and $|t| \leq \frac{\eps}{2^k}$ we now have
$|nt| \leq \eps$, $n\alpha(t) \in \frac{n}{2^k} U_1 \subeq U_1$, and
$$ \exp_G(n\alpha(t))= \gamma(t)^n = \gamma(nt) = \exp_G(\alpha(nt)). $$
Therefore the injectivity of $\exp_G$ on $U_1$ yields
\begin{equation}
  \label{eq:2.2.4}
\alpha(nt) = n \alpha(t) \quad \hbox{ for } \quad n \leq 2^k, |t| \leq \frac{\eps}{2^k}.
\end{equation}
Combining (\ref{eq:2.2.3}) and (\ref{eq:2.2.4}), leads to
$$ \alpha(\textstyle \frac{n}{2^k}t) = \textstyle \frac{n}{2^k}\alpha(t) \quad \hbox{ for } \quad |t| \leq \eps,
k \in \N, |n| \leq 2^k. $$
Since the set of all numbers $\frac{nt}{2^k}$, $n \in \Z$, $k \in \N$,
$|n| \leq 2^k$, is dense
in the interval $[-t,t]$, the continuity of $\alpha$ implies that
$$ \alpha(t) = \frac{t}{\eps} \alpha(\eps) \quad \hbox{ for } \quad  |t| \leq \eps. $$
In particular, $\alpha$ is smooth and of the form $\alpha(t) = tx$ for some
$x \in \g$. Hence $\gamma(t) =\exp_G(tx)$ for $|t| \leq \eps$, but then
$\gamma(nt) = \exp_G(ntx)$ for $n \in \N$ leads to $\gamma(t) = \exp_G(tx)$ for each $t \in \R$.
\end{prf}

\begin{prop} \mlabel{prop:5.4.2} Let $G$ be a locally exponential Lie group. 
For $x, y \in \L(G)$ the 
\index{Trotter Product Formula} 
{\emph Trotter Product Formula} 
$$ \exp_G(x + y) = \lim_{n \to \infty} \Big(\exp_G\big(\frac{x}{n}\big)
\exp_G\big(\frac{y}{n}\big)\Big)^n$$
and the {\emph Commutator Formula} \index{Commutator Formula} 
$$ \exp_G([x,y]) = \lim_{n \to \infty} \Big(
\exp_G\big(\frac{x}{n}\big)\exp_G\big(\frac{y}{n}\big)
\exp_G\big(-\frac{x}{n}\big)\exp_G\big(-\frac{y}{n}\big)\Big)^{n^2} $$ 
hold. 
\end{prop}

\begin{prf} 
To verify the product formula, we first observe that, for $n$ sufficiently 
large, we have  
\[\Big(\exp_G\big(\frac{x}{n}\big)\exp_G\big(\frac{y}{n}\big)\Big)^n
= \exp_G\Big(n\Big(\frac{x}{n} * \frac{y}{n}\Big)\Big),\]
where $m(x,y) := x * y := \exp_G^{-1}(\exp_G x \exp_G y)$ denotes the local 
multiplication in exponential coordinates. Hence it suffices to show that 
$$ n\cdot \frac{x}{n} * \frac{y}{n} \to x + y, $$
but this follows from 
$\dd m(0,0)(x,y) = x + y$ 
(Remark~\ref{rem:brack-taylor}). 

For the commutator formula, we likewise observe that 
$$\Big(
\exp_G\big(\frac{x}{n}\big)\exp_G\big(\frac{y}{n}\big)
\exp_G\big(-\frac{x}{n}\big)\exp_G\big(-\frac{y}{n}\big)\Big)^{n^2}
\!\!\!\!\!= \exp_G\Big(n^2  \Big(\frac{x}{n} * \frac{y}{n} * 
\frac{-x}{n} * \frac{-y}{n}\Big)\Big), $$
so that we have to verify that 
$$ [x,y] = \lim_{n \to \infty} n^2 \Big(\frac{x}{n} * \frac{y}{n} * 
\frac{-x}{n} * \frac{-y}{n}\Big). $$
This follows from the fact that the second order Taylor polynomial 
of the commutator map $x * y * (-x)* (-y)$ coincides with the 
Lie bracket $[x,y]$ (Remark~\ref{rem:brack-taylor}(c)). 
\end{prf}

Before we come to the Automatic Smoothness Theorem, 
we state a convenient criterion for smoothness of group 
homomorphisms. 

\begin{prop} \mlabel{prop:smooth-crit}
Let $G_1$ and $G_2$ be Lie groups with an exponential function 
and $\phi \: G_1 \to G_2$ be a group homomorphism.
If $G_1$ is locally exponential, then the following are equivalent:
\begin{enumerate}
\item[\rm(a)] $\phi$ is smooth in an identity neighborhood of
$G_1$.
\item[\rm(b)] $\phi$ is smooth.
\item[\rm(c)] There exists a
continuous linear map $\psi \: \L(G_1) \to \L(G_2)$ satisfying
\begin{equation}
  \label{eq:exp-rel2}
 \exp_{G_2} \circ \psi = \phi \circ \exp_{G_1}.
\end{equation}
\end{enumerate}
\end{prop}

\begin{prf} (a) $\Rarrow$ (b): (cf.\ Exercise~\ref{exer:3.1.4b}) 
Let $U$ be an open $\be$-neighborhood
of $G_1$ such that $\phi\res_U$
is smooth. Since each left translation $\lambda_g$ is a diffeomorphism,
$\lambda_g(U) = gU$ is an open neighborhood of $g$, and we have
$$ \phi(gx) = \phi(g)\phi(x), \quad \hbox{i.e.,} \quad
\phi \circ \lambda_g = \lambda_{\phi(g)} \circ \phi. $$
Hence the smoothness of $\phi$ on $U$ implies the smoothness of
$\phi$ on $gU$, and therefore that $\phi$ is smooth.

(b) $\Rarrow$ (c): If $\phi$ is smooth, then $\psi := \L(\phi)$ satisfies
\eqref{eq:exp-rel2} by Proposition~\ref{prop:exp-diag}. 

(c) $\Rarrow$ (a): If $\psi$ is a continuous linear map
satisfying \eqref{eq:exp-rel2}, then the fact that 
$\exp_{G_1}$ is a local diffeomorphism and
the smoothness of $\psi$ imply (a).
\end{prf}

\begin{thm} [Automatic Smoothness Theorem] \mlabel{thm:5.4.4} 
\index{Automatic Smoothness Theorem} 
Each continuous homomorphism $\phi \: G \to H$ 
of locally exponential Lie groups is smooth. 
\end{thm}

\begin{prf} For $x \in \L(G)$, put $\gamma_x(t) := \exp_G(tx)$. Then 
$\phi \circ \gamma_x \: \R \to H$ is a continuous one-parameter group,
hence can be written as 
$$ \phi \circ \gamma_x = \gamma_{\psi(x)}\quad \mbox{ for some } \quad \psi(x) \in \L(H) $$
(Lemma \ref{lem:5.4.1}). 
We conclude that 
\begin{equation}
  \label{eq:exp-hom}
\phi \circ \exp_G = \exp_H \circ \psi. 
\end{equation}
Using the Trotter Product Formula (Proposition~\ref{prop:5.4.2}), 
we see that $\psi$ is linear. 
Since $\exp_G$ and $\exp_H$ are local diffeomorphisms and 
$\phi$ is continuous, $\psi \: \L(G) \to \L(H)$ is continuous. 
Now Proposition~\ref{prop:smooth-crit} implies 
that $\phi$ is smooth, and $\psi = \L(\phi)$ follows from 
$T_0(\exp_{G_j}) = \id_{\g_j}$ for $j =1,2$. 
\end{prf}

For Lie groups which are not locally exponential, we still have 
that any $C^1$-homomorphism is smooth (Exercise~\ref{ex:c1smooth}), 
but we do not have any example of a continuous homomorphism between 
Lie groups which is not smooth. 

\begin{cor} \mlabel{cor:5.4.5} 
A topological group $G$ carries at most one locally exponential 
Lie group structure.   
\end{cor}

\begin{defn}
  \mlabel{def:locexptopgrp}
In view of the preceding corollary, we may think of locally exponential 
Lie groups as a special class of topological groups. More precisely, 
the category of locally exponential 
Lie groups is a full subcategory of the category of topological groups. 
In this sense, we call a topological group 
{\it locally exponential}  \index{Lie algebra!locally exponential}
if it is a locally exponential Lie group. 
\end{defn}

We conclude this subsection with a criterion for a morphism 
of locally exponential Lie groups to be open. 

\begin{prop} \mlabel{prop:hom-prop}
For a morphism $\phi\: G_1 \to G_2$ of locally exponential Lie groups, the following
assertions hold:
\begin{description}[\textup{(3)}]
\item[\textup{(1)}] $\ker\L(\phi) = \{ x \in\L(G_1) \: \exp_{G_1}(\R x)
\subeq \ker\phi\}$.
\item[\textup{(2)}] $\phi$ is an open map if and only if $\L(\phi)$ is an 
open map, i.e., a linear quotient map. 
\item[\textup{(3)}] If $\L(\phi)$ is a topological isomorphism and $\phi$
is bijective, then $\phi$ is an isomorphism of Lie groups.
\end{description}
\end{prop}

\begin{prf} (1) The condition $x \in\ker\L(\phi)$ is equivalent to
\[
\{\be\} = \exp_{G_2}\bigl(\R\L(\phi)x\bigr) = \phi\bigl(\exp_{G_1}(\R x)\bigr).
\]

(2) Suppose first that $\phi$ is an open map. Since $\exp_{G_i}$, $i=1,2$,
are local diffeomorphisms,
\begin{equation}
\label{eq:ex-rel}
\exp_{G_2} \circ\L(\phi) = \phi\circ\exp_{G_1}
\end{equation}
implies that there exists some $0$-neighborhood in
$\L(G_1)$ on which $\L(\phi)$ is an open map, hence that $\L(\phi)$
is open. 

If, conversely, $\L(\phi)$ is an open map, 
then \eqref{eq:ex-rel} implies that there exists an 
open $\be$-neighborhood $U_1$ in $G_1$
such that $\phi\res_{U_1}$ is an open map. We claim that this implies
that $\phi$ is an open map. In fact, suppose that $O \subeq G_1$ is open
and $g \in O$. Then there exists an open $\be$-neighborhood $U_2$ of $G_1$
with $g U_2 \subeq O$ and $U_2 \subeq U_1$. Then
$\phi(O) \supeq\phi(g U_2) = \phi(g) \phi(U_2),$ 
and since $\phi(U_2)$ is open in $G_2$, we see that $\phi(O)$ is a
neighborhood of $\phi(g)$, hence that $\phi(O)$ is
open because $g \in O$ was arbitrary.

(3) From the relation $\exp_{G_2} \circ\L(\phi) = \phi\circ\exp_{G_1}$
and the bijectivity of $\phi$, we derive that the group homomorphism
$\phi^{-1}$ satisfies
\[
\phi^{-1} \circ\exp_{G_2} = \exp_{G_1} \circ\L(\phi)^{-1},
\]
so that Proposition~\ref{prop:smooth-crit} implies that $\phi^{-1}$ is also smooth.
\end{prf}

\begin{thm} \mlabel{thm4.12} 
If $\phi \: G_1 \to G_2$ is a morphism of locally exponential Fr\'echet--Lie 
groups for which $\L(\phi) \: \g_1 \to \g_2$ is surjective, then 
$\phi$ has a continuous local section. 
\end{thm}

\begin{prf} Let $\exp_{G_j} \: \g_j \to G_j$ denote the exponential
functions and $U_j \subeq \g_j$ be an open convex identity neighborhood for 
which $\exp_{G_j}\res_{U_j}$ is a diffeomorphism onto an open subset
of $G_j$. 
Then Michael's Theorem (\cite[Ch.\,II, \S4.7, Prop.\,12]{Bou87}) 
implies that there exists a continuous
section 
$\sigma_L \: \g_2 \to \g_1$. Hence 
$$ \sigma_G := \exp_{G_1} \circ \sigma_L\res_{U_2} \circ
(\exp_{G_2}\res_{U_2})^{-1} \: \exp_{G_2}(U_2) \to G_1 $$
is a continuous section of $\phi$ on $\exp_{G_2}(U_2)$. 
\end{prf}

\subsection{Integrating Lie algebra morphisms}

Note that we do not assume in the following theorem that the target group 
$H$ is locally 
exponential, we only assume that it has an  
exponential function. 

\begin{thm} [Integration Theorem] \mlabel{thm:int-thm} 
\index{Integration Theorem!loc.\ exp.\ to group with exp.~fct}
Let $H$ be a Lie group 
with an exponential function, 
$G$ a $1$-connected locally exponential Lie group, 
and $\psi \: \L(G) \to \L(H)$ a continuous homomorphism of Lie algebras. 
Then there exists a unique smooth homomorphism $\eta \: G \to H$
with $\L(\eta) = \psi$. 
\end{thm}

\begin{prf} Let $U \subeq \L(G)$ be an open $0$-neighborhood for which 
$\exp_G\res_U$ is a diffeomorphism onto an open subset $W := \exp_G U$ of $G$. We consider 
the smooth function 
$$ \phi := \exp_H \circ \psi \circ (\exp_G\res_U)^{-1}\: W \to H, $$ 
satisfying 
$$ \phi \circ \exp_G\res_U = \exp_H \circ \psi. $$
For the left logarithmic derivative 
$\delta(\phi) = \phi^* \kappa_H$ of $\phi$, this identity leads
on $U$ to 
$$ \exp_G^*\phi^*\kappa_H 
= \psi^*\exp_H^*\kappa_H 
= \psi^*\delta(\exp_H). $$
For $x \in U$, we further get with Exercise~\ref{ex:adint-compat}
\begin{eqnarray*}
&& \big(\psi^*\delta(\exp_H)\big)_x 
= \delta(\exp_H)_{\psi(x)} \circ \psi 
= \int_0^1 e^{-t \ad \psi(x)}\, dt \circ \psi \\
&=& \int_0^1 e^{-t \ad \psi(x)} \circ \psi\, dt 
= \int_0^1 \psi \circ e^{-t \ad x}\, dt 
= \psi \circ \int_0^1 e^{-t \ad x}\, dt \\
&=& \psi \circ \delta(\exp_G)_x 
= \psi \circ (\exp_G^* \kappa_G)_x= (\exp_G^* (\psi \circ \kappa_G))_x.
\end{eqnarray*}
We thus obtain on $U$ the relation 
$\exp_G^*\phi^*\kappa_H = \exp_G^* (\psi \circ \kappa_G).$
Since $\exp_G\res_U$ is a diffeomorphism onto $\exp_G(U)$, this leads
to 
$$ \delta(\phi) = \phi^*\kappa_H =\psi \circ \kappa_G$$
on $\exp_G(U)$. Now Lemma~\ref{lem:c.15} 
implies the existence of an open symmetric \break $\be$-neighborhood $W_1 \subeq W$ 
with $\phi(xy) = \phi(x)\phi(y)$ for $x,y \in W_1$. Since 
$G$ is $1$-connected, 
Proposition~\ref{prop:lochomo-ext} yields an extension $\eta$ of $\phi\res_{W_1}$ 
to a homomorphism $G \to H$. We then have $\L(\eta) = \delta(\phi)_\be = \psi$, 
which proves existence. The uniqueness of $\eta$ with $\L(\eta) = \psi$ 
follows from the connectedness of~$G$ (Proposition~\ref{propc.14}). 
\end{prf} 

\begin{cor} \mlabel{cor:locexp-unique} 
If $G_1$ and $G_2$ are $1$-connected locally exponential Lie groups with isomorphic 
Lie algebras, then $G_1 \cong G_2$. 
\end{cor}

\begin{prf} Apply Theorem~\ref{thm:int-thm} to an isomorphism 
$\psi \: \L(G_1) \to \L(G_2)$ and its inverse 
$\psi^{-1} \: \L(G_2) \to \L(G_1)$ to obtain 
morphisms $\phi_1 \: G_1 \to G_2$ and $\phi_2 \: G_2 \to G_1$ with 
$\L(\phi_1) = \psi$ and $\L(\phi_2) = \psi^{-1}$. 
Then $\phi_1 \circ \phi_2= \id_{G_2}$ and 
$\phi_2 \circ \phi_1= \id_{G_1}$ follow from Proposition~\ref{propc.14}, 
so that $\phi_1$ is an isomorphism. 
\end{prf}

\begin{cor} \mlabel{cor:autgrp} If $G$ is a $1$-connected locally 
exponential Lie group, then the map 
$$ \L \: \Aut(G) \to \Aut(\g) $$
is an isomorphism of groups. 
\end{cor}

\begin{thm}[Locally exponential groups with the same Lie algebra]
\mlabel{thm:5.2.11} 
\index{Lie groups!locally exponential, with ame Lie algebra} 
Let $\g$ be the Lie algebra of a locally exponential Lie group. 
Then there exists a $1$-connected Lie group $G$ with Lie algebra~$\g$ 
which is unique up to isomorphism. 
Every connected locally exponential Lie group with Lie algebra~$\g$ is isomorphic 
to a quotient $G/\Gamma$, where $\Gamma\subeq Z(G)$ is a discrete subgroup. 
For two discrete subgroups $\Gamma_1, \Gamma_2 \subeq Z(G)$, 
the groups $G/\Gamma_1$ and $G/\Gamma_2$ are isomorphic if and only 
if there exists an automorphism $\phi \in \Aut(G)\cong \Aut(\g)$ with 
$\phi(\Gamma_1) = \Gamma_2$. 
\end{thm}

\begin{prf} Proposition~\ref{prop:4.1.5} implies the existence of a 
$1$-connected locally exponential Lie group $G$ with Lie algebra~$\g$ 
and the uniqueness up to isomorphism follows  from Corollary~\ref{cor:locexp-unique}. 
For any connected locally exponential Lie group~$H$ whose Lie algebra 
is isomorphic to $\g$, it follows that the universal covering group 
$\tilde H$ is isomorphic to $G$. 
Now the remaining assertions follow from Theorem~\ref{thm:3.3.6}. 
\end{prf}

\begin{rem}
If $G_0$ and a discrete group $\Gamma$ are given, then the determination of all Lie groups 
$G$ with identity component $G_0$ and component group $\pi_0(G)\cong \Gamma$ 
corresponds to the classification of all Lie group extensions 
\[  \1 \to G_0 \into G \onto \Gamma \to \1 \] 
(cf.\ Definition~\ref{def:liegrp-ext}), 
i.e., to a description of the set $\Ext(\Gamma,G_0)$. 
Extension problems of this type are discussed in \cite[Ch.~18]{HiNe12}; 
see also Definition~\ref{def:liegrp-ext} below. 
\end{rem}

\begin{rem} Theorem \ref{thm:int-thm} ensures for each continuous homomorphism 
$\psi \: \L(G) \to \L(H)$ of the Lie algebra of a locally exponential 
connected Lie group $G$ to the Lie algebra of a Lie group $H$ with exponential 
function, the existence of a unique morphism of 
Lie groups  $\tilde\eta \: \tilde G \to H$ with 
$\L(\tilde\eta) = \psi$. Here we identify the Lie algebras of $G$ and $\tilde G$, 
so that  the universal covering 
morphism $q_G \: \tilde G \to G$ satisfies $\L(q_G) = \id_{\L(G)}$. 
Then the period homomorphism 
$$ \per_\psi := \tilde\eta\res_{\pi_1(G)} \: \pi_1(G) \to H $$
is an obstruction to the existence of $\eta \: G \to H$ with 
$\L(\eta) = \psi$. In fact, if $\eta$ exists, then 
$\tilde\eta = \eta\circ q_G$ implies that $\per_\psi$ vanishes. 
If, conversely, $\per_\psi$ vanishes, then $\tilde\eta$ factors 
through some $\eta$ with $\eta \circ q_G = \tilde\eta$, so that 
$\L(q_G) = \id_{\L(G)}$ yields $\L(\eta) = \psi$. 
\end{rem}

\subsection{Locally $k$-exponential Lie groups} 

The Integration Theorem~\ref{thm:int-thm} extends naturally to a larger class of 
groups which are not locally exponential. 
\begin{defn} A Lie group $G$ with an exponential function 
is called {\it locally $k$-exponential} \index{Lie group!locally $k$-exponential} 
if its Lie algebra 
$\g$ is a topological direct sum 
\[ \g = \g_1\oplus \cdots \oplus \g_k\]
 of closed subspaces 
such that the map 
\[ \Phi \: \g \to G, \quad x = x_1 + \ldots + x_k \mapsto 
\exp_G(x_1) \cdots \exp_G(x_k) \] 
is a local diffeomorphism in $0$.   
\end{defn}

\begin{ex} In Example~\ref{ex:bad-expfct} we have 
seen that Lie groups of the form 
$G = E \rtimes_\alpha \R$, where $\alpha \: \R \to \GL(E)$ defines a smooth 
$\R$-action on the locally convex space $E$ need not be locally 
exponential. However, they are certainly locally-$2$-exponential. 
More generally, any semidirect product Lie group 
$N \rtimes_\alpha G$ of two locally exponential Lie groups 
$N$ and $G$ is locally $2$-exponential. 
\end{ex}

\begin{thm} 
[Integration Theorem--Locally $k$-exponential case] \mlabel{thm:int-thm-lock} 
\index{Integration Theorem!locally $k$-exponential case} 
Let $H$ be a Lie group 
with an exponential function and $G$ be a $1$-connected locally $k$-exponential Lie group, 
and $\psi \: \L(G) \to \L(H)$ a continuous homomorphism of Lie algebras. 
Then there exists a unique smooth homomorphism $\eta \: G \to H$
with $\L(\eta) = \psi$. 
\end{thm}

\begin{prf} Write $\g = \L(G)$ as a topological direct sum 
$\g = \g_1\oplus \cdots \oplus \g_k$ of closed subspaces 
such that the map 
\[ \Phi \: U \to G, \quad x = x_1 + \ldots + x_k \mapsto 
\exp_G(x_1) \cdots \exp_G(x_k) \] 
maps an open connected $0$-neighborhood of the form 
$U = U_1 + \cdots + U_k$ with open symmetric $0$-neighborhoods 
$U_j \subeq \g_j$, diffeomorphically onto an open subset $W \subeq G$. 
We consider the smooth function 
\[ \phi \: W \to H,\quad 
\exp_G(x_1) \cdots \exp_G(x_k) \mapsto 
\exp_H(\psi(x_1)) \cdots \exp_H(\psi(x_k)).\] 
For $x \in U$, we get with Exercise~\ref{ex:adint-compat} 
\begin{eqnarray*}
&&\delta(\exp_H)_{\psi(x)} \circ \psi 
= \int_0^1 e^{-t \ad \psi(x)}\, dt \circ \psi = \int_0^1 e^{-t \ad \psi(x)} \circ \psi\, dt \\
&=& \int_0^1 \psi \circ e^{-t \ad x}\, dt 
= \psi \circ \int_0^1 e^{-t \ad x}\, dt = \psi \circ \delta(\exp_G)_x.
\end{eqnarray*}
For the left logarithmic derivative of $\phi \circ \Phi$, we thus obtain with the product rule 
\begin{align*}
 \delta(\phi\circ \Phi)_x v 
&= \delta(\exp_H)_{\psi(x_k)}\psi(v_k)
+ e^{-\ad \psi(x_k)} \delta(\exp_H)_{\psi(x_{k-1})} \psi(v_{k-1}) + \cdots \\ 
&\qquad \qquad \qquad + e^{-\ad \psi(x_k)}\cdots e^{-\ad \psi(x_2)} \delta(\exp_H)_{\psi(x_1)} \psi(v_1)\\
&= \psi\delta(\exp_G)_{x_k} v_k 
+ e^{-\ad \psi(x_k} \psi\delta(\exp_G)_{x_{k-1}} v_{k-1} + \cdots \\ 
&\qquad \qquad \qquad + e^{-\ad \psi(x_k)}\cdots e^{-\ad \psi(x_2)} \psi\delta(\exp_G)_{x_1} v_1 \\ 
&= \psi\Big(\delta(\exp_G)_{x_k} v_k 
+ e^{-\ad x_k} \delta(\exp_G)_{x_{k-1}} v_{k-1} + \cdots \\ 
&\qquad \qquad \qquad + e^{-\ad x_k}\cdots e^{- \ad x_2} \delta(\exp_G)_{x_1} v_1\Big)\\
&= \big(\psi \circ \delta(\Phi)_x\big)v. 
\end{align*}
This means that 
\[ \Phi^* \delta(\phi) = \delta(\phi \circ \Phi) 
= \psi \circ \delta(\Phi) = \Phi^*(\psi \circ \kappa_G),\] 
and this in turn leads to 
\[  \delta(\phi) = \psi \circ \kappa_G.\] 
Now Lemma~\ref{lem:c.15} 
implies the existence of an open symmetric \break $\be$-neighborhood $W_1 \subeq W$ 
with $\phi(xy) = \phi(x)\phi(y)$ for $x,y \in W_1$. Since 
$G$ is $1$-connected, 
Proposition~\ref{prop:lochomo-ext} yields an extension $\eta$ of $\phi\res_{W_1}$ 
to a homomorphism $G \to H$. We then have $\L(\eta) = \delta(\phi)_\be = \psi$, 
which proves existence. The uniqueness of $\eta$ with $\L(\eta) = \psi$ 
follows from the connectedness of~$G$ (Proposition~\ref{propc.14}). 
\end{prf}

\subsection{Topological groups with Lie algebra} 

We have seen above that a topological group $G$ carries at most one 
structure of a locally exponential Lie group. It is therefore a natural 
question how the ``Lie structure'' can be recovered from the 
topological group $G$. The following proposition recovers 
the locally convex space underlying the Lie algebra. 

For a topological group $G$, we write $\fL(G) := \Hom_c(\R,G)$ for 
the set of one-parameter groups, 
endowed with the compact open topology. If $G$ is a locally 
exponential Lie group, we have seen in Lemma~\ref{lem:5.4.1} that the map 
$$ \Gamma \: \L(G) \to \fL(G), 
\quad x \mapsto \gamma_x, \quad \gamma_x(t) = \exp_G(tx) $$
is a bijection. We now show that it is also compatible with the topology 
on both sides. 

\begin{prop} For any locally exponential Lie group $G$, the map 
$\Gamma$ is a homeomorphism. 
\end{prop}

\begin{prf} (1) First we observe that the map 
$$ (\exp_G)_* \: C(\R,\L(G))_{c.o.} \to C(\R,G)_{c.o.}, \quad \xi \mapsto \exp_G \circ \xi $$
is continuous with respect to the compact open topology. In fact, for each 
compact subset $K \subeq G$ and each open subset $O \subeq G$, the set 
\[  (\exp_G)_*^{-1}(\lfloor K,O\rfloor) = \lfloor K,\exp_G^{-1}(O)\rfloor  \] 
is open in $C(\R,\L(G))_{c.o.}$. 

(2) Next we observe that, for each topological vector space space $E$, the map 
$$ \ev \: \fL(E,+) \to E, \quad \gamma \mapsto \gamma(1) $$
is a homeomorphism whose inverse is given by $\ev^{-1}(x)(t) = tx$. 

The continuity of the evaluation map follows from 
$\ev^{-1}(O) \supeq \lfloor \{\be\},O\rfloor$ for each open 
subset $O \subeq E$. For the bijectivity, we refer to Lemma~\ref{lem:5.4.1}. 
To see that $\ev^{-1}$ is continuous, i.e., that $\ev$ is an open map, 
note that, for each compact subset 
$K \subeq \R$ and each open subset $O \subeq E$, we have 
$$ \ev(\lfloor K,O\rfloor) = \{ x \in E \: K \cdot x \subeq O\}, $$
which is an open subset of $E$ because the scalar multiplication 
$\R \times E \to E$ is continuous and $K$ is compact. 

Combining (1) and (2), we see that $\Gamma= (\exp_G)_* \circ \ev^{-1}$ 
is continuous for every Lie group with an exponential function. 

(3) To see that the inverse map $\Gamma^{-1} \: \gamma \mapsto \gamma'(0)$ is also continuous, 
we need that $G$ is locally exponential. Let 
$U_G \subeq G$ be an open identity neighborhood 
for which there exists a balanced  
open convex $0$-neighborhood $U_\g \subeq \g$, such that 
$\exp_G \res_{2U_\g}\: 2U_\g\to U_G$ is a diffeomorphism. 

Consider $\gamma = \gamma_x \in \fL(G)$. Then there exists some 
$\eps > 0$ with $[0,\eps] x \subeq U_\g$, which implies that 
$\gamma_x$ is contained in the 
open subset $\lfloor [0,\eps],\exp_G(U_\g)\rfloor$ of $\fL(G)$. 

For any $\gamma_y \in \lfloor [0,\eps],\exp_G(U_\g)\rfloor$, we have 
$\exp_G([0,\eps]y) \subeq \exp_G(U_\g)$. 
We claim that $\eps y \in U_\g$. If this is not the case, then 
there exists some positive $t < \eps$ with $ty \in 2 U_\g\setminus U_\g$. 
Then 
$$ \exp_G(ty) \in \exp_G([0,\eps]y) \subeq \exp_G(U_\g), $$
together with the injectivity of $\exp_G$ on $2 U_\g$ leads to a contradiction; 
hence our claim. 
We conclude that $y \in \eps^{-1} U_\g,$ so that 
$$ \Gamma^{-1}(\lfloor [0,\eps], \exp_G(U_\g)\rfloor) \subeq \eps^{-1} U_\g. $$
On $\lfloor [0,\eps], \exp_G(U_\g)\rfloor$ 
the map $\Gamma^{-1}$ is given by the continuous map
\[ \gamma_z \mapsto z = \frac{1}{\eps}\exp_G^{-1}(\gamma_z(\eps)),\] 
which shows that $\gamma_y \to \gamma_x$ uniformly on compact sets implies that 
$y \to x$. This proves that $\Gamma^{-1}$ is continuous in~$\gamma_x$.
\end{prf}

Identifying $\L(G)$ with $\fL(G)= \Hom_c(\R,G)$, the 
scalar multiplication of $\L(G)$ can be written as 
\begin{eqnarray}
  \label{eq:4.1.3b}
(\lambda\alpha)(t) 
:= \alpha(\lambda t), \quad \lambda \in \R, \alpha \in \Hom_c(\R,G), 
\end{eqnarray}
and, in view of Proposition~\ref{prop:5.4.2}, 
addition and Lie bracket may be written on the  level 
of one-parameter groups by 
\begin{eqnarray}
  \label{eq:4.1.4b}
(\alpha + \beta)(t) := \lim_{n \to \infty} 
\Big(\alpha(\frac{t}{n})\beta(\frac{t}{n})\Big)^n
\end{eqnarray}
and 
\begin{eqnarray}
  \label{eq:4.1.5}
[\alpha,\beta](t^2) := \lim_{n \to \infty} \Big(
\alpha(\frac{t}{n})\beta(\frac{t}{ n})
\alpha(-\frac{t}{ n})\beta(-\frac{t}{ n})\Big)^{n^2}. 
\end{eqnarray}

In \cite{HoM05, HoM07}, Hofmann and Morris use 
(\ref{eq:4.1.3b})-(\ref{eq:4.1.5}) 
as the starting 
point in the investigation of a remarkable class of topological groups: 

\begin{defn}
Let $G$ be a topological group 
and ${\mathfrak L}(G) := \Hom_c(\R,G)$ be the set of one-parameter groups, 
endowed with the compact open topology. 
Then $G$ is said to be \index{topological group with Lie algebra} 
{\it a topological group with Lie algebra} if 
the limits in (\ref{eq:4.1.3b})-(\ref{eq:4.1.5}) 
exist for $\alpha,\beta \in \Hom_c(\R,G)$, 
define elements of ${\mathfrak L}(G)$, addition and bracket are continuous maps 
${\mathfrak L}(G) \times {\mathfrak L}(G) \to {\mathfrak L}(G),$
and ${\mathfrak L}(G)$ is a real Lie algebra
with respect to the scalar multiplication (\ref{eq:4.1.3b}), 
the addition (\ref{eq:4.1.4b}), 
and the bracket (\ref{eq:4.1.5}). This implies 
that $\fL(G)$ is a topological Lie algebra. 
The exponential function of $G$ is defined by 
$$ \exp_G \: \fL(G) \to G, \quad \gamma \mapsto \gamma(1). $$
\end{defn}

In \cite[Sect.~1.5]{BCR81}, Boseck, Czichowski and Rudolph define smooth 
functions on a topological group in terms of restrictions to  
one-parameter groups, which leads them to (\ref{eq:4.1.4b})-(\ref{eq:4.1.5}), 
together 
with the assumption that $\fL(G)$ can be identified with the set of 
derivations of the algebra of germs of smooth functions in $\be$. 

In these terms, we have just seen that any locally exponential Lie group 
is a topological group with Lie algebra. Since $\R$ is connected, a topological 
group $G$ has a Lie algebra if and only if its identity component $G_0$ does. 
In \cite[Thm.~2.3]{HoM05} it is also 
observed that any abelian topological group is a group with Lie algebra,
where the addition on ${\mathfrak L}(G)$ is pointwise multiplication and the bracket 
is trivial. Here is a slight generalization: 

\begin{thm}   Each $2$-step nilpotent topological group has a Lie 
algebra. 
\end{thm}

\begin{prf}  (Sketch) The commutator map $c \: G \times G \to Z(G)$ is an 
alternating bihomomorphism. Then direct calculations lead to the formulas 
$$ (\alpha + \beta)(t) = \alpha(t) \beta(t) c(\alpha(t), \beta(-{\textstyle\frac{t}{ 2}})) 
\quad \hbox{ and } \quad  [\alpha,\beta](t) = c(\alpha(1), \beta(t)), $$
which can be used to verify all requirements. 
\end{prf}

\begin{rem} More generally, one can actually show that 
every nilpotent topological group $G$ has a Lie algebra. 
Here Hall's formula (\cite[Ch.~2]{Bou89}) is the key, 
when combined with the fact that 
nilpotent groups generated by finitely many one-parameter groups 
are Lie groups. The latter result is proved in \cite{MS75}, 
although the arguments in this paper are not easy to digest. 
\end{rem}

\begin{small}
\subsection*{Exercises for Section~\ref{sec:5.5}}

\begin{exer} \mlabel{ex:c1smooth} 
Let $G$ and $H$ be Lie groups and 
$\phi \: G \to H$ be a homomorphism which is a $C^1$-map. 
Show that $\phi$ is smooth. \\ 
Hint: Use induction to show that, if $\phi$ is $C^k$, then it is 
$C^{k+1}$: Show that, on the tangent bundle $TG$, we have, for 
$g \in G$ and $x \in \g$, the relation  
\[ T(\phi)(g.x) = \phi(g).T_\be(\phi)x.\] 
Derive from this relation that, if $\phi$ is $C^k$, then $T(\phi)$ is 
$C^k$, and hence that $\phi$ is $C^{k+1}$.   
\end{exer}

\begin{exer}  \mlabel{exer:5.1.1} 
Let $G$ be a locally exponential Lie group 
and $\gamma \: [0,\eps] \to G$ be a curve with $\gamma(0) = \be$ which 
is differentiable in $0$. Show that 
\[ \exp(\gamma'(0)) = \lim_{n \to \infty} \gamma(1/n)^n. \]  
Use this relation to obtain a direct proof of the Trotter Product Formula 
not using the local $*$-multiplication. 
\end{exer}

\begin{exer} Let $G$ be a finite-dimensional Lie group. Show that there exist 
$x_1, \ldots, x_n \in \L(G)$ for which 
the map 
\[ \Phi \: \R^n \to G, \quad \Phi(t_1,\ldots, t_n) := \exp_G(tx_1) \cdots \exp_G(t_n x_n) \] 
is a local diffeomorphism in $0$. Use this fact to derive directly from Lemma~\ref{lem:5.4.1} 
that any continuous homomorphism $\phi \: G \to H$ into any Lie group $H$ is smooth. 
\end{exer}

\begin{exer} Let $G$ be a connected locally exponential Lie group with Lie algebra 
$\g$ and 
$q_G \: \tilde G \to G$ be a universal covering map with $\L(q_G) = \id_\g$ (Corollary~\ref{cor:unicov}). 
We identify the groups $\Aut(\g)$ and $\Aut(\tilde G)$ as in 
Corollary~\ref{cor:autgrp}. 
\begin{enumerate}
\item[\rm(1)] Explain how to associate to each automorphism 
$\alpha \in \Aut(\g)$ a period homomorphism 
$$ \per_\alpha \: \pi_1(G) \to Z(G) $$
with the property that $\per_\alpha = \be$ (the constant homomorphism) is 
equivalent to $\alpha = \L(\phi)$ for some $\phi \in \Aut(G)$. 
\item[\rm(2)] Show that 
$\Aut(G) \cong \{ \phi \in \Aut(\tilde G) \: \phi(\ker q_G) = \ker q_G\}.$
\end{enumerate}
\end{exer}
  
\begin{exer} \mlabel{ex:adint-compat} 
Let $\psi \: \g \to \fh$ be a continuous homomorphism of $\ad$-integrable 
Lie algebras (cf.\ Definition~\ref{def4.1.5}). Then each $x \in \g$ satisfies the relation 
\[ \psi \circ e^{\ad x} = e^{\ad \psi(x)} \circ \psi.\] 
Hint: Show that the curves $\gamma \: \R \to \fh$, 
$\gamma(t) := e^{-t\ad \psi(x)} \circ \psi \circ e^{t\ad x}y$ are constant. 
\end{exer}

\end{small}

\section{The Baker--Campbell--Dynkin--Hausdorff formula} \mlabel{sec:5.4} 

In finite-dimensional and Banach--Lie theory, the formalism related to the 
Hausdorff series plays an important role because it 
provides a very direct way to pass from the Lie algebra, the infinitesimal 
level, to the local level. From the existence of non-analytic exponential 
Lie groups, it follows that, in general, we need more refined tools 
for this passage. 

In this section, we  turn to the  
Hausdorff series and its applications. 
We first show that, for any nilpotent Lie algebra $\g$, 
any corresponding exponential group structure $(\g,*)$ is given by the 
Hausdorff series 
$$ x * y = x + y + \frac{1}{2}[x,y]+\frac{1}{12}[x,[x,y]]+\frac{1}{12}[y,[y,x]]+\ldots. $$
It follows in particular that it is uniquely determined by the Lie algebra 
structure and that it actually is polynomial, where the degree is bounded by the 
nilpotent length of $\g$. 

In Subsection~\ref{subsec:5.4.2}, we show the converse, namely that the 
Hausdorff series 
defines on each locally convex nilpotent Lie algebra a polynomial group structure. 
In the infinite-dimensional context, pro-nilpotent Lie algebras form a natural 
generalization of the nilpotent ones and our results on nilpotent Lie algebras 
immediately imply that all pro-nilpotent ones are also exponential. 
Typical examples of pro-nilpotent Lie algebras and groups arise from 
groups of formal diffeomorphisms, resp., vector fields, which are discussed in some 
detail in Subsection~\ref{subsec:5.4.3}. 

In Subsection~\ref{subsec:5.4.4}  below, we apply the information 
obtained on the Hausdorff series 
in showing that all Banach--Lie algebras 
are in fact locally exponential because the Hausdorff series converges locally 
and defines an associative multiplication on sufficiently small $0$-neighborhoods. 

We conclude this section with the result that, although the multiplication 
in an exponential local Lie group need not be analytic, its Taylor series 
in $0$ always coincides with the Hausdorff series. 

\subsection
[Exponential nilpotent Lie groups]
{The BCH formula for exponential nilpotent Lie groups}
 \mlabel{subsec:5.4.1}  

We start with some remarks on ``nilpotent'' functional calculus, which is 
quite elementary and powerful 
enough to deal with nilpotent Lie algebras. Our overall 
strategy is to reduce the general case to the nilpotent one. 

Let $\cA$ be a unital associative $\K$-algebra, where $\K$ is any field of 
characteristic zero. 
We write 
\[ \K[X] := \{ a_0 + a_1 X + \ldots + a_d X^d \: d \in \N_0, a_j \in \K\} \] 
for the polynomial ring in one indeterminate $X$ over $\K$, 
and 
\[ \K[[X]]:= \Big\{ \sum_{k = 0}^\infty a_k X^k \: a_k \in \K\Big\} \] 
for the ring of formal power series in $X$. 
Then we have for each $a \in \cA$ a 
unique insertion homomorphism 
$$ \eta_a \: \K[X] \to \cA, \quad  f(X) \mapsto f(a). $$
If, in addition, $a$ is nilpotent, this insertion homomorphism extends to a 
homomorphism 
$$ \eta_a \: \K[[X]] \to \cA, \quad  \sum_{i = 0}^\infty \lambda_i X^i \mapsto 
\sum_{i=0}^\infty \lambda_i a^i, $$
which is defined because sufficiently large powers of $a$ vanish. 
Next we recall the composition map 
$$ \K[[X]] \times X \K[[X]] \to \K[[X]], \quad 
(f,g) \mapsto f \circ g, $$
defined for $f(X) = \sum_{i = 0}^\infty \lambda_i X^i$ by 
$$ (f \circ g)(X) 
= \sum_{i = 0}^\infty \lambda_i g(X)^i, $$
which makes sense because $g(X)^i$ contains only terms of degree $\geq i$.

\begin{lem} \mlabel{lem:nilp-compo} 
If $a \in \cA$ is nilpotent, 
$f \in \K[[X]]$ and $g \in X\K[[X]]$, then $g(a)$ is also nilpotent, and 
$$ f(g(a)) = (f \circ g)(a). $$ 
\end{lem}

\begin{prf} Suppose that $a^N = 0$. 
Writing $g(X) = X h(X)$ with $h(X) \in \K[[X]]$, we see that 
$g(a)^N = a^N h(a)^N$ vanishes, so that $g(a)$ is nilpotent. 
Hence $f(g(a))$ is defined. 
For $f(X) = \sum_{i = 0}^\infty \lambda_i X^i$, we have 
\[  f(g(a)) 
= \sum_{i = 0}^\infty \lambda_i g(a)^i 
= \Big(\sum_{i = 0}^\infty \lambda_i g^i\Big)(a) 
= (f \circ g)(a), \] 
where only finitely many summands are non-zero. 
\end{prf}

\begin{lem} \mlabel{lem:Phi-Psi} For the power series  
$$\Phi(X)=\frac{1-e^{-X}}{X} 
:= \sum_{k = 1}^\infty (-1)^{k-1}\frac{X^{k-1}}{k!} 
= \sum_{k = 0}^\infty (-1)^k \frac{X^k}{(k+1)!} \in \Q[[X]] $$
and 
$$\Psi(X)=\frac{(X+1)\log (X+1)}{X} 
:= (X+1)\sum_{k=0}^\infty \frac{(-1)^{k}}{k+1}X^k \in \Q[[X]], $$
the relation 
$$ \Psi(e^X-1)\Phi(X)=1 \quad \mbox{ holds in } \quad 
\Q[[X]]. $$
\end{lem}

\begin{prf} Since the power series $\Phi$ and $\Psi$ define complex analytic 
functions in a neighborhood of $0$ in $\C$, it suffices to 
verify the asserted relation for the corresponding functions. 
If $|z|<\log 2$, then $|e^z-1|<1$ and we obtain from 
$\log(e^z) = z$ for the holomorphic logarithm function 
on the unit circle around $1$, defined by the log-power series, the relation 
\[\Psi(e^z-1) \Phi(z) = \frac{e^z z}{e^z-1}\frac{1-e^{-z}}{z}=1.
\qedhere\] 
\end{prf}

\begin{prop} \mlabel{prop:nil-rel} 
For any nilpotent element $a$ of an associative unital algebra $\cA$, 
the element $e^a - \be$ is nilpotent, 
$\Phi(a)$ is unipotent, hence invertible, and we have 
$$ \Phi(a)^{-1} = \Psi(e^a-1). $$
\end{prop}

We now derive an explicit formula for any exponential Lie group structure 
on a nilpotent Lie algebra, called the 
{\it Baker--Campbell--Dynkin--Hausdorff formula}. 
\index{Baker--Campbell--Dynkin--Hausdorff formula} 

\begin{thm} \mlabel{thm:nilp-bch}
Let $(\g,*)$ be an exponential Lie group for which $\g$ is a nilpotent Lie algebra. 
Then the product is given by the BCH formula 
\begin{eqnarray*}
 x * y 
&&=x + \sum_{\substack{k,m\ge 0 \\ p_i+q_i>0}}(-1)^k \frac{
(\ad x)^{p_1}(\ad y)^{q_1}\ldots
(\ad x)^{p_k}(\ad y)^{q_k}(\ad x)^m}{ 
(k+1)(q_1+\ldots+q_k+1)p_1!q_1!\ldots p_k!q_k!m!}y \\
&&=x+y+\frac{1}{2}[x,y]+\frac{1}{12}[x,[x,y]]+\frac{1}{12}[y,[y,x]]+\ldots 
\end{eqnarray*}
\end{thm}

The power series in the preceding theorem is called the 
{\it Haus\-dorff series}. 
\index{Hausdorff series}

\begin{prf} Before 
we turn to the explicit formula, we start with an integral formula which 
already contains all the information.  
We claim that, for $x,y \in \g$, we have 
  \begin{eqnarray}
    \label{eq:pre-bch}
 x * y = x + \int_0^1 \Psi(e^{\ad x}e^{t\ad y}-\be) y\, dt.
  \end{eqnarray}

For $t \in \R$, consider 
$\gamma(t) := x * ty \in \g.$
We are interested in an explicit formula for $\gamma(1) = x* y$. We obtain 
with Propositions~\ref{prop:nil-rel} and Theorem~\ref{thm:exp-logder} (applied 
to $\exp_\g = \id_\g$): 
$$ \gamma'(t) 
= T_0(\lambda_{\gamma(t)})y 
= \kappa_\g(\gamma(t))^{-1} y 
= \Phi(\ad \gamma(t))^{-1} y 
= \Psi(e^{\ad \gamma(t)}-\be) y. $$
Using 
$$ e^{\ad \gamma(t)} = e^{\ad(x * ty)} 
= \Ad(x * ty) = \Ad(x)\Ad(ty) = e^{\ad x} e^{t\ad y} $$ 
(cf.\ Definition~\ref{def:kappa-liealg}), we thus get 
\begin{eqnarray*}
x * y 
&=& \gamma(1) 
= \gamma(0) + \int_0^1 \gamma'(t)\, dt 
= x  + \int_0^1 \Psi(e^{\ad x}e^{t\ad y}-\be) y\, dt.
\end{eqnarray*}

It only remains to expand the integrand of (\ref{eq:pre-bch}) as a 
series (which actually is a finite sum because of the nilpotency of $\g$): 
  \begin{align*}
&\ \ \ \ \int_0^1\Psi\big(\exp(\ad x)\exp(\ad ty)-\be\big)y\, dt\cr
&=\int_0^1\sum_{k=0}^\infty \frac{(-1)^k\big(\exp(\ad x)\exp(\ad ty)-\id\big)^k}
{(k+1)}\big(\exp(\ad x)\exp(\ad ty)\big)y\, dt\\ 
&=\int_0^1\sum_{\substack{k\ge 0\\ p_i+q_i>0}}\frac{(-1)^k}{(k+1)}
\frac{(\ad x)^{p_1}(\ad ty)^{q_1}\ldots
(\ad x)^{p_k}(\ad ty)^{q_k}}{p_1!q_1!\ldots p_k!q_k!}\exp(\ad x)y\, dt\cr
&=\sum_{\substack{k,m\ge 0 \\ p_i+q_i>0}}\!\! \frac{(-1)^k}{(k+1)}
\frac{(\ad x)^{p_1}(\ad y)^{q_1}\ldots (\ad x)^{p_k}(\ad y)^{q_k}(\ad x)^m y}
{p_1!q_1!\ldots p_k!q_k!m!}\int_0^1 \!\!\! t^{q_1+\ldots + q_k}\, dt\cr
&=\sum_{\substack{k,m\ge 0\\ p_i+q_i>0}}\frac{(-1)^k (\ad x)^{p_1}(\ad y)^{q_1}\ldots
(\ad x)^{p_k}(\ad y)^{q_k}(\ad x)^m y}
{(k+1)(q_1+\ldots+q_k+1) p_1!q_1!\ldots p_k!q_k!m!}. 
\qedhere \end{align*}
\end{prf}

\begin{cor} \mlabel{cor:pro-nil} If $(\g,*) = \prolim (\g_j,*)$ is a projective 
limit of nilpotent exponential Lie groups, 
then the product $*$ on $\g$ is given by the Hausdorff series. 
\end{cor} 

\subsection{Associativity of the BCH multiplication} \mlabel{subsec:5.4.2}

Now we turn to the converse problem: On each nilpotent rational Lie algebra~$\g$, 
the Hausdorff series defines a polynomial product 
$x*y$ satisfying $x* y = x + y$ if $[x,y]=0$. From that it follows in 
particular that 
\[ tx * s x = (t+s)x \quad \mbox{ for } \quad t,s \in \Q, \] 
but it is not clear from the explicit form
of the power series that the product is actually associative. In this 
section we provide this information without going too much into formal 
power series arguments. Instead we use the uniqueness result from the 
preceding subsection. 

\begin{defn} Let $\cA$ be a unital associative algebra. 
A {\it filtration of $\cA$} \index{filtration of $\cA$}
is a sequence $(\cF_n)_{n \in \N_0}$ of subspaces of $\cA$ 
with 
$$\cF_0 = \cA, \quad \bigcap_n \cF_n = \{0\} 
 \quad \mbox{ and } \quad \cF_n \cF_m \subeq \cF_{n+m}. $$
In particular, each $\cF_n$ is an ideal of $\cA$. 
\end{defn}

\begin{prop} \mlabel{prop:lim-filter} Let $(\cF_n)_{n \in \N}$ be a filtration of the 
algebra $\cA$ for which all quotient algebras $\cA_n := \cA/\cF_n$ are finite-dimensional 
and endow the projective limit algebra $\hat \cA := \prolim \cA_n \subeq \prod_n \cA_n$ 
with the projective limit topology. 
Put $\hat \cF_1 := \prolim \cF_1/\cF_n \subeq \hat \cA$. Then 
the closed affine subspace $U := \1 + \hat \cF_1$ is a subgroup of $\hat \cA^\times$ 
and a Lie group with respect to the induced multiplication. 
Its exponential map 
$$ \exp_U \: \hat \cF_1 \to \1 + \hat \cF_1, \quad x \mapsto \sum_{n = 0}^\infty \frac{x^n}{n!} $$
is a diffeomorphism and, for $x,y \in \hat \cF_1$, the product 
$$ x * y = \exp_U^{-1}(\exp_U(x) \exp_U(y)) $$
is given by the Hausdorff series. 
\end{prop}

\begin{prf} The ideal $\cF_1/\cF_n \subeq \cA_n$ is nilpotent because all $n$-fold products 
of elements in this ideal vanish. For $x \in \cF_1/\cF_n$, we therefore obtain an 
inverse of $\1 - x$ by $\1 + x + \ldots + x^{n-1} = \sum_{k = 0}^\infty x^k$. 
Hence the unit group of $\cA_n$ contains $\1 + \cF_1/\cF_n$. 
For any $x = (x_n)_{n \in \N} \in \hat \cF_1$, the element 
$\1 - x_n \in \1 + \cF_1/\cF_n \subeq \cA_n^\times$ is invertible, and its inverse is given by 
the geometric series. This implies that, in $\hat \cA$, we have 
$$ \1 -x \in \hat \cA^\times \quad \mbox{ with } \quad 
(\1-x)^{-1} = \sum_{k = 0}^\infty x^k, $$
where the convergence of the series in $\hat \cA$ follows from its finiteness in 
each quotient algebra $\cA_n$. We conclude that the affine subspace 
$U$ is contained in $\cA^\times$, and since $\hat \cF_1$ is an ideal 
of $\hat \cA$, namely the kernel of the canonical homomorphism $\hat \cA \to \cA_1$, 
we see that $U$ is a subgroup of $\hat \cA^\times$. 

That $U$ actually is a Lie group is a consequence of the smoothness of the 
(bilinear) multiplication and the inversion, which is given by the geometric series,  
whose smoothness follows from the smoothness of its compositions with the projections 
to the algebras $\cA_n$, which are polynomial 
(Lemma~\ref{lemPL}). 

We likewise see that $\exp_U$ defines a smooth map whose $\cA_n$-components 
are polynomials of degree $\leq n$. Its inverse is given by 
$$ \log_U \: U \to \hat \cF_1, \quad \1 - x \mapsto 
-\sum_{n = 1}^\infty \frac{x^n}{n}. $$
This function is also smooth because all its $\cA_n$-components are polynomial. 
The relations $\exp_U \circ \log_U = \id_U$ 
and $\log_U \circ \exp_U = \id_{\hat \cF_1}$ follows from 
the corresponding relations for $\exp$ and $\log$ on nilpotent, resp., unipotent 
elements of associative algebras (Exercise~\ref{exer:5.4.1}). 

We conclude that $U$ is an exponential Lie group, and since it is a 
projective limit of the nilpotent Lie groups $\1 + \cF_n/\cF_1$, 
the $\cA_n$-component of $x*y$ is given by the Hausdorff series 
by Theorem~\ref{thm:nilp-bch}, and this implies that $x * y$ is given by 
the Hausdorff series in the projective limit algebra $\hat \cF_1$. 
\end{prf}

\begin{ex} Let $\cA := \cF_{\rm ass}(x_1,\ldots, x_N)$ be the 
free associative algebra in $n$ generators. 
Then $\cA$ has a natural filtration 
$$\cF_n := \Spann\{ s_1\cdots s_m \: s_i \in \{x_1,\ldots, x_N\}, m \geq n \}. $$
Each quotient $\cA_m := \cA/\cF_m$ is finite-dimensional, so that the preceding 
Proposition~\ref{prop:lim-filter} applies. The algebra 
$$\hat \cF_{\rm ass}(x_1,\ldots, x_N) := \prolim \cA_n $$
can be identified with the algebra of non-commutative formal power series in 
the generators $x_1, \ldots, x_N$. Further, $\hat \cF_1$ is a hyperplane ideal, 
topologically generated by $x_1,\ldots, x_N$. 

From Proposition~\ref{prop:lim-filter}, we derive that 
the Hausdorff series defines on $\hat \cF_1$ an exponential Lie group structure. 
\end{ex}

Using the Poincar\'e--Birkhoff--Witt Theorem, it is easy to see that the 
closed subspace generated by the Lie polynomials in $x_1,\ldots, x_N$ 
is the complete free Lie algebra $\hat \cF_{\rm Lie}(x_1,\ldots, x_N)$, 
generated by $x_1,\ldots, x_N$ (Exercise~\ref{exer:5.4.2}). 
We thus obtain:

\begin{thm} \mlabel{thm:BCH-free} The 
Hausdorff series defines on the complete free Lie algebra 
$\hat \cF_{\rm Lie}(x_1,\ldots, x_N)$, generated by  $x_1,\ldots, x_N$, an exponential 
Lie group structure. 
\end{thm}

The following theorem complements Theorem~\ref{thm:nilp-bch} by showing that it actually 
applies to {\it all} nilpontet Lie groups.

\begin{thm} \mlabel{thm:BCH-nil} 
On each nilpotent locally convex Lie algebra $\g$, the Hausdorff series defines an 
exponential Lie group structure. 
\end{thm}

\begin{prf} For $a,b \in \g$, 
we define $a * b$ by the Hausdorff series (Theorem~\ref{thm:nilp-bch}), so that we obtain a polynomial 
binary operation on $\g$. We then have 
$$ ta * sa = (t+s) a \quad \mbox{ for } \quad a \in \g, t,s \in \R. $$
In particular, $0 * a = a * 0 = a$ and $a * (-a) = (-a) * a = 0$. 
For $a,b,c \in \g$, there exists a unique morphism of Lie algebras 
\[  \phi \: \hat \cF_{\rm Lie}(x_1, x_2, x_3) \to \g 
\quad \mbox{ with } \quad x_1 \mapsto a, x_2 \mapsto b, x_3 \mapsto c.\]
 Now the associativity of 
the product on $\g$ follows from Theorem~\ref{thm:BCH-free}: 
\begin{eqnarray*}
(a*b) * c 
&=& (\phi(x_1) * \phi(x_2)) * \phi(x_3) 
= \phi(x_1 * x_2) * \phi(x_3) = \phi((x_1*x_2)*x_3) \\
&=& \phi(x_1*(x_2*x_3)) = \phi(x_1) * (\phi(x_2)*\phi(x_3)) = a * (b*c). 
\end{eqnarray*}
We conclude that $(\g,*)$ is an exponential Lie group. 
\end{prf}

\begin{cor} If $\g$ is any nilpotent $\Q$-Lie algebra, then the BCH multiplication 
defines a group structure on $\g$, satisfying 
$$ t x * sx = (t+s)x \quad \mbox{ for } \quad x \in \g,s,t \in \Q.  $$
\end{cor}

\begin{prf} Since each finite subset of a nilpotent Lie algebra generates a 
finite-dimensional subalgebra (Exercise~\ref{exer:2.6.10}), 
it suffices to verify the assertion for finite-dimensional Lie algebras. 
If $\dim \g < \infty$, we 
embed $\g$ into the real (locally convex) Lie algebra $\g \otimes_\Q \R$ 
and apply Theorem~\ref{thm:BCH-nil} to see that the multiplication defined 
by  the Hausdorff series is associative.   
\end{prf}

\subsection{Pro-nilpotent Lie algebras} \mlabel{subsec:5.4.3} 

As an application of the general theory developed above, we briefly 
discuss projective limits of exponential Lie algebras. 

\begin{prop} \mlabel{prop:prolim-exp} Projective 
limits of exponential Lie algebras are exponential. 
\end{prop}

\begin{prf} (a) First we consider the case of a direct product, which is a special 
projective limit. The corresponding directed set is the set of finite subsets of~$J$. 
Let $(\g_j)_{j \in J}$ be a family of exponential Lie algebras 
and $\g := \prod_{j \in J} \g_j$ be their topological product, which is a locally convex 
Lie algebra. To see that $\g$ is exponential, let 
$m_j \: \g_j \times \g_j \to\g_j$ define an exponential Lie group 
structure on $\g_j$.  Then the product map 
$$m \: \g \times \g \to \g, \quad 
(x,y) = ((x_j), (y_j))\mapsto (m_j(x_j, y_j))_{j \in J} $$
defines an exponential Lie group structure on $\g$. 

(b) Now let $\g = \prolim \g_j$ be an arbitrary projective limit of a family of 
exponential Lie algebras $(\g_j)_{j \in J}$. 
Then Proposition~\ref{prop:5.2.13b} implies that the corresponding morphisms 
of Lie algebras $\phi_{ij} \: \g_j \to \g_i$ 
are also morphisms for the corresponding group structures. 
Hence 
$\g$ is a $*$-subgroup of the direct product 
Lie group $\prod_{j \in J} (\g_j,*)$. We thus obtain on 
$\g$ an exponential Lie group structure with $\L(\g,*)= \g$. 
\end{prf}

Combining the preceding proposition with Theorem~\ref{thm:BCH-nil}, we get: 

\begin{cor} \mlabel{cor:pro-nilb}
Pro-nilpotent Lie algebras, i.e., projective limits of nilpotent Lie algebras, 
are exponential. 
\end{cor}

\begin{ex} \label{ex:5.2.formal} 
The preceding corollary applies in many interesting situations, where it implies 
that certain pro-nilpotent groups carry natural Lie group structures. 

An important class of examples arises as follows. Let $E$ be a finite-dimensional $\K$-vector space, 
let $P_d(E)$ denote the space of all homogeneous polynomial functions $E \to E$ of degree $d$.  
Then, for each $n \geq 2$, the space $\g_n := \oplus_{k = 2}^n P_k(E)$ 
carries a natural Lie algebra structure, given for 
$f \in P_i(E)$ and $g \in P_j(E)$ by 
$$ [f,g](x) := \left\{ 
       \begin{array}{cl} 
           \dd f(x)g(x) -\dd g(x)f(x) & \mbox{for $i + j -1 \leq n$} \cr 
           0 & \mbox{for $i + j - 1 > n$}. 
\end{array} \right. $$
This is a modification of the natural Lie bracket on the space $C^\infty(E,E)$ 
obtained by changing signs and omitting all terms of degree $> n$. From 
$$ [P_i(E), P_j(E)] \subeq P_{i+j-1}(E), $$
it immediately follows that each $\g_n$ is a nilpotent Lie algebra. 
For $n < m$, we have a natural projection 
$$ \phi_{nm} \: \g_m \to \g_n, $$
which is a homomorphism of Lie algebras. The projective limit 
Lie algebra 
$\g := \prolim \g_n$ can be identified with the space of $E$-valued formal 
power series in $\dim(E)$ variables starting in degree~$2$.

To obtain a natural Lie group corresponding to $\g_n$, we consider the set 
$G_n$ of all polynomial maps $f = \id_E + \xi \: E \to E$ with $\xi\in \g_n$. 
On $G_n$ we define a group structure 
by composition and then omitting all terms of order $> n$: 
\begin{eqnarray}
  \label{eq:mult}
f \# g := (f \circ g)_{\leq n} \quad \mbox{ with } \quad 
g_{\leq n} := g_2 + g_3 + \ldots + g_n.
\end{eqnarray}

Clearly, $\id_E$ is a neutral element with respect to this multiplication. 
From the Chain Rule for Taylor polynomials (Proposition~\ref{chainrtay}), we obtain for 
$f,g,h \in G_n$ the associativity condition: 
\begin{eqnarray*}
f \# (g \# h) 
&=& (f \circ (g \circ h)_{\leq n})_{\leq n} 
= (f \circ (g \circ h))_{\leq n} \\
&=& ((f \circ g) \circ h)_{\leq n} 
= ((f \circ g)_{\leq n} \circ h)_{\leq n} 
= (f\#g) \#h. 
\end{eqnarray*}
This shows that $(G_n,\#,\id_E)$ is a monoid. For $k < n$, omitting terms of order 
$> k$ leads to monoid homomorphisms $\phi_{kn} \: G_n \to G_k.$

To see that $(G_n,\#,\id_E)$ actually is a group, let 
$f \in G_n$ and assume that $g \in G_n$ is such that 
$g \circ f - \id_E$
vanishes up to terms of order $k$. 
Modulo terms of degree $k+1$, we find 
\begin{eqnarray*}
(g \circ f)_{\leq k} 
&=& (g_k \circ f)_{\leq k} + \sum_{j = 1}^{k-1} (g_j \circ f)_{\leq k} 
= g_k   + \sum_{j = 1}^{k-1} (g_j \circ f)_{\leq k}\\
&=& g_k + \id_E + \sum_{j = 2}^{k-1} (g_j \circ f)_{k}. 
\end{eqnarray*}
Therefore the requirement that $g \circ f - \id_E$ vanishes up to terms of order 
$k+1$ leads to the condition 
$$ g_k = - \sum_{j = 2}^{k-1} (g_j \circ f)_k. $$
Inductively, we thus find an element $g \in G_n$ with 
$g \# f = \id_E$, i.e., a left inverse of $f$.  

Next we construct a right inverse $h$  of $f$. 
Again, we assume that $k > 1$ and that 
$f \circ h - \id_E$
vanishes up to terms of order $k$. 
Modulo terms of degree $\geq k+1$, we then have 
\begin{eqnarray*}
(f \circ h)_{\leq k} 
&=& h_{\leq k} + \sum_{j = 2}^{k} (f_j \circ h)_{\leq k} 
= h_{\leq k} + \sum_{j = 2}^{k} 
(f_j \circ h_{\leq k-1})_{\leq k} \cr
&=& \id_E + h_k + \sum_{j = 2}^{k} (f_j \circ h_{\leq k-1})_k.
\end{eqnarray*}
Thus 
$$ h_k := - \sum_{j = 2}^{k} \big(f_j \circ h_{\leq k-1}\big)_k$$
is equivalent to 
$(f \circ h)_{\leq k} = \id_E.$
Inductively, we thus obtain a right inverse $h$ of $f$ in $G_n$. 
We now have 
$$ h = (g\#f)\#h = g \# (f \# h) = g, $$
i.e., left and right inverse coincide, which shows that $G_n$ is a group. 
Since the preceding construction shows that the inversion is a smooth map, 
we see that $G_n$ actually is a Lie group. 
Recall that the smoothness of the inversion 
in the finite-dimensional group $G_n$ follows from the Inverse 
Function Theorem (Exercise~\ref{exer:3.1.7}). 

For $\eta \in \g_n$ and $\xi \in P_n(E)$ we have 
\begin{eqnarray*}
(\id_E + \eta) \# (\id_E + \xi) 
&=& \id_E + \xi + \big(\eta \circ  (\id_E + \xi)\big)_{\leq n} \\
&=& \id_E + \xi + \eta 
= (\id_E + \xi) \# (\id_E + \eta), 
\end{eqnarray*}
so that $Z_n := \id_E + P_n(E)$ is central in $G_n$.  
Therefore $G_n$ is a central extension of 
$G_{n-1}$ by $Z_n$, so that we inductively see that each group $G_n$ 
is nilpotent. 

To see that $\g_n$ is the Lie algebra of $G_n$, 
we first make the composition in $G_n$ more explicit by using the 
Taylor formula:  
$$ \xi(x + \eta(x)) = \xi(x) + \dd \xi(x)(\eta(x)) + \ldots 
+ \frac{1}{n!} (\dd^n \xi)(x)(\eta(x), \ldots, \eta(x)). $$
For $f = \id_E + \xi$ and $g = \id_E + \eta$, this leads to 
\begin{eqnarray}\label{eq:5.2.prod} 
&& (f \circ g)(x)  
= x + \eta(x) + \xi(x + \eta(x)) \\
&=& x + \eta(x) + \xi(x) + \dd\xi(x)(\eta(x)) + \sum_{k = 2}^n \frac{1}{k!}  
(\dd^k \xi)(x)(\eta(x), \ldots, \eta(x)).\notag 
\end{eqnarray}
From the preceding formula, we already see that the left invariant 
vector fields on $G_n$ are given in 
$f = \id_E+ \xi$ by 
$$ \eta_l(f) = f.\eta = T_{\id_E}(\lambda_f)\eta 
= \dd f(\cdot)\eta = \eta + \dd\xi(\cdot)\eta. $$
Expanding $\xi = \sum_{j = 2}^n \xi_j$ and 
$\eta = \sum_{i=2}^n \eta_i$ into homogeneous summands of degree $j$, 
we obtain summands of the form 
$$(\dd^k \xi_j)(x)(\eta_{i_1}(x), \ldots, \eta_{i_k}(x)), $$
which are homogeneous of degree $j-k + i_1 + \ldots + i_k \geq j + k \geq 4$. 
From (\ref{eq:5.2.prod}) it immediately follows that the 
Taylor polynomial of degree $2$ of the multiplication in $(\id_E, \id_E)$ is given by 
$$ (\xi,\eta) \mapsto \xi + \eta + \dd\xi(\cdot) \eta. $$
This leads to the Lie bracket 
$$ [\xi,\eta](x) = \dd\xi(x)\eta(x) - \dd\eta(x)\xi(x) $$ 
(Remark~\ref{rem:brack-taylor}). 

We now discuss the exponential function 
$$ \exp_{G_n} \: \g_n \to G_n, $$
which is given by ``integrating'' a left invariant vector field $\xi_l$ on $G_n$ 
modulo terms of order 
$> n$. Since $G_n$ is a finite-dimensional nilpotent $1$-connected Lie group, 
we know a priory that its exponential function is a diffeomorphism because its 
universal covering group is $(\g_n,*)$, 
where $*$ is the BCH multiplication on $\g_n$ (Corollary~\ref{cor:to-existhom}). 
For $\xi \in \g_n$, we want to solve the initial value problem 
$$ \gamma'(t) = \gamma(t)\xi = \dd \gamma(t)(\cdot) \xi, \quad \gamma(0) = \id_E, $$
which leads to the Picard iteration: 
$$ \gamma_0 = \id_E, \quad 
\gamma_{k+1}(t) = \id_E + \int_0^t \dd \gamma_k(s)(\cdot)\xi\ ds. $$
Since the right hand side is linear in $\gamma_k$, it follows inductively that 
$\gamma_k(t)$ is a polynomial of degree $\leq k$ in $t$. 
For $\gamma_k(t) = \sum_{j = 0}^k a_j^k t^j$, we obtain $a_0^{k}= \id_E$ for each $k$, and 
the recurrence relation 
\[ a_{j+1}^{k+1} = \frac{1}{j+1} \dd a_j^k(\cdot) \xi \quad \mbox{ for } \quad 
j \leq k, \] 
which shows that $a_j^k$ is obtained directly from $a_0^{k-j} = \id_E$. In particular, 
we see that 
$a_j^k = a_j^j$ for $k \geq j$. Omitting terms of order $\geq n$, we further see that 
Picard iteration leads to a fixed point after $n$ steps: 
$\gamma(t) = \sum_{j = 0}^n a_j t^j$ with 
\[ a_0 = \id_E, \ a_1 = \xi, \ a_2 = \shalf\dd\xi(\cdot)  \xi \quad \mbox{ and } \quad 
a_{j+1} = \frac{1}{j+1}\dd a_j(\cdot)  \xi, \quad j < n.\]

We now form the projective limit group $G := \prolim G_n \cong \id_E + \g$ 
whose manifold structure is obtained from the structure as an affine space 
with translation group $\g$. Since the exponential functions of the groups $G_n$ are 
compatible with the limiting process, we see that $G$ is an exponential 
Lie group with a pro-nilpotent Lie algebra. The group $G$ can be identified 
with the set of all formal diffeomorphisms of $E$ fixing $0$ and with 
first order term given by $\id_E$. Likewise $\g$ can be identified with 
the Lie algebra of formal vector fields on~$E$ with vanishing constant term.
In the terminology of Subsection~\ref{subsec:sternberg}, we have 
$G \cong \Gf_n(\R)$ for $E \cong \R^n$. 
\end{ex}

\begin{small}
\subsection*{Exercises for Section~\ref{sec:5.4}} 

\begin{exer} \mlabel{exer:3.6.10} Show that each nilpotent Lie algebra 
$\g$ is locally finite, i.e., each finite subset of $\g$ generates 
a finite-dimensional Lie subalgebra. 
\end{exer}

\begin{exer} \mlabel{exer:5.4.1} Let $\cA$ be a unital real algebra, 
$\cA_n$ the set of nilpotent elements of $\cA$ and 
$\cA_u := \1 + \cA_n$ be the set of unipotent elements of $\cA$. Consider the maps 
$$ \Exp \: \cA_n \to \cA, \ \  x\mapsto \sum_{n = 0}^\infty \frac{x^n}{n!}, \quad 
\Log \: \cA_u \to \cA, \ \  \1 - x\mapsto -\sum_{n = 1}^\infty \frac{x^n}{n}. $$
Show that: 
\begin{enumerate}
\item[\rm(1)] $\Exp(\cA_n) \subeq \cA_u$ and $\Log(\cA_u) \subeq \cA_n$: 
Hint: $\Exp(x) - \1 = xy$ and \break $\Log(\1 -x) = x z$ for some $y$ and $z$, commuting 
with $x$. 
\item[\rm(2)] $\Exp \circ \Log = \id_{\cA_u}$. 
Hint: Lemma~\ref{lem:nilp-compo} and $\exp(\log X)) = X$ in $\R[[X]]$.  
\item[\rm(3)] $\Log \circ \Exp = \id_{\cA_n}$. 
Hint: Lemma~\ref{lem:nilp-compo} and $\log(\1 + (\exp X - \1)) = X$ in $\R[[X]]$.  
\end{enumerate}
\end{exer} 

\begin{exer} \mlabel{exer:5.4.2} Let $X$ be a set and $\cA(X)$ denote the free associative 
$\K$-algebra over $X$, i.e., the tensor algebra over the free vector space $\K^{(X)}$. 
We write $\delta\: X \to \cA$ for the canonical map. Show that 
the Lie subalgebra $L(X)\subeq \cA(X)$ generated by $\delta(X)$ is a free 
Lie algebra over $X$, i.e., for each map $\alpha \: X \to \g$ into a Lie algebra 
$\g$, there exists a unique homomorphism of Lie algebras $\hat\alpha \: L(X) \to \g$ 
with $\hat \alpha \circ \delta = \alpha$.
Hint: Use the embedding $\g \into U(\g)$ into the enveloping algebra 
(Poincar\'e--Birkhoff--Witt Theorem). 
\end{exer} 

\begin{exer} \mlabel{exer:5.4.3} Each Banach--Lie algebra $\g$ possesses a compatible 
norm.\\ Hint: If $\|[x,y]\| \leq C \|x\|\|y\|$ for some $C> 0$. 
\end{exer}

\begin{exer} \mlabel{exer:5.6.4} We consider the special case of 
Example~\ref{ex:5.2.formal}, where $E = \R$. 
We then identify $\g_n$ with the space of polynomials 
$\g_n = \Spann \{ x^2,\ldots, x^n\}$ and the group $G_n$ with the affine space 
$G_n := x + \g_n$. The group structure on $G_n$ is given by composition 
and then omitting all terms of order $> n$. 
For $n = 2,3,4$ and $5$: 
\begin{enumerate}
\item[\rm(1)]  Give an explicit formula for the composition in the group $G_n$. 
\item[\rm(2)]  Give an explicit formula for the exponential function 
$\exp_{G_n} \: \g_n \to G_n$. 
\item[\rm(3)]  Show that the exponential function of $G_n$ can be obtained 
by  truncated Picard iteration. 
\end{enumerate}
\end{exer} 

\begin{exer} \mlabel{exer:2.6.10} 
Show that every finite subset of a nilpotent Lie algebra generates 
a finite-dimensional subalgebra (cf.\ Exercise~\ref{ex:e.9}). 
\end{exer}
 
\end{small}


\section{Notes and comments on Chapter~\ref{ch:5}} 

{\bf\ref{sec:5.1}:} The concept of a locally exponential Lie group generalizes in a natural 
way the concept of a Banach--Lie group and, as we have seen in this chapter, 
many important properties of Banach--Lie groups carry over to locally exponential 
Lie groups modeled on locally convex spaces. 
In Robart's paper \cite{Rb97}, locally exponential Lie groups 
are called groups of {\it type one}, 
but we prefer the more instructive term ``locally exponential.''

\nin{\bf\ref{sec:5.5}:} The concept of a locally $k$-exponential Lie group (cf.\ Theorem~\ref{thm:int-thm-lock}) 
is inspired by J.~Mical's papers \cite{Mi92, Mi94}, where he even considers infinite 
products of exponential factors; see also \cite[\S 4.1]{Rob02}. 

For (local) Banach--Lie groups, the Automatic Smoothness Theorem 
can already be found in \cite{Bir38}, 
and for BCH--Lie groups (Definition~\ref{def:VI.1.9})
 in \cite{Gl02c} (see also \cite{Mil84} 
for the statement without proof). 
The special case of one-parameter groups $\R \to \cA^\times$ 
of the unit group of a Banach algebra $\cA$ 
is due to Nagumo \cite{Nag36} and Nathan \cite{Nat35}. 
In the context of topological groups, there exist 
``automatic continuity'' results, deriving continuity of group 
homomorphisms from measurability (see \cite{KR07, Ro09, Ro19}, \cite{NP00})

The Integration Theorem is due to Maissen for Banach--Lie groups 
\cite[Satz~10.3]{Ms62}. 

As the non-locally exponential semidirect products 
$E \rtimes_\alpha \R$ from Example~\ref{ex:bad-expfct} 
show, local exponentiality is not an extension property. 
However, one can show that local exponentiality is at least preserved 
under central extensions (Section~\ref{sec:11.3}). 

\nin{\bf\ref{sec:5.4}:} In \cite{SchF90a}, F.~Schur derived recursion formulas for the summands of 
the series describing the multiplication of a Lie group in canonical 
coordinates (i.e., in an exponential chart). He also proved the local convergence 
of the series given by this recursion relations, which can in turn be used to obtain 
the associativity of the BCH multiplication. His approach is quite 
close to our treatment of locally exponential Lie algebras in the sense that 
he derived the series from the Maurer--Cartan form by integration of 
a partial differential equation of the form 
$f^*\kappa_\g = \kappa_\g$ with $f(0) = x$, whose unique solution 
is the left multiplication $f = \lambda_x$ in the local group. 

The Hausdorff series was made more explicit by {  Campbell} in \cite{Cam97, 
Cam98} 
and in \cite{Hau06} F.~Hausdorff approached the Hausdorff series on a formal level, 
showing that the formal expansion of $\log(e^x e^y)$ 
can be expressed in terms of Lie polynomials. Parts of his results 
had been obtained earlier by Baker \cite{Bak01, Bak05}. See \cite{Ei68} for a more 
recent short argument that all terms in the Hausdorff series are Lie brackets. 
For an extensive discussion of various aspects of the 
Hausdorff series and the BCH formula, we refer to the nice 
monogrpaph \cite{BF12} by Bonfiglioli and Fulci. 
For a new approach to the $y$-homogeneous terms in $\log(e^x e^{2y} e^x)$, 
we refer to \cite{ML18}. 

Our estimates in the discussion of the exponential local Lie group 
associated to a Banach--Lie algebra in Theorem~\ref{thm:loc-bangrp} are not 
optimal. In \cite[Ch.~2, \S 7.2, Prop.~1/2]{Bou89}, it is shown that triple 
products are associative under the weaker condition 
$\|x\| + \|y\| + \|z\| < \log\frac{3}{2}$.

\chapter{Local theory of Lie groups} 
\mlabel{ch:7}

To develop the theory of locally exponential Lie groups 
further, we have to take a closer look at 
local Lie groups. This leads in particular to the concept 
of an {\it exponential local Lie group}. \index{local Lie group!exponential} 
The global and the local concepts of ``local exponentiality'' 
are related by the fact that each 
locally exponential Lie group contains arbitrarily small symmetric 
$\be$-neighborhoods which are 
exponential local Lie groups. The Lie algebras of exponential local Lie groups 
are called {\it locally exponential}, and these are the natural candidates 
for  Lie algebras of locally exponential groups. The local convergence 
of the Hausdorff series in a Banach--Lie algebra 
implies that every Banach--Lie algebra is locally exponential. 

The natural integrability problem 
in the context of locally exponential groups is: 
Which locally exponential Lie algebra is the Lie algebra 
of a locally exponential Lie group? We thus call it {\it enlargeable}. 
Since \index{enlargeable}  
each locally exponential Lie algebra belongs to some exponential local 
Lie group, the preceding question concerns a passage from a local 
object to a global one. In this sense it rather is an enlargeability  
than an integrability problem, hence considerably more tractable. 
In the present chapter, we  
show in particular that, for each locally exponential Lie algebra $\g$, 
the quotient $\g/\z(\g)$ is enlargeable. 
Based on this result, we shall derive in Section~\ref{sec:11.4} 
a general criterion 
for the integrability of a locally exponential Lie algebra in terms of 
period groups.

\section{Local Lie groups} \mlabel{sec:3.7} 

To any Lie group $G$ we have attached its Lie algebra $\L(G)$, which is 
an infinitesimal structure attached to the global object $G$. 
In numerous situations it turns out to be fruitful to have an intermediate 
local level between the infinitesimal and the global one. 
This leads to the concept of a local Lie group. Every open symmetric 
$1$-neighborhood of a Lie group is a local Lie group, but there 
are local Lie groups which do not embed into global 
Lie groups (see Section~\ref{sec:6.4}). 
Nevertheless, we may associate to any local Lie group $G$ 
a Lie algebra $\L(G)$ in a functorial fashion. The passage between a 
local Lie group and its Lie algebra will be of particular importance 
for the special class of exponential local Lie groups, the local variants 
of locally exponential Lie groups (cf.\ Chapter~\ref{ch:5}). 

\subsection{Local Lie groups and their Lie algebras}

First we recall the abstract concept of a local group. 

\begin{defn} \mlabel{def:3.1.1b} 
Let $G$ be a set, $D \subeq G \times G$ a subset, $\be \in G$, and 
$$m_G \: D \to G, (x,y) \mapsto x * y$$ 
a map. We say that the product $x*y$ is {\it defined} if $(x,y) \in D$. 
We call $G$, resp.\ $(G,D,m_G,\be)$, a {\it local group} if the following conditions are
satisfied: \index{local group}
\begin{description}
\item[\rm(Loc1)] Assume that $(x,y),(y,z) \in D$. If $(x*y,z)$ or
$(x,y*z) \in D$, then both are contained in $D$ and 
$(x*y)*z = x*(y*z).$ 
\item[\rm(Loc2)] For each $x \in G$, we have $(x,\be),(\be,x) \in D$ and 
$x*\be = \be*x =x$. 
\item[\rm(Loc3)] For each $x \in G$, there is a element $x^{-1} \in G$ with 
$(x,x^{-1}),(x^{-1},x)\in D$ and $x*x^{-1} = x^{-1} * x = \be$. 
\item[\rm(Loc4)] If $(x,y)\in D$, then $(y^{-1}, x^{-1}) \in D$. 
\end{description}
\end{defn}

Many of the familiar properties of groups carry over to local Lie groups. 

\begin{rem} \mlabel{rem:locgrp0} Let $(G,D,m_G,\be)$ be a local group. \\
(a) The inverse element required in (Loc3) is unique: 
Suppose that $x,y,z \in G$ satisfy 
$(x,y),(z,x) \in D$ with 
$x * y = \be = z * x$, i.e., $y$ is a right inverse of $x$ and 
$z$ a left inverse. 
Then 
$$z = z * \be = z * (x * y) = (z * x) * y = \be * y = y$$ 
shows that $y = z$ is the unique inverse of $x$.  

(b) Suppose that $(x,y) \in D$. 
We always have $(y,y^{-1}) \in D$ and \break $(x,y * y^{-1}) = (x,\be) \in D$, so that 
$$ (x*y)*y^{-1} = x*(y*y^{-1}) = x * \be = x $$
follows from (Loc1). This further leads to 
$$ (x*y)*(y^{-1}*x^{-1}) 
= ((x*y)*y^{-1}) * x^{-1}
= x* x^{-1} = \be, $$
so that (a) entails 
\begin{eqnarray}
  \label{eq:4.1.1}
(x *y)^{-1} = y^{-1} * x^{-1}. 
 \end{eqnarray}

(c) Uniqueness of idempotents: Let $x \in G$ with 
$(x,x) \in D$ and $x*x = x$. Then 
$$ \be = x * x^{-1} = (x*x)*x^{-1} = x * (x* x^{-1}) = x * \be = x, $$
so that $\be$ is the only idempotent in $G$. 
\end{rem}

\begin{defn} A {\it local Lie group} \index{local Lie group} 
is a local group $(G,D,m_G,\be)$, for which 
$G$ is a locally convex manifold, $D$ is an open subset of $G \times G$, 
and multiplication $m_G \: D \to G$ and inversion 
$\eta_G \: G \to G, x \mapsto x^{-1}$ are smooth. 

If $(G,D_G,m_G,\be_G)$ and $(H,D_H,m_H,\be_H)$ are local Lie groups, then 
a smooth map $\phi \: G \to H$ is called a 
\index{ morphism of local Lie groups} 
{\it morphism of local Lie groups} if 
$\phi(\be_G) = \be_H$, $(\phi \times \phi)(D_G)\subeq D_H$, and 
$m_H \circ (\phi \times \phi) = \phi \circ m_G$, i.e.,  
$$ \phi(x*y) = \phi(x) * \phi(y) \quad \mbox{ for } \quad (x,y) \in D_G. $$ 
\end{defn} 

The following remark justifies the terminology: 

\begin{rem} \mlabel{rem:locsubgroup} (a) Let $G$ be a group and $U \subeq G$ a symmetric subset
containing the identity element $\be$. Then 
$(U,D_U,m_U,\be)$ is a local group if we put 
$$ D_U := \{ (x,y) \in U \times U \: x * y \in U\} 
\quad \mbox{ and } \quad m_U := m_G\res_{D_U}.$$ 

If, more generally, $(G,D,m_G,\be)$ is a local group and $U \subeq G$ a symmetric subset
containing $\be$, then 
$(U,D_U,m_U,\be)$ is a local group for 
$$ D_U := \{ (x,y) \in D \cap (U \times U) \: x * y \in U\} 
\quad \mbox{ and } \quad m_U := m_G\res_{D_U}$$ 
(Exercise~\ref{exer:3.7.1}). 

(b) If, in addition, $G$ is a Lie group and $U \subeq G$ is open, 
then $(U,D,m_U, \be)$ is a local Lie group. 
\end{rem} 

From Theorem~\ref{thm:locglob} we derive immediately the following result 
on the existence of a Lie group structures on groups containing local Lie groups. 

\begin{thm} \mlabel{thm:locinglob} 
Let $(U,D_U,m_U,\be)$ be a local Lie group and suppose that 
$G$ is a group for which we have an injection 
$\iota \: U \to G$ of local groups such that 
\[ D_U = \{ (x,y) \in U \: \iota(x)\iota(y) \in \iota(U)\}.\] 
Then $G$ carries a unique Lie group structure for which 
the inclusion $U \into G$ is a diffeomorphism onto an open subset of $G$. 
\end{thm}

We now turn to the definition of the Lie algebra of a local Lie group. Although 
the local concepts require more caution than the global ones because the 
multiplication is not everywhere defined, the key facts from the 
global context carry over to the local situation. 

\begin{defn} {\rm(The Lie algebra of a local Lie group)} 
Let $G$ be a local Lie group. For $x \in G$ define the left multiplication 
$$ \lambda_x \: D_x := \{ y \in G \: (x,y) \in D\} \to G, \quad y \mapsto x * y $$
and the right multiplication 
$$ \rho_y \: D^y := \{ x \in G \: (x,y) \in D\} \to G, \quad x \mapsto x * y. $$
The domains $D_x$, resp., $D^y$ are open subsets of 
$G$ containing $\be$ and $x^{-1}$, resp.,~$y^{-1}$. 

Suppose that $(x,y) \in D$, so that $x*y$ is defined. 
From (Loc1) we derive that if $z \in D_{x*y} \cap D_y$, then 
$\lambda_{x*y}(z) = \lambda_x(\lambda_y(z))$. 
In particular we have $\lambda_{x*y}= \lambda_x \circ \lambda_y$ on a neighborhood of $\be$, 
which leads to  
\begin{eqnarray}
  \label{eq:4.1.2}
T_\be(\lambda_{x*y}) = T_y(\lambda_x) \circ T_\be(\lambda_y). 
\end{eqnarray}
We likewise obtain 
\begin{eqnarray}
  \label{eq:4.1.3}
T_\be(\rho_{y*x}) =  T_y(\rho_x) \circ T_\be(\rho_y). 
\end{eqnarray}

For each $x \in \g = T_\be(G)$, we now define the corresponding 
\index{left invariant vector field, on local Lie group} 
{\it left invariant vector field} $x_l \in {\cal V}(G)$ by 
$$ x_l(g) := T_\be(\lambda_g)(x). $$
For $h \in D_g$ we then have 
$$ x_l(gh) 
= T_\be(\lambda_{gh})(x)
= T_h(\lambda_{g})(x)T_\be(\lambda_{h})(x)
= T_h(\lambda_{g})(x) x_l(h), $$
i.e., 
$$ x_l \circ \lambda_g = T(\lambda_g) \circ x_l\res_{D_g}. $$
This means that the vector field $x_l$ is $\lambda_g$-related to $x_l\res_{D_g}$. 
This implies that for $x,y \in \g$ the Lie bracket $[x_l, y_l]$ is also 
$\lambda_g$-related to $[x_l,y_l]\res_{D_g}$ (Lemma~\ref{larelglob}). 
It follows in particular that 
this vector field is left invariant:
\[  [x_l,y_l](g) = T_\be(\lambda_g)[x_l,y_l](\be) \quad \mbox{ for } \quad g \in G.\] 
We conclude that the set 
${\cal V}(G)^l$ of left invariant vector fields on $G$ is a Lie subalgebra of the 
Lie algebra ${\cal V}(G)$ of all smooth vector fields on $G$. Since the map 
$$ \g \to {\cal V}(G)^l, \quad x \mapsto x_l $$
is a linear isomorphism with inverse $\ev_\be\: X \mapsto X(\be)$, we obtain a Lie bracket on $\g$ by 
$$ [x,y] := [x_l, y_l](\be) $$
(cf.\ Exercise~\ref{exer:4.1.2}). 
By definition, this Lie bracket satisfies $[x_l, y_l] = [x,y]_l$, 
and we shall see in Remark~\ref{rem:locgrp} below that it is continuous. 
The Lie algebra 
$$ \L(G) := \L(G,D,m_G,\be)  := (\g,[\cdot,\cdot]) $$
is called the 
{\it Lie algebra of the local Lie group $G$}. 
\index{Lie algebra!of a local Lie group}
\end{defn}

\begin{rem} {\rm(The Lie bracket and Taylor expansions)} \mlabel{rem:locgrp} 

Let $E$ be a locally convex space and $G \subeq E$ an open 
$0$-neighborhood carrying the structure of a local Lie group 
$(G,D_G,m_G,0)$. For $(x,y) \in D$ we write $x * y = m_G(x,y)$. 

We consider the Taylor expansion of the multiplication in $(0,0) 
\in E \times E$. In view of $0 * 0 =0$, the constant term vanishes, so that 
\[  x * y = b_1(x,y) + b_2(x,y) + b_3(x,y) + \ldots, \quad \mbox{ where } \quad 
b_k = \frac{1}{k!} \delta^k_{(0,0)} m_G \]
are continuous homogeneous polynomials of degree $k$. 
With exactly the same arguments as in Remark~\ref{rem:brack-taylor}, 
we see that $b_1(x,y)= x+y$, and that $b_2$ is bilinear with 
$$ [x,y] = b_2(x,y) - b_2(y,x). $$

For the Taylor expansion $\eta(x) = x^{-1} = s_1(x) + s_2(x) +\cdots$ of 
the inversion map in $0$, we also derive 
\begin{eqnarray}
  \label{eq:4.l.4}
s_1(x) = -x \quad \hbox{ and } \quad s_2(x) = -b_2(x,s_1(x)) =
b_2(x,x) 
\end{eqnarray}
with the same reasoning as in Remark~\ref{rem:brack-taylor}. 
\end{rem}

We now turn to homomorphisms of local Lie groups.

\begin{prop} If $\phi \: (G,D_G,m_G,\be) \to (H,D_H,m_H,\be)$ is a morphism 
of local Lie groups, then 
$$ \L(\phi) := T_{\be}(\phi) \: \L(G) \to \L(H) $$
is a continuous homomorphism of Lie algebras. 
\end{prop}

\begin{prf} Let $x \in \L(G)$ and $x_l(g) = T_\be(\lambda_g)(x)$ be 
the corresponding left invariant 
vector field on $G$. We have $\phi((D_G)_g) \subeq (D_H)_{\phi(g)}$, and 
on $(D_G)_g$ we have 
\[ \phi \circ \lambda_g = \lambda_{\phi(g)} \circ \phi \quad \mbox{ for } \quad g \in G.\] 
This implies that 
$$ T_g(\phi)x_l(g) = T_g(\phi) T_\be(\lambda_g)x  
= T_\be(\lambda_{\phi(g)}) T_\be(\phi)(x)
= \big(\L(\phi)(x)\big)_l(\phi(g)), $$
so that the vector fields $x_l$ on $D_g$ 
and $\big(\L(\phi)x\big)_l$ on $D_{\phi(g)}$ are $\phi$-related. 
We conclude with Lemma~\ref{larelglob} 
that for $x,y \in \L(G)$ the vector fields 
$[x_l,y_l]$ and 
$[(\L(\phi)x)_l,(\L(\phi)y)_l]$ 
are also $\phi$-related on $D_g$. Evaluating in $\be$, we obtain 
$$ [\L(\phi)x,\L(\phi)y] = [\L(\phi)x,\L(\phi)y]_l(\be) 
= T_\be(\phi)[x_l,y_l](\be) = \L(\phi)[x,y]. $$
Hence $\L(\phi)$ is a homomorphism of Lie algebras. Its continuity follows from 
the smoothness of $\phi$. 
\end{prf} 

We record that the Lie algebra of a local Lie group only depends on its ``germ in $\be$'': 
\begin{cor} \mlabel{cor:4.l.1} If $(G,D_G,m_G,\be)$ is a local Lie group, 
  $U \subeq G$ an open symmetric $\be$-neighborhood,
  $D_U := m_G^{-1}(U) \cap (U \times U)$, and $m_U := m_G\res_{D_U}$, then 
$(U,D_U,m_G,\be)$ is a local Lie group and the inclusion 
$i_U \: U \to G$ induces an isomorphism 
$$ \L(i_U) \: \L(U) \to \L(G) $$
of locally convex Lie algebras. 
\end{cor}

\begin{rem} As a direct consequence of the Chain Rule, we get 
$$\L(\id_G)  = \id_{\L(G)} \quad \mbox{ and } \quad 
\L(\phi_1 \circ \phi_2) = \L(\phi_1) \circ \L(\phi_2)$$ 
for two composable morphisms $\phi_1$ and $\phi_2$ of local Lie groups. Therefore 
$\L$ defines a functor from the category of local Lie groups to 
the category of locally convex Lie algebras. 
\end{rem}

\begin{defn}
We call a locally convex Lie algebra $\g$ 
\index{Lie algebra!locally  integrable} 
{\it locally  integrable} if 
there exists a local Lie group $G$ with $\L(G) \cong \g$. 
\end{defn}

The concept of local integrability splits the integrability problem into 
two parts: to decide if a Lie algebra is integrable, we may first ask 
for local integrability and then for an embedding of a corresponding local 
group into a global one (the {\it enlargibility problem}) 
(cf.\ Theorem~\ref{thm:locinglob}).

\begin{defn}\mlabel{defn:adjrep-locgrp} {\rm(The adjoint representation)} Let $G$ be a local Lie group. 
For $(x,y)$, $(y,x^{-1})$ and $(xy,x^{-1}) \in D$, we put 
$$ c_x(y) := (xy)x^{-1} = x(yx^{-1}).$$ 
Then $c_x$ is a smooth map defined on some $\be$-neighborhood of $\be$ and 
$c_x(\be) = \be$. 

We claim that if $(x,y) \in D$, then 
\begin{eqnarray}\label{eq:4.l.5} 
c_x \circ c_y = c_{xy} 
\end{eqnarray}
holds on a neighborhood of $\be$. 
In fact, if  $z \in G$ is sufficiently close to $\be$, we have 
$(x,yz),(y,z),(z,y^{-1}),(yz,y^{-1}), (x, c_y(z)), (c_y(z),x^{-1}), 
(xc_y(z), x^{-1}) \in D$, and this leads to 
\begin{eqnarray*}
c_x(c_y(z)) 
&=& \big(x((yz)y^{-1})\big)x^{-1} = \big((x(yz)) y^{-1}\big)x^{-1} 
=\big(((xy)z) y^{-1}\big)x^{-1} \\
&=& ((xy)z)(y^{-1}x^{-1}) =((xy)z)(xy)^{-1} = c_{xy}(z).
\end{eqnarray*}

Fix $x \in G$. If $z_1, z_2$ are sufficiently close to $\be$, 
both $c_x(z_1)$ and $c_x(z_2)$ are defined, and we have 
\begin{eqnarray*}
 c_x(z_1 z_2) 
&=& (x(z_1z_2))x^{-1} 
= ((xz_1)z_2)x^{-1} 
= ((xz_1)(x^{-1}(x z_2)))x^{-1} \\
&=& \big(c_x(z_1) (xz_2)\big)x^{-1} 
= c_x(z_1) c_x(z_2). 
\end{eqnarray*}
Hence there exists some symmetric $\be$-neighborhood $U_x \subeq G$, such that 
$c_x \: U_x \to G$ is a morphism of local Lie groups, and 
Corollary~\ref{cor:4.l.1} implies that $\L(c_x) \: \g \cong \L(U_x) \to \g$ 
is an automorphism of Lie algebras.

\index{adjoint representation} 
The {\it adjoint representation of $G$} is defined by 
$$ \Ad \: G \to \Aut(\g), \quad g \mapsto \L(c_g) 
= T_\be(\rho_g)^{-1} \circ T_\be(\lambda_g) 
= T_{g^{-1}}(\lambda_g) \circ T_\be(\rho_{g^{-1}}).  $$ 
Note that (\ref{eq:4.l.5}) implies $\Ad(gh) = \Ad(g)\Ad(h)$ for 
$(g,h) \in D$. Since we also have $\Ad(\be) = \id_\g$, 
$\Ad$ is homomorphism of local groups, if we consider 
the group $\Aut(\g)$ as a local group. 
\end{defn}

\begin{rem} Let $(G,D,m,0)$ be a local Lie group, 
where $G$ is an open $0$-neighborhood in some 
locally convex space. 

For the Taylor polynomial of order $2$ of the conjugation map, we obtain 
with the Chain Rule for Taylor Polynomials (Proposition~\ref{chainrtay}) 
and Remark~\ref{rem:locgrp}: 
\begin{eqnarray*} \label{eq:conjug-loc}
(x * y) * x^{-1} 
&=& (x + y + b_2(x,y) + \cdots) * (-x + b_2(x,x) + \cdots) \\ 
&=& (x + y -x) + b_2(x,y) + b_2(x,x) + b_2(x + y,-x) + \cdots \\ 
&=& y + b_2(x,y) - b_2(y,x) + \cdots \\ 
&=& y + [x,y] + \cdots 
\end{eqnarray*}
For the commutator we therefore get 
\begin{eqnarray*} \label{eq:commut-loc}
(x * y * x^{-1}) * y^{-1} 
&&= (y + [x,y] + \cdots) * (-y + b_2(y,y) + \cdots) \cr
&&=  [x,y] + b_2(y,y) + b_2(y,-y)  + \cdots \cr
&&=  [x,y] + \cdots   
\end{eqnarray*}
\end{rem}

\subsection{The Maurer--Cartan form and its applications} 

\begin{defn} \mlabel{defn:mcart-locgrp} Let $(G,D,m_G,\be)$ be a local Lie group. 
We define the \index{Maurer--Cartan form} 
{\it Maurer--Cartan form} $\kappa_G \in \Omega^1(G,\g)$ by 
$$ (\kappa_G)_g 
:= T_g(\lambda_{g^{-1}})= T_\be(\lambda_g)^{-1} $$
(cf.\ (\ref{eq:4.1.2})).
This means in particular, that for each left invariant 
vector field $x_l$, the function $\kappa_G(x_l) = x$ is constant.

Let $M$ be a smooth connected manifold and 
$f \: M \to U$ be a smooth function. We define the \index{logarithmic derivative} 
{\it logarithmic derivative} of $f$ 
as in Definition~\ref{def:c.11} by 
$$\delta(f) := f^*\kappa_G \in \Omega^1(M,\g). $$ 
Pointwise this means that 
$$ \delta(f)_m = (\kappa_G)_{f(m)} \circ T_m(f) =
 T_\be(\lambda_{f(m)})^{-1} \circ T_m(f). $$ 
\end{defn}

\begin{lem} \mlabel{lem:4.l.14} If $G$ is a local Lie group, $g \in G$, and 
$\lambda_g \: D_g \to G$ the corresponding left multiplication, 
then we have on $D_g$: 
$$ \lambda_g^*\kappa_G = \delta(\lambda_g) = \kappa_G\res_{D_g}. $$
\end{lem}

\begin{prf} For $h \in D_g$ we get with (\ref{eq:4.1.2}) 
$$ \delta(\lambda_g)_h 
= T_{\be}(\lambda_{gh})^{-1} T_h(\lambda_g) 
= T_{\be}(\lambda_{h})^{-1} T_h(\lambda_{g})^{-1} T_h(\lambda_g) 
= T_{\be}(\lambda_{h})^{-1}  = (\kappa_G)_h. $$
\end{prf}

\begin{lem} \mlabel{lem:4.l.15} 
Let $(G,D_G,m_G,\be)$ be a local Lie group with Lie algebra $\g$ 
and $M$ a smooth connected manifold. 
For $f,g \in C^\infty(M,G)$ we put  
\[ M_{f,g} := (f, g)^{-1}(D_G)\quad \mbox{  and define } \quad 
f \cdot g := m_G \circ (f, g) \: M_{f,g} \to G.\] 
For $\alpha \in \Omega^1(M,\g)$ we define 
$\Ad(f)\alpha \in \Omega^1(M,\g)$ by  
\[ (\Ad(f)\alpha)_m := \Ad(f(m)) \circ \alpha_m \quad \mbox{  for } \quad m \in M.\] 
Then the following assertions hold: 
\begin{enumerate}
\item[\rm(i)] $\delta(f \cdot g) = \delta(g) + \Ad(g)^{-1}\delta(f)$ on 
$M_{f,g}$ {\em(Product Rule).} \index{Product Rule} \index{Quotient Rule} 
\item[\rm(ii)] $\delta(f^{-1}) = - \Ad(f)\delta(f)$ {\em(Quotient Rule).} 
\end{enumerate}
\end{lem} 

\begin{prf} (i) From the Chain Rule and \eqref{eq:4.1.2} we get 
\begin{eqnarray*}
&& \delta(f \cdot g)_m 
= T_\be(\lambda_{f(m)g(m)})^{-1} T_m(f\cdot g) \\
&=& T_\be(\lambda_{g(m)})^{-1} T_{g(m)}(\lambda_{f(m)})^{-1} 
\big(T_{g(m)}(\lambda_{f(m)}) T_m(g) + T_{f(m)}(\rho_{g(m)}) T_m(f)\big) \\ 
&=& T_\be(\lambda_{g(m)})^{-1} T_m(g) + 
T_\be(\lambda_{g(m)})^{-1} T_{g(m)}(\lambda_{f(m)})^{-1} T_{f(m)}(\rho_{g(m)})T_m(f) \\ 
&=& \delta(g)_m 
+ T_\be(\lambda_{g(m)})^{-1} T_\be(\rho_{g(m)}) T_\be(\lambda_{f(m)})^{-1} \circ T_m(f) \\ 
&=& \delta(g)_m 
+ \Ad(g(m))^{-1} \circ \delta(f)_m. 
\end{eqnarray*}

(ii) In view of $M_{f,f^{-1}} = M$, (i) implies 
$$ 0 = \delta(\be) = \delta(f \cdot f^{-1})
= \delta(f^{-1}) + \Ad(f)\delta(f), $$
which leads to 
$\delta(f^{-1}) = - \Ad(f)\delta(f).$
\end{prf}

\begin{lem} \mlabel{lem:uni-lemma} {\rm(Uniqueness Lemma)} 
Let $(G,D_G,m_G,\be)$ be a local Lie group with Lie algebra $\g$ 
and $M$ a smooth connected manifold. 
If $f,g \in C^\infty(M,G)$ satisfy 
$(g,f^{-1})(M) \subeq D_G$ and $\delta(f) = \delta(g)$, 
then $g =  \lambda_x \circ f$ holds for some $x \in G$. 
\end{lem}

\begin{prf} By assumption, we have $M_{g,f^{-1}} = M$, so that the 
two formulas from Lemma~\ref{lem:4.l.15} yield 
$$ \delta(g\cdot f^{-1}) 
= \delta(f^{-1}) +  \Ad(f)\delta(g) 
= - \Ad(f)\delta(f) +  \Ad(f)\delta(f)= 0. $$
Hence $T_m(g \cdot f^{-1}) = 0$ in each $m \in M$, so that $g \cdot f^{-1}$ is constant 
equal to some $x \in G$. This means that $g(m) = x \cdot f(m)$ for all $m \in M$. 
\end{prf}

\begin{prop} \mlabel{prop:lochom} Let 
$(G_1, D_1, m_1, \be)$ and 
$(G_2, D_2, m_2, \be)$ be local Lie groups, 
let $\psi \: \L(G_1) \to \L(G_2)$ be a continuous homomorphism of Lie algebras and 
$f \: G_1 \to G_2$ a smooth map with $\delta(f) = \psi \circ \kappa_{G_1}$ 
and $f(\be) = \be$. 
Then there exists an open $\be$-neighborhood $U \subeq G_1$ with 
$U \times U \subeq D_1$, and 
$$ f(x * y) = f(x) * f(y) \quad \mbox{ for } \quad x,y \in U.$$
\end{prop}

\begin{prf} Let $E \subeq G_2$ be an open symmetric $\be$-neighborhood 
so small that $E \times E \subeq D_2$ and 
$U \subeq G_1$ a connected open $\be$-neighborhood 
with $U \times U \subeq D_1$ and $f(U * U) \subeq E$. 

Fix $x \in U$. We consider the function 
$$ h := \lambda_{f(x)}^{-1} \circ f \circ \lambda_x \: U \to E, \quad 
y \mapsto f(x)^{-1}f(x *y) $$
and claim that $h = f\res_U$. 
In fact, we get with Lemma~\ref{lem:4.l.14} 
$$ \delta(h) = h^*\kappa_{G_2} 
= \lambda_x^*f^*\lambda_{f(x)^{-1}}^*\kappa_{G_2} 
= \lambda_x^*f^*\kappa_{G_2} 
= \lambda_x^*(\psi \circ \kappa_{G_2})
= \psi \circ \kappa_{G_2} = \delta(f). $$
Since $U$ is connected and $f(\be) = h(\be) = \be$, we derive from the 
Uniqueness Lemma~\ref{lem:uni-lemma} $h = f\res_U$, and this leads to 
$f(x*y) = f(x)*f(y)$ for $x,y \in U$. 
\end{prf}

\subsection{Equivariant differential forms on a local Lie group} 

\begin{defn} \mlabel{def:equiv-form} Let $(G,D,m_G,\be)$ be a local Lie group 
and $E$ a locally convex space. We further assume that 
we are given a homomorphism of local groups 
$\rho_E \: G \to \GL(E)$ such that the map 
$G \times E \to E, (g,v) \mapsto g.v := \rho_E(g)v$ is smooth. 
Then we call $E$ a {\it smooth $G$-module}. 
\index{smooth $G$-module} 

We call a $p$-form $\alpha \in \Omega^p(G,E)$ 
\index{equivariant differential form, on Lie group} 
if we have {\it equivariant} if we have 
for each $g \in G$ the relation 
$$ \lambda_g^*\alpha = \rho_E(g) \circ \alpha\res_{D_g}. $$
We write $\Omega^p(G,E)^G \subeq \Omega^p(G,E)$ for the subspace of equivariant $p$-forms. 
If $E$ is a trivial module, then an equivariant form is 
a left invariant $E$-valued form on $G$. 
An equivariant $p$-form $\alpha$ is uniquely determined by the 
corresponding element 
$\alpha_\be \in C^p_c(\g,E) = \Alt^p(\g,E)$: 
\begin{eqnarray}
  \label{eq:4.l.6}
\alpha_g(T_\be(\lambda_g)x_1, \ldots, T_\be(\lambda_g)x_p) 
= \rho_E(g) \circ \alpha_\be(x_1, \ldots, x_p) 
\end{eqnarray}
for $g \in G, x_i \in \g$.  

Conversely, (\ref{eq:4.l.6}) defines for each $\omega \in
C^p_c(\g,E)$ a unique 
equivariant $p$-form $\omega^{\rm eq}$ on $G$ with $\omega^{\rm eq}_\be = \omega$.
\end{defn}

\begin{prop} \mlabel{prop:CE}
For each $\omega \in C^p_c(\g,E)$,  $p \in \N_0$, 
we have the relation $\dd(\omega^{\rm eq}) = (\dd_\g\omega)^{\rm eq}$. 
In particular, the evaluation map 
$$ \ev_\be \: \Omega^p(G,E)^G \to C^p_c(\g,E), \quad 
\omega \mapsto \omega_\be $$
defines an isomorphism from the chain complex $(\Omega^\bullet(G,E)^G,\dd)$ 
of equivariant
$E$-valued differential forms on $G$ to the continuous $E$-valued Lie algebra
complex $(C^\bullet_c(\g,E),\dd_\g)$
discussed in {\rm Section~\ref{app:liealg-cohom}}. 
\end{prop}

\begin{prf} For $g \in G$, we have on the open subset $D_g$ of $G$: 
$$ \lambda_g^* \dd\omega^{\rm eq}
= \dd\lambda_g^* \omega^{\rm eq}
= \dd(\rho_E(g) \circ \omega^{\rm eq})
= \rho_E(g) \circ (\dd\omega^{\rm eq}), $$
showing that $\dd\omega^{\rm eq}$ is equivariant. 

For $x \in \g$, we write $x_l$ for the corresponding
left invariant vector field on~$G$. 
In view of (\ref{eq:4.l.6}), it 
suffices to calculate the value of $\dd\omega^{\rm eq}$ on $(p+1)$-tuples of 
left invariant vector fields in the identity element. From 
$$ \omega^{\rm eq}(x_{1,l},\ldots, x_{p,l})(g) 
= g.\omega(x_1, \ldots, x_p), $$
we obtain 
$$ \big(x_{0,l}.\omega^{\rm eq}(x_{1,l},\ldots, x_{p,l})\big)(\be) 
= x_0.\omega(x_1, \ldots, x_p), $$
and therefore 
\begin{eqnarray*}
&&\ \ \ \ \big(\dd\omega^{\rm eq}(x_{0,l}, \ldots, x_{p,l}\big)\big)(\be) \\
&&=\sum_{i=0}^{p} (-1)^i \big(x_{i,l}.\omega^{\rm eq}(x_{0,l}, \ldots, \hat{x_{i,l}}, 
\ldots, x_{p,l})\big)(\be) \cr
&&\ \ \ \ + \sum_{i < j}(-1)^{i+j} \omega^{\rm eq}([x_{i,l}, x_{j,l}], x_{0,l}, 
\ldots,\hat{x_{i,l}}, \ldots, \hat{x_{j,l}}, \ldots, x_{p,l})(\be)\cr           
&&=\sum_{i=0}^{p} (-1)^i x_i.\omega(x_0, \ldots, \hat{x_i}, 
\ldots, x_p) \cr
&&\ \ \ \ + \sum_{i < j}(-1)^{i+j} \omega([x_i, x_j], x_0, 
\ldots,\hat{x_i}, \ldots, \hat{x_j}, \ldots, x_p) \cr
&&=(\dd_\g\omega)(x_0, \ldots, x_p). 
\end{eqnarray*}
This proves our assertion. 
\end{prf}

\section{Exponential local Lie groups} \mlabel{sec:5.2} 

Before we turn to the deeper global Lie theory of locally exponential 
Lie groups, we need to take a closer look at their local 
versions, the exponential local Lie groups. We shall see later that 
suitable identity neighborhoods in locally exponential Lie groups 
provide exponential local Lie groups, but the converse is not true 
in general. Not every exponential local Lie group is 
enlargeable in the sense that it is obtained 
from a global locally exponential one. 

The most important results in this section are the following:  
The first one is that any continuous homomorphisms of locally exponential 
Lie algebras induces a homomorphisms of the corresponding exponential 
local Lie groups. In the classical situation of Banach--Lie algebras,  
this is a trivial consequence of the fact that the local group 
structure is given by the Hausdorff series. Since, in general, 
the local multiplication 
is not analytic, we have to use different methods to achieve this result. 
The crucial tool is the explicit formula for the Maurer--Cartan form of an 
exponential local Lie group and the Uniqueness Lemma~\ref{lem:c.12b}. 
An important consequence of the fact that Lie algebra homomorphisms 
lead to group homomorphisms is that the germ of the corresponding 
exponential local Lie groups is unique, i.e., any two local multiplications 
corresponding to the same Lie bracket coincide on some $0$-neighborhood. 
This is quite remarkable because it means that the multiplication 
is completely determined by its second derivative in the origin. 
From this uniqueness result, we derive that the Taylor 
series of the local multiplication is the Hausdorff series. 

In the following we shall always write the multiplication in a local Lie group 
$(G, D, m_G,\be)$ as $x * y := m_G(x,y)$ for $(x,y) \in D$. 

\subsection{Basic facts and some examples} 

\begin{defn} \mlabel{def:5.2.1} 
\index{local Lie group!exponential} 
An {\it exponential local Lie group} is a 
local Lie group of the form 
$(G,D,m,0)$, where $G$ is an open balanced $0$-neighborhood 
in a locally convex space $\g$, and  
for $x \in U$ and $t,s,t+s \in [-1,1]$ we have $(tx,sx) \in D$ with 
\begin{eqnarray}
  \label{eq:5.2.1} 
tx * sx = (t+s)x. 
\end{eqnarray}
The Lie algebra $\L(G) = \L(G,D,m,0)$ of an exponential local Lie group 
is called a 
{\it locally exponential Lie algebra}. 
\index{Lie algebra!locally exponential} 
\end{defn} 

\begin{lem} \mlabel{lem:5.2.2} Let $G$ be a locally exponential Lie group and 
$U \subeq \L(G)$ be a balanced $0$-neighborhood for which 
$\exp_G\res_U$ is a diffeomorphism onto an open subset of $G$. Then 
$(U,m_U,D,0)$ is an exponential  local Lie group with respect to 
$$ D := \{(x,y) \in U \times U \: \exp_G x \exp_G y \in \exp_G(U)\} $$
and 
$$ x*y := m_U(x,y) := (\exp_G\res_U)^{-1}(\exp_G(x)\exp_G(y)). $$
\end{lem}

\begin{prf} Clearly $(U,m_U,D,0)$ is a local group isomorphic to 
the local subgroup $\exp_G(U)$ of $G$ (Remark~\ref{rem:locsubgroup}). 

For each $x \in U$ and $t,s,t+s \in [-1,1]$, we have 
$tx,sx,(t+s)x \in U$ since $U$ is balanced. Now 
$$ \exp_G(tx)\exp_G(sx) = \exp_G((t+s)x) \in \exp_G(U) $$
implies $(tx,sx) \in D$ and $(tx) * (sx) = (t+s)x$. 
\end{prf}

\begin{exs} \mlabel{exs:5.2.3} 
(a) Let $\g$ be a Banach--Lie algebra. In Theorem~\ref{thm:loc-bangrp} below, we shall see 
that, for any norm $\|\cdot\|$ on $\g$ defining the topology and satisfying 
$\|[x,y]\| \leq \|x\|\cdot\|y\|$ for $x,y \in \g$, the set 
$$ U := \Big\{ x \in \g \: \|x\| < \frac{\log 2}{8}  \Big\}$$
carries the structure of an exponential local Lie group, 
where the multiplication is given by the Hausdorff series 
$$ x * y = x + y + \frac{1}{2}[x,y]+\frac{1}{12}[x,[x,y]]+\frac{1}{12}[y,[y,x]]+\ldots, $$
converging on $U \times U$. 
We conclude that Banach--Lie algebras are locally exponential. 

(b) Finite-dimensional Lie algebras are in particular Banach and therefore 
locally exponential. 

(c) For each real $\eps > 0$, 
the interval $G := ]-\eps,\eps[ \subeq \R$ is an  
exponential local Lie group with respect to addition and 
$D := \{ (x,y) \in G^2 \: x+y \in G\}. $

More generally, every balanced subset $G$ of a locally convex space 
$E$ is an exponential local Lie group with respect to addition. 
\end{exs}

\begin{rem} \mlabel{rem:5.2.4} Let $(G,D,m,0)$ be an exponential local Lie group. 

(a) For $x \in G$ we have 
$x * (-x) = (-x) * x = (1-1)x = 0$
for all $x \in \g$, so that inversion is given by 
\begin{eqnarray}
  \label{eq:5.2.2}
x^{-1} = -x.
\end{eqnarray}

(b) In view of Remark~\ref{rem:locgrp}, the second order term 
in the Taylor series of the inversion is $b_2(x,x)$, but since 
$x^{-1} = -x$ is linear, we obtain $b_2(x,x) = 0$. Hence 
$b_2(x,x) = 0$, which means that $b_2$ is skew-symmetric, and therefore 
\begin{eqnarray}
  \label{eq:5.2.3}
b_2(x,y) = \frac{1}{2} [x,y]. 
\end{eqnarray}

(c) For $x \in G$, let $S := \sup \{ t \in \R \: tx \in G\} \in [0,\infty]$. 
We claim that  we have 
\begin{eqnarray}
  \label{eq:5.2.4}
tx * s x = (t+s)x \quad \mbox{ for } \quad |t|,|s|,|t+s| < S.   
\end{eqnarray}
In fact, if $|t|,|s|,|t+s| < S$, then there exists a $T \in \R$ with 
$$ |t|,|s|,|t+s| < T < S. $$
Since $G$ is balanced, $y := Tx \in G$, and we have 
$$ tx * sx 
= \frac{t}{T} Tx * \frac{s}{T} Tx 
= \frac{t}{T} y * \frac{s}{T} y 
= \Big(\frac{t}{T} + \frac{s}{T}\Big) y  
= \frac{t+s}{T} y  = (t+s) x. $$
Clearly, the conditions $|t|,|s|,|t+s| < S$ are also necessary for the validity 
of (\ref{eq:5.2.4}) because otherwise $tx,sx$ or $(t+s)x$ is not contained in $G$. 

(d) For $n \in \N$ and $x \in \frac{1}{n} G$ we claim that 
$$ x^n := x * (x * (x * \ldots )) $$ 
simply equals $nx$. 
We argue by induction that $x^k = kx$ for $k \leq n$. 
For $k =1$, there is nothing to show. 
Suppose that $k < n$ and $x^k = kx$. Then 
$kx \in \frac{k}{n} G$ and $k + 1 \leq n$, so that (c) implies that 
$$ x * x^k = x * (kx) = (k+1)x. $$
\end{rem} 

\begin{ex} \mlabel{ex:6.2.5} Let $M$ be a compact smooth manifold (with boundary) 
and $(G,D,m,0)$ an exponential local 
Lie group with Lie algebra $\g$. Then, for each $r \in \N \cup \{\infty\}$, 
$C^r(M,\g)$ is a locally convex 
topological Lie algebra,  
when endowed with the compact open $C^r$-topology. 
The subset $C^r(M,G)$ is open in $C^r(M,\g)$, and it carries a local group 
structure given by pointwise multiplication 
$$ (f * g)(x) := f(x) * g(x). $$
This multiplication is smooth (cf.\ Theorem~\ref{thm:mapgro-Lie}) 
and defines on $C^r(M,G)$ the structure of an exponential local Lie group 
with Lie algebra $C^r(M,\g)$ (see also Exercise~\ref{exer:5.6.2}). 
\end{ex}

\begin{prop} \mlabel{prop:5.2.6} Let $(G,D,m_G,0)$ be an exponential 
local Lie group 
and $x, y \in \g = \L(G)$. Then we have the 
\index{Trotter Product Formula} 
{\em Trotter Product Formula} 
$$ x + y = \lim_{n \to \infty} n \Big(\frac{x}{n} * \frac{y}{n}\Big)$$
and the 
{\em Commutator Formula} 
\index{Commutator Formula} 
$$ [x,y] = \lim_{n \to \infty} n^2 \Big(\frac{x}{n} * \frac{y}{n} * 
\frac{-x}{n} * \frac{-y}{n}\Big). $$
\end{prop}

\begin{prf} The Trotter Product Formula  follows from 
\[  x + y 
= T_{(0,0)}(m_G)(x,y)
= \lim_{n \to \infty} n \cdot m_G\Big(\frac{x}{n}, \frac{y}{n}\Big)
= \lim_{n \to \infty} n \cdot \Big(\frac{x}{n} * \frac{y}{n}\Big). \] 
For the Commutator Formula, we recall from  
(Remark~\ref{rem:locgrp}) that the Taylor expansion in 
$(0,0)$ of the commutator satisfies 
$$ x * y * (-x)* (-y) 
= [x,y] + \cdots.  $$
Now the second order Taylor Formula yields 
\[ \lim_{n \to \infty} n^2\Big(\frac{x}{n} 
* \frac{y}{n} * \frac{-x}{n} * \frac{-y}{n}\Big) = [x,y]. \qedhere\] 
\end{prf}

\subsection
[Local group homomorphisms] 
{Local group homomorphisms versus Lie algebra homomorphisms} 

In this subsection we first prove the Linearity Theorem, 
that any continuous homomorphism of 
exponential local Lie groups is the restriction of a Lie algebra homomorphism, 
hence in particular smooth. We then address the converse and show that 
Lie algebra homomorphisms also define homomorphisms of local Lie groups. 

The following lemma deals 
with the one-parameter case of the Linearity Theorem. 
It is a local version of Lemma~\ref{lem:5.4.1}. 

\begin{lem} \mlabel{lem:5.2.7} Let $(G,D,m,0)$ be an exponential local Lie group, 
$\eps > 0$, and 
$$ \gamma \: ]-\eps,\eps[ \to G \subeq \g = \L(G) $$
a continuous map with 
$$ \gamma(t + s) = \gamma(t) * \gamma(s) \quad \mbox{ for } \quad |t|,|s|,|t+s| < \eps. $$
Then there exists a unique element $x\in \g$ with 
$$ \gamma(t) = t x \quad \mbox{ for } \quad |t| < \eps. $$
\end{lem}

\begin{prf} First we observe that $\gamma(0) = \gamma(0) * \gamma(0)$ 
implies $\gamma(0) = 0$ (Remark~\ref{rem:locgrp0}(c)). 
For $|t| < \frac{\eps}{n}$ and 
$\gamma(t) \in \frac{1}{n} G$, we have 
\begin{eqnarray}
  \label{eq:5.2.5}
n\gamma(t) = \gamma(t)^n = \gamma(nt)
\end{eqnarray}
(Remark~\ref{rem:5.2.4}(d)), 
where the last equality follows by induction as in Remark~\ref{rem:5.2.4}(d). 

Since $\gamma$ is continuous in $0$, there exists some 
$T  \in ]0,\eps[$ such that $\gamma(t) \in \frac{1}{2}G$ for $|t| \leq T$, 
so that $\gamma(2t) = \gamma(t)*\gamma(t) = 2 \gamma(t)$. 
Let $s \in \R$ and $N \in \N$ with $2^{-N} |s| \leq T$. 
Then 
$$ 2^k \gamma(2^{-N-k}s) 
= 2^{k-1} \gamma(2^{-N-(k-1)}s) = \ldots 
= 2 \gamma(2^{-N-1}s) = \gamma(2^{-N}s) $$
for $k \in \N$ (Remark~\ref{rem:5.2.4}(d)) shows that the limit 
$$ \delta(s) := \lim_{n \to \infty} 2^n \gamma(2^{-n}s) $$ 
exists for each $s \in \R$ because the sequence is eventually constant. 
For ${|s| \leq 2^N T}$, we have 
$$ \delta(s) = 2^N \gamma(2^{-N}s). $$ 
Hence $\delta \: \R \to \g$ is a continuous function 
with $\delta(s) = \gamma(s)$ for $|s| \leq T$. 

Now $\delta(2s) = 2\delta(s)$ holds for each $s \in \R$ 
by definition. 
For $k \in \N$ and $s \in \R$, we further derive with (\ref{eq:5.2.5}) 
$$ \delta(ks) 
= \lim_{n \to \infty} 2^n \gamma(2^{-n}ks) 
= \lim_{n \to \infty} 2^n k\gamma(2^{-n}s) = k \delta(s). $$
This implies that 
$$ \delta(ts) = t \delta(s) \quad \hbox{ for } \quad 
t \in 2^\Z \N. $$
Since $2^\Z \N$ is dense in $]0,\infty[$, and $\delta$ is continuous,
we obtain $\delta(t) = t \delta(1)$ for $t \geq 0$. 
Further $\gamma(-t) * \gamma(t) = \gamma(0) = 0 = \gamma(t)*\gamma(-t)$ 
leads to $\gamma(-t) = \gamma(t)^{-1} = - \gamma(t)$. 
We thus obtain 
$\delta(-t) = - \delta(t)$ for all $t \in \R$. This proves that 
$\gamma(t) = t \delta(1)$ for $|t| \leq T$. 

To see that this formula is valid for each $t \in ]-\eps,\eps[$, 
consider the connected 
component $I$ of $0$ in the set of all $t\in [0,\eps[$ with $\gamma(t) = tx$. 
Its closedness follows from the continuity of $s \mapsto \gamma(s)-  sx$. 
We claim that it is also an open subset of $[0,\eps[$. 
Pick $t \in I$. If $t = 0$, then $t$ is an interior point of $I$ in $[0,\eps[$. 
If $t > 0$, then $\gamma(t) = tx$, and we choose 
$0 < h < \min(\eps - t, t)$ in such a way that $(t + h) x \in G$. We now 
obtain for $s < h$ with Remark~\ref{rem:5.2.4}(c): 
$$ \gamma(t+s) = \gamma(t) * \gamma(s) = tx * sx = (t+s)x. $$
Hence $I$ is open in $[0,\eps[$, and the connectedness of $[0,\eps[$ 
yields $I = [0,\eps[$. We likewise obtain 
$\gamma(t) = tx$ for all $t \in ]-\eps,0]$. 

The uniqueness of $x$ follows from $\gamma'(0) = x$. 
\end{prf}

The following theorem shows that continuous homomorphisms of 
exponential local groups are restrictions of linear maps, hence in 
particular smooth. 
\begin{thm} [Linearity Theorem] \mlabel{thm:5.2.8} Let $(G,D_G,m_G,0)$ and 
\break $(H,D_H,m_H,0)$ be 
exponential local Lie groups and 
$\phi \: G \to H$ be a continuous homomorphism of local groups. 
Then there exists a continuous homomorphism of Lie algebras $\psi \:\g \to\h$ 
with $\phi = \psi\res_G$. In particular, $\phi$ is smooth. 
\end{thm} 

\begin{prf} Let $x \in \g$ and choose $\eps > 0$ such that 
$t x \in G$ holds for $|t| < \eps$. 
In view of
\[  \phi(tx) * \phi(sx) = \phi(tx * sx) = \phi((t + s)x) \quad \mbox{ for }\quad 
|t|, |s|, |t+s| < \eps,\] 
Lemma~\ref{lem:5.2.7} implies the existence of a unique element $\psi(x) \in \h$ with 
$\phi(tx) = t \psi(x)$ for $|t| < \eps$. 
For $x \in G$, we may choose $\eps > 1$, which leads to 
$\psi(x) = \phi(x)$. 

We further have 
$\psi(x) = \lim_{t \to 0} \frac{1}{t} \phi(tx),$
which immediately leads to $\psi(0) = 0$ and, for $\lambda \in \R^\times$, to 
\begin{eqnarray}
  \label{eq:5.2.6}
\psi(\lambda x) 
&=& \lim_{t \to 0} \frac{1}{t}\phi(t \lambda x) 
= \lambda \lim_{t \to 0} \frac{1}{t\lambda} \phi(t \lambda x) 
= \lambda \psi(x). 
\end{eqnarray}
Since $\psi\res_{G} = \phi$ is continuous, (\ref{eq:5.2.6}) implies that 
$\psi$ is continuous on $\g$. 

From the Trotter Product Formula in $G$, resp., $H$ 
(Proposition~\ref{prop:5.2.6}), we now get 
\begin{eqnarray*}
\psi(x + y) 
&=& \lim_{n \to \infty} \psi\Big(n\big(\frac{x}{n} * \frac{y}{n}\big)\Big)
= \lim_{n \to \infty} n \psi\Big(\frac{x}{n} * \frac{y}{n}\Big)
= \lim_{n \to \infty} n \phi\Big(\frac{x}{n} * \frac{y}{n}\Big)\\
&=& \lim_{n \to \infty} n \Big(\phi\big(\frac{x}{n}\big) * \phi\big(\frac{y}{n}\big)\Big)
= \lim_{n \to \infty} n \Big(\psi\big(\frac{x}{n}\big) * \psi\big(\frac{y}{n}\big)\Big)\\
&=& \lim_{n \to \infty} n \Big(\frac{\psi(x)}{n} * \frac{\psi(y)}{n}\Big)
= \psi(x) + \psi(y). 
\end{eqnarray*}
This shows that $\psi$ is linear. Further, the Commutator Formula leads to 
\begin{align*}
\psi([x,y]) 
&= \lim_{n \to \infty} \psi\Big(n^2\big(\frac{x}{n} * \frac{y}{n} 
* \frac{-x}{n} * \frac{-y}{n}\big)\Big)\\
&= \lim_{n \to \infty} n^2 \psi\Big(\frac{x}{n} * \frac{y}{n} 
* \frac{-x}{n} * \frac{-y}{n}\Big) \\
&= \lim_{n \to \infty} n^2 \phi\Big(\frac{x}{n} * \frac{y}{n} 
* \frac{-x}{n} * \frac{-y}{n}\Big) \\
&= \lim_{n \to \infty} n^2 \Big(\phi\big(\frac{x}{n}\big) * \phi\big(\frac{y}{n}\big) 
* \phi\big(\frac{-x}{n}\big) * \phi\big(\frac{-y}{n}\big)\Big) \\
&= \lim_{n \to \infty} n^2 \Big(\psi\big(\frac{x}{n}\big) * \psi\big(\frac{y}{n}\big) 
* \psi\big(\frac{-x}{n}\big) * \psi\big(\frac{-y}{n}\big)\Big) \\
&= \lim_{n \to \infty} n^2 \Big(\frac{\psi(x)}{n} * \frac{\psi(y)}{n}
* \frac{-\psi(x)}{n} * \frac{-\psi(y)}{n}\Big) 
= [\psi(x),\psi(y)]. \qedhere 
\end{align*}
\end{prf}

\begin{lem} \mlabel{lem:5.2.9} If $(G,D,m,0)$ is an exponential local Lie group, 
then its Lie algebra $\g$ is $\ad$-integrable, i.e., for each $x \in \g$ 
the linear operator $\ad x$ generates a smooth $\R$-action 
$\alpha \: \R \to \GL(\g), t \mapsto e^{t \ad x}$ with $\alpha'(0)= \ad x$. 
For each $x \in G$, 
$$ \Ad(x) = e^{\ad x}. $$
\end{lem} 

\begin{prf} Let 
$\Ad \: G \to \Aut(\g), g\mapsto \L(c_g)$
be the adjoint representation of the local group $G$ on its Lie algebra 
$\g$ (cf.\ Definition~\ref{defn:adjrep-locgrp}). 
In view of Theorem~\ref{thm:5.2.8}, we have 
$$ \Ad(x)(y) = (x*y)*x^{-1}= c_x(y) $$
if $y$ is sufficiently close to $0$. 
We now 
obtain a well-defined group homomorphism 
$$ \alpha_x \: \R \to \Aut(\g), \quad 
t \mapsto \lim_{n \to \infty} \Ad\Big(\frac{tx}{n}\Big)^n $$
with $\alpha_x(t) = \Ad(tx)$ for $|t| \leq 1$ 
(Exercise~\ref{exer:e.10}). 

To see that $\alpha_x$ defines a smooth action, it suffices to see that it is smooth 
in a $0$-neighborhood because each $\alpha_x(t)$ is a continuous linear automorphism 
of $\g$, hence in particular a smooth map. 
In view of the linearity of the maps $\alpha_x(t)$, it further suffices to show that the 
map 
$$ \R \times \g \to \g \times \g, \quad (t,y) \mapsto 
(\alpha_x(t)y, \alpha_x(t)^{-1}y) 
= (\alpha_x(t)y, \alpha_x(-t)y)  $$
is smooth in a $0$-neighborhood, which follows from 
$\alpha_x(t)y = (tx * y) * (- tx)$
for $t$ and $y$ sufficiently small. 

Below we shall use that the linearity of $c_x(y)$ in $y$ implies that the 
$n$th order terms of the Taylor expansion of $c_x(y)$ in $(0,0)$  
are of order $n-1$ in $x$ and that 
\[ c_x(y) = x * y * (-x) 
=  y + [x,y] + \cdots + \frac{1}{n!} (\ad x)^n y + \cdots \] 
(cf.~Remark~\ref{rem:locgrp}). 
As $\alpha_x(\R) \subeq \Aut(\g)$, the derivative 
$\alpha_x'(0)$ is a continuous derivation of $\g$ 
(Proposition~\ref{prop:aut-der}). For $y$ sufficiently small, we find that 
\begin{eqnarray*}
\alpha_x'(0)y 
&&= \lim_{n \to \infty} n\big(\alpha_x(\frac{1}{n}) - \alpha_x(0)\big)y  
= \lim_{n \to \infty} n\big( \frac{x}{n} * y * \frac{-x}{n} - y\big) \\
&&= \lim_{n \to \infty} n\big(y + \frac{1}{n} [x,y]  + \cdots - y\big) = [x,y].
\end{eqnarray*}
This means that $\alpha_x'(0) = \ad x$ (pointwise on $\g$), 
so that $\g$ is $\ad$-integrable. 
\end{prf}

The following proposition is the key to the 
Lie algebra structure of a locally exponential 
Lie algebra determines the germ of the local multiplication in $(0,0)$. 
It shows that the differential of the left multiplications is 
determined by the Lie bracket. Here the amazing point is that the multiplication 
is only assumed to be smooth, and in general the germ of a smooth function is not determined 
by its Taylor series in any point. 
This result is a Lie algebra version of Theorem~\ref{thm:exp-logder}  
on the derivative of the exponential function. 
We recall from Definition~\ref{defn:mcart-locgrp} the Maurer--Cartan form $\kappa_G$ of a 
local Lie group $(G,D,m,\be)$ defined by $(\kappa_G)_g = T_\be(\lambda_g)^{-1}$. 

\begin{prop} \mlabel{prop:5.2.10} If $(G,D,m_G,0)$ is an exponential local Lie group 
and $x \in G$, then the Maurer--Cartan form $\kappa_G$ is given by the operator-valued 
integral  
$$ (\kappa_G)_x = \int_0^1 e^{-t \ad x}\, dt, $$
interpreted in the pointwise sense. 
\end{prop}

\begin{prf} Fix $t,s \in \R$ with $|t|, |s| < 1$. 
Then 
$$f_{t+s} \: G \to \g, \quad x \mapsto (t + s) x = tx * sx$$
is a smooth function which can be written as 
$f_{t+s} = f_t \cdot f_s$ for $f_t(x) := tx$ and $f_s(x) = sx$. 
For $x \in G$ we consider the smooth curve 
$$ \psi \: [0,1] \to \g, \quad t \mapsto (\kappa_G)_{tx}(ty) = t (\kappa_G)_{tx}(y) 
= \delta(f_t)_x(y), $$
satisfying $\psi(0) = 0$ and 
$\psi'(0) = \lim_{t \to 0} t^{-1} \psi(t) = (\kappa_G)_0(y) = y.$
With the product formula for left logarithmic derivatives 
(Lemma~\ref{lem:4.l.15}), we now get with 
Remark~\ref{lem:5.2.9}  
\begin{eqnarray*}
\psi(t+s) 
&=& \delta(f_{t+s})_x(y)
= \delta(f_t * f_s)_x(y)\\  
&=& \delta(f_s)_x(y) + \Ad(f_s(x))^{-1} \delta(f_t)_x(y) 
= \psi(s) + \Ad(sx)^{-1} \psi(t) \\
&=&  \psi(s) + e^{-s \ad x} \psi(t).
\end{eqnarray*}
Taking derivatives in $t = 0$, this leads to the differential equation 
$$ \psi'(s) = e^{-s\ad x} \psi'(0) = e^{-s\ad x}y, $$
which has the unique solution 
$$ \psi(s) 
= \psi(0) + \int_0^s e^{-t \ad x}y\, dt 
= \int_0^s e^{-t \ad x}y\, dt. $$
Now the assertion follows from $\psi(1) = (\kappa_G)_xy$. 
\end{prf} 

\begin{rem} \label{rem:5.2.11} (a) Let $(G,D,m_G,0)$ be an exponential local Lie 
group. We consider the multiplication map $m_G \: D \to G$ as a smooth function 
on the manifold $D$ with values in the local Lie group $G$. 
Let 
$$p_j \: G \times G \to G, \quad (x_1,x_2) \mapsto x_j, \quad j=1,2,$$ 
denote the projections. 
We consider $p_1, p_2$ as smooth functions $D \to G$ and note that 
$(p_1,p_2)(D) = D$. Therefore the Product Rule for logarithmic derivatives 
(Lemma~\ref{lem:4.l.15}) implies that 
$$ \delta(m_G) = \delta(p_2) + \Ad(p_2)^{-1}\delta(p_1). $$
With Proposition~\ref{prop:5.2.10} we get 
$$ \delta(p_1)_{(x,y)}
= (\kappa_G)_x \circ p_1
= \int_0^1 e^{-t\ad x}\, dt \circ p_1 $$
and 
$$ \delta(p_2)_{(x,y)}= \int_0^1 e^{-t\ad y}\, dt \circ p_2, $$
which leads to 
$$ \delta(m_G)_{(x,y)} = 
\Big(e^{-\ad y} \int_0^1 e^{-t\ad x}\, dt, \int_0^1 e^{-t\ad y}\, dt\Big) 
\in {\cal L}(\g,\g)^2 \cong {\cal L}(\g \times \g,\g). $$

(b) Let $G$ be a Lie group with an exponential function 
$\exp_G \: \g \to G$. We want to compute the logarithmic derivative 
of the function 
$$ m_G \circ (\exp_G \times \exp_G) \: \g \times \g \to G. $$
With the notation from (a), this function is the product of 
$\exp_G \circ p_1$ with $\exp_G \circ p_2$, so that 
the Product Rule and 
$\delta(\exp_G)_x = \int_0^1 e^{-t\ad x}\, dt$ (Theorem~\ref{thm:exp-logder}) lead to 
\begin{eqnarray*}
&& \delta\big(m_G \circ (\exp_G \times \exp_G)\big)
=  \delta((\exp_G \circ p_1) \cdot (\exp_G \circ p_2)) \\
&=&  \delta(\exp_G \circ p_2) + \Ad(\exp_G \circ p_2)^{-1}\delta(\exp_G \circ p_1).  
\end{eqnarray*}
Evaluating in $(x,y) \in \g \times \g$, thus leads to 
\begin{eqnarray*}
&& \delta\big(m_G \circ (\exp_G \times \exp_G)\big)_{(x,y)} \\
&=& \Big(e^{-\ad y} \int_0^1 e^{-t\ad x}\, dt, \int_0^1 e^{-t\ad y}\, dt\Big) 
\in {\cal L}(\g,\g)^2 \cong {\cal L}(\g \times \g,\g).
\end{eqnarray*}
\end{rem}

\begin{thm} \mlabel{thm:5.2.11c} Let 
$(G,D_G,m_G,0)$ and $(H,D_H,m_H,0)$ be exponential local Lie groups 
and $\psi \: \g = \L(G) \to \h= \L(H)$ be 
a continuous homomorphism of Lie algebras. 
Let $D_0$ denote the connected component of $D_G \cap (\psi \times \psi)^{-1}(D_H)$ 
containing $(0,0)$. Then 
$$ \psi(x*y) = \psi(x) * \psi(y) \quad \mbox{ for } \quad (x,y) \in D_0. $$

For $G = \g$ and $H = \h$ we have in particular $D_0 = \g \times\g$, so that 
\break $\psi \: (\g,*) \to (\h,*)$ is a homomorphism of exponential Lie groups. 
\end{thm}

\begin{prf} Let $D := D_G \cap (\psi \times \psi)^{-1}(D_H)$. 
On $D$ we have two smooth functions 
$$ f_1 \: D \to H, (x,y) \mapsto \psi(x*y)
\quad \mbox{ and } \quad 
f_2 \: D \to H, (x,y) \mapsto \psi(x)*\psi(y) $$
satisying $f_1(0,0) = 0 = f_2(0,0)$. 

We claim that $\delta(f_1) = \delta(f_2)$. We first observe 
that on $G \cap \psi^{-1}(H)$ we have 
$$ \delta(\psi)_x 
= (\kappa_H)_{\psi(x)} \circ \psi 
= \int_0^1 e^{-t\ad \psi(x)}\, dt \circ \psi 
= \psi \circ \int_0^1 e^{-t\ad x}\, dt  = \psi \circ (\kappa_G)_x $$
(Exercise~\ref{ex:adint-compat} and Proposition~\ref{prop:5.2.10}),
i.e., $\delta(\psi) = \psi \circ \kappa_G$. 
This leads to 
$$ \delta(f_1) 
= \delta(\psi \circ m_G) 
= m_G^*\delta(\psi) = m_G^*(\psi \circ \kappa_G) = \psi \circ (m_G^*\kappa_G) 
=\psi \circ \delta(m_G). $$
On the other hand, we have 
$$ \delta(f_2) 
= \delta(m_H \circ (\psi \times \psi)) 
= (\psi \times \psi)^* \delta(m_H). $$
Using Remark~\ref{rem:5.2.11}, evaluating this in $(x,y) \in D$ leads to 
\begin{eqnarray*}
\delta(f_2)_{(x,y)} 
&=& \delta(m_H)_{(\psi(x),\psi(y))} \circ (\psi \times \psi) \\
&=& \Big(e^{-\ad \psi(y)} \int_0^1 e^{-t\ad \psi(x)}\, dt, \int_0^1 e^{-t\ad \psi(y)}\, dt\Big) 
\circ (\psi \times \psi) \\
&=& \psi \circ \Big(e^{-\ad y} \int_0^1 e^{-t\ad x}\, dt, \int_0^1 e^{-t\ad y}\, dt\Big) \\ 
&=& \psi \circ \delta(m_G)_{(x,y)} = \delta(f_1)_{(x,y)}. 
\end{eqnarray*}
Now the Uniqueness 
Lemma~\ref{lem:uni-lemma} implies that $f_1$ and $f_2$ coincide 
on the connected component $D_0$ of $D$ containing $(0,0)$. 
For $G = \g$ and $H = \h$ we have $D = \g \times \g$, so that $D_0 = \g \times \g$. 
\end{prf} 

There is also a ``semilocal'' version of the preceding theorem: 
\begin{thm} \mlabel{thm:5.2.11b} Let 
$(G,D_G,m_G,0)$ be an exponential local Lie group with Lie algebra $\g$,  
$H$ a Lie group with an exponential function \break $\exp_H \: \h \to H$ 
and $\psi \: \g  \to \h$ be a continuous homomorphism of Lie algebras. 
Let $D_0$ denote the connected component of $D_G$ 
containing $(0,0)$. Then 
$$ \exp_H(\psi(x*y)) = \exp_H(\psi(x))\exp_H(\psi(y)) \quad \mbox{ for } \quad (x,y) \in D_0. $$
\end{thm}

\begin{prf} With $\phi := \exp_H \circ \psi \: G \to H$ we obtain 
two smooth function 
$$ f_1 \: D_G \to H, (x,y) \mapsto \phi(x*y)
\quad \mbox{ and } \quad 
f_2 \: D_G \to H, (x,y) \mapsto \phi(x)\phi(y) $$
satisfying $f_1(0,0) = \be = f_2(0,0)$. 

To see that $\delta(f_1) = \delta(f_2)$, we use Theorem~\ref{thm:exp-logder} to get 
$$ \delta(\phi)_x 
= \delta(\exp_H)_{\psi(x)} \circ \psi 
= \int_0^1 e^{-t\ad \psi(x)}\, dt \circ \psi 
= \psi \circ \int_0^1 e^{-t\ad x}\, dt  = \psi \circ (\kappa_G)_x $$
(Exercise~\ref{ex:adint-compat} and Proposition~\ref{prop:5.2.10}), 
i.e., $\delta(\phi) = \psi \circ \kappa_G$. 
This leads to 
$$ \delta(f_1) 
= \delta(\phi \circ m_G) 
= m_G^*\delta(\phi) = m_G^*(\psi \circ \kappa_G) = \psi \circ (m_G^*\kappa_G) 
=\psi \circ \delta(m_G). $$
On the other hand, we have 
\begin{eqnarray*}
\delta(f_2) 
&=& \delta(m_H \circ (\phi \times \phi)) 
= \delta(m_H \circ (\exp_H \times \exp_H) \circ (\psi \times \psi)) \cr
&=& (\psi \times \psi)^* \delta(m_H \circ (\exp_H \times \exp_H)). 
\end{eqnarray*}
Using Remark~\ref{rem:5.2.11}(b) and evaluating this in $(x,y) \in D$, leads to 
\begin{eqnarray*}
\delta(f_2)_{(x,y)} 
&=& \delta(m_H \circ (\exp_H \times \exp_H))_{(\psi(x),\psi(y))} \circ (\psi \times \psi) \\
&=& \Big(e^{-\ad \psi(y)} \int_0^1 e^{-t\ad \psi(x)}\, dt, \int_0^1 e^{-t\ad \psi(y)}\, dt\Big) 
\circ (\psi \times \psi) \\
&=& \psi \circ \Big(e^{-\ad y} \int_0^1 e^{-t\ad x}\, dt, \int_0^1 e^{-t\ad y}\, dt\Big) \\ 
&=& \psi \circ \delta(m_G)_{(x,y)} = \delta(f_1)_{(x,y)}. 
\end{eqnarray*}
Now the Uniqueness Lemma~\ref{lem:uni-lemma} implies that $f_1$ and $f_2$ coincide 
on the connected component $D_0$ of $D_G$ containing $(0,0)$. 
\end{prf} 

As an important consequence of Theorem~\ref{thm:5.2.11c}, 
the germ of the multiplication in $(0,0)$ is determined by the Lie bracket. This is 
quite remarkable because the $*$-multiplication is not assumed to be analytic. 

\begin{thm} [Uniqueness of germs of exponential local groups] 
 \mlabel{thm:5.2.12} Let 
$(G,D_G, m_G,0)$ and $(H,D_H,m_H,0)$ be two exponential local Lie 
groups with the same Lie algebra $\g \supeq G, H$ and $D_0$ 
the connected component of 
$D_G \cap D_H$ containing $(0,0)$. Then 
$$ m_G(x,y) = m_H(x,y) \quad \mbox{ for all } \quad x,y \in D_0. $$
\end{thm}

\begin{prf} We apply the preceding theorem with $\psi = \id_\g$, which leads to 
\break 
$D_G \cap (\psi \times \psi)^{-1}(D_H) = D_G \cap D_H$.   
\end{prf}

\begin{cor}  \mlabel{cor:5.2.14} 
Let $\g$ be a nilpotent locally convex Lie algebra 
and $G \subeq \g$ be an exponential local group with Lie algebra~$\g$. 
Then, for each pair $(x,y)$ in the connected component of $D_G$ containing 
$(0,0),$ the product $x*y$ is given by the Hausdorff series. 
\end{cor} 

\begin{prf} This follows by applying Theorem~\ref{thm:5.2.12} 
to $\psi = \id_\g$ and $H = (\g,*)$, where $*$ is the BCH product 
(Theorem~\ref{thm:BCH-nil}). 
In this case, 
\[D_G \cap (\psi \times \psi)^{-1}(D_H) = D_G. \qedhere \] 
\end{prf}

\section{Locally exponential Lie algebras} \mlabel{sec:5.3} 

In this section, we change our perspective from 
the local to the infinitesimal, i.e., from the exponential local 
Lie groups to their Lie algebras, called locally exponential Lie algebras. 
The first main result of this section is a characterization of those 
ideals $\fn$ of a locally exponential Lie algebra $\g$ for which 
$\g/\fn$ is locally exponential. As we shall see in Section~\ref{sec:6.4}, 
this characterization is crucial to understand which 
normal subgroups $N$ of a locally exponential Lie group $G$ lead to 
Lie group quotients $G/N$. The second main result of the present 
section is the Adjoint Enlargeability 
Theorem: For any locally exponential Lie algebra $\g$, the 
quotient $\g/\z(\g)$ is enlargeable, which is shown by endowing the subgroup 
of $\Aut(\g)$ generated by $e^{\ad \g}$ with a locally exponential Lie 
group structure. 

\begin{defn} \mlabel{defn:5.3.1} A Lie algebra $\g$ is said to be 
\index{Lie algebra!locally exponential} 
{\it locally exponential} 
if it is the Lie algebra of an exponential local Lie group, i.e., 
if there exists an open symmetric $0$-neighborhood $G \subeq \g$ carrying the structure of 
an exponential local Lie group $G$ with $\L(G) = \g$. 
\end{defn}

The following lemma provides a characterization of closed locally 
exponential subalgebras in terms of the local multiplication. 

\begin{lem} \mlabel{lem:locexpsub-crit} Let $\h$ be a Lie subalgebra of the locally exponential 
Lie algebra~$\g$. 
If $\h$ is locally exponential, then there exists an exponential local Lie group $G \subeq \g$ with 
\begin{eqnarray}
  \label{eq:locexp-subalg} 
  m_G(D_G \cap (\h \times \h)) \subeq \h. 
\end{eqnarray}
  If, conversely, this condition is satisfied and $\h$ is closed, then $\h$ is locally exponential. 
\end{lem}

\begin{prf} (a) We first assume that $\h$ is locally exponential. Let $H \subeq \h$ and $G_1 \subeq \g$ 
be exponential local Lie groups. Then Theorem~\ref{thm:5.2.11c} implies that for 
$(x,y) \in D_0$, the connected component of $(0,0)$ in 
$D_H \cap D_{G_1}$, which is an open subset of $\h \times \h$, we have 
$x *_H y = x *_{G_1} y.$
Pick an exponential local Lie subgroup $G \subeq G_1$ satisfying 
$$ (G \times G) \cap (\h \times \h) \subeq D_0. $$
For $(x,y) \in D_G \cap (\h \times \h)$,  we then have 
$x *_G y = x*_{G_1} y = x *_H y \subeq \h.$

(b) If $G$ satisfies (\ref{eq:locexp-subalg}), then 
$H := G \cap \h$ is an open balanced $0$-neighborhood in $\fh$ and 
$D_H := D_G \cap (\h \times \h)$ is mapped by 
$m_H := m_G\res_{D_H}$ into~$\h$. 
If $\h$ is closed, then $m_H$ is a smooth map, and now it is easy to see that 
$(H,D_H,m_H,0)$ is an exponential local Lie group in $\h$. 
\end{prf}

\begin{prop}{\rm(Inverse images of locally exponential subalgebras)} 
 \mlabel{prop:invim} Let 
$\psi \: \g \to \h$ be a morphism of 
locally exponential 
Lie algebras and $\h_1 \subeq \h$ be a closed locally exponential subalgebra. 
Then $\psi^{-1}(\h_1)$ is locally exponential. 
In particular, $\ker \psi$ is locally exponential. 
\end{prop}

\begin{prf} First we use Lemma~\ref{lem:locexpsub-crit} 
to find an exponential 
local Lie group $H \subeq \h$ with 
$$m_H(D_H \cap (\h_1 \times \h_1)) \subeq \h_1.$$
Next we choose an exponential local Lie group 
$G \subeq \g$ with $\psi(G) \subeq H$, $\psi(D_G) \subeq D_H$ and 
$\psi(x*y) = \psi(x)*\psi(y)$ for $(x,y) \in D_G$ 
(Theorem~\ref{thm:5.2.11c}). 
We then obtain for the inverse image $\g_1 := \psi^{-1}(\fh_1)$ 
and $(x,y) \in D_G \cap (\g_1 \times \g_1)$ the relation 
$$ \psi(x * y) = \psi(x) * \psi(y) \in \h_1, $$
so that $x * y \in \g_1$. Since $\g_1$ is closed, 
Lemma~\ref{lem:locexpsub-crit} now shows that $\g_1$ is locally exponential. 
\end{prf}

\begin{cor} \mlabel{cor:equalizer} Let 
$\psi_j \: \g_j \to \h$, $j =1,2$, be  morphisms of 
locally exponential 
Lie algebras. Then the equalizer 
\[ \fe := \{ (x_1, x_2) \in \g_1 \times \g_2 \: \psi_1(x_1) = \psi_2(x_2)\} \] 
is locally exponential. 
\end{cor}

\begin{prf} We apply the preceding proposition to the morphism of locally exponential 
Lie algebras  
$\psi_1 \times \psi_2 \: \g_1 \times \g_2 \to \fh \times \fh$ and 
observe that the diagonal $\Delta_\fh \subeq \fh \times \fh$ is a 
closed locally exponential subalgebra. 
\end{prf}

One of the most serious problems Lie theory one has to cope with beyond the realm of  
Banach--Lie groups is that there exist smooth actions of 
$\R$ on locally convex spaces $E$ whose infinitesimal generator 
leaves a closed subspace $F$ invariant, but the whole one-parameter group 
does not (cf.\ Example~\ref{ex:e.2.9}). In the context of Lie algebras, 
this leads to the concept of stability: 

\begin{defn} \mlabel{def:5.3.3} We call a subalgebra $\h \leq \g$ of an 
$\ad$-integrable Lie algebra $\g$ (cf.\ Lemma~\ref{lem:5.2.9}) 
\index{Lie algebra! stable} 
\index{ideal! stable} 
{\it stable} if 
$$ e^{\ad x}\h = \h \quad \mbox{ for all } x \in \h. $$
An ideal $\fn \trile \g$ is called {\it stable} if 
$$ e^{\ad x}\fn = \fn \quad \mbox{ for all } x \in \g. $$
\end{defn}

\begin{lem} \mlabel{lem:5.3.4} If 
$\alpha \: \g \to \h$ is a morphism of $\ad$-integrable Lie algebras, 
then $\im(\alpha)$ is a stable subalgebra of $\fh$ and 
$\ker(\alpha)$ is a stable ideal. 

If, in addition, $G \subeq \g$ and 
$H \subeq \h$ are exponential local Lie groups with $\alpha(G) \subeq H$, then 
$$ (\kappa_G)_x(\ker(\alpha)) = \ker(\alpha) $$
for all $x \in G$. 
\end{lem}

\begin{prf} Since $\g$ and $\h$ are $\ad$-integrable, 
$e^{\ad_\g \g}$ and $e^{\ad_\h \h}$ are defined, and for each $x \in \g$, 
we have 
$$ \alpha \circ e^{\ad_\g x} = e^{\ad_\h \alpha(x)} \circ \alpha $$
(Exercise~\ref{ex:adint-compat}), 
which directly implies that $\im(\alpha)$ is a stable subalgebra and 
that $\ker \alpha$ is a stable ideal. 

Now let $G \subeq \g$ be an exponential local Lie group. 
From $(\kappa_G)_x = \kappa_\g(x)$ (Proposition~\ref{prop:5.2.10}) we derive   
$\alpha \circ (\kappa_G)_x = (\kappa_H)_{\alpha(x)} \circ \alpha.$
Since $(\kappa_G)_x$ and $(\kappa_H)_{\alpha(x)}$ are invertible, 
it follows that 
\[\ker \alpha = \ker(\alpha \circ (\kappa_G)_x) = (\kappa_G)_x^{-1}(\ker \alpha). \qedhere\] 
\end{prf}

\begin{ex} \mlabel{ex:5.3.4b} We recall the exponential Fr\'echet--Lie group 
$G = E \rtimes_\alpha \R$ from Example~\ref{ex:5.1.11}, where 
$E$ is the space of smooth $1$-periodic functions $\R \to \C$, and 
$(\alpha_t f)(x) = e^{\mu t} f(x+t)$ is the smooth $\R$-action on 
$E$ with infinitesimal generator $Df = \mu f + f'$. 

Let $M \subeq [0,1]$ be an open subset which is not dense and put 
$$E_M := \{ f \in E \: f\res_M = 0\}.$$ 
Then $E_M$ is a closed subspace of 
$E$ with $DE_M \subeq E_M$ but 
$\alpha(t)(E_M) \not\subeq E_M$ for some $t \in \R$. 
Therefore $\h_M := E_M \rtimes_D \R \leq \g = E \rtimes_D \R$ is a closed subalgebra of 
the exponential Fr\'echet--Lie algebra $\g$ which is not stable. 
\end{ex}

\begin{ex}
  \mlabel{ex:V.2.15} Let $V := C^\infty(\R,\R)$ and 
consider the one-parameter group \break $\alpha \: \R \to \GL(V)$, given by 
$\alpha_t(f)(x) = f(x+t)$. Then $\R$ acts smoothly on $V$, so that we may 
form the corresponding semidirect product group 
$$ G := V \rtimes_\alpha \R. $$
This is a Lie group with a smooth exponential function given by 
$$ \exp(v,t) 
= \Big( \int_0^1 \alpha_{st}v\, ds, t\Big), $$
where 
$$  \Big(\int_0^1 \alpha_{st}v\, ds\Big)(x) = \int_0^1 v(x + st)\, ds. $$

The Lie algebra $\g$ has the corresponding semidirect product structure 
$\g = V \rtimes_D \R$ with $Dv = v'$, i.e., 
$$ [(f,t), (g,s)] = (tg'-sf', 0). $$

In $\g \cong V \rtimes \R$, we now consider the subalgebra $\h := V_{[0,1]} \rtimes \R$, 
where 
$$ V_{[0,1]} := \{ f \in V \: \supp(f) \subeq [0,1]\}. $$
Then $\h$ clearly is a closed subalgebra of $\g$. It is not stable because 
$\alpha_{-t} V_{[0,1]} = V_{[t,t+1]}.$
The subgroup of $G$ generated by $\exp \h$ contains $\{0\} \rtimes \R$, 
$V_{[0,1]}$, and hence all intervales $V_{[t,t+1]}$, which implies that 
$\la \exp \h \ra = C^\infty_c(\R) \rtimes \R$. 

Lemma~\ref{lem:5.3.4} 
implies that the inclusion $\fh \into \g$ does not 
integrate to a homomorphism $\phi \: H \to G$ of Lie group with an exponential 
function, for which $\Lie(\phi)$ is the inclusion $\fh \into \g$. 
\end{ex}

\index{Quotient Theorem for Locally Exponential Lie Algebras}
\begin{thm}[Quotient Theorem for Locally Exponential Lie Algebras] 
\mlabel{thm:5.3.5} 
Let $\fn \trile \g$ be a closed ideal of the locally exponential Lie algebra $\g$. 
Then the quotient Lie algebra $\g/\fn$ is locally exponential if and only if 
$\fn$ is stable and there exists a $0$-neighborhood $U \subeq \g$ such that 
$\kappa_\g(x)(\fn)~=~\fn$ for all $x \in U$. 
\end{thm}

\begin{prf} Let $q \: \g \to \g/\fn$ denote the quotient morphism. If 
$\g/\fn$ is locally exponential, then Lemma \ref{lem:5.3.4} implies that 
$\fn$ is a stable ideal with $\kappa_\g(x)(\fn) = \fn$ 
whenever $x\in\g$ is  sufficiently close to~$0$.

Suppose, conversely, that 
$G \subeq \g$ is an exponential local Lie group which is a 
convex subset of $\g$, 
and that $\fn$ is stable with $\kappa_G(x)(\fn) = \fn$ for all $x \in G$. 
Let $U \subeq G$ be a convex symmetric $0$-neighborhood such that 
all $4$-fold products of elements of $U$ are defined in $G$ and associative. 

Pick $x,x',x'' \in U$ and assume that $y := x'' - x' \in \fn$. 
For the curve 
$$ \gamma \: [0,1] \to \g, \quad t \mapsto x * (x' + ty) $$
we then obtain with 
$T_b(\lambda_a) T_0(\lambda_b) = T_0(\lambda_{a * b})$ for $(a,b) \in D_G$ and 
the invariance condition on $\fn$: 
\begin{eqnarray*}
\gamma'(t) 
&&= T_{(x' + ty)}(\lambda_x) y 
= T_0(\lambda_{x*(x' + ty)}) T_0(\lambda_{x' + ty})^{-1} y \\
&&= \kappa_G(x*(x' + ty))^{-1} \kappa_G(x' + ty)(y) \in \fn. 
\end{eqnarray*}
We conclude that 
$$ x * x'' = x * (x' + y) = \gamma(1) \in \gamma(0) + \fn = x * x' + \fn. $$
This further leads to 
\begin{eqnarray*}
x'' * x 
&&= (-x) * (x * x'') * x = e^{-\ad x}(x * x'') 
\in e^{-\ad x}(x * x' + \fn) \\
&&= x' * x + e^{-\ad x}\fn 
= x' * x + \fn 
\end{eqnarray*}
because $\fn$ is a stable ideal. 
Therefore the smooth map 
\begin{eqnarray}
  \label{eq:5.2.10}
U  \times U \to \g/\fn, \quad (x,x') \mapsto x * x' + \fn 
\end{eqnarray}
is constant on the intersection of the convex set $U \times U$ with the cosets of the 
subspace $\fn \times \fn$. Let $q \: \g \to \g/\fn$ denote the quotient map 
and $Q := q(U)$, which is an open symmetric convex $0$-neighborhood in $\g/\fn$. 
We put 
$$ D_Q := \{ (q(x),q(y)) \in Q \times Q \: q(x * y) \in Q\}, $$
and obtain a well-defined multiplication map 
$$ m_Q \: D_Q \to Q, \quad (q(x), q(y)) \mapsto q(x) * q(y) := q(x*y). $$
It is smooth because the map 
\eqref{eq:5.2.10} is smooth by Lemma~\ref{lemquot}.
Since triple products in $U$ are defined and associative, it easily follows that 
$(Q,D_Q, m_Q,0)$ is an exponential local Lie group with 
Lie algebra $\g/\fn$. 
\end{prf} 

Combining Lemma~\ref{lem:5.3.4} with Theorem~\ref{thm:5.3.5}, we obtain: 

\begin{cor} [Canonical factorization] If $\alpha \: \g \to \h$ is a morphism of 
locally exponential Lie algebras, then 
$\oline\g := \g/\ker \alpha$ is locally exponential and 
$\alpha$ factors through an injective morphism 
$\oline\alpha \: \oline\g \to \h$. We thus obtain the factorization 
$\alpha = \oline\alpha \circ q,$
where $q \: \g \to \oline\g$ is the quotient map and $\oline\alpha$ is injective.   
\end{cor}

\begin{cor} \mlabel{cor:5.3.6}  If $\g$ is locally exponential and 
$\fn \trile \g$ a closed ideal, then $\g/\fn$ is locally exponential if 
\begin{enumerate}
\item[\rm(i)] $\fn \supeq [\g,\g]$, 
\item[\rm(ii)] $\fn \subeq \z(\g)$, or 
\item[\rm(iii)] $\g/\fn$ is a Banach--Lie algebra. 
\end{enumerate}
\end{cor}

\begin{prf} (i) If $\fn$ be a closed ideal containing all commutators, 
then $\g/\fn$ is abelian and therefore locally exponential with 
$x * y :=  x + y$. 

 (ii) If $\fn \subeq \z(\g)$, then we consider, for $z \in \z(\g)$ 
and $x \in \g$, the curve 
$\gamma(t) := e^{t \ad x}z$ that satisfies 
$$ \gamma(0) = z, \quad \gamma'(t) = e^{t \ad x}[x,z] = 0, $$
which implies $e^{\ad x}z = z$. From that it follows that 
each central ideal $\fn$ is stable and satisfies $\kappa_\g(x)(\fn) = \fn$ for 
each $x \in \g$. Now we apply Theorem~\ref{thm:5.3.5}.  

(iii) This follows from Example~\ref{exs:5.2.3}(a). 
\end{prf}

\begin{cor} \mlabel{cor:5.3.7} 
Let $\fn \trile \g$ be a closed ideal of the locally exponential Lie algebra $\g$. 
Suppose that $\fn$ is stable and that there exists a $0$-neighborhood $U \subeq \g$ such that 
$\kappa_\g(x)(\fn)~=~\fn$ for all $x \in U$. Then $\fn$ is locally exponential.  
\end{cor}

\begin{prf} In view of Theorem~\ref{thm:5.3.5}, $\g/\fn$ is locally exponential, 
so that the quotient map $q \: \g \to \g/\fn$ is a morphism of locally exponential 
Lie algebras and the assertion 
follows from Proposition~\ref{prop:invim}. 
\end{prf}

The following theorem is a key result in the theory of locally exponential 
Lie algebras, resp., Lie groups. For non-Banach--Lie groups, the adjoint 
groups have no reason to be actual Lie groups, but for locally 
exponential groups, resp., Lie algebras, the theorem implies in particular 
that there is no pathology involved 
with the adjoint group. 

\begin{thm} [Adjoint Enlargeability Theorem] \mlabel{thm:5.3.8} 
For each locally exponential Lie algebra 
$\g$, the group $G_{\rm ad} := \la e^{\ad \g} \ra \subeq \Aut(\g)$ 
carries the structure of a locally exponential Lie group 
with Lie algebra $\g_{\ad} := \g/\z(\g)$ and 
the exponential function 
$$ \Exp \: \g_{\ad} 
= \g/\z(\g) \to G_{\rm ad}, \quad x + \z(\g) \mapsto e^{\ad x}.$$ 
In particular, $0$ is isolated in $\Exp^{-1}(\be)$ and
$\z(\g)$ is an open subset of \break $\{ x \in \g \: e^{\ad x} = \id_\g\}$. 
Moreover, the canonical action of $G_{\rm ad}$ on $\g$ is smooth. 
\end{thm}

\begin{prf} From Corollary~\ref{cor:5.3.6} we know 
already that the Lie algebra $\g_{\ad}$ is locally  exponential. 
According to Lemma~\ref{lem:5.2.9}, $\g$ is $\ad$-integrable, 
so that we can consider the group 
$G_{\ad} := \la e^{\ad \g} \ra \subeq \GL(\g)$. The map 
$$ \Exp \: \g_{\ad} = \g/\z(\g) \to G_{\ad}, \quad \oline x := x + \z(\g) \mapsto e^{\ad x} $$
is well-defined since $\z(\g) = \ker \ad$. 

Let $G \subeq \g$ be an exponential local Lie group with Lie algebra $\g$. 
Further, let $U \subeq G$ be a convex symmetric open $0$-neighborhood such that 
all $4$-fold products of elements in $U$ are defined and associative. 
Let $q \: \g \to \g_{\rm ad}$ be the quotient homomorphism. 
We have seen in the proof of Theorem~\ref{thm:5.3.5} that 
$V := q(U) \subeq \g_{\rm ad}$ is an exponential local Lie group with 
Lie algebra $\g_{\rm ad}$ if we define the product of 
$\oline x := q(x)$ and $\oline y := q(y)$ by 
$\oline x * \oline y := \oline{x*y}$. 

For $x,y \in U$ we now get with $\Ad(x) = e^{\ad x}$ (Lemma~\ref{lem:5.2.9}) 
$$ \Exp(\oline x) \Exp(\oline y) 
= e^{\ad x} e^{\ad y} = \Ad(x) \Ad(y) = \Ad(x*y) 
= e^{\ad(x*y)} = \Exp(\oline x * \oline y), $$
so that $\Exp$ restricts to a homomorphism $V \to G_{\rm ad}$ of 
local groups. 

We want to show that $\Exp$ is injective on a suitable subset of $V$. 
We first observe that, for $x \in G$ with 
$\Ad(x) = e^{\ad x} = \id_\g$, the invertibility 
of $\kappa_\g(x) = \int_0^1 e^{-t\ad x}\, dt$, together with 
the relation 
$$ \kappa_\g(x) \circ \ad x  = e^{\ad x}  - \id_\g = 0 $$
(cf. \eqref{eq:kappa-comm} 
in Definition~\ref{def:kappa-liealg}), leads to $\ad x= 0$. 
Since $e^{\ad \z(\g)} = \{\id_\g\}$ holds trivially, we conclude that 
$$ \Ad^{-1}(\id_\g) = G \cap \z(\g). $$

If $\oline x,\oline y \in V$ with 
$(-x) * y \in U$ satisfy $\Exp(\oline x) = \Exp(\oline y)$, then 
$$ e^{\ad((-x)*y)} 
= \Exp((-\oline x) * \oline y) = \Exp(-\oline x) \Exp(\oline y) = \Exp(\oline x)^{-1} \Exp(\oline y)
 = \id_\g $$
and $z := (-x) * y \in U \cap \Ad^{-1}(\id_\g) \subeq \z(\g)$. 
Next we note that $x * z = x * ((-x) * y) = y$. 

For each $a \in G$ and $z \in \z(\g)$, 
we have $e^{\ad a}z = z$, hence 
$\kappa_G(a)z = z$. Applying this to 
the curve 
$\eta(t) = a * tz$ for $t \in [0,1]$ and $a,z \in U$, we obtain
$$ \eta'(t) 
= T_0(\lambda_{a * tz})z = \kappa_G(a * tz)^{-1}z= z, $$ 
which leads to $a * z  = \eta(1) = a + z$. 
In particular, we obtain $y = x * z = x + z$ and hence 
$\oline x = \oline y$. 

Let $V_1 \subeq V$ be a balanced open $0$-neighborhood with 
$$V_1 \times V_1 \subeq D_V := \{ (\oline x, \oline y) \in V^2 \: x * y \in U\}.$$ 
Then the preceding argument implies that 
$\Exp\res_{V_1}$ is injective, so that 
it may be identified with a subset of the group $G_{\rm ad}$. 
Since $G_{\ad}$ is generated by 
$\Exp(V_1)$, the smoothness of multiplication and inversion on $V_1$ 
imply the existence of a unique Lie group structure on $G_{\rm ad}$ for which 
$\Exp$ defines a chart in an identity neighborhood
(Theorem~\ref{thm:locglob}). It follows immediately from 
the construction that $\Exp \: \g_{\rm ad} \to G_{\rm ad}$ is an exponential 
function for this group and hence that $G_{\rm ad}$ is locally exponential. 

We finally shows that the action of $G_{\rm ad}$ on $\g$ is smooth. 
The map 
\[ \g \times \g \to \g, \qquad 
(x,y) \mapsto x * y * (-x) = e^{\ad x}y = \exp_{G_{\rm ad}}(\ad x)y \] 
(Lemma~\ref{lem:5.2.9}) is smooth on a neighborhood of $(0,0)$ in 
$\g \times \g$. Therefore Lemma~\ref{lemquot} shows  that 
the corresponding map $\g_{\rm ad} \times \g \to \g$ 
is smooth in a $0$-neighborhood. Since $G_{\rm ad}$ 
is locally exponential, Exercise~\ref{exer:smooth-crit} now 
implies that the canonical action $\alpha$ of $G_{\rm ad}$ on $\g$ is smooth. 
\end{prf}

\begin{cor} \mlabel{cor:5.3.9} Each 
locally exponential Lie algebra $\g$ with trivial center is enlargeable. 
\end{cor}

From the relation $\Ad(\exp_G x) = e^{\ad x}$ and 
Proposition~\ref{prop:smooth-crit}, we immediately derive: 

\begin{cor} \mlabel{cor:adj-smooth} 
For every locally exponential Lie group, the adjoint 
representation defines a morphism of locally exponential  
Lie groups $\Ad \: G \to G_{\rm ad}$. 
\end{cor}

We shall see in Section~\ref{sec:11.4} 
that Theorem~\ref{thm:5.3.8} has the interesting consequence 
that one can characterize the integrability of a locally exponential Lie algebra $\g$ 
in terms of the discreteness of a certain subgroup of its center; 
the so-called period group. 
In this context, the main point of the preceding theorem is that 
it reduces the enlargeability problem 
for locally exponential Lie algebras to the enlargeability 
of the central extension $\g$ of $\g_{\rm ad}$.

\begin{prop} \mlabel{prop:semidir-Ad} Let $\g$ be a locally exponential 
Lie algebra, $|\g|$ the underlying locally convex space, and let 
$G_{\rm ad}$ be a simply connected locally exponential 
Lie group with Lie algebra $\g_{\rm ad} := \g/\fz(\g)$. 
Then $|\g| \rtimes G_{\rm ad}$ is a locally exponential Lie group. 
\end{prop}

\begin{prf} First we recall the existence of $G_{\rm ad}$ 
from the Adjoint Enlargeability Theorem~\ref{thm:5.3.8} 
which already asserts that the action of this group on $\g$ 
is smooth. Therefore $|\g| \rtimes_\alpha G_{\rm ad}$ is a Lie group. 

From Proposition~\ref{prop:5.2.10} we know that the 
Maurer--Cartan form of any exponential local Lie group $U_\g \subeq \g$ is given by 
\begin{equation}
  \label{eq:kappa-int1}
 \kappa_x = \int_0^1 e^{-t \ad x}\, dt.
\end{equation}
In particular, the map 
\[ U_\g \times \g \to \g, \quad 
(x,y) \mapsto \kappa_x y = T_0(\lambda_x^*)^{-1} y  
= T_x(\lambda_{-x}^*) y \] 
is smooth and the map 
\[ U_\g \times \g \to \g, \quad 
(x,y) \mapsto \kappa_x^{-1} y = T_0(\lambda_x^*) y \] 
is smooth as well. 
With the same argument as in the proof of 
Proposition~\ref{prop:semdir-exp} it now follows that 
\[  \exp(y,\ad x) := (\kappa_{-x}y, e^{\ad x}) \] 
is an exponential function of the semidirect product 
group $|\g| \rtimes_\alpha G_{\rm ad}$. Here we do not have
to assume that $\g$ is Mackey complete because the 
integrals \eqref{eq:kappa-int1} automatically exists pointwise in~$\g$ 
(see also Lemma~\ref{lem:Texp} for the exponential 
function of $TG \cong |\g| \rtimes G$ if $G$ exists). 
As $G_{\rm ad}$ is locally exponential and the map 
$(y,\ad x) \mapsto (\kappa_{-x}^{-1}y, \ad x)$ is smooth,
the assertion follows. 
\end{prf}

That the semidirect product Lie algebra $\fh := |\g| \rtimes_{\ad } \g_{\rm ad}$ 
is locally exponential also follows from the fact that 
$\fh \cong (|\g| \rtimes_{\ad} \g)/\{0\} \times \fz(\g)$ 
and that $T\g \cong |\g| \rtimes_{\ad}\g$ is locally exponential 
(Lemma~\ref{tang-locexp} and Corollary~\ref{cor:5.3.6}).

\section{Hausdorff series and locally exponential Lie algebras} 

This section is devoted to two results concerning the Hausdorff series. 
In the first subsection we show that, for a Banach--Lie algebra $\g$, 
it converges on a $0$-neighborhood and that it defines an analytic 
local group (Theorem~\ref{thm:loc-bangrp}). 
In the second subsection we show that, for any exponential local Lie group, 
the Taylor series of the multiplication in $(0,0)$ coincides with the 
Hausdorff series, regardless of its convergence.

\subsection{Banach--Lie algebras} \mlabel{subsec:5.4.4}

Let $\g$ be a Banach--Lie algebra. A norm $\|\cdot\|$ on $\g$ is said to 
be {\it compatible} \index{compatible norm, on Banach--Lie algebra} 
if it induces the given topology on $\g$ and satisfies 
$$\|[x,y]\| \leq \|x\|\cdot\|y\| \quad \mbox{ for all } x,y \in \g.$$ 
From now on we fix a compatible norm on $\g$ (cf.\ Exercise~\ref{exer:5.4.3}) 
and note that this implies that in the Banach algebra ${\cal L}(\g)$, we have 
$\|\ad x\| \leq \|x\|$. 

We want to show that the Hausdorff series converges on the set of all pairs 
$(x,y) \in \g \times \g$ with $\|x \| + \|y\| < \log 2$. 

\nin{\bf Step 1.} We start with the observation that for elements 
$a,b$ in an (associative) Banach algebra $\cA$ with a submultiplicative norm, we have 
$$ \|e^a - \1\| 
= \|\sum_{k = 1}^\infty \frac{a^k}{k!}\| 
\leq \sum_{k = 1}^\infty \frac{\|a\|^k}{k!} = e^{\|a\|} - 1. $$
This in turn leads to 
$$ \|e^a e^b - \1\| 
\leq \| (e^a - \1)e^b\| + \|e^b - \1\| 
\leq (e^{\|a\|} -1)e^{\|b\|} + e^{\|b\|} -1
=  e^{\|a\|+\|b\|} -1. $$
We conclude in particular, that 
$$\|a\| + \|b\| < \log 2 \quad \Rarrow \quad \|e^a e^b - \1\| < 1. $$

\nin{\bf Step 2.} Now we turn to the power series 
$$\Psi(X)=\frac{(X+1)\log (X+1)}{X} 
= (X+1)\sum_{k=0}^\infty \frac{(-1)^{k}}{k+1}X^k \in \R[[X]], $$
whose radius of convergence is $1$. For each element $a$ in a unital Banach algebra 
$\cA$ with $\|a\| < 1$, we then put 
$$ \Psi(a) := (a+1)\sum_{k=0}^\infty \frac{(-1)^{k}}{k+1}a^k  $$ 
and observe that 
$$ \|\Psi(a)\| \leq \|a+1\| \sum_{k=0}^\infty \frac{1}{k+1}\|a\|^k 
= (\|a\|+1) \frac{-\log(1 - \|a\|)}{\|a\|}.$$
The function on the right is increasing on the interval $[0,1[$ because it is 
defined by a power series with non-negative coefficients. 

\nin{\bf Step 3.} Returning to the Banach--Lie algebra $\g$, we assume that 
\break $\|x\| + \|y\| < \log 2$. Then Step 1 implies that 
$\|e^{\ad x} e^{t\ad y} -\1\|< 1$ for all $t \in [0,1]$, so that 
$\Psi(e^{\ad x} e^{t\ad y} -1)$ is defined. We may thus define 
$$ x * y := x + \int_0^1 \Psi(e^{\ad x}e^{t\ad y}-\1) y\, dt $$
because the integrand is a continuous function in $t$. 
Inspection of the second half of the proof of Theorem~\ref{thm:nilp-bch} shows that 
we may expand the absolutely convergent power series $\Psi$ and $e^{\ad x}$ etc.\ on the right hand side, 
which leads precisely 
to the Hausdorff series. We conclude that, for $\|x\| + \|y\| < \log 2$, 
this series converges to $x * y$. 

\nin{\bf Step 4.} Now we assume, in addition, that $\|x\|+\|y\| < \log\frac{3}{2}$. We then obtain 
for $t \in [0,1]$: 
$$ \|e^{\ad x}e^{t\ad y}-\1\| \leq e^{\|x\| + \|y\|} -1 \leq 
\textstyle{\frac{3}{2}} - 1 = \shalf, $$
and hence
$$ \|\Psi(e^{\ad x}e^{t\ad y}-\1)\| \leq 
 \frac{3}{2} \frac{\big(-\log(1 - \frac{1}{2})\big)}{\frac{1}{2}}= 3 \log 2. $$
We conclude that 
$\|x*y\| \leq \|x\| + 3\log 2 \cdot \|y\|.$
For $\|x\|, \|y\| < \frac{1}{8} \log 2$, this leads to 
$$ \|x*y\| < \textstyle{\frac{1}{8}}\log 2 + \frac{3}{8}\log 2 \log 2 \leq \frac{1}{2} \log 2. $$
In view of Step 3, this shows that, for 
\[ \|x\|, \|y\|, \|z\| < \frac{\log 2}{8}, \]
both products 
$(x*y) *z$ and $x*(y*z)$ are defined by convergent power series in $x$, $y$ and $z$. 

From Theorem~\ref{thm:BCH-free}, we know that the Hausdorff series defines on 
the completed free 
Lie algebra in three variables an exponential Lie group structure, so that 
projecting to the 
quotient by terms of order $n+1$, we see that the power series defining 
$(x*y)*z$ and $x*(y*z)$ coincide up to terms of order $n+1$. Since $n$ is arbitrary 
and the series converge, both products coincide. 

We collect the result of this discussion in the following theorem: 

\begin{thm} \mlabel{thm:loc-bangrp} Let $\g$ be a Banach--Lie algebra with the compatible norm $\|\cdot\|$ and 
$$ U := \Big\{ x \in \g \: \|x\| < \frac{\log 2}{8}\Big\}. $$
Then the Hausdorff series defines an analytic multiplication  on 
$$ \{ (x,y) \in \g \times \g \: \|x\| + \|y\| < \log 2\}, $$
and, for $x,y,z \in U$, both 
iterated products $(x*y)*z$ and $x*(y*z)$ are defined by 
the convergent Hausdorff series and coincide. With 
$$ D := \{ (x,y) \in U \: x*y \in U \} $$
we thus obtain an exponential local Lie group with Lie algebra $\g$, whose inversion is given by $x^{-1} = -x$. 
\end{thm}

\begin{prf} It only remains to prove the last assertion on the local group. 

Since $x*y = x + y$ holds for all pairs with $[x,y] =0$, we obtain in particular 
$x*(-x) = 0$ and 
\[  tx * sx = (t+s) x \quad \mbox{ for } \quad |t|,|s| \leq 1.\]

From that, (Loc1)-(Loc3) follow directly. To verify (Loc4), let 
$(x,y) \in D$. To see that $(-y,-x) \in D$, we claim that 
$(-y)*(-x) = -(x*y)$, so that (Loc4) follows from $U=-U$. 
The relation 
$$(-y)*(-x) = -(x*y) $$
holds for all pairs $(x,y) \in \g\times \g$ with $\|x\| + \|y\| < \log 2$ 
because it holds for the BCH multiplication in the exponential group structure on the 
completed free Lie algebra of two generators, which implies that the homogeneous terms of the 
series $-(x*y)$ and $(-x)*(-y)$ mutually coincide. 

In view of these remarks, the remaining assertions follow easily. 
\end{prf}

\begin{defn} \mlabel{def:bcdh-Lie-alg} 
A {\it BCH--Lie algebra} \index{BCH--Lie algebra} 
is a locally convex Lie algebra 
with the property that the Hausdorff series 
$x * y = \sum_{n = 1}^\infty H_n(x,y)$ converges 
on $U \times U$ for some open $0$-neighborhood $U \subeq \g$
and defines an analytic exponential local group structure on some 
$0$-neighborhood in~$U$. 
\end{defn}

\begin{rem} \mlabel{rem:bch-liealg}
(a) BCH--Lie algebras are in particular locally exponential. 

(b) The preceding theorem implies in particular that Banach--Lie algebras are BCH.

(c) If $\g$ is a projective limit of nilpotent Lie algebras, 
then it is BCH. In fact, we know already that $\g$ is exponential 
(Corollary~\ref{cor:pro-nilb}) 
and that the global multiplication $*$ is given by the convergent 
Hausdorff series. The same is true for the complexification 
$\g_\C$, so that $* \: \g_\C \times \g_\C \to \g_\C$ 
defines a complex analytic map. Restricting to $\g \times\g$, 
we thus obtain a real analytic map. 
This proves that $\g$ is BCH.
\end{rem}

\subsection{Locally exponential Lie algebras} \mlabel{subsec:5.4.5}

Let $(G,D,m,0)$ be an exponential local Lie group and 
consider the Taylor series 
\begin{eqnarray}
  \label{eq:5.2.8}
x * y = \sum_{n = 1}^\infty b_n(x,y) 
\end{eqnarray}
of its multiplication map in $(0,0)$. We know already that 
$$b_1(x,y) = x+ y \quad \mbox{ and } \quad b_2(x,y) = \frac{1}{2} [x,y] $$
(Remark~\ref{rem:5.2.4}).  
The main result of this subsection states that (\ref{eq:5.2.8}) is the 
Hausdorff series.

\begin{lem} \mlabel{tang-locexp} Let $(G,D,m_G,0)$ be an exponential local 
Lie group and write $|\g|$ for the space $\g = \L(G)$, endowed with the 
zero Lie bracket.  Then $(TG,TD, Tm_G,0_0)$ is an exponential local 
Lie group with Lie algebra 
$T\g := |\g|\rtimes_{\ad} \g$. 
\end{lem}

\begin{prf} Applying the tangent functor $T$ to the local Lie group 
$(G,D,m_G,0)$, we get on the domain 
$TD \subeq TG \times TG \cong T(G \times G)$ a smooth multiplication 
$Tm_G \: TD \to TG$. First we show that we obtain a local Lie group 
(cf.\ Exercise~\ref{exer:4.1.3}). 

On the domain 
$$D^1 := \{ (x,y,z) \in G^3 \: (x,y), (y,z), (x*y,z) \in D\}, $$
the associativity condition (Loc1) reads 
$$ m_G \circ (m_G \times \id_G) = m_G \circ (\id_G \times m_G) \: D^1 \to G. $$
Applying $T$, we get 
$$ Tm_G \circ (Tm_G \times \id_{TG}) = Tm_G \circ (\id_{TG} \times Tm_G) $$
on 
$$ T(D^1) = \{ (x,y,z) \in T(G)^3 \: (x,y), (y,z), (x*y,z) \in T(D)\} = (TD)^1. $$
Together with the observation that, for $(x,y) \in T(D)$, the relation 
$$(Tm_G(x,y),z) \in T(D)$$ implies $(x, Tm_G(y,z)) \in T(D)$, and vice versa, 
we obtain (Loc1) for $TG$. 

For each $x \in TG$, we clearly have 
$(x,0_0), (0_0,x) \in TD$ with 
$$ x * 0_0 = (Tm_G)(x,0_0) = x = (Tm_G)(0_0,x) = 0_0 * x, $$
which is (Loc2). 

To verify (Loc3), we write $\eta_G \: G \to G$ for the inversion map, 
which satisfies 
$$m_G \circ (\eta_G \times \id_G) = m_G \circ (\id_G \times \eta_G) = 0. $$
Applying $T$ leads to 
$$ Tm_G \circ (T\eta_G \times \id_{TG}) = Tm_G \circ (\id_{TG} \times T\eta_G) = 0_0, $$
which implies (Loc3) for the inversion map $T\eta_G = \eta_{TG}$ of $TG$. 

The last condition (Loc4) means that $D$ is invariant under the map 
$$ F \: G \times G \to G \times G, \quad (x,y) \mapsto (\eta_G(y), \eta_G(x)), $$
which implies that $TD$ is invariant under 
$$ TF \: TG \times TG \to TG \times TG, \quad (x,y) \mapsto (T\eta_G(y), T\eta_G(x)) 
= (\eta_{TG}(y), \eta_{TG}(x)), $$
and this is (Loc4) for $TG$.

To see that $Tm_G$ satisfies the exponentiality condition
\begin{eqnarray}
  \label{eq:5.1.2b}
tx * sx = (t+s) x\quad \mbox{ for} \quad |s|,|t|, |s+t| \leq 1, 
\end{eqnarray}
we write $\mu_s$ for the multiplication by 
$s \in \R$ and $\Delta_G(x) := (x,x)$ for the diagonal map. Then 
\eqref{eq:5.1.2b} reads 
$$ m_G \circ (\mu_t \times \mu_s) \circ \Delta_G = \mu_{t + s} \quad \mbox{for}  \quad 
|s|, |t|, |s+t| \leq 1. $$
Applying the tangent functor $T$, leads to 
\begin{eqnarray}
  \label{eq:5.2.9}
Tm_G \circ (T\mu_t \times T\mu_s) \circ T\Delta_G = T\mu_{t + s}. 
\end{eqnarray}
Next we observe that $(T\mu_s)(x,v) = (sx,sv)$ 
is the scalar multiplication by $s$ on the vector space 
$T\g \cong \g \times \g$. Further $T\Delta_G = \Delta_{TG}$, so that (\ref{eq:5.2.9}) 
is \eqref{eq:5.1.2b} for the local group $(TG,TD,Tm_G,0_0)$. 

To determine the Lie algebra of the local group $TG$, 
we recall that the product on $TG$ is given by 
$$ (x,v) * (x',v') = (x * x', T_{x'}(\lambda_x)v' + T_x(\rho_{x'})v ). $$
Since $T_{x'}(\lambda_x)v'$ is linear in $v'$, constant $v'$ for $x = 0$, and 
for $x' = 0$ we have $T_0(\lambda_x)v' = v' + \frac{1}{2}[x,v'] + \cdots$, 
the Taylor polynomial of order $2$ of $Tm_G$ in $(0,0)$ is given by 
$$ \Big(x + x' + \textstyle{\frac{1}{2}}[x,x'], v + v' + \textstyle{\frac{1}{2}}[x,v'] 
- \textstyle{\frac{1}{2}} [x',v]\Big). $$
Hence the Lie bracket on $T\g$ is given by 
$$ [(x,v), (x',v')] = ([x,x'], [x,v'] - [x',v]). $$
This means that $T\g \cong \g \ltimes_{\ad} |\g|$. 
\end{prf}

To iterate the construction of the preceding lemma, 
we first write the Lie algebra $T\g := \L(TG) \cong |\g| \rtimes_{\ad} \g$ differently. 
The {\it dual numbers} \index{dual numbers} 
$\R[\eps]$ form  the $2$-dimensional real algebra 
$\R \oplus \R \eps$ with $\eps^2 = 0$. With this commutative algebra, we obtain a 
description of $T\g$ by scalar extension from $\R$ to $\R[\eps]$: 
$$ T\g \cong \g \otimes_\R \R[\eps] = \eps \g \rtimes \g \cong |\g|\rtimes_{\ad}\g. $$

Iterating this procedure, we obtain a sequence of exponential local Lie groups 
$$ (T^n G, T^n D, T^n m_G,0) \quad \mbox{ with } \quad 
\L(T^n G) \cong T^n \g \cong \g \otimes_\R T^n(\R), $$
where 
$$ T^n(\R) := \R[\eps_1, \ldots, \eps_n],\quad \eps_i \eps_j = \eps_j \eps_i, 
\quad \eps_i^2 = 0. $$ 

On the higher tangent bundles $T^n G$, we have a natural action of 
the symmetric group $S_n$ whose fixed point set is the 
\index{jet bundle}
{\it $n$th jet bundle}
$J^n G$. On the Lie algebra level, we have 
$$ J^n \g := (T^n \g)^{S_n} \cong (\g \otimes_\R T^n(\R))^{S_n} 
\cong \g \otimes_\R T^n(\R)^{S_n}, $$
where $S_n$ acts on $T^n(\R)$ by permuting the generators $\eps_1, \ldots, \eps_n$.
The fixed point set for this action is the unital subalgebra spanned by the elements 
$$ \delta^{(k)} := \sum_{i_1 < \ldots < i_k} \eps_{i_1} \cdots \eps_{i_k}, $$
satisfying 
$$ \delta^{(k)} \delta^{(\ell)} = 
   \left\{
     \begin{array}{cl} 
      {k + \ell \choose k} \delta^{(k+ \ell)} 
         & \mbox{ for $k + \ell \leq n$} \\ 
      0  & \mbox{ for $k + \ell >n.$} 
     \end{array} \right. $$
Restricting to $J^n G$, naturality of the $S_n$-action implies that 
$J^n\g$ is a subalgebra of $T^n\g$ and that $m_G$ induces a smooth multiplication 
$$ J^n(m_G) \: J^n(D) \to J^n(G), $$
turning $J^n(G)$ into an exponential local Lie group. 

\begin{thm} \mlabel{bch-locexp} Let $(G,D,m_G,0)$ be an exponential local Lie group. 
Then the Taylor series of the multiplication $m_G$ in $(0,0)$ coincides with the 
Hausdorff series. 
\end{thm}

\begin{prf} It clearly  suffices to show  that the $n$-jet of $m_G$ in $(0,0)$ 
is the sum of all terms up to order $n$ of the Hausdorff series 
$\sum_{r,s \leq n} H_{r,s}(x,y)$, where $H_{r,s}(x,y)$ denotes the summand 
homogeneous of degree $r$ in $x$ and of degree $s$ in $y$.

In terms of the terminology introduced above, the $n$-jet of a smooth function 
$f$ in a point $x$ satisfies for $\delta := \eps_1 + \ldots +\eps_n$: 
\begin{eqnarray*}
(J^n f)(x + \delta h) 
&=& (T^n f)(x + \eps_1 h + \cdots + \eps_n h) \\
&=& \sum_{k = 0}^n \delta^{(k)} (\dd^k f)(x)(h,\ldots, h) 
= \sum_{k = 0}^n \frac{\delta^k}{k!} (\dd^k f)(x)(h,\ldots, h).
\end{eqnarray*}
This permits us to obtain the $n$th order Taylor polynomial of $f$ in $x$ by 
evaluating $J^n f$ in $x + \delta h$. We therefore have to show that 
$$ (J^n m_G)(\delta x, \delta y) 
= \sum_{k = 0}^n \delta^{k} \sum_{r+s = k} H_{r,s}(x,y). $$

To verify this relation, let $\eps_{T^n(\R)} \: {T^n(\R)} \to \R$ 
be the augmentation homomorphism. 
This is the unique algebra homomorphism ${T^n(\R)} \to \R$ whose kernel 
is the hyperplane ideal generated by the elements $\eps_1, \ldots, \eps_n$. 
Tensoring with $\id_\g$, this leads to a homomorphism 
$E := \id_\g \otimes \eps_{T^n(\R)} \: T^n \g \cong \g \otimes_\R {T^n(\R)} \to \g$ of Lie algebras. 
We then have 
$$ \ker E \cong \g \otimes_\R (\ker \eps_{T^n(\R)}), $$
and since 
$\ker \eps_{T^n(\R)} = \sum_{i =1 }^n \eps_i {T^n(\R)}$
is a nilpotent commutative algebra, $\ker E$ is a nilpotent 
Lie algebra. Therefore Theorem~\ref{thm:5.2.13} implies that the 
multiplication on $\ker E$ is given by 
the Hausdorff series. We thus obtain in particular in $J^n\g$ the relation 
\begin{eqnarray*}
&& (J^n m_G)(\delta x,\delta y) = \delta x * \delta y
= \sum_{r,s \leq n} H_{r,s}(\delta x,\delta y)
= \sum_{r,s \leq n} \delta^{r+s} H_{r,s}(x,y), 
\end{eqnarray*}
and this completes the proof. 
\end{prf}

\subsection*{Open problems in Chapter~\ref{ch:7}} 

\begin{probl} \mlabel{prob:5.5.1} (Assumptions of the Quotient Theorem) 
Suppose that \break $\fn \trile \g$ is a closed ideal of the locally 
exponential Lie algebra $\g$. 
\begin{enumerate}
\item[(1)] Is the stability of $\fn$ necessary for $\fn$ being locally exponential? 
\item[(2)] Does the stability of $\fn$ imply that 
$\kappa_\g(x)(\fn) = \fn$ holds for all $x$ in a $0$-neighborhood? 
\item[(3)] Suppose that $\fn$ is stable and locally exponential. Does this imply 
$\kappa_\g(x)(\fn) = \fn$ for $x$ close to $0$? 
\end{enumerate}
\end{probl}

The following problem describes the functional analytic essence of the 
preceding problem. 
\begin{probl} \mlabel{prob:5.5.2} Let $E$ be a Mackey complete 
space, $\alpha \: \R \to \GL(E)$ a smooth action of $\R$ on $E$ 
and $D = \alpha'(0)$ its infinitesimal generator. 
Let 
$$\beta(t) := \int_0^1 \alpha(st)\, dt. $$
We assume that $\beta(t) \in \GL(E)$ for $|t| \leq \eps$, and that 
the corresponding map $\tilde\beta \: [-\eps,\eps] \times E \to E \times E$ is smooth. 
Suppose that $W \subeq E$ is a closed subspace invariant under $\alpha(\R)$, 
which clearly implies that $\beta(t)(W) \subeq W$ for all $t \in \R$. 
Does this imply that 
\begin{eqnarray}
  \label{eq:prob}
\beta(t)(W) = W \quad \hbox{ for } \quad |t| \leq \eps? 
\end{eqnarray}

Note that $W$ is a closed ideal of 
$\g = E \rtimes_D \R$, and that the condition $\alpha(\R)(W) \subeq W$ is equivalent 
to the stability of this ideal. The condition $\kappa_\g(x)(W) = W$ for 
$x$ close to $0$ is equivalent to (\ref{eq:prob}). 

For every closed stable ideal $\fn$ of a locally exponential 
Lie algebra $\g$ and $x \in \g$ we obtain a smooth $\R$-action on $E := \g$ 
by $\alpha(t) := e^{t\ad x}$, and the corresponding operator 
$\beta(t)$ equals $\kappa_\g(-tx)$. This correspondence shows that 
the Lie theoretic Problem~\ref{prob:5.5.1}(2)  
 is closely related to the 
operator theoretic one which only deals with a very simple type of 
$2$-step solvable Lie algebras. 

Since, for every smooth $\R$-module $E$, we have an equivariant embedding 
$E \into C^\infty(\R,E)$, the smooth modules of the form 
$E = C^\infty(\R,V)$ with $(\alpha_t f)(x) = f(t+x)$ and $Df = f'$ 
are of central importance. 
\end{probl}

\begin{probl} \mlabel{prob:5.1} (One-parameter groups and local 
exponentiality) 
Let \break $\alpha \: \R \to \GL(E)$ define a smooth action of 
$\R$ on the Mackey complete 
locally convex space and $D := \alpha'(0)$ be its infinitesimal 
generator. Characterize the local exponentiality of 
$G := E \rtimes_\alpha \R$ in terms of the infinitesimal generator~$D$. 

Writing the exponential function as 
$\exp_G(v,t) = (\beta(t)v, t)$ with
$\beta(t) = \int_0^1 \alpha(st)\, ds$ (Proposition~\ref{prop:semdir-exp}), 
we obtain the curve $\beta \: \R \to {\cal L}(E)$. 
We are looking for a characterization of those operators 
$D$ for which there exists some $T > 0$ such that 
\begin{enumerate}
\item[\rm(1)] $\beta(]{-T},T[) \subeq \GL(E)$, and 
\item[\rm(2)] $\tilde\beta \: ]{-T},T[\times E \to E, (t,v) 
\mapsto \beta(t)^{-1}v$ is smooth. 
\end{enumerate}

\nin Note that $(t,v) \mapsto \beta(t)v$ is always smooth. 
If $E$ is a Banach space, then $G$ is a Banach--Lie group, 
hence locally exponential. In this case $D$ is a bounded operator 
and we have for each $t \not=0$: 
\[  \beta(t) 
= \frac{1}{ t} \int_0^t e^{sD}\, ds 
= \frac{1}{ t} \frac{e^{tD} - \1}{D} 
= \frac{e^{tD} - \1}{tD} 
= \sum_{k = 0}^\infty \frac{t^k}{(k+1)!} D^k. \]
As $\beta \: \R\to {\cal L}(E)$ is continuous w.r.t.\ to the operator norm on 
${\cal L}(E)$, $\beta(0) = \be$, and $\GL(E)$ is open, 
conditions (1) and (2) follow immediately. 
Moreover, the Spectral Theorem implies that 
$$ \Spec(\beta(t)) = 
\Big\{ \frac{e^{t\lambda} - 1}{t\lambda} \: \lambda \in \Spec(D)\Big\}, $$
which means that $\beta(t)$ is invertible if $t$ is sufficiently small. 
\end{probl}

\begin{probl} \mlabel{prob:5.1b} 
Characterize closed subalgebras $\h$ of a locally exponential 
Lie algebra $\g$ which are locally exponential. 
In view of Lemma~\ref{lem:5.3.4}, 
stability is a necessary condition. 
We also need $\kappa_\g(x) \h = \h$ for $x \in \h$ close to~$0$. Are these two 
conditions sufficient?   
\end{probl}


\begin{small}
\subsection*{Exercises for Chapter~\ref{ch:7}}

\begin{exer} \mlabel{exer:3.7.1} 
Let $(G,D,m_G,\be)$ be a local group and $U \subeq G$ a symmetric subset
containing $\be$. Show that 
$(U,D_U,m_U,\be)$ is a local group for 
\[ D_U := \{ (x,y) \in D \cap (U \times U) \: xy \in U\} \quad \mbox{ and } \quad 
m_U := m_G\res_{D_U}.\] 
\end{exer} 

\begin{exer} \mlabel{exer:4.1.2} 
Let $E$ be a vector space, $\g$ a Lie algebra and 
$\phi \: E \to \g$ a linear isomorphism. Show that 
$$ [x,y] := \phi^{-1}([\phi(x),\phi(y)]) $$
defines on $E$ a Lie algebra structure for which $\phi \: E \to \g$ 
is an isomorphism of Lie algebras. 
\end{exer}

\begin{exer} \mlabel{exer:4.1.3} 
Show that on the tangent bundle $T(G)$ of a local Lie group $G$, 
there exists a local Lie group structure  $(TG,TD,Tm_G,0_\be)$. 
For this Lie group structure the zero section 
$\sigma \: G \to TG$ 
and the bundle projection 
$\pi_G \: TG \to G$
are morphism of local Lie groups satisfying $\pi_G \circ \sigma_G = \id_G$. 
\end{exer}

\begin{exer} \mlabel{exer:4.l.4} 
Let $(G,D,m_G,\be)$ be a local Lie group and 
$D_0$ the connected component of $(\be,\be)$ in $D$. Assume that 
$G$ is connected. Show that 
the restriction $m_G^0 := m_G \res_{D_0}$ satisfies (Loc2), (Loc3) and (Loc4). 
\end{exer}

\begin{exer} \mlabel{exer:4.l.5} (Homomorphisms of local groups) 
Let $(G,D_G,m_G,\be_G)$ and \break $(H,D_H,m_H,\be_H)$ be local groups 
and $\phi \: G \to H$ a map satisfying 
\[ (\phi \times \phi)(D_G)\subeq D_H \]  and 
$$ \phi(x*y) = \phi(x) * \phi(y) \quad \mbox{ for } \quad (x,y) \in D_G. $$ 
Show that $\phi(\be_G) = \be_H$. 
\end{exer} 

\begin{exer} \mlabel{exer:4.l.6} Let $(G,D, m_G,\be)$ be a local group. 
One of the most subtle points of local groups is the validity of the associative 
law for iterated products. Let $(G_n)_{n \in \N}$ be a decreasing sequence of 
local subgroups of $G$ with 
$$G_0 = G, \quad G_n \times G_n \subeq D \quad \mbox{ and } \quad G_n * 
G_n \subeq G_{n-1} \quad \mbox{ for } \quad n \in \N. $$
Show that: 
\begin{enumerate}
\item[\rm(1)]  For $x,y,z \in G_2$ the products $(x * y) * z$ and $x*(y*z)$ are defined and equal. 
\item[\rm(2)]  For elements of $G_n$ all $n$-fold products (in any order) are defined and associative. 
Hint: For $x_1, \ldots, x_n \in G_n$ define 
$$ x_1 * x_2 * \cdots * x_n := (\cdots ((x_1 * x_2) * x_3) * \cdots x_n) $$
and proceed by induction. 
\end{enumerate}
\end{exer} 

\begin{exer} \mlabel{exer:4.l.7} Let $(G,D, m_G,0)$ be a local Lie group, where 
$G$ is an open subset of the locally convex space $\g$. 
If $b(x,y)$ is the second order term in the Taylor polynomial 
$x*y = x + y + b(x,y)$ of order $2$, we put 
$$ [x,y] := b(x,y)- b(y,x)\quad \mbox{ for } \quad x,y \in \g. $$
Show directly, without using vector fields, that: 
\begin{enumerate}
\item[\rm(1)] The bracket $[\cdot,\cdot]$ on $\g$ is continuous, bilinear and alternating. 
\item[\rm(2)] The second order Taylor expansion in $0$ of the inversion 
is 
$$ x^{-1} = -x + b(x,x) + \ldots. $$
\item[\rm(3)] The second order Taylor expansion in $(0,0)$ of the conjugation map is 
$$ c_x(y) := (x * y) * x^{-1} = x * (y * x^{-1}) = y + [x,y]  + \ldots. $$
\item[\rm(4)] The second order Taylor expansion in $(0,0)$ of the commutator map is 
$$ (x,y) := x * y * x^{-1} * y^{-1} = [x,y]  + \ldots. $$
\item[\rm(5)] We have the third order Taylor expansion in $(0,0,0)$: 
$$ ((x,y), c_y(z)) = \big[[x,y],z\big] + \ldots  $$
Hint: Chain Rule for Taylor polynomials 
(Proposition~\ref{chainrtay}), (1) and (2). 
\item[\rm(6)] For $x,y,z$ sufficiently close to $0$ we have Hall's identity 
$$ ((x,y), c_y(z)) * ((y,z), c_z(x)) * ((z,x), c_x(y)) = \be. $$
\item[\rm(7)] The bracket $[\cdot,\cdot]$ satisfies the Jacobi identity. 
\end{enumerate}
\end{exer} 


\begin{exer} \mlabel{exer:5.6.2} 
Let $M$ be a compact smooth manifold (with boundary) 
and $G$ a local Lie group. 
Show that, for each $r \in \N \cup \{\infty\}$,  the group $C^r(M,G)$, 
endowed with the compact open $C^r$-topology, 
carries a natural local Lie group structure, given by pointwise 
multiplication. 
\end{exer} 
  
\begin{exer} \mlabel{exer:e.10} Let $V$ be a locally convex space and 
$D \in \gl(V)$ a continuous linear operator. Assume that for some $\eps > 0$ the curve 
$$ \alpha \: [-\eps,\eps] \to \GL(V) $$
is smooth in the sense that the map 
$$ \hat\alpha \: [-\eps,\eps] \times V \to V \times V, \quad (t,v) \mapsto (\alpha(t)(v), 
\alpha(t)^{-1}(v)) $$
is smooth, and that $\alpha$ satisfies the IVP 
$$ \alpha'(t) = \alpha(t) \circ D, \quad \alpha(0) = \id_V. $$
Show that: 
\begin{enumerate}
\item[\rm(1)] $\alpha(t+s) = \alpha(t)\alpha(s)$ for $|t|,|s|,|t+s|\leq \eps$. 
\item[\rm(2)] $\alpha(\frac{t}{n})^n = \alpha(t)$ for $|t| \leq \eps$. 
\item[\rm(3)] $\alpha(t) := \lim_{n \to \infty} \alpha\big(\frac{t}{n}\big)^n$
defines an extension of $\alpha$ to a smooth $\R$-action on $V$. 
\end{enumerate}
\end{exer}

\end{small}

%
%
%
%


\chapter{Subgroups and quotients of Lie groups}  \mlabel{ch:6}
It is one of the fundamental problems of Lie theory to understand 
to which extent subgroups of Lie groups carry natural Lie group 
structures. Similarly, one wants to 
understand under which conditions quotients of Lie groups by 
closed normal subgroups are Lie groups. Clearly, the closedness 
of the normal subgroup is necessary for the quotient to be Hausdorff 
and for finite-dimensional groups this condition is already sufficient.
However, already for Banach--Lie groups, the situation is more complicated. 


This chapter is devoted to these two problems. 
It consists of four sections. The first one contains general material 
on subgroups and quotients of Lie groups. It provides 
definitions and concepts and derives some of the results 
familiar from the finite-dimensional theory in the infinite-dimensional 
context. This often means to find the appropriate set of assumptions which 
make the constructions work. Here we introduce the term 
{\it Lie subgroup} 
\index{Lie subgroup} 
\index{Lie subgroup!split} 
for a subgroup $H \leq G$ which is a submanifold of $G$, hence in particular 
a Lie group in its own right. We also discuss sufficient 
conditions on a Lie subgroup $H$ ensuring 
that the quotient $G/H$ 
carries a natural manifold structure; then we call $H$ a 
{\it split Lie subgroup}.

In Sections~\ref{sec:6.2} and \ref{sec:6.3} we focus on subgroups and quotients of 
locally exponential Lie groups. This is justified by the fact that the 
assumption of local exponentiality provides the necessary framework to 
develop powerful general results on subgroups and quotients. 

There are two natural classes of subgroups of locally exponential 
Lie groups. The first class, the integral subgroups, consists of those subgroups
$H \leq G$ for which there exists an injective morphism of locally exponential 
Lie algebras $\h \into \L(G)$ for which $H = \la \exp_G \h \ra$ 
carries a locally exponential Lie group structure with $\L(H) = \h$. These groups are natural 
generalizations of the integral subgroups (sometimes called 
``analytic subgroups'') of finite-dimensional Lie groups. 
It is an important point that, for the theory of integral subgroups,  
it suffices to a large extent that $G$ has an exponential function and that 
$\h$ is locally exponential, which permits us to push some of the main theorems 
beyond the class of locally exponential groups. 
The second natural class of subgroups are those which are locally exponential 
Lie groups when endowed with the topology inherited from~$G$. Since any 
topological group carries at most one structure of a locally exponential 
Lie group, we call these groups 
\index{Lie subgroup!locally exponential} 
{\it locally exponential Lie subgroups}. 
Any closed locally exponential Lie subgroup is a submanifold, 
hence a Lie subgroup in the sense defined above, and this condition is 
equivalent to the closedness of its Lie algebra (which is automatic in 
finite dimensions). 

Quotients are discussed in Section~\ref{sec:6.4}, where we show in particular 
that, for a closed normal subgroup $N$ of a Banach--Lie group $G$, the quotient 
$G/N$ is a Banach--Lie group if and only if $N$ is a Banach--Lie subgroup. 
If $G$ is only locally exponential, we have to assume, in addition, that 
$\L(G)/\L(N)$ is a locally exponential Lie algebra to get the Lie group 
structure on $G/N$. 

We conclude this chapter with a discussion of open problems related 
to manifold structures on quotients, orbits of Lie group  actions and 
the information available in this context. 

 {\bf Prerequisites:} Sections~\ref{seccuint}-\ref{secCk} and 
\ref{secTay}, and Chapter~\ref{chapmanif} (Sections 1-4), 
Chapter~\ref{ch:3} (Sections 1-3), Chapter~\ref{ch:4} (Section 1) and 
Chapter~\ref{ch:5} (Section 2).

\section{Subgroups of general Lie groups} \mlabel{sec:6.1}

To develop a flexible infinite-dimensional Lie theory, 
it is important to set up the concepts in the least restrictive way, still 
preserving the crucial part of the interesting structure. In this spirit, we 
start with a discussion of the Lie algebra of a {\it general} subgroup $H$ of a 
Lie group $G$. On a very general level, we can also ask if a subgroup $H$ 
of a group $G$ carries a Lie group structure which is determined by 
the manifold structure of $G$ in the sense that a map $M \to H$ from a smooth manifold 
$M$ into $H$ is smooth if and 
only it is smooth as a map $M \to G$. We call these Lie groups 
\index{Lie subgroup!initial} 
{\it initial Lie subgroups}. Since the definition easily implies that each 
subgroup carries at most one initial Lie group structure, an important point 
is to find sufficient criteria for a given subgroup to carry an initial 
Lie subgroup structure. It is surprising that the answer to this question 
is positive in an amazing generality. Conversely, it is quite hard to 
prove the non-existence of initial Lie group structures on subgroups. 

We conclude this section with a general construction of manifold 
structures on homogeneous spaces $G/H$.

\subsection{The differential Lie algebra of a subgroup} 

A first step to a Lie group structure on a subgroup $H$ of 
a Lie group $G$ is to find a natural candidate for the Lie algebra 
of~$H$. In the following, we write $C^1_*([0,1],G)$ for the set of 
all $C^1$-curves $\gamma \: [0,1] \to G$ with $\gamma(0) = \be$. 
We shall now verify that the following definition provides a Lie algebra for {\em any} subgroup 
$H \leq G$:  
\begin{equation}
  \label{eq:Ld}
\L^d(H) := \{ \alpha'(0) \in \L(G) = T_\be(G) \: \alpha \in C^1_*([0,1],G),\ \im(\alpha) 
\subeq H\}. 
\end{equation}

For the proof that this always is a Lie algebra, we need the following proposition 
and its corollary. 

\begin{prop} \mlabel{prop:diff-sqrt} Let $E,F$ be locally convex spaces,  
let $U\subeq E$ be an open $0$-neighborhood 
and $f \in C^3(U,F)$ with $f(0) = 0$ and $\delta_0 f = 0$. Then, for each 
$\eta \in C^1([0,1],U)$ with $\eta(0) = 0$, the prescription 
$\gamma(t) := f(\eta(\sqrt t))$ defines a $C^1$-curve $[0,1] \to F$ with 
\begin{eqnarray}
  \label{eq:diff0}
\gamma'(0) 
=  \frac{1}{2} \delta_0^2 f(\eta'(0)). 
\end{eqnarray}
\end{prop}

\begin{prf} Clearly, $\gamma$ is $C^1$ in $]0,1]$ and it remains 
to show that it is differentiable in $0$ and that $\gamma'$ is continuous. 
From the Taylor expansion of $f$ in $0$, we obtain 
\[  f(h) = \frac{1}{2}\delta_0^2 f(h) + \frac{1}{2} 
\int_0^1 (1-t)^2 \delta^3_{th} f(h)\, dt \]
(Theorem~\ref{thmtay}). 
We therefore have 
\[  \frac{f(\eta(s))}{s^2} 
= \frac{1}{2}\delta_0^2 f\Big(\frac{\eta(s)}{s}\Big) + \frac{1}{2} 
\int_0^1 (1-t)^2 s\delta^3_{t\eta(s)} f\Big(\frac{\eta(s)}{s}\Big)\, dt. \]
Since $\tilde\eta(0) := \eta'(0)$ yields a continuous extension of the curve 
$\tilde\eta(s) := \frac{\eta(s)}{s}$ to $[0,1]$, we conclude that 
\begin{eqnarray*}
\lim_{s \to 0} \frac{f(\eta(s))}{s^2} 
&=& \frac{1}{2}\delta_0^2f(\eta'(0)).  
\end{eqnarray*}
Hence $\gamma$ is differentiable in $0$ with $\gamma'(0)$ given by \eqref{eq:diff0}. 
In $t \not=0$ we have 
$$ \gamma'(t) = \dd f\big(\eta(\sqrt t), \eta'(\sqrt t)\big) 
\cdot \frac{1}{2 \sqrt t}, $$
which is continuous on $]0,1]$. 
Writing 
\[ \gamma'(t) = h(1) \quad \mbox{ with  } \quad 
h(s) = \dd f\big(s\eta(\sqrt t), \eta'(\sqrt t)\big) 
\cdot \frac{1}{2 \sqrt t},\] 
we obtain 
\begin{eqnarray*} 
\gamma'(t) 
&=& \int_0^1 h'(s)\, ds 
= 
\frac{1}{2} \int_0^1 d^{(2)} f\Big(s \eta(\sqrt t), 
\eta'(\sqrt t), \frac{\eta(\sqrt t)}{\sqrt t}\Big)
\, ds.
\end{eqnarray*}
Since $\eta$ is differentiable in $0$ and $d^{(2)}f$ is continuous, 
the integrand extends to a continuous function on $[0,1]^2$, which leads to 
\[\lim_{t \to 0} \gamma'(t) 
= \frac{1}{2} \int_0^1 d^{(2)} f(0,\eta'(0), \eta'(0))\, ds
= \frac{1}{2} \delta^2_0 f(\eta'(0)) = \gamma'(0). 
\qedhere\] 
\end{prf}

\begin{cor} \mlabel{cor:sqrt} Let $\alpha, \beta \in C^1_*([0,1],G)$ and 
$$\gamma(t) := \alpha(\sqrt t)\beta(\sqrt t)\alpha(\sqrt t)^{-1}\beta(\sqrt t)^{-1}. $$
Then $\gamma \in C^1_*([0,1],G)$ with $\gamma'(0) = [\alpha'(0), \beta'(0)]$.   
\end{cor}

\begin{prf} Expressing everything in a local chart $(\phi,U)$ of $G$ with 
$\phi(\be) = 0$ and $T_\be(\phi) = \id_{\L(G)}$, we apply 
Proposition~\ref{prop:diff-sqrt} with 
\[ \eta(t) = (\alpha(t), \beta(t)) \quad \mbox{ and } \quad 
\gamma(t) = f(\eta(\sqrt t)) \quad \mbox{ for } \quad 
f(x,y) = x* y * x^{-1} * y^{-1}.\]
We know from 
Remark~\ref{rem:brack-taylor} that $\dd f(0,0) = 0$ and 
$\frac{1}{2}\delta_0^2 f(x,y) = [x,y],$
so that the assertion follows. 
\end{prf}

\begin{prop} \mlabel{prop:diff-liealg} Let $H \leq G$ be a subgroup of the Lie group $G$. Then 
$$ \L^d(H) = \{ \alpha'(0) \in \L(G)\: \alpha \in C^1_*([0,1],G),\ \im(\alpha) \subeq H\} $$
is a Lie subalgebra of $\L(G)$. 
\end{prop}

\begin{prf} If $\alpha,\beta \in C^1_*([0,1],G)$, then 
$$ (\alpha\beta)'(0) = \alpha'(0) + \beta'(0), \quad 
(\alpha^{-1})'(0) = -\alpha(0), $$
and, for $0 \leq \lambda \leq 1$, the curve $\alpha_\lambda(t) := \alpha(\lambda t)$ 
satisfies $\alpha_\lambda'(0) = \lambda \alpha'(0)$. Therefore 
$\L^d(H)$ is a real linear subspace of $\L(G)$. 

Using Corollary~\ref{cor:sqrt}, 
we see that the curve 
$$\gamma(t) := \alpha(\sqrt t)\beta(\sqrt t)\alpha(\sqrt t)^{-1}\beta(\sqrt t)^{-1} $$
with $\gamma(0) = \be$ is $C^1$ with $\gamma'(0) 
= [\alpha'(0), \beta'(0)]$. Hence $\L^d(H)$ is a Lie subalgebra of $\L(G)$. 
\end{prf}

We put the superscript $d$ for ``differentiable'' 
to distinguish $\L^d(H)$ from the Lie algebra of a Lie group. Later, we shall encounter 
another way to define the Lie algebra of a subgroup which works well for closed 
subgroups of locally exponential Lie groups. 

\begin{rem} For any $C^1$-curve 
$\alpha \: [0,1] \to H$ into a subgroup $H$ of the Lie group $G$, the continuous 
curve $\delta(\alpha) \in C^0([0,1],\L(G))$ 
has values in $\L^d(H)$. In fact, for any $t_0 \in [0,1]$, 
the $C^1$-curve $\beta(t) := \alpha(t_0)^{-1}\alpha(t_0+t)$ satisfies 
$\beta(0) = \be$ and $\beta'(0) = \delta(\alpha)_{t_0}$. 
Note that $\beta$ is defined on $[-t_0, 1-t_0]$ and not on $[0,1]$, 
but this can be adjusted easily by reparameterization. 

Conversely, if $H$ is a subgroup of $G$ and 
\[ \alpha \in C^1_*([0,1],G) \quad\mbox{ with } \quad
  \delta(\alpha) \in C^0([0,1],\L^d(H)),\] 
it is not clear which natural additional condition we have to impose 
on $H$ to conclude that $\alpha(1) \in H$. 
We shall see below  that the concept of regularity helps to address  this problem. 
\end{rem}

\subsection{Initial Lie subgroups} 

In this subsection,  
we briefly discuss the rather weak concept of an initial Lie subgroup. As a consequence 
of the universal property built into its definition, such a structure is unique 
whenever it exists. We prove its existence in two cases: if each $C^1$-arc in $H$ is 
constant and if $H$ is the center and $G$ satisfies some mild extra conditions.

\begin{defn} (Initial Lie subgroups) \mlabel{def-init-Lie} 
An injective morphism $\iota \: H \to G$ of Lie groups 
is called {\it $C^k$-initial} \index{Lie subgroup!$C^k$-initial}
($k \in \N \cup \{\infty\}$) 
if, 
for each $C^k$-map  $f \: M \to G$ from a $C^k$-manifold $M$ 
to $G$ with $\im(f) \subeq \iota(H)$, the corresponding map 
$\iota^{-1} \circ f \: M \to H$ is $C^k$. 

\index{Lie subgroup!strictly initial}
We then call $H$ a {\it $C^k$-initial Lie subgroup}, 
and for $k = \infty$, we simply 
call it {\it initial}. 
We say that $H$ is {\it strictly initial} if, in addition, 
$\L(\iota)$ is injective. Recall that this condition is automatically
satisfied for $k < \infty$. 

In the following, we shall mostly identify $H$ with 
the subgroup $\iota(H)$ of~$G$. In this sense we shall think of the injective smooth 
map $\iota$ as 
a Lie group structure on $H$ and write $H^L$ for the subgroup $H$ of $G$, 
endowed with this Lie group structure. In view of Lemma~\ref{lem:init-unique} below, 
this structure is unique. 
\end{defn} 


\begin{defn}
The following lemma shows that the existence of an initial Lie group structure 
only depends on the subgroup $H$, considered as a subset of $G$. 
It therefore makes sense to call a subgroup $H \subeq G$ {\it initial} if 
\index{Lie subgroup!initial}
it carries an initial subgroup structure. 
\end{defn}

\begin{lem} \mlabel{lem:init-unique} Any subgroup $H$ of a Lie group $G$ carries 
at most one structure of an initial Lie subgroup.
\end{lem}

\begin{prf} If $\iota' \: H' \into G$ is an 
initial Lie subgroup with the same range as $\iota \: H \to G$,  then 
$\iota^{-1} \circ \iota' \: H' \to H$ and 
$\iota'^{-1} \circ \iota \: H \to H'$ are smooth 
morphisms of Lie groups, so that $H$ and $H'$ are isomorphic. 
\end{prf}

\begin{ex} (A non-initial subgroup; \cite[Thm.~2.3]{Gl08a}) 
We consider the locally convex space $G := \R^\N$, endowed with 
the product topology, as an abelian Lie group. 
We claim that the subgroup $H = \ell^\infty$ of bounded sequences, 
carries no initial Lie subgroup structure. 

Let $\iota \: N \to G$ be a $C^\infty$-initial Lie subgroup with $\iota(N) = H$. 
We endow $H$ with the subspace topology, inherited from $G$, which turns 
it into a locally convex space, hence into a Lie group 
for which the inclusion map $j \: H \to G$ is smooth and a topological embedding. 
By the universal property of $\iota$, the map 
$k := \iota^{-1} \circ j \: H \to N$ is a morphism of Lie groups. 
Next we observe that the corestriction $\zeta := \iota|^H \: N \to H$ is continuous,
so that $k = \zeta^{-1}$ implies that $\zeta$ is a homeomorphism, resp., 
that $\iota$ is a topological embedding. 

We may thus identify $H$ and $N$ as topological groups. 
The subgroup $H = \bigcup_{n \in \N} [-n,n]^\N$ is $\sigma$-compact 
because the cubes $[-n,n]^\N$ are compact by Tychonoff's Theorem. 
As a subspace of a metric space, $H$ is also regular. 
For the smooth manifold $N$, the existence of a chart around $0$  
thus implies that its tangent space 
$T_0(N)$ contains a $\sigma$-compact $0$-neighborhood, 
hence that the locally convex space $T_0(N)$ is $\sigma$-compact. 

Now we show that the tangent map $T_0(\iota) \:  T_0(N) \to \R^\N$ is surjective. 
Pick $y = (y_n)_{n \in \N} \in \R^\N$ and consider the smooth curve 
\[  g \: \R \to \R^\N, \quad g(t) = (\sin(t y_n))_{n \in \N}.\] 
As $g$ is componentwise smooth, it is a smooth curve in $\R^\N$. 
As its range is contained in $H$, there exists a smooth curve
$h \: \R \to N$ with $g = \iota \circ h$, so that 
$y = g'(0) \in T_0(\iota)$. 
The surjectivity of $T_0(\iota)$ now implies 
that $\R^\N$ is $\sigma$-compact, so that Baire's Theorem, applied 
to the complete metric space $\R^\N$, leads to the assertion that 
$\R^\N$ is locally compact, which is absurd. 
\end{ex}

\begin{lem} \mlabel{lem:reg-subgrp} 
Let $H$ be a regular Lie group and $j \: H \to G$ 
be an injective morphism of Lie groups for which $\L(j)$ is a topological 
embedding with closed range. Suppose that, for every smooth curve 
$\gamma \: [a,b] \to G$ with $\gamma([a,b]) \subeq j(H)$, 
the curve $\gamma^H := j^{-1} \circ \gamma$ is smooth. 
Then $j$ defines an initial Lie subgroup of $G$. 
\end{lem}

\begin{prf} Let $M$ be a smooth manifold and $\phi \in C^\infty(M,G)$ with 
$\im(\phi) \subeq j(H)$. We have to show that $j^{-1} \circ \phi \: M 
\to H$ is smooth. Since this assertion is local, it suffices to assume that 
$U := M$ is an open circular $0$-neighborhood in a locally convex space~$E$ 
and that $\phi(0) = \be$. We write $\fh := \im(\L(j))$ for the Mackey complete 
subalgebra of $\g$ obtained by embedding $\L(H)$. 

Our assumption implies that the $\g$-valued $1$-form 
$\delta(\phi)$ takes values in $\fh$ and since $\fh$ is closed, 
it defines an element of $\Omega^1(U,\fh)$. 
As it satisfies the Maurer--Cartan equation, $U$ is simply connected, 
and $H$ is regular, there exists a smooth map 
$\phi^H \: U \to H$ with $\delta(\phi^H) = \delta(\phi)$ and 
$\phi^H(0) = \be$ (Theorem~\ref{thm-fundamental}). 
The Uniqueness Lemma~\ref{lem:c.12b} now implies that 
$j \circ \phi^H = \phi$, so that $j^{-1} \circ \phi = \phi^H$ is smooth.
\end{prf}

\begin{lem} If 
$\iota_H \: H \to G$ is a $C^k$-initial Lie subgroup for some 
$k \in \N\cup \{ \infty\}$, then the following assertions hold: 
\begin{enumerate}
\item[\rm(i)] For $k =1$, we have $\L^d(H) = \im(\L(\iota_H))$. 
\item[\rm(ii)] If $H$ is normal, then  the conjugation action of $G$ on $H$ is~$C^k$.
\end{enumerate}
\end{lem}

\begin{prf} (i) If $H$ is $C^1$-initial, then 
$$C^1_*([0,1],H) = \{ \alpha \in C^1_*([0,1],G) \: \im(\alpha) \subeq H\}$$ 
implies that $\L(\iota_H)(\L(H)) = \L^d(H)$. 

(ii) The conjugation map $G \times H \to G, (g,h) \mapsto ghg^{-1}$, is a smooth 
map whose values lie in $H$, hence it is smooth as a map $G \times H \to H$. 
\end{prf}

\begin{lem} \mlabel{lem:liealg-cen} $\L^d(Z(G)) \subeq \z(\L(G))$. 
\end{lem}

\begin{prf} Each $\gamma \in C^1_*([0,1],G)$ with $\im(\gamma) \subeq Z(G)$ 
satisfies $\Ad(\gamma(t)) = \id_{\L(G)}$ for each $t$, and hence 
$$ 0 = \derat0  \Ad(\gamma(t))x = [\gamma'(0), x] $$
for each $x \in \L(G)$ (Proposition~\ref{prop:der-Ad}), thus 
$\gamma'(0) \in \z(\L(G))$. 
\end{prf}

\begin{prop} \mlabel{prop:initial-disc1} 
Let $G$ be a Lie group and $H \subeq G$ be a subgroup with 
the property that all $C^1$-arcs in $H$ are constant. Then 
the discrete topology turns $H$ into a $C^k$-initial Lie subgroup of~$G$
for any $k \in \N \cup \{\infty\}$. 
\end{prop} 

\begin{prf} Write $H_d$ for the group $H$, endowed with the discrete topology. 
If $M$ is a smooth manifold and $f \: M \to G$ a $C^k$-map with 
$f(M) \subeq H$, then $f$ is locally constant because each $C^1$-arc in $H$ 
is constant. 
Hence the map $f \: M \to H_d$ is locally constant and therefore smooth. 
\end{prf} 

In Corollary~\ref{cor:c1-disc} we shall see that, 
in a closed subgroup $H$ of a locally exponential Lie group~$G$, 
all $C^1$-arcs are constant if all one-parameter subgroups of $H$ 
are trivial. This leads to the concrete example of the subgroup 
$H := L^2([0,1],\Z)$ of the Banach--Lie group 
$G := L^2([0,1],\R)$ in which all $C^1$-curves are constant 
(Example~\ref{ex:L2-z}).

\begin{prop}\mlabel{prop:initial-disc2}  
Let $G$ be a Lie group with exponential function for which the subgroup 
\begin{equation}
  \label{eq:gamma-z}
\Gamma_Z := \{ x \in \z(\g) \: \exp_G x = \be\} \subeq \z(\g)
\end{equation}
is discrete and $\z(\g)$ is Mackey complete. 
Then $Z(G)$ is an initial Lie subgroup for which 
$\exp\res_{\z(\g)}$ induces an isomorphism $\z(\g)/\Gamma_Z \to Z(G)_0$. 
\end{prop}

\begin{prf}  First we observe that the exponential function 
$$ \exp_Z := \exp_G\res_{\z(\g)} \: \z(\g) \to Z(G) $$
is a group homomorphism and our assumption implies that 
$Z := \z(\g)/\Gamma_Z$ carries a natural Lie group 
structure (Example~\ref{exs:ab-quot}). The smooth function 
$\exp_Z \: \z(\g) \to G$ 
factors through an injective morphism of Lie groups 
\break $\iota \: Z \to G$
whose image is the subgroup $\exp_G \z(\g)$ of $Z(G)$ and whose differential 
$\L(\iota)$ is the embedding $\z(\g) \into \g$.
We now endow $Z(G)$ with the Lie group structure for which 
$\exp_G \z(\g)$ is an open subgroup isomorphic to $Z$ and write 
$Z(G)^L$ for this Lie group (Corollary~\ref{cor:open-lie}). 
Then the inclusion map 
$$ \iota \: Z(G)^L \to G $$
is an injective morphism of Lie groups with injective differential. 

To see that this defines an initial Lie group structure on $Z(G)$, let 
\break $f \: M \to G$ be a smooth map with $f(M) \subeq Z(G)$ 
and fix a point $m_0 \in M$. To show that $f \: M \to Z(G)^L$ 
is smooth in a neighborhood of $m_0$, 
we may w.l.o.g.\ assume that $M$ is an open convex subset of a locally convex space 
and $f(m_0)=\be$. 
Then the left logarithmic derivative 
$\delta(f) \in \Omega^1(M,\z(\g))$ 
is a smooth $\z(\g)$-valued $1$-form on $M$ 
(Lemma~\ref{lem:liealg-cen}). 
Since $\delta(f)$ satisfies the Maurer--Cartan equation 
$$ \dd(\delta(f)) + \frac{1}{2} [\delta(f), \delta(f)] = 0, $$
(Lemma~\ref{lem:MC}), the fact that $\z(\g)$ is abelian implies that 
$\delta(f)$ is closed. 
As $\z(\g)$ is Mackey complete and 
$M$ is simply connected, there exists a smooth function 
$f_1 \: M \to \z(\g)$
with $\dd f_1 = \delta(f)$ and $\exp_Z(f_1(m_0))= f(m_0)$ 
(Fundamental Theorem~\ref{thm-fundamental}). 
Then the smooth function 
$$ f_2 := \exp_Z \circ f_1 \: M  \to Z(G)^L  $$
satisfies $\delta(f_2) = \dd f_1 = \delta(f)$ and $f_2(m_0) = f(m_0),$
so that $f_2 = f$, and this implies that $f \: M \to Z(G)^L$ is smooth. 
\end{prf}

The following example shows that, in general, the center 
$Z(N)$, considered as a topological subgroup, is not a Lie group because 
it may be a totally disconnected non-discrete subgroup.  
Since all Lie groups modeled on Mackey complete spaces 
known to the authors 
satisfy the requirements from Proposition~\ref{prop:initial-disc2},  
for all these Lie groups $G$, the center 
$Z(G)$ is an initial Lie subgroup. 

\begin{ex} \mlabel{ex:nonlie-cen}
We consider the Lie group 
$$ G := \C^\N \rtimes \R^\N \quad \hbox{ with } \quad 
((z_n), (t_n)) * ((z_n'), (t_n')) := ((z_n + e^{t_n}z_n'), (t_n + t_n')), $$
where the locally convex structure on $\C^\N$ and $\R^\N$ is given by the 
product topology. 
We observe that $G$ is an infinite topological product of groups isomorphic to 
$\C \rtimes \R$, endowed with the multiplication 
$$ (z,t) * (z',t') = (z + e^t z', t + t'). $$
Therefore the center of $G$ is 
$$ Z(G) = \{0\} \times (2\pi i \Z)^\N \cong \Z^\N. $$
This group is totally disconnected, but not discrete, 
hence not a Lie group with respect to the subspace topology. 
Therefore the center of a Lie group, although closed, is not always a 
Lie group. 

Since $Z(G)$ is totally disconnected, $Z(G)$ is initial with respect to 
the discrete topology. 
\end{ex}

\begin{rem} If $G$ has an exponential function which is injective 
in some $0$-neighborhood, then $0$ is isolated in $\Gamma_Z$ from 
\eqref{eq:gamma-z}, hence 
$\Gamma_Z$ is discrete. Therefore the preceding proposition applies 
in particular to locally exponential Lie groups and in particular 
to Banach--Lie groups. 
\end{rem}

\begin{probl} Find a Lie group $G$ with exponential function 
for which the closed subgroup $\Gamma_Z = \ker(\exp\res_{\z(\g)}) 
\subeq \z(\g)$ is not discrete. 
\end{probl}

\subsection{Lie subgroups} 

Throughout this book, we shall encounter several classes of ``nice'' subgroups 
of Lie groups. One of these classes consists of those subgroups $H$ which 
are submanifolds of $G$: 

\begin{defn} \mlabel{defn:liesubgroup} A subgroup $H$ of the Lie group $G$ 
with Lie algebra $\g$ is called a 
{\it Lie subgroup} 
\index{Lie subgroup} 
if $H$ is a submanifold of $G$, i.e., there exists a closed subspace 
$\h \subeq \g$ and a $\g$-chart $(\phi,U)$ of 
$G$ with $\be \in U$ and $\phi(U \cap H) = \phi(U) \cap \h$. 
\end{defn} 

Since every submanifold of a manifold carries a natural manifold structure, 
turning it into an initial submanifold (Definition~\ref{definitmfd}), we derive 
immediately: 

\begin{prop} Each Lie subgroup $H$ of a Lie group $G$ is a Lie group with respect 
to its natural manifold structure. Moreover, it is $C^1$-initial. 
\end{prop}

\begin{prf} The multiplication map $m_H \: H \times H \to H$ is smooth because 
$H$ is an initial submanifold of $G$ and $m_H$ is the restriction of the smooth map 
$m_G \: G \times G \to G$ to $H \times H$. 
We likewise see that the inversion $\eta_H$ is smooth. 
Therefore $H$ is a Lie group 
with respect to its natural manifold structure. 

To describe this manifold structure, 
pick  $\h \subeq \g$ and a $\g$-chart $(\phi,U)$ of 
$G$ with $\be \in U$ with$\phi(U \cap H) = \phi(U) \cap \h$. 
Then the charts of $H$ are given by the maps 
$$\psi_h \: h(U \cap H) = (hU) \cap H \to \phi(U) \cap \h, \quad 
g \mapsto \phi(h^{-1}g). $$

Next we show that $H$ is $C^1$-initial. Let 
$f \: M \to G$ be a $C^1$-map with values in $H$. 
To show that the corestriction $f \: M \to H$ is $C^1$, 
it suffices to consider small neighborhoods 
of some element $m_0 \in M$, and after replacing $f$ by the function 
$m \mapsto f(m_0)^{-1} f(m)$, we may w.l.o.g.\ assume that 
$f(M) \subeq U$ (where $U$ is as in the definition above). 
Then $\hat f:= \phi \circ f \: M \to \L(G)$ 
is $C^1$ and $\h$-valued, so that the closedness of $\h$ implies that its 
corestriction to $\h$ is $C^1$. 
Therefore $f = \psi_\be^{-1} \circ \hat f \: M \to H$ is a $C^1$-map. 
\end{prf}

\begin{rem} (a) Let $H \leq G$ be a Lie subgroup and 
$(\phi,U)$ be a chart of $G$ as in Definition~\ref{defn:liesubgroup}. 
Then the local multiplication $x*y := \phi(\phi^{-1}(x)\phi^{-1}(y))$ on 
$$D := \{(x,y) \in \phi(U) \times \phi(U) \: \phi^{-1}(x)\phi^{-1}(y) \in U\}$$ 
satisfies 
\begin{eqnarray}
  \label{eq:3.2.x}
x*y \in \h \quad \hbox{ for } \quad (x,y) \in D \cap (\h \times \h) 
\end{eqnarray}
and 
\begin{eqnarray}
  \label{eq:3.2.y}
x^{-1} \in \h \quad \hbox{ for } \quad x \in \h \cap \phi(U). 
\end{eqnarray}
In view of Remark~\ref{rem:brack-taylor}, 
this implies that $\h$ is a closed Lie subalgebra of $\g$, isomorphic to 
$\L(H)$. 

If, conversely, $\h \leq \L(G)$ is a closed Lie subalgebra for which there 
is a chart $(\phi,U)$, satisfying (\ref{eq:3.2.x}) and (\ref{eq:3.2.y}), then 
we apply Theorem~\ref{thm:locglob} to  
the embedding $\phi^{-1} \: \phi(U) \cap \h \to G$. This 
leads to a Lie group 
structure on the subgroup $H := \la \phi^{-1}(\phi(U) \cap \h)\ra$ of $G$, generated by 
$\phi^{-1}(\phi(U) \cap \h)$. 
As we shall see in Example~\ref{ex:dense-wind} below, this does 
not always lead to a submanifold of~$G$. 

(b) A weaker concept of a ``Lie subgroup'' is obtained by requiring only that 
$H \leq G$ is a subgroup, for  which there exists an identity neighborhood $U^H$ 
whose $C^\infty$-arc-component $U^H_0$ containing $\be$ is a submanifold of $G$ 
(cf.\ \cite[p.~45]{KYM85}). Then we can use Theorem~\ref{thm:locglob} to obtain 
a Lie group structure on $H$ for which some identity neighborhood is diffeomorphic to 
an identity neighborhood in $U^H_0$. 
\end{rem}

\subsection{Quotients of Lie groups} 

In this subsection we collect some general material on Lie group structures on
(normal) subgroups and quotient groups. 
In particular, we give a condition on a normal subgroup
$N \trile G$ which ensures that the quotient group $G/N$ is a manifold, 
for which the quotient map $q \: G \to G/N$ defines on $G$
the structure of a smooth $N$-principal bundle. 

\begin{defn} (Split Lie subgroups) \mlabel{defc.4} 
Let $G$ be a Lie group. A subgroup $H$ of $G$
is called a {\it split Lie subgroup} \index{Lie subgroup!split} 
if it carries a Lie group
structure for which the inclusion $i \: H \into G$ is an injective morphism 
of Lie groups and there exists an open subset $U$ of some locally convex space $E$ 
and a smooth map $\sigma \: U \to G$ with $\sigma(0) = \be$ such that the map 
$$ \mu \: U \times H \to G, \quad (x,h) \mapsto \sigma(x) i(h) $$
is a diffeomorphism onto an open subset of $G$. Then the set $\mu(U\times H)= \sigma(U)H$ is a 
tube around the subgroup $H$ of $G$ and $\mu$ provides
 product coordinates on this 
tube. 
\end{defn}

\begin{lem} Each split Lie subgroup $H$ of the Lie group $G$ is in particular 
a submanifold, hence a Lie subgroup. 
\end{lem}

\begin{prf} This is immediate from the fact that, for each 
$x \in U$ (as in the definition), the set $\{x\} \times H$ is a submanifold of 
$U \times H$, hence that $H$ is a submanifold of $G$. 
\end{prf}

\begin{prop} \mlabel{prop:split-Lie} For a Lie group 
structure on a subgroup $H$ of the Lie group $G$,  
the  following are equivalent: 
\begin{enumerate}
\item[\rm(i)] $H$ is a split Lie subgroup. 
\item[\rm(ii)] The canonical right action of 
$H$ on  $G$ defined by $(g,h) \mapsto gh$ 
defines a smooth $H$-principal bundle for a suitable manifold structure 
on the quotient space $G/H$. 
\end{enumerate}
\end{prop}

\begin{prf} (i) $\Rarrow$ (ii): Our assumption implies in particular 
that $i$ is injective and that 
$H \cong i(H)$ is a closed subgroup of $G$. We therefore identify $H$ with the 
image of~$i$. Let $M := G/H$ denote the topological quotient space of $G$ by 
the closed subgroup $H$ and observe that the quotient topology turns $M$ into a 
Hausdorff space (Exercise~\ref{exer:6.1.1}). 
We have to show that $M$ is a smooth manifold and that the quotient map 
$q \: G \to M, g \mapsto gH$ defines on $G$ the structure of a smooth $H$-principal bundle. 

Let $U$ be an open subset of some locally convex space $E$ 
and $\sigma \: U \to G$ be a smooth map with $\sigma(0) = \be$ such that the map 
$$ \mu \: U \times H \to G, \quad (x,h) \mapsto \sigma(x) i(h) $$
is a diffeomorphism onto an open subset of $G$. 
Let $W := \mu(U \times H) = \sigma(U)H$ and write 
$p_U \: W \to U$ for the smooth map with 
$w  \in \sigma(p_U(w))H$ for all $w \in W$. 
Then $q(W) = q\circ \sigma(U)$ is an open subset of $G/H$ and 
$q(W) = W/H$ is homeomorphic to $U$ because $\mu$ is a homeomorphism. 
Therefore the map 
$q \circ \sigma \: U \to q(W)$ is a homeomorphism. We put 
\[ \psi := (q \circ \sigma)^{-1} \: q(W) \to U.\]
Then, for $g \in W$, we have 
$$\psi(gH) = \psi(\sigma(p_U(g)) H) = p_U(g). $$
Hence, for each $g \in G$, the map 
\[ \psi_g \: V_g := gq(W) \to U, \quad gq(\sigma(x)) \mapsto x \] 
also is a homeomorphism. We claim that this 
collection of homeomorphisms is a smooth atlas of $G/H$. 

Let $g, h \in G$ and note that $V_g \cap V_h = g(V_\be \cap V_{g^{-1}h}).$
For 
$$x \in U' := \psi(V_\be \cap V_{g^{-1}h}) = \{ x \in U \: 
h^{-1} g \sigma(x) \in W \}, $$
we then have 
$$ \psi_h \circ \psi_g^{-1}(x) 
=  \psi(h^{-1}q(g \sigma(x)))
=  \psi(h^{-1}g \sigma(x)H)
=  p_U(h^{-1}g \sigma(x)), $$
showing that $\psi_h\circ \psi_g^{-1} \res_{U'} \: U' \to U$ is a smooth map. 
This proves that the maps $(\psi_g, V_g)_{g \in G}$ form a smooth atlas of $G/H$.

On $V_g = gq(W)$, the map 
$$g \sigma(x)H \mapsto 
\sigma(\psi_g(g \sigma(x)H)) 
= \sigma(\psi(\sigma(x)H))= \sigma(p_U(\sigma(x))) = \sigma(x)$$ 
is a smooth section of the quotient map $q \: G \to M = G/H$. Moreover, the map 
$$ \psi_g \circ q \: g W \mapsto U, \quad gw \mapsto 
\psi_g(gwH) = p_U(w) $$ 
is smooth, which shows that $q$ is smooth. We conclude 
that $q$ defines a smooth $H$-principal bundle over $M$, hence 
that $H$ is a split Lie subgroup. 

(ii) $\Rarrow$ (i): Let $(\phi,U)$ be an $E$-chart of the manifold $G/H$ and assume 
that $\sigma \: U \to G$ is a smooth section of the quotient map $q \: G \to G/H$. 
Then (ii) implies that the map 
$$ \Phi \: U \times H \to G, \quad (x,h) \mapsto \sigma(x)h $$
is a diffeomorphism onto the open subset $\sigma(U)H$ of~$G$. 
\end{prf} 

\begin{ex}
Since the Lie algebra $\h$ of a Lie subgroup $H$ of a Lie group 
$G$ need not have a closed complement in 
$\g$, not every Lie subgroup is split. A simple example is the
closed subspace $H := c_0(\N,\R)$ in the Banach space 
$G := \ell^\infty(\N,\C)$ (\cite[Satz IV.6.5]{We95}). 
\end{ex}

\begin{cor} \mlabel{cor:quot-splitnorm} 
Let $G$ be a Lie group and $N \trile G$ a split normal 
Lie subgroup. Then the quotient group $G/N$ carries a natural Lie group 
structure for which the quotient map 
$q \: G \to G/N$ defines on $G$ the structure of a principal $N$-bundle. 
\end{cor}

\begin{prf} We endow $G/N$ with the manifold structure from 
Proposition~\ref{prop:split-Lie}, which already provides the $N$-principal 
bundle structure on $G$. Then the product map $q \times q \: G \times G 
\to G/N \times G/N$ defines an $(N \times N)$-principal 
bundle structure on $G \times G$. Since the map 
$$ m_{G/N} \circ (q \times q) = q \circ m_G \: 
G \times G \to  G/N, \quad (g,h) \mapsto ghN $$ 
is smooth, the multiplication map $m_{G/N}$ on 
$G/N$ is smooth. It likewise follows from 
$\eta_{G/N} \circ q = q \circ \eta_G$ that the inversion $\eta_{G/N}$ of $G/N$ 
is smooth, hence that $G/N$ is a Lie group. 
\end{prf}

\begin{prop} \mlabel{prop:split-crit-ban} A Lie subgroup $H$ of the 
Banach--Lie group $G$ is split if and only if 
its Lie algebra $\L(H)$ has a closed complement. 
\end{prop}

\begin{prf} Suppose that $H$ is split and that 
$E \subeq \L(G)$ is a closed complement of 
$\L(H)$ in $\L(G)$. Then it follows immediately from the definition that 
the tangent map $T_{(0,0)}(\mu) \: E \times \L(H) \to \L(G)$ is a linear topological 
isomorphism, hence that $\L(H)$ is split. 

Suppose, conversely, that there exists a closed subspace $E \subeq \L(G)$ for which 
the summation map $E \times \L(H) \to \L(G), (x,y) \mapsto x + y$ is a homeomorphism. 
Then we consider the smooth map 
$$ \mu \: E \times H \to G, \quad (x,y) \mapsto \exp_G(x) h, $$
for which $T_{(0,\be)}(\mu)(x,y)= x + y$ is a linear isomorphism. Hence the 
Inverse Function Theorem implies the existence of an open $0$-neighborhood 
$U_E \subeq E$ and an open $\be$-neighborhood $U_H \subeq H$ for which 
the restriction 
$\mu\res_{U_E\times U_H}$ is a diffeomorphism onto an open subset of $G$. 

From $\mu(x,yh) = \mu(x,y)h$ for $h \in H$, we now derive that 
$T_{(x,y)}(\mu)$ is invertible for each pair $(x,y) \in U_E \times H$, so that 
$\mu$ defines a local diffeomorphism $U_E \times H \to G$. 

Next we choose an open $\be$-neighborhood $U_G$ in $G$ with 
$U_G \cap H = U_H$, and an open $0$-neighborhood $V_E \subeq U_E$ such that 
$\exp_G(-V_E) \exp_G(V_E) \subeq U_G$. For two pairs 
$(x,h), (x',h') \in V_E \times H$ mapped to the same element of $G$ under $\mu$, 
we then have 
$$ \exp_G(-x')\exp_G(x) = h'h^{-1} \in U_G \cap H = U_H. $$
This leads to $\mu(x,\be) = \mu(x',h'h^{-1})$, so that the injectivity of $\mu$ on 
$U_E \times U_H$ implies $x = x'$ and $h' = h$. Therefore 
$\mu$ is injective on $V_E \times H$, hence a diffeomorphism 
onto an open subset of $G$. 
\end{prf}

\section{Locally exponential integral subgroups} \mlabel{sec:6.2} 

In this section we give a criterion for the integrability 
of a continuous injection $\h \to \L(G)$ of a locally exponential Lie algebra $\h$ 
into the Lie algebra $\L(G)$ of a Lie group $G$ with an 
exponential function. 

\subsection{Integral subgroups} 

\begin{defn} \mlabel{defn:5.5.1} 
Let $G$ be a Lie group. An 
\index{integral subgroup, of a Lie group} 
{\it integral subgroup} 
is an injective morphism $\iota \: H \to G$ of 
Lie groups for which $H$ is connected and the differential $\L(\iota) \: \L(H) \to \L(G)$ 
is injective. 
\end{defn} 

We think of these subgroups as arising by ``integration'' of an injective 
morphism of Lie algebras, hence the terminology. 

Integral subgroups are submanifolds in a rather weak sense. 
We shall see below that, in general, they are not initial 
submanifolds (Example~\ref{ex:non-separ}), which leads to the 
problem to determine when an integral subgroup actually is initial.  
In Theorem~\ref{thm:sep-subgroup2} we show that this is the case for a large class 
of integral subgroups. 

\begin{lem}\mlabel{lem:integ-subgrp} Let $\iota \: H \to G$ be a locally exponential 
integral subgroup, $\h := \im(\L(\iota))$, $\g =\L(G)$,  
and assume that $G$ has an exponential function. 
Then the following assertions hold: 
\begin{enumerate}
  \item[\rm(i)] $\h$ is a stable subalgebra of $\g$, i.e., 
$e^{\ad x}\h = \h$ for $x \in \h$. 
  \item[\rm(ii)] $\iota(H) = \la \exp_G(\h) \ra$. 
\end{enumerate}
\end{lem}

\begin{prf} (i) Since $G$ and $H$ have exponential functions, 
for each $x \in \L(H)$, the operators $\ad_{\L(H)} x$ and $\ad_\g(\L(\iota)x)$ satisfy 
$$ \L(\iota) \circ \ad_{\L(H)} x = \ad_\g(\L(\iota)x) \circ \L(\iota) $$
and generate smooth one-parameter groups of automorphisms of $\h$, resp., $\g$. 
Then 
$$ \L(\iota) \circ e^{\ad_{\L(H)} x}  = e^{\ad_\g(\L(\iota)x)} \circ \L(\iota) $$
follows from (the proof of) Lemma~\ref{lem4.1.9} in the appendix. 
From that we immediately obtain (i). 

If $G$ is locally exponential, we can also argue with Lemma~\ref{lem:5.3.4}. 

(ii) follows from $H = \la \exp_H \L(H) \ra$ and 
$\exp_G \circ \L(\iota)= \iota \circ \exp_H$. 
\end{prf} 

In view of the preceding lemma, we may consider locally exponential 
integral subgroups as obtained by endowing the subgroup 
$H := \la \exp_G \h \ra$ of $G$, where $\h \into \L(G)$ is an inclusion of a locally exponential 
Lie algebra, with some Lie group structure for which $\L(H) = \h$. 

\begin{lem} \mlabel{intsubgrp-discrete} Let 
$G$ be a group with an exponential function and Lie algebra $\g$, 
let $\h$ be a locally exponential Lie algebra and 
$\alpha \: \h \to \g$ a continuous injective homomorphism of Lie algebras. 
Then the following are equivalent: 
\begin{enumerate}
\item[\rm(i)] $0$ is isolated in the subset $(\exp_G \circ \alpha)^{-1}(\be)$ of $\h$. 
\item[\rm(ii)] $\Gamma_\alpha := \{ x \in \z(\h) \: \exp_G(\alpha(x)) = \be\}$ is a 
discrete subgroup of $\z(\h)$. 
\end{enumerate} 
\end{lem}

\begin{prf} (i) $\Rarrow$ (ii): Condition (i) implies that $0$ is isolated in 
$\Gamma_\alpha$, but this is equivalent to $\Gamma_\alpha$ 
being discrete (Exercise~\ref{exer:3.3.8}). 

(ii) $\Rarrow$ (i): 
Let $H_{\rm loc} \subeq \h$ be an exponential local Lie group 
with the property that 
$$ \{ x \in H_{\rm loc} \: e^{\ad x} = \be\} \subeq \z(\h) $$
(cf.\ the Adjoint Enlargeability Theorem~\ref{thm:5.3.8}). 
If $\exp_G(\alpha(x))= \be$ 
holds for $x \in H_{\rm loc}$, then 
$$ e^{\ad_\g \alpha(x)} = \Ad(\exp_G(\alpha(x))) = \id_\g $$
leads to 
$$ \alpha \circ e^{\ad_\h x} = e^{\ad_\g \alpha(x)} \circ \alpha = \alpha, $$
so that the injectivity of $\alpha$ yields $e^{\ad_\h x} = \be$. We therefore 
have 
$$ (\exp_G \circ \alpha)^{-1}(\be) \cap H_{\rm loc} \subeq \z(\h), $$
so that (ii) implies (i). 
\end{prf}

\begin{thm} [Integral Subgroup Theorem] \mlabel{thm:5.5.3} 
Let $G$ be a Lie group with an exponential function, Lie algebra 
$\g$ and $\alpha \: \h \to \g$ an injective continuous morphism of locally convex Lie algebras, 
where $\h$ is locally exponential. 
If the subgroup 
$$\Gamma_\alpha = \{ z \in \z(\h) \: \exp_G(\alpha(z)) = \be\} $$ 
of $\z(\h)$ is discrete, then 
there exists a locally exponential integral subgroup $\iota \: H \to G$ with 
$\L(H) = \h$ and $\L(\iota) = \alpha$. 

The subgroup $\Gamma_\alpha$ is discrete if
\begin{enumerate}
\item[\rm(a)] $G$ is locally exponential or  
\item[\rm(b)] $\dim \z(\h) < \infty$. 
\end{enumerate}
\end{thm}

\begin{prf} We consider the subgroup 
$H := \la \exp_G \h \ra  \subeq G.$
Let 
$H_{\rm loc} \subeq \h$ be an exponential local Lie group with Lie algebra $\h$ 
for which 
\begin{eqnarray}
  \label{eq:loc-cond}
(\exp_G \circ \alpha)^{-1}(\be) \cap H_{\rm loc} = \{0\} 
\end{eqnarray}
(Lemma~\ref{intsubgrp-discrete}). 
Then Theorem~\ref{thm:5.2.11b} implies the existence of an open balanced 
$0$-neighborhood $V$ in $H_{\rm loc}$ 
with $V \times V \subeq D_{H_{\rm loc}}$ for which the map 
$$ \psi \: H_{\rm loc} \to G, \quad x \mapsto \exp_G(\alpha(x)) $$
satisfies 
\begin{eqnarray}
  \label{eq:mulrel}
\psi(x*y) = \psi(x) \psi(y) \quad \mbox{ for all } \quad x,y \in V. 
\end{eqnarray}
This implies in particular $\psi(x)\psi(-x) = \be$, so that $\psi(-x) = \psi(x)^{-1}$. 

To see that $\psi$ is injective, suppose that 
$\psi(x) = \psi(y)$ for $x,y \in V$. Then 
$$\psi((-x)*y) = \psi(-x)\psi(y) = \psi(x)^{-1}\psi(y) = \be $$ 
and (\ref{eq:loc-cond}) lead to 
$(-x)*y \in \psi^{-1}(\be) \cap H_{\rm loc} = \{0\}$, 
so that the injectivity of $\alpha$ yields $(-x) * y = 0$. 
We conclude that 
$$x = x * 0 = x * ((-x) * y) = (x * (-x)) * y = 0 * y = y $$
and hence that $\psi$ is injective. 

This means that $V \cong \psi(V) \subeq H$ carries a natural manifold structure 
and that multiplication and inversion are smooth on an identity neighborhood in $V$. 
Since $H$ is generated by $V$, it carries a unique Lie group 
structure  for which $\psi$ defines a chart in a neighborhood of $0$ 
(cf.\ Theorem~\ref{thm:locglob}). 

The curves $\gamma_x \: \R \to H, \gamma_x(t) := \exp_G(\alpha(tx))$ define {one-parameter} groups, and since 
$\gamma_x(t) = \psi(tx)$ holds 
for $t$ close to $0$, these one-parameter groups are 
smooth with $\gamma_x'(0) = x$, 
where we identify $\h$ via $T_0(\psi)$ with $T_\be(H)$. 
Therefore $H$ has an exponential function $\exp_H : \h \to H$ which coincides on 
a $0$-neighborhood with $\psi$. In particular,  $H$ is locally exponential, 
and the inclusion map $\iota \: H \to G$ is smooth in 
an identity neighborhood, hence a morphism of Lie groups. 

This completes the first part of the proof. We now turn to the sufficient conditions 
for $\Gamma_\alpha$ to be discrete. 

(a) If $G$ is locally exponential, then $0$ is isolated in $\exp_G^{-1}(\be)$, 
so that $0$ is isolated in 
$\alpha^{-1}(\exp_G^{-1}(\be))
= (\exp_G \circ \alpha)^{-1}(\be),$ which is condition 
(i) in Lemma~\ref{intsubgrp-discrete}. 

(b) 
Since $\z(\h)$ is abelian, the restriction of $\exp_G \circ \alpha$ to $\z(\h)$ is a smooth 
group homomorphism $\z(\h) \to G$ with kernel $\Gamma_\alpha$ 
(Lemma~\ref{lem:4.1.4}). 
Hence $\Gamma_\alpha$ is a closed subgroup of $\z(\h)$. Since 
$\L(\exp_G\circ \alpha) = \alpha\res_{\z(\h)}$ is injective, 
the derivative of each $C^1$-curve 
$\gamma \: [0,1] \to \Gamma_\alpha$ lies in $\ker \alpha = \{0\}$. 
Hence each $C^1$-curve in $\Gamma_\alpha$ is constant, 
so that $\Gamma_\alpha$ contains no non-trivial vector subspace, hence is discrete 
if $\fz$ is finite-dimensional (Exercise~\ref{exer:3.5.6}). 
\end{prf}

\begin{rem}
The condition, that a closed subalgebra $\h \subeq \g$ is locally exponential,  
is quite subtle. First of all, stability of $\h$ is necessary 
(Lemma~\ref{lem:5.3.4}). 

In view of Theorem~\ref{thm:5.2.11c}, it is equivalent to 
$x*y \in \h$ for $x,y \in \h$ sufficiently close to $0$. 
To verify this condition, one would like 
to show that for $x,y \in \h$ 
the integral curve $\gamma(t) := x * ty$ of the left 
invariant vector field $y_l$ through $x$ does not leave the closed subspace 
$\h$ of $\g$. This leads to the necessary condition 
$T_0(\lambda_x)(\h) \subeq \h$. Under the assumption that $\h$ is stable, 
this means that the operator $\kappa_\g(x) = \int_0^1 e^{-t\ad x}\, dt$ 
satisfies $\kappa_\g(x)(\h) = \h$ for $x \in \h$ sufficiently close to $0$ 
(whenever the integral is defined in ${\cal L}(\g)$, which is always the 
case if $\g$ is Mackey complete). 
We have already seen in Corollary~\ref{cor:5.3.7} 
that, for a closed stable ideal $\fn \trile \g$
of a locally exponential Lie algebra $\g$, the corresponding 
condition $\kappa_\g(x)(\fn) = \fn$ for all $x\in\g$ 
sufficiently close $0$ is sufficient 
for $\fn$ to be locally exponential. 
\end{rem}

\begin{cor}[Integral Subgroup Theorem--Banach version] \mlabel{cor:intsub-ban}
Let $G$ be a locally exponential Lie group and 
$\h$ be a Banach--Lie algebra. Then any continuous injective morphism 
$\alpha \: \h \to \L(G)$ of Lie algebras integrates to an integral subgroup of~$G$. 
\end{cor}

\begin{prf} In view of Theorem~\ref{thm:5.5.3}, this follows from the fact that all 
Banach--Lie algebras are locally exponential (Example~\ref{exs:5.2.3}). 
\end{prf}

\begin{cor} \mlabel{cor:5.5.6}  
If $\h$ is a Banach--Lie algebra for which there exists an injective 
continuous Lie algebra homomorphism into the Lie algebra 
of a Banach--Lie group $G$, then $\h$ is enlargeable. 
\end{cor}

We also record the following formulation of the Integral Subgroup Theorem in 
terms of integrability of morphisms of Lie algebras. 

\begin{cor} \mlabel{cor:5.5.6a}
An injective morphism $\alpha \: \h \to \g$ of locally exponential Lie algebras 
is integrable to a morphism of locally exponential Lie groups if and only if 
$\g$ is enlargeable. 
\end{cor}

\begin{prop} \mlabel{prop:5.5.4} If $\g$ is a BCH-Lie algebra, 
then each closed subalgebra $\h \subeq \g$ is locally exponential. 
\end{prop}

\begin{prf} Let $G \subeq \g$ be an exponential local Lie group on which 
the multiplication $m_G$ is analytic and given by the Hausdorff series. 
Since $\h$ is closed, the explicit form of the Hausdorff series implies that 
$m_G(D_G \cap (\h \times \h)) \subeq \h,$
and hence that $\h$ is BCH and therefore locally exponential. 
\end{prf}

The following corollary shows that a seemingly weaker assumption 
already implies that a locally exponential Lie algebra is enlargeable. 

\begin{cor} \mlabel{cor:5.5.5b} 
If the locally exponential Lie algebra $\g$ is integrable 
to a group with an exponential function $\exp_G \: \g \to G$ 
for which the subgroup 
$\Gamma_\g := \{ x \in \z(\g) \: \exp_G(x) = \be\}$ is discrete, 
then $\g$ is enlargeable. 
\end{cor}

\begin{rem} Suppose that $\g$ is a locally exponential Lie algebra  
and that $G_1$ is a Lie group with Lie algebra $\g$ and an exponential 
function $\exp \: \g \to G_1$ which is injective on some $0$-neighborhood. 
Then Corollary~\ref{cor:5.5.5b} implies the existence of a locally exponential 
Lie group $G$ and an injective morphism of Lie groups 
$\iota \: G \to G_1$ with $\L(\iota) = \id_\g$. 
Although we do not know of any example where $\iota$ is not an isomorphism of 
Lie groups, we need additional assumptions on $G_1$ to ensure that 
$\iota$ is an isomorphism. 
If the locally exponential Lie group $G$ is regular, 
then the same holds for its 
universal covering $\tilde G$ (see Theorem~\ref{thm:reg-ext-prop} below), 
and we obtain with Theorem~\ref{thm3.2.11} a morphism  of Lie groups 
$\phi \: \tilde G_1 \to \tilde G$
with $\L(\phi) = \id_{\L(G)}$. We can also lift $\iota$ to a Lie group morphism 
$\tilde\iota \: \tilde G \to \tilde G_1$ with 
$\L(\tilde\iota) = \id_{\L(G)}$. Then the uniqueness of Lie group morphisms corresponding 
to a given morphism of Lie algebras yields 
$$ \phi \circ \tilde\iota = \id_{\tilde G} \quad \mbox{ and } \quad 
\tilde\iota \circ \phi = \id_{\tilde G_1}, $$
so that $\tilde\iota$ is an isomorphism of Lie groups. Therefore $\iota$ is 
a local isomorphism, hence also surjective, and thus an isomorphism of Lie groups. 
\end{rem}

\begin{ex} Let $\cA$ be a Mackey complete 
continuous inverse algebra and $G := \GL_n(\cA)_0$ 
be the identity component of the unit group of the matrix 
algebra $M_n(\cA)$, which also is a continuous inverse algebra 
(Proposition~\ref{prop:8.1.4}). 
Then $G$ is a BCH--Lie group (Theorem~\ref{thm:IV.1.11}), 
hence in particular locally exponential. 

To obtain a Lie group associated to the Lie algebra 
$$\fsl_n(\cA) := \oline{[\gl_n(\cA), \gl_n(\cA)]}, $$ 
we use the Integral Subgroup Theorem~\ref{thm:5.5.3} to obtain an 
integral subgroup $\iota \: S \to G$ for which 
$\L(S) = \fsl_n(\cA)$ and $\L(\iota)$ is the inclusion map. 

If $\cA$ is commutative, then $\det \: \GL_n(\cA) \to \cA^\times$ 
is a morphism of Lie groups, so that 
$$ \SL_n(\cA) := \ker(\det) $$
is a locally exponential Lie subgroup with Lie algebra 
$\ker \L(\det) = \ker \tr = \fsl_n(\cA)$ 
(Exercise~\ref{exer:5.6.13}). 
In this case $\SL_n(\cA)_0 = \la \exp \fsl_n(\cA) \ra = S$ 
is the integral subgroup from above. 
Therefore the Lie group $S$ is an analog of the identity component 
$\SL_n(\cA)_0$ for non-commutative continuous inverse algebras. 

In general, there is no determinant function on $\GL_n(\cA)$, so that it 
is not clear how $\SL_n(\cA)$ should be defined, but there 
are natural replacements. The non-commutative version 
of the trace is the Lie algebra homomorphism 
$$ \Tr \: \g = \gl_n(\cA) \to \cA/\oline{[\cA,\cA]}, \quad 
(a_{ij}) \mapsto \Big[\sum_{j = 1}^n a_{jj}\Big]. $$ 
Let $q_G \: \tilde G \to G$ denote the universal covering group. 
Then $\Tr$ integrates to a homomorphism of locally exponential Lie groups 
$$ \alpha \: \tilde G \to \cA/\oline{[\cA,\cA]} $$
(Theorem~\ref{thm:int-thm}).  
Now the group $\hat S := \ker \alpha$ 
is a locally exponential Lie subgroup of $\tilde G$ whose Lie algebra is 
\[ \L(\hat S) = \L(\ker \alpha) = \ker \L(\alpha)= \ker \Tr= \fsl_n(\cA).\] 

The universal covering map $q_G \: \tilde G \to G$ clearly maps 
$\hat S_0$ to $S$, so that $\hat S_0$ is a covering group of $S$, 
but $q_G(\hat S)$ need not be connected. This happens if 
$\pi_1(G) \cong \ker q_G \not\subeq \hat S$, i.e., if 
$$ \per(\Tr) \: \pi_1(\GL_n(\cA)) \to \cA/\oline{[\cA, \cA]} $$
is non-zero. 

If $\cA$ is commutative, then $\alpha$ is the unique homomorphic 
lift of the determinant homomorphism $\det \: \GL_n(\cA)_0 \to \cA^\times_0 
\cong \cA/(\ker \exp_\cA)$. In this case 
$$ q_G^{-1}(\SL_n(\cA)) = \alpha^{-1}(\ker \exp_\cA), $$
so that 
$$ \SL_n(\cA) \cap \GL_n(\cA)_0 = q_G\big(\alpha^{-1}(\ker \exp_\cA)). $$
\end{ex}

\subsection{Separability and the Lie algebra of a subgroup}

We now turn to the question if an integral subgroup determines 
its ``Lie algebra'', considered as a subalgebra of~$\L(G)$. 

\begin{thm} \mlabel{thm:sep-subgroup1} Let $G$ be a Lie group with 
Lie algebra $\g$ and an exponential 
function $\exp_G \: \g \to G$ which is injective on some $0$-neighborhood. 
Further, let $\iota_H \: H \into G$ be a locally exponential integral 
subgroup whose Lie algebra is separable. 
Then 
$$ \im(\L(\iota_H)) = \{ x \in \g \: \exp_G(\R x) \subeq \iota_H(H)\}. $$
In particular, the Lie algebra of $H$ is determined by the 
subgroup $i_H(H)$ of~$G$. 
\end{thm}

\begin{prf} Put $\h := \im(\L(\iota_H))$. 
Since the exponential function $\exp_H$ of $H$ satisfies 
$\iota_H \circ \exp_H = \exp_G \circ \L(\iota_H)$, we have 
$\exp_G(\h) \subeq \iota_H(H)$. Let $x \in \g\setminus \h$. 
We have to show that $\exp_G(\R x) \not\subeq \iota_H(H)$. 

Let $U_\g \subeq \g$ be a symmetric open $0$-neighborhood on which 
$\exp_G$ is injective. 
We further choose an exponential local Lie group $U_\h \subeq \L(\iota_H)^{-1}(U_\g)$ in 
$\L(H)$ such that 
$\exp_H\res_{U_\h}$ is a diffeomorphism onto an open subset of $H$ and 
$V_\h = - V_\h \subeq U_\h$ with $\exp_H(V_\h)\exp_H(V_\h) \subeq \exp_H(U_\h)$. 

Pick $\eps > 0$ such that $tx \in U_\g$ holds for $|t| \leq \eps$. 
For any such $t\not=0$ we have $\exp_G(tx) \not \in \exp_G(U_\h)$ 
because $tx \not\in \h$ and $\exp_G\res_{U_\g}$ is injective. 
For $|t|,|s|< \frac{\eps}{2}$ with $t \not=s$, this leads to 
$$ \exp_G((t-s)x) \not\in \exp_G(\L(\iota_H)V_\h)\exp_G(\L(\iota_H)V_\h) $$
because the set on the right is contained in $\exp_G(\L(\iota_H)U_\h)$. 
We conclude that the sets $\exp_G(tx) \exp_G(\L(\iota_H)V_\h)$, $|t| \leq \frac{\eps}{2}$, 
are pairwise disjoint. On the other hand, the assumption that 
$\L(H)$ is separable and $H = \la \exp_H \L(H) \ra$ imply that 
$H$ is separable, so that it cannot contain uncountably many pairwise 
disjoint open subsets. Hence $\exp_G(tx)$ cannot be contained in 
$\iota_H(H)$ for all~$t \in \R$.   
\end{prf}

\begin{ex} \mlabel{ex:non-separ} For non-separable subalgebras $\h$, the conclusion 
of the preceding theorem is no longer valid, 
as the following counterexample shows (\cite[p.~157]{HoM98}): We consider
the abelian Lie group $\g := \ell^1(\R,\R) \times \R$, where 
the group structure is given by addition. We write $(e_r)_{r
\in \R}$ for the canonical topological basis elements of
$\ell^1(\R,\R)$. Then the additive subgroup $\Gamma\subeq \g$ generated by the pairs 
$(e_r, -r)$, $r \in \R$, is  discrete (Exercise~\ref{exer:5.6.7}), so that 
$G := \g/\Gamma$ is an abelian Lie group with universal covering group $\tilde G = \g$ 
and $\pi_1(G) \cong \Gamma$. Now we consider the closed
subalgebra $\h := \ell^1(\R,\R)$ of $\g$. As $\h + \Gamma = \g$, we have 
$\exp_G(\h) = G$. For the Banach--Lie group 
$H := \h/(\Gamma \cap \fh)$, we thus obtain 
a bijective morphism $i_H \: H \to G$ of Banach--Lie groups 
for which $\L(i_H)$ is an embedding of a closed hyperplane. 

To see that $i_H$ is not an isomorphism, we note that the 
element $(0,1) \in \g \setminus \h$ generates a 
one-parameter group $\gamma \: \R \to G,t \mapsto \exp_G(0,t)$ 
which is smooth as a map into $G$, but 
$i_H^{-1} \circ \gamma$ is not continuous as a one-parameter group 
of $H$. 

We thus obtain an example for the following pathologies: 
\begin{itemize}
\item[$\bullet$] A bijective morphisms of Banach--Lie groups 
$H \to G$ which is not an isomorphism, hence in particular not open. 
\item[$\bullet$] An integral subgroup $H \into G$ whose range coincides with $G$, 
but whose Lie algebra $\fh$ is a proper hyperplane in $\g$.
\item[$\bullet$] A continuous one-parameter group 
$\gamma \:\R \to G$ and a surjective continuous homomorphism 
$i_H \: H \to G$ for which there exists no continuous one-parameter 
group $\tilde\gamma \: \R \to H$ with $\gamma = i_H \circ \tilde\gamma$, 
i.e., $\gamma$ does not lift to~$H$.
\end{itemize}

The third can also be studied from a different perspective. The closed subgroup 
\[ C := \{ (t,h) \in \R \times H \: \gamma(t) = i_H(h) \} \subeq \R \times H \] 
admits a continuous bijection $q \: C \to \R, (t,h) \mapsto t$ which is 
not a homeomorphism. It can also be written as 
\begin{align*}
 C 
&= \{ (t,\exp_H(e_t)) \: t \in \R  \} 
= \{ (t, e_t)  + (\Gamma \cap \fh) \: t \in \R  \} \\
&  \subeq \R \times H \cong \R \times \big(\fh/(\Gamma \cap \fh)\big),
\end{align*}
which describes $C$ as the image of the discrete subset 
\[ \{ (t,e_t) \: t \in \R\} \subeq \R \times \ell^1(\R) \]  under 
the map $\id_\R \times \exp_H$. 
\end{ex}

\subsection{Initiality of integral subgroups} 

We now turn to sufficient conditions for locally exponential 
subgroups to be initial, 
even if they are not closed. 

\begin{thm} \mlabel{thm:sep-subgroup2} Let $G$ be a locally exponential Lie group with 
Lie algebra~$\g$. Further, let $\iota_H \: H \into G$ be a locally exponential integral 
subgroup whose Lie algebra $\h$ is a closed separable subalgebra of $\g$. 
If $H$ is $C^{k-1}$-regular, then $\iota_H \: H \to G$ defines a $C^k$-initial Lie subgroup of $G$. 
\end{thm}

\begin{prf} Let $\phi \: M \to G$ be a $C^k$-map with 
$\phi(M) \subeq H$. We have to show that the map 
$\iota_H^{-1} \circ \phi \: M \to H$ is $C^k$. 

Let $m_0 \in M$. It suffices to verify the $C^k$-property in a neighborhood of $m_0$. 
Multiplying with $\phi(m_0)$ on the left, we may w.l.o.g.\ assume that 
$\phi(m_0) = \be$. 
Let $U_\g\subeq \g$ be an open symmetric $0$-neighborhood for which 
$\exp_G\res_{U_\g}$ is a diffeomorphism onto its image 
and $V_\g \subeq U_\g$ an open symmetric $0$-neighborhood with 
$\exp_G(V_\g)\exp_G(V_\g) \subeq \exp_G(U_\g)$, so that the 
local multiplication $x*y := (\exp_G\res_{U_\g})^{-1}(\exp_G(x)\exp_G(y))$ 
in $U_\g$ is defined on $V_\g \times V_\g$. We may further assume that 
$V_\h := V_\g \cap \h$ satisfies 
$V_\h * V_\h \subeq \h$. 
Then $\phi^{-1}(\exp_G(U_\g))$ is an open neighborhood of $m_0$, so that we may 
further assume that $\phi(M) \subeq \exp_G(U_\g)$ and that 
$M$ is connected. 
Then $\phi_\g := (\exp_G\res_{U_\g})^{-1} \circ \phi \: M \to \g$ 
is a $C^k$-map. If we can show that $\phi_\g(M) \subeq \h$, then 
$\phi_\g$ is $C^k$ as a map into $\h$, and the relation 
$$ \iota_H^{-1} \circ \phi = \exp_H \circ \phi_\g $$
implies that $\iota_H^{-1} \circ \phi$ is $C^k$. 

Since $M$ is connected, it is connected by smooth arcs, so that 
we may w.l.o.g.\ assume that $M = I = [0,1]$ with $m_0 = 0$. 
Next we show that $\gamma := \phi_\g \: I \to \g$ satisfies 
$\delta(\gamma)_t \in \h$ for each $t\in I$.  
Suppose that this is not the case for some $t_0 \in I$. 
Replacing $\gamma(t)$ by $\gamma(t_0)^{-1}\gamma(t)$, 
we may thus assume that $\gamma(t_0) = 0$ and $\delta(\gamma)_{t_0}\not\in \h$. 

With the Hahn--Banach Theorem, we find a continuous linear functional 
$\lambda \in \g'$ vanishing on $\h$ and satisfying $\lambda(\gamma'(t_0)) = 1$.  
We claim that 
\[ (-\gamma(t)) * \gamma(s) \not\in \fh \] 
for 
$t\not=s$ in some $\eps$-neighborhood of $t_0$. 
We consider the smooth map 
$$ F \: I^2 \to \R, \quad 
F(s,t) := \lambda((-\gamma(t)) * \gamma(s)). $$
Then $F(t,t) = 0$ for each $t \in I$ and 
$$ \frac{\partial F}{\partial s}(t_0, t_0) = 1 = - 
 \frac{\partial F}{\partial t}(t_0, t_0). $$
Hence there exists an $\eps > 0$ such that 
$$ \frac{\partial F}{\partial s} - \frac{\partial F}{\partial t}\geq 1 
\quad \mbox{ on } \quad [t_0-\eps,t_0+\eps]^2. $$
We conclude that, for $|t-t_0|, |s-t_0| \leq \eps$, 
$$ F(s,t) > 0 \quad \mbox{ for } \quad s > t. $$
It follows in particular that 
$(-\gamma(t)) * \gamma(s) \not\in \h,$
and hence that 
$$\gamma(s) \exp_H(V_\h) \cap \gamma(t) \exp_H(V_\h) = \eset. $$
This contradicts the separability of $H$. 
We therefore have $\delta(\gamma)_t \in \h$ for each $t\in I$. 

As $H$ is $C^{k-1}$-regular, there exists a $C^k$-curve 
$\eta \: I \to H$ with $\delta(\eta) = \delta(\gamma)$ and 
$\eta(t_0) = \gamma(t_0)$. Then $\gamma = \eta$ follows from 
Lemma~\ref{lem:c.12b}, and we conclude that $\gamma$ is $C^k$ as a map $I \to H$. 
This completes the proof. 
\end{prf}

\begin{lem}\mlabel{lem:init-crit-arc}  Let $H$ be a subgroup of the Lie group $G$ 
and $H_a$ be the smooth arc-component of the identity in~$H$. 
If $H_a$ is an initial Lie subgroup of $G$, then so is~$H$. 
\end{lem}

\begin{prf} Suppose that $j_a \: H_a \to G$ defines an initial Lie subgroup 
structure on $H_a$. For every $h \in H \subeq N_G(H_a)$, the conjugation map 
$c_h \: H_a \to H_a$ defines a smooth map $c_h \circ j_a \: H_a \to G$, 
hence a smooth map $H_a \to H_a$.  Therefore 
$H$ carries a Lie group structure for which $H_a$ is an 
open subgroup (Corollary~\ref{cor:open-lie}). 
We claim that the inclusion $j \: H \to G$ is initial. 
To this end, let $f \: M \to G$ be a smooth map, where 
$M$ is a connected manifold and $f(M) \subeq H$. 
For any $m_0 \in M$, we then have $f(m_0)^{-1}f(M) \subeq H_a$, so that the 
initiality of $H_a$ implies that  the map 
$M \to H_a, m \mapsto f(m_0)^{-1} \cdot f$ is smooth, but then $f$ 
also defines a smooth map $M \to H$. 
\end{prf}

\begin{thm} \mlabel{thm:findim-initial} Any subgroup of a 
finite-dimensional Lie group carries an initial Lie subgroup structure. 
\end{thm}

\begin{prf} (cf.\ \cite[Thm.~9.6.13]{HiNe12}) 
Let $H_a \subeq H$ be the arc-component of $\be$ in~$H$, 
viewed as a topological subgroup of $G$. 
According to Yamabe's Theorem (\cite{Go69}, \cite[Sect.~9.6]{HiNe12}), 
the arc-component $H_a$ of $G$ is an integral subgroup, 
i.e., there exists a Lie subalgebra $\fh \subeq \L(G)$ such that 
$H_a = \la \exp_G \h \ra$. In particular, $H_a$ coincides with the 
smooth arc component of the identity in~$H$. 
Therefore it suffices to show that $H_a$ is initial 
(Lemma~\ref{lem:init-crit-arc}), and this follows from 
Theorem~\ref{thm:sep-subgroup2}. 
\end{prf} 

\subsection{Infinite-dimensional tori in Lie groups} 

Let $G$ be a Lie group with a smooth exponential function 
and $\fa \subeq \g$ be an abelian Lie 
subalgebra. Then the exponential function 
\[ \exp_\fa := \exp_G\res_\fa \: \fa \to G \] 
is a smooth group homomorphism by Lemma~\ref{lem:4.1.4}. 
As we have seen in the 
Integral Subgroup Theorem~\ref{thm:5.5.3}, 
it is important to understand the nature of  its kernel 
\[  \Gamma_\fa := \ker(\exp_\fa) = \exp_\fa^{-1}(\be).\] 

As $\exp_\fa$ is continuous, $\Gamma_\fa$ is a closed subgroup 
of the locally convex space~$\fa$. In the proof of the 
Integral Subgroup Theorem~\ref{thm:5.5.3} we already discussed 
such subgroups. In particular, we have seen that: 
\begin{itemize}
\item[$\bullet$] If $G$ is locally exponential, then $0$ is isolated in 
$\exp_G^{-1}(\be)$, 
so that $0$ is isolated in $\exp_\fa^{-1}(\be)$ and therefore 
$\Gamma_\fa$ is discrete. 
\item[$\bullet$] As $\L(\exp_\fa) \:  \fa \to \fg$ is injective, 
the derivative of each $C^1$-curve \break 
$\gamma \: [0,1] \to \Gamma_\fa$ is trivial, so that 
each $C^1$-curve in $\Gamma_\fa$ is constant. 
In particular $\Gamma_\fa$ contains no non-trivial vector subspace. 
If $\fa$ is finite-dimensional, this implies that $\Gamma_\fa$ is 
discrete (Exercise~\ref{exer:3.5.6}). 
\item[$\bullet$] 
In the proof of Theorem~\ref{thm:findim-discrete} we have also seen that  
$\Gamma_\fa$ is discrete if $\fa$ is an abelian subalgebra consisting 
of complete vector fields in the Fr\'echet--Lie algebra $\cV(M)$ of smooth 
vector fields on a $\sigma$-compact finite-dimensional manifold~$M$ 
with finitely many connected components. This is due to the fact that the 
dimension of a torus acting faithfully on $M$ is bounded from  above. 
\end{itemize}

If $\Gamma_\fa$ is discrete, then 
$A := \fz/\Gamma_\fa \cong \exp_\fa(\fa)$ carries a natural 
Lie group structure. Therefore the main question is to find sufficient 
conditions for this to be the case. 
Presently we are not aware of any example for which $\Gamma_\fa$ is 
not discrete. As the following theorem shows, such examples would lead 
to large ``non-Lie'' subgroups of $G$. 

\begin{thm} {\rm(Subtorus Theorem)}  \mlabel{thm:comp-tor}
Let $G$ be a Fr\'echet--Lie group and 
$\fa \subeq \g$ be a closed abelian subalgebra for which 
the subgroup $\Gamma_\fa$ is not discrete. Then 
the closure of $\exp_G(\fa)$ contains an infinite-dimensional torus 
subgroup~$T \cong \T^\N$. 
\end{thm}

\begin{prf} We fix a basis $(U_n)_{n \in \N}$ of absolutely convex $0$-neighborhoods 
in $\fa$ satisfying 
\[ U_n + U_n \subeq U_{n-1} \quad \mbox{ for } \quad n > 1.\] 
Inductively, this hypothesis leads for $m > n$ to  
\begin{equation}
  \label{eq:u-sum}
U_n + \cdots + U_m \subeq U_{n-1}.
\end{equation}

That $\Gamma_\fa$ is not discrete is equivalent to the existence 
of a sequence of non-zero elements $x_n \in \fa$ with $x_n \to 0$ 
and $\exp_G x_n= \be$. Passing to a subsequence, we may assume that 
the $x_n$ are linearly independent and that 
$x_n \in U_n$ holds for every $n \in \N$. Then 
\eqref{eq:u-sum} implies that, for 
$0 \leq a_n \leq 1$, the series $\sum_{n = 1}^\infty a_n x_n$ converges 
by the Cauchy Criterion because $\fa$ is complete. 

In particular, we obtain a well-defined map 
\[ S \: [0,1]^\N \to \fa, \quad 
\psi((a_n)_{n \in \N}) := \sum_{n=1}^\infty a_n x_n \] 
whose continuity with respect to the product topology on $[0,1]^\N$ 
also follows immediately from \eqref{eq:u-sum}. 
Note that 
\[ q \: [0,1]^\N \to \T^\N, \quad q((a_n)_{n \in \N}) = (e^{2\pi i a_n})_{n\in \N} \] 
is a continuous surjective map. Now the compactness of $[0,1]^\N$ 
implies that $q$ is a quotient map. We have 
$q(a)= q(b)$ if and only if $b_n - a_n \in \Z$ and 
then 
\begin{align*}
\exp_\fa\Big(\sum_{n = 1}^\infty a_n x_n\Big) 
&= \lim_{N \to \infty} \exp_\fa\Big(\sum_{n=1}^N a_n x_n\Big) \\
&= \lim_{N \to \infty} \exp_\fa\Big(\sum_{n=1}^N b_n x_n\Big) 
=  \exp_\fa\Big(\sum_{n = 1}^\infty b_n x_n\Big).
\end{align*}
Therefore $\exp_\fa \circ S \:  [0,1]^\N \to G$ factors through a 
continuous map 
\[ \psi \: \T^\N \to G, \qquad 
(e^{2\pi i a_n})_{n \in\N} \mapsto 
\lim_{N \to \infty} \exp_\fa\Big(\sum_{n = 1}^N a_n x_n\Big) \] 
which is easily seen to be a group homomorphism. 

Therefore $T := \psi(\T^\N) \cong \T^\N/\ker \psi$ is a compact abelian 
subgroup of $G$ whose character group $\hat T := \Hom(T,\T)$ is isomorphic 
to 
\[ \{ \chi \in \hat{\T^\N} \cong \Z^{(\N)}
\: \ker \psi \subeq \ker \chi\}.\]  
As subgroups of free abelian groups are free, $\hat T$ is a 
free abelian group and therefore $T \cong \hats{T}$ is a torus 
(\cite[Thm.\ 7.63]{HoM98}). 
That it is infinite-dimensional follows from the fact that 
it contains all tori 
\[ T_n := \exp(\R x_1 + \cdots + \R x_n)\] whose dimension 
equals $n$ because the $x_n$ are linearly independent. 
\end{prf} 

\begin{rem} As above, let $\fa \subeq \g$ be an abelian 
Lie subalgebra of the Lie algebra $\g$ of the Fr\'echet--Lie group~$G$.
Suppose that $x_n \in \fa$ is a sequence of non-zero 
elements in $\Gamma_\fa$ converging to $0$ and 
consider the subspaces $\fa_n := \Spann \{ x_1, \ldots, x_n\}$. 
We consider the 
subgroup $\Gamma := \sum_{n \in \Z} \Z x_n\subeq \Gamma_\fa$ 
generated by the~$x_n$. If $\Gamma$ is closed, then 
it is a complete metric space, hence Baire. As 
$\Gamma$ is the union of the subgroups 
$\Gamma_n := \Z x_1 + \ldots + \Z x_n$, Baire's Theorem 
implies that some $\Gamma_n$ contains an open subset of 
$\Gamma$. But the subspaces $\fa_n$ are finite-dimensional, so that 
$\Gamma_{\fa_n}$ is discrete. This implies that $\Gamma$ is discrete, 
contradicting $x_n \to 0$.  Therefore $\Gamma$ is not closed. 
\end{rem}

\begin{ex} Let $\cH$ be an infinite-dimensional complex Hilbert space 
with orthonormal basis $(e_n)_{n\in\N}$. 
Then its unitary group $\U(\cH)$ is a Lie group with respect to the norm 
topology. This group contains an infinite-dimensional torus 
\[ T := \{ U \in \U(\cH) \: (\forall n \in \N) \ U e_n \in \T e_n\}.\] 
As an abstract group $T$ is isomorphic to $\T^\N$, but the topology $T$ inherits 
from the operator norm is much finer than the norm topology. It corresponds 
to the metric 
\[ d(U,V) = \sup \{ |u_n - v_n| \: n \in \N\}, 
\quad \mbox{ where } \quad 
U_n e_n = u_n e_n, 
V_n e_n = v_n e_n.\] 

The elements $X_n \in \fu(\cH) = \{ X \in \cL(\cH) \: X^* = - X\}$ 
given by $X_n e_m = 2\pi i \delta_{nm} e_m$ satisfy 
$\exp X_n = \be$ in $\U(\cH)$ but $\|X_n\| = 2\pi$, so that 
the $X_n$ do not converge to $0$. However, 
$X_n \to 0$ in the strong operator topology, i.e., in the topology 
of pointwise convergence on $\cL(\cH)$. As the sequence $(X_n)_{n \in \N}$ 
is uniformly bounded, it follows that $X_n \to 0$ also holds in the 
topology of uniform convergence on compact subsets of $\cH$. 
This means that the corresponding linear vector fields $\hat X_n$ converge 
in $\cV(\cH)_c$ (the Lie algebra of smooth vector fields on $\cH$, 
endowed with the compact-open topology) to $0$.
\end{ex}

\begin{rem}
An important special case of a closed abelian Lie subalgebra 
is the center $\fz := \fz(\g)$. 
In Proposition~\ref{prop:initial-disc2}, we have seen that the 
discreteness of the subgroup $\Gamma_\fz$ implies that 
$Z(G)$ is an initial Lie subgroup of $G$. 
However, it may still happen that the Lie group topology 
on $Z(G)$ is strictly finer than the subspace topology inherited from 
$G$ (Example~\ref{ex:nonlie-cen}). 
If $\g$ is Fr\'echet and $\Gamma_\fz$ is not discrete, then 
Theorem~\ref{thm:comp-tor} implies the existence of an 
infinite-dimensional torus in $Z(G)$. 
\end{rem}

\section{Locally exponential subgroups} \mlabel{sec:6.3}

In this subsection we introduce the important concept of a locally exponential 
Lie subgroup of a locally exponential Lie group. 
It generalizes in a natural way the situation one encounters in finite-dimensional 
Lie groups for closed subgroups. In particular, we shall see below that 
each locally compact subgroup of a locally exponential Lie group 
is a finite-dimensional Lie group. 
Another important point is that 
inverse images of locally exponential Lie subgroups under smooth homomorphisms 
do also have this property, which implies in particular that all kernels 
are locally exponential Lie subgroups. 
In the next section we shall reverse this perspective in the sense that we 
characterize those normal Lie subgroups $N$ for which the quotient $G/N$ 
is a locally exponential Lie group. 

\subsection{Locally exponential subgroups and their Lie algebras}

\begin{defn} \mlabel{def:5.5.7} 
Let $G$ be a locally exponential Lie group. A subgroup \break $H \subeq G$ 
is called a 
\index{subgroup!locally exponential} 
\index{Lie subgroup!locally exponential} 
{\it locally exponential (Lie) subgroup} if it is a locally exponential 
Lie group with respect to the topology inherited from $G$ (cf.\ Corollary~\ref{cor:5.4.5}). 

\index{Lie algebra!of locally exponential subgroup} 
To define the {\it Lie algebra of a locally exponential Lie subgroup}, 
we first recall that each continuous homomorphism $\R \to G$ is smooth, 
so that 
$$\gamma \mapsto \gamma'(0), \quad \Hom(\R,G) \to \L(G)$$ 
yields an 
identification of $\Hom(\R,G)$ with the Lie algebra of $G$ 
(cf.\ Theorem~\ref{thm:5.4.4}). Since the inclusion $H \into G$ is a 
topological embedding, loc.\ cit.\ also implies that $\iota \: H \to G$ is a 
morphism of Lie groups, so that we have a morphism of topological Lie algebras 
$$ \L(\iota) \: \L(H) \to \L(G). $$
Identifying $\L(H)$ likewise with $\Hom(\R,H)$, we see that 
the injectivity of $\iota$ implies that $\L(\iota)$ is injective, so that 
$$ \im(\L(\iota)) = \{ x \in \L(G) \: \exp_G(\R x)\subeq H \} =: \L^e(H) $$
can be identified with the Lie algebra of $H$. In this sense, we call 
$\L^e(H)$ the {\it Lie algebra of the subgroup $H$}. 
\index{Lie algebra!of subgroup $H$, $L^e(H)$} 
\end{defn}

\begin{lem} \mlabel{lem:locexp-idcomp} For a topological group $G$, 
the following are equivalent: 
  \begin{enumerate}
  \item[\rm(i)] $G$ is a locally exponential Lie group. 
  \item[\rm(ii)] The identity component $G_0$ is a locally exponential Lie group. 
  \end{enumerate}
\end{lem}

\begin{prf} Clearly, (i) implies (ii). If, conversely, (ii) is satisfied, then 
each conjugation $c_g$, $g \in G$, defines a topological automorphism of 
$G_0$, hence a Lie automorphism. Now the assertion follows from 
Corollary~\ref{cor:open-lie}. 
\end{prf}

The preceding lemma reduces the problem to decide whether a 
subgroup is locally exponential to the case of connected subgroups. 
To get a better picture of the structure of a subgroup close to the 
identity, we first record 
the following observation concerning arbitrary closed subgroups 
of locally exponential Lie groups: 
\begin{prop} \mlabel{prop:5.4.3} 
Let $G$ be a locally exponential Lie group and $H \subeq G$ 
be a closed subgroup. 
\begin{enumerate}
\item[\rm(a)] The subset 
$$ \L^e(H) := \{ x \in \L(G) \: \exp_G(\R x) \subeq H \} $$
is a closed Lie subalgebra of $\L(G)$, satisfying 
$$ \Ad(h)\L^e(H) = \L^e(H) \quad \hbox{ for all } \quad h \in H. $$
In particular, $\L^e(H)$ is a stable subalgebra of $\L(G)$. 
\item[\rm(b)] If, in addition, $H$ is a normal subgroup, then 
$$ \Ad(g)\L^e(H) = \L^e(H) \quad \hbox{ for all } \quad g \in G. $$
In particular, $\L^e(H) \trile \L(G)$ is a stable ideal. 
\item[\rm(c)] $\L^e(H)= \L^d(H)$. 
\end{enumerate}
\end{prop}

\begin{prf} (a) Since $H$ is closed, the fact that $\L^e(H)$ is a closed Lie subalgebra of $\L(G)$ 
is a direct consequence of the product and the commutator
formula (Proposition \ref{prop:5.4.2}). Since 
$\L(c_g) = \Ad(g)$, the invariance of $H$ under the automorphisms
$c_h$, $h \in H$, implies that $\L^e(H)$ is invariant under $\Ad(H)$. 

(b) follows from the invariance of $H$ under all inner
automorphisms $c_g$, $g \in G$, and $\L(c_g) = \Ad(g)$. 

(c) Clearly, $\L^e(H) \subeq \L^d(H)$. For the converse, 
consider a $C^1$-curve $\gamma \: [0,1] \to G$ with 
$\gamma([0,1]) \subeq H$, $\gamma(0) = \be$ and $\gamma'(0) = x$. 
Then 
$$ \alpha(t) := \lim_{n \to \infty} \gamma\Big(\frac{t}{n}\Big)^n $$
exists and equals $\exp_G(tx)$ (Exercise~\ref{exer:5.1.1}). In fact, 
let $U_\g \subeq \g$ be a star-like open $0$-neighborhood 
for which $\exp_G\res_{U_\g}$ is a diffeomorphism onto the 
open $\be$-neighborhood $U_G := \exp_G(U_\g)$ in $G$. 
Suppose that $\gamma([0,\eps]) \subeq U_G$, so that there 
exists a smooth curve $\eta \: [0,\eps] \to U_\g$ with 
$\exp_G \circ \eta = \gamma$. We then have, for $t \leq n \eps$, 
$$ \gamma\Big(\frac{t}{n}\Big)^n 
= \exp_G\Big(\eta\Big(\frac{t}{n}\Big)\Big)^n
= \exp_G\Big(n\eta\Big(\frac{t}{n}\Big)\Big)
\to \exp_G(t \eta'(0)) = \exp_G(tx). $$
Now the closedness of $H$ implies that 
$\exp_G(\R x) \subeq H$, so that $x \in \L^e(H)$. 
\end{prf}

From (c) above, we further obtain: 

\begin{cor}\mlabel{cor:c1-disc} 
In a closed subgroup $H$ of a locally exponential Lie group~$G$, 
all $C^1$-arcs in $H$ are constant if and only if 
$\L^e(H) = \{0\}$.   
\end{cor}

\begin{ex} \mlabel{ex:L2-z} In the Hilbert space 
$E = L^2([0,1],\R)$, we consider the closed subgroup 
$H := L^2([0,1],\Z)$ of all $\Z$-valued functions (up to sets of measure zero). 
Then $\L^e(H) = \{0\}$ follows from the fact that 
$\R f \subeq H$ implies $f = 0$. 
Now Corollary~\ref{cor:c1-disc} implies that all $C^1$-curves in $H$ are constant. 
\end{ex}

\begin{prop} \mlabel{prop:5.4.6} 
Let $\phi \: G \to H$ be a morphism of the locally exponential 
Lie groups $G$ into the Lie group $H$ with a smooth exponential function. 
Then
    \begin{eqnarray}
      \label{eq:5.4.2}
\L^e(\ker \phi) = \ker \L(\phi).        
    \end{eqnarray}
\end{prop}

\begin{prf} For $x \in \L(G)$, the relation 
$x \in \L^e(\ker \phi)$ is equivalent to 
$\exp_G(\R x) \subeq \ker \phi$, i.e., 
$$ \phi(\exp_G(tx)) = \exp_H(t \L(\phi)x) = \be \quad \mbox{ for all } \quad 
t \in \R. $$
This is equivalent to $\L(\phi)x =0$. 
%
\end{prf}

\begin{lem} \mlabel{lem:liealg-emb} 
For any locally exponential Lie subgroup $H$ of the locally exponential 
Lie group $G$, the canonical map 
$\L(H) \into \L(G)$ is a topological embedding with image $\L^e(H)$. 
\end{lem}

\begin{prf} Let $i_H \: H \to G$ denote the topological embedding of $H$ in $G$. 
Then the relation 
$i_H \circ \exp_H = \exp_G \circ \L(\iota_H)$
implies that there exists an open  \break $0$-neighborhood $U \subeq \L(H)$ such that 
$\L(\iota_H) \res_U$ is a topological embedding into $\L(G)$, but this implies that 
$\L(\iota_H)$ is a topological embedding because it is linear. 
\end{prf}

\begin{prop} \mlabel{prop:7.3.8} If $G = (E,+)$ is  the additive group of a 
locally convex space, then $G$ is an exponential Lie group 
with $\exp_E = \id_E$. A subgroup $H$ of $E$ is 
locally exponential if and only if the maximal linear 
subspace contained in $H$ is an open subgroup of~$H$. 
\end{prop}

\begin{prf} If $H \leq G$ is a locally exponential Lie subgroup, then 
the preceding lemma implies that 
$\L^e(H) = \{ x \in E \: \R x \subeq H\},$
the largest linear subspace contained in $H$, is isomorphic to $\L(H)$. 
Therefore $\L^e(H)$ is the identity component of $H$, hence an 
open subgroup. 

If, conversely, $\L^e(H)$ is an open subgroup of $H$, 
then the subspace topology turns $H$ into a topological group whose identity 
component is $(\L^e(H),+)$, and this implies that $H$ is locally exponential.   
\end{prf}

Since linear subspaces need not be closed, the preceding proposition shows in particular
that there are many locally exponential Lie subgroups which are not closed. The following 
lemma provides a criterian for a, not necessarily closed, subgroup $H$ to be locally exponential. 
It will be refined below for closed subgroups.

\begin{lem} \mlabel{lem:equiv-locexp-Liesub} For a subgroup $H$ 
of a locally exponential Lie group $G$, the 
following are equivalent: 
\begin{enumerate}
\item[\rm(i)] $H$ is a locally exponential Lie subgroup. 
\item[\rm(ii)] $\L^e(H)$ is a locally exponential subalgebra of $\L(G)$ and 
there exists an open $0$-neighborhood $U \subeq \L(G)$ 
for which $\exp_G\res_U$ is a diffeomorphism onto an open subset of $G$ and 
$$ \exp_G(U) \cap H = \exp_G(U \cap \L^e(H)). $$
\end{enumerate}
\end{lem}

\begin{prf} (i) $\Rarrow$ (ii): Let $V \subeq \L(G)$ be an open convex 
$0$-neighborhood  for 
which $\exp_G\res_V$ is a diffeomorphism onto an open subset of $G$. 
Since the inclusion $H \into G$ is a topological embedding, there exists a 
$0$-neighborhood $U_\h \subeq \L(H)$ such that 
$\exp_H\res_{U_\h}$ is a diffeomorphism onto an open subset of $H$ and 
$U_\h \subeq V$, if we identify $\L(H)$ with the Lie subalgebra 
$\L^e(H)$ of $\L(G)$ 
(cf.\ Lemma~\ref{lem:liealg-emb}). Since $H$ is a topological subgroup of $G$, 
there exists an open subset $U_\g \subeq V$ with 
$$ \exp_G(U_\g) \cap H = \exp_H(U_\h) = \exp_G(U_\h), $$
which leads to $U_\h = U_\g \cap \L^e(H)$, and (ii) follows. 

(ii) $\Rarrow$ (i): In view of the Integral Subgroup Theorem~\ref{thm:5.5.3}, 
the inclusion of locally exponential Lie algebras 
$\alpha \: \L^e(H) \to \L(G)$ integrates to a morphism of locally 
exponential Lie groups $\beta \: H^L \to G$ with $\L(H^L) = \L^e(H)$. Then 
$\beta(H^L) = \la \exp_G \L^e(H)\ra$, so that (ii) implies that 
$\beta(H^L)$ coincides with the identity component $H_0$ of $H$. 
As $\L(\beta) = \alpha$ is a topological 
embedding, and $H^L$ and $G$ are locally exponential, 
condition (ii) further implies that $\beta$ is an embedding, hence a topological 
isomorphism onto $H_0$. Thus $H_0$ is locally exponential and 
Lemma~\ref{lem:locexp-idcomp} shows that $H$ is locally exponential. 
\end{prf}

The following proposition is a useful criterion for 
subgroups of locally exponential Lie groups to be locally exponential 
Lie groups. 

\begin{prop} \mlabel{prop:locexp-subgroup-crit} 
Let $G$ be a locally exponential Lie group and $H \subeq G$ be 
a closed subgroup 
for which there exists an open $0$-neighborhood $V \subeq \L(G)$ such that 
$\exp_G\res_V$ is a diffeomorphism onto an open subset of $G$ and 
\begin{eqnarray}
  \label{eq:liesub}
\exp_G(V \cap \L^e(H)) = (\exp_G V) \cap H.
\end{eqnarray}
Then $H$ is a locally exponential Lie subgroup and a submanifold of $G$. 
\end{prop}

\begin{prf} We first recall from Proposition~\ref{prop:5.4.3} that 
$\h := \L^e(H)$ is a closed Lie subalgebra of $\g := \L(G)$. 
Let $V \subeq \g$ be an open symmetric $0$-neigh\-bor\-hood 
for which $\exp_G\res_V \: V \to \exp_G(V)$ is a diffeomorphism onto an
open subset of $G$ and (\ref{eq:liesub}) holds. 
Using the chart $((\exp_G\res_V)^{-1},\exp_G(V))$ of $G$, we see that 
$H$ is a Lie subgroup of $G$ in the sense of Definition~\ref{defn:liesubgroup}, 
hence in particular a submanifold and a Lie group with respect to the inherited 
manifold structure. 

Next we show that $\exp_H := \exp_G\res_\h \: \h \to H$ is a smooth 
exponential function for $H$, endowed with the Lie group structure from above. 
In fact, for each $x \in \h$, the curve 
$\gamma_x \: \R \to H, t \mapsto \exp_H(tx)$ 
is a group homomorphism which is smooth in a $0$-neighborhood, hence smooth, 
and satisfies $\gamma_x'(0) = x$. Therefore $\exp_H$ is an exponential 
function for $H$. Its smoothness follows from its smoothness on 
a $0$-neighborhood and Lemma~\ref{lem:5.1.2}. 
Since $\exp_H$ coincides on a $0$-neighborhood of $\h$ with 
$(\exp_G\res_V)\res_\h$, which defines the charts of $H$, 
the group $H$ is locally exponential and hence $\h = \L^e(H)$ is locally 
exponential. 
\end{prf} 

\begin{prop} {\rm(Inverse images of locally exponential subgroups)} 
\mlabel{prop:5.5.9} 
If $\phi \: G_1 \to G_2$ is a morphism of
locally exponential Lie groups and $H \subeq G_2$ is a closed 
locally exponential subgroup, 
then $H_1 := \phi^{-1}(H)$ is a locally exponential subgroup with 
Lie algebra 
$$\L^e(H_1) = \L(\phi)^{-1}\big(\L^e(H)\big).$$ 
In particular,  
$\ker \phi$ is a locally exponential subgroup of $G_1$ with Lie algebra 
$\ker \L(\phi)$. 
\end{prop}

\begin{prf} We choose an open symmetric 
$0$-neighborhood $V \subeq \L(G_2)$ such that $\exp_{G_2}\res_V$ is a
diffeomorphism onto an open subset of ${G_2}$, 
and $$\exp_{G_2}(V \cap \L^e(H))= \exp_{G_2}(V) \cap H$$ 
(Lemma~\ref{lem:equiv-locexp-Liesub}). 
Then we choose an open symmetric $0$-neighborhood 
$V_1 \subeq \L(\phi)^{-1}(V)$ such that $\exp_{G_1}\res_{V_1}$ is a
diffeomorphism onto the open subset $\exp_{G_1}(V_1)$. 

Let $x \in V_1$ with $\exp_{G_1} x \in H_1$. 
Then 
$$ \phi(\exp_{G_1} x) = \exp_{G_2}(\L(\phi)x) \in \exp_{{G_2}}(V) \cap H $$
with $\L(\phi)x \in V$. Hence $\L(\phi)x \in \L^e(H)$ and therefore 
$$x \in \L^e(H_1) = \L(\phi)^{-1}(\L^e(H)). $$ 
In view of the latter equality,  Proposition~\ref{prop:invim} 
implies that $\L^e(H_1)$ is locally exponential. 
The preceding arguments show that 
$\exp_{G_1}(V_1) \cap H_1 \subeq \exp_{G_1}(V_1 \cap \L^e(H_1))$, 
and therefore equality, 
because the converse inclusion is trivial. 
Now Lemma~\ref{lem:equiv-locexp-Liesub} implies that 
$H_1$ is locally exponential. 

Since $\{\be\}$ is a  locally exponential subgroup of ${G_2}$, $\ker \phi$ is a  locally exponential subgroup of $G_1$.
\end{prf}

\begin{cor} \mlabel{cor:5.5.10} 
If $N \trile G$ is a closed normal subgroup such that the
quotient group $G/N$ carries a locally exponential 
Lie group structure for which the
quotient map $q \: G \to G/N$ is a morphism of Lie groups, then 
$N$ is a closed locally exponential subgroup with $\L^e(N) = \ker(\L(q))$. 
\end{cor}

\begin{prop} \mlabel{prop:5.5.10b} 
If $G$ is a connected locally exponential Lie group with Lie algebra~$\g$, then 
the following assertions hold: 
\begin{enumerate}
\item[\rm(i)] $Z(G)$ is a locally exponential subgroup of $G$ with Lie algebra $\z(\g)$, and 
\item[\rm(ii)] $G/Z(G)$ is a locally exponential Lie group.  
\end{enumerate}
\end{prop}

\begin{prf} From Corollary~\ref{cor:adj-smooth} we know 
that the adjoint representation defines a morphism 
$\Ad \: G \to G_{\rm ad}$ of locally exponential Lie groups. 
Therefore $Z(G) = \ker \Ad$ is a closed locally exponential 
subgroup with Lie algebra 
$\L^e(\ker \phi) = \ker \L(\phi) = \z(\g)$ (Proposition~\ref{prop:5.5.9}). 

Since $\L(\Ad) \: \g \to \L(G_{\rm ad}) = \g/\z(\g)$ is a 
quotient map, $\Ad$ is open. Hence it factors through a 
topological isomorphism $\oline\Ad \: G/Z(G) \to G_{\rm ad}$. 
This means that the topological group $G/Z(G)$ carries a locally 
exponential Lie group structure. 
\end{prf}

\subsection{Locally compact subgroups}

\begin{thm} \mlabel{thm:loccpt-subgrp} Each locally compact subgroup of a locally exponential 
group is a finite-dimensional Lie subgroup. 
\end{thm}

\begin{prf} Let $G$ be a locally exponential Lie group and 
$H \leq G$ be a locally compact subgroup. 
We fix an open convex symmetric $0$-neighborhood 
$U_\g \subeq \g$ such that $\exp_G\res_{U_\g}$ is a diffeomorphism onto 
an open subset of $G$. 

Since $H$ is locally compact and $G$ is Hausdorff, $H$ is locally closed, 
hence closed (Exercise~\ref{exer:6.1.3}) 
and $\h := \L^e(H)$ is a closed Lie subalgebra of $\g := \L(G)$. 
First we show that $\h$ is finite-dimensional. 
The set $U_H :=  (\exp_G\res_{U_\g})^{-1}(H)$ is a locally compact 
subset of $U_\g$, so that $\h \cap U_H$ is a locally compact 
subset of $\h$. On the other hand, $\exp_G(\h) \subeq H$ implies that 
$U_H$ is a $0$-neighborhood in $\h$, so that $\h$ is a locally compact 
topological vector space, hence finite-dimensional 
by Proposition~\ref{lcpfindim}. 
Since each finite-dimensional 
Lie algebra is locally exponential (Corollary~\ref{exs:5.2.3}), 
we may shrink $U_\g$ in such a way 
that $U_\h := U_\g \cap \h$ satisfies $U_\h * U_\h \subeq \h$ (Theorem~\ref{thm:5.2.11c}). 

As $\h$ is finite-dimensional, the Hahn--Banach Theorem implies 
that the linear map $\id_\h \: \h \to \h$ extends to a continuous linear projection 
$p_\h \: \g \to \h$, so that $E := \ker p_\h$ is a closed subspace of $\g$ for 
which the map 
$$ \g \to E \times \h, \quad x \mapsto (x - p_\h(x),p_\h(x)) $$
is a topological isomorphism. 

Let $V \subeq U_\g$ be an open $0$-neighborhood with 
$\exp_G(V)\exp_G(V) \subeq \exp_G(U_\g)$, and put $V_\h := V \cap \h$. 
Then 
$$ F \: V \times V_\h \to \h, \quad (x,y) \mapsto p_\h(x*y), $$
is a smooth map. We put $F_x(y) := F(x,y)$. 
Then $T_0(F_0) = \id_\h$, so that 
the Implicit Function Theorem in \cite[Th.~2.3(i)(c)]{Gl04a}, 
implies the existence of open $0$-neighborhoods 
$W \subeq V$ and $W_\h \subeq V_\h$ such that 
$$U := \bigcup_{x \in W} \{x\} \times F_x(W_\h)$$ 
is an open subset 
of $\g \times \h$, the map 
$$ W \times W_\h \to U, \quad (x,y) \mapsto (x, F_x(y)) $$ 
is a diffeomorphism, and the map 
$U \to W_\h, (x,y) \mapsto F_x^{-1}(y)$ 
is smooth. In particular, the map 
$$ g \: U' := U \cap (\g \times \{0\}) \to W_\h, \quad x \mapsto F_x^{-1}(0) $$
is smooth and satisfies 
$x * g(x) \in E = \ker p_\h.$ Note that $U' \subeq W \subeq V$. 
Hence each $x \in U'$ has a decomposition 
$$ x = (x * g(x)) * (-g(x)) \quad \mbox{ with } \quad 
x * g(x) \in E, g(x) \in \h. $$
For $x \in U' \cap E$ we have $F_x(0) =0$, so that $g(x) = 0$, 
and for $x \in U' \cap U_\h$ the relation $F_x(-x) = 0$ leads to $g(x) = -x$. 

We want to show that $U'_\h := U' \cap \h$ is open in $U'_H := U' \cap U_H \subeq V$.
If $x \in U'_H$, then $g(x) \in W_\h \subeq V_\h$ implies that 
$$\exp_G(x * g(x)) = \exp_G(x) \exp_G(g(x))\in H,$$ 
so that 
$x * g(x) \in U_H \cap E$. If $U'_\h$ is not open in $U_H'$ and 
$x_i \in U_H' \setminus U_\h'$ converges to some $x \in U'_\h$, 
then $x_i * g(x_i) \to x * g(x) = x * -x = 0$ 
with $x_i * g(x_i) \not= 0$. Hence it suffices to show that 
$0$ is isolated in $U_H \cap E$.

Since $U_H$ is locally compact, there exists a continuous seminorm 
$p$ on $\g$ such that 
$p^{-1}([0,2]) \subeq U_\g$ and $\Gamma_H := p^{-1}([0,2]) \cap U_H$ is compact. 
We claim that $0 \not= x \in U_H \cap E$ implies $p(x) > 0$. In fact, 
if $p(x) = 0$, then $\R x \subeq p^{-1}([0,2]) \subeq U$,  
and the local group property of $U_H$ implies that 
$\Z x \subeq p^{-1}(0) \cap U_H$, contradicting the compactness of 
$p^{-1}(0) \cap U_H$. 
Now $$U_H^1 := \{ x \in U_H \cap E \: p(x) = 1\}$$ is a compact subset of $\g$. 
Arguing by contradiction, we assume that $0$ is not isolated in 
$U_H \cap E$, so that there exists a net $x_i \in U_H\cap E \setminus \{0\}$,  
converging to $0$. Then $p(x_i) \to 0$ and $p(x_i) > 0$ by the preceding 
argument. For $r \in \R$, let 
$$\lfloor r \rfloor := \max \{m \in \Z \: m \leq r\}. $$
We put $n_i := \lfloor\frac{1}{p(x_i)}\rfloor$. 
Then $p(n_i x_i) = n_i p(x_i) \leq 1$ 
and $n_i x_i = x_i * \ldots * x_i \in U_H$. 
Furthermore, $x_i \to 0$ and $\frac{1}{p(x_i)} - n_i \leq 1$ lead to 
$$ n_i x_i - \frac{x_i}{p(x_i)} = \Big(n_i - \frac{1}{p(x_i)}\Big) x_i \to 0. $$
In particular, $n_i p(x_i) = p(n_i x_i)\to 1$. 

The closed subset 
$$\Gamma_H' := \big\{ x \in U_H \: 
\shalf \leq p(x) \leq 1\big\} \subeq \Gamma_H $$
of $\Gamma_H$ is also compact, and we may w.l.o.g.\ assume that $n_i x_i \in \Gamma_H'$ for 
each~$i$. Let $y \in \Gamma_H'$ be a cluster point of the net 
$(n_i x_i)$. Passing to a subnet, we may assume that 
$n_i x_i \to y$. 

Fix $m \in \N$. We write $n_i = q_i m + r_i$ with 
$0 \leq r_i \leq m-1$ and $r_i, q_i \in \Z$. Then 
$$ q_i x_i - \frac{y}{m}
= \frac{1}{m}(m q_i x_i - y) 
= \frac{1}{m}(\underbrace{n_i x_i - y}_{\to 0} - \underbrace{r_i x_i}_{\to 0}) 
\to 0. $$
From $q_i \leq n_i$ and $x_i, 2 x_i, \ldots, n_i x_i \in U_H$, we 
derive $q_i x_i \in U_H$, which leads to 
$\frac{1}{m}y \in U_H$ because the closedness of $H$ implies the 
closedness of $U_H$ in $U_\g$. As $m$ was arbitrary, the local group property 
of $U_H$ implies that 
$$ (\Q \cap [0,1])y \subeq U_H, $$
and the closedness of $U_H$ in $U_\g$ entails 
$[0,1]y \subeq U_H$. This in turn leads to 
$\exp_G([0,1]y) \subeq H$, and hence to $y \in \L(H) = \h$, 
contradicting $y \in E\setminus\{0\}$. 
\end{prf}

\begin{cor} Each closed subgroup $H$ of a finite-dimensional Lie group is a Lie subgroup. 
\end{cor}

\begin{cor} \mlabel{cor:7.3.16} {\rm(Canonical factorization for 
finite-dimensional Lie groups)} 
If $\phi \: G \to H$ is a morphism of Lie groups and $G$ 
is finite-dimensional, 
then $N := \ker \phi$ is a Lie subgroup of $G$, $G/N$ is a Lie group, and 
$\phi$ induces a morphism of Lie groups $\oline\phi \: G/N \to H$. 
\end{cor}

\begin{prf} The preceding corollary implies that $N$ is a Lie subgroup of $G$, 
and Proposition~\ref{prop:split-crit-ban} implies that $N$ is split, so that 
$G/N$ carries a natural Lie group structure for which the quotient map 
$q \: G \to G/N$ defines on $G$ the structure of an $N$-principal bundle. 
We conclude in particular that the smooth map $\phi$, 
which is constant on the 
$N$-cosets, factors through a smooth map $\oline\phi \: G/N \to H$. 
\end{prf}

\begin{cor} \mlabel{cor:diff-inj} If $\phi \: G \to H$ 
is a morphism of Lie groups for which 
$G$ is finite-dimensional and $\L(\phi)$ is injective, 
then $\ker \phi$ is discrete. 
\end{cor}

\begin{prf} Since $N := \ker \phi$ is a closed subgroup of the finite-dimensional 
Lie group $G$, it is a Lie subgroup. From 
$\L(N) = \ker \L(\phi) = \{0\}$ we thus derive that $N$ is discrete. 
\end{prf}

\begin{ex} Theorem~\ref{thm:loccpt-subgrp} generalizes the result that 
closed subgroups of finite-dimensional Lie groups are Lie subgroups 
(see \cite{vN29} for closed subgroups of $\GL_n(\R)$) 
and that locally compact subgroups of Banach--Lie groups are 
finite-dimensional Lie groups (cf.\ \cite[Th.\ 5.41(vi)]{HoM98}). 

How bad closed subgroups can be is
illustrated by the following example due to K.~H.~Hofmann: 
We consider the real Hilbert space $G := L^2([0,1],\R)$ as a
Banach--Lie group. Then the subgroup $H := L^2([0,1],\Z)$ of all those
functions which almost everywhere take values in $\Z$ is a closed
subgroup. Since the one-parameter subgroups of $G$ are of the form 
$t \mapsto tf$, $f \in G$, we have $\L^e(H) = \{0\}$ 
(cf.\ Example~\ref{ex:L2-z}). 
On the other hand, the group $H$ is arcwise connected, and even
contractible, because the map $F \: [0,1] \times H \to H$ given by 
$$ F(t,f)(x) := 
\left\{ 
  \begin{array}{cl} 
f(x) & \mbox{for $0 \leq x \leq t$} \\ 
0& \mbox{for $t < x \leq 1$} 
\end{array} \right. $$
is continuous with $F(1,f) = f$ and $F(0,f) = 0$. 

This pathology can be avoided by the assumption that the subgroup 
is connected by $C^1$ arcs (cf.\ Theorem~\ref{thm:smootharc} below). 
In view of Corollary~\ref{cor:c1-disc}, all $C^1$-arcs in 
$H$ are constant because $\L^e(H) = \{0\}$. 
\end{ex}

\begin{ex} \mlabel{ex:subgroup-circle} (Closed subgroups of $\T$) 
Let 
\[ H \subeq \T = \{ z \in \C^\times \: |z| = 1\} \] 
be a closed proper subgroup 
and $q = \exp_\T \: \R \to \T, t \mapsto e^{2\pi i t}$, be the exponential 
function of $\T$. 
Then $q^{-1}(H)$ is a closed proper subgroup of $\R$, hence cyclic 
(cf.\ Exercise~\ref{exer:3.5.1}). Since it contains $1$, its generator 
is a rational number. We conclude that $H = q(q^{-1}(H))$ is one of the groups 
$$ C_n := \{ z \in \T \: z^n = 1\} $$
of $n$th roots of unity. 
\end{ex}

The discussion of the following example describes in particular 
the prototype of a situation where the image of a 
morphism of Lie groups (here a one-parameter group) is not a submanifold. 

\begin{ex} (Subgroups of $\T^2$) \mlabel{ex:dense-wind}
(a) Let $H \subeq \T^2$ be a closed proper subgroup. Then 
$\L(H) \not=\L(\T^2)$ implies $\dim H < \dim \T^2  = 2$. 
Further, $H$ is compact, so that the  
group $\pi_0(H)$ of connected components of $H$ is finite. 

If $\dim H = 0$, then $H$ is finite, and for $n := |H|$ it is 
contained in a subgroup of the form 
$C_n \times C_n$, where $C_n \subeq \T$ is the subgroup of $n$th 
roots of unity 
(cf.\ Example~\ref{ex:subgroup-circle}). 

If $\dim H = 1$, then $H_0$ is a compact connected $1$-dimensional Lie group, 
hence isomorphic to $\T$ (Exercise~\ref{exer:3.5.3}). Therefore 
$H_0 = \exp_{\T^2}(\R x)$ for some $x \in \L(H)$ with $\exp_{\T^2}(x) 
= (e^{2\pi i x_1}, e^{2\pi i x_2})= (1,1)$, 
which is equivalent to $x \in \Z^2$. We conclude that 
the Lie algebras of the closed one-dimensional subgroups 
are of the form $\L(H) = \R x$ for some $x \in \Z^2$. 

(b) For each $\theta \in \R\setminus \Q$, the image of the 
$1$-parameter group 
$$ \gamma \: \R \to \T^2, \quad t \mapsto (e^{i\theta t}, e^{it}) $$
is not closed because $\gamma$ is injective. Hence the closure of 
$\gamma(\R)$ is a closed subgroup of dimension at least $2$, which shows that 
$\gamma(\R)$ is dense in $\T^2$.  
\end{ex}

\subsection{The commutator subgroup}

\begin{lem} \mlabel{lem:5.5.10c} 
For a connected locally exponential Lie group 
$G$ with Lie algebra $\g$, the commutator group $[G,G]$ satisfies 
$$ [G,G] \subeq \la \exp \oline {[\g,\g]}\ra \subeq \oline{[G,G]} $$
and thus $\oline{[G,G]} = \oline{\la \exp {[\g,\g]} \ra}$. 
\end{lem}

\begin{prf} From the Commutator Formula, we derive that 
$$\exp_G([x,y]) \in \oline{[G,G]} \quad \mbox{ for } \quad x,y \in \g.$$ 
Therefore ${[\g,\g]} \subeq \L(\oline{[G,G]})$ (Proposition~\ref{prop:5.4.3}) leads to 
$\exp \oline {[\g,\g]} \subeq \oline{[G,G]}.$

For the converse inclusion, we consider the simply connected covering group 
$q_G \: \tilde G \to G$. Let $\Ab(\g)$ denote the 
quotient space $\g/\oline {[\g,\g]}$. 
Then there exists a continuous homomorphism 
$\phi \: \tilde G \to \Ab(\g)$ for which the differential 
$\L(\phi) \: \g \to \Ab(\g)$ is the quotient map (Theorem~\ref{thm:int-thm}). 
Since $\ker \phi$ is a locally exponential Lie subgroup (Proposition \ref{prop:5.5.9}), 
$$ [\tilde G, \tilde G] \subeq (\ker \phi)_0 = \la \exp \ker \L(\phi) \ra 
= \la \exp_{\tilde G} \oline {[\g,\g]} \ra. $$
Hence 
$$ [G,G] = q_G([\tilde G, \tilde G]) \subeq 
q_G(\la \exp_{\tilde G} \oline {[\g,\g]}\ra)
= \la \exp_{G} \oline {[\g,\g]} \ra. $$
This completes the proof.
\end{prf}

\begin{rem} If $G$ is a connected 
subgroup of some $\GL_n(\R)$, then its commutator group $[G,G]$ is 
closed because it is a semidirect product of a unipotent normal subgroup 
and a semisimple group (\cite[Cor.~14.5.7]{HiNe12}), 
hence in particular a Lie subgroup (Theorem~\ref{thm:loccpt-subgrp}). 
In general, the commutator subgroup of a finite-dimensional Lie group is not closed 
(Exercise~\ref{exer:5.6.8}).
\end{rem}

\subsection{Initiality of closed subgroups} 

We have already seen that all subgroups of finite-dimensional Lie groups 
are initial (Theorem~\ref{thm:findim-initial}). We do not know if this is still 
true for Banach--Lie groups, or even locally exponential Lie group, 
but the results in this subsection 
cover at least all closed subgroups of Banach--Lie groups. 
The following theorem is a variant of Theorem~\ref{thm:sep-subgroup2}, where 
we assume that the Lie algebra $\fh$ is closed, but not that it is separable.  

\begin{thm} \mlabel{thm:closed-initial} Let $G$ be a locally exponential Lie group and 
$H \subeq G$ be a closed subgroup for which 
$\fh := \L^e(H)$ is integrable to a regular locally exponential Lie group. 
Then the subgroup $H$ is an initial Lie subgroup.
\end{thm}

\begin{prf} Let $j \:  H_1 \into G$ be an integral subgroup 
with Lie algebra $\fh$ 
(Theorem~\ref{thm:5.5.3}). Our assumption implies that 
$H_1$ is regular. Then $j(H_1) = \la \exp_G \fh \ra \subeq H$. 

Let $\phi \: M \to G$ be a smooth map with values in $H$. We have to show 
that $\phi$ defines a smooth map $M \to H$. Since this assertion is local, 
we may w.l.o.g.\ assume that $M$ is $1$-connected and 
that there exists an $m_0 \in M$ with ${\phi(m_0)=\be}$. 
For every $v \in T(M)$, we have 
$\delta(\gamma)v \in \L^d(H) = \fh$ 
(Proposition~\ref{prop:5.4.3}).  As $\fh$ is closed, $\delta(\gamma) 
\in \Omega^1(M,\fh)$. 
The regularity of $H_1$ then implies the existence of a smooth 
map $\eta \: M \to H_1$ with 
$\delta(\eta) = \delta(\gamma)$ and $j(\eta(m_0)) = \phi(m_0)= \be$. Then 
$j \circ \eta = \phi$ by the Uniqueness Lemma~\ref{lem:c.12b}. 
We conclude that $\phi$ defines a smooth map $M \to H$. 
This proves that $H$ is initial. 
\end{prf}

\begin{exs} \mlabel{ex:closed-initial} 
Typical examples satisfying the assumptions of the preceding theorem are:

\nin (a) $G = (E,+)$ a Mackey complete locally convex space and $H$ any 
closed subgroup. 

\nin (b) $G$ is a Banach--Lie group and $H$ any closed subgroup. 
Here $\L^e(H)$ is a Banach--Lie algebra of a closed subgroup, hence 
enlargeable (Corollary~\ref{cor:intsub-ban}) 
and every Banach--Lie group is regular (Theorem~\ref{thm:banach-regular}). 
\end{exs}

For closed subgroups of Banach--Lie groups, 
we even have the following refinement. 

\begin{thm} \mlabel{thm:smootharc} Let 
$G$ be a Banach--Lie group and $H \subeq G$ be a closed subgroup. 
Then $H$ carries the structure of an initial Lie subgroup, whose 
$C^1$-arc component is an integral subgroup for the closed Lie subalgebra 
$\L^e(H) \subeq \L(G)$. 
\end{thm}

\begin{prf} We have aready seen in Example~\ref{ex:closed-initial}(b) 
that $H$ carries the structure of an initial Lie subgroup whose 
identity component is the integral subgroup 
$i \: H_a \to G$ 
corresponding to the closed Lie subalgebra $\fh := \L^e(H)$ of $\g= \L(G)$. 

If $\gamma \: I \to G$ is a $C^1$-path with 
$\gamma(0) = \be$, then 
$\gamma(I) \subeq H$ implies 
that $\gamma'(0) \in \L^e(H) = \fh$ (Proposition~\ref{prop:5.4.3}). 
This shows that,  for any $C^1$-curve $\gamma \: I \to H \subeq G$, we have 
$\delta(\gamma)_t \in \fh$ for each $t \in I$. 
For any $C^1$-path $\gamma \: I \to G$ with values in $H$,  
we then obtain a continuous path $\delta(\gamma) \: I \to \fh$, 
and since the Banach--Lie group $H_a$ is $C^0$-regular 
(Theorem~\ref{thm:banach-regular}), there exists a $C^1$-path 
$\eta \: I \to H_a$ with $\delta(\eta) = \delta(\gamma)$ 
and $\eta(0)= \be$. Then 
the paths $\iota \circ \eta$ and $\gamma$ coincide 
by the Uniqueness Lemma~\ref{lem:c.12b}. 
We conclude that $i(H_a) = \la \exp \fh \ra$ 
is the $C^1$-arc component of the identity in $H$. 
\end{prf}

\subsection{Small subgroups of locally exponential groups} 

\begin{defn} We say that a topological group $G$ has 
\index{small subgroups} 
{\it no small subgroups} 
if there exists an identity neighborhood $U \subeq G$ such that 
$\{\be\}$ is the only subgroup of $G$ contained in~$U$.
\end{defn}

It is a classical result that Banach--Lie groups do not contain 
small subgroups. This is no longer true for locally convex vector groups, such 
as $G = \R^\N$, endowed with the product topology. Then each $0$-neighborhood 
contains non-zero vector subspaces, so that $G$ has small subgroups. 
For a locally convex space $E$, the non-existence of small subgroups 
is equivalent to the existence of a continuous norm on $E$ 
(Lemma~\ref{lem:8.11.2}) and we shall see that a locally exponential Lie group 
$G$ has small subgroups if and only if its Lie algebra 
$\g$ carries no continuous norm (Proposition~\ref{prop:8.11.3}).
Since any real vector space is torsion free, this implies that 
no locally exponential Lie group contains small torsion subgroups. 
For strong ILB--Lie groups, it is also known that they do not contain 
small subgroups \cite[Thm.~III.2.3]{Ne06a}, and this implies in particular 
that for each compact manifold $M$ the group $\Diff(M)$ 
does not contain small subgroups. 
We also know that direct limits of finite-dimensional Lie groups do not contain 
small subgroups (Theorem~\ref{thm:VII.1.3}). 

\begin{lem} \mlabel{lem:8.11.2}
  For a locally convex space $E$, the following are equivalent: 
  \begin{description}
  \item[\rm(a)] $E$ has no small subgroups, i.e., there exists a 
$0$-neighborhood not containing any non-trivial subgroup. 
\item[\rm(b)] There exists a $0$-neighborhood containing no non-zero 
linear subspace. 
\item[\rm(c)] $E$ possesses a continuous norm. 
  \end{description}
\end{lem}

\begin{prf} (a) $\Rarrow$ (b) is trivial. 

(b) $\Rarrow$ (c): Let $U \subeq E$ be a closed $0$-neighborhood 
not containing any non-zero subspace. We may w.l.o.g.\ assume that 
$U$ is absolutely convex. Then its Minkowski funtional 
\[ \mu_U(x)\,:=\,\inf\, \{ t>0\!: x\in tU\} \] 
is a continuous seminorm on $E$ (Proposition~\ref{Minkowsk}). 
As $U$ does not contain any non-trivial linear subspace, 
$\mu_U$ actually is a norm. 

(c) $\Rarrow$ (a): If $\|\cdot\|$ is a continuous norm 
on $E$, then its unit ball contains no non-trivial subgroup. 
\end{prf}

\begin{prop} \mlabel{prop:8.11.3} 
A locally exponential Lie group $G$ has no small subgroups 
if and only if its Lie algebra $\g$ possesses a continuous norm. 
\end{prop}

\begin{prf} Let $U \subeq \g$ be an open absolutely convex $0$-neighborhood 
for which $\exp\res_{2U}$ is a diffeomorphism onto an open subset of $G$. 
Let $H \subeq \exp(U)$ be a subgroup. We consider the subset 
\[ \tilde H := \exp^{-1}(H) \cap U = (\exp\res_U)^{-1}(H).\] 
If $x \in \tilde H$, then $\exp(nx) = \exp(x)^n \in H$ holds 
for each $n \in \Z$. Assume that $\R x \not\subeq U$. 
Then $\R x \cap U = ]-T,T[\cdot x$ for some $T > 1$. 
Let $n_0 \in \N$ be minimal $\geq T$ with $n_0 x\not\in U$. Then 
$n_0 x \in 2U \setminus U$. Now 
$\exp(n_0 x) \in H \subeq \exp(U)$ contradicts the injectivity 
of $\exp$ on $2U$. This implies that $\R x \subeq U$. 
We conclude that  $\tilde H$ is contained in the largest linear 
subspace $E$ of $\g$ contained in $U$. 

If $\g$ possesses a continuous norm, we may choose $U$ in such a way 
that it does not contain any non-zero subspace. 
Then $\tilde H = \{0\}$ implies $H = \{\be\}$, so that 
$G$ has no small subgroups. 

Suppose, conversely, that $G$ has no small subgroups. Then we may assume that 
$\exp(U)$ contains no non-trivial subgroup. If $\R x \subeq U$, 
then $\exp(\R x)$ is a subgroup contained in $\exp(U)$, so that $x = 0$. 
Now Lemma~\ref{lem:8.11.2} implies that $\g$ possesses a continuous norm.
\end{prf}

\begin{rem} The preceding argument shows that, for any locally 
exponential Lie group and $U$ as in the proof of Proposition~\ref{prop:8.11.3}, 
any subgroup $H \subeq \exp(U)$ is contained in $\exp(E)$, 
where $E \subeq U$ is the maximal linear subspace contained in $U$. 
As $\exp\res_E \: E \to \exp(E)$ is a homeomorphism, this shows that 
$H$ contains no non-trivial compact subgroup. In particular, it contains 
no torsion elements.
\end{rem}

\section{Quotients of locally exponential Lie groups} 
\mlabel{sec:6.4} 

In Section~\ref{sec:6.3} we introduced the concept of a locally exponential 
subgroup and we have seen that kernels of morphisms of 
locally exponential Lie groups are closed locally exponential subgroups. 
In this section, we turn to the converse situation, i.e., 
we study under which conditions on a closed normal subgroup $N$ 
of a locally exponential Lie group $G$, 
the quotient group $G/N$ is a locally exponential Lie group. 
The key result is the Quotient Theorem~\ref{thm:5.5.11}. 
We shall also discuss 
an interesting connection between the Lie subgroup property of integral 
subgroups corresponding to a closed ideal $\fn$ and the enlargeability of the 
quotient Lie algebra $\g/\fn$.

\subsection{Quotient groups} 

\index{Quotient Theorem for locally exponential groups}
\begin{thm}[Quotient Theorem for locally exponential groups]  \mlabel{thm:5.5.11} 
Let $G$ be a locally exponential Lie group with Lie algebra $\g$, 
$N \trile G$ a closed normal locally exponential Lie subgroup and 
$\fn$ its Lie algebra.
Then the quotient group $G/N$ is a locally exponential Lie group if and only 
if the Lie algebra $\g/\fn$ is locally exponential. 
\end{thm}

\begin{prf} If the topological group $G/N$ is a locally exponential Lie group, 
then the quotient map $q \: G \to G/N$ is a smooth morphism of Lie groups and 
$\L(q) \: \g = \L(G) \to \L(G/N)$ also is a quotient morphism of Lie algebras with 
kernel $\fn = \L(N)$ (Theorem~\ref{thm:5.4.4}). Therefore 
$\L(G/N) \cong \g/\fn$, and since $G/N$ was assumed to be locally exponential, 
the Lie algebra $\g/\fn$ is locally exponential. 

Suppose, conversely, that $Q \subeq \g/\fn$ is an exponential local Lie group 
and write $q \: \g \to \g/\fn$ for the quotient map. 
Further, let $U \subeq \g$ be an exponential local Lie group 
for which $\exp_G\res_U$ is a diffeomorphism onto an 
open subset of~$G$, 
$$ \exp_G(U) \cap N = \exp_G (U \cap \fn), $$
and $q(x*y) = q(x) * q(y)$ for $x,y \in U$ (Theorem~\ref{thm:5.2.11c}). 
In view of Theorem~\ref{thm:5.3.5}, we may further assume that 
$\kappa_\g(x)\fn = (\kappa_U)_x(\fn) = \fn$ holds for all $x \in U$. 

Let $U_1 \subeq U$ be a convex balanced open $0$-neighborhood with 
$U_1 \times U_1 \subeq D_U$. 
For $x,x' \in U_1$ with $x' - x \in \fn$, we then have 
$$ q((-x)*x') = q(-x)*q(x') = (-q(x)) * q(x) = 0, $$
so that $(-x)*x' \in \ker q = \fn$. 
This implies that 
$$\exp_G(x)^{-1}\exp_G(x') = \exp_G((-x)*x') \in \exp_G \fn \subeq N,$$ 
and hence that $\exp_G x' \in \exp_G x\cdot N$. 

We thus obtain a well-defined map 
$$ \psi \: V_1 := q(U_1) \to G/N, \quad  
q(x) \mapsto (\exp_G x) N $$
satisfying 
\begin{align*}
\psi(q(x)*q(y)) 
&= \psi(q(x*y)) = \exp_G(x*y) N \\
&= \exp_G x \exp_G y N = \psi(q(x)) \psi(q(y)),
\end{align*}
i.e., $\psi \: (V_1, *) \to G/N$ is a homomorphism of local groups. 

Moreover, $\psi$ is open because this holds for $\exp_G\res_{U_1}$ and 
the quotient map $G \to G/N$. 
If $\psi(q(x)) = \psi(q(y))$ holds for $x,y \in U_1$, then 
$$ \exp_G((-y)*x) = \exp_G(y)^{-1} \exp_G x \in N \cap \exp_G(U) \subeq \exp_G(\fn), $$
which leads to $(-y)*x \in \fn$ and hence to 
$$ 0 = q((-y)*x) = q(-y)*q(x) = (-q(y)) * q(x). $$
Thus $q(y) = q(x)$ 
follows from the uniqueness of inverses in the local group $Q$.
This implies that $\psi$ is open and injective, hence a homeomorphism onto its image. 

Multiplication and inversion in $V_1 \subeq Q$ are smooth. 
Moreover, we have, for each $g \in G$,  
on a sufficiently small neighborhood of $0$ in $\g$ the relation  
$$ c_{gN} \circ \psi = \psi \circ \Ad_{\g/\fn}(g), $$
which implies that the conjugation with $x$ in $V_1$ is given in a neighborhood of 
$0$ by the continuous linear map $\Ad_{\g/\fn}(\exp_G x)$, hence smooth. 
Now Theorem~\ref{thm:locglob} implies the existence of a unique Lie group 
structure on $G/N$ for which $\psi$ defines a chart in a neighborhood of $\be$. 

For each $x \in \g$, the one-parameter group 
\[ \gamma_{q(x)} \: \R \to G/N, \quad t \mapsto \exp_G(tx) N \] 
satisfies for sufficiently small $t$ the relation 
$\gamma_{q(x)}(t) = \psi(t q(x)),$
hence is smooth with $\gamma_{q(x)}'(0) = q(x)$, where we identify 
$\g/\fn$ with $T_\be(G/N)$ via $T_0(\psi)$. 
This shows that the Lie group $G/N$ has an exponential function  
given by 
$$ \exp_{G/N}(q(x)) = \exp_G(x)N. $$
Since $\psi \circ q  = \exp_{G/N} \circ q$ holds on $U_1$, 
$\exp_{G/N}$ is smooth and the group $G/N$ is locally exponential. 
\end{prf}

\begin{cor} \mlabel{cor:5.5.12} 
Let $N \trile G$ be a closed 
normal locally exponential Lie subgroup and $\fn$ be its Lie algebra.
If $N$ is central or contains the commutator group, then 
the quotient group $G/N$ is a locally exponential Lie group. 
\end{cor}

\begin{prf} In view of the Quotient Theorem~\ref{thm:5.5.11}, we only have to 
show  that $\g/\fn$ is locally exponential if $N$ is central or contains the 
commutator algebra. If $N$ is central, then $\fn \subeq \z(\g)$ 
and $\g/\fn$ is locally exponential by Corollary~\ref{cor:5.3.6}. 
If $N$ contains the commutator group, then it also contains its closure, 
so that $[\g,\g] \subeq \fn$ by Lemma~\ref{lem:5.5.10c}, and we can apply 
Corollary~\ref{cor:5.3.6} again. 
\end{prf}

\begin{cor} \mlabel{cor:5.5.13b} 
Let $G$ be a Banach--Lie group and $N \trile G$ a closed normal subgroup. 
Then the quotient group $G/N$ carries a Banach--Lie group structure if and only if 
$N$ is a Banach--Lie subgroup. 
\end{cor}

\begin{prf} This follows immediately from Theorem~\ref{thm:5.5.11} because 
the quotient Lie algebra $\g/\fn$ is Banach, hence automatically locally exponential. 
\end{prf}

In analogy to the preceding corollary, 
it would be nice to have a similar criterion for general subgroups:

\begin{conj} \label{conj}
A closed subgroup $H$ of a Banach--Lie group $G$ is a Banach--Lie
subgroup if and only if $G/H$ carries the structure of a Banach
manifold for which the quotient map $q \: G \to G/H, g \mapsto gH$
has a surjective differential at each point and $G$ acts smoothly on
$G/H$.
\end{conj}

The main difficulty arises from the possible absence of closed
complements to the closed subspace $\L(H)$ in $\L(G)$. If $H$ is a
Banach--Lie subgroup of $G$ and $\L(H)$ has a closed complement,
then the classical Inverse Function Theorem provides natural charts
on the quotient space $G/H$, hence a natural manifold structure with
all desired properties (Propositions~\ref{prop:split-Lie} and 
\ref{prop:split-crit-ban}). Without a closed
complement for $\L(H)$ it is not known how to construct charts of
$G/H$. However, the natural model space is the quotient space
$\L(G)/\L(H)$. With the version of the Implicit Function Theorem 
obtained in \cite{AnN09},  we can now prove one half of the conjecture above:

\begin{thm} \label{thm:liesub} Let $G$ be a Banach--Lie group and
$\sigma \: G \times M \to M$ a smooth action of $G$ on a Banach
manifold $M$, $p\in M$. Suppose that for the derived action
$$ \dot\sigma \: \L(G) \to {\cal V}(M), \quad
\dot\sigma(x)(m) := -d\sigma(\be,m)(x,0) $$ the subspace
$\dot\sigma(\L(G))(p)$ of $T_p(M)$ is closed. Then the stabilizer
$$ G_p := \{ g \in G \: g.p = p\} $$
is a Lie subgroup.
\end{thm}

\begin{prf} First we note that $G_p$ is a closed subgroup of
$G$ with
$$ \L(G_p) = \{ x \in \L(G) \: \dot\sigma(x)(p) =0\}. $$
We consider the smooth maps
$$ \psi \: G \to M, \quad g \mapsto g.p
\quad \mbox{ and } \quad
\phi \: \L(G_p) \to G_p, \quad x \mapsto \exp_G(x). $$
Then $G_p = \psi^{-1}(p)$, $\phi(0)= \be$,
$$ \im(d\phi(0)) = \L(G_p) = \ker(d\psi(\be)). $$
Since $\im(d\phi(0))=\dot\sigma(\L(G))(p)$ is closed, 
the Implicit Function Theorem \cite[Thm.~2.1]{AnN09} 
implies that there exists an open 
$0$-neighborhood $U \subeq \L(G_p)$, which can be chosen to be
arbitrarily small, such that $\exp_G(U)$ is a $\be$-neighborhood in
$G_p$. By Lemma~\ref{lem:equiv-locexp-Liesub}, 
$G_p$ is a Banach--Lie subgroup of~$G$.
\end{prf}

\begin{cor} Let $G$ be a Banach--Lie group and let 
$H \subeq G$ be a closed subgroup for which $G/H$ carries the structure
of a Banach manifold. Suppose the quotient map $q \: G \to G/H, g
\mapsto gH$ has a surjective differential at some point $g_0\in G$
and $G$ acts smoothly on $G/H$. Then $H$ is a Banach--Lie subgroup
of~$G$.
\end{cor}

\begin{prf} The action of $G$ on $G/H$ is given by
$\sigma(g,xH) = gxH = q(gx),$
so that
$$ \dot\sigma(\L(G))(gH) = dq(g)(\L(G)) = T_{gH}(G/H) $$
and the stabilizer of $p := gH = q(g) \in G/H$ is $gHg^{-1}$. By
Theorem~\ref{thm:liesub}, $g_0 Hg_0^{-1}$ is a Banach--Lie subgroup. So
$H$ is a Banach--Lie subgroup.
\end{prf}

\begin{rem} \label{rem:more-gen}
The preceding theorem has a natural generalization to the following
setting. Let $H$ be a closed subgroup of a Banach--Lie group $G$.
Suppose that there exists an open $\be$-neighborhood $U_G \subeq G$
and a $C^1$-map $F \: U_G \to M$ to some Banach manifold $M$ such
that $T_\be(F)$ is surjective and that 
\[ F^{-1}(F(\be)) \cap U_G = H \cap U_G,\] 
then $H$ is a Banach--Lie subgroup of $G$.
\end{rem}

\begin{probl} Suppose that the Banach--Lie group $G$ acts smoothly
on the Banach manifold $M$ and $p \in M$. Is it always true that the
stabilizer $G_p$ is a Lie subgroup of $G$?

From the Algebraic Subgroup Theorem (Theorem~\ref{thm:5.5.14}), 
it follows that, 
for each Banach space $Y$, for each $y \in Y$ and each closed subspace
$X \subeq Y$, the subgroups
\[  \GL(Y)_y := \{ g \in \GL(Y) \: g(y)=y\},  \quad 
 \GL(Y)_X := \{ g \in \GL(Y) \: g(X)=X\} \] 
are Banach--Lie subgroups. Since inverse images of Lie subgroups are
Lie subgroups (Proposition~\ref{prop:5.5.9}),   the problem
has a positive solution for linear actions and for stabilizers of
closed subspaces.

If $X$ is not complemented in $Y$, then it is not clear how to turn
the Grassmannian $\Gr_X(Y) := \{ g.X \: g \in \GL(Y) \}$ into a
smooth manifold on which $\GL(Y)$ acts. We are therefore in a
situation where we know that the closed subgroup $\GL(X)_Y$ is a
Banach--Lie subgroup, but nothing is known about the quotient
$\Gr_X(Y) \cong \GL(Y)/\GL(Y)_X$. Therefore a solution of
Conjecture~\ref{conj} would in particular lead to a natural Banach
manifold structure on the Grassmannian of all closed subspaces of
$Y$.
\end{probl}

\subsection{Integrability of quotient Lie algebras} 

In this subsection we combine the Quotient Theorem and the 
fact that kernels of morphisms of locally exponential Lie groups 
are locally exponential subgroups to derive that for a closed ideal 
$\fn$ in the Lie algebra 
$\g$ of a $1$-connected locally exponential Lie group $G$,  
the quotient $\g/\fn$ (if assumed to be locally exponential) 
is enlargeable if and only if $\exp_G\fn$ generates a closed 
locally exponential Lie subgroup of $G$. This observation has a variety of applications, 
one of which is the Integrability Theorem for Quotients stated below. 
This topic will be taken up in the more extensive discussion of 
integrability in Section~\ref{sec:11.4}.

\begin{prop} \mlabel{prop:int-quot-locexp} 
Let $G$ be a $1$-connected locally exponential Lie group with Lie algebra 
$\g$ and $\fn \trile \g$ a closed ideal for which 
the quotient Lie algebra $\g/\fn$ is locally 
exponential. Then $\g/\fn$ is enlargeable 
if and only if the normal subgroup $N := \la \exp_G\fn\ra \trile G$ is a 
locally exponential Lie subgroup.   
\end{prop}

\begin{prf} If $\fq := \g/\fn$ is integrable and 
$Q$ a corresponding connected Lie group, then there exists a 
homomorphism $q \: G \to H$ for which $\L(q) \: \g \to \fq$ 
is the quotient map (Theorem~\ref{thm:5.4.4}). 
Hence $\ker q$ is a normal locally exponential Lie subgroup 
of $G$ with Lie algebra $\fn$ (Proposition~\ref{prop:5.5.9}). In particular 
$N := \la \exp \fn \ra = (\ker q)_0$ is a closed 
locally exponential Lie subgroup of $G$. 

If, conversely, $N$ is a closed locally exponential 
Lie subgroup of $G$, then the Quotient 
Theorem~\ref{thm:5.5.11} 
implies that the topological group $G/N$ carries a Lie group 
structure with $\L(G/N) \cong \g/\fn$.   
\end{prf}

\begin{cor} \mlabel{cor:int-quot-ban} Let $G$ be a $1$-connected Banach--Lie group 
and $\fn \trile \g$ a closed ideal. Then the 
quotient Lie algebra $\g/\fn$ is enlargeable if and only if 
the normal subgroup $N := \la \exp_G\fn\ra \trile G$ is a 
closed locally exponential Lie subgroup.   
\end{cor}

\begin{prf} We only have to recall that each quotient $\g/\fn$ 
is Banach, hence in particular locally exponential 
(Proposition~\ref{prop:5.1.3}). Then 
Proposition~\ref{prop:int-quot-locexp} applies. 
\end{prf}

\begin{thm} [Enlargeability Theorem for Quotients] \mlabel{thm:quot-alg} 
Let $G$ be a $1$-connected locally exponential 
Lie group with Lie algebra $\g$ and $\fn \trile \g$ a closed 
ideal for which the quotient Lie algebra $\fq := \g/\fn$ is locally exponential. 
We write $N := \la \exp_G\fn \ra$ for the integral subgroup 
generated by $\exp_G\fn$ and 
$$Z(G,\fn) := \{ g \in G \: (\Ad(g) - \be)(\g) \subeq \fn\} $$ 
for the kernel of the natural representation of $G$ on $\fq$. 
Then $Z(G,\fn) \trile G$ is a normal locally exponential Lie subgroup with 
Lie algebra 
$$\z(\g,\fn) := \{ x \in \g \: [x,\g] \subeq \fn\}$$ 
and the Lie algebra homomorphism 
$q \: \z(\g,\fn) \to \z(\fq)$ defines a period homomorphism 
$$ P_\fq \: \pi_1(Z(G,\fn)) \to \z(\fq), \quad 
P_\fq([\gamma]) = \int_0^1 q(\delta(\gamma)_t)\, dt, $$
where $\gamma \:[0,1] \to Z(G,\fn)$ is a piecewise smooth loop. 
The Lie algebra $\fq$ is enlargeable if and only if the image of $P_\fq$ is a discrete 
subgroup of $\z(\fq)$. 
\end{thm}

\begin{prf} First we use the Integrability Theorem~\ref{thm:5.3.8}  
to see that the adjoint group $Q_{\rm ad} := \la e^{\ad \fq} \ra \subeq \Aut(\fq)$ 
carries the structure of a locally exponential Lie group 
with Lie algebra $\fq_{\ad} := \fq/\z(\fq) \cong \g/\z(\g,\fn)$.  
Let $\Ad_\fq \: G \to Q_{\rm ad}$ denote the natural representation of 
$G$ on $\fq$ and note that 
$$ \Ad_\fq \circ \exp_G = \exp_{Q_{\rm ad}} \circ \ad_{\fq} $$
implies that $\Ad_{\fq}$ is a morphism of locally exponential Lie groups. 
Therefore Proposition~\ref{prop:5.5.9} implies that 
$Z(G,\fn) = \ker \Ad_{\fq}$ is a normal locally exponential Lie subgroup of $G$ 
whose Lie algebra is $\ker \ad_{\fq} = \z(\g,\fn)$. 

Integrating the Lie algebra homomorphism 
$\z(\g,\fn) \to \z(\fq)$ to a group homomorphism 
$\tilde P_\fq \: Z(G,\fn)\,\tilde{}_0\to \z(\fq)$ of the corresponding universal 
covering group, we obtain the homomorphism 
$$ P_\fq \: \pi_1(Z(G,\fn)) \to \z(\fq) $$
by restricting to the discrete central subgroup $\pi_1(Z(G,\fn))$ 
(Theorem~\ref{thm:5.4.4}). The explicit formula for $P_\fq$ can be verified 
as follows: Let $\tilde \gamma \: [0,1] \to Z(G,\fn)\,\tilde{}_0$ be a 
lift of a loop $\gamma$ with $\tilde\gamma(0) = \be$. Then 
$$ P_\fq([\gamma]) 
= \tilde P_\fq(\tilde\gamma(1)) 
= \int_0^1 (\tilde P_\fq \circ \tilde\gamma)'(t)\, dt 
= \int_0^1 q(\delta(\tilde\gamma)_t)\, dt 
= \int_0^1 q(\delta(\gamma)_t)\, dt $$
(cf.\ Theorem~\ref{thm-fundamental}). 

If $\Gamma := \im(P_\fq)$ is a discrete subgroup of 
$\z(\fq)$, then the quotient group \break 
$Z := \z(\fq)/\Gamma$ carries a natural Lie group structure 
and the homomorphism 
$$Z(G,\fn)\,\tilde{}_0 \to \z(\fq)$$ factors through 
a homomorphism 
$\beta \: Z(G,\fn)_0 \cong Z(G,\fn)\,\tilde{}_0/\pi_1(Z(G,\fn)) \to Z$ for which 
$\L(\beta) \: \z(\g,\fn) \to \z(\fq)$ is the quotient map 
with $\ker(\L(\beta)) = \fn$. We conclude that 
$\ker \beta$ is a normal locally exponential Lie subgroup of $Z(G,\fn)_0$ with Lie 
algebra $\fn$, and hence that $N$ is a closed locally exponential Lie subgroup of $G$. 
Now $\fq$ is enlargeable by Proposition~\ref{prop:int-quot-locexp}. 

Suppose, conversely, that $\fq$ is  enlargeable, and hence 
that $N$ is a closed locally exponential 
subgroup of $G$ (Proposition~\ref{prop:int-quot-locexp}), 
which implies that $N$ is a locally exponential subgroup of $Z(G,\fn)$. 
In view of the Quotient Theorem~\ref{thm:5.5.11}, 
$Z := Z(G,\fn)_0/N$ is a locally exponential Lie group with Lie algebra 
$\L(Z) = \z(\g,\fn)/\fn \cong \z(\fq)$, hence a quotient 
$\z(\fq)/\Gamma$ for some discrete subgroup 
$\Gamma \cong \pi_1(Z)$ of $\z(\fq)$ (Proposition~\ref{prop:5.1.3b}). 
The homomorphism of 
Lie groups $Z(G,\fn)_0 \to Z$ induces 
a homomorphism $\pi_1(Z(G,\fn)) \to \pi_1(Z)\cong \Gamma$ which 
coincides with $P_\fq$. Hence $\im(P_\fq) \subeq \Gamma$ is a discrete 
subgroup of $\z(\fq)$. 
\end{prf}

\begin{rem} We shall see in Section~\ref{sec:11.4} that, 
for each Mackey complete locally exponential 
Lie algebra $\fq$ (for which $\fq_{\rm ad}$ is also Mackey complete), 
there is a natural homomorphism 
$$ \per_\fq \: \pi_2(Q_{\rm ad}) \to \z(\fq) $$
whose image is discrete if and only if $\fq$ is enlargeable. 

Assuming that $\g$ is Fr\'echet, we have 
in the situation of the preceding theorem  
$Q_{\rm ad} \cong G/Z(G,\fn)$, and with Michael's Selection Theorem 
(\cite{MicE59}, \cite{Bou87}) it follows that $G$ is a topological $Z(G,\fn)$-principal bundle. 
We therefore have a connecting homomorphism 
\[ \delta \: \pi_2(Q_{\rm ad}) \to \pi_1(Z(G,\fn)).  \] 
The surjectivity of $\delta$ follows from the $1$-connectedness of $G$ 
and the exactness of the long exact homotopy sequence of the bundle 
$G \to Q_{\rm ad}$ (Theorem~\ref{homseq-princ}). 
Then it is not hard to see that 
$$ P_\fq \circ \delta  = \per_{\fq}, $$
which makes the two pictures fit. 
\end{rem}

\section{Regularity of Lie group extensions} \mlabel{sec:6.5}

The main result of this section is 
Theorem~\ref{thm:reg-ext-prop}, asserting that an extension 
$\hat G$ of a Lie group $G$ by a Lie group $N$ is regular 
if and only  if $G$ and $N$ are regular. This regularity criterion 
is an important tool to verify regularity of Lie groups.

\subsection{Regularity of initial Lie subgroups} 

Before we turn to the proof, we need some criteria for subgroups 
of regular Lie groups to be regular:

\begin{prop} \mlabel{prop:reg-crit1} Let 
$G$ be a ${C^k}$-regular Lie group with Lie algebra $\g$ and 
$H \leq G$ an initial Lie subgroup with Lie algebra $\h \subeq \g$, for which 
there exists an open 
identity neighborhood $U \subeq G$ and a smooth function 
$f \: U \to E$ into some locally convex space $E$, such that 
$f$ is constant on $U \cap gH$ for each $g \in U$, and $H \cap U = f^{-1}(0) \cap U$. Then 
$H$ is $C^k$-regular.   
\end{prop}

\begin{prf} 
The ${C^k}$-regularity of $G$ implies the existence of a smooth evolution map 
$$ \evol_H \: {C^k}(I,\fh) \to G, \quad \xi \mapsto 
\gamma_\xi(1), $$
and since $H$ is initial, it suffices to see that the range 
of this map lies in $H$. 

If $\xi \in {C^k}(I,\fh)$ such that 
$\im(\gamma_\xi) \subeq U$, then  for every $t\in I$, 
$$ (f \circ \gamma_\xi)'(t) 
= \dd f({\gamma_\xi(t)})\gamma_\xi'(t) = 0 $$
because $\gamma_\xi'(t) = \dd\ell_{\gamma_\xi(t)}(\be)\xi(t)$ is 
the derivative of a curve in the set $U \cap \gamma_\xi(t)H$, on which 
$f$ is constant. 
Therefore $f \circ \gamma_\xi$ is constant, which leads to 
$\im(\gamma_\xi) \subeq f^{-1}(f(\be)) = f^{-1}(0) = U \cap H$. 

If $\xi \in {C^k}(I,\fh)$ is arbitrary, we apply the 
preceding argument to the curves 
$t\mapsto \gamma_\xi(t_0)^{-1} \gamma_\xi(t)$
on sufficiently small intervals $[t_0, t_0 + \eps]$ and see that 
$\im(\gamma_\xi)$ is contained in~$H$. 
\end{prf}

\begin{cor} \mlabel{cor:reg-crit2} Let $\sigma \: G \times M \to M$ 
be a smooth action of the regular Lie group $G$ on the 
smooth manifold $M$. Suppose that the stabilizer subgroup 
$G_m$ of some $m \in M$ is an initial Lie subgroup of $G$. Then it is a regular 
Lie group.
\end{cor}

\begin{prf} Let $(\phi,U_M)$ be an $E$-chart of $M$ with $m \in U_M$. 
We then have the smooth map 
$$ f := \phi \circ \sigma^m \: U_G := \{ g \in G \: g.m \in U_M\}
\to E, \quad g \mapsto \phi(g.m), $$
which is constant on the left cosets of $G_m$ in $U_G$. Hence 
Proposition~\ref{prop:reg-crit1} applies. 
\end{prf}

\subsection{Regularity is an extension property}

\begin{lem}\mlabel{lem.e.2} 
Let $M$ be a compact connected smooth 
manifold (possibly with boundary)  and $G$ be a Lie group 
with Lie algebra~$\g$. 
Then the prescription 
\[  \alpha * f := \delta(f) + \Ad(f)^{-1}.\alpha \]
defines a smooth affine right action of the group 
$C^\infty(M,G)$ on $\Omega^1(M,\g)$ and of the groups 
$C^r(M,G)$ on $\Omega^1_{C^{r-1}}(M,\g)$ for $r \in \N$. 
\end{lem}

\begin{prf} 
That $\alpha * f$ defines an action of the group $C^\infty(M,G)$ follows from 
\begin{eqnarray*}
(\alpha * f_1) * f_2 
&=& \delta(f_2) + \Ad(f_2)^{-1}.(\delta(f_1) + \Ad(f_1)^{-1}.\alpha) \\
&=& \delta(f_1f_2) + \Ad(f_1 f_2)^{-1}.\alpha = \alpha * (f_1 f_2). 
\end{eqnarray*}
The smoothness of the action follows from the smoothness of 
$\delta$ (Proposition~\ref{prop:smooth-logder0}) 
and the the corresponding linear action (Lemma~\ref{lem:smoothprop}). 
\end{prf} 

\begin{rem} \mlabel{rem:7.5.7} If $\xi = \delta(\eta)$ for some smooth function 
$h \: M \to G$, then the Product Rule implies that 
$\delta(\eta\gamma) = \xi * \gamma,$ 
so that the action from above corresponds to the right multiplication action 
on the level of group-valued functions. 
\end{rem}

\begin{defn} \mlabel{def:liegrp-ext} 
(Lie group extensions) (a) Let $N, \hat G$ and $G$ be Lie groups. 
A pair $(i,q)$ of morphisms $i \: N \to \hat G$ and $q \: \hat G \to G$ is called a 
\index{Lie group extension} 
{\it Lie group extension} of $G$ by $N$ if 
$i$ is injective, $\ker q = i(N)$, and $q$ is surjective. 
We further require that there exists an open $\be$-neighborhood 
$U\subeq G$ and a smooth map $\sigma \: U \to \hat G$ with $q\circ \sigma = \id_U$ such that 
$U \times N \to \hat G, (g,n) \mapsto \sigma(g) i(n)$ is a diffeomorphism onto an open 
subset of $\hat G$. This implies that 
$N$ is a split Lie subgroup of $G$ in the sense of Definition~\ref{defc.4}. 

\index{Lie group extension!abelian} 
\index{Lie group extension!central} 
\index{Lie group extension!equivalence} 
(b) A Lie group extension is called {\it abelian} if $N$ is abelian and 
{\it central} if $i(N) \subeq Z(\hat G)$. 

(c) We call two extensions 
$N \sssmapright{i_1} \hat G_1 \sssmapright{q_1} G$ and 
$N \sssmapright{i_2} \hat G_2 \sssmapright{q_2} G$ 
of the Lie group $G$ by the 
Lie group $N$ {\it equivalent} if there exists a Lie group morphism 
$\phi \: \hat G_1 \to \hat G_2$ such that the following diagram commutes: 
\[  \begin{matrix} 
 N & \sssmapright{i_1} & \hat G_1 &  \sssmapright{q_1}  & G \cr 
\mapdown{\id_N} & & \mapdown{\phi} & & \mapdown{\id_G} \cr 
 N &  \sssmapright{i_2} & \hat G_2 &  \sssmapright{q_2}  & G. \end{matrix}\] 
It is easy to see that any such $\phi$ is an isomorphism of 
 groups and that its inverse is smooth (Exercise~\ref{exer:equiv-ext}).
Thus $\phi$ is an isomorphism
 of Lie groups, and we obtain indeed an equivalence relation. 
We write $\Ext(G,N)$ for the set of equivalence classes of
Lie groups extensions of $G$ by~$N$. 
\end{defn}

\begin{thm} \mlabel{thm:reg-ext-prop} 
{\rm($C^k$-regularity is an extension property)} Let $q \: \hat G \to G$ 
be an extension of the Lie group $G$ by the Lie group $N$ and $k \in \N_0$. Then 
$\hat G$ is ${C^k}$-regular if and only if $N$ and $G$ are ${C^k}$-regular. 
\end{thm}

\begin{prf} 
\textbf{Step 1.}
We assume that $G$ and $N$ are $C^k$-regular and show that this implies 
the $C^k$-regularity of $\hat G$. 
Since $G$ is $C^k$-regular, the evolution map 
$$\Evol_G \: {C^k}(I,\g) \to C^{k+1}_*(I,G)$$ 
is smooth (for $k < \infty$ this is shown in  \cite[Thm A]{Gl12b} and for 
$k = \infty$ in Proposition~\ref{prop:reg-Mack}). 
Let $U_G \subeq G$ be an open $\be$-neighborhood for which we have a smooth 
section $\sigma \: U_G \to \hat G$ with $\sigma(\be_G) = \be_{\hat G}$ and 
$$ p \: {C^k}(I,\hat\g) \to {C^k}(I,\g), \quad \xi \mapsto q_\g \circ \xi $$
be the projection map. Then 
$ V := p^{-1}(\Evol_G^{-1}(C^{k+1}(I,U_G))) $
is an open $0$-neighborhood in ${C^k}(I,\hat \g)$. 
Further, by Lemma~\ref{ClmapCk}, 
the map 
$$ \Phi \: V \to C^{k+1}_*(I,\hat G), \quad \xi \mapsto \sigma \circ \Evol_G(q_\g \circ \xi) $$
is smooth. For $\xi \in {C^k}(I,\hat\g)$, we find 
\begin{align*}
&\quad q_\g(\xi * \Phi(\xi)^{-1}) 
= q_\g\Big(\Ad(\Phi(\xi))(\xi - \delta(\Phi(\xi)))\Big) \\
&= \Ad(q_G(\Phi(\xi)))\Big( q_\g \circ \xi - \delta(q_G(\Phi(\xi)))\Big) 
= \Ad(q_G(\Phi(\xi)))\big( q_\g \circ \xi - q_\g \circ \xi\big) = 0.
\end{align*}
This means that $\xi * \Phi(\xi)^{-1} \in {C^k}(I,\fn)$. 
Now Lemma~\ref{lem.e.2}, applied to the action of $C^{k+1}(I,\hat G)$ 
on $\Omega^1_{{C^k}}(I,\hat \g) \cong {C^k}(I,\hat\g)$, 
 shows that it depends smoothly on $\xi$. 
We thus obtain a smooth map 
$$ \widehat E \: V \to \hat G, \quad \widehat E(\xi) := 
\evol_N(\xi * \Phi(\xi)^{-1}) \cdot \Phi(\xi)(1). $$
The curve $\gamma := \Evol_N(\xi * \Phi(\xi)^{-1})\Phi(\xi)$ 
satisfies 
$$ \delta(\gamma) 
= \delta(\Phi(\xi)) + \Ad(\Phi(\xi))^{-1} 
\big(\delta(\Phi(\xi)^{-1}) + \Ad(\Phi(\xi))\xi\big) 
= \xi. $$
Therefore $\widehat E$ is a smooth local evolution map for $\hat G$, 
and Lemma~\ref{lem:reg-crit} implies that $\hat G$ is 
$C^k$-regular.

\textbf{Step 2.} Conversely, we show that 
the ${C^k}$-regularity of $\hat G$ implies the ${C^k}$-regularity of $N$ and~$G$. 
To see that $N$ is ${C^k}$-regular, we choose a chart 
$(\phi_G,U_G)$ of $G$ and consider the map 
$$ f := \phi_G \circ q \: q^{-1}(U_G) \to \phi_G(U_G), $$
which is constant on the left cosets of $N$, lying in this set. 
Therefore Proposition~\ref{prop:reg-crit1} implies that $N$ is ${C^k}$-regular because 
it is a submanifold of $\hat G$, hence in particular an initial submanifold. 

To see that $G$ is ${C^k}$-regular, we first choose a continuous linear section 
$\sigma \: \g \to \hat \g$, which induces a continuous linear section 
$$ \sigma_* \: {C^k}(I,\g) \to {C^k}(I,\hat\g), \quad 
\xi \mapsto \sigma \circ \xi. $$
Then, for each $\xi \in {C^k}([0,1],\g)$, the curve 
$\gamma_\xi := q \circ \Evol_{\hat G}(\sigma\circ \xi)$ satisfies 
$$ \delta(\gamma_\xi) = \L(q) \circ \sigma \circ \xi = \xi, $$
so that 
$\evol_G = q \circ \evol_{\hat G} \circ \sigma$ 
is a composition of smooth maps, hence smooth. 
\end{prf}

\begin{exs} (a) If $G$ is a regular Lie group, then its tangent bundle 
$T(G)$ is an extension of $G$ by $\g = \L(G)$. 
In fact, it is a semidirect product 
$T(G) \cong \g \rtimes_{\Ad} G$. If $G$ is regular, then $\g$ is 
Mackey complete (Proposition~\ref{prop:reg-Mack}), 
hence a regular Lie group (Propositions~\ref{prop:reg-Mack} and 
\ref{prop:mc-reg}). Therefore 
the preceding theorem implies that $T(G)$ is also regular. 

(b) If $G$ and $N$ are regular Lie groups and $S \: G \to \Aut(N)$ defines a 
smooth action of $G$ on $N$, then the preceding theorem implies that the 
semidirect product Lie group 
$N \rtimes_S G$ is also regular. 

(c) Let $q \: G_1 \to G_2$ be a covering of Lie groups. Then 
$G_1$ is regular if and only if $G_2$ is regular because 
$G_2 \cong G_1/\ker q$, and the discrete group $\ker q$ is trivially regular. 
\end{exs}

\subsection{Smoothness of group actions on regular Lie groups} 

We now take a closer look at functions depending smoothly on a parameter. 
The main application we are aiming at is the 
smoothness of group actions on regular Lie groups. 

Suppose that $\Lambda$ is a smooth manifold. Then a family 
$(\alpha_\lambda)_{\lambda \in \Lambda}$ of elements of $\Omega^1(M,\g)$ 
is said to be {\it smooth} if the map 
$$ \Lambda \times TM \to \g, \quad (\lambda, v) \mapsto \alpha_\lambda(v) $$
is smooth. 
\index{smooth family of $1$-forms}

\begin{prop}\mlabel{prop:para} 
Suppose that $G$ is regular. 
Let $(\alpha_\lambda)_{\lambda \in \Lambda}$ be a smooth family of 
elements of $\Omega^1(M,\g)$ satisfying the Maurer--Cartan equation and 
suppose that $M$ is $1$-connected.  Pick a point $m_0 \in M$. 
Let $f_\lambda \: M \to G$ be the uniquely determined smooth function 
with $f_\lambda(m_0) = \be$ and $\delta(f_\lambda) = \alpha_\lambda$. 
Then the function 
$$ F \: \Lambda \times M \to G, \quad (\lambda, m) \mapsto f_\lambda(m) $$
is smooth. 
\end{prop}

\begin{prf}
First we observe that we may w.l.o.g.\ assume that $M$ is an 
open convex $0$-neighborhood of 
a locally convex space. We then have to show that the function 
$$F \: \Lambda \times M \to G, \quad (\lambda, m) \mapsto f_\lambda(m) $$
is smooth. 

For $x \in M$ we have $f_\lambda(x) = \evol_G(\xi_{\lambda,x})$ for 
$$ \xi_{\lambda,x}(t) 
= f_\lambda(tx)^{-1}\frac{d}{dt} f_\lambda(tx) 
= f_\lambda(tx)^{-1}\dd f_\lambda(tx)(x) = \alpha_\lambda(tx)(x). $$
Since $(\alpha_\lambda)_{\lambda \in \Lambda}$ is a smooth family of $1$-forms, the map 
$$ \Lambda \times M \to C^\infty([0,1],\g), \quad 
(\lambda,x) \mapsto \xi_{\lambda,x} $$
is smooth, because the map 
$$ \Lambda \times M \times [0,1] \to \g, \quad (\lambda, x,t) \mapsto 
\xi_{\lambda,x}(t) = \alpha_\lambda(tx)(x) $$
is smooth (Proposition~\ref{prop:cartes-closed}). Now 
$F(\lambda,x) = \evol_G(\xi_{\lambda,x})$ 
and the smoothness of $\evol_G$ imply that $F$ is smooth.
\end{prf}

\begin{cor}
  \mlabel{cor:para} 
Let $\alpha \: H \times \g \to \g$ be a smooth action of the Lie group 
$H$ on the Lie algebra $\g$ of the connected regular Lie group $G$ by automorphisms. 
Suppose further that each automorphism $\alpha_h$ integrates to an automorphism 
$\alpha_h^G$ of $G$. Then the corresponding action 
$$ \alpha^G \: H \times G \to G $$
is smooth. 
\end{cor}

\begin{prf}
We consider the action $\alpha^G$ as a smooth family of automorphisms of $G$, 
where $\Lambda = H$ and $M = G$. 
For each $h \in H$ we have 
$$ \delta(\alpha^G_h) = \alpha_h \circ \kappa_G \in \Omega^1(G,\g), $$
where $\kappa_G$ is the left Maurer--Cartan form of $G$. In view of the smoothness of the 
action of $H$ on $\g$, the map 
$$ H \times TG  \to \g, \quad (h,gx) \mapsto (\alpha_h \circ \kappa_G)(gx) 
= \alpha_h(x) = \alpha(h,x) $$
is smooth, so that $(\alpha_h \circ \kappa_G)_{h \in H}$ is a smooth family in $\Omega^1(G,\g)$. 
Now the assertion follows from Proposition~\ref{prop:para}. 
\end{prf}

\section*{Open problems in Chapter~\ref{ch:6}}

In Corollary~\ref{cor:5.5.13b} we have seen that a quotient $G/N$ of a Banach--Lie group 
$G$ by a closed normal subgroup is a Banach--Lie group if and only if 
$N$ is a Banach--Lie subgroup. To obtain a natural manifold structure on 
a homogeneous space $G/H$, we have to assume that $H$ is a split Lie subgroup 
(Proposition~\ref{prop:split-Lie}), a condition which apparently is not necessary 
in the case of normal subgroups. It is an important open problem to understand 
when a homogeneous space $G/H$ of a Banach--Lie group carries a natural manifold 
structure. 

\begin{probl} (Banach transformation groups) 
Let $\sigma \: G \times M \to M$ be a smooth action of a 
Banach--Lie group on the Banach manifold $M$. 
\begin{enumerate}
\item[\rm(1)] Is the stabilizer group 
$G_m := \{ g \in G \: g.m = m\}$ of a point $m \in M$ always a Banach--Lie 
subgroup?  This is true for linear actions by 
Theorem~\ref{thm:5.5.14}. 
\item[\rm(2)] Do the orbits of $G$ in $M$ carry natural manifold structures?  
\end{enumerate}
\end{probl}

\begin{probl} Let $G$ be a Banach--Lie group and $H \subeq G$ be a closed subgroup. 
\begin{description}
\item[\rm(a)] Suppose that $G/H$ is a smooth manifold with submersive 
$q$ and a smooth action of $G$. Does this imply that $H$ is a Lie subgroup of $G$? 
\item[\rm(b)] Characterize those 
Lie subgroups $H$ for which $G/H$ is a smooth manifold. 
\item[\rm(c)] Let $H \subeq G$ be a closed subgroup and $\fh := \L^e(H)$ its 
Lie algebra. 
Then the normalizer $N_G(\fh)$ of $\fh$ is a Lie subgroup 
(Proposition~\ref{prop:5.5.9}  and  Proposition~\ref{prop:5.5.15}). 
Is it true that $\Ad(G)\fh \cong G/N_G(\fh)$ carries a natural manifold structure? 
If $H$ is connected, it is a normal subgroup of $N_G(\fh)$. If $H$ is a Lie subgroup, 
this implies that $N_G(\fh)/H$ carries a Lie group structure and 
therefore a manifold structure.
\end{description}
\end{probl}

\begin{probl} (cf.\ Corollary~\ref{cor:7.3.16}) Show that for a morphism 
$\phi \: G \to H$ from a locally exponential Lie group 
$G$ to a Lie group $H$ with exponential function, the quotient 
$G/\ker \phi$ is a locally exponential Lie group and 
$\phi$ induces a morphism of Lie groups $G/\ker \phi \into H$. 
\end{probl}

\begin{probl} (cf.\ Corollary~\ref{cor:diff-inj}) Show that, if $\phi \: G \to H$ 
is a morphism of Lie groups for which 
$G$ is locally exponential and $\L(\phi)$ is injective, 
then $\ker \phi$ is discrete. 

This problem can be reduced to the case where $G$ is the 
additive group of a locally convex space: 
First we observe that we may w.l.o.g.\ assume that $G$ is connected and put 
$N := \ker \phi$. From Proposition~\ref{prop:5.4.6} we derive that 
$\L^e(N) = \{0\}$, so that all $C^1$-arcs in $N$ are constant. 
Hence, for every $n \in \N$, the map 
$G \to N, g \mapsto c_g(n) = gng^{-1}$, is constant, which means that 
$N$ is central. From Proposition~\ref{prop:5.5.10b} we know that 
$Z(G)$ is a locally exponential Lie group, which implies that 
$Z(G)_0 \cong \z(\g)/\ker(\exp_G\res_{\z(\g)})$. 
Therefore $N$ is discrete if and only if 
$\exp_G^{-1}(N) \cap \z(\g)$ is discrete, 
Hence it suffices to prove the assertion for the morphism 
$$ \phi \circ \exp_G\res_{\z(\g)} \: \z(\g) \to H. $$
\end{probl}

\begin{probl} Suppose that $G$ is a locally exponential 
Lie group and $H \subeq G$ is a closed subgroup. 
Is $\L^e(H)$ always locally exponential? 
\end{probl}

\begin{small}

\subsection*{Exercises for Chapter~\ref{ch:6}}

\begin{exer}
  \mlabel{exer:equiv-ext} 
Let $N \into \hat G_1 \onto G$ and 
$N \into \hat G_2 \onto G$ be Lie group extensions of $G$ by $N$. 
Show that every equivalence $\phi \: \hat G_1 \to \hat G_2$ 
(Definition~\ref{def:liegrp-ext}) 
is an isomorphism of Lie groups, i.e., that $\phi$ is bijective and that 
$\phi^{-1}$ is also smooth. 
\end{exer}

\begin{exer} \mlabel{exer:6.1.2} (Split Lie subgroups) Show that the following condition implies that 
the Lie subgroup $H$ of the Lie group $G$ is a split Lie subgroup: 
There exists an open subset $U$ of some locally convex space $E$, 
a smooth map $\sigma \: U \to G$, and 
an open $\be$-neighborhood $U_H \subeq H$ such that the map 
$$ \mu_0 \: U \times U_H \to G, \quad (x,h) \mapsto \sigma(x) i(h) $$
is a diffeomorphism onto an open subset of $G$. 
Hint: Show first that the corresponding map $\mu \: U \times H \to G$ is a 
local diffeomorphism and then that, for a suitable subset $V \subeq U$ 
(we need $\sigma(V)^{-1}\sigma(V) \subeq U_G$ for 
some $\be$-neighborhood $U_G$ of $G$ 
satisfying $U_G \cap H \subeq  U_H$), the 
restriction of $\mu$ to $V \times H$ is injective. 
\end{exer}

\begin{exer} \mlabel{exer:6.1.4} Every Lie subgroup of a Lie group is closed. 
Hint: Exercise~\ref{exer:6.1.3}. 
\end{exer}


\begin{exer} \mlabel{exer:5.6.7} 
Let $J$ be a set. For a tuple 
$x = (x_j)_{j \in J} \in (\R^+)^J$, we define 
$$ \sum_{j \in J} x_j := \sup\Big\{ \sum_{j \in F} x_j \: F \subeq J\hbox{\rm\ finite}\Big\}. $$
Show that 
$$ \ell^1(J,\R) := \Big\{ x= (x_j)_{j \in J} \: \sum_{j \in J} |x_j| <
\infty\Big\} $$
is a Banach space with respect to the norm  
$\|x\|_1 := \sum_{j \in J} |x_j|.$
Define $e_j \in \ell^1(J,\R)$ by $(e_j)_i = \delta_{ij}$. Show that the
subgroup $\Gamma$ generated by $\{e_j \: j \in J\}$ is discrete. 
\end{exer}

\begin{exer} \mlabel{exer:5.6.8} 
We define a Lie bracket on 
$\g := \R^4$ by 
$$ [x,y] := (x_3 y_4 - x_4 y_3, \sqrt{2}(x_3 y_4 - x_4 y_3), 0,0). $$
Show that: 
\begin{enumerate}
\item[\rm(i)] $\g$ is a $2$-step nilpotent Lie algebra with commutator algebra 
$\R(1,\sqrt{2})$. 
\item[\rm(ii)] $\Gamma := \Z e_1 + \Z e_2$ is a discrete subgroup of the group $(\g,*)$, 
where $x*y := x + y + \frac{1}{2}[x,y]$ denotes the BCH product. 
\item[\rm(iii)] In the Lie group $G := (\g,*)/\Gamma$, the commutator group 
$[G,G]$ is $1$-dimensional, but its closure is a $2$-dimensional torus. 
\end{enumerate}
\end{exer}

The following exercise shows that non-closed commutator groups may be simple 
Lie groups. 
\begin{exer} \mlabel{exer:5.6.9} 
Let $H := \tilde\SL_2(\R)$ denote the universal covering group of $\SL_2(\R)$ 
and recall that $Z(H) \cong \Z$. 
Let $\gamma \: Z(H) \to \T^2$ be an injective homomorphism with dense range and form the group 
$$ G := (\T^2 \times H)/D, $$
where $D := \{ (\gamma(z),z) \: z \in Z(H)\}$ is the graph of $\gamma$. 
Show that: 
\begin{enumerate}
\item[\rm(i)] $G$ is a $5$-dimensional Lie group with Lie algebra 
$\g \cong \R^2 \times \fsl_2(\R)$. 
\item[\rm(ii)] The commutator group $[G,G]$ 
of $G$ is $3$-dimensional and dense. 
\item[\rm(iii)] $[G,G] \cong H \cong \tilde\SL_2(\R)$. 
\end{enumerate}
\end{exer}
  
\end{small}

\section{Notes and comments on Chapter~\ref{ch:6}} 

\nin{\bf \ref{sec:6.1}:} The definition of the differential Lie algebra $\L^d(H)$ of a subgroup 
is modeled after \cite{Lau56} and \cite{vN29}. 

Homogeneous spaces $G/H$ of Lie groups play a particular important 
role in the context of Banach--Lie groups \cite{Bel06}, and they 
arise naturally as orbits of Banach--Lie groups 
in spaces of idempotents in Banach algebras 
(generalized Grassmannians) and state spaces 
(see f.i.~\cite{CIJM19}). 

\nin{\bf \ref{sec:6.2}:} In \cite{Rb97}, Robart gives a criterion for the existence of locally exponential 
integral subgroups 
for a prescribed injective morphism $\alpha \: \h \to \g = \L(G)$. Translated 
into our terminology, it says that 
the Lie algebra morphism $\alpha$ can be integrated to an integral subgroup if 
and only if $\h$ is $\ad$-integrable and 
$\h/\z(\h)$ is integrable to a local exponential Lie group 
isomorphic to $H_{\rm ad} :=\la e^{\ad \h} \ra \subeq \Aut(\h)$. 
To see how Robart's result follows from ours, we observe that if 
$\h$ is locally exponential, then $\h$ is stable and 
the integrability of $\h/\z(\h)$ to $H_{\rm ad}$ follows from 
the Integrability Theorem (Theorem~\ref{thm:5.3.8}). Conversely, 
Robart's assumption on $\h$ implies that $\h/\z(\h) = \L(H_{\rm ad})$ is locally exponential,  
and therefore its central extension $\h$ is locally exponential by 
Theorem~\ref{centext-thm}. 
Thus Robart's assumption is equivalent to the local exponentiality of $\h$ 
that we have used to characterize integral subgroups in Theorem~\ref{thm:5.5.3}. 

That Banach--Lie algebras injecting into enlargeable ones are enlargeable 
(Corollary~\ref{cor:5.5.6}) is due to van Est and Korthagen (\cite{EK64}). 
The existence of integral subgroups for closed Lie subalgebras of 
$\L(G)$ for a Banach--Lie group is due to Maissen (\cite[Satz~12.3]{Ms62}). 

The result that if $G$ is a Banach--Lie group 
and $\h \subeq \g := \L(G)$ a closed separable subalgebra, then 
the integral subgroup $H := \la \exp \h \ra \subeq G$ satisfies 
$$ \L(H) = \{ x \in \g \: \exp(\R x) \subeq H\} = \h, $$
i.e., $\exp \R x \subeq H$ implies $x \in \h$ 
is \cite[Thm.~5.52]{HoM98}. Here it is a special case of 
Theorem~\ref{thm:sep-subgroup1} which also applies to non-closed subalgebras 
and arbitrary locally exponential Lie groups.

\nin {\bf \ref{sec:6.3}:} In the finite-dimensional context, one even has a local version 
of Theorem~\ref{thm:loccpt-subgrp} due to E.~Cartan (\cite[Thm.~3.1.7]{Ta14}). 
Tao's monograph \cite{Ta14} is a very nice reference for the circle 
of ideas concerning Hilbert's Fifth Problem, asserting essentially that 
topological groups which are finite-dimensional topological manifolds are 
Lie groups. Closely related is the result that a locally compact group 
without small subgroups is a Lie group (cf.\ \cite{HoM03}). 
It is remarkable that similar results hold for subgroups 
of $\Homeo_+(\bS^1)$. By \cite[Thm.~4.7]{Gh01}, any locally compact 
subgroup of this group is a Lie group.

\nin {\bf \ref{sec:6.4}:} The Quotient Theorem for locally exponential Lie groups 
has grown out of its Banach version (see Corollary~\ref{cor:5.5.13b}  and \cite{GN03}) 
and its generalization 
to BCH Lie groups (\cite{Gl02a}). 

\nin{\bf \ref{sec:6.5}:} 
Proposition~\ref{prop:reg-crit1} is modeled after \cite[38.7]{KM97}. 
The $C^k$-version also appears in Appendix B of \cite{NS13}. 
%
%
%


\chapter{Continuous inverse algebras and linear 
Lie groups} \mlabel{ch:lingrp} 

Many features of finite-dimensional Lie theory can already 
be studied and developed for the class of matrix groups 
without much loss of generality. Accordingly, the full matrix groups 
$\GL_n(\R)$ are key examples of finite-dimensional Lie groups and many 
general results can be derived from the corresponding 
properties of subgroups of $\GL_n(\R)$ (cf.~\cite{HiNe12} for this approach 
to finite-dimensional Lie groups). 
Matrix groups are much easier to understand than abstract Lie groups 
because they are more concrete. This is due to several simplifications 
provided by the matrix context: the group sits in an associative algebra, 
the Lie bracket is a commutator bracket and the exponential function is given by 
the exponential series of matrices. Moreover, the manifold structure 
on $\GL_n(\R)$ is immediately inherited from $M_n(\R)$ because it is an open subset. 

The natural class of locally convex algebras for  
infinite-dimensional Lie theory are continuous inverse algebras 
(Definition~\ref{defncia}). 
A unital locally convex algebra $\cA$ is called a continuous inverse algebra 
(cia  for short) if its unit group $\cA^\times$ is open and the inversion
$\cA^\times \to \cA, a \mapsto a^{-1}$ is continuous. 

For a non-unital locally convex algebra $\cA$,
there is an analogous concept of a continuous quasi-inverse algebra
(see Definition~\ref{defn:cqia}).
We shall see that
$\cA$ has the latter property if and only
if its
unitization 
$\cA_+ = \cA \oplus \K \1$
(as in Exercise~\ref{exer:3.1.4}) is a cia
(see Proposition~\ref{cqia-necessities}).

Returning to a (unital) cia,
we already know from Corollary~\ref{invsmoocia} that the inversion 
on the unit group $\cA^\times$ defines a smooth map, so that $\cA^\times$ 
carries a natural Lie group structure. 
In this chapter we take a closer look at unit groups of cias. In particular, 
we shall show that 
\begin{itemize}
\item[$\bullet$] The inversion of a real cia $\cA$ is real 
analytic, so that $\cA^\times$ actually is an 
analytic Lie group (Theorem~\ref{thm:IV.1.11}). 
This is first verified for 
complex cias and then derived from the fact that the complexification 
$\cA_\C$ of a 
real cia is a complex cia. 
\item[$\bullet$] If $\cA$ is Mackey complete, then a suitable generalization of 
the holomorphic functional calculus for Banach algebras 
provides an analytic exponential and logarithm function. We derive 
that $\cA^\times$ is locally exponential and even a BCH--Lie group 
(Theorem~\ref{thm:IV.1.11}). We also discuss sufficient conditions for 
regularity. 
\end{itemize}
After this inspection of the Lie theoretic properties 
of unit groups, we turn to various types of examples of cias. 
That matrix groups over cias are Lie groups follows from the basic 
observation that, for any  cia $\cA$, the matrix algebra $M_n(\cA)$ 
also is a cia (Proposition~\ref{prop:12.2.3}). 
We conclude this chapter with a discussion of the subalgebra 
$\cA^\infty$ of smooth vectors for a continuous action of a Lie group 
$G$ on a cia $\cA$. The main results are sufficient criteria for 
$\cA^\infty$ to be a cia and for the regularity of its unit group 
(Section~\ref{sec:cia-smoothvec}).\medskip


\noindent
{\bf Conventions for Chapter~\ref{ch:lingrp}.}
In this chapter, we consider algebras
over a ground field $\K\in\{\R,\C\}$.
All algebras considered are assumed to be
associative $\K$-algebras.
They need not have a unit
element, unless they are are referred to as unital
algebras.

\section{Basic facts on cias} 
\mlabel{sec:cia-1}

The following lemma, which follows from Exercise~\ref{exer:I.2.}, 
sometimes simplifies the verification that a 
given locally convex algebra $\cA$ is a cia. 
 
\begin{lem} \mlabel{lem:10.1.1} 
If $\cA$ is a unital locally convex algebra for which 
$\cA^\times$ is a neighborhood of the unit element $\1$ and 
inversion is continuous in $\1$, then $\cA$ is a cia. 
\end{lem}

\begin{ex} (Banach algebras are cias) Every Banach algebra $\cA$ is a cia. It suffices to verify this claim 
for unital Banach algebras, because $\cA_+$ is a Banach algebra if $\cA$ is.
If $\cA$ is a unital 
Banach algebra and $x \in \cA$ satisfies $\|x\| < 1$, then 
$\1- x \in \cA^\times$ with $(\1-x)^{-1} = \sum_{n = 0}^\infty x^n$. Therefore 
$\cA^\times$ is a $\1$-neighborhood, hence open. Further, the above 
formula shows that inversion is continuous in $\1$. This implies that 
$\cA$ is a cia (Lemma~\ref{lem:10.1.1}).   
\end{ex}

If $\cA$ is a Fr\'echet space, then it is 
even easier to verify the cia property. In this case the openness of the unit 
group $\cA^\times$ already implies the continuity of the inversion. 
Since the open subset $\cA^\times$ 
of the Fr\'echet space is completely metrizable, this 
follows from the following lemma. 

\begin{lem} \mlabel{lem:autocont} 
Let $G$ be a group which is a complete metric space and 
assume that the multiplication is separately continuous. 
Then the inversion in $G$ is continuous. 
\end{lem}

\begin{prf} (\cite[p.~115]{Wa71})  Let $d$ be the metric on $G$. 
We show that 
$x_n \to \be$ implies the existence of a subsequence 
with $x_{n_k}^{-1} \to\be$. This subsequence will be chosen by induction.
Pick any $n_1\in\N$.
Let $k\in \N$ and assume that $n_1, \ldots, n_k$ have been chosen. 
Put $y_j := x_{n_j}$ and $p_k := y_1 \cdots y_k$. Since $x_n \to \be$, 
we can choose $n_{k+1}>n_k$ in such a way that 
\[
d(p_k, p_{k+1}) < 2^{-k-1} \,\mbox{ and }\; 
d(p_ky_s^{-1}, p_{k+1}y_s^{-1}) < 2^{-k-1} 
\mbox{ for all } s\in \{1,\ldots, k\},
\] 
where we put $y_{k+1} := x_{n_{k+1}}$ and $p_{k+1} := p_k y_{k+1}$. 
For each $k$, only finitely many conditions have to be satisfied, 
and these conditions hold if $n_{k+1}$ is sufficiently large. 
Here we use that multiplication in $G$ is separately continuous. 

Let $p := \lim_k p_k$ and $q_s := \lim_k p_k y_s^{-1}$. The existence 
of these limits follows from the completeness of $G$ because our choices 
ensure that both sequences are Cauchy. For all $s,k\in\N$ with $s\leq k$, we now have 
\[
d(p,p_k) < 2^{-k} \quad \mbox{ and } \quad 
d(q_s, p_k y_s^{-1}) < 2^{-k}.
\]
Moreover, 
\begin{align*}
d(p,q_s) 
&\leq d(p, p_s y_s^{-1}) + d(p_s y_s^{-1}, q_s) \\
&= d(p, p_{s-1}) + d(p_s y_s^{-1}, q_s) 
< 2^{-s+1} + 2^{-s} < 2^{-s+2},
\end{align*}
so that $q_s \to p$ for $s \to \infty$. 
Note that $q_s = py_s^{-1}$ by separate continuity of  
multiplication. This shows that 
$y_s^{-1} = p^{-1}(p y_s^{-1}) = p^{-1} q_s \to p^{-1} p = \be.$
\end{prf}

\begin{prop}
Let $\cA$ be a unital topological algebra
which is a Fr\'{e}chet space.
If the unit group $\cA^\times$ is open, 
then the inversion map $\iota\colon \cA^\times\to\cA$ is 
continuous and thus~$\cA$ is a cia.
\end{prop}
\begin{prf} Since $\cA^\times$ is an open subset of the Fr\'echet space $\cA$,
its topology can be defined by a complete metric. 
Hence the assertion follows from Lemma~\ref{lem:autocont}.
\end{prf}
\begin{defn}
Let $\cA$ be a cia and write $\cA[[t]]$ for 
the algebra of 
\index{formal power series} 
{\it formal power series} 
with coefficients in $\cA$. We topologize $\cA[[t]]$ by the bijection 
$\cA[[t]] \to \cA^{\N_0},  \sum_{n= 0}^\infty a_n t^n \mapsto (a_n)_{n \in \N}.$
Then $\cA[[t]]$ is a locally convex topological 
algebra with respect to the product 
\[ \Big(\sum_n a_n t^n\Big) * \Big(\sum_m b_m t^m\Big) 
:= \sum_{n = 0}^\infty \Big(\sum_{j+k= n} a_j b_k\Big) t^n.\]
\end{defn}

\begin{prop} $\cA[[t]]$ is a cia with the unit group 
$\cA[[t]]^\times = \cA^\times + t \cA[[t]]. $  
\end{prop}

\begin{prf} Obviously $\cA^\times$ is contained in the unit group of $\cA[[t]]$. 
To see that each element of the form $a + t b$, $a \in \cA^\times$, is invertible, 
we observe that $a + tb = a(\1 + ta^{-1} b)$, and that for each $c \in \cA$ 
the geometric series $\sum_{n = 0}^\infty t^n c^n$ converges 
to an inverse of $\1 - tc$. This implies that the unit group of 
$\cA[[t]]$ is open. 

To see that the inversion map $\eta_\cA \: \cA[[t]]^\times \to \cA[[t]]$ is continuous, 
we write it as $\eta_\cA(a) = \sum_n \eta_\cA(a)_n t^n$. We have to 
show that each map $$\eta_\cA(a)_n \: \cA[[t]]^\times \to \cA$$ is continuous. 
We argue by induction: the map $\eta_\cA(a)_0 = a_0^{-1}$ is continuous because 
$\cA$ is a cia. 
Now let $n > 0$ and assume that $\eta_\cA(a)_m$ is continuous for $m < n$. Then 
we use $\eta_\cA(a)a = \1$ to derive that 
$$ \eta_\cA(a)_n a_0 + \eta_\cA(a)_{n-1} a_1 + \ldots + \eta_\cA(a)_0 a_n = 0, $$
which leads to 
$$ \eta_\cA(a)_n = -(\eta_\cA(a)_{n-1} a_1 + \ldots + \eta_\cA(a)_0 a_n) a_0^{-1}, $$
showing that $\eta_\cA(a)_n$ is continuous. Therefore $\cA[[t]]$ is a cia. 
\end{prf}

\begin{rem} The subgroup $U := \1 + t \cA[[t]] \subeq \cA[[t]]$ is pro-nilpotent. 
It is the projective limit of the nilpotent groups 
$U_n := U/(\1 + t^n \cA[[t]]) \subeq \cA[t]/(t^n)$. 
\end{rem}

A more general situation arises as follows. Let 
$\cB$ be a unital associative algebra, endowed with a filtration 
$(F_n)_{n \in \N}$ satisfying 
$$F_0 = \cB \quad \mbox{ and } \quad \bigcap_n F_n = \{0\}. $$
Then the relation $F_n F_m \subeq F_{n+m}$ implies that each $F_n$ is an ideal of $\cB$. 

\begin{prop} \mlabel{prop:lim-cia} If all the quotient algebras 
$\cB_n := \cB/F_n$ are cias for which the connecting maps 
$\phi_{m,n} \:\cB_n \to \cB_m$, $m \leq n$, are quotient maps, 
then $\hat \cB := \prolim \cB/F_n$ is a cia. 
\end{prop}

\begin{prf}
The ideal $F_1/F_n \subeq \cB_n$ is nilpotent because all $n$-fold products 
of elements in this ideal vanish. Therefore 
the unit group of $\cB_n$ contains $F_1/F_n$. 
Note that $\ker(\phi_{m,n}) = F_m/F_n$. 
Then $\phi_{1,n}$ induces a homomorphism \break 
{$(\cB_n)^\times  \to (\cB_1)^\times$} with 
kernel $\1 + F_1/F_n$. As $\phi_{1,n}$ is a quotient map, it 
maps $\cB_n^\times$ onto an open subgroup of $\cB_1^\times$, hence 
$(\cB_n^\times)_0$ to $(\cB_1^\times)_0$. 
We conclude that 
$$ \phi_{1,n}^{-1}((\cB_1^\times)_0) = (\cB_n^\times)_0. $$

Let $(b_n) \in \hat \cB$ with $b_1 \in (\cB_1^\times)_0$. Then 
$\phi_{1,n}(b_n) = b_1$ implies that $b_n \in (\cB_n^\times)_0$ is invertible, 
and hence that the sequence $(b_n)$ is invertible in the product 
algebra $\prod_{n \in \N_0} \cB_n$. 
Then $\phi_{m,n}(b_n^{-1}) = b_m^{-1}$ implies that 
$(b_n^{-1}) \in \hat \cB$, and hence that 
$(b_n)$ is invertible in $\hat \cB$. We conclude that 
$$ \hat \cB^\times \supeq \phi_1^{-1}((\cB_1^\times)_0), $$
so that $\hat \cB^\times$ is a neighborhood of $\1$, hence open. 
The continuity of the inversion in $\hat \cB$ follows from the 
continuity of the inversion in the algebras $\cB_n$. 
\end{prf}

\begin{ex} Let $\cB= \K[x_1,\ldots, x_n]$, $\K \in \{\R,\C\}$, 
be the free associative algebra in $n$ generators. 
Then $\cB$ has a natural filtration 
$$F_n := \Spann\{ s_1\cdots s_m \: s_i \in \{x_1,\ldots, x_n\}, m \geq n \}, $$
defined by the degree graduation of $\cB$. 
Each quotient $\cB_m := \cB/F_m$ is finite-dimensional, hence a cia. 
Therefore the algebra 
$\hat \cB := \prolim \cB_n$, which can be identified with the algebra 
$\K[[x_1,\ldots, x_n]]$ of non-commutative formal power series in 
the generators $x_1, \ldots, x_n$, is a cia (Proposition~\ref{prop:lim-cia}). 
\end{ex}

\begin{prop} \mlabel{prop:spec-compact} 
If $\cA$ is a Mackey complete complex 
cia and $a \in \cA$, 
then 
\[ \Spec(a)  := \{ \lambda \in \C \:  a- \lambda \1\not\in\cA^\times \} \] 
is a non-empty compact subset of $\C$ and the resolvent 
\[ r_a \: \C \setminus \Spec(a) \to \cA, \quad 
\lambda \mapsto (\lambda \1-a)^{-1} \] 
is holomorphic. 
\end{prop}

\begin{prf} As $\cA^\times$ is open, for 
$\lambda\in \C$ sufficiently large, 
$\lambda^{-1} a - \1 \in \cA^\times$, so that $\lambda \not\in \Spec(a)$. 
This shows that $\Spec(a)$ is bounded. 
That $\Spec(a)$ is closed follows from the closedness of
$\cA \setminus \cA^\times$ and the continuity of the map 
${\C \to \cA,\lambda \mapsto a - \lambda \1}$. 
It remain to verify that $\Spec(x) \not=\eset$. Let us assume that 
$\Spec(x)$ is empty. Then the resolvent 
\[ r_a  \:  \C \to \cA, \quad r_a(\lambda) := (\lambda\1-a)^{-1} \] 
is a holomorphic map by 
Proposition~\ref{invalongCk}. For $\lambda\not=0$, we have 
$r_a(\lambda) = \lambda^{-1}(\1-\lambda^{-1}a)^{-1}$, where 
the second factor tends to $\1$ for $\lambda \to \infty$. 
Therefore $r$ is bounded, so that Liouville's Theorem leads to 
the absurd conclusion that $r$ is constant. 
\end{prf}

\begin{lem}\label{control-spec}
Let $\cA$ be a complex cia, $a\in \cA$ and $U\sub\C$ be a neighborhood of
$\Spec(a)$. Then there exists a neighbourhood $B$ of~$a$ in~$\cA$
such that $\Spec(b)\sub U$ for all $b\in B$.
\end{lem}
\begin{prf}
We have $\Spec(a)\sub \bD_r^0$ for some $r>0$.
Then 
\[ \Spec(za)=z\Spec(a)\sub \bD_1^0 
\quad \mbox{ for all } \quad z\in \bD_{1/r},\] 
whence $\1-za\in\cA^\times$.
Thus
\[
\{(z\1-a)^{-1}=z^{-1}(\1-z^{-1}a)^{-1}\colon |z|\geq r\}
\sub \bD_{1/r}\{(\1-za)^{-1}\colon z\in\bD_{1/r}\}
\]
is a relatively compact subset of~$\cA$. Moreover, $\{(z\1-a)^{-1}\colon z\in (\C\setminus U)\cap\bD_r\}$
is compact. Hence
\[
C:=\{(z\1-a)^{-1}\colon z\in \C\setminus U\}
\]
is relatively compact in~$\cA$. There is a $0$-neighborhood $Y\sub\cA$
such that $CY\sub \1-\cA^\times$ and thus $\1-CY\sub\cA^\times$.
Set $B:=a+Y$. If $b\in B$, then $b=a+y$ for some $y\in Y$.
If $z\in \C\setminus U$, then
\[ z\1-b=z\1-a-y=(z\1-a)(\1-(z\1-a)^{-1}y)\in\cA^\times (\1-CY)\sub\cA^\times,\]
whence $z\not\in\Spec(b)$. Thus $\Spec(b)\sub U$.
\end{prf}

\begin{prop}\label{autocpx}
If~$\cA$ is a complex
continuous inverse algebra, then inversion
$\iota\!: \cA^\times \to \cA^\times, a \mapsto a^{-1}$ is complex analytic. 
\end{prop}
\begin{proof} If~$\cA$ is a complex continuous
inverse algebra, then
$\iota$ is smooth (Corollary~\ref{invsmoocia}), with
\begin{equation}
  \label{eq:inv-differ}
d\iota(x;y)=-x^{-1}\,y\,x^{-1}\quad \mbox{ for } \quad x\in \cA^\times, y\in \cA,
\end{equation}
showing that $d\iota(x;\cdot)$ is complex linear.
By Theorem~\ref{charcxcompl}, the mapping~$\iota$ is complex analytic.
\end{proof}

\begin{prop}\label{autocxan}
If $\cA$ is a real continuous inverse algebra,
then also its complexification $\cA_\C$ is a continuous
inverse algebra and the inversion ${\iota\!: \cA^\times \to \cA^\times}$ is real analytic.
\end{prop}
\begin{proof} 
Since $\cA$ is a continuous inverse algebra, we find
an open identity neighborhood $U\sub \cA^\times$
and an open zero-neighborhood $V\sub \cA$
such that $\1+(a^{-1}b)^2\in \cA^\times$
for all $a\in U$, $b\in V$.
For $a,b$ as before, we have
\[ a+ib=a(\1+ia^{-1}b) \quad \mbox{ in } \quad \cA_\C,\] 
where $a$ is invertible. We claim that $\1+ia^{-1}b$ is also invertible. 
Abbreviating $c:=a^{-1}b$, we observe that $\1+ic$ and $\1-ic$ commute,
whence $\1+ic$ also commutes with $\1+c^2=(\1+ic)(\1-ic)$
and its inverse.
Now 
\[ (\1+ic)(\1+c^2)^{-1}(\1-ic)=(\1+c^2)^{-1}(\1-ic)(\1+ic)=\1 \] 
shows that $(\1+c^2)^{-1}(\1-ic)$ is the inverse of $\1+ic$.
We have shown that the open set $U+iV$ is contained in $\cA_\C^\times$, 
and
\[
(a+ib)^{-1}=(\1+(a^{-1}b)^2)^{-1}(\1-ia^{-1}b)a^{-1}
\;\;\;\;
\mbox{for all $\;(a,b)\in U\times V$,}
\]
which depends continuously on~$(a,b)$.
By Lemma~\ref{lem:10.1.1}, $\cA_\C$ is a continuous
inverse algebra and so inversion on $(\cA_\C)^\times$ is complex analytic 
by Proposition~\ref{autocpx}. Thus inversion $\iota\!: \cA^\times\to \cA^\times$
has a complex analytic extension to a mapping between open subsets
of $\cA_\C$, which means that~$\iota$ is a real analytic map.
\end{proof}
\begin{lem}\label{Neumann} 
Let $\cA$ be a complex cia and $W\sub \cA$  be a balanced open $0$-neighborhood 
such that $\1 + W \subeq \cA^\times$. 
Then $(\1-x)^{-1}=\sum_{n=0}^\infty x^n$ for all $x\in W$. 
\end{lem}
\begin{proof}
Given $x\in W$, we have $B_r(0)x\sub W$ for some $r>1$,
where $B_r(0)\subeq \C$ is the open disk of radius~$r$ centered at~$0$.
We consider the complex analytic function
\[f\!: B_r(0)\to \cA,\;\;\; f(z):=(\1-zx)^{-1}.\]
From \eqref{eq:inv-differ} we derive 
$f'(z)=-(\1-zx)^{-1}(-x)(\1-zx)^{-1}=xf(z)^2$.
A simple induction gives $f^{(n)}(z)=n!\cdot x^n \cdot f(z)^{n+1}$
for all $n\in \N_0$. 
As we may consider~$f$ as a complex analytic function to the
completion of~$\cA$, Corollary~\ref{analytsingle}
shows that
\[
(\forall z\in B_r(0))\qquad  f(z)=\sum_{n=0}^\infty
\frac{z^n}{n!}f^{(n)}(0)=\sum_{n=0}^\infty z^n x^n,\]
whence $(\1-x)^{-1}=f(1)=\sum_{n=0}^\infty x^n$ in particular.
\end{proof}


\begin{thm}
  \mlabel{thm:IV.1.11}  If $\cA$ is a Mackey complete cia, 
then its unit group $\cA^\times$ is an analytic Lie group 
whose exponential function is locally bianalytic. In particular, 
$\cA^\times$ is a BCH--Lie group. 
\end{thm}

\begin{prf} That the multiplication 
on  $\cA^\times$ is analytic follows from its bilinearity on $\cA$  
and the analyticity of the inversion follows from 
Proposition~\ref{autocxan}. 
This in turn leads to the expansion by the Neumann series 
$(\1 - x)^{-1} = \sum_{n = 0}^\infty x^n$ 
(Lemma~\ref{Neumann}). 
This shows that $\cA^\times$ is an analytic Lie group. 

In view of Example~\ref{exs:5.1.4}, 
$\cA^\times$ is locally exponential with an 
analytic exponential function $\exp \: \cA \to \cA^\times$. 
Likewise, functional calculus 
provides on the open star-like subset 
\[ U := \{ a \in \cA \: \Spec(a) \cap ]{-\infty},0] = \eset\} \subeq \cA^\times \] 
an analytic logarithm function $\log\: U \to \cA$ 
(cf.\ Definition~\ref{def:exp-log}). 
Finally, Lemma~\ref{lem:logexp} implies that $\exp_\cA$ is locally bianalytic. 

In view of Theorem~\ref{thm:IV.1.8}, this implies that $\cA^\times$ is~BCH. 
\end{prf}

\begin{rem} We have seen in \cite{Ne08} that, for any cia $\cA$, 
the set $\Idem(\cA)$ of idempotents carries a 
natural manifold structure. The set $\Gr(\cA) := \{ p\cA \: p \in \Idem(\cA)\}$ 
of corresponding right ideals is called the 
\index{Grassmannian of cia} 
{\it Grassmannian of $\cA$} 
(cf.~\cite{DEG98, DG01}). That it carries a natural manifold structure if $\cA$ 
is a Banach algebra has been shown in \cite{DG01}, and the general case of a cia 
follows from \cite{BerN05}. 

If $E$ is a right $\cA$-module, then $\Emb(E,\cA) \subeq \Hom_\cA(E,\cA)$ is called 
the corresponding 
\index{Stiefel manifold}  
{\it Stiefel manifold}. 
If $E = p\cA$ is a projective right ideal, then 
$\Hom_\cA(E,\cA) \cong \Hom_\cA(p\cA,\cA) \cong \cA p$, and 
$\Emb(E,\cA)$ corresponds to the set of elements $x \in \cA p$ for which 
$\lambda_x \: p\cA \to \cA$ is an embedding, i.e., 
$xa = xpa = 0$ implies $a \in (\1-p)\cA$. 
\end{rem}

\section{Algebras of mappings and algebras of germs} 
\mlabel{sec:cia-2}
 
The algebra $C(K,\C)$ of continuous functions
on a compact topological space $K\not=\emptyset$
is a prime example of a Banach algebra
when endowed with the supremum norm.
Likewise, many function spaces with values in $\C$ or
a given cia provide examples of continuous inverse algebras.
Related examples are algebras of function germs
around a point or a compact set.\medskip

In the next lemma, we endow the algebra $C(X,\cA)$ of continuous $\cA$-valued functions on
$X$ with the compact-open topology;
it coincides with the topology of uniform convergence on compact subsets of
$X$ (see Lemma~\ref{sammelsu}(b)).
\begin{lem}
  \mlabel{lem:A.3}
Let $\cA$ be a locally convex algebra and $X$ a topological space. 
Then the following holds: 
\begin{description}[(D)]
\item[\rm(a)] $C(X,\cA)$ is a locally convex algebra. 
\item[\rm(b)] If $\cA$ is unital and inversion in $\cA$ is continuous, then $C(X,\cA)^\times
=\;\;\;\;$\linebreak
$C(X,\cA^\times)$, and inversion in $C(X,\cA)$ is continuous.
\item[\rm(c)] If $X$ is compact and $\cA$ a cia, then $C(X,\cA)$ is a cia.
\end{description}
\end{lem}

\begin{prf}
(a) By Lemma~\ref{sammelsu}(g), $C(X,\cA)$ is a locally convex space.
Let\linebreak
$\mu\colon \cA\times\cA\to\cA$ be the
algebra multiplication, which is continuous.
For $j\in\{1,2\}$,
let $\pr_j\colon\cA\times\cA\to \cA$, $(x_1,x_2)\mto x_j$ be the projection onto the $j$th component.
By Lemma~\ref{cotopprod}, the map 
\[ \Phi\colon C(X,\cA\times\cA)\to C(X,\cA)\times
C(X,\cA), \quad f\mto (\pr_1\circ f,\pr_2\circ f) \] 
is a homeomorphism.
The algebra multiplication on $C(X,\cA)$ is given by $C(X,\mu)\circ\Phi^{-1}$
and hence continuous, by Lemma~\ref{covsuppo}.
%
%

(b) If $f \in C(X,\cA)^\times$, then each value $f(x)$ is a
unit in $\cA$. If, conversely, $f(X) \subeq \cA^\times$, then the
continuity of the inversion in $\cA$ implies that the function 
$X \to \cA, x \mapsto f(x)^{-1}$ is continuous. Thus 
$C(X,\cA)^\times = C(X,\cA^\times)$. 
Since $\iota\colon\cA^\times\to\cA$, $a\mto a^{-1}$
is continuous, also the map
\[
C(X,\iota)\colon C(X,\cA^\times)\to C(C,\cA), \;\; f\mto \iota\circ f
\]
is continuous, by Lemma~\ref{covsuppo}; this is the inversion map.
%
%

(c) By (b), inversion in $C(X,\cA)$ is continuous. If $X$ is compact, then  
$C(X,\cA^\times)=\lfloor X,\cA^\times\rfloor$ is open in $C(X,\cA)$ and
thus $C(X,\cA)$ is a cia.
\end{prf}

\begin{lem}
  \mlabel{lem:A.5} Let $M$ be a locally convex manifold, $r \in \N_0 
\cup \{\infty\}$ and $\cA$ be a locally convex algebra. 
Then the following assertions hold: 
\begin{description}
\item[\rm(a)] $C^r(M,\cA)$ is a locally convex algebra.
\item[\rm(b)] If $\cA$ is unital and inversion in $\cA$ is continuous, then 
\[ C^r(M,\cA)^\times =\{ f \in C^r(M,\cA) \: f(M) \subeq \cA^\times \} \] 
and inversion in $C^r(M,\cA)$ is continuous.
\item[\rm(c)] If $M$ is compact and $\cA$ a cia, then so is $C^r(M,\cA)$.
\end{description}
\end{lem}

\begin{prf}
(a) By (f) and (g) in Lemma~\ref{sammelsu}, $C^r(M,\cA)$ is a locally convex space.
Let $\mu\colon \cA\times\cA\to\cA$ be the
algebra multiplication, which is continuous
and bilinear and thus smooth.
For $j\in\{1,2\}$,
let $\pr_j\colon\cA\times\cA\to \cA$, $(x_1,x_2)\mto x_j$ be the projection onto the $j$th component.
By Lemma~\ref{prodmapsp}, the map $\Phi\colon C^r(M,\cA\times\cA)\to C^r(M,\cA)\times
C^r(M,\cA)$, $f\mto (\pr_1\circ f,\pr_2\circ f)$ is a homeomorphism.
The algebra multiplication on $C^r(M,\cA)$ is given by $C^r(M,\mu)\circ\Phi^{-1}$
and hence continuous, by Lemma~\ref{Cktoppu}(b).

(b) If $f \in C^r(M,\cA)^\times$, then each value $f(x)$ is a
unit in $\cA$. If, conversely, $f(M) \subeq \cA^\times$, then the
smoothness of the inversion in $\cA$ implies that the function 
$M \to \cA, x \mapsto f(x)^{-1}$ is $C^r$. Thus 
$C^r(M,\cA)^\times = C^r(M,\cA^\times)$. 
Since $\iota\colon\cA^\times\to\cA$, $a\mto a^{-1}$
is smooth, the map $C^r(M,\iota)\colon C^r(M,\cA^\times)\to C^r(M,\cA)$, $f\mto \iota\circ f$
is continuous, by Lemma~\ref{Cktoppu}(b); this is the inversion map.

(c) By (b), inversion in $C^r(M,\cA)$ is continuous. If $M$ is compact, then  
$C^r(M,\cA^\times)=\lfloor M,\cA^\times\rfloor$ is an open subset of $C^r(M,\cA)$
(whose topology is finer than the compact-open topology) and
thus $C^r(M,\cA)$ is a cia.
\end{prf}   

\begin{defn} Let $\cA$ and $\cB$ be unital algebras. We call a unital algebra 
homomorphism 
$\phi \: \cA \to \cB$ 
\index{isospectral homomorphism of algebras} 
{\it isospectral} if 
\[ \phi^{-1}(\cB^\times) = \cA^\times.\] 
\end{defn}

Note that $\cA^\times \subeq \phi^{-1}(\cB^\times)$ does always hold. 
Continuous isospectral morphisms are sometimes useful to verify that 
the unit group $\cA^\times$ is open: If $\phi$ is continuous and isospectral 
and $\cB^\times$ is open, then so is $\cA^\times$. 

\begin{ex} \mlabel{ex:9.2.6a} 
Let $X$ be a locally compact space and consider the continuous inclusion 
\[ \eta \: C_c(M,\C) \to C_0(M,\C)\] 
which obviously exteds to a continuous homomorphism 
\[ \eta_+ \: C_c(M,\C)_+ = C_c(M,\C) + \C 1   \to C_0(M,\C)_+ = C_0(M,\C) + \C 1.\] 
Then $\eta_+$ is isospectral and therefore 
$C_c(M,\C)_+^\times = \eta_+^{-1}(C_0(M,\C)_+^\times)$ is open. 
\end{ex}

\begin{ex}
Let $M$ be a $\sigma$-compact smooth manifold. 
If $M$ is not compact, then it carries unbounded smooth functions~$f$. 
Then $f$ has unbounded spectrum in $C^\infty(M,\R)$, so that 
$C^\infty(M,\R)$ is not a cia (Proposition~\ref{prop:spec-compact}). 
However, we claim that 
\[ C^\infty_0(M,\R) := C^\infty(M,\R) \cap C_0(M,\R) \] 
is a continuous quasi-inverse algebra with 
respect to the topology obtained from the embedding into the product 
space $C^\infty(M,\R) \times C_0(M,\R)$. 

Indeed, we first use Lemma~\ref{lem:A.5}(2) to see that 
 the inclusion 
\[C_0^\infty(M,\R)_+ \into C_0(M,\R)_+ \] 
 is isospectral. 
Hence the openness of the unit group in $C_0^\infty(M,\R)_+$ 
follows from the continuity of the inclusion map. 
The continuity of the inversion is a consequence of the corresponding continuity 
in $C^\infty(M,\R)$ and $C_0(M,\R)$ (Lemma~\ref{lem:A.5}). 
Finally, the continuity of the multiplication 
follows from the continuity of the multiplication in $C_0(M,\R)$ and 
$C^\infty(M,\R)$. 
\end{ex}

\begin{lem}\label{semi-via-ball}
For a seminorm $q\colon \cA\to [0,\infty[$
on an algebra $\cA$,
the following conditions are equivalent:
\begin{description}[(D)]
\item[\rm(a)]
$q$ is submultiplicative,
i.e., $q(xy)\leq q(x)q(y)$ for all $x,y\in \cA$;
\item[\rm(b)]
$\wb{B}^q_1(0)\wb{B}^q_1(0)\sub \wb{B}^q_1(0)$.
\end{description}
\end{lem}
\begin{prf}
The implication (a)$\Rightarrow$(b) is obvious.
If (b) holds, then also (a), as a consequence of Lemma~\ref{like-op-no-mult}.
\end{prf}
If seminorms $q_1$ and $q_2$ on $\cA$ are submultiplicative,
then also their pointwise maximum $q\colon \cA\to[0,\infty[$,
$x\mto\max\{q_1(x),q_2(x)\}$ is submultiplicative.
%
%
%
\begin{defn}
A not necessarily Hausdorff, locally convex
topological algebra $\cA$ is called \emph{locally $m$-convex}
(or locally multiplicatively convex)
if its locally convex vector topology can be defined
by a set $\Gamma$ of submultiplicative
seminorms on~$\cA$.
\end{defn}
%
%
%
\begin{rem}
Locally $m$-convex unital topological algebras need not
have an open unit group.
For example, $C(\R,\C)$ is locally $m$-convex
(the seminorm $f\mto \|f|_K\|_\infty$ being submultiplicative
for each compact subset $K\sub \R$) but $C(\R,\C)$
does not have an open unit group.
\end{rem}
Local $m$-convexity passes to locally convex direct limits of countable
direct systems of Banach algebras (see \cite{AkN96}, \cite{DW97}).
We prove a special case.
\begin{thm}\label{union-cia}
Let $(\cA_n,\|\cdot\|_n)_{n\in\N}$ be a sequence
of unital Banach algebras such that $\cA_n\sub \cA_{n+1}$
for all $n\in \N$ and the inclusion map $\cA_n\to\cA_{n+1}$
is a homomorphism of unital algebras and continuous.
Endow $\cA:=\bigcup_{n\in \N}\cA_n$ with the unique unital
algebra structure making each inclusion map $\cA_n\to\cA$ a
homomorphism of unital algebras.
Then the following holds.
\begin{description}[(D)]
\item[\rm(a)]
The locally convex direct limit topology $\cO$
on $\cA=\dl\cA_n$\vspace{-.7mm} makes the latter a not necessarily
Hausdorff, locally convex topological algebra with an open unit group
and a continuous inversion map $\iota\colon \cA^\times \to\cA$.
\item[\rm(b)]
If $\cO$ is Hausdorff, then $(\cA,\cO)$ is a cia.
\item[\rm(c)]
If each of the inclusion maps $\cA_n\to\cA_{n+1}$
has operator norm $\leq 1$,
then $(\cA,\cO)$ is locally $m$-convex.
\end{description}
\end{thm}
\begin{prf}
For $n\leq m$, let $c_{m,n}$ be the operator norm of
the inclusion $\cA_n\to\cA_m$.

(a)
To see that the algebra multiplication of $\cA$ is continuous
at $(0,0)$ and hence continuous, let $U\sub \cA$ be a $0$-neighborhood.
Then
\[
\bigcup_{n\in \N} \sum_{k=1}^n B^{\cA_k}_{s_k}(0)\sub U
\]
for some sequence $(s_k)_{k\in\N}$ of positive real numbers~$s_k$.
For $2\leq k\in \N$, let $M_k$ be the set of all $(a,b)\in\N^2$ such that
$\max\{a,b\}=k$ and $a\not=b$.
Recursively, we find a sequence $(r_k)_{k\in\N}$ of positive real numbers such that,
for each $k\in \N$,
\begin{description}[(D)]
\item[(i)]
$r_k^2<s_k/2$ holds and
\item[(ii)]
If $k\geq 2$, the $r_k\sum_{(a,b)\in M_k}c_{k,\min\{a,b\}} r_{\min\{a,b\}}<s_k/2$.
\end{description}
Then $V:=\bigcup_{n\in\N}\sum_{k=1}^n B^{\cA_k}_{r_k}(0)$ is a $0$-neighborhood
in~$\cA$ and we show that $VV\sub U$. Let $n\in \N$ and
$x_k,y_k\in B^{\cA_k}_{r_k}(0)$ for $k\in\{1,\ldots,n\}$.
Then
\[
(x_1+\cdots+ x_n)(y_1+\cdots+ y_n)=\sum_{k=1}^n \sum_{\max\{(a,b)=k}x_ay_b,
\]
with $(a,b)\in\N^2$. Let us write $w_k$ for inner sum. Then $w_k\in B^{\cA_k}_{s_k}(0)$
(as required). In fact, $\|w_1\|_1=\|x_1y_1\|_1\leq \|x_1\|_1\|y_1\|_1\leq r_1^2<s_1/2\leq s_1$.
For $k\in\{2,\ldots, n\}$,
\[
\|w_k\|_k\leq \|x_ky_k\|_k+\sum_{(a,b)\in M_k}\|x_ay_b\|_k< s_k/2+s_k/2=s_k
\]
by (i) and (ii), since $\|x_k\|_k\|y_k\|_k\leq r_k^2$ and
\[ \|x_ay_b\|_k\leq
\|x_a\|_k\|y_b\|_k\leq c_{k,\min\{a,b\}}\|x_a\|_a\|y_b\|_b
\leq r_k\,c_{k,\min\{a,b\}} r_{\min\{a,b\}}.\]
Thus
\[
(x_1+\cdots+x_n)(y_1+\cdots+ y_n)\in \sum_{k=1}^n B^{\cA_k}_{s_k}(0)\sub U.
\]
We now show that $\cA^\times$ is a $\1$-neighborhood in~$\cA$ and $\iota$ is continuous
at~$\1$. Then $\cA$ is open in~$\cA$ and $\iota$ is continuous, by
Lemma~\ref{lem:10.1.1}.
To this end, let $U\sub \cA$ be a $\1$-neighborhood.
Then $\1+\bigcup_{n\in\N}\sum_{k=1}^n B^{\cA_k}_{s_k}(0)\sub U$
for some sequence $(s_n)_{n\in\N}$ of positive real numbers.
For $n\in\N$, define
\[
C_n:=1+\sum_{k=1}^{n-1}c_{n,k}s_k.
\]
Let $\ve_n:=\frac{s_n}{C_n}$
and $\delta_n:=\frac{\ve_n}{1+\ve_n}$;
thus
\[
\ve_nC_n\leq s_n\;\;\,\mbox{and}\;\;\,
\frac{1}{1-\delta_n}-1\leq \ve_n.
\]
Finally, let $r_n:=\delta_n/C_n$
and
\[
V_n:=\1 + B^{\cA_1}_{r_1}(0)+\cdots + B^{\cA_n}_{r_n}(0).
\]
Then $C_n\geq 1$, $\delta_n\in \,]0,1[$, and $r_n\in \,]0,1[$.
We claim that
\[
V_n\sub \cA_n^\times
\]
for each $n\in \N$ and
\begin{equation}\label{thsbdd}
V_n^{-1}\sub \1+\sum_{k=1}^n B^{\cA_k}_{s_k}(0).
\end{equation}
If this is true, then
\[
V:=\bigcup_{n\in\N}V_n=
\1+\bigcup_{n\in\N}\sum_{k=1}^nB^{\cA_n}_{r_n}(0)
\]
is an open $\1$-neighborhood in~$\cA$ which is contained in~$\cA^\times$
and $V^{-1}\sub U$, whence $\cA^\times$ is a $\1$-neighborhood
and $\iota$ is continuous at~$\1$.

We prove the claim by induction. If $n=1$,
then $C_1=1$, $\ve_1=s_1$, and $r_1=\delta_1$.
Since $r_1\in \,]0,1]$, for each $a\in V_1\sub B^{\cA_1}_{r_1}(0)\sub B^{cA_1}_1(0)$
we have $a\in \cA_1^\times$ and
\[
\|a^{-1}-\1\|_1\leq \sum_{k=1}^\infty \|a^k\|_1< \frac{1}{1-r_1}-1=\ve_1=s_1
\]
(using Neumann's series), whence (\ref{thsbdd}) holds.
If $n\geq 2$ and the assertion holds for $n-1$ in place of~$n$,
let $x_k\in B^{\cA_k}_{r_k}(0)$ for $k\in\{1,\ldots, n\}$.
Then
\[
a:=\1 + x_1+\cdots+x_{n-1}\in V_{n-1}\sub \cA_{n-1}^\times
\]
and $a^{-1}\in \1+\sum_{k=1}^{n-1} B^{\cA_k}_{s_k}(0)$,
whence $\|a^{-1}\|_n< C_n$.
Thus
\[ \|a^{-1}x_n\|_n< C_n r_n\leq \delta_n<1,\]
showing that $\1-a^{-1} x_n\in\cA_n^\times$
and
\[
\|(\1-a^{-1}x_n)^{-1}-\1\|_n\leq \frac{1}{1-\|a^{-1}x_n\|_n}-1
< \frac{1}{1-\delta_n}-1\leq \ve_n=\frac{s_n}{C_n}.
\]
As a consequence, $a+x_n=a(\1+a^{-1}x_n)\in\cA_n^\times$
and
\begin{eqnarray*}
(a+x_n)^{-1}  &= & (\1-a^{-1}x_n)^{-1}a^{-1}\\
&=&  ((\1-a^{-1}x_n)^{-1}-\1)a^{-1}+a^{-1}\in B^{\cA_n}_{s_n}(0)+V_{n-1}=V_n
\end{eqnarray*}
as $a^{-1}\in V_{n-1}$ by induction and
$\|((\1-a^{-1}x_n)^{-1}-\1)a^{-1}\|_n < \frac{s_n}{C_n}C_n=s_n$.

(b) is an immediate consequence of~(a).

(c) If $U\sub \cA$ is a $0$-neighborhood, then
$V:=\bigcup_{n\in\N}\sum_{k=1}^n B^{\cA_k}_{s_k}(0)\sub U$
for some sequence $(s_k)_{k\in \N}$
of positive real numbers $s_k$.
Then $V$ is an absolutely convex, open $0$-neighborhood in~$\cA$.
After shrinking~$V$, we may assume that
the sequence $(s_k)_{k\in\N}$
is monotonically decreasing and $\sum_{k=1}^\infty s_k\leq \frac{1}{2}$.
We show that $VV\sub V$.
As a consequence, the Minkowski functional
$\mu_V$ of $V$ is submultiplicative (see Lemma~\ref{semi-via-ball})
and we deduce that the set of submultiplicative continuous
seminorms on~$\cA$ defines its vector topology.
Let $n\in\N$ and $x_k,y_k\in B^{\cA_k}_{s_k}(0)$
for $k\in \{1,\ldots,n\}$.
Then
\[
(x_1+\cdots+x_n)(y_1+\cdots+y_n)=\sum_{k=1}^n w_k
\]
with $w_k:=\sum_{\max\{a,n)\}=k}x_ay_b$.
Now $w_k\in B^{\cA_k}_{s_k}(0)$,
as
\[
\|w_k\|_k\leq\sum_{\max\{a,b\}=k}\|x_a\|_k\|y_b\|_k
< s_k\sum_{b=1}^k s_b +s_k\sum_{a=1}^{k-1} s_a < s_k.
\]
This completes the proof.
%
%
%
\end{prf}
We apply Theorem~\ref{union-cia}
to an algebra of germs.
\begin{numba}
Let $E$ be a complex Banach space
and $K\sub E$ be a compact, non-empty subset.
If $U$ and $V$ are open neighborhoods of~$K$ in~$E$,
we call holomorphic functions $f\colon U\to \C$
and $g\colon V\to \C$ \emph{equivalent} if
$f|_W=g|_W$ for some neighborhood $W$ of $K$ in~$U\cap V$.
The equivalence classes $[f]$ are called \emph{germs}
around~$K$ and the algebraic operations of the set $\cO(K)$
of germs are inherited from operations on representatives $f$
(restricted to a common domain).
\end{numba}
\begin{prop} \mlabel{prop:germcia}
For each non-empty compact subset $K$
of a complex Banach space~$E$,
the algebra
$\cO(K)$ of germs of $\C$-analytic complex-valued functions
around $K$ in~$E$
is a cia.
\end{prop}
\begin{prf}
Let $\|\cdot\|$ be a norm on~$E$
defining its topology and $B_r(x)\sub E$ the open ball
of radius $r>0$ around $x\in E$.
Pick a sequence $r_1>r_2>\cdots$
converging to~$0$.
For $n\in \N$, define
\[
U_n:=K+B_{r_n}(0)=\bigcup_{x\in K}B_{r_n}(x).
\]
The supremum norm $\|\cdot\|_\infty$ turns the algebra $BC(U_n,\C)$
of continuous, bounded functions $U_n\to \C$
into a Banach algebra. Let $\cA_n$
be the subalgebra of all $f\in BC(U_n,\C)$
which are complex analytic.
Then $\cA_n$ is closed in $BC(U_n,\C)$ and hence
a Banach algebra. By construction, the restriction maps $\rho_{m,n}\colon
\cA_n\to \cA_m$
have operator norm $\leq 1$ for all $m\geq n$.
If $f\in \cA_n$ and $x\in K$, then $f|_{B_{r_n}(x)}$
is uniquely determined by $f|_{B_{r_m}(x)}$,
by the Identity Theorem. Hence $\rho_{m,n}$ is injective.
The mappings
\[
\rho_n\colon \cA_n\to \cO(K),\;\; f\mto [f]
\]
are injective and they turn $\cO(K)$ into the direct limit of the
vector spaces $\cA_n$ because each $f \in \cO(K)$ is bounded on some~$U_n$. 
We endow $\cO(K)$
with the (not necessarily Hausdorff) locally convex direct limit topology.
Identifying $\cA_n$ with $\rho_n(\cA_n)$,
we can interpret $\cO(K)$ as the union of an ascending sequence
of Banach algebras and Proposition~\ref{union-cia} will complete
the proof once we know that $\cO(K)$ is Hausdorff.
Consider a non-zero element in $\cO(K)$,
say $\rho_n(f)$ with $f\not=0$.
By the Identity Theorem, $f$ must have a G\^{a}teaux derivative
$\delta^k_xf(y)\not=0$
for certain $x\in K$, $k\in \N_0$, and $y\in E$.
The linear map
\[
\lambda_n\colon \cA_n\to\C,\;\; g\mto \delta^k_xg(y)
\]
is continuous for each $n\in \N$.
In fact, there exists $r>0$ such that $x+z y\in U_n$
for each $z \in \C$ with $|z|\leq r$.
By Corollary~\ref{corcxcurve1}(c), we have
\[
\lambda_n(g)=\delta^k_xg(y)=\frac{k!}{2\pi i}\int_{|\zeta|=1} \frac{g(x+\zeta y)}{\zeta^{k+1}}\, d\zeta,
\]
which is a linear $\C$-valued function in $g\in \cA_n$ and continuous
with operator norm
$\leq \frac{k!}{r^k}$.
Since $\lambda_m\circ \rho_{mn}=\lambda_n$
for $m\geq n$, we get a well-defined linear map
\[
\lambda\colon \cO(K)\to \C,\;\; [g]\mto \delta^k_xg(y)
\]
which is continuous as a consequence of Remark~\ref{firstremlcxDL}(d).
Since $\lambda([f])=\delta^k_x f(y)\not=0$,
we can separate $[f]$ and the $0$-element of $\cO(K)$,
whence the topological vector space $\cO(K)$ is Hausdorff.
\end{prf}
\begin{rem}
If $E=\C^n$ in Proposition~\ref{prop:germcia},
then each restriction map $\cA_n\to \cA_{n+1}$
is a compact operator, whence $\cO(K)$
is a Silva space. The Hausdorff property of
$\cO(K)$ is then guaranteed by Proposition~\ref{silvahaveDL}(b)
and we
might use Proposition~\ref{diff-silva}
instead of Proposition~\ref{union-cia}
to establish continuity of the algebra multiplication
and the inversion map.
In fact, Montel's Theorem shows that the unit ball $B$ in $\cA_n$
is relatively compact in the locally convex space $\cO(U_n)$
of holomorphic functions on~$U_n$, endowed with
the compact-open topology.\footnote{By hand: It is a uniformly
bounded set of functions
and equicontinuous due to the Cauchy estimates
for directional derivatives and Lemma~\ref{lipviaprime}(c).
Thus Ascoli's Theorem applies.}
The restriction map $\rho\colon \cO(U_n)\to \cA_{n+1}$
being continuous, we deduce that $\rho_{n+1,n}(B)=\rho(B)$
is relatively compact in $\cA_{n+1}$.
\end{rem}
\begin{defn}
Let $\cA$ be an associative algebra
over $\K\in\{\R,\C\}$
and $x\in \cA$. An element $y\in \cA$ is called
\emph{quasi-inverse} to~$x$ if
$xy=yx$ and $x+y-xy=0$.
If $x$ admits a quasi-inverse, then
$x$ is called \emph{quasi-invertible}.
We write $Q(\cA)$ for the set
of all quasi-invertible elements in~$\cA$.
\end{defn}
Quasi-inverses are unique if they exist,
by the following lemma.
We can therefore write $q(x)$ for the quasi-inverse
of $x\in Q(\cA)$.

In the following two lemmas,
we consider $\cA_+:=\cA\times \K$
as the direct sum $\cA\oplus \K \1 $, writing $\1 := (0,1)$
and identifying $\cA$ with $\cA\times \{0\}$.
Then the multiplication $(a,t)(b,s):=(tb+sa+ab,ts)$ makes
$\cA_+$ a unital associative $\K$-algebra
(see Exercise~\ref{exer:I.1}).
\begin{lem}\label{quasi-rel-inv}
Let $\cA$ be an associative $\K$-algebra.
For $x\in \cA$, we have:
\begin{description}[(D)]
\item[\rm(a)]
If $\cA$ is unital, then $x\in Q(\cA)$ if and only if $\1-x\in \cA^\times$.
In this case, $y\in \cA$ is quasi-inverse to $x$ if and only if
$\1-y=(\1-x)^{-1}$.
\item[\rm(b)]
In any case,
we have $x\in Q(\cA)$ if and only $\1 - x\in (\cA_+)^\times$.
Then ${y :=\1 - (\1 -x)^{-1}}$
is in $\cA$ and this is the unique quasi-inverse of~$x$ in~$\cA$.
\end{description}
\end{lem}
\begin{prf}
(a) For $x,y\in \cA$,
we have
\begin{eqnarray*}
(\1- x)(\1 -y) &= & \1 -x-y+xy\;\;\mbox{and}\\
(\1 -y)(\1 -x) &=& \1 -x-y+yx.\quad
\end{eqnarray*}
Thus $y$ is quasi-inverse to $x$ if and only
if both $(\1-x)(\1 -y)$ and $(\1-y)(\1-x)$ equal $\1$,
i.e., if and only if $\1-x$ is invertible with inverse $\1-y$. 

(b) We readily check that
$Q(\cA)\sub \cA\cap Q(\cA_+)$.
If $x\in  \cA\cap Q(\cA_+)$,
then $x$ has a quasi-inverse $t\1 +y$ in $\cA_+ $,
with $t\in \K$, $y\in \cA$. Then
\[
0=x + t\1  +y-x(t\1 +y)=t\1 +(x+y-tx-xy)
\]
implies that $t=0$ and $x+y-xy=0$.
Since $x$ and $t \1 +y=y$ commute in $\cA_+$,
they commute in $\cA$. Thus $y=q(x)$ in $\cA$.
The assertions follow.
\end{prf}
\begin{defn}\label{defn:cqia}
A locally convex associative topological algebra $\cA$
(which need not be unital) is called a \emph{continuous
quasi-inverse algebra} if $Q(\cA)$ is open in $\cA$
and the quasi-inversion map $q\colon Q(\cA)\to \cA$
is continuous.
\end{defn}
\begin{rem}\label{remciacqia}
Every cia $\cA$ also is a continuous quasi-inverse
algebra, as $Q(\cA)=\1-\cA^\times$ is open in $\cA$
and $q(x)=\1 - (\1 -x)^{-1}$ a continuous function of $x\in Q(\cA)$.
\end{rem}
\begin{prop}\label{cqia-necessities}
Let $\cA$ be an associative $\K$-algebra
and $q\colon Q(\cA)\to\cA$ be its quasi-inversion map.
Then $(\cA_+)^\times =\K^\times (\1 -Q(\cA))$
and
\begin{equation}\label{inv-via-quasi}
(t\1 +a)^{-1}= \frac{1}{t}\1 -\frac{1}{t}\, q\Big({-\frac{1}{t}a}\Big)
\end{equation}
for all $(a,t)\in \cA \times \K$ such that $t\1 +a\in (\cA_+)^\times$.
If $\cA$ as before is a locally convex
topological algebra, the the following holds:
\begin{description}[(D)]
\item[\rm(a)]
$Q(\cA)$ is open in $\cA$ if and only if $(\cA_+)^\times$
is open in $\cA_+$.
\item[\rm(b)]
$q$ is continuous if and only if the inversion map
$j\colon (\cA_+)^\times \to\cA_+$ is continuous.
\item[\rm(c)]
If $\cA$ is a continuous quasi-inverse algebra, then $q$
is $\K$-analytic.
\end{description}
\end{prop}
\begin{prf}
By Lemma~\ref{quasi-rel-inv}(b),
we have $\1-Q(\cA)\sub (\cA_+)^\times$.
As $\K^\times\1\sub (\cA_+)^\times$,
we deduce that $\K^\times (\1-Q(\cA))=(\K^\times \1)(\1-Q(\cA))\sub (\cA_+)^\times$.
If $t\1 +a\in (\cA_+)^\times$
with $t\in \K$ and $a\in \cA$,
then $t\not=0$
(if $t=0$, then $t\1 +a$ cannot be invertible in $\cA_+$,
as $\cA_+ \cA\sub \cA$).
Then $\1-(-(1/t)a)=(1/t)(t\1+a)\in (\cA_+)^\times$,
whence $-(1/t)a\in Q(\cA)$ by Lemma~\ref{quasi-rel-inv}(b)
and 
\[ (\1-(-(1/t)a))^{-1}=\1-q(-(1/t)a).\]
 Thus
$(t\1+a)^{-1}=(1/t)(\1-q(-(1/t)a))$ satisfies~(\ref{inv-via-quasi}).

(a) If $(\cA_+)^\times$ is open in $\cA_+$, then
$Q(\cA)=\{x\in \cA\colon \1-x\in (\cA_+)^\times\}$
is open in $\cA$. If, conversely, $Q(\cA)$ is open,
then it contains an open, balanced \break {$0$-neighborhood} $U\sub\cA$.
For each $t\in \K$ such that $|t|>1/2$,
we have $-t Q(\cA)\supseteq -t U\supseteq (1/2)U$,
whence
\[
(\cA_+)^\times=\K^\times (\1-Q(\cA))\supseteq ((1/2)U)\times \{t\in\K\colon |t|>1/2\}.
\]
Thus $(\cA_+)^\times$ is a neighborhood of $\1=(0,1)$ in $\cA_+=\cA\times \K$
and hence open, by Exercise~\ref{exer:I.2.}(2).

(b) If $q$ is continuous, then also $\iota$, by (\ref{inv-via-quasi}).
If $\iota$ is continuous, then $q(x)=\1-(\1-x)^{-1}$ is continuous in~$x$.

(c) If $\cA$ is a continuous quasi-inverse algebra,
then $\cA_+$ is a cia by (a) and~(b).
Hence $\iota$ is $\K$-analytic
(see Propositions~\ref{autocpx} and~\ref{autocxan}),
whence $q(x)=\1-(\1-x)^{-1}$ is $\K$-analytic in~$x$.
\end{prf}
For further discussions, see Exercise~\ref{cia-vs-quasi}.
\begin{prop}\mlabel{prop:tefu-cia}
Let
$M$ be $\sigma$-compact locally compact $C^r$-manifold over~$\R$
with $r\in \N_0\cup\{\infty\}$
and $\cA$ be a cia $($or a continuous quasi-inverse algebra$)$.
Then $C^r_c(M,\cA)$ is a continuous quasi-inverse algebra.
Notably, $C^r_c(M,\C)$ is a continuous quasi-inverse algebra.
\end{prop}
\begin{prf}
In view of Remark~\ref{remciacqia},
it suffices to consider a continuous quasi-inverse algebra~$\cA$.
Let $Q(\cA)\sub\cA$ be the open subset of quasi-invertible
elements and $q\colon Q(\cA)\to \cA$ be the quasi-inversion map,
which is smooth by Lemma~\ref{cqia-necessities}(c).
By Example~\ref{tefu-algs}(b),
$C^r_c(M,\cA)$ is a topological algebra.
Its subset $C^r_c(M, Q(\cA))$ is open by Corollary~\ref{open-in-tefu},
and we readily check that $Q(C^r_c(M,\cA))=C^r_c(M,Q(\cA))$
with quasi-inversion-map $C^r_c(M,q)$,
which is smooth (and hence continuous) by Proposition~\ref{variant-pushforwards}(a).
\end{prf}

\section{Matrix algebras over cias} 
\mlabel{sec:cia-matrix}

In this section we show that, for any (unital) cia $\cA$, the 
corresponding matrix algebra $M_n(\cA)$ also is a cia.

\begin{prop} \mlabel{prop:12.2.3} 
If $\cA$ is a  cia, then the same holds for the matrix algebras 
$M_n(\cA)$, $n \in \N$.  In particular, $\GL_n(\cA) = M_n(\cA)^\times$ is a Lie group. 
\end{prop}

\begin{prf} We endow $M_n(\cA)$ with the product topology from $\cA^{n^2}$. Then it is 
clear that $M_n(\cA)$ is a locally convex unital algebra. We prove the assertion 
by induction on $n \in \N$. For $n = 1$ there is nothing to show. 

First we show that $\GL_n(\cA)$ is open in $M_n(\cA)$. 
Let $x \in M_n(\cA)$ and write 
\[x = \begin{pmatrix} \alpha & \beta \\ \gamma & \delta \end{pmatrix} 
\ \ \\ \mbox{with}\ \ 
\alpha \in \cA, \beta \in M_{1,n-1}(\cA),
\gamma \in M_{n-1,1}(\cA), \delta \in M_{n-1}(\cA).\] 
Then $x$ is contained in the image of the smooth map 
\begin{eqnarray*}
\Phi &&\: \cA^\times \times M_{1,n-1}(\cA) \times M_{n-1,1}(\cA) \times \GL_{n-1}(\cA) 
\to M_n(\cA), \cr
(a,b,c,d) &&\mapsto 
 \begin{pmatrix}\1 & b \cr 0 & \1 \end{pmatrix} 
\begin{pmatrix}a & 0 \cr 0 & d \end{pmatrix} 
\begin{pmatrix}\1 & 0 \cr c & \1 \end{pmatrix} 
= \begin{pmatrix}a  + b d c & b d \cr d c & d \end{pmatrix}
\end{eqnarray*} 
if and only if $\delta \in \GL_{n-1}(\cA)$ and 
$\alpha - \beta \delta^{-1} \gamma \in \cA^\times$. 
Since $\cA$ and $M_{n-1}(\cA)$ are cias, the image of $\Phi$ is open, 
and from the explicit formula for its inverse 
$$ \Phi^{-1}\begin{pmatrix} \alpha & \beta \\ \gamma & \delta \end{pmatrix} 
= (\alpha - \beta \delta^{-1}\gamma, \beta \delta^{-1}, \delta^{-1}\gamma, \delta) $$ 
is follows that $\Phi$ is a diffeomorphism onto an open subset of 
$M_n(\cA)$ contained in $\GL_n(\cA)$. We conclude in particular that $\GL_n(\cA)$ is open.

From 
$$ \Phi(a,b,c,d)^{-1} 
= \begin{pmatrix}a^{-1} & a^{-1}b \cr -ca^{-1} & d^{-1} + c a^{-1}b \end{pmatrix}$$ 
and the continuity of the inversion in $\cA$ and 
$M_{n-1}(\cA)$, 
we further derive that the inversion in $\GL_n(\cA)$ is continuous in an identity 
neighborhood. In view of Lemma~\ref{lem:10.1.1}, 
this proves that $M_n(\cA)$ is a cia and hence that $\GL_n(\cA)$ is a Lie group.
\end{prf}

\begin{rem} \mlabel{rem-cia} 
For $n \in \N$ we obtain in particular 
that 
\begin{eqnarray*}
C^\infty(M,\GL_n(\cA)) 
&&= C^\infty(M,M_n(\cA)^\times) 
= C^\infty(M,M_n(\cA))^\times  \\
&&= M_n(C^\infty(M,\cA))^\times  = \GL_n(C^\infty(M,\cA))
\end{eqnarray*}
is a Lie group (Lemma~\ref{lem:A.5}(3)). 
\end{rem}

\begin{prop} \mlabel{prop:12.1.1} 
If $\cA$ is a  cia and $M \subeq \cA$ is a subset, then 
its commutant 
\[ M' := \{ a \in \cA \:  (\forall m \in M)\, am=ma \} \] 
of $M$ in $\cA$ is also a cia with 
$(M')^\times = \cA^\times \cap M'.$ 
\end{prop} 

\begin{prf} Since $a \in \cA^\times$ commutes with 
$M$ if and only if $a^{-1}$ does, we have 
$(M')^\times = \cA^\times \cap M'$. 
This implies that the unit group $(M')^\times$ is open, and 
the continuity of the inversion follows by restriction. 
\end{prf}

\begin{lem} \mlabel{lem:12.2.7} If $\cA$ is a finite-dimensional unital $\K$-algebra and 
$\rho_b(a):= ab$ its right multiplications, 
then 
\[ \cA \cong \rho_\cA' := \{ A \in \End_\K(\cA) \: 
(\forall b \in \cA) A \rho_b = \rho_b A\}. \]
\end{lem} 

\begin{lem} \mlabel{lem:12.2.8} If $\cA$ and $\cB$ are unital $\K$-algebras and 
$M \subeq \cB$, then 
$$ (\1 \otimes M)' = \cA \otimes M' $$
holds for the commutants of $M$ in $\cA \otimes \cB$, resp., $\cB$. 
\end{lem}

\begin{prop} \mlabel{prop:8.1.4} If $\cA$ is a unital 
cia and $\cB$ is a finite-dimensional unital algebra, 
then $\cA \otimes \cB$ is a cia. 
\end{prop} 

\begin{prf} For $n := \dim \cB$, the regular representation 
$\lambda_b(x) := bx$ of $\cB$ 
yields an embedding $\lambda \: \cB \into \End_\K(\cB) \cong M_n(\K)$. 
This leads to an embedding 
\[ \cA \otimes \cB \to \cA \otimes \End_\K(\cB)
\cong \cA \otimes M_n(\K) \cong M_n(\cA), \] 
and since $\cB$ is a commutant in 
$\End_\K(\cB)$ (Lemma~\ref{lem:12.2.7}), the same holds for 
$\cA \otimes \cB$ in $M_n(\cA)$ (Lemma~\ref{lem:12.2.8}). 
As $M_n(\cA)$ is a cia by Proposition~\ref{prop:12.2.3}, 
the assertion follows from Proposition~\ref{prop:12.1.1}. 
\end{prf}

The following corollary is particularly useful for spectral theoretic 
arguments in real cias. In particular, it provides a second argument 
for the stability of the cia property under complexification 
(Proposition~\ref{autocxan}). 

\begin{cor} \mlabel{cor:8.1.5} 
For any real cia $\cA$, the complexification $\cA_\C\cong \cA \otimes_\R \C$ is a 
complex cia. 
\end{cor}  

\begin{defn} Let $\cA$ be a real unital associative algebra and $\cA_\C$ be 
its complexification. 
Complex conjugation provides an involutive automorphism $\sigma(a + ib) = a - i b$ 
of $\cA_\C$ with 
$(\cA_\C)^\sigma = \cA$. If $a \in \cA$ is invertible in $\cA_\C$, then $a^{-1} \in (\cA_\C)^\sigma = \cA$
implies that 
\begin{equation}
  \label{eq:complexif-equispec}
\cA \cap \cA_\C^\times = \cA^\times. 
\end{equation}
For $a \in \cA$ we define 
$$ \Spec(a) := \Spec_\cA(a) := 
\{ \lambda \in \C \: a - \lambda \1 \not\in \cA_\C^\times\}. $$

If $\lambda = \alpha + i \beta \in \Spec(a)$, then 
\begin{eqnarray}
  \label{eq:comp-spec}
(a - \lambda \1)(a - \oline\lambda\1) = a^2 - 2\alpha a + |\lambda|^2 \1 
\end{eqnarray}
is an element of $\cA$, not invertible in $\cA_\C$ (Exercise~\ref{ex:inv-prod}), 
hence not invertible in $\cA$ by \eqref{eq:complexif-equispec}. 
If, conversely, (\ref{eq:comp-spec}) is not invertible 
in $\cA$, it is not invertible in $\cA_\C$, which implies that 
$a - \lambda \1$ or $a - \oline \lambda \1$ is not invertible. 
In view of 
\[ \sigma(a - \lambda \1) = a - \oline\lambda \1,\]
 it follows 
that $\lambda \in \Spec(a)$. We therefore have the following 
description of the spectrum in terms of the real algebra $\cA$: 
$$ \Spec_\cA(a) = \{ \lambda \in \C \: (a - \lambda \1)(a -\oline \lambda \1) \not\in \cA^\times\}. $$
\end{defn}

\section{Algebraic subgroups} 
\mlabel{subsec:alg-grp}

We will now discuss a very convenient criterion for a subgroup of a 
Banach--Lie algebra to be a Lie group which applies easily 
in many concrete situations. 
First, we need the concept of a polynomial function and of an algebraic subgroup. 

\begin{defn} \mlabel{defn:5.5.13} Let $\cA$ be a Banach algebra. 
A subgroup $G \subeq \cA^\times$ is called 
\index{algebraic subgroup, of Banach algebra}
{\it algebraic} if there exists a
$d \in \N_0$ and a set 
${\cal F}$ of Banach space valued polynomial functions on $\cA\times \cA$ 
of degree $\leq d$ such that 
$$ G = \{ g \in \cA^\times \: (\forall f \in {\cal F})\ f(g,g^{-1}) = 0\}. $$
\end{defn}

To prove the Algebraic Subgroup Theorem below, we first have
to verify the following special case. 

\begin{prop} \mlabel{prop:5.5.15} 
Let $E$ be a Banach space and $F \subeq E$ be a closed
subspace. Then 
$$ \GL(E)_F := \{ g \in \GL(E) \: gF = F\} $$
is  a Banach--Lie subgroup of $\GL(E)$ with 
$$ \L^e(\GL(E)_F) = \{ X \in \cL(E) \: XF \subeq F\}. $$
\end{prop}

\begin{prf}
Let $V \subeq
\gl(E) = {\cal L}(E)$ be an open $0$-neighborhood such that $\exp\res_V \: V \to
\exp V$ is a diffeomorphism and $\|\exp x - \1\| < 1$ for all $x \in V$. 
Then the inverse function 
$$ \log \:= (\exp\res_V)^{-1} \: \exp V \to \gl(E) $$
is given by the convergent power series 
$$ \log(g)  = \sum_{n=1}^\infty \frac{(-1)^{n+1}}{n} (g-\1)^n $$
(cf.\ Example~\ref{ex:abgrp}(b)). 
For $g = \exp x\in (\exp V) \cap H$, we then
obtain $x(F) \subeq F$ directly from the power series of $\log(g)$. 

Conversely, if $XF \subeq F$ and $t \in \R$, then $\exp(tX)F \subeq
F$ follows directly from the power series expansion of $\exp$.
In particular, we have ${\exp(-tX)F \subeq F}$, too, so that 
$\exp(tX)F = F$. 
\end{prf}

\index{Algebraic Subgroup Theorem} 
\begin{thm}[Algebraic Subgroup Theorem] \mlabel{thm:5.5.14} 
Every algebraic subgroup $G$ of the unit group $\cA^\times$ 
of a unital Banach algebra $\cA$ is a Banach--Lie subgroup. 
\end{thm}

\begin{prf} In view of the Hahn--Banach Theorem, we may assume that 
$${\cal F} \subeq P := \bigoplus_{k=0}^d P_k(\cA \times \cA, \K), $$
the Banach space of scalar-valued continuous polynomials on $\cA \times
\cA$ of degree $\leq d$. 
The space $P$ carries a natural Banach space structure such that the
action of $\cA^\times$ on $P$ given by 
$$ \big(\pi(g)f\big)(x,y) := f(xg, g^{-1}y) $$
yields a continuous homomorphism 
$\pi \: \cA^\times \to \GL(P)$ (Exercise~\ref{exer:oden-iv.6}(d)). 
Replacing ${\cal F}$ by 
$$ F := \{ f \in P \:(\forall g \in G)\, f(g,g^{-1}) = 0\}, $$
we may assume that ${\cal F} = F$. The space $F$ is a closed
subspace of $P$. We claim that 
\begin{equation}
  \label{eq:oden-4.1}
G = \{ g \in \cA^\times \: \pi(g)F \subeq F \}.
\end{equation}
In fact, if $g,x \in G$ and $f \in F$, then 
$$ \big(\pi(g)f\big)(x,x^{-1}) = f(xg, g^{-1}x^{-1}) = 0, $$
showing that $\pi(g)f \in F$. If, conversely, $g \not\in G$, then
there exists an $f \in F$ with 
$$ 0 \not= f(g,g^{-1}) = \big(\pi(g)f\big)(\1,\1). $$
It follows in particular that $\pi(g)f \not\in F$. Next, since 
$G$ is a group, \eqref{eq:oden-4.1} implies 
$$ G = \{ g \in \cA^\times \: \pi(g)F = F \}. $$
We conclude that 
$$ G = \pi^{-1}(\{ g \in \GL(P) \: \pi(g)F = F \}), $$
so that the assertion follows from Proposition~\ref{prop:5.5.9} 
on inverse images of Lie subgroups, 
combined with Proposition~\ref{prop:5.5.15}.
\end{prf}

\begin{exs} \mlabel{ex:IV.15-ne04} 
(a) If $(\cA,*)$ is a unital Banach-* algebra, then its 
unitary group 
\[  \U(\cA) := \{ a \in \cA \: aa^* = a^* a = \1 \} \] 
is a real Banach--Lie subgroup. This follows from 
Theorem~\ref{thm:5.5.14} if we consider $\cA$ as a real Banach algebra 
and $\U(\cA)$ as an algebraic subgroup specified by 
$f(g,g^{-1}) = 0$, where $f(g,h) = g^* - h$. 

(b) If $\cA$ is a Banach space and $m \: \cA \times \cA \to \cA$ 
a continuous bilinear map, then the corresponding automorphism group 
\[ \Aut(\cA,m) := \{ g \in \GL(\cA) \: (\forall a,b \in \cA)\ m(ga,gb) =
gm(a,b)\} \]
is a Banach--Lie subgroup whose Lie
algebra is the space 
$$ \der(\cA,m) := \{ X \in {\cal L}(\cA) \: (\forall a,b \in \cA)\, Xm(a,b) =
m(Xa,b) + m(a,Xb)\} $$
of derivations of $(\cA,m)$ (Exercise~\ref{exer:5.6.6}). 

(c) If $V$ and $W$ are Banach spaces and $\beta \: V \times  V\to W$ is a 
continuous bilinear map, then 
\[ \OO(V,\beta) := \{ g\in \GL(V) \: (\forall v,v' \in V)\, 
\beta(gv,gv') = \beta(v,v') \} \] 
is a Banach--Lie subgroup of $\GL(V)$ because it is algebraic in the 
unital Banach algebra $\cL(V)$. Each pair $v,v' \in V$ specifies 
a polynomial function $f_{v,v'}(g) := \beta(gv,gv') - \beta(v,v')$, and 
$\OO(V,\beta)$ is the simultaneous zero sets of these polynomials in $\GL(V)$. 
It is easy to see that the Lie algebra of this groups is 
\[ \fo(V,\beta) = \{ X \in \cL(V) \: (\forall v,v' \in V)\, 
\beta(Xv,v') + \beta(v,Xv') = 0\}.\] 
\end{exs}

\begin{exs}
  \mlabel{ex:IV.15-ne04b} (a) If $\cA$ is a unital Banach algebra and 
$M_n(\cA)$ is the algebra of $(n \times n)$-matrices with entries in
$A$, then $M_n(\cA)$ also is a Banach algebra. In fact, on the space 
$\cA^n$ we consider the norm given by 
\[  \|x\| := \max\{\|x_1\|, \ldots, \|x_n\|\}. \] 
Then $\cA^n$ is a Banach space and we have a natural embedding 
$M_n(\cA)  \into \cL(\cA^n)$ 
which we use to define a norm on $M_n(\cA)$. It is not hard to verify 
that $M_n(\cA)$ is closed in $\cL(\cA^n)$, hence a Banach algebra. 
We write $\GL_n(\cA) := M_n(\cA)^\times$ for the unit group of this Banach 
algebra (see also Proposition~\ref{prop:12.2.3}). 

(b) As we will see below, it sometimes is convenient to
refine the preceding construction as follows. Let $\cJ \trile \cA$ be an
ideal which is a Banach algebra in its own right such that the
multiplication map 
\[  \cA \times \cJ \to \cJ, \quad (a,b) \mapsto ab \] 
is continuous, i.e., there exists a $C > 0$ with 
\[ \|ab\|_J \leq C \|a\|_\cA \|b\|_\cJ \quad \mbox{  for } \quad a \in \cA, b \in \cJ.\] 
After replacing the norm on $A$ by the equivalent norm 
\[  \|a\|' := \max(\|a\|_\cA,\sup\{ \|ab\|_\cJ\: \|b\|_\cJ \leq 1\}) 
\leq \max(1,C)\|a\|_\cA, \]
we may assume that $\|ab\|_J \leq \|a\|_A\|b\|_J$ holds for 
$a \in \cA$ and $b \in \cJ$. 

We consider the algebra 
\[  M_2(\cA,\cJ) := \Big\{\pmat{ a & b \\  c & d} \: 
a,d \in \cA, b,c \in \cJ \Big\} \] 
endowed with the norm 
\[  \Big\|\pmat{ a & b \\ c & d}\Big\| :=
2 \max\{\|a\|_\cA,\|b\|_\cJ,\|c\|_\cJ, \|d\|_\cA\}. \] 
Then $\|xy\| \leq \|x\|\|y\|$ holds for $x,y \in M_2(\cA,\cJ)$ 
(Exercise!), 
so that $M_2(\cA,\cJ)$ is a Banach algebra. 

A similar construction works for $(n \times n)$-matrices, where one
defines 
\[  M_n(\cA,\cJ) := \{ (x_{ij})_{i,j=1,\ldots,n} \in M_n(\cA) \: 
i \not= j \Rarrow x_{ij} \in J \} \] 
and 
\[  \|(x_{ij})\| := n \max \{ \|x_{ij}\|_\cJ, i\not= j; \|x_{ii}\|_\cA,
i=1,\ldots,n\}. \] 
We write $\GL_n(\cA,\cJ)$ for the unit group of this Banach algebra. 

(c) If $\cA$ is a Banach algebra over $\K \in \{\R,\C\}$ without a unit element, then
we endow the space $\cA_+ := \cA \oplus \K$ 
with the Banach algebra structure given by 
\[  \|(a,z)\| := \|a\| + |z| \quad \hbox{ and } \quad 
(a,z)(a',z') := (aa' + z a' + z' a, zz'). \] 
Then $\cA \cong \cA \times \{0\}$ is a closed ideal in $\cA_+$, and we 
have an algebra homomorphism $\eps \: \cA_+ \to \K$ given by  
$\eps(a,z) = z$. We define 
\[ \cA^\times := \eps^{-1}(1) \cap \cA_+^\times.\] 
This is a closed subgroup of $\cA_+^\times$ and 
$\{ (a,1) \: \|a\| < 1 \}$
is an open $\1$-neigh\-bor\-hood in $\cA^\times$. Therefore 
$\cA^\times$ is a Lie subgroup of $\cA_+^\times$ with the Lie algebra $(\cA, [\cdot,\cdot])$ and the exponential function 
\[ \exp \: \cA \to \cA^\times, \quad x \mapsto e^x = (e^x - \1,\1). \] 

(d) If $\cJ \trile \cA$ is an ideal of the unital Banach algebra $\cA$, then 
\[  J^\times := \cA^\times \cap (\1 + J) \] 
is the kernel of the homomorphism 
\[ \cA^\times \to (\cA/\cJ)^\times, \quad g \mapsto g + \cJ, \] 
where $(\cA/\cJ)^\times$ denotes the unit group of the unital algebra $\cA/\cJ$
which is not required to carry a natural Banach space structure 
because we did not assume that $\cJ$ is closed. 
We have for $\cJ_+ := \K \1 + \cJ \subeq \cA$ the relation 
\[ \cJ_+^\times = \cA^\times \cap \cJ_+ = \K^\times \cdot \cJ^\times. \] 

If, in addition, $\cJ$ is closed, then $J^\times$ is a closed subgroup of 
$\cA^\times$ and a Lie group with Lie algebra $\cJ$. Further, 
$\cA/\cJ$ carries a natural Banach algebra structure with the norm 
\[  \|a + \cJ\| := \inf \{ \|a + b\| \: b \in \cJ \}, \] 
and the quotient map $\cA \to \cA/\cJ$ is continuous, so that we have an 
exact sequence of Lie groups 
\[  \{\1\} \to \cJ^\times \into \cA^\times \to (\cA/\cJ)^\times. \] 
Here the map on the right hand side need not be surjective. 
A typical example is $\cA = \cL(\cH)$ for an infinite-dimensional Hilbert space $\cH$ 
and $\cJ = K(\cH)$ (the ideal of compact operators). In fact, there 
are Fredholm operators with non-vanishing index (cf.\ \cite{Ru91}). 
\end{exs}


\section{Cias of Smooth Vectors} 
\mlabel{sec:cia-smoothvec}

In this section we show that for a Lie group action 
on a cia~$\cA$, the subalgebra $\cA^\infty$ also is a cia. 

\begin{prop}  \mlabel{prop:gramsch} 
Let $\cA$ is a  cia and $M$ be a smooth manifold. 
We consider a map $\alpha \: M \to \Aut(\cA), m \mapsto \alpha_m$ and 
put $\alpha^a(m) := \alpha_m(a)$. Then the subspace 
\[ \cA^\infty := \{ a \in \cA \: \alpha^a \in C^\infty(M,\cA)\} \] 
is a cia with respect to the subspace topology 
coming from the smooth compact open topology on 
$C^\infty(M,\cA)$.
\end{prop}

\begin{prf} We have already seen in Lemma~\ref{lem:A.5} that 
$C^\infty(M,\cA)$ is a locally convex algebra with continuous inversion and 
unit group $C^\infty(M,\cA^\times)$. We consider the inclusion 
\[ \eta \:  \cA^\infty \into C^\infty(M,\cA), \quad \eta(a) := \alpha^a.\]
Since every $\alpha_m$ is an automorphism of $\cA$, $\eta$ is an algebra 
homomorphism. By definition, $\eta$ is a topological embedding. 
Therefore $\cA^\infty$ is a locally convex algebra with continuous inversion. 
It remains to see that its unit group is open. 

We may w.l.o.g.\ assume that $M \not=\eset$. Pick $m \in M$. 
If $a \in \cA^\infty$ is invertible, then 
$\alpha_m(a) = \alpha^a(m)$ is invertible, and therefore 
$a \in \cA^\times$ because $\alpha_m$ is an automorphism. 
If, conversely, $a \in \cA^\times$, 
then $\alpha^a(M) \subeq \cA^\times$ implies that $\alpha^a$ is invertible 
in $C^\infty(M,\cA)$ with inverse $(\alpha^a)^{-1} = \alpha^{a^{-1}}$, 
and therfore $a \in (\cA^\infty)^\times$. 
This shows that 
\[ (\cA^\infty)^\times = \cA^\infty \cap \cA^\times = 
(\ev_{m}\circ \eta)^{-1}(\cA^\times).\] 
As $\ev_{m} \circ \eta$ is continuous, $(\cA^\infty)^\times$ is open and 
therefore $\cA^\infty$ is a cia. 
\end{prf}

In non-commutative differential geometry, one needs derivations of an algebra, 
because they are expected to describe the infinitesimal structure in a similar 
fashion as vector fields describe the infinitesimal structure of a manifold. 
Therefore one is forced to consider actions of a Lie group 
which lead to an action of its Lie algebra by derivations. 
In this context, the following proposition  is very useful tool.

\begin{prop}
  \mlabel{prop:VI.7} {\rm($\cA^\infty$ as a cia)} Let $\cA$ be a cia,  
$G$ be a Lie group  and 
$\alpha \: G \to \Aut(\cA), g \mapsto \alpha_g$ define 
an action of $G$ on $\cA$. 
Then the subspace $\cA^\infty$ of smooth vectors is a 
cia with respect to the  topology inherited from the
embedding $\cA^\infty \into C^\infty(G,\cA), a \mapsto \alpha^a$, 
$\alpha^a(g) := \alpha_g(a)$. 

The action of $G$ on $\cA$ leaves $\cA^\infty$ invariant and 
induces on $\cA^\infty$ an action by algebra automorphisms, with 
smooth orbit maps, and 
for which the derived action $\g \times \cA^\infty \to \cA^\infty$ 
is separately continuous. If the action of $G$ 
and the derived action of $\g$ on $\cA^\infty$ are continuous, then the $G$-action 
on $\cA^\infty$ is smooth. 
\end{prop} 

\begin{prf} Proposition~\ref{prop:gramsch} implies that 
$\cA^\infty$ is a cia. For $\alpha^a(g) =: g.a$ and $h \in G$ we have 
\begin{equation}
  \label{eq:6.1x}
\alpha^{h.a}(g) = (gh).a = \alpha^a(gh), 
\end{equation}
so that, for each $a \in \cA^\infty$, the orbit map 
$G \to \cA^\infty, h \mapsto h.a$ corresponds to the map 
$$ G \to C^\infty(G,\cA), \quad h \mapsto \alpha^a \circ \rho_h, $$
and the smoothness of this map follows from the smoothness of the map 
$$ G \times G \to \cA, \quad (g,h) \mapsto \alpha^a(\rho_h(g)) = \alpha^a(gh) $$
and Proposition~\ref{prop:cartes-closed}. Therefore the action of $G$ on $\cA^\infty$ has smooth 
orbit maps. 

From \eqref{eq:6.1x} we derive that the action of $G$ 
on $\cA^\infty$ corresponds to the right translation action of $G$ on $C^\infty(G,\cA)$. 
Now 
\[  x.a := T_\be(\alpha^a)x \quad \mbox{ resp. } \quad 
\alpha^{x.a} = x^l \alpha^a  \] 
(cf.\ Proposition~\ref{prop:3e.1.5}) 
defines a representation of the Lie algebra $\g$ of $G$, which is a 
subrepresentation of the representation of $\g$ on $C^\infty(G,\cA)$, 
defined by left invariant vector fields. 
Clearly, $\g$ acts by continuous operators on $\cA^\infty$ and, 
for each $a \in \cA$, the map $x \mapsto x.a$ is continuous linear. 
This means that the action map $\g \times \cA \to \cA$ is separately continuous. 

Now we assume that the $G$-action 
$\alpha_g(a) =:g.a$ on $\cB := \cA^\infty$ and its derived 
action $\alpha^\g \:\g \times \cB \to \cB, (x,a) \mapsto x.a$ are continuous. 
First we note that we have a continuous tangent action 
\[  T\alpha \: TG \times T\cB \to T\cB, \quad 
(g.x).(a,w) := (g.a, g.w + g.(x.a)), \] 
and this implies that $\alpha$ is a $C^1$-action. The formula for $T\alpha$ now 
implies that $T\alpha$ is also $C^1$. By iteration we conclude that 
the action of $G$ on $\cB$ is smooth.  
\end{prf}

The preceding proposition applies in particular to the important class of actions of 
finite-dimensional Lie groups on Banach algebras. 

\begin{ex}
  \mlabel{ex:VI.8} Smooth quantum tori form an important interesting class of cias 
(\cite{BEG89}), they are defined as follows.  We consider a matrix 
$\bq := (q_{jk}) \in M_d(\C)$ satisfying 
\begin{description}
\item[\rm(Q1)] $q_{jj} = 1$ for each $j$. 
\item[\rm(Q2)] $q_{jk} = q_{kj}^{-1}$ for all $j,k$. 
\item[\rm(Q3)] $|q_{jk}| = 1$ for all $j,k$. 
\end{description}
Then there is a universal $C^*$-algebra $\cC$ generated by 
$d$ unitary elements $u_1,\ldots, u_d$ satisfying the commutation relations 
\[ u_j u_k = q_{jk} u_k u_j \quad \hbox{for} \quad 1 \leq j,k\leq d. \] 
Then 
\[ \alpha \: \T^d \to \Aut(\cC), \quad 
\alpha(\bt) u_j = t_j u_j, \qquad \bt = (t_1,\ldots, t_d) \] 
defines a continuous action of the torus group $\T^d$ by automorphisms on 
$\cC$. The corresponding smooth quantum torus is the Fr\'echet cia 
$\cA_\bq := \cC^\infty$ 
of smooth vectors for this action (Proposition~\ref{prop:VI.7}). 
This implies that each element of $\cA_\bq$ can be represented by a convergent series 
\[  a= \sum_{\alpha \in \Z^n} a_\alpha u^\alpha\quad \mbox{ with } \quad 
a_\alpha \in \C, u^\alpha := u_1^{\alpha_1} \cdots u_d^{\alpha_d}, \alpha \in \Z^d \]
(cf.\ Theorem~\ref{thm:peter-weyl}). More precisely, we have 
\begin{align*}
\cA_\bq 
&= \Bigg\{ \sum_{\alpha \in \Z^d} a_\alpha u^\alpha \: 
(\forall k \in \N) \sum_\alpha  |\alpha|^k |a_\alpha| < \infty\Bigg\}\\
&= \Bigg\{ \sum_{\alpha \in \Z^d} a_\alpha u^\alpha \: 
(\forall k \in \N) \sup |\alpha|^k |a_\alpha| < \infty\Bigg\}. 
\end{align*}
If $q_{jk} = 1$ for all $j,k$, then $\cC \cong C(\T^d)$ and 
$\cA_\bq \cong C^\infty(\T^d)$ is the algebra of smooth functions on~$\T^d$. 
\end{ex}

\begin{rem} \mlabel{rem:VI.9} 
If $G$ is a finite-dimensional Lie group and 
$\cA$ is Mackey complete, then, for each 
$f \in C^\infty_c(G,\R)$, we obtain a continuous linear operator on 
$\cA$ by 
\[  \alpha_f(a) := \int_G\ f(g) \alpha_g(a)\, d\mu_G(g), \] 
where $\mu_G$ is a left Haar measure on $G$. For any sequence of functions
$(f_n)$ 
with $\int_G f_n\ d\mu_G = 1$ and support converging to $\be$, we then have 
$\alpha(f_n)a \to a$. Since $\im(\alpha_f)$ consists of smooth vectors, 
$\cA^\infty$ is dense in  $\cA$ (G\aa{}rding's Theorem \cite{Ga47}). 
\end{rem}

\begin{rem}
  \mlabel{rem:VI.10}Above we have endowed the algebra $\cA^\infty$ with the topology 
inherited from the space $C^\infty(G,\cA)$. There is an alternative 
way to obtain refined topologies via derivations. 

So let $\cB$ be a subalgebra of the locally convex algebra $\cA$ 
and consider a derivation $D \: \cB \to \cA$. Consider the 
tangent algebra 
$T\cA \cong \cA \times \cA$ 
with the multiplication 
\[ ((a,v),(a',v')) \mapsto (aa', av' + va').\] 
Then 
\[ \eta_D \: \cB \to T\cA, \quad a \mapsto (a,Da) \] 
is an injective algebra homomorphism (Exercise~\ref{ex:9.1.x}). 
Therefore the topology on $\cB$ for which 
$\eta_D$ is an embedding defines on $\cB$ the structure of a locally 
convex algebra. If $\cA$ is a cia, then so is $T\cA$, with unit 
group $(T\cA)^\times = \cA^\times \times \cA$. In particular, 
$\eta_D^{-1}((T\cA)^\times) = \cB \cap \cA^\times$, is open and therefore 
$\cB^\times$ is open. 

A typical example is $\cA = C(\T,\R)$ and $\cB := C^1(\T,\R)$, topologized 
as a subalgebra of the Banach algebra $\cA$. 
We then have the derivation 
\[ D \: \cB \to \cA, \quad D f = f' \] 
and the above process leads to the norm $\|f\|_1 := \|f\| + \|f'\|$ on $\cB$, 
turning $\cB$ into a Banach algebra. 
\end{rem}

\section{Regularity of unit groups} 
\mlabel{sec:cia-regular}

In this section we describe a  sufficient criterion for the 
regularity of the unit group $\cA^\times$ of a cia. 
In view of Proposition~\ref{prop:reg-Mack}, the Mackey 
completeness of $\cA$ is a necessary requirement. 
Since Banach--Lie groups are regular, it is also clear that 
$\cA^\times$ is regular if $\cA$ is a Banach algebra, but beyond 
that it is not so easy to find sufficient conditions. The conditions 
we described have different flavors: 
\begin{itemize}
\item[$\bullet$] A condition in terms of series involving constants 
$c_n(p,q)$, where $p$ and $q$ are continuous seminorms on $\cA$. 
This encodes some quantitative information on the continuity 
of iterated multiplication maps. 
\item[$\bullet$] We prove Turpin's Theorem asserting that commutative 
cias are locally $m$-convex, so that we have seminorms satisfying 
$c_n(p,p) \leq 1$, but in the commutative case, regularity of $\cA^\times$ 
follows easily from the Mackey completeness of~$\cA$. 
Since we also show that the completion of any commutative 
cia $\cA$ is a cia, this completely settles the commutative case. 
\item[$\bullet$] We also  show that regularity is inherited 
from $\cA$ to algebras $\cA^\infty$ of smooth vectors 
for automorphic actions of Banach--Lie groups on $\cA$. 
\end{itemize}

\subsection{Sufficient conditions}

Given continuous seminorms $p,q$ on $\cA$, we put 
\[ c_n(p,q) := \sup \{ p(a_1 \cdots a_n) \: 
q(a_j) \leq 1, j =1,\ldots, n\} \in [0,\infty].\] 
In \cite{GN12} the following criterion for the regularity of $\cA^\times$ is obtained: 

\begin{thm} Let $\cA$ be a Mackey complete cia such that, for every 
continuous seminorm $p$ on $\cA$, there exists a continuous 
seminorm $q$ on $\cA$ and an $r > 0$ such that 
\begin{equation}
  \label{eq:suff-cond-reg}
\sum_{n = 0}^\infty r^n c_n(p,q) < \infty.
\end{equation}
Then $\cA^\times$ is $C^1$-regular and if $\cA$ is sequentially complete, 
then $\cA^\times$ even is $C^0$-regular. 
\end{thm}

The main idea for the proof is to show convergence of the 
Picard iteration. Suppose that $\cA$ is a sequentially complete cia. 
For $u \in C([0,1],\cA)$ we want to solve the linear initial value problem 
\begin{equation}
  \label{eq:4.1.2b}
\gamma(0) = \1,\quad \gamma'(t) = \gamma(t)u(t). 
\end{equation}
To this end, we consider the Picard iteration: 
\[  \gamma_0(t) := \1, \quad 
\gamma_{n+1}(t) := \1 + \int_0^t \gamma_n(\tau)u(\tau)\, d\tau, \] 
which leads to 
\[  \gamma_n(t) =  \1 + \sum_{k = 1}^n \int_0^t \int_0^{\tau_n} \cdots \int_0^{\tau_2} 
u(\tau_1)u(\tau_2)\cdots u(\tau_n)\ d\tau_1\, d\tau_2 \cdots d\tau_n. \] 
If $q$ is a continuous seminorm with $q(u(t)) \leq r$ for $0 \leq t \leq 1$, then 
\begin{align*}
& p\Big(\int_0^t \int_0^{\tau_n} \cdots \int_0^{\tau_2} 
u(\tau_1)u(\tau_2)\cdots u(\tau_n)\ d\tau_1\, d\tau_2 \cdots d\tau_n\Big) \\
&\leq c_n(p,q) r^n \int_0^t \int_0^{\tau_n} \cdots \int_0^{\tau_2} 
\ d\tau_1\, d\tau_2 \cdots d\tau_n 
= c_n(p,q)r^n \frac{t^n}{n!}.
\end{align*}
Therefore Picard iteration leads to a uniformly convergent series 
\[ \Evol(u) := \lim_{n \to \infty} \gamma_n \] 
and it remains to show that this limit defines a smooth evolution map 
\[ \evol \: C([0,1],\cA)\to \cA^\times, \quad u \mapsto \Evol(u)(1).\] 
If $\cA$ is only Mackey complete, then the ingrals are only defined if 
$u$ is a $C^1$-curve.

\begin{rem} (a) Presently we do not know whether \eqref{eq:suff-cond-reg} 
is really needed or if all sequentially complete cias have $C^0$-regular 
unit groups.

If $\cA^\times$ is $C^0$-regular,
then $\cA^\times$ is locally $\mu$-convex in the sense of \cite{Gl12b},
as shown in~\cite{Hn22}.
This implies that, for each continuous seminorm~$p$ on~$\cA$,
there exists a continuous seminorm~$q$ on~$\cA$ and $a>0$ such that
\begin{equation}\label{weaker}
(\forall n\in\N) \;\, c_n(p,q)\leq a^nn^n,
\end{equation}
as we shall presently verify.
The estimate~(\ref{weaker}) is not strong enough to get~\eqref{eq:suff-cond-reg}.
In fact, the condition in~\eqref{eq:suff-cond-reg}
is equivalent to the existence of some $b>0$ such that
\[
(\forall n\in \N)\;\, c_n(p,q)\leq b^n,
\]
which is not a consequence of~(\ref{weaker}).\\[2.3mm]
[To prove (\ref{weaker}), let~$p$ be a continuous seminorm on~$\cA$.
By local $\mu$-convexity, there exists a continuous seminorm~$q$ on~$\cA$
and $\ve>0$ such that
\begin{eqnarray}
\lefteqn{(\forall n\in\N)\,(\forall x_1,\ldots, x_n\in\cA)}\quad\qquad \\ \notag
& & \hspace*{-8mm}\sum_{j=1}^n q(x_j)\leq\ve\;\Rightarrow \;p((\one+x_1)(\one+x_2)\cdots (\one+x_n)-\one)
\leq\sum_{j=1}^nq(x_j).\label{alg-muconv}
\end{eqnarray}
We may assume that $\ve\leq 1$.
Given $n\in\N$, $j\in\{1,\ldots,n\}$ and $y\in\cA$, consider the difference operator~$\Delta_{j,y}$
acting on functions $g\colon \cA^n\to\cA$
via
\[
\Delta_{j,y}(g)(a_1,\ldots,a_n)=g(a_1,\ldots, a_{j-1},a_j+y,a_{j+1},\ldots, a_n)
-g(a_1,\ldots, a_n)
\]
for $a_1,\ldots,a_n\in\cA$.
Consider the function 
\[ f_n\colon \cA^n\to \cA, \quad f_n(x_1,\ldots,x_n):=
(\one+x_1)\cdots(\one+x_n)-\one.\]
Given $y_1,\ldots, y_n\in\cA$ with $\sum_{j=1}^n q(y_j)\leq \ve$, we then have
\begin{eqnarray*}
y_1y_2\cdots y_n&=&(\Delta_{n,y_n}\cdots\Delta_{1,y_1}f)(0,\ldots,0)\\
&=&\sum_{\alpha_1,\ldots,\alpha_n=0}^1(-1)^{n-\alpha_1-\cdots-\alpha_n}
f_n(\alpha_1y_1,\ldots,\alpha_ny_n)
\end{eqnarray*}
and hence $p(y_1,\ldots,y_n)\leq 2^n\ve$ (using (\ref{weaker})).
Taking any $y_1,\ldots, y_n\in\cA$ such that $q(y_j)\leq \ve/n$ for all $j\in\{1,\ldots,n\}$,
we deduce that
\[
c_n(p,q)\leq \frac{2^n\ve n^n}{\ve^n}\leq a^nn^n
\]
with $a:=2/\ve$.]

(b) If the topology on $\cA$ can be defined by a family of submultiplicative 
seminorms $p$, i.e., $\cA$  is 
\index{cia!locally $m$-convex} 
{\it locally $m$-convex}, 
then $c_n(p,p) \leq 1$ and \eqref{eq:suff-cond-reg} is satisfied for $r < 1$. 
As we shall see below, 
this covers in particular all commutative cias 
(Theorem~\ref{turpi}, \cite{Tu70}), but in the latter case we simply have 
\[ \evol(u) = \exp\Big(\int_0^1 u(\tau)\, d\tau\Big).\] 
\end{rem}

\subsection{Turpin's Theorem concerning local m-convexity}

In this subsection we show Turpin's Theorem asserting that commutative 
cias are locally $m$-convex, which provides optimal control over the 
constants~$c(p,q)$ from the preceding subsection. 

\begin{prop}[Turpin's Theorem]\label{turpi}
Every commutative cia~$\cA$ is locally $m$-convex.
\end{prop}
Our proof uses two lemmas.
\begin{lem}\label{exc-esti-mono}
If $\cA$ is a complex cia, then 
for each continuous seminorm~$p$ on~$\cA$ with $p(\one)\leq 1$, there exists a continuous seminorm~$q$
on~$\cA$ such that
\begin{equation}\label{estihompolmult}
(\forall n\in\N_0)(\forall x\in \wb{B}^q_1(0))\quad p(x^n)\leq 1.
\end{equation}
\end{lem}
\begin{prf}
Given~$p$, let $q$ be a continuous seminorm on~$\cA$
such that $V:=\wb{B}^q_2(\1)\sub \cA^\times$ and $V^{-1}\sub \wb{B}^p_1(\one)\sub\wb{B}^p_2(0)$.
Let $\iota\colon \cA^\times \to \cA$, $x\mto x^{-1}$ be the inversion map,
which is complex analytic by Corollary~\ref{invsmoocia}. By Exercise~\ref{excTaycia}, the G\^{a}teaux differentials of~$\iota$
at~$\one$
are given by
\[
\delta^n_{\1}\iota=(-1)^n n!\mu_n,\quad 
\mbox{ where } \quad \mu_n(x)=x^n \quad \mbox{ for } \quad 
n \in \N, x\in \cA.\]
Hence Exercise~\ref{excgenesti} implies (\ref{estihompolmult}) for $n\in\N$,
and (\ref{estihompolmult}) also holds for $n=0$ since 
$p(\one)\leq 1$.
\end{prf}
\begin{lem}\label{convmulti}
If $\cA$ is an associative algebra and $S\sub \cA$ is a subset which is
closed under multiplication, then also the convex hull $\conv(S)$
is closed under multiplication.
\end{lem}
\begin{prf} Simply observe that products 
$(\sum_j \lambda_jx_j)(\sum_k \mu_k y_k) 
= \sum_{j,k} \lambda_j \mu_k x_j y_k$ of convex combinations are 
convex combinations of products. 
%
\end{prf}
%
%
%
\noindent
\emph{Proof of Proposition}~\ref{turpi}.
We may assume that the ground field is~$\C$ (if $\cA$ is a real cia,
discuss~$\cA_\C$; cf.~Proposition~\ref{autocxan}). 
Let $p$ be a continuous seminorm on~$\cA$;
after replacing $p$ with $rp$ for suitable $r>0$,
we may assume that $p(\one)\leq 1$.
For $n\in\N_0$,
consider the continuous symmetric $n$-linear map
\[
\beta_n\colon \cA^n\to \cA,\quad (x_1,\ldots, x_n)\mto x_1x_2\cdots x_n
\]
and the associated continuous homogeneous polynomial
\[
\mu_n\colon \cA\to \cA,\quad \mu_n(x):=\beta_n(x,x,\ldots,x)=x^n.
\]
By Lemma~\ref{exc-esti-mono}, there is a continuous seminorm~$q$
on~$\cA$ such that
\[
(\forall n\in\N_0)\;\; \sup\{p(x^n) \colon x\in \wb{B}^q_1(0)\}\leq 1.
\]
Using the Polarisation Formula~\ref{proppolarvar},
we deduce that
\[
\sup\{p(\beta_n(x_1,\ldots,x_n))\colon x_1,\ldots,x_n\in \wb{B}^q_1(0)\}\leq
\frac{n^n}{n!}
\]
(cf.\ Exercise~\ref{excnormpolvsmulti}).
Since $\sqrt[n]{n^n/n!}=n/\sqrt[n]{n!}\to e$ for $n\to\infty$ by Stirling's Formula, the power series
$\sum_{n=0}^\infty \frac{n^nt^n}{n!}$ has positive radius
of convergence, whence
we find $t>0$ such that
\[
(\forall n\in \N_0)\quad \frac{n^nt^n}{n!}\leq 1.
\]
Set $W:=\wb{B}^q_t(0)$.
For all $n\in \N$, we deduce that
\[
(\forall x_1,\ldots, x_n\in W)\quad p(x_1\cdots x_n)\leq \frac{n^nt^n}{n!}\leq 1.
\]
Hence
\[
S:=\bigcup_{n\in\N} \underbrace{WW\cdots W}_{n}
\]
is a balanced {$0$-neighbourhood} of~$\cA$ which is closed under multiplication
and contained in $\wb{B}^p_1(0)$.
Let $C:=\wb{\conv(S)}$ be the closed convex hull of~$S$.
Then $C$ is a closed, absolutely convex $0$-neighborhood in~$\cA$.
Let $\mu_C$ be its Minkowski functional.
Then $\mu_C$ is a continuous seminorm with unit ball $C\sub \wb{B}^p_1(0)$
and thus $\mu_C\geq p$. Since $CC\sub C$ by Lemma~\ref{convmulti},
$\mu_C$ is submultiplicative and so
submultiplicative seminorms define the topology of~$\cA$.\qed

\subsection{Completions of commutative cias}
Let $\cA$ be a locally convex topological algebra
and $\wt{\cA}$ be its completion as a locally convex space.
Then the continuous bilinear algebra multiplication $\cA\times \cA\to\cA$
extends to a continuous bilinear map $\wt{\cA}\times\wt{\cA}\to\wt{\cA}$,
which turns~$\wt{\cA}$ into a topological algebra. If~$\cA$
is unital, then also~$\wt{\cA}$ is unital.
\begin{probl}
Is the completion $\wt{A}$ of a cia also a cia~$\cA$?
\end{probl}
If $\cA$ is a commutative cia, then the answer is affirmative,
as we show now. The general case is open.
We begin with three lemmas.
\begin{lem}\label{spec-hom}
Let $\phi\colon\cA\to\cB$ be a unital algebra homomorphism
between unital complex algebras $\cA$ and $\cB$.
Given $a\in\cA$, let $\sigma_\cA(a)$ be its spectrum
and $\sigma_\cB(\phi(a))$ be the spectrum of $\phi(x)\in\cB$.
Then $\sigma_\cB(\phi(a))\sub \sigma_\cA(a)$.
\end{lem}
\begin{prf}
If $\lambda\in\sigma_\cB(\phi(a))$, then $\lambda \one-\phi(a)=\phi(\lambda\one-a)\not\in\cB^\times$,
entailing that $\lambda\one-a\not\in\cA^\times$ and thus $\lambda\in\sigma_\cA(a)$.
\end{prf}
\begin{lem}\label{spec-PL}
Let $((\cA_i)_{i\in I},(q_{ij})_{i\leq j})$ be a projective system of unital algebras,
and $\cA:=\pl \cA_i\sub\prod_{i\in I}\cA_i$,
with the limit maps $q_i\colon \cA\to\cA_i$.\vspace{-1.3mm}
Then
\[
\cA^\times=\pl\cA_i^\times.\vspace{-1.3mm}
\]
Moreover,
$\sigma_\cA(a)=\bigcup_{i\in I}\sigma_{\cA_i}(q_i(a))$ for each $a\in\cA$.
\end{lem}
\begin{prf}
If $a=(a_i)_{i\in I}\in\cA^\times$, then $a_i=q_i(a)\in\cA_i^\times$ for each $i\in I$.
Moreover,\vspace{-1mm}
\[
q_{ij}\Big|_{\cA_j^\times}^{\cA_i^\times}(a_j)=q_{ij}(a_j)=a_i\vspace{-1.3mm}
\]
whenever $i\leq j$, whence $a=(a_i)_{i\in I}\in\pl \cA_i^\times$.\vspace{-1.3mm}

If $a=(a_i)_{i\in I}\in\pl\cA_i^\times$,
then $a_i\in\cA_i^\times$ and $q_{ij}(a_j)=a_i$ whenever $i\leq j$,
whence $q_{ij}(a_j^{-1})=a_i^{-1}$ and hence $b:=(a_i^{-1})_{i\in I}\in\pl\cA_i=\cA$.\vspace{-1.3mm}
Since $ab=ba=\one$, we see that $a\in\cA^\times$.

To prove the final assertion, let $a=(a_i)_{i\in I}\in\cA$ 
and $\lambda\in\C$. By the preceding,
\[
\lambda\in\sigma_\cA(a)\aeq \lambda\one-a \not\in\cA^\times
\aeq (\exists i\in I)\; \lambda\one-a_i\not\in\cA_i^\times\aeq (\exists i\in I)\;\lambda\in\sigma_{\cA_i}(a_i).
\]
Hence $\sigma_\cA(a)=\bigcup_{i\in I}\sigma_{\cA_i}(a_i)$, as required.
\end{prf}
\begin{lem}\label{specr-comm}
For each unital commutative complex Banach algebra $\cA$,
the map $r\colon\cA\to[0,\infty[$, taking $a\in \cA$
to its spectral radius $r(a)$, is continuous.
\end{lem}
\begin{prf}
Let $\cA=\Hom(\cA,\C)$ be the Gelfand spectrum of~$\cA$
(the set of unital algebra homomorphisms to~$\C$)
and $\gamma_\cA\colon\cA\to C(\wh{\cA},\C)$
be the Gelfand homomorphism given by $\gamma_\cA(a)=\hat{a}$
with $\hat{a}(h)=h(a)$ for $a\in\cA$ and $h\in\wh{\cA}$.
It is well known from Gelfand theory
that $r(a)=\|\hat{a}\|_\infty$ for each $a\in\cA$.
As both $\|\cdot\|_\infty$ and $\gamma_\cA$ are continuous,
the continuity of $r=\|\cdot\|_\infty\circ\gamma_\cA$ follows.
\end{prf}
\begin{prop}
For any commutative cia~$\cA$, its completion~$\wt{\cA}$
is~a~cia.
\end{prop}
\begin{prf}
Write $\beta\colon \cA\times\cA\to\cA$ for the continuous bilinear
algebra multiplication.
Let $I$ be the set of all submultiplicative continuous seminorms on~$\cA$.
Given $p,q\in I$, write $p\leq q$ if $p(a)\leq q(a)$ for all~$a$.
As $a\mto \max\{p(a),q(a)\}$ is submultiplicative for all $p,q\in I$,
we see that $(I,\leq)$ is directed. Moreover, $rp\in I$ for all $p\in I$ and $r\geq 1$. 
Given $p\in I$, let $\cA_p:=\cA/p^{-1}(\{0\})$
be the corresponding normed space, $\alpha_p\colon \cA\to\cA_p$
be the canonical quotient map and $(\wt{\cA}_p,\|\cdot\|_p)$
be the corresponding Banach space. If $a\in\cA$ satisfies 
$p(a)=0$, then $p(ab)=p(ba)=0$ for all $b\in\cA$ follows from the 
submultiplicativity of~$p$. Hence $p^{-1}(\{0\})$
is an ideal in~$\cA$ and thus~$\cA_p$ is an algebra, with multiplication $\beta_p\colon\cA_p
\times\cA_p\to\cA_p$. Since $\beta_p\circ (\alpha_p\times\alpha_p)
=\alpha_p\circ \beta$
is continuous and $\alpha:= p\times\alpha_p$ 
a quotient map, we see that $\beta_p$
is continuous and hence extends to a continuous algebra multiplication on~$\wt{\cA}_p$.
Thus $\wt{\cA}_p$ is a unital Banach algebra. For $p\leq q$,
the canonical maps $\cA_q\to\cA_p$ extend to continuous unital algebra homomorphisms
$\alpha_{pq}\colon \wt{\cA}_q\to\wt{\cA}_p$. Then 
$((\wt{\cA})_{p\in I},(\alpha_{pq})_{p\leq q})$
is a projective system of 
Banach algebras and we can form its projective limit~$\wt{\cA}$.
For $p\in I$, let $\pi_p\colon\wt{\cA}\to\wt{\cA}_p$ be the limit map.
The canonical map
\[
\phi\colon \cA\to\pl \wt{\cA}_p=\wt{\cA},\quad a\mto (\alpha_p(a))_{p\in I}\vspace{-1.3mm}
\]
is a topological embedding with dense image, by 
Proposition~\ref{hand-on-PL}(a) and~(b).
The projective limit being complete, we see that $(\wt{\cA},\phi)$ is a completion of~$\cA$
(as the notation $\wt{\cA}$ for the projective limit already suggested).
Being a projective limit of unital Banach algebras, $\wt{\cA}$ has a continuous inversion map\vspace{.5mm}
\[
\iota_{\wt{\cA}}\colon\wt{\cA}^\times\to\wt{\cA},\vspace{-1mm}
\]
since
$\pi_p\circ \iota_{\wt{\cA}}=\iota_{\wt{A}_p}\circ \pi_p$
is continuous for each $p\in I$.

To see that $\wt{\cA}^\times$ is a $\one$-neighborhood in~$\wt{\cA}$ (and hence open),
let $U\sub\cA$ be an open $0$-neighborhood such that
\[
r_\cA(a)\leq\frac{1}{2}\quad \mbox{for all $\,a\in U$}
\]
(cf.\ Lemma~\ref{control-spec}).
Since $\phi$ is a topological embedding with dense image, the closure $\wb{\phi(U)}$
has $\phi(U)$ in its interior. Thus $\wb{\phi(U)}$ is a $0$-neighborhood in~$\wt{\cA}$
and the proof will be complete if we can show that
\begin{equation}\label{goal-now}
\one+\wb{\phi(U)}\, \sub\, \wt{\cA}^\times.
\end{equation}
Now
$r_{\wt{\cA}_p}(\alpha_p(a))\leq r_\cA(a)\leq\frac{1}{2}$ for all $a\in U$,
by Lemma~\ref{spec-hom}.
Thus $r_{\wt{\cA}_p}(\pi_p(a))\leq\frac{1}{2}$ for all $a\in \phi(U)$
(as $\pi_p\circ\phi=\alpha_p$) and hence
\[
r_{\wt{\cA}_p}(\pi_p(a))\leq\frac{1}{2}\quad\mbox{for all $\,a\in\wb{\phi(U)}$,}
\]
by continuity of the maps~$r_{\wt{\cA}_p}$ 
(see Lemma~\ref{specr-comm}) and~$\pi_p$.
Thus $\Spec_{\wt{\cA}_p}(\pi_p(a))\sub \wb{B}^\C_{1/2}(0)$ 
for all $a\in\wb{\phi(U)}$
and hence
\[
\Spec_{\wt{\cA}_p}(\one+\pi_p(a))\, \sub\, \wb{B}^\C_{1/2}(1)\, \sub \, \C\setminus\{0\},
\]
whence $\pi_p(\one+a)=\one+\pi_p(a)\in\wt{\cA}_p^\times$ for all $p\in I$.
Thus $\one+a\in \wt{\cA}^\times$, by Lemma~\ref{spec-PL},
and~(\ref{goal-now}) is established.
\end{prf}
\subsection{Regularity of algebras of smooth vectors}

Suppose that $\cA$ is a cia for which the unit group $\cA^\times$ is a regular 
Lie group and that $G$ is a Lie group acting on $\cA$ via 
$\alpha \: G \to \Aut(\cA), g \mapsto \alpha_g$, $g.a = \alpha_g(a)$ 
by topological automorphisms. 

Let $\cA^\infty \subeq \cA$ be the subspace of smooth vectors for this action 
and recall from the preceding subsection that it carries a natural 
cia structure with 
\[  (\cA^\infty)^\times = \cA^\times \cap \cA^\infty \] 
and a continuous inclusion 
$\iota \: \cA^\infty \to \cA$ 
(Proposition~\ref{prop:VI.7}). 

\begin{lem}
  \mlabel{lem:VI.13}If $\alpha \: G \times E \to E$ is a smooth action of $G$ 
on the locally convex space $E$ and $M$ is a finite-dimensional manifold (possibly with 
boundary), then the induced action $(g.f)(m) := g.f(m)$ of $G$ on the locally convex 
space $C^\infty(M,E)$ is smooth. 
\end{lem}

\begin{prf}
We recall from Lemma~\ref{lem:smooth-eval-vec} that the evaluation map 
$$ \ev \: C^\infty(M,E) \times M \to E $$
is smooth. Therefore the map 
$$ G \times C^\infty(M,E) \times M \to E, \quad 
(g,f,m) \mapsto (g.f)(m) = g.f(m) = g.\ev(f,m) $$
is smooth, as a composition of the evaluation map and the action map. 
Now  Proposition~\ref{prop:cartes-closed} shows that the action map 
$$ G \times C^\infty(M,E) \to C^\infty(M,E), \quad (g,f) \mapsto g.f $$
is also smooth. 
\end{prf}

\begin{thm}
  \mlabel{thm:VI.14} Assume that $\cA^\times$ is a regular Lie group and 
that the induced $G$-action on $\cA^\infty$ is smooth. Then 
the unit group of $\cA^\infty$ is also regular. 
\end{thm}

\begin{prf}
   Let $\xi \in C^\infty(I,\cA^\infty)$ be a smooth map. 
Then $\xi$ also is a smooth map to $\cA$, and the regularity of $\cA^\times$ 
implies the existence of a smooth path $\gamma \:I \to \cA^\times$ solving the 
linear initial value problem
$$ \gamma(0) = \1, \qquad \gamma'(t)= \gamma(t)\cdot \xi(t)\quad \mbox{ for } \quad 
t \in [0,1]. $$
To verify that $\gamma(I) \subeq \cA^\infty$, we first observe that 
the evolution map 
$$ \evol_\cA \: C^\infty(I,\cA) \to \cA^\times $$
is $\Aut(\cA)$-equivariant with respect to the natural actions, hence in 
particular $G$-equivariant. 

From Lemma~\ref{lem:VI.13} we know that the induced action of $G$ on  
$C^\infty(I,\cA^\infty)$ is smooth. Therefore the equivariance of 
$\evol_\cA$ implies that, for each $t \in I$, we have 
$\gamma(t) \in \cA^\infty$ (here the smoothness of the $G$-action on $\cA^\infty$ is used). 
The smoothness of the action of $G$ on $\cA^\infty$ implies that 
the map 
$$ G \times I \to \cA^\infty, \quad (g,t) \mapsto g.\xi(t) $$
is smooth. Applying $\evol_\cA$, we conclude that 
$$ G \times I \to \cA, \quad (g,t) \mapsto g.\gamma(t) = \alpha^{\gamma(t)}(g) $$
is smooth. Hence the map 
\[  I \to C^\infty(G,\cA), \quad t \mapsto \alpha^{\gamma(t)} \] 
is smooth (Proposition~\ref{prop:cartes-closed}), and this means that $\gamma \: I \to \cA^\infty$ is a smooth curve. 
Therefore the evolution map 
$$ \evol_{\cA^\infty} \: C^\infty(I,\cA^\infty) \to (\cA^\infty)^\times $$
is defined. To see that it is smooth, we may use Proposition~\ref{prop:cartes-closed} again, so that 
it suffices to see that the map  
\[ G \times C^\infty(I,\cA^\infty) \to \cA^\times, \quad 
(g,\xi) \mapsto g.\evol_{\cA}(\xi) = \evol_{\cA}(g.\xi)\] 
is smooth, but this follows from the smoothness of the action of $G$ on 
$C^\infty(I,\cA^\infty)$ and the smoothness of the composition 
\[ C^\infty(I,\cA^\infty) \into C^\infty(I,\cA) \to \cA^\times, \quad 
\xi \mapsto \evol_\cA(\xi).\qedhere\] 
\end{prf}

\begin{cor}
  \mlabel{cor:VI.15}  If $G$ is a finite-dimensional Lie group acting on the 
cia $\cA$ whose unit group $\cA^\times$ is a regular Lie group, then 
the unit group of the cia $\cA^\infty$ is also regular. 
\end{cor}

\begin{prf}
If $G$ is a finite-dimensional Lie group, then the 
right translation action of $G$ on the space $C^\infty(G,\cA)$ is smooth. 
In fact, the evaluation map $C^\infty(G,\cA) \times G \to \cA$ is smooth, and 
this implies that the map 
$$ G \times C^\infty(G,\cA) \times G \to \cA, \quad 
(g,f,x) \mapsto (g.f)(x) = f(xg) = \ev(f,xg) $$
is smooth, so that Proposition~\ref{prop:cartes-closed} shows that the $G$-action on 
$C^\infty(G,\cA)$ is smooth. Now Theorem~\ref{thm:VI.14} applies. 
\end{prf}

\begin{rem} If $G$ is a Banach--Lie group, then the smoothness of 
the $G$-action on $\cA^\infty$ follows from \cite[Thm.~4.4]{Ne10}, so that 
Theorem~\ref{thm:VI.14} also applies in this case. 
\end{rem}

\begin{cor} 
  \mlabel{cor:VI.16} The unit groups of smooth quantum tori 
{\rm(Example~\ref{ex:VI.8})} are regular Lie groups. 
\end{cor}

\begin{prop}
  \mlabel{prop:VI.17} Suppose that the unit group $\cA^\times$ of the 
cia $\cA$ is regular  
and $p \in \Idem(\cA)$. Then the unit group of the cia $p\cA p$ is also regular. 
\end{prop}

\begin{prf}
We consider the involution $\alpha_p := 2p - \1 \in \cA^\times$. 
Then the direct sum decomposition of 
\[ \cA = p\cA p \oplus p\cA(\1-p) \oplus (\1 -p)\cA p \oplus (\1-p)\cA(\1-p) \] 
implies that the centralizer of $\alpha_p$ in $\cA$ is the closed subspace 
$$ p\cA p \oplus (\1-p)\cA (\1-p) = \cA^{\alpha_p}$$
and the regularity of $\cA^\times$ easily implies the regularity of the product 
group 
$$ (p\cA p)^\times \oplus \big((\1-p)\cA (\1-p))^\times $$
from which one directly derives the regularity of the $(p\cA p)^\times$.
\end{prf}

\begin{probl}
Does the regularity of $\cA^\times$ imply the regularity of $\GL_n(\cA)$ 
for each $n \in \N$? \\ 
\end{probl}

\section{Appendix: Spectral theory of cias} 

\mlabel{sec:cia-1-app}

In this section, we develop the holomorphic functional
calculus for Mackey complete complex  cias 
with a view towards Lie-theoretic applications such as the 
BCH property of the group $\cA^\times$. 
Throughout the remainder of this section,
$\Omega$ denotes a non-empty, open subset of~$\C$ and 
$\cA$ a Mackey complete unital complex cia.
\subsection{Spectral calculus} 

\begin{defn}
We let $\cO(\Omega)$ denote
the space of holomorphic functions $f\!:\Omega\to \C$,
equipped with the topology of uniform convergence
on compact sets. Define
\[
\cA_\Omega:=\{x\in \cA\!:\Spec(x)\sub \Omega\}.
\]
\end{defn}
It is easy to see that
$\cO(\Omega)$ is a Fr\'echet algebra under pointwise operations 
because it is a closed subalgebra of the Fr\'echet algebra 
$C(\Omega,\C)$ (Lemma~\ref{sammelsu}(b),(c),(e)).
By Lemma~\ref{control-spec}, $\cA_\Omega$ is open in~$\cA$.
\begin{defn}
Recall that a formal finite $\Z$-linear combination
$\Gamma =\sum_\gamma n_\gamma\, \gamma$ (where almost all $n_\gamma=0$)
of smooth curves $\gamma\!:[a_\gamma,b_\gamma]\to \C$
is called a 
\index{contour} 
\index{contour!surrounding $K$} 
{\em contour\/} if 
$\sum \{ n_\gamma \: \gamma(a_\gamma) =z \} 
= \sum \{ n_\gamma \: \gamma(b_\gamma) = z\}$ for all $z\in\C$.
If $K\sub \C$ is a compact set, $U$ an open neighborhood of~$K$ in~$\C$,
and $\Gamma$ as before, we say that
$\Gamma$ is a contour {\em surrounding~$K$ in~$U$\/}
if $\im(\gamma)\sub U\,\setminus\, K$ whenever
$n_\gamma\not=0$, and
\[
\Ind_\Gamma(z)=\frac{1}{2\pi i}\oint_\Gamma \frac{d\zeta}{\zeta-z}
=
\begin{cases}
1 & \text{ if } z\in K \\ 
0 & \text{ if } z\in \C\,\setminus\, U.\end{cases}
\]
It is easy to see that a contour surrounding $K$ in~$U$ always exists (cf.~\cite{Ru91}).
\end{defn}
\begin{defn}\label{defncalculus}
Given
$f\in \cO(\Omega)$ and $x\in \cA_\Omega$,
we define
\begin{align*}
f[x]&:=\frac{1}{2\pi i}\oint_\Gamma f(\zeta)\,(\zeta  \1-x)^{-1}\, d\zeta \\
&:=\frac{1}{2\pi i} \sum_{\gamma} 
n_\gamma \cdot \int_{a_\gamma}^{b_\gamma} f(\gamma(t))\, \gamma'(t)
\, (\gamma(t)\1-x)^{-1}\, dt,
\end{align*}
where $\Gamma$ is a contour surrounding $\Spec(x)$ in~$\Omega$
(the integrals exist because $\cA$ is assumed to be Mackey complete). 
\end{defn}
Applying continuous linear functionals,
we deduce from the scalar-valued Cauchy Theorem 
that $f[x]$ is independent of the chosen contour.
\begin{prop}\label{eltyprops}
For fixed $x\in \cA_\Omega$,
the mapping
\[
\Phi\!: \cO(\Omega)\to \cA,\;\;\; f\mto f[x]\]
is a homomorphism of unital
algebras which takes
$J:=\id_\C|_{\Omega}$ to~$x$.
\end{prop}
\begin{proof}
It is obvious that $\Phi$ is linear.\\
\nin{\bf Step 1:} {\em We have $\Phi(J^k)=x^k$ for each $k\in \N_0$.}
Let $W\sub \cA^\times$ be a balanced open zero-neighborhood
in~$\cA$. There is $r>0$ such that $r^{-1} x\in W$
and $\Spec(x)\sub B_r(0)$.
The mapping $J^k \in \cO(\Omega)$
extends to the mapping $z\mto z^k$ in $\cO(\C)$;
in view of the independence of the chosen contour
in Definition~\ref{defncalculus},
it therefore suffices to assume $\Omega=\C$
in the present part of the proof,
enabling us to choose a positive
parametrization~$\gamma$ of the circle of radius~$2r$
around~$0$ as the contour surrounding $\Spec(x)$
used to define $\Phi(J^k)=J^k[x]$.
For fixed $\lambda\in \cA'$,
consider the complex analytic function
\[
f\!: B_{\frac{1}{r}}(0)\to \C,\;\;\; f(z)=\lambda((\1-zx)^{-1})
=\sum_{n=0}^\infty z^n\lambda(x^n)\]
(Proposition~\ref{autocpx}, Lemma~\ref{Neumann}).
The radius of convergence
of the power series being at least $\frac{1}{r}$,
the series converges absolutely and
uniformly on $\wb{B_{\frac{1}{2r}}(0)}$,
entailing that
\begin{eqnarray*}
\lambda\left(\int_\gamma \zeta^k (\zeta \1-x)^{-1}\, d\zeta\right)
& = & \int_\gamma \zeta^{k-1}\lambda((\1-\zeta^{-1}x)^{-1})\, d\zeta\\
& = & \int_\gamma \sum_{n=0}^\infty \zeta^{k-1-n}\lambda(x^n)\, d\zeta\\
& = & \sum_{n=0}^\infty \Big(\int_\gamma \zeta^{k-1-n}\, d\zeta\Big)\,\lambda(x^n)
=2\pi i \, \lambda(x^k)
\end{eqnarray*}
using that $\int_\gamma\zeta^{k-1-n}\, d\zeta=2\pi i \, \delta_{k,n}$.
The continuous linear functionals separating points
on~$\cA$, we deduce that $\Phi(J^k)=x^k$.

\nin{\bf Step 2:} {\em $\Phi$ is a homomorphism of unital algebras.\/}
As we already know that~$\Phi$ is linear and takes
$\1=J^0$ to $x^0=\1$, it only remains to show
that $\Phi$ is multiplicative. To see this, let $f,g\in \cO(\Omega)$.
There exists a relatively compact open neighborhood
$U$ of $\Spec(x)$ in~$\C$ whose closure
is contained in~$\Omega$.
We choose a contour $\Gamma_1$ surrounding $\Spec(x)$
in~$U$, and a contour~$\Gamma_2$ surrounding
$\wb{U}$ in~$\Omega$.
Then we use the First Resolvent Identity
\begin{equation}\label{resid}
r_x(z)r_x(\zeta)=(\zeta-z)^{-1}(r_x(z)-r_x(\zeta)) 
\end{equation}
to obtain 
\begin{eqnarray*}
f[x]g[x] & = & \left( \frac{1}{2\pi i}\int_{\Gamma_1}
f(z)r_x(z)\,dz\right) \left( \frac{1}{2\pi i}
\int_{\Gamma_2}g(\zeta)r_x(\zeta)\, d\zeta \right)\\
 & = & \frac{1}{(2\pi i)^2}
\int_{\Gamma_1}\int_{\Gamma_2}f(z)g(\zeta)r_x(z)r_x(\zeta)\,d\zeta dz\\
& = & 
\frac{1}{(2\pi i)^2}
\int_{\Gamma_1} f(z)
\left( \int_{\Gamma_2}\frac{g(\zeta)}{\zeta-z}\,d\zeta \right) r_x(z)\,dz\\
& &
+ \frac{1}{(2\pi i)^2}
\int_{\Gamma_2}g(\zeta)\left(
\int_{\Gamma_1}\frac{f(z)}{z-\zeta}\, dz\right) r_x(\zeta)\, d\zeta,
\end{eqnarray*}
where we use Exercise~\ref{excmorewint} twice
to obtain the second equality and 
Fubini's Theorem to obtain the last line.
Since $\int_{\Gamma_2}\frac{g(\zeta)}{\zeta-z}\, d\zeta=2\pi i \,g(z)$
by Cauchy's Integral Formula
(as $\Ind_{\Gamma_2}(z)=1$)
and $\int_{\Gamma_1}\frac{f(z)}{z-\zeta}\, dz=0$
as $\Ind_{\Gamma_1}(\zeta)=0$,
we obtain 
\[ f[x]g[x]=\frac{1}{2\pi i}\int_{\Gamma_1}f(z)g(z)r_x(z)\, dz=
(f\cdot g)[x]. \qedhere \] 
\end{proof}
\begin{thm}\label{calcanalyt}
For every non-empty open subset $\Omega\sub \C$,
the mapping
\[
\Xi\!: \cO(\Omega)\times \cA_\Omega\to \cA,\;\;\; \Xi(f,x):=
f[x] \] 
is complex analytic. In particular,
the mapping $\cA_\Omega\to \cA$, $x\mto f[x]$ is complex
analytic  for each $f\in \cO(\Omega)$, and, for any $x\in \cA$,
the homomorphism~$\Phi$ {\rm(as in Proposition~\ref{eltyprops})} 
is continuous.
\end{thm}
\begin{proof}
Let
$(f,x)\in \cO(\Omega)\times \cA_\Omega$.
There is a relatively compact open neighborhood~$U$
of $\Spec(x)$ in~$\C$, such that $\wb{U}\sub \Omega$.
Then $\cA_U$ is an open neighborhood of~$x$ in~$\cA$. 
Fix a contour $\Gamma=\sum_\gamma n_\gamma \gamma$ surrounding $\wb{U}$
in~$\Omega$ and note that, for all $(f,x) \in \cO(\Omega) \times \cA_U$, we have  
\[ \Xi(f,x)
= \frac{1}{2\pi i}\oint_\Gamma f(\zeta)\,(\zeta  \1-x)^{-1}\, d\zeta. \]
Therefore it suffices to show that, for every path 
$\gamma \: [a,b] \to \Omega$ occuring in $\Gamma$, the map 
\[ \cO(\Omega) \times \cA_U \to \cA,\ 
(f,x) \mapsto \int_a^b \Psi(t,f,x)\, dt, \quad 
\Psi(t,f,x) := f(\gamma(t))(\gamma(t)\1 - x)^{-1} \] 
is complex analytic.

First we observe that the function 
\[ \cO(\Omega) \to C^1([a,b],\C), \quad f \mapsto f\circ \gamma \] 
is continuous because $\gamma([a,b])$ is a compact subset of $\Omega$. 
Further, continuity of inversion in $\cA$ implies that 
\[ [a,b] \times \cA_U \to \cA, \quad (t,x) \mapsto 
(\gamma(t)\1 - x)^{-1} \] 
is continuous. Now the continuity of multiplication in $\cA$ 
implies that $\Psi$ is continuous. For each fixed $t \in [a,b]$, 
the function $\Psi_t(f,x) :=  \Psi(t,f,x)$ is complex analytic with the differential 
\begin{align}\label{eq:derivPhi}
& d\Phi_t((f,x),(g,y)) \notag \\ 
&= g(\gamma(t))(\gamma(t)\1 - x)^{-1} 
+ f(\gamma(t))(\gamma(t)\1 - x)^{-1} y (\gamma(t)\1 - x)^{-1}.
\end{align}
Since all these functions are $C^1$ in $t$ and $\cA$ is Mackey 
complete, Proposition~\ref{diffpar}, combined with 
Theorem~\ref{charcxcompl}, now shows that $\Xi$ is 
complex analytic on the open neighborhood $C(\Omega) \times \cA_U$ of $(f,x)$ 
because $\cA$ is Mackey complete. 
As $(f,x)$ was arbitrary in $\cO(\Omega) \times \cA_\Omega$, 
this proves the theorem. 
\end{proof}
\begin{cor}\label{pseries}
Let $f\!: B_r(z_0)\to\C$ be a holomorphic function,
defined on a disk $B_r(z_0)\sub \C$,
with power series expansion
$f(z)=\sum_{n=0}^\infty a_n (z-z_0)^n$.
Then
\[
(\forall x\in \cA_{B_r(z_0)})\;\;\;\;\;
f(x)=\sum_{n=0}^\infty a_n (x-z_0\1)^n\;.
\]
\end{cor}
\begin{proof}
We have $f_k \to f$ in $\cO(B_r(z_0))$,
where 
$f_k(z):=\sum_{j=0}^k a_j\, (z-z_0)^j.$ 
Thus $f[x]=\lim_{k\to\infty} f_k[x]$
for every $x\in \cA_{B_r(z_0)}$,
by
Theorem~\ref{calcanalyt}.
It remains to observe that
$f_k[x]=\sum_{j=0}^k a_j\, (x-z_0e)^j$,
by Proposition~\ref{eltyprops}.
\end{proof}
The familiar Spectral Mapping Theorem carries
over to Mackey complete cias. 
\begin{prop}\label{specmap}
Let $x\in \cA_\Omega$ and $f\in \cO(\Omega)$. Then
\begin{description}
\item[\rm(a)]
$f[x]$ is invertible in~$\cA$ if and only if
$f(z)\not=0$ for all $z\in \Spec(x)$.
\item[\rm(b)]
$\Spec(f[x])=f(\Spec(x))$.
\end{description}
\end{prop}
\begin{proof}
In view of Proposition~\ref{eltyprops}
above, the proof of the Banach case
as formulated in~\cite[Thm.~10.28]{Ru91} can be repeated verbatim.
\end{proof}
Concerning compositions, the usual proof of the Banach case
(see \cite[Thm.~10.29]{Ru91}) applies without changes. We obtain:
\begin{prop}\label{compos}
Let $x\in \cA_\Omega$, $f\in \cO(\Omega)$,
$\Omega_1$ be an open neighborhood~of
$f(\Spec(x))$ in~$\C$, and $g\in \cO(\Omega_1)$.
Set $h:=g\circ (f|_{\Omega_0}^{\Omega_1})$,
where $\Omega_0:=f^{-1}(\Omega_1)$.
Then $f[x]\in \cA_{\Omega_1}$, and $h[x]=g[f[x]]$.
\end{prop}
For the final proposition in this section,
let $\cA$ be a Mackey complete,
{\em real\/} continuous inverse
algebra. Given $x\in \cA$,
we let $\Spec(x)$ be the spectrum of~$x$,
considered as an element of~$\cA_\C$.
Let $\kappa\!: \C\to\C$
and
\[ \tau\!: \cA_\C\to \cA_\C,\quad
  \tau(x+iy)=\wb{x+iy}:=x-iy \quad \mbox{ for } \quad x,y\in \cA \] 
denote the respective complex conjugation.
Given an open non-empty subset $\Omega$
of~$\C$ invariant under complex conjugation~$\kappa$,
we set $\cA_\Omega:=(\cA_\C)_\Omega\cap \cA$
and define $f[x]\in \cA_\C$ using the holomorphic functional
calculus for~$\cA_\C$ when $f\in \cO(\Omega)$ and $x\in \cA_\Omega$.
Clearly $f^*(z):=\wb{f(\wb{z})}$
({\em i.e.}, $f^*:=\kappa\circ f\circ \kappa|_\Omega^\Omega$)
defines a continuous involution
$*\!: \cO(\Omega)\to \cO(\Omega)$
making $\cO(\Omega)$ a $*$-algebra.
\begin{prop}\label{realcase}
For $x\in (\cA_\C)_\Omega$ and $f\in \cO(\Omega)$,
we have 
\[ \tau(f[x])=f^*[\tau(x)].\] 
In particular,
$f[x]\in \cA$ for all $x\in \cA_\Omega$
and hermitian elements $f=f^* \in \cO(\Omega)$.
\end{prop}
\begin{proof}
It is easily verified that $\tau$ is a complex antilinear, continuous
real unital algebra automorphism of~$\cA_\C$.
We readily deduce that
$\Spec(\tau(x))=\wb{\Spec(x)}\sub \Omega$
for
$x\in (\cA_\C)_\Omega$.
For $z\in \rho(x)$,
we have 
\[ \tau(r(x,z))=\tau((z\1-x)^{-1})=(\wb{z}\1-\tau(x))^{-1}=
r(\tau(x),\wb{z}).\]
Let $\Gamma=\sum_\gamma n_\gamma \gamma$
be a contour surrounding $\Spec(x)$ in~$\Omega$.
Then
\[ \Gamma_2:=-\sum_\gamma n_\gamma \,\kappa\circ \gamma \] 
is a contour surrounding
$\wb{\Spec(x)}=\Spec(\tau(x))$ in~$\Omega$.
Thus
\begin{eqnarray*}
&&\tau(f[x]) = 
\sum_\gamma \frac{-n_\gamma}{2\pi i} \int_{a_\gamma}^{b_\gamma}
\wb{f(\gamma(t))}\,\wb{\gamma'(t)}\,r(\tau(x),\wb{\gamma(t)})\, dt\\
& = &
\sum_\gamma \frac{-n_\gamma}{2\pi i} \int_{a_\gamma}^{b_\gamma}
f^*((\kappa\circ \gamma)(t))\,(\kappa\circ\gamma)'(t)\,
r(\tau(x),(\kappa\circ \gamma)(t))\, dt=f^*[\tau(x)],
\end{eqnarray*}
as asserted. The remainder is an immediate consequence.
\end{proof}
Note that $\cO(\Omega)=\Herm(\cO(\Omega))_\C$ in the preceding situation,
where $\Herm(\cO(\Omega)):=\{f\in \cO(\Omega)\!: f=f^*\}$
is the real vector subspace of hermitian elements in $\cO(\Omega)$.
Since $\Psi\!:\Herm(\cO(\Omega))\times \cA\to \cA$,
$(f,x)\mto f[x]$ extends to a complex analytic mapping
$(\Herm(\cO(\Omega))\times \cA)_\C=\cO(\Omega)\times \cA_\C\to \cA_\C$
(Theorem~\ref{calcanalyt}),
the mapping~$\Psi$ is real analytic. 
\subsection{Exponential and logarithm function} 
\label{exp}
\begin{defn}\mlabel{def:exp-log}
Using the exponential function
\[ \exp\!:\C\to\C, \qquad \exp(z) := \sum_{k=0}^\infty \frac{z^k}{k!} \] 
and logarithm
\[ \log\!: B_1(1)\to \C, \qquad \log(z):=\sum_{k=1}^\infty (-1)^{k+1}\frac{(z-1)^k}{k},\]
we define mappings
$\exp_\cA\!: \cA\to \cA$, $\exp_\cA(x):=\exp[x]$
and $\log_\cA\!: \cA_{B_1(1)}\to \cA$,
$\log_\cA(x):=\log[x]$
using the holomorphic functional
calculus described in the preceding section 
(Proposition~\ref{realcase}). 
\end{defn}
As a consequence of Theorem~\ref{calcanalyt},
in the case $\K=\C$, 
$\log_\cA$ is a complex analytic mapping, and $\exp_\cA$ is complex analytic
on the open subset $\cA_{B_r(0)}$ of $\cA$
for each $r>0$ and thus complex analytic. 
If $\K=\R$, then $\exp_{\cA_\C}$ and $\log_{\cA_\C}$
are complex analytic mappings between open subsets
of the respective complexified algebras (Proposition~\ref{autocxan}) 
which extend $\exp_\cA$, resp., $\log_\cA$,
and thus $\exp_\cA$ and $\log_\cA$ are real analytic mappings.
Corollary~\ref{pseries}
shows that
\begin{equation}\label{f1}
\!\!\!\!(\forall x\in \cA)\;\;\;
\exp_\cA(x)=\sum_{k=0}^\infty \frac{1}{k!}\, x^k\;\;\;\;\;
\mbox{and}
\end{equation}
\vspace{-2 mm}
\begin{equation}\label{f2}
(\forall x\in \cA_{B_1(1)})\;\;\;
\log_\cA(x)=\sum_{k=1}^\infty \frac{(-1)^{k+1}}{k}\, (x-\1)^k.
\end{equation}
\begin{lem} \mlabel{lem:logexp}
We have
\begin{enumerate}
\item[\rm(a)]
$\exp_\cA(\log_\cA(x))=x$ for all $x\in \cA_{B_1(1)}$ and
\item[\rm(b)]
$\log_\cA(\exp_\cA(x))=x$ for all $x\in \cA_{B_{\log(2)}(0)}$.
\item[\rm(c)]
Let $V:=\log_\cA^{-1}(\cA_{B_{\log(2)}(0)})$
and $U:=\{x\in \cA_{B_{\log(2)}(0)}\!:\exp_\cA(x)\in V\}$.
Then $V$ is an open identity neighborhood
in~$\cA$, $U$ is an open zero-neighborhood
in~$\cA$, 
and $\exp_\cA|_U^V\!: U\to V$ is a
$\K$-analytic diffeomorphism, with inverse $\log_\cA|_V^U$.
\end{enumerate}
\end{lem}
\begin{proof}
(a) It is well-known that $\exp(\log(z))=z$ for all
$z\in B_1(1)\sub \C$. Proposition~\ref{compos}
and Proposition~\ref{eltyprops}
entail that $\exp_\cA(\log_\cA(x))=x$ for all
$x\in \cA_{B_1(1)}$.

(b) It is well-known that $|\exp(z)-1|<1$ for all
$z\in B_{\log(2)}(0)$,
and $\log(\exp(z))=z$ for such~$z$.
Proposition~\ref{specmap}
shows that $\exp_\cA(\cA_{B_{\log(2)}(0)})\sub \cA_{B_1(1)}$.
Thus $\log_\cA(\exp_\cA(x))$ $=x$
for all $x\in \cA_{B_{\log(2)}(0)}$,
by Proposition~\ref{compos}
and Proposition~\ref{eltyprops}.

(c) $\log_\cA$ and $\exp_\cA|_{\cA_{B_{\log(2)}(0)}}$
being continuous mappings defined
on open subsets of~$\cA$, clearly $U$ and~$V$ are open.
Since $\log_\cA(\1)=0$ and $\exp_\cA(0)=\1$ in view of (\ref{f1})
and (\ref{f2}), $\1\in V$ and $0\in U$ hold.
Given $x\in V$, we have
$\exp_\cA(\log_\cA(x))=x\in V$
by (a), where
$\log_\cA(x)\in \cA_{B_{\log(2)}(0)}$
by definition of~$V$.
Thus $\log_\cA(x)\in U$, by definition
of~$U$. We have proved that
$V\sub \exp_\cA(U)$.
Since $\exp_\cA(U)\sub V$ by definition of~$U$,
we deduce that ${\exp_\cA(U)=V}$.
We have also shown that $\log_\cA(V)\sub U$,
and
\begin{equation}\label{f3}
\exp_\cA|_U^V\circ \log_\cA|_V^U=\id_V.
\end{equation}
By the preceding, $\exp_\cA|_U^V\!:U\to V$ is surjective.
By (b), $\exp_\cA|_U^V$ is also injective
and thus a bijection. Composing with $(\exp_\cA|_U^V)^{-1}$
in (\ref{f3}), we find that
$\log_\cA|_V^U=(\exp_\cA|_U^V)^{-1}$. Both $\exp_\cA|_U^V$ and
$\log_\cA|_V^U$ being $\K$-analytic maps,
$\exp_\cA|_U^V$ is a diffeomorphism
of $\K$-analytic manifolds.
\end{proof}
\begin{lem}
We have $\im(\exp_\cA)\sub \cA^\times$ 
and the co-restriction \break $\exp_\cA|^{\cA^\times}\!: \cA\to \cA^\times$ is the exponential
function of the Lie group $\cA^\times$,
i.e., for every $x\in \cA$,
the mapping $\gamma_x\!:\R\to \cA^\times$, $\gamma_x(t):=\exp_\cA(tx)$
is a smooth homomorphism such that $\gamma_x'(0)=x$.
\end{lem}
\begin{proof}
Let $x\in \cA$.
The mapping $\exp_\cA$ being $\K$-analytic and thus smooth,
$\gamma_x$ is a smooth curve; we have $\gamma_x(0)=\1$.
Given $r\in \R$,
define
$\mu_r\!:\C\to\C$, $\mu_r(z):=r\cdot z$
and $e_r:=\exp\circ \mu_r\!:\C\to \C$,
$e_r(z)=\exp(rz)$.
Then $\mu_r[x]=r\cdot x$ by Proposition~\ref{eltyprops}
and thus, using Proposition~\ref{compos},
$e_r[x]=\exp[\mu_r[x]]=\exp_\cA(rx)$.

Given $r,s\in \R$, we have $e_r\cdot e_s=e_{r+s}$
and thus 
\[ \gamma_x(r+s)=\exp_\cA((r+s)x)=e_{r+s}[x]=e_r[x]e_s[x]
=\gamma_x(r)\gamma_x(s) \]  by Proposition~\ref{eltyprops}.
In particular, 
\[ \gamma_x(r)\gamma_x(-r)=\gamma_x(-r)\gamma_x(r)=\gamma_x(0)=\1,\]
showing that $\gamma_x(r)\in \cA^\times$
and $\gamma_x(r)^{-1}=\gamma_x(-r)$. Thus $\gamma_x\!: \R\to
\cA^\times$ is a homomorphism.
Using Equation~\ref{eq:derivPhi}, we see that 
\[ \gamma_x'(0)=d\exp_\cA(0;x)=\frac{1}{2\pi i}
\oint_\Gamma \frac{\exp(z)}{z^2}x\, dz=\exp'(0)x=x,\]
where $\Gamma$ is a contour
surrounding $\{0\}$.
\end{proof}
\begin{lem}
  \mlabel{lem:a.1-ne04}
Let $\cA$ be a complex unital Banach algebra and  consider the right 
open half plane $\Omega := \{ z \in \C \: \Re z > 0\}.$ 
Then we obtain a holomorphic map 
\[ \log \: \cA_\Omega \to \cA.\] 
If $*$ is an antilinear antiautomorphism of 
$\cA$, then $\log(x^*) = \log(x)^*$ for $x \in~\cA_\Omega$. 
\end{lem}

\begin{prf} By Lemma~\ref{control-spec}, $\cA_\Omega$ is open, 
and the holomorphy of $\log$ from Theorem~\ref{calcanalyt}. 
The relation $\log(a^*) = \log(a)^*$ is obtained as in 
the proof of Proposition~\ref{realcase}.
\end{prf}

\begin{small}
\subsection*{Exercises for Chapter~\ref{ch:lingrp}}

\begin{exer}
Let $\cA$ be a cia and $\cB \subeq \cA$ be a unital 
subalgebra. If $\cB^\times = \cA^\times \cap \cB$, then $\cB$ is also a 
cia.   
\end{exer}

\begin{exer} \mlabel{exer:I.1} For an associative 
$\K$-algebra $\cA$, we write $\cA_+$ for the algebra $\cA \times \K$ 
with the multiplication 
$$ (a,s)(b,t) := (ab + s b + ta,st). $$
\begin{description}
\item[\rm(1)] Verify that $\cA_+$ is a unital algebra with unit $\1 = (0,1)$. 
\item[\rm(2)] Show that $\GL_1(\cA) := \cA_+^\times \cap (\cA \times \{1\})$ is a group. 
\item[\rm(3)] If $e \in \cA$ is an identity element, then $\cA_+$ is isomorphic to the 
direct product algebra $\cA \times \K$ with the product 
$(a,s)(b,t) = (ab,st)$. 
\end{description}
\end{exer}

\begin{exer}\mlabel{exer:I.2.} A 
\index{topological ring} 
{\it topological ring} is a ring $\cR$ endowed with a topology 
for which addition and multiplication are continuous. 
Let $\cR$ be a unital topological ring. Show that:
\begin{description}
\item[\rm(1)] For $x \in \cR^\times$, the left and right multiplications 
$\lambda_x(y) := xy$ and $\rho_x(y) := yx$ are homeomorphisms of $\cR$. 
\item[\rm(2)] The unit group $\cR^\times$ is open if and only if it is a neighborhood of $\1$. 
\item[\rm(3)] The inversion $\cR^\times \to \cR$ is continuous, i.e., 
$(\cR^\times, \cdot)$ is a topological group, if it is continuous in~$\1$. 
\end{description}
\end{exer} 

\begin{exer}\label{cia-vs-quasi}
Let $\cA$ be a locally convex topological $\K$-algebra
over $\K\in \{\R,\C\}$ (which may not be unital)
and $q\colon Q(\cA)\to\cA$ be its quasi-inversion map.
Show the following:
\begin{description}[(D)]
\item[(a)]
$Q(\cA)$ is a $0$-neighborhood in $\cA$ if and only if
$(\cA_+)^\times$ is a $\1$-neighborhood on $\cA_+$.
\item[(b)]
$q$ is continuous at $0$ if and only if $\iota\colon (\cA_+)^\times\to \cA_+$,
$a\mto a^{-1}$ is continuous at~$\1$.
\end{description}
\end{exer}

\begin{exer} \mlabel{exer:I.3} Let $\cR$ be a unital ring with $2 \in \cR^\times$, 
$n \in \N$ and 
$M_n(\cR)$ the ring of all $(n \times n)$-matrices with entries in $\cR$. 
In the following, we write elements $x \in M_n(\cR)$ as 
$$ x = \pmat{ a & b \cr c & d \cr} 
\in M_n(\cR) = \pmat{ M_{n-1}(\cR) & M_{n-1,1}(\cR) \cr M_{1,n-1}(\cR) & \cR \cr}
= \pmat{ M_{n-1}(\cR) & \cR^{n-1}\cr (\cR^{n-1})^\top & \cR \cr}. $$
\begin{description}
\item[\rm(1)] Show that a matrix $x$ is of the form 
$$ \pmat{\1 & \beta \cr 0 & \1} \pmat{\alpha & 0 \cr 0 & \delta} 
\pmat{\1 & 0 \cr \gamma & \1} \quad \hbox{ with } 
\quad \alpha \in \GL_{n-1}(\cR), \beta, \gamma^\top \in \cR^{n-1},\delta \in \cR^\times $$
if and only if $d \in \cR^\times, a - bd^{-1}c \in \GL_{n-1}(\cR)$, and that in this case 
$$ \delta = d, \quad \beta = bd^{-1}, \quad \gamma = d^{-1}c, 
\quad \alpha  = a - bd^{-1}c. $$
\item[\rm(2)] Assume, in addition, that $\cR$ is a topological ring with open unit group and 
continuous inversion. Show by induction on $n$ that 
\begin{description}
\item[\rm(a)] $\GL_n(\cR)$ is open in $M_n(\cR)$. 
\item[\rm(b)] Inversion in $\GL_n(\cR)$ is continuous, i.e., 
$\GL_n(\cR)$ is a topological group. 
\end{description}
\end{description}
\end{exer} 

\begin{exer}\mlabel{exer:I.4} Let $\cR$ be a unital ring and consider the 
right $\cR$-module $\cR^n$, where the module structure is given by 
$(x_1,\ldots, x_n)r := (x_1r,\ldots, x_nr)$. 
Let $M$ be a right $\cR$-module, $\sigma \: r \mapsto r^\sigma$ 
an involution on $\cR$, i.e., an involutive anti-automorphism and 
$\eps \in \{\pm 1\}$. 
A biadditive map $\beta \: M \times M \to \cR$ is called 
\index{$\sigma$-sesquilinear map} 
{\it $\sigma$-sesquilinear} if 
$$ \beta(xr,ys) = r^\sigma \beta(x,y) s \quad \hbox{ for } \quad 
x,y \in M, r,s \in \cR. $$
It is called 
\index{$\sigma$-$\eps$-hermitian map} 
\index{$\sigma$-hermitian map} 
\index{$\sigma$-antihermitian map} 
\index{unitary group} 
{\it $\sigma$-$\eps$-hermitian} if, in addition, 
$$ \sigma(x,y)^\sigma = \eps \sigma(y,x)\quad \hbox{ for }  \quad x,y \in M. $$ 
For $\eps = 1$, we call the form {\it $\sigma$-hermitian}, and 
{\it $\sigma$-antihermitian} for $\eps = -1$. 
For a $\sigma$-$\eps$-hermitian form $\beta$ on $M$,  
\[  \U(M,\beta) := \{ \phi \in \Aut_\cR(M) \: 
(\forall x,y \in M)\ \beta(\phi(x),\phi(y)) = \beta(x,y)\} \]
is called the corresponding {\it unitary group}.  
Show that: 
\begin{description}
\item[\rm(1)] $\End_\cR(\cR^n) \cong M_n(\cR)$, where $M_n(\cR)$ acts 
by left multiplication on the space $\cR^n$ of column vectors. 
\item[\rm(2)] $\Aut_\cR(\cR^n) \cong \GL_n(\cR)$. 
\item[\rm(3)] $\beta(x,y) := \sum_{i=1}^n x_i^\sigma y_i$ is a 
$\sigma$-hermitian form on $\cR^n$. Describe the corres\-ponding 
unitary group in terms of matrices. 
\item[\rm(4)] $\beta(x,y) := \sum_{i=1}^n x_i^\sigma y_{n+i} - x_{n+i}^\sigma y_i$ is a 
$\sigma$-antihermitian form on $\cR^{2n}$. Describe the corresponding 
unitary group in terms of matrices. 
\end{description}
\end{exer} 

\begin{exer}\mlabel{exer:6.3.1} 
Show that, in general, for a closed subspace $F$ of a Banach space $E$, 
the set $\{ g \in \GL(E) \: gF \subeq F \}$
is not a subgroup. 
\end{exer}

\begin{exer} \mlabel{exer:5.6.13} Let $\cA$ be a commutative continuous inverse algebra and 
$$ \det \: M_n(\cA) \to \cA $$
the determinant function. Show that $\det$ restricts to a smooth morphism of 
Lie groups $\det \: \GL_n(\cA) \to \cA^\times$ with 
$\L(\det) = \tr$. 
\end{exer} 

\begin{exer} \mlabel{exer:5.6.5} Let $E$ be a Banach space and $v \in E$. 
Show that 
 $$\GL(E)_v := \{ g \in \GL(E) \: gv = v\}$$ 
is a Lie subgroup with 
$$ \L^e(\GL(E)_v) = \{ Y \in {\cal L}(E) \: Yv = 0\}. $$
\end{exer} 

\begin{exer} \mlabel{exer:5.6.6} 
Let $\cA$ be a Banach space and $m \: \cA \times \cA \to \cA$
a continuous bilinear map. Then the group 
$$ \Aut(\cA,m) := \{ g \in \GL(\cA) \: (\forall a,b \in \cA)\ m(ga,gb) =
gm(a,b)\} $$
of automorphisms of the (not necessarily associative) algebra $(\cA,m)$ is a Lie group whose Lie
algebra is the space 
$$ \der(\cA,m) := \{ X \in {\cal L}(\cA) \: (\forall a,b \in \cA)\, Xm(a,b) =
m(Xa,b) + m(a,Xb)\} $$
of derivations of $(\cA,m)$.\\
 Hint: Algebraic Subgroup Theorem 
(cf.\ also Proposition~\ref{prop:aut-der}). 
\end{exer}

\begin{exer}
  \mlabel{exer:oden-iv.6}
Let $E$ and $F$ be Banach spaces and 
$\cL^k(E,F)$ be the space of continuous $k$-linear maps $E^k \to F$. 
Show that: 
\begin{enumerate}
\item[\rm(a)] $\cL^k(E,F)$ is a Banach space with respect to the norm 
$$ \|f\| := \sup \{ \|f(x_1, \ldots, x_k)\| \: x_i \in E, \|x_1\|, \ldots,
\|x_k\| \leq 1\}. $$
\item[\rm(b)] The assignment 
$$ \big(\pi(g)f\big)(x_1, \ldots, x_k) := f(g^{-1}x_1, \ldots,
g^{-1}x_k) $$
defines a continuous homomorphism $\pi \: \GL(E) \to \GL(\cL^k(E,F))$. \\
Hint: The map $\eta \: \cL(E)^k \to \cL(\cL^k(E,F))$ with 
$$\big(\eta(A_1,\ldots, A_k)f\big)(x_1, \ldots, x_k) 
:= f(A_1x_1, \ldots,A_kx_k) $$
is a continuous $k$-linear map. 
\item[\rm(c)] Calculate the derived representation 
$\L(\pi) \: \gl(E) \to \gl(\cL^k(E,F))$. 
\item[\rm(d)] We identify the space $P_k(E,F)$ of $F$-valued continuous
polynomial functions of degree $k$ on $E$ with the closed subspace 
$\Sym^k(E,F) \subeq \cL^k(E,F)$. Then a compatible norm on this space 
is given by
$$ \|f\| = \sup \{ \|f(x)\| \: \|x\| \leq 1\} $$
and the assignment 
$\big(\pi(g)f\big)(x) := f(g^{-1}x)$
defines a continuous homomorphism $\pi \: \GL(E) \to \GL(P_k(E,F))$. 
\end{enumerate}
\end{exer}

\begin{exer}
  \mlabel{exer:oden-iv.7} 
  \begin{enumerate}
\item[\rm(a)] Let $\cH$ be a complex Hilbert space. Show that there
exists an antilinear isometric map $I \: \cH \to \cH$ with $I^2 = \1$. 

\item[\rm(b)] If $I_1$ and $I_2$ are two such maps, then there exists a
unitary operator $g \in \U(\cH)$ with $I_2 = g I_1 g^{-1}$. 

\item[\rm(c)] Show that, for a complex Hilbert space $\cH$, all groups
\[ \OO(\cH,I) := \{ g \in \GL(\cH) \: g^* = Ig^{-1}I\} \] 
are isomorphic. 
  \end{enumerate}
\end{exer} 

\begin{exer}
  \mlabel{exer:oden-iv.8}
  \begin{enumerate}
  \item[\rm(a)] Let $\cH$ be an infinite-dimensional or
even-dimensional complex Hilbert space. Show that there
exists an antilinear isometric map \break 
$I \: \cH \to \cH$ with $I^2 = -\1$. 

\item[\rm(b)] If $I_1$ and $I_2$ are two such maps, then there exists a
unitary operator $g \in \U(\cH)$ with $I_2 = g I_1 g^{-1}$. 

\item[\rm(c)] Show that for a complex Hilbert space $\cH$, all groups
\[ \Sp(\cH,I) := \{ g \in \GL(\cH) \: g^* = Ig^{-1}I^{-1}\} \] 
are isomorphic. 
  \end{enumerate}
\end{exer}

\begin{exer}
  \mlabel{exer:oden-iv.9}
Let $\cH$ be a complex Hilbert space and $I$ an
antilinear isometry with $I^2 = \pm \1$. We consider the complex
bilinear form 
$$ \beta(v,w) := \la v, Iw \ra. $$
\begin{enumerate}
\item[\rm(a)] $\beta$ is symmetric (skew-symmetric) if $I^2 = \1$ ($I^2
= -\1$). 
\item[\rm(b)] For $I^2 = \1$, 
$$ \OO(\cH,I) = \{ g \in \GL(\cH) \: (\forall v, w \in \cH)\, \beta(gv, gw) 
= \beta(v,w)\}, $$
and the Lie algebra of this subgroup is 
\begin{align*}
\fo(\cH,I) 
&= \{ X \in \cL(\cH) \: IX^*I^{-1} + X = 0\} \cr
&= \{ X \in \cL(\cH) \: (\forall v, w \in \cH)\, \beta(Xv,w) +\beta(v,Xw)=0\}. 
\end{align*}
\item[\rm(c)] For $I^2 = -\1$, 
$$ \Sp(\cH,I) = \{ g \in \GL(\cH) \: (\forall v, w \in \cH)\, \beta(gv, gw) 
= \beta(v,w)\}, $$
and the Lie algebra of this subgroup is 
\begin{align*}
\fsp(\cH,I) 
&= \{ X \in \cL(\cH) \: IX^*I^{-1} + X = 0\} \cr
&= \{ X \in \cL(\cH) \: (\forall v, w \in \cH)\, \beta(Xv,w) +
\beta(v,Xw)=0\}. 
\end{align*}
\item[\rm(d)] For $I^2 = \1$ and $\dim \cH = \infty$, there exists  
an orthonormal basis $(e_j)_{j \in 2J}$ of $\cH$ with 
$Ie_j = e_{-j}$ for  $j \in 2J := J \dot\cup - J$. Then 
$$\cH \cong \ell^2(2J,\C) 
\cong \ell^2(J,\C) \oplus \ell^2(-J,\C)\cong \ell^2(J,\C) \oplus \ell^2(J,\C), $$
and with respect to this decomposition, we write elements of 
$\cL(\cH)$ as $2 \times 2$-block matrices. For  
$Q(v,w) = (w,v)$, we then have 
$$ \OO(\cH,I) = \{ g \in \GL(\cH) \: g^{-1} = Q g^\top Q^{-1} \} $$
and, for 
$g = \pmat{ a & b \cr c & d \cr}$, 
this means that 
$$ c b^\top + da^\top = \1, \quad 
cd^\top + dc^\top = 0 \quad \hbox{ and } \quad 
ab^\top + ba^\top = 0. $$
\item[\rm(e)] If $I^2 = -\1$, then there exists  
an orthonormal basis $(e_j)_{j \in 2J}$ of $\cH$ with 
\[Ie_j = 
\begin{cases}
e_{-j}, & \text{for } j \in J, \cr
-e_{-j}, & \text{for } j \in -J.
\end{cases}\] 
Then 
\[  \Sp(\cH,I) 
= \{ g \in \GL(\cH) \: g^{-1} = Q g^\top Q^{-1} \} \quad \hbox{ with } \quad 
Q = \pmat{ \0 & -\1 \cr \1 & \0}, \]
and for 
$g = \pmat{ a & b \cr c & d \cr}$
this means that 
\[ c^\top a = a^\top c, \quad d^\top b = b^\top d \quad \hbox{ and }
\quad a^\top  d- c^\top b = \1. \] 
\end{enumerate}
\end{exer}

\begin{exer} \mlabel{ex:inv-prod} Let 
$\cA$ be a unital associative algebra and $\cA^\times$ is unit group. 
If $a,b \in \cA$ commute and $ab\in \cA^\times$, then $a,b \in \cA^\times$. 
\end{exer}

\begin{exer} \mlabel{ex:9.1.x} 
Let $T(\cA) := \cA \oplus \cA$ denote the tangent algebra of 
$\cA$. Then $D \in \cL(\cA)$ is a derivation if and only if 
\[ \phi_D \: \cA \to T(\cA), \quad a \mapsto (a,Da) \] 
is an algebra homomorphism. 
\end{exer}

\end{small}

\section{Notes and comments on Chapter~\ref{ch:lingrp}}


Continuous inverse algebras can be regarded as generalizations
of Banach algebras which retain many of the essential properties of
those.
See \cite{Wa54a, Wa54b, Wa54c}
for early works concerning continuous inverse algebras
and functional calculus therein,
including generalizations when the algebra
multiplication is only separately continuous.
Compare also \cite{Ze60, Ze85}
for related studies.
Turpin's Theorem was taken from \cite{Tu70};
we give a new proof
which does not presuppose completeness properties or existence
of integrals. Local $m$-convexity of
topological algebras was introduced in \cite{MicE52}.
Our presentation of functional calculus in cias is based on~\cite{Gl02b},
where the Lie-theoretic properties of their unit groups
were elaborated. Unit groups of Banach-Lie algebras
are classical examples of Banach-Lie groups
(see, e.g., \cite{Bou89} and older works),
and their Lie subgroups
(as in \cite{dlH72}).
Further examples of cias
are suitable algebras of weighted mappings
with values in a Banach algebra
(see~\cite{Wa12}).
Compare also~\cite{Gra84} for examples;
implicitly, this source
also makes use of isospectral
homomorphisms (for inclusions
of subalgebras).
Proposition~\ref{prop:gramsch} 
generalizes \cite[Satz~5.11]{Gra84} from Banach algebras to cias. 
The above results concerning regularity of unit groups
of cias are taken from~\cite{GN12}.
Unit groups of some
projective limits of Banach algebras were already considered as regular
Fr\'{e}chet--Lie groups in~\cite{MR95}.
For multi-variable functional
calculus in cias, see~\cite{Bil07};
corresponding $*$-algebras were discussed in \cite{Bil04, BelN10}.
The Algebraic Subgroup Theorem~\ref{thm:5.5.14}  is due to 
Harris and Kaup \cite{HK77}. 

For a generalization of the Noether--Index as a means to classify 
connected components of certain groups of operators on Banach spaces 
to other types of pseudo-unitary groups, we refer to \cite{SV17}. 
These results are closely related to the Maslov index 
in the infinite-dimensional context; 
see \cite{Fu04, NO06}. 


\chapter{Mapping Groups} \mlabel{ch:mapgrp} 

%

\red{This chapter is not yet finished. Some proofs remain
  to be added and some sections have to be written. 
Below we marked in red where pieces are missing.} \\

Lie groups of mappings
and related Lie groups like gauge groups
of principal bundles are among the most important
examples of infinite-dimensional Lie groups.
In this section,
we discuss the Lie group $C^\ell(M,K)$
of $C^k$-maps from a compact manifold~$M$
(which may have a rough boundary)
to a given Lie group~$K$
modeled on a locally convex space.
For a $\sigma$-compact finite-dimensional
$C^\ell$-manifold $M$, we also discuss the
Lie group $C^\ell_c(M,K)$
of all $C^\ell$-maps $\gamma\colon M\to K$
which are compactly supported in the sense
that $\gamma(x)=\be$
for all $x$ off some compact subset of~$M$.
As a tool for the study of the latter,
we discuss fine box products of sequences of Lie groups.
For each countable locally finite family $(M_i)_{i\in I}$
of compact, full submanifolds $M_i\sub M$,
we shall see that $C^\ell_c(M,K)$ is isomorphic to a closed
Lie subgroup of the box product
\[
\bigoplus_{i\in I}C^\ell(M_i,K),
\]
which allows properties of $C^\ell_c(M,K)$
to be deduced from those of the $C^\ell(M_i,K)$.
\section{Mapping groups on compact manifolds}
In this section we start with the most regular class of mapping
  groups, namely those of maps $M \to K$, where $M$ is a compact smooth manifold
and $K$ is a Lie group.
\begin{prop}\label{gpCkmps}\mlabel{prop:cartes-closed-lie-a}
Let $K$ be a Lie group with Lie algebra~$\fk:=\L(K)$.
Let $\ell\in \N_0\cup\{\infty\}$
and $M$ be a compact $C^\ell$-manifold which may have a rough boundary.
Then the following holds:
\begin{description}[(D)]
\item[\rm(a)]
$C^\ell(M,K)$, endowed with
the canonical smooth manifold structure
as in Theorem~\emph{\ref{thmmfdmps}},
makes $C^\ell(M,K)$ a Lie group
modeled on $C^\ell(M,\fk)$.
\item[\rm(b)]
If $k\in \N_0\cup\{\infty\}$ and $L$ is a $C^k$-manifold which may have a rough boundary,
then a map $f\colon L\to C^\ell(M,K)$ is $C^k$ if and only if
\[ f^\wedge\colon L\times M\to K, \quad (x,y)\mto f(x)(y) \]
is~$C^{k,\ell}.$ 
\end{description}
In particular, if~$L$ is a smooth manifold, then a map
$f\colon L\to C^\infty(M,K)$ is smooth if and only if $f^\wedge\colon L\times M\to K$
is smooth.
\end{prop}
\begin{prf}
Recall from Proposition~\ref{exliegploa}
that $K$ admits a local addition,
whence Theorem~\ref{thmmfdmps}
applies.
Let $\mu\colon K\times K\to K$ be the smooth group multiplication
and $\eta\colon K\to K$ be the smooth inversion map
of the Lie group~$K$. The inversion map of $C^\ell(M,K)$
is the map $C^\ell(M,\eta)\colon C^\ell(M,K)\to C^\ell(M,K)$ and hence smooth,
by Lemma~\ref{ClmapCk}.
Likewise, the multiplication map is smooth as it can be identified
with $C^\ell(M,\mu)\colon C^\ell(M,K\times K)\to C^\ell(M,K)$
if we identify $C^\ell(M,K)\times C^\ell(M,K)$ with $C^\ell(M,K\times K)$
as in Lemma~\ref{prodasmfd}. The assertions concerning $f$ and $f^\wedge$
are subsumed by the Exponential Law in Theorem~\ref{thmmfdmps}(a).
\end{prf}
\begin{rem}\label{mapgp-intrinsic}
(a)
In Proposition~\ref{gpCkmps}, let
$\psi\colon U_\psi\to V_\psi$ be a $C^\infty$-diffeomorphism
from an open $e$-neighborhood $U_\psi\sub K$
onto an open subset $V_\psi\sub\fk$ (a $\fk$-chart of~$K$ around~$e$).
Then $C^\ell(M,U_\psi)$ is open in $C^\ell(M,K)$ (by Lemma~\ref{inCkopensub})
and the map $C^\ell(M,\psi)\colon C^\ell(M,U_\psi)\to
C^\ell(M,V_\psi)$ is a $C^\infty$-diffeomorphism (and hence a chart for $C^\ell(M,K)$
around its neutral element), as a consequence of Lemma~\ref{ClmapCk}.\medskip

\noindent
(b) If $C^\ell(M,K)$ is endowed with a Lie group
structure modeled on $C^\ell(M,\fk)$
which makes $C^\ell(M,\psi)$ a chart for $C^\ell(M,K)$
for some $\psi$
as in (a), then the Lie group structure coincides with the one constructed above,
as it induces the same smooth manifold structure on the
open $\be$-neighborhood $C^\ell(M,U_\psi)$.
\end{rem}
In the situation of
Proposition~\ref{prop:cartes-closed-lie-a},
the following holds:
\begin{prop}\label{mapgp-the-exp}
\begin{description}[(D)]
\item[\rm(a)] The map
\[
\Theta\colon \L(C^\ell(M,K))\to C^\ell(M,\fk),\quad
v\mto ((T_{\be}\ve_x)(v))_{x\in M}
\]
is an isomorphism of topological Lie algebras.
\item[\rm(b)]
If $K$ has an exponential function,
then also $C^\ell(M,K)$ has an exponential function.
Identifying $\L(C^\ell(M,K))$
with $C^k(M,\fk)$ as in {\rm(a)},
the latter is the mapping
\[
C^\ell(M,\exp_K)\colon C^\ell(M,\fk)\to C^\ell(M,K),\quad \gamma\mto \exp_K\circ \, \gamma.
\]
\item[\rm(c)]
If $K$ is locally exponential,
then $C^\ell(M,K)$ is locally exponential.
\item[\rm(d)]
If $k\in \N_0\cup\{\infty\}$
and $K$ is $C^k$-regular,
then also the Lie group $C^\ell(M,K)$
is $C^k$-regular.
\end{description}
\end{prop}
\begin{prf}
(a)
Let $f\colon M\to K$ be the constant function $x\mto\be$,
which is the neutral element of $C^\ell(M,K)$.
Then $\Gamma_f$, defined as in~\ref{special-gamma}, is given by
\[
\Gamma_f=C^\ell(M,\fk)
\]
and hence a Lie algebra with the pointwise Lie bracket.
By Theorem~\ref{thmmfdmps}(d),
$\Theta$ is an isomorphism of topological vector spaces.
Since $\ve_x\colon C^\ell(M,K)\to K$
is a morphism of Lie groups,
we see that
\[
T_{\be}(\ve_x)=\L(\ve_x)\colon L(C^\ell(M,K))\to \fk
\]
is a Lie algebra homomorphism for each $x\in M$,
entailing that $\Theta$ is a Lie algebra homomorphism.\smallskip

(b) By Lemma~\ref{ClmapCk}, the map $h:=C^\ell(M,\exp_K)$
is smooth. Given a map $\gamma \in C^\ell(M,\fk)$,
consider the smooth curve $\eta\colon [0,1]\to C^\ell(M,K)$,
$t\mto h(t\gamma)$. We have to show that $\Theta\circ \delta(\eta)$
is the constant curve $t\mto \gamma$.
For each $x\in M$,
we have that $\ve_x\circ\eta=(t\mto \exp_K(t \gamma(x)))$,
whence
\[
\Theta(\delta(\eta)(t))(x)
=\L(\ve_x)(\delta(\eta)(t))=\delta(\ve_x\circ \eta)(t)=\delta(t\mto\exp_K(t\gamma(x)))
\]
is the constant path $t\mto \gamma(x)$. Thus $\Theta\circ \delta(\eta)$
is the constant path $t\mto\gamma$.\smallskip

(c) If $K$ is locally exponential, then there exists an open $0$-neighborhood
$V\sub \fk$ such that $U:=\exp_K(V)$ is open in~$K$
and $\exp_K|_V\colon V\to U$ is a $C^\infty$-diffeomorphism.
Then $\psi:=(\exp_K|_V)^{-1}\colon U\to V$ is a $C^\infty$-diffeomorphism,
whence
\[
C^\ell(M,\psi)\colon C^\ell(M,U)\to C^\ell(M,V),\quad \gamma\mto \psi\circ\gamma
\]
is a $C^\infty$-diffeomorphism, see Remark~\ref{mapgp-intrinsic}(a).
It only remains to observe that the latter map is the inverse of
the exponential map $C^\ell(M,\exp_K)$, restricted to
a mapping from $C^\ell(M,V)$ to $C^\ell(M,U)$.\smallskip

(d) Let $I:=[0,1]$. By Proposition~\ref{prop:cartes-closed}
and Corollary~\ref{CkellCellk}, the map
\begin{equation}\label{change-sides}
\Phi\colon C^k(I,C^\ell(M,\fk))\to C^\ell(M,C^k(I,\fk))
\end{equation}
determined by $\Phi(\gamma)(x)(t)=\gamma(t)(x)$
is an isomorphism of topological vector spaces.
By Theorem~\ref{thmmfdmps}(a) and Corollary~\ref{CkellCellk},
the analogous map
\[
\Psi\colon C^{k+1}(I,C^\ell(M,K))\to C^\ell(M,C^{k+1}(I,K))
\]
is a bijection and hence an isomorphism of groups.
Let $\phi\colon U\to V$ be a $C^\infty$-diffeomorphism
from an open $\be$-neighborhood $U\sub K$
onto an open $0$-neighborhood $V\sub \fk$.
The map
\[
C^{k+1}(I,C^\ell(M,\fk))\to C^\ell(M,C^{k+1},\fk))
\]
analogous to (\ref{change-sides})
is an isomorphism of topological vector spaces
and
restricts to a $C^\infty$-diffeomorphism
\[
\Xi\colon C^{k+1}(I,C^\ell(M,V))\to C^\ell(M,C^{k+1}(I,V)
\]
between open $0$-neighborhoods as stated.
Note that the group isomorphism $\Psi$ restricts to the $C^\infty$-diffeomorphism
\[
C^\ell(M,C^{k+1}(I,\phi))^{-1}\circ \Xi\circ C^{k+1}(I,C^\ell(M,\phi))
\]
between open identity neighborhoods. Hence $\Psi$
is an isomorphism of Lie groups.
Thus $h:=\Psi^{-1}\circ \hspace*{.2mm}C^\ell(M,\Evol_K)\circ \hspace*{.2mm}\Phi$ defines a $C^\infty$-map
\[
h\colon C^k(I,C^\ell(M,\fk))\to C^{k+1}(I,C^\ell(M,K)).
\]
The proof will be complete if we can show that
$\delta(h(\gamma))=\Theta^{-1}\circ \gamma$
for each map $\gamma\in C^k(I,C^\ell(M,\fk))$;
the evolution map of $G:=C^\ell(M,K)$
will then be given by $\Evol_G=h\circ\Theta$,
which is a smooth map.
Let $\eta:=h(\gamma)\in C^{k+1}(I,C^\ell(M,K))$.
The family $(\ve_x)_{x\in M}$ of point evaluations
has the property that the $\L(\ve_x)$
separate points on $\L(C^\ell(M,K))$,
by~(a). For each $x\in M$, we have
\[
\ve_x\circ h(\gamma)=\Evol_K(t\mto \gamma(t)(x))
\]
and thus $\ev_x\circ \, \gamma=\delta(\ve_x\circ h(\gamma))=\L(\ve_x)\circ \delta(h(\gamma))=
\ev_x\circ \, \Theta\circ \delta(h(\gamma))$,
using the continuous linear evaluation maps $\ev_x\colon C^\ell(M,\fk)\to \fk$, $\theta\mto \theta(x)$.
Hence $\gamma=\Theta\circ \delta(h(\gamma))$,
which finishes the proof.
\end{prf}
\begin{prop}\label{mapgp-analyt}
Let $\ell\in \N_0\cup\{\infty\}$,
$M$ be a compact $C^\ell$-manifold $($which may have a rough boundary$)$,
$\K\in\{\R,\C\}$
and $K$ be a $\K$-analytic Lie group with Lie algebra~$\fk$.
Then the following holds:
\begin{description}[(D)]
\item[\rm(a)]
There is a $\K$-analytic manifold structure
on $C^\ell(M,K)$ modeled on $C^\ell(M,\fk)$
making it a $\K$-analytic Lie group,
which is uniquely determined by the property that
$C^\ell(M,U)$ is open in $C^\ell(M,K)$
and
\[
C^\ell(M,\phi)\colon C^\ell(M,U)\to C^\ell(M,V),\quad \gamma\mto \phi\circ\gamma
\]
is a diffeomorphism of $\K$-analytic manifolds
for each $\K$-analytic diffeomorphism $\phi\colon U\to V$
from an open symmetric $\be$-neighborhood $U\sub K$
onto an open $0$-neighborhood $V\sub \fk$.
\item[\rm(b)]
The smooth manifold structure underlying the $\K$-analytic Lie group
$C^\ell(M,K)$ from~{\rm(a)}
is the canonical structure described in
Proposition~\emph{\ref{gpCkmps}}.
\item[\rm(c)]
If $K$ has an exponential function
$\exp_K\colon \fk\to K$ which is $\K$-analytic,
then also the exponential function $C^\ell(M,\exp_K)$
of $C^\ell(M,K)$ is $\K$-analytic.
\item[\rm(d)]
If $K$ is a $\K$-analytic BCH-Lie group,
then also $C^\ell(M,K)$ is a $\K$-analytic BCH-Lie group.
\end{description}
\end{prop}
\begin{prf}
(a) and (b). Let $\phi\colon U\to V$ be a $\K$-analytic diffeomorphism from
an open, symmetric $\be$-neighborhood $U\sub K$ onto an open $0$-neighborhood $V\sub \fk$.
Then
\[
D_U:=\{(x,y)\in U\times U\colon xy\in U\}
\]
is open in $U\times U$ and contains $(U\times\{\be\})\cup(\{\be\}\times U)$.
Its image
\[
D_V:=(\phi\times \phi)(D_U)
\]
is open in $V\times V$.
We give $C^\ell(M,U)$ the $\K$-analytic manifold structure modeled on $C^\ell(M,\fk)$
making $\Phi:=C^\ell(M,\phi)$ a $\K$-analytic diffeomorphism;
its underlying smooth manifold structure makes it an open submanifold
of $C^\ell(M,K)$ with its canonical smooth manifold
structure (see Remark~\ref{mapgp-intrinsic}(a)).
We identify $C^\ell(M,U\times U)$ with $C^\ell(M,U)\times C^\ell(M,U)$
and $C^\ell(M,V\times V)$ with $C^\ell(M,V)\times C^\ell(M,V)$.
Then
$C^\ell(M,D_U)$
is open in $C^\ell(M,U\times U)$
and $\Phi\times \Phi$ restricts to a $\K$-analytic diffeomorphism
\[
\Psi\colon C^\ell(M,D_U)\to C^\ell(M,D_V).
\]
The group multiplication of~$K$ restricts to a $\K$-analytic map $m\colon D_U\to U$;
the group inversion restricts to a $\K$-analytic map $\eta_U\colon U\to U$.
Then also the mappings $\mu\colon D_V\to V$ and $j\colon V\to V$
given by
\[
\mu(x,y):=\phi(\phi^{-1}(x)\phi^{-1}(y))\quad\mbox{and}\quad
\iota(z):=\phi(\eta(\phi^{-1}zx)))
\]
for $(x,y)\in D_V$ and $z\in V$ are $\K$-analytic.
Thus
\[
C^\ell(M,\mu)\colon C^\ell(M,D_V)\to C^\ell(M,V),\quad (\gamma,\eta)\mto \mu\circ (\gamma,\eta)
\]
and $C^\ell(M,\iota)\colon C^\ell(M,V)\to C^\ell(M,V)$ are $\K$-analytic, by
Proposition~\ref{variant-pushforwards}(a).
As a consequence, also
\[
C^\ell(M,m)\colon C^\ell(M,D_U)\to C^\ell(M,U),\quad (\gamma,\eta)\mto m\circ (\gamma,\eta)=\gamma\eta
\]
and $C^\ell(M,\eta)\colon C^\ell(M,U)\to C^\ell(M,U)$, $\gamma\mto\eta_U\circ\gamma=\gamma^{-1}$
are $\K$-analytic.
Consider the continuous mapping
\[
c\colon K\times U\to K,\quad (g,y)\mto gyg^{-1}.
\]
If $\zeta\in C^\ell(M,K)$,
then $\zeta(M)$ is a compact subset of~$K$ and we also have
$c(\zeta(M)\times\{\be\})=\{\be\}\sub U$.
The Wallace Lemma provides an open $\be$-neighborhood $P\sub U$
and an open subset $N\sub K$ with $\zeta(M)\sub N$
such that $c(N\times P)\sub U$. Let $Q:=\phi(P)$.
Then
\[
I_\zeta\colon C^\ell(M,K)\to C^\ell(M,K),\quad \gamma\mto \zeta\gamma\zeta^{-1}
\]
maps $C^\ell(M,P)$ into $C^\ell(M,U)$.
Consider the mapping
\[
f\colon M\times Q\to V,\quad (x,y)\mto \phi(\zeta(x)\phi^{-1}(y)\zeta(x)^{-1}).
\]
Then
\[
h\colon N\times Q\to V,\quad (g,x)\mto \phi(g\phi^{-1}(x)g^{-1}) 
\]
is a $\K$-analytic mapping and $f(x,y)=h(\zeta(x),y)$ for all $(x,y)\in M\times Q$.
Thus
\[
f_*\colon C^\ell(M,Q)\to C^\ell(M,V),\quad \gamma\mto(\id_M,\gamma)
\]
is $\K$-analytic, by Proposition~\ref{variant-pushforwards}(b).
Then also
\[
I_\zeta|_{C^\ell(M,P)}=\Phi^{-1} \circ f_*\circ \Phi|_{C^\ell(M,P)}
\]
is $\K$-analytic.
Now the $\K$-analytic version of 
Theorem~\ref{thm:locglob}
%
%
provides a unique $\K$-analytic manifold structure
on $C^\ell(M,K)$ making it a Lie group~$G$
and turning $C^\ell(M,U)$ into an open $\K$-analytic
submanifold. By Remark~\ref{mapgp-intrinsic}(b), the smooth manifold
structure underlying the latter
is the canonical smooth manifold structure.
For the uniqueness statement, consider a $\K$-analytic diffeomorphism
$\phi'\colon U'\to V'$
analogous to $\phi$ and make $C \ell(M,K)$ a Lie group~$G'$
as before, using~$\phi'$ in place of~$\phi$.
Then $U\cap U'$ is open in both $U$ and $U'$.
Moreover, $W:=\phi(U\cap U')$
is open in $V$ and $W':=\phi'(U\cap U')$
is open in~$V'$. The set $C^\ell(M,U\cap U')$
is open in both $C^\ell(M,U)$ and $C^\ell(M,U')$.
Abbreviating $\Phi':=C^\ell(M,\phi')$, we have that
\[
\Phi'\circ \Phi^{-1}|_{C^\ell(M,W)}\colon C^\ell(M,W)\to C^\ell(M,W')
\]
is the map $C^\ell(M,(\phi'\circ \phi^{-1})|_W^{W'})$ and hence
a $\K$-analytic diffeomorphism. Hence $G$ and $G'$
induce the same $\K$-analytic manifold structure on the open set $C^\ell(M,U\cap U')$.
As a consequence, $\id\colon G\to G'$ is $\K$-analytic
and also $\id\colon G'\to G$, whence $G=G'$.

(c) Let $\phi\colon U\to V$ be as in (a). Then $\exp_K^{-1}(U)$
is an open $0$-neighborhood in~$\fk$
we choose an open, balanced $0$-neighborhood $W\sub \exp_K^ {-1}(U)$.
Then
\[
C^\ell(M,\phi)\circ C^\ell(M,\exp_K|_W)=C^\ell(M,\phi\circ \exp_K|_U)
\]
is $\K$-analytic, by Proposition~\ref{variant-pushforwards}(a). Hence $C^\ell(M,\exp_K)$ is $\K$-analytic
on the open set $C^\ell(M,W)$.
Now the open sets $C^\ell(M,nW)$ cover $C^\ell(M,\fk)$
for $n\in \N$. Since
\[
\exp_G(x)=\exp_({\textstyle\frac{1}{n}}x)^n,
\]
we have that $C^\ell(M,\exp_K)(\gamma)=(C^\ell(M,\exp_K|_W)({\textstyle\frac{1}{n}}\gamma))^n$
is a $\K$-analytic function of $\gamma\in C^\ell(M, nW)$. Thus $C^\ell(M,\exp_K)\colon C^\ell(M,\fk)\to C^\ell(M,K)$
is a $\K$-analytic map. By Proposition~\ref{mapgp-the-exp}(b), the latter map can be identified with the exponential
map of $C^\ell(M,K)$, up to composition with
an isomorphism of locally convex spaces (which is a $\K$-analytic diffeomorphism).

(d) By hypothesis, there exists an open $0$-neighborhood $V\sub \fk$
such that $U:=\exp_K(V)$ is open in $G$ and $\exp_K|_V\colon V\to U$
is a $\K$-analytic diffeomorphism.
After replacing $V$ with $(\exp_K|_V)^{-1}(U\cap U^{-1})$,
we may assume that $U$ is symmetric. Then (a) can be applied with
\[
\phi:=(\exp_K|_V)^{-1}\colon U\to V,
\]
showing that the map $\Phi:=C^\ell(M,\phi)\colon C^\ell(M,U)\to C^\ell(M,V)$
is a $\K$-analytic diffeomorphism
from an open identity neighborhood in $C^\ell(M,K)$
onto an open $0$-neighborhood in $C^\ell(M,\fk)$.
Hence
\[
C^\ell(M,\exp_K|_V)=\Phi^{-1}
\]
is a $\K$-analytic diffeomorphism. Recall that the exponential function $\exp_G$
of $G:=C^\ell(M,K)$ is given by $C^\ell(M,\exp_K)\circ\Theta$,
where $\Theta\colon \L(C^\ell(M,K))\to C^\ell(M,\fk)$
is as in Proposition~\ref{mapgp-the-exp}(a). Then $W:=\Theta^{-1}(C^\ell(M,V))$
is an open $0$-neighborhood in $\L(C^\ell(M,K))$
and $\exp_G|_W=C^\ell(M,\exp_K|_V)\circ \Theta|_W$
is a $\K$-analytic diffeomorphism onto the open
identity neighborhood $C^\ell(M,U)$.
\end{prf}
\section{Box products of Lie groups}
%

One way to deal with mapping groups on non-compact manifolds $M$, is
  to embed them into box products of mapping groups on compact
  submanifolds $M_i$ with boundary, exhausting~$M$. We therefore
  start this section by introducing box products of Lie groups. 
If $(G_i)_{i\in I}$
is a family of groups and $G:=\prod_{i\in I}G_i$
their direct product, then
\[
\bigoplus_{i\in I}G_i:=\{(g_i)_{i\in I}\in G\colon\,
\mbox{$g_i\not=e$ for only finitely many $i\in I$}\}
\]
is a normal subgroup of~$G$.
In the following, we construct Lie group
structures on such box products $\bigoplus_{i\in I}G_i$.
\begin{thm}\label{box-prod-lie}
Let $I$ be a countable set. For $i\in I$,
let $G_i$ be
a Lie group
modeled on a locally convex space~$E_i$.
Then there is a unique Lie group
structure on $G:=\bigoplus_{i\in I}G_i$
which is modeled on the locally convex direct sum $\bigoplus_{i\in I}E_i$
and satisfies the following condition:
\begin{description}[D]
\item[\rm$(*)$]
For every family $(\phi_i)_{i\in I}$
of charts $\phi_i\colon U_i\to V_i\sub E_i$
of $G_i$ such that $e\in U_i$ and $\phi_i(e)=0$,
the set $\bigoplus_{i\in I}U_i:= G\cap \prod_{i\in I}U_i$
is an open identity neighborhood in~$G$
and the following map is a chart for~$G$:
\[
\oplus_{i\in I}\phi_i\colon \bigoplus_{i \in I}U_i\to \bigoplus_{i\in I}V_i,\;\,
(x_i)_{i\in I}\mto (\phi_i(x_i))_{i\in I}.
\]
\end{description}
An analogous conclusion holds
when Lie groups are replaced with $\K$-analytic
Lie groups.
\end{thm}
We also show:
\begin{prop}\label{box-prod-regular}
If, for some $k\in \N_0$,
each of the Lie groups $G_i$ in Theorem~{\rm\ref{box-prod-lie}}
is $C^k$-regular, then also $\bigoplus_{i\in I}G_i$
is $C^k$-regular.
\end{prop}
We shall deduce Theorem~\ref{box-prod-lie}
from the following lemma.
\begin{lem}\label{fine-box-lie}
Let
$I$ be a countable set and $G_i$
be a Lie group $($resp., a $\K$-analytic Lie group$)$
for each $i\in I$.
Then the pointwise group operations
turn $G:=\prod_{i\in I}^{\fbx}G_i$
into a Lie group $($resp., a $\K$-analytic Lie group$)$.
\end{lem}
\begin{prf}
For each $i\in I$,
the group multiplication $\mu_i\colon G_i\times G_i\to G_i$
and group inversion $\eta_i\colon G_i\to G_i$, $g\mto g^{-1}$
are smooth (resp., $\K$-analytic) maps.
By Lemma~\ref{maps-between-box}, also the map
$\prod_{i\in I}\eta_i\colon G\to G$ is smooth (resp.,
$\K$-analytic), which is the group inversion in~$G$.
Identifying $G\times G$ with $\prod_{\in I}^{\fbx}(G_i\times G_i)$
as in Lemma~\ref{finebox-product}, the group multiplication of~$G$
turns into the map $\oplus_{i\in I}\mu_i$,
which is smooth (resp., $\K$-analytic) by Lemma~\ref{maps-between-box}.
\end{prf}
\noindent
\emph{Proof of Theorem}~\ref{box-prod-lie}.
For $\phi:=(\phi_i)_{i\in I}$
as in ($*$), $\bigoplus_{i\in I}U_i$
coincides with the set $U_\phi$ in~\ref{the-fine-top}
(by (\ref{U-phi-concrete})), whence $\bigoplus_{i\in I}U_i$
is open in the smooth (resp., $\K$-analytic$)$
Lie group $P:=\prod_{i\in I}^{\fbx}G_i$.
Being an identity neighborhood, the subgroup
$\bigoplus_{i\in I}G_i$ is open in~$P$
and hence a smooth (resp., $\K$-analytic) Lie group
using the induced manifold structure.
The mapping described in ($*$)
is the chart $\Phi_\phi$ of $P$ described in~\ref{basics-fine-box};
it also is a chart for $\bigoplus_{i\in I}G_i$.\qed\\

The following five lemmas prepare the proof
of Proposition~\ref{box-prod-regular}.
\begin{lem}\label{shape-lie}
In the situation of Theorem~{\rm\ref{box-prod-lie}}, we have:
\begin{description}[(D)]
\item[\rm(a)]
For each $i\in I$, the map
$\lambda_i\colon G_i\to G$
determined by $\lambda_i(g)_i=g$, $\lambda_i(g)_j=e$ for $j\in I\setminus\{i\}$
is a smooth  $($resp., $\K$-analytic$)$
group homomorphism.
\item[\rm(b)]
The component-wise Lie bracket makes $\bigoplus_{i\in I}\L(G_i)$
a topological Lie algebra and the map
$\Xi\colon \bigoplus_{i\in I}\L(G_i)\to \L(G)$, $(v_i)_{i\in I}\mto \sum_{i\in I}\L(\lambda_i)(v_i)$
is an isomorphism of topological Lie algebras.
\end{description}
\end{lem}
\begin{prf}
For $i\in I$, let $\phi_i\colon U_i \to V_i\sub E_i$ be a chart for~$G_i$
suhthat $e\in U_i$ and $\phi(e)=0$.
Let $\phi:=(\phi_i)_{i\in I}$
and $\Phi_\phi:=\oplus_{i\in }\phi_i
\colon \bigoplus_{i \in I}U_i\to\bigoplus_{i\in I}V_i\sub \bigoplus_{i\in I}E_i=:E$
be the associated chart of~$G$.
For $i\in I$, let $\alpha_i\colon E_i\to E$
be the continuous linear map with $\alpha_i(v)_i=v$,
$\alpha_i(v)_j=0$ for $j\in I\setminus \{i\}$.
Then
\[
\Phi_\phi\circ \lambda_i\circ \phi_i^{-1}=\alpha_i|_{V_i},
\]
whence $d\Phi_\phi|_{T_eG}\circ T_e\lambda_i\circ T\phi_i^{-1}(0,\cdot)=\alpha_i$.
The map $d\Phi_\phi|_{T_eG}$ being linear,
\[
d\Phi_\phi|_{T_eG}\circ \Xi\circ \oplus_{i\in I}T\phi_i^{-1}(0,\cdot)=\id_E
\]
follows. Since $d\Phi_\phi|_{T_eG}$ and $\oplus_{i\in I}T\phi_i^{-1}(0,\cdot)$
are isomorphsms of topological vector spaces, also $\Xi$
is an isomorphism of topological vector spaces.
For every non-empty finite subset $J$ of~$I$,
consider the map
\[
\lambda_J\colon \prod_{j\in J}G_j\to G
\]
determined by $\lambda_J((x_j)_{j\in J})_i=x_i$ if $i\in J$,
$\lambda_J((x_j)_{j\in J})_i=e$ if $i\in I\setminus J$.
Then $\lambda_J$
is a smooth group homomorphism.
Thus $\L(\lambda_J)$ is a Lie algebra homomorphism.
We now identify $\L(\prod_{j\in J}G_j)$ with $\prod_{j\in J}\L(G_j)$
as in Lemma~\ref{tanprod}(b); moreover,
we identify $\prod_{j\in J}\L(G_i)$ with its image in $\bigoplus_{i\in I}\L(G_i)$
under the map which is analogous to $\lambda_J$.
Then $\Xi|_{\prod_{j\in J}\L(G_j)}=\L(\lambda_J)$
is a Lie algebra homomorphism.
As $(\prod_{j\in J}\L(G_j))_{J\sub I\,\text{finite}}$
is an upward directed family
of Lie subalgebras of
$\bigoplus_{i\in I}\L(G_i)$ with union $\bigoplus_{i\in I}\L(G_i)$,
we deduce that~$\Xi$ is a Lie algebra homomorphism.
\end{prf}
\begin{rem}\label{idenfy-sum}
By Lemma~\ref{shape-lie}(b), we may consider $\L(G)$ as the locally convex direct sum $\bigoplus_{i\in I}\L(G_i)$,
together with the continuous linear mappings $\L(\lambda_i)\colon \L(G_i)\to \L(G)$.
\end{rem}
\begin{lem}\label{maps-to-sum}
Let $(E_i)_{i\in I}$
be a countable family of locally convex spaces
with locally convex direct sum $E:=\bigoplus_{i\in I}E_i$.
For $i\in I$, let
$\alpha_i\colon E_i\to E$
be the continuous linear map
determined by $\alpha_i(v)_i=v$, $\alpha_i(v)_j=0$ for $j\in I\setminus\{i\}$.
\begin{description}[(D)]
\item[\rm(a)]
If $K$ is a compact topological space,
give $C(K,E_i)$ the compact-open topology
for $i\in I$ and form the locally convex direct sum
$\bigoplus_{i\in I}C(K,E_i)$.
Then the following map is an isomorphism
of topological vector spaces:
\[
\alpha\colon \bigoplus_{i\in I}C(K,E_i)\to C(K,E),\;\,
(\gamma_i)_{i\in I}\mto \sum_{i\in I} (\alpha_i\circ \gamma_i).
\]
\item[\rm(b)]
For $k\in \N_0$, give $C^k([0,1],E)$
the compact-open $C^k$-topology
and form the locally convex direct sum $\bigoplus_{i\in I}C^k([0,1],E_i)$.
Then the mapping\linebreak
$\beta_k\colon \bigoplus_{i\in I}C^k([0,1],E_i)\to C^k([0,1],E)$,
$(\gamma_i)_{i\in I}\mto \sum_{i\in I} (\alpha_i\circ \gamma_i)$
is an isomorphism of topological vector spaces.
\end{description}
\end{lem}
\begin{prf}
(a) We first note that $\alpha$ is linear and injective.
For each $i\in I$, the linear map $C(K,\alpha_i)\colon C(K,E_i)\to C(K,E)$
is continuous by Lemma~\ref{covsuppo}. Hence $\alpha$ is continuous, by Lemma~\ref{firstlasum}(a).
If $\gamma\in C(K,E)$, then
\[ \gamma(K)\sub \prod_{i\in J}E_i=:F \]
for some finite subset $J\sub I$, by Lemma~\ref{firstlasum}(e).
For each $j\in J$, we set $\gamma_j:=\pr_j\circ \, \gamma$,
using the continuous linear projection $\pr_j\colon F\to E_j$
onto the $j$th component. For $i\in I\setminus J$, set $\gamma_i:=0\in C(K,E_i)$.
Then $\alpha((\gamma_i)_{i\in I})=\gamma$,
whence $\alpha$ is surjective and hence an isomorphism of vector spaces.
Every $0$-neighborhood in $C(K,E)$
contains a $0$-neighborhood of the form $C(K,U)$
for some $0$-neighborhood $U\sub E$ (see Lemma~\ref{sammelsu}(b)).
After shrinking $U$, we may assume that $U=\bigoplus_{\in I}U_i$
with open $0$-neighborhoods $U_i\sub E_i$, see Remark~\ref{firstremsums}(a).
Then $W:=\bigoplus_{i\in I}C(K,U_i)$ is an open $0$-neighborhood
in $\bigoplus_{i\in I}C(K,E_i)$ and $\alpha(W)\sub C(K,U)$,
whence $\alpha^{-1}$ is continuous at~$0$ and hence continuous.

(b) $\beta_k$ is linear and continuous.
As in (a), each $\gamma\in C^k([0,1],E)$ has image in
a finite partial product; since $\gamma_i:=\pr_i\circ \,\gamma$ is $C^k$,
we see that $\beta_k$ is surjective.
Now
\[
\lambda_i\colon C^k([0,1],E_i)\to\prod_{j=0}^kC([0,1],E_i),\;\,
\gamma\mto(\gamma^{(j)})_{j=0}^k
\]
is a topological embedding, for each $i\in I$
(cf.\ Exercise~\ref{exc-the-topf}).
Hence
\[
\oplus_{i\in I}\lambda_i\colon
\bigoplus_{i\in I} C^k([0,1],E_i)\to\bigoplus_{i\in I} \prod_{j=0}^kC([0,1],E_i)
\cong \prod_{j=0}^k\bigoplus_{i \in I} C([0,1],E_i)
\]
is a topological embedding (cf.\ (d) and (c) in Lemma~\ref{firstlasum}). The topology
$\cO$ on
$\bigoplus_{i\in I}C^k([0,1],E_i)$
is therefore initial with respect to the maps
\[
\delta_j\colon
\bigoplus_{i\in I}C^k([0,1],E_i)\to \bigoplus_{i\in I}C([0,1],E_i),\;\,
(\gamma_i)_{i\in I}\mto (\gamma_i^{(j)})_{i\in I}
\]
for $j\in \{0,1,\ldots, k\}$, and hence
also with respect to the maps $\alpha\circ \delta_j$
with $\alpha\colon \bigoplus_{i\in I}C([0,1],E_i)\to C([0,1],E)$
as in~(a).
The topology on $C^k([0,1],E)$ is initial with respect to the maps
$D_j\colon C^k([0,1],E)\to C([0,1],E)$, $\gamma\mto \gamma^{(j)}$
for $j\in\{0,1,\ldots, k\}$ (cf.\ Exercise~\ref{exc-the-topf}).
By transitivity of initial topologies (Lemma~\ref{transinit}),
the topology $\cT$ on $\bigoplus_{i\in I}C^k([0,1],E_i)$
making $\beta_k$ a topological
embedding is therefore initial
with respect to the maps $D_j\circ \beta_k$
for $j\in \{0,1\ldots, k\}$.
Since $D_j\circ \beta_k=
\alpha\circ\delta_j$, we see that $\cO=\cT$.
\end{prf}
\begin{rem}
Taking $k=\infty$ in Lemma~\ref{maps-to-sum}(b),
the proof still shows that $\beta_\infty$
is linear, bijective, and continuous.
However, $\beta_\infty$ need not be a homeomorphism.
For example, if $I:=\N$ and $E_n:=\R$ for $n\in \N$,
then
\[
U_n:=\{\gamma \in C^\infty([0,1],\R)\colon (\forall j\in\{0,1,\ldots,n\})\;\,
\|\gamma^{(j)}\|_\infty<1\}
\]
is an open $0$-neighborhood in $C^\infty([0,1],\R)$
but $\beta_\infty(\bigoplus_{n\in\N}U_n)$
is not a $0$-neighborhood
in $C^\infty([0,1],E)$,
the topology of which is initial with respect to the inclusion maps
$C^\infty([0,1],E)\to C^k([0,1],E)$ for $k\in \N_0$.
\end{rem}
\begin{lem}\label{cp-in-box-group}
In the situation of Theorem~{\rm\ref{box-prod-lie}},
every compact subset $K\sub G$ is contained in
$\prod_{j\in J}G_j$ for some finite subset
$J\sub I$.
\end{lem}
\begin{prf}
For $i\in I$, let $\phi_i\colon U_i \to V_i\sub E_i$ be a chart for~$G_i$
such that $e\in U_i$ and $\phi_i(e)=0$.
Let $\phi:=(\phi_i)_{i\in I}$
and
\[ \Phi_\phi:=\oplus_{i\in }\phi_i
\colon U_\phi\to\bigoplus_{i\in I}V_i\sub \bigoplus_{i\in I}E_i=:E \] 
be the corresponding chart of~$G$,
with $U_\phi:=\bigoplus_{i\in I}U_i$.
As~$G$, being a topological group,
is a regular topological space,
$U_\phi$ contains am identity neighborhood~$A\sub G$
which is closed in~$G$.
Then $\Phi_\phi(A\cap g^{-1}K)$ is a compact subset of $\bigoplus_{i \in I}E_i$
for each $g\in G$ and thus
$\Phi_\phi(A\cap g^{-1}K)\sub \prod_{j\in J_g}E_j$
for some finite subset $J_g\sub I$, whence
\begin{equation}\label{fi-fin}
K\cap gA\sub g\prod_{j\in J_g}G_j.
\end{equation}
Since $K$ is compact, there exist $g_1,\ldots, g_n\in K$ such that
\begin{equation}\label{se-fin}
K\sub \bigcup_{k=1}^ng_kA^0.
\end{equation}
For each $k\in \{1,\ldots, n\}$,
we have $g_k\in \prod_{j\in I_k}G_j$
for some finite subset $I_k$ of~$J$.
Then $J:=\bigcup_{k=1}^n (I_k\cup J_{g_k})$
is a finite subset of~$I$;
by (\ref{fi-fin}) and (\ref{se-fin}), we have $K\sub \prod_{j\in J}G_j$.
\end{prf}
\begin{lem}\label{maps-to-box}
In the situation of Theorem~{\rm\ref{box-prod-lie}},
let $k\in \N_0$ and consider the map
\[
\Theta\colon \bigoplus_{i\in I}C^k([0,1],G_i)\to C^k([0,1],G)
\]
determined by $\Theta((\gamma_i)_{i\in I})(t):=(\gamma_i(t))_{i\in I}$.
Then $\Theta$ is an isomorphism of Lie
groups $($resp., of $\K$-analytic Lie groups$)$.
\end{lem}
\begin{prf}
In view of Lemma~\ref{shape-lie}(a), Lemma~\ref{ClmapCk} implies that $\Theta((\gamma_i)_{i\in I})$
is an element of $C^k([0,1],G)$
for each $(\gamma_i)_{i\in I}$.
Abbreviate $\g_i:=L(G_i)$ for $i\in I$ and identify $\g:=\L(G)$
with $\bigoplus_{i\in I}\g_i$
as in Remark~\ref{idenfy-sum}.

Consider the isomorphism $\beta_k\colon \bigoplus_{i\in I}C^k([0,1],\g_i)\to C^k([0,1],\g)$
of topological vector spaces as in Lemma~\ref{maps-to-sum}(b).
We easily see that~$\Theta$ is a homomorphism of groups and injective.
Lemma~\ref{cp-in-box-group} implies that $\Theta$ is surjective,
whence $\Theta$ is an isomorphism of groups.
For $i\in I$, let $\phi_i\colon U_i \to V_i\sub \g_i$ be a $\g_i$-chart for~$G_i$
such that $e\in U_i$ and $\phi_i(e)=0$.
Let $\phi:=(\phi_i)_{i\in I}$
and $\Phi_\phi:=\oplus_{i\in }\phi_i
\colon \bigoplus_{i \in I}U_i\to\bigoplus_{i\in I}V_i\sub \bigoplus_{i\in I}\g_i=\g$
be the corresponding $\g$-chart of~$G$.
By Remark~\ref{mapgp-intrinsic}(a), the map
\[
(\phi_i)_*\colon C^k([0,1],U_i)\to C^k([0,1],V_i)\sub C^k([0,1],E_i),\;\,
\gamma\mto \phi\circ \gamma
\]
is a chart for $C^k([0,1],G)$. We set $\psi:=((\phi_i)_*)_{i \in I}$
and consider the associated chart
$\Phi_\psi:=\oplus_{i\in I}(\phi_i)_*\colon
U_\psi\to V_\psi$
with $U_\psi:=\bigoplus_{i\in I}C^k([0,1],U_i)$
and $V_\psi:=\bigoplus_{i\in I}C^k([0,1],V_i)\sub \bigoplus_{i\in I}C^k([0,1],\g_i)$.
Likewise, we have the chart
\[
(\Phi_\phi)_*\colon C^k([0,1],U_\phi)\to C^k([0,1],V_\phi),\;\, \gamma\mto \Phi_\phi\circ\gamma
\]
of $C^k(M,G)$, with $U_\phi:=\bigoplus_{i\in I}U_i$
and $V_\phi:=\bigoplus_{i\in I}V_i$.
Now
\[
(\Phi_\phi)_*\circ \Theta|_{U_\psi}\circ (\Phi_\psi)^{-1}=\beta_k|_{V_\psi},
\]
where $\beta_k|_{V_\psi}\colon V_\psi\to C^k([0,1], V_\phi)$
is a $C^\infty$-diffeomorphism (resp., a $C^\omega_\K$-diffeomorphism).
As a consequence, the group isomorphism $\Theta$
is an isomorphism of Lie groups (resp., of $\K$-analytic Lie groups).
\end{prf}
\noindent
\emph{Proof of Proposition}~\ref{box-prod-regular}.
Let $\Xi$ be the isomorphism
of topological Lie algebras introduced in Lemma~\ref{shape-lie}.
For each $\gamma\in C^k([0,1],\L(G))$,
we have $\Xi^{-1}(\gamma)\in \prod_{j\in J}C^k([0,1],\L(G_j))$
for some non-empty finite subset $J\sub I$.
Write $\zeta:=\Xi^{-1}(\gamma)=(\zeta_j)_{j\in J}$
with $\zeta_j\in C^k([0,1],G_j)$.
Each $G_j$ being $C^k$-regular,
also $G_J:=\prod_{j\in J}G_j$
is $C^k$-regular and
\[
\Evol_{G_J}(\zeta)(t)=(\Evol_{G_j}(\zeta_j)(t))_{j\in J}\;\,\mbox{for all $t\in [0,1]$,}
\]
identifying $\L(G_J)$ with $\prod_{j\in J}\L(G_j)$
as in Remark~\ref{lie-produ} (see Lemma~\ref{prod-reg}).
For $\lambda_J\colon G_J\to G$ as in the proof of Lemma~\ref{shape-lie},
we have left logarithmic derivatives
\[
\delta(\lambda_J\circ \Evol_{G_j}(\zeta))=\L(\lambda_J)\circ\delta(\Evol_{G_J}(\zeta))
=\L(\lambda_J)\circ \zeta=\gamma,
\]
whence $\Evol_G(\gamma)=\lambda_J\circ \Evol_{G_J}(\zeta)$.
Abbreviate $\g:=\L(G)$, $\g_i:=\L(G_i)$
and ${\mathbb I}:=[0,1]$.
By the preceding, $\Evol_G$ can be written as the composition
\[
C^k({\mathbb I},\g )\stackrel{\Xi^{-1}}{\longrightarrow}
\bigoplus_{i\in I}C^k({\mathbb I},\g_i)
\lmapright{\oplus_{i\in I}\Evol_{G_i}}
\bigoplus_{i\in I}C^k({\mathbb I},G_i)\stackrel{\Theta}{\longrightarrow}
C^k({\mathbb I},G)
\]
using the map $\oplus_{i\in I}\Evol_{G_i}$
(which is $C^\infty$ by Lemma~\ref{maps-between-box})
and the isomorphism
$\Theta$ of Lie groups from Lemma~\ref{maps-to-box}.
Hence $\Evol_G$ is smooth.\qed
\section{Mapping groups on non-compact manifolds}
We now consider Lie groups
of compactly supported Lie group-valued mappings on a non-compact
manifold~$M$. The case of a finite-dimensional manifold~$M$
without boundary may be of highest interest,
but the construction works without changes
if $M$ has a boundary, or rough boundary.\medskip

Given a Lie group~$K$ and $\ell\in \N_0\cup\{\infty\}$,
we let $C^\ell_c(M,K)$
be the group of all $C^\ell$-functions $\gamma\colon M\to K$
for which there exists a compact subset $L\sub G$
such that $\gamma(x)=\be$ for all $x\in M\setminus L$.
In particular,
given a locally convex space $E$ over $\K\in \{\R,\C\}$

We now describe a Lie group structure on $C^\ell_c(M,K)$,
modeled on the locally convex space $C^\ell_c(M,\fk)$
(cf.~Lemma~\ref{details-gamma-nonc}). 
More generally, we are able to turn $C^\ell(M,K)$
into a Lie group having $C^\ell_c(M,K)$
as an open subgroup. However, our main interest is in $C^\ell_c(M,K)$
itself.
If $M$ is not compact,
the resulting topology on $C^\ell_c(M,K)$
can be finer than the compact-open $C^\ell$-topology
and may be of limited usefulness.
%
%
%
\begin{thm}\label{thm-tefu-gp}
Let $K$ be a Lie group with Lie algebra $\fk:=\L(K)$.
Let $\ell\in \N_0\cup\{\infty\}$
and $M$ be a compact $C^\ell$-manifold which may have a rough boundary.
Then the following holds:
\begin{description}[(D)]
\item[\rm(a)]
The smooth manifold structure
described in Theorem~\ref{thmmfdmps} 
makes $C^\ell(M,K)$ a Lie group
modeled on $C^\ell(M,\fk)$.
\end{description}
Using this manifold structure, we have:
\begin{description}[(D)]
\item[\rm(b)]
For each $C^\infty$-diffeomorphism $\phi\colon U\to V$
from an open, symmetric $\be$-neighborhood $U\sub K$
onto an open $0$-neighborhood $V\sub \fk$
such that $\phi(\be)=0$, the set
\[
C^\ell_c(M,U):=\{\gamma\in C^\ell_c(M,K)\colon \gamma(M)\sub U\}
\]
is open in $C^\ell(M,K)$ and the map $C^\ell_c(M,\phi)\colon C^\ell_c(M,U)\to
C^\ell_c(M,V)$ is a $C^\infty$-diffeomorphism onto the open
subset $C^\ell_c(M,V)\sub C^\ell_c(M,\fk)$.
\item[\rm(c)]
The normal subgroup $C^\ell_c(M,K)$ is open in $C^\ell(M,K)$.
\item[\rm(d)]
For each $x\in M$, the map $\ve_x\colon C^\ell(M,K)\to K$, $\gamma\mto\gamma(x)$
is a smooth group homomorphism. Moreover, $\Theta(v):=(\L(\ve_x)(v))_{x\in M}\in C^\ell_c(M,\fk)$
for each $v\in \L(C^\ell(M,K))$ and the map $\Theta\colon \L(C^\ell(M,K))\to C^\ell_c(M,\fk)$
is an isomorphism of topological Lie algebras.
\item[\rm(e)]
If $K$ has a smooth exponential function $\exp_K\colon \fk\to K$,
then
\[ C^\ell(M,\exp_K)\circ \Theta\colon \L(C^\ell(M,K))\to C^\ell(M,K),\quad
v\mto \exp_K\circ \Theta(v).\] 
\item[\rm(f)]
If $K$ is locally exponential, then also $C^\ell(M,K)$
is locally exponential.
\item[\rm(g)]
If $K$ is $C^k$-regular for some finite $k\in \N_0$,
then also $C^\ell(M,K)$ is $C^k$-regular.
\end{description}
\end{thm}
\begin{prf} \bred{Proof to be written}
\end{prf}
\begin{rem}
(i) We can replace $C^\ell(M,K)$ with its open subgroup $C^\ell_c(M,K)$
in parts (b), (d), (e), (f) and (g)
of Theorem~\ref{thm-tefu-gp}.
\medskip

(ii) Let $C^\ell(M,K)$ (or $C^\ell_c(M,K)$)
be endowed with a smooth manifold structure modeled
on $C^\ell(M,\fk)$ making it a Lie group $G'$.
Assume that, for some $C^\infty$-diffeomorphism $\phi\colon U\to V$
as  in Theorem~\ref{thm-tefu-gp}(b),
the set $C^\ell_c(M,U)$ is open in $G'$ and $C^\ell_c(M,\phi)\colon C^\ell_c(M,U)\to C^\ell_c(M,V)\sub C^\ell_c(M,\fk)$
is a $C^\infty$-diffeomorphism. Then $G'$ coincides with
the Lie group $G=C^\ell(M,K)$ (resp., $G=C^\ell_c(M,K)$)
described in Theorem~\ref{thm-tefu-gp}(a). In fact, the group homomorphisms
$\id\colon G\to G'$ and $\id\colon G'\to G$ are smooth on an identity neighborhood and hence
smooth; thus $G=G'$.
\end{rem}
We now consider Lie groups
of compactly supported Lie group-valued mappings on a non-compact
manifold~$M$. The case of a finite-dimensional manifold~$M$
without boundary may be of highest interest,
but the construction works without changes
if $M$ has a boundary, or rough boundary.\medskip

Let $r\in \{\infty,\omega\}$. If $r=\infty$, let $\K=\R$.
If $r=\omega$, let $\K\in \{\R,\C\}$.
The manifolds~$M$ we consider are manifolds
over the ground field of real numbers.
Given a Lie group~$K$ and $\ell\in \N_0\cup\{\infty\}$,
we let $C^\ell_c(M,K)$
be the group of all $C^\ell$-functions $\gamma\colon M\to K$
for which there exists a compact subset $L\sub G$
such that $\gamma(x)=\be$ for all $x\in M\setminus L$.
\begin{thm} 
Let $\ell\in \N_0\cup\{\infty\}$,
$M$ be a $\sigma$-compact, locally compact
$C^\ell$-manifold $($which may have a rough boundary$)$
and $K$ be a $C^r_\K$-Lie group modeled
on a locally convex topological $\K$-vector space,
with Lie algebra $\fk$.
Then the following holds:
\begin{description}[(D)]
\item[\rm(a)]
There is a $C^r_\K$-manifold structure
on $C^\ell_c(M,K)$ modeled on $C^\ell_c(M,\fk)$
making it a $C^r_\K$-Lie group,
which is uniquely determined by the property that
$C^\ell_c(M,U):=\{\gamma\in C^r_c(M,K)\colon \gamma(M)\sub U\}$
is open in $C^\ell_c(M,K)$
and
\[
C^\ell_c(M,\phi)\colon C^\ell_c(M,U)\to C^\ell_c(M,V),\quad \gamma\mto \phi\circ\gamma
\]
is a diffeomorphism of $C^r_\K$-manifolds
for each $C^r_\K$-diffeomorphism \break $\phi\colon U\to V$
from an open symmetric $\be$-neighborhood $U\sub K$
onto an open $0$-neighborhood $V\sub \fk$, such that $\phi(\be)=0$.
%
%
\item[\rm(b)]
If $K$ has an exponential function
$\exp_K\colon \fk\to K$ which is $\K$-analytic,
then also the exponential function $C^\ell(M,\exp_K)$
of $C^\ell(M,K)$ is $\K$-analytic.
\item[\rm(c)]
If $K$ is a $\K$-analytic BCH-Lie group,
then also $C^\ell(M,K)$ is a $\K$-analytic BCH-Lie group.
\end{description}
\end{thm}
%
%
%
%
%
%

%


%
\section{Gauge groups of principal bundles}
\bred{This section has to be written. Here is the main result
with a sketch of its proof.} 

\begin{thm}
  \mlabel{thm:IV.1.12}  If $K$ is a locally exponential Lie group with
  Lie algebra~$\fk$ and 
$q \: P \to M$ a smooth $K$-principal bundle over the $\sigma$-compact finite-dimensional 
manifold $M$, then the group $\Gau_c(P)$ of compactly supported gauge transformations 
is a locally exponential Lie group. 
In particular, the Lie group $C^\infty_c(M,K)$ is locally exponential. 

If, in addition, $K$ is regular brauche endliches k falls M nicht kompakt, then $\Gau_c(P)$ is regular and 
if $K$ is BCH, then so is $\Gau_c(P)$. 
\end{thm}

\begin{prf}  (Sketch) 
Let $\exp_K \: \fk \to K$ be the exponential function of $K$ 
and realize $\Gau(P)$ as the subgroup 
$C^\infty(P,K)^K$ of $K$-fixed points in $C^\infty(P,K)$ with respect to the 
$K$-action given by $(k.f)(p) := k f(p.k)k^{-1}$. 
Then we put 
\begin{align*}
\gau(P) 
&:= C^\infty(P,\fk)^K \\
&= \{ \xi \in C^\infty(P,\fk) \: 
(\forall p \in P)(\forall k \in K)\, \Ad(k)\xi(p.k) = \xi(p) \},
\end{align*}
and observe that, for the group $G := \Gau_c(P)$, the map 
$$ \exp_G \:\g := \gau_c(P) \to G, \quad \xi \mapsto \exp_K \circ \xi $$
is a local homeomorphism in $0$. Using Theorem~\ref{thm:locglob}, 
this can be used to define a Lie group structure on $G$. 
Then $\exp_G$ is an exponential function of~$G$, and, by construction, 
it is a local diffeomorphism in $0$. 
\end{prf}


%


\chapter{Direct Limit Lie Groups} \mlabel{ch:dirlim} 

\vspace{5mm}
Frequently, a group $G$ of interest can be written
as the union $G=\bigcup_{n\in\N}G_n$
of Lie groups $G_1\sub G_2\sub\cdots$,
such that all inclusion maps $G_n\to G_{n+1}$
are smooth group homomorphisms.
In this chapter, we first consider
the case that each $G_n$ is a finite-dimensonal Lie group.
We show that the smooth manifold structure
on $G$ constructed in Section~\ref{dl-construction}
makes $G$ a $C^0$-regular Lie group
whose Lie algebra is an ascending
union of finite-dimensional Lie algebras
isomorphic to $L(G_n)$.
Every ascending union of finite-dimensional Lie algebras
arises in this way (see Section~\ref{sec:VII.1}).
In Section~\ref{more-findim-DL},
we survey some additional results
concerning such locally finite-dimensional Lie groups
$G=\bigcup_{n\in \N}G_n$ without proofs, but with references to the literature.
In Section~\ref{sec:VII.2},
we then briefly consider the case that the $G_n$
may be infinite-dimensional.
In many cases, a Lie group structure on the union $G=\bigcup_{n\in \N}G_n$
is available by ad hoc arguments suitable for a given class of examples.
For example, we already encountered
test function groups $C^\infty_c(M,K)$
and Lie groups $\Diff_c(M)$ of compactly supported diffeomorphisms,
for a $\sigma$-compact finite-dimensional smooth manifold~$M$
and a Lie group~$K$.
For a compact exhaustion $(L_n)_{n\in\N}$ of~$M$,
the Lie group $C^\infty(M,K)$ is the union of the ascending sequence
$C^\infty_{K_n}(M,K)$ of smooth $K$-valued functions supported
in~$K_n$; likewise, $\Diff_c(M)=\bigcup_{n\in \N}\Diff_{K_n}(M)$.
We survey some results concerning
direct limit properties of such ascending unions.
\section{Direct limits of finite-dimensional Lie groups}\mlabel{sec:VII.1} 
We discuss Lie group structures on ascending unions of finite-dimensional Lie groups.
In this section, Lie groups are real Lie groups.
Notably, we show:
\begin{thm}\label{dl-lie}
Let $G_1\sub G_2\sub\cdots$
be a sequence of finite-dimensional Lie groups
such that the inclusion map $\lambda_{n,m}\colon G_m\to G_n$
is a smooth group homomorphism and an immersion
for all $n,m\in\N$ with $n\geq m$.
Give $G:=\bigcup_{n\in \N}G_n$ the unique group structure making $G_n$ a subgroup for each $n\in \N$.
Then the smooth manifold structure on $G:=\bigcup_{n\in \N}G_n$
provided by Theorem~{\rm\ref{dl-mfd}}
makes $G$ a Lie group.
The inclusion map $\lambda_n\colon G_n\to G$
is a smooth group homomorphism for each $n\in \N$.
If $H$ is a Lie group, then a group homomorphism $f\colon G\to H$
is smooth if and only if $f|_{G_n}$ is smooth for each $n\in\N$.
\end{thm}
\begin{ex}
  \mlabel{ex:VII.1.4} (a) One of the most famous examples of a direct limit Lie group 
is the group 
$$ \GL_\infty(\R) := \indlim \GL_n(\R) $$
with the connecting maps 
\[ \lambda_{n+1,n} \: \GL_n(\R) \to \GL_{n+1}(\R), \quad a \mapsto \pmat{ a & 0 \\ 0 & 1}.\]
Its Lie algebra is the Lie algebra 
$\gl_\infty(\R)$ of all $(\N \times \N)$-matrices with 
only finitely many non-zero entries.
Direct limit Lie groups $\SL_\infty(\R)$,
$O_\infty$ and $U_\infty$
can be formed
using the analogous connecting maps $\SL_n(\R)\to\SL_{n+1}(\R)$,
$O_n\to O_{n+1}$ and $U_n\to U_{n+1}$,
respectively.

(b) In the context of $C^*$-algebras, direct limits of finite-dimensional ones are particularly 
interesting objects. On the level of unit groups one encounters in particular 
groups of the form 
\[  G := \indlim \GL_{2^n}(\C), \quad \lambda_{n+1,n}(a) = \pmat{ a & 0 \\ 0 & a}.\] 

(c) Let $E := \C^{(\N)}$ be the free vector space with basis 
$(e_n)_{n \in \N}$ and $D \in {\cal L}(E)$ be defined by 
$D(e_n) = 2\pi i n e_n$. 
Then the Lie algebra $\g := E \rtimes_D \R$ is locally finite and 
we obtain a corresponding locally finite-dimensional Lie group 
$G = E \rtimes_\alpha \R$, where $\alpha(t) = e^{tD}$. 
Since the sequence $(0,{1\over n})_{n \in \N}$ consists of singular points 
of the exponential function, the Lie algebra $\g$ is not locally exponential 
(cf.\ Proposition~\ref{prop:exp-noninj}). 
\end{ex}
We first have a look at the relevant Lie algebras.
\begin{lem}\label{do-not-dodge}
Let $\g_1\sub \g_2\sub\cdots$
be a sequence of finite-dimensional
Lie algebras such that $\g_n$ is a subalgebra of
$\g_{n+1}$ for each $n\in\N$.
Give $\g:=\bigcup_{n\in\N}\g_n$
the Lie bracket making each inclusion map
$\g_n\to \g$ a Lie algebra homomorphism.
Then the
locally convex direct limit topology on
$\g$ makes it a topological Lie algebra;
it concides with the finest locally convex topology.
\end{lem}
\begin{prf}
The locally convex direct limit topology on $\g$
concides with the finest locally convex vector topology (see Exercise~\ref{exc-when-finest})
and makes it the direct limit $\dl\g_n$\vspace{-.3mm}
as a topological space (cf.\
Corollary \ref{cp-DL-TVS} and Example~\ref{finest-lcx-vec}).
Each $\g_n$ being locally compact, the product topology on $\g\times \g$
makes it the direct limit topological space\vspace{-.3mm} $\dl(\g_n\times\g_n)$
(see Lemma~\ref{DLcompaprod}). The Lie bracket\vspace{-.3mm}
$\beta_n\colon \g_n\times\g_n\to\g_n$ is continuous, whence also
the Lie bracket
\[ \dl\beta_n\colon \dl(\g_n\times\g_n)\to \dl \g_n \]
of $\g$ is continuous.
\end{prf}
\noindent
\emph{Proof of Theorem}~\ref{dl-lie}.
The group multiplication
$\mu_n\colon G_n\times G_n\to G_n$
is smooth for each $n\in\N$ and also the inversion map
$\eta_n\colon G_n\to G_n$, $g\mto g^{-1}$.
The inversion map $\eta\colon G\to G$
equals\vspace{-.3mm} $\dl\,\eta_n$, whence $\eta$
is smooth by Lemma~\ref{dl-of-ck-maps}.
Also the group multiplication is smooth:
Using the $C^\infty$-diffeomorphism $\Psi\colon
\dl\,(G_n\times G_n)\to G\times G$
from Lemma~\ref{prod-of-dl-mfd}, we have\vspace{-.7mm} $\mu=
\big(\dl \mu_n)\big)\circ \Psi^{-1}$,
where $\dl\mu_n\colon \dl(G_n\times G_n)\to \,\dl G_n$
is smooth by Lemma~\ref{dl-of-ck-maps}.\vspace{.7mm}\qed

We recall an identity for left logarithmic derivatives:
\begin{equation}\label{identity-left-log}
\delta(\eta\xi^{-1})(t)=(\Ad \xi(t))(\delta(\eta)(t)-\delta(\xi)(t))\;\,\mbox{for all $t\in [0,1]$,}
\end {equation}
for each Lie group $G$ and all $\eta,\xi\in C^1([0,1],G)$.
\begin{prop}\label{dl-lie-regular}
Let $G_1\sub G_2\sub\cdots$
be an ascending sequence of Lie groups such that
the inclusion maps $\lambda_{m,n}\colon G_n\to G_m$
are smooth group homomorphisms for all $m\geq n$ in~$\N$.
For the Lie group $G=\bigcup_{n\in\N}G_n$
as in Theorem~{\rm\ref{dl-lie}}, let $\lambda_n\colon G\to G_n$
be the inclusion map. Then we have:
\begin{description}[(D)]
\item[\rm(a)]
$\L(G)$, with the limit maps $\L(\lambda_n)\colon \L(G_n)\to \L(G)$,
is the direct limit of $((\L(G_n))_{n\in \N},(\L(\lambda_{m,n})_{n\leq m})$
as a locally convex space, topological space, set, and as a topological Lie algebra.
\item[\rm(b)]
$G$ is $C^0$-regular.
\end{description}
\end{prop}
\begin{prf}
For each $n\in\N$, we identify $\L(G_n)$ with its image in $\L(G_{n+1})$
under the injective continuous linear map
$\L(\lambda_{n+1,n})$. We endow $\g:=\bigcup_{n\in \N}\L(G_n)$
with the Lie algebra structure making each of the inclusion maps
$\L(G_n)\to\g$ a Lie algebra homomorphism.
We give $\g$ the topology making
it the locally convex direct limit $\dl\, \L(G_n)$,\vspace{-.9mm}
which has the asserted properties by Lemma~\ref{do-not-dodge}.

(a) For $x_0=e$, we have $e\in G_{n_0}$ with $n_0:=1$.
Let $d$ be the supremum of $\dim(G_n)$ for $n\in\N$
and $\phi=\dl\,\phi_n\colon U\to V\sub\R^d$\vspace{-.3mm}
be a chart as in the proof of Theorem~\ref{dl-mfd},
with $\phi_n\colon U_n\to V_n\sub \R^{d_n}$.
Thus $U=\bigcup_{n\in \N}U_n$ and $V=\bigcup_{n\in \N}V_n$.
Moreover, $G_n\cap U=U_n$ and $E_n\cap V=V_n$ for all $n\in \N$.
For $n\in \N$, let $i_n\colon \R^{d_n}\to\R^{d_n}\times \R^{d-d_n}$
be the map $v\mto (v,0)$, which is continuous and linear.
Then $\phi\circ \lambda_n|_{U_n}\circ \phi_n^{-1}=i_n|_{V_n}$, whence
\[
d(\phi\circ \lambda_n\circ \phi_n^{-1})(0,\cdot)=i_n.
\]
Thus $\R^d$ is the union of the images of the maps
$d(\phi\circ\lambda_n)|_{\L(G_n)}$,
whence $\L(G)$ is the union of the images of the maps
$\L(\lambda_n)$. Since $i_n$ is injective, also $\L(\lambda_n)$
is injective, entailing that the linear map
\begin{equation}\label{leftunion}
\psi\colon \g=\bigcup_{n\in \N}\L(G_n)\to \L(G)\vspace{-.5mm}
\end{equation}
determined by $\psi|_{\L(G_n)}=\L(\lambda_n)$
is injective. By the preceding, $\psi$ is surjective.
As both $\g$ and $\L(G)\cong \R^d$ are endowed
with the finest locally convex vector topology,
the bijective linear map $\psi$ is an isomorphism
of topological vector spaces (see Exercise~\ref{exc-maps-between-finest}).
Since $\psi|_{\L(G_n)}$
is a Lie algebra homomorphism
for each $n\in \N$,
$\psi$ is a Lie algebra homomorphism
and hence an isomorphism of Lie algebras.
Since $\L(G)$ is a topological Lie algebra, so is~$\g$.

(b) If $\gamma\in C([0,1],\L(G))$,
then $\psi^{-1}\circ \gamma\colon [0,1]\to \g$
is continuous. As $\psi^{-1}(\gamma([0,1])\sub \g$
is compact, there exists $n\in\N$
such that $\psi^{-1}(\gamma([0,1])\sub \L(G_n)$, see Proposition~\ref{henceSilvacoreg}.
The inclusion map $\L(G_n)\to \g$ being a topological
embedding (cf.\ Proposition~\ref{onlyyou}), we deduce that
\[ \zeta:=(\psi^{-1}\circ\gamma)|^{\L(G_n)}\colon [0,1]\to \L(G_n) \] 
is continuous. The finite-dimensional Lie group $G_n$
being $C^0$-regular (Theorem~\ref{thm:banach-regular}),
there exists a $C^1$-function $\eta\colon [0,1]\to G_n$ with $\eta(0)=e$
and left logarithmic derivative $\delta(\eta)=\zeta$.
Then $\lambda_n\circ \eta\colon [0,1]\to G$
is $C^1$, $\lambda_n(\eta(0))=e$, and $\delta(\lambda_n\circ \eta)=\L(\lambda_n)\circ
\delta(\eta)=\psi|_{\L(G_n)}\circ (\psi^{-1}\circ\gamma)|^{\L(G_n)}=\gamma$.
Thus $G$ is $C^0$-semiregular.
To see that $G$ is $C^0$-regular,
it only remains to show that $\Evol\colon C([0,1],\L(G))\to C([0,1],G)$
is continuous at~$0$ (see Theorem~\ref{thm-cts-0-reg}).
By the preceding discussion, we can identify $\L(G_n)$
with its image in $\L(G)$ under the injective Lie algebra
homomorphism $\L(\lambda_n)$.
Then $\L(G)=\bigcup_{n\in\N} \L(G_n)$.
We set $F_1:=\L(G_1)$ and write $\L(G_n)=\L(G_{n-1})\oplus F_n$
for integers $n\geq 2$
using a complementary vector subspace $F_n$ for $\L(G_{n-1})$ in $\L(G_n)$.
Then
\[
\L(G)=\bigoplus_{n\in \N}F_n
\]
as a locally convex space.
Let $U\sub C([0,1],G)$
be an identity neighborhood.
After shrinking~$U$,
we may assume that $U=C([0,1],V)$
for some identity neighborhood $V\sub G$ (see Lemma~\ref{sammelsu}(b)).
Setting $V_0:=V$,
there are open identity neighborhoods $V_n$
in $G$ for $n\in\N$ such that $V_nV_n\sub V_{n-1}$.
Then
\[
V_n\cdots V_1\sub V_0\;\,\mbox{for all $n\in\N$.}
\]
After shrinking $V_n$, we may assume that $V_n=\bigcup_{j\in \N}V_{n,j}$
for a
sequence
\[
V_{n,1}\sub V_{n,2}\sub \cdots
\]
of relatively compact identity neighborhoods $V_{n,j}\sub G_j$, see Lemma~\ref{la-rel-cp-nbds}.
We claim that there exist open $0$-neighborhoods $W_n\sub F_n$
for $n\in \N$ with
\begin{equation}\label{gives-reg-DL}
\hspace*{-4.2mm}\Evol_{G_n}(\gamma)\in C([0,1], V_{n,n}\cdots V_{1,1}) \;
\mbox{for all $\gamma\in C([0,1],W_1\oplus\cdots\oplus W_n)$.}\hspace*{-3.2mm}
\end{equation}
Then $W:=\bigoplus_{n\in \N}W_n$ is an open $0$-neighborhood in $\L(G)$.
For each $\gamma\in C([0,1],W)$,
we have $\gamma([0,1])\in \bigoplus_{j=1}^nF_j$
for some $n\in \N$ by Lemma~\ref{firstlasum}(e)
and thus $\gamma\in C([0,1],W_1\oplus\cdots\oplus W_n)\sub C([0,1], \L(G_n))$,
whence
\[
\Evol_G(\gamma)=\Evol_{G_n}(\gamma)\in C([0,1],V_{n,n}\cdots V_{1,1})
\sub C([0,1],V_n\cdots V_1)\sub U.
\]
Thus $\Evol_G(W)\sub U$, whence $\Evol_G$ is continuous at~$0$.
It only remains to prove the claim.
For $n\in \N$, $\Evol_{G_n}^{-1}(C([0,1], V_{n,n}))$
is an open $0$-neighborhood in $C([0,1],\L(G_n))$,
since $\Evol_{G_n}\colon C([0,1],\L(G_n))\to C([0,1],G_n)$
is continuous. Then
$C([0,1],Q_n)\sub \Evol_{G_n}^{-1}(C([0,1], V_{n,n}))$
for an open $0$-neighborhood $Q_n\sub \L(G_n)$, by Lemma~\ref{sammelsu}(b).
If $n=1$, we set $W_1:=Q_1$.
If $n\geq 2$, we use that restriction
\[
h_n \colon G_n\times F_n \to \L(G_n),\;\,
(g,v)\mto \Ad_{G_n}(g)(v)
\]
of the adjoint action is continuous and $h_n((\wb{V_{n-1,n-1}}\cdots \wb{V_{1,1}})\times \{0\})=\{0\}\sub Q_n$; here $\wb{V_{j,j}}$ denotes the closure of
$V_{j,j}$ in $G_j$ for $j\in \{1,\ldots, n-1\}$,
which is a compact set.
By the Wallace Lemma, there exists an open $0$-neghborhood
$W_n\sub F_n$ such that $h_n((\wb{V_{n-1,n-1}}\cdots\wb{V_{1,1}})\times W_n)\sub Q_n$,
whence
\begin{equation}\label{givesreg-2}
\Evol_{G_n}(t\mto \Ad_{G_n}(\alpha(t))(\beta(t)))\in C([0,1],V_{n,n})
\end{equation}
for all $\alpha\in C([0,1], V_{n-1,n-1}\cdots V_{1,1})$
and $\beta\in C([0,1], W_n)$.
We now show by induction that (\ref{gives-reg-DL})
holds for all $n\in\N$.
If $n=1$, then
$\Evol_{G_1}(\gamma)\in C([0,1],V_{1,1})$
for all $\gamma\in C([0,1],W_1)$, by definition of~$W_1$.
If $n\geq 2$ and (\ref{gives-reg-DL})
holds for $n-1$ in place of~$n$,
let $\gamma=\sum_{j=1}^n\gamma_j$
with $\gamma_j\in C([0,1],W_j)\sub C([0,1], F_j)$.
Thus $\gamma=\theta+\gamma_n$
with
\[ \theta:=\sum_{j=1}^{n-1}\in C([0,1],W_1\oplus\cdots\oplus W_{n+1})\sub C([0,1], \L(G_{n-1})),\]
whence $\Evol_{G_{n-1}}(\theta)\in C([0,1],V_{n-1,n-1})$
by the inductive hypothesis.
Now
\begin{eqnarray*}
\Evol_{G_n}(\gamma) &=&\Evol_{G_n}(\gamma)(\Evol_{G_n}(\theta))^{-1}\Evol_{G_n}(\theta)\\
&=& \Evol_{G_n}(\gamma)(\Evol_{G_n}(\theta))^{-1}\Evol_{G_{n-1}}(\theta)
\end{eqnarray*}
where $\Evol_{G_{n-1}}(\theta)\in C([0,1],V_{n-1,n-1}\cdots V_1)$
by induction. By (\ref{identity-left-log}) and (\ref{givesreg-2}),
$\Evol_{G_n}(\gamma)(t)(\Evol_{G_n}(\theta)(t))^{-1} =
\Evol_{G_n}(\Ad_{G_m}(\Evol_{G_n}(\theta))(\gamma-\theta))(t)$
$=\Evol_{G_n}(\Ad_{G_n}(\Evol_{G_n}(\theta))(\gamma_n))(t)$ is in
\[
h_n((\wb{V_{n-1,n-1}}\cdots\wb{V_{1,1}})\times W_n)\sub C([0,1],V_{n,n})
\]
for all $t\in[0,1]$, whence
$\Evol_{G_n}(\gamma)\in C([0,1],V_{n,n})C([0,1],V_{n-1,n-1}\cdots V_{1,1})$
$\sub C([0,1], V_{n+1,n+1}\cdots V_{1,1})$.
\end{prf}
\begin{prop}\label{homog-of-DL}
If $H$ is a closed subgroup of $G=\bigcup_{n\in \N}G_n$
in the situation of Theorem~{\rm\ref{dl-lie}},
then $H$ is a split Lie subgroup of~$G$
and $G/H$ admits a smooth manifold structure
turning the canonical map $q\colon G\to G/H$, $g\mto gH$
into a submersion. For this smooth manifold structure, we have
$G/H=\dl\, G_n/(H\cap G_n)$ as a smooth manifold.
\end{prop}
\begin{prf}
For each $n\in \N$, $H_n:=H\cap G_n$ is a closed
subgroup of~$G_n$ and hence a submanifold.
We give $M_n:=G_n/H_n$
the unique smooth manifold structure making $q_n\colon G_n\to M_n$,
$g\mto gH_n$ a submersion. Then $\ker T_e(q_n)=T_e(H_n)$.
For all $m\geq n$ in~$\N$, the map $q_{m,n}\colon M_n\to M_m$, $gH_n\mto g H_m$
is injective,
and smooth since $q_{m,n}\circ q_n=q_m$ is smooth
(see Proposition~\ref{subm-sec}(b)).
The tangent map $T_{H_n}(q_{m,n})$ is injective and thus
$q_{n,m}$ an immersion, since $q_{n,m}$ is equivariant for the
left action of $G_n$ on $M_n$ and $M_m$.
In fact, if $v\in \ker(T_{H_n}q_{m,n})$, then
$v=Tq_n(w)$
for some $w\in T_eG$ and $0=Tq_{m,n}(Tq_n(w))=Tq_m(w)$.
whence $w\in T_eG_n\cap T_eH_m=T_eH_n$ and thus $v=0$.
Consider the canonical map $q\colon G\to G/H$, $g\mto gH$
and the injective mappings $j_n\colon M_n\to G/H$, $gH_n\mto gH$.
The topology on $G$ being final with respect to
the inclusion maps $\lambda_n\colon G_n\to G$,
the quotient topology on $G/H$ is final with respect to
the mappings $q\circ \lambda_n=j_n\circ q_n$,
and hence with respect to the mappings $j_n$
(using the transitivity of final topologies).
Hence $G/H=\dl\, M_n$ as a topological space.
In the following, we identify $gH_n\in M_n$ with $gH\in G/H$.
Now Theorem~\ref{dl-mfd} provides a smooth manifold structure
on $G/H$ making it the direct limit of the smooth manifolds
$M_n$. Then $q=\dl\, q_n$ is smooth.
We show that $q$ admits smooth local sections.
Then $q$ is a submersion (by Lemma~\ref{GH-principal-sub})
and $G$ with $q$ a smooth $H$-principal bundle.
Let $g\in G$ and $n_0\in \N$ such that $g\in G_{n_0}$.
Let $\phi_n\colon U\to V$
be a chart for $\dl\,M_n=G/H$ around $gH$
as in the proof of Theorem~\ref{dl-mfd}.
Thus $U=\bigcup_{n\geq n_0}U_n$
with relatively compact, smoothly contractible
open subsets $U_n\sub M_n$
such that $U_n$ is a submanifold of $U_{n+1}$
and there exists a smooth retraction $r_n\colon U_{n+1}\to U_n$.
For each $n\geq n_0$, $q_n\colon G_n\to M_n$ restrict to
an $H_n$-principal bundle
$q_n^{-1}(U_n)\to U_n$ which is trivalizable as $U_n$
is smoothly contractible. Let $\theta_n\colon q_n^{-1}(U_n)\to U_n\times H_n$
be a local trivialization with second component $\theta_{n,2}\colon q_n^{-1}(U_n)\to H_n$.
Then a sequence $(\sigma_n)_{n\geq n_0}$
of smooth local sections $\sigma_n\colon U_n\to G_n$
can be found for $n\geq n_0$
such that $\sigma_n|_{U_{n-1}}=q_{n,n-1}\circ \sigma_{n-1}$
for all $n > n_0$.
We start with $\sigma_{n_0}(x):=\theta_n^{-1}(x,e)$.
If $\sigma_{n_0},\ldots,\sigma_n$ have been found,
we can take
\[
\sigma_{n+1}(x):=\theta_{n+1}^{-1}(x,(\theta_{n+1,2}\circ \sigma_n\circ r_n)(x)).
\]
Then $\sigma:=\dl_{n\geq n_0}\sigma_n\colon \dl_{n\geq n_0}U_n\to \dl_{n\geq n_0}G_n$
is a smooth map from $U$ to $G$, using~Lemma~\ref{dl-of-ck-maps}.
By construction, $q\circ \sigma=\id_U$.
Since $q$ is a submersion, $H=q^{-1}(\{H\})$
is a submanifold of~$G$ (see Proposition~\ref{subm-sec}(c)).
As $q\colon G\to G/H$ is an $H$-principal bundle,
$H$ actually is split Lie subgroup of~$G$ (cf.\ Lemma~\ref{GH-principal-sub}).
\end{prf}
Consider a direct sequence $((G_n)_{n\in\N}, (\lambda_{m,n})_{m\geq n})$
of finite-dimensional Lie groups and smooth group
homomorphisms $\lambda_{m,n}\colon G_n\to G_m$.
We can then associate a direct sequence of Lie groups with injective bonding maps.
\begin{numba}\label{fi-set-non-inj}
In fact,
$(\ker(\lambda_{m,n}))_{m\geq n}$
is an ascending sequence of closed normal subgroups of $G_n$
for each $n\in \N$, whence
\[
K_n:=\bigcup_{m\geq n} \ker(\lambda_{m,n})
\]
is a normal subgroup of $G_n$.
Let $H:=\dl \,G_n$\vspace{-.3mm} as a set and a group,
with limit maps $\Lambda_n\colon G_n\to H$.
We endow $H$ with the final topology with respect to the
family $(\Lambda_n)_{n\in\N}$,
making it the direct limit $\dl\, G_n$\vspace{-.3mm}
as a topological space.
If $h \in H$, then $h=\Lambda_{n_0}(g)$
for some $n_0\in \N$ and $G\in G_n$.
For $n\geq n_0$, the left translations $\lambda^{G_n}_g\colon G_n\to G_n$, $x\mto gx$
are continuous.
Now the left translation $\lambda_h^H\colon H\to H$, $x\mto hx$
is the map $\dl_{n\geq n_0}\lambda_g^{G_n}$ and hence
continuous (cf.\ Definition~\ref{DLmaps}).
Likewise, the right translation $\rho_h^H\colon H\to H$
is continuous and\vspace{-.3mm} the inversion map $H\to H$, $h\mto h^{-1}$.
As a consequence, the closure $N:=\wb{\{e\}}$
of the trivial subgroup in~$H$
is a closed normal subgroup of~$H$.
Thus
\[
N_n:=\Lambda_n^{-1}(N)
\]
is a closed normal subgroup of $G_n$.
Endow $\wb{G}_n:= G_n/N_n$
with the unique Lie group structure turning
the canonical map $q_n\colon G_n\to \wb{G}_n$, $g\mto gN_n$
into a submersion.
By construction, we have
\begin{equation}\label{this-gives-inj}
N_n=\lambda_{m,n}^{-1}(N_m)\;\;\,\mbox{for all $\,m\geq n$.}
\end{equation}
Thus
$\lambda_{m,n}(N_n)\sub N_m$, whence
\[
\bar{\lambda}_{m,n}\colon \wb{G}_n\to\wb{G}_m,\;\;
gN_n\mto gN_m
\]
is a well-defined group homomorphism
and smooth since $\bar{\lambda}_{m,n}\circ q_n=q_m$
is smooth (see Proposition~\ref{subm-sec}(b)).
By construction, $\bar{\lambda}_{m,n}$ is injective
for all $m\geq n$ in~$\N$.
Then $((\wb{G}_n)_{n\in \N}(\bar{\lambda}_{m,n})_{m\geq n})$
is a direct sequence of finite-dimensional Lie groups;
stimulated by \cite{NRW91},
we call it the associated \emph{injective quotient system}.
\end{numba}
\begin{numba}
The maps $\lambda_n\colon \wb{G}_n\to H/N$, $gN_n\mto gN$
are injective group homomorphism for $n\in \N$
and $G/N$ is the ascending union of the images.
Hence
\[
H/N=\dl\, \wb{G}_n\vspace{-.3mm}
\]
as a set and as a group, with limit maps $\lambda_n$.
We give $G:=H/N$ the quotient topology with respect to~$q$.
By transitivity of final topologies,
this topology is final with respect to the maps $q\circ \Lambda_n=\lambda_n\circ q_n$
for $n\in\N$ and hence final with respect to the maps $\lambda_n$.
Thus $G=\dl\, \wb{G}_n$\vspace{-.3mm} as a topological space.
\end{numba}
\begin{prop}\label{dl-lie-non-injective}
Let $((G_n)_{n\in \N},(\lambda_{m,n})_{m\geq n})$
be a direct sequence of Lie groups
and $((\wb{G}_n)_{n\in \N}, (\wb{\lambda}_{m,n})_{m\geq n})$
be the associated injective quotient system,
with canonical maps $q_n\colon G_n\to\wb{G}_n$
for $n\in\N$.
Let $G:=\dl\,\wb{G}_n$\vspace{-.5mm}
be its direct limit Lie group as in Theorem~{\rm\ref{dl-lie}},
with limit maps $\lambda_n\colon \wb{G}_n\to G$.
Then the following holds:
\begin{description}[(D)]
\item[\rm(a)]
$G$, with the limit maps $\lambda_n\circ \,q_n$,
is the direct limit of $((G_n)_{n\in \N},(\lambda_{m,n})_{m\geq n})$
as a Lie group, Hausdorff topological space, Hausdorff topological
group and also in the category of $C^r$-manifolds with rough boundary
and $C^r$-maps between them,
for each $r\in \N_0\cup\{\infty\}$.
\item[\rm(b)]
$\L(G)\hspace*{-.4mm}=\hspace*{-.4mm}\dl \hspace*{.05mm}\L(\wb{G}_n)$\vspace{-.3mm} with bonding maps $\L(\lambda_{m,n})$
and limit maps $\L(\lambda_n\hspace*{-.2mm}\circ\hspace*{-.1mm}q_n)$.
\item[\rm(c)]
If $G_n$ induces the discrete topology on
$K_n$ for each $n\in\N$,
then $N_n=K_n$
and $\L(G)=\dl \hspace*{.2mm}\L(G_n)$.
\end{description}
\end{prop}
\begin{prf}
(b) holds by Theorem~\ref{dl-lie-regular}.

(a) Let $Y$ be a Hausdorff topological space and
$(f_n)_{n\in\N}$
be a sequence of continuous maps $f_n\colon G_n\to Y$
such that $f_m\circ \lambda_{m,n}=f_n$
for all $m\geq n$ in~$\N$.
Then there is a unique continuous map $f\colon H\to Y$
such that $f\circ\Lambda_n=f_n$ for all $n\in\N$.
For each $x\in H$, the singleton $\{f(x)\}$
is closed in the Hausdorff space~$Y$,
whence $x N=\wb{\{x\}}\sub f^{-1}(\{f(x)\})$.
Thus $g\colon G\to Y$, $q(x)\mto f(x)$
is a well defined map and continuous
since $g\circ q=f$ is continuous.
If $Y$ is a Hausdorff group and each $f_n$ a continuous
group homomorphism, then so is~$g$.
If $Y$ is a $C^r$-manifold (possibly with rough boundary$)$
and each $f_n$ a \break $C^r$-map,
then $g\circ \lambda_n\circ q_n=g\circ q\circ \circ \Lambda_n=f_n$
is $C^r$ for each $n\in\N$, whence $g\circ\lambda_n$
is $C^r$ by Proposition~\ref{subm-sec}(b)
(noting Remark~\ref{subm-sec-bd}).
Since $G=\dl\,\wb{G}_n$
as a $C^r$-manifold
by Theorem~\ref{dl-lie},
we deduce that $g$ is $C^r$.

(c) If each $K_n$ is discrete,
we define $\wb{Q}_n:=G_n/K_n$ for $n\in \N$
and let $p_n\colon G_n\to Q_n$
be the canonical map.
For $M\geq n$ in $\N$,
we let $p_{m,n}\colon \wb{Q}_n\to\wb{Q}_m$
be the smooth group homomorphism
determined by $p_{m,n}\circ p_n=p_m$.
Then $p_{m,n}$
is injective and we can form the direct limit Lie group
$Q=\dl\, Q_n$\vspace{-.3mm}
of the direct sequence $((Q_n)_{n\in\N},(p_{m,n})_{m\geq n})$,
with limit maps $j_n\colon Q_n\to Q$ which are injective.
For $H$ as above, there exists a continuous
map $j\colon H\to Q$ such that $j\circ\Lambda_n=j_n\circ p_n$,
entailing that $N_n=K_n$.
\end{prf}
\begin{ex}
Let $G_n:=\bS$ be the circle group for each
$n\in \N$ and $\lambda_{m,n}\colon \bS\to\bS$ be the map $z\mto z^{2^{m-n}}$.
Then $K_n$ is the group of all $z\in \bS$ such that $z^{2^k}=1$
for some $k\in \N$.
Thus $K_n$ is dense in $G_n$ and hence $\wb{K_n}=G_n\sub \ker(\lambda_n\circ q_n)$,
showing that $(\lambda_n\circ q_n)(g)=e$ for all $g\in G_n$.
As a consequence, $G:=\dl \,G_n=\{e\}$\vspace{-.3mm} as a Lie group and $\L(G)=\{0\}
\not=\dl \, \L(G_n)\cong\R$. As an abstract group, $\dl \,G_n=\dl \,G_n/K_n$,\vspace{-.3mm}
which is uncountable and hence not isomorphic to~$G$.
\end{ex}
A Lie algebra is {\em locally finite\/}
if every finite subset generates a finite-dimensional
subalgebra.
\begin{thm}\label{intthm}
Let $\g$ be a
locally finite Lie algebra
of countable dimension.
Then the finest locally convex topology
makes $\g$ a topological Lie algebra.
Then there exists a simply connected, connected, $C^0$-regular Lie group $G$
such that $\L(G)\cong \g$ and $G=\dl\hspace*{.2mm}G_n$
for some ascending sequence $G_1\sub G_2\cdots$
of finite-dimensonal Lie groups.
\end{thm}
The following fact is useful for the proof.
\begin{lem}\label{vectgp}
Let $G$ be a simply connected, connected,
finite-dimensional real Lie group.
Let $Z(G)$ be the centre of $G$
and $Z(G)_0$ its connected component of the identity.
Then
$Z(G)_0$ is a vector group,
i.e., $Z(G)_0\cong \R^m$ for some $m\in \N_0$.
\end{lem}
\begin{prf}
By L\'{e}vi's Theorem,
$\L(G)=\fr\rtimes \fs$
internally,
where $\fr$ is the radical of~$\L(G)$
and $\fs$ a L\'{e}vi
complement
(see
\cite[Ch.\,I, \S6.8, Thm.\,5]{Bou89}).
Let
$R$ and $S$ be the
integral subgroups of $G$ corresponding to
$\fr$ and $\fs$, respectively.
Then $R$ and $S$ are simply connected,
$R$ is a closed normal subgroup,
$S$ a closed
subgroup of~$G$, and $G=R\rtimes S$ internally
(see \cite[Prop.\,11.1.19]{HiNe12}).
Now consider the identity component
$Z(G)_0$ of the centre $Z(G)$ of~$G$.
Let \break $\pi \colon  G\to S$ be the projection
onto~$S$, with kernel $R$.
Then $\pi(Z(G)_0)\sub Z(S)_0=\{e\}$,
entailing that $Z(G)_0\sub R$.
Thus $Z(G)_0$ is an integral subgroup
of the simply connected solvable Lie group~$R$,
whence $Z(G)_0$ is simply connected
\cite[Prop.\,11.2.15]{HiNe12}.
Being a simply connected abelian
Lie group, $Z(G)_0$ is a vector group.
\end{prf}
\noindent
\emph{Proof of Theorem}~\ref{intthm}.
As $\g$ is locally finite and
of countable dimension, $\g$ is the union of
some ascending sequence
$\g_1\sub \g_2\sub\cdots$
of finite-dimensional Lie subalgebras
of~$\g$.
The finest locally convex topology on $\g$
makes it the direct limit $\g=\dl \g_n$
as a locally convex space, and a topological
Lie algebra (see Lemma~\ref{do-not-dodge}).
For $n\in \N$, let
$G_n$ be a simply connected, connected Lie group
with Lie algebra $\L(G_n)\isom \g_n$;
fix an isomorphism $\kappa_n\colon \L(G_n)\to \g_n$.
If $m\geq n$, then
the injective Lie algebra homomorphisms
\[ \kappa_{m,n}:=\kappa_m^{-1}\circ \kappa_n\colon
\L(G_n)\to \L(G_m) \] 
make $(\L(G_n),(\kappa_{m,n})_{m\geq n})$
a directed sequence of Lie algebras.
Lemma~\ref{do-not-dodge} permits us to
form the direct limit $\dl\hspace*{.2mm} \L(G_n)$
as a topological Lie algebra;
the map $\kappa:=\dl\kappa_m\colon \g=\dl\hspace*{.2mm}\g_n\to \dl \hspace*{.2mm}\L(G_n)$
is an isomorphism of topological Lie algebras.
For all $m\geq n$,
there exists a unique smooth group homomorphism
$\lambda_{m,n}\colon G_n\to G_m$
such that $\L(\lambda_{m,n})=\kappa_{m,n}$
(see Theorem~\ref{thm3.2.11}).
Since
$\L(\lambda_{k,m}\circ \lambda_{m,n})=\L(\lambda_{k,n})$,
we have
\begin{equation}\label{hence-ascending}
\lambda_{k,m}\circ\lambda_{m,n}=\lambda_{k,n}\;\;\,\mbox{for all $\,k\geq m\geq n$,}
\end{equation}
whence $((G_n)_{n\in\N},(\lambda_{m,n})_{m\geq n})$ is a
direct sequence of finite-dimensional Lie groups.
Let $K_n$ be as in \ref{fi-set-non-inj}.
We claim that $K_n$ is discrete
for each $n\in\N$,
whence $\wb{G}_n:=G_n/K_n$ admits a Lie group
structure making $q_n\colon G_n\to \wb{G}_n$, $g\mto gK_n$
a local $C^\infty$-diffeomorphism.
By parts (a) and (c) of Proposition~\ref{dl-lie-non-injective},
we can form the direct limit $G=\dl\, G_n/K_n$
in the category of Lie groups
and smooth group homomorphisms, which is $C^0$-regular
by Proposition~\ref{dl-lie-regular},
and we also have $G=\dl\,G_n$ as a Lie group.
Let $\lambda_n\colon \bar{G}_n\to G$ be the limit map.
Moreover, $\L(G)\cong\dl \L(\wb{G}_n)$,
by Proposition~\ref{dl-lie-non-injective}(b),
using the bonding maps $\L(\lambda_{m,n})$.
As also $\dl \L(q_n)\colon \dl \L(G_n)\to\dl \L(\wb{G}_n)$
is an isomorphism of topological Lie algebras, we have that
$\g\cong \dl \L(G_n)\cong \dl \L(\wb{G}_n)\cong \L(G)$
as a topological Lie algebra.
Let $\wt{G}$
be the universal covering Lie group of~$G$
and $\pi\colon \wt{G}\to G$
be a surjective group homomorphism
which is a local $C^\infty$-diffeomorphism.
Since $G$ is regular, also $\wt{G}$ is regular (see Proposition~\ref{prop:4.1.5}).
For each $n\in \N$, the continuous Lie algebra
homomorphism $\L(\pi)\circ \L(\lambda_n)$
integrates to a smooth group homomorphism
$j_n\colon \wb{G}_n\to \wt{G}$ (see Theorem~\ref{thm3.2.11}).
Then $j_m\circ\lambda_{m,n}=j_n$,
whence there exists a smooth group homomorphsm
$j\colon G=\dl\,\wb{G}_n\to\wt{G}$
such that $j\circ \lambda_n=j_n$ for all $n\in \N$.
Then $\L(\pi\circ j)=\id_{\L(G)}$ and $\L(j\circ \pi)=\id_{\L(\wt{G})}$
imply
that $\pi\circ j=\id_G$ and $j\circ \pi=\id_{\wt{G}}$.
Thus $G\cong \wt{G}$ is simply connected.

It remains to show that each $K_n$ is discrete.
We show that the closure
$N_n:=\wb{K_n}\sub G_n$ is discrete.
The homomorphism $\lambda_{m,n}$ has discrete
kernel for all $m,n\in \N$ with $m\geq n$,
because $\L(\lambda_{m,n})=\kappa_{m,n}$ is
injective. Now $\ker\lambda_{m,n}$ being a discrete
normal subgroup of the connected group $G_n$,
it is central. Thus $N_n\sub Z(G_n)$,
for each~$n$,
whence $(N_n)_0\sub Z(G_n)_0$ is a vector group
being a connected closed subgroup
of a vector group (Lemma~\ref{vectgp}).
As a consequence of (\ref{hence-ascending}),
we have
$\lambda_{m,n}(N_n)\sub N_m$ for all $m\geq n$,
whence $\lambda_{m,n}((N_n)_0)\sub (N_m)_0$.
Being a continuous homomorphism between
vector groups, $\psi_{m,n}:=\lambda_{m,n}|_{(N_n)_0}\colon
(N_n)_0\to (N_m)_0$
is real linear.
Hence, being a linear map with
discrete kernel, $\psi_{m,n}$ is
injective.
Thus $(N_n)_0=\wb{\bigcup_{m\geq n}\ker \psi_{m,n}}=\{e\}$,
whence~$N_n$ is discrete. \qed
\begin{cor}
Every connected regular Lie group $G$ whose Lie algebra is countably 
dimensional, locally finite and carries the finest locally convex topology is a direct limit 
of finite-dimensional Lie groups.
\end{cor}
\begin{prf}
By Theorem~\ref{intthm} and Proposition~\ref{dl-lie-regular},
there exists a direct limit $H:=\dl\hspace*{.2mm}G_n$
of finite-dimensional Lie groups $G_n$ with $\L(H)\cong \L(G)$,
such that~$H$ is connected, simply connected and
a regular Lie group. The universal covering group~$\wt{G}$
being regular, we have $H\cong \wt{G}$ as Lie groups as
a consequence of Theorem~\ref{thm3.2.11}.
Hence $G\cong H/N$ for a discrete normal subgroup $N\sub H$.
By Proposition~\ref{homog-of-DL}, $H/N\cong G_n/(H\cap G_n)$.
\end{prf}
\section[Further facts concerning locally finite-dimensional Lie groups]{Further facts concerning locally finite-dimensional\\
\hspace*{10.37mm}Lie groups}\label{more-findim-DL}
In the following, we shall call the class of Lie groups described by
Theorem~\ref{dl-lie}
{\it locally finite-dimensional} (regular) Lie groups. 
We mention several facts without proof,
with references to the literature.
Recall that a topological group~$G$
is said to have \emph{no small subgroups}
if there exists an identity neighborhood $U\sub G$
which does not contain a subgroup
of~$G$ as a subset except for the trivial subgroup.
\begin{thm}
  \mlabel{thm:VII.1.3} {\rm(\cite{Gl05c, Gl07b})} 
Let $G$ be a locally finite-dimensional Lie group. Then the 
following assertions hold: 
\begin{description}[(D)]
\item[\rm(a)] Every subalgebra $\h \subeq \L(G)$ integrates to an integral subgroup.
\item[\rm(b)] Every locally compact subgroup $H \subeq G$ 
is a finite-dimensional Lie group. 
\item[\rm(c)] $G$ does not contain small subgroups. 
\end{description}
\end{thm}
%
%
\begin{thm}{\rm\cite{Gl05c}}
  \mlabel{thm:VII.1.5}  Every continuous homomorphism between locally 
finite-dimensional Lie groups is smooth. 
\end{thm}
In the same way as 
the corresponding result for locally exponential Lie groups 
(Theorem~\ref{thm:5.4.4}), the preceding theorem implies that 
locally finite-dimensional Lie groups form a full sub-category 
of the category of topological groups. 
We even have the following stronger version of the preceding theorem: 
\begin{thm}
  \mlabel{thm:VII.1.6} Let $G = \indlim G_n$ be a locally finite-dimensional 
Lie group and $H$ be a Lie group. 
\begin{description}
\item[\rm(a)] A group homomorphism $\phi \: G \to H$ is smooth if and only if the corresponding 
homomorphisms $\phi_n \: G_n \to H$ are smooth. 
\item[\rm(b)] If $H$ has a smooth exponential map, then every 
continuous homomorphism $\phi \: G \to H$ is smooth. 
\end{description}
\end{thm}

\begin{prf}
 (a) is contained in \cite{Gl05c}. 

In view of (a), part (b) follows 
from the finite-dimensional case, which in turn follows from the existence of 
local coordinates of the second kind: 
$(t_1,\ldots, t_n)$ $\mapsto 
\prod_{i = 1}^n \exp_G(t_i x_i).$ 
\end{prf}
\begin{rem}
  \mlabel{rem:VII.1.2} (a) Beyond countable directed systems, several serious obstacles arise. 
For vector spaces of countable dimension, the finest locally 
convex topology coincides with the finest topology for which all 
inclusions of finite-dimensional subspaces are continuous. This is crucial 
for many arguments in this context. If $E$ is not of countable dimension, 
the addition on $E$ is not continuous for the latter topology. Similar problems 
occur for uncountable direct limits of topological groups: in many cases the direct 
limit topology does not lead to a continuous multiplication 
(cf.\ \cite{Gl03a} for more details). 

(b) Any vector space $E$ of countable dimension, endowed with the finest locally convex topology,
can be considered as a direct limit space of finite-dimensional subspaces $E_n$ of
$\dim E_n = n$, see Example~\ref{finest-lcx-vec}.
Since each $E_n$ is a closed subspace which is Banach, and all 
inclusions $E_n \to E_{n+1}$ are compact operators, $E$ is an LF space and a Silva space 
at the same time.
\end{rem}
%
%
%
\section{Direct limits of infinite-dimensional Lie groups} 
\mlabel{sec:VII.2}
Direct limit constructions also play an important role when applied to sequences of 
infinite-dimensional Lie groups. On the level of Banach-, resp., Fr\'echet spaces, 
different types of directed systems lead to the important classes of 
LF-spaces and Silva spaces.\\[1mm]
Many interesting Lie groups $G$ like mapping groups and other classes embed
naturally into box products (also
called direct sums,
restricted direct products, or weak direct products),
as constructed in~Theorem~\ref{box-prod-lie}.
This often gives a nice atlas for~$G$
and can reduce differentiability questions of mappings to $G$
to questions concerning mappings to box products,
which may be easier to access.
We recall from \cite{Gl03b}:
\begin{thm}
  \mlabel{thm:VII.2.1}  {\rm(\cite{Gl03c})} 
If $(G_j)_{j \in J}$ is a countable family of 
Lie groups, then their box product
\[  G : =\bigoplus_{j \in J} G_j := \Big\{ (g_j)_{j \in J} \in \prod_{j \in J} G_j \: 
|\{i \: g_j \not=\1\}| < \infty\Big\} \] 
carries a natural Lie group structure, where 
$\L(G) \cong \bigoplus_{j \in I} \L(G_j)$ carries the
box topology making it the locally convex direct sum.
If each $G_j$ is locally exponential, then also~$G$.
\end{thm}
An analogous result is available for
uncountable families,
using a suitable ``fine box product''
modeled on the locally convex direct sum (see loc.\,cit.,
in different terminology).
For $J=\N$,
$\bigoplus_{n\in \N}G_n$ is the direct limit
of the partial products $\prod_{j=1}^nG_j$
in the category of smooth Lie groups
and smooth group homomorphisms,
but not necessarily as a smooth manifold
(see \cite{Gl07a}).\\[1mm]
An example of major interest are diffeomorphism groups.
Consider the Lie group $\Diff_c(M)$
of compactly supported smooth diffeomorphsms
of a $\sigma$-compact finite-dimensional smooth manifold~$M$,
as in Chapter~\ref{ch:diffeo} 
(and going back to \cite{Mr80}).
It is modeled on the LF-space $\cV_c(M)$
of compactly supported smooth vector fields.
\begin{thm}
  \mlabel{thm:VII.2.2}  {\rm(\cite{Gl07a})} 
Let $M$ be a $\sigma$-compact, non-compact smooth manifold
of  positive dimension and $(L_n)_{n\in \N}$
be a compact exhaustion of~$M$. 
Then the Lie group $\Diff_c(M)$ of compactly supported diffeomorphisms 
is not a direct limit of the subgroups $\Diff_{L_n}(M)$
in the category of smooth manifolds, but a homomorphism 
$\Diff_c(M) \to H$ to a Lie group $H$ is smooth if and only if it is smooth on 
each subgroup $\Diff_{L_n}(M)$. 
\end{thm}
Thus $\Diff_c(M)$ is the direct limit of the $\Diff_{L_n}(M)$
in the category of smooth Lie groups and smooth group homomorphisms.
A crucial tool for the proof of the theorem is the following lemma
from~\cite{Gl07a},
which exploits box products of Lie groups.
\begin{lem}
  \mlabel{lem:VII.2.3} {\rm(Fragmentation Lemma)} Let $M$ be a $\sigma$-compact 
manifold of finite dimension. 
Then there exists a locally finite cover $(K_n)_{n \in \N}$ of $M$ by compact 
sets, an open identity neighborhood $U \subeq \Diff_c(M)$ and a smooth mapping 
$\Phi \: U \to \bigoplus_{n \in \N} \Diff_{K_n}(M)$
which satisfies 
$\gamma = \Phi(\gamma)_1 \circ \ldots \circ \Phi(\gamma)_n$ 
for each $\gamma \in U$. 
\end{lem}
Let us retain
a $\sigma$-compact finite-dimensional manifold~$M$.
Let $K$ be a Lie group and $C^\infty_c(M,K)$ be the Lie group
of compactly supported 
smooth maps $M \to K$,
as discussed in Theorem~\ref{thm-tefu-gp}
(see also \cite{Gl02c}). 
Then $C^\infty_c(M,K)$
is the union of the directed system of subgroups 
$C^\infty_L(M,K) := \{ f \in C^\infty(M,K) \: \supp(f) \subeq L\}$
for compact subsets $L\sub M$.
If $K$ is a Banach--Lie group or Fr\'{e}chet--Lie group,
then $C^\infty_L(M,K)$ is
a Fr\'{e}chet--Lie group
and hence $C^\infty_c(M,K)$ an LF-Lie group.
But the construction of a Lie group structure works for general~$K$.
\begin{thm}
  \mlabel{thm:VII.2.4} {\rm(\cite{Gl07a})} For a $\sigma$-compact, non-compact manifold $M$ 
of positive dimension 
and a finite-dimensional Lie group $K$ of positive dimension, 
the Lie group $C^\infty_c(M,K)$ of compactly supported $K$-valued 
smooth functions, endowed with the Lie group structure modeled on the direct limit 
space $C^\infty_c(M,K) = \indlim C^\infty_{M_n}(M,K)$, $(M_n)_{n \in \N}$ an exhaustion of $M$, 
is not a direct limit of the subgroups $C^\infty_{M_n}(M,K)$ 
in the category of smooth manifolds, but a homomorphism 
$C^\infty_c(M,K) \to H$ to a Lie group $H$ is smooth if and only if it is smooth on 
each subgroup $C^\infty_{M_n}(M,K)$. 
\end{thm}
Again, a fragmentation map to a suitable box product
is behind the final positive assertion.
\begin{rem}
Uncountable direct systems (without cofinal subsequences)
are encountered when passing to uncountable index sets,
or paracompact finite-dimensional smooth manifolds.
Diffeomorphism groups can be constructed nonetheless (see \cite{Mr80}),
as well as test function groups (cf.\ \cite{Gl12b}, \cite{Gl15b}, \cite{Gl19a},
\cite{Gl21a}).
While direct limit properties in the category of all
Lie groups are unknown, the following can be shown
for a paracompact finite-dimensional smooth manifold~$M$
(see \cite{Gl21a}):\medskip

(a)
$\Diff_c(M)=\dl\, \Diff_L(M)$
in the category of \emph{regular}
Lie groups and smooth group homomorphisms
for $L$ in the directed set
of compact subsets of~$M$.\medskip

(b) $C^\infty_c(M,K)=\dl C^\infty_L(M,K)$
as a regular Lie group,
for each Lie group~$K$ which is $C^k$-regular for some $k\in \N_0$.
\end{rem}
\begin{rem}
We mention that some open problems from \cite{Ne06a}
concerning direct limits and related Lie groups
could be solved in the meantime:\medskip

(a) Locally finite-dimensional Lie groups 
are topological groups with Lie algebra in the sense of \cite{HoM04a}
(see \cite{Gl08a}).
Thus, the compact open topology
on ${\frak L}(G) \cong \Hom(\R,G)$ coincides with the given one on 
$\L(G)$.\medskip

(b) Every subgroup $H$ of a locally finite-dimensional Lie group $G$ 
is an initial Lie subgroup (see \cite{Gl08a}).\medskip

(c) For a Banach space~$E$ over~$\K$ and non-empty compact subset $K\sub E$,
there is a Lie group of germs of $\K$-analytic diffeomorphisms $\phi$ around $K$
such that $\phi|_K=\id_K$.
This Lie group is regular (see \cite{Da11}).
For $K$ a singleton in $\C^n$,
such Lie groups were first considered in \cite{Pis77}.\\[1mm]
Thus Problem VI.5, Problem VII.2 and Problem VII.3
from~\cite{Ne06a} were solved in the affirmative.
\end{rem}
\noindent
The following former problem has a negative answer
if $G$ is regular.\medskip

\noindent
%
\emph{Let $M$ be a locally convex manifold and 
$\g \subeq {\cal V}(M)$ be a countably dimensional locally finite-dimensional 
subalgebra consisting of complete vector fields. 
Does the inclusion $\g \to {\cal V}(M)$ integrate to a smooth action 
of a corresponding Lie group $G$ with $\L(G) = \g$}?\medskip
%
%

\noindent
In fact, consider the Fr\'{e}chet space $M=\R^\N$
and the abelian Lie group $(\R^\infty,+)$.
We have a left action
\[
\sigma\colon \GL_\infty(\R)\times \R^\N\to \R^\N,\quad
((a_{ij})_{i,j\in\N},(x_j)_{j\in \N})\mto
\Bigg(\sum_{j=1}^\infty a_{ij}x_j\Bigg)_{i\in \N}.
\]
For $(a_k)_{k\in\N}\in \R^\infty$,
consider the matrix $h((a_k)_{k\in \N})\in\GL_\infty(\R)$
with the \break $(i,j)$-entry $\delta_{i,j}$ if $i\geq 2$ or $j=1$, while $h((a_k)_{k\in\N})_{1,j}:=a_{j-1}$
for integers $j\geq 2$.
Then $h\colon \R^\infty\to\GL_\infty(\R)$
is a group homomorphism and thus
\[
\tau\colon \R^\infty\times\R^\N\to\R^\N,\quad (a,x)\mto \sigma(h(a),x)
\]
is a left action of $(\R^\infty, +)$ on $\R^\N$.
Note that $\tau$ is not continuous,
as the first component of $\tau(a,x)\in\R^\N$
is given by
\[
x_1+\sum_{k=1}^\infty a_kx_{k+1}
\]
and hence discontinuous (see Exercise~\ref{exc-product-not-k}).
For $n\in \N$, the restriction of $\tau$ to an action
\[
\tau_n\colon (\R^n,+)\times\R^\N\to\R^\N
\]
is smooth and hence gives rise to
a Lie algebra homomorphism
$\psi_n\colon \R^n\cong L(\R^n)\to \cV(\R^\N)$
using the negative of the Lie bracket of vector fields
on $\cV(\R^\N)$, which is determined
by
\[
\psi(v)(x):=\frac{d}{dt}\Big|_{t=0}
\tau(\exp_{\R^n}(tv),x)
=\frac{d}{dt}\Big|_{t=0}
\tau(tv,x)
=\sum_{k=1}^\infty v_kx_{k+1}
\]
and hence injective. We can therefore identify $\g=\bigcup_{n\in\N}\im(\psi_n)$
with $\R^\infty$.
The vector fields in $\g$ are complete because we have the $\R^n$-action.
Every regular, connected Lie group $G$
with Lie algebra $\g$ would be quotient of
$\R^\infty$ by a discrete subgroup
and hence act discontinuously as well.
%
%
%
\begin{rem}\mlabel{prob:VII.5}
We mention that the methods developed in \cite{Gl03a} and \cite{Gl05c} 
for the analysis of direct limit 
Lie groups have potential to apply to more general classes of 
Lie groups $G$ which are direct limits of finite-dimensional manifolds 
$(M_n)_{n \in \N}$, 
with the property that for $n,m \in \N$ there exist $c(n,m)$ and $d(n)$ with 
$$ M_n \cdot M_m \subeq M_{c(n,m)} \quad \hbox{ and } \quad M_n^{-1} \subeq M_{d(n)}, $$
a situation which occurs in free constructions.
See \cite[Lemma~14.5]{Gl07a} for one result
of this form.
Similar situations, with infinite-dimensional $M$, occur in the ind-variety 
description of Kac--Moody groups (cf.\ \cite{Kum02}, \cite{BiPi02}). 
\end{rem}
\section{Open Problems for Chapter~\ref{ch:dirlim}}

\begin{probl} \mlabel{prob:VII.1} Is every Lie group $G$ regular whose Lie algebra 
$\g$ is countably dimensional, 
locally finite, and endowed with the finest locally convex 
topology?
\end{probl}

\section{Notes and comments on Chapter~\ref{ch:dirlim}} 
The systematic study of Lie group structures on direct limit 
Lie groups $G = \indlim G_n$ was started in the 1990s by J.~Wolf  
and his coauthors \cite{NRW91, NRW93}. 
They used certain conditions on the groups $G_n$ and the connecting 
maps $G_n \to G_{n+1}$ to ensure that the direct limit 
groups are locally exponential.  
Since not all direct limit groups are locally exponential 
(Example~\ref{ex:VII.1.4}(c)), 
their approach does not cover all cases. The picture for countable 
direct limits of finite-dimensional Lie groups was completed in 
\cite{Gl03b, Gl05c}, where it
was shown that 
arbitrary countable direct limits of finite-dimensional 
Lie groups exist. The key to these results were general 
construction principles for direct limits of finite-dimensional manifolds,
as presented here in Section~\ref{dl-construction}. 
The approach of Wolf et al.\
was later adapted to ascending unions of Banach--Lie groups
in well-behaved situations
and regularity established for such (see \cite{Da14}).
In \cite[Thm.~47.9]{KM97}, 
it is shown that every subalgebra $\h$ of $\gl_\infty(\R)$ is integrable 
to an integral subgroup, 
which is a special case of Theorem~\ref{thm:VII.1.3}. Here $\h$ is even BCH. 
For $\dim K < \infty$, the groups $C^\infty_c(M,K)$ are also discussed in 
\cite{NRW94} as direct limit Lie groups which are BCH;
they also occur in \cite{AHM93}.


\chapter{Groups of Diffeomorphisms} \mlabel{ch:diffeo} 

%
%
\red{This chapter is still unfinished. Some proofs remain
  to be added and some sections have to be written.
Below we marked in red where pieces are missing.} \\

\section{Diffeomorphism groups of compact manifolds}
%

%
If $K$ is a compact smooth manifold (without boundary)
and $\ell\in\N_0\cup\{\infty\}$,
we write $\Diff_{C^\ell}(K)$ for the group
of all $C^\ell$-diffeomorphisms $\gamma \colon K\to K$
(with composition as the group operation and neutral element $\id_K$).
We abbreviate $\Diff(K):=\Diff_{C^\infty}(K)$.
As an application of the techniques developed, we shall
establish the following theorem:
\begin{thm}\label{diffKsmooth}
If $K$ is a compact smooth manifold without boundary,
then $\Diff(K)$ is an open subset of $C^\infty(K,K)$
and the smooth manifold structure as an open submanifold
turns $\Diff(K)$ into a Lie group.
\end{thm}
This theorem is subsumed by Proposition~\ref{detail-diff}
below, which provides additional information.
To prepare the result, et us
first record differentiability properties of the composition map.
\begin{prop}\label{compallmaps}
Let $K$ be a compact smooth manifold $($possibly with rough boundary$)$,
$L$ be a compact smooth manifold $($without boundary$)$,
$M$ be a smooth manifold admitting a local addition, and $k,\ell\in \N_0\cup\{\infty\}$.
Then the following holds:
\begin{description}[(D)]
\item[\rm(a)] The mapping
\[
\Gamma \colon C^{\ell+k}(L,M)\times C^\ell(K,L)\to C^\ell(K,M),\;\, (\gamma,\eta)
\mto\gamma\circ\eta
\]
is $C^k$-map.
%
%
Notably, the map
$C^\infty(L,M)\times C^\infty(K,L)\to C^\infty(K,M)$, $(\gamma,\eta)\mto\gamma\circ\eta$
is smooth.
\item[\rm(b)]
For each $\gamma\in C^\ell(K,L)$, the map $\rho_\gamma\colon C^\ell(L,M)\to C^\ell(K,M)$, $\eta\mto\eta\circ\gamma$
is smooth.
\end{description}
\end{prop}
\begin{prf}
(a) We have $\Gamma^\wedge(\gamma,\eta,x)=\gamma(\eta(x))=\ve(\gamma,\wt{\ve}(\eta,x)$
in terms of the evaluation maps $\ve\colon C^{\ell+k}(L,M)\times L\to M$
and $\wt{\ve}\colon C^\ell(K,L)\times K\to L$.
Since $\wt{\ve}$ is $C^{\infty,\ell}$ (and thus $C^{k,\ell}$)
and the projection $C^{\ell+k}(L,M)\times C^\ell(K,L)\to C^\ell(K,L)$
onto the second component is~$C^k$, the map
\[ h \colon (C^{\ell+k}(L,M)\times C^\ell(K,L))\times K\to L, \quad
  (\gamma,\eta,x)\mto\eta(x) \] 
is~$C^{k,\ell}$, by the Chain Rule~(Proposition~\ref{chainC1}).
Now
\[ g\colon (C^{\ell+k}(L,M)\times C^\ell(K,L))\times K\to C^{\ell+l}(L,M),
  \quad (\gamma,\eta,x)\mto x \] 
is a smooth map and hence $C^{k,\ell}$ as well.
Since $\ve$ is $C^{\infty,l+\ell}$ and thus $C^{k+\ell}$,
the Chain Rule~(Proposition~\ref{chainC1})
shows that $\Gamma=\ve\circ (g,h)$ is $C^{k,\ell}$.

(b) \red{Sketch:} The idea is that
\begin{align*}
  (\phi_{\eta\circ\gamma}\circ \rho_\gamma\circ \phi_\eta^{-1})(\sigma)
  &=\Sigma\circ (\theta^{-1}\circ (\eta,\sigma)) \circ \gamma\\
&  =\pr_2\circ\theta\circ\theta^{-1}\circ (\eta,\sigma)\circ\gamma=\sigma\circ\gamma
\end{align*}
but first must show that the leftmost term is defined everywhere.
We need that
$\rho_\gamma\colon \Gamma \eta^*(E)\to \Gamma((\eta\circ\gamma)^*(E))$ is continuous linear.
\end{prf}
The following five lemmas prepare the discussion of diffeomorphism groups.
The first one is devoted to immersions.
\begin{defn}
Let $\ell\in\N\cup\{\infty\}$.
A $C^\ell$-map $f\colon M\to N$ from a finite-dimensional $C^\ell$-manifold~$N$
to a $C^\ell$-manifold~$N$ (both without boundary) is called
an \emph{immersion} (or $C^\ell$-immersion)
if $T_xf\colon T_xM\to T_xN$ is injective for all $x\in M$.
We write $\Imm_{C^\ell}(M,N)$ for the set of all $C^\ell$-immersions
$\gamma\colon M\to N$, and abreviate $\Imm(M,N):=\Imm_{C^\infty}(M,N)$.
\end{defn}
\begin{rem}
The reader is referred to \red{(suitable exercises)} for
a discussion of immersions, embeddings and submersions between general
$C^\ell$-manifolds (modeled on arbitrary locally convex spaces).
\end{rem}
\begin{lem}\label{immlocinjsubset}
Let $n\in\N$ and $\gamma\colon V\to\R^n$ be a $C^1$-map
on an open subset $V\sub\R^n$ such that $\gamma'(x_0)=\id_{\R^n}$
for some $x_0\in V$.
Then there exists a \break $\gamma$-neighborhood $P \sub C^1(V,F)$
and an open $x_0$-neighborhood $Q\sub V$
such that $\eta|_Q$ is an injective $C^1$-immersion for each $\eta\in P$.
\end{lem}
\begin{prf}
Note that $\gamma'\colon V\to\cL(\R^n)$ is continuous
as $\gamma'(x)(y)=d\gamma(x,y)\in\R^n$ depends continuously on $x\in V$
for each $y\in\R^n$, and the unique vector topology on $\cL(\R^n)\cong \R^{n^2}$
coincides with the topology of pointwise convergence.
We fix a norm~$\|\cdot\|$ on~$\R^n$ and write $\|\cdot\|_{\op}$
for the corresponding operator norm on $\cL(\R^n)$. Also balls shall refer to the given norm.
Using the continuity of $\gamma'$, we find $r>0$ with $K:=\wb{B}_r(x_0)\sub U$
such that
\begin{equation}\label{estigamm}
\sup_{x\in \wb{B}_r(x_0)}\|\gamma'(x)-\id_{\R^n}\|_{\op}<\frac{1}{3}.
\end{equation}
We set $Q:=B_r(0)$ and let $P$ be the set of all $\eta\in C^1(V,F)$
such that
\begin{equation}\label{estiet2}
\sup_{x\in \wb{B}_r(x_0)}\|\eta'(x)-\id_{\R^n}\|_{\op}<\frac{1}{3},
\end{equation}
which entails that $\eta'(x)\in \GL(\R^n)$ for each $x\in\wb{B}_r(x_0)$.
In particular, $\eta|_Q$ is an immersion.
Moreover,
\[ \|\eta'(x_0)^{-1}\|_{\op}=\Bigg\|\sum_{k=0}^\infty(\id_{\R^n}-\eta'(x))^k\Bigg\|_{\op}
  <\frac{1}{1-\frac{1}{3}}=3/2\]
using Neumann's series and the geometric series, whence
$\frac{1}{\|\eta'(x_0)^{-1}\|_{\op}}> 2/3.$ 
Since $\|\eta'(x)-\eta'(x_0)\|_{\op}\leq\|\eta'(x)-\id_{\R^n}\|_{\op}+\|\eta'(x_0)-\id_{\R^n}\|_{\op}
<\frac{3}{3}$ we have $\Lip(\eta-\eta'(x_0))\leq 2/3 <\frac{1}{\|\eta'(x_0)^{-1}\|_{\op}}$,
whence $\eta|_{B_r(0)}$ is injective, Theorem~\ref{lip-inv-fct}.
It only remains to note that~$P$ is an open $\gamma$-neighborhood in $C^1(V,\R^n)$
since
\[
D\colon C^1(V,\R^n)\to C(V,\cL(\R^n)), \;\, \eta\mto\eta'
\]
is a continuous linear map (as $D\eta=(d\eta)^\vee$
and Lemma~\ref{df-vs-Tf}, Proposition~\ref{ctsexp},
as well as Lemma~\ref{ctsemb} apply); moreover,
\[ C(V,\cL(\R^n))\to [0,\infty[, \quad
  \zeta\mto \sup_{x\in K}\|\zeta(x)\|_{\op}\] 
is a continuous seminorm on $C(V,\cL(\R^n)$ (cf.\ Lemma~\ref{sammelsu}(b)).
\end{prf}
\begin{lem}\label{immlocinj}
Let $M$ be a finite-dimensional $C^1$-manifold,
$\gamma\colon M\to N$ be a $C^1$-immersion
to a manifold~$N$, and $x_0\in M$.
Then there exists a \break $\gamma$-neighborhood $W \sub C^1(M,N)$
and an open $x_0$-neighborhood $U\sub M$
such that $\eta|_U$ is an injective $C^1$-immersion for each $\eta\in W$.
\end{lem}
\begin{prf}
Let $F$ be the modeling space of~$N$ and
$\psi\colon U_\psi\to V_\psi\sub F$ be a chart of~$N$ around~$\gamma(x)$.
Let $n$ be the dimension of~$M$ and 
$\phi\colon U_\phi\to V_\phi\sub\R^n=:E$
be a chart of~$M$ around~$x_0$ such that $\gamma(U_\phi)\sub U_\psi$ and $\phi(x_0)=0$.
Then $f:=\psi\circ \gamma| \circ \phi^{-1}\colon V_\phi\to F$ is~$C^1$
and the vectors $y_j:=dh(x_0,e_j)$ are linearly independent in~$F$ for $j\in\{1,\ldots,n\}$.
Set $Y:=\Spann\{y_1,\ldots,y_n\}$ and let $y_1^*,\ldots,y_n^*\in Y'$ be the dual basis
determined by $y_j^*(y_k)=\delta_{j,k}$ for all $j,k\in\{1,\ldots, n\}$.
By the Hahn-Banach Extension Theorem, for each $j\in\{1,\ldots, n\}$
there exists $\lambda_j\in F'$ such that $\lambda_j|_Y=y_j^*$.
Let $K\sub U_\phi$ be a compact $x_0$-neighborhood.
Then
\[
\Omega:=\{\eta\in C^1(M,N)\colon \eta(K)\sub U_\psi\}=C^1(M,N)\cap\lfloor K,U_\psi\rfloor
\]
is an open $\gamma$-neighborhood in $C^1(M,N)$. Then $V:=\phi(K^0)$
is an open subset of~$\R^n$ and
\[
\kappa
\colon \colon \Omega\to C^1(V,\R^n),\quad \eta\mto f\circ \psi\circ \eta \circ \phi^{-1}|_V
\]
is continuous, by Lemma~\ref{Cktoppu}(b),
Remark~\ref{earlyremCk}, and Lemma~\ref{Cktoppu}(a).
By Lemma~\ref{immlocinjsubset}, there exists an open $\phi(x_0)$-neighborhood
$Q\sub V$ and a $\kappa(\gamma)$-neighborhood $P\sub C^1(V,\R^n)$
such that $\eta|_Q$ is an injective immersion for each $\eta\in P$.
Then $U:=\phi^{-1}(Q)$ is an open $x$-neighborhood in~$M$
and $W:=\kappa^{-1}(P)$ is a $\gamma$-neighborhood in $C^1(M,N)$
such that $f\circ\psi\circ \eta|_U$ (and hence also $\eta|_U$)
is an injective immersion for each $\eta\in W$.
\end{prf}
\begin{lem}
If $K$ is a compact smooth manifold,
$M$ a smooth manifold admitting a local addition $($both without boundary$)$
and $\ell\in \N\cup\{\infty\}$, then $\Imm_{C^\ell}(K,M)$
is $C^1$-open in $C^\ell(K,M)$.
\end{lem}
\begin{prf}
Since $\Imm_{C^\ell}(K,M)=\Imm_{C^1}(K,M)\cap C^\ell(K,M)$,
we may assume that $\ell=1$. Let $\gamma\in \Imm_{C^1}(K,M)$.
For each $x\in M$, Lemma~\ref{immlocinj} provides an open $\gamma$-neighborhood $W_x\sub
C^1(K,M)$ and an open $x$-neighborhood $U_x\sub K$
such that $\eta|_{U_x}$ is an immersion for each $\eta\in W_x$.
Since~$K$ is compact, we have $K=\bigcup_{x\in\Phi}U_x$ for
a finite subset $\Phi\sub K$. Then $W:=\bigcap_{x\in\Phi}W_x$ is an open
$\gamma$-neighborhood such that $W\sub\Imm_{C^1}(K,N)$.
\end{prf}
\begin{lem}\label{global-invpar}
Let $P$ be a $C^k$-manifold $($possibly with rough boundary$)$
modeled on a locally convex space~$E$
and $M$ be a $C^k$-manifold modeled on a Banach space~$F$,
where $k\in \N\cup\{\infty\}$.
Let $f\colon P\times M\to M$ be a $C^k$-map such that
$f_p:=f(p,\cdot)\colon M\to M$ is a $C^k$-diffeomorphism for
each $p\in P$. If $F$ is infinite-dimensional, assume that $k\geq 2$. Then also
\[
g\colon P\times M\to M,\quad (p,y)\mto :=g_p^{-1}(y)
\]
is a $C^k$-map.
\end{lem}
\begin{prf}
Given $(p,y)\in P\times M$, let $\kappa\colon U_\kappa\to V_\kappa\sub F$
be a chart of~$M$ around~$y$. There exist a chart $\phi\colon U_\phi\to V_\phi\sub E$
of~$P$ around~$p$ and a chart
\[ \psi\colon U_\psi\to V_\psi\sub F \]  of~$M$ around
$x:=f_p^{-1}(y)$ such that $f(U_\phi\times U_\psi)\sub U_\kappa$.
Then
\[ \kappa\circ (\phi^{-1}\times \psi^{-1})\colon
V_\phi\times V_\psi\to V_\kappa \]  satisfies the hypotheses
of \red{(reference/theorem to be added)}. As a consequence, there exists an open $p$-neighborhood $Q\sub U_\phi\sub P$
and an open neighborhood $W\sub U_\kappa\sub M$ such that $g|_{Q\times W}$
is~$C^k$.
\end{prf}
\begin{lem}\label{C1metriza}
If $K$ is a compact smooth manifold,
then $C^1(K,K)$ is metrizable.
\end{lem}
\begin{prf}
By Whitney's Embedding Theorem \cite{Wh36},
$K$ embeds in $\R^d$ for some $d\in\N$. 
Hence $TK$ embeds in $\R^{2d}$. Thus $TK$
is a $\sigma$-compact locally compact space and metrizable,
whence $C(TK,\R^{2d})$ is metrizable (by Lemma~\ref{sammelsu}(c))
and hence also $C(TK,TK)$ (using Lemma~\ref{ctsemb}).
As a consequence, also $C^1(K,K)$ is metrizable, as it embeds in $C(TK,TK)$ (see Lemma~\ref{onlyT}).
\end{prf}
\begin{prop}\label{detail-diff}
  Let $\ell\in\N\cup\{\infty\}$ and $K$ be a compact smooth manifold.
Then
$\Diff_{C^\ell}(K)$ is $C^1$-open in $C^\ell(K,K)$. We endow $\Diff_{C^\ell}(K)$ with the
smooth manifold structure as an open subset of $C^\ell(K,K)$.
Then the following holds for all $k\in \N_0\cup\{\infty\}$:
\begin{description}[(D)]
\item[\rm(a)]
The inversion map $\Diff_{C^{\ell+k}}(K)\to \Diff_{C^\ell}(K)$ is $C^k$;
\item[\rm(b)]
The composition map
$\Diff_{C^{\ell+k}}(K)\times \Diff_{C^\ell}(K) \to \Diff_{C^\ell}(K)$ is $C^k$;
\item[\rm(c)]
  For each $\gamma\in \Diff_{C^\ell}(K)$, the right translation
  \[ \rho_\gamma\colon\Diff_{C^\ell}(K)\to\Diff_{C^\ell}(K),
\quad \eta\mto\eta\circ\gamma \]  is a $C^\infty$-map;
\item[\rm(d)]
$\Diff_{C^\ell}(K)$ is a topological group; and
\item[\rm(e)]
$\Diff(K)$ is a Lie group.
\end{description}
\end{prop}
\begin{prf}
If $\Diff_{C^\ell}(K)$ was not $C^1$-open in $C^\ell(K,K)$,
we could find $\gamma\in\Diff(K)$
and a sequence $(\gamma_n)_{n\in \N}$ in $C^\ell(K,K)\setminus \Diff_{C^\ell}(K)$
which converges to~$\gamma$ in the $C^1$-topology
(cf.\ Lemma~\ref{C1metriza}).
For each $x\in M$, Lemma~\ref{immlocinj} provides an open $\gamma$-neighborhood $W_x\sub
C^1(K,M)$ and an open $x$-neighborhood $U_x\sub K$
such that $\eta|_{U_x}$ is an injective immersion for each $\eta\in W_x$.
Since~$K$ is compact, we have $K=\bigcup_{x\in\Phi}U_x$ for
a finite subset $\Phi\sub K$. Then $W:=\bigcap_{x\in\Phi}W_x$ is an open
$\gamma$-neighborhood such that $W\sub\Imm_{C^1}(K,N)$.

Given $x\in K$, let $V(x)$ be the connected component of~$K$
containing~$x$ (which is both open and compact). Pick $x_1,\ldots, x_m\in K$ such that $K=\bigcup_{j=1}^m V(x_j)$.
There exists $n_0\in\N$ such that
\[
\gamma_n\in \lfloor \{x_j\},V(\gamma(x_j))\rfloor\;\,\mbox{for all $j\in\{1,\ldots,m\}$,}
\]
for all $n\geq n_0$. Thus $\gamma_n(x_j)\in V(\gamma(x_j))$ and thus $\gamma_n(V(x_j))\sub V(\gamma(x_j))$.
Since $\gamma_n$ is \'{e}tale and hence an open map, we deduce that
$\gamma_n(V(x_j))$ is both open and compact. Thus $\gamma_n(V(x_j))=V(\gamma(x))$,
by connectedness, entailing that $\gamma_n(K)=K$. After omitting $\gamma_1,\ldots,\gamma_{n_0-1}$,
we may assume that each $\gamma_n$ is surjective. Now, not being a diffeomorphism,
$\gamma_n$ must not be injective. We therefore find $x_n\not=y_n$ in~$K$ such that $\gamma_n(x_n)=\gamma_n(y_n)$.
Since~$K$ is compact, after passing to subsequences we may assume that both $(x_n)_{n\in\N}$
and $(y_n)_{n\in\N}$ are convergent sequences, with limits $x$ and $y$, say.
The evaluation map $\ve\colon C(K,K)\times K\to K$ being continuous,
we deduce that
\[
\gamma(x)=\lim_{n\to\infty}\ve(\gamma_n,x)=\lim_{n\to\infty}\ve(\gamma_n,y)=\gamma(y).
\]
Hence $x=y$, as $\gamma$ is a diffeomorphism. Now $x\in U_z$ for some $z\in\Phi$.
Since $x_n\to x$ and $y_n\to x$, there exists $N\in\N$ such that $x_n,y_n\in U_z$
for all $n\geq N$. After increasing~$N$, we may assume that, moreover,
$\gamma_n\in W$ for all $n\geq N$. As $\gamma_N\in W\sub W_z$,
the restriction $\eta_N|_{W_z}$ is injective. Since $x_N,y_N\in U_z$, we now
deduce from
$\gamma_N(x_N)=\gamma_N(y_N)$ that $x_N=y_N$, contradicting the choice of these
elements. Thus $\Diff_{C^\ell}(K)$ is $C^1$-open in $C^\ell(K,K)$.

(a) To see that the inversion map $\eta\colon \colon\Diff_{C^{\ell+k}}(K)\to\Diff_{C^\ell}(K)$
is~$C^k$, note that the evaluation map
\[
f \colon \Diff_{C^{k+l}}(K)\times K\to K,\;\, (\gamma,x)\mto \gamma(x)
\]
is $C^{k+\ell}$ (see Theorem~\ref{thmmfdmps}(c)) and $f(\gamma,\cdot)=\gamma$ is a diffeomorphism
for each $\gamma\in\Diff_{C^{\ell+k}}(K)$. By Lemma~\ref{global-invpar},
the map
\[
g\colon \Diff_{C^{k+\ell}}(K)\times K\to K;\quad (\gamma,y)\mto g(\gamma,\cdot)^{-1}(y)=\gamma^{-1}(y)
\]
is~$C^{k+\ell}$ and thus $C^{k,\ell}$. As $\eta^\wedge=g$,
the exponential law (as in Theorem~\ref{thmmfdmps}) entails that $\eta$ is~$C^k$.
 
(b) and (c) follow from Proposition~\ref{compallmaps}.

(d) follows from~(a) and (b), applied with $k=0$.

(e) follows from~(a) and~(b), applied with $k=\ell=\infty$.
\end{prf}
\begin{rem}
We mention that the compact-open topology makes
\[ \Diff_{C^0}(K)=\Homeo(K) \]
a topological group,
for each compact smooth manifold~$K$ (and more generally:
see Exercise~\ref{exer:I.5}).
\end{rem}
identify Lie algebra

prove regularity
\section{Diffeomorphism groups of non-compact manifolds}
\red{This section has to be written. It should cover the following
topics: 
\begin{enumerate}
\item[\rm(a)] $\Diff_c(M)$ for $M$ non-compact and semidirect products 
$\Diff_c(M) \rtimes G$, where $G$ acts smoothly on $M$. 
\item[\rm(b)] Higher order jet groups as extensions of 
$\Diff(M)$ by gauge groups of higher order tangent bundles 
(see the appendix in \cite{BiNe08}). 
\end{enumerate}}

For every locally compact space $X$, the group 
$\Homeo(X)$ of homeomorphisms of $X$ carries a natural group structure, 
which is determined by the property that the map 
\[ \Homeo(X) \into C(X,X)^2, \quad 
\phi \mapsto (\phi, \phi^{-1}) \] 
is a topological embedding (Exercise~\ref{exer:I.5}). 

If $M$ is a finite-dimensional smooth manifold, then all iterated tangent 
bundles $T^k(M)$ are locally compact, so that we obtain a natural homomorphic 
embedding 
\[ \Diff(M) \into \prod_{k \in \N_0} \Homeo(T^k(M)), \quad 
\phi \mapsto (T^k(\phi)).\] 
This leads to a natural group topology on $\Diff(M)$, but in general one 
cannot hope for this topology to be compatible with a Lie group 
structure. 

\begin{ex} If $M = \coprod_{k \in \N} \bS^1$ is a countable union of circles, 
then the arc-component of the identity if $\Diff(M)$ is the group 
$\Diff_+(\bS^1)^\N$ which is not locally contractible, hence does not carry 
a manifold structure. 
\end{ex}

\section{Diff($\bS^1$) is not locally exponential} 

Below we show that the exponential function
$$ \exp \: {\cal V}(\bS^1) \to \Diff(\bS^1)$$ 
is not a local diffeomorphism by proving that every identity neighborhood of 
$\Diff(\bS^1)$ contains elements which do not lie on any one-parameter group. 

Let $G := \Diff_+(\bS^1)$ denote the group of orientation preserving diffeomorphisms of 
$\bS^1 \cong \R/2pi \Z$, i.e., the identity component of $\Diff(\bS^1)$. 
To get a better picture of this group, we first construct 
its universal covering group $\tilde G$. Let 
$$ \tilde G := \{ \phi \in \Diff(\R) \: (\forall x \in \R)\ \phi(x + 2\pi) 
= \phi(x) + 2 \pi, \phi' > 0\}. $$
We consider the map 
$$ q \: \R \to \bS^1 := \R/2\pi \Z, \quad x \mapsto x + 2 \pi \Z $$
as the universal covering map of $\bS^1$. Then every orientation preserving 
diffeomorphism $\psi \in \Diff_+(\bS^1)$ 
lifts to a diffeomorphism $\tilde\psi$ of $\R$, commuting with the translation 
action of the group $2\pi \Z \cong \pi_1(\bS^1)$, which means that 
$\tilde\psi(x + 2\pi) = \tilde\psi(x) + 2\pi$ for each $x \in \R$. 
The diffeomorphism $\tilde\psi$ is uniquely determined by the choice of an 
element in $q^{-1}(\psi(q(0)))$. That $\psi$ is orientation preserving means that 
$\tilde\psi' > 0$. Hence we have a surjective homomorphism 
$$ q_G  \: \tilde G \to G, \quad q_G(\phi)(q(x)) := q(\phi(x)) $$
with kernel isomorphic to $\Z$. 

The Lie group structure of $\tilde G$ is rather simple. It can be defined by a 
global chart. Let $C^\infty_{2\pi}(\R,\R)$ denote the Fr\'echet space of 
$2\pi$-periodic smooth functions on $\R$, which is considered as a closed subspace 
of the Fr\'echet space $C^\infty(\R,\R)$. In this space,  
$$ U := \{ \phi \in C^\infty_{2\pi}(\R,\R) \: \phi' > -1\} $$
is an open convex subset, and the map 
\[  \Phi \: U \to \tilde G, \quad \Phi(f) = \id_\R + f\] 
is a bijection. 
In fact, let $f \in U$. Then $\Phi(f)(x + 2\pi) = \Phi(f)(x) + 2\pi$ 
follows directly from the requirement that $f$ is $2\pi$-periodic, and 
$\Phi(f)' > 0$ follows from $f' > -1$. Therefore $\Phi(f)$ is strictly increasing, 
hence a diffeomorphism of $\R$ onto the interval $\Phi(f)(\R)$. As the latter 
interval is invariant under translation by $2\pi$, $\Phi(f)$ is 
surjective and therefore $\Phi(f) \in \tilde G$. 
Clearly, $\Phi^{-1}(\psi) = \psi - \id_\R$ maps $\tilde G$ to $U$. 
We define the manifold structure on $\tilde G$ by declaring 
$\Phi$ to be a global chart. With respect to this chart, the group operations 
in $\tilde G\cong U$ are given by 
$$ m(f,g)(x) := f(g(x) + x) - x 
\quad \hbox{ and } \quad  \eta(f)(x) = (f + \id_\R)^{-1}(x) - x, $$
which can be shown directly to be smooth maps. We thus obtain on 
$\tilde G$ the structure of a Lie group such that 
$\Phi \: U \to \tilde G$ is a diffeomorphism. In particular, $\tilde G$ 
is contractible and therefore simply connected, so that the map 
$q_G \:\tilde G \to G$ turns out to be the universal covering map of $G$. 

\msk 

\begin{thm}
  \mlabel{thm:V.1.7}  Every identity neighborhood in $\Diff(\bS^1)$ contains elements 
not contained in the image of the exponential function. 
\end{thm}

\begin{prf} First we construct certain elements in $\tilde G$ which are close to the identity. 
For $0 < \eps < {1\over n}$, we consider the function 
$$ f \: \R \to \R, \quad x \mapsto x + {\pi \over n} + \eps \sin^2(n x) $$
and observe that $f \in \tilde G$ follows from 
\[ f'(x) = 1 + 2 \eps n \sin(nx) \cos(nx) = 1 + \eps n\sin(2nx)> 0.\]

\nin {\bf Step 1.} For $n$ large fixed and $\eps \to 0$, we get elements in $\tilde G$ which are 
arbitrarily close to $\id_\R$.

\nin {\bf Step 2.} $q_G(f)$ has a unique periodic orbit of order $2n$ on $\bS^1$: 
Under $q_G(f)$, the point $q(0) \in \bS^1$ is mapped to 
${\pi \over n}$ etc., so that we obtain the orbit 
$$ q(0) \to q({\textstyle{\pi \over n}}) \to q({\textstyle{2\pi \over n}}) 
\to \ldots \to q({\textstyle{(2n-1)\pi\over n}}) 
\to q(0). $$

For $0 < x_0 < {\pi \over n}$,  we have for $x_1 := f(x_0)$: 
$$ x_0 + {\pi \over n} < x_1 < {2\pi \over n}, $$
and for $x_n := f(x_{n-1})$, the relations 
$$ 0 < x_0 < x_1 -{{\textstyle\pi \over n}} <  x_2 -{{\textstyle 2\pi \over n}} 
<  \cdots < {{\textstyle \pi \over n}}. $$
Therefore $x_k - x_0 \not\in 2 \pi \Z$ for each $k \in \N$, and hence the orbit of 
$q(x_0)$ under $q_G(f)$ is not finite. This proves that $q_G(f)$ has a unique periodic 
orbit and that the order of this orbit is $2n$.

\nin {\bf Step 3.} 
$q_G(f) \not= g^2$ for all $g \in \Diff(\bS^1)$: We analyze the 
periodic orbits. Every periodic point of $g$ is a periodic point of $g^2$ and vice versa. 
If the period of $x$ under $g$ is odd, then the period of $x$ under 
$g$ and $g^2$ is the same. If the period of $x$ under $g$ 
is $2m$, then its $g$-orbit under $g$ breaks 
up into two orbits under $g^2$, each of order $m$. Therefore $g^2$ can never have a single 
periodic orbit of even order, and this proves that $q_G(f)$ has no square root in 
$\Diff(\bS^1)$. It follows in particular that $q_G(f)$ does not lie on any one-parameter 
subgroup, i.e., $q_G(f) \not = \exp X$ for each $X \in {\cal V}(M)$. 
\end{prf} 

\begin{rem}
  \mlabel{rem:V.1.8} (a) If $M$ is a compact manifold, then one can show that the identity 
component $\Diff(M)_0$ of $\Diff(M)$ is a simple group 
(Epstein, Hermann and Thurston; see \cite{Ep70}). 
Being normal in $\Diff(M)_0$, the subgroup $\la \exp {\cal V}(M) \ra$ 
coincides with $\Diff(M)_0$. Hence every 
diffeomorphism homotopic to the identity is a finite product of 
exponentials. 

  (b) Although $\Diff(M)_0$ is a simple Lie group, its Lie algebra 
${\cal V}(M)$ is far from being simple. For each open subset 
$U \subeq M$, the set ${\cal V}(M)_U$ of all vector fields vanishing in $U$ 
is an ideal which is proper if $U$ is not dense. 
\end{rem}

\begin{ex}
Identifying the Lie algebra $\g := {\cal V}(\bS^1)$ of $\Diff(\bS^1)$ 
with smooth $2\pi$-periodic functions on~$\R$, the Lie bracket corresponds to 
$$ [f,g] = fg'-f'g. $$
For the constant function $c_0 = 1$, 
$c_n(t) := \cos(nt)$ and $s_n(t) = \sin(nt)$, this leads to 
$$ [c_0, s_n] = n c_n \quad \hbox{ and }  \quad [c_0, c_n] = -n s_n, $$
so that $\Spann\{c_0, s_n, c_n\} \leq {\cal V}(\bS^1)$ is a Lie subalgebra 
isomorphic to $\fsl_2(\R)$. It further follows that 
$((\ad c_0)^2 + n^2\be) s_n = 0$, so that 
$\kappa_\g({2\pi\over n}c_0)s_n = 0$ implies that $\exp_G$ is 
not injective in any neighborhood of ${2\pi\over n}c_0$ 
(Proposition~\ref{prop:exp-noninj}) (cf.\ \cite{Mil82}, Ex. 6.6). 
Therefore $\exp_{\Diff(\bS^1)}$ is neither locally surjective nor injective. 
\end{ex}

\section{Lie transformation groups} 

One of the fundamental references on topological transformation groups 
is the monograph \cite{MZ55} by Montgomery and Zippin. 
Since it also deals with differentiability properties of transformation groups on manifolds, 
some of the techniques described there have interesting applications 
in the context of infinite-dimensional Lie theory. 

\subsection{Smooth Lie group actions} 

\begin{thm}
  \mlabel{thm:IX.1.1} {\rm(\cite{BoMo45}, \cite[p.212]{MZ55})} 
Any continuous action \break {$\sigma\: G \times M \to M$}
of a finite-dimensional Lie group on a finite-dimensional 
smooth manifold $M$ by diffeomorphisms 
is smooth. 
\end{thm}

For compact manifolds we obtain the following ``automatic smoothness'' result 
on homomorphisms of Lie groups (see also \cite{CM70} for one-parameter groups; and 
\cite{Gl08b} for the non-compact case). 

\begin{cor}
  \mlabel{cor:IX.1.2} If $M$ is a $\sigma$-compact finite-dimensional 
manifold and $G$ a finite-di\-men\-sio\-nal Lie group, 
then every continuous homomorphism \break $\phi \: G \to \Diff_c(M)$ is smooth. 
\end{cor}

The following result answers question (FP9) in \cite{Ne06a},
asking for which Lie groups all locally compact subgroups are
finite-dimensional   Lie groups,
in the affirmative  for diffeomorphism groups. 

\begin{thm} \mlabel{thm:IX.1.3} {\rm(\cite[Th.~5.2.2, p.~208]{MZ55})} If 
a locally compact group $G$ acts faithfully on a smooth finite-dimensional manifold $M$ 
by diffeomorphisms, then $G$ is a finite-dimensional Lie group. 
If $M$ is compact, then each locally compact subgroup of 
$\Diff(M)$ is a Lie group. 
\end{thm}

The preceding results take care of the actions of locally compact groups 
on manifolds. As the work of {de la Harpe} and { Omori} (\cite{OdH71, OdH72}) 
shows, the situation for Banach--Lie groups is more subtle. 
If $\alpha \: \g \to {\cal V}(M)$ is an injective homomorphism, then for each 
$p \in M$ the subspace 
$$\g_p := \{ x \in \g \: \alpha(x)(p) =0\}$$ 
is a finite-codimensional subalgebra 
with $\bigcap_p \g_p = \{0\}$. Therefore the existence of many finite-codimensional 
subalgebras is necessary for Lie algebras to be realizable by vector fields on 
a finite-dimensional manifold. 

\begin{thm}
  \mlabel{thm:IX.1.4} {\rm(\cite{OdH72})} Let $G$ be a Banach--Lie group. 
If $\L(G)$ has no proper finite-codimensional closed ideals, then 
$\L(G)$ has no proper finite-codimensional closed subalgebra and each smooth action  
of $G$ on a finite-dimensional manifold is trivial. 
\end{thm}

\begin{thm}
  \mlabel{thm:IX.1.5} {\rm(\cite{OdH72})} If a Banach--Lie group $G$ acts smoothly, effectively, 
amply (for each $m \in M$ the evaluation map $\g \to T_m(M)$ is surjective),  
and primitively (it leaves no foliation invariant) on a 
finite-dimensional manifold $M$, then it is finite-dimensional. 
\end{thm}

\begin{thm}
  \mlabel{thm:IX.1.6}  {\rm(\cite[Thms.~B/C]{Omo78})} 
Let $G$ be a connected Banach--Lie group with Lie algebra $\g$ 
acting smoothly, effectively 
and transitively on a finite-dimensional manifold $M$. 
\begin{description}
\item[\rm(1)] If $M$ is compact, then $G$  is finite-dimensional. 
\item[\rm(2)] If $M$ is non-compact, then $\g$ contains a finite-codimensional 
closed solvable ideal. 
\end{description}
\end{thm}

Since $\Diff(M)$ acts smoothly, effectively and transitively on $M$, this implies: 

\begin{cor}
  \mlabel{cor:IX.1.7} If $M$ is a compact manifold, then $\Diff(M)$ cannot be 
given a Banach--Lie group structure for which the natural action on~$M$ is smooth. 
\end{cor}

In Section 4 of \cite{OdH72}, Omori and de la Harpe construct 
an example of a Banach--Lie group $G$ acting smoothly and amply, but not primitively 
on~$\R^2$. 

The preceding discussion implies in particular that Banach--Lie groups 
rarely act on finite-dimensional manifolds. As the gauge groups of 
principal bundles $q \: P \to M$ over compact manifolds $M$ show, the 
situation is different for locally exponential Lie groups 
(cf.\ Theorem~\ref{thm:IV.1.12}). 
Therefore it is of some interest to have good criteria for 
the integrability of infinitesimal actions of locally exponential 
Lie algebras on finite-dimensional manifolds (cf.\ (FP7) in \cite{Ne06a}).  

We start with a more general setup for infinite-dimensional manifolds which 
need extra smoothness assumptions (\cite{AbNe08}): 

\begin{thm}[Integration of locally exponential Lie algebras of vector fields]
  \mlabel{thm:IX.1.8} 
Let 
$M$ be a smooth manifold modeled on a locally convex space, 
$\g$ a locally exponential Lie algebra and 
$\alpha \: \g \to {\cal V}(M)$ be a homomorphism of Lie algebras whose 
range consists of complete vector fields. 
Suppose further that the map 
$$ \Exp \: \g \to \Diff(M), \quad x \mapsto \Phi^{\alpha(x)}_1 $$
is smooth in the sense of {\rm Definition~\ref{def:e.1.1a}} and that 
$0$ is isolated in 
\[ \z(\g) \cap \Exp^{-1}(\id_M).\] Then there exists a 
locally exponential Lie group $G$ and a smooth action 
$\sigma \: G \times M \to M$ whose derived action 
$\dot\sigma \: \g \to {\cal V}(M)$ coincides with $\alpha$. 
\end{thm}

In the finite-dimensional case, the smoothness assumptions in Theorem~\ref{thm:IX.1.8} 
follows from the 
smooth dependence of solutions of ODEs on parameters and initial values, and 
the condition on the exponential function can be verified with methods 
to be found in \cite{MZ55}. This leads to the following less technical 
generalization of the Lie--Palais Theorem which subsumes in particular 
{ Omori}'s corresponding results for Banach--Lie algebras \cite[Thm.~A]{Omo80}, 
\cite[Thm.~4.4]{Pe95b}). 

\begin{thm}
  \mlabel{thm:IX.1.9} Let $M$ be a smooth finite-dimensional 
manifold, $\g$ a locally exponential Lie algebra and 
$\alpha \: \g \to {\cal V}(M)$ a continuous homomorphism of Lie algebras 
whose range consists of complete vector fields. Then there exists a 
locally exponential Lie group $G$ and a smooth action 
$\sigma \: G \times M \to M$ with $\dot\sigma = \alpha$. 
\end{thm}

The following result is a generalization of Palais' Theorem in another direction. 
Since $\Diff(M)$ is $\mu$-regular (\cite[Thm.~III.3.1]{Ne06a}), 
it also follows from \cite[Thm.~III.2.8]{Ne06a}. 


\begin{thm}
  \mlabel{thm:IX.1.10} {\rm(\cite{Les68})} If $M$ is compact, then 
a subalgebra $\g \subeq {\cal V}(M)$ is integrable to an integral subgroup 
if $\g$ is finite-dimensional or closed and finite-codimensional. 
\end{thm}

\subsection{Groups of diffeomorphisms as automorphism groups} 

In this subsection, we collect some results stating that 
automorphism groups of certain algebras, Lie algebras or groups associated to 
geometric structure on manifolds are not larger than the natural candidates. 
Most of the results 
formulated below for automorphisms of geometric structures attached to a manifold $M$ 
generalize to results saying that if $M_1$ and $M_2$ are two manifolds and 
two objects of the same kind attached to $M_1$ and $M_2$ are isomorphic, 
then this isomorphism can be implemented by a diffeomorphism $M_1 \to M_2$, compatible 
with the geometric structures under consideration. 

\begin{thm}
  \mlabel{thm:IX.2.1} Let $M$ be a $\sigma$-compact finite-dimensional 
smooth manifold. Then the following assertions hold: 
\begin{description}
\item[\rm(1)] For the Fr\'echet algebra $C^\infty(M,\R)$, each homomorphism 
to $\R$ is a point evaluation. 
\item[\rm(2)] $\Aut(C^\infty(M,\R)) \cong \Diff(M)$. 
\item[\rm(3)] $\Aut({\cal V}_c(M)) \cong \Aut({\cal V}(M)) \cong \Diff(M)$. 
\item[\rm(4)] If $M$ is complex and ${\cal V}^{(1,0)}(M)\subeq {\cal V}(M)_\C$ 
is the Lie algebra of complex vector fields of type $(1,0)$, then 
$\Aut({\cal V}^{(1,0)}(M)) \cong \Aut_{\cal O}(M)$ is the group of biholomorphic 
automorphisms of $M$. 
\item[\rm(5)] For each finite-dimensional $\sigma$-compact manifold $M$ and each 
simple (real or complex) finite-dimensional Lie algebra $\fk$, the natural homomorphism 
$$ C^\infty(M,\Aut_\K(\fk)) \rtimes \Diff(M) \to \Aut_\K(C^\infty(M,\fk)) $$
is bijective. 
\item[\rm(6)] If $M$ is a Stein manifold and $\fk$ is a finite-dimensional 
complex simple Lie algebra, then 
$\Aut({\cal O}(M,\fk)) \cong {\cal O}(M,\Aut(\fk)) \rtimes \Aut_{\cal O}(M)$, where 
$\Aut_{{\cal O}}(M)$ denotes the group of biholomorphic diffeomorphisms of $M$. 
\item[\rm(7)] If $K \subeq \C^n$ is a polyhedral domain and 
${\cal O}(K,\C)$ the algebra of germs of holomorphic 
$\C$-valued functions in $K$, then the  
group $\Aut({\cal O}(K,\C))$ consists of the germs of biholomorphic maps 
of $K$ and $\der({\cal O}(K,\C))$ consists of the germs holomorphic vector fields on $K$. 
\end{description}
\end{thm}

\begin{prf}
(1) (cf.\ \cite{My54} for the compact case; \cite{Pu52}; \cite{Co94}). 
(2) follows easily from (1) because each automorphism of the algebra $C^\infty(M,\R)$ 
acts on $\Hom(C^\infty(M,\R),\R) \cong M$. 

(3) The representability of each automorphism of ${\cal V}_c(M)$ by a diffeomorphism 
is due to Pursell and Shanks (\cite{PuSh54}), and the other assertion 
follows from \cite[Thm.~2]{Ame75}. It is based on the fact that 
the maximal proper subalgebras of finite codimension are 
all of the form 
\[ {\cal V}(M)_m := \{ X \in {\cal V}(M) \: X(m) = 0\} \] 
for some $m \in M$. They are permuted by each automorphism. 
Alternatively, one can use that 
all maximal ideals consist of all vector fields whose jet vanishes in some 
$m \in M$.

(4) follows from \cite[Thm.~1]{Ame75}. 

(5) (cf.~\cite[Prop.~3.4.2]{PS86}) A central point is that every non-zero endomorphism of 
$\fk$ is an automorphism. Further, it is used that 
$[\fk, C^\infty(M,\fk)] = C^\infty(M,\fk)]$ and that distributions supported 
by one point are of finite order. 

(6) \cite{NeW08b} 

(7) This is \cite[Thm.~III]{vHo52b}, where it is first shown that the maximal ideals 
in the Silva cia ${\cal O}(K,\C)$ \cite[Ex.VIII.3(d)]{Ne06a} 
are the kernels of the point evaluations \cite[Thm.I]{vHo52b}). 
\end{prf}

\begin{rem}
  \mlabel{rem:IX.2.2} Let $K \subeq \C^n$ be a compact subset and 
$\Aut_{{\cal O}}(K)$ the group of germs of bihomolorphic maps, defined on some neighborhood 
of $K$, mapping $K$ onto itself. 
In \cite{vHo52a}, van Hove introduces a group topology on 
this group as the topology for which the map 
$$ \Aut_{{\cal O}}(K) \to {\cal O}(K,\C^n) \times {\cal O}(K,\C^n), \quad 
g \mapsto (g,g^{-1}) $$
is an embedding. He shows that, under certain geometric conditions
 on the set $K$, this group is complete and contains no small subgroups. 
Moreover, its natural action on ${\cal O}(K,\C)$ is continuous. 
\end{rem}

\subsection{Automorphism groups of differential forms} 

\begin{thm}
  \mlabel{thm:IX.2.3} {\rm(\cite[\S 10]{Omo74})} 
Let $M$ be a $\sigma$-compact finite-dimensional 
smooth manifold. For a differential form $\alpha$ on $M$ we put 
\[ {\cal V}(M,\alpha) := \{ X \in {\cal V}(M) \: {\cal L}_X\alpha = 0\}.\] 
Then the following assertions hold: 
\begin{description}
\item[\rm(1)] If $\mu$ is a volume form or a symplectic form on $M$, 
then every automorphism of the real Lie algebra  
${\cal V}(M,\mu)$ is induced by an element of the group 
$$\{ \phi \in \Diff(M) \: \phi^*\mu \in \R \mu\}.$$
\item[\rm(2)] If $\alpha$ is a contact $1$-form on $M$, then every
  automorphism of the real Lie algebra 
${\cal V}(M,\alpha)$ is induced by an element of the group 
\[ \{ \phi \in \Diff(M) \: \phi^*\alpha \in C^\infty(M,\R^\times)\cdot \alpha\}. \] 
\end{description}
\end{thm}

\subsection{Automorphism groups of Banach symmetric spaces}

\begin{defn}
Let $M$ be a smooth manifold. We say that $(M,\mu)$
is a {\it symmetric space} (in the sense of Loos) \cite{Lo69} if 
$\mu \: M \times M \to M, (x,y) \mapsto x \cdot y$ 
is a smooth map with the following properties: 
\begin{description}
\item[\rm(S1)] $x \cdot x$ for all $x \in M$. 
\item[\rm(S2)] $x \cdot (x \cdot y) =y$ for all $x,y \in M$. 
\item[\rm(S3)] $x \cdot (y \cdot z) = (x \cdot y) \cdot (x \cdot z)$ 
for all $x,y \in M$. 
\item[\rm(S4)] $T_x(\mu_x) = -\id_{T_x(M)}$ for $\mu_x(y) := \mu(x,y)$ and each $x \in M$. 
\end{description}
\end{defn}


\begin{thm} {\rm(M.~Klotz, \cite{Kl11, Kl12})}
  If $(M,\mu)$ is a connected Banach symmetric space, then
  \[ \Aut(M,\mu) := \{ \phi \in \Diff(M) \:
    \mu \circ (\phi \times \phi) = \phi \circ \mu \} \]
  carries a natural Banach--Lie group structure.
\end{thm}

\subsection{Automorphism groups of singularities} 

In \cite{Omo80}, one finds another interesting result of this type. 
Let $V$ be a germ of an affine variety in $0 \in \C^n$. Two such germs $V$ and $V'$ 
are said to 
be {\it biholomorphically equivalent} if there exists an element 
$\phi \in \Gh_n(\C)$ of the group of germs of biholomorphic maps fixing $0$ 
(as in \cite[Ex.~VI.2.12]{Ne06a}), such that 
$\phi(V) = \phi(V')$. On the infinitesimal level the automorphisms of a germ 
$V$ are given by the Lie algebra 
$$ \g(V) := \{ X \in \gh_n(\C) \: X.J(V) \subeq J(V)\}, $$
where $J(V) \subeq {\cal O}(0,\C)$ (the germs of holomorphic functions in $0$) 
is the annihilator ideal of $V$. Let $\g(V)_k \trile \g(V)$ denote the ideal 
consisting of all vector fields vanishing of order $k$ in $0$ and form the 
projective limit Lie algebra 
$$ \oline\g(V) := \prolim \g(V)/\g(V)_k, $$
which can be viewed as a Fr\'echet completion of $\g(V)$. 

An element $X \in \gh_n(\C)$ is called {\it semi-expansive} 
if it is $\Gh_n(\C)$-conjugate to a linear diagonalizable vector field for which all 
eigenvalues lie in some open halfplane. 
The germ $V$ is called an {\it expansive singularity} if 
$\g(V)$ contains an expansive vector field. 

\begin{thm}
  \mlabel{thm:IX.2.4} Two expansive singularities $V$ and $V'$ are biholomorphically equivalent 
if and only if the pro-finite 
Lie algebras $\oline\g(V)$ and $\oline\g(V')$ are isomorphic. 
Moreover, $\Aut(\oline\g(V))$ can be identified with the stabilizer $\Gh_n(\C)_V$ of 
$V$ in the group $\Gh_n(\C)$. 
\end{thm}

On the group level, we have the following analog of Theorem~\ref{thm:IX.2.3},  
see \cite{Fil82} for (1) 
and \cite[Thm.~7.1.4/5/6]{Ban97} for (2)-(4).  

\begin{thm}
  \mlabel{thm:IX.2.5}  Let $M$ be a $\sigma$-compact connected finite-dimensional 
smooth manifold. Then the following assertions hold: 
\begin{description}
\item[\rm(1)] Every automorphism of the discrete group $\Diff(M)$ is inner. 
\item[\rm(2)] If $\alpha$ is a contact $1$-form on $M$, then every
  automorphism of the discrete group 
$\Diff(M,\alpha)$ is conjugation with an element of the group 
\[ \{ \phi \in \Diff(M) \: \phi^*\alpha \in C^\infty(M,\R^\times)\cdot \alpha\}.\] 
\item[\rm(3)] If $\omega$ is a symplectic form and $M$ is compact of dimension $\geq 2$, 
  then every automorphism of the discrete group
  $\Diff(M,\omega)$ is conjugation by an element of the group 
$$\{ \phi \in \Diff(M) \: \phi^*\omega \in \R \omega\}.$$
\item[\rm(4)] If $\mu$ is a volume form and $M$ is of dimension $\geq 2$, 
  then every automorphism of the discrete group
  $\Diff(M,\mu)$ is conjugation by an element of the group 
$\{ \phi \in \Diff(M) \: \phi^*\mu \in \R \mu\}.$
\end{description}
\end{thm}

\subsection*{Open Problems} 

\begin{probl}
  \mlabel{prob:OX.1} Let ${\cal V}(M)_{\rm cp}$ denote the set of complete 
vector fields on the finite-dimensional manifold $M$. Then we have an exponential 
function 
$$ \Exp \: {\cal V}(M)_{\rm cp} \to \Diff(M), \quad X \mapsto \Phi^X_1. $$
Is it true that $0$ is isolated in $\Exp^{-1}(\id_M)$ with respect to the natural 
Fr\'echet topology on ${\cal V}(M)$? 

That this is true for compact manifolds follows from  
Newman's Theorem \cite[Thm.~2]{Dr69}.
For the proof of Theorem~\ref{thm:IX.1.9}, we show 
for each continuous homomorphism $\alpha \: \g \to {\cal V}(M)$ of a 
locally exponential Lie algebra $\g$ to ${\cal V}(M)$ with 
range in ${\cal V}(M)_{\rm cp}$ that $0$ is isolated in 
$(\Exp \circ \alpha)^{-1}(\id_M)$, which is a weaker statement. 

Since the set $\Exp^{-1}(\id_M)$ is in one-to-one correspondence with the 
smooth $\T$-actions on $M$, the problem is to show that the trivial 
action is isolated in this ``space'' of all smooth $\T$-actions on $M$. 
\end{probl}

\begin{ex}
If $M$ is the real Hilbert space $\ell^2(\N,\R)$ with the Hilbert basis $e_n$, $n \in \N$, 
then we have linear vector fields $X_n(v) := 2\pi i \la v, e_n \ra e_n$ with 
$\exp(X_n) = \id_M$ and $X_n \to 0$ uniformly on compact subsets of $E$. Hence 
the finite-dimensionality of $M$ is crucial. 
\end{ex}

\begin{probl}
  \mlabel{prob:IX.4} Show that, for each compact subset $K \subeq \C^n$, the group 
$\Aut_{{\cal O}}(K)$ from Remark~\ref{rem:IX.2.2} 
is a Lie group with respect to the manifold structure 
inherited from the embedding into ${\cal O}(K,\C^n)$. 
\end{probl}

\begin{probl}
  \mlabel{prob:IX.5} (Automorphisms of gauge algebras) 
Let $q \: P \to M$ be a smooth $K$-principal bundle over the (compact) manifold $M$. 
Determine the group $\Aut(\gau(P))$ of automorphisms of the gauge Lie algebra 
(Example~\ref{exs:5.1.4}(c)). 
Does it coincide with the automorphism group $\Aut(\ad(P))$ of the adjoint bundle 
$\ad(P)$,  whose Lie algebra of sections is isomorphic to $\gau(P)$? If $K$ is a simple complex Lie group, 
then the results in \cite[Thm.~16]{Lec80} provide a local description of the 
automorphisms of this Lie algebra in terms of diffeomorphisms of $M$ and 
sections of the automorphism bundle $\Aut(\ad(P))$ (see also Theorem~\ref{thm:IX.2.1}). 
\end{probl}

\begin{probl}
  \mlabel{prob:IX.6} Determine the automorphism groups of the Lie algebras 
$\gf_n(\K)$, $\gs_n(\K)$ and $\gh_n(\C)$. 
\end{probl}

\begin{probl}
  \mlabel{prob:IX.7} Characterize connected Banach--Lie groups acting 
\break smoothly, effectively and transitively on a finite-dimensional manifold. 
In view of Theorem~\ref{thm:IX.1.6}, for each Banach--Lie group $G$, 
the Lie algebra $\g$ contains a finite-codimensional closed solvable ideal. 
If, conversely, $\g$ is a Banach--Lie algebra with a finite-codimensional 
closed solvable ideal, then Theorem VI.1.19 implies that $\g$ is enlargeable. 
Under which conditions do the corresponding groups $G$ act effectively on 
some finite-dimensional homogeneous space? (see also the corresponding 
discussion in \cite{Omo97}). 
\end{probl}

\begin{probl}
  \mlabel{prob:IX.8} 
Let $G$ be a Banach--Lie group and $H \subeq G$ be a closed subgroup for which 
$\L^e(H)$ has finite codimension. Does this imply that $G/H$ is a manifold? 
\end{probl}

\section{Circle actions on manifolds}

In this subsection we show that, for the diffeomorphism group $\Diff(M)$ 
of a compact smooth manifold, 
$0$ is isolated in the set of complete vector fields $X$ satisfying 
$\Phi^X_1 = \id_M$ (cf.\ Appendix~\ref{app:nonlie}). 
This implies in particular that each continuous 
morphism $\alpha \: \h\to {\cal V}(M)$, $\h$ locally exponential, satisfies the 
assumptions of the Integral Subgroup Theorem~\ref{thm:5.5.3}. 

Let $M$ be a finite-dimensional manifold. 
For any smooth action $\sigma$ of the circle group 
$\T \cong \R/\Z$, the infinitesimal generator is a vector field 
satisfying $\Phi^X_1 = \id_M$. Conversely, every complete vector field 
$X \in {\cal V}(M)$ with $\Phi^X_1 = \id_M$ defines a 
smooth $\T$-action on $M$. 

\begin{thm} \mlabel{thm:diff-1isol} If $M$ is a compact manifold, then there exists a 
$0$-neigh\-bor\-hood $U$ in ${\cal V}(M)$ such that 
$\Phi^X_1 = \id_M$ for $X \in U$ implies $X = 0$. 
\end{thm}

\begin{prf} Let $g$ be a Riemannian metric on $M$ and $d_M \: M \times M \to \R$ the 
corresponding distance function. We then 
obtain a continuous norm on ${\cal V}(M)$ by 
$$ \|X\|_\infty := \sup_{m \in M} \|X(m)\|= \sup_{m \in M} \sqrt{g_m(X(m), X(m))}. $$
For $\|X\|_\infty \leq r$ and $0 \leq t \leq 1$, we then have 
$d_M(\Phi^X_t(p), p) \leq r$. 

According to Newman's Theorem (\cite[Thm.~2]{Dr69}), there exists, for each 
smooth manifold $M$ with a metric $d_M$,  an $\eps > 0$ 
such that every action of a finite group $G$ on $M$ satisfying 
$$d_M(g.x,x) < \eps \quad \mbox{ for } \quad g \in G, x \in M, $$
is trivial. 

Applying this to the smooth $\T$-action obtained from $X \in {\cal V}(M)$ 
with $\Phi^X_1 = \id_M$, it follows that $\|X\|_\infty < \eps$ implies that 
$\Phi^X_{\frac{1}{n}}$ defines an action of $\Z/n\Z$ on $M$, and 
that this action is trivial for each $n \in \N$. We conclude that 
$X = 0$. 
\end{prf}

\begin{ex} We give a direct argument for Theorem~\ref{thm:diff-1isol} in the case $M = \T = \R/\Z$. 

Let $X \in {\cal V}(\T)$ be a non-zero vector field with $\Phi^X_1 = \id_M$. 
Then the corresponding $\T$-action has a non-trivial orbit, which 
is a compact $1$-dimensional submanifold of $M$, hence also open. 
This implies that the action is transitive and hence that $X$ has no zeros. 

Let $X_0 := \frac{\partial}{\partial\theta}$ be the standard generator 
of the rotation action on $\T \cong \R/\Z$. Then 
$X(\theta) = f(\theta) X_0$ for some smooth $1$-periodic function $f \: \R \to \R$ 
without zeros. 

For $\phi \in \Diff(\T)$, we then have 
$\Ad(\phi)X  = f_\phi X_0,$
where 
$$ f_\phi(\theta) = f(\phi^{-1}(\theta)) \phi'(\phi^{-1}(\theta)), $$
resp., 
$f_\phi \circ \phi = f \cdot \phi'.$
Let 
$$c := \frac{1}{\int_0^1 \frac{dt}{f(t)}} 
\quad \mbox{ and } \quad 
\phi(\theta) := \int_0^\theta \frac{c}{f(t)}\, dt. $$
Then $\phi(1) = 1$ and 
$$ \phi(\theta + 1) 
= \phi(\theta) + \int_\theta^{\theta+1} \frac{c}{f(t)}\, dt
= \phi(\theta) + \int_0^{1} \frac{c}{f(t)}\, dt
= \phi(\theta) + 1. $$
Hence $\phi \: \R \to \R$ defines an element of $\Diff(\T)$. 
This diffeomorphism satisfies 
$\phi' \cdot f = c,$
which leads to 
$\Ad(\phi)X = c \cdot X_0.$

Now $\Phi^X_1 = \id_\T$ implies that $\Phi^{X_0}_c = \id_\T$, which is equivalent 
to $c \in \Z$. We conclude that 
\[  \int_0^1 \frac{dt}{|f(t)|} \in \Big\{ \frac{1}{n} \: n \in \N\Big\} \] 
whenever $X$ is non-zero with $\Phi^X_1 = \id_\T$. 
For $|f| \leq \frac{1}{2}$ we already obtain 
$\int_0^1 \frac{dt}{|f(t)|} \geq 2,$
so that $\Phi^X_1 \not= \id_\T$ in this case. 
\end{ex}

\begin{rem} It is also instructive so see that, for vector fields with 
zeros, one can give a more direct argument for the conclusion 
of Theorem~\ref{thm:diff-1isol}. This applies in particular 
to compact orientable manifolds with non-zero Euler characteristic $\chi(M)$ 
on which every vector field has a zero (\cite[Satz~13.11]{tD00}). 

Let $g$ be a Riemannian metric on the compact connected smooth manifold $M$. 
If $p \in M$ is a zero of the vector field $X \in {\cal V}(M)$, 
then the vector field $T(X) \in {\cal V}(T(M))$ restricts to a linear 
vector field $T_p(X)$ on the submanifold $T_p(M)$ of $T(M)$. 
Suppose that $\|T_p(X)\| < 2\pi$, where the norm is the operator norm 
corresponding to a Riemannian metric on $T_p(M)$. 
If $\Phi^X_1 = \id_M$, then $\Phi^{T(X)}_1 = \id_{T(M)}$, and in particular 
$\exp_{\GL(T_p(M))} T_p(X) = \id_{T_p(M)}$ holds in $\GL(T_p(M))$. 
This implies that $T_p(X)$ is semisimple with 
eigenvalues in $2\pi i \Z$. Now the assumption $\|T_p(X)\| < 2\pi$ 
leads to $T_p(X) = 0$. 

Let $\sigma \: \T \to \Diff(M), \sigma(t) := \Phi^X_t$ 
denote the smooth $\T$-action on $M$ generated by $X$ 
and $\tilde g$ a $\T$-invariant Riemannian metric on $M$. 
Then the corresponding exponential function $\Exp \: T(M) \to M$ is 
$\T$-equivariant. From $T_p(X) = 0$ it follows that $\T$ acts trivially 
on $T_p(M)$, hence on $\Exp(T_p(M))$, which contains a neighborhood of $p$. 
This argument implies that the set $M^\T$ of fixed points is open. As it is 
also closed, the assumption $M^\T \not=\eset$ leads to $M^\T = M$, and hence 
to $X = 0$. 
\end{rem}




\begin{small}

\subsection*{Exercises for Chapter~\ref{ch:diffeo}} 

\begin{exer}\mlabel{exer:I.5} (cf.\ \cite{Str06}) 
Let $X$ be a topological space and endow the set 
$C(X,X)$ of continuous self-maps of $X$ with the compact open topology. 
We endow the group $\Homeo(X)$ with the initial topology with respect to the map 
$$\Homeo(X) \to C(X,X)^2,\quad  \phi \mapsto (\phi,\phi^{-1}).$$ 
Show that, if $X$ is locally compact, then this topology turns 
$\Homeo(X)$ into a topological group. 
Hint: If $f\circ g \in \lfloor K,O \rfloor$ choose 
a compact subset $K'$ and an open subset $O'$ with 
$g(K) \subeq O' \subeq K' \subeq f^{-1}(O)$ 
\end{exer}

\begin{exer} Let $E$ and $F$ be locally convex spaces. Then 
the inclusion $\cL(E, F)_{c.o.} \into C^\infty(E,F)$ 
is a topological embedding. 
\end{exer}

\end{small}

\section{Notes and comments on Chapter~\ref{ch:diffeo}} 

That, for 
$G := \Diff(\bS^1)$, the image of the exponential 
function is not a neighborhood of $\be$ (Theorem~\ref{thm:V.1.7}) 
can already be found in  \cite{Ham82} and \cite{Mil84}; 
see also \cite[Ex.~43.2]{KM97}. 
More generally, for any compact manifold $M$, 
the image of the exponential function of $\Diff(M)$ does not contain any 
identity neighborhood (cf.\ \cite{Gr88}, \cite{Pali68, Pali74}, 
and \cite{Fre68} for some $2$-dimensional cases). 

Lie algebras associated to the group of homeomorphisms of the circle 
are constructed by Malikov and Penner in \cite{MP98}; in particular 
they study vector fields which are piecewise~$\fsl_2$. For  
more on topological aspects of groups acting on the circle, 
see \cite{Gh01}. This paper identifies $\SO_2(\R)$ as a maximal 
compact subgroup of $\Homeo_+(\bS^1)$ and contains various classification 
results concerning finite-dimensional subgroups. 
This is also related to the classification results on 
finite-dimensional Lie subalgebras of $\cV(\bS^1)$ in \cite{Sp17}. 
In this context, an important result that can be found in 
\cite[pp.~208ff]{MZ55} (cf.~Theorem~\ref{thm:IX.1.1}) 
is that any action of a finite-dimensional Lie group~$G$ 
by $C^k$-diffeomorphism on a finite-dimensional manifold $M$ corresponds to a 
$C^k$-map $G \times M \to M$. If $M$ is compact and $k = \infty$, we thus 
obtain Lie group homomorphisms $G \to \Diff(M)$. 

A recent construction of the Virasoro--Bott group, a non-trivial 
central extension of $\Diff(\bS^1)$ in the context of 
Chern--Simons Theory and \break $3$-dimensional gravity has been given 
in \cite{FK21}.

For interesting uniqueness results for topologies on groups 
of diffeomorphisms and homeomorphisms, see \cite{Kal86}.

The regularity of the composition of $H^s$-diffeomorphisms of 
a compact manifold $M$ for $s > \frac{1}{2}\dim M + 1$ has 
recently been discussed in detail in \cite{IKT13}.

Results on the non-existence of conjugation invariant norms on 
diffeomorphism groups can be found in \cite{BIP08}. 

Various types of groups of diffeomorphisms of $\R^n$ 
are analyzed  with a perspective of applications 
to shape spaces in \cite{MM13}. 

%
%
\def\date{22.11.21} 

\chapter{Homotopy groups of infinite-dimensional 
Lie groups} 
\mlabel{ch:top}

In this chapter we address several aspects of the topology of infinite-di\-men\-sio\-nal 
Lie groups. For a Lie group $G$, the most important topological information is contained 
in the first three homotopy groups  $\pi_0(G)$ (the group of connected components), 
$\pi_1(G)$ (the fundamental group), and the second homotopy group $\pi_2(G)$. 
The importance of $\pi_0(G)$ is clear 
because one often needs to know whether a concretely 
given group is connected and the influence of the Lie algebra 
does not reach beyond the identity component (cf.\ Proposition~\ref{propc.14}). 
Information on the fundamental group 
is important for the integration of Lie 
algebra homomorphisms to group homomorphisms and hence in particular for 
representation theory (Theorem~\ref{thm3.2.11}). It also shows up in the integration of 
$1$-forms satisfying the Maurer--Cartan equation (Theorem~\ref{thm-fundamental}). 
The interest in $\pi_2(G)$ stems 
from the crucial role this group plays for enlargibility of Lie algebras and for extensions of $G$ 
(Subsection~\ref{subsec:per-grp}). 

We start this chapter with a section introducing higher homotopy 
groups of pointed spaces $(X,x_0)$ and the long exact homotopy sequence 
of topological principal bundles. We then develop some methods (originally 
due to R.~Palais)  
that are often helpful to calculate homotopy groups of open subsets of locally 
convex spaces. In Section~\ref{sec:10.2} we briefly recall some of the key 
results on homotopy groups of finite-dimensional Lie groups. 
The same is done in Section~\ref{sec:top-diffeo} for groups of 
diffeomorphisms. 

Homotopy groups of groups of operators on Hilbert spaces are discussed 
in some detail in the next three sections. 
Here Kuiper's Theorem on the contractibility 
of the unitary group of an infinite-dimensional Hilbert space 
(real, complex or quaternionic) is a cornerstone (Section~\ref{sec10.4}). 
In Section~\ref{sec:top-dirlim} we then discuss direct limit groups 
and this is applied in Section~\ref{sec:top-schatten} to groups associated 
to the Schatten ideals.

To treat mapping groups, such as $C^r_c(M,K)$, for $r \in \N \cup \{\infty\}$, 
we start with some 
approximation results in Section~\ref{sec:top.4}, asserting that  
inclusions such as 
$C^r_c(M,K) \into C_c(M,K) \into C_*(M_\infty,K)$ are weak homotopy equivalences. 
This permits us in Section~\ref{sec10.6} to actually calculate homotopy groups 
of $C_*(M,K)$ in terms of the homotopy groups of $K$ in the case where 
$M$ is a sphere, a torus or a $2$-dimensional surface. 
In the final Section~\ref{sec10.6} we take a closer look at mapping groups on 
products of spheres, in particular tori, and on compact 
surfaces (orientable or not) with boundary.

\section{Basics of homotopy groups} \mlabel{sec:10.1}

In this first section we start with an introduction to 
higher homotopy groups and the long exact homotopy sequence of fiber 
bundles and discuss some applications to homogeneous spaces 
(Subsection~\ref{subsec:14.1.3}). 
We then continue in Subsection~\ref{subsec:14.1.2} with some generalities 
on homotopy groups of open subsets in (locally convex) direct limits. 
Here Palais' Theorem on weak homotopy equivalences 
provides an important tool for concrete calculations of homotopy groups 
(Subsection~\ref{subsec:14.1.3b}).

\subsection{Homotopy groups} 
\mlabel{subsec:14.1.3}

In this subsection, we present some basic facts and definitions 
concerning homotopy groups of spaces and, in particular, of topological 
groups. 

\begin{defn}
  \mlabel{def-homot}
(a) Let $X$ and $Y$ be topological spaces. Two maps 
$f, g \: X \to Y$ are said to be 
{\it homotopic} \index{homotopic maps} 
(denoted $f \sim g$) 
if there exists a continuous map 
$H \: [0,1] \times X \to Y$, a so-called 
\index{homotopy} 
{\it homotopy}, with 
$$ H_0(x) := H(0,x) = f(x) 
\quad \mbox{ and } \quad 
 H_1(x) := H(1,x) = g(x) $$
for all $x \in X$. We write $[X,Y]$ for the set of 
homotopy classes of continuous maps $X \to Y$ and 
$[f]$ for the homotopy class of a map~$f$. 

(b) An important refinement of this concept is obtained as follows. 
A {\it space pair} \index{space pair} 
is a pair $(X, X_0)$ of topological spaces, 
where $X_0 \subeq X$ is a subspace. For two space 
pairs $(X, X_0)$ and $(Y, Y_0)$ we write 
$$ C((X, X_0), (Y, Y_0)) := \{ f \in C(X, Y) \: 
f(X_0) \subeq Y_0\} $$
for the set of continuous maps $X \to Y$ mapping $X_0 \to Y_0$. 
Two maps $f,g \: (X, X_0) \to (Y, Y_0)$ are 
called {\it homotopic} \index{homotopic maps of space pairs} 
if there exists a homotopy 
$H \: [0,1] \times X \to Y$ of $f$ and $g$ for which 
all maps $H_t := H(t,\cdot)$ map $X_0 \to Y_0$. We write 
$[(X, X_0), (Y, Y_0)]$ for the set of homotopy classes of 
maps of pairs. 

(c) A particularly important special case arises if $X_0 := \{x_0\}$ 
and $Y_0 := \{y_0\}$. Then the pairs $(X, x_0)$, $(Y, y_0)$ 
are called {\it pointed spaces} \index{pointed space} 
and we write 
$$ C_*(X,Y) := C((X, x_0), (Y,y_0)) $$
for the set of continuous base point preserving maps. 
We likewise write $[X,Y]_*$ for the homotopy classes of 
base point preserving maps. 

(d) There is a refinement of these concepts to space 
triples $(X, X_0, X_1)$ with $X_1 \subeq X_0 \subeq X$ and 
$$ C((X, X_0, X_1), (Y, Y_0, Y_1)) := \{ f \in C(X, Y) \: 
f(X_0) \subeq Y_0, f(X_1) \subeq Y_1 \}.$$
Accordingly, we write $[(X, X_0, X_1), (Y, Y_0, Y_1)]$ for the set of 
homotopy classes of maps of space triples. 
\end{defn}

Sets of homotopy classes form important topological invariants of 
topological spaces, but a priori they are only sets not carrying any  
algebraic structure. In the following, we shall discuss two important situations
where $[X,Y]$ carries a group 
structure, namely if $Y$ is a topological group and if $(X, x_0)$ is a pointed 
sphere $\bS^n$.

\begin{defn}
  \mlabel{def-homotgrp} 
We now define 
\index{relative homotopy groups of pairs} 
{\it relative homotopy groups of pairs}. 
Let $I^n := [0,1]^n$ denote the $n$-dimensional cube. Then the
boundary $\partial I^n$ of $I^n$ can be written as 
$I^{n-1} \cup J^{n-1}$, where $I^{n-1} := \{ x \in I^n \: x_1 = 1\}$ 
is called the {\it initial face} 
and $J^{n-1}$ is the union of all other faces of $I^n$. 

Let $(X, X_0, \{x_0\})$ be a space triple. We write 
$$F^n(X,X_0,x_0) := C((I^n, \partial I^n, J^{n-1}), (X,X_0, \{x_0\})) $$
for the continuous maps of space triples and 
$$ \pi_n(X,X_0,x_0) := [(I^n, \partial I^{n}, J^{n-1}), (X,X_0, \{x_0\})] $$
for the set of homotopy classes of such maps. 
If $X_0 = \{x_0\}$, then we simplify the notation to 
$$F^n(X,x_0) := C((I^n, \partial I^{n}), (X,x_0)) $$
and 
\[  \pi_n(X,x_0) := [(I^n, \partial I^{n}), (X,x_0)] 
= [I^n/\partial I^n, X]_*
= [\bS^n, X]_*.\] 
Here we use the image of $\partial I^n$ as the base point of the quotient 
$I^n/\partial I^n \cong \bS^n$. 

For $f,g \in F^n(X,X_0, x_0)$, $n \geq 2$, we define a product
$$ (f * g)(t_1, \ldots, t_n) 
:= \left\{   \begin{array}{cl}
f(t_1, t_2, \ldots, 2t_n) & \mbox{for $0 \leq t_n \leq \frac{1}{2}$} \\ 
g(t_1, t_2, \ldots, 2t_n-1) & \mbox{for $\frac{1}{2} < t_n \leq 1$} \\ 
  \end{array}
\right.  $$ 
and note that $f*g$ is continuous because $f(x) = g(x) = x_0$ for $x_n = 1$.
For $n \geq 1$ and $X_0 = \{x_0\}$, this product also yields 
an element of $F^n(X,x_0)$. 

Now one verifies that (Exercise~\ref{exer:10.1.2}): 
\begin{enumerate}
\item[\rm(1)] $f \sim f', g \sim g' \Rarrow f * g \sim f' * g'$, so that 
we obtain a well-defined product 
$$ [f] * [g] := [f * g] $$
of homotopy classes.  
\item[\rm(2)] (Neutral element) 
If $x_0$ denotes the constant map $I^n \to X, t \mapsto x_0$, then 
$$ [x_0] * [f] = [f] * [x_0] = [f]. $$
\item[\rm(3)] (Inverse) Let $\oline f(t_1, \ldots, t_n) := f(t_1, \ldots, t_{n-1}, 1-t_n)$. 
Then $[f] * [\oline f] = [\oline f] * [f] = [x_0]$. 
\item[\rm(4)] (Associativity) $[f * g] * [h] = [f] * [g * h]$ 
for $f,g,h \in F^n(X,X_0,x_0)$, $n \geq 2$, or $n = 1$ and $X_0 = \{x_0\}$. 
\item[\rm(5)] (Commutativity) $[f] * [g] = [g] * [f]$ 
for $f,g \in F^n(X,X_0,x_0)$, $n \geq 3$, or $n = 2$ and $X_0 = \{x_0\}$. 
\item[\rm(6)] (Functoriality) For $f,g \in F^n(X,X_0,x_0)$ 
and $\phi \in C((X,X_0,x_0), (Y,Y_0,y_0))$ we have 
$$ [(\phi \circ f) * (\phi \circ g)] = [\phi \circ (f * g)], $$ 
so that $\phi$ induces a homomorphism 
$$ \phi_* = \pi_n(\phi) \: \pi_n(X,X_0,x_0) \to \pi_n(Y,Y_0, y_0). $$
\end{enumerate}
The so obtained groups 
$$ \pi_n(X,X_0, x_0), \quad n \geq 2 \quad \mbox{ and } \quad
\pi_n(X,x_0), \quad n \geq 1, $$
are called the 
\index{homotopy group of the pair} 
\index{homotopy group} 
{\it homotopy groups of the pair $(X,X_0)$ with respect 
to the base point $x_0$}, resp., the 
 {\it homotopy groups of $X$ with respect 
to the base point $x_0$}. 
\end{defn}

\begin{rem} \mlabel{rem:10.1.3} 
Let $G$ be a topological group and consider the identity element 
$\be$ as a base point. Then the sets 
$F^n(G,\be)$, $n \in \N_0$, 
also carry a natural group structure given by the pointwise 
product $(f \cdot g)(x) := f(x)g(x)$ and one can show that 
(Exercise~\ref{exer:10.1.3}): 
\begin{enumerate}
\item[\rm(1)] $f\sim f', g \sim g' \Rarrow f \cdot g \sim f' \cdot g'$, 
so that we obtain a well-defined product 
$$ [f] \cdot [g] := [f \cdot g] $$
of homotopy classes for all $n \in \N_0$.  
\item[\rm(2)] $f \sim g \Longleftrightarrow fg^{-1} \sim \be$, 
the constant map. Here $g^{-1}$ is the pointwise inverse.  
\item[\rm(3)] (Commutativity) $[f] \cdot [g] = [g] \cdot [f]$ 
for $f,g \in F^n(G,\be)$, $n \geq 1$. 
\item[\rm(4)] (Consistency) $[f] \cdot [g] = [f] * [g]$ 
for $f,g \in F^n(G,\be)$, $n \geq 1$. 
\end{enumerate}
As a consequence of (4), we can calculate the product of homotopy 
classes in 
$$ \pi_n(G) := \pi_n(G,\be) $$
via pointwise products of representatives. Note also that 
$\pi_0(G)$ can be identified with the group $G/G_a$ of arc-components 
of $G$, where $G_a = [\be] \trile G$ is the arc-component of the identity, 
which is a normal subgroup. 
\end{rem}

A proof of the following theorem can be found in 
\cite[\S VII.5]{Bre93}.  
\index{Long exact homotopy  sequence of a pair}
\begin{thm} [The long exact homotopy  sequence of a pair] \mlabel{pair-seq}
For any pointed space pair $(X,Y,y_0)$,  we have a natural map 
$$ \partial \: \pi_n(X,Y,y_0) \to \pi_{n-1}(Y,y_0), \quad 
[f] \mapsto [f\res_{I^{n-1}}] $$
which is a group homomorphism for $n \geq 2$. Together 
with the homomorphisms induced by the inclusion $Y \into X$, 
we thus obtain a long exact homotopy sequence 
\begin{align*}
\cdots  &\to \pi_k(X,y_0) \to \pi_k(X,Y,y_0) \sssmapright{\partial} \pi_{k-1}(Y,y_0) 
\to \ldots \\
 \ldots &\to \pi_1(X,y_0) \to \pi_1(X,Y,y_0) 
 \sssmapright{\partial} \pi_0(Y,y_0) \to \pi_0(X,y_0), 
\end{align*}
where the last three maps are not  group homomorphisms, and in this 
case exactness means that the range equals the inverse image of the base 
point, i.e., the homotopy class of the constant map $y_0$. 
\end{thm}

In the following, we shall frequently use an important consequence 
of the preceding theorem for principal bundles. 

\begin{defn} Let $G$ be a topological group and 
$X$ a topological space. A 
\index{principal bundle} 
{\it $G$-principal bundle over $X$} 
is a right $G$-space $P$ together with a surjective continuous map 
$q \: P \to X$ with the following properties: 
\begin{enumerate}
\item[(1)] The fibers of $q$ coincide with the $G$-orbits in $P$. 
\item[(2)] (Local triviality) For each $x \in X$, there exists an open 
neighborhood $U$ and a  homeomorphism 
$\phi \: U \times G \to q^{-1}(U),$
satisfying 
\[  q \circ \phi = \pr_U \quad \mbox{ and } \quad 
\phi(u,hg) = \phi(u,h).g, \ u \in U,h,g \in G. \]
\end{enumerate}  
\end{defn}

\begin{rem} (a) For each $G$-principal bundle $q \: P \to X$, 
the group $G$ acts freely on $P$ and we may identify $X$ with the 
set $P/G$ of $G$-orbits in $P$, endowed with the quotient topology 
(Exercise~\ref{exer:10.1.4}). 

(b) Let $G$ be a topological group and $H$ a closed subgroup, so that 
the space $X := G/H$ of $H$ left-cosets carries a natural Hausdorff topology 
(the quotient topology). We then have a continuous open map 
$q \: G \to X$ whose fibers coincide with the orbits of $H$ under the right 
action of $H$ on $G$ given by $g.h := gh$. Now the local triviality requirement 
is equivalent to the existence of some open subset $U \subeq G/H$ 
on which there exists a continuous section $\sigma \: U \to G$ of the 
quotient map $q$. In view of the equivariance of $q$ with respect to the 
$G$-left action, 
this implies the local triviality (Exercise~\ref{exer:10.1.5}), so that 
$G$ carries the structure of an $H$-principal bundle if local 
continuous sections exist. 

The local triviality condition is always satisfied if $G$ is a finite-dimensional 
Lie group or if $G$ is a Banach--Lie group 
and the closed subspace $\Lie(H)$ of $\Lie(G)$ is complemented 
(Exercise~\ref{exer:10.1.6}). More generally, it holds if $H$ is a split 
Lie subgroup (cf.\ Proposition~\ref{prop:split-Lie}). 
\end{rem}

\begin{rem} (Principal bundles and coverings) If 
$G$ is a discrete group, then a $G$-principal bundle 
$q \: P \to X$ is a covering map for which the group 
of {\it deck transformations} acts simply transitively on each fiber. 
\index{deck transformations} 

Suppose, conversely, that $q \: P \to X$ is a covering map which is 
regular in the sense that the group 
$$ \Delta := \Deck(P,q) := \{ \phi \in \Homeo(P) \: q \circ \phi = q\} $$
of deck transformations acts transitively on each fiber.  
Assume, in addition, that $P$ is connected and $X$ is locally connected. 
If a deck transformations $\phi$ has a fixed point 
$p \in P$, then we pick a continuous section 
$\sigma \: U \to P$ of a connected neighborhood $U$ of $x := q(p)$ with 
$\sigma(x) = p$ (the existence follows directly from the covering property). 
In addition, we may assume that 
$q^{-1}(U)$ is a disjoint union of open subsets of $P$ homeomorphic to 
$U$. Then $\sigma(U)$ is a connected component of $q^{-1}(U)$, 
and since $\phi \circ \sigma$ and $\sigma$ are two sections 
$U \to P$ mapping $x$ to the same point $p \in \sigma(U)$, it follows 
that $\sigma = \phi \circ \sigma$, so that $\phi$ fixes each point 
in $\sigma(U)$. We conclude that the set $P^\phi$ of points fixed by 
$\phi$ is an open closed subset of $P$, hence the assumption that 
$P$ is connected implies $\phi = \id_P$. Therefore  the 
group $\Delta$ acts simply transitively on each fiber, which implies 
that $q \: P \to X$ is a principal $\Delta$-bundle. 
\end{rem}

For the proof of the following theorem, we refer to \cite[Cor.\ 17.2]{Ste51}. 

\index{Long exact homotopy sequence of a principal bundle}
\begin{thm} [The long exact homotopy sequence of a principal bundle] 
\mlabel{homseq-princ} 
Let $q \: P \to B$ be a $G$-principal bundle, 
$y_0 \in P$ a base point, $b_0 := q(y_0)$ and identity 
$G$ with the fiber $F := q^{-1}(b_0) = y_0.G \cong G$. 
Then the maps 
$$ q_* \: \pi_k(P,F,y_0) \to \pi_k(B,b_0), 
\quad [f] \mapsto [q \circ f] $$
are bijective, so that we obtain 
connecting homomorphisms 
$$ \delta := \partial \circ (q_*)^{-1} \: \pi_k(B,b_0) \to 
\pi_{k-1}(F,y_0) \cong \pi_{k-1}(G). $$  
The so obtained sequence 
\begin{align*}
& \ldots \to \pi_k(P,y_0) \to \pi_k(B,b_0) \to \pi_{k-1}(G) \to \ldots \\
& \ldots \to \pi_1(P,y_0) \to \pi_1(B,b_0) \to \pi_0(G) \to 
\pi_0(P,p_0) \onto \pi_0(B,b_0) 
\end{align*}
is exact. The last two maps are not group homomorphisms;  
exactness in $\pi_0(P,y_0)$ 
means that the orbits of $\pi_0(G)$ in the set $\pi_0(P,y_0)$ are the 
fibers of the map $\pi_0(P,p_0) \to \pi_0(B,b_0)$, 
and the surjectivity of this map 
is the exactness in $\pi_0(B,b_0)$. 
\end{thm}

\begin{cor} \mlabel{lem10.1.2}
If $X$ is a semilocally simply connected arcwise connected space and 
$q_X \: \tilde X \to X$ is a universal covering of $X$, then 
$q_X$ induces isomorphisms 
$$ \pi_k(q_X) \: \pi_k(\tilde X) \to \pi_k(X), \quad k \geq 2. $$
\end{cor}

\begin{prf}  Pick $x_0 \in X$. 
We consider $q_X \: \tilde X \to X$ as a principal bundle for the 
discrete group $G := \pi_1(X,x_0)$ and apply the exact homotopy sequence 
(Theorem~\ref{homseq-princ}).
Since $G$ is discrete, we have $\pi_k(G) = \{\be\}$ for $k \geq 1$, and the assertion 
follows from the exactness of the homotopy sequence. 
\end{prf}

\begin{cor} \mlabel{cor10.1.2b}
Let $G$ be a topological group and $H \leq G$ a closed subgroup for which the 
natural map $q \: G \to G/H$ defines an $H$-principal bundle. 
Let $x_0 := \be H$ be the canonical base point in $G/H$. 
Then we have an exact sequence 
\begin{align*}
& \ldots \to \pi_k(G) \to \pi_k(G/H,x_0) \to \pi_{k-1}(H) \to \ldots \\
& \ldots \to \pi_1(G) \to \pi_1(G/H,x_0) \to \pi_0(H) \to \pi_0(G) \onto \pi_0(G/H,x_0). 
\end{align*}
\end{cor}

In view of \cite{MicE59}, the preceding corollary applies in 
particular for quotient morphisms $q \: G \to G/N$ of Banach--Lie groups. 

Sometimes one encounters situations where $q \: G \to G/H$ does not 
admit continuous local sections, but one would nevertheless 
like to use at least pieces of the long exact homotopy sequence. 
Such a situation arises from quotient morphisms of locally exponential 
Lie groups. 

\begin{itemize}
\item Let $G$ be a locally exponential Lie group,
\item   $N \subeq G$ be a closed normal locally exponential Lie subgroup, and
\item   assume that the quotient group $Q := G/N$ is also locally exponential. 
\end{itemize}

We write $q \: G\to Q$ for the quotient morphism and 
$j \: N \to G$ for the inclusion. 

We then have natural homomorphism 
\[  \pi_0(N) \to \pi_0(G) \to \pi_0(G/N) \quad\mbox{ and } \quad 
 \pi_1(N) \to \pi_1(G) \to \pi_1(G/N),\] 
but if $q$ has no local continuous sections, then it is not clear how to 
obtain a connecting map $\delta \: \pi_1(G/N) \to \pi_0(N)$. 
We now explain how this can nevertheless be done, and that the so 
obtained sequence 
\[  \pi_1(N) \to \pi_1(G) \to \pi_1(G/N) \sssmapright{\delta} 
 \pi_0(N) \to \pi_0(G) \to \pi_0(G/N) \to \1 \] 
is exact (Proposition~\ref{prop:homot-quot}). 

\begin{lem} Let $N \subeq G$ be a closed normal Lie subgroup of the 
locally exponential Lie group~$G$. Then the following 
assertions hold: 
\begin{description}
\item[\rm(a)] The natural map $G/N_0 \to G/N = Q, gN_0 \mapsto gN$ 
is a covering morphism. 
\item[\rm(b)] If $G$ is $1$-connected and $N$ is connected, then $G/N$ 
is $1$-connected. 
\item[\rm(c)] If $G$ is $1$-connected, then $\pi_1(G/N) \cong \pi_0(N)$.
\end{description}
\end{lem}

By Exercise~\ref{exer:10.1.5b}, (a) remains true for any 
Lie subgroup $N \subeq G$ 
because $N_0$ is a normal subgroup of~$N$.

\begin{prf} (a) The map $p \: G/N_0 \to G/N$ is a quotient map whose 
kernel is the discrete subgroup $N/N_0 \cong \pi_0(N)$ of $G/N_0$, 
hence $p$ is a covering. 

(b) Now we assume that $G$ is $1$-connected. We have to show that 
every loop $\gamma \: [0,1] \to G/N$ in $\be$ is contractible. 
To this end, we may w.l.o.g.\ assume that $\gamma(t) = \be$ for 
$\frac{1}{2} \leq t \leq 1$. 
Using an exponential chart, we can deform $\gamma$ to a loop $\eta$ for which 
there exists a subdivison 
\[ t_0 = 0 < t_1 < \ldots < t_k  < 1 \quad \mbox{ and } \quad 
x_1, \ldots, x_k \in \fq = \L(Q),\] 
such that 
\[ \eta(t) = \eta(t_j) \exp_Q((t-t_j)x_{j+1}) \quad \mbox{ for } \quad 
t_j \leq t < t_{j+1}, \quad j < k, \] 
and $\eta(t) = \be$ for $t \geq t_k$. 
Pick $y_j \in \g$ with $\L(q)x_j = y_j$. Then 
\[ \tilde\eta(t) = \tilde\eta(t_j) \exp_G((t-t_j)y_{j+1}) \quad \mbox{ for } \quad 
t_j \leq t < t_{j+1}, \quad j < k, \] 
defines a continuous lift $\tilde\eta \: [0,t_k] \to G$ of 
$\eta\res_{[0,t_k]}$. Since  $\eta(t_k) = \be$, $\tilde\eta(t_k) \in N$. 
As $N$ is connected, we can extend $\tilde\eta$ on $[t_k,1]$ by a path 
in $N$ from $\tilde\eta(t_k)$ to $\be$. 
Then $\eta = q \circ \tilde\eta$, and since the loop $\tilde\eta$ in 
$G$ is contractible, the same holds for $\eta$ in $Q$. 

(c) If $G$ is $1$-connected and $N$ is not connected, then 
(b) implies that $G/N_0$ is simply connected, and since 
$p \:  G/N_0 \to G/N$ is a covering for which 
$\pi_0(N)$ acts as a group of deck transformations on the fibers, 
it follows that $\pi_0(N) \cong \pi_1(G/N)$. 
\end{prf}

In the general context, we write 
$q_G \: \tilde G_0 \to G_0$ for the universal covering morphism of the identity 
component $G_0$, 
$\hat N \subeq \tilde G_0$ for the integral subgroup with Lie algebra 
$\fn$ and $N^\sharp := q_G^{-1}(N)$. Then 
$\hat N = (N^\sharp)_0$, so that the preceding lemma implies that 
$\tilde G_0/\hat N$ is simply connected. As 
\[ Q_0 = G_0/(N\cap G_0)  \cong \tilde G_0/N^\sharp,\] 
it follows that 
\begin{equation}
  \label{eq:pi1Q}
\pi_1(Q) = \pi_1(Q_0) \cong \pi_1(\tilde G_0/N^\sharp) 
\cong \pi_0(N^\sharp).
\end{equation}
Therefore the natural map $\pi_0(N^\sharp) \to \pi_0(N)$ leads to a connecting 
homomorphism 
\[  \delta \: \pi_1(Q) \to \pi_0(N).\] 
It maps a homotopy class $[\gamma]$ represented by a loop with a 
continuous lift 
$\tilde \gamma \: [0,1] \to G$ starting in $\be$ to the connected component 
of its endpoint $\tilde\gamma(1)$. 

\begin{prop} \mlabel{prop:homot-quot} 
Let $q \: G \to Q$ be a quotient morphism of 
locally exponential Lie groups with kernel~$N$. 
Then the sequence 
\[  \pi_1(N) \to \pi_1(G) \to \pi_1(G/N) \sssmapright{\delta} 
 \pi_0(N) \to \pi_0(G) \to \pi_0(G/N) \to \1 \] 
is exact. 
\end{prop}

\begin{prf} 
{\bf Exactness in $\pi_0(G)$ and $\pi_0(G/N)$:} 
The quotient group $G/N$ is a Lie group in which $G_0/N$ is an open connected 
subgroup. Therefore $Q_0 \cong \tilde G_0/N$ and 
\[ \pi_0(Q) = \pi_0(G/N) \cong G/(G_0 N).\] 
This implies that $\pi_0(q) \: \pi_0(G) \to \pi_0(Q)$ is surjective. 
We also see that a connected component $g G_0$ is mapped to $Q_0$ 
if and only if $g \in G_0 N$, which means that $[g] \in \pi_0(G)$ is contained 
in the image of $\pi_0(j)$. 

\nin {\bf Exactness in $\pi_0(N)$:} The image of $\delta$ in $\pi_0(N)$ 
is the image of $\pi_0(N^\sharp)$, which is the set of all connected 
components of $N$ contained in $G_0$. Therefore $\im(\delta)$ is the kernel 
of the map $\pi_0(N) \to \pi_0(G)$. 

\nin {\bf Exactness in $\pi_1(G/N)$:} That $[\gamma] \in \pi_1(G/N)$ 
is contained in $\ker\delta$ means that $\tilde\gamma(1) \in N_0$, and then 
$\tilde\gamma$ can be extended to a loop in $G$. Therefore 
it is contained in the image of $\pi_1(G)$. If, conversely, 
$[\gamma]$ is in the image of $\pi_1(G)$, i.e., ,
$\gamma = q \circ \tilde\gamma$ for a loop in $\tilde\gamma$ in $G$, 
then $\tilde\gamma(1) = \be \in N$. This implies that 
the image of $\pi_1(G)$ coincides with the kernel of $\delta$.  

\nin {\bf Exactness in $\pi_1(G)$:} From $q \circ j = \be$ it follows 
that $\im(\pi_1(j)) \subeq \ker(\pi_1(q))$. 
Suppose, conversely, that $[\gamma] \in \pi_1(G)$ is contained in 
the kernel of $\pi_1(q)$. Identifying 
$\pi_1(Q)$ with $\pi_0(N^\sharp)$ as in \eqref{eq:pi1Q} and 
$\pi_1(G)$ with the subgroup $\ker q_G$ of $\tilde G_0$, this means that 
\[ \ker(\pi_1(q)) 
\cong (\ker q_G) \cap (N^\sharp)_0 
\cong (\ker q_G) \cap \hat N,\] 
and this is the image of $\pi_1(N) \cong \ker q_N$ under the natural 
homomorphism $\tilde N_0 \to \tilde G_0$ whose range is $\hat N$. 
\end{prf}

\begin{cor} If $N$ is $1$-connected, 
then $\pi_1(G) \to \pi_1(G/N)$ is an isomorphism. 
\end{cor}

\subsection{Direct limits} 
\mlabel{subsec:14.1.2}

In this section, we collect some material on homotopy groups of direct limit 
spaces (cf.~Section~\ref{app-basic-DL} in Appendix~\ref{appA}). 
It is mainly be applied to direct limit Lie groups 
(Chapter~\ref{ch:dirlim}) and also to certain 
other Lie groups containing direct limit groups as dense subgroups. 
For the latter application we prove Palais' Theorem, one version of which 
states that, for a continuous map $f \: X \to Y$ between locally convex spaces 
with dense range and an open subset $O \subeq Y$, the induced map 
$f^{-1}(O) \to O$ is a weak homotopy equivalence.

\begin{lem} \mlabel{lem10.1.3}
Let $X$ be a Hausdorff space which carries the direct
 limit topology with respect to the increasing sequence $(X_n)_{n \in \N}$ 
of subspaces with $\bigcup_n X_n = X$, i.e., $O \subeq X$ is open if and only if 
all sets $X_n \cap O$ are open in $X_n$. Pick a base point $x_0 \in X_1$. Then 
\[  \pi_k(X,x_0)\cong \indlim \pi_k(X_n,x_0) \quad \mbox{ for every } \quad k \in \N_0.\]
\end{lem}

\begin{prf}  We claim that each compact subset $K \subeq X$ is contained in some $X_n$. 
If this is not so, then we pick $x_n \in K \setminus X_n$ for each $n \in \N$. 
The set $M := \{ x_n \: n~\in~\N\}$ satisfies 
$M \cap X_m \subeq \{ x_1, \ldots, x_{m-1} \}$. Therefore $M\cap X_m$ is closed 
for each $m$, so that $M$ is a closed subset of $X$. Thus $M \subeq K$ implies 
that $M$ is compact. The same argument applies to the subsets 
$M_m := \{ x_m, x_{m+1}, \ldots \} \subeq M$. Now $\bigcap_{m \in \N}
M_m \not=\eset$ follows from the compactness of the 
sets $M_m$. But $M_{m+1} \cap X_m = \eset$ implies that 
$$ \bigcap_{m \in \N} M_m \subeq X \setminus \bigcup_{m \in \N} X_m = \eset. 
$$
This contradiction shows that there exists an $n \in \N$ with $K
\subeq X_n$.

Now let $\gamma \: \bS^k \to X$ be a continuous map. Then 
$\gamma(\bS^k)$ is a compact subset of $X$, hence contained in some
$X_n$, and since $X_n \into X$ is an embedding, the corestriction 
$\gamma\res^{X_n} \: \bS^k \to X_n$ is continuous. Therefore the natural map 
$$ \indlim \pi_k(X_n,x_0)  \to \pi_k(X,x_0) $$
is surjective. To verify injectivity, we apply the same argument
to the range of a homotopy of two continuous maps 
$\gamma_1 \: \bS^k \to X_{n_1}$ and $\gamma_2 \: \bS^k \to X_{n_2}$. 
We find that there exists $n_3 > \max(n_1, n_2)$ such that 
$\pi_k(\phi_{n_3, n_1})([\gamma_1]) = \pi_k(\phi_{n_3,
n_2})([\gamma_2])$, 
where $\phi_{n_3, n_1} \: X_1 \to X_3$ and 
$\phi_{n_3, n_2} \: X_2 \to X_3$ are the embeddings. Therefore the 
images of $[\gamma_1]$ and $[\gamma_2]$ in $\pi_k(X,x_0)$ coincide. 
\end{prf}

\begin{lem} \mlabel{lem-pal0} Let $E$ be a real vector space endowed with the 
finite topology, i.e., the direct
limit topology with respect to its finite-dimensional subspaces. Then
the following assertions hold: 
\begin{enumerate}
\item[\rm(1)] Each linearly independent subset of $E$ is closed and discrete. 
\item[\rm(2)] Each compact subset of $E$ is contained a finite-dimensional subspace. 
\item[\rm(3)] For each subset $U \subeq E$ and $u_0 \in U$, we have 
$$ \pi_k(U,u_0) \cong \indlim \pi_k(U \cap F,u_0), $$
where the limit is taken over the directed set of all finite-dimensional
subspaces $F \subeq E$ containing $u_0$. 
\item[\rm(4)] If $U \subeq E$ is a subset for which the intersection
with all finite-dimensional subspaces are open, then the subspace
topology on $U$ coincides with the direct limit topology with respect
to the sets $U \cap F$, $F \subeq E$ a finite-dimensional subspace. 
\end{enumerate}
\end{lem}

\begin{prf}  (1) (\cite[Lemma 10.2]{Pa66}) Let $S \subeq E$ be a linearly
independent subset. Then, for each finite-dimensional subspace $F
\subeq E$, the intersection $S \cap F$ is closed, so that $S$ is
closed in~$E$. The same argument implies that each subset of $S$ is
also closed in $E$. It follows in particular that $S$ is a discrete
topological space. 

 (2) (\cite[Lemma~10.3]{Pa66}) Let $K \subeq E$ be a compact subset and $S \subeq K$ a
maximal linearly independent subset. Then $K \subeq \Spann S$. In view of (1), 
$S$ is closed,
hence compact. On the other hand, $S$ is discrete and therefore
finite. 

 (3) Let $Y$ be a compact space with base point $y_0$ and 
$f \: Y \to U$ a
continuous map with $f(y_0) = u_0$. Then $f(Y)$ is a compact subset of $E$, hence
by (2) contained in a finite-dimensional subspace $F$, and we clearly have 
$u_0 = f(y_0) \in F$. 

For $Y = \bS^k$, it follows that the natural homomorphism 
$$ \eta \: \indlim\pi_k(U \cap F,u_0) \to \pi_k(U,u_0) $$
is surjective. To see that it is also injective, suppose that 
$f \: \bS^k \to U \cap F$ is a continuous map which in $U$ is
homotopic to the constant map $\bS^k \to \{u_0\}$. 
Let $H \: [0,1] \times \bS^k \to U$ be a homotopy with 
$H(0,x)= f(x)$ and $H(1,x) = u_0$. Then the compact set 
$\im(H)$ is contained in a
finite-dimensional subspace $F \subeq E$, and therefore the homotopy
class of $f$ in $\pi_k(U \cap F, u_0)$ is trivial. This implies that
$\eta$ is injective. 

 (4) Our assumption implies that $U$ is open in the finite topology. 
Hence a subset $V \subeq U$ is open in the finite topology 
if and only if it is open in $E$, which means that 
all intersections $V \cap F$ with a finite-dimensional subspace $F \subeq E$,
are open in $F$. This in turn shows that $U$ carries the direct limit topology with 
respect to the subspaces $U \cap F$. 
\end{prf}

\subsection{Palais' Theorem on weak homotopy equivalences} 
\mlabel{subsec:14.1.3b}

In this subsection, we turn to the circle of ideas related to Palais' Theorem. 
Throughout $V$ denotes a locally convex space, 
$E \subeq V$ is a dense subspace endowed with its finite topology, i.e., the direct 
limit topology with respect to all finite-dimensional subspaces of $E$, and 
$O \subeq V$ is an open subset. We write $O_E := O \cap E$ for its intersection 
with $E$, endowed with the topology inherited from $E$. One of our goals in this section is to show that the inclusion map $O_E \to O$ is a weak homotopy equivalence.  

\begin{lem} \mlabel{lem-pal1} Let $O \subeq V$ be an open subset of the locally convex 
space~$V$ and $K \subeq O$ compact. If $W$ is a sufficiently small $0$-neighborhood in 
$V$ and $k_1, \ldots, k_n \in K$ are such that 
$$ K \cap \bigcap_{i = 1}^n (k_i + W) \not= \eset, \quad 
\mbox{ then } \quad \conv\Big(\bigcup_{i = 1}^n (k_i + W)\Big) \subeq O.$$ 
\end{lem}

\begin{prf} For each $k \in K$, let $U_k$ be a $0$-neighborhood in $V$ with 
$k + U_k \subeq O$ and $W_k$ be a neighborhood of zero with $W_k + W_k \subeq U_k$. 
Then there exist finitely many elements $\ell_1, \ldots, \ell_m \in K$ with 
$K \subeq \bigcup_{j = 1}^m \ell_i + W_{\ell_i}.$
Let $W_1$ be a convex $0$-neighborhood in $V$ contained in 
$\bigcap_{j = 1}^m W_{\ell_j}$. Then we have 
$$ \ell_i + W_{\ell_i} + W_1 \subeq \ell_i + W_{\ell_i} + W_{\ell_i} \subeq \ell_i + U_{\ell_i} \subeq O $$
for each $i$ and therefore $K + W_1 \subeq O$. 

Now let $W$ be a convex $0$-neighborhood 
with $W + (-W) \subeq W_1$ and suppose that there exist elements 
$k, k_1,\ldots, k_n \in K$ with $k \in \bigcap_{i = 1}^n (k_i + W)$. 
Since each of the sets $k_i + W$ is convex, we have 
$$ \conv\Big(\bigcup_{i = 1}^n k_i + W\Big) 
= \Big\{ \sum_{i =1}^n \lambda_i (k_i + w_i) \: w_i \in W, \lambda_i \geq 0, 
\sum_i \lambda_i = 1
\Big\}. $$
In view of $k + W_1 \subeq O$, it suffices to show that 
$$ \Big(\sum_i \lambda_i (k_i + w_i)\Big) - k 
= \sum_i \lambda_i (k_i + w_i- k)  \in W_1 $$
for $w_i \in W, \lambda_i \geq 0$ and $\sum_i \lambda_i= 1$. 
Since $W_1$ is convex, it even suffices to show that $k_i + w_i - k \in W_1$. 
In view of $k - k_i \in W$, this follows from $W - W \subeq W_1$. 
\end{prf}

In the following, we call a family 
$(h_t)_{0 \leq t \leq 1}$ of continuous maps $X \to Y$ a 
\index{ deformation, of a map} 
{\it deformation} if the map 
$h \: [0,1] \times X \to Y,  (t,x) \mapsto h_t(x)$
is continuous. 

\begin{lem} \mlabel{lem-pal2} 
Let $K \subeq O$ be a compact subset, $E \subeq V$ a dense subspace and $L \subeq K \cap E$ 
be a compact subset contained in a 
finite-dimensional subspace. Then there is a deformation $h_t \: K \to O$ of $K$ in $O$ 
such that $h_t\res_L = \id_L$ and $\im(h_1)$ is contained in a finite-dimensional 
subspace of 
$E$. 
\end{lem}

\begin{prf} Let $E_1 \subeq E$ be a finite-dimensional subspace containing $L$. 
By the Tietze Extension Theorem (\cite[Thm.~I.10.4]{Bre93}) 
the identity map $L \to E_1$ extends to a continuous map 
$g \: K \to E_1$. Let $W$ be as in Lemma~\ref{lem-pal1} and put 
$$ U_0 := \{ k \in K \: g(k) - k \in W \}. $$
Then $U_0$ is a neighborhood of $L$ in $K$. Let $U_1, \ldots, U_n$ be an open cover 
of $K \setminus U_0$ by sets disjoint from $L$ and of the form 
$U_i = K \cap (k_i + W_i)$, where $W_i \subeq W$ is a $0$-neighborhood in $V$. 
Let $\phi_0, \ldots, \phi_n$ be a partition of unity for $K$ with 
$\supp(\phi_i) \subeq U_i$ and define 
$$ h_t(k) := (1 - t)k + t \Big( \phi_0(k) g(k) + \sum_{i =1 }^n \phi_i(k) e_i\Big), $$
where $e_i \in (k_i + W_i) \cap E \subeq (k_i + W) \cap E$ (recall that $E$ is dense in $V$). 
Then clearly $h_t$ is a deformation of $K$ in $V$. 
If $k \in L$, then $\phi_i(k) = 0$ for $i \geq 1$, since $U_i \cap L = \eset$, 
so $\phi_0(k)= 1$ and $h_t(k) = (1-t)k + t g(k) = (1-t)k + t k = k$. 
Let $E_2 := \Spann(e_1,\ldots, e_n) + E_1$. Then $h_1$ maps $K$ continuously into 
$E_2$, and since $E_2$ is a topological subspace of $E$, it maps $K$ continuously into $E$. 

It remains to show that $h_t(K) \subeq O$. For a given $k\in K$, 
arrange the indices so that 
$\phi_i(k) > 0$ for $i = 1,\ldots, m$ and $\phi_i(k) = 0$ for $i > m$. Then 
$$ k \in (k + W) \cap \bigcap_{i=1}^m U_i \subeq (k + W) \cap \bigcap_{i = 1}^m (k_i + W). $$
Now 
$(1 -t) + t \sum_{i = 0}^m \phi_i(k) = 1$, so by the preceding lemma, it will suffice to 
show that $k$ and the $e_i$ belong to 
$(k + W) \cup \bigcup_{i = 1}^m (k_i + W)$, and that if 
$\phi_0(k) > 0$, so does $g(k)$. But $k \in k + W$, $e_i \in k_i + W$ and if $\phi_0(k) > 0$, 
then $k \in U_0$ and so $g(k) \in k + W$. 
\end{prf}

\begin{lem} \mlabel{lem-pal3} Let $A$ be a closed subspace of a compact space $X$ and 
let $f_0 \: X \to O$ be a continuous map such that $f_0\res_A$ is a continuous map 
of $A$ into~$O_E$. 
Then there is a homotopy $f_t \: X \to O$ of $f_0$ such that 
$f_1$ is a continuous map of $X$ into $O_E$ and 
$f_t \res_A = f_0\res_A$ for $0 \leq t \leq 1$. 
\end{lem}

\begin{prf} In Lemma \ref{lem-pal2}, take $K := f_0(X)$ and $L := f_0(A)$ and put 
$f_t := h_t \circ f_0$.   
\end{prf}

\begin{defn} \mlabel{def:weak-homo-equiv} 
Let $(X,x_0)$ and $(Y,y_0)$ be pointed topological spaces. A~map 
$f \in C_*(X,Y)= \{ h \in C(X,Y) \: h(x_0) = y_0\}$ 
is called a 
\index{homotopy equivalence!weak} 
\index{homotopy equivalence} 
{\it weak homotopy equivalence} if all induced maps 
$\pi_k(f) \: \pi_k(X,x_0) \to \pi_k(Y,y_0)$, $[\gamma] \mapsto [f \circ \gamma]$, are bijections. 

A map $f \in C_*(X,Y)$ is called a {\it homotopy equivalence} if 
there exists a $g \in C_*(Y,X)$ such that $f\circ g$, resp., $g\circ f$ 
are homotopic to 
$\id_Y$, resp., $\id_X$ in $C_*(Y,Y)$, resp., $C_*(X,X)$. 
The functoriality of the higher homotopy groups immediately 
implies that a homotopy equivalence is a weak homotopy equivalence. 
\end{defn}

\index{Palais' First Theorem on weak homotopy equivalences} 
\begin{thm} [Palais' First Theorem on weak homotopy equivalences] \mlabel{thm-pal1} 
    Let $V$ be a locally convex space and 
    $E \subeq V$ a dense linear subspace. We endow $E$ with the finite topology. For each open subset 
    $O \subeq V$, the continuous map $O_E = O \cap E\to O$ 
is a weak homotopy equivalence if 
    $O_E$ is considered as a topological subspace of $E$. 
\end{thm}

\begin{prf} Let $f_0 \: \bS^n \to O$ be a continuous map and 
$f_t \: \bS^n \to O$ a homotopy of $f_0$ with the continuous map 
$f_1 \: \bS^n \to O_E$ (Lemma~\ref{lem-pal3}). Then the inclusion map 
$i \: O_E \to O$ satisfies 
$\pi_n(i)[f_1] = [f_0],$
showing that $\pi_n(i)$ is surjective. 

Now suppose $h_0, h_1  \: \bS^n \to O_E$ are continuous maps for which 
$i \circ h_0$ and $i \circ h_1$ are homotopic. 
Then there exists a continuous map 
$H \: [0,1] \times \bS^n \to O$ with 
$H_t = i \circ h_t$ for $t =0,1$. Again by Lemma~\ref{lem-pal3}, 
there exists a homotopy 
$F_s \: [0,1] \times \bS^n \to O$, $0 \leq s \leq 1$, with $F_0 = H$, $F_1 
\: [0,1] \times \bS^n \to O_E$ continuous and 
$F_s \res_{\{t\} \times \bS^n} = h_t$ for $t = 0,1$ and all $s$. 
Then $F_1$ is a homotopy of $h_0$ and $h_1$ as continuous maps 
$\bS^n \to O_E$. Therefore $\pi_n(i)$ is also injective. 
\end{prf}

The following theorem is quite useful to 
calculate homotopy groups: 

\index{Palais' Second Theorem on weak homotopy equivalences}
\begin{thm} [Palais' Second Theorem on weak homotopy equivalences]
\mlabel{thm-pal2}
Let $V_1$ and $V_2$ be locally convex spaces and 
$f \: V_1 \to V_2$ a continuous linear map with dense range. Let 
$O \subeq V_2$ be an open subset and put $\tilde O := f^{-1}(O)$ and 
$\tilde f := f\res_{\tilde O}$. 
Then $\tilde f \: \tilde O \to O$ is a weak homotopy equivalence. 
\end{thm}

\begin{prf} Let $W := \ker f$, $\pi \: V_1 \to V_1/W$  the quotient map, 
$f = g \circ \pi$ the canonical factorization of $f$, and $O^* := g^{-1}(O) \subeq V_1/W$. 
Then it will suffice to show that 
$$\pi\res_{\tilde O} \: \tilde O \to O^*\quad \mbox{ and } \quad 
g\res_{O^*} \: O^* \to O $$
are weak homotopy equivalences. Since $g$ is injective and $\pi$ is a quotient 
map, it will suffice to prove the theorem in the following two cases. 

\nin {\bf Case 1. $f$ is injective:} Without loss of generality we can assume that 
$V_1$ is a dense linear subspace of $V_2$ with a finer topology and  that 
$f$ is the inclusion map. So $\tilde O = O \cap V_1$ as a subset of $V_1$. 
Let $E := V_1$ with the finite topology and let 
$O_E := O \cap E$ as a subspace of $E$. Let $i \: O_E \to \tilde O$ and 
$j \: O_E \to O$ be the inclusion maps which by Theorem~\ref{thm-pal1} are 
weak homotopy equivalences. Then $f \circ i = j$ implies that 
$f$ is a weak homotopy equivalence. 

\nin {\bf Case 2. $f$ is surjective:} Let $E_j := V_j$, $j =1,2$, be 
endowed with the finite 
topology, $O_E:= O$ as a subspace of $V_2$, $\tilde O_E := \tilde O$ as a subspace of 
$E_1$ and define $g \: \tilde O_E \to O_E$ by $g := f\res_{\tilde O_E}$. Since a linear 
map between vector spaces, each endowed with the finite topology, is clearly 
continuous, $g$ is continuous, and $\tilde f \circ i = j \circ g$, where 
the identity maps $i$ and $j$ are weak homotopy equivalences by Theorem~\ref{thm-pal1}. 
Thus it suffices to prove that $g$ is a weak homotopy equivalence. 

Let $h \: \bS^n \to O_E$ be a continuous map. Then $\im(h)$ is contained in a 
finite-dimensional subspace $V$ of $E_2$ (Lemma~\ref{lem-pal0}). 
Let $e_1, \ldots, e_n$ be a basis of~$V$. Since $f$ is surjective, we can find 
elements $\tilde e_i \in E_1$ with $f(\tilde e_i) = e_i$. Let 
$h(q) := \sum_{i = 1}^n h_i(q) e_i$ be the coordinate representation of the 
continuous map $h \: \bS^n \to V$ and define $\tilde h \: \bS^n \to E_1$ by 
$$ \tilde h(q) := \sum_{i = 1}^n h_i(q) \tilde e_i. $$
Then $g \circ \tilde h = h$ and we obtain $[h] = \pi_n(g) [\tilde h]$, 
so that $\pi_n(g)$ is surjective. 

Next let $k_0, k_1 \: \bS^n \to \tilde O_E$ be continuous maps with 
$g \circ k_i$ homotopic and $K \: [0,1] \times \bS^n \to O_E$ be 
a continuous map with $K_t = g \circ k_t$ for $t = 0,1$. 
Then $\im(k_0) \cup \im(k_1)$ is contained in a finite-dimensional subspace 
$Y$ of $E_1$ of which we may assume, in addition, that $f(Y)$ contains the image 
of $K$ (Lemma~\ref{lem-pal0}(2)). Let 
$$\tilde O_Y := \tilde O \cap Y = f^{-1}(O) \cap Y = (f\res_Y)^{-1}(O). $$
With $O_Y := f(Y) \cap O$ we then obtain 
$$ \tilde O_Y \cong O_Y \times \ker (f\res_Y), $$
which is homotopy equivalent to $O_Y$. Hence the maps 
$\bS^n \to \tilde O_Y \subeq \tilde O_E$ obtained by corestriction from 
$k_0$ and $k_1$ are homotopic as maps to $\tilde O_Y$. 
Therefore $k_0$ and $k_1$ are homotopic. This proves that 
$\pi_n(g)$ is injective. 
\end{prf}

\begin{defn} \mlabel{def10.1.12} 
Let $X$ be a vector space which is the union of an increasing sequence of linear  
subspaces $X_n$, $n \in \N$. Suppose that $\tau_n$ is a topology on $X_n$ turning it 
into a locally convex space such that 
$\tau_{n+1} \res_{X_n} = \tau_n$ for each $n \in \N$. Let 
$\tau$ denote the finest locally convex topology on $X$ for which all the inclusion maps 
$X_n \into X$ are continuous. This means that a seminorm $q \: X \to \R$ is continuous 
if and only if all restrictions $q\res_{X_n}$ are continuous. 
We call $(X,\tau)$ the 
\index{inductive limit!strict, of l.c. spaces}
{\it strict inductive limit} of the locally convex spaces 
$(X_n,\tau_n)$. 
\end{defn}

\begin{lem} \mlabel{ext-lem} Suppose that the locally convex space $X$ 
is the strict inductive limit of the locally convex spaces 
$X_n$. Then the following assertions hold: 
\begin{description}
\item[\rm(a)] Every continuous seminorm $p_n$ on $X_n$ is the restriction of a 
continuous seminorm $p_{n+1}$ on $X_{n+1}$, hence also of a continuous seminorm 
on $X$. 
\item[\rm(b)] If each $X_n$ is closed in $X_{n+1}$, 
$y \in X_{n+1} \setminus X_n$ and $c \in \R$, 
then each continuous seminorm $p_n$ on $X_n$ 
is the restriction of a continuous seminorm $p_{n+1}$ on $X_{n+1}$ with 
$p_{n+1}(y) \geq  c$. 
\end{description}
\end{lem}

\begin{prf} (cf.\ \cite[Prop.~1.10.2]{He89}) (a) Let $p_n$ be a continuous seminorm on $X_n$ 
and $U_0 := \{ x \in X_n \: p_n(x) \leq 1\}$. As a $0$-neighborhood in $X_n$, 
$U_0$ is the intersection with $X_n$ of a $0$-neighborhood in $X_{n+1}$, 
so that there exists a continuous seminorm $p'$ on $X_{n+1}$ with 
$p_n \leq p'\res_{X_n}$. 

The set 
$$ A := \{ \alpha x + (1- \alpha)x' 
\: x \in X_n, p_n(x) \leq 1; x' \in X_{n+1}, p'(x') \leq 1, 0 \leq \alpha \leq 1\} $$
is convex, balanced and absorbing, so that its gauge 
$$ p_{n+1}(x) := \inf \{ t > 0 \: x \in t A \} $$
is a seminorm (Proposition~\ref{Minkowsk}). 
Since $A$ also is a neighborhood of the origin in $X_{n+1},$ the seminorm 
$p_{n+1}$ is continuous. Further 
$X_n \cap A = U_0$ implies $p_{n+1} \res_{X_n} = p_n$. 

Iterating this construction, we get a sequence $(p_k)_{k \leq n}$ of continuous seminorms 
$p_k$ on $X_k, k \geq n$, with $p_{k+1}\res_{X_{k}} = p_k$. Setting 
$q(x) := p_k(x)$ for $x \in X_{k}$ now defines a continuous seminorm on $X$ 
whose restriction to $X_n$ is $p_n$. 

(b) Since $y + X_n$ is a closed subset of $X_{n+1}$ that does not contain the origin, 
the seminorm $p'$ in the proof of (a) can be chosen so that $p' \geq c$ on 
$y + X_n$. Then $y \in tA$ for $y = \alpha x + (1- \alpha) x'$ as above 
implies 
$$t \geq t(1- \alpha) p(x') = p'((1-\alpha)t x') = p'(y - \alpha t x) \geq c$$ 
and therefore $p_{n+1}(y) \geq c$. 
\end{prf}

\begin{lem} \mlabel{lem-herve} 
If the locally convex space $X$ is the strict inductive limit of the spaces $X_n$, then each bounded subset 
of $X$ is contained in some $X_n$ and bounded as a subset of $X_n$. 
\end{lem}

\begin{prf} (cf.~\cite[Prop.\ 1.10.3]{He89}) Let $B \subeq X$ be bounded and assume that 
$B$ is not contained in any $X_n$. Then there  exists an increasing sequence 
$(n_k)_{k \in \N}$ of integers and 
$y_k \in B \cap (X_{n_{k+1}} \setminus X_{n_k})$ for all $k \in \N$. 
Lemma~\ref{ext-lem}(b) yields a sequence of continuous 
seminorms $p_n \: X \to \R$ with $p_{n+1} \res_{X_{n+1}} = p_n$ and 
$p_{n_{k+1}}(y_k) > k$ for all $k \in \N$. By $q(x) := p_n(x)$ for $x \in X_n$, 
we then obtain a continuous seminorm $q$ on $X$ which is not bounded on $B$, 
contradicting the boundedness of $B$ in $X$. 

Hence there is some $n$ with $B \subeq X_n$. That $B$ is bounded as a subset of 
$X_n$ follows from the fact that each continuous seminorm on $X_n$ extends 
to a continuous seminorm on $X$ (Lemma~\ref{ext-lem}). 
\end{prf}

\begin{prop} \mlabel{prop10.1.13}
Assume that 
$X$ is the strict inductive limit of the spaces $X_n$ and that 
each $X_n$ is closed in $X_{n+1}$. Let $O \subeq X$ be an open subset 
and write $O_\infty$ for the set $O$, endowed with the direct limit topology of the subsets $O_n := O \cap X_n$ of $X_n$. 
Then the continuous map $O_\infty \to O$ is a weak homotopy 
equivalence. 
\end{prop}

\begin{prf}  It follows directly from the definitions that the map 
$\eta \: O_\infty \to O$ is continuous. We claim that each continuous map 
$f \: K \to O$, $K$ a compact space, factors through a continuous map 
$\tilde f \: K \to O_\infty$, i.e., $f$ is also continuous as a map to $O_\infty$. 

According to Lemma~\ref{lem-herve}, each bounded subset of $X$ is contained in some 
$X_n$. This implies in particular that each compact subset of $O$ is contained in some 
$O_n$. Therefore we find for each continuous map $f \: K \to O$ an $n \in \N$ with 
$f(K) \subeq O_n$. Then $f \: K \to O_n$ is continuous 
because the inclusion $O_n \to O$ is an embedding (Lemma~\ref{ext-lem}(a)), 
and this implies that the map $\tilde f \: K \to O_\infty$ is also continuous. 

Since the preceding observation 
applies to homotopies 
$h \: [0,1] \times K \to O$ between two maps $K \to O$, we 
conclude that the map $\eta$ induces a bijection of pointed homotopy
classes $[K,O_\infty]_* \to [K,O]_*$, 
and hence in particular that $\eta$ is a weak homotopy equivalence. 
\end{prf}

\index{Palais' Third Theorem on weak homotopy equivalences} 
\begin{thm} [Palais' Third Theorem on weak homotopy equivalences] \mlabel{thm8.6} 
Let $X$ be a locally convex space and $(X_n)_{n \in \N}$ be an
increasing sequence of closed subspaces whose union $X_\infty :=
\bigcup_{n \in \N} X_n$ is dense in $X$. Let $O \subeq X$ be an open subset 
and write $O_\infty$ for the set $O \cap X_\infty$, endowed with the
direct limit topology of the open subsets 
$O_n := O \cap X_n$ of~$X_n$. 
Then the continuous map $O_\infty \to O$ is a weak homotopy 
equivalence. 
\end{thm}

\begin{prf}  We endow $X_\infty$ with the 
strict inductive limit topology with respect to the sequence 
$(X_n)_{n \in \N}$ (Definition~\ref{def10.1.12}). 
Let $O_0 := O \cap X_\infty$ (considered as a subspace of $X_\infty$). 
Then we have continuous maps 
$$ f_0 \: O_0 \to O, \quad 
f_{\infty,0} \: O_\infty \to O_0, \quad \hbox{ and } \quad 
f_{\infty} \: O_\infty \to O $$ 
with $f_\infty = f_0 \circ f_{\infty,0}$. 

In view of Theorem~\ref{thm-pal2}, the map $f_0$ is a weak homotopy 
equivalence and Proposition~\ref{prop10.1.13} implies that $f_{\infty,0}$ is a weak
homotopy equivalence, so that $f_\infty = f_0 \circ f_{\infty,0}$ also
is a weak homotopy equivalence. 
\end{prf}

The following criterion for contractibility can be found in 
\cite[Cor.\ to Th.\ 15]{Pa66}. We shall use it in the proof of 
Kuiper's Theorem~\ref{kuiper-thm} below. 
\begin{thm}\mlabel{pal-contr} 
Let $(M,m_0)$ be a pointed metrizable manifold. 
Then $M$ is contractible if and only if 
all homotopy groups $\pi_n(M,m_0)$, $n \in \N_0$, are trivial. 
\end{thm}

\section{Finite-dimensional Lie groups} \mlabel{sec:10.2}

In this section, we briefly recall some facts on the topology of
finite-di\-men\-sio\-nal 
Lie groups. Since the proofs of these facts require quite elaborate structure 
theory, we refer to the literature for detailed arguments. 

Let $G$ be a connected finite-dimensional Lie group with 
finitely many connected components and $K \subeq G$ a maximal compact 
subgroup. Then there exists a diffeomorphism $G \to K \times \R^d$ 
for some $d \in \N_0$ (\cite[Thm.~14.3.11]{HiNe12}). 
This implies in particular that the inclusion map 
$K\into G$ is a homotopy equivalence, hence induces isomorphisms 
$\pi_k(K) \to \pi_k(G)$ for each $k \in \N_0$. 
This reduces all questions on the topology 
of finite-dimensional Lie groups to compact groups. 

\begin{rem} \mlabel{rem1.3} 
(a) We recall some results on the homotopy groups of
a compact Lie group~$K$. 
Cartan's Theorem asserts that 
$\pi_2(K)$ vanishes (\cite[Th.\ 3.7]{Mim95}) and Bott's Theorem states 
that $\pi_3(K) \cong \Z$ for every compact connected simple Lie group~$K$ 
(\cite[Th.\ 3.9]{Mim95}). A generator of $\pi_3(K)$ can
be obtained from a homomorphism $\eta \: \SU_2(\C) \cong \bS^3 \to K$ 
(cf.\ Exercise~\ref{exer:10.2.1}). 
More precisely, let $\alpha$ be a long root in the root system $\Delta_\fk$ of $\fk$ 
and 
$\fk(\alpha) \subeq \fk$ the corresponding $\su_2(\C)$-subalgebra. 
Then the corresponding homomorphic inclusion $\SU_2(\C) \cong
\bS^3 \to K$ represents a generator of $\pi_3(K)$ (\cite{Bo58}). 

In \cite[pp.\ 969-970]{Mim95} one finds a table with $\pi_k(K)$ up to
$k = 15$, showing that 
\[  \pi_4(K) \cong 
\begin{cases}
\Z_2 \oplus \Z_2 & \text{for } K = \SO_4(\R) \\
\Z_2 & \text{for }  K = \U_{n}(\bH), \SU_{2}(\C), \SO_{3}(\R), \SO_{5}(\R) \\ 
\1  & \text{for } K = \SU_n(\C), n \geq 3, \text{ and } \SO_n(\R), n \geq 6 \\ 
\1  & \text{for } K = G_2, F_4, E_6, E_7, E_8.  
\end{cases}\] 
\[  \pi_5(K) \cong 
\begin{cases}
\Z_2 \oplus \Z_2 & \text{for } K = \SO_4(\R) \\ 
\Z_2 & \text{for } K = \U_n(\bH), \SU_2(\C), \SO_3(\R), \SO_5(\R) \\ 
\Z  & \text{for } K = \SU_n(\C), n \geq 3, \text{ and } \SO_6(\C) \\ 
\1  & \text{for} K= \SO_n(\R), n \geq 7,\quad  G_2, F_4, E_6, E_7, E_8. 
\end{cases}
\] 

(b) Let $K$ be a connected compact Lie group, 
$K_1, \ldots, K_m$ be 
the connected simple normal subgroups of $K$, and $Z(K)$ be its center. 
Then the multiplication map 
$$ Z(K)_0 \times K_1 \times \ldots \times K_m \to K $$
has finite central kernel $\Gamma$, hence is a covering map. 
Since $Z(K)_0$ is a torus, its universal covering group 
$\tilde Z(K)_0 \cong \fz(\fk)$ is contractible, so that 
all its homotopy groups of degree $\geq 2$ vanish. We therefore 
obtain with Lemma~\ref{lem10.1.2} 
\[  \pi_k(K)  \cong \prod_{j = 1}^m \pi_k(K_j) \quad \mbox{ for }\quad 
k > 1.\] 
Corollary~\ref{cor10.1.2b} further yields a short exact sequence 
\[ \1\to \pi_1(Z(K)) \times \prod_{j = 1}^m \pi_1(K_j) \to \pi_1(K) \to \Gamma \to \1. \]
\end{rem} 

\begin{rem} \mlabel{rem1.4} 
A fundamental result in topology states that the 
spheres $\bS^d$ carry a Lie group structure if and only if $d \in \{0,1,3\}$ 
(cf.\ \cite{Hz91}): 
\[ \bS^0 \cong \OO_1(\R), \qquad \bS^1 \cong \T \cong \SO_2(\R) \quad \mbox{ and } \quad 
\bS^3 \cong \SU_2(\C).\] 
That $d \in \{0,1,3\}$ is necessary for the existence of a Lie group 
structure on $\bS^d$ can already be derived from Remark~\ref{rem1.3} 
as follows. Let $K$ be a compact Lie group homeomorphic to~$\bS^d$. 
If $K$ is abelian, then it is a torus (Example~\ref{exs:ab-quot}(a)) 
and therefore not simply connected, so that $d = 1$. 
If $K$ is not abelian, then it contains a simple compact subgroup, 
so that $\pi_3(K)$ is non-zero and thus $d \leq 3$. The number 
$d = 2$ is ruled out because $\pi_2(K)$ vanishes (Remark~\ref{rem1.3}). 

For finite-dimensional Lie groups $G$, the Lie group structures on $\bS^1$ and $\bS^3$ 
actually leads to the interesting fact that, for $d = 1,3$,  
every homotopy class in $[\bS^d,K]$ can be represented by a 
group homomorphism. For $d = 3$ this follows from Remark~\ref{rem1.3}, 
and for $d = 1$ it follows from 
the fact that, for a maximal torus $T \subeq K$, the homomorphism 
\begin{equation}
  \label{eq:maxtorus}
\Hom(\T,T) \cong \pi_1(T) \to \pi_1(K) 
\end{equation}
is surjective (\cite[Cor.~14.2.10]{HiNe12}).  
Since the long exact homotopy sequence of the $T$-principal bundle 
$K$ over $K/T$ leads to an exact sequence 
\[ \1 = \pi_2(K) \to \pi_2(K/T) \to \pi_1(T) \to \pi_1(K) \to \pi_1(K/T) \to \pi_0(T) 
= \1 \]
(Corollary~\ref{cor10.1.2b}),  
the surjectivity of~\eqref{eq:maxtorus} implies that 
$\pi_1(K/T)= \1$. 

In Example~\ref{ex:VII.9-mn03} below we shall see that, for infinite-dimensional 
groups, elements of $\pi_1(G)$ are in general not represented by homomorphisms~$\T \to K$. 
\end{rem}

  \begin{rem} \mlabel{rem1.5} 
For a topological group $G$ and $k \geq 1$, the groups 
$\pi_k(G)$ are abelian (Remark~\ref{rem:10.1.3}). The groups 
$$\pi_k^\Q(G) := \Q \otimes \pi_k(G)$$ are called 
\index{homotopy group!rational} 
the {\it rational homotopy groups of $G$}. For most purposes, 
including applications to the period maps 
arising for abelian extensions \cite{Ne02a, Ne04a} (see also 
Section~\ref{sec:11.4}), it suffices to know the rational homotopy groups because 
each homomorphism from $\pi_k(G)$ to a rational vector space factors through the natural 
map $\eta_k \: \pi_k(G) \to \pi_k^\Q(G)$, killing the torsion subgroup of $\pi_k(G)$.

We have seen above that a finite-dimensional connected Lie group is homotopy 
 equivalent to a compact connected Lie group, hence, up to a finite covering 
 to a product of a torus and finitely many compact simple Lie groups. 
For a simply connected simple compact Lie group its rational homotopy 
groups are the same as those of a product of odd-dimensional spheres 
whose dimensions can be 
computed from the corresponding root system (Remark~\ref{rem1.3}). 
The rational homotopy groups of spheres are 
$$ \pi_k^{\Q}(\bS^{2d+1}) \cong 
\left\{   \begin{array}{cl}
\Q & \mbox{for $k = 2d + 1$} \\ 
\0 & \mbox{otherwise }  \\ 
  \end{array}
\right. 
\ \ \ \hbox{ and } \ \ \ 
\pi_k^{\Q}(\bS^{2d}) \cong 
\left\{   \begin{array}{cl}
\Q & \mbox{for $k \in \{2d,4d-1\}$} \\ 
\0 & \mbox{otherwise.} \\ 
\end{array}
\right. $$ 
We therefore have complete information on the rational homotopy groups of 
finite-dimensional Lie groups~$G$. In particular, $\pi_2(G)$ vanishes and 
all groups $\pi_{2k}(G)$ are torsion groups. 
\end{rem}

For the following proposition we recall that a locally exponential 
Lie group $G$ has no small subgroups if and only if its Lie algebra 
$\g$ carries a continuous norm (Proposition~\ref{prop:8.11.3}). 

\begin{prop}
  \mlabel{prop:mn03-lemma-VII.8}
Let $K$ be a locally exponential Lie group with no small subgroups 
and $(X,x_0)$ be a connected pointed  topological space. 
Then the constant map $\be$ is the only element of finite order in 
the group~$C_*(X,K)$.
\end{prop}

\begin{prf}
Assume that $f^k = \be$ holds for some $f \in C_*(X,K)$. 
Let $U \subeq K$ be an identity neighborhood 
containing no small subgroups and $V \subeq U$ be an open 
identity neighborhood  with $V^k \subeq U$. 
Then the only element of order $k$ in $V$ is $\be$ because 
otherwise $U$ would contain a non-trivial subgroup of $K$. 
Therefore $f^{-1}(V)$ is an open subset of $X$ which coincides with 
$f^{-1}(\{\be\})$, hence is also closed. As $f$ preserves base points, this 
set is non-empty, and the connectedness of $X$ implies that $f$ is constant $e$. 
\end{prf}

\begin{ex}
  \mlabel{ex:VII.9-mn03} 
Let $K$ be a compact connected simple 
Lie group. 

(a) For $G := C_*(\bS^1,K)$, the group 
$\pi_0(G) \cong \pi_1(K)$ is finite (Example~\ref{ex:spheres}) and 
Proposition~\ref{prop:mn03-lemma-VII.8} implies that 
the short exact sequence $G_0 \into G \onto \pi_0(G)$
does not split because $G$ contains no element of finite order. 

(b) For $G := C_*(\bS^2,K)$, we have 
$\pi_1(G) \cong \pi_3(K) \cong \Z$ (Remark~\ref{rem1.3}, Example~\ref{ex:spheres}). 
Since $G$ contains no elements of finite order 
(Proposition~\ref{prop:mn03-lemma-VII.8}), all homomorphisms 
$\T \to G$ are trivial: 
\[ \Hom(\T,G) \cong \exp_G^{-1}(\be) = \{0\}.\] 
Therefore no element of $\pi_1(G)$ can be realized 
by a homomorphism $\T \to G$. This contrasts the corresponding situation 
for finite-dimensional Lie groups.  

As $\pi_2(K)$ vanishes, the group $G$ is connected. 
Let $q_G \: \tilde G \to G$ denote its universal covering group. 
As all homomorphisms $\T \to G$ are trivial, it follows that 
\[ \pi_1(G) \cong \ker q_G \not\subeq \exp_{\tilde G}(\g).\] 
This never happens for connected finite-dimensional Lie 
groups, for which $Z(G)$ is always contained in $\exp(\g)$ \cite[Thm.~14.2.1]{HiNe12}. 
\end{ex}

\section{Diffeomorphism groups} 
\mlabel{sec:top-diffeo}

In this section, we briefly 
collect some results on the homotopy types of diffeomorphism groups 
of compact manifolds with boundary. 
For more details we refer to \cite{Mil84, Ha12}. 

\subsection{Spheres and balls} 

For $M = \bS^1$ we have already seen in Chapter~\ref{ch:diffeo} 
that the universal covering group 
of $\Diff_+(\bS^1)$ is contractible. This implies that 
$$ \pi_k(\Diff(\bS^1)) \cong 
\left\{   \begin{array}{cl}
\Z_2 & \mbox{for $k = 0$} \\ 
\Z & \mbox{for $k = 1$} \\ 
\0 & \mbox{otherwise.} \\ 
 \end{array}
\right.  $$ 

For $M = \bS^d$, $d \geq 2$, 
the situation is more complicated. Of course we have a natural 
inclusion $\OO_{d+1}(\R)\into \Diff(\bS^d)$ and one may ask for which dimension $d$ this 
inclusion is a homotopy equivalence. For $d =1,2$ this has been proved by S.\ Smale 
in 1959 and conjectured for $d= 3$. This conjecture was 
proved in 1983 by Hatcher \cite{Ha83}. 
For $d = 4$ the answer is not known to, 
and for $d > 4$ the inclusion is not 
a homotopy equivalence (\cite[p.~1053]{Mil84}). 

For $d =2$ this leads to the following information on the homotopy groups. 
As $\OO_3(\R) \cong \SO_3(\R) \times \Z_2$ and the universal covering group 
$\SU_2(\C)$ of $\SO_3(\R)$ is homeomorphic to $\bS^3$, we obtain from 
Remark~\ref{rem1.3}: 
$$ \pi_k(\Diff(\bS^2)) \cong 
\left\{   \begin{array}{cl}
\Z_2 & \mbox{ \ for\   $k = 0$} \\ 
\Z_2 & \mbox{ \ for\   $k = 1$} \\ 
\0 & \mbox{ \ for\   $k = 2$} \\ 
\Z & \mbox{ \ for\   $k = 3$} \\ 
\Z_2 & \mbox{ \ for\   $k = 4$} \\ 
\Z_2 & \mbox{ \ for\   $k = 5$.} 
\end{array}
\right.  $$ 

For $d = 3$ we have $\OO_4(\R) \cong \SO_4(\R) \rtimes \Z_2$ and the 
universal covering group of $\SO_4(\R)$ is a two-fold covering by 
$\SU_2(\C)^2$. This leads to 
$$ \pi_k(\Diff(\bS^3)) \cong 
\left\{   \begin{array}{cl}
\Z_2 & \mbox{\ for\  $k = 0$} \cr 
\Z_2 & \mbox{\ for\  $k = 1$} \cr 
\0 & \mbox{\ for\  $k = 2$} \cr 
\Z^2 & \mbox{\ for\  $k = 3$} \cr
\Z_2^2 & \mbox{\ for\  $k = 4$} \cr 
\Z_2^2 & \mbox{\ for\  $k = 5$} \cr
\end{array}
\right.  $$ 

The group $\pi_0(\Diff_+(\bS^d))$, which is finite for $d \geq 5$, has a remarkable 
differential geometric interpretation. Its elements correspond to 
oriented diffeomorphism classes of smooth $(d+1)$-dimensional manifolds homeomorphic to 
$\bS^{d+1}$. According to Kervaire and Milnor \cite{KM63}, all these groups are finite 
and abelian. A table with the orders of these groups for $d = 5$ to $17$ can be found in 
\cite{Mil84}. 

\begin{rem} (\cite{Ha12}) As the sphere $\bS^{n-1}$ is the boundary 
$\partial\bD^n$ of the closed unit ball $\bD^n$ in $\R^n$, one has a 
natural fibration 
\[ \Diff_\partial(\bD^n) \to \Diff(\bD^n) \into \Diff(\bS^{n-1}) \] 
obtained by restricting to the boundary. Here 
$\Diff_\partial(M)$ \index{$\Diff_\partial(M)$} 
denotes for a manifold $M$ with boundary the group 
of diffeomorphisms fixing all boundary points. 
 
One also has 
\[ \Diff(\bS^n)\sim \OO_{n+1}(\R) \times \Diff_\partial(\bD^n),\]
where $\sim$ denotes homotopy equivalence,  \index{$\sim$; homotopy equivalence} 
so that the homotopy type of $\Diff(\bS^n)$ can also be described in terms 
of diffeomorphisms of the ball. For $n = 1,2,3$, the group 
$\Diff_\partial(\bD^n)$ is weakly contractible 
(\cite[Thm.~1.3.2]{Ku18}) and, for $n \geq 5$, 
$\pi_0(\Diff_\partial(\bD^n)) \cong \Theta_{n+1}$, where 
$\Theta_{n}$ is the group of homotopy $n$-spheres under 
connected sum (\cite{Ku18}). 
\end{rem}

\subsection{Products of spheres and intervals} 

\begin{exs} (Diffeomorphism groups of surfaces; \cite{Ha12}) 
For surfaces we have: 
  \begin{enumerate}
  \item[\rm(a)] $\pi_0(\Diff(\T^2)) \cong \GL_2(\Z)$ and 
$\Diff(\T^2)_0 \sim \T^2$ (\cite{EE67, EE69}).
  \item[\rm(b)] $\pi_0(\Diff(\bS^1 \times [0,1])) \cong \Z_2 \times \Z_2$ 
and $\Diff(\bS^1 \times I)_0 \sim \T$ (\cite{ES70})
  \item[\rm(c)] $\Diff(\bP(\R^3))_0 \sim \SO_3(\R) \sim \Diff(\bS^2)_0$. 
For the Klein bottle $K$, the disc $\bD^2$, the annulus 
$\bS^1 \times I$ and the M\"obius strip $M$, the identity component 
of the diffeomorphism 
group has $\SO_2(\R) \cong \T$ as a strong deformation retract 
(\cite{ES70}). 
  \item[\rm(d)] $\Diff(\Sigma)_0$ is contractible for all 
other surfaces~$\Sigma$ (\cite{EE67, EE69}). The group $\pi_0(\Diff(\Sigma))$ is 
the {\it mapping class group}. \index{mapping class group} 
  \end{enumerate}
\end{exs}

\begin{ex} (Diffeomorphism groups of some $3$-manifolds; \cite{Ha12})
  \begin{enumerate}
  \item[\rm(a)] The group $\Diff(\bS^1 \times \bS^2)$ is homotopy equivalent to 
the product group $\OO_2(\R) \times \OO_3(\R) \times \Omega(\SO_3(\R))$ 
(\cite{Ha81}).
  \item[\rm(b)] $\Diff(\T^3)_0 \sim \T^3$. 
  \item[\rm(c)] $\Diff(\T^2 \times [0,1])_0 \sim \T^2$. 
  \end{enumerate}
\end{ex}

\begin{rem} (a) For $n \geq 5$, $\pi_0(\Diff(\T^n))$ contains a copy of 
$\Z_2^\infty$ (\cite{Ha12}).   

(b) \cite[Cor.~4.1]{Yag10} implies that, for any connected, 
$\sigma$-compact smooth manifold $M$ with boundary, 
endowed with a volume form $\mu$, 
the group $\Diff(M,\mu)$ of volume preserving diffeomorphisms is a strong 
deformation retract of the group $\Diff(M)_+$ of orientation 
preserving diffeomorphisms. This reduces all questions on the 
homotopy type of $\Diff(M)$ to the group $\Diff(M,\mu)$. 
For compact manifolds, a corresponding result can be 
derived from the Nash--Moser Inverse Function Theorem 
(\cite[Thm.~2.5.3]{Ham82}). This exhibits the group 
$\Diff(M)$ as a smooth principal bundle over the space 
of smooth probability measures on~$M$.
\end{rem}

\begin{rem} For a symplectic manifold $(M,\omega)$, 
not much is known about the topology of the group 
$\Sp(M,\omega)$ of symplectomorphisms and its normal subgroup 
$\Ham(M,\omega)$ of hamiltonian diffeomorphisms 
(\cite[\S 7.2]{Pol01}): 
\begin{itemize}
\item[$\bullet$] $\Ham(\R^{2n},\omega)$ is contractible for 
$n =1,2$, but nothing is known for $n > 2$. 
\item[$\bullet$] For $M = \bS^2$, the inclusion 
$\OO_3(\R) \into \Ham(\bS^2,\omega) = \Sp(\bS^2,\omega)$ is a 
homotopy equivalence, so that $\pi_1(\Ham(\bS^2,\omega)) \cong \Z/2\Z$. 
\item[$\bullet$] For a compact symplectic surface $(\Sigma, \omega)$ 
of genus $g \geq 1$, the group $\Ham(\Sigma,\omega)$ is simply 
connected. 
\end{itemize}
\end{rem}


\section{Groups of operators on Hilbert spaces}  
\mlabel{sec10.4}

In this section, we first discuss Kuiper's Theorem. 
It asserts that the unitary group (with the operator norm)  
for a Hilbert space $\cH$ over $\K \in \{\R,\C,\bH\}$ is 
contractible.\begin{footnote}
{Here a Hilbert space 
over $\K$ means a $\K$-vector space endowed with a positive definite hermitian form 
$\la \cdot, \cdot \ra \: \cH \times \cH \to \K,$ 
linear in the first argument, which is complete with respect to 
the norm $\|v\| := \sqrt{\la v,v\ra}$.}
\end{footnote}
We then use 
this result to prove that various classical groups of operators on 
Hilbert spaces such as 
$$ \GL(\cH,I) = \{ g \in \GL(\cH) \: Ig^*I^{-1} = g^{-1}\} $$
($I$ an antilinear isometry on $\cH$ with $I^2 \in \{\pm \1\}$) 
are contractible. 
Then we turn, for a set $J$, 
to the direct limit groups $\GL_J(\K)$ of those invertible 
$J \times J$-matrices $g$ for which $g - \1$ has only finitely many non-zero 
entries. We will see that, for any infinite set $J$, this group is weakly 
homotopy equivalent to $\GL_{\N}(\K)$, i.e., to the direct limit of the groups 
$\GL_n(\K)$. Combining these insights with the general results of Palais from 
Section~\ref{sec:10.1}, we compute the homotopy groups of 
the congruence groups 
\[ \GL_p(\cH) := \GL(\cH) \cap (\1 + \cL_p(\cH)) \] 
of the Schatten ideals $\cL_p(\cH) \subeq \cL(\cH)$ 
and of the groups 
\[ \GL_p(\cH,I) := \GL(\cH,I) \cap \GL_p(\cH) 
\quad \mbox{ for } \quad 1 \leq p \leq \infty.\]

\subsection{Kuiper's Theorem} \mlabel{kuipersec} 

In this subsection, we explain how Kuiper's Theorem that the group
$\GL_\K(\cH)$ of $\K$-linear continuous operators on an 
infinite-dimensional separable \break $\K$-Hilbert space $\cH$ 
is contractible (\cite{Ku65}) can be obtained in a quite elementary
way for inseparable Hilbert spaces. 
It  is based on the following lemma, which is a refinement
of \cite[Th.~14.10]{vN50}.  

\begin{lem} \mlabel{kuiperlem} 
Let $\cH$ be a Hilbert space over $\K \in \{\R, \C, \bH\}$ and 
${\cal M} \subeq \cL_\K(\cH)$ a separable set of operators. Then there exists
an orthogonal decomposition $\cH \cong \hat\bigoplus_{j \in J} \cH_j$ into
${\cal M}$-invariant subspaces such that each $\cH_j$ is separable. 
If $\cH$ is infinite-dimensional, then the spaces 
$\cH_j$ can be chosen in such a way that they are all
infinite-dimensional, hence isomorphic to $\cH_s := \ell^2(\N,\K)$, and
thus $\cH \cong \ell^2(J,\cH_s)$. 
\end{lem}

\begin{prf} Since the closed $*$-subalgebra of $\cL(\cH)$ generated by ${\cal
M}$ is separable, we may assume that ${\cal M}$ is $*$-invariant with
$\1 \in {\cal M}$, otherwise we replace it by ${\cal M} \cup {\cal
M}^* \cup \{\1\}$. Now the assertion follows by a standard application of Zorn's Lemma:  
Let $(\cH_j)_{j \in J}$ be a maximal set of non-zero 
closed ${\cal M}$-invariant separable subspaces of $\cH$ 
such that the sum $\sum_{j\in J} \cH_j$ is orthogonal. 
Set $\cH_0 := \oline{\sum_{j\in J} \cH_j}$. Then $\cH_0^\bot$ is ${\cal
M}$-invariant because ${\cal M}$ is $*$-invariant. Assume that 
$\cH_0 \not= \cH$. For $0 \not= v \in \cH_0^\bot$ the subspace 
$\cH_v := \oline{\Spann({\cal M}v)}$ is separable and orthogonal 
to all the spaces $\cH_j$, contradicting the 
maximality of the family $(\cH_j)_{j \in J}$. 
This proves the first assertion.

To prove the second part, let us assume that $\cH$ is
infinite-dimensional and consider a decomposition 
$\cH \cong \hat\bigoplus_{j \in J} \cH_j$ as above. Let 
$$ J_0 := \{ j \in J \: \dim \cH_j < \infty\}.$$

\nin {\bf Case $1$:} If $J_0$ is finite, then there exists a $j_0 \in J \setminus
J_0$. Replacing $\cH_{j_0}$ by the larger space  
$\cH_{j_0} + \sum_{i \in J_0} \cH_i$, we obtain the desired decomposition. 

\nin {\bf Case $2$:} If $J_0$ is infinite, then the equality of cardinals 
$|J_0 \times \N| = |J_0|$ (\cite[App.~2]{La93})
implies that $J_0$ can be partitioned into infinite countable subsets 
$I_i$, $i \in J_0$. Then all the subspaces 
$\cK_i := \oline{\sum_{j \in I_i} \cH_j}$
are infinite-dimensional and separable, and we have the desired orthogonal
decomposition of $\cH$: 
\[  \cH = \hat\bigoplus_{j \in J \setminus J_0} \cH_j 
\oplus \hat\bigoplus_{i \in J_0} \cK_i. \qedhere \] 
\end{prf}

The proof of the following proposition is inspired by the setting in Mityagin's 
paper~\cite{Mit70}. 

\begin{prop} \mlabel{prop-kuip}  Let $Y$ be a separable topological space and 
$\cH$ an inseparable Hilbert space over $\K \in \{\R, \C,\bH\}$. 
Then each continuous map $f \: Y \to \GL_\K(\cH)$ is homotopic to a
constant map. 
\end{prop}

\begin{prf} Since $f(Y)$ is a separable set of operators, Lemma~\ref{kuiperlem} implies that
there exists a set $J$ and an isomorphism 
$\cH \to \ell^2(J,\cH_s)$ with $\cH_s := \ell^2(\N,\K)$ 
such that the operators in $f(Y)$ are diagonal operators on $\ell^2(J,\cH_s)$. 

\nin{\bf Step 1:} Since $\cH$ is not separable, the set $J$ is (uncountably) infinite. 
First we consider a decomposition $J = J_1 \dot\cup J_2$ with 
$|J_1| = |J_2| = |J|$. This leads to an orthogonal decomposition 
$\cH \cong \cH \oplus \cH$, and we consider operators on $\cH$ accordingly as 
block $2 \times 2$-matrices. Then $f \: Y \to \GL_\K(\cH)$ can be written as 
\[  f(y) = \pmat{g_1(y) & 0 \\ 0 & g_2(y)},\]
where $g_j \: Y \to \GL_\K(\cH)$ are continuous maps. We claim that 
$f$ is homotopic to the map 
\begin{equation}
  \label{eq:top.3.2}
f_1(y) = \begin{pmatrix} g_1(y) g_2(y) & 0 \\  0 & \1 \end{pmatrix}.
\end{equation}
It suffices to show that 
$$ f(y)^{-1} f_1(y)  = \begin{pmatrix} g_2(y) & 0 \\ 0 & g_2(y)^{-1} \end{pmatrix}$$
is homotopic to a constant map. This is implemented by the homotopy 
$$ H(t,a) := 
\begin{pmatrix} 1 & 0 \cr t(a^{-1} - 1) & 1\end{pmatrix} 
\begin{pmatrix} 1 & t \cr 0  & 1\end{pmatrix} 
\begin{pmatrix} 1 & 0 \cr t(a - 1) & 1\end{pmatrix} 
\begin{pmatrix} 1 & - t a^{-1}\cr 0  & 1\end{pmatrix}  $$
which satisfies 
$$ H(1,a) = \begin{pmatrix} a & 0 \cr 0 & a^{-1} \end{pmatrix} 
\quad \hbox{ and } \quad H(0,a) = \begin{pmatrix} \1 & 0 \cr 0 & \1 \end{pmatrix}. $$

\nin{\bf Step 2:} In view of Step $1$, we may assume that $f_1 \: Y \to
\GL_\K(\cH)$ has the form \eqref{eq:top.3.2}. 
Next we observe that $\cH \cong \ell^2(\N,\cH)$
because $|J| = |\N \times J|$ (\cite[App.~2]{La93}). Therefore 
$$ \cH \cong \cH \oplus \ell^2(\N, \cH) \cong \ell^2(\N_0, \cH), $$
and we may assume that 
$$ f_1(y) = \diag( g(y), \1, \1, \ldots ),  $$
where $d = \diag(d_1,d_2, \ldots)$ is the operator defined by 
$d.(v_n)_{n \in \N_0} = (d_n v_n)_{n \in \N}$. 
Partitioning $\N$ into odd numbers $\N_{\rm odd}$
and even numbers $\N_{\rm even}$, and writing accordingly 
$$ \ell^2(\N,\cH) \cong \ell^2(\N_{\rm odd},\cH) \oplus \ell^2(\N_{\rm even},\cH), $$
it follows from Step 1 that the constant map $Y \to \GL_\K(\ell^2(\N,\cH)), y \mapsto \1$ is
homotopic to the map 
$$ y \mapsto \diag(g(y)^{-1}, g(y), g(y)^{-1}, g(y), \ldots ). $$
Therefore $f_1$ is homotopic to 
$$ f_2(y) = \diag( g(y), g(y)^{-1}, g(y), g(y)^{-1}, \ldots ). $$
Applying the same argument again to the decomposition 
$$ \cH 
\cong \cH \oplus \ell^2(\N,\cH) 
\cong \ell^2(\{0\} \cup \N_{\rm even}, \cH) \oplus \ell^2(\N_{\rm odd},\cH), $$
we see that $f_2$ is homotopic to a constant map. 
\end{prf}

\index{Kuiper's Theorem for general Hilbert spaces} 
\begin{thm} [Kuiper's Theorem for general Hilbert spaces] \mlabel{kuiper-thm}  
If $\cH$ is an infinite-dimensional Hilbert space over $\K \in \{\R, \C, \bH\}$,
then the group $\GL_\K(\cH)$ is contractible. 
\end{thm}

\begin{prf} In view of Theorem \ref{pal-contr}, it suffices to show that all
homotopy groups of $\GL_\K(\cH)$ vanish. 
In \cite{Ku65} (see also \cite[Thm.~3.15]{Bel06}), 
this assertion is proved for infinite-dimensional
separable Hilbert spaces. For inseparable Hilbert spaces, this 
follows from Proposition \ref{prop-kuip} because the spheres 
$\bS^k$, $k \in \N_0$, are separable. 
\end{prf}

\subsection{Applications of Kuiper's Theorem} 

Kuiper's Theorem can be used to prove that many ``classical'' groups of operators on 
a Hilbert space are contractible. Below we briefly discuss these applications. 

\begin{defn} \mlabel{def10.4.2}
(a) If $\cH$ is a Hilbert space over $\K \in
\{\R,\bC,\bH\}$, then we define 
$$ \U_\K(\cH) := \{ g \in \GL_\K(\cH) \: g^* g = gg^* = \1\} $$
as the unitary part of this group. We also write 
$$ \OO(\cH) := \U_\R(\cH), \quad 
\U(\cH) := \U_\C(\cH) \quad \hbox{ and } \quad 
\Symp(\cH) := \U_\bH(\cH). $$

 (b) Let $\cH$ be a complex Hilbert space and $I$  be an
antilinear isometry with $I^2 \in \{\pm \1\}$. Then 
$$ \GL(\cH,I) := \{ g \in \GL(\cH) \: Ig^*I^{-1} = g^{-1}\} $$
is a complex Lie subgroup of $\GL(\cH)$ (Exercise~\ref{exer:oden-iv.9})). 
For $I^2 = \1$ we then have 
\[  \U(\cH,I) := \U(\cH) \cap \GL(\cH,I) \cong \OO(\cH_\R) \quad \hbox{ with } \quad 
\cH_\R := \{ x \in \cH \: Ix = x\}, \]
and for $I^2 = -\1$ we have 
$$ \U(\cH,I) \cong \U_\bH(\cH) \cong \Symp(\cH), $$
where the quaternionic structure on $\cH$ is given by the subalgebra 
$\C \1 + \C I \cong \bH$ of $\cL_\R(\cH)$, the real linear endomorphisms of $\cH$. 

 (c) (Hermitian groups) Let $\cH$ be a complex Hilbert space and
$$\cH = \cH_+ \oplus \cH_-$$ 
be an orthogonal decomposition. Further let 
$T = T^* \in \cL(\cH)$ with $\cH_\pm = \ker(T \mp \1)$. We define the
corresponding 
{\it pseudo-unitary group}   
\index{pseudo-unitary group}   
$$ \U(\cH_+, \cH_-) := \{ g \in \GL(\cH) \: Tg^*T^{-1} = g^{-1}\}. $$

We further obtain a real bilinear skew-symmetric form 
$\Omega(x,y) := \Im \la x, y \ra$ and write $\cH^\R$ for the
real Hilbert space underlying $\cH$. Then 
$$ \Symp(\cH,\Omega) := \{ g \in \GL_\R(\cH^\R) \: (\forall v,w \in \cH^\R) 
\, \Omega(gv, gw) = \Omega(v,w)\} $$
is called the 
\index{symplectic group of Hilbert space}
{\it symplectic group of $\cH$}. 
If we start with the real Hilbert space $\cH^\R$ and consider an
isometric complex structure $I$ on $\cH^\R$, then we can define 
$$ \Omega(x,y) := - \la Ix, y \ra =  \la x, Iy \ra $$
and put 
$$ \Symp(\cH^\R,I) := \{ g \in \GL_\R(\cH^\R) \: (\forall v,w \in \cH^\R) 
\, \Omega(gv, gw) = \Omega(v,w)\}. $$
It is easy to see that both constructions lead to isomorphic groups 
$\Symp(\cH^\R,I) \cong \Symp(\cH,\Omega)$. 

If $I$ is a conjugation on the complex Hilbert space $\cH$
and $\cH_+ \subeq \cH$ be a subspace for which we get an orthogonal
decomposition $\cH = \cH_+ \oplus \cH_-$ with $\cH_- := I\cH_+$. Then we define 
$$ \OO^*(\cH,I) := \U(\cH,I) \cap \U(\cH_+, \cH_-). $$
\end{defn}

\begin{prop}
  \mlabel{prop:a.2-ne04} If $\cH$ is a Hilbert space over $\K \in \{\R,\C,\bH\}$, 
then the polar map 
\[  p \: \U_\K(\cH) \times \Herm_\K(\cH) \to \GL_\K(\cH), \quad 
(u,X) \mapsto u e^X \]
is a diffeomorphism. 
\end{prop}

\begin{prf} We first consider the case $\K =\C$. 
Let $g \in \GL(\cH)$. Then $g^*g$ is a positive hermitian operator, so
that the continuous functional calculus provides a positive hermitian
operator 
$X := {\textstyle{1\over 2}}\log(g^*g).$ 
Let $u := g e^{-X}$. Then 
\[  uu^* = g e^{-2X} g^* 
= g (g^*g)^{-1} g^* = \1 \] 
and 
\[  u^* u 
= e^{-X} g^* g e^{-X} 
= e^{-X} e^{2X} e^{-X} = \1. \] 
We conclude that every operator $g \in \GL(\cH)$ has a unique
decomposition $g = u e^X$ with $X \in \Herm(\cH)$. This means that $p$
is bijective. It is also clear that $p$ is a smooth map. 
That $p^{-1}$ is also smooth follows from the smoothness of the function 
\[  \log \: \{ g \in \GL(\cH) \cap \Herm(\cH) \: \inf \Spec(g) > 0 \} \to
\Herm(\cH) \] 
(Lemma~\ref{lem:a.1-ne04}). 

If $\K = \R$, we consider the complex Hilbert space $\cH_\C$, on 
which we obtain an antilinear isometric involution 
by $I(x + i y) := x - iy$. Then $\tau(A) := IAI^{-1} = IAI$ defines an 
antilinear involutive automorphism of the Banach-$*$ algebra $\cL(\cH_\C)$ and 
\[ \GL_\C(\cH_\C)^\tau \cong \GL_\R(\cH), \quad 
\U_\C(\cH_\C)^\tau \cong \U_\R(\cH), \quad 
\Herm_\C(\cH_\C)^\tau \cong \Herm_\R(\cH),\] 
where the isomorphisms are implemented by restriction to the real 
subspace $\cH = (\cH_\C)^I$ of $I$-fixed points. 
Writing $g \in \GL_\C(\cH_\C)$ according to the polar decomposition 
as $g = u e^X$ with $u \in \U_\C(\cH_\C)$ and $X \in \Herm_\C(\cH_\C)$, the relation 
$\tau(g) = \tau(u) e^{\tau(X)}$ and the uniqueness of the polar decomposition 
implies that $\tau(g) = g$ is equivalent to 
$\tau(u) = u$ and $\tau(X) =X$. Therefore the polar map 
$p$ for $\GL_\C(\cH_\C)$ restricts to a diffeomorphism 
\[ \U_\R(\cH) \times \Herm_\R(\cH) \to \GL_\R(\cH), \quad 
(u,X) \mapsto u e^X.\] 

If $\K = \bH$, we consider the complex Hilbert space $\cH^\C$ obtained 
by restricting the scalars from $\bH = \C + J \C$ to the subalgebra $\C$ 
and the scalar product to its $\C$-component. 
Then $\tau(A) := JAJ^{-1} = -JAJ$ defines an 
antilinear involutive automorphism of the Banach-$*$ algebra $\cL(\cH^\C)$ and 
\[ \GL_\bH(\cH^\C)^\tau = \GL_\bH(\cH), \quad 
\U_\C(\cH^\C)^\tau = \U_\bH(\cH), \quad 
\Herm_\C(\cH_\C)^\tau = \Herm_\bH(\cH).\] 
Now the same arguments as in the real case show that the 
polar decomposition of $\GL_\C(\cH^\C)$ restricts to a diffeomorphism 
\[ \U_\bH(\cH) \times \Herm_\bH(\cH) \to \GL_\bH(\cH).\qedhere\] 
\end{prf}

\begin{rem}
  \mlabel{rem:a.3-ne04} Our proof for the polar decomposition works also for
abstract $C^*$-algebras, where it provides a diffeomorphism 
\[  p \: \U(\cA) \times \Herm(\cA) \to \cA^\times, \quad 
(u,x) \mapsto u e^x, \] 
where 
\[  \U(\cA) = \{ a \in \cA \: aa^* = a^*a = \1\}\] 
is the unitary group of $\cA$. 
For commutative algebras $\cA = C(X,\C)$, $X$ a compact space, this is
the trivial decomposition 
\[  C(X,\T) \times C(X,\R) \to C(X,\C^\times), \quad 
(u,f) \mapsto u e^f.\qedhere \] 
\end{rem}

\begin{thm} \mlabel{thm10.4.3} 
If $\cH$ is an infinite-dimensional Hilbert space over 
$\K \in \{\R,\bC,\bH\}$, then the following groups are contractible: 
\begin{enumerate}
\item[\rm(i)] the group of $\K$-linear automorphisms $\GL_\K(\cH)$. 
\item[\rm(ii)] the group of isometric $\K$-linear automorphisms
$\U_\K(\cH)$. 
\item[\rm(iii)] the group $\GL(\cH,I)$ if $\cH$ is complex and $I$ an
antilinear isometry with $I^2 \in \{\pm \1\}$. Moreover, $\GL(\cH,I)$ has
a smooth polar decomposition. 
\item[\rm(iv)] the hermitian groups $\U(\cH_+, \cH_-)$, where $\cH = \cH_+ \oplus
\cH_-$ is an orthogonal decomposition with two infinite-dimensional
summands, $\Symp(\cH,\Omega)$, and $\OO^*(\cH,I)$. 
\end{enumerate}
\end{thm}

\begin{prf}  (i) is Theorem \ref{kuiper-thm}. 

 (ii) follows from (i) and the polar decomposition 
(Proposition~\ref{prop:a.2-ne04}).

 (iii) In view of Definition \ref{def10.4.2}, the group $\U(\cH,I)$ is
contractible, because it is one of the groups in (ii). Hence the
assertion follows from the polar decomposition 
of $\GL(\cH,I)$ which can be obtained as follows. 
We consider the automorphism $\tau(g) := I(g^*)^{-1}I^{-1}$ of $\GL(\cH)$ and 
write $\tau_\g(x) := - Ix^*I^{-1}$ for the corresponding antilinear automorphism of its 
Lie algebra $\gl(\cH)$. Then 
$$ \GL(\cH,I) = \GL(\cH)^\tau := \{ g \in \GL(\cH) \: \tau(g) = g \}. $$
Let $g = u e^x$ be the polar decomposition of $g \in \GL(\cH)$. 
Then $\tau(g) = \tau(u) e^{\tau_\g(x)}$ is the polar 
decomposition of $\tau(g)$, so that 
the uniqueness of this decomposition implies that $\tau(g) = g$  
is equivalent to 
$\tau(u) = u$ and $\tau_\g(x) = x$, i.e., 
$u \in \U(\cH,I)$ and $x \in \Herm(\cH,I)$. 

 (iv) For the hermitian groups we shall show that they 
have polar decompositions with 
$$ \U(\cH_+, \cH_-) \cap \U(\cH) \cong \U(\cH_+) \times \U(\cH_-), 
\quad \Symp(\cH,\Omega) \cap \OO(\cH^\R) \cong \U(\cH)  $$
and 
$$ \OO^*(\cH,I) \cap \U(\cH) \cong \U(\cH_+), $$
where $\cH \cong \cH_+ \oplus I\cH_+$ as in Definition~\ref{def10.4.2}(c). 
Therefore (ii) implies that all these groups are contractible. 

To derive the polar decomposition of $\U(\cH_+, \cH_-)$, let 
$g \in \GL(\cH)$ with polar decomposition 
$g = u e^x$, $u \in \U(\cH)$ and $x=x^*$. 
For $T$ as in Definition~\ref{def10.4.2}(c), 
 we consider the automorphism $\tau(g) := T(g^*)^{-1}T^{-1}$ of $\GL(\cH)$ and 
write $\tau_\g(x) := - Tx^*T^{-1}$ for the corresponding antilinear automorphism of its 
Lie algebra $\gl(\cH)$. Then 
$\tau(g) = \tau(u) e^{\tau_\g(x)}$ is the polar decomposition of $\tau(g)$, so that 
the uniqueness of this decomposition implies that $\tau(g) = g$ is equivalent to 
$\tau(u) = u$ and $\tau_\g(x) = x$. Therefore $g \in \U(\cH_+,\cH_-)$ if 
and only if 
$$u \in \U(\cH_+,\cH_-) \cap \U(\cH) \cong \U(\cH_+) \times \U(\cH_-) 
\quad \hbox{ and } \quad x \in \fu(\cH_+,\cH_-).$$ 

To see that $\Symp(\cH, \Omega)$ is adapted to the polar decomposition, we 
observe that 
$$ \Omega(x,y) = \Im \la x, y \ra = \Re \la x, iy \ra = (x,Jy), $$
where $(\cdot, \cdot) := \Re \la \cdot, \cdot \ra$ denotes the real 
scalar product on $\cH^\R$. 
Therefore $g \in \Symp(\cH, \Omega)$ is equivalent to 
$g^\top J g = J$, i.e., $g = \tau(g) :=  J(g^\top)^{-1}J^{-1}$. 
Then $\tau$ is an involutive automorphism of $\GL_\R(\cH^\R)$ and 
$\tau_\g(x) := -Jx^{\top}J^{-1}$ is the corresponding Lie algebra automorphism.
Let $g = u e^x$ be the polar decomposition of $g \in \GL_\R(\cH^\R)$, where 
$u \in \OO(\cH^\R)$ and $x^\top = x$. Then 
$\tau(g) = \tau(u) e^{\tau_\g(x)}$ is the polar decomposition of $\tau(g)$ because 
$ue^{-x}$ is the polar decomposition of $(g^\top)^{-1}$. 
Therefore $g \in \Symp(\cH, \Omega)$ is equivalent to 
$\tau(u) = u$, i.e., $u \in \U(\cH)$, and to $Jx = - xJ$, i.e., $x$ is antilinear. 

The argument for the group $\OO^*(\cH,I)$ is similar. 
\end{prf}

\section{Direct limit groups} 
\mlabel{sec:top-dirlim}

In this section we turn to homotopy groups of direct limits 
of finite-dimensional Lie groups. A central result in this 
context is Bott's Periodicity Theorem (Theorem~\ref{thm2.12}). 
It provides complete information on all homotopy groups of 
direct limits related to the inclusions 
$\GL_n(\K) \into \GL_{n+1}(\K)$ for $\K = \R,\C,\bH$, and 
also for the associated symmetric spaces. 
Similar techniques apply to the groups 
$\bigcup_n \GL_n(\cA)$, where $\cA$ is a unitary Banach algebra. 
Here we only discuss the identity of its homotopy groups 
with those of its natural completion.

\subsection{Direct limits of finite-dimensional matrix groups} 

\begin{defn} \mlabel{def:2.10} For an infinite set $J$, we view a function 
${m \: J \times J \to \K}$ as a matrix with entries $m(i,j)$. 
In this sense, we write $M_J(\K)$ for the set of all $J \times
J$-matrices with finitely many non-zero entries in $\K$. 
Then $M_J(\K)$ is a real algebra with respect to matrix
multiplication. It has a unit if and only if $J$ is finite. 
We write $\1 = (\delta_{ij})_{i,j \in J}$ for the identity matrix. 
Then $\1 + M_J(\K)$ is a multiplicative monoid, and we define 
$\GL_J(\K)$ to be its group of units. We endow $\GL_J(\K)$
with the direct limit topology with respect to the subgroups 
$\GL_F(\K) := \GL_J(\K) \cap (\1 + M_F(\K))$, 
where $F \subeq J$ is a finite subset. It follows directly from the
constructions that the left and right multiplications in the group 
$\GL_J(\K)$ are continuous, but if $J$ is uncountable, then the
multiplication is not jointly continuous (\cite[Th.~7.1]{Gl03a}). 
Here we identify 
$M_F(\K)$ in a natural way with a subset of $M_J(\K)$ and likewise 
$\GL_F(\K)$ with a subset of $\GL_J(\K)$. 
\end{defn}

\begin{prop} \mlabel{prop2.11}
Let $J$ be an infinite set.  
Then, for each injective map $\N \into J$, the cor\-res\-pond\-ing map 
$\GL_{\N}(\K) \to \GL_J(\K)$ is a weak homotopy equivalence. 
\end{prop}

\begin{prf}  We may w.l.o.g.\ assume that $\N \subeq J$. 
Let $\eta \: \GL_{\N}(\K) \into \GL_J(\K)$
be the corresponding embedding of groups. 

Let $Y$ be a compact space and $f \: Y \to \GL_J(\K)$ be a continuous
map. Then there exists a finite subset $F \subeq J$ with $f(Y) \subeq
\GL_F(\K)$ (Lemma~\ref{lem-pal0}(2)). If $F' \subeq J$ is finite with $F \cap F' = \eset$ and 
$|F| = |F'| = n$, then $\GL(F \cup F',\K) \cong \GL_{2n}(\K)$, where we
identify $F$ with $\{1,\ldots, n\}$ and 
$F'$ with $\{n+1,\ldots, 2n\}$. Then $f$ is a map of the form 
$$ f(y) = \begin{pmatrix} g(y) & \0 \cr \0 & \1 \end{pmatrix}, $$
where we write the elements of $\GL_{2n}(\K)$ as block $2 \times
2$-matrices with entries in $M(n,\K)$, and $g \: Y \to \GL_n(\K)$ is a continuous map. 

We consider the map $H \: [0,1] \times \GL_n(\K) \to \GL_{2n}(\K)$ given
by 
$$ H(t,a) := 
\begin{pmatrix} \1 & 0 \\ t(a^{-1} - \1) & \1 \end{pmatrix}
\begin{pmatrix} \1 & t \\ 0  & \1 \end{pmatrix}
\begin{pmatrix} \1 & 0 \\ t(a - \1) & \1\end{pmatrix}
\begin{pmatrix} \1 & - t a^{-1}\\ 0  & \1 \end{pmatrix}  $$
which satisfies 
$$ H(1,a) = \begin{pmatrix}a & 0 \\ 0 & a^{-1} \end{pmatrix}
\quad \hbox{ and } \quad H(0,a) = \begin{pmatrix}\1 & 0 \\ 0 & \1 \end{pmatrix}. $$
Then 
$$\tilde H \: [0,1] \times Y \to \GL_{2n}(\K), \quad (t,y) \mapsto
f(y) H(t,g(y))^{-1} $$
is continuous with $\tilde H(0,y) = f(y)$ and 
$$ \tilde H(1,y) = \begin{pmatrix}\1 & \0 \\ \0 & g(y)^{-1} \end{pmatrix}. $$
This construction shows that every continuous map 
$f \: Y \to \GL_F(\K)$ is homotopic in $\GL_J(\K)$ to a continuous map 
$f' \: Y \to \GL_{F'}(\K)$. 

In particular, we see that for each continuous map 
$f \: Y \to \GL_J(\K)$ there exists a finite subset $E \subeq \N$ such
that $f$ is homotopic to a continuous map 
$\tilde f \: Y \to \GL(E,\K)$. In fact, with $F$ as above, we simply choose $E
\subeq \N$ such that $|E| = |F|$ and $E \cap F = \eset$. 
This argument shows that the natural homomorphism 
$\pi_k(\eta) \: \pi_k(\GL_{\N}(\K)) \to \pi_k(\GL_J(\K))$ is surjective. 

To see that $\pi_k(\eta)$ is injective, suppose that 
$f \: Y \to \GL_n(\K) \subeq \GL_{\N}(\K)$ is in $\GL_J(\K)$ homotopic to a constant
map. Let $H \: [0,1] \times Y \to \GL_J(\K)$ be a homotopy with 
$H(0,y) = \1$ and $H(1,y) = f(y)$ for all $y \in Y$. Then there exists
a finite subset $F \subeq J$ with $\im(\cH) \subeq \GL_F(\K)$. Then 
we may assume that $F \supeq \{1,\ldots, n\}$, and since $\GL_{|F|}(\K) \cong \GL_F(\K)$,
we see that the homotopy class of $f$ vanishes in $\GL_{|F|}(\K)
\subeq \GL_{\N}(\K)$. In particular, $f$ is homotopic to a constant map in $\GL_{\N}(\K)$. 
\end{prf}

The homotopy groups for $\GL_{\N}(\K)$ and hence for all groups $\GL_J(\K)$, where $J$ is 
an infinite set (Proposition~\ref{prop2.11}), are given by the following theorem. 

\index{Bott Periodicity Theorem} 
\begin{thm} [Bott Periodicity Theorem] \mlabel{thm2.12}
Let $\K \in
\{\R,\C,\bH\}$ and $d := \dim_\R \K$. 
Then, for $k \leq d (n+1) - 3$ and $q \in \N$, the maps 
$$ \pi_k(\GL_n(\K)) \to \pi_k(\GL_{n+q}(\K)) $$
are isomorphisms, so that 
$$ \pi_k(\GL_{\N}(\K)) \cong \pi_k(\GL_n(\K)). $$
Moreover, we have the periodicity relations 
$$ \pi_{n+2}(\GL_{\N}(\C)) \cong 
\pi_{n}(\GL_{\N}(\C)), \quad \pi_{n+4}(\GL_{\N}(\R)) \cong 
\pi_{n}(\GL_{\N}(\bH)), $$
$$ \pi_{n+4}(\GL_{\N}(\bH)) \cong 
\pi_{n}(\GL_{\N}(\R)), $$
$$ \pi_{n}(\GL_{\N}(\bH)/\GL_{\N}(\C)) \cong 
\pi_{n+1}(\GL_{\N}(\bH)), $$
$$ \pi_{n}(\GL_{\N}(\C)/\GL_{\N}(\R)) 
:= \lim_{m \to \infty} \pi_{n}(\GL_m(\C)/\GL_m(\R)) 
\cong \pi_{n+2}(\GL_{\N}(\bH)), $$
$$ \pi_{n}(\GL_{\N}(\R)/\GL_{\N}(\C)) := \lim_{m \to \infty} \pi_{n}(\GL_{2m}(\R)/\GL_m(\C)) 
\cong \pi_{n+1}(\GL_{\N}(\R)), $$
$$ \pi_{n}(\GL_{\N}(\C)/\GL_{\N}(\bH)) 
:= \lim_{m \to \infty} \pi_{n}(\GL_{2m}(\C)/\GL_m(\bH)) 
\cong \pi_{n+2}(\GL_{\N}(\R)). $$
In particular, the homotopy groups of $\GL_{\N}( \K)$ are determined by the
following table: 

\noindent
\center{\begin{tabular}{||l||l|l|l||}\hline
{} & \ \  $\GL_{\N}(\R)$ \ \ &\ \  $\GL_{\N}(\C)$\ \  &\ \ $\GL_{\N}(\bH)$\ \ \\ \hline\hline
$\pi_0$ & \hfil{$\Z/2\Z$}\hfil & \hfil $\0$ \hfil & \hfil $\0$ \hfil \\ \hline
$\pi_1$ & \hfil$\Z/2\Z$\hfil &$ \hfil\Z$\hfil & \hfil$\0$\hfil  \\ \hline
$\pi_2$ & \hfil$\0$\hfil & \hfil$\0$\hfil & \hfil$\0$\hfil \\ \hline
$\pi_3$ & \hfil$\Z$\hfil & \hfil$\Z$\hfil & \hfil$\Z$\hfil \\ \hline
\end{tabular}
}
\end{thm}

\begin{prf}  The first easy part is \cite[Th.~8.4.1]{Hu94} and the remaining
assertions can be found in \cite[pp.314ff]{Bo59} (cf.\ also \cite[Cor.~9.10.2]{Hu94}). 

For the sake of completeness, we include a proof of the first part. 
Using the polar decomposition, we may consider the corresponding maps 
of the unitary groups $\U_n(\K)$. To understand the homomorphisms 
induced by the inclusion maps 
$\eta_n \: \U_n(\K) \to \U_{n+1}(\K)$ on the homotopy groups, we consider the transitive 
action of $\U_{n+1}(\K)$ on the sphere $\bS^{d(n+1)-1}$ which leads to
a locally trivial principal bundle 
\[  \U_n(\K) \into \U_{n+1}(\K) \to \bS^{d(n+1)-1}. \] 
The exact homotopy sequence of this bundle (Theorem~\ref{homseq-princ}) contains 
\[  \pi_{k+1}(\bS^{d(n+1)-1}) \to \pi_k(\U_n(\K)) \smapright{\pi_k(\eta_n)} 
\pi_k(\U_{n+1}(\K)) \to 
\pi_k(\bS^{d(n+1)-1}) \]
For $k < d(n+1)-1$, the group 
$\pi_k(\bS^{d(n+1)-1})$ vanishes (this follows by smoothing and Sard's Theorem), 
so that $\pi_k(\eta_n)$ is surjective. 
If, in addition, $k \leq d(n+1) - 3$, then 
$k + 1 < d(n+1) - 1$ implies that also 
$\pi_{k+1}(\bS^{d(n+1)-1})$ vanishes, so that the injectivity of 
$\pi_k(\eta_n)$ follows. 
\end{prf}


\subsection{A Banach--Lie completion of $\GL_\infty(\cA)$}  \mlabel{sec10.9} 

In this subsection, we discuss a natural Banach--Lie completion of the
group $\GL_\infty(\cA) := \indlim \GL_n(\cA)$ for a unital Banach algebra $\cA$. 

Let $\cA$ be a unital Banach algebra and endow the Banach algebras $M_n(\cA)$, $n \in \N$, 
with the operator norms obtained from the action on $\cA^n$, endowed with the 
sup-norm. This norm has the property that it is submultiplicative. 
We have isometric embeddings $M_n(\cA) \into M_{n+1}(\cA)$, and we thus obtain 
a norm on the direct limit algebra 
$M_\infty(\cA) = \indlim M_n(\cA)$ such that all the embeddings 
$M_n(\cA) \into M_\infty(\cA)$ are isometric. Let 
$\hat M_\infty(\cA)$ denote the Banach completion of this normed algebra. 
This is a Banach algebra isometrically contained in $\cL(c_0(\N,\cA))$, where 
$c_0(\N,\cA)$ is the completion of $\bigcup_{n \in \N} \cA^n \cong \cA^{(\N)}$ 
in the sup-norm. 
Each element $x \in M_\infty(\cA)$ is a finite matrix, hence annihilates elements of 
norm one in $c_0(\N,\cA)$, so that $\|x - \1 \| \geq 1$. This implies that 
$\hat M_\infty(\cA)$ has no unit, and we also consider the Banach algebra 
\[ \hat M_\infty(\cA)_+ := \hat M_\infty(\cA) + \C \1 \subeq \cL(c_0(\N,\cA)).\] 

Let 
\[ \hat M_\infty(\cA)^\times := \hat M_\infty(\cA)_+^\times \cap 
(\1 + \hat M_\infty(\cA)) \] 
(a congruence subgroup) and observe that this is a Banach--Lie 
group with Lie algebra $\hat M_\infty(\cA)$. It contains the direct limit group 
$\GL_\infty(\cA) = \indlim \GL_n(\cA)$, the unit group of the submonoid $\1 + M_\infty(\cA) \subeq \hat M_\infty(\cA)_+$. 

\begin{lem} \mlabel{lem8.2} 
For each $n \in \N$, we have $\hat M_\infty(\cA)^\times \cap (\1 + M_n(\cA)) =
\GL_n(\cA)$. In particular,  
$$ \hat M_\infty(\cA)^\times \cap (\1 + M_\infty(\cA)) = \GL_\infty(\cA). $$
\end{lem}

\begin{prf}  The inclusion ``$\supeq$'' is trivial. So let $x \in M_n(\cA)$ and suppose that 
$\1 + x$ is invertible in $\hat M_\infty(\cA)_+^\times$, i.e., contained in 
$\hat M_\infty(\cA)^\times$. 
Let $\1 + y := (\1 + x)^{-1}$. Then 
$$ \1 = (\1 + y)(\1 + x) = (\1 + x)(\1 + y) $$
is equivalent to 
$$ xy = yx \quad \hbox{ and } \quad x + y + xy = 0. $$
In the matrix $xy$, at most the first $n$ rows are non-zero, whereas in the matrix 
$yx$ at most the first $n$  columns are non-zero. Therefore 
$xy =yx \in M_n(\cA)$. This implies that $y = -x - xy \in M_n(\cA)$, so that 
\[  (\1 + M_n(\cA)) \cap \hat M_\infty(\cA)^\times = \GL_n(\cA), \] 
and this implies the assertion. 
\end{prf}

Since $\hat M_\infty(\cA)^\times$ is an open subset of $\1 + \hat M_\infty(\cA)$ and $M_\infty(\cA)$ is dense in $\hat M_\infty(\cA)$, it 
follows that 
$\GL_\infty(\cA)$ is dense in $\hat M_\infty(\cA)^\times$, 
so that we may view $\hat M_\infty(\cA)^\times$ as a Banach--Lie group completion of 
$\GL_\infty(\cA)$. 

\begin{cor} \mlabel{cor8.7}
Let $\cA$ be the isometric direct limit of the Banach
algebras $(\cA_n)_{n \in \N}$ 
and assume that $\cA^\times \cap (\cA_n)_+ = \cA_n^\times$ holds for each
$n \in \N$. Then the continuous map 
$\indlim \cA_n^\times \to \cA^\times$
is a weak homotopy equivalence. 
\end{cor}

\begin{prf}  Our assumption implies that the $\cA_n$ are closed subspaces of $\cA$
whose union is dense. We recall that 
$\cA^\times = \cA_+^\times \cap (\1 + \cA)$ 
and put \break $U := \cA^\times  - \1 \subeq \cA$. This is an open subset of 
$\cA$, and, in view of the assumption that 
$\cA^\times \cap (\cA_n)_+ = \cA_n^\times$ for each $n \in \N$, we obtain from 
Theorem~\ref{thm8.6} that the map 
$\indlim \big(\cA_n^\times - \1\big) \to \cA^\times - \1$
is a weak homotopy equivalence. Since addition of $\1$ is continuous
on both sides, the assertion follows. 
\end{prf}

\begin{thm} \mlabel{thm8.8}
 Let $\cA$ be a Banach algebra and 
$\hat M_\infty(\cA)$ be the completion of $M_\infty(\cA) = \bigcup_{n \in \N} M_n(\cA)$ with
respect to a norm for which each $M_n(\cA)$ is a closed subalgebra of 
$\hat M_\infty(\cA)$. Then the continuous map 
\[  \GL_\infty(\cA) := \indlim \GL_n(\cA) \to \hat M_\infty(\cA)^\times \] 
is a weak homotopy equivalence. In particular, 
for each $k \in \N_0$, we have 
\[  \pi_k(\GL_\infty(\cA)) \cong \indlim \pi_k(\GL_n(\cA)) \cong 
\pi_k(\hat M_\infty(\cA)^\times). \] 
\end{thm}

\begin{prf} In view of Theorem~\ref{thm8.6}, we only have to use 
Lemma~\ref{lem8.2} to see 
that, for each $n \in \N$, we have $\GL_n(\cA) = \hat M_\infty(\cA)^\times 
\cap (\1 + M_n(\cA))$. 
\end{prf}

The main point in the preceding theorem is that it describes a
Banach--Lie group completion $\hat M_\infty(\cA)^\times$ 
of the direct limit topological group $\GL_\infty(\cA)$
which has the same homotopy groups. This leads to a closer connection 
between Lie theory and K-theory, which deals with the homotopy groups 
of $\GL_\infty(\cA)$ (\cite{At67, Ka81}).

\section{Congruence subgroups for Schatten ideals} 
\mlabel{sec:top-schatten}

In this section we turn to Lie groups of operators 
associated to the Schatten ideals $\cL_p(\cH) \subeq \cL(\cH)$, 
$1 \leq p \leq \infty$, consisting for $p < \infty$  
of those operators $X$ for which $(X^*X)^{p/2}$ is trace class, 
and for $p = \infty$ of all compact 
operators. We determine the homotopy groups of some Lie groups 
of operators specified in terms of the Schatten ideal.

Let $\cH$ be a Hilbert space over $\K \in \{\R,\C,\bH\}$. For $x, y \in \cH$ 
we define $P_{x,y}(v) := \la v, y \ra x$ and put $P_x := P_{x,x}$. 
Note that $P_{x,y} \in \cL_1(\cH)$. 

\begin{lem} \mlabel{lem2.13}
Let $\cH$ be a $\K$-Hilbert space and $(e_j)_{j \in J}$ be an
orthonormal basis. Then $\cL_0(\cH) := \Spann \{ P_{e_j, e_k} \: j,k \in J\}$ is
dense in each of the spaces $\cL_p(\cH)$, $1 \leq p \leq \infty$. 
\end{lem}

\begin{prf}  For each $p \in [1,\infty]$, we have 
$\|x\|_\infty \leq \|x\|_p \leq \|x\|_1$ and accordingly 
$\cL_1(\cH) \subeq \cL_p(\cH) \subeq K(\cH) = \cL_\infty(\cH)$. 

Since $\cL_0(\cH)$ and $\cL_p(\cH)$ are $*$-invariant, it suffices to
see that each hermitian operator in $\cL_p(\cH)$ is contained in the
closure of $\cL_0(\cH)$. The Spectral Theorem for Compact Hermitian
Operators (\cite{Ru91}) directly implies that the ideal $\cL_{\rm fin}(\cH)$ of 
continuous maps with finite-dimensional image is dense in $\cL_p(\cH)$,
and hence that $\cL_1(\cH)$ is dense in $\cL_p(\cH)$. Therefore it suffices to
see that $\cL_0(\cH)$ is dense in $\cL_1(\cH)$ because $\|x\|_p \leq \|x\|_1$. 

In view of the Hahn--Banach Theorem, we have to show that each
continuous linear functional $f \in \cL_1(\cH)'$ vanishing on $\cL_0(\cH)$ is
zero. As $\cL_1(\cH)' \cong \cL(\cH)$ (\cite[Prop.\ A.I.10(vi)]{Ne00}), 
the functional $f$ can be written as $f(X) = \tr(AX)$ for some 
$A \in \cL(\cH)$. Hence 
$$ f(P_{e_j,e_k}) = \tr(A P_{e_j, e_k}) = \tr(P_{Ae_j, e_k}) 
= \la Ae_j, e_k \ra. $$
If $f$ vanishes on $\cL_0(\cH)$, then the matrix of $A$ with respect to
the orthonormal basis $(e_j)_{j \in J}$ vanishes, and this means that
$A = 0$. 
\end{prf}

The following theorem is well known for the case of separable Hilbert
spaces (cf.\ \cite{Pa65} and \cite[p.~II.29]{dlH72}). The results on direct limit
groups obtained in the preceding subsection easily permit us to extend
it to general Hilbert spaces. 

\begin{thm} \mlabel{thm2.14}
Let $\cH$ be an infinite-dimensional Hilbert space over
$\K \in \{\R, \C,\bH\}$ and $p \in [1,\infty]$. Then the following assertions hold: 
\begin{enumerate}
\item[\rm(i)] For every $k \in \N_0$, we have 
$$\pi_k(\GL_p(\cH)) \cong \pi_k(\GL_{\N}(\K)) 
\cong \indlim \pi_k(\GL_n(\K)\big).$$
\item[\rm(ii)] If $\cH_s \subeq \cH$ is an infinite-dimensional separable subspace, then
the inclusion map $\GL_p(\cH_s) \into \GL_p(\cH)$ is a weak homotopy
equivalence. 
\item[\rm(iii)] For $1 \leq p \leq q \leq \infty$, the inclusion map 
$\GL_p(\cH) \into \GL_q(\cH)$ is a weak homotopy equivalence. 
\end{enumerate}
\end{thm}

\begin{prf}  (i) Let $(e_j)_{j \in J}$ be an orthonormal basis of $\cH$. 
Then Lemma~\ref{lem2.13} 
above shows that $\cL_0(\cH) = \Spann \{ P_{e_j, e_k} \: j,k \in J\}$ is
dense $\cL_p(\cH)$. We endow $\cL_0(\cH)$ with the direct limit topology with
respect to the directed set of finite-dimensional subspaces of
$\cL_0(\cH)$. 

For the open subset $U := \GL_p(\cH) - \1 \subeq \cL_p(\cH)$ we have 
$$ U \cap \cL_0(\cH) = \GL_J(\K) - \1, $$
where $\GL_J(\K)$ is the set of those elements $g \in \GL(\cH)$ for which 
the matrix of $g- \1$ with respect to $(e_j)_{j \in J}$ 
has only finitely many non-zero entries, i.e., $g$ and $g^*$ fix all but
finitely many $e_j$. It easily follows from $(g^*)^{-1} = (g^{-1})^*$
that $g^{-1}$ has the same property. Therefore the natural
identification of $\cL_0(\cH)$ with the matrix algebra $M_J(\K)$ leads to
an identification of the group $\1 + (U \cap \cL_0(\cH))$ with 
$\GL_J(\K)$ as in Definition~\ref{def:2.10}. 

Theorem~\ref{thm-pal1} implies that if we endow $\GL_J(\K)$ with the final 
topology with respect to the subgroups $\GL_F(\K)$, $F \subeq J$ a
finite subset, then the inclusion map 
$\GL_J(\K) \to \GL_p(\cH)$ is a weak homotopy equivalence. 
Further Proposition~\ref{prop2.11} shows that we have a weak homotopy equivalence
$\GL_{\N}(\K) \into \GL_J(\K)$. 
Composition of these two weak homotopy equivalences yields a weak homotopy
equivalence. In view of Lemma~\ref{lem-pal0}(3), this proves~(i). 

 (ii) We first choose an orthonormal basis
$(e_n)_{n \in \N}$ of $\cH_s$ and then complete it to an orthonormal 
basis $(e_j)_{j \in J}$ of $\cH$. We consider the corresponding maps 
$$ \phi_1 \: \GL_{\N}(\K) \to \GL_p(\cH_s), \quad 
\phi_2 \: \GL_J(\K) \to \GL_p(\cH),  $$
$$ \phi_3 \: \GL_{\N}(\K) \to \GL_J(\K) \quad \mbox{ and } \quad 
\phi_4 \: \GL_p(\cH_s) \to \GL_p(\cH) $$
with $\phi_4 \circ \phi_1 = \phi_2 \circ \phi_3$. Since 
$\phi_1$ and $\phi_2$ are weak homotopy equivalences by the first part
of the proof, and $\phi_3$ is a weak homotopy equivalence by
Proposition~\ref{prop2.11}, $\phi_4$ also is a weak homotopy
equivalence. 

 (iii) From the elementary inclusion $\ell^p(\N, \R) \subeq \ell^q(\N, \R)$,  
we derive that $\cL_p(\cH) \subeq \cL_q(\cH)$, and Lemma~\ref{lem2.13} implies that 
$\cL_p(\cH)$ is a dense subspace. Therefore (iii) follows by applying 
Theorem~\ref{thm-pal2} to the open subset $U := \GL_q(\cH) -\1 \subeq \cL_q(\cH)$ 
which satisfies $U \cap \cL_p(\cH) = \GL_p(\cH)-\1$. 
\end{prf}

\begin{prop}
  \mlabel{prop:a.4-ne04} For every $p \in [1,\infty]$ the polar map 
\[  p \: \U_p(\cH) \times \Herm_p(\cH) \to \GL_p(\cH), \quad 
(u,X) \mapsto u e^X \] 
is a diffeomorphism. 
\end{prop}

\begin{prf} We consider the Banach-$*$ subalgebra 
$\cL_p(\cH)_+ = \C \1 + \cL_p(\cH) \subeq B(\cH).$ 
From Example~\ref{ex:IV.15-ne04} it follows that 
\[ \cL_p(\cH)_+ \cap \GL(\cH) = \cL_p(\cH)_+^\times, \] 
so that the spectrum $\Spec_p(X)$ of an element $X \in \cL_p(\cH)_+$
coincides with the spectrum $\Spec(X)$ of $X$ as an element of
$\cL(\cH)$. Therefore Lemma~\ref{lem:a.1-ne04}
 implies that, for $g \in \GL_p(\cH)$, we
have $\log(g^*g) \in \cL_p(\cH)$, and that the map 
$$ \log \: \{ g \in \GL_p(\cH) \cap \Herm_p(\cH) \: \inf \Spec(g) > 0 \} \to
\Herm_p(\cH) $$
is smooth. This implies the assertion. 
\end{prf}

\begin{cor} \mlabel{cor2.15} 
If $\cH$ is an infinite-dimensional complex Hilbert space, \break 
$1 \leq p \leq \infty$, and 
$$\GL_p(\cH,I) := \GL(\cH,I) \cap \GL_p(\cH),$$ 
then the following assertions hold: 
\begin{enumerate}
\item[\rm(1)] $\pi_k(\GL_p(\cH,I)) \cong 
\left\{   \begin{array}{cl}
\pi_k(\GL_{\N}(\R)) & \mbox{for $I^2 = \1$} \\ 
\pi_k(\GL_{\N}(\bH)) & \mbox{for $I^2 = -\1$}
\end{array}\right. $ 
\item[\rm(2)] If $\cH_s \subeq \cH$ is 
an infinite-dimensional separable $I$-invariant subspace, then
the inclusion map $\GL_p(\cH_s,I\res_{\cH_s}) \into \GL_p(\cH,I)$ is a weak homotopy
equivalence. 
\item[\rm(3)] For $1 \leq p \leq q \leq \infty$, the inclusion map 
$\GL_p(\cH,I) \into \GL_q(\cH,I)$ is a homotopy equivalence. 
\end{enumerate}
\end{cor}

\begin{prf}  We first observe that the polar decomposition of $\GL_p(\cH)$ 
(Proposition~\ref{prop:a.4-ne04}) implies that its intersection with 
$\GL_p(\cH,I)$ also has a 
polar decomposition (see the proof of Theorem~\ref{thm10.4.3}), 
hence is homotopy equivalent to 
$\U_p(\cH,I) := \U_p(\cH) \cap \GL(\cH,I)$. For $I^2 = -\1$ we have 
$\U_p(\cH,I) \cong \U_{\bH,p}(\cH)$, and for $I^2 = \1$ we get 
$\U_p(\cH,I) \cong \U_{\R,p}(\cH_\R)$, where $\cH_\R = {\{x\in \cH \: Ix = x\}}$. 
Since the group $\U_{\K,p}(\cH)$ is homotopy equivalent to $\GL_{\K,p}(\cH)$
by Theorem~\ref{thm10.4.3}, the assertions on the groups 
$\GL_p(\cH,I)$ follow from 
Theorem~\ref{thm2.14} and the existence of polar decompositions. 
\end{prf}

\section{Approximation of continuous maps by smooth ones} 
\mlabel{sec:top.4}

The goal of this section is to show that continuous 
functions $f \: M \to N$, $M$ a finite-dimensional compact manifold
and $N$ an infinite-dimensional manifold, can be approximated by smooth
functions. In particular, we shall see that every homotopy class 
 has a smooth representative. 

\subsection{General approximation results} 

In the following, $C(M,N)_c$ denotes the space $C(M,N)$ of continuous
maps $M \to N$ endowed with the compact open topology. 

\begin{thm} \mlabel{thma.1} 
Let $s \in \N \cup \{\infty\}$ and  
$M$ be a finite-dimensional $\sigma$-compact $C^s$-manifold with boundary. 
Further, let $V$ be a locally convex space and $O \subeq V$ an open subset. 
Then the subset $C^\infty(M,O)$ is dense in $C(M,O)_c$. 

If $0 \in O$, $f \in C(M,O)$ has compact support and $U$ is an open neighborhood
of $\supp(f)$, then each neighborhood of $f$ in $C(M,O)_c$ contains a
smooth function whose support is contained in $U$. 
\end{thm}

\begin{prf} (based on \cite[Th.~2.2]{Hr76}). First we observe that $M$ permits
\linebreak $C^s$-partitions of unity. The local convexity of $V$ is
crucial for the following arguments with partitions of unity. 

Let $f \in C(M,V)$ and $(U_\alpha)_{\alpha \in A}$ be a locally finite
open cover of $M$ by relatively compact subsets. Then, for each $\alpha$, 
the set $f(\oline{U_\alpha}) \subeq O$ is compact, so that there exists an 
open convex $0$-neighborhood $W_\alpha$ in $V$ with 
$$f(\oline{U_\alpha}) + W_\alpha~\subeq~O.$$ 
In view of the paracompactness of $M$, 
there exists a subordinated locally finite 
open cover $(U_i)_{i \in I}$ 
of $M$ and constant maps $g_i \: M \to V$ such that for all $y \in
U_i \cap U_\alpha$ we have 
$$ g_i(y) - f(y) \in W_\alpha. $$
Let $(\lambda_i)_{i \in I}$ be a partition of unity subordinated to
$(U_i)_{i \in I}$ and define  \break 
$g \:= \sum_i \lambda_i g_i \: M \to V.$ 
Since this sum if locally finite, it defines a \break $C^s$-function $M \to V$,  
and on $U_\alpha$ we have $g - f \in W_\alpha$ because $W_\alpha$ is convex. 
This implies that 
$$ g(U_\alpha) \subeq f(U_\alpha) + W_\alpha \subeq O. $$

If, in addition, $\supp(f)$ is compact and contained in the open set
$U$, then we may assume that each set $U_i$ is either contained
in $U$ or satisfies $U_i \cap \supp(f) = \eset$. 
For $U_i \cap \supp(f) = \eset$ we then put $g_i = 0$, and the
assertion follows. 
\end{prf} 

\begin{thm} \mlabel{thma.3} 
Let  $s\in \N \cup \{\infty\}$ and $M$ and $N$ be $C^s$-manifolds with $\dim M <
\infty$. Then $C^s(M,N)$ is dense in $C(M,N)_c$. 
Let $f \in C(M,N)$ and $K \subeq M$ closed such that $f$ is smooth on $M \setminus K$. 
Then there exists for each neighborhood $U_f$ of $f$ in $C(M,N)_c$ and each 
open neighborhood $U$ of $K$ in $M$ 
a $C^s$-function $g \: M \to N$ in $U_f$ with $f = g$ on $M \setminus U$.
\end{thm}

\begin{prf} First we need a refinement of \cite[Thm.~2.5]{Hr76}. 
Let $U \subeq \R^n$ be open, $V$ be a locally convex space,  
$O \subeq V$ open and $f \: U \to O$ a 
continuous map. Further let $K \subeq U$ be closed and $W \subeq U$ open
such that $f$ is $C^s$ on a neighborhood of the closed subset $K
\setminus W$. Then the set of all functions  
$h \in C(U,O)$ which are $C^s$ on a neighborhood of $K$ and coincide
with $f$ on $U\setminus W$ intersects every neighborhood of $f$ in 
$C(U,O)_c$. For the proof we may w.l.o.g.\ assume that $O = V$, so
that Theorem \ref{thma.1} can be used. The 
remaining arguments can be copied from \cite[Thm.~2.5]{Hr76}. 

To conclude the proof, one uses that
$M$ has a countable open cover, and then an inductive argument as in 
\cite[Thm.~2.6]{Hr76}. The argument given in \cite{Hr76} shows in particular that 
if $f$ is smooth outside of a compact subset $K$ of $M$ and $U$ an open neighborhood 
of $K$, then we find in each neighborhood 
of $f$ a smooth function $g$ which coincides with $f$ on $M \setminus U$. 
\end{prf}

\subsection{Mapping groups on compact manifolds} 

\begin{rem} \mlabel{rema.4} 
If $F$ is a locally convex space and $X$ a compact
space, then $C(X,F)$ is a locally convex space with respect to the topology
of uniform convergence (Lemma~\ref{sammelsu}(g)). 

If $U \subeq F$ is an open subset, then 
$C(X,U) = \lfloor X,U \rfloor$ is an open subset of $C(X,F)$. Now let 
$U_j \subeq F_j$, $j =1,2$, be open subsets of 
sequentially complete locally convex spaces and 
$\phi \: U_1 \to U_2$ be a smooth map. We consider the map 
$$ \phi_X \: C(X,U_1) \to C(X,U_2), \quad \gamma \mapsto \phi \circ
\gamma. $$
Then $\phi_X$ is smooth by Proposition~\ref{superctisCk}. 
\end{rem}

\begin{defn} \mlabel{defa.5}
If $K$ is a Lie group and $X$ is a compact space,
then $C(X,K)$, endowed with the topology of uniform convergence is a
Lie group with Lie algebra $C(X,\fk)$ (this follows from the proof of 
Theorem~\ref{thm:mapgro-Lie}, specialized to $r = 0$). 
The main point is to verify that inversion and 
multiplication in the canonical local charts are smooth, 
and this follows from Remark~\ref{rema.4}. 
\end{defn}

\begin{lem} \mlabel{lema.6} Let $M$ be a compact manifold and $K$ a Lie group. 
Then $C^\infty(M,K)$ is dense in $C(M,K)_c$. In particular,  
every connected component of the Lie group $C(M,K)_c$ 
contains a smooth map. Moreover, we have 
\begin{eqnarray} \label{eq:10.4.3} 
 C^\infty(M,K) \cap C(M,K)_0 =  C^\infty(M,K)_0. 
\end{eqnarray}
\end{lem}

\begin{prf} As $K$ is a topological group, the compact open topology on
$C(M,K)$ coincides with the topology of uniform convergence (Exercise~\ref{exer:4.1.1}) 
which turns $C(M,K)$ into a Lie group with Lie algebra 
$C(M,\fk)$ (Theorem~\ref{thm:mapgro-Lie}). In
particular $C(M,K)$ is locally arcwise connected, so that the
first assertion follows immediately from Theorem~\ref{thma.3}. 

To verify (\ref{eq:10.4.3}), we first observe that every smooth 
map $f \: M \to K$ which is sufficiently close to the identity is 
homotopic to the identity in $C^\infty(M,K)$ because its range lies in an 
open identity neighborhood diffeomorphic to an open convex set. 
Now homogeneity implies that $C^\infty(M,K)$ is locally connected 
with respect to the compact open topology, and hence that its 
connected components are also open in the coarser compact open topology. 
This implies that the connected components of $C^\infty(M,K)$ are closed 
in the compact open topology, and therefore that the closure of 
$C^\infty(M,K)_0$ is open and closed in $C(M,K)$, hence coincides with 
$C(M,K)_0$ and satisfies (\ref{eq:10.4.3}). 
\end{prf}

\subsection{Groups of compactly supported functions} 

Throughout this subsection $M$ is a paracompact finite-dimensional smooth manifold 
and $K$ is a Lie group with Lie algebra~$\fk$. 

\begin{lem} \mlabel{lemc2a.1} For 
each compact subset $Y$ of $C^\infty_c(M,K)$, there exists a compact 
subset $X \subeq M$ with 
\[ Y \subeq C^\infty_X(M,K) := \{f \in C^\infty(M,K) \: f(M\setminus X) \subeq \{\be\}\}.\]
\end{lem}

\begin{prf} Let $U \subeq \fk$ be an open $0$-neighborhood, and 
$\phi \: U \to \phi(U)$ a chart with $\phi(0) = \be$. 
Then there exists an open $0$-neighborhood $U_0 \subeq U$ such that we 
obtain a local chart for $G := C^\infty_c(M,K)$ by 
$\phi_G(f) := \phi  \circ f$ (Theorem~\ref{thm:mapgro-Lie}). Let 
$$V := \{ f \in C^\infty_c(M;\fk) \: f(M) \subeq U_0\}$$ and observe that 
\[ \phi_G(V) = \{ f \in C^\infty_c(M,K) 
\: f(M) \subeq \phi(U_0)\}. \]
Then, for each $f \in G$, the set $f \phi_G(V)$ is an open neighborhood, and the map 
$$ \phi_f \: V \to f \phi_G(V), \quad \xi \mapsto f \phi_G(\xi) $$
is a diffeomorphism. Let $W \subeq V$ be a closed $0$-neighborhood 
such that $\phi_G(W) \phi_G(W) \subeq \phi_G(V)$. Since 
$\oline{\phi_G(W)}$ is the intersection of all sets $\phi_G(W) N$,
where $N$ is an identity neighborhood in $C^\infty_c(M,K)$, 
$\oline{\phi_G(W)} \subeq \phi_G(V)$, so that the closedness of $W$
implies that $\phi_G(W)$ is closed. 

Since the compact set $Y$ is covered by the open 
sets $f \phi_G(W^0)$, $f \in Y$, there exist $f_1, \ldots, f_n \in Y$
with 
$$ Y \subeq f_1 \phi_G(W^0) \cup \ldots \cup f_n \phi_G(W^0). $$
The closedness of $\phi_G(W)$ implies that each set 
$Y \cap f_j \phi_G(W)$ is compact, so that, for each $j$, the closed set 
$$ \phi_{f_j}^{-1}(Y \cap f_j \phi_G(W))  = W \cap \phi_{f_j}^{-1}(Y) 
\subeq C^\infty_c(M;\fk) = \indlim C^\infty_X(M;\fk)$$
is compact, which implies the existence of a compact subset $X_j \subeq M$ with 
$\phi_{f_j}^{-1}(Y \cap f_j \phi_G(W)) \subeq C^\infty_{X_j}(M;\fk)$ 
(Lemma~\ref{lem-herve}). 
Let 
$$X := X_1 \cup \ldots \cup X_n \cup \supp(f_1) \cup \ldots \cup \supp(f_n).$$ 
Then $X$ is compact and $Y \subeq C^\infty_X(M,K)$. 
\end{prf}

\begin{lem} \mlabel{lemc2a.2} 
Let $Y$ be a compact space and $f \: Y \to
C^\infty_c(M,K)$ 
a continuous map. Then there 
exists a compact subset $X \subeq M$ and a continuous map 
$f_X \: Y \to C^\infty_X(M,K)$ such that 
$f = \eta_X \circ f_X$ holds for the inclusion map 
$\eta_X \: C^\infty_X(M,K) \to C^\infty_c(M,K)$. 
\end{lem}

\begin{prf} From the fact that $C^\infty_c(M;\fk)$ is a strict 
inductive limit of spaces 
$C^\infty_X(M; \fk)$, 
$X \subeq M$ compact, and the description of the natural charts of the
Lie group $C^\infty_c(M,K)$, we derive that, for each compact subset 
$X \subeq M$, the inclusion map 
$C^\infty_X(M,K) \into C^\infty_c(M,K)$ is a topological embedding. 

Since $C^\infty_c(M,K)$ is Hausdorff, the set $f(Y)$ is compact. 
In view of Lemma~\ref{lemc2a.1}, there exists a compact subset 
$X \subeq M$ with $f(Y) \subeq C^\infty_X(M,K)$. 
Let $f_X \: Y \to C^\infty_X(M,K)$ denote the corestriction of $f$ to $C^\infty_X(M,K)$. 
Since $\eta_X$ is a topological embedding, the map 
$f_X$ is continuous. It obviously satisfies $f = \eta_X \circ f_X$. 
\end{prf}

\begin{prop} \mlabel{propc2a.3} 
For $(X_n)_{n \in \N}$ be an exhaustion of $M$, i.e., 
$X_n \subeq M$ is  compact with 
$X_n \subeq X_{n+1}^0$ and $M = \bigcup_n X_n$. 
Then the map 
$$ \indlim C^\infty_{X_n}(M,K) \to C^\infty_c(M,K) $$
is a weak homotopy equivalence. In particular,  
\[ \pi_m(C^\infty_c(M,K)) \cong \indlim \pi_m(C^\infty_{X_n}(M,K)) \quad \mbox{ 
for every } \quad m \in \N_0.\] 
\end{prop}

\begin{prf} Lemma~\ref{lemc2a.2} implies that each continuous map 
$f \: \bS^m \to C^\infty_c(M,K)$ factors through some inclusion 
$C^\infty_{X_n}(M,K) \to C^\infty_c(M,K)$. 
If two such maps $f_1, f_2$ are homotopic, then each homotopy 
$h \: \bS^m \times [0,1] \to C^\infty_c(M,K)$ also factors through 
some group $C^\infty_{X_k}(M,K)$. 
This implies that the natural map 
$$ \indlim \pi_m(C^\infty_{X_n}(M,K))
 \cong \pi_m(\indlim C^\infty_{X_n}(M,K)) \to \pi_m(C^\infty_c(M,K)) $$
is bijective, i.e., that the continuous map
$\indlim C^\infty_{X_n}(M,K) \to C^\infty_c(M,K)$ is a weak homotopy equivalence. 
\end{prf}

\begin{rem} \mlabel{remc2a.4} 
A similar argument as in the proof of 
Proposition~\ref{propc2a.3} shows that the map 
$$ \indlim C_{X_n}(M,K) \to C_c(M,K) $$
is a weak homotopy equivalence. 
\end{rem}

Next we show that the map 
\[C^\infty_c(M,K) \to C_c(M,K) \] 
is a weak homotopy equivalence, so that the homotopy groups of $C^\infty_c(M,K)$ 
are the limits of the corresponding homotopy groups of 
$C_X(M,K)$, where $X \subeq M$ is a compact submanifold with boundary. 
These groups are more approachable since they are isomorphic to 
$C_*(X/\partial X,K)$, where $X/\partial X$ is a compact space, with 
the image of $\partial X$ as the base point. 

If $M$ is a compact manifold with boundary, then the homotopy 
groups $\pi_m(C_*(M/\partial M,K))$ might be well accessible 
(cf.~Section~\ref{sec10.6}). 
If $\partial M$ is empty, then we identify the group 
$C_*(M/\partial M,K)$ with $C(M,K)$. 

\begin{lem} \mlabel{lemc2a.5} Let $X_1, X_2 \subeq M$ be compact subsets with 
$X_1 \subeq X_2^0$ and $f \in C_{X_1}(M,K)$. Then every neighborhood 
of $f$ contains a map $f'$ in $C^\infty_{X_2}(M,K)$. The image of the homomorphism 
$$ \eta \: \pi_0(C^\infty_{X_2}(M,K)) \to \pi_0(C_{X_2}(M,K)) $$
contains the image of $\pi_0(C_{X_1}(M,K))$. 
Moreover, if $f$ is contained in the identity component 
$C_{X_1}(M,K)_0$, then we may choose 
$f' \in C^\infty_{X_2}(M,K)_0$. 
\end{lem}

\begin{prf} The first assertion follows from Theorem~\ref{thma.3}. 
Since the groups $C_X(M,K)$ and $C^\infty_X(M,K)$ are Lie groups, 
their connected components are open, so that   
every connected component of $C_{X_2}(M,K)$ meeting 
$C_{X_1}(M,K)$ contains a smooth element. 

If the map $f \in C_{X_1}(M,K)$ is sufficiently close to $\be$ in the sense that 
$f(M) \subeq V$ for some $\be$-neighborhood $V \subeq K$
diffeomorphic to an open convex set, we find 
$f_1 \in C^\infty_{X_2}(M,K)$ with $f_1(M) \subeq V$. 
Any two smooth maps $f_1, f_2 \in  C^\infty_{X_2}(M,K)$ with $f_j(M) \subeq V$ 
are smoothly homotopic, hence contained in the same connected component of 
$C^\infty_{X_2}(M,K)$. 

If $f \in C_{X_1}(M,K)$ is contained in the 
identity component, then there exists a continuous curve 
$\gamma \: [0,1] \to C_{X_1}(M,K)$ with $\gamma(0) = \be$ and $\gamma(1) = f$. 
For a sufficiently fine subdivision 
$0 = t_0 < t_1 < \ldots < t_N = 1$, we now find smooth maps 
$f_j \in C^\infty_{X_2}(M,K)$ close to $\gamma(t_j)$ in the sense that
$(f_j^{-1} \cdot \gamma(t_i))(M) \subeq V$, where for $j < N$ the maps 
$f_j$ and $f_{j+1}$ are smoothly homotopic. Hence 
$f_N$ is contained in the identity component of $C^\infty_{X_2}(M,K)$. 
\end{prf}

\begin{lem} \mlabel{lemc2a.6} 
The map $\iota \: C^\infty_c(M,K) \to C_c(M,K)$ induces an isomorphism 
$$ \pi_0(\iota) \: \pi_0(C^\infty_c(M,K)) \to \pi_0(C_c(M,K)). $$
\end{lem}

\begin{prf} The surjectivity of $\pi_0(\iota)$ follows directly from 
Lemma~\ref{lemc2a.5}. If $f \in C^\infty_c(M,K)$ satisfies $[f] \in \ker \pi_0(\iota)$, 
then there exists a compact subset $X \subeq M$ and a continuous map 
$\gamma \: [0,1] \to C_X(M,K)$ with $\gamma(0) = \be$ and $\gamma(1) =
f$ (Lemma~\ref{lemc2a.2}). 
Let $Y \subeq M$ be a compact subset with $X \subeq Y^0$. Then 
Lemma~\ref{lemc2a.5} implies that we can approximate $f$ by smooth functions 
$f'$ in the identity component of $C^\infty_Y(M,K)$. It follows in particular 
that $f$ is contained in the identity component of $C^\infty_Y(M,K)$, hence also in 
the identity component of $C^\infty_c(M,K)$. This shows that $\pi_0(\iota)$ is 
injective. 
\end{prf}

Fix a base point $x_M \in M$ 
and in any Lie group $K$ 
we consider the unit element $\be$ as the base point. We
write $C^\infty_*(M,K) \subeq C^\infty(M,K)$ for the subgroup of base
point preserving maps and observe that 
\begin{equation}
  \label{eq:semidir-map}
C^\infty(M,K) \cong C^\infty_*(M,K) \rtimes K 
\end{equation}
as Lie groups, where we identify $K$ with the subgroup of
constant maps. This relation already leads to 
\begin{eqnarray} \label{eq:10.4.1}
 \pi_k(C^\infty(M,K)) \cong \pi_k(C^\infty_*(M,K)) \times \pi_k(K)
\quad \mbox{ for } \quad k \geq 1 
\end{eqnarray}
and 
\begin{eqnarray} \label{eq:10.4.1b}
 \pi_0(C^\infty(M,K)) \cong \pi_0(C^\infty_*(M,K)) \rtimes \pi_0(K).
\end{eqnarray}

\begin{rem} If $K$ is a topological group, then the group 
$\pi_0(K)$ of arc-components acts naturally by automorphisms 
on the groups $\pi_k(K)$, $k \in \N_0$, via 
\[ [g].[\sigma] := [c_g \circ \sigma], \quad \mbox{ where } \quad 
c_g(h) = ghg^{-1}.\] 

A typical example where this action is non-trivial is the group 
$K = \OO_2(\R) \cong \SO_2(\R) \rtimes \{\be,\sigma\}$, where 
$\sigma$ is a linear reflection. Then 
$\pi_0(K) = \{ \be, [\sigma]\}$ and $[\sigma]$ acts 
on $\pi_1(K) \cong \Z$ by inversion. 
\end{rem}

\begin{rem} (Action of $\pi_0(K)$ on higher homotopy groups) 
For each topological group $K$, the arc-component $K_a$ 
of the identity, and each $k \in \N$, we have
\begin{eqnarray} \label{eq:10.4.2} 
\pi_k(K) 
\cong \pi_0(C_*(\bS^k,K)) = \pi_0(C_*(\bS^k,K_a)) = \pi_0(C(\bS^k,K_a)). 
\end{eqnarray} 
In Lemma~\ref{lem-trivact} we shall see that 
$\pi_0(C_*(\bS^k,K))$ acts trivially 
on $\pi_\ell(C_*(\bS^k,K))$ for $k \geq 1$ and $\ell \in \N_0$. Therefore 
\eqref{eq:10.4.2} implies that the corresponding action of 
$\pi_k(K)$ on $\pi_m(K)$ is trivial for $m \geq k \geq 1$. 
\end{rem}

The following theorem is one of the two main results of this section. It provides
a valuable tool to determine the homotopy groups of groups of 
compactly supported smooth
maps in terms of the corresponding groups of continuous maps. 

\begin{thm} \mlabel{thmc2a.7} 
If $M$ is a connected $\sigma$-compact finite-dimensional 
manifold and $K$ a Lie group, then the inclusion 
$C^\infty_c(M,K) \to C_c(M,K)$ is a weak homotopy
equivalence. If $M$ is compact and $x_M \in M$ is a base point, then the inclusion 
\begin{eqnarray} \label{eq:10.4.3b} 
C^\infty_*(M,K) \to C_*(M,K) := \{ f \in C(M,K) \: f(x_M) = \be\} 
\end{eqnarray}
is a weak homotopy equivalence.
\end{thm}

\begin{prf} We have to show
that the inclusion induces an isomorphism 
\[  \pi_k(C^\infty_c(M,K)) \to \pi_k(C_c(M,K)) \quad \mbox{ for each } \quad 
k \in \N_0.\] 
For $k = 0$ this is Lemma~\ref{lemc2a.6}. 
If $M$ is compact, then 
\[  \pi_0(C^\infty_c(M,K)) = \pi_0(C^\infty(M,K)) \cong
\pi_0(C_*^\infty(M,K)) \rtimes \pi_0(K) \] 
and 
\[ \pi_0(C_c(M,K)) = \pi_0(C(M,K)) \cong \pi_0(C_*(M,K)) \rtimes \pi_0(K), \] 
so that (\ref{eq:10.4.3b}) follows from Lemma~\ref{lemc2a.6}. We only observe that, if 
$(f_t)_{0 \leq t \leq 1}$ 
is a homotopy between $f_0$ and $f_1$ in $C^\infty_c(M,K)$ and 
$x_M \in M$ is a base point, then
$f_t(x) f_t(x_M)^{-1}$ is a homotopy between $f_0$ and $f_1$ in
$C^\infty_*(M,K)$. 

Next we assume that $k \geq 1$ and observe that the inclusions 
\begin{eqnarray*}
&&\ \ \ \ C_*(\bS^k, C^\infty_c(M,K)) = C_*(\bS^k, C^\infty_c(M,K)_a) 
\into C(\bS^k, C^\infty_c(M,K)_a) \\
&& \into 
C(\bS^k, C_c(M,K)_a)
 \into C(\bS^k, C_c(M,K)) \cong C_c(\bS^k \times M,K)    
\end{eqnarray*}
(cf.\ Lemma~\ref{lemc2a.1}) 
are continuous homomorphisms of Lie groups, where 
$$ C(\bS^k, C_c(M,K)_a) \into C(\bS^k, C_c(M,K)) $$
is an open embedding. For the group of connected
components, we obtain for $k \geq 1$ with (\ref{eq:10.4.2}) the homomorphisms 
\begin{eqnarray*}
 \pi_k(C^\infty_c(M,K)) 
&&\cong \pi_0\big(C_*(\bS^k, C^\infty_c(M,K))\big)
\cong \pi_0\big(C(\bS^k, C^\infty_c(M,K)_a)\big)\cr
&&\to \pi_0\big(C(\bS^k, C_c(M,K)_a)\big)\cong
\pi_k\big(C_c(M,K)\big).
\end{eqnarray*}

If $f \: \bS^k \times M \to K$ is a continuous map with compact
support corresponding to an element of 
$C_*(\bS^k; C_c(M,K)_a)$, then 
Lemma~\ref{lemc2a.5} first implies that every neighborhood of 
$f$ contains a smooth map with compact support. Thus 
every connected component of $C_c(\bS^k \times M,K)$ contains an element
of $C(\bS^k, C^\infty_c(M,K))_a$ by the openness argument from above. 
This means that the homomorphism 
$\pi_k(C^\infty_c(M,K)) \to \pi_k(C_c(M,K))$ is surjective. 

To see that it is injective, suppose that 
$\sigma \in C\big(\bS^k, C^\infty_c(M,K)_a\big)$ satisfies 
$\sigma \in C\big(\bS^k, C_c(M,K)_a\big)_a \cong C_c(\bS^k \times
M,K)_a$. From Lemma~\ref{lemc2a.6}, we obtain 
$$ C^\infty_c(\bS^k \times M,K) \cap C_c(\bS^k \times M,K)_a \subeq
C^\infty_c(\bS^k \times M,K)_a, $$
so that approximating $\sigma$ by elements in 
$C^\infty_c(\bS^k \times M,K)$ (Lemma~\ref{lemc2a.5}), we may even
approximate it by elements in $C^\infty_c(\bS^k \times M,K)_a$, which
implies that $\sigma$ lies in the identity component of 
$C\big(\bS^k, C^\infty_c(M,K)_a\big)$. This proves that the
homomorphisms 
$\pi_k(C^\infty_c(M,K)) \to \pi_k(C_c(M,K)),$ $k \in \N_0,$
are isomorphisms. 
\end{prf}

\begin{rem}
Theorem~\ref{thmc2a.7} can also be extended to non-connected 
$\sigma$-compact manifolds $M$ as
follows. Let $M = \bigcup_{j \in J} M_j$ be the decomposition of $M$
into connected components $M_j$. Here one can use 
\begin{equation}
  \label{eq:comp-dir-sum}
C_c(M,K) = \bigoplus_{j \in J} C_c(M_j,K), 
\end{equation}
and, for each compact subset $X \subeq M$, we have the finite sum
decomposition 
$$ C_X(M,K) = \bigoplus_{X \cap M_j \not= \eset} C_{X \cap M_j}(M_j,K). $$
If $M$ has infinitely many (and thus at most countably many) 
connected components, then the right hand side 
of \eqref{eq:comp-dir-sum} carries the Lie group topology 
corresponding to the locally convex direct sum topology on $C_c(M;\fk)$ 
(cf.\ Theorem~\ref{thm:IV.1.12}). 
\end{rem}

\begin{rem} If 
$K$ is a connected Fr\'echet--Lie group and $M$ a compact manifold, 
then $C^\infty(M,K)$ also is a 
Fr\'echet--Lie group (Theorem~\ref{thm:mapgro-Lie}), 
so that combining Theorem~\ref{pal-contr} 
with Theorem \ref{thmc2a.7} implies that the inclusion
$C^\infty(M,K) \into  C(M,K)$ is a homotopy equivalence. 
\end{rem}

Lemma~\ref{lemc2a.8} and Proposition~\ref{propc2a.9} 
below provide additional information on the homotopy type of the groups of continuous maps. 

\begin{lem} \mlabel{lemc2a.8} 
If $M$ is a locally compact space, then 
the inclusion $$\eta \: C_c(M,K) \to C_0(M,K)$$ 
induces an isomorphism 
$\pi_0(C_c(M,K)) \to \pi_0(C_0(M,K))$. 
\end{lem}

\begin{prf} Let $f \in C_0(M,K)$. Then there exists a compact subset $X
\subeq M$ such that $f(M\setminus X)$ is contained in an identity
neighborhood $U_K \subeq K$ for which there exists a convex open 
$0$-neighborhood $U \subeq \fk$ and a diffeomorphism 
$\phi \: U \to U_K$ with $\phi(0)= \be$. 
Using a continuous function $h \in C_c(M;\R)$ which is $1$ on $X$ and 
satisfies $h(M) \subeq [0,1]$, we define a function 
$f_h\in C_c(M,K)$ by $f_h = f$ on $X$ and 
$f_h(m) = \phi\big(h(m) \phi^{-1}(f(m))\big)$ 
on $M\setminus X$. Then 
\[  F \: M \times [0,1] \to K, \quad 
F(x,t) := f_{t + (1-t)h}(x) \] 
is a homotopy between $f = f_1$ and $f_h$, and we see that 
$\pi_0(\eta)$ is surjective. 

A similar argument shows that, for $f,g \in C_c(M,K)$, any path joining
$f$ and $g$ in $C_0(M,K)$ can be deformed to a path lying completely inside of
$C_X(M,K)$ for a compact subset $X$ of $M$. Therefore $\pi_0(\eta)$ is
injective. 
\end{prf}

In the following, we write $M_\infty= M \cup \{ \infty\}$ for the 
\index{one-point compactification} 
{\it one-point compactification} of $M$. 
This is the uniquely determined compact space 
containing $M$ as an open subspace which is the complement of one point, 
denoted $\infty$. Note that 
\begin{equation}
  \label{eq:13.7}
 C_0(M,K) \cong C_*(M_\infty, K) 
\end{equation}
for every topological group~$K$. 

\begin{prop} \mlabel{propc2a.9} 
If $M$ is a locally compact space, then 
the inclusion $\eta \: C_c(M,K) \to C_0(M,K)$ is a weak homotopy
equivalence. 
\end{prop}

\begin{prf} For every compact space $X$, we have an embedding of
topological groups  
$$ C(X,C_0(M,K)) \cong C(X,C_*(M_\infty,K)) \into C(X,C(M_\infty,K))
\cong C(X \times M_\infty,K), $$
which leads to the isomorphism 
\[ C(X,C_0(M,K)) \cong C_0(X \times M,K). \] 

In view of Lemma~\ref{lemc2a.8}, there exists for each 
$f \in C_0(X \times M,K)$ some compact subset $Y \subeq M$ and a
continuous map $f_Y \in C(X,C_Y(M,K)) \subeq C(X \times Y,K)$
homotopic to $f$. The same argument applies to $[0,1] \times X$
instead of $X$, so that we see that the inclusion 
$C_c(M,K) \to C_0(M,K)$ induces a bijection 
$[X,C_c(M,K)] \to [X,C_0(M,K)]$ on the level of homotopy classes. 

Applying this to $X := \bS^k$, $k \in \N$, we obtain with Lemma~\ref{lemc2a.8} that
the natural map 
\begin{align*}
\pi_k(C_c(M,K)) 
&\cong [\bS^k,C_c(M,K)]_* 
\cong [\bS^k,C_c(M,K)_a] \\ 
&\to 
[\bS^k,C_0(M,K)_a] \cong [\bS^k,C_0(M,K)]_* 
\cong \pi_k(C_0(M,K)) 
\end{align*}
is bijective, hence an isomorphism of groups. 
\end{prf}

Combining Proposition~\ref{propc2a.9} with Theorem~\ref{thmc2a.7}, 
we obtain: 

\begin{thm} \mlabel{thmc2a.10} 
For each $\sigma$-compact connected finite-dimensional
manifold $M$ and each Lie group $K$, the inclusion map 
$$ C^\infty_c(M,K) \to C_0(M,K) \cong C_*(M_\infty,K) $$ 
is a weak homotopy equivalence. 
\end{thm}

\begin{ex} \mlabel{exc2a.11} 
For $M = \R^n$ with $M_\infty \cong \bS^n$, 
we obtain with Theorem~\ref{thmc2a.10} for each
$k \in \N_0$: 
$$ \pi_k(C^\infty_c(\R^n,K)) 
\cong \pi_k(C_*(\R^n_\infty,K))
\cong \pi_k(C_*(\bS^n,K))
\cong \pi_{k+n}(K). $$
\end{ex}

\begin{lem} \mlabel{lemc2a.12} 
Let $M$ and $N$ be finite-dimensional smooth manifolds and 
{$\phi \: N \to M$} be a smooth proper map. 
\begin{description}
\item[\rm(i)] The map 
$$ \phi_K \: C^\infty_c(M,K) \to C^\infty_c(N,K), \quad f \mapsto f \circ \phi $$
is a morphism of Lie groups. 
\item[\rm(ii)] Let $\phi_\infty \: M_\infty \to N_\infty$ denote the continuous 
extension of $\phi$ to the one-point compactifications. 
Then, for each $k \in \N_0$, the map 
$$ \pi_k(\phi_K) \: \pi_k(C^\infty_c(M,K)) \to \pi_k(C^\infty_c(N,K)) $$
only depends on the homotopy class of $\phi_\infty$ in 
the set $[M_\infty, N_\infty]_*$ of pointed homotopy classes. 
\end{description}
\end{lem}

\begin{prf} (i) It is clear that $\phi_K$ maps $C^\infty_c(N,K)$ into $C^\infty_c(M,K)$ and that 
it is a group homomorphism. It therefore suffices to show smoothness in some 
identity neighborhood. 

Let $U \subeq K$ be an open identity neighborhood and 
$\psi \: U \to W$ be a chart of $K$, where $W \subeq \fk$ is 
an open subset and $\psi(\be) = 0$. Then there exists an open $0$-neighborhood 
$V \subeq W$ such that 
\[ C^\infty_c(N,\psi^{-1}(V)) := \{ f \in C^\infty_c(N,K) \: f(N) \subeq \psi^{-1}(V)\}\] 
is an open subset of 
$C^\infty_c(N,K)$ (Theorem~\ref{thm-tefu-gp}). 
 Now it suffices to see that the map 
\[  C^\infty_c(M,V) \to C^\infty_c(N,V), \quad f \mapsto f \circ \phi  \] 
is smooth. As this map is the restriction of a linear map, we only have to 
show that it is continuous. 

For each compact subset $X \subeq M$, we have 
$$ C^\infty_X(M,K) \circ \phi \subeq C^\infty_{\phi^{-1}(X)}(M,K), $$
so that the assertion follows from the observation that, for each $n \in \N$, 
the map $d^n(f \circ \phi)$ depends continuously on $f$, 
when considered as an element 
of $C(T^n(N), \fk)_c$
(cf.~Proposition~\ref{Crpuba}). 

 (ii) Let 
$\eta_M \: C^\infty_c(M,K) \to C_*(M_\infty,K)$ denote the natural inclusion. 
Then $\eta_N \circ \phi_K = \tilde\phi_K \circ \eta_M$ holds with 
\[ \tilde\phi_K \: C_*(M_\infty,K) \to C_*(N_\infty,K), \quad f \mapsto f \circ \phi_\infty. \]

We know from Theorem~\ref{thmc2a.10} that the maps $\eta_M$ and $\eta_N$ are weak homotopy 
equivalences. 
Therefore it suffices to show that the maps $\pi_k(\tilde\phi_K)$ only depend on 
the homotopy class of $\phi$. If $\phi, \psi \: M \to N$ are proper and smooth such that 
$\phi_\infty$ and $\psi_\infty$ are homotopic, then it is easy to see that the maps 
$\tilde\phi_K$ and $\tilde\psi_K$ are homotopic, hence induce the same homomorphisms 
on homotopy groups. 
\end{prf}

\subsection{Holomorphic mapping groups} 

Let $M$ be a {\it Stein manifold}, \index{Stein manifold}  
i.e., a complex manifold which can be realized as a closed 
submanifold of some $\C^n$. 
If $G$ is a complex Banach--Lie group, then the groups 
$C(M,G)$ and $\Hol(M,G)$ are metrizable 
topological groups with respect to the topology of 
uniform convergence on compact subsets of~$M$. 
In general, these groups are not 
Lie groups and it is an interesting open problem to characterize those 
Stein manifolds $M$ for which they are. We have a natural inclusion map 
\[  \eta \: \Hol(M,G) \into C(M,G)_{\rm c.o.},  \] 
and one can show that this inclusion is a weak homotopy equivalence. 
This can be reduced to 
the Oka Principle \cite[Satz I, p.~268]{Gr58a} 
which asserts that the inclusion $\eta$ induces a bijection on the level 
of connected components. Further, one uses that, for each $k \in \N$, the group 
$C(\bS^k,G)$ is also a complex Banach--Lie group, so that Oka's 
Principle applies to the 
topological group 
\[ \Hol(M,C(\bS^k,G)) \cong C(\bS^k,\Hol(M,G)).\] 

These results are of particular interest if $M = \Sigma \setminus P$, where 
$\Sigma$ is a compact Riemann surface and $P \subeq \Sigma$ is a finite set 
(cf.\ Subsection~\ref{subsec:13.8.4} below).

\section{Concrete mapping groups} \mlabel{sec10.6} 

In this section we take a closer look at mapping groups on 
products of spheres, in particular tori, and on compact 
surfaces with boundary. 

Let $K$ be a connected Lie group and $M$ a compact
connected manifold. We write $C^\infty(M,K)$ for the
corresponding 
\index{ mapping group} 
{\it mapping group}. In the preceding subsection, we have seen that the 
inclusion 
\[  C^\infty(M,K) \into C(M,K) \]
is a weak homotopy equivalence (Theorem~\ref{thmc2a.7}), 
so that the determination of homotopy groups of 
smooth mapping groups can be carried out for groups of continuous maps, 
and this is much easier. For compactly supported smooth maps, we likewise 
have a weak homotopy equivalence 
\[ C^\infty_c(M,K) \into C_*(M_\infty,K). \] 
Motivated by these reductions, we study in 
this section the homotopy groups of topological groups 
of the form $C_*(X,K)$, where $K$ is a topological group and $(X,x_0)$ 
is a pointed compact space. 

\subsection{Groups of connected components} 
\mlabel{subsec:14.8.1} 

We start this section with an investigation of the group 
\[ [X,K]_* = \pi_0(C_*(X,K)) \] 
of arc-components and write $K_a$ for 
the arc-component of the identity in $K$. 
For $k \in \N$,  we then have 
\begin{eqnarray} \label{eq:10.6.3} 
 \pi_k(K) = [\bS^k,K]_* 
\cong \pi_0(C_*(\bS^k,K)) \cong \pi_0(C_*(\bS^k,K_a)).   
\end{eqnarray} 

For two pointed spaces $(X,x_0)$ and $(Y,y_0)$, we define their 
\index{smash product, of pointed spaces}
\index{wedge sum, of pointed spaces}
{\it smash product} by 
\[ X \wedge Y := (X \times Y)/(X \vee Y),\] 
where 
\begin{equation}
  \label{eq:vee-product}
X \vee Y := (X \times \{y_0\}) \cup (\{x_0\} \times Y) 
\end{equation}
is their {\it wedge sum}. 
For $X,Y$ locally compact, we then obtain a homeomorphism 
\[ C_*(X, C_*(Y, Z)) \cong C_*(X \wedge Y, Z) \] 
from the homeomorphism 
\[ C(X, C(Y, Z)) \cong C(X \times Y, Z) \] 
(Proposition~\ref{ctsexp}). 

For a pointed topological space $(X,x_0)$, let 
$SX := \bS^1 \wedge X$ denote the 
\index{suspension} 
{\it suspension of $X$}. Iterating 
this process $k$ times leads to 
\[ S^k X = (\bS^1 \wedge \cdots \wedge \bS^1) \wedge X \cong \bS^k \wedge X \]  
(\cite{Bre93}). In view of (\ref{eq:10.6.3}), we have 
\begin{equation}
  \label{eq:13-8}
\pi_k(C_*(X,K)) 
\cong \pi_0\big(C_*(\bS^k,C_*(X,K))\big) \cong 
[\bS^k \wedge X, K]_*.
\end{equation}

This discussion already shows that 
it is important to have more information on spaces 
which are products of a sphere with another space. 
The following lemma provides some insights in this situation. 

\begin{lem} \mlabel{lem-trivact} Let $K$ be a topological group, 
$(X,x_0)$ be a pointed compact space, and $p \geq 1$. 
Then , for $G := C_*(\bS^p, K)$, the conjugation actions of the group 
$\pi_0(G) = [\bS^p,K]_*$ on $[X,G]_*$ and $[X,G]$ defined by 
\[ [\sigma].[f] := [\sigma * f], \quad 
(\sigma * f)(x,y) := \sigma(y) f(x,y) \sigma(y)^{-1} \] 
for $\sigma \in C_*(\bS^p,K), f \in C(X,G) \subeq C(X \times \bS^p,K),$ 
are trivial. 
\end{lem}

\begin{prf} Let $\sigma \in G$ and $f \: X \to G = C_*(\bS^p,K)$ 
be a continuous map. 
We identify maps $X \to C_*(\bS^p,K)$ with maps 
$X \times \bS^p \to K$ vanishing on 
$X \times \{y_0\}$, where $y_0\in \bS^p$ is the base point. 

As $p \geq 1$, we may choose a representative of $[\sigma] \in [\bS^p,K]$ supported in 
the upper hemisphere of $\bS^p$ and 
a representative of $[f] \in [X,G]$ supported in the product of $X$ and 
the lower hemisphere of $\bS^p$. Then 
$\sigma * f = f$, and this implies the triviality of the action. 
\end{prf}

The following theorem provides valuable information on the groups $[X,K]$ 
(cf.\ \cite[Th.~X.3.6]{Whi78}): 

\begin{thm} \mlabel{thm:whitehead} 
Let $X= \bS^{p_1} \times \cdots \times \bS^{p_n}$, $p_j > 0$, 
be a product of $n$ spheres and $K$ an arcwise 
connected topological group. Then the following assertions hold: 
\begin{itemize}
\item[\rm(1)]\ \ $[X,K] \cong \pi_0(C_*(X,K))$ is nilpotent of class $< n+1$. 
\item[\rm(2)]\ \ $[X,K]$ acts trivially on $\pi_m(C_*(X,K))$ for 
$m \geq 1$.  
\end{itemize}
\end{thm}

\begin{prf} We proceed by induction on the number $n$  of factors in $X$. 

(a) First we assume that $X = \bS^p$ is a $p$-dimensional sphere. 
Since $K$ is arcwise connected, $[X,K] = [X,K]_* \cong \pi_p(K)$ is abelian, 
i.e., nilpotent of length $< 2$. Moreover, for $p > 0$, the group 
$[X,K] \cong \pi_p(K) = \pi_0(C_*(\bS^p,K))$ acts trivially on 
$\pi_m(C_*(\bS^p,K)) \cong [\bS^m,C_*(\bS^p,K)]_*$ by Lemma~\ref{lem-trivact}. 

(b) Now we assume that the assertion holds for all products 
$X$ of $n$ spheres and consider the group 
$[X \times \bS^m, K] = [X \times \bS^m, K]_*$. 
Let $x_0 \in X$ and $y_0 \in \bS^m$ be base points, and choose 
$(x_0, y_0)$ as the base point in $X \times \bS^m$. 

Let $p_1\: X \times \bS^m \to X$ and $p_2 \: X \times\bS^m \to \bS^m$ 
denote the projection maps. Then we have continuous group homomorphisms 
$$ p_1^* \: C_*(X,K) \to C_*(X \times \bS^m,K), \quad 
p_2^* \: C_*(\bS^m,K) \to C_*(X \times \bS^m,K). $$
On the other hand, the inclusion maps 
$$ q_1 \: X \to X \times \bS^m, \quad x \mapsto (x,y_0), \quad 
q_2 \: \bS^m \to X \times \bS^m, \quad y \mapsto (x_0,y) $$
yield continuous homomorphisms 
$$ q_1^* \: C_*(X \times \bS^m,K) \to C_*(X,K), \quad 
q_2^* \: C_*(X \times \bS^m,K) \to C_*(\bS^m,K), $$
and these maps satisfy the relations 
\[  q_1^* \circ p_1^* = \id, \quad 
q_2^* \circ p_2^* = \id, \quad 
q_1^* \circ p_2^* = \be, \quad q_2^* \circ p_1^* = \be. \] 

We thus have a short exact sequence of topological groups 
\[ \1  \to C_*(X \wedge \bS^m, K) \into C_*(X \times \bS^m,K) 
\smapright{(q_1^*, q_2^*)} C_*(X,K) \times C_*(\bS^m,K) \to \1 \] 
with a continuous section, given by $p_1^* \cdot p_2^*$. 
On the level of connected components, this topologically split short 
exact sequence of groups induces 
a short exact sequence 
\begin{eqnarray}
  \label{eq:ex}
\1 \to [X \wedge \bS^m, K] \into [X \times \bS^m,K]  
\smapright{(q_1^*, q_2^*)} [X,K] \times [\bS^m,K] \to \1.\quad 
\end{eqnarray}
The groups $[X \wedge \bS^m, K] \cong \pi_m(C_*(X,K))$ and 
$[\bS^m,K] \cong \pi_m(K)$ are abelian, 
and the induction hypothesis implies that  
$[X,K]$ is nilpotent of length ${< n+1}$. 
Our induction hypothesis also implies that 
$[X,K]$ acts trivially by conjugation on $[X \wedge \bS^m,K]$ and, 
furthermore, $\pi_m(K) = [\bS^m,K]_*$ acts trivially on 
$[X \wedge \bS^m,K] \cong [X, C_*(\bS^m,K)]_*$ by 
Lemma~\ref{lem-trivact}. 
Therefore the extension (\ref{eq:ex}) is central, so that 
$[X \times \bS^m,K]$ is nilpotent of length $< n+2$. 
\end{prf}

If $K$ is arcwise connected, then 
\[ \pi_0(C_*(X,K)) = [X,K]_* = [X,K] = \pi_0(C(X,K)),\] 
and this group also acts naturally on the higher homotopy groups 
\[ \pi_m(C(X,K)) \cong \pi_m(C_*(X,K)) \rtimes \pi_m(K).\] 
To analyze this action, we need the following concept: 

\begin{defn} (Samelson products) 
Let $(X,x_0)$ and $(Y,y_0)$ be pointed compact spaces and $K$ 
be a topological group. For $f \in C_*(X,K)$ and $g \in C_*(Y,K)$, the continuous map 
\[ [f,g] \in C(X \times Y,K), \qquad [f,g](x,y) := f(x)g(y) f(x)^{-1} g(y)^{-1} \] 
maps $X \vee Y$ to $\be$, hence defines a continuous map 
\[ [f,g] \in C_*(X \wedge Y,K).\] 
We thus obtain a map 
\[ C^{X,Y} \: [X,K]_* \times [Y,K]_* \to [X \wedge Y,K]_*, \quad 
[[f],[g]] \mapsto [[f,g]] \] 
called the 
\index{Samelson product} 
{\it Samelson product}  (\cite{Whi78}, \cite{Sa53}). 
\end{defn}

\begin{rem} \mlabel{rem:13.8.4} The Samelson product provides in particular a map 
\[ C^{\bS^k,Y} \: [\bS^k,K]_* \times [Y,K]_* 
= \pi_k(K) \times [Y,K]_* \to [\bS^k \wedge Y,K]_* 
\cong\pi_k(C_*(Y,K)),\] 
and, for $Y = \bS^m$ it specializes to a map 
\begin{equation}
  \label{eq:triang}
 C^{\bS^k,\bS^m} \: \pi_k(K) \times \pi_m(K) \to 
[\bS^k \wedge \bS^m,K]_* = [\bS^{k+m},K]_* = \pi_{k+m}(K).
\end{equation}
\end{rem}

\begin{rem} For $X = \bS^p \times \bS^q$, if follows from the preceding theorem 
that the group $[X,K]$ is nilpotent of length $\leq 2$. More concretely, 
(\ref{eq:ex}) specializes to a central extension 
\begin{eqnarray}
  \label{eq:ex2}
\1 \to \pi_{p+q}(K) \cong [\bS^p\wedge \bS^q, K] \into [X,K]  
\sssmapright{q} \pi_p(K) \times \pi_q(K)\to \1. 
\end{eqnarray}
The commutator map of $[X,K]$ factors through a biadditive map 
$$ C \: \pi_p(K) \times \pi_q(K) \to \pi_{p+q}(K), \quad 
([\alpha],[\beta]) \mapsto [p_1^*(\alpha)p_2^*(\beta)p_1^*(\alpha^{-1})
p_2^*(\beta^{-1})] $$
(Exercise~\ref{exer:2-step}).
Writing the right hand side explicitly, we see that $C([\alpha], [\beta])$ coincides with the 
Samelson product for $X = \bS^p$ and $Y = \bS^q$.

It is not particularly easy to find non-trivial Samelson products. 
For the unitary groups $\U_n(\C)$, the situation has been studied in detail 
by Bott in \cite{Bo60}. In this case 
$$ \pi_{2n}(\U_n(\C)) \cong \Z/n!\Z, $$
and Bott shows that, for $r + s = n-1$, we have for suitable generators 
$\alpha \in \pi_{2r+1}(\U_n(\C)) \cong \Z$ and 
$\beta \in \pi_{2s+1}(\U_n(\C)) \cong \Z$ the relation 
$$ C(\alpha,\beta) = r!s! \gamma $$
for some generator $\gamma$ of $\pi_{2n}(\U_n(\C))$. 
\end{rem} 

\begin{ex}
For the direct limit group $\U_\infty(\C) = \indlim \U_n(\C)$, all 
even degree homotopy groups vanish and all odd degree homotopy 
groups are isomorphic to 
$\Z$ (Theorem~\ref{thm2.12}). 
In particular, all Samelson products for this group are trivial 
(the sum of two odd numbers is even). 
For the other classical groups $\U_\infty(\bH)$ and $\OO_\infty(\R)$,  
infinite homotopy groups only arise in degrees congruent $3$ modulo $4$ 
(Theorem~\ref{thm2.12}). Therefore 
all Samelson products for these groups factor through finite groups. 
\end{ex}

\subsection{Reduction by suspension} 

From now on, $K$ will be a topological group and $X$ a compact
space. We further assume that $X$ is arcwise connected. 
In $X$ we fix a base point $x_0$
and in any group we consider the unit element $\be$ as the base point. We
write $C_*(X,K) \subeq C(X,K)$ for the subgroup of base
point preserving maps and observe that 
\begin{eqnarray} \label{eq:10.6.1} 
 C(X,K) \cong C_*(X,K) \rtimes K 
\end{eqnarray} 
as topological groups, where we identify $K$ with the subgroup of
constant maps. This leads to 
\begin{eqnarray} \label{eq:10.6.2} 
 \pi_k(C(X,K)) \cong \pi_k(C_*(X,K)) \times \pi_k(K)\quad \mbox{ for } \quad 
k \in \N 
\end{eqnarray} 
and 
\begin{eqnarray} \label{eq:10.6.2b} 
 \pi_0(C(X,K)) = [X,K] \cong [X,K]_* \rtimes \pi_0(K).
\end{eqnarray} 

\begin{ex} [Spheres]  \mlabel{ex:spheres} 
If $M = \bS^m$ is an $m$-dimensional sphere, then we obtain for every $k$ 
\begin{eqnarray} \label{10.6.6} 
 \pi_k(C_*(\bS^m,K)) \cong [\bS^k \wedge \bS^m,K]_* \cong [\bS^{k+m},K]_* 
\cong \pi_{k+m}(K) 
\end{eqnarray}
and therefore with (\ref{eq:10.6.2}) 
\[  \pi_k(C(\bS^m,K)) \cong \pi_{k+m}(K) \times \pi_k(K) \quad \mbox{ for } \quad  
k > 0 \] 
and 
\[  \pi_0(C(\bS^m,K)) = [\bS^m,K] \cong \pi_m(K) \rtimes \pi_0(K).\] 
\end{ex} 

\begin{ex} [Tori] \mlabel{ex:tori}
Let $M = \T^m$ be the $m$-dimensional torus. 
Then 
\begin{equation}
  \label{eq:m-tori}
C(\T^m,K) 
\cong C(\T, C(\T^{m-1},K)) 
\cong C_*\big(\T,C(\T^{m-1},K)\big) \rtimes C(\T^{m-1},K) 
\end{equation}
implies that 
\[ \pi_k(C(\T^m,K)) 
\cong \pi_{k+1}(C(\T^{m-1},K)) \oplus \pi_k(C(\T^{m-1},K)) 
\quad \mbox{ for } \quad k \geq 1.\] 
By induction, we obtain 
\begin{eqnarray} \label{10.6.7} 
\pi_k(C(\T^m,K)) \cong \sum_{j = 0}^m \pi_{k+j}(K)^{m \choose j} 
\quad \mbox{ for } \quad k \geq 1. 
\end{eqnarray}
For $m = 2$, we get in particular 
\[  \pi_k(C(\T^2,K)) \cong \pi_k(K) \oplus \pi_{k+1}(K)^2 \oplus
\pi_{k+2}(K) \quad \mbox{ for } \quad k \geq 1,\] 
which also follows from the calculations for compact pointed surfaces below 
(cf.\ Theorem~\ref{thm-orient}). 
For $k =0$, we obtain from \eqref{eq:m-tori} that 
\[ [\T^m, K] = \pi_0(C(\T^m,K)) \cong \pi_1(C(\T^{m-1},K)) \rtimes \pi_0(C(\T^{m-1},K)),\] 
so that the determination of the group structure on $\pi_0(C(\T^m,K))$ also involves 
the action of the groups $\pi_0(C(\T^d,K))$ on the groups 
$\pi_k(C(\T^d,K))$ for $d < m$. Since the action of 
$[\T^{m-1},K] $ on $\pi_1(C_*(\T^{m-1},K))$ is trivial by 
Theorem~\ref{thm:whitehead}, the action on 
\[ \pi_1(C(\T^{m-1},K)) \cong 
 \pi_1(C_*(\T^{m-1},K)) \rtimes \pi_1(K) \] 
takes the form 
\[ [\sigma].([\alpha], [\beta]) = ([\alpha] + C([\sigma],[\beta]), [\beta])\] 
because $[\T^{m-1},K] = [\T^{m-1},K]_*$ acts also trivially modulo the 
normal subgroup $\pi_1(C_*(\T^{m-1},K))$. It now follows directly from the 
construction 
that 
\[ C \: [\bS^1,K]_*  \times [\T^{m-1},K]_* \to 
\pi_1(C_*(\T^{m-1},K))  
\cong [\bS^1 \wedge \T^{m-1}, K]_* \] 
is the Whitehead product. The simplest non-trivial case arises for 
$m =2$, where 
\[ [\T^2,K] \cong \big(\pi_2(K) \oplus  \pi_1(K)\big) \rtimes \pi_1(K), \] 
and the group structure involves the Whitehead product 
$\pi_1(K) \times \pi_1(K) \to \pi_2(K)$. 
\end{ex}

\subsection{Reduction by retraction} 

Let $x_0 \in X$ be a base point and $Y \subeq X$ a closed subset
containing $x_0$. We write $X/Y$ for the quotient space obtained by
collapsing $Y$ to a point, and view $Y$ as a base point in $X/Y$. Then
we have a natural group homomorphism 
\[  q^* \: C(X/Y,K) \to C(X,K), f \mapsto f \circ q \quad \hbox{ with  }
\quad q^*\big(C_*(X/Y,K)\big) \subeq C_*(X,K). \] 
We further have a restriction map 
$$ R \: C(X,K) \to C(Y,K) $$
which leads to an exact sequence 
\begin{eqnarray} \label{eq:10.6.8} 
\1 \to C_*(X/Y,K) \to C(X,K) \to C(Y,K). 
\end{eqnarray}

\begin{rem} \mlabel{rem10.6.1} 
(Cofibrations) An inclusion $i \: Y \into X$ is called a
{\it cofibration} \index{cofibration} 
if the pair $(X,Y)$ has the 
\index{homotopy extension property} 
{\it homotopy extension property} 
for all spaces $Z$ (cf. \cite[Sect.\ VII.1]{Bre93}). 
This means that, if $f \: X \to Z$ is
a continuous map and $h_Y \in C(Y,Z)$ 
with $f\res_{Y} \sim h_Y$, then there exists a
continuous extension $h \: X \to Z$ of $h_Y$ 
and an extension of the homotopy 
of $f\res_Y$ to $h_Y$ to a homotopy of $f$ to $h$. 

If $Y \into X$ is a cofibration, then it follows from 
\cite[Cor.\ VII.10.3]{Bre93} that, for any space $Z$, the sequence of maps 
\[  [X/Y,Z]_* \to [X,Z]_* \to [Y,Z]_* \] 
induced on pointed homotopy classes is exact in the sense that the inverse
image of the constant class in $[Y,Z]_*$ consists of those classes
coming from maps $f \:X/Y \to Z$, i.e., maps $f \: X \to Z$ with $f(X)
= \{*\}$. 

This can be seen directly as follows. It is clear that the image of
$[X/Y,Z]_*$ consists of classes whose restriction to $Y$ is constant 
equal to the base point. 
If, conversely, $f\res_Y$ is homotopic to the constant map $*$, then
the homotopy extension property yields a homotopic map 
$\tilde f \: X \to Z$ with $f(Y) = \{*\}$. 

In particular, we obtain for every topological group $K$ exact
sequences of group homomorphisms 
\[ [X/Y,K]_* \to [X,K]_* \to [Y,K]_* \quad \mbox{ and } \quad 
[X/Y,K]_* \to [X,K] \to [Y,K]. \]
\end{rem}

\begin{rem} \mlabel{rem4.2} 
Let $(X,x_0)$ be an arcwise connected pointed topological space. 
Suppose that $Z$ is an arcwise connected pointed space 
which has a universal covering
space $q_Z \: \tilde Z \to Z$, i.e., $Z$ is semilocally
simply connected (\cite{Bre93}). Then we can identify 
$C_*(X,\tilde Z)$ with the subset of $C_*(X,Z)$ consisting
of all those maps which can be lifted to pointed maps $X \to \tilde Z$. 
Recall that a map $f \: X \to Z$ lifts to $\tilde Z$ if and only if 
$\pi_1(f) \: \pi_1(X,x_0) \to \pi_1(Z,z_0)$ vanishes 
(\cite[Satz~6.12]{tD91}). 
As $X$ is
arcwise connected, base point preserving lifts are unique whenever they
exist. Therefore $C_*(X,\tilde Z) = \delta^{-1}(\be)$ for the map 
$$ \delta \: C_*(X,Z) \to \Hom(\pi_1(X), \pi_1(Z)), \quad f \mapsto \pi_1(f). $$
On the level of homotopy classes we thus obtain maps 
$$ [X,\tilde Z]_* \to [X,Z]_* \to \Hom(\pi_1(X),\pi_1(Z)). $$
Since each covering has the homotopy lifting property 
(\cite[Satz 8.3]{tD00}), the map  $[X,\tilde Z]_* \to [X,Z]_*$
is injective. For the case where $Z := K$ is a topological group, we
therefore obtain an exact sequence of group homomorphisms
\begin{eqnarray} \label{eq:10.6.9} 
\1 \to [X,\tilde K]_* \to [X,K]_* \to
\Hom(\pi_1(X),\pi_1(K)).
\end{eqnarray}
\end{rem}

\begin{lem} \mlabel{lem:hom-class-pi1}
Let $E$ be  a vector space, $\Gamma \subeq E$ be a discrete 
subgroup, $K := E/\Gamma$ and $X$ be a paracompact semilocally 
simply connected space. Then the natural homomorphism 
\[ \Phi \: [X,K]_* \to \Hom(\pi_1(X),\pi_1(K)) = \Hom(\pi_1(X),\Gamma) \] 
is an isomorphism. 
\end{lem}

\begin{prf} Since $\tilde K \cong E$ is contractible, $[X,\tilde K]_*$ is trivial, 
so that $\Phi$ is injective by the preceding discussion. It remains to show that 
$\Phi$ is also surjective. 

The paracompactness of $X$ implies that every 
principal $E$-bundle over $X$ has a continuous 
section (an easy argument using partitions of unity), hence is trivial. Now let 
$\gamma \in \Hom(\pi_1(X), \pi_1(K))$ and construct the $E$-principal 
bundle given by 
\[  P_\gamma := (\tilde X \times E)/\pi_1(X), \] 
where $\pi_1(X)$ acts on $X \times E$ from the right by 
\[  (x,g).d := (x.d, g-\gamma(d)). \] 
Then continuous sections $\sigma \: X \to P_\gamma$ correspond via 
$$ \sigma(q_X(x)) = [x, f(x)] $$
to continuous maps $f \: \tilde X \to E$ satisfying 
\[ f(x.d) = f(x) - \gamma(d)\quad \mbox{ for } \quad x \in \tilde X, d \in \pi_1(X). \] 
Then $\oline f(q_X(x)) := [f(x)]$ defines a continuous map $X \to K$ 
with ${\Phi([\oline f]) = \gamma}$. 
\end{prf}

\begin{rem} \mlabel{rem10.6.2}  (Retracts) (a) We want to use the 
exact sequence (\ref{eq:10.6.8}) 
to reduce the problem to calculate the
homotopy groups of $C(X,K)$ to the same problem for the conceivably
simpler groups $C(X/Y,K)$ and $C(Y,K)$. To make this work, we first assume
that there exists a 
\index{retraction} 
{\it retraction} $r \: X \to Y$, i.e., a 
continuous map with $r\res_Y = \id_Y$. Then 
\[  r^* \: C(Y,K) \to C(X,K) \] 
is a homomorphic cross-section of the restriction homomorphism, so that
\[  C(X,K) \cong C_*(X/Y,K) \rtimes C(Y,K) \] 
is a semidirect product of topological groups. 

In general, a retraction from $X$ to $Y$ may  not exist. 
It is much more natural to assume that $r\res_Y \sim
\id_Y$ ($r\res_Y$ is homotopic to the identity). So let us make this
weaker assumption. 

First we consider the map
$$ \Phi \: C(Y,K) \times C_*(X/Y,K) \to C(X,K), \quad 
(\alpha, \beta) \mapsto \beta \cdot (r^*\alpha). $$
Then $R(\Phi(\alpha,\beta)) = r^*\alpha\res_Y \sim \alpha$ shows 
that the induced homomorphism 
$$ [X,K] \to [Y,K] $$
is surjective. Its kernel consists of those maps whose restriction
to $Y$ is homotopic to the constant map $\be$. To obtain a good description of
this subgroup, we have to assume that the
inclusion map $Y \into X$ is a cofibration. Then Remark~\ref{rem10.6.1} leads to
the longer exact sequence 
$$ [X/Y,K]_* \to [X,K] \to [Y,K] \to \be, $$
where the homomorphism 
$$ [\sigma] \: [Y,K] \to [X,K], \quad [\alpha] \mapsto [\alpha \circ
r] $$
splits the homomorphism onto $[Y,K]$. 

 (b) For a locally contractible arcwise connected space $X$, one can use
the isomorphisms 
\begin{equation}
  \label{eq:coloops}
\Phi \: [X,\T]_* \to \Hom(\pi_1(X),\Z), \quad [f] \mapsto \pi_1(f) 
\end{equation}
from Lemma~\ref{lem:hom-class-pi1} 
to represent homomorphisms $\pi_1(X) \to \Z$ as $\pi_1(f)$ 
for a continuous map $f \: X \to \T$. 

Suppose that we have a subset $Y \cong \bS^1 \subeq X$ such that the
inclusion is a cofibration, and a continuous map $r \: X \to \bS^1
\cong Y$ with $r\res_Y \sim \id_Y$. This means that there exists a
homomorphism $\pi_1(X) \to \Z$ mapping the homotopy class of the 
inclusion map $\bS^1 \to Y
\subeq X$ to $1$. Then the considerations in (a) apply to this situation. 
\end{rem}

\subsection{Closed orientable surfaces} 

In this subsection, $\Sigma$ denotes an orientable compact
surface of {\it genus $g$} \index{genus of closed surface} 
and $K$ is an arbitrary topological 
group. 

\begin{rem} \mlabel{rem4.1} 
We recall that $\Sigma$ can be described as a CW complex by starting
with a bouquet  
\[ A_{2g} \cong \underbrace{\bS^1 \vee \bS^1 \vee \ldots \vee \bS^1}_{2g} \] 
of $2g$-circles (cf.\ \eqref{eq:vee-product}). We write 
$a_1, b_1, \ldots, a_g, b_g \: \bS^1 \to A_{2g}$ for the corresponding generators of the
fundamental group of $A_{2g}$ which is a free group on $2g$
generators. Then we consider the continuous 
map $\gamma \: \bS^1 \to A_{2g}$  corresponding to the commutator product 
$$ [a_1, b_1] \cdots [a_g, b_g] \in \pi_1(A_{2g}), $$
where $[x,y] := xyx^{-1}y^{-1}$ denotes a commutator. Now 
$\Sigma$ is homeomorphic to the space obtained by identifying  
the points in $\partial \bB^2 \cong \bS^1$ with their
images in $A_{2g}$ under $\gamma$, i.e., 
$$ \Sigma \cong A_{2g} \cup_\gamma \bB^2. $$ 
In this sense, we can identify
$A_{2g}$ with a subset of $\Sigma$. It is instructive to 
view $\bB^2$ as a regular $4g$-gon, 
where we identify certain points on the boundary edges such that, in
counterclockwise order, the sequence of edges corresponds to the loop 
\[  a_1 b_1 a_1^{-1} b_1^{-1} a_2 \cdots  a_g^{-1} b_g^{-1}. \] 
Now $A_{2g}$ corresponds to the boundary polygon modulo these identifications. 

This procedure shows that a continuous
map $f \: A_{2g} \to Z$ into a topological space $Z$ 
extends to a map $\Sigma \to Z$ if and only if
the corresponding map $\partial \bB^2 \to Z$ extends to the
interior of $\bB^2$, which in turn means that it is a contractible curve. 
Finally, this can be expressed by the condition that 
$$ \pi_1(f) \: \pi_1(A_{2g}) \cong \Z^{*2g} 
= \underbrace{\Z * \Z * \cdots * \Z}_{2g} \to \pi_1(Z)
$$
annihilates the commutator 
$a_1 b_1 a_1^{-1} b_1^{-1} a_2 \cdots  a_g^{-1} b_g^{-1},$
hence factors to a homomorphism 
$\pi_1(\Sigma) \to \pi_1(Z).$

Conversely, if such a homomorphism is given, then we can lift it to a
homomorphism $\pi_1(A_{2g}) \to \pi_1(Z)$ which can be represented by
a continuous map $A_{2g} \to Z$. As we have seen above, this map
extends to $\Sigma$, showing that the map 
\begin{eqnarray} \label{eq:10.6.10} 
 C_*(\Sigma,Z) \to \Hom(\pi_1(\Sigma), \pi_1(Z)), \quad f \mapsto \pi_1(f) 
\end{eqnarray}
is surjective for any pointed space $Z$. 
\end{rem}

\begin{thm} \mlabel{thm-orient} 
For each orientable surface $\Sigma$ of genus
$g$ and each topological group $K$,  we have a homeomorphism 
\begin{equation}
  \label{eq:homeo-sigma}
 C_*(\Sigma,K) \cong C_*(\bS^2,K) \times C_*(\bS^1,K)^{2g}
\end{equation}
and 
\begin{equation}
  \label{eq:13.26}
 \pi_k(C_*(\Sigma,K)) \cong \pi_{k+2}(K) \times \pi_{k+1}(K)^{2g}
\quad \hbox{ for all} \quad k \in \N. 
\end{equation}
\end{thm}

\begin{prf} Let $(\gamma_1, \ldots, \gamma_{2g})$ be the natural generators of
$\pi_1(\Sigma)$ coming from the maps 
$\bS^1 \to A_{2g} \into \Sigma$ given by 
$a_1, b_1, \ldots, a_g, b_g$. 
From (\ref{eq:10.6.10}), we obtain for $Z = \bS^1\subeq \C$ 
with $\pi_1(Z) \cong \Z$ and 
base point $1$,  pointed continuous maps 
$\chi_1, \ldots, \chi_{2g} \: \Sigma \to \bS^1$ with 
$$ \pi_1(\chi_j)(\gamma_i) \cong [\chi_j \circ \gamma_i] = \delta_{ij}. $$
We even get maps with 
\begin{equation}
  \label{eq:split-rel}
\chi_j\circ \gamma_i = 
\left\{   \begin{array}{cl}
1 & \mbox{for $i \not=j$} \cr
\id_{\bS^1} & \mbox{for $i = j$} 
 \end{array}
\right.  
\end{equation}
if we start with the continuous maps $\chi_j^0 \: A_{2g} \to \bS^1$ with the
required property and observe that all these maps extend continuously
to $\Sigma$ because $\pi_1(\chi_j^0) \: \pi_1(A_{2g}) \to \pi_1(\bS^1)
\cong \Z$ annihilates all commutators since $\Z$ is abelian. 

Now we obtain for each topological group $K$ a nice splitting of the restriction 
homomorphism 
\[  R \: C_*(\Sigma, K) \to C_*(A_{2g}, K) \cong C_*(\bS^1, K)^{2g} \] 
by the extension homomorphism 
$$ E \: C_*(\bS^1, K)^{2g} \to C_*(\Sigma, K), \quad 
(\alpha_1, \ldots, \alpha_{2g}) \mapsto (\alpha_1 \circ \chi_1)\cdots
(\alpha_{2g} \circ \chi_{2g}). $$
Then $RE = \id$ follows directly from \eqref{eq:split-rel}. We conclude that 
$$ C_*(\Sigma,K) \to  \ker(R) \times C_*(\bS^1,K)^{2g}, \quad 
f \mapsto (f E(R(f))^{-1}, R(f)) $$
is a homeomorphism whose inverse is given by 
$(\alpha, \beta) \mapsto \alpha E(\beta)$. Next we observe that 
$$ \ker R \cong C_*(\Sigma/A_{2g},K) \cong C_*(\bS^2,K), $$
so that we obtain the homeomorphism~\eqref{eq:homeo-sigma} which in 
turn leads to \eqref{eq:13.26}. 
\end{prf} 

\begin{rem} \mlabel{rem4.7} 
Suppose that $\Sigma$ is an orientable surface of genus ${g \geq 1}$. 
Then the universal covering
space $\tilde \Sigma$ of $\Sigma$ is contractible, showing that the
only non-trivial homotopy group of $\Sigma$ is $\pi_1(\Sigma)$. 
Theorem~\ref{thm-orient} and Remark~\ref{rem4.1} yield a homomorphism 
\[  [\Sigma, K]_* \sssmapright{\alpha}
 \Hom(\pi_1(\Sigma), \pi_1(K)) \cong \pi_1(K)^{2g}, \] 
where $\alpha$ 
has a kernel isomorphic to $\pi_2(K) = [\bS^2,K]_*$, hence is not injective in 
general. This 
means that the homotopy classes of maps $\Sigma \to K$ are not 
classified by the sequence of homomorphisms $\pi_k(\Sigma) \to
\pi_k(K)$, $k \in \N_0$. If $K$ is simply connected, then we simply
have $[\Sigma,K]_* \cong \pi_2(K)$ by \eqref{eq:13.26}. 
\end{rem}

For the study of central group extensions, it is crucial to have a good
description of the first two homotopy groups 
(cf.~Subsection~\ref{subsec:per-grp}). These are given for 
$C_*(\Sigma,K)$ by 
\[  \pi_1(C_*(\Sigma,K)) 
\cong \pi_3(K) \times \pi_2(K)^{2g} \quad \mbox{ and } \quad 
\pi_2(C_*(\Sigma,K)) \cong \pi_4(K) \times \pi_3(K)^{2g}.\] 
If $K$ is a finite-dimensional Lie group, then $\pi_2(K)$ vanishes and
$\pi_4(K)$ is a torsion group (Remark~\ref{rem1.5}). 

\subsection{Orientable surfaces with boundary}
\mlabel{subsec:13.8.4}

After the discussion of mapping groups on closed surfaces, we now turn 
to orientable surfaces with boundary. 

Let $X$ be a topological space, $K$ an arcwise connected topological group and 
$Y \subeq X$ a closed subset. Then we write 
$$ C(X,Y,K) := \{ f \in C(X,K) \: f(Y) = \{\be\}\}. $$

\begin{prop}  \mlabel{prop4.4} 
Let $\Sigma$ be an orientable surface of genus $g$, $\eset\not=P 
\subeq \Sigma$ be a finite subset,  and 
$$ C(\Sigma,P,K) := \{ f \in C(\Sigma,K) \: f(P) = \{\be\}\}. $$
Then 
\begin{equation}
  \label{eq:prop4.4}
 \pi_k(C(\Sigma,P,K)) \cong \pi_{k+1}(K)^{2g+|P|-1} \times \pi_{k+2}(K) 
\quad \mbox{ for } \quad k\in \N. 
\end{equation}
\end{prop}

\begin{prf} As above, we describe $\Sigma$ by a polygon whose boundary is, after
gluing, homeomorphic to $A_{2g}$ (Remark~\ref{rem4.1}). 
We then may assume that the base point
is contained in $P$ and that all other points in $P$ are lying on 
one loop $\alpha_1 \subeq A_{2g}$. Then the homeomorphism 
$$ C_*(\Sigma,K) \cong C_*(\bS^1,K)^{2g} \times C_*(\bS^2,K) $$
(Theorem~\ref{thm-orient}) is adapted to the situation in the sense that 
$C_*(\bS^2,K) \subeq C(\Sigma,P,K)$, so that we have homeomorphisms 
\begin{eqnarray*}
C(\Sigma,P,K) 
&&\cong C(\bS^1,P,K) \times C_*(\bS^1,K)^{2g-1} \times C_*(\bS^2,K) \\ 
&&\cong C_*(\bS^1,K)^{P} \times C_*(\bS^1,K)^{2g-1} \times C_*(\bS^2,K) \\ 
&&\cong C_*(\bS^1,K)^{2g+|P|-1} \times C_*(\bS^2,K).
\end{eqnarray*}
This implies~\eqref{eq:prop4.4} (cf.~Remark~\ref{rem:13.8.4}). 
\end{prf} 

\begin{cor} \mlabel{cor4.5} 
Let $\Sigma$ be a compact oriented surface with 
$p$ boundary components and such that collapsing all boundary
components to points leads to a surface of genus $g$. Then 
\[  \pi_k\big(C(\Sigma,\partial \Sigma,K)\big) 
\cong \pi_{k+1}(K)^{2g+|P|-1} \times \pi_{k+2}(K) \quad \mbox{ for } \quad  k \in \N. \]
\end{cor}

\begin{prf} Let $q \: \Sigma \to M$ be the quotient map collapsing the
boundary components of $\Sigma$ to points and put $P := q(\partial
\Sigma)$. Then $M$ is an orientable surface of genus $g$ and 
$C(\Sigma, \partial \Sigma,K) \cong C(M,P,K),$
so that the assertion follows from Proposition~\ref{prop4.4}. 
\end{prf}

\subsection{Closed non-orientable surfaces} 

The result from the preceding subsection can even be extended to 
mapping groups on non-orientable surfaces. 
In this subsection, $\Sigma$ denotes a non-orientable compact
surface of genus $g$ and $K$ is an arbitrary topological 
group.

\begin{rem} \mlabel{rem4.8} 
We recall that $\Sigma$ can be described as a CW-complex by starting
with a bouquet 
$$A_g \cong \underbrace{\bS^1 \vee \bS^1 \vee \ldots \vee \bS^1}_{g} $$ 
of $g$-circles. We write 
$a_1, \ldots, a_g$ for the corresponding generators of the
fundamental group of $A_g$ which is a free group on $g$
generators. Then we consider the continuous 
map $\gamma \: \bS^1 \to A_g$  corresponding to the element 
\[  a_1^2 \cdots a_g^2 \in \pi_1(A_g). \] 
Now $\Sigma$ is homeomorphic to the space obtained by identifying  
the points in $\partial \bB^2 \cong \bS^1$ with their
images in $A_g$ under $\gamma$, i.e., 
\[  \Sigma \cong A_g \cup_\gamma \bB^2. \] 
In this sense, we can identify
$A_g$ with a subset of $\Sigma$. The most instructive picture is to
view $\bB^2$ as a $2g$-gon, 
where we identify certain points on the boundary edges such that in
counterclockwise order the sequence of edges corresponds to the loop 
$a_1^2 \cdots a_g^2.$
Now $A_g$ corresponds to the polygon modulo these identifications on
the boundary. 

With similar arguments as in Remark~\ref{rem4.1}, we see that the map 
\begin{eqnarray} \label{eq:10.6.12} 
 C_*(\Sigma,Z) \to \Hom(\pi_1(\Sigma), \pi_1(Z)) 
\end{eqnarray}
is surjective for any pointed space $Z$. 
\end{rem}

\begin{thm} \mlabel{thm4.9} 
For every topological group $K$, we have a homeomorphism 
\begin{equation}
  \label{eq:homeo-4.9}
 C_*(\Sigma,K) \cong C_*(P_2(\R),K) \times C_*(\bS^1,K)^{g-1} 
\end{equation}
and accordingly 
\[  \pi_k(C_*(\Sigma,K)) \cong \pi_k(C_*(P_2(\R),K)) \times \pi_{k+1}(K)^{g-1}
\quad \hbox{ for all} \quad k \geq 1.\] 
\end{thm}

\begin{prf} Let $(\gamma_1, \ldots, \gamma_{g})$ be the natural generators of
$\pi_1(\Sigma)$ coming from the maps 
$\bS^1 \to A_g \into \Sigma$ given by 
$a_1, \ldots, a_g$. 
From (\ref{eq:10.6.12}), we obtain for $Z = \bS^1$ with $\pi_1(Z) \cong \Z$ 
pointed continuous maps 
$\chi_1, \ldots, \chi_{g-1} \: \Sigma \to \bS^1$ with 
$$ \pi_1(\chi_j)(\gamma_i) \cong [\chi_j \circ \gamma_i] = \delta_{ij}
\quad \hbox{ for } \quad i,j=1,\ldots, g-1$$
and 
$$ \pi_1(\chi_j)(\gamma_g) = -1 \quad \hbox{ for } \quad j=1,\ldots, g-1.$$
The maps $\chi_j$ can be chosen in such a way that 
\begin{equation}
  \label{eq:triang2}
\chi_j\circ \gamma_i = 
\left\{   \begin{array}{cl}
1 & \mbox{for $i \not=j, i,j < g$} \\ 
\id_{\bS^1} & \mbox{for $i = j < g$} \\ 
z \mapsto \oline z& \mbox{for $j < i = g$, $z \in \bT \cong \bS^1$}  
  \end{array}
\right. 
\end{equation}
if we start with the continuous maps $\chi_j^0 \: A_g \to \bS^1$ with the
required property and observe that all these maps extend continuously
to $\Sigma$ because the homomorphism 
$\pi_1(\chi_j^0) \: \pi_1(A_g) \to \pi_1(\bS^1)
\cong \Z$ annihilates $a_1^2 \cdots a_g^2$. 

Let $A_{g_1} \subeq A_g$ denote the bouquet of the $g-1$ circles
corresponding to $a_1, \ldots, a_{g-1}$, resp., to the maps 
$\gamma_1,\ldots, \gamma_{g-1}$. 
Then we obtain a topological splitting 
of the restriction homomorphism 
\[  R \: C_*(\Sigma, K) \to C_*(A_{g-1}, K) \cong C_*(\bS^1, K)^{g-1} \] 
by the extension map 
$$ E \: C_*(\bS^1, K)^{g-1} \to C_*(\Sigma, K), \quad 
(\alpha_1, \ldots, \alpha_{g-1}) \mapsto (\alpha_1 \circ \chi_1)\cdots
(\alpha_{g-1} \circ \chi_{g-1}). $$
Then $RE = \id$ follows from \eqref{eq:triang2}. This leads to a homeomorphism 
\[  C_*(\Sigma,K) \to  \ker(R) \times C_*(\bS^1,K)^{g-1}, \quad 
f \mapsto (f E(R(f))^{-1}, R(f)) \] 
whose inverse is given by 
$(\alpha, \beta) \mapsto \alpha E(\beta)$. Next we observe that 
$$ \ker R \cong C_*(\Sigma/A_{g-1},K) \cong C_*(P_2(\R),K), $$
so that we obtain the homeomorphism~\eqref{eq:homeo-4.9}. 
This implies the theorem. 
\end{prf} 

\subsection{Non-orientable surfaces with boundary}

Finally, we also turn to 
mapping groups on non-orientable surfaces with boundary.

\begin{prop} \mlabel{prop4.10} 
Let $\Sigma$ be a non-orientable surface of genus
$g > 1$ and $\eset\not=P \subeq \Sigma$ be a finite subset. Then 
$$ \pi_k(C(\Sigma,P,K)) 
\cong \pi_{k+1}(K)^{g+|P|-2} \times \pi_k(C_*(P_2(\R),K)), \quad k
\in \N. $$
\end{prop}

\begin{prf} As above, we describe $\Sigma$ by a polygon whose boundary is, after
gluing, homeomorphic to $A_{g}$. We then may assume that the base point
is contained in $P$ and that all other points are lying on 
one loop $\gamma_1 \subeq A_{g}$ (see the proof of Theorem \ref{thm4.9}). 
Then the homeomorphism 
\[  C_*(\Sigma,K) \cong C_*(P_2(\R), K) \times C_*(\bS^1,K)^{g-1} \] 
(Theorem \ref{thm4.9}) is adapted to the situation in the sense that 
$C_*(P_2(\R),K) \subeq C(\Sigma,P,K)$, so that we have homeomorphisms 
\begin{eqnarray*}
C(\Sigma,P,K) 
&&\cong C(\bS^1,P,K) \times C_*(\bS^1,K)^{g-2} \times C_*(P_2(\R),K) \cr
&&\cong C_*(\bS^1,K)^{P} \times C_*(\bS^1,K)^{g-2} \times C_*(P_2(\R),K) \cr
&&\cong C_*(\bS^1,K)^{g+|P|-2} \times C_*(P_2(\R),K).
\end{eqnarray*}
This implies the assertion.
\end{prf}

\begin{cor} \mlabel{cor4.11} 
Let $\Sigma$ be a compact non-orientable surface with 
$p$ boundary components such that collapsing all boundary
components to points leads to a surface of genus $g > 1$. Then 
\[  \pi_k(C(\Sigma,\partial \Sigma,K)) 
\cong \pi_{k+1}(K)^{g+|P|-2} \times \pi_k(C_*(P_2(\R),K)) \quad \mbox{ for } \quad k
\in \N. \] 
\end{cor}

\begin{prf} As in the proof of Corollary~\ref{cor4.5}, 
this can be derived from Proposition~\ref{prop4.10}.
\end{prf}

It remains to determine the homotopy groups for the case $g = 1$,
where $\Sigma \cong P_2(\R)$ is the real projective plane. 
In the following, we will show in
particular that, for each Fr\'echet--Lie group $K$, all homotopy groups 
$\pi_k(C_*(P_2(\R),K))$ are $2$-groups. More precisely, they are
extensions of an elementary abelian $2$-group by another one, which
implies in particular that the order of an element does not exceed $4$.

\begin{lem} \mlabel{lem-mn1} 
Let $(M_i,x_i)$, $i =1,2$, be two pointed compact spaces and 
$\alpha_1, \alpha_2 \: M_1 \to M_2$ be two pointed homotopic 
maps. Then the Lie group homomorphisms 
$$ \alpha_{j,K} \: C_*(M_2,K) \to C_*(M_1,K),
\quad f \mapsto f \circ \alpha_j $$
satisfy $\pi_m(\alpha_{1,K}) = \pi_m(\alpha_{2,K})$ for
each $m \in \N_0$. 
\end{lem}

\begin{prf} Let $F \: [1,2] \times M_1 \to M_2$ be a base point preserving 
homotopy with $F_1 = \alpha_1$ and $F_2 = \alpha_2$. Then the map 
$$ \Phi \: [1,2] \times C_*(M_2,K) \to C_*(M_1,K), \quad 
\Phi(t,f)(s) := f(F(t,s)) $$
is continuous because the map 
$$ \hat\Phi \: [1,2] \times C_*(M_2,K) \times M_1\to K, \quad 
\hat\Phi(t,f,s) := f(F(t,s)) = \ev(f, F(t,s)) $$
is continuous, which in turn follows from the continuity of the
evaluation map 
$$ \ev \: C_*(M_2,K) \times M_2 \to K. $$
We conclude that the two maps $\Phi_1, \Phi_2 \: C_*(M_2,K) \to C_*(M_1,K)$
are homotopic, hence induce the same map 
$\pi_m(C_*(M_2,K)) \to \pi_m(C_*(M_1,K))$ for each $m \in \N_0$. 
\end{prf}

\begin{lem} \mlabel{lem:map-degree} 
For $h \in C(\bS^1, \bS^1)$ and $m \in \N_0$, the homomorphism 
\[ \pi_m(h_K) \: \pi_m(C_*(\bS^1,K)) \to \pi_m(C_*(\bS^1,K))\] 
is given by 
$$ \pi_m(h_K)([\sigma]) = \deg(h) \cdot [\sigma], $$
where $\deg(h) = [h] \in \pi_1(\bS^1) \cong \Z$ is the mapping degree of $h$. 
\end{lem}

\begin{prf} We realize $\bS^1$ as $\R/\Z$, so that continuous functions $\bS^1 \to
K$ correspond to continuous $1$-periodic functions $\R \to K$. 
In view of Lemma~\ref{lem-mn1}, $\pi_m(h_K)$ only depends on the
homotopy class of $h$, so that we may assume that 
$h(z) = n z$ for  $n = \deg(h) \in \Z$. 

This leads to the maps 
$$ \phi_n \: C_*(\bS^1,K)  \to C_*(\bS^1,K), \quad 
\phi_n(f)(t) = f(nt). $$
We claim that $\phi_n$ is homotopy equivalent to the map $\psi_n(f) := f^n$. 

We assume that $n > 0$. The case $n = 0$ is trivial and the case $n < 0$ 
is treated similarly. For each 
interval $[\frac{i}{n}, \frac{i+1}{n}]$, $i =0,\ldots, n-1$, we define a
continuous map 
$$ \alpha_i \: C_*(\bS^1,K)  \to C_*(\bS^1,K), 
\quad \alpha_i(f)(t) := f(\tilde\alpha_i(t)), \qquad 0 \leq t \leq 1, $$
where 
$$ \tilde\alpha_i \: [0,1] \to [0,1], \quad t \mapsto 
\left\{ \begin{array}{cl} 
0 & \mbox{for $t \leq \frac{i}{n}$} \cr 
n t-i & \mbox{for $\frac{i}{n} \leq t \leq \frac{i+1}{n}$} \cr 
1 & \mbox{for $\frac{i+1}{n} \leq t \leq 1$}. 
\end{array} \right. $$
This means that the functions $\alpha_i(f)$ are ``supported'' by the
$\Z$-translates of the interval $[\frac{i}{n}, \frac{i + 1}{n}]$.
Then each map $\tilde\alpha_i$ is homotopic to the identity of
$[0,1]$ with fixed endpoints, and the same carries over to $\alpha_i$. 
Now 
$$ \phi_n(f) = \alpha_1(f) \cdot \alpha_2(f) \cdots\alpha_n(f)$$
is a pointwise product because the supports of the factors are
disjoint. As each map $\alpha_i$ is homotopic to $\id_{C_*(\bS^1,K)}$, the
map $\phi_n$ is homotopic to the $n$th power map. 

The $n$th power map on $C_*(\bS^1,K)$ induces the $n$th power map on
the corresponding homotopy groups, where the multiplication is induced
by pointwise multiplication in $K$, and we conclude that 
$$\pi_m(\phi_n) \: \pi_m(C_*(\bS^1,K)) \to \pi_m(C_*(\bS^1,K)) $$
is the $n$th power map in the abelian group $\pi_m(C_*(\bS^1,K))$. 
\end{prf}

If $A$ is an abelian group and $n \in \N$, then we write 
\[  A[n] := \{ a \in A \: a^n = \be\}, \qquad n \in \N\] 
for the 
\index{$n$-torsion subgroup of abelian group $A$, $A[n]$}
{\it $n$-torsion subgroup.}

\begin{thm} \mlabel{thm4.13} 
Let $K$ be a connected Fr\'echet--Lie group with Lie algebra~$\fk$ and 
$\eset \not= P \subeq \Sigma \cong P_2(\R)$ be a finite subset. Then, for
each $k \in \N_0$, there exists a short exact sequence 
$$ \pi_{k+2}(K)/2\pi_{k+2}(K)
\into \pi_k(C(P_2(\R),P,K)) \onto \pi_{k+1}(K)[2] \times \pi_{k+1}(K)^{|P|-1} . $$
\end{thm}

\begin{prf} Let $p := |P|$, and 
$q_\Sigma \: \bB^2 \to P_2(\R)$ be the quotient map for which 
$\gamma := q_\Sigma\res_{\partial \bB^2} \: \bS^1 \to q_\Sigma(\bS^1) = A_1$ is
the two-fold covering map with $\gamma(-x) = \gamma(x)$ for $x \in
\partial \bB^2$. We may w.l.o.g.\ assume that $P \subeq A_1$. 
Let $\tilde P := q_\Sigma^{-1}(P) = \gamma^{-1}(P)$. 

Since $K$ is a Fr\'echet--Lie group, the same holds for all the groups 
$C(X,Y,K)$, where $X$ is a compact space and $Y \subeq X$ a subset. 
According to  Theorem \ref{thm4.12}, 
the restriction map 
$$R_B\: C(\bB^2, K) \to C(\bS^1, K)$$ 
has continuous local sections because the restriction map 
$$\L(R_B) \: C(\bB^2, \fk) \to C(\bS^1, \fk) $$
is surjective. The image of the restriction map consists of all 
contractible loops, hence with the identity component 
$C(\bS^1, K)_0$ of $C(\bS^1, K)$. We therefore have a topological 
principal bundle 
\[  C_*(\bS^2, K) \into C_*(\bB^2, K) \onto C_*(\bS^1, K)_0 \]
which defines a long exact sequence of homotopy groups 
\[   \ldots \to \pi_{k+1}(C_*(\bS^1, K)) 
\smapright{\delta_{k+1}} \pi_k(C_*(\bS^2,K)) 
\to \pi_k(C_*(\bB^2, K)) \to \ldots \]
As the group $C_*(\bB^2,K)$ is contractible (Exercise~\ref{exer:contract}), 
all its homotopy groups vanish, and the connecting maps 
$$  \delta_{k+1} \: \pi_{k+1}(C_*(\bS^1,K)) \to \pi_k(C_*(\bS^2,K)) $$
are isomorphisms. 

Next we consider the restriction map 
$R_P\: C(\Sigma,P,K) \to C(A_1,P,K)$, where $A_1 = q_\Sigma(\bS^1)
\cong~\bS^1$. 
We have $\pi_1(\Sigma) \cong \Z_2$, and a continuous map 
$f \in C(A_1,P,K)$ extends to $\Sigma$ if and only if 
$\pi_1(f)(2) = \pi_1(f)(1)^2 \in \pi_1(K)$ vanishes. Let 
$$ C(\bS^1,P,K)_2 := \{ f \in C(\bS^1,P,K) \: \pi_1(f)(2) = 0\}. $$
Then $\im(R_P) = C(\bS^1,P,K)_2$. In particular, on the Lie
algebra level the map 
$$ \L(R_P) \: C(\Sigma,P,\fk) \to C(A_1,P,\fk) \cong C(\bS^1,P,\fk) $$
is surjective, so that Theorem \ref{thm4.12} leads to the topological 
principal bundle 
$$ C_*(\bS^2,K) \into C(P_2(\R),P,K) \onto C(\bS^1,P,K)_2. $$
The pullback map $\gamma^* \: f \mapsto f\circ \gamma$ 
satisfies $R_P \circ \gamma^* = \gamma^* \circ R_B$,
so that we obtain a commutative diagram 
\[ \begin{matrix} 
C_*(\bS^2,K) & \into & C(P_2(\R),P,K) & \onto & C(\bS^1,P,K)_2 \cr 
\mapdown{\cong} & & \mapdown{q_\Sigma^*} &  & \mapdown{\gamma^*} \cr 
C_*(\bS^2,K) & \into & C_*(\bB^2,K) & \onto & C_*(\bS^1,K)_0. \end{matrix} \]

In view of $C(\bS^1,P,K) \cong C_*(\bS^1,K)^P$, 
the pull-back map
\[ \gamma^* \: C(\bS^1,P,K)_2 \to C_*(\bS^1,K)_0 \]  can be viewed as
a composition map for loops which maps 
$(f_1,\ldots, f_p)$ to a loop with the homotopy class of 
$$ 2 [f_1] + \ldots + 2 [f_p]$$
because $\gamma^* \: C_*(\bS^1,K)_2 \to C_*(\bS^1,K)_0$ 
is homotopic to the squaring map $f \mapsto f^2$ 
by Lemma~\ref{lem:map-degree}. 
 
We also need that $\pi_1(K) \cong \pi_0(C_*(\bS^1,K))$, 
that $C_*(\bS^1,K)_2$ contains the identity component 
of $C_*(\bS^1,K)$, and that 
$$ \pi_0(C_*(\bS^1,K)_2) \cong \{ g \in \pi_1(K) \: g^2 = \be\} =
\pi_1(K)[2]. $$

Now the naturality of the long exact homotopy sequence leads to
the commutative diagrams 
$$ \begin{matrix} 
\ldots & \to & \pi_{k+1}(C(\bS^1,P,K)) & \sssmapright{\delta_{k+1}^P} 
& \pi_k(C_*(\bS^2,K)) & \to &
\pi_k(C(P_2(\R),P,K)) & \to & \ldots \cr 
& & \mapdown{\pi_{k+1}(\gamma^*)} & & \mapdown{\cong} &  &
\mapdown{\pi_k(q_\Sigma^*)} & & \cr 
\ldots & \sssmapright{0} & \pi_{k+1}(C_*(\bS^1,K)) &
\sssmapright{\delta_{k+1}^\bB}
 & \pi_k(C_*(\bS^2,K)) & \to & \pi_k(C_*(\bB^2,K)) = \0& \to & \ldots \end{matrix} $$
The commutativity of the diagram above implies that the connecting map 
$\delta_{k+1}^P$ is equivalent to the  map 
\begin{eqnarray*}
\pi_{k+1}(C_*(\bS^1,K)^P) \cong \pi_{k+2}(K)^{P} 
&\to & \pi_{k}(C_*(\bS^2,K)) \cong \pi_{k+2}(K) \\ 
(a_1,\ldots, a_p) &\mapsto& 2(a_1 + \ldots + a_p), 
\end{eqnarray*}
so that 
$$ \im(\delta_{k+1}^P) = 2 \pi_{k+2}(K)
\quad \hbox{ and } \quad \ker(\delta_{k+1}^P) \cong \pi_{k+2}(K)^{p-1}
\times \pi_{k+2}(K)[2]. $$
From the exactness of the long exact homotopy sequence of $R_P$, we
now obtain short exact sequences 
$$  \pi_{k+2}(K)/2\pi_{k+2}(K) 
\into \pi_k(C(P_2(\R),P,K)) \onto 
\pi_{k+1}(K)^{p-1} \times \pi_{k+1}(K)[2],$$ 
for $k \in \N_0.$ 
This completes the proof. 
\end{prf}

\begin{cor} \mlabel{cor4.14} Let $K$ be a connected Fr\'echet--Lie group. Then 
there exist short exact sequences 
\[  \pi_{k+2}(K)/ 2\pi_{k+2}(K) \into \pi_k(C(P_2(\R),K)) \onto
\pi_{k+1}(K)[2] \quad \mbox{ for } \quad k \in \N. \]
\end{cor}

\begin{rem} \mlabel{rem4.15} A priori we do not know if the group 
$\pi_0(C_*(P_2(\R),K)) = [P_2(\R),K]_*$ is abelian. This group is an
extension of the $2$-group $\pi_1(K)[2]$ by the $2$-group 
$\pi_2(K)/2\pi_2(K)$. 

For an abelian group $A$ we write $A_\Q := A \otimes \Q$. 
For $k \in \N$, the groups $\pi_k(C_*(P_2(\R),K))$ are abelian
$2$-groups, and it follows in particular that 
$$ \pi_k(C_*(P_2(\R),K))_\Q =\0.$$
Therefore we obtain for each non-oriented compact surface of genus $g$
for the rational homotopy groups 
\[  \pi_k(C(\Sigma,K))_\Q \cong\pi_{k+1}(K)^{g-1}_\Q 
\times \pi_k(K)_\Q \quad \mbox{ for } k \in \N \] 
and in particular 
\[  \pi_2(C(\Sigma,K))_\Q \cong\pi_3(K)^{g-1}_\Q. \] 
\end{rem}

\begin{small}
\subsection*{Exercises for Chapter~\ref{ch:top}} 

\begin{exer} \mlabel{exer:10.1.1} 
Show that the quotient space $I^n/\partial I^n$ obtained by collapsing the 
boundary of the cube $I^n$ to one point is homeomorphic to the $n$-dimensional 
sphere~$\bS^n$.\\
Hint: The interior of $I^n$ is homeomorphic to $\R^n$ and 
$I^n/\partial I^n$ is a one-point compactification. 
\end{exer}

\begin{exer} \mlabel{exer:10.1.2} Show that for a space triple $(X,X_0,\{x_0\})$ the 
product on \break $F^n(X,X_0,x_0)$ satisfies: 
\begin{enumerate}
\item[\rm(1)] $f \sim f', g \sim g' \Rarrow f * g \sim f' * g'$. 
\item[\rm(2)] If $x_0$ denotes the constant map $I^n \to X$, then 
$x_0 * f \sim f * x_0 \sim f.$
\item[\rm(3)] (Inverse) Let $\oline f(t_1, \ldots, t_n) := f(t_1, \ldots, t_{n-1}, 1-t_n)$. 
Then $[f] * [\oline f] = [\oline f] * [f] = [x_0]$. 
\item[\rm(4)] (Associativity) $(f * g) * h \sim f * (g * h)$  
for $f,g,h \in F^n(X,X_0,x_0)$, $n \geq 2,$ or $n = 1$ and $X_0 = \{x_0\}$. 
\item[\rm(5)] (Commutativity) $f * g \sim g* f$ 
for $f,g \in F^n(X,X_0,x_0)$, $n \geq 3,$ or $n = 2$ and $X_0 = \{x_0\}$. \\
Hint: Define a second product 
$$ (f \sharp g)(x_1, \ldots, x_n) 
:= \left\{   \begin{array}{cl}
f(x_1, x_2, \ldots, 2x_{n-1}, x_n) & \mbox{for $0 \leq x_{n-1} \leq \frac{1}{2}$} \\ 
g(x_1, x_2, \ldots, 2x_{n-1}-1, x_n) & \mbox{for$\frac{1}{2} < x_{n-1} \leq 1$}. \\ 
  \end{array} 
\right.  $$ 
and then show that 
$$ f * g \sim (f \sharp x_0) * (x_0 \sharp g) \sim f \sharp g 
\sim (x_0 * f) \sharp (g * x_0) \sim g * f. $$
\item[\rm(6)] (Functoriality) For $f,g \in F^n(X,X_0,x_0)$ 
and $\phi \in C((X,X_0,x_0), (Y,Y_0,y_0))$ we have 
$(\phi \circ f) * (\phi \circ g) = \phi \circ (f * g).$ 
\end{enumerate}
\end{exer}

\begin{exer} \mlabel{exer:10.1.3} 
Let $G$ be a topological group and 
$f,f',g,g',h \in F^n(G,\be), {n \geq 0}$. Show that:  
\begin{enumerate}
\item[\rm(1)] $f\sim f', g \sim g' \Rarrow f \cdot g \sim f' \cdot g'$. 
\item[\rm(2)] $f \sim g \Longleftrightarrow fg^{-1} \sim \be$, 
the constant map. 
\item[\rm(3)] (Commutativity) $f \cdot g \sim g \cdot f$ for $n \geq 1$. \\
Hint: Show that 
$f \cdot g \sim (f * \be) \cdot (\be * g) = (\be * g)  \cdot (f * \be) 
\sim g \cdot f.$
\item[\rm(4)] (Consistency) $f \cdot g \sim f * g$ for $n \geq 1$.   
\end{enumerate}
\end{exer} 

\begin{exer} \mlabel{exer:10.1.4} Let $G$ be a topological group and 
$q \: P \to X$ a principal $G$-bundle. Show that the canonical map 
$P/G \to X$ is a homeomorphism. 
\end{exer} 

\begin{exer} \mlabel{exer:10.1.5} Let $G$ be a topological group and 
$H \leq G$ be a closed subgroup. Show that the fiber bundle 
$q \: G \to G/H$ is locally trivial if and only if there exists 
an open subset $U \subeq G/H$ 
on which there exists a continuous section $\sigma \: U \to G$ of the 
quotient map~$q$. 
\end{exer} 

\begin{exer} \mlabel{exer:10.1.5b} Let $G$ be a topological group and 
$H \leq G$ be a closed subgroup. Further, let $K \subeq H$ be an open 
subgroup. We then have a natural map 
\[ q \: G/K \to G/H, \quad g K \mapsto gH.\] 
We write $q_K \: G \to G/K$ and $q_H \: G \to G/H$ 
for the quotient maps. Show that: 
\begin{enumerate}
\item[\rm(a)] $q$ is continuous and open. 
\item[\rm(b)] If $U \subeq G$ is an  open symmetric $\be$-neighborhood for which 
$U^2 K \cap H = K$ (existence follows from $K$ being open in $H$), then, 
for each $h \in H$, the map  
$q^h := q\res_{q_K(Uh)} \: q_K(Uh) \to q_H(U)$ is bijective, 
continuous and open, hence a 
homeomorphism. 
\item[\rm(c)] Fix a set $(h_j)_{j \in J}$ of representatives of the 
left $K$-cosets in $H$. Then 
\[ q_K(UH) = \dot\bigcup_{j \in J} q_K(Uh_j), \] 
and if $K$ is normal in $H$, then these cosets are mutually disjoint. 
\item[\rm(d)] If $K \trile H$ is normal in $G$, then $q$ is a covering 
map with $\Deck(q)\cong H/K$. 
\end{enumerate}
\end{exer} 

\begin{exer} \mlabel{exer:10.1.5c} Let $G$ be a connected locally 
exponential Lie group and $H \subeq G$ 
be a closed connected Lie subgroup with inclusion map $\iota_H \: H \into G$. 
Further, let $q_G \: \tilde G \to G$ be a universal covering group 
with $\L(q_G) = \id_{\L(G)}$, 
and \break $\hat H := \la \exp_{\tilde G} \L(H) \ra$ be the integral 
subgroup corresponding to the Lie subalgebra $\L(H) \subeq \L(G) \cong \L(\tilde G)$. 
Show that: 
\begin{enumerate}
\item[\rm(a)] $\hat H$ is the identity component of the Lie subgroup 
$H^\sharp := q_G^{-1}(H) \subeq \tilde G$. 
\item[\rm(b)] The identification $\pi_1(G) \cong \ker q_G \subeq H^\sharp$ 
leads to an exact sequence 
\[ \pi_1(H) \smapright{\pi_1(\iota_H)} \pi_1(G) \onto \pi_0(H^\sharp).\] 
\item[\rm(c)] $\pi_1(\hat H) \cong \ker(\pi_1(\iota_H)).$ 
\item[\rm(c)] $\hat H$ is $1$-connected if and only if 
$\pi_1(\iota_H)$ is an isomorphism. 
\end{enumerate}
\end{exer}

\begin{exer} \mlabel{exer:10.1.6} Let $G$ be a Banach--Lie group 
and $H \leq G$ a closed Lie subgroup and which $\L(H)$ has a closed 
complement $\fq \subeq \L(G)$. Show that there exists 
an open subset $U \subeq G/H$ 
on which there exists a continuous section $\sigma \: U \to G$ of the 
quotient map $q \: G \to G/H$.\\ 
Hint: Apply the Inverse Function Theorem to the map 
$\fq \times H \to G, (x,h) \mapsto \exp(x)h$. 
\end{exer}

\begin{exer} \mlabel{exer:10.2.1} Show that group 
\[  \SU_2(\C) = \Big\{ 
\begin{pmatrix} a & b \\ -\oline b & \oline a \end{pmatrix} 
\in M_2(\C)\: |a|^2 + |b|^2 = 1\Big\} \] 
is diffeomorphic to $\bS^3$.
\end{exer} 
  
\begin{exer} \mlabel{exer:10.4.1} 
Let $\cH$ be a complex Hilbert space and write 
$Iv = iv$ for the complex structure on $\cH$. 
Then 
\[  \Im \la x, y \ra  = \Re\la x, Iy\ra \quad \mbox{ for all } \quad 
x,y\in \cH. \]
\end{exer} 

\begin{exer} \label{exer:2-step} 
Let $N$ be a $2$-step nilpotent group and $Z := Z(N)$ its center. Then the 
commutator map $(x,y) \mapsto xyx^{-1}y^{-1}$ factors through a biadditive map 
$$ N/Z \times N/Z \to Z, \quad (xZ,yZ) \mapsto xyx^{-1}y^{-1}. $$
\end{exer}

\begin{exer} \mlabel{exer:contract} 
Let $(X,x_0)$ be a compact contractible pointed space and $G$ a 
topological group. Show that the group $C_*(X,G)$ is contractible. 
\end{exer}

\begin{exer} \mlabel{exer:contr-sphere} 
Let $\cH$ be an infinite-dimensional Hilbert space 
and 
\[ \bS := \{ v \in \cH \: \|v\| = 1\} \quad \mbox{ the sphere in } 
\quad \cH.\]
We discuss a few arguments, showing that $\bS$ is contractible. 
Show that: 
\begin{enumerate}
\item[\rm(a)] For $\cH = L^2([0,1])$, the constant function $1$ is 
contained in $\bS$. We consider the map 
\[  F(t,f)(x) := \frac{f_t}{\|f_t\|} \quad \mbox{ with }\quad 
f_t(x) := 
\begin{cases}
f(x) & \text{for } 0 \leq x \leq t \\ 
1& \text{for } t < x \leq 1. 
\end{cases}\] 
Then $F \: [0,1] \times \bS \to \bS$ 
is continuous with $F(1,f) = f$ and $F(0,f) = 1$. 
Note that every separable Hilbert space is unitarily isomorphic 
to  $L^2([0,1])$. 
\item[\rm(b)] For every Hilbert space (infinite-dimensional or not), 
the complement $\bS \setminus \{p\}$ of any point $p\in \bS$ 
is homeomorphic to the hyperplane $p^\bot$, hence contractible. \\ 
Hint: Stereographic projection. 
\item[\rm(c)] Let $T \: \cH \to \cH$ be an isometry which is not 
surjective and for which $-1$ is not an eigenvalue. Then 
\[  F(t,f)(x) := \frac{f_t}{\|f_t\|} \quad \mbox{ with }\quad 
f_t := t f + (1-t) Tf \] 
is a homotopy from $T \: \bS \to \bS$ to $\id_\bS$. 
Since $T(\bS)$ is a proper subset, $T$ is homotopic to a constant map 
by~(b). 
\end{enumerate}
\end{exer}

\end{small}

\section{Notes and comments on Chapter~\ref{ch:top}} 

{\bf Section~\ref{sec10.4}:} The proof of Theorem~\ref{kuiper-thm} is taken from 
\cite{Ne02c}. It is based on a hint in a
footnote in Kuiper's paper. 
For complex Hilbert spaces Kuiper's Theorem can also be derived 
from more general results of Br\"uning and Willgerodt on the
contractibility of unit groups of von Neumann algebras of infinite
type (\cite{BW76}). For the unitary group $\U(\cH)$ of a separable 
Hilbert space, endowed with the strong topology, the contractibility 
is shown in \cite{DiD63}. 

As we have seen in Section~\ref{sec10.4}, many natural groups of operators 
have polar decompositions, and these provide important tools to understand their 
topology. In finite dimensions, the Iwasawa decomposition and its generalizations 
are equally important, but they are much less robust 
with respect to infinite-dimensional generalizations. A 
very concrete result in this direction is \cite[Prop.~1.1]{Bel09} 
which implies that the Iwasawa decomposition fails for $\GL(\cH)$ on the Hilbert 
space $\cH = \ell^2(\N)$, but this paper also contains a series 
of positive results for other contexts. For loop groups consisting 
of maps with absolutely convergent Fourier series, 
V.~Balan and J.~Dorfmeister obtained in \cite{BD01} natural 
versions of the Birkhoff and the Iwasawa decomposition. Here the 
summability of the Fourier series is needed to obtain a splitting 
into positive and negative Fourier modes on the Lie algebra level. \\

\nin {\bf Section~\ref{sec:top-schatten}:} 
From the results in Sections~\ref{sec:top-dirlim} and 
\ref{sec:top-schatten} it follows that the restricted orthogonal 
group $\OO_2(\cH) = \OO(\cH)\cap \cL_2(\cH)$ of an infinite-dimensional 
real Hilbert space $\cH$ satisfies 
\[ \pi_0(\OO_2(\cH)) \cong \pi_1(\OO_2(\cH)) \cong \Z/2\Z.\] 
Hence its identity component has a $1$-connected 
simply connected covering $\Spin_2(\cH)$. For a detailed 
discussion of this group we refer to \cite{PS75} and \cite{Pl94}. 

\nin {\bf Section~\ref{sec:top.4}} is mostly taken from 
 Appendix~A.3 in \cite{Ne02a} for compact manifolds, and from Appendix~A in 
\cite{Ne04c} for non-compact manifolds. 

\begin{probl}
  \mlabel{prob:VI.7} For a linear 
quotient map $q \: E \to Q$ of Fr\'echet spaces, we may use 
\cite{MicE59} to find a continuous linear cross section $\sigma \: Q \to E$, which implies in 
particular that $q$ defines a topologically trivial fiber bundle. For more general 
locally convex spaces, cross sections might not exist, but it would still be 
interesting if quotient maps of locally convex spaces are Serre fibrations, 
i.e., have the homotopy lifting property for cubes (cf.\ \cite{Bre93} 
and Section~\ref{sec:10.1}). 
Here the main point is to show that continuous curves can be lifted. 
If this is the 
case, the long exact homotopy sequence (Theorem~\ref{homseq-princ}) 
would also be available for quotient maps 
of locally exponential Lie groups, which would 
be an important tool to calculate homotopy 
groups of such Lie groups. 
\end{probl}

\chapter{Selected topics} 
\mlabel{ch:sectop} 

In this chapter we discuss a variety of special topics concerning 
infinite-dimensional Lie groups. 

We start in Section~\ref{sec:3.6} by showing that a connected Lie group 
$G$ is nilpotent, resp., solvable if and only if its Lie algebra $\g$ 
has this property. We also show that, for solvable Lie groups, the 
existence of an exponential function and the Mackey completeness implies 
regularity. 

In Section~\ref{sec:polgrp} we study Lie groups defined by a polynomial 
multiplication of a locally convex space $G$ for which the degrees of the 
iterated multiplication maps are bounded. The main result is that 
such a Lie group $G$ has a polynomial exponential map 
$\Exp \: \g \to G$ with a global polynomial inverse and that 
its Lie algebra $\g$ is nilpotent. As a consequence, 
every polynomial Lie group $G$ is polynomially isomorphic to a 
Lie group of the form $(\g,*)$, where $\g$ is a nilpotent Lie algebra,
endowed with the BCH multiplication. 

In Section~\ref{sec:formgrp} we define formal Lie groups as group 
objects in the category of locally convex spaces whose morphisms are 
formal power series without constant term. 
We show that every formal Lie group $G$ has a formal exponential 
function that implements an isomorphism to the formal group 
$\Ch(\g)$, on which the formal multiplication is given by the Hausdorff 
series in terms of the Lie bracket on~$\g$.

Section~\ref{sec:8.3} is devoted to 
Lie group extensions $\hat G$ of $G$ by a regular Lie group~$N$. 
In particular, we show that, if $G$ has an exponential function, 
then so does $\hat G$. 
We also obtain a relative Fundamental Theorem on the 
integrability of Lie algebra-valued $1$-forms 
satisfying the Maurer--Cartan equation, even if the 
group $\hat G$ is not known to be regular. Then 
integrability modulo $N$ implies integrability for~$\hat G$. 
We also provide criteria for the existence of abelian extensions  
with global smooth cocycles and develop 
methods to deal with central extensions without any smooth local 
sections. On the Lie algebra level this leads to the concept of a 
generalized central extension. 

In Section~\ref{sec:11.4} we characterize 
enlargeable locally exponential Lie algebras~$\g$ 
by the discreteness of an additive subgroup 
$\Pi(\g) \subeq \fz(\g)$, called the period group. 
We also describe several methods to calculate these groups 
and construct examples of non-enlargeable locally exponential 
Lie algebras. 

In Section~\ref{sec:11.5} we discuss integrability of a 
locally convex Lie algebra $\g$, i.e., the existence of a 
Lie group $G$ with Lie algebra $\g$, without assuming any additional properties, 
such as local exponentiality. In this context only rather 
 isolated results are available. In the context of analytic groups, we 
provide several non-integrability results for complex locally convex 
Lie algebras. For a continuous linear operator $D$ on a locally convex space 
$E$, the integrability of the Lie algebra $E \rtimes_D \R$ 
can be characterized in terms of the integrability of the operator $D$, 
i.e., by the existence of a smooth one-parameter groups $\R \to \GL(E)$ 
with infinitesimal generator~$D$. 

The analytic version of locally exponential Lie groups, BCH--Lie groups, 
are briefly discussed in Section~\ref{sec:BCH}. 

We already encountered Lie groups whose exponential function 
is a global diffeomorphism in the context of nilpotent and pro-nilpotent 
Lie groups. In Section~\ref{sec:banach-reg-exp} we discuss the global 
regularity properties of the exponential function for Banach--Lie groups. 
This leads in particular to a Banach version of the Dixmier--Saito Theorem 
characterizing finite-dimensional exponential Lie groups. 
For complex Banach--Lie algebras we prove W.~Wojty\'nski's remarkable 
theorem asserting that quasi-nilpotence is equivalent to exponentiality 
and we also provide a characterization of exponential projective limits 
of finite-dimensional Lie algebras. 

The final Section~\ref{sec:14.9} discusses the failure of local 
exponentiality of semidirect products of locally exponential Lie groups.

\section{Nilpotent and solvable Lie groups} \mlabel{sec:3.6} 

In Subsection~\ref{subsec:8.1.1}  we show that a connected Lie group 
$G$ is nilpotent, resp., solvable if and only if its Lie algebra $\g$ 
has this property. Moreover, the degree of nilpotency, resp., solvability 
is the same for the group and its Lie algebra. 
In Subsection~\ref{subsec:9.1.2} we show that, for solvable Lie groups, the 
existence of an exponential function and the Mackey completeness of $\g$ 
implies regularity.

\subsection{Characterization in terms of the Lie algebra} 
\mlabel{subsec:8.1.1} 

For two subsets $A,B$ of a group $G$ we write 
$(A,B)$ for the subgroup generated by all commutators 
$(a,b) := aba^{-1}b^{-1}$, $a \in A$, $b \in B$. 
We define the 
\index{derived series of $G$, $D^n(G)$} 
\index{central series of $G$, $C^n(G)$} 
\index{group!solvable} 
\index{solvable length} 
\index{nilpotent length} 
\index{group!nilpotent} 
{\it derived series of $G$} inductively by 
$$ D^0(G) := G, \quad D^{n+1}(G) := (D^n(G), D^n(G)), n \in \N_0, $$
and the {\it lower (or descending) central series} by 
$$ C^1(G) := G, \quad C^{n+1}(G) := (G, C^n(G)), n \in \N_0. $$
A group $G$ is said to be {\it solvable} if $D^N(G) = \{\be\}$ for 
some $N \in \N_0$. The minimal $N$ with this property is called 
the {\it solvable length of $G$}. 
A group $G$ is called {\it nilpotent} if 
$C^{N+1}(G) = \{\be\}$ for some $N \in \N_0$. The minimal such $N$ is called 
its {\it nilpotent length}. 

\index{derived series of Lie algebra $\g$, $D^n(\g)$} 
\index{central series of Lie algebra $\g$, $C^n(\g)$} 
\index{Lie algebra!solvable} 
\index{Lie algebra!nilpotent} 
For a Lie algebra $\g$ we likewise define 
the {\it derived series} inductively by 
$$ D^0(\g) := \g, \quad D^{n+1}(\g) 
:= [D^n(\g), D^n(\g)], n \in \N_0, $$
and the {\it lower (or descending) central series} by 
$$  C^1(\g) := \g, \quad  C^{n+1}(\g) := 
[\g,C^n(\g)], n \in \N_0. $$
The Lie algebra $\g$ is said to be {\it solvable} if $D^N(\g) = \{0\}$ for 
some $N \in \N$, and the minimal such $N$ is called its 
{\it solvable length}. Likewise we say that $\g$ is {\it nilpotent 
(of length $N$)} if $C^{N+1}(\g) = \{0\}$ 
for some $N \in \N_0$ (cf.\ Exercises~\ref{ex:e.8}, \ref{ex:e.9}). 
If, in addition, $\g$ is a topological Lie algebra, we also put 
$$ \oline D^n(\g) := \oline{D^n(\g)} \quad \mbox{ and } \quad 
\oline C^n(\g) := \oline{C^n(\g)} $$
and observe that these are closed ideals invariant under all 
automorphisms. 

In this section we collect some results on nilpotent and solvable Lie groups. 
Although nilpotent Lie groups behave in all respects very much 
like in the finite-dimensional case, the theory of solvable Lie groups 
faces more analytic obstacles, which is partly due to the 
fact that their universal covering group need not be exponential. 

Throughout this section $G$ denotes a connected Lie group and 
$\g = \L(G)$ its Lie algebra. Our first goal is to show that 
the group $G$ is nilpotent, resp., solvable if and only if its 
Lie algebra is nilpotent, resp., solvable.

\begin{thm} \mlabel{thm:nilpo-liealg} A connected Lie group $G$ is nilpotent if and only if 
its Lie algebra $\g$ is nilpotent. 
\end{thm}

\begin{prf}
Inductively, we define a sequence of smooth functions 
$f_n \: G^{n+1} \to G$ by 
\[  f_1(g_1, g_2) := g_1 g_2 g_1^{-1} g_2^{-1},\]
which is the {\it group commutator function},  \index{commutator!group}
and 
\[ f_{n+1}(g_0,\ldots, g_{n+1}) := f_1(g_0, f_n(g_1,\ldots, g_{n+1})),\]
which is the $(n+1)$-st {\it commutator function on $G$}. 
A group $G$ is nilpotent with vanishing $C^{N+1}(G)$ 
if and only if the function $f_N$ is constant $\be$ (Exercise~\ref{ex:e.5}). 

Clearly the differential $\dd f_1(\be,\be)$ vanishes, so that its 
$2$-jet is a  well-defined quadratic map 
$\g \times \g \to \g$, given by 
\[  \frac{1}{2} \delta^2_{(\be,\be)} f_1(x_1,x_2) =  [x_1,x_2] \]
(Remark~\ref{rem:brack-taylor}). 
With the Chain Rule for Taylor polynomials 
(Proposition~\ref{chainrtay}) and induction, one easily derives that 
the $n$-jet of $f_n$ in $(\be,\ldots, \be)$ vanishes, and that its 
$(n+1)$-jet is given by 
\begin{eqnarray}
  \label{eq:e.1}
\frac{1}{(n+1)!} \delta^{n+1}_{(\be,\ldots, \be)} f_n (x_0,x_1, \ldots, x_n) 
= [x_0,[x_1,[x_2, \ldots ,[x_{n-1},x_n]]]]. 
\end{eqnarray}
From (\ref{eq:e.1}) we immediately derive that if $G$ is nilpotent, 
then $\g$ is nilpotent: $f_N = \be$ implies that 
all $(N+1)$-fold iterated brackets of elements of $\g$ vanish, 
which means that $\g$ is a nilpotent Lie algebra with 
$C^{N+1}(\g) = \{0\}$ (cf.\ Exercise~\ref{ex:e.8}). 

Now we assume, conversely, that $\g$ is nilpotent with $C^N(\g) = \{0\}$. 
We consider the adjoint representation $\Ad \: G \to \Aut(\g)$, 
which leads to a short exact sequence of groups 
\[  \1 \to Z(G) \to G \to \Ad(G) \to \1 \] 
(Corollary~\ref{cor:3.2.8}). 
Since central extensions of nilpotent groups are nilpotent 
(Exercise~\ref{ex:e.4}(2)), it 
suffices to show that $\Ad(G)$ is a nilpotent group. This will 
be accomplished by showing that $\Ad(G)$ is contained in the 
unipotent group of a finite flag of subspaces of $\g$ (cf.\ Exercise~\ref{ex:e.3}). 

Clearly, each subspace $\oline C^n(\g)$ is a closed ideal of 
$\g$ which is invariant 
under each automorphism, and in particular under the group 
$\Ad(G)$. Therefore the adjoint action defines for each 
$n \in \N$ a smooth representation of $G$ on the 
space $\oline C^{n}(\g)/\oline C^{n+1}(\g)$ (Definition~\ref{def:3.1.2} and 
Lemma~\ref{lemquot}), and this action is trivial 
because $[\g,\oline C^{n}(\g)] \subeq \oline C^{n+1}(\g)$ 
follows immediately from the 
continuity of the bracket on $\g$ (Proposition~\ref{prop:e.2.6}). 
We conclude that $(\Ad(G) - \id_\g)\oline C^{n}(\g) \subeq \oline C^{n+1}(\g)$ 
holds for each $n \in \N$, so that $\Ad(G)$ is contained in the 
unipotent group $U({\cal F})$ of the flag 
$$ {\cal F} \: \quad 
\0 = \oline C^N(\g) \subeq \oline C^{N-1}(\g) 
\subeq \ldots \subeq \oline C^2(\g) \subeq \oline C^1(\g) = \g,$$
and thus $C^{N-1}(\Ad(G)) = \{\1\}$. 
This implies that $C^N(G) =\{\be\}$, and therefore $G$ is nilpotent. 
\end{prf}

\begin{thm} \mlabel{thm:solv-grpalg} A connected Lie group $G$ is solvable if and only if 
its Lie algebra $\g$ is solvable. 
\end{thm}

\begin{prf} The first half of the proof follows the corresponding 
argument in the nilpotent case. 
We consider the sequence of smooth functions 
$h_n \: G^{2^n}~\to~G$ defined by 
$$ h_1(g_1, g_2) :=g_1 g_2 g_1^{-1} g_2^{-1}$$
and 
\[  h_{n+1}(g_1,\ldots, g_{2^{n+1}}) := h_1\big(
h_n(g_1,\ldots, g_{2^n}), h_n(g_{2^n+1}, \ldots,  g_{2^{n+1}})\big). \]
A Lie group $G$ is solvable with $D^N(G) = \{\be\}$ 
if and only if the 
function $h_N$ is identically $\be$ (Exercise~\ref{ex:e.6}). 
Again, we use the Chain Rule for Taylor polynomials 
(Proposition~\ref{chainrtay}) to derive that 
the $(2^n-1)$-jet of $h_n$ in $(\be,\ldots, \be)$ vanishes and that 
(in any local chart) 
\begin{eqnarray}
  \label{eq:e.2}
 \frac{1}{2^n!} \delta_{(\be,\ldots, \be)}^{2^n} h_n \: \g^{2^n} \to \g 
\end{eqnarray}
equals the function $b_n$, defined inductively by 
$b_1(x_1, x_2) := [x_1, x_2]$ and 
$$ b_n(x_1, \ldots, x_{2^n}) 
:= [b_{n-1}(x_1,\ldots, x_{2^{n-1}}), b_{n-1}(x_{2^{n-1}+1}, \ldots,  x_{2^{n}})]. $$
With \eqref{eq:e.2} it follows that $D^N(G) = \{\be\}$ implies that 
$b_N$ vanishes, i.e., that $\g$ is solvable with $D^N(\g) = \{0\}$  
(Exercise~\ref{ex:e.9}). 

Assume, conversely, that $\g$ is a solvable Lie algebra. 
Let $I = [0,1]$ denote the unit interval and consider two 
smooth curves $\alpha,\beta \: I \to G$ starting in $\be$. 
We are interested in their commutator 
$$ \gamma \: I \to G, \quad t \mapsto \alpha(t)\beta(t) \alpha(t)^{-1}\beta(t)^{-1}. $$
In view of 
$$ \delta(\alpha\beta) = \delta(\beta) + \Ad(\beta)^{-1}\delta(\alpha) 
\quad \mbox{ and } \quad 
\delta(\alpha^{-1}) = - \Ad(\alpha)\delta(\alpha) $$
(Lemma~\ref{lem:c.12}), we have 
$$ \delta(\alpha^{-1}\beta^{-1}) 
= \delta(\beta^{-1}) + \Ad(\beta)\delta(\alpha^{-1})
= -\Ad(\beta)\delta(\beta) - \Ad(\beta\alpha)\delta(\alpha). $$
This leads to 
\begin{eqnarray*}
 \delta(\gamma) 
&&= \delta(\alpha^{-1}\beta^{-1}) + \Ad(\beta\alpha) \delta(\alpha\beta) \\ 
&&=  -\Ad(\beta)\delta(\beta) - \Ad(\beta\alpha)\delta(\alpha) 
+ \Ad(\beta\alpha)\big(\delta(\beta) + \Ad(\beta)^{-1}.\delta(\alpha)\big) \\
&&=  -\Ad(\beta)\delta(\beta) - \Ad(\beta\alpha)\delta(\alpha) 
+ \Ad(\beta\alpha)\delta(\beta) + \Ad(\beta\alpha\beta^{-1})\delta(\alpha) \\
&&= \Ad(\beta)\big(\Ad(\alpha)\delta(\beta)-\delta(\beta)\big)  
+ \Ad(\beta\alpha\beta^{-1})\big(\delta(\alpha)-\Ad(\beta)\delta(\alpha)\big). 
\end{eqnarray*}

Suppose that $\delta(\alpha)$ and $\delta(\beta)$ have values in $\oline D^n(\g)$. 
For a fixed $t \in I$ we then consider the smooth curve 
$f \: I \to \oline D^n(\g)$ given by 
$$ f(s) := \Ad(\alpha(st))\delta(\beta)(t)-\delta(\beta)(t) $$
and note that 
$$ f'(s) = t \Ad(\alpha(st))[\delta(\alpha)(st),\delta(\beta)(t)]
\in \oline D^{n+1}(\g). $$
In view of $f(0) = 0$, this implies that 
$$ \Ad(\alpha(t))\delta(\beta)(t)-\delta(\beta)(t)
= f(1) \in \oline D^{n+1}(\g). $$
We likewise find that 
$\Ad(\beta)\delta(\alpha) - \delta(\alpha)$ has values in 
$\oline D^{n+1}(\g)$, 
and hence that $\delta(\gamma)$ has values in $\oline D^{n+1}(\g)$. 

Inductively, it now follows that if 
$\alpha_i \: I \to G$, $i =1,\ldots, 2^n$, are smooth curves, then 
the curve 
\[  \gamma := h_n \circ (\alpha_1, \ldots, \alpha_{2^n}) \] 
satisfies $\im(\delta(\gamma))\subeq \oline D^n(\g)$. 
For $n = N$, we derive in 
particular from \break $D^N(\g) = \{0\}$ 
that $\gamma$ is constant. Since $G$ is connected, 
it follows that $h_N$ is constant $\be$ (Lemma~\ref{lem:c.12b}), 
so that $G$ is solvable with $D^N(G)= \{\be\}$. 
\end{prf}

\begin{rem} For finite-dimensional solvable Lie algebras, the fact that the 
solvability of $\g$ implies the solvability of $G$ can be proved very much 
in the spirit of our argument in the nilpotent case: Since solvability is 
an extension property, it suffices to show that $\Ad(G)$ is solvable. 
We may w.l.o.g.\ assume that the Lie algebra $\g$ is complex. 
Then Lie's Theorem implies the existence of a 
\index{complete flag} 
{\em complete} 
invariant flag 
$$ {\cal F} \: F_0 = \0 \subeq F_1 \subeq \ldots \subeq F_N = \g
\quad \mbox{ with } \quad \dim F_n = n,$$ 
consisting of ideals of $\g$. 
Hence 
$$ \Ad(G) \subeq 
\Aut({\cal F}) = \{ \phi \in \GL(\g) \: (\forall n)\ \phi(F_n) = F_{n}\}, $$
and since $\Aut({\cal F})$ is solvable (Exercise~\ref{ex:e.3}(2)), the assertion follows. 
\end{rem} 

\begin{ex} \mlabel{ex:der-non-nilp} 
Let $E := C^\infty(\R,\C)$ denote the space of smooth 
complex valued functions on $\R$, endowed with its natural Fr\'echet 
topology (Definition~\ref{def:smooth-co-top}). Further, let 
$\h := \Spann_\C \{p,q,z\}$ denote the complex $3$-dimensional Heisenberg algebra 
whose bracket is given by 
$$ [p,q] = z \quad \mbox{ and }\quad [p,z] = [q,z] =0. $$ 
We then have a natural representation $\alpha$ of $\h$ on $E$ given by 
$$ \alpha(z) = \id_E, \quad 
\alpha(p) = \frac{d}{dt} \quad \mbox{and}\quad 
(\alpha(q)f)(t)  = t f(t). $$
Now the semidirect sum 
$$ \g := E \rtimes_\alpha \h $$
is a locally convex Lie algebra whose derived series is given by 
$$ D^1(\g) = E \rtimes \C z, \quad 
 D^2(\g) = E \quad \mbox{ and } \quad 
 D^3(\g) = \0. $$
It follows in particular that $\g$ is solvable. 
We observe that $D^1(\g)$ is not nilpotent, 
which never happens for finite-dimensional Lie algebras. 
The fact that the commutator $z = [p,q]$ does not act in a nilpotent 
fashion on $\g$ further implies that there is no 
``complete'' flag of $\g$ consisting of ideals. Moreover, 
the derived algebra $D^1(\g)$ is not nilpotent, which is always the case 
for finite-dimensional solvable Lie algebras.   
\end{ex}

\subsection{Regularity of solvable Lie groups} 
\mlabel{subsec:9.1.2}

In this section, we prove that for each solvable Lie group $G$ 
whose Lie algebra $\g$ is Mackey complete, the existence of a smooth 
exponential function implies regularity. 
For the following, we recall the definition of the derived series 
$\oline D^k(\g) := \oline{D^k(\g)}$ from Subsection~\ref{subsec:8.1.1}. 

\begin{lem} \mlabel{lem:stab-derser} Let 
$\g$ be a locally convex Lie algebra. 
Then each $\oline D^k(\g)$ is a closed characteristic ideal of $\g$, i.e., 
invariant under all topological automorphisms of $\g$. 
\end{lem}

\begin{lem} \mlabel{lem:d1-red} Let $G$ be a Lie group with an exponential 
function and assume that $\g = \L(G)$ is Mackey complete. 
Further let $\fn \trile \g$ be an $\Ad(G)$-invariant closed ideal. 
Then the following assertions hold: 
\begin{description}
\item[\rm(i)] For each $k \in \N \cup \{\infty\}$, the map 
$$\Phi \: C^k([0,1],\g) \to C^{k+1}([0,1],G), \quad 
\Phi(\xi)(t) := \exp_G\Big(\int_0^t \xi(s)\, ds\Big) $$
is smooth. 
\item[\rm(ii)] For $\xi \in C^k([0,1],\fn)$ and $\eta := \Phi(\xi)$, 
the map $\xi_1 := \Ad(\eta)(\xi - \delta(\eta))$ is a $C^{k}$-curve with 
values in $\oline D^1(\fn)$. 
\item[\rm(iii)] A curve $\eta_1 \in C^1([0,1],G)$ satisfies 
$\delta(\eta_1\eta) = \xi \in C^1([0,1],\fn)$ if and only if 
$$ \delta(\eta_1) = \Ad(\eta)(\xi - \delta(\eta)) = \xi_1 \quad \mbox{ in } \quad 
C([0,1], \oline D^1(\fn)). $$
\item[\rm(iv)] If $\evol_1 \: C^k([0,1],\oline D^1(\fn)) \to G$ is a smooth evolution map, 
then 
$$ \evol_G(\xi) 
:= \evol_1(\xi_1)\cdot \Phi(\xi)(1) 
= \evol_1(\xi_1)\cdot \exp_G\Big(\int_0^1 \xi(s)\, ds\Big) $$
defines a smooth evolution map 
$C^k([0,1],\fn) \to G$. 
\end{description}
\end{lem}

\begin{prf} (i) Since $\exp_G$ is smooth, composition with 
$\exp_G$ induces a smooth map $C^{k+1}([0,1],\g) \to C^{k+1}([0,1],G)$. 
In fact, since it suffices to see that this map is smooth on some $0$-neighborhood 
of the form $C^{k+1}([0,1],U_\g)$, $U_\g \subeq \g$ an open $0$-neighborhood 
(Lemma~\ref{lem:5.1.2}), 
the assertion follows from Proposition~\ref{superctisCk}.
On the other hand, for $k \geq 1$, the integration map 
$C^k([0,1],\g)\to C^{k+1}([0,1],\g)$ is continuous 
and linear, hence smooth. 

(ii) Let $q \: \fn \to \fn/\oline D^1(\fn)$ denote the quotient map and 
$\tilde\xi(t) := \int_0^t \xi(s)\, ds$. Then 
$q$ is a morphism onto an abelian Lie algebra, so that 
$q \circ \ad x= 0$ for each $x\in \fn$. 
If $\gamma \in C^1([0,1],G)$ satisfies $\delta(\gamma)([0,1]) \subeq \fn$ 
and $\gamma(0) = \be$, we obtain for each $x \in \fn$ the relation 
$$ \frac{d}{dt} \Ad(\gamma(t))x 
= \Ad(\gamma(t))[\delta(\gamma)_t, x] 
= [\Ad(\gamma(t))\delta(\gamma)_t, \Ad(\gamma(t))x] $$
(Proposition~\ref{prop:der-Ad}), which leads to 
$\frac{d}{dt} q\big(\Ad(\gamma(t))x\big) = 0.$
Hence the curve $q\big(\Ad(\gamma(t))x\big)$ in $\fn/\oline D^1(\fn)$ is constant for 
each $x$, and this means that \break $q \circ \Ad(\gamma(t))= q$ for each $t$. 
In particular, 
$q \circ e^{\ad x} = q$ for $x \in \fn$ and thus 
$$ q \circ \kappa_\fn(x) 
= q \circ \int_0^1 e^{-t\ad x}\, dt 
= \int_0^1 q\circ e^{-t\ad x}\, dt 
= \int_0^1 q\, dt = q $$
for each $x \in \fn$. Here we use the notation from 
Definition~\ref{def:kappa-liealg}, resp., Definition~\ref{def4.1.5}. 
For the $C^1$-curve $\eta = \exp_G \circ \tilde\xi$, we obtain 
with Theorem~\ref{thm:exp-logder} 
$$ \delta(\eta)_t 
= (\tilde\xi^*\delta(\exp_G))_t 
= \delta(\exp_G)_{\tilde\xi(t)}\tilde\xi'(t) 
= \kappa_\g(\tilde\xi(t))\xi(t) 
= \kappa_\fn(\tilde\xi(t))\xi(t) \in \fn. $$ 
In view of the preceding discussion, this 
implies that $q \circ \Ad(\eta(t)) = q$ for each~$t$, so 
that we arrive at 
\begin{eqnarray*}
q(\xi_1(t)) 
&=& q\big(\xi(t) - \delta(\eta)_t\big)
= q(\xi(t)) - q\big(\kappa_\fn(\tilde\xi(t))\xi(t)\big) \\ 
&=& q(\xi(t)) - q(\xi(t)) = 0. 
\end{eqnarray*}

(iii) This follows immediately from 
$\delta(\eta_1 \cdot \eta) 
= \delta(\eta) + \Ad(\eta)^{-1}.\delta(\eta_1).$

(iv) In view of (iii), it remains to show that the map 
$$C^k([0,1],\fn) \to C^k([0,1],\oline D^1(\fn)), \quad 
\xi \mapsto \xi_1$$ is smooth.
This follows from the smoothness of the adjoint action of 
the group $C^k([0,1],G)$ (Lemma~\ref{lem:smoothprop}) 
and the smoothness of the logarithmic derivative 
(Proposition~\ref{prop:smooth-logder0}). 
\end{prf}

\begin{prop} \mlabel{prop:solv-regular} Let $G$ be a solvable Lie group with 
exponential function 
whose Lie algebra $\g$ is Mackey complete. Then $G$ is $C^1$-regular. 
If, in addition, $\g$ is sequentially complete, then $G$ is 
$C^0$-regular. 
\end{prop}

\begin{prf} In view of Theorem~\ref{thm:solv-grpalg}, the 
Lie algebra $\g$ is solvable. Let $N$ be its solvable length. 
We argue by induction on $n = N,N-1,\ldots, 0$ that there exists a smooth evolution map 
$$ \evol_G \: C^1([0,1],\oline D^n(\g)) \to G. $$
For $n = N$ there is nothing to show. If the assertion holds 
for some $n > 0$, then Lemma~\ref{lem:d1-red}(iv) implies that 
it also holds for $n-1$ because the ideal $\oline D^{n-1}(\g)$ of $\g$ is closed 
and $\Ad(G)$-invariant (Lemma~\ref{lem:stab-derser}). The assertion follows for $n = 0$. 

If $\g$ is sequentially complete, the integrals of continuous curves exist, 
so that $k = 0$ can be permitted in Lemma~\ref{lem:stab-derser} in this case. 
\end{prf}

In Theorem~\ref{thm:reg-ext-prop} 
we have seen that, in the context of Lie group extensions (which are 
locally trivial smooth fiber bundles), regularity 
is an extension property. 
Since we do not want to assume that all ideals 
$\oline D^n(\g)$ are topologically split, the preceding proposition 
does not follow from Theorem~\ref{thm:reg-ext-prop}. 

\begin{ex} We give an example of a solvable Banach--Lie algebra $\g$, 
where the ideal $\oline D^1(\g)$ does not split topologically. 
Consider the space $\g := \ell^\infty(\N,\R)$, endowed with the sup-norm. 
We write $\g = \g_0 \oplus \R x$, where $\g_0$ is a closed subspace containing 
the non-split subspace $c_0(\N,\R)$ of sequences converging to $0$. 
We define a continuous linear operator 
\[ T \: \g \to \g, \quad T(x)_n = \lambda_n x_n,\] 
where $\lim_{n \to \infty} \lambda_n = 0$ but $\lambda_n \not=0$ for each $n$. 
Note that $T(\g) \subeq \g_0$ implies 
that $T$ defines a continuous linear operator on $\g_0$. We now consider the 
Lie algebra structure on $\g$ given by 
$\g := \g_0 \rtimes_T \R$. Then $[\g,\g] = T(\g_0) \subeq c_0(\N,\R)$ 
contains $e_n = \lambda_n^{-1} T(e_n)$, which leads to 
$\oline D^1(\g) = c_0(\N,\R)$. 
In view of $\Spec(T) \subeq \R$, the Lie algebra $\g$ is an exponential 
Banach--Lie algebra of solvable length $2$ 
(cf.\ Proposition~\ref{prop:semdir-exp} and Example~\ref{ex:bad-expfct}). 
\end{ex}

\begin{small}
\subsection*{Exercises for Section~\ref{sec:3.6}} 

\begin{exer} \mlabel{ex:e.3a} Let $E$ be a vector space and 
$$ {\cal F} \: F_0 = \0 \subeq F_1 \subeq \ldots \subeq F_N = E $$
a finite flag of linear subspaces. Show that
$$ \fu({\cal F}) := \{ \phi \in \gl(E) \: (\forall n)\ \phi(F_n) \subeq 
F_{n-1}\} $$
is a nilpotent Lie algebra.\\
 Hint: Show inductively that 
$C^j(\fu({\cal F}))(F_k) \subeq F_{k-j}$. 
\end{exer}

\begin{exer} \mlabel{ex:e.3} Let $E$ be a vector space and 
$$ {\cal F} \: F_0 = \0 \subeq F_1 \subeq \ldots \subeq F_N = E $$
a finite flag of subspaces. We write 
$$ \Aut({\cal F}) := \{ \phi \in \GL(E) \: (\forall n)\ \phi(F_n) = F_{n}\} $$
for the group of 
\index{automorphisms of flag} 
{\em automorphisms of this flag}. 
Show that:
\begin{description}
\item[(1)] The set 
$$ U({\cal F}) := \{ \phi \in \Aut({\cal F}) \: (\forall n)\ (\phi-\id_E)(F_n) \subeq 
F_{n-1}\} $$ \index{group!unipotent} 
is a nilpotent group with $C^N(\U(\cal F))= \{\1\}$. 
It is called the {\it unipotent group of ${\cal F}$}. \\
Hint: If $(\phi - \id_E)(F_j) \subeq F_{j-m}$ and 
$(\psi - \id_E)(F_j) \subeq F_{j-n}$ for all $j$, then the commutator 
$(\phi,\psi) = \phi\psi\phi^{-1}\psi^{-1}$ 
satisfies $((\phi,\psi)-\id_E)(F_j) \subeq F_{j-n-m}$. Use that 
\begin{align*}
(\phi,\psi) -\id_E = \phi\psi\phi^{-1}\psi^{-1} -\id_E 
&= (\phi\psi - \psi\phi)\circ \phi^{-1}\psi^{-1}
= [\phi,\psi] \circ \phi^{-1}\psi^{-1} \\
&= [\phi-\id_E,\psi-\id_E] \circ \phi^{-1}\psi^{-1} 
\end{align*}
and Exercise~\ref{ex:e.3a}.
\item[(2)] If $E$ is finite-dimensional and ${\cal F}$ is complete, 
i.e., $\dim F_j = j$ for each $j$, then the group $\Aut({\cal F})$ is solvable. 
\end{description}
\end{exer}

\begin{exer} \mlabel{ex:e.4} Let $G$ be a group. Show that: 
\begin{description}
\item[(1)] If $G$ is nilpotent, then each subgroup and each quotient group of $G$ is nilpotent. 
\item[(2)] Any central extension of a nilpotent group is nilpotent. 
\item[(3)] If $G$ is solvable, then each subgroup and each quotient group of $G$ 
is solvable. 
\item[(4)] Any extension of a solvable group by a solvable group is solvable 
(solvability is an extension property). 
\item[(5)] Find an example of an extension of an abelian group by an abelian group 
which is not nilpotent (nilpotency is not an extension property). 
\end{description}
\end{exer}

\begin{exer} \mlabel{ex:e.5} Show that, for a group $G$, the following 
are equivalent: 
\begin{description}
\item[(1)] $G$ is nilpotent, i.e., $C^N(G) = \{\be\}$ for some $N \in \N$. 
\item[(2)] There exists an $N \in \N$ such that all iterated commutators of the form 
$(g_1,(g_2,\ldots, (g_{N-1},g_N)\ldots ))$ are trivial. 
\end{description}
\end{exer}

\begin{exer} \mlabel{ex:e.6} Show that, for a group $G$, the following 
are equivalent: 
\begin{description}
\item[(1)] $G$ is solvable, i.e., $D^N(G) = \{\be\}$ for some $N \in \N$. 
\item[(2)] There exists an $N \in \N$ such that all binary iterated commutators of the form 
$$((\cdots ((g_1, g_2), (g_3, g_4)) \cdots),(\cdots (g_{2^N-1},g_{2^N}) \cdots )) $$
are trivial. 
\end{description}
\end{exer}

\begin{exer} \mlabel{ex:e.7} Let $\g$ be a Lie algebra. Show that: 
\begin{description}
\item[(1)] If $\g$ is nilpotent, then each subalgebra and each quotient algebra of 
$\g$ is nilpotent. 
\item[(2)] Any central extension of a nilpotent Lie algebra is nilpotent. 
\item[(3)] If $\g$ is solvable, then each subalgebra
 and each quotient algebra of $\g$ 
is solvable. 
\item[(4)] Any extension of a solvable Lie algebra by a solvable Lie algebra is solvable 
(solvability is an extension property). 
\item[(5)] Find an example of an extension of an abelian Lie algebra by an abelian 
Lie algebra which is not nilpotent (nilpotency is not an extension property). 
\end{description}
\end{exer}

\begin{exer} \mlabel{ex:e.8} Show that, for a Lie algebra $\g$, 
the following conditions 
are equivalent: 
\begin{description}
\item[(1)] $\g$ is nilpotent, i.e., $C^{N+1}(\g) = \{0\}$ for some $N \in \N$. 
\item[(2)] There exists an $N \in \N$ such that all 
iterated commutators of the form 
$[x_1,[x_2,\ldots, [x_{N-1},x_N]\ldots ]]$ are trivial. 
\end{description}
\end{exer}

\begin{exer} \mlabel{ex:e.9} Show that, for a Lie algebra $\g$, 
the following conditions are equivalent: 
\begin{description}
\item[(1)] $\g$ is solvable, i.e., $D^N(\g) = \{0\}$ for some $N \in \N$. 
\item[(2)] There exists an $N \in \N$ such that all binary iterated commutators of the form 
$$[[\cdots [[x_1, x_2], [x_3, x_4]] \cdots],[\cdots [x_{2^N-1},x_{2^N}] \cdots ]] $$
vanish. 
\end{description}
\end{exer}
  
\end{small}

\nin{\bf Notes on Section~\ref{sec:3.6}:} \\ 
\nin The exercises in this section show that nilpotent Lie algebras 
arise naturally from finite flags $\cF$ because the Lie algebra 
$\fu(\cF)$ is nilpotent (Exercise~\ref{ex:e.3a}). On the other 
hand, a Lie algebra $\g$ is nilpotent if and only if 
the adjoint representation maps into $\fu(\cF)$ for a finite flag 
$\cF$ in~$\g$. This leads to the question whether any nilpotent 
locally convex Lie algebra has a faithful representation in some 
$\fu(\cF)$. As nilpotent Lie algebras have non-zero center, the 
adjoint representation is not faithful. However, an affirmative answer 
to the above question has been obtained by 
D.\ and I.~Belti\c{t}\u{a} in \cite{BB15}. They construct 
for every nilpotent locally convex Lie algebra a faithful 
representation on certain spaces of polynomial maps of finite degree 
on $\g \to \R$ 
whose range is contained in $\fu(\cF)$ for a suitable flag 
(determined by the degree). For a nilpotent Banach--Lie algebra, 
their construction provides a faithful Banach representation 
in some $\fu(\cF)$.

\section{Polynomial Lie groups} \mlabel{sec:polgrp}

In this section we study Lie groups defined by a polynomial 
multiplication of a locally convex space $G$ for which the degrees of the 
iterated multiplication maps are bounded. The main result is that 
such a Lie group $G$ has a polynomial exponential map 
$\Exp \: \g \to G$ with a global polynomial inverse and that 
its Lie algebra $\g$ is nilpotent. As a consequence, 
every polynomial Lie group $G$ is polynomially isomorphic to a 
Lie group of the form $(\g,*)$, where $\g$ is a nilpotent Lie algebra,
endowed with the BCH multiplication. 

\index{Lie group!real polynomial}
\begin{defn} A {\it real polynomial Lie 
group}\begin{footnote} {Although we shall discuss polynomial groups only 
as examples 
of Lie groups, the theory of polynomial groups can be developed over 
any rational vector space and no topology is needed.}
\end{footnote}is a real locally convex 
space~$G$, with a continuous polynomial group structure 
$m_G \: G \times G \to G$, for which $0$ is the identity element 
and all iterated product maps 
$$ m_G^{(j)} \: G^j \to G, \quad (x_1,\ldots, x_j) \mapsto 
x_1\cdots x_j \quad \mbox{ for } \quad j \in \N,$$ 
are polynomial maps whose degree is bounded by some $d \in \N$. 
\end{defn} 

We shall see below that the above requirements imply that 
the inversion map $\eta_G \: G \to G$ is also polynomial of 
degree $\leq d$ (Corollary~\ref{cor:inv-pol}). 

\begin{exs} \mlabel{ex:polygrp} (a) Let $\cA$ be a nilpotent
  associative locally convex 
algebra with $a_1\cdots a_d = 0$ for 
$a_1,\ldots, a_d \in \cA$. Then 
$m_\cA(x,y) := x + y + xy$ defines a group structure on $\cA$ for which inversion 
is given by 
$$ \eta_\cA(x) := \sum_{i = 0}^{d-1} (-1)^i x^i, $$
as is easily verified directly by identifying $\cA$ with the affine subspace 
$\cA + \1$ of the algebra $\cA \times \R$ with unit $\1 = (0,1)$, endowed 
with the 
multiplication $(a,t)(a',t') = (aa' + ta'+t'a,tt')$ 
(Exercise~\ref{exer:3.1.4}). 

Then $(\cA,m_\cA)$ is a polynomial Lie group 
and the degree of all the maps $m_\cA^{(j)}$ is bounded by $d-1$ because 
$d$-fold products of elements of $\cA$ vanish. 

(b) If $G$ is a Lie group with Lie algebra $\g$ 
and $d \in \N$, then we know already 
that the iterated tangent bundle $T^d(G)$ is a Lie group for which the 
canonical projection $\pi \: T^d(G) \to G$ is a morphism of Lie groups 
(Proposition~\ref{prop:tangentgrp}). 
The fiber $T^d_\be(G) := \ker \pi$ also carries a natural group 
structure, and by using the right trivializations, it also inherits the 
structure of a locally convex space, isomorphic to $\g^{2^d-1}$. 
One can show that the group structure on $T^d_\be(G)$ is polynomial 
of degree $\leq d$. For $d = 1$, we obtain the abelian group  
$T_\be(G) \cong (\g,+)$, and for $d = 2$ we get the $2$-step nilpotent group 
$$ T^2_\be(G) \cong \g^3, \quad 
(x_{12}, x_2, x_1) * (x_{12}', x_2', x_1') 
= (x_{12} + x_{12}' + [x_1, x_2'], x_2+x_2', x_1+x_1'). $$

A good way to understand the group structure on $T_\be^d(G)$ is to 
use a $\g$-chart $(\phi,U)$ of $G$ with $\phi(\be) = 0$ and 
$T_\be(\phi) = \id_\g$. The Lie algebra 
$\L(T^d(G))$ is isomorphic to 
\[ T^d(\g) := \g \otimes_\R T^d(\R), \quad \mbox{ where } \quad 
T^d(\R) := \R[\eps_1,\ldots, \eps_d]\]
is the $2^d$-dimensional commutative 
algebra with generators $\eps_i$ satisfying $\eps_i^2~=~ 0$. 
Writing $\eps^\alpha := \eps_1^{\alpha_1} \cdots \eps_d^{\alpha_d}$ 
for $\alpha \in \{0,1\}^d$, we then have 
$$ T^d_\be(G) \cong \bigoplus_{\alpha\not=0} \g\otimes \eps^\alpha. $$
The multiplication on the group $T_\be^d(G)$ is now given by 
the Taylor expansion $\mu_1 + \cdots + \mu_d$ 
of degree $d$ of the local multiplication 
$$ \mu(x,y) := \phi(\phi^{-1}(x)\phi^{-1}(y)). $$
Here each $\mu_j$ is a homogeneous polynomial of degree $j$ from 
$\g \times \g \to \g$, hence of the form 
$$ \mu_j(x) = \tilde \mu_j(x,\ldots, x), \quad \mbox{ where } \quad 
\tilde \mu_j \: (\g \times \g)^j \to \g $$
is a symmetric $j$-linear map.  By scalar extension from $\R$ to the 
$\R$-algebra $T^d(\R)$, the $j$-linear map 
$\tilde \mu_j$ induces a unique $T^d(\R)$-$j$-linear map 
$$ (T^d(\g) \times T^d(\g))^j \to T^d(\g), $$
inducing a homogeneous polynomial 
$T^d(\mu_j) \: T^d(\g) \times T^d(\g) \to T^d(\g)$,  
mapping  $T^d_\be(\g) \times T^d_\be(\g) \to T^d_\be(\g)$.  
This leads to the multiplication 
$$ T^d_\be(\mu) = \sum_{j=1}^d T^d_\be(\mu_j) \: T^d_\be(\g) 
\times T^d_\be(\g) \to T^d_\be(\g), $$
and it is clear from the construction that the order of the 
iterated multiplication maps does not exceed $d$ because 
the algebra $T^d_\be(\R)$ is nilpotent and all products of order 
$d+1$ vanish. 

(c) If $\g$ is a nilpotent Lie algebra, then the BCH multiplication 
defines on $\g$ the structure of a polynomial Lie group $(\g,*)$ 
whose inversion is given by $x^{-1} = -x$. 

(d) Let $\cA$ be any (associative) locally convex algebra, 
$M_n(\cA)$ the algebra of $(n\times n)$-matrices with entries in $\cA$, 
and 
\[ N_n(\cA) := \{(g_{ij}) \in M_n(\cA) \: (i> j \Rarrow g_{ij} = 0)\,\ 
(\forall i)\ g_{ii} = \1\} \]
be the group of strictly upper triangular matrices. Then 
$G := N_n(\cA) - \1$ is a polynomial group with respect to the multiplication 
\[ m_G(x,y) = x + y + xy \] 
(cf.~Exercise~\ref{exer:3.1.4}). 
Since $G$ can also be identified with the nilpotent 
associative subalgebra $\fn_n(\cA)$ of $M_n(\cA)$ 
consisting of strictly upper triangular matrices, 
it follows easily from (a) that all maps $m_G^{(j)}$ are polynomial of degree 
$\leq n$. 
\end{exs}

\subsection{Binomial polynomial functions} 

To deal with one-parameter groups of polynomial groups, we need 
binomial polynomial functions with values in an abelian group.

\begin{defn}
Let $A$ be an abelian group and $D = \{ n \in \Z \: n \geq n_0\}$ for some 
$n_0 \in \Z$. We call an $A$-valued function 
$f \: D \to A$  \index{polynomial!values in abelian group}
an {\it $A$-valued binomial polynomial of degree $\leq n$} if it 
is of the form 
$$ f(k) = \sum_{i = 0}^n {k \choose i} v_i \quad \mbox{ with } 
\quad v_i \in A. $$
\end{defn}

To each function $f \: D \to A$, we associate the {\it difference function} 
\[ \Delta f \: D \to A, \quad n \mapsto f(n+1) - f(n).\] 
The $\Z$-valued binomial functions 
$b_i \: \Z \to \Z, b_i(k) := {k \choose i}$, $i \in \N_0,$ satisfy 
\[  (\Delta b_i)(k) = {k+1 \choose i} - {k\choose i} = {k \choose i-1} = b_{i-1}(k) 
\quad \mbox{ for } \quad i \geq 1 \] 
and $\Delta b_0 = 0$. 
For $f$ as above, we therefore get 
\begin{equation}
  \label{eq:ndelta}
(\Delta^n f)(k) = v_n \quad \mbox{ and } \quad \Delta^{n+1} f = 0.
\end{equation}
It follows inductively that all $v_i$ are uniquely determined by~$f$. 

\begin{lem} \mlabel{lem:bipol-crit} A function $f \: D \to A$ 
is a binomial polynomial function 
of degree $\leq n$ if and only if $\Delta^{n+1} f = 0$. 
\end{lem}

\begin{prf} We have already seen above that $\Delta^{n+1}f = 0$ holds for 
each binomial polynomial function of degree $\leq n$. 

Suppose, conversely, that $\Delta^{n+1} f = 0$. 
We argue by induction that $\Delta^{n-r} f$ is a binomial polynomial function 
of degree $\leq r$. First we observe that $\Delta^n f$ is constant, 
i.e., of degree $0$. Suppose that 
$\Delta^{n-r} f$ is a binomial polynomial 
of degree $\leq r$, hence of the 
form $\sum_{i = 1}^{r+1} {k\choose i-1} v_i$. Now 
$(\Delta^{n-r-1}f)(k) - \sum_{i = 1}^{r+1}{k \choose i} v_i$ is constant, 
which implies that $\Delta^{n-r-1}f$ is a binomial polynomial of degree 
$\leq r+1$. This proves that, for each $r \in \{0,\ldots, n\}$, the function 
$\Delta^{n-r}f$ is a binomial polynomial of degree $\leq r$. 
For $r = n$, we obtain the assertion. 
\end{prf}

\subsection{Constructing one-parameter groups} 

Let $(G,m_G)$ be a polynomial Lie group. We define the power functions 
$p_k \: G \to G$, $k \in \N$, inductively by 
\[ p_0(x) \quad 
\quad \mbox{ and } \quad 
p_{k+1}(x) := m_G(x,p_k(x)) \quad \mbox{ for } \quad 
k \in \N_0. \] 
Note that the definition of a polynomial Lie group implies in particular, 
that all these functions are polynomial of degree $\leq d$ because they are 
restrictions of the functions $m_G^{(k)}$ to the diagonal in $G^k$.  

The following theorem is a variant of a result of Lazard 
(\cite[p.114]{Sr92}): 
\index{Lazard's Theorem!polynomial version}
\begin{thm}[Lazard's Theorem---polynomial version] \mlabel{thm:lazard-pol} 
For $k =1,\ldots, d$, there exist polynomial maps 
$\psi^k \: G \to G$, which are uniquely determined 
by the following requirements: 
\begin{description}
\item[\rm(i)] $\psi^1(x) = x$. 
\item[\rm(ii)] $\psi^k(x)$ is of order $\geq k$, i.e., the homogeneous 
component $\psi^k_j$ vanishes for $j < k$. 
\item[\rm(iii)] For all $k \in \N$, we have 
$p_k(x) = \sum_{i = 1}^d {k \choose i}\psi^i(x)$. 
\end{description}
Each $\psi^k$ is continuous. 
\end{thm}

\begin{prf} {\bf Uniqueness} follows from equation (iii) which 
asserts that the map $k \mapsto p_k(x)$ is a $G$-valued binomial polynomial 
function with coefficients $\psi^i(x)$. These coefficients 
are uniquely determined  by~\eqref{eq:ndelta}. 

\nin{\bf Existence:} Write $P_n(G,G)$ for the abelian group of 
continuous polynomial functions $G \to G$. 
From our discussion of binomial polynomial functions above, 
we have to show that, for each $n \in \N$, the function 
\[  f_n \: \N_0 \to P_n(G,G), \quad k \mapsto (p_k)_n \] 
is binomial polynomial of degree $\leq n$.  
We prove this by induction on $n$. For $n = 0$ we have 
\[  (p_k)_0(x) = p_k(0) = 0. \] 
Assume $n > 0$ and that 
the assertion holds for the maps $f_j$, $j \leq n-1$. 
To see that $f_n$ is a binomial polynomial of 
degree $\leq n$, it suffices to show that 
$\Delta f_n$ is a binomial polynomial of degree $\leq n-1$ 
(Lemma~\ref{lem:bipol-crit}). We derive from 
$p_{k+1}(x) = m_G(x,p_k(x)) = x + p_k(x) + \cdots$ that 
$$ (\Delta f_n)(k)(x) =  (p_{k+1})_n(x) - (p_k)_n(x) 
= \sum_{\substack{i + j \leq n \\ 1 \leq i,j}} m_{ij}(x,p_k(x))_n, $$
where $m_{ij}$ denotes the component of $m_G(x,y)$ which is homogeneous of 
degree $i$ in $x$ and of degree $j$ in $y$. In each term 
$m_{ij}(x,p_k(x))$, any contribution to the $n$th order term 
comes from the terms in $p_k(x)$ of order $< n$, which, 
according to our induction hypothesis, can be written as 
$$ (p_k)_{\leq n-1} = \sum_{i \leq n-1} {k\choose i} \psi^i_{n-1}, $$
a binomial polynomial of degree $\leq n-1$. 
Since a product of two $\Z$-valued binomial polynomials of degree $i$ and $j$ 
is a $\Z$-valued binomial polynomial of degree $\leq i + j$  
(Exercise~\ref{exer:bin-mult}), 
$\Delta f_x$ is a binomial polynomial of degree $< n$. 
This implies that $f_x$ is a binomial polynomial of degree $\leq n$, 
and this completes the proof. 
\end{prf}

\begin{ex} If $\cA$ is a nilpotent associative algebra and 
  \[ m_\cA(x,y) =  x + y + xy \]
  (Example~\ref{ex:polygrp}), 
then 
$p_k(x) = \sum_{i = 1}^k {k \choose i} x^i$
implies that $\psi^i(x) = x^i$ are the power functions of the associative algebra $\cA$. 
\end{ex}

We now follow the arguments in \cite{Ber08} 
(cf.\ also \cite[Ch.~III, \S 5]{Bou89}). 
For each $x \in G$, we obtain a 
polynomial curve 
$$ \zeta_x \: \R \to G, \quad 
\zeta_x(t) := x^t := \sum_{j = 1}^d {t \choose j} \psi^j(x). $$
For $s,t \in \N$, we then have 
\begin{eqnarray}
\label{eq:additiv}   
\zeta_x(t+s) = p_{t+s}(x) = p_t(x)p_s(x) = \zeta_x(t)\zeta_x(s).
\end{eqnarray}
This is an identity for the values of a polynomial map $\R^2 \to G$ 
on the subset $\N^2$, and since each polynomial on 
$\R^2$ is uniquely determined by its values on $\N^2$, 
(\ref{eq:additiv}) 
also holds for all $s,t \in \R$.
We conclude that $\zeta_x$ is a one-parameter group with $\zeta_x(1) = x$. 
In particular, we obtain: 
\begin{cor} \mlabel{cor:inv-pol} The 
inversion of $G$ is the continuous polynomial map given by  
\[  \eta_G(x) 
= \zeta_x(-1) = \sum_{j = 1}^d {-1 \choose j} \psi^j(x) 
= \sum_{j = 1}^d (-1)^j \psi^j(x) \quad \mbox{ for } \quad x \in G. \] 
\end{cor}

The derivative of $\zeta_x$ in $0$ is 
$$ \zeta_x'(0) = \sum_{j = 1}^d \frac{(-1)^{j+1}}{j} \psi^j(x) $$
because 
\[  \frac{d}{dt}\Big|_{t=0} {t \choose j} = \frac{(-1)(-2)\cdots (-j+1)}{j!} = \frac{(-1)^{j-1}}{j}. \]
We consider the polynomial map 
$$ \log_G \:G \to G, \quad \log_G(x) := \zeta_x'(0) 
= \sum_{j = 1}^d \frac{(-1)^{j+1}}{j} \psi^j(x). $$
To define an exponential function of the group $(G,m_G)$, 
we have to invert the polynomial map $\log_G$ (the logarithm function of $G$). 
This is achieved by  the following theorem (\cite[Thm.~PG.6]{Ber08}): 

\begin{thm} Each real polynomial Lie group $(G,m_G)$ has a   
polynomial exponential map $\exp_G \: G \to G$ 
which is bijective. Its inverse $\log_G \: G \to G$ 
is also a continuous polynomial. If $\psi^p_m$ is the homogeneous term of degree 
$m$ in $\psi^p$, 
we further have 
\begin{equation}
  \label{eq:exp-pol}
\exp_G(x) = \sum_{p=1}^d \frac{1}{p!} \psi^p_p(x) 
\quad \mbox{ and } \quad 
\log_G(x) = \sum_{j = 1}^d \frac{(-1)^{j+1}}{j} \psi^j(x).
\end{equation}
\end{thm}

\begin{prf} First we establish the functional equation for 
$\log_G$. By definition, $\log_G(x) = \zeta_x'(0)$. 
In view $\zeta_{x^n}(m) = (x^n)^m = x^{nm} = \zeta_x(nm)$ for 
$m \in \N$, we have 
$\zeta_{x^n}(t) = \zeta_x(nt)$ for all $t \in \R$, and hence 
\begin{equation}
  \label{eq:log}
\log_G(x^n) = \zeta_{x^n}'(0) = n \zeta_x'(0) = n \log_G(x), \quad 
x \in G, n \in \Z. 
\end{equation}
Hence 
\begin{eqnarray}
  \label{eq:logexp1}
\log_G(x^t) = t \log_G(x) \quad \mbox{ for} \quad x \in G, t \in \R,  
\end{eqnarray}
as both sides of \eqref{eq:log} are polynomial functions in $t$ and coincide 
on~$\N$.

We now define $\exp_G$ by \eqref{eq:exp-pol} and claim that 
$\log_G \circ \exp_G = \id_G$. We define homogeneous 
polynomials $\phi_{r,m}$ of degree $m$ in $x$  by 
\begin{eqnarray*}
x^t = \zeta_x(t) = \sum_{j = 1}^d {t \choose j}\psi^j(x)
= \sum_{r,m} t^r \phi_{r,m}(x). 
\end{eqnarray*}
For $t \in \R^\times$, we then have 
\begin{eqnarray}
  \label{eq:logexp2}
\log_G(tx) = t \log_G((tx)^{t^{-1}}) 
= t\log_G\Big(\sum_{r \leq m} t^{m-r} \phi_{r,m}(x)\Big).
\end{eqnarray}
From the definition of $x^t$ and the components $\psi^p_m$, we get in 
particular 
$$ \phi_{1,m}(x) = \sum_{p \leq m} \frac{(-1)^{p+1}}{p} \psi^p_m(x) \quad 
\mbox{and} \quad \phi_{m,m}(x) = \frac{1}{m!} \psi^m_{m}(x). $$
We also recall that $\psi^p_{m} = 0$ for $m < p$. 
Passing in (\ref{eq:logexp2}) to the derivative in $t = 0$ and using 
$\psi^1(x) = x$, we get 
$$ x
= \log_G\Big(\sum_m \phi_{m,m}(x)\Big) 
= \log_G\Big(\sum_m \frac{1}{m!} \psi^m_m(x)\Big) 
= \log_G(\exp_G(x)). $$
It follows in particular that $\log_G$ is surjective. 

Next we show that $\exp_G \circ \log_G = \id_G$. For $x \in G$ and 
$y := \log_G(x)$, we have 
$\zeta_{\exp_G(y)}'(0) = \log_G(\exp_G(y)) = y$. 
We also know that $\zeta_x$ is a polynomial one-parameter group with 
$\zeta_x'(0) = \log_G(x) = y$. Since one-parameter groups are determined by 
their derivative in $0$, it follows that 
$\zeta_x = \zeta_{\exp_G(y)}$, and hence that 
$$ x= \zeta_x(1) = \zeta_{\exp_G(y)}(1) = \exp_G(y) = \exp_G(\log_G(x)). $$

Applying $\exp_G$ to (\ref{eq:logexp1}), we find that 
$$ x^t = \exp_G(t \log_G(x)), $$
showing that the curves $t \mapsto \exp_G(ty)$ are one-parameter groups. 
Moreover, $\exp_G'(0) = \psi_1^1 = \id_G$. 
Hence $\exp_G$ is an exponential function of $G$. 
\end{prf}

\begin{prop} Each polynomial Lie group $(G,m_G)$ is nilpotent. 
\end{prop}

\begin{prf} Inductively, we define a sequence of smooth functions 
$f_n \: G^{n+1} \to G$ by 
\[  f_1(g_0, g_1) := g_0 g_1 g_0^{-1} g_1^{-1}\]
and 
$$ f_{n+1}(g_0,\ldots, g_{n+1}) := f_1(g_0, f_n(g_1,\ldots, g_{n+1})).$$
We have to show that, for some $N \in \N$, the 
function $f_N$ vanishes (cf.\ Exercise~\ref{ex:e.5}). 

As in the proof of Theorem~\ref{thm:nilpo-liealg}, we see that 
the $n$-jet of $f_n$ in $(0,\ldots, 0)$ is constant $0$, and that 
\begin{eqnarray}
  \label{eq:e.1b}
d^{(n+1)} f_n(0,0,\ldots, 0)(x_0,x_1, \ldots, x_n) 
= [x_0,[x_1,[x_2, \ldots ,[x_{n-1},x_n]]]] 
\end{eqnarray}
is the first non-trivial term in the Taylor expansion of $f_n$ in 
$(0,\ldots, 0)$. 

Each function $f_n$ can be written in the form 
$$ f_n(g_0,\ldots, g_n) = F_n(g_0,\ldots, g_n, g_0^{-1},\ldots, g_n^{-1}), $$
where $F_n$ is a restriction of some iterated multiplication map, 
hence polynomial with degree bounded by some $d \in \N$. 
Since the inversion $\eta_G$ is also polynomial of degree $\leq d$ 
(Corollary~\ref{cor:inv-pol}), 
it follows that $f_n$ is polynomial of degree $\leq d^2$. 
Combining with the information from above, we see that 
$f_n = 0$ for $n \geq d^2$ and hence $d^{(n+1)} f_n =0$. 
\end{prf}

\begin{small}
\subsection*{Exercises for Section~\ref{sec:polgrp}}  

\begin{exer}  \mlabel{exer:bin-mult} Let 
$D = \{ n \in \Z \: n \geq n_0\}$ for some $n_0 \in \Z$ and 
$f,g \: D \to R$ be functions with values in the ring~$R$. 
Show that 
\[ \Delta(fg) = \Delta(f)\tau(g) + f\Delta(g), \quad \mbox{ where } \quad 
\tau(g)(n) = g(n+1).\] 
Conclude that, if 
$f$ is polynomial of degree $\leq n$ and 
$g$ is polynomial of degree $\leq m$, then $fg$ is polynomial of degree $\leq nm$. 
\end{exer}

\begin{exer} \mlabel{exer:bin-multb} Show that, for $m,n,i,j \in \N_0$, we have 
$$ \begin{pmatrix} m \\ i \end{pmatrix} 
\begin{pmatrix} n \\ j \end{pmatrix} 
 = \sum_{i,j \leq k \leq i+j} \frac{k!}{(k-i)!(k-j)!(i+j-k)!} \
\begin{pmatrix}
m \\ k 
\end{pmatrix}. $$
Hint: Calculate $(1 + X)^k(1 + Y)^k \in \Z[X,Y]$ in two different ways.  
Use this to show that the $\Z$-span of the binomial polynomial functions 
$p_k \: \Z \to \Z, p_k(t) := \begin{pmatrix}
  t \\ k
\end{pmatrix}$
is a ring. 
\end{exer}
\end{small}

\section{Formal Lie groups} 
\mlabel{sec:formgrp} 

In the preceding section we have seen that polynomial Lie groups 
are isomorphic to nilpotent Lie groups of the form 
$(\g,*)$, where $*$ is the BCH multiplication on a nilpotent Lie algebra. 
In this section we show that this result extends from the polynomial 
context to the formal context, but this requires some extra conceptual 
background. We start with the introduction of the 
formal locally convex category whose objects are locally convex 
spaces $E$ and whose morphisms are ``formal power series'' $\phi 
= \sum_{n = 1}^\infty \phi_n \: E \to F$, 
represented by sequences $(\phi_n)_{n\in \N}$ of continuous polynomial 
maps $\phi_n$ of degree~$n$. Replacing $E$ by the space 
$E_0[[t]] \cong E^\N$ of formal power series without constant term, every 
formal map $E \to F$ implements a smooth map $E_0[[t]] \to F_0[[t]]$, 
so that the formal locally convex category can actually be realized as a 
subcategory of the category of smooth locally convex manifolds.

Group objects in this category are 
called formal Lie group. In a similar fashion as for polynomial 
group, we show that every formal Lie group $G$ has a formal exponential 
map $\Exp \:  \g \to  G$ (where $G$ and $\g$ have the same underlying 
locally convex spaces) and that this map is a formal diffeomorphism. 
This implies that every formal Lie group is formally isomorphic 
to  one for which $G = \g$ and $\Exp = \id_\g$. This in turn implies that 
the formal group structure on $\g$ is given by  the Hausdorff series 
defined by the Lie bracket on~$\g$. 
As a consequence, the category of 
formal Lie groups is equivalent to the category of locally convex Lie 
algebras.

\subsection{Formal maps} 

First we introduce the 
\index{formal locally convex category} 
{\em formal locally convex category} $\flc$. 
Its objects are locally convex spaces and a 
\index{formal morphism} 
{\em formal morphism} 
$\phi \: E  \to F$ is given by a sequence 
$(\phi_n)_{n \in \N}$ of homogeneous continuous $F$-valued polynomials 
$\phi_n \: E \to F$ of degree $n$. We think of $\phi$ as corresponding to 
the {\it formal series} \index{formal series} 
$\sum_{n=1}^\infty \phi_n$, which leads to the 
following rule for the composition: For formal morphisms 
$\phi = (\phi_n) \: E \to F$ and $\psi = (\psi_n) \: F \to G$,  we 
define  
$(\psi \circ \phi)_n$ as the term of degree $n$ in the Taylor expansion of 
the polynomial 
$$ \sum_{i \leq n} \psi_i \circ \Big(\sum_{j \leq n} \phi_j\Big). $$
In the following, we shall also write 
$$ \phi_{\leq n} := \sum_{j \leq n} \phi_j $$
for the polynomial map of degree $\leq n$ obtained by truncation from~$\phi$. 

From the Chain Rule for Taylor polynomials (Proposition~\ref{chainrtay}) 
it easily follows that 
this composition is associative, and it is clear that $(\phi_n)$ with 
$$ \phi_1 = \id_E \quad \mbox{ and } \quad \phi_n = 0 \quad \mbox{ for } \quad n > 1 $$
can be identified with $\id_E$ in $\flc$. Restricing morphisms to those 
for which only finitely many polynomials $\phi_n$ are non-zero, 
we obtain the subcategory of polynomial maps between locally convex spaces 
preserving the origin. In particular, we write $0$ for the constant map 
$E \to F$ with value zero, 
and identify this map with the corresponding morphism in $\flc$. 

\begin{rem} Although it is not covered by the categorical setup, it is clear 
that we may also admit non-zero constant terms $\phi_0$ in formal maps, 
but then the composition $(\phi_n) \circ (\psi_n)$ is only defined if $\psi_0 = 0$. 
\end{rem}

\begin{rem} \label{rem:form2} There is a natural functor from the category of 
locally convex spaces $E$ whose morphisms are germs of smooth map 
in $0 \in E$ to a locally convex space $F$ mapping $0$ to $0$.
It is simply given by 
assigning to any representing smooth map $f \: U \to F$, $U \subeq E$ a 
$0$-neighborhood, 
its Taylor series $T_0^\infty(f)$ in the origin. The functoriality of this 
assignment follows from the Chain Rule 
for Taylor series (Proposition~\ref{chainrtay}). 
\end{rem}

There is also a converse, namely a realization of $\flc$ as a subcategory 
of the category of locally convex spaces whose morphisms are smooth maps: 

\begin{rem}
For each locally convex space we write 
$E_0[[t]]$ for the set of formal $E$-valued power series of the 
form $\sum_{i =  1}^\infty a_i t^i$. Each polynomial 
map $\phi \: E \to F$ which is homogeneous of degree $d$, 
hence of the form $\phi(x) := \beta(x,x,\ldots, x)$ for some symmetric 
$d$-linear map 
$\beta \: E^d \to F$ extends to a map 
$E_0[[t]] \to F_0[[t]]$ by 
\[  \hat \phi\Big(\sum_{i = 1}^\infty a_i t^i\Big) 
:= \beta\Big(\sum_{i = 1}^\infty a_i t^i, \ldots, \sum_{i = 1}^\infty a_i t^i\Big) 
= \sum_{i_1,\ldots, i_d = 1}^\infty \beta(a_{i_1}, \ldots, a_{i_d}) 
t^{i_1 + \ldots + i_d}, \]
which is well-defined because in each degree $n$, we have only finitely many summands 
contributing to the coefficient of $t^n$. We also observe that the image of 
$\hat\phi$ is contained in $\sum_{j \geq d} F t^j$. 
The preceding observation implies that if $\phi = (\phi_n) \: E \to F$ is a
formal morphism, then we can also extend it to a  map 
$$ \hat\phi \: E_0[[t]] \to F_0[[t]] 
\quad \mbox{ by } \quad \hat\phi\Big(\sum_{i = 1}^\infty a_i t^i\Big) 
:= \sum_{n = 1}^\infty \hat\phi_n\Big(\sum_{i = 1}^\infty a_i t^i\Big) $$
because in each degree $n$, only finitely many $\hat\phi_n$ contribute to the coefficient 
of $t^n$.

We topologize $E_0[[t]]$ as the direct product space $E^\N$ and observe that with respect 
to this topology, each map $\hat\phi$ can be written as 
$\hat\phi = \sum_{n = 1}^\infty \hat\phi^n t^n$, where the maps 
$\hat\phi^n \: E_0[[t]] \to F$ are polynomial, factoring through some quotient space 
$E_0[[t]]/(t^m) \cong E^{m-1}$. We conclude in particular that each map $\hat\phi$ 
is smooth. 

We further note that for $h \in E$ we have 
$$ \hat\phi(ht) 
= \sum_{n = 1}^\infty \hat\phi_n(h t) 
= \sum_{n = 1}^\infty \phi_n(h) t^n, $$
so that we recover $\phi$ as the Taylor series of the restriction of 
$\hat\phi$ to the subspace $E\cdot t \subeq E_0[[t]]$. 

The assignment $E \mapsto E_0[[t]]$, $\phi \mapsto \hat\phi$ defines a functor to 
the smooth category. 
\end{rem}

\subsection{Formal Lie groups} 

It is obvious how to define products in $\flc$, so that we are lead to 
the concept of a formal group in this setting:  

\begin{defn} A 
\index{formal Lie group} 
{\em formal (locally convex) Lie group} is a group object in the category 
$\flc$, i.e., a locally convex space $E$, 
with two morphisms $m = (m_n) \: E \times E \to E$ 
and $\eta = (\eta_n) \: E \to E$, satisfying 
\begin{description}
\item[\rm(FG1)] $m \circ (m \times \id_E) = m \circ (\id_E \times m)$ (associativity). 
\item[\rm(FG2)] $m \circ (0, \id_E) = m \circ (\id_E, 0) = \id_E$ (neutral element).
\item[\rm(FG3)] $m \circ (\eta, \id_E) = m \circ (\id_E, \eta) = 0$ (inverse). 
\end{description}
\end{defn}

Note that each of these identities is equivalent to a sequence of identities 
between the polynomial functions $m_n \: E \times E \to E$ and 
$\eta_n \: E \to E$. 
F.i., (2) is equivalent to 
$$ m_n(x,0) = m_n(0,x) = x $$
for each $n \in \N$. This implies in particular that the linear terms is given by 
$$ m_1(x,y) = m_1(x,0) + m_1(0,y) = x + y $$
and that the quadratic term $m_2$ is a bilinear map $E \times E \to E$. 
Analyzing first order terms in (3) now leads to 
$m_1(\eta_1(x),x) = 0$, hence to 
$$ \eta_1(x) = -x. $$

Putting 
$$ [x,y] := m_2(x,y) - m_2(y,x), $$
we obtain a continuous alternating bilinear map on $E$, and 
using Hall's identity (of order 42!) (Exercise~\ref{exer:4.l.7}; \cite{Sr92}), 
it follows that this is a Lie bracket. 

Conversely, for every locally convex Lie algebra $\g$, 
the BCH series defines a formal Lie group $\Ch(\g)$  on $\g$ 
with $\eta(x) = -x$. 

Below we shall see that $\eta$ is redundant. Its existence follows from (FG1) and~(FG2). 

\begin{rem} In \cite{Ri50} Ritt defines formal groups in infinitely many variables 
as the data given by a sequence of formal power series $f^i(u,v)$ in the variables 
$(u_i)_{i \in \N}$ and $(v_i)_{i \in \N}$, where it is assumed that each 
power series $f^i$ contains only a finite number of terms of a given degree. 
This is equivalent to deal with formal 
groups on the locally convex space $E := \R^\N$. 

Ritt's axioms contain some redundancy, such as 
$$ f^i_0 = 0 \quad \mbox{ and } \quad f^i_1(u,v) = u_i + v_i. $$
\end{rem}

\begin{rem} If $(G,D_G, m_G,0)$ is a local Lie group, where 
$G$ is a $0$-neighborhood in a locally convex space $E$, then the 
functoriality of the Taylor expansion implies that 
$(E,T^\infty_{(0,0)}(m_G))$ is a formal Lie group 
(cf.\ Remark~\ref{rem:form2}).
\end{rem}

\begin{thm} [Formal Inverse Function Theorem] \mlabel{thm:inv-fct} 
A formal map $\phi = (\phi_n)_{n \in \N} \:  E \to F$ is an isomorphism 
if and only if $\phi_1 \: E \to F$ is a linear isomorphism. 
\end{thm}

\begin{prf} If $\psi \: F \to E$ is a formal map with 
$\psi \circ \phi = \id_E$ and $\phi \circ \psi = \id_F,$
then the first order terms satisfy 
$$ \psi_1 \circ \phi_1 = \id_E \quad \mbox{ and } \phi_1 \circ \psi_1 = \id_F, $$
so that $\phi_1$ is a linear isomorphism. 

Suppose, conversely, that $\phi_1$ is a linear isomorphism. Then it suffices 
to show that $\phi_1^{-1} \circ \phi \: E \to E$ is an automorphism. 
We may therefore assume that $E = F$. Then $\Hom(E,E)$ is a monoid, 
and we have to show that $\phi$ is invertible if $\phi_1 \in \GL(E)$. 

We first construct inductively a left inverse $(\psi_n)$ of $\phi$ in $\Hom(E,E)$. 
We start with $\psi_1 := \phi_1^{-1}$. Assume that $k > 1$ and that 
$\psi_j$ is already defined for $j \leq k-1$ in such a way that 
$\psi \circ \phi - \id_E$
vanishes up to terms of order $k$. 
Modulo terms of degree $k+1$, we then have 
\begin{eqnarray*}
(\psi \circ \phi)_{\leq k} 
&=& (\psi_k \circ \phi)_{\leq k} + \sum_{j = 1}^{k-1} (\psi_j \circ \phi)_{\leq k} 
= \psi_k \circ \phi_1  + \sum_{j = 1}^{k-1} (\psi_j \circ \phi)_{\leq k}\\
&=& \psi_k \circ \phi_1  + \id_E + \sum_{j = 2}^{k-1} (\psi_j \circ \phi)_{k}. 
\end{eqnarray*}
We thus put 
\[ \psi_k := - \Bigg(\sum_{j = 2}^{k-1} (\psi_j \circ \phi)_k\Bigg) \circ \phi_1^{-1} \] 
to obtain 
$(\psi \circ \phi)_{\leq k} = \id_E.$
Inductively, we get a formal map $(\psi_n)_{n\in \N}$ with 
$\psi \circ \phi = \id_E$. This proves that $\phi$ has a left inverse. 

Next we construct a right inverse $(\eta_n)$ of $\phi$. 
We start again with $\eta_1 := \phi_1^{-1}$. Assume that $k > 1$ and that 
$\eta_j$ is already defined for $j \leq k-1$ in such a way that 
$\phi \circ \eta - \id_E$
vanishes up to terms of order $k$. 
Modulo terms of degree $\geq k+1$, we then have 
\begin{eqnarray*}
(\phi \circ \eta)_{\leq k} 
&=& (\phi_1 \circ \eta)_{\leq k} + \sum_{j = 2}^{k} (\phi_j \circ \eta)_{\leq k} 
= \phi_1 \circ \eta_{\leq k} + \sum_{j = 2}^{k} 
(\phi_j \circ \eta_{\leq k-1})_{\leq k} \cr
&=& \id_E + \phi_1 \circ \eta_k + \sum_{j = 2}^{k} (\phi_j \circ \eta_{\leq k-1})_k.
\end{eqnarray*}
With 
$$ \eta_k := - \phi_1^{-1} \circ 
\Big(\sum_{j = 2}^{k} \big(\phi_j \circ \eta_{\leq k-1}\big)_k\Big).$$
we now obtain 
$(\phi \circ \eta)_{\leq k} = \id_E.$
Inductively, we thus obtain a right inverse $\eta$ of $\phi$. 
We now have 
$\eta = (\psi \circ \phi)  \circ \eta 
= \psi \circ (\phi \circ \eta) = \psi,$
and hence that $\phi$ is invertible. 
\end{prf}

\begin{thm} \label{thm:form-monoid} Each formal monoid is 
a formal group. 
\end{thm}

\begin{prf} Let $m \: E \times E \to E$ be a formal associative multiplication with 
$$ m(0,x) = m(x,0) = x,  $$
i.e., 
\begin{description}
\item[\rm(1)] $m \circ (m \times \id_E) = m \circ (\id_E \times m)$ (associativity), 
\item[\rm(2)] $m \circ (\eps_E, \id_E) = m \circ (\id_E, \eps_E) = \id_E$ (unit element), 
where $\eps_E \: E \to E$ is the map which is constant~$0$.
\end{description}
Then the same arguments as in Remark~\ref{rem:brack-taylor} show that 
$$ m_1(x,y) = x + y $$
is the first order term. We consider the map 
$$ F \: E \times E \to E \times E, \quad 
(x,y) \mapsto (x,xy). $$
Then $F_1(x,y) = (x,x+y)$ is an invertible continuous linear map, 
so that Theorem~\ref{thm:inv-fct} implies that $F$ is invertible. 
Let $G \: E \times E \to E \times E$ be the unique inverse of $F$. 
We write 
$$ G(x,y) = (G_1(x,y), G_2(x,y)). $$
Then  
$$ (x,y) = (F \circ G)(x,y) = (G_1(x,y), m(G_1(x,y),G_2(x,y))) $$
implies that $G_1(x,y) = x$, and 
$$ y = m(x,G_2(x,y)). $$
We put 
$$ \eta_r(x) := G_2(x,0) $$
and obtain 
$$ 0 = m(x,\eta_r(x)). $$

Likewise, a consideration of the map $F'(x,y) := (x,yx)$ leads to a 
morphism $\eta_l \: E \to E$ with 
$$ 0 = m(\eta_l(x),x). $$
We thus obtain 
\begin{eqnarray*}
\eta_r(x) 
&=& m(0, \eta_r(x)) 
= m(m(\eta_l(x),x),\eta_r(x))\\
&=& m(\eta_l(x),m(x,\eta_r(x)))
= m(\eta_l(x),0) = \eta_l(x). 
\end{eqnarray*}
Therefore $\eta := \eta_r = \eta_l$ satisfies 
$m \circ (\eta_E, \id_E) = m \circ (\id_E, \eta_E) = \eps_E,$
which means that $m$ defines a formal group structure on $E$. 
\end{prf}

\subsection{Formal integral curves of formal vector fields} 

Next we construct the exponential map for a formal group law 
$m$ on $E$. To this end, we have to discuss formal integral curves 
of formal vector fields. 

\index{formal vector field}
\index{formal (integral) curve} 
A {\it formal vector field on $E$} is a formal map 
$X \: E \to E$, where we do {\sl not} require that $X_0 = 0$. 

A {\it formal curve in $E$} is a morphism 
$\gamma = (\gamma_k)_{k \in \N} \: \R \to E$. 
Then $\gamma_k(t) = t^k c_k$ for some $c_k \in E$, so that we may 
identify formal curves with the corresponding sequence 
$(c_k)_{k \in E}$ in $E$ and write $\gamma_k$ for $c_k$. 
For any formal curve 
\[ \gamma(t) = \sum_{k = 1}^\infty t^k \gamma_k, \] 
we define its derivative by 
$$\gamma'(t) = \sum_{k = 0}^\infty (k+1) t^k \gamma_{k+1}. $$
A {\it formal integral curve of $X$} is a morphism 
$\gamma \: \R \to E$
with 
\begin{eqnarray}
  \label{eq:ode}
\gamma'(t) = X(\gamma(t)), \quad \mbox{ i.e.}, \quad 
\gamma' = X \circ \gamma. 
\end{eqnarray}
Note that $X \circ \gamma$ is defined because $\gamma_0 = 0$, so that 
(\ref{eq:ode}) makes sense. 

\begin{prop} Each formal vector fields $X$ has a unique formal 
integral curve $\gamma_X$.   
\end{prop}

\begin{prf} Up to terms of degree $n+1$, the equation $\gamma' = X \circ \gamma$ 
implies that 
$$ \sum_{k = 0}^n (k+1) t^k \gamma_{k+1} 
= (X\circ \gamma)_{\leq n} 
= X_0 + \sum_{k = 1}^n (X_k \circ \gamma_{\leq n})_{\leq n}. $$
Inductively, we thus obtain 
$$ \gamma_1 = X_0. $$
Further, 
$ \sum_{k = 0}^{n-1} (k+1) t^k \gamma_{k+1} 
= (X\circ \gamma)_{\leq n-1}, $
implies that 
$$ (n+1) t^n \gamma_{n+1} 
= (X\circ \gamma)_n 
= \sum_{k = 1}^n (X_k \circ \gamma_{\leq n})_n, $$
so that 
\begin{eqnarray}
  \label{eq:solut}
\gamma_{n+1} 
= \frac{1}{n+1}\sum_{k = 1}^n (X_k \circ \gamma_{\leq n})_n(1) 
= \frac{1}{n+1} (X \circ \gamma_{\leq n})_n(1). 
\end{eqnarray}
From this formula, it is clear that there exists at most one 
formal integral curve for each formal vector field. 

Conversely, its existence follows from the inductive definition: 
$\gamma_1 := X_0$ and 
$$ \gamma_{n+1} := \frac{1}{n+1}\sum_{k = 1}^n (X_k \circ \gamma_{\leq n})_n(1), \quad 
n \in \N. $$ 
We then have 
\begin{eqnarray*}
\sum_{k = 0}^n (k+1) t^k \gamma_{k+1} 
&=& X_0 + \sum_{k = 1}^n \sum_{j = 1}^k (X_j \circ \gamma_{\leq k})_k(t) \\
&=& X_0 + \sum_{k = 1}^n (X \circ \gamma_{\leq k})_k(t) 
= X_0 + \sum_{k = 1}^n (X \circ \gamma)_k(t) \\
&=& (X \circ \gamma)_{\leq n}(t). 
\end{eqnarray*}
Since this holds for each $n \in \N$, we have 
$\gamma' = X \circ \gamma$. 
\end{prf}

\begin{lem} \mlabel{lem:intcurve-trafo} Let $\phi \: E \to F$ be a formal diffeomorphism, 
$X$ a formal vector field on $E$, and 
$$ \tilde X := T(\phi) \circ X \circ \phi^{-1} $$
the corresponding formal vector field on $F$. Then 
$\phi \circ \gamma_X = \gamma_{\tilde X}.$
\end{lem}

\begin{prf} The formal curve $\tilde \gamma_X := \phi \circ \gamma_X$ satisfies 
$$ \tilde\gamma_X'
= T\phi  \circ \gamma'
= T\phi \circ X \circ \gamma
= \tilde X \circ \phi \circ \gamma = \tilde X \circ \tilde\gamma_X. $$
Hence $\tilde\gamma_X$ is an integral curve of $\tilde X$, and the assertion follows from 
the uniqueness of formal integral curves. 
\end{prf}

\begin{rem} (Formal Picard Iteration) Let $X$ be a formal vector field on $E$. 
We then define a sequence of polynomials by Picard iteration: 
$$ \gamma^0 := 0, \quad 
\gamma^{n+1}(t) := \int_0^t (X \circ \gamma^n)_{\leq n}(s)\, ds. $$
This leads to 
$$ \gamma^1(t) = \int_0^t X(0)\, ds = t X_0 $$
and 
$$ \gamma^2(t) 
= \int_0^t X_{\leq 1}(s X_0)\, ds 
= \int_0^t X_0 + s X_1(X_0)\, ds 
= t X_0 + \frac{t^2}{2} X_1(X_0). $$
Moreover, 
\begin{eqnarray*}
\gamma^3(t) 
&=& \int_0^t (X \circ \gamma^2)_{\leq 2}(s)\, ds  \\
&=& \int_0^t X_0 + s X_1(X_0) + \frac{s^2}{2} X_1^2(X_0) + s^2 X_2(X_0)\, ds \\
&=& t X_0 + \frac{t^2}{2} X_1(X_0) 
+ \frac{t^3}{3!} X_1^2(X_0) 
+ \frac{t^3}{3} X_2(X_0). 
\end{eqnarray*}
By definition, we then have 
$$ (\gamma^{n+1})' = (X \circ \gamma^n)_{\leq n}. $$

On the other hand, the unique integral curve $\gamma_X$ of $X$ satisfies 
$$ \gamma_X' = X \circ \gamma_X, $$
which leads to 
$$ (\gamma_X')_{\leq n} 
= (X \circ \gamma_X)_{\leq n}. $$
Clearly, $\gamma^1 = (\gamma_X)_1$. So assume that 
$\gamma^n = (\gamma_X)_{\leq n}$. Then 
$$ (\gamma_X')_{\leq n} 
= (X \circ (\gamma_X)_{\leq n})_{\leq n} 
= (X \circ \gamma^n)_{\leq n}
= (\gamma^{n+1})', $$
and therefore 
$$ (\gamma_X)_{\leq n+1} = \gamma^{n+1}. $$

We conclude that $\gamma_X$ can be obtained by Picard iteration as above.  
\end{rem}

From (\ref{eq:solut}), it follows that the dependence of the integral curve 
$\gamma_X$ from the vector field $X$ is such that, for each 
$n \in \N$, the element $(\gamma_X)_n$, depends polynomially on $X$, 
where the degree of the polynomial is bounded by $n$. F.i., $\gamma_1 = X_0$ 
depends linearly on $X$, and 
$$ \gamma_2 
= \frac{1}{2} (X_1 \circ \gamma_1)(1) = \frac{1}{2} X_1(X_0) $$
is quadratic in the components of $X$. Further 
$$ \gamma_3 = \frac{1}{3}
\Big(\frac{1}{2} X_1^2(X_0) + X_2(X_0)\Big). $$

\begin{lem} On the formal level, we have 
$$ \gamma_{sX}(t) = \gamma_X(st), \quad \mbox{ for } \quad s,t \in \R. $$
In particular, $(\gamma_X)_n$ is homogeneous 
of degree $n$ in the components of $X$. 
\end{lem}

\begin{prf} The curve 
$\eta(t) := \gamma_X(st)$ satisfies 
$$ \eta'(t) = s\gamma_X'(st) = s X(\gamma_X(st)) = (sX)(\eta(t)), $$
so that the uniqueness of integral curves implies that 
$\eta = \gamma_{sX}$, i.e., 
$$ \gamma_{sX}(t) = \gamma_X(st). $$
This relation implies imediately that $(\gamma_X)_n$ is homogeneous 
of order $n$ in~$X$. 
\end{prf}

Now we return to our original problem of constructing the 
exponential function of a formal Lie group $(E, m, \eta)$. 
First we associate to each $x \in E$ a formal vector field by 
\[  x_l(y) := \derat0 m(y,tx),  \] 
which means that we first consider the formal map 
$$ \R \times E \times E \to E, \quad 
(t,x,y) \mapsto m(y,tx) $$
and then pick out those terms which are linear in $t$. 
From $m_1(y,tx) = y + tx$ it follows in particular, that $m_1$ 
contributes the linear term $x$ to $x_l(y)$. More generally, 
we obtain from $m_k(y,tx)$ a term of order $1$ in $x$ and order~$k-1$ in~$y$. Hence 
\[  (x_l)_k = \derat0 m_{k+1}(y,tx). \] 

Let $\gamma_x \: \R \to E$ denote the unique formal integral curve 
of the formal vector field $x_l$. This vector field depends linearly 
on $x$, so that $(\gamma_x)_n$ is a homogeneous polynomial of order $n$ 
in $x$. We thus obtain a formal map 
$$ \Exp \: E \to E, \quad \Exp_n(x) := (\gamma_x)_n, $$
called the 
\index{formal group!exponential function} 
{\it exponential function of the formal group $(E,m_E)$}. 
Since 
$$ \Exp_1(x) = (\gamma_x)_1 = x, $$
we have $\Exp_1= \id_E$, and the Formal Inverse Function Theorem implies that 
$\Exp \: E \to E$ is a formal diffeomorphism of $E$. 

We now consider the new formal group structure 
$$ \tilde m := \Exp^{-1} \circ m \circ (\Exp \times \Exp). $$
We want to show that this new structure coincides with $\Ch(\g)$, 
i.e., $\tilde m$ is given by the BCH series. 

Let 
$$ \tilde x_l(y) := \derat0 \tilde m(y,tx) = T_{(y,0)}(\tilde m)(0,x). $$
Then we have 
\begin{eqnarray*}
\tilde x_l(y) &=& T(\Exp)^{-1} T(m) T_{(y,0)}(\Exp \times \Exp)(0,x) \\
&=& T(\Exp)^{-1} T_{(\Exp y,0)}(m) (0,T_0(\Exp)x) \\
&=& T(\Exp)^{-1} T_{(\Exp y,0)}(m) (0,x) 
=  \big(T(\Exp)^{-1} \circ x_l \circ \Exp\big)(y), 
\end{eqnarray*}
so that Lemma~\ref{lem:intcurve-trafo} implies that 
$$ \tilde\gamma_x := \Exp^{-1} \circ \gamma_x $$
is the unique formal integral curve of $\tilde x_l$. Therefore 
$$ \gamma_x(t) = \gamma_{tx}(1) = \Exp(tx) $$
implies that $\tilde\gamma_x(t) = tx.$
Therefore the transformation by $\Exp$ transforms the formal group law 
$m$ into the formal group law $\tilde m$ which has the addition property that the 
corresponding exponential  
$$ \tilde \Exp = \id_E $$
is the identity map. 

\index{formal group!canonical} 
In the following we call a formal Lie group $(E, m)$ {\it canonical} if 
the corresponding exponential map is $\id_E$. We have thus show that: 

\begin{prop} Every formal Lie group $(E,m)$ is isomorphic to a canonical one, 
i.e., to one with $\exp_E = \id_E$. 
\end{prop}

Let $\g = (E,[\cdot,\cdot])$ denote the Lie algebra of the formal Lie group 
$(E,m)$. Clearly $\Ch(\g)$ is a canonical Lie group with the Lie algebra $\g$ 
(cf.~Theorem~\ref{thm:BCH-free}), and we want  to show that this is the only one. 

\begin{lem} If the formal Lie group $(E,m)$ is canonical, then the same holds for 
its tangent bundle $(T(E), T(m))$. 
\end{lem}

\begin{prf} If $(F,m_F)$ is a formal Lie group and $\Exp_F$ its exponential function, 
then $T(\Exp_F)$ is the exponential function of the formal Lie group $T(F)$. In fact, 
for each $x \in F$, we have the relation 
$$ \frac{d}{dt} \Exp_F(tx) = \deras0 \Exp_F(tx) * sx. $$
Applying the tangent functor, this leads for 
$(x,v) \in T(F) \cong F \times F$ to the relation 
$$ \frac{d}{dt} T(\Exp_F)(tx,tv) = \deras0 T(\Exp_F)(tx,tv) * (sx,sv), $$
which means that $t \mapsto T(\Exp_F)(tx,tv)$ is the unique formal integral curve of the 
vector field defined by 
$$(x,v)_l(y,w) := \deras0 (y,w) * (sx,sv). $$

The preceding argument shows in particular that $\Exp_F = \id_F$ implies 
$\Exp_{T(F)} = T(\Exp_F) = T(\id_F) = \id_{T(F)}$, so that $(T(F), T(m_F))$ is canonical if 
$(F,m_F)$ has this property. 
\end{prf}

Any formal Lie group $(E,m_E)$ is uniquely determined by the maps 
$(m_E)_n$, which can also be obtained from the iterated tangent maps 
$$T^n(m_E) \: T^n(E) \times T^n(E) \to T^n(E),$$ restricted to the 
fiber $T^n(E)_0$ over $0$, which is a nilpotent group, endowed with a polynomial 
multiplication. 

To show that any canonical formal Lie group $(E,m_E)$ is of the form 
$\Ch(\g)$, it therefore suffices to verify this under the additional assumption that 
it is nilpotent and that $m_E$ is polynomial. Then the formal group structure is an actual 
group structure, and we have to show that for any nilpotent polynomial Lie group 
structure on a locally convex space $E$, the assumption that it is canonical as a formal 
group implies that $m_E$ is the BCH multiplication. 

In view of Theorem~\ref{bch-locexp}, it suffices to show that 
\[  t x * s x = (t+s) x \quad \mbox{ for } \quad t,s \in \R, x \in E.\] 
This condition is equivalent to the curve 
$\eta_x(t) := tx$ being an integral curve of the left invariant vector field $x_l$. 
This does {\bf not} follow immediately from the assumption that $E$ is canonical as a formal  
group, which only implies that $\eta_x$ represents the Taylor series of any integral 
curve of the vector field $x_l$ in $0$. 
We shall close this gap in the following section, using a 
method due to Lazard that we have already seen in the polynomial 
context in Section~\ref{sec:polgrp}. 

\subsection{One-parameter groups} 

Let $(E,m)$ be a formal Lie group. We define the power functions 
$p_k$, $k \in \N$, inductively by 
$$ p_1(x) = x \quad 
\quad \mbox{ and } \quad 
p_{k+1}(x) := m(x,p_k(x)). $$

The following theorem is another variant of a result of Lazard 
(cf.\ \cite{Sr92}, p.114) that we have already seen in a polynomial 
version in Theorem~\ref{thm:lazard-pol}: 
\index{Lazard's Theorem!formal version}
\begin{thm}[Lazard's Theorem---formal version] There exist formal maps 
$\psi^k \: E \to E$, $k \in \N$, which are uniquely determined by the following 
requirements: 
\begin{description}
\item[\rm(1)] $\psi^1(x) = x$. 
\item[\rm(2)] $\psi^k(x)$ is of order $\geq k$, i.e., $\psi^k_j = 0$ for $j < k$. 
\item[\rm(3)] For all $k \in \Z$ we have 
$p_k(x) = \sum_{i = 1}^\infty {k\choose i}  \psi^i(x).$ 
\end{description}
\end{thm}

\begin{prf} {\bf Uniqueness:} From (3) we derive for the terms of 
order $n$ with (2): 
$$ (p_k)_n = \sum_{i \leq n} {k \choose i} \psi^i_n, $$
which is an $E$-valued binomial polynomial function, and the corresponding 
coefficients $\psi^i_n$ are uniquely determined and hence homogeneous 
of degree $n$. 

\nin{\bf Existence:} From our discussion of binomial polynomial functions above, 
we have to show that for each $n \in \N$, the functions 
$$ f_x \: \Z \to E, \quad k \mapsto (p_k)_n(x) $$
are binomial polynomials of degree $\leq n$. 

We prove this by induction on $n$. For $n = 1$ we have 
$$ (p_k)_1(x) = kx, $$
as an easy consequence of $m_1(x,y) = x+ y$. 
Assume $n > 1$ and that 
the assertion holds for the maps $(p_k)_j$, $j \leq n-1$. 
To see that $f_x(n) := (p_k)_n(x)$ is a binomial polynomial of 
degree $\leq n$, it suffices to show that 
$\Delta f_x$ is a binomial polynomial of degree $\leq n-1$ 
(Lemma~\ref{lem:bipol-crit}). We derive from 
$p_{k+1}(x) = m(x,p_k(x)) = x + p_k(x) + ...$ that 
$$ \Delta f_x(k) =  (p_{k+1})_n(x) - (p_k)_n(x) 
= \sum_{i + j \leq n} m_{ij}(x,p_k(x))_n, $$
where $m_{ij}$ denotes the component of $m(x,y)$ which is homogeneous of 
degree $i$ in $x$ and of degree $j$ in $y$. In each term 
$m_{ij}(x,p_k(x))$, any contribution to the $n$-th order term 
comes from the terms in $p_k(x)$ of order $< n$, which can be written as 
\[  (p_k)_{\leq n-1} = \sum_{i \leq k \leq n-1} {k \choose i} \psi^i_k, \] 
a binomial polynomial of degree $\leq n-1$. 
Since a product of two $\Z$-valued binomial polynomials of degree $i$ and $j$ 
is a $\Z$-valued binomial polynomial of degree $\leq i + j$  
(cf.\ Exercise~\ref{exer:bin-mult}), 
it follows that $\Delta f_x$ is a binomial polynomial of degree $< n$. 
This implies that $f_x$ is a binomial polynomial of degree $\leq n$, 
and this completes the proof. 
\end{prf}

\begin{ex} If $\cA$ is an associative algebra and 
$$ m(x,y) = x + y + xy, $$
then 
$p_k(x) = \sum_{i = 1}^k {k \choose i} x^i$
implies that $\psi_i(x) = x^i$ are the power functions of the associative algebra $\cA$. 
\end{ex}

As in Section~\ref{sec:polgrp}, we obtain for each $x \in E$ a formal curve 
$$ \zeta_x \: \R \to E, \quad 
\zeta_x(t) := \sum_{j = 1}^\infty {t \choose j}\psi^j(x), $$
and note that this is a well-defined formal map, because we have 
only finitely many summands in each degree. 

For $s,t \in \N$, we now have 
$$ \zeta_x(t+s) = p_{t+s}(x) = p_t(x)p_s(x) = \zeta_x(t)\zeta_x(s), $$
as formal functions $E \to E$. 
For each fixed degree of homogeneity in $x$, this is an identity for the values of 
a polynomial map $\R^2 \to E$ on the subset $\N^2$, and since each polynomial on 
$\R^2$ is uniquely determined by its values on $\N^2$, 
$$ \zeta_x(t+s) = p_{t+s}(x) = p_t(x)p_s(x) = \zeta_x(t)\zeta_x(s) $$
also holds for all $s,t \in \R$.
We conclude that $\zeta_x$ is a formal one-parameter group. Its derivative in $0$ is given by 
$$ \zeta_x'(0) = \sum_{j = 1}^\infty \frac{(-1)^{j+1}}{j} \psi^j(x) $$
because 
\[  \frac{d}{dt}\Big|_{t = 0} {t \choose j} = \frac{(-1)(-2)\cdots (-j+1)}{j!} = \frac{(-1)^{j-1}}{j}. \]
We consider the formal map 
$$ L \:E \to E, \quad L(x) := \sum_{j = 1}^\infty \frac{(-1)^{j+1}}{j} \psi^j(x). $$
Then $L_1(x) = \psi^1(x) = x$ implies that $L$ is a formal diffeomorphism. 

One the other hand, the fact that $\zeta_x$ is a one-parameter group implies that 
$\zeta_x$ is the unique integral curve of the corresponding vector field 
$L(x)_l$, which means that 
$$ \Exp(t L(x)) = \gamma_{L(x)}(t) = \zeta_x(t), $$
and for $t = 1$ we see that 
$\Exp(L(x)) = x.$
We thus obtain 
$\Exp = L^{-1}$
and see that all curves 
$t \mapsto \Exp(tx)$
are formal $1$-parameter groups.

With Theorem~\ref{bch-locexp}, applied to the exponential 
Lie group structure on $E_0[[t]]$ we now obtain: :

\begin{thm} For any canonical formal Lie group $G$, the multiplication 
$m = (m_n)_{n \in \N}$ is given by the Hausdorff series, i.e., 
$(G,m) = \Ch(\g)$. 
\end{thm}

\section{Some results on Lie group extensions} 
\mlabel{sec:8.3}

In this section we collect some miscellaneous results on extensions 
of Lie groups. 
In Subsection~\ref{sec:8.4b} we start with a criterion 
for a Lie group extension $\hat G$ of $G$ by $N$ to have an exponential 
function: it suffices that $N$ is regular and that $G$ has an exponential function. Following this philosophy on extending nice properties 
from Lie groups to extensions by regular Lie groups, we develop 
in Subsection~\ref{subsec:relfundtheo} 
a tool that can be used to integrate Lie algebra-valued $1$-forms 
satisfying the Maurer--Cartan equation, even if the corresponding 
Lie group is not known to be regular. 
The relative version of the Fundamental Theorem (Theorem~\ref{rel-funda}) 
ensures integrability for Lie group extensions $\hat G$ of $G$ by $N$ 
if we have integrability in $G$, and $N$ is a regular Lie group. 

In Subsection~\ref{sec7.6} we study abelian extensions of a 
Lie group $G$ by a smooth $G$-module $\fa$ which is assumed to be 
Mackey complete. Our main results concern the existence of smooth global 
group cocycles for given Lie algebra cocycles. In particular, we show that 
they always exist if $G$ is smoothly contractible. 

In Subsection~\ref{sec:11.3} we address the problem 
  that many central extensions of locally exponential Lie groups
  may not even have smoothly local sections. Typical examples arise from 
  $G$, considered as a central extension of the adjoint group if 
the closed subspace  $\fz(\g) \subeq \g$ is not split. 
To deal with such 
  central extensions we introduce generalized central extensions and 
develop tools based on embeddings of such extensions 
into central extensions defined by cocycles.

\subsection{The exponential function of a Lie group extension} 
\mlabel{sec:8.4b}

For the following theorem, we recall the concept of a Lie group 
extension from Definition~\ref{def:liegrp-ext}. 

\begin{thm} \mlabel{thm:exp-ext} 
Let $q \: \hat G \to G$ be a Lie group extension of $G$ by 
the regular Lie group $N = \ker q$. 
If $G$ has a smooth exponential function, then 
also $\wh{G}$ has a smooth exponential 
function.
\end{thm}

The assumption of the preceding theorem is in particular satisfied 
if $N =(\fz ,+)$ for a Mackey complete locally convex space~$\fz$. 

\begin{prf}
Let $\sigma\colon U\to \wh{G}$ be a smooth local section for~$q$,
defined on an open identity neighborhood $U\subeq G$, such that $\sigma({\bf e})={\bf e}$.
Let $V\subeq \fg:=\L(G)$ be a balanced open $0$-neighborhood such that
\[
\exp_G(V)\subeq U.
\]
Then $W:=\L(q)^{-1}(V)$ is a balanced open $0$-neighborhood in $\wh{\fg}:=\L(\wh{G})$.
If $w\in W$, then
\[
\oline\gamma_w\colon [0,1]\to U,\quad \oline\gamma_w(t):=\exp_G(t\L(q)(w))
\]
is a smooth curve in~$U$ such that $\gamma_w(0)={\bf e}$ and $\delta^\ell(\gamma_w)=\L(q)(w)$.
Hence
\[
\eta_w:= \sigma\circ \oline\gamma_w
\]
is a smooth curve in $\wh{G}$
such that $\eta_w(0)={\bf e}$.
Since $q\circ\eta_w=\oline\gamma_w$, we have
\[
\L(q)\circ\delta^\ell(\eta_w)=\delta^\ell(q\circ \eta_w)=\delta^\ell(\oline\gamma_w)=\L(q)(w),
\]
whence
\begin{equation}\label{inn}
\delta^\ell(\eta_w)_t -w\in\ker \L(q)=\L(N)=: \fn 
\quad \mbox{ for each }\quad t \in [0,1]. 
\end{equation}
We are looking for a smooth curve
$\theta_w\colon [0,1]\to N$ such that 
\begin{equation}\label{need}
\theta_w(0)=\be \quad \mbox{ and } \quad 
\delta^\ell(\theta_w\eta_w)=w.
\end{equation}
Since
\[
\delta^\ell(\theta_w\eta_w)=\delta^\ell(\eta_w)+\Ad(\eta_w^{-1})\delta^\ell(\theta_w),
\]
the identity (\ref{need}) is equivalent to
$\delta^\ell(\eta_w)+\Ad(\eta_w^{-1})\delta^\ell(\theta_w)=w.$
Equivalently,
\[
\delta^\ell(\theta_w)=\Ad(\eta_w)(w-\delta^\ell(\eta_w)),
\]
where the right hand side is a smooth curve in $\fn$ as a consequence of (\ref{inn}).
The preceding equation has the unique solution
\[
\theta_w:=\Evol_N\big(\Ad(\eta_w)(w-\delta^\ell(\eta_w))\big).
\]
Then $\zeta_w:=\theta_w\eta_w$ is a smooth curve $[0,1]\to\wh{G}$
taking $0$ to ${\bf e}$ such that $\delta^\ell(\zeta_w)=w$. 
In view of Exercise~\ref{exer:local-one-para}, 
$\zeta_w$ admits a unique extension to a smooth homomorphism
$\gamma_w\colon \R\to \wh{G}$ which is uniquely determined by 
the initial condition $\gamma_w'(0)= w$. 
For each $s\in\R\setminus\{0\}$, the curve $t\mapsto \gamma_w(st)$ 
is the unique smooth one-parameter group $\gamma_{sw}$
with derivative $sw$ at~$0$.
We therefore have $\gamma_w\colon \R\to G$ at our disposal
for all $w\in \wh{\fg}$.
We define
\[
\exp\colon \wh{\fg}\to \wh{G},\quad \exp(w):=\gamma_w(1).
\]
If we can show that $\exp|_W$ is smooth, then $\exp$ will be smooth by 
Lemma~\ref{lem:5.1.2}. Note first that the map
\[
f\colon W\to C^\infty([0,1],\wh{G}),\quad w\mapsto \eta_w
\]
is smooth because
\[
\wh{f}\colon W\times [0,1]\to \wh{G},\quad \wh{f}(w,t):=\eta_w(t)
=\sigma\big(\exp_G(t\L(q)w)\big)
\]
is smooth (Proposition~\ref{prop:cartes-closed-lie-a}). 
Since $\delta^\ell\colon C^\infty([0,1],\wh{G})\to C^\infty([0,1],\wh{\fg})$ is smooth 
(Proposition~\ref{prop:smooth-logder0}), 
we deduce that $\delta^\ell\circ f$ is smooth and hence also the map
\[
h:=\wh{\delta^\ell\circ f}\colon W\times [0,1]\to \wh{\fg},\quad h(w,t):=\delta^\ell(\eta_w)(t).
\]
We deduce that
\[
\phi\colon W\to C^\infty([0,1],\fn),\quad \phi(w):=\Ad(\eta_w)(w-\delta^\ell(\eta_w))
\]
is smooth because
\[
\wh{\phi}\colon W\times [0,1]\to \fn,\quad \wh{\phi}(w,t):=\wh{\Ad}(\wh{f}(w,t),w-h(w,t))
\]
is smooth.
As a consequence, also
\[
\psi\colon W\to C^\infty([0,1],N),\quad w\mapsto \theta_w=\Evol_N(\phi(w))
\]
is smooth.
Now, for $w\in W$, we have
\[
\exp(w)=\gamma_w(1)=\zeta_w(1)=\theta_w(1)\eta_w(1)=\psi(w)(1)f(w)(1)=\wh{\psi}(w,1)\wh{f}(w,1),
\]
which is a smooth function of $w\in W$.
\end{prf}

\begin{cor} \mlabel{cor:8.4.2} If, in {\rm Theorem~\ref{thm:exp-ext}}, 
$G$ is locally exponential and $N$ is central in $\hat G$, 
then $\hat G$ is also locally exponential. 
\end{cor}

\begin{prf} Let $U \subeq \g$ be an open convex $0$-neighborhood 
for which $(\exp\res_U)^{-1}$ is a chart of $G$. Then 
$\sigma(\exp_G x) := \exp_{\hat G} x$ defines a smooth section 
$\exp_G(U) \to \hat G$ of the quotient map and we thus obtain the inverse of a 
chart of $\hat G$ by 
\[ \fz \times U \to \hat G, \quad 
(z,x) \mapsto \exp_Z(z) \exp_{\hat G} x 
= \exp_{\hat G} z \exp_{\hat G} x 
= \exp_{\hat G}(z + x)\] 
(Lemma~\ref{lem:4.1.4}). 
This implies that $\hat G$ is locally exponential. 
\end{prf}

\subsection{A relative Fundamental Theorem} 
\mlabel{subsec:relfundtheo}

Consider a Lie group extension 
$q \: \hat G \to G$ of $G$ by $N$. If $\hat G$ is regular and 
$\alpha \in \Omega^1(M,\hat\g)$ satisfies the Maurer--Cartan 
equation, then it is locally integrable by the 
Fundamental Theorem~\ref{thm-fundamental}. 
In situations where $\hat G$ is not known to be regular, 
one would like to have tools 
to bypass it in concrete situations. 
The relative version of the Fundamental Theorem presented 
in this subsection is such a tool. It asserts that,
if $N$ is regular, then the local integrability of 
$\alpha$ follows from the local integrability of 
$\L(q)\alpha \in \Omega^1(M,\g)$. 

If $M = \hat\g$ and 
$\kappa_{\hat\g}(x) := \int_0^1 e^{-t\ad x}\, dt$ 
exists and defines a smooth $1$-form in $\Omega^1(\hat\g,\hat\g)$, 
then the local integrability of $\kappa_{\hat\g}$ 
is equivalent to the existence of a smooth exponential 
function of $\hat G$ and the integrability of $\L(q)\kappa_{\hat\g} 
= \kappa_\g \circ \L(q)$ follows from the existence of a smooth exponential 
function of $G$. In this sense the local Fundamental Theorem 
would imply Theorem~\ref{thm:exp-ext} if we already had $\kappa_{\hat\g}$, 
which we do not want to assume a priori.

\begin{defn} Let $M$ be a smooth manifold 
and $G$ a Lie group with Lie algebra $\g$. 
We write 
$$ Z^1_{\rm dR}(M,\g) := \big\{ \alpha \in \Omega^1(M,\g) \: \dd \alpha 
+ \shalf [\alpha,\alpha] = 0\big\} $$
for the set of solutions of the Maurer--Cartan equation. 
If $\g$ is abelian, then $Z^1_{\rm dR}(M,\g)$ is the space of closed 
$\g$-valued $1$-forms. 
\end{defn}

If $G$ is regular, $M$ is connected and $m_0 \in M$, 
then we obtain from  the Fundamental Theorem~\ref{thm-fundamental} a sequence of maps 
\[  G \into C^\infty(M,G) \sssmapright{\delta} Z^1_{\rm dR}(M,\g) 
\smapright{\per} \Hom(\pi_1(M,m_0), G), \]
which is exact as a sequence of pointed sets, i.e., at each point, 
the inverse image of the base point is the range of the preceding map.

Let $\g$ be a locally convex Lie algebra and $M$ a 
smooth manifold. Further, let $\alpha \in \Omega^1(M,\g)$ be a smooth $\g$-valued $1$-form and 
$\gamma \: M \to \Aut(\g)$ be a map which is smooth in the sense that 
$$ \tilde \gamma\: M \times \g 
\to \g \times \g, \quad (m,x) \mapsto (\gamma(m)x, \gamma(m)^{-1}x) $$
is smooth (cf.\ Section~\ref{sec:e.2}). 
If $\gamma$ is smooth, then 
\[ \delta(\gamma)_m := \gamma(m)^{-1} \circ T_m(\gamma) \: T_m(M) \to \der(\g) \] 
is such that the corresponding map $TM \times \g \to \g$ is smooth 
(Proposition~\ref{prop:aut-der}). 
We write $\gamma.\alpha \in \Omega^1(M,\g)$ for the smooth $1$-form defined pointwise by 
$(\gamma.\alpha)_m = \gamma(m) \circ \alpha_m$ for $m \in M$. 
For any smooth vector field 
$X \in {\cal V}(M)$, we have 
$$ (X\gamma)(m) := T_m(\gamma) X(m) 
= \gamma(m)\big(\delta(\gamma)_m X(m)\big), $$
i.e., 
$$ X\gamma = \la \gamma.\delta(\gamma), X \ra. $$ 
We now obtain for vector fields $X,Y \in {\cal V}(M)$: 
\begin{eqnarray*}
&& \dd(\gamma.\alpha)(X,Y) 
= X\big((\gamma.\alpha)(Y)\big) - Y\big((\gamma.\alpha)(X)\big) 
- (\gamma.\alpha)([X,Y]) \\ 
&=& (X\gamma).\alpha(Y) + \gamma.(X\alpha(Y)) 
- (Y\gamma).\alpha(X) + \gamma.(Y\alpha(X)) 
- (\gamma.\alpha)([X,Y]) \\ 
&=& (\gamma.\dd \alpha)(X,Y) + (X\gamma).\alpha(Y) - (Y\gamma).\alpha(X) \\
&=& (\gamma.\dd \alpha)(X,Y) 
+ \big(\gamma.\delta(\gamma)\big)(X)(\alpha(Y)) 
- \big(\gamma.\delta(\gamma)\big)(Y)(\alpha(X))  \\
&=& \gamma.\Big(\dd \alpha(X,Y) + \delta(\gamma)(X)(\alpha(Y))- \delta(\gamma)(Y)(\alpha(X))\Big).
\end{eqnarray*}
With 
$$ (\delta(\gamma)\wedge \alpha)_m(v,w) := 
\delta(\gamma)_m(v)(\alpha_m(w)) -\delta(\gamma)_m(w)(\alpha_m(v)), $$
this can be written as 
\begin{eqnarray}
  \label{eq:diff-x}
\dd(\gamma.\alpha) = \gamma.(\dd \alpha + \delta(\gamma) \wedge \alpha). 
\end{eqnarray}
From this formula we immediately derive 
\begin{eqnarray}
  \label{eq:diff-y}
&& \dd(\gamma.\alpha) +\shalf[\gamma.\alpha, \gamma.\alpha] 
= \gamma.(\dd \alpha + \delta(\gamma) \wedge \alpha) + \shalf \gamma.[\alpha,\alpha]\notag\\
&=& \gamma.\big(\dd \alpha +\shalf[\alpha,\alpha] 
+ \delta(\gamma) \wedge \alpha\big). 
\end{eqnarray}

\begin{lem}\mlabel{lem.e.1} Let $G$ be a Lie group with Lie algebra $\g$.  
For $f \in C^\infty(M,G)$ and $\alpha \in \Omega^1(M,\g)$, we have 
$$ \dd(\Ad(f).\alpha) = \Ad(f).\big(\dd \alpha + [\delta(f), \alpha]\big).$$  
\end{lem}

\begin{prf} We define a map $\gamma \: M \to \Aut(\g)$ by 
$\gamma := \Ad \circ f$. Then $\delta(\gamma) = \ad \circ \delta(f)$ and 
formula~\eqref{eq:diff-x} lead to 
\begin{align*}
\dd(\Ad(f).\alpha) 
&= \Ad(f).(\dd \alpha + \delta(\Ad(f)) \wedge \alpha) 
= \Ad(f).(\dd \alpha + \ad(\delta(f))\wedge \alpha) \\
&= \Ad(f).(\dd \alpha + [\delta(f),\alpha]). 
\qedhere
\end{align*}
\end{prf}

\begin{lem}\mlabel{lem.e.2b} 
Let $M$ be a connected smooth manifold and $G$ be a Lie group 
with Lie algebra~$\g$. 
Then the set $Z^1_{\rm dR}(M,\g)$ of solutions of the Maurer--Cartan equation is invariant 
under the smooth affine right action 
\[ \alpha * f := \delta(f) + \Ad(f)^{-1}.\alpha\] 
of $C^\infty(M,G)$ on $\Omega^1(M,\g)$. 
\end{lem}

\begin{prf} We recall from Lemma~\ref{lem.e.2} that $*$ defines a smooth 
affine action. To verify the invariance of $Z^1_{\rm dR}(M,\g)$, 
we put $R(\alpha) := \dd \alpha + \frac{1}{2}[\alpha,\alpha]$ and recall that 
$\delta(f)$ satisfies the Maurer--Cartan equation. Then 
Lemma~\ref{lem.e.1} implies that 
\begin{eqnarray*}
&& R(\alpha * f) = \dd(\alpha * f) + \shalf [\alpha * f, \alpha * f] \\
&=& \dd(\delta(f)) 
+ \Ad(f)^{-1}.\big(\dd \alpha + [\delta(f^{-1}), \alpha]\big) 
+ \shalf [\delta(f), \delta(f)] \cr
&&\ \ \ \ + [\delta(f), \Ad(f)^{-1}.\alpha] 
+ \shalf [\Ad(f)^{-1}.\alpha, \Ad(f)^{-1}.\alpha] \cr
&=& \Ad(f)^{-1}.\big(\dd \alpha + \shalf[\alpha,\alpha]\big) 
+ [\Ad(f^{-1}).\delta(f^{-1}), \Ad(f^{-1}).\alpha]\big) \\
&&\ \ \ \ + [\delta(f), \Ad(f)^{-1}.\alpha] \\
&=& \Ad(f)^{-1}.\big(\dd \alpha + \shalf[\alpha,\alpha]\big) 
= \Ad(f)^{-1}.R(\alpha). 
\end{eqnarray*}
Therefore $\alpha * f$ satisfies the Maurer--Cartan equation if and only if 
$\alpha$ does. 
\end{prf} 

The following theorem is particularly useful to deal with Lie group extensions 
$q \: \hat G \to G$ (Definition~\ref{def:liegrp-ext}), where we do not know anything about 
the regularity of the group $G$. In this context one sometimes wants to solve 
equations of the form 
$\delta(f) = \alpha$ for $\alpha \in Z^1_{\rm dR}(M,\hat\g)$, where 
we already have a smooth function $f_G \in C^\infty(M,G)$ with 
$\delta(f_G) = \L(q)\circ \alpha$. 

\begin{thm} [Relative Fundamental Theorem for Lie group extensions] \mlabel{rel-funda}
Let $q \: \hat G \to G$ be an extension of the Lie group 
$G$ by the regular Lie group $N$ and 
$q_\g \: \hat\g \to \g$ be the corresponding extension of the Lie algebra $\g$ of $G$ 
by the Lie algebra $\fn$ of $N$. 

Further, let $M$ be a smooth 
manifold and $\alpha \in \Omega^1(M,\hat\g)$ 
be a $1$-form satisfying the Maurer--Cartan equation. 

If the $\g$-valued $1$-form $q_\g \circ \alpha$ is locally integrable, then 
$\alpha$ is locally integrable.
\end{thm}

\begin{prf} We have to show that 
each point $m \in M$ has an open neighborhood $U$ for which there is a 
smooth function $f_U \: U \to \hat G$ with $\delta(f_U) = \alpha\res_U$. 
This is a local problem. In 
view of the local integrability of $q_\g \circ \alpha$, 
we may therefore assume that 
$q_\g \circ \alpha = \delta(f_G)$ for some smooth function $f_G \: M \to G$. 

For $m \in M$ we choose a $1$-connected open neighborhood $U$ for 
which $f_G(U)$ is contained in an open subset $U_G$ of $G$ on which 
we have a smooth section $\sigma \: U \to \hat G$ of the quotient map 
$q \: \hat G \to G$. In terms of the right action of the group 
$C^\infty(U,\hat G)$ on $\Omega^1(U,\hat\g)$ (Lemma~\ref{lem.e.2}), we now define 
$$ \beta := \alpha * (\sigma \circ f_G)^{-1} 
= \delta((f \circ \sigma_G)^{-1}) + \Ad(\sigma \circ f_G).\alpha 
= \Ad(\sigma \circ f_G).(\alpha - \delta(\sigma \circ f_G)). $$
Since $\alpha$ satisfies the Maurer--Cartan equation, the same holds for $\beta$ 
by (Lemma~\ref{lem.e.2b}). Moreover, 
$$ q_\g \circ \beta 
= \Ad_G(\sigma \circ f_G).(q_\g \circ \alpha - \delta(q \circ \sigma \circ f_G)) 
= \Ad_G(\sigma \circ f_G).(q_\g \circ \alpha - \delta(f_G)) =0, $$
so that $\beta \in \Omega^1(U,\fn)$. 

Now the Fundamental Theorem~\ref{thm-fundamental} implies the existence of a smooth 
function $f_N \: U \to N$ with $\delta(f_N) = \beta$.
With $f := f_N \cdot (\sigma \circ f_G)$ we now get on $U$ the relation  
$\delta(f) = \delta(f_N) * (\sigma \circ f_G) 
= \beta * (\sigma \circ f_G) = \alpha.$
\end{prf}

Since locally integrable $1$-forms on $1$-connected 
manifolds are integrable, we immediately derive with 
the Fundamental Theorem~\ref{thm-fundamental}(ii)(a): 

\begin{cor} \mlabel{cor:4.1.8} 
Let $q \: \hat G \to G$ be an extension of the Lie group 
$G$ by the regular Lie group $N$ and 
$q_\g \: \hat\g \to \g$ be the corresponding extension of $\g$ by~$\fn$. 
Further, let $M$ be a $1$-connected 
manifold and $\alpha \in Z^1_{\rm dR}(M,\g)$. If 
$q_\g \circ \alpha$ is integrable, then $\alpha$ is integrable. 
\end{cor}

\begin{cor}
Let $q \: \hat G \to G$ be a $1$-connected extension of the 
Lie group $G$ by the regular Lie group $N$. Suppose further that 
\[ \psi \in \Aut(\hat\g,\hat\fn) := \{ \phi \in \Aut(\hat\g) \:  \phi(\hat\fn) = \hat\fn \} \] is such that the induced automorphism $\psi_\g$ of $\g$ 
integrates to an automorphism $\psi_G \in \Aut(G)$. 
Then there exists a unique automorphism $\psi_{\hat G} \in \Aut(\hat G)$ with 
$\L(\psi_{\hat G}) = \psi$. 
\end{cor}

\begin{prf}
We put $M := \hat G$ and $\alpha :=  \psi \circ \kappa_{\hat G}$, where 
$\kappa_{\hat G}$ is the left Maurer--Cartan form of $\hat G$. 
For the function $f_G := \psi_G \circ q \: \hat G \to G$, we then obtain 
$$ q_\g \circ \alpha 
=  q_\g \circ \psi \circ \kappa_{\hat G}
=  \psi_\g \circ q_\g \circ \kappa_{\hat G}
=  \delta(\psi_G \circ q) = \delta(f_G), $$
so that the assumptions of Theorem~\ref{rel-funda} are satisfied, and we obtain a 
unique smooth function $\psi_{\hat G} \: \hat G \to \hat G$ with 
$$ \delta(\psi_{\hat G}) = \psi \circ \kappa_{\hat G} 
\quad \hbox{ and } \quad \psi_{\hat G}(\be) = \be. $$
Now Proposition~\ref{prop:homocrit} implies that $\psi_{\hat G}$ is 
a homomorphism of Lie groups with $\L(\psi_{\hat G}) = \psi$. 
Applying the same reasoning to $\psi^{-1}$, 
the fact that homomorphisms are uniquely determined by their derivative in 
$\be$ (Proposition~\ref{propc.14}) 
implies that $\psi_{\hat G}$ is indeed an automorphism. 
\end{prf}

Another application of the Relative Fundamental Theorem is that 
equivalence of simply connected Lie group extensions by regular 
Lie groups follows already from equivalence of the corresponding Lie algebra 
extension. 

\begin{cor}
Let 
$N \sssmapright{i_1} \hat G_1 \sssmapright{q_1} G$ and 
$N \sssmapright{i_2} \hat G_2 \sssmapright{q_2} G$ 
be extensions of~$G$ by the regular Lie group 
$N$ for which the corresponding Lie algebra extensions are equivalent 
in the sense that there exists a Lie algebra isomorphism $\psi \: \hat\g_1 \to \hat \g_2$ 
such that the diagram 
\[  \begin{matrix} 
 \fn & \smapright{\L(i_1)} & \hat \fg_1 &  \smapright{\L(q_1)}  & \g \\ 
\mapdown{\id_\fn} & & \mapdown{\psi} & & \mapdown{\id_\g} \cr 
 \fn &  \smapright{\L(i_2)} & \hat \g_2 &  \smapright{\L(q_2)}  & \g \end{matrix}\] 
commutes. If $\hat G_1$ is simply connected, then there 
exists a unique isomorphism $\hat\phi \: \hat G_1 \to \hat G_2$ with 
$\L(\hat\phi) = \psi$. 
\end{cor}

\begin{prf}
We apply Corollary~\ref{cor:4.1.8} with $M = \hat G_1$, 
$\alpha := \psi \circ \kappa_{\hat G_1}$ and ${f_G := \id_G}$. 
\end{prf}

\subsection{Abelian extensions with smooth global sections} 
\label{sec7.6} 

Let $G$ be a Lie group, let $\fa$ be a smooth $G$-module and 
$\omega \in Z^2(\g,\fa)$ be a continuous $2$-cocycle, so that we can form the 
abelian Lie algebra extension $\fa \oplus_\omega \g$ by the bracket 
\[ [(a,x),(a',x')] := (\omega(x,x') + x.a' - x'.a, [x,x']) \quad \mbox{ on } \quad 
\fa \oplus \g.\] 
In this subsection, we discuss sufficient conditions for the integrability 
of the Lie algebra $\fa \oplus_\omega \g$, provided we already have a Lie group 
$G$ with Lie algebra $\g$. 
In particular, we shall see that our construction always works 
if $G$ is smoothly contractible, or if at least the de Rham cohomology 
groups $H^k_{\rm dR}(G,\fa)$ vanish for $k \in \{1,2\}$. 

\begin{defn}\mlabel{def:smooth-mod} 
Let $(A,+)$ be an abelian Lie group (written additively),  $G$ be  a Lie group 
and $\alpha \:  G \to \Aut(A)$ be a homomorphism defining a smooth 
$G$-action $(g,a) \mapsto \alpha_ga =: g.a$ on $A$. 
We then call $A$ a {\it smooth $G$-module}. 
\index{smooth $G$-module} 
\index{$2$-cocycle} 

A function $f \: G \times G \to A$ 
satisfying $f(g,\be) = f(\be,g) = 0$ and 
\[ f(g,g') + f(gg',gg') = g.f(g',g'') + f(g,g'g'') \quad \mbox{ for } \quad 
g,g',g''\in G \] 
is called a {\it $2$-cocycle}. 
We write $Z^2_s(G,A)$ for those cocycles which are smooth in an 
open neighborhood of $(\be,\be)$ in $G^2$.
\end{defn}

\begin{lem} \label{lem7.1.1} 
Let $G$ be a group, $A$ a smooth $G$-module and $f \: G^2 \to A$ 
be a $2$-cocycle. 
Then we obtain a group $A \times_f G$ by endowing the product set 
$A \times G$ with the multiplication 
\begin{equation}
  \label{eq:7.1.1}
(a,g) (a',g') := (a + g.a' + f(g,g'), gg'). 
\end{equation}
The unit element of this group is $(0,\be)$, inversion is given by 
\begin{equation}
  \label{eq:7.1.2}
(a,g)^{-1} = (-g^{-1}.(a + f(g,g^{-1})),g^{-1}), 
\end{equation}
and conjugation by the formula 
\begin{align} 
  \label{eq:7.1.3}
& (a,g)(a',g') (a,g)^{-1} \notag \\
&= \big(a + g.a' - gg'g^{-1}.a + f(g,g')
-f(gg'g^{-1},g),gg'g^{-1}\big).
\end{align}
The map 
\[ q \: A \times_f G \to G, \quad (a,g) \mapsto g \] is a surjective
homomorphism whose kernel $A \times \{\be\}$ is isomorphic to $A$. 
The conjugation action of $A \times_f G$ on the normal subgroup $A$
factors through the original action of $G$ on $A$. 
\end{lem}

\begin{prf}The condition $f(\be,g) = f(g,\be) = 0$ implies that 
$(0,\be)$ is an identity element in $A \times_f G$, and the associativity of
the multiplication is equivalent to the cocycle condition. 
The formula for the inversion is easily verified. 
Conjugation in $A \times_f G$ is given by 
\begin{eqnarray*}
&&\ \ \ \  (a,g)(a',g') (a,g)^{-1} \\ 
&&= \big(a + g.a' + f(g,g'), gg'\big) \big(-g^{-1}.(a + f(g,g^{-1})),g^{-1}\big)\cr
&&= \big(a + g.a' + f(g,g')- gg'g^{-1}.(a + f(g,g^{-1})) +
f(gg',g^{-1}),gg'g^{-1}\big).
\end{eqnarray*}
To simplify this expression, we use 
$$ f(g,g^{-1}) = f(g,g^{-1}) + f(\be,g) = f(g,\be) + g.f(g^{-1},g) = g.f(g^{-1},g) $$
and 
$$ f(gg',g^{-1}) + f(gg'g^{-1},g) = f(gg',\be) + gg'.f(g^{-1},g) = gg'.f(g^{-1},g) $$
to obtain 
\begin{eqnarray*}
&&\ \ \ \ (a,g)(a',g') (a,g)^{-1} \\ 
&&= \big(a + g.a' + f(g,g')- gg'g^{-1}.a - gg'g^{-1}.f(g,g^{-1}) +f(gg',g^{-1}),gg'g^{-1}\big)\cr
&&= \big(a + g.a' + f(g,g')- gg'g^{-1}.a - gg'.f(g^{-1}.g) +f(gg',g^{-1}),gg'g^{-1}\big)\cr
&&= \big(a + g.a' + f(g,g')- gg'g^{-1}.a -
f(gg'g^{-1},g),gg'g^{-1}\big).  
\end{eqnarray*}
In particular, we obtain 
$$ (0,g)(a,\be) (0,g)^{-1} 
= (g.a,\be).  $$
This means that the action of $G$ on $A$ given by 
$q(g).a := gag^{-1}$ for $g \in A \times_f G$ 
coincides with the given action of $G$ on $A$. 
\end{prf}

The following lemma will be helpful in the proof of
Proposition~\ref{prop7.6.2} below. 

\begin{lem} \label{lem7.6.1} 
Let $G$ be a connected Lie group, $A$ be a smooth
$G$-module and $f \in Z^2_s(G,A)$ be such that all functions 
$f_g \: G \to A, x \mapsto f(g,x)$ are smooth. Then 
$f \: G \times G \to A$ is a smooth function. 
\end{lem}

\begin{prf} We write the cocycle condition as 
\[  f(xy,z) = f(x,yz) + \alpha_xf(y,z)- f(x,y), \quad x,y,z \in G. \] 
For $x$ fixed, this function is smooth as a function of the pair
$(y,z)$ in a neighborhood of $(\be,\be)$. This implies that 
$f$ is smooth on a neighborhood of the points $(x,\be)$, $x \in G$. 
Fixing $x$ and $z$ shows that there exists an $\be$-neighborhood $V
\subeq G$ (independent of $x$) such that the functions $f^z :=f(\cdot, z)$, 
$z \in V$, are smooth
in a neighborhood of $x$. Since $x \in G$ was arbitrary, we conclude
that the functions $f^z$, $z \in V$, are smooth.  
Now 
\[  f^{yz} = f^z(\cdot y) - \alpha(\cdot)f(y,z) + f^y\] 
shows that the same holds for the functions $f^u$, $u \in V^2$. 
Iterating this process, using $G = \bigcup_{n \in \N} V^n$, we
derive that $f^x$ is smooth for every $x \in G$. 
Finally we see that the function 
\[  (x,y) \mapsto  f(x,yz)  = f(xy,z) - \alpha_xf(y,z) + f(x,y) \] 
is smooth in a neighborhood of each point $(x_0, \be)$, hence that $f$
is smooth in each point $(x_0, z_0)$. This proves that $f$ is
smooth on $G \times G$. 
\end{prf}

For the following proposition we recall equivariant differential 
forms introduced in Definition~\ref{def:equiv-form}. 

\begin{prop} \mlabel{prop7.6.2} 
Let $G$ be a 
connected Lie group, $\fa$ a Mackey complete locally convex 
smooth $G$-module, $\omega \in Z^2_c(\g,\fa)$ a
continuous $2$-cocycle, and $\omega^{\rm eq} \in \Omega^2(G,\fa)$ 
the corresponding equivariant $2$-form on $G$ with $\omega^{\rm eq}_\be = \omega$. 
We assume that 
\begin{description}
\item[\rm(1)] $\omega^{\rm eq} = d\theta$ for some $\theta \in \Omega^1(G,\fa)$, 
and 
\item[\rm(2)] for each $g \in G$, the closed $1$-form $\lambda_g^* \theta - \alpha_g
\circ \theta$ is exact. 
\end{description} 
Both assumptions are satisfied if $G$ is smoothly contractible. 

Then there exists a smooth $2$-cocycle $f \: G^2 \to\fa$ 
defining a Lie group structure on the product manifold 
$\hat G := \fa \times G$
\[  (a,g) (a',g') := (a + g.a' + f(g,g'),gg') \] 
and the Lie algebra $\hat\g$ of $\hat G$ is isomorphic to 
$\fa \oplus_\omega\g$. 
\end{prop}

\begin{prf} If $G$ is smoothly contractible, then the Poincar\'e 
Lemma~\ref{poincarelemma} implies that closed $\fa$-valued forms of 
positive degree are exact. Note that $\omega^{\rm eq}$ is closed by 
Proposition~\ref{prop:CE} and that, 
for each $g \in G$, the relation 
$\alpha_g \circ \omega^{\rm eq} = \lambda_g^*\omega^{\rm eq}$ leads to 
\[  d\big(\alpha_g \circ \theta - \lambda_g^*\theta\big) = 
\alpha_g \circ \omega^{\rm eq} - \lambda_g^*\omega^{\rm eq} = 0.\] 
Therefore (1) and (2) are satisfied or every smoothly contractible 
Lie group. 

From (2), we obtain 
 for each $g \in G$ a smooth function $f_g \: G \to \fa$
with $f_g(\be) = 0$ and 
\[ df_g = \lambda_g^*\theta - \alpha_g  \circ \theta. \] 
Observe that $f_\be = 0$. For $g, h \in G$ this leads to 
\begin{eqnarray*}
df_{gh} 
&&= \lambda_{gh}^*\theta - \alpha_{gh} \circ\theta 
= \lambda_h^*(\lambda_{g}^*\theta - \alpha_g \circ\theta) + \lambda_h^*(\alpha_g\circ\theta) - \alpha_{gh} \circ\theta \cr 
&&= \lambda_h^*df_g + \alpha_g (\lambda_h^*\theta - \alpha_h \circ\theta) 
= \lambda_h^*df_g + \alpha_g  \circ df_h \\
&&= d(f_g \circ \lambda_h + \alpha_g  \circ f_h). 
\end{eqnarray*}
Comparing values of both functions in $\be$, we get 
\begin{equation}
  \label{eq:7.6.1}
f_{gh} = f_g \circ \lambda_h + \alpha_g \circ f_h - f_g(h). 
\end{equation}
Now we define $f \: G \times G \to \fa$ by $f(x,y) := f_x(y)$. Then 
\eqref{eq:7.6.1}  means that 
$$ f(gh,u) = f(g,hu) + \alpha_g f(h,u) - f(g,h), \quad g,h,u \in G, $$
so that $f$ is a group cocycle. 

The concrete local formula for $f_x$ in the 
Poincar\'e Lemma (Theorem~\ref{poincarelemma}) 
and the smooth dependence of the integral on $x$ imply that 
$f$ is smooth on a neighborhood of $(\be,\be)$, so that Lemma~\ref{lem7.6.1}
implies that $f \: G \times G \to \fa$ is a smooth function. 
With Lemma~\ref{lem7.1.1} we therefore obtain on the space 
$\hat G := \fa \times G$
a Lie group structure with the multiplication given by 
\[  (a,g) (a',g') := (a + \alpha_ga' + f(g,g'), gg') \] 
and the corresponding Lie bracket is given by 
\[  [(a,x), (a',x')] 
= \big(x.a' - x'.a + \shalf\delta^2_{(\be,\be)}f (x,x') - 
\shalf\delta^2_{(\be,\be)}f(x',x), 
[x,x']\big).\]
Here $\delta^2_{(\be,\be)} f$  is well-defined and bilinear because 
$f(\be,x) = f(y,\be) = 0$ for $x,y \in G$, 
which implies in particular $\dd f(\be,\be) = 0$
(cf.~Exercise~\ref{exer:2.x.x}). 
Now we relate this formula to the Lie algebra cocycle $\omega$. 
The relation $df_g = \lambda_g^* \theta - \alpha_g \circ\theta$ leads to 
$$ df(g,\be)(0,y) =  df_g(\be)y = (\lambda_g^* \theta - \alpha_g \circ \theta)_\be(y) 
= \la \theta, y_l \ra(g) - \alpha_g \theta_\be(y), $$
where $y_l$ denotes the left invariant vector field with $y_l(\be) =
y$. Taking second derivatives, we further obtain for $x \in \g$: 
\begin{eqnarray*}
\frac{1}{2} \delta^2_{(\be,\be)} f(x, y)  &&=  x_l(\la \theta, y_l \ra)(\be) - x\theta_\be(y) \\
&&=  (d\theta)(x_l, y_l)(\be) + y_l(\la \theta, x_l \ra)(\be) 
+ \theta([x_l,y_l])(\be)- x\theta_\be(y) \cr
&&=  \omega(x, y) + y_l(\la \theta, x_l \ra)(\be) + \theta_\be([x,y])- x\theta_\be(y), 
\end{eqnarray*}
Subtracting 
$\frac{1}{2} \delta^2_{(\be,\be)} f(y,x) 
 =  y_l(\la \theta, x_l \ra)(\be) - y\theta_\be(x),$
leads to the $2$-cocycle 
\[ \eta(x,y) := 
 \omega(x, y) + \theta_\be([x,y]) - x\theta_\be(y) + y\theta_\be(x) 
=  \omega(x, y) - (d\theta_\be)(x,y). \] 
Since this cocycle is equivalent to $\omega$, $\L(\hat G) \cong 
\fa \oplus_\omega \g$. 
\end{prf}

\begin{cor} \label{cor7.6.3} 
If $G$ is simply connected and 
$\omega^{\rm eq}$ is exact, then there exists a smooth cocycle $f \: G\times G
\to \fa$, so that $\hat G := \fa \times_f G$ is a Lie group with Lie
algebra $\hat\g = \fa \oplus_\omega \g$. 
\end{cor}

\begin{prf}Since $\pi_1(G)$ is trivial and $\fa$ is Mackey complete, 
$H^1_{\rm dR}(G,\fa)$ vanishes 
by the Fundamental Theorem~\ref{thm-fundamental}, so that condition (2) in
Proposition~\ref{prop7.6.2} is automatically satisfied. 
\end{prf}

\subsection{Generalized 
central extensions of locally exponential Lie algebras} 
\mlabel{sec:11.3}

In this subsection, we discuss the concept of a generalized central extension 
of a locally convex Lie algebra. It generalizes central extensions 
$q \: \hat\g \to \g$ 
which are quotient maps with central kernel to maps with dense range. 
It is an extremely useful tool which can be used to show the 
local exponentiality of central extensions of locally exponential Lie algebras 
for which $\ker q$ is not split. 
Here the key point is that any central extension can be embedded as a closed 
ideal into a topologically split central extension 
(Proposition~\ref{propD.3.5}).

\index{generalized central extension} 
\begin{defn} \mlabel{def:F.3.1} 
Let $\g$ and $\hat\g$ be locally convex Lie algebras. 
A continuous Lie algebra homomorphism $q \: \hat\g \to \g$ with dense range 
is called a {\it generalized central extension} if there exists a continuous 
bilinear map 
$b \: \g \times \g \to \hat\g$ with 
\begin{eqnarray}
  \label{eq:f.4}
b(q(x), q(y)) = [x,y] \quad \hbox{ for } \quad x,y \in \hat\g. 
\end{eqnarray}
Since $q$ has dense range, the map 
$b$ is uniquely determined by (\ref{eq:f.4}) and (\ref{eq:f.4}) implies that $\ker
q$ is central in $\hat\g$. 
\end{defn} 

\begin{exs} \mlabel{exs:remD.3.2} (a) If $q \: \hat\g \to \g$ is a quotient homomorphism of 
locally convex Lie algebras with central kernel, i.e., a 
\index{central extension!of Lie algebra} 
{\it central extension}, 
then $$q \times q \: \hat\g \times \hat\g \to \g \times \g$$ also is a 
quotient map. Therefore the Lie bracket of $\hat\g$ factors through 
a continuous bilinear map $b \: \g\times \g \to \hat\g$ 
with $b(q(x), q(y)) = [x,y]$ for $x,y \in \hat\g$, 
showing that $q$ is a generalized central extension.

(b) Let $\omega \in Z^2_c(\g,\z)$ be a continuous $2$-cocycle, i.e., 
an alternating bilinear map satisfying 
\begin{equation}
  \label{eq:2-coc}
  \omega([x,y],z) +   \omega([y,z],x) + \omega([z,x],y) = 0 \quad 
\mbox{ for } \quad x,y,z \in \g.
\end{equation}
Then we obtain a central extension 
$p \: \fz \oplus_\omega \g \to \g, (z,x) \mapsto x$, where 
$\fz \oplus\g$ denotes the topological product space $\fz \times \g$, 
endowed with the Lie bracket 
\[ [(z,x),(z',x')] := (\omega(x,x'), [x,x']), \quad z,z'\in \fz, 
x,x'\in\g.\] 
Suppose that $\hat\g \subeq \z \oplus_\omega \g$ is a closed subalgebra 
for which $p(\hat\g)$ is dense in $\g$. We claim that 
$q := p\res_{\hat\g} : \hat \g \to \g$ is a generalized central extension with 
$b(x,y) = (\omega(x,y),[x,y])$ for $x,y\in\g$.
In fact, the range of $q$ is dense by the assumption that 
$p(\hat\g)$ is dense in $\g$. It is also clear that 
$b \circ (p \times p)$ is the bracket on $\fz \oplus_\omega \g$, 
but it remains to show that $\im(b) \subeq \hat\g$. 
For $x = q(x'), y = q(y')$ in $\im(q) = p(\hat\g)$ we have 
$$ b(x,y) = b(q(x'), q(y')) = [x',y'] \in \hat\g. $$
Now the continuity of $b$, the density of $\im(q)$ in $\g$, and the closedness of 
$\hat\g$ imply that $\im(b) \subeq\hat\g$. 
\end{exs}

\begin{lem} \mlabel{lemD.3.4} 
For a generalized central extension
$q \: \hat\g\to\g$ with the corresponding map $b$, 
 the following assertions hold: 
\begin{description}
\item[\rm(1)] $[x,y] = q(b(x,y))$ for all $x,y \in \g$. 
\item[\rm(2)] $[\g,\g] \subeq \im(q)$ and $\ker q \subeq \z(\hat\g)$. 
\item[\rm(3)] $b \in Z^2_c(\g,\hat\g)$, i.e., 
$b([x,y],z) + b([y,z],x) +b([z,x],y) =0$ for $x,y, z \in \g$.
\item[\rm(4)]  For $x \in \g$ we define 
$$ \hat\ad(x) \: \hat\g \to \hat \g, \quad y \mapsto b(x,q(y)). $$
Then $\hat\ad$ defines a continuous representation of $\g$ on $\hat\g$
by derivations for which $q$ is equivariant with respect to the
adjoint representation of~$\g$. 
\item[\rm(5)] If $\hat\g$ is topologically perfect, then 
$q^{-1}(\z(\g)) = \z(\hat\g)$. 
\end{description}
\end{lem}

\begin{prf} (1) If $x = q(a)$ and $y = q(b)$ holds for $a, b \in \hat\g$, then 
  \begin{equation}
    \label{eq:8.3.1}
[x,y] = [q(a), q(b)] = q([a,b]) = q(b(x,y)).
  \end{equation}
Therefore the Lie bracket on $\g$ coincides on the dense subset 
$\im(q) \times \im(q)$ of $\g \times \g$ with the continuous map 
$q \circ b$, so that (1) follows from the continuity of both maps. 

\nin (2) The first assertion follows from (1) and the second assertion 
from \eqref{eq:f.4}. 

\nin (3) In view of (\ref{eq:f.4}), the Jacobi identity in $\hat\g$ leads to 
\begin{eqnarray*}
0 &=& [[x,y],z] + [[y,z],x] + [[z,x],y] \\
&=& b(q([x,y]), q(z)) + b(q([y,z]), q(x)) + b(q([z,x]), q(y)) \\
&=& b([q(x),q(y)], q(z)) + b([q(y), q(z)], q(x)) + b([q(z), q(x)], q(y)).
\end{eqnarray*}
Therefore the restriction of $b$ to $\im(q)$ is a Lie algebra 
cocycle, and since $\im(q)$ is dense and $b$ is continuous, 
$b$ is a Lie algebra cocycle on $\g$. 

\nin (4) First we observe that the bilinear 
map $\g \times \hat\g \to \hat\g, (x,y) \mapsto b(x,q(y))$ is continuous. 
Moreover, (1) implies 
$$ q(\hat\ad(x)y) = q(b(x,q(y))) = [x,q(y)],   $$
i.e., $q \circ \hat\ad(x) = \ad x \circ q$. 

From the cocycle identity 
$$ b([x,y],z) + b([y,z],x) + b([z,x],y) = 0, \quad x,y,z \in \g, $$
we derive in particular for $x \in \g$ and $y,z \in \hat\g$: 
\begin{eqnarray*}
0 
&=& b([x,q(y)],q(z)) + b([q(y), q(z)],x) + b([q(z),x],q(y)) \\ 
&=& b(q(\hat\ad(x)y), q(z)) + b(q([y,z]),x) -  b(q(\hat\ad(x)z),q(y)) \\ 
&=& [\hat\ad(x)y,z] -  \hat\ad(x)[y,z] -  [\hat\ad(x)z,y]. 
\end{eqnarray*}
Therefore each $\hat\ad(x)$ is a derivation of $\hat\g$. 
On the other hand, the cocycle identity for $b$ leads for $x,y \in \g$
and $z \in \hat\g$ to 
\begin{eqnarray*}
0 &=& b([x,y],q(z)) + b([y,q(z)],x) + b([q(z),x],y)\\ 
&=& \hat\ad([x,y])z + b(q(\hat\ad(y)z),x) - b(q(\hat\ad(x)z),y) \\
&=& \hat\ad([x,y])z - \hat\ad(x)\hat\ad(y)z + \hat\ad(y)\hat\ad(x)z,   
\end{eqnarray*}
so that 
$\hat\ad \:\g \to \der(\hat\g)$
is a representation of $\g$ by derivations of $\hat\g$, and
the map $q$ is equivariant with respect to the adjoint representation
of~$\g$. 

\nin (5) Let $\hat\z(\g) := q^{-1}(\z(\g))$. We first observe that 
$[\hat \z(\g), \hat \g]$ is contained in $\ker q \subeq \z(\hat \g)$ because 
$q([\hat \z(\g), \hat \g]) \subeq [\z(\g), \g] = \{0\}.$
This leads to 
$$ [\hat \z(\g), [\hat \g, \hat \g]] 
\subeq [\hat \g, [\hat \z(\g), \hat \g]] \subeq [\hat \g, \ker q] = \{0\}. $$
If $\hat \g$ is topologically perfect, we obtain 
$\hat \z(\g) \subeq \z(\hat \g)$. The other inclusion follows from 
the density of the image of $q$. 
\end{prf}

The following proposition shows that generalized central extensions 
can be characterized as certain closed subalgebras of central extensions 
defined by cocycles (cf.\ Example~\ref{exs:remD.3.2}(b)). 

\begin{prop} \mlabel{propD.3.5} 
Write $|\g|$ for the locally convex space underlying 
$\g$, considered as an abelian Lie algebra. 
If $q \: \hat\g \to \g$ is a generalized central 
extension and  $b \: \g \times \g \to \hat\g$ the corresponding cocycle, then 
the map 
$$ \psi \: \hat\g \to |\hat\g| \oplus_b \g, \quad x \mapsto (x,q(x)) $$
is a topological embedding of $\hat\g$ onto a closed 
ideal of $|\hat\g| \oplus_b \g$, containing the commutator algebra. 
The map 
$$ \eta \: |\hat\g| \oplus_b \g \to |\g|, \quad (x,y) \mapsto y - q(x)  $$
is a quotient morphism of Lie algebras whose kernel is $\im(\psi) \cong \hat\g$. 
\end{prop}

\begin{prf} The bracket in $|\hat\g| \oplus_b \g$ 
is given by 
\begin{equation}
  \label{eq:8.3.2}
[(x,y),(x',y')] = (b(y,y'),[y,y']). 
\end{equation}
Now 
\begin{eqnarray*}
[\psi(x), \psi(x')]  
&=& [(x,q(x)),(x', q(x'))] = (b(q(x),q(x')),[q(x),q(x')]) \\ 
&=& ([x,x'], q([x,x'])) = \psi([x,x'])  
\end{eqnarray*}
implies that the continuous linear map $\psi$ is a morphism of Lie algebras. 
As the graph of the continuous linear map $q$, the image of $\psi$ is a
closed subspace of $|\hat\g| \oplus_b \g$, and the projection onto the
second factor is a continuous linear map. Therefore 
$\psi$ is a topological embedding onto a closed subalgebra. 

Moreover, \eqref{eq:8.3.2} and 
$q(b(x,x')) = [x,x']$ show that $\im(\psi)$ contains all brackets, hence is an ideal. Therefore 
the map $\eta \: |\hat\g| \oplus_b \g\to |\g|$ whose kernel is $\im(\psi)$ is a morphism of Lie algebra. 
That it is a quotient map follows from the fact that its restriction to the subspace $\g$ is a topological 
isomorphism. 
\end{prf}


\begin{thm} \mlabel{centext-thm} If $q \: \hat\g \to \g$ is a Mackey complete 
generalized central extension of 
the locally exponential Lie algebra $\g$, then $\hat\g$ is also locally 
exponential. 
\end{thm}

\begin{prf} {\bf Step 1:} First we consider the special case of 
a topologically split central extension 
$\hat\g = \z \oplus_\omega \g$ of $\g$ by a Mackey complete space $\fz$, 
where $\omega \in Z^2(\g,\fz)$ is a continuous $2$-cocycle. 
 
Let $(\g, *, D)$ be the data of a locally exponential Lie algebra on $\g$ 
and $U_n$ as in Definition~\ref{def:5.2.1}. 
We then define a left invariant $\fz$-valued 
$2$-form $\omega^l$ on $U_1$ by 
\[  \omega^l(x)(v,w) := \omega(d\lambda_{-x}^*(x)v, d\lambda_{-x}^*(x)w). \] 
Since $\omega$ is a $2$-cocycle, $\omega^l$ is closed 
by Proposition~\ref{prop:CE}. 
Since $U_1$ is convex and 
$\z$ is Mackey complete, there exists a $\z$-valued $1$-form 
$\theta$ on $U_1$ with $\omega^l = d\theta$ 
(Corollary~\ref{cor:poincare-lemma}). 

A corresponding smooth local group cocycle can be constructed as follows. 
We connect two points $x,y \in U_2$ by the path 
\[ \alpha_{x,y} \: [0,1] \to U, \quad t \mapsto x * t((-x) * y).\] 
For $x,y, z, z* x, z * y \in U_2$, we then have 
\[  \alpha_{z*x, z*y} = \lambda_z^* \circ \alpha_{x,y}
\quad \mbox{ and } \quad 
\alpha_{0,y}(t) = ty. \]
Recall the standard $2$-simplex 
$\Delta_2 := \{(s,t) \in \R^2 \: 0 \leq s,t; s + t \leq 1\}$.
For $x,y,z \in U_2$, we obtain a singular smooth $2$-cycle 
$$ \alpha_{x,y,z} := \alpha_{x,y} + \alpha_{y,z} - \alpha_{x,z}, $$
which corresponds to the piecewise smooth map 
$\alpha_{x,y,z} \: \partial\Delta_2 \to \g$ given by 
$$ \alpha_{x,y,z}(s,t) = 
\left\{ 
  \begin{array}{cl}
\alpha_{x,y}(s), & \mbox{for $t = 0$} \\ 
\alpha_{y,z}(1-s), & \mbox{for $s + t = 1$} \\ 
\alpha_{x,z}(t), & \mbox{for $s = 0$}. 
  \end{array}
\right.  $$ 
For $x,y,z \in U_3$ we now consider the smooth map 
\begin{align}\label{eq:gamma}
 \gamma_{x,y,z} &\: \Delta_2 \to U, \notag \\
(s,t) &\mapsto x * \big( s\cdot \alpha_{(-x)*y,(-x)*z}(t) 
+ t \cdot \alpha_{(-x)*y,(-x)* z}(1-s)\big) 
\end{align}
whose restriction to the boundary $\partial \Delta_2$ coincides with 
$\alpha_{x,y,z}$. For $a,x,y,z \in U_4$ we further have 
\begin{eqnarray}
  \label{eq:5.1.5}
\gamma_{a*x, a*y, a * z} = \lambda_a^* \circ \gamma_{x,y,z}.
\end{eqnarray}
Now we consider the smooth function 
\begin{align} \label{eq:F(x,y,z)}
F &\: U_3 \times U_3  \times U_3 \to \z, \notag \\  
F(x,y,z) 
&:= \int_{\gamma_{x,y,z}} \omega^l 
= \int_{\gamma_{x,y,z}} d\theta 
= \int_{\alpha_{x,y}} \theta + \int_{\alpha_{y,z}} \theta 
- \int_{\alpha_{x,z}} \theta. 
\end{align}
From the left invariance of $\omega^l$ and (\ref{eq:5.1.5}) we obtain 
$$ F(a*x, a*y, a*z) = F(x,y,z)\quad \mbox{for} \quad a,x,y,z \in U_4. $$
For the function 
$$ f\: U_3 \times U_3 \to \z, \quad (x,y) \mapsto F(0,x,x*y) $$
this leads for $x,y,z \in U_4$ to the cocycle relation 
\begin{eqnarray*} 
&&\ \ \ \  f(x,y) + f(x*y,z) - f(x, y * z) - f(y,z) \\
&&= F(0,x,x*y) + F(0, x*y, x*y*z) - F(0, x, x * y * z) - F(0,y, y * z) \\ 
&&= F(0,x,x*y) + F(0, x*y, x*y*z) - F(0, x, x * y * z) \\
&&\ \ - F(x,x*y, x* y * z) \\ 
&&= \int_{\alpha_{0,x}} \theta + \int_{\alpha_{x,x*y}} \theta - \int_{\alpha_{0,x*y}}\theta 
+ \int_{\alpha_{0,x*y}} \theta + \int_{\alpha_{x*y,x*y*z}} \theta 
- \int_{\alpha_{0,x*y*z}}\theta \\ 
&&\ \ \ \ -\int_{\alpha_{0,x}} \theta - \int_{\alpha_{x,x*y*z}} \theta 
+ \int_{\alpha_{0,x*y*z}}\theta 
- \int_{\alpha_{x,x*y}} \theta -\int_{\alpha_{x*y,x*y*z}} \theta 
+ \int_{\alpha_{x,x*y*z}}\theta \\ 
&&= 0. 
\end{eqnarray*}
We now define a smooth associative multiplication on 
$\hat U := \z \times U_4 \subeq \hat\g$
by 
\begin{eqnarray}
  \label{eq:5.1.6}
(z,x) * (z', x') := (z + z' + f(x,x'), x * x').
 \end{eqnarray}

If $x' \in \R x$, then the image of the map 
$\gamma_{0,x,x*x'}$ lies in the line $\R x$, and therefore 
the integral \eqref{eq:F(x,y,z)} defining $f(x,x')$ vanishes, so that 
$$ (z,x) * (z', x') = (z + z', x +x'), $$
and this implies \eqref{eq:5.2.1}. 
Next we observe that the cocycle 
relation for $f$ implies that the multiplication (\ref{eq:5.1.6}) is 
associative on all triples in $\fz \times U_5$ on which it is defined. 
We recall from Remark~\ref{rem:locgrp} that 
the second order term 
$\hat b_2$ of the Taylor series of the product (\ref{eq:5.1.6})  
defines a Lie bracket on 
$\hat\g$ which if of the form 
\[ \hat b_2((z,x),(z',x'))
- \hat b_2((z',x'),(z,x))
 = (\eta(x,x'), [x,x']), \]
where $\eta \: \g \times \g \to \fz$ is alternating and bilinear. 
It remains to show that $\eta = \omega$, i.e., 
\begin{equation}
  \label{eq:d2omega}
\frac{1}{2}\delta^2_{(0,0)} f(x,y) - \frac{1}{2}\delta^2_{(0,0)} f(y,x) = \omega(x,y).
\end{equation}
This implies that $\hat\g$ is the Lie algebra of the exponential 
local Lie group $(U_4, *)$ that we constructed, hence locally exponential. 

To verify \eqref{eq:d2omega}, we put 
$\gamma_{x,y} := \gamma_{0,x,x*y}$, so that 
\[ f(x,y) = \int_{\gamma_{x,y}}\omega^l 
= \int_{\Delta^2} \omega^l(\gamma_{x,y}(s,t))
\Big(\frac{\partial \gamma_{x,y}}{\partial s}(s,t), 
\frac{\partial \gamma_{x,y}}{\partial t}(s,t)\Big)\, d(s,t). \]
From 
\begin{equation}
  \label{eq:gamma2}
\gamma_{x,0}(s,t) = (s+t)x  \quad \mbox{ and }\quad \gamma_{0,y}(s,t) = ty 
\end{equation}
it follows in particular that 
$\frac{\partial \gamma_{0,0}}{\partial s} 
= \frac{\partial \gamma_{0,0}}{\partial t} = 0$ 
and that $\frac{\partial}{\partial s} \gamma_{x,y} \wedge \frac{\partial}{\partial t} 
\gamma_{x,y} = 0$ for $x = 0$ or $y = 0$. We conclude that 
$f(x,0) = f(0,y) = 0$, so that the $1$-jet of $f$ in $(0,0)$ 
vanish and the quadratic map $\delta^2_{(0,0)} f$ is bilinear 
(Lemma~\ref{lem:second-order}). 
We further derive from \eqref{eq:gamma2} that the first order term of 
$\frac{\partial \gamma_{x,y}}{\partial s}$ in $(0,0)$ is $x$ 
and the first order term of 
$\frac{\partial \gamma_{x,y}}{\partial t}$ in $(0,0)$ is $x+y$. 
For the integrand we thus obtain the second order term as 
\[ \omega^l(\gamma_{0,0}(s,t))(x,x+y) 
= \omega(x,y), \]
so that 
\[  \frac{1}{2} \delta^2_{(0,0)} f(x,y) 
= \int_{\Delta_2} \, ds\, dt \cdot \omega(x,y) = \frac{1}{2} \omega(x,y). \]
This proves~\eqref{eq:d2omega}. 

\nin{\bf Step 2:} Let $b \in Z^2_c(\g,\hat\g)$ be the cocycle of the generalized 
central extension (Definition~\ref{def:F.3.1}) and 
$$ \g_1 := |\hat\g| \oplus_b \g, \quad q_1 \: \g_1 \to \g, \quad (z,x) \mapsto x $$
denote the corresponding topologically split 
central extension of $\g$ by the trivial $\g$-module 
$\hat\g$, which is denoted by $|\hat\g|$ to distinguish it from the Lie algebra 
$\hat\g$. Then the map 
$$ \alpha \: \g_1 \to |\g|, \quad(z,x) \mapsto x - q(z) $$
is a homomorphism of Lie algebras whose kernel 
$$ \ker \alpha = \{ (z,q(z)) \: z \in \hat\g \} $$
is isomorphic to $\hat\g$ as a Lie algebra (Proposition~\ref{propD.3.5}) 
and it is a closed ideal containing all commutators. 

Since every abelian Lie algebra is locally exponential 
and by Step 1, the topologically split central extension $\g_1$ of 
$\g$ is locally exponential, the ideal 
$\hat\g$ of $\g_1$ is locally exponential 
by Proposition~\ref{prop:invim}.
\end{prf}

The following corollary provides a criterion for the enlargeability 
of a locally exponential Lie algebra $\g$ if it already known that there exists a 
Lie group $G$ with Lie algebra $\g$ which possesses an exponential function. 

\begin{cor} \mlabel{cor:locexp-centext-enl-crit} 
Let $G$ be a Lie group 
with Mackey complete Lie algebra $\g$ and 
and an exponential function 
$\exp \:\g \to G$. Suppose further that 
$\fz \subeq \g$ is a closed central subspace 
for which $\g/\fz$ is locally exponential and enlargeable. 
Then $\g$ is locally exponential and, if the additive subgroup 
\[ \Gamma :=\{ z \in \fz \: \exp_G z = \be \} \] 
of $\fz$ is discrete, then $\g$ is enlargeable. 
\end{cor}

\begin{prf} First we use Theorem~\ref{centext-thm}  to see that 
$\g$ is locally exponential. Replacing $G$ by the universal 
covering group of its identity component, we may assume that 
$G$ is $1$-connected. In view of Corollary~\ref{cor:5.5.5b}, it 
suffices to show that 
\[ \Gamma_\g := \{ z\in \fz(\g) \: \exp_G z = \be \} \] 
is discrete. To derive this from our assumption, 
let $Q$ be a connected locally exponential Lie group 
with Lie algebra $\fq := \g/\fz$. By the Integration Theorem~\ref{thm:int-thm}, 
there exists a morphism of Lie groups 
$\phi \:  G \to Q$ for which $q := \L(\phi) \: \g \to \fq$ is the quotient map. 
Let $U \subeq\g$ be a convex $0$-neighborhood for which 
$\exp_Q$ is injective on $q(U)$ and $U \cap \Gamma = \{0\}$. 
For $z \in \Gamma_\g \cap U$ we 
then have $\exp_Q q(z) = \phi(\exp z) = \be = \exp_Q 0$ and therefore 
$q(z) = 0$, i.e., $z \in \fz$. We conclude that $z \in \Gamma \cap U = \{0\}$. 
Therefore $\Gamma_\g$ is discrete and the assertion follows from 
Corollary~\ref{cor:5.5.5b}. 
\end{prf}

\begin{prop}
  \mlabel{prop:VI.1.6b}  {\rm(Enlargeability of generalized central extensions)} Let 
$G$ be a smoothly contractible locally exponential Lie group with Lie algebra $\g$ and 
$q \: \hat\g\to \g$ be a generalized central extension for which $\hat\g$ 
is Mackey complete. Then $\hat\g$ is enlargeable. 
\end{prop}  

\begin{prf} Let $b \: \g \times \g \to \hat\g$ be the corresponding 
cocycle defined by 
\[ b(q(x),q(y))= [x,y]\quad \mbox{ for }  \quad x,y \in \hat\g.\] 
First we consider the corresponding topologically split central extension 
\[ \g_1 := |\hat\g| \oplus_b \g. \] 
Theorem~\ref{centext-thm} shows that $\hat\g$ and $\g_1$ are locally exponential. 
Further \[ \hat\g \to \g_1,\quad  x \mapsto (x,q(x))\]  realizes $\hat\g$ as a closed 
ideal of $\g_1$. By the Integral Subgroup Theorem~\ref{thm:5.5.3}, 
it therefore suffices to prove the enlargeability 
of $\g_1$. With Proposition~\ref{prop7.6.2} we obtain 
a central Lie group extension $G_1$ with Lie algebra $\g_1$ 
which is diffeomorphic to $\hat\g \times G$. Finally 
Corollary~\ref{cor:8.4.2} implies that $G_1$ is locally exponential.
\end{prf}

\begin{rem} \mlabel{rem:8.3.15} 
The preceding proposition can be generalized to the situation, 
where $G$ is only contractible, not necessarily smoothly. 
Then one can use \cite[Thm.~7.12]{Ne02a} 
instead of Proposition~\ref{prop7.6.2} to integrate
the topologically split central extension $\g_1$ 
to a Lie group extension. The other arguments are the same. 
\end{rem}

\nin{\bf Notes on Section~\ref{sec:8.3}:}\\ 
\nin {\bf \ref{sec7.6}:} For central extensions of finite-dimensional Lie groups, 
the construction described in Proposition~\ref{prop7.6.2} is  
due to \'E.~Cartan, 
who used it to construct a central extension of a simply connected
finite-dimensional Lie group $G$ by the additive group $\fz$. Since in this case 
\[ H^2_{\rm dR}(G,\fz) \cong \Hom(\pi_2(G),\fz)= \{0\} 
\ \  \mbox{ and }\ \ 
H^1_{\rm dR}(G,\fz) \cong \Hom(\pi_1(G),\fz) = \{ 0\},\] 
(cf.\ \cite{God71}; see also Section~\ref{sec:10.2}), 
the requirements of the construction are satisfied for
every Lie algebra cocycle $\omega \in Z^2_c(\g,\fz)$. 
Appying this to the adjoint group $G_{\rm ad}$ of a finite-dimensional 
Lie algebra $\g$ and the central extension $\g$ of $\g_{\rm ad}$ by $\fz(\g)$, 
this argument implies the integrability of $\g$; see also 
Corollary~\ref{cor:lie3} below. 

\nin{\bf \ref{sec:11.3}}: An interesting aspect of generalized 
central extensions which is not exploited here, is that 
topologically perfect 
locally convex Lie algebras have universal generalized central 
extensions. This is remarkable because 
universal central extensions do not always exist, not even for 
topologically perfect Banach--Lie algebras. We refer to 
\cite{Ne03} for a more detailed discussion of this issue. 

For the case of central extensions of Banach--Lie algebras, 
part of Proposition~\ref{propD.3.5} can already be found in a  
footnote in \cite[p.~58]{vES73}.

\section{Enlargeability of locally exponential Lie algebras} 
\mlabel{sec:11.4}  

Let $\g$ be a Mackey complete 
locally exponential Lie algebra for which 
\[ \g_{\rm ad} = \g/\fz(\g) \] is also Mackey complete. This condition is satisfied 
if either $\g$ is Fr\'echet or $\g$ is Mackey complete with 
$\fz(\g)$ a split subspace. 
We assign to $\g$ an additive subgroup $\Pi(\g)$ of its center $\fz(\g)$, called its 
period group. The period group has the remarkable property 
that it detects enlargeability: $\Pi(\g)$ is discrete if and only if 
$\g$ is enlargeable. If this is the case 
and $G$ is a $1$-connected locally exponential Lie group with Lie algebra $\g$, 
then 
\[ \Pi(\g) \cong\pi_1(Z(G)).\] 
Since Banach--Lie algebras are locally exponential, 
this leads in particular a convenient enlargeability criterion
 for Banach--Lie algebras. 
%

\subsection{Path groups of Lie groups} 
\mlabel{subsec:11.4.1} 

\begin{defn} \mlabel{def4.4.1} 
Let $G$ be a connected Lie group with Lie algebra $\g$ 
and $r \in \N_0 \cup \{\infty\}$. We write 
\[ P^{(r)}(G) := \{ \gamma \in C^r([0,1],G)  : \; \gamma(0) = \be\} \]
for the {\it $C^r$-path group of $G$}, 
\index{path group} 
\index{loop group} 
where the multiplication on $P^{(r)}(G)$ is
defined pointwise. This is a Lie group 
by Theorem~\ref{thm:mapgro-Lie} 
and 
\[ P^{(r)}(\g) := \{ \gamma \in C^r(\, [0,1],\g) : \; \gamma(0) = 0\} \] 
is the Lie algebra of $P^{(r)}(G)$. The evaluation map
\[  \ev_1 \!: P^{(r)}(G) \to G, \;\;\;\; \gamma \mapsto \gamma(1) \] 
is a surjective morphism of Lie groups. Its kernel 
\[ \Omega^{(r)}(G) = \{ \gamma \in C^r([0,1],G) \: \gamma(0) = \gamma(1) =\be \} \]
is called the {\it $C^r$-loop group of $G$}.
It is a Lie group with Lie algebra 
\[ \Omega^{(r)}(\g) = \{ \gamma \in C^r([0,1],\g) \: \gamma(0) = \gamma(1) =0 \} \]
(Exercise~\ref{exer:8.5.1}). 
We also write $P(G) := P^{(\infty)}(G)$ and $\Omega(G) := \Omega^{(\infty)}(G)$ 
for the smooth path, resp., loop group of~$G$ and 
$P(\g)$, resp., $\Omega(\g)$ for their Lie algebras.
\end{defn} 

\begin{lem} The Lie group morphism 
$\ev_1 \:  P^{(r)}(G) \to G$ has smooth local sections. In particular,
$P^{(r)}(G)$ is a Lie group extension of $G$ by $\Omega^{(r)}(G)$ and 
$G \cong P^{(r)}(G)/\Omega^{(r)}(G)$.
\end{lem}

\begin{prf} For any chart $\phi \: U \to G$, where $U \subeq \g$ 
is a convex open neighborhood of $0$ 
and $\phi(0) = \be$, the map 
$\sigma \: \phi(U) \to P^{(r)}(G), \sigma(\phi(x))(t) := \phi(tx)$ 
is a smooth section of $\ev_1$ (Exercise~\ref{exer:8.6.1.1}). 
\end{prf}

On the Lie algebra level 
$\{ \alpha \in P^{(r)}(\g) \!: (\forall t) \; \alpha(t) = t \,\alpha(1)\}$
is a natural vector space complement to
$\Omega^{(r)}(\g)$ in $P^{(r)}(\g)$, but 
this subspace is not a Lie subalgebra unless $\g$ is abelian. 

\begin{prop} \mlabel{prop:8.6.3} For every $r \in \N_0 \cup \{\infty\}$, the group 
$P^{(r)}(G)$ is contractible, and for $r= \infty$ it is smoothly contractible. 
\end{prop}

\begin{prf} For $r = \infty$, the map 
\[ H \: [0,1] \times P^{(r)}(G) \to P^{(r)}(G), \quad H(t,\gamma)(s) := \gamma(ts) \] 
is smooth. In fact, as $\ev$ is multiplicative in the second component, 
it suffices to verify smoothness in an identity neighborhood, where 
it reduces to the smoothness of the evaluation map 
$[0,1] \times P^{(r)}\g \to \g$, which 
in turn follows from Proposition~\ref{prop:smooth-eval}. 

For $r < \infty$, the continuity of $H$ follows similarly from the continuity 
of the evaluation map $[0,1] \times P^{(r)}(\g) \to \g$.
\end{prf}

\begin{rem} \mlabel{rem:homot-loop} 
By Proposition~\ref{prop:8.6.3}, $P^{(r)}(G)$ is simply connected,
so that the natural homomorphism 
$P(G)/\Omega(G)_0 \to G$ is a simply connected covering 
group (Corollary~\ref{cor10.1.2b}). It follows in particular that 
\begin{equation}
  \label{eq:pi0-pi1}
\pi_1(G) \cong \pi_0(\Omega(G))
\end{equation}
(cf.\ Subsection~\ref{subsec:14.8.1}).
More generally, the long exact homotopy sequence of the 
principal fiber bundle $\ev_1 \:  P^{(r)}(G) \to G$ with fiber 
$\Omega^{(r)}(G)$ (Theorem~\ref{homseq-princ}) leads to isomorphisms 
\[ \delta_k \: \pi_k(G) \to \pi_{k-1}(\Omega^{(r)}(G)) \quad \mbox{ for } 
\quad k \in \N, r \in \N_0 
\cup \{\infty\}.\] 
\end{rem}

\subsection{Enlargeability of path Lie algebras} 

In this subsection $\g$ is a locally exponential Lie algebra. 

\begin{lem} $P^{(r)}(\g)$ and $\Omega^{(r)}(\g)$ are locally exponential Lie algebras. 
\end{lem}

\begin{prf} For every compact manifold $M$ with boundary, 
the Lie algebra $C^r(M,\g)$ is locally exponential 
(Example~\ref{ex:6.2.5}).  
As $\ev_0 \: C^r([0,1],\g) \to \g, \xi \mapsto \xi(0)$ 
is a morphism of locally exponential Lie algebras, 
$P^{(r)}(\g) = \ker(\ev_0)$ is locally exponential 
(Proposition~\ref{prop:invim}). 
Likewise $\Omega^{(r)}(\g)$ is locally exponential 
because it is the kernel of 
the homomorphism $\ev_1 \: P^{(r)}(\g) \to \g$. 
\end{prf}

\index{Enlargeability Theorem for Path Lie Algebras} 
\begin{thm}[Enlargeability Theorem for Path Lie Algebras]  \mlabel{thm:enl-path}
If $\g$ is locally exponential and $r \geq 1$, then the locally exponential 
Lie algebra $P^{(r)}(\g)$ is enlargeable. Any $1$-connected locally 
exponential Lie group with Lie algebra $P^{(r)}(\g)$ is contractible 
and its center is simply connected 
\end{thm}

\begin{prf} (a) We consider the locally exponential 
Lie algebra $\fm := C^{r-1}([0,1],\g)$ (Exercise~\ref{exer:5.6.2}). 
According to Proposition~\ref{prop:semidir-Ad}, the Lie algebra 
\[ \hat\fm := |\fm| \rtimes \fm_{\rm ad} \] 
is locally exponential and enlargeable. A corresponding Lie group is 
the semidirect product $\hat M := |\fm| \rtimes M_{\rm ad}$, where $M_{\rm ad}$ is the $1$-connected
locally exponential Lie group with Lie algebra $\fm$ 
(Theorem~\ref{thm:5.3.8}). 
Now the adjoint action of the subalgebra $P^{(r)}(\g)$ on $\fm$ leads to a 
homomorphism 
\[ \alpha \: P^{(r)}(\g) \to \hat\fm = |\fm| \rtimes \fm_{\rm ad}, \quad 
\xi \mapsto (\xi', \ad \xi)  \] 
which is injective because $\xi' = 0$ implies that $\xi = 0$ because 
$\xi(0) = 0$. Next we observe that, for  
$\xi\in \fz(P^{(r)}(\g)) = P^{(r)}(\fz(\g))$, we have 
$\alpha(\xi) = (\xi',0)$, so that 
\[ \exp_{\hat M}(\alpha(\xi)) = (\xi',\be) = (0,\be) \] 
leads to $\xi = 0$. 
Therefore the Integral Subgroup Theorem~\ref{thm:5.5.3} 
implies the existence of a locally exponential Lie group 
with Lie algebra $P^{(r)}(\g)$. 
The center of this integral subgroup $\exp\big(P^{(r)}(\fz(\g))\big)$  
is contained in $|\fm| \times \{0\}$, hence simply connected. 

(b) Now let $H$ be a $1$-connected locally exponential Lie group 
with Lie algebra $P^{(r)}(\g)$. For each 
$t \in [0,1]$, we obtain a continuous endomorphism 
\[ f_t \: P^{(r)}(\g) \to P^{(r)}(\g), \quad 
f_t(\xi)(s) := \xi(ts).\] 
Let $F_t \: H \to H$ be the corresponding morphism of Lie groups 
with ${\L(F_t)=f_t}$. Then $F_0 = \be$ is constant and $F_1 = \id_H$. 
It remains to see that the map $F \: [0,1] \times H \to H, (t,h) \mapsto F_t(h)$ is continuous. 
Let $U \subeq \g$ be an open $0$-neighborhood mapped diffeomorphically 
to $\exp(U)$ by the exponential function. 
From 
\[ F_t(\exp \xi) = \exp(f_t(\xi))\] 
and the continuity of the map 
$[0,1] \times P^{(r)}(\g) \to P^{(r)}(\g), (t,\xi) \mapsto f_t(\xi)$ 
(Proposition~\ref{prop:8.6.3}), it follows that 
$F \res_{[0,1] \times \exp(U)}$ is continuous. 
As $F$ is multiplicative in the second argument and 
$H = \la \exp U \ra$, all the curves 
$t \mapsto F_t(h)$, $h \in H$, are continuous. 
In view of $F_t(h_0 h) = F_t(h_0) F_t(h)$, this implies that 
$F$ is continuous in every pair $(t_0, h_0) \in [0,1] \times H$. 
\end{prf}

In Proposition~\ref{prop:path-int} below we shall see that 
the continuous path Lie algebra 
$P^{(0)}(\g)$ is also enlargeable, but this requires different arguments, 
based on the fact that its adjoint group is contractible.

\begin{rem} Suppose that $\fz \subeq \fz(\g)$ is a central ideal 
for which $\fq := \g/\fz$ is enlargeable. 
If $\fz$ is not topologically split, then it is not clear that the natural 
map $P(\fg) \to P(\fq)$ is surjective. Nonetheless 
Corollary~\ref{cor:5.3.6} implies that the quotient 
$P(\g)/P(\fz)$ is locally exponential. 

If $P_\g$ is the $1$-connected Lie group with Lie algebra $P(\g)$, 
then its center is simply connected, so that 
$\exp \: P(\fz(\g)) \to Z(P_\g)_0$ is a homeomorphism. 
In particular $P_\fz := \exp(P(\fz))$ is a closed Lie subgroup and 
since $P(\g)/P(\fz)$ is locally exponential, 
$P_\g/P_\fz$ is a locally  exponential Lie group with Lie algebra 
$P(\g)/P(\fz)$ (Theorem~\ref{thm:5.5.11}). 

Since the contraction $F \: [0,1] \times P_\g\to P_\g$ 
is continuous and each $F_t$ preserves the subgroup $P_\fz$, 
we obtain a well-defined map 
\[ \oline F \: [0,1] \times (P_\g/P_\fz) \to P_\g/P_\fz, \quad 
\oline F_t([g]) := [F_t(g)].\]
Since the map $[0,1] \times P_\g \to [0,1] \times (P_\g/P_\fz)$ is open, its is a 
quotient map, and therefore $\oline F$ is continuous. 
This implies that $P_\g/P_\fz$ is also contractible. 
\end{rem}

\begin{rem} \mlabel{rem:pre-integrable}
If $G$ is a Lie group with Lie algebra $\g$, then 
$P(G)$ is a Lie group with Lie algebra $P(\g)$ 
(Definition~\ref{def4.4.1}).  
The logarithmic derivative 
$\delta \: P(G) \to C^\infty([0,1],\g)$ 
is a smooth map satisfying 
\[ \delta(\alpha\beta) = \delta(\beta) + \Ad(\beta)^{-1}.\delta(\alpha)
\quad \mbox{ and }\quad T_\be(\delta)(\xi) = \xi' \] 
 (Proposition~\ref{prop:smooth-logder0}). 
As $[\xi,\eta]' = [\xi', \eta] + [\xi,\eta']$, it follows that 
\[ T_\be(\delta) \: P(\g) \to C^\infty([0,1],\g), \quad \xi \mapsto \xi'\] 
becomes a topological isomorphism of Lie algebras if 
$C^\infty([0,1],\g)$ is endowed with the bracket 
\begin{equation}
  \label{eq:6.1.2}
[\xi,\eta](t) := \Big[\xi(t), \int_0^t \eta(\tau)\, d\tau\Big] + 
\Big[\int_0^t \xi(\tau)\, d\tau, \eta(t)\Big].
\end{equation}
The evaluation map $\ev_1 \: P(\g) \to \g$ corresponds to the quotient map 
$$ C^\infty([0,1],\g) \to \g, \quad \xi \mapsto \int_0^1 \xi(\tau)\, d\tau. $$

If, in addition, $G$ is regular, then $\delta$ is a diffeomorphism, 
and it follows that $C^\infty([0,1],\g)$, endowed with the bracket \eqref{eq:6.1.2}, is 
integrable. Since this property is clearly necessary for the regular integrability of 
$\g$, Lie algebras with this property are called 
{\it pre-integrable} \index{Lie algebra!pre-integrable} 
in \cite{RK97} (see also \cite{Les93}). 
\end{rem}

\subsection{The period group} 
\mlabel{subsec:per-grp}

Let $\g$ be a Mackey complete locally exponential Lie algebra, 
$\z := \z(\g)$ be its center and $\g_{\rm ad} = \g/\z$. 
We assume that $\g_{\rm ad}$ is also Mackey complete, which is in 
particular the case if $\g$ is Fr\'echet or if $\fz(\g)$ is a split subspace. 
Theorem~\ref{thm:5.3.8} 
implies that $\g_{\rm ad}$ is integrable to a connected Lie group 
$G_{\rm ad} \subeq \Aut(\g)$.
The central extension 
$$ \z \into \g \onto \g_{\rm ad} $$
can be pulled back via the evaluation map 
$\ev_1 \: P(\g_{\rm ad}) \to \g_{\rm ad}$ to a central extension 
\[  \z \into \hat P(\g) \onto P(\g_{\rm \ad}),   
\quad \hbox{ where } \quad 
\hat P(\g) := \{ (\xi,x) \in P(\g_{\rm ad}) \times \g :
\, \xi(1) = \ad x \}. \]
Then $\hat P(\g)$ also is locally exponential 
by Corollary~\ref{cor:equalizer}. 

In view of Proposition~\ref{prop:8.6.3}, 
the group $P(G_{\rm ad})$ is smoothly contractible. 
Its Lie algebra $P(\g_{\rm ad})$ is Mackey complete 
because $\g_{\rm ad}$ is Mackey complete 
(Exercise~\ref{exer:pathspace-mcomp})
and $\g$ is Mackey complete, so that the closed subspace 
$\hat P(\g) \subeq P(\g_{\rm ad}) \times \g$ is 
Mackey complete as well. Therefore 
Proposition~\ref{prop:VI.1.6b} implies that 
$\hat P(\g)$ is enlargeable (cf.~Example~\ref{exs:remD.3.2}). 
Let $\hat P_\g$ be a $1$-connected 
locally exponential Lie group with Lie algebra $\hat P(\g)$. 
In the proof of Proposition~\ref{prop:VI.1.6b}, 
a Lie group with Lie algebra $\hat P(\g)$ is constructed 
from the embedding 
\[ \hat P(\g) \into \g_1 := |\hat P(\g)| \oplus_b P(\g_{\rm ad}), \quad 
x \mapsto (x, q(x))\] 
and since the group integrating $\g_1$ is such that 
$|\hat P(\g)|$ is a simply connected subgroup, it follows in particular 
that $\exp_{\hat P_\g}\res_{\fz}$ is injective. We therefore have a central 
extension 
\begin{equation}
  \label{eq:path-center}
\fz \into \hat P_\g \sssmapright{q} P(G_{\rm ad}) 
\end{equation}
because $\hat P_\g/\fz$ is the unique simply connected locally 
exponential Lie group with Lie algebra $P(\g_{\rm ad})$. 
If the subspace $\fz \subeq \g$ is not split, 
this need not be a Lie group extension, i.e., smooth local 
sections need not exist. 

Consider the homomorphism
\[ \hat\Ad_1 := \ev_1 \circ q : \hat P_\g \to G_{\rm ad} \quad \mbox{ with } \quad 
\L(\hat\Ad_1)(\xi, x) := \xi(1).\] 
Then 
\[ \hat\Omega_\g := \ker \hat\Ad_1 = q^{-1}(\Omega(G_{\rm ad})) \] 
is a locally exponential Lie subgroup of $\hat P_\g$ with Lie algebra 
\begin{equation}
  \label{eq:8.27}
\hat\Omega(\g)  := \ker \L(\hat\Ad_1) = \Omega(\g_{\rm ad}) \times \z
\end{equation}
(Proposition~\ref{prop:5.5.9}), and 
\begin{equation}
  \label{eq:8.28}
\hat\Omega_\g/\fz \cong \Omega(G_{\rm ad}).\
\end{equation}
Let $q_\Omega \: \tilde\Omega(G_{\rm ad}) \to \Omega(G_{\rm ad})_0$ 
be the universal covering homomorphism of the identity component  
$\Omega(G_{\rm ad})_0$ with $\L(q_\Omega) = \id_{\Omega(\g_{\rm ad})}$. 
The inclusion homomorphism $\Omega(\g_{\rm ad}) \to \hat\Omega(\g)$ 
now integrates to a morphism of Lie groups 
\begin{equation}
  \label{eq:9.34}
f \: \tilde\Omega(G_{\rm ad}) \to \hat\Omega_\g\quad \mbox{ with } \quad 
q \circ f = q_\Omega.  
\end{equation}
In particular $f(\ker q_\Omega) \subeq \fz$. 
Identifying $\pi_1(\Omega(G_{\rm ad}))$ with $\ker q_\Omega$ 
(Theorem~\ref{thm:3.3.6}), we thus obtain 
a homomorphism 
\[ \per_\g := f\res_{\pi_1(\Omega(G_{\rm ad}))} \: \pi_1(\Omega(G_{\rm ad})) \to \fz.\] 

\begin{defn} \index{period homomorphism of Lie algebra $\g$, $\per_\g$} 
We call $\per_\g$ the {\it period homomorphism of $\g$} and its image
\[ \Pi(\g) := \im(\per_\g) \subeq \z(\g) \] 
the {\it period group of $\g$}. 
\end{defn}
\index{period homomorphism of Lie algebra $\g$, $\per_\g$} 
\index{period group of Lie algebra $\g$, $\Pi(\g)$} 

We recall that $\Omega(\g_{\rm ad}) \trile P(\g_{\rm ad})$ is the kernel of the evaluation 
morphism $P(\g_{\rm ad}) \to \g_{\rm ad}$, which is a morphism of locally exponential 
Lie algebras. Hence $\Omega(\g_{\rm ad})$ is locally exponential 
by Proposition~\ref{prop:invim}, 
so that the inclusion $\Omega(\g_{\rm ad}) \into \hat P(\g)$ 
integrates to an integral subgroup $N \into \hat P_\g$. The following 
lemma connects the behavior of this integral subgroup to the period group of~$\g$. 

\begin{lem} \mlabel{lem4.4.6} 
Let $N \subeq \hat P_\g$ be the connected integral subgroup
corresponding to the closed locally exponential subalgebra $\Omega(\g_{\rm ad})
\subeq \hat P(\g)$. Then $N \cap \, \z = \Pi(\g)$, and $N$ is a Lie
subgroup if and only if $\Pi(\g)$ is a discrete subgroup of~$\z$. 
\end{lem}

\begin{prf} The description of $(\hat\Omega_\g)_0$ as the quotient 
\[ \tilde\Omega(G_{\rm ad}) \times \z \to  (\hat\Omega_\g)_0, \quad 
(g,z) \mapsto f(g) z \] 
(see \eqref{eq:9.34}) shows that 
$N \cap \z \cong \im(\per_\g) = \Pi(\g)$ 
because $N$ is the image of $\tilde\Omega(G_{\rm ad})$ in~$\hat\Omega_\g$.

That the normal subgroup $N \subeq \hat P_\g$ is a Lie subgroup is
equivalent to $N$ being a Lie subgroup of $\hat\Omega_\g$. 
The Lie algebra $\hat\Omega(\g)$ is a direct product $\Omega(\g_{\rm
ad}) \times \z$. Therefore $N$ is a Lie subgroup if and only if 
there exists a \break {$0$-neighborhood} $U$ in $\z$ with 
$N \cap U = \{0\}$ (Lemma~\ref{lem:equiv-locexp-Liesub}), 
which is equivalent to $\Pi(\g)$ being discrete. 
\end{prf}

\begin{thm} [Enlargeability of locally exponential Lie algebras] \mlabel{thm4.4.7} 
Let $\g$ be a Mackey complete locally exponential Lie algebra 
for which $\g_{\rm ad}$ is also Mackey complete. Then $\g$ is enlargeable 
if and only if $\Pi(\g)$ is discrete. 
\end{thm}

\begin{prf} We have seen in the 
construction of $\Pi(\g)$ that there exists a group extension 
$\hat\Omega_\g \into \hat P_\g \onto G_{\rm ad},$ 
where $\hat P_\g$ is a simply connected locally exponential group with Lie algebra
$\hat P(\g)$.

If $\g$ is the Lie algebra of a locally exponential Lie group $G$, 
then we may assume w.l.o.g.\ that $G$ is $1$-connected. 
The $1$-connectedness of $\hat P_\g$ permits us to
integrate the homomorphism $\hat P(\g) \to \g$ of locally exponential Lie algebras 
to a Lie group homomorphism $p \! : \hat P_\g \to G$ with 
$$ \L(\ker p) = \ker \L(p) = \Omega(\g_{\rm ad}). $$
In view of Theorem~\ref{thm:int-thm}, we then have 
$G \cong \hat P_\g/\ker p$, where $\ker p$ is connected because $G$
is simply connected. Thus $\ker p$ coincides with the
connected intgral subgroup $N$ corresponding to
the Lie subalgebra $\fn := \Omega(\g_{\rm ad})$ of $\hat P(\g)$. 
In view of Lemma~\ref{lem4.4.6}, this implies that $\Pi(\g)$ is discrete.

If, conversely, $\Pi(\g)$ is discrete, then $N$ is a Lie subgroup,
and since $\hat P(\g)/\fn = \hat P(\g)/\Omega(\g_{\rm ad}) \cong \g$ is locally 
exponential, Theorem~\ref{thm:5.5.11}  implies that 
$\hat P_\g/N$ is a locally exponential Lie group with Lie algebra 
$\hat P(\g)/\fn \cong~\g$. 
\end{prf}

\begin{cor} \mlabel{cor:pi2-enl} 
If $\pi_2(G_{\rm ad})$ is trivial, then $\g$ is enlargeable. 
\end{cor}

\begin{prf} This follows from the isomorphism 
$\delta_2 \: \pi_2(G_{\rm ad}) \to \pi_1(\Omega(G_{\rm ad}))$ 
(Remark~\ref{rem:homot-loop}). 
\end{prf}

Since every Banach--Lie group is automatically locally exponential, 
the preceding theorem can be formulated in a simpler fashion for Banach--Lie algebras. 

\begin{cor} [Enlargeability Theorem for Banach--Lie algebras] \mlabel{cor4.4.7a} 
A Banach--Lie algebra $\g$ is integrable if and only if $\Pi(\g)$ is discrete. 
\end{cor}

\begin{cor}\mlabel{cor:lie3} 
{\rm(Lie's Third Theorem)}  Every finite-dimensional Lie algebra 
$\g$ is enlargeable.   
\end{cor}

\begin{prf} If $\g$ is finite-dimensional, then $G_{\rm ad}$ is also finite-dimensional, 
and therefore $\pi_2(G_{\rm ad})\cong \pi_1(\Omega(G_{\rm ad}))$ 
vanishes (Remark~\ref{rem1.3}). 
Hence the period group $\Pi(\g)$ is trivial and therefore discrete. 
\end{prf}

\begin{rem}
  If $\g$ is a finite-dimensional Lie algebra, then another tool to verify
  enlargeability is Ado's Theorem: There exists an $n \in \N$
  and  an embedding $\g \into \gl_n(\R)$. Then we can use
  the Integral Subgroup Theorem~\ref{thm:5.5.3}. 
We refer to \cite[Thm.~9.4.11]{HiNe12} for more details on this approach.
\end{rem}

\begin{lem} \mlabel{lem4.4.3} {\rm(Naturality of the period group)}
Let $\phi \! : \g \to \h$ be a homomorphism of locally exponential Lie 
algebras with $\phi(\z(\g)) \subeq \z(\h)$ and $\phi^{G_{\rm ad}} \:
\tilde G_{\rm ad} \to \tilde H_{\rm ad}$ the group homomorphism induced by $\phi$. Then 
\[ \phi(\Pi(\g)) \subeq \Pi(\h),\]  and the following diagram 
commutes: 
\begin{equation}
  \label{eq:diag1}
 \begin{matrix}  
\pi_2(G_{\rm ad}) \cong \pi_2(\tilde G_{\rm ad}) &
\smapright{\pi_2(\phi_{\rm ad}^G)} & \pi_2(H_{\rm ad}) \cong \pi_2(\tilde H_{\rm ad}) \cr 
\mapdown{\per_\g}  & {}  & \mapdown{\per_\h} \cr 
\z(\g)  & \smapright{\phi} & \z(\h). 
\end{matrix} 
\end{equation}
\end{lem}

\begin{prf} Since $\phi$ maps $\z(\g)$ to $\z(\h)$, it induces a
homomorphism 
\[ \phi_{\rm ad} \! : \g_{\rm ad} \to \h_{\rm ad}, \quad 
\ad x \mapsto \ad \phi(x) \] 
and hence a homomorphism 
$\hat P(\phi) \! : \hat P(\g) \to \hat P(\h), (\xi,x) \mapsto (\phi \circ \xi, 
\phi(x))$ with 
$\hat P(\phi)\big(\hat\Omega(\g)\big) \subeq \hat\Omega(\h)$. 
Integration to the simply connected group $\hat P_\g$ further leads to
a group homomorphism 
$$ \hat P_\phi \! : \hat P_\g \to \hat P_\fh 
\quad \hbox{ with } \quad 
\L(\hat P_\phi) = \hat P(\phi). $$
It is clear that this homomorphism maps the subgroup 
$\hat \Omega_\g$ to $\hat \Omega_\fh$, hence induces a homomorphism 
$\pi_1(\hat\Omega_\g) \to \pi_1(\hat\Omega_\fh)$. 
This means that the induced map 
\[  (\hat\Omega_\g)\,\tilde{}\ \cong \tilde\Omega(G_{\rm ad}) \times \z(\g) \to  
\tilde\Omega(H_{\rm ad}) \times \z(\h) \cong  (\hat\Omega_\fh)\,\tilde{} \]
of the simply connected covering groups 
maps the graph of $\per_\g$ into the graph of $\per_\h$. This translates into 
\begin{equation}
  \label{eq:4.4.1}
\phi\res_{\z(\g)} \circ \per_\g
= \per_\h \circ \, \pi_1\big(\Omega(\phi_{\rm ad}^G)\big), 
\end{equation}
where $\phi_{\rm ad}^G \! : \tilde G_{\rm ad} \to \tilde H_{\rm ad}$ is the
homomorphism with $\L(\phi_{\rm ad}^G) = \phi_{\rm ad}$, and 
$\Omega(\phi_{\rm ad}^G)$ is the corresponding homomorphism 
\[ \Omega(\tilde G_{\rm ad}) \cong \Omega(G_{\rm ad})_0 \to 
\Omega(\tilde H_{\rm ad}) \cong \Omega(H_{\rm ad})_0, \quad 
\gamma \mapsto \phi_{\rm ad}^G \circ \gamma\] 
(cf.\ Remark~\ref{rem:homot-loop}). 
Next we observe that the naturality of the connecting homomorphism 
$\delta_2^\g \: \pi_2(G_{\rm ad}) \to \pi_1(\Omega(G_{\rm \ad}))$ 
(Remark~\ref{rem:homot-loop}) 
leads to the commutative diagram 
\[  \begin{matrix} 
 \pi_2(G_{\rm ad}) & \mapright{\pi_2(\phi_{\rm ad}^G)} & \pi_2(H_{\rm ad})  \\ 
\mapdown{\delta_2^\g} & & \mapdown{\delta_2^\fh}  \\ 
 \pi_1(\Omega(G_{\rm ad})) & \mapright{\pi_1(\Omega(\phi_{\rm ad}^G))} & 
\pi_1(\Omega(H_{\rm ad})) \\ 
\mapdown{\per_\g} & & \mapdown{\per_\fh}  \\ 
 \fz(\g)  & \mapright{\phi} & \fz(\fh). \end{matrix}\] 
Identifying 
$\pi_2(G_{\rm ad})$ with $\pi_1(\Omega (G_{\rm ad}))$ 
via $\delta_2^G$ now implies the commutativity of the diagram 
\eqref{eq:diag1} and hence in particular that $\phi(\Pi(\g)) \subeq \Pi(\h)$. 
\end{prf}

\begin{cor} \mlabel{cor4.4.4} 
  \begin{description}
  \item[\rm(a)] If $\phi \!: \g \to \h$ is a homomorphism of locally exponential  
Lie algebras with $\phi(\z(\g)) \subeq \z(\h)$ for which the induced
homomorphism $\pi_2(G_{\rm ad}) \to \pi_2(H_{\rm ad})$ is surjective, then 
$\phi(\Pi(\g))= \Pi(\h)$. 
  \item[\rm(b)] If $\phi \!: \g \to \h$ is a quotient homomorphism of
locally exponential Lie algebras with $\phi(\z(\g)) = \z(\h)$  and 
$\ker \phi \subeq \z(\g)$, then 
$\phi(\Pi(\g))= \Pi(\h)$. 
  \end{description}
\end{cor}

\begin{prf} (a) This is an immediate consequence of Lemma~\ref{lem4.4.3}.

\par\nin (b) Our assumption implies that 
$\g_{\rm ad} = \g/\z(\g) \cong \h/\phi(\z(\g)) \cong \h_{\rm ad}.$
Therefore the induced map $\phi_{\rm ad}^G \: G_{\rm ad} \to H_{\rm ad}$ is 
an isomorphism, and the assertion follows from~(a).
\end{prf}

\begin{cor}
  \mlabel{cor:VI.1.13} If $\g_1, \g_2$ are Mackey complete 
locally exponential Lie algebras whose adjoint algebras are also 
Mackey complete, then 
\[ \Pi(\g_1 \times \g_2) = \Pi(\g_1) \times \Pi(\g_2). \] 
\end{cor} 

\begin{rem} \mlabel{rem4.4.5} 
Let 
${\bf LZ}$ denote the category whose objects are 
Mackey complete locally exponential Lie algebras whose adjoint quotients 
are also Mackey complete, and whose morphisms are continuous Lie algebra
homomorphisms  mapping center to center. Then Lemma~\ref{lem4.4.3} 
shows that $\Pi \: \g \mapsto \Pi(\g)$ defines 
a functor from ${\bf LZ}$ to the category of abelian
topological groups. 
\end{rem}

\begin{prop} \mlabel{prop4.4.8} 
If $G$ is a simply connected locally exponential Lie group with Lie
algebra $\g$ for which $\g$ and $\g_{\rm ad}$ are Mackey complete, then 
$$Z(G)_0 \cong \z/ \Pi(\g) \quad \hbox{ and } \quad 
\pi_1(Z(G)) \cong \Pi(\g) = \ker(\exp_G\res_\z). $$
\end{prop}

\begin{prf} As in the proof of Theorem~\ref{thm4.4.7},
we write $G$ as $\hat P_\g/\ker p$. Since $G$ is simply connected, the group $\ker p$ is
connected. Moreover, $\z(\g) \cong \hat\Omega(\g)/\Omega(\g_{\rm ad})$
implies that 
\[ Z(G)_0 \cong \hat\Omega_\g/\ker p, \] 
so that 
$\hat\Omega_\g\cong (\tilde\Omega(G_{\rm ad}) \times \z)/ \Gamma(-\per_\g)$ 
yields $Z(G)_0 \cong \z/\im(\per_\g) =
\z/\Pi(\g)$. 
\end{prf}

\begin{rem}
  \mlabel{rem:VI.1.14} (Non-enlargeable Lie algebras) 
Suppose that $\g$ is a Mackey complete locally exponential Lie algebra 
with finite-dimensional center  with 
$\Pi(\g) = \Z d$ and $d \not=0$. Let $\theta \in \R \setminus \Q$. 
Then $\z := \{ (x,\theta x) \: x \in \z(\g)\}$ is a central ideal of 
$\g \times \g$, so that $\h := (\g \times \g)/\z$ is locally exponential 
(Corollary~\ref{cor:5.3.6}). 
Note that $\fh$ is Mackey complete since $\g \times \g$ is so, and $\fz$ is 
finite-dimensional. Corollary~\ref{cor:VI.1.13} and Lemma~\ref{lem4.4.3} 
now imply that 
\[ \Pi(\h) \cong \Z[(d,0)] + \Z[(0,d)] = (\Z + \Z \theta) [(d,0)], \] 
which is not discrete. Hence $\fh$ is not enlargeable. 
\end{rem} 

\begin{ex}
  \mlabel{ex:VI.1.16} (Non-enlargeable Banach--Lie algebras) 
(a) The first example of a non-enlargeable Banach--Lie algebra 
was given by van Est and Korthagen with the method described in 
Remark~\ref{rem:VI.1.14} (\cite{EK64}). 
It is the central extension 
$\g = \R \oplus_\omega C^1(\bS^1, \su_2(\C))$ of the Banach--Lie algebra 
$C^1(\bS^1, \su_2(\C))$ defined by the cocycle 
$$ \omega(f,g) := \int_0^1 \tr(f(t)g'(t))\, dt, $$
where we identify functions on $\bS^1 \cong \R/\Z$ with $1$-periodic functions on $\R$. 
Then $\g_{\rm ad} \cong C^1(\bS^1, \su_2(\C))$ and 
$G_{\rm ad} \cong C^1(\bS^1, \SU_2(\C))$ leads to 
\[ \pi_2(G_{\rm ad}) \cong \pi_3(\SU_2(\C)) \cong \pi_3(\bS^3) \cong \Z.\]  
Now one shows that $\per_\g = \per_\omega$ is non-trivial to verify that 
$\Pi(\g) \cong \Z$ (see also \cite{MN03} for this calculation). 

(b) Douady and Lazard found a simpler example (\cite{DL66}) 
by observing that the $1$-connectedness of the unitary group $\U(\cH)$ of an 
infinite-dimensional complex Hilbert space $\cH$ 
(Kuiper's Theorem~\ref{kuiper-thm}) 
implies that its Lie algebra 
$\fu(\cH) := \{ X \in {\cal L}(\cH) \: X^* = - X\}$ satisfies 
$$ \Pi(\fu(\cH)) \cong \pi_1(Z(\U(\cH))) = \pi_1(\T) \cong \Z $$
(Proposition~\ref{prop4.4.8}). 

(c) Based on the fact that $\U(\cH)$ is $1$-connected, one can also give the following 
direct argument. For any irrational $\theta \in \R \setminus \Q$, the line 
$\fn := \R i (\1,\theta \1)$ generates a dense subgroup of the center 
$Z(\U(\cH) \times \U(\cH)) \cong \T^2$ of the \break 
{$1$-connected} group $\U(\cH) \times \U(\cH)$, 
so that Theorem~\ref{thm:quot-alg}  implies that the quotient Lie algebra 
$\fq := (\fu(\cH) \times \fu(\cH))/\fn$ is not enlargeable. 
As in Remark~\ref{rem:VI.1.14}, one can show that, 
\[ \Pi(\fq) = (\Z + \Z \theta) \subeq \fz(\fq) \cong \R, \] 
which is not discrete. 
\end{ex}

\begin{rem} With the refined methods developed in \cite{Ne02a} for 
central extensions of infinite-dimensional Lie groups, 
one obtains for every Mackey complete 
locally exponential Lie algebra $\g$ a period map 
\[ \per_b \: \pi_2(G_{\rm ad}) \to |\g| \]
for the cocycle $b \: \g_{\rm ad} \times \g_{\rm ad} \to \g$ defined by 
$b(\ad x,\ad y) := [x,y]$. This does not require the Mackey completeness 
of $\g_{\rm ad}$. This period map is obtained by integration 
of the left invariant closed $2$-form $\omega \in \Omega^2(G_{\rm ad}, |\g|)$ 
with $\omega_\be = b$. 

As $\ad \circ b$ is the Lie bracket of $\g_{\rm ad}$, it is the  
Lie algebra coboundary $- \dd_\g(\id_{\g_{\rm ad}})$ 
(cf.~Appendix~\ref{app:liealg-cohom}). 
Hence the corresponding $2$-form 
$\ad \circ \omega$ is exact, so that its periods vanish by Stoke's Theorem. 
This implies that $\im(\per_b) \subeq \ker \ad = \z(\g)$, so that 
we obtain a homomorphism 
\[ \per_b \: \pi_2(G_{\rm ad}) \to \fz(\g) \] 
without the assumption that $\g_{\rm ad}$ is Mackey complete. 
\end{rem}

\begin{prop} \mlabel{prop:path-int} 
Let $\g$ be a Mackey complete locally exponential Lie algebra 
whose adjoint algebra is also Mackey complete. 
Then the continuous path algebra $P^{(0)}(\g)$ is enlargeable. 
In particular, $P^{(0)}(\g)$ is enlargeable for every Fr\'echet--Lie algebra. 
\end{prop}

\begin{prf} {\bf First proof:} First we observe that the path Lie algebra 
$P^{(0)}(\g_{\rm ad})$ is the Lie algebra of the contractible 
locally exponential group $P(G_{\rm ad})$ 
(Proposition~\ref{prop:8.6.3}). 
Next we show that 
\[ q := P_{\rm ad} \: P^{(0)}(\g) \to P^{(0)}(\g_{\rm ad}), \quad 
\xi \mapsto \ad \circ \xi \] 
is a generalized central extension (Definition~\ref{def:F.3.1}). 
Since the polygonal paths in
 $P^{(0)}(\g_{\rm ad})$ lift to elements in $P^{(0)}(\g)$, the image of 
$q$ is dense. Let 
\[ b_\g \: \g_{\rm ad} \times \g_{\rm ad} \to \g, \quad 
(\ad x, \ad y) \mapsto [x,y] \] 
be the canonical cocycle of the generalized central extension 
$\ad \:  \g \to \g_{\rm ad}$. Then 
\[ b \: P^{(0)}(\g_{\rm ad}) \times P^{(0)}(\g_{\rm ad}) \to 
P^{(0)}(\g), \quad 
b(\xi_1, \xi_2)(t) := b_\g([\xi_1(t), \xi_2(t)]) \] 
is a continuous bilinear map with 
\[ b(q(\xi_1), q(\xi_2)) = [\xi_1, \xi_2] \quad \mbox{ for } \quad 
\xi_1, \xi_2 \in P^{(0)}(\g). \] 
This shows that $q$ defines a generalized central extension. 

The Mackey completeness of $\g$ implies the 
Mackey completeness of $P^{(0)}(\g)$ (Exercise~\ref{exer:pathspace-mcomp}). 
Therefore we derive the enlargeability of $P^{(0)}(\g)$ 
from Remark~\ref{rem:8.3.15}.\begin{footnote}
{Note that Proposition~\ref{prop7.6.2} does not apply here because 
we do not know if the cohomology spaces $H^k_{\rm dR}(P^{(0)}(G_{\rm ad}),|P^{(0)}(\g)|)$, 
$k =1,2$, vanish. This would follow if $P^{(0)}(G_{\rm ad})$ is smoothly 
contractible, but we only know that it is topologically contractible. }  
\end{footnote}

Since $P^{(0)}(\g)$ naturally embeds into 
$\fh \cong |P^{(0)}(\g)| \oplus_b P^{(0)}(\g_{\rm ad})$ 
(Proposition~\ref{propD.3.5}), it follows 
from the Integral Subgroup Theorem~\ref{thm:5.5.3} that 
$P^{(0)}(\g)$ is enlargeable. 

\nin{\bf Second proof:} An alternative proof of the enlargibility 
of $\fp := P^{(0)}(\g)$ can be derived from Corollary~\ref{cor:pi2-enl} 
by showing that $\pi_2(P_{\rm ad})$ vanishes. 
We have $\fp_{\rm ad} = P^{(0)}(\g)/P^{(0)}(\fz(\g))$ which is the subalgebra 
of $P^{(0)}(\g_{\rm ad})$ consisting of all elements with a continuous lift to~$\g$. 
Let 
\[ F \: [0,1] \times P^{(0)}(\g_{\rm ad}) \to P^{(0)}(\g_{\rm ad}), \quad 
F_t(\xi)(s) := \xi(ts) \] 
denote the canonical contraction. 
Then each $F_t$ maps $\fp_{\rm ad}$ into itself, so that we obtain a map 
\[ \oline F\: [0,1] \times \fp_{\rm ad} \to \fp_{\rm ad}, \quad 
\oline F_t (\ad \xi)(s) = \ad(\xi(st)).   \] 
This map is continuous because $[0,1] \times \fp \to [0,1]\to \fp_{\rm ad}$ 
is an open map and the map 
\[[0,1] \times \fp \to \fp_{\rm ad}, \quad 
(t,\xi) \mapsto (s \mapsto \ad(\xi(st))) \] 
is continuous. Therefore Exercise~\ref{exer:8.6.1.3} implies 
that the simply connected covering group $\tilde P_{\rm ad}$ of 
$P_{\rm ad}$ is contractible. Hence the group 
$\pi_2(P_{\rm ad}) \cong \pi_2(\tilde P_{\rm ad})$ 
(cf.~Lemma~\ref{lem10.1.2}) is trivial. 
\end{prf}

\begin{rem} \mlabel{rem:periods-r finite}
With Proposition~\ref{prop:path-int}  and Theorem~\ref{thm:enl-path}, 
we know that, 
for $\g$ and $\g_{\rm ad}$ Mackey complete, all path algebras 
$P^{(r)}(\g)$ are enlargeable. 
Since the groups $P^{(r)}(G_{\rm ad})$ are contractible, 
we can use Remark~\ref{rem:8.3.15} to 
see that the corresponding central extensions 
$\hat P^{(r)}(\g)$ are enlargeable to central extensions 
\[ \fz \into \hat P^{(r)}_\g \onto P^{(r)}(G_{\rm ad}).\] 
We thus obtain natural homomorphisms 
\[ f^{(r)} \: \tilde\Omega^{(r)}(G_{\rm ad}) \to \hat\Omega_\g^{(r)} 
= \ker(\hat\Ad_1^{(r)})\] 
which lead to period homomorphisms 
\[ \per_\g^{(r)} \: \pi_1(\Omega^{(r)}(G_{\rm ad}))  \to \fz.\] 
Since the inclusion maps 
$\Omega(G_{\rm ad}) \into \Omega^{(r)}(G_{\rm ad})$ 
are weak homotopy equivalences (Theorem~\ref{thma.3}), 
we obtain isomorphisms 
\[ \eta_r \: \pi_1(\Omega(G_{\rm ad})) \to \pi_1(\Omega^{(r)}(G_{\rm ad})) 
\quad \mbox{ with } \quad 
\per_\g^{(r)} \circ \eta_r = \per_\g, \] 
so that all period maps $\per_\g^{(r)}$ 
have the same range $\Pi(\g)$. 
\end{rem}

\subsection{A period group for non-Mackey complete algebras} 

In this subsection we also introduce a period group
$\Pi'(\g)\subeq \fz(\g)$ of a locally exponential Lie algebra~$\g$ 
which does not have to be Mackey complete and which can also be used 
to characterize enlargeability by its discreteness. 
If $\g$ is Fr\'echet or $\fz(\g)$ is topologically split, 
then both period groups $\Pi(\g)$ and $\Pi'(\g)$ coincide 
(Proposition~\ref{prop:periods-same}). 
This alternative approach is of interest because, for a locally exponential 
Lie algebra $\g$, it is not clear if the completion or the Mackey 
completion are locally exponential as well. 

We start with the smooth path Lie algebra 
$P(\g) := P^\infty(\g)$ 
and recall from Theorem~\ref{thm:enl-path} that $P(\g)$ is an enlargible 
locally exponential Lie algebra. Let $P_\g$ be a corresponding $1$-connected 
locally exponential Lie group. Then the evaluation map 
$\ev_1 \: P(\g) \to \g, \xi \mapsto \xi(1)$ is a quotient homomorphism 
with kernel $\Omega(\g)$, so that $\g \cong P(\g)/\Omega(\g)$. 
Now Proposition~\ref{prop:int-quot-locexp} implies that 
$\g$ is enlargible if and only if the integral subgroup 
$\la \exp \Omega(\g) \ra \subeq P_\g$ is a locally exponential Lie subgroup. 
We want to use the Enlargeability Theorem for Quotients~\ref{thm:quot-alg} to get a 
more direct characterization of the enlargeability of~$\g$. 

The action of $P(G)$ on $\g$ corresponds to a morphims 
\[ \Ad_1 \: P_\g \to G_{\rm ad} \quad \mbox{ with } \quad 
\L(\Ad_1)  \: P(\g) \to \g_{\rm ad}, \quad \xi \mapsto \ad(\xi(1)). \]
We consider the kernel 
\[ N_\g := \ker(\Ad_1) 
\quad \mbox{ with Lie algebra } \quad 
\fn_\g := \{ \xi  \in P(\g) \: \xi(1) \in \fz(\g) \}.\] 
Since $\L(\Ad_1)$ is a quotient map, $\g_{\rm ad} \cong P(\g)/\fn_\g.$ 
Accordingly $N_\g$ is a closed normal subgroup with 
$G_{\rm ad} \cong P_\g/N_\g$. Unfortunately, we don't know without any further 
assumptions on $\g$ whether $\Ad_1$ has continuous local sections. 
The problem is the existence of continuous local sections for the 
linear quotient map $\g \to \g_{\rm ad}$. Now 
\[ \ev_1 \:  \fn_\g \to \fz(\g),\quad \xi \mapsto \xi(1) \] 
is a homomorphism of Lie algebras, i.e., a continuous linear map vanishing on all brackets. This leads to a period homomorphism 
\[ \per_\g' \: \pi_1(N_\g) \to \fz(\g), \quad 
[\gamma] \mapsto \int_0^1 \delta(\gamma)_t(1)\, dt,\] 
where $\gamma \in \Omega(N_\g)$ is a smooth loop. 

\begin{defn} We call the additive subgroup
  \[ \Pi'(\g) := \per_\g'(\pi_1(N_\g)) \subeq \fz(\g) \]  
the {\it second period group of $\g$}. \index{second period group of $\g$} 
\end{defn}

From Theorem~\ref{thm:quot-alg} we  obtain: 
\begin{thm} \mlabel{thm4.4.7b} 
The locally exponential Lie algebra $\g$ is enlargeable 
if and only if $\Pi'(\g)$ is discrete. 
\end{thm}

In the remainder of this subsection we discuss the naturality 
of $\Pi'(\g)$ and the coincidence with $\Pi(\g)$ if both are defined. 

\begin{defn} We call $\g$ 
\index{Lie algebra!topologically $\ad$-split} 
{\it topologically $\ad$-split} if 
there exists a \break {$0$-neighorhood} $U \subeq \g_{\rm ad}$ and a continuous 
map $\sigma \: U \to \g$ with $\ad \circ \sigma = \id_U$. 
\end{defn}

\begin{lem} If $\g$ is Fr\'echet or $\fz(\g)$ is a split subspace of $\g$,
then $\g$ is topologically $\ad$-split. 
\end{lem}

\begin{prf} That $\fz(\g)$ is a split subspace means that there exists a 
continuous linear map $\sigma \: \g_{\rm ad} \to \g$ with 
$\ad \circ \sigma = \id_{\g_{\rm ad}}$. 

If $\g$ is Fr\'echet, then the existence of continuous local sections follows 
from Michael's Selection Theorem; see \cite{MicE59} for the Banach case and 
\cite[Ch.\,II, \S4.7, Prop.\,12]{Bou87}
for the general~result.
\end{prf}

The following proposition draws some consequences. 

\begin{prop} \mlabel{prop:iso-homtop-grp}
Suppose that $\g$ is topologically $\ad$-split. 
Then the homomorphism $\Ad_1 \: P_\g \to G_{\rm ad}$ has continuous local sections, 
so that we obtain on $P_\g$ the structure of a locally trivial $N_\g$-principal 
bundle over $G_{\rm ad}$. The corresponding connecting homomorphism 
\[ \delta_2 \: \pi_2(G_{\rm ad}) \to \pi_1(N_\g) \] 
is an isomorphism. 
\end{prop}

\begin{prf} Since $\ev_1 \: P(\g) \to \g$ has a linear continuous section, 
our assumption implies that 
$\ad_1 \: P(\g) \to \g_{\rm ad}$ has continuous local sections, and 
since $P_\g$ and $G_{\rm ad}$ are locally exponential, the homomorphism 
$\Ad_1$ has continuous local sections. Therefore $P_\g$ is a principal 
$N_\g$-bundle, so that we have an exact sequence 
\[ \cdots \to \pi_2(N_\g) \to \pi_2(P_\g) \to \pi_2(G_{\rm ad}) 
\sssmapright{\delta_2} \pi_1(N_\g) \to \cdots \] 
(Theorem~\ref{homseq-princ}). As $P_\g$ is contractible 
(Theorem~\ref{thm:enl-path}), all its homotopy 
groups are trivial, and therefore $\delta_2$ is an isomorphism. 
\end{prf}

\begin{prop} \mlabel{prop:periods-same} 
If $\g$ is Fr\'echet or Mackey complete with $\fz(\g)$ split, 
then 
\[ \per_\g' \circ \delta_2 = \per_\g \: \pi_2(G_{\rm ad}) \to \fz.\] 
In particular, $\Pi(\g)$ coincides with $\Pi'(\g)$. 
\end{prop}

\begin{prf} Our assumptions imply that $\g$ and $\g_{\rm ad}$ are 
Mackey complete, so that 
$\per_\g \: \pi_1(N_\g) \to \fz(\g)$ is defined. 

\nin(a): From Remark~\ref{rem:homot-loop} we know that 
the connecting  homomorphism 
\[  \delta_2^\g \: \pi_2(G_{\rm ad}) \to \pi_1(\Omega(G_{\rm ad})) \] 
is an isomorphism. We recall that 
$\per_\g \: \pi_1(\Omega(G_{\rm ad})) \to \fz(\g)$ was obtained 
by restring the natural homomorphism 
\[ f \: \tilde\Omega(G_{\rm ad}) \to \hat\Omega_\g = \ker(\Ad_1)  
\subeq \hat P_\g\] 
integrating the inclusion 
$\Omega(\g_{\rm ad}) \into \hat\Omega(\g) \cong \Omega(\g_{\rm ad}) \times \fz$. 

In Remark~\ref{rem:periods-r finite} we have seen how to obtain period 
homomorphisms $\per_\g^{(r)} \: \pi_1(\Omega^{(r)}(G_{\rm ad})) \to \fz(\g)$ 
and natural isomorphism 
\[ \eta_r \: \pi_1(\Omega)(G_{\rm ad})) \to \pi_1(\Omega^{(r)}(G_{\rm ad})) \] 
satisfying $\per_\g^{(r)} \circ \eta_r = \per_\g$. 

\nin (b): The other period homomorphism $\per_\g' \: \pi_1(N_\g) \to \g$ is 
obtained from the natural homomorphism 
$f' \: \tilde{N_{\g,0}} \to \fz(\g)$ by restriction  to $\pi_1(N_\g)$. 

For $r < \infty$, our assumptions imply that the natural map 
\[ P^{(r)}(\g) \to P^{(r)}(\g_{\rm ad}) \] 
is a quotient map with kernel $P^{(r)}(\fz(\g))$ 
(Exercise~\ref{exer:pathspace-mcomp-2}). Accordingly, we get a 
quotient map 
\[ q^{(r)} \: P^{(r)}(\g) \to \hat P^{(r)}(\g), \quad \xi \mapsto (\ad \xi, \xi(1)) \] 
with kernel $\Omega^{(r)}(\fz(\g))$, so that 
$\hat P^{(r)}(\g) \cong P^{(r)}(\g)/\Omega^{(r)}(\fz(\g))$ 
and, for the corresponding locally exponential groups, 
\[ q^{(r)}_G \: P_\g^{(r)} \to \hat P^{(r)}_\g\quad \mbox{ yields } \quad 
\hat P^{(r)}_\g \cong P_\g^{(r)}/\Omega^{(r)}(\fz(\g)).\] 
If $\fz(\g)$ is split, then this even holds for $r = \infty$. 

The homomorphisms 
$\ad_1^{(r)} \:  P^{(r)}(\g) \to \g_{\rm ad}$
and $\hat\ad_1^{(r)} \:  \hat P^{(r)}(\g) \to \g_{\rm ad}$ 
satisfy $\hat\ad_1^{(r)} \circ q_r = \ad_1^{(r)}$. 
On the group level this leads to a commutative diagram with exact rows and columns
\[  \begin{matrix} 
\Omega^{(r)}(\fz(\g))  & \ssmapright{\id} & \Omega^{(r)}(\fz(\g))  && \\ 
\mapdown{}  & & \mapdown{} & &  \\ 
 N_\g^{(r)} & \ssmapright{i_1} & P_\g^{(r)} &  \smapright{\Ad_1^{(r)}}  & G_{\rm ad} \\ 
\mapdown{}  & & \mapdown{q_G^{(r)}} & & \mapdown{\id_{G_{\rm ad}}} \\ 
 \hat\Omega_\g^{(r)} &  \ssmapright{i_2} & \hat P_\g^{(r)} &  \smapright{\hat\Ad_1^{(r)}}  
& G_{\rm ad}. \end{matrix}\] 
Here the leftmost column is a topologically split quotient morphism 
with contractible kernel, hence it induces isomorphisms 
\[ \pi_1(N_\g^{(r)}) \to \pi_1(\hat\Omega_\g^{(r)}) \cong \pi_1(\Omega^{(r)}(G_{\rm ad}))\] 
(Theorem~\ref{homseq-princ}). 
Next we note that the diagram 
\[  \begin{matrix} 
\fn_\g^{(r)} & \smapright{\ev_1} & \fz(\g) \\ 
\mapdown{q^{(r)}}  & & \mapdown{\id_{\fz(\g)}}  \\ 
 \hat\Omega^{(r)}(\g) 
&  \smapright{\pr_{\fz(\g)}}  & \fz(\g) \end{matrix}\] 
commutes, and this implies that 
\[ \per^{(r)}_\g \circ \pi_1(q_G^{(r)}) = \per'^{(r)}_\g.\] 
Combined with  (a) for the passage between $r$ finite and infinite, 
the assertion now follows. 
\end{prf}


\begin{lem} \mlabel{lem4.4.3b} 
Let $\phi : \g \to \h$ be a homomorphism of locally exponential Lie 
algebras with $\phi(\fz(\g)) \subeq \fz(\fh)$. 
Then $\phi$ induces by composition 
a homomorphism $P(\g) \to P(\fh)$ and further a morphism of 
Lie groups $P_\phi \: P_\g \to P_\fh$ which restricts to a continuous 
homomorphism 
$\phi_N  \: (N_\g)_0 \to (N_\fh)_0$. 
Then $\phi(\Pi'(\g)) \subeq \Pi'(\h)$ and the following diagram 
commutes: 
\[  \begin{matrix}  
\pi_1(N_\g) &
\smapright{\pi_1(\phi_N)} & \pi_1(N_\fh) \cr 
\mapdown{\per_\g'}  & {}  & \mapdown{\per_\h'} \cr 
\z(\g)  & \smapright{\phi} & \z(\h). 
\end{matrix}  \]
\end{lem}

\begin{prf} The homomorphism 
\[ P(\phi) \: P(\g) \to P(\fh), \quad \xi \mapsto \phi \circ \xi \] 
integrates to a group homomorphism 
$P_\phi \: P_\g \to P_\fh.$ 
As 
$\phi(\fz(\g)) \subeq \fz(\fh)$, 
$\phi$ induces a continuous homomorphism 
$\phi^{\rm ad} \: \g_{\rm ad} \to \fh_{\rm ad}, 
\ad x \mapsto \ad \phi(x)$; let $\phi^{\rm ad}_G \: \tilde G_{\rm ad} \to 
\tilde H_{\rm ad}$ be the corresponding morphism of Lie groups. 
Let 
\[ \ad^\g_1 \: P(\g) \to \g_{\rm ad}, \qquad \xi \mapsto \ad(\xi(1)) 
\quad \mbox{ and } \quad 
\Ad_1^\g \: P_\g \to G_{\rm ad} \] 
be the corresponding group homomorphism. We then have 
\[ \ad_1^\fh \big(P(\phi)\xi\big) = \ad(\phi(\xi(1))) 
= \phi^{\rm ad} \ad(\xi(1)) 
= \phi^{\rm ad} \circ \ad_1^\g(\xi)\quad \mbox{ for } \quad \xi \in P(\g),\] 
i.e., 
\[ \ad_1^\fh  \circ P(\phi) = \phi^{\rm ad} \circ \ad_1^\g \: 
P(\g) \to \fh_{\rm ad}.\] 
On the groups level this leads to 
\[ \tilde\Ad_1 \circ P_\phi = \phi_G^{\rm ad} \circ \tilde\Ad_1^G \: 
P_\g \to \tilde H_{\rm ad}.\] 
We conculde that 
$P_G(\phi)((N_G)_0) = P_\phi(\ker \tilde\Ad_1^G)_0  \subeq 
(\ker \tilde\Ad_1^H)_0 = (N_\fh)_0$, so that 
we obtain a continuous morphism 
$\phi_N \: (N_\g)_0 \to (N_\fh)_0$. 
Further, $\ev_1^\g \: \fn_\g \to \fz(\g)$ satisfies 
\[ \phi \circ \ev_1^\g = \ev_1^\fh \circ P(\phi)\res_{\fn_\g},\] 
which leads to 
\[ \phi \circ \per_\g' = \per_\fh' \circ \pi_1(\phi_N).\] 
This relation shows in particular that 
$\phi(\Pi'(\g)) \subeq \Pi'(\fh)$. 
\end{prf}

\begin{small}
\subsection*{Exercises for Section~\ref{sec:11.4}}

\begin{exer}\mlabel{exer:pathspace-mcomp} 
Let $E$ be a Mackey complete space and 
$M$ be a compact manifold (possibly with corners). Show that, 
for $r \in \N_0 \cup \{\infty\}$, the space 
$C^r(M,E)$ is also Mackey complete. \\ 
Hint: For a smooth curve $\gamma \: [0,1] \to C^r(M,E)$, let 
$\Gamma(m) := \int_0^1 \gamma(t)(m)\, dt$ and show that 
$\Gamma \in C^r(M,E)$ is a weak integral of $\gamma$. 
\end{exer}

\begin{exer}\mlabel{exer:pathspace-mcomp-2} 
Let $E$ be a Fr\'echet space, $F \subeq E$ be a closed subspace and 
${q \: E \to E/F}$ the quotient map. 
Show that, for every $r < \infty$, the map 
\[ P^{(r)}(q) \: P^{(r)}(E) \to P^{(r)}(E/F), \quad 
\gamma \mapsto q \circ \gamma \] 
is a quotient map and in particular surjective. \\ 
Hint: Use Michael's Selection Theorem for $r = 0$ and use induction 
on $r$.
\end{exer}

\begin{exer}
  \mlabel{exer:8.5.0} 
Let $\phi \: G \to H$ be a morphism of connected Lie groups 
and $\tilde \phi \: G \to \tilde H$ a lift with 
$q_H \circ \tilde\phi = \phi$, where 
$q_H \: \tilde H \to H$ is the universal covering morphism. 
Show that $(\ker \phi)_0 = (\ker \tilde\phi)_0$. 
\end{exer}

\begin{exer}
  \mlabel{exer:8.5.1} 
Let $M$ be a compact manifold with boundary, $K$ be a locally exponential 
Lie group and $r \in \N_0 \cup \{\infty\}$. For a closed subset 
$S \subeq M$, we consider the subgroup
\[ C^r_S(M,K) := \{ f \in C^r(M,K) \: (\forall s \in S)\, f(s) = \be\}.\] 
Show that $C^r_S(M,K)$ is a locally exponential Lie subgroup 
of $C^r(M,K)$ (Theorem~\ref{thm:mapgro-Lie}). \\ 
Hint: Consider the restriction map $C^r(M,K) \to C(S,K)$ 
and use Proposition~\ref{prop:5.5.9}. 
\end{exer}
  
\begin{exer}\mlabel{exer:8.6.1.1}
For $G$ be a Lie group and $r \in \N_0 \cup \{\infty\}$ and consider the evaluation morphism 
\[ \ev_1 \: P^{(r)}(G) \to G, \quad f \mapsto f(1).\] 
Show that, for any chart $\phi \: U \to G$, where $U \subeq \g$ is a convex $0$-neighborhood 
and $\phi(0) = \be$, the map 
\[  \sigma \: \phi(U) \to P^{(r)}G, \quad \sigma(\phi(x))(t) := \phi(tx) \]
is a smooth section of $\ev_1$.   
\end{exer}

\begin{exer} \mlabel{exer:8.6.1.3} 
Let $G$ be a $1$-connected locally exponential Lie group 
with Lie algebra $\g$. Suppose that there exists a continuous 
function 
$F \: [0,1] \times \g \to \g$ with 
$F_0 = \id_\g$ and $F_1 = 0$ such that every $F_t$ is a homomorphism 
of Lie algebras. Show that $G$ is contractible. \\ 
Hint: Mimik the arguments in the proof of Theorem~\ref{thm:enl-path}(b) 
and use open $0$-neighborhoods $U_1, U_2 \subeq \g$ such that 
$\exp\res_{U_2}$ is a diffeomorphism onto an open subset of $G$ and 
$F_t(U_1) \subeq U_2$ for every $t \in [0,1]$.   
\end{exer}

\begin{exer} Let $(V,\omega)$ be a locally convex symplectic vector space 
and $\g := \heis(V,\omega)$ be the corresponding 
\index{Heisenberg algebra} 
{\it Heisenberg algebra}, i.e., 
$\g = \R \times V$ with the bracket 
\[ [(z,v), (z',v')] := (\omega(v,v'), 0).\] 
Show that $f_t(z,v) := (t^2z, t v)$ is a family of Lie endomorphisms 
of $\g$ satisfying the assumptions of Exercise~\ref{exer:8.6.1.3}.   
\end{exer} 

\begin{exer} Let $G$ be a Lie group and 
$q_G \:\tilde G \to G_0$ be the simply connected covering group of its 
identify component. Show that: 
\begin{description}
\item[\rm(a)] For each $r \in \N_0 \cup \{\infty\}$, the map 
$P^{(r)}(\tilde G) \to P^{(r)}(G), \gamma \mapsto q_G \circ \gamma$ is an isomorphism 
of Lie groups. 
\item[\rm(b)] $\tilde G \cong P^{(r)}(G)/\Omega^{(r)}(G)_0$
\item[\rm(c)] $\pi_1(G) \cong \pi_0(\Omega^{(r)}(G))$. 
\end{description}
\end{exer}

\begin{exer} Let $\phi \:  \g \to \fh$ be a continuous surjective homomorphism 
of locally convex Lie algebras. Show that: 
\begin{description}
\item[\rm(a)] $\phi(\fz(\g)) \subeq \fz(\fh)$. 
\item[\rm(b)] For every continuous derivation $D \in \der(\g)$, there is at most 
one continuous derivation $D' \in \der(\fh)$ with 
$D' \circ \phi = \phi \circ D$. Further, 
\[ \fd := \{ D \in \der(\g) \: (\exists D' \in \der(\fh))\ \phi \circ D 
= D' \circ \phi\} 
\] 
is a Lie subalgebra and $\fd \to \der(\fh), D \mapsto D'$ 
is a homomorphism of Lie algebras. 
\item[\rm(c)] For every automorphism $\alpha \in \Aut(\g)$, there is at most 
one automorphism $\alpha' \in \Aut(\fh)$ with 
$\alpha' \circ \phi = \phi \circ \alpha$. Further, 
\[ A := \{ \alpha \in \Aut(\g) \: (\exists \alpha' \in \Aut(\fh)) \ 
\phi \circ \alpha = \alpha' \circ \phi\} \] 
is a subgroup and $A \to \Aut(\fh), \alpha \mapsto \alpha'$ 
is a group homomorphism. 
\item[\rm(d)] If, in addition, $\g$ and $\fh$ are locally exponential, 
then there exists a continuous homomorphism $\phi_{\rm ad} \: G_{\rm ad} \to H_{\rm ad}$ 
with $\phi_{\rm ad}(e^{\ad x}) = e^{\ad \phi(x)}$ for $x \in \g$. 
\end{description}
\end{exer}

\begin{exer} \mlabel{exer:cent-uh} Show that, for any complex Hilbert 
space~$\cH$, 
the center $\fz(\fu(\cH))$ of the Lie algebra $\fu(\cH)$ of 
skew-hermitian bounded operators coincides with $i\R \1$. \\
Hint: Evaluate the condition that any central element commutes with the 
rank-$1$-projections $P_v(w) := \la w, v \ra v$, where $v \in \cH$ 
is a unit vector. 
\end{exer}

\end{small}

\subsection*{Notes on Section~\ref{sec:11.4}} 

The first systematic investigation of enlargeability 
for Banach--Lie algebras is due to van Est and Korthagen 
\cite{EK64}. In particular, 
they introduce for every central extension 
$\hat\g$ of an enlargeable Banach--Lie algebra $\g$ by $\fz$ a 
subgroup $\Pi \subeq \fz$ and show that its discreteness is equivalent to 
the enlargeability of~$\hat\g$. Their 
period homomorphism arises as an element of 
$H^2_{\rm sing}(G,\z) \cong \Hom(H_2(G),\z)$ 
obtained from the enlargeability 
theory of local groups \cite{Est62}. Here $G$ is a $1$-connected Banach--Lie 
group with Lie algebra~$\g$. 
Then the Hurewicz homomorphism $\pi_2(G) \to H_2(G)$ 
is an isomorphism, so that their period homomorphism also is a 
homomorphism $\pi_2(G) \to \fz$, and one can even show that, 
up to sign, it is the one introduced above. 
The definition of the period homomorphism 
in terms of integration of differential forms as in \cite{Ne02a} 
makes it much more accessible than the implicit construction in \cite{EK64}. 

The results of the present section generalize 
\cite{EK64} to locally exponential Lie algebras and simplify 
the approach to non-enlargeability considerably. 
Our main point is that that the 
Enlargeability Theorem~\ref{thm4.4.7} 
can be obtained rather directly from the 
enlargeability of the path Lie algebra $P(\g)$ and 
the Enlargeability Theoren for Quotients~\ref{thm:quot-alg}.

Proposition~\ref{prop:path-int} implies in particular 
the integrability of $P^{(0)}(\g)$ for every Banach--Lie algebra, which 
was proved by by Swierczkowski \cite{Sw71}. 
The method to study enlargeability of a locally exponential 
Fr\'echet--Lie algebra $\g$ by writing it as a quotient of an enlargeable 
Lie algebra, such as $P(\g)$, already appears in 
Robart's work (cf.~\cite[Thm.~5]{Rob02}). 
This method is the starting point of a theory 
of certain topological groups which are more general than Lie groups, namely 
quotients of Lie groups. This leads to the concept of a 
{\it scheme of Lie groups}, or {\it $S$-Lie groups} (cf.\ \cite{Sr92}, 
\cite{DL66} and \cite{Est84}). 
\index{scheme of Lie groups} \index{$S$-Lie groups} 

In \cite{Woj06}, Wojty\'nski describes a variant of our approach 
for Banach--Lie algebras. 
Instead of the Banach--Lie algebra $P^{(0)}(\g)$ or the Fr\'echet--Lie 
algebra $P(\g)$, he considers 
analytic paths $\gamma(t) := \sum_{n = 1}^\infty a_n t^n$ in $\g$, 
for which 
$\|\gamma\|_1 := \sum_{n = 1}^\infty \|a_n\|$ is finite. Identifying these curves 
with their coefficient sequences, we denote this space by 
$\ell^1(\g) := \ell^1(\N,\g)$. The Lie bracket on this sequence space is 
given by 
\begin{equation}
  \label{eq:6.1.1}
 [(a_n), (b_n)] = (c_n)\quad \hbox{ with} \quad c_n = \sum_{j = 1}^{n-1} [a_j, b_{n-j}]. 
\end{equation}

With the same Lie bracket, we also turn the full sequence space 
$\g^\N$ into a pro-nilpotent Fr\'echet--Lie algebra, which is exponential 
(Proposition~\ref{prop:prolim-exp}). 
Since the Banach--Lie algebra $\ell^1(\g)$ injects into the 
exponential Lie algebra~$\g^\N$, it is enlargeable 
by the Integral Subgroup Theorem~\ref{thm:5.5.3} 
(one can also use the enlargeability of $P(\g)$). 
We have an evaluation map 
\[  q \: \ell^1(\g) \to \g, \quad (a_n) \mapsto \sum_{n = 1}^\infty a_n,  \] 
which is a quotient morphism of Lie algebras and since 
the subgroup 
$\la \exp \ell^1(\g)\ra$ is contractible \cite{Woj06}, one may proceed with 
the Enlargeability Theorem for Quotients~\ref{thm:quot-alg} 
in a similar fashion as for $P^{(0)}(\g)$. 

In \cite{Pe93a, Pe95a}, Pestov shows that if $E$ is a Banach space of 
$\dim E > 1$, then the free Banach--Lie 
algebra $F(E)$ over $E$ has trivial center. As a consequence, every Banach--Lie algebra 
$\g$ of dimension $> 1$ is a quotient of a centerless Banach--Lie algebra 
$F(\g)$, the free Banach--Lie algebra over the Banach space $\g$, which is enlargeable 
because its center is trivial (Corollary~\ref{cor:5.3.9}). Again, we can proceed 
with Theorem~\ref{thm:quot-alg} to derive enlargeability criteria. 

One also finds interesting variations on this theme in the literature, 
such as the following result of Swierczkowski, where one assumes 
enlargeability of an idean $\fn$ and the quotient $\g/\fn$  
(\cite[Thm., Sect. 12]{Sw67}).  
\begin{thm}
  \mlabel{thm:VI.1.19} Suppose that $\g$ is a Banach--Lie algebra 
and $\fn \trile \g$ is a closed enlargeable ideal for which $\fq := \g/\fn$ is enlargeable 
to some group $Q$. If $\pi_2(Q)$ vanishes, then $\g$ is enlargeable. 
\end{thm}

Note that the assumption is always satisfied if $\fq$ is abelian 
or finite-dimensional. It generalizes Theorem~\ref{thm:solv-int} below. 

\begin{defn}
  \mlabel{def:VI.1.20}  A 
Banach--Lie algebra is said to be {\it lower solvable} 
\index{Banach--Lie algebra!lower solvable} 
if there exists 
an ordinal number $\alpha$ and an ascending chain of closed subalgebras 
$$ \{0\} = \g_0 \subeq \g_1 \subeq \g_2 \subeq \ldots \subeq \g_\beta \subeq 
\g_{\beta+1} \subeq \ldots \subeq \g_\alpha = \g $$
such that the following condition is satisfied:  
\begin{description}
\item[\rm(a)] If $\beta \leq \alpha$ is not a limit ordinal, then 
$\g_{\beta-1}$ is an ideal of $\g_\beta$ containing all commutators. 
\item[\rm(b)] If $\beta \leq \alpha$ is a limit ordinal, then 
$\g_{\beta}$ is the closure of $\bigcup_{\gamma < \beta} \g_\gamma$. 
\end{description}
\end{defn}

The following theorem {\rm(\cite[Thm.~2]{Sw65})} 
is an immediate consequence of Theorem~\ref{thm:VI.1.19}, 
applied to the situation where $\fq$ is abelian. 
It implies in particular Theorem~\ref{thm:solv-int} below. 

\begin{thm}
  \mlabel{thm:VI.1.21} Each lower solvable 
Banach--Lie algebra is enlargeable.
\end{thm}

We call a norm $\|\cdot\|$ on a Lie algebra $\g$ 
{\it submultiplicative} if 
\index{submultiplicative norm on Lie algebra} 
\index{Banach--Lie algebra!contractive} 
$\|[x,y]\| \leq \|x\|\|y\|$ for all $x,y \in \g$. A Banach--Lie algebra 
$(\g,\|\cdot\|)$ is called {\it contractive} if its norm is submultiplicative. 
For any contractive Banach--Lie algebra, we define 
$$\delta_\g := \inf \{ \|x\|\: 0 \not= x \in \Pi(\g)\} \in [0,\infty] $$
and note that $\g$ is enlargeable if and only if $\delta_\g > 0$, which is equivalent 
to the discreteness of the period group $\Pi(\g)$ 
(Theorem~\ref{thm4.4.7}). The following theorem is originally due to 
Pestov who proved it with non-standard methods. A ``standard'' proof 
has been given by D.~Belti\c{t}\u{a} in \cite{Bel04}. 

\begin{thm}
  \mlabel{thm:VI.1.24} {\rm(Pestov's Local Theorem on Enlargeability)} 
A contractive Banach--Lie algebra $\g$ is enlargeable if and only if 
there exists a directed family ${\cal H}$ of closed subalgebras of $\g$ 
whose union is dense in $\g$ and 
\[ \inf \{ \delta_\h \: \h \in {\cal H}\} > 0.\] 
\end{thm}

\begin{probl}
  \mlabel{prob:VI.6} Extend Pestov's Theorem (Theorem~\ref{thm:VI.1.24}) 
to locally exponential Lie algebras.  
\end{probl}

Since, for each finite-dimensional Lie algebra $\g$, the period group is 
trivial, we have $\delta_\g = \infty$, and the preceding theorem, applied to the 
directed family of finite-dimensional subalgebras of $\g$ leads to:

\begin{cor}
  \mlabel{cor:VI.1.25} {\rm(\cite{Pe92, Bel04})} If $\g$ is a Banach--Lie 
algebra containing 
a locally finite-dimensional dense subalgebra, then $\g$ is enlargeable. 
\end{cor}

Along similar lines, Pestov also obtained the following criterion 
for enlargeability (\cite[Thm.~7]{Pe93b}). 
\begin{thm}
  \mlabel{thm:VI.1.26} A Banach--Lie algebra $\g$ is 
enlargeable if and only if all its separable closed subalgebras are. 
\end{thm}

\nin{\bf Linear Banach--Lie algebras:} Let us call a Banach--Lie algebra $\g$ 
\index{Banach--Lie algebra!linear} 
{\it linear} if it has a faithful homomorphism into some Banach algebra~$\cA$. 
According to Ado's Theorem, each finite-dimensional Lie algebra is linear, 
but the situation becomes more interesting, and also more complicated, 
for Banach--Lie algebras. 
In view of the Integral Subgroup Theorem (Theorem~\ref{thm:5.5.3}), 
 enlargeability is necessary for linearity, but it 
is not sufficient. 
In fact, if the Lie algebra $\g$ of a $1$-connected Banach--Lie group $G$ 
contains elements $p,q$ for which $[p,q]$ is a non-zero central element  
with $\exp_G([p,q]) = \1$, then $\g$ is not linear, because 
any morphism $\g \to \gl(E)$, $E$ Banach, would lead to a linear representation of the 
quotient $H/Z$ of the $3$-dimensional Heisenberg 
group~$H$ modulo a cyclic central subgroup~$Z$. 
Such elements exist in the Lie algebra $\hat\g$ of the central 
extension of the Banach--Lie algebra $C^1(\bS^1, \su_2(\C))$ by $\R$ 
(Example~\ref{ex:VI.1.16}, \cite{vES73}). 
Since, for each Banach--Lie algebra $\g$ the quotient $\g_{\rm ad} = \g/\z(\g)$ 
is linear, the intersection $\fn$ of all kernels of linear representations of $\g$ 
is a central ideal of $\g$. This links the linearity problem intimately with 
central extensions: When is a central extension of a linear Banach--Lie algebra 
linear? As the enlargeability is necessary, the discreteness of the corresponding 
period group is necessary, but what else?  

In \cite{vES73}, W.T.\ van Est and S.\ \'Swierczkowski 
describe a condition on the cohomology class of a central extension which 
is sufficient for linearity. They apply this in particular to show that, 
under some cohomological condition involving the center,  
for a Banach--Lie algebra $\g$, the $C^1$-path algebra $P^{(1)}(\g)$ 
is linear. It is remarkable 
that their argument does not work for $C^0$-curves. 
Closely related to this  circle of ideas is van Est's proof of 
Ado's theorem (\cite{Est66}), 
based on the vanishing of $\pi_2(G)$ for each finite-dimensional 
Lie group~$G$ (cf.~Remark~\ref{rem1.5}). 

It is also interesting to note that, for a real Banach--Lie algebra $\g$,  
linearity implies the linearity of the complexification $\g_\C$, 
which in turn implies that 
$\g_\C$ is enlargeable, and this is crucial for the existence of universal 
complexifications of the corresponding groups (cf.\ Theorem~\ref{thm:VI.1.23} below). 
In view of the Integral Subgroup Theorem (Theorem~\ref{thm:5.5.3}), 
we thus have the implications 
\begin{equation}
  \label{eq:14.6.s}
 \g \quad \mbox{ linear } \quad 
\Rarrow \quad \g_\C \quad \mbox{\ enlargeable\ } \quad \Rarrow \quad \g 
\quad \mbox{ enlargeable}.
\end{equation}

\section{Integrability} 
\mlabel{sec:11.5}

We recall that a locally convex Lie algebra $\g$ is said to be 
integrable if 
there exists some Lie group $G$ with Lie algebra $\g$. 
We have already seen many results on 
the enlargeability of locally exponentiall 
Lie algebra. For Lie algebras which are not 
locally exponential only rather isolated results are available.

\subsection{Analytic groups} 

\begin{prop}
  \mlabel{prop:monas-VI.2}
Let $G$ be a connected real analytic Lie group with Lie algebra~$\g$. Then 
each closed ideal of $\g$ is invariant under $\Ad(G)$. 
\end{prop}

\begin{prf}
Let $\fn \trile \g = \L(G)$ be a closed ideal. 
Since $G$ is connected, it suffices to show that there exists 
an $\be$-neighborhood $U_G \subeq G$ with $\Ad(U_G)\fn \subeq \fn$. 
Let $U \subeq \g$ be an open convex $0$-neighborhood and 
$\phi \: U \to G$ be a bianalytic diffeomorphism onto an open subset 
of $G$ with $\phi(0) = \be$. 
Let $x \in U$, $w_0 \in \fn$ and $w(t) := \Ad(\phi(tx))w_0$ for $t \in [0,1]$. 
We have to show that 
$w(1) \in \fn$. For $p(t) := \phi(tx)$, the right logarithmic derivative 
$v :=  \Ad(p)\delta(p)\: [0,1] \to \g$ of $p$ 
satisfies the differential equation 
\begin{align} 
  \label{eq:6.1}
 w'(t) 
&= \Ad(p(t))[p^{-1}(t)p'(t), w_0] \notag\\
&= \Ad(p(t))[\delta(p)_t, w_0] = [v(t), w(t)]. 
\end{align}
Since $v$ and $w$ are analytic, their Taylor expansions converge for 
$t$ close to $0$:
\[ v(t) = \sum_{n=0}^\infty v_n t^n \quad \mbox{ and }\quad  
 w(t) = \sum_{n=0}^\infty w_n t^n \quad \mbox{ in } \quad \g.\]  
Then the differential equation \eqref{eq:6.1} for $w$ can be written as 
\[  \sum_{n=0}^\infty (n+1) w_{n+1} t^n = w'(t) = [v(t), w(t)]
= \sum_{n=0}^\infty t^n \sum_{k=0}^n [v_k, w_{n-k}]. \]
Comparing coefficients now leads to
\[  w_{n+1} = \frac{1}{ n+1} \sum_{k=0}^n [v_k, w_{n-k}], \]
so that we obtain inductively $w_n \in \fn$ for each $n \in \N$. 
Since $\fn$ is closed, we get $w(t) \in \fn$ for $t$ close to $0$. 
Applying the same argument in other points $t_0 \in [0,1]$, we see that 
the set $w^{-1}(\fn)$ is an open closed subset of $[0,1]$, and therefore that 
$w(1) \in \fn$ because $w(0) \in \fn$ and $[0,1]$ is connected. 
\end{prf} 

\begin{cor}
  \mlabel{cor:monas-VI.2}
Let $G$ be a connected complex Lie group with Mackey complete 
Lie algebra~$\g$. Then 
each closed ideal of $\g$ is invariant under $\Ad(G)$. 
\end{cor}

\begin{prf} We only have to observe that, by Theorem~\ref{charcxcompl} 
and the Mackey completeness of $\g$, the complex Lie group $G$ 
is complex analytic, hence real analytic.
\end{prf}

For the following corollary, we recall that an ideal 
$\fn \trile \g$ is called stable \index{stable} 
if it is invariant under $e^{\ad \g}$ (Definition~\ref{def:5.3.3}). 

\begin{cor}
  \mlabel{cor:monas-VI.3}
If $\g$ is a complex Fr\'echet--Lie algebra containing 
a closed ideal which is not stable, 
then $\g$ is not integrable to a 
complex Lie group. 
\end{cor}

\subsection{Lie algebras of vector fields} 

Proposition~\ref{prop:monas-VI.2} has also direct implications for 
Lie algebras of vector fields. 

\begin{cor}
\mlabel{cor:monas-VI.4}
If $M$ is a compact smooth manifold, then 
the Lie group $\Diff(M)$ does not possess an analytic Lie group structure. 
\end{cor} 

\begin{prf} For each non-dense open subset $U \subeq M$, the subspace 
\[  {\cal V}(M)_U := \{ X \in {\cal V}(M) \: X\res_U = 0\} \]
is a closed ideal of the Lie algebra 
${\cal V}(M)$ not invariant under $\Diff(M)$ because 
$\Ad(\phi){\cal V}(M)_U = {\cal V}(M)_{\phi(U)}$  for 
$\phi \in \Diff(M)$. Therefore Proposition~\ref{prop:monas-VI.2} applies. 
\end{prf}

\begin{thm}
  \mlabel{thm:monas-VI.5}If $M$ is a compact manifold, 
then the Fr\'echet--Lie algebra ${\cal V}(M)_\C$ is not integrable to a 
regular complex Lie group. 
\end{thm} 

\begin{prf} (Sketch; see \cite{Mil84}) Let $\g := {\cal V}(M)_\C$ and $U \subeq M$ be 
an open non-empty 
subset of $M$ which is not dense. Then  
\[  \fn_U := \{ x \in \g \: x\res_U = 0\} \] 
is a closed ideal of $\g$. 
Let $G$ be a regular complex Lie group with Lie algebra $\g$ and 
let $q \: \tilde\Diff(M) \to \Diff(M)_0$ denote the 
universal covering homomorphism of $\Diff(M)_0$. 
Then the 
inclusion ${\cal V}(M) \into \g$ can be integrated to a Lie group homomorphism 
$\phi \: \tilde\Diff(M) \to G$. For $g \in \tilde\Diff(M)$, we then have 
\[  \Ad(\phi(g))\fn_U = \fn_{\phi(g)(U)}, \] 
contradicting the invariance of $\fn_U$ under $\Ad(G)$ 
(Proposition~\ref{prop:monas-VI.2}). 
\end{prf}

The following non-integrability result by L.~Lempert (\cite{Lem97}) 
is much stronger 
because no regularity is assumed. It does not even 
assume the existence of an exponential function. 
Its outcome is that complexifications 
of Lie algebras of vector fields are rarely integrable. 

\index{Lempert's No-Complexification Theorem}
\begin{thm}[Lempert's No-Complexification Theorem] 
  \mlabel{thm:VI.2.5} 
Let $M$ be a compact manifold of positive dimension. 
Then the complexifications $\g_\C$ of the following Lie algebras $\g$ are not integrable: 
\begin{description}
\item[\rm(1)] The Lie algebra ${\cal V}(M)$ of smooth vector fields on $M$. 
\item[\rm(2)] If $M$ is analytic, the Lie algebra ${\cal V}^\omega(M)$ 
of analytic vector fields on $M$. 
\item[\rm(3)] If $\Omega$ is a symplectic $2$-form on $M$, the Lie algebra 
\[ {\cal V}(M,\Omega) := \{ X \in {\cal V}(M) \: {\cal L}_X\Omega = 0\} \] 
of symplectic vector fields on $M$. 
\item[\rm(4)] If $M$ is analytic and $\Omega$ is an analytic 
symplectic $2$-form on $M$, the Lie algebra 
${\cal V}^\omega(M,\Omega)$ of analytic symplectic vector fields on $M$. 
\end{description}
\end{thm}

\begin{prf}  (Sketch) Lempert's 
proof is based on the following result, which is obtained by PDE methods: 
If $\xi \: \R \to \g_\C$ is a smooth curve such that, for each 
$x \in \g_\C$, the initial value problem 
\[ \gamma(0) = x, \quad \dot\gamma(t) = [\xi(t),\gamma(t)] \] 
has a smooth solution, then $\xi(0) \in \g$. 

Lempert gives another argument, based on the fact that 
\begin{equation}
  \label{eq:6.2.1}
 \Aut({\cal V}(M)_\C) 
\cong \Aut({\cal V}(M)) \rtimes \{\id,\sigma\} 
\cong \Diff(M) \rtimes \{\id,\sigma\}, 
\end{equation}
where $\sigma$ denotes the complex conjugation on ${\cal V}(M)_\C$. 
The first isomorphism is obtained in \cite{Lem97}, using \cite{Ame75}, and 
the second is an older result of Pursell and Shanks \cite{PuSh54} 
(cf.~\cite[Thm.~IX.2.1]{Ne06a}). 

Clearly \eqref{eq:6.2.1} implies that, for any 
connected Lie group $G$ with Lie algebra $\L(G) = {\cal V}(M)_\C$, the group 
$\Ad(G) \subeq \Aut({\cal V}(M)_\C)$ preserves the real 
subspace ${\cal V}(M)$. Taking derivatives of orbit maps, this 
leads to the contradiction 
$[{\cal V}(M)_\C, {\cal V}(M)] \subeq {\cal V}(M)$.

Another key point of Lempert's proof is that, 
for any $X \in {\cal V}(M)$, 
the one-parameter group $\exp(\R X)$ acts on $\g$ precisely 
as the corresponding one-parameter 
group of $\Diff(M)$. This argument requires a uniqeness 
lemma for ``smooth'' maps with values in 
$\Aut(\g)$, which is far from being a Lie group (cf.~Lemma~\ref{lem:e.1.3}). 
\end{prf}

Analytic vector fields do no longer form a Fr\'echet space because 
they have a natural direct limit structure. They arise as Lie algebras 
of groups of analytic diffeomorphisms.

\begin{thm}
  \mlabel{thm:VI.2.4} {\rm(\cite{Les82, Les83})} Let $M$ be a compact analytic manifold and 
${\cal V}^\omega(M)$ the Lie algebra of analytic vector fields on $M$. 
Then ${\cal V}^\omega(M)$ carries a natural Silva space structure, turning it into a 
topological Lie algebra, and the group $\Diff^\omega(M)$ of 
analytic diffeomorphisms 
carries a smooth Lie group structure 
for which ${\cal V}^\omega(M)^{\rm op}$ is its Lie algebra. 
\end{thm}

It was shown by Tagnoli \cite{Ta68} that the group $\Diff^\omega(M)$, 
$M$ a compact analytic manifold, carries no analytic Lie group structure 
(cf.\ \cite[Ex.~6.12]{Mil82}). That there is no analytic Lie group with an 
analytic exponential function and Lie algebra ${\cal V}^\omega(M)$ 
can be seen by verifying that the map $(X,Y) \mapsto \Ad(\Phi^X_1)Y$ 
is not analytic on a $0$-neighborhood in ${\cal V}^\omega(M) \times 
{\cal V}^\omega(M)$ (cf.\ \cite[Ex.~6.17]{Mil82}). 

\begin{thm}
  \mlabel{thm:VI.2.6} {\rm(\cite{Omo81})} For any non-compact $\sigma$-compact smooth 
manifold $M$ of positive dimension, 
the Lie algebra ${\cal V}(M)$ is not integrable to any Lie group 
with an exponential function. 
\end{thm}

\begin{prf}
(Sketch) If $G$ is a Lie group with Lie algebra $\L(G) = {\cal V}(M)$ 
and an exponential function, then we obtain 
 for each 
$X \in {\cal V}(M)$ a smooth  {one-parameter} group 
$t \mapsto \Ad(\exp_G(tX))$ of automorphisms of ${\cal V}(M)$ 
with generator $\ad X$. By \cite[Thm.~2]{Ame75}, 
$\Aut({\cal V}(M)) \cong \Diff(M)$, so that we obtain a 
one-parameter group $\gamma_X$ of $\Diff(M)$ which then is shown to coincide 
with the flow generated by $X$ (Proposition~\ref{prop:e.1.4}; 
\cite[Sect.~3.4]{KYM85}). 
This contradicts the existence of
non-complete vector fields on $M$. 
\end{prf}

\subsection{Integrability of linear operators} 

In this subsection we discuss several interesting links between 
Lie theory and the integrability of continuous operators 
$D \: E \to E$ on a locally convex space in the sense that 
they are generators of smooth linear $\R$-actions on~$E$. 
In particular, we show that the solvable Lie algebra 
$E \rtimes_D \R$ is integrable if and only if $D$ is integrable.

\begin{defn} (cf.\ Definition~\ref{def:integrable-op}) 
Let $E$ be a locally convex space and 
$D \in \cL(E)$ be a continuous linear operator on $E$. 
We say that $D$ is {\it integrable} if 
\index{endormorphism!integrable} 
there exists a homomorphism $\gamma_D \: \R \to \GL(E)$ which is 
{\it smooth}  \index{smooth homomorphism to $\GL(E)$} 
in the sense that the corresponding action $\R \times E \to E, 
(t,v) \mapsto \gamma_D(t)v$ is smooth and 
\[ \frac{d}{dt}\Big|_{t= 0} \gamma_D(t)v = Dv \quad \mbox{ for every } \quad v 
\in E.\] 
Note that $\gamma_D$ is uniquely determined whenever it exists 
(Lemma~\ref{lem:e.1.3}). 
\end{defn}

\begin{ex} (a) If $E$ is a Banach space, then each bounded operator 
$D\in \cL(E)$ integrates to  
a homomorphism $\alpha$ which is continuous with respect to the norm topology on 
$\GL(E)$ and given by the convergent exponential series 
$\gamma_D(t) := \sum_{k = 0}^\infty \frac{t^k}{k!} D^k = e^{tD}.$

If, conversely, $\gamma \: \R\to \GL(E)$ defines a smooth $\R$-action 
on $E$, then its infinitesimal generator $D := \gamma'(0)$ is a closed operator 
which is everywhere defined, hence bounded by the Closed Graph Theorem. 

(b) If $\gamma \: \R \to \GL(E)$ 
is a one-parameter group on a locally convex space~$E$ 
defining a continuous action, then the subspace $E^\infty \subeq E$ of smooth vectors 
carries a natural topology inherited from the embedding 
\[ E^\infty \into C^\infty(\R,E),\quad  v \mapsto \gamma^v, \quad 
\gamma^v(t) = \gamma(t)v. \]
It can also be obtained from the seminorms 
$p_n(v) := p(D^nv)$ on $E^\infty$, where $p$ is a continuous seminorm 
on $E$ and $D := \gamma'(0)$. Then the induced one-parameter group 
$\gamma^\infty \: \R \to \GL(E^\infty)$ defines a smooth action 
by \cite[Thm.~4.4]{Ne10}, so that the operator 
$D$ on $E^\infty$ is integrable. 
In this sense each continuous linear action leads to an integrable 
operator on the space of smooth vectors. 
\end{ex}

\begin{ex} If $M$ is a finite-dimensional $\sigma$-compact manifold 
and $X \in {\cal V}(M)$ a vector field, then the corresponding derivation of the Fr\'echet 
algebra $C^\infty(M,\R)$ is integrable if and only if the vector field 
$X$ is complete. 
\end{ex}

\begin{ex} Let $E$ be a real pre-Hilbert space 
and $\cH$ be the Hilbert space completion of $(E, \la \cdot,\cdot\ra)$. 
We consider an operator $D \in \cL(E)$ which is skew-symmetric, i.e., 
\[ \la Dv,w \ra = - \la v, D w\ra \quad \mbox{ for } \quad v,w \in E,\] 
and assume that the topology on $E$ can be defined by the seminorms 
\[ p_n(v) := \|D^nv\|, \quad n \in \N_0.\] 
Then the inclusion $E \into \cH$ is continuous. 

Suppose first that $D$ is integrable. Then 
$\gamma_D \: \R \to \GL(E)$ preserves the scalar product because 
\[ \frac{d}{dt}  \la \gamma_D(t)v, \gamma_D(t)w \ra 
= \la D\gamma_D(t)v, \gamma_D(t)w \ra + \la \gamma_D(t)v, D\gamma_D(t)w \ra  
= 0.\] 
Therefore $\gamma_D$ extends to a strongly continuous orthogonal one-parameter group 
$U_t$ on~$\cH$. Since $U_t$ leaves the subspace $E$ invariant, 
Stone's Theorem \cite[Thm.~VIII.10]{RS75} implies that $A = \gamma_D'(0)$ 
is the closure of $D$, i.e., the graph $\Gamma(A)\subeq \cH \times \cH$ 
is the closure of the graph $\Gamma(D)$. 

Suppose, conversely, that the closure $A$ of $D$ is skew-adjoint, 
i.e., $\oline D = - D^*$ (in the sense 
of unbounded operators; see \cite{Ru91}). Then Stone's 
Theorem \cite[Thm.~VIII.10]{RS75} implies 
the existence of a strongly continuous one-parameter 
group $\gamma \: \R \to \OO(\cH)$ with infinitesimal generator $A$, 
i.e., the domain $\cD(A)$ of $A$ is the set of elements for which 
$t \mapsto \gamma(t)v$ is $C^1$ and $A v= \derat0 \gamma(t)v$ for 
every $v \in \cD(A)$. Suppose, in addition, that each $\gamma(t)$ leaves 
$E$ invariant. We put $\gamma_D(t) := \gamma(t)\res_E$. 
Now $E$ consists of smooth vectors for $\gamma$, i.e., 
$E \subeq \bigcap_{n \in \N} \cD(A^n) = \cH^\infty$ and we recall from \cite{Ne10} 
that the $\R$-action on $\cH^\infty$ is smooth. As the Fr\'echet topology on $\cH^\infty$ is defined by the seminorms 
$p_n(v) = \|A^nv\|$, $n \in \N_0$, the inclusion 
$E \into \cH^\infty$ is a topological embedding. 
Therefore the smoothness of the action on $\cH^\infty$ implies the smoothness 
of the action on $E$. 

This discussion shows that 
$D$ is integrable if and only if its closure is skew-adjoint 
and the corresponding one-parameter group of isometries of $\cH$ leaves 
$E$ invariant.  The latter condition is automatically satisfied 
if $E$ is complete, because this implies that $E = \cH^\infty$. 
\end{ex}

The preceding example shows that the integrability problem 
for linear operators on Fr\'echet spaces contains in particular the problem 
to decide for a symmetric operator $A \in \cL(\cD)$ on a dense 
subspace $\cD$ of a Hilbert space whether it is selfadjoint or not. 

\begin{lem} \mlabel{lem:ad-int-exp} 
If $\g$ is the Lie algebra of a Lie group $G$ with 
exponential function, then each operator $\ad x \in \cL(\g)$ 
is integrable. 
\end{lem}

\begin{prf}
The one-parameter group 
$\gamma(t) := \Ad(\exp tx)$ defines a smooth 
$\R$-action on $\g$ and $\gamma'(0) = \ad x$ by 
Proposition~\ref{prop:der-Ad}. 
\end{prf}

Since the Hausdorff series can be used to define a Lie group structure on 
any nilpotent locally convex Lie algebra, all these Lie algebras are integrable 
(Theorem~\ref{thm:BCH-nil}). For the class of solvable 
Lie algebras, the integrability problem contains in particular the 
integrability problem for continuous linear operators on locally convex spaces:  

\begin{thm}
  \mlabel{thm:VI.2.7} Let $E$ be a locally convex space and $D \in \gl(E)$. 
Then the Lie algebra 
\[ \g := E \rtimes_D \R \quad \mbox{ with the bracket }\quad 
 [(v,t), (v',t')] := (tDv' - t'Dv, 0) \] 
is integrable 
if and only if $D$ is integrable. 
\end{thm}

\begin{prf}
If $D$ is integrable to a smooth representation 
$\alpha \: \R \to \GL(E)$ with $\alpha'(0) = D$, 
then  the semi-direct product $G := E \rtimes_\alpha \R$ is a Lie group with 
the Lie algebra $\g$ (Proposition~\ref{prop:liealg-semdir}). 

Suppose, conversely, that $G$ is a connected Lie group with Lie algebra~$\g$. 
Replacing $G$ by its universal covering group, we may 
assume that $G$ is \break {$1$-connected}. Then the regularity of the additive 
group $(\R,+)$ implies the existence of a smooth homomorphism 
$\chi \: G \to \R$ with $\L(\chi) = q$, where 
$q(v,t) = t$ is the projection $\g = E \rtimes_D \R \to \R$
(Theorem~\ref{thm3.2.11}). 
Let $\eps > 0$ and $\gamma \: [-\eps,\eps] \to G$ 
be any smooth curve with $\gamma(0) = \be$ 
and $q(\gamma'(0))~=~1$. 
Then $\eta := \chi \circ \gamma \: [-\eps,\eps] \to \R$ is a smooth curve 
with $\eta(0) = 0$ and $\eta'(0) = 1$. 
Rescaling $\gamma$, we may even assume that 
$\eta'(t)  > 0$ holds for each $|t| \leq \eps$. 
Hence the one-dimensional 
Inverse Function Theorem implies the smoothness of the inverse function 
$\eta^{-1} \: \im(\eta) = [\eta(-\eps),\eta(\eps)] \to [-\eps,\eps]$. 
Then $\beta := \gamma \circ \eta^{-1}$ is a smooth curve in $G$ 
with $\chi \circ \beta = \id_{\im(\eta)}$. 

Next we observe that $[\g,\g] \subeq E$ implies that 
$E$ is $\Ad(G)$-invariant, so that $\Ad_E(g) := \Ad(g)\res_E$ defines 
a smooth action of $G$ on $E$ whose derived representation is given by $\ad_E(x,t) = tD$. 
We now put 
$\alpha(t) := \Ad_E(\beta(t))$ and observe that 
\[  \delta(\alpha)_t = \ad_E(\delta(\beta)_t) 
= q(\delta(\beta)_t)\cdot D 
= \delta(\chi \circ \beta)_t\cdot D = D. \] 
It easily follows that, for 
$|t| \leq \eps$, we have $\alpha(\frac{t}{n})^n = \alpha(t)$, so that 
\[  \alpha(t) := \lim_{n \to \infty} \alpha\Big(\frac{t}{n}\Big)^n \] 
defines an extension of $\alpha$ to a smooth 
$\R$-action on $E$ (Exercise~\ref{exer:e.10}). 
Hence $D$ is integrable. 
\end{prf}

The following example is also of some interest for the integrability of 
Lie algebras of formal vector fields (Exercise~\ref{ex:VI.2.8}). 

\begin{ex} \label{ex:3.13} 
We consider the space 
$E := \R[[x]]$ of formal power series $\sum_{n = 0}^\infty a_n x^n$ in one variable 
and endow it with the Fr\'echet topology for which the map 
$\R^{\N_0} \to \R[[x]], (a_n) \mapsto \sum_{n = 0}^\infty a_n x^n$ is a topological isomorphism. 
Then the derivative 
\[ D\Big(\sum_{n=0}^\infty a_n x^n\Big) := \sum_{n = 0}^\infty a_{n+1} (n+1) x^n \] 
defines a continuous linear 
operator on $E$. We claim that this operator is not integrable. 

We argue by contradiction, and assume that $\alpha \: \R \to \GL(E)$ 
is a smooth $\R$-action of $E$ with $\alpha'(0) = D$. For each 
$n \in \N$, the curve $\gamma \: \R \to E, \gamma(t) := (x+t)^n$, 
satisfies 
$\dot\gamma(t) = n (x+t)^{n-1} = D\gamma(t)$, so that 
Lemma~\ref{unique-intcurve} 
implies that $\alpha(t)x^n = (x+t)^n$ for all $t \in \R$. 
Then we obtain
\[ \alpha(1)x^n = 1 + n x+ \ldots.\] 
In view of $\lim_{n \to \infty} x^n \to 0$ in $E$, this contradicts the 
continuity of the operator~$\alpha(1)$. Therefore $D$ is not integrable. 
\end{ex}

\begin{ex}
  \mlabel{ex:VI.2.8} Let 
\[ \g := \gf_n(\R)_{-1} := \R^n[[x_1,\ldots, x_n]] \] 
denote 
the space of all $\R^n$-valued formal power series in $n$ variables, 
considered as the 
\index{Lie algebra of formal vector fields}  
{\it Lie algebra of formal vector fields}, 
endowed with the bracket 
$$ [f,g](x) := dg(x)f(x) - df(x)g(x), $$
which makes sense on the formal level because if $f$ is homogeneous of degree $p$ 
and $g$ is homogeneous of degree $q$, then $[f,g]$ is homogeneous 
of degree $p+q-1$. 

We have already seen in Subsection~\ref{subsec:sternberg} that the subalgebra 
$\gf_n(\R)$ of all elements with vanishing constant term is the Lie algebra 
of the Fr\'echet--Lie group $\Gf_n(\R)$ of formal diffeomorphisms of $\R^n$ 
fixing $0$. If $\fa \subeq \g$ denotes the abelian Lie subalgebra 
of constant vector fields, then 
$\g$ is the vector space direct sum of the Lie subalgebras $\fa$ and~$\gf_n(\R)$. 

We claim that the Lie algebra $\g$ is not integrable to any Lie group 
with an exponential function. This strengthens a statement in 
\cite[p.~80]{KYM85}, 
that it is not integrable to a $\mu$-regular Fr\'echet--Lie group. 
Suppose that $\g$ is the Lie algebra of a Lie group with 
exponential function. Then Lemma~\ref{lem:ad-int-exp} implies that 
all operators $\ad x$ are integrable. 
For $x = e_1$, we have 
$[e_1, g] = \frac{\partial g}{\partial x_1},$
and we can now argue as 
in Example~\ref{ex:3.13} that $\ad e_1$ is not integrable, 
hence that $\g$ is not integrable to any Lie group with an exponential function. 
\end{ex}

The preceding example shows that the constant terms create problems in integrating 
Lie algebras of formal vector fields. This 
is very natural because the formal completion 
distinguishes the point $0 \in \R^n$. 
A similar phenomenon arises in the context of groups of germs of local diffeomorphisms. 
For germs of functions in $0$, the non-integrability of vector fields with 
non-zero constant term follows from the fact that all automorphisms preserve the 
unique maximal ideal of functions vanishing in $0$. 

Let $\gs_n(\R)_{-1}$ \index{$\gs_n(\R)_{-1}$} 
denote the space of germs of smooth maps $\R^n \to \R^n$ in $0$, 
identified with germs of vector fields in $0$. 
According to \cite[Sect.~5.2]{RK97}, 
this space carries a natural Silva structure, turning it into a 
locally convex Lie algebra. Let $\gs_n(\R)$ denote the subspace of all germs vanishing in $0$ 
and $\gs_n(\R)_1$ the set of germs with vanishing $1$-jet in~$0$. 

\begin{thm}
  \mlabel{thm:VI.2.9} {\rm(\cite[Thm.~3]{RK97})} The 
group $\Gs_n(\R)$ of germs of diffeomorphism of $\R^n$ in $0$ fixing $0$ 
carries a Lie group structure for which the Lie algebra is the space $\gs_n(\R)$ of 
germs of vector fields vanishing in $0$. 

We have a semidirect product decomposition $\Gs_n(\R) \cong \Gs_n(\R)_1 \rtimes\GL_n(\R)$, 
where $\Gs_n(\R)_1$ is the normal subgroup of those germs $[\phi]$ for which 
$\phi - \id_{\R^n}$ has vanishing $1$-jet. The map 
$$ \Phi \: \gs_n(\R)_{1} \to \Gs_n(\R)_{1}, \quad \xi \mapsto \id + \xi $$
is a global diffeomorphism. 
\end{thm}

In view of the preceding theorem, it is a natural problem 
to integrate Lie algebras of germs of vector fields 
vanishing in the base point to Lie groups of germs of diffeomorphisms. 
This program has been  carried out by N.~Kamran and T.~Robart 
in several papers (cf.\ \cite{RK97}, \cite{KR01, KR04}, \cite{Rob02}). 
It provides several interesting classes 
of Silva--Lie groups of germs of smooth and also analytic local diffeomorphisms, 
where the corresponding Silva--Lie algebras depend on certain 
parameters which are used to obtain a good topology. 

The following proposition (cf.\ \cite{Gl08a}) 
is a variant of E.~Borel's theorem on 
the Taylor series of smooth functions. It provides an interesting connection 
between the smooth global and the formal perspective on diffeomorphism groups. 

\begin{prop}
  \mlabel{prop:VI.2.11} Let 
$M$ be a smooth finite-dimensional manifold, 
$m_0 \in M$ and write $\Diff_c(M)_{m_0}$ for the stabilizer of $m_0$. 
For each $\phi \in \Diff_c(M)_{m_0}$, let 
$T^\infty_{m_0}(\phi) \in \Gf_n(\R)$ denote the Taylor series of $\phi$ in $m_0$ 
with respect to some local chart. 
Then the map 
$$ T^\infty_{m_0} \: \Diff_c(M)_{m_0,0} \to \Gf_n(\R)_0 $$
is a surjective homomorphism of Lie groups. 
Here the identity component $\Gf_n(\R)_0\subeq \Gf_n(\R)$ 
is the subgroup of index $2$ consisting of those formal 
diffeomorphisms $\psi$ with $\det(T_0(\psi)) > 0$. 
\end{prop}

\begin{ex}
  \mlabel{ex:VI.2.12} (Holomorphic germs)  Let $\gh_n(\C)$ denote the 
space of germs of holomorphic maps $f \: \C^n \to \C^n$ in $0$ satisfying $f(0) = 0$.  
We endow this space with the locally convex 
direct limit topology of the Banach spaces $E_k$ of holomorphic 
functions on the closed unit disc of radius $\frac{1}{ k}$ 
in $\C^n$ (with respect to any norm). Thinking of the elements of $\gh_n(\C)$ 
as germs of vector fields in $0$ leads to the Lie bracket 
\[ [f,g](z) := dg(z)f(z) - df(z) g(z), \] 
which turns $\gh_n(\C)$ into a topological Lie algebra.
This follows easily from
Proposition~\ref{prop:germcia}, applied with $K = \{0\}$.

The set $\Gh_n(\C)$ of all germs $[f]$ with $\det(f'(0)) \not=0$ is an open subset 
of $\gh_n(\C)$ which is a group with respect to composition 
$[f][g] := [f \circ g]$. In \cite{Pis77}, D.~Pisanelli 
shows that composition and inversion in 
$\Gh_n(\C)$ are holomorphic, so that $\Gh_n(\C)$ is a complex Lie group 
with respect to the manifold structure it inherits as an open subset of 
$\gh_n(\C)$. This Lie group has a holomorphic exponential function 
which is not locally surjective, where the latter fact can be obtained by 
adapting Sternberg's example 
\[ f(z) = e^{\frac{2\pi i}{m}} z + p z^{m+1}\quad \mbox{ for } \quad m > 0 \] 
from Subsection~\ref{subsec:sternberg} (see also \cite{Pis76}). 

Note that  $\Gh_n(\C) \cong \Gh_n(\C)_1 \rtimes \GL_n(\C)$, 
where $\Gh_n(\C)_1$ is the subgroup of all diffeomorphisms with linear term 
$\id_{\C^n}$. 
\end{ex}

\subsection*{Recently solved problems}

In \cite{Eyn12}, J.-M.~Eyni obtained the following result: 

\begin{thm} \mlabel{thm:cobanach} {\rm(Banach Complemented Subgroup Theorem)} 
Let $G$ be a Lie group and 
$\fh \subeq \L(G)$ be a closed Lie subalgebra such that there exists a 
Banach subspace $E \subeq \L(G)$ for which the addition map 
$\fh \oplus E \to \L(G)$ is a topological isomorphism. 
Then there exists an integral subgroup $H \subeq G$ with Lie algebra~$\fh$.
\end{thm}

\begin{cor} \mlabel{cor:cofin} {\rm(Finite Codimension Subgroup Theorem)} 
Let $G$ be a Lie group and 
$\fh \subeq \L(G)$ be a closed Lie subalgebra of finite codimension. 
Then there exists an integral subgroup $H \subeq G$ with Lie algebra~$\fh$.
\end{cor}

For $\mu$-regular groups, this corollary follows from 
\cite[Thm.~4.2]{OMY83a}.

In \cite{Da11}, R.~Dahmen shows in particular that: 

\begin{prop} The group $\Gh_n(\C)$ of germs of complex analytic diffeomorphisms 
around $0 \in \C^n$ {\rm(Example~\ref{ex:VI.2.12})} 
is a regular Lie group.   
\end{prop}

\subsection*{Some open problems}

\begin{probl}
\mlabel{prob:VI.2} ($2$-step solvable Lie algebras) 
Theorem~\ref{thm:VI.2.7} gives an integrability criterion 
for solvable Lie algebras of the type $\g = E \rtimes_D \R$. 

Since abelian Lie algebras are integrable for trivial reasons, it is natural 
to address the integrability problem for solvable Lie algebras by first 
restricting to algebras of {\it solvable length~$2$}, i.e., 
$D^1(\g) = \oline{[\g,\g]}$ is an abelian ideal of $\g$. Clearly, the adjoint 
action defines a natural topological module structure for the abelian Lie algebra 
$W := \g/D^1(\g)$ on $E := D^1(\g)$. Here are some problems concerning this situation: 
\begin{description}
\item[\rm(1)] Does the integrability of $\g$ imply that the Lie algebra module  
structure for $W$ on $E$ integrates to a smooth action of the Lie group 
$(W,+)$ on $E$?  If $W$ is finite-dimensional, this can be proved by 
an argument similar to the proof of Theorem~\ref{thm:VI.2.7}. 
\item[\rm(2)] Assume that the Lie algebra module structure for $W$ on $E$ integrates 
to a smooth action of $(W,+)$. Does this imply that $\g$ is integrable? 
\end{description}

If $\g \cong E \rtimes W$ is a semidirect product, the latter is obvious, but 
if $\g$ is a non-trivial extension of $W$ by $V$, the situation is more complicated. 
Note that all solvable Banach--Lie algebras are integrable by Theorem~\ref{thm:VI.1.21}. 
\end{probl}

\begin{probl}
  \mlabel{prob:VI.3} 
Is the group $\Gs_n(\R)_1$ of germs of diffeomorphisms $\phi$ of $\R^n$ fixing $0$, 
for which the linear term of $\phi - \id_{\R^n}$ vanishes, (locally) exponential? 
\end{probl}

\begin{probl}
  \mlabel{prob:VI.5} 
Is the subgroup $\Gh_n(\C)_1$ (Example~\ref{ex:VI.2.12}) 
an exponential Lie group? (cf.\ Problem~\ref{prob:VI.3}) 
\end{probl}

\section*{Notes on Section~\ref{sec:11.5}} 

Proposition~\ref{prop:monas-VI.2} is Lemma~9.1 in  \cite{Mil84}. 
For complexifications of Lie algebras of ILB--Lie groups, 
Omori describes results similar to Theorem~\ref{thm:VI.2.5} 
in \cite[Cor.~4.4]{Omo97}. 

\section{BCH--Lie groups} 
\mlabel{sec:BCH}

We recall from Definition~\ref{def:bcdh-Lie-alg} 
that a BCH--Lie algebra is a locally convex Lie algebra 
with the property that the Hausdorff series 
\[ x * y = \sum_{n = 1}^\infty H_n(x,y) 
= x + y + \frac{1}{2}[x,y] + \cdots \] 
converges 
for $x,y$ in some open $0$-neighborhood $U \subeq \g$
and defines an analytic exponential local group structure on some 
$0$-neighborhood in~$U$.

On the global level we have the following result whose proof requires 
the uniqueness assertion from Theorem~\ref{bch-locexp}: 

\begin{thm}
  \mlabel{thm:IV.1.8} For a Lie group $G$ with Lie algebra $\g$, 
the following are equivalent: 
  \begin{description}
\item[\rm(1)] $G$ is analytic with an analytic 
exponential function which is a local analytic 
diffeomorphism in $0$. 
\item[\rm(2)] $G$ is locally exponential and $\g$ is BCH. 
  \end{description}
\end{thm}

\begin{prf} 
(1) $\Rarrow$ (2): Let $U \subeq 0$ be an open $0$-neighborhood 
such that $\exp\res_U$ is an analytic diffeomorphism onto an open subset 
of $G$. Then 
\[ x * y := \exp\res_U^{-1}(\exp x \exp y) \] 
defines an analytic exponential local Lie group on a $0$-neighborhood 
of~$\g$. Therefore $\g$ is BCH by Theorem~\ref{bch-locexp}. 

(2) $\Rarrow$ (1): If $G$ is locally exponential, 
we can construct an atlas as in Theorem~\ref{thm:locglob} 
from an exponential chart around $\be$. Then the BCH property 
implies that the chart changes are actually analytic, so that 
we obtain an analytic Lie group structure on $G$. 
\end{prf}

\begin{defn}
  \mlabel{def:VI.1.9}  A 
group satisfying the equivalent conditions of Theorem~\ref{thm:IV.1.8} is 
called a {\it BCH--Lie group}. \index{BCH--Lie group}
\end{defn}

In Example~\ref{ex:5.1.10} we describe an analytic Lie group 
$G$ with an analytic exponential 
function which is a smooth diffeomorphism, but $\g$ is not BCH. This is a 
negative answer to a question raised in \cite[p.~31]{Mil84}. 

\begin{exs} \mlabel{ex:8.6.4} 
(a) For every pro-nilpotent Lie algebra $\g$, the corresponding group 
$(\g,*)$ defined by the Hausdorff series 
is a BCH--Lie group (Remark~\ref{rem:bch-liealg}(c)). 

(b) Every  Banach--Lie group is BCH (Theorem~\ref{thm:loc-bangrp}). 

(c) If $\cA$ is a unital Mackey complete cia, then its unit group 
$\cA^\times$ is BCH (Theorem~\ref{thm:IV.1.11}). Since 
$M_n(\cA)$ also is a Mackey complete cia (Proposition~\ref{prop:12.2.3}), 
the group $\GL_n(\cA) = M_n(\cA)^\times$ is also BCH. 

(d) Let $K$ be a BCH--Lie group and 
$q \: P \to M$ a smooth $K$-principal bundle over the 
$\sigma$-compact finite-dimensional 
manifold $M$. Then 
the group $\Gau_c(P)$ of compactly supported gauge transformations 
is BCH (Examples~\ref{exs:5.1.4}(c)). 
\end{exs}

The following result is quite useful to show that certain Lie algebras are not BCH: 

\begin{thm}
  \mlabel{thm:IV.1.7} {\rm(Robart's Criterion; \cite{Rob04})} If $\g$ is a 
sequentially complete BCH--Lie algebra, then there exists a 
$0$-neighborhood $U \subeq \g$ such that 
$f(x,y) := \sum_{n = 0}^\infty (\ad x)^n y$
converges and defines an analytic function on $U \times \g$. 
\end{thm}

\begin{ex}
  \mlabel{ex:VI.2.10} The formal analog of the Lie algebra $\gs_1(\R)_1$ 
of germs of diffeomorphisms 
of the form $x + a_2 x^2 + \cdots$ in $0$ 
is the Lie algebra $\gf_1(\R)_1$ of formal diffeomorphisms of 
this form, which is pro-nilpotent, 
hence in particular BCH (Example~\ref{ex:8.6.4}(a)). 
In contrast to this fact, Robart observed 
that $\gs_1(\R)_1$ is not BCH. In fact, for the elements 
$\xi(x) = ax^2$, $a \in \R,$ and $\eta(x) = x^3$, we have 
$$\sum_{n = 0}^\infty \big((\ad \xi)^n \eta\big)(x) 
= x^3 \sum_{n = 0}^\infty a^n n! x^n, $$
which converges for no $x \not=0$ if $a \not=0$. 
With Floret's results from \cite[p.~155]{Fl71}, it follows that 
this series does not converge in the Silva space $\gs_1(\R)_1$, so that 
Theorem~\ref{thm:IV.1.7} shows that $\gs_1(\R)_1$ is not BCH. 
\end{ex}

To discuss the following examples, we recall the 
concept of $\mu$-regularity from \cite{OMY82, OMY83a}: 
\begin{defn} \mlabel{def:Ne06-III.2.4} Let $G$ be a Lie group with Lie algebra $\g$ and 
\[ \Delta := \{t_0, \ldots, t_m\} \]  a division of the real 
interval $J := [a,b]$ with 
$a = t_0$ and $b = t_m$. We write 
\[ |\Delta| := \max\{t_{j+1} -t_j \; j = 0,\ldots, m-1\}. \] 
For $|\Delta| \leq \eps$, a pair $(h, \Delta)$ is called a 
{\it step function} \index{step function} 
on $[0,\eps] \times J$ 
if the map $h \: [0,\eps] \times J \to G$ satisfies 
\begin{enumerate}
\item[\rm(1)] $h(0,t) = \be$ for all $t \in J$ and all maps 
$h^t(s) := h(s,t)$ are $C^1$. 
\item[\rm(2)] $h(s,t) = h(s,t_j)$ for $t_j \leq t < t_{j+1}$. 
\end{enumerate}

\nin For a step function $(h,\Delta)$, we define the 
\index{product integral} 
{\it product integral} 
$\prod_a^t (h,\Delta) \in G$ by 
$$ \prod_a^t (h,\Delta) := h(t-t_k, t_k) h(t_k-t_{k-1}, t_{k-1})\cdots h(t_1 - t_0, t_0)
\quad \hbox{ for } \quad t_k \leq t < t_{k+1}. $$ 

Now let $(h_n,\Delta_n)$ be a sequence of step functions with $|\Delta_n| \to 0$ 
for which the sequence $(h_n, {\partial h_n \over \partial s})$ converges 
uniformly to a pair $(h, {\partial h \over \partial s})$ for a 
function $h \: [0,\eps] \times J \to G$. 
Then the limit function $h$ is a 
\index{$C^1$-hair} 
{\it $C^1$-hair in $\be$}, 
i.e., it is continuous, differentiable with respect to $s$, and 
${\partial h \over \partial s}$ is continuous on $[0,\eps] \times J$. 

The Lie group $G$ is called 
\index{Lie group!$\mu$-regular}
{\it $\mu$-regular}
\begin{footnote}{$\mu$ stands for ``multiplicative''.}  \end{footnote} 
(called ``regular'' in 
\cite{OMY82, OMY83a}), if the product integrals $\prod_a^t (h_n,\Delta_n)$ converges  
uniformly on $J$ for each sequence $(h_n, \Delta_n)$ converging in the sense explained 
above to some $C^1$-hair in $\be$. 
Then the limit is denoted $\prod_a^t (h, d\tau)$ 
and called the {\it product integral of $h$}.
\end{defn}

\begin{ex}
  \mlabel{ex:IV.1.14}  (Formal diffeomorphisms) (a) Important examples of 
pro-nilpotent BCH Lie groups are the Fr\'echet--Lie groups 
$\Gf_n(\K)$ of formal diffeomorphisms (cf.~Subsection~\ref{subsec:sternberg}, 
Section~\ref{sec:formgrp} and Example~\ref{ex:VI.2.8}).

These groups are $\mu$-regular Lie groups \cite{Omo80}): 
In view of the semidirect decomposition and the fact that 
$\mu$-regularity is an extension property (\cite[Thm.~V.1.8]{Ne06a}), 
it suffices to observe that pro-nilpotent Lie groups are $\mu$-regular, which follows 
by an easy projective limit argument. 

(b) For $\K = \R$, the identity component 
$\Gf_1(\R)_0$ is exponential and analytic, but not BCH. 
For $n \geq 2$, the group $\Gf_n(\R)$ is analytic, but not locally exponential. 
If a subgroup $H \subeq \GL_n(\R)$ consists of matrices with real eigenvalues, 
then the subgroup $\Gf_n(\R)_1 \rtimes H \subeq \Gf_n(\R)$ 
is locally exponential (\cite[Thms.~6/7]{Rob02}).
\end{ex}

\begin{rem}
  \mlabel{rem:VI.2.13} Let $\g(A)$ be a symmetrizable Kac--Moody Lie algebra. 
In \cite{Rod89}, E.~Rod\-ri\-guez--Carrington describes certain 
Fr\'echet completions of $\g(A)$, 
including a smooth version 
$\g^\infty(A)$ and an analytic version~$\g^\omega(A)$, 
which are BCH--Lie algebras (\cite[Prop.~1]{Rod89}). 
Corresponding groups are constructed for the unitary real forms by unitary 
highest weight modules of $\g(A)$, as integral subgroups of the 
unitary group of a Hilbert space (cf.~Theorem~\ref{thm:5.5.3}). 
In \cite{Su88}, Suto obtains closely related results, 
but no Lie group structures. 

In a different direction, Leslie describes in \cite{Les90} 
a certain completion $\oline\g(A)$ of $\g(A)$ which leads to a Lie group 
structure on the path space $C^\infty([0,1],\oline\g(A))$, corresponding to the 
pointwise Lie algebra structure on this space. One thus obtains an integrable 
Lie algebra extension of $\oline\g(A)$ in the spirit of pre-integrable Lie algebras 
(see Remark~\ref{rem:pre-integrable}). 
For an approach to Kac--Moody groups in the context of diffeological 
groups, we refer to \cite{Les03} (see also~\cite{So84}) and for an approach 
in the context of algebraic groups we refer to S.~Kumar's monograph~\cite{Kum02}.
\end{rem}

\subsection*{Complexifications} 

The setting of BCH--Lie groups is the natural one for complexifications 
because if $\g$ is a locally exponential Lie algebra for which $\g_\C$ 
is locally exponential as a complex Lie algebra, then the local multiplication 
is complex analytic. This implies that $\g_\C$ is BCH which in turns entails that 
$\g$ is BCH.

\begin{defn}  
If $G$ is a real BCH--Lie group,  then a morphism of real Lie groups 
$\zeta_G \: G \to G_\C$ to a complex BCH--Lie group $G_\C$ is called a 
\index{universal complexification, of BCH--Lie group} 
{\it universal complexification} if for each morphism 
$\alpha \: G \to H$ to a complex BCH--Lie group $H$, there exists a unique morphism 
$\beta \: G_\C \to H$ with $\alpha = \beta \circ \zeta_G$. 
\end{defn}

It is 
well known that for finite-dimensional Lie groups 
a universal complexification always exists  
(cf.\ \cite{Ho65} , \cite[Th.~XIII.5.6]{Ne00}), but it need not be locally injective, 
so that it may happen 
that $\dim_\C G_\C < \dim_\R G$. The following theorem shows 
that, due to the existence of non-enlargeable Lie algebras,  the situation 
becomes more complicated in infinite dimensions.
 
\begin{thm}
  \mlabel{thm:VI.1.23} {\rm(Existence of universal complexifications; \cite{GN03, Gl02c})} 
Given a real BCH--Lie group~$G$, 
let $N_G\trile G$ be the intersection of all kernels
of smooth homomorphisms from~$G$ to complex BCH--Lie groups. 
Then~$G$ has a universal complexification 
 if and only if $N_G$ is a BCH--Lie subgroup of~$G$ and the 
complexification of the quotient $\L(G)/\L^e(N_G)$ is enlargeable. 
\end{thm}

From Theorem~\ref{thm:quot-alg} one can derive that, if $G$ is $1$-connected, then 
the existence of a universal complexification is equivalent to the enlargeability 
of $\L(G)/\L^e(N_G)$. 
In \cite{GN03}, one finds an example of a Banach--Lie group for which
$N_G$ fails to be a Lie subgroup 
and also examples where $N_G=\{\be\}$ but $\L(G)_\C$ is not enlargeable. 

\subsection*{Open problems}

\begin{probl}
  \mlabel{prob:VI.8}  Prove an appropriate version of Theorem~\ref{thm:VI.1.23} 
on the existence of a 
universal complexification for locally exponential Lie algebras. 

This already is an interesting issue on the level of Lie algebras because the complexification 
of a locally exponential Lie algebra need not be locally exponential. 
In fact, in Example~\ref{ex:5.3.4b} we have seen an exponential Lie algebra $\g$ 
containing an unstable closed subalgebra $\h$. If $\g_\C$ is locally exponential, 
as a complex Lie algebra, then the local multiplication in $\g_\C$ is holomorphic, 
so that $\g$ is BCH, contradicting the existence of unstable closed subalgebras. 
\end{probl}

\section*{Notes on Section~\ref{sec:BCH}} 

The Lie group concept used in \cite{BCR81} 
is stronger than our concept of a BCH--Lie group 
because additional properties of the Lie algebra 
are required, namely that it is a so-called AE--Lie algebra, a property 
which encodes the existence of certain seminorms, compatible with the Lie bracket. 

\section{Lie algebras with regular exponential function} 
\mlabel{sec:banach-reg-exp}

For finite-dimensional Lie groups, the class of those 
for which the exponential function $\exp \: \g \to G$ is a 
global diffeomorphism is well understood. By 
the Dixmier--Saito Theorem, this is the case 
if and only if $G$ is $1$-connected, $\g$ is solvable 
and the operators $\ad x$, $x \in \g$, have no purely imaginary 
eigenvalues. If, in addition, $\g$ is complex, then 
$G$ must be a $1$-connected nilpotent Lie group. 

In this section we discuss various generalizations of these
results to Banach--Lie groups. 
From Remark~\ref{rem:ban-dexp} we know that 
$T_x(\exp)$ is invertible for every $x \in \g$ if and only if 
the operators $\ad x$ have no non-zero purely imaginary spectral 
values; below we call this condition (SC). 
In Subsection~\ref{subsec:9.8.1} we inspect condition (SC) 
and take a closer look at solvable Banach--Lie groups. 
In particular we show that all solvable Banach--Lie algebras 
are enlargeable to contractible Lie groups 
and these groups are exponential if and only if $\g$ satisfies~(SC) 
(Theorem~\ref{thm:dix}). 

In Subsection~\ref{subsec:9.8.2} we then turn to the condition 
that $\Spec(\ad x) = \{0\}$ for every $x \in \g$, i.e., all 
operators $\ad x$ are quasi-nilpotent; we then say that $\g$ is 
quasi-nilpotent. In the finite-dimensional 
case this is equivalent to the nilpotence of $\g$ by Engel's Theorem. 
For complex Banach--Lie algebras we prove W.~Wojty\'nski's remarkable 
theorem asserting that quasi-nilpotence is equivalent 
 to the convergence of the 
Hausdorff series on all pairs in \break {$\g \times \g$}, and also 
to exponentiality. In the real case the first two conditions 
are still equivalent and imply exponentiality, but there 
exist exponential real Lie algebras which are not nilpotent. 

We conclude this section with a brief discussion of exponentiality 
of projective limits of finite-dimensional Lie algebras 
in Subsection~\ref{subsec:9.8.3}. These Lie algebras are of particular 
importance because they arise in many formal constructions 
(see Example~\ref{ex:VI.2.8}).

\subsection{Solvable Banach--Lie algebras} 
\mlabel{subsec:9.8.1}

According to Theorem~\ref{thm:BCH-nil}, 
nilpotent Lie algebras are always exponential, but the converse 
is not true. In finite-dimensional Lie theory, the 
Dixmier--Saito Theorem \cite[Thm.~14.4.8]{HiNe12} 
characterizes exponential Lie algebras as 
those satisfying the 
{\it spectral condition}: \index{spectral condition}  
\begin{description}
 \item[\rm(SC)] \qquad $(\forall x \in \g) \quad \Spec(\ad x) \cap i \R \subeq \{0\}.$ 
\end{description}
This is why we start this section with a discussion of (SC) for 
Banach--Lie algebras. 

We shall need the following fact from spectral theory: 

\begin{thm} \mlabel{thm-spec-inh} {\rm(\cite[Thm.~10.18]{Ru91})} Let $\cB$ be a 
unital complex Banach algebra 
and $\cA \subeq \cB$ a closed subalgebra containing the identity. 
Then, for each $x \in \cA$, the spectrum $\Spec_\cA(x)$ with respect to $\cA$ 
is a union of $\Spec_\cB(x)$ and a collection of bounded components 
of $\C \setminus \Spec_\cB(x)$. 
\end{thm}

As an application, we obtain: 

\begin{lem} \mlabel{spec-subspace} Let $X$ be a Banach space and $Y \subeq X$ a closed subspace. 
Further, let $A \in {\cal L}(X)$ be a continuous linear map with 
$A(Y) \subeq Y$. Then 
$$ \Spec(A) \cap i \R \subeq \{0\} \quad \Rarrow \quad 
  \Spec(A\res_Y) \cap i \R \subeq \{0\}. $$
\end{lem}

\begin{prf} First we observe that 
$${\cal L}(X)_Y := \{ \phi \in {\cal L}(X) \: \phi(Y) \subeq Y \} $$
is a closed unital Banach subalgebra of ${\cal L}(X)$. 
The condition 
\[ \Spec(A) \cap i \R \subeq \{0\} \]  means that 
$\Spec_{{\cal L}(X)}(A)$ is contained in $\C \setminus i\R^\times$. 
Then $i\R^\times_+$ and $i\R^\times_-$ are contained in the unbounded component 
of the resolvent set of $A$, so that Theorem~\ref{thm-spec-inh} implies that 
$\Spec_{{\cal L}(X)_Y}(A)$ is also contained in $\C \setminus i\R^\times$. 

The map 
${\cal L}(X)_Y \to {\cal L}(Y), \phi \mapsto \phi\res_Y$
is a homomorphism of unital Banach algebras, so that 
$\Spec(A\res_Y) \subeq \Spec_{{\cal L}(X)_Y}(A) \subeq \C \setminus i\R^\times.$
\end{prf}

\begin{prop} \mlabel{prop:spec-cond} For a Banach--Lie algebra $\g$ 
satisfying {\rm(SC)}, the following assertions hold: 
\begin{description}
\item[\rm(1)] Each closed subalgebra of $\g$ satisfies {\rm(SC)}. 
\item[\rm(2)] Each quotient of $\g$ satisfies {\rm(SC)}. 
\item[\rm(3)] If $\g$ is finite-dimensional, then $\g$ is solvable. 
\end{description}
\end{prop}

\begin{prf} (1) Let $\h \subeq \g$ be a closed subalgebra and $x \in \h$. 
Then the operator $\ad x$ on $\g$ preserves $\h$, 
so that Lemma~\ref{spec-subspace} implies (1). 

(2) Let $x \in \g$ and $\fn \trile \g$ be a closed ideal. 
Then for each $\lambda \in \R^\times$ the operator 
$\ad x + i \lambda \be$ is invertible on $\g$ and preserves $\fn$. 
Moreover, (1) implies that its restriction to $\fn$ is invertible, so that 
it is also invertible on $\g/\fn$.   

(3) Suppose that $\g$ is a finite-dimensional Lie algebra which is 
not solvable. Then the Levi decomposition implies that 
it contains a non-trivial simple subalgebra $\fs$. 
Let $\fs = \fk \oplus \fp$ be a Cartan decomposition of $\fs$ 
(\cite[Thm.~13.1.7]{HiNe12})
and $0 \not= x \in \fk$. Then $\ad x$ is semisimple with 
$\Spec(\ad x) \subeq i\R$, so that our assumption leads to $\ad x = 0$. 
As $\fz(\fs) = \{0\}$, this contradicts $x \not=0$. 
\end{prf}

\begin{prop} \mlabel{prop:derived-homot} Let $G$ be a $1$-connected 
Banach--Lie group, $n \in \N$ and 
$\oline D^n(\g) := \oline{D^n(\g)}$ the closed derived series of $\g$. 
Then 
$\oline D^n(G) := \la \exp \oline D^n(\g)\ra$ is a Lie subgroup of $G$ and the inclusion map 
$\oline D^n(G) \into G$ is a weak homotopy equivalence. In particular, each 
subgroup $\oline D^n(G)$ is $1$-connected. 
\end{prop}

\begin{prf} We argue by induction on $n$. In view of 
$\oline D^n(\g) = \oline D^1(\oline D^{n-1}(\g))$ and the fact that a Lie subgroup 
of a Lie subgroup is a Lie subgroup, it suffices to deal with the case $n = 1$. 

First we note that $\fa := \g/\oline D^1(\g)$ is an abelian Banach--Lie algebra, 
which we consider as an abelian Lie group. Since $G$ is $1$-connected, 
there exists a unique Lie group 
homomorphism $q \: G \to \fa$ for which $\L(q) \: \g \to \fa$ is the quotient map. 
Now $N := \ker q$ is a normal Lie subgroup of $G$ with 
Lie algebra $\oline D^1(\g),$ and $G/N \cong (\fa,+)$ (Proposition~\ref{prop:5.5.9}). 
Since $\fa$ is simply connected, $N$ is connected, which implies that $N = \oline D^1(G)$. 

It follows in particular that $G/\oline D^1(G)$ is contractible, and since 
$q$ defines a topologically locally trivial fiber bundle 
by Michael's Selection Theorem~\cite{MicE59}, 
the long exact homotopy sequence (Theorem~\ref{homseq-princ}) 
implies 
that all homomorphisms $\pi_k(\oline D^1(G)) \to \pi_k(G)$ are isomorphisms. 
This means that that the inclusion of $\oline D^1(G)$ is a weak homotopy equivalence.
In particular, $\pi_1(\oline D^1(G)) \cong \pi_1(G)$ vanishes, so that 
$\oline D^1(G)$ is $1$-connected. 
\end{prf}

\begin{cor} \mlabel{cor:solvcontr} All 
homotopy groups of a $1$-connected solvable Banach--Lie group 
vanish.   
\end{cor}

\begin{prf} Let $G$ be a simply connected solvable Banach--Lie group. 
In view of Proposition~\ref{prop:derived-homot}, 
for each $n \in \N$, the inclusion $\oline D^n(G) \to G$ is a weak homotopy  
equivalence. Since $G$ is solvable, we have $\oline D^N(G) = \{\be\}$ for 
some $N$, and the assertion follows. 
\end{prf}

\begin{rem}
In view of Palais' Theorem~\ref{pal-contr}, 
the preceding proposition even implies that 
all $1$-connected solvable Banach--Lie groups are contractible. 
\end{rem}

\begin{thm} \mlabel{thm:solv-int} Each solvable Banach Lie algebra $\g$ is enlargeable. 
\end{thm}

\begin{prf} Let $\tilde G_{\rm ad}$ be a $1$-connected Lie group with Lie algebra $\g/\z(\g)$ 
(Theorem~\ref{thm:5.3.8}). 
In view of the solvability of $\g/\z(\g)$, Corollary~\ref{cor:solvcontr} 
implies that $\pi_2(G_{\rm ad})$ vanishes, and the assertion follows from 
Corollary~\ref{cor:pi2-enl}. 
\end{prf}

\begin{prop} \mlabel{prop:solv-exp-sur} 
If $\g$ is a solvable Banach--Lie algebra satisfying 
{\rm(SC)} and $G$ a connected Lie group with Lie algebra $\g$, 
then $\exp_G$ is surjective. 
\end{prop}

\begin{prf} We argue by induction on 
the solvable length of $\g$. It clearly suffices to assume 
that $G$ is simply connected. Let $\fa := D^k(\g) \trile \g$ be the last non-zero 
term of the derived series of $\g$. 
Then $\fa \trile \g$ is an abelian ideal of $\g$  
and $\g/\fa$ is a solvable Lie algebra of length $k$, whereas the length of $\g$ is 
$k+1$. 

If $\g$ is abelian, then $\tilde G \cong (\g,+)$ implies that the 
exponential function of $G$ is surjective. We now argue inductively 
by the solvable length of $\g$. The abelian case corresponds to solvable length~$1$. 
Using Theorem~\ref{thm:solv-int}, we obtain a $1$-connected Lie group 
$H$ with $\L(H) = \g/\fa$. Since the quotient $\g/\fa$ also satisfies 
the spectral condition (Proposition~\ref{prop:spec-cond}), 
our induction hypothesis implies that the exponential function of 
$H$ is surjective. Next we recall that there exists a unique homomorphism 
of Lie groups $\phi \: G \to H$ for which $\L(\phi) \: \g \to \g/\fa$ is the 
quotient homomorphism. Its kernel $A := \ker \phi$ is a normal Lie 
subgroup with $H \cong G/A$. From the surjectivity of $\exp_H$ we now obtain 
for each $g \in G$ and $x \in \g$ with $gA = \exp_H(x + \fa) = \exp_G(x)A$. 
This means that $g \in \exp_G(x)A$, so that it remains to show that 
the subgroup $\exp_G(\R x) A$ of $G$ is contained in the image of the exponential 
function. 

As $G/A \cong H$ is $1$-connected, the subgroup $A$ is connected. 
Moreover, it is abelian, so that $\exp_A \: \fa \to A$ is surjective. 
We may therefore assume that $g \not\in A$, which implies that 
$x \not\in \fa$. Then $\fb := \fa + \R x$ is a closed subalgebra of $\g$, 
isomorphic to $\fa \rtimes_D \R$, where $D := \ad x\res_\fa$. 
We conclude that with $\alpha(t) := e^{tD}$ we obtain a
simply connected group 
$B := \fa \rtimes_\alpha \R$
with Lie algebra $\fb$. The exponential function of $B$ is given explicitly by 
$$ \exp_B(v,t) = (\beta(t)v, t), 
\quad \mbox{ where } \quad \beta(t) = \int_0^1 e^{-st D}\,ds 
= \frac{\1- e^{-t\ad x}}{t \ad x}$$
(Proposition~\ref{prop:semdir-exp}). 
From the spectral condition and Proposition~\ref{prop:spec-cond} 
we conclude that the operator $\beta(t)$ is invertible 
for each $t \in \R$, which implies the surjectivity of $\exp_B$. 

Now let $j_B \: B \to G$ be the unique Lie group homomorphism for which 
$\L(j_B) \: \fb \to \g$ is the natural inclusion map. 
Then 
$$g\in A \exp_G(\R x) = j_B(B) = j_B(\exp_B(\fb)) = \exp_G(\fb)$$ 
implies that $g \in \exp_G(\g)$. 
\end{prf}

\begin{rem} Let $G$ be a simply connected Banach--Lie group whose Lie 
algebra $\g$ satisfies (SC). Then 
$\kappa_\g(x)$ is invertible for each $x \in \g$ 
(Remark~\ref{rem:ban-dexp}), which implies that 
$\exp_G \: \g \to G$ is a local diffeomorphism in each $x \in\g$. 
Hence $\im(\exp_G)$ is an open subset of~$G$. 

Moreover, $\exp_G(x) = \exp_G(y)$ implies $\exp_G(x-y) = 0$ 
(Proposition~\ref{prop:exp-noninj}), which 
leads to $x - y \in \z(\g)$ because of the spectral condition. 
Let 
$$\Gamma_\z := \{ x \in \z(\g) \: \exp x = \be\} \cong \pi_1(Z(G)).$$ 
Then $\Gamma_\z$ is a discrete subgroup of $\z(\g)$ and the exponential 
map 
$$ \exp_G \: \g \to \im(\exp_G) $$
is a covering map with $\Gamma_\z$ acting additively as deck transformations. 
\end{rem}

The following theorem is a Banach variant of the classical Dixmier--Saito 
Theorem for finite-dimensional Lie groups \cite[Thm.~14.4.8]{HiNe12}. 

\begin{thm} \mlabel{thm:dix}
Let $G$ be a connected 
Banach--Lie group with Lie algebra~$\g$. For the assertions 
\begin{description}
\item[\rm(1)] For each $x \in \g$ we have $\Spec(\ad x) \cap i\R = \{0\}$. 
\item[\rm(2)] The exponential 
function $\exp_G \: \g \to G$ has no singular points. 
\item[\rm(3)] $\g$ is exponential. 
\item[\rm(4)] The exponential function $\exp_{\tilde G} \: \g \to \tilde G$ is a 
diffeomorphism. 
\end{description}
we have the implications 
\[  (1) \Leftrightarrow (2) \Leftarrow   (3) \Leftrightarrow (4).  \] 
If $\g$ is finite-dimensional or solvable, then all assertions are equivalent.
\end{thm}

\begin{prf} (1) $\Leftrightarrow$ (2): 
We know from Proposition~\ref{thm:exp-logder} that the logarithmic derivative 
of the exponential function is given by 
\[ \delta(\exp)_x = \kappa_\g(x) = 
\int_0^1 e^{-t\ad x}\, dt = \frac{\1 - e^{-\ad x}}{\ad x}.\] 
By the Spectral Mapping Theorem, 
this operator is invertible if and only if $\Spec(\ad x) \cap 2\pi i \Z \subeq \{0\}$. 
Therefore the exponential function of $G$ has no singular points if and only if 
(1) is satisfied. 

(3) $\Leftrightarrow$ (4): Since two $1$-connected Banach--Lie groups 
with Lie algebra $\g$ are isomorphic, (3) and (4) are equivalent. 

(4) $\Rarrow$ (2) is trivial. 

(2) $\Rarrow$ (4): If $\g$ is finite-dimensional, then 
Propoition~\ref{prop:spec-cond}(3) implies that $\g$ is solvable. 

Suppose that $\g$ is a solvable Banach--Lie algebra. 
Proposition~\ref{prop:solv-exp-sur} implies that $\exp_G$ is surjective. 
Replacing $G$ by $\tilde G$, we may assume that $G$ is $1$-connected. 
It remains to see that $\exp_G$ is injective. In view of the regularity of 
$\exp_G$, the relation $\exp_G x = \exp_G y$ implies $\exp_G(x -y) =\be$ 
(Proposition~\ref{prop:exp-noninj}), 
so that it suffices to show that $\exp_G^{-1}(\be) = \{0\}$. 

Let $x \in \g$ with $\exp_G x = \be$. 
If $x \not=0$, then there exists a maximal $n$ with $x \in \oline D^n(\g)$. 
Then the subgroup $\oline D^n(G) \subeq G$ is $1$-connected by 
Proposition~\ref{prop:derived-homot}. Replacing $G$ with $\oline D^n(G)$, we may 
therefore assume that $n = 0$, i.e., $x \not\in \oline{[\g,\g]}$. 
Consider the quotient map 
$p \: \g \to \fq := \g/\oline{[\g,\g]}$. 
As the group $G$ is $1$-connected, there is a corresponding group 
homomorphism $p_G \: G \to (\fq,+)$. This leads us to 
$0 = p_G(\exp_G x) = \exp_\fq q(x) = q(x)$, contradicting the assumption 
$x\not\in\oline D^1(\g)$.
\end{prf}

\begin{cor} A finite-dimensional 
complex Lie algebra is exponential if and only it is nilpotent. 
\end{cor}

\begin{prf} We have already seen that any nilpotent Lie algebra is exponential. 
According to Engel's Theorem \cite{HiNe12}, it suffices to show that each 
$\ad x$, $x \in \g$, is nilpotent. If this is not the case, 
then there exists an $x \in \g$ such that $\ad x$ has a non-zero eigenvalue 
$\lambda$. Then $i \in \Spec(\ad (i\lambda^{-1}x)) \cap i\R$. 
In view of the Dixmier--Saito Theorem, this contradicts the exponentiality of $\g$. 
\end{prf}

\subsection{Quasi-nilpotent Banach--Lie algebras} 
\mlabel{subsec:9.8.2}

The main goal of this subsection is to prove 
Wojty\'nski's Theorem characterizing Banach--Lie algebras 
for which the Hausdorff series converges everywhere in terms 
of the quasi-nilpotence of all operators $\ad x$. 

We start with some preliminary remarks on plurisubharmonic functions. 

\index{plurisubharmonic function}
\begin{defn} (cf.~\cite[Def.~4.1.3]{He89}) Let $X$ be a complex locally convex space 
and $\Omega \subeq X$ be an open subset. 
A function $u \: \Omega \to [-\infty,\infty[$ is called 
{\it plurisubharmonic}\begin{footnote}{In \cite{He89} these functions are called plurihypoharmonic 
and such a function $u$ is called plurisubharmonic if $u^{-1}(-\infty)$ has no 
interior points.}   
\end{footnote}
if it satisfies the following conditions: 
\begin{description}
\item[\rm(PS1)] $u$ is upper semicontinuous, i.e., for any $t \in \R$, 
the set 
$\{ z \in \Omega \: u(z) < t \}$ is open. 
\item[\rm(PS2)] For $a \in \Omega$ and $b \in X$ with $a + \zeta b \in \Omega$ for 
$|\zeta| \leq 1$, 
\[ u(a) \leq \frac{1}{2\pi} \int_{-\pi}^\pi u(a + e^{i\theta}b)\, d\theta.\] 
\end{description}
\end{defn}

\begin{rem} \mlabel{rem:all-minusinf} 
If $u \: \Omega \to [-\infty,\infty[$ is plurisubharmonic, 
$\Omega$ is connected and $u^{-1}(-\infty)$ has an interior point, then 
$u$ is constant $-\infty$ by \cite[Prop.~4.1.2(a)]{He89}. 
\end{rem}

\begin{lem} \mlabel{lem:psh-vanish} 
Let $f \: \C^2 \to [0,\infty[$ be homogeneous in the sense that  
there exists a $\sigma \in \R$ with 
\[ f(\lambda z) = |\lambda|^\sigma f(z) \quad \mbox{ for } \quad z \in \C^2, 
\lambda \in \C.\] 
If $f$ is plurisubharmonic and vanishes on $\R^2$, then 
$f = 0$. 
\end{lem}

\begin{prf} Any element $z \in \C^2$ can be written in the form 
$z = a + ib$ with $a,b \in \R^2$. If $a,b$ are linearly dependent, 
then $z \in \C \cdot \R^2$ and thus $f(z) =0$ because $f$ is homogeneous. 
If $a,b$ are linearly independent, we choose them as a new basis, so that 
it suffices to show that $f(1,i) = 0$. 
For $|\zeta| = 1$, 
\[ (1,i) + \zeta(1,-i) = (1 + \zeta, i(1 - \zeta)) \in \C \cdot \R^2 \] 
because $i(1-\zeta)/(1 + \zeta) \in \R$ for $\zeta \not=-1$. 
Therefore the subharmonic function $h \: \C \to [0,\infty[$ defined by 
$h(\zeta) := f(1 + \zeta, i(1-\zeta))$ vanishes on the unit circle, 
so that $h(\zeta) = 0$ for $|\zeta| \leq 1$ follows from the mean value 
property (PS2) and $h \geq 0$. We obtain in particular 
that $f(1,i) = h(0)= 0$. 
\end{prf}

We shall need the following result due to E.~Vesentini 
(cf.~\cite[Thm.~2]{Ves68}). 

\index{Vesentini Theorem} 
\begin{thm}[Vesentini] \mlabel{thm:vesentini} If $\cA$ is a unital complex Banach algebra, 
then the spectral radius 
\[ \rho \: \cA \to [0,\infty[, \qquad 
\rho(a) = \lim_{n \to \infty} \|a^n\|^{1/n} \] 
and the function $\log \rho$ are plurisubharmonic. 
\end{thm}

\begin{prf} In view of \cite[Ex.~4.1.5(d)]{He89}, it suffices to show that 
$\log \rho$ is plurisubharmonic. If $\rho(a) = 0$, then 
$\rho(a + t \1) = |t|$ for every $t \in \C$. 
This shows (PS3). 
Further, (PS1) follows from the 
lower semicontinity of the spectrum \cite[Thm.~5.2.3]{HP57}. 

Note that 
\[ \log(\rho(a)) =\lim_{n \to \infty} \frac{1}{2^n} \log \|a^{2^n}\|\] 
and that the functions $f_n(a) := \frac{1}{2^n} \log \|a^{2^n}\|$ 
are plurisubharmonic because the power functions $a \mapsto a^{2^n}$ 
are holomorphic (\cite[Ex.~4.1.5(c)]{He89}). We also have 
$f_{n+1} \leq f_n$, so that $\log \rho = \lim_n f_n$ is plurisubharmonic  
by (\cite[Prop.~4.1.2(d)]{He89}). 
\end{prf}

\begin{defn} A Banach--Lie algebra $\g$ is called 
\index{Banach--Lie algebra!quasi-nilpotent} 
if 
\[ \Spec(\ad x) = \{0\} \quad \mbox{ for every } \quad x \in\g,\]
 i.e., if all operators $\ad x$ are quasi-nilpotent.
\end{defn}

\index{Wojty\'nski's Theorem} 
\begin{thm}[Wojty\'nski's Theorem] \mlabel{thm:wowo} 
For a complex Banach--Lie algebra 
$\g$, the following are equivalent: 
\begin{description}
\item[\rm(a)] $\g$ is quasi-nilpotent. 
\item[\rm(b)] All operators $\kappa_\g(x) = \int_0^1 e^{-t\ad x}\, dt 
= \frac{\1- e^{-\ad x}}{\ad x}$ are invertible. 
\item[\rm(c)] $\g$ satisfies (SC). 
\item[\rm(d)] $\g$ is exponential. 
\item[\rm(e)] The Hausdorff series 
$x * y = x + y + \frac{1}{2}[x,y] + \cdots 
= \sum_{n = 1}^\infty H_n(x,y),$ 
where $H_n(x,y)$ is the homogeneous term of order $n$, 
converges in any two elements $x,y \in \g$. 
\end{description}
\end{thm}

\begin{prf} (a) $\Leftrightarrow$ (b)  $\Leftrightarrow$ (c):  
In Remark~\ref{rem:ban-dexp} we have seen that $\kappa_\g(x)$ 
is not invertible if and only if $\Spec(\ad x)$ contains an element of the 
form $2\pi i n$ with $n \in \Z$ non-zero. 
Since $\g$ is complex, 
this never happens for any $x \in \g$ if and only if $\Spec(\ad x) = \{0\}$ 
for every $x \in \g$ and this in turn is equivalent to~(SC). 

(b) $\Rarrow$ (d): For $x,y \in \g$, consider the holomorphic function 
\[ F \: \C \to \cL(\g), \quad F(z) := e^{z\ad x}e^{z\ad y} - \1.\] 
If $z$ is sufficiently small, then 
\[ F(z) = e^{(z \ad x) * (z \ad y)} - \1 = e^{\ad(z x * zy)} - \1,\] 
and since $\g$ is quasi-nilpotent, $\Spec(F(z)) = \{0\}$ for $z$ close to $0$. 
Since $F$ is holomorphic, the function 
\[ f \: \C \to [-\infty,\infty[, \quad 
f(z) := \log\rho(F(z)) \] 
is plurisubharmonic by Theorem~\ref{thm:vesentini} and the fact that 
composing a plurisubharmonic function with a holomorphic function yields 
a plurisubharmonic function (\cite[Thm.~4.3.1]{He89}). 
We have already seen that $f^{-1}(-\infty)$ is a neighborhood of $0$, 
so that Remark~\ref{rem:all-minusinf} implies that $f$ is constant 
$-\infty$. This shows that $\Spec(F(z)) = \{0\}$ for every $z$, so that 
$\Spec(e^{\ad x}e^{\ad y}) = \{1\}$ for $x,y \in \g$. 

In Subsection~\ref{subsec:5.4.4} we have seen that, for the real power series  
\[ \Psi(X)=\frac{(X+1)\log (X+1)}{X} 
= (X+1)\sum_{k=0}^\infty \frac{(-1)^{k}}{k+1}X^k \in \R[[X]] \] 
with radius of convergence $1$, and for $\|x\| + \|y\| < \log 2$, we have  
\begin{equation}
  \label{eq:bch-quasinil}
 x * y  = x + \int_0^1 \Psi(e^{\ad x}e^{t\ad y}-\1) y\, dt.
\end{equation}
Since 
\[ \g \times \g \to \Omega := \{ a \in \cL(\g) \: \rho(a) < 1 \}, 
\quad (x,y) \mapsto e^{\ad x} e^{\ad y} - \1 \] 
is a holomorphic map and insertion of $a$ in $\Psi(a)$ is also holomorphic 
as a map $\Omega \to \cL(\g)$ 
(Theorem~\ref{calcanalyt}), 
we see that the right hand side of \eqref{eq:bch-quasinil} 
defines a holomorphic function $\g \times \g \to \g$. 
We thus obtain an exponential complex Banach--Lie group $(\g,*)$. 

(d) $\Rarrow$ (e): The Taylor series of the holomorphic multiplication 
$x * y$ in~$0$, which is the Hausdorff series by Theorem~\ref{bch-locexp}, 
converges everywhere by Corollary~\ref{corcxmap}. 

(e) $\Rarrow$ (b): If the Hausdorff series converges in every pair 
$(x,y) \in \g^2$, then \cite[Thm.~5.2]{BS71b} implies that the so obtained 
function $x * y$ is analytic, so that we obtain an exponential 
Lie group $(\g,*)$. 

(d) $\Rarrow$ (c) If $(\g,*)$ is an exponential Lie group, 
then Proposition~\ref{prop:5.2.10} implies that 
$\kappa_\g \in \Omega^1(\g,\g)$ is the corresponding Maurer--Cartan 
form, hence invertible in every $x \in \g$. 
\end{prf}

\begin{thm} \mlabel{thm:real-wowo}
For a real Banach--Lie algebra $\g$, the following are equivalent: 
\begin{description}
\item[\rm(a)] $\g$ is quasi-nilpotent. 
\item[\rm(b)] $\g_\C$ is quasi-nilpotent. 
\item[\rm(c)] The Hausdorff series 
$x * y = x + y + \frac{1}{2}[x,y] + \cdots 
= \sum_{n = 1}^\infty H_n(x,y)$ 
converges in any two elements $x,y \in \g$. 
\end{description}
\end{thm}

\begin{prf} (a) $\Leftrightarrow$ (b): Clearly, (b) implies (a). To see the converse, 
we consider the function 
\[ F \: \g_\C \to [0,\infty[, \quad F(x) := \rho(\ad x)  
= \lim_{n \to \infty} \|(\ad x)^n\|^{1/n} \] 
which is plurisubharmonic by Vesentini's Theorem~\ref{thm:vesentini}.
For $x,y \in \g$, the function 
\[ f \: \C^2 \to [0,\infty[, \quad 
f(z,w) := F(z x + w y) \] 
is plurisubharmonic by \cite[Thm.~4.3.1]{He89}. 
By (a) it vanishes on $\R^2$, and since it is homogeneous, 
Lemma~\ref{lem:psh-vanish} implies that $f = 0$. This shows that 
$\Spec(\ad(x + i y)) = \Spec(\ad x + i \ad y) = \{0\}$ for all 
$x,y \in \g$, and this means that $\g_\C$ is quasi-nilpotent. 

(b) $\Rarrow$ (c) follows from Theorem~\ref{thm:wowo}. 

(c) $\Rarrow$ (b): If the Hausdorff series 
converges in every pair $(x,y) \in \g^2$, 
then \cite[Prop.~4.1]{BS71b} implies that the convergence is uniform 
on any compact subset contained in a finite-dimensional subspace. 
By \cite[Lemma~4.1]{BS71a}, this property is inherited by the natural 
extension of the homogeneous polynomials $H_n \: \g \times \g \to \g$ to 
maps $\g_\C \times \g_\C \to \g_\C$. Therefore the Hausdorff series 
also converges in every pair $(x,y) \in \g_\C^2$, and the assertion 
follows from Wojty\'nski's 
Theorem~\ref{thm:wowo}.\begin{footnote}
{Alternatively, one may argue directly with \cite[Prop.~8.2]{BS71b}.
}
\end{footnote}  
\end{prf}

\begin{rem} \mlabel{rem:quasi-nil} 
(a) If a real Banach--Lie algebra $\g$ is quasi-nilpotent, then 
the same argument as in the proof of Theorem~\ref{thm:wowo} implies that 
it is exponential and BCH. But there are real exponential Banach--Lie algebras 
which are not quasi-nilpotent. The smallest example is the two-dimensional 
non-abelian Lie algebra $\g = \aff(\R)$ of the affine group of the real 
line (Example~\ref{ex:5.1.9}). 

(b) If $\cA$ is a Banach algebra consisting of quasi-nilpotent elements, then 
the underlying Banach--Lie algebra is also quasi-nilpotent 
because, for each $a \in \cA$, we have 
$\Spec(\ad a) \subeq \Spec(a) - \Spec(a) = \{0\}$ 
(\cite[Prop.~3.3]{BelN08}). 

In this case $\cA^\times = \1 + \cA$ is an affine subspace of the unital 
Banach algebra $\cA_+$ on which the logarithm function is defined,  
and an inverse to the exponential function. 
Therefore $\cA^\times$ is an exponential Banach--Lie group. 
If $\cA$ is complex, Wojty\'nski's Theorem~\ref{thm:wowo} thus 
implies that the underlying Lie algebra is quasi-nilpotent. 
\end{rem}

\begin{ex} (a) On the Hilbert space 
$\cH = \ell^2(\N,\R)$ with the orthonormal basis 
$(e_n)_{n \in \N}$, we consider the Banach algebra $\cA$ 
of all compact operators $A$ with a strictly lower triangular 
matrix, i.e., $\la A e_m, e_n \ra = 0$ for $n \leq m$. 
Since all spectral values of $A$ are eigenvalues by the Spectral 
Theorem for Compact Operators, it follows that 
$\Spec(A) = \{0\}$ for $A \in \cA$. Therefore the underlying Banach--Lie algebra 
is quasi-nilpotent by Remark~\ref{rem:quasi-nil}. 

(b) On the Hilbert space 
$\cH = L^2([0,1],\C)$, we consider the Banach algebra $\cA$ generated by 
all integral operators of the form 
\[ (A_K f)(x) = \int_0^1 K(x,y)f(y)\, dy,\] 
where $K \in C([0,1]^2)$ and $K(x,y) = 0$ for $x \leq y$. 
Since all these operators are Hilbert--Schmidt with 
\[ \|A\|_2^2 = \int_0^1 \int_0^1 |K(x,y)|^2\, dx\, dy,\] 
they are in particular compact, so that $\cA$ consists of compact operators. 

If there exists an $\eps > 0$ such that $K(x,y)$ vanishes if $x \leq y- \eps$, 
then $(A_K f)(x)$ only depends on the values $f(y)$ for $y < x + \eps$, 
and this implies that $A_K$ is nilpotent, hence in particular quasi-nilpotent. 
Therefore $\cA$ contains a dense subspace of quasi-nilpotent elements, 
hence is quasi-nilpotent because the spectrum is continuous 
on the set of compact operator. This follows from the fact that 
their spectra are totally disconnected 
(\cite{New51}; see also \cite{CM79} for more general results). 
\end{ex}

\subsection{Exponential pro-finite Lie algebras} 
\mlabel{subsec:9.8.3}

In this subsection we briefly discuss 
exponentiality of projective limits of finite-dimensional Lie algebras 
(pro-finite Lie algebras). These Lie algebras are of particular 
importance because they arise many formal constructions 
(see Example~\ref{ex:VI.2.8}). 
 
\begin{prop} Let $\g$ be a pro-finite Lie algebra. Then 
$\g$ is exponential if and only if it is the projective limit 
of a system of finite-dimensional exponential Lie algebras. 

In particular, pro-nilpotent Lie algebras are exponential. 
\end{prop}

\begin{prf} We have already seen in Proposition~\ref{prop:prolim-exp} 
that projective limits of exponential Lie algebras are exponential. 

Assume, conversely, that $\g$ is exponential and that 
$\g = \prolim \g_j$, where each $\g_j$ is finite-dimensional. 
Since we may w.l.o.g.\ assume that the corresponding maps 
$\alpha_j \: \g \to \g_j$ have dense range, they are surjective. 
We claim that each Lie algebra $\g_j$ is exponential. 
For each $x \in \g$, the operator $\kappa_\g(x)$ is invertible, and we have 
\[  \alpha_j \circ \kappa_\g(x) = \kappa_{\g_j}(\alpha_j(x)) \circ \alpha_j 
\quad \mbox{ for } \quad  x\in \g.\]
Therefore $\kappa_{\g_j}(\alpha_j(x))$ is surjective for each 
$x \in \g$, and since $\g_j$ is finite-dimensional, 
$\kappa_{\g_j}(z)$ is invertible for each $z \in \g_j$. 
In view of the Theorem~\ref{thm:dix}, 
this means that $\g_j$ is exponential. 
\end{prf} 

In the preceding proposition we encounter the situation where we have 
a quotient map $q \: \g \to \h$ onto a finite-dimensional Lie algebra $\h$. 
Our argument above shows that $\h$ is exponential if $\g$ is exponential 
(Proposition~\ref{prop:derived-homot}). 
This in turn implies that $q$ is a morphism of exponential Lie algebras, 
and hence that $\fn := \ker q$ is an exponential ideal of finite codimension 
(Proposition~\ref{prop:invim}). 

The following example shows that if, conversely, $\fn \trile \g$ is a finite-codimensional 
exponential ideal of the pro-finite Lie algebra $\g$, 
then $\g$ need not be locally exponential. 

\begin{ex} (a) First we consider the profinite Lie algebra 
$$ \fs := (\C \rtimes_{\alpha_n} \R)^\N, $$
where $\alpha_n = 1 + in \in \End(\C)$. This means that 
$$ [(z_n,t_n), (z_n',t_n')] = (t_n \alpha_n z_n' - t_n' \alpha_n' z_n, 0). $$
Since all the Lie algebras $\C \rtimes_{\alpha_n} \R$ are exponential 
by Theorem~\ref{thm:dix}, their product $\fs$ is exponential. 

Now we form $\g := \fs \rtimes_D \R$, where $D \in \der(\fs)$ is given by 
$D(z_n,t_n) = (z_n, 0).$
Then 
$$ \g \cong \C^\N \rtimes \R^{\N_0}, $$
and, for each $n \in \N$, we have 
\[  [(0,t e_n -te_0), (e_n, 0)]  =\big( (t(1 + in) -t)e_n, 0\big) = (itn e_n,0), \]
showing that $(e_n,0)$ is an eigenvector for the eigenvalue $2\pi i$ for 
$\frac{2\pi}{n}(0,e_n - e_0)$. As this sequence converges to $0$ in $\g$, 
each $0$-neighborhood contains singular points of the exponential function, 
so that the Lie algebra $\g$ is not locally exponential. 

(b) We consider the profinite Lie algebra 
$$ \g := \C^\N \rtimes_D \R, $$
where $D$ acts on $\C^\N$ by 
$D e_n = in e_n$. Then $\C^\N$ is an abelian, hence exponential ideal of codimension $1$, 
 and in each point $(0,\frac{2\pi}{n}) \in \g$ the operator $\kappa_\g$ is singular. 
Therefore $\g$ is not locally exponential. 
\end{ex}

\section*{Notes on Section~\ref{sec:banach-reg-exp}}

{\bf Subsection~\ref{subsec:9.8.1}:} 
In Proposition~\ref{prop:solv-exp-sur} we encountered the difficulty 
to verify the surjectivity of the exponential map for certain solvable 
Banach--Lie groups. To understand the related difficulties, one considers 
for a connected Lie group $G$ the subsets 
\[ E^n(G) := (\exp \g)^n = 
\{ \exp X_1 \cdots \exp X_n \: X_1, \ldots, X_n \in \g\}, n \in \N, \] 
consisting of $n$-fold products of elements in the image of the exponential 
function. For finite-dimensional Lie groups, it is well-known that 
$G = E^1(G) = \exp \g$ for connected Lie groups whose Lie algebra 
is compact or nilpotent. But one easily finds examples of solvable 
connected Lie groups or simple Lie groups, such as $\SL_2(\R)$, for which 
the exponential function is not surjective 
(cf.~\cite{Wue03, Wue05}). Therefore it is a remarkable fact that 
$G = E^2(G)$ holds for all connected finite-dimensional Lie group 
(\cite{MS03}), i.e., every element is a product of two exponentials. 
The key method to obtain this result are decompositions similar to the 
Iwasawa decomposition which implies the assertion for reductive Lie groups. 
For infinite-dimensional Lie groups no general results of this kind are known. 
Already for the natural 
analogs of compact groups, such as the orthogonal group $\OO(\cH)$ 
of a real Hilbert space, may show pathological behavior, im particular, 
the exponential function of $\OO(\cH)$ is not surjective 
(see the notes on Chapter~\ref{ch:3}). 

\nin {\bf Subsection~\ref{subsec:9.8.2}:} 
Wojty\'nski's Theorem~\ref{thm:wowo} appears in~\cite{Woj98}, 
where it is also stated that it leads to Theorem~\ref{thm:real-wowo}. 

For a finite-dimensional Lie algebra $\g$, the commutator algebra 
$[\g,\g]$ is always nilpotent and $\ad x$ is nilpotent for 
$x \in [\g,\g]$. Therefore one can ask for similar results 
concerning quasi-nilpotency for Banach--Lie algebras. 
Interesting results in this direction have been obtained by 
Yu.~V.~Turovskii in \cite{Tu85}. In particular, he extends the spectral radius 
to bounded subsets of a Banach algebra and obtains generalizations of 
Vesentini's Theorem (Theorem~\ref{thm:vesentini}) in this context. 
Much of this theory is concerned with ``triagonalization'' of 
Banach--Lie algebras of operators, and thus generalizing Lie's Theorem. 
Extending the finite-dimensional representation theory of 
solvable Lie algebras to the Banach context, is one of the central 
themes of the monograph \cite{BS01} by D.~Belti\c t\u a and M.~Sabac. 

\section{Local exponentiality of semidirect products} 
\mlabel{sec:14.9}

In general semidirect products of locally exponential Lie groups 
are not locally exponential. In this section we describe various natural  
examples of semidirect products where this happens. 

\begin{ex} \mlabel{ex:15.9.b}
Let $E$ be a Mackey complete space and $\alpha$ a smooth $\R$-action on $E$.
We consider a semidirect product Lie group 
\[  G = E \rtimes_\alpha \R  \quad \mbox{ with the exponential function } \quad 
\exp_G(v,t) = (\beta(t)v, t) \] 
and $D = \alpha'(0)$. 
The multiplication in this group is given by 
\[  (v,t)(v',t') = (v + \alpha(t)v',t + t'), \quad \mbox{ so that }\quad 
T_\be(\lambda_{v,t})(v',t') = (\alpha(t)v', t'). \] 
We also have 
$$ T_{(v,t)}(\exp_G)(v',t') = (\beta(t)v' + \beta'(t)t', t'). $$
Therefore $\kappa_\g(v,t)$ is injective if and only if $\beta(t)$ is injective. 

The relation$\exp_G(v,t) = \exp_G(v',t')$ is equivalent to 
$t = t'$ and \break{$\beta(t)(v'-v) = 0$.} Therefore $\exp_G$ is injective in 
some $0$-neighborhood if and only if 
$\beta(t)$ is injective for sufficiently small $t$. 

In view of Proposition~\ref{prop:kerexp-locconv}, the kernel of 
$\beta(t)$ is topologically generated by the subspaces 
$$ \ker(t^2 D^2 + n^2 4 \pi^2\1), \quad n \in \N. $$
We write $\Spec_p(D)$ for the set of all eigenvalues of $D$. 
Therefore the condition $\beta(t)$ being injective is  
$$ \frac{-n^2 4 \pi^2}{t^2} \not\in \Spec_p(D^2), \quad \mbox{ resp., } \quad 
\frac{n 2\pi i }{t} \not\in \Spec_p(D), $$
for all $n \in \N$. This holds for $|t| \leq \eps$ if and only if 
\begin{equation}
  \label{eq:15.1.a}
  \Spec_p(D) \cap 2\pi i \big[\eps^{-1},\infty\big[ = \eset. 
\end{equation}

Since $\exp_G^{-1}(\be) = \{0\}$, an element $x = (v,t) \in \g$ is isolated in 
its fiber of $\exp_G$ if and only if $\kappa_\g(x)$ is injective, which is 
equivalent to the injectivity of $\beta(t)$ 
(Proposition~\ref{prop:exp-noninj}). 
\end{ex}

\begin{prop} 
Let $G$ be a Lie group with a smooth exponential function and 
assume that 
\begin{description}
\item[\rm(1)] $\L(G)$ is locally exponential, and 
\item[\rm(2)] the homomorphism $\exp_\z :=\exp_G\res_{\z(\g)} \: \z(\g) \to G$
has discrete kernel. 
\end{description}
Then there exists a $0$-neighborhood 
$U \subeq \g$ such that $\exp_G\res_U$ is injective. 
\end{prop}

\begin{prf} Let $H \subeq \g$ be an exponential local Lie group with Lie algebra 
$\g$. Then for each $x \in H$ the operator $\kappa_\g(x) = (\kappa_H)_x$ is invertible, 
which implies that $H \subeq \Reg_p(\g)$. 

Let $U \subeq \g$ be a balanced $0$-neighborhood with $U + U \subeq H$. 
If $\exp_G x = \exp_G y$ holds for $x,y \in U$, then Proposition~\ref{prop:exp-noninj} 
implies that $x -y \in \z(\g)$ with $\exp_G(x - y) = \be$. 
Therefore (2) implies that $\exp_G\res_U$ is injective. 
\end{prf}

\begin{probl} Let $G$ be a Lie group whose Lie algebra $\L(G)$ is locally exponential and 
which has a smooth exponential function $\exp_G \: \L(G) \to G$. 
Is the kernel of the homomorphism $\exp_Z := \exp_G\res_{\z(\L(G))} \: \z(\L(G)) \to G$ 
discrete?   
\end{probl}

\begin{prop} Let $M$ be a smooth finite-dimensional manifold and 
$X \in {\cal V}(M)$ a non-zero complete vector field with $\exp X = \id_M$. 
Then the Lie algebra $\g := C^\infty(M,\R) \rtimes_{{\cal L}_X} \R$ is 
not locally exponential. 
\end{prop}

\begin{prf} There exists a non-zero integer $k \in \Z$ and $f \in 
C^\infty(M,\C)$ with 
$\cL_Xf = 2\pi i k  f$. Then $\cL_X f^n = 2\pi i nk f$ for each $n \in \N$, so that 
$(0,\frac{1}{n}) \in \g$ is $\exp$-singular. 
\end{prf}
 
\begin{prop} Let $(M,\Omega)$ be a connected smooth finite-dimensional symplectic manifold and 
$X \in {\cal V}(M,\Omega)$ a non-zero complete symplectic vector field with $\exp X = \id_M$. 
Then the Lie algebra $\g := \ham(M,\omega) \rtimes_{{\cal L}_X} \R$ is 
not locally exponential. 
\end{prop}

\begin{prf} This follows by the same argument as above, because 
$\ham(M,\omega) \cong C^\infty(M,\R)/\R \1$. 
\end{prf}

\begin{probl} Let $M$ be a compact manifold, $X \in {\cal V}(M)$ and 
$\gamma_X(t) := \exp(tX)$ be the corresponding smooth $\R$-action on~$M$. 
This leads to a smooth action $\alpha$ on $E := C^\infty(M,\R)$. 
When is the group $G := E \rtimes_\alpha \R$ locally exponential? 

If $f \in E_\C$ is an $\cL_X$-eigenvector for some purely imaginary 
eigenvalue $i \mu$, then all non-zero positive multiples of 
$i\mu$ are eigenvalues of $\cL_X$ on $E_\C$, and this implies that 
$E \rtimes_{\cL_X} \R$ is not locally exponential. 
Therefore the main point is to show that the derivation 
$\cL_X$ on $E_\C$ always has a purely imaginary eigenvalue if 
$M$ is compact. 
%
\end{probl}

\begin{prop} If $\fk$ is a finite-dimensional Lie algebra, then 
the Lie algebra $\g := C^\infty(\R,\fk) \rtimes_D \R$ 
with $Df = f'$ is not locally exponential. 
\end{prop}

\begin{prf} For $x \in \fk$, consider $\xi \in \g$, defined by 
$\xi(t) := \cos(\omega t) x$. Then 
\[ D^2 \xi = - \omega^2 \xi,\]
 so that the assertion follows from the 
discussion in Example~\ref{ex:15.9.b}.   
\end{prf}

\begin{ex} 
If $\R$ acts ``linearly'' on the torus $\T^n$ with a dense orbit, 
then the representation of $\R$ on $C^\infty(\T^n,\R)$ 
extends to a representation 
of the torus $\T^n$, which implies the existence of 
$\lambda_1, \ldots, \lambda_n \in \R$ such that the subgroup 
\[ \Gamma := \Z \lambda_1 + \ldots + \Z \lambda_n \] of $\R$ has rank $n$ 
and $i\Gamma$ is the set of eigenvalues of $\cL_X$ on $C^\infty(\T^n,\R)$. 
\end{ex}

\appendix

\chapter{Tools from point set topology}\label{appA}
\chaptermark{Tools from point set topology\hfill \copyright{} H.
Gl\"{o}ckner and K.-H. Neeb}
We assume that the reader is familiar with the most basic concepts
of topology, like the notion of a (not necessarily open)
neighborhood, continuity of mappings,
and compactness.
In this appendix, we compile further topological concepts and results
which are important for our purposes.
Sections~\ref{secnets} to \ref{appcotop}
are fundamental for many parts of the book,
both for the theory and examples.
Section~\ref{app-basic-DL} on direct limits
is only used for specific examples
(like ascending unions of finite-dimensional Lie groups)
and for a related functional analytic tool (Proposition~\ref{silvahaveDL}).
\section{Nets}\label{secnets}
The notion of a net generalizes the notion
of a sequence. While topological questions concerning
metric spaces can usually be phrased
in the language of sequences, to deal with
arbitrary
topological spaces one has to use the
more general concept of a net.
\begin{defn}\label{defdirset}
Recall that a {\em directed set\/} is a non-empty
partially ordered set $(A,\leq)$
such that, for all $\alpha,\beta\in A$, there exists
$\gamma\in A$ such that $\alpha,\beta \leq
\gamma$.
\end{defn}
\begin{rem}
If $(A,\leq)$ is a directed set and
$F\sub A$ a finite subset,
then there exists $\beta \in A$ such that
$\alpha \leq \beta$ for all $\alpha\in F$,
as can be shown by a simple induction.
\end{rem}
\begin{defn}\label{defnnet}
A {\em net\/} in a set~$X$
is a family $(x_\alpha)_{\alpha\in A}$
of elements of~$X$, where $(A,\leq)$ is
a directed set.
\end{defn}
\begin{ex}
$(\N,\leq)$ is a directed set;
every sequence $(x_n)_{n\in\N}$ of elements $x_n\in X$
is a net in the given set~$X$.
\end{ex}
Thus, nets are generalizations of sequences,
obtained by replacing~$\N$
with more general directed sets.
Generalizing the definition
of convergence of sequences, one defines
convergence
of nets as follows:
\begin{defn}\label{defnetconv}
A net $(x_\alpha)_{\alpha\in A}$ in a topological space~$X$
is said to
{\em converge\/} to an element $x\in X$
(in which case we write $x_\alpha\to x$)
if for every neighborhood~$U$ of~$x$,
there exists $\alpha\in A$ such that
\[
x_\beta\in U \quad \mbox{for
all $\beta\in A$ such that $\beta\geq \alpha$.}
\]
\end{defn}
If~$X$ is Hausdorff ({\em i.e.\/},
if all elements $x\not=y$ in~$X$ have disjoint neighborhoods), then
a net $(x_\alpha)$ in~$X$ converges to at most
one point $x\in X$.\\[3mm]
Important examples of directed sets
are neighborhood filters in topological spaces.
We recall the relevant definitions.
\begin{defn}\label{deffilter}
A set
$\cF$ of subsets of a set~$X$
is called a {\em filter\/}
if the following conditions are satisfied:
\begin{description}[(F3)]
\item[(F1)]
$\,\emptyset\not\in \cF$;
\item[(F2)]
$\,U\in \cF$ and $U\sub V\sub X$
entail $V\in \cF$; and
\item[(F3)]
$\,U\cap V\in \cF$ for all
$U,V \in \cF$.
\end{description}
\end{defn}
\begin{defn}\label{defnbfil}
If $X$ is a topological space and $x\in X$,
then the set
$\cU_x(X)$ of all
(not necessarily open) neighborhoods of~$x$ in~$X$
is called the \emph{neighborhood filter} of~$x$.
A set $\cB$
of neighborhoods
of~$x$ is called a
\emph{basis of $x$-neighborhoods} if each neighborhood
of~$x$ contains some $U\in \cB$.
\end{defn}
It is clear
that $\cU_x(X)$ is a filter.
\begin{lem}\label{lemfi}
The neighborhood filter
$\cU_x(X)$
is a directed
set under inverse inclusion,
i.e., under the partial order defined
via
$U\leq V$ $:\aeq$ $V\sub U$.
\end{lem}
\begin{prf}
Clearly ``$\leq$''
is a partial order on $\cU_x(X)$.
If $U, V \in \cU_x(X)$,
then $W:=U\cap V \in \cU_x(X)$
with
$W \sub U, V$
and thus $U, V \leq W$.
Since furthermore
$\cU_x(X)$ is non-empty (as $X\in \cU_x(X)$),
we see that $(\cU_x(X),\leq)$ is a directed set.
\end{prf}
Here are two typical applications of nets in topology.
\begin{prop}\label{exemplify}
Let $X$, $Y$ be topological spaces,
$f\!:X\to Y$ be a mapping,
and $x\in X$.
Then~$f$ is continuous at~$x$
if and only if
$f(x_\alpha)\to f(x)$ for every net $(x_\alpha)_{\alpha\in A}$
in~$X$ which converges to~$x$.
\end{prop}
\begin{prf}
If $f$ is continuous at~$x$, let
$(x_\alpha)_{\alpha\in A}$ be a net in~$X$ converging to~$x$.
If $U$ is a neighborhood of~$f(x)$,
then $f^{-1}(U)$ is a neighborhood of~$x$;
since $x_\alpha\to x$,
there exists $\alpha\in A$ with $x_\beta\in f^{-1}(U)$ for
all $\beta\geq \alpha$. Thus $f(x_\beta)\in U$
for all $\beta\geq \alpha$. We have proved that
$f(x_\alpha)\to f(x)$.

If, conversely, $f$ is not continuous at~$x$,
we find a neighborhood~$U$ of~$f(x)$ in~$Y$
such that $f^{-1}(U)$ is not a neighborhood
of~$x$ in~$X$. Hence, for every neighborhood~$V$
of~$x$ in~$X$, we have $V\not\sub f^{-1}(U)$,
whence there exists $x_V\in V$ such that $f(x_V)\not\in U$.
Considering $\cU_x(X)$
as a directed set as in Lemma~\ref{lemfi},
$(x_V)_{V\in \cU_x(X)}$ becomes a net in~$X$,
and clearly $x_V\to x$,
since given any \break{$W\in \cU_x(X)$,}
we have $x_V\in V\sub W$ for all $V\geq W$.
Because $f(x_V)\not\in U$ for all~$V$,
the net $(f(x_V))_{V\in \cU_x(X)}$ does not converge to~$f(x)$.
The proof is complete.
\end{prf}
\begin{prop}\label{netclosure}
Let~$X$ be a topological space, $M$ be a subset
of~$X$, and $x\in X$. Then $x\in \wb{M}$ $($the closure
of~$M)$ if and only if there exists a net
$(x_\alpha)_{\alpha\in A}$ in~$M$ which converges
to~$x$ in~$X$.
\end{prop}
\begin{prf}
Let $x\in \wb{M}$.
For every neighborhood~$U$ of~$x$ in~$X$,
there exists an element $x_U\in U\cap M$.
Then $(x_U)_{U\in \cU_x(X)}$ is a net
in~$M$ which apparently converges to~$x$ in~$X$.

Conversely, let $x\in X$ and assume
that there is a net $(x_\alpha)_{\alpha\in A}$
in~$M$ which converges to~$x$. Then, for every
neighborhood~$U$ of~$x$, there exists $\beta\in A$
with $x_\alpha\in U$ for all $\alpha\geq \beta$.
In particular, $x_\beta\in U$. Since $x_\beta\in M$,
we have $U\cap M\not=\emptyset$.
This entails $x\in \wb{M}$.
\end{prf}
\begin{defn}
Let $(x_\alpha)_{\alpha\in A}$ be a net in a set~$X$, where $(A,\leq)$
is a directed set. A \emph{subnet} of $(x_\alpha)_{\alpha\in A}$
is a net of the form $(x_{f(j)})_{j\in J}$
where $(J,\leq)$ is a directed set and $f\colon J\to A$
a monotone map with the property that, for each $\alpha\in A$, there exists
$j_0\in J$ such that $\alpha\leq f(j_0)$. 
\end{defn}
If a net $(x_\alpha)_{\alpha\in A}$ in a topological space~$X$
converges to $x\in X$, then also each subnet of $(x_\alpha)_{\alpha\in A}$
converges to~$x$.
\section{Initial topologies; direct products}\label{secinitop}
\begin{defn}\label{definitop}
Let $X$ be a set and $(f_i)_{i\in I}$
be a family of mappings $f_i\colon X\to Y_i$
into topological spaces~$Y_i$.
Then there is a coarsest topology~$\cO$
on~$X$ making each $f_i$ continuous;
it is called the {\em initial topology
on~$X$ with respect to the family $(f_i)_{i\in I}$}.
\end{defn}
Indeed, the sets of the form
$\bigcap_{i\in F}f_i^{-1}(U_i)$, where $F\sub I$ is finite
and $U_i\sub Y_i $ is open,
form a basis for a topology~$\cO$
on~$X$ which has the desired
properties.
For later use, observe that we still get a basis
for the topology on~$X$ if we choose $U_i$ only in
a basis for the topology on~$Y_i$.
\begin{ex}\label{indutop}
If~$Y$ is a
topological space and $X\sub Y$
a subset, then the {\em induced topology\/}
on~$X$ is defined as the initial topology
with respect to the inclusion map $f\colon X\to Y$,
$f(x):=x$.
\end{ex}
\begin{ex}\label{exambije}
If $f\colon X\to Y$ is a bijection
and $Y$ is a topological space,
then the initial topology on~$X$
with respect to~$f$
makes $f$ a homeomorphism.
\end{ex}
\begin{ex}\label{defprodu}
Let $X:=\prod_{i\in I} X_i$
be the cartesian product
of a family $(X_i)_{i\in I}$
of topological
spaces,
with the canonical projections
$\pr_j\colon X\to X_j$,
$(x_i)_{i\in I}\mto x_j$
for $j\in I$.
The {\em product topology\/}
on~$X$ is defined as the initial topology
on~$X$ with respect to the family $(\pr_i)_{i\in I}$.
We shall always equip direct products
with the product topology.
\end{ex}
\begin{lem}\label{lempropsinit}
Let $X$ be equipped with the initial
topology with respect to a family
$(f_i)_{i\in I}$ of maps $f_i\colon X\to Y_i$
into topological spaces~$Y_i$.
Then the following holds:
\begin{description}[(D)]
\item[\rm (a)]
Let $(x_\alpha)$ be a net in~$X$ and $x\in X$.
Then $x_\alpha\to x$ in~$X$ if and only if
$f_i(x_\alpha)\to f_i(x)$
for each $i\in I$.
\item[\rm (b)]
A map $f\colon Z\to X$ from a topological
space~$Z$ to~$X$ is continuous
if and only if $f_i\circ f\colon Z\to Y_i$
is continuous for each $i\in I$.
\end{description}
\end{lem}
\begin{prf}
(a) If $x_\alpha\to x$, then $f_i(x_\alpha)\to f_i(x)$
by continuity of~$f_i$.
Conversely, assume that $f_i(x_\alpha)\to f_i(x)$
for each $i\in I$. Every neighborhood of~$x$
contains a neighborhood of the form
$V:=\bigcap_{i\in F}f_i^{-1}(U_i)$,
where $F\sub I$ is finite and $U_i\sub Y_i$
an open neighborhood of $f_i(x)$.
For each $i\in I$, there exists $\alpha_i$
such that $f_i(x_\alpha)\in U_i$ for all $\alpha\geq \alpha_i$.
Choose $\beta$ such that $\beta\geq \alpha_i$ for each
$i\in F$. Then $x_\alpha \in V$ for all $\alpha\geq \beta$.
Thus $x_\alpha\to x$.

(b) If $x_\alpha\to x$ in~$Z$,
then $f(x_\alpha)\to f(x)$
if and only if $f_i(f(x_\alpha))\to f_i(f(x))$
for all $i\in I$, by (a).
In view of Proposition~\ref{exemplify},
the assertion follows.
\end{prf}
Applying Lemma~\ref{lempropsinit}
to cartesian products,
we obtain:
\begin{lem}\label{ctsinprod}
Let $(X_i)_{i\in I}$
be a family of
topological spaces and
$X:=\prod_{i\in I}X_i$ be their cartesian
product,
with canonical projections
$\pr_i\colon X\to X_i$.
Then the following holds:
\begin{description}[(D)]
\item[\rm (a)]
A net $(x_\alpha)_{\alpha\in A}$ in $\prod_{i\in I} X_i$
converges if and only if,
for each $i\in I$,
the net $(\pr_i(x_\alpha))_{\alpha\in A}$
converges in~$X_i$.
If $\pr_i(x_\alpha)\to x_i$,
then $x_\alpha\to (x_i)_{i\in I}$.
\item[\rm (b)]
A map $f \colon Z\to \prod_{i\in I}X_i$
is continuous if and only if all of its components
$f_i:=\pr_i\circ f$ are continuous.\qed
\end{description}
\end{lem}
\begin{lem}[{\bf Transitivity of initial topologies}]\label{transinit}
Let $(f_i)_{i\in I}$ be a family of maps
$f_i\colon X\to Y_i$ into topological spaces
$Y_i$ whose topology is initial
with respect to a family $(g_j)_{j\in J_i}$
of mappings $g_j\colon Y_i\to Z_j$
into topological spaces~$Z_j$.
Set $K:=\{(i,j)\colon i\in I, j\in J_i\}$.
Then the initial topology on~$X$ with respect to the
family $(f_i)_{i\in I}$
and the initial topology with respect to
$(g_j\circ f_i)_{(i,j)\in K}$
coincide.
\end{lem}
\begin{prf}
A basis for both topologies
is given by the
sets $\bigcap_{(i,j)\in P}f_i^{-1}(g_j^{-1}(U_{i,j}))$,
where $P\sub K$ is a finite subset
and each $U_{i,j}$ an open subset of $Z_j$.
\end{prf}
\begin{rem}\label{reintinit}
Let $X$ be a set and $(f_i)_{i\in I}$
be a family of maps $f_i\colon X\to Y_i$ to topological
spaces~$Y_i$.
Lemma~\ref{transinit} entails
that the initial topology on $X$ with respect to
the family $(f_i)_{i\in I}$
coincides with the initial
topology with respect to the map
$f:=(f_i)_{i\in I}\colon X\to\prod_{i\in I } Y_i$
with components $\pr_i\circ f=f_i$.
\end{rem}
\begin{defn}\label{defembed}
A map $f\colon X\to Y$
between topological spaces is called
a (topological)
{\em embedding\/}
if~$f$ is injective and the topology
on $X$ is initial with respect to~$f$.
\end{defn}
\begin{rem}\label{embhomeoim}
By Example~\ref{exambije}
and Lemma~\ref{transinit},
an injective map $f\colon X\to Y$ between
topological spaces is an embedding if
and only if the corestriction
$f|^{f(X)}\colon X\to f(X)$
is a homeomorphism,
where $f(X)$ is equipped with the
topology induced by~$Y$.
\end{rem}
\section{Final topologies and quotient maps}\label{secfinaltop}
\begin{defn}\label{defnfinal}
Let $X$ be a set and $(f_i)_{i\in I}$
be a family of mappings $f_i\colon Y_i\to X$
from topological spaces~$Y_i$ to~$X$.
Then there is a finest topology~$\cO$ on~$X$
making each $f_i$ continuous;
it is called the
{\em final topology with respect to the family $(f_i)_{i\in I}$.}
\end{defn}
To see that~$\cO$ exists, call $U\sub X$ open
if and only if $f_i^{-1}(U)$ is open in~$Y_i$ for each $i\in I$.
Then $\cO$ is a topology on~$X$
with the desired properties.
\begin{lem}\label{fctnonfinal}
Suppose that $X$ is equipped with the final topology with
respect to a family $(f_i)_{i\in I}$
of mappings $Y_i\to X$, as in Definition~{\rm \ref{defnfinal}}.
Then a map $f\colon X\to Z$ from~$X$ to a topological space~$Z$
is continuous if and only if $f\circ f_i$ is continuous
for each $i\in I$.
\end{lem}
\begin{prf}
The map $f$ is continuous if and only if
$f^{-1}(U)$ is open in~$X$
for each open subset $U\sub Z$.
Here $f^{-1}(U)$ is open
if and only if
$f_i^{-1}(f^{-1}(U))=(f\circ f_i)^{-1}(U)$
is open for each~$i$.
The latter holds for each $U$ and each~$i$
if and only if all of the maps
$f\circ f_i$ are continuous.
\end{prf}
\begin{lem}[{\bf Transitivity of final
topologies}]\label{transfinal}
Let $(f_i)_{i\in I}$ be a family of mappings
$f_i\colon Y_i\to X$
from topological spaces~$Y_i$
to~$X$, such that the topology on $Y_i$
is final with respect to a family
$(g_j)_{j\in J_i}$ of mappings
$g_j\colon Z_j\to Y_i$ from
certain topological spaces $Z_j$ to~$Y_i$.
Define $K:=\{(i,j)\colon i\in I, j\in J_i\}$.
Then the final topology $\cS$ on~$X$ with respect to
the family $(f_i)_{i\in I}$ and
the final topology~$\cT$
with respect to $(f_i\circ g_j)_{(i,j)\in K}$ coincide.
\end{lem}
\begin{prf}
Let $U\sub X$.
Then $U\in \cS$ if and only if $f_i^{-1}(U)$
is open in~$Y_i$ for each $i\in I$,
which holds if and only if
$g_j^{-1}(f_i^{-1}(U))=(f_i\circ g_j)^{-1}(U)$
is open in~$Z_j$ for each $j\in J_i$.
But this condition is equivalent to $U\in \cT$.
\end{prf}
\begin{defn}\label{defquotmap}
A map
$q\colon X\to Q$ between
topological spaces
is called a {\em quotient map\/}
if $q$ is surjective
and a subset $U\sub Q$ is open in~$Q$
if and only if $q^{-1}(U)$ is
open in~$X$.
In this case, the topology on~$Q$ is called
the {\em quotient topology\/}
with respect to~$q$.
\end{defn}
Note that every quotient map is
continuous.
Clearly,
a map $q\colon X\to Q$ between topological
spaces is a quotient map
if and only if~$q$ is surjective
and $Q$ is equipped with the
final topology with respect to~$q$.
We therefore obtain as a special case
of Lemma~\ref{fctnonfinal}:
\begin{lem}\label{ctsonquot}
If $q\colon X\to Q$ is a quotient map
and $f\colon Q\to Y$ a map
to a topological space~$Y$,
then $f$ is continuous if and only if
$f\circ q\colon X\to Y$ is continuous.\qed
\end{lem}
\begin{lem}\label{opencquot}
Let $f\colon X\to Q$ be a surjective
continuous map between topological spaces.
If $f$ is closed or open,
then $f$ is a quotient map.
\end{lem}
\begin{prf}
Assume that $f$ is an open map,
i.e., $f(W)$ is open in~$Q$ for
each open subset $W\sub X$.
Let $U\sub Q$.
If $U$ is open,
then $f^{-1}(U)$ is open by continuity
of~$f$.
If $f^{-1}(U)\sub X$ is open,
then $U=f(f^{-1}(U))$
since~$f$ is surjective,
and this set is open
as $f$ is an open map.
Hence $f$ is a quotient map.
If $f$ is a closed
and $A\sub Q$,
then $f^{-1}(A)$ is closed in~$X$
if and only if $A$ is closed
in~$Q$,
by an analogous argument.
Hence $f$ is a quotient map also
in this case.
\end{prf}
\subsection*{Exercises for Section~\ref{secfinaltop}}

\begin{small}

\begin{exer}\label{topforbun}
Let $E$ be a topological space whose topology is final
with respect to a family $(f_j)_{j\in J}$ of injective mappings
$f_j\colon X_j\to E$ whose domains are topological spaces~$X_j$.
Assume that $f_i^{-1}(f_j(X_j))=f_i^{-1}(f_i(X_i)\cap f_j(X_j))$ is open in $X_i$
and
\[
f_j^{-1}\circ f_i\colon f_i^{-1}(f_i(X_i)\cap f_j(X_j))\to f_j^{-1}(f_i(X_i)\cap f_j(X_j))
\]
is a homeomorphism for all $i,j\in J$.
\begin{description}[(D)]
\item[(a)]
Show that $f_j(X_j)$ is open in~$E$
and $f_j$ is a homeomorphism onto its open image~$f_j(X_j)$,
for each $j\in J$.
\item[(b)]
If $X_j$ is Hausdorff for some $j\in J$, show that $\phi(X_j)$ is Hausdorff and
all points $x\not=y$ in~$f_j(X_j)$
have disjoint open neighborhoods in~$X$.
\item[(c)]
Assume that $X_j$ is Hausdorff for all $j\in J$
and assume that there exists a continuous map $f\colon E\to M$ to a Hausdorff space~$M$
such that $f(x)=f(y)$ for $x,y\in E$ imples that $x,y\in f_j(X_j)$
for some $j\in J$. Show that $E$ is Hausdorff.
\end{description}
[Such observations are useful for some constructions of vector bundles].
\end{exer}
\end{small}
\section{Compactness and locally compact spaces}
We collect facts concerning compact topological spaces and locally compact
spaces.

Recall that a topological space $K$ is called \emph{compact}
if it is Hausdorff and every open cover of~$K$ has a finite subcover.
Thus, if $(V_j)_{j\in J}$ is a family of open subsets $V_j$
of $K$ such that
\[
K=\bigcup_{j\in J} V_j,
\]
then there exists a finite subset $F\sub J$ such that $K=\bigcup_{j\in F} V_j$.
A subset $K\sub X$ of a topological space~$X$ is called compact
if it is compact in the induced topology.
This is the case if and only if $K$ is Hausdorff and
for each family $(W_j)_{j\in J}$ of open subsets $W_j$ of $X$
with
\[
K\sub \bigcup_{j\in J}W_j,
\]
there is a finite subset $F\sub J$ such that $K\sub\bigcup_{j\in F} W_j$.

Also other properties of topological spaces play a role.
\begin{defn}
Let $X$ be a topological space.
\begin{description}[(D)]
\item[(a)]
$X$ is called \emph{locally compact}
if~$X$ is Hausdorff and each neighborhood~$U$
of a point $x\in X$ contains a compact $x$-neighborhood.
\item[(b)]
$X$ is called \emph{regular}
if $X$ is Hausdorff and each neighborhood $U$
of a point $x\in X$ contains an $x$-neighborhood which is closed in~$X$.
\item[(c)]
$X$ is called \emph{normal}
if $X$ is Hausdorff and, for all closed subsets $A$ and $B$
of~$X$ such that $A\cap B=\emptyset$,
there exist open subsets $U$ and $V$ of~$X$
such that $A\sub U$, $B\sub V$ and $U\cap V=\emptyset$.
\item[(d)]
If $X=\bigcup_{n\in \N}K_n$ for a sequence $(K_n)_{n\in\N}$
of compact subsets $K_n$ of~$X$,
then $X$ is called \emph{$\sigma$-compact}.
\end{description}
\end{defn}
\begin{rem}
Note that every locally compact topological space
is regular, as compact subsets of a Hausdorff space~$X$ are
closed in~$X$.
\end{rem}
\begin{lem}[{\bf Wallace Lemma}]\label{Wallla}
Let $X_1$ and $X_2$ be topological spaces,
$K_1\sub X_1$ and $K_2\sub X_2$ be compact subsets and $U\sub X_1\times X_2$
be an open subset such that
\[
K_1\times K_2\sub U.
\]
Then there exist open subsets $P \sub X_1$ and $Q \sub X_2$ such that
\[
K_1\times K_2\sub P\times Q\sub U.
\]
\end{lem}
\begin{prf}
For all $x\in K_1$ and $y\in K_2$, there are open neighborhoods $P_{x,y}\sub X_1$ and $Q_{x,y}\sub X_2$
of $x$ and $y$, respectively, such that
\[
P_{x,y}\times Q_{x,y}\sub U.
\]
Let $x\in X_1$. Then $(Q_{x,y})_{y\in K_2}$ is a family of open subsets of $X_2$ such that
\[
K_2\sub \bigcup_{y\in K_2}Q_{x,y}.
\]
By compactness of $K_2$, there is a finite subset $F_x\sub K_2$ such that
\[
K_2\sub\bigcup_{y\in F_x}Q_{x,y}=:Q_x.
\]
Then $P_x:=\bigcap_{y\in F_x}P_{x,y}$ is an open neighborhood of~$x$ in $X_1$
such that
\[
P_x\times Q_x\sub U.
\]
Now $(P_x)_{x\in K_1}$ is a family of open subsets of $X_1$ such that
$K_1\sub\bigcup_{x\in K_1}P_x$.
By compactness of $K_1$, there is a finite subset $F\sub K_1$ such that
\[
K_1\sub \bigcup_{x\in F} P_x=: P.
\]
Now $Q\!:=\!\bigcap_{x\in F} Q_x$ is open in $X_2$ and $K_2\!\sub \!Q$.
By construction, $P\!\times\! Q \sub  U$.
\end{prf}
%
%
%
%
\begin{lem}\label{lemcpnorm}
Every compact topological space is normal.
\end{lem}
\begin{prf}
Let $K$ be a compact topological space;
let $A$ and $B$ be disjoint closed subsets of~$K$.
As $K$ is Hausdorff,
$U:=\{(x,y)\in K\times K\colon x\not=y\}$ is an open subset
of $K\times K$. That $A$ and $B$ are disjoint is equivalent to
$A\times B\sub U$. Using the Wallace Lemma,
we obtain open subsets $P,Q\sub K$ such that $A\sub P$, $B\sub Q$
and $P\times Q\sub U$, whence $P\cap Q=\emptyset$.
\end{prf}
\begin{rem}\label{betternormal}
In the preceding lemma, one can always achieve that even the closures $\wb{P}$
and $\wb{Q}$ are disjoint.
If fact, if $P$ and $Q_1$ are disjoint open subsets of $K$ with $A\sub P$ and $B\sub Q_1$,
then $P$ (and hence also $\wb{P}$)
is contained in the closed set $K\setminus Q_1$. Hence $\wb{P}$
and $B$ are disjoint closed sets, and using the lemma again we find
disjoint open subsets $P_1,Q\sub K$
such that $\wb{P}\sub P_1$ and $B\sub Q$.
Then $\wb{Q}\sub K\setminus P_1$,
whence $\wb{Q}\cap \wb{P}=\emptyset$.
\end{rem}
\begin{rem}\label{trivi-lcp}
If $X$ is a locally compact topological space, $K\sub X$ a compact subset and
$U\sub X$ an open subset such that $K\sub U$,
then there exists a compact subset $L\sub U$ which contains~$K$
in its interior.\\[2.3mm]
In fact, we can cover $K$
with finitely many open subsets of~$U$
which are relatively compact in~$U$, and define $L$
as the union of their closures.
\end{rem}
\begin{defn}\label{defn-exhaust}
Let $X$ be a locally compact topological space.
A sequence $(K_n)_{n\in \N}$
of compact subsets $K_n$ of~$X$ is called a \emph{compact exhaustion of~$X$}
if $X=\bigcup_{n\in \N}K_n=X$
and $K_n\sub K_{n+1}^0$ for each $n\in\N$.
\end{defn}
\begin{lem}\label{exhaustions-exist}
Every $\sigma$-compact, locally compact topological space~$X$
admits a compact exhaustion.
\end{lem}
\begin{prf}
There exists a sequence $(K_n)_{n\in \N}$
of compact subsets of~$X$ such that $X=\bigcup_{n\in \N}K_n$.
Let $L_1:=K_1$. If compact subsets $L_1,\ldots, L_n$
of~$X$ have already been constructed such that $K_j\sub L_j$
for all $j\in \{1,\ldots, n\}$ and \break $L_j\sub L_{j+1}^0$ if $j<n$,
then $L_n\cup K_{n+1}$ is a compact subset of~$X$.
By Remark~\ref{trivi-lcp},
there exist a compact subset $L_{n+1}$ of~$X$ such that $L_n\cup K_{n+1}\sub L_{n+1}^0$.
Then $(L_n)_{n\in\N}$ is a compact exhaustion.
\end{prf}
\section{Paracompactness}
We recall the definition of a paracompact topological space
and record some consequences which are useful
in infinite-dimensional calculus.
\begin{numba}\label{subordinate}
If $(U_i)_{i\in I}$
is an open cover of a topological space~$X$,
let us say that a family $(Y_j)_{j\in J}$
of subsets of~$X$ is \emph{subordinate}
to $(U_i)_{i\in I}$ if, for each $j\in J$,
there exists $i\in I$ such that $Y_j\sub U_i$.
A family $(Y_j)_{j\in J}$ of subsets of~$X$
is called \emph{locally finite}
if each $x\in X$ has a neighborhood $W\sub X$ such that
$\{j\in J\colon Y_j\cap W\not=\emptyset\}$ is finite.
\end{numba}
\begin{defn}\label{defn-paracompact}
A topological space $X$ is called \emph{paracompact}
if it is Hausdorff and has the following property:
For
each open cover $(U_i)_{i \in I}$,
there exists a locally finite open cover $(V_j)_{j\in J}$
of~$M$ which is subordinate to $(U_i)_{i\in I}$.
\end{defn}
\begin{lem}\label{ops-locfin}
Let $X$ be a topological space and $(A_j)_{j\in J}$
be a locally finite famly of subsets of~$X$.
Then the following holds:
\begin{description}[(D)]
\item[\rm(a)]
Also the family $(\wb{A_j})_{j\in J}$
of closures is locally finite.
\item[\rm(b)]
For every compact subset $K\sub X$, the set
$\{j\in J\colon K\cap A_j\not=\emptyset\}$ is finite.
\item[\rm(c)]
If $A_j$ is closed in $X$ for each $j\in J$,
then $A:=\bigcup_{j\in J}A_j$ is closed in~$X$.
\end{description}
\end{lem}
\begin{prf}
(a) Each $x\in X$ has an open neighborhood $W=W(x)$ in $X$ such that
\[
J_x:=\{j\in J\colon A_j\cap W\not=\emptyset\}
\]
is finite.
For $j\in J\setminus J_0$,
we have $A_j\sub X\setminus W$. The latter set being closed,
$\wb{A_j}\sub X\setminus W$
follows, whence
$\wb{A_j}\cap W=\emptyset$.
Thus $\{j\in J\colon \wb{A_j}\cap W\not=\emptyset\}
=J_x$ is finite.\\[1mm]
(b) The subset $K$ being compact, we find
$x_1,\ldots, x_n\in K$ such that
\[
K\sub W(x_1)\cup\cdots\cup W(x_n).
\]
For each
$j\in J$, the condition $A_j\cap K\not=\emptyset$
implies $A_j\cap W(x_k)\not=\emptyset$
for some $k\in\{1\ldots, n\}$, whence
$j$ is contained in the finite set $J_{x_1}\cup \cdots\cup J_{x_n}$.\\[1mm]
(c) For each
$x\in X\setminus A$,
there exists an open $x$-neighborhood $W\sub X$ such that
$J_0:=\{j\in J\colon W\cap A_j\not=\emptyset\}$
is finite. Then
$W\setminus A=W\setminus \bigcup_{j\in J_0}A_j$
is an open $x$-neighborhood in~$X$, whence
$x\in (X\setminus A)^0$. Thus $X\setminus A$
is open.
\end{prf}
Every paracompact topological space is normal. We only need:
%
\begin{prop}\label{paracomp-reg}
Every paracompact topological space is regular.
\end{prop}
\begin{prf}
Let $x\in X$ and $U$ be an open neighborhood of~$x$ in~$X$.
For each $y\in X\setminus \{x\}$,
there exist disjoint open neighborhoods $P_y$ of~$x$
and $Q_y$ of~$y$ in~$X$.
Then $X\setminus\{x\}=\bigcup_{y\not=x}Q_y$
is open and thus $\{x\}$ is closed in~$X$.
Let $Q_x:=U$. 
As $X$ is assumed paracompact,
we can find a locally finite open cover $(V_j)_{j\in J}$
of~$X$ subordinate to $(Q_y)_{y\in X}$.
Thus, for each $j\in J$, there exists $x(j)\in X$ such that $V_j\sub Q_{x(j)}$.
For each
$j$ in $J_1:=\{j\in J\colon x(j)\not=x\}$,
we have $V_j\sub Q_{x(j)}$, whence $V_j\sub X\setminus P_{x(j)}$
and thus also $\overline{V_j}\sub X\setminus P_{x(j)}$,
entailing that $x\not\in \overline{V_j}$.
Since $(\overline{V_j})_{j\in J_1}$ is locally finite
by Lemma~\ref{ops-locfin}(a),
$A:=\bigcup_{j\in J_1}\overline{V_j}$
is closed in~$X$ by Lemma~\ref{ops-locfin}(b).
Since $x\not\in A$, we see that
$W:=X\setminus A$ is an open $x$-neighborhood in~$X$.
Now $X=U\cup \bigcup_{j\in J_1}V_j$ entails that
\[
X\setminus U\sub \bigcup_{j\in J_1}V_j\sub A,
\]
whence $U\supseteq X\setminus \bigcup_{j\in J_1}V_j\supseteq
X\setminus A=W$. The second set being closed in~$X$,
we deduce that $\overline{W}\sub X\setminus
\bigcup_{j\in J_1}V_j\sub U$.
\end{prf}
Let us turn to locally compact spaces.
%
\begin{lem}\label{lem-rings}
Let $X$ be a locally compact topological space
and $A\sub X$ be a $\sigma$-compact closed subset.
Then there are compact subsets $L_j\sub A$
and open subsets $O_j\sub X$ with $L_j\sub O_j$ for $j\in\N$
such that $\bigcup_{j\in \N}L_j=A$ holds,
\begin{eqnarray}
O_i\cap O_j\not=\emptyset &\Rightarrow & |i-j|\leq 1, \quad \mbox{and}\label{henceleq2}\\
L_i\cap L_j\not=\emptyset &\Rightarrow  & |i-j|\leq 1.\label{henceleq1}
\end{eqnarray}
\end{lem}
\begin{prf}
Let $(C_j)_{j\in\N}$
be a sequence of compact subsets of~$A$ with union~$A$.
Then there exists a sequence $(K_j)_{j\in\N}$
of compact subsets of~$X$ such that, setting $K_{-1}:=K_0:=\emptyset$,
\[
K_{j-1}\sub K_j^0\quad \mbox{and}\quad  C_j\sub K_j
\]
holds for all $j\in\N$, using
the interior $K_j^0$ of $K_j$
as a subset of~$X$.
In fact, if $j\in\N$ and $K_i$ has been constructed for $i\in \{0,\ldots,j-1\}$,
then the compact set $K_{j-1}\cup C_j$
is contained in $K_j^0$ for some compact subset $K_j$
of~$X$, by Remark~\ref{trivi-lcp}.
We now define
\[
L_j:= A\cap (K_j\setminus K_{j-1}^0)
\mbox{ and }U_j:=K_{j+1}^0\setminus K_{j-2} 
\]
for $j\in\N$.
Then $L_j$ is compact, $U_j$ an open subset of~$X$ and $L_j\sub U_j$,
since $K_j\sub K_{j+1}^0$ and $K_{j-2}\sub K_{j-1}^0$.
If $1\leq i\leq j-2$, then $L_i\sub K_i\sub K_{j-2}\sub K_{j-1}^0$
and thus $L_i\cap L_j=\emptyset$,
whence (\ref{henceleq1}) holds.
Similarly, $1\leq i\leq j-3$ implies that $U_i\sub K_{i+1}\sub K_{j-2}$
and thus $U_i\cap U_j=\emptyset$.
As a consequence,
\[
(\forall i,j\in\N)\;\,
U_i\cap U_j\not=\emptyset\Rightarrow |i-j|\leq 2,
\]
entailing that $(U_j)_{j\in \N}$
(hence also $(L_j)_{j\in\N}$)
is locally finite in the open subset $U:=\bigcup_{j\in\N}U_j$ of~$X$.
Consequently, $A_J:=\bigcup_{j\in \N\setminus J}L_j$ is closed in~$U$
for each subset $J\sub \N_0$.
Let $O_{-1}:=O_0:=\emptyset$.
Recursively,
using Remark~\ref{trivi-lcp}
we find
relatively compact, open subsets $O_j$ of~$U$
for $j\in \N$ such that $L_j\sub O_j$ 
and\vspace{-1mm}
\begin{equation}\label{nice-in}
\wb{O_j}\sub (U\setminus A_{\{j-1,j,j+1\}})\setminus\bigcup_{i=-1}^{j-2}\wb{O_i},\vspace{-1mm}
\end{equation}
as the right-hand side is open and contains~$L_j$.
By (\ref{nice-in}), we have (\ref{henceleq2}).
\end{prf}
\begin{lem}\label{pre-lcp-para}
Every $\sigma$-compact locally compact
space~$X$
is paracompact.
\end{lem}
\begin{prf}
%
Let $(L_n)_{n\in \N}$ and $(O_n)_{n\in\N}$ be as in
Lemma~\ref{lem-rings}.applied with $A:=X$.
Let $(U_i)_{i\in I}$ be an open cover of~$X$.
For each $x\in X$, there exists $i(x)\in I$
such that $x\in U_{i(x)}$.
For each $n\in \N$, the open sets
\[
V_{(n,x)}:= O_n\cap U_{i(x)}
\]
cover $L_n$ for $x\in L_n$, whence $L_n\sub\bigcup_{x\in \Phi_n}V_{(n,x)}$
for some finite subset \break $\Phi_n\sub L_n$.
We define $J:=\bigcup_{n\in\N}\{n\}\times\Phi_n$
and obtain the open cover $(V_j)_{j\in J}$ of~$X$
which is subordinate to $(U_i)_{i\in I}$
as $V_{(n,x)}\sub U_{i(x)}$
and locally finite as each $x\in X$ is contained in $O_n$
for some $n\in\N$ and $O_n\cap V_{(m,y)}$
implies $m\in \{n-1,n,n+1\}$ with $y$ in the finite set $\Phi_m$.
\end{prf}
We need topological sums as a tool.
\begin{numba}
Let
$(X_j,\cO_j)_{j\in J}$
be a family of topological spaces such that \break
$X_i\cap X_j=\emptyset$
for all $i\not=j$ in~$J$.
The union
$X:=\bigcup_{j\in J}X_j$
is called the \emph{topological sum}
of the topological spaces $X_j$,
if it is endowed with the final topology
with respect to the inclusion maps
$\lambda_j\colon X_j\to X$, $x\mto x$.
\end{numba}
By Exercise~\ref{topforbun}, each $X_j$ is open in~$X$
and a subset
$U\sub X_j$
is open in $(X_j,\cO_j)$ if and only if it is open in
$(X,\cO)$.
If each $X_j$ is Hausdorff, then also the topological sum~$X$.
In fact, distinct points within $X_j$ can be separated by open sets in $X_j$.
If $x\in X_i$ and $y \in X_j$
with $i\not=j$, then
$X_i$ and $X_j$ are disjoint neighborhoods
of $x$ and $y$, respectively.
\begin{prop}\label{lcp-parac}
A locally compact topological space is paracompact
if and only if it is a topological sum of $\sigma$-compact
locally compact
spaces.
\end{prop}
\begin{prf}
If $X$ is paracompact, then each $x\in X$ has a relatively compact,
open neighborhood $U(x)$ in~$X$.
By paracompactness, there exists a locally finite open cover $(V_j)_{j\in J}$
of~$X$ subordinate to $(U(x))_{x\in X}$.
Then $V_j$ is relatively compact for each $j\in J$,
as $V_j\sub U(x_j)$ for some $x_j\in X$.
For $x,y\in X$, we write $x\sim y$
if there exist $n\in\N$ and $x_0,x_1,\ldots, x_n\in X$
such that
$x_0=x$, $x_n=y$ and,
for all $k\in\{0,\ldots, n-1\}$, there exists $j\in J$ with
\[
\{x_k,x_{k+1}\}\sub \wb{V_j}.
\]
We readily check that $\sim$ is an equivalence relation
on~$X$.
%
%
We now show that the equivalence classes $[x]$
are open in~$X$ (whence they are locally compact)
and $\sigma$-compact.
If $y\in [x]$,
there exists $j\in J$
such that $y\in V_j$.
Now $\{y,z\}\sub V_j$ for every
$z\in V_j$, whence $z\sim y$ and thus $V_j\sub [x]$.
Hence $[x]$ is a neighborhood of $y$ and hence open,
as $y\in [x]$ was arbitrary.
For $x\in X$ let $C_n(x)$
be the set of all $y\in [x]$, for which we can use $x=x_0,\ldots, x_n=y$
as above with the given~$n$.
The family $(\wb{V_j})_{j\in J}$
being locally finite by Lemma~\ref{ops-locfin}(a),
$I_1:=\{j\in J\colon x\in\wb{V_j}\}$
is finite and hence
\[
C_1(x)=\bigcup_{j\in I_1}\wb{V_j}
\]
compact.
If we know that
$C_n(x)=\bigcup_{j\in I_n}\wb{V_j}$
for a finite subset
$I_n\sub J$,
then
$J_j:=\{i\in J\colon \wb{V_i}\cap C_n(x) \not=\emptyset\}$
is a finite set for each $j\in I_n$, by (a) and (b) in Lemma~\ref{ops-locfin}.
Setting
\[
I_{n+1}:=\bigcup_{j\in I_n}J_j,
\]
we deduce that $C_{n+1}(x)=\bigcup_{j\in I_{n+1}}\wb{V_j}$
is compact. Thus $[x]=\bigcup_{n\in\N}C_n(x)$ is $\sigma$-compact.\\[1mm]
If, conversely, $X$ is a topological sum
of a family $(X_a)_{a\in A}$ of $\sigma$-compact
locally compact topological spaces,
then each $X_a$ is paracompact, by Lemma~\ref{pre-lcp-para}.
If $(U_i)_{i\in I}$
is an open cover of~$X$, then
there exists a locally finite open cover $(V_{(a,b)})_{b\in B_a}$
of $X_a$ subordinate to the open cover $(U_i\cap X_a)_{i\in I}$
of~$X_a$. Let $J:=\bigcup_{a\in A}\{a\}\times B_a$.
Then $(V_j)_{j\in J}$ is a locally finite open cover of $X=\bigcup_{a\in A}X_a$
which is subordinate to $(U_i)_{i\in I}$.
\end{prf}
We frequently know that a map is
a local homeomorphism. The following theorem
will ensure it is a homeomorphism
on a large open set.
\begin{thm}\label{godement-globalize}
Let $f\colon Y\to Z$ be a map between
topological spaces and $X\sub Y$ be a subset.
We assume:
\begin{description}[(D)]
\item[\rm(a)]
$f$ is a local homeomorphism;
\item[\rm(b)]
$f(X)$ is closed in~$Z$
and $f|_X\colon X\to f(X)$
is a homeomorphism;
\item[\rm(c)]
$Z$ is paracompact.
\end{description}
Then there exists an open subset
$U\sub Y$ with $X\sub U$ such that $f(U)$
is open in~$Z$ and $f|_U\colon U\to f(U)$
is a homeomorphism.
\end{thm}
The proof of Theorem~\ref{godement-globalize}
uses the following lemma,
which is of independent interest (see Proposition~\ref{global-analyt-extension}).
\begin{lem}\label{pre-godement}
Let $Y$ and $Z$ be topological spaces,
$A\sub Z$ be a closed subset and
$(V_i)_{i\in I}$ be a family of open subsets of~$Z$
such that $A\sub\bigcup_{i\in I} V_i$.
For each $i\in I$, let $g_i\colon V_i\to Y$
be a continuous function. We assume:
\begin{description}[(D)]
\item[\rm(a)]
For all $i,j\in I$ and $a\in V_i\cap V_j\cap A$,
there exists an open $a$-neighborhood
$O=O(a,i,j)\sub V_i\cap V_j$
such that $g_i|_O=g_j|_O$;
\item[\rm(b)]
$Z$ is paracompact.
\end{description}
Then there exists an open subset
$V\sub \bigcup_{i\in I}V_i$
with $A\sub V$ and a continuous function
$g\colon V\to Y$
such that each
$x\in V$
has an open neighborhood $W\sub V$
such that $W\sub V_i$ for some $i\in I$ and
$g|_W=g_i|_W$.
\end{lem}
\begin{prf}
For each $z\in A$, there exists
$\iota(z)\in I$
such that $z\in V_{\iota(z)}$.
Since $Z$ is paracompact and thus regular (see Proposition~\ref{paracomp-reg}),
we can find a closed neighborhood~$Q_z$ of~$z$ in~$Z$
such that $Q_z\sub V_{\iota(z)}$.
For $z\in Z\setminus A$, let $Q_z:=Z\setminus A$.
By paracompactness of~$Z$,
there exists a locally finite open cover
$(W_j)_{j\in J}$
of~$Z$ which is subordinate to
$((Q_z)^0)_{z\in Z}$.
Thus, for each
$j\in J$, there exists
$z(j)\in Z$ such that $W_j\sub (Q_{z(j)})^0$.
Let $J_0:=\{j\in J\colon z(j)\in A\}$.
Then $(W_j)_{j\in J_0}$ is a locally finite family of open subsets
of~$Z$ and $A\sub\bigcup_{j\in J_0}W_j$.
For each $j\in J_0$, we have
\[
\overline{W_j}\sub Q_{z(j)}\sub V_{i(j)}
\quad \mbox{ with } \quad i(j):=\iota(z(j)).\]
By Lemma~\ref{ops-locfin}(a),
each $a\in A$ has an open neighborhood~$U_a$ in~$Z$
such that
\[
\Phi(a):=\{j\in J_0\colon U_a\cap \overline{W_j}\not=\emptyset\}
\]
is finite.
If $\Psi:=\{j\in\Phi(a)\colon a\not\in \overline{W_j}\}
\not=\emptyset$,
we can replace $U_a$ with its open subset
\[
U_a\setminus \bigcup_{j\in\Psi} \overline{W_j};
\]
we may therefore assume that
$a\in\wb{W}_j$ for all $j\in \Phi(a)$.
For all $j,k\in \Phi(a)$, we have
$a\in V_{i(j)}\cap V_{i(k)}\cap A$,
and thus $a\in O(a, i(j),i(k))$;
after replacing $U_a$ with
\[
U_a\cap\bigcap_{j,k\in\Phi(a)}O(a,i(j),i(k)),
\]
we may assume that $U_a\sub O(a,i(j),i(k))$
for all $j,k\in\Phi(a)$.
Shrinking the set further, we may assume that
$U_a\sub W_{j(a)}$ for some $j(a)\in J_0$;
then $j(a)\in \Phi(a)$.
Now
\[
V:=\bigcup_{a\in A} U_a
\]
is an open subset of~$Z$
with $A\sub Z$
and $g\colon V\to Y$, $z\mto g_{i(j(a))}(z)$
for $z\in U_a$ is well defined.
In fact, if $z\in U_a\cap U_b$, then $z\in U_b\sub W_{j(b)}$,
whence $j(b)\in \Phi(a)$ and thus
$U_a\sub O(a,i(j(a)),i(j(b))=:O$;
hence
\[
g_{i(j(a))}(z)=g_{i(j(a))}|_{O}(z)=g_{i(j(b))}|_O(z)=g_{i(j(b))}(z).
\]
Notably,
$g|_{U_a}=g_{i(j(a)}|_{U_a}$ for each $a\in A$,
whence $g|_{U_a}$ is continuous. Hence $g$ is continuous.
\end{prf}
\noindent
\emph{Proof of Theorem}~\ref{godement-globalize}.
For $z\in f(X)$, there is a unique $x\in X$
with $f(x)=z$.
By hypothesis,
$x$ has an open neighborhood
$U_z$ in~$Y$
such that
$V_z:=f(U_z)$ is open in~$Z$
and $f|_{U_z}\colon U_z\to V_z$
is a homeomorphism.
We let
\[
g_z:=(f|_{U_z})^{-1}\colon V_z\to U_z.
\]
As $f(X\cap U_z)$ is relatively open in $f(X)$,
there exists an open subset
$P$ of~$Z$
such that $f(X\cap U_z)=f(X)\cap P$.
After replacing $U_z$ with $U_z\cap f^{-1}(P)$,
we may assume that
\[
f(X)\cap V_z=f(X\cap U_z);
\]
then $g_z|_{f(X)\cap V_z}=(f|_{X\cap U_z})^{-1}=(f|_X)^{-1}|_{f(X)\cap V_z}$.
If $a\in f(X)\cap V_z\cap V_w$
with $z,w\in f(X)$,
then $(f|_X)^{-1}(a)\in X\cap U_z\cap U_w$.
Now
\[ O=O(a,z,w):=f(U_z\cap U_w) =f|_{U_z}(U_z\cap U_w) \] 
is an open $a$-neighborhood in $V_z\cap V_w$
and $g_z|_O=(f|_{U_z\cap U_w})^{-1}=g_w|_O.$ 
By Lemma~\ref{pre-godement},
there exists an open subset
$V\sub \bigcup_{z\in f(X)}V_z$
with $f(X)\sub V$
and a continuous function
$g\colon V\to Y$ such that each
$x\in V$
has an open neighborhood
$W\sub V$ such that $W\sub V_z$ for some
$z\in f(X)$
and $g|_W=g_z|_W$.
Hence $g$ is a local homeomorphism
and thus $U:=g(V)$ is open in~$Y$.
Moreover,
$(f\circ g)|_W=(f\circ g_z)|_W=\id_W$,
whence $f|_U\circ g=\id_V$;
thus the maps $f|_U\colon U\to V$ and $g\colon V\to U$
are mutually inverse homeomorphisms.
If $x_0\in X$, then $x:=f(x_0)\in V$
und for $W$ and $z$ as above we have
$g(x)=g_z(x)=(f|_X)^{-1}(x)=x_0$,
whence $X\sub U$. \qed
\section{The compact-open topology}
\label{appcotop}
In this section, we give a self-contained introduction to the compact-open topology
on function spaces.
Most of the results are classical, and no originality is claimed
(cf.\ \cite[Chapter~X, \S1--\S3]{Bou66}, for example).
However, the sources treating
the topic do so with a different thrust.
Here, we compile precisely those results which are the foundation
for the study of non-linear mappings between function spaces,
and their differentiability properties.

We recall: If $X$ and $Y$ are Hausdorff topological spaces,
then the \emph{compact-open topology} on $C(X,Y)$
is the topology given by the subbasis of open subsets
\begin{equation}\label{thesubba}
\lfloor K,U\rfloor :=\{\gamma\in C(X,Y)\colon \gamma(K)\sub U\},
\end{equation}
for $K$ ranging through the set $\cK(X)$ of all compact subsets
of $X$ and $U$ through the set of open subsets of~$Y$.
In other words, finite intersections
of sets as in (\ref{thesubba}) form a basis for the compact-open topology on $C(X,Y)$.
We always endow $C(X,Y)$ with the compact-open topology (unless the contrary is stated).
\begin{rem}
The compact-open topology on $C(X,Y)$ makes the point evaluation
\[
\ev_x\colon C(X,Y)\to Y,\quad \gamma\mto \gamma(x)
\]
continuous, for each $x\in X$.
(If $U\sub Y$ is open, then $\ev_x^{-1}(U)=\lfloor \{x\}, U\rfloor$
is open in $C(X,Y)$).
As the maps $\ev_x$ separate points on $C(X,Y)$ for $x\in X$,
it follows that
the compact-open topology on $C(X,Y)$
is Hausdorff.
\end{rem}
\begin{lem}\label{subbaca}
Let $X$ and $Y$ be Hausdorff topological spaces and
$\cS$ be a subbasis of open subsets of~$Y$.
Then
\[
\cV:=\{\lfloor K,U\rfloor \colon K\in \cK(X), U\in \cS\}
\]
is a subbasis for the compact-open topology on $C(X,Y)$.
\end{lem}
\begin{prf}
Each of the sets $V\in \cV$ is open in the compact-open topology.
To see that $\cV$ is a subbasis for the latter,
it suffices to show that
$\lfloor K,U\rfloor$
is open in the topology generated by $\cV$,
for each $K\in \cK(X)$ and open subset $U\sub Y$.
To check this, let $\gamma\in \lfloor K,U\rfloor$.
For each $x\in K$, there is a finite subset $F_x\sub \cS$ such that
\[
\gamma(x)\in \bigcap_{S\in F_x}S
\sub U.
\]
Let $S_x:=\bigcap_{S\in F_x}S$.
Since~$K$ is compact and hence locally compact,
each $x\in K$ has a compact neighborhood $K_x\sub \gamma^{-1}(S_x)$ in~$K$.
By compactness, there is a finite subset $\Phi\sub K$
such that $K=\bigcup_{x\in \Phi}K_x^0$ (where $K_x^0$
denotes the interior of $K_x$ relative $K$).
Then $V:=\bigcap_{x\in \Phi}\lfloor K_x,S_x\rfloor=
\bigcap_{x\in \Phi}\bigcap_{S\in F_x}\lfloor K_x,S\rfloor$
is open in the topology generated by~$\cV$,
and $V\sub \lfloor K,U\rfloor$.
In fact, let $\gamma\in V$.
For each $y\in K$, there is $x\in \Phi$ such that $y\in K_x$.
Hence $\gamma(y)\in \gamma(K_x)\sub S_x\sub U$
and thus $\gamma\in \lfloor K,U\rfloor$.
\end{prf}
Mappings between function spaces
of the form $C(X,f)$ (so-called superposition operators)
are used frequently.
\begin{lem}\label{covsuppo}
Let $X$, $Y_1$ and $Y_2$ be Hausdorff topological spaces.
If a map $f\colon Y_1\to Y_2$ is continuous, then also
the following map is continuous:
\[
C(X,f)\colon C(X,Y_1)\to C(X,Y_2), \quad \gamma\mto f\circ \gamma\, .
\]
\end{lem}
\begin{prf}
The map $C(X,f)$ will be continuous if pre-images of subbasic open subsets are open.
To this end, let $K\in \cK(X)$ and $U\sub Y_2$ be open.
Then $C(X,f)^{-1}(\lfloor K,U\rfloor)=\lfloor K,f^{-1}(U)\rfloor$ is open in $C(X,Y_1)$
indeed.
\end{prf}
\begin{lem}\label{inipush}
Let $X$ and $Y$ be Hausdorff topological spaces.
Assume that the topology on $Y$ is initial with respect to
a family $(f_j)_{j\in J}$ of maps \break $f_j\colon Y\to Y_j$
to Hausdorff topological spaces~$Y_j\!$.\hspace*{-.5mm}
Then the compact-open topology on $C(X,Y)$ is initial with respect to the
family $(C(X,f_j))_{j\in J}$ of the mappings 
$C(X,f_j)\colon C(X,Y)\to C(X,Y_j)$.
\end{lem}
\begin{prf}
Let $\cS$ be the set of all subsets of $Y$ of the form
$f_j^{-1}(W)$, with $j\in J$ and $W$ an open subset of $Y_j$.
By hypothesis, $\cS$ is a subbasis for the topology of~$Y$.
Hence, by Lemma~\ref{subbaca},
the sets $\lfloor K, f_j^{-1}(W)\rfloor$
form a subbasis for
the compact-open topology on $C(X,Y)$,
for $j,W$ as before and $K\in \cK(X)$.
But $\lfloor K, f_j^{-1}(W)\rfloor=C(X,f_j)^{-1}(\lfloor K,W\rfloor)$
(since $\gamma(K)\sub f_j^{-1}(W)$ is equivalent to $f_j(\gamma(K))\sub W$, 
and this is equivalent to $(f_j\circ \gamma)(K)\sub W$),
and these sets form a subbasis for the initial topology on $C(X,Y)$ with respect to the family $(C(X,f_j))_{j\in J}$.
\end{prf}
We mention three consequences of Lemma~\ref{inipush}.
\begin{lem}\label{ctsemb}
If $X$ is a Hausdorff topological space and $f\colon Y_1\to Y_2$ is a\linebreak
topological embedding of Hausdorff topological spaces,
then also the mapping\linebreak
$C(X,f)\colon C(X,Y_1)\to C(X,Y_2)$
is a topological embedding.
\end{lem}
\begin{prf}
Since $f$ is injective, also $C(X,f)$ is injective.
By hypothesis, the topology on $Y_1$ is initial with respect to $f$.
Hence the topology on $C(X,Y_1)$ is initial with respect to $C(X,f)$
(by Lemma~\ref{inipush}).
As a consequence, the injective map $C(X,f)$ is a topological embedding.
\end{prf}
\begin{rem}\label{reminduco}
In particular, if $X$ and $Y_2$ are Hausdorff topological spaces and
$Y_1\sub Y_2$ is a subset, endowed with the induced topology,
then the topology induced by $C(X,Y_2)$ on $C(X,Y_1)$ coincides with the
compact-open topology on $C(X,Y_1)$.
\end{rem}
Next, we deduce that $C(X,\prod_{j\in J} Y_j)=\prod_{j\in J}C(X,Y_j)$.
\begin{lem}\label{cotopprod}
Let $X$ be a Hausdorff topological space and $(Y_j)_{j\in J}$
be a family of Hausdorff topological spaces, with the cartesian product
$Y:=\prod_{j\in J}Y_j$ $($endowed with the product topology$)$
and the projections
$\pr_j\colon Y\to Y_j$. Then the natural map
\[
\Phi:=(C(X,\pr_j))_{j\in J}\colon C(X,Y)\to \prod_{j\in J}C(X,Y_j)
\]
is a homeomorphism.
\end{lem}
\begin{prf}
The map $\Phi$ is a bijection because
a map $f\colon X\to Y$ is continuous if and only if all of its components
$\pr_j\circ f$ are continuous. The topology on~$Y$
being initial with respect to the family $(\pr_j)_{j\in J}$,
the topology on $C(X,Y)$ is initial with respect to the family
$(C(X,\pr_j))_{j\in J}$ (Lemma~\ref{inipush}).
Thus $\Phi$ is a topological embedding
and hence a homeomorphism (being bijective).
\end{prf}
It is often useful that $C(X,\pl \, Y_j)=\pl \, C(X,Y_j)$.
\begin{lem}
Let $(J,\leq)$ be a directed set, $((Y_j)_{j\in J}, (\phi_{j,k})_{j\leq k})$
be a projective system of Hausdorff topological
spaces\footnote{Thus $\phi_{j,k}\colon Y_k\to Y_j$ is a continuous
map for $j,k\in J$ such that $j\leq k$,
with $\phi_{j,j}=\id_{Y_j}$
and $\phi_{j,k}\circ \phi_{k,\ell}=\phi_{j,\ell}$ if $j\leq k\leq \ell$.}
and $Y$ be its projective limit,
with the limit maps $\phi_j\colon Y\to Y_j$.
Then the topological space $C(X,Y)$ is the projective limit of
$((C(X,Y_j))_{j\in J}, (C(X,\phi_{j,k}))_{j\leq k})$,
together with the limit maps $C(X,\phi_j)\colon C(X,Y)\to C(X,Y_j)$.
\end{lem}
\begin{prf}
Let $P:=\prod_{j\in J}Y_j$ and $\pr_j\colon P\to Y_j$
be the projection onto the $j$th component.
The map
\[
\phi:=(\phi_j)_{j\in J}\colon Y\to P
\]
is a topological embedding with image 
\[ \phi(Y)=\{(x_j)_{j\in J}\in P\colon(\forall j,k\in J)\ 
j\leq k \impl x_j=\phi_{j,k}(x_k)\}.\]
Thus $C(X,\phi)$ is a topological embedding (see Lemma~\ref{ctsemb})
and hence so is
\[
\Phi\circ C(X,\phi) \colon C(X,Y)\to\prod_{j\in J}C(X,Y_j),
\]
using the homeomorphism $\Phi:=(C(X,\pr_j))_{j\in J}\colon C(X,P)\to
\prod_{j\in J}C(X,Y_j)$ from Lemma~\ref{cotopprod}.
The image of $\Phi\circ C(X,\phi)$ is contained in the projective limit
\[
L:=\Big\{(f_j)_{j\in J}\in\prod_{j\in J}C(X,Y_j)\colon
(\forall j,k\in J)\, j\leq k \impl f_j=\phi_{j,k}\circ f_k\Big\}\]
of the spaces $C(X,Y_j)$.
If $(f_j)_{j\in J} \in L$ and $x\in X$, then
$\phi_{j,k}(f_k(x))=f_j(x)$ for all $j\leq k$, whence
there exists $f(x)\in Y$ with $\phi_j(f(x))=f_j(x)$.
The topology on~$Y$ being initial with respect to the maps $\phi_j$,
we deduce from the continuity of the maps $\phi_j\circ f= f_j$
that $f\colon X\to Y$ is continuous.
Now $(\Phi\circ C(X,\phi))(f)=(f_j)_{j\in J}$.
Hence $\Phi\circ C(X,\phi)$ is a homeomorphism from
$C(X,Y)$ onto~$L$, and the assertions follow.
\end{prf}
Also composition operators (or pullbacks) $C(f,Y)$
are essential tools.
\begin{lem}\label{pubas}
Let $X_1$, $X_2$ and $Y$ be Hausdorff topological spaces.
If a map $f\colon X_1\to X_2$ is continuous, then also
the map
\[
C(f,Y)\colon C(X_2,Y)\to C(X_1,Y), \quad \gamma\mto \gamma\circ f
\]
is continuous.
\end{lem}
\begin{prf}
If $K\sub X_1$ is compact and $U\sub Y$ an open subset,
then $f(K)\sub X_2$ is compact and $C(f,Y)^{-1}(\lfloor K,U\rfloor )
= \lfloor f(K),U\rfloor$, 
since $(\gamma\circ f)(K)\sub U$ is equivalent to $\gamma (f(K))\sub U$
for $\gamma\in C(X_2,Y)$.
\end{prf}
\begin{rem}\label{resop}
If $X_2$ and $Y$ are Hausdorff topological spaces and $X_1\sub X_2$ a subset,
then the restriction map
\[
\rho\colon C(X_2,Y)\to C(X_1,Y),\quad \gamma\mto \gamma|_{X_1}
\]
is continuous. In fact, since $\gamma|_{X_1}=\gamma\circ f$
with the continuous inclusion map $f\colon X_1\to X_2$, $x\mto x$,
we have $\rho=C(f,Y)$ and Lemma~\ref{pubas} applies.
\end{rem}
\begin{lem}\label{coveremb}
Let $X$ and $Y$ be Hausdorff topological spaces
and $(X_j)_{j\in J}$ be a family of subsets of $X$ whose interiors $X_j^0$
cover~$X$.
Then the map
\[
\rho\colon C(X,Y)\to \prod_{j\in J}C(X_j,Y),\quad \gamma\mto (\gamma|_{X_j})_{j\in J}
\]
is a topological embedding with closed image.
\end{lem}
\begin{prf}
It is clear that the map $\rho$ is injective,
and it is continuous since each of its components is continuous,
by the preceding remark.
The image of $\rho$ consists of all $(\gamma_j)_{j\in J}\in
\prod_{j\in J}C(X_j,Y)$ such that
\[
(\forall j,k\in J)\,(\forall x\in X_j\cap X_k)\quad
\gamma_j(x)=\gamma_k(x)\,.
\]
The point evaluation $\ev_x\colon C(X_j,Y)\to Y$
and the corresponding one on $C(X_k,Y)$ are continuous if $x\in X_j\cap X_k$.
Since $Y$ is Hausdorff (and thus the diagonal is closed in $Y\times Y$),
it follows that $\im(\rho)$ is closed.
Since $\rho$ is injective, it will be an open map onto its image
if it takes the elements of a subbasis to relatively open subsets.
To verify this property,
let $K\sub X$ be compact and $U\sub Y$ be open.
Each $x\in K$ is contained in the interior $X_{j_x}^0$ (relative $X$) for
some $j_x\in J$.
Because $K$ is locally compact, $x$ has
a compact neighborhood $K_x$ in $K$
such that $K_x\sub X_{j_x}^0$.
By compactness, there is a finite subset $\Phi\sub K$
such that $K=\bigcup_{x\in \Phi}K_x$.
Then $W_x:=\lfloor K_x,U\rfloor\sub C(X_{j_x},U)$ is an open subset of $C(X_{j_x},Y)$
for all $x\in \Phi$. Hence
\[ W:=\Big\{(\gamma_j)_{j\in J}\in \prod_{j\in J}C(X_j,Y) \colon (\forall x\in \Phi)\;\gamma_{j_x}\in W_x\Big\} \] 
is open in
$\prod_{j\in J}C(X_j,Y)$.
Since $\rho(\lfloor K,U\rfloor)
=W\cap\im(\rho)$,
we see that $\rho(\lfloor K,U\rfloor)$ is relatively open in $\im(\rho)$,
which completes the proof.
\end{prf}
\begin{lem}\label{laeval}
If $X$ is a locally compact space and $Y$ a Hausdorff topological space,
then the evaluation map
\[
\ve\colon C(X,Y)\times X\to Y,\quad (\gamma,x)\mto \gamma(x)
\]
is continuous.
\end{lem}
\begin{prf}
Let $U\sub Y$ be open and $(\gamma,x)\in \ve^{-1}(U)$.
By local compactness, there exists a compact neighborhood $K\sub X$ of $x$
such that $K\sub \gamma^{-1}(U)$ .
Then $\lfloor K,U\rfloor\times K$
is a neighborhood of $(\gamma,x)$
in $C(X,Y)\times X$ and $\lfloor K,U\rfloor\times K\sub \ve^{-1}(U)$.
As a consequence, $\ve^{-1}(U)$
is open and thus $\ve$ is continuous.
\end{prf}
\begin{lem}
Let $X$, $Y$ and $Z$ be Hausdorff topological spaces.
If $Y$ is locally compact, then the composition map
\[
\Gamma\colon C(Y,Z)\times C(X,Y)\to C(X,Z),\quad (\gamma,\eta)\mto \gamma\circ \eta
\]
is continuous.
\end{lem}
\begin{prf}
Let $K\sub X$ be compact and $U\sub Z$ be an open subset.
Let $(\gamma,\eta)\in \Gamma^{-1}(\lfloor K,U\rfloor)$.
For each $x\in K$, there is a compact neighborhood~$L_x$
of $\eta(x)$ in~$Y$ such that $L_x\sub \gamma^{-1}(U)$.
We choose a compact neighborhood $M_x$
of $\eta(x)$ in~$Y$ such that $M_x\sub L_x^0$.
Then $K_x:=K\cap\eta^{-1}(M_x)$ is a compact
neighborhood of~$x$ in~$K$. Hence,
there is a finite subset $\Phi\sub K$ such that
$K=\bigcup_{x\in \Phi}K_x$.
Then $W:=\bigcap_{x\in \Phi}\lfloor K_x,L_x^0\rfloor$
is an open neighborhood of~$\eta$ in $C(X,Y)$
and $V:=\bigcap_{x\in \Phi}\lfloor L_x,U\rfloor$
is an open neighborhood of~$\gamma$ in $C(Y,Z)$.
We claim that $V\times W\sub \Gamma^{-1}(\lfloor K,U\rfloor)$.
If this is true, then $\Gamma^{-1}(U)$ is open and hence $\Gamma$
is continuous. To prove the claim, let $\sigma\in V$ and $\tau\in W$.
If $y\in K$, then $y\in K_x$ for some $x\in \Phi$.
Hence $\tau(y)\in \tau(K_x)\sub L_x^0$
and thus $\sigma(\tau(y))\in \sigma(L_x)\sub U$,
showing that $\sigma\circ \tau\in \lfloor K,U\rfloor$ indeed.
\end{prf}
\begin{lem}\label{fewercp}
Let $X$ and $Y$ be Hausdorff topological spaces
and $\cL$ be a set of compact subsets of~$X$.
Assume that,
for each compact subset $K\sub X$
and open subset $V\sub X$ with $K\sub V$,
there exist $n\in \N$ and $K_1,\ldots, K_n\in \cL$
such that $K\sub \bigcup_{i=1}^nK_i\sub V$.
Then the sets
\[
\lfloor K,U\rfloor,\quad\mbox{for $K\in \cL$ and open subsets $\,U\sub Y$,}
\]
form a basis for the compact-open topology on $C(X,Y)$.
\end{lem}
\begin{proof}
Given a compact subset $K\sub X$ and open subset
$U\sub Y$, let\linebreak
$\gamma\in \lfloor K,U\rfloor$.
Then $V:=\gamma^{-1}(U)$
is an open subset of $X$ that contains~$K$,
whence we can find $K_1,\ldots,K_n\in \cL$ as described in the hypotheses.
Then $\gamma\in \bigcap_{i=1}^n \lfloor K_i,U\rfloor\sub \lfloor K,U\rfloor$.
The assertion follows.
\end{proof}
\begin{rem}\label{defnkkR}
(a) Recall that a topological space~$X$
is called a \emph{$k$-space}
if~$X$ is Hausdorff and its topology is final with
respect to the inclusion maps $K\to X$,
for $K$ ranging through the set $\cK(X)$ of compact subsets of~$X$.
More explicitly, this means that
a subset $A\sub X$ is closed if and only if $A\cap K$ is closed
for each $K\in\cK(X)$.
Or, equivalently, a subset $U\sub X$ is open if and only if $U\cap K$
is (relatively) open in~$K$ for each $K\in\cK(X)$.
A Hausdorff space is a $k$-space if and only if mappings
$f\colon X\to Y$ from~$X$ to a topological space~$Y$
are continuous if and only if $f|_K$ is continuous for each compact subset
$K\sub X$ (see Exercise~\ref{exc-kviamap}).\medskip

\noindent
(b) A topological space~$X$ is called \emph{completely regular}
if $X$ is Hausdorff and its topology is initial with respect
to the set $C(X,\R)$ of all continuous
real-valued functions on~$X$.
It is easy to see that every locally convex topological vector space is completely
regular (Exercise~\ref{exc-lcxcoreg}).
We mention that,
more generally, every Hausdorff topological
group is completely regular~(see \cite[Thm.~8.4]{HR79}).\medskip

\noindent
(c) A topological space $X$ is called a \emph{$k_\R$-space}
if~$X$ is Hausdorff and real-valued functions $f\colon X\to\R$
on~$X$ are continuous if and only if $f|_K$ is continuous for each compact subset $K\sub X$.
If $X$ is a $k_\R$-space, then a map $f\colon X\to Y$ from~$X$ to a completely regular
topological space~$Y$ is continuous if and only if $f|_K$ is continuous
for each compact subset $K\sub X$ (by (b), 
this is clear from Lemma~\ref{lempropsinit}(b)).
\end{rem}
\begin{rem}
(a) Note that every $k$-space is a $k_\R$-space, in particular.
For each uncountable set~$I$ whose cardinality is at least that of~$\R$, the power $X:=\R^I$ of real lines
(with the product topology) is a $k_\R$-space (see~\cite{Nb70}
or also~\cite{GM20})
but not a $k$-space. In fact, $(X,+)$ has a non-discrete, closed subgroup~$A$
such that every compact subset of~$A$ is finite (cf.\ \cite{FHS17}).
If~$X$ was a $k$-space, then also the closed subset~$A$
would be a $k$-space (see Exercise~\ref{exc-kprem}).
Hence~$A$ would be discrete (since every compact subset of~$A$ is finite
and hence a discrete topological space), contradiction.\medskip

(b) If $X$ and $Y$ are $k$-spaces, then the product topology need not make $X\times Y$
a $k$-space (see Exercise~\ref{exc-product-not-k} for a counterexample).
\end{rem}
\begin{prop}\label{ctsexp}
Let $X$, $Y$ and $Z$ be Hausdorff topological spaces.
If $f\colon X\times Y\to Z$ is continuous, then also
\begin{equation}\label{deffv}
f^\vee\colon X\to C(Y,Z), \quad f^\vee(x):=f(x,\cdot)
\end{equation}
is continuous. Moreover, the map
\begin{equation}\label{expiso}
\Phi \colon C(X\times Y,Z)\to C(X,C(Y,Z)),\quad f\mto f^\vee
\end{equation}
is a topological embedding.
If $Y$ is locally compact, or $X\times Y$ is a $k$-space,
or $X\times Y$ is a $k_\R$-space and $Z$ is completely regular,
then $\Phi$ is a homeomorphism.
\end{prop}
\begin{prf}
Let $\gamma\in C(X\times Y,Z)$.
To see that $\gamma^\vee$ is continuous, let $K\sub Y$ be compact and $U\sub Z$ be an open subset.
If $x\in X$ such that $\gamma^\vee(x)\in \lfloor K,U\rfloor$,
then $\gamma(\{x\}\times K)\sub U$
and thus the product $\{x\}\times K$ of compact sets is contained in the open subset
$\gamma^{-1}(U)$.
By the Wallace Lemma~\ref{Wallla}, 
there is an open subset $V\sub X$ with $\{x\}\sub V$
and $V \times K\sub \gamma^{-1}(U)$.
Then $V\sub (\gamma^\vee)^{-1}(\lfloor K,U\rfloor)$,
showing that  $(\gamma^\vee)^{-1}(\lfloor K,U\rfloor)$
is a neighborhood of~$x$ and hence open (as $x$ was arbitrary).
Thus $\gamma^\vee$ is continuous.

To see that $\Phi$ is continuous, recall that
the sets $\lfloor L,U\rfloor$, with $L\in \cK(Y)$ and open subsets $U\sub Z$,
form a subbasis of the compact-open topology on $C(Y,Z)$.
Hence,
by Lemma~\ref{subbaca},
the sets $\lfloor K,\lfloor L,U\rfloor\rfloor$
form a subbasis of the compact-open topology on $C(X,C(Y,Z))$,
for $L$ and $U$ as before and $K\in \cK(X)$.
We clearly have 
\begin{equation}\label{hlp}
\Phi^{-1}(\lfloor K,\lfloor L,U\rfloor\rfloor)=
\lfloor K\times L, U\rfloor, 
\end{equation}
so that $\Phi$ is continuous.
Because $\Phi$ is (obviously) injective,
it will be open onto its image if it takes
open subsets in a subbasis to relatively open subsets.
Let $\cL\sub \cK(X\times Y)$ be the set of all products $K\times L$,
where $K\in\cK(X)$ and $L\in \cK(Y)$.
We claim that $\cL$ satisfies
the hypotheses of Lemma~\ref{fewercp}
(for the function space $C(X\times Y,Z)$).
If this is so, then the sets $\lfloor K\times L,U\rfloor$
with $K\in \cK(X)$ and $L\in \cK(Y)$
form a subbasis for the compact-open topology on $C(X\times Y,Z)$,
and now~(\ref{hlp}) shows that $\Phi$ is open onto its image
(and hence a topological embedding).

To verify the claim,
let $\pi_1\colon X\times Y\to X$ and
$\pi_2\colon X\times Y\to Y$ be the projection onto the first and second component, respectively.
If $M\sub X\times Y$ is compact and $V\sub X\times Y$ an open subset
such that $K\sub V$, then $M\sub M_1\times M_2$,
where $M_1:=\pi_1(M)$ and $M_2:=\pi_2(M)$ are compact.
For each $(x,y)\in M$,
there is an open neighborhhod $U_{x,y}\sub X$ of $x$ and
and open neighborhood $V_{x,y}\sub Y$ of~$y$ such that
$U_{x,y}\times V_{x,y}\sub V$.
Because $M_1$ and $M_2$ are locally compact, there exist
compact neighborhoods $K_{x,y}\sub M_1$ of $x$
and $L_{x,y}\sub M_2$ of~$y$
such that $K_{x,y}\sub U_{x,y}$ and $L_{x,y}\sub V_{x,y}$.
Then
\[ K_{x,y}\times L_{x,y}
\sub U_{x,y}\times V_{x,y}\sub V.\]
By compactness, $K$ is covered by finitely many of the sets
$K_{x,y}\times L_{x,y}$. Thus all hypotheses
of Lemma~\ref{fewercp} are indeed satisfied.

Now assume that $Y$ is locally compact.
Let $\eta\in C(X,C(Y,Z))$.
Because the evaluation map $\ve\colon C(Y,Z)\times Y\to Z$
is continuous by Lemma~\ref{laeval},
the map $\gamma:=\eta^\wedge:= \ve\circ (\eta\times \id_Y)\colon X\times Y\to Z$,
$(x,y)\mto \eta(x)(y)$ is continuous.
Since $\Phi(\gamma)=(\eta^\wedge)^\vee=\eta$,
we see that $\Phi$ is surjective and hence (being an embedding) a homeomorphism.

Finally, assume that $X\times Y$ is a $k$-space (resp., that $X\times Y$
is a $k_\R$-space and $Z$ is completely regular).
Again, we only need to show that $\Phi$ is surjective.
Let $\eta\in C(X,C(Y,Z))$
and define $\gamma:=\eta^\wedge\colon X\times Y\to Z$, $(x,y)\mto \eta(x)(y)$.
If we can show that~$\gamma$ is continuous, then $\Phi(\gamma)=\eta$
(as required).
Because $X\times Y$ is a $k$-space (resp., $X\times Y$ is a $k_\R$-space and $Z$ is completely
regular), it suffices to show
that $\gamma|_K$ is continuous for each compact subset
$K\sub X\times Y$ (Remark~\ref{defnkkR}(c)).
This will follow if $\gamma|_{K_1\times K_2}$
is continuous for all compact subsets $K_1\sub X$ and $K_2\sub Y$
(given $K$ as before, $K$ is contained in the compact set $K_1\times K_2$
with $K_1 :=\pi_1(K)$, $K_2:= \pi_2(K)$).
We now use that the restriction map
$\rho\colon C(Y,Z)\to C(K_2,Z)$
is continuous (see Remark~\ref{resop})
and hence also the map
\[
\zeta := \rho\circ \eta|_{K_1} \colon K_1\to C(K_2,Z).
\]
Since $K_2$ is compact,
$\zeta^\wedge\colon K_1\times K_2\to Z$ is continuous
(by the preceding part of the proof).
Because $\gamma|_{K_1\times K_2}=\zeta^\wedge$,
the continuity of $\gamma$
follows.
\end{prf}
The $k_\R$-version of the exponential law goes back to \cite{Hu71}.
\begin{rem}\label{rem-cts-exp}
(a) Given Hausdorff topological spaces $X$, $Y$, $Z$ and a map
$\eta\colon X\to C(Y,Z)$, define
\[
\eta^\wedge\colon X\times Y\to Z,\quad (x,y)\mto \eta(x)(y).
\]
The preceding proposition entails:
If $\eta^\wedge$ is continuous, then $\eta=(\eta^\wedge)^\vee$
is continuous.
If $Y$ is locally compact, or $X\times Y$ is a $k$-space, or $X\times Y$ is a $k_\R$-space
and $Z$ is completely regular,
then $\eta=(\eta^\wedge)^\vee$ is continuous if and only if
$\eta^\wedge$ is continuous.\medskip

(b) The map $\Phi$ in Proposition~\ref{ctsexp} need not be surjective (see Exercise~\ref{exc-no-explaw}).
\end{rem}
\begin{lem}\label{tocon}
Let $X$ and $Y$ be Hausdorff topological spaces.
For $y\in Y$, let $c_y\colon X\to Y$ be the constant map $x\mto y$.
Then the map
\[
c\colon Y\to C(X,Y),\quad y\mto c_y
\]
is continuous
and $($if $X\not=\emptyset)$
in fact a topological embedding.
\end{lem}
\begin{prf}
If $K\sub X$ is compact and $U\sub Y$ an open subset,
then $c^{-1}(\lfloor K, U\rfloor)=U$ (if $K\not=\emptyset$)
and $c^{-1}(\lfloor K, U\rfloor)=Y$ (if $K=\emptyset$).
In either case, the preimage is open, entailing that $c$ is continuous.
If $X$ is not empty, we pick $x\in X$.
Then the point evaluation $\ev_x\colon C(X,Y)\to Y$, $\gamma\mto \gamma(x)$
is continuous and $\ev_x\circ c=\id_Y$,
entailing that $c$ is a topological embedding.
\end{prf}
By the next lemma, so-called pushforwards
are continuous.
\begin{lem}\label{ctspfapp}
Let $X$, $Y$ and $Z$ be Hausdorff topological spaces
and\linebreak
$f\colon X\times Y\to Z$ be continuous.
Then also the following map is continuous:
\[
f_*\colon C(X,Y)\to C(X,Z),\quad \gamma\mto f\circ (\id_X,\gamma).
\]
\end{lem}
Note that $f_*(\gamma)(x)=f(x,\gamma(x))$. 
\begin{prf}
Identifying $C(X,X)\times C(X,Y)$ with $C(X,X\times Y)$ as in
Lemma~\ref{cotopprod},
we can write
$f_*(\gamma)=C(X,f)(\id_X,\gamma)$.
Since $C(X,f)$ is continuous by Lemma~\ref{covsuppo},
the continuity of $f_*$ follows.
\end{prf}
\noindent
We also have two versions with parameters.
\begin{lem}\label{pushpar}
Let $X$, $Y$, $Z$ and $P$ be Hausdorff topological spaces
and\linebreak
$f\colon X\times Y\times P\to Z$ be a continuous map.
For $p\in P$, abbreviate
\[ f^p  :=f(\cdot,p)\colon X\times Y\to Z.\] 
Then also the following map is continuous:
\[
C(X,Y)\times P\to C(X,Z),\quad (\gamma,p)\mto f^p\circ (\id_X,\gamma).
\]
\end{lem}
\begin{prf}
Let $c\colon P\to C(X,P)$, $p\mto c_p$ be the continuous map discussed in
Lemma~\ref{tocon}.
The map $C(X,f)\colon C(X,X\times Y\times P)\to C(X,Z)$ is continuous
(by Lemma~\ref{covsuppo}).
Identifying $C(X,X)\times C(X,Y)\times C(X,P)$ with\linebreak
$C(X,X\times Y\times P)$,
we have $f^p\circ (\id_X,\gamma)=C(X,f)(\id_X,\gamma,c_p)$,
which is continuous in $(\gamma,p)$.
\end{prf}
\begin{lem}\label{parcov}
Let $X$, $Y$, $Z$ and $P$ be Hausdorff topological spaces
and\linebreak
$f\colon Y\times P\to Z$ be a continuous map.
For $p\in P$, set $f^p:=f(\cdot,p)\colon Y\to Z$.
Then also the following map is continuous:
\begin{equation}\label{phi}
\psi\colon C(X,Y)\times P\to C(X,Z),\quad (\gamma,p)\mto f^p\circ \gamma.
\end{equation}
\end{lem}
\begin{prf}
Since $g\colon X\times Y\times P \to Z$, $g(x,y,p):=f(y,p)$
is continuous, so is
$\psi\colon C(X,Y)\times P\to C(X,Z)$, $(\gamma,p)\mto g^p\circ (\id_X,\gamma)$,
by Lemma~\ref{pushpar}.
\end{prf}
\begin{numba}\label{hemic}
Recall that a topological space~$X$ is called \emph{$\sigma$-compact}
if it can be written as a union $X=\bigcup_{n\in\N} K_n$ of compact subsets $K_n\sub X$.
A topological space~$X$ is called
\emph{hemicompact}
if there exists a sequence $(K_n)_{n\in\N}$ of compact subsets
of~$X$ such that each compact subset of~$X$ is contained in some~$K_n$
(then $X$ is $\sigma$-compact in particular, as $X=\bigcup_{n\in\N}K_n$).
\end{numba}
For example, every $\sigma$-compact locally compact
space is hemicompact (see Exercise~\ref{lcp-hemi}).
\begin{lem}\label{sammelsu}
If $X$ is a Hausdorff topological space and $G$ a topological group,
then the following holds:
\begin{description}[(D)]
\item[{\rm(a)}]
$C(X,G)$, endowed with the compact open topology, 
is a topological group with respect to the pointwise
group 
operations.
\item[{\rm(b)}]
Let $\cL$ be a set of compact subsets of $X$ such that each $K\in \cK(X)$
is contained in some $L\in \cL$.
Also, let $\cU$ be a basis of open identity neighborhoods in~$G$.
For $\gamma\in C(X,G)$,
the sets $\gamma\cdot \lfloor K, U\rfloor$
then form a basis of open neighborhoods of $\gamma$ in $C(X,G)$
$($for $K\in \cL$ and $U\in \cU)$, and so do the sets
$\lfloor K, U\rfloor\cdot \gamma$.
The compact-open topology therefore coincides with the topology of
uniform convergence on compact sets,
both with respect to the left and also the right uniformity on~$G$.
\item[{\rm(c)}]
If $X$ is hemicompact and $G$ is metrizable, then
$C(X,G)$ is metrizable.
\item[{\rm(d)}]
If $X$ is a $k_\R$-space $($e.g., if $X$ is a $k$-space$)$ and $G$ is complete,
then $C(X,G)$ is complete.
\item[{\rm(e)}]
If $X$ is a $k_\R$-space and $G$ is sequentially complete,
then $C(X,G)$ is\linebreak
sequentially complete.
\item[{\rm(f)}]
If $E$ is a topological vector space over
$\K\in\{\R,\C\}$, then the pointwise operations make
$C(X,E)$ a topological $\K$-vector space.
Moreover, $C(X,\K)$ is a topological $\K$-algebra and
$C(X,E)$ is a topological $C(X,\K)$-module.
\item[{\rm(g)}]
If $\K \in \{\R,\C\}$ and $E$ is a locally convex topological $\K$-vector space,
then also $C(X,E)$ is locally convex.
\item[{\rm(h)}]
If $X$ is a $k_\R$-space and $E$ is a Mackey complete locally convex space,
then also $C(X,E)$ is Mackey complete.
\item[{\rm(i)}]
If $X$ is a $k_\R$-space and $E$ is a quasi-complete locally convex space,
then also $C(X,E)$ is quasi-complete.
\end{description}
\end{lem}
\begin{prf}
(a) Since $G$ is a topological group, the group multiplication\linebreak
$\mu\colon G\times G\to G$, $(x,y)\mto xy$ and the inversion map $\iota\colon G\to G$, $x\mto x^{-1}$
are continuous. Identifying $C(X,G)\times C(X,G)$ with $C(X,G\times G)$ (as in Lemma~\ref{cotopprod}),
the group multiplication of $C(X,G)$ is the mapping 
\[ C(X,\mu)\colon C(X,G\times G)\to C(X,G),\] 
which is continuous by Lemma~\ref{covsuppo}.
The group inversion is the map $C(X,\iota)\colon C(X,G)\to C(X,G)$
and hence continuous as well. Thus $C(X,G)$ is a topological group.

(b) If $K_1,\ldots, K_n\in \cK(X)$ and $U_1,\ldots, U_n\sub G$ are open identity neighborhoods,
we find $U\in \cU$ such that $U\sub U_1\cap\cdots\cap U_n$,
and $K\in \cL$ such that $K_1\cup\cdots\cup K_n\sub K$.
Then $\bigcap_{j=1}^n\lfloor K_j,U_j\rfloor\supseteq \lfloor K,U\rfloor$.
Therefore the sets $\lfloor K,U\rfloor$ with $K\in \cL$ and $U\in \cU$ form a basis
of open identity neighborhoods for $C(X,G)$.
Since left and right translations in the topological group $C(X,G)$ are homeomorphisms,
the remainder of (b) follows.

(c) Let $K_1\sub K_2\sub \cdots$ be an ascending sequence of compact subsets of~$X$, with union $X$,
such that each compact subset of $X$ is contained in some $K_n$.
Also, let $U_1\supseteq U_2\supseteq\cdots$ be a descending sequence of open identity neighborhoods in $G$
which gives a basis for the filter of identity neighborhoods.
By~(b), the sets $\lfloor K_n, U_n\rfloor$, for $n\in \N$, provide a countable basis
of identity neighborhoods in $C(X,G)$.
As a consequence, the topological group
$C(X,G)$ is metrizable \cite[Thm.~8.3]{HR79}. 

(d) Let $(\gamma_a)_{a\in A}$ be a Cauchy net in $C(X,G)$.
For each $x\in X$, the point evaluation $\ev_x\colon C(X,G)\to G$ is a continuous homomorphism.
Hence $(\gamma_a(x))_{a\in A}$ is a Cauchy net in~$G$
and hence convergent to some element $\gamma(x)\in G$
(as  $G$ is assumed complete).
For each compact subset $K\sub X$, the restriction map $\rho_K\colon C(X,G)\to C(K,G)$ is a continuous
homomorphism (cf.\ Remark~\ref{resop}), whence $(\gamma_a|_K)_{a\in A}$
is a Cauchy net in $C(K,G)$.
We claim that $C(K,G)$ is complete. If this is true,
then $\gamma_a|_K\to \gamma_K$
for some continuous map $\gamma_K\in C(K,G)$.
Since $\gamma(x)=\gamma_K(x)$ for each $x\in K$, we see that $\gamma|_K=\gamma_K$
is continuous. Because $X$ is a $k_\R$-space and~$G$ (being a topological group) is completely regular, this implies that~$\gamma$
is continuous. If $K\in \cK(X)$ and $U\sub G$ is an identity neighborhood,
then $(\gamma_K)^{-1}\gamma_a|_K\in \lfloor K,U\rfloor$ inside $C(K,G)$
for sufficiently large~$a$ (see (b))
and hence $\gamma^{-1}\gamma_a\in \lfloor K,U\rfloor$ inside $C(X,G)$,
showing that $\gamma_a\to \gamma$.

To prove the claim, we may assume that $X=K$ is compact.
Let $U\sub G$ be an open identity neighborhood
and $V\sub U$ be an open identity neighborhood
with closure $\overline{V}\sub U$.
There is $a\in A$ such that $\gamma_b^{-1}\gamma_c\in \lfloor K, V\rfloor$
for all $b,c\geq a$ in $A$. Thus, for all $x\in K$, $\gamma_b(x)^{-1}\gamma_c(x)\in V$.
Passing to the limit in $c$, we obtain $\gamma_b(x)^{-1}\gamma(x)\in \overline{V}\sub U$
and thus
\begin{equation}\label{sireu}
\gamma(x)\in \gamma_b(x)U, \quad \mbox{for all $\,x\in K$ and $b\geq a$.}
\end{equation}
If $W$ is any identity neighborhood in $G$,
we can choose $U$ from before so small that $U^{-1}UU\sub W$.
By continuity of $\gamma_a$, each $x\in K$ has a neighborhood $L\sub K$
such that $\gamma_a(y)^{-1}\gamma_a(x)\in U$ for all $y\in L$.
Combining this with~(\ref{sireu}), we see that
\[
\gamma^{-1}(y)\gamma(x)\in U^{-1}\gamma_a(y)^{-1}\gamma_a(x)U\sub U^{-1}UU\sub W
\]
for all $y\in L$. Thus $\gamma$ is continuous at~$x$ and hence continuous.
Finally, we have $\gamma_b^{-1}\gamma\in \lfloor K,U\rfloor$ for all $b\geq a$,
by (\ref{sireu}). Hence $\gamma_b\to \gamma$ in $C(K,G)$.

(e) Proceed as in (d), replacing the Cauchy net with a Cauchy sequence.

(f) By (a), $C(X,E)$ is a topological group.
Let $\sigma\colon \K\times E\to E$, $(z,v)\mto zv$ be multiplication with scalars.
Then the map $\K\times C(X,E)\to C(X,E)$, $(z,\gamma)\mto \sigma(z,\cdot) \circ \gamma$
is continuous (by Lemma~\ref{parcov}).
As this is the multiplication by scalars in $C(X,E)$,
the latter is a topological vector space.
Let $\mu\colon \K\times \K\to\K$ be the multiplication
in the field $\K$. Identifying \break $C(X,\K)\times C(X,\K)$
with $C(X,\K\times\K)$,
the algebra multiplication of $C(X,\K)$ is the map
$C(X,\mu)\colon C(X,\K\times\K)\to C(X,\K)$,
which is continuous by Lemma~\ref{covsuppo}.
Hence $C(X,\K)$ is a topological algebra.
Identifying \break $C(X,\K)\times C(X,E)$
with $C(X,\K\times E)$,
the $C(X,\K)$-module multiplication on $C(X,E)$ is the map
$C(X,\sigma)\colon C(X,\K\times E)\to C(X, E)$,
which is continuous by Lemma~\ref{covsuppo}.
Hence $C(X,E)$ is a topological $C(X,\K)$-module.

(g) By hypothesis, the open convex $0$-neighborhoods $U\sub E$ form a basis 
of $0$-neighborhoods in $E$. For each $K\in \cK(X)$,
the set $\lfloor K,U\rfloor$ then is an open $0$-neighborhood
in $C(X,E)$. 
Part\,(b) implies that these sets form a basis of $0$-neighborhoods
in $C(X,E)$. Thus $C(X,E)$ is locally convex.

(h) Let $\gamma\colon [0,1]\to C(X,E)$ be a Lipschitz curve.
If $\wt{E}$ is a completion of~$E$ such that $E\sub\wt{E}$, then $C(X,\wt{E})$
is complete (by~(d)). The weak integral $\theta:=\int_0^1\gamma(t)\,dt$ therefore
exists in $C(X,\wt{E})$ (see Proposition~\ref{rieman}).
For $x\in X$, let $\ve_x\colon C(X,E)\to E$
and $\wt{\ve}_x\colon C(X,\wt{E})\to\wt{E}$
be the respective point evaluation $\eta\mto\eta(x)$
at~$x$, which is continuous and linear.
Then $\ve_x\circ \gamma\colon [0,1]\to E$ is
a Lipschitz curve, entailing that $\zeta(x):=\int_0^1(\ve_x\circ \gamma)(t)\,dt$
exists in~$E$. Using
Exercise~\ref{excmorewint} (a) and (b),
we find that
\[
\theta(x)=\wt{\ve}_x(\theta)=\int_0^1 (\wt{\ve}_x\circ\gamma)(t)\,dt =\zeta(x).
\]
As a consequence, $\theta\in C(X,E)$ and hence $\theta=\int_0^1\gamma(t)\, dt$ in $C(X,E)$,
using Exercise~\ref{excmorewint}(b) again.

(i) Let $\wt{E}$ be a completion of~$E$ with $E\sub\wt{E}$.
If $(\gamma_\alpha)_{\alpha\in A}$ is a bounded Cauchy net in $C(X,E)$,
then $\gamma_\alpha$ converges in $C(X,\wt{E})$ to a continuous function
$\gamma\colon X\to\wt{E}$, by~(d). For each $x\in X$, the point evaluation
$C(X,E)\to E$, $\eta\mto \eta(x)$ is continuous and linear,
whence $(f_\alpha(x))_{\alpha\in A}$ is a bounded
Cauchy net in~$E$ and hence convergent, entailing that $f(x)\in E$.
Then $f\in C(X,E)$ and $f_\alpha\to f$ also in $C(X,E)$.
\end{prf}
\begin{prop}
 \label{prop:B.6} Let $X$ and $Y$ be Hausdorff topological spaces. 
On $C(X,Y)$, the compact-open topology coincides with the 
\index{graph topology}  
{\it graph topology}, 
i.e., the topology generated by the sets of the form 
$$C(X,Y)_{U,K} := \{ f \in C(X,Y) \: \Gamma(f\res_K) \subeq U\}, $$
where $U \subeq X \times Y$ is open, $K \subeq X$ is compact, 
and $\Gamma(f) \subeq X \times Y$ is the graph of $f$. 

If, in addition, $X$ is compact, then a basis for the graph topology is given by 
the sets 
$$C(X,Y)_{U} := \{ f \in C(X,Y) \: \Gamma(f) \subeq U\}, $$
where $U \subeq X \times Y$ is open. 
\end{prop}

\begin{prf} Let $f \in C(X,Y)$, $K \subeq X$ compact and $U \supeq \Gamma(f\res_K)$ be 
an open subset 
of $X \times Y$. 
Then there exists for each $x \in X$ a compact neighborhood 
$K_x$ of $x$ in $K$ and an open neighborhood $U_{f(x)}$ of $f(x)$ in $Y$ with 
$K_x \times U_{f(x)} \subeq U$ and $f(K_x) \subeq U_{f(x)}$. 
Covering $K$ with finitely many sets $K_{x_i}$, $i =1,\ldots, n$, we see that 
$$ \bigcap_{i = 1}^n \lfloor K_{x_i}, U_{f(x_i)}\rfloor  \subeq C(X,Y)_{U,K}. $$
This implies that each set $C(X,Y)_{U,K}$ is open in the compact-open topology. 

Conversely, let $K \subeq X$ be compact and $O \subeq Y$ open. 
Then 
$$
\lfloor K,O\rfloor  = \{ f \in C(X,Y) \: \Gamma(f\res_K) \subeq X \times O\} 
= C(X,Y)_{X \times O,K} $$
is open in the graph topology. We conclude that the graph topology 
coincides with the compact-open topology. 

Assume, in addition, that $X$ is compact. The system of the 
sets $C(X,Y)_U$ is stable under intersections, hence a basis for the 
topology it generates. Each set 
$C(X,Y)_U = C(X,Y)_{U,X}$ is open in the graph topology. 
If, conversely, $K \subeq X$ is compact and $U \subeq X \times Y$ 
is open with $f \in C(X,Y)_{U,K}$, then 
$V := \big((X \setminus K) \times Y\big) \cup U$ is an open subset 
of $X \times Y$ with 
$f \in C(X,Y)_V \subeq C(X,Y)_{U,K}$. This completes the proof. 
\end{prf} 
\subsection*{Exercises for Section~\ref{appcotop}}

\begin{small}

\begin{exer}\label{exc-Kelleyfic}
Let $X$ be a Hausdorff topological space,
$\cO$ be its topology and $\cK(X)$ be the set of
all compact subsets $K\sub X$.
Let $\cO_k$ be the final topology on~$X$
with respect to the inclusion maps $K\to X$
for $K\in\cK(X)$.
Show that $(X,\cO_k)$ is a $k$-space
and $\cO\sub \cO_k$.
Moreover, $\cO_k\sub\cT$ for every topology $\cT$
on~$X$ such that $\cO\sub\cT$ and $(X,\cT)$ is a $k$-space.
\end{exer}

\begin{exer}\label{exc-kviamap}
Let $(X,\cO)$ be a Hausdorff topological space.
Show that $X$ is a $k$-space if and only if $X$ has the following property:
a map $f\colon X\to Y$ to a topological space~$Y$ is continuous
if and only if $f|_K$ is continuous for each compact set $K\sub X$.\\[2mm]
[Hint: Choose $Y:=X$, endowed with the topology
$\cO_k$ as in Exercise~\ref{exc-Kelleyfic}.]
\end{exer}

\begin{exer}\label{exc-lcxcoreg}
Show that every locally convex space~$E$ is
completely regular as a topological space.\\[2mm]
[Hint: The translates of open unit balls $B^q_1(0)=q^{-1}(]{-1},1[)$
of continuous seminorms $q\colon E\to\R$ form a basis for the
topology of~$E$.]
\end{exer}

\begin{exer}\label{exc-altcreg}
\begin{description}[(D)]
\item[(a)]
Show that a Hausdorff topological space is completely regular if and only if
it has the following property: For each $x\in X$ and neighborhood~$U$
of~$x$ in~$X$, there is a continuous function
$f\colon X\to[0,1]$ such that $f(x)=1$ and $\Supp(f):=\wb{\{y\in X\colon f(y)\not=0\}}\sub U$.
\item[(b)]
Show that one can even achieve that $f|_V=1$ for a neighborhood $V\sub U$ of~$x$.
\end{description}
[Hint: Let $X$ be completely regular. After shrinking $U$, we may assume that $U=\bigcap_{j=1}^nf_j^{-1}(]a_j-\ve,a_j+\ve[)$
for certain $f_1,\ldots, f_n\in C(X,\R)$, $a_1,\ldots, a_n\in\R$ and some $\ve>0$.
Then $g:=\sum_{j=1}^n(f_j-a_j)^2/\ve^2\geq 0$ and
$g^{-1}(]{-1},1[)\sub U$. Now build further functions (and, eventually, a suitable~$f$) with the help of~$g$.]
\end{exer}

\begin{exer}\label{exc-kprem}
Let $X$ be a $k$-space and $Y\sub X$ be a subset which is closed or open.
Show that $Y$ is a $k$-space in the topology induced by~$X$.\\[2mm]
[Hint: If $Y$ is open in~$X$ and $K\sub X$ a compact subset, then every $x\in K\cap Y$
has a compact neighborhood~$L$ in~$K$ which is contained in the (relatively) open subset
$K\cap Y$ of~$K$.]
\end{exer}

\begin{exer}
Let $X$ be a completely regular topological space.
Then every subset $Y\sub X$ is completely regular in the topology induced by~$X$.
\end{exer}

\begin{exer}
Let $X$ be a topological space which is a $k_\R$-space and completely regular.
Show that every open subset $U\sub X$ is $k_\R$ in the induced topology.\\[2mm]
[Hint: Exercise~\ref{exc-altcreg}(b) may be useful.]
\end{exer}

\begin{exer}
Let $X$ and $Y$ be Hausdorff topological spaces
whose underlying sets are non-empty.
Show that if $X\times Y$ is a $k_\R$-space, then $X$ and $Y$
are $k_\R$-spaces.
\end{exer}

\begin{exer}
Let $X$ be a $k$-space (resp., a $k_\R$-space), $Y$ be a Hausdorff space and $q\colon X\to Y$ be a quotient map.
Show that~$Y$ is a $k$-space (resp., a $k_\R$-space).
\end{exer}

\begin{exer}\label{lcp-hemi}
Show that every
locally compact, $\sigma$-compact topological space~$X$
has an \emph{exhaustion by compact sets},
i.e.,  a sequence $(K_n)_{n\in\N}$
of compact subsets of~$X$ such that $\bigcup_{n\in\N}K_n=X$
and $K_n\sub K_{n+1}^0$ for each $n\in\N$.
Deduce that~$X$ is hemicompact.
\end{exer}

\end{small}
\section{Direct limits of topological spaces}\label{app-basic-DL}
We now discuss a special type of final topology,
the direct limit topology on an ascending union of topological spaces.
\begin{defn}\label{defnDLtop}
Let $X_1\sub X_2\sub\cdots$ be an ascending sequence of topological
spaces $(X_n,\cO_n)$ such that the inclusion map
$j_{m,n}\colon X_n\to X_m$ is continuous for all $n\leq m$ in~$\N$.
We then call $X_1\sub X_2\sub\cdots$ a \emph{direct sequence of topological spaces}.
The final topology $\cO_{\DL}$ on $X:=\bigcup_{n\in\N}X_n$ with respect to
the inclusion maps $j_n\colon X_n\to X$
is called the \emph{direct limit topology} on~$X$.
If~$X$ is endowed with the direct limit topology, then we say that
$X$ is the \emph{direct limit} of the topological spaces~$X_n$, or also that
\[
X=\dl X_n\vspace{-.5mm}
\]
as a topological space.
If $j_{n+1,n}$ (and hence each $j_{m,n}$) is a topological embedding
for each $n\in\N$,
then $X_1\sub X_2\sub\cdots$ is called a \emph{strict} direct sequence.
The direct sequence is called \emph{compactifying}
if, for each $n\in \N$ and $x\in X_n$, there exists a neighborhood~$V$ of~$x$ in~$X_n$
and $m\geq n$ such that the closure $\wb{V}$ of~$V$ in $(X_m,\cO_m)$ is compact. 
\end{defn}
\begin{rem}\label{firstremDL}
(a) By definition of the final topology, a subset $M\sub X$ is open in $(X,\cO_{\DL})$
if and only if $M\cap X_n$ is open in~$(X_n,\cO_n)$
for each $n\in\N$.
Likewise, a subset $M\sub X$ is closed in $(X,\cO_{\DL})$ if and only if $M\cap X_n$ is closed in~$(X_n,\cO_n)$
for each $n\in\N$.

(b)
Each of the inclusion maps $j_n\colon (X_n,\cO_n)\to (X,\cO_{\DL})$ is continuous,
by definition of the final topology $\cO_{\DL}$,
and $\cO_{\DL}$ is the finest topology with this property.

(c) If each $X_n$ is compact, then $X_1\sub X_2\sub\cdots$ is strict.
If each $X_n$ is locally compact, then $X_1\sub X_2\sub\cdots$ is compactifying.
Further examples of compactifying direct sequences
will be encountered in Section~\ref{sec-appDLvec} in connection with so-called Silva spaces.
In this case, each~$E_n$ is a Banach space and each inclusion map $E_n\to E_{n+1}$
is a compact operator.

(d)
For each subsequence $(X_{n_k})_{k\in\N}$
of $(X_n)_{n\in\N}$, we have
\[
\dl X_n=\dl X_{n_k}.\vspace{-1mm}
\]
In fact, if $V\sub \dl X_n$\vspace{-.5mm} is open,
then $V\cap X_{n_k}$ is open for each~$k$, whence $V$ is open in $\dl X_{n_k}$.
If, conversely, $V\sub \dl X_{n_k}$\vspace{-.5mm} is open,
given $n\in\N$ we find $k\in\N$ such that $n_k\geq n$.
Then $V\cap X_{n_k}$ is open in $X_{n_k}$, entailing that
$V\cap X_n=j_{n_k,n}^{-1}(V\cap X_{n_k})$ is open in~$X_n$. Thus~$V$
is open in $\dl X_n$.\vspace{-.5mm}
\end{rem}
\begin{lem}\label{ctsonsteps}
Let $X_1\sub X_2\sub\cdots$ be a direct sequence of topological spaces
and $\cO$ be a topology on $X=\bigcup_{n\in\N} X_n$.
Then $(X,\cO)=\dl X_n$\vspace{-.5mm} as a topological space
if and only if the following two conditions are satisfied:
\begin{description}[(DD)]
\item[\rm(i)]
The inclusion map $j_n\colon (X_n,\cO_n)\to (X,\cO)$
is continuous for all~$n\in\N$;
\item[\rm(ii)]
For every
map $f\colon X\to Y$ to a topological space, $f$ is continuous on $(X,\cO)$
if and only if $f|_{X_n}$ is continuous on $(X_n,\cO_n)$
for all $n\in\N$.
\end{description}
\end{lem}
\begin{prf}
The direct limit topology $\cO_{\DL}$ satisfies both conditions
as it is the final topology with respect to $(j_n)_{n\in\N}$
(see Definition~\ref{defnfinal} and Lemma~\ref{fctnonfinal}).

If, conversely, $\cO$ satisfies (i) and (ii), then $\cO\sub \cO_{\DL}$ (see Remark~\ref{firstremDL}(b)).
If $h\colon (X,\cO)\to(X,\cO_{\DL})$,
$h(x):=x$, then $h\circ j_n=j_n\colon X_n\to (X,\cO_{\DL})$
is continuous and hence so is~$h$, by (ii) for~$\cO$. Thus $\cO_{\DL}\sub \cO$
and $\cO=\cO_{\DL}$.
\end{prf}
Some observations concerning compactifying sequences will be useful.
\begin{lem}\label{compactify-la}
Let $X_1\sub X_2\sub\cdots$ be a compactifying
direct sequence of Hausdorff topological spaces
and $n\in\N$. Then the following holds.
\begin{description}[(DD)]
\item[\rm(a)]
For each compact subset $K\sub X_n$,
there exists an open subset $V\sub X_n$ and $m\geq n$ such that
$K\sub V$ and $V$ has compact closure $\oline V$ in~$X_m$.
If $U$ is an open subset of $X_m$ such that $K\sub U$,
we can achieve that, moreover, $\wb{V}\sub U$
$($after shrinking $V$ if necessary$)$.
\item[\rm(b)]
If $K$ and $L$ are disjoint compact subsets of~$X_n$,
then there exist open subsets $V$ and $W$ in $X_n$
and $m\geq n$ such that $K\sub V$ and $L\sub W$,
and moreover the closures $\wb{V}$ and $\wb{W}$ in~$X_m$
are compact and disjoint.
\end{description}
\end{lem}
\begin{prf}
(a) Omitting a trivial case, we may assume that $K\not=\emptyset$.
For each $x\in K$, there exist an open neighborhood $W_x$ of~$x$ in~$X_n$
and $m_x\geq n$ such that $W_x$ has compact closure $\wb{W_x}$ in $X_{m_x}$.
Since~$K$ is compact, we have $K\sub\bigcup_{x\in F}W_x=:W$ for a finite subset $F\sub K$.
Set $m:=\max\{m_x\colon x\in F\}$. Since $\wb{W_x}$ is also compact
in~$X_m$ for each $x\in F$, we see that the closure of~$W$ in $X_m$ is given by
$\wb{W}=\bigcup_{x\in F}\wb{W_x}$ and hence compact.
Now assume that $U\sub X_m$ is an open subset such that $K\sub U$ (e.g., $U=X_m$).
Since $K$ and $\wb{W}\setminus U$ are disjoint compact subsets of $\wb{W}$,
there is an open subset $P\sub \wb{W}$ such that $K\sub P$
and $\wb{P}\cap (\wb{W}\setminus U)=\emptyset$ (see Remark~\ref{betternormal}).
Then $V:=P\cap W$ is an open subset of~$W$ (and hence of~$X_n$) which contains~$K$, and
$\wb{V}\sub \wb{P}$ is compact and contained in $\wb{W}\setminus (\wb{W}\setminus U) \sub U$.

(b) Using~(a), we find
open subsets $Y,Z\sub X_n$
and $m\geq n$ such that $K\sub Y$, $L \sub Z$, and moreover
the closures $\wb{Y}$ and $\wb{Z}$ in $X_m$
are compact. Then $C:=\wb{Y}\cup\wb{Z}$ is a compact subset of~$X_m$.
Since $K$ and $L$ are disjoint compact subsets of~$C$,
there are open subsets $P,Q\sub C$ such that $K\sub P$,
$L\sub Q$, and such that the closures $\wb{P}$ and $\wb{Q}$ in~$C$
are disjoint (see Remark~\ref{betternormal}). Then $V:=Y\cap P$ and $W:=Z\cap Q$
have the desired properties.
\end{prf}
The following property is useful in the context of direct limits.
\begin{defn}\label{defncpregu}
Let $X_1\sub X_2\sub\cdots$ be a direct sequence of topological spaces $(X_n,\cO_n)$.
Assume that $X:=\bigcup_{n\in\N}X_n$ is 
endowed with a topology making all inclusion maps
$X_n\to X$ continuous.
We call $(X,(X_n)_{n\in\N})$ \emph{compactly regular}
if, for each compact subset~$K$ of~$X$,
there exists $n\in\N$ such that~$K$ is contained $X_n$
and $(X_n,\cO_n)$ induces the given compact topology on~$K$.
\end{defn}
We now record various basic properties of direct limit topologies.
Recall that a topological space $X$ is said to satisfy the separation axiom $T_1$
if the singleton $\{x\}$ is closed in~$X$ for each $x\in X$.
\begin{lem}\label{basicDL}
Let $X_1\sub X_2\sub\cdots$ be a direct sequence of topological spaces
and $X=\dl X_n$\vspace{-.5mm} be the direct limit topological space.
\begin{description}[(DD)]
\item[\rm(a)]
If $V_1\sub V_2\sub\cdots$ is an ascending sequence of open subsets
$V_n\sub X_n$,
then $V:=\bigcup_{n\in\N}V_n$ is open in $X=\dl X_n$.\vspace{-.5mm}
\item[\rm(b)]
In \emph{(a)}, the topology induced by~$X$ on~$V$ coincides with the direct limit topology on
$V=\dl V_n$.\vspace{-.5mm}
\item[\rm(c)]
If $M\sub X$ is open or closed,
then the topology induced by $(X,\cO_{\DL})$ on~$M$
coincides with the direct limit topology on
$M=\dl(M\cap X_n)$,\vspace{-1.3mm}
using the topology induced
by $(X_n,\cO_n)$ on $M\cap X_n$. 
\item[\rm(d)]
If each $X_n$ is $T_1$, then also $X$ is $T_1$.
\item[\rm(e)]
If each $X_n$ is $T_1$, then every compact subset
$K\sub X$ is contained in $X_n$ for some $n\in\N$.
\item[\rm(f)]
If $X_1\sub X_2\sub\cdots$ is a strict direct sequence,
then the topology induced by $(X,\cO_{\DL})$ on~$X_n$ coincides
with $\cO_n$, for each $n\in\N$.
\item[\rm(g)]
If $X_1\sub X_2\sub\cdots$ is a strict direct sequence of $T_1$-spaces,
then $(X,(X_n)_{n\in\N})$ is compactly regular.
\item[\rm(h)]
If each $X_n$ is Hausdorff and the direct sequence $X_1\sub X_2\sub\cdots$
is compactifying, then also $\dl X_n$\vspace{-.5mm}
is Hausdorff.
\end{description}
\end{lem}
\begin{prf}
(a) For each $n\in\N$, the preimage $j_n^{-1}(V)=X_n\cap \bigcup_{m\geq n}V_m=\bigcup_{m\geq n}j_{m,n}^{-1}(V_m)$
is open in $(X_n,\cO_n)$, whence~$V$ is open in $(X,\cO_{\DL})$.\vspace{.4mm}

(b) A subset $W\sub V$ is open in $\dl V_n$\vspace{-0.5mm} if and only if
$W\cap V_n$ is open in each~$V_n$, endowed with the topology induced by~$X_n$.
Since $V_n$ is open in~$X_n$, the latter holds if and only if $W\cap V_n$ is open in~$X_n$.
This property is satisfied if~$W$ is open in $(X,\cO_{\DL})$
(then $W\cap X_n$ is open in $X_n$, entailing that $W\cap V_n=(W\cap X_n)\cap V_n$
is open in~$X_n$, being an intersection of two open subsets).
If, conversely, $W\cap V_n$ is open in~$X_n$ for each~$n$, then $W=\bigcup_{n\in\N}(W\cap V_n)$
is open in $(X,\cO_{\DL})$, by~(a).

(c) Since $M$ is open (resp., closed) in $(X,\cO_{\DL})$,
the intersection $M\cap X_n$ is open (resp., closed) in $(X_n,\cO_n)$
for each $n\in\N$. A subset $N\sub M$ is open (resp., closed) in $\dl(M\cap X_n)$\vspace{-.5mm}
if and only if $N\cap X_n=N\cap M\cap X_n$ is open (resp., closed) in $M\cap X_n$.
The latter holds if and only if $N\cap X_n$ is open (resp., closed) in~$X_n$,
as $M\cap X_n$
is an open (resp., closed) subset of~$X_n$. But this condition defines
openness (resp., closedness) of~$N$ in~$X$.

(d) If $x\in X$, then $\{x\}\cap X_n$ is either $\{x\}$ or the empty set
and hence closed in~$X_n$, for each $n\in\N$. Hence $\{x\}$ is closed
in $(X,\cO_{\DL})$.

(e) If not, then we could find some $x_n\in K\setminus X_n$ for each $n\in\N$.
Then $D:=\{x_n\colon n\in\N\}$ would be an infinite set.
For each $n\in\N$, the intersection $D\cap X_n\sub\{x_1,\ldots, x_{n-1}\}$
is a finite subset of the $T_1$-space $X_n$ and hence
closed, entailing that $D$ is a closed subset of $(X,\cO_{\DL})$
and hence compact (as~$D\sub K$).
However, the $T_1$-space $X_n$ induces the discrete topology
on its finite subset $D\cap X_n$ (Exercise~\ref{exc-fininT1}).
Since $D=\dl(D\cap X_n)$\vspace{-.5mm} by (c),
we deduce that~$D$ is discrete (Exercise~\ref{exc-DL-disc}).
Being both compact and discrete, $D$ has to be finite
(contradiction).

(f) By Remark~\ref{firstremDL}(d), we may assume $n=1$.
The inclusion map $j_1\colon (X_1,\cO_1)\to (X,\cO_{\DL})$ is continuous.
To see that~$j_1$ is a topological embedding, let $W_1$ be an open subset
of $(X_1,\cO_1)$; we show that $W_1=j_1^{-1}(W)$ for an open subset~$W$
of $(X,\cO_{\DL})$. Recursively, for each $m\in\N$,
we find an open subset $W_{m+1}$ of $(X_{m+1},\cO_{m+1})$ such that $X_m\cap W_{m+1}=j_{m+1,m}^{-1}(W_{m+1})=W_m$,
as $j_{m+1,m}$ is a topological embedding.
By~(a), $W=\bigcup_{m\in\N}W_m$ is open in $(X,\cO_{\DL})$.
Since $W_m\cap X_1=W_1$ for all~$m$,
we have $j_1^{-1}(W)=W\cap X_1=W_1$.

(g) is immediate from (e) and (f).

(h) Given $x\not=y$ in~$X$,
we have $x,y\in X_{n_0}$ for some $n_0\in\N$.
Then $K_0:=\{x\}$ and $L_0:=\{y\}$ are compact subsets of $X_{n_0}$.
Using Lemma~\ref{compactify-la}(b) and induction,
we find integers $n_i>n_{i-1}$
and open subsets $V_i,W_i\sub X_{n_{i-1}}$ for $i\in\N$
such that $K_{i-1}\sub V_i$, $L_{i-1}\sub W_i$, and such that
the closures $K_i:=\wb{V_i}$ and $L_i:=\wb{W_i}$ in $X_{n_i}$
are compact and disjoint.
Then $V:=\bigcup_{i\in \N}V_i$ and $W:=\bigcup_{i\in\N}W_i$ are disjoint
open neighborhoods
of~$x$ and~$y$, respectively, in $\dl X_{n_i}=\dl X_n$\vspace{-.5mm}
(see (a) and Remark~\ref{firstremDL}(d)). Note that $V_i \subeq V_{i+1}$ 
and $W_i \subeq W_{i+1}$ hold by construction. 
\end{prf}
\begin{rem}
If $X_1\sub X_2\sub\cdots$ is compactifying, then also $V_1\sub V_2\sub\cdots$
is compactifying in Lemma~\ref{basicDL}(a).
In fact, if $x\in V_n$, then $x$ has an open neighborhood~$P$ in~$X_n$
whose closure $\wb{P}$ in $X_m$ is compact for some $m\geq n$.
Let~$Q$ be a compact neighborhood of $x$ in~$\wb{P}$ such that $Q\sub V_m$,
and $Q^0$ be its interior relative~$\wb{P}$. Then $W:=P\cap Q^0\cap V_n$
is an open neighborhood of~$x$ in~$V_n$ such that $\wb{W}\sub Q\sub V_m$
is compact.
\end{rem}
\begin{defn}\label{DLmaps}
Let $X_1\sub X_2\sub\cdots$ and $Y_1\sub Y_2\sub\cdots$
be direct sequences of topological spaces.
Give $X:=\bigcup_{n\in\N} X_n$ and  $Y:=\bigcup_{n\in\N} Y_n$
the direct limit topology and let
$j_n\colon Y_n\to Y$ be inclusion.
If $(f_n)_{n\in\N}$ is a sequence of continuous maps
$f_n\colon X_n\to Y_n$ such that $f_m|_{X_n}=f_n$ for all $n\leq m$,
then
\[
\dl f_n\colon X\to Y,\quad x\mto f_n(x)\;\;\; \mbox{if $x\in X_n$}\vspace{-.5mm}
\]
is a well-defined map.
Since $(\dl f_n)|_{X_n}=j_n\circ f_n$\vspace{-.5mm}
is continuous, also
\[
\dl f_n\colon \dl X_n\to \dl Y_n\vspace{-.5mm}
\]
is continuous,  by Lemma~\ref{ctsonsteps}.
\end{defn}
Let $X_1\sub X_2\sub\cdots$ and $Y_1\sub Y_2\sub\cdots$
be direct sequences of topological spaces,
$X:=\bigcup_{n\in\N} X_n$ and  $Y:=\bigcup_{n\in\N} Y_n$.
Then we have two natural topologies on $X\times Y$:
On the one hand, we can consider $X\times Y$ as the product
of $\dl X_n$\vspace{-.5mm} and $\dl Y_n$, and endow it with the product
topology~$\cO$.
On the other hand, we can consider $X\times Y$
as the union
\[
X\times Y=\bigcup_{n\in\N} X_n\times Y_n\vspace{-1mm}
\]
and give it the topology $\cO_{\DL}$ making it the direct limit
\[
X\times Y=\dl(X_n\times Y_n)\vspace{-.8mm}
\]
of the topological spaces $X_n\times Y_n$ (using the product topology on~$X_n\times Y_n$).
Then
$\cO\sub\cO_{\DL}$ (by the next lemma);
$\cO\not=\cO_{\DL}$ can occur
(Exercise~\ref{exc-product-not-k}).
\begin{lem}\label{DLcompaprod}
If $X_1\sub X_2\sub\cdots$ and $Y_1\sub Y_2\sub\cdots$
are direct sequences of topological spaces, then the identity map\vspace{-1mm}
\[
h\colon \dl(X_n\times Y_n)\to\left(\dl\, X_n\right)\times \left(\dl\, Y_n\right),\quad (x,y)\mto(x,y)\vspace{-1mm}
\]
is continuous. If both $X_1\sub X_2\sub\cdots$ and $Y_1\sub Y_2\sub\cdots$
are compactifying and all of the spaces $X_n$ and $Y_n$ are Hausdorff, then~$h$ is a homeomorphism.
\end{lem}
\begin{prf}
For $n\in\N$, let $\pr_{1,n}\colon X_n\times Y_n\to X_n$ and $\pr_{2,n}\colon X_n\times Y_n\to Y_n$
be the projection onto the first and second component.
Then $\pr_{1,n}$ and $\pr_{2,n}$ are continuous and hence also
$h=\big(\dl \pr_{1,n},\, \dl \pr_{2,n}\big)$.
Thus $\cO\sub\cO_{\DL}$.

Now assume the $X_n$ and $Y_n$ are Hausdorff and both direct sequences
are compactifying.
We show that $\cO_{\DL}\sub \cO$ (which completes the proof).
To this end, let $U\sub \dl(X_n\times Y_n)$\vspace{-.5mm} be an open subset.
It suffices to show that $U$ is a neighborhood of each $(x,y)\in U$
with respect to the product topology. We have $(x,y)\in X_{n_0}\times Y_{n_0}$
for some $n_0\in\N$. Then $K_0:=\{x\}$ and $L_0:=\{y\}$ are
compact subsets of $X_{n_0}$ and $Y_{n_0}$, respectively,
and $K_0\times L_0$
is contained in the open subset $U\cap (X_{n_0}\times Y_{n_0})$
of $X_{n_0}\times Y_{n_0}$.
By Lemma~\ref{compactify-la}(a),
we find $n_1>n_0$ and open subsets $V_1\sub X_{n_0}$ as well as $W_1\sub Y_{n_0}$
such that $K_0\sub V_1$, $L_0\sub W_1$ and moreover the closures~$\wb{V_1}\sub X_{n_1}$
and $\wb{W_1}\sub Y_{n_1}$ are compact.
The Wallace Lemma yields open subsets $P\sub X_{n_1}$ and $Q\sub Y_{n_1}$
such that
\[
\mbox{$K_0\sub P$, $\,L_0\sub Q$, \,and $\,P\times Q\sub U\cap(X_{n_1}\times Y_{n_1})$.}
\]
By Lemma~\ref{compactify-la}(a),
after shrinking~$V_1$ and~$W_1$, we may assume that $K_1:=\wb{V_1}\sub P$ and $L_1:=\wb{W_1}\sub Q$,
entailing that $K_1\times L_1\sub U\cap (X_{n_1}\times Y_{n_1})$.

Proceeding in this way, we obtain $n_j>n_{j-1}$ for $j\in \N$
and open subsets $V_j\sub X_{n_{j-1}}$ and $W_j\sub Y_{n_{j-1}}$
whose closures $K_j:=\wb{V_j}\sub X_{n_j}$ and $L_j:=\wb{W_j}\sub Y_{n_j}$
are compact, and such that $K_{j-1}\sub V_j$, $L_{j-1}\sub W_j$
and $K_j\times L_j\sub U\cap (X_{n_j}\times Y_{n_j})$.
Then $V:=\bigcup_{j\in\N}V_j$ and $W:=\bigcup_{j\in\N}W_j$
are open subsets of $\dl X_n$\vspace{-.5mm} and $\dl Y_n$, respectively (see Remark~\ref{firstremDL}(d)
and Lemma~\ref{basicDL}(a)), and $(x,y)\in V\times W\sub U$.
\end{prf}
\begin{defn}\label{dirseqTG}
(a) An ascending sequence $G_1\sub G_2\sub\cdots$ of topological groups $G_n$
is called a \emph{direct sequence of topological groups} if all inclusion maps
$G_n\to G_m$ (for $n\leq m$) are continuous group homomorphisms.
Then $G:=\bigcup_{n\in\N}G_n$ can be considered as a group:
Let $\mu_{G_n}\colon G_n\times G_n\to G_n$
be the group multiplication of~$G_n$.
The group multiplication $\mu_G\colon G\times G\to G$
is
\[
\mu_G(g,h):=\mu_{G_n}(g,h)
\]
for $g,h\in G$, where $n\in\N$ is chosen such that $g,h\in G_n$.
As all inclusion maps $G_n\to G_m$ (for $n\leq m$) are group homomorphisms,
the product is well defined, independent of the choice of~$n$.
For the preceding group structure on~$G$, each of the inclusion maps $G_n\to G$
is a group homomorphism, and it is uniquely determined
by this property.

(b) An
ascending sequence $E_1\sub E_2\sub\cdots$ of topological vector spaces
is called a \emph{direct sequence of topological vector spaces} if all inclusion maps
$E_n\to E_m$ (for $n\leq m$) are continuous and linear.
Then $E:=\bigcup_{n\in\N}E_n$ admits a unique vector space structure
making each~$E_n$ a vector subspace.
\end{defn}
We now investigate conditions ensuring that the direct limit topology $\cO_{\DL}$ on $G$
(or $E$) turns the latter into a topological group (resp., a topological vector space).
This is not always so (see Exercise~\ref{exc-DLnotGP});
the group multiplication need not be continuous as a map
$(G,\cO_{\DL})\times (G,\cO_{\DL})\to (G,\cO_{\DL})$.
\begin{prop}\label{DLofTG}
Let $G_1\sub G_2\sub \cdots$ be a direct sequence
of topological groups. Give $G:=\bigcup_{n\in\N}G_n$
the group multiplication turning each inclusion map $G_n\to G$
into a homomorphism, and endow $G$ with the direct limit topology.
Then the inversion map $\eta_G\colon G\to G$, $x\mto x^{-1}$
is continuous, and the left translations $\lambda_g^G\colon G\to G$, $x\mto gx$
and right translations $\rho^G_g\colon G\to G$,
$x\mto xg$ are continuous for each $g\in G$.
If the direct sequence $G_1\sub G_2\sub\cdots$ is compactifying,
then $G$ is a topological group
and the sets
\[
\bigcup_{n\in\N} U_1U_2\cdots U_n\vspace{-1mm}
\]
form a basis $\cB$ of open identity neighborhoods in~$G$,
when the $U_n$ range through a basis $\cB_n$ of open
identity neighborhoods in~$G_n$.
\end{prop}
\begin{prf}
Let $\mu_{G_n}\colon G_n\times G_n\to G_n$
be the continuous group multiplication of~$G_n$
and $\eta_{G_n}\colon G_n\to G_n$ be the continuous inversion map.
Then $\eta_G=\dl\eta_{G_n}$\vspace{-.6mm}
is continuous. Let $g\in G$. To see that $\lambda^G_g$ and $\rho^G_g$
are continuous, we may assume that $g\in G_1$ (see Remark~\ref{firstremDL}(d)).
Then
\[
\lambda^G_g=\dl\lambda^{G_n}_g\quad\mbox{and}\quad
\rho^G_g=\dl\rho^{G_n}_g\vspace{-.5mm}
\]
in terms of the continuous left (resp., right) translations in the topological groups~$G_n$,
and so $\lambda^G_g$ and $\rho^G_g$ are continuous.
If $G_1\sub G_2\sub\cdots$ is compactifying, then Lemma~\ref{basicDL}(h) shows that $G$ is Hausdorff
(recalling that we assume by definition that all topological groups are Hausdorff).
Moreover,
\[
G\times G=\dl (G_n\times G_n)\vspace{-.5mm}
\]
as a topological space (see Lemma~\ref{DLcompaprod}), entailing that $\mu_G=\dl\mu_{G_n}$\vspace{-.5mm}
is continuous. Thus $G$ is a topological group.

To prove the final assertion, let $W\in \cB$.
There are $V_n\in \cB_n$ for $n\in\N$ such that $W=\bigcup_{n\in\N}W_n$ with
\[
W_n:=V_1V_2\cdots V_n.
\]
Now $W_n$ is open in~$G_n$ (being a union of translates of the open subset~$V_n$).
Since $W_1\sub W_2\sub\cdots$, Lemma~\ref{basicDL}(a)
shows that $W=\bigcup_{n\in\N}W_n$ is open in~$G$.
Hence each $W\in\cB$ is an open identity neighborhood.
To see that $\cB$ is a basis, let $P_0\sub G$ be an open identity neighborhood.
Since $G$ is a topological group, we recursively find open identity neighborhoods
$P_n\sub G$ such that $P_nP_n\sub P_{n-1}$, for all $n\in\N$.
Then $P_1P_2\cdots P_n\sub P_0$ for all $n\in\N$.
Since $\cB_n$ is a basis of identity neighborhoods in~$G_n$,
for each $n\in\N$ we find $U_n\in \cB_n$ such that $U_n\sub P_n\cap G_n$.
Then
\[
U:=\bigcup_{n\in\N}U_1\ldots U_n\in\cB\vspace{-1.3mm}
\]
and $U\sub \bigcup_{n\in\N}P_1\cdots P_n\sub P_0$.
\end{prf}
\begin{cor}\label{cp-DL-TVS}
Given a direct sequence $E_1\sub E_2\sub\cdots$ of topological vector spaces,
make $E:=\bigcup_{n\in\N}E_n$ a vector space in the natural way.
Let $\cO_{\DL}$ be the direct limit topology on~$E$.
If $E_1\sub E_2\sub\cdots$ is compactifying, then $(E,\cO_{\DL})$
is a topological vector space. If, moreover, each~$E_n$ is locally convex, then
also $(E,\cO_{\DL})$ is locally convex.
\end{cor}
\begin{prf}
By Proposition~\ref{DLofTG}, the direct limit topology makes $(E,+)$ a topological group.
Let $m_{E_n}\colon\K\times E_n\to E_n$ and $m_E\colon \K\times E\to E$
be multiplication with scalars.
Since $\K\sub \K\sub\cdots$ and $E_1\sub E_2\sub\cdots$ are
compactifying,
\[
\K\times (E,\cO_{\DL})=\dl (\K\times E_n)\vspace{-1.3mm}
\]
as a topological space. Hence $m_E=\dl m_{E_n}$\vspace{-.5mm}
is continuous and thus~$E$ is a topological vector space.
If each~$E_n$ is locally convex, let $\cB_n$ be the set of all absolutely convex, open
$0$-neighborhoods in~$E_n$. If $U_n\in\cB_n$ for $n\in\N$,
then
\[
U:=\bigcup_{n\in\N}(U_1+\cdots+ U_n)\sub E\vspace{-1.3mm}
\]
is an ascending union of absolutely convex sets and hence absolutely convex.
As such sets~$U$ form a basis~$\cB$ of $0$-neighborhoods in~$E$ by Proposition~\ref{DLofTG},
we see that the topological vector space~$E$ is locally convex.
\end{prf}
\begin{ex}\label{exDLgp}
If $G_1\sub G_2\sub\cdots$ is a direct sequence of locally
compact topological groups,
then the direct limit topology makes $\bigcup_{n\in \N}G_n$
a (Hausdorff) topological group
(see Proposition~\ref{DLofTG} and
Remark~\ref{firstremDL}(c)).
\end{ex}
\begin{ex}\label{theKinfty}
Let $\K^{(\N)}\sub\K^\N$ be the vector subspace of all sequences
$(x_n)_{n\in\N}$ of scalars such that $x_n=0$
for all but finitely many~$n$. Identifying
$(x_1,\ldots, x_n)\in \K^n$ with $(x_1,\ldots, x_n,0,0,\ldots)\in\K^{(\N)}$,
we can write
\[
\K^{(\N)}=\bigcup_{n\in\N}\K^n\vspace{-1mm}
\]
and give $\K^{(\N)}$ the direct limit topology~$\cO_{\DL}$.
By Corollary~\ref{cp-DL-TVS}, $(\K^{(\N)},\cO_{\DL})$
is a locally convex topological vector space.
\end{ex}
\begin{rem}\label{rem-cof-subs}
(a) Let $(A,\leq)$ be a directed set and $(X_\alpha)_{\alpha\in A}$
be a family of topological spaces such that $X_\alpha\sub X_\beta$
and the inclusion map $j_{\beta,\alpha}\colon X_\alpha\to X_\beta$ is continuous
whenever $\alpha\leq\beta$. Then $(X_\alpha)_{\alpha\in A}$
is called a \emph{direct system of topological spaces}.
The final topology $\cO_{\DL}$ on $X:=\bigcup_{\alpha\in A}X_\alpha$
with respect to the inclusion maps $X_\alpha\to X$
is called the \emph{direct limit topology} on~$X$,
and we say that $(X,\cO_{\DL})=\dl X_\alpha$\vspace{-.8mm} as a topological space.

(b) An ascending sequence $\alpha_1\leq\alpha_2\leq\cdots$ in~$A$ is called
\emph{cofinal} if, for each $\beta\in A$, there exists $n\in \N$ such that $\beta\leq \alpha_n$.
Then $(X_{\alpha_n})_{n\in\N}$ is a direct sequence of topological spaces,
a so-called \emph{cofinal subsequence} of $(X_\alpha)_{\alpha\in A}$.
If $\cO$ is the direct limit topology on $X=\dl X_{\alpha_n}$\vspace{-.5mm},
then $\cO_{\DL}\sub \cO$ by definition. If $V\in\cO$ and $\beta\in A$, find $n\in\N$ such that
$\beta\leq\alpha_n$. Then
\[V\cap X_\beta=(V\cap X_{\alpha_n})\cap X_\beta
  =j_{\alpha_n,\beta}^{-1}(V\cap X_{\alpha_n})\] 
is open in $X_\beta$ and thus $V\in \cO_{\DL}$. Hence $\cO=\cO_{\DL}$ and thus
\[
\dl X_\alpha=\dl X_{\alpha_n}\vspace{-1mm}
\]
for each cofinal subsequence. All of the direct systems encountered in this book
will have cofinal subsequences and hence can be reduced to the case of direct sequences, as
treated in this appendix.

(c) More generally, one can consider direct systems $((X_{\alpha})_{\alpha\in A},(j_{\alpha,\beta})_{\alpha\geq \beta})$
and associated direct limits when $(A,\leq)$ is a directed set, each $X_\alpha$ a topological space and $j_{\alpha,\beta}\colon X_\beta\to X_\alpha$
a continuous map for $\alpha\geq \beta$ in~$A$ such that $j_{\alpha,\alpha}=\id_{X_\alpha}$
and $j_{\alpha, \beta}\circ j_{\beta,\gamma}=j_{\alpha,\gamma}$
for all $\alpha\geq\beta\geq\gamma$.
We can avoid them here, and the more technical notation that would go along with them.
\end{rem}
\begin{ex}\label{finest-lcx-vec}
Let $E$ be an infinite-dimensional vector space.
Then the set $A$ of all finite-dimensional vector subspaces~$F$ of~$E$
is directed under inclusion. Endow each $F\in A$ with the unique topology making it
a topological vector space (see Proposition~\ref{onlyyou}).
Then $(F)_{F\in A}$ is a direct system of locally convex topological vector spaces.
Let $\cO_{\DL}$ be the direct limit topology on $E=\bigcup_{F\in A} F$;
this topology is sometimes called the \emph{finite topology}
in the literature. Now assume that~$E$ has countable dimension.
Let $(v_n)_{n\in\N}$ be a basis for~$E$ and set
\[
F_n:=\K v_1\oplus\cdots\oplus \K v_n\sub E.
\]
Then $(F_n)_{n\in\N}$ is a cofinal subsequence in $(F)_{F\in A}$ and thus
$\dl F=\dl F_n$. Hence $(E,\cO_{\DL})$ is a locally convex topological
vector space (see Corollary~\ref{cp-DL-TVS}).
If $\cO$ is any vector topology on~$E$,
then $(E,\cO)$ induces the unique vector topology on each $F\in A$
and hence makes the inclusion map $F\to (E,\cO)$ continuous,
whence $\cO\sub \cO_{\DL}$.
Thus $\cO_{\DL}$ is the \emph{finest vector topology} on~$E$,
and also the \emph{finest locally convex vector topology}.\footnote{On a vector space of uncountable dimension,
the finite topology, the finest vector topology and the finest locally convex vector topology
all disagree, see \cite{KK63} and \cite{Bisg93}.}
Since
\[
f_n\colon \K^n\to F_n,\quad (x_1,\ldots, x_n)\mto x_1 v_1+\cdots+ x_n v_n
\]
is an isomorphism of topological vector spaces for each $n$ and $f_{n+1}|_{F_n}=f_n$,
\[
\dl f_n\colon \dl \K^n\to \dl F_n\vspace{-1mm}
\]
is an isomorphism of topological vector spaces (with inverse $\dl f_n^{-1}$\vspace{-1.2mm}) and thus
$E\cong \K^{(\N)}$ as a topological vector space.
\end{ex}
\begin{ex}
If $X$ is a Hausdorff topological space, then the set $\cK(X)$ of compact subsets of~$X$
is directed under inclusion (as $K_1\cup K_2$ is compact if so are $K_1$ and $K_2$).
Then $(K)_{K\in\cK(X)}$ is a direct system of topological spaces and~$X$ is a $k$-space if and only if
$X=\dl K$.\vspace{-1.2mm}
\end{ex}
\begin{defn}
We say that a topological space~$X$ is a $k_\omega$-space
if~$X$ is a $k$-space and hemicompact (as in \ref{hemic}).
Thus $X$ has sequence of compact subsets $K_1\sub K_2\sub\cdots$
which is cofinal in $\cK(X)$ (viz., each compact subset $K\sub X$ is contained in some~$K_n$).
In this context, $(K_n)_{n\in\N}$ is called a \break \emph{$k_\omega$-sequence}
for~$X$.
A topological space~$X$ is called a \emph{locally $k_\omega$-space}
if it is Hausdorff and every point $x\in X$ has an open neighborhood~$V$ in~$X$ which is a $k_\omega$-space
in the induced topology.
Then every neighborhood of $x$ in~$X$ contains such a neighborhood~$V$ (cf.\
Exercise~\ref{exc-k-omeg}).
\end{defn}
\begin{rem}
(a) Every $k_\omega$-space is locally $k_\omega$.
Every $\sigma$-compact locally compact space is a $k_\omega$-space;
every locally compact space is locally $k_\omega$.

(b) A topological space $X$ is $k_\omega$ if and only if $X=\dl K_n$\vspace{-.5mm}
for a direct sequence $K_1\sub K_2\cdots$ of compact topological spaces
(cf.\ Lemma~\ref{basicDL}(g)).

(c) Every locally $k_\omega$-space
is a $k$-space
(cf.\ Exercise~\ref{cover-k}).
\end{rem}
\begin{prop}\label{henceSilvacoreg}
If a direct sequence $X_1\sub X_2\sub\cdots$ of Hausdorff topological spaces is compactifying, then $X:=\dl X_n$\vspace{-1mm}
is locally $k_\omega$ and $(X,(X_n)_{n\in\N})$
is compactly regular.
\end{prop}
\begin{prf}
If $x\in X_{n_1}$, set $K_1:=\{x\}$ and find integers $n_1<n_2<\cdots$ and open subsets $V_j$
of $X_{n_j}$ with $K_j\sub V_j$ such that $K_{j+1}:=\wb{V_j}$ is compact in $X_{n_{j+1}}$ 
(Lemma~\ref{compactify-la}).
Then $V=\dl V_j$ is open in $\dl X_n=\dl X_{n_j}$,\vspace{-.5mm} by Lemma~\ref{basicDL}(a).
Moreover, $V=\dl K_j$\vspace{-.7mm}
(see Lemma~\ref{basicDL}(b) and Exercise~\ref{exc-easycofin}), whence $V$ is a
$k_\omega$-space.\vspace{.3mm}

If $\emptyset\not=K\sub X$ is compact, then each $x\in K$ has an open neighborhood~$V_x$ in~$X$
which is a $k_\omega$-space and admits a $k_\omega$-sequence $(K_{x,j})_{j\in\N}$
such that each $K_{x,j}$ is a compact subset of~$X_n$ for some $n\in\N$ (as just shown).
Now~$x$ has a compact neighborhood~$L_x$ in~$K$ such that $L_x\sub V_x$. By compactness of~$K$,
\[
K=\bigcup_{x\in F} L_x\vspace{-.8mm}
\]
for some finite subset $F\sub K$. Since $(V_x,(K_{x,j})_{j\in\N})$
is compactly regular by Lemma~\ref{basicDL}(g), we find $j_x\in\N$ for $x\in F$ such that~$L_x$
is compact in $K_{x,j_x}$ and hence also in $X_{m_x}$ for some $m_x\in \N$.
Thus~$K$ is contained in~$X_m$ for $m:=\max\{m_x\colon x\in F\}$ and since~$X_m$ induces
the given compact topology on~$L_x$ for each $x\in F$, also $K=\bigcup_{x\in F}L_x$
is compact in~$X_m$ in the induced topology,
which (being coarser) coincides with the given topology on~$K$.
\end{prf}
\begin{rem}\label{prod-of-k-omeg}
(a) As shown in \cite{GHK10}, $\dl X_n$\vspace{-.5mm}
is $k_\omega$ (resp., locally $k_\omega$)
for each direct sequence $X_1\sub X_2\sub\cdots$
of $k_\omega$-spaces (resp., locally $k_\omega$-spaces).
In particular, the direct limit is Hausdorff.

(b) Let $X_1\sub X_2\sub\cdots$ and $Y_1\sub Y_2\sub\cdots$ be direct sequences of \break 
$k_\omega$-spaces (resp., locally $k_\omega$-spaces)
and $X:=\dl X_n$, $Y:=\dl Y_n$.\vspace{-.8mm}
Then $X\times Y$ is $k_\omega$ (resp., locally $k_\omega$) and
\[
X\times Y=\dl (X_n\times Y_n)\vspace{-.8mm}
\]
as a topological space (see \cite{GHK10}).
\end{rem}
\begin{lem}\label{la-rel-cp-nbds}
Consider an ascending sequence $X_1\sub X_2\sub\cdots$
of locally compact topological spaces such that all inclusion maps $X_n\to X_{n+1}$
are continuous. Endow $X:=\bigcup_{n\in\N}X_n$ with the direct limit topology.
If $x\in X_1$ and $U\sub X$ is an open $x$-neighborhood,
then there exists a sequence $(V_n)_{n\in \N}$ open $x$-neighborhoods $V_n\sub X_n$
with compact closure $\wb{V_n}\sub X_n \cap U$
such that $\wb{V_n}\sub V_{n+1}$ for each $n\in \N$.
\end{lem}
\begin{prf}
The open $x$-neighborhood $X_1\cap U$ in $X_1$
contains a compact $x$-neighborhood~$K_1$,
by local compactness of $X_1$; we let $V_1$ be the interior of~$K_1$ in~$X_1$.
If relatively compact open $x$-neighborhoods $V_j$ in $X_j\cap U$ have been found for $j\in\{1,\ldots, n\}$
such that $\wb{V_j}\sub V_{j+1}$ for all $j<n$,
then the compact subset $\wb{V_n}$ of $X_n\cap U$
is also compact in the open subset $X_{n+1}\cap U$ of $X_{n+1}$.
By Remark~\ref{trivi-lcp},
there is a compact subset $K_{n+1}\sub X_{n+1}\cap U$ whose interior
$V_{n+1}:=K_{n+1}^0$ in~$X_{n+1}$
contains $\wb{V_n}$.
\end{prf}
\subsection*{Exercises for Section~\ref{app-basic-DL}}

\begin{small}

\begin{exer}\label{exc-fininT1}
Show that every finite subset of a $T_1$-space~$X$
is closed in~$X$.
\end{exer}

\begin{exer}\label{exc-DL-disc}
Let $X_1\sub X_2\sub\cdots$ be a direct sequence
of topological spaces. Show: If each $X_n$ is discrete,
then also $X=\dl X_n$\vspace{-.5mm} is discrete.
\end{exer}

\begin{exer}\label{exc-perm-k}
Show that if $X_1\sub X_2\sub\cdots$ is a direct sequence of $k$-spaces,
then also $\dl X_n$\vspace{-.5mm} is a $k$-space. Likewise for $k_\R$-spaces.
\end{exer}

\begin{exer}\label{exc-product-not-k}
Consider the bilinear map\vspace{-1mm}
\[
\beta\colon \R^\N\times\R^{(\N)}\to\R,\quad
\beta(x,y):=\sum_{n=1}^\infty x_ny_n\vspace{-1mm}
\]
for $x=(x_n)_{n\in\N}\in\R^\N$ and $y=(y_n)_{n\in\N}\in\R^{(\N)}$.
We endow $\R^\N$ with the product topology (which is metrizable by~\ref{furex}(a))
and hence makes $\R^\N$ a $k$-space) and give $\R^{(\N)}$
the topology making it the direct limit $\R^{(\N)}=\dl \R^n$\vspace{-.5mm}
of its vector subspaces $\R^n\times \{0\}$, which we identify with~$\R^n$.
By Exercise~\ref{exc-perm-k}, $\R^{(\N)}$ is a $k$-space.
\begin{description}[(D)]
\item[(a)]
Show that $\beta|_{\R^\N\times\R^n}$ is continuous for each $n\in\N$.
\item[(b)]
Show that every compact subset $K$ of $\R^\N\times \R^{(\N)}$
is contained in $\R^\N\times\R^n$ for some $n\in \N$.
Hence $\beta|_K$ is continuous.
\item[(c)]
Show that $\beta$ is not continuous at $(0,0)$.
Hence $\R^\N\times\R^{(\N)}$ is not a $k$-space (nor a $k_\R$-space)
although $\R^\N$ and $\R^{(\N)}$ are $k$-spaces.
Moreover, the product topology on $\R^\N\times \R^{(\N)}$
has to be properly coarser than the direct limit topology on $\dl(\R^\N\times \R^n)$\vspace{-.5mm}
(as the latter makes~$\beta$ continuous).
\end{description}
\end{exer}

\begin{exer}\label{exc-when-DL-GP}
Let $G_1\sub G_2\sub\cdots$ and $H_1\sub H_2\sub\cdots$ be direct sequences of topological groups;
endow the groups $G:=\bigcup_{n\in\N}G_n$ and $H:=\bigcup_{n\in\N}H_n$
with the direct limit topology.
\begin{description}[(D)]
\item[(a)]
Using that the maps $j_n\colon G_n\to G_n\times H_n$, $x\mto(x,\be)$ are continuous,
show that
\[
j_G\colon G\to \dl(G_n\times H_n),\quad x\mto (x,\be)
\]
is continuous. Likewise, $j_H\colon H\to\dl (G_n\times H_n)$, $y\mto(\be,y)$ is continuous.
\item[(b)]
Show:
If $\dl(G_n\times H_n)$\vspace{-.5mm} is a topological group,
then it coincides with $G\times H$, endowed with the product topology.\\[1mm]
[Hint: The map $G\times H\to \dl(G_n\times H_n)$, $(x,y)\mto(x,y)$\vspace{-.8mm}
is continuous since $(x,y)=j_G(x)j_H(y)$.]
\end{description}
\end{exer}

\begin{exer}\label{exc-DLnotGP}
Using Exercise~\ref{exc-product-not-k} and Exercise~\ref{exc-when-DL-GP}, deduce
that the direct limit topology on $\dl(\R^\N\times\R^n)$
does not make the latter a topological group.
\end{exer}

\begin{exer}\label{exc-easycofin}
Let $V_1\sub V_2\sub\cdots$ be a direct sequence of topological spaces
and $L_n\sub V_n$ be a subset for $n\in\N$ such that $V_n\sub L_{n+1}$.
Give $L_n$ the topology induced by~$V_n$. Show that $\dl L_n=\dl V_n$\vspace{-.2mm}
as a topological space.\\[1mm]
[Both $(L_n)_{n\in\N}$ and $(V_n)_{n\in\N}$ are subsequences of $L_1\sub V_1\sub L_2\sub V_2\sub\cdots\,$]
\end{exer}

\begin{exer}\label{exc-k-omeg}
Let $X$ be a $k_\omega$-space.
\begin{description}[(D)]
\item[(a)]
Show that every closed subset $A\sub X$ is a $k_\omega$-space
in the induced topology.
\item[(b)]
Let $x\in X$. Show that every open neighborhood $U$ of $x$ in $X$ contains an open neighborhood $V$ of~$x$
which is a $k_\omega$-space.\\[1mm]
[Hint: Let $(K_n)_{n\in\N}$ be a $k_\omega$-sequence for~$X$. Set $L_0:=\{x\}$ and,
for $n\in \N$, choose a compact subset $L_n\sub K_n\cap U$ which has $L_{n-1}$ in its interior $L_n^0$ (relative $K_n$).
Now take $V:=\dl L_n=\dl L_n^0$.]
\end{description}
\end{exer}

\begin{exer}\label{cover-k}
Let $X$ be a Hausdorff topological space.
Show: If $X$ has an open cover $(U_j)_{j\in J}$
such that each $X_j$ is a $k$-space,
then $X$ is a $k$-space.
\end{exer}

\end{small}
\section{Notes and comments on Appendix~\ref{appA}}\label{notes-appA}
Much of the material in this appendix is standard,
whence we refrain from references for most results;
our main references for general topology are~\cite{Bou88}, \cite{Ek89}, also~\cite{Sch75}
and~\cite{Kel75}.
For our purposes, it suffices to consider nets whose directed index sets
are ordered (not only quasi-ordered).
For background concerning $k$-spaces, we refer to~\cite{Kel75}
and, with a view towards algebraic topology,~\cite{Ste67}.
Exercise~\ref{exc-product-not-k} draws on \cite[Remark~17.9]{Bcz91}.
For further information on $k_\omega$-spaces,
see~\cite{Mic68}, \cite{Mr56}, \cite{FS77}, and~\cite{GHK10};
the latter work also discusses spaces which are locally~$k_\omega$.
Background concering $k_\R$-spaces can be found in \cite{Hu71}
and~\cite{Nb70}; in contrast to the latter,
we do not impose complete reguarity on such spaces
(as this property does not play a role for our application
to exponential laws. Of course, all locally convex
spaces are completely regular, as well as subsets thereof).
The notation
$\lfloor K, U\rfloor$ for subbasic sets of the compact-open topology
has been adopted from~\cite{Str06}.
As to basic properties of direct limits of topological spaces,
cf.\ \cite{NRW91}, \cite{NRW93}, \cite{Gl03a}, and~\cite{Gl05c}
for Lemma~\ref{basicDL} (a)--(d) and~(f),
\cite{Han71} for Lemma~\ref{basicDL} (e)~and~(g).
Direct limits of topological groups and those of the underlying topological spaces
were compared in~\cite{TSH98}, \cite{HST01}, and~\cite{Ym98},
also~\cite{Gl03a}, \cite{Gl07a}, \cite{Gl11}, \cite{Gl13a}, and (with a focus on $k_\omega$-spaces)
in~\cite{GHK10}. The topic of compactifying direct sequences will be taken up
in Appendix~\ref{sec-appDLvec}, when we discuss so-called Silva spaces,
which are direct limits of Banach spaces whose bonding maps are compact operators.
Also the topic of compact regularity will be resumed in the context
of locally convex direct limits. For applications of compact regularity
in infinite-dimensional Lie theory, see, e.g.,~\cite{Da14}, \cite{Gl11}, and~\cite{Gl12b};
it can be helpful in proofs of Lie-theoretic regularity for ascending unions of Lie groups.
Moreover, it can simplify the calculation of homotopy groups of such, see~\cite{Gl08c}.

For further information concerning paracompactness,
see, e.g., \cite{Ek89}.
Theorem~\ref{godement-globalize} and Lemma~\ref{pre-godement}
vary \cite[p.\,150]{Gd58}
and the discussion of tubular neighborhoods in \cite{La99}.

\chapter[Basic facts concerning locally convex spaces]{Basic facts concerning
locally convex spaces}\label{chaplcx}
\chaptermark{Basics of locally convex spaces\hfill \copyright{}
H. Gl\"{o}ckner and K.-H. Neeb}
In this appendix,
we recall the definition of
locally convex topological vector spaces,
describe various basic constructions,
and prove all results required for our purposes.
Throughout the following,
$\K\in \{\R,\C\}$.
All vector spaces will be $\K$-vector
spaces and all linear maps
are $\K$-linear, unless the contrary
is stated.
We assume that the reader is familiar
with some basic facts concerning
point set topology and
normed spaces (like the uniform boundedness principle).
\section{Basic definitions and first properties}\label{decdefnlcx}
\begin{defn}\label{defnlcx}
A {\em topological vector space\/}
is a vector space~$E$, equipped
with a topology
turning both the addition map
\[
E\times E\to E\,,\quad (x,y)\mto x+y
\]
and scalar multiplication
\[
\K\times E\to E\,,\quad (z,x)\mto zx
\]
into continuous mappings.\footnote{Here $\K$ is equipped
with its usual topology, and the products
are equipped with the product topology.}
Such a topology is also called a {\em
vector topology}.
A topological vector space
is called {\em locally convex\/}
if every $0$-neighborhood
contains a convex $0$-neighborhood.\footnote{A subset $U$ of a
vector space is called {\em
convex\/} if it contains with any two
points $x,y\in U$ also the line segment joining them,
viz.\ $tx+(1-t)y\in U$ for all $t\in [0,1]$.}
We use the term \emph{locally convex space}
as a shorthand for ``locally convex topological
vector space.''
We shall always assume that topological
vector spaces are Hausdorff, unless the contrary is stated.
\end{defn}
Given subsets $U$ and $V$ of a vector space~$E$
and $t\in \K$,
we write
\[ U+V:=\{x+y\colon x\in U, y\in V\} \quad \mbox{ and } \quad
  tU:=\{tx\colon x\in U\}.\]
We set
$\bD_r := \{ z\in \K\colon |z|\leq r\}$
for $r>0$ and $\bD:=\bD_1$.
\begin{lem}\label{verybase}
For each topological vector space~$E$,
the following holds:
\begin{description}
\item[\rm (a)]
For each $x\in E$, the translation
$\tau_x\colon E\to E$, $\tau_x(y):=x+y$
is a homeomorphism.
\item[\rm (b)]
For each $t\in \K^\times$,
the map $h_t\colon E\to E$, $h_t(x):=tx$
is a homeomorphism.
\item[\rm (c)]
Every $0$-neighborhood $U\sub E$
contains a $0$-neighborhood~$V$
which is ``balanced,'' i.e.,
$tV\sub V$ for each $t\in \K$ such that
$|t|\leq 1$. Then $V={-V}$ in particular.
\item[\rm (d)]
For each $0$-neighborhood~$U\sub E$,
there exists a $0$-neighborhood
$V\sub E$ such that
$V+V\sub U$.
\item[\rm (e)]
Each $0$-neighborhood
$U\sub E$ is ``absorbing,''
i.e., $E=\bigcup_{r>0} r U$.
\item[\rm (f)]
If $\cB$ is a basis of $0$-neighborhoods
in~$E$,
then for each $x\in E$ the set
$\{x+U\colon U\in \cB\}$ is a basis
of neighborhoods of~$x$ in~$E$.
\item[\rm (g)]
Each $0$-neighborhood contains
a closed $0$-neighborhood.
\item[\rm (h)]
If $K\sub E$ is compact and
$U\sub E$ an open subset such
that $K \sub U$,
then there exists
an open $0$-neighborhood $W\sub E$ such that
$K+W\sub U$.
\item[\rm(i)]
If $M\sub E$ is a subset
and $\cB$ a basis of $0$-neighborhoods
in~$E$, then the closure of~$M$ is given by
$\wb{M}=\bigcap_{U\in \cB}(M+U)$.
\end{description}
\end{lem}
\begin{prf}
(a) The addition map
$\alpha\!: E\times E\to E$, $\alpha(u,v):=u+v$ and
the map $i_x\!: E\to E\times E$,
$i_x(y):=(x,y)$ being continuous,
so is $\tau_x=\alpha\circ i_x$.
Clearly $\tau_x$ is invertible,
with inverse $\tau_{-x}$
which is continuous. Hence $\tau_x$ is a homeomorphism.

(b) Scalar multiplication $\mu\colon \K\times E\to E$
is continuous and hence so are its partial maps
$h_t=\mu(t,\cdot)$
and $(h_t)^{-1}=h_{t^{-1}}$.

(c) By continuity of the scalar multiplication~$\mu$,
the preimage $\mu^{-1}(U)$ is open in $\K\times E$,
entailing that there exists $\ve>0$ and
a $0$-neighborhood $W\sub E$ such that
$\bD_\ve \times W \sub \mu^{-1}(U)$
and thus $V:= \bD_\ve W=\mu(\bD_\ve \times W)\sub U$.
Then $V\sub U$ and $\bD V=V$
because $\bD \bD_\ve=\bD_\ve$.

(d) Since addition $\alpha\colon E\times E\to E$
is continuous and $\alpha(0,0)=0$,
the preimage $\alpha^{-1}(U)$ is a $(0,0)$-neighborhood.
As $E\times E$ is equipped with the
product topology, there exist $0$-neighborhoods
$V_1,V_2\sub E$ such that $V_1\times V_2\sub \alpha^{-1}(U)$.
Then
$V :=V_1\cap V_2$ has the desired property.

(e) Given $x\in E$, the map $f\colon \K\to E$,
$f(t):=tx$ is continuous, and
$f(0)=0$. Hence $f^{-1}(U)$ is a $0$-neighborhood in~$\K$,
and so there
exists $n\in \N$ such that $\frac{1}{n}\in f^{-1}(U)$.
Then $\frac{1}{n}x\in U$ and thus $x\in nU$.

(f) follows
from the fact that $\tau_x$ is a homeomorphism
(by (a)).

(g) Given a $0$-neighborhood $U\sub E$,
let $V\sub E$ be a $0$-neighborhood
such that $V-V\sub U$ (see (d) and (c)).
Given $w\in \wb{V}$, by\,(f) the set $w+V$ is a neighborhood
of~$w$ in~$E$.
Hence, there exists $v_1\in V$ such that
$v_1\in w+V$, whence in turn there
is $v_2\in V$ such that $v_1=w+v_2$.
Then $w=v_1-v_2\in U$.
Thus $\wb{V} \sub U$.

(h)
By (f), for each $x\in K$ there exists
a $0$-neighborhood $V_x\sub V$ such that
$x+V_x\sub U$. By (d), there is an open $0$-neighborhood
$W_x\sub E$ such that $W_x+W_x\sub V_x$.
Then $(x+W_x)_{x\in K}$
is an open cover of~$K$.
By compactness of~$K$, there is a finite subset
$F\sub K$ such that $K\sub \bigcup_{x\in F}(x +W_x)$.
Then $W:=\bigcap_{x\in F}W_x$ is an open $0$-neighborhood.
Given $y\in K$, there exists $x\in F$ such that
$y\in x+W_x$. Then $y+W\sub x+W_x+W\sub x+V_x\sub U$.
As $y$ was arbitrary, $K+W=\bigcup_{y\in K}(y+W)\sub U$
follows.

(i) Because the map $E\to E$, $x\mto -x$ is a homeomorphism,
also the sets $-U$ form a basis of $0$-neighborhoods,
for $U\in \cB$.
Let $x\in E$. Then $x\in \wb{M}$ if and only if $V\cap M\not=\emptyset$
for all~$V$ in a basis of neighborhoods of~$x$.
Thus $x\in \wb{M}$ if and only if
\begin{equation}\label{soonchg}
(x-U)\cap M\not=\emptyset
\end{equation}
for all $U\in \cB$ (using (f)).
Since (\ref{soonchg})
is equivalent to $x\in M+U$, the assertion follows.
\end{prf}
\begin{lem}\label{baseconvex}
Let $E$ be a topological vector space
and $U\sub E$ be convex.
Then the following holds:
\begin{description}
\item[\rm (a)]
The sets $tU$ and $U+x$ are convex,
for each $t\in \K$ and $x\in E$.
\item[\rm (b)]
The interior $U^0$ and the closure
$\wb{U}$ of~$U$ are convex.
Moreover, we have $tx+(1-t)U^0\sub U^0$ for all $x\in U$
and $t\in [0,1[$.
\item[\rm (c)]
If $U^0\not=\emptyset$, then $U^0$ is dense
in~$U$.
\item[\rm(d)]
If $U$ is a $0$-neighborhood, $x\in U$
and there exists $r>1$ such that $rx\in U$, then $x\in U^0$.
\end{description}
\end{lem}
\begin{prf}
(a) Is trivial.

(b) Given $t\in [0,1]$,
consider $\phi\colon E\times E\to E$,
$\phi(x,y):=tx+(1-t)y$.
Since $U$ is convex, we have $\phi(U\times U)\sub U$
and hence $\phi(\wb{U}\times \wb{U})\sub \wb{U}$,
by continuity of~$\phi$.
Hence $\wb{U}$ is convex.
Given $x\in U$
and $t\in [0,1[$,
the set $tx+(1-t)U^0\sub U$
open in~$E$ by Lemma~\ref{verybase}\,(a) and\,(b),
and thus $tx+(1-t)U^0\sub U^0$.
Hence $tU^0+(1-t)U^0\sub U^0$
for each $t\in \;]0,1[$,
and thus $U^0$ is convex.

(c) Given $x\in U$ and
$y\in U^0$, we have
$u_t:=tx+(1-t)y\in tx+(1-t)U^0\sub U^0$
for each $t\in [0,1[$ as just shown,
where $u_t\to x$ as $t\to 1$. Thus $x\in \wb{U^0}$.

(d) By (b), we have $x=\frac{1}{r} rx\in \frac{1}{r}(rx)+(1-\frac{1}{r})U^0\sub U^0$.
\end{prf}
\begin{rem}\label{remlikewise}
Likewise,
$tU$, $x+U$ and $\wb{U}$
are balanced if so is $U\sub E$.
If $U\sub E$ is a balanced
$0$-neighborhood, then so is~$U^0$.
\end{rem}
Recall that a subset $U$
of a vector space~$E$ is called {\em absolutely convex\/}
if $U$ is both balanced and convex.
\begin{lem}\label{lemabscon}
Let $U$ be a subset of a vector space~$E$.
Then
\begin{description}
\item[\rm (a)]
There exists a smallest
convex subset $\cnv(U)\sub E$
containing~$U$.
\item[\rm (b)]
An element $x\in E$ belongs to~$\cnv(U)$
if and only if there exists
$k\in \N$, elements $x_1,\ldots, x_k\in U$
and $t_1,\ldots, t_k\in \;] 0,\infty[$
such that $\sum_{j=1}^k t_j=1$
and $x=\sum_{j=1}^k t_j x_j$.
\item[\rm (c)]
$\abscnv(U):=\cnv(\bD U)$
is the smallest absolutely convex subset
of~$E$ which contains~$U$.
\end{description}
\end{lem}
\begin{prf}
Let $C$ be the set of all
elements $x=\sum_{j=1}^k t_j x_j$ as described in~(b).
We claim that if $D \sub E$ is convex
and $U\sub D$, then $C\sub D$.
To see this, we show by induction on~$k$
that each element of~$C$ of the form
$x=\sum_{j=1}^k t_j x_j$
is in~$D$.
The case $k=1$ is trivial.
If the assertion holds
for~$k$ and $x=\sum_{j=1}^{k+1}t_j x_j\in C$,
then $x= (1-t_{k+1}) y+t_{k+1}x_{k+1}$
where $y:=\sum_{j=1}^k \frac{t_j}{1-t_{k+1}}x_j\in D$
by induction and $x_{k+1}\in U\sub D$
and thus $x\in D$, the set $D$ being convex.
Thus~$C\sub D$.

To see that $C$ is convex,
note that we may assume that
the elements $x_j$ are pairwise
distinct in the formula $x=\sum_{j=1}^k t_jx_j$.
Furthermore, clearly we get
the same set~$C$
if we allow $t_j=0$.
Hence, given $x,y\in C$, we may assume
that $x=\sum_{j=1}^k t_j x_j$
and $y=\sum_{j=1}^k s_j x_j$,
using the same
elements $x_1,\ldots, x_k\in C$
and certain $s_j,t_j\in [0,\infty[$
such that $\sum_{j=1}^k t_j=\sum_{j=1}^k s_j=1$.
Given $t\in [0,1]$, we then have
$tx+(1-t)y=\sum_{j=1}^k r_j x_j$
where $r_j:=t t_j+(1-t)s_j\geq 0$ and
$\sum_{j=1}^k r_j=t\sum_{j=1}^kt_j+(1-t)\sum_{j=1}^k s_j=t+(1-t)=1$
and thus $tx+(1-t)y\in C$. Thus $C$ is convex
and hence $C=\cnv(U)$.

Since every absolutely convex set $D$ containing~$U$
also contains $\bD U$,
we must have $\cnv(\bD U)\sub D$.
On the other hand, it is clear
from the formula from~(b)
that $\conv(\bD U)$ is balanced
and thus absolutely convex.
Hence
$\cnv(\bD U)$ is
the smallest absolutely convex subset of~$E$
containing~$U$.
\end{prf}
\begin{lem}\label{basisnice}
If $E$ is a locally convex space,
then every
$0$-neighborhood $U\sub E$ contains an open,
absolutely convex $0$-neighborhood
and a closed,
absolutely
convex $0$-neighborhood.
Furthermore, each neighborhood
of a point $x\in E$ contains an
open $x$-neighborhood and
a closed, convex $x$-neighborhood.
\end{lem}
\begin{prf}
By Lemma~\ref{verybase}(g),
there exists a closed $0$-neighborhood
$V\sub E$ such that $V\sub U$.
By definition of a locally convex space, there
exists a convex $0$-neighborhood $C\sub E$
such that $C\sub V$.
By Lemma~\ref{verybase}(c), there exists
a balanced $0$-neighborhood $B\sub C$.
Then $Q:=\abscnv(B)=\cnv(\bD B)=\cnv(B)\sub C$
is absolutely convex and a $0$-neighborhood.
By Lemma~\ref{baseconvex}(b)
and Remark~\ref{remlikewise},
the sets $Q^0$ and $\wb{Q}\sub V\sub U$
are absolutely convex $0$-neighborhoods
of~$E$ contained in~$U$ which are
open and closed, respectively.
Hence the set $\cB$ of all
open (resp., closed) absolutely convex $0$-neighborhoods
is a basis of $0$-neighborhoods.
Given $x\in E$, the set $\{x+W \colon W\in \cB\}$
is a basis of $x$-neighborhoods,
by Lemma~\ref{verybase}(f).
Here $x+W$ is open (resp., closed)
since translation by~$x$ is a homeomorphism
(Lemma~\ref{verybase}(a)).
Also, $W$ being convex, so is
$x+W$ (Lemma~\ref{baseconvex}(a)).
\end{prf}
\begin{lem}\label{TVSHausdorff}
Let $E$ be a topological vector space
which we do not assume Hausdorff. If $\{0\}$
is closed in~$E$, then $E$ is Hausdorff.
\end{lem}
\begin{prf}
Translations being homeomorphisms, $\{x-y\}$ is closed
whenever $x\not=y$ are distinct elements of~$E$
(cf.\ Lemma~\ref{verybase}\,(a)).
Hence $U:=E\setminus\{x-y\}$ is a $0$-neighborhood.
Using Lemma~\ref{verybase}~(c) and~(d), we find a $0$-neighborhood $V\sub E$
such that $V-V\sub U$.
Then $x+V$ and $y+V$ are disjoint neighborhoods
of~$x$ and~$y$, respectively.
In fact, if $x+v=y+w$ with $v,w\in V$,
we would obtain the contradiction $x-y=w-v\in V-V\sub U$.
\end{prf}
\subsection*{Exercises for Section~\ref{decdefnlcx}.}

\begin{small}
\begin{exer}\label{ac-kernel}
Let $E$ be a topological vector space
and $U\sub E$ be a convex, open $0$-neighborhood.
\begin{description}[(D)]
\item[(a)]
Show that the set $V:=\{x\in U\colon \bD x\sub U\}$
is absolutely convex. Using the Wallace Lemma (Lemma~\ref{Wallla}), show that~$V$ is open.
Deduce that~$V$ is the largest absolutely convex $0$-neighborhood contained in~$U$.
\item[(b)]
Let $\alpha\colon E\to E_1$ be a linear map to a topological vector space~$E_1$
and $U_1\sub E_1$ be an open, convex $0$-neighborhood
such that $\alpha(U)\sub U_1$. Show that $\alpha(V)\sub V_1$ if we define
$V_1:=\{y\in E_1\colon\bD y\sub U_1\}$.
\end{description}
\end{exer}
\end{small}
\section{Continuity and openness of
linear maps}\label{secconopen}
\begin{prop}\label{ctsopen}
Let $f\colon E\to F$ be a linear
map between topological vector spaces.
Then the following holds:
\begin{description}
\item[\rm (a)]
$f$ is continuous if and only if~$f$
is continuous at~$0$.
\item[\rm (b)]
$f$ is open if and only if
$f(U)$ is a $0$-neighborhood in~$F$,
for each $0$-neighborhood
$U\sub E$.
\end{description}
\end{prop}
\begin{prf}
In each case, the condition
is clearly necessary.
Let us prove sufficiency.

(a) Assume that $f$ is continuous at~$0$;
given $x\in E$, we show that $f$ is continuous
at~$x$. Using the translations
$\tau_x$ and $\tau_{f(x)}$
on~$E$ and $F$, respectively,
we have
$f\circ \tau_x = \tau_{f(x)}\circ f$
as $f(x+y)=f(x)+f(y)$.
Hence
\begin{equation}\label{move}
f\; = \; \tau_{f(x)}\circ f
\circ\tau_x^{-1}\,.  
\end{equation}
Now $\tau_x^{-1}$ being continuous
with $\tau_x^{-1}(x)=0$,
the right hand side of (\ref{move})
is continuous at~$x$ and hence so is the left hand
side, $f$.

(b) Let $U\sub E$ be an open subset.
For each $x\in U$, the set $U-x$ is a $0$-neighborhood
in~$E$ (see Lemma~\ref{verybase}(a))
and hence $f(U-x)=f(U)-f(x)$
is a $0$-neighborhood in~$F$, by the hypothesis.
As a consequence, $f(U)=f(U-x)+f(x)$
is a neighborhood of $f(x)$ in~$F$,
and thus $f(x)\in f(U)^0$.
Hence $f(U)=f(U)^0$, whence $f(U)$ is open.
\end{prf}
Of course, it suffices to consider
$U$ in a basis of $0$-neighborhoods
of~$E$ in Proposition~\ref{ctsopen}(b).
\begin{rem}\label{precompare}
Let $E$ be a vector space and $\cO_1$, $\cO_2$
be vector topologies on~$E$.
Applying Proposition~\ref{ctsopen}\,(a)
to the identity map $(E,\cO_1)\to (E,\cO_2)$,
we see that $\cO_2\sub \cO_1$
if and only if every $0$-neighborhood of $(E,\cO_2)$
contains a $0$-neighborhood of $(E,\cO_1)$.
\end{rem}
\section{Vector subspaces, quotients and direct products}\label{secsubquot}
\begin{prop}\label{baconlcx}
Let $E$ be a topological vector space
and $(E_i)_{i\in I}$ be a family
of topological vector spaces.
\begin{description}[(D)]
\item[\rm (a)]
If $F\sub E$ is a vector subspace,
then the induced topology
makes $F$ a topological vector space.
Furthermore, the closure $\wb{F}$
of~$F$ is a vector subspace of~$E$.
If $E$ is locally convex, then so
is~$F$.
\item[\rm (b)]
The cartesian product $P:=\prod_{i\in I}E_i$,
equipped with the product topology
and componentwise addition and scalar
multiplication
is a topological vector space.
If each $E_i$ is locally convex, then
so is~$P$.
\item[\rm (c)]
If $F\sub E$ is a closed vector subspace,
then the quotient topology
with respect to $q\colon E\to E/F$,
$q(x):=x+F$
makes the
quotient vector space $E/F$
a topological vector space, which is
locally convex if so is~$E$.
The quotient map $q\colon E\to E/F$ is
open.
\end{description}
\end{prop}
\begin{prf}
(a) Since addition $\alpha\colon E\times E \to E$
is continuous, so is its restriction $\alpha|_{F\times F}\colon
F\times F\to E$
and the co-restriction $\alpha|_{F\times F}^F\colon
F\times F\to F$, which is the addition on~$F$.
Similarly, continuity of scalar multiplication on~$F$
is inherited from that on~$E$.
Thus~$F$ is a topological vector space.
Since $\alpha$ is continuous
and $F+F=\alpha(F\times F)\sub F$,
we deduce that $\wb{F}+\wb{F}=\alpha(\wb{F}\times \wb{F})
\sub \wb{F}$.
Likewise, $\K \, \wb{F}\sub \wb{F}$,
and thus $\wb{F}$ is a vector subspace
of~$E$.
If $E$ is locally convex and $V\sub F$
is a $0$-neighborhood,
by definition of the induced topology
there exists a $0$-neighborhood
$U\sub E$ such that $U\cap F\sub V$.
There exists a convex $0$-neighborhood
$W\sub E$ such that $W\sub U$.
Then $F\cap W$ is convex
(as an intersection of convex sets),
and it is a $0$-neighborhood of~$F$ contained
in~$V$. Hence~$F$ is locally convex.

(b) Let $\pr_i\!: P\to E_i$ be the coordinate
projection for $i\in I$.
The addition $\alpha \colon P\times P\to P$
is given by $\alpha((x_i)_{i\in I},(y_i)_{i\in I})
=(\alpha_i(x_i,y_i))_{i\in I}$
in terms of the addition $\alpha_i$
of~$E_i$. Hence $\pr_i\circ \, \alpha=\alpha_i\circ
(\pr_i\times \pr_i)$
is continuous for each $i\in I$,
whence
$\alpha$ is continuous
by Lemma~\ref{ctsinprod}(b).
Similarly, we see that scalar multiplication on~$P$
is continuous.
If each $E_i$ is locally convex
and $U\sub P$ is a $0$-neighborhood,
there exists a finite set $J\sub I$
and $0$-neighborhoods
$U_i\sub E_i$ such that $V:=\bigcap_{i\in J}\pr_i^{-1}(U_i)
\sub U$, by definition of the product
topology.
After shrinking~$U_i$,
we may assume that each~$U_i$ is convex.
Then also the preimage $\pr_i^{-1}(U_i)$
under the linear map $\pr_i$
is convex (exercise),
and hence the intersection~$V$
of convex sets is convex,
whence $V$ is a convex $0$-neighborhood
contained in~$U$.

(c) We equip $Q:=E/F$ with the quotient topology.
If $U\sub E$ is open, then also
$U+x$ is open in~$E$ for each $x\in F$
(cf.\ Lemma~\ref{verybase}\,(a))
and hence $q^{-1}(q(U))=U+F=\bigcup_{x\in F} (U+x)$ is open
as a union of open subsets.
Thus $q(U)$ is open in~$Q$, by definition of
the quotient topology.
Now $q$ being open and surjective, so is
$q\times q\!: E\times E\to Q\times Q$,
$(x,y)\mto (q(x),q(y))$,
whence also $q\times q$
is a quotient map.
For the addition maps $\alpha$ on~$E$ and $\beta$
on~$Q$, we have
\[
\beta\circ (q\times q)\; =\; q\circ \alpha\,.
\]
Here the right hand side is continuous and hence
so is $\beta\circ (q\times q)$.
Since $q\times q$ is a quotient map,
this entails that~$\beta$ is continuous
(see Lemma~\ref{ctsonquot}).
Similarly, we find that scalar multiplication
is continuous on $\K\times Q$.
Hence~$Q$ is a not necessarily Hausdorff topological vector space.
Now $\{0\}$ is closed in~$Q$,
because $q^{-1}(\{0\})=F$ is closed in~$E$ and~$q$
a quotient map.
Hence~$Q$ is Hausdorff (see Lemma~\ref{TVSHausdorff})
and thus~$Q$ is a topological vector space.

If, finally, $E$ is locally convex and $R\sub Q$
a $0$-neighborhood,
then $q^{-1}(R)$ is a $0$-neighborhood
in~$E$ and hence contains a convex, open
$0$-neighborhood~$S$ of~$E$.
Then $q(S)\sub R$
is open in~$Q$
and is convex, being the image
of a convex set under a linear map (exercise).
\end{prf}
A special type of vector subspaces
is of interest.
\begin{defn}\label{defcplsub}
A vector subspace $X$ of a topological
vector space~$E$ is called \emph{complemented}
(in~$E$) or \emph{split} if there exists a vector subspace
$Y\sub E$ (a ``complement'') such that
the addition map
\[
X\times Y\to E\,,\quad (x,y)\mto x+y
\]
is an isomorphism of topological vector spaces.
We then write
$E=X\oplus Y$.
\end{defn}
\begin{rem}
We compile some elementary facts
concerning complemented vector subspaces.
The details are left to the reader
as an exercise.
\begin{description}[(D)]
\item[(a)]
If $F$ is complemented in~$E$,
then $F$ is closed in~$E$.
\item[(b)]
If~$E$ is a Banach space
and $F\sub E$ a closed vector subspace,
then a vector subspace $H\sub E$ is a complement
for~$F$ if and only if $H$ is closed,
$F+H=E$, and $F\cap H=\{0\}$.
This follows from the open mapping theorem.
\item[(c)]
Every closed vector subspace
of a Hilbert space is complemented.
\end{description}
\end{rem}
\begin{small}
\subsection*{Exercises for Section~\ref{secsubquot}}
\begin{exer}\label{exc-inittvs}
Let $E$ be a vector space
and $(f_j)_{j\in J}$ be a family of linear maps
$f_j\colon E\to E_j$ to topological vector spaces~$E_j$
which separates points on~$E$.
Then the initial topology $\cO$ on~$E$ with respect to $(f_j)_{j\in J}$
turns the linear map
\[
f:=(f_j)_{j\in J}\colon E\to\prod_{j\in J}E_j
\]
into a homeomorphism onto the image.
Deduce that $(E,\cO)$ is a topological vector space.
If each $E_j$ is locally convex, then also $(E,\cO)$ is locally convex.
\end{exer}
\end{small}
\section{Neighborhoods of $0$
and continuous seminorms}\label{secsemino}
Although a locally convex space
need not be normable (i.e.,
there need not be a norm on it
which defines its topology),
locally convex spaces are
not too far away from
normed spaces:
At least,
their topology
can always be defined by a \emph{set}
of \emph{semi}norms,
as we recall now.
\begin{defn}\label{defnsemino}
Let $E$ be a vector space.
A map $p\!: E\to [0,\infty[$
is called a {\em seminorm\/} on~$E$
if it satisfies the following conditions:
\begin{description}
\item[(a)] $\,p(x+y)\leq p(x)+p(y)$ for all $x,y\in E$
(subadditivity);
\item[(b)]
$\,p(zx)=|z|p(x)$ for all $z\in \K$ and $x\in E$.
\end{description}
Given $x\in E$ and $\ve>0$,
we write
$B_\ve^p(x):=\{y\in E\colon p(y-x)<\ve\}$
and $\wb{B}_\ve^p(x):=\{y\in E\colon p(y-x)\leq \ve\}$.
If $(E,\|\cdot\|)$ is a normed space,
we usually write $B^E_r(x)$ (or simply $B_r(x)$)
for $B^{\|\cdot\|}_r(x)$,
and $\wb{B}^E_r(x)$ (or $\wb{B}_r(x)$)
for $\wb{B}^{\|\cdot\|}_r(x)$.
\end{defn}
\begin{numba}\label{modout}
In contrast to the case of a norm,
there may be non-zero vectors $x\in E$
such that $p(x)=0$.
It follows from (a) and (b)
that $N_p:=\{x\in E\!: p(x)=0\}$
is a vector subspace of~$E$, for each
seminorm~$p$ on~$E$
(for example, if $x,y\in N_p$,
we have $0\leq
p(x+y)\leq p(x)+p(y)=0$, whence $p(x+y)=0$ and thus
$x+y\in N_p$).
We let $E_p:=E/N_p$ be the quotient vector
space, formed by the cosets $x+N_p$
with $x\in E$.
Then
\[
\|\cdot\|_p\!: E_p\to [0,\infty[\,,\quad
\|x+N_p\|_p\,:=\, p(x)
\]
is well defined
(if $x+N_p=y+N_p$,
then $x=y+z$ for some $z\in N_p$
and thus $p(x)\leq p(y)+p(z)=p(y)$;
likewise $p(y)\leq p(x)$).
Furthermore, since~$p$
is a seminorm, so is $\| .\|_p$,
and it is in fact a norm on~$E_p$
because $0=\|x+N_p\|_p=p(x)$ entails
$x\in N_p$.
We let $\alpha_p\!: E\to E_p$,
$\alpha_p(x):=x+N_p$ be the canonical quotient map.\footnote{Occasionally,
we shall also write $\|\cdot\|_p$ instead of~$p$
for a seminorm on~$E$, to increase
the readability. No confusion with
$\|\cdot\|_p$ on $E_p$ just defined is likely.}
\end{numba}
Convex $0$-neighborhoods
and continuous seminorms are closely related.
\begin{prop}\label{Minkowsk}
Let $E$ be a topological vector space
and $U\sub E$ be an
absolutely convex $0$-neighborhood.
Then the Minkowski functional
\[
\mu_U\colon E\to [0,\infty[\,,\quad
\mu_U(x)\,:=\,
\inf\, \{ t>0\!: x\in tU\}
\]
is a continuous seminorm on~$E$.
If $U$ is closed,
then~$U$
coincides with the closed
unit ball $\{x\in E\!: \mu_U(x)\leq 1\}$
of~$\mu_U$.
Conversely, for every continuous seminorm
$p\colon E\to [0,\infty[$, its closed
unit ball $B:=\wb{B}^p_1(0)$
is an absolutely convex $0$-neighborhood in~$E$
and $\mu_B=p$.
\end{prop}
\begin{prf}
The final assertion is clear.
To prove the others, let $U\sub E$ be an absolutely
convex $0$-neighborhood.
Since $U$ is absorbing
by Lemma~\ref{verybase}(e),
we have $\mu_U(x)<\infty$ for each~$x\in E$.
Clearly $\mu_U(0)=0$ and
thus $\mu_U(zx)=0$ if $z=0$.
Now suppose that
$z\in \K^\times$.
If $t>0$ such that $zx\in tU$,
then $x\in tz^{-1}U=t|z|^{-1}U$
and hence $t|z|^{-1}\geq \mu_U(x)$.
Passing to the infimum, we obtain
$\mu_U(zx)|z|^{-1}\geq \mu_U(x)$
and thus
$\mu_U(zx)\geq |z|\mu_U(x)$.
Similarly,
$\mu_U(x)=\mu_U(z^{-1}zx)\geq |z|^{-1}\mu_U(zx)$
and thus $\mu_U(zx)\leq |z|\mu_U(x)$,
whence equality~holds.

If $x,y\in E$ and $t,s>0$ such that
$x\in tU$ and $y\in sU$, then 
\[ x+y\in tU+sU=(t+s) \big(
\frac{t}{t+s}U+\frac{s}{t+s}U\big)=(t+s)U,\]
entailing that $\mu_U(x)+\mu_U(y)\geq \mu_U(x+y)$.
Hence $\mu_U$ is a seminorm.

We now assume that~$U$ is closed.
If $x\in U$, then $x\in tU$
for $t=1$ and thus
$\mu_U(x)\leq 1$.
Conversely, let $x\in E$ such that
$\mu_U(x)\leq 1$. If $\mu_U(x)<1$,
then there is $t\in \;]0,1[$ such that
$x\in tU\sub U$ and thus $x\in U$.
If $\mu_U(x)=1$,
then $\mu_U(tx)=t<1$ for each $t\in \;]0,1[$
and hence $tx\in U$, by what has just been
shown. Letting $t\to 1$,
we find that $x\in \wb{U}=U$.
Hence indeed $U$ is the closed unit ball of~$\mu_U$.

As $\mu_U^{-1}([0,\ve])\supseteq \ve U$
for each $\ve>0$, the map $\mu_U$ is continuous
at~$0$.
Given $x\in E$, we have
$|\mu_U(x)-\mu_U(y)|\leq |\mu_U(x-y)|$
for all $x,y\in E$,
as a consequence of the
subadditivity of~$\mu_U$.
Hence
$\mu_U(x+\ve U)\sub [\mu_U(x)-\ve,\mu_U(x)+\ve]$
for each $\ve>0$, entailing that $\mu_U$ is continuous
at~$x$.
\end{prf}
\begin{defn}\label{defndefntop}
Let $E$ be a vector space.
\begin{description}
\item[(a)]
A set $\cP$
of seminorms on~$E$ is called
{\em separating\/}
if $p(x)=0$
for all $p\in \cP$
implies $x=0$.
\item[(b)]
If $\cP$ is a separating set of
seminorms on~$E$, then the map
\[
\alpha:=(\alpha_p)_{p\in \cP}\colon
E\to \prod_{p\in \cP}\, E_p\,,
\quad \alpha(x):=(\alpha_p(x))_{p\in \cP}
\]
is injective and linear (where
$E_p$ and $\alpha_p$ are as in~\ref{modout}).
We equip the product of normed
spaces with the product topology
(see Proposition~\ref{baconlcx}(b)).
There is a uniquely
determined topology~$\cO$ on~$E$ which makes
$\alpha$ a homeomorphism onto $\alpha(E)$,
equipped with the induced topology;
it is called the {\em locally convex topology on~$E$
defined by the set of seminorms~$\cP$.}
\end{description}
\end{defn}
Thus $\cO$ is the initial topology on~$E$
with respect to the family $(\alpha_p)_{p\in \cP}$.
\begin{rem}
Since $\alpha(E)$ is a locally
convex space (being a vector subspace
of a locally convex space)
and the co-restriction
$\alpha|^{\alpha(E)}\!: E\to \alpha(E)$
is an isomorphism of vector spaces,
it readily follows that $\cO$ makes~$E$
a locally convex space.
Furthermore,
each $p\in \cP$ is
a {\em continuous\/} seminorm on $(E,\cO)$,
as it can be written
as the composition
$p=\|\cdot\|_p\circ \pr_p\circ \alpha$
of continuous maps (where
$\pr_p$ is the projection from the cartesian
product onto
the factor $E_p$).
\end{rem}
\begin{rem}\label{basisgen}
Note that, in the situation
of Definition~\ref{defndefntop}(b),
a basis of $0$-neighborhoods
for the locally convex topology on~$E$ defined
by $\cP$ is given by the sets of the form
\begin{equation}\label{typicnbd}
B_\ve^{p_1}(0)
\cap \cdots \cap B_\ve^{p_n}(0)\,,
\end{equation}
where $\ve>0$, $n\in \N$ and
$p_1,\ldots, p_n\in \cP$.
In fact, clearly each of the sets in (\ref{typicnbd})
is an open $0$-neighborhood in~$E$,
as the locally convex topology defined
by~$\cP$ makes each $p\in \cP$
a continuous seminorm.
Given a $0$-neighborhood $U\sub E$,
there exists a finite subset
$F\sub \cP$
and $0$-neighborhoods $U_p\sub E_p$
for $p\in F$
such that $\alpha(U)\supseteq
\alpha(E)\cap
\bigcap_{p\in F}\pr_p^{-1}(U_p)$
(by definition of the product topology
and the induced topology
on~$\alpha(E)$).
Now, $E_p$ being normed,
we find $\ve_p$ such that $\{x\in E_p\colon \|x\|_p<\ve_p\}\sub U_p$.
Set $\ve:=\min \{\ve_p\colon p\in F\}$.
Then $V:=\bigcap_{p\in F}B_\ve^p(0)\sub U$
because $\alpha(V)\sub \alpha(E)\cap
\bigcap_{p\in F}\pr_p^{-1}(U_p)\sub \alpha(U)$,
which in turn holds because
$\|\pr_p(\alpha(x))\|_p=\|\alpha_p(x)\|_p=p(x)<\ve$
for each $x\in V=\bigcap_{p\in F}B_\ve^p(0)$
and each $p\in F$.
\end{rem}
The following observation is clear
from Remarks~\ref{precompare} and \ref{basisgen}:
\begin{lem}\label{comparedtop}
Let $E$ be a vector space and
$\cP_1$, $\cP_2$ be separating sets
of seminorms on~$E$.
For $j\in \{1,2\}$,
let $\cO_j$ be the vector topology on~$E$
defined by~$\cP_j$.
Then the following conditions are equivalent:
\begin{description}[(D)]
\item[\rm (a)]
For each
$p\in \cP_2$, there exist $n\in \N$,
$p_1,\ldots, p_n\in \cP_1$ and
$r_1,\ldots, r_n>0 $
such that $p\leq r_1p_1+\ldots+r_np_n$.
\item[\rm (b)]
$\cO_2\sub \cO_1$ holds, i.e.,
$\cO_2$ is coarser than $\cO_1$.\qed
\end{description}
\end{lem}
In particular,
Lemma~\ref{comparedtop}
tells us when two sets of seminorms
define the same vector topology.
The following fact follows from
Proposition~\ref{Minkowsk}:
\begin{prop}\label{embedprod}
If $E$ is a locally convex topological vector space,
then the set $\cP$ of all
continuous seminorms on~$E$ is separating.
For each $p\in \cP$,
the map $\alpha_p\!: E\to E_p$ is continuous,
and
\[
\alpha:=(\alpha_p)_{p\in \cP}
\colon E\to \prod_{p\in \cP}\, E_p \, =:\, P
\]
is linear and a topological embedding
$($i.e., a homeomorphism onto its image$)$.
In other words, $\cP$
defines the given locally convex topology on~$E$.
\end{prop}
\begin{prf}
For each $p\in \cP$,
the linear map $\alpha_p$ is continuous
as $\alpha^{-1}(B_\ve^{\|\cdot\|_p}(0))=B_\ve^p(0)$
is a $0$-neighborhood for each $\ve>0$.
Let $\pr_p\colon P\to E_p$
be the canonical projection.
Since $\pr_p\circ \alpha=\alpha_p$ is continuous
for each $p\in \cP$, the map $\alpha$ is continuous
(see Lemma~\ref{ctsinprod}(b)).
Furthermore,
clearly $\alpha$ is linear since so is each~$\alpha_p$.
Given $0\not=x\in E$,
there exists a $0$-neighborhood
$U\sub E$ such that $x\not\in U$.
By Lemma~\ref{basisnice},
there exists a closed, absolutely
convex $0$-neighborhood $B\sub E$ such that
$B\sub U$. By Proposition~\ref{Minkowsk},
there is a continuous seminorm $p\in \cP$
with closed unit ball~$B$. Since $x\not \in B$,
we deduce that $1<p(x)=\|\alpha_p(x)\|_p$,
whence $\alpha_p(x)\not=0$ and
hence $\alpha(x)\not=0$.
Therefore $\alpha$ is injective.
It remains to show
that $\alpha$ is open onto its image.
By Proposition~\ref{ctsopen}\,(b), we only
need to show that $\alpha(V)$ is
a $0$-neighborhood in $\alpha(E)$,
for each $0$-neighborhood $V\sub E$.
As before, we see that $V$
contains the closed unit ball~$Q$
of some continuous seminorm $q\in \cP$.
Given $y \in \alpha(E)$, we have $y=\alpha(x)$
for some $x\in E$;
due to the injectivity of~$\alpha$,
\[
y=\alpha(x)\in \alpha(Q)
\]
holds if and only if $y\in Q$.
The latter holds if and only if
\[ 1\geq q(x)
=\|\alpha_q(x)\|_q=\|\pr_q(y)\|_q,\]
i.e., if and only if $y\in \pr_q^{-1}(R)$,
where~$R\sub E_q$ is the closed unit ball
of~$\|\cdot\|_q$.
Hence $\alpha(V)\supseteq \alpha(Q)
=\alpha(E)\cap \pr_q^{-1}(R)$,
which is a $0$-neighborhood.
\end{prf}
\begin{rem}\label{nulittlerem}
In the situation of
Proposition~\ref{embedprod},
let $\cQ\sub \cP$ be a set of continuous
seminorms such that, for each $p\in \cP$,
there exist $q_1,\ldots, q_n\in \cQ$
and $r_1,\ldots, r_n>0$ such that
$p\leq r_1q_1+\cdots+ r_nq_n$.
Then $\cQ$ defines the
given locally
convex vector topology on~$E$
(as follows from
Proposition~\ref{embedprod} and
Lemma~\ref{comparedtop}).
Hence
$(\alpha_q)_{q\in \cQ}\colon E\to \prod_{q\in \cQ}E_q$
is a topological embedding.
\end{rem}
\begin{rem}\label{littlerem}
Let $E$ be a vector space
and $p,q\colon E\to [0,\infty[$
be seminorms such that $p\leq q$
(pointwise).
Then $q(x)=0$ entails
$p(x)=0$, whence $N_q\sub N_p$.
By the homomorphism theorem,
there is a unique linear map
\[
\alpha_{pq}\colon E_q\to E_p
\]
such that $\alpha_{pq}\circ\alpha_q=\alpha_p$.
Given $y\in E_q$, we have $y=\alpha_q(x)$
for some $x\in E$.
Then
\[ \|\alpha_{pq}(y)\|_p=\|\alpha_{pq}(\alpha_q(x))\|_p=\|\alpha_p(x)\|_p
=p(x)\leq q(x)=\|\alpha_q(x)\|_q=\|y\|_q \] 
shows that
$\alpha_{pq}$ is a continuous
linear map, of operator norm $\|\alpha_{pq}\|\leq 1$.
\end{rem}
\begin{ex}
Let $r\in \N_0\cup\{\infty\}$.
For $k\in \N_0$ such that
$k\leq r$, the map
\[
p_k\!: C^r([0,1],\R)\to [0,\infty[\,,
\quad
p_k(\gamma)\, :=\, \sup\{|\gamma^{(k)}(t)|\!: t\in [0,1]\}
\]
is a seminorm on the space
$C^r([0,1],\R)$
of real-valued $C^r$-maps on
$[0,1]$.
The set $\cP$
of all~$p_k$ as before
is a separating set of
seminorms on $C^r([0,1],\R)$.
We give $C^r([0,1],\R)$
the locally convex topology defined by~$\cP$.
\end{ex}
Recall that a topological space~$X$
is called {\em metrizable\/}
if there exists a metric~$d$
on~$X$ whose associated topology
coincides with the given one.
\begin{cor}\label{charmetriz}
For a locally
convex
topological vector space~$E$,
the following conditions
are equivalent:
\begin{description}
\item[\rm (a)]
$E$ is metrizable;
\item[\rm (b)]
$E$ is first countable, i.e.,
for each $x\in E$ there exists a countable
basis of $x$-neighborhoods;
\item[\rm (c)]
There exists a countable basis
of $0$-neighborhoods in~$E$;
\item[\rm (d)]
There exists a countable set
$\cP$ of continuous seminorms
on~$E$ which defines its topology.
\end{description}
\end{cor}
\begin{prf}
(a)$\impl$(b): Let $d$ be a metric on~$E$ defining its
topology. Then, for each $x\in E$, the balls
$B_{1/n}(x):=\{y\in E\colon d(x,y)<1/n\}$
(indexed by $n\in \N$)
form a countable basis of $x$-neighborhoods.

(b)$\impl$(c) is trivial.

(c)$\impl$(d):
Let $\cU$ be a countable basis
of $0$-neighborhoods in~$E$. For each $U\in \cU$,
there exists a closed absolutely convex
$0$-neighborhood $B_U\sub U$
(see Lemma~\ref{basisnice}) and
a continuous seminorm $p_U$ on~$E$
with closed unit ball $B_U$
(see Proposition~\ref{Minkowsk}).
Then $\cP:=\{p_U\colon U\in \cU\}$
is a countable set of continuous
seminorms which is separating as $\bigcap_{U\in \cU}B_U=\{0\}$.
Set $\cB:=\{B_U\colon U\in \cU\}$.
For each finite sequence $U_1,\ldots, U_n$ in~$\cU$
and $r>0$,
the set $r B_{U_1}\cap\cdots\cap rB_{U_n}$
is a $0$-neighborhood in~$E$
and hence contains some $U\in \cU$
and hence also $B_U$.
Therefore
Remark~\ref{basisgen} implies that
$\cB$ is a basis of $0$-neighborhoods
for the locally convex topology on~$E$
defined by~$\cP$.
But $\cU$ also is a basis
of $0$-neighborhoods for the original
topology, entailing that the two
vector topologies coincide.

(d)$\impl$(a): Suppose that
$(p_n)_{n\in \N}$ be a sequence
of continuous seminorms on~$E$
which
defines the locally convex topology~$\cO$ on~$E$.
Then
\[
d\colon E\times E\to [0,\infty[\, ,
\quad d(x,y) \,
:=\, \sum_{n=1}^\infty \frac{2^{-n}\, p_n(x-y)}{1+p_n(x-y)}
\]
is a continuous function,
because the series converges uniformly
on~$E\times E$.
Furthermore,
using that
\[
h\colon [0,\infty[\;\to [0,1[\,,\qquad
h(t)\, :=\, \frac{t}{1+t}\,=\, 1-\frac{1}{t+t}
\]
is monotonically increasing
and $h(s+t)\leq h(s)+h(t)$ for all $s,t\in [0,\infty[$,
it is easy to see that $d$
is a metric (exercise).
Since $d$ is continuous on $(E,\cO)$,
the topology $\cT$ determined by~$d$
is coarser then the given one,
$\cT\sub \cO$.
Given $x\in E$ and a neighborhood
$U$ of $x$ in~$(E,\cO)$,
using Lemma~\ref{verybase}(f)
and Remark~\ref{basisgen},
we find
$\ve>0$ and $n\in \N$
such that $B:=\bigcap_{k=1}^n B_\ve^{p_k}(x)\sub U$.
Set $\delta := 2^{-n} h(\ve)$.
If $y\in E$ such that $d(x,y)<\delta$,
then
$2^{-k} h(p_k(x-y))< \delta$
for each $k\in \{1,\ldots, n\}$,
whence $h(p_k(x-y))<2^k \delta \leq 2^n \delta=h(\ve)$
and hence $p_k(x-y)<\ve$
as $h$ is strictly monotonically
increasing.
Therefore $\{y\in E\colon d(x,y)<\delta\}
\sub B\sub U$.
Since, for each $x\in E$,
every
$x$-neighborhood in $(E,\cO)$
also is an $x$-neighborhood
in $(E,\cT)$,
we deduce that
$\cO \sub \cT$.
Hence both topologies coincide.
\end{prf}
\begin{defn}\label{normable}
A locally convex space $E$ is called
{\em normable\/}
if its locally convex vector topology
can be defined using a single norm.
\end{defn}
\begin{numba}[{\bf Further examples}]\label{furex}
Here are further examples
of locally convex spaces (which are not normable):
\begin{description}[(D)]
\item[\rm (a)]
By Corollary~\ref{charmetriz},
the product topology
makes $\R^\N:=\prod_{n\in \N}\R$
a metrizable locally convex space
(as the countable set of the seminorms $\R^\N\to[0,\infty[$, $f\mto |f(n)|$
with $n\in\N$ is separating and defines the product topology
by Definition~\ref{defndefntop}(b)). 
More generally, every countable
direct product of normed spaces
is metrizable.
\item[\rm (b)]
We equip $C(\R,\R)$, the space of continuous
real-valued functions on~$\R$,
with the locally convex topology
defined by the set of all
seminorms of the form
\[
p_K \colon C(\R,\R)\to  [0,\infty[\,,\quad
p_K(\gamma):=\sup\{|\gamma(x)|\colon x\in K\}\,,
\]
where $K$ ranges through the set
of all compact subsets~$K$ of~$\R$.
The (same) locally convex topology
on $C(\R,\R)$ is defined
by the countable set $\{p_{[-n,n]}\colon n\in \N\}$
of seminorms (exercise),
showing that $C(\R,\R)$ is metrizable.
Following the same pattern,
$C(X,\K)$ can be made a locally convex topological $\K$-vector space,
for each topological space~$X$.
\item[\rm (c)]
Given an open subset $\Omega\sub \C$,
the set $\Hol(\Omega)$ of all holomorphic
functions $\gamma\colon \Omega\to\C$
is a vector subspace of the
locally convex space $C(\Omega,\C)$
(defined as in (b)),
which is metrizable since we can find
a sequence $K_1\sub K_2\sub\cdots$
of compact subsets of~$\Omega$
such that $\Omega=\bigcup_{n\in \N}K_n$
and each compact subset of~$K\sub \Omega$
is contained in some $K_n$ (see Exercise~\ref{lcp-hemi}).
The induced topology makes $\Hol(\Omega)$
a metrizable locally convex space.
\end{description}
\end{numba}
\begin{small}
\section*{Exercises for Section~\ref{secsemino}}

\begin{exer}\label{exc-mink-op}
Show that
$\{x\in E\colon\mu_U(x)<1\}\sub U\sub \{x\in E\colon \mu_U(x)\leq 1\}$
in the situation of Proposition~\ref{Minkowsk}.
If~$U$ is open, then $U=\{x\in E\colon \mu_U(x)<1\}$.
\end{exer}
\end{small}
\section{Linear functionals and the Hahn--Banach Theorem}\label{sechbanach}
Given a topological $\K$-vector space~$E$, we let
$E'$ be the set of all continuous linear
maps (``functionals'') $\lambda\colon E\to\K$.
Since sums and scalar multiples of
continuous linear functionals
are such,
$E'$ is a vector subspace of the space $\K^E$ of
all $\K$-valued functions on~$E$.
We call $E'$ the {\em dual space\/} of~$E$.

The Hahn--Banach theorems
provide continuous linear functionals
with specific properties.
Notably, they ensure
that the continuous
linear functionals on a locally convex space~$E$
separate the points of~$E$ (Theorem~\ref{dualsep}).
The Hahn--Banach theorems are powerful tools
of functional analysis.
We shall follow a
geometric, calculation-free approach
to the Hahn--Banach theorems.\\[3mm]
Two preliminary lemmas will be used:
\begin{lem}\label{1dim1}
Let $E$ be a real topological vector space.
If $E\setminus \{0\}$ is disconnected,
then $E$ is $1$-dimensional.
\end{lem}
\begin{prf}
If $\dim(E)>1$,
let $x,y\in E\setminus \{0\}$.
If $x$ and $y$ are linearly independent,
then $[0,1]\to E$, $t\mto x+t(y-x)$
is a path from $x$ to~$y$ in $E\setminus \{0\}$.
If $x$ and $y$ are linearly dependent,
we choose $z \in E\setminus \Spann \{x,y\}$.
Then $x$ and $z$
as well as $z$ and $y$ 
are linearly independent,
whence $x$ and $z$
(resp., $z$ and $y$)
can be joined by a path $\gamma$ (resp., $\eta$)
in $E\setminus\{0\}$, as just shown.
Concatenating the two paths we obtain a path from
$x$ to~$y$.
We have shown that $E\setminus \{0\}$
is path connected and hence connected.
\end{prf}
\begin{lem}\label{1dim2}
Let $E$ be a real topological
vector space.
If there exists an open, non-empty
convex subset
$U\sub E\setminus \{0\}$ such that
$U \cap \R x \not=\emptyset$
for each $x\in E\setminus\{0\}$,
then $E$ is $1$-dimensional.
\end{lem}
\begin{prf}
Abbreviate $P:=\;]0,\infty[$.
The set $PU$ is stable under
multiplication with positive scalars
(as $PP=P$).
Furthermore, $PU+PU\sub PU$
because $av+bw=(a+b)\big( \frac{a}{a+b}v+\frac{b}{a+b}w\big)
\in (a+b)U\sub PU$
for all $a,b\in P$
and $v,w\in U$, exploiting the convexity
of~$U$.
Hence $PU$ is what one calls a
``cone''
(a set stable under addition and multiplication with
positive scalars).
As $PU=\bigcup_{r>0}rU$,
the set $PU$ is open.
Furthermore, $0\not\in PU$,
since $0\not\in U$.
We also have $PU\cap (-PU)=\emptyset$,
because $au=-bv$
with $a,b>0$, $u,v\in U$
would entail that $0=au+bv\in PU+PU\sub PU$,
which we just excluded.
Since $U$ meets $\R x$ for each $0\not= x\in E$,
we have
\begin{equation}\label{crucdisj}
E\;=\; PU\cup (-PU)\cup \{0\}\,.
\end{equation}
By the preceding, the union in (\ref{crucdisj})
is
disjoint.
Therefore
$E\setminus \{0\}=PU\cup (-PU)$
as a disjoint union, where
both $PU$ and $-PU$
are open and non-empty.
Hence
$E\setminus\{0\}$
is not connected,
whence $E$ is $1$-dimensional
by Lemma~\ref{1dim1}.
\end{prf}
\begin{thm}[{\bf Hahn--Banach Theorem, real case}]\label{geomHBreal}
If $E$ is a real topological vector space
and $U\sub E$ an open, convex subset
such that $0\not\in U$,
then there exists $\lambda\in E'$ such that
\[
\lambda(U)\; \sub \;\, ]0,\infty[\,.
\]
\end{thm}
\begin{prf}
Without loss of generality $U\not=\emptyset$
(otherwise, take $\lambda:=0$).
We consider the set $\cA$ of all vector subspaces
$F\sub E$ such that $F\cap U=\emptyset$.
Inclusion of sets provides a partial  order~on~$\cA$
which makes $\cA$
an inductive set.
Indeed, let $\Gamma\sub \cA$ be a chain
(i.e., a totally ordered subset).
If $\Gamma=\emptyset$, then the trivial
vector space $\{0\}\in \cA$ is an upper bound
for~$\Gamma$. Otherwise, clearly $\bigcup\Gamma\in \cA$,
and this union is an upper bound for~$\Gamma$.
Hence, by Zorn's Lemma,
there exists a maximal element $H\in \cA$.
We claim that $H$ is a closed hyperplane,
i.e., a closed vector subspace
such that $\dim(E/H)=1$.

{\em $H$ is closed}: Since the closure $\wb{H}$
is a vector subspace and $\wb{H}\sub E\setminus U$
(as the latter set is closed and contains~$H$),
we deduce $H=\wb{H}$ from the maximality of~$H$.

{\em $E/H$ is $1$-dimensional.}
To see this, let $q\colon E\to E/H$ be
the quotient map and set $V:=q(U)$.
Since $q$ is linear and open, the set $V$
is an open, convex subset of~$E/H$.
Note that $V$ does not contain
the $0$-element~$H$
of the quotient space, because $H\in V$
would entail
\[ \emptyset\not=
q^{-1}(\{H\}\cap V)=
q^{-1}(\{H\})\cap q^{-1}(V)=
H\cap (U+H) \] 
whence there exist $h_1,h_2\in H$ and $u\in U$
such that $h_1=u+h_2$
and therefore $u=h_1-h_2\in U\cap H$,
which is a contradiction. 
We now show that every vector subspace
$K\sub E/H$ such that $K\cap V=\emptyset$
is trivial. To see this, note that
$q^{-1}(K)$ is a vector subspace of~$E$
such that $H\subseteq q^{-1}(K)$.
Furthermore, $U\cap q^{-1}(K)\sub q^{-1}(q(U))\cap q^{-1}(K)=q^{-1}(V)\cap
q^{-1}(K)=q^{-1}(V\cap K)=\emptyset$.
Hence $H=q^{-1}(K)$, by maximality of~$H$,
and thus $K=q(q^{-1}(K))=\{H\}$
is the trivial vector subspace.
Using Lemma~\ref{1dim2},
we see that the locally convex space
$E/H$ is $1$-dimensional
and hence isomorphic to~$\R$
as a consequence of Corollary~\ref{onlyyou}.
We choose an isomorphism of topological
vector spaces $\phi\colon E/H\to \R$
and set $\lambda:=\phi \circ q\colon
E\to \R$.
Then $\lambda$ is a continuous linear
functional such that $\ker\, \lambda=H$
does not meet~$U$,
whence $0\not\in \lambda(U)$.
Being the image of a convex set under
a linear map, $\lambda(U)$ is a convex
subset of~$\R$.
As $0\not\in \lambda(U)$,
we must have $\lambda(U)\sub \;]0,\infty[$
or $\lambda(U)\sub \;]{-\infty}, 0[$.
In the first case, we are home;
in the second case, we replace $\lambda$
with ${-\lambda}$.
\end{prf}
To transfer the Hahn--Banach Theorem from the real case to
the complex case,
we first establish a
one-to-one correspondence between real linear
and complex linear functionals
on a complex topological vector space.
\begin{lem}\label{lemrealcx}
Let $E$ be a complex topological
vector space.
If $\lambda\colon E\to \C$ is
a continuous complex linear
functional on~$E$,
then $\Repart\circ \lambda\colon E\to \R$
is a continuous real linear functional.
Conversely, for each continuous
real linear functional $u\colon E\to \R$,
there is a unique continuous
complex linear functional
$\lambda\colon E\to\C$ such that
$u=\Repart\circ \lambda$;
it is given by the formula
\begin{equation}\label{reptocx}
\lambda(x)\; =\; u(x)-iu(ix)\qquad\mbox{for all $\, x\in E$.}
\end{equation}
\end{lem}
\begin{prf}
It is obvious that $u:=\Repart\circ \lambda$
is continuous real linear if $\lambda\in E'$.
Conversely, let $u\colon E\to \R$
be a continuous real linear functional.
If there exists $\lambda\in E'$ such that
$u=\Repart\circ \lambda$,
then $\lambda=u+iv$ with
$v:=\Impart\circ \lambda$.
For each $x\in E$, we have
\[
u(ix)+iv(ix)\;=\; \lambda(ix)\; =\; i\lambda(x)\;=\; i(u(x)+iv(x))
\; =\; iu(x)-v(x)\,;
\]
comparing the real and imaginary parts,
we find that $v(x)=-u(ix)$.
Hence $\lambda(x)=u(x)-iu(ix)$
indeed, showing in particular that
$\lambda$ is uniquely
determined by its real part~$u$.
To see that a complex linear functional
with real part~$u$ exists,
we simply define a map
$\lambda\colon E\to \C$ via~(\ref{reptocx}).
Then $\lambda$ is continuous,
real linear, and $u=\Repart\circ \lambda$.
Being real linear, clearly
$\lambda$ will be complex linear
if we can show that $\lambda(ix)=i\lambda(x)$.
This is the case:
We have $\lambda(ix)=u(ix)-iu(i(ix))=i(u(x)-iu(ix))=i\lambda(x)$.
\end{prf}
\begin{thm}[{\bf Hahn--Banach Theorem, complex case}]\label{complanalog}
If $E$ is a complex topological vector space,
and $U\sub E$ an open, convex subset
such that $0\not\in U$,
then there exists $\lambda\in E'$ such that
\[
\Repart(\lambda(U))\; \sub \;\, ]0,\infty[\,.
\]
\end{thm}
\begin{prf}
By Theorem~\ref{geomHBreal},
there exists a continuous real linear
functional $u\colon E\to \R$
such that $u(U)\sub \;]0,\infty[$.
Lemma~\ref{lemrealcx}
provides $\lambda\in E'$ such that
$\Repart\circ \lambda =u$.
Then $\Repart(\lambda(U))=u(U)\sub \;]0,\infty[$.
\end{prf}
\begin{thm}\label{dualsep}
If $E$ is a locally convex topological
$\K$-vector space, then $E'$ separates
points on~$E$,
i.e., for all $x,y\in E$ such that $x\not=y$,
there exists $\lambda\in E'$
such that $\lambda(x)\not=\lambda(y)$.
\end{thm}
\begin{prf}
Since $x-y\not=0$, there exists an open,
convex $0$-neighborhood $U\sub E$ such that
$x-y\not\in U$.
Then $U-(x-y)$ is an open, convex subset of~$E$
such that $0\not\in U-(x-y)$.
The Hahn--Banach theorem provides
$\lambda\in E'$ such that $\Repart(\lambda(U-(x-y)))\sub \; ]0,\infty[$.
Then
\[ \Repart(\lambda(y))-\Repart(\lambda(x))
=\Repart(\lambda(-(x-y)))\in \;]0,\infty[ \] 
and thus $\Repart(\lambda(x))\not=\Repart(\lambda(y))$,
whence also $\lambda(x)\not=\lambda(y)$.
\end{prf}
The following observation is useful:
\begin{lem}\label{obsopen}
If $E$ is a topological $\K$-vector space
and $\lambda\colon E\to \K$
a linear functional such that $\lambda\not=0$,
then $\lambda$ is an open map.
\end{lem}
\begin{prf}
By Proposition~\ref{ctsopen}\,(b),
we only need to show that $\lambda(U)$
is a $0$-neighborhood in~$\K$
for each $0$-neighborhood~$U$ in~$E$.
After shrinking~$U$,
we may assume that $U$ is balanced.
Then also $\lambda(U)\sub \K$ is balanced.
Since $\lambda\not=0$,
there exists $x_0\in E$ such that $\lambda(x_0)\not=0$;
after replacing $x_0$ with $tx_0$
for sufficiently small $0\not= t\in \K$,
we may assume that $x_0\in U$.
Setting $r:=|\lambda(x_0)| >0$,
we then have
$\bD_r=\bD \lambda(x_0)\sub \lambda(U)$,
showing that indeed $\lambda(U)$
is a $0$-neighborhood.
\end{prf}
\begin{thm}[{\bf Hahn--Banach Separation Theorem}]\label{HBsep}
Let
$E$ be a real or complex locally convex space,
$A\sub E$ be a closed convex subset and
$K\sub E$ be a compact, convex non-empty set
such that $K\cap A=\emptyset$.
Then
there exists $\lambda\in E'$ and $r\in \R$ such that
\begin{equation}\label{sparation}
\Repart(\lambda(A))\;\sub \; \, ]{-\infty}, r[\quad
\mbox{and}\quad \Repart(\lambda(K))\sub  \;\, ]r,\infty[\,.
\end{equation}
\end{thm}
\begin{prf}
Since $K\sub E\setminus A$,
where $E\setminus A$ is open
(as $A$ is closed),
using Lemma~\ref{verybase}\,(h)
we find an
absolutely convex, open
$0$-neighborhood $W\sub E$ such that $K+W\sub E\setminus A$.
Then $K\cap (A+W)=\emptyset$ (otherwise,
$x=a+w$ for suitable elements $x\in K$, $a\in A$ and
$w\in W$
and hence $a=x -w \in A\cap (K+W)$,
contradiction).
Thus, $U:=A+W-K$ is an open, convex subset
of~$E$ such that $0\not\in U$.
Theorem~\ref{geomHBreal}
(resp., Theorem~\ref{complanalog})
provide $\mu \in E'$ such that
\begin{equation}\label{yieldssepa}
u(A)+u(W)-u(K)\;=\; u(U)\;\sub \, \; ]0,\infty[\,,
\end{equation}
where $u:=\Repart\circ \mu\colon E\to\R$.
Since $u(W)$ is an open $0$-neighborhood
in~$\R$ (by Lemma~\ref{obsopen}),
there is $w\in W$ such that $u(w)<0$.
The set $K$ being compact, we have
$\max u(K)=u(x_0)$ for some $x_0\in K$.
Now (\ref{yieldssepa})
shows that $u(A)\sub \; ]u(x_0)-u(w),\infty[$,
whence $u(A)\sub \; ]r,\infty[$
and $u(x_0)\in \;]{-\infty},r[$
with $r:=u(x_0)-\frac{1}{2}u(w)$.
Then also $u(K)\sub \;]{-\infty},r[$.
Now set $\lambda:=-\mu$.
\end{prf}
Thus, in the real case,
$A$ and $K$ are
contained in the disjoint open half-spaces
$\lambda^{-1}(]{-\infty}, r[)$ and
$\lambda^{-1}(]r,\infty[)$, respectively.
In other words, $A$ and $K$ can be separated
by the closed hyperplane $\lambda^{-1}(\{r\})$.\\[2.5mm]
Note that Theorem~\ref{HBsep}
applies in particular if $K:=\{x_0\}$,
where $x_0\in E\setminus A$.
\begin{thm}[{\bf Bipolar Theorem}]\label{bipolar}
Let $E$ be a real or complex locally convex space,
$A\sub E$ be an absolutely convex,
closed, non-empty subset, and $x_0\in E$ such that $x_0\not\in A$.
Then there is $\lambda\in E'$ such that
$|\lambda(x)|\leq 1$ for all $x\in A$,
but $|\lambda(x_0)|>1$.
\end{thm}
\begin{prf}
By the Hahn--Banach Separation Theorem,
there exists $\mu\in E'$ and $r\in \R$
such that $\Repart(\mu(A))\sub \;]{-\infty},r[$
and $\Repart(\mu(x_0))>r$.
Let $x \in A$.
There exists $z\in \K$ such that $|z|=1$
and $\mu(zx)=z\mu(x)\in [0,\infty[$.
Now $zx\in A$ as $A$
is balanced.
Hence
$0\leq |\mu(x)|=z\mu(x)=\mu(zx)=\Repart(\mu(zx))\in
\;]{-\infty},r[$,
entailing that $r>0$ and $|\mu(x)|<r$.
For $\lambda:=\frac{1}{r}\mu$,
we then have $|\lambda(x)|<1$
for each $x\in A$, and furthermore
$\Repart(\lambda(x_0))>1$ and thus $|\lambda(x_0)|>1$.
\end{prf}
The name of the Bipolar Theorem
has the following origin.
\begin{rem}\label{rembipol}
Given a subset $A$ of a locally convex space~$E$,
one calls
\[
A^\circ\; :=\; \{\lambda\in E'\!:
\mbox{$|\lambda(x)|\leq 1$ for all $x\in A$}\}
\]
the {\em polar of $A$ in~$E'$}.
Similarly, given $B\sub E'$
one defines the {\em polar of $B$ in~$E$\/}
as the set
$B_\circ:=\{x\in E\!: \mbox{$|\lambda(x)|\leq 1$
for all $\lambda\in B$}\}$.
It is obvious from
these definitions that
$A\sub (A^\circ)_\circ$ for each subset $A\sub E$.
The Bipolar Theorem entails that
\begin{equation}\label{formbipol}
A\;=\; (A^\circ)_\circ\,,
\end{equation}
for each absolutely convex,
closed, non-empty subset $A\sub E$ (exercise).
\end{rem}
\begin{thm}[{\bf Hahn--Banach Extension Theorem}]\label{HBext}
Let $E$ be a real vector space, $q\colon E\to [0,\infty[$ be a seminorm,
$E_0\sub E$ be a vector subspace and $\lambda\colon E_0\to\R$
be a linear functional such that
\[
|\lambda(x)|\leq q(x)\quad\mbox{for all $x\in E_0$.}
\]
Then there exists a linear functional $\Lambda\colon E\to\R$
such that $\Lambda|_{E_0}=\lambda$ and
\[
|\Lambda(x)|\leq q(x)\quad\mbox{for all $x\in E$.}
\]
If $E$ is a locally convex space and $q$ is continuous, then $\Lambda$ is continuous.
\end{thm}
\begin{prf}
Final statement:
If $q$ is continuous, then $\Lambda^{-1}(]{-\ve},\ve[)\supseteq q^{-1}(]{-\ve},\ve[)$
is a $0$-neighborhood for each $\ve>0$ and thus $\Lambda$ is continuous.\\[2.3mm]
To prove the existence of $\Lambda$,
endow $E$ with any Hausdorff locally convex vector topology which makes $q$
continuous. Then
\[
U:=\{(x,t)\in E\times \R\colon t>q(x)\}
\]
is an open convex subset of $E\times\R$
such that
\[
\text{graph}(\lambda)\cap U=\emptyset.
\]
The Hahn--Banach Separation Theorem (in whose proof we may consider vector
subspaces containing $\text{graph}(\lambda)$) yields a continuous linear functional
\[
\theta\colon E\times\R\to\R
\]
such that $\text{graph}(\lambda)\sub \ker(\theta)$ and $\ker(\theta)\cap U=\emptyset$.
We can write
\[
\theta(x,t)=\Lambda(x)+bt
\]
with $\Lambda(x):=\theta(x,0)$ and $b:=\theta(0,1)$.
If we suppose that $b=0$, then $\{0\}\times \R\sub\ker(\theta)$
and hence $\{0\}\times \;]0,\infty[\; \sub U\cap \ker(\theta)$,
contradiction. Hence $b\not=0$ and after replacing $\theta$ with $-\frac{1}{b}\theta$,
we may assume that $b=-1$. Thus
\[
\ker(\theta)=\{(x,t)\in E\times \R\colon\Lambda(x)-t=0\}=\text{graph}(\Lambda).
\]
Notably, $\text{graph}(\lambda)\sub\text{graph}(\Lambda)$ and thus $\Lambda|_{E_0}=\lambda$.
For each $x\in E$, we have $(x,\Lambda(x))\in\ker(\theta)$ and thus $(x,\Lambda(x))\not\in U$.
Thus
\[
\lambda(x)\leq q(x).
\]
Likewise, $-\lambda(x)=\lambda(-x)\leq q(-x)=q(x)$ and thus $|\lambda(x)|\leq q(x)$.
\end{prf}
\begin{cor}\label{corHBext}
Let $E$ be a real locally convex space
and $q\colon E\to[0,\infty[$ be a continuous seminorm.
For each $z\in E$, there exists a continuous linear functional
$\lambda\colon E\to\R$ such that
\[
\mbox{$\lambda(z)=q(z)$ and $|\lambda(x)|\leq q(x)$ for all $x\in E$.}
\]
\end{cor}
\begin{prf}
The linear functional $\lambda_0\colon \R z\to\R$, $tz\mto tq(z)$ for $t\in \R$
satisfies $\lambda_0(z)=q(z)$ and $|\lambda_0(tz)|=|t|q(z)=q(tz)\leq q(tz)$
for all $t$, i.e.,
\[ |\lambda_0(x)|\leq q(x) \quad \mbox{  for all } \quad x\in \R z.\] 
The Hahn--Banach Extension Theorem (Theorem~\ref{HBext}) provides
$\lambda\in E'$ such that $|\lambda(x)|\leq q(x)$ for all $x\in E$ and $\lambda|_{\R z}=\lambda_0|_{\R z}$,
whence $\lambda(z)=\lambda_0(z)=q(z)$ in particular.
\end{prf}
\subsection*{Exercises for Appendix~\ref{sechbanach}}
\begin{small}
\begin{exer}
Prove formula (\ref{formbipol}) in Remark~\ref{rembipol}.
\end{exer}
\begin{exer}\label{exc-annih}
Let $E$ be a locally convex space. Given a subset $M\sub E$, we call
$M^\perp:=\{\lambda\in E'\colon (\forall x\in M)\;\lambda(x)=0\}$
its {\em annihilator} in~$E'$.
Given a subset $N\sub E'$, we call $N_\perp:=\{x\in E\colon (\forall \lambda\in N)\;\lambda(x)=0\}$
its annihilator in~$E$.
\begin{description}[(D)]
\item[(a)]
Show that $M^\perp$ is a vector subspace of~$E'$.
Show that $N_\perp$ is a closed vector subspace of~$E$.
\item[(b)]
Show that $M^\perp=M^\circ$ if $M$ is a vector subspace of~$E$.
Likewise, $N_\perp=N_\circ$ if $N$ is a vector subspace of~$E'$.
\item[(c)]
Using (\ref{formbipol}) in Remark~\ref{rembipol},
deduce that $M=(M^\perp)_\perp$
for each closed vector subspace $M\sub E$.
\end{description}
\end{exer}
\end{small}
\section{Completeness and sequential completeness}\label{seccomplete}
We now study completeness properties
of topological vector spaces, which are useful
for analysis as they can frequently be used
to show that certain limits of interest exist.
Although the results
concerning sequences are
usually sufficient for our
purposes, it is more natural
to discuss convergence of nets
as well.
The definitions of nets
and convergence of nets
can be looked up in Appendix~A.
\begin{defn}\label{defnCauchy}
Let $E$ be a topological $\K$-vector space.
\begin{description}[(D)]
\item[\rm (a)]
A sequence $(x_n)_{n\in \N}$ in
$E$ is called a {\em Cauchy sequence\/}
if, for every zero-neighborhood
$U$ in~$E$, there exists $N\in \N$
such that
\[
x_n-x_m\, \in \, U\qquad \mbox{for all $\, n,m\geq N$.}
\]
\item[\rm (b)]
The space~$E$ is called
{\em sequentially complete\/} if every Cauchy sequence $(x_n)_{n\in \N}$
in~$E$ converges.
\item[\rm (c)]
A net $(x_\alpha)_{\alpha\in A}$ in~$E$ is
called a {\em Cauchy net\/}
if, for every $0$-neighborhood~$U$ in~$E$,
there exists $\alpha_0\in A$
such that $x_\alpha -x_\beta\in U$
for all $\alpha,\beta\geq \alpha_0$.
\item[\rm (d)]
$E$ is called {\em complete\/}
if every Cauchy net in~$E$ converges.
\end{description}
\end{defn}
Note that
a sequence $(x_n)_{n\in \N}$ in $E$
is a Cauchy sequence if and only if
it is a Cauchy net.
Consequently, every complete
topological vector space~$E$ also is
sequentially complete.

If $E$ is a topological vector space, $F\sub E$ a vector subspace
and $(x_\alpha)_{\alpha\in A}$ a net in~$F$,
then $(x_\alpha)_{\alpha\in A}$ is a Cauchy net in~$F$
if and only if it is a Cauchy net in~$E$ (as is clear from the definition of the induced topology).
We infer:
\begin{lem}\label{compl-clo}
If $E$ is a topological vector space and $F\sub E$ a vector subspace
which is complete in the induced topology,
then~$F$ is closed in~$E$.
\end{lem}
\begin{prf}
If $x\in E$ is in the closure of~$F$,
then there exists a net $(x_\alpha)_{\alpha\in A}$ in~$F$
which converges to~$x$.
Thus $(x_\alpha)_{\alpha\in A}$ is a Cauchy net in~$E$
and hence also in~$F$, whence it converges in~$F$ (and hence in~$E$) to some
$y\in F$. By uniqueness of limits of nets in Hausdorff spaces,
we have $x=y\in F$.
\end{prf}
\begin{prop}\label{complmetr}
A metrizable topological vector space $E$
is complete if and only if it is sequentially complete.
\end{prop}
\begin{prf}
As the converse is clear, we only need to consider a metrizable, sequentially
complete topological vector space~$E$ and show that it is complete.\\[3mm]
Since $E$ is metrizable, there exists
a sequence $(U_n)_{n\in \N}$
of $0$-neighborhoods $U_n\sub E$
which form a basis of $0$-neighborhoods.
After replacing $U_n$ with $\bigcap_{k\leq n}U_k$,
we may assume that
$U_1\supseteq U_2\supseteq\cdots$.\\[3mm]
Let $(x_\alpha)_{\alpha\in A}$
be a Cauchy net in~$E$.
Since~$(x_\alpha)$ is
a Cauchy net,
there exists
an element $\alpha_1\in A$
with $x_\alpha-x_\beta\in U_1$
for all $\alpha,\beta\geq \alpha_1$.
Next, we find $\alpha_2\in A$
such that
$x_\alpha-x_\beta\in U_2$
for all $\alpha,\beta\geq \alpha_2$.
Note that, since $A$ is directed, there is $\alpha'_2\in A$
such that $\alpha_2'\geq \alpha_1$ and $\alpha_2'\geq \alpha_2$.
Replacing $\alpha_2$ with $\alpha_2'$,
we may assume that $\alpha_2\geq \alpha_1$.
Proceeding in this way, we find an increasing sequence
$\alpha_1\leq \alpha_2\leq \alpha_3\leq \cdots$
in $A$ such that, for each $n\in \N$,
\[
x_\alpha-x_\beta \, \in \, U_n\qquad
\mbox{for all $\, \alpha,\beta\geq \alpha_n$.}
\]
For all $k,m\geq n$, we have
$\alpha_k,\alpha_m\geq \alpha_n$
and hence
$x_{\alpha_k}-x_{\alpha_m}\in U_n$;
thus $(x_{\alpha_n})_{n\in\N}$ is a Cauchy sequence
and therefore convergent to some~$x\in E$ by hypothesis.
Let~$U$ be any neighborhood of~$x$ in~$E$.
Then
there exists some $n\in \N$ such that
$U_n +x \sub U$.
As a consequence of Lemma~\ref{verybase}(d),
there exists
some $m\in \N$ with $U_m+U_m\sub U_n$.
The sequence $(x_{\alpha_k})_{k\in\N}$
being convergent to~$x$, there exists
a natural number $\ell \geq m$ with
$x_{\alpha_k}\in U_m+x$ for all $k\geq \ell$.
For every $\alpha\geq \alpha_\ell$,
we obtain
\[
x_\alpha \; =\; (x_\alpha-x_{\alpha_\ell})+ x_{\alpha_{\ell}}
\; \in \; U_\ell+U_m+x\;\sub\;
U_m+U_m+x\; \sub \; U_n+x\sub U \, .
\]
Thus the net $(x_\alpha)$ converges to~$x$.
\end{prf}
For example, every Banach space is a complete
locally convex space by the preceding proposition.
Complete, metrizable locally
convex spaces are encountered
frequently and deserve their own name.
\begin{defn}\label{defnFrech}
A complete, metrizable locally convex
space is called a {\em Fr\'{e}chet space.}
\end{defn}
\begin{lem}\label{lemcauch}
If $\lambda\!: E\to F$ is a continuous linear
map between topological vector spaces and $(x_\alpha)_{\alpha\in A}$
a Cauchy net in $E$, then $(\lambda(x_\alpha))_{\alpha\in A}$
is a Cauchy net in~$F$.
\end{lem}
\begin{prf}
In fact, given a $0$-neighborhood
$U\sub F$, due to the continuity of~$\lambda$
the pre-image $\lambda^{-1}(U)$ is a $0$-neighborhood
in~$E$. Now $(x_\alpha)_{\alpha\in A}$
being a Cauchy net, there exists $\alpha_0\in A$
such that $x_\alpha-x_\beta\in \lambda^{-1}(U)$
for all $\alpha,\beta\geq \alpha_0$.
Then, for all $\alpha,\beta\geq \alpha_0$, we have
\[
U\; \supseteq \; \lambda(\lambda^{-1}(U))\; \ni\;
\lambda(x_\alpha-x_\beta)\; =\; \lambda(x_\alpha)-\lambda(x_\beta)\,,
\]
whence indeed $(\lambda(x_\alpha))_{\alpha\in A}$ is a Cauchy net.
\end{prf}
\begin{defn}\label{defseqclo}
Let $X$ be a Hausdorff topological space.
A subset $A\sub X$ is called
{\em sequentially
closed\/} if $\lim_{n\to\infty} a_n\in A$,
for each sequence $(a_n)_{n\in \N}$ in~$A$ which converges
in~$X$.
\end{defn}
Note that if
$A\sub X$ is closed, then $A$ is sequentially closed.\\[3mm]
The following constructions preserve completeness.
\begin{prop}\label{newcompl}
\begin{description}[(D)]
\item[\rm (a)]
If $(E_i)_{i\in I}$
is a family of complete $($resp., sequentially complete$)$
topological vector spaces, then also the direct product
$P:=\prod_{i\in I} E_i$ is complete
$($resp., sequentially complete$)$.
\item[\rm (b)]
Let $E$ be a complete $($resp., sequentially complete$)$
topological vector space
and $F\sub E$ be a closed $($resp., sequentially
closed$)$ vector subspace.
Then also $F$ is complete
$($resp., sequentially complete$)$.
\end{description}
\end{prop}
\begin{prf}
We only discuss completeness;
to prove the analogues
concerning sequential completeness,
simply replace nets by sequences.

(a) Let $(x_\alpha)_{\alpha\in A}$
be a Cauchy net in $P$. Given $i\in I$,
the canonical projection
$\pr_i\!: P\to E_i$ is continuous linear
map, and thus
$x_{i,\alpha}:=\pr_i(x_\alpha)$
is a Cauchy net in $E_i$ (Lemma~\ref{lemcauch})
and thus convergent to some $y_i\in E_i$
as $E_i$ is assumed complete.
Then $(x_\alpha)_{\alpha\in A}$
converges to
$y:=(x_i)_{i\in I}\in P$,
since $x_{i,\alpha}\to y_i$
for each $i$ (see Lemma~\ref{ctsinprod}(b)).

(b) If $(x_\alpha)_{\alpha\in A}$ is a Cauchy net
in $F$, then $(x_\alpha)_{\alpha\in A}$
also is a Cauchy net in $E$ (apply Lemma~\ref{lemcauch}
to the continuous linear inclusion map $F\to E$)
and hence converges in~$E$, to $x\in E$
say. Now $F$ being closed in $E$,
we deduce
from $\{x_\alpha\!: \alpha\in A\} \sub F$
that $x\in \wb{F}=F$ (see Proposition~\ref{netclosure}).
Then $x_\alpha\to x$ in $F$,
equipped with the topology induced by~$E$ (exercise).
\end{prf}
\begin{prop}\label{onlyyou}
On each finite-dimensional $\K$-vector space~$E$,
there is exactly one topology~$\cO$ making it a $($Hausdorff$)$ topological
$\K$-vector space. Every norm on~$E$ defines
this vector topology, and $(E,\cO)$ is complete.
\end{prop}
\begin{prf}
The proof is by induction on the dimension~$n$ of~$E$. Let $\cO$
be a vector topology on~$E$. It suffices to find an isomorphism
\[
\phi\colon \K^n\to (E,\cO)
\]
of topological vector spaces, where $\K^n$ is endowed with the product topology.
In fact, if also $\cO'$ is a vector topology on~$E$ and $\psi\colon \K^n\to (E,\cO')$
is an isomorphism of topological vector spaces, then $\psi^{-1}\circ \id_E\circ \phi$
is an automorphism of the vector space $\K^n$ and hence a homeomorphism
(as is known from elementary calculus). As a consequence, $\id_E\colon (E,\cO)\to (E,\cO')$
is a homeomorphism and thus $\cO=\cO'$.

If $n=1$, choose $0\not=v\in E$.
Then the map $h\colon\K\to E$, $z\mto zv$ is an isomorphism of vector spaces
and continuous.
If $r>0$, then there is a balanced open $0$-neighborhood $W\sub E$ such that $rv\not\in W$.
Hence $h^{-1}(W)$ is a balanced open subset of $\K$
such that $r\not\in h^{-1}(W)$, and thus $h^{-1}(W)=\bD_s^0$
for some $s\in\;]0,r]$.
As a consequence, $h(\bD_r^0)\supseteq W$, whence $h(\bD_r^0)$ is a $0$-neighborhood.
Thus~$h$ is an open map (by Proposition~\ref{ctsopen}(b)) and hence~$h$
is a homeomorphism.

If $n>1$, let $F\sub E$ be a vector subspace of dimension $n-1$. Give~$F$ the topology induced by~$E$.
By induction, $F$ is complete and we find an isomorphism $\phi\colon \K^{n-1}\to F$ of topological vector spaces.
Being complete, $F$ is closed in~$E$, whence $E/F$ is a (Hausdorff) topological vector space.
Let $q\colon E\to E/F$, $x\mto x+F$ be the canonical quotient map.
Since $E/F$ has dimension~$1$, we find an isomorphism $\theta\colon\K\to E/F$
of topological vector spaces. Choose $v\in E$ such that $q(v)=\theta(1)$.
By uniqueness of the vector topology on $1$-dimensional vector spaces, the vector space isomorphism
$q|_{\K v}\colon \K v\to E/F$ is a homeomorphism. Now
\[
\psi\colon\K^n\to E,\quad (z_1,\ldots, z_n)\mto \phi(z_1,\ldots, z_{n-1})+z_n v=:w
\]
is continuous and an isomorphism of vector spaces. Since $z_n=\theta^{-1}(q(w))$
and $(z_1,\ldots, z_{n-1})=\phi^{-1}(w-z_nv)$ depend continuously on~$w$,
we see that $\psi^{-1}$ is continuous. Hence $\psi$ is an isomorphism of
topological vector spaces.
\end{prf}
\begin{prop}\label{lcpfindim}
Every locally compact topological vector space
has
finite dimension.
\end{prop}
\begin{prf}
If $E$ is locally compact, let $K\sub E$ be a compact $0$-neighborhood.
Since $2 K$ is covered by translates of $K^0$,
we have
\begin{equation}\label{givesall}
2 K\sub \Phi+ K
\end{equation}
for a finite subset $\Phi\sub E$. Then $F:=\Spann(\Phi)$ has finite dimension.
By Proposition~\ref{onlyyou}, $F$ is complete and
hence closed in~$E$, whence $E/F$ is Hausdorff.
Let $q\colon E\to E/F$ be the quotient map. By (\ref{givesall}), we have
\[
\mbox{$q(2K)\sub q(K)\;$ and thus $\;(\forall n\in\N_0)\;\;2^n q(K)\sub q(K)\sub 2^{-n}q(K)$.}
\]
Hence
$E/F=\bigcup_{n\in\N_0}2^nq(K)=q(K)=\bigcap_{n\in\N_0}2^{-n}q(K)=\{0\}$,
using that $q(K)$ is absorbing and bounded. Thus $E/F=\{0\}$ and hence~$E=F$.
\end{prf}
\begin{rem}
If $E$ is a complete locally convex space and $F\sub E$ a closed vector subspace,
then $E/F$ need not be complete (see Example~\ref{bad-quot}
for an explicit counterexample). The pathology cannot occur if~$E$
is metrizable (see Lemma~\ref{extprp}(c)).
\end{rem}
\section{The completion of a locally convex space}\label{seccomption}
It is useful that every locally convex space
can be completed.
\begin{defn}\label{defcomption}
A {\em completion\/} of a locally convex space
$E$ is a complete locally convex space
$\tilde{E}$, together with a linear
topological embedding\linebreak
$\kappa_E\!: E\to \tilde{E}$
with dense image.
\end{defn}
If a completion exists,
then standard set theoretic arguments
can be used to manufacture a completion $(\tilde{E},\kappa_E)$
such that~$E$ is a vector subspace of~$\tilde{E}$
and $\kappa_E$ the inclusion map, $x\mto x$.
\begin{prop}\label{cmplexists}
Every locally convex space~$E$ has a completion.
\end{prop}
\begin{prf}
Let $\cP$ be the set of all
continuous seminorms on~$E$.
For each $p\in \cP$,
we let
$\wt{E_p}$ be the Banach space
obtained by completing the normed
space $(E_p,\|\cdot\|_p)$,
where $E_p:=E/N_p$ with $N_p:=\{x\in E\!: p(x)=0\}$
and the norm is defined via
$\|x+N_p\|_p:=p(x)$ for $x\in E$.
We may assume that $E_p\sub \wt{E_p}$.
For each $p\in \cP$,
let $\alpha_p\!: E\to \wt{E_p}$, $\alpha_p(x):=x+N_p$
be the natural map.
We set
$\alpha:=(\alpha_p)_{p\in \cP}\!: E\to\prod_{p\in \cP}\wt{E_p}=:P$.
Being a direct product of complete locally
convex spaces, $P$ is complete
(Proposition~\ref{newcompl}(a)).
Then also the closure
$\tilde{E}:=\wb{\alpha(E)}$
of the image of~$\alpha$ is a
complete locally convex space,
being a closed vector subspace of the complete
locally convex space~$P$
(Proposition~\ref{newcompl}\,(b)).
We let $\kappa_E:=\alpha|^{\tilde{E}}\!: E\to \tilde{E}$
be the co-restriction of $\alpha$ to $\tilde{E}$.
As a consequence of
Proposition~\ref{embedprod},
the linear map~$\alpha$,
and hence also $\kappa_E$, is a topological
embedding.
The image of $\kappa_E$ being dense in $\tilde{E}$ by construction,
$(\tilde{E},\kappa_E)$ is a completion of~$E$.
\end{prf}
Completions
are unique up to canonical isomorphism
(see Remark~\ref{remuniqcomp}\,(c)).
\section{Continuous extension of linear maps}\label{secctsex}
\begin{prop}\label{extendlin}
Let $E$ be a locally convex space,
$E_0\sub E$ be a dense vector subspace,
$F$ a complete locally convex space, and
$\lambda\!: E_0\to F$ be a continuous linear map.
Then there exists a unique continuous linear map
\[
\Lambda\!: E\to F
\]
which extends $\lambda$, i.e., $\Lambda|_{E_0}=\lambda$.
\end{prop}
\begin{prf}
{\bf Step~1:} {\em Definition of $\Lambda$.}
Given $x\in E$,
Proposition~\ref{netclosure}
provides a net $(x_\alpha)_{\alpha\in A}$
in $E_0$ converging to~$x$ in~$E$, since $\wb{E_0}=E$.
By Lemma~\ref{lemcauch}, $(\lambda(x_\alpha))_{\alpha\in A}$
is a Cauchy net in $F$ and thus
convergent, as $F$ is assumed complete.
We want to define
\begin{equation}\label{badlydef}
\Lambda(x)\, :=\, \lim\, \lambda (x_\alpha)
\end{equation}
as the limit of $(\lambda (x_\alpha))_{\alpha\in A}$.
To see that this definition is independent of the choice
of net $(x_\alpha)_{\alpha\in A}$, let
also $(y_\beta)_{\beta\in B}$ be a net in $E_0$
converging to~$x$ in~$E$. Let $W\sub F$ be an open $0$-neighborhood.
Then $\lambda^{-1}(W)$ is an open
$0$-neighborhood in~$E_0$.
As $E_0$ is equipped with the induced
topology, we find an open neighborhood
$V\sub E$ of $0$ such that $V\cap E_0= \lambda^{-1}(W)$.
Let $U\sub E$ be a $0$-neighborhood
such that $U-U\sub V$. We find $\alpha_0\in A$
and $\beta_0\in B$ such that $x-x_\alpha\in U$
and $x-y_\beta\in U$, for all $A\ni
\alpha\geq \alpha_0$
and $B\ni \beta\geq \beta_0$.
Noting that $E_0$
is a vector subspace, we obtain
\[
x_\alpha-y_\beta\,=\, x-y_\beta-(x-x_\alpha)\, \in \, (U-U)\cap E_0
\, \sub \, V\cap E_0\, =\, \lambda^{-1}(W)
\]
for any such $\alpha$ and $\beta$
and thus $\lambda(x_\alpha)-\lambda(y_\beta)=\lambda(x_\alpha
-y_\beta)\in W$. Passing to the limit in $\alpha$ first
and then to the limit in $\beta$, we deduce that
\[
\lim\, \lambda(x_\alpha)-\lim\, \lambda(y_\beta)\, \in \, \wb{W}\,.
\]
Since sets of the form $\wb{W}$ constitute a basis
of $0$-neighborhoods and~$F$ is Hausdorff,
we deduce that
$\lim\, \lambda(x_\alpha)-\lim\, \lambda(y_\beta)=0$,
whence a mapping\linebreak
$\Lambda\!: E\to F$ is well defined by (\ref{badlydef}).\\[3mm]
{\bf Step~2:} {\em $\Lambda$ is linear.}
In fact, given $x,y\in E$ and $z\in \K$, let
$(x_\alpha)_{\alpha\in A}$ and
$(y_\beta)_{\beta\in B}$ be nets in~$E_0$
converging to $x$ and $y$ in~$E$, respectively.
Then $C:=A\times B$ becomes a directed set by declaring
$(\alpha_1,\beta_1)\leq (\alpha_2,\beta_2)$
for $\alpha_1,\alpha_2\in A$, $\beta_1,\beta_2\in B$
if and only if $\alpha_1\leq \alpha_2$ and $\beta_1\leq \beta_2$
(exercise). Furthermore, the net $(x_\alpha)_{(\alpha,\beta)\in C}$
converges to~$x$ and $(y_\beta)_{(\alpha,\beta)\in C}$
converges to~$y$ (exercise).
Using these nets instead of $(x_\alpha)_{\alpha\in A}$
and $(y_\beta)_{\beta\in B}$,
we may assume without loss of generality that $(A,\leq)=(B,\leq)$.
Then $(x_\alpha+zy_\alpha)_{\alpha\in A}$ is easily
seen to be a Cauchy net in~$E_0$
converging to $x+zy$ in~$E$, and thus
\begin{eqnarray*}
\Lambda(x+zy) & = & \lim\,\lambda(x_\alpha+zy_\alpha)\, =\,
\lim\, (\lambda(x_\alpha)+z\lambda(y_\alpha))\\
& = &\lim\, \lambda(x_\alpha)+z \,\lim\,\lambda(y_\alpha)
\, =\, \Lambda(x)+z\Lambda(y)\, ,
\end{eqnarray*}
using Proposition~\ref{exemplify} and
Lemma~\ref{ctsinprod}\,(a)
(exploiting that addition and scalar multiplication by $z$
are continuous maps).
Thus $\Lambda$ is indeed linear.\\[3mm]
{\bf Step~3:} {\em $\Lambda$ is continuous.}
Given any $0$-neighborhood $Q\sub F$, there exists
an open $0$-neighborhood
$W\sub F$ whose closure $\wb{W}$ is contained in~$Q$.
Then $V:=\lambda^{-1}(W)$ is an open $0$-neighborhood
in~$E_0$. Let $U\sub E$ be an open $0$-neighborhood
in $E$ such that $U\cap E_0=V$.
Since $E_0$ is dense is $E$, the set
$U\cap E_0=V$ is dense in $U$, and thus
$\wb{U}=\wb{V}$, the closure of $V$ in~$E$.
Given $x\in U$, we have $x\in \wb{V}$,
whence there exists a net $(x_\alpha)_{\alpha\in A}$
in~$V$ such that $x_\alpha\to x$ in~$E$.
Since $\lambda(x_\alpha)\in W\sub \wb{W}$
for each $\alpha$, we deduce from Proposition~\ref{netclosure}
that
\[
\Lambda(x)\; =\; \lim\, \lambda(x_\alpha)\; \in \; \wb{W}\; \sub\; Q\,.
\]
Thus $\Lambda(U)\sub Q$, showing that the linear map
$\Lambda$ is continuous at~$0$ and thus
continuous. Since $E_0$ is dense in~$E$, the continuous map
$\Lambda$ extending $\lambda$ is uniquely determined.
\end{prf}
\begin{rem}\label{remuniqcomp}
Let $E$ be a locally convex space and $\wt{E}$
be a completion of~$E$, together
with the canonical embedding $\kappa\colon E\to \wt{E}$.
We may assume that $E\sub \wt{E}$ and $\kappa(x)=x$.
\begin{description}[(D)]
\item[(a)]
As $\K$ is a complete locally convex space,
Proposition~\ref{extendlin} shows
that every continuous linear functional
$\lambda\colon E\to\K$
extends uniquely to a continuous
linear functional $\wt{E}\to \K$.
Conversely, $\kappa'(\lambda):=\lambda\circ \kappa=\lambda|_E$ is a
continuous linear functional
on~$E$, for each $\lambda\in (\wt{E})'$.
Hence, as it is also linear, the bijection
\[
\kappa'\colon (\wt{E})'\to E'\,,\qquad
\lambda\mto \kappa'(\lambda)=\lambda\circ \kappa
\]
is an isomorphism of vector spaces.
\item[(b)]
Likewise, each continuous linear map
$E\to F$ from~$E$ to a complete
locally convex space~$F$ extends uniquely
to a continuous linear map $\wt{E}\to F$.
This means
that the natural linear map
\[
\cL(\kappa, F)\colon \cL(\wt{E}, F)\to \cL(E,F)\,,\quad
\lambda\mto \lambda\circ \kappa
\]
between the indicated
spaces of continuous linear maps
is an isomorphism of vector spaces.
\item[(c)]
In particular, for each completion
$\gamma \colon E\to F$ of~$E$
there are unique continuous linear maps
$\gamma_1 \colon \wt{E}\to F$
and $\kappa_1 \colon F\to \wt{E}$
such that $\gamma_1 \circ  \kappa=\gamma$
and $\kappa_1 \circ \gamma =\kappa$.
Then $(\gamma_1 \circ \kappa_1 )\circ \gamma
=\gamma_1 \circ \kappa =\gamma=\id_F \circ \gamma$
implies that $\gamma_1\circ\kappa_1 =\id_F$,
and likewise $\kappa_1\circ \gamma_1=\id_{\wt{E}}$.
Hence $\gamma_1 \colon \wt{E}\to F$
and $\kappa_1 \colon F\to \wt{E}$
are mutually inverse isomorphisms
of topological vector spaces. Thus
$\wt{E}$ and $F$ are naturally isomorphic.
\end{description}
\end{rem}
\subsection*{Exercises for Appendix~\ref{secctsex}}
\begin{small}
\begin{exer}\label{sep-points-failure}
Let $E$ be a locally convex space and
$\wt{E}$ be a completion of~$E$ such that $E\sub\wt{E}$.
Then each $\lambda\in E'$ has a unique extension to
a continuous linear functonal $\tilde{\lambda}$ on~$\wt{E}$.
We show that if a subset $\Lambda\sub E$ separates points on~$E$,
then $\{\tilde{\lambda}\colon \lambda\in \Lambda\}\sub (\wt{E})'$
need not separate points on~$E$.\\[1mm]
To get an example, let $F$ be a complete
locally convex space which has a proper dense vector
subspace~$E$; let $0\not=v\in F\setminus E$ (e.g., $F:=\ell^1$, $E\sub F$
the subspace of finite sequences, $v=(1/n^2)_{n\in\N}$).
Thus~$F$ is a completion of~$E$.
Then $M:=\R v$ is a closed vector subspace of~$F$ (cf.\ Proposition~\ref{onlyyou}),
whence $(M^\perp)_\perp=M$ (see Exercise~\ref{exc-annih}). Show that $\Lambda:=\{\lambda|_E\colon\lambda\in M^\perp\}$\vspace{-.3mm}
separates points on~$E$ but $M^\perp=\{\tilde{\lambda}\colon \lambda\in\Lambda\}$
fails to separate points on $F=\wt{E}$.
\end{exer}

\begin{exer}\label{fct-space-completion}
Let $X$ be a Hausdorff topological space,
$E$ be a locally convex space, and $\tilde{E}$ be a completion
of~$E$ such that $E\sub\tilde{E}$.
Write $k(X)$ for~$X$, endowed with the final topology with respect
to the inclusion maps $i_K\colon K\to X$ for $K$ in the set of compact subsets of~$X$.
Then the map $k(X)\to X$, $x\mto x$ is continuous.
\begin{description}[(D)]
\item[(a)]
Show that a subset $K\sub X$ is compact in~$X$
if and only if it is compact in~$k(X)$. Deduce that $k(X)$
is a $k$-space.
For each Hausdorff topological space~$Y$,
show that the inclusion map $C(X,Y)\to C(k(X),Y)$
is a topological embedding.
\item[(b)]
Show that the closure of $C(X,E)$ in $C(k(X),\tilde{E})$
is a completion $C(X,E)\,\tilde{\;}$ of $C(X,E)$.
\item[(c)]
The point evaluations $\ve_x\colon C(k(X),\tilde{E})$
at elements $x\in X$ separate points on $C(X,E)\,\tilde{\;}\sub C(k(X),\tilde{E})$,
and hence so do the maps $\lambda\circ\ve_x$ with $x\in X$ and $\lambda\in (\tilde{E})'$.
The pathology described in example~\ref{sep-points-failure}
therefore does not occur for
$\Lambda:=\{\lambda\circ \ve_x\colon x\in X,\lambda\in E'\}\sub C(X,E)'$
writing $\ve_x$ also for the point evaluation on $C(X,E)$;
this can be useful in connection with $C(X,E)$-valued weak integrals.
\end{description}
\end{exer}
\end{small}
\section{Bounded sets and the weak topology}\label{secbdd}
\begin{defn}\label{defnbdd}
A subset $B\sub E$ of a topological vector space
$E$ is called {\em bounded\/}
if it
is absorbed by each $0$-neighborhood,
i.e., for each neighborhood 
$U\sub E$ of $0$, there exists $r>0$ such that $B\sub rU$.
\end{defn}
\begin{lem}\label{lemboundd}
If $\alpha\colon E\to F$ is a continuous
linear map between topological
vector spaces and $B\sub E$ is bounded,
then $\alpha(B)\sub F$ is bounded.
\end{lem}
\begin{prf}
If $U$ is a $0$-neighborhood in~$F$,
then $\alpha^{-1}(U)$ is a $0$-neighborhood in~$E$.
Since~$B$ is bounded, we find
$r>0$ such that $B\sub r\alpha^{-1}(U)=\alpha^{-1}(rU)$.
Consequently, $\alpha(B) \sub rU$.
\end{prf}
\begin{lem}\label{lembdsub}
If $E$ is a topological vector space,
$F\sub E$ a vector subspace
and $B\sub F$, then $B$ is bounded
in~$F$ if and only if $B$ is bounded in~$E$.
\end{lem}
\begin{prf}
If $B$ is bounded in~$F$,
then $B=\lambda(B)$ is bounded in~$E$,
the inclusion map $\lambda\colon F\to E$, $\lambda(x):=x$
being continuous linear (Lemma~\ref{lemboundd}).
Conversely, assume that $B$ is bounded in~$E$,
and let $U\sub F$ be a $0$-neighborhood.
Then there exists a $0$-neighborhood
$V\sub E$ such that $F\cap V \sub U$.
The set~$B$ being bounded in~$E$, we find $r>0$ such
that $B\sub rV$. Thus $B\sub F\cap rV=r(F\cap V)\sub rU$,
showing that $B$ is bounded in~$F$.
\end{prf}
\begin{lem}\label{lembdprod}
Let $(E_j)_{j\in J}$
be a family of topological vector
spaces and $E:=\prod_{j\in J}E_j$
be their direct
product, with canonical projections
$\pr_j\colon E\to E_j$.
Then a subset $B\sub E$ is bounded
if and only if $\pr_j(B)$ is bounded
in~$E_j$ for each $j\in J$.
\end{lem}
\begin{prf}
If $B$ is bounded, then so is $\pr_j(B)$
for each $j\in J$, by Lemma~\ref{lemboundd}
(the projection $\pr_j$ being continuous linear).
Conversely, suppose that $\pr_j(B)$ is bounded
for each~$j\in J$.
If $U\sub E$ is a $0$-neighborhood,
then there exists a finite subset
$F\sub J$ and $0$-neighborhoods
$U_j\sub E_j$ such that
\[\bigcap_{j\in F}\pr_j^{-1}(U_j)\sub U.\] 
After shrinking $U_j$,
we may assume that each $U_j$
is balanced.
Since $\pr_j(B)$ is bounded,
we find $r_j>0$ such that $\pr_j(B)\sub r_jU_j$.
Set $r:=\max\{r_j\colon j\in F\}$.
Then $\pr_j(B)\sub r_jU_j\sub rU_j$
and hence $B\sub \pr_j^{-1}(rU_j)=r\pr_j^{-1}(U_j)$
for each $j\in F$,
using that $U_j$
is balanced.
Thus $B\sub r \bigcap_{j\in F}\pr_j^{-1}(U_j)\sub r U$,
showing that $B$ is bounded.
\end{prf}
\begin{lem}\label{lembdsmn}
If $E$ is a locally convex space,
then a subset $B\sub E$ is bounded
if and only if $p(B)$ is a bounded
subset of~$\K$,
for each continuous seminorm~$p$ on~$E$.
\end{lem}
\begin{prf}
If $B$ is bounded and
$p\colon E\to [0,\infty[$ is a continuous seminorm,
then the ball $U:=\{x\in E\colon p(x)\leq 1\}$
is a $0$-neighborhood,
whence $B\sub rU$
for some $r>0$ and thus $p(B)\sub p(rU)=rp(U)\sub [0,r]$.
Conversely, suppose that $p(B)$ is bounded
for each continuous seminorm~$p$.
If $U\sub E$ is a $0$-neighborhood,
Lemma~\ref{basisnice} and Proposition~\ref{Minkowsk}
provide a continuous seminorm $q$
such that $V:=\{x\in E\colon q(x)\leq 1\}\sub U$.
By hypotheses, $q(B)\sub [0,r]$ for some
$r>0$, whence $q(r^{-1}B)\sub [0,1]$
and thus $r^{-1}B\sub V$
and so $B\sub rV\sub rU$.
\end{prf}
\begin{defn}\label{defweaktop}
The {\em weak topology\/}
on a locally convex space~$E$ is the
initial topology
with respect to the set $E'$
of all continuous linear functionals,
i.e., the topology $\cO_w$
which makes
\begin{equation}\label{embedweak}
\eta\colon (E,\cO_w)\to\prod_{\lambda\in E'}\K\; =\; \K^{E'}\,,
\qquad
x\mto (\lambda(x))_{\lambda\in E'}
\end{equation}
a topological embedding, where
the right hand side is equipped with the
product topology.
We write $E_w$ for~$E$,
equipped with the weak topology $\cO_w$.
The product on the right hand side of (\ref{embedweak})
being a locally
convex space, so is $E_w$.
A subset of $E$ is called {\em weakly open\/}
({\em weakly closed\/},
resp., {\em weakly bounded\/}) if it is open (closed, resp., bounded)
in~$E_w$.
\end{defn}
Note that each $\lambda\in E'$
is continuous on $E_w$,
by definition of the weak topology;
hence $E'\sub (E_w)'$.
The weak topology being coarser
than the original topology,
we also have
$E' \supseteq (E_w)'$
and thus
\begin{equation}\label{weakdual}
E'\; =\; (E_w)'\, .
\end{equation}
\begin{lem}\label{lemweakbd}
Let $E$ be a locally convex space.
Then a subset $B\sub E$ is weakly
bounded if and only if $\lambda(B)$
is a bounded subset of~$\K$,
for each $\lambda \in E'$.
\end{lem}
\begin{prf}
If $B\sub E_w$ is bounded,
then $\lambda(B)\sub \K$
is bounded for each $\lambda\in E'=(E_w)'$,
by Lemma~\ref{lemboundd}.
If, conversely,
$\lambda(B)$ is bounded in~$\K$
for each $\lambda\in E'$,
then $\ev_\lambda (\eta(B))=\lambda(B)$
is bounded for each $\lambda\in E'$
(where $\eta\colon E\to \K^{E'}$
is as in (\ref{embedweak}),
and $\ev_\lambda=\pr_\lambda\colon
\K^{E'}\to \K$,
$\xi\mto \xi(\lambda)$
the evaluation map).
Therefore $\eta(B)$ is bounded in $\K^{E'}$
by Lemma~\ref{lembdprod}.
Since $\eta(B)\sub \eta(E)$,
this entails that $\eta(B)$ is bounded
in~$\eta(E)$.
Therefore $B=(\eta|^{\eta(E)})^{-1}(\eta(B))$
is bounded in~$E_w$ (by Lemma~\ref{lemboundd}).
\end{prf}
\begin{thm}[{\bf Mackey's Theorem}]\label{Mackey}
A subset $B$ of a locally convex space~$E$
is bounded if and only if it is weakly bounded.
\end{thm}
\begin{prf}
The map $\Lambda\colon E \to E_w$, $x\mto x$ being
continuous, boundedness of $B$ in~$E$ entails
boundedness of $B=\Lambda(B)$ in $E_w$.
For the converse, assume that $B\sub E$ is weakly
bounded.
Let $\cP$ be the set of all continuous seminorms
on~$E$. By Lemma~\ref{lembdsmn}, $B$ will
be bounded if we can show that
$p(B)$ is bounded in~$\R$ for each
$p\in \cP$.
We define
$E_p:=E/N_p$
and let $\| .\|_p$
be the norm on $E_p$
determined by $\|\alpha_p(x)\|_p=p(x)$,
as in~\ref{modout}.
We choose a Banach space
$(\wt{E_p},\|\cdot\|_p^{\wt{\;}})$
which is a completion
of the normed space $(E_p,\|\cdot\|_p)$,
i.e., we are also given a linear isometry
$\gamma_p\colon E_p\to \wt{E_p}$
with dense image.
Then $\beta_p:=\gamma_p\circ \alpha_p\colon E\to \wt{E_p}$ 
is a continuous linear map
such that
$\|\beta_p(x)\|_p^{\wt{\;}}=p(x)$
for each $x\in E$.
Hence $B$ will be bounded
if we can show that $\|\beta_p(B)\|_p^{\wt{\;}}=p(B)$
is bounded, for each $p\in \cP$.
Note that, for each continuous linear
functional $\lambda\colon \wt{E}_p \to \K$,
the composition $\lambda\circ \beta_p$
is a continuous linear functional on~$E$,
whence $\lambda(\beta_p(B))=(\lambda\circ \beta_p)(B)\sub \K$
is bounded.
Hence $\beta_p(B)$ is weakly bounded
in $\wt{E_p}$,
and it only remains to show that $\beta_p(B)$
is bounded in $\wt{E_p}$.
We may therefore assume now without loss
of generality that $E$ is a Banach space,
and that $B\sub E$ is weakly bounded.
As we know from functional
analysis on normed spaces,
$E'$ becomes a Banach space when equipped
with the operator norm.
Repeating this procedure,
also the bi-dual $E''$ (the dual of $E'$)
becomes a Banach space. It is well known that
the evaluation homomorphism
\[
\eta\colon E\to E''
\]
taking $x\in E$ to
$\eta(x)\colon E'\to\K$, $\lambda\mto \lambda(x)$
is a linear isometry.\footnote{This also follows
from the Bipolar Theorem (Theorem~\ref{bipolar}),
noting that the polar $U^\circ$ of the closed
unit ball $U\sub E$ is the closed unit ball
in~$E'$, whence the polar
$U^{\circ\circ}$ of $U^\circ$ in $E''$
is the closed
unit ball in~$E''$. By the Bipolar Theorem, $\eta(U)=U^{\circ\circ}\cap\eta(E)$.}
As a consequence of Lemma~\ref{lembdsub},
$B$ will be bounded if we can show
that $\eta(B)$ is bounded in $E''$.
Given $\lambda\in E'$, let
$\ev_\lambda\colon E''\to\K$,
$\ev_\lambda(\xi)=\xi(\lambda)$
be the point evaluation.
Since $\ev_\lambda(\eta(x))=\eta(x)(\lambda)=\lambda(x)$
for each $x\in E$, we obtain
$\ev_\lambda(\eta(B))=\lambda(B)$,
which is bounded in~$\K$.
Hence, applying the Uniform Boundedness
Principle
(see \cite[Thm.\,5.8]{Ru87};
cf.\ \cite[Thm.\,2.6]{Ru91})
to the set $\eta(B)$
of continuous linear maps $E'\to\K$
on the Banach space~$E'$,
we find that $\eta(B)$
is bounded in $E''$,
as required.
\end{prf}
\begin{prop}\label{prop-normable}
A locally convex space $E$ is normable
if and only if it has a bounded $0$-neighborhood.
\end{prop}
\begin{prf}
If the topology on~$E$ comes from a norm, then the unit ball
$B_1^E(0)$ is a bounded $0$-neighborhood.
Conversely, assume that $U\sub E$ is a bounded
$0$-neighborhood.
After shrinking~$U$ if necessary, we may assume that~$U$ is absolutely convex and closed
(see Lemma~\ref{basisnice}). Let $p:=\mu_U$ be the Minkowski functional
associated with~$U$.
Then~$p$ is a continuous seminorm on~$E$.
If $V\sub E$ is any $0$-neighborhood, then $U\sub rV$
for some $r>0$ and thus $\wb{B}_{1/r}^p(0)=\frac{1}{r}\wb{B}^p_1(0)=\frac{1}{r}U\sub V$.
The vector topology~$\cT$ arising from~$p$
is therefore finer than the given topology.
Since~$\cT$ is also coarser, the two topologies
coincide. Hence $\cT$ is Hausdorff,
entailing that the seminorm~$p$ is a norm.
\end{prf}
\begin{lem}\label{extprp}
Let $E$ be a locally convex space,
$F\sub E$ be a closed vector subspace
and
$q\colon E\to E/F$ be the quotient map.
Then the following holds:
\begin{description}
\item[\rm(a)]
$E$ is metrizable if and only if both $F$ and $E/F$
are metrizable.
\item[\rm(b)]
$E$ is normable if and only if both $F$ and $E/F$ are normable.
\item[\rm(c)]
$E$ is a Fr\'{e}chet space
if and only if both $F$ and $E/F$
are Fr\'{e}chet spaces.
\item[\rm(d)]
$E$ is a Banach space
if and only if both $F$ and $E/F$
are Banach spaces.
\end{description}
\end{lem}
\begin{prf}
Set $Q:=E/F$.

(a) If $E$ is metrizable, then also
$F$. Moreover, $E$ has a countable
basis $\{U_n\colon n\in \N\}$
of $0$-neighborhoods.
Then $\{q(U_n)\colon n\in\N\}$
is countable and a basis of $0$-neighborhoods for~$Q$
(cf.\ Proposition~\ref{baconlcx}(c)),
whence $Q$ is metrizable, by Corollary~\ref{charmetriz}.
If, conversely, $F$ and~$Q$ are metrizable,
pick countable bases
$\{V_n\colon n\in \N\}$ and $\{W_n\colon n\in \N\}$
of $0$-neighborhoods for~$F$ and~$Q$,
respectively. For each $n\in \N$,
there exists a $0$-neighborhood $U_n\sub E$
such that $F\cap U_n\sub V_n$ and $U_n\sub q^{-1}(W_n)$.
With the help of Lemma~\ref{verybase} (d)~and~(c), after shrinking $U_2 , U_3, \ldots$,
we may assume that
\begin{equation}\label{easyco}
(\forall n\in \N)\quad U_{n+1}-U_{n+1}\sub U_n\,.
\end{equation}
Then
\begin{equation}\label{hol}
(\forall n\in \N)\quad
U_{n+1}\sub E \setminus (U_{n+1}+(F\setminus U_n))\,.
\end{equation}
In fact, otherwise we would have $x=y+z$ for certain
$x,y\in U_{n+1}$ and $z\in F\setminus U_n$,
entailing the contradiction
$z=x-y\in U_{n+1}-U_{n+1}\sub U_n$.
We show that $\{U_n\colon n\in \N\}$
is a basis of $0$-neighborhoods for~$E$
(whence~$E$ is metrizable by Corollary~\ref{charmetriz}).
To this end, let $R\sub E$ be any $0$-neighborhood,
and $S\sub E$ be a $0$-neighborhood
with $S+S\sub R$.
There exists $n\in\N$ such that $F\cap U_n \sub F\cap S$.
Also, there is $k>n$ such that $W_k\sub q(S\cap U_{n+1})$
and thus
\begin{equation}\label{fir}
U_k\sub (S\cap U_{n+1})+F\,.
\end{equation}
Using~(\ref{hol}), we deduce that
\[
U_k\sub U_{n+1}\sub
E\setminus (U_{n+1}+(F\setminus U_n)) \sub
E \setminus ((S\cap U_{n+1})+(F\setminus U_n)).
\]
Combining this inclusion with (\ref{fir}), we obtain
\[ U_k\sub (S\cap U_{n+1})+(F\cap U_n)
\sub S+ (F\cap S)\sub R.\]

(b) If $U\sub E$ is a bounded $0$-neighborhood,
then $U\cap F$ is a bounded $0$-neighborhood in~$F$
(Lemma~\ref{lembdsub})
and $q(U)$ a bounded $0$-neighborhood in~$Q$ (see Proposition~\ref{baconlcx}
and Lemma~\ref{lemboundd}), 
entailing that~$F$ and~$Q$ are normable (Proposition~\ref{prop-normable}).
Conversely, assume that both~$F$ and~$Q$
are normable.
Let $V\sub F$ and $W\sub Q$ be bounded $0$-neighborhoods
(see Proposition~\ref{prop-normable}).
There exists an absolutely
convex $0$-neighborhood $U\sub E$ such that $F\cap U\sub V$
and $q(U)\sub W$.
Define $U_n:=2^{-n}U$ for $n\in \N$.
Then the sets $F\cap U_n\sub 2^{-n} V$
and $q(U_n)\sub 2^{-n}W$
(for $n\in \N$)
form bases of $0$-neighborhoods
for~$F$ and~$Q$, respectively.
Moreover, (\ref{easyco}) is satisfied.
Hence $\{U_n\colon n\in \N\}$
is a basis of $0$-neighborhoods for~$E$ (by the proof of~(a)),
entailing that~$U$ is bounded and hence~$E$ normable.

(c) We let $U_1\supseteq U_2\supseteq\cdots$
be a basis of absolutely convex,
closed neighborhoods of $0$ for~$E$,
and $p_1\leq p_2\leq\cdots$ be the corresponding Minkowski functionals.

If $E$ is a Fr\'{e}chet space,
then also its closed vector subspace~$F$
is complete (see Proposition~\ref{newcompl}\,(b))
and hence Fr\'{e}chet. Also completeness of~$Q$ follows:
If $(y_n)_{n\in\N}$ is a Cauchy sequence in~$Q$,
it will converge if we can find a convergent subsequence.
After passage to a subsequence, we may hence
assume that $y_{n+1}-y_n\in 2^{-n}q(U_n)$
for all $n\in \N$, and thus $y_{n+1}-y_n=q(z_n)$
for some $z_n\in 2^{-n}U_n$.
Choose $x_1\in E$ such that $q(x_1)=y_1$.
Recursively, define $x_{n+1}:=x_n+z_n$
for $n\in \N$. Then
$q(x_n)=y_n$, by induction: In fact,
$q(x_{n+1})=q(x_n)+q(z_n)=y_n+(y_{n+1}-y_n)=y_{n+1}$.
Moreover, $x_{n+1}-x_n=z_n \in 2^{-n}U_n$.
Then $p_n(x_{n+m}-x_n)\leq 2^{-n}+\cdots+2^{-n-m+1}\leq 2^{-n+1}\leq 1$
for all $n,m\in \N$, whence $x_{n+m}-x_n\in U_n$
and thus $x_k-x_\ell\in 2 U_n$ for all $k,\ell\geq n$.
Hence $(x_n)_{n\in\N}$ is a Cauchy sequence in~$E$
and thus convergent to some $x\in E$. Then $(y_n)_{n\in\N}$ converges
to $q(x)$.

Conversely, assume that $F$ and $Q$ are Fr\'{e}chet spaces.
We have to show that every Cauchy sequence
$(x_n)_{n\in \N}$ in~$E$ converges.
Again, we may assume that $x_{n+1}-x_n\in 2^{-n}U_n$
for all $n\in\N$.
Then
\begin{equation}\label{prstp}
(\forall n,m\in\N)\quad
p_n(x_{n+m}-x_n)\leq (2^{-n}+\cdots+2^{-n-m+1})\leq 2^{-n+1}\leq 1
\end{equation}
and hence $q(x_{n+m})-q(x_n)\in q(U_n)$.
Thus
$(q(x_n))_{n\in \N}$ is a Cauchy sequence
in~$Q$ and hence convergent, to $y\in Q$ say.
Pick $a\in E$ such that $q(a)=y$.
As a consequence of (\ref{prstp}),
we have $q(x_{n+m})-q(x_n)\in 2^{-n+1}q(U_n)$
for all $n,m\in\N$ and therefore, letting $m\to\infty$,
\[
y-q(x_n)\in 2^{-n+1}\wb{q(U_n)}\,.
\]
Hence $y-q(x_n)\in 2^{-n+2}q(U_n)$,
using that
$\wb{q(U_n)}\sub q(U_n)+q(U_n)=2q(U_n)$
(see Lemma~\ref{verybase}\,(i)).
Thus, we can find $z_n\in F$
such that
\begin{equation}\label{tozer}
a+z_n-x_n\in 2^{-n+2}U_n.
\end{equation}
Then
\begin{eqnarray*}
p_n(z_{n+1}-z_n)&=&p_n((a+z_{n+1}-x_{n+1})-(a+z_n-x_n)+(x_{n+1}-x_n))\\
&\leq&
2^{-n+1}+2^{-n+2}+2^{-n}\leq 2^{-n+3}
\end{eqnarray*}
for each $n\in \N$, whence
$(z_n)_{n\in \N}$ is a Cauchy sequence in~$F$
and hence convergent to some $z\in F$.
Then $a+z_n\to a+z$ as $n\to\infty$
and hence also $x_n\to a+z$, since $(a+z_n)-x_n\to 0$
by (\ref{tozer}).

(d) Follows directly from (b) and (c).
\end{prf}
\begin{defn}
If $E$ is a locally convex space,
we write $E'_\sigma$ for $E'$, endowed
with the initial topology $\cO_\sigma$
with respect to the point evaluations
\[
\ve_x\colon E'\to\K,\quad \lambda\mto\lambda(x)
\]
for $x\in E$ (the so-called \emph{weak $*$-topology}).
Then $(E',\cO_\sigma)$ is a locally convex space,
as the $\ve_x$ are linear and separate points on~$E'$.
\end{defn}
For later use, we mention:
\begin{lem}\label{dualofweakstar}
Let $E$ be a locally convex space
and $\theta\in (E'_\sigma)'$.
Then $\theta=\ve_x$ for some $x\in E$.
\end{lem}
\begin{prf}
Since $\theta^{-1}(\bD)\sub E'_\sigma$ is a $0$-neighborhood,
we find $x_1,\ldots, x_n\in E$ and $r_1,\ldots, r_n>0$ such that
\[
\{\lambda\in E'\colon (\forall j\in\{1,\ldots, n\})\;\; |\lambda(x_j)|\leq r_j\}\sub\theta^{-1}(\bD)
\]
and hence $F:=\ker(\ve_{x_1})\cap\cdots\cap\ker(\ve_{x_n})\sub \ker\theta$.
Let $\beta\colon E\to\K^n$ be the linear map with components $\ve_{x_1},\ldots,\ve_{x_n}$.
Then $\ker\beta=F\sub\ker\theta$, whence
there is a unique linear map
\[
\wb{\theta}\colon \beta(E)\to\K
\]
with $\wb{\theta}\circ \beta=\theta$.
Extend $\wb{\theta}$ to a linear map $\Theta\colon\K^n\to \K$.
Let $e_1,\ldots, e_n$ be the standard basis vectors of~$\K^n$ and $a_j:=\Theta(e_j)\in\K$
for $j\in\{1,\ldots, n\}$. Then
\[
\Theta(t)=\sum_{j=1}^na_jt_j
\]
for all $t=(t_1,\ldots, t_n)\in\K^n$.
Hence
\[
\theta(\lambda)=\wb{\theta}(\beta(\lambda))=\Theta(\beta(\lambda))=\sum_{j=1}^n a_j\ve_{x_j}(\lambda)
=\lambda(x)=\ve_x(\lambda)
\]
holds for $x:=\sum_{j=1}^n a_jx_j$.
\end{prf}
\begin{thm}[Banach-Alaoglu]\label{alaoglu}
Let $E$ be a locally convex space and $U\sub E$ be
a $0$-neighborhood. Then the polar~$U^\circ$ is compact
in~$E'_\sigma$.
\end{thm}
\begin{prf}
If $\lambda\colon E\to\K$ is linear and $\lambda(U)\sub \bD$,
then $\lambda$ is continuous.
Hence
\[
U^\circ=\{\lambda\in\Hom_\K(E,\K)\colon (\forall x\in U)\;\lambda(x)\in\bD\},
\]
which is a closed subset of $\Hom_\K(E,\K)$ and hence also closed in $\K^E$,
endowed with the product topology (Exercise~\ref{exc-linclo}).
For each $x\in E$, there is $r_x>0$ such that $x\in r_xU$.
Then $\lambda(x)\in r_x\lambda(U)\sub r_x\bD$,
entailing that $U^\circ$ is a closed subset of
\[
\prod_{x\in E} r_x\bD\sub \K^E.
\]
As products of compact sets are compact by Tychonoff's Theorem,
we see that~$U^\circ$ is compact.
\end{prf}
\begin{defn}
Let $E$ be a locally convex space.
We write $E'_b$ for $E'$, endowed with the locally convex vector topology defined by the seminorms
\[
\|\cdot\|_B\colon E'\to[0,\infty[,\quad \|\lambda\|_B:=\sup_{x\in B}|\lambda(x)|,
\]
for $B$ ranging through the set of all
bounded subsets of~$E$.
We write $E'_c$ for $E'$, endowed with the locally convex vector topology $\cO_c$
defined by $\|\cdot\|_K$ for~$K$ ranging through the set of all compact subsets
of~$E$.
\end{defn}
\begin{rem}\label{remcvsb}
(a) Since every compact set $K\sub E$ is bounded, the map
\[
E'_b\to E'_c,\quad x\mto x
\]
is continuous.

(b) $\cO_c$ coincides with the compact-open topology induced by $C(E,\K)$
(cf.\ Lemma~\ref{sammelsu}(b)).

(c) For each $x\in E$, the point evaluation $\ve_x\colon E'_b\to\K$, $\lambda\mto\lambda(x)$
is linear and continuous (since $|\ve_x(\lambda)|=|\lambda(x)|\leq \|\lambda\|_{\{x\}}$).
Likewise for~$E'_c$.

(d) Since $r\|\cdot\|_B=\|\cdot\|_{rB}$ for each bounded subset $B\sub E$ and $r>0$,
it follows that the closed unit balls
\[
\{\lambda\in E'\colon \|\lambda\|_B\leq 1\}=B^\circ
\]
(which coincide with the polars) form a basis of $0$-neighborhoods
in $E'_b$, for~$B$ ranging through the set of bounded subsets of~$E$.

(e) If $U\sub E$ is a $0$-neighborhood, then its polar $U^\circ$ is bounded
in~$E'_b$. To see this, let $V\sub E'_b$ be a $0$-neighborhood.
Then $B^\circ\sub V$ for some bounded subset $B\sub E$.
Now $B\sub r U$ for some $r>0$,
whence $V\supseteq B^\circ\supseteq (rU)^\circ=\frac{1}{r}U^\circ$.
\end{rem}
\begin{prop}\label{opspacecompl}
Let $E$ be a locally convex space.
\begin{description}[(D)]
\item[\rm(a)]
If $E$ is a $k_\R$-space, then $E'_c$ is complete.
\item[\rm(b)]
If $E'_c$ is complete, then also $E'_b$ is complete.
\end{description}
\end{prop}
\begin{prf}
(a) Since $E'$ is closed in $C(E,\K)$ (cf.\ Exercise~\ref{exc-linclo})
and $C(E,\K)$ is complete (see Lemma~\ref{sammelsu}(d)),
we see that~$E'_c$ is complete.

(b) If $(\lambda_j)_{j\in J}$
is a Cauchy net in~$E'_b$, then also in~$E'_c$.
In~$E_c'$, it therefore converges to some $\lambda\in E'_c$.
If $B\sub E$ is a bounded set and $\ve>0$,
we find $j_0\in J$ such that
\[
\|\lambda_i-\lambda_j\|_B\leq\ve
\]
for all $i,j\in J$ such that $i\geq j_0$ and $j\geq j_0$.
Thus
\[
|\lambda_i(x)-\lambda_j(x)|\leq\ve
\]
for all $x\in B$ and $i,j$ as before. Passing to the limit in~$i$, we find that
\[
|\lambda(x)-\lambda_j(x)|\leq \ve\quad
\mbox{for all $x\in B$ and $j\geq j_0$.}
\]
Thus $\|\lambda-\lambda_j\|_B\leq \ve$
for all $j\geq j_0$ and thus $\lambda_j\to\lambda$ in~$E'_b$.
\end{prf}

\begin{rem}\label{linearlykR}
Let $E$ be a locally convex space.
The conclusion of Proposition~\ref{opspacecompl}(a)
remains valid if
linear maps $\lambda\colon E\to\K$
are continuous if and only if $\lambda|_K$ is continuous for each compact subset $K\sub E$.
This property is inherited by locally convex direct limits (cf.\ Section~\ref{sec-appDLvec});
in particular, all (LF)-spaces enjoy the property (and
hence their dual spaces are complete).
\end{rem}
More generally, for locally convex spaces $E$ and $F$,
we can consider topologies on the space
$\cL(E,F)$ of all continuous linear mappings $A\colon E\to F$.
\begin{defn}\label{bd-conv-semi}
We write $\cL(E,F)_b$ for $\cL(E,F)$, endowed with the locally convex vector topology~$\cO_b$
defined by the seminorms
\[
\|\cdot\|_{B,q}\colon \cL(E,F)\to[0,\infty[,\;\;
\|A\|_{B,q}:=\sup_{x\in B} q(Ax),
\]
for $B$ ranging through the set of bounded subsets of~$E$ and~$q$
in the set of continuous seminorms on~$F$
(the \emph{topology of bounded convergence}).
The seminorms $\|\cdot\|_{K,q}$ for $K$ in the set
of compact subsets of~$E$
define the compact-open topology (cf.\ Lemma~\ref{sammelsu}(b)); when it is used, we write~$\cL(E,F)_c$.
\end{defn}

\begin{small}
\section*{Exercises for Section~\ref{secbdd}}

\begin{exer}\label{exc-linclo}
Show that the vector subspace
$\Hom_\K(E,\K)$ of linear maps is closed in $\K^E$
(with respect to the product topology/topology of pointwise
convergence), for each vector space~$E$.
\end{exer}

\begin{exer}\label{exer-semin-norm}
Let $(E,\|\cdot\|)$ be a normed space and $q\colon E\to[0,\infty[$ be a continuous
seminorm. Show that there is $r>0$ such that $q(x)\leq r \|x\|$ for all $x\in E$.\\[1mm]
[Since $\wb{B}^{\|\cdot\|}_1(0)$ is bounded, we have $\wb{B}^{\|\cdot\|}_1(0)\sub r\wb{B}^q_1(0)=\wb{B}^q_r(0)$ for some $r>0$. For $0\not=x\in E$, this implies
$q(x)=\|x\|q(\frac{1}{\|x\|}x)\leq r\|x\|$.]
\end{exer}

\begin{exer}\label{exc-bvsopnorm}
Show: If $E$ and $F$ are normed spaces, then the topology $\cO_b$
of bounded convergence on $\cL(E,F)$
coincides with the vector topology defined by the operator norm.
\end{exer}

\begin{exer}\label{exc-on-fin}
Let $E$ and~$F$ be locally convex spaces.
\begin{description}[(D)]
\item[(a)]
Fix $\lambda\in E'$. Show that the linear mapping $F\to\cL(E,F)_b$, $y\mto \lambda\otimes y$
with $(\lambda\otimes y)(x):=\lambda(x)y$ is continuous.
\item[(b)]
If $E$ is finite-dimensional, let $v_1,\ldots, v_n$ be a basis for~$E$
and $v_1^*,\ldots, v_n^*\in E'$ be the dual basis determined by $v_i^*(v_j)=\delta_{ij}$
for $i,j\in\{1,\ldots, n\}$. Show that the map $\Phi\colon\cL(E,F)_b\to F^n$,
$\alpha\mto (\alpha(v_1),\ldots,\alpha(v_n))$
is an isomorphism of topological vector spaces.\\[1mm]
[Check that $\Phi^{-1}(y_1,\ldots,y_n)=\sum_{j=1}^n v_j^*\tensor y_j$.]
\end{description}
\end{exer}

\begin{exer}\label{exc-opscomplete}
Show that $\cL(E,F)_b$ is complete whenever $E$ is a $k_\R$-space
and $F$ is complete.
\end{exer}

\begin{exer}\label{exc-allseq}
By Example~\ref{furex}(a) and Proposition~\ref{newcompl}(a),
the space $\R^\N=\prod_{n\in\N}\R$
of all real sequences is a Fr\'{e}chet space.
Show that $\R^\N$ is not normable.\\[2.3mm]
[Hint: Every $0$-neighborhood in $\R^\N$ contains
$\{0\}^n\times \R^{\{n+1,n+2,\ldots\}}$ for some $n\in\N$,
a non-trivial vector subspace.]
\end{exer}

\begin{exer}\label{exc-no-explaw}
Show that the map $\Phi\colon C(E'_c\times E,\R)\to C(E'_c,C(E,\R))$, $f\mto f^\vee$
(as in Proposition~\ref{ctsexp})
fails to be surjective for every non-normable real locally convex space~$E$.\\[2mm]
[Hint: The inclusion map $j\colon E'_c\to C(E,\R)$, $\lambda\mto\lambda$
is an element of $C(E'_c,C(E,\R))$ but $j^\wedge=\ve$
is the evaluation map $E'_c\times E\to\R$ which fails to be continuous
by Proposition~\ref{evalctsnormable}.]
\end{exer}

\end{small}
\section{Reflexive spaces and barreled spaces}\label{secreflex}
In this section, we briefly discuss reflexivity
of locally convex spaces.
If $E$ is a locally convex space,
then the point evaluation $\ve_x\colon E'_b\to\K$
is an element of $(E'_b)'$ (see Remark~\ref{remcvsb}(c)).
We therefore obtain a linear map
\[
\eta_E\colon E\to (E'_b)'_b,\quad x\mto\ve_x
\]
which is injective as $E'$ separates points on~$E$.
\begin{defn}
If $\eta_E$ is an isomorphism of topological vector spaces,
then~$E$ is called \emph{reflexive}.
\end{defn}
\begin{rem}\label{needforsilva}
Reflexivity is a useful tool to establish completeness properties:
If $E$ is a reflexive locally convex space and $E'_b$
is a $k_\R$-space (or has the weakened property
described in Remark~\ref{linearlykR}),
then $E\cong (E'_b)'_b$ is complete,
by Proposition~\ref{opspacecompl}.
\end{rem}
\begin{lem}\label{openontoim}
For every locally convex space~$E$,
the map $\eta_E$ is open onto its image.
\end{lem}
\begin{prf}
If $U\sub E$ is an absolutely
convex, closed $0$-neighborhood, then $U^\circ$
is bounded in $E'_b$ (see Remark~\ref{remcvsb}(e)), whence $(U^\circ)^\circ$ is a $0$-neighborhood
in $(E'_b)'_b$. Since $\eta_E(E)\cap (U^\circ)^\circ=\eta_E((U^\circ)_\circ)=\eta_E(U)$
(as $U=(U^\circ)_\circ$ by the Bipolar Theorem),
we see that $\eta_E$ is open as a map to $\eta_E(E)$.
\end{prf}
\begin{defn}
A subset $W$ of a locally convex space~$E$ is called
a \emph{barrel} if $W$ is absolutely convex, closed
and absorbing.
\end{defn}
For example,
every absolutely convex, closed $0$-neighborhood
is a barrel.
\begin{defn}
A locally convex space~$E$
is called \emph{barreled} if every barrel
in~$E$ is a $0$-neighborhood.
\end{defn}
\begin{prop}\label{frechbarr}
Every Fr\'{e}chet space~$F$ is barreled.
\end{prop}
\begin{prf}
If $B\sub F$ is a barrel, then
$F=\bigcup_{n\in\N}nB$ with $nB$ closed,
whence the interior $(nB)^0=nB^0$ is non-empty
for some $n\in\N$ (by Baire's Category Theorem) and hence
$B^0\not=\emptyset$.
Let $x\in  B^0$; since $B^0=-B^0$,
we see that $0\in \frac{1}{2}x+\frac{1}{2}B^0\sub B$,
whence $B$ is a $0$-neighborhood.
\end{prf}
\begin{prop}\label{charrefl}
A locally convex space $E$ is reflexive
if and only if $E$ is barreled and every
absolutely convex, bounded, closed subset $B\sub E$ is weakly compact $($i.e.,
compact in~$E_w)$.
\end{prop}
\begin{prf}
If $E$ is reflexive and $W\sub E$ a barrel,
then $\eta_E(W)$ (like~$W$) is absorbing,
whence the polar $W^\circ$ is weakly bounded in $E'_b$
and hence bounded. Hence $(W^\circ)^\circ$ is a $0$-neighborhood
in~$(E'_b)'_b$ and hence $W=\eta_E^{-1}((W^\circ)^\circ)$
is a $0$-neighborhood in~$E$
(with equality by the Bipolar Theorem).

If~$E$ is reflexive and $B\sub E$ an absolutely convex, bounded, closed
subset, then $B=\eta_E^{-1}((B^\circ)^\circ)$
by the Bipolar Theorem. Note that
\[
B^\circ=\{\lambda\in E'\colon \|\lambda\|_B\leq 1\}
\]
is a $0$-neighborhood in~$E'_b$. Hence $(B^\circ)^\circ$ is compact
in $(E'_b)'$, endowed with the weak $*$-topology,
by the Banach--Alaoglu Theorem (Theorem~\ref{alaoglu}). Since $\eta_E\colon E_w\to (E'_b)'$ is
a homeomorphism for this topology,
we deduce that $B=\eta_E^{-1}((B^\circ)^\circ)$ is weakly compact.

If $E$ is barreled and $U\sub (E'_b)'_b$ a $0$-neighborhood, then
$U\supseteq B^\circ$ for some bounded subset
$B\sub E'_b$ (see Remark~\ref{remcvsb}(d)). Then~$B$ is also bounded in the weak $*$-topology,
whence $\ve_x(B)\sub\K$ is bounded for each $x\in E$ and hence $x\in n B_\circ$ for
some $n\in\N$. Hence $B_\circ$ is absorbing in~$E$,
whence $B_\circ$ is a barrel in~$E$ and hence a $0$-neighborhood.
Since $\eta_E^{-1}(U)\supseteq \eta_E^{-1}(B^\circ)=B_\circ$,
we deduce that $\eta_E^{-1}(U)$ is a $0$-neighborhood in~$E$
and so~$\eta_E$ is continuous. Hence~$\eta_E$ is a homeomorphism onto
its image, by Lemma~\ref{openontoim}.

Now assume that absolutely convex, bounded, closed subsets of~$E$ are weakly compact.
Suppose that $\eta_E$ is not surjective; we shall derive a contradiction.
Let $\alpha\in (E'_b)'\setminus \eta_E(E)$.
Then the polar $\{\alpha\}_\circ$ is a $0$-neighborhood in~$E'_b$.
We therefore find a bounded subset $B\sub E$ such that $B^\circ\sub \{\alpha\}_\circ$.
After replacing $B$ with its closed absolutely convex hull,
we may assume that~$B$ is, moreover,
absolutely convex and closed. Hence $B$ is weakly compact (by hypothesis)
and thus $\eta_E(B)$ is a compact, absolutely convex subset in $(E'_b)'$, endowed with the weak $*$-topology.
Since $\alpha\not\in\eta_E(B)$,
the Bipolar Theorem provides $\theta \in ((E'_b)'_\sigma)'$
such that $\theta(\eta_E(B))\sub \bD$ and $|\theta(\alpha)|>1$.
By Lemma~\ref{dualofweakstar}, we have $\theta=\ve_\lambda$ for some $\lambda\in E'$,
where $\ve_\lambda\colon (E'_b)'\to \K$ is evaluation at~$\lambda$.
Thus
\[
\lambda(B)\sub\bD
\]
(whence $\lambda\in B^\circ$) and $|\alpha(\lambda)|>1$
(whence $\lambda\not\in \{\alpha\}_\circ$). But
$B^\circ\sub \{\alpha\}_\circ$.
\end{prf}
\section{Equicontinuity, bilinear and multilinear mappings}\label{sec-equic}
In this section, we record useful basic facts concerning
multilinear mappings, and define topologies on spaces of such.
We then
recall the notion of equicontinuity and apply it to linear mappings.
Afterwards, we discuss weakened continuity properties of bilinear maps.
The latter topic can be skipped on a first reading.
\begin{lem}\label{like-op-no-mult}
Given real vector spaces $E_1,\ldots, E_n$ and $F$,
consider an $n$-linear map
$\beta\colon E_1\times\cdots\times E_n\to F$. 
Let $q_1,\ldots, q_n$ and~$p$ be seminorms on $E_1,\ldots, E_n$ and~$F$, respectively,
and $r_1,\ldots, r_n,s>0$. If
$p(\beta(x_1,\ldots,x_n))\leq s$
for all $(x_1,\ldots, x_n)\in B^{q_1}_{r_1}(0)\times\cdots\times B^{q_n}_{r_n}(0)$,
then
\[
p(\beta(x_1,\ldots, x_n))\,\leq\, \frac{s}{r_1\cdots r_n}\,q_1(x_1)\cdots q_n(x_n)
\]
for all $(x_1,\ldots,x_n) \in E_1\times\cdots\times E_n$.
\end{lem}
\begin{prf}
If $t_j>q_j(x_j)$ for $j\in\{1,\ldots, n\}$,
then $q_j(r_jx_j/t_j)<r_j$ and thus
\begin{eqnarray*}
p(\beta(x_1,\ldots, x_n))&=&(t_1/r_1)\cdots (t_n/r_n)\beta(r_1x_1/t_1,\ldots,
r_nx_n/t_n)\\
&\leq &s(t_1/r_1)\cdots(t_n/r_n).
\end{eqnarray*}
Letting $(t_1,\ldots, t_n)\to (q_1(x_1),\ldots, q_n(x_n))$ in~$\R^n$, the conclusion follows.
\end{prf} 
Taking $q_j:=q$, $r_j:=r$, $x_j:=x$ and $t_1=\cdots=t_n$
here, we also get
\begin{lem}\label{like-op-no}
$\!$Let $E$ and $F$ be real vector spaces,
$\!\beta\colon \!E^n\!\!\to \!F$ be $n$-linear and
$f\colon E\to F$, $x\mto \beta(x,\ldots, x)$
the corresponding homogeneous polynomial.
Let $q$ and $p$ be seminorms on~$E$ and~$F$, respectively,
and $r,s>0$. If
\[
p(f(x))\leq s\;\;\mbox{for all $\,x\in B^q_r(0)$,}
\]
then $p(f(x))\leq \frac{s}{r^n}(q(x))^n$ for all $x\in E$.\qed
\end{lem}
\begin{lem}\label{bf-multi-cts}
Let $E_1,\ldots, E_n$ and~$F$ be locally convex spaces
and let $\beta\colon E_1\times\cdots\times E_n\to F$ be an $n$-linear map.
The following conditions are equivalent:
\begin{description}[(D)]
\item[\rm(a)]
$\beta$ is continuous;
\item[\rm(b)]
$\beta$ is continuous at~$0$;
\item[\rm(c)]
For each continuous seminorm~$p$ on~$F$, there exist continuous seminorms $q_k$ on~$E_k$ for $k\in\{1,\ldots, n\}$ such that
\begin{equation}\label{est-as-i}
p(\beta(x_1,\ldots,x_n))\leq q_1(x_1)\cdots q_n(x_n)\;\,\mbox{for all $\, (x_1,\ldots, x_n)\in
E_1\times\cdots\times E_n$.}
\end{equation}
\end{description}
\end{lem}
\begin{prf}
Trivially, (a) implies~(b).
If~(b) holds and~$p$ is a continuous seminorm on~$F$,
then $\beta^{-1}(B^p_1(0))$ is a $0$-neighborhood in $E_1\times\cdots\times E_n$,
whence we find continuous seminorms $q_k$ on~$E_k$ for $k\in\{1,\ldots, n\}$
such that
\[
B^{q_1}_1(0)\times\cdots\times B^{q_n}_1(0)\sub \beta^{-1}(B^p_1(0)).
\]
Then (\ref{est-as-i}) holds, by Lemma~\ref{like-op-no-mult}.

If (c) holds, let us show that~$\beta$ is continuous at each
$x=(x_1,\ldots,x_n)\in E_1\times\cdots\times E_n$. Let~$p$ be a continuous seminorm on~$F$ and
$q_1,\ldots, q_n$ be as in~(b). For $y=(y_1,\ldots,y_n)\in E_1\times\cdots\times E_n$, we have
\begin{eqnarray*}
p(\beta(y)-\beta(x))&\leq&\sum_{k=1}^n p(\beta(x_1,\ldots,x_{k-1},y_k-x_k,y_{k+1},\ldots,y_n))\\
&\leq & \sum_{k=1}^n q_1(x_1)\cdots q_{k-1}(x_{k-1})q_k(y_k-x_k)q_{k+1}(y_{k+1})\cdots q_n(y_n),
\end{eqnarray*}
which tends to~$0$ as $y\to x$.
\end{prf}
\begin{defn}\label{defn-top-multi}
Given $n\in\N$, let $E_1,\ldots, E_n$ as well as~$F$ be locally convex spaces.
We write $\cL^n(E_1,\ldots,E_n;F)$
for the space of all continuous $n$-linear mappings
$\beta\colon E_1\times\cdots\times E_n\to F$.
We write $\cL^n(E_1,\ldots,E_n;F)_b$ for $\cL^n(E_1,\ldots, E_n;F)$, endowed with the locally convex vector topology defined by the seminorms
\[
\|\cdot\|_{B,p}\colon \cL^n(E_1,\ldots, E_n;F)\to[0,\infty[,\quad \|\beta\|_{B,p}
:=\sup_{x\in B}p(\beta(x)),
\]
for $B$ ranging through the set of all
bounded subsets of~$E_1\times\cdots\times E_n$
and~$p$ in the set of all continuous seminorms on~$F$.
We write $\cL^n(E_1,\ldots, E_n;F)_c$ for $\cL^n(E_1,\ldots, E_n;F)$, endowed with the locally convex vector topology $\cO_c$
defined by $\|\cdot\|_{K,p}$ for~$K$ ranging through the set of all compact subsets
of $E_1\times\cdots\times E_n$ and continuous seminorms~$p$ as before.
Thus $\cO_c$ is the topology
induced by the compact-open topology on $C(E_1\times\cdots\times E_n,F)$
(cf.\ Lemma~\ref{sammelsu}(b)).
\end{defn}
Note that $\|\beta\|_{B,p}<\infty$. To see that, let $B_k$ be the projection
of $B$ onto the $k$th component, for $k\in\{1,\ldots,n\}$. Then~$B_k$ is bounded in~$E_k$.
If $q_1,\ldots,q_n$ are as in Lemma~\ref{bf-multi-cts}(c),
then $\|\beta\|_{B,p}\leq\sup q_1(B_1)\cdots \sup q_n(B_n)<\infty$.

\begin{defn}\label{deftopalg}
(a) A \emph{$\K$-algebra} (or ``algebra,''
for short)
is a $\K$-vector space~$\cA$,
equipped with a
bilinear map
\[
\beta\colon  \cA\times \cA\to \cA\,.
\]
We call $\beta$ the \emph{algebra multiplication},
and write $xy:=\beta(x,y)$.
If the associative law $x(yz)=(xy)z$ holds in~$\cA$,
we call~$\cA$ \emph{associative}.
If, furthermore, there exists an element $\one \in \cA$ with $\1\not=0$
such that $\one x=x \one$ for each $x\in \cA$,
we call~$\cA$ a \emph{unital} associative
algebra. In this case,
we let
\[
\cA^\times\; :=\; \{x\in \cA\colon
\mbox{$(\exists x^{-1}\in \cA)\, x^{-1}x=xx^{-1}=\one $}\}
\]
be the set of all invertible
elements (``units'')
of~$\cA$. Then $\cA^\times$ is a group
under multiplication,
called the \emph{unit group}
of~$\cA$.\medskip

\noindent
(b) A $\K$-algebra $\g$ with bilinear multiplication $\g\times\g\to \g$,
$(x,y)\mto [x,y]$ is called a \emph{Lie algebra} if $[y,x]=-[x,y]$ for all $x,y\in \g$
and the Jacobi identity holds:
\[
[x,[y,z]]+[y,[z,x]]+[z,[x,y]]=0\;\,\mbox{for all $\,x,y,z\in \g$.}
\]
(c) Unless we speak about Lie algebras, we shall always assume that algebras are associative,
without mention.\medskip

\noindent
(d) A \emph{topological algebra}
is a topological $\K$-vector space~$\cA$,
equipped with a continuous bilinear map
$\beta\colon \cA\times \cA\to \cA$ making it an associative algebra.
If the topological vector space~$\cA$ is locally convex, a Fr\'{e}chet space,
or a Banach space, we say that $\cA$ is a \emph{locally convex algebra}, \emph{Fr\'{e}chet algebra},
or \emph{Banach algebra}, respectively.\medskip

\noindent
(e) A \emph{topological Lie algebra} is a topological $\K$-vector space~$\g$,
equipped with a continuous bilinear map
$\g\times \g\to \g$, $(x,y)\mto[x,y]$ making it a Lie algebra.
If the topological vector space~$\g$ is locally convex, a Fr\'{e}chet space,
or a Banach space, we say that $\g$ is a \emph{locally convex Lie algebra}, 
\emph{Fr\'{e}chet--Lie algebra},
or \emph{Banach--Lie algebra}, respectively.\medskip

\noindent
(f) As usual, a map $\alpha\colon \cA_1\to\cA_2$ between (not necessarily associative)
algebras is called an \emph{algebra homomorphism} if it is linear and $\alpha(xy)=\alpha(x)\alpha(y)$
for all $x,y\in\cA_1$. If $\cA_1$ and $\cA_2$ are unital associative algebras
and $\alpha\colon\cA_1\to\cA_2$ is an algebra homomorphism such that $\alpha(\1)=\1$,
then $\alpha$ is called a \emph{unital algebra homomorphism}.
Continuous Lie algebra homomorphisms between topological Lie algebras are also called
\emph{morphisms of topological Lie algebras}.
\end{defn}
\begin{rem}\label{whatisban}
The norm $\|\cdot\|$ defining the topology on a unital Banach algebra~$\cA$
can always be chosen such that the following holds
(see Exercise~\ref{ban-compatible}):
\begin{description}[(D)]
\item[(a)]
$\|\1\|=1$;
\item[(b)]
$\|\cdot\|$ is \emph{sub-multiplicative},
i.e., $\|xy\|\leq\|x\|\, \|y\|$ for all $x,y\in \cA$.
\end{description}
We shall always assume that the norm on a unital Banach algebra
satisfies (a) and~(b), without saying so explicitly.
\end{rem}
For example, $\cA=\cL(E)$ with the operator norm $\|\cdot\|_{\op}$
is a unital Banach algebra for each Banach space $(E,\|\cdot\|)$ with $E\not=\{0\}$.
\begin{defn}
Let $X$ be a set, $E$ be a topological vector space
and $\Gamma\sub E^X$ be a set of $E$-valued functions on~$X$.
Given $x\in X$, the set~$\Gamma$ is called \emph{equicontinuous at~$x$}
if, for each $0$-neighborhood $V\sub E$, there exists a neighborhood~$U$ of~$x$ in~$X$
such that
\[
\gamma(y)-\gamma(x)\in V\;\;\mbox{for all $\gamma\in\Gamma$ and $y\in U$.}
\]
If $\Gamma$ is equicontinuous at each $x\in X$,
then~$\Gamma$ is called \emph{equicontinuous}.
\end{defn}
If $\Gamma$ is equicontinuous, then $\Gamma\sub C(X,E)$.
The proofs of the following simple facts are left to the reader.
\begin{lem}\label{basics-equi}
Let $X$ and $Y$ be topological spaces,
$E$ and $F$ be topological vector spaces
and $\Gamma\sub E^Y$.
\begin{description}[(D)]
\item[{\rm(a)}]
If $f\colon X\to Y$ is continuous at a point $x\in X$
and $\Gamma$ is equicontinuous at~$f(x)$, then
$\{\gamma\circ f\colon \gamma\in\Gamma\}\sub E^X$
is equicontinuous at~$x$.
\item[{\rm(b)}]
If $\Gamma$ is equicontinuous at $y\in Y$
and $g\colon E\to F$ is a continuous affine-linear map,
then also $\{g\circ \gamma\colon \gamma\in\Gamma\}\sub F^Y$ is equicontinuous at~$y$.
The same conclusion holds if $g\colon E\to F$ is any uniformly continuous
map.\footnote{$g\colon E\to F$ is called \emph{uniformly continuous}
if, for each $0$-neighborhood $V\sub F$, there exists a $0$-neighborhood
$W\sub E$ such that $g(z)-g(y)\in V$ for all $y\in E$ and $z\in y+W$.}
\item[{\rm(c)}]
If $E=\prod_{j\in J}E_j$ with topological vector spaces~$E_j$
and canonical projections $\pr_j\colon E\to E_j$,
then $\Gamma$ is equicontinuous at $y\in Y$
if and only if $\{\pr_j\circ\, \gamma\colon \gamma\in\Gamma\}\sub E_j^Y$ is equicontinuous at~$y$
for each $j\in J$.
\end{description}
\end{lem}
\begin{lem}\label{lin-equi}
If $E$ and~$F$ are topological vector spaces, then a set
$\Gamma\sub F^E$ of linear mappings
is equicontinuous if and only if~$\Gamma$ is equicontinuous~at~$0$.
\end{lem}
\begin{prf}
If $V\sub F$ and $U\sub E$ are $0$-neighborhoods such that
$\gamma(U)\sub V$ for all $\gamma\in\Gamma$, then
\[
(\forall x\in E)\,(\forall \gamma\in \Gamma)\quad
\gamma(x+U)-\gamma(x)=\gamma(U)\sub V.
\]
Thus equicontinuity of $\Gamma$ at~$0$ entails its equicontinuity.
\end{prf}
We mention a version of the Banach--Steinhaus Theorem.
\begin{prop}\label{bar-BS}
Let $E$ and~$F$ be locally convex spaces and $\Gamma\sub \cL(E,F)$
be a set of continuous linear mappings which is pointwise bounded in the sense
that
$\Gamma(x):=\{\gamma(x)\colon\gamma\in\Gamma\}$
is bounded in~$F$ for each $x\in E$.
If~$E$ is barreled, then~$\Gamma$ is equicontinuous.
\end{prop}
\begin{prf}
Let $V\sub F$ be a $0$-neighborhood. After shrinking~$V$, we may assume that~$V$
is absolutely convex and closed. Then
\[
U:=\bigcap_{\gamma\in\Gamma}\gamma^{-1}(V)\vspace{-1.3mm}
\]
is an absolutely convex, closed subset of~$E$.
To see that $U$ is absorbing (and hence a barrel),
let $x\in E$. By hypothesis, $\Gamma(rx)=r\Gamma(x)\sub V$
for some $r>0$, whence $rx\in U$. Since~$E$ is barreled,
$U$ is a $0$-neighborhood. Since $\gamma(U)\sub V$ for each $\gamma\in\Gamma$
by construction, $\Gamma$ is equicontinuous at~$0$ and hence equicontinuous,
by Lemma~\ref{lin-equi}.
\end{prf}
\begin{rem}\label{equi-bd}
If $E$ and $F$ are locally convex spaces and $\Gamma\sub\cL(E,F)$
is equicontinuous, then $\Gamma$ is a bounded subset of $\cL(E,F)_b$
(and hence also of $\cL(E,F)_c$).\\[1mm]
[Given a bounded subset $B\sub E$ and a continuous seminorm~$p$ on~$F$,
there exists a $0$-neighborhood $U\sub E$ such that $\gamma(U)\sub B^p_1(0)$
for all $\gamma\in\Gamma$. Pick $r>0$ such that $B\sub rU$. Then $\gamma(B)\sub
rB^p_1(0)=B^p_r(0)$ for all $\gamma\in\Gamma$
and thus $\sup\{\|\gamma\|_{B,p}\colon\gamma\in\Gamma\}\leq r$.
Thus~$\Gamma$ is bounded, using Lemma~\ref{lembdsmn}.]
\end{rem}
Many natural examples of bilinear mappings
(like the following ones) fail to be continuous.
This motivates the study of weakened continuity
properties.
\begin{prop}\label{evalctsnormable}
Let $E$ be a locally convex space.
If $E'$ admits a vector topology $\cO$ making the evaluation map
\[
\ve\colon (E',\cO)\times E\to\K,\quad (\lambda,x)\mto\lambda(x)
\]
continuous, then $E$ is normable.
\end{prop}
\begin{prf}
There are open $0$-neighborhoods $U\sub E'$ and $V\sub E$
such that
\[
\ve(U\times V)\sub \bD.
\]
Since $U$ is absorbing, for each $\lambda\in E'$ we find $r>0$
such that $\lambda\in rU$, whence $\lambda(V)\sub r \bD$.
Hence $V$ is weakly bounded and thus
bounded, by Theorem~\ref{Mackey}.
Hence~$E$ is normable, by Proposition~\ref{prop-normable}. 
\end{prf}
Let $E_1$, $E_2$, and $F$ be topological vector spaces.
A bilinear mapping\linebreak
$\beta\colon E_1\times E_2\to F$
is called \emph{separately continuous} if $\beta(x,\cdot)\colon E_2\to F$
is continuous for each $x\in E_1$ and $\beta(\cdot,y)\colon E_1\to F$
is continuous for each $y\in E_2$.
A separately continuous bilinear map between locally convex spaces is called
\emph{hypocontinuous} (with respect to the second argument)
if it satisfies one (and hence all) of the following
equivalent conditions.
\begin{prop}\label{three-hypo}
Let $E_1$, $E_2$, and $F$ be locally convex spaces.
The following three conditions are equivalent for a separately continuous bilinear mapping
$\beta\colon E_1\times E_2\to F$:
\begin{description}[(D)]
\item[{\rm(a)}]
For each bounded subset $B\sub E_2$,
the restriction $\beta|_{E_1\times B}\colon E_1\times B\to F$ is continuous.
\item[{\rm(b)}]
For each bounded subset $B\sub E_2$ and each $0$-neighborhood $V\sub F$,
there exists a $0$-neighborhood $U\sub E_1$ such that $\beta(U\times B)\sub V$.
\item[{\rm(c)}]
For each bounded subset $B\sub E_2$,
the set $\{\beta(\cdot,y)\colon y\in B\}\sub \cL(E_1,F)$
is equicontinuous.
\end{description}
\end{prop} 
\begin{prf}
(c) the equicontinuity of the sets in (c) is equivalent to their equicontinuity
at~$0$ (see Lemma~\ref{lin-equi}); the latter is condition~(b).

(a)$\impl$(b): After replacing $B$ with its absolutely convex hull,
we may assume that $B$ is absolutely convex.
Let $V\sub F$ be a $0$-neighborhood.
Since $\beta|_{E\times B}$ is continuous at~$(0,0)$,
there exist $0$-neighborhoods $U\sub E_1$ and $W\sub E_2$ with
\[
\beta(U\times (B\cap W))\sub V.
\]
Since~$B$ is bounded, we find
$r\in\;]0,1]$ such that $rB\sub W$.
Then $rU$ is a $0$-neighborhood in~$E_1$ and
we deduce from $rB\sub B\cap W$ that
\[ \beta((rU)\times B)=
  \beta(U\times (rB))\sub\beta(U\times (B\cap W))\sub V.\]

(b)$\impl$(a): Let $B\sub E_2$ be a bounded set and $(x,y)\in E\times B$.
Let $W\sub F$ be a $0$-neighborhood; after shrinking~$W$,
we may assume that $W$ is absolutely convex.
Since $\beta(x,\cdot)$ and $\beta(\cdot,y)$ are continuous,
we find $0$-neighborhoods $U\sub E_1$ and $V\sub E_2$
such that
\[
\beta(U\times\{y\})\sub \frac{1}{3}W\quad\mbox{and}\quad
\beta(\{x\}\times V)\sub \frac{1}{3}W.
\]
As $B-B$ is bounded, after shrinking~$U$ we may assume that, moreover,
\[
\beta(U\times (B-B))\sub\frac{1}{3}W,
\]
by~(b).
For all $(x_1,y_2)\in (x+U)\times (B\cap (y+V))$,
\begin{eqnarray*}
\beta(x_1,y_1)-\beta(x,y)&=&
\beta(x_1-x,y_1-y)+\beta(x_1-x,y)+\beta(x,y_1-y)\\
&\in&
\frac{1}{3}W+\frac{1}{3}W+\frac{1}{3}W\sub W
\end{eqnarray*}
follows. Hence $\beta|_{E\times B}$ is continuous and thus~(b) entails~(a).
\end{prf}
\begin{ex}\label{exa-eval}
If $E$ and $F$ are locally convex spaces,
endow $\cL(E,F)$ with the topology of bounded
convergence.
Then the evaluation map
\[
\ve\colon \cL(E,F)_b\times E\to F,\quad (A,x)\mto A(x)
\]
is a separately continuous bilinear map. If $B\sub E$ is a bounded subset
and $V\sub F$ a $0$-neighborhood, there exists a continuous seminorm
$q\colon F\to [0,\infty[$ on~$F$ such that $B^q_1(0)\sub V$.
Then
\[
U:=\{ A \in\cL(E,F)\colon \|A\|_{B,q}<1\}
\]
(with notation as in Definition~\ref{bd-conv-semi}) is a $0$-neighborhood in $\cL(E,F)_b$
and $\ve(U\times B)\sub B^q_1(0)\sub V$.
Hence $\ve$ is hypocontinuous, by Proposition~\ref{three-hypo}(b).
Yet, if $F=\K$ and $E$ is not normable, then $\ve$ is not continuous
(see Proposition~\ref{evalctsnormable}).
\end{ex}
In the following two corollaries, $E_1$, $E_2$, and $F$
are locally convex spaces.
\begin{cor}\label{hypo-seq}
Every hypocontinuous bilinear map $\beta\colon E_1\times E_2\to F$
is sequentially continuous.
\end{cor}
\begin{prf}
Let $(x_n,y_n)_{n\in\N}$ be a convergent sequence in $E_1\times E_2$,
with limit $(x,y)$. Then $B:=\{y_n\colon n\in\N\}\cup\{y\}$
is a compact subset of~$E_2$ and hence bounded.
As $\beta$ is hypocontinuous,
$f:=\beta|_{E_1\times B}$
is continuous by Proposition~\ref{three-hypo}(a).
Thus $\beta(x_n,y_n)=f(x_n,y_n)\to f(x,y)=\beta(x,y)$
as $n\to\infty$.
\end{prf}
\begin{cor}\label{bar-hypo}
If $E_1$ is barreled, then each separately continuous bi\-linear map
$\beta\colon E_1\times E_2\to F$ is hypocontinuous.
\end{cor}
\begin{prf}
If $B\sub E_2$ is a bounded subset, then $\{\beta(\cdot,y)\colon y\in B\}\sub\cL(E_1,F)$
is pointwise bounded as $\{\beta(x,y)\colon y\in B\}\sub F$ is bounded for each $x\in E_1$,
using that $\beta(x,\cdot)\colon E_2\to F$ is continuous linear.
Thus $\{\beta(\cdot,y)\colon y\in B\}$ is equicontinuous, by Proposition~\ref{bar-BS}.
Hence $\beta$ is hypocontinuous, by Proposition~\ref{three-hypo}. 
\end{prf}

\begin{small}
\section*{Exercises for Section~\ref{sec-equic}}

\begin{exer}\label{ban-compatible}
Let $\cA$ be a unital topological algebra
whose underlying topological vector space is normable.
Let $\|\cdot\|$ be a norm defining the topology of~$\cA$.
By Lemma~\ref{bf-multi-cts} and Exercise~\ref{exer-semin-norm},
there exists $C\geq 0$ such that $\|ab\|\leq C\|a\|\,\|b\|$ for all $a,b\in \cA$.
We define $L_a\in\cL(\cA)$ via $L_a(b):=ab$ for $a,b\in\cA$.
Then $\|L_a\|_{\op}\leq C\|a\|$ for all $a\in\cA$.
\begin{description}[(D)]
\item[(a)]
Show that the map $L\colon \cA\to\cL(\cA)$ is a homomorphism of unital algebras
and continuous (as $\|L\|_{\op}\leq C$).
\item[(b)]
Check that $\ve_{\1}\colon\cL(\cA)\to\cA$, $S\mto S(\1)$ is a continuous linear map
and deduce from $\ve_{\1}\circ L=\id_{\cA}$ that~$L$ is a topological embedding.
\item[(c)]
Show that $\|a\|':=\|L_a\|_{\op}$ is a norm on~$\cA$ which is equivalent to~$\|\cdot\|$.
Moreover, $\|\1\|'=1$ and~$\|\cdot\|'$ is submultiplicative.
\end{description}
\end{exer}

\begin{exer}\label{bil-bd-zero}
Let $E_1$, $E_2$, and $F$ be locally convex spaces,
$\beta\colon E_1\times E_2\to F$ be a continuous
bilinear map,
$(x_n)_{n\in\N}$ be a sequence in~$E_1$
such that $x_n\to 0$ as $n\to\infty$,
and $(y_n)_{n\in\N}$ be a sequence in~$E_2$ which is bounded
(in the sense that $\{y_n\colon n\in\N\}$ is bounded in~$E_2$).
Deduce from Lemma~\ref{bf-multi-cts}(c)
that $\beta(x_n,y_n)\to 0$ as $n\to\infty$.
Show that the conclusion remains valid if~$\beta$
is only assumed hypocontinuous in the second argument.
\end{exer}

\begin{exer}\label{exc-equimult}
Let $F$ be a locally convex space,
$E:=E_1\times\cdots\times E_n$
be a product of locally convex spaces with $n\in\N$ and $\Gamma\sub F^E$
be a set of $n$-linear mappings. Show that the following
conditions are equivalent:
\begin{description}[(D)]
\item[(a)]
$\Gamma$ is equicontinuous;
\item[(b)]
$\Gamma$ is equicontinuous at $0$;
\item[(c)]
For each continuous seminorm $p$ on~$F$,
there exist continuous seminorms $q_i$ on~$E_i$ for $i\in\{1,\ldots, n\}$
such that
\[
(\forall \beta\in\Gamma)\,(\forall (x_1,\ldots,x_n)\in E)\;
p(\beta(x_1,\ldots,x_n))\leq q_1(x_1)\cdots q_n(x_n).
\]
\end{description}
\end{exer}

\begin{exer}\label{exer-wopbd}
Let $(E,\|\cdot\|_E)$ and $(F,\|\cdot\|_F)$ be normed
spaces and $B$ be a bounded subset in $(\cL(E,F),\|\cdot\|_{\op})$;
thus $\sup\{\|A\|_{\op} \colon A\in B\}<\infty$.
For $y\in E$, let
\[
\ve_y\colon B \to F,\quad A\mto A(y)
\]
be evaluation at~$y$.
Let $\cO$ be the initial topology on~$B$ with respect to $(\ve_y)_{y\in E}$
and $\cO_D$ be the initial topology on~$B$
with respect to $(\ve_y)_{y\in D}$,
where $D\sub E$ is a dense vector subspace.
\begin{description}[(D)]
\item[(a)]
Show that a net $(A_j)_{j\in J}$ in~$B$ converges to
$A\in B$ with respect to~$\cO$ if and only if it does so in $(B,\cO_D)$.
Deduce that $\cO=\cO_D$.
\item[(b)]
Show that the evaluation map
\[
\ve\colon (B,\cO)\times E\to F,\quad (A,y)\mto A(y)
\]
is continuous.
\item[(c)]
Let $X$ be a topological space and $f\colon X\times E\to F$
be a map such that $f^\vee(x):=f(x,\cdot)\in B$ for each $x\in X$.
Show that $f$ is continuous if and only if
\[ f(\cdot,y)\colon X\to F, \quad x\mto f(x,y) \] 
is continuous for each $y\in D$.\\[2mm]
[Hint: We have $f(x,y)=\ve(f^\vee(x),y)$ for $x\in X$ and $y\in E$.]
\end{description}
\end{exer}

\begin{exer}\label{exc-polar-co}
Let $E$ be a locally convex space.
\begin{description}[(D)]
\item[(a)]
Show that a subset $\Gamma\sub E'$
is equicontinuous if and only if $\Gamma\sub U^\circ$
for some $0$-neighborhood $U\sub E$.
Notably, $U^\circ$ is equicontinuous
for each $0$-neighborhood.
\item[(b)]
Show that on an equicontinuous subset
$\Gamma\sub E'$, the topologies induced by $E'_\sigma$
and $E'_c$ concide.
\end{description}
Combining (a) and (b) with the Theorem of Banach-Alaoglu,
we see that the polar $U^\circ$ is a compact subset of $E'_c$
for each open $0$-neighborhood $U\sub E$.
\end{exer}

\begin{exer}\label{exc-eval-co}
Let $E$ be a real locally convex space. Show that the evaluation map
$\ve\colon E'_c\times E\to\R$, $(\lambda,x)\mto\lambda(x)$
is continuous if and only if~$E$ has finite dimension.\\[1mm]
[If $\ve$ is continuous, there exist $0$-neighborhoods $U\sub E$
and $V\sub E'_c$ such that $\ve(V\times U)\sub[{-1},1]$
and thus $V\sub U^\circ$, whence $U^\circ$ is a compact
$0$-neighborhood in $E'_c$ (see Exercise~\ref{exc-polar-co})
and thus $\dim(E')<\infty$, whence $\dim(E)<\infty$.\,]
\end{exer} 
\end{small}
\section{Projective limits}\label{secprolim}
We now define projective limits
of topological vector spaces,
which are useful tools in various chapters of the
book.
\begin{defn}\label{defprosys}
Let $(I,\leq)$ be a directed set
(see Definition~\ref{defdirset}).
A~{\em projective system
of topological vector spaces\/} (over~$I$)
is a family $(E_i)_{i\in I}$
of topological vector spaces
indexed by~$I$,
together with a family $(q_{ij})_{i\leq j}$
of continuous linear maps $q_{ij}\!: E_j\to E_i$,
indexed by pairs $(i,j)\in I\times I$ such that
$i\leq j$, with the property that
\begin{description}
\item[(a)]
$q_{ii}=\id_{E_i}$ for each $i\in I$ and
\item[(b)]
$q_{ij}\circ q_{jk}=q_{ik}$
for all $i,j,k\in I$ such that $i\leq j\leq k$.
\end{description}
Given a projective system $\cS=((E_i)_{i\in I},
(q_{ij})_{i\leq j})$ of topological vector
spaces,
let $P:=\prod_{i\in I} E_i$ be their cartesian
product (equipped with the product topology)
and define
\begin{equation}\label{stan-PL}
E\;:=\; \Big\{(x_i)_{i\in I}\in P \colon
\mbox{$q_{ij}(x_j)=x_i$ for all $i,j\in I$ such that $i\leq j$}\Big\}\,.
\end{equation}
Then $E$ is
a closed vector subspace of~$P$
(exercise); equipped with the induced topology,
it becomes a topological vector
space called the {\em projective limit\/}
of the projective system~$\cS$.
The mappings $\pi_i:=\pr_i|_E\!: E\to E_i$
(where $\pr_i\!: P\to E_i$ is the
respective coordinate projection)
are called the {\em limit maps.} 
More generally,
if $F$ is a topological vector space
isomorphic to~$E$
and $\phi\!: F\to E$
an isomorphism of topological vector
spaces,
we call~$F$, together with the maps
$q_i:=\pi_i\circ \phi$,
a projective limit of~$\cS$ and also write
$\pl\, E_i:=F$.
It will always be clear from the context whether we are working with the standard projective limit
$(E,(\pi_i)_{i\in I})$
or a projective limit $(F,(q_i)_{i\in I})$.
\end{defn}
\begin{rem}\label{first-PL}
(a) Clearly, projective limits
of locally convex spaces are locally convex. If each $E_i$ is complete,
then $\pl E_i$\vspace{-.8mm} is complete,
being isomorphic to a closed vector subspace
of the direct product $\prod_{i\in I}E_i$.\smallskip

(b) If each $\cA_i :=E_i$ is a locally convex topological algebra in
the situation of Definition~\ref{defprosys} and each $q_{ij}$
a continuous algebra homomorphism, then $((\cA_i)_{i\in I}, (q_{ij})_{i\leq j})$
is called a \emph{projective system of locally convex topological algebras}
and $\cA:=E$ from~(\ref{stan-PL}) is a subalgebra of~$P$
and hence a topological algebra. If each $\cA_i$ is unital and $q_{ij}(\one)=\one$
for all $i\leq j$, then~$\cA$~is unital.\smallskip

(c) An analogous definition applies to topological spaces:
Let $(I,\leq)$ be a directed set, $(X_i)_{i\in I}$
be a family of topological spaces and $(q_{ij})_{i\leq j}$
a family of continuous mappings $q_{ij}\colon X_j\to X_i$ for $i\leq j$ in~$I$
such that $q_{ii}=\id_{X_i}$ for all $i\in I$ and condition~(b)
of Definition~\ref{defprosys} is satisfied.
Then
\[
\pl X_i\,:=\, \Big\{(x_i)_{i\in I}\in \prod_{i\in I}X_i \colon (\forall i\leq j)\;q_{ij}(x_j)
=x_i\Big\},\vspace{-.8mm}
\]
endowed with the topology induced by the direct product,
is called the projective limit of the given projective system
of topological spaces.
We write $\pi_i$ (or $q_i$) for the projection from the projective limit
onto the $i$th component,~$X_i$.
\end{rem}
The next lemma includes projective limits
of locally convex spaces.
\begin{lem}\label{basisPL}
Consider the projective limit~$X$
of a projective system $((X_i)_{i\in I}, (q_{ij})_{i\leq j})$
of topological spaces, with limit maps $q_i\colon X\to X_i$.
Let $x\in X$. Then the sets $q_i^{-1}(U)$ form a basis of $x$-neighborhoods in~$X$,
for $i$ ranging through $I$ and $U$ ranging through the neighborhoods of~$q_i(x)$ in~$X_i$.
\end{lem}
\begin{prf}
Every neighborhood~$V$ of~$x$ in $X\sub\prod_{i \in I}X_i$
contains a neighborhood of the form $X\cap\prod_{i\in I}V_i$,
where $V_i$ is a neighborhood of $x_i:=q_i(x)$ in~$X_i$
and $F:=\{i\in I\colon V_i\not=X_i\}$
is finite. Let $i_0\in I$ such that $i_0\geq i$ for all $i\in F$.
Then
\[
U:=\{y\in X_{i_0}\colon (\forall i\in F)\;q_{ii_0}(y)\in V_i\}
\]
is a neighborhood of~$x_{i_0}$ in~$X_{i_0}$
and $q_{i_0}^{-1}(U)\sub V$ since $z\in q_{i_0}^{-1}(U)$
entails that
$q_i(z)=q_{ii_0}(q_{i_0}(z))\in q_{ii_0}(U)\sub V_i$ for all $i\in F$.
\end{prf}
The following fact frequently allows a given locally convex space
to be identified as a projective limit.
\begin{prop}\label{hand-on-PL}
Let $((E_i)_{i\in I},(q_{ij})_{i\leq j})$ be a projective system
of topological vector spaces and $E\sub\prod_{i\in I}E_i$
be its projective limit with limit maps $q_i\colon E\to E_i$.
Let $F$ be a topological vector space and $\phi_i\colon F\to E_i$
be continuous linear mappings for $i\in I$
such that
\begin{equation}\label{plcone}
q_{ij}\circ \phi_j=\phi_i\quad\mbox{for all $i,j\in I$ such that $i\leq j$.}
\end{equation}
Then there is a unique map $\phi\colon F\to E$ such that
\begin{equation}\label{induPL}
q_i\circ\phi=\phi_i\quad\mbox{for all $\in I$,}
\end{equation}
and $\phi$ is continuous and linear.
The following holds:
\begin{description}[(D)]
\item[{\rm(a)}]
If $\phi_i$ has dense image in~$E_i$ for all $i\in I$,
then $\phi$ has dense image.
\item[{\rm(b)}]
If the sets $\phi_i^{-1}(U)$
form a basis of $0$-neighborhoods in~$F$
for $i$ ranging through~$I$ and $U$ ranging through the set of
$0$-neighborhoods in~$E_i$,
then $\phi$ is a topological embedding.
\item[{\rm(c)}]
If $F$ is complete and the hypotheses of \emph{(a)}
and \emph{(b)}
are satisfied, then $\phi\colon F\to E$ is an isomorphism
of topological vector spaces.
\end{description}
\end{prop}
\begin{prf}
The map $(\phi_i)_{i\in I}\colon F\to \prod_{i\in I}E_i$
is continuous and linear. Its image is contained in~$E$,
by~(\ref{plcone}). Now
$\phi:=(\phi_i)_{i\in I}|^E$ has the desired properties.

(a) Let $x\in E$ and $V$ be a neighborhood of~$x$ in~$E$.
By Lemma~\ref{basisPL}, there exists $i\in I$ and a neighborhood~$U$ of
$q_i(x)$ in~$E_i$ such that $q_i^{-1}(U)\sub V$.
By hypothesis, there exists $y\in F$ such that $\phi_i(y)\in U$.
Then $q_i(\phi(y))=\phi_i(y)\in U$ shows that
$\phi(y)\in q_i^{-1}(U)\sub V$. Thus $\phi(F)$ is dense in~$E$.

(b) As~$F$ is Hausdorff, the hypothesis of~(b) entails that the maps
$\phi_i$ separate points on~$F$ for $i\in I$. As a consequence, $\phi$
is injective. We already know that~$\phi$ is continuous and linear.
To see that~$\phi$ is open onto its image,
let $V\sub F$ be a $0$-neighborhood. By hypothesis,
there exists $i\in I$ and a $0$-neighborhood $U\sub E_i$
such that $\phi_i^{-1}(U)\sub V$.
Then $q_i^{-1}(U)\cap \phi(F)$ is a $0$-neighborhood in $\phi(F)$
in the topology induced by~$E$. Since $q_i(\phi(x))=\phi_i(x)$, we have
\[
q_i^{-1}(U)\cap\phi(F)=\{\phi(x)\colon x\in F\;\mbox{with $q_i(\phi(x))\in U$}\}
=\phi(\phi_i^{-1}(U))\sub\phi(V),
\]
showing that also $\phi(V)$ is a $0$-neighborhood in~$\phi(E)$.
Hence $\phi$ is an open map (and therefore a topological embedding),
by Proposition~\ref{ctsopen}(b).

(c) By (a) and~(b), $\phi$ is a topological embedding with dense
image. Since $\phi(F)$ (like $F$) is complete,
it is closed in~$E$, by
Lemma~\ref{compl-clo}.
Hence $E=\wb{\phi(F)}=\phi(F)$
and ~$\phi$ is an isomorphism of topological vector spaces.
\end{prf}
\begin{ex}
$\N_0$, equipped with its usual order,
is a directed set.
The sequence $(C^i([0,1],\R))_{i\in \N_0}$,
together with the inclusion maps
\begin{equation}\label{defincl}
q_{ij}\colon C^j([0,1],\R)\to
C^i([0,1],\R)\,,\quad q_{ij}(\gamma)\, :=\, \gamma
\end{equation}
for $i\leq j$,
apparently is a projective system.
We claim that
\begin{equation}\label{expl}
C^\infty([0,1],\R)\; \cong \; \pl C^k([0,1],\R)\,.\vspace{-.8mm}
\end{equation}
To see this,
let $P:=\prod_{i\in \N_0} C^i([0,1 ],\R)$.
As $q_{ij}(\gamma_j)=\gamma_i$
holds for $\gamma_j\in C^j([0,1],\R)$
and $\gamma_i\in C^i([0,1],\R)$
if and only if $\gamma_i=\gamma_j$
(see (\ref{defincl})),
we get
\[
E\;:=\; \big\{(\gamma_i)_{i\in \N_0}\in P\!:
\mbox{$q_{ij}(\gamma_j)
=\gamma_i$ if $i\leq j$} \big\}
\;=\; \{(\gamma)_{i\in \N_0}\!: \gamma\in C^\infty([0,1],\R)\}\,.
\]
Thus $\phi\!: C^\infty([0,1],\R)\to E$,
$\phi(\gamma):=(\gamma)_{i\in \N_0}=(\gamma,\gamma,\ldots)$
is an isomorphism of vector spaces.
Then $q_i:=\pr_i\circ \,\phi\!: C^\infty([0,1 ],\R)\to
C^i([0,1],\R)$ is the inclusion map
for each $i\in \N_0$,
which is continuous by definition of the topologies.
Each component $q_i$ of~$\phi$ being continuous,
$\phi$ is continuous.
Proposition~\ref{hand-on-PL}(b)
shows that $\phi$ is an open map.
Hence
$\phi$ is an isomorphism of
topological vector spaces,
and thus (\ref{expl}) holds.
\end{ex}
\begin{ex}
Every Fr\'{e}chet space~$E$
is a projective limit of normed
spaces.
Indeed, by Proposition~\ref{charmetriz},
there exists a sequence $(p_n)_{n\in \N}$
of continuous seminorms on~$E$ defining its topology.
After replacing each $p_n$ by the pointwise maximum $\max\{p_1,\ldots, p_n\}$,
we may assume that $p_1\leq p_2\leq\cdots$.
For $n,m\in \N$ such that $n\leq m$,
Remark~\ref{littlerem} provides a unique
continuous linear map $\alpha_{n,m}\colon E_{p_m}\to E_{p_n}$
such that $\alpha_{n,m}\circ \alpha_{p_m}=\alpha_{p_n}$.
For each $n\in \N$, we complete $E_{p_n}$
to a Banach space~$F_n$, such that $E_{p_n}\sub F_n$.
We write $q_n$ for $\alpha_{p_n}$, considered
as a map to~$F_n$.
For $n\leq m$, we let $q_{n,m}\colon F_m\to F_n$
be the unique continuous extension of $\alpha_{n,m}$
provided by Proposition~\ref{extendlin}.
Then $\cS:=((F_n)_{n\in \N},(q_{n,m})_{n\leq m})$
is a projective system
and $E$, together with the maps $q_n$,
is a projective limit for~$\cS$
(as a special case of Proposition~\ref{hand-on-PL}(c)).
A similar argument shows that every
complete locally convex space~$E$ is
a projective limit of Banach spaces;
the projective system is $((\wt{E}_p)_{p\in \cP},
(\wt{\alpha}_{p,q})_{p\leq q})$,
where $\cP$ is the set of all
continuous seminorms on~$E$ and
$\wt{\alpha}_{p,q}\colon \wt{E}_q\to\wt{E}_p$
the unique continuous
extension of $\alpha_{p,q}\colon E_q\to E_p$
to a map between the completions.
\end{ex}
\section{Direct sums and locally convex direct limits}\label{sec-appDLvec}
We now consider locally convex direct sums, as these are an important
tool for the study of spaces of compactly supported smooth functions (and spaces of compactly supported smooth sections in vector bundles).
Afterwards, we consider locally convex direct limits $E=\bigcup_{n\in\N}E_n$
for ascending sequences $E_1\sub E_2\sub\cdots$
of locally convex spaces.
We compile facts concerning two important cases
(which are only needed for specific examples in the book):
\begin{description}[(D)]
\item[(i)]
Strict direct limits (when $E_{n+1}$ induces the given topology on~$E_n$);\footnote{Although spaces
of compactly supported functions and sections are examples of strict
direct limits, their essential properties are more easily established via
embeddings into locally convex direct sums.}
\item[(ii)]
Silva spaces~$E$ (when each $E_n$ is a Banach space and each inclusion map $E_n\to E_{n+1}$
a compact operator).
\end{description}
\subsection*{Locally convex direct sums}
For a family $(E_j)_{j\in J}$ of locally convex spaces
with countable index set~$J$, we consider the vector subspace
\[
\bigoplus_{j\in J}E_j:=\Big\{(x_j)_{j\in J}\in\prod_{j\in J}E_j\colon \mbox{$x_j=0$ for all but finitely many $j\,$}\Big\}\sub \prod_{j\in J}E_j.
\]
For $j\in J$, let $\cP_j$ be the set of all continuous
seminorms on~$E_j$. Then
\begin{equation}\label{refsupsemi}
p \colon
\bigoplus_{j\in J} E_j\to[0,\infty[,\quad (x_j)_{j\in J}\mto \sup\{p_j(x_j)\colon j\in J\}
\end{equation}
is a seminorm, for all $(p_j)_{j\in J}\in \prod_{j\in J} \cP_j$.
Let $\cO_{\lcx}$ be the locally convex vector topology on $\bigoplus_{j\in J} E_j$
defined by all seminorms~$p$ as in~(\ref{refsupsemi}).
When endowed with this topology, $\bigoplus_{j\in J}E_j$
is called the \emph{locally convex direct sum} of the spaces $E_j$,
for $j\in J$.
If all spaces $E_j$ are equal, say $E_j=E$, we abbreviate
\[
E^{(J)}:=\bigoplus_{j\in J}E.
\]
\begin{rem}\label{firstremsums}
(a) Note that
\[
B^p_r(x)=\bigoplus_{j\in J} B^{p_j}_r(x_j)
\]
for all $p$ as in (\ref{refsupsemi}), $x=(x_j)_{j\in J}\in \bigoplus_{j\in J} E_j$ and $r>0$, writing
\[
\bigoplus_{j\in J} V_j:=\bigoplus_{j\in J}E_j\cap\prod_{j\in J} V_j
\]
for all subsets $V_j\sub E_j$ such that $0\in V_j$ for all but finitely many~$j$.
Since seminorms can be replaced with positive multiples, this implies that the sets of the form
\[
\bigoplus_{j\in J}B^{q_j}_{r_j}(x_j)
\]
with $p_j\in \cP_j$ and $r_j>0$ form a basis of open neighborhoods of~$x$
and hence also the sets $\bigoplus_{j\in J} V_j$, if $V_j$ ranges through
a basis of open neighborhoods of $x_j$ in~$E_j$ for each~$j$.

(b) It is clear from (a) that the inclusion map $\bigoplus_{j\in J}E_j\to\prod_{j\in J}E_j$
is continuous.

(c) For a finite subset $M\sub J$, consider the linear map
\[
\lambda_M\colon \prod_{j\in M}E_j\to \bigoplus_{j\in J}E_j
\]
taking $(x_j)_{j\in M}$ to $(y_j)_{j\in J}$ with $y_j:=x_j$ if $j\in M$, $y_j:=0$ otherwise.
It is clear from (a) that $\lambda_M$ is a topological embedding.
Using $\lambda_M$, we identify $\prod_{j\in M} E_j$ with $\im(\lambda_M)$.
In particular, we can identify $E_j$ with $\lambda_{\{j\}}(E_j)$.
If~$J$ is finite, then $\bigoplus_{j\in J}E_j=\prod_{j\in J}E_j$
with the product topology.
\end{rem}
\begin{lem}\label{firstlasum}
Let $(E_j)_{j\in J}$
be a countable family of locally convex spaces.
\begin{description}
\item[\rm(a)]
If $\alpha\colon \bigoplus_{j\in J}E_j\to F$ is a linear map to a topological vector space~$F$,
then $\alpha$ is continuous if and only if $\alpha|_{E_j}$ is continuous
for each $j\in J$.
\item[\rm(b)]
If $(J_i)_{i\in I}$ is a partition of~$J$ into non-empty disjoint subsets $J_i$,
then the vector space isomorphism
\[
\bigoplus_{j\in J} E_j\to \bigoplus_{i\in I}\Big(\bigoplus_{j\in J_i}E_j\Big),\quad
(x_j)_{j\in J}\mto ((x_j)_{j\in J_i})_{i\in I}
\]
is an isomorphism of topological vector spaces.
\item[\rm(c)]
If $(E_j)_{j\in J}$ and $(F_j)_{j\in J}$ are countable families of locally convex spaces, then
the following map is an isomorphism of topological vector spaces:
\[
\bigoplus_{j\in J}(E_j\times F_j)\to \Big(\bigoplus_{j\in J} E_j\Big)\times\Big(\bigoplus_{j\in J}F_j\Big),\;\;
(x_j,y_j)_{j\in J}\mto ((x_j)_{j\in J}, (y_j)_{j\in J}).
\]
\item[\rm(d)]
If $F_j\sub E_j$ is a vector subspace for each $j\in J$, then the topology induced by $\bigoplus_{j\in J}E_j$
makes $\bigoplus_{j\in J} F_j$ the locally convex direct sum of the spaces~$F_j$, for $j\in J$.
\item[\rm(e)]
Every bounded subset of $E:=\bigoplus_{j\in J}E_j$ is contained in $\prod_{j\in M} E_j$
for some finite subset $M\sub J$.
\item[\rm(f)]
If each of the spaces $E_j$ is complete $($resp., quasi-complete, resp., sequentially complete, resp.,
Mackey complete$)$, then also $\bigoplus_{j\in J}E_j$ is complete $($resp., quasi-complete, resp., sequentially complete, resp.,
Mackey complete$)$.
\end{description}
\end{lem}
\begin{prf}
(a) For the non-trivial implication, let $Q_0\sub F$ be a $0$-neighborhood;
for $k\in\N$, find a $0$-neighborhood $Q_k\sub F$ such that $Q_k+Q_k\sub Q_{k-1}$.
Then $Q_1+ \cdots+Q_k\sub Q_0$ for all $k$. Let $J\to\N$, $j\mto n_j$ be an injective map.
Now
\[
W:=\bigoplus_{j\in J} (\alpha|_{E_j})^{-1}(Q_{n_j})
\]
is a $0$-neighborhood in $\bigoplus_{j\in J}E_j$ such that $\alpha(W)\sub\bigcup_{k\in \N}(Q_1+\cdots+Q_k)
\sub Q_0$.

(b) This follows from the fact that sets of the form
\[
\bigoplus_{i\in I}\Big(\bigoplus_{j\in J_i}B^{p_j}_r(0)\Big)
\]
form a basis of $0$-neighborhoods in $\bigoplus_{i\in I}\left(\bigoplus_{j\in J_i}E_j\right)$
(see Remark~\ref{firstremsums}(a)). 

(c) follows from (b), applied to two partitions of the index set $J\times \{1,2\}$.

(d) The sets $\bigoplus_{j\in J} (F_j\cap V_j)$ form a basis of $0$-neighborhoods
in $\bigoplus_{j\in J}F_j$ for both topologies,
for $V_j$ in the set of all open $0$-neighborhoods in~$E_j$.

(e) Suppose $B\sub E$ was bounded but not contained in a finite partial product;
we derive a contradiction.
We may assume that $J=\N$. There are positive integers $m_1<m_2<\cdots$
and elements $x^{(n)}=(x^{(n)}_k)_{k\in\N}
\in B$ for $n\in\N$ such that $x^{(n)}_{m_n}\not=0$,
constructed as follows: Let $m_0:=0$. Assume that, for some $n\in\N$, elements $m_j\in\N$ and  $x^{(j)}\in B$
with $x^{(j)}_{m_j}\not=0$
have already been found for $j\in\{1,\ldots, n-1\}$
such that $m_{i-1}\leq m_i$ for all $i\in\{1,\ldots,n-1\}$
(a vacuous condition if $n=1$). Since $B\not\sub \prod_{k=1}^{m_{n-1}}E_k$,
there exists $x^{(n)}\in B\setminus \prod_{k=1}^{m_{n-1}}E_k$.\vspace{-.3mm}
Thus $x^{(n)}_{m_n}\not=0$ for some $m_n>m_{n-1}$,
which completes the construction.

Choose a balanced open $0$-neighborhood $V_n\sub E_{m_n}$ with $x^{(n)}_{m_n}\not\in V_n$.
Then
\[
W:=\Big\{(y_n)_{n\in\N}\colon (\forall n\in \N)\;\; y_{m_n}\in\frac{1}{n} V_n\Big\}
\]
is a $0$-neighborhood in~$E$ such that $x^{(n)}\not\in nW$ and thus
$B\not\sub n W$ for all $n\in\N$, contradicting the boundedness of~$B$.

(f) Since every bounded subset of~$E$ is contained in a finite partial product and~$E$
induces the product topology on the latter,
the assertions concerning quasi-completeness, sequential completeness and Mackey completeness
hold because they are well known for finite direct products.
Now assume that $E_j$ is complete for each $j\in J$, and let $(x^{(\alpha)})_{\alpha\in A}$
be a Cauchy net in $\bigoplus_{j\in J}E_j$, where $x^{(\alpha)}=(x^{(\alpha)}_j)_{j\in J }$.
By Remark~\ref{firstremsums}(b),
$(x_\alpha)_{\alpha\in A}$ also is a Cauchy net in $\prod_{j\in J}E_j$, whence
\[
x_\alpha\to x\quad\mbox{in $\,\prod_{j\in J}E_j$}
\]
for some $x=(x_j)_{j\in J}\in\prod_{j\in J} E_j$. To see that $x\in \bigoplus_{j\in J} E_j$, assume to the contrary
that $I:=\{j\in J\colon x_j\not=0\}$ is infinite. For $j\in I$, choose a closed $0$-neighborhood $V_j\sub E_j$
such that $-x_j\not\in V_j$. For $j\in J\setminus I$, set $V_j:=E_j$.
There exists $\alpha_0\in A$ such that
\[
x^{(\alpha)}-x^{(\beta)}\in \bigoplus_{j\in J}V_j
\]
for all $\alpha,\beta\geq\alpha_0$. For $j\in J$, we therefore have
\[
x^{(\alpha)}_j-x^{(\beta)}_j\in V_j
\]
for all $\alpha,\beta\geq\alpha_0$.
Passing to the limit, we see that
$x^{(\alpha)}_j-x_j\in V_j$ for all $\alpha\geq\alpha_0$. Thus
\[
x^{(\alpha_0)}_j\in x_j+V_j 
\]
and thus $x^{(\alpha_0)}_j\not=0$ for all $j\in J$, contrary to $x^{(\alpha_0)}\in \bigoplus_{j\in J}E_j$.

To see that $x^{(\alpha)}\to x$ in $\bigoplus_{j\in J} E_j=E$, let $U\sub E$ be a $0$-neighborhood;
then~$U$ contains a $0$-neighborhood~$W$ of the form
\[
W=\bigoplus_{j\in J} W_j
\]
with closed $0$-neighborhoods $W_j\sub E_j$.
There exists $\alpha_0\in A$ such that\linebreak
$x^{(\alpha)}- x^{(\beta)}\in W$
for all $\alpha,\beta\geq\alpha_0$. For each $j\in J$, we then have
\[
x^{(\alpha)}_j-x^{(\beta)}_j\in W_j\quad\mbox{for all $\alpha,\beta\geq\alpha_0$.}
\]
Passing to the limit in $\beta$, we see that $x^{(\alpha)}_j-x_j\in W_j$ for all $\alpha\geq\alpha_0$ and $j\in J$,
whence $x^{(\alpha)}-x\in E\cap \prod_{j\in J} W_j=\bigoplus_{j\in J} W_j=W\sub U$.
\end{prf}
\begin{ex}\label{Rinftysum}
By Example~\ref{theKinfty},
the direct limit topology (as a topological space) makes
$\R^{(\N)}=\dl\R^n$\vspace{-1.3mm} a locally convex space.
It is clear that
\[
\R^{(\N)}=\bigoplus_{n\in\N}\R
\]
as a vector space. Let us verify that $\R^{(\N)}=\bigoplus_{n\in\N}\R$
as a locally convex space; this will entail, e.g., that $\R^{(\N)}$ is
complete (using Lemma~\ref{firstlasum}(f)).

We have to show that the identity map $\phi\colon\R^{(\N)}\to\bigoplus_{n\in\N}\R$
is a homeo\-morphism. As the restriction  $\phi|_{\R^n}=\lambda_{\{1,\ldots,n\}}$ is continuous
for each $n\in\N$
by Remark~\ref{firstremsums}(c),
the map~$\phi$ is continuous by
Lemma~\ref{ctsonsteps}(ii). Let us write $\R e_n=\lambda_{\{n\}}(\R)$ for the $n$th
summand in $\bigoplus_{n\in\N}\R$.
Now $\phi^{-1}|_{\R e_n}$ is continuous as it can be written as the composition
\[
\R e_n\to \R^n\to \R^{(\N)}
\]
of continuous linear inclusion maps.
Since $\phi^{-1}$ is linear, Lemma~\ref{firstlasum}(a)
shows that~$\phi^{-1}$ is continuous. Thus~$\phi$ is a homeomorphism.
\end{ex}
\subsection*{Locally convex direct limits}
Many vector spaces of interest can be written as an ascending union of
locally convex spaces.
\begin{defn}
An ascending sequence $E_1\sub E_2\sub\cdots$ of locally convex spaces
$(E_n,\cO_n)$, for $n\in\N$,
is called a \emph{direct sequence of locally convex spaces}
if each inclusion map $E_n\to E_{n+1}$ is continuous and linear.
We give $E:=\bigcup_{n\in\N}E_n$ its natural vector space structure (see
Definition~\ref{dirseqTG}).
Let~$P$ be the set of all seminorms $p\colon E\to[0,\infty[$ such that
\[
p|_{E_n}\colon (E_n,\cO_n)\to [0,\infty[
\]
is continuous for each $n\in\N$.
The (not necessarily Hausdorff) locally convex vector topology $\cO_{\lcx}$ on~$E$
defined by the seminorms $p\in P$ is called the \emph{locally convex direct limit topology}
on~$E$ and $(E,\cO_{\lcx})$ is called the \emph{locally convex direct limit} of the sequence.
If we say that $E=\dl E_n$\vspace{-.5mm} \emph{as a locally convex space}, we refer
to the topology~$\cO_{\lcx}$ on~$E$.
If each $E_n$ is a Banach space (resp., a Fr\'{e}chet space), then $(E,\cO_{\lcx})$ is called
an (LB)-space (resp., an (LF)-space).
\end{defn}
\begin{rem}\label{firstremlcxDL}
(a) Note that each inclusion map $j_n\colon (E_n,\cO_n) \to (E,\cO_{\lcx})$ is continuous,
since $p(j_n(x))=(p|_{E_n})(x)$ for each $p\in P$ and $x\in E_n$.

(b)
In fact, $\cO_{\lcx}$ is the finest locally convex topology on~$E$ making each inclusion map
$j_n\colon E_n\to E$ continuous.\footnote{If $\cO$ is a locally convex vector topology on~$E$
making each $j_n$ continuous and $p$ is a continuous seminorm on $(E,\cO)$,
then $p|_{E_n}=p\circ j_n$ is continuous for each $n\in\N$ and thus $p\in P$,
whence $\cO\sub \cO_{\lcx}$.}

(c) If $p$ is a seminorm on~$E$ and $p|_{E_n}$ is continuous for some~$n$, then $p|_{E_m}$ is continuous
for all positive integers $m\leq n$. Hence, if $n_1< n_2<\cdots$ are positive integers, then
$p$ is continuous on each $E_n$ if and only if $p|_{E_{n_j}}$
is continuous for each $j$, whence the locally convex direct limits agree:
\[
\dl E_n=\dl E_{n_j}.\vspace{-1mm}
\]

(d) If $\lambda\colon E\to F$ is a linear map to a (not necessarily Hausdorff) locally convex space~$F$
such that $\lambda|_{E_n}\colon \!(E_n,\cO_n)\!\to \! F$ is continuous for all $n\in\N$, then
$\lambda\colon (E,\cO_{\lcx})\to F$ is continuous. To see this, let $V\sub F$ be a $0$-neighborhood. We find a continuous seminorm
$q$ on~$F$ such that $B^q_1(0)\sub V$. Then
$q\circ\lambda|_{E_n}$ is a continuous seminorm on $(E_n,\cO_n)$ for each $n\in\N$
and thus $p:=q\circ \lambda \in P$. Since $\lambda(B^p_1(0))\sub B^q_1(0)\sub V$
by construction, $\lambda$ is continuous at~$0$ and hence continuous.

(e) If there exists an injective linear map $\lambda\colon E\to F$ to a Hausdorff locally convex
space~$F$ such that $\lambda|_{E_n}\colon (E_n,\cO_n)\to F$ is continuous for each~$n$,
then $\lambda$ is continuous (by (d)) and thus $(E,\cO_{\lcx})$ is Hausdorff.
This simple criterion applies in many situations of interest.
Other sufficient conditions for the Hausdorff property will be
given in Propositions \ref{strictDlprops}(b) and \ref{silvahaveDL}(b).

(f) On $E=\bigcup_{n\in\N} E_n$, there also is the topology $\cO_{\DL}$
making it the direct limit topological space (see Definition~\ref{defnDLtop}).
Then $\cO_{\lcx}\sub \cO_{\DL}$ (as the map
\[
(E,\cO_{\DL})\to (E,\cO_{\lcx}),\quad x\mto x
\]
is continuous by (a) and Lemma~\ref{ctsonsteps}).
We warn the reader that $\cO_{\DL}\not=\cO_{\lcx}$ in many cases
(an example can be found in Exercise~\ref{exc-DLnotGP}). A special situation where
$\cO_{\DL}=\cO_{\lcx}$ will be
encountered in Proposition~\ref{silvahaveDL}, and is useful
for differential calculus (see Section~\ref{sec-silva-calc}).
\end{rem} 
The following lemma sheds some light on $0$-neighborhoods in locally\linebreak
convex direct limits.
\begin{lem}\label{lemnbhdlcxDL}
If $E_1\sub E_2\sub \cdots$ is a direct sequence of locally convex spaces $(E_n,\cO_n)$
and $E=\dl E_n$\vspace{-.5mm} its locally convex direct limit, then we have:
\begin{description}[(D)]
\item[\rm(a)]
A convex subset $V\sub E$ with $0\in V$ is a $0$-neighborhood in $(E,\cO_{\lcx})$
if and only if $V_n:=V\cap E_n$ is a $0$-neighborhood in $(E_n,\cO_n)$ for each $n\in\N$.
\item[\rm(b)]
A convex subset $V\sub E$ is open in $(E,\cO_{\lcx})$
if and only if $V_n:=V\cap E_n$ is open in $(E_n,\cO_n)$ for each $n\in\N$.
\item[\rm(c)]
If $(V_n)_{n\in\N}$ is a sequence of open $0$-neighborhoods $V_n$ in $(E_n,\cO_n)$,
then
\[
V:=\sum_{n\in\N}V_n:=\bigcup_{n\in\N}(V_1+\cdots+V_n)\vspace{-.8mm}
\]
is an open $0$-neighborhood in~$(E,\cO_{\lcx})$.
Such sets~$V$ form a basis of $0$-neighborhoods in $(V,\cO_{\lcx})$.
\item[\rm(d)] $\cO_{\lcx}$ is the quotient topology under the summation map
\[
q\colon \bigoplus_{n\in\N} E_n\to E,\quad (x_n)_{n\in\N}\mto \sum_{n\in \N} x_n,\vspace{-.8mm}
\]
whose domain is the locally convex direct sum.
\item[\rm(e)]
If $E_n$ is barreled for each $n\in\N$, then also $(E,\cO_{\lcx})$ is barreled.
\end{description}
\end{lem}
\begin{prf}
(d) Let $\cO_q$ be the quotient topology on~$E$.
The restriction $q|_{E_n}$ to a summand coincides with the continuous inclusion map
$(E_n,\cO_n) \to (E,\cO_{\lcx})$. Hence~$q$ is continuous, by Lemma~\ref{firstlasum}(a), and thus $\cO_{\lcx}\sub \cO_q$.
Since $q|_{E_n}\colon (E_n,\cO_n) \to (E,\cO_q)$ is continuous for each~$n$ and $\cO_q$ is locally convex,
we have $\cO_q\sub \cO_{\lcx}$ (see Remark~\ref{firstremlcxDL}(b)) and hence $\cO_q=\cO_{\lcx}$.

(c) Since $q$ is continuous and open by (d) and the sets $W:=\bigoplus_{n\in \N}V_n$ form
a basis of open $0$-neighborhoods in $\bigoplus_{n\in\N} E_n$, it follows that the
sets $V=q(W)$ form a basis of $0$-neighborhoods in~$E$.

(a) Follows from (c) since $\sum_{n\in\N} 2^{-n}V_n\sub V$.

(b) Assume each $V_n$ is open. Let $x\in V$, say $x\in E_m$. After passing to a subsequence, we may assume that $m=1$.
Then $(V-x)\cap E_n=V_n-x$ is a $0$-neighborhood in~$E_n$ for each~$n$.
By (a), $V-x$ is a $0$-neighborhood, so $V$ is a neighborhood of~$x$.

(e) If $V\sub E$ is a barrel, then $V\cap E_n$ is a barrel in~$V_n$
for each $n\in\N$. Hence $V\cap E_n$ is a $0$-neighborhood in $E_n$ (as
$E_n$ is assumed barreled), whence~$V$ is a $0$-neighborhood in~$E$, by~(a).
\end{prf}
\subsection*{Strict direct limits of locally convex spaces}
A direct sequence $E_1\sub E_2\sub\cdots$ of locally convex spaces
is called \emph{strict} if $E_{n+1}$ induces the given topology on~$E_n$
for each $n\in\N$ (see Definition~\ref{defnDLtop});
the locally convex direct limit $\dl E_n$\vspace{-.5mm} is
called a \emph{strict direct limit} in this case.
We now record some properties of these.
The following concept
is useful:
\begin{defn}
Let $E_1\sub E_2\sub\cdots$ be a direct sequence of
locally convex spaces $(E_n,\cO_n)$, with locally convex direct limit $E=\dl E_n$.\vspace{-.5mm}
We say that $(E,(E_n)_{n\in\N})$ is \emph{boundedly regular}\footnote{Often,
such direct limits are simply called \emph{regular} in the literature.}
if for every bounded subset $B\sub E$, there exists $n\in\N$
such that $B$ is contained and bounded in $(E_n,\cO_n)$.
\end{defn}
For the related notion of compact regularity,
see Definition~\ref{defncpregu}. We recall a useful fact:
\begin{lem}\label{extendnbhd}
Let $E$ be a locally convex space, $F\sub E$ be a vector subspace
and $V\sub F$ be an open, absolutely convex $0$-neighborhood. Then there exists
an absolutely convex, open $0$-neighborhood $Q\sub E$ such that $V=F\cap Q$.
\end{lem}
\begin{prf}
By the definition of the induced topology, there exists an open $0$-neighborhood $P\sub E$ such that
$F\cap P=V$. Let $W\sub P$ be an open, absolutely convex $0$-neighborhood in~$E$
and set
\[
Q:=\conv(V\cup W).
\]
Then $Q$ is absolutely convex. If $x\in Q$, then $x=tv+(1-t)w$ for some $t\in[0,1]$, $v\in V$ and $w\in W$
(see Exercise~\ref{convex-of-convex}).
Now assume that $x\in F\cap Q$. If $t=1$, then $x=v\in V$.
If $1-t\not=0$, we deduce that $w=(x-tv)/(1-t)\in F$. Thus $w\in F\cap W\sub F\cap P=V$,
whence $x=tv+(1-t)w\in V$ (as $V$ is convex). Hence $F\cap Q\sub V$ and thus $F\cap Q=V$
(as $V\sub Q$). Since $W\sub Q$, we know that $Q$ is a $0$-neighborhood in~$E$.
If $x\in Q$, write $x=tv+(1-t)w$ as before. Since~$V$ is open in~$F$ and $W$ is open in~$E$,
we find $r>1$ such that $rv\in V$ and $rw\in W$.
Then $rx=t(rv)+(1-t)(rw)\in Q$, whence $x\in Q^0$ by Lemma~\ref{baseconvex}(d).
Thus $Q\sub Q^0$, showing that $Q=Q^0$ is open.
\end{prf}
\begin{prop}\label{strictDlprops}
If $E_1\sub E_2\sub\cdots$ is a strict direct sequence of locally convex spaces
and $E:=\dl E_n$\vspace{-.5mm} its locally convex direct limit,
then we have:
\begin{description}[(D)]
\item[\rm(a)]
$E$ induces the given topology on $E_n$ for each $n\in\N$;
\item[\rm(b)]
$E$ is Hausdorff;
\item[\rm(c)]
If $E_n$ is closed in $E_{n+1}$ for each $n\in\N$, then
$(E,(E_n)_{n\in\N})$ is boundedly regular
and compactly regular;
\item[\rm(d)]
If $E_n$ is complete for each $n\in\N$, then also~$E$ is complete;
\item[\rm(e)]
If $E_n$ is closed in $E_{n+1}$
and $E_n$ is sequentially complete $($resp., quasi-complete, resp, Mackey complete$)$
for each $n\in \N$, then also $E$ is sequentially complete $($resp., quasi-complete,
resp., Mackey complete$)$.
\end{description}
\end{prop}
\begin{prf}
(a) The inclusion map $j_n\colon E_n\to E$ is continuous.
To see that~$j_n$ is a topological embedding, we may assume that $n=1$. Let $V_1\sub E_1$ be
an absolutely convex, open $0$-neighborhood. By Lemma~\ref{extendnbhd}, for each $n\in\N$,
there exists an absolutely convex open $0$-neighborhood $V_{n+1}\sub E_{n+1}$
such that $V_{n+1}\cap E_n=V_n$. Using Lemma~\ref{lemnbhdlcxDL}(b),
we see that $V:=\bigcup_{n\in\N}V_n$
is an open $0$-neighborhood in~$E$. By construction, $V\cap E_1=V_1$.
Hence $j_1$ is a topological embedding.

(b) By Lemma~\ref{TVSHausdorff}, it suffices to show that $\{0\}$ is closed in~$E$. Let $0\not=x\in E$, say $x\in E_n$;
we show that $x\not\in\wb{\{0\}}$. Let $V\sub E_n$ be an open $0$-neighborhood such that $x\not\in V$.
By (a), there exists an open $0$-neighborhood $W\sub E$ such that $E_n\cap W=V$ and thus $x\not\in W$.
Then $x-W$ is an open neighborhood of~$x$ in~$E$ such that $0\not\in x-W$
and thus $\wb{\{0\}}\sub E\setminus (x-W)$.

(c) By (a), we need only show that every bounded subset $B\sub E$ is contained in~$E_n$ for some~$n$.
If not, then $B\setminus E_n\not=\emptyset$ for each~$n$.
Pick $x_1\in B\setminus E_1$ and $m_1\in\N$ such that $x_1\in E_{m_1}$.
Recursively, find $x_{n+1}\in B\setminus E_{m_n}$ and $m_{n+1}>m_n$ such that
$x_{n+1}\in E_{m_{n+1}}$. After passing to the subsequence $E_1, E_{m_1}, E_{m_2},\ldots$,
we may assume that $m_n=n+1$ for all $n\in\N$.
Let $V_1\sub E_1$ be an open, absolutely convex $0$-neighborhood.
Inductively, for each $n\in\N$ we find an open, absolutely convex $0$-neighborhood $V_{n+1}\sub E_{n+1}$
with
\begin{equation}\label{ints-prop}
V_{n+1}\cap E_n=V_n
\end{equation}
and
\begin{equation}\label{notconta}
x_n\not\in (n+1)V_{n+1}.
\end{equation}
In fact, using Lemma~\ref{extendnbhd}, we find a balanced open $0$-neighborhood $V_{n+1}\sub E_{n+1}$
such that (\ref{ints-prop}) holds.
Since $E_n$ is closed in $E_{n+1}$ and $x_{n+1}\not\in E_n$,
there is an open, absolutely convex $0$-neighborhood $W\sub
E_{n+1}$ such that $(x_n+(n+1)W)\cap E_n=\emptyset$
and thus
\begin{equation}\label{henceouts}
x_n\not\in E_n+(n+1)W=(n+1)(E_n+W).
\end{equation}
After replacing $V_{n+1}$ with $V_{n+1}\cap (E_n+ W)$,
we may assume that $V_{n+1}\sub E_n+ W$.
Hence~(\ref{notconta}) holds, by~(\ref{henceouts}).

Now $V:=\bigcup_{n\in \N} V_n$ is an open $0$-neighborhood in~$E$ such that $V\cap E_{n+1}=V_{n+1}$
for each $n\in\N$ and hence $x_n\not\in (n+1)V$, entailing that~$B$ is not a subset of $(n+1)V$.
This contradicts the boundedness of~$B$.

(d) Let $(x_\alpha)_{\alpha\in A}$
be a Cauchy net in~$E$. Let~$\cU$ be the set of all $0$-neighborhoods in~$E$.
We claim that there exists $m\in\N$
such that, for all $\alpha\in A$ and $W\in\cU$,
there is $\beta\geq\alpha$ such that
\[
x_\beta\in E_m+W.
\]
If this is true, then
\[
A_W:=\{\alpha\in A\colon x_\alpha\in E_m+W\}
\]
is cofinal in~$A$ and thus
\[
M:=\{(W,\alpha)\in\cU\times A\colon \alpha\in A_W\}
\]
becomes a directed set if we write $(W_1,\alpha_1)\leq (W_2,\alpha_2)$ if and only if
$W_2\sub W_1$ and $\alpha_1\leq \alpha_2$.
For $a=(W,\alpha)\in M$, pick $y_a\in E_m$ and $w_a\in W$ such that
\begin{equation}\label{thusconvsu}
x_\alpha=y_a+w_a.
\end{equation}
Then $(y_a)_{a\in M}$ is a Cauchy net in~$E_m$.
In fact, if $U\sub E_m$ is an open $0$-neighborhood, (a) provides $Q\in\cU$ such that $U=E_m\cap Q$.
Let $P\in\cU$ be such that $P-P+P\sub Q$ and $\gamma\in A$ be such that
\[
x_\alpha-x_\beta\in P\quad\mbox{for all $\,\alpha,\beta\geq \gamma$.}
\]
We may assume that $\gamma\in A_P$. For all $a,b\geq (P,\gamma)$ in~$M$, say $a=(V,\alpha)$ and $b=(W,\beta)$,
we then have  
\[
y_a-y_b=x_\alpha-x_\beta-w_a+w_b\in P-V+W\sub Q
\]
and thus $y_a-y_b\in E_m\cap Q=U$. Let $y\in E_m$ be the limit of $(y_a)_{a\in M}$.
Then $y_a\to y$ also in~$E$. Given $W\in\cU$, let $\beta\in A_W$.
Since $w_a\in W$ for $a\geq (W,\beta)$, the net $(w_a)_{a\in M}$ converges to~$0$ in~$E$.
In view of (\ref{thusconvsu}), the subnet
\[
(x_\alpha)_{(\alpha,V)\in M}
\]
of $(x_\alpha)_{\alpha\in A}$
(and hence also the Cauchy net $(x_\alpha)_{\alpha\in A}$) therefore converges to~$y$.

To prove the claim, assume it was wrong.
Then, for each $m\in\N$,
there exist $\alpha_m\in A$ and $W_m\in\cU$
such that
\begin{equation}\label{no-inters}
x_\alpha\not \in E_m+W_m\quad\mbox{for all $\alpha\geq\alpha_m$.}
\end{equation}
We may assume that $W_m$ is absolutely convex for each~$m$
and $W_1\supseteq W_2\supseteq\cdots$.
Let~$W$ be the convex hull of
\[
\bigcup_{n\in\N} (W_n\cap E_n).
\]
Then
\begin{equation}\label{set-smaller}
E_m+W\sub E_m+W_m
\end{equation}
for each  $m\in\N$. In fact, if $x\in W$, then there is $(x_n)_{n\in\N}\in\prod_{n\in\N} (W_n\cap E_n)$ and
a sequence $(r_n)_{n\in\N}\in
\R^{(\N)}$ of real numbers $r_n\geq 0$ with sum~$1$ such that
\[
x=\sum_{n=1}^\infty r_n x_n=\sum_{n=1}^mr_nx_n+\sum_{n=m+1}^\infty r_n x_n
\]
(see Exercise~\ref{convex-of-convex}(b)).
Since $E_n\sub E_m$ for $n\in\{1,\ldots, m\}$, the first partial sum is an element of~$E_m$.
Since $W_n\sub W_m$ if $n\geq m$ and~$W_m$ is absolutely convex, the second partial sum
is in~$W_m$. Hence $x\in E_m+W_m$ indeed.

By definition of a Cauchy net, we find $\gamma\in A$ such that
\begin{equation}\label{diff-in}
x_\alpha- x_\beta\in W
\end{equation}
for all $\alpha,\beta\geq \gamma$. Now $x_\gamma\in E_{m_0}$ for some $m_0\in\N$.
Since $A$ is directed, we find $\alpha\in A$ such that $\alpha\geq \gamma$ and $\alpha\geq\alpha_{m_0}$.
Using (\ref{diff-in}), we obtain
\[
x_\alpha =x_\gamma+ (x_\alpha-x_\gamma)\in E_{m_0}+W.
\]
But $x_\alpha\not\in E_{m_0}+W_{m_0}$ by (\ref{no-inters}) and thus
$x_\alpha\not\in E_{m_0}+W$ (by (\ref{set-smaller})),
which is absurd.

(e) is immediate from (c) and (a).
\end{prf}
\begin{rem}
A family $(E_\alpha)_{\alpha\in A}$ of locally convex spaces $(E_\alpha,\cO_\alpha)$ is called
a \emph{direct system} if $A$ is a directed set,
$E_\alpha\sub E_\beta$ for all $\alpha,\beta\in A$ with $\alpha\leq\beta$
and the inclusion map $j_{\beta,\alpha}\colon E_\alpha\to E_\beta$ is continuous linear.
If each $j_{\beta,\alpha}$ is, moreover, a topological embedding,
then the direct system is called \emph{strict}.
Give $E:=\bigcup_{\alpha\in A}E_\alpha$ the vector space
structure which turns each $E_\alpha$ into a vector subspace
and let $P$ be the set of all seminorms~$q$ on~$E$ such that $q|_{E_\alpha}$ is continuous
on $(E_\alpha,\cO_\alpha)$ for each $\alpha\in A$.
The locally convex vector topology $\cO_{\lcx}$ on~$E$ defined by the seminorms $p\in P$
is called the \emph{locally convex direct limit topology}, and we say that
$E=\dl E_\alpha$\vspace{-.5mm} as a locally convex space when
this topology is used.

If $q$ is a seminorm on~$E$ and $q|_{E_\beta}$ is continuous for some $\beta\in A$,
then $q|_{E_\alpha}=q|_{E_\beta}\circ j_{\beta,\alpha}$ is continuous for all $\alpha\leq\beta$.
Hence, if $(E_\alpha)_{\alpha\in A}$ admits a cofinal subsequence $(E_{\alpha_n})_{n\in\N}$
(as in Remark~\ref{rem-cof-subs}(b)), then
\[
\dl E_\alpha=\dl E_{\alpha_n}\vspace{-1mm}
\]
as locally convex spaces.
In this book, we shall only work with direct systems admitting cofinal subsequences,
to avoid additional complications.
\end{rem}
\begin{ex}
For each compact subset $K\sub\R$, let $C^\infty_K(\R)\sub C^\infty(\R)$
be the closed vector subspace of all smooth functions $f\colon\R\to\R$ with
support $\supp(f)\sub K$.
Then the $C_{[-n,n]}^\infty(\R)$ form a cofinal subsequence
for $n\in\N$, whence the locally convex direct limit
\[
C^\infty_c(\R):=\dl C^\infty_K(\R)=\dl C^\infty_{[-n,n]}(\R)
\]
(the so-called space of test functions) is a strict (LF)-space
and hence enjoys all of the properties described in
Proposition~\ref{strictDlprops}.
\end{ex}
\subsection*{Compact operators and Silva spaces}
Another class of locally convex direct limits is important for
analysis and infinite-dimensional Lie theory,
the
Silva spaces.
Recall that a linear map
\[
\alpha\colon E\to F
\]
between Banach spaces is called a \emph{compact operator}
if the closure $\wb{\alpha(B)}$ is compact in~$F$
for each bounded subset $B\sub E$.
Choosing $B:=B^E_1(0)$, we see that any compact operator
is bounded and hence continuous.
\begin{defn}
A locally convex space~$E$ is called a \emph{Silva space}
if $E$ is the locally convex direct limit $E=\dl E_n$\vspace{-.5mm}
of some direct sequence $E_1\sub E_2\sub\cdots$ of Banach spaces
such that each inclusion map $E_n\to E_{n+1}$ is a compact operator.
\end{defn}
The following proposition compiles properties of Silva spaces.
\begin{prop}\label{silvahaveDL}
Let $E_1\sub E_2\sub\cdots$ be a direct sequence of Banach spaces
such that all inclusion maps $E_n\to E_{n+1}$ are compact
operators. Consider the locally convex direct limit $E=\dl E_n$,\vspace{-.5mm}
endowed with the locally convex direct limit topology~$\cO_{\lcx}$.
Then the following holds:
\begin{description}[(D)]
\item[\rm(a)]
$\cO_{\lcx}$ coincides with the direct limit topology $\cO_{\DL}$ on
$E=\bigcup_{n\in\N}E_n$;
\item[\rm(b)]
$(E,(E_n)_{n\in\N})$ is Hausdorff and compactly regular;
\item[\rm(c)]
$(E,(E_n)_{n\in\N})$ is boundedly regular;
\item[\rm(d)]
Every bounded subset of~$E$ is relatively compact;
\item[\rm(e)]
$E$ is reflexive;
\item[\rm(f)]
$E'_b$ is a Fr\'{e}chet space and $E$ is complete.
\item[\rm(g)]
$E$ is a $k_\omega$-space.
\end{description}
\end{prop}
\begin{prf}
(a)  Since each inclusion map $E_n\to (E,\cO_{\lcx})$ is continuous, we have $\cO_{\lcx}\sub \cO_{\DL}$.
The direct sequence $E_1\sub E_2\sub\cdots$
is compactifying (in the sense of Definition~\ref{defnDLtop})
as each inclusion map $E_n\to E_{n+1}$ is a compact operator and hence takes balls
to relatively compact sets.
Hence $\cO_{\DL}$ is a locally convex vector topology
(see Corollary~\ref{cp-DL-TVS}).
As the inclusion maps $E_n\to (E,\cO_{\DL})$ are continuous, we deduce
that $\cO_{\DL}\sub\cO_{\lcx}$
and hence $\cO_{\lcx}=\cO_{\DL}$.

(b) By Proposition~\ref{henceSilvacoreg}, $(E,(E_n)_{n\in\N})$ is compactly regular
and Hausdorff (being locally $k_\omega$).

(c) After replacing the norm $\|\cdot\|_n$ on each~$E_n$ by a suitable multiple,
we may assume that
\begin{equation}\label{inflate}
\|x\|_{n+1}\leq \|x\|_n\quad\mbox{for all $n\in\N$ and $x\in E_n$.}
\end{equation}
Let $B\sub E$ be bounded. Then $B\sub E_n$ for some $n\in\N$,
because otherwise we could find an element $x_n\in B\setminus E_n$
for each $n\in\N$. Since $B$ is bounded,
\[
\frac{1}{n}x_n\to 0\quad\mbox{as $\;n\to\infty$,}
\]
entailing that $K:=\{0\}\cup\{\frac{x_n}{n}\colon n\in\N\}$ is compact in~$E$
and hence a compact subset of some~$E_n$, by~(b). This contradicts $x_n\not\in E_n$.

Thus $B\sub E_n$ and we may assume that $B\sub E_1$. Then $B$ is bounded in~$E_n$
for some $n\in\N$, because otherwise we can find $x_n\in B$ such that
\begin{equation}\label{secstepunb}
\|x_n\|_n\geq n,
\end{equation}
for each $n\in\N$. Then $L:=\{0\}\cup\{x_n/\sqrt{n}\colon n\in\N\}$ is a compact subset of~$E$
and hence a compact subset of~$E_n$ for some $n\in\N$.
Hence~$L$ is bounded in~$E_n$. However, for each $m\geq n$, by (\ref{inflate}) and (\ref{secstepunb})
we have
\[
\|x_m/\sqrt{m}\|_n\geq \|x_m/\sqrt{m}\|_m=\|x_m\|_m/\sqrt{m}\geq \sqrt{m},
\]
and these numbers form an unbounded set
(contradiction).

(d) Let $B\sub E$ be bounded. By (c), $B$ is a bounded subset of~$E_n$ for some $n\in \N$.
As the inclusion map $E_n\to E_{n+1}$ is a compact operator, this implies that the closure~$K$ of $B$ in
$E_{n+1}$ is compact. Then $K$ is also compact in~$E$.

(e) Each $E_n$ is barreled (see Proposition~\ref{frechbarr}), whence~$E$ is barreled by Lemma~\ref{lemnbhdlcxDL}(e).
By~(d), the closure~$K$ of a bounded subset $B\sub E$ in~$E$ is compact.
Hence~$K$ is also compact in~$E_w$ and thus~$E$ is reflexive, by Proposition~\ref{charrefl}.

(f) and (g). By (d), we have $E'_b=E'_c$. Using norms as in the proof of~(c), we have $B^{E_n}_r(0)\sub B^{E_{n+1}}_r(0)$
for each $n\in\N$ and $r>0$. Let $K_n$ be the closure of $B^{E_n}_n(0)$ in $E_{n+1}$.
Then each~$K_n$ is compact and $K_1\sub K_2\sub\cdots$. If $K\sub E$ is compact, then $K$ is contained
in $B^{E_n}_k(0)$ for some $n,k\in\N$ (by (b)), and thus $K\sub K_m$ with $m:=\max\{n,k\}$.
Hence $E$ is hemicompact (as in \ref{hemic}),
whence $E'_b=E'_c$ is metrizable (by Lemma~\ref{sammelsu}(c)).
Therefore $(E'_b)'_b$ is complete (by Proposition~\ref{opspacecompl}) and hence also~$E$, as $\eta_E\colon E\to (E'_b)'_b$
is an isomorphism of  topological vector spaces.
Since $(E,\cO_{\DL})$ is a direct limit of the Banach spaces~$E_n$ which are $k$-spaces,
also~$E$ is a $k$-space (see Exercise~\ref{exc-perm-k}), entailing that~$E'_b$ is complete (see Proposition~\ref{opspacecompl})
and hence a Fr\'{e}chet space. Being a hemicompact $k$-space, $E$ is a $k_\omega$-space.
\end{prf}
\begin{rem}
Using that the dual $\alpha'\colon Y'\to X'$ of a compact operator $\alpha\colon X\to Y$
is again compact, one can show that the dual $E'_b$ of a Silva space is a so-called
Fr\'{e}chet--Schwartz space, i.e., a projective limit of a projective sequence
of Banach spaces whose bonding maps are compact operators.
Hence $E$ is isomorphic to the dual of a Fr\'{e}chet--Schwartz space
(and this property characterizes Silva spaces).
Silva spaces are therefore called (DFS)-spaces in part of the literature.
\end{rem}
\subsection*{Incomplete quotients}
We describe an (LF)-space
which is not complete. Yet, it is a quotient of a complete locally convex space.
\begin{numba}\label{setnonc}
Consider the Fr\'{e}chet space
$F:=\R^\N$ of all sequences of real numbers.
Then
\[
E_n:=\{(x_n)_{n\in\N}\in F^\N\colon (\forall m\geq n)\,x_m\in\ell^1\}\cong F^{n-1}\times (\ell^1)^{\{n,n+1,\ldots\}}
\]
is a Fr\'{e}chet space for each $n\in \N$.
The locally convex direct limit
\[
E:=\dl\,E_n=\bigcup_{n\in\N}E_n\vspace{-1.3mm}
\]
is Hausdorff, as the inclusion map
\[
\phi\colon E\to F^\N
\]
is continuous (being linear, and continuous on~$E_n$ for each $n\in\N$).
\end{numba}
\begin{lem}\label{pre-badq}
In the situation of~\emph{\ref{setnonc}},
$\phi$ is a topological embedding.
Hence $E$ is a proper, dense vector subspace
of $F^\N$ and carries the induced topology,
whence~$E$ is a non-complete metrizable locally convex space.
\end{lem}
\begin{prf}
We have to show that each $0$-neighborhood $U$
in the locally convex direct limit~$E$ also
is a $0$-neighborhood in the topology induced by $F^\N$ on~$E$.
After shrinking~$U$, we may assume that~$U$ is closed in~$E$
and absolutely convex.
Since $U\cap E_1$
is a $0$-neighborhood in $E_1=(\ell^1)^{\N}$,
there exists $n\in\N$ and $r>0$
such that
\[
V:=(B^{\ell^1}_r(0))^n\times (\ell^1)^{\{n+1, \ldots\}} \sub U.
\]
For each $m>n$, the intersection $U\cap E_m$
is closed in~$E_m$ and hence contains the closure
of~$V$ in~$E_m$. The latter contains
\[
V_m:=(B^{\ell^1}_r(0))^n\times F^{m-n-1}\times (\ell^1)^{\{m, \ldots\}}, 
\]
as $\ell^1$ is dense in~$F$.
Thus
\[
P:= (B^{\ell^1}_r(0))^n\times E=\bigcup_{m>n}V_m\sub U.\vspace{-.8mm}
\]
Since $U\cap E_{n+1}$ is a $0$-neighborhood in~$E_{n+1}$,
there exist $s>0$, $k>n$ and a $0$-neighborhood $W\sub F$
such that
\[
Q:=W^n\times (B^{\ell^1}_s(0))^{k-n}\times (\ell^1)^{\{k+1, \ldots\}}\sub U\cap E_{n+1}\sub U.
\]
Then
\[
U\supseteq \frac{1}{2}P+\frac{1}{2}Q\supseteq \left(\frac{1}{2}W^n\right)\times E
=E\cap \left(\left(\frac{1}{2}W^n\right)\times F^\N\right),
\]
showing that~$U$ is a $0$-neighborhood in the topology induced by $F^{\N}$
on~$E$.
\end{prf}
\begin{ex}\label{bad-quot}
In the situation of~\ref{setnonc},
$S:=\bigoplus_{n\in\N}E_n$ is a complete locally convex space
(by Lemma~\ref{firstlasum}(f) or Lemma~\ref{strictDlprops}(d)).
The map
\[
q\colon \bigoplus_{n\in\N}E_n\to E,\quad (x_n)_{n\in\N}\mto\sum_{n\in\N}x_n\vspace{-.8mm}
\]
is continuous linear and~$E$ carries the quotient topology with respect to~$q$,
by Lemma~\ref{lemnbhdlcxDL}(d). Thus $\ker(q)$ is a closed vector subspace of the complete
locally convex space~$S$ (which is a strict (LF)-space)
such that $S/\ker(q)\cong E$ fails to be complete.
Moreover, $S/\ker(q)\cong E$ is metrizable and barreled,
by Lemmas~\ref{pre-badq} and \ref{lemnbhdlcxDL}(e).
\end{ex}
\subsection*{Exercises for Section~\ref{sec-appDLvec}}

\begin{small}

\begin{exer}\label{convex-of-convex}
Let $E$ be a vector space. 
\begin{description}[(D)]
\item[(a)]
Show that $\conv(A\cup B)=\{ta+(1-t)b\colon t\in [0,1], a\in A, b\in B\}$
holds for all
non-empty convex subsets $A,B\sub E$.
\item[(b)]
Let $(W_j)_{j\in J}$ be a family of
non-empty,  convex subsets of~$E$.
Show that the convex hull of $\bigcup_{j\in J} W_j$ coincides with the set of all
elements $x\in E$ of the form
\[
x=\sum_{j\in J} r_jw_j,
\]
where $w_j\in W_j$ for $j\in J$ and $(r_j)_{j\in J}\in \R^{(J)}$ such that
$0\leq r_j$ for all $j\in J$ and $\sum_{j\in J} r_j=1$.
\end{description}
\end{exer}

\begin{exer}\label{conv-modular}
Let $E$ be a vector space, $F\sub E$ be a vector subspace and
$V\sub F$ as well as $W\sub E$ be convex subsets
such that $W\cap F\sub V$. Show that $F\cap \cnv(V\cup W)=V$.
\end{exer}

\begin{exer}\label{exc-dlnotmetr}
Let $E_1\sub E_2\sub\cdots$ be a strict direct sequence
of Fr\'{e}chet spaces such that $E_n\not=E_{n+1}$ for all $n\in \N$.
Show that the locally convex direct limit $E:=\dl E_n$\vspace{-.7mm}
is not metrizable.\\[2mm]
[Hint: If there was a metric $d$ on~$E$ defining its topology,
we pick $x_n\in E_{n+1}\setminus E_n$ for $n\in\N$
and may assume that $d(x_n,0)<\frac{1}{n}$, after replacing
$x_n$ with a suitable scalar multiple. Then $x_n\to 0$ as $n\to\infty$,
whence $\{x_n\colon n\in\N\}$ is a bounded subset of~$E$.
Derive a contradiction.]
\end{exer}

\begin{exer}\label{exc-when-finest}
Using Proposition~\ref{onlyyou},
show:
\begin{description}[(D)]
\item[(a)]
If $F_1\sub F_2\sub\cdots$
is an ascending sequence of finite-dimensional
vector spaces such that each inclusion map
$F_n\to F_{n+1}$ is linear,
then the locally convex direct limit topology on
$\dl\,F_n=\bigcup_{n\in\N}F_n$
is the finest locally convex vector topology
(as in Example~\ref{finest-lcx-vec}).
\item[(b)]
If $(E_n)_{n\in\N}$
is any sequence of finite-dimensional vector spaces,
then the locally convex direct sum $\bigoplus_{n\in\N}E_n$
is endowed with the finest locally convex vector topology.
\end{description}
\end{exer}

\begin{exer}\label{exc-maps-between-finest}
Let $E$ and $F$ be vector spaces of countable dimension.
Endow $E$ and $F$ with the finest locally convex vector topology.
Show that every linear map $\lambda\colon E\to F$ is continuous.
Every bijective linear map $\lambda\colon E\to F$
is an isomorphism of topological vector spaces.
\end{exer}
\end{small}
\section{Notes and comments on Appendix~\ref{chaplcx}}\label{notes-appB}
Most of the results of this appendix
are well known and can
be found in textbooks and monographs like
\cite{Bou87}, \cite{Jr81}, \cite{Ru91}, and \cite{ShW99};
we refrain from individual references for such facts,
and do not claim originality.
In particular, we profited from~\cite{Bou87}.

The geometric approach to the Hahn--Banach
theorems goes back to Mazur~\cite{Maz33}.
Our discussion of three-space properties
in Lemma~\ref{extprp}
was stimulated by the discussion of extension properties
of topological groups in \cite[5.38(e)]{HR79}
(cf.\ also \cite{Grv50}). We mention a related general fact:
If $N$ is a closed normal subgroup of a topological
group~$G$ and both~$N$ and~$G/N$ are complete,
then~$G$ is complete (cf.\ \cite[Thm.~12.3(a)]{RD81}).

The discontinuity of the evaluation map
$\ve\colon E'\times E\to\K$
for non-normable~$E$ (Proposition~\ref{evalctsnormable})
is taken from~\cite{KM97}, whose authors
refer to~\cite{Ms63}.
See \cite{Bou87} and
\cite{Gr73} for further discussions of hypocontinuity
(cf.\ also \cite{Tr67}).
As every hypocontinuous bilinear map is sequentially continuous,
it behaves to some extent like a differentiable map
(see, e.g., Exercise~\ref{seq-cts-Leib},
also \cite{Th96}, \cite{Gl07e}, and~\cite{Gl08b}).
The fairly technical Exercise~\ref{exer-wopbd}
is used in Exercice~\ref{exc-C1notFC1} to obtain an example
of a mapping $f\colon L^1[0,1]\to L^1[0,1]$
which is $C^1$ but fails to be continuously Fr\'{e}chet
differentiable in the classical sense.

Our discussion of locally convex direct
sums and strict direct limits is based on~\cite{Bou87}.
A prime example of a strict direct limit is
the space $C^\infty_c(\R)$ of test functions on~$\R$
used in the theory of distributions~\cite{Schw57}.
In the context of infinite-dimensional Lie theory,
certain generalizations are needed like
spaces of vector-valued compactly supported smooth functions
on $\sigma$-compact finite-dimensional smooth manifolds,
and spaces of compactly supported smooth vector fields thereon.
Any such
admits a linear topological embedding onto a closed
(and complemented) vector subspace
of a suitable locally convex direct sum, which can be exploited
in differential calculus (see, e.g., Chapter~\ref{ch:dirlim}). 
For this reason,
the specialized, but simple
concept of a locally convex direct sum
is of interest.

Silva spaces were first studied by
J. Sebast\~{a}o e Silva \cite{SeS55}
and Ra\u{\i}kov \cite{Rai59},
cf.\ also \cite{Gr54}); we profited, moreover,
from~\cite{Fl71}.
Examples of Silva spaces include spaces of germs of
$\K$-valued analytic functions around $0$ (or any non-empty compact subset)
in~$\K^n$, as well as
spaces of real analytic functions or
real analytic vector fields on a compact real analytic manifold.
Further information on compact regularity and
bounded regularity of (LF)-spaces
can be found in~\cite{We03}.

The example of a non-complete, metrizable (LF)-space
discussed in \ref{setnonc} and Lemma~\ref{pre-badq}
is based on~\cite[Ex.~1]{SN81}.
Further results concerning barreled locally convex spaces
can be found in~\cite{PCB87};
for refined results concerning
Fr\'{e}chet spaces, see~\cite{MV97}.

%
\def\date{20.3.19}

\chapter{Topological Groups} \mlabel{ch:topgrp} 

\vspace{5mm}
In this appendix we collect some very basic material on topological 
groups which is used to develop the theory of Lie groups. 


\section{Group topologies from local data} \mlabel{sec:a.2}

\begin{defn} \mlabel{def:t.1.1} 
A {\it topological group} $G$ is a Hausdorff space $G$ 
 endowed with a group structure such that the multiplication 
map and the inversion map are continuous. 

If $G$ and $H$ are topological groups, then a group homomorphism 
\break $\phi \: G \to H$ is called a {\it morphism of topological groups} 
if $\phi$ is continuous. 
\end{defn}

The following lemma describes how to construct a group topology on a group from 
a filter basis of subsets which then becomes a filter basis of 
identity neighborhoods for the group topology. 

\begin{defn} \mlabel{def:4.2.1}
Let $X$ be a set. A subset ${\cal F} \subeq \bP(X)$ 
(the set of all subsets of $X$) is called a {\it filter
basis}\index{filter!basis} if the following conditions are satisfied: 
\begin{description}
\item[(F1)] ${\cal  F} \not= \eset$.
\item[(F2)] Each set $F \in {\cal F}$ is nonempty.
\item[(F3)] $A, B \in {\cal F}\Rarrow (\exists C \in {\cal F})\ C \subeq A \cap B.$
\end{description}
\end{defn}

\begin{lem} \mlabel{lem:loctop} 
Let $G$ be a group and ${\cal F}$ be a filter basis of
subsets of $G$ satisfying $\bigcap {\cal F} = \{\be\}$ and 
\begin{description}
\item[\rm(U1)] $(\forall U \in {\cal F})(\exists V \in {\cal F})\ VV \subeq U. $
\item[\rm(U2)] $(\forall U \in {\cal F})(\exists V \in {\cal F})\ V^{-1} \subeq U. $
\item[\rm(U3)] $(\forall U \in {\cal F})(\forall g \in G) (\exists V \in
{\cal F})\ gVg^{-1} \subeq U. $
\end{description}

\nin Then there exists a unique group topology on $G$ such that
${\cal F}$ is a basis of $\be$-neighborhoods in $G$, i.e., 
a subset $U \subeq G$ is a $\be$-neighborhood if and only if it is a 
superset of an element of $\cF$. 
This topology is
given by 
$$\{ U \subeq G \: (\forall g \in U) (\exists V \in {\cal F})\ gV
\subeq U\}.$$
\end{lem}

\begin{prf} Let 
$$\tau  := \{ U \subeq G \: (\forall g \in U) (\exists V \in {\cal F})\ gV
\subeq U\}.$$
First we show that $\tau$ is a topology. Clearly $\eset, G \in \tau$. 
Let $(U_j)_{j \in J}$ be a family of elements of $\tau$ and 
$U := \bigcup_{j \in J} U_j$. For each $g \in U$, there exists a $j_0 \in J$ 
with $g \in U_{j_0}$ and a $V \in {\cal F}$ with $gV \subeq U_{j_0} \subeq U$. 
Thus $U \in \tau$ and we see that $\tau$ is stable under arbitrary 
unions. 

If $U_1, U_2 \in \tau$ and $g \in U_1 \cap U_2$, there exist 
$V_1, V_2 \in {\cal F}$ with $g V_i \subeq U_i$, $i=1,2$. 
Since $\cF$ is a filter 
basis, there exists $V_3 \in \cF$ with $V_3 \subeq V_1 \cap V_2$, and 
then $gV_3 \subeq U_1 \cap U_2$. We conclude that $U_1 \cap U_2 \in \tau$, 
and hence that $\tau$ is a topology on $G$. 

We claim that the interior $U^0$ of a subset $U \subeq G$ is given by 
$$ U_1 := \{ u \in U \: (\exists V \in \cF)\ u V \subeq U\}. $$
In fact, if $u \in U_1$ and there exists a $V \in \cF$ with 
$uV \subeq U$, then we pick a $W \in \cF$ with $WW \subeq V$ and obtain 
$uW W \subeq U$, so that $uW \subeq U_1$. Hence $U_1$ is open, and it 
clearly is the largest open subset contained in $U$, i.e., $U_1 = U^0$. 
It follows in particular that $U$ is a neighborhood of $g$ if and only if 
$g \in U_1$. It follows in particular that $\cF$ is a basis of the 
neighborhood filter of $\be$. The property $\bigcap \cF = \{\be\}$ 
implies that for $x \not =y$ there exists $U \in \cF$ with 
$y^{-1} x \not \in U$. For $V \in \cF$ with $VV \subeq U$ and 
$W \in \cF$ with $W^{-1} \subeq V$ we then obtain 
$y^{-1} x \not \in VW^{-1}$, i.e., 
$xW \cap y V = \eset$. Thus $(G,\tau)$ is a Hausdorff space. 

To see that $G$ is a topological group, we have to verify that the 
map 
$$ f \: G \times G \to G, \quad (x,y) \mapsto xy^{-1} $$
is continuous. So let $x,y \in G$, $U \in \cF$ and pick 
$V \in \cF$ with $yVy^{-1} \subeq U$ and 
$W \in \cF$ with $WW^{-1} \subeq V$. Then 
$$ f(xW, yW) = xW W^{-1} y^{-1} = xy^{-1} y(WW^{-1}) y^{-1} 
\subeq xy^{-1} yV y^{-1} \subeq xy^{-1} U, $$
implies that $f$ is continuous in $(x,y)$. 
\end{prf}

Before we turn to Lie group structures, it is illuminating to first 
consider the topological variant of the following theorem. 

\begin{lem} \mlabel{lem:toplocglob} 
Let $G$ be a group and 
$U = U^{-1}$ a symmetric subset containing~$\be$. 
We further assume that $U$ carries a topology for which 
\begin{description}
\item[\rm(T1)] $D := \{ (x,y) \in U \times U \: xy\in U\}$ is an open subset 
and the group multiplication $m_U \: D \to U, (x,y) \mapsto xy$ is continuous, 
\item[\rm(T2)] the inversion map $\eta_U \: U \to U, u \mapsto u^{-1}$ is
continuous, and,  
\item[\rm(T3)] for each $g \in G$, there exists an open 
$\be$-neighborhood $U_g$ in $U$ with $c_g(U_g) \subeq U$, such that the conjugation map 
$$ c_g \: U_g \to U, \quad x \mapsto gxg^{-1} $$
is continuous. 
\end{description}
Then there exists a unique group topology on $G$ 
for which the inclusion map $U \into G$ is a homeomorphism 
onto an open subset of $G$. 

If, in addition, $U$ generates $G$, then {\rm(T1/2)} imply {\rm(T3)}.
\end{lem}

\begin{prf} First we consider the filter basis 
${\cal F}$ of $\be$-neighborhoods in $U$. 
Then (T1) implies (U1), (T2) implies (U2), and (T3) implies (U3). 
Moreover, the assumption that $U$ is Hausdorff implies that 
$\bigcap {\cal F} = \{\be\}$. Therefore Lemma \ref{lem:loctop} implies that $G$
carries a unique structure of a (Hausdorff) topological group for which 
${\cal F}$ is a basis of $\be$-neighborhoods. 

We claim that the inclusion map $U \to G$ is an open embedding. 
So let $x \in U$. Then 
$$U_x := U \cap x^{-1}U = \{ y\in U \: (x,y) \in D \} $$
is open in $U$ and $\lambda_x$ restricts to a continuous map 
$U_x \to U$ with image $U_{x^{-1}}$. Its 
inverse is also continuous. Hence 
$\lambda_x^U \: U_x \to U_{x^{-1}}$ is a homeomorphism. 
We conclude that a basis for the neighborhood filter of 
$x \in U$ consists of sets of the form $xV$, $V$ a neighborhood of $\be$, 
and hence that the inclusion map $U \into G$ is an open embedding. 

Suppose, in addition, that $G$ is generated by $U$. 
For each $g \in U$, there exists an open $\be$-neighborhood $U_g$ 
with $gU_g \times \{g^{-1}\} \subeq D$. 
Then $c_g(U_g) \subeq U$, and the continuity of $m_U$ 
implies that $c_g\res_{U_g} \: U_g\to U$ is continuous. 
Hence, for each $g \in U$, the conjugation $c_g$ is continuous 
in a neighborhood of $\be$. 
Since the set of all these
$g$ is a submonoid of $G$ containing $U$, it contains $U^n$ for each
$n \in \N$, hence all of $G$ because $G$ is generated by $U = U^{-1}$. 
Therefore (T3) follows from (T1) and (T2). 
\end{prf}

\section{Extending local homomorphisms of topological groups} \mlabel{sec:a.3}

In the following we call a connected topological group $G$ {\it simply connected} 
if each connected group covering of $G$ is trivial. Although we won't need it 
in the context of Lie groups, this concept has the advantage that it also  
applies to groups which are not necessarily locally arcwise connected. 
If $G$ is a Lie group, then Corollary~\ref{cor:unicov} implies that 
$G$ is simply connected as a topological group if and only  if $\pi_1(G)$ 
vanishes, i.e., if and only if $G$ is simply connected as a topological space. 

\begin{prop} \mlabel{prop:lochomo-ext} 
Let $G$ be a connected simply connected topological group and 
$H$ a group. Let $V$  be an open symmetric connected identity neighborhood in
$G$ and $f \: V \to H$ a function with 
$$ f(xy) = f(x)f(y) \quad \hbox{ for } \quad x, y, xy \in V. $$
Then there exists a unique group homomorphism extending $f$. 
If, in addition, $H$ is a topological group and $f$ is continuous,
then its extension is also continuous. 
\end{prop}

\begin{prf} We consider the group $G \times H$ and the subgroup 
$S \subeq G \times H$ generated by the subset $U := \{ (x,f(x)) \: x \in V \}$. 
We endow $U$ with the topology for which 
$x \mapsto (x,f(x)), V \to U$ is a homeomorphism. 
Note that 
$f(\be)^2 = f(\be^2) = f(\be)$
implies $f(\be) = \be$, which further leads to 
$\be = f(xx^{-1}) = f(x)f(x^{-1})$, so that $f(x^{-1}) = f(x)^{-1}$. Hence 
$U = U^{-1}$. 

Since $S$ is generated by $U$ and (T1/2) in Lemma~\ref{lem:toplocglob} 
follow from from the corresponding 
properties of $V$, we obtain a group topology on $S$. 
The connectedness of $S = \la U \ra$ follows from the connectedness 
of $U$ (Exercise~\ref{exer:3.1.3}(7)). The projection 
$p_G \: G \times H \to G$ induces a covering homomorphism 
$q \: S \to G$ (Exercise~\ref{exer:3.3.1}(3)), 
so that the connectedness of $S$ and the simple
connectedness of $G$ imply that $q$ is a homeomorphism. Now 
$F := p_H \circ q^{-1} \: G \to H$ provides the required extension of
$f$. In fact, for $x \in U$ we have $q^{-1}(x) = (x,f(x))$, and
therefore $F(x) = f(x)$. 
\end{prf}

\begin{small}
\subsection*{Exercises for 
Appendix~\ref{ch:topgrp}}

\begin{exer} \mlabel{exer:3.3.8}
Let $G$ be a topological group and $\Gamma \subeq G$ be a subgroup. 
Then $\Gamma$ is discrete (with respect to the subspace topology) 
if and only if there exists a $\be$-neighborhood $U$ in $G$ with 
$U \cap \Gamma = \{\be\}$. 
\end{exer}

\begin{exer} \mlabel{exer:3.1.3} 
Let $G$ be a topological group. Show that:
\begin{description}
\item[\rm(1)] If $U \subeq G$ is an open set and $V \subeq G$ arbitrary, 
then 
$$V \cdot U := \{ vu \: v \in V, u \in U\}$$ is open in $G$. 
\item[\rm(2)] If $U, V \subeq G$ are compact, then $U\cdot V$ is compact. 
\item[\rm(3)] If $H \subeq G$ is a subgroup containing an identity 
neighborhhood, then it is open and closed. In particular, 
every open subgroup is closed. Hint: Consider the cosets of $H$ in $G$. 
\item[\rm(4)] If $\be \in U = U^{-1}\subeq G$ is an open subset, then 
each set $U^n := U^{n-1} \cdot U$ defines an ascending sequence of open 
subsets of $G$ and $\bigcup_{n \in \N} U^n$ is an open closed subgroup. 
\item[\rm(5)] If $\be \in U = U^{-1}\subeq G$ is a compact identity 
neighborhood, then each set $U^n$ is compact and $\bigcup_n U^n$ contains 
the identity component of $G$. 
\item[\rm(6)] Each locally compact topological group $G$ 
is paracompact, i.e., all its connected components are $\sigma$-compact. 
\item[\rm(7)] If $G$ is generated by a (arcwise) connected subset, then 
$G$ is (arcwise) connected. 
\end{description}
\end{exer}

\begin{exer} \mlabel{exer:3.1.5} 
Show that in a topological group $G$, the closure of a subset 
$S$ is given by 
$$ \oline S = \bigcap_{U \in \cU(\be)} SU. $$
Hint: $g \in \oline S$ is equivalent to $S \cap gU \not=\eset$ 
for every $\be$-neighborhood in $G$.   
\end{exer}

\begin{exer} \mlabel{exer:3.1.6} 
Let $G$ be a group, endowed with a topology for which multiplication and 
inversion are continuous maps. Show that $G$ is Hausdorff if 
$\{\be\}$ is a closed subset. 
\end{exer}

\begin{exer} \mlabel{exer:3.3.2} (Refining Lemma~\ref{lem:toplocglob}) 
Show that the conclusion of \break Lemma~\ref{lem:toplocglob} is still valid if 
the assumption (T1) is weakened as follows: 
There exists an open subset $D \subeq U \times U$ with 
$xy\in U$ for all $(x,y) \in D$, 
containing all pairs $(x,x^{-1})$, $(x,\be)$, $(\be,x)$ for $x \in U$, such that 
the group multiplication 
$m \: D \to U$ is continuous. 
\end{exer}

\begin{exer} \mlabel{exer:4.1.1} Let 
$X$ be a topological Hausdorff space and $G$ a topological 
group. Show that the compact open topology turns $C(X,G)$ into a topological group. 
Hint: Show first that the set ${\cal F}$ of all sets $W(K,U)$ with 
$\be \in U$, $U$ open in $G$, $K$ compact in $X$, is a filter basis of subsets 
of $C(X,G)$, satisfying all conditions in Lemma~\ref{lem:loctop}. 
To verify (U3), find for $g \in C(X,G)$ and $W(K,U)$ an open 
$\be$-neighborhood $U' \subeq G$ with 
$g(x)U'g(x)^{-1} \subeq U$ for all $x \in K$. Now Lemma~\ref{lem:loctop} 
yields a group topology for which ${\cal F}$ is a local basis. 
It is easy to see that all sets $W(K,U)$, $K$ compact in $X$ and $U$ open in $G$ 
are open in this group because they are contained in $g{\cal F}$ for some constant map 
$g$. Conversely, one has to show that any set of the form $g\cdot W(K,U)$, 
$\be \in U$, is open in the compact open topology. 
Let $V \subeq U$ be an open symmetric $\be$-neighborhood with $VV \subeq U$. 
Then show that there exist finitely many points $x_1,\ldots, x_n \in K$ 
and a covering $K_1,\ldots, K_n$ of $K$ by compact subsets of $K$ with 
$g(K_i) \subeq g(x_i)V$ for $i =1,\ldots,n$. Now 
$f \in \bigcap_{i = 1}^n W(K_i, g(x_i)V)$ 
implies for $x \in K_i$ the relation 
$g(x)^{-1}f(x) \in g(x)^{-1}g(x_i)V \subeq VV \subeq U$, so that 
$$ g \in \bigcap_{i = 1}^n W(K_i, g(x_i)V) \subeq gW(K,U). $$
\end{exer}

\begin{exer} \mlabel{exer:6.1.1} Let $G$ be a topological group and $H \leq G$ a subgroup. 
We endow the set $G/H = \{ gH \: g \in G \}$ of left cosets of $H$ with the 
quotient topology. Show that $G/H$ is Hausdorff if and only if $H$ is closed 
in $G$. Hint: If $H$ is closed and $g \not\in H$, then there exists an open symmetric 
$\be$-neighborhood $U \subeq G$ with $UU \cap gH = \eset$. 
\end{exer}

\begin{exer} \mlabel{exer:6.1.3} Show that a subgroup $H$ of the topological group $G$ is closed if 
and only if there exists a $\be$-neighborhood $U$ in $G$ for which 
$U \cap H$ is closed in $U$. Hint: For $x \in \oline H$, observe that $x \in HU^{-1}$ and pick 
$u \in U$ with $y := xu \in H$. Then 
$u = x^{-1} y \in U$ lies in the closure of $H \cap U$
\end{exer}

\begin{exer} \mlabel{exer:3.3.1} Let $\phi \: G \to H$ be a surjective morphism of topological 
groups. Show that the following conditions are equivalent: 
\begin{description}
\item[\rm(1)] $\phi$ is open with discrete kernel. 
\item[\rm(2)] $\phi$ is a covering map, i.e., 
 each $h \in H$ has an open neighborhood $U$ such that 
$\phi^{-1}(U) = \bigcup_{i \in I} U_i$ is a disjoint union of open subsets 
$U_i$ for which all restrictions $\phi\res_{U_i} \: U_i \to U$ are 
homeomorphisms. 
\item[\rm(3)] There exists an open $\be$-neighborhood $U \subeq G$ 
such that $\phi\res_U \: U \to \phi(U)$ is a homeomorphism onto an open 
subset of~$H$. 
\end{description}
Hint: For the implication (1) $\Rarrow$ (2), pick an open symmetric 
$\be$-neighborhood $U_G \subeq G$ with $U_G^3 \cap \ker(\phi) = \{\be\}$ 
and show that $U_G$ is mapped homeomorphically onto its image. 
\end{exer}

\begin{exer} \mlabel{exer:3.3.4}
Let $G$ be a connected topological group and 
$\Gamma \trile G$ a discrete normal subgroup. Then $\Gamma$ is central. 
\end{exer}

\begin{exer} \mlabel{exer:3.3.5}
Let $G$ be a simply connected topological group, \break 
$\alpha \: H_1 \to H_2$ a covering morphism and 
$f \: G \to H_2$ a morphism of topological groups. 
Show that there exists a morphism of topological groups 
$\tilde f \: G \to H_1$ with $\alpha \circ \tilde f = f$. 
Hint: Consider the group 
$$ \Gamma  := \{ (g, h_1) \in G \times H_1 \: f(g) = \alpha(h_1) \} $$
and show that the projection $p_{G} \: \Gamma \to G$ is a covering 
morphism if $\Gamma$ is endowed with the subspace topology of 
$G \times H_1$. 
\end{exer}
  
\end{small}

\section*{Notes and Comments on Appendix~\ref{ch:topgrp}}

Proposition~\ref{prop:lochomo-ext} 
on extensions of local homomorphisms can also be found in 
\cite[Cor.\ A.2.26]{HoM98}. 


\chapter{Smooth maps into non-Lie groups} \mlabel{app:nonlie}

\vspace{5mm}
In this appendix, we develop a framework for dealing with 
smooth maps into groups such as 
diffeomorphism groups of in\-fi\-nite-di\-men\-sio\-nal manifolds which 
are not Lie groups, but for which 
the concept of a smooth map can be made sense of as follows: 
Let $M$ be a smooth manifold modeled on a locally convex space and  let 
$\Diff(M)$ be the group of diffeomorphisms of $M$. Further, let $N$ be  a
smooth manifold. 
Although $\Diff(M)$ has no
natural Lie group structure, we call a map 
$\phi \: N \to \Diff(M)$ {\it smooth} if the map 
\[  \tilde \phi \: N \times M \to M \times M, \quad 
(n,m) \mapsto (\phi(n)(m), \phi(n)^{-1}(m)) \]
is smooth. Even in this general context, we have 
left and right logarithmic derivatives which are 
${\cal V}(M)$-valued $1$-forms on $N$. 
If $N$ is an interval in $\R$, these are simply time-dependent 
vector fields on~$M$. 

A key tool is the Uniqueness Lemma which states that on a connected manifold 
$N$ a smooth map $\phi \: N\to \Diff(M)$ is uniquely determined by its left
or right logarithmic derivative and the value in any point of~$N$. 
This implies in particular that solutions to certain initial value problems
are unique whenever they exist. In this generality, this is quite 
remarkable because there are ordinary linear differential 
equations with constant coefficients on Fr\'echet spaces $E$ for which
solutions are not unique. Nevertheless, the Uniqueness Lemma implies 
that solutions of the corresponding operator-valued 
initial value problems on the group $\GL(E)$ are unique whenever they exist. 
Here we simply consider the group $\GL(E)$ as a subgroup of
$\Diff(E)$. 

The Uniqueness Lemma implies also that a smooth action of 
a connected Lie group $G$ on $M$ is uniquely determined by 
the corresponding homomorphism of Lie algebras 
$\Lie(G) \to {\cal V}(M)$, and that any smooth representation 
$\pi \: G \to \GL(E)$ is uniquely determined by its derived representation 
$d\pi \: \Lie(G) \to \gl(E)$.  

Another consequence of the Uniqueness Lemma is that we may define 
complete vector fields $X$ on $M$ as those which are the 
derivative $X = \gamma_X'(0)$ 
of a smooth one-parameter subgroup $\gamma_X \: \R \to \Diff(M)$. 
The complete vector fields constitute the domain of the exponential
function $\exp(X) := \gamma_X(1)$ 
of $\Diff(M)$. Likewise, the domain of the exponential function of $\GL(E)$ is
the set of all continuous 
linear operators $D$ on $E$ for which the corresponding
linear vector field $v \mapsto Dv$ is complete, i.e., there exists a 
smooth representation $\alpha \: \R \to \GL(E)$ with 
$\alpha'(0) = D$. 
We may also define, for each
Lie group $G$, the domain of the exponential function of $G$ as those
elements of the Lie algebra for which the corresponding left invariant
vector field is complete. 

We conclude this appendix with a discussion of automatic 
smoothness result of the following type. If $\phi \: N \to \Diff(M)$ 
is a map for which the map $\tilde\phi_1 \: N \times M \to M, (n,m) 
\mapsto \phi(n)(m)$ is smooth and if $M$ is a Banach manifold, 
then the second component 
$\tilde\phi_2(n,m) := \phi(n)^{-1}(m)$ is automatically smooth. 
In general, we have to require at least that $\tilde\phi_2$ is $C^1$. 
If $M = E$ is a vector space and $\phi(N) \subeq \GL(E)$, then 
the continuity of $\tilde \phi_2$ implies its smoothness.  

\section{Smooth maps into diffeomorphism groups} \mlabel{sec:e.1}

In this section, we study ``smooth'' maps with values in the group 
$\Diff(M)$ of diffeomorphisms of a, possibly infinite-dimensional, smooth 
manifold $M$. One of the main points is that the group $\Diff(M)$ behaves 
in many respects like a Lie group, although it is far from being one if 
$M$ is not compact. We start with a discussion of the group structure 
on the ``tangent bundle'' of $\Diff(M)$, which permits us to define 
left and right logarithmic derivatives of a smooth map $N \to \Diff(M)$ in such a 
way that logarithmic derivatives are 
smooth $1$-forms on $N$ with values in the locally convex space 
${\cal V}(M)$. Here ${\cal V}(M)$ is endowed with a natural 
locally convex topology inherited from the $C^\infty$-topology on 
$C^\infty(M,T(M))$. We also discuss the adjoint representation of 
$\Diff(M)$ on ${\cal V}(M)$ whose ``derived representation'' is given 
by the negative Lie bracket. We also prove the Uniqueness Lemma 
and explain how it can be applied to the ``exponential function'' of $\Diff(M)$.

\subsection*{The tangent bundle of the diffeomorphism group} 

\begin{defn} \mlabel{def:e.1.1} (The tangent bundle of $\Diff(M)$) 

(a) Let $M$ be a manifold modeled on a locally convex space. 
We write $\Diff(M)$ for the group
of all diffeomorphisms of
$M$ and $\cV(M)$ for the Lie algebra of smooth vector fields on $M$,
i.e., the set of all smooth maps $X \: M \to TM$ satisfying 
$\pi_{TM} \circ X = \id_M$, where $\pi_{TM} \: TM \to M$ is the bundle
projection. 

Although we we do not define a differentiable structure on 
the group $\Diff(M)$ (if $M$ is not compact), we nevertheless think of the set 
$$ T(\Diff(M)) := \{ X \in C^\infty(M,TM) \: \pi_{TM} \circ X \in
\Diff(M)\} $$
as the {\it tangent bundle of $\Diff(M)$}, where the map 
$$ \pi \: T(\Diff(M)) \to \Diff(M), \quad X \mapsto \pi_{TM}
\circ X $$
is the bundle projection, and 
$T_\phi(\Diff(M)) := \pi^{-1}(\phi)$
is considered as the {\it tangent space} in $\phi \in \Diff(M)$. 
The vector space structure comes from the fact that, 
for $X,Y \in T_\phi(\Diff(M))$, the values $X(m), Y(m) \in T_{\phi(m)}(M)$ 
can be added. 

(b) We have natural left and right
actions of $\Diff(M)$ on $T(\Diff(M))$ by 
$$ \phi.X = T(\phi) \circ X \quad \mbox{ and } \quad X.\phi := X \circ \phi. $$
The action 
$$ \Ad \: \Diff(M) \times \cV(M) \to {\cal V}(M), \quad 
(\phi, X) \mapsto \Ad(\phi)X := T(\phi) \circ X \circ \phi^{-1} $$
is called the {\it adjoint action of $\Diff(M)$ on $\cV(M)$}. 

(c) We obtain a natural locally convex topology on 
the space $\cV(M)$ by the embedding 
$$ \cV(M) \into \prod_{k \in \N_0} C(T^k(M), T^{k+1}(M))_{c.o.}, \quad  
X \mapsto (T^k(X))_{k \in \N_0},  $$
where each factor of the product on the right carries the compact open 
topology. 
\end{defn} 

\begin{lem} \mlabel{lem:tang-mult} 
  \begin{description}
\item[\rm(a)] The 
set $T(\Diff(M))$ carries a group structure 
given by 
$$ X * Y := T(\phi) \circ Y + X \circ \psi 
   \quad \mbox{ for } \quad X \in T_\phi(\Diff(M)), Y \in T_\psi(\Diff(M)). $$
Inversion in this group is given by 
$$ X^{-1} = - T(\phi)^{-1} \circ X \circ \phi^{-1} \quad   \mbox{ for } \quad 
X \in T_\phi(\Diff(M)). $$
The maps $\pi \: T(\Diff(M)) \to \Diff(M)$ and 
$\iota \: \Diff(M) \to T(\Diff(M))$ (the zero section) are group homomorphisms 
with $\pi \circ \iota = \id_{\Diff(M)}$. We thus obtain an isomorphism of groups 
\[  ({\cal V}(M),+) \rtimes_{\Ad} \Diff(M) \to T(\Diff(M)), \quad 
(X,\phi) \mapsto X * \iota(\phi) = X \circ \phi. \] 

\item[\rm(b)] The map 
$T^\sharp \: T(\Diff(M)) \to \Diff(T(M))$
defined by 
$$ T^\sharp(X)(v) := T(\phi)v + X(p) \quad \mbox{ for } \quad 
X \in T_\phi(\Diff(M)) $$
defines an action of the group $T(\Diff(M))$ on $T(M)$. 
  \end{description}
\end{lem}

\begin{prf} (a) First we observe that we indeed have 
$X * Y \in T_{\phi \circ \psi}(\Diff(M))$. 

To verify the associativity, we note that, for 
$X \in T_\phi(\Diff(M))$, 
$Y \in T_\psi(\Diff(M))$ and $Z \in T_\zeta(\Diff(M))$, we have 
\begin{eqnarray*}
(X * Y) * Z 
&=& (T(\phi) \circ Y + X \circ \psi) * Z \\
&=& T(\phi\circ \psi)\circ Z + (T(\phi) \circ Y + X \circ \psi) \circ \zeta \\ 
&=& T(\phi) \circ T(\psi) \circ Z + T(\phi) \circ Y \circ \zeta + X \circ \psi \circ \zeta \\ 
&=& T(\phi) \circ (T(\psi) \circ Z + Y \circ \zeta) + X \circ (\psi \circ \zeta)
= X * (Y * Z). 
\end{eqnarray*}

Obviously, $\be := \iota(\id_M)$ is an identity in $T(\Diff(M))$. 
That each element $X \in T_\phi(\Diff(M))$ is invertible with inverse 
$X^{-1} = - T(\phi)^{-1} \circ X \circ \phi^{-1}$ follows 
by direct calculation. We conclude that $T(\Diff(T(M)))$ is a group 
and that $\pi$ and $\iota$ are group homomorphism. The last assertion 
follows from 
$$\iota(\phi) * X * \iota(\phi^{-1}) = T(\phi) \circ X \circ \phi^{-1} = \Ad(\phi)X. $$

(b) Since $T(\phi) \in \Diff(T(M))$, each map $T^\sharp(X)$ is a diffeomorphism of $T(M)$. 
Moreover, $T^\sharp$ is a group homomorphism because 
we have, for  $X \in T_\phi(\Diff(M))$, $Y \in T_\psi(\Diff(M))$ and 
$v \in T_p(M)$: 
\begin{align*}
T^\sharp(X * Y)v 
&= T(\phi \circ \psi)v + (X * Y)(p) 
= T(\phi) T(\psi)v + T(\phi)Y(p) + X(\psi(p)) \\ 
&= T(\phi)\big(T(\psi)v + Y(p)\big) + X(\psi(p)) 
= T^\sharp(\phi)T^\sharp(\psi)v.\\[-12mm]
\end{align*}
\end{prf}

\begin{defn} \mlabel{def:e.1.1a} 
(a) Let $N$ be a smooth manifold. A map $\phi \: N 
\to \Diff(M)$ is called {\it smooth} if the
corresponding map 
\[ \tilde \phi \: N \times M \to M \times M, \quad (n,m) \mapsto
(\phi(n)(m), \phi(n)^{-1}(m)) \]
is smooth. 
We also write 
\[ \tilde\phi_1(n,m) := \phi(n)(m) \quad 
\mbox{  and } \quad \tilde\phi_2(n,m) := \phi(n)^{-1}(m)\]
the two components of $\tilde\phi$. 

(b) If $\tilde\phi_1$ is smooth, then we have a smooth tangent map 
$$ T(\tilde\phi_1) \: T(N \times M) \cong T(N) \times T(M) \to T(M) , $$
and, for each $v \in T_p(N)$, the partial map 
$$ T_p(\phi)v \: M \to T(M), \quad 
m \mapsto T_{(p,m)}(\tilde\phi_1)(v,0) $$
is an element of $T_{\phi(p)}(\Diff(M))$. We thus obtain a {\it tangent map} 
\[ T(\phi) \: T(N) \to T(\Diff(M)), \quad T_p(N) \ni v \to T_p(\phi)v. \] 

We now define the {\it left logarithmic derivative of $\phi$ in $p$} by 
$$ \delta^l(\phi)_p \: T_p(N) \to \cV(M), \quad 
v \mapsto \phi(p)^{-1}.T_p(\phi)(v) = T(\phi(p))^{-1}\circ T_p(\phi)(v). $$
Postponing smoothness considerations to Proposition~\ref{prop:smooth-logder} 
below, we interprete $\delta^l(\phi)$ as a 
$\cV(M)$-valued $1$-form on $N$.
Similar arguments apply to the {\it right logarithmic derivative}, defined by 
$$ \delta^r(\phi)_p \: T_p(N) \to \cV(M), \quad 
v \mapsto T_p(\phi)(v).\phi(p)^{-1} = T_p(\phi)(v) \circ \phi(p)^{-1}. $$

If $N = J$ is an interval in $\R$, then we also write 
\[  \phi'(t) := T_t(\phi)(1) \in T_{\phi(t)}(\Diff(M))\] 
and note that the logarithmic derivative 
$\delta^l(\phi)_t := \phi(t)^{-1}.\phi'(t)$
describes a time-dependent vector field on $M$. 
\end{defn}

\begin{lem} \mlabel{lem:tang-mult-c}
If $\phi \: N \to \Diff(M)$ is smooth, then 
$$ T^\sharp \phi \: T(N) \to \Diff(T(M)),\quad 
v \mapsto T^\sharp(T(\phi)v) $$ 
is also smooth. Restricting to the zero section, we 
get the smooth map 
\[  N \to \Diff(T(M)), \quad n \mapsto T(\phi(n)). \] 
\end{lem}

\begin{prf} First we observe that, for $v \in T_p(N)$ and $w \in T_m(M)$, we have 
\begin{eqnarray*}
(T^\sharp\phi)(v)w 
&=& (T_p(\phi)v)^\sharp w 
= T(\phi(p))w + (T_p(\phi)v)(m)\\
&=& T_{(p,m)}(\tilde\phi_1)(0,w) + T_{(p,m)}(\tilde\phi_1)(v,0) 
= T_{(p,m)}(\tilde\phi_1)(v,w),   
\end{eqnarray*}
so that 
$$ (T^\sharp\phi)\,\tilde{}_1 = T(\tilde\phi_1) \: T(N) \times T(M) \to T(M). $$
In view of Lemma~\ref{lem:tang-mult}(a), we also have 
\begin{align*}
&(T^\sharp\phi)(v)^{-1}w 
= ((T_p(\phi)v)^{-1})^\sharp w
= \big( - T(\phi(p))^{-1} \circ (T_p(\phi)v) \circ \phi(p)^{-1}\big)^\sharp w \\
&= T(\phi(p))^{-1} w - T(\phi(p))^{-1} \circ (T_p(\phi)v)\big(\phi(p)^{-1}(m))  
= T_{(p,m)}(\tilde\phi_2)(v,w),   
\end{align*} 
which leads to 
$(T^\sharp\phi)\,\tilde{}_2 = T(\tilde\phi_2) \: T(N) \times T(M) \to T(M).$
This shows that $T^\sharp \phi$ is smooth with 
$\tilde{T^\sharp\phi} = T(\tilde\phi).$
\end{prf}

\begin{rem} \mlabel{rem:smooth-act} If $\sigma \: G \times M \to M$ is an action 
of the Lie group $G$ 
on the smooth manifold $M$, then $\sigma$ is smooth if and only if the corresponding 
map 
$$ \tilde \sigma \: G \to \Diff(M), \quad g \mapsto \sigma_g, \quad 
\sigma_g(m) := \sigma(g,m) $$
is smooth. In fact, if $\sigma$ is smooth, 
then the map $(g,m) \mapsto \sigma_g^{-1}(m) = \sigma(g^{-1},m)$ is automatically 
smooth. The converse is clear. 
\end{rem}

\begin{prop} \mlabel{prop:vec-topliealg}
{\rm(a)} If $M$ is finite-dimensional, then 
the Lie bracket on ${\cal V}(M)$ is continuous 
with respect to the natural topology on ${\cal V}(M)$. 

{\rm(b)} For each $\phi \in \Diff(M)$, the map 
$\Ad(\phi)$ is a linear topological automorphism of the Lie algebra ${\cal V}(M)$. 
\end{prop}

\begin{prf} (a) Let ${\cal A} = \{ (\phi_i, U_i) \: i \in I\}$ be an $E$-atlas of $M$. 
Then the restriction maps define a topological embedding 
\begin{eqnarray}
  \label{eq:vec-emb}
{\cal V}(M) \to \prod_{i \in I} {\cal V}(U_i) \cong \prod_{i \in I} C^\infty(U_i, E),  
\end{eqnarray}
so that it suffices to show that, for each $i$, the map 
$$ {\cal V}(M) \times {\cal V}(M) \to {\cal V}(U_i), \quad 
(X,Y) \mapsto [X,Y]\res_{U_i} $$
is continuous. Since the restriction map 
${\cal V}(M) \to {\cal V}(U_i)$ is continuous 
and $[X,Y]\res_{U_i} = [X\res_{U_i}, Y\res_{U_i}]$, 
it suffices to show that the Lie bracket on 
each ${\cal V}(U_i)$ is continuous.  

Hence we may assume that $M$ is an open subset 
of some finite-dimensional vector space~$E$. Then 
we write vector fields as $X(m) = (m,\tilde X(m))$ with 
$\tilde X \in C^\infty(M,E)$. We then have 
$$ [X,Y]\,\tilde{}\, (m) = d\tilde Y(m,\tilde X(m)) - d\tilde X(m, \tilde Y(m)), $$
and it is an easy consequence of the Chain Rule that this bilinear 
operation on ${\cal V}(M) \cong C^\infty(M,E)$ is continuous. 

(b) From $\Ad(\phi)X = T(\phi) \circ X \circ \phi^{-1}$ and the Chain Rule,  
we derive for each $k \in \N_0$: 
$$ T^k(\Ad(\phi)X) = T^{k+1}(\phi) \circ T^k(X) \circ T^k(\phi)^{-1}. $$
Since the map 
\[  C(T^k(M), T^{k+1}(M))\to C(T^k(M), T^{k+1}(M)), \quad 
f \mapsto T^{k+1}(\phi) \circ f \circ T^k(\phi)^{-1} \]
is continuous with respect to the compact open topology 
(Lemma~\ref{covsuppo}, Lemma~\ref{pubas}), 
the continuity of $\Ad(\phi)$ follows from the continuity of the 
map ${\cal V}(M) \to C(T^k(M), T^{k+1}(M)), X \mapsto T^k(X)$. 

Since each map $\Ad(\phi)$ is continuous, the relation $\Ad(\phi)^{-1} = \Ad(\phi^{-1})$ 
implies that each $\Ad(\phi)$ is an isomorphism of topological 
vector spaces. For $X, Y \in {\cal V}(M)$, the vector fields 
$X$ and $\Ad(\phi)(X)$ are $\phi$-related 
$$ T(\phi) \circ X = (\Ad(\phi)X) \circ \phi. $$
Likewise $Y$ and $\Ad(\phi)Y$ are $\phi$-related, and therefore the vector fields 
$[X,Y]$ and $[\Ad(\phi)X,\Ad(\phi)Y]$ are $\phi$-related  (Lemma~\ref{larelglob}), which 
means that $\Ad(\phi)[X,Y] = [\Ad(\phi)X,\Ad(\phi)Y]$, 
i.e., $\Ad(\phi)$ is Lie algebra automorphism.
\end{prf}

\begin{rem} If $M$ is infinite-dimensional, 
we cannot expect the bracket on ${\cal V}(M)$ 
to be continuous. To see this, we assume that 
$M= U$ is an open subset of a locally convex space $E$ 
and consider the subalgebra 
$\aff(E) \cong E \rtimes \gl(E)$ 
of affine vector fields 
$X_{A,b}$ with $X_{A,b}(v) = Av + b$, $A \in {\cal L}(E)$, $b \in E$. It is easy to see that 
the natural topology on ${\cal V}(U)$ induces on $\aff(E)$ 
the product topology of the original topology on $E$ 
and the compact open topology on $\gl(E) := ({\cal L}(E),[\cdot,\cdot])$. 
In view of 
$$ [X_{A,b}, X_{A',b'}] = X_{[A',A], A'b - Ab'}, $$
it therefore suffices to show that the bilinear evaluation map 
${\cal L}(E)_{c.o.} \times E \to E$ 
is not continuous if $\dim E = \infty$. 
Pick $0 \not = v \in E$ and embed $E'_c \into {\cal L}(E)_{c.o.}$ by assigning to 
$\alpha \in E'$ the operator $v \otimes \alpha \: x \mapsto \alpha(x) v$. 
Hence it suffices to see that the evaluation map 
$$ \ev \: E'_c \times E \to \R, \quad (\alpha,v) \mapsto \alpha(v) $$ 
is not continuous. Basic neighborhoods of $(0,0)$ in $E'_c \times E$ are 
of the form $\hat K\times U_E$, where 
$U_E \subeq U$ is a $0$-neighborhood, $K \subeq E$ is compact, 
and 
\[ \hat K := \{ f \in E' \: (\forall k \in K)\, |f(k)| \leq 1\} \] 
 is the polar set of~$K$. 
On $\hat K \times U_E$ the evaluation map is bounded if and only if 
$U_E$ is contained in some multiple of the bipolar $\hats K$, 
which, according to the Bipolar Theorem~\ref{bipolar}, coincides with the balanced convex closure 
of $K$, which is pre-compact (\cite[Prop.~7.11]{Tr67}). 
Then $\hats K$ is a pre-compact $0$-neighborhood in $E$, so that 
$E$ is finite-dimensional (cf.\ \cite[proof.~of Thm.~1.22]{Ru91}). 
A similar argument shows that, if we endow $E'$ with the finer 
topology of uniform convergence on bounded subsets of $E$, 
then the evaluation map is continuous if and only if $E$ is normable, which is 
equivalent to the existence of a (weakly) bounded $0$-neighborhood 
(Proposition~\ref{normable}). 
\end{rem} 

\begin{prop} \mlabel{prop:smooth-logder} 
Let $M$ and $N$ be smooth manifolds. 
\begin{description}
\item[\rm(1)] A map $f \: N \to \cV(M)$ is smooth if the corresponding map 
$$ \hat f \: N \times M \to T(M), \quad (n,m) \mapsto f(n)(m) $$
is smooth. In particular, each fiberwise linear 
map $\alpha \: T(N) \to {\cal V}(M)$ for which 
$\hat\alpha \: T(N) \times M \to T(M), (v,m) \mapsto \alpha(v)(m)$ is smooth 
defines a smooth ${\cal V}(M)$-valued $1$-form on $N$. 
\item[\rm(2)] If $\alpha \in \Omega^1(N,{\cal V}(M))$ is such that 
$\hat\alpha$ is smooth and $f \: N \to \Diff(M)$ is smooth, then 
$$ \Ad(f).\alpha \: T(N) \to {\cal V}(M), \quad \big(\Ad(f).\alpha\big)_p 
:= \Ad(f(p)) \circ \alpha_p $$
defines a smooth ${\cal V}(M)$-valued $1$-form on~$N$.
\item[\rm(3)] If $\phi \: N \to \Diff(M)$ is a smooth map, then 
the left and right logarithmic derivatives of $\phi$ 
are smooth ${\cal V}(M)$-valued $1$-forms, for 
which the functions $\delta^l(\phi)\,\hat{}$ and $\delta^r(\phi)\,\hat{}$ are smooth. 
\end{description}
\end{prop} 

\begin{prf} (1) In view of (\ref{eq:vec-emb}), 
it suffices to show that, for each chart $(\phi_i,U_i)$ of~$M$, the map 
$$ n \mapsto f(n)\res_{U_i} \: N \to {\cal V}(U_i)$$
is smooth. We may therefore assume that $M$ is an open subset 
of some locally convex space $E$. Then 
${\cal V}(M) \cong C^\infty(M,E)$, and 
\[ \hat f(n,m) = (f(n), f^\sharp(n,m)), \]
where $f^\sharp \: N \times M \to E$ is smooth by assumption. 
Now the smoothness of the corresponding map 
$N \to C^\infty(M,E) \cong {\cal V}(M)$
follows from Proposition~\ref{prop:cartes-closed}. 

(2) In view of (1), this follows from 
\begin{eqnarray*}
&&\ \ \ \ \Ad(f(p)).\alpha_p(v)(m) 
= T(f(p)).\alpha_p(v)\big(f(p)^{-1}(m)\big)  \\
&=& T_{(p,f_2(p,m))}(\hat f)\Big(0, \alpha_p(v)(f(p)^{-1}(m))\Big) 
= T_{(p,f_2(p,m))}(\hat f)\Big(0, \tilde\alpha(v, \tilde f_2(p,m))\Big) 
\end{eqnarray*}
for $v \in T_p(N)$ because the maps 
$T(\hat f)$, $\tilde\alpha$ and $\tilde f_2$ are smooth. 

(3) For each $m \in M$ and $v \in T_p(N)$ we have 
$$ (\delta^l(\phi)v)(m) = T(\tilde\phi_2)(0_p, T(\tilde\phi_1)(v,0_m)), $$
so that the smoothness of $\delta^l(\phi)$ follows similarly as in (2). 
\end{prf}

\subsection*{The Uniqueness Lemma and its applications} 
We now turn to the rules for logarithmic derivatives and the Uniqueness Lemma; 
both are key tools throughout this book. 
The following lemma generalizes the corresponding result for Lie algebra-valued 
forms from Lemma~\ref{lem:c.12}, used for Lie group-valued maps, to $\Diff(M)$-valued 
maps. 

\begin{lem} [Product and Quotient Rule] \mlabel{lem:e.1.2}
For two smooth maps $f, g \:  N \to \Diff(M)$,  
the following assertions hold: 
\begin{description}
\item[\rm(1)] {\rm(Product Rule)} The map $(fg)(p) := f(p) \circ g(p)$ is smooth with 
$$ \delta^l(fg) 
= \delta^l(g) + \Ad(g^{-1}).\delta^l(f) 
\quad \hbox{ and } \quad 
\delta^r(fg) = \delta^r(f) + \Ad(f).\delta^r(g). $$  
\item[\rm(2)] {\rm(Quotient Rule)} The map $fg^{-1} \: p \mapsto f(p) \circ g(p)^{-1}$ 
is smooth with 
$$ \delta^l(fg^{-1}) = \Ad(g).(\delta^l(f) - \delta^l(g)). $$
In particular, we have 
$$ \delta^l(f^{-1}) = - \Ad(f).\delta^l(f) = - \delta^r(f).$$
\end{description}
\end{lem}

\begin{prf} The smoothness of $f^{-1}$ follows directly from the
definitions. For the smoothness of the product map $fg$, we
observe that the maps 
\[  N^2 \times M \to M, \quad 
(t,s,x) \mapsto f(t)\big(g(s)(x)\big) = \big(f(t) \circ g(s)\big)(x)\] 
and 
\[   N^2 \times M \to M, \quad (t,s,x) \mapsto \big(f(t) \circ g(s)\big)^{-1}(x) = 
g(s)^{-1}\big(f(t)^{-1}(x)\big) \]
are smooth, as compositions of smooth maps. 

The Chain Rule now implies 
\[  T_p(fg)v = T(f(p))(T_p(g)v) + (T_p(f)v) \circ g(p), \] 
which leads to 
\begin{eqnarray*}
\delta^l(fg)_p 
&=& 
g(p)^{-1}.T_p(g) + g(p)^{-1}f(p)^{-1}.\big(T_p(f) \circ g(p)\big) \\
&=& \delta^l(g)_p + \Ad(g(p)^{-1}).\delta^l(f)_p,
\end{eqnarray*}
and this further leads to 
$$ \delta^r(fg) = 
\Ad(fg).\delta^l(fg) 
= \Ad(f).\delta^r(g) + \delta^r(f). $$  
For $g = f^{-1}$, we obtain in particular 
$$ 0 = \delta^l(ff^{-1}) 
= \delta^l(f^{-1}) + \Ad(f).\delta^l(f).  $$
Combining this with (1), we obtain the Quotient Rule. 
\end{prf}

\begin{lem} [Uniqueness Lemma] \mlabel{lem:e.1.3}
Suppose that $N$ is connected and that $f, g \:  N \to \Diff(M)$ are smooth. 
\begin{description}
\item[\rm(1)] The relation $\delta^l(f)= \delta^l(g)$ is equivalent to the
existence of $\phi \in \Diff(M)$ with 
$g(p) = \phi \circ f(p)$ for all $p \in N$. 
\item[\rm(2)] The relation $\delta^r(f)= \delta^r(g)$ is equivalent to the
existence of $\phi \in \Diff(M)$ with 
$g(p) = f(p) \circ \phi$ for all $p \in N$. 
\end{description}
\end{lem} 

\begin{prf} (1) If $g(p) = \phi \circ f(p)$ for each $p \in N$, then 
$T_p(g)  = \phi(p).T_p(f)$, and therefore 
$\delta^l(g) = \delta^l(f)$. 

If, conversely, $\delta^l(g) = \delta^l(f)$, 
then Product and Quotient Rule lead to 
\[  \delta^l(gf^{-1}) 
= \delta^l(f^{-1}) + \Ad(f).\delta^l(g) 
= \delta^l(f^{-1}) + \Ad(f).\delta^l(f) 
= \delta^l(ff^{-1}) = 0. \]
Let $\gamma := gf^{-1}$. Then $\delta^l(\gamma) = 0$ implies 
$T_p(\gamma) = 0$ for each $p \in N$, so that, for each $x \in M$, the
map $p \mapsto \gamma(p)(x)$ has vanishing derivative, hence is
locally constant because $M$ is modeled on a locally convex space 
(Lemma~\ref{locconst}). 
Since $N$ is connected, $\gamma$ is constant. 
We conclude that $g f^{-1}$ is constant, which
means that $g = \phi \circ f$ for some $\phi \in \Diff(M)$. 

(2) follows with similar arguments as (1). 
\end{prf}

\begin{prop} [Uniqueness for time-dependent vector fields] \mlabel{prop:e.1.4}
Let $J \subeq \R$ be an interval, $t_0 \in J$, $\phi_0 \in \Diff(M)$,
 and $X \: J\to \cV(M)$ be a smooth curve. 
Then there exists at most one smooth curve 
$\gamma \: J \to \Diff(M)$ solving the initial value problem 
\begin{eqnarray}
  \label{eq:e.2a}
\gamma'(t) = \gamma(t).X(t) = T(\gamma(t)) \circ X(t), \quad \gamma(t_0) =
\phi_0 
\end{eqnarray}
or 
\begin{eqnarray}
  \label{eq:e.2b}
\gamma'(t) = X(t).\gamma(t) = X(t) \circ \gamma(t), \quad \gamma(t_0) =
\phi_0. 
\end{eqnarray}
\end{prop}

\begin{prf} The curve $\gamma$ solves (\ref{eq:e.2a}) if and only if 
$\delta^l(\gamma) = X$. If $\eta$ is a second solution of (\ref{eq:e.2a}), then
$\eta(t_0) = \gamma(t_0)$, so that the Uniqueness Lemma~\ref{lem:e.1.3} 
implies $\phi = \psi$. We likewise argue for (\ref{eq:e.2b}) with 
$\delta^r(\gamma) = X$. 
\end{prf}

We can even sharpen the assertion of the preceding proposition: 
Whenever there is a smooth curve 
$\gamma \: J \to \Diff(M)$ satisfying the initial value problem 
\begin{eqnarray}
  \label{eq:2.3.2}
\gamma(0) = \id_M \quad \hbox{ and } \quad 
\gamma'(t) = X_t \circ \gamma(t) 
\end{eqnarray}
for a time-dependent vector field $X \: J \to {\cal V}(M)$, then all integral 
curves of $X$ on $M$ are of the form 
\begin{eqnarray}
  \label{eq:2.3.3b}
\eta(t) = \gamma(t)(m),
\end{eqnarray}
hence unique. In particular, the existence 
of multiple integral curves of $X$ implies that (\ref{eq:2.3.2}) has no solution. 
Below we shall see examples where this situation arises, even for linear 
differential equations. 

\begin{lem} \label{unique-intcurve} 
Let $J \subeq \R$ be an interval containing $0$,  
$\gamma \: J \to \Diff(M)$ 
be a smooth curve with $\gamma(0) = \id_M$ and $m_0 \in M$. 
Let $X_t := \delta^r(\gamma)_t$ be the corresponding time-dependent 
vector field on $M$ with $X_t \circ \gamma(t) = \gamma'(t)$ 
and assume that $\eta \: J \to M$ 
is a solution of the initial value problem: 
$$ \eta(0)= m_0 \quad \hbox{ and } \quad 
\eta'(t) = X_t(\eta(t)) \quad \hbox{ for } \quad t \in J. $$
Then $\eta(t) = \gamma(t)(m_0)$ holds for all $t \in J$. 
\end{lem}

\begin{prf} The smooth curve $\alpha \: J \to M, t \mapsto \gamma(t)^{-1}(\eta(t))$ 
satisfies $\alpha(0) = m_0$ and 
\begin{eqnarray*}
\alpha'(t) 
&=& (\gamma^{-1})'(\eta(t)) + T(\gamma(t)^{-1})(\eta'(t)) 
= T(\gamma(t)^{-1})\Big(\delta(\gamma^{-1})_t(\eta(t)) + \eta'(t)\Big) \\ 
&=& T(\gamma(t)^{-1})\big(-\delta^r(\gamma)_t(\eta(t)) + \eta'(t)\big) \\
&=& T(\gamma(t)^{-1})\big(-X_t(\eta(t)) + \eta'(t)\big) =0.
\end{eqnarray*}
Hence $\alpha$ is constant $m_0$, and the assertion follows. 
\end{prf}

\begin{prop}  \mlabel{prop:e.1.5} Let $G$ be a Lie group and $M$ a smooth manifold. 
  \begin{description}
  \item[\rm(1)] If $\sigma \: M \times G \to M$ is a smooth right action of $G$
on $M$, then 
$$ \dot\sigma\:  \g \to {\cal V}(M), \qquad 
\dot\sigma(x) :=  T_{(p,\be)}(\sigma)(0,x) $$ 
defines a homomorphism of Lie algebras. 
\item[\rm(2)] If $G$ is connected and $\alpha \: \g \to {\cal V}(M)$ a homomorphism 
of Lie algebras, then there is at most one smooth right action 
$\sigma$ of $G$ on $M$ with $\dot\sigma = \alpha$. 
  \end{description}
\end{prop}

\begin{prf} (1) We pick $p \in M$ and write 
$\sigma^p \: G \to M, g \mapsto p.g := \sigma(p,g)$ for the smooth orbit map of $p$. 
Then $\sigma^p \circ \lambda_g = \sigma^{p.g}$ leads to 
$$ T(\sigma^p)(x_l(g)) = T_\be(\sigma^{p.g})x = \dot\sigma(x)_{p.g} 
= (\dot\sigma(x) \circ \sigma^p)(g), $$
which means that the vector fields $x_l$ and $\dot\sigma(x)$ are 
$\sigma^p$-related. We conclude that, for $x,y \in \g$, the vector fields 
$[x_l, y_l] = [x,y]_l$ and $[\dot\sigma(x), \dot\sigma(y)]$ are also 
$\sigma^p$-related 
(Lemma~\ref{larelglob}), 
which leads to 
$$ [\dot\sigma(x), \dot\sigma(y)](p) 
= T_\be(\sigma^p)[x,y]_l(\be) 
= T_{(p,\be)}(\sigma)(0,[x,y]) 
= \dot\sigma([x,y])(p). $$

(2) If $\sigma$ is a smooth right action with $\dot\sigma = \alpha$ and 
$\sigma_g(p) := \sigma(p,g)$, then $\sigma(gh) = \sigma(h) \sigma(g)$ 
implies that the map 
$\sigma^\vee \: G \to \Diff(M)$, $\sigma^\vee(g) := \sigma_g$ 
satisfies 
\[ T_g(\sigma^\vee) \circ T_\be(\lambda_g) = T_\be(\sigma^\vee)(\cdot) \circ \sigma(g), \] 
so that 
\[  \delta^r(\sigma^\vee)_g =  T_g(\sigma^\vee)(\cdot) \circ \sigma(g)^{-1}  
= T_\be(\sigma^\vee) \circ T_\be(\lambda_g)^{-1} 
= \dot\sigma \circ T_\be(\lambda_g)^{-1}. \]
Hence (2) follows from the Uniqueness Lemma. 
\end{prf}

\subsection*{The derivative of the adjoint representation} 

In the following lemma, we calculate the derivative of the adjoint representation 
of $\Diff(M)$ on the Lie algebra 
${\cal V}(M)$. The ``wrong'' sign we obtain here is due to the fact that the 
Lie bracket on ${\cal V}(M)$ turns it into the ``Lie algebra'' of the opposite 
group $\Diff(M)^{\rm op}$ with the product $\phi\psi := \psi \circ \phi$. 
This group acts on $M$ from the right via $(m,\phi) \mapsto \phi(m)$ 
(cf.~Proposition~\ref{prop:e.1.5}). 

\begin{lem} \mlabel{lem:e.1.6} 
Let $J \subeq \R$ be an interval and $\gamma \: J \to \Diff(M)$, \break 
$\alpha \: J \to \cV(M)$ be smooth curves. Then the curve 
$$ \beta\: J \to {\cal V}(M), \quad t \mapsto \Ad(\gamma(t))\alpha(t) $$
is smooth and satisfies the differential equation 
$$ \beta'(t) 
= \Ad(\gamma(t))\Big([\alpha(t),\delta^l(\gamma)_t]+ \alpha'(t)\Big). $$
\end{lem}

\begin{prf} First we show that $\beta$ is smooth. 
With the same argument as in the proof of Proposition~\ref{prop:smooth-logder}, 
we see that it suffices to show that the map 
$$ \hat \beta \: J \times M \to T(M), \quad 
(t,m) \mapsto \beta(t)(m) $$
is smooth. This follows from 
\[  \beta(t)(m) = T(\gamma(t))\alpha(t)\big(\gamma(t)^{-1}(m)\big) 
= T(\hat \gamma)\big(0, \alpha(t)(\tilde\gamma_2(t,m))\big). \]

Let $t_0 \in J$. From the continuity of the linear map 
$\Ad(\gamma(t_0))$ (Proposition~\ref{prop:vec-topliealg}(b)), 
we obtain 
\[  \frac{d}{dt} \Ad(\gamma(t))\alpha(t) 
= \Ad(\gamma(t_0))\frac{d}{dt}\Big(\Ad\big(\gamma(t_0)^{-1}\gamma(t)\big)\alpha(t)\Big). \]
Now $\eta(t) := \gamma(t_0)^{-1}\gamma(t)$ satisfies 
$\eta(t_0) = \id_M$ and we put 
$X := \alpha(t_0)$, $Y := \eta'(t_0) = \delta^l(\gamma)_{t_0}$ and 
$Z := \alpha'(t_0)$. Then it remains to show that 
$\beta(t) := \Ad(\eta(t))\alpha(t)$ satisfies 
$\beta'(t_0) = [X,Y] + Z$. 

Differentiating the relation $\eta(t)(\eta(t)^{-1}(m)) = m$, we get 
$$ (\eta^{-1})'(t) = -T(\eta(t))^{-1}\circ \eta'(t) \circ \eta(t)^{-1}. $$
Therefore 
$\beta(t) = T(\eta(t)) \circ \alpha(t) \circ \eta(t)^{-1}$
leads, for each $m \in M$, to  
\begin{eqnarray}\label{eq:d.1.1}
 \beta'(t_0)(m) 
&=& (T \circ \eta)'(t_0)\big(\alpha(t_0)(\eta(t_0)^{-1}m)\big) 
+ T^2(\eta(t_0))(\alpha'(t_0)(\eta(t_0)^{-1}m)) \notag\\ 
&& + T^2(\eta(t_0))T(\alpha(t_0))(\eta^{-1})'(t_0)(m)\notag \\
&=& (T \circ \eta)'(t_0)X(m) + Z(m) - T(X)Y(m). 
\end{eqnarray}
We now calculate this expression in local coordinates, 
where we have $X(m) = (m,  X^\sharp(m))$, $Y(m) = (m, Y^\sharp(m))$ and 
\[  T(\eta(t))(m,v) = \big(\eta(t)(m), d(\eta(t))(m)v\big), \]
which leads to 
\begin{eqnarray} \label{eq:d.1.2}
(T \circ \eta)'(t_0)X(m) 
&=& \Big((m, X^\sharp(m)), \big(\eta'(t_0)(m), d(\eta'(t_0))(m)X^\sharp(m)\big)
\Big)\notag \\
&=& \Big((m,X^\sharp(m)), (Y^\sharp(m), dY^\sharp(m)X^\sharp(m))\Big) 
\end{eqnarray}
and 
\begin{equation}
  \label{eq:d.1.3}
T(X)Y(m) = \big((m, X^\sharp(m)), (Y^\sharp(m), dX^\sharp(m)Y^\sharp(m))\big).
\end{equation}
Subtracting \eqref{eq:d.1.3} from \eqref{eq:d.1.2}, we obtain with \eqref{eq:d.1.1} 
\[  \beta'(t_0)(m) - Z(m)  
= ((m,0), (0,dY^\sharp(m)X^\sharp(m) - d X^\sharp(m)Y^\sharp(m))\big) 
= [X,Y](m), \]
where we identify $T_m(M)$ with the tangent space $T_v(T_m(M)) 
\subeq T(T(M))$ in 
any point~$v \in T_m(M)$. 
\end{prf} 

\begin{prop} \mlabel{prop:e.1.8} 
Let $J \subeq \R$ be an interval, $\gamma \: J \to \Diff(M)$ be a smooth curve, $t_0
\in J$, and $\phi \in \Diff(M)$ be such that 
\begin{description}
\item[\rm(1)] $\gamma(t_0)$ commutes with $\phi$, and 
\item[\rm(2)]$\Ad(\phi)\delta^l(\gamma)_t = \delta^l(\gamma)_t$ for each $t \in J$. 
\end{description}
Then $\phi$ commutes with each $\gamma(t)$, $t \in J$. 
\end{prop}

\begin{prf} Let $\beta(t) := \phi \circ \gamma(t) \circ \phi^{-1}$. 
Then $\beta \: J \to \Diff(M)$ is a smooth curve with 
$\beta(t_0) = \gamma(t_0)$. Moreover, 
\begin{eqnarray*}
\delta^l(\beta)_t 
&=&  \big(\phi \circ \gamma(t)^{-1} \circ \phi^{-1}\big).((\phi.\gamma'(t)) \circ \phi^{-1}) \\
&=&  (\phi.\delta^l(\gamma)_t) \circ \phi^{-1} 
=  \Ad(\phi)\delta^l(\gamma)_t = \delta^l(\gamma)_t. 
\end{eqnarray*}
Now the assertion follows from the Uniqueness Lemma \ref{lem:e.1.3}.
\end{prf}

The following proposition provides a criterion for two curves in 
$\Diff(M)$ to commute. 
\begin{prop} \mlabel{prop:e.1.9} 
Let $J \subeq \R$ be an interval, $\gamma,\eta \: J \to \Diff(M)$ smooth curves 
and $X \in \cV(M)$. Then the following assertions hold: 
\begin{description}
\item[\rm(a)] If $[X, \delta^l(\gamma)_t] = 0$ or $[X, \delta^r(\gamma)_t] = 0$ 
for all $t \in J$, then $\Ad(\gamma(t))X = X$ for all $t \in J$. 
\item[\rm(b)] If $[\delta^l(\gamma)_t, \delta^l(\eta)_s] = 0$ for all $t,s \in J$, then 
$\gamma(t) \eta(s) = \eta(s) \gamma(t)$ for all $t,s \in J$, and 
$\delta^l(\gamma\eta) = \delta^l(\gamma) + \delta^l(\eta).$
\end{description}
\end{prop}

\begin{prf} (a) If $X$ commutes with $\delta^l(\gamma)$, then the assertion 
follows from the fact that the curve $\beta(t) := \Ad(\gamma(t))X$ in 
$\cV(M)$ is constant by Lemma~\ref{lem:e.1.6}.  
If $X$ commutes with $\delta^r(\gamma) = - \delta^l(\gamma^{-1})$ 
(Lemma~\ref{lem:e.1.2}), then we apply the preceding argument to 
$\beta(t) := \Ad(\gamma(t))^{-1}X$. 

(b) First we use (a) to see that $\Ad(\gamma(t))\delta^l(\eta)_s 
= \delta^l(\eta)_s$ holds for all $t,s \in J,$ and then 
Proposition~\ref{prop:e.1.8} to see that each $\gamma(t)$ commutes with 
each $\eta(s)$. Finally Lemma \ref{lem:e.1.2}(2) implies 
$ \delta^l(\gamma\eta) = \delta^l(\gamma) + \delta^l(\eta).$
\end{prf}

\subsection*{The exponential function of Diff(M)} 

\begin{defn} \mlabel{def:e.1.10} 
Let $M$ be a locally convex manifold. 
We call a vector field $X \in \cV(M)$ {\it complete} if there exists a
smooth curve $\gamma_X \: \R \to \Diff(M)$ with 
\begin{eqnarray}
  \label{eq:e.6}
\delta^l(\gamma_X)_t = X \quad \hbox{ for all} \quad t \in \R
\quad \hbox{ and } \quad 
\gamma_X(0) = \id_M. 
\end{eqnarray}
In view of Proposition~\ref{prop:e.1.4}, for each vector field $X$, the initial value
problem (\ref{eq:e.6}) has at most one solution, and the completeness of a 
vector field means that a smooth solution exists. 
We shall see in Lemma~\ref{lem:e.1.16}(3) below that completeness of a 
vector field $X$ is equivalent to an existence of a smooth $\R$-action 
on $M$ with infinitesimal generator~$X$. 

We write ${\cal D}_{\exp}(\Diff(M)) := \cV(M)_{\rm cp} \subeq \cV(M)$ 
for the set of complete vector
fields and define the {\it exponential function of $\Diff(M)$} by 
$$ \exp \: \cV(M)_{\rm cp} \to \Diff(M), \quad 
X \mapsto \gamma_X(1). $$
\end{defn}

\begin{rem}
If $M$ is a compact manifold, then 
${\cal D}_{\exp}(\Diff(M)) 
= \cV(M)$ because each smooth vector field on $M$ is complete (\cite[Cor.~IV.2.4]{La99}), 
but if $M$ is not compact, then ${\cal D}_{\exp}(\Diff(M))$ is strictly smaller than 
$\cV(M)$. The situation is even worse if $M$ is infinite-dimensional 
(cf.\ Example~\ref{ex:e.2.3} below). 
\end{rem}

\begin{lem} \mlabel{lem:e.1.16} 
The domain $\cV(M)_{\rm cp}$ of the exponential function of $\Diff(M)$ satisfies: 
\begin{description}
\item[\rm(1)] $0 \in \cV(M)_{\rm cp}$. 
\item[\rm(2)] $\R \cdot \cV(M)_{\rm cp} = \cV(M)_{\rm cp}$. 
\item[\rm(3)] For $X \in \cV(M)_{\rm cp}$, we have 
$\exp((t+s)X) = \exp(tX) \exp(sX)$ for $s,t \in \R$, i.e., 
$\gamma_X(t) := \exp(tX)$ is a group homomorphism 
$\gamma_X \: \R \to \Diff(M)$. 
\item[\rm(4)] $\Ad(\phi) \cV(M)_{\rm cp} = \cV(M)_{\rm cp}$ for all
$\phi \in \Diff(M)$ and $\exp \circ \Ad(\phi) = c_\phi \circ \exp$ on
$\cV(M)_{\rm cp}$. 
\item[\rm(5)] If $X, Y \in \cD_{\exp}$ commute, i.e., $[X,Y] =0$, then 
  \begin{description}
  \item[\rm(a)] $X + Y \in \cV(M)_{\rm cp}$. 
  \item[\rm(b)] $\exp(X + Y)=\exp X \exp Y = \exp Y \exp X$. 
  \item[\rm(c)] $\Ad(\exp X)Y = Y$. 
  \end{description}
\end{description}
\end{lem}

\begin{prf} (1) is trivial. 

 (2) Let $X \in \cV(M)_{\rm cp}$ and $s \in \R$. Consider 
the curve $\eta \: \R \to \Diff(M)$ defined by 
$\eta(t) := \gamma_X(st)$. It solves the initial value problem 
\[  \eta'(t) = s \eta(t).X = \eta(t).(sX), \quad \eta(0) = \id_M, \] 
which implies $\eta = \gamma_{sX}$ by uniqueness. Therefore 
$sX \in \cV(M)_{\rm cp}$. 

 (3) Fix $s \in \R$. Then the curve $\eta(t) := \gamma_X(s)^{-1}\gamma_X(t+s)$
solves the initial value problem 
\[  \eta'(t) = \gamma_X(s)^{-1}.\gamma_X'(t+s) =
\gamma_X(s)^{-1}.(\gamma_X(t+s).X) = \eta(t).X, \quad 
\eta(0) = \id_M,  \]
so that the uniqueness of the solutions (Proposition~\ref{prop:e.1.4}) implies 
$\eta(t) = \gamma(t)$ for all $t \in \R$, which yields (3). 

 (4) Let $X\in \cV(M)_{\rm cp}$. We consider the smooth curve 
$$ \gamma \: \R \to \Diff(M), \quad t \mapsto c_\phi(\exp tX) = \phi  \circ (\exp tX) \circ 
\phi^{-1}. $$
Then 
\begin{eqnarray*}
\gamma'(t) 
&=& \big(\phi \circ (\exp t X)\big).(X \circ \phi^{-1}) \\
&=& \big(\phi \circ (\exp t X) \circ \phi^{-1} \circ \phi\big).(X \circ \phi^{-1}) 
= \gamma(t).\Ad(\phi)X, 
\end{eqnarray*}
so that $\delta^l(\gamma) = \Ad(\phi)X$ is constant. This implies that
$\Ad(\phi)X\in \cV(M)_{\rm cp}$ with 
$$ \exp(\Ad(\phi)X) = \phi \circ (\exp X) \circ \phi^{-1}. $$

(5) Assertions (a) and (b) follow from Proposition \ref{prop:e.1.9}(b). 
For each $t \in \R$, we obtain from (4) and (5) that 
\[ \exp(tY) = \exp(X)\exp(tY) \exp(-X) = \exp(t\Ad(\exp X)Y), \] 
so that we simply have to take the derivative in $t = 0$ to obtain~(c). 
\end{prf}

The following lemma is very useful if one considers families of 
complete vector fields depending on a parameter. 
It provides a suitable differential of the exponential 
function in a very general context. 

\begin{lem} \label{lem:dexp-vf} Let $X, Y \in {\cal V}(M)$ be two smooth vector fields 
with the property that $E := \R X + \R Y$ consists of complete vector fields. 
If the function $\Exp := \exp_{\Diff(M)}\res_E$ is smooth, then we have in the locally convex space~${\cal V}(M)$: 
$$ \delta^l(\Exp)_X(Y) = \int_0^1 \Ad(\Exp(-tX))Y\, dt \quad \mbox{ for } \quad X,Y \in E. $$
\end{lem}

\begin{prf} We claim that the smooth function 
$$ \psi \: I \to {\cal V}(M), \quad 
\psi(t) := \delta^l(\Exp)_{tX}(tY), $$
satisfies the functional equation  
\begin{eqnarray}
  \label{eq:psi-eq}
\psi(t+s) = \Ad(\Exp(-sX))\psi(t) + \psi(s). 
\end{eqnarray}
We consider the three smooth functions 
$F,G,H \: E \to \Diff(M)$, given by 
$$ F(X) := \Exp((t+s)X), \quad 
G(X) := \Exp(tX) \quad \mbox{and} \quad H(X) := \Exp(sX),  $$
satisfying $F = G \cdot H$ pointwise. The 
Product Rule (Lemma~\ref{lem:e.1.2}) implies that 
$$ \delta^l(F) = \delta^l(H) + \Ad(H)^{-1}.\delta^l(G). $$
Now (\ref{eq:psi-eq}) follows from 
$$ \delta^l(G)_X(Y) = \psi(t), \quad 
\delta^l(H)_X(Y) = \psi(s) \quad \mbox{ and } \quad 
\delta^l(F)_X(Y) = \psi(s+t). $$

Clearly $\psi(0) = 0$ and 
$\psi'(0) = \lim_{t \to 0} \delta^l(\Exp)_{tX}(Y) = \delta^l(\Exp)_0(Y) = Y.$
Taking derivatives with respect to $t$ in $0$, (\ref{eq:psi-eq}) thus leads to 
$$ \psi'(s) = \Ad(\Exp(-sX))Y, $$
which implies the assertion by integration: 
\[  \delta^l(\Exp)_X(Y) = \psi(1) = \int_0^1 \psi'(s)\, ds 
= \int_0^1 \Ad(\Exp(-sX)).Y\, ds. \qedhere \] 
\end{prf}

\begin{small}
\subsection*{Exercises for Section~\ref{sec:e.1}} 

\begin{exer} Show that if $\sigma \: M \times G \to M$ is a smooth right action of 
the Lie group $G$ on $M$, then $T\sigma \: TM \times TG \to TM$ is a
smooth action of $TG$ on $TM$. \\
Hint: Express $\sigma(\sigma(p,g_1), g_2)= \sigma(p,g_1g_2)$ as 
$$ \sigma \circ (\sigma \times \id_G) = \sigma\circ (\id_M \times m_G) $$
and apply the tangent functor $T$. 
\end{exer}
  
\end{small}

\section{Smooth maps into general linear groups} \mlabel{sec:e.2}

Let $E$ be a locally convex space. 
Then the {\it general linear group} $\GL(E)$ is 
a subgroup of $\Diff(E)$, so that the results of the 
preceding section apply in particular to maps with values in $\GL(E)$. 
An important class of such maps 
are solutions of ordinary linear differential equations and 
smooth representations of Lie groups on $E$. In this context,  
we also discuss smooth representations of compact Lie groups 
and smooth representations of $\R$. 

\subsection*{Linear ordinary differential equations}

The following lemma provides a characterization of $\GL(E)$ 
as the centralizer of the scalar multiplications. Combined with the 
Uniqueness Lemma, it provides a bridge to the material in Section~\ref{sec:e.1}. 

\begin{lem} \mlabel{lem:e.2.1} For each $\lambda \in \R^\times$, we consider the diffeomorphism 
$\mu_\lambda \in \Diff(E)$ given by $\mu_\lambda(v) = \lambda v$. 
Then the following assertions hold: 
\begin{description}
\item[\rm(1)] $\GL(E) = \{ \phi \in \Diff(E) \: (\forall \lambda \in \R^\times)\ 
\phi \circ \mu_\lambda = \mu_\lambda \circ \phi\}$. 
\item[\rm(2)] $\{ X \in {\cal V}(E) \: (\forall \lambda \in \R^\times)\ 
\Ad(\mu_\lambda)X = X\}$ is the set of linear vector fields. 
\end{description}
\end{lem}

\begin{prf} (1) If $\phi \in \GL(E)$, then $\phi$ commutes with all scalar
multiplications $\mu_\lambda$. If, conversely, $\phi \in \Diff(E)$ commutes with
all scalar multiplications, then 
\[  \lambda \phi(0) = \phi(\mu_\lambda(0)) = \phi(0) \quad \mbox{ for all } \quad 
\lambda \in \R \] 
implies that $\phi(0) = 0$. We now obtain for each $\lambda\not=0$: 
$$ \phi(x) 
= \lambda^{-1} \phi(\lambda x) 
= \lambda^{-1} (\phi(0+ \lambda x) - \phi(0)), $$
and by passing to the limit $\lambda \to 0$, we obtain 
$\phi(x) = T_0(\phi)x,$
and hence that $\phi$ is linear. 

(2) Since $T(\mu_\lambda)$ is the multiplication with $\lambda$ on 
$T(E) \cong E \times E$, for $X \in {\cal V}(E) \cong C^\infty(E,E)$,  
the relation $\Ad(\mu_\lambda)X = X$ is equivalent to 
$$ X(\lambda v) = \lambda X(v) \quad \mbox{ for all } \quad \lambda \in \R, v \in E. $$
This condition is obviously satisfied for linear vector fields. 
Conversely, it implies that 
$$ X(v) = \lim_{t \to 0} \frac{1}{t} X(tv) = T_0(X)v, $$
showing that $X$ is linear. 
\end{prf}

\begin{defn} If $N$ is a smooth manifold, then we call a map \break 
$\phi \: N \to \GL(E)$ 
{\it smooth} if it is smooth in the
sense of Definition~\ref{def:e.1.1a}. 
\end{defn}

\begin{rem} \mlabel{rem:e.2.2} (a) We consider the Lie algebra 
$\gl(E) := ({\cal L}(E), [\cdot,\cdot])$ of continuous endomorphisms of $E$,  
endowed with the commutator bracket $[A,B] = AB - BA$ 
and observe that, for linear vector fields 
$X_A(v) := Av$, we have $[X_A, X_B] = - X_{[A,B]}$ 
(Exercise). 

For $\phi \in \GL(E) \subeq \Diff(E)$, we have 
$$ \Ad(\phi)X_A = T(\phi) \circ X_A \circ \phi^{-1} = X_{\phi \circ A \circ \phi^{-1}}, $$
so that the adjoint action of $\GL(E) \leq \Diff(E)$ 
on $\gl(E) \subeq {\cal V}(E)$ coincides with 
the conjugation action. 

(b) For each curve $u \: J \to \gl(E)$, a smooth solution 
$\gamma \: J \to \GL(E)$ of the initial value problem 
$$ \gamma'(t) = \gamma(t) \circ u(t), \quad \gamma(t_0)= \gamma_0, $$
can be viewed as solution of the equation 
$\delta^l(\gamma) = u$ on $\Diff(E)$. Therefore Proposition~\ref{prop:e.1.4} 
implies that it has at most one solution. 

(c) Assume that the smooth curve $\gamma \: J \to \Diff(E)$ satisfies 
$$ \gamma'(t) = \gamma(t) \circ u(t), \quad \gamma(t_0) \in \GL(E). $$
Then $\gamma(t_0)$ commutes with all scalar multiplications 
$\mu_\lambda$ and the vector fields $u(t) = \delta^l(\gamma)(t)$ on $E$ are
linear, hence invariant under each $\mu_\lambda$ (Lemma~\ref{lem:e.2.1}). 
Now Proposition~\ref{prop:e.1.8} implies that each $\gamma(t)$ commutes with each
$\mu_\lambda$, and Lemma~\ref{lem:e.2.1} yields $\gamma(J) \subeq \GL(E)$. 

This means that, if $u(t)$ consists of linear vector fields, 
then any solution of the corresponding IVP in the group $\Diff(E)$, 
starting in $\GL(E)$, remains in $\GL(E)$. \
It follows in particular that solvability of the corresponding 
problems on $\GL(E)$ and $\Diff(E)$ are equivalent. 

(d) From (c), it follows that, for each 
$D \in {\cal L}(E)$, the initial value problem 
$$ \gamma'(t) = \gamma(t) \circ D, \quad \gamma(0)= \id_E, $$
has at most one solution. If $E$ is a Banach space, then solutions
always exist, but this is not true in general, as Example \ref{ex:e.2.3} below shows. 

(e) The relations from Lemma~\ref{lem:e.1.16} also hold for the exponential function 
$\exp_{\GL(E)} \: {\cal D}_{\exp}(\GL(E))\to \GL(E)$. 
\end{rem}

\begin{defn} \mlabel{def:integrable-op} If $E$ is a locally convex space, then the 
domain of the exponential function of $\GL(E)$ is 
the set of all continuous linear operators $D$ on $E$ for which the corresponding
linear vector field $X_D(v) = Dv$ is complete, i.e., there exists a 
smooth representation $\alpha \: \R \to \GL(E)$ with 
$\alpha'(0) = D$. We call these operators $D$ {\it integrable}. 
\end{defn}

\begin{ex} \mlabel{ex:e.2.3} 
We consider the Fr\'echet space $E :=
C^\infty(]0,1[,\R)$ (Definition~\ref{def:smooth-co-top}) 
and the continuous 
linear operator $Df := f'$ on this space. Consider the initial value problem 
\begin{eqnarray}
  \label{eq:e.3}
\dot\gamma(t) = D \circ \gamma(t)= \gamma(t)', \quad \gamma(0) = \gamma_0. 
\end{eqnarray}
Suppose that $\gamma \: J \to E$ is a smooth solution. Then 
$$ \hat\gamma \: J \times ]0,1[ \to \R, \quad (t,s) \mapsto \gamma(t)(s) $$
is a smooth function satisfying the differential equation 
$$ \frac{\partial \hat\gamma}{\partial t} = \frac{\partial \hat\gamma}{\partial s}, $$ 
which implies that $\hat\gamma$ is constant on the intersection of the lines 
$t + s = c$ with $J \times ]0,1[$. If $\eps > 0$ and $[0,\eps] \subeq J$, it
follows that 
$t \mapsto \gamma(\eps)(t- \eps)$
defines a smooth extension of $\gamma_0$ to the interval $]-1,1 +
\eps[$ and likewise for solutions on $]-\eps,0[$. Hence the existence
of a global solution $\gamma \: \R \to E$ implies the existence of a
smooth extension of $\gamma_0$ to a function $\R \to \R$. 

This proves that there are functions $\gamma_0 \in E$ for which 
(\ref{eq:e.3}) has no solution on any non-trivial interval $J$ containing
$\{0\}$, and hence  there exists no non-trivial interval $J \subeq \R$ on
which the initial value problem 
\[ \dot\gamma(t) = \gamma(t) \circ D, \quad \gamma(0) = \id_E \]
on the group $\GL(E)$ has a smooth solution. 
\end{ex}

\begin{prop} \mlabel{prop:int-ad} If 
$M$ is a smooth finite-dimensional manifold and $X \in {\cal V}(M)$ a complete 
vector field, then the continuous linear operator 
$\ad X$ on ${\cal V}(M)$ generates a smooth $\R$-action 
given by $\alpha(t) := \Ad(\exp(-tX))$. 
\end{prop}

\begin{prf} Since $\gamma_X(t) := \exp(tX)$ defines a homomorphism
$\R \to \Diff(M)$ (Lemma~\ref{lem:e.1.16}(3)), 
$\alpha$ defines a linear action of $\R$ on 
${\cal V}(M)$. In view of Proposition~\ref{prop:smooth-logder} and 
Remark~\ref{rem:smooth-act}, 
the smoothness of this action follows from the smoothness of the map 
\begin{eqnarray*}
\R \times \cV(M) \times M &\to& T(M), \\ 
(t,Y,m) &\mapsto& \big(\Ad(\gamma_X(t))Y\big)(m)
= T(\gamma_X(t)) Y\big(\gamma_X(t)^{-1}(m)\big), 
\end{eqnarray*}
which in turn follows from the smoothness of 
\[ \R \times T(M) \to T(M), \quad (t,v) \mapsto T(\gamma_X(t))v \] 
(Lemma~\ref{lem:tang-mult-c}) and the smoothness of the evaluation map 
$$ \ev \: \cV(M) \times M \to T(M), \quad (Z,m) \mapsto Z(m), $$
which in turn follows from the smoothness of the evaluation map 
$$ \ev \: C^\infty(U,\R^n) \times U \to \R^n $$
for any open subset $U \subeq \R^n$ (Lemma~\ref{lem:smooth-eval-vec}). 
Finally, $\alpha'(0) = - \ad X$ follows from Lemma~\ref{lem:e.1.6}. 
\end{prf}

\subsection*{Smooth representations of Lie groups} 

\begin{defn} Let $E$ be a locally convex space and $G$ a Lie group. 
A~{\it smooth representation} of $G$ on $E$ is 
a homomorphism $\pi \: G \to \GL(E)$ which is smooth in the sense that the 
linear action $\hat\pi \: G \times V \to V$ is smooth. 
We mostly write $(\pi,E)$ for a representation $\pi$ of $G$ on~$E$. 
\end{defn} 

\begin{prop} [The derived representation] \mlabel{prop:e.2.6} 
Let $(E, \pi)$ be a smooth representation of the Lie group $G$ with 
Lie algebra $\g$. Then 
$$ d\pi \: \g \to \gl(E), \quad d\pi(X)v := Xv := T_\be(\pi)(x)(v) $$
is a homomorphims of Lie algebras. 
If $G$ is connected and $\alpha \: \g \to \gl(E)$ is a representation of 
$\g$ on $E$ defining on $E$ the structure of a topological $\g$-module, 
 then there is at most one smooth representation of 
$G$ on $E$ with $d\pi = \alpha$. 
\end{prop}

\begin{prf} To see that $d\pi$ is a homomorphism of Lie algebras, we 
apply Proposition~\ref{prop:e.1.5}(1) 
to the smooth right action $\sigma(v,g) := \pi(g^{-1})v$, 
and recall from Remark~\ref{rem:e.2.2} that the natural inclusion 
$\gl(E) \to {\cal V}(E), A \mapsto X_A$ is an antihomomorphism. The uniqueness assertion 
follows from Proposition~\ref{prop:e.1.5}(2). 
\end{prf}

The following lemma permits us to apply 
the Peter--Weyl Theorem in a context where the usual completeness requirements 
are not satisfied. 
\begin{lem} \mlabel{lem:rep-ext} Let $K$ be a compact topological group and 
$$\sigma \: K \times E \to~E, \quad (k,v) \mapsto \sigma_k(v)$$ 
a continuous action of $K$ on the locally 
convex space $E$. 
Then $\sigma$ extends uniquely to a continuous action $\hat\sigma$ of $K$ on 
the completion $\hat E$ of~$E$. 
\end{lem}

\begin{prf} First we claim that the topology of $E$ can be defined by 
$K$-invariant continuous seminorms. To see this, let 
$U \subeq E$ be a convex balanced $0$-neighborhood. Then 
$$ U_K := \bigcap_{k \in K} \sigma_k(U) $$
also is zero neighborhood in $E$ because 
$\sigma^{-1}(U) \subeq K \times E$ is an open neighborhood of the 
compact subset $K \times \{0\}$. By the Wallace Lemma~\ref{Wallla}, 
there exists an open $0$-neighborhood $W \subeq E$ 
with $K \times W \subeq \sigma^{-1}(U)$, and this 
means that $W \subeq U_K$. Clearly $U_K$ is $K$-invariant, so that its 
gauge functional 
$$ p_{U_K}(v)  := \inf \{ t > 0 \: v \in t U_K \} $$
is a $K$-invariant continuous seminorm on $E$ (Proposition~\ref{Minkowsk}). 
This shows that $E$ possesses a basis of $K$-invariant convex balanced $0$-neighborhoods, 
and hence that the topology on $E$ is defined by a collection 
${\cal P}$ of $K$-invariant seminorms. 

For $p \in {\cal P}$, we write $q_p \: E \to E_p$ for the quotient map 
into the corresponding normed space $E_p$. Since $p$ is $K$-invariant, the subspace 
$p^{-1}(0)$ is $K$-invariant, and we obtain an action 
$\sigma_p \: K \times E_p \to E_p$. From the $K$-invariance of $p$, we derive 
that this action is isometric. Moreover, 
$$ \sigma_p \circ (\id_K \times q_p) = q_p \circ \sigma, $$
and since $\id_K \times q_p$ is a quotient map, 
$\sigma_p$ is continuous. Let $\hat E_p$ denote the completion of $E_p$,  
which is a Banach space. Since each $\sigma_{p,k}$ is isometric, it extends 
uniquely to an isometry $\hat\sigma_{p,k}$ of $\hat E_p$. 
We thus obtain an action $\hat\sigma_p$ of $K$ on $\hat E_p$. 
The continuity of this action follows from the fact that $K$ 
acts by isometries and the set of elements of $\hat E_p$ with continuous orbits maps 
contains $E_p$, hence is dense, which implies the continuity of $\hat\sigma_p$ 
(Exercise~\ref{exer:cont-crit}). 

Forming the product space, we obtain a 
continuous action 
$\hat\sigma_{{\cal P}}$ of $K$ on $\prod_{p \in {\cal P}} \hat E_p$. 
Let 
$$q_{\cal P} \: E \to \prod_{p \in {\cal P}} E_p, \quad v \mapsto (q_p(v))_{p \in {\cal P}} $$
denote the canonical embedding. Then the completion $\hat E$ can be identified 
with the closure  
\[  \hat E \cong \oline{q_{\cal P}(E)} \subeq \prod_{p \in {\cal P}} \hat E_p, \] 
because products of Banach spaces are complete (Exercise). 
Since $q_{\cal P}$ is $K$-equivariant, $\oline{q_{\cal P}(E)}$ is $K$-invariant, 
and by restriction we obtain a continuous action of $K$  on $\hat E$. 
It is uniquely determined 
because each operator $\sigma_k \in \GL(E)$ extends uniquely to $\hat E$. 
\end{prf}

\begin{thm} [Peter--Weyl Theorem] \mlabel{thm:peter-weyl} 
Let $K$ be a compact Lie group acting smoothly 
on the Mackey complete locally convex space. Then the space of 
$K$-finite vectors is dense in $E$. 

Let $\hat K$ denote the set of irreducible characters of $K$. 
For each $\chi \in \hat K$, we then have a continuous projection 
$p_\chi \: E \to E_\chi$ onto the isotypic component $E_\chi$ of type $\chi$, 
and, for each $v \in E$, its Fourier series 
$$ v = \sum_{\chi \in \hat K} p_\chi(v) $$
converges in $E$. 
\end{thm}

\begin{prf} With Lemma~\ref{lem:rep-ext} we first extend the smooth action 
$\sigma$ of $K$ on $E$ to a continuous action $\hat\sigma$ of $K$ on the completion 
$\hat E$ on $E$. Then the Peter--Weyl Theorem \cite[Thm.~3.5.1]{HoM98} 
implies that the space $\hat E^{[K]}$ of $K$-finite vectors is dense in $\hat E$. 
This space is a direct sum 
$$\hat E^{[K]} = \bigoplus_{\chi \in \hat K} \hat E_{\chi}, $$
where 
$\hat E_{\chi}$ is the isotypic component corresponding to the 
irreducible character $\chi \in \hat K$. 

We further have continuous projections onto the isotypic components given by 
$$ p_\chi(v) = \chi(\be) \int_K \oline{\chi(k)} \sigma(k,v)\, d\mu_K(k), $$
where $\mu_K$ is the normalized Haar measure of $K$. 
Since $E$ is dense in $\hat E$, the sum of the subspaces 
$p_\chi(E)$ is also dense in $\hat E$. 
On the other hand, the Mackey completeness of $E$ implies that, for 
each $v \in E$, the element $p_\chi(v)$, which is given by an integral of a smooth function, 
is also contained in $E$. We conclude that 
$\sum_\chi p_\chi(E) \subeq E^{[K]}$ is dense in $E$. 

The convergence of the Fourier series of smooth vectors in $\hat E$ 
follows from Harish--Chandra's Theorem (\cite[Thm.~4.4.2.1]{Wa72}), but since 
$E$ consists of smooth vectors, the assertion follows. 
\end{prf}

\subsection*{Smooth representations of $\R$} 

Of particular importance are the smooth representation of 
the group $(\R,+)$. Since they correspond to smooth one-parameter subgroups 
of $\GL(E)$, any such representation $\alpha \: \R \to \GL(E)$ is uniquely 
determined by it {\it infinitesimal generator} $D := \alpha'(0)$, which 
is an element of ${\cal D}_{\rm exp}(\GL(E))$. Conversely, 
each element in the domain of the exponential function of $\GL(E)$ 
is the infinitesimal 
generator of a smooth representation of $\R$. 
Using the ``exponential notation'', we also write 
$e^{tD} := \exp_{\GL(E)}(tD) = \alpha(t)$ in this case. 

The following example shows that, for smooth representations of $\R$, 
invariance of a subspace under the derived action does not imply its invariance. 
\begin{ex} \mlabel{ex:e.2.9} 
Let $E := C^\infty(\R,\R)$ be the Fr\'echet space of smooth functions on $\R$ and 
consider the closed subspaces 
\[  W(a,b) := \{ f \in E \: \supp(f) \subeq [a,b]\}, \qquad a < b. \] 
These subspaces are invariant under the continuous linear operator 
$$ D\: E \to E, \quad Df = f'. $$
The curve 
$$ \alpha \: \R \to \GL(E), \quad \alpha(t)(f)(x) = f(t + x) $$
is smooth because the action 
$$ \R \times E \to E, \quad (t,f) \mapsto \alpha(t)f $$
is smooth (Remark~\ref{rem:smooth-act}), 
which in turn follows from the smoothness of the map 
$$ \R \times E \times \R \to \R, \quad (t,f,x) \mapsto f(x + t) $$
(Proposition~\ref{prop:cartes-closed}(ii)). In view of 
\[  \alpha'(t) = \alpha(t) \circ D \quad \mbox{ for }\quad  t \in \R, \] 
we have 
\[  D \in {\cal D}_{\exp}(\GL(E)) \quad 
\mbox{ and } \quad \exp(tD) =  \alpha(t). \]

The interesting point of this example is that the subspaces $W(a,b)$ are invariant 
under the infinitesimal generator $D$ of the one-parameter group $\alpha$, but not 
under $\alpha$ itself. In fact, we have 
$$ \alpha(t)W(a,b) = W(a-t,b-t). $$
\end{ex}

\subsection*{Singularities of the exponential function}

The following proposition is one of our main tools for the analysis of the 
singularities of the exponential function of a Lie group. 

\begin{prop} \mlabel{prop:e.2.10} Let $E$ be a Mackey complete and 
$\alpha \: \R \to \GL(E)$ a smooth representation with 
$\alpha'(0) = D$. Then the following assertions hold: 
\begin{description}
\item[\rm(1)] The operator-valued integral 
$$ \beta(t) := \int_0^1 \alpha(st)\, ds = \int_0^1 e^{stD}\, ds \in \gl(E) $$ 
exists pointwise. 
\item[\rm(2)] The map 
$\hat\beta \: \R \times E \to E, (t,v) \mapsto \beta(t)v$
is smooth. 
\item[\rm(3)] For each $t \in \R^\times$,  we have 
$D \circ \beta(t) = \beta(t) \circ D = \frac{\alpha(t) - \1}{t}.$
\item[\rm(4)] Assume that $E$ is a Banach space. 
  \begin{description}
  \item[\rm(a)] Then $\beta \: \R \to {\cal L}(E)$ is real analytic and given by 
the power series 
$$ \beta(t) 
= \frac{e^{tD} - \1}{tD} 
= \sum_{k = 0}^\infty \frac{t^k}{(k+1)!}  D^k. $$ 
  \item[\rm(b)] $\beta(1)$ is invertible if and only if $\Spec(D) \cap 2\pi i \Z \subeq \{0\}$.
If $\|D\| < 2\pi$, this condition is satisfied and $\beta(1)^{-1}$ is given by the 
power series 
$$ g(D) := \frac{D}{e^D - 1} := \sum_{n = 0}^\infty \frac{B_n}{n!} D^n,  $$
where $(B_n)_{n \in \N_0}$ are the {\it Bernoulli numbers}. 
\end{description}  
\end{description}  
\end{prop}

\begin{prf} (1) follows from the fact that, for each $v \in E$, the curve 
$s \mapsto \alpha(st)v$ is smooth and $E$ is Mackey complete. 

(2) follows from $\hat\beta(t)v = \int_0^1 \alpha(st)v\, ds$ and 
the smoothness of the map 
$$ I \times \R \times E \to E, \quad (s,t,v) \mapsto \alpha(st)v $$
(Lemma~\ref{lem:e.1.16}). 

(3) From $D = \Ad(\alpha(t))D = \alpha(t)D \alpha(t)^{-1}$ for each $t \in \R$  
(Lemma~\ref{lem:e.1.16}(5)), we obtain for each $v \in E$ 
\[  D t\beta(t)v 
= \int_0^1 D t \alpha(st)v\, ds 
= \int_0^t D \alpha(s)v\, ds 
= \int_0^t \alpha'(s)v\, ds = \alpha(t)v - v, \] 
and further 
$$  t\beta(t)Dv 
= \int_0^1 t \alpha(st)Dv\, ds 
= \int_0^t D\alpha(s)v\, ds = t D \beta(t)v. $$

(4) (a) Let $f(z) := \frac{e^z - 1}{z}$ and observe that $f$ defines an entire 
function. Therefore the holomorphic functional calculus (Section~\ref{sec:cia-1}) 
provides 
a holomorphic map 
$$ \C \to {\cal L}(E)_\C \cong {\cal L}(E_\C), \quad 
z \mapsto f(zD), $$
whose restriction to the real line is~$\beta$. 

(4)(b) follows directly from the Spectral Mapping Theorem 
(\cite[Thm.~10.28]{Ru91}) and the fact that $f^{-1}(0)= 2 \pi i \Z \setminus \{0\}$. 
If, in addition, $\|D\| < 2\pi$, then 
the relation $g(z)f(z) = 1$ for $z \in \C$ with $f(z)\not=0$ 
implies that the power series $g(D)$ converges to an operator whose inverse if $f(D)$. 
\end{prf}

\begin{rem} The {\it Bernoulli numbers} $(B_n)_{n \in \N_0}$ are defined by the 
power series 
$$ g(z) := \frac{z}{e^z - 1} = \sum_{n = 0}^\infty \frac{B_n}{n!} z^n,  $$
which defines a holomorphic function on 
$\C \setminus \{ k 2\pi i \: 0 \not=k \in \Z\}$. 
They arise in particular (for $p \in \N$) in the summation formula 
$$ s_p(n) := \sum_{k = 1}^{n-1} k^p 
= \frac{1}{p+1} \sum_{k = 0}^{p} B_k \begin{pmatrix}p+1 \\ k\end{pmatrix} n^{p-k+1} 
= p! \sum_{k = 0}^{p} \frac{B_k}{k!} \frac{n^{p-k+1}}{(p-k+1)!}. $$

Below we explain how this formula can be derived from 
Proposition~\ref{prop:e.2.10}. 
Let $E \subeq C^\infty(\R,\R)$ denote the finite-dimensional 
space of polynomials of 
degree $\leq N$. We then obtain a smooth linear action of $\R$ on 
this space by 
$$ (\alpha(t)f)(x) = f(t+x). $$
Its infinitesimal generator is given by $Df = f'$, so that $D^{N+1} = 0$.  
Proposition~\ref{prop:e.2.10}(3) implies that the operator 
\[  A := \sum_{k = 0}^\infty \frac{1}{(k+1)!} D^k = \beta(1) \] 
satisfies 
\begin{eqnarray}
  \label{eq:fund-rel}
 D \circ A = \Delta := \alpha(1) - \id_E = e^D - \id_E. 
\end{eqnarray} 
Explicitly, we have 
\[ (\Delta f)(x) = f(x+1) - f(x) 
\quad \mbox{ and } \quad 
(Af)(x) = \int_0^1 f(x + t)\, dt.\] 
Since $D$ is nilpotent, $A$ is unipotent with 
$A^{-1} = \sum_{k = 0}^\infty \frac{B_k}{k!} D^k.$

Suppose that $p < N$. Then the monomial $x^p$ is contained in the image of 
the differentiation operator $D$, and (\ref{eq:fund-rel}) implies that 
$x^p = \Delta F$ for some polynomial $F \in E$ for which we may 
assume that $\int_0^1 F(s)\, ds= 0$. In view of $1 = \dim(\ker(D)) 
= \dim(\ker \Delta)$, the latter condition determines $F$ uniquely. Now $F$ satisfies 
$F(x+1) - F(x) = x^p$, which immediately leads to 
$$ s_p(n) 
= \sum_{k = 0}^{n-1} k^p 
= \sum_{k = 0}^{n-1} F(k+1) - F(k) = F(n) - F(0). $$

On the other hand, $x^p = (\Delta F)(x) =  (DAF)(x)$ 
and $(AF)(0) = \int_0^1 F(s)\, ds = 0$ by assumption, which leads to 
$AF = \frac{x^{p+1}}{p+1}$, and hence to 
$$ F(x) = A^{-1}\frac{x^{p+1}}{p+1}
= \sum_{k = 0}^\infty \frac{B_k}{k!} D^k\frac{x^{p+1}}{p+1} 
= \frac{1}{p+1} \sum_{k = 0}^{p+1} B_k {p+1 \choose k}x^{p-k+1}, $$
which implies that 
$$ s_p(n) = F(n) - F(0) 
= \frac{1}{p+1} \sum_{k = 0}^{p} B_k {p+1 \choose k} n^{p-k+1}. $$
\end{rem}

Although a holomorphic functional calculus is not available for continuous 
operators on locally convex spaces, the following proposition provides a 
quite handy criterion for the injectivity of certain operators related 
to the differential of the exponential function of a Lie group. 

\begin{prop} \mlabel{prop:kerexp-locconv}
Let $E$ be a Mackey complete locally convex space, \break 
$D \in \cD_{\exp}(\GL(E))$ and $A := \int_0^1 e^{tD}\, dt$. 
Then 
\begin{eqnarray}
  \label{eq:ev1}
\ker(e^D - \1) = \oline{\bigoplus_{n \in \N_0} \ker(D^2 + n^2 4 \pi^2 \1)}
\end{eqnarray}
and 
\begin{eqnarray}
  \label{eq:ev2}
\ker(A) = \oline{\bigoplus_{n \in \N} \ker(D^2 + n^2 4 \pi^2 \1)}. 
\end{eqnarray}
If $E$ is complex and $D$ is complex linear, then 
\begin{eqnarray}
  \label{eq:ev3}
\ker(e^D - \1) = \oline{\bigoplus_{n \in \Z}  \ker(D - n 2\pi i \1)}
\end{eqnarray}
and 
\begin{eqnarray}
  \label{eq:ev4}
\ker(A) = \oline{\bigoplus_{n \in \Z \setminus \{0\}} \ker(D - n 2\pi i \1)}. 
\end{eqnarray}
In particular, $e^D = \1$ and $\ker(D- n2\pi i) = \{0\}$ for each 
$n \in \Z \setminus \{0\}$ imply $D = 0$. 
\end{prop} 

\begin{prf} Since $e^D$ is a continuous endomorphism of $E$, its fixed point space 
$W := \ker(e^D - \1)$ is closed, hence also Mackey complete. 
By 
$$ \sigma \: \T \times W \to W, \quad (e^{2\pi it},w) \mapsto e^{tD}.w $$
we obtain a smooth action of the circle group $\T \cong\R/\Z$ on $W$. 
Since the group $\T$ is compact, the 
Peter--Weyl Theorem~\ref{thm:peter-weyl} implies that the $\T$-finite 
vectors span a dense subspace of $W$ on which the representation of $\T$ is semisimple, 
i.e., a direct sum of simple modules. 

If $E$ is complex, then the simple $\T$-submodules consist of eigenvectors. 
As a $\T$-eigenvector $w \in W$ corresponds to the character $\chi_n(z) = z^n$ of 
$\T$ if and only if $Dw = 2 \pi i n w$, the description of 
$\ker(e^D - \1)$ follows. 

If $E$ is real, then the preceding argument applies to the complexification 
$E_\C$ of $E$, and the assertion follows from 
\[  \ker(D - n 2\pi i \1) \oplus \ker(D + n 2\pi i \1) 
= \ker(D^2 + n^2 4 \pi^2 \1). \]

From Proposition~\ref{prop:e.2.10}, we get 
$D A = e^D - \1$, which leads to $\ker(A) \subeq \ker(e^D - \1)$. 
For $Dw = 2 \pi i n w$ we have 
$$ Aw = \int_0^1 e^{t2\pi i n} w\, dt = 
\left\{ \begin{array}{cl} 
0 & \mbox{ for $n \not= 0$} \cr 
w & \mbox{ for $n = 0.$} 
\end{array} \right.$$ 
Since $A$ commutes with the action of the circle group $\T$ on the space $W$, 
its kernel is a $\T$-invariant closed subspace, hence topologically 
generated by the 
$\T$-eigenvectors. This implies (\ref{eq:ev4}) 
if $E$ is complex and $D$ is complex linear. The real case is obtained immediately 
from the complex case by considering the complexification $E_\C$ of $E$.
\end{prf} 

\begin{cor} \mlabel{cor:kerexp-cia} 
Let $\cA$ be a Mackey complete complex cia and $x \in \cA$. 
\begin{description}
\item[\rm(1)] $\ker (\lambda_{e^x} - \1) = \bigoplus_{n \in \Z} \ker(\lambda_x 
- n 2\pi i \1)$, where the sum on the right is finite. 
\item[\rm(2)] If $e^x = \1$, then $x$ is a linear combination of mutually orthogonal idempotents 
with coefficients in $2\pi i \Z$. 
\end{description}
\end{cor}

\begin{prf} (1) For each element $x \in \cA$, the linear action 
$\R \times \cA \to \cA, (t, a) \mapsto e^{tx} a$
is smooth with infinitesimal generator $D = \lambda_x$. Hence 
$\lambda_x \in \cD_{\exp}(\GL(\cA))$, and the assertion follows from 
Proposition~\ref{prop:kerexp-locconv} and the boundedness 
of the $\Spec(x)$ (Proposition~\ref{prop:spec-compact}). 

(2) is a direct consequence of (1). 
\end{prf} 

\begin{ex} (a) Let $V$ be a locally convex space. We consider the space 
$E := C^\infty(\R,V)$ with the operator 
$(D \phi)(x) := ix \phi(x)$ which generates a smooth $\R$-action given by 
$$ (\alpha(t) \phi)(x) =  e^{itx} \phi(x). $$
The operator $A$ is given by $A \phi = f \phi$ for 
$f(x) := \frac{e^{ix} - 1}{ix}$. 
Hence $A$ is injective because the zeros of $f$ are isolated, 
lying in the points $2\pi n$, $0 \not=n \in \Z$. 
Therefore each function in $\im(A)$ vanishes in 
the zero set of $f$. 

Conversely, let $\phi \: \R \to V$ be a smooth function vanishing in the 
zeros of $f$. Let $x_0 := 2 \pi n$ be such a zero. Then 
$$ \frac{\phi(x_0 + h)}{h} = \frac{\phi(x_0 + h) - \phi(x_0)}{h} 
= \int_0^1 \phi'(x_0 + th)\, dt $$
is a smooth function of $h$. 
Therefore  the function 
$x \mapsto \phi(x)f(x)^{-1}$ on the complement of $f^{-1}(0)$ extends to a smooth 
function $h$ on $\R$ with $fh = \phi$. 
We conclude that 
$$ \im(A) = \{ \phi \in C^\infty(\R,V) \: (\forall n \in \Z \setminus \{0\})\ \phi(2\pi n) =0\}. $$

(b) Let $E := {\cal O}(\C)$ denote the space of entire functions on $\C$. 
On this space we consider the multiplication operator 
$(D \phi)(z) := z\phi(z)$ which generates a smooth $\R$-action given by 
$$ (\alpha(t) \phi)(z) =  e^{tz} \phi(z). $$
Then $A\phi = f \phi$ for 
$f(z) := \frac{e^z - 1}{z}$, and $A$ is injective. 
Each function in $\im(A)$ vanishes in the zero set 
$2 \pi i \Z \setminus \{0\}$
of $f$. For any entire function $\phi \in E$ vanishing in $f^{-1}(0)$, the function 
$f^{-1} \cdot \phi$ is meromorphic on $\C \setminus f^{-1}(0)$, and in 
each point $2\pi i n$, $0 \not= n \in \Z$, it has a removable singularity. 
Therefore it is holomorphic, and 
$$ \im(A) = \{ \phi \in {\cal O}(\C) \: (\forall n \in \Z \setminus \{0\})\ \phi(2\pi i n) =0\}. $$

Recall the Euler product expansion 
\[  f(z) = \frac{e^z - 1}{z} 
= e^{z/2} \frac{e^{z/2} - e^{-z/2}}{z}
= e^{z/2} \frac{\pi}{2} \prod_{n = 1}^\infty \Big(1 + \frac{z^2}{4 \pi^2 n^2}\Big) \] 
converging in ${\cal O}(\C)$ (\cite[\S 1.3]{Re95}).  This product expansion 
might also be of some interest in the analysis of 
$A = f(D)$ in more general situations. In particular, we may always consider the 
expression 
\[  \Bigg(\prod_{n = 1}^\infty \Big(1 + \frac{D^2}{4 \pi^2 n^2}\Big)\Bigg)\phi 
\quad \mbox{ for } \quad \phi \in E, \]
and ask if the sequence corresponding to the finite products converges in $E$. 
This depends on the behavior of the sequence $(D^{2n}\phi)_{n \in \N}$. 
\end{ex} 

\subsection*{ad-integrable Lie algebras} 

\begin{defn} \mlabel{def4.1.5}  A locally convex Lie algebra $\g$ is called 
{\it $\ad$-integrable} if $\ad \g \subeq {\cal D}_{\exp}(\GL(\g))$, i.e., the
exponential function of the group $\GL(\g)$ is defined on the Lie
subalgebra $\ad \g$ of $\gl(\g)$. In this case, we write 
$$ e^{\ad x} := \exp(\ad x) \quad \hbox{ for } \quad x \in \g.  $$
If the integral 
\begin{equation}
  \label{eq:kappag}
\kappa_\g(x)(y) := \int_0^1 e^{-t\ad x}y\, dt 
\end{equation}
exists for  $x$ and $y$ in $\g$ (which is the case if $\g$ is Mackey complete), 
then we also write 
$$ \kappa_\g(x) := \int_0^1 e^{-t\ad x}\, dt. $$

Note that this integral also exists if 
there exists a Lie group $G$ with Lie algebra $\g$ and a smooth 
exponential function $\exp_G \: \g \to G$. In the latter case, we have 
$$ \delta(\exp_G)_x = \kappa_\g(x) $$
for each $x \in \g$ (Theorem~\ref{thm:exp-logder}). 

We define the {\it point-regular set of $\g$} by 
$$ \Reg_p(\g) := \{ x \in \g \: \ker \kappa_\g(x) = \{0\}\}, $$
i.e., the set of all points $x \in \g$ for which $0$ is not contained in the 
point spectrum of $\kappa_\g(x)$. 
\end{defn} 

\begin{rem} If $\g$ is Mackey complete, then Proposition~\ref{prop:kerexp-locconv} 
implies that 
$$ \Reg_p(\g) := \{ x \in \g \: (\forall n \in \N)\ \ker((\ad x)^2 + 4 \pi^2 n^2\1)
 = \{0\}\}.$$ 
If $\Spec_p$ denotes the set of all eigenvalues of the complexification of an 
operator, we obtain 
$$ \{ x \in \g \: \sup|\Im(\Spec_p(\ad x))| < 2 \pi \} \subeq \Reg_p(\g). $$ 
For a Banach--Lie algebra $\g$ for which the norm satisfies 
$\|\ad x \| \leq \|x\|$ we thus obtain 
$$ B_{2\pi}(0) := \{ x \in \g \: \|x\| < 2 \pi \} \subeq \Reg_p(\g). $$
\end{rem}

\begin{lem} \mlabel{lem4.1.9} If $\psi \: \g \to \h$ is a morphism of $\ad$-integrable
Lie algebras, then 
$$ \psi \circ e^{\ad x} = e^{\ad \psi(x)} \circ \psi 
\quad \hbox{ for each} \quad x \in \g $$
and 
$$ \psi \circ \kappa_\g(x) = \kappa_\h(\psi(x)) \circ \psi $$
if $\kappa_\g(x)$ and $\kappa_\h(\psi(x))$ are defined by \eqref{eq:kappag}. 
\end{lem}

\begin{prf} We consider the continuous endomorphism 
$\tilde\psi$ of $\g \times \h$ given by 
$\tilde \psi(x,y) = (0, \psi(x))$. For a pair 
$(\alpha,\beta) \in {\cal L}(\g) \times {\cal L}(\h) \subeq {\cal L}(\g \times \h)$,  
we then have 
$$ \tilde \psi(\alpha(x), \beta(y)) = (0, \psi(\alpha(x))) 
\quad \hbox{ and } \quad 
(\alpha,\beta)\tilde\psi(x,y) = (0, \beta(\psi(x))). $$
Therefore $\tilde\psi$ commutes with $(\alpha,\beta)$ if and only if 
$\beta \circ \psi = \psi \circ \alpha$. 

Now we consider the smooth curve 
\begin{eqnarray*}
\gamma &\:& \R \to \GL(\g) \times \GL(\h) \subeq \GL(\g \times \h) \\ 
\gamma(t) &:=& (e^{t \ad x}, e^{t \ad \psi(x)}) 
= \exp\big(t(\ad x,\ad \psi(x))\big). 
\end{eqnarray*}
Then $\gamma(0) = (\id_\g, \id_\h) = \id_{\g \times \h}$ commutes with 
$\id_{\g \times \h} + \tilde\psi \in \GL(\g \times \h)$, and 
$\delta^l(\gamma) = (\ad x, \ad \psi(x))$ satisfies 
\begin{eqnarray*}
& & \Ad(\1 + \tilde\psi)(\ad x, \ad \psi(x)) 
= (\1 + \tilde\psi)(\ad x, \ad \psi(x)) (\1 + \tilde\psi)^{-1} \\ 
&=& (\1 + \tilde\psi)(\ad x, \ad \psi(x)) (\1 - \tilde\psi)
= (\ad x, \ad \psi(x)) + [\tilde\psi, (\ad x, \ad \psi(x))]. 
\end{eqnarray*}
Since $\psi \circ \ad x = (\ad \psi(x)) \circ \psi$ follows from the
fact that $\psi$ is a morphism of Lie algebras, the map 
$\tilde\psi$ commutes with $(\ad x, \ad \psi(x))$, and we obtain 
$$ \Ad(\1 + \tilde\psi)(\ad x, \ad \psi(x)) 
= (\ad x, \ad \psi(x)). $$
Now Proposition~\ref{prop:e.1.8} 
shows that $\1 + \tilde\psi$, and hence
$\tilde\psi$, commutes with $\gamma(1)$, which proves the first assertion of the lemma. 
The second one follows immediately by integration. 
\end{prf}

\begin{rem} Let $\g$ be an $\ad$-integrable Lie algebra for which 
\[ \kappa_x(y) := \int_0^1 e^{-t \ad x}y\, dt \] 
defines a smooth $\g$-valued $1$-form on $\g$. 
We further assume that the function 
\[ E \:  \g \to \Aut(\g), \quad E(x) := e^{\ad x} \] 
is smooth (at least on $2$-dimensional subspaces), 
so that Lemma~\ref{lem:dexp-vf} implies that 
\[ \delta(E)_x(y) 
= \int_0^1 e^{- t\ad_{\ad x}} \ad y\, dt 
= \int_0^1 e^{- t\ad x} \circ \ad y \circ e^{t\ad x} \, dt.\] 

We now show that $\kappa$ always 
satisfies the Maurer--Cartan equation, which in this case can be written 
\begin{equation}
  \label{eq:mc-x}
\partial_y\big(\kappa_x(z)\big) - \partial_z\big(\kappa_x(y)\big) 
= [\kappa_x(z), \kappa_x(y)]\quad \mbox{ for }\quad 
x,y,z \in \g,
\end{equation}
where $\partial_y$ is the partial derivative with respect to $x$ in 
the direction of $y \in \g$. The right hand side is 
\[ [\kappa_x(z), \kappa_x(y)] 
= \int_0^1  \int_0^1 [e^{-t \ad x}z, e^{-s \ad x}y]\, dt\, ds. \] 
We further obtain 
\begin{align*}
&\partial_y \kappa_x(z) 
= \int_0^1 e^{-t \ad x} \delta(E)_{-tx}(-t y)z\, dt \\ 
&= \int_0^1 \int_0^1 e^{-t \ad x} e^{ st\ad x}(-t \ad y) e^{-ts\ad x}z \, ds\, dt \\ 
&= \int_0^1 \int_0^t e^{-t \ad x} e^{ s\ad x}(-\ad y) e^{-s\ad x}z \, ds\, dt \\ 
&= \int_0^1 \int_0^t e^{-(t-s) \ad x}(-\ad y) e^{-s\ad x}z \, ds\, dt \\
&= \int_0^1 \int_0^t e^{-(t-s) \ad x}[e^{-s\ad x}z,y] \, ds\, dt \\ 
&= \int_0^1 \int_0^t [e^{-t \ad x}z,e^{-(t-s) \ad x}y] \, ds\, dt 
= \int_0^1 \int_0^t [e^{-t \ad x}z,e^{-s \ad x}y] \, ds\, dt \\ 
&= \int_{0 \leq s \leq t\leq 1} [e^{-t \ad x}z,e^{-s \ad x}y] \, ds\, dt. 
\end{align*}
We likewise get 
\begin{align*}
-\partial_z \kappa_x(y) 
&= -\int_{0 \leq s \leq t\leq 1} [e^{-t \ad x}y,e^{-s \ad x}z] \, ds\, dt\\
&= \int_{0 \leq t \leq s\leq 1} [e^{-t \ad x}z, e^{-s \ad x}y] \, ds\, dt. 
\end{align*}
From these calculations \eqref{eq:mc-x} follows.
\end{rem}

\subsection*{Automorphisms and derivations} 

\begin{lem} [Linear Quotient Rule] \mlabel{lem:lin-quot}  
If $E$ is a locally convex space, $J \subeq \R$ an interval 
and $\gamma \: J \to \GL(E)$ smooth, then 
\[  (\gamma^{-1})'(t)  = - \gamma(t)^{-1} \circ \gamma'(t) \circ \gamma(t)^{-1} 
\quad \mbox{ for } \quad  t \in J.\] 
\end{lem}

\begin{prf} Differentiating the relation $\gamma(t)\gamma(t)^{-1}v = v$ leads for 
each $v \in E$ to 
$$ \gamma'(t) \gamma(t)^{-1}v + \gamma(t) (\gamma^{-1})'(t)v = 0, $$
and the assertion follows. 
\end{prf}

\begin{prop} \mlabel{prop:aut-der} Let $\cA$ be a locally convex space, endowed with a continuous 
bilinear multiplication $m_\cA \: \cA \times \cA \to \cA, (a,b) \mapsto ab$. We write 
$$ \der(\cA,m_\cA) := \{ D \in \gl(\cA) \: (\forall a,b \in \cA)\, D(ab) = D(a)b + a D(b)\} $$
for the Lie algebra of derivations of $(\cA,m_\cA)$ 
and 
$$ \Aut(\cA,m_\cA) := \{ \phi \in \GL(\cA) \: (\forall a,b \in \cA)\, \phi(ab) = \phi(a)\phi(b)\}$$
for its automorphism group. 
Then the following assertions hold for any interval $J \subeq \R$: 
\begin{description}
\item[\rm(1)] If $\gamma \: J \to \GL(\cA)$ is a smooth curve with values in 
$\Aut(\cA,m_\cA)$, then $\delta^l(\gamma)$ has values in $\der(\cA,m_\cA)$. 
\item[\rm(2)] If $\gamma \: J \to \GL(\cA)$ is a smooth curve for which 
$\delta^l(\gamma)$ has values in $\der(\cA,m_\cA)$ and $\gamma(t_0)\in \Aut(\cA,m_\cA)$ for some 
$t_0 \in J$, then 
$\gamma(t) \in \Aut(\cA,m_\cA)$
for all $t \in J$. 
\end{description}
\end{prop}

\begin{prf} (1) Let $a,b \in \cA$. For each $t \in J$, we then have 
$$ \gamma(t)(ab) = \gamma(t)(a) \cdot \gamma(t)(b), $$
and taking derivatives leads to 
$$ \gamma'(t)(ab) = \gamma'(t)(a) \cdot \gamma(t)(b) + \gamma(t)(a) \cdot \gamma'(t)(b), $$
so that 
$$  \delta^l(\gamma)_t(ab) = \gamma(t)^{-1}\circ \gamma'(t)(ab) 
=\delta^l(\gamma)_t(a) \cdot b + a \cdot \delta^l(\gamma)_t(b). $$

(2) Now we assume that $\delta^l(\gamma)$ has values in $\der(\cA,m_\cA)$. We consider the 
smooth curve $\alpha(t) := \gamma(t)\big(\gamma(t)^{-1}(a)\cdot\gamma(t)^{-1}(b)\big)$ with 
$\alpha(t_0) = ab$. We have to show that $\alpha$ is constant, i.e., that 
$\alpha'(t) = 0$. Using Lemma~\ref{lem:lin-quot}, this follows from 
\begin{eqnarray*}
\alpha'(t) 
&=&  \gamma'(t)\big(\gamma(t)^{-1}(a)\cdot \gamma(t)^{-1}(b)\big)
- \gamma(t)\Big( \gamma(t)^{-1}\gamma'(t)\gamma(t)^{-1}(a)\cdot \gamma(t)^{-1}(b)\Big)\\
&&\qquad \qquad \ \ \ \ - \gamma(t)\Big( \gamma(t)^{-1}(a) \cdot \gamma(t)^{-1}\gamma'(t) \gamma(t)^{-1}(b)\Big) \\
&=&  \gamma(t)\Big( \delta^l(\gamma)_t\big(\gamma(t)^{-1}(a)\cdot\gamma(t)^{-1}(b)\big)
- \delta^l(\gamma)_t\big(\gamma(t)^{-1}(a)\big)\cdot \gamma(t)^{-1}(b) \\
&&\qquad \qquad\ \ \ \ - \gamma(t)^{-1}(a) \cdot \delta^l(\gamma)_t\big(\gamma(t)^{-1}(b)\big)\Big)=0.
\quad\qquad\qquad\qquad\qquad\qquad\qedhere\end{eqnarray*}
\end{prf}

\begin{cor} Let $(\cA,m_\cA)$ be a locally convex algebra, 
$M$ a connected smooth manifold and $f \: M \to \GL(\cA)$ be a smooth map. 
Then the following assertions hold:
\begin{description}
\item[\rm(1)] If $f(M) \subeq \Aut(\cA,m_\cA)$, then $\delta^l(f)$ has values in 
$\der(\cA,m_\cA)$. 
\item[\rm(2)] If $\delta^l(f)$ has values in $\der(\cA,m_\cA)$ and $f(m_0)\in \Aut(\cA,m_\cA)$ for one 
$m_0 \in M$, then $f(M) \subeq \Aut(\cA,m_\cA)$. 
\end{description}
\end{cor}

\begin{small}
\subsection*{Exercises for Section~\ref{sec:e.2}} 

\begin{exer} Let $m_\cA \: \cA \times \cA \to \cA$ be a continuous bilinear map on the 
locally convex space $\cA$, considered as a multiplication. We consider the diffeomorphism 
\[  \phi_\cA \: \cA^3 \to \cA^3, \quad \phi(a,b,c) := (a,b,c + m_\cA(a,b)). \] 
Show that: 
\begin{description}
\item[(1)] $\Aut(\cA,m_\cA) = \{ \psi \in \GL(\cA) \: \phi_\cA \circ \psi^{\times 3} 
= \psi^{\times 3} \circ \phi_\cA\}$. 
\item[(2)] For $D \in \gl(\cA)$, let $D^{\times 3}$ denote the linear operator on 
$\cA^3$, given by $D^{\times 3}(a,b,c) := (Da,Db,Dc)$, which we also 
consider as a linear vector field. 
Then 
$$ \der(\cA,m_\cA) = \{ D \in \gl(\cA) \: \Ad(\phi_\cA)(D^{\times 3}) = D^{\times 3}\}. $$
\item[(3)] Use (1), (2) and Proposition~\ref{prop:e.1.8} to give another proof of 
Proposition~\ref{prop:aut-der}. 
\end{description}
\end{exer}

\begin{exer} \mlabel{ex:e.1} Let $V$ be a smooth module of the connected 
Lie group $G$. Show that $V$ is a trivial $G$ module if and only if 
$V$ is a trivial $\g$-module. \\
Hint: The ``smooth'' function $\rho_V \: G \to \GL(V)$ satisfies 
$\delta(\rho_V) = \dot\rho_V \circ \kappa_G.$
Now use the Uniqueness Lemma. 
\end{exer}

\begin{exer} \mlabel{ex:e.2} Let $V$ be a smooth module of the Lie group $G$ and 
$W \subeq V$ a closed submodule. Show that the 
natural action of 
$G$ on the quotient space $V/W$ induces on $V/W$ a smooth 
$G$-module structure.\\
 Hint: Observe first that the natural map 
$G \times V \to V/W$ is smooth and then use Lemma~\ref{lemquot} to derive the assertion. 
\end{exer}

\begin{exer} \mlabel{exer:rep-embed} Let $\alpha \: G \to \GL(E)$ define a smooth representation 
of the Lie group $G$ on the locally convex space $E$ and 
consider for each $v \in E$ the smooth function 
$\alpha^v \:G \to E, g \mapsto \alpha(g)v$. 
\begin{description}
\item[(1)] $\alpha^{\alpha(g)v} = \alpha^v \circ \rho_g$ for $v \in E$, $g \in G$. 
\item[(2)] The map $\Phi \: E \to C^\infty(G,E), v \mapsto \alpha^v$ is continuous and injective. 
\item[(3)] The evaluation map $\ev_\be \: C^\infty(G,E) \to E, f \mapsto f(\be)$ satisfies $\ev_\be \circ \Phi = \id_E$. 
\item[(4)] $\Phi$ is a topological embedding. 
\end{description}
\end{exer} 
  
\begin{exer}
\mlabel{exer:cont-crit}
Let $G$ be a topological group, let $E$ be a normed space and 
$\pi \: G \to \Isom(E)$ be a group homomorphism. Show that the corresponding 
action of $G$ on $E$ is continuous if the subset $E_c \subeq E$ of elements 
with continuous orbit maps is dense. 
\end{exer}

\begin{exer} \mlabel{exer:appe-contract} Let $E$ be a locally convex space and 
$\phi \in \Diff(E)$ with $\phi(0)=0$. Show that 
\[ \psi(t)(x) :=
\begin{cases}
  \frac{1}{t} \phi(tx) & \text{ for }\ t > 0 \\
  T_0(\phi)x  & \text{ for }\ t = 0 
\end{cases} \] 
defines a smooth curve $\psi \: [0,1] \to \Diff(E)$. It is a smooth 
homotopy from $\phi$ to a linear map. 
\end{exer}
\end{small}

\section{Smooth maps into Lie groups} \mlabel{sec:e.3}

In this short section, we describe the bridge between the general setting 
for smooth maps $N \to \Diff(M)$ and smooth maps into Lie groups that we have 
already discussed in Chapter~\ref{ch:4}. 
The main point is that each Lie group $G$ embeds naturally into $\Diff(G)$ 
via $g \mapsto \lambda_g$. Its image is the centralizer of the right 
multiplications $\rho_g$, $g \in G$.

\begin{ex} \mlabel{ex:e.3.1} 
Let $G$ be a Lie group with Lie algebra $\g$. 
Then we embed $G$ into
$\Diff(G)$ via $g \mapsto \lambda_g$, where $\lambda_g(x) = gx$ is the left multiplication. 
If $J \subeq \R$ is an interval and 
$u \: J \to \g$ is a smooth function, then it also defines
a smooth function $\tilde u \: J \to \cV(G)$ with 
$\tilde u(t)(g) :=
 u(t).g$, whose values consist of right invariant
vector fields. Each smooth solution $\gamma \: J \to G$ of the initial
value problem 
\begin{eqnarray}
  \label{eq:e.4}
\gamma'(t) = u(t).\gamma(t), \quad \gamma(t_0) = \gamma_0
\end{eqnarray}
defines a smooth solution 
$\eta \: J \to \Diff(G), t \mapsto \lambda_{\gamma(t)}$ of the
initial value problem 
\begin{eqnarray}
  \label{eq:e.5}
\eta'(t) = \tilde u(t).\eta(t), \quad \eta(t_0) =
\lambda_{\gamma_0}. 
\end{eqnarray}
Therefore Proposition \ref{prop:e.1.4} implies uniqueness of~solutions. 

Suppose, conversely, that $\eta \: J \to \Diff(G)$ satisfies (\ref{eq:e.5}). 
Let $g \in G$ and write $\rho_g(x) := xg$ for the corresponding right multiplication map. 
Then $\eta(t_0) = \lambda_{\gamma_0}$ commutes with every~$\rho_g$ and each vector field 
$\tilde u(t)$ is right invariant, so that Proposition~\ref{prop:e.1.8} implies 
that each $\eta(t)$ commutes with every $\rho_g$. Hence 
$\eta(t) = \lambda_{\gamma(t)}$ holds for the map 
$\gamma \: J \to G, t \mapsto \eta(g)(\be)$ satisfying (\ref{eq:e.4}).
This proves that (\ref{eq:e.4}) has a solution in $G$ if and only if the
corresponding problem (\ref{eq:e.5}) has a solution in $\Diff(G)$. 
\end{ex}

\begin{defn} \mlabel{def:exp-lie} 
Let $G$ be a Lie group with Lie algebra $\g$. We write 
$gv$ for the action of $g \in G$ on a tangent vector $v \in TG$. 
We define ${\cal D}_{\exp} := {\cal D}_{\exp}(G) \subeq \g$ as the set of all
elements $x \in \Lie(G)$ for which the initial value problem 
\begin{eqnarray}
  \label{eq:exp-ivp}
 \gamma'(t) = \gamma(t)x, \quad \gamma(0) = \be 
\end{eqnarray}
has a smooth solution on $\R$. According to the Uniqueness Lemma~\ref{lem:e.1.3}, 
the solution is unique whenever it exists. We write 
$\gamma_x \: \R \to G$
for the corresponding integral curve. It follows from 
(\ref{eq:exp-ivp}) that the curve $\gamma_x$ is smooth. The map 
$$ \exp_G \: {\cal D}_{\exp}(G) \to G, \quad x \mapsto \gamma_x(1) $$
is called the {\it exponential function of $G$}. 
\end{defn}

\begin{lem} \mlabel{lem:exp-lie} The domain ${\cal D}_{\exp}(G)$ of the exponential function of $G$ 
has the following properties: 
  \begin{description} 
  \item[\rm(1)] $0 \in \cD_{\exp}(G)$. 
  \item[\rm(2)] $\R \cdot \cD_{\exp}(G) = \cD_{\exp}(G)$. 
  \item[\rm(3)] For $x \in \cD_{\exp}(G)$ we have 
    $\exp_G((t+s)x ) = \exp_G(tx) \exp_G(sx)$ for $s,t \in \R,$ i.e., 
$\gamma_x \: \R \to G$ is a homomorphism of Lie groups. 
  \item[\rm(4)] $\Ad(g) \cD_{\exp}(G) = \cD_{\exp}(G)$ for all
    $g \in G$ and 
    $$\exp_G \circ \Ad(g)  = c_g \circ \exp_G, 
    \quad \hbox{ where } \quad c_g(h) = ghg^{-1}. $$
  \item[\rm(5)] If $x, y \in \cD_{\exp}(G)$ commute, then 
    \begin{description}
    \item[\rm(a)] $x + y \in \cD_{\exp}(G)$. 
    \item[\rm(b)] $\exp_G(x + y) = \exp_G(x)\exp_G(y) = \exp_G(y)\exp_G(x)$. 
    \item[\rm(c)] $\Ad(\exp_G(x))y = y$. 
    \end{description}  
  \item[\rm(6)] If $x \in \cD_{\exp}(G)$, then $\ad x \in {\cal D}_{\exp}(\GL(\g))$. 
\end{description}
\end{lem}

\begin{prf} In view of Example~\ref{ex:e.3.1}, 
(1)-(5) follow immediately from the corresponding assertions on $\Diff(G)$ in 
Lemma~\ref{lem:e.1.16}. 

(6) Since the adjoint action of $G$ on $\g$ is smooth, 
$\alpha(t) := \Ad(\exp_G(tx))$ defines a smooth $\R$-action on $\g$. 
To calculate its infinitesimal generator, we translate 
the action on $\g$ to the space of left invariant vector fields. 
From $\alpha(t) = \Lie(c_{\exp_G(tx)})$ and 
Exercise~\ref{exer:e3.2}, we obtain for $y \in \g$ 
\begin{eqnarray*}
(\alpha(t)y)_l 
&=& \Ad(c_{\exp_G(tx)})y_l 
= \Ad(\rho_{\exp_G(-tx)})\Ad(\lambda_{\exp_G(tx)})y_l \\
&=& \Ad(\rho_{\exp_G(-tx)})y_l. 
\end{eqnarray*}
As $\rho_{\exp_G(-tx)}$ describes the flow on $G$ generated by the 
left invariant vector field $-x_l$ (Exercise~\ref{exer:e3.1}), 
we can use Lemma~\ref{lem:e.1.6} to obtain 
$$ \alpha'(0)y = [y_l, - x_l](\be) = [x,y], $$
i.e., $\alpha'(0) = \ad x$. 
\end{prf}

\begin{lem} \mlabel{lem:hom-lie} If $\phi \: G \to H$ is a morphism of Lie groups, then 
$$ \Lie(\phi)(\cD_{\exp}(G)) \subeq \cD_{\exp}(H) 
\quad \mbox{ and } \quad 
\phi \circ \exp_G = \exp_H \circ \Lie(\phi). $$
\end{lem}

\begin{prf} If $x \in \Lie(G)$, then $\phi \circ \gamma_x \: \R \to H$ is a smooth 
one-parameter subgroup with $(\phi \circ \gamma_x)'(0) = \Lie(\phi)x$. Thus 
\[ \Lie(\phi)x  \in \cD_{\exp}(H) \quad \mbox{ and } \quad
  \exp_H \circ \Lie(\phi) = \phi \circ \exp_G.\qedhere\]   
\end{prf}

\begin{lem} \mlabel{lem:cov-exp} If $G$ is a connected Lie group with 
$\L(G) = \g$ and $q_G \: \tilde G \to G$ is a covering morphism with 
$\Lie(q_G) = \id_{\g}$, then 
$$ \cD_{\exp}(G) = \cD_{\exp}(\tilde G) 
\quad \mbox{ and } \quad \exp_{\tilde G} = \exp_G \circ q_G. $$
If $\cD_{\exp}(G) = \g$ and $\exp_G \: \g \to G$ is smooth, then the same holds for 
the exponential function $\exp_{\tilde G} \: \g \to \tilde G$. 
\end{lem}

\begin{prf} Lemma \ref{lem:hom-lie} implies that 
$\cD_{\exp}(\tilde G) \subeq \cD_{\exp}(G)$ and 
$\exp_G = q_G \circ \exp_{\tilde G}$. 

To see that $\cD_{\exp}(G) \subeq \cD_{\exp}(\tilde G)$, let 
$x \in \g$ and $\gamma_x \: \R \to G$ be a smooth one-parameter subgroup 
with $\gamma_x'(0) = x$. As $\R$ is simply connected, there exists a 
unique continuous lift $\tilde \gamma_x \: \R \to \tilde G$ with 
$\tilde\gamma_x(0) = \be$ and $q_G \circ \tilde\gamma_x = \gamma_x$. 
Since $q_G$ is a local diffeomorphism, the curve $\tilde\gamma_x$ is smooth, and we 
have 
\[  T_{\tilde\gamma_x(t)}(q_G)\tilde\gamma_x'(t) = \gamma_x'(t) = \gamma_x(t)x 
=  T_{\tilde\gamma_x(t)}(q_G)\tilde\gamma_x(t)x, \] 
which leads to 
$\tilde\gamma_x'(t) =\tilde\gamma_x(t)x,$
and hence to $x \in \cD_{\exp}(\tilde G)$. 

Now we assume that $\cD_{\exp}(G) = \g$ and that $\exp_G$ is smooth. 
Let $V \subeq G$ be an open connected $\be$-neighborhood for which 
$q_G^{-1}(V) \cong V \times J$ holds for a discrete set $J$ and let 
$\tilde V \subeq \tilde G$ denote the identity component of $q_G^{-1}(V)$. 
Then $q_G\res_{\tilde V} \: \tilde V \to V$ is a diffeomorphism. 
Now let $U \subeq \g$ be a convex $0$-neighborhood with $\exp_G(U) \subeq V$. 
For $x \in U$ we then have 
$$ \exp_{\tilde G}(tx) \in q_G^{-1}(V) $$
for $0 \leq t \leq 1$ and since the curve $t \mapsto \exp_{\tilde G}(tx)$ 
is continuous, it follows that $\exp_{\tilde G}(x) \in \tilde V$, and hence that 
$\exp_{\tilde G}(U) \subeq \tilde V$. Furthermore 
$$ \exp_{\tilde G}\res_U  = (q_G\res_{\tilde V})^{-1} \circ \exp_G \res_U, $$
which implies that $\exp_{\tilde G}\res_U$ is smooth. Now the smoothness of 
$\exp_{\tilde G}$ follows from  
$\exp_{\tilde G}(x) = \exp_{\tilde G}(\frac{x}{n})^n$
for every $n \in \N$. 
\end{prf}

\begin{small}
\subsection*{Exercises for Section~\ref{sec:e.3}} 

\begin{exer} \mlabel{exer:e3.1} Let $G$ be a Lie group and $x \in \Lie(G)$. Then the corresponding 
left invariant vector field 
$x_l \in {\cal V}(G)$ with $x_l(\be) = x$ is complete if and only if 
$x \in {\cal D}_{\exp}(G)$. If this is the case, then the corresponding smooth 
$\R$-action on $G$ is given by $\sigma(t,g) := g \exp_G(tx) = \rho_{\exp_G(tx)}(g)$. 
\end{exer}

\begin{exer} \mlabel{exer:e3.2} Let $G$ be a Lie group, $\phi \in \Aut(G)$ 
and $x \in \Lie(G)$. Then $\Ad(\phi)x_l = T(\phi) \circ x_l \circ \phi^{-1} = (\Lie(\phi)x)_l$. 
\end{exer}
  
\end{small}

\section{Automatic smoothness of inversions} \mlabel{sec:e.4}

We have defined smooth maps 
$\phi \: N \to \Diff(M)$ as those for which the two maps 
\[ \tilde\phi_j \: N \times M \to M,\qquad 
\tilde\phi_1(n,m) = \phi(n)(m), \quad 
\tilde\phi_2(n,m) = \phi(n)^{-1}(m) \] 
are smooth. In many situations, it is much easier to verify the smoothness of 
$\tilde\phi_1$, whereas a verification of the smoothness of the 
second map is usually harder. In this section, we provide 
several tools to show the smoothness of~$\tilde \phi_2$. 
In particular, we shall see that, in the linear case, i.e., 
if $M= E$ is a locally convex space and $\phi(N) \subeq \GL(E)$, the smoothness of 
$\tilde\phi_1$ and the continuity of $\tilde \phi_2$ already imply 
the smoothness 
of $\tilde\phi_2$. In non-linear settings, it is necessary to assume that 
$\tilde f_2$ is $C^1$ and if $M$ is Banach, no assumption on 
$\tilde f_2$ is needed at all. 

\subsection*{The linear case} 

\begin{thm} \mlabel{thm:autsmooth-lin} Let $E$ be a locally convex space, $N$ be 
a smooth manifold and $k \in \N$. 
Further, let $f \: N \to \GL(E)$ be a map for which 
\[  \tilde f_1 = f^\wedge\: N \times E \to E, \quad (m,v) \mapsto f(m)v \] 
is $C^k$ and 
\[ \tilde f_2  \: N \times E \to E, \quad (m,v) \mapsto f(m)^{-1}v \] 
is continuous. Then the map $\tilde f_2$ also is $C^k$.   
\end{thm}

\begin{prf} We may w.l.o.g.\ assume that 
$N$ is an open subset of a locally convex space $X$. 
We argue by induction. For $k=1$, we have
\begin{eqnarray*}
&&\ \ \  
\tilde f_2((x,v)+s(h_1,h_2))-\tilde f_2(x,v) \\ 
&& =
f(x+sh_1)^{-1}(v+sh_2)- f(x)^{-1}v \\ 
&&=
s f(x+sh_1)^{-1}h_2 + (f(x+sh_1)^{-1} - f(x)^{-1}) v \\ 
&& =s \tilde f_2(x+sh_1, h_2) +f(x)^{-1}\big(f(x)-f(x+sh_1)\big) \tilde f_2(x+sh_1, v)\\ 
&&=s \tilde f_2(x+sh_1, h_2) \\
&& \ \ \ \ +
f(x)^{-1} \Big(\tilde f_1(x,\tilde f_2(x+sh_1,v)-
\tilde f_1\big(x+sh_1,\tilde f_2(x+sh_1, v)\big)\Big) \\ 
&&=s \tilde f_2(x+sh_1, h_2) +
s f(x)^{-1} \tilde f_1^{[1]} \big((x,\tilde f_2(x+sh_1, v)),(h_1,0),s\big) 
\end{eqnarray*}
(see~\eqref{defnU1} for the notation). 
Dividing by $s$, we get 
\begin{eqnarray}
  \label{eq:x.1a}
&& (\tilde f_2)^{[1]} ((x,v),(h_1,h_2),s) \\
&=&
 \tilde f_2(x+sh_1, h_2)+
 f(x)^{-1} \tilde f_1^{[1]}((x,\tilde f_2(x+sh_1,v),(h_1,0),s) \notag \\ 
&=&\tilde f_2(x+sh_1,h_2)
+\tilde f_2\Big(x, \tilde f_1^{[1]}\big((x, \tilde f_2(x+s h_1,v)),
(h_1,0),s\big)\Big),\notag
\end{eqnarray}
which,
according to our assumptions, is a continuous map. 
It follows that $\tilde f_2$ is $C^1$. For $s=0$, we get
\begin{equation}
  \label{eq:x.1}
T_{(x,v)}(\tilde f_2)(h_1,h_2)=
f(x)^{-1} h_2 - f(x)^{-1} T_{(x,f(x)^{-1}v)}(\tilde f_1)(h_1,0).
\end{equation}
Using Equation~(\ref{eq:x.1}) together with the Chain Rule, we
can iterate this argument, and it follows that $\tilde f_2$
is $C^k$ if $\tilde f_1$ is~$C^k$.
\end{prf}

In view of the preceding theorem, the crucial point verifying that a 
map ${f \: N \to \GL(E)}$ is smooth is to verify that the map $\tilde f_2$ is 
continuous, which might also be difficult. The following 
observations describe cases where it can be simplified. 

\begin{prop} \mlabel{prop:autsmooth-ban} Let $M$ be a smooth manifold, 
let $E$ be a Banach space, and let $f \: M \to \GL(E)$ be a map. 
Then the following are equivalent: 
\begin{description}
\item[\rm(1)] $f$ is a smooth map into the Banach--Lie group $\GL(E)$. 
\item[\rm(2)] $\tilde f \: M \times E \to E \times E, (m,v) \mapsto (f(m)v, f(m)^{-1}v)$ is smooth. 
\item[\rm(3)] $\tilde f_1 \: M \times E \to E, (m,v) \mapsto f(m)v$ is smooth. 
\end{description}
\end{prop}

\begin{prf} (1) $\Rarrow$ (2) $\Rarrow$ (3) are trivial. 

(3) $\Rarrow$ (1): We may w.l.o.g.\ assume that $M$ is an open subset of a 
locally convex space $F$. Assume that the map $\tilde f_1$ is smooth. Then 
\cite[Lemma~2.2]{Gl06a} implies that the map 
\[  d\tilde f_1 \: M \times E \to {\cal L}(F \times E, E)_b \]
is smooth, where ${\cal L}(F \times E, E)_b$ is the set of continuous linear maps 
$F \times E \to E$, endowed with the topology of uniform convergence on bounded subsets. 
Since we have a natural embedding ${\cal L}(E) \into {\cal L}(F \times E, E)$, as those 
linear maps vanishing on $F \times \{0\}$, it follows that the map 
\[ f \: M \to {\cal L}(E), \quad f(m)v = d\tilde f_1(m,0)(0,v) \]
is smooth. Since its values lie in the open subset $\GL(E)$, (1) follows. 
\end{prf}

\begin{lem} Suppose that $E = \prolim E_j$ is a projective limit 
of the Banach spaces $E_j$ and that 
$f \: M \to \GL(E)$ is a map with the property that, for each~$j$, 
the operator $f(m)$ induces an invertible continuous linear operator 
$f_j(m)$ on $E_j$. If all the maps $f_j \: M \to \GL(E_j)$ are smooth, 
then $f$ is smooth. 
\end{lem}

\begin{prf} Let $q_j \: E \to E_j$ denote the canonical projection. 
Our assumption implies that the map $q = (q_j)_{j \in J} \: E \to \prod_{j \in J} E_j$
is a topological embedding. Therefore the map $\tilde f$ is smooth if and only 
if all its components 
\begin{eqnarray*} \tilde f_j \: M \times E 
&\to& E_j \times E_j, \\ 
(m,v) &\mapsto& (q_j(f(m)v),q_j(f(m)^{-1}v)) = (f_j(m)q_j(v), f_j(m)^{-1}q_j(v)) 
\end{eqnarray*}
are smooth. These maps factor through the maps 
$$ M \times E_j \to E_j \times E_j, \quad 
(m,v_j) \mapsto (f_j(m)v_j, f_j(m)^{-1}v_j) $$
which are smooth because all the maps $f_j \: M \to \GL(E_j)$ are smooth 
 (Proposition~\ref{prop:autsmooth-ban}). Therefore $\tilde f$ is smooth. 
\end{prf}

\subsection*{Smoothness of the inversion in groups} 

Let $G$ be a group with a manifold structure for which the multiplication 
map $m_G \: G \times G \to G$ is smooth. What can be said about the smoothness 
of the inversion map $\eta_G$? Does it follow from weaker assumptions? 
If $G$ is a Banach manifold, then the smoothness of $\eta_G$ follows from the 
Inverse Function Theorem, applied to the map 
$$ G \times G \to G \times G, \quad (x,y) \mapsto (x,xy). $$

On the other hand, there exist locally convex algebras $\cA$ 
which are fields, 
so that  their unit group $\cA^\times = \cA \setminus \{0\}$ is open, but 
the inversion map is not  continuous. 
Then $\cA^\times$ carries a natural 
manifold structure, the group multiplication on $\cA^\times$ is smooth, but 
the inversion is not continuous. 

Therefore we cannot expect to get the smoothness of $\eta_G$ for free. 

\begin{prop} \mlabel{prop:aut-smooth-inv} 
Let $G$ be a group with a manifold structure for which the multiplication 
map $m_G \: G \times G \to G$ is smooth. If $\eta_G$ is $C^1$, then it is smooth. 
If $G$ is a Banach manifold, then $\eta_G$ is always smooth. 
\end{prop}

  \begin{prf} 
Let $T(m_G) \: TG \times TG \to TG$ denote the multiplication map on the tangent bundle, 
which turns it into a monoid with a smooth multiplication 
(Proposition~\ref{prop:tangentgrp}). 
Since $\eta_G$ is assumed to be $C^1$, we have a continuous tangent map 
$T(\eta_G) \: TG \to TG$, and the usual arguments (Proposition~\ref{prop:tangentgrp}) 
imply that 
$TG$ actually is a group with continuous inversion map $\eta_{TG} = T(\eta_G)$. 
Identifying $G$ with the zero section in $TG$, which is a subgroup, 
we obtain by differentiating the relation $\eta_G(g) g = \be$ the formula 
$$ T_g(\eta_G)v  = -g^{-1} v g^{-1}, $$
which leads to 
$$ T(\eta_G)v = - \eta_G(q_{TG}(v)) \cdot v \cdot \eta_G(q_{TG}(v)), $$
if $q_{TG} \: TG \to G$ denotes the quotient map. 
Hence the smoothness of the multiplication in $TG$ and the assumption that 
$\eta_G$ is $C^1$ imply that $T(\eta_G)$ is $C^1$, i.e., that $\eta_G$ is $C^2$. 
Iterating the argument, it follows that $\eta_G$ is smooth. 

If, in addition, $G$ is a Banach manifold, then the Inverse Function Theorem applies 
to the map 
$$ \phi \: G \times G \to G \times G, \quad (x,y) \mapsto (x,xy) $$
and shows that 
$\phi^{-1} \: G \times G \to G \times G, (x,y) \mapsto (x,xy^{-1})$
is smooth, which implies that $\eta_G$ is smooth. 
\end{prf}

\subsection*{Smoothness of the inversion in diffeomorphism groups} 

After the linear and the group case, we now turn to the general case of 
maps into diffeomorphism groups. 

We start with a lemma generalizing Theorem~\ref{thm:autsmooth-lin} 
to a non-linear setting. 
For two locally convex spaces $E$ and $W$, we write $\Iso(E,W)$ for the set 
of topological isomorphisms $E \to W$. 

\begin{lem} \mlabel{lem:e.4.3} Let $N$ and $M$ be smooth manifolds, 
$g_j \: N \to M$, $j =1,2$, two $C^k$ maps, 
$$\bE_j :=  g_j^*T(M) := \{ (n,v) \in N \times TM \: v \in T_{g_j(n)}(M)\} $$ 
the corresponding pull-backs of the tangent bundle of $TM$ to $N$ 
(which are $C^k$-manifolds), 
$f(n) \in \Iso(T_{g_1(n)}(M), T_{g_2(n)}(M))$ such that the map 
$$ F_1 \: \bE_1 \to TM, \quad (\bE_1)_n \ni v \mapsto f(n)(v) $$
is $C^k$ and the map 
\[ F_2 \: \bE_2 \to TM, \quad (\bE_2)_n \ni v \mapsto f(n)^{-1}(v) \] 
is continuous. Then $F_2$ is $C^k$. 
\end{lem} 

\begin{prf}
Let $X$ be the model space of $N$ and $Y$ the model space of $M$. 
Pick $n_0 \in N$. Replacing $N$ by a small open neighborhood of $n_0$ and 
$M$ by the union of open neighborhoods of $g_1(n_0)$ and $g_2(n_0)$, we may 
w.l.o.g.\ assume that $N$ is an open subset of $X$ and $M$ an open subset of $Y$. 

Then we have natural isomorphisms $\bE_j \cong N \times Y$, and under this 
identification, $f$ is a map $N \to \GL(Y)$, and the maps $F_j$ are given by 
$$ F_1(n,v) = f(n)v \quad \hbox{ and } \quad 
F_2(n,v) = f(n)^{-1}v. $$
Hence Theorem~\ref{thm:autsmooth-lin} shows that $F_2$ is $C^k$ if $F_1$ is $C^k$ and $F_2$ is 
continuous. 
\end{prf}

\begin{thm} \mlabel{thm:aut-smooth} Let 
$M$ and $N$ be manifolds and $f \: N \to \Diff(M)$ a map for which 
$\tilde f_1$ is smooth and $\tilde f_2$ is $C^1$. Then $\tilde f_2$ is smooth.
\end{thm}

\begin{prf} Our assumption implies that the tangent map 
$$T(\tilde f_2) \: TN \times TM \to TM $$
exists and is continuous. The major point of the proof is to show that, for 
$k \in \N_0$, the map $T(\tilde f_2)$ is $C^{k+1}$ if it is $C^k$. By assumption, it is 
$C^0$ and $\tilde f_2$ is $C^1$, so that we may argue by induction. We assume that 
$\tilde f_2$ is $C^{k+1}$, so that $T(\tilde f_2)$ is $C^k$. 

In view of 
\[ T_{(n,m)}(\tilde f_2)(v,w) = 
T_{(n,m)}(\tilde f_2)(v,0) + T_{(n,m)}(\tilde f_2)(0,w), \] 
$T(\tilde f_2)$ is the fiberwise sum of two maps which can be treated separately. 

We start with the second summand. 
\begin{eqnarray}
  \label{eq:e.x}
   T_{(n,m)}(\tilde f_2)(0,w) = T_m(f(n)^{-1})w = T_{\tilde f_2(n,m)}(f(n))^{-1}w. 
\end{eqnarray}
With the two $C^{k+1}$-maps 
\[  g_1 := \tilde f_2 \: N \times M \to M \quad \mbox{ and } \quad 
g_2 \:  N \times M \to M, \ (n,m) \mapsto m, \] 
we obtain, for each pair $(n,m)$, a topological isomorphism 
\[ F(n,m) := T_{\tilde f_2(n,m)}(f(n)) \: T_{g_1(n,m)}(M) \to T_{g_2(n,m)}(M) 
= T_m(M),\] 
and our assumption 
implies that the corresponding map 
$$ g_1^*T(M) \to T(M), \quad (n,m,v) \mapsto F(n,m)v $$ 
is $C^{k+1}$. Therefore the preceding Lemma~\ref{lem:e.4.3} implies that the map 
$$ g_2^*T(M) \to T(M), \quad (n,m,v) \mapsto F(n,m)^{-1}v$$ 
is $C^{k+1}$, which implies that the map in (\ref{eq:e.x}) is $C^{k+1}$. 

Now we turn to the first summand. 
$$ T_{(n,m)}(\tilde f_2)(v,0) = \big(T_n(f^{-1})v\big)(m), $$
where we identify $T(\Diff(M))$ with a subset of $C^\infty(M,TM)$. 
Differentiating the relation $f \circ f^{-1} = \id_M$ pointwise leads to 
$$ T_n(f^{-1})v = -T(f(n)^{-1}) \circ T_n(f)(v) \circ f(n)^{-1} \in C^\infty(M,TM) $$
(Quotient Rule). 
From the preceding paragraph, we know already that the map 
$(n,m,w) \mapsto T_m(f(n)^{-1})w$ is $C^{k+1}$, so that the assertion follows from the fact that 
the map 
\begin{eqnarray*}
TN \times M &\to& TM,\\ 
T_p(N) \times M \ni  (v,m) 
&\mapsto& T_n(f)(v)(f(n)^{-1}(m)) = T_n(f)(v)(\tilde f_2(n,m)) \\
&& = T_{(n,\tilde f_2(n,m))}(\tilde f_1)(v,0) 
\end{eqnarray*}
is $C^{k+1}$, which follows from our assumptions and the Chain Rule. 

We conclude that $T(\tilde f_2)$ is $C^k$ for each $k$, and this means that 
$\tilde f_2$ is smooth. 
\end{prf}

\begin{thm} \mlabel{thm:autsmooth-bandiff} Let 
$M$ be a Banach manifold, $N$ a locally convex manifold, 
and $f \: N \to \Diff(M)$ be a map for which 
$\tilde f_1$ is smooth. Then $\tilde f_2$ is smooth as well.
\end{thm}

\begin{prf} We may w.l.o.g.\ assume that 
$N$ is an open subset of some locally convex space $E$. 
Let $F$ denote the model space of $M$, which is a Banach space. 

We fix $p_0 \in N$, $m_0 \in M$, and 
open neighborhoods $U$ of $p_0$, $V$ of $m_0$ and $V'$ of $m_0' := f(p_0)(m_0)$ 
such that 
$\tilde f_1(U \times V) \subeq V'$ and there exist charts 
$(\phi,V)$ and $(\phi',V')$ of $M$. 

Now 
$$ \Phi \:= 
\phi' \circ \tilde f_1 \circ (\id_U \times \phi^{-1}) \: U \times \phi(V) \to \phi'(V') $$
is a smooth map from an open subset of $E \times F$ to $F$, 
and the differential of the 
map $\Phi_{p_0} \: \phi(V) \to \phi'(V'), x \mapsto \Phi(p_0,x)$ 
is invertible in each $x \in \phi(V)$ because $\Phi_{p_0} = \phi' \circ f(p_0) 
\circ \phi^{-1}$ is a diffeomorphism onto an open subset of~$F$. 
In view of \cite[Th.~2.3(i)(c)]{Gl04b}, there exist open neighborhoods 
$U_1$ of $p_0$ and $V_1$ of $m_0$ such that 
$W := \bigcup_{p_1 \in U_1} \{p_1\} \times \Phi_{p_1}(\phi(V_1))$ is an open subset 
of $E \times F$, the map 
$U_1 \times \phi(V_1) \to W, (p_1,x) \mapsto (p_1, \Phi_{p_1}(x))$ 
is a diffeomorphism, and the map 
$$W \to \phi(V_1), \quad (p_1,v) \mapsto \Phi_{p_1}^{-1}(v) $$
is smooth. Hence the map 
$\tilde f_2(p,x) = \phi^{-1}(\Phi_p^{-1}(\phi'(x)))$ 
is smooth in a neighborhood of $(p_0, m_0')$, and since 
$p_0$ and $m_0$ are arbitrary, $\tilde f_2$ is smooth. 
\end{prf}

\begin{cor} \mlabel{cor:d.4.9} Let $M$ be a Banach manifold, 
$G$ a Lie group and let $\sigma \: G \to \Diff(M)$ be a homomorphism.  
Then the map $\sigma$ is smooth if and only if the corresponding 
action of $G$ on $M$ is smooth. 
\end{cor}

\begin{prf} Since the action is the map $\tilde\sigma_1$, 
this follows from Theorem~\ref{thm:autsmooth-bandiff}. 
\end{prf}

We conclude this section with the observation that the general smoothness 
concept discussed here coincides with the natural one in all cases where 
the target group has a natural Lie group structure. 

\begin{cor} Let $M$ be a compact smooth manifold  
and $N$ be any smooth manifold. 
Then a map $f \: N \to \Diff(M)$ 
is smooth with respect to the Lie group structure on $\Diff(M)$ if and only 
if it is smooth in the sense defined above, resp., if and only if 
$\tilde f_1$ is smooth. 
\end{cor} 

\begin{prf} For the Lie group structure on $\Diff(M)$, the 
evaluation map 
\[ \ev \: \Diff(M) \times M \to M, \quad (\phi, m) \mapsto \phi(m)  \] 
is smooth by Theorem~\ref{thmmfdmps}(c) because 
$\Diff(M)$ is an open submanifold of $C^\infty(M,M)$. 
If $f$ is smooth with respect to the Lie group structure on $\Diff(M)$, then 
the smoothness of the evaluation map implies that $\tilde f_1$ is smooth, 
and Theorem~\ref{thm:autsmooth-bandiff} implies that $\tilde f$ is smooth. 

If, conversely, $\tilde f_1 = f^\wedge$ is a smooth map, then 
Theorem~\ref{thmmfdmps}(a) implies 
that $f$ is smooth as a map into the Lie group $\Diff(M) 
\subeq C^\infty(M,M)$. 
\end{prf}

From Proposition~\ref{prop:autsmooth-ban} we derive: 
\begin{prop} 
If $E$ is a Banach space, then $\GL(E)$ carries a natural Lie 
group structure and a map $f \: E \to \GL(E)$ is smooth as a map into the Lie group 
$\GL(E)$ if and only if $\tilde f_1$ is smooth. 
\end{prop}

\begin{prf} Since $E$ is a Banach space, the group $\GL(E)$ carries a natural Lie 
group structure, and it follows from Proposition~\ref{prop:autsmooth-ban} 
that~$f$ is smooth as a map into the Lie group 
$\GL(E)$ if and only if $\tilde f_1$ is smooth. 
\end{prf}

\begin{prop}
If $G$ is a Lie group, then $f \: M \to G$ is smooth 
if and only if the map 
$\phi \: M \to \Diff(G), m \mapsto \lambda_{f(m)}$ is smooth.
\end{prop}

\begin{prf} The map 
$\phi \: M \to \Diff(G), m \mapsto \lambda_{f(m)}$ is smooth if and only if 
\[  \tilde\phi \: M \times G \to G \times G, 
\quad (m,g) \mapsto (f(m)g, f(m)^{-1}g) \] 
is smooth. As $G$ is a Lie group, this is equivalent to the smoothness 
of~$f$.
\end{prf}

\subsection*{Open problems} 

The preceding results solve 
a special case of the following problem whose solution 
would be most valuable to verify smoothness of maps into diffeomorphism groups. 

\begin{probl} Let $N$ and $M$ be smooth manifolds. 
Further let \break 
 $\phi \: N \to \Diff(M)$ be a map for which 
$\tilde\phi_1$ is smooth and 
$\tilde\phi_2$ is continuous. 
Does this imply that $\tilde\phi_2$ is smooth? 
\end{probl}

\begin{probl} Let $G$ be a smooth manifold with a group structure for which 
the multiplication map $m_G$ is smooth. 
Does the continuity of the inversion $\eta_G$ imply that it is $C^1$ and therefore smooth?   
If $G$ is Banach, the smoothness of $\eta_G$ follows from the Inverse Function Theorem 
without assuming its continuity. Therefore it would be interesting to see 
if there is some Fr\'echet group with smooth multiplication 
for which the inversion is continuous but not smooth. Note that this pathology 
does not occur for unit groups of cias, where the continuity of the inversion 
implies its smoothness (Corollary~\ref{invsmoocia}). 
\end{probl}

\section{Locally exponential Lie algebras of vector fields} 

R.~Palais' Integrability Theorem \cite{Pa57} asserts that 
any fi\-nite-di\-men\-sio\-nal Lie algebra $\g$ of complete vector fields on a 
finite-dimensional smooth manifold $M$ can be integrated to a global action of a Lie group $G$ on $M$. 
The purpose of this section is to extend this result to 
locally exponential Lie algebras of complete vector fields 
on locally convex manifolds (cf.\ \cite{AbNe08}). 

Let $M$ be a locally convex manifold, 
$\g$ be a locally exponential Lie algebra and 
$\alpha \: \g \to {\cal V}(M)$ a morphism of Lie algebras. We shall work with the 
following assumptions: 
\begin{description}
\item[(A1)] Each vector field $\alpha(x)$, $x \in \g$, is complete. 
\item[(A2)] $\Exp_\alpha \: \g \to \Diff(M), x \mapsto \Exp(-\alpha(x))$ is smooth. 
\end{description}

Note that (A2) implies in particular that the linear 
map $-T_0(\Exp_\alpha) = \alpha \: \g \to {\cal V}(M)$ is continuous. 

If $M$ is infinite-dimensional, then we cannot expect the action of 
$\Diff(M)$ on ${\cal V}(M)$ to be smooth in any sense. Even worse, 
we do not know if, for a smooth curve 
$\phi \: \R \to \Diff(M)$, the corresponding map 
\[ \R \times {\cal V}(M) \to {\cal V}(M),\quad (t,X) \mapsto \Ad(\phi(t))X \] 
is smooth. Hence one-parameter groups of diffeomorphisms of infinite-dimensional 
manifolds have to be handled with particular care. However, in our context, the 
following lemma provides what we need.

\begin{lem} \label{lem:ad-exp-vect} For each $x \in \g$ we have in 
${\cal L}(\g,{\cal V}(M))$: 
  \begin{eqnarray}
    \label{eq:int.1}
\Ad(\Exp_\alpha(x)) \circ \alpha = \alpha \circ e^{\ad x}. 
  \end{eqnarray}
\end{lem}

\begin{prf} We recall that 
$\phi(t) := \exp(t\alpha(x))= \Exp_\alpha(-tx)$ is a smooth curve in $\Diff(M)$ 
and $\beta(t) := \alpha(e^{t\ad x}y)$ is a smooth curve in ${\cal V}(M)$ 
by (A2), so 
that $\psi(t) := \Ad(\phi(t))\beta(t)$ is smooth with 
\begin{eqnarray*}
\psi'(t) 
&=& \Ad(\phi(t))([\beta(t), \delta(\phi)_t] + \beta'(t))
= \Ad(\phi(t))\big([\beta(t), \alpha(x)] + \beta'(t)\big)
\end{eqnarray*}
(Lemma~\ref{lem:e.1.6}). 
Since our assumptions imply that $\alpha$ is a continuous linear map, 
the Chain Rule leads to 
\begin{eqnarray*}
\beta'(t) &=& \alpha([x,e^{t\ad x}y]) 
= [\alpha(x), \alpha(e^{t\ad x}y)]
= [\alpha(x), \beta(t)], 
\end{eqnarray*}
showing that $\psi'=0$, and hence that $\psi$ is constant, which leads to 
$$\Ad(\Exp_\alpha(-x))\alpha(e^{\ad x}y) 
=\Ad(\exp(\alpha(x))) \alpha(e^{\ad x}y) = \psi(1) = \psi(0) = \alpha(y), $$
so that (\ref{eq:int.1}) holds.
\end{prf}

For the following proposition, we recall from Section~\ref{sec:5.2} the concept 
of an exponential local Lie group. 

\begin{prop} \label{prop-mult} Let $G \subeq \g$ be an exponential local Lie group 
with Lie algebra $\g$ and $D_G^0 \subeq D_G$ the connected component of 
$(0,0)$. Then we have for all $(x,y) \in D_G^0$: 
$$ \Exp_\alpha(x*y) = \Exp_\alpha(y) \circ \Exp_\alpha(x). $$
\end{prop}

\begin{prf} We have two smooth functions 
  \begin{eqnarray*}
f_1 &\:& D_G \to \Diff(M), (x,y) \mapsto \Exp_\alpha(x*y), \\
f_2 &\:& D_G \to \Diff(M), (x,y) \mapsto \Exp_\alpha(x)\Exp_\alpha(y) 
\end{eqnarray*}
satisfying $f_1(0,0) = \id_M = f_2(0,0)$. 
Since $D_G^0$ is connected, in view of the Uniqueness Lemma, it suffices to 
show that $\delta^l(f_1) = \delta^l(f_2)$. We use Lemmas~\ref{lem:dexp-vf} 
and \ref{lem:ad-exp-vect} to get 
\begin{eqnarray*}
\delta^l(\Exp_\alpha)_x y  
&=& \int_0^1 \Ad(\Exp_\alpha(-tx))(-\alpha(y))\, dt 
= -\alpha\Big(\int_0^1 e^{-t\ad x}y\, dt\Big) \\
&=& -\alpha \circ (\kappa_G)_x(y), 
\end{eqnarray*}
i.e., $\delta^l(\Exp_\alpha) = -\alpha \circ \kappa_G$. 
This leads to 
\begin{eqnarray*}
\delta^l(f_1) 
&=& \delta^l(\Exp_\alpha \circ m_G) 
= m_G^*\delta^l(\Exp_\alpha) = -m_G^*(\alpha \circ \kappa_G) = -\alpha \circ (m_G^*\kappa_G)\\ 
&=&-\alpha \circ \delta^l(m_G). 
\end{eqnarray*}
Let $p_{1/2} \: D_G \to \g$ denote the projections 
$p_1(x,y) = x$ and $p_2(x,y) = y$. Then 
$f_2 = (\Exp_\alpha \circ p_1)(\Exp_\alpha \circ p_2)$, so that we 
can use the Product Rule (Lemma~\ref{lem:e.1.2}), 
Lemma~\ref{lem:ad-exp-vect}  and Remark~\ref{rem:5.2.11} 
to obtain 
\begin{eqnarray*}
\delta^l(f_2) 
&=& \delta^l(\Exp_\alpha \circ p_2) + \Ad(\Exp_\alpha \circ p_2)^{-1}.\delta^l(\Exp_\alpha \circ p_1)\\ 
&=& p_2^*\delta^l(\Exp_\alpha) + \Ad(\Exp_\alpha \circ p_2)^{-1}.p_1^*\delta^l(\Exp_\alpha) \\ 
&=& -p_2^*(\alpha \circ \kappa_G)- \Ad(\Exp_\alpha \circ p_2)^{-1}.p_1^*(\alpha \circ \kappa_G) \\
&=& -\alpha \circ \delta^l(p_2) - \Ad(\Exp_\alpha \circ p_2)^{-1}.(\alpha \circ \delta^l(p_1))\\
&=& -\alpha \circ \Big(\delta^l(p_2) + \Ad_G(p_2)^{-1}.\delta^l(p_1)\Big) 
= -\alpha \circ \delta^l(m_G) = \delta^l(f_1).\qquad \qedhere
\end{eqnarray*}
\end{prf}

\begin{thm} \label{thm:gen-palais} Let $\g$ be a locally exponential Lie algebra, 
$M$ a smooth locally convex manifold  and  
$\alpha \: \g \to {\cal V}(M)$ be an injective morphism of Lie algebras satisfying 
\begin{description}
\item[\rm(A1)] Each vector field $\alpha(x)$, $x \in \g$, is complete. 
\item[\rm(A2)] $\Exp_\alpha \: \g \to \Diff(M), x \mapsto \exp(-\alpha(x))$ is smooth. 
\item[\rm(A3)] The subgroup $\Gamma_\alpha := \Exp_\alpha^{-1}(\id_M) \cap \z(\g)$ is discrete. 
\end{description}
Then there exists a connected locally exponential Lie group $G$ and a 
map $\beta \: G \to \Diff(M)$ defining an effective smooth $G$-action on $M$ 
with $\dot\beta = \alpha$. 
\end{thm}

\begin{prf} We consider the subgroup $G := \la \Exp_\alpha(\g) \ra  \subeq \Diff(M).$
Let $G_{\rm loc} \subeq \g$ be an exponential local Lie group with Lie algebra $\g$. 
If $x \in G_{\rm loc}$ satisfies $\Exp_\alpha(x) = \id_M$, then Lemma~\ref{lem:ad-exp-vect} 
implies that $\alpha \circ e^{\ad x} = \alpha$, and hence that 
$e^{\ad x} = \id_\g$ because $\alpha$ is injective. 

From the Adjoint Integrability Theorem~\ref{thm:5.3.8} we know that 
$\z(\g)$ is open in $\{ x \in \g \: e^{\ad x} = \id_\g\}$, so that 
we may shrink $G_{\rm loc}$ in such a way that, for 
$x \in G_{\rm loc}$, the relation $e^{\ad x} = \id_\g$ implies $x \in \z(\g)$. 
In view of the discreteness of $\Gamma_\alpha$, we may further shrink 
$G_{\rm loc}$ so that 
$$ \Exp_\alpha^{-1}(\id_M) \cap G_{\rm loc} = \{0\}. $$
Then Proposition~\ref{prop-mult} implies the existence of an open balanced 
$0$-neigh\-bor\-hood $V$ in $G_{\rm loc}$ 
with $V \times V \subeq D_{G_{\rm loc}}$ for which  
$$ \Exp_\alpha(x*y) = \Exp_\alpha(x) \Exp_\alpha(y) \quad \mbox{ for all } \quad x,y \in V.  $$

To see that $\Exp_\alpha\res_V$ is injective, suppose that 
$\Exp_\alpha(x) = \Exp_\alpha(y)$ for $x,y \in V$. Then 
$$\Exp_\alpha((-x)*y) = \Exp_\alpha(-x)\Exp_\alpha(y) = \Exp_\alpha(x)^{-1}\Exp_\alpha(y) = \id_M $$ 
leads to 
$(-x)*y \in \Exp_\alpha^{-1}(\id_M) \cap G_{\rm loc} = \{0\}$. 
We conclude that 
$$x = x * 0 = x * ((-x) * y) = (x * (-x)) * y = 0 * y = y $$
and hence that $\Exp_\alpha\res_V$ is injective. 

This means that $\Exp_\alpha(V) \subeq G$ carries a manifold structure 
for which $\Exp_\alpha\res_V$ is a diffeomorphism, 
and that multiplication and inversion are smooth in an identity neighborhood in $V$. 
Since $G$ is generated by the image of every $0$-neighborhood in $V$, 
it carries a unique Lie group 
structure  for which $\Exp_\alpha$ defines a chart in a neighborhood of $0$ 
(Theorem~\ref{thm:locglob}). 

The curves $\gamma_x \: \R \to G$ defined by 
$\gamma_x(t) := \Exp_\alpha(tx)$ define smooth one-parameter groups 
with $\gamma_x'(0) = x$, where we identify $\g$ via $T_0(\Exp_\alpha)$ with $T_\be(G)$. 
Therefore $\Exp_\alpha$ is an exponential function of $G$. 
In particular,  $G$ is locally exponential  
and the inclusion map $\beta \: G \to \Diff(M)$ is smooth in 
an identity neighborhood, hence everywhere, and thus it 
defines a smooth action of $G$ on $M$ 
with $\dot\beta = -T_0(\Exp_\alpha) = \alpha$. 
\end{prf}

\begin{cor} Let $G$ be a $1$-connected locally exponential Lie group, 
$M$ a smooth locally convex manifold  and  
$\alpha \: \Lie(G) \to {\cal V}(M)$ be 
a morphism of Lie algebras satisfying {\rm(A1)-(A3)}. 
Then there exists a unique map $\beta \: G \to \Diff(M)$ defining an 
effective smooth $G$-action on $M$ with $\dot\beta = \alpha$. 
\end{cor}

Fortunately, the conditions (A2) and (A3) are automatically satisfied 
if $M$ is finite-dimensional connected and $\sigma$-compact 
(which, in this case, is equivalent 
to paracompactness). This leads to the following easier integrability 
criterion: 

\begin{thm} \mlabel{thm:findim-discrete} 
Suppose that $M$ is finite-dimensional and $\sigma$-compact 
with finitely many connected components 
and that $\alpha \: \g \to {\cal V}(M)$ is a continuous homomorphism 
of Lie algebras, where $\g$ is locally exponential and 
$\im(\alpha)$ consists of complete vector fields. 
Then there exists a connected locally exponential Lie group $G$ and a 
map $\beta \: G \to \Diff(M)$ defining an effective smooth $G$-action on $M$ 
with $\dot\beta = \alpha$. 
\end{thm}

\begin{prf} In view of Theorem~\ref{thm:gen-palais}, we have to verify 
(A2) and (A3). 

(A2): Let $E := \im(\alpha)$, endowed with the subspace topology, 
inherited from ${\cal V}(M)$. Then $\alpha \: \g \to E$ is continuous linear, 
hence smooth, so that it suffices to show that the map 
$\Phi \: E \times M \to M, (X,p) \mapsto \exp(X)(p)$
is smooth. 

Since $M$ is $\sigma$-compact, ${\cal V}(M)$ is a Fr\'echet space. 
It follows in particular, that 
$E$ is a metrizable locally convex space. In view of \cite[Thm.~12.4]{BGN04}, 
it suffices to show that, for each smooth map 
$f \: \R^k \to E$, the map 
$$ \Phi  \circ (f \times \id_M) \: \R^k \times M \to M $$
is smooth. 

For any such $f$, the smoothness of the evaluation map ${\cal V}(M) \times M \to T(M)$ 
implies that the map $\R^k \times M \to T(M), (x,m) \mapsto f(x)(m)$ is smooth. 
We thus obtain on the manifold $\hat M := \R^k \times M$ a smooth  vector field 
$X(x,m) := (0, f(x)(m))$ with 
\[ \Phi^X_1(x,m) = (x,\Phi^{f(x)}_1(m)) = \big(x, \exp(f(x))(m)\big), \]
so that the assertion follows from the smoothness of the map 
$\Phi^X_1$ on $\R^k \times M$ (\cite[Thm.~IV.2.6]{La99}).   

(A3): Suppose that $\Gamma_\alpha$ is not discrete, i.e., that $0$ is 
not isolated in $\Gamma_\alpha$. Then $\alpha(\Gamma_\alpha)$ is not discrete, 
and since ${\cal V}(M)$ is metrizable, there exists a sequence 
$(x_n)_{n \in \N}$ in $\g$ with $\alpha(x_n) \to 0$ and 
$\exp(\alpha(x_n)) = \id_M$ for each $n \in \N$. 

For each each finite subset $F \subeq \N$, the smooth vector fields 
$\alpha(x_n)$, $n \in F$, commute and define smooth actions of the circle 
group $\T = \R/\Z$ on $M$, so that they combine to a smooth action of the 
torus $\T^F$ on $M$. 
If $K$ is a compact group acting faithfully on 
$M$, then $M$ possesses a $K$-invariant 
Riemannian metric $g$, which leads for $d := \dim M$ and 
the number $c$ of connected components of $M$ to the estimate 
\[ \dim K \leq \dim(\Iso(M,g)) \leq c\frac{d(d+1)}{2} \] 
(cf.\ \cite[Thm.~II.1.3]{Ko95}). 

Applying this to the action of the torus $\T^F$ from above, it follows that, 
for each finite subset $F \subeq \N$, we have 
\[  \dim \Spann \{ \alpha(x_n) \: n \in F \} \leq c\frac{d(d+1)}{2}, \] 
but this implies that $\Spann \{ \alpha(x_n) \: n \in \N\}$ is finite-dimensional. 
As $\alpha$ is injective, $\fa := \Spann \{ x_n \: n \in \N\}$ is a 
finite-dimensional central subalgebra of $\g$ and 
$\sigma(x,m) := \Exp(x)(m)$ defines a smooth action of the 
finite-dimensional Lie group $A := (\fa,+)$ on $M$. 
In particular, $\Gamma_\alpha = \{ x \in \fa \: \Exp(x) = \id_M\} = \ker \sigma$ 
is a closed normal subgroup of $A$ whose Lie algebra 
is $\ker \alpha \cap \fa = \{0\}$. Hence $\Gamma_\alpha$ is discrete. 
\end{prf}

\begin{cor} Let $G$ be a $1$-connected locally exponential Lie group, 
$M$ a compact smooth manifold  and  
$\alpha \: \Lie(G) \to {\cal V}(M)$ a continuous morphism of topological Lie algebras. 
Then there exists a unique morphism of Lie groups 
$\beta \: G \to \Diff(M)$ with $\Lie(\beta) = -\alpha$. 
\end{cor}

The preceding corollary can also be obtained by different methods, exploiting the fact that, 
for any compact manifold $M$, the group $\Diff(M)$ is a regular Lie group 
(Chapter~\ref{ch:diffeo}), 
so that the general machinery developed in Chapter~\ref{ch:4} applies. 

\begin{ex} The assumption that the number of connected components 
of $M$ is finite is crucial. In fact, let 
$M := \bigcup_{n \in \N} M_n$ be a disjoint union of infinitely 
many circles $M_n \cong \bS^1$. Identifying $M_n$ with $\R/\Z$, 
we obtain a vector field $X_n$ on $M$ vanishing on $M\setminus M_n$, 
and on $M_n$ it is given by $\frac{\partial}{\partial x}$. 
The Lie algebra $\g := \R^\N$ is abelian, hence locally exponential, 
and we obtain a continuous homomorphism 
$$\alpha \: \g \to {\cal V}(M), \quad 
\alpha((t_n)) := \sum_{n= 1}^\infty t_n X_n $$
because on each connected component, the series on the right consists of 
only one summand.  
Clearly $\exp((t_n)) = \id_M$ is equivalent to $t_n \in \Z$ for 
each $n \in \N$, so that $\Gamma_\alpha= \Z^\N$, which is not discrete. 
\end{ex}

\begin{probl} Find a locally exponential Lie algebra $\g$, a 
smooth connected manifold 
$M$ and a continuous homomorphism 
$\alpha \: \g \to {\cal V}(M)$ for which $\Gamma_\alpha$ is not discrete. 
For that we may w.l.o.g.\ assume that $\g$ is abelian. As we have already seen that 
$\Gamma_\alpha$ is always discrete if $M$ is finite-dimensional 
connected and $\sigma$-compact, for any  such example 
the manifold $M$ must be infinite-dimensional or not paracompact. 
\end{probl}

\begin{probl} Find a Lie group $G$ with exponential function 
and an abelian subalgebra $\fa \subeq \g$ for which 
$\ker(\exp\res_\fa)$ is not discrete. Clearly, $G$ cannot be locally 
exponential. 
\end{probl}

\subsection*{Notes on Appendix~\ref{app:nonlie}} 

In \cite{OMY82}, one finds the particular version of Lemma~\ref{unique-intcurve}
dealing with solutions of the initial value problem 
$$ \eta'(t) = [\eta(t),\xi(t)] + \eta(t), \quad \eta(0) = x $$
in the Lie algebra $\g$ of a regular Lie group $G$ (cf.~Lemma~\ref{lem:e.1.6}), 
see also \cite[2.5/2/6]{KYM85}. 

More examples of linear ODEs with constant coefficients on Fr\'echet spaces for which
solutions are not unique can be found in \cite{Mil84}. 

One-parameter semigroups of operators on locally convex spaces 
have already been studied in the 1960s \cite{Yo63}, 
\cite{Ko64}, \cite{Ko68}. In particular in has been shown in a rather 
general context that the ``infinitesimal generator'' determines 
the one-parameter group uniquely under mild assumptions, 
such as equicontinuity of the semigroup.



\begin{thebibliography}{MMMa99}

\bibitem[AbNe08]{AbNe08} Abouqateb, A., and K.-H. Neeb, {\it Integration of locally exponential 
Lie algebras of vector fields}, 
Ann.\ Global Anal.\ Geom.\ {\bf 33:1} (2008), 89--100. 
%
\bibitem[ARS85]{ARS85} Adams, M., T. Ratiu, and R. Schmid, {\it 
The Lie group structure of diffeomorphism groups and 
invertible Fourier integral operators, with applications}, 
pp.~1--69 in: V. Kac (ed),
``Infinite-dimensional Groups with Applications'' 
(Berkeley, Calif., 1984), Math.\ Sci.\ Res.\ Inst.\ Publ.\ {\bf 4}, 
Springer, New York, 1985.
%
\bibitem[ARS86a]{ARS86a} ---, {\it A Lie group structure for pseudodifferential operators}, 
Math.\ Ann.\ {\bf 273:4} (1986), 529--551.
%
\bibitem[ARS86b]{ARS86b} ---, {\it A Lie group structure for Fourier integral operators}, 
Math.\ Ann.\ {\bf 276:1} (1986), 19--41. 
%
\bibitem[AkN96]{AkN96}
Akkar, M. and C. Nacir, C.,
\emph{Structure $m$-convexe d'une alg\`{e}bre limite inductive localement convexe
d'alg\`{e}bres de Banach},
Rend.\ Sem.\ Mat.\ Univ.\ Padova {\bf 95} (1996), 107--126.
%
%
\bibitem[AHM93]{AHM93}
Albeverio, S., J. H\o{}egh-Krohn, J.~A. Marion,
D.~H. Testard, and B.~S. Torr\'{e}sani, ``Noncommutative Distributions'', 
Marcel-Dekker, New York, 1993. 
%
\bibitem[AHV83]{AHV83}
Albeverio, S., R. Hoegh-Krohn, and A.~Vershik, {\it Factorial representations of path 
groups}, J. Funct.\ Anal.\ \textbf{51}  (1983), 115--131.  
%
\bibitem[Al65]{Al65} Allan, G.~R., {\it A spectral theory for locally convex algebras}, 
Proc.\ London Math.\ Soc.\ (3) {\bf 15} (1965), 399--421. 
%
\bibitem[AAB97]{AAB97} Allison, B. N., S. Azam, S.~Berman, Y.~Gao, and
A.~Pianzola, ``Extended Affine Lie Algebras and Their Root Systems,'' 
Mem.\ Am.\ Math.\ Soc.\ {\bf 603}, Providence R.~I., 1997.
%
\bibitem[ABG00]{ABG00} Allison, B., G. Benkart, and Y. Gao, {\it Central 
extensions of Lie algebras graded by finite-root
systems}, Math.\ Ann.\ \textbf{316} (1993), 499--527.
%
\bibitem[Alz13]{Alz13}
Alzaareer, H., ``Lie Groups of Mappings on Non-Compact Spaces and Manifolds",
doctoral dissertation,
Universit\"{a}t Paderborn, 2013;
{\tt http://nbn-resolving.de/urn:n:de:hbz:466:2-11572}
%
\bibitem[Alz19]{Alz19} ---,
\emph{Differential calculus on multiple products},
Indag.\ Math.\ {\bf 30} (2019), 1036--1060.
%
\bibitem[Alz21]{Alz21} ---,
\emph{Lie groups of $C^k$-maps on non-compact manifolds and the
fundamental theorem for Lie group-valued mappings},
J. Group Theory {\bf 24} (2021), 1099--1134.
%
\bibitem[AlS15]{AlS15}
Alzaareer, H., and A. Schmeding, \emph{Mappings on products with different degrees
of differentiability in the two factors}, Expo.\ Math.\ \textbf{33} (2015), 184--222.
%
\bibitem[Ame75]{Ame75} Amemiya, I., {\it Lie algebra of vector fields and complex structure}, 
J.~Math.\ Soc.\ Japan {\bf 27:4} (1975), 545--549. 
%
\bibitem[AGS20]{AGS20}
Amiri, H., H. Gl\"{o}ckner, and A. Schmeding,
\emph{Lie groupoids of mappings taking values in a Lie groupoid},
Arch.\ Math.\ (Brno) {\bf 56} (2020),
307--356.
%
%
\bibitem[AnN09]{AnN09} An, J., Neeb, K.-H., {\it An implicit function 
theorem for Banach spaces and some applications}, Math. Z. {\bf 262:3} (2009), 
627--643.

\bibitem[AW08]{AW08} An, J., and Z. Wang, {\it Curve selection lemma 
for semianalytic sets and conjugacy classes of finite order in Lie groups}, 
Sci. China Ser. A {\bf 51:3} (2008), 383--388.
%
\bibitem[ADG94]{ADG94} Ancel, F.D., Dobrowolski, T., and J.\ Grabowski, 
{\it Closed subgroups in Banach spaces}, Studia Math. {\bf 109:3} (1994), 
277--290 
%
\bibitem[AB68]{AB68} 
Anderson, R.~D., and R.~H. Bing, {\it A complete elementary
proof that Hilbert space is homeomorphic to the countable infinite
product of lines}, Bull.\ Am.\ Math.\ Soc. {\bf 74} (1968), 771--792. 
%
\bibitem[ADM20]{ADM20} Ando, H., M. Doucha, and Y. Matsuzawa, 
{\it Large scale geometry of Banach--Lie groups}, 
arXiv:math.OA:2011.10376  
%
\bibitem[AW64]{AW64}
Araki, H., and W. Wyss, {\it Representations of canonical
anticommutation relations} Helv.\ Phys.\ Acta \textbf{37} (1964), 136--159. 
%
\bibitem[ASS71]{ASS71} 
Araki, H., M.-S. Bae Smith, and L. Smith, {\it On the
homotopical significance of the type of von Neumann algebra factors},  
Commun.\ Math.\ Phys.\ {\bf 22} (1971), 71--88. 
%
\bibitem[Arn66]{Arn66} Arnold, V.~I., {\it Sur la g\'eom\'etrie diff\'erentielle des groupes 
de Lie de dimension infinie et ses applications \`a l'hydrodynamique des 
fluides parfaits}, Ann.\ Inst.\ Fourier {\bf 16:1} (1966), 319--361.  
%
\bibitem[AK98]{AK98} 
Arnold, V.~I., and B.~A. Khesin, ``Topological Methods in Hydrodynamics'',  
Springer-Verlag, Berlin, 1998.
%
\bibitem[At67]{At67} 
Atiyah, M., ``K-Theory'', Benjamin, Amsterdam, 1969. 
%
\bibitem[Au99]{Au99} 
Au\ss{}enhofer, L., ``Contributions to the Duality Theory of
Abelian Topological Groups and to the Theory of Nuclear Groups'', 
Diss.\ Math.\ {\bf 384}, Warsaw, 1999. 
%
\bibitem[AS68]{AS68}
Averbukh, V.~I., and O.~G. Smolyanov, 
{\it The various definitions of the derivative in linear
topological spaces}, Russ.\ Math.\ Surv.\ \textbf{23} (1968), 67--113. 
%
\bibitem[AI95]{AI95}
de Azcarraga, J.~A., and J.~M. Izquierdo, ``Lie Groups, Lie
Algebras, Cohomology and some Applications in Physics, '' 
Cambridge Monographs on Math.\ Physics,
Cambridge Univ.\ Press, Cambridge, 1995. 
%
\bibitem[BSZ92]{BSZ92}
Baez, J.~C., I. Segal, and Z. Zhou, ``Introduction to Algebraic  
and Constructive Quantum Field Theory'', 
Princeton Univ.\ Press, Princeton, 1992. 
%
\bibitem[Bak01]{Bak01} Baker, H.~F., {\it On the exponential theorem for a simply 
transitive continuous group, and the calculation of the finite equations from 
the constants of structure}, 
J. London Math.\ Soc.\ {\bf 34} (1901), 91--127. 
%
\bibitem[Bak05]{Bak05} ---, {\it On the calculation of the finite equations of a 
continuous group}, London M.~S. Proc.\ {\bf 35} (1903), 332--333. 
%
\bibitem[BD01]{BD01} Balan, V., and J. Dorfmeister, {\it Birkhoff decompositions 
and Iwasawa decompositions for loop groups}, 
Tohoku Math. J. (2) {\bf 53:4} (2001), 593--615
%
\bibitem[BMSY11]{BMSY11} Banakh, T., K. Mine, K. Sakai, and T.~Yagasaki, 
{\it Homeomorphism and diffeomorphism groups of non-compact manifolds 
with the Whitney topology}, Topol.\ Proc.\ {\bf 37} (2011), 61--93.
%
\bibitem[BR10]{BR10}
Banakh, T. and D. Repov\v{s},
\emph{Topological structure of direct limits in the category of uniform spaces},
Topology Appl.\ {\bf 157}, (2010), 1091--1100. 
%
\bibitem[BY15]{BY15} Banakh, T., and T.~Yagasaki, 
{\it Diffeomorphism groups of non-compact manifolds endowed with the 
Whitney $C^\infty$-topology}, 
Topology Appl. {\bf 179} (2015), 51--61. 
%
\bibitem[Bcz91]{Bcz91}
Banaszczyk, W.,
``Additive Subgroups of Topological Vector Spaces,''
Lecture Notes in Math.\ {\bf 1466}, Springer-Verlag, Berlin, 1991.
%
\bibitem[Ban97]{Ban97} Banyaga, A., ``The Structure of Classical Diffeomorphism
Groups,'' Kluwer Academic Publishers, 1997. 
%
\bibitem[Ban02]{Ban02} ---, {\it On the geometry of locally conformal 
symplectic manifolds}, 
in ``Infinite Dimensional Lie Groups in Geometry and Representation 
Theory'', Banyaga, A., J. A. Leslie, and T. Robart, eds., 
World Scientific, River Edge, NJ, 2002; 79--104. 
%
\bibitem[BLR02]{BLR02} Banyaga, A., J. A. Leslie, and T. Robart, eds., 
``Infinite Dimensional Lie Groups in Geometry and Representation 
Theory'', World Scientific, River Edge, NJ, 2002 
%
\bibitem[Ba64]{Ba64}
Bastiani, A., {\it Applications diff\'erentiables et vari\'et\'es diff\'erentiables 
de dimension infinie}, J. Anal.\ Math.\ \textbf{13} (1964), 1--114. 

\bibitem[BHM25]{BHM23} Bauer, M., P.~Harms and P.W.~Michor,
  {\it Regularity and completeness of half-Lie groups},
  J. Europ. Math. Soc., DOI 10.4171/JEMS/1587;
 arXiv:2302.01631 
%

\bibitem[Bg87]{Bg87}
Beggs, E., {\it De Rham's theorem for infinite-dimensional
manifolds}, Quart.\ J. Math.\ \textbf{38}  (1987), 131--154. 
%
\bibitem[Bel04]{Bel04} 
Belti\c{t}\u{a}, D., {\it Asymptotic 
products and enlargibility of Banach-Lie algebras},
J. Lie Theory {\bf 14} (2004), 215--226. 
%
\bibitem[Bel06]{Bel06} ---, ``Smooth Homogeneous Structures 
in Operator Theory,''  Chapman and Hall, CRC Monographs and 
Surveys in Pure and Applied mathematics, 2006. 

\bibitem[Bel09]{Bel09} ---, {\it Iwasawa decompositions of some 
infinite-dimensional Lie groups}, 
Trans. Amer. Math. Soc. {\bf 361:12} (2009), 6613--6644

\bibitem[BB15]{BB15} Belti\c{t}\u{a}, D., and I.~Belti\c{t}\u{a}, 
{\it Faithful representations of infinite-dimensional nilpotent Lie algebras}, 
Forum Math. {\bf 27:1} (2015), 255--267

\bibitem[BelN08]{BelN08}
Belti\c t\u a, D., and K.-H. Neeb, 
{\it Finite-dimensional Lie subalgebras of algebras with continuous inversion},
Studia Math.\ \textbf{185} (2008), 249--262.

\bibitem[BelN10]{BelN10}
---,
\emph{Geometric characterization of Hermitian algebras with continuous inversion},
Bull.\ Aust.\ Math.\ Soc.\ {\bf 81} (2010), 96--113.

\bibitem[BNi15]{BNi15}
Belti\c{t}\u{a}, D., and M. Nicolae,
{\it On universal enveloping algebras in a topological setting},
Stud.\ Math.\ {\bf 230:1} (2015), 1--29.
%
\bibitem[BP07]{BP07}
Belti\c{t}\u{a}, D., and B. Prunaru, {\it Amenability, completely bounded 
projections, dynamical systems and smooth orbits},
Integral Equations Operator Theory \textbf{57} (2007), 1--17.
%
\bibitem[BR05]{BR05}
Belti\c t\u a, D., and T.~S. Ratiu,
{\it Symplectic leaves in real Banach Lie-Poisson spaces}, 
Geom.\ Funct.\ Anal.\ \textbf{15} (2005), 753--779. 
%
\bibitem[BR07]{BR07}
---, {\it  Geometric representation theory for unitary groups of operator algebras}, 
Adv.\ Math.\ \textbf{208} (2007), 299--317. 
%
\bibitem[BS01]{BS01} Belti\c t\u a, D., and M.~Saba\c c, ``Lie Algebras of Bounded 
Operators,'' Operator Theory: Advances and Applications {\bf 120}, 
Birkh\"auser Verlag, Basel, 2001. 
%
\bibitem[BZ96]{BZ96} Benkart, G., and E. Zelmanov, {\it Lie algebras graded by
finite root systems and intersection matrix algebras}, 
Invent.\ Math.\ \textbf{126} (1996), 1--45. 
%
\bibitem[Be74]{Be74}
Berezin, F.~A., {\it Quantization in complex symmetric spaces}, 
Math.\ USSR-Izv.\ \textbf{9} (1974), 341--379. 
%
\bibitem[BGV04]{BGV04} Berline, N., E. Getzler, and M. Vergne, ``Heat Kernels and 
Dirac Operators,'' Springer-Verlag, 2004.
%
\bibitem[BM92]{BM92} Berman, S., and R.~V. Moody, {\it Lie algebras
graded by finite root systems and the intersection matrix algebras
of Slodowy}, Invent.\ Math.\ \textbf{108} (1992), 323--347.
%
\bibitem[Ber00]{Ber00}
Bertram, W., ``The Geometry of Jordan and Lie Structures'', 
Lecture Notes in Math.\ {\bf 1754}, Springer-Verlag, Berlin, 2000. 
%
\bibitem[Ber01]{Ber01} ---, {\it Generalized projective geometries: general
theory and equivalence with Jordan structures}, 
Adv.\ Geom.\ \textbf{2}  (2001), 329--369. 
%
\bibitem[Ber08]{Ber08}
---, ``Differential Geometry, Lie Groups and Symmetric Spaces
over General Base Fields and Rings'',
Mem.\ Am.\ Math.\ Soc.\ \textbf{192}, 2008. 
%
\bibitem[BGN04]{BGN04}
Bertram, W., H.~Gl\"{o}ckner, and
K.-H.~Neeb, {\it Differential calculus over general base fields and rings}, 
Expo.\ Math.\ \textbf{22} (2004), 213--282. 
%
\bibitem[BerN04]{BerN04} Bertram, W., and K.-H. Neeb, {\it Projective completions of Jordan pairs, Part I. 
The generalized projective geometry of a Lie algebra}, 
J.~Algebra {\bf 277:2} (2004), 474--519.
%
\bibitem[BerN05]{BerN05} ---, {\it Projective completions of Jordan pairs, Part II}, 
Geom.\ Dedicata {\bf 112:1} (2005), 75--115. 
%
\bibitem[BP75]{BP75}
Bessaga, C. and A. Pe\l{}czy\'{n}ski,
``Selected Topics in Infinite-Dimensional Topology,'' PWN, Warsaw, 1975.
%
\bibitem[Bi88]{Bi88}
Bierstedt, K.-D., {\it An introduction to
locally convex inductive limits}, pp.~35--133
in: Hogbe-Nlend, H. (ed),
``Functional Analysis and its
Applications'', (Nice, 1986), 
World Sci.\ Publ., Singapore, 1988. 
%
\bibitem[BBP94]{BBP94}
Bierstedt, K.-D., J. Bonet, and A. Peris, {\it Vector-valued 
holomorphic germs on Fr\'{e}chet-Schwartz
spaces}, Proc.\ Roy.\ Irish Acad Sect.\ A \textbf{94}  (1994), 31--46. 
%
\bibitem[BM77]{BM77}
Bierstedt, K.-D., and R. Meise, {\it Nuclearity and the Schwartz property in the theory
of holomorphic functions on metrizable locally convex spaces}, pp.~93--129
in: Matos, M.~C. (ed),
``Infinite Dimensional Holomorphy and Applications'', 
North-Holland, Amsterdam, 1977. 
%
\bibitem[Bil02]{Bil02}
Biller, H., {\it The exponential law
for smooth functions}, Manuscript, TU Darmstadt, July 2002. 
%
\bibitem[Bil04]{Bil04}
---, {\it Continuous inverse algebras with involution}, 
Forum Math.\ \textbf{22} (2010), 1033--1059.
%
\bibitem[Bil07]{Bil07}
---,
\emph{Analyticity and naturality of the multi-variable functional calculus},
Expo.\ Math.\ {\bf 25} (2007), 131--163.
%
\bibitem[Bil03]{Bil03} 
Billig, Y., {\it Abelian extensions of the group of
diffeomorphisms of a torus}, Lett.\ Math.\ Phys.\ \textbf{64:2} (2003), 155--169.
%
\bibitem[BiNe08]{BiNe08} Billig, Y., and K.-H. Neeb, {\it  
On the cohomology of vector fields on parallelizable manifolds}, 
Ann.\ Inst.\ Fourier {\bf 58} (2008), 1937--1982
%
\bibitem[BiPi02]{BiPi02} Billig, Y., and A. Pianzola, {\it Free Kac-Moody groups and their Lie algebras}, 
Algebr.\ Represent.\ Theory {\bf 5:2} (2002), 115--136. 
%
\bibitem[Bir36]{Bir36} Birkhoff, G., {\it Continuous groups and linear spaces}, 
Mat. Sbornik {\bf 1} (1936), 635--642 

\bibitem[Bir38]{Bir38} ---, {\it Analytic groups}, 
Trans. Amer. Math. Soc. {\bf 43} (1938), 61--101 

\bibitem[BG14]{BG14}
Birth, L., and H. Gl\"{o}ckner,
\emph{Continuity of convolution of test functions on Lie groups},
Can.\ J. Math.\ {\bf 66:1} (2014), 102--140.
%
\bibitem[Bisg93]{Bisg93}
Bisgaard, T.~M., {\it The topology
of finitely open sets is not a vector space topology}, Arch.\ Math.\ \textbf{60}  (1993), 
546--553. 
%
%
\bibitem[Bl98]{Bl98}
Blackadar, B., ``K-Theory for Operator Algebras'', 2nd edition, 
Cambridge Univ.\ Press, Cambridge, 1998. 
%
\bibitem[BS71a]{BS71a}
Bochnak, J., and J. Siciak, {\it Polynomials and multilinear mappings in topological vector spaces}, Studia Math.\ \textbf{39} (1971), 59--76.
%
\bibitem[BS71b]{BS71b} ---, {\it Analytic functions in topological vector spaces}, 
Studia Math.\ \textbf{39} (1971), 77--112. 
%
\bibitem[Bo46]{Bo46} Bochner, S., {\it Formal Lie groups},
Ann.\ of Math.\ {\bf 47} (1946), 192--201.   
%
\bibitem[BoMo45]{BoMo45} Bochner, S., and D. Montgomery, {\it Groups 
of differentiable and real or complex analytic transformations}, 
Ann.\ of Math.\ (2) {\bf 46} (1945), 685--694.  
%
\bibitem[BDS17]{BDS17} Bogfjellmo, G., R. Dahmen, and A. Schmeding,
  {\it Overview of (pro-)Lie group structures on Hopf algebra
    character groups}, Preprint, arXiv:1711.05963v1 [math.GR] 

\bibitem[Bm67]{Bm67} 
Boman, J., {\it Differentiability of a function and of
its compositions with functions of one variable}, 
Math.\ Scand.\ \textbf{20}  (1967), 249--268.
%
\bibitem[BDM94]{BDM94}
Bonet, J., P. Domanski, and
J. Mujica, {\it Complete spaces of vector-valued
holomorphic germs}, Math.\ Scand.\ \textbf{75} (1994), 150--160. 
%
\bibitem[BF12]{BF12} Bonfiglioli, A., and R. Fulci, 
``Topics in Noncommutative Algebra: The Theorem of Campbell, Baker, Hausdorff and Dynkin,'' 
Lecture Notes in Math.\ {\bf 2034}, Springer, 2012.
%
\bibitem[BD73]{BD73}  Bonsall, F.~F., and J. Duncan, ``Complete Normed Algebras,''
Ergeb.\ Mathem.\ und ihrer Grenzgebiete {\bf 80}, Springer-Verlag, New York, Heidelberg, 1973. 
%
\bibitem[BB08]{BB08}
Borzellino, J.~E., and V. Brunsden,
\emph{A manifold structure for the group of orbifold diffeomorphisms of a smooth orbifold},
J. Lie Theory {\bf 18} (2008), 979--1007. 
%
\bibitem[BCR81]{BCR81}
Boseck, H., G. Czichowski, and K.-P. Rudolph,
`` Analysis on Topological Groups -- General Lie Theory'', Teubner Texte zur Mathematik 
\textbf{137}, Teubner Verlag, Leipzig, 1981. 
%
\bibitem[Bos90]{Bos90}
Bost, J.-B., {\it Principe d'Oka, $K$-theorie et syst\`emes dynamiques 
non-commu\-ta\-tifs}, Invent.\ Math.\ \textbf{101} (1990), 261--333. 
%
\bibitem[Bo58]{Bo58}
Bott, R., {\it The space of loops on a Lie group}, Michigan Math.\ J.  
\textbf{5} (1958), 35--61. 
%
\bibitem[Bo59]{Bo59} ---, {\it The stable homotopy groups of the classical
groups}, Ann.\ Math.\ \textbf{70:2} (1959), 313--337. 
%
\bibitem[Bo60]{Bo60} ---, {\it A note on the Samelson product in the 
classical groups}, Comment.\ Math.\ Helv.\ {\bf 34} (1960), 249--256. 
%
\bibitem[Bo77]{Bo77} ---, {\it On the characteristic classes of groups of diffeomorphisms}, 
Enseign.\ Math.\ (2) {\bf 23:3-4} (1977), 209--220.  
%
\bibitem[Bou67]{Bou67}
Bourbaki, N., ``Vari\'{e}t\'{e}s diff\'{e}rentielles
et analytiques. Fascicule de r\'{e}sultats'',  
Hermann, Paris, 1967. 
%
\bibitem[Bou87]{Bou87} ---, ``Topological Vector Spaces (Chapters~1--5)'', 
Springer-Verlag, Berlin, 1987. 
%
\bibitem[Bou88]{Bou88}
---, ``General Topology (Chapters 1--6)'', Springer-Verlag, Berlin, 1988. 
%
\bibitem[Bou66]{Bou66} ---,  ``General Topology, Part 2,''
Addison-Wesley, Reading and Hermann, Paris, 1966.
%
\bibitem[Bou89]{Bou89}
---, ``Lie Groups and Lie Algebras
(Chapters 1--3)'', Springer-Verlag, Berlin, 1989. 
%
\bibitem[Bo80]{Bo80}
Boyer, R., {\it Representations of the Hilbert Lie group
$U({\cal H})_2$}, Duke Math.\ J. \textbf{47}  (1980), 325--344. 
%
\bibitem[Bo83]{Bo83}
---, {\it Infinite traces of AF-Algebras and characters
of $U(\infty)$}, Operator Theory {\bf 9} (1983), 205--236. 
%
\bibitem[Bo92]{Bo92}
---, {\it Characters and factor representations of the
infinite-dimensional classical groups}, J. Oper.\ Theory \textbf{29} (1992),
281--307.
%
\bibitem[BIP08]{BIP08} Burago, D., S.\, Ivanov, and L.\, Polterovich, {\it 
Conjugation-invariant norms on groups of geometric origin}, in 
``Groups of diffeomorphisms,'' 221–250, 
Adv. Stud. Pure Math. {\bf 52}, Math. Soc. Japan, Tokyo, 2008

\bibitem[BEG89]{BEG89} Bratteli, O., 
G. A. Elliott, F.\ M. Goodman, and P.\ E.\ T.\ Jorgensen, 
{\it Smooth Lie group actions on non-commutative tori}, 
Nonlinearity {\bf 2} (1989),  271--286 
%
\bibitem[Bre72]{Bre72} Bredon, G. E., ``Introduction to Compact Transformation Groups'',  
Academic Press, 1972.
%
\bibitem[Bre93]{Bre93} ---, ``Topology and Geometry'', Grad.\ Texts Math.\ \textbf{139}, Springer-Verlag, Berlin, 1993.
%
\bibitem[Br70]{Br70} 
Breuer, M., {\it On the homotopy type of the groups of regular
elements of semifinite von Neumann algebras}, Math.\ Ann.\ {\bf 185}  (1970), 61--74. 
%
\bibitem[BJ73]{BJ73} 
Br\"ocker, Th., and K. J\"anich, ``Einf\"uhrung in die Differentialtopologie'', 
Springer-Verlag, Berlin, 1973. 
%
\bibitem[Bro82]{Bro82} 
Brown, K. S., ``Cohomology of Groups'', Grad.\ Texts Math.\  
\textbf{87}, Springer-Verlag, Berlin, 1982. 
%
\bibitem[BW76]{BW76} 
Br\"uning, J., and W. Willgerodt, {\it Eine Verallgemeinerung
eines Satzes von N.~Kuiper}, Math.\ Ann.\ \textbf{220} (1976), 47--58. 
%
\bibitem[BrW59]{BrW59}Bruhat, F., and H. Whitney, H.,
{\it Quelques propri\'{e}t\'{e}s fondamentales des ensembles analytique-r\'{e}els},
Comment.\ Math.\ Helv.\ {\bf 33} (1959), 132--160.
%
\bibitem[BBM14]{BBM14} Bauer, M., Bruveris, M., and P.W. Michor, 
{\it  Overview of the geometries of shape spaces and diffeomorphism groups}, 
J. Math. Imaging Vision {\bf 50:1-2} (2014), 60--97. 
%
\bibitem[BBM16]{BBM16} ---, {\it Why use Sobolev metrics 
on the space of curves}, 
in ``Riemannian computing in computer vision,'' 233--255, Springer, 
Cham, 2016. 
%
\bibitem[Bry93]{Bry93}
Brylinski, J.-L., ``Loop Spaces, Characteristic Classes and
Geometric Quantization'', Progress in Math.\ \textbf{107}, Birkh\"auser
Verlag, Boston, 1993. 
%
\bibitem[CaP19]{CaP14}
Cabau, P., and F. Pelletier,
{\it Integrability on direct limits of Banach manifolds},
Ann.\ Fac.\ Sci.\ Toulouse, Math.\ (6)
\textbf{28:5} (2019), 909--956. 
%
\bibitem[Cal70]{Cal70} Calabi, E., {\it On the group of automorphisms of a symplectic manifold}, 
pp.~1-26 in: R.~C. Gunning (ed),
``Problems in Analysis,'' (Lectures at the Sympos.\ in honor of Salomon Bochner), 
Princeton Univ.\ Press, Princeton, N.J., 1970. 
%
\bibitem[Cam97]{Cam97} Campbell, J.~E., {\it On a law of combination of operators bearing on the theory 
of continuous transformation groups}, Proc.\ London Math.\ 
Soc.\ {\bf 28} (1897), 381--390.  
%
\bibitem[Cam98]{Cam98} ---, {\it On a law of combination of operators. (Second paper)}, 
Proc.\ London Math.\ 
Soc.\ {\bf 28} (1897), 381--390.  
%
\bibitem[Ca18]{Ca18}  Caprace, P.-E., {\it Non-discrete simple locally compact 
groups}, Preprint, 2018 
%
\bibitem[CP89a]{CP89a}
Carey, A.. and J. Palmer, {\it Infinite complex spin groups}, 
J. Funct.\ Anal.\ \textbf{83} (1989), 1--43. 
%
\bibitem[CP89b]{CP89b}
---, {\it Gauge anomalies on $\bS^2$ and group 
extensions}, J. Math.\ Phys.\ \textbf{30:9} (1989), 2181--2191. 
%
\bibitem[CaE98]{CaE98} Cartan, E., {\it Les groupes bilin\'eaires et les syst\`emes 
de nombres complexes}, Ann.\ Fac.\ Sci.\ Toulouse Sci.\ Math.\ Sci.\ Phys.\ 
{\bf 12:1} (1898), B1--B64. 
%
\bibitem[CaE01]{CaE01} ---, {\it L'int\'egration des syst\`emes d'\'equations
aux differentielles totales},
Ann.\ Sci.\ Ecol.\ Norm.\ Sup.\ (3) {\bf 18} 
(1901), 241--311. 
%
\bibitem[CaE04]{CaE04} ---, {\it Sur la structure des groups infinies des transformations}, 
Ann.\ Sci.\ Ecol.\ Norm.\ Sup.\ {\bf 21} (1904), 153--206; 
{\bf 22} (1905), 219--308.  
%
\bibitem[CaE08]{CaE08} ---, {\it Les sous-groupes des groupes continus de transformations}, 
Ann.\ Sci.\ Ecol.\ Norm.\ Sup.\  {\bf 25} (1908), 57--194.  
%
\bibitem[CaE09]{CaE09} ---, {\it Les groupes de transformations continus, infinis, simples}, 
Ann.\ Sci.\ Ecol.\ Norm.\ Sup.\ {\bf 26} (1909), 93--161.  
%
\bibitem[CaE30]{CaE30} ---, {\it Le troisi\`eme th\'eor\`eme fondamental de Lie}, 
C. R. Acad.\ Sci.\ Paris {\bf 190} (1930), 914--916, 1005--1007.  
%
\bibitem[CaE36]{CaE36} ---, {\it La topologie des groupes de Lie. (Expos\'es de g\'eom\'etrie Nr. 8.)}, 
Actual.\ sci.\ industr.\ {\bf 358} (1936), 28 p.
%
\bibitem[CaE52a]{Ca52a} ---, 
{\it Le troisi\`eme th\'eor\`eme fondamental de Lie}, pp.\ 1143--1148
in: ``Oeuvres I'', Gauthier--Villars, Paris , 1952.  
%
\bibitem[CaE52b]{Ca52b} ---, {\it La topologie des espaces repr\'esentatifs de groupes de
Lie}, pp.\ 1307--1330 in: ``Oeuvres I'', Gauthier--Villars, Paris, 1952. 
%
\bibitem[CaE52c]{Ca52c} ---, {\it Les repr\'esentations lin\'eaires des groupes
de Lie}, pp.\ 1339--1350 in: ``Oeuvres I'', Gauthier--Villars, Paris, 1952. 
%
\bibitem[CaH67]{CaH67}
Cartan, H., ``Calcul Diff\'{e}rentiel'', Hermann, Paris, 1967. 
%
\bibitem[CVL98]{CVL98} 
Cassinelli, G., E. de Vito, P. Lahti, and A. Levrero, {\it 
Symmetries of the quantum state space and group representations},  
Reviews in Math.\ Physics \textbf{10:7} (1998), 893--924. 
%
\bibitem[Ce61]{Ce61}
Cerf, J.,
\emph{Topologie de certains \'{e}spaces de plongements},
Bull.\ Soc.\ Math.\ France {\bf 89} (1961), 227--380.
%
\bibitem[CP86]{CP86}
Chari, V., and A. Pressley, {\it New unitary representations
of loop groups}, Math.\ Ann.\ \textbf{275}  (1986), 87--104. 
%
\bibitem[CM70]{CM70} Chernoff, P., and J. Marsden, 
{\it On continuity and smoothness of group actions}, Bull.\ Amer.\ Math.\ Soc.\ {\bf 76} (1970), 
1044--1049.
%
\bibitem[Ch46]{Ch46}
Chevalley, C., ``Theory of Lie Groups I'', Princeton Univ.\ Press, Princeton, 1946. 
%
\bibitem[CE48]{CE48}
Chevalley, C., and S. Eilenberg, {\it Cohomology theory of Lie groups and Lie 
algebras},  Trans.\ Am.\ Math.\ Soc.\ \textbf{63} (1948), 85--124. 
%
\bibitem[CDD85]{CDD85}
Choquet-Bruhat, Y.,
C. DeWitt-Morette, and M. Dillard-Bleick,
``Analysis, Manifolds and Physics,''
North-Holland (Elsevier), ${}^2$1985.
%
\bibitem[CH82]{CH82}
Chow, S.~N., and J.~K. Hale, ``Methods of Bifurcation Theory,'' Springer-Verlag, New York, 1982.
%
\bibitem[CIJM19]{CIJM19} Ciaglia, F.M., A.~Ibort, J.~Jost, and G.~Marmo, {\it 
Manifolds of classical probability distributions and quantum density 
operators in infinite dimensions}, Inf. Geom. {\bf 2:2} (2019), 231--271.
%
\bibitem[ClH12]{ClH12}
Clark, D. E., and J. Houssineau,
{\it Fa\'{a} di Bruno's formula for G\^{a}teaux differentials and interacting stochastic population processes},
preprint, 2012, {\tt arXiv:1202.0264}.
%
\bibitem[Co85]{Co85}
Connes, A., {\it Noncommutative Differential Geometry}, 
Publ.\ Math.\ IHES \textbf{62} (1985), 41--144. 
%
\bibitem[Co94]{Co94} Connes, A., ``Non-commutative  Geometry,'' Academic Press, 1994.  
%
\bibitem[CM79]{CM79} Conway, J. B., and B. B. Morrel, {\it Operators that are 
points of spectral continuity},
Integral Equations Oper.\ Theory {\bf 2:2} (1979), 174--198. 
%
\bibitem[CGM90]{CGM90} Cuenca Mira, J.\ A., A.\ Garcia Martin, and C.\ Martin
Gonzalez, {\it Structure theory of $L^*$-algebras}, Math.\ Proc.\
Camb.\ Phil.\ Soc.\ {\bf 107} (1990), 361--365. 
%
\bibitem[Cu04]{Cu04}
Cuntz, J., {\it Cyclic theory and the bivariant Chern--Connes character}, 
pp.~73--135 in:
Doplicher, S. et al.\ (eds) ``Noncommutative Geometry,''
Lect.\ Notes Math.\ {\bf 1831}, 2004.
%
\bibitem[Da10]{Da10}
Dahmen, R., \emph{Analytic mappings between LB-spaces and applications in infinite-dimensional Lie theory},
Math.\ Z. {\bf 266} (2010), 115--140.
 %
\bibitem[Da11]{Da11}
---,
``Direct limit constructions in infinite-dimensional Lie theory,''
Ph.D.-thesis, University of Paderborn, 2011.
%
\bibitem[Da14]{Da14}
---, {\it Regularity in Milnor’s sense for ascending unions of Banach-Lie groups},
J. Lie Theory {\bf 24:2} (2014), 545--560.
%
\bibitem[Da15]{Da15}
---,
\emph{Peetre's theorem in the locally convex setting},
Forum Math.\ {\bf 27} (2015),
2957--2979. 
%
\bibitem[DS15]{DS15}
Dahmen, R., and A. Schmeding,
{\it The Lie group of real analytic diffeomorphisms is not real analytic},
Stud.\ Math.\ {\bf 229:2} (2015), 141--172.
%
\bibitem[DGS14]{DGS14}
Dahmen, R., H. Gl\"{o}ckner, and A. Schmeding,
{\it Complexifications of infinite-dimensional manifolds and
new constructions of infinite-dimensional Lie groups},
preprint, 2014; arXiv:1410.6468.
%
\bibitem[DP03]{DP03} Dai, J., and D. Pickrell, {\it The orbit method and the 
Virasoro extension of $\Diff_+(S^1)$. I. Orbital integrals}, 
J.~Geom.\ Phys.\ {\bf 44:4} (2003), 623--653.  
%
\bibitem[Da94]{Da94}
Dazord,  P., {\it Lie groups and algebras in infinite dimension: a new approach}, 
pp.~17--44 in:
``Symplectic Geometry and Quantization,'' 
Contemp.\ Math.\ {\bf 179}, Amer.\ Math.\ Soc., Providence, RI, 1994. 
%
\bibitem[Dei77]{Dei77}
Deimling, K.,
``Ordinary Differential Equations in Banach Spaces,''
Lecture Notes in Math.\ {\bf 596}, Springer Verlag, Berlin, 1977.
%
\bibitem[De32]{De32}
Delsartes, J., ``Les groups de transformations lin\'eaires dans l'espace 
de Hilbert,'' M\'emoirs des Sciences Math\'ematiques, Fasc.\ {\bf 57}, Paris. 
%
\bibitem[DWL82]{DWL82} De Wilde, M., and P. B. A. Lecomte, {\it Isomorphisms of Lie algebras of 
vector fields},
Commentat.\ Math.\ Univ.\ Carol.\ {\bf 23} (1982), 513--523.  
%
\bibitem[DW97]{DW97}
Dierolf, S., and J. Wengenroth, {\it Inductive limits of topological algebras}, 
Linear Topol.\ Spaces Complex Anal.\ \textbf{3} (1997), 45--49. 
%
\bibitem[Di44]{Di44}
Dieudonn\'{e}, J.,       
\emph{Sur la completion des groupes topologiques},
C. R. Acad.\ Sci.\ Paris {\bf 218} (1944), 774--776.
%
%
\bibitem[Di60]{Di60}
Dieudonn\'{e}, J., ``Foundations of Modern Analysis'', Academic Press, New York, 1960. 
%
\bibitem[DiPe99]{DiPe99} Dimitrov, I., and I.~Penkov, {\it Weight modules of direct
limit Lie algebras},
International Math.\ Res.\ Notices {\bf 5} (1999), 
223--249.
%
\bibitem[DJNV21]{DJNV21} Diez, T., Janssens, B., Neeb, K.-H., and C.~Vizman, {\it 
Central extensions of Lie groups preserving a differential form}, 
Int. Math. Res. Not. IMRN {\bf 5} 2021, 3794--3821

\bibitem[DR19]{DR19} Diez, T., and G.~Rudolph, {\it Slice theorem and orbit type 
stratification in infinite dimensions}, Differential Geom. Appl. {\bf 65} (2019), 176--211.

\bibitem[Dn81]{Dn81}
Dineen, S., {\it Holomorphic germs on compact subsets of 
locally convex spaces}, pp.~247--263 in:
``Functional Analysis, Holomorphy, and Approximation Theory'' 
(Rio de Janeiro, 1978), 
Lecture Notes in Math., {\bf 843}, 
Springer-Verlag, Berlin, 1981. 
%
\bibitem[dCa92]{dCa92}
do Carmo, M. P., ``Riemannian Geometry'', Birkh\"{a}user, Boston, 1992. 
%
\bibitem[DiD63]{DiD63} Dixmier, J., and A. Douady, {\it Champs continus 
d'espaces hilbertiens et de $C^*$-alg\`ebres}, Bull. Soc. Math. France 
{\bf 91} (1963), 227--284
%
\bibitem[DGV16]{DGV16}
Dodson, C.~T.~J., G. Galanis, and E.~Vassiliou,
``Geometry in a Fr\'{e}chet Context. A Projective Limit Approach,''
London Math.\ Soc.\ Lecture Note Series {\bf 428},
Cambridge University Press, Cambridge, 2016.
%
\bibitem[DI85]{DI85}
Donato, P., and P. Iglesias, {\it Examples de groupes diff\'eologiques: 
flots irrationnels sur le tore}, C. R. Acad.\ Sci.\ Paris, ser.\ 1, 
\textbf{301} (1985), 127--130. 
%
\bibitem[Dou61]{Dou61}
Douady, A.,
\emph{Vari\'{e}t\'{e}s \`{a} bord anguleux et voisinages tubulaires},
S\'{e}minaire Henri Cartan, 1961/62, Exp.\,1, 11 pp.
%
\bibitem[Dou66]{Dou66} Douady, A., {\it Un espace de Banach dont le groupe lin\'eaire n'est pas 
connexe}, Nederl.\ Akad.\ Wetensch.\ Proc.\ Ser.\ A 68 = Indag.\ Math.\ {\bf 27} (1965), 
787--789. 
%
\bibitem[DL66]{DL66}
Douady, A., and M. Lazard, {\it Espaces fibr\'es en alg\`ebres
de Lie et en groupes}, Invent.\ Math.\ {\bf 1} (1966), 133--151. 
%
\bibitem[Dr69]{Dr69} Dress, A., {\it Newman's Theorem 
on transformation groups}, Topology {\bf 8} (1969), 203--207.  
%
\bibitem[DG01]{DG01} Dupre, M. J., and J. F. Glazebrook, {\it The Stiefel bundle of a Banach 
algebra},
Integr.\ equ.\ oper.\ theory {\bf 41} (2001), 264--287.  
%
\bibitem[DEG98]{DEG98} Dupre, M. J., J.-C.~Evard, 
and J. F. Glazebrook, {\it Smooth parametrization of subspaces 
in a Banach space}, Revista de la Uni\'on Mathem\'atica 
Argentina {\bf 41:2} (1998), 1--13. 

\bibitem[Dy47]{Dy47} Dynkin, E. B., {\it Calculation
of the coefficients in the Campbell--Hausdorff
formula} (Russian), Doklady Akad.\ Nauk.\
SSSR (N.S.) {\bf 57} (1947), 323--326.  
%
\bibitem[Dy53]{Dy53} ---, ``Normed Lie Algebras and Analytic Groups,'' 
Amer.\ Math.\ Soc.\ Translation 1953, no. 97, 66 pp.  
%
\bibitem[Dz92]{Dz92} 
Dzhumadildaev, A., {\it Central extensions of 
infinite-dimensional Lie algebras}, Funct.\ Anal.\ Appl.\ \textbf{26:4}  (1992), 247--253. 
%
\bibitem[EE67]{EE67} Earle, C. J., and J.~Eells, {\it The 
diffeomorphism group of a compact Riemann surface}, 
Bull. Amer. Math. Soc. {\bf 73} (1967), 557--559

\bibitem[EE69]{EE69} ---, {\it 
A fibre bundle description of Teichm\"uller theory}, 
J. Differential Geometry {\bf 3} (1969), 19--43 

\bibitem[ES70]{ES70} Earle, C. J., and A.~H.\ Schatz, 
{\it Teichm\"uller theory for surfaces with boundary}, 
J. Differential Geometry {\bf 4} (1970), 169--185 

\bibitem[Eb70]{Eb70} Ebin, D. G., {\it The manifold of riemannian metrics}, 
in ``Global Analysis'' (Berkeley, Calif., 1968), Amer.\ Math.\ Soc.\ 
Proc.\ Symp.\ Pure Math.\ {\bf 15} (1970), 11--40.  
%
\bibitem[EM69]{EM69} Ebin, D. G., and J. E. Marsden, {\it Groups of diffeomorphisms and 
the solution of the classical Euler equations for a perfect fluid}, 
Bull.\ Amer.\ Math.\ Soc.\ {\bf 75} (1969), 962--967. 
%
\bibitem[EM70]{EM70} ---, {\it Groups of diffeomorphisms and the motion of 
an incompressible fluid}, Ann.\ of Math.\ {\bf 92} (1970), 102--163.  
%
\bibitem[EMi99]{EMi99} Ebin, D. G., and G. Misiolek, {\it The exponential map on 
${\cal D}^s_\mu$},
in ``The Arnoldfest (Toronto, ON, 1997),'' 
Fields Inst.\ Commun.\ {\bf 24},
Amer.\ Math.\ Soc., Providence, RI, 1999, 153--163. 
%
\bibitem[Ee58]{Ee58} Eells, J. Jr., {\it On the geometry of function spaces}, in 
``International Symposium on Algebraic Topology,''
pp.~303--308; 
Universidad Nacional Autonoma de M\'exico and UNESCO, Mexico City. 
%
\bibitem[Ee66]{Ee66} ---, {\it A setting for global analysis}, Bull.\ Amer.\
Math.\ Soc.\ {\bf 72} (1966), 751--807.  
%
\bibitem[Ed99]{Ed99} Edamatsu, T., {\it On the bamboo-shoot topology of
certain inductive limits of topological groups}, 
J. Math.\ Kyoto Univ.\ \textbf{39} (1999), 715--724. 
%
\bibitem[EW17]{EW17}
Egeileh, M., and T. Wurzbacher,
{\it Infinite-dimensional manifolds as ringed spaces},
Publ.\ Res.\ Inst.\ Math.\ Sci.\ {\bf 53:1} (2017), 187--209.
%
\bibitem[Ei07]{Ei07}
Eichhorn, J.,
``Global Analysis on Open Manifolds,''
Nova Science Publishers, New York, 2007.
%
\bibitem[ES96]{ES96} Eichhorn, J., and R. Schmid, {\it Form 
preserving diffeomorphisms on open manifolds}, 
Ann.\ Global Anal.\ Geom.\ {\bf 14:2} (1996), 147--176. 
%
\bibitem[ES01]{ES01} ---, {\it Lie groups of Fourier integral operators on open manifolds}, 
Comm.\ Anal.\ Geom.\ {\bf 9:5} (2001), 983--1040. 
%
\bibitem[Ei68]{Ei68} Eichler, M., {\it A new proof of the Baker--Campbell--Hausdorff formula}, 
J.\ Math.\ Soc.\ Japan {\bf 20} (1968), 23--25.  
%
\bibitem[El67]{El67} Eliasson, H., {\it Geometry of manifolds of maps}, J.\ Diff.\
Geom.\ {\bf 1} (1967), 169--194.  
%
\bibitem[Ek89]{Ek89}
Engelking, R.,  ``General Topology'', Heldermann Verlag, Berlin, 1989. 
%
\bibitem[Ep70]{Ep70}
Epstein, D.~B.~A., {\it The simplicity of certain groups of homeomorphisms}, 
Comp.\ Math.\ \textbf{22} (1970), 165--173. 
%
\bibitem[Ep84]{Ep84} ---, {\it Commutators of $C^\infty$-diffeomorphisms}, 
Comment.\ Math.\ Helv.\ {\bf 59} (1984), 111--122.  
%
\bibitem[Est55]{Est55}
van Est, W.~T., {\it On the algebraic cohomology concepts in
Lie groups I,II}, Nederl.\ Akad.\ 
Wet., Proc.\ \textbf{A 58}=Indag.\ math.\ {\bf 17}  (1955), 225--233; 286--294. 
%
\bibitem[Est62]{Est62}
---, {\it Local and global groups}, 
Ned.\ Akad.\ Wet., Proc.\ {\bf A 65}=Indag.\ math.\ {\bf 24} (1962), 391--425. 
%
\bibitem[Est66]{Est66} ---, {\it On Ado's theorem},
Proc.\ Kon.\ Ned.\ Akad.\ v.\ Wet.\ Series {\bf A 69}
=Indag.\ Math.\ {\bf 28} (1966), 176--191.  
%
\bibitem[Est84]{Est84} ---, {\it Rapport sur les S-atlas}, Ast\'erisque {\bf 116} (1984), 235--292.  
%
\bibitem[Est88]{Est88} ---, {\it Une d\'emonstration de E.\ Cartan du
troisi\`eme th\'eor\`eme de Lie}, pp.~83--96 
in: Dazord P. et al.\ (eds) ``Seminaire sud-rhodanien de g\'eom\'etrie
VIII: Actions hamiltoniennes de groupes; troisi\`eme th\'eor\`eme de
Lie'', Hermann, Paris, 1988. 
%
\bibitem[EK64]{EK64}
van Est, W.~T., and Th.\ J. Korthagen, {\it Non-enlargible Lie
algebras}, 
Ned.\ Acad.\ v.\ Wet., Proc.\ {\bf A 67}=Indag.\ Math.\ {\bf 26} (1964), 15--31.
%
\bibitem[ES73]{ES73} van Est, W. T., and S. \'Swierczkowski, {\it The path functor and 
faithful representability of Banach Lie algebras},
in ``Collection of articles dedicated to the memory of Hannare Neumann, I.''
J. Austral.\ Math.\ Soc.\ {\bf 16} (1973), 54--69. 
%
\bibitem[vES73]{vES73} van Est, W. T., and S.~\'Swierczkowski, {\it The 
path functor and faithful representability of Banach Lie algebras}, 
in ``Collection of articles dedicated to the memory of Hannare Neumann, I,'' 
J. Austral. Math. Soc. {\bf 16} (1973), 54--69

\bibitem[EF94]{EF94}
Etingof, P.~I., and I. B. Frenkel, {\it Central extensions of current 
groups in two dimensions}, Commun.\ Math.\ Phys.\ \textbf{165} (1994), 429--444. 
%
\bibitem[EK98]{EK98} 
Evans, D.~E., and Y. Kawahigashi, ``Quantum Symmetries on
Operator Algebras'', Oxford Math.\ Monographs, Clarendon Press, 1998. 
%
\bibitem[Eyn12]{Eyn12}
Eyni, J.~M., {\it Frobeniuss\"{a}tze f\"{u}r Vektordistributionen auf unendlich-dimensionalen
Mannigfaltigkeiten}, Master's thesis, Universit\"{a}t Paderborn, 2012; 
(advisor: H. Gl\"{o}ckner)
%
\bibitem[Eyn14]{Eyn14}
---,
{\it The Frobenius theorem for Banach distributions on infinite-dimensional manifolds
and applications in infinite-dimensional Lie theory}, preprint, 2014,
arXiv:1407.3166
%
\bibitem[Eyn15]{Eyn15}
---, {\it A Lie group structure on the group of real analytic diffeomorphisms of a
compact real analytic manifold with corners}, prepint, 2015, arXiv:1512.01506
%
\bibitem[Eyn16]{Eyn16}
---, ``New Examples and Constructions in Infinite-Dimensional Lie Theory,''
Doctoral Thesis, Universit\"{a}t Paderborn, 2016 (advisor: H. Gl\"{o}ckner);
see nbn-resolving.de/urn:nbn:de:hbz:466:2-27103
%
%
\bibitem[Fef07]{Fef07}
Fefferman, C.,
\emph{$C^m$-extension by linear operators},
Ann. Math. {\bf 166} (2007), 779--835.
%
\bibitem[Fe88]{Fe88}
Feigin, B.~L., {\it On the cohomology of the Lie algebra of
vector fields and the current algebra}, Sel.\ Math.\ Sov.\ \textbf{7:1} (1988), 49--62.
%
\bibitem[FF00]{FF00} Feigin, B.~L., and D. B. Fuks, {\it Cohomologies of Lie groups and Lie algebras}, 
In: I.~L.~Onishchik and E.~B. Vinberg (eds), ``Lie Groups and Lie Algebras II'', 
Encyclopaedia of Math.\ Sciences \textbf{21}, Springer-Verlag, Berlin, 2000. 
%
\bibitem[FHS17]{FHS17}
Ferrer, M.~V., S. Hern\'{a}ndez, and D. Shakhmatov,
\emph{A countable free closed non-reflexive subgroup of $\Z^{\mathfrak c}$},
Proc.\ Am.\ Math.\ Soc.\ {\bf 145:8} (2017), 3599--3605.
%
\bibitem[Fil82]{Fil82} Filipkiewicz, R. P., {\it Isomorphisms between diffeomorphism groups}, 
Ergodic Theory Dynamical Systems {\bf 2:2} (1983), 159--171. 
%
\bibitem[Fl71]{Fl71}
Floret, K., {\it Lokalkonvexe Sequenzen mit kompakten Abbildungen}, 
J.~Reine Ang.\ Math.\ \textbf{247} (1971), 155--195.
%
\bibitem[Fl80]{Fl80} ---, {\it Some aspects of the theory of
locally convex inductive limits}, pp.~205--237  
in: K.-D. Bierstedt and B. Fuchssteiner (eds)
``Functional Analysis: Surveys and Results II'', 
North Holland, Amsterdam, 1980.
%
\bibitem[FS77]{FS77}
Franklin, S.~P., and B.~V.~ S. Thomas, \emph{A survey of $k_\omega$-spaces},
Topol.\ Proc.\ {\bf 2} (1977), 111--124.
%
\bibitem[Fre25]{Fre25}
Fr\'echet, M., {\it La notion de diff\'{e}rentielle
dans l'analyse g\'{e}n\'{e}rale}, Ann.\ Sci.\ Ec.\ Norm.\ Sup.\ III 
\textbf{42} (1925), 293--323. 
%
\bibitem[Fr84]{Fr84}
Frenkel, I.~B., {\it Orbital theory for affine Lie algebras}, 
Invent.\ Math.\ \textbf{77} (1984), 301--352. 
%
\bibitem[FK21]{FK21}
Frenkel, I.~B., and H.K. Kim, 
{\it Three dimensional construction of the Virasoro--Bott group}, 
arXiv:math.GT:2107.11693
%
\bibitem[Fre68]{Fre68} Freifeld, Ch., {\it One-parameter subgroups do not fill a 
neighborhood of the identity in an infinite-dimensional Lie (pseudo-) group},
1968, Battelle Rencontres. 1967 Lectures in Mathematics and Physics; 538--543, 
Benjamin, New York. 
%
\bibitem[Frk07]{Frk07}
Frerick, L.,
\emph{Extension operators for spaces of infinite differentiable Whitney jets},
J. Reine Angew.\ Math.\ {\bf 602} (2007), 123--154.
%
\bibitem[FJW11]{FJW11}
Frerick, L., E. Jord\'{a}, and J. Wengenroth,
{\it Tame linear extension operators for smooth Whitney functions},
J. Funct.\ Anal.\ {\bf 261:3} (2011), 591--603.
%
\bibitem[FJW16]{FJW16}
---,
\emph{Whitney extension operators without loss of derivatives},
Rev.\ Mat.\ Iberoam.\ {\bf 32} (2016), 377--390.
%
\bibitem[FB66]{FB66} Fr\"olicher, A., and W. Bucher, ``Calculus in Vector Spaces without Norm,'' 
Lecture Notes in Math.\ {\bf 30}, Springer-Verlag, Berlin, 1966.  
%
\bibitem[FK88]{FK88}
Fr\"olicher, A., and A. Kriegl, ``Linear Spaces and Differentiation Theory'', 
Wiley Interscience, Chichester, 1988. 
%
\bibitem[Fu70]{Fu70}
Fuchs, L., ``Infinite Abelian Groups, I'', Acad.\ Press, New York, 1970. 
%
\bibitem[Fu73]{Fu73} 
---, ``Infinite Abelian Groups, II'' Acad.\ Press, New York, 1973. 
%
\bibitem[Fu86]{Fu86}
Fuks, D.~B., ``Cohomology of Infinite Dimensional Lie Algebras'',
Contemp.\ Sov.\ Math., Consultants Bureau, New York, London, 1986. 
%
\bibitem[Fu04]{Fu04} Furutani, K., {\it Fredholm-Lagrangian--Grassmannian and 
the Maslov index}, J. Geom. Phys. {\bf 51} (2004), 269--331
%
\bibitem[Ga96]{Ga96} Galanis, G., {\it Projective limits of Banach-Lie groups}, 
Period.\ Math.\ Hungar.\ {\bf 32:3} (1996), 179--191. 
%
\bibitem[Ga97]{Ga97} ---, {\it On a type of linear differential equations in Fr\'echet spaces}, 
Ann.\ Scuola Norm.\ Sup.\ Pisa Cl.\ Sci.\ (4) {\bf 24:3} (1997), 501--510. 
%
\bibitem[Ga47]{Ga47} G\aa{}rding, L., {\it Note on continuous 
representations of Lie groups}, 
Proc.\ Nat.\ Acad.\ Sci.\ U.S.A. {\bf 33} (1947), 331--332. 
%
\bibitem[GF68]{GF68}
Gelfand, I., and D. B. Fuks, {\it Cohomology of the Lie algebra of vector fields 
on the circle}, Funct.\ Anal.\ Appl.\ \textbf{2:4}  (1968), 342--343. 
%
\bibitem[Gh01]{Gh01} Ghys, E., {\it Groups acting on the circle}, 
L'Enseignement Math\'ematique {\bf 47} (2001), 329--407
%
\bibitem[GG61]{GG61} Glashow, S. L., and M. Gell-Mann, {\it Gauge theories of vector particles}, 
Ann.\ Phys.\ {\bf 15} (1961), 437--460.
%
\bibitem[Gl02a]{Gl02a}
Gl\"{o}ckner, H., {\it Infinite-dimensional Lie groups without completeness
restrictions}, pp.~43--59 in: A. Strasburger, W. Wojtynski, J. Hilgert, and K.-H. Neeb (eds),
``Geometry and Analysis on Finite and Infinite-dimensional Lie Groups'',
Banach Center Publications \textbf{55}, Warsaw, 2002. 
%
\bibitem[Gl02b]{Gl02b}
---, {\it Algebras whose groups of units are Lie groups}, Studia Math.\ \textbf{153} (2002), 
147--177. 
%
\bibitem[Gl02c]{Gl02c}
---, {\it Lie group structures on quotient groups
and universal complexifications for infinite-dimensional Lie groups}, 
J. Funct.\ Anal.\ \textbf{194} (2002), 347--409. 
%
\bibitem[Gl03a]{Gl03a}
---, {\it Direct limit Lie groups and manifolds}, 
J. Math.\ Kyoto Univ.\ \textbf{43} (2003), 
1--26. 
%
\bibitem[Gl03b]{Gl03b}
---,  {\it Lie groups of measurable mappings}, Canadian J.\ Math.\ \textbf{55} (2003), 
969--999.
%
\bibitem[Gl03c]{Gl03c}
---, {\it Patched locally convex spaces,
almost local mappings, and diffeomorphism groups of
non-compact manifolds}, manuscript in preparation, 2003.
%
\bibitem[Gl04a]{Gl04a}
---, {\it Lie groups of germs of analytic mappings}, pp.~1--16 in: 
T. Wurzbacher (ed), ``Infinite Dimensional Groups and Manifolds 
in Mathematics and Quantum Physics,''
IRMA Lecture Notes in Mathematics
and Math.\ Physics, de Gruyter, 2004.
%
\bibitem[Gl04b]{Gl04b}
---, {\it Lie groups over non-discrete topological fields}, preprint, 2004,
arXiv:math.GR/0408008
%
\bibitem[Gl05a]{Gl05a}
---, {\it H\"older continuous homomorphisms between infinite-dimensional 
Lie groups are smooth}, J. Funct.\ Anal.\ {\bf 228:2} (2005), 419--444.
%
\bibitem[Gl05b]{Gl05b}
---, \emph{Conveniently H\"{o}lder homomorphisms are smooth in the convenient sense},
Ann.\ Global Anal.\ Geom.\ {\bf 27} (2005), 227--255. 
%
\bibitem[Gl05c]{Gl05c}
 ---, {\it Fundamentals of direct limit Lie theory}, 
Compositio Math.\ {\bf 141} (2005), 1551--1577.
%
\bibitem[Gl05d]{Gl05d}
---, {\it $\Diff(\R^n)$ as a Milnor--Lie group},
Math.\ Nachr.\ {\bf 278:9} (2005), 1025--1032. 
%
\bibitem[Gl06a]{Gl06a}
---, {\it Implicit functions from topological vector spaces 
to Banach spaces},
Israel J. Math.\ {\bf 155} (2006), 205--252.
%
\bibitem[Gl06b]{Gl06b}
---, {\it Discontinuous non-linear mappings on locally convex
direct limits}, Publ.\ Math.\ Debrecen {\bf 69} (2006), 1--13. 
%
\bibitem[Gl06c]{Gl06c}
 ---, {\it Fundamental problems in the theory of infinite-dimensional Lie groups},
J. Geom.\ Symm.\ Phys.\ {\bf 5} (2006), 24--35.
%
\bibitem[Gl06d]{Gl06d}
---, \emph{Implicit functions from topological vector spaces to Fr\'{e}chet spaces in the presence of metric estimates},
preprint, 2006, arXiv:math/0612673.
%
\bibitem[Gl07a]{Gl07a}
---, {\it Direct limits of infinite-dimensional Lie groups
compared to direct limits in related categories}, 
J. Funct.\ Anal.\ {\bf 245} (2007), 19--61. 
%
\bibitem[Gl07b]{Gl07b}
 ---, {\it Direct 
limit groups do not have small subgroups}, 
Topol.\ Appl.\ {\bf 154} (2007), 1126--1133. 
%
\bibitem[Gl07c]{Gl07c}
---, \emph{Instructive examples of smooth, complex differentiable and complex
analytic mappings into locally convex spaces}, J. Math.\ Kyoto Univ.\ {\bf 47}
(2007), 631--642. 
%
\bibitem[Gl07d]{Gl07d}
---, \emph{Finite order differentiability properties, fixed points and implicit functions over valued fields},
preprint, 2007, math/0511218.
%
\bibitem[Gl07e]{Gl07e}
---, \emph{Aspects of differential calculus related to infinite-dimensional vector bundles and
Poisson vector spaces}, manuscript in preparation, 2007.
%
\bibitem[Gl08a]{Gl08a}
---, \emph{Solutions to open problems in Neeb's recent survey on infinite-dimensional Lie groups},
Geom.\ Dedicata {\bf 135} (2008), 71--86.
%
\bibitem[Gl08b]{Gl08b}
---, \emph{Applications of hypocontinuous bilinear maps in infinite-dimensional differential calculus},
pp.~171--186 in: S. Silvestrov, E. Paal, V. Abramov and A. Stolin (eds),
``Generalized Lie Theory in Mathematics, Physics and Beyond,'' Springer Verlag,
2008.
%
\bibitem[Gl08c]{Gl08c}
---, \emph{Homotopy groups of ascending unions of infinite-dimensional manifolds},
to appear in Ann.\ Inst.\ Fourier (cf.\ arXiv:0812.4713v1)
%
\bibitem[Gl11]{Gl11}
---, \emph{Direct limits of infinite-dimensional Lie groups}, pp.~243--280
in: K.-H. Neeb and A. Pianzola (eds), ``Developments and Trends in Infinite-Dimensional Lie Theory,'' Birkh\"{a}user Progr.\ Math.\ 288, 2011.
%
\bibitem[Gl12a]{Gl12a}
---, \emph{Continuity of bilinear maps on direct sums of topological vector spaces},
J. Funct.\ Anal.\ {\bf 262:5} (2012), 2013--2030.
%
\bibitem[Gl12b]{Gl12b} ---,  {\it Regularity properties of infinite-dimensional Lie groups, and semiregularity}, 
arXiv:math.FA:1208.0715 (Version from Feb.~6, 2016)
%
\bibitem[Gl13a]{Gl13a}
---, {\it Continuity of LF-algebra representations associated to representations of Lie
groups}, Kyoto J. Math.\ {\bf 53:3} (2013), 567--595. 
%
\bibitem[Gl13b]{Gl13b}
---, {\it Differentiable mappings between spaces of sections}, 
preprint, 2013, arXiv:1308.1172.
%
\bibitem[Gl15a]{Gl15a}
---, {\it Fundamentals of submersions and immersions between infinite-dimensional manifolds},
preprint, 2015, arXiv:1502.05795.
%
\bibitem[Gl15b]{Gl15b}
---, {\it Measurable regularity properties of infinite-dimensional Lie groups},
preprint, 2015, arXiv:1601.02568.
%
\bibitem[Gl17]{Gl17}
---, {\it Completeness of locally $k_\omega$-groups and related infinite-dimensional
Lie groups}, Topol.\ Appl.\ {\bf 228} (2017), 277--284.
%
\bibitem[Gl19a]{Gl19a}
---, {\it Completeness of infinite-dimensional Lie groups in their left uniformity},
Canadian J. Math.\ {\bf 71} (2019),
131--152.
%
\bibitem[Gl20a]{Gl20a}
---,
\emph{Smoothing operators for vector-valued functions and extension operators},
preprint, 2020, arXiv:2006.00254.
%
\bibitem[Gl20b]{Gl20b}
---,
\emph{Diffeomorphism groups of real-analytic manifolds are $L^1$-regular},
preprint, 2020, arXiv:2007.15611.
%
\bibitem[Gl21a]{Gl21a}
---,
\emph{Direct limits of regular Lie groups}, Math.\ Nachr.\ {\bf 294} (2021),
74--81.
%
\bibitem[Gl21b]{Gl21b}
---,
\emph{Manifolds of mappings on rough manifolds and manifold
structures on box products}, in preparation.
%
\bibitem[Gl21c]{Gl21c}
---,
\emph{Semidirect products involving Lie groups of compactly supported
diffeomorphisms}, in preparation.
%
\bibitem[Gl23]{Gl23}
---,
{\it Diffeomorphism groups of convex polytopes},
J. Convex Anal.\ {\bf 30} (2023), 343--358.
%
\bibitem[GHK10]{GHK10}
Gl\"{o}ckner, H., T. Hartnick, and R. K\"{o}hn,
\emph{Final group topologies, Kac--Moody groups and Pontryagin duality},
Israel J. Math.\ {\bf 177} (2010), 49--102.
%
\bibitem[GH23]{GH23}
Gl\"{o}ckner, H. and J. Hilgert,
\emph{Aspects of control theory on infinite-dimensional Lie groups and $G$-manifolds},
J. Differ.\ Equations {\bf 343} (2023), 186--232. 
%
\bibitem[GL12]{GL12}
Gl\"{o}ckner, H., and B. Langkamp,
\emph{Topological algebras of rapidly decreasing matrices and generalizations},
Topol.\ Appl.\ {\bf 159} (2012), 2420--2422. 
%
%
\bibitem[GM20]{GM20}
Gl\"{o}ckner, H., and N. Masbough,
{\it Products of regular locally compact spaces are $k_\R$-spaces},
Top.\ Proc.\ {\bf 55} (2020), 35--38.
%
\bibitem[GN03]{GN03}
Gl\"{o}ckner, H., and K.-H. Neeb, {\it Banach-Lie quotients, enlargibility,
and universal complexifications}, 
J. Reine Angew.\ Math.\ \textbf{560} (2003), 1--28. 
%
\bibitem[GN12]{GN12}
---, {\it When unit groups of continuous inverse algebras are regular 
Lie groups}, Studia Math. {\bf 211:2} (2012), 95--109.
%
\bibitem[GN17]{GN17}
---, {\it Diffeomorphism groups of compact convex sets},
Indag.\ Math.\ {\bf 28} (2017), 760--783. 
%
\bibitem[GS22]{GS22}
Gl\"{o}ckner, H. and A. Schmeding,
\emph{Manifolds of mappings on Cartesian products},
Ann.\ Global Anal.\ Geom.\ {\bf 61} (2022), 359--398. 
%
\bibitem[GSS24]{GSS23}
Gl\"{o}ckner, H., A. Schmeding, and A. Suri,
\emph{Manifolds of continuous BV-functions and vector measure regularity of Banach--Lie groups}
preprint, 2024, arXiv:2407.05190.
%
\bibitem[GSu23]{GSu23}
Gl\"{o}ckner, H. and A. Suri,
\emph{Strong ILB-Lie groups are $L^1$-regular},
manuscript, 2023.
%
\bibitem[God71]{God71}
Godbillon, C., ``El\'ements de Topologie Alg\'ebrique'', Hermann, Paris, 
1971. 
%
\bibitem[Gd58]{Gd58} Godement, R., ``Topologie alg\'{e}brique et th\'{e}orie des faisceaux'',
Hermann, Paris, 1958.
%
\bibitem[Go04]{Go04} Goldin, G. A., {\it Lectures on 
diffeomorphism groups in quantum physics},
Contemporary problems in mathematical physics, 3--93, World Sci.\ Publ.,
Hackensack, NJ, 2004. 
%
\bibitem[GW84]{GW84}
Goodman, R., and N. R. Wallach, {\it Structure and unitary
cocycle representations of loop groups and the group of
diffeomorphisms of the circle}, J. Reine Ang.\ Math.\ \textbf{347} (1984), 
69--133.
%
\bibitem[GW85]{GW85}
---, {\it Projective unitary positive energy
representations of} $Diff(\bS^1)$, J. Funct.\ Anal.\ \textbf{63} (1985), 
299--321. 
%
\bibitem[GOV97]{GOV97} Gorbatsevich, V.\, V., A.L.\, Onishchik, 
and E.B.\, Vinberg, ``Foundations of Lie Theory and Lie Transformation 
Groups,'' Springer, 1997 
%
\bibitem[Go69]{Go69} Goto, M., {\it On an arcwise connected subgroup of a Lie group}, Proc.\ Amer.\ 
Math.\ Soc.\ {\bf 20} (1969), 157--162.   
%
\bibitem[Grb88]{Gr88} Grabowski, J., {\it Free subgroups of diffeomorphism groups}, 
Fund.\ Math.\ {\bf 131:2} (1988), 103--121.
%
\bibitem[Grb93]{Gr93}
Grabowski, J., {\it Derivative of the exponential mapping for
infinite-dimensional Lie groups}, 
Ann.\ Global Anal.\ Geom.\ \textbf{11:3} (1993), 213--220. 
%
\bibitem[GVF01]{GVF01} Gracia-Bondia, J. M., J.C. Vasilly, and H. Figueroa,
``Elements of Non-com\-mu\-ta\-tive Geometry,'' Birkh\"auser Advanced Texts,
Birkh\"auser Verlag, Basel, 2001.
%
\bibitem[Grv50]{Grv50}
Graev, M. I., {\it Theory of topological groups.} I.\ {\it Norms and metrics on groups.
Complete groups. Free topological groups},
Usp.\ Mat.\ Nauk\ N.S. 5, {\bf 2:36} (1950), 3--56.
%
\bibitem[Gra84]{Gra84} Gramsch, B., 
{\it Relative Inversion in der St\"{o}rungstheorie
von Operatoren und $\Psi$-Algebren}, Math.\ Ann.\ \textbf{269} (1984), 27--71.
%
\bibitem[Gr57]{Gr57} Grauert, H., {\it Holomorphe 
Funktionen mit Werten in komplexen Lieschen
Gruppen}, Math.\ Ann.\ {\bf 133} (1957), 450--472. 
%
\bibitem[Gr58a]{Gr58a} ---, {\it Analytische 
Faserungen \"uber holomorph-vollst\"andigen R\"aumen}, 
Math.\ Ann.\ {\bf 135} (1958), 263--273. 
%
\bibitem[Gr58b]{Gr58b}
---, {\it On Levi’s problem and the imbedding of real-analytic manifolds},
Ann.\ Math.\ (2) {\bf 68} (1958), 460--472.
%
\bibitem[GW72]{GW72}
Greenfield, S.~J., and N. R. Wallach, {\it Automorphism 
groups in Banach spaces}, 
Trans.\ Amer.\ Math.\ Soc.\ \textbf{166}  (1972), 45--57. 
%
\bibitem[GE92]{GE92}
Gro\ss{}e-Erdmann, K.-G., {\it The Borel-Okada theorem revisited}, 
Habilitationsschrift,
Fernuniversit\"{a}t Hagen,  1992. 
%
\bibitem[GE04]{GE04}
---, {\it A weak criterion for vector-valued holomorphy},
Math.\ Proc.\ Camb.\ Philos.\ Soc.\ {\bf 136:2} (2004), 399--411.
%
\bibitem[Gr54]{Gr54}
Grothendieck, A.,
{\it Sur les espaces} (F) {\it et} (DF),
Summa Brasil.\ Math.\ {\bf 3} (1954), 57--121.
%
\bibitem[Gr55]{Gr55}
---, ``Produits tensoriels topologiques et espaces
nucl\'eaires'', Mem.\ Amer.\ Math.\ Soc.\ \textbf{16}, 1955.
%
\bibitem[Gr73]{Gr73}
---, ``Topological Vector Spaces,'' Gordon and Breach,
New York, 1973.

\bibitem[Gu08]{Gu08} G\"undogan, H., {\it Lie algebras of Smooth Sections}, 
Preprint,  arXiv:math.RA:0803.2800v4 

\bibitem[Gu11]{Gu11} G\"undogan, H., ``Classification and Structure Theory
of Lie Algebras of Smooth Sections,'' 
PhD Thesis, FAU Erlangen--Nuremberg, October 2011 
%
\bibitem[GR07]{GR07}
Guieu, L., and C. Roger, ``L'Alg\`ebre et le Groupe de Virasoro,''
Les Publications CRM, Montreal, QC, 2007. 
%
\bibitem[Gu80]{Gu80}
Gurarie, D., {\it Banach uniformly continuous representations
of Lie groups and algebras}, J.\ Funct.\ Anal.\ \textbf{36} (1980), 401--407.
%
\bibitem[Gu77]{Gu77} Gutknecht, J., {\it Die $C^\infty_\Gamma$-Struktur auf der 
Diffeomorphismengruppe einer kompakten Mannigfaltigkeit}, 
Ph.D. Thesis, Eidgen\"ossische Technische Hochschule Z\"urich, Diss.\ No.\ {\bf 5879}, 
Juris Druck + Verlag, Zurich, 1977. 
%
\bibitem[HV04]{HV04} 
Haller, S., and C. Vizman, {\it Nonlinear Grassmannians as coadjoint orbits},   
Math.\ Ann.\ \textbf{329} (2004), 771--785. 
%
\bibitem[Ham82]{Ham82}
Hamilton, R., {\it The inverse function theorem of Nash and
Moser},  Bull.\ Amer.\ Math.\ Soc.\ \textbf{7} (1982), 65--222. 
%
\bibitem[Han71]{Han71}
Hansen, V.~L., {\it Some theorems on direct limits
of expanding systems of manifolds}, Math.\ Scand.\ \textbf{29} (1971), 5--36. 
%
\bibitem[Han83]{Han83} ---, {\it The homotopy groups of a space of maps
between oriented closed surfaces}, Bull.\ London Math.\ Soc.\ {\bf 15}  (1983), 
360--364.
%
\bibitem[Hn19]{Hn19}
Hanusch, M.,
{\it Differentiability of the evolution map and Mackey continuity},
Forum Math.\ {\bf 31} (2019), 1139--1177.

\bibitem[Hn20]{Hn20}
---, {\it The strong Trotter property for locally $\mu$-convex groups},
J. Lie Theory {\bf 30} (2020), 25--32. 
%
\bibitem[Hn22]{Hn22}
---, {\it Regularity of Lie groups},
Commun.\ Anal.\ Geom.\ \textbf{30} (2022), 53--152.
%
%
\bibitem[Hn23]{Hn23}  
---,
{\it A $C^k$-Seeley-Extension-Theorem for
Bastiani's differential calculus},
Can.\ J. Math.\ \textbf{75} (2023), 170--201. 
%
\bibitem[dlH72]{dlH72}
de la Harpe, P., ``Classical Banach--Lie Algebras and
Banach--Lie Groups of Operators in Hilbert Space'', 
Lecture Notes in Math.\ \textbf{285}, Springer-Verlag, Berlin, 1972. 
%
\bibitem[dlH79]{dlH79}
---, {\it Les extensions de $\gl(H)$ par un noyau de dimension finie 
sont triviales}, J.\ Funct.\ Anal.\ \textbf{33}  (1972), 362--373. 
%
\bibitem[dlH82]{dlH82} ---, {\it Classical groups and classical Lie algebras of operators}, 
Proc.\ Symp.\ Pure Math.\ {\bf 38:1} (1982), 477--513. 
%
\bibitem[Ha74]{Ha74}
Harris, L.~A., {\it Bounded symmetric homogeneous domains in
infinite dimensional spaces}, pp.~13--38 in: T.~L. Hayden and T.~J. Suffridge (eds),
``Proceedings on Infinite Dimensional
Holomorphy'',  Lecture Notes in Math.\ \textbf{364},
Springer-Verlag, Berlin, 1974.
%
\bibitem[HK77]{HK77} Harris, L. A., and W. Kaup,
{\it Linear algebraic groups in infinite dimensions}, 
Illinois J. Math.\ {\bf 21} (1977), 666--674.
%
\bibitem[Ha81]{Ha81} Hatcher, A., {\it 
On the diffeomorphism group of $\bS^1 \times \bS^2$}, 
Proc. Amer. Math. Soc. {\bf 83:2} (1981), 427--430
%
\bibitem[Ha83]{Ha83} ---, {\it 
A proof of the Smale conjecture, $\Diff(\bS^3) \cong \OO(4)$}, 
Ann. of Math. (2) {\bf 117:3} (1983), 553--607
%
\bibitem[Ha02]{Ha02} ---, ``Algebraic Topology,''
Cambridge University Press, 2002.
%
\bibitem[Ha12]{Ha12} Hatcher, A., ``A 50-year view of diffeomorphism groups,'' 
A talk at the 50th Cornell Topology Festival in May 2012; 
http://pi.math.cornell.edu/~hatcher/

\bibitem[Hau06]{Hau06} Hausdorff, F., {\it Die symbolische Exponentialformel in der 
Gruppentheorie}, Leip\-zi\-ger Berichte {\bf 58} (1906), 19--48.  
%
\bibitem[HS68]{HS68} Hausner, M., and J. T. Schwartz, ``Lie Groups; Lie Algebras,''
Gordon and Breach, New York, London, Paris, 1968.  
%
\bibitem[HeMa02]{HeMa02} Hector, G., and E. Mac\'{\i}as-Virg\'{o}s, {\it Diffeological groups},
Research Exp.\ Math.\ {\bf 25} (2002), 247--260. 
%
\bibitem[Hel93]{Hel93} Helemskii, A. Ya., ``Banach and Locally Convex Algebras,'' 
Oxford Science Publications,  Oxford University Press, New York, 1993. 
%
\bibitem[Hel78]{Hel78}
Helgason, S., ``Differential Geometry, Lie Groups, and
Symmetric Spa\-ces'', Acad.\ Press, London, 1978. 
%
\bibitem[He82]{He82}
Helton, W., {\it Non-euclidean functional analysis and
electronics}, Bull.\ Amer.\ Math.\ Soc.\ {\bf 7}  (1982), 1--63. 
%
\bibitem[He71]{He71} Herv\'{e}, M., 
``Analytic and Plurisubharmonic Functions'', 
Lecture Notes in Math.\ \textbf{198}, Springer-Verlag, Berlin, 1971. 
%
\bibitem[He89]{He89} ---, ``Analyticity in Infinite Dimensional Spaces'', 
de Gruyter, Berlin,  1989.
%
\bibitem[HR79]{HR79}
Hewitt, E., and K. A. Ross, ``Abstract Harmonic Analysis I'', 
Springer-Verlag, New York, 1979. 
%
\bibitem[HG27]{HG27}
Hildebrandt, T.~H., and L. M.  Graves, {\it Implicit functions and their differentials
in general analysis}, Trans.\ Amer.\ Math.\ Soc.\ 
\textbf{29} (1927), 127--153. 
%
\bibitem[HiNe12]{HiNe12}
Hilgert, J., and K.-H. Neeb, ``Structure and Geometry of Lie Groups,'' 
Springer-Verlag, New York, 2012.
%
\bibitem[HP57]{HP57} Hille, E., and R. S. Phillips, ``Functional Analysis and 
Semigroups,'' Amer.\ Math.\ Soc., 1957. 
%
\bibitem[Hi99]{Hi99}
Hiltunen, S., {\it Implicit functions from locally convex
spaces to Banach spaces}, Studia Math.\ \textbf{134}  (1999), 235--250. 
%
\bibitem[Hi00]{Hi00} ---, {\it A Frobenius theorem for locally convex global
analysis}, Monatsh.\ Math.\ \textbf{129}  (2000), 109--117. 
%
\bibitem[HST01]{HST01}
Hirai, T., H. Shimomura,
N. Tatsuuma, and E. Hirai, {\it Inductive limits
of topologies, their direct products, and problems
related to algebraic structures}, 
J. Math.\ Kyoto Univ.\ \textbf{41} (2001), 475--505. 
%
\bibitem[Hr76]{Hr76} Hirsch, M.~W.. ``Differential Topology'', Grad.\ Texts Math.\ 
{\bf 33}, Springer-Verlag, New York, 1976.
%
\bibitem[HrP70]{HrP70}
Hirsch, M. W., and C.~C. Pugh,
{\it Stable manifolds and hyperbolic sets}, pp.~133--163 in:
S.-S. Chern and S. Smale (eds),
``Global Analysis,''
Proc.\ Sympos.\ Pure Math.\ {\bf 14}, 1970.
%
\bibitem[Hz91]{Hz91} Hirzebruch, F., {\it Division algebras and topology}, 
pp.~281--302 in: ``Numbers,'' Graduate Texts in Mathematics {\bf 123}, 
Springer-Verlag, New York, 1991. 
%
\bibitem[HS17]{HS17}
Hjelle, E. O. and A. Schmeding,
\emph{Strong topologies for spaces of smooth maps with infinite-dimensional target},
Expo.\ Math.\ {\bf 35} (2017), 13--53.
%
\bibitem[Ho51]{Ho51}
Hochschild, G., {\it Group extensions of Lie groups I, II},  
Ann.\ of Math.\ \textbf{54} (1951), 96--109; 537--551. 
%
\bibitem[Ho65]{Ho65} ---, ``The Structure of Lie Groups,'' Holden Day, San 
Francisco, 1965.
%
\bibitem[HS53a]{HS53a} 
Hochschild, G., and J.-P. Serre, {\it Cohomology of group
extensions}, Trans.\ Amer.\ Math.\ Soc.\ \textbf{74} (1953), 
110--134.
%
\bibitem[HS53b]{HS53b}
---, {\it Cohomology of Lie
algebras}, Ann.\ of Math.\ \textbf{57:2} (1953), 591--603. 
%
\bibitem[Hof68]{Hof68} Hofmann, K. H., ``Introduction to the Theory of Compact Groups. Part I,'' 
Dept.\ Math.\ Tulane Univ., New Orleans, LA, 1968. 
%
\bibitem[Hof72]{Hof72} ---, {\it Die Formel von Campbell, Hausdorff und Dynkin 
und die Definition Liescher Gruppen},
in ``Theory 
of Sets and Topology'' (in honour of Felix Hausdorff, 1868--1942); 
251--264. VEB Deutsch, Verlag Wissensch., Berlin, 1972. 
%
\bibitem[Hof75]{Hof75} ---, {\it Analytic groups without analysis}, Symposia Mathematica {\bf 16} 
(Convegno sui Gruppi Topologici e Gruppi di Lie, INDAM, Rome, 1974); 357--374, 
Academic Press, London, 1975. 
%
%
\bibitem[HoM98]{HoM98} Hofmann, K.~H., and S. A. Morris, ``The Structure
of Compact Groups'', Studies in Math.\ \textbf{25},
de Gruyter, Berlin, New York, 1998. 
%
\bibitem[HoM03]{HoM03} ---, {\it Projective limits of 
finite-dimensional Lie groups},
Proc.\ London Math.\ Soc.\ (3) {\bf 87:3} (2003), 647--676.
%
\bibitem[HoM04a]{HoM04a} ---, {\it Lie theory and the structure of pro-Lie groups and pro-Lie algebras}, 
Top.\ Proc.\ {\bf 28} (2004), 541--567. 
%
\bibitem[HoM04b]{HoM04b} ---, {\it The structure of abelian pro-Lie groups}, Math.\ Z. {\bf 248} (2004), 
867--891.
%
\bibitem[HoM05]{HoM05} ---, {\it Sophus Lie's third fundamental theorem and the adjoint functor theorem}, 
J. Group Theory {\bf 8} (2005), 115--123. 
%
\bibitem[HoM07]{HoM07}
---, ``The Lie Theory of Connected Pro-Lie Groups,''
EMS Tracts in Mathematics {\bf 2},
European Mathematical Society, Z\"urich, 2007.
%
\bibitem[HMP04]{HMP04} Hofmann, K. H., S. A. Morris, and D. Poguntke, {\it 
The exponential function of locally connected compact abelian groups}, 
Forum Math.\ {\bf 16:1} (2004), 1--16.
%
\bibitem[HN77]{HN77} Hogbe-Nlend, H., ``Bornologies and Functional
Analysis'', North-Holland, Amsterdam, 1977. 
%
\bibitem[Hop42]{Hop42} Hopf, H., {\it Fundamentalgruppe und zweite Bettische Gruppe}, 
Comment.\ Math.\ Helv.\ {\bf 14} (1942), 257--309.
%
\bibitem[vHo52a]{vHo52a} van Hove, L., {\it Topologie des espaces fonctionnels analytiques, et des groups 
infinis des transformations}, Bull.\ de la Classe de Sc., Acad.\ Roy.\ de Belgique, 
S\'er. 5, {\bf 38} (1952), 333--351.
%
\bibitem[vHo52b]{vHo52b} ---, {\it L'ensemble des fonctions analytiques sur un compact en tant 
qu'alg\`ebre topologique}, Bull.\ Soc.\ Math.\ Belgique 1952, (1953), 8--17. 
%
\bibitem[HM74]{HM74}
Hunt, D. C., and S. A. Morris, {\it Free subgroups of free topological groups},
pp.~377--387 in:
M.~F. Newman (ed),
``Proceedings of the 2nd international conference on the theory of groups. Australian National University, August 13--24, 1973,''
Lecture Notes Math.\ {\bf 372},
Springer-Verlag, Berlin, 1974.
%
\bibitem[Hu71]{Hu71}
Hu\v{s}ek, M.,
{\it Products of quotients and of $k'$-spaces},
Comment.\ Math.\ Univ.\ Carolinae {\bf 12:1} (1971), 61--68.
%
\bibitem[Hu94]{Hu94} Husemoller, D., ``Fibre Bundles,'' 
Graduate Texts in Math., Springer-Verlag, New York, 1994. 
%
\bibitem[Il03]{Il03}
Illman, S., \emph{The very-strong $C^\infty$-topology on $C^\infty(M,N)$
and $K$-equivariant maps},
Osaka J. Math.\ {\bf 40} (2003), 409--428. 
%
\bibitem[IKT13]{IKT13} Inci, H., T.~Kappeler, and P.~Topalov, 
``On the regularity of the composition of diffeomorphisms,'' 
Mem. Amer. Math. Soc. {\bf 226} (2013), no. 1062, vi+60 pp.

\bibitem[Ing19]{Ing19}
Ingrisch, S., ``Gemischte Differenzierbarkeit von Fixpunkten und impliziten Funktionen'',
Master's thesis, University of Paderborn, 2019 (advisor: H. Gl\"{o}ckner).
%
\bibitem[Ir72]{Ir72}
Irwin, M. C.,
\emph{On the smoothness of the composition map},
Q. J. Math.\ (1972), 113--133.
%
\bibitem[Ir80]{Ir80}
---,
\emph{A new proof of the pseudo-stable manifold theorem},
J. Lond.\ Math.\ Soc., II. Ser.\
{\bf 21}, (1980), 557--566.
%
\bibitem[Is90]{Is90} Ismagilov, R.~S., 
{\it Representations of certain infinite
dimensional groups}, in:
A.~I. Vershik and D.~P. Zhelobenko (eds), 
``Representations of Lie Groups and Related
Topics'', Gordan and
Breach Science Publ., New York, 1990.
%
\bibitem[Is96]{Is96} 
---, ``Representations of Infinite-Dimensional
Groups'',
Translations of Math.\ Monographs \textbf{152}, Amer.\ Math.\ Soc., 
1996.
%
\bibitem[Ja10]{Ja10} Janssens, B., ``Transformation \& Uncertainty. 
Some Thoughts on Quantum Probability Theory, Quantum Statistics, 
and Natural Bundles,'' PhD Thesis, 2010; arXiv:math.DG:1011.3035 

\bibitem[JLM17]{JLM17}
Jaramillo, J.~A., S.~Lajara, and \'{O}. Madiedo,
{\it Inversion of nonsmooth maps between Banach spaces},
preprint, 2017, arXiv:1712.00565.
%
\bibitem[Jr81]{Jr81}
Jarchow, H.,
``Locally Convex Spaces,'' B.~G. Teubner, Stuttgart, 1981. 

\bibitem[Jo82]{Jo82} Joris, H., {\it Une $C^\infty$-application non-immersive 
qui poss\`ede la propri\'et\'e universelle des immersions}, 
Arch. Math. (Basel) {\bf 39:3} (1982), 269--277

%
\bibitem[Ka85]{Ka85} Kac, V. G., {\it Constructing groups associated to infinite-dimensional 
Lie algebras},
in:  V.~G. Kac (ed), ``Infinite-Dimensional Groups with Applications,'' 
MSRI Publications {\bf 4}, 
Springer-Verlag, Berlin, Heidelberg, New York, 1985. 
%
\bibitem[Ka90]{Ka90} ---, ``Infinite-dimensional Lie Algebras", Cambridge University 
Press, 1990.
%
\bibitem[KP83]{KP83} Kac, V.~G., and D. H. Peterson,  {\it Regular functions on
certain infinite-dimensional groups},
in: M. Artin and J. Tate (eds)
``Arithmetic and Geometry, Vol.\ 2'', 
Birkh\"auser, Boston, 1983. 
%
\bibitem[KP84]{KP84}
---, {\it Unitary structure in
representations of infinite dimensional groups and a convexity
theorem}, Invent.\ Math.\ \textbf{76} (1984), 1--14. 
%
\bibitem[KP87]{KP87} ---, {\it On geometric invariant theory 
for infinite-di\-men\-sio\-nal groups}, pp.\ 109--142 in: 
Cohen A.~M. et al.\ (eds) 
``Algebraic Groups'' (Utrecht 1986),
Lecture Notes in Math.\ \textbf{1271}, Springer-Verlag, Berlin, 1987. 
%
\bibitem[KR87]{KR87}
Kac, V.~G., and A. K. Raina, ``Bombay Lectures on
Highest Weight Representations of 
Infinite Dimensional Lie Algebras'', Advanced Series in Math.\ Physics \textbf{2},
World Scientific, Singapore, 1987. 
%
\bibitem[Ka52]{Ka52}
Kadison, R.~V., {\it Infinite unitary groups}, 
Trans.\ Amer.\ Math.\ Soc.\ \textbf{72} (1952), 386--399. 
%
\bibitem[Ka54]{Ka54}
---, {\it Infinite general linear groups}, 
Trans.\ Amer.\ Math.\ Soc.\ \textbf{76} (1954), 66--91. 
%
\bibitem[Ka55]{Ka55}
---, {\it The general linear group of infinite factors}, 
Duke Math.\ J. \textbf{22} (1955), 110--122. 
%
\bibitem[KK63]{KK63}
Kakutani, S., and V. Klee, {\it The finite topology
of a linear space}, Arch.\ Math.\ \textbf{14} (1963), 55--58.  
%
\bibitem[Kal86]{Kal86} Kallmann, R. R., {\it Uniqueness results 
for homeomorphism groups}, Transactions of the Amer. Math. Soc. 
{\bf 295:1} (1986), 389--396
%
\bibitem[KR01]{KR01}
Kamran, N., and T. Robart, {\it A manifold structure 
for analytic Lie pseudogroups of infinite type}, 
J.~Lie Theory \textbf{11} (2001), 57--80. 
%
\bibitem[KR04]{KR04} ---, {\it An infinite-dimensional manifold structure 
for analytic Lie pseudogroups of infinite type}, 
Internat.\ Math.\ Res.\ Notices {\bf 34} (2004), 1761--1783. 
%
\bibitem[KA64]{KA64}
Kantorivich, L. ~V., and G.~P. Akilov,
``Functional Analysis in Normed Spaces,'' Pergamon Press, Oxford, 1964.
%
\bibitem[Kp52]{Kp52}
Kaplan, S., {\it Cartesian products of reals}, 
Amer.\ J.~Math.\ \textbf{74} (1952), 936--954. 
%
\bibitem[Ka78]{Ka78} 
Karoubi, M., ``$K$-theory'',
Grundlehren math.\ Wiss.\ {\bf 226}, Springer-Verlag, Berlin, 1978.
%
\bibitem[Ka81]{Ka81}
Kaup, W., {\it \"Uber die Klassifikation
der symmetrischen
her\-mi\-te\-schen Mannigfaltigkeiten unendlicher Dimension I,II}, 
Math.\ Ann.\ \textbf{257} (1981), 463--486; \textbf{262} (1983), 57--75. 
%
\bibitem[Ka83]{Ka83}
---, {\it A Riemann mapping theorem for bounded symmetric
domains in complex Banach spaces}, 
Math.\ Z. \textbf{183} (1983), 503--529. 
%
\bibitem[KR07]{KR07} Kechris, A.S., and C. Rosendal, 
{\it Turbulence, amalgamation, and generic automorphisms of 
homogeneous structures}, 
Proc. Lond. Math. Soc. (3) {\bf 94:2} (2007), 302--350

\bibitem[KKM06]{KKM06}
K\c{e}dra, J., D. Kotchick, and S. Morita, {\it Crossed flux homomorphisms 
and vanishing theorems for flux groups},
Geom.\ Funct.\ Anal.\ {\bf 16} (2006), 1246--1273. 
%
\bibitem[Kr74]{Kr74} Keller, H.~H., ``Differential Calculus
in Locally Convex Spaces'', Springer-Verlag,
Berlin, Heidelberg, New York, 1974.
%
\bibitem[Ke73]{Ke73} 
Keesling, J., {\it Topological groups whose underlying spaces
are separable Fr\'echet manifolds},  Pacific J. Math.\ {\bf 44:1}  (1973), 181--189. 
%
\bibitem[Kel75]{Kel75}
Kelley, J.~L., ``General Topology,''
Graduate Texts in Math.\ {\bf 27}, Springer-Verlag,
New York, ${}^2$1975. 
%
\bibitem[KM63]{KM63} Kervaire, M. A., and J. W. Milnor, {\it 
Groups of homotopy spheres. I}, Ann.\ of Math.\ (2) {\bf 77} (1963), 504--537.
%
\bibitem[KW09]{KW09}
Khesin, B., and R. Wendt, ``The Geometry of Infinite-Dimensional Groups,''
Springer-Verlag, Berlin, 2009.
%
\bibitem[KMM24]{KMM24} Khesin, B., G. Misiolek, and K. Modin,
  {\it Information geometry of diffeomorphism groups},
  Preprint, arXiv:2411.03265v1 [math.DG]

  
\bibitem[Kh83]{Kh83} Kholina, G. G., {\it One-parameter subgroups 
of differentiable locally convex groupuscles}, 
Mat.\ Zametki {\bf 33:6} (1983), 823--832. 
%
\bibitem[Ki87]{Ki87} Kirillov, A., {\it K\"ahler structure on the $k$-orbits of the group 
of diffeomorphisms of a circle}, Funct.\ Anal.\ Appl.\ {\bf 21} (1987), 122--125.  
%
\bibitem[Ki05]{Ki05} ---, {\it The orbit method beyond Lie groups. Infinite-dimensional 
groups},
Surveys in modern mathematics, 292--304; 
London Math.\ Soc.\ Lecture Note Ser.\ {\bf 321}, 
Cambridge Univ.\ Press, Cambridge, 2005. 
%
\bibitem[KY87]{KY87} Kirillov, A.~A., and D. V. Yuriev, 
{\it K\"ahler geometry of
the infinite-dimensional homogeneous space} $M =
\Diff_+(\bS^1)/\Rot(\bS^1)$, Funct.\ Anal.\ Appl.\ \textbf{21}  (1987), 284--294. 
%
\bibitem[Kir01]{Kir01}
Kirwan, P.,
{\it Complexification of multilinear mappings and polynomials},
Math.\ Nachr.\ {\bf 231} (2001), 39--68.
%
\bibitem[Kl53]{Kl53}
Klee, V.~L., {\it Convex bodies and periodic heomeomorphisms in Hilbert space},
Trans.\ Amer.\ Math.\ Soc.\ {\bf 74} (1953), 10--43.
%
\bibitem[Kl1872]{Kl1872} Klein, F., 
``Vergleichende Betrachtungen über neuere geometrische Forschungen,'' 
('A comparative review of recent researches in geometry'), Math. Annalen {\bf 43} 
 (1893), 63--100

\bibitem[Kl11]{Kl11} Klotz, M., {\it The automorphism group of a Banach 
principal bundle with {$\1$}-structure}, Geom.\ Dedicata {\bf 154} (2011), 161--182. 
%
\bibitem[Kl12]{Kl12} ---, {\it An integrability criterion for Banach--Lie triple 
systems}, J. Lie Theory {\bf 22:1} (2012), 205--244.
%
\bibitem[Kn30]{Kn30} Kneser, H., {\it Die kanonische Parametergruppe},
  Jahresber. Deutsch. Math. Verein {\bf 39} (1930), 72--78
  
\bibitem[KYM85]{KYM85} Kobayashi, O., 
A. Yoshioka, Y. Maeda, and H. Omori, {\it The 
theory of infinite-dimensional Lie groups and its applications}, 
Acta Appl.\ Math.\ {\bf 3:1}  (1985), 71--106. 
%
\bibitem[Ko95]{Ko95} Kobayashi, S., ``Transformation Groups in Differential 
Geometry,'' Reprint of the 1972 edition, 
Classics in Mathematics, Springer-Verlag, Berlin, 1995.
%
\bibitem[Ko69]{Ko69}
Koecher, M., {\it Gruppen und Lie-Algebren von rationalen
Funktionen}, Math.\ Z. \textbf{109} (1969), 349--392. 
%
\bibitem[K\"o69]{Koe69} K\"othe, G., ``Topological Vector Spaces I,'' Grundlehren der
Math. Wis\-sen\-schaf\-ten {\bf 159}, Springer-Verlag, Berlin etc., 1969. 
%
\bibitem[Kop70]{Kop70} Kopell, N., {\it Commuting diffeomorphisms},
Proc.\ Symp.\ Pure Math.\ 
{\bf 14}, Amer.\ Math.\ Soc., 1970, 165--184. 
%
\bibitem[KS77]{KS77}
Kolomytsev, V.~I., and Y. S. Samoilenko, {\it Irreducible
representations of inductive limits of groups}, Ukrain.\ Mat.\ Zh.\
\textbf{29:4}  (1977), 526--531; Engl.\ transl., 402--405.
%
\bibitem[Ko64]{Ko64} Komatsu, H., {\it Semigroups of operators in locally 
convex spaces}, J.~Math.\ Soc.\ Japan {\bf 16:3} (1964), 230--262.
%
\bibitem[Ko68]{Ko68} Komura, T., {\it Semigroups of operators in locally 
convex spaces}, J.~Funct.\ Anal.\ {\bf 2} (1968), 258--296.
%
\bibitem[Kos70]{Kos70} Kostant, B., {\it Quantization and unitary representations}, 
in: ``Lectures in Modern Analysis and Applications III'',  
Springer Lecture Notes Math.\ {\bf 170} (1970), 87--208. 
%
\bibitem[KP02a]{KP02a} Krantz, S.~G., and H. R. Parks, ``The Implicit Function Theorem'',
Birkh\"{a}user-Verlag, Basel, 2002.
%
\bibitem[KP02b]{KP02b}
---,
``A Primer of Real Analytic Functions,''
Birkh\"{a}user, Basel, ${}^2$2002.
%
\bibitem[Kri80]{Kri80}
Kriegl, A., ``Eine Theorie glatter Mannigfaltigkeiten
und Vektorb\"{u}ndel'',
Doctoral Dissertation, University of Vienna, 1980. 
%
\bibitem[KM97]{KM97}
Kriegl, A., and P. W. Michor, `` The Convenient Setting of
Global Analysis'', Math.\ Surveys and Monographs {\bf 53}, Amer.\
Math.\ Soc., Providence, 1997.
%
\bibitem[KM97b]{KM97b}
---, {\it  Regular infinite-dimensional Lie groups}, 
J. Lie Theory \textbf{7}  (1997), 61--99. 
%
\bibitem[KMR16]{KMR16}
Kriegl, A., P. W. Michor, and A. Rainer,
{\it The exponential law for spaces of test functions and diffeomorphism groups},
Indag.\ Math.\ {\bf 27:1} (2016), 225--265.
%
\bibitem[Kr76]{Kr76}
Krikorian, N., {\it Invariance of domain in Banach spaces},
Oric.\ Amer.\ Math.\ Soc.\ {\bf 60} (1976), 367--368.
%
\bibitem[Ku65]{Ku65}
Kuiper, N.~H., {\it The homotopy type of the unitary group of
Hilbert spaces}, Topology \textbf{3} (1965), 19--30. 
%
\bibitem[Kum02]{Kum02} Kumar, S., ``Kac-Moody Groups, their Flag Varieties and Representation Theory,'' 
Progress in Math. {\bf 204}, Birkh\"auser, Boston, 2002.
%
\bibitem[Ku18]{Ku18} Kupers, A., 
``Lectures on Diffeomorphism Groups of Manifolds,'' 
Lectures Notes, Version Feb 22, 2019 \\
https://people.math.harvard.edu/~kupers/teaching/272x/book.pdf
%
\bibitem[Kur59]{Kur59} Kuranishi, M., 
{\it On the local theory of continuous infinite pseudo groups I}, 
Nagoya Math. J. {\bf 15} (1959), 225--260 

\bibitem[LMP98]{LMP98} Lalonde, F., D. McDuff, and L. Polterovich, 
{\it On the flux conjectures},
in: ``Geometry, topology, and dynamics'' 
(Montreal, PQ, 1995); 69--85, CRM Proc.\ Lecture Notes {\bf  15}, 
Amer.\ Math.\ Soc., Providence, RI, 1998.
%
\bibitem[La75]{La75}
Lang, S.,
``$\SL_2(\R)$,'' Addison-Wesley, Reading, Mass.\, 1975.
%
\bibitem[La93]{La93} 
Lang, S., ``Algebra'', Addison Wesley Publ.\ Comp., London, 1993. 
%
\bibitem[La99]{La99}
---, ``Fundamentals of Differential Geometry'', 
Grad.\ Texts Math.\ \textbf{191},
Springer-Verlag, Berlin, 1999. 
%
\bibitem[Lar99]{Lar99} Laredo, V.\ T., {\it Integration of unitary representations of
infinite dimensional Lie groups}, J.\ Funct.\ Anal.\ {\bf 161:2}
(1999), 478--508. 
%
\bibitem[Las57]{Las57} Lashof, R. K., {\it Lie algebras of locally compact groups}, Pacific J. Math.\ {\bf 7} 
(1957), 1145--1162.
%
\bibitem[Lars99]{Lars99} 
Larsson, T. A., {\it Lowest-energy representations of non-centrally extended
diffeomorphism algebras}, Commun.\ Math.\ Phys.\ \textbf{201} (1999), 461--470. 
%
\bibitem[Lau55]{Lau55} Laugwitz, D., {\it \"Uber unendliche kontinuierliche Gruppen.} I. {\it Grundlagen der 
Theorie; Untergruppen}, Math.\ Ann.\ {\bf 130} (1955), 337--350. 
%
\bibitem[Lau56]{Lau56} ---, {\it \"Uber unendliche kontinuierliche Gruppen.} II. 
{\it Strukturtheorie lokal Banachscher Gruppen}, Bayer.\ Akad.\ Wiss.\ Math.-Nat.\ 
Kl.\ S.-B. 1956, (1957), 261--286.
%
\bibitem[Lau73]{Lau73} ---, {\it Ist Differentialrechnung ohne
Grenzwertbegriff m\"oglich}\,?, 
Math.\ Phys.\ Semesterberichte \textbf{20} (1973), 182--201. 
%
\bibitem[LaTi66]{LaTi66} Lazard, M., and J. Tits, 
{\it Domaines d'injectivit\'e de l'application exponentielle}, Topology {\bf 4} (1966), 315--322.  
%
\bibitem[Lea61]{Lea61}
Leach, E.~B.,
{\it A note on inverse function theorems},
Proc.\ Am.\ Math.\ Soc.\ {\bf 12} (1961), 694--697.
%
\bibitem[Lec80]{Lec80} Lecomte, P., {\it Sur 
l'alg\`ebre de Lie des sections d'un fibr\'e en alg\`ebre de Lie}, 
Ann.\ Inst.\ Fourier {\bf 30} (1980), 35--50.  
%
\bibitem[Lec81a]{Lec81a} ---, {\it Sur l'alg\`ebre de Lie des sections d'un fibr\'e en alg\`ebre de Lie}, 
Ann.\ Inst.\ Fourier {\bf 30} (1980), 35--50. 
%
\bibitem[Lec81b]{Lec81b} ---, {\it On the infinitesimal automorphisms of a vector bundle}, 
J. Math.\ Pures Appl., IX. S\'er. {\bf 60} (1981), 229--239. 
%
\bibitem[Lec85]{Lec85} ---, {\it  Sur la suite exacte canonique associ\'ee \`a un fibr\'e 
principal}, Bull.\ Soc.\ Math.\ Fr.\ {\bf 13} (1985), 259--271. 
%
\bibitem[Lem95]{Lem95} Lempert, L., {\it The Virasoro group as a complex manifold}, 
Math.\ Res.\ Letters {\bf 2} (1995), 479--495. 
%
\bibitem[Lem97]{Lem97} ---, {\it The problem of complexifying
a Lie group}, pp.~169--176 in: P.~D. Cordaro et al.\ (eds),
``Multidimensional complex analysis and partial differential
equations'', Contemp.\ Math.\ {\bf 205}, American Math.\ Soc.,
Providence, 1997.
%
\bibitem[Les66]{Lec66} Leslie, J. A., {\it On a theorem of E. Cartan}, 
Ann.\ Mat.\ Pura Appl.\ (4) {\bf 74} (1966), 173--177. 
%
\bibitem[Les67]{Les67}
---, {\it On a differential structure
for the group of diffeomorphisms}, Topology~\textbf{6} (1967), 263--271. 
%
\bibitem[Les68]{Les68} ---, {\it Some Frobenius theorems in global analysis}, 
J. Diff.\ Geom.\ {\bf 2} (1968), 279--297.
%
\bibitem[Les82]{Les82} ---, {\it On the group of real analytic
diffeomorphisms of a compact real analytic manifold}, 
Trans.\ Amer.\ Math.\ Soc.\ {\bf 274} (1982), 651--669. 
%
\bibitem[Les83]{Les83} ---, {\it A Lie group structure for the group of analytic diffeomorphisms}, 
Boll.\ Un.\ Mat.\ Ital.\ A (6) {\bf 2:1} (1983), 29--37.
%
\bibitem[Les90]{Les90} ---, {\it A path functor for Kac-Moody Lie algebras},
in: ``Lie 
Theory, Differential Equations and Representation Theory (Montreal, PQ, 1989),''
265--270; Univ.\ Montreal, QC, 1990.
%
\bibitem[Les92]{Les92} ---, {\it Some integrable subalgebras of infinite-dimensional Lie groups}, 
Trans.\ Amer.\ Math.\ Soc.\ {\bf 333} (1992), 423--443. 
%
\bibitem[Les93]{Les93} ---, {\it On the integrability of some infinite dimensional Lie algebras}, 
Preprint, Howard University, 1993  
%
\bibitem[Les03]{Les03} ---, {\it On a diffeological group realization of
certain generalized symmetrizable Kac-Moody Lie
algebras}, J.~Lie Theory {\bf 13} (2003), 427--442.
%
\bibitem[Lew39]{Lew39} Lewis, D., {\it Formal power series transformations}, Duke Math.\ J. 
{\bf 5} (1939), 794--805. 
%
\bibitem[Li74]{Li74} 
Lichnerowicz, A., {\it Alg\`ebre de Lie des automorphismes infinit\'esimaux 
d'une structure unimodulaire}, Ann.\ Inst.\ Fourier {\bf 24}  (1974), 219--266. 
%
\bibitem[Lie80]{Lie80} Lie, S., {\it Theorie der Transformationsgruppen I},  Math.\ Ann.\ 
{\bf 16:4} (1880), 441--528.
%
\bibitem[Lie95]{Lie95} ---, {\it Unendliche kontinuierliche Gruppen}, Abhandlungen, S\"achsische 
Gesell\-schaft der Wissenschaften {\bf 21} (1895), 43--150.
%
\bibitem[LE93]{LE93} Lie, S., and F. Engel, 
``Theorie der Transformationsgruppen III,'' Teubner, Leipzig, 1893; 

\bibitem[Lod98]{Lod98}
Loday, J.-L., ``Cyclic Homology'', 
Grundlehren math.\ Wiss.\ \textbf{301}, Springer-Verlag, Berlin, 1998. 
%
\bibitem[Lo69]{Lo69} Loos, O., ``Symmetric Spaces I: General Theory'',
Benjamin, New York, 1969. 
%
\bibitem[Lo75]{Lo75} ---, ``Jordan Pairs'', 
Lecture Notes in Mathematics \textbf{460}, Springer-Verlag, Berlin, 1975. 
%
\bibitem[Lk92]{Lk92}
Losik, M.~V., {\it Fr\'{e}chet manifolds as diffeologic spaces}, Russ.\ Math.\ \textbf{36} 
 (1992), 31--37. 
%
\bibitem[LV94]{LV94} Luminet, D., and A. Valette, {\it Faithful uniformly continuous 
representations of Lie groups}, J. London Math.\ Soc.\ (2) {\bf 49} (1994), 100--108. 
%
\bibitem[Ma01]{Ma01}
Ma, T.~W., {\it Inverse mapping theorem on
coordinate spaces}, Bull.\ London Math.\ Soc.\ \textbf{33}  (2001), 473--482. 
%
\bibitem[ML63]{ML63} MacLane, S., ``Homology'',
Grundlehren math.\ Wiss.\
{\bf 114}, Springer-Verlag, Berlin, 1963.
%
\bibitem[ML71]{ML71} ---, ``Categories for the
Working Mathematician'', Springer-Verlag, New York, 1971. 
%
\bibitem[ML78]{ML78}  ---, {\it Origins of the cohomology of groups},
Enseign.\ math.\ {\bf 24:2} (1978), 1--29. 
%
\bibitem[MOK85]{MOK85} Maeda, Y., H. Omori, O. Kobayashi, and A. Yoshioka, {\it On 
regular Fr\'echet--Lie groups. VIII. Primordial operators and Fourier integral operators}, 
Tokyo J. Math.\ {\bf 8:1} (1985), 1--47.
%
\bibitem[MS75]{MS75}
Magnin, L., and J. Simon, {\it Lie algebras associated with topological nilpotent groups}, 
Rep.\ Math.\ Phys.\ {\bf 8} (1975), 171--180.
%
\bibitem[Mai02]{Mai02}
Maier, P., {\it Central extensions of topological
current algebras}, pp.\ 61--76 in: A. Strasburger et al.\ (eds), 
``Geometry and Analysis on Finite-
and Infinite-Dimensional Lie Groups,''
Banach Center Publications {\bf 55}, Warszawa, 2002.
%
\bibitem[MN03]{MN03} 
Maier, P., and K.-H. Neeb, {\it Central extensions of current groups},  
Math.\ Ann.\ \textbf{326:2} (2003), 367--415. 
%
\bibitem[Ms62]{Ms62}
Maissen, B., {\it Lie-Gruppen mit 
Banachr\"{a}umen als Parameterr\"{a}ume}, Acta Math.\ \textbf{108} (1962), 229--270. 
%
\bibitem[Ms63]{Ms63}
---, {\it \"{U}ber Topologien im Endomorphismenraum
eines topologischen Vektorraums}, Math.\ Ann.\ \textbf{151} (1963), 283--285. 
%
\bibitem[MP98]{MP98} Malikov, F., and R. C. Renner, {\it 
The Lie algebra of homeomorphisms of the circle}, 
Adv. Math. {\bf 140:2} (1998), 282--322. %
%
\bibitem[MO92]{MO92}
Margalef-Roig, J., and E. Outerelo Dominguez, ``Differential Topology,''
North-Holland, Amsterdam, 1992.
%
\bibitem[MR95]{MR95} Marion, J., and T. Robart, {\it Regular 
Fr\'{e}chet Lie groups of invertibe elements
in some inverse limits of unital involutive Banach algebras}, 
Georgian Math.\ J. \textbf{2} (1995), 425--444. 

\bibitem[Ma18]{Ma18} Marquis, T., ``An Introduction to Kac-Moody Groups over 
Fields,'' EMS Textbooks in Mathematics, 
 European Mathematical Society (EMS), Z\"urich, 2018 
%

\bibitem[MNe18]{MNe18}
Marquis, T., and K.-H. Neeb,
{\it Half Lie groups}, 
Transf. Groups \textbf{23} (2018), 801--840 

%
\bibitem[Mar67]{Mar67} Marsden, J. E., {\it Hamiltonian one parameter groups: A mathematical 
exposition of infinite dimensional Hamiltonian systems with applications 
in classical and quantum mechanics}, Arch.\ Rational Mech.\ Anal.\ {\bf 28} (1968), 
362--396.
%
\bibitem[MA70]{MA70}
Marsden, J., and R. Abraham, {\it Hamiltonian mechanics on Lie groups and 
Hydrodynamics}, pp.~237--243 in: S. S. Chern and S. Smale (eds), 
``Global Analysis'', Proc.\ Symp.\ Pure Math.\ \textbf{16}, 
Amer.\ Math.\ Soc., Providence, 1970.
%
\bibitem[Mat74]{Mat74} Mather, J. N., {\it Commutators of
diffeomorphisms}, Comment.\ Math.\ Helv.\ {\bf 49} (1974), 512--528. 
%
\bibitem[Mat75]{Mat75} ---, {\it Commutators of diffeomorphisms. II}, 
Comment.\ Math.\ Helv.\ {\bf 50} (1975), 33--40. 
%
\bibitem[MaTh35]{MaTh35} Mayer, W., and T. Y. Thomas, {\it Foundations of the theory of Lie groups}, 
Ann.\ of Math.\ {\bf 36:3} (1935), 770--822. 
%
\bibitem[Mau88]{Mau88} Maurer, L., {\it \"Uber allgemeinere Invarianten-Systeme}, 
M\"unchner Berichte {\bf 43} (1888), 103--150.
%
\bibitem[Maz33]{Maz33}
Mazur, S.,
{\it \"{U}ber konvexe Mengen in linearen normierten R\"{a}umen},
Stud.\ Math.\ {\bf 4} (1933), 70--84.
%
\bibitem[MDS98]{MDS98}
McDuff, D., and D. Salamon, ``Introduction to Symplectic
Topology'', Oxford Math.\ Monographs, Clarendon Press, Oxford, 1998. 
%
\bibitem[MD05]{MD05} McDuff, D., {\it Enlarging the Hamiltonian group}, 
Conference on Symplectic Topology,
J. Symplectic Geom.\ {\bf 3} (2005), 481--530. 
%
\bibitem[MV97]{MV97}
Meise, R., and D. Vogt,
``Introduction to Functional Analysis,'' Clarendon Press, Oxford, 1997.
%
\bibitem[Mi92]{Mi92} Mical, J., {\it Coordinates of the second kind on groups of currents}, 
Bulletin of the Polish Acad.\ Sciences, Mathematics {\bf 40:1} (1992), 31--39 
%
\bibitem[Mi94]{Mi94} ---, {\it Applications exponentielles pour les groupes 
de courants et la d\'ecomposition de Birkhoff pour les groupes des noeuds}, 
Dissertationes Math. CCCXXXII, Polska Akad. Nauk. Inst. Matmatycny, 
Warszawa, 1994 
%
\bibitem[MicE52]{MicE52} Michael, E.~A., ``Locally Multiplicatively Convex Topological
Algebras'', Mem.\ Amer.\ Math.\ Soc.\ \textbf{11}, Amer.\ Math.\ Soc., 1952. 
%
\bibitem[MicE59]{MicE59} ---, {\it Convex structures and continuous selections}, 
Can.\ J. Math.\ \textbf{11} (1959), 556--575.
%
\bibitem[Mic68]{Mic68} Michael, E., {\it Bi-quotient maps and Cartesian products
of quotient maps}, Ann.\ Inst.\ Fourier (Grenoble) {\bf 18} (1968), 287--302.
%
\bibitem[MicA38]{MicA38}
Michal, A.~D., {\it Differential calculus in linear topological spaces}, 
Proc.\ Nac.\ Acad.\ Sci.\ USA \textbf{24}  (1938), 340--342. 
%
\bibitem[MicA40]{MicA40} ---, {\it Differential of functions with arguments and values in 
topological abelian groups}, 
Proc.\ Nat.\ Acad.\ Sci.\ USA {\bf 26} (1940), 356--359. 
%
\bibitem[MicA45]{MicA45} ---, {\it The total differential equation for the exponential 
function in non-com\-mu\-ta\-tive normed linear rings}, 
Proc.\ Nat.\ Acad.\ Sci.\ USA  {\bf 31} (1945),  315--317.
%
\bibitem[MicA48]{MicA48} ---, {\it Differentiable infinite continuous groups in abstract 
spaces}, Revista Ci., Lima {\bf 50} (1948), 131--140.
%
\bibitem[MiE37]{MiE37}  Michal, A. D., and V. Elconin, {\it Differential 
properties of abstract transformation groups
with abstract parameters}, Amer.\ J. Math.\ {\bf 59} (1937), 129--143.
%
\bibitem[Mr80]{Mr80}
Michor, P.~W.,  ``Manifolds of Differentiable
Mappings'', Shiva Publishing, Orpington, Kent (U.K.), 1980. 
%
\bibitem[Mi84]{Mi84} ---, {\it A convenient setting for differential geometry
and global analysis I, II}, Cah.\ Topologie G\'{e}om.\ Differ.\ {\bf 25} (1984), 
63--109, 113--178.
%
\bibitem[Mi87]{Mi87} ---, {\it The cohomology of the diffeomorphism group of a manifold is a 
Gelfand-Fuks cohomology}, in ``Proc.\ of the 14th Winter School on Abstr.\ Analysis, 
Srni, 1986,'' Suppl.\ Rend.\ del Circ.\ Mat.\ Palermo II {\bf 14} (1987), 235--246. 
%
\bibitem[Mi91]{Mi91} ---, ``Gauge Theory for Fiber Bundles,'' Bibliopolis, ed.\ di fil.\ sci., 
Napoli, 1991. 
%
\bibitem[Mi20]{Mi20} ---,
\emph{Manifolds of mappings for continuum mechanics},
pp.\,3--75 in:
R. Segev et al.\ (eds.),
``Geometric Continuum Mechanics,'' Birkh\"{a}user, Cham, 2020.
%
\bibitem[MM13]{MM13} Michor, P.\ W., and D. Mumford, {\it A zoo of 
diffeomorphism groups on $\R^n$}, Ann. Global Anal. Geom. {\bf 44:4} 
(2013), 529--54
%
\bibitem[MT99]{MT99}
Michor, P.~W., and J. Teichmann, {\it Description of infinite
dimensional abelian regular Lie groups}, J. Lie Theory \textbf{9} (1999), 487--489. 
%
\bibitem[Mick87]{Mick87}
Mickelsson, J., {\it Kac-Moody groups, topology of the Dirac determinant
bundle, and fermionization}, Commun.\ Math.\ Phys.\ \textbf{110} (1987), 173--183. 
%
\bibitem[Mick89]{Mick89} ---, ``Current Algebras and Groups'', 
Plenum Press, New York, 1989. 
%
\bibitem[Mil82]{Mil82}
Milnor, J., {\it On infinite-dimensional
Lie groups}, Preprint, Institute for Advanced Study, 
Princeton, 1982. 
%
\bibitem[Mil84]{Mil84}
---, {\it Remarks on infinite-dimensional Lie groups}, pp.~1007--1057 
in:  B. DeWitt and R. Stora (eds), 
``Relativit\'{e}, groupes et topologie II
(Les Houches, 1983), North Holland, Amsterdam, 1984.
%
\bibitem[Mim95]{Mim95}
Mimura, M., {\it Homotopy theory of Lie groups}, 
in: I.~M. James (ed), ``Handbook
of Algebraic Topology'', North Holland, Amsterdam, 1995. 
%
\bibitem[Mit61]{Mit61} Mityagin, B. S.,
\emph{Approximative dimension and bases in nuclear spaces},
Russ.\ Math.\ Surv.\ {\bf 16} (1961), 59--128.
%
\bibitem[Mit70]{Mit70} ---, {\it The homotopy structure of the linear
group of a Banach space}, Uspekhi Mat.\ Nauk {\bf 25:5} (1970), 63--106.
%
\bibitem[ML18]{ML18} Moodie, J.C., and M.W. Long, {\it 
An exact power series representation of the 
Baker--Campbell--Hausdorff formula}, arXiv:math-ph:1807.07884v1 20 Jul 2018 
%
\bibitem[MR91]{MR91} Moerdijk, I., and G.~E. Reyes, ``Models for Smooth
Infinitesimal Analysis'', Springer-Verlag, Berlin, 1991.
%
\bibitem[MZ55]{MZ55}
Montgomery, D., and L. Zippin, ``Topological Transformation Groups'', 
Interscience, New York, 1955.
%
\bibitem[Mo64]{Mo64} 
Moore, C.~C., {\it Extensions 
and low dimensional cohomology theory of locally compact 
groups. I,~II}, Trans.\ Amer.\ Math.\ Soc.\ {\bf 113} (1964), 40--63; 63---86. 
%
\bibitem[Mo76]{Mo76}  ---, {\it Group extensions and cohomology of locally compact groups. III, IV},  
Trans.\ Amer.\ Math.\ Soc.\ {\bf 221:1} (1976), 1--33; 35--58.  
%
\bibitem[Mr56]{Mr56}
Morita, K., {\it On decompositon spaces of locally compact spaces},
Proc.\ Japan Acad.\ {\bf 32} (1956), 544--548.
%
\bibitem[Mo61]{Mo61} Moser, J., {\it A new technique for the construction of solutions of nonlinear 
differential equations}, Proc.\ Nat.\ Acad.\ Sci.\ USA {\bf 47} (1961), 1824--1831. 
%
\bibitem[MS03]{MS03} Mostkowitz, M., and R. Sacksteder, {\it 
The exponential map and differential equations on real Lie groups}, 
J. Lie Theory {\bf 13} (2003), 291--306 

\bibitem[Mu77]{Mu77}
Mujica, J., {\it Holomorphic germs on infinite-dimensional
spaces}, pp.\ 313--321 in: M.~C. Matos (ed),
``Infinite Dimensional Holomorphy and Applications'', 
North-Holland, Amsterdam, 1977. 
%
\bibitem[Mu79]{Mu79}
---, {\it Spaces of germs of holomorphic functions}, Studies in Analysis, Adv.\ in Math.\
Suppl.\ Stud.\ \textbf{4} (1979), 1--41. 
%
\bibitem[Mu83]{Mu83}
---, \emph{Spaces of continuous functions with values in an inductive
limit}, pp.~359--367 in: G. I. Zapata (ed), ``Functional Analysis, Holomorphy,
and Approximation Theory,'' Lect.\ Notes Pure Appl.\ Math.\ {\bf 83}, Marcel
Dekker, New York, 1983.
%
\bibitem[Mu84]{Mu84}
---, {\it A completeness criterion
for inductive limits of Banach spaces}, pp.\ 319--329 
in:
``Functional Analysis, Holomorphy and Approximation Theory, II'', 
(Rio de Janeiro, 1981), North Holland, Amsterdam, 1984. 
%
\bibitem[MW09]{MW09}
M\"{u}ller, C., and C. Wockel, \emph{Equivalences of smooth and continuous principal
bundles with infinite-dimensional structure group},
Adv.\ Geom.\ {\bf 9} (2009), 605--626. 
%
\bibitem[Mu08]{Mu08}
M\"uller, O.,
\emph{A metric approach to Fr\'{e}chet geometry},
J. Geom.\ Phys.\ {\bf 58} (2008), 1477--1500.
%
\bibitem[MST99]{MST99}
Mu\~{n}oz, G.~A., Y. Sarantopoulos, and A. Tonge,
\emph{Complexifications of real Banach spaces, polynomials and multilinear maps},
Stud.\ Math.\ {\bf 134:1} (1999), 1--33.
%
\bibitem[Mu88]{Mu88}
Murray, M. K., {\it Another construction of the central extension 
of the loop group}, Commun.\ Math.\ Phys.\ \textbf{116} (1988), 73--80. 
%
\bibitem[My54]{My54} Myers, S. B., {\it Algebras of differentiable functions}, 
Proc.\ Amer.\ Math.\ Soc.\ {\bf 5} (1954), 917--922.
%
\bibitem[Na69]{Na69}
Nachbin, L., ``Topology on Spaces of Holomorphic Mappings'', 
Er\-geb\-nis\-se der Math.\ \textbf{47}, Springer-Verlag, Berlin, 1969. 
%
\bibitem[Nag36]{Nag36} Nagumo, M., {\it Einige analytische Untersuchungen in linearen metrischen Ringen}, 
Jap.\ J. Math.\ {\bf 13} (1936), 61--80. 

\bibitem[Nar25]{Nar25} Nariman, S., {\it Continuous cohomology of diffeomorphism
    groups}, Preprint


\bibitem[Nat94]{Nat94}
Natarajan, L., {\it Unitary highest weight-modules of inductive
limit Lie algebras and groups}, J. Algebra \textbf{167} (1994), 9--28. 
%
\bibitem[NRW91]{NRW91} Natarajan, L., E. Rodriguez-Carrington, and
J.~A. Wolf, {\it Differentiable structure for direct limit groups},
Letters Math.\ Phys.\ \textbf{23} (1991), 99--109. 
%
\bibitem[NRW93]{NRW93}
---, {\it Locally convex Lie groups}, Nova J. Algebra Geom.\ \textbf{2:1} (1993), 59--87. 
%
\bibitem[NRW94]{NRW94}
---, {\it New classes of infinite dimensional Lie groups},
Proc.\ Symp.\ Pure Math.\ \textbf{56:2} (1994), 377--392. 
%
\bibitem[NRW01]{NRW01}
---, {\it The Bott-Borel Weil Theorem for direct limit
groups}, Trans.\ Amer.\ Math.\ Soc.\ \textbf{353} (2001), 4583--4622. 
%
\bibitem[Nat35]{Nat35} Nathan, D. S., {\it One-parameter groups of transformations in abstract vector spaces}, 
Duke Math.\ J. {\bf 1:4} (1935), 518--526.
%
\bibitem[Ne98]{Ne98} Neeb, K.--H., {\it Holomorphic highest weight representations
of infinite dimensional complex classical groups}, 
J.\ Reine Angew. Math.\ {\bf 497} (1998), 171--222  

\bibitem[Ne00]{Ne00} Neeb, K.-H., 
``Holomorphy and Convexity in Lie Theory,'' 
Expositions in Mathematics {\bf 28}, de Gruyter Verlag, Berlin, 2000.  
%
\bibitem[Ne01a]{Ne01a}
---, {\it Infinite-dimensional groups and their representations}, pp.~131--178
in: A. Huckleberry and T. Wurzbacher (eds), ``Infinite Dimensional K\"ahler Manifolds'',  
DMV-Seminar {\bf 31}, Birkh\"auser Verlag, Basel, 2001. 
%
\bibitem[Ne01b]{Ne01b} ---, {\it Borel-Weil theory for loop groups}, pp.\ 179--229 
in: A.~Huckleberry and T.~Wurzbacher (eds), 
``Infinite Dimensional K\"ahler Manifolds,'' DMV-Seminar \textbf{31}, 
Birkh\"auser Verlag, Basel, 2001.
%
\bibitem[Ne01c]{Ne01c} ---, {\it Locally finite Lie algebras with unitary highest weight
representations},
manu\-scrip\-ta mathe\-matica {\bf 104:3} (2001), 343--358.
%
\bibitem[Ne02a]{Ne02a}
---, {\it Central extensions of infinite-dimensional
Lie groups}, Ann.\ Inst.\ Fourier \textbf{52:5} (2002), 1365--1442. 
%
\bibitem[Ne02b]{Ne02b}
---, {\it Universal central extensions of Lie groups},  
Acta Appl.\ Math.\ \textbf{73} (2002), 175--219. 
%
\bibitem[Ne02c]{Ne02c} ---, {\it Classical Hilbert--Lie groups, their extensions and their
homotopy groups}, pp.~87--151 in: A. Strasburger et al.\ (eds), 
``Geometry and Analysis on Finite- and Infinite-dimensional Lie Groups,'' 
Banach Center Publications {\bf 55}, Warsaw, 2002. 
%
\bibitem[Ne02d]{Ne02d} ---, {\it A Cartan--Hadamard Theorem for
Banach--Finsler manifolds},
Geometriae Dedicata {\bf 95} (2002), 115--156. 
%
\bibitem[Ne03]{Ne03}  
---, {\it Locally convex root graded Lie algebras},   
Travaux math.\ {\bf 14} (2003), 25--120.  
%
\bibitem[Ne04a]{Ne04a}  ---, {\it Abelian extensions of infinite-dimensional Lie
groups}, Travaux math.\ {\bf 15} (2004), 69--194. 
%
\bibitem[Ne04b]{Ne04b}
---, {\it Infinite-dimensional Lie groups and their
representations}, pp.~213--328 in:
J.~P. Anker and B. \O{}rsted (eds),
  ``Lie Theory: Lie Algebras and Representations,'' 
Progress in Math.\ {\bf 228}, 
Birkh\"auser Verlag, 2004.
%
\bibitem[Ne04c]{Ne04c}
 ---,  {\it Current groups for non-compact manifolds and their 
central extensions}, pp.~109--183 in: T. Wurzbacher (ed),
``Infinite Dimensional Groups and Manifolds'', 
IRMA Lectures in Mathematics and
Theoretical Physics {\bf 5}, de Gruyter Verlag, Berlin, 2004.
%
\bibitem[Ne05]{Ne05}
---, {\it Exact sequences for Lie group cohomology with non-abelian coefficients}, 
in preparation. 
%
\bibitem[Ne06a]{Ne06a} 
---, {\it Towards a Lie theory of locally convex groups},
Jpn.\ J. Math.\ {\bf 1} (2006), 291--468. 
%
\bibitem[Ne06b]{Ne06b} ---, {\it Lie algebra extensions and higher order cocycles},  
J. Geom.\ Symmetry Phys.\ {\bf 5} (2006), 48--74.
%
\bibitem[Ne06c]{Ne06c} ---, {\it On the period group of a continuous inverse algebra}, in preparation.  
%
\bibitem[Ne06d]{Ne06d}
---, {\it Non-abelian extensions of topological Lie algebras}, 
Comm.\ Algebra {\bf 34} (2006), 991--1041. 
%
\bibitem[Ne07]{Ne07}
---, Non-abelian extensions of infinite-dimensional Lie groups, 
Ann.\ Inst.\ Fourier (Grenoble) {\bf 57} (2007), 209--271.
%
\bibitem[Ne08]{Ne08} ---, {\it Lie 
group extensions associated to projective modules 
of continuous inverse algebras}, 
in ``Proceedings of the Winter School 
on Geometry and Physics in Srni, 2008'', 
Archivum Mathematicum (Brno) {\bf 44} (2008), 339--363. 
%
\bibitem[Ne10]{Ne10} ---,  {\it 
On differentiable vectors for representations of infinite dimensional 
Lie groups},
J.~Funct.\ Anal.\ {\bf 259} (2010), 2814--2855. 
%
\bibitem[Ne11]{Ne11} ---, 
{\it On analytic vectors for 
unitary representations of infinite dimensional Lie groups}, 
Ann.\ Inst.\ Fourier (Grenoble) {\bf 61:5} (2011), 1441--1476. 
%
\bibitem[Ne14]{Ne14} ---, {\it Unitary representations of Unitary Groups}, 
in ``Lie theory workshops'', Eds. G.~Mason, 
I. Penkov, J. Wolf, Developments in Math. {\bf  37}, 
Springer, 2014, 197--243
%
\bibitem[Ne17]{Ne17} ---,  {\it Bounded and Semi-bounded 
Representations of Infinite Dimensional Lie Groups}, 
in ``Representation Theory---Current Trends and Perspectives,'' 
Eds. P. Littelmann et al, European Math. Society, 
EMS Ser. Congr. Rep., Eur. Math. Soc., Z\"urich, 2017, pp. 541--563, 
arXiv:math.RT:1510.08695. 
%
\bibitem[N\O{}06]{NO06} Neeb, K.-H., and B.~\O{}rsted, {\it A topological 
Maslov index for $3$-graded Lie groups}, 
Journal of Funct. Anal. {\bf 233} (2006), 426--477. 

\bibitem[NP00]{NP00} Neeb, K.-H., and D. Pickrell, 
{\it Supplements on the papers 
entitled ``On a Theorem of S.\ Banach'' and 
``The separable representations of $U(H)$''}, 
J. Lie Theory {\bf 10} (2000), 107--109 

\bibitem[NS13]{NS13} Neeb, K.-H., and H. Salmasian, 
{\it Differentiable vectors and unitary representations 
of Frechet--Lie supergroups}, Math.\ Z.
{\bf 275:1} (2013), 419--451.
%
\bibitem[NSZ14]{NSZ14} 
Neeb, K.-H., H. Salmasian, and C. Zellner,
\emph{On an invariance property of the space of smooth vectors},
Kyoto~J. Math.\ {\bf 55:3} (2015), 501--515. 

\bibitem[NSZ17]{NSZ17} ---, 
{\it Smoothing operators and $C^*$-algebras for 
infinite dimensional Lie groups}, 
International Journal of Mathematics {\bf 28:5} (2017), 
DOI 10.1142/S0129167X1750042. 

%
\bibitem[NS01]{NS01} Neeb, K.-H., and N. Stumme, {\it On the classification of 
locally finite split simple Lie algebras}, J.\ reine angew.\ Math.\ {\bf 533} (2001),
25--53.
%
\bibitem[NTS15]{NTS15} Neeb, K.-H., T. Thiemann, and H. Sahlmann, 
{\it Weak Poisson structures on infinite dimensional manifolds and hamiltonian actions}, 
in ``Lie Theory and Its Applications in Physics,''  
V.~Dobrev, ed., Springer Proceedings in Mathematics \& Statistics, 
Springer, 2015; 105-136. 

\bibitem[NV03]{NV03}
Neeb, K.-H., and C. Vizman, {\it Flux homomorphisms and principal bundles over
infinite-dimensional manifolds}, Monatsh.\ Math.\ {\bf 139} (2003), 309--333. 
%
\bibitem[NeW08a]{NeW08a} Neeb, K.-H., and F. Wagemann,
{\it The second cohomology of current algebras 
of general Lie algebras},
Canadian J. Math.\ {\bf 60:4} (2008), 892--922.  
%
\bibitem[NeW08b]{NeW08b} ---, {\it Lie group structures on groups of 
smooth and holomorphic maps}, Geom.\ Dedicata {\bf 134} (2008), 17--60. 
%
\bibitem[NeWo09]{NeWo09} Neeb, K.-H., and C.~Wockel, {\it 
Central extensions of groups of sections},
Annals of Global Analysis and Geometry {\bf 36:4} (2009), 381--418. 
%
\bibitem[Neh93]{Neh93} Neher, E., {\it Generators and relations for $3$-graded Lie
algebras}, J. Algebra {\bf 155} (1993), 1--35. 
%
\bibitem[Neh96]{Neh96} ---, {\it Lie algebras graded by  $3$-graded root systems and 
Jordan pairs covered by grids}, Amer.\ J. Math.\ {\bf 118} (1996), 439--491. 
%
\bibitem[NR17]{NR17}
Nenning, D.~N. and A. Rainer, {\it The Trouv\'{e} group for spaces
of test functions}, preprint, 2017, arXiv:1711.01196.
%
\bibitem[Ner83]{Ner83}
Neretin, Yu.\ A., {\it Boson representation of the diffeomorphisms
of the circle},  Sov.\ Math.\ Dokl.\ \textbf{28} (1983), 411--414. 
%
\bibitem[Ner87a]{Ner87a}
---, {\it On spinor representations of $\OO(\infty, \C)$}, 
Sov.\ Math.\ Dokl.\ \textbf{34:1} (1983), 71--74.
%
\bibitem[Ner87b]{Ner87b}
---, {\it A complex semigroup containing the group
of diffeomorphisms of the circle}, Funct.\ Anal.\ Appl.\ \textbf{21}  (1987), 82--83. 
%
\bibitem[Ner90]{Ner90}
---, {\it Holomorphic extensions of representations of the 
group of diffeomorphisms of the circle}, Math.\ USSR Sbornik 
\textbf{67} (1990), 75--97. 
%
\bibitem[Ner93]{Ner93} ---, {\it Categories of bistochastic measures, and 
representations of some infinite dimensional groups}, 
Russ.\ Acad.\ Sci.\ Sb.\ Math.\ \textbf{75} (1993), 197--219. 
%
\bibitem[Ner96]{Ner96} ---, ``Categories of Symmetries and
Infinite-Dimensional Groups'', London Math.\ Soc.\ Monographs, N.S. \textbf{16}, 
Clarendon Press, Oxford, 1996.
%
\bibitem[NO97]{NO97} Neretin, Yu.~A., and G.~I. Ol'shanski\u\i{}, 
{\it Boundary values
of holomorphic functions, singular unitary representations of
$O(p,q)$, and their limits as $q \to \infty$}, 
J. Math.\ Sci.\ (New York) \textbf{87} (1997), 3983--4035. 
%
\bibitem[Ne87]{Ne87} Nessonov, N.~I., 
{\it A complete classification of the
representations of $\GL(\infty)$ containing the identity
representation of the unitary subgroup}, 
Math.\ USSR Sbornik \textbf{58:1} (1987), 127--147.
%
\bibitem[vN29]{vN29} von Neumann, J., {\it \"Uber die analytischen Eigenschaften von Gruppen 
linearer Transformationen}, Math.\ Z. {\bf 30:1} (1929), 3--42.
%
\bibitem[vN50]{vN50}
von Neumann, J., {\it Functional Operators II: The Geometry of
Orthogonal Spaces}, Annals of Math.\ Studies {\bf 22}, Princeton
Univ.\ Press, Princeton, 1950. 
%
\bibitem[New51]{New51} Newburgh, J. D., {\it The variation of spectra}, Duke Math.\ J.
{\bf 18} (1951), 165--176
%
\bibitem[NW56]{NW56}
Newns, W. F.
and A. G. Walker,
\emph{Tangent planes to a differentiable manifold},
J. Lond.\ Math.\ Soc.\ {\bf 31} (1956), 400--406.
%
\bibitem[Nij74]{Nij74}
Nijenhuis, A.,
{\it Strong derivatives and inverse mappings},
Am.\ Math.\ Mon.\ {\bf 81} (1974), 969--980.
\bibitem[Lec79]{Lec79} Lecomte, P., {\it Sur l'alg\`ebre de Lie des automorphismes infinitesimaux 
du fibre tangent}, C. R. Acad.\ Sci.\ Paris, S\'er.\ A {\bf 288} (1979), 661--663. 
%
\bibitem[Nik15]{Nik15}
Nikitin, N.,
{\it Exponential laws for weighted function spaces and regularity of weighted mapping groups},
preprint, 2015, arXiv:1512.07211
%
\bibitem[Nik21]{Nik21}
---,
``Regularity Properties of Infinite-Dimensional Lie Groups
and Exponential Laws,'' doctoral thesis, Paderborn University,
2021.
%
%
\bibitem[Nb70]{Nb70}
Noble, N.,
{\it The continuity of functions on cartesian products},
Trans.\ Am.\ Math.\ Soc.\ {\bf 149} (1970), 187--198.
%
\bibitem[OMS81]{OMS81}
Okamoto, K., H. Matsushima, and T. Sakurai, {\it On a certain class of
irreducible unitary representations of the infinite dimensional
rotation group I}, Hiroshima Math.\ J. \textbf{11:1} (1981), 181--193. 
%
\bibitem[OS82a]{OS82a} Okamoto, K., and T. Sakurai, 
{\it On a certain class of irreducible unitary representations of the
infinite dimensional rotation group II}, Hiroshima Math.\ J. \textbf{12} (1982), 385--397. 
%
\bibitem[OS82b]{OS82b}
---, {\it An analogue of Peter-Weyl
Theorem for the infinite dimensional unitary group},  
Hiroshima Math.\ J. \textbf{12} (1982), 529--541.
%
\bibitem[Ol84]{OL84}
Ol'shanski\u\i{}, G.~I., {\it Infinite dimensional classical
groups of finite $\R$-rank: Description of representations and
asymptotic theory}, Funkts.\ Anal.\ Prilozhen \textbf{18:1} (1984), 28--42. 
%
\bibitem[Ol86]{Ol86} ---, {\it Unitary representation of the
group $SO_o(\infty, \infty)$ as limits of unitary representations of
the groups $SO_o(n, \infty)$ as $n \to \infty$},  
Funkts.\ Anal.\ Prilozhen \textbf{20:4} (1986), 46--57. 
%
\bibitem[Ol88]{Ol88}
---, {\it Method of holomorphic extension in the
theory of unitary representations of infinite dimensional classical
groups}, Funct.\ Anal.\ Appl.\ \textbf{22:4} (1988), 273--285. 
%
\bibitem[Ol90]{Ol90} ---, {\it Unitary representations of infinite
dimensional pairs $(G,K)$ and the formalism of R.\ Howe}, pp.\ 269--463 
in: A. Vershik and D. Zhelobenko (eds), ``Representations of 
Lie Groups and Related Topics'', Gordon \& Breach, New York, 1990. 
%
\bibitem[Ol91a]{Ol91a}
---, {\it Representations of infinite dimensional 
classical groups, limits of envelopping algebras, and Yangians}, pp.\ 1-66 
in:
``Topics in Representation Theory'', Advances Sov.\ Math.\ {\bf 2},
Amer.\ Math.\ Soc., Providence, 1991. 
%
\bibitem[Ol91b]{Ol91b}
---, {\it On semigroups related to infinite 
dimensional groups}, pp.\ 67--101 in: 
``Topics in Representation Theory'', Advances Sov.\ Math.\ \textbf{2}, 
Amer.\ Math.\ Soc., Providence, 1991. 
%
\bibitem[Ol91c]{Ol91c} ---, 
{\it Caract\`eres generalis\'es de
$U(\infty)$ et fonctions int\'erieurs},  C. R. Acad.\ Sci.\ Paris \textbf{313} 
(1991), 9--12.
%
\bibitem[OV96]{OV96}
Ol'shanski\u\i, G.~I. and A. Vershik,  {\it Ergodic unitarily
invariant measures on the space of infinite hermitian matrices}, 
Amer.\ Math.\ Soc.\ Transl.\ (2) \textbf{175} (1996), 137--175. 
%
\bibitem[Olv93]{Olv93} Olver, P.~J., ``Applications of Lie Groups to Differential Equations.''  
Graduate Texts in Math.\ {\bf 107}, 
Springer-Verlag, New York, ${}^2$1993.
%
\bibitem[Omo70]{Omo70} Omori, H., {\it On the group of diffeomorphisms on a compact manifold}, 
in
``Global Analysis'' (Proc. Sympos. Pure Math., Vol. XV, Berkeley, Calif., 1968), 
Amer.\ Math.\ Soc., Providence, R. I., 1970, 167--183.
%
\bibitem[Omo72a]{Omo72a} Omori, H., {\it Local structures of groups of diffeomorphisms}, 
J. Math.\ Soc.\ Japan {\bf 24} (1972), 60--88.
%
\bibitem[Omo72b]{Omo72b} ---, {\it On smooth extension theorems}, 
J. Math.\ Soc.\ Japan {\bf 24} (1972), 405--432.
%
\bibitem[Omo73]{Omo73} ---, {\it Groups of diffeomorphisms and their subgroups}, 
Trans.\ Amer.\ Math.\ Soc.\ {\bf 179} (1973), 85--122.
%
\bibitem[Omo74]{Omo74} ---, ``Infinite Dimensional Lie Transformation Groups,'' 
Lect.\ Notes Math.\ {\bf 427}, Springer-Verlag, Berlin-New York, 1974.
%
\bibitem[Omo78]{Omo78} ---, {\it On Banach-Lie groups 
acting on finite-dimensional manifolds}, T\^{o}hoku Math.\ J. {\bf 30} (1978), 223--250. 
%
\bibitem[Omo80]{Omo80} ---, {\it A method of classifying expansive singularities}, 
J.\ Diff.\ Geom.\ {\bf 15} (1980), 493--512. 
%
\bibitem[Omo81]{Omo81} ---, {\it A remark on non-enlargible Lie algebras}, 
J. Math.\ Soc.\ Japan {\bf 33:4} (1981), 707--710. 
%
\bibitem[Omo97]{Omo97} ---, ``Infinite-Dimensional Lie Groups'', Translations
of Math.\ Monographs \textbf{158}, Amer.\ Math.\ Soc.,
Providence, 1997. 
%
\bibitem[OdH71]{OdH71} Omori, H., and P. de la Harpe, {\it Op\'eration de groupes de Lie 
banachiques sur les vari\'et\'es diff\'erentielles de dimension finie}, 
C. R. Acad.\ Sci.\ Paris S\'er. A-B {\bf 273} (1971), A395--A397.
%
\bibitem[OdH72]{OdH72} ---, {\it About interactions between 
Banach-Lie groups and finite dimensional mani\-folds}, 
J. Math.\ Kyoto Univ.\ {\bf 12} (1972), 543--570.
%
\bibitem[OMY80]{OMY80} Omori, H., Y. Maeda, and A. Yoshioka, 
{\it On regular Fr\'echet--Lie groups}. I.
{\it Some differential geometrical expressions of Fourier 
integral operators on a Riemannian manifold}, Tokyo J. Math.\ {\bf 3:2} (1980), 353--390.
%
\bibitem[OM81a]{OM81a} ---, {\it On regular Fr\'echet--Lie groups}. II. {\it Composition 
rules of Fourier-integral operators on a Riemannian manifold}, 
Tokyo J. Math.\ {\bf 4:2}  (1981), 221--253. 
%
\bibitem[OM81b]{OM81b} Omori, H., Y. Maeda, A. Yoshioka, and O. Kobayashi, 
{\it On regular Fr\'echet--Lie groups}. III. {\it A second cohomology class related 
to the Lie algebra of pseudodifferential operators of order one},  
Tokyo J. Math.\ {\bf 4:2}  (1981), 255--277. 
%
\bibitem[OMY82]{OMY82} ---, {\it On regular Fr\'{e}chet-Lie groups}.
IV. {\it Definition and fundamental theorems},
Tokyo J. Math.\ {\bf 5} (1982), 365--398.
%
\bibitem[OMY83a]{OMY83a} ---, {\it On regular Fr\'echet--Lie groups. V. 
Several basic properties}, Tokyo J. Math.\ {\bf 6:1} (1983), 39--64. 
%
\bibitem[OMY83b]{OMY83b} ---, {\it On regular Fr\'echet--Lie groups}.
VI. {\it Infinite-dimensional 
Lie groups which appear in general relativity},
Tokyo J. Math.\ {\bf 6:2} (1983), 217--246. 
%
\bibitem[Oer80]{Oer80}
\O{}rsted, B., {\it A model for an interacting quantum field}, 
J. Funct.\ Anal.\ \textbf{36} (1980), 53--71. 
%
\bibitem[Ono06]{Ono06}
Ono, K., {\it Floer-Novikov cohomology and the flux conjecture},
Geom.\ Funct.\ Anal.\ {\bf 16} (2006), 981--1020. 
%
\bibitem[Ot95]{Ot95} Ottesen, J.~T., ``Infinite Dimensional
Groups and Algebras in Quantum
Physics'',  Lecture Notes in Physics, N.S. \textbf{27}, 
Springer-Verlag, Berlin, 1995.   
%
\bibitem[OR98]{OR98}
Ovsienko, V., and C. Roger, {\it Generalizations of Virasoro
group and Virasoro algebra through extensions by modules of tensor 
densities on $S^1$}, Indag.\ Math.\ N. S. {\bf 9:2} (1998), 277--288. 
%
\bibitem[Pa57]{Pa57} Palais, R. S., ``A Global Formulation of the Lie Theory of Transformation Groups,'' 
Mem.\ Amer.\ Math.\ Soc.\ {\bf 22}, Amer.\ Math.\ Soc., 1957. 
%
\bibitem[Pa61]{Pa61} ---, {\it Equivalence of nearby differentiable 
actions of a compact Lie group}, 
Bull. Amer. Math. Soc. {\bf 67} (1961), 362--364 
%
\bibitem[Pa65]{Pa65} ---, {\it On the homotopy type of certain groups of 
operators}, Topology \textbf{3} (1965), 271--279. 
%
\bibitem[Pa66]{Pa66} ---, {\it Homotopy theory of infinite dimensional manifolds}, 
Topology \textbf{5} (1966), 1--16. 
%
\bibitem[Pa68]{Pa68} ---, ``Foundations of Global Nonlinear Analysis,'' Benjamin, N. Y., 1968.  
%
%
\bibitem[Pali68]{Pali68} Palis, J., {\it On Morse-Smale dynamical systems}, Topology {\bf 8} (1968), 385--404.
%
\bibitem[Pali74]{Pali74} ---, {\it Vector fields generate few diffeomorphisms}, Bull.\ Amer.\ Math.\ 
Soc.\ {\bf 80:3} (1974), 503--505.
%
\bibitem[Pat83]{Pat83}
Paterson, A.~L.~T., {\it Harmonic analysis on unitary groups}, 
J. Funct.\ Anal.\ \textbf{53} (1983), 203--223. 
%
\bibitem[PS09]{PS09}
Patyi, I., and S.~B. Simon,
{\it On real analytic Banach manifolds},
J. Lond.\ Math.\ Soc.\ (2) {\bf 80:2} (2009), 375--387.
%
\bibitem[PCB87]{PCB87}
P\'{e}rez Carreras, P., and J.  Bonet,
``Barrelled Locally Convex Spaces,''
North-Holland, Amsterdam, 1987.
%
\bibitem[Pe92]{Pe92} Pestov, V.~G., {\it Nonstandard hulls
of Banach-Lie groups and algebras}, 
Nova J.\ Algebra Geom.\ \textbf{1} (1992), 371--381. 
%
\bibitem[Pe93a]{Pe93a} ---, {\it Free Banach--Lie algebras, couniversal Banach-Lie groups, 
and more}, Pacific J. Math.\ {\bf 157:1} (1993), 137--144. 
%
\bibitem[Pe93b]{Pe93b} ---, {\it Enlargible Banach--Lie algebras
and free topological groups},
Bull.\ Aust.\ Math.\ Soc.\ {\bf 48} (1993), 13--22.
%
\bibitem[Pe95a]{Pe95a}  ---, {\it Correction to ``Free Banach--Lie algebras, couniversal Banach-Lie groups, 
and more''}, Pacific J. Math.\ {\bf 171:2} (1995), 585--588. 
%
\bibitem[Pe95b]{Pe95b} ---, {\it Regular Lie groups and a theorem of Lie-Palais}, 
J. Lie Theory {\bf 5:2} (1995), 173--178.
%
\bibitem[PU21]{PU21} Pestov, V.~G., and V.V.~Uspenskij, 
{\it On epimorphisms in some categories of infinite-dimensional 
Lie groups}, J. Lie Theory {\bf 31:3} (2021), 871--884. 

\bibitem[Ptr59]{Ptr59}
Peetre, J.,
\emph{Une caract\'{e}risation abstraite des op\'{e}rateurs diff\'{e}rentiels},
Math.\ Scand.\ {\bf 7} (1959), 211--218.
%
\bibitem[Ph91]{Ph91}
Phillips, N.~C., {\it $K$-theory for Fr\'{e}chet
algebras}, Int.\ J. Math.\ \textbf{2} (1991), 77--129.
%
\bibitem[Pi88]{Pi88}
Pickrell, D., {\it The separable representations of $U(\cH)$}, 
Proc.\ Amer.\ Math.\ Soc.\ \textbf{102} (1988), 416--420.
%
\bibitem[Pi89]{Pi89} ---, {\it On the Mickelsson-Faddeev
Extension and Unitary Representations},
Comm.\ Math.\ Phys.\ \textbf{123} (1989), 617--625. 
%
\bibitem[Pi90]{Pi90}
---, {\it On $U(\infty)$-invariant measures}, 
unpublished preprint, 1990. 
%
\bibitem[Pi96]{Pi96}
---, {\it On $YM_2$ measures and area-preserving
diffeomorphisms}, J. Geom.\ Physics \textbf{19} (1996), 315--367.
%
\bibitem[Pic00a]{Pic00a} Pickrell, D., ``Invariant measures for unitary groups 
associated to Kac-Moody Lie algebras,''
Mem.\ Amer.\ Math.\ Soc.\  {\bf 146/693}, 2000.
%
\bibitem[Pic00b]{Pic00b} ---, {\it On the action of the group of 
diffeomorphisms of a surface on sections of the determinant line bundle}, 
Pacific J. Math.\ {\bf 193:1} (2000), 177--199.
%
\bibitem[Pi72]{Pi72}
Pietsch, A., ``Nuclear Locally Convex Spaces'', 
Springer-Verlag, New York, 1972.
%
\bibitem[Pis76]{Pis76} Pisanelli, D., {\it An extension of the exponential of a matrix 
and a counter example to the inversion theorem in a space $H(K)$},
Rendiconti di Matematica (6) (3) {\bf 9:3} (1976), 465--475.
%
\bibitem[Pis77]{Pis77} ---, {\it An example of an infinite
Lie group}, Proc.\ Amer.\ Math.\ Soc.\ {\bf 62} (1977), 156--160.
%
\bibitem[Pis79]{Pis79} ---, {\it The second Lie theorem in the group $Gh(n,\C)$},
in: J.~A. Barroso (ed), ``Advances in Holomorphy,'' North Holland Publ., 1979. 
%
\bibitem[Pl94]{Pl94} Plymen, R.J., ``Spinors in Hilbert space,'' 
Cambridge Tracts in Math. {\bf 114}, Cambridge Univ. Press, Cambridge, 1994 

\bibitem[PS75]{PS75} Plymen, R.J., and R.F.\ Streater, {\it 
A model of the universal covering group of $\SO_2(\cH)$}, 
Bull. London Math. Soc. {\bf 7} (1975), 283--288

\bibitem[Pol01]{Pol01} Polterovich, L.,``The geometry of the group of 
symplectic diffeomorphisms,'' Lectures in Mathematics, ETH Z\"urich, 
Birkh\"auser Verlag, Basel, 2001.
%
\bibitem[Po39]{Po39} Pontrjagin, L., ``Topological Groups,'' 
Princeton Mathematical Series {\bf 2}, 
Princeton University Press, Princeton, 1939.
%
\bibitem[PBN03]{PBN03} Popovych, R.O., V.M.\,Boyko, M.O.\,Nesterenko, 
and M.W.\, Lutfullin, {\it Realizations of real low-dimensional 
Lie algebras}, J. Phys. A {\bf 36:26} (2003), 7337--7360.
%
\bibitem[Pop99]{Pop99}
Poppenberg, M.,
{\it An inverse function theorem for Fr\'{e}chet spaces satisfying a smoothing property and}
(DN). Math.\ Nachr.\ {\bf 206} (1999), 123--145.
%
\bibitem[Pr80]{Pr80}
Pressley, A., {\it Decompositions of the space of loops on a
Lie group}, Topology \textbf{19} (1980), 65--79.
%
\bibitem[PS86]{PS86}
Pressley, A., and G. Segal, ``Loop Groups'', Oxford University Press, Oxford, 
1986.
%
\bibitem[PzS21]{PzS21}
Prinz, D. and A. Schmeding,
\emph{Lie theory for asymptotic symmetries in general relativity:
the NU group},
preprint, arXiv:2109.11476.
%
\bibitem[PSS25]{PSS25}
Prinz, D., A. Schmeding, and P.K. Schwartz, 
\emph{Group contractions via infinite-dimensional Lie theory},
Preprint, arXiv:2512.18530

\bibitem[Pu52]{Pu52} Pursell, M. E., ``Algebraic structures associated with smooth manifolds,'' 
Thesis, Purdue Univ., 1952.
%
\bibitem[PuSh54]{PuSh54} Pursell, M. E., and M. E. Shanks, {\it The Lie algebra of a smooth manifold}, 
Proc.\ Amer.\ Math.\ Soc.\ {\bf 5} (1954), 468--472. 
%
\bibitem[PW52]{PW52} Putnam, C.~R., and A.\ Wintner, {\it The orthogonal group in
Hilbert space}, Amer.\ J.\ Math.\ {\bf 74} (1952), 52--7.8
%
\bibitem[Rae77]{Rae77} Raeburn, I., {\it The relationship between a commutative Banach
algebra and its maximal ideal space}, J. Funct.\ Anal.\ {\bf 25} (1977), 366--390. 
%
\bibitem[Rai59]{Rai59}
Ra\u{\i}kov, D.\,A., {\it A criterion of completeness
of locally convex spaces}, Uspehi Mat.\ Nauk {\bf 14:85} (1959), 223--229.
%
\bibitem[RS81]{RS81} Ratiu, T., and R. Schmid, {\it The differentiable structure
of three remarkable diffeomorphism groups}, Math.\ Z. {\bf 177} (1981), 81--100.
%
\bibitem[RO03]{RO03} Ratiu, T., and A. Odzijewicz, {\it Banach Lie--Poisson spaces and reduction}, 
Comm.\ Math.\ Phys.\ {\bf 243:1} (2003),  1--54.
%
\bibitem[RO04]{RO04} ---, {\it Extensions of Banach Lie--Poisson spaces}, 
 J. Funct.\ Anal.\ {\bf 217:1} (2004), 103--125.
%
\bibitem[RS75]{RS75} Reed, S., and B. Simon, 
``Methods of Modern Mathematical Physics I: Functional Analysis,'' 
Academic Press, New York, 1973.
%
\bibitem[Re95]{Re95} Remmert, R., ``Funktionentheorie II'', Springer Lehrbuch, 
Springer-Verlag, Berlin, 1995.   
%
\bibitem[Ri50]{Ri50} Ritt, J. F., {\it Differential groups
and formal Lie theory for an infinite number
of parameters}, Ann.\ of Math.\ {\bf 52} (1950), 708--726. 
%
\bibitem[Rs53]{Rs53}
Riss, J.
{\it \'{E}l\'{e}ments de calcul diff\'{e}rentiel et th\'{e}orie des distributions sur les groupes ab\'{e}liens localement compacts},
Acta Math.\ {\bf 89} (1953), 45--105.
%
\bibitem[Rob96]{Rob96} Robart, T., 
{\it Groupes de Lie de dimension infinie. Second et troisi\`eme 
th\'eor\`emes de Lie. I. Groupes de premi\`ere esp\`ece}, 
C. R. Acad.\ Sci.\ Paris S\'er.\ I Math.\ {\bf 322:11} (1996), 1071--1074. 
%
\bibitem[Rb97]{Rb97}
---, {\it Sur l'int\'{e}grabilit\'{e}
des sous-alg\`{e}bres de Lie en dimension infinie}, 
Canad.\ J. Math.\ \textbf{49} (1997), 820--839. 
%
\bibitem[Rob02]{Rob02} ---, {\it Around the exponential mapping}, pp.\ 11-30
in ``Infinite Dimensional Lie Groups in Geometry and Representation 
Theory'', Banyaga, A., J. A. Leslie, and T. Robart, eds., 
World Scientific, River Edge, NJ, 2002; 11--30.  
%
\bibitem[Rob04]{Rob04}  ---, {\it On Milnor's regularity and the path-functor for the 
class of infinite dimensional Lie algebras of CBH type}, 
Algebras, Groups and Geometries {\bf 21} (2004), 367--386.
%
\bibitem[RK97]{RK97} Robart, T., and N. Kamran, {\it Sur la th\'eorie locale des pseudogroupes 
de transformations continus infinis. I}, Math.\ Ann.\ {\bf 308:4} (1997), 593--613.
%
\bibitem[RS18]{RS18}
Roberts, D. M. and A. Schmeding,
\emph{Extending Whitney's extension theorem: nonlinear function spaces},
to appear in Annales de l'Institut Fourier
(cf.\ arXiv:1801.04126).
%
\bibitem[Ro68]{Ro68}
Robbin, J.~W.,
\emph{On the existence theorem for differential equations},
Proc.\ Amer.\ Math.\ Soc.\ {\bf 19} (1968), 1005--1006. 
%
\bibitem[Ry63]{Ry63}
Roby, N., {\it Lois polyn\^omes et lois formelles
en th\'eorie de modules}, Ann.\ Sci.\ Ec.\ Norm.\ Sup.\ 3e.\ s\'erie 
\textbf{80} (1963), 213--348.
%
\bibitem[Rod89]{Rod89} Rodriguez-Carrington, E., 
{\it Lie groups associated to Kac--Moody Lie algebras: an analytic approach}, 
in ``Infinite-dimensional Lie Algebras and Groups,'' (Luminy-Marseille, 1988), Adv.\
Ser.\ Math.\ Phys.\ {\bf 7}, World Scientific Publ., Teaneck, NJ, 1989; 
57--69.
%
\bibitem[RD81]{RD81}
Roelcke, W., and S. Dierolf, ``Uniform Structures on Topological Groups and their
Quotients,'' McGraw-Hill, New York, 1981.
%
\bibitem[RLL00]{RLL00} 
R\o{}rdam, M., F. Larsen, and N.~J. Laustsen, 
``An Introduction to K-theory for $C^*$-algebras'', London Math.\ Soc.,
Student Texts {\bf 49}, Cambridge Univ.\ Press, 2000. 
%
\bibitem[Ro95]{Ro95}
Roger, C., {\it Extensions centrales d'alg\`ebres et de groupes
de Lie de dimension infinie, alg\`ebres de Virasoro et
g\'en\'eralisations}, Reports Math.\ Phys.\ \textbf{35} (1995), 225--266.
%
\bibitem[Ros94]{Ros94} Rosenberg, J., ``Algebraic K-theory and its Applications,'' 
Grad.\ Texts  Math.\ {\bf 147}, Springer-Verlag, 1994. 
%

\bibitem[Ro09]{Ro09} Rosendal, C., {\it 
Automatic continuity of group homomorphisms}, 
Bull. Symbolic Logic {\bf 15:2} (2009), 184--214.

\bibitem[Ro13]{Ro13} ---,  {\it Global and local boundedness 
of Polish groups}, Indiana Univ. Math. J. {\bf 62:5} (2013), 1621--1678

\bibitem[Ro18]{Ro18} ---, ``Coarse Geometry of Topological Groups,'' 
http://homepages.math.uic.edu/~rosendal/PapersWebsite/Coarse-Geometry-Book17.pdf

\bibitem[Ro19]{Ro19} ---,  {\it Continuity of universally measurable 
homomorphisms}, Forum Math. Pi {\bf 7} (2019), e5, 20pp. 

\bibitem[Ru87]{Ru87} Rudin, W. , ``Real and Complex Analysis'', 
McGraw-Hill, New York, 1987.
%
\bibitem[Ru91]{Ru91} ---, ``Functional Analysis'', McGraw-Hill, New York, 1991.
%
\bibitem[Sa53]{Sa53} Samelson, H., {\it A connection between the Whitehead product 
and the Pontryagin product}, Amer.\ J. Math.\ {\bf 75} (1953), 744--752. 
%
\bibitem[SN81]{SN81}
Saxon, S. A., and P.~P. Narayanaswami,
{\it Metrizable} (LF)-{\it spaces}, (db)-{\it spaces, and the separable quotient problem},
Bull.\ Aust.\ Math.\ Soc.\ {\bf 23} (1981), 65--80.
%
\bibitem[ShW99]{ShW99}
Schaefer, H.~H., and M.~P. Wolff, ``Topological Vector
Spaces'', Springer-Verlag, Berlin, 1999.
%
\bibitem[Sr29]{Sr29}
Schauder, J., {\it Invarianz des Gebiets in Funktionalr\"{a}umen},
Studia Math.\ {\bf 1} (1929), 123--139.
%
\bibitem[Scr97]{Scr97}
Schechter, E.,
``Handbook of Analysis and its Foundations,''
Academic Press, San Diego, 1997.
%
\bibitem[Sf84]{Sf84} Schikhof, W.~H., ``Ultrametric Calculus -- An Introduction to
$p$-adic Analysis'', Cambridge University Press, Cambridge, 1984. 
%
\bibitem[Sme15]{Sme15}
Schmeding, A.,
``The Diffeomorphism Group of a Non-Compact Orbifold,''
Diss.\ Math.\ {\bf 507} (2015), 179 pp.
%
\bibitem[Sme23]{Sme23}
---, ``An Introduction to Infinite-Dimensional DIfferential Geometry,''
Cambridge University Press, Cambridge, 2023.
%
\bibitem[Sch87]{Sch87} Schmid, R., ``Infinite-dimensional Hamiltonian systems,'' 
Monographs and Textbooks in Physical Science, Lecture Notes {\bf 3}, 
Bibliopolis, Naples, 1987.
%
\bibitem[Sch02]{Sch02} ---, {\it The Lie group of Fourier integral 
operators on open manifolds}, 
in ``Infinite Dimensional Lie Groups in Geometry and Representation 
Theory'', Banyaga, A., J. A. Leslie, and T. Robart, eds., 
World Scientific, River Edge, NJ, 2002; 44--64. 

\bibitem[Sch04]{Sch04} ---, {\it Infinite dimensional Lie groups with applications 
to mathematical physics}, J. Geom.\ Symm.\ Phys.\ {\bf 1} (2004), 54--120.
%
\bibitem[SAR84]{SAR84} Schmid, R., M. Adams, and T. Ratiu, {\it The group of Fourier 
integral operators as symmetry group},
in ``XIIIth international colloquium 
on group theoretical methods in physics'' (College Park, Md., 1984), 
World Sci.\ Publishing, Singapore, 1984, 246--249. 
%
\bibitem[Sch75]{Sch75} 
Schubert, H., ``Topologie'', Teubner Verlag, Stuttgart, 1975. 
%
\bibitem[Sc60]{Sc60} Schue, J.~R., {\it Hilbert space methods in the theory of Lie
algebras}, Trans.\ Amer.\ Math.\ Soc.\ {\bf 95} (1960),
69--80.
%
\bibitem[Sc61]{Sc61} ---, {\it Cartan decompositions for $L^{*}$-algebras}, Trans.\
Amer.\ Math.\ Soc.\ {\bf 98} (1961), 334--349.
%
\bibitem[SJ13]{SJ13}
Sch\"{u}tt, J., {\it Symmetry groups of principal bundles over non-compact bases},
preprint, 2013, arXiv:1310.8538
%
\bibitem[SV17]{SV17}  Schulz-Baldes, H., and C.~Villegas-Blas, {\it 
Signatures for J-hermitians and J-unitaries on Krein spaces 
with real structures}, Math. Nachr. {\bf 290:11-12} (2017),  1840--1858

\bibitem[ScF90a]{SchF90a} Schur, F., {\it Neue Begr\"undung der 
Theorie der endlichen Transformationsgruppen}, 
Math.\ Ann.\ {\bf 35} (1890), 161--197. 
%
\bibitem[ScF90b]{SchF90b} ---,  {\it Beweis f\"ur die Darstellbarkeit 
der infinitesimalen Transformationen aller transitiven endlichen Gruppen durch 
Quotienten best\"andig convergenter Potenz\-rei\-hen}, Leipz.\ Ber.\ {\bf 42} (1890), 1--7. 
%
\bibitem[Schw57]{Schw57}
Schwartz, L., ``Th\'{e}orie des Distributions'', Tome~I. Hermann, Paris, 1957. 
%
\bibitem[SeS55]{SeS55}
Sebasti\~{a}o e Silva, J.,
{\it Su certe classi di spazi localmente convessi importanti per le applicazioni},
Rend.\ Mat.\ Appl.\ {\bf 14} (1955), 388--410.
%
\bibitem[See64]{See64}
Seeley, R.,
{\it Extension of $C^\infty$ functions defined in a half space},
Proc.\ Am.\ Math.\ Soc.\ {\textbf 15} (1964), 625--626.
%
\bibitem[SeG68]{SeG68}
Segal, G., {\it Classifying spaces and spectral sequences}, 
Publ.\ Math.\ I.H.E.S. Paris \textbf{34} (1968), 105--112.
%
\bibitem[SeG70]{SeG70}
---, {\it Cohomology of topological groups},  
Sympos.\ Math.\ Roma,  \textbf{4},
Teoria Numeri Dic.\ 1968, e Algebra, Marzo 1969,
\textbf{1970}, 377--387. 
%
\bibitem[SeG75]{SeG75}
---, {\it Classifying space of a topological group in the
Gelfand-Fuks sense}, Funkt.\ Anal.\ Prilozhen \textbf{9:2} (1975), 131--133. 
%
\bibitem[SeG81]{SeG81} ---, {\it Unitary representations of some
infinite-dimensional groups}, Comm.\ Math.\ Phys.\ \textbf{80} (1981), 301--342.
%
\bibitem[SeI78]{SeI78}
Segal, I.~E., {\it The complex-wave representation of the free
boson field}, pp.\ 321--343 in:
``Topics in Funct.\ Anal.'', Adv.\ in
Math.\ Suppl.\ Stud.\ \textbf{3} (1978), Academic Press, New York
%
\bibitem[Sei72]{Sei72} Seip, U., ``Kompakt erzeugte Vektorr\"{a}ume und Analysis",
Lecture Notes in Math.\ {\bf 273}, Springer, Berlin, 1972.
%
\bibitem[Sr92]{Sr92}
Serre, J.-P., ``Lie Algebras and Lie Groups'', 
Lecture Notes Math.\ {\bf 1500}, Springer-Verlag, Berlin, 1992.
%
\bibitem[Sh49]{Sh49}
Shapiro, A., {\it Group extensions of compact Lie groups}, 
Ann.\ of Math.\ \textbf{50:3} (1949), 581--586.
%
\bibitem[Sh77]{Sh77} Shereshevskij, I.~A., 
{\it Quantization based on infinite-dimensional
Hermitian symmetric spaces}, 
Moscow Univ.\ Math.\ Bull.\ \textbf{32} (1977), 28--36. 
%
\bibitem[Si52]{Si52}
Singer, I.~M., {\it Uniformly continuous representations of Lie
groups}, Annals Math.\ \textbf{56:2} (1952), 242--247.
%
\bibitem[SiSt65]{SiSt65}
Singer, I. M., and S. Sternberg, {\it The infinite groups of Lie and Cartan}. I. {\it The transitive 
groups},
J. d'Anal.\ Math.\ {\bf 15} (1965), 1--114. 
%
\bibitem[Sz88]{Sz88}
Slez\'{a}k, B., {\it On the inverse function theorem
and implicit function theorem in Banach spaces}, 
pp.~186--190 in: J. Musielak (ed) ``Function Spaces'', (Pozna\'{n}, 1986).
Teubner-Verlag, Leipzig, 1988. 
%
\bibitem[Sl85]{Sl85} Slodowy, P., {\it An adjoint quotient for certain groups
attached to Kac--Moody algebras}, pp.\ 307--333  
in: V. Kac (ed),
``Infinite Dimensional groups with Applications'', 
Math.\ Sciences Research Institute Publications \textbf{4}, 
Springer-Ver\-lag, Berlin, Heidelberg, New York, 1985.
%
\bibitem[Sl86]{Sl86} ---, {\it Beyond Kac-Moody algebras, and inside}, 
Canad.\ Math.\ Soc.\ Conf.\ Proceedings,
\textbf{5} (1986), 361--371. 
%
\bibitem[Son98]{Son98}
Sontag, E.~D.,
``Mathematical Control Theory. Deterministic Finite Dimensional Systems,''
Texts in Applied Math.\ {\bf 6}, Springer Verlag, New York, ${}^2$ 1998.
%
\bibitem[So84]{So84}
Souriau, J.-M., {\it Groupes diff\'{e}rentiels
de physique math\'{e}matique}, pp.~73--119 in:
P. Dazord and N. Desolneux-Moulis (eds)
``Feuilletages et quantification g\'{e}ometrique,
Journ. lyonnaises Soc.\ math.\ France 1983,
S\'{e}min.\
sud-rhodanien de g\'{e}om.\ II'', 
Hermann, Paris, 1984. 
%
\bibitem[So85]{So85}
---, {\it Un algorithme g\'en\'erateur de structures
quantiques}, Soc.\ Math.\ Fr., Ast\'erisque, hors s\'erie 
 (1985), 341--399.  
%
\bibitem[Sp66]{Sp66}
Spanier, E. H., ``Algebraic Topology'',  
Mc Graw Hill Series in Higher Math., Univ.\ of Cal., Berkeley, 
1966. 
%
\bibitem[Sp17]{Sp17} Spichak, S.~V., {\it Classification of 
realizations of Lie algebras of vector fields on a circle}, 
arXiv:math.RT:1304.2241V2  26 Sep 2017 
%
\bibitem[Ste51]{Ste51} 
Steenrod, N.~E., ``The Topology of Fibre Bundles'',  Princeton
University Press, Princeton, New Jersey 1951. 
%
\bibitem[Ste67]{Ste67} ---,
{\it A convenient category of topological spaces},
Mich.\ Math.\ J. {\bf 14} (1967), 133--152.
%
\bibitem[Stn70]{Stn70}
Stein, E. M.,
\emph{Singular integrals and differentiability properties of functions},
Princeton Univ.\ Press, Princeton, 1970.
%
\bibitem[St61]{St61} Sternberg, S., {\it Infinite Lie groups and the formal aspects 
of dynamical systems}, J. Math.\ Mech.\ {\bf 10} (1961), 451--474. 
%
\bibitem[Str06]{Str06} Stroppel, M., ``Locally Compact Groups,'' 
EMS Textbooks in Mathematics,  European Mathematical Society (EMS), 
Z\"urich, 2006. 
%
\bibitem[St99]{St99} Stumme, N., ``The Structure of Locally Finite Split Lie 
algebras,'' Ph.\ D.\ thesis, Darmstadt University of Technology, 1999. 
%
\bibitem[Sus73]{Sus73} Sussmann, H. J., {\it Orbits of families of vector fields and 
integrability of distributions}, Trans.\ Amer.\ Math.\ Soc.\ {\bf 180} (1973), 171--188.
%
\bibitem[Su88]{Su88} Suto, K., 
{\it Groups associated with unitary forms of
Kac--Moody algebras}, 
J.\ Math.\ Soc.\ Japan \textbf{40:1}  (1988), 85--104. 
%
\bibitem[Su97]{Su97}
---, {\it Borel-Weil-type Theorem for the flag manifold of
a generalized Kac-Moody algebra}, 
J. Algebra \textbf{193}  (1997), 529--551. 
%
\bibitem[Swa62]{Swa62} Swan, R. G., {\it Vector bundles and projective modules}, 
Trans.\ Amer.\ Math.\ Soc.\ {\bf 105} (1962), 264--277.
%
\bibitem[Swa77]{Swa77} ---, {\it Topological examples of
projective modules}, 
Trans.\ Amer.\ Math.\ Soc.\ \textbf{230}  (1977), 201--234. 
%
\bibitem[Sw65]{Sw65} Swierczkowski, S., {\it  Embedding theorems for local analytic groups}, 
Acta Math.\ {\bf 114} (1965), 207--235. 
%
\bibitem[Sw67]{Sw67}
---, {\it Cohomology of local groups extensions}, 
Trans.\ Amer.\ Math.\ Soc.\ \textbf{128}  (1967), 291--320. 
%
\bibitem[Sw71]{Sw71} ---, {\it The path-functor on Banach Lie
algebras}, 
Nederl.\ Akad.\ Wet., Proc.\ {\bf A 74}=Indag.\ Math.\ {\bf 33} (1971), 235--239. 
%
\bibitem[Swi71b]{Swi71b} ---, {\it On the cohomology of local groups}, J. Austral.\ Math.\ Soc.\ 
{\bf 12} (1971), 249--255. 
%
\bibitem[Swi71c]{Swi71c} ---, {\it Cohomology of group germs and Lie algebras},
Pac.\ J. Math.\ {\bf 39} (1971), 471--482. 
%
\bibitem[Ta68]{Ta68} Tagnoli, A., {\it La variet\`a analitiche reali come 
spazi omogenei}, Boll.\ Un.\ Mat.\ Ital.\ (4) {\bf 1} (1968), 422--426.
%
\bibitem[Ta14]{Ta14} Tao, T., ``Hilbert's Fifth Problem and Related Topics,'' 
Graduate Studies in Math.\ {\bf 153}, Amer.\ Math.\ Soc., 2014.
%
\bibitem[Ta15]{Ta15} ---, ``Expansion in Finite Simple Groups 
of Lie Type,'' Graduate Studies in Math.\ {\bf 164}, Amer.\ Math.\ Soc., 2015.
%
\bibitem[TSH98]{TSH98}
Tatsuuma, N., H. Shimomura, and T. Hirai, 
{\it On group topologies and unitary representations of inductive limits of
topological groups and the case of the group of diffeomorphisms}, 
J. Math.\ Kyoto Univ.\ \textbf{38} (1998), 551--578.
%
\bibitem[Te99]{Te99} Teichmann, J., 
``Infinite Dimensional Lie Theory from the Point
of View of Functional Analysis'',  Ph.D. Thesis, Vienna, 1999. 
%
\bibitem[Te01a]{Te01a} ---, {\it Regularity of infinite-dimensional Lie groups by metric space
methods}, Tokyo J. Math.\ {\bf 24} (2001), 39--58. 
\bibitem[Te01b]{Te01b}
---, {\it A Frobenius theorem on convenient manifolds},
Monatsh.\ Math.\ {\bf 134:2} (2001), 159--167.
%
\bibitem[Te02]{Te02}
---, {\it Inheritance properties for Lipschitz-metrizable Fr\"olicher 
groups}, in ``Infinite Dimensional Lie Groups in Geometry and Representation 
Theory'', Banyaga, A., J. A. Leslie, and T. Robart, eds., 
World Scientific, 2002; 1--10  

\bibitem[Th95]{Th95}
Thomas, E.~G.~F., {\it Vector fields as derivations on nuclear
manifolds}, Math.\ Nachr.\ \textbf{176} (1995), 277--286.
%
\bibitem[Th96]{Th96}
---, {\it Calculus on locally convex spaces}, 
Preprint \textbf{W-9604}, Univ.\ of Groningen, 1996. 
%
\bibitem[Th74]{Th74} Thurston, W. P., {\it Foliations and groups of diffeomorphisms},
Bull.\ Amer.\ Math.\ Soc.\ {\bf 80} (1974), 304--307.
%
\bibitem[Tn79]{Tn79}
Tidten, M., \emph{Fortsetzungen von $C^\infty$-Funktionen, welche auf einer abgeschlossenen
Menge in $\R^n$
definiert sind},
Manuscr.\ Math.\ {\bf 27} (1979), 291--312.
%
\bibitem[Ti83]{Ti83}
Tits, J., ``Liesche Gruppen und Algebren'', 
Springer-Verlag, Berlin, New York, Heidelberg, 1983. 
%
\bibitem[Tr67]{Tr67}
Treves, F., ``Topological Vector Spaces, Distributions, and
Kernels'', Academic Press, New York, 1967. 
%
\bibitem[TL99]{TL99} Toledano Laredo, V.,
{\it Integrating unitary representations of infinite-dimensional Lie groups}, 
J. Funct.\ Anal.\ \textbf{161} (1999), 478--508. 
%
\bibitem[tD91]{tD91}
tom Dieck, T., ``Topologie'', de Gruyter, Berlin, New York, 1991. 
%
\bibitem[tD00]{tD00}
---, ``Topologie'', 2nd ed.,  de Gruyter, Berlin, New York, 2000. 
%
\bibitem[Tr72]{Tr72}
Tromba, A.~J., {\it Some theorems on Fredholm maps}, Proc.\ Amer.\ Math.\ Soc.\
(1972), 578--585.
%
\bibitem[Tsu81]{Tsu81} Tsujishita, T., ``Continuous cohomology of the Lie algebra of vector fields,'' 
Memoirs Amer.\ Math.\ Soc.\ {\bf 253}, 1981.
%
\bibitem[Tu85]{Tu85} Turovskii, Yu. V., {\it Spectral properties of certain Lie 
subalgebras and the spectral radius of subsets of a Banach algebra} (Russian), 
in Spectral theory of operators and its applications,  
``\`Elm'', Baku,  {\bf 6} (1985), 144--181 
%
\bibitem[Tu70]{Tu70} Turpin, P., 
{\it Une remarque sur les alg\`{e}bres \'{a} inverse continu}, 
C.~R.\ Acad.\ Sci.\ Paris \textbf{270} (1970), 1686--1689. 
%
\bibitem[Tu95]{Tu95} Tuynman, G.~M. 
{\it An elementary proof of Lie's Third Theorem}, Unpublished note, 
1995. 
%
\bibitem[TW87]{TW87} Tuynman, G.~M., and W.~A.~J.~J. Wiegerinck, 
{\it Central extensions and physics}, J. Geom.\ Physics \textbf{4:2} (1987), 207--258.
%
\bibitem[Up85]{Up85}
Upmeier, H., ``Symmetric Banach Manifolds and Jordan
$C^*$-algebras'', North-Holland Math.\ Studies \textbf{104},
North-Holland, Amsterdam, 1985.
%
\bibitem[Va84]{Va84} Varadarajan, V. S., ``Lie Groups, Lie Algebras, and Their Representations,'' 
Graduate Texts in Math.\ {\bf 102}, Springer--Verlag, 1984. 
%
\bibitem[Va85]{Va85}
---, ``Geometry of Quantum Theory'',  
Springer-Verlag, New York, ${}^2$1985. 
%
\bibitem[Ver82]{Ver82}
Verdier, J.-L., {\it Les representations des algebres de Lie
affines: applications a quelques problemes de physique}, 
Sem.\ Bourbaki \textbf{596} (1982), 365--377.  
%
\bibitem[VK82]{VK82}
Ver\v s{}ik, A.\ M., and S.~V. Kerov, {\it Characters and
factor representations of the infinite unitary group}, 
Soviet Math.\ Dokl.\ \textbf{26} (1982), 570--574.
%
\bibitem[Ves68]{Ves68} Vesentini, E., {\it On the subharmonicity of the spectral 
radius}, Boll.\ Un.\ Mat.\ Ital.\ {\bf 1} (1968), 427--429.
%
\bibitem[Vi02]{Vi02} 
Vizman, C., {\it Geodesics on extensions of the Lie algebra of
vector fields on the circle}, pp.~165--172 in:  A. Strasburger et al.\ (eds), 
``Geometry and Analysis on Finite-and Infinite-Dimensional Lie Groups'',
Banach Center Publications {\bf 55}, Warszawa. 
%
\bibitem[Vo87]{Vo87} Vogt, D., {\it On the functors $\Ext^1(E,F)$ for Fr\'echet spaces}, Studia 
Math.\ {\bf 85} (1987), 163--197. 
%
\bibitem[Vo76]{Vo76}
Voiculescu, D., {\it Repr\'esentations factorielles de type
$II_1$ de $U(\infty)$}, J.~Math.\ pures appl.\ \textbf{55} (1976), 1--20.
%
\bibitem[Voi92]{Voi92}
Voigt, J.,
\emph{On the convex compactness property for the strong operator topology},
Note Mat.\ {\bf 12}, 259--269.
%
\bibitem[Wa54a]{Wa54a}
Waelbroeck, L., {\it Les alg\`{e}bres \`{a} inverse continu}, 
C.~R. Acad.\ Sci.\ Paris \textbf{238} (1954), 640--641.
%
\bibitem[Wa54b]{Wa54b}
---, {\it Le calcul symbolique dans les
alg\`{e}bres commutatives}, J. Math.\ 
Pures Appl.\ \textbf{33} (1954), 147--186.
%
\bibitem[Wa54c]{Wa54c} ---, {\it Structure des alg\`{e}bres \`{a} inverse continu}, 
C.~R. Acad.\ Sci.\ Paris {\bf 238} (1954), 762--764. 
%
\bibitem[Wa71]{Wa71} ---, ``Topological Vector Spaces and Algebras'',
Springer-Verlag, Berlin, Heidelberg, New York, 1971. 
%
\bibitem[Wa75]{Wa75}
---, {\it The holomorphic functional calculus
and non-Banach algebras}, pp.\ 187--251 in: J.~H. Williamson (ed), 
``Algebras in Analysis'', Academic Press, London, New York, 1975.  
%
\bibitem[Wa82]{Wa82}
---, {\it The holomorphic functional calculus
as an operational calculus}, pp.\ 512--552
in:
``Spectral Theory'', Banach Center Publications,
\textbf{8}, Warsaw, 1982. 
%
\bibitem[Wa99]{Wa99} 
Wagemann, F., ``Sur la cohomologie de Gelfand-Fuks des champs
de vecteurs holomorphes'', Ph.\ D. Thesis, Univ.\ Claude Bernard--Lyon~I, 
1999.
%
\bibitem[Wa12]{Wa12}
Walter, B., \emph{Weighted diffeomorphism groups of Banach spaces and weighted mapping groups},
Diss.\ Math.\ {\bf 484}, 2012.
%
\bibitem[WHO18]{WHO18}
Walther, H.~O.,
{\it Fr\'{e}chet differentiability in Fr\'{e}chet spaces, and differential equations with unbounded variable delay},
preprint, 2018, arXiv:1801.09213.
%
\bibitem[Wa72]{Wa72} Warner, G., ``Harmonic Analysis on Semisimple Lie Groups I'',  
Springer-Verlag, Berlin, Heidelberg, New York, 1972. 
%
\bibitem[We95]{We95} Weibel, C.~A.,
``An Introduction to Homological Algebra'', 
Cambridge Studies in Advanced Math.\ \textbf{38}, Cambridge Univ.\ Press,
Cambridge, 1995.
%
\bibitem[Wei69]{Wei69} Weinstein, A., {\it  Symplectic structures on Banach manifolds}, 
Bull.\ Amer.\ Math.\ Soc.\ {\bf 75} (1969), 1040--1041.
%
\bibitem[Wei89]{Wei89} ---, {\it Cohomology of symplectomorphism groups and critical values 
of hamiltonians}, Math.\ Z. {\bf 201} (1989), 75--82. 
%
\bibitem[Wsg52]{Wsg52}
Weissinger, J., \emph{Zur Theorie und Anwendung des Iterationsverfahrens},
Math.\ Nachr.\ {\bf 8} (1952), 193--212.
%
\bibitem[Wz12]{Wz12}
Weizs\"{a}cker, H. von,
\emph{In which spaces is every curve Lebesgue--Pettis integrable}?,
preprint, 2012, {\tt arXiv:1207.6034}.
%
\bibitem[We76]{We76}
Wells, J.~C.,
\emph{Invariant manifolds of nonlinear operators},
Pac.\ J. Math.\ {\bf 62} (1976), 285--293.
%
\bibitem[We03]{We03}
Wengenroth, J., ``Derived Functors in Functional Analysis,''
Lecture Notes in Mathematics~{\bf 1810}, Springer-Verlag, Berlin, 2003.
%
\bibitem[Wer95]{Wer95}
Werner, D., ``Funktionalanalysis'', Springer-Verlag, Berlin, Heidelberg, 1995. 
%
\bibitem[Whi78]{Whi78} Whitehead, G. W., ``Elements of Homotopy Theory,'' 
Graduate Texts in Mathematics {\bf 61}, Springer-Verlag, New York-Berlin, 1978.
%
\bibitem[Wh34]{Wh34}
Whitney, H.,
\emph{Analytic extensions of differentiable functions defined in closed sets},
Trans.\ Am.\ Math.\ Soc.\ {\bf 36} (1934), 63--89.
%
\bibitem[Wh36]{Wh36} ---,
\emph{Differentiable manifolds},
Ann.\ Math.\
{\bf 37} (1936), 645--680.
%
%
%
\bibitem[Wi88]{Wi88}
Wi\c{e}s\l{}aw, W., ``Topological Fields'', Marcel Dekker, New York and Basel, 1988.
%
\bibitem[Wie49]{Wie49} Wielandt, H., {\it \"Uber die Unbeschr\"anktheit der Operatoren 
der Quantenmechanik},  Math.\ Ann.\ {\bf 121} (1949), 21.
%
\bibitem[Wo03]{Wo03} Wockel, Chr., ``Central Extensions of Gauge Groups,'' 
Thesis (Dipl.), TU Darmstadt, Sept.\ 2003.  
%
%
\bibitem[Wo06]{Wo06}
---, {\it Smooth extensions and spaces of smooth and holomorphic 
mappings},
J. Geom.\ Symmetry Phys.\ {\bf 5} (2006), 118--126.
%
\bibitem[Wo07a]{Wo07a}
---, \emph{Lie group structures on symmetry groups of principal bundles},
J. Funct.\ Anal.\ {\bf 251} (2007), 254--288. 
%
\bibitem[Wo07b]{Wo07b}
---, {\it The Samelson product and rational homotopy for gauge groups},
Abh.\ Math.\ Semin.\ Univ.\ Hamb.\ {\bf 77} (2007), 219--228.
%
\bibitem[Wo09]{Wo09}
---, {\it A generalization of Steenrod's approximation theorem},
Arch.\ Math.\ Brno {\bf 45:2} (2009), 95--104.
%
\bibitem[Wol05]{Wol05} Wolf, J. A., {\it Principal series representations
of direct limit groups}, Compos.\ Math.\ {\bf 141:6} (2005), 1504--1530.
%
\bibitem[Woj98]{Woj98} Wojty\'nski, W., {\it Quasinilpotent Banach--Lie algebras are 
Baker--Campbell--Hausdorff}, J.~Funct.\ Anal.\ {\bf 153} (1998), 405--413.  
%
\bibitem[Woj06]{Woj06} ---, {\it Effective integration of Lie algebras}, 
J. Lie Theory {\bf 16} (2006), 601--620. 
%
\bibitem[Wue03]{Wue03} W\"ustner, M., {\it Supplements on the theory of exponential Lie groups}, 
J. Algebra {\bf 265:1} (2003), 148--170.
%
\bibitem[Wue05]{Wue05} ---, {\it The classification of all simple Lie groups with surjective 
exponential map},  J. Lie Theory {\bf 15:1} (2005), 269--278.
%
\bibitem[Wur01]{Wur01} Wurzbacher, T., {\it Fermionic second quantization and the
geometry of the restricted Grassmannian}, pp.~287--375
in: A.~Huckleberry and T.~Wurzbacher (eds),
``Infinite Dimensional K\"ahler Manifolds,'' 
DMV-Seminar {\bf 31}, 
Birkh\"auser Verlag, 2001.
%
\bibitem[Yag01]{Yag01} Yagasaki, T., {\it 
Homotopy types of diffeomorphism groups of noncompact $2$-manifolds}, 
arXiv:math.GT:0109183

\bibitem[Yag10]{Yag10} Yagasaki, T., {\it 
Groups of volume-preserving diffeomorphisms of noncompact manifolds and
mass flow toward ends}, 
Trans. Amer. Math. Soc. {\bf 362:11} (2010), 5745--5770.

\bibitem[Ya85]{Ya85} Yamasaki, Y., `` Measures on Infinite Dimensional
Spaces'',  Series in Pure Math., {\bf 5}, World Scientific, Singapore, 1985. 
%
\bibitem[Ym98]{Ym98} Yamasaki, A., {\it Inductive limits of general linear
groups}, J. Math.\ Kyoto Univ.\ \textbf{38} (1998), 769--779.
%
\bibitem[Yo36]{Yo36} Yosida, K., {\it On the groups embedded in the metrical complete ring}, 
Japanese J. Math.\ {\bf 13} (1936), 7-26.
%
\bibitem[Yo63]{Yo63} Yosida, K., {\it Holomorphic semi-groups in a locally convex 
linear topological space}, Osaka Math.\ J. {\bf 15} (1963), 51--57.
%
\bibitem[Ze60]{Ze60} \.{Z}elazko, W., 
{\it On the locally bounded and m-convex topological algebras}, 
Studia Math.\ \textbf{19} (1960), 333--356.
%
\bibitem[Ze85]{Ze85} ---, {\it A non-m-convex
algebra on which operate all entire functions}, 
Ann.\ Pol.\ Math.\ \textbf{46}  (1985), 389--394. 
%
\end{thebibliography}
\end{document}